\documentclass[12pt]{book}
\usepackage[margin=0.8in]{geometry}

\usepackage[htt]{hyphenat}

\newcommand{\R}{\mathbb{R}}
\newcommand{\Z}{\mathbb{Z}}
\newcommand{\N}{\mathbb{N}}
\newcommand{\Q}{\mathbb{Q}}
\newcommand{\cOMPL}{\mathbb{C}}

\renewcommand{\epsilon}{\varepsilon}

\newcommand{\e}{\varepsilon}

\usepackage[ukrainian, english]{babel}
\usepackage[alphabetic]{amsrefs}
\usepackage{amsmath,amssymb,amsfonts,amsthm,enumerate}
\usepackage{mathtools}
\usepackage{hyperref}
\usepackage{graphicx,epstopdf,color}
\usepackage{newtxtext}
\usepackage[varvw]{newtxmath}
\usepackage{calligra}
\usepackage[T1]{fontenc}

\numberwithin{equation}{section}
\numberwithin{figure}{section}

\newtheorem{theorem}{Theorem}[section]
\newtheorem{lemma}[theorem]{Lemma}
\newtheorem{definition}[theorem]{Definition}

\newtheorem{proposition}[theorem]{Proposition}
\newtheorem{corollary}[theorem]{Corollary}

\renewcommand{\leq}{\leqslant}
\renewcommand{\le}{\leqslant}
\renewcommand{\geq}{\geqslant}
\renewcommand{\ge}{\geqslant}

\DeclareMathOperator{\DIV}{div}
\DeclareMathOperator{\curl}{curl}
\renewcommand{\div}{\DIV}
\DeclareMathOperator\cAPAC{Cap}
\DeclareMathOperator{\sign}{sign}
\allowdisplaybreaks

\usepackage{imakeidx}
\makeindex

\begin{document}

\author{Serena Dipierro \& Enrico Valdinoci\thanks{University of Western Australia,
35 Stirling Highway,
Crawley WA 6009, Australia. {\tt serena.dipierro@uwa.edu.au}, 
{\tt enrico.valdinoci@uwa.edu.au}}}

\title{Elliptic partial differential equations \\
from an elementary viewpoint}
\maketitle

\tableofcontents

\chapter*{Preface}

\begin{center}
\begin{quote}
{\color{blue}\large\calligra{``There\raisebox{-1.9ex}{\kern-.35em'\kern.35em}s
more to life than mathematics'',
Joan said. ``But not much more''.\\

{Greg Egan}, Glory.}}
\end{quote}
\end{center}\bigskip\bigskip\bigskip\bigskip

These notes are the outcome of some courses taught to undergraduate and graduate students
from the University of Western Australia,
the Pontif\'{\i}cia Universidade Cat\'olica do Rio de Janeiro, 
the Indian Institute of Technology Gandhinagar, the
Ukrainian Catholic University ({\foreignlanguage{ukrainian}{Український Католицький Університет}}) and the Politecnico di Milano in 2021--2023.
\medskip

Far from aiming at being all-encompassing, the following pages wish to
shed some light on a number of selected topics in the theory
of elliptic partial differential equations with a style that should be
accessible to third-year undergraduate students, possibly under an
inspired mentorship, but might
also provide some interest to more advanced students
and possibly professional researchers. While all the
topics presented are of classical flavor,
the exposition and the way the material are organized is
perhaps rather original, with the intention of addressing several quite
difficult points with a style that is as self-contained as possible,
rigorous as well as intuitive,
and approachable without major prerequisites (indeed,
we only assume prior knowledge of the ``basic'' analysis, giving also
reference to useful results and theorems whenever the less advanced readers
may need to consolidate their backgrounds).
\medskip

We stress again that the list of topics covered here is far from being exhaustive,
since we mainly dealt with classical problems related to the Laplace operator
and their natural counterpart for equations in nondivergence form
(the divergence structure case is not really covered here,
though the Divergence Theorem is obviously utilized ubiquitously, the variational structure of the Laplace-Beltrami
operator is discussed quite in detail
and several pages are devoted to the theory of capacity which possesses a variational essence).
Also, we do not address here nonlinear, singular/degenerate elliptic operators, fully nonlinear equations
and fractional/nonlocal elliptic equations.
We do not even linger too much on explicit solutions
and cheap tricks to find them (e.g., we do not repeat over and over very specific methods such as the separation of variables),
since, all things considered, it is very unlikely that one can solve explicitly\footnote{Moreover, explicit solutions
may provide a handy resource to develop an initial theory, but they usually do not exhaust the complexity
of the problem. For example, according to~\cite[page~452]{MR1908418},
it is quite common to run into ``one of the serious problems with such exact solutions [...]: namely, they often do not determine all possible solutions
and indeed, may not even give the most relevant one''.}
a partial differential equation,
instead it is often more useful to understand the qualitative and quantitative properties
of the solution without solving explicitly the corresponding equation.

The many parts of the theory of elliptic partial differential equations which are missing from these pages
are by no means less important than the topics covered: the issue is that we had to make
a rather harsh selection of topics just to be able to collect all the material into finitely many pages
(arguably, Jorge Luis Borges would remark that a comprehensive treatment of elliptic equations
can be found in the Library of Babel~\cite{BORGES}).
\medskip

Some additional fun for the interested reader could come from
the treatment of a few natural problems with geometric flavor, such as
the Soap Bubble Theorem and Serrin's overdetermined problem. To treat them,
but this is also a general approach used everywhere in the book,
we tried to provide different perspectives and different
points of view, also providing physical motivations
to trigger creativity and imagination, and favoring a dynamic interplay between geometric and analytic arguments.\medskip

These notes also collect a number of classical, albeit not always well-known,
topics, which we believe can serve as an excellent training camp to develop
some familiarity with elliptic equations.\medskip

Here are
a couple of features of these notes that we hope can be appealing for many readers. First, all is built in an artisanal way, perhaps not by following the most elegant, concise or general approach, but aiming at a possibly slow, but conceptually clear, sequence of strategic steps. Second, this set of notes is planned in a way that one can read starting almost any place and freely jump from one topic to another just by following a personal stream of thoughts.
\medskip

Though a list of classical and modern references is given at the end,
this is not a book focused on the history
of mathematics, hence we will not address topics such as priority
of discoveries or progress of subject over the course of time.
\medskip

As usual, these notes may contain errors or inconsistencies: if you find
any, please let us know. In general, we will be happy to receive
comments and criticisms to possibly improve this work.
After all, scientific knowledge is based on a dynamic flow
of information and we
will certainly cherish and treasure readers' feedback and advice.
\medskip

All right, enough chitchat, please fasten your PFD,
it's time to start our journey together.

\begin{center}
Serena and Enrico
\end{center}

\vfill

\begin{figure}[h]
  \centering
  \includegraphics[width=.65\linewidth]{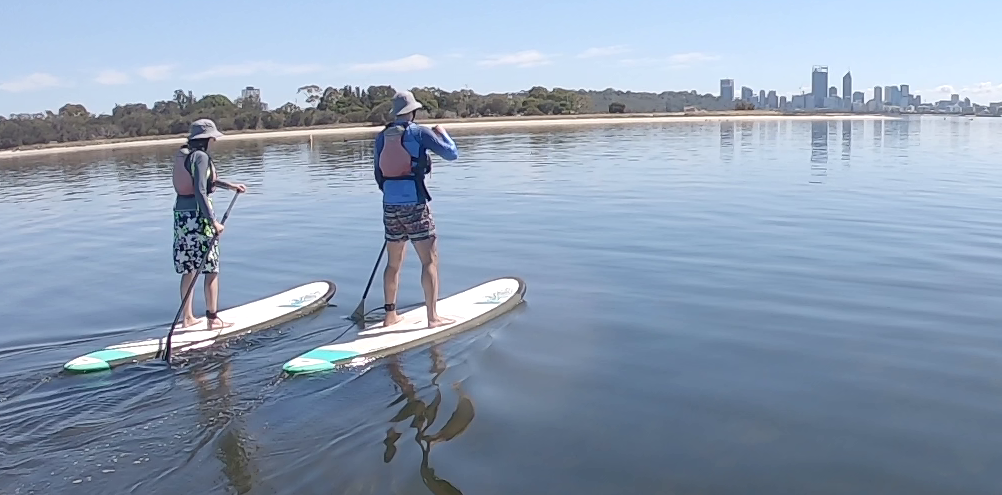}
  \caption{\sl Courtesy of Burbuqe Shaqiri.}
\end{figure}

\chapter{Motivations}

\section{The Laplacian coming onto the scene}

Without aiming at being exhaustive, we present here some classical circumstances in which partial differential
equations, and especially (but not only)
those relying on the Laplace operator, naturally arise to model interesting natural phenomena.
The objective is not much to use mathematics  to find the Answer to the Ultimate Question of Life, the Universe, and Everything
(which is well-known to be~$42$ anyway) but rather to develop some familiarity with some of the basic features
of partial differential equations, to see them in action and to build an intuition and an instinctive feeling
of the problems,
which will turn out handy in our voyager through the mighty jungle of technical and seemingly abstract mathematics.

\subsection{The heat equation}\label{HEATSEC:02w}

\begin{figure}
                \centering
                \includegraphics[width=.65\linewidth]{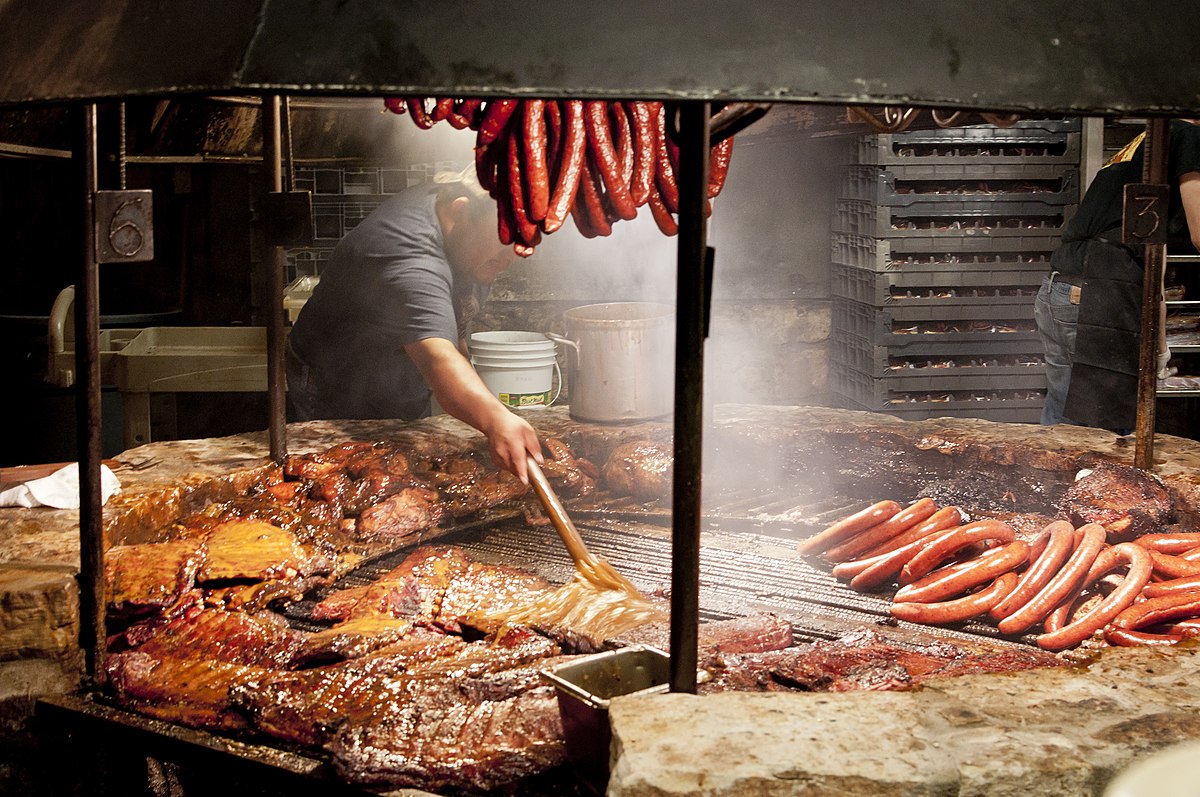}
        \caption{\sl The heat equation applied in a restaurant in Texas (photo by vxla, image from
        Wikipedia,
        licensed under the Creative Commons Attribution 2.0 Generic license).}\label{HAFOBAR89GIJ7solFUMHDNOJHNFOJED}
\end{figure}

Among the many occurrences in which elliptic equations naturally arise in nature, one of the most widely popularized comes from the stationary state of the heat equation (see Figure~\ref{HAFOBAR89GIJ7solFUMHDNOJHNFOJED} for a popular application of the heat equation).

\begin{figure}
                \centering
                \includegraphics[width=.35\linewidth]{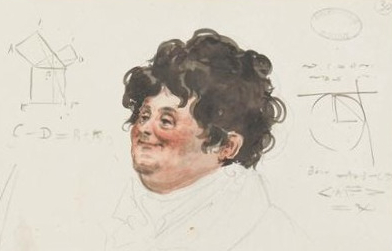}
        \caption{\sl Caricature of Joseph Fourier (Public Domain image from
        Wikipedia).}\label{HAFOUMSldRGIRA4AXELHARLROA7789GIJ7solFUMHDNOJHNFOJED}
\end{figure}

\begin{figure}
                \centering
                \includegraphics[width=.55\linewidth]{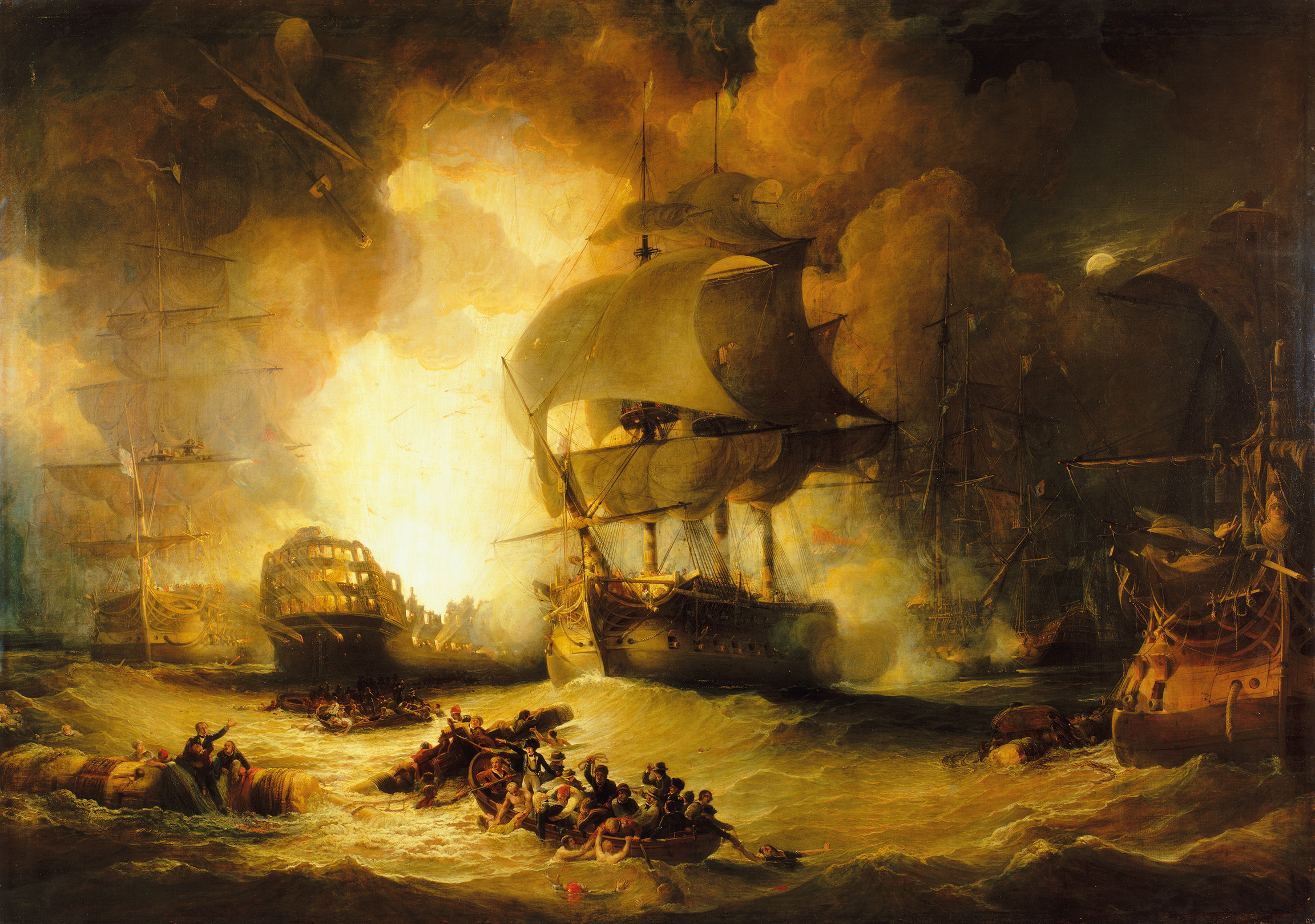}
        \caption{\sl The explosion of ``L'Orient'' at the Battle of the Nile, as depicted by artist George Arnald (Public Domain image from Wikipedia).}\label{HAFLORIENTF347-3-D}
\end{figure}

In its most basic formulation, the setting of this equation is due to
Jean-Baptiste Joseph Fourier who\footnote{Interestingly, while changing the course of mathematics and physics, Fourier was acting as a full time politician, since in~1801 Napol\'eon Bonaparte decided to appoint him as Prefect of the Department of Is\`ere, in the Alps. And perhaps it was not a completely trivial task to upset the First Consul (Emperor of the French from 1804) and dedicate oneself to mathematical investigations; in any case, playing around
with the heat equation at that point was not really part of Fourier's workload model.

See Figure~\ref{HAFOUMSldRGIRA4AXELHARLROA7789GIJ7solFUMHDNOJHNFOJED}
for a watercolor caricature of Fourier (by Julien-L\'eopold Boilly).

By the way, Napoleon had already included Fourier as the leader of the ``Legion of Culture'' which was sent to Egypt in~1798. The 150 renowned scholars of this selected group had the (maybe debatable) honor to travel to Egypt with Napoleon himself on the flagship of the French fleet ``L'Orient'' and to entertain him every evening with a scientific lecture. Allegedly, Fourier enjoyed the weather of the desert and, when he returned to Paris, he kept his rooms unreasonably hot and used to wrap up with multiple layers of very warm clothes.

The fate of the flagship ``L'Orient'' was instead unfortunate, since it exploded spectacularly at 10:30 PM, on August 1, 1798, during the Battle of the Nile, see Figure~\ref{HAFLORIENTF347-3-D}. Someone reported that the ship was repainted slightly before the battle, which was not a very wise choice, since the inflammable contents of the fresh painting acted as an energizer for the fire produced on the ship by Horatio Nelson's cannonballs.} presented for the first time a partial differential equation to describe conductive diffusion of heat (and, to study this equation, he also introduced a marvelous instrument that was going to change forever the History of the Universe, namely, the Fourier Series).

The ansatz\footnote{``Ansatz'' (plural: ``Ans\"atze'') is a masculine German name.
As quite customary in German, this name is obtained by a preposition (``an'',
meaning ``at'', or ``on'', or ``to'', or ``against'', or pretty much whatever one likes) and another name
(``Satz'', plural: ``S\"atze'', meaning ``sentence'', or ``statement''). And as quite customary in German,
a precise translation goes well beyond the expertise of modest mathematicians of our caliber.
The ``ansatz'' could be literally a ``starting point''. In jargon, it is frequently used to denote
an ``educated guess'' about a problem,
typically an additional simplifying assumption set at the beginning of an argument,
possibly to be verified later on (if it produces any result worthy of further consideration).

Maybe mathematics is not quite a deductive discipline, after all. Perhaps it is mostly the art of making good guesses.
Of course, anybody can make guesses. But for a good guess one needs to become acquainted with the problem
and to have developed that sort of familiarity that permits to separate the essential aspects from
the minor details (and, in trying to do this, sometimes
a little bit of luck doesn't hurt).} made by Fourier is that, at some\footnote{The space variable~$x$ here is supposed to be in~$\R^n$. Though the models presented refer to physical spaces, in which one can focus on the case~$n=3$ (or even~$n=2$, when dealing with plates, or~$n=1$, when dealing with ropes of negligible thickness), whenever possible we prefer to work in (Euclidean) spaces of arbitrary dimension. When one does that, the notation becomes usually more effective and the essential treats of the main ideas are very often more transparent. Also, when it will be needed to focus our attention on the special case of particular dimensions, one will be able to understand why and to what extent these dimensions are different from the general case. In dealing with the arbitrary dimension case, we can't help quoting the first ``general working principle'' in the Preface of~\cite{MR1625845}: ``PDE theory is (mostly) not restricted to two independent variables. Many texts describe PDE as if functions of the two variables~$(x, y)$ or~$(x, t)$ were all that matter. This emphasis seems to me misleading, as modern discoveries concerning many types of equations, both linear and nonlinear, have allowed for the rigorous treatment of these in any number of dimensions. I also find it unsatisfactory to classify partial differential equations: this is possible in two variables, but creates the false impression that there is some kind of general and useful classification scheme available in general''.}
point~$x\in\R^n$, at an instant of time~$t\in\R$, the variation of temperature~$u(x,t)$ of some body is determined by the heat flow around the point under consideration, plus possibly additional heat\footnote{For the sake of precision, we should probably say here ``heating or cooling sources'',
e.g., fireplaces, bonfires, stoves, radiators, gas heaters, infrared heaters, but also refrigerators, air conditioners, swamp boxes, chilled beams, or whatever heating or cooling device Fourier might have possessed at that time.
With this understanding, positive values of~$f$ would correspond to heating systems and negative values to cooling systems.}
sources. Let us denote by~$B(x,t)$ the heat
flux vector and by~$f(x,t,u(x,t))$ the scalar intensity of the external\footnote{To make things simpler, one can suppose that~$f$ depends only on~$x$, as a heat source which is permanently switched on or off, or on~$x$ and~$t$, in case the source is turned on and off or its intensity is modified over time. The additional dependence on~$u$ itself models the interesting case of a thermostat controlling a heating or cooling system.} heat sources.

A precise measure of the temperature precisely at a point~$x$ is certainly rather difficult from a practical point of view, so it is useful sometimes to consider instead an average\footnote{As usual, we use
here (and repeatedly over this set of notes) the integral notation for averages\index{average}, that is
$$ \fint_\Omega u(x,t)\,dx:=\frac1{|\Omega|}\,\int_\Omega u(x,t)\,dx,$$
being~$|\Omega|$ the Lebesgue measure of~$\Omega$.}
temperature, measured in some region of the space~$\Omega$,
say 
$$ U(t;\Omega):=\fint_\Omega u(x,t)\,dx.$$
In this setting, up to physical constants that we omit, we can equate the variation of~$U$
with the heat flux through~$\partial\Omega$, possibly adding to it the effect of the heat sources (if any) in~$\Omega$,
namely we write that
\begin{equation}\label{DAGB-ADkrVoiweLL4re2346ytmngrrUj}
\partial_t U(t;\Omega)=-\frac1{|\Omega|}\int_{\partial\Omega} B(x,t)\cdot\nu(x)\,d{\mathcal{H}}^{n-1}_x+\fint_\Omega f(x,t,u(x,t))\,dx.\end{equation}
In terms of notation, we are denoting here (and we will mostly keep this notation throughout the all set of notes)
by~$\nu(x)$ the unit normal vector at~$x\in\partial\Omega$ pointing outwards
\label{71KJSDNc9wov03ifgu9430eirqh89h32Xyr8tgyh4btr2} from~$\Omega$ and
by~$d{\mathcal{H}}^{n-1}_x$ the surface element\footnote{In many textbooks, the surface
element\index{surface element} is denoted by~$ds$, or~$dS$, or~$d\Sigma$. The
notation~$d{\mathcal{H}}^{n-1}_x$ comes from the fact that we are denoting by~${\mathcal{H}}^k$ the $k$-dimensional
Hausdorff measure, and, for smooth $(n-1)$-dimensional surfaces, the surface measure itself is precisely
equal to~${\mathcal{H}}^{n-1}$. See e.g.~\cite[Chapter~2]{MR3409135} for a thorough presentation of the Hausdorff measure.

The notation~$d{\mathcal{H}}^{n-1}_x$ used here is perhaps a bit heavier than~$ds$, or~$dS$, or~$d\Sigma$, or similar ones, but it has the benefit of being clearer and more explicit. Also, it allows us to consider surface integrals along less regular objects without having to introduce a new notation on a case-by-case basis.} on~$\partial\Omega$. The minus sign appearing on the right hand side of~\eqref{DAGB-ADkrVoiweLL4re2346ytmngrrUj} is due to the fact that the normal~$\nu$ points towards the exterior, while we are computing there the flux coming into the region~$\Omega$.

By differentiating under the integral sign\footnote{In this chapter about motivations, the arguments are developed at a formal level, we feel free to exchange derivatives and integrals, we do not keep track of lower order terms in the expansions, we do not discuss convergence issues and existence of limits, etc. This is a rather customary approach when one deals with providing convincing, but not necessarily circumstantial, motivations for a problem with the objective of developing some intuition about it. We promise to try to be more rigorous from the next chapter on.} in~\eqref{DAGB-ADkrVoiweLL4re2346ytmngrrUj}, we find that
\begin{equation*}\fint_\Omega\partial_t u(x,t)\,dx=-\frac1{|\Omega|}\int_{\partial\Omega} B(x,t)\cdot\nu(x)\,d{\mathcal{H}}^{n-1}_x+\fint_\Omega f(x,t,u(x,t))\,dx.\end{equation*}
Hence, by the Divergence\footnote{Typically we
reserve the notation ``$\nabla$'' for the vectors of the derivatives with respect to the space variables and~``$\div$'' for the corresponding divergence. No confusion should arise with respect to derivatives with respect to the time variable, which is usually denoted here by~``$\partial_t$''.}
Theorem,
\begin{equation}\label{DAGB-ADkrVoiweLL4re2346ytmngrrUj2}\fint_\Omega\partial_t u(x,t)\,dx=-\fint_{\Omega} \div B(x,t)\,dx+\fint_\Omega f(x,t,u(x,t))\,dx.\end{equation}

{F}rom this, Fourier took a further step by realizing that to make this identity manageable one needs to
take a constitutive law about the flux vector~$B$, relating it to the temperature. Fourier's new ansatz was to suppose that the 
flux of heat between two adjacent (infinitesimal) regions is proportional to the (infinitesimal) difference of their temperatures, namely
\begin{equation}\label{DAGB-ADkrVoiweLL4re2346ytmngrrUj3}
B(x,t)=-\kappa(x,t)\,\nabla u(x,t).
\end{equation}
The proportional coefficient~$\kappa$ is sometimes called\footnote{See Section~\ref{Whenyouareamathematician} for its counterpart in electrostatics.} ``heat conduction coefficient''\index{heat conduction coefficient}. We take~$\kappa$ to be positive:
with respect to this, we stress that the minus sign in~\eqref{DAGB-ADkrVoiweLL4re2346ytmngrrUj3} is motivated by the fact that the heat flows from hotter regions to colder ones (hence in the opposite direction of the growth of the temperature function~$u$).
The relation in~\eqref{DAGB-ADkrVoiweLL4re2346ytmngrrUj3} is also sometimes\footnote{Actually, the same presentation
here could have been done to describe the transport of mass through diffusive means.
In this setting, one argue that the flux goes from regions of high concentration to regions of low concentration.
The ansatz corresponding to~\eqref{DAGB-ADkrVoiweLL4re2346ytmngrrUj3} would be that the
magnitude of the flux is proportional to the concentration gradient.
In the context of mass diffusion, $\kappa$ is sometimes called ``diffusion coefficient''
and~\eqref{DAGB-ADkrVoiweLL4re2346ytmngrrUj3} is sometimes refereed to with the name of
``Fick's Law''\index{Fick's Law} after the German physician and physiologist
Adolf Eugen Fick.

Alternative choices for the constitutive relation are possible. For instance, one can replace~\eqref{DAGB-ADkrVoiweLL4re2346ytmngrrUj3} with a more general equation in which the flux depends possibly in a nonlinear way
from the gradient of the temperature (in the case of Fourier's model, or of the transported mass in the case of Fick's model),
that is one could suppose that~$B=-\Phi(\nabla u)$, for some function~$\Phi:\R^n\to\R^n$. A typical example is the case in which~$\Phi(\nabla u)=|\nabla u|^{p-2}\nabla u$, that is~$B=-\kappa|\nabla u|^{p-2}\nabla u$,
meaning that the flux is proportional to a power of the gradient of~$u$ (more precisely, the vector~$B$ has the same direction of the gradient of~$u$ and its magnitude is proportional to a power of the magnitude of the gradient of~$u$):
this setting would lead to the so-called ``$p$-Laplace equation''.

Another setting of interest is the one in which the flux is related to a pressure drop
(this framework is usually related to the so-called Darcy's Law\index{Darcy's Law}), say~$B=-\kappa\nabla P$, where~$P$ has some physical meaning of pressure. In this sense, the case~$P=u$ reduces to~\eqref{DAGB-ADkrVoiweLL4re2346ytmngrrUj3},
but other situations are of interest. For instance, one could assume that~$P$ and~$u$ are related by
a state equation of the form~$P(x)=p(u(x))$. A natural possibility is to take the function~$p$ to be a
power of~$u$, e.g.~$p=u^m$. This choice would lead to the so-called porous medium equation.

These notes are of elementary nature, hence we will not address the cases of the $p$-Laplace equation,
or of the porous medium equation, or other types of ``anomalous diffusions'', such as the type of diffusion that takes into
account mass transfer from remote regions due to long-range interactions.
The reader interested in these more advanced topics may look e.g. to~\cite{MR781350, MR1230384, MR2305115, MR2286292, MR3469920} and the references therein.

For the readers interested in using Darcy's Law to have a proper cup of coffee, see e.g.~\cite{2008AmJPh}
and the references therein.} called ``Fourier's Law''.

\begin{figure}
                \centering
                \includegraphics[width=.3\linewidth]{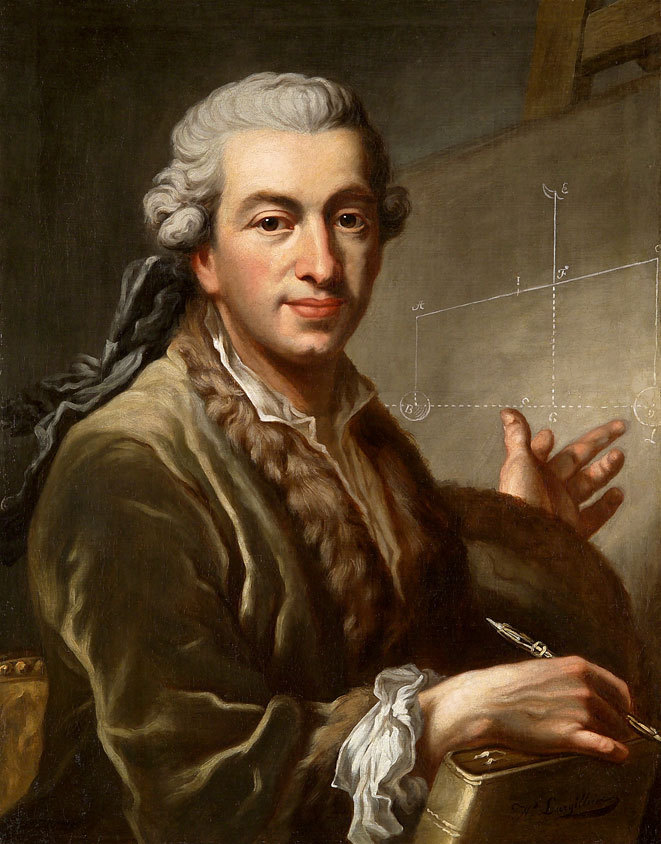}
        \caption{\sl Portrait of Pierre-Simon Laplace by Johann Ernst Heinsius (image from
        Wikipedia,
        licensed under the Creative Commons Attribution-Share Alike 4.0 International license).}\label{2HAFODAKEDFUMSPierre-SimonLaplacldRGIRA4AXELEUMDJOMNFHARLROA7789GIJ7solFUMHDNOJHNFOJED231}
\end{figure}

By substituting~\eqref{DAGB-ADkrVoiweLL4re2346ytmngrrUj3} into~\eqref{DAGB-ADkrVoiweLL4re2346ytmngrrUj2} we find that
\begin{equation}\label{DAGB-ADkrVoiweLL4re2346ytmngrrUj4}
\fint_\Omega\partial_t u(x,t)\,dx=\fint_{\Omega} \div \big(
\kappa(x,t)\,\nabla u(x,t)
\big)\,d{\mathcal{H}}^{n-1}_x+\fint_\Omega f(x,t,u(x,t))\,dx.\end{equation}
Since this identity holds true for all regions~$\Omega$, we find a pointwise counterpart of~\eqref{DAGB-ADkrVoiweLL4re2346ytmngrrUj4} by writing
\begin{equation}\label{DAGB-ADkrVoiweLL4re2346ytmngrrUj6}
{\partial_t u(x,t)}= \div \big(
\kappa(x,t)\,\nabla u(x,t)
\big)+f(x,t,u(x,t)).\end{equation}
The case of homogeneous media in which~$\kappa$ is constant (say, equal to~$1$ up to a renormalization of units of measure)
is of particular interest: in this case~\eqref{DAGB-ADkrVoiweLL4re2346ytmngrrUj6} boils down\footnote{We exploit the
standard notation for the Laplacian\index{Laplacian} given by
$$ \Delta u:=\div(\nabla u)=\sum_{j=1}^n \frac{\partial^2 u}{\partial x_j} .$$
In some textbook, the use of~``$\Delta$'' is replaced by~``$\nabla^2$''
or by~``$|\nabla|^2$''. Of course all notations are good. Personally, we have a preference for
the $\Delta$ notation since we find it simpler to read. Also, in a sense, it highlights the fact that the Laplace operator
has a ``dignity'' which is ``independent from the one of the gradient''.
This ``philosophical'' point will be perhaps clarified by the geometric interpretation of the Laplacian that we will discuss
in the forthcoming Theorem~\ref{KAHAR2}.

The name of the Laplace operator comes from Pierre-Simon, marquis de Laplace, who introduced it
while studying celestial mechanics, in connection with the gravitational potential (this approach will
be exploited here when dealing with the fundamental solution in Section~\ref{lfundsP-S}.
See Figure~\ref{2HAFODAKEDFUMSPierre-SimonLaplacldRGIRA4AXELEUMDJOMNFHARLROA7789GIJ7solFUMHDNOJHNFOJED231}
for a portrait of Laplace (by Johann Ernst Heinsius) when he was 35 years old, with an easel that quite resembles
a blackboard,  working instrument and sign of distinction of every passionate mathematician.

Let us remark that a slightly different approach to the derivation of the heat equation is possible,
by considering ``heat'' instead of ``temperature'' as the main building block of the equation.
The disadvantage of having heat, rather than temperature, in a pivotal role is that it is perhaps a more vague and less intuitive concept 
than temperature (just because temperature seems a notion we are so familiar with, e.g. in view of weather forecasts).
However, the notion of heat relates more directly to thermal energy and constitutes an extensive (i.e., additive)
property. In this perspective,
the heat equation reflects an energy budget
in which the heat change in time is equal to heat produced (source) plus heat entering through the boundary (flux).
Namely, if one is willing to take heat, instead of temperature, as the building block for the heat equation,
the balance equation in~\eqref{DAGB-ADkrVoiweLL4re2346ytmngrrUj} reads at the level of
variation of heat, instead of averaged temperature
(and note that the constitutive relation in the Fourier's Law~\eqref{DAGB-ADkrVoiweLL4re2346ytmngrrUj3}
is already a relation between the notions of heat and temperature, since it assumes the
heat flux to be proportional to the gradient of temperature).

The final heat equation in~\eqref{DAGB-ADkrVoiweLL4re2346ytmngrrUj7} would remain essentially
unchanged since the changes of heat content~$Q$ directly relates to the changes in temperature
via the relation~$\partial_t Q= c \rho \partial_tu$, where~$\rho$ is mass density and~$c$ is the specific heat capacity (a constant that depends on the internal properties of the material, which can be tabulated through experiments).}
to
\begin{equation}\label{DAGB-ADkrVoiweLL4re2346ytmngrrUj7}
\partial_t u(x,t) = \Delta u(x,t)+f(x,t,u(x,t)).\end{equation}
Equation~\eqref{DAGB-ADkrVoiweLL4re2346ytmngrrUj7} (or variations of it) is typically\footnote{Often,
equations as in~\eqref{DAGB-ADkrVoiweLL4re2346ytmngrrUj7} are referred to as ``parabolic''.
The name comes from the following classification method.
A general linear differential operator of second order in the variables~$(X_1,\dots,X_N)$ has the form
$$ \sum _{i,j=1}^{N} a_{i,j}{\frac{\partial ^{2}}{\partial X_{i}\partial X_{j}}}.$$
Then one classifies the operator, and the corresponding partial differential equation,
depending on the sign of the eigenvalues of the coefficient matrix~$a_{i,j}$.
More specifically, when all the eigenvalues of~$a_{i,j}$ have the same sign (either are all strictly
positive or all strictly negative), the equation is named ``elliptic''\index{elliptic}.
When one eigenvalue is zero and all the other eigenvalues have the same sign
(either all strictly
positive or all strictly negative), the equation is referred to as ``parabolic''\index{parabolic}.
When all the eigenvalues have the same sign except one which has the opposite sign
(i.e. if either there is only one negative eigenvalue and all the rest are positive,
or there is only one positive eigenvalue and all the rest are negative),
the equation is named ``hyperbolic''\index{hyperbolic}.
This classification is not exhaustive: the remaining cases not classified here are typically quite
hard to study, require specific techniques and no general theory is available for those. \label{CLASSIFICATIONFOOTN}

In the case of the heat equation in~\eqref{DAGB-ADkrVoiweLL4re2346ytmngrrUj7}, one has~$N=n+1$, $X=(x,t)$ and
$$ a_{ij}=\begin{dcases}1 & {\mbox{ if }}i=j\in\{1,\dots,n\},\\
0&{\mbox{ otherwise }}.\end{dcases}$$
The corresponding eigenvalues are therefore~$0$, with multiplicity~$1$, and~$1$, with multiplicity~$n$,
hence the equation in~\eqref{DAGB-ADkrVoiweLL4re2346ytmngrrUj7} is parabolic.}
called ``the heat equation''\index{heat equation}.\medskip

Suppose now that~$f$ in~\eqref{DAGB-ADkrVoiweLL4re2346ytmngrrUj7} is independent of time, i.e.~$f=
f(x,u(x,t))$. In this framework,
of particular interest are certainly stationary solutions of~\eqref{DAGB-ADkrVoiweLL4re2346ytmngrrUj7},
i.e. solutions~$u=u(x)$ which are also independent of time: these solutions are actually equilibria of~\eqref{DAGB-ADkrVoiweLL4re2346ytmngrrUj7}, since they provide solutions of~\eqref{DAGB-ADkrVoiweLL4re2346ytmngrrUj7}
which remain the same at every instant of time. Of course, when~$f$ and~$u$ are independent of time,
equation~\eqref{DAGB-ADkrVoiweLL4re2346ytmngrrUj7} boils down to\footnote{With respect to the classification
presented in footnote~\ref{CLASSIFICATIONFOOTN}, equation~\eqref{DAGB-ADkrVoiweLL4re2346ytmngrrUj8},
\eqref{DAGB-ADkrVoiweLL4re2346ytmngrrUj9} and~\eqref{DAGB-ADkrVoiweLL4re2346ytmngrrUj10} are of elliptic type
(actually, they provide us with the ``paradigmatic'' form of elliptic equations).
To check their ellipticity, one exploits the setting of footnote~\ref{CLASSIFICATIONFOOTN}
with~$N=n$, $X=x$ and
$$a_{ij}=\delta_{ij}:=\begin{dcases}
1 & {\mbox{ if }}i=j,\\
0 &{\mbox{ if }}i\ne j.
\end{dcases}$$
In this scenario, all the eigenvalues of the coefficient matrix
are equal to~$1$, hence strictly positive.

Further insight on the notion of ellipticity will be given on pages~\pageref{CLASSIFICATIONFOOTN2}
and~\pageref{CLASSIFICATIONFOOTN3}.

The above definition of~$\delta_{ij}$ will be also used throughout all these notes (sometimes \label{Kronecker-notation}
this definition of~$\delta_{ij}$ is referred to with the name of Kronecker notation\index{Kronecker notation}).}
\begin{equation}\label{DAGB-ADkrVoiweLL4re2346ytmngrrUj8}
\Delta u(x)+f(x,u(x))=0.\end{equation}
In particular, when~$f$ is independent of~$u$, equation~\eqref{DAGB-ADkrVoiweLL4re2346ytmngrrUj8} reduces to
\begin{equation}\label{DAGB-ADkrVoiweLL4re2346ytmngrrUj9}
\Delta u(x)+f(x)=0,\end{equation}
which is sometimes called Poisson's equation.

The case in which~$f$ in~\eqref{DAGB-ADkrVoiweLL4re2346ytmngrrUj9} vanishes identically produces
\begin{equation}\label{DAGB-ADkrVoiweLL4re2346ytmngrrUj10}
\Delta u(x)=0,\end{equation}
which is called Laplace's equation\index{Laplace's equation}.

Solutions of~\eqref{DAGB-ADkrVoiweLL4re2346ytmngrrUj10} are named ``harmonic functions''\index{harmonic function}.
The study of harmonic functions is essential, one way or another, in virtually any aspect of mathematics
(including mathematical analysis, mathematical physics, complex analysis, geometry, probability, finance, statistics, you name it). Also, the theory of harmonic functions is a clear example of adamantine beauty and elegance
in which human creativity has reached one of its highest peaks ever.

In the forthcoming pages we will do our best to present some bits of this theory, but we definitely invite
the reader to look to as many other sources as possible to collect the broadest possible amount of information and
develop their own point of view on this topic of paramount importance.

\subsection{Population dynamics, chemotaxis and random walks}\label{CHEMOTX}

\begin{figure}
                \centering
                \includegraphics[width=.75\linewidth]{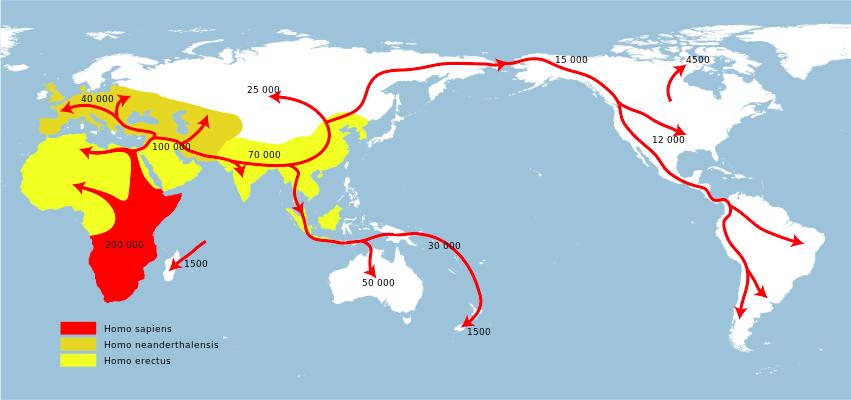}
        \caption{\sl Spreading humanity all over the world (Public Domain image by NordNordWest from
        Wikipedia).}\label{HAFOBAomujnerfOJED}
\end{figure}

Another fascinating situation in which partial differential equations surface very often
occurs in the description of biological populations. This is a very captivating field of study,
which is truly crossdisciplinary in spirit and collects fundamental questions from (at least) mathematics, physics,
biology, ethology and social sciences.
For the purpose of these pages, we limit ourselves to a simple description of a 
biological species exhibiting at a point~$x\in\R^n$ and at a time~$t\in\R$ a population density of the form~$u(x,t)$.
We assume that a fraction~$\rho\in[0,1]$ of
the population has the tendency of moving randomly, while the rest of the population
(corresponding to a fraction~$\mu:=1-\rho$ of the totality of individuals) is driven
towards regions of higher concentration of a given attractant (e.g. a chemical signal)
distributed according to a function~$w(x,t)$ (the
movement of an organism in response to a chemical stimulus and specifically
in dependence of the increasing or decreasing concentration of a particular substance
is called in jargon ``chemotaxis''\index{chemotaxis}).
Additionally, we suppose that the population is subject to a drift\index{drift} (e.g. due to the wind, or the stream, or the tide) \label{DRIFPAG}
in a given direction~$b(x,t)$.
The simpler cases in which there is no chemotactic effect (corresponding to~$\mu:=0$),
or no random motion (corresponding to~$\rho:=0$), or no drift (corresponding to~$b:=0$)
are interesting special situations of the complex phenomena for which we now present a mathematical model
(the random component of this discussion will also be further generalized on page~\pageref{0uojf29249-45kpkfdSmd11493839429efv}
to more elaborated environments). See Figure~\ref{HAFOBAomujnerfOJED}
for a description of how our ancestors spread over the world. \label{RANDOW}

In outline, we consider a small spatial scale~$h$ and a small time scale~$\tau$.
We will choose in what follows a suitable relation between space and time scales to provide
coherent asymptotics in the formal limit.
In this setting,
at every unit of time the variation of the density of
the population is influenced by the previous factors in a rather explicit way. For the moment, let us focus on the simple case in which there is no chemotaxis and no drift. In this situation, $\rho:=1$, $\mu:=0$, $b$ vanishes identically
and the only drive for the population comes from a random walk.
To get to the bottom of this case,
we consider, at a given point~$x$, at time~$t+\tau$, the new random population~$\rho u(x,t+\tau)$
(which actually coincides with the whole population~$u(x,t+\tau)$, but let us keep the fraction~$\rho$
for future use, even if for the moment~$\rho=1$). Then, we have that~$\rho u(x,t+\tau)$
is produced by the random population at the previous time~$t$ which is located at some point~$x+h\omega$, for a given direction~$\omega\in\partial B_1$, times the probability
that these individuals move from~$x+he$ to~$x$ in the unit of time. If all directions are equally probable, this says that the new random population~$\rho u(x,t+\tau)$ is produced by a term of the form~$\rho\fint_{\partial B_1} u(x+h\omega,t)\,d{\mathcal{H}}^{n-1}_\omega$, that is
\begin{equation}\label{JOSN-Oohpiklfoigheigohbo-gherpoguev}
\rho u(x,t+\tau)
=\rho\fint_{\partial B_1} u(x+h\omega,t)\,d{\mathcal{H}}^{n-1}_\omega.
\end{equation}
Subtracting~$\rho u(x,t)$ from both sides, dividing by~$\tau$ and choosing
\begin{equation}\label{SCAKLI-UJs903th}
h:=\sqrt\tau,\end{equation}
we find that
\begin{equation}\label{7tnfesewuxttau-x21s46}\begin{split}
\rho \frac{u(x,t+\tau)-u(x,t)}\tau\,&
=\,\frac1{\tau}\left[\rho\fint_{\partial B_1} u(x+h\omega,t)\,d{\mathcal{H}}^{n-1}_\omega-\rho u(x,t)\right]\\&
=\,\frac{\rho}{h^2}\fint_{\partial B_1} \Big(u(x+h\omega,t)-u(x,t)\Big)\,d{\mathcal{H}}^{n-1}_\omega.\end{split}
\end{equation}
Now we aim at showing that, for small~$h$,
\begin{equation}\label{KMS:RUIDKVODIG}
\fint_{\partial B_1} \Big(u(x+h\omega,t)-u(x,t)\Big)\,d{\mathcal{H}}^{n-1}_\omega=
\frac{c h^2}2 \Delta u(x,t) +o(h^2),
\end{equation}
for some~$c>0$ depending on the dimension~$n$ (see formula~\eqref{CDEGLDDEFC} below).
For this, we
consider the formal Taylor expansion
$$ u(x+h\omega,t)=u(x,t)+h\nabla u(x,t)\cdot \omega+\frac{h^2}2\,D^2u(x,t) \omega\cdot \omega+o(h^2).$$
We observe, by symmetry, that
$$ \int_{\partial B_1} \nabla u(x,t)\cdot \omega\,d{\mathcal{H}}^{n-1}_\omega=0,$$
because\footnote{Spotting these types of cancellations will be often an essential ingredient of our arguments.
After all, one way to obtain significant information about a complex problem is to understand
which quantities do not play much of a role since they ``average out''.
At this stage, we discuss these simplifications at an intuitive level. Later on, see e.g.~\eqref{22345XXse5q6u5at3i332o2n},
we will give more rigorous arguments to justify them, but for the moment the main goal is to try to visualize the
possibility of these cancellations and, especially, to appreciate their importance.}
any (positive or negative) contribution to such surface integral coming from a point~$\omega\in\partial B_1$
is canceled precisely by the (negative or positive) contribution coming from~$-\omega$.
{F}rom these observations we arrive at
\begin{eqnarray*} &&\fint_{\partial B_1} \Big(u(x+h\omega,t)-u(x,t)\Big)\,d{\mathcal{H}}^{n-1}_\omega
=\fint_{\partial B_1} \left(h\nabla u(x,t)\cdot \omega+\frac{h^2}2\,D^2u(x,t) \omega\cdot \omega\right)\,d{\mathcal{H}}^{n-1}_\omega+o(h^2)\\&&\qquad=
\frac{h^2}2\,\fint_{\partial B_1}D^2u(x,t) \omega\cdot \omega\,d{\mathcal{H}}^{n-1}_\omega+o(h^2)=\frac{h^2}2
\sum_{i,j=1}^n\fint_{\partial B_1}\partial_{ij} u(x,t) \omega_i\omega_j\,d{\mathcal{H}}^{n-1}_\omega+o(h^2)
.\end{eqnarray*}
We now identify an additional simplification by observing that if~$i\ne j$ then
$$ \fint_{\partial B_1}\partial_{ij} u(x,t) \omega_i\omega_j\,d{\mathcal{H}}^{n-1}_\omega=0,$$
because whatever contribution comes from~$\omega_i$ is canceled precisely by that coming from~$-\omega_i$
(notice how helpful was for all these cancellations that we are integrating over such a symmetric domain as~$\partial B_1$,
which remains invariant under all these reflections). 

Therefore,
$$ \fint_{\partial B_1} \Big(u(x+h\omega,t)-u(x,t)\Big)\,d{\mathcal{H}}^{n-1}_\omega=\frac{h^2}2
\sum_{i=1}^n\fint_{\partial B_1}\partial_{ii} u(x,t) \omega_i^2\,d{\mathcal{H}}^{n-1}_\omega+o(h^2).
$$
Now, once again, we understand the symmetries of the problem under consideration. For this, we notice that,
for every~$i$, $j\in\{1,\dots,n\}$
$$ \fint_{\partial B_1} \omega_i^2\,d{\mathcal{H}}^{n-1}_\omega
=\fint_{\partial B_1} \omega_j^2\,d{\mathcal{H}}^{n-1}_\omega,$$
because\footnote{Once again, we aim here at developing arguments in the most intuitive way possible.
This idea with be more rigorously retaken later on, see~\eqref{H66}.}
the role played by the~$i$th coordinate in the sphere is precisely the same as the one played by the~$j$th coordinate.
Therefore, we can define
\begin{equation}\label{CDEGLDDEFC} c:=\fint_{\partial B_1} \omega_i^2\,d{\mathcal{H}}^{n-1}_\omega\end{equation}
and we stress that this quantity does not depend on~$i$.

Consequently,
$$ \fint_{\partial B_1} \Big(u(x+h\omega,t)-u(x,t)\Big)\,d{\mathcal{H}}^{n-1}_\omega=\frac{c h^2}2
\sum_{i=1}^n \partial_{ii} u(x,t) +o(h^2),
$$
from which we obtain~\eqref{KMS:RUIDKVODIG}, as desired.

Hence we insert~\eqref{KMS:RUIDKVODIG} into~\eqref{7tnfesewuxttau-x21s46} to find that
\begin{equation}\label{y8oqiwhfesjd8923wfyug3iuweLuhonwfdpinrgNqikedOKS}\begin{split}&
\rho \frac{u(x,t+\tau)-u(x,t)}\tau 
=\frac{c\rho}{2}\Delta u(x,t) +o(1).
\end{split}
\end{equation}
By sending~$\tau\searrow0$ (and thus~$h\searrow0$), we thereby find that
\begin{equation}\label{HJA:A89MAOL8789892}
\rho \partial_t u(x,t)
= \frac{c\rho}{2}\Delta u(x,t) ,
\end{equation}
which coincides, up to constant, with the heat equation\footnote{See
e.g.~\cite{MR2732325}
for more information about the strong connections between
the random walk and the heat equation.} presented in~\eqref{DAGB-ADkrVoiweLL4re2346ytmngrrUj7}
(here, with no source term).

\begin{figure}
                \centering
                \includegraphics[width=.4\linewidth]{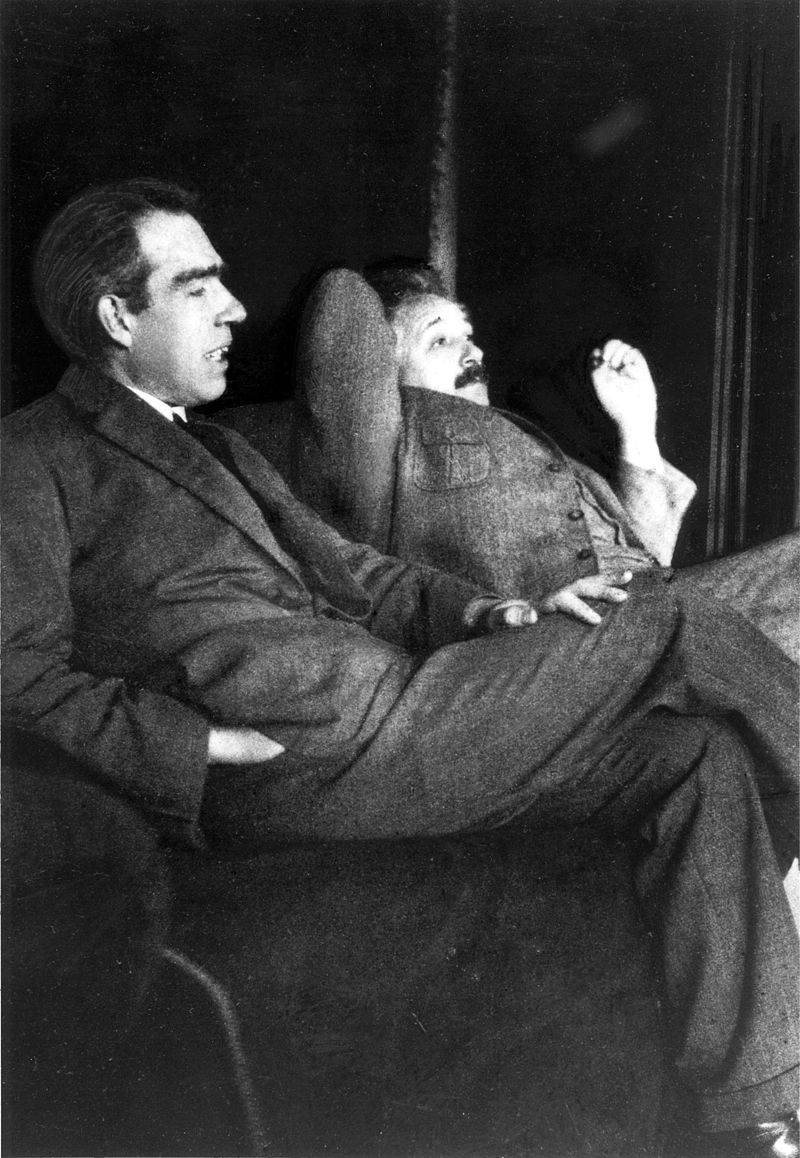}
        \caption{\sl Hard-working scientists (Public Domain image from
        Wikipedia).}\label{h2HAFODAKEDFUMSldRGIRA4AXELEUMDJOMNFHARLROA7789GIJ7solFUMHDNOJHNFOJED231EIN}
\end{figure}

\begin{figure}
                \centering
                \includegraphics[width=.58\linewidth]{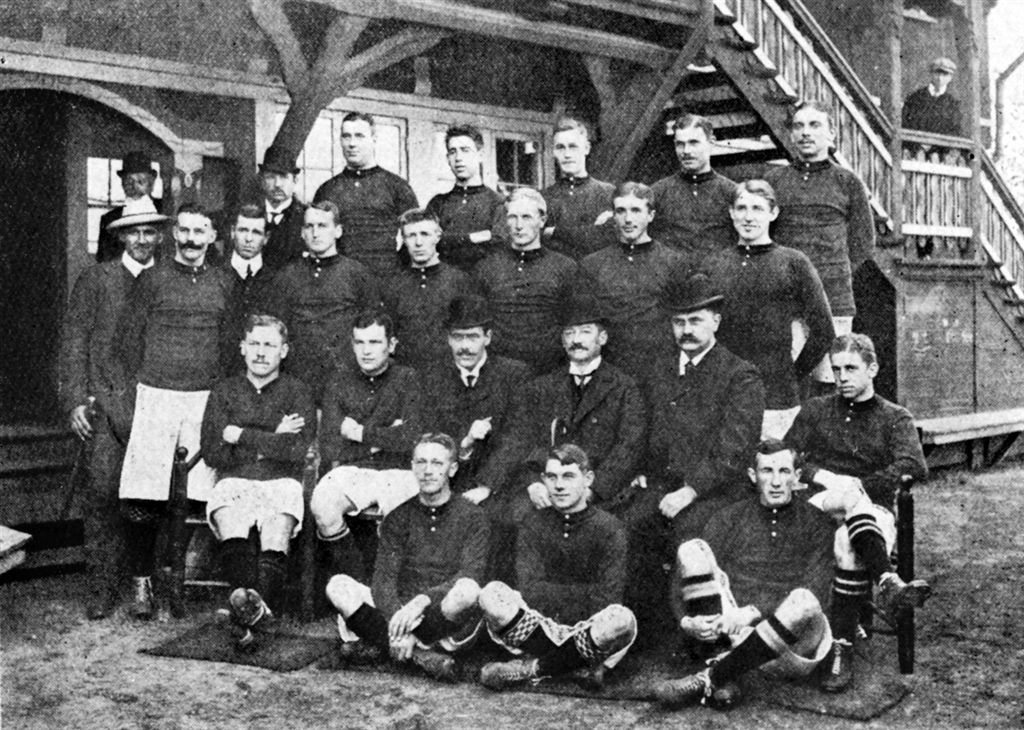}
        \caption{\sl Danish soccer team at the 1908 Olympic Games (Public Domain image from
        Wikipedia).}\label{2HAFODAKEDFUMSldRGIRA4AXELEUMDJOMNFHARLROA7789GIJ7solFUMBOHr32HDNOJHNFOJED231}
\end{figure}

This computation was not only instructive from the technical point of view, but it also unveiled one of the mainsprings
of the mathematical theory of diffusion, revealing the strong conceptual connection
of a substance (or a population) moving randomly according to a Brownian motion\index{Brownian motion} and the way in which
heat dissipates: after all, the spreading of temperature all over the region subject to diffusion is, in some sense, nothing else\footnote{The link between the random motion of molecules and the
macroscopic phenomenon of diffusion is indeed a deep feature. While the Brownian motion was discovered
in~1827 by botanist Robert Brown by looking at dust grains floating in water, Albert Einstein
in his ``annus mirabilis''~1905
gave a quantitative model for the motion of floating particles as being moved by individual water molecules~\cite{zbMATH02652222}. Not only this article
founded the statistical physics analysis of Brownian motion, but also it
provides a way to determine the mass and the dimensions of the atoms involved in the
process, thus changing atomic theory from a controversial set of conjectures
into an established fact of science. We will come back to the theory proposed by Einstein in Section~\ref{EIN:SEC}.

See Figure~\ref{h2HAFODAKEDFUMSldRGIRA4AXELEUMDJOMNFHARLROA7789GIJ7solFUMHDNOJHNFOJED231EIN}
depicting Albert Einstein together with Niels Bohr.

By the way, Niels Bohr's brother was Harald Bohr, who was a mathematician, pioneering almost periodic functions,
also author, together with Edmund Landau, of an elegant result
related to the Riemann Conjecture (see~\cite[Section~9.6]{MR1854455}), and a soccer player. The Denmark national soccer team in which he was playing took part in the  1908 Summer Olympics, where football was an official event for the first time. Harald Bohr scored two goals in the first match and played in the semifinal (Denmark 17 - France 1, which remains the Olympic record of the most goals scored by one team). After that, the Denmark team lost the final and won the silver medal. In all the subsequent editions of the Olympic Games, Denmark soccer team won two more silver medals (in 1912 and 1960) and a bronze (in 1948); hence Harald Bohr's team also holds the record (ex aequo) of best placement for a Denmark soccer team in the Olympics.

See Figure~\ref{2HAFODAKEDFUMSldRGIRA4AXELEUMDJOMNFHARLROA7789GIJ7solFUMBOHr32HDNOJHNFOJED231}
for a photo of the Danish soccer team in the 1908 Olympic Games: Harald Bohr is the second player from the left in the top row.

Niels Bohr was also a passionate soccer player. He played goalkeeper in the Copenhagen-based team Akademisk Boldklub (at the time, one of the best clubs in Denmark; actually, the two brothers played several matches together in this team). One can wonder however why Niels never made it to the national team. Well, according to {\tt https://www.theguardian.com/football/2005/jul/27/theknowledge.panathinaikos}
in a match against a German team, one of the midfielders of the opposite team launched a very long shot and Niels, leaning against the post, did not react at all, missing an easy save and letting the German team score. After the game, Niels admitted that on that occasion he had been distracted by a mathematical problem he was thinking about. Strangely enough, Niels did not play for Akademisk Boldklub after that season.}
but the tendency of random movements to distribute mass in average from regions of high density
towards regions of small density. In both cases, the process
has the strong tendency to ``balance out'' differences:
as time goes, hot spots lose their temperature in favor of cold spots which are heated up by
the temperature of their neighbors, as well as under a random walk
the regions with low density get occupied by the population
coming from the highly populated regions. 

This ``democratic'' tendency of averaging out differences is typical of the Laplace operator and it \label{DEEFFNVST}
will be the underlying feature of all the ``regularity theories'' for elliptic equations that we will present in the forthcoming
pages.

Having well understood the case in which the biological species is subject solely on a random motion,
we now retake the interesting case in which a chemotactic agent
also comes into play, for a fraction~$\mu$ of the population (for the moment, we still suppose that there is no additional drift, hence~$b$ vanishes identically). In this situation,
the fraction~$\rho$ of the population performing the random walk would be still described
by~\eqref{y8oqiwhfesjd8923wfyug3iuweLuhonwfdpinrgNqikedOKS} (that was the reason to include~$\rho$ in the previous
computation, even if before~$\rho$ was just equal to~$1$) and we focus now on
the fraction~$\mu$ of the population following the attractant~$w$.
To appreciate the effect of the chemotactic factor, we suppose that such a fraction~$\mu$ of the population with density~$u$ does not move completely randomly, but, at each time step, picks a direction of motion with a probability~$\lambda$ that is proportional to values of the attractant having density~$w$. 

More explicitly, given~$\omega\in\partial B_1$,
we suppose for the chemotactic population that a jump from the point~$x$ to the point~$x+h\omega$
occurs with probability
\begin{equation}\label{RISUGIBNSUPJDMD} \lambda(x,\omega,t):=\frac1{{\mathcal{H}}^{n-1}(\partial B_1)}+w(x+h\omega,t)-w(x-h\omega,t).\end{equation}
In this setting, we are assuming that the oscillations of~$w$ are sufficiently small\footnote{More precisely,
one could rewrite~\eqref{RISUGIBNSUPJDMD} as
\[ \lambda(x,\omega,t):=\frac1{{\mathcal{H}}^{n-1}(\partial B_1)}+\e_0\Big(w(x+h\omega,t)-w(x-h\omega,t)\Big),\]
where~$\e_0>0$ is a parameter that takes into account the ``sensibility'' of the chemotactic population
to the chemical attractant. Just to be consistent with the probability scenario, it is convenient
to think that~$\lambda\ge0$ in view of the smallness of~$\e_0$ with respect to the oscillations
of the attractor's density~$w$.}
to make the above quantity positive.
In this sense, the probability~$\lambda$ differs from that of the classical random walk (corresponding to~$w$
being constant) since it increases in favor of the directions~$\omega$ in which the concentration of the chemical
attractant is higher. We stress that the above setting provides a normalized probability,
since the total probability of
jumping to the sphere of radius~$h$ around the given point~$x$ is~$1$, because
\begin{equation*}
\int_{\partial B_1}\lambda(x,\omega,t)\,d{\mathcal{H}}_\omega^{n-1}
=1+\int_{\partial B_1}\Big(w(x+h\omega,t)-w(x-h\omega,t)\Big)\,d{\mathcal{H}}_\omega^{n-1}=1,
\end{equation*}
due to odd symmetry cancellations.

We also point out that, in the above setting, the probability of a jump from a point~$a$ to
a point~$b$ with~$|a-b|=h$ is given by
\begin{eqnarray*} \lambda\left(a,\frac{b-a}{h},t\right)&=&\frac1{{\mathcal{H}}^{n-1}(\partial B_1)}+w\left(a+h\frac{b-a}{h},t\right)-w\left(a-h\frac{b-a}{h},t\right)\\&=&
\frac1{{\mathcal{H}}^{n-1}(\partial B_1)}+w(b,t)-w(2a-b,t).\end{eqnarray*}
In consequence, given~$\omega\in\partial B_1$, 
taking~$a:=x+h\omega$ and~$b:=x$, we find that
the jump from~$x+h\omega$ to~$x$ occurs with probability
$$ \lambda (x+h\omega,-\omega,t)=
\frac1{{\mathcal{H}}^{n-1}(\partial B_1)}+w(x,t)-w(x+2h\omega,t).$$
This being so, we can detect the chemotactic counterpart of~\eqref{JOSN-Oohpiklfoigheigohbo-gherpoguev}.
Indeed, since the chemotactic population~$\mu u(x,t+\tau)$ is produced by a term of the form~$\mu\int_{\partial B_1} \lambda (x+h\omega,-\omega,t) u(x+h\omega,t)\,d{\mathcal{H}}^{n-1}_\omega$, we have that
\begin{equation*} \begin{split}
\mu u(x,t+\tau)&=\mu\int_{\partial B_1} \lambda (x+h\omega,-\omega,t) u(x+h\omega,t)\,d{\mathcal{H}}^{n-1}_\omega\\
&=\mu\int_{\partial B_1} \left(\frac1{{\mathcal{H}}^{n-1}(\partial B_1)}+w(x,t)-w(x+2h\omega,t)\right)u(x+h\omega,t)\,d{\mathcal{H}}^{n-1}_\omega\\&=\mu\fint_{\partial B_1} u(x+h\omega,t)\,d{\mathcal{H}}^{n-1}_\omega
+\mu\int_{\partial B_1} \Big(w(x,t)-w(x+2h\omega,t)\Big)u(x+h\omega,t)\,d{\mathcal{H}}^{n-1}_\omega
.\end{split}\end{equation*}
Thus, we can subtract to both sides~$\mu u(x,t)$, divide by~$\tau=h^2$ and recall~\eqref{KMS:RUIDKVODIG}
and~\eqref{CDEGLDDEFC} to find that
\begin{eqnarray*}
&&\mu \frac{u(x,t+\tau)-u(x,t)}\tau\\ &=&
\frac1{h^2}\left[\mu\fint_{\partial B_1}\Big( u(x+h\omega,t)-u(x,t)\Big)\,d{\mathcal{H}}^{n-1}_\omega+
\mu\int_{\partial B_1} \Big(w(x,t)-w(x+2h\omega,t)\Big)u(x+h\omega,t)\,d{\mathcal{H}}^{n-1}_\omega
\right]\\&=& \frac{c\mu}{2}\Delta u(x,t) +
\frac{\mu}{h^2}\int_{\partial B_1} \Big(w(x,t)-w(x+2h\omega,t)\Big)u(x+h\omega,t)\,d{\mathcal{H}}^{n-1}_\omega
+o(1)\\&=& \frac{c\mu}{2}\Delta u(x,t)-
\frac{\mu}{h^2}\int_{\partial B_1} \Big(
2h\nabla w(x,t)\cdot\omega
+2h^2D^2w(x,t)\omega\cdot\omega
\Big)\Big(u(x,t)+h\nabla u(x,t)\cdot\omega\Big)\,d{\mathcal{H}}^{n-1}_\omega\\&&\qquad\qquad\qquad
+o(1)\\&=& \frac{c\mu}{2}\Delta u(x,t) -
\frac{\mu}{h^2}\int_{\partial B_1} \Big(
2h^2\nabla w(x,t)\cdot\omega\,\nabla u(x,t)\cdot\omega
+2h^2D^2w(x,t)\omega\cdot\omega\, u(x,t)\Big)\,d{\mathcal{H}}^{n-1}_\omega
+o(1)\\&=& \frac{c\mu}{2}\Delta u(x,t) -
2\mu \int_{\partial B_1} \Big(
\nabla w(x,t)\cdot\omega\,\nabla u(x,t)\cdot\omega
+D^2w(x,t)\omega\cdot\omega\, u(x,t)\Big)\,d{\mathcal{H}}^{n-1}_\omega\\&&\qquad\qquad\qquad
+o(1)\\&=& \frac{c\mu}{2}\Delta u(x,t) -
2\mu \sum_{i=1}^n\int_{\partial B_1} \Big(\partial_i w(x,t)\partial_i u(x,t)\omega_i^2
+\partial_{ii}w(x,t) u(x,t)\omega_i^2\Big)\,d{\mathcal{H}}^{n-1}_\omega
+o(1)
\\&=& \frac{c\mu}{2}\Delta u(x,t) -
2\mu \sum_{i=1}^n\int_{\partial B_1} \partial_i\big(\partial_i w(x,t) u(x,t)\big)\omega_i^2
\,d{\mathcal{H}}^{n-1}_\omega
+o(1)\\&=& \frac{c\mu}{2}\Delta u(x,t) -
\widetilde{c}\mu \sum_{i=1}^n \partial_i\big(\partial_i w(x,t) u(x,t)\big)
+o(1)\\&=& \frac{c\mu}{2}\Delta u(x,t) -
\widetilde{c}\mu \div\big(u(x,t) \nabla w(x,t)\big)
+o(1)
\end{eqnarray*}
for a suitable~$
\widetilde{c}>0$.

Notice that here above we have repeatedly taken advantage of odd symmetry cancellations.
Gathering this and~\eqref{y8oqiwhfesjd8923wfyug3iuweLuhonwfdpinrgNqikedOKS} we thereby find that
\begin{equation}\label{DRIFPAG2}\begin{split}&
\frac{u(x,t+\tau)-u(x,t)}\tau=
\rho \frac{u(x,t+\tau)-u(x,t)}\tau +\mu \frac{u(x,t+\tau)-u(x,t)}\tau\\&\qquad
=\frac{c\rho}{2}\Delta u(x,t) +\frac{c\mu}{2}\Delta u(x,t) -
\widetilde{c}\mu \div\big(u(x,t) \nabla w(x,t)\big)+o(1)\\&\qquad
=\frac{c}{2}\Delta u(x,t) -
\widetilde{c}\mu \div\big(u(x,t) \nabla w(x,t)\big)+o(1).\end{split}\end{equation}
Sending~$\tau\searrow0$ (and thus~$h\searrow0$) we obtain
\begin{equation}\label{CHEMLKSDEQ}
\partial_t u(x,t)=\frac{c}{2}\Delta u(x,t)  -
\widetilde{c}\mu \div\big(u(x,t) \nabla w(x,t)\big).\end{equation}
For constant distributions~$w$ of the chemical attractant, equation~\eqref{CHEMLKSDEQ}
can be seen as a special case (up to normalizing constants) of the heat equation in~\eqref{DAGB-ADkrVoiweLL4re2346ytmngrrUj7}.

Stationary solutions of~\eqref{CHEMLKSDEQ} solve
\begin{equation}\label{OJS-PJDN-0IHGDOIUGDBV02ujrf}
\frac{c}{2}\Delta u(x,t)  -
\widetilde{c}\mu \div\big(u(x,t) \nabla w(x,t)\big)=0.\end{equation}

This is the equation solved by the equilibria corresponding to a biological population in the presence of
a chemotactic factor.

Let us now introduce one last complication in the model by accounting for a possible drift
(as mentioned on page~\pageref{DRIFPAG}). In this situation, we may suppose that, in the time unit,
the population is also moved by the velocity vector~$b$ by a space displacement~$b\tau$.
Regarding this, the population density~$u(x,t+2\tau)$
% that we have previously computed as located at~$x$ at time~$t+\tau$
has in fact moved to the location~$x+b\tau$: correspondingly, to take into account this drift,
one can observe that~$ u(x,t+2\tau)=u(x-b(x,t)\tau,t+\tau)$, and accordingly, recalling~\eqref{DRIFPAG2}
and using that
$$u\big(x-b(x,t)\tau,t+\tau\big)=u(x,t )
-\tau b(x,t)\cdot\nabla u (x,t )+\tau\partial_t u(x,t)
+o(\tau),$$ we have that
\begin{eqnarray*}
&& \frac{u(x,t+2\tau)-u(x,t)}\tau= \frac{u(x,t+2\tau)-u(x,t+\tau)}\tau+ \frac{u(x,t+\tau)-u(x,t)}\tau
\\&&\qquad
=\frac{u(x-b(x,t)\tau,t+\tau)-u(x,t+\tau)}\tau+
\frac{c}{2}\Delta u(x,t) -
\widetilde{c}\mu \div\big(u(x,t) \nabla w(x,t)\big)+o(1)\\&&\qquad= 
\frac{u(x,t )
-\tau b(x,t)\cdot\nabla u (x,t )+\tau\partial_t u(x,t)-
u(x,t)-\tau\partial_t u(x,t) +o(\tau)
}\tau
\\&&\qquad\qquad+
\frac{c}{2}\Delta u(x,t) -
\widetilde{c}\mu \div\big(u(x,t) \nabla w(x,t)\big)+o(1)\\&&\qquad=
-b(x,t)\cdot\nabla u (x,t )+
\frac{c}{2}\Delta u(x,t) -
\widetilde{c}\mu \div\big(u(x,t) \nabla w(x,t)\big)+o(1)
\end{eqnarray*}
Thus, we can pass to the limit and obtain
\begin{equation}\label{OJS-PJDN-0IHGDOIUGDBV02ujrf220}
\partial_t u(x,t)=\frac{c}{2}\Delta u(x,t) -
\widetilde{c}\mu \div\big(u(x,t) \nabla w(x,t)\big)-b(x,t)\cdot\nabla u(x,t).\end{equation}
Stationary solutions of~\eqref{OJS-PJDN-0IHGDOIUGDBV02ujrf220} give rise to the equation
\begin{equation}\label{OJS-PJDN-0IHGDOIUGDBV02ujrf22}
\Delta u(x,t) -
\div\big(u(x,t) \nabla w(x,t)\big)-b(x,t)\cdot\nabla u(x,t)=0,\end{equation}
up to normalizing constants that we have omitted for the sake of simplicity.
In this spirit, \eqref{OJS-PJDN-0IHGDOIUGDBV02ujrf22} models\footnote{The equation
in~\eqref{OJS-PJDN-0IHGDOIUGDBV02ujrf22}, together with the diffusion of the attractant,
reveals the complexity
of the model under consideration which, depending on different parameter thresholds, exhibit a variety of
different stable geometrical patterns, colony formation, blowup mechanisms, flocculation and aggregation phenomena.

All these traits play a fundamental role for life as we know it, since these configurations can arise
from equations such as the ones presented here and possess quite a precise counterpart in nature, when
the corresponding behaviors occur, for example, in response to external stress,
or to lack of resources (sometimes the chemical attractant is produced by the biological species itself,
e.g. by amoebae in case of food scarceness) or in response to predation
(autoaggregation provides protection and defense against predators). These aggregation processes often produce slimy yet durable coatings, called biofilms. Even blowup phenomena
can be exploited by bacteria to enhance their possibility of surviving in hostile environments
(e.g., the individual on the top of the tower formed by the blowup have more chances to be picked up, say, by the wind
or by an external factor, and be deposited possibly in a less hostile environment), so all in all the variety
of patterns exhibited by the solutions of a relatively simple differential equation is in some sense the mathematical counterpart
of the variety of ways in which ``life finds a way'' (as Ian Malcolm utters in ``Jurassic Park'').} the equilibria reached by a biological species possibly subject
to chemotaxis and to an external drift.

Notice that~\eqref{OJS-PJDN-0IHGDOIUGDBV02ujrf22} reduces to~\eqref{OJS-PJDN-0IHGDOIUGDBV02ujrf}
in the absence of drift. 

We have written~\eqref{OJS-PJDN-0IHGDOIUGDBV02ujrf22} in a form which emphasizes its ``divergence structure'':
if instead one wants to highlight the presence of a second order
differential operator, it suffices to notice that~$\div(u \nabla w)=u\Delta w+\nabla w\cdot\nabla u$ to recast~\eqref{OJS-PJDN-0IHGDOIUGDBV02ujrf22} into
\begin{equation}\label{OJS-PJDN-0IHGDOIUGDBV02ujrf22NDF}
\Delta u(x,t) -\big( b(x,t)+\nabla w(x,t)\big)\cdot\nabla u(x,t)- u(x,t)\Delta w(x,t)=0,\end{equation}
which is also a telling expression since it is revealing that the effect of the chemotaxis
is to alter the diffusion with an additional drift term~$\nabla w$
which is favoring the movements towards regions with higher concentrations
of the the attractant.

We refer to~\cite{MR3925816, MR1908418, MR1952568, MR2448428} and the references therein for further information about chemotaxis and its biological implications.
See also e.g.~\cite{MR0210469}, \cite{MR0421700}, \cite[Section~2.2]{MR3469920} and~\cite{MR4316246} for applications of the theory of random walks to game theory.

\subsection{Pattern formation, or how the leopard gets its spots}

\begin{figure}
  \centering
  \includegraphics[width=.6\linewidth]{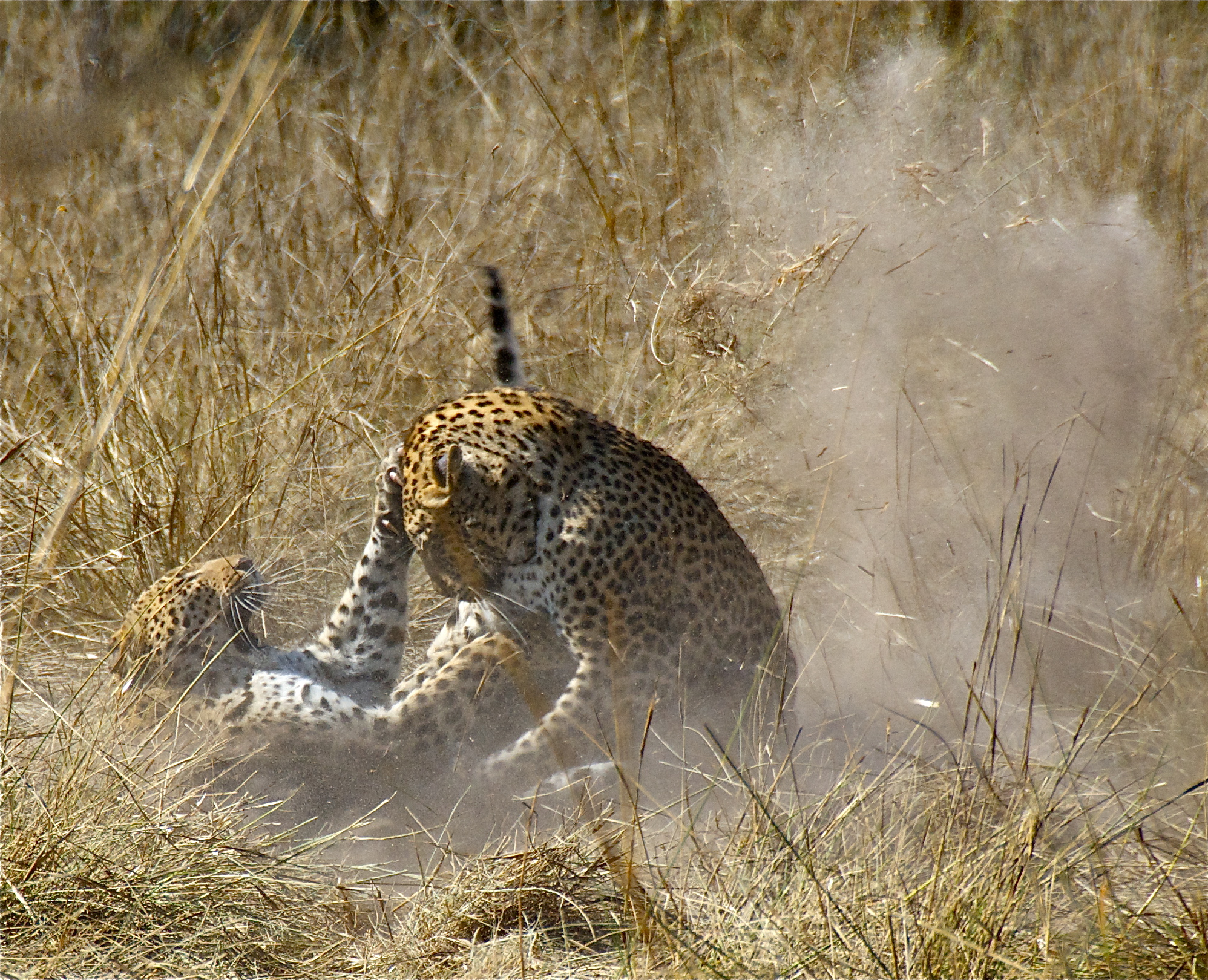}
 \caption{\sl  Leopard mating dance (image from
 Wikipedia, licensed under the Creative Commons Attribution 2.0 Generic license;
see Steve Jurvetson's website {\tt https://www.flickr.com/photos/jurvetson/5913330010/}
for the complete series of shots).}\label{leoP2DAILACRPOGGSTENSIWASTREDItangeFIlASELEO}
\end{figure}

\begin{figure}
  \centering
  \includegraphics[height=.35\linewidth]{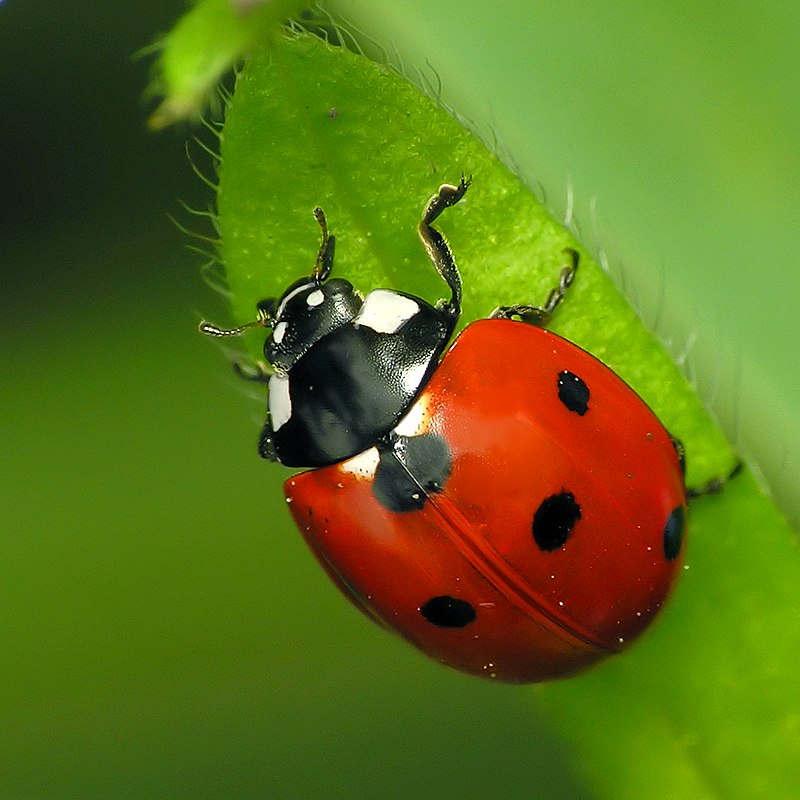}$\quad$
    \includegraphics[height=.35\linewidth]{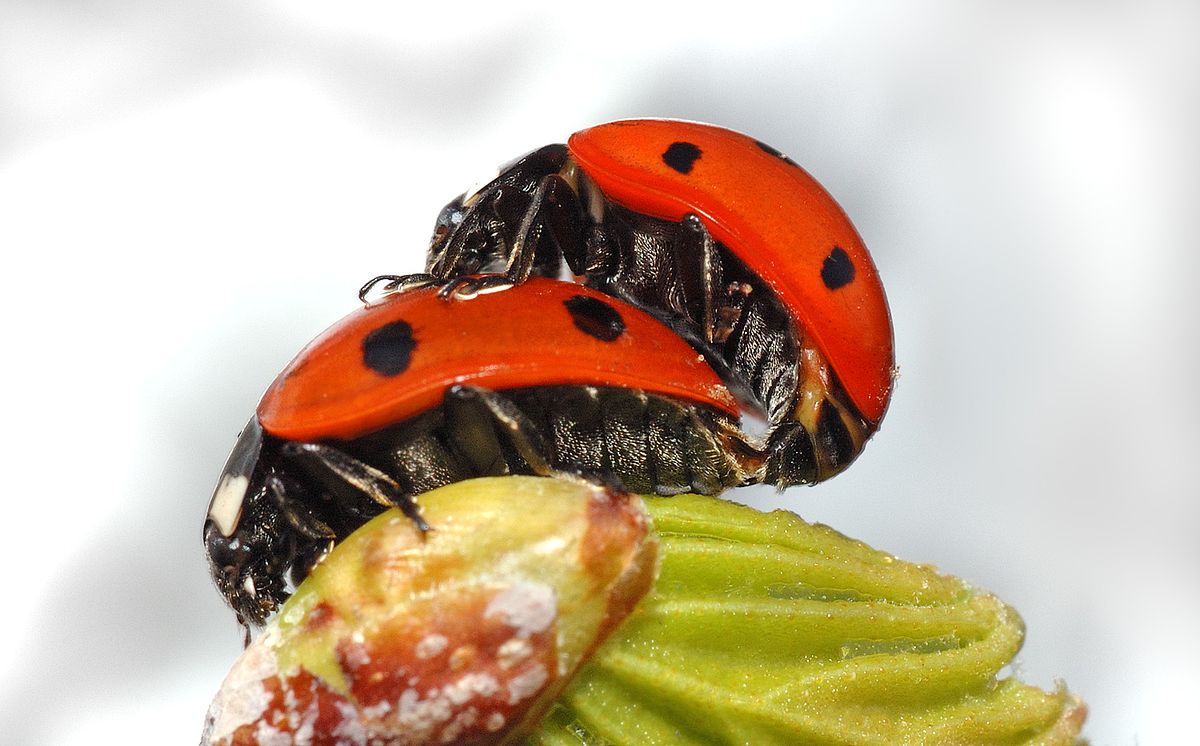}
 \caption{\sl Examples of seven-spotted ladybugs (left:
 photo by Dominik Stodulski; image from
 Wikipedia, licensed under the Creative Commons Attribution-Share Alike 3.0 Unported license; right: photo by Andr\'e Karwath; image from
 Wikipedia, licensed under the Creative Commons Attribution-Share Alike 2.5 Generic license).}\label{P2DAILACRP67890-00OGGSTENSIWASTREDItange3FIlAALANSE}
\end{figure}

\begin{figure}
  \centering
  \includegraphics[height=.33\linewidth]{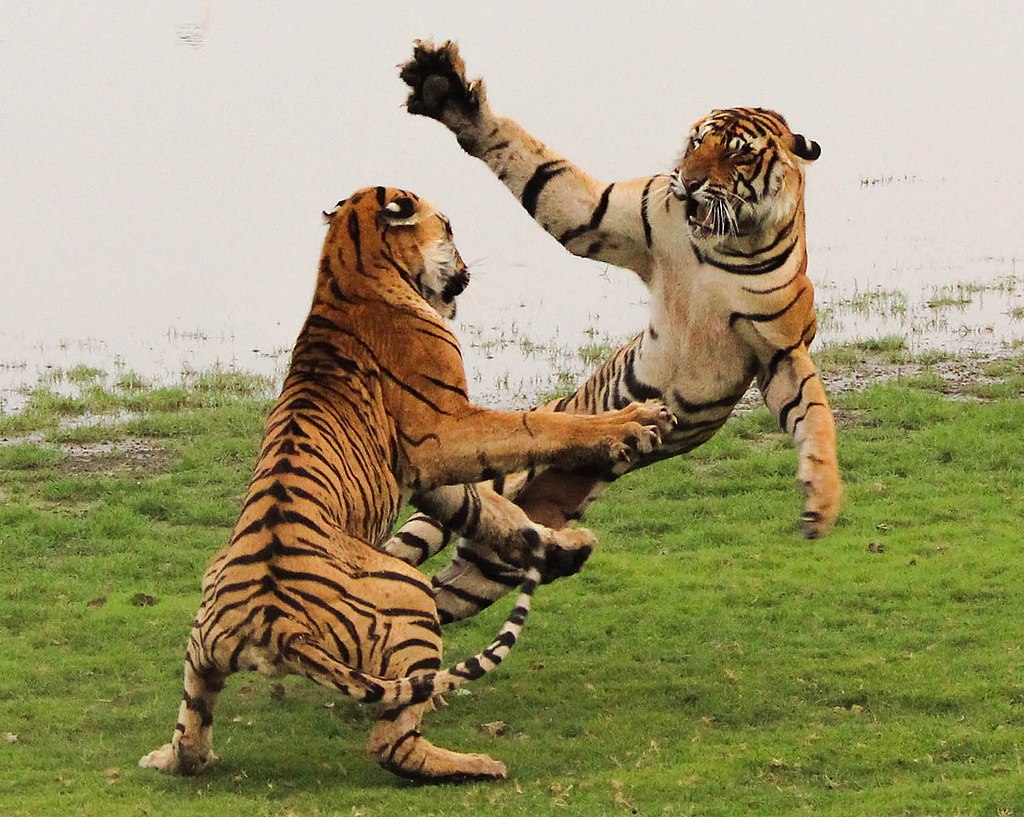}$\quad$
  \includegraphics[height=.33\linewidth]{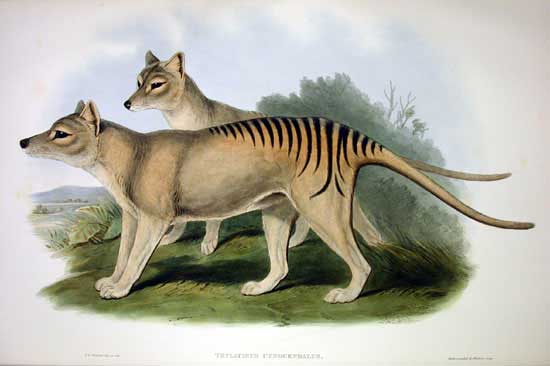}
 \caption{\sl Left: two young female tigers in a playful mood (photo by Vedang Vadalkar; image from
 Wikipedia, licensed under the Creative Commons Attribution-Share Alike 4.0 International license). Right: lithographic plate by John Gould's representing two examples of thylacines, 
 the extinct (or maybe not?) carnivorous marsupials a.k.a. Tasmanian tigers (Public Domain image from
        Wikipedia).}\label{P2DAILACRPOGGSTENSIWASTREDIta888iu34-20o3rjtgrnXnge3FIlAALANSE}
\end{figure}

We know well that the natural world offers a variety of amazing
regular patterns, such as symmetries, tessellations, spirals, spots, etc.,
see e.g. Figures~\ref{leoP2DAILACRPOGGSTENSIWASTREDItangeFIlASELEO}, \ref{P2DAILACRP67890-00OGGSTENSIWASTREDItange3FIlAALANSE}
and~\ref{P2DAILACRPOGGSTENSIWASTREDIta888iu34-20o3rjtgrnXnge3FIlAALANSE}.
These patterns visible in nature
are determinant factors in the processes of natural and sexual selection,
providing organisms with structures recognizable by their conspecifics (thus favoring social
interactions and reproduction) and often realizing optimal configurations for mobility, hunt or camouflage (see e.g.~\cite{MR1398883} and the references therein).

The biological process leading to the development of the specific shape
of an organism is dubbed ``morphogenesis'' \index{morphogenesis}
and a scientific investigation of the formation of patterns in nature
is relatively recent.
One of the pioneer scientists interested in the mathematical analysis of the growth, form and evolution of plants and animals
was Sir D'Arcy Wentworth Thompson: his 1116 pages book~\cite{MR0006348}
(793 pages in the first edition of 1917) combined classical natural philosophy,
biology and mathematics to give insight on a number of biological shapes
and analyzed the differences in the forms in nature in the light of mathematical transformations.

\begin{figure}
  \centering
  \includegraphics[width=.45\linewidth]{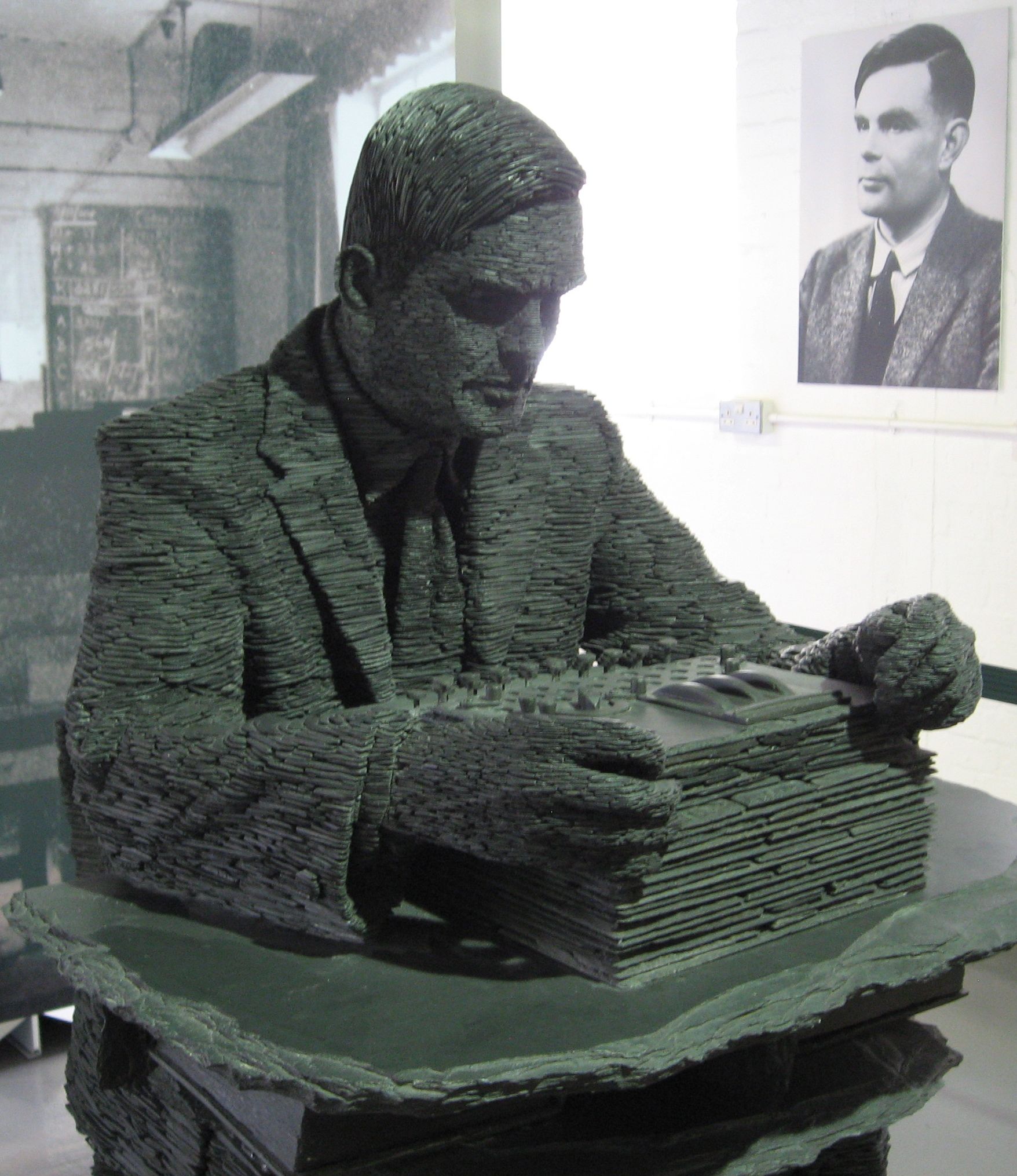}
 \caption{\sl Statue of Alan Turing (photo by Jon Callas; image from
 Wikipedia, licensed under the Creative Commons Attribution 2.0 Generic license).}\label{P2DAILACRPOGGSTENSIWASTREDItange3FIlAALANSE}
\end{figure}

A modern mathematical treatment of morphogenesis was initiated
by\footnote{Not only the father of mathematical morphogenesis,
computer science and artificial intelligence, Alan Turing was also an exceptional long-distance runner, 
capable of world-class marathon standards. See {\tt https://kottke.org/18/04/alan-turing-was-an-excellent-runner}
for a nice photo of Turing during one of his running performances.

During World War II, Turing worked as a code breaker for the British Intelligence,
devising an electromechanical machine, named ``Bombe'',
to decipher the encrypted messages created by the German cipher device ``Enigma''.
The headquarter of the British code breakers was located in a
mansion named Bletchley Park which hosts today a statue of Turing
made with slates (by Stephen Kettle). In this statue, Turing is depicted seated and looking at a German Enigma machine,
see Figure~\ref{P2DAILACRPOGGSTENSIWASTREDItange3FIlAALANSE}.

Turing was~39 years old when he was charged with ``gross indecency'' on the basis of a homosexual relationship
under Section~11 of the Criminal Law Amendment Act~1885.
The case, Regina v Turing and Murray, was brought to trial on 31 March 1952
and Turing pleaded guilty, insisting that he saw nothing wrong with his actions.

He was convicted and given a choice between imprisonment and chemical castration via synthetic estrogens.
As a result of choosing the second option, Turing suffered of severe mutations of his body,
was barred from his work for Government Communication Headquarters and
was denied entry into the United States.
On June 8 1954, Turing was found dead of cyanide poisoning with a half-eaten apple
next to his bed lay. The inquest into his death recorded a verdict of suicide, but
the apple was not tested for cyanide.

Turing's favorite fairy tale was Snow White and the Seven Dwarfs.}
Alan Mathison Turing in~1952, see~\cite{MR3363444}.
Turing's brilliant idea is that simple physical laws were sufficient to justify
the shaping of complex patterns, such as
animal markings and arrangement of leaves and florets in plants. 

In a nutshell, Turing proposed that the process of morphogenesis
is regulated by the interaction of chemical substances, called morphogens, which
diffuse through a tissue: these substances could be
hormones, or genes, or any other essence which can act and react in presence of another one.
In practice, see e.g.~\cite{zbMATH07203210}, one of the morphogens may act as an
``activator'', which is self-sustaining and introduces positive feedback, while the other
may play the role of an ``inhibitor'' which tends to suppress the self-amplification of the activator. 
In this interplay, the pattern may be created by the different speeds of diffusion
of the two substances:
namely, the faster diffusion of the inhibitor
can catch up with the activator's  self-replication (that is, roughly speaking, on the one side the activator's capacity of self-replicating could be strong enough to produce
local patches, but the predominant speed of the inhibitor could avoid that these patches grow
incessantly). All in all, the whole process can be thereby considered as a ``diffusion driven instability''. See~\cite[page~76]{MR1952568} for a very clear explanation of the roles of activators and inhibitors, via an example with sweating grasshoppers.

The mathematical formulation of Turing's idea combines the notion of diffusion,
as modeled for instance by the heat equation~\eqref{HJA:A89MAOL8789892},
and that of reaction, taking into account that each morphogen can chemically react with the others
and the effect of this interaction can depend on the concentration of the diffusing substances in the tissue. The combination of reaction and diffusion in this type
of systems of partial differential equations justified the name \label{REDI6NFHARLROA7KOLM789GIJ7solFUMHDNOJHNFOJED231}
of ``reaction-diffusion
equations''. \index{reaction-diffusion equation}

For instance, one can consider the case of two morphogens with density~$u$ and~$v$ respectively
and a system of reaction-diffusion equations of the form
\begin{equation}\label{SYS:READI}
\begin{dcases}
&\partial_t u=\mu\Delta u+f(u,v),\\
&\partial_t v=\nu\Delta v+g(u,v),
\end{dcases}
\end{equation}
where~$\mu$ and~$\nu$ are positive coefficients which model the speed of diffusion
of the chemical substances with density~$u$ and~$v$ respectively.

Now, the formation of patterns as an outcome of~\eqref{SYS:READI}
is, in a sense, not completely intuitive: on the contrary, given the
``democratic'' tendency of the Laplace operator (as discussed on page~\pageref{DEEFFNVST})
one may imagine that the solution~$(u(t),v(t))$, as~$t\to+\infty$, will evolve spontaneously towards
some constant values~$(u_0,v_0)$, which are just the common zeros of~$f$ and~$g$
(that is, such that~$f(u_0,v_0)=g(u_0,v_0)=0$). And this is indeed one of the possible
destination for the solutions of~\eqref{SYS:READI}. However, there is also a more intriguing possibility:
the constant state~$(u_0,v_0)$ may well exist (and it may also be ``stable'' for small
perturbations when~$\mu:=0$ and~$\nu:=0$), but it may be triggered off by random disturbances.
In this situation, when~$\mu>0$ and~$\nu>0$, the solutions may end up drifting away
from the constant~$(u_0,v_0)$. Thus (unless for some reasons the solution diverges)
several interesting patterns may arise, such as oscillations between equilibria,
stationary waves, moving fronts, etc.

Though a full understanding of Turing's theory of morphogenesis goes well beyond the scopes of this
set of notes, following~\cite{MR3363444} one can at least grasp some of the ideas involved. For instance,
one can consider a simple subcase of~\eqref{SYS:READI} in which~$f$ and~$g$ are linear functions
(this simplification can also be inspiring to treat the general case, since, in the vicinity of
the constant equilibrium~$(u_0,v_0)$,
one can try to ``linearize the equation'' to obtain information on its dynamics).
That is, let us consider the system of equations
\begin{equation}\label{52568-345256890hro3}
\begin{dcases}
&\partial_t u=\mu\Delta u+a(u-u_0)+b(v-v_0),\\
&\partial_t v=\nu\Delta v+c(u-u_0)+d(v-v_0),
\end{dcases}
\end{equation}
for some~$a$, $b$, $c$, $d\in\R$. For simplicity, let us also suppose that the problem is set on a circle,
namely~$x\in\R$ and~$u$ and~$v$ are periodic functions of period~$2\pi$.
One can thus look for solutions of~\eqref{52568-345256890hro3} in Fourier Series of the form
\begin{equation}\label{52568-345256890hro4} u(x,t)=u_0+\sum_{j\in\Z} U_j(t)\,e^{ijx}\qquad{\mbox{and}}\qquad
v(x,t)=v_0+\sum_{j\in\Z} V_j(t)\,e^{ijx},\end{equation}
with~$U_j$ and~$V_j$ to be determined.

Substituting~\eqref{52568-345256890hro4} into~\eqref{52568-345256890hro3} we have that
\begin{eqnarray*}&&
\sum_{j\in\Z} \dot U_j(t)\,e^{ijx}=
\partial_t u=\mu\Delta u+a(u-u_0)+b(v-v_0)\\&&\qquad
= -\mu\sum_{j\in\Z} j^2U_j(t)\,e^{ijx}
+a\sum_{j\in\Z} U_j(t)\,e^{ijx}+b\sum_{j\in\Z} V_j(t)\,e^{ijx}
\end{eqnarray*}
and similarly
\[ \sum_{j\in\Z} \dot V_j(t)\,e^{ijx}=
-\nu\sum_{j\in\Z} j^2V_j(t)\,e^{ijx}
+c\sum_{j\in\Z} U_j(t)\,e^{ijx}+d\sum_{j\in\Z} V_j(t)\,e^{ijx}.\]
{F}rom these equations we arrive at
\begin{equation}\label{52568-345256890hro4-1} \left(\begin{matrix} \dot{U}_j(t)\\ \dot{V}_j(t)
\end{matrix}\right)=
\left( \begin{matrix}a-j^2\mu & b\\
c& d-j^2\nu\end{matrix}\right)
\left(\begin{matrix} {U}_j(t)\\ {V}_j(t)
\end{matrix}\right).
\end{equation}
This is a first order ordinary differential equation with constant coefficients.
Hence, we suppose for simplicity that the (possibly complex) eigenvalues of the matrix~$\left( \begin{matrix}a-j^2\mu & b\\
c& d-j^2\nu\end{matrix}\right)$, which we denote by~$\lambda_j$ and~$\Lambda_j$, are distinct.
We also denote by~$w_j$ and~$W_j\in \cOMPL^2$ the corresponding eigenvectors.
With this notation,
we find (see e.g.~\cite[Theorem~3.6]{MR3450069}) that the solutions of~\eqref{52568-345256890hro4-1}
are of the form
\begin{equation}\label{52568-345256890hro4-2} \left(\begin{matrix} {U}_j(t)\\ {V}_j(t)\end{matrix}\right)
=\xi_j \,e^{\lambda_j t} w_j+\Xi_j\, e^{\Lambda_j t} W_j
\end{equation}
for some~$\xi_j$, $\Xi_j\in\cOMPL$.

\begin{figure}
                \centering
                \includegraphics[width=.35\linewidth]{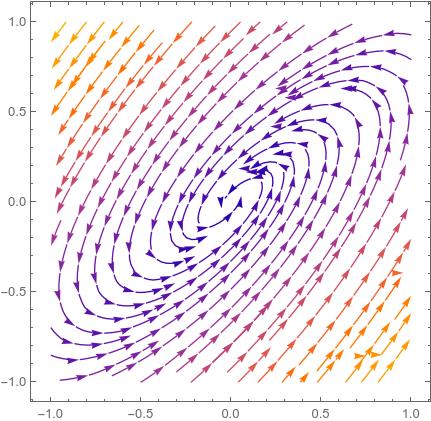}
        \caption{\sl Stream plot of the system of ordinary differential equations in~\eqref{STREAMPLO}.}\label{52568-345256890hro3ODEstre}
\end{figure}

It is convenient to use the vector notation
$$ w_j=\left(\begin{matrix} w_{j1}\\ w_{j2}\end{matrix}\right)\qquad{\mbox{and}}\qquad
W_j=\left(\begin{matrix} W_{j1}\\ W_{j2}\end{matrix}\right)$$
and let
$$ A_j:=\xi_j w_{j1},\qquad B_j:=\Xi_j W_{j1},\qquad
C_j:=\xi_j w_{j2}\qquad{\mbox{and}}\qquad D_j:=\Xi_j W_{j2}.$$
In this way, \eqref{52568-345256890hro4-2} yields that
\begin{equation}\label{UJ:lao-myND9ona9nyuton-0201}
U_j(t)=A_j e^{\lambda_j t}+B_j e^{\Lambda_j t}\qquad{\mbox{and}\qquad}
V_j(t)=C_j e^{\lambda_j t}+D_j e^{\Lambda_j t}.\end{equation}
Moreover,
\begin{equation*}\begin{split}&\left( \begin{matrix}a-j^2\mu & b\\c& d-j^2\nu\end{matrix}\right)
\left( \begin{matrix}A_j\\ C_j\end{matrix}\right)=
\left( \begin{matrix}a-j^2\mu & b\\c& d-j^2\nu\end{matrix}\right)
\left( \begin{matrix}\xi_jw_{j1} \\ \xi_jw_{j2}\end{matrix}\right)
\\&\qquad=\xi_j\left( \begin{matrix}a-j^2\mu & b\\c& d-j^2\nu\end{matrix}\right)w_j=
\xi_j\lambda_jw_j=\lambda_j\left( \begin{matrix}A_j \\ C_j\end{matrix}\right),
\end{split}\end{equation*}
leading to
\begin{equation}\label{UJ:lao-myND9ona9nyuton-0202}
(a-j^2\mu-\lambda_j)A_j+b C_j=0.
\end{equation}
Similarly,
\begin{equation}\label{UJ:lao-myND9ona9nyuton-0203}
(a-j^2\nu-\Lambda_j)B_j+b D_j=0.\end{equation}
In jargon, \eqref{UJ:lao-myND9ona9nyuton-0202} and~\eqref{UJ:lao-myND9ona9nyuton-0203}
are sometimes dubbed ``dispersion relations''\index{dispersion relation}:
their interest lies in the fact that they relate the speed of oscillation in the time variable
(quantified in~\eqref{UJ:lao-myND9ona9nyuton-0201}
by the eigenvalues~$\lambda_j$ and~$\Lambda_j$) with the spatial periodicity of the medium
(characterized by the eigenvalues~$-j^2$ of the one-dimensional Laplacian and modulated
by the speeds of diffusion~$\mu$ and~$\nu$).

To recap briefly, from the system of reaction-diffusion equations in~\eqref{52568-345256890hro3}
one arrives at the solutions introduced in~\eqref{52568-345256890hro4},
with~$U_j$ and~$V_j$ as in~\eqref{UJ:lao-myND9ona9nyuton-0201},
where
\begin{equation}\label{o2re2s3p4on21di21ng214ly}
{\mbox{$\lambda_j$ and~$\Lambda_j$ are (distinct)
complex eigenvalues of the matrix }}\left( \begin{matrix}a-j^2\mu & b\\
c& d-j^2\nu\end{matrix}\right),\end{equation}
and with the parameters satisfying~\eqref{UJ:lao-myND9ona9nyuton-0202}
and~\eqref{UJ:lao-myND9ona9nyuton-0203}.

As detailed in~\cite[Sections~7 and~8]{MR3363444}
and~\cite[Chapters~2 and~3]{MR1952568}, this explicit mathematical construction
has a number of important biological consequences and presents a sufficiently rich structure
to account for many patterns visible in nature. To see these features, one may focus
on the case in which one of the eigenvalues has the largest real part (roughly speaking,
one expects that the other modes are dominated by this one). Also,
it is convenient to distinguish between the case in which the dominant eigenvalue is real
from the one in which it is complex and with nonzero imaginary part: indeed, real eigenvalues will be related to stationary states
and complex eigenvalues to oscillatory cases.

More specifically, suppose that
\begin{equation}\label{CMOPSDEArajeke24w}
{\mbox{$\Lambda_{j_0}$ is the eigenvalue with largest real part.}}\end{equation}
We notice that also~$\Lambda_{-j_0}=\Lambda_{j_0}$, since 
\begin{equation}\label{INVAMAJJ}
{\mbox{the matrix }}\left( \begin{matrix}a-j^2\mu & b\\ c& d-j^2\nu\end{matrix}\right) {\mbox{ remains the same if we exchange~$j$ with~$-j$}}.\end{equation}
Hence, dropping the higher order terms, we can assume that the dynamics of the solutions in~\eqref{52568-345256890hro4} is governed by the following long-time asymptotics:
\begin{equation}\label{INVAMAJJ2}
\begin{split}&
u(x,t)\simeq u_0+ U_{j_0}(t)\,e^{i{j_0}x}+U_{-j_0}(t)\,e^{-i{j_0}x}
\simeq u_0+ e^{\Lambda_{j_0} t}\Big( B_{j_0} e^{i{j_0}x}+B_{-j_0} e^{-i{j_0}x}\Big)\\
{\mbox{and}}\qquad&
v(x,t)\simeq v_0+ e^{\Lambda_{j_0} t}\Big( D_{j_0} e^{i{j_0}x}+D_{-j_0} e^{-i{j_0}x}\Big).\end{split}
\end{equation}
The invariance in~\eqref{INVAMAJJ} also suggests that if~$A_j$, $B_j$, $C_j$ and~$D_j$ are solutions of~\eqref{UJ:lao-myND9ona9nyuton-0202} and~\eqref{UJ:lao-myND9ona9nyuton-0203}, then so are~$A_{-j}$, $B_{-j}$, $C_{-j}$ and~$D_{-j}$: for this reason, we can suppose that~$B_{-j_0}=B_{j_0}$
and~$D_{-j_0}=D_{j_0}$
in~\eqref{INVAMAJJ2}, obtaining that
\begin{equation}\label{INVAMAJJ3}
\begin{split}&
u(x,t)\simeq u_0+ B_{j_0} e^{\Lambda_{j_0} t}\Big( e^{i{j_0}x}+e^{-i{j_0}x}\Big)=
u_0+ 2B_{j_0} e^{\Lambda_{j_0} t}\cos( j_0x)\\
{\mbox{and}}\qquad&
v(x,t)\simeq v_0+2D_{j_0} e^{\Lambda_{j_0} t}\cos( j_0x).\end{split}
\end{equation}

Without aiming at exhausting all the possible patterns included in~\eqref{INVAMAJJ3},
let us now show a concrete case of interest.

\begin{figure}
                \centering
                \includegraphics[width=.45\linewidth]{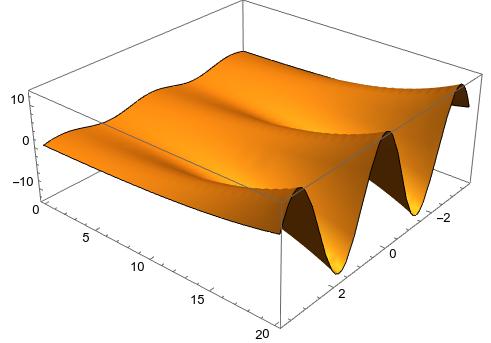}
        \caption{\sl Plot of the function~$(x,t)\mapsto\exp\left(\frac{264\left(7 -4\sqrt{3}\right)}{11\left(13\sqrt{3}-9\right)}\,t\right)\cos (2 x)$.}\label{23COSInMHDNOJHNFOJED231}
\end{figure}

For instance, let~$\vartheta\in\N$, \begin{equation}\label{G1A2MbdRbnU2U8n4ngVA9}\gamma:=
\frac{33\,\vartheta^2}{13\sqrt{3}-9}
\end{equation} and also
\begin{equation}\label{MSUOJS5UIHJ6OR6YUFH8YTFUYHGJFDTYG}
u_0:=0,\qquad v_0:=0, \qquad a:=\gamma, \qquad b:=-2\gamma,\qquad c:=2\gamma\qquad {\mbox{and}}\qquad d:=-2\gamma. \end{equation}
In this case,
the system in~\eqref{52568-345256890hro3} describes an activator with density~$u$
which is an autocatalytic activator: that is, such a substance stimulates the production
of itself (since~$a>0$) and also activates the production of the substance with density~$v$ (since also~$c>0$). Also, the substance with density~$v$ corresponds to a self-degrading inhibitor:
indeed, higher concentrations of this reactant are noxious for itself (since~$d<0$)
and for the activator with density~$u$ (since~$b<0$).

Interestingly, the origin, which corresponds to the equilibrium~$(u_0,v_0)$, is a stable sink for the system of ordinary differential equations corresponding to~\eqref{52568-345256890hro3} when~$\mu:=0$ and~$\nu:=0$. See Figure~\ref{52568-345256890hro3ODEstre} for a sketch of the trajectories
of
\begin{equation}\label{STREAMPLO}
\begin{dcases}
&\partial_t u= \gamma u -2\gamma v,\\
&\partial_t v=2\gamma u-2\gamma v.
\end{dcases}
\end{equation}

Quite remarkably, as discovered by Turing, the stability of~\eqref{STREAMPLO} can be destroyed by random fluctuations arising from the diffusivity of the chemical reactant. To appreciate this, given~$u_0$, $v_0$, $a$, $b$, $c$ and~$d$ as in~\eqref{MSUOJS5UIHJ6OR6YUFH8YTFUYHGJFDTYG},
we take~$\mu:=1$ and~$\nu:=12$ in~\eqref{52568-345256890hro3}. Note that this corresponds to a situation in which the diffusion of the inhibitor is faster than the one of the activator.
This scenario gives that
$$ \left( \begin{matrix}a-j^2\mu & b\\ c& d-j^2\nu\end{matrix}\right)=
\left( \begin{matrix}\gamma-j^2 & -2\gamma\\ 2\gamma& -2\gamma- 12j^2\end{matrix}\right),$$
which possesses eigenvalues of the form
\begin{equation}\label{CMOPSDEArajeke24w2}-\frac12\left(\gamma+13 j^2  \pm \sqrt{(11j^2+7\gamma)(11j^2-\gamma)}\right).\end{equation}
By~\eqref{CMOPSDEArajeke24w}, we aim at detecting the greatest possible real part in~\eqref{CMOPSDEArajeke24w2}. To this end, note that when~$(11j^2+7\gamma)(11j^2-\gamma)\le0$
then the real part in~\eqref{CMOPSDEArajeke24w2} is equal to~$-\frac12(\gamma+13 j^2)<0$.
Instead, if~$(11j^2+7\gamma)(11j^2-\gamma)>0$ the largest possible real part in~\eqref{CMOPSDEArajeke24w2} is equal to
\begin{equation}\label{LESOCHEDIBMECOHowrj-0}\begin{split}&
\sup_{j\in\Z}\frac12\left(\sqrt{(11j^2+7\gamma)(11j^2-\gamma)}
-\gamma-13 j^2\right)\\=\,&\frac\gamma2\sup_{j\in\Z}\left(
\sqrt{\left(\frac{11j^2}\gamma+7\right)\left(\frac{11j^2}\gamma-1\right)
}-1-\frac{13 j^2}\gamma\right)\\=\,&
\frac\gamma2\sup_{j\in\Z}\Phi\left(\frac{j^2}\gamma\right),\end{split}\end{equation}
where
$$ \Phi(\tau):=\sqrt{(11\tau+7)(11\tau-1)}-1-13 \tau.$$
Using elementary calculus, one checks that
$$ \max_{\tau\ge0}\Phi(\tau)= \frac4{11} \left(7 -4\sqrt{3}\right)
=\Phi\left(\frac{13\sqrt{3}-9}{33}\right)=\Phi\left(\frac{\vartheta^2}\gamma\right),$$
thanks to~\eqref{G1A2MbdRbnU2U8n4ngVA9}.

\begin{figure}
                \centering
                \includegraphics[width=.27\linewidth]{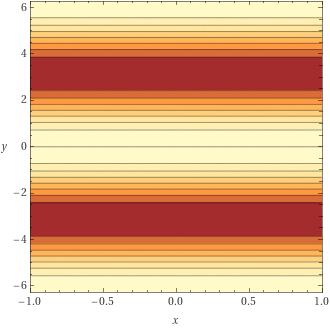} $\quad$ \includegraphics[width=.27\linewidth]{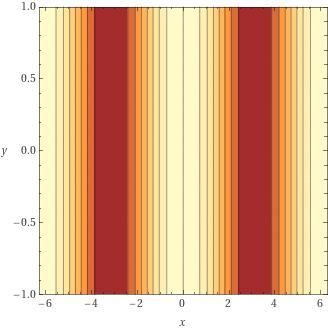} $\quad$
                \includegraphics[width=.27\linewidth]{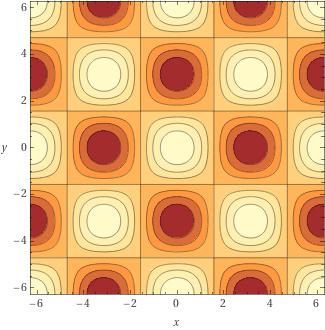}
        \caption{\sl Level sets of~$(x_1,x_2)\mapsto2\cos x_2$, $(x_1,x_2)\mapsto2\cos x_1$ and~$(x_1,x_2)\mapsto2\cos x_1\cos x_2$.}\label{Cos2HAFODAKEDFUMSldRGIRA4AXELEUMDJOMNFHARLROA7789GIJ7solFUMHDNOJHNFOJED23COSEN1}
\end{figure}

This observation and~\eqref{LESOCHEDIBMECOHowrj-0} give that the eigenvalues with largest possible real part in~\eqref{CMOPSDEArajeke24w2} correspond to the choice~$j:=\vartheta$ and are of the form
$$ \frac{2\gamma}{11} \left(7 -4\sqrt{3}\right)
=\frac{66\left(7 -4\sqrt{3}\right)\,\vartheta^2}{11\left(13\sqrt{3}-9\right)}
.$$
Hence, up to constants, the corresponding setting in~\eqref{INVAMAJJ3} takes the form
\begin{equation}\label{LESOCHEDIBMECOHowrj}
\begin{split}&
u(x,t)\simeq  \exp\left(\frac{66\left(7 -4\sqrt{3}\right)\,\vartheta^2}{11\left(13\sqrt{3}-9\right)}\,t\right)\cos(\vartheta x)\\
{\mbox{and}}\qquad&
v(x,t)\simeq \exp\left(\frac{66\left(7 -4\sqrt{3}\right)\,\vartheta^2}{11\left(13\sqrt{3}-9\right)}\,t\right)\cos (\vartheta x) .\end{split}
\end{equation}
See Figure~\ref{23COSInMHDNOJHNFOJED231}
for an example with~$\vartheta:=2$.

\begin{figure}
                \centering
                \includegraphics[width=.55\linewidth]{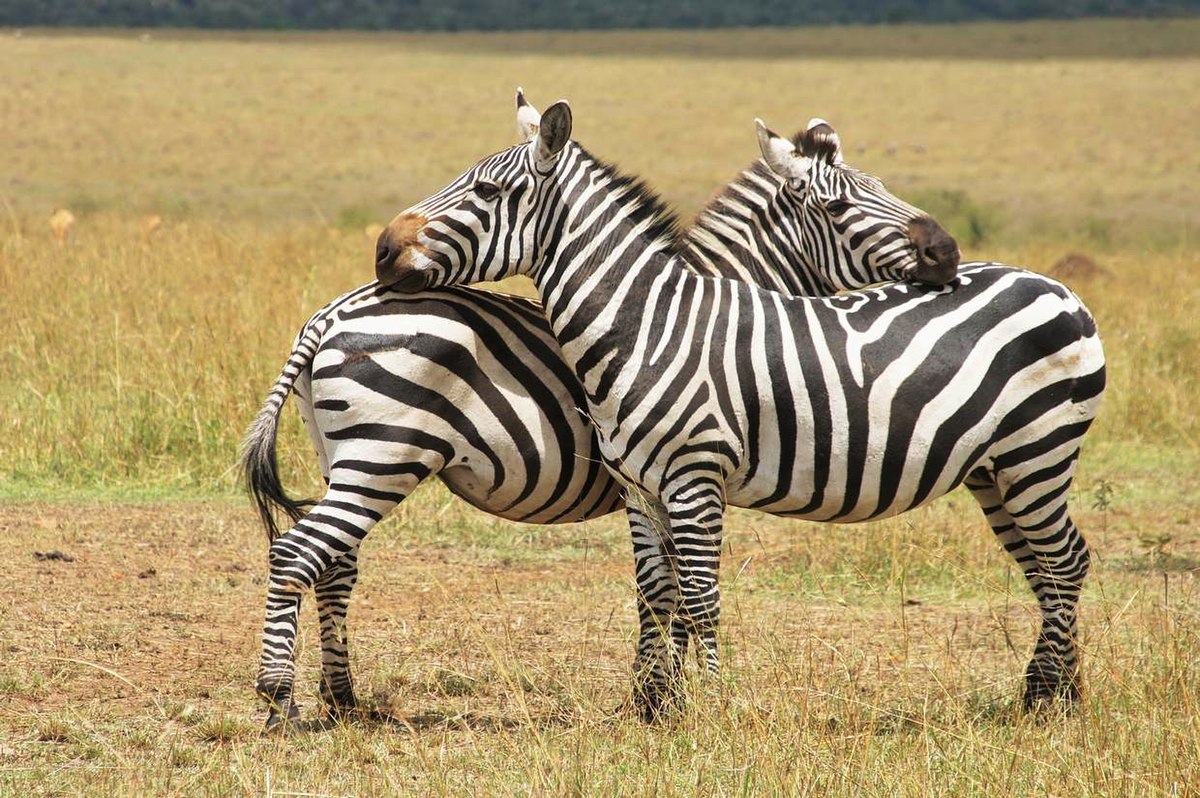}
        \caption{\sl Mutually grooming zebras (photo by Duvignau Alain from
        Wikipedia, available
        under the Creative Commons CC0 1.0 Universal Public Domain Dedication).}\label{ZEB2HAFODAKEDFUMSldRGIRA4AXELEUMDJOMNFHARLROA7789GIJ7solFUMHDNOJHNFOJED231}
\end{figure}

Of course, the asymptotics in~\eqref{LESOCHEDIBMECOHowrj} are divergent as~$t\to+\infty$,
which would correspond to the chemical substances to reach infinite density, which is certainly unfeasible in practice,
hence the meaning of~\eqref{LESOCHEDIBMECOHowrj} has to be understood only in the vicinity of the equilibrium~$(u_0,v_0)$,
which was set to be the origin for simplicity. Indeed, in practice the linear system in~\eqref{52568-345256890hro3}
must be considered as an efficient linearization only in the vicinity of the equilibrium~$(u_0,v_0)$,
while the general situation is more accurately described by a nonlinear system as in~\eqref{SYS:READI}.
For practical purposes, the nonlinear sources~$f$ and~$g$ would force a bound on the densities~$u$ and~$g$
and possibly favor the convergence of the solution for large times to steady solutions~$u=u(x)$ and~$v=v(x)$ of
\begin{equation}\label{SYS:READI:STEAD}
\begin{dcases}
&\mu\Delta u(x)+f(u(x),v(x))=0,\\
&\nu\Delta v(x)+g(u(x),v(x))=0.
\end{dcases}
\end{equation}
That is, for small times, the linear mechanism identified in~\eqref{LESOCHEDIBMECOHowrj} is helpful to detect how diffusion can drive
instability and place the kinetics of the system out of the ``trivial'' state~$(u_0,v_0)$; then, at a longer time scale,
the nonlinear structure of~\eqref{SYS:READI} becomes instrumental to confine the solution and lead it towards spatially
inhomogeneous patterns, as described by the steady solutions of~\eqref{SYS:READI:STEAD}
(as another option, the evolution of the equation may be stopped after a certain time in case
the release of the chemical substances stops; this could be the case in which
the pattern is formed at an embryonic stage for the animal due to chemical substances that are released only
during specific periods of the early stage development of an organism).

We also remark that, for the sake of simplicity, here we confined ourselves to the case in which the spatial domain
is a circle (i.e., the real line with periodic assumptions in~$x$): in general, if one considers more complicated domains
(say, closer to biological situations of specific interest) then it is convenient to replace~$e^{ijx}$ in~\eqref{52568-345256890hro4}
with the eigenfunctions
of the Laplacian in the domain of interest (with the corresponding boundary conditions).
In like manner, the dispersion relations in~\eqref{UJ:lao-myND9ona9nyuton-0202} and~\eqref{UJ:lao-myND9ona9nyuton-0203}
must take into account the corresponding eigenvalues in the place of~$-j^2$.
These eigenvalues replace~$-j^2$ in the matrix in~\eqref{o2re2s3p4on21di21ng214ly} too, and the diffusion eigenvalues~$\lambda_j$ and~$\Lambda_j$ must be modified accordingly.
The diffusion eigenvalues corresponding to diffusive instability still correspond to the ones with positive real part,
and one can focus for concreteness on the diffusive eigenvalue with largest real part in~\eqref{CMOPSDEArajeke24w}.
The structural difference in the general case is however that the excited modes, that is the diffusive
eigenvalues with positive real part, depend on the domain. Since their region of positive and negative values
create a visible pattern in the dynamics of the solution (compare with the positive and negative values of the solution
depicted in Figure~\ref{23COSInMHDNOJHNFOJED231}), it is conceivable that these regions have a connection with the visible
patterns in nature. For instance, two-dimensional
rectangular domains present eigenfunctions of the type~$\cos\left(\frac{2\pi j_1 x_1}{\ell_1}\right)\cos\left(\frac{2\pi j_2 x_2}{\ell_2}\right)$, where~$\ell_1$ and~$\ell_2$ account for the lengths of the side of the rectangle
and~$j_1$, $j_2\in\N$. In this scenario, excited states corresponding to~$j_1=0$ give rise to horizontal stripes,
the ones corresponding to~$j_2=0$ to vertical stripes, and the ones with~$j_1\ne0$ and~$j_2\ne0$ to maculate patterns, see
e.g. Figure~\ref{Cos2HAFODAKEDFUMSldRGIRA4AXELEUMDJOMNFHARLROA7789GIJ7solFUMHDNOJHNFOJED23COSEN1}.

With respect to this, notice also that different excite modes correspond to different size of the pattern (such as width of the
stripes
or possible elongation of the spots). Of course, the analysis of simple rectangular regions is insufficient to capture the whole complexity of animal patterns: yet, it is suggestive to ``approximate'' an animal's coat with rectangular regions and
observe how for instance the orientation of stripes ``locally'' follow the proportion of the approximating rectangles,
see e.g. Figure~\ref{ZEB2HAFODAKEDFUMSldRGIRA4AXELEUMDJOMNFHARLROA7789GIJ7solFUMHDNOJHNFOJED231}
which clearly shows a change between vertical and horizontal stripes patterns in zebras
at the junctions between the body and the legs and between the body and the tail.

\begin{figure}
                \centering
                \includegraphics[width=.85\linewidth]{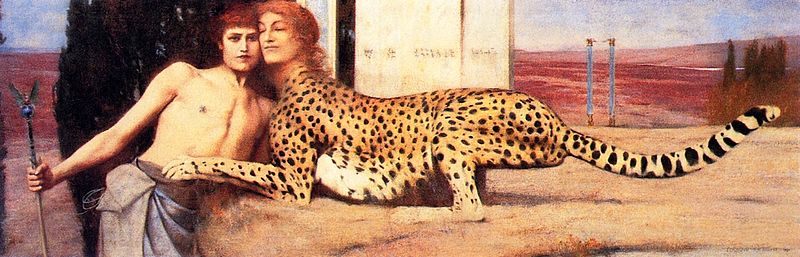}
        \caption{\sl Caress of the Sphinx, by Fernand Khnopff (Public Domain image from
        Wikipedia).}\label{23HAFODAK2EDka4nFUw1MSldRGSPHIIJ7solFUMHDNOJHNFOJED231}
\end{figure}

\begin{figure}
                \centering
                \includegraphics[width=0.98\linewidth]{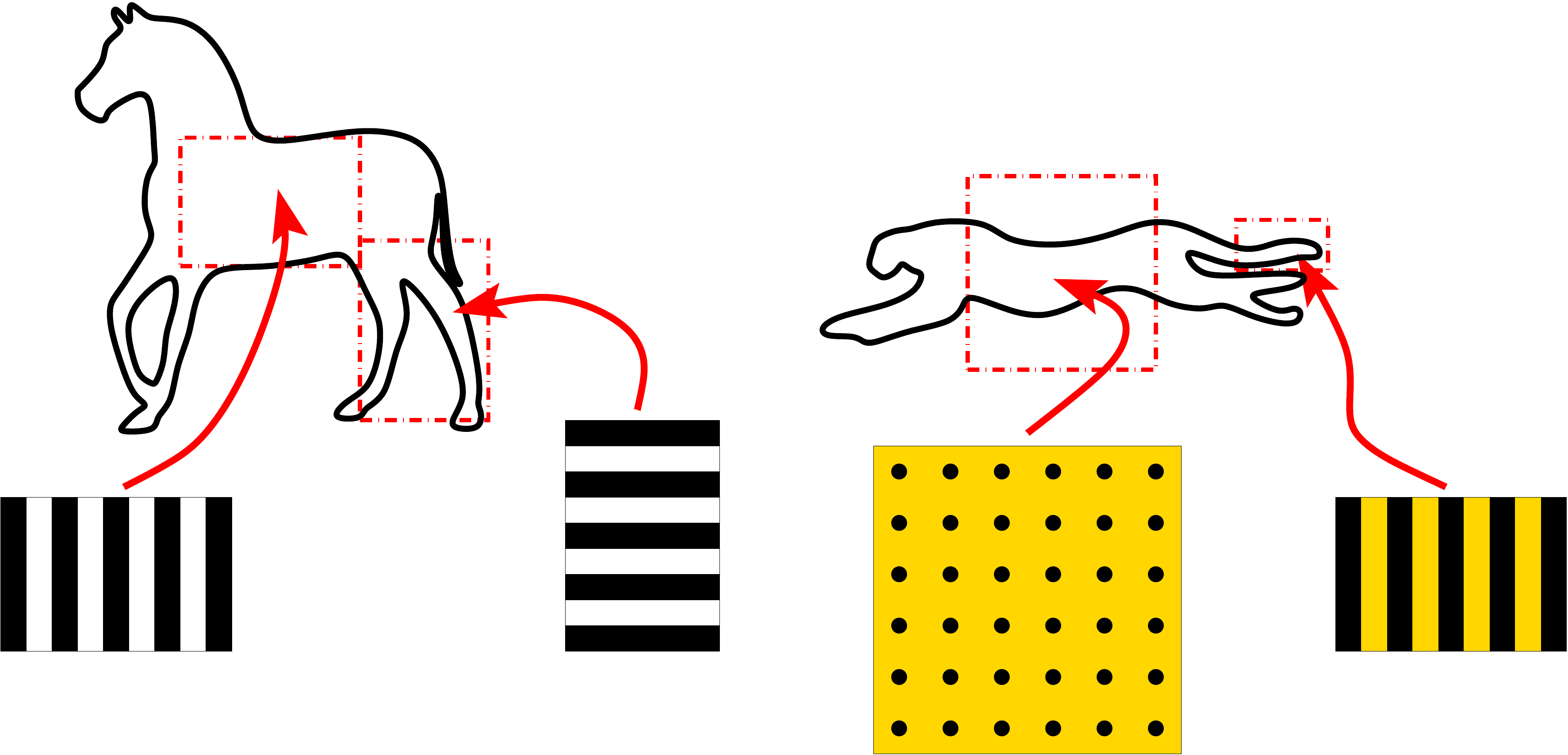}
        \caption{\sl Change of visual markers in different areas of zebras and cheetahs.}\label{23HAFChamd:lFUMHDNOJHNFOJED231}
\end{figure}

Animal patterns may also bifurcate from stripes to spots in different regions\footnote{It is indeed quite tempting
to divide the whole fur of zebras and cheetahs into ``rectangular'' regions to approximately reproduce the
eigenfunction patterns depicted in Figure~\ref{Cos2HAFODAKEDFUMSldRGIRA4AXELEUMDJOMNFHARLROA7789GIJ7solFUMHDNOJHNFOJED23COSEN1}.
Of course, such a reduction is a bit simplistic, since these rectangles do overlap
a bit in regions which cannot be accurately described by the simple eigenfunction patterns
(equivalently, to understand the nodal lines of eigenfunctions, it is too crude
to approximate a complicated domain with nonoverlapping rectangles,
not only because these rectangles may not shade faithfully the precise contour of the domain, but also,
and most importantly, because a rigid rectangular subdivision imposes new artificial boundary conditions): in any case,
an attempt to perform such a rectangular reduction is sketched in Figure~\ref{23HAFChamd:lFUMHDNOJHNFOJED231}.
See also Figures~3.2, 3.6 and~3.8 in~\cite{MR1952568}
for accurate simulations accounting for realistic body surfaces.

On a different note, a common curiosity question is:
are zebras white animals with black stripes, or black animals with white stripes?
Well, the most reasonable answer is that arguably zebras are black with white stripes.
Indeed, zebra's fur grows from follicles that contain melanocyte cells, which generate the pigment that gives color to hair.
In this context, zebra's white fur represents an absence of melanin.
That is, white is not really the pigment of the stripe
and white stripes only exist because pigment is denied.}
of the organism, possibly in response of a variation
of the ``local geometry'', compare e.g. with the tail markings of the cheetah in which the typical spots of the animal
become stripes at the tip of the tail, see Figure~\ref{23HAFODAK2EDka4nFUw1MSldRGSPHIIJ7solFUMHDNOJHNFOJED231}.\medskip

After~\cite{MR3363444},
Turing's ideas about morphogenesis have become a cornerstone
in mathematical biology and led to a number of fantastic accomplishments, see e.g.
the mathematical reconstruction of very sophisticated animal patterns
presented in~\cite{LEOPARD, MR1952568}.
Also, the name of Turing pattern\index{Turing pattern}
is nowadays commonly used to denote biological structures
and specifically those arising by the imbalances between diffusion rates of different chemical agents which make a stable system sensitive to perturbations.

See~\cite[Chapters~2 and~3]{MR1952568}, \cite[Section~12.5]{MR2573923}, \cite{MR3408563}
and the references therein for a thorough analysis of diffusion driven pattern formation,
with several remarkable examples inspired by concrete natural phenomena.
See also~\cite{MR688146} and the reference therein for more information about
reaction-diffusion equations.

\subsection{Space invaders}

\begin{figure}
                \centering
                \includegraphics[width=.17\linewidth]{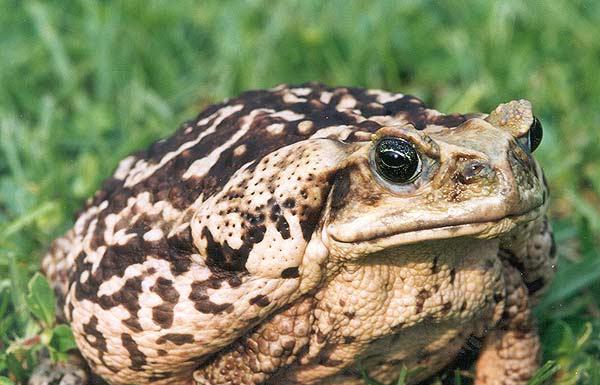}$\quad$\includegraphics[width=.17\linewidth]{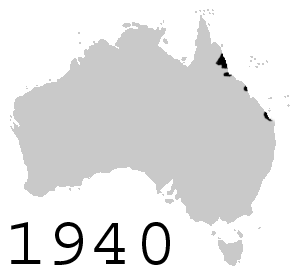}
                $\quad$\includegraphics[width=.17\linewidth]{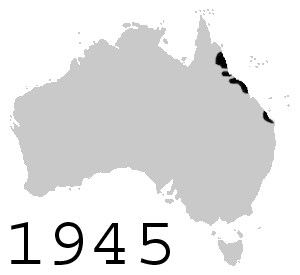}$\quad$\includegraphics[width=.17\linewidth]{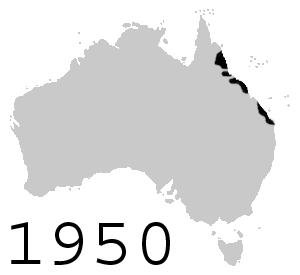}
                $\quad$\includegraphics[width=.17\linewidth]{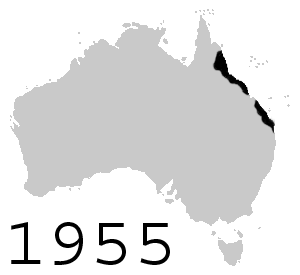}\\ \includegraphics[width=.17\linewidth]{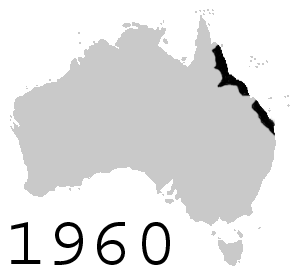}$\quad$
                \includegraphics[width=.17\linewidth]{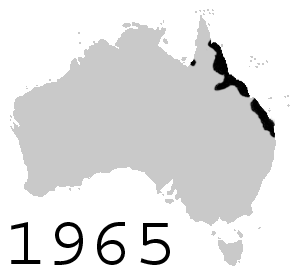}$\quad$\includegraphics[width=.17\linewidth]{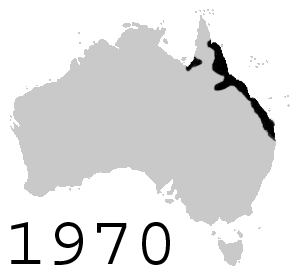}
                $\quad$\includegraphics[width=.17\linewidth]{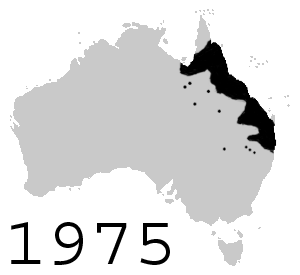}$\quad$\includegraphics[width=.17\linewidth]{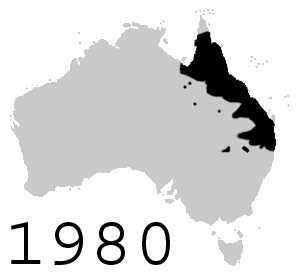}
        \caption{\sl A cane toad (bufo marinus) and the spread of cane toads in Australia from 1940 to 1980 in five-year intervals
        (photo by Bill Waller and animated map by  Froggydarb; images from
        Wikipedia, licensed under the Creative Commons Attribution-Share Alike 3.0 Unported license).}\label{2HAFORANE7solFUMHDNOJHNFOJED231}
\end{figure}

\begin{figure}
                \centering
                \includegraphics[width=.17\linewidth]{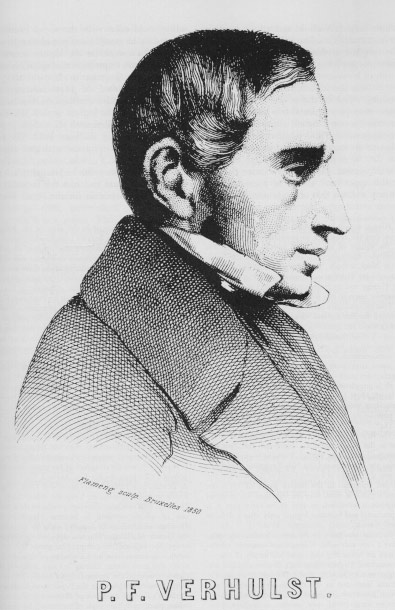}
        \caption{\sl Pierre Verhulst (Public Domain image from
        Wikipedia).}\label{2HAFORABEVErDNOJHNFOJED231}
\end{figure}

A topical argument in mathematical biology consists in the study of invasive species and in their territorial colonization ability,
see e.g. Figure~\ref{2HAFORANE7solFUMHDNOJHNFOJED231}.
To introduce ourselves to this subject, it is first opportune to understand the so-called \index{logistic equation}
logistic equation. This model was introduced by Pierre Fran\c{c}ois Verhulst,
see Figure~\ref{2HAFORABEVErDNOJHNFOJED231} for his portrait,
and describes the evolution in time of a biological population.
The number of individuals~$N(t)$ is supposed to grow at an intrinsic growth rate, parameterized by a given~$\rho>0$,
which somewhat accounts for the ideal birth rate of the population ($N$ individuals would give rise to~$\rho N$ newborns in the unit of time). Additionally, the population undergoes an intraspecific competition which causes the death of some individuals
due to possible overcrowding (assuming, for instance,
that the environmental resources only allow a maximum number of individuals~$K$). The combination of these effects
lead to the logistic ordinary differential equation
\begin{equation}\label{HJA:A89MAOL8789892BIOEQLOG}
\frac{dN}{dt}(t)=\rho N(t)\left(1-\frac{N(t)}K\right).
\end{equation}

This equation considers all the population located basically at the same place, but, for many practical purposes,
it is also convenient to describe biological individuals in a spatial environment, in which case the function~$N$
depends on the time variable~$t$ and also on the space variable~$x$. If the population performs some kind of random walk,
as discussed in Section~\ref{CHEMOTX} (see in particular equation~\eqref{HJA:A89MAOL8789892}), combining~\eqref{HJA:A89MAOL8789892BIOEQLOG} with the random diffusive tendency of the population we arrive at the
partial differential equation
\begin{equation}\label{HJA:A89MAOL8789892BIOEQLOG-2}
\partial_t N(x,t)=c\Delta N(x,t)+\rho N(x,t)\left(1-\frac{N(x,t)}K\right),
\end{equation}
for some diffusion coefficient~$c>0$.

\begin{figure}
  \centering
  \includegraphics[width=.6\linewidth]{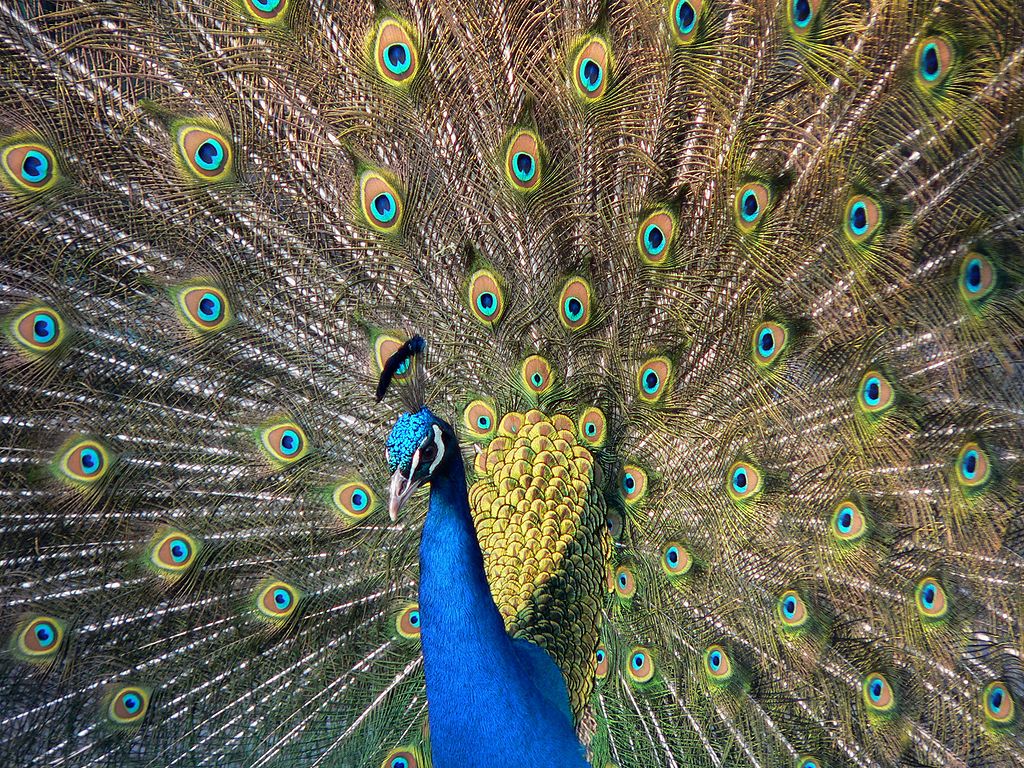}
 \caption{\sl A blue peacock presenting its feathers (image by from
 Wikipedia, available under the Creative Commons Attribution-Share Alike 3.0 Unported).}\label{P1232DAILACRpeacWASTREDItangeFIlASE}
\end{figure}

We can actually simplify~\eqref{HJA:A89MAOL8789892BIOEQLOG-2} via a simple rescaling: namely, defining
\begin{equation}\label{HJA:A89MAOL8789892BIOEQLOG-4}
u(x,t):=\frac1K\,N\left(\sqrt{\frac{c}\rho}x,\frac{t}{\rho}\right),\end{equation}
equation~\eqref{HJA:A89MAOL8789892BIOEQLOG-2} reduces to
\begin{equation}\label{HJA:A89MAOL8789892BIOEQLOG-3}
\partial_t u=\Delta u+u(1-u),
\end{equation}
which is dubbed\footnote{Equation~\eqref{HJA:A89MAOL8789892BIOEQLOG-3}
is named after Sir Ronald Aylmer Fisher~\cite{zbMATH02523582}, 
Andrey Nikolaevich Kolmogorov, Ivan Georgievich Petrovsky and Nikolai Semenovich Piskunov~\cite[pages 248--270]{MR1175399}.

Though biological invasions and
the spread of a genetic trait were
the prime motivations for the study of~\eqref{HJA:A89MAOL8789892BIOEQLOG-3},
nowadays this equation is also intensively studied for its applications
in combustion and flame propagation: that's the unifying power of
the good mathematics!

Besides the introduction of equation~\eqref{HJA:A89MAOL8789892BIOEQLOG-3} and some questionable
studies in eugenics, Fisher is considered one of the founders of population genetics and of modern statistical science.
His name is also linked to the so-called
Fisherian runaway sexual selection mechanism, aiming at an explanation of
some exaggerated, costly and apparently maladaptive male ornamentation in nature (seemingly in conflict with natural selection)
as dictated by persistent female choice (possibly triggered initially by the
ornament signaling greater potential fitness, hence likelihood of leaving more descendants):
the classical example of Fisherian runaway is the elaborate peacock plumage, see Figure~\ref{P1232DAILACRpeacWASTREDItangeFIlASE}.

A classical story about Fisher's contribution towards the standardization of statistical experiments is the ``lady tasting tea''~\cite{zbMATH02532791}. In a nutshell, once Fisher offered a cup of tea to phycologist Blanche Muriel Bristol.
She politely declined it, saying that she preferred the flavor when the milk was poured into the cup before (not after) the tea. So we have an intriguing dilemma: can the order of pouring milk affect the flavor of the tea? and, specifically, can Bristol notice the difference?

To test this, Fisher provided Bristol with eight randomly ordered cups of tea: four of which were prepared by first pouring the tea and then adding milk, and four by first pouring the milk and then adding the tea. Bristol had to identify the different cups. 
Fisher had to interpret the data to distinguish the case in which Bristol could actually tell the difference in the pouring order
and the one in which she could identify correctly some cups, but just by chance.

For this,
with a bit of combinatorics, Fisher computed that there are~${{8}\choose{4}}=\frac{8!}{4! (8-4)!}=70$ possible combinations of cups.
Since Bristol was aware that there were four cups of each type, her answer would have included four of each.

There is of course only one possibility of making precisely the right identification: hence achieving a complete identification
just by chance occurs with probability~$\frac{1}{70}=0.01428571428...$;
there are instead~16 possibilities of making exactly one error (indeed, one error comes from a single swap of two cups of different type, say the $i$th cup of the first type with the~$j$th cup of the second type, with~$i$, $j\in\{1,2,3,4\}$).
Hence, there are~$1+16=17$ possibilities of making at most one mistake
and therefore, if we adopt a conventional probability criterion of success to be below the~$5\%$ threshold, the ability of Bristol to properly categorize the cups of tea would have been confirmed if, and only if, she was able to correctly identify all of them, without making any mistake.

This story is relevant since it became an example of a randomized experiment for checking a ``null hypothesis'' in statistics and it contributed to a scientific establishment of randomization analysis of experimental data.

For the record, Muriel Bristol allegedly succeeded in classifying all eight cups correctly.

See~\cite{zbMATH03151158, MR1815390} for further reading on the lady tasting tea experiment.

See also Figure~\ref{2HAFODAKEDFUMSPierre-SimonLaplacldRGIRA4AXELEUMDJOMNFHARLROA7KOLM789GIJ7solFUMHDNOJHNFOJED231}
for a picture of Kolmogorov and footnote~\ref{KOLMONOTE}
on page~\pageref{KOLMONOTE} for more information about the prominent figure of Kolmogorov.}
in jargon the Fisher-Kolmogorov-Petrovsky-Piskunov equation
(or Fisher-KPP equation for short). \index{Fisher-Kolmogorov-Petrovsky-Piskunov equation}
Interestingly, this equation is a reaction-diffusion equation \index{reaction-diffusion equation}
in the setting introduced on page~\pageref{REDI6NFHARLROA7KOLM789GIJ7solFUMHDNOJHNFOJED231}.

\begin{figure}
                \centering
                \includegraphics[height=.21\textheight]{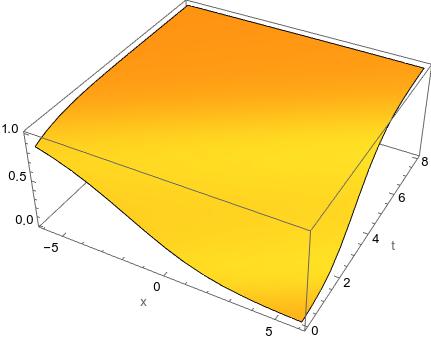}$\quad$\includegraphics[height=.21\textheight]{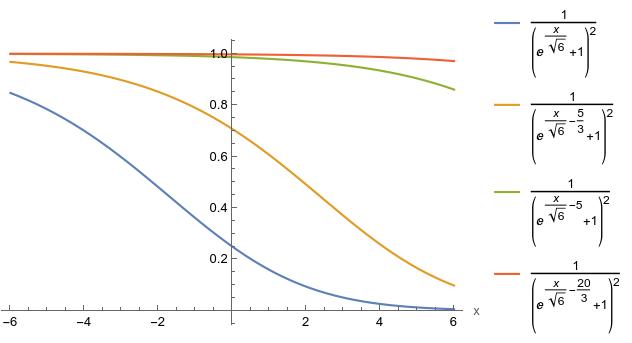}
        \caption{\sl A traveling wave of a species invading the spatial domain from left to right.}\label{2HAFODAKEDFUMSldRGIRA4AXTRAVEJHNFOJED231}
\end{figure}

The question of determining the invasivity of biological species is very related to~\eqref{HJA:A89MAOL8789892BIOEQLOG-3}.
In this setting, in light of~\eqref{HJA:A89MAOL8789892BIOEQLOG-4}, the level~$u=0$ corresponds to the absence of the population
and the level~$u=1$ to the maximal density permitted by environmental resources.
Remarkably, for any given~$v\in[2,+\infty)$
\begin{equation}\label{BEVEPAVKCiuneA3OMAkac}
\begin{split}&
{\mbox{equation~\eqref{HJA:A89MAOL8789892BIOEQLOG-3} admits a solution~$u$ describing the full environmental invasion}}\\ &
{\mbox{of the biological species with velocity~$v$.}}\end{split}\end{equation}
More precisely, given any~$v\in[2,+\infty)$ and any~$\omega\in \partial B_1$, there exists a smooth function~$U:\R\to(0,1)$, with
\begin{equation}\label{BEVEPAVKCiuneA3OMAkac22} \lim_{\tau\to+\infty}U(\tau)=0\qquad{\mbox{and}}\qquad\lim_{\tau\to-\infty}U(\tau)=1,\end{equation}
and a solution~$u$ of~\eqref{HJA:A89MAOL8789892BIOEQLOG-3} of the form
\begin{equation}\label{SOLIWA} u(x,t):=U(\omega\cdot x-vt).\end{equation}
The solutions\footnote{The concrete interest of the very special solutions in~\eqref{SOLIWA} is given by the fact that
under suitable initial conditions, ``reasonable'' solutions of~\eqref{HJA:A89MAOL8789892BIOEQLOG-3} evolve into a solitary wave,
see~\cite{MR1175399, MR0427837, MR511740}. Roughly speaking, wells in the initial conditions
are quickly reabsorbed by the evolving dynamics, since regions with~$u$ close to~$1$ tend to grow
and eat up adjacent regions with~$u$ close to~$0$, producing, for large times, a monotonic traveling wave with constant speed.
This phenomenon makes solitary waves an essential building block for understanding the dynamics of the
Fisher-Kolmogorov-Petrovsky-Piskunov equation.} of this form are often called ``solitary waves''. \index{solitary wave}
See Figure~\ref{2HAFODAKEDFUMSldRGIRA4AXTRAVEJHNFOJED231} for a representation of a solitary wave.

To establish~\eqref{BEVEPAVKCiuneA3OMAkac},
we consider\footnote{As remarked in~\cite{MR639998}, when~$v:=\frac{5}{\sqrt {6}}=2.04124145232...$, 
one can find explicit solutions of~\eqref{BEVEPAVKCiuneA3OMAkac2} in the form
$$ U(\tau)=\frac1{\displaystyle\left(1+C\exp\left(\frac\tau{\sqrt{6}}\right)\right)^{2}},$$
for all~$C>0$.

Figure~\ref{2HAFODAKEDFUMSldRGIRA4AXTRAVEJHNFOJED231}
was obtained with this function with~$C:=1$ and the frames of the second image there are obtained
by taking times~$t\in\{0,2,4,6,8\}$.

See also~\cite[Section~13.4]{MR1908418} for more information about this explicit solution.
We are not aware of other solutions of~\eqref{BEVEPAVKCiuneA3OMAkac2} which
can be written explicitly in a nice and closed form.}
the ordinary differential equation
\begin{equation}\label{BEVEPAVKCiuneA3OMAkac2}
\frac{d^2U}{d\tau^2}(\tau)+v\frac{dU}{d\tau}(\tau)+U(\tau)\big(1-U(\tau)\big)=0.\end{equation}
By setting~$V(\tau):=\frac{dU}{d\tau}(\tau)$, we obtain the system of equations
\begin{equation}\label{BEVEPAVKCiuneA3OMAkac3b}
\left(\begin{matrix}
\displaystyle\frac{dU}{d\tau} \\
\\
\displaystyle\frac{dV}{d\tau}
\end{matrix}\right)=\left(\begin{matrix}
V \\
\\
-vV-U(1-U)
\end{matrix}\right).\end{equation}
The equilibria of~\eqref{BEVEPAVKCiuneA3OMAkac3b} are~$E_0:=(0,0)$ and~$E_1:=(1,0)$.
The linearized dynamics in the vicinity of~$E_1$ is described by the matrix~$M_1:=\left(
\begin{matrix}
0 & 1\\
1 & -v
\end{matrix}\right)$, which possesses eigenvalues~$\frac12\left(-v \pm \sqrt{v^2 + 4}\right)$
with eigenvectors~$\left(\frac12 \left(v\pm\sqrt{v^2+4}\right), 1\right)$.

\begin{figure}
  \centering
  \includegraphics[height=4.3cm, width=6cm]{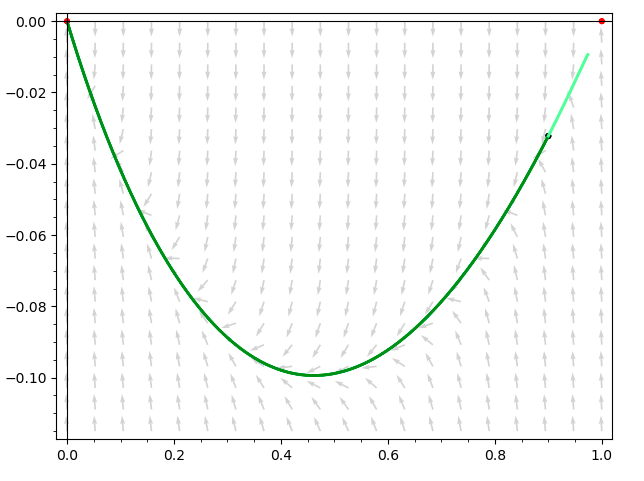}$\quad$\includegraphics[height=4.3cm, width=6cm]{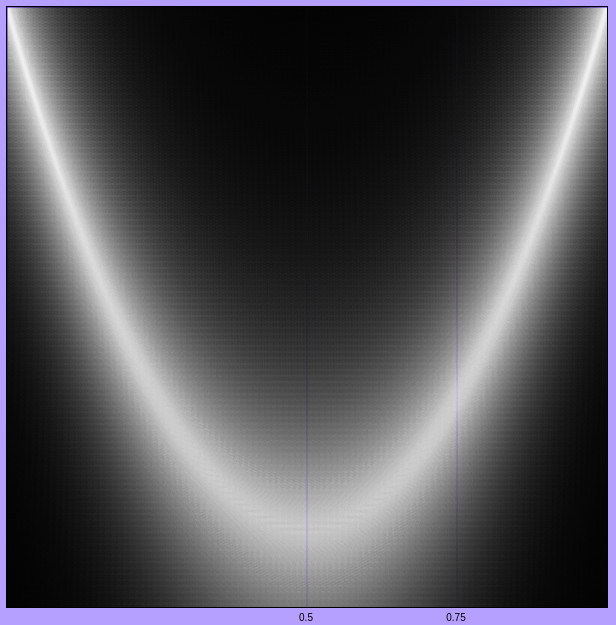}
 \caption{\sl Phase portrait and heteroclinc connection for the dynamical system in~\eqref{BEVEPAVKCiuneA3OMAkac3b} with~$v:=5/2$.}\label{12FINEST}
\end{figure}

In particular, $E_1$ is a hyperbolic saddle for all values of~$v\in\R$. So we consider the unstable manifold
(see e.g.~\cite[Section~2.7]{MR1801796}) of~$E_1$ in the direction~$-\left(\frac12 \left(v+\sqrt{v^2+4}\right), 1\right)$.
This provides us with an orbit~$(U(\tau),V(\tau))$ such that
\begin{equation}\label{JMNSRAPPORA1}
\lim_{\tau\to-\infty}(U(\tau),V(\tau))=E_1\end{equation}
and
\begin{equation}\label{JMNSRAPPORA}
\lim_{\tau\to-\infty}\frac{1}{\lambda e^{\lambda t}}
\left(\frac{dU}{d\tau} (\tau),\frac{dV}{d\tau}(\tau)\right)= -
\left(\frac12 \left(v+\sqrt{v^2+4}\right), 1\right),\end{equation}with~$\lambda:=\frac{\sqrt{v^2+4}-v}{2}$,
see Figure~\ref{12FINEST}.

Given~$v\in[2,+\infty)$, we now consider the planar region
$$ {\mathcal{R}}:=\left\{{\mbox{$(U,V)\in\R^2$ \;s.t.\; $U<1$, \;$V\in\left(-\displaystyle\frac1v,0\right)$\;
and \; $\displaystyle\frac{v+\sqrt{v^2-4}}{2}\,U+V>0$}}\right\}.$$

Notice that, in light of~\eqref{JMNSRAPPORA1} and~\eqref{JMNSRAPPORA},
there exists~$\tau_\star\in\R$ such that~$(U(\tau),V(\tau))\in{\mathcal{R}}$ for all~$\tau\in(-\infty,\tau_\star)$.

We claim that
\begin{equation}\label{CONGFIuynU6mnemt}
{\mbox{$(U(\tau),V(\tau))\in{\mathcal{R}}$\; for all \;$\tau\in(-\infty,+\infty)$.}}\end{equation}
Indeed, suppose not. Then, there exists~$\tau_0\ge \tau_\star$ such that~$(U(\tau),V(\tau))\in{\mathcal{R}}$ for all~$\tau\in(-\infty,\tau_0)$,
but either~$V(\tau_0)=-1/v$, or~$V(\tau_0)=0$
or~$U(\tau_0)=1$, or~$\frac{v+\sqrt{v^2-4}}{2}\,U(\tau_0)+V(\tau_0)=0$. We exclude all these cases to prove~\eqref{CONGFIuynU6mnemt}.

First, we show that~$V(\tau_0)=-1/v$ cannot hold. This is because otherwise~$
\frac{dV}{d\tau}(\tau_0)\le0$, but we know from~\eqref{BEVEPAVKCiuneA3OMAkac3b}
that
$$ \frac{dV}{d\tau}(\tau_0)=-vV(\tau_0)-U(\tau_0)\big(1-U(\tau_0)\big)\ge-vV(\tau_0)-\frac14=1-\frac14>0.$$

Similarly, $V(\tau_0)=0$ cannot happen, otherwise
$$ 0\le\frac{dV}{d\tau}(\tau_0)=-vV(\tau_0)-U(\tau_0)\big(1-U(\tau_0)\big)\le-vV(\tau_0)=0.$$
This entails that either~$U(\tau_0)=0$ or~$U(\tau_0)=1$, whence either~$(U(\tau_0),V(\tau_0))=E_0$
or~$(U(\tau_0),V(\tau_0))=E_1$, which is impossible, being both~$E_0$ and~$E_1$ equilibria for the system.

Also~$U(\tau_0)=1$ cannot hold true, otherwise, by~\eqref{BEVEPAVKCiuneA3OMAkac3b},
$$ 0\le\frac{dU}{d\tau}(\tau_0)=V(\tau_0)\le0,$$
that gives that~$(U(\tau_0),V(\tau_0))=E_1$, which is impossible.

Finally, the case~$\frac{v+\sqrt{v^2-4}}{2}\,U(\tau_0)+V(\tau_0)=0$ cannot occur as well, otherwise we argue as follows.
Let
$$ F(\tau):=\frac{v+\sqrt{v^2-4}}{2}\,U(\tau)+V(\tau).$$
We note that~$F(\tau)>0$ for~$\tau\in(-\infty,\tau_0)$ and~$F(\tau_0)=0$. As a result, from~\eqref{BEVEPAVKCiuneA3OMAkac3b} we arrive at
\begin{eqnarray*} 0&\ge&\frac{dF}{d\tau}(\tau_0)\\&=&
\frac{v+\sqrt{v^2-4}}{2}\,\frac{dU}{d\tau}(\tau_0)+\frac{dV}{d\tau}(\tau_0)\\&=&
\frac{v+\sqrt{v^2-4}}{2}\,V(\tau_0)-vV(\tau_0)-U(\tau_0)+U^2(\tau_0)\\&=&
\frac{\sqrt{v^2-4}-v}{2}\,V(\tau_0)-U(\tau_0)+U^2(\tau_0)\\&=&-
\frac{\big(\sqrt{v^2-4}-v\big)\big(v+\sqrt{v^2-4}\big)}{4}\,U(\tau_0)-U(\tau_0)+U^2(\tau_0)\\&=&-
\frac{(v^2-4)-v^2}{4}\,U(\tau_0)-U(\tau_0)+U^2(\tau_0)\\&=&
U^2(\tau_0).
\end{eqnarray*}
This gives that~$U(\tau_0)=0$, whence~$V(\tau_0)=0$, but this contradicts the fact that~$E_0$ is an equilibrium for the system.

The proof of~\eqref{CONGFIuynU6mnemt} is thereby complete.

As a byproduct of~\eqref{CONGFIuynU6mnemt}, we also have that, for all~$t\in\R$,
\begin{equation}\label{CONGFIuynU6mnemtLI9}
1>U(\tau)>-\frac{2}{v+\sqrt{v^2-4}}\,V>0.
\end{equation}

Now we claim that
\begin{equation}\label{CONGFIuynU6mnemtLI}
\lim_{\tau\to+\infty}(U(\tau),V(\tau))=E_0.
\end{equation}
Indeed, by~\eqref{CONGFIuynU6mnemt} we know that~$V(\tau)<0$ for all~$\tau\in\R$, therefore by~\eqref{BEVEPAVKCiuneA3OMAkac3b}, $U(\tau)$ is decreasing.
Since also by~\eqref{CONGFIuynU6mnemt} we have that~$U(\tau)\in(0,1)$ for all~$\tau\in\R$, we infer that
there exists~$U_\infty\in[0,1)$ such that~$U(\tau)\to U_\infty$ as~$\tau\to+\infty$.

Now, let~$\delta>0$, to be taken as small as we wish here below. Using the ``dots'' to denote derivatives, recalling~\eqref{BEVEPAVKCiuneA3OMAkac3b} and~\eqref{CONGFIuynU6mnemt} we observe that
\begin{equation}\label{KAQUINLtyY93l2}\begin{split}
|V(\tau)|\,&=\frac1\delta\left|\int_\tau^{\tau+\delta}\big(V(\tau)-V(\theta)+V(\theta)\big)\,d\theta\right|\\&
\le \frac1\delta\left|\int_\tau^{\tau+\delta}\big(V(\tau)-V(\theta)\big)\,d\theta\right|+
\frac1\delta\left|\int_\tau^{\tau+\delta} V(\theta)\,d\theta\right|\\&=
\frac1\delta\left|\int_\tau^{\tau+\delta}\left(\int_\theta^\tau\dot{V}(\zeta)\,d\zeta\right)\,d\theta\right|+
\frac1\delta\left|\int_\tau^{\tau+\delta} \dot{U}(\theta)\,d\theta\right|\\&=
\frac1\delta\left|\int_\tau^{\tau+\delta}\left(\int^\theta_\tau vV(\zeta)+U(\zeta)\big(1-U(\zeta)\big)\,d\zeta\right)\,d\theta\right|+
\frac1\delta\left|U(\tau+\delta)-U(\tau)\right|\\&\le
\frac2\delta\int_\tau^{\tau+\delta}  \left(\int^\theta_\tau \,d\zeta\right)\,d\theta+
\frac1\delta\left|U(\tau+\delta)-U(\tau)\right|\\&=
\delta+
\frac1\delta\left|U(\tau+\delta)-U(\tau)\right|.\end{split}
\end{equation}
{F}rom this we arrive at
$$ \limsup_{\tau\to+\infty}|V(\tau)|\le\limsup_{\tau\to+\infty}
\delta+
\frac1\delta\left|U(\tau+\delta)-U(\tau)\right|=\delta+
\frac1\delta\left|U_\infty-U_\infty\right|=\delta$$
and accordingly, by choosing~$\delta$ arbitrarily small,
$$ \lim_{\tau+\infty} V(\tau)=0.$$
As a result, we have that~$(U(\tau),V(\tau))\to(U_\infty,0)\in[0,1)\times\{0\}$.
Since the only equilibrium of~\eqref{BEVEPAVKCiuneA3OMAkac3b} on~$[0,1)\times\{0\}$
is~$E_0$, the proof of~\eqref{CONGFIuynU6mnemtLI} is complete.

Hence, in view of~\eqref{BEVEPAVKCiuneA3OMAkac3b}, \eqref{JMNSRAPPORA1}, \eqref{CONGFIuynU6mnemtLI9} and~\eqref{CONGFIuynU6mnemtLI},
the trajectory~$U(\tau)$ provides a solution of~\eqref{BEVEPAVKCiuneA3OMAkac2} satisfying the asymptotics in~\eqref{BEVEPAVKCiuneA3OMAkac22}. Then, by~\eqref{SOLIWA},
\begin{eqnarray*}&&
-\partial_t u(x,t)+\Delta u(x,t)+u(x,t)\big(1-u(x,t)\big)\\&&\qquad=
v \frac{dU}{d\tau}(\omega\cdot x-vt)
+\frac{d^2U}{d\tau^2}U(\omega\cdot x-vt)
+U(\omega\cdot x-vt)\big(1-U(\omega\cdot x-vt)\big)=0.
\end{eqnarray*}
This shows that~$u$ is the desired solitary wave solution of the Fisher-Kolmogorov-Petrovsky-Piskunov equation in~\eqref{HJA:A89MAOL8789892BIOEQLOG-3}.\medskip

For completeness, we show also that biologically interesting solitary waves of the Fisher-Kolmogorov-Petrovsky-Piskunov equation
do not exist when~$v\in[0,2)$, since 
\begin{equation}\label{1343g43678EPAVKCiuneA3OMAkI536}
\begin{split}&
{\mbox{when~$v\in[0,2)$ any solution~$U:\R\to(-\infty,1]$ of~\eqref{BEVEPAVKCiuneA3OMAkac22} and~\eqref{BEVEPAVKCiuneA3OMAkac2}}}\\&{\mbox{necessarily attains also negative values.}}\end{split}\end{equation} We recall indeed that biologically interesting values of the solutions are the ones modeling a population density,
hence they are confined in the interval~$[0,1]$, due to~\eqref{HJA:A89MAOL8789892BIOEQLOG-4}.

To check the claim in~\eqref{1343g43678EPAVKCiuneA3OMAkI536}, 
suppose by contradiction that a solution~$U:\R\to(0,1]$ of~\eqref{BEVEPAVKCiuneA3OMAkac22} and~\eqref{BEVEPAVKCiuneA3OMAkac2}
exists for some~$v\in[0,2)$ and let~$V:=\frac{dU}{d\tau}$ to find the system of equations in~\eqref{BEVEPAVKCiuneA3OMAkac3b}.

We claim that, for all~$\tau\ge0$,
\begin{equation}\label{KAQUINLtyY93l}
|V(\tau)|< |V(0)|+\frac1v.
\end{equation}
Indeed, suppose not. Then, there exists~$\underline\tau>0$ such that~$|V(\tau)|<|V(0)|+\frac1v$
for all~$\tau\in[0,\underline\tau)$ and~$|V(\underline\tau)|=|V(0)|+\frac1v$. Two cases can occur,
either~$V(\underline\tau)\ge0$ or~$V(\underline\tau)<0$. In the first case, we have that
$$ 0\le\frac{dV}{d\tau}(\underline\tau)=-vV(\underline\tau)-U(\underline\tau)\big(1-U(\underline\tau)\big)
=-v\left(|V(0)|+\frac1v\right)-U(\underline\tau)\big(1-U(\underline\tau)\big)
\le-1,$$
which is a contradiction.

If instead~$V(\underline\tau)<0$, then
$$ 0\ge\frac{dV}{d\tau}(\underline\tau)=-vV(\underline\tau)-U(\underline\tau)\big(1-U(\underline\tau)\big)=
v\left(|V(0)|+\frac1v\right)-U(\underline\tau)\big(1-U(\underline\tau)\big)\ge1-\frac14>0,
$$
which is a contradiction too, hence~\eqref{KAQUINLtyY93l} is proved.

Using this and computing as in~\eqref{KAQUINLtyY93l2}, for all~$\tau\ge0$ and~$\delta>0$ we arrive at
\begin{eqnarray*}
|V(\tau)|&\le&
\frac1\delta\left|\int_\tau^{\tau+\delta}\left(\int^\theta_\tau vV(\zeta)+U(\zeta)\big(1-U(\zeta)\big)\,d\zeta\right)\,d\theta\right|+
\frac1\delta\left|U(\tau+\delta)-U(\tau)\right|\\&\le&
C\delta +
\frac1\delta\left|U(\tau+\delta)-U(\tau)\right|,
\end{eqnarray*}
with~$C:=v\left(|V(0)|+\frac1v\right)+1$.
Therefore,
$$\limsup_{\tau\to+\infty}
|V(\tau)|\le C\delta +
\frac1\delta\left|0-0\right|=C\delta$$
and accordingly, taking~$\delta$ as small as we wish, we find that~$V(\tau)\to0$ as~$\tau\to+\infty$.

In this way, we have established that~$(U(\tau),V(\tau))\to(0,0)$ as~$\tau\to+\infty$.
Hence, we pick~$\eta>0$, to be chosen conveniently small in what follows and we find~$\tau_\eta\in\R$ such that~$|(U(\tau),V(\tau))|\le\eta$ for all~$\tau\ge\tau_\eta$. Up to a time shift, i.e. up to replacing~$(U(\tau),V(\tau))$ with~$(U(\tau-\tau_\eta),V(\tau-\tau_\eta))$, we can assume from now on that~$\tau_\eta=0$ and write that
\begin{equation}\label{L-X2343545y678we34eS}
\big|(U(\tau),V(\tau))\big|\le\eta\qquad{\mbox{for all }}\,\tau\ge0.
\end{equation}

Now we perform an argument related to the nonlinear stability of equilibria in planar dynamical systems.
Namely, we take~$\varrho(\tau)>0$ and~$\varphi(\tau)\in\R$ such that
$$ U(\tau)=\varrho(\tau)\,\cos(\varphi(\tau))\qquad {\mbox{and}}\qquad
V(\tau)=\varrho(\tau)\,\sin(\varphi(\tau)).$$
Actually, since~$U(\tau)>0$, we have that
\begin{equation}\label{922i3k2e944racrhfo}
{\mbox{we can pick~$\varphi(\tau)$ in the interval }}\,\left(-\frac\pi2,\frac\pi2\right).\end{equation}
Thus, denoting by the dot the derivative with respect to~$\tau$,
$$ \dot U=\dot\varrho\cos\varphi- \varrho\dot\varphi\sin\varphi
\qquad {\mbox{and}}\qquad
\dot V=\dot\varrho\sin\varphi+\varrho\dot\varphi\cos\varphi.$$
{F}rom this we arrive at
\begin{eqnarray*}
\dot U\sin\varphi-\dot V\cos\varphi&=&
\big(\dot\varrho\cos\varphi- \varrho\dot\varphi\sin\varphi\big)\sin\varphi-
\big(\dot\varrho\sin\varphi+\varrho\dot\varphi\cos\varphi\big)\cos\varphi
\\&=&
\dot\varrho\sin\varphi\cos\varphi- \varrho\dot\varphi\sin^2\varphi-
\dot\varrho\sin\varphi\cos\varphi-\varrho\dot\varphi\cos^2\varphi
\\&=&-\varrho\dot\varphi.
\end{eqnarray*}
Since, by~\eqref{BEVEPAVKCiuneA3OMAkac3b},
\begin{eqnarray*}
\dot U\sin\varphi-\dot V\cos\varphi&=&
V\sin\varphi+\big(vV+U-U^2\big)\cos\varphi\\
&=&\varrho\sin^2\varphi+v\varrho\sin\varphi\cos\varphi+
\varrho\cos^2\varphi-\varrho^2\cos^3\varphi\\&=&\varrho+\frac{v\varrho}2\sin(2\varphi)
-\varrho^2\cos^3\varphi,
\end{eqnarray*}
we conclude that
\[ -\varrho\dot\varphi=\varrho+\frac{v\varrho}2\sin(2\varphi)
-\varrho^2\cos^3\varphi\]
and accordingly
\[ -\dot\varphi=1+\frac{v}2\sin(2\varphi)
-\varrho\cos^3\varphi.\]
In particular, using~\eqref{L-X2343545y678we34eS} and the assumption that here~$v\in[0,2)$,
\[ -\dot\varphi\geq 1-\frac{v}2-\varrho\ge 1-\frac{v}2-\eta
\ge\frac{1-\displaystyle\frac{v}2}{2}>0,\]
as long as~$\eta$ is small enough.
This gives that there exists~$\tau_\star>0$ such that~$\varphi(\tau_\star)<-\frac{\pi}{2}$,
in contradiction with~\eqref{922i3k2e944racrhfo}.
This completes\footnote{Here is an alternative argument
to prove~\eqref{1343g43678EPAVKCiuneA3OMAkI536}, relying on the additional simplifying assumptions that
$$ \lim_{\tau\to+\infty}\dot{U}(\tau)=0$$
and the following limit exists
$$ \ell:=\lim_{\tau\to+\infty} \frac{\ddot{U}(\tau)}{\dot{U}(\tau)}.$$
{F}rom this, if~$U(\tau)\in(0,1)$ for all~$\tau\in\R$ and~$U(\tau)\to0$ as~$\tau\to+\infty$, we can apply L'H\^opital's Rule and deduce that
$$ \lim_{\tau\to+\infty} \frac{\dot{U}(\tau)}{U(\tau)}=
\lim_{\tau\to+\infty} \frac{\ddot{U}(\tau)}{\dot{U}(\tau)}=\ell.$$
As a result, dividing~\eqref{BEVEPAVKCiuneA3OMAkac2} by~$U(\tau)$ and taking the limit we find that
\begin{eqnarray*}
0&=&\lim_{\tau\to+\infty}\frac{\ddot{U}(\tau)}{U(\tau)}+v\frac{\dot{U}(\tau)}{U(\tau)}+1-U(\tau)\\
&=&\lim_{\tau\to+\infty}\frac{\ddot{U}(\tau)}{\dot U(\tau)}\;
\frac{\dot{U}(\tau)}{U(\tau)}+v\frac{\dot{U}(\tau)}{U(\tau)}+1\\
&=&\ell^2+v\ell+1\\
&=&\left( \ell+\frac{v}2\right)^2+1-\frac{v^2}4
\\&\ge&1-\frac{v^2}4,
\end{eqnarray*}
which gives that~$v^2\ge4$, hence (in the convention that~$v\ge0$) necessarily~$v\ge2$ to allow invading fronts.

These computations are also interesting since for more general equations such as~$\ddot{U}+v\dot{U}+f(U)=0$ with~$f$ smooth, $f(r)>0$ for all~$r\in(0,1)$, $f(0)=f(1)=0$ and~$f'(1)<0$, a calculation as above leads to
\begin{eqnarray*}&& 0=\ell^2+v\ell+\lim_{\tau\to+\infty}\frac{f(U(\tau))}{U(\tau)}=
\ell^2+v\ell+f'(0)=\left( \ell+\frac{v}2\right)^2+f'(0)-\frac{v^2}4\ge f'(0)-\frac{v^2}4,
\end{eqnarray*}
leading to the necessary condition for invasion~$v\ge 2\sqrt{f'(0)}$.

This is quite interesting, since it also points out that the speed of the invasive front is dictated by its low density fringes not by the high density regions (namely, the velocity of invasion is completely determined by the forcing term~$f$
and specifically by the values of~$f$ in the vicinity of~$0$, not in the vicinity of~$1$).

A function~$f$ with the above structure is often called a \index{bistable nonlinearity} bistable nonlinearity.

For further details on bistable nonlinearities and plane wave solutions, see~\cite[Section~4]{MR511740}.}
the proof of~\eqref{1343g43678EPAVKCiuneA3OMAkI536}.\medskip

For additional information on the Fisher-Kolmogorov-Petrovsky-Piskunov, see~\cite[Chapter~13]{MR1908418},
\cite{MR2081104}, \cite{MR3408563}
and the references therein.

\subsection{Equations from hydrodynamics}\label{EUMSD-OS-32456i7DNSSE234R:SEC}

A classical arena for partial differential equations is provided by fluid mechanics and hydrodynamics.
Let us present some of the situations that naturally occur in this framework.
Let us consider a fluid with density~$\rho(x,t)$, vectorial velocity~$v(x,t)$ and pressure~$p(x,t)$.
The position of a parcel of fluid at time~$t$ is given by~$x(t)$ and therefore such a parcel
travels according to the velocity prescription
\begin{equation}\label{0uojfSANMFOSJFEL}
\dot x(t)=v(x(t),t).
\end{equation}

\begin{figure}
  \centering
  \includegraphics[width=.35\linewidth]{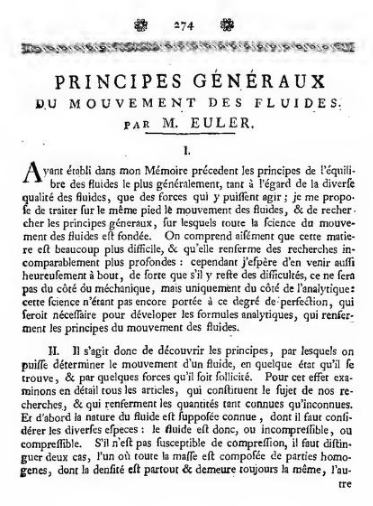}
 \caption{\sl The title page to Euler's original article ``Principes g\'en\'eraux du mouvement des fluides'',
 published in M\'emoires de l'Acad\'emie des Sciences de Berlin in 1757 (Public Domain source from
 Internet Archive).}\label{GR-EULEROFIDItangeFI}
\end{figure}

Also, the total mass of a fluid element occupying a region~$\Omega$ of the space at a given time~$t_0$
is given by the quantity~$\int_{\Omega} \rho(x,t_0)\,dx$. This amount of fluid will move around in a small time~$\tau$
possibly occupying a different region, that we name~$\Omega_\tau$ which collects all the evolution trajectories~$x(t_0+\tau)$
of the fluid parcels such that~$x(t_0)\in\Omega$. Since we are assuming that matter is neither created nor disappear in this process
we deduce a conservation law given by the equation
$$ \int_{\Omega} \rho(x,t_0)\,dx=\int_{\Omega_\tau} \rho(x,t_0+\tau)\,dx$$
or equivalently
\begin{equation}\label{0ojCOMAS} \frac{d}{d\tau}\int_{\Omega_\tau} \rho(x,t_0+\tau)\,dx=0.\end{equation}
Notice that we can change variable in the previous integral simply by flowing back the fluid.
Specifically, given a point~$q\in\R^n$ we denote by~$\Phi^t(q)$ the evolution according to the fluid vector field~$v$
with initial position at time~$t_0$ equal to~$q$, or more formally we define~$\Phi^t(q)$ as the solution of the following
Cauchy problem for ordinary differential equations:
$$ \begin{dcases}
\frac{d}{dt} \Phi^t(q)=v\big(\Phi^t(q),t),\\
\Phi^{t_0}(q)=q.
\end{dcases}$$
This is nothing else as a rewriting of~\eqref{0uojfSANMFOSJFEL}, but we are keeping track here explicitly of the initial position.
With this notation, we have that~$\Omega_\tau=\{\Phi^\tau(y)$, $y\in\Omega\}$.
As a result, one can perform the change of variable
$$x=\Phi^\tau(y)=y+\tau v(y,t_0)+o(\tau),$$
which produces
$$ |\det \partial_y x|=1+\tau \div v(y,t_0)+o(\tau)$$
and therefore
\begin{eqnarray*}&&
\int_{\Omega_\tau} \rho(x,t_0+\tau)\,dx\\&=&
\int_{\Omega} \rho(y+\tau v(y,t_0),t_0+\tau)\,\big(1+\tau \div v(y,t_0)\big)\,dy+o(\tau)
\\&=&\int_{\Omega} 
\Big(\rho(y,t_0)+\tau\nabla\rho(y,t_0)\cdot v(y,t_0)
+\tau\partial_t \rho(y,t_0)\Big)
\,\big(1+\tau \div v(y,t_0)\big)\,dy+o(\tau)\\&=&\int_{\Omega} 
\Big(\rho(y,t_0)+\tau\nabla\rho(y,t_0)\cdot v(y,t_0)
+\tau\partial_t \rho(y,t_0)+
\tau\rho(y,t_0)\div v(y,t_0)\Big)\,dy+o(\tau)\\&=&\int_{\Omega} 
\Big(\rho(y,t_0)
+\tau\partial_t \rho(y,t_0)+
\tau\div\big(\rho(y,t_0)\, v(y,t_0)
\big)\Big)\,dy+o(\tau).
\end{eqnarray*}
This and~\eqref{0ojCOMAS} formally lead to
\begin{equation}\label{D23q4wy5OPACLPASIJMDNIAKS279uyt3M}
\begin{split}0=\left.
\frac{d}{d\tau}\int_{\Omega_\tau} \rho(x,t_0+\tau)\,dx\right|_{\tau=0}=
\int_{\Omega} 
\Big(\partial_t \rho(y,t_0)+\div\big(\rho(y,t_0)\, v(y,t_0)\big)\Big)\,dy.\end{split}\end{equation}
Since this holds true for all domains~$\Omega$ and all times~$t_0$, we have obtained the equation
\begin{equation}\label{OJS-PJDN-0IHGDOIUGDBV02ujrfMTE}
\partial_t \rho+\div(\rho v)=0.
\end{equation}
This is called\footnote{It is interesting to observe that in the framework of equation~\eqref{OJS-PJDN-0IHGDOIUGDBV02ujrfMTE}
one can reinterpret the divergence term in~\eqref{OJS-PJDN-0IHGDOIUGDBV02ujrf}
as a transport term with velocity~$\nabla w$ (that is, the effect of chemotaxis
is to produce a transport of the biological population with a velocity proportional to the gradient of the chemical attractant).

Moreover, in chemistry one often models the spreading of some substrate with concentration~$\rho$
in a fluid with velocity~$v$: in this context equation~\eqref{OJS-PJDN-0IHGDOIUGDBV02ujrfMTE}
is often called ``convection equation''\index{convection equation}.}
in jargon ``mass transport equation''\index{mass transport equation}
(or simply ``transport equation'') \index{transport equation}
 or ``continuity equation''\index{continuity equation}.

Of course, transport of mass, corresponding to the conservation of matter, is not the only ingredient to accurately
describe the dynamics of a fluid, hence one has to complement~\eqref{OJS-PJDN-0IHGDOIUGDBV02ujrfMTE} with additional information. In particular, one can take into account Newton's Law regarding momentum conservation.
In this setting, the change of fluid momentum must correspond to the forces acting on the fluid
(for simplicity, we assume here that the fluid is only subject to gravity and to its own pressure force
with no other external forces to be exerted).
Notice that the momentum corresponding to the mass of the fluid in a region of space~$\Omega$
is given by the quantity~$\int_\Omega \rho(x,t) v(x,t)\,dx$.
The corresponding gravity force is provided by~$-ge_n\int_\Omega \rho(x,t)\,dx$, being~$g$ the gravity acceleration constant
(and~$e_n=(0,\dots,0,1)$, supposing the gravity acting downwards in the vertical direction).
As for the pressure force, it arises from~$p$ in the normal direction
along the surface~$\partial\Omega$, with a conventional
minus sign (so that a resistance to velocity increasing occurs when moving towards regions with high pressures, while
low pressure regions produce a sucking effect). These considerations translate 
Newton's Law about momentum conservation into the formula
\begin{equation} \label{90uyoihfs83wEUSlNOE}
\frac{\partial}{\partial t}\int_\Omega \rho(x,t) v(x,t)\,dx=
-ge_n\int_\Omega \rho(x,t)\,dx-\int_{\partial\Omega} p(x,t) \nu(x)\,d{\mathcal{H}}^{n-1}_x,\end{equation}
being, as usual,~$\nu(x)$ the unit external normal at~$x\in\partial\Omega$. 

Furthermore, in light of~\eqref{D23q4wy5OPACLPASIJMDNIAKS279uyt3M}
(applied here to~$\rho v_j$ instead of~$\rho$) we have that, for all~$j\in\{1,\dots,n\}$,
\begin{eqnarray*}
&&\frac{\partial}{\partial t}\int_\Omega \rho(x,t) v_j(x,t)\,dx\\
&=& \int_{\Omega} 
\Big(\partial_t (\rho v_j)(x,t)+\div\big((\rho v_j)(x,t)\, v(x,t)\big)\Big)\,dx\\&=&
\int_{\Omega} 
\Big(\rho (x,t)\partial_t v_j(x,t)+\rho(x,t) v(x,t)\cdot\nabla v_j(x,t)\Big)\,dx,
\end{eqnarray*}
where~\eqref{OJS-PJDN-0IHGDOIUGDBV02ujrfMTE} has been exploited in the latter step.

\begin{figure}
                \centering
                \includegraphics[width=.35\linewidth]{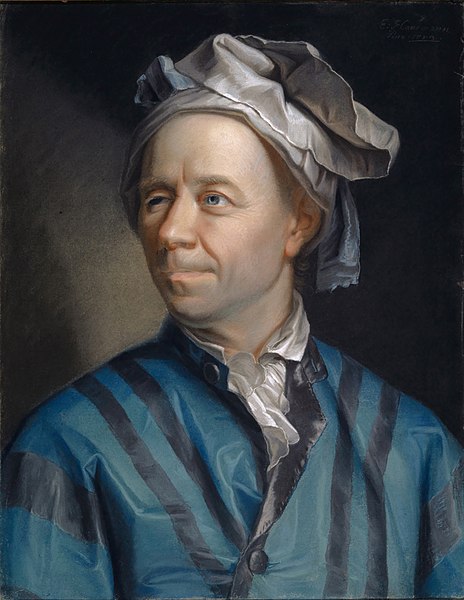}
        \caption{\sl Portrait of Leonhard Euler (Public Domain image from
        Wikipedia).}\label{2HAFODAKEDFUMSldRGIRA4AXELEUMDJOMNFHARLROA7789GIJ7solFUMHDNOJHNFOJED231}
\end{figure}

Thus, using the notation~$v\cdot\nabla$ to denote the operator~$\displaystyle\sum_{i=1} v_i\partial_i$,
possibly applied to a vector valued function componentwise, we can write that
\begin{eqnarray*}
&&\frac{\partial}{\partial t}\int_\Omega \rho(x,t) v(x,t)\,dx=
\int_{\Omega} 
\Big(\rho (x,t)\partial_t v(x,t)+\rho(x,t) \big(v(x,t)\cdot\nabla\big) v(x,t)\Big)\,dx.
\end{eqnarray*}
{F}rom this and~\eqref{90uyoihfs83wEUSlNOE}
we obtain that
\begin{equation} \label{90uyoihfs83wEUSlNOE-VXBNC}\begin{split}&
\int_{\Omega} 
\Big(\rho (x,t)\partial_t v(x,t)+\rho(x,t) \big(v(x,t)\cdot\nabla\big) v(x,t)\Big)\,dx\\&\qquad\qquad=
-ge_n\int_\Omega \rho(x,t)\,dx-\int_{\partial\Omega} p(x,t) \nu(x)\,d{\mathcal{H}}^{n-1}_x.\end{split}\end{equation}
Also, from the Divergence Theorem, for all~$j\in\{1,\dots,n\}$,
$$ \int_{ \Omega} \partial_j p(x,t)\,dx=
\int_{ \Omega} \div\big( p(x,t) e_j\big)\,dx
=\int_{\partial\Omega} p(x,t) e_j\cdot\nu(x)\,d{\mathcal{H}}^{n-1}_x$$
and in consequence
\[ \int_{ \Omega} \nabla p(x,t)\,dx=\int_{\partial\Omega} p(x,t) \nu(x)\,d{\mathcal{H}}^{n-1}_x.\]
{F}rom these considerations and~\eqref{90uyoihfs83wEUSlNOE-VXBNC}
we arrive at
$$\int_{\Omega} 
\Big(\rho (x,t)\partial_t v(x,t)+\rho(x,t) \big(v(x,t)\cdot\nabla\big) v(x,t)\Big)\,dx=
-ge_n\int_\Omega \rho(x,t)\,dx-\int_{ \Omega} \nabla p(x,t)\,dx$$
and accordingly
\begin{equation}\label{EUMSD-OS-D}
\rho \partial_t v+\rho (v\cdot\nabla) v=
-ge_n \rho-\nabla p.
\end{equation}
The system containing the mass conservation equation in~\eqref{OJS-PJDN-0IHGDOIUGDBV02ujrfMTE}
and the momentum balance equation in~\eqref{EUMSD-OS-D} (plus possibly a constitutive relation
linking pressure and density) constitutes what in jargon\footnote{The fountainhead of these equations
came indeed from the work of Leonhard Euler,
see Figure~\ref{GR-EULEROFIDItangeFI}.

See also Figure~\ref{2HAFODAKEDFUMSldRGIRA4AXELEUMDJOMNFHARLROA7789GIJ7solFUMHDNOJHNFOJED231}
for Euler's iconic portrait (by Jakob Emanuel Handmann).
Many of us wonder why Euler is depicted in a pijama with a towel on his head.
Arguably, it was not a pijama but a fashionable banyan, not a towel but a silk cloth
and possibly at that time it was believed that loose and informal
dresses (the banyan, no wig)
contribute to exercise of the faculties of the mind (and indeed we also find it much more comfortable
to do mathematics wearing very informal clothes).

Note that Euler is portrayed as facing left,
possibly to mitigate the impact on the spectator of a right eyelid ptosis and a divergent strabismus.
At age 31, Euler became almost blind in his right eye~\cite{158}
(and this might be a possible explanation for which Euler was nicknamed ``the cyclops''
by Frederick II, King of Prussia).

This partial loss of vision did not discourage Euler from producing the finest possible
mathematics (actually, he allegedly stated ``Now I will have fewer distractions'').
At age~59, a surgical restoration for a cataract in his left eye rendered Euler
almost totally blind. Yet, Euler's passion and talent for mathematics remained undefeated
and he kept producing a vast number of impressive works which changed the course of mathematics
till he died at age~76. In his eulogy, Marquis de Condorcet wrote
``il cessa de calculer et de vivre'' (French for ``he ceased to calculate and to live'').}
are called ``Euler fluid flow equations''\index{Euler equation}.\medskip

In a sense, these equations are obtained for a rather ideal situation and several modifications
of the previous setting can be performed to account for more complex models. Among these
modifications, one of the most popular one consists in accounting for some friction between the fluid molecules
which creates a ``viscosity''\index{viscosity} resisting to the change of fluid velocity. More specifically,
let us quantify\footnote{This idea of confronting the pointwise value of a function
with the average nearby will be extensively used in Chapter~\ref{CHA-2} and it will be
one of the leitmotifs of the study of harmonic functions.}
the difference of velocity between a fluid parcel at~$x$ and that of the parcels nearby
by the difference of~$v(x,t)$ and its average in a small ball, say of radius~$h$ centered at~$x$, namely
$$ {\mathcal{D}}(x,t):= v(x,t)-\fint_{B_h(x)}v(y,t)\,dy.$$
By exploiting~\eqref{KMS:RUIDKVODIG} and polar coordinates, we find that
\begin{equation*}
\begin{split}
{\mathcal{D}}(x,t)&= \fint_{B_h(x)}\big( v(x,t)-v(y,t)\big)\,dy\\
&=-\frac{1}{|B_h|} \int_0^h \left( \int_{\partial B_1}
r^{n-1} \big(v(x+r\omega,t)-v(x,t)\big)\,d{\mathcal{H}}^{n-1}_\omega\right)\,dr\\
&=-\frac{|B_1|}{|B_h|} \int_0^h
r^{n-1} \left(
\frac{c r^2}2 \Delta v(x,t) +o(h^2)
\right)\,dr\\
&=-\mu h^2  \Delta v(x,t) +o(h^2),
\end{split}\end{equation*}
for some~$\mu>0$.

Then, we suppose that the velocity of the fluid is reduced when~${\mathcal{D}}$ is positive (since the parcel at~$x$
is faster than the ones in its vicinity, and therefore it is dragged back by them)
and is enhanced when~${\mathcal{D}}$ is negative (since the parcel at~$x$
is slower than the ones in its vicinity, and therefore it is pulled forward by them).
For simplicity, we can therefore assume that such an additional acceleration term is proportional to~${\mathcal{D}}$
(say, to make it a finite quantity in the limit, proportional to~${\mathcal{D}}/h^2$). By taking~$h\searrow0$,
we thus obtain an additional acceleration (or deceleration) given by a term of the form~$\mu \Delta v$.
By incorporating this correction into~\eqref{EUMSD-OS-D} we thus find the equation
\begin{equation}\label{EUMSD-OS-DNS}
\rho \partial_t v+\rho (v\cdot\nabla) v=\mu\Delta v
-ge_n \rho-\nabla p,
\end{equation}
which is called in jargon\footnote{Equation~\eqref{EUMSD-OS-DNS} (possibly complemented with
other structural equations of the fluid)
come from Claude Louis Marie Henri Navier and Sir George Gabriel Stokes, 1st Baronet.
Our experience with the Laplacian as a ``democratic'' operator (recall the discussion on page~\pageref{DEEFFNVST}) \label{DEEFFNVSTBIS}
may suggest that the Laplacian in~\eqref{EUMSD-OS-DNS} is providing a smoothing effect.
This may be the case, but our present understanding of the solutions of~\eqref{EUMSD-OS-DNS}
is unsatisfactory and highly incomplete. In particular, it is not yet known if smooth (say,
globally defined and meeting some natural conditions)
solutions always exist. If the reader finds an answer to this question,
they do not only earn immortal fame among the circle of PDEs enthusiasts, but also obtain a
substantial amount of money, since the problem is listed in the Millennium Prize Problems \label{106106}
and a~$10^6$ US \$ award is offered to the first person providing a solution
(the case of  incompressible fluids in~$\R^3$ is hard enough). See e.g.~\cite{MR2238274, MR2246251} for further details about the most difficult way to become rich.}
the Navier-Stokes equation\index{Navier-Stokes equation}.\medskip

We stress that~\eqref{EUMSD-OS-DNS} is a vectorial equation (or, equivalently, a system of equations in each component of~$v$). There are also interesting situations in which~\eqref{EUMSD-OS-DNS} reduces to a scalar equation,
such as the one in which the fluid is constrained within an infinite pipe with a given direction. For instance,
to make things as simple as possible, 
let us suppose that the fluid is incompressible\index{incompressibility}: hence, we have that the density of the fluid remains constant in the flow
(namely, a change in the density over time would imply that the fluid had either compressed or expanded). That is, for incompressible fluids the quantity~$\rho(x(t),t)$ must be constant in time, leading to
\begin{equation*} 0=\frac{d}{dt}\rho(x(t),t)=\nabla\rho(x(t),t)\cdot \dot x(t)+\partial_t\rho(x(t),t)=
\nabla\rho(x(t),t)\cdot v(x(t),t)+\partial_t\rho(x(t),t).\end{equation*}
Combining this with~\eqref{OJS-PJDN-0IHGDOIUGDBV02ujrfMTE},
\begin{equation*} \rho\div v=\div(\rho v)-\nabla\rho\cdot v=\div(\rho v)+\partial_t\rho(x(t),t)=0,\end{equation*}
whence
\begin{equation}\label{INCOMPRE}
\div v=0,\end{equation}
that is the velocity is a divergence-free vector field.

\begin{figure}
  \centering
  \includegraphics[width=.7\linewidth]{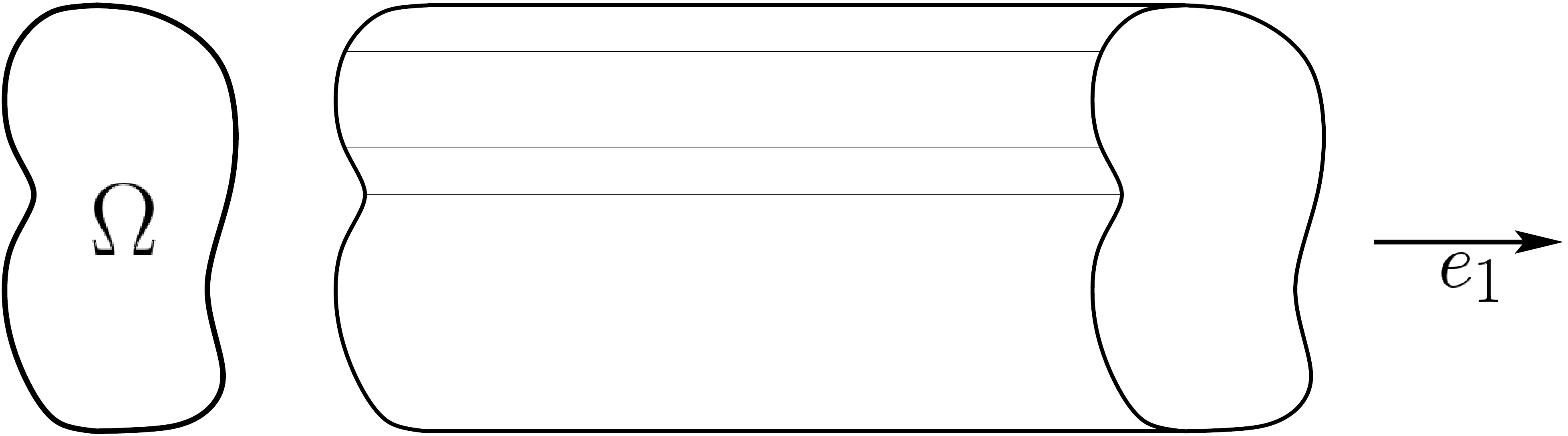}
 \caption{\sl A horizontal pipe in direction~$e_1$ with cross section~$\Omega$.}\label{PIPEDItangeFI}
\end{figure}

For this reason, recalling~\eqref{EUMSD-OS-DNS}, stationary solutions of incompressible viscous fluids
are solutions of
\begin{equation}\label{EUMSD-OS-DNSSER}
\begin{dcases}&
\mu\Delta v-\rho (v\cdot\nabla)v-ge_n \rho-\nabla p=0,\\&
\div v=0.
\end{dcases}
\end{equation}
Suppose now that the above fluid moves in parallel streamlines
through a horizontal straight pipe of a given cross sectional form.
Up to rotations which leave invariant the vertical directions, we can assume that the pipe has the form~$\R\times\Omega$
for some~$\Omega\subset\R^{n-1}$ and~$v=(u,0,\dots,0)$ for some scalar function~$u$ (that is, the pipe directing the fluid is oriented towards the~$e_1$-direction, see Figure~\ref{PIPEDItangeFI}). In this situation, 
we have that~$\Delta v=\Delta u \,e_1$,
$(v\cdot\nabla)v=u\partial_1u\,e_1$ and~$\div v=\partial_1u$.
Consistently with these observations,
\eqref{EUMSD-OS-DNSSER} reduces to
\begin{equation}\label{EUMSD-OS-DNSSE234R}
\begin{dcases}&
\mu\Delta u\,e_1-\rho u\partial_1u\,e_1-ge_n \rho-\nabla p=0,\\&
\partial_1u=0.
\end{dcases}
\end{equation}
Taking the components~$\{2,\dots,n-1\}$ of the first equation in~\eqref{EUMSD-OS-DNSSE234R} we see that~$\partial_i p=0$
for each~$i\in\{2,\dots,n-1\}$, giving that the fluid pressure is constant in these directions.
Taking the last component instead, we find~$\partial_n p=-g\rho$.

The most interesting information however comes from the first component of the first equation in~\eqref{EUMSD-OS-DNSSE234R} since
this provides an elliptic equation for~$u:\Omega\to\R$. Indeed, focusing on this aspect, we infer from~\eqref{EUMSD-OS-DNSSE234R} that~$u$ is a solution of
\begin{equation*}
\begin{dcases}&
\mu\Delta u-\rho u\partial_1u-\partial_1 p=0,\\&
\partial_1u=0.
\end{dcases}
\end{equation*}
Hence, plugging the second equation into the first one,
\begin{equation}\label{EUMSD-OS-32456i7DNSSE234R}
\mu\Delta u=\partial_1 p \qquad{\mbox{ in }}\,\Omega.
\end{equation}
A special case of interest, which constitutes one of the main motivations for the theory that will be presented
in Sections~\ref{SEC:OVER-SERRINTheorem} and~\ref{BACKSE}, is when:
\begin{equation}\label{NK090SOJDKHJNZ}\begin{split}
&{\mbox{the pressure rate is constant,}}\\
&{\mbox{the velocity on the boundary of the pipe is constant}}\\
&{\mbox{and the tangential stress on the pipe is constant.}}
\end{split}\end{equation}
The problem discussed in detail in Sections~\ref{SEC:OVER-SERRINTheorem} and~\ref{BACKSE}, and dating back
to James Serrin~\cite{MR333220}, is precisely
to identify the possible shape of the pipe that allows such special configurations (that is, to determine all the cross sections~$\Omega$ that give rise to solutions of~\eqref{EUMSD-OS-32456i7DNSSE234R}
which are compatible with the prescriptions in~\eqref{NK090SOJDKHJNZ}):
this is certainly a problem of great practical importance, since ideally one would like to project pipes which
do not produce excessive tangential stress, to reduce the wearing of the material and the consequential formation of defects, holes
and licking (for this, maintaining a constant stress along the pipe could be a nice way to balance
the wearing somewhat uniformly along the whole surface of the pipe).

To efficiently address a question of this type, we first need to translate the physical
requests described in~\eqref{NK090SOJDKHJNZ}
into a mathematical framework. In order to achieve this,
we point out that the prescription that the pressure rate is constant simply says that~$\partial_1 p=c_1$, for some~$c_1\in\R$.

Also, the assumption that the velocity on the boundary of the pipe is constant means that~$u=c_2$ along~$\partial\Omega$,
for some~$c_2\in\R$. As a matter of fact, since equation~\eqref{EUMSD-OS-32456i7DNSSE234R} remains invariant
if we replace~$u$ with~$u-c_2$, we can simply assume\footnote{Equivalently, one can argue that the invariance of the form of the description of physical problems
among mutually translating reference frames allows us to choose an inertial frame of reference
moving at constant speed~$c_2$. In this frame the previous prescription reduces to that
the velocity on the boundary of the pipe is equal to zero.}
that~$c_2=0$, hence~$u=0$ along~$\partial\Omega$.

As for the stress on the pipe, we can suppose that this is due to the fact that fluid parcels near the pipe may move faster or slower
that the ones directly adhering to the pipe (that is, in the normalization above, the fluid parcel in the vicinity of the pipe
could have strictly positive or strictly negative velocity). We can additionally suppose that the wearing of the pipe is directly
influenced by this change of velocity (say, some particles of the pipe get removed by the fluid parcels moving either forwards or backwards). In this spirit, one can model the stress along the pipe as proportional to the normal derivative of the fluid velocity.
With this, the prescription of having a constant stress along the pipe is translated into~$\partial_\nu u=c_3$ along~$\partial\Omega$,
for some~$c_3\in\R$.

In consideration of this, and normalizing constants for simplicity,
we write~\eqref{EUMSD-OS-32456i7DNSSE234R} and~\eqref{NK090SOJDKHJNZ} as
\begin{equation}\label{EUMSD-OS-32456i7DNSSE234R-BNM}
\begin{dcases}
\Delta u=1&{\mbox{ in }}\Omega,\\
u=0&{\mbox{ on }}\partial\Omega,\\
\partial_\nu u=c&{\mbox{ on }}\partial\Omega.
\end{dcases}\end{equation}
This equation will be indeed\footnote{We observe that~\eqref{EUMSD-OS-32456i7DNSSE234R-BNM}
does possess a solution when~$\Omega$ is a ball. Specifically, given~$x_0\in\R^n$ and~$r>0$, the function~$u(x):=\frac{|x-x_0|^2-r^2}{2n}$ is a solution of~\eqref{EUMSD-OS-32456i7DNSSE234R-BNM}
in the ball~$B_r(x_0)$
(and note that this solution corresponds to our physical intuition of the problem,
since it parallels to a fluid reaching its maximal speed at the center of a circular pipe, with the velocity of the fluid being minimal at the boundary of the pipe, due to the viscous effects). \label{OJSLLAFOr5O3OSMEr3r}

We will discuss in Sections~\ref{SEC:OVER-SERRINTheorem} and~\ref{BACKSE}
whether this is the only type of possible solution
or whether there are domains different from the ball which allow different kinds of solutions.}
the starting point of the topics presented in
Sections~\ref{SEC:OVER-SERRINTheorem} and~\ref{BACKSE}. A different model leading to the same mathematical
problem will be also presented in Section~\ref{RODSBACKSE}.\medskip

We end this section by recalling a simple, but helpful, use of the ``principle of inertia'' in the fluid dynamics setting. \label{PRINERY}
Namely, a case of special interest is provided by a rigid body (say, with a given shape described by a set~$\Omega$) moving in a fluid with constant speed. In this situation, the equations presented in this section hold outside the moving domain described, for instance, by~$\Omega(t):=\{x\in\R^3$ s.t. $x+v_0 t\in\Omega\}$, for a given constant vector~$v_0\in\R^3$.
Since fluid dynamics is already very difficult in a given domain, it is usually more agreeable to change the inertial frame to describe the motion by following the moving body (and we expect by Newton's First Law that this new choice of coordinates does not alter the physical description of the phenomena).

More explicitly, it comes in handy to define
\begin{equation}\label{ISUJONPRENsyhfngTNAbaujnfRGbhndfIKAR}
\widetilde v(x,t):=v(x-v_0 t,t)+v_0,\qquad
\widetilde p(x,t):=p(x-v_0 t,t)\qquad{\mbox{and}}\qquad\widetilde \rho(x,t):=\rho(x-v_0 t,t)\end{equation} and observe that
equations~\eqref{OJS-PJDN-0IHGDOIUGDBV02ujrfMTE}, \eqref{EUMSD-OS-D}, \eqref{EUMSD-OS-DNS}
and~\eqref{INCOMPRE} are all preserved by this transformation, since~$\div\widetilde v=\div v$,
\begin{eqnarray*}
\partial_t \widetilde\rho+\div(\widetilde\rho \widetilde v)&=&
\partial_t \rho-v_0\cdot\nabla\rho+\div\big(\rho( v+v_0)\big)\\&=&\partial_t \rho+\div(\rho v)
\end{eqnarray*} and
\begin{eqnarray*}&&
\widetilde\rho \partial_t \widetilde v+\widetilde\rho (\widetilde v\cdot\nabla) \widetilde v-\mu\Delta \widetilde v+ge_n \widetilde\rho+\nabla \widetilde p\\&=&
\rho \big(\partial_t v-v_0\cdot\nabla v\big)
+\rho \big((v+v_0)\cdot\nabla) v-\mu\Delta v+ge_n \rho+\nabla p\\&=&
\rho \partial_t v
+\rho v\cdot\nabla v-\mu\Delta v+ge_n \rho+\nabla p.
\end{eqnarray*}
Also, the new domain of reference is now simply the complement of~$\Omega$.

Notice that the transformation in~\eqref{ISUJONPRENsyhfngTNAbaujnfRGbhndfIKAR}
is simply replacing the situation of a rigid body in constant speed motion in a fluid at rest at infinity with the one of a still body in a fluid with constant speed at infinity. This is precisely the idea leading to the construction of wind tunnels: to replicate the aerodynamic interactions between a moving object and the surrounding air, it is common to construct large tubes with air blowing through them against a static model of the object.

In the forthcoming Sections~\ref{WRItvHS-SE1} and~\ref{AIRFVBD01ogt9u432gtb}, where we will describe rigid bodies moving at a constant speed, we will tacitly assume to have performed the transformation in~\eqref{ISUJONPRENsyhfngTNAbaujnfRGbhndfIKAR} to reduce ourselves to the case of still objects and domains that do not vary with time.

\subsection{Irrotational fluids}\label{BARPotA}

Among all possible fluids, a class deserving some special attention is provided by those which present ``no vortex''.
In light of Stokes' Theorem, this notion is made precise by the mathematical concept\footnote{While most of the arguments
in other sections are valid in~$\R^n$, when dealing with the curl we reduce our attention to the case of~$\R^3$.
The case of~$\R^2$ is also included, simply by identifying~$\R^2$ with~$\R^2\times\{0\}\subset\R^3$.
The curl of the velocity field is often named ``vorticity''. \index{vorticity}}
of curl, that is
we say that a fluid is irrotational if the curl of its velocity field vanishes.

A natural question in this setting is whether ``vortexes can be produced out of nothing'':
for instance if the initial conditions of a fluid present no vortexes, is it possible that vortexes arise at a later stage?

To increase our familiarity with multivariate calculus and partial differential equations,
we show that this is impossible, at least for inviscid, incompressible and barotropic flows \label{BAROFL}
(in this setting, a flow is said to be barotropic if its density is a function of the pressure only, say~$\rho=g(p)$
for some positive function~$g$). More precisely, we consider a barotropic solution in the absence of gravity of
the momentum balance equation in~\eqref{EUMSD-OS-D}
satisfying the incompressibility condition in~\eqref{INCOMPRE}, that is
\begin{equation}\label{CURLDESE1}
\begin{dcases}
\rho \partial_t v+\rho (v\cdot\nabla) v=-\nabla p,\\
\div v=0,\\
\rho=g(p)
\end{dcases}\end{equation}
and we show that
\begin{equation}\label{CURLDESE}
{\mbox{if~$\curl v=0$ at time~$t=0$, then~$\curl v=0$ at every time~$t$.}}\end{equation}
To achieve this goal, we first recall the vector calculus identity, valid for all smooth vector fields~$V=(V_1,V_2,V_3):\R^3\to\R^3$,
\begin{equation}\label{VCISTPOK}
(V\cdot\nabla )V=\frac{1}{2}\nabla |V|^2-
V\times (\curl V),\end{equation}
having denoted by~$\times$ the vector product operation.
To prove~\eqref{VCISTPOK}, up to exchanging the order of the coordinates, we can concentrate our attention on the first
coordinate: hence we aim at proving that
\begin{equation}\label{VCISTPOK-7}
\sum_{j=1}^3 V_j\partial_j V_1=\frac{1}{2}\partial_1 |V|^2-
\big(V\times (\curl V)\big)_1.\end{equation}
To check this, we calculate:
\begin{eqnarray*}&&
\big(V\times (\curl V)\big)_1-\frac{1}{2}\partial_1 |V|^2\\&=&V_2(\curl V)_3-V_3(\curl V)_2-
V\cdot\partial_1 V\\
&=& V_2\big( \partial_1 V_2-\partial_2 V_1\big) - V_3 \big( \partial_3 V_1-\partial_1 V_3\big)-V_1\partial_1 V_1
-V_2\partial_1 V_2-V_3\partial_1 V_3\\
&=& - V_2\partial_2 V_1 - V_3 \partial_3 V_1-V_1\partial_1 V_1,
\end{eqnarray*}
which proves~\eqref{VCISTPOK-7} and thus~\eqref{VCISTPOK}.

We also recall another vector calculus identity, valid for all smooth vector fields~$V=(V_1,V_2,V_3)$, $W=(W_1,W_2,W_3):\R^3\to\R^3$:
\begin{equation}\label{nhfhgtfredwM3S54E4y5C236terdhrCHE}
\curl(V\times W)
= \div W\;V- \div V \;W+(W\cdot\nabla)\, V- (V\cdot\nabla)\,W.\end{equation}
To prove this, we can write~$V=V_1e_1+V_2e_2+V_3 e_3$ and notice that, since~\eqref{nhfhgtfredwM3S54E4y5C236terdhrCHE} is linear in~$V$, it suffices to prove it when~$V$
reduces to each of the components~$V_j e_j$.

Thus, up to reordering the coordinates, we can focus on the case in which~$V$ is actually~$V_1 e_1$.
Hence we calculate
\begin{eqnarray*}
&&\curl((V_1e_1)\times W)+ \div( V_1e_1) \;W-(W\cdot\nabla)\, V_1e_1+( (V_1e_1)\cdot\nabla)\,W
\\&=&
\curl(0,-V_1W_3,V_1W_2)
+ \partial_1 V_1\,W-\sum_{j=1}^3W_j\partial_j V_1e_1+V_1\partial_1 W\\
&=&\big( \partial_2(V_1W_2)+\partial_3(V_1W_3),\,-\partial_1(V_1W_2),\,-\partial_1(V_1W_3)\big)
+ \partial_1 V_1\,W-\sum_{j=1}^3W_j\partial_j V_1e_1+V_1\partial_1 W
\\&=&
\big(
\partial_2V_1\,W_2+V_1\partial_2 W_2+\partial_3 V_1\,W_3+V_1\partial_3 W_3+\partial_1 V_1\,W_1
+V_1\partial_1W_1,\,
0,\,0\big)-\sum_{j=1}^3W_j\partial_j V_1e_1\\&=&
\big(V_1\partial_2 W_2+V_1\partial_3 W_3
+V_1\partial_1W_1,\,
0,\,0\big)
\\&=& \div W\;V_1e_1
\end{eqnarray*}
and this completes the proof of~\eqref{nhfhgtfredwM3S54E4y5C236terdhrCHE}.

\begin{figure}
  \centering
  \includegraphics[width=.5\linewidth]{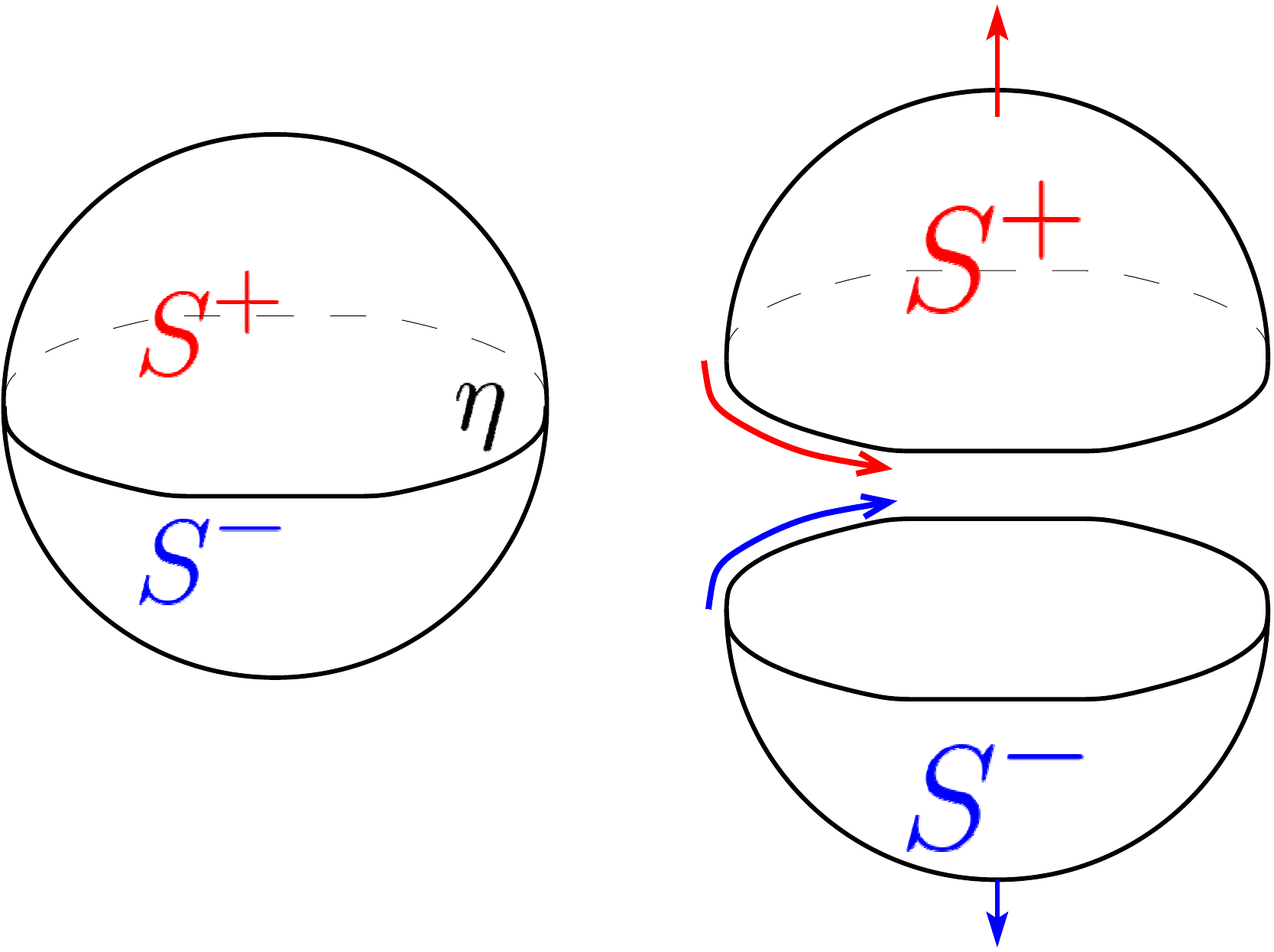}
 \caption{\sl The geometric argument used to prove~\eqref{YY08ojew9oik903hyto3298ythgfjb8ihyri-1urX}.}\label{P2DAILATENSIWASTREDItangeFIlAFINASDJGIOSE}
\end{figure}

It is also useful to recall the vector calculus identity
\begin{equation}\label{YY08ojew9oik903hyto3298ythgfjb8ihyri-1urX} \div(\curl V)=0.\end{equation}
For this, we take any ball~$B\subset\R^3$ and we consider the spherical surface~$S:=\partial B$.
We split~$S$ into two hemispheres~$S^+$ and~$S^-$
with a common equator~$\eta$, which is the intersection between the boundaries of the surfaces~$S^+$ and~$S^-$.
Looking at~$\eta$ as a curve, the natural direction of travel along~$\eta$ when considered as the boundary of~$S^+$
is opposite to the one obtained when considering it as the boundary of~$S^-$, see Figure~\ref{P2DAILATENSIWASTREDItangeFIlAFINASDJGIOSE}
(here, we are endowing~$S^+$ and~$S^-$ with the external unit vector field of~$\partial B$).
As a result, the circulation of a vector field~$V$ along~$\eta$ considered as the boundary of~$S^+$ (that we denote by~${\mathcal{C}}^+(V)$)
is opposite to the circulation of~$V$ along~$\eta$ considered as the boundary of~$S^-$
(that we denote by~${\mathcal{C}}^-(V)$), that is
$$ {\mathcal{C}}^+(V)=-{\mathcal{C}}^-(V).$$
In addition, by Stokes' Theorem, we know that~${\mathcal{C}}^\pm(V)$ agrees with the flux
of~$\curl V$ through the surface~$S^+$. 

Using these bits of information and the Divergence Theorem, we conclude that
\begin{eqnarray*}
&&0={\mathcal{C}}^+(V)+{\mathcal{C}}^-(V)=
\int_{S^+}\curl V\cdot\nu+\int_{S^-}\curl V\cdot\nu=\int_{\partial B}\curl V\cdot\nu=
\int_{B}\div(\curl V).
\end{eqnarray*}
{F}rom the arbitrariness of the ball~$B$, we obtain~\eqref{YY08ojew9oik903hyto3298ythgfjb8ihyri-1urX}, as desired.

We also point out that, for any smooth function~$\psi$,
\begin{equation}\label{08ojew9oik903hyto3298ythgfjb8ihyri-1urX} \curl\nabla\psi=0.\end{equation}
This can be checked by employing Stokes' Theorem. Indeed, if~$\Sigma$ is any smooth oriented element of surface in~$\R^3$
with boundary~$\partial \Sigma $, the flux of~$\curl\nabla\psi$ through the surface~$\Sigma$
is equal to the circulation of~$\nabla\psi$ over the loop~$\partial\Sigma$. The latter object, if we describe~$\partial\Sigma$
by a closed curve~$\gamma:[0,L]\to\R^3$ with arc length parameterization, can be written as
$$ \int_0^L \nabla\psi(\gamma(\tau))\cdot\dot\gamma(\tau)\,d\tau=
\int_0^L \frac{d}{d\tau}\psi(\gamma(\tau))\,d\tau=\psi(\gamma(L))-\psi(\gamma(0))=0.$$
Since this is valid for an arbitrary surface~$\Sigma$, the proof of~\eqref{08ojew9oik903hyto3298ythgfjb8ihyri-1urX} is complete.

Now, to prove~\eqref{CURLDESE}, we let~$\omega:=\curl v$, we use the momentum balance equation 
and the barotropicity in~\eqref{CURLDESE1}
in combination with~\eqref{VCISTPOK} and we see that
\begin{equation}\label{dui39urteifjewgtlie4y9548796575903549873204327957608497685490549030594}
\begin{split}
0\,&=\curl\left( \partial_t v+ (v\cdot\nabla) v+\frac{\nabla p}{\rho}\right)\\&=\curl\left( \partial_t v+ \frac{1}{2}\nabla |v|^2-
v\times (\curl v)+\frac{\nabla p}{g(p)}\right)\\
&=\partial_t \omega+\curl\left(  \frac{1}{2}\nabla |v|^2-
v\times \omega+\frac{\nabla p}{g(p)}\right).
\end{split}\end{equation}
Now we observe that
\begin{eqnarray*} &&
\curl\left(\frac{\nabla p}{g(p)}\right)= \left(\partial_2 \left(\frac{\partial_3 p}{g(p)}\right)-
\partial_3 \left(\frac{\partial_2 p}{g(p)}\right),\partial_3 \left(\frac{\partial_1 p}{g(p)}\right)-
\partial_1 \left(\frac{\partial_3 p}{g(p)}\right),
\partial_1 \left(\frac{\partial_2 p}{g(p)}\right)-
\partial_2 \left(\frac{\partial_1 p}{g(p)}\right)
\right)\\&&\qquad=\left( \frac{\partial_{23} p}{g(p)}-\frac{\partial_3p\,g'(p)\,\partial_2p}{g^2(p)}
-\frac{\partial_{23} p}{g(p)} + \frac{\partial_2 p\,g'(p)\,\partial_3p}{g^2(p)},\right.\\ &&\qquad\qquad\qquad\left.
\frac{\partial_{13} p}{g(p)}-\frac{\partial_1 p\,g'(p)\,\partial_3p}{g^2(p)}
-\frac{\partial_{13} p}{g(p)}+\frac{\partial_3 p\,g'(p)\,\partial_1p}{g^2(p)},\right.\\ &&\qquad\qquad\qquad\left.
\frac{\partial_{12} p}{g(p)}-\frac{\partial_2 p\,g'(p)\,\partial_1p}{g^2(p)}
-\frac{\partial_{12} p}{g(p)}+\frac{\partial_1 p\,g'(p)\,\partial_2p}{g^2(p)}
\right)\\&&\qquad=(0,0,0).
\end{eqnarray*}
Hence, from this, \eqref{08ojew9oik903hyto3298ythgfjb8ihyri-1urX} and~\eqref{dui39urteifjewgtlie4y9548796575903549873204327957608497685490549030594},
\begin{equation*}
\begin{split}
0=\partial_t \omega-\curl\left( v\times \omega\right).
\end{split}\end{equation*}
This and~\eqref{nhfhgtfredwM3S54E4y5C236terdhrCHE} lead to
\begin{equation*}
\partial_t \omega=
\curl(v\times \omega)
= \div \omega\;v- \div v \; \omega+( \omega\cdot\nabla)\, v- (v\cdot\nabla)\, \omega.\end{equation*}
Thus, combining the incompressibility condition in~\eqref{CURLDESE1}
and~\eqref{YY08ojew9oik903hyto3298ythgfjb8ihyri-1urX},
\begin{equation*}
\partial_t \omega=( \omega\cdot\nabla)\, v- (v\cdot\nabla)\, \omega\end{equation*}
and therefore, following the fluid parcel~$x(t)$,
\begin{equation*}\begin{split}&
\frac{d}{dt} \omega(x(t),t)=
\partial_t \omega(x(t),t)+(\dot x(t)\cdot\nabla)\omega(x(t),t)\\&\qquad
=\partial_t \omega(x(t),t)+\big(v( x(t),t)\cdot\nabla\big)\omega(x(t),t)
=\big( \omega(x(t),t)\cdot\nabla\big)\, v(x(t),t).
\end{split}\end{equation*}
As a result, if~$Z(t)=(Z_1(t),Z_2(t),Z_3(t)):=\omega(x(t),t)$, we have that~$Z$ is a solution of the Cauchy problem
for ordinary differential equations
$$ \begin{dcases}\displaystyle
\dot Z(t)=\sum_{k=1}^3 Z_k(t)\, \partial_k v(x(t),t),\\
Z(0)=0.
\end{dcases}$$
By the uniqueness result for solutions of ordinary differential equations we thereby infer that~$Z(t)=0$ for all~$t\in\R$
and this completes\footnote{As a technical detail, we point out that we are freely assuming that if, for a given function~$f$,
we know that~$f(x(t),t)=0$ for all times~$t$, then we can infer that~$f(x,t)=0$ \label{LAFOPRI}
for all spatial positions~$x$ and times~$t$. This is because we are implicitly assuming that the flow
of the fluid exists for all times: hence, given any position~$x$ we can ``flow~$x$ backwards for a time~$t$'',
that is consider the fluid parcel evolution starting at some point~$x_0$ at time~$0$ arriving at~$x$ at time~$t$
(this gives us the possibility of choosing~$x(t)=x$).}
the proof of~\eqref{CURLDESE}.

\subsection{Propagation of sound waves}\label{PROSAWA0qiodjwfl0owfej23t98y4823rty}

Now we retake the set of Euler fluid flow equations for the velocity~$v$ of an inviscid fluid,
namely we look at the mass conservation equation in~\eqref{OJS-PJDN-0IHGDOIUGDBV02ujrfMTE}, at
the momentum balance equation in~\eqref{EUMSD-OS-D} and at a constitutive relation
linking the pressure~$p$ and the density~$\rho$, namely (neglecting the gravity effects)
\begin{equation}\label{EUMSD-OS-D-TUT}
\begin{dcases}&
\partial_t \rho+\div(\rho v)=0,\\&
\rho \partial_t v+\rho (v\cdot\nabla) v=-\nabla p,\\&
p=f(\rho),
\end{dcases}
\end{equation}
for some function~$f$. In practice, it is useful to take~$f$ to be strictly increasing (the higher the density of the fluid, the higher the pressure produced; this would also give that the fluid is barotropic in the setting of Section~\ref{BARPotA}).
We can think of~\eqref{EUMSD-OS-D-TUT} as a very simple model for a gas
and we aim here at understanding the propagation of sound waves in such a medium.

For this objective, we first point out that a solution of~\eqref{EUMSD-OS-D-TUT}
is provided by~$(v,p,\rho)=(0,p_0,\rho_0)$, for every~$\rho_0\in(0,+\infty)$ and~$p_0:=f(\rho_0)$.
This configuration corresponds to the physical situation of a gas at rest, with constant density and pressure.

Suppose now that we perturb this configuration by creating a small variation of the gas pressure,
for instance by singing, or by playing\footnote{To play guitar, see the forthcoming Section~\ref{GUITAR}.}
a guitar. The idea is thus to look for (at least approximate) solutions of the form
$$ v(x,t)=\e v_1(x,t),\qquad p(x,t)=p_0+\e p_1(x,t)\qquad{\mbox{and}}\qquad \rho(x,t)=\rho_0+\e\rho_1(x,t).$$
In this setting, $\e$ is a small parameter and our goal is to determine the functions~$(v_1,p_1,\rho_1)$ in order to satisfy the set of equations in~\eqref{EUMSD-OS-D-TUT}, formally up to negligible errors in~$\e$. To this end, we observe that
$$ p_0+\e p_1=p=f(\rho)=f(\rho_0+\e \rho_1)=f(\rho_0)+\e f'(\rho_0) \rho_1+o(\e)=p_0+\e f'(\rho_0) \rho_1+o(\e),$$
leading to the choice
\begin{equation*}
p_1=f'(\rho_0) \rho_1,
\end{equation*}
up to higher orders in~$\e$ that we here sloppily disregard.

\begin{figure}
  \centering
  \includegraphics[width=.35\linewidth]{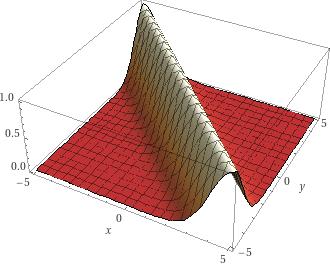}$\quad$
  \includegraphics[width=.35\linewidth]{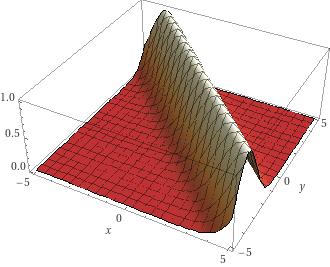}\\
   \includegraphics[width=.35\linewidth]{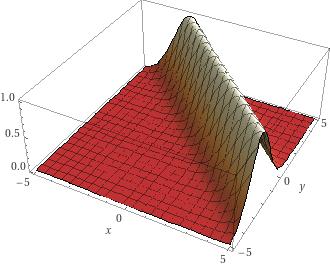}$\quad$
  \includegraphics[width=.35\linewidth]{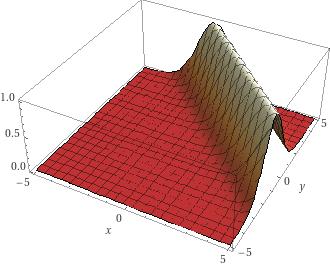}
 \caption{\sl A plane wave.}\label{PLAWASTREDItangeFI}
\end{figure}

Moreover,
\begin{equation*}
0=\partial_t \rho+\div(\rho v)=\e\partial_t \rho_1+\e\rho_0\div v_1+o(\e)
\end{equation*}
and
\begin{equation*}
0= \rho \partial_t v+\rho (v\cdot\nabla) v+\nabla p=\e\rho_0 \partial_t v_1+\e\nabla p_1+o(\e)=
\e\rho_0 \partial_t v_1+\e f'(\rho_0)\nabla \rho_1+o(\e).
\end{equation*}
These observations give
\begin{equation*}
\begin{dcases}
&\partial_t \rho_1+ \rho_0\div v_1=0,\\
&\rho_0 \partial_t v_1+f'(\rho_0)\nabla \rho_1=0.
\end{dcases}\end{equation*}
As a consequence,
\begin{equation*}
\partial_{tt} \rho_1=-\partial_t\big(\rho_0\div v_1\big)
=-\div \big(\rho_0\partial_t v_1\big)=
\div\big(
f'(\rho_0)\nabla \rho_1
\big)
=f'(\rho_0)\Delta \rho_1
.\end{equation*}
That is, the perturbed density~$\rho_1$
satisfies the equation
\begin{equation}\label{WAYEBJJD121t416JH098327uyrhdhc832nbcM}
\partial_{tt}u=c^2 \Delta u
\end{equation}
with~$c:=\sqrt{f'(\rho_0)}>0$.
Equation~\eqref{WAYEBJJD121t416JH098327uyrhdhc832nbcM} is called\footnote{In terms of the classification
mentioned in footnote~\ref{CLASSIFICATIONFOOTN}
on page~\pageref{CLASSIFICATIONFOOTN}, equation~\eqref{WAYEBJJD121t416JH098327uyrhdhc832nbcM} is hyperbolic.
Indeed, we can take here~$N=n+1$, $X=(x,t)$ and
$$ a_{ij}=\begin{dcases}
c^2 & {\mbox{ if }}i=j\in\{1,\dots,n\},\\
-1 & {\mbox{ if }}i=j=n+1,\\
0 & {\mbox{ otherwise,}}
\end{dcases}$$
thus producing $n$ strictly positive, and one strictly negative, eigenvalues.}
in jargon the ``wave equation''\index{wave equation}.

To have a feeling of the ``propagation of waves'' encoded in equation~\eqref{WAYEBJJD121t416JH098327uyrhdhc832nbcM},
we can consider a smooth function~$u_\star:\R\to\R$ and a direction~$\omega\in\partial B_1$,
and thus define~$u_\omega(x,t):=u_\star(\omega\cdot x-ct)$. Notice that~$u_\omega$ is indeed a solution of~\eqref{WAYEBJJD121t416JH098327uyrhdhc832nbcM},
physically corresponding to a traveling plane wave\index{plane wave}. \label{01ojehniTNasIplamnewa}
Indeed, its evolution in time corresponds to a translation of~$u_\star$ along the direction~$\omega$
with speed~$c$. Also, for~$\kappa\in\R$, the parallel hyperplanes~$\{\omega\cdot x-ct=\kappa\}$
(again, traveling with constant speed~$c$ in direction~$\omega$)
correspond to the level sets~$\{ u_\omega=\kappa_\star\}$, where~$\kappa_\star:=u_\star(\kappa)$
(the level sets of a wave are often called ``wavefronts''
and they have special importance since they correspond to the surfaces on which,
at a given moment of time, all particles of the medium undergo the same motion).

\begin{figure}
  \centering
  \includegraphics[width=.25\linewidth]{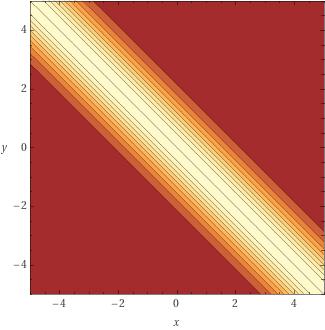}$\quad$
  \includegraphics[width=.25\linewidth]{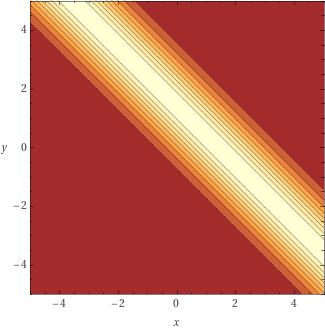}\\
   \includegraphics[width=.25\linewidth]{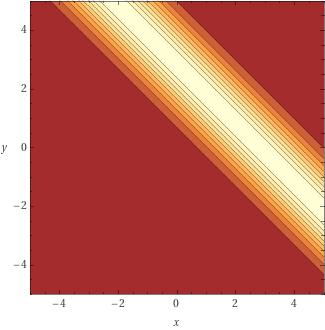}$\quad$
  \includegraphics[width=.25\linewidth]{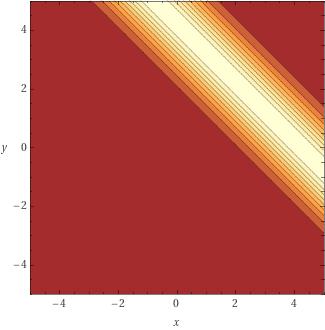}
 \caption{\sl Wavefronts of the plane wave in Figure~\ref{PLAWASTREDItangeFI}.}\label{PLAWASTREDItangeFIlASE}
\end{figure}

See Figure~\ref{PLAWASTREDItangeFI} for a series of frames of a traveling plane wave
(the picture is obtained by choosing~$u_\star(r):=e^{-r^2}$, $\omega:=\left(\frac{\sqrt2}{2},\frac{\sqrt2}{2}\right)$, $c=1$
and instants of time~$t\in\{0,1,2,3\}$).
See also Figure~\ref{PLAWASTREDItangeFIlASE}
for the pictures of the corresponding level sets of the wave
(i.e., of the corresponding wavefronts).

\subsection{Hydrodynamics doesn't always work right}\label{WRItvHS-SE1}

Now that we built some confidence in the equations of hydrodynamics,
it is time to question and challenge our own knowledge, to appreciate its power as well as its limitation.
After all, it's not so important to know, but rather to know of not knowing.

In particular, the great Euler (``Master of us all'' according to Laplace) allowing,
we will recall here a rather striking, and perhaps surprising, flaw of the hydrodynamic equations presented in Section~\ref{EUMSD-OS-32456i7DNSSE234R:SEC} in the description of aerodynamic forces.

We will focus here on the notion of ``drag''\index{drag}, which is the
force acting in the opposite direction
to an object moving in a fluid (say, the air resistance, for example).
This is not just a mathematical curiosity: for instance, according to Wikipedia
{\tt https://en.wikipedia.org/wiki/Drag\_(physics)},
``induced drag tends to be the most important component for airplanes during take-off or landing flight''
and, more importantly,
``in the physics of sports, the drag force is necessary to explain the motion of balls, javelins, arrows and frisbees and the performance of runners and swimmers''.\medskip

\begin{figure}
                \centering
                \includegraphics[width=.35\linewidth]{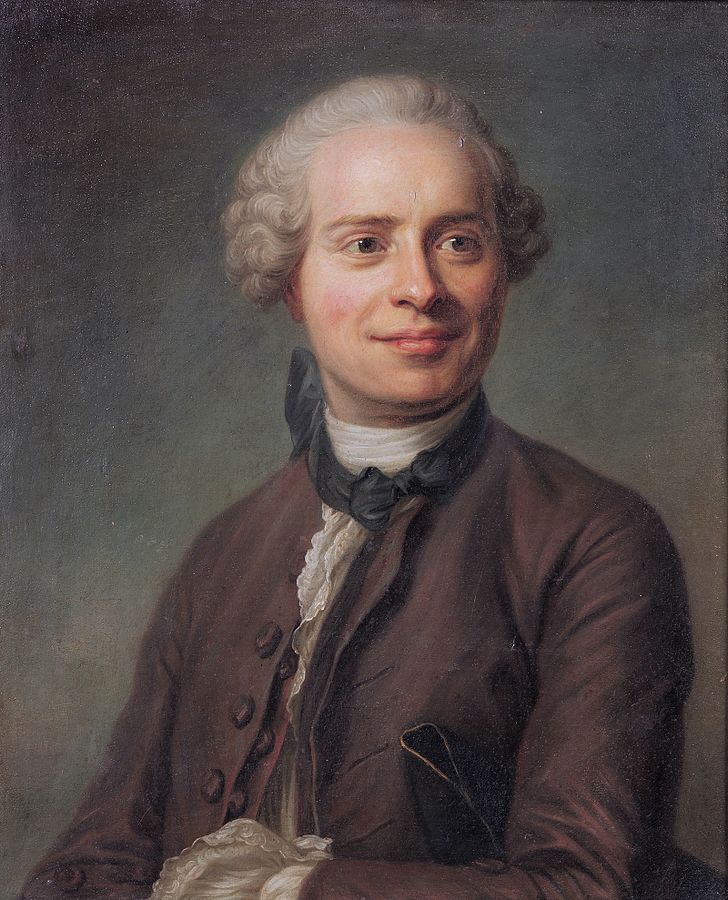}
        \caption{\sl Portrait of Jean Le Rond d'Alembert (Public Domain image from
        Wikipedia).}\label{HAFODAKEDFUMSldRGIRA4AXELHARLROA7789GIJ7solFUMHDNOJHNFOJED231}
\end{figure}

The issue about the drag that we discuss here was discovered by Jean le Rond d'Alembert in~1752.
Often, this remarkable discovery is called with the name of
``d'Alembert's paradox''\index{d'Alembert's paradox};
see Figure~\ref{HAFODAKEDFUMSldRGIRA4AXELHARLROA7789GIJ7solFUMHDNOJHNFOJED231}
for a portrait of d'Alembert's intense glance (by Maurice Quentin de La Tour). In a nutshell, we will prove that
\begin{equation}\label{PARADDA}
\begin{split}
&{\mbox{if a body is moving with constant velocity}}\\&{\mbox{in an inviscid, irrotational and steady flow
with constant density,}}\\
&{\mbox{then the corresponding drag force is zero.}}
\end{split}
\end{equation}
Of course, we may wonder how such a statement can be coherent with our own
everyday experience, in which the effect of air resistance is usually apparent and quite decisive
(at least when we catch a flight or we kick a ball). Though the full understanding
of hydrodynamics is likely way above the present possibilities of science,
one common explanation of the discrepancy between the claim in~\eqref{PARADDA} and real life situations
relies on remarking that the fluid's description in~\eqref{PARADDA} is too idealistic.
In particular, the occurrence of the paradox is usually attributed
to the neglected effects of viscosity\footnote{The importance of viscosity actually reflects the importance of the Laplace operator in the model: compare equations~\eqref{EUMSD-OS-D} and~\eqref{EUMSD-OS-DNS}.

Though we do not explicitly use this here, for completeness we point out that
an irrotational fluid with constant density is automatically incompressible
(because in this case~$\div v=\frac{1}\rho\div(\rho v)=\frac{1}\rho\big(\partial_t\rho+\div(\rho v)\big)=0$
by the continuity equation
in~\eqref{OJS-PJDN-0IHGDOIUGDBV02ujrfMTE}), and
an incompressible and  irrotational fluid is automatically inviscid
(hence, the lack of viscosity is also a byproduct of two assumptions, namely the incompressibility and
the lack of vortexes, which, together, entail very restrictive byproducts): indeed, using the vector calculus
identity
$$ \nabla (\div V)-\curl(\curl V)=\Delta V,$$
(see e.g.~\eqref{KMSJUNDNSMDEJDN2} for a proof) we have that for incompressible and irrotational fluid it holds that
$$ \Delta v=\nabla (\div v)-\curl(\curl v)=0-0=0,$$
hence the Navier-Stokes equation in~\eqref{EUMSD-OS-DNS}
boils down to the Euler equation in~\eqref{EUMSD-OS-D}.}
(roughly speaking, one of the effects of viscosity~\cite{zbMATH02651679} is to produce
thin boundary layers\index{boundary layer} near the surface of the moving object, which may entail
friction, flow separation and a low-pressure wake, leading to pressure drag).
However, the official resolution of the question is possibly debatable in its full generality:
for instance, in relation to the paradoxical features of~\eqref{PARADDA},
Garrett Birkhoff~\cite{MR0038180} states\footnote{Garrett Birkhoff's statements
were in turn criticized by James J. Stoker~\cite{MR1565341}
(who was concerned about the possibility that readers were
``very likely to get wrong ideas about some of the important
and useful achievements in hydrodynamics''
being misled by the ``negative aspects of the theory''). Some of the original statements in~\cite{MR0038180}
were revised in the second edition~\cite{MR0122193}.

All in all, much more knowledge has to be built before we reach a complete understanding of
fluid mechanics and hydrodynamics.}
that: ``I think that to attribute them all to the neglect of viscosity is an unwarranted oversimplification. The root lies deeper, in lack of precisely that deductive rigor whose importance is so commonly minimized by physicists and engineers''.
\medskip

Now we work out the mathematical details needed to prove the claim in~\eqref{PARADDA}.
To achieve this goal, we recall
the Euler equation in~\eqref{EUMSD-OS-D}
and we complement it with the mass transport equation in~\eqref{OJS-PJDN-0IHGDOIUGDBV02ujrfMTE},
the assumption that the density is constant
and an irrotational condition. We also assume that the flow is ``steady'':
this word usually refers to a situation in
which the fluid properties at a point in the system do not change over time
(they only depend on the point); concretely, to avoid ambiguities, we will suppose here
that the relevant quantities of the fluid, such as velocity and pressure, follow the body moving at constant
velocity~$v_0$ (that is, they are not functions of~$x$ and~$t$ separately, but of~$x-v_0t$).
Equivalently,
according to the principle of inertia
discussed on page~\pageref{PRINERY}, if we take a system of reference
moving together with the body at constant speed~$v_0$, then
the relevant quantities of the fluid would be constant in time (in particular, $\partial_tv=0$).

Hence, in this system, 
by~\eqref{OJS-PJDN-0IHGDOIUGDBV02ujrfMTE} and~\eqref{EUMSD-OS-D},
the velocity~$v$
of a steady and irrotational fluid, subject to no gravity effects, with constant density (say~$\rho:=1$)
can be described by the set of equations (valid outside the body~$\Omega$, which is supposed to be
still in this reference system)
\begin{equation}\label{VCISTPOK-2}
\begin{dcases}
(v\cdot\nabla) v=-\nabla p,\\
\div v=0,\\
\curl v=0.
\end{dcases}
\end{equation}
Now, combining~\eqref{VCISTPOK} with the irrotationality condition in~\eqref{VCISTPOK-2} we deduce that
\begin{equation}\label{VCISTPOK-3}
(v\cdot\nabla )v=\frac{1}{2}\nabla |v|^2.\end{equation}
We also consider the differential form~$\omega:=v_1\,dx_1+v_2\,dx_2+v_3\,dx_3$ and we observe that
\begin{eqnarray*} d\omega&=&-\partial_2 v_1\,dx_1\wedge dx_2-
\partial_3 v_1\,dx_1\wedge dx_3+\partial_1 v_2\,dx_1\wedge dx_2\\&&\qquad
-\partial_3 v_2\,dx_2\wedge dx_3+\partial_1v_3\,dx_1\wedge dx_3
+\partial_2v_3\,dx_2\wedge dx_3\\&=&(\curl v)_3\,dx_1\wedge dx_2-(\curl v)_2\,dx_1\wedge dx_3+
(\curl v)_1\,dx_1\wedge dx_2\\&=&0,
\end{eqnarray*}
thanks again to the irrotationality condition in~\eqref{VCISTPOK-2}.

Accordingly, the differential form~$\omega$ is closed, and therefore exact, thanks to the
Poincar\'e Lemma (see e.g.~\cite[Theorem~8.3.8]{MR1209437}). This gives that there exists a velocity potential function~$\varphi$
satisfying~$\omega=d\varphi$, that is
\begin{equation}\label{jiaknfujsuE0ujkjHAHLoefvujs98ou38jhPD98jhdjhfX8}
v=\nabla \varphi. \end{equation}
We stress that~$\varphi$ is harmonic outside~$\Omega$, since~$0=\div v=\div\nabla\varphi=\Delta\varphi$.
Moreover, the first equation in~\eqref{VCISTPOK-2},
\eqref{VCISTPOK-3} 
and~\eqref{jiaknfujsuE0ujkjHAHLoefvujs98ou38jhPD98jhdjhfX8} entail\footnote{In Section~\ref{AIRFVBD01ogt9u432gtb}
we will have a closer
look at expressions such as the one in~\eqref{kSUDSpREdfskljhREHGkozWI23f3aziw2}
and we will interpret them as a sort of \index{Bernoulli's Principle} Bernoulli's Principle.} that
\begin{equation}\label{kSUDSpREdfskljhREHGkozWI23f3aziw2}\begin{split}&
0= (v\cdot\nabla) v+ \nabla p=
\frac{1}{2}\nabla |v|^2+ {\nabla p} =
\frac{1}{2}\nabla |\nabla\varphi|^2+ {\nabla p}=
\nabla\left( \frac{1}{2}|\nabla\varphi|^2\right)+ {\nabla p}.\end{split}
\end{equation}

With this, we can now compute the drag that the air exerts on the moving body.
For this, we recall that the force~$F$ of the air is produced by the pressure, acting normally on the surface
of the body (see~\eqref{90uyoihfs83wEUSlNOE}). Hence, calling the moving object~$\Omega$ and~$\nu$ its outer normal, we have that
\begin{equation}\label{FEUYT90hfbvjjxJ} F=-\int_{\partial\Omega} p\nu.\end{equation}
The minus sign here above is just recalling that positive pressures are assumed to go in the opposite direction
with respect to the exterior normal.
Since the drag~$D$ is the force acting in the direction of motion of the body, using the Divergence Theorem we thus conclude that
\begin{equation}\label{EXDIVTHMQU} D=F\cdot\frac{v_0}{|v_0|}
=-\int_{\partial\Omega} p\frac{v_0\cdot\nu}{|v_0|}=\int_{\R^3\setminus\Omega}
\div\left(p\frac{v_0}{|v_0|}\right)=\int_{\R^3\setminus\Omega} \frac{\nabla p\cdot v_0}{|v_0|}
.\end{equation}
Notice that to apply the Divergence Theorem here to the infinite region~$\R^3\setminus\Omega$
(that corresponds to the region occupied by the fluid) we are implicitly supposing that the pressure decays sufficiently fast at infinity
(that is, the pressure disturbance is essentially localized in the vicinity of the moving object, which is a reasonable\footnote{Once again, in this chapter, if an argument sounds convincing we buy it! We will be much more skeptical about heuristic reasonings from next chapter on. For instance, when we need to integrate on exterior domains in Section~\ref{Capacitance-and-capacity}, the reader will appreciate how careful we will be in checking explicitly the appropriate decay properties of the functions involved (which is perhaps a delicate and annoying, albeit absolutely necessary, detail to be taken care of to avoid nonsensical computations).

In any case, let us mention that the possibility of utilizing the Divergence Theorem in an exterior domain
in~\eqref{EXDIVTHMQU} is quite a delicate matter (see e.g. the thoughtful argument in~\cite[formulas~(6.4.2) and~(6.4.3)]{MR1744638}) \label{DECA2DSFGBD}
and it is a byproduct of the decay at infinity of the fluid velocity and pressure which
is modeled on that of harmonic functions.

Arguably, in the setting of these notes, alternative arguments to the ones in~\cite[formulas~(6.4.2) and~(6.4.3)]{MR1744638}
could be provided by using
the Green's Representation Formula in~\eqref{300G}
(say, in the domain~$B_R\setminus\overline{\Omega}$,
sending~$R\to+\infty$): for this one would need to assume that the velocity potential
has a limit at infinity (thus, from Cauchy's Estimates
in Theorem~\ref{CAUESTIMTH}, one can also bound the derivative of the potential along~$\partial B_R$).

In any case, the decay estimates of the velocity fields rely
on the fact that the fluid cannot penetrate inside the object (as formalized in~\eqref{PSLmoLMSxHNS-AR}),
since this information allows one to get rid of one order of magnitude in the corresponding estimates. This is an interesting
physical feature since it reveals the ``global'' character of fluid dynamics, in which some features 
in the proximity of the moving body have a significant influence on the velocity field at infinity, and vice versa.

In general, the decay analysis for solutions of fluid mechanics equations
is a delicate matter, see e.g.~\cite{MR0481652} and the references therein.}
physical assumption, though quantifying it precisely and checking it rigorously would require some technical skills).

We also point out that, for all~$k\in\{1,2,3\}$,
\begin{eqnarray*}&&
\partial_k\left(\frac{1}{2}|\nabla\varphi|^2\right)=
\sum_{j=1}^3 \partial_j\varphi\,\partial_{jk}\varphi
=\sum_{j=1}^3 \partial_j\big(\partial_j \varphi\,\partial_{k}\varphi\big)-\Delta\varphi\partial_k\varphi
=\div\big(\partial_{k}\varphi\nabla\varphi\big)-0.\end{eqnarray*}
Consequently,
\begin{eqnarray*}&&
v_0\cdot\nabla\left(\frac{1}{2}|\nabla\varphi|^2\right)=
\sum_{k=1}^3 v_{0,k}\partial_k\left(\frac{1}{2}|\nabla\varphi|^2\right)=
\sum_{k=1}^3 v_{0,k}\div\big(\partial_{k}\varphi\nabla\varphi\big)
= \div\Big((v_0\cdot\nabla\varphi)\nabla\varphi\Big).
\end{eqnarray*}
This entails that
\begin{equation*}\begin{split}&
\int_{\R^3\setminus\Omega} \frac{ v_0}{|v_0|}\cdot\nabla\left(
\frac{1}{2}|\nabla\varphi|^2\right)=\frac{1}{|v_0|}
\int_{\R^3\setminus\Omega} \div\Big((v_0\cdot\nabla\varphi)\nabla\varphi\Big)\\&\qquad
=-\frac{1}{|v_0|}
\int_{\partial\Omega} (v_0\cdot\nabla\varphi)(\nabla\varphi\cdot\nu)
=-\frac{1}{|v_0|}
\int_{\partial\Omega} (v_0\cdot v)(v\cdot\nu).\end{split}
\end{equation*}
Hence, recalling~\eqref{kSUDSpREdfskljhREHGkozWI23f3aziw2} and~\eqref{EXDIVTHMQU},
\begin{equation}\label{KSM-dso02or3kf03ptjhohwgi02034}
D=\int_{\R^3\setminus\Omega} \frac{\nabla p\cdot v_0}{|v_0|}=
-\int_{\R^3\setminus\Omega} \frac{v_0}{|v_0|}\cdot\nabla\left( \frac{1}{2}|\nabla\varphi|^2\right)=
\frac{1}{|v_0|}
\int_{\partial\Omega} (v_0\cdot v)(v\cdot\nu).
\end{equation}

We now stress that since the air cannot penetrate inside the object, necessarily
\begin{equation}\label{PSLmoLMSxHNS-AR}
{\mbox{the normal component
of the fluid velocity vanish along~$\partial\Omega$,}}\end{equation} namely~$v\cdot\nu=0$
on~$\partial\Omega$.
This and~\eqref{KSM-dso02or3kf03ptjhohwgi02034} yield that~$D=0$,
which completes the proof of d'Alembert's paradox in~\eqref{PARADDA}.

See e.g.~\cite{MR631832, MR1218879, MR1744638, MR3837529} and the references therein for additional information
on d'Alembert's paradox.

\begin{figure}
  \centering
  \includegraphics[width=.6\linewidth]{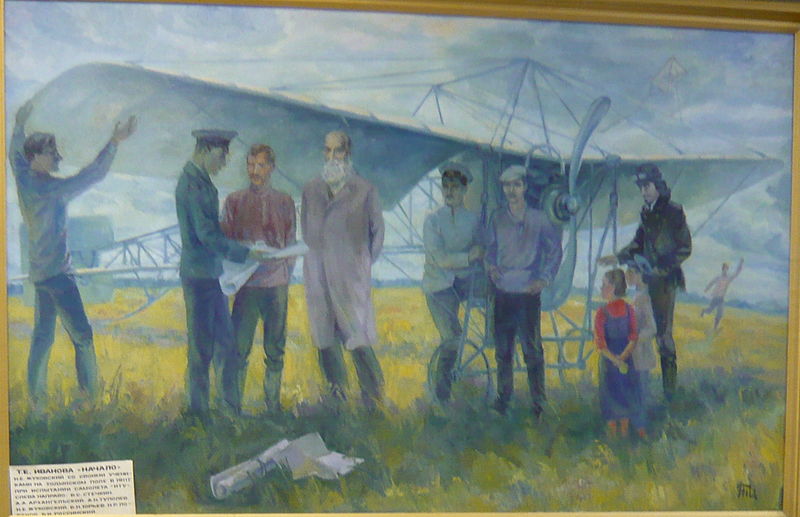}
 \caption{\sl Portrait of Nikolai Yegorovich Joukowski,
Museum of Moscow Aviation Institute (image by Just from
 Wikipedia, licensed under the Creative Commons Attribution-Share Alike 4.0 International license).}\label{P2DAILACRPOGGSTENSIWASTREDItangeFIlASE}
\end{figure}

\subsection{Lift of an airfoil}\label{AIRFVBD01ogt9u432gtb}

To safely recover from the shock received in Section~\ref{WRItvHS-SE1}
by d'Alembert's paradox, we give now a positive result
in terms of the aerodynamic lift of an airfoil\index{airfoil}. This result was established independently
by Martin Kutta and Nikolai Yegorovich Joukowski in the early twentieth century and, in its simplest formulation,
can be stated\footnote{Though the theory that
we present in these pages is way too crude to
comprise aerodynamics in its full complexity, the importance of a statement like~\eqref{JTTA} in the theory of
flight appears to be paramount, as confirmed by the solid recognition bestowed to the scientists involved. For instance,
the Russian Air Force Academy and one of the airports in Moscow
were named after Joukowski. See also Figure~\ref{P2DAILACRPOGGSTENSIWASTREDItangeFIlASE}.

In terms of real world applications, a common belief is indeed that the setting in~\eqref{JTTA}
is too idealized to detect the intricate patterns produced by real fluids in the vicinity of a traveling object: however,
it is also believed that the circulation detected in~\eqref{JTTA} does reflect significant physical
information since in many concrete situations
the flow around a thin airfoil is composed of a narrow viscous region near the body
(a sort of ``boundary layer'') outside which the idealized 
and inviscid description of the flow in~\eqref{JTTA} turns out to be sufficiently realistic. 

The gist is then to apply the setting in~\eqref{JTTA} not quite to the traveling body but to the aggregate of the body and its own boundary layer:
in particular, the loop to compute the circulation in~\eqref{JTTA} must be chosen outside this boundary layer.

With this, the boundary layer gives any traveling object an ``effective shape'' that may be different from its physical shape by accounting for the region in which the velocity changes from zero at the surface to the stream value away from the surface: the effectiveness of~\eqref{JTTA} is thus to incorporate, as much as possible, the effects of turbulence and skin friction into the effective shape described by this boundary layer, remaining with a more mathematically treatable description away from it.

In this sense, the success of
the Kutta-Joukowski theory \index{Kutta-Joukowski theory}
consists in providing a simple, but not trivial, approach to aerodynamics which incorporates some aspects of viscous effects, while neglecting others.}
as follows:
\begin{equation}\label{JTTA}
\begin{split}&
{\mbox{an airfoil in relative motion with constant velocity~$-v_0$}}\\&
{\mbox{to an ambient inviscid homogeneous, irrotational fluid}}\\&{\mbox{has a lift force (that is the component of the force perpendicular to~$v_0$)}}\\&{\mbox{of magnitude }}
\rho \,|v_0| \Gamma, \\&{\mbox{where~$\Gamma$ is a circulation of~$v_0$ along the cross section of the airfoil.}}\end{split}\end{equation}
It is interesting to observe that~\eqref{JTTA} explains for instance the generation of a lift on a wing as a result
of the contribution of the circulation~$\Gamma$ of the velocity field around the wing.
Remarkably, the fact that~\eqref{JTTA} takes into account the component of the force perpendicular to~$v_0$
clearly states the importance of such a result for the theory of flight (that is, this lift force is what
allows, in principle
at least, an airplane to take off, provided they manage to create
a sufficiently large circulation~$\Gamma$).

\begin{figure}
  \centering
  \includegraphics[width=.6\linewidth]{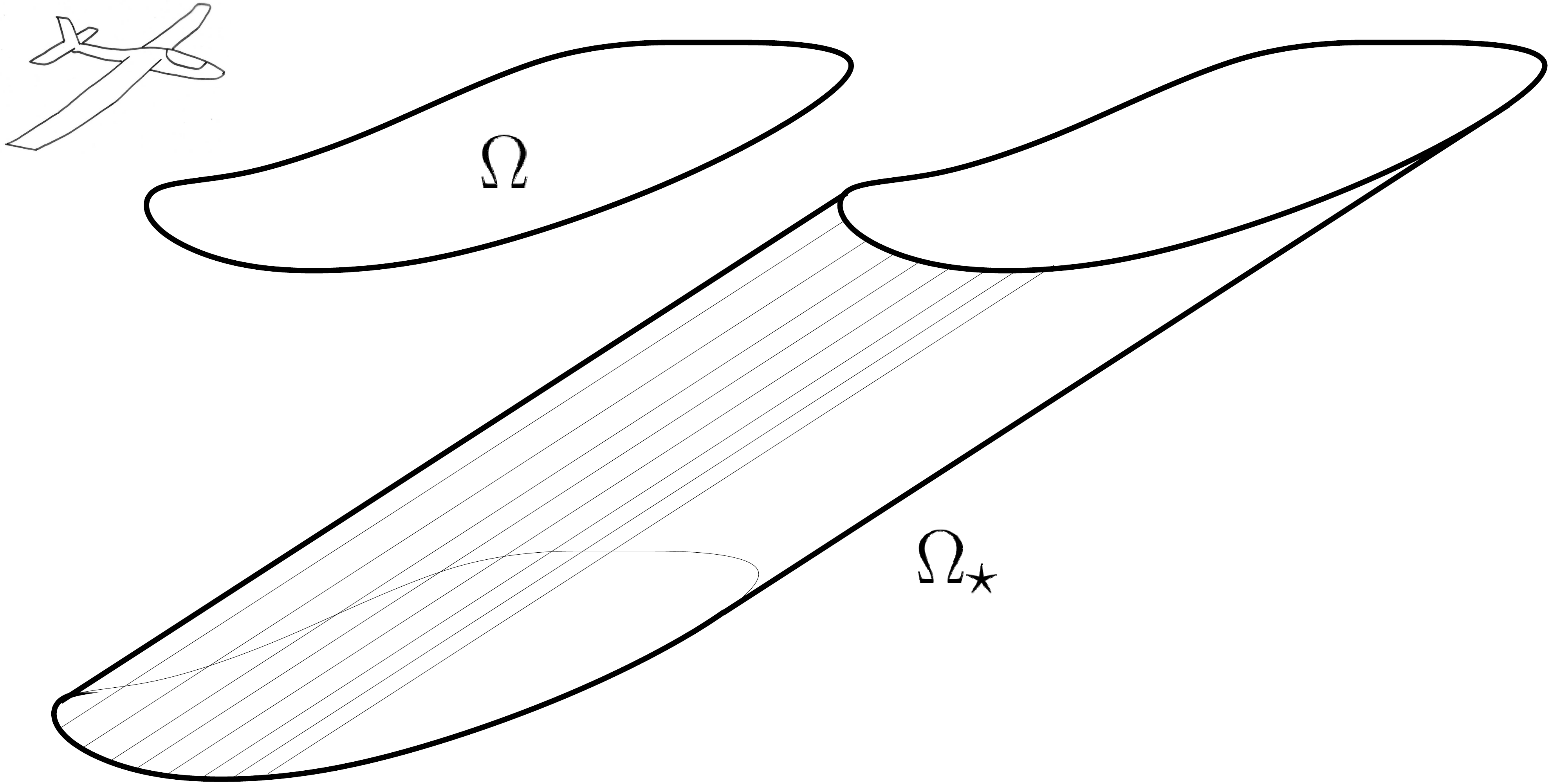}
 \caption{\sl Two-dimensional model of an airfoil, as inspired by a sailplane.}\label{P2DAILATENSIWASTREDItangeFIlASE}
\end{figure}

To model an airfoil, that is a ``long wing'' (actually, infinitely long for simplicity),
we consider a (nice, contractible) planar domain~$\Omega\subset\R^2$
and the wing given by~$\Omega_\star:=\Omega\times\R$, see Figure~\ref{P2DAILATENSIWASTREDItangeFIlASE}.
We reconsider the Euler equation in~\eqref{EUMSD-OS-D} (in the absence of gravity) and the incompressibility condition in~\eqref{INCOMPRE} outside~$\Omega_\star$
assuming that, by symmetry, the fluid parameters~$v$ and $p$
are actually independent of~$x_3$ and also on time (leading to a completely
steady solution of the problem). 
In this way, we may consider
the following set of equations for~$v=v(x,y)$ and~$p=p(x,y)$, with~$X=(x,y)\in\R^2\setminus\Omega$,
for a given constant~$\rho\in(0,+\infty)$
(the constancy of the density being the mathematical translation
that the fluid is homogeneous):
\begin{equation}\label{DFRE}
\begin{dcases}&\rho (v\cdot\nabla) v=-\nabla p,\\
&\div v=0,\\&\curl v=0.\end{dcases}\end{equation}

Notice that in this setting (thinking the vector field~$(v_1(x,y),v_2(x,y),0)$ to be three-dimensional
but with a trivial last entry to compute the curl), we have that
\begin{equation}\label{CUR25t34yuD2}
0=\curl v=\big(\partial_xv_2-\partial_y v_1\big) e_3.
\end{equation}
Moreover, since, according to the principle of inertia
discussed on page~\pageref{PRINERY},
we are taking an inertial frame in which the body appears to be still,
\begin{equation}\label{ALLIVO0}
\lim_{|(x,y)|\to+\infty}v(x,y)=v_0.\end{equation}
Also, the Euler equation in~\eqref{DFRE} and the assumption that the density is constant give that
\begin{eqnarray*}&& v\cdot\nabla\left(\frac{|v|^2}2+\frac{p}\rho\right)=
v\cdot\left(\sum_{k=1}^2 v_k\nabla v_k+\frac{\nabla p}\rho\right)=
v\cdot\left(\sum_{k=1}^2 v_k\nabla v_k-(v\cdot\nabla) v\right)
\\&&\qquad
=\sum_{j,k=1}^2 v_j v_k\partial_j v_k-
v\cdot\left(\sum_{j=1}^2 v_j\partial_j v\right)
=\sum_{j,k=1}^2 v_j v_k\partial_j v_k-
\sum_{j,k=1}^2 v_k v_j\partial_j v_k=0.
\end{eqnarray*}
As a consequence, if~$\beta(X):= \frac{|v(X)|^2}2+\frac{p(X)}\rho$,
$$ v\cdot \nabla\beta=0,$$
thus
\begin{equation*}\frac{d}{dt}\beta(X(t))=\nabla\beta(X(t))\cdot \dot X(t)=
\nabla\beta(X(t))\cdot v( X(t))=0.
\end{equation*}
As a result, $\beta(X(t))=\beta(X(0))$ for every time~$t$, that is
\begin{equation}\label{BERNEQl}
\frac{|v(X(t))|^2}{2}+\frac{p(X(t))}{\rho}=
\frac{|v(X(0))|^2}{2}+\frac{p(X(0))}{\rho}.
\end{equation}
Equations such as~\eqref{BERNEQl} are often referred to with the name of Bernoulli's Principle. \index{Bernoulli's Principle}

\begin{figure}
  \centering
  \includegraphics[width=.5\linewidth]{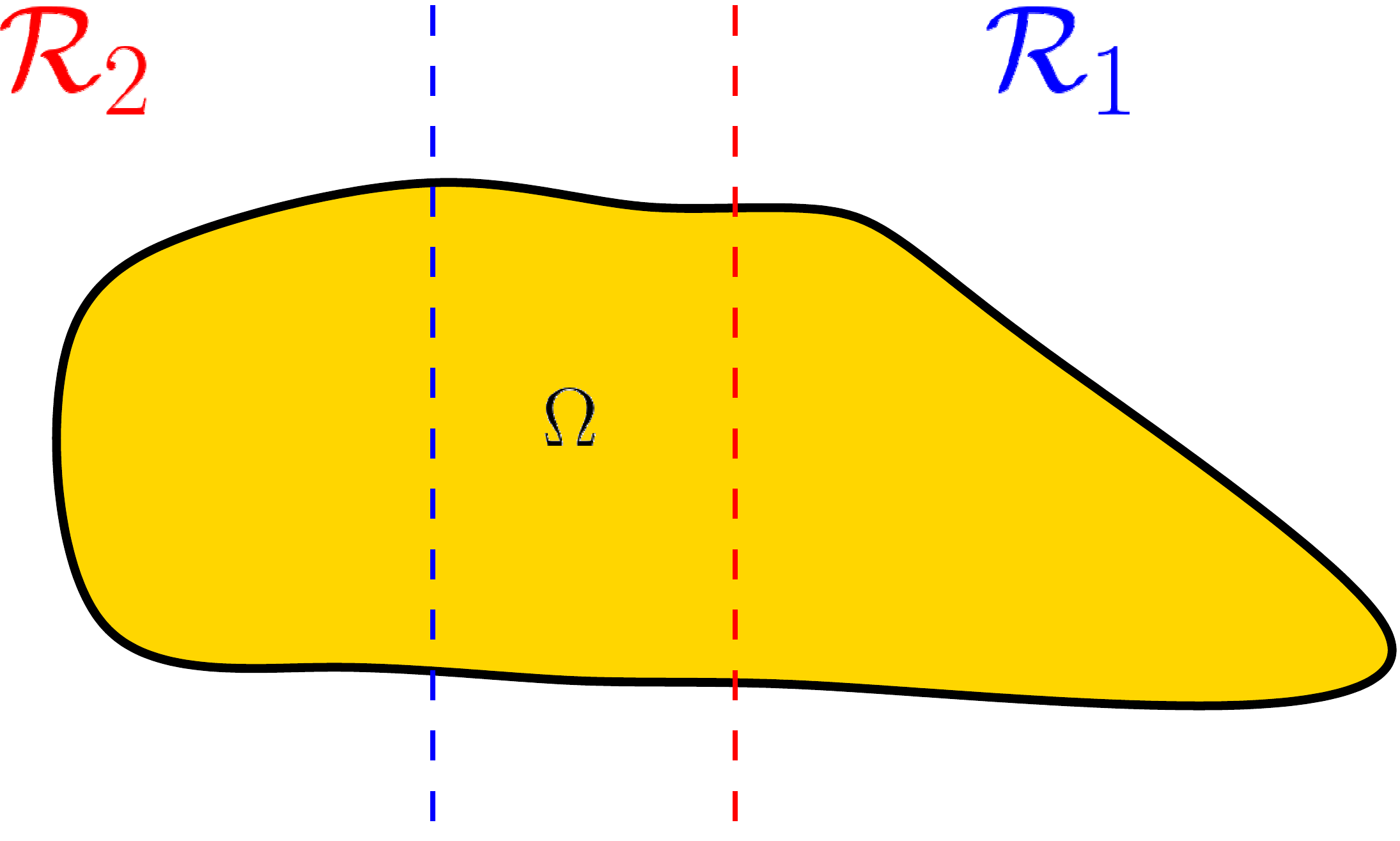}
 \caption{\sl A topological argument used to take care of the fact that~$\R^2\setminus\Omega$
 has a hole.}\label{P2DAILATENSIWASTREDItangeFIlAHOLEWFINASDJGIOSE}
\end{figure}

Now we rely on the two-dimensional structure of the problem combined with the divergence free condition
in~\eqref{DFRE} to see that the differential form~$\zeta:= v_2\,dx-v_1\,dy$ is closed,
since \begin{equation}\label{0pod-45ia7n7b9e0r}
d\zeta=-(\partial_y v_2+\partial_x v_1)\,dx\wedge dy=-\div v\,dx\wedge dy=0.\end{equation}
As a matter of fact, the exterior of~$\Omega$ is not simply connected in~$\R^2$, however we can consider two
regions~${\mathcal{R}}_1$ and~${\mathcal{R}}_2$ which are simply connected and their union is the exterior of~$\Omega$, see Figure~\ref{P2DAILATENSIWASTREDItangeFIlAHOLEWFINASDJGIOSE}.
In this setting, we can apply the Poincar\'e Lemma (see e.g.~\cite[Theorem~8.3.8]{MR1209437})
to each of the regions~${\mathcal{R}}_1$ and~${\mathcal{R}}_2$, finding two potentials~$\psi_j:{\mathcal{R}}_j\to\R$ such that~$d\psi_j=\zeta$ in~${\mathcal{R}}_j$, with~$j\in\{1,2\}$.
In particular,
\begin{equation}\label{56STRE78} \partial_x\psi_j=v_2\qquad{\mbox{and}}\qquad\partial_y\psi_j=-v_1.\end{equation}
Note that~${\mathcal{R}}_1\cap{\mathcal{R}}_2$ consists of two connected components which reach infinity,
say~${\mathcal{S}}_1$ and~${\mathcal{S}}_2$. We observe that
\[ \nabla(\psi_1-\psi_2)=
\nabla\psi_1-\nabla\psi_2=(v_2,-v_1)-(v_2,-v_1)=0,\]
in~${\mathcal{S}}_1$, hence there exists~$c_1\in\R$ such that~$\psi_1-\psi_2=c_1$ in~${\mathcal{S}}_1$.

Similarly, there exists~$c_2\in\R$ such that~$\psi_1-\psi_2=c_2$ in~${\mathcal{S}}_2$.
Up to incorporating the constant into~$\psi_2$ (which leaves~\eqref{56STRE78} invariant), we can suppose that~$c_2=0$.

We claim that
\begin{equation}\label{Cik231ma3223y0-2tkhojnhoiuwhg8732gtpiuywfng2uahHBSne}
c_1=0.
\end{equation}
For this, let~$A\in(\partial\Omega)\cap\overline{{\mathcal{S}}_1}$ and~$B\in(\partial\Omega)\cap\overline{{\mathcal{S}}_2}$
and let~$\vartheta:[0,\ell]\to(\partial\Omega)\cap\overline{{\mathcal{R}}_1}$ be a curve, parameterized by its arc length, lying in~$(\partial\Omega)\cap\overline{{\mathcal{R}}_1}$ and
joining~$A$ to~$B$, say with~$\vartheta(0)=A$ and~$\vartheta(\ell)=B$.
It follows from the impenetrability condition in~\eqref{PSLmoLMSxHNS-AR} 
that the velocity field is tangential to~$\partial\Omega$ and consequently
\begin{equation}\label{DEFSTTHiow9805e04938} v(\vartheta(\tau))
=\alpha(\tau)\,\dot\vartheta(\tau)\qquad{\mbox{ for all }}\tau\in[0,\ell],\end{equation}
for some scalar function~$\alpha(\tau)$.

Thus, since~$\vartheta$ lies in~$\overline{{\mathcal{R}}_1}$, we can utilize~\eqref{56STRE78} with~$j:=1$
and write that
$$ \Big( -\partial_y\psi_1\big(\vartheta(\tau)\big),\,\partial_x\psi_1\big(\vartheta(\tau)\big)\Big)=
\alpha(\tau)\,\dot\vartheta(\tau)$$
and therefore
$$ \nabla\psi_1\big(\vartheta(\tau)\big)=
\alpha(\tau)\,\big(\dot\vartheta_2(\tau),\,-\dot\vartheta_1(\tau)\big).$$
{F}rom this, we obtain that
$$ \frac{d}{d\tau}\psi_1\big(\vartheta(\tau)\big)=
\nabla\psi_1\big(\vartheta(\tau)\big)\cdot\dot\vartheta(\tau)
=
\alpha(\tau)\,\big(\dot\vartheta_2(\tau),\,-\dot\vartheta_1(\tau)\big)\cdot
\big(\dot\vartheta_1(\tau),\,\dot\vartheta_2(\tau)\big)=0
$$
and consequently
$$ \psi_1(B)=\psi_1\big(\vartheta(\ell)\big)=\psi_1\big(\vartheta(0)\big)=\psi_1(A).$$
Similarly, by considering a curve lying in~$(\partial\Omega)\cap\overline{{\mathcal{R}}_2}$, one obtains that
$$ \psi_2(B)=\psi_2(A).$$
As a result,
\begin{eqnarray*}
0=c_2=\psi_1(B)-\psi_2(B)=\psi_1(A)-\psi_2(A)=c_1
\end{eqnarray*}
and the proof of~\eqref{Cik231ma3223y0-2tkhojnhoiuwhg8732gtpiuywfng2uahHBSne} is complete.

Hence, we can now define, in the exterior of~$\Omega$,
$$ \psi:=
\begin{dcases}
\psi_1 & {\mbox{ in }}{\mathcal{R}}_1,\\
\psi_2 & {\mbox{ in }}{\mathcal{R}}_2\\
\end{dcases}$$
and we stress that this is a fair definition, since~$\psi_1=\psi_2$ in~${\mathcal{R}}_1\cap{\mathcal{R}}_2$,
owing to~\eqref{Cik231ma3223y0-2tkhojnhoiuwhg8732gtpiuywfng2uahHBSne}.

Furthermore, by~\eqref{56STRE78}, in the exterior of~$\Omega$,
\begin{equation}\label{9uojhexterimega}
\partial_x\psi =v_2\qquad{\mbox{and}}\qquad\partial_y\psi =-v_1.\end{equation}
Note in particular that~$|\nabla\psi|=v$. Moreover, $\psi$ is harmonic, thanks to the irrotationality
situation pointed out in~\eqref{CUR25t34yuD2}.

The function~$\psi$ is often referred to with the name of
``stream function''. 
This name comes from a physical intuition. Indeed, the streamlines, i.e. the trajectories of fluid particles, are described by the level sets of the stream function, because
\begin{equation}\label{STERTEAMSDN} \frac{d}{dt}\psi(X(t))=\nabla\psi(X(t))\cdot \dot X(t)=\nabla\psi(X(t))\cdot v(X(t))=0,\end{equation}
thanks to~\eqref{9uojhexterimega}.

Let us now reassess Bernoulli's Principle in~\eqref{BERNEQl}. \index{Bernoulli's Principle}
In our setting we need a refinement of~\eqref{BERNEQl} along the boundary of our moving object, namely we claim that, for every~$x\in\partial\Omega$,
\begin{equation}\label{BERNEQ}
|v(x)|^2+\frac{2p(x)}{\rho}=p_0,
\end{equation}
for some constant~$p_0\in\R$.

To prove this, we take the simplifying assumption that the velocity field possesses at most finitely many zeros along~$\partial\Omega$. We denote by~$Z_1,\dots,Z_N$ these zeros (if any).
We pick~$X_0\in(\partial\Omega)\setminus\{Z_1,\dots,Z_N\}$. Up to a rigid motion, we assume that~$X_0=0$
and that the tangent vector of~$\partial\Omega$ at~$0$ is horizontal.
Say, we describe~$\Omega$ near~$0$ as the subgraph of a smooth function~$f:\R\to\R$ with~$f(0)=0$ and~$f'(0)=0$.
We can also replace~$\psi$ by~$\psi-\psi(0)$ and thus assume additionally that~$\psi(0)=0$.
In this setting, from the impenetrability condition in~\eqref{PSLmoLMSxHNS-AR}, we have that the fluid vector field at the origin
is horizontal and thus, by~\eqref{9uojhexterimega}, we have that~$\nabla\psi(0)$ is vertical and different from zero.
Let us suppose that
\begin{equation}\label{deltastreamline}
\partial_y\psi(0)>0,\end{equation}
the case~$\partial_y\psi(0)<0$ being similar, just considering level sets
of~$\psi$ corresponding to negative, instead of positive, values.

In this framework, near~$0$, the level sets of~$\psi$ can be written as graph of smooth functions: more precisely,
near the origin, for small~$\delta>0$, the level set~$\{\psi=\delta\}$ can be identified with the graph~$\{y=f_\delta(x)\}$
for a suitable~$f_\delta$. It follows from~\eqref{deltastreamline} that~$f_{\delta'}\ge f_\delta$ if~$\delta'\ge\delta$.
Also, since, by~\eqref{STERTEAMSDN},~$\{\psi=\delta\}$ describes the streamlines of the fluid,
we have that~$\{\psi=\delta\}$ is contained in the complement of~$\Omega$ for all~$\delta>0$ and consequently,
near the origin,~$\partial\Omega$ lies below the graph of~$f_\delta$ for all~$\delta>0$.

\begin{figure}
  \centering
  \includegraphics[width=.5\linewidth]{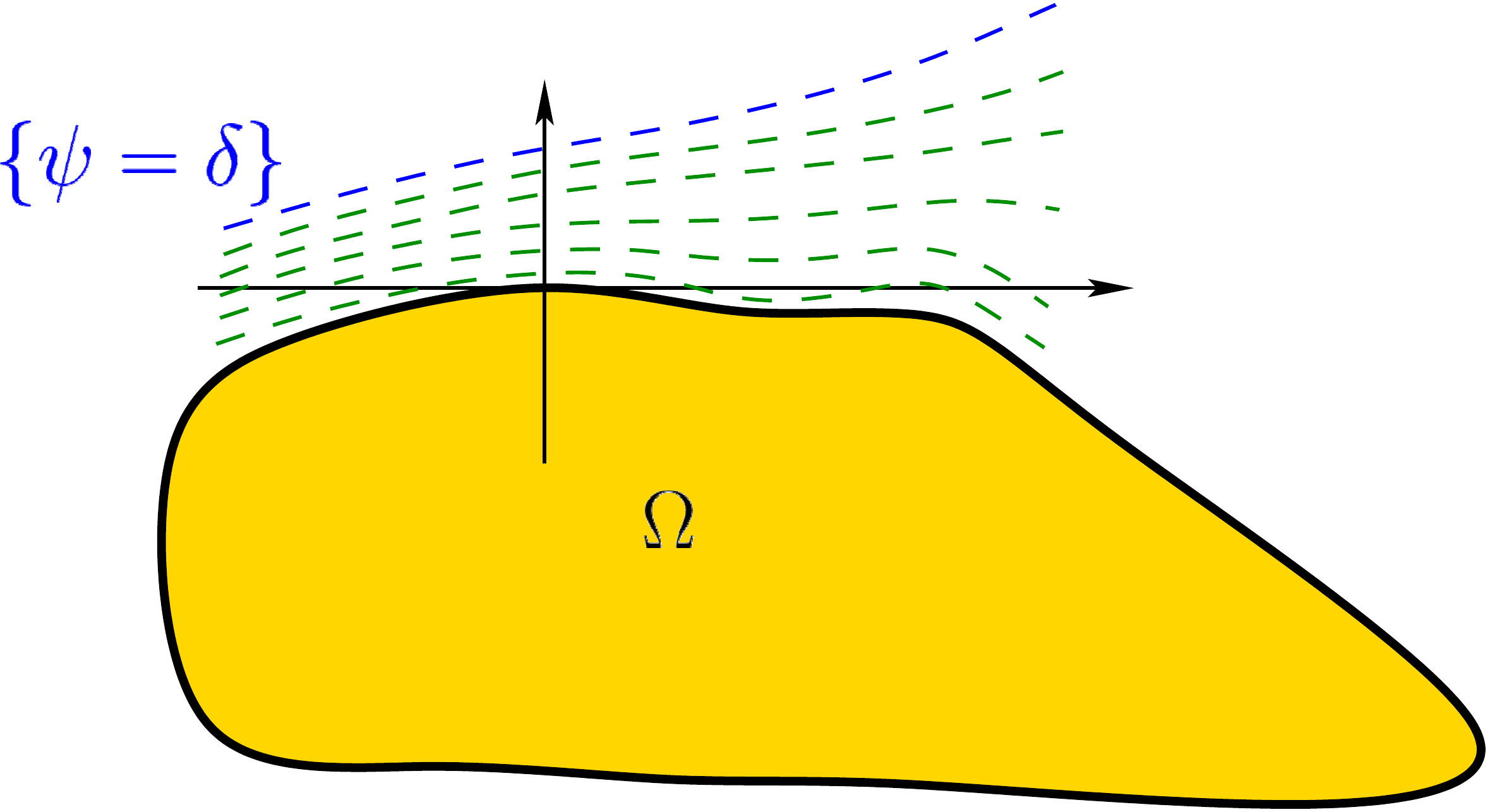}
 \caption{\sl Foliation of streamlines in the vicinity of an airfoil.}\label{PTOP12DAILATENSIWASTREDItangeFIlAHOLEWFINASDJGIOSE}
\end{figure}

Thus, given the monotonicity of~$f_\delta$, near the origin we can define
$$ f_0(x):=\lim_{\delta\searrow0}f_\delta(x)$$
and we have that~$f\le f_0$. Moreover, since
$$ \psi(0,f_0(0))=\lim_{\delta\searrow0}\psi(0,f_\delta(0))=\lim_{\delta\searrow0}\delta=0=\psi(0,0)=\psi(0,f(0)),$$
we deduce from~\eqref{deltastreamline} that~$f_0(0)=f(0)$.
Accordingly, $f_0=f$ (otherwise different fluid trajectories would meet at the origin with the same velocity,
in contradition with the uniqueness results for ordinary differential equations), see Figure~\ref{PTOP12DAILATENSIWASTREDItangeFIlAHOLEWFINASDJGIOSE}.

This shows that the streamline emanating from the origin (that is, from every point of~$\partial\Omega$ with nonzero velocity field)
remains on~$\partial\Omega$.

Thus, if the velocity field does not vanish on~$\partial\Omega$, then~$\partial\Omega$ consists of a streamline
and then~\eqref{BERNEQ} follows from~\eqref{BERNEQl}. 

\begin{figure}[b]
  \centering
  \includegraphics[width=.5\linewidth]{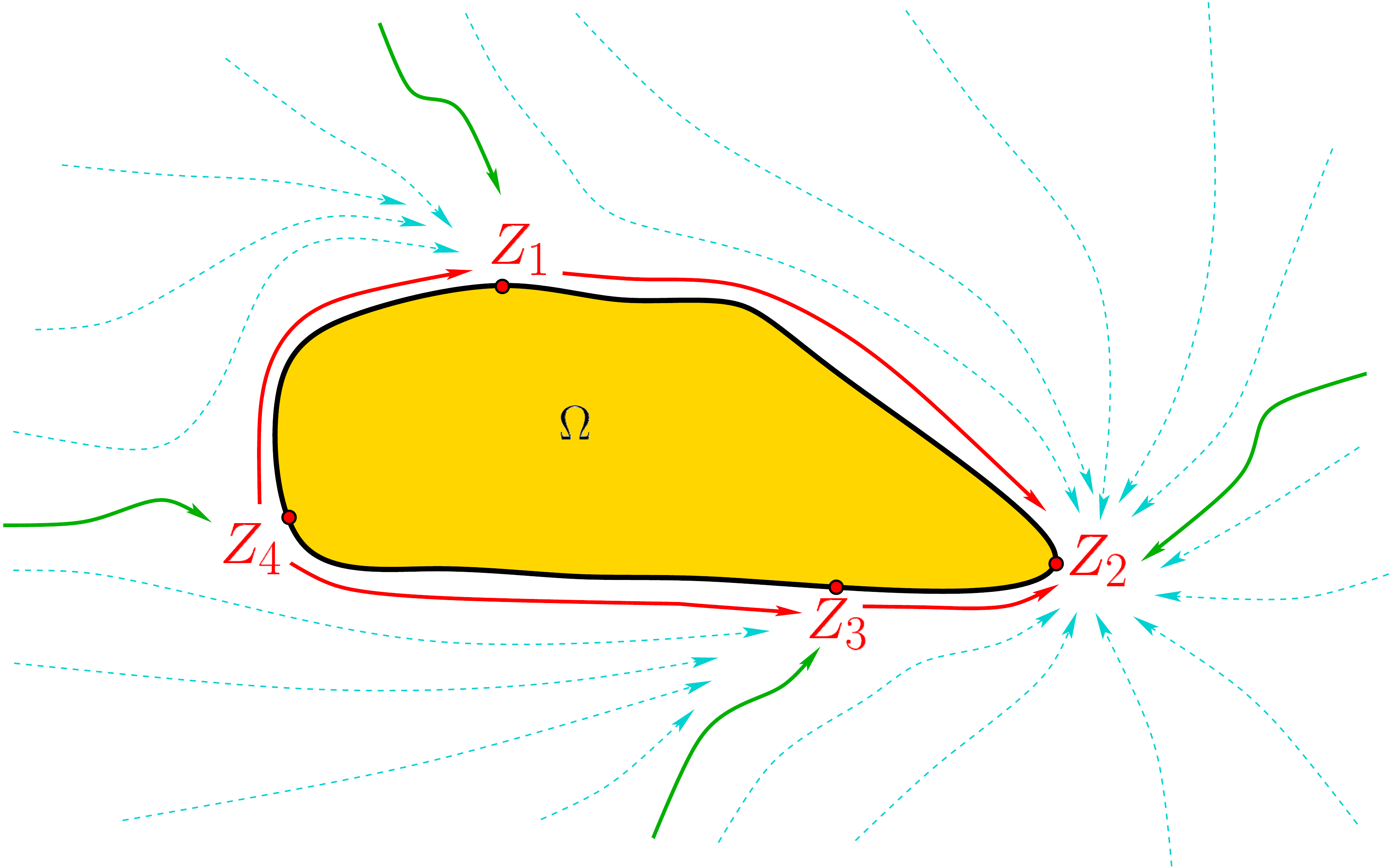}
 \caption{\sl Partitioning~$\partial\Omega$ as in~\eqref{PTOP12DAILATENSIWASTbni9buj9unooREDItangeFIlAHOLEWFINASDJGIOSE2}
 (not necessarily a realistic picture).}\label{PTOP12DAILATENSIWASTbni9buj9unooREDItangeFIlAHOLEWFINASDJGIOSE}
\end{figure}

If instead
the velocity field vanishes at~$Z_1,\dots,Z_N\in\partial\Omega$, we can decompose~$\partial\Omega$
into open arcs~$\Lambda_1,\dots,\Lambda_N$ such that
\begin{equation}\label{PTOP12DAILATENSIWASTbni9buj9unooREDItangeFIlAHOLEWFINASDJGIOSE2}\partial\Omega=\Lambda_1\cup\dots\cup\Lambda_N\cup\{Z_1,\dots,Z_N\},\end{equation}
in such a way that the velocity field does not vanish on each~$\Lambda_j$ and~$\Lambda_j$
is a streamline connecting~$Z_j$ to~$Z_{j+1}$ in infinite time (with~$Z_{N+1}:=Z_1$), see Figure~\ref{PTOP12DAILATENSIWASTbni9buj9unooREDItangeFIlAHOLEWFINASDJGIOSE}.
In this setting, we can consider a trajectory~$X(t)$ starting at a given point of~$\Lambda_j$ and deduce from~\eqref{BERNEQl} that,
for every~$x\in\Lambda_j$,
$$ |v(x)|^2+\frac{2p(x)}{\rho}=p_j,$$
for some~$p_j\in\R$, and
$$ p_j=\lim_{t\to\pm\infty}|v(X(t))|^2+\frac{2p(X(t))}{\rho}$$
which leads to
$$ p_j=|v(Z_j)|^2+\frac{2p(Z_j)}{\rho}=|v(Z_{j+1})|^2+\frac{2p(Z_{j+1})}{\rho}.$$
This actually says that~$p_1=\dots=p_N$ and thus~\eqref{BERNEQ} plainly follows, as desired.

Now we describe~$\partial\Omega$ as a smooth and closed curve, traveled counterclockwise and parameterized by its arc length, say~$\gamma:[0,L]\to\cOMPL$, for some~$L>0$ (in particular, notice that~$\gamma(0)=\gamma(L)$).
We  identify points of~$\R^2$ with complex numbers.
In this way, given~$\tau\in[0,L]$, a unit tangent vector to~$\partial\Omega$ at~$\gamma(\tau)$
is given by~$\dot\gamma(\tau)$ and the unit normal vector~$\nu(\gamma(\tau))$ pointing outwards is given by a clockwise rotation
of~$\dot\gamma(\tau)$ by an angle~$\frac\pi2$, corresponding to the complex multiplication by~$e^{-i\pi/2}=-i$
(that is, the outward unit normal can be written as~$-i\dot\gamma(\tau)$).

For this reason, in light of~\eqref{FEUYT90hfbvjjxJ}, we have that
the force acting on the airfoil (after identifying~$\R^2$ with~$\cOMPL$) takes the form
\begin{equation*} F=-\int_{\partial\Omega} p\nu
=-\int_0^L p(\gamma(\tau))\,\nu(\gamma(\tau))\,d\tau=i\int_0^L p(\gamma(\tau))\,\dot\gamma(\tau)\,d\tau.
\end{equation*}
{F}rom this and~\eqref{BERNEQ} we arrive at
\begin{equation}\label{BERNEQ2} \frac{2F}\rho
=i\int_0^L \frac{2p(\gamma(\tau))}\rho\,\dot\gamma(\tau)\,d\tau
=i\int_0^L 
\left(p_0-|\nabla\psi(\gamma(\tau))|^2\right)
\,\dot\gamma(\tau)\,d\tau.
\end{equation}
Since
$$ \int_0^L \dot\gamma(\tau)\,d\tau=\gamma(L)-\gamma(0)=0,$$
we deduce from~\eqref{BERNEQ2} that
\begin{equation}\label{BERNEQ3} \frac{2F}\rho
=-i\int_0^L |\nabla\psi(\gamma(\tau))|^2\,\dot\gamma(\tau)\,d\tau.
\end{equation}
We also recall that, in light of the impenetrability condition in~\eqref{PSLmoLMSxHNS-AR}
combined with~\eqref{9uojhexterimega},
we have that~$\nabla\psi$ is orthogonal to
the unit tangent vector~$\dot\gamma$.
Thus, we write that
\begin{equation}\label{BERNEQ4} \nabla\psi(\gamma(\tau))=\pm i
|\nabla\psi(\gamma(\tau))|\,\dot\gamma(\tau).\end{equation}
Also, in the complex variable setting, that is using the notation~$z=x+iy$ and identifying~$x+iy$ with~$(x,y)$, we set
\begin{equation}\label{ENGBSC2O232M4P3AND5D} w(z)=w(x+iy):=\partial_x\psi(x,y)-i\partial_y\psi(x,y),\end{equation}
and it follows that~$w$ is the complex conjugated of~$\nabla\psi$ and~$|w|=|\nabla\psi|$.

Hence, taking the complex conjugation in~\eqref{BERNEQ4},
\begin{equation*} w(\gamma(\tau))=
\pm \overline{i|w(\gamma(\tau))|\,\dot\gamma(\tau)}=
\mp i|w(\gamma(\tau))|\,\overline{\dot\gamma(\tau)}.\end{equation*}
Taking the square of this identity
\begin{equation*} w^2(\gamma(\tau))=-|w(\gamma(\tau))|^2\,(\overline{\dot\gamma(\tau)})^2
\end{equation*}
and as a result
\begin{equation*} \begin{split}&w^2(\gamma(\tau)) \dot\gamma(\tau)=-|w(\gamma(\tau))|^2\,(\overline{\dot\gamma(\tau)})^2{\dot\gamma(\tau)}
=-|w(\gamma(\tau))|^2\,\overline{\dot\gamma(\tau)}\,|\dot\gamma(\tau)|^2\\&\quad=-|w(\gamma(\tau))|^2\,\overline{\dot\gamma(\tau)}=-
|\nabla\psi(\gamma(\tau))|^2\,\overline{\dot\gamma(\tau)}.\end{split}
\end{equation*}
We thus combine this information with the complex conjugation of~\eqref{BERNEQ3} and we conclude that
\begin{equation}\label{wdotgammatau}
\begin{split}&\overline{\frac{2i F}\rho}
=\overline{\int_0^L |\nabla\psi(\gamma(\tau))|^2\,\dot\gamma(\tau)\,d\tau}=
\int_0^L |\nabla\psi(\gamma(\tau))|^2\,\overline{\dot\gamma(\tau)}\,d\tau\\&\quad\qquad\qquad\qquad=-\int_0^L w^2(\gamma(\tau)) \dot\gamma(\tau)\,d\tau=-
\oint_\gamma w^2(z)\,dz.\end{split}
\end{equation}
Also, using the harmonicity of~$\psi$, we observe that
$$ \partial_x(\Re w)=\partial_{xx}\psi=-\partial_{yy}\psi=\partial_y(\Im w)$$
and moreover
$$ \partial_y(\Re w)=\partial_{xy}\psi(x,y)=-\partial_x(\Im w).$$
These observations give that~$w$ satisfies the Cauchy-Riemann equations
and therefore~$w$ is holomorphic outside~$\Omega$.

We can therefore take~$R>0$ large enough such that~$\Omega\Subset B_R$
and employ Laurent's Theorem (see e.g.~\cite[Theorem~11.1]{MR3839273}, applied here with~$R_1:=R$, $R_2:=+\infty$
and~$z_0:=0$): in this way we
obtain that, for all~$z$ in the exterior of~$B_R$, 
the following series representation holds
true for suitable~$a_k$, $b_k\in \cOMPL$:
$$ w(z)=w_0(z)+w_\star(z)\qquad{\mbox{with}}\qquad
w_0(z):=\sum_{k=0}^{+\infty} a_k z^k\qquad{\mbox{and}}\qquad
w_\star(z):=\sum_{k=1}^{+\infty}\frac{b_k}{ z^k},$$
where the power series representing~$w_0$ converges in~$\cOMPL$ and the one representing~$w_\star$
converges in the exterior of~$B_R$.

In particular, $w_0$ is holomorphic in the whole of~$\cOMPL$
and the convergence of the power series defining~$w_\star$ entails (see e.g.~\cite[Theorems~3.19 and~3.23]{MR3839273}) that
$$ \limsup_{k\to+\infty} |b_k|^{1/k}\le{2R}.$$
As a result, we can take~$k_0\in\N\cap[2,+\infty)$ sufficiently large such that~$|b_k|\le(3R)^k$ for every~$k> k_0$
and therefore
\begin{eqnarray*}&& \limsup_{|z|\to+\infty}|w_\star(z)|\le
\limsup_{|z|\to+\infty}\sum_{k=1}^{k_0}\frac{|b_k|}{| z|^k}+\limsup_{|z|\to+\infty}\sum_{k=k_0+1}^{+\infty}\frac{|b_k|}{| z|^k}
\le0+\limsup_{|z|\to+\infty}\sum_{k=k_0+1}^{+\infty}\left(\frac{3R}{| z|}\right)^k\\&&\qquad\qquad\qquad\qquad
=\limsup_{|z|\to+\infty}
\left(\frac{3R}{| z|}\right)^{k_0+1}\frac{1}{1-\displaystyle\frac{3R}{| z|}}=0.
\end{eqnarray*}
For this reason, recalling~\eqref{ALLIVO0},
$$ |v_0|=\limsup_{|(x,y)|\to+\infty}|\nabla\psi (x,y)|=\limsup_{|z|\to+\infty}|w(z)|=
\limsup_{|z|\to+\infty}|w_0(z)+w_\star(z)|=\limsup_{|z|\to+\infty}|w_0(z)|.$$
As a consequence, $w_0$ is bounded in the whole of~$\cOMPL$ and therefore, by Liouville's Theorem
(see e.g.~\cite[Theorem~10.8]{MR3839273}) we have that~$w_0$ is constant.

That is, by~\eqref{ALLIVO0}, $w_0(z)=i \overline{v_0}$ for all~$z\in\cOMPL$, whence, for all~$z$ in the exterior of~$B_R$,
\begin{equation}\label{EXTE6B-pjrmfR} w(z)=i \overline{v_0}+w_\star(z)=
i \overline{v_0}+\sum_{k=1}^{+\infty}\frac{b_k}{ z^k}.\end{equation}
%%Now we define
%%$$ W(z):=v_0 z+\sum_{k=1}^{+\infty}\frac{b_k}{(1-k) z^{k-1}}$$
%%and we point out that, for every~$z$ in the exterior of~$\Omega$,
%%$$ W'(z)=w(z).$$
Moreover, \begin{eqnarray*}&&\limsup_{M\to+\infty}\left|\oint_{\partial B_M}\sum_{k=k_0+1}^{+\infty}\frac{b_k}{ z^k}\,dz\right|\le
\limsup_{M\to+\infty}\sum_{k=k_0+1}^{+\infty}\frac{2\pi|b_k|}{M^{k-2}}\\&&\quad
\le\limsup_{M\to+\infty}\sum_{k=k_0+1}^{+\infty}\frac{2\pi M^2(3R)^k}{M^{k}}
=\limsup_{M\to+\infty}
2\pi M^2\left(\frac{3R}{M}\right)^{k_0+1}\frac{1}{1-\displaystyle\frac{3R}{M}}=0.\end{eqnarray*}
Thus, by~\eqref{EXTE6B-pjrmfR} and Cauchy's Theorem (see e.g.~\cite[Theorem~9.13]{MR3839273}),
\begin{equation}\label{SunmgeDediamJofeMSy0r987wqugf7wegbDEDS}\begin{split}&
\oint_\gamma w(z)\,dz=
\lim_{M\to+\infty}\oint_{\partial B_M} w(z)\,dz\\&\quad=
\lim_{M\to+\infty}\oint_{\partial B_M} \left(i \overline{v_0}+\sum_{k=1}^{+\infty}\frac{b_k}{ z^k}\right)\,dz=
\lim_{M\to+\infty}\oint_{\partial B_M} \left(i \overline{v_0}+\sum_{k=1}^{k_0}\frac{b_k}{ z^k}\right)\,dz.\end{split}
\end{equation}
We stress that the function in the latter integrand is meromorphic in~$B_M$ and we can thereby employ
Cauchy's Residue Theorem (see e.g.~\cite[Theorem~12.3]{MR3839273}), thus  deducing from~\eqref{SunmgeDediamJofeMSy0r987wqugf7wegbDEDS} that
\begin{equation}\label{SunmgeDediamJofeMSy0r987wqugf7wegbDEDS-2}
\oint_\gamma w(z)\,dz=
\lim_{M\to+\infty}\oint_{\partial B_M} w(z)\,dz=2\pi i b_1.
\end{equation}
In addition,
\begin{equation*}
\begin{split}&
\oint_\gamma w(z)\,dz=\int_0^L
\Big(\partial_x\psi(\gamma(\tau))-i\partial_y\psi(\gamma(\tau))\Big)
\Big( \dot\gamma_1(\tau)+i\dot\gamma_2(\tau)\Big)\,d\tau\\&\quad=
\int_0^L
\Big(\partial_x\psi(\gamma(\tau)) \dot\gamma_1(\tau) +\partial_y\psi(\gamma(\tau))\dot\gamma_2(\tau)-i\partial_y\psi(\gamma(\tau))\dot\gamma_1(\tau)
+i\partial_x\psi(\gamma(\tau))\dot\gamma_2(\tau)
\Big)\,d\tau\\&\quad=
\int_0^L
\left(  \frac{d}{d\tau}\psi(\gamma(\tau)) -i\partial_y\psi(\gamma(\tau))\dot\gamma_1(\tau)
+i\partial_x\psi(\gamma(\tau))\dot\gamma_2(\tau)
\right)\,d\tau\\&\quad=\psi(\gamma(L))-\psi(\gamma(0))-i
\int_0^L
\left(  \partial_y\psi(\gamma(\tau))\dot\gamma_1(\tau)-\partial_x\psi(\gamma(\tau))\dot\gamma_2(\tau)
\right)\,d\tau\\&\quad=
-i
\int_0^L
\left(  \partial_y\psi(\gamma(\tau))\dot\gamma_1(\tau)-\partial_x\psi(\gamma(\tau))\dot\gamma_2(\tau)
\right)\,d\tau.
\end{split}\end{equation*}
Therefore, recalling the definition of~$\Gamma$ in~\eqref{JTTA} and the relations in~\eqref{9uojhexterimega},
\begin{equation*}
\begin{split}&
\oint_\gamma w(z)\,dz=i
\int_0^L
\left( v_1(\gamma(\tau))\dot\gamma_1(\tau)+v_2(\gamma(\tau))\dot\gamma_2(\tau)
\right)\,d\tau=i\Gamma.
\end{split}\end{equation*}
As a result of this and~\eqref{SunmgeDediamJofeMSy0r987wqugf7wegbDEDS-2} we deduce that
$$ b_1=\frac{\Gamma}{2\pi }.$$
This and~\eqref{EXTE6B-pjrmfR} lead to
$$ w^2(z)=\left(
i \overline{v_0}+\frac{\Gamma}{2\pi z}+\sum_{k=2}^{+\infty}\frac{b_k}{ z^k}
\right)^2=-\overline{v_0}^2+\frac{i \overline{v_0} \Gamma}{\pi z}+\sum_{k=2}^{+\infty}\frac{c_k}{ z^k},$$
for suitable~$c_k\in\cOMPL$.

This, Cauchy's Residue Theorem (see e.g.~\cite[Theorem~12.3]{MR3839273}) and~\eqref{wdotgammatau} give that
\begin{equation*}
\overline{\frac{2i F}\rho}
=2  \overline{v_0} \Gamma
\end{equation*}
and therefore
\begin{equation*}
{\frac{2i F}\rho}
=\overline{2\overline{v_0} \Gamma}=2v_0 \Gamma.
\end{equation*}
Consequently,
$$ F=-i\rho v_0 \Gamma$$
and this establishes the claim in~\eqref{JTTA}.\medskip

\begin{figure}
  \centering
  \includegraphics[width=.35\linewidth]{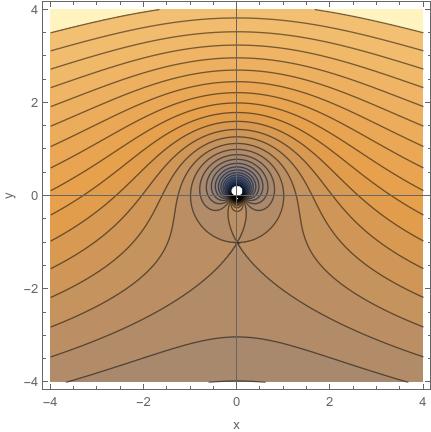}
 \caption{\sl The level sets of the potential function in~\eqref{HARGIRA4AXELHARLROA7789GIJ7soloDItangeFI2}.}\label{HARGIRA4AXELHARLROA7789GIJ7soloDItangeFI}
\end{figure}

\begin{figure}
  \centering
  \includegraphics[width=.35\linewidth]{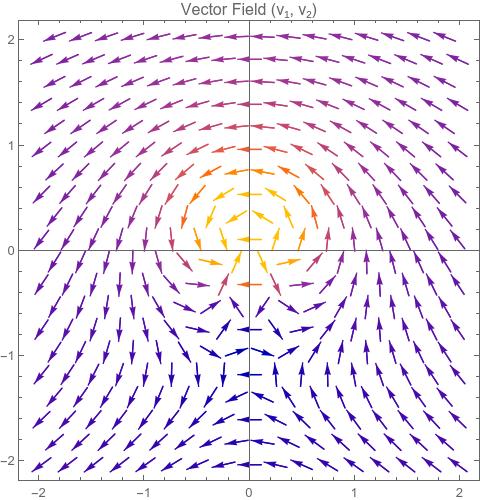}
 \caption{\sl The vector field in~\eqref{HARGIRA4AXELHARLROAGIJ7soloDItangeFI2}.}\label{HARGIRA4AXELHARLROAGIJ7soloDItangeFI}
\end{figure}

A natural question now is however if it is possible to produce a nonzero circulation in~\eqref{JTTA}.
This is in general quite an intriguing problem: here we just provide a simple and explicit\footnote{Though
we do not really exploit this fact here, the example is inspired by a rotating cylinder producing
a Magnus effect, that is what soccer players use to make the ball curve during flight. See {\tt https://www.youtube.com/watch?v=XdL7EDKr\_rk} for a famous application of 
the Magnus effect in soccer.} example.
For this, we let \begin{equation}\label{HARGIRA4AXELHARLROA7789GIJ7soloDItangeFI2A}\Omega:=B_1\end{equation} and
\begin{equation}\label{HARGIRA4AXELHARLROA7789GIJ7soloDItangeFI2} \psi(x,y):= y\left(1-\frac1{x^2+y^2}\right)+\ln(x^2+y^2).\end{equation}
See Figure~\ref{HARGIRA4AXELHARLROA7789GIJ7soloDItangeFI} for the level sets of~$\psi$.
These level sets will correspond to streamlines, since we will consider the vector field~$v=(v_1,v_2)$ with
\begin{equation}\label{HARGIRA4AXELHARLROAGIJ7soloDItangeFI2}
\begin{split}
&v_1(x,y):=-\partial_y\psi(x,y)=
\frac{2 x^2}{(x^2 + y^2)^2 }-\frac{1 + 2 y}{x^2 + y^2}-1
\\ {\mbox{and}}\qquad&
v_2(x,y):=\partial_x\psi(x,y)=\frac{2 x}{x^2 + y^2} +\frac{2 x y}{(x^2 + y^2)^2}.\end{split}
\end{equation}
See Figure~\ref{HARGIRA4AXELHARLROAGIJ7soloDItangeFI} for a sketch
of the vector field~$v$.

\begin{figure}[b]
  \centering
  \includegraphics[width=.45\linewidth]{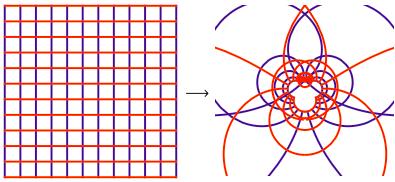}
 \caption{\sl The complex map~$w$ in equation~\eqref{COMPLEXPOTE}.}\label{HARGIRA4AXELHARLROAGIJ7soloDItangeFINOTa345}
\end{figure}

The complex
analysis enthusiasts will also enjoy the fact that the setting in~\eqref{ENGBSC2O232M4P3AND5D}
is also available by posing \label{HARGIRA4AXELHARLROAGIJ7soloDItangeFINOTa345N}
\begin{equation}\label{COMPLEXPOTE}w:=\frac2z + i \left(\frac1{z^2} - 1\right).\end{equation}
Notice that~$w$ is holomorphic away of the origin. And the knowledge of this complex
map is pretty much all that is needed to produce the fluid velocity field in~\eqref{HARGIRA4AXELHARLROAGIJ7soloDItangeFI2}.
The action of the above map~$w$ is sketched in Figure~\ref{HARGIRA4AXELHARLROAGIJ7soloDItangeFINOTa345}.

We observe that~$v$ is tangential to~$\partial\Omega$ since
$$ v(\cos\theta,\sin\theta)\cdot(\cos\theta,\sin\theta)=0$$
and this is consistent with the impenetrability condition in~\eqref{PSLmoLMSxHNS-AR}.
Additionally, we have that~$\Delta\psi=0$ outside~$\Omega$, which leads that~$\div v=0$
outside~$\Omega$ as well.

Furthermore, recalling~\eqref{CUR25t34yuD2},
\begin{equation}\label{8ihTGBAlFer3ONmdfZ345tyIjdmvREBSLv02}
\curl v=\big(\partial_xv_2-\partial_y v_1\big) e_3=\big(\partial_x(\partial_x\psi)+\partial_y( \partial_y\psi )\big) e_3=0
\end{equation}
and
\begin{eqnarray*}
\lim_{|(x,y)|\to+\infty} v(x,y)=
\lim_{|(x,y)|\to+\infty}\left(
\frac{2 x^2}{(x^2 + y^2)^2 }-\frac{1 + 2 y}{x^2 + y^2}-1,
\;\frac{2 x}{x^2 + y^2} +\frac{2 x y}{(x^2 + y^2)^2}\right)
=(-1,0)=:v_0.\end{eqnarray*}

Thus, to check that this example is indeed a solution of our fluid dynamics problem, it remains
to show that the first equation in~\eqref{DFRE} holds true.
For this, we choose~$\rho:=1$ and, inspired by the Bernoulli's Principle \index{Bernoulli's Principle}
in~\eqref{BERNEQl},
$$ p(x,y):=\frac{|v_0|^2-|v(x,y)|^2}{2}=\frac{1-|v(x,y)|^2}{2}$$
and we combine the vectorial identity in~\eqref{VCISTPOK} with~\eqref{8ihTGBAlFer3ONmdfZ345tyIjdmvREBSLv02}
to calculate that \begin{eqnarray*}&&\rho (v\cdot\nabla) v+\nabla p=
(v\cdot\nabla) v-\nabla\frac{|v|^2}{2}=-
v\times (\curl v)
=0.
\end{eqnarray*}
This shows that the vector field in~\eqref{HARGIRA4AXELHARLROAGIJ7soloDItangeFI2}
is consistent with the fluid dynamics setting in~\eqref{JTTA}. To show the interest of~\eqref{JTTA}
it thus remains to check that this vector field provides a nontrivial circulation~$\Gamma$: to this end, we compute that
\begin{eqnarray*}
\Gamma&=&\int_0^{2\pi} v(\cos\theta,\sin\theta)\cdot(-\sin\theta,\cos\theta)\,d\theta\\
&=&2\int_0^{2\pi} (\cos^2\theta-\sin\theta-1,\;\cos\theta +\sin\theta\cos\theta)\cdot(-\sin\theta,\cos\theta)\,d\theta\\&=&2\int_0^{2\pi} (1+\sin\theta)\,d\theta\\&=&4\pi
\\&\ne&0,
\end{eqnarray*}
as desired.\medskip

\begin{figure}
  \centering
  \includegraphics[width=.45\linewidth]{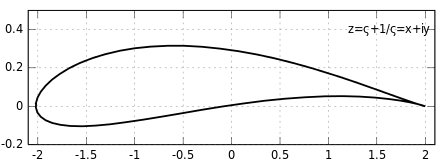}
 \caption{\sl Joukowski airfoil
 (image by Krishnavedala from
 Wikipedia, licensed under the Creative Commons Attribution-Share Alike 4.0 International license).}\label{123GR12edHARGIRA4AXELHARLROA7789GIJ7soloDItangeFI1}
\end{figure}

But hold on a sec, the skeptical reader (who is always very welcome) will complain we've been cheating on them:
we have promised this section was devoted to airfoils, but then we tried to sell in~\eqref{HARGIRA4AXELHARLROA7789GIJ7soloDItangeFI2A} the disk as an example of airfoils.
Come on, airfoils are objects as the ones in Figure~\ref{123GR12edHARGIRA4AXELHARLROA7789GIJ7soloDItangeFI1}, nothing in Figures~\ref{HARGIRA4AXELHARLROA7789GIJ7soloDItangeFI}, \ref{HARGIRA4AXELHARLROAGIJ7soloDItangeFI} or~\ref{HARGIRA4AXELHARLROAGIJ7soloDItangeFINOTa345}
looks like an airfoil, just because the disk is not an airfoil at all!

Well, in fact, to a certain degree, it is: this is the beauty of mathematics, and of complex analysis in particular.
Indeed, one of the brilliant ideas of the pioneers of aerodynamics is that the potential in~\eqref{COMPLEXPOTE}
describes all the fluid dynamics outside the disk of the example that we have discussed explicitly in detail,
therefore other examples can be constructed via conformal transformations.
In particular, using dilations, translations and maps of the form~$z+\frac1z$, one can transform a disk
into an airfoil shaped as in Figure~\ref{123GR12edHARGIRA4AXELHARLROA7789GIJ7soloDItangeFI1}.

The set of conformal transformations linking circles and airfoils can be visualized for instance via the Wolfram Demonstrations Projects

\noindent
{\tt https://demonstrations.wolfram.com/TheJoukowskiMappingAirfoilsFromCircles/}

\noindent {\tt https://demonstrations.wolfram.com/JoukowskiAirfoilGeometry/}

\noindent or via the GeoGebra application

\noindent {\tt https://www.geogebra.org/m/XwmqSR49}

\noindent or likely via a number of IT resources.

\begin{figure}[t]
  \centering
  \includegraphics[width=.3\linewidth]{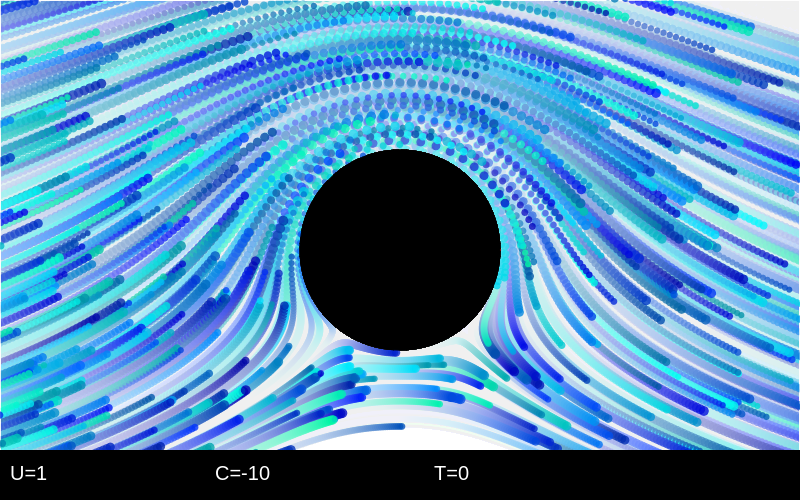}$\qquad$
  \includegraphics[width=.3\linewidth]{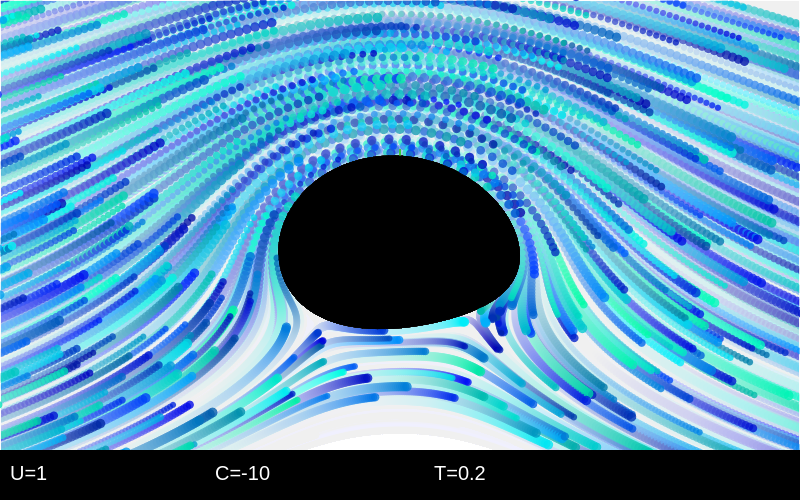}\\ $\qquad$
  \\
  \includegraphics[width=.3\linewidth]{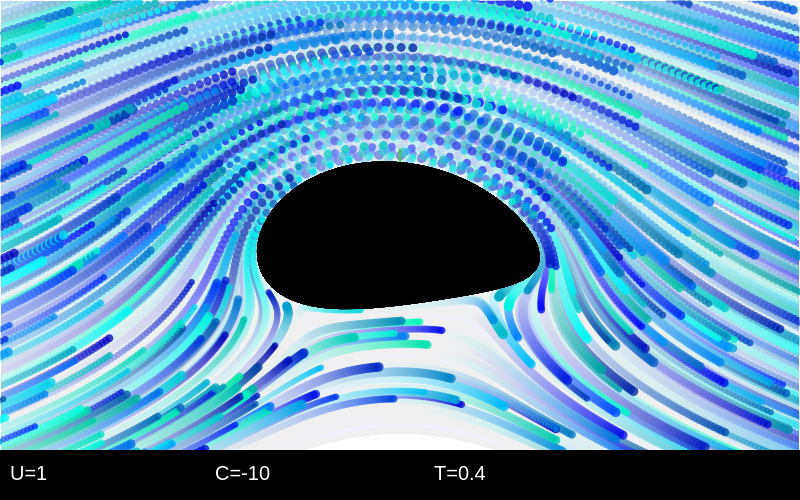}$\qquad$
  \includegraphics[width=.3\linewidth]{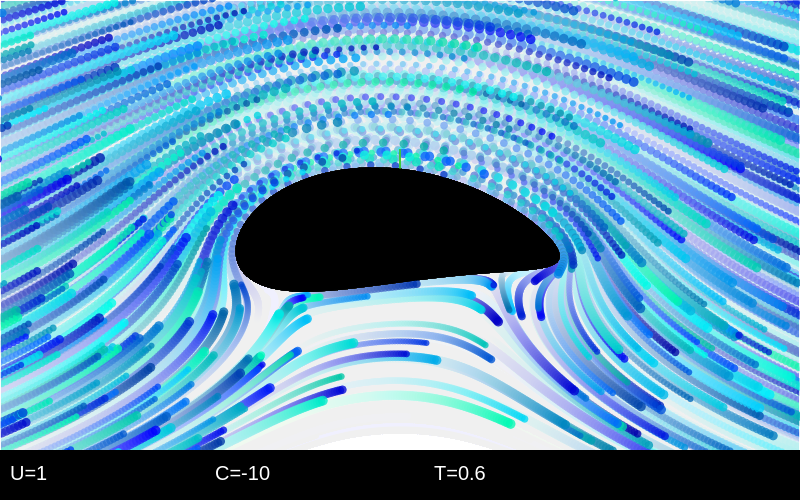} \\ $\qquad$
  \\
  \includegraphics[width=.3\linewidth]{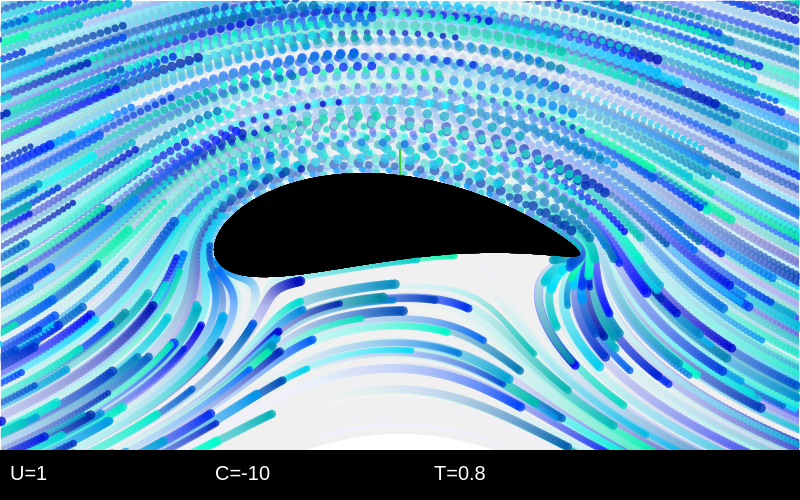}$\qquad$
  \includegraphics[width=.3\linewidth]{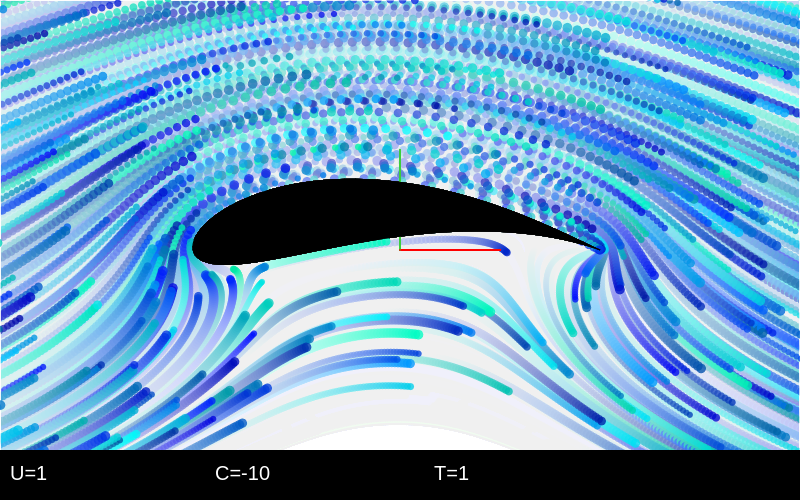} 
 \caption{\sl {F}rom a disk to an airfoil, and the corresponding fluid flow (images produced
 by the online application in
 {\tt https://complex-analysis.com/content/joukowsky\_airfoil.html}).}\label{GwedfgREEywuqwerthygfeigwegrhthrigheasfsdrirehg745NFIDItangeFI}
\end{figure}

Interestingly, on
{\tt https://complex-analysis.com/content/joukowsky\_airfoil.html}
one can also plot the resulting flow around the airfoil: the result that we have obtained
using this application is reported in Figure~\ref{GwedfgREEywuqwerthygfeigwegrhthrigheasfsdrirehg745NFIDItangeFI},
which nicely shows how a disk is mapped into an airfoil.

For further readings about the Kutta-Joukowski theory and the analysis of airfoils see e.g.~\cite{MR0112435, MR2650049, MR3503185}
and the references therein.
See also~\cite[pages~227--238]{MR0348082} for an extensive treatment of fluid dynamics from a complex analysis perspective.\medskip

\begin{figure}
                \centering
                \includegraphics[width=.45\linewidth]{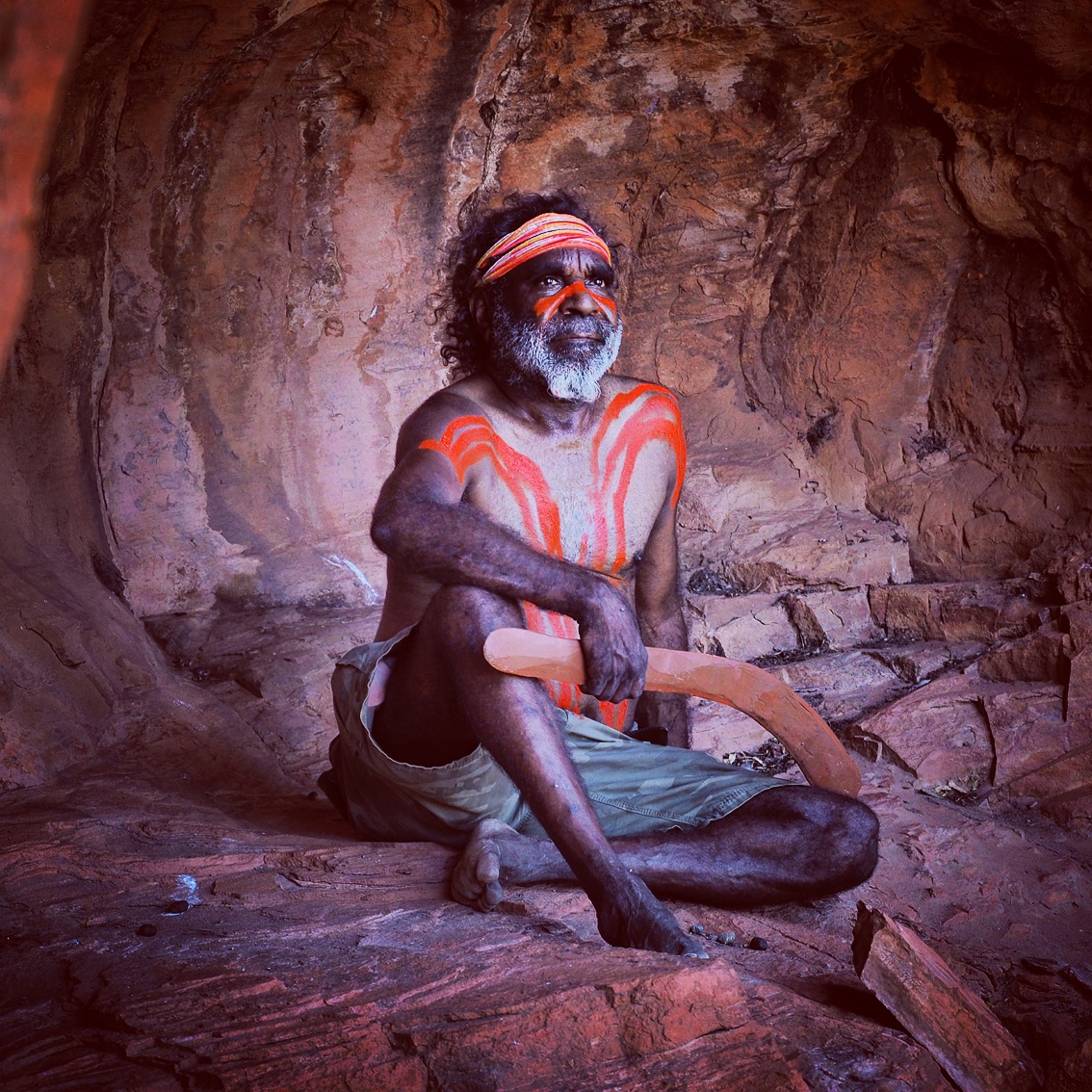}
        \caption{\sl Otto Jungarryi Sims with a boomerang sitting inside a cave in the Northern Territory of Australia
        (photo by Ed Gold; image from Wikipedia for free use under ticket \#2020101210010454).}\label{HARGIRA4AXELHARLROA7789GIJ7solFUMHDNOJHNFOJED7jmNThj-03}
\end{figure}

We cannot avoid mentioning that one of the most spectacular applications of fluid dynamics is probably showcased
by the motion\footnote{The origin of the name boomerang is a bit uncertain: some references link it to an extinct Aboriginal language of New South Wales, others to the language of the Turuwal people (a sub-group of the Darug) of the Georges River (Tucoerah River).
Besides the well-known traditional employment by some Aboriginal Australian peoples for hunting, see Figure~\ref{HARGIRA4AXELHARLROA7789GIJ7solFUMHDNOJHNFOJED7jmNThj-03}, it seems that boomerangs have been used also in ancient Europe, Egypt, and North America (ancient Egyptian boomerangs have been tested and seemed to work well as returning boomerangs and a boomerang discovered in the Carpathian Mountains in Poland dated back to about $30\times10^3$ years ago.

Interestingly, a boomerang
was used to set the Guinness World Record 
for the longest throw of any object by a human: namely, in 2005, at Murarrie Recreation Ground, in Queensland,
David Schummy 
performed a throw of 427.2 metres, see {\tt https://www.youtube.com/watch?v=ly3nCEbcQig}

Currently, long distance boomerangs are mostly 
shaped as a question mark and often have
a beveled edge, to facilitate the pitch
and lower the drift, since the boomerang in this case
is usually thrown almost horizontally (indeed, in the above mentioned
record, the objective
was not to make the boomerang come back to the throw location;
actually the boomerang ended up on a tree).} of a boomerang.

\begin{figure}
                \centering
                \includegraphics[width=.35\linewidth]{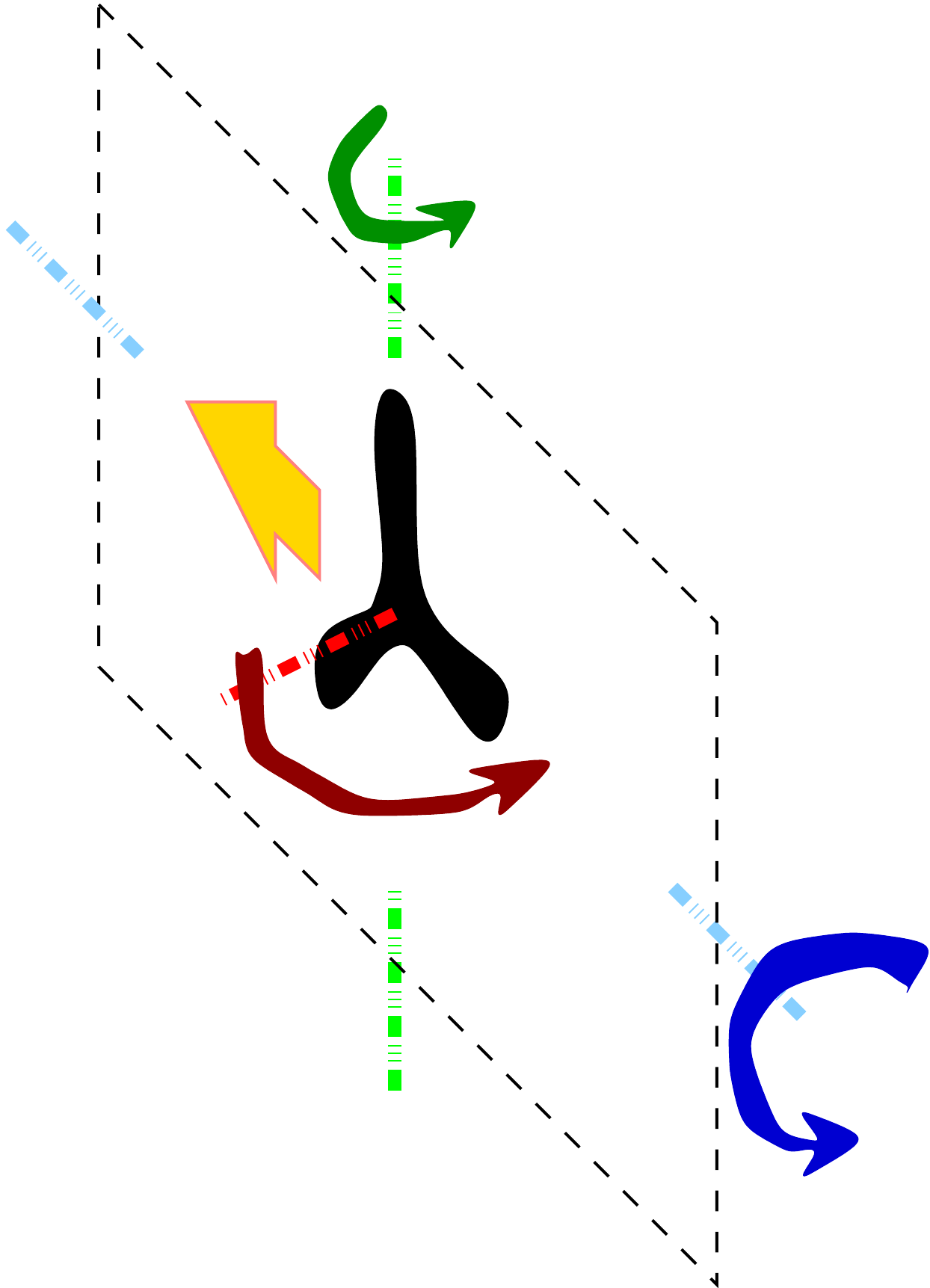}
        \caption{\sl Why do boomerangs return?}\label{BomHARGIRA4AXELHARLROA7789GIJ7solFUMHDNOJHNFOJED7jmNThj-03}
\end{figure}

Roughly speaking, each wing of a boomerang is shaped as an airfoil section,
allowing the airflow over the wings to create a significant lift.
If the boomerang is thrown nearly upright, the rotating blades generate more lift at the top than the bottom,
because at the top the speed of the rotation adds up to the forward speed, while at the bottom
the speed of the rotation subtracts from the forward speed, see Figure~\ref{BomHARGIRA4AXELHARLROA7789GIJ7solFUMHDNOJHNFOJED7jmNThj-03}
in which the forward speed is represented by the yellow arrow.
This additional lift from the top produces a torque (represented by the blue arrow in Figure~\ref{BomHARGIRA4AXELHARLROA7789GIJ7solFUMHDNOJHNFOJED7jmNThj-03})
whose effect is to make the rotation plane turn around: we point out that this torque is typically not sufficient
to tilt the boomerang around its axis of travel, given its high
angular momentum (the spinning of the boomerang being represented by the red arrow in Figure~\ref{BomHARGIRA4AXELHARLROA7789GIJ7solFUMHDNOJHNFOJED7jmNThj-03}):
hence, the stability of the rotating plane ensured by the gyroscopic precession combines with the aerodynamic torque
and leads to the curved trajectory sketched by the green arrow in Figure~\ref{BomHARGIRA4AXELHARLROA7789GIJ7solFUMHDNOJHNFOJED7jmNThj-03}.

\subsection{Surfing the waves}

A topical problem in the dynamics of fluids consists in the description of waves in shallow waters,
since this is a typical case arising in the proximity of land. Though a large number of different models are available for this goal,
in these pages we recall the classical equation

\begin{figure}
                \centering
                \includegraphics[width=.6\linewidth]{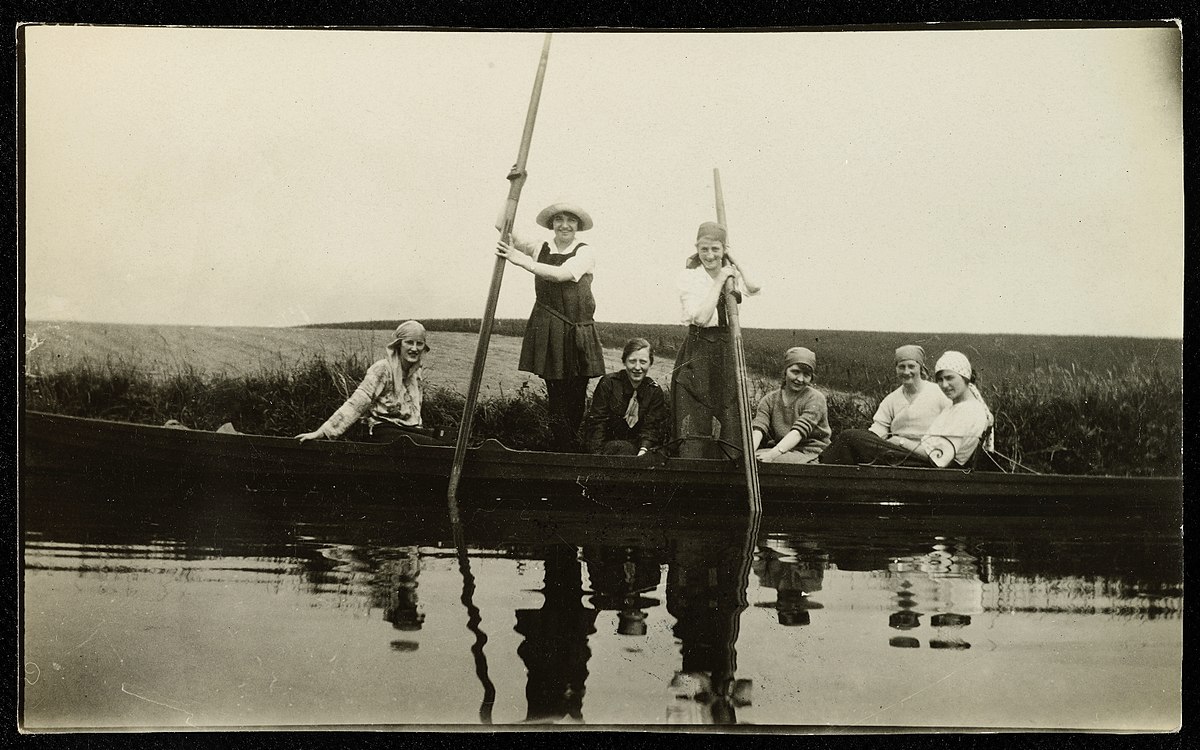}
        \caption{\sl Students picnic on the Union Canal in 1922
        (copyright The University of Edinburgh, available for public use; image from
        Wikipedia,
        licensed under the Creative Commons Attribution-Share Alike 3.0 Unported license).}\label{CROS12HAFODUNIOKsdxOLM789GIJ7solFUuuMHDNOJHNFOJED231}
\end{figure}

\begin{equation}\label{KDV-or}
a\partial _{t}\phi +b\partial _{x}^{3}\phi +c\,\phi \,\partial _{x}\phi =0
\end{equation}
for~$a$, $b$, $c\in\R\setminus\{0\}$, with~$\phi=\phi(x,t)$, $x\in\R$ and~$t\in(0,+\infty)$.
Equation~\eqref{KDV-or} is called\footnote{Equation~\eqref{KDV-or}
is named after Diederik Johannes Korteweg and Gustav de Vries~\cite{MR3363408}, though it was
introduced by Joseph Valentin Boussinesq~\cite{zbMATH02713295}.
De Vries completed his PhD under Korteweg's supervision and then worked all his life
as a high school teacher in Haarlem in the Netherlands.

The strong interest in the formation and propagation of waves in canals was possibly the outcome of a direct observation by
John Scott Russell~\cite{RUSSE} which took place in the Union Canal
(a canal in Scotland, running from Falkirk to Edinburgh).
Russell's own words have been repeated in virtually all papers and books which even remotely discuss water wave problems and we have no intention of breaking this consolidated tradition, hence here is his report on the astounding experience of meeting a traveling wave for the first time: ``I was observing the motion of a boat which was rapidly drawn along a narrow channel by a pair of horses, when the boat suddenly stopped -- not so the mass of water in the channel which it had put in motion; it accumulated round the prow of the vessel in a state of violent agitation, then suddenly leaving it behind, rolled forward with great velocity, assuming the form of a large solitary elevation, a rounded, smooth and well-defined heap of water, which continued its course along the channel apparently without change of form or diminution of speed. [...] Such, in the month of August 1834, was my first chance interview with that singular and beautiful phenomenon which I have called the Wave of Translation''.

See Figure~\ref{CROS12HAFODUNIOKsdxOLM789GIJ7solFUuuMHDNOJHNFOJED231} for a historical picture
of the Union Canal.} the
Korteweg-de Vries equation (or KdV equation for short). \index{Korteweg-de Vries equation|(}

The physical content of equation~\eqref{KDV-or} is that~$\phi$ models the shape of a wave in a canal.
The position on the canal is given by~$x\in\R$ and~$t$ stands for the time variable.
The amplitude of the wave and the depth of the channel are supposed to be small
(as will be discussed below in the approximations leading to the derivation of the equation),
hence the validity of~\eqref{KDV-or} has to be limited to small amplitude waves in shallow waters.

Though equation~\eqref{KDV-or} should be considered as a very simplified model, not able to capture the complexity
of oscillatory phenomena and waves in the real world, it is interesting to notice that the Korteweg-de Vries equation
can describe some interesting features such as traveling waves
and provide interesting information on the shape of the waves, as we now discuss.

\begin{figure}
  \centering
  \includegraphics[width=.45\linewidth]{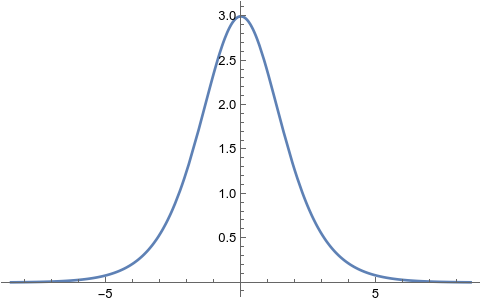}
 \caption{\sl The traveling wave solution for the Korteweg-de Vries equation, as found in
 in~\eqref{DIKDVSOLIFI-EQ}.}  \label{DIKDVSOLIFI}
\end{figure}

First of all (recalling the presentation on page~\pageref{01ojehniTNasIplamnewa}) one can seek for
traveling wave solutions of~\eqref{KDV-or} in the form
\begin{equation}\label{BSNEGIAaon}
\phi(x,t)=\phi_0(x-vt)\end{equation}
for some velocity~$v$.
In this setting, equation~\eqref{KDV-or} reduces to
\begin{equation}\label{KAlyhchtvdhr4htFAn}
0=-av\phi_0' +b\phi'''_0 +c\,\phi_0 \,\phi_0' =-av\phi_0' +b\phi'''_0 +\frac{c}2\,(\phi_0^2)'=\left(-av\phi_0 +b\phi''_0 +\frac{c}2\,\phi_0^2\right)'.
\end{equation}
We integrate this equation by assuming that at infinity~$\phi_0$ 
and its derivatives converge to zero, finding that
\begin{equation*}0=
-av\phi_0 +b\phi''_0 +\frac{c}2\,\phi_0^2.
\end{equation*}
As a result,
\begin{equation*}
0=\left( -av\phi_0 +b\phi''_0 +\frac{c}2\,\phi_0^2\right)\phi_0'=
\left(-\frac{av}2\phi_0^2 +\frac{b}2 (\phi'_0)^2 +\frac{c}{6}\,\phi_0^3 \right)'.
\end{equation*}
Integrating this as above, we deduce that
\begin{equation*} 0
=-\frac{av}2\phi_0^2 +\frac{b}2 (\phi'_0)^2+\frac{c}{6}\,\phi_0^3.\end{equation*}
This ordinary differential equation
can be solved by separation of variables (assuming~$\phi_0\ge0$, $c>0$ and~$abv>0$) leading,
up to a translation, to
\begin{equation}\label{DIKDVSOLIFI-EQ} \phi_0(r)=
\frac{6 a v}{c \left( \cosh\left(r \sqrt{\frac{a v}b} \right)+ 1\right)}.\end{equation}
See Figure~\ref{DIKDVSOLIFI} for a sketch\footnote{See also {\tt http://lie.math.brocku.ca/$\sim$sanco/solitons/kdv\_solitons.php}
for several animations of (possibly interacting) traveling waves
of the Korteweg-de Vries equation.} of this traveling wave (recall~\eqref{BSNEGIAaon})
when~$a:=1$, $b:=1$, $c:=1$ and~$v:=1$.

It is also interesting to observe that the highest crest of the wave in~\eqref{DIKDVSOLIFI-EQ} is attained for~$r:=0$
and it is equal to
\begin{equation}\label{d3ae9ikjm3f82KSMpf334}
\frac{3 a v}{c}.\end{equation} This suggests that the fastest waves (i.e., waves with largest velocity~$v$) correspond to the
highest ones.\medskip

\begin{figure}
                \centering
                \includegraphics[width=.55\linewidth]{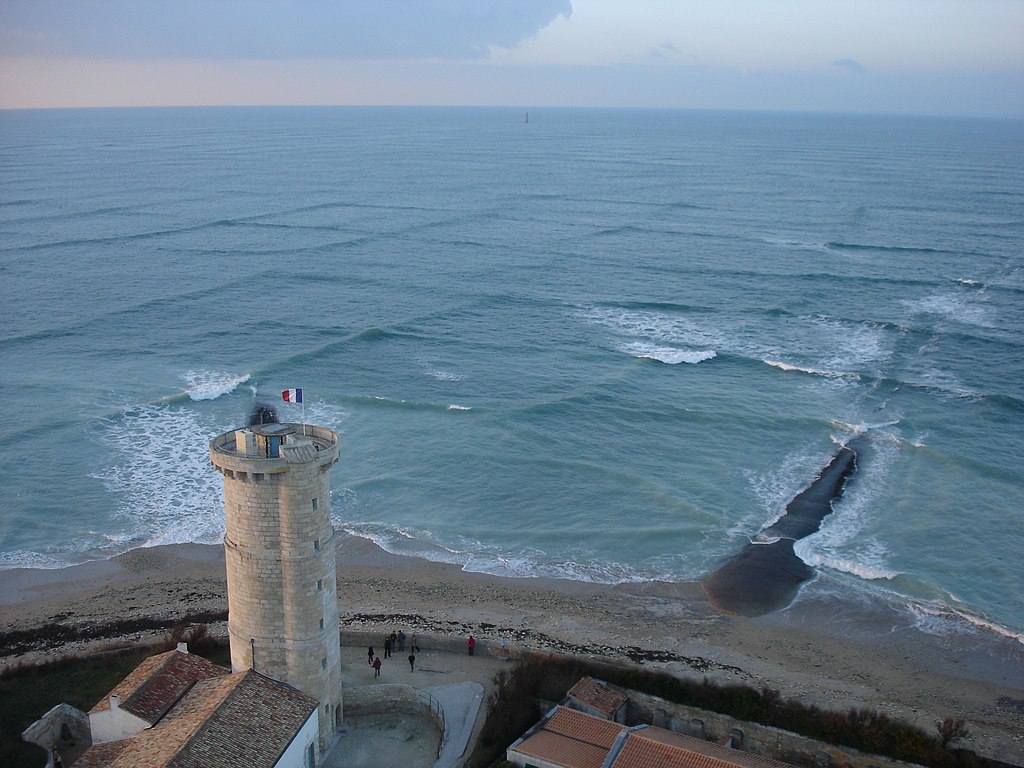}
        \caption{\sl Wave trains crossing in front of \^{I}le de R\'e, in the Atlantic Ocean
        (photo by Michel Griffon; image from
        Wikipedia,
        licensed under the Creative Commons Attribution-Share Alike 3.0 Unported license).}\label{CROS12HAFODAKEDFUMSPierre-SimonLaplacldRGIRA4AXELEUMDJOMNFHARLROA7KOLM789GIJ7solFUMHDNOJHNFOJED231}
\end{figure}

\begin{figure}
                \centering
                \includegraphics[width=.85\linewidth]{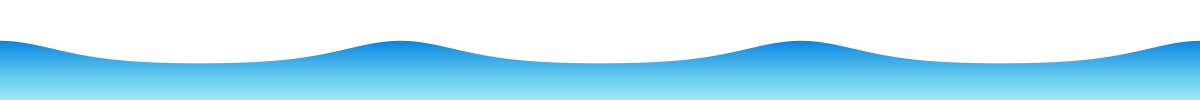}
        \caption{\sl Cnoidal wave solution to the Korteweg-de Vries equation (image by Kraaiennest from
        Wikipedia,
        licensed under the Creative Commons Attribution-Share Alike 3.0 Unported license).}\label{2HAFODAKEDFUMSPierre-SimonLaPALSusM789GIJ7solFUMHDNOJHNFOJED231}
\end{figure}

\begin{figure}
                \centering
                \includegraphics[width=.5\linewidth]{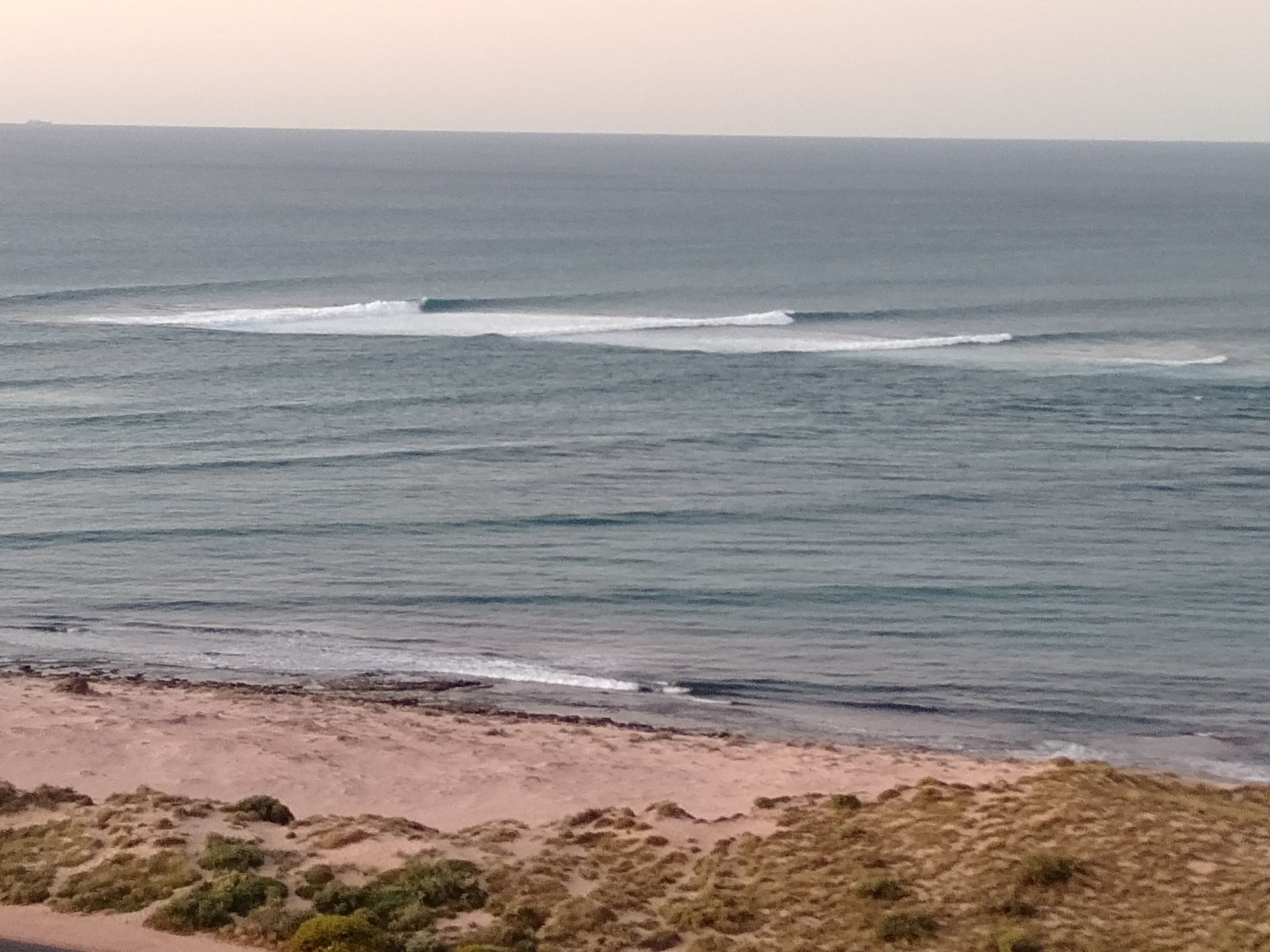}
        \caption{\sl Cnoidal looking waves seen from North West Cape.}\label{2HAFODAKEDFUMSPierre-SimonLaPALSusM789GIJ7solFUMHDNOJHNFOJED231-99-0987654rthKSMSj}
\end{figure}

Another nice feature captured by  the Korteweg-de Vries equation in~\eqref{KDV-or}
is a sufficiently accurate description of the shape of a set of waves: it is indeed a common experience that
``real'' waves exhibit sharper crests and flatter troughs than those of sine and cosine functions,
see Figures~\ref{CROS12HAFODAKEDFUMSPierre-SimonLaplacldRGIRA4AXELEUMDJOMNFHARLROA7KOLM789GIJ7solFUMHDNOJHNFOJED231}, \ref{2HAFODAKEDFUMSPierre-SimonLaPALSusM789GIJ7solFUMHDNOJHNFOJED231}
and~\ref{2HAFODAKEDFUMSPierre-SimonLaPALSusM789GIJ7solFUMHDNOJHNFOJED231-99-0987654rthKSMSj}.

To understand the shape of these waves, it is convenient to define, given~$\varphi\in\R$ and~$m\in[0,1)$,
\begin{equation}\label{NokmsLPSo8s0} U(\varphi,m):=\int _{0}^{\varphi }{\frac{d\theta }{\sqrt {1-m\sin ^{2}\theta }}}.\end{equation}
We observe that~$\frac{\partial U}{\partial \varphi}>0$, hence, given~$m\in(0,1)$, the function~$\varphi\mapsto U(\varphi,m)$ is invertible
and, with a slight abuse of notation, we call~$\varphi(U,m)$ the inverse function.
In this setting, one defines the elliptic cosine (in Latin, cosinus amplitudinis) by \index{elliptic cosine}
\begin{equation}\label{NokmsLPSo8s0BIS} \operatorname{cn} (U,m):=\cos (\varphi(U,m)).\end{equation}
Now, to find periodic solutions of the Korteweg-de Vries equation in~\eqref{KDV-or}
and to compare their shapes with Figure~\ref{2HAFODAKEDFUMSPierre-SimonLaPALSusM789GIJ7solFUMHDNOJHNFOJED231} we go back to~\eqref{KAlyhchtvdhr4htFAn}
and we integrate it, now not assuming that~$\phi_0$ decays at infinity: in this way, we have that
$$ -av\phi_0 +b\phi''_0 +\frac{c}2\,\phi_0^2=\kappa_1,$$
for some constant~$\kappa_1\in\R$.

As a consequence,
$$0=\Big(-av\phi_0 +b\phi''_0 +\frac{c}2\,\phi_0^2-\kappa_1\Big)\phi_0'=
\left(-\frac{av}2\phi_0^2 +\frac{b}2 (\phi'_0)^2 +\frac{c}6\,\phi_0^3-\kappa_1\phi_0\right)'
,$$
from which a further integration produces
\begin{equation}\label{MS-pqewkrfciaub} -\frac{av}2\phi_0^2 +\frac{b}2 (\phi'_0)^2 +\frac{c}6\,\phi_0^3-\kappa_1\phi_0=\kappa_2,\end{equation}
for some constant~$\kappa_2\in\R$.

We observe that we can write~\eqref{MS-pqewkrfciaub} as
\begin{equation}\label{MS-pqewkrfciaub2} (\phi'_0)^2={\mathcal{P}}(\phi_0),\end{equation}
for a suitable polynomial~${\mathcal{P}}$ of degree~$3$. The setting that we consider now
is that in which~${\mathcal{P}}$ possesses three distinct real roots~$\beta_1>\beta_2>\beta_3$.
With this, we rewrite~\eqref{MS-pqewkrfciaub2} as
\begin{equation}\label{MS-pqewkrfciaub-x}
(\phi'_0)^2=-\kappa\,(\phi_0-\beta_1)(\phi_0-\beta_2)(\phi_0-\beta_3),\end{equation}
where~$\kappa:=c/(3b)$, that we assume to be positive.

Now it turns out to be useful to seek solutions in the form
\begin{equation}\label{MS-pqewkrfciaub-x98ytr} \phi_0(r)=\beta_1\cos^2 (\Phi(r))+\beta_2\sin^2 (\Phi(r)), \end{equation}
where the function~$\Phi$ is supposed to be increasing and has to be determined. 

We observe that
$$ (\phi_0')^2=\Big( -2\beta_1\cos\Phi\,\sin\Phi\, \Phi'+2\beta_2\sin\Phi \,\cos\Phi\,\Phi'\Big)^2=
4(\beta_1-\beta_2)^2 \sin^2\Phi\,\cos^2\Phi\,(\Phi')^2.$$
Additionally,
\begin{eqnarray*}
\phi_0-\beta_1=\beta_1\big( \cos^2 \Phi -1\big)+\beta_2\sin^2 \Phi=-(\beta_1-\beta_2)\sin^2 \Phi
\end{eqnarray*}
and
\begin{eqnarray*}
\phi_0-\beta_2=\beta_1 \cos^2 \Phi +\beta_2\big(\sin^2 \Phi-1)=(\beta_1-\beta_2)\cos^2 \Phi.
\end{eqnarray*}
These observations and~\eqref{MS-pqewkrfciaub-x} yield that
\begin{equation*}
4(\Phi')^2=\kappa\,(\phi_0-\beta_3).
\end{equation*}
Thus, since
$$ \phi_0-\beta_3=\beta_1\Big(1-\sin^2 \Phi\Big)+\beta_2\sin^2 \Phi-\beta_3=(\beta_1-\beta_3)-
(\beta_1-\beta_2)\sin^2\Phi,
$$
we conclude that
\begin{equation}\label{MS-pqewkrfciaub-x89} \Phi'=\frac{\sqrt\kappa}2\,\sqrt{ (\beta_1-\beta_3)-(\beta_1-\beta_2)\sin^2\Phi}=
\kappa_0\,\sqrt{ 1-m\sin^2\Phi}
,\end{equation}
where
$$ \kappa_0:=\frac{\sqrt{\kappa\,(\beta_1-\beta_3)}}2\qquad{\mbox{and}}\qquad
m:=\frac{\beta_1-\beta_2}{\beta_1-\beta_3}.$$
We can thus utilize the separation of variable method in~\eqref{MS-pqewkrfciaub-x89}
and, in view of~\eqref{NokmsLPSo8s0}, obtain that
$$ \kappa_0 r=\int_{0}^{\Phi(r)}\frac{d\vartheta}{\sqrt{ 1-m\sin^2\vartheta}}=U\big(\Phi(r),m\big),
$$
where we normalized the picture so that~$\Phi(0)=0$.

{F}rom this and~\eqref{NokmsLPSo8s0BIS} it follows that
$$ \operatorname{cn}(\kappa_0 r)=\operatorname{cn}\big( U\big(\Phi(r),m\big),m\big)=\cos(\Phi(r)),$$
where we have used the short setting~$\operatorname{cn}(\cdot):=
\operatorname{cn}(\cdot,m)$ to ease notation.

\begin{figure}
                \centering
                \includegraphics[width=.21\linewidth]{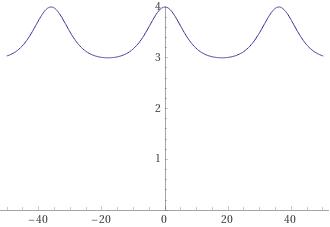}$\quad$
                \includegraphics[width=.21\linewidth]{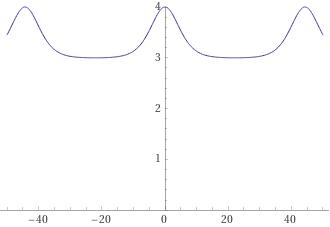}$\quad$
                                \includegraphics[width=.21\linewidth]{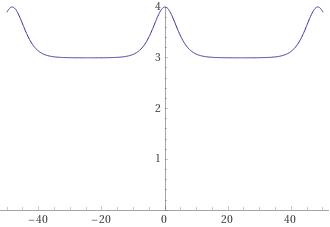}$\quad$
                                                \includegraphics[width=.21\linewidth]{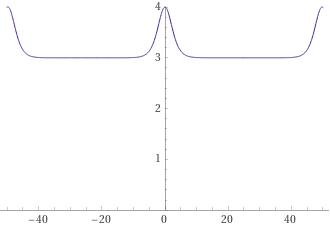}
        \caption{\sl The cnoidal waves in~\eqref{EQcnoidal-waves}
        with~$A:=1$, $B:=3$ and~$(\kappa_0,m)\in\big\{ (1/7,0.9),\,(1/6,0.99),\,(1/5,0.999),\,(1/3,0.999999)
        \big\}$.}\label{CROS12HAFODAKEDFUMSPieCNLEUMDJOMNFHARLROA7KOLM789GIJ7solFUMHDNOJHNFOJED231}
\end{figure}

Going back to~\eqref{MS-pqewkrfciaub-x98ytr} we thereby conclude that
\begin{equation}\label{EQcnoidal-waves} \begin{split}&\phi_0(r)=\beta_1\cos^2 (\Phi(r))+\beta_2\big(1-\cos^2 (\Phi(r))\big)
=(\beta_1-\beta_2)\cos^2 (\Phi(r))+\beta_2\\&\qquad\quad\qquad=(\beta_1-\beta_2)\operatorname{cn}^2(\kappa_0 r )+\beta_2=
A\operatorname{cn}^2(\kappa_0 r)+B,\end{split}
 \end{equation}
where~$A:=\beta_1-\beta_2$ and~$B:=\beta_2$.
The profile of these functions, for suitable choices of parameters, resembles the shape of waves
in Figures~\ref{CROS12HAFODAKEDFUMSPierre-SimonLaplacldRGIRA4AXELEUMDJOMNFHARLROA7KOLM789GIJ7solFUMHDNOJHNFOJED231} and~\ref{2HAFODAKEDFUMSPierre-SimonLaPALSusM789GIJ7solFUMHDNOJHNFOJED231} (which are often
called ``cnoidal waves'' \index{cnoidal wave} due to the presence of the elliptic cosine~$\operatorname{cn}$
in their expression). To have a feeling on how the different parameters
change the shape of a cnoidal waves see Figure~\ref{CROS12HAFODAKEDFUMSPieCNLEUMDJOMNFHARLROA7KOLM789GIJ7solFUMHDNOJHNFOJED231}.\medskip

\begin{figure}
                \centering
                \includegraphics[width=.22\linewidth]{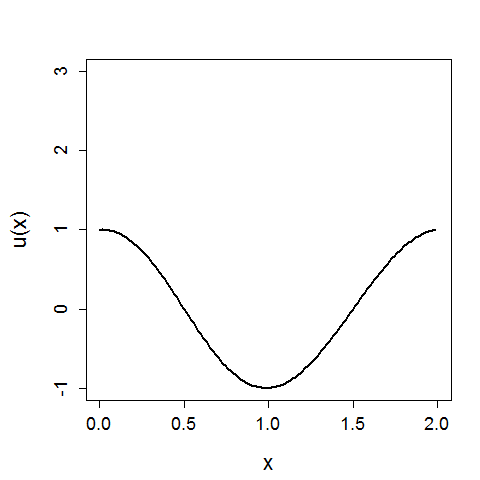} $\quad$
                \includegraphics[width=.22\linewidth]{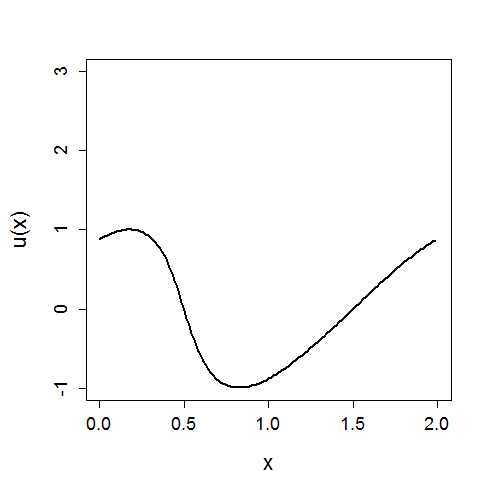} $\quad$
                \includegraphics[width=.22\linewidth]{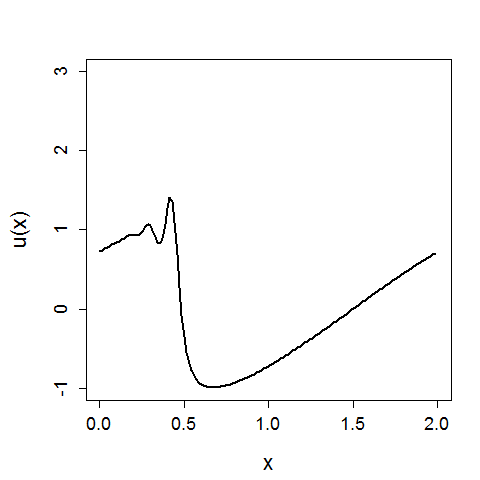} $\quad$
                \includegraphics[width=.22\linewidth]{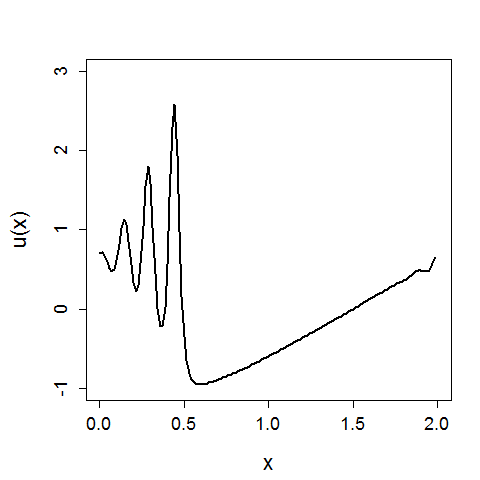} \\
                   \includegraphics[width=.22\linewidth]{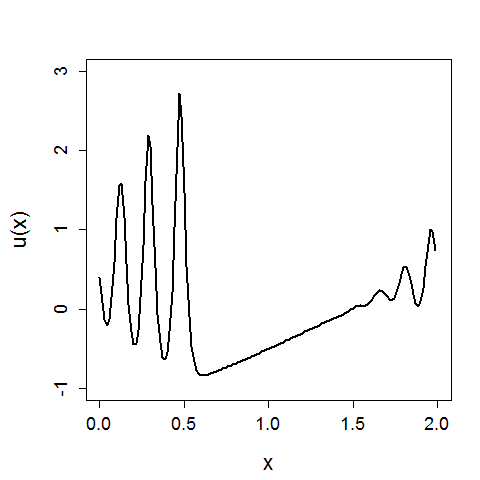} $\quad$
                \includegraphics[width=.22\linewidth]{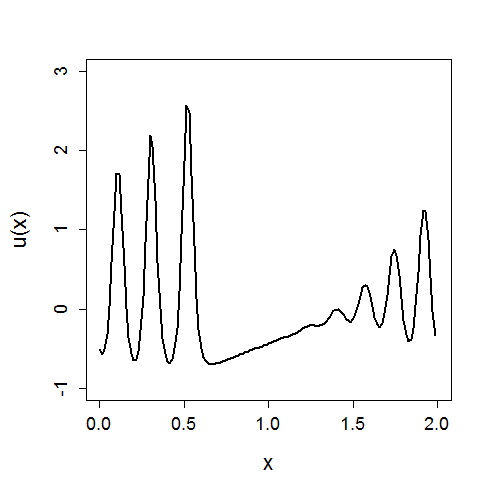} $\quad$
                \includegraphics[width=.22\linewidth]{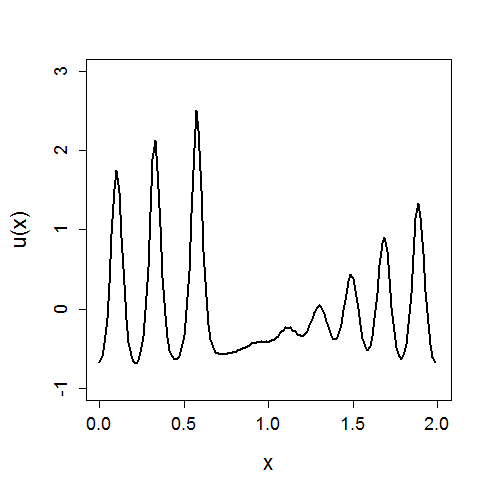} $\quad$
                \includegraphics[width=.22\linewidth]{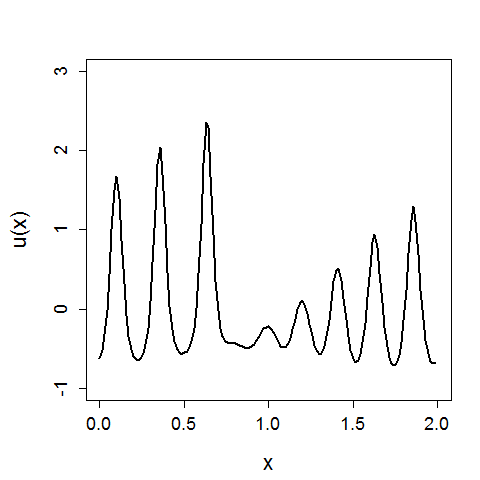}
        \caption{\sl Numerical solution of the KdV equation~$\partial_t u + u\partial_x u + (0.022)^2\partial_x^3u = 0$
with initial condition~$ u(x, 0) = \cos(\pi x)$.
The initial cosine wave evolves into a train of solitary-type waves~\cite{zbMATH05826808}
(by Ta2o; image from
        Wikipedia,
        licensed under the Creative Commons Attribution-Share Alike 3.0 Unported license).}\label{NUMEKSDUNIOKsdKDxOLM789GIJ7solFUuuMHDNOJHNFOJED231}
\end{figure}

It is also interesting to take note of the fact that, in spite of its (relative) simplicity, the Korteweg-de Vries equation
in~\eqref{KDV-or} is capable of capturing complex phenomena such as the formation of oscillatory waves, with crests
traveling at different speeds according to their heights: these complicated situations may even arise from very simple
initial data, see Figure~\ref{NUMEKSDUNIOKsdKDxOLM789GIJ7solFUuuMHDNOJHNFOJED231}.\medskip

\begin{figure}
                \centering
                \includegraphics[height=.24\textheight]{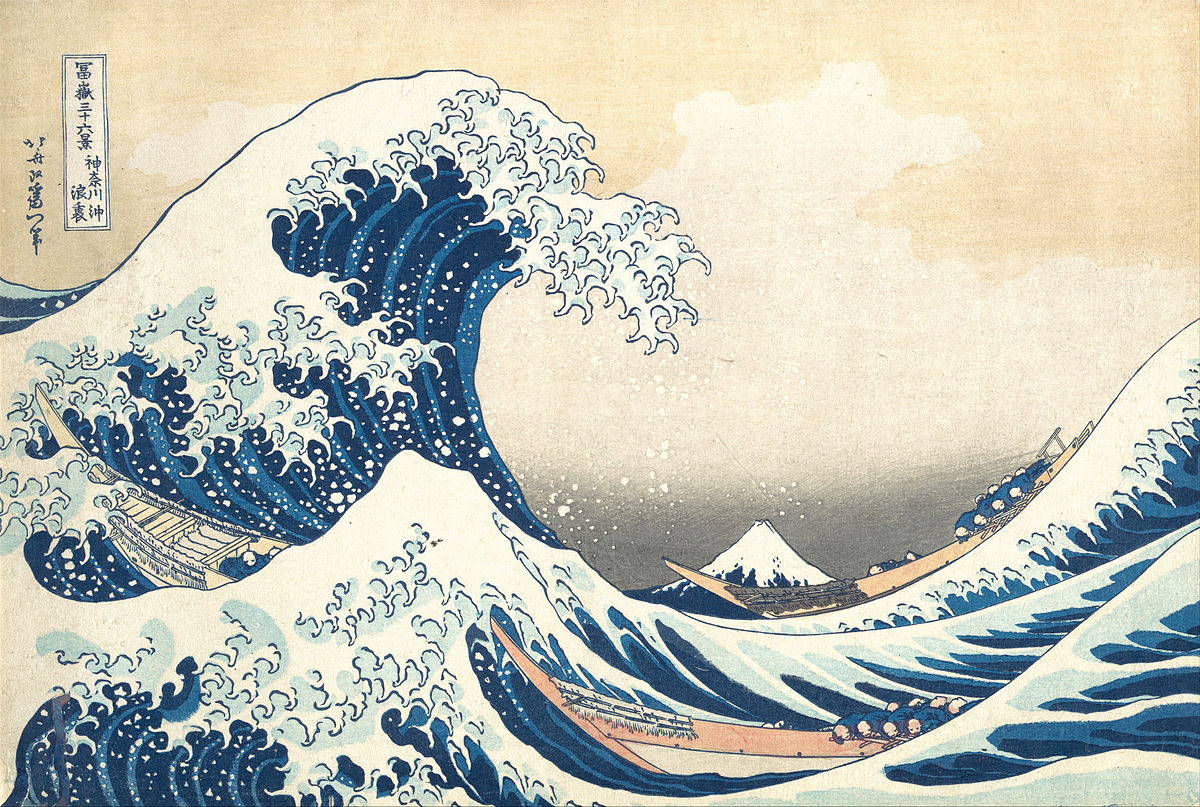} $\quad$
                \includegraphics[height=.24\textheight]{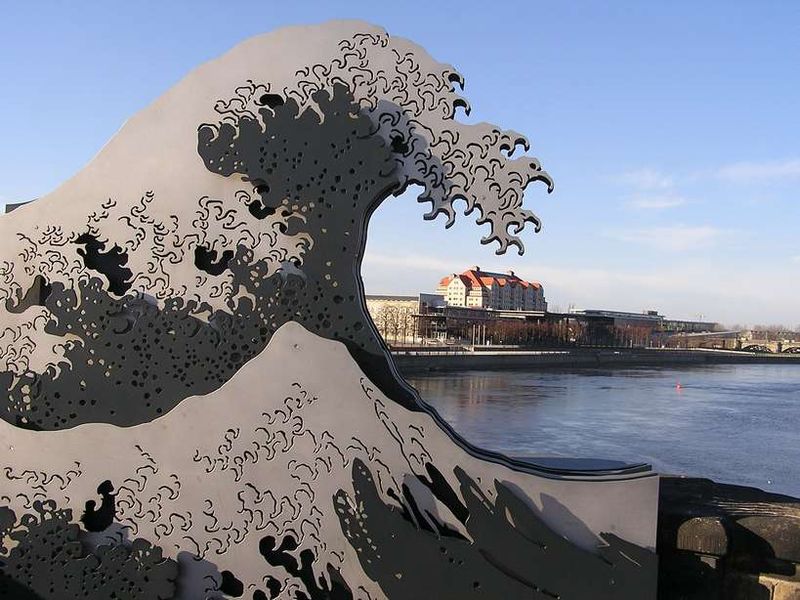}
        \caption{\sl Left:
        The Great Wave off Kanagawa, famous woodblock print by the ukiyo-e artist Katsushika Hokusai
(H. O. Havemeyer Collection, Bequest of Mrs. H. O. Havemeyer, 1929, Metropolitan Museum of Art; Public Domain
image from Wikipedia). Right: Die Woge, sculpture by Tobias Stengel in Dresden, Germany (photo by Christoph M\"unch,
Press Section of the city of Dresden's official homepage; image from
        Wikipedia,
        licensed under the Creative Commons Attribution-Share Alike 3.0 Unported license).}\label{WOGEOLM789GIJ7solFUuuMHDNOJHNFOJED231}
\end{figure}

Furthermore, we point out that the Korteweg-de Vries equation often turns in handy
in the description of natural phenomena that seem to be well beyond its range of applicability.
For example (see e.g.~\cite{TSUNA} and the references therein), the Korteweg-de Vries equation
has been exploited in the description\footnote{Although it is debatable whether the artists
intended their works to be interpreted in that way, it is customary to
refer to the artworks in Figure~\ref{WOGEOLM789GIJ7solFUuuMHDNOJHNFOJED231}
in connection with tsunamis.}
of tsunamis \index{tsunami}
even though the extreme depth of oceans may appear to be incompatible with an equation which was originally
designed for shallow waters. The reason for the success of the Korteweg-de Vries equation even in this setting
is, at least, twofold. First of all, though oceans are very deep, tsunami waves can reach a spatial
extent of even larger size: therefore, when the depth of the ocean is sufficiently small with respect to the breadth of the wave
the shallow water theory may still be applied, since what counts is the relative (and not absolute) sizes of the length parameters.
Secondly, while tsunamis in the deep ocean consists typically of very long waves with quite small amplitudes, as they approach shallow waters the nonlinear effects become predominant and they change significantly
the shape and velocity of the leading waves, for which the Korteweg-de Vries equation can also provide an approximate, but often
effective, model.
\medskip

Now we give a brief motivation for the Korteweg-de Vries equation
in~\eqref{KDV-or}. While its expression may appear rather mysterious
at first sight, especially in view of the third derivative appearing as a leading term, we will
see now that equation~\eqref{KDV-or} arises naturally from the equations of fluid dynamics
presented in Section~\ref{EUMSD-OS-32456i7DNSSE234R:SEC}.

More specifically, we consider the Euler fluid flow equation in~\eqref{EUMSD-OS-D}.
Here we take coordinates~$(x,z)\in\R^2$ where~$x\in\R$ represents the direction of a canal
and~$z\in\R$ is the vertical displacement. The fluid velocity~$v$ will be written in components as~$v=(\overline u,\overline w)$,
with~$\overline u$ and~$\overline w$ scalar functions. Hence, being~$\rho$ the density of the fluid (assumed to be constant)
and~$p$ the pressure, equation~\eqref{EUMSD-OS-D} reads
\begin{equation}\label{Cm-PAPMam}
\rho \partial_t (\overline u,\overline w)+\rho ((\overline u,\overline w)\cdot(\partial_x,\partial_z)) (\overline u,\overline w)=
-\rho g(0,1) -(\partial_x,\partial_z) p.
\end{equation}
Hence, setting
\begin{equation}\label{Cm-PAPMam20} \overline P:=\frac{p}\rho+gz,
\end{equation}
we rewrite~\eqref{Cm-PAPMam} as the following system of two equations:
\begin{equation}\label{Cm-PAPMam2}
\begin{dcases}
\partial_t \overline u+
\overline u\partial_x\overline u+\overline w \partial_z \overline u=- \partial_x\overline P,
\\
\partial_t \overline w+
\overline u \partial_x\overline w
+\overline w \partial_z \overline w=
-\partial_z\overline P.
\end{dcases}
\end{equation}

We suppose that the bottom of the canal is~$z=0$ and the fluid at rest has height~$h_0\in(0,+\infty)$.
The wave is thus described as a perturbation of the rest fluid surface~$z=h_0$
described by the graph of a function~$\overline\eta=\overline\eta(x,t)$, modulated by a small amplitude~$\alpha\in(0,+\infty)$:
namely, at every instant of time~$t$ the fluid in the channel corresponds to the region of the points~$(x,z)\in\R^2$
with~$0<z<h_0+\alpha\overline\eta(x,t)$. That is, the equations in~\eqref{Cm-PAPMam2} are set in the domain
\begin{equation}\label{Cm-PAPMam30} \Omega^{(t)}:=\Big\{ (x,z)\in\R^2{\mbox{ s.t. }}x\in\R{\mbox{ and }}z\in\big(0,h_0+\alpha\overline\eta(x,t)\big)\Big\}.\end{equation}
We also assume that the fluid is incompressible, hence, by~\eqref{INCOMPRE},
\begin{equation}\label{Cm-PAPMam3}
0=\div (\overline u,\overline w)=\partial_x\overline u+\partial_z\overline w.\end{equation}
The equations in~\eqref{Cm-PAPMam2} and~\eqref{Cm-PAPMam3} are complemented by boundary conditions
corresponding to the bottom of the canal and the free surface of the fluid.
Indeed, at the bottom~$z=0$ we suppose that the floor of the canal is impenetrable, hence the vertical component
of the fluid velocity vanishes, which reads
\begin{equation}\label{Cm-PAPMam4}
\overline w(x,0,t)=0.\end{equation}
As for the free surface, we assume that fluid parcels staying at the top remain at the top:
hence if~$(x(t),z(t))$ represents the trajectory of a fluid particle on the upper surface of the fluid at time~$t$ we have that
$$ z(t)=h_0+\alpha\overline\eta(x(t),t).$$
Thus, taking the time derivative of this expression and recalling that the fluid parcel velocity~$(\dot x(t),\dot z(t))$ agrees with~$(\overline u,\overline w)$,
\begin{equation}\label{Cm-PAPMam44}\begin{split}&
\overline w\big(x(t),h_0+\alpha\overline\eta(x(t),t),t\big)=\dot z(t)=\alpha\partial_x\overline\eta(x(t),t)\dot x(t)+\alpha\partial_t\overline\eta(x(t),t)\\&\quad=\alpha\partial_x\overline\eta(x(t),t)
\overline u\big(x(t),h_0+\alpha\overline\eta(x(t),t),t\big)
+\alpha\partial_t\overline\eta(x(t),t).\end{split}\end{equation}

The free surface~$z=h_0+\alpha\overline\eta$ also presents an additional boundary condition related to pressure.
Indeed, the pressure on the top surface of the fluid is balanced by the atmospheric pressure, that we denote by~$p_0$
and assume to be constant. Hence, by~\eqref{Cm-PAPMam20},
\begin{equation}\label{Cm-PAPMam45} \overline P\big(x,h_0+\alpha\overline\eta(x,t),t\big)=\frac{p_0}\rho+g\,\big(h_0+\alpha\overline\eta(x,t)\big).\end{equation}

We now consider a regime of shallow waters and small amplitudes.
That is, we pick~$\Theta>0$ and\footnote{To simplify the computation, one can take simply~$\Theta:=1$.
The role of the ``arbitrary'' constant~$\Theta$ is just to introduce a quantity
with the dimension of time. In this way, we can think that~$\overline\eta$, $\e$ and~$\delta$
are dimensionless, while~$\alpha$ and~$h_0$ have the dimension of space.
In the end, the auxiliary dimensional constant~$\Theta$ will simplify and will not appear
any longer in the final equation~\eqref{UHJNZzekmdLSyhn8ijP235}.}
we consider the quantities
\begin{equation}\label{EPSIDELKDV}
\e:=\frac{\alpha}{h_0}\qquad{\mbox{and}}\qquad\delta:=\frac{h_0}{\Theta^2 \,g}.\end{equation}
We suppose that~$\e$ and~$\delta$ are small and also, up to normalizing constants, that
\begin{equation}\label{sawqewet54b5494pppppplllllkkkk}
\delta=\e.\end{equation}
We define
\begin{equation}\label{RTArasmsImPNSd-2}\begin{split}& u^\star(x,z,t):=\frac{\overline u\big(\sqrt{gh_0} x,h_0z,t\big)}{\sqrt{gh_0}},\qquad
w^\star(x,z,t):=\frac{\Theta\,\overline w\big(\sqrt{gh_0} x,h_0z,t\big)}{h_0},\\&\qquad
\eta^\star(x,t):= {\overline\eta(\sqrt{gh_0} x,t)}
\qquad{\mbox{and}}\qquad
P^\star(x,z,t):=\frac{\overline P\big(\sqrt{gh_0} x,h_0z,t\big)}{{gh_0}}.
\end{split}\end{equation}
In view of~\eqref{Cm-PAPMam30}, the above functions are defined in the region
\begin{equation}\label{NEWPD-do1}
\begin{split}
\Omega_{ t }^\star&:=
\Big\{ (x,z)\in\R^2{\mbox{ s.t. }}(\sqrt{gh_0} x,h_0z)\in\Omega^{(t)}\Big\}\\&
=\Big\{ (x,z)\in\R^2{\mbox{ s.t. }}x\in\R{\mbox{ and }} z\in\big(0,1+\e\eta^\star( x,t)\big)\Big\}.\end{split}
\end{equation}
We thus collect the equations in~\eqref{Cm-PAPMam2}, \eqref{Cm-PAPMam3}, \eqref{Cm-PAPMam4}, \eqref{Cm-PAPMam44} and~\eqref{Cm-PAPMam45}, exploiting also~\eqref{EPSIDELKDV} and~\eqref{sawqewet54b5494pppppplllllkkkk}, by writing
\begin{equation}\label{w78uhdbc-384irton2rve7g}
\begin{dcases}
\partial_t u^\star+
u^\star\partial_xu^\star+\displaystyle\frac{w^\star \partial_z u^\star}\Theta=- \partial_x P^\star
& {\mbox{ in }}\Omega_{ t }^\star,
\\
\e\big( \Theta\partial_t w^\star+\Theta
u^\star \partial_x w^\star
+ w^\star \partial_z w^\star\big)=
-\partial_z P^\star & {\mbox{ in }}\Omega_{ t }^\star,\\ \Theta\partial_xu^\star+\partial_zw^\star=0& {\mbox{ in }}\Omega_{ t }^\star,\\
P^\star=\displaystyle\frac{p_0}{\rho gh_0}+{1+\e\eta^\star }
& {\mbox{ on }} z=1+\e\eta^\star ,\\
w^\star=\e\Theta\,(\partial_t\eta^\star+u^\star\partial_x\eta^\star)
& {\mbox{ on }} z=1+\e\eta^\star ,\\
w^\star=0& {\mbox{ on }} z=0.
\end{dcases}
\end{equation}
We let~$c_0:=\frac{p_0}{\rho gh_0}+1$.
It is now convenient to ``surf the wave'' and look at the new variable~$\xi:=x-t$ \label{SUERF:PAHNDFGE}
(roughly speaking, using this variable one is moving with a wave that travels at unit\footnote{Interestingly,
moving at unit speed for~$\eta^\star$ corresponds to moving at speed~$\sqrt{gh_0}$
for the original profile~$\overline\eta$, thanks to the change of spatial variables in~\eqref{RTArasmsImPNSd-2}.
The quantity~$\sqrt{gh_0}$ happens indeed to be one of the characteristic velocities of small waves
in shallow waters (relating the velocity to~$\sqrt{gh_0}$
also tells us that the speed of these waves decreases when the height of the canal
is smaller and, conversely, their speed increases when moving from very shallow to slightly deeper water).} speed).
It is also appropriate to look at large times by setting~$\tau:=\e t$. Thus, we rephrase our physical quantities with
respect to these new variables by defining
\begin{equation}\label{RTArasmsImPNSd-1}
\begin{split}& 
u(\xi,z,\tau):=u^\star\left(\xi+\frac\tau\e,z,\frac\tau\e\right),\qquad w(\xi,z,\tau)=
w^\star\left(\xi+\frac\tau\e,z,\frac\tau\e\right),\\&\qquad\eta(\xi,\tau)=
\eta^\star\left(\xi+\frac\tau\e,\frac\tau\e\right)
\qquad{\mbox{and}}\qquad P(\xi,z,\tau):=
P^\star\left(\xi+\frac\tau\e,z,\frac\tau\e\right).
\end{split}\end{equation}
By~\eqref{NEWPD-do1}, these functions are defined in the region
\begin{eqnarray*} \Omega_\tau&:=&\Big\{ (\xi,z)\in\R^2{\mbox{ s.t. $\xi=x-t$, $\tau=\e t$ and~$(x,z)\in\Omega_{ t }^\star$}}\Big\}\\&
=&\Big\{ (\xi,z)\in\R^2{\mbox{ s.t. }}\xi\in\R{\mbox{ and }} z\in\big(0,1+\e\eta( \xi,\tau)\big)\Big\}
.\end{eqnarray*}
In this setting, one can rewrite~\eqref{w78uhdbc-384irton2rve7g} in the form
\begin{equation}\label{KDV-0m-001}
\begin{dcases}
-\partial_\xi u+\e\partial_\tau u+u\partial_\xi u+\frac{w\partial_z u}\Theta=-\partial_\xi P&{\mbox{ in }}\Omega_\tau,\\
\e\big( -\Theta\partial_\xi w+\e\Theta\partial_\tau w+\Theta u\partial_\xi w+w\partial_z w\big)=-\partial_z P&{\mbox{ in }}\Omega_\tau,\\
\Theta\partial_\xi u+\partial_z w=0&{\mbox{ in }}\Omega_\tau,\\
P=c_0+\e\eta & {\mbox{ on }}z=1+\e\eta,\\
w=\e\Theta(\e\partial_\tau\eta+(u-1)\partial_\xi\eta)& {\mbox{ on }}z=1+\e\eta,\\
w=0 & {\mbox{ on }}z=0.
\end{dcases}
\end{equation}

It is now convenient to consider formal Taylor
expansions of the physical quantities involved in powers of~$\e$ by writing
\begin{eqnarray*}&& u(\xi,z,\tau)=\sum_{k=0}^{+\infty}\e^k u_k(\xi,z,\tau),\qquad
w(\xi,z,\tau)=\sum_{k=0}^{+\infty}\e^k w_k(\xi,z,\tau)\\&&
P(\xi,z,\tau)=\sum_{k=0}^{+\infty}\e^k P_k(\xi,z,\tau)
\qquad{\mbox{and}}\qquad
\eta(\xi,\tau)=\sum_{k=0}^{+\infty}\e^k \eta_k(\xi,\tau).\end{eqnarray*}
Substituting these expressions into~\eqref{KDV-0m-001} we obtain that
\begin{equation}\label{KDV-0m-002}
\begin{dcases}
-\displaystyle\sum_{k=0}^{+\infty}\e^k\partial_\xi  u_k+\displaystyle\sum_{k=0}^{+\infty}\e^{k+1}\partial_\tau u_k+
\displaystyle\sum_{k,m=0}^{+\infty}\e^{k+m}
u_k\partial_\xi u_m+\displaystyle\frac1\Theta\sum_{k,m=0}^{+\infty}\e^{k+m}
w_k\partial_z u_m=-\displaystyle\sum_{k=0}^{+\infty}\e^k\partial_\xi P_k&{\mbox{ in }}\Omega_\tau,\\
\begin{matrix}
-\Theta\displaystyle\sum_{k=0}^{+\infty}\e^{k+1}\partial_\xi w_k+
\Theta\displaystyle\sum_{k=0}^{+\infty}\e^{k+2}\partial_\tau w_k+\Theta
\displaystyle\sum_{k,m=0}^{+\infty}\e^{k+m+1}u_k\partial_\xi w_m+
\displaystyle\sum_{k,m=0}^{+\infty}\e^{k+m+1}w_k\partial_z w_m\\
=-\displaystyle\sum_{k=0}^{+\infty}\e^k\partial_z P_k\end{matrix}&{\mbox{ in }}\Omega_\tau,\\
\displaystyle\sum_{k=0}^{+\infty}\e^k(\Theta\partial_\xi u_k+\partial_z w_k)=0&{\mbox{ in }}\Omega_\tau,\\
P\left( \xi, 1+\displaystyle\sum_{k=0}^{+\infty}\e^{k+1}\eta_k,\tau\right)=c_0+\displaystyle\sum_{k=0}^{+\infty}\e^{k+1}\eta_k, & \\ \frac1\Theta\,
w\left( \xi, 1+\displaystyle\sum_{k=0}^{+\infty}\e^{k+1}\eta_k,\tau\right)=
\displaystyle\sum_{k=0}^{+\infty}\e^{k+2}\partial_\tau\eta_k+
\displaystyle\sum_{k,m=0}^{+\infty}\e^{k+m+1}u_k\partial_\xi\eta_m-
\displaystyle\sum_{k=0}^{+\infty}\e^{k+1}\partial_\xi\eta_k,& \\
w(\xi,0,\tau)=0. &
\end{dcases}
\end{equation}

The strategy is now to approximate~\eqref{KDV-0m-002} by neglecting the terms of size~$\e^k$ with~$k\ge3$.
For this, we consider the orders in~$\e^k$ in~\eqref{KDV-0m-002} with~$k\in\{0,1,2\}$ and we formally equate the corresponding
coefficients. Namely, by formally considering the coefficients related to~$\e^0$ in~\eqref{KDV-0m-002} we see that
\begin{equation*}
\begin{dcases}
-\partial_\xi  u_0+u_0\partial_\xi u_0+\frac{w_0\partial_z u_0}\Theta=-\partial_\xi P_0&{\mbox{ in }}\Omega_\tau,\\
0=\partial_z P_0&{\mbox{ in }}\Omega_\tau,\\
\Theta\partial_\xi u_0+\partial_z w_0=0&{\mbox{ in }}\Omega_\tau,\\
P_0( \xi, 1,\tau)=c_0, & \\
w_0( \xi, 1,\tau)=0,& \\
w_0(\xi,0,\tau)=0. &
\end{dcases}
\end{equation*}
Since we are interested in small perturbations of steady states, we also take~$u_0:=0$ and~$w_0:=0$, therefore
\begin{equation}\label{KDV-0m-003}
\begin{dcases}
0=-\partial_\xi P_0&{\mbox{ in }}\Omega_\tau,\\
0=\partial_z P_0&{\mbox{ in }}\Omega_\tau,\\
P_0( \xi, 1,\tau)=c_0. & 
\end{dcases}
\end{equation}
{F}rom the second and third lines in~\eqref{KDV-0m-003} it follows that
\begin{equation}\label{KDV-0m-005}P_0(\xi,z,\tau)=c_0.\end{equation}

Now we consider the coefficients related to~$\e^1$ in~\eqref{KDV-0m-002} and we obtain that
\begin{equation}\label{KDV-0m-004}
\begin{dcases}
-\partial_\xi  u_1=-\partial_\xi P_1&{\mbox{ in }}\Omega_\tau,\\
0=-\partial_z P_1&{\mbox{ in }}\Omega_\tau,\\
\Theta\partial_\xi u_1+\partial_z w_1=0&{\mbox{ in }}\Omega_\tau,\\
P_1(\xi,1,\tau)+\partial_z P_0(\xi,1,\tau)\eta_0=\eta_0, & \\
w_1(\xi,1,\tau)=-\Theta\partial_\xi\eta_0,& \\
w_1(\xi,0,\tau)=0. &
\end{dcases}
\end{equation} Owing to~\eqref{KDV-0m-005},
the fourth line in~\eqref{KDV-0m-004} can be written as~$P_1(\xi,1,\tau)=\eta_0(\xi,\tau)$, which
combined with the second line gives that
\begin{equation}\label{KDV-0m-009}
P_1(\xi,z,\tau)=\eta_0(\xi,\tau).\end{equation} Hence, from the first line,
\begin{equation}\label{TH:EVE}
\partial_\xi u_1=\partial_\xi\eta_0.\end{equation} This and the third line give that
$$ \partial_z \left(\frac{w_1}\Theta+z\partial_\xi\eta_0\right)=\frac{\partial_z w_1}\Theta+\partial_\xi\eta_0=-\partial_\xi u_1+\partial_\xi\eta_0=0
$$
and therefore, by the last line in~\eqref{KDV-0m-004}, we infer that
\begin{equation}\label{KDV-0m-006}\frac{w_1}\Theta=-z\partial_\xi\eta_0,\end{equation}
which is also in agreement with the second last line in~\eqref{KDV-0m-004}.

\begin{figure}
                \centering
                \includegraphics[height=5.6cm]{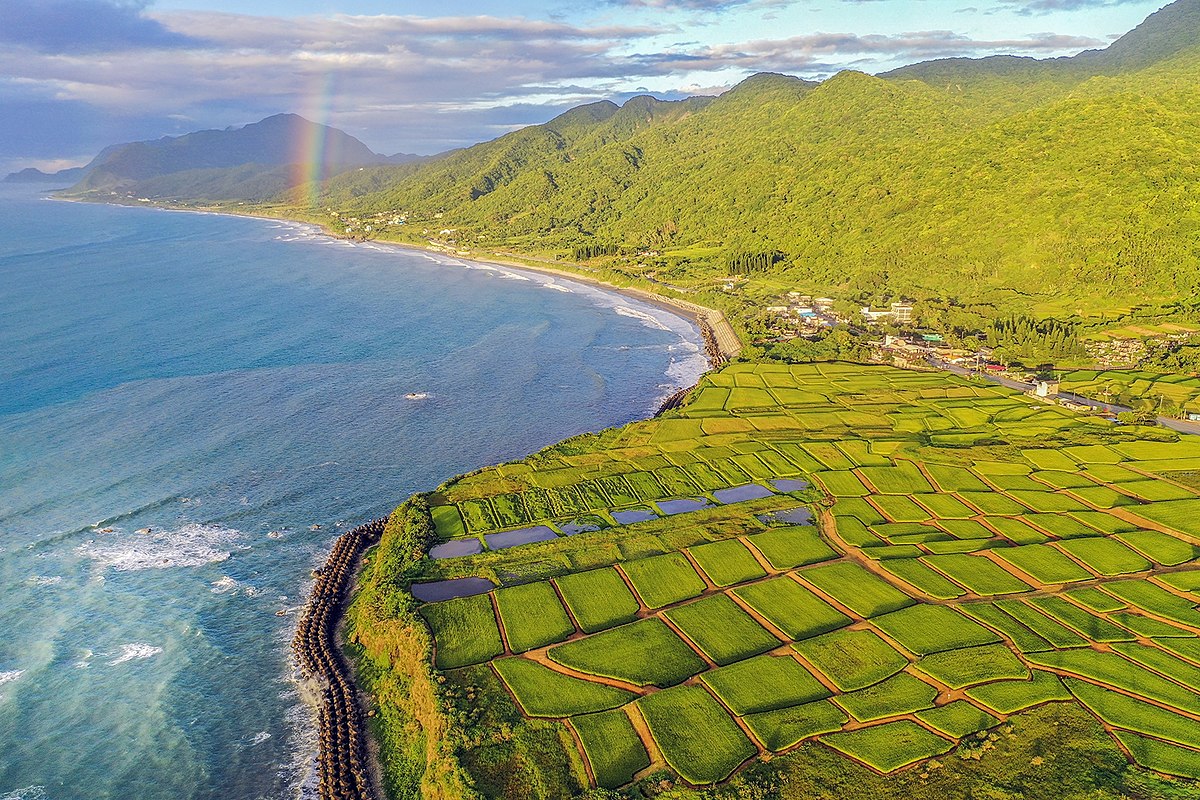}$\quad$
                                \includegraphics[height=5.6cm]{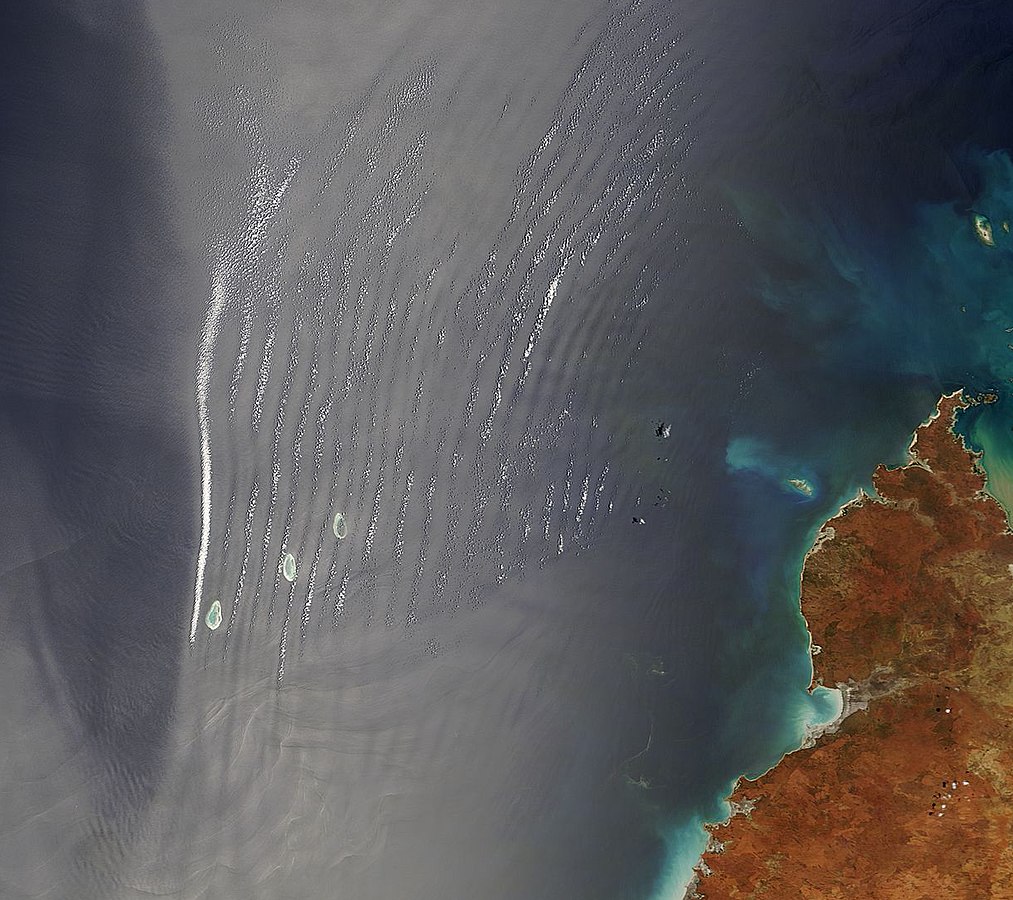}
        \caption{\sl Left: the coast of Fengbin, Taiwan; right: satellite view of the
Shark Bay, Western Australia
(images from Wikipedia, photo by Wu Pei Hsuan
licensed under the Creative Commons Attribution-Share Alike 4.0 International license for the first,
Public Domain 
photo by Jeff Schmaltz, NASA, 
for the second).}\label{H8942kkXcoARGIRA4AXELHARLROA7789GIJ7solFUMHDNOJHNFOJEDappro}
\end{figure}

Additionally, from~\eqref{TH:EVE}, we have that~$u_1(\xi,\tau,z)=\eta_0(\xi,\tau)+\upsilon(\tau,z)$,
for an auxiliary function~$\upsilon$.
As a matter of fact, if we have a solution of~\eqref{KDV-0m-004}
and we add to~$u_1$ a function of~$(\tau,z)$ we obtain another solution. Therefore,
we choose the ``simplest possible'' solution by picking~$\upsilon$ identically equal to zero, thus obtaining
\begin{equation}\label{TH:EVE:2}
u_1(\xi,\tau,z)=\eta_0(\xi,\tau).
\end{equation}

Now we consider the order related to~$\e^2$ in the expansion in~\eqref{KDV-0m-002}
(the higher orders being formally neglected). In this way, we find that
\begin{equation}\label{KDV-0m-008}
\begin{dcases}
-\partial_\xi  u_2+\partial_\tau u_1+
u_1\partial_\xi u_1+\frac{w_1\partial_z u_1}\Theta=-\partial_\xi P_2&{\mbox{ in }}\Omega_\tau,\\
-\Theta\partial_\xi w_1=-\partial_z P_2&{\mbox{ in }}\Omega_\tau,\\
\Theta\partial_\xi u_2+\partial_z w_2=0&{\mbox{ in }}\Omega_\tau,\\
\partial_zP_0(\xi,1,\tau)\eta_1+
\displaystyle\frac12\partial_{zz}P_0(\xi,1,\tau)\eta_0^2+\partial_zP_1(\xi,1,\tau)\eta_0+P_2(\xi,1,\tau)=
\eta_1, & \\
\partial_zw_1(\xi,1,\tau)\eta_0+w_2(\xi,1,\tau)=
\Theta\partial_\tau\eta_0+\Theta u_1\partial_\xi\eta_0-\Theta\partial_\xi\eta_1,& \\
w_2(\xi,0,\tau)=0. &
\end{dcases}
\end{equation}
Recalling~\eqref{KDV-0m-005} and~\eqref{KDV-0m-009}, the fourth line in~\eqref{KDV-0m-008} becomes
$$ P_2(\xi,1,\tau)=
\eta_1(\xi,\tau).$$

Also, from the second line in~\eqref{KDV-0m-008} and~\eqref{KDV-0m-006},
$$ \partial_zP_2=-\Theta^2\,\partial_\xi(z\partial_\xi\eta_0)=-z\Theta^2\,\partial_{\xi \xi}\eta_0.
$$
These observations entail that, for all~$z\in(0,1)$,
\begin{eqnarray*}&&\eta_1(\xi,\tau)-P_2(\xi,z,\tau)=
P_2(\xi,1,\tau)-P_2(\xi,z,\tau)=
\int_z^1\partial_zP_2(\xi,z',\tau)\,dz'\\&&\qquad=-\Theta^2\int_z^1 z'\partial_{\xi \xi}\eta_0(\xi,\tau)\,dz'=-\frac{\Theta^2\,(1-z^2)}2\partial_{\xi \xi}\eta_0(\xi,\tau)
\end{eqnarray*}
and therefore
$$ P_2(\xi,z,\tau)=\eta_1(\xi,\tau)+\frac{\Theta^2\,(1-z^2)}2\partial_{\xi \xi}\eta_0(\xi,\tau).$$
{F}rom this and the first equation in~\eqref{KDV-0m-008} we infer that
\begin{eqnarray*}
\partial_\xi  u_2&=&
\partial_\tau u_1+u_1\partial_\xi u_1+\frac{w_1\partial_z u_1}\Theta+\partial_\xi P_2\\
&=&
\partial_\tau u_1+u_1\partial_\xi u_1+\frac{w_1\partial_z u_1}\Theta+\partial_\xi \eta_1+\frac{\Theta^2\,(1-z^2)}2\partial_{\xi\xi \xi}\eta_0.
\end{eqnarray*}
Notice the appearance of a third derivative for~$\eta_0$. Thus, recalling
the third equation in~\eqref{KDV-0m-008},
\begin{eqnarray*}
-\frac{\partial_z w_2}\Theta=
\partial_\xi u_2=
\partial_\tau u_1+u_1\partial_\xi u_1+\frac{w_1\partial_z u_1}\Theta+\partial_\xi \eta_1+\frac{\Theta^2\,(1-z^2)}2\partial_{\xi\xi \xi}\eta_0.
\end{eqnarray*}

\begin{figure}
                \centering
                \includegraphics[height=5.6cm]{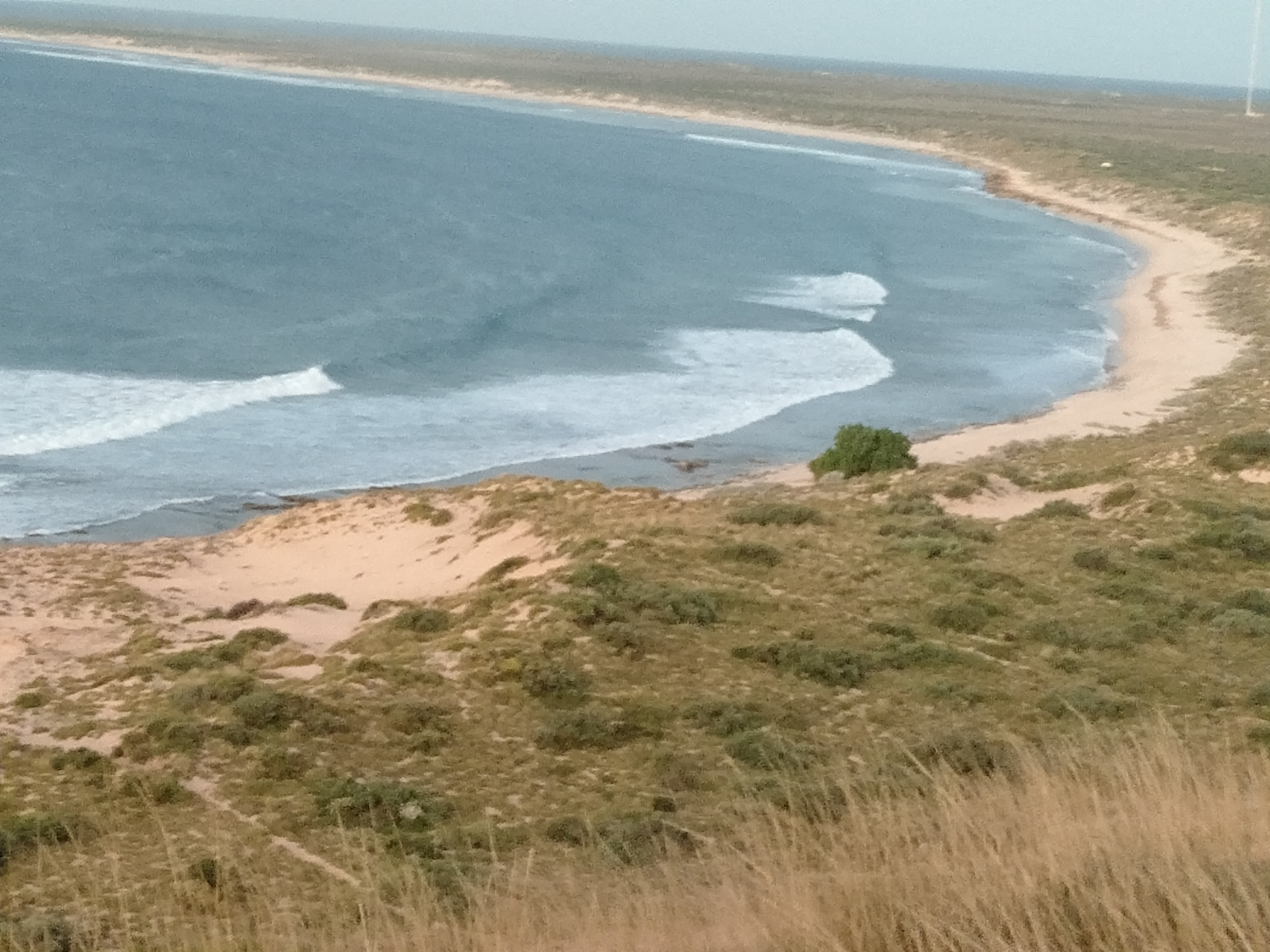}
        \caption{\sl Wave fronts approaching the shore at North West Cape.}\label{H8942kkXcoARGIRA4AXELHARLROA7789GIJ7solFUMHDNOJHNFOJEDappro-99}
\end{figure}

Combining this with the last two equations in~\eqref{KDV-0m-008} we conclude that
\begin{eqnarray*}&&
\partial_\tau\eta_0(\xi,\tau)+u_1(\xi,1,\tau)\partial_\xi\eta_0(\xi,\tau)-\partial_\xi\eta_1(\xi,\tau)-\frac{\partial_zw_1(\xi,1,\tau)\eta_0(\xi,\tau)}\Theta\\&=&\frac{w_2(\xi,1,\tau)}\Theta\\
&=& \frac{w_2(\xi,1,\tau)-w_2(\xi,0,\tau)}\Theta
\\&=&\frac1\Theta\int_0^1 \partial_zw_2(\xi,z,\tau)\,dz\\&=&
-\int_0^1 \Bigg(
\partial_\tau u_1(\xi,z,\tau)+u_1(\xi,z,\tau)\partial_\xi u_1(\xi,z,\tau)+\frac{w_1(\xi,z,\tau)\partial_z u_1(\xi,z,\tau)}\Theta\\&&
\qquad\qquad+\partial_\xi \eta_1(\xi,\tau)+\frac{\Theta\,(1-z^2)}2\partial_{\xi\xi \xi}\eta_0(\xi,\tau)
\Bigg)\,dz
.\end{eqnarray*}
This and~\eqref{KDV-0m-006} yield that
\begin{eqnarray*}&&
\partial_\tau\eta_0(\xi,\tau)+u_1(\xi,1,\tau)\partial_\xi\eta_0(\xi,\tau)-\partial_\xi\eta_1(\xi,\tau)+\partial_\xi\eta_0(\xi,\tau)\eta_0(\xi,\tau)\\&=&
-\int_0^1 \left(
\partial_\tau u_1(\xi,z,\tau)+u_1(\xi,z,\tau)\partial_\xi u_1(\xi,z,\tau)-z\partial_\xi\eta_0(\xi,\tau)\partial_z u_1(\xi,z,\tau)
\right)\,dz\\&&\qquad\qquad-\partial_\xi \eta_1(\xi,\tau)-\frac{\Theta^2}3\partial_{\xi\xi \xi}\eta_0(\xi,\tau).
\end{eqnarray*}

\begin{figure}
                \centering
                \includegraphics[height=7.6cm]{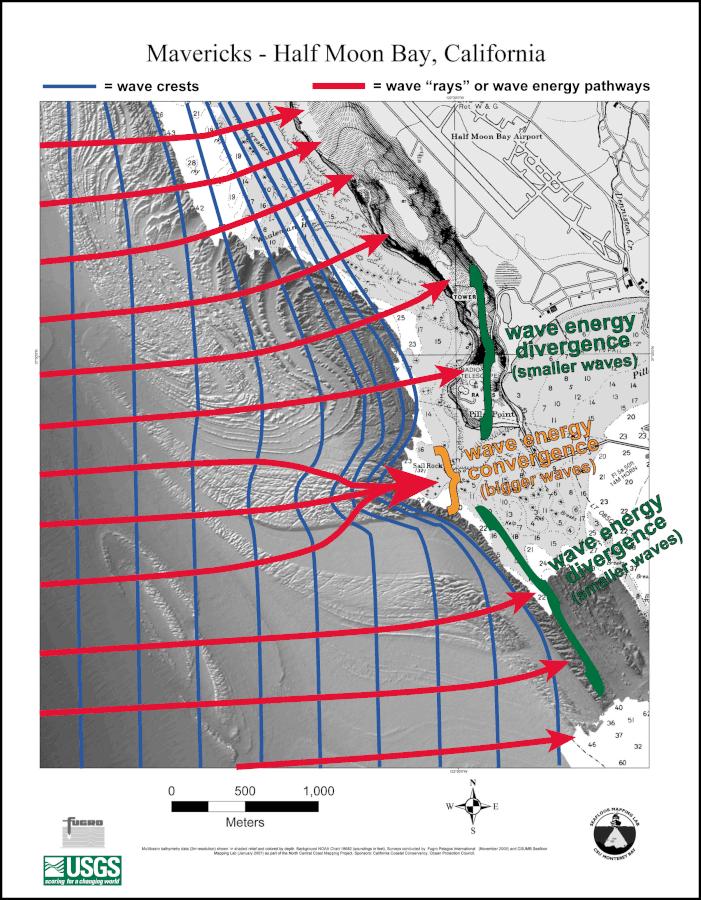}
        \caption{\sl Wavefronts at Mavericks, California
        (Public Domain image from
        Wikipedia).}\label{PDH8942kkXcoARGIRA4AXELHARLROA7789GIJ7solFUMHDNOJHNFOJEDappro}
\end{figure}

Then, we can simplify the term~$-\partial_\xi\eta_1(\xi,\tau)$ and exploit~\eqref{TH:EVE:2}, finding that
\begin{eqnarray*}&&
\partial_\tau\eta_0(\xi,\tau)+\eta_0(\xi,\tau)\partial_\xi\eta_0(\xi,\tau)+\partial_\xi\eta_0(\xi,\tau)\eta_0(\xi,\tau)\\&=&
-\int_0^1 \left(
\partial_\tau \eta_0(\xi,\tau)+\eta_0(\xi,\tau)\partial_\xi \eta_0(\xi,\tau)
\right)\,dz-\frac{\Theta^2}3\partial_{\xi\xi \xi}\eta_0(\xi,\tau)\\&=&
-\partial_\tau \eta_0(\xi,\tau)-\eta_0(\xi,\tau)\partial_\xi \eta_0(\xi,\tau)
-\frac{\Theta^2}3\partial_{\xi\xi \xi}\eta_0(\xi,\tau),
\end{eqnarray*}
which is the Korteweg-de Vries equation
\begin{equation}\label{VJlaVDa8jr3} \partial_\tau\eta_0+ \frac32\partial_\xi\eta_0\eta_0+\frac{\Theta^2}6\partial_{\xi\xi \xi}\eta_0=0,\end{equation}
to be compared with~\eqref{KDV-or}, as desired.

\begin{figure}
                \centering
                \includegraphics[width=.49\linewidth]{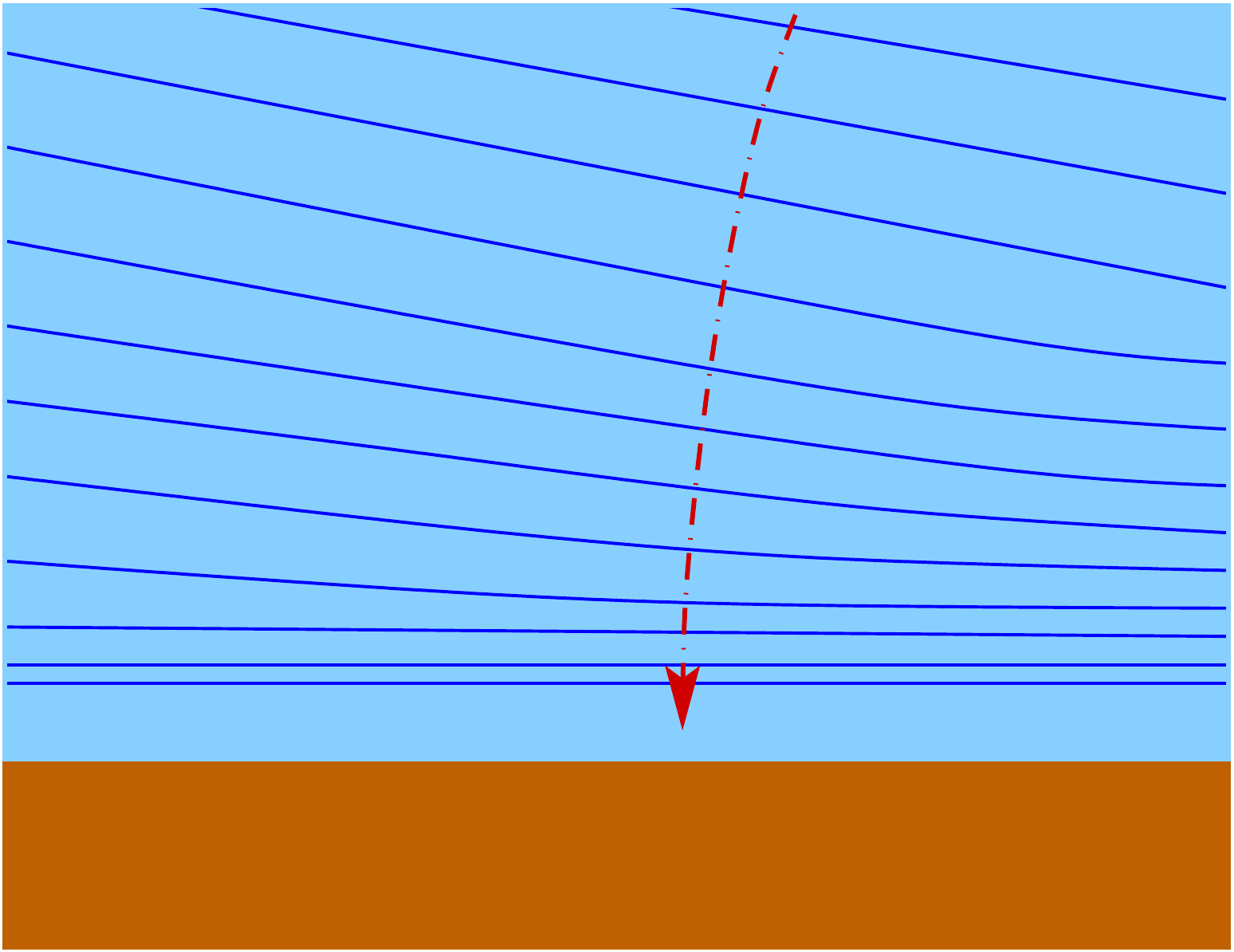}
        \caption{\sl Waterfronts approaching the coastline.}\label{HcoARGIRA4AXELHARLROA7789GIJ7solFUMHDNOJHNFOJEDappro}
\end{figure}

It is now instructive to recover equation~\eqref{VJlaVDa8jr3} in terms of the original variables. For this, owing to~\eqref{RTArasmsImPNSd-2} and~\eqref{RTArasmsImPNSd-1},
\begin{equation}\label{VJlaVDa8jr3dra6mo9}\eta(\xi,\tau)=
\eta^\star\left(\xi+\frac\tau\e,\frac\tau\e\right)=\overline\eta\left(\sqrt{gh_0} \left(\xi+\frac\tau\e\right),\frac\tau\e\right).\end{equation}
Moreover, since~$\overline\eta(x,t)$ was the profile of the free surface of the fluid
in the canal and~$\eta_0(\xi,\tau)$ the leading order of this profile in the new variables~$(\xi,\tau)$
in a shallow water and low amplitude approximation, in light of~\eqref{VJlaVDa8jr3dra6mo9} we can consider the leading order of
the original profile~$\overline\eta_0$ as defined implicitly by the relation
$$ \eta_0(\xi,\tau)=\overline\eta_0\left(\sqrt{gh_0} \left(\xi+\frac\tau\e\right),\frac\tau\e\right).$$
Therefore~\eqref{VJlaVDa8jr3} becomes
\begin{eqnarray*}\frac{
\sqrt{gh_0}}\e\partial_x\overline\eta_0+\frac1\e\partial_t\overline\eta_0+ \frac{3\sqrt{gh_0}}2\partial_x\overline\eta_0\overline\eta_0+\frac{ (gh_0)^{3/2}\,\Theta^2}6\partial_{x}^3\overline\eta_0=0
\end{eqnarray*}
and thus, recalling~\eqref{EPSIDELKDV},
\begin{equation}\label{BEsF3647cehHNMMauJDR6tEmwfr9w0}
\partial_t\overline\eta_0+\frac{ \sqrt{g} h_0^{5/2}}{6}\partial_{x}^3\overline\eta_0
+ \frac{3h_0^{3/2}}{2\Theta^2\,\sqrt{g}}\partial_x\overline\eta_0\overline\eta_0+\sqrt{gh_0} \partial_x\overline\eta_0=0.
\end{equation}
In this expression, one can get rid of the term~$\partial_x\overline\eta_0$ by a vertical translation of~$\overline\eta_0$,
that is, if
\begin{equation}\label{ZETATRAVEKDG32}\zeta_0(x,t):=\alpha\left(\overline\eta_0(x,t)+\frac{2\Theta^2\,g}{3h_0}\right)
=\frac{h_0^2}{\Theta^2\,g} \left(\overline\eta_0(x,t)+\frac{2\Theta^2\,g}{3h_0}\right),\end{equation} it follows that
\begin{equation}\label{UHJNZzekmdLSyhn8ijP235}
\partial_t\zeta_0+\frac{ \sqrt{g} h_0^{5/2}}{6}\partial_{x}^3\zeta_0
+ \frac{3\sqrt{g}}{2\sqrt{h_0}}\partial_x\zeta_0\zeta_0=0,
\end{equation}
which corresponds\footnote{Retaking~\eqref{d3ae9ikjm3f82KSMpf334}, we know that
the velocity of the traveling wave~$\phi_0$ is proportional to~$\frac{c}a$ times the amplitude of the wave. \label{BEsF3647cehHNMMauJDR6tEmwfr9w0-PA}
It is suggestive to compare this information with the numerology of the coefficients
and the relation with the physical shape of the wave. Namely, in this case we have that~$\frac{c}a=c=\frac{3\sqrt{g}}{2\sqrt{h_0}}$
and, by~\eqref{ZETATRAVEKDG32}, the amplitude of the traveling wave would formally correspond to~$
\frac{h_0^2}{\Theta^2\,g} \left(\max\overline\eta_0+\frac{2\Theta^2\,g}{3h_0}\right)=\frac{h_0^2}{\Theta^2\,g} \max\overline\eta_0+\frac{2h_0}{3}$.
The velocity of the traveling wave would thus be proportional to~$\frac{3\sqrt{g}}{2\sqrt{h_0}}\left(
\frac{h_0^2}{\Theta^2\,g} \max\overline\eta_0+\frac{2h_0}{3}\right)=
\frac{3h_0^{3/2}}{2\Theta^2\,\sqrt{g}} \max\overline\eta_0+{\sqrt{gh_0}}
$. Once again, this suggests that the speed of small amplitude waves in a shallow canal is dictated by a quantity proportional to~${\sqrt{gh_0}}$,
plus an increasing term in the amplitude of the wave.

The leading term proportional to~${\sqrt{gh_0}}$ for the speed of the wave, which gets smaller in
shallower water near the coast, also suggests that waves have the
tendency to slow down as they approach the shore because of the force exerted on them by the seabed. 

This phenomenon is somewhat related to the fact that, in our experience,
the waves coming to the shore have the strong tendency to follow the shape of the coastline, see e.g. Figures~\ref{H8942kkXcoARGIRA4AXELHARLROA7789GIJ7solFUMHDNOJHNFOJEDappro}
and~\ref{H8942kkXcoARGIRA4AXELHARLROA7789GIJ7solFUMHDNOJHNFOJEDappro-99}.

This effect can be explained in terms of the slowing process induced by shallow waters.
Indeed, the waves offshore can well be oblique with respect to the coastline, but
when approaching the shoreline the part of the wave closest to the coast is
affected sooner by the slowing effect of shallow waters, while the farthest wavefront continues its run for a while with the same speed: this relative variation of speed causes the ridges of the waves to gradually arrange themselves in an
almost parallel fashion to the coastline, see Figures~\ref{PDH8942kkXcoARGIRA4AXELHARLROA7789GIJ7solFUMHDNOJHNFOJEDappro}
and~\ref{HcoARGIRA4AXELHARLROA7789GIJ7solFUMHDNOJHNFOJEDappro}.

Of course, we do not aim here at exhausting the study of the
complicated topic of water waves' speed. See e.g.~\cite{MR3744219, w13212990}
and the references therein for further details on this subject.}
to~\eqref{KDV-or} with~$a:=1$, $b:=\frac{ \sqrt{g} h_0^{5/2}}{6}$ and~$c:=\frac{3\sqrt{g}}{2\sqrt{h_0}}$.\medskip

It is also worth reconsidering the solitary waves
of the Korteweg-de Vries equation found in~\eqref{DIKDVSOLIFI-EQ}
and note that these correspond, in principle, to waves with arbitrarily large velocity~$v$ and arbitrarily large amplitude (proportional to~$v$,
see~\eqref{d3ae9ikjm3f82KSMpf334}): the large amplitude solutions however are, strictly speaking, not consistent with the technical derivation of the model, which relies on expansions that are valid only for small amplitude waves (recall the setting in~\eqref{EPSIDELKDV}).
We also remark that, since the coefficient~$\sqrt{\frac{av}{b}}$ in~\eqref{DIKDVSOLIFI-EQ}
is monotone increasing in~$v$, and since this coefficient modulates the hyperbolic cosine at the denominator,
faster and higher waves in~\eqref{DIKDVSOLIFI-EQ} are also narrower.

The fact that higher waves (at least in this setting)
need to be narrower can be guessed also by the fact that,
at least for waves ``vanishing at infinity'' and when~$a\ne0$, the total mass of the wave~$m(t):=\int_\R \phi(x,t)\,dx$ is a conserved quantity for the
Korteweg-de Vries equation~\eqref{KDV-or}, since
\begin{eqnarray*}
\dot m(t)&=&\int_\R \partial_t\phi(x,t)\,dx\\&=&
-\frac{b}a\int_\R\partial _{x}^{3}\phi(x,t)\,dx-\frac{c}{a}\int_\R \phi(x,t) \,\partial _{x}\phi (x,t)\,dx\\&=&
-\frac{b}a\int_\R\partial_x\Big(\partial _{x}^{2}\phi(x,t)\Big)\,dx-\frac{c}{2a}\int_\R \partial _{x}
\Big(\phi^2 (x,t)\Big)\,dx\\&=&
-\frac{b}a\Big( \partial _{x}^{2}\phi(+\infty,t)-\partial _{x}^{2}\phi(-\infty,t)\Big) -\frac{c}{2a}
\Big(\phi^2 (+\infty,t)-\phi^2 (-\infty,t)\Big)\\&=&0.
\end{eqnarray*}

For further information\footnote{For simplicity, we did not discuss the role of the surface tension
in the derivation of the Korteweg-de Vries equation, but we mention that this additional
term can be also easily included in the coefficient~$b$ in~\eqref{KDV-or}.} about the Korteweg-de Vries equation
and related topics see e.g.~\cite{MR271526, MR0352741, zbMATH03735376, MR856894, MR985322, MR1462745, MR1629555}
and the references therein. \medskip

\begin{figure}
                \centering
                \includegraphics[width=.92\linewidth]{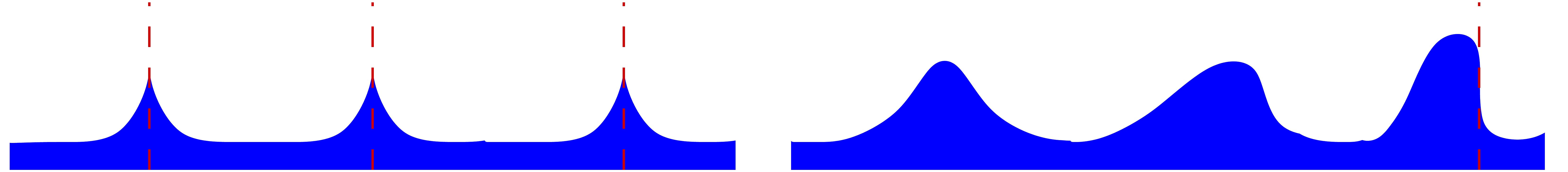}
        \caption{\sl Breaking waves: 
        sharp crests versus waves which break into bores.}\label{C12435tHIMSINghItangeFI}
\end{figure}

Since we are speaking here about waves in a shallow canal,
let us spend a few, certainly not exhaustive, words about the phenomenon of wave breaking, \index{breaking waves}
that is the case in which a solution of a certain equation somewhat modeling water waves
remains bounded\footnote{On a different note, we mention that
other possible singularities
are the ones in which crests exhibit a corner. These may
occur for a set of equations
similar to, but structurally different from, the ones presented here.
The theory of some of these equations go back to Stokes
and several remarkable properties hold true about the
maximal periodic deep-water waves arising in this context,
e.g. the wave height is about one-seventh of the wavelength
and the maximum height waves have a~$120^\circ$ corner.
See~\cite{MR513927, MR642535} for further details about these types of problems.

Traveling waves with corners on their top edge
(that is, with discontinuous first derivatives) are sometimes called \index{peakon}
``peakons''.

Another type of equation which develops peakons and singularities
(i.e., moving waves with corners as well as waves whose slope
becomes vertical in finite time, provided it is initially
sufficiently negative) is the so-called
Camassa-Holm equation, \index{Camassa-Holm equation}
see~\cite{MR1234453}.}
but its slope becomes infinite in finite time.
In particular, one would aim at detecting breaking phenomena\footnote{The classification of breaking waves is not homogeneous across the scientific literature
and other types of wave breaking have been classified as well,
see e.g.~\cite[page~42]{zbMATH00054189} and~\cite[page~60]{lin2008numerical}.}
that
correspond to ``sharp crests'' (namely, waves with cuspidal tops)
or to waves which ``break into bores'' (namely, waves~$u$ which travel, say, from left to right,
starting from a configuration with~$u'$ sufficiently negative at some point in the space,
whose profile steepens until~$u'$ becomes~$-\infty$ in finite time), see Figure~\ref{C12435tHIMSINghItangeFI}.
However, the existence and the properties of breaking waves
heavily depend on the specific partial differential equation chosen to model the phenomenon.
For instance, the Korteweg-de Vries equation does not allow breaking waves,
as proved in~\cite{MR1824796, MR1969209}. \label{GITARFKMAOSTGABafOKS02o3t}
The heuristic reason for which the Korteweg-de Vries equation prevents wave braking
is the third derivative term in~\eqref{KDV-or}, which provides\footnote{Quite surprisingly,this
smoothing effect will play a fundamental role in the evolution phenomena
of chain of oscillators, see the forthcoming discussion
on page~\pageref{GITARFKMAOSTGABafOKS02o3t-2} for full details.}
a smoothing effect. As a matter of fact,
by dropping this term,
one obtains another equation, namely the
Bateman-Burgers equation, \index{Bateman-Burgers equation}
which will be introduced below in~\eqref{VERBURG2},
and which develops singularities in finite time: however,
as clarified in~\cite[page~457]{MR1699025},
the theory of the Bateman-Burgers equation
(that is, of the Korteweg-de Vries equation without third derivative term)
``goes too far: it predicts that all waves carrying an increase of
elevation break. Observations have long since established that some waves
do not break. So an invalid theory seems to be right sometimes and wrong
sometimes!'' \index{Korteweg-de Vries equation|)}\medskip

To develop a mathematical theory of breaking waves,
several modifications of the Korteweg-de Vries equation have been introduced
in the literature (for instance by replacing the third derivative term
with a nonlocal or fractional term, see~\cite{MR208903, MR3291137, MR3682673, MR3749383}).
Yet, the matter appears to be not completely settled.
For instance, several attempts have been made to
determine the breaking wave height~$\alpha$
with respect to the water depth~$h$ and\footnote{The ratio between the height of the breaking wave
and the depth of the bottom is sometimes referred to with the name of \index{breaker index}
``breaker index''.}
the ``numerology''
does not always match conjectures, theory and experiments.
As an example, in~\cite[page~14]{SHORT-b},
one finds that the water depth corresponds to~$1.5$ times the breaking wave height, giving~$\frac\alpha{h}=\frac23=0.\overline{6}$.
The problem has an ancient tradition:
already in~\cite[equation~(52)]{doi:10.1080/14786449108621390},
John McCowan, a pioneer in the study of the fluid mechanics behind surfing,
calculated~$\frac\alpha{h}\simeq0.9$, concluding that ``the wave will break for an elevation rather less than
the mean depth'' believing this condition for breaking ``not be far wrong. Scott Russell's
experiments confirm this: he found that the wave broke when
the elevation was about equal to the depth; but from some
experiments of my own I am inclined to think that~$3/4$ is a closer approximation for the elevation at the breaking-point''. This experiment, determining~$\frac\alpha{h}\simeq0.75$,
was retaken from the theoretical point of view in~\cite[equation~(34)]{doi:10.1080/14786449408620643}
proposing at this stage the value~$\frac\alpha{h}\simeq0.78$.
This value has been commented about in~\cite[page~39]{zbMATH00054189} as
``appropriate for estimating the breaker height in mild-sloped surf-zones''.

The reliability of these numbers is possibly a bit controversial, see~\cite{KisacikTroch2014},
according to which
``\cite{GODAA} mentions  that  a 
value of~$0.8261$ is more accurate~\cite{YAMADA}. {F}rom the field observations, 
it is found to be between~$0.78$ and~$0.86$''.
But~\cite[page~29]{IppenKulin1954} lists
five different theoretical findings for the ratio between breaking wave elevation
and depth of the bottom, ranging from~$0.73$ to~$1.03$.
Numerous laboratory experiments have been performed to
empirically quantify the wave breaking phenomenon,
also by installing offshore structures resembling the
actual coastal terrain, but most of the formulas found depend significantly
on the specific laboratory experiment data,
see e.g.~\cite[Table~1]{j.1749-6632.1949.tb27281.x} for the effect of the beach characteristics on break height.
See also~\cite{jmse9070731, jmse10010050} for a machine learning viewpoint to the problem.
Other approaches to describe wave breaking relate to equations with variable\footnote{It is possible
that the steepness of the bottom also plays a role in the determination of the height
of a breaking wave: according to~\cite{10.211205-0566.1},
``most formulas show quite good predictions for cases including gentle slopes,
however, the predictions are typically not satisfactory for breaking waves on steep slopes''.
This paper also presents a comparison of different approaches, methodologies and formulas
on the topic of breaking waves.}
depth, instead of flat bottoms, see e.g.~\cite{MR2655953}.
In general, a common belief is that the phenomena
related to wave breaking are very complex and are linked to many
additional properties of the environment, such as the
shape of the land (sometimes referred to with the name of
``topography'') and the specific features of the bottom
(such as impermeability, presence of irregularities, slope oscillations).
\medskip

In essence, about two centuries after the solitary wave of translation observed by John Scott Russell, a unified and all-encompassing theory of water waves with a outright match between rigorous results and experimental findings is, perhaps, still not available. 

\subsection{Plasma physics}\label{PLASECT-0uorwehfg}

\begin{figure}
                \centering
                \includegraphics[width=.35\linewidth]{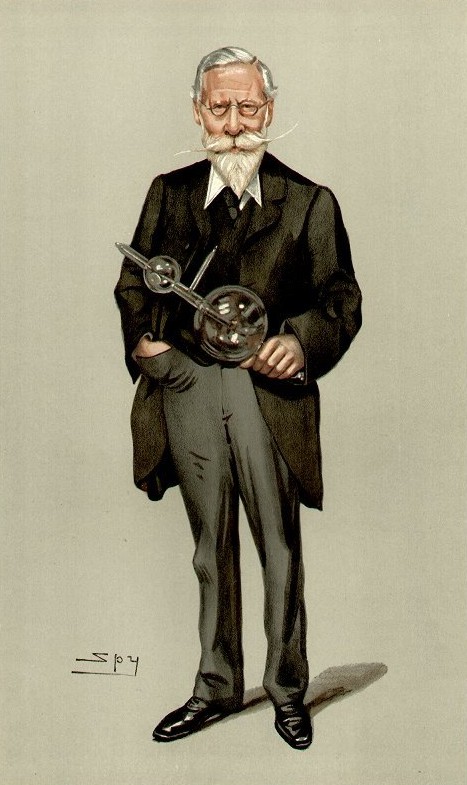}
        \caption{\sl Caricature of Sir William Crookes (Public Domain image from
        Wikipedia).}\label{HARGIRA4AXELHARLROA7789GIJ7solFUMHDNOJHNFOJED}
\end{figure}

A state of matter is one of the distinct forms in which matter can exist in nature.
Classically, nature was classified
into three states of matter: solid, liquid, and gas.
Nowadays it is common\footnote{Also, a large number of
intermediate states are nowadays known to exist, such as 
liquid crystals,
Bose-Einstein condensates, 
Fermionic condensates,
neutron-degenerate matter, quark-gluon plasma, superfluids, superconductors, etc., but
we do not investigate these states here.

{F}rom the historical point of view,
plasma was first identified in laboratory by Sir William Crookes in 1879
(though more systematic studies of plasma only began in the twentieth century).

Besides discovering plasma, Crookes is also credited with discovering
thallium, inventing a radiometer and, quite conveniently, some~$100\%$ 
ultraviolet blocking sunglass lens.

He was also interested in spiritualism and became president of the Society for Psychical Research,
whose purpose is to understand psychic and paranormal events.

See Figure~\ref{HARGIRA4AXELHARLROA7789GIJ7solFUMHDNOJHNFOJED}
for a caricature of Sir William Crookes (by caricaturist Sir Leslie Ward).}
to consider also a fourth state
of matter, namely plasma\index{plasma}.

Though our everyday experience with plasma appears to be rather limited
(we are experienced to plasma mostly only by
neon tubes, plasma lamps and lightnings),
plasma is allegedly the most abundant form of ordinary matter in the universe
(excluding dark matter and dark energy) and likely amounts to about~$
99\%$ of the total mass of the visible universe. 
It consists of a gas of ions and free electrons, usually produced by
very high temperatures ($10^4$ degrees or more)
which make electrons leave
their orbit around the corresponding nuclei. In this way,
if the initial status of the gas is overall in an electrical neutral state,
the plasma consists in a mixture of charged particles, namely the free electrons
and the corresponding ions.\medskip

The mathematical description of a plasma is a highly advanced subject but,
in a nutshell, there is a hierarchy of models accounting for the
the evolution of a plasma. 

The most intuitive one would be to describe the plasma
by detecting positions and velocities of all its particles (ions and electrons) at a given time.

This approach is very precise but often impractical, given the high number of particles involved,
and therefore intermediate
models (often called ``kinetic models'') have been introduced
to describe the plasma ``in average'' through the statistical analysis of the particle
distribution.

Moreover, at a large scale, efficient models for a plasma leverage the knowledge of
fluid dynamics equations (somewhat close in the spirit to the ones that we 
presented in Sections~\ref{EUMSD-OS-32456i7DNSSE234R:SEC} and~\ref{BARPotA}):
namely, at a macroscopic level, close to thermodynamic equilibrium,
one can identify
each species of particles of a plasma (ions and electrons) with a fluid
and describe the corresponding density,
velocity and energy via a set of\footnote{To confirm the importance of partial differential
equations in our understanding of plasma, let us mention for instance that the
laboratory specialized in plasma and located in Toulouse, France,
is named LAboratoire PLAsma et Conversion d'Energie (LAPLACE).} partial differential equations.\medskip

To make this model concrete, we denote by an index~$j\in\{1,2\}$
the species of particles of a given plasma (e.g., $j=1$ corresponding to ions and~$j=2$ corresponding to
electrons). We suppose that the plasma, being neutral in average
but composed at a small scale by charged particles, is subject to its own magnetic and electric fields (denoted by~$B$
and~$E$, respectively).
We recall that charged particles are influenced by both the electric and the magnetic fields: more specifically,
the so-called
Lorentz force acting on a single charged particle (say, with charge~$q_j$ and velocity~$v_j$)
is of the form~$q_j(E+v_j\times B)$. Therefore the total contribution of the Lorentz force acting on a portion of plasma with
particle density~$\mu_j$ takes the form
\begin{equation}\label{TOTAFO-1}
\mu_jq_j(E+v_j\times B).\end{equation}
Additionally, the possible variations of density exert a force on the plasma similar to a pressure. Thus, we denote by \begin{equation}\label{TOTAFO-2}
-\nabla p_j\end{equation} this type of pressure (the minus sign
manifesting the fact that higher densities oppose the motion).

Differently from the case of classical fluids, in the description of a plasma we have also to account for
the effect of possible particle collisions. Since we are looking here at a macroscopic description
of the plasma, we suppose that the total momentum of the particles of species~$j$
comes from the average over all particles of that species and therefore
like particle collisions do not change the total momentum (say, if two electrons collide,
one may slow down after collision and the other can be accelerated, but the total momentum
is preserved, thanks to Newton's Law; same if two ions collide).
We have therefore to focus on the collisions between unlike particles,
which allow momentum to be exchanged between the species;
in this situation, when an electron and a ion collide, the momentum loss of one particle entails
a corresponding momentum gain of the other particle:
that is if the first particle experiences a variation of momentum given by~$P_1$,
the second particle experiences a variation of momentum given by~$P_2=-P_1$.

Computing the momentum exchange of colliding particles can be in general a rather complicated task.
A simplifying assumption in this framework is to assume that the particles are reduced to a point.
An additional simplification arises if we are willing to suppose that the collisions are perfectly elastic (no dissipation
of energy due to particle collisions).
Furthermore, we can also reduce the problem to a simpler one if we are willing to use the fact that the mass
of the electrons is typically much smaller than the mass of the ions and take therefore
the mass of the ions as the leading term of a more complex approximation.
To make these ideas work, at least in this simplified setting, we recall that the elastic collision
of two point masses, say with masses~$m_1$ and~$m_2$,
initial (before collision) velocity~$v_{1i}$ and~$v_{2i}$
and final (after collision) velocity~$v_{1f}$ and~$v_{2f}$
is fully described by the equations
$$
\begin{dcases}v_{1f}=\displaystyle{\frac{(m_{1}-m_{2})v_{1i}+2m_{2}v_{2i}}{m_{1}+m_{2}}},\\v_{2f}=\displaystyle{\frac{(m_{2}-m_{1})v_{2i}+2m_{1}v_{1i}}{m_{1}+m_{2}}},\end{dcases}
$$
see e.g.~\cite[equation~(17.2)]{MR0120782}.

Hence, since we consider the index~$j=1$ as corresponding to the electrons and we assume that~$m_1\ll m_2$,
we have that the momentum exchange for the colliding electron is
\begin{eqnarray*}
m_1(v_{1f}-v_{1i})&=&
m_1\left(\frac{(m_{1}-m_{2})v_{1i}+2m_{2}v_{2i}}{m_{1}+m_{2}}-v_{1i}\right)\\&=&
\frac{2m_1m_2(v_{2i}-v_{1i})}{m_{1}+m_{2}}\\
&\simeq&2m_1(v_{2i}-v_{1i}).
\end{eqnarray*}

In the framework of plasma collisions, this would give that~$P_1=2\mu_1 m_1(v_2-v_1)$
and accordingly~$P_2=2\mu_1 m_1(v_1-v_2)$.
Thus, if we denote by~$\phi$ the collision frequency between the two species
(for simplicity, we are disregarding here the collisions between ions, or electrons,
with the remaining neutral gas),
the species corresponding to~$j=1$ changes its momentum by a quantity~$\phi\mu_1 m_1(v_2-v_1)$
and the species corresponding to~$j=2$ changes its momentum by a quantity~$\phi\mu_2 m_2(v_1-v_2)$
(where constants are omitted for the sake of simplicity).

\begin{figure}
                \centering
                \includegraphics[width=.35\linewidth]{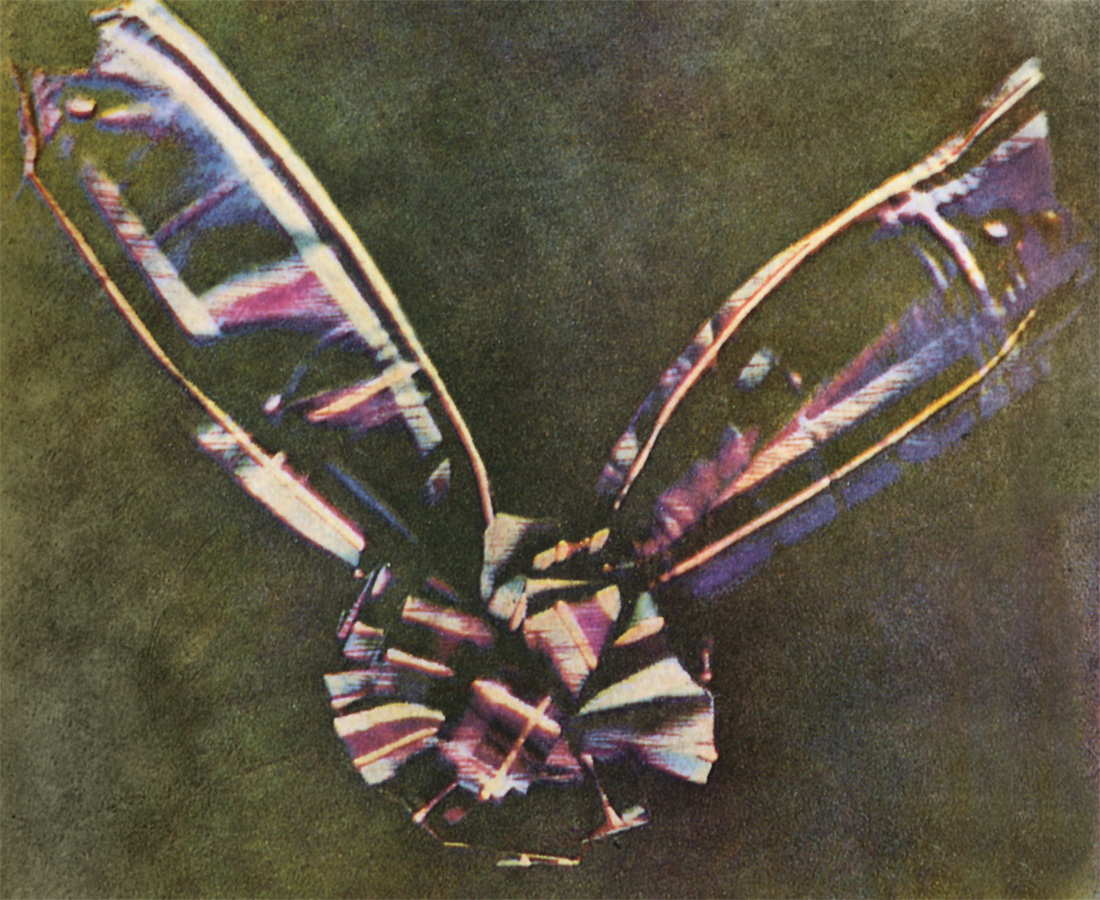}
        \caption{\sl The first durable color photographic image (Public Domain image from Wikipedia).}\label{GALAHAFOUMSldRGIRA4AXELHARLROA7789GIJ7solFUMHDNOJHNFOJEDMAXW}
\end{figure}

{F}rom this, \eqref{TOTAFO-1} and~\eqref{TOTAFO-2}, by Newton's Law of momentum balance have that
\begin{eqnarray*} \mu_j m_j\partial_t v_j +\mu_j m_j(v_j\cdot\nabla)v_j&=& \mu_jm_j\frac{d}{dt} (v_j)\\&=&
\mu_jq_j(E+v_j\times B)-\nabla p_j+\phi\mu_j m_j(v_{k_j}-v_j),
\end{eqnarray*}
where
$$ k_j:=\begin{dcases}
1 & {\mbox{ if }}j=2,\\
2 & {\mbox{ if }}j=1.
\end{dcases}$$
This equation is usually complemented by a continuity equation (corresponding to the conservation
of the total number of particles, which can be seen as the counterpart here of the mass transport equation
in~\eqref{OJS-PJDN-0IHGDOIUGDBV02ujrfMTE}),
an equation of state (that is a constitutive law relating pressure and particle density,
which can be seen as the counterpart of the barotropic flow description given on page~\pageref{BAROFL},
usually taken of the form~$p_j =\kappa \mu_j^\gamma$, for suitable positive constants~$\kappa$ and~$\gamma$)
and the classical
Maxwell's equations\footnote{James Clerk Maxwell
founded the theory of electromagnetism, proving that electric and magnetic fields travel through space as waves moving at the speed of light and argued that light itself is an undulation causing electric and magnetic phenomena. This great conceptual unification of light and electrical phenomena made Maxwell one of the nineteenth century scientists having the greatest impact on the subsequent developments of relativity and quantum field theory.

Maxwell also worked on the kinetic theory of gases and pioneered the early days of control theory. Also, he established that the rings of Saturn were made of numerous small particles and he is credited for presenting the first durable color photograph. For this, he had the idea of superimposing on a screen shots taken with red, green and blue filters: his result, depicting a tartan ribbon, is shown in Figure~\ref{GALAHAFOUMSldRGIRA4AXELHARLROA7789GIJ7solFUMHDNOJHNFOJEDMAXW}.
The photographic plates available at the time were almost insensitive to red and barely sensitive to green, hence the results in themselves were perhaps far from perfect, but Maxwell clearly knew this was after all just a rather minor detail in the Big Picture of Science
and he wrote
``by finding photographic materials more sensitive to the less refrangible rays, the representation of the colors of objects might be greatly improved''.} for the electric and magnetic fields:
omitting structural constants, this set of equations thus takes the form\footnote{Sometimes the setting in~\eqref{MAXWEUL}
is called in jargon ``Euler-Maxwell's equations''\index{Euler-Maxwell's equations}, since it mixes the classical Euler's formulation
of fluid dynamics presented in Section~\ref{EUMSD-OS-32456i7DNSSE234R:SEC} (after suitable corrections
to deal with particle collision) with the Maxwell's equations\index{Maxwell's equations} for electromagnetism.}
\begin{equation}\label{MAXWEUL}
\begin{dcases}
\mu_j m_j\partial_t v_j +\mu_j m_j(v_j\cdot\nabla)v_j=\mu_jq_j(E+v_j\times B)-\nabla p_j+\phi\mu_j m_j(v_k-v_j),\\
\partial_t\mu_j+\div(\mu_jv_j)=0,\\
p_j =\kappa \mu_j^\gamma,\\
\div B=0,\\
\div E=\mu_1q_1+\mu_2q_2,\\
\curl E=-\partial_t E,\\
\curl B=\mu_1q_1v_1+\mu_2q_2v_2+\partial_t E,
\end{dcases}
\end{equation}
with~$j\in\{1,2\}$.

For additional information about plasma physics, see e.g.~\cite{MR1368631, MR1960956, MR2779346, MR3184808} and the references therein.

\subsection{Galaxy dynamics}

\begin{figure}
                \centering
                \includegraphics[width=.65\linewidth]{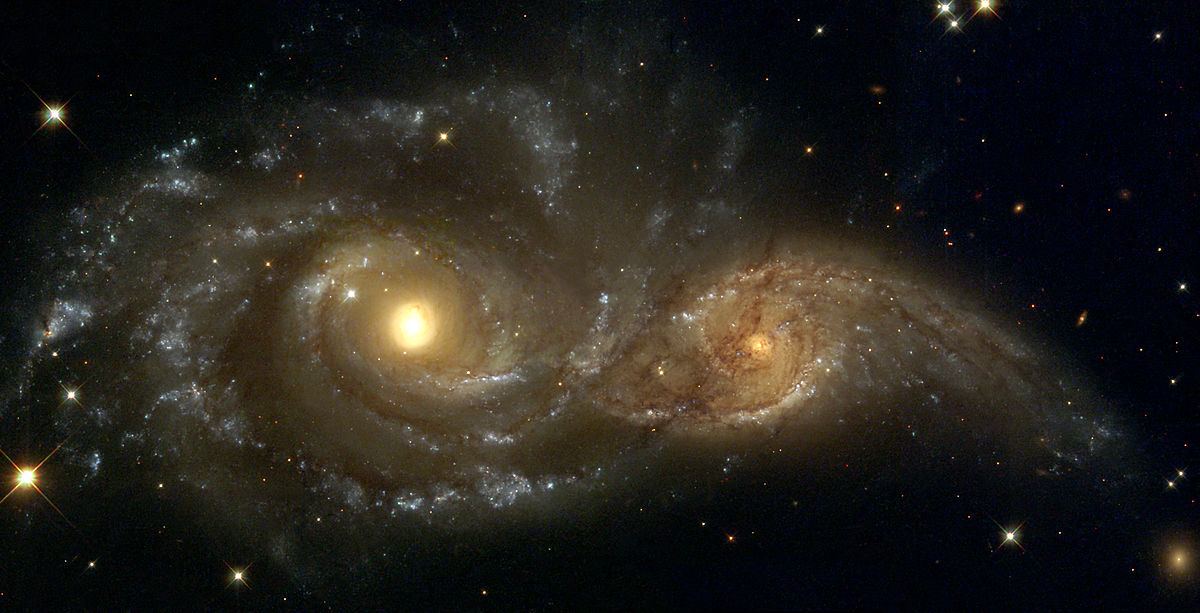}
        \caption{\sl Hubble Space Telescope picture of the galaxies NGC 2207 (left) and IC 2163 (right) (Public Domain image from Wikipedia).}\label{GALAHAFOUMSldRGIRA4AXELHARLROA7789GIJ7solFUMHDNOJHNFOJED}
\end{figure}

Interestingly, a suitable modification of the model used to describe the motion of fluids
presented in Section~\ref{EUMSD-OS-32456i7DNSSE234R:SEC}
can represent the dynamics of galaxies and globular clusters
(see e.g. Figure~\ref{GALAHAFOUMSldRGIRA4AXELHARLROA7789GIJ7solFUMHDNOJHNFOJED} for
a fascinating picture).
The ansatz of the model that we describe here is that the stars forming the galaxies interact only
by the gravitational field that they create collectively
and they can be described by parcels of a self-gravitating gas.
Neglecting collisional and relativistic effects, and disregarding the physical and chemical reactions
that continuously modify the internal structure of the stars, one
can provide a suitable set of equations describing the large scale galaxy motion.
These are named Jeans equations\index{Jeans equations} after Sir James Hopwood Jeans
(though were probably first derived by James Clerk Maxwell
and can be seen as a collisionless version of the Boltzmann equations
that describe the statistical behavior of a thermodynamic system not in equilibrium)
and are as follows:
\begin{equation}\label{STARS}
\begin{dcases}
\partial_t f+v\cdot\nabla_x f-\nabla u\cdot\nabla_v f=0,\\
\Delta u=\rho,
\end{dcases}
\end{equation}
where~$f=f(x,t,v)$ denotes the phase-space\footnote{By
phase-space here we mean the collection of position and velocity variables.}
density of stars, that is the density of stars corresponding to time~$t$, position~$x\in\R^3$ and velocity~$v\in\R^3$,
$u=u(x,t)$ is the gravitational potential induced collectively by the stars,
and
$$ \rho=\rho(x,t):=\int_{\R^3}f(x,t,v)\,dv$$
is the star density corresponding to time~$t$ and position~$x$.

To understand the rationale of equation~\eqref{STARS} we argue as follows.
In analogy with the dynamics of fluid parcels described in~\eqref{0uojfSANMFOSJFEL},
we assume that star dust trajectories follow the law~$\dot x(t)=v(x(t),t)$
and that their acceleration follows the gravitational law~$\ddot x(t)=-\nabla u(x(t),t)$.

One can therefore obtain an analogous of the continuity equation in~\eqref{OJS-PJDN-0IHGDOIUGDBV02ujrfMTE}
by assuming that, in absence of encounters, collisions
and collapses, the amount of stars is preserved by the flow.
Namely, the total quantity of stars occupying a region~$Z$ of the phase-space~$\R^3\times\R^3$ at a given time~$t_0$
is given by the quantity~$\int_{Z} f(x,t_0,v)\,dx\,dv$
and, in a small time~$\tau$,
possibly occupying a different region, that we name~$Z_\tau$ which collects all the evolution trajectories
in phase space given by
\begin{eqnarray*}
&& \Big(x(t_0+\tau),v\big(x(t_0+\tau),t_0+\tau\big)\Big)=
\big(x(t_0+\tau),\dot x(t_0+\tau)\big)=
(x(t_0),\dot x(t_0))+\tau (\dot x(t_0),\ddot x(t_0))+o(\tau)\\&&\qquad=
\Big(x(t_0),v\big(x(t_0),t_0\big)\Big)+\tau \Big(v\big( x(t_0),t_0\big),-\nabla u(x(t_0),t_0)\Big)+o(\tau).\end{eqnarray*}
Thus, if~$(x_0,v_0)\in Z$ and~$(x(\tau),z(\tau))=\Phi^\tau(x_0,v_0)$ denote the integral flow starting at~$(x_0,v_0)$, it follows that
$$ \left|\det D_{(x_0,v_0)} \Phi^\tau(x_0,v_0)\right|=
\left| \det D_{(x_0,v_0)}\Big( (x_0,v_0)+\tau
\big(v_0,-\nabla u(x_0,t_0)\big)\Big)\right|+o(\tau)=1+o(\tau),
$$
leading to
\begin{eqnarray*}
\int_{Z_\tau}f(x,t_0+\tau,v)\,dx\,dv=
\int_{Z}f\big(x_0+\tau v_0,t_0+\tau,v_0-\tau\nabla u(x_0,t_0)\big)\,dx_0\,dv_0+o(\tau).
\end{eqnarray*}
Hence, the assumption that stars are conserved thus translates into
\begin{eqnarray*}
0&=&\left.\frac{d}{d\tau}\int_{Z_\tau} f(x(t_0+\tau),t_0+\tau,v(t_0+\tau))\,dx\,dv\right|_{\tau=0}\\
&=&\left.
\frac{d}{d\tau}\left(
\int_{Z}f\big(x_0+\tau v_0,t_0+\tau,v_0-\tau\nabla u(x_0,t_0)\big)\,dx_0\,dv_0+o(\tau)\right)\right|_{\tau=0}\\
&=& \int_{Z}
\Big[ \nabla_x f(x_0,t_0,v_0)\cdot v(x_0,t_0)+\partial_t f(x_0,t_0,v_0)
-\nabla_v f(x_0,t_0,v_0)\cdot\nabla u(x_0,t_0)
 \Big]\,dx_0\,dv_0.
\end{eqnarray*}
The arbitrariness of the phase-space domain~$Z$ thus leads to
$$ \nabla_x f(x,t_0,v)\cdot v(x,t_0)+\partial_t f(x,t_0,v)
-\nabla_v f(x,t_0,v)\cdot\nabla u(x,t_0)=0,$$
which is the first equation in~\eqref{STARS}.

Additionally, by Gau{\ss}'s Flux Law for gravity, 
we know that the flux of the gravitational field through a surfaces
balances the total mass enclosed within the surface: namely,
neglecting dimensional constants,
for every domain~$\Omega\subset\R^3$,
$$ \int_{\partial\Omega} \nabla u(x)\cdot\nu(x)\,d{\mathcal{H}}^{n-1}_x=\int_\Omega \rho(x,t)\,dx.$$
{F}rom this and the Divergence Theorem we arrive at
$$ \int_{\Omega} \Delta u(x)\,dx=\int_\Omega \rho(x,t)\,dx,$$
and consequently, since~$\Omega$ is arbitrary,
we find that~$\Delta u=\rho$, which 
is the second equation in~\eqref{STARS}.

For additional information about galaxy dynamics see e.g.~\cite{BINNEY} and the references therein.

\subsection{How to count what we cannot see}\label{EIN:SEC}

\begin{figure}
  \centering
  \includegraphics[width=.25\linewidth]{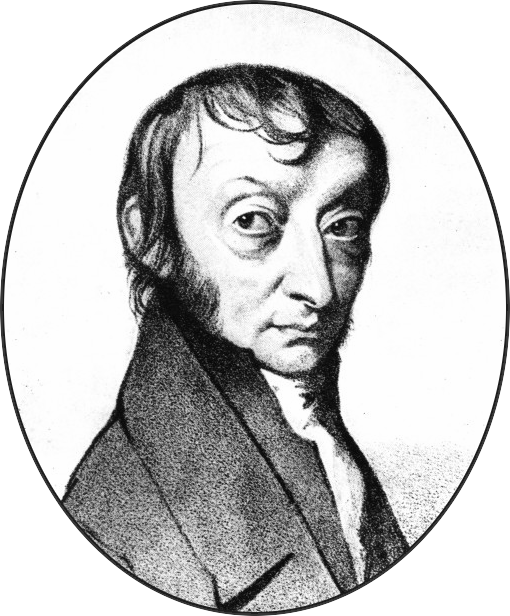}
 \caption{\sl Portrait of Amedeo Avogadro (Public Domain image from
        Wikipedia).}\label{AVPGRa-EULEROFIDItangeFI}
\end{figure}

In 1812, Lorenzo Romano Amedeo Carlo Avogadro, Count of Quaregna and Cerreto,
hypothesized that the volume of a gas at a given pressure and temperature is proportional to the number of atoms or molecules.
That is, equal volumes of gases at the same temperature and pressure have the same number of molecules,
regardless of the nature of the gas.
In particular,
by Avogadro's Law, the number of molecules or atoms in a given volume of ideal gas is independent of their size.
While this prescription only holds true for ideal gases,
and real gases show instead small deviations from this ideal configuration,
Avogadro's Law is often a very useful approximation and it always provides a great conceptual tool,
since it detects a universal quantity that only depends on reasonable macroscopic parameters, such as
temperature and pressure, and is independent of the specific situation under consideration, namely the type of gas
that one is taking into account. 
See Figure~\ref{AVPGRa-EULEROFIDItangeFI} for a portrait of Avogadro (by C. Sentier).

It is perhaps worth stressing that Avogadro's Law is not so intuitive.
For instance, if one has three units of volumes of hydrogen
and one unit of volume of nitrogen, after they combine together and produce
ammonia, how many units of volume of ammonia do we expect (assuming that
pressure and temperature are maintained constant)?
One unit? Three units? Four units?

\begin{figure}
                \centering
                \includegraphics[width=.69\linewidth]{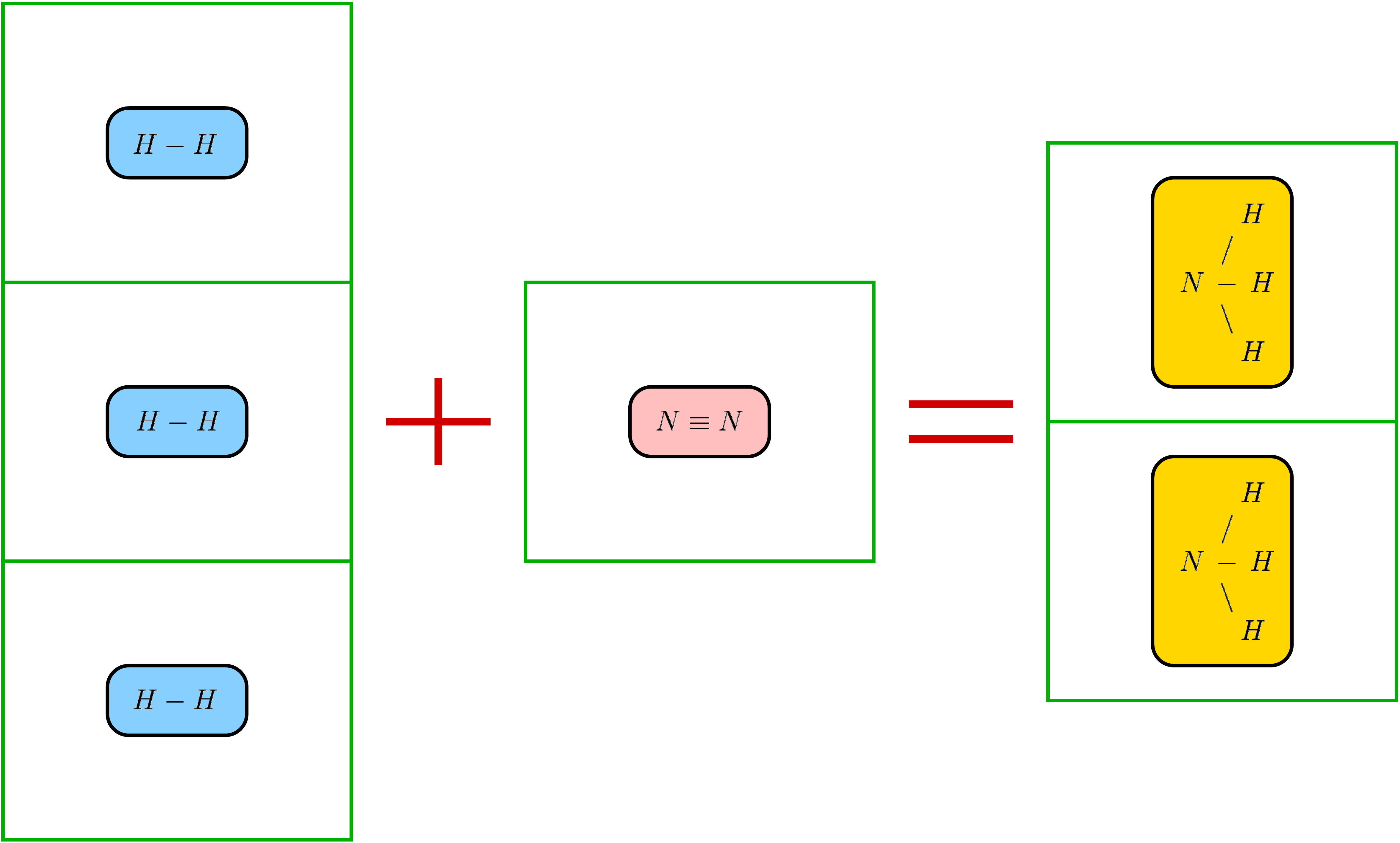}
        \caption{\sl The algebra in~\eqref{REACCHI}
with the corresponding gas volumes according to Avogadro's Law.}\label{CHIMIPIPEDItangeFI}
\end{figure}

The correct answer according to Avogadro's Law is {\em two} units of volume.
This is because three volumes of hydrogen contain~$3k$ molecules of hydrogen~$H_2$ (for some~$k\in\N$)
and one volume of nitrogen contains~$1k$ molecules of nitrogen~$N_2$
(and the proportionality factor~$k$ is the same for both hydrogen
and nitrogen, by Avogadro's Law).
Thus, from the reaction
\begin{equation}\label{REACCHI} 3H_2+1N_2=2NH_3,\end{equation}
the combination of~$3k$ molecules of hydrogen
and~$1k$ molecules of nitrogen produces~$2k$ molecules
of ammonia~$NH_3$. And again by Avogadro's Law
this corresponds to two volumes of ammonia. See Figure~\ref{CHIMIPIPEDItangeFI}
for a sketch of this phenomenon.\medskip

Also, in chemistry, an obvious practical issue is given by the fact that molecules are small
while in the lab one has to deal with macroscopic quantities. 
Thus, while it is desirable from the theoretical point of view to found the concept of
``amount of a given substance'' on the number of ``elementary entities'', i.e., atoms or molecules, present in the substance
(because these elementary entities are the ones which will play a role in the chemical reactions),
for practical purposes it is often convenient to relate the notion of ``amount of a given substance''
to the ratio of measured macroscopic quantities (because measuring a huge number of microscopic entities
is typically unfeasible).

We remark that the equivalence between these two approaches is a notable consequence of
Avogadro's Law. Indeed, if, at a given pressure and temperature, we take, say, one gram of
atomic hydrogen~$H$ (or, more simply, two grams of the common molecule of hydrogen~$H_2$)
and we measure its volume, and then we consider the same volume of another gas, such as the carbon dioxide~$CO_2$,
we know from Avogadro's Law that the latter contains the same number of molecules of~$CO_2$
as the number of atoms of hydrogen~$H$ (or the number of molecules of~$H_2$) in the previous container.
If we measure the weight of the carbon dioxide~$CO_2$ in the second container,
we find that it amounts to~$44$ grams: we can accordingly conclude that
the mass of a molecule of carbon dioxide~$CO_2$
is~$44$ times bigger than the mass of an atom of hydrogen~$H$
(or~$22$ times bigger than the mass of a molecule of hydrogen~$H_2$).
That is, we have measured, with a minimal effort, the microscopical weight of any molecule of any gases
in terms of a given ``unit of measure'', such as the atomic weight of~$H$ (or the molecular weight of~$H_2$).

Hence, it has become a common practice to adopt
the notion of ``mole'' to denote the mass of a given substance
which contains the same number of molecules (or atoms in case of
pure atomic elements) as one gram of atomic hydrogen~$H$.
In the previous example, since the mass of a molecule of carbon dioxide~$CO_2$
happens to be~$44$ times bigger than the mass of an atom of hydrogen~$H$,
we infer that
\begin{equation}\label{LAMOLE}
{\mbox{a mole of carbon dioxide weights~$44$ grams.}}\end{equation}

Moreover, by Avogadro's Law,
a mole of any gas at a given temperature and pressure occupies the same volume, making it a quite convenient
way to measure the (relative) amount of substance just by comparing volumes of gases.

The natural question is thus, given a certain chemical compound
(say, an ideal gas), how many molecules
of it correspond to a mole? How does this number depend on the specific chemical compound?
The answer is easy: by construction this number is the same for all (ideal)
substances and equals
the number of hydrogen atoms~$H$ necessary to form one gram of hydrogen, and this
is called\footnote{We stress that
Avogadro's Number is not a ``pure number'' but a physical constant of dimension
``one over moles''. See also~\cite{AVOG} for additional information on Avogadro's Number.}
Avogadro's Number~${\mathcal{N}}_{\!\mathcal{A}}$. \index{Avogadro's Number}
\medskip

Notice that the knowledge of~${\mathcal{N}}_{\!\mathcal{A}}$ also permits to ``pass from the ratios
to the real quantities''. For instance from~\eqref{LAMOLE} we obtain that each molecule of
carbon dioxide weights~$\frac{44}{ {\mathcal{N}}_{\!\mathcal{A}} }$ grams:
in general, the knowledge of Avogadro's Number entails the knowledge of atomic and molecular masses
as a byproduct of macroscopic measurements.

Therefore, it is highly desirable to have a precise\footnote{It turns
out (possibly with some mild approximation)
that~${\mathcal{N}}_{\!\mathcal{A}}=6.02214076\times 10^{23}$. Several people noticed
a similarity between this number and~$2^{79}$ (which is half of the so-called yobibyte~$2^{80}$):
indeed, it is a rather surprising coincidence that 
$$ \frac{{\mathcal{N}}_{\!\mathcal{A}}-2^{79}}{{\mathcal{N}}_{\!\mathcal{A}}+2^{79}}$$
is almost zero.
There is nothing special in this, coincidences happen.
For instance, it is a mere, and quite remarkable, coincidence that the most important physical constant, namely the speed
of light, is almost equal to~$300.000.000$ m/s.
It is also a coincidence that
the angular diameter of the Sun seen from Earth is quite close to that of the Moon 
making it possible to see a full solar eclipse.
It is also a coincidence that
$$ \left(\frac{\pi^e-e^\pi}{\pi^e+e^\pi}\right)^{\pi+e}$$
turns out to be almost equal to zero.}
quantification of~${\mathcal{N}}_{\!\mathcal{A}}$. This is not a mathematical problem per se, since
to quantify~${\mathcal{N}}_{\!\mathcal{A}}$ one can only measure it through experiments.
However, one needs a brilliant piece of mathematics to contrive
a practical procedure to experimentally calculate~${\mathcal{N}}_{\!\mathcal{A}}$
and we recall here the method devised by Albert Einstein in~\cite{zbMATH02652222}
(at that time, he was working as a patent clerk in Bern, Switzerland,
and the very notion of atoms and molecules was still a subject of controversy and intense scientific debate, since individual
atoms and molecules were, for the instruments available in~1905, simply
``too small'' and ``too fast'').

To measure Avogadro's Number, Einstein utilized a molecular theory of heat
more or less in the lines of the Brownian motion described here
in Section~\ref{CHEMOTX} and leading to the heat equation in~\eqref{HJA:A89MAOL8789892}.
For this, Einstein thought of some ideal spherical particles of radius~$a$
subject to a Brownian motion in an ideal (for simplicity, one-dimensional and infinitely long) pipe
containing some fluid with viscosity coefficient equal to~$\mu$.
In this setting, according to Einstein's calculation,
\begin{equation}\label{EINAVOGEAFG}
{\mathcal{N}}_{\!\mathcal{A}}=\frac{RTt}{3\pi \mu a\,X^2(t)},
\end{equation}
where~$R$ is the universal\footnote{Note that in principle one can measure~$R$
via the ideal gas prescription~$PV=nRT$ in terms of macroscopic quantities (though more precise measures
of~$R$ typically rely on a finer analysis). Thus, for the purposes
of this notes, $R$ is just a ``known'' constant.}
gas constant,
$T$ denotes the temperature (or, more precisely, the thermodynamic temperature, that is the temperature measured in Kelvin), $t$ is any given time and~$X^2(t)$ is the averaged square distance traveled by a Brownian particle in an interval
of time equal to~$t$.

Equation~\eqref{EINAVOGEAFG} is quite remarkable since it equates
Avogadro's Number with quantities which can be measured experimentally, thus leading
to a precise determination of~${\mathcal{N}}_{\!\mathcal{A}}$
(that's the power of mathematics!).

\begin{figure}
  \centering
  \includegraphics[width=.25\linewidth]{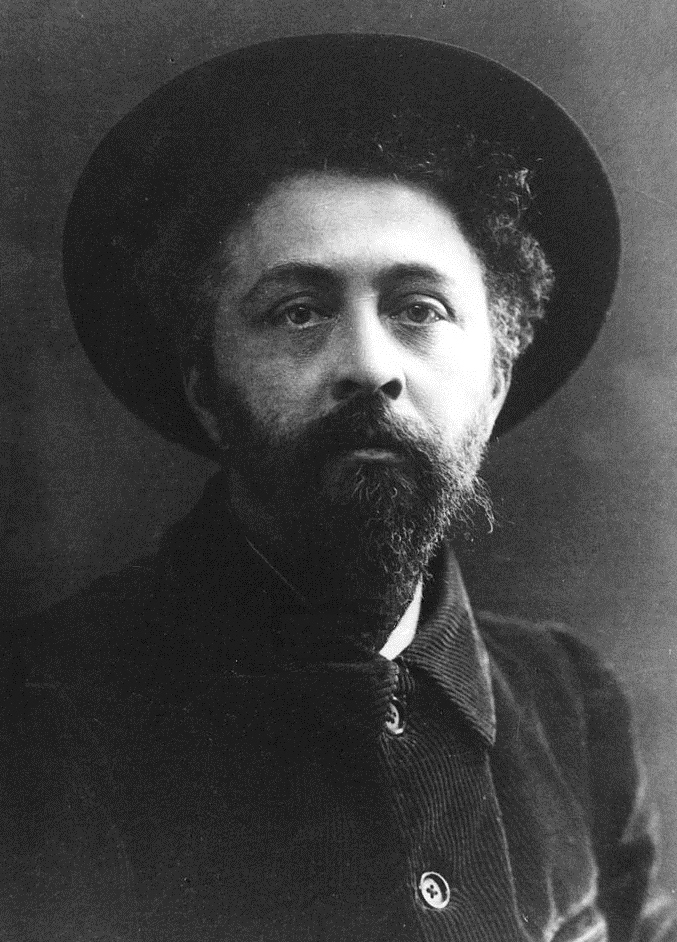}
 \caption{\sl Jean Perrin (Public Domain image from
 Wikipedia).}\label{HENITONFJeanPerrinIDItangeFI}
\end{figure}

The experiment suggested by Einstein was indeed carried through by Jean Baptiste Perrin (see Figure~\ref{HENITONFJeanPerrinIDItangeFI})
and his team of research students
in~1909. The set up of the test consisted in a camera lucida with a microscope to
observe and record the motion of suspended gamboge\footnote{Gamboge
is a yellowish pigment.
The word gamboge comes from gambogium, the Latin word for the pigment, which in turns derives from Gambogia, the Latin word for Cambodia. 

Besides its strong laxative properties, which did not play a role in Perrin's experiment, 
gamboge also possessed the convenient feature of producing, with appropriate alcoholic additives and after a selective centrifuge,
almost perfect spherules with equal radius, thus providing the ideal suspension particles for the experiment
proposed by Einstein (this was not a cheap preparation,
about one kilo of gamboge and several months of work produced few decigrams of useful spherules).}
particles in
a liquid of a given viscosity and constant temperature.
Marking the particle's position on a piece of
graph paper at timed intervals, and repeating the test under different conditions
(such as different viscosity, temperature, radius of the gamboge grains, etc.), Perrin
obtained consistent results about the right hand side of~\eqref{EINAVOGEAFG},
thus producing a rather accurate measurement of Avogadro's Number~${\mathcal{N}}_{\!\mathcal{A}}$
through a wealth of measurements that could not be contested.

This experiment also provided experimental confirmation of Einstein's equation~\eqref{EINAVOGEAFG} and
raised atoms and molecules from the status of hypothetical objects to concrete and observable entities,
thus offering the ultimate confirmation to the atomic nature of matter
and concluding the struggle regarding the physical reality of molecules.
For these achievements, Perrin was awarded with the Nobel Prize for Physics in 1926.\medskip

It is now time to go back to Einstein's equation~\eqref{EINAVOGEAFG}
to understand its roots in the theory of partial differential equations.
To this end, we recall that the first ingredient towards~\eqref{EINAVOGEAFG} is Stokes' Law\index{Stokes' Law}
about the viscosity force acting on a small sphere moving through a viscous fluid.
Namely, studying the incompressible steady flow of the Navier-Stokes equation (see Section~\ref{EUMSD-OS-32456i7DNSSE234R:SEC}), George Gabriel Stokes in~1851 had proposed the formula
\begin{equation}\label{SYO-EDCSJDKA235athcma}
F=6\pi \mu a v\end{equation}
to describe the frictional force~$F$
acting on a spherical particle of radius~$a$ with relative velocity~$v$ in a fluid with viscous coefficient equal to~$\mu$.

\begin{figure}
  \centering
  \includegraphics[width=.25\linewidth]{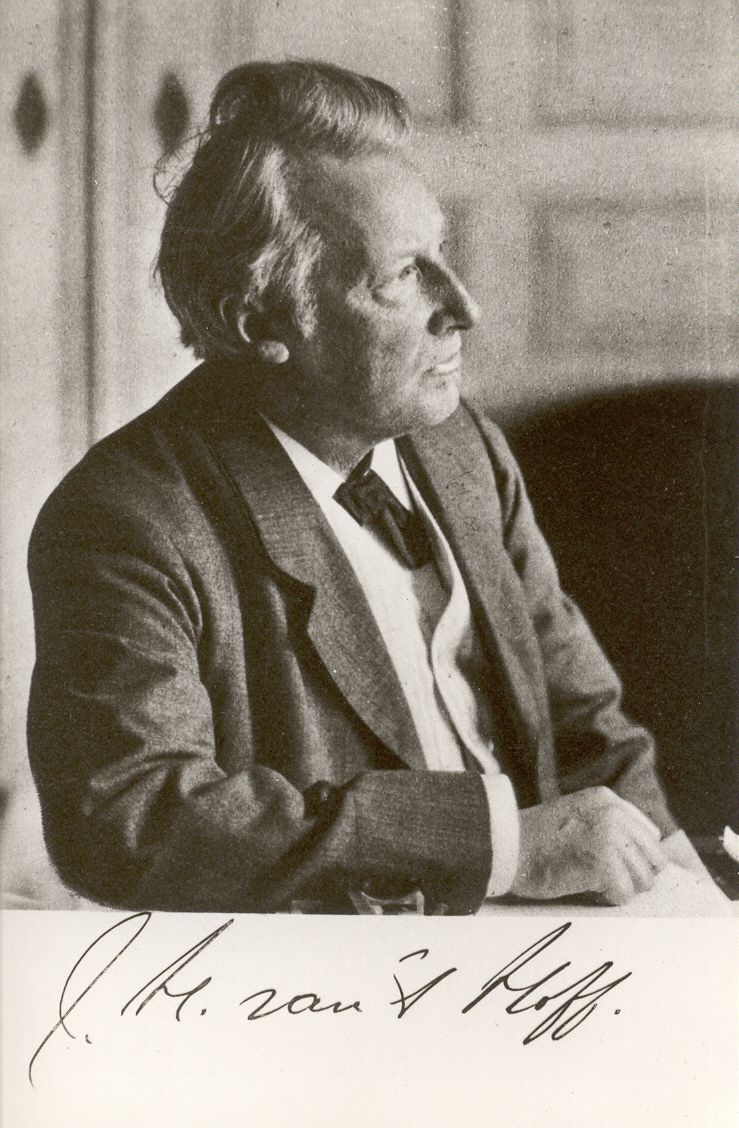}
 \caption{\sl Henry van 't Hoff (Public Domain image from
 Wikipedia).}\label{HENITONFIDItangeFI}
\end{figure}

The second ingredient towards~\eqref{EINAVOGEAFG} is
the van 't Hoff equation\index{van 't Hoff equation} for\footnote{Equation~\eqref{BVNBSstVASNnsrNA}
is named after Jacobus Henricus van 't Hoff Jr.,
first winner of the Nobel Prize in Chemistry, see Figure~\ref{HENITONFIDItangeFI}. Allegedly, van 't Hoff
chose to study chemistry
against the wishes of his father. Also, it seems he came up with equation~\eqref{BVNBSstVASNnsrNA}
after a chance encounter and a quick chat with a botanist friend during a walk in a park in Amsterdam,
thinking about an analogy to the law for ideal gases.} the osmotic pressure~$p$
of a solution containing~$  n$  moles of solute particles in a solution of volume~$  V $,
given by
\begin{equation}\label{BVNBSstVASNnsrNA}
p=\frac{nRT}V,\end{equation}
in which~$R$ is the universal gas constant.

The third ingredient towards~\eqref{EINAVOGEAFG} is the assumption that
the Brownian motion of the suspension particles can be effectively described by the heat equation,
as discussed in Section~\ref{CHEMOTX}. Hence, we suppose that the density~$u$ of the particles fulfills the equation
\begin{equation}\label{EINHEAT}
\partial_t u= c\partial_{xx} u,
\end{equation}
for some constant diffusion coefficient~$c>0$ (here, we are assuming that the motion of the particles is confined
to a linear pipe, hence the variable~$x$ is one-dimensional).

We can now retake van 't Hoff equation~\eqref{BVNBSstVASNnsrNA}
by considering that the number of floating particles~$N$ over the volume~$V$
equals the density~$u$ and, by definition, a mole contains~${\mathcal{N}}_{\!\mathcal{A}}$ particles.
These observations lead to
\begin{equation}\label{EDCSJDKA235S6IMfujNMimfr0}
p=\frac{NRT}{{\mathcal{N}}_{\!\mathcal{A}}\,V}=\frac{RT u}{{\mathcal{N}}_{\!\mathcal{A}}}.
\end{equation}
Also, the force~$F$ acting on each floating particle arises from a pressure gradient, namely the gradient pressure force with respect to the unit of mass is~$\frac{\partial_x p}u$.
This remark and~\eqref{EDCSJDKA235S6IMfujNMimfr0} give that
\begin{equation*}
F u=\partial_xp= \frac{RT\,\partial_x u }{{\mathcal{N}}_{\!\mathcal{A}}}.
\end{equation*}
At equilibrium, the force must correspond to the frictional force in~\eqref{SYO-EDCSJDKA235athcma}, therefore we obtain that
\begin{equation}\label{EDCSJDKA235S6IMfujNMimfr02}
6\pi \mu a v u=\frac{RT\,\partial_x u}{{\mathcal{N}}_{\!\mathcal{A}}}.\end{equation}

Now we compute the flux of floating particles through an ideal cross section of the pipe
by using two possible strategies. On the one hand, we can consider the model of particles traveling with velocity~$v$,
hence we can quantify the flux by the product between~$v$ and the particle density~$u$.
On the other hand, the flux at some cross section~$x$ can also be obtained by the variation in time
of the mass preceding~$x$,
that is~$\partial_t \int_{-\infty}^x u(\xi,t)\,d\xi$. These observations and~\eqref{EINHEAT} give that
$$ v(x,t) u(x,t)=\partial_t \int_{-\infty}^x u(\xi,t)\,d\xi=
\int_{-\infty}^x \partial_t u(\xi,t)\,d\xi=c\int_{-\infty}^x \partial_{xx} u(\xi,t)\,d\xi=
c\partial_x u(x,t),$$
where, as usual, we have performed formal calculations and assumed some convenient decay of the solution at~$-\infty$.

{F}rom this and~\eqref{EDCSJDKA235S6IMfujNMimfr02} we arrive at
\begin{equation}\label{EDCSJDKA235S6IMfujNMimfr03}
{\mathcal{N}}_{\!\mathcal{A}}=\frac{RT\,\partial_x u}{6\pi \mu a v u}=\frac{RT}{6\pi \mu a c}
.\end{equation}

Now we observe that
$$g (x,t):=\frac{1}{\sqrt {4\pi ct}}\exp \left(-{\frac{x^{2}}{4ct}}\right)$$
is a solution of~\eqref{EINHEAT} concentrating at the origin when~$t=0$, hence
we can suppose that the general solution~$u$ of~\eqref{EINHEAT} occurs as a superposition\footnote{The formalization
of this superposition method is the core of the notion of fundamental solution and, in the elliptic setting, it will
play a decisive role in Section~\ref{lfundsP-S} and in particular in the forthcoming
Proposition~\ref{ESI:NICE:FB}.}
of traslations of~$g$ of the form~$g(x-x_0,t)$ (notice that each of these translations
would correspond to a particle located at~$x_0$ when~$t=0$).

The average squared displacement at time~$t$ of the particle starting at the origin is therefore provided by the quantity
$$ \int_\R x^2 g(x,t)\,dx=\int_\R\frac{x^2}{\sqrt {4\pi ct}}\exp \left(-{\frac{x^{2}}{4ct}}\right)\,dx=2ct.
$$
Since this displacement should not depend on the original position, we infer that it is the same for all
possibilities of initial location of the particles, namely
\begin{equation}\label{EINAVOGEAFG2} X^2(t)=2ct.\end{equation}
This is conceptually an important result since it suggests that velocities (namely ratios of~$X(t)$ over~$t$)
of Brownian motions behave sort of~$\frac{\sqrt{2ct}}{t}=\frac{\sqrt{2c}}{\sqrt{t}}$,
hence quite irregularly\footnote{For a formalization of this heuristic and imprecise discussion see e.g.~\cite[Section~3.4]{MR3154922}.} for short times. That is,
measuring the velocity of the floating particles in an extremely short interval of time
just produces a result that approaches infinity (in a sense,
attempted experiments in this direction would simply end up measuring the wrong quantity).
One of the merit of Einstein's approach is thus to have bypassed this hindrance by detecting
the relevant physical quantities which were both consistent from the mathematical point of view and experimentally
measurable in laboratories. Notice indeed that by plugging~\eqref{EINAVOGEAFG2}
into~\eqref{EDCSJDKA235S6IMfujNMimfr03} we obtain~\eqref{EINAVOGEAFG}, as desired.

\subsection{Good vibrations}\label{GUITAR}

\begin{figure}
                \centering
                \includegraphics[width=.6\linewidth]{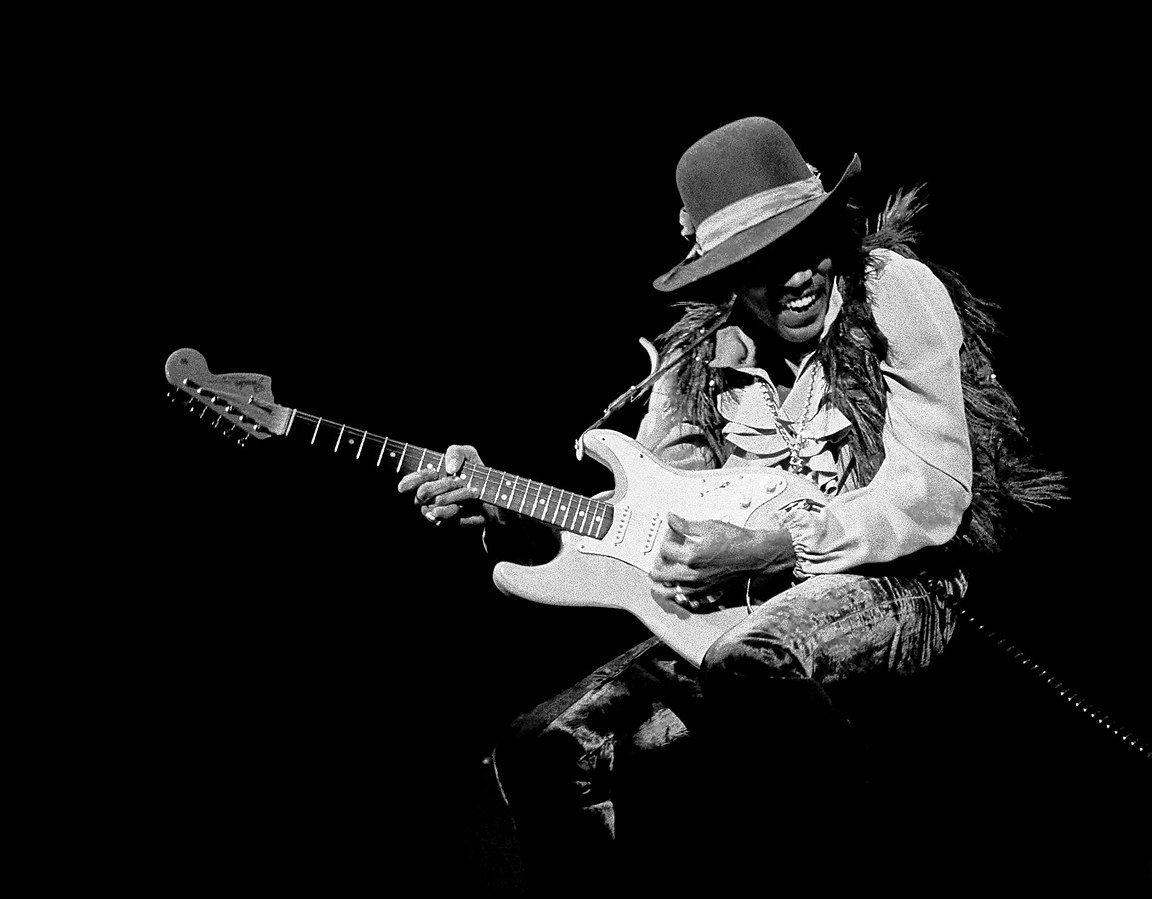}
        \caption{\sl Jimi Hendrix (photo by Steve Banks, image from
        Wikipedia,
        licensed under the Creative Commons Attribution-Share Alike 4.0 International license).}\label{2HAFODAKEDFUMSPierre-SimonLaplacldRGIRA4AXELEUMDJOMNFHARLROA7789GIJ7solFUMHDNOJHNFOJED231HEND}
\end{figure}

Another classical occurrence for partial differential equations arises in the description
of vibrating strings\index{vibrating string} (e.g., the string of a guitar, see
Figure~\ref{2HAFODAKEDFUMSPierre-SimonLaplacldRGIRA4AXELEUMDJOMNFHARLROA7789GIJ7solFUMHDNOJHNFOJED231HEND}). Different models
are possible, taking into account different physical assumptions (see e.g.~\cite[Remark~10.5]{MR2164768}),
but we focus here on the simplest possible (indeed linear) equation.
In our description, the string at rest is modeled by the segment~$[0,L]\times\{0\}$
for some~$L>0$ which represents the length of the string in the absence of further external forces.
One can assume that the string is constrained at the extrema
and can be deformed into the graph of a function~$u$: namely, the position
of the vibrating string at each moment of time~$t$ is described by the graph~$(x,u(x,t))\in[0,L]\times\R$
and the constraints for the string correspond to the boundary prescription~$u(0)=u(L)=0$.

We suppose that the string is subject to gravity (acting downwards in the vertical direction)
and to its own ``tension''\index{tension} (that is the internal force of the string acting between its elements
providing an elastic reaction to the external forces). The magnitude of this tension
at each point of the string~$(x,u(x,t))$ will be denoted by~$T(x,t)$.

\begin{figure}
  \centering
  \includegraphics[width=.6\linewidth]{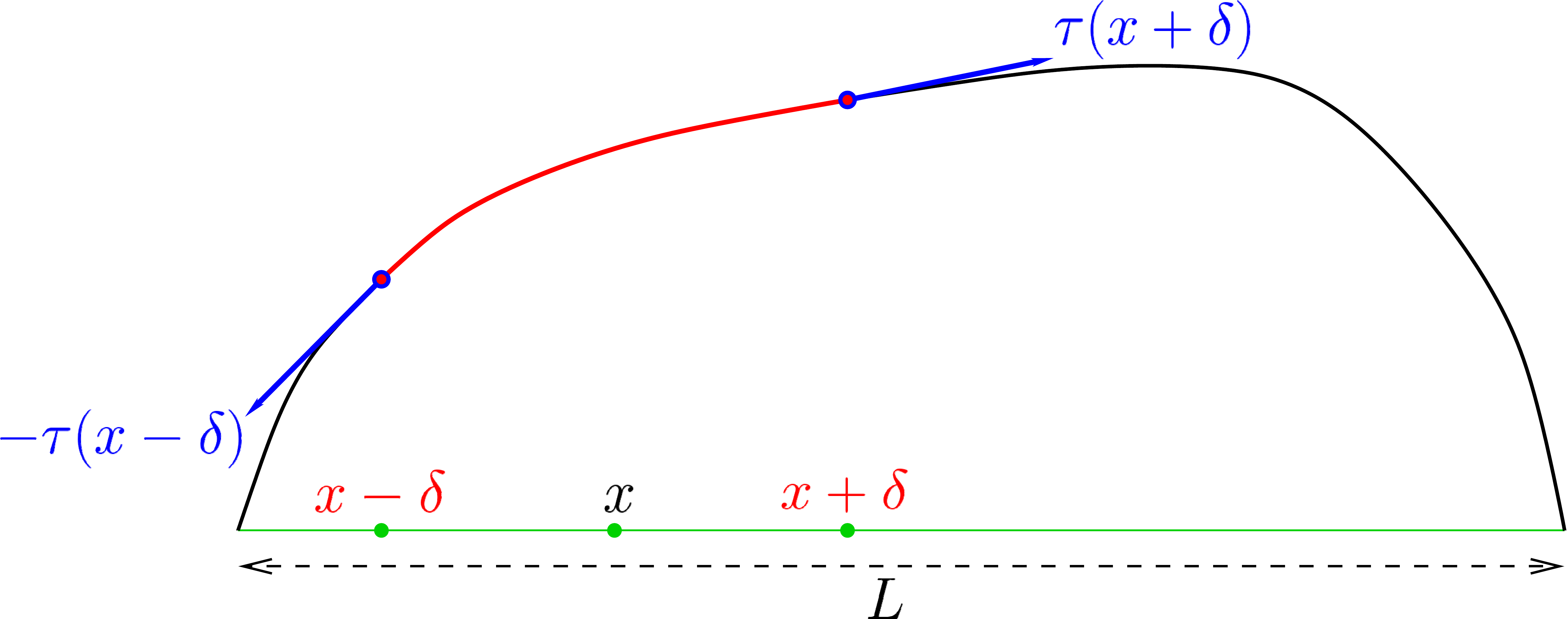}
 \caption{\sl A string subject to its own tension.}\label{PLATENSIWASTREDItangeFIlASE}
\end{figure}

We also consider the unit tangent vector to the string, as given by
$$ \tau(x,t):=\frac{\big(1,\partial_x u(x,t)\big)}{\sqrt{1+\big(\partial_x u(x,t)\big)^2}}.$$
The vibrating string model then consists in the assumption that, given a small~$\delta>0$,
the tension forces acting at the string point~$(x,u(x,t))$ are the byproduct of
a force~$F_-$ acting on~$(x-\delta,u(x-\delta,t))$ and a force~$F_+$ acting on~$(x+\delta,u(x+\delta,t))$
whose magnitude is proportional to the tension~$T(x,t)$ and opposite tangential directions
(namely, the direction of~$F_-$ is~$-\tau(x-\delta,t)$ and
the direction of~$F_+$ is~$\tau(x+\delta,t)$), see Figure~\ref{PLATENSIWASTREDItangeFIlASE}.

Therefore, the total tension force acting on the string at the point~$(x,u(x,t))$ takes the form
\begin{equation}\label{TOTFOSTR}\begin{split}&F_++F_-=\kappa T(x)\big(\tau(x+\delta,t)-\tau(x-\delta,t)\big)=
2\kappa \delta T(x)\partial_x\tau(x,t)+o(\delta)\\&\qquad\qquad\qquad=
2\kappa \delta T(x)\partial_x\left(
\frac{\big(1,\partial_x u(x,t)\big)}{\sqrt{1+\big(\partial_x u(x,t)\big)^2}}\right)
+o(\delta).\end{split}\end{equation}
The constant~$\kappa>0$ depends on the material of which the string is made
and accounts for the elastic properties of the string.

The total force acting on the string, taking into account gravity, is therefore
\begin{equation*} F=(F_1,F_2):=F_++F_- -(0,m_\delta g)=
2\kappa\delta T(x)\partial_x\left(
\frac{\big(1,\partial_x u(x,t)\big)}{\sqrt{1+\big(\partial_x u(x,t)\big)^2}}\right)-(0,m_\delta g)
+o(\delta),\end{equation*}
where~$m_\delta$ is the mass of the string located between~$x-\delta$ and~$x+\delta$ and~$g$ denotes
the gravity acceleration.

If we suppose that the density of the string is constant,  say equal to some~$\rho>0$, the mass~$m_\delta$
is the product of~$\rho$ and the length of the string located between~$x-\delta$ and~$x+\delta$, namely
\begin{equation}\label{MEPSI} m_\delta= \rho \int_{x-\delta}^{x+\delta} \sqrt{1+\big(\partial_x u(\zeta,t)\big)^2}\,d\zeta
=2\rho\delta \sqrt{1+\big(\partial_x u(x,t)\big)^2}+o(\delta).\end{equation}

Furthermore, by Newton's Second Law, the vertical component of the force is equal to the mass times
the vertical acceleration of the string, that is
$$ F_2=m_\delta\, \partial_{tt} u(x,t).$$
In light of these considerations, we find that\footnote{{F}rom the geometric point
of view, it is interesting to observe that the term~$\partial_x\left(
\frac{ \partial_x u}{\sqrt{1+\big(\partial_x u\big)^2}}\right)$ in~\eqref{CUSTRJMSONGJMSD}
detects the curvature of the string
(remarks of this type will be formalized in a general setting
in Theorem~\ref{INCOO}).}
for small~$\delta$
\begin{equation}\label{CUSTRJMSONGJMSD}
\begin{split}
\partial_{tt} u\,&= \frac{F_2}{m_\delta}
\\&=\frac1{m_\delta}\left[
2\kappa\delta T\partial_x\left(
\frac{ \partial_x u}{\sqrt{1+\big(\partial_x u\big)^2}}\right)- m_\delta g+o(\delta)
\right]\\&=
\frac{2\kappa\delta T}{m_\delta}
\partial_x\left(
\frac{ \partial_x u}{\sqrt{1+\big(\partial_x u\big)^2}}\right)- g+o(1)\\&=
\frac{\kappa T}{
\rho \sqrt{1+\big(\partial_x u\big)^2}+o(1)}
\partial_x\left(
\frac{ \partial_x u}{\sqrt{1+\big(\partial_x u\big)^2}}\right)- g+o(1)\\
&=
\frac{\kappa T}{
\rho \sqrt{1+\big(\partial_x u \big)^2}}
\partial_x\left(
\frac{ \partial_x u}{\sqrt{1+\big(\partial_x u\big)^2}}\right)- g+o(1).
\end{split}\end{equation}
By formally sending~$\delta\searrow0$ we thus obtain the partial differential equation
\begin{equation}\label{CUSTRJMSONGJMSD-2}
\partial_{tt} u=
\frac{\kappa T}{
\rho \sqrt{1+\big(\partial_x u\big)^2}}
\partial_x\left(
\frac{ \partial_x u}{\sqrt{1+\big(\partial_x u\big)^2}}\right)- g.
\end{equation}
Several different models can be consider to describe explicitly the tension~$T$
and thus find an expression for~\eqref{CUSTRJMSONGJMSD-2}
solely depending on the shape of the string, on its density and elastic constant
and on the gravity acceleration.
One commonly accepted model is to consider the string as ``inextensible''
and~$T$ as a constant (the tension being distributed uniformly along the whole string).
Another possibility is to consider the string as a superposition of infinitesimal elastic springs and thus take~$T$ as proportional to the string's elongation.
Both these models however share the common treat that for small elongations~$T$
is constant. Therefore, in the small elongation approximation, we take~$T:=1$
and we disregard quadratic terms in~$u$ in~\eqref{CUSTRJMSONGJMSD-2}
(the ansatz being that for small~$u$ the linear terms
will prevail against the quadratic ones). With this simplification in mind,
one can reduce~\eqref{CUSTRJMSONGJMSD-2} to
\begin{equation}\label{CUSTRJMSONGJMSD-3}
\partial_{tt} u=
\frac{\kappa }{
\rho }
\partial_{xx}u- g.
\end{equation}
In particular, in the absence of gravity, the linear approximation of the vibrating string can be seen
as a particular case of the wave equation in~\eqref{WAYEBJJD121t416JH098327uyrhdhc832nbcM}.
Interestingly, the speed of propagation~$c$ in~\eqref{WAYEBJJD121t416JH098327uyrhdhc832nbcM}
corresponds here to~$\sqrt{\frac{\kappa }{
\rho }}$. Namely the higher the elastic parameter of the string, the higher the
speed of propagation; the higher the density (equivalently, the heavier the string),
the slower the speed of propagation.

The dependence of the solution of~\eqref{CUSTRJMSONGJMSD-3} upon the structural parameters
of the string is not a mere mathematical curiosity: on the contrary it has a deep impact on music,
since these parameters precisely determine the pitch of the sound that the string generates.
Let us see, for instance, how understanding partial differential equations may turn out 
to be useful when tuning a guitar. For accomplishing this goal, we neglect the gravity effect and summarize~\eqref{CUSTRJMSONGJMSD-3} and the boundary conditions of the string
into the following mathematical setting:
\begin{equation}\label{CwUSTRJMSONGJMSD-4}
\begin{dcases}
\partial_{tt} u(x,t)=\displaystyle{\frac{\kappa }{
\rho }}\partial_{xx}u(x,t)\qquad{\mbox{for all }}(x,t)\in(0,L)\times(0,+\infty),\\
u(x,0)=\eta\displaystyle\sin\frac{\pi x}{L},\\
\partial_t u(x,0)=0,\\
u(0,t)=u(L,t)=0\qquad{\mbox{for all }}t\in(0,+\infty).
\end{dcases}
\end{equation}
The prescriptions in~\eqref{CwUSTRJMSONGJMSD-4} model a guitar string
with rest length~$L$ and null initial velocity which is initially displaced (say by 
an expert fingerpicking) as a sinusoidal graph (for the sake of simplicity, though more complicated
initial situations could be taken into account). The small parameter~$\eta>0$
has been introduced here above just to be consistent with the small elongation approximation.

\begin{figure}
  \centering
  \includegraphics[width=.40\linewidth]{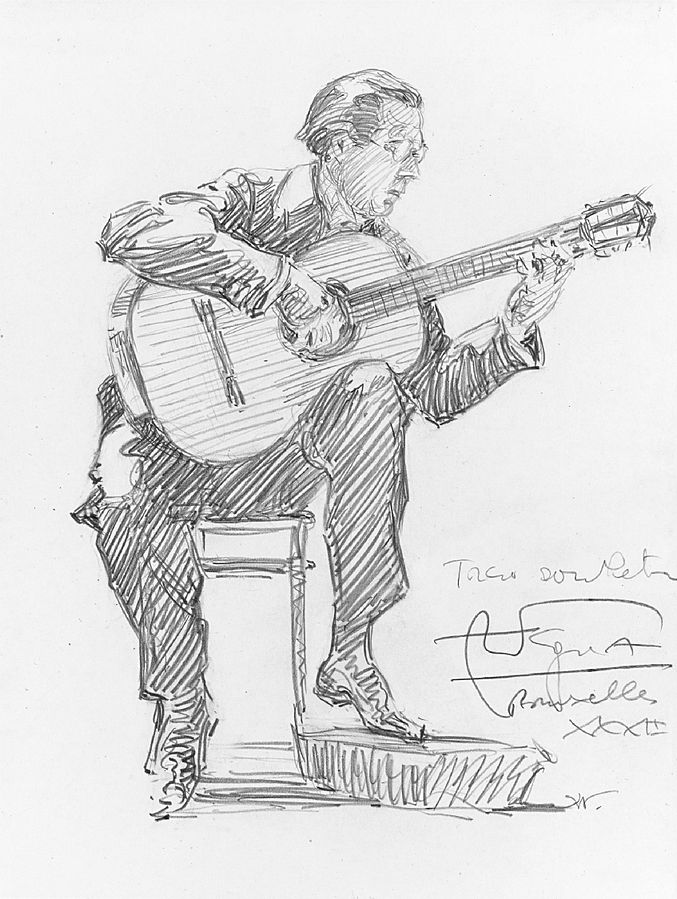}
 \caption{\sl 
Andr\'es Segovia Torres, 1st Marquis of Salobre\~{n}a
during a recital in Brussels on December 15, 1932
 (Public Domain image from
        Wikipedia, source Sotheby's).}\label{9ijj0SIegib Ov801tangeFI}
\end{figure}

It can be readily checked that the function
\begin{equation}\label{CwUSTRJMSONGJMSD-5} u(x,t)=\eta\cos\frac{\sqrt\kappa\,\pi t}{\sqrt\rho\,L}
\sin\frac{\pi x}{L}\end{equation}
solves~\eqref{CwUSTRJMSONGJMSD-4}.
This solution is quite telling regarding the pitch played by the string.
Indeed, in light of~\eqref{CwUSTRJMSONGJMSD-5},
the string vibrates with a frequency
\begin{equation}\label{CwUSTRJMSONGJMSD-6}
\omega:=
\frac{\sqrt\kappa\,\pi}{\sqrt\rho\,L}.\end{equation} In particular:
\begin{itemize}\item
Strings with a higher elastic coefficient produce a higher pitch (and we already know this by experience, because tightening the tuning pegs of a guitar enhances the string's tension and correspondingly raises the pitch),
\item Longer strings produce lower pitches and shorter strings produce higher pitches (and we already know this by experience too, since, in playing, one changes the length of the vibrating string by holding it firmly against the fingerboard with a finger and shortening the string, that is stopping it on a higher fret, gives higher pitch),
\item Strings with higher density produce lower pitches, that is more massive strings vibrate more slowly (and this is also in agreement with our experience since, for instance, on classical guitars, the low density nylon strings are used for the high pitches and the higher density wire-wound strings are employed for low pitches).\end{itemize}
To appreciate even more the impact of mathematics on guitar playing,
we recall that an octave\footnote{The name ``octave'' comes from the Latin adjective ``octava'', meaning ``eighth''. This etymology however is possibly a bit confusing, since it only represents the interval between one musical pitch and another with double its frequency: hence, an octave is more related to the number~$2$ than to the number~$8$. 
Also, in the common western musical composition, the octave is composed by~$7$ notes, or~$12$ semitones, hence the octave seems to relate more to the numbers~$7$ and~$12$ rather than~$8$. The fact is that when talking about intervals of notes the tradition is to count from the first note to the final note included (which makes~$8$ notes, somewhat justifying the name of octave).}\index{octave} \label{JS-PGAGPA}
is the distance between one pitch and another with half or double its frequency:
thus, from~\eqref{CwUSTRJMSONGJMSD-6},
holding the string against the fingerboard halving its length has the effect of doubling the frequency of the root note (consistently with the fact that each fret on the guitar's fretboard corresponds\footnote{The choice of
dividing the octave in twelve semitones\index{semitone} is also mathematically well-grounded, since it corresponds
to frequency ratios of the form~$2^{\frac{j}{12}}$ for~$j\in\{0,\dots,12\}$.
On the one hand, these are irrational numbers (except for~$j=0$ and~$j=12$).
On the other, notes sound harmonious to our ear if the frequency of the notes is close to a simple interval:
for instance a frequency ratios such as~$\frac32$ (the ``perfect fifth''
in musical jargon), $\frac43$ (the ``fourth'') and~$\frac54$ (the ``major third'').
These two observations seem contradictory, till we realize
how ``close to rationals'' the irrational numbers of the form~$2^{\frac{j}{12}}$ are
(for~$j\in\{1,\dots,11\}$ these numbers are ``almost''~$\frac{16}{15}$,
$\frac98$, $\frac65$, $\frac54$, $\frac43$, $\frac75$, $\frac32$, $\frac85$, $\frac53$, $\frac{16}9$
and~$\frac{15}8$).
For example, 
$$ \frac32-2^{\frac{7}{12}}=1.5-1.49830707688...=0.00169292312...$$
showing how close the seventh semitone is to the perfect fifth.

While other tuning scales are possible and are indeed employed in several contexts,
the equal temperament one (i.e., the one splitting the scale into equal intervals)
presents several advantages, such as the possibility of transposing a tune into other keys
thus allowing instruments to play in all keys
(or singers to sing a song in a key which is more congenial to their voice)
with minimal flaws in intonation.
The twelve-tone equal temperament based on
the division of the octave into twelve equally spaced parts
on a logarithmic scale appears the most widespread system in music today.
Likely, the common adoption of this specific temperament
also influenced the composition of music in order to accommodate the system and minimize dissonance
coming from irrational approximation of rational ratios.}
to one semitone,
the distance of twelve semitones is the octave, and the position of the twelfth fret
corresponds to a string half as long as the original one).\medskip

\begin{figure}
  \centering
  \includegraphics[height=.15\textheight]{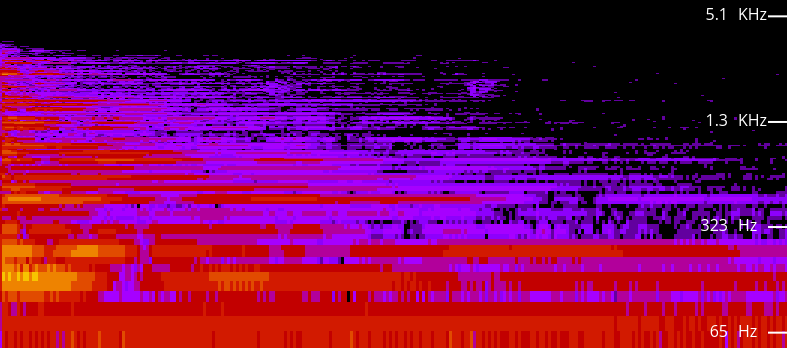}$\,$
    \includegraphics[height=.15\textheight]{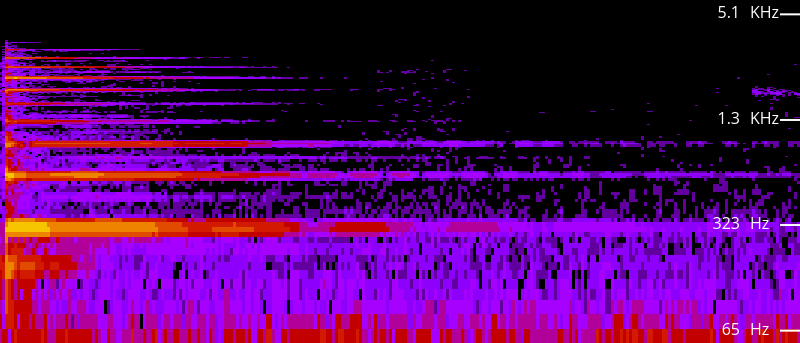}
 \caption{\sl Frequency plot, over time, of a metal E2 string (left)
 and a  synthetic E4 string (right).}\label{SIMSTR1kngGngeFI}
\end{figure}

Let us also mention a typical question asked by many beginners of
(say, classical) guitar:
\begin{equation}\label{CHANGE:STRI}
{\mbox{can I replace the steel strings with nylon strings?}}
\end{equation}
This question is often motivated by a practical reason:
\begin{equation}\label{sidofgOfRTYHJtt6P}
{\mbox{nylon strings are usually ``softer'' compared to steel strings,}}\end{equation}
therefore they are less likely to cause pain on beginners' fingers
(professional guitarists develop callouses on their fingertips with practice, so they
do not perceive this problem): instead,
steel strings are generally perceived by beginners
more uncomfortable and may require more effort to play
in terms of muscle tension in the hand and wrist.

Well, from the musical theory point of view,
there are technical issues about how the timbre
of the instrument would be affected by the change proposed in~\eqref{CHANGE:STRI}:
many guitarists would argue that messing
up with the proper setting of the guitar may well compromise its unique\footnote{And
here we confine
ourselves to a rough dichotomy between steel and nylon strings,
not going into finer distinctions and specific materials, such as
nickel, brass and bronze, on the metals side, catgut, on the natural fibers side,
and fluorocarbon polymers, on the synthetic materials side. Interestingly,
an important part of the development of modern nylon strings occurred in
late forties, thanks to the close interactions between luthier Albert Augustine,
chemical company DuPont, and virtuoso classical guitarist Andr\'es Segovia
(see Figure~\ref{9ijj0SIegib Ov801tangeFI} for a portrait of the latter, by Hilda Wiener).}
tone quality
(nylon strings tend to have a ``rounder'' and ``mellower'' tone,
inappropriate materials may make
the strings buzz excessively and slow down the playing action,
some strings lend themselves better to certain techniques than others,
e.g. strumming with a pick may work not too well with a nylon string,
which are instead often more appreciated by experienced fingerstylists,
while ordinarily finger picking the steel strings is considered to be easier, etc.).

But if
you are going to tolerate these performance issues, in principle the answer to~\eqref{CHANGE:STRI} is
yes, you can: {\em
placing nylon strings on a guitar built for steel strings will (more or less) work}.
\medskip

But, beware, not the other way around: {\em
steel strings on a guitar built for nylon strings can break the neck of the instrument due to the increased tension}.

Indeed, by the discussions after~\eqref{CwUSTRJMSONGJMSD-6},
we know that the higher the density, the lower the pitch.
Thus, putting a steel (denser) string in place of a nylon (less dense) one would take down the pitch:
to compensate the pitch drop, again by the discussions after~\eqref{CwUSTRJMSONGJMSD-6},
one will need to increase the tension of the string
(since we know that this would raise the pitch). So, all in all, to adjust the tuning, {\em replacing nylon with steel strings would require
an increased tension} which, on the long run, can damage the instrument.
\medskip

This fact, showcasing that nylon strings have
lower tension than their steel counterpart, explains also the phenomenon described
in~\eqref{sidofgOfRTYHJtt6P}, since, due to their reduced tension,
nylon strings are softer and more gentle on players' fingers.\medskip

The difference of tension between nylon and steel strings also justifies why
they affect differently the timbre of the instruments.
Roughly speaking, on the one hand,
with less tension, the nylon strings may move more freely, thus producing
more overtones, by
dislodging some energy on adjacent or resonant modes, while the
higher tension of steel strings may tend to produce sounds with a greater clarity
and a brighter effect.
On the other hand, denser strings have in principle
more potential energy, which gets converted into kinetic energy by string releasing:
this can also enhance some high frequency harmonics
in a steel string, providing its characteristic ``metallic''
sound (and, if one relates high frequency modes with more irregular
behaviors of the function that they describe, this also
justifies the fact that steel strings often
produce a vibration which may sound less ``round'',
or ``harder'', than their nylon counterpart). For a diagram comparing
the sound spectrum of steel and nylon strings see e.g. Figure~3.10 in~\cite{PANOS}.
\medskip

It is also instructive to compare the frequency diagrams simply obtained
by two pick strokes on a cheap guitar, one on a metal E2 string and one on a synthetic E4 string, see Figure~\ref{SIMSTR1kngGngeFI}. In this setting, both the strings correspond to the same note~E,
but the synthetic string plays two octave higher than the metal one (specifically, E4
should correspond approximately to~329.6276 Hz
while E2 should have frequency of approximately~82.40689 Hz, and note the ratio of approximately~$4$ between
the two frequencies, in agreement with their distance of two octaves).
We notice from Figure~\ref{SIMSTR1kngGngeFI} that indeed the
metal strings not only tend
to excite additional harmonics of higher frequencies
(and possibly too many of them, indicating
a poor quality of the instrument and an excessively forceful
pitch on the string),
but also their vibration is more persistent
in time: this fact is also in agreement with our intuition, since heavier metal strings
possess more (potential, hence kinetic) energy than lighter synthetic ones,
whence they keep oscillating longer since longer is the time needed to dissipate energy due to
friction.\medskip

Moreover,
it is experienced by some guitarists\footnote{Other musical instruments can
also be understood through the lens of mathematics, see e.g. {\tt http://newt.phys.unsw.edu.au/jw/flutes.v.clarinets.html}
for flutes and clarinets.}
that stiffer strings play ``faster'', namely steel strings with
high tension
have a ``quicker attack'', while their nylon counterpart a slower one (the attack being
the time it takes from the pick of the string to the reaching of the loudest peak of the sound, i.e.
the time one begins to hear the note after the string is released).
\medskip

\begin{figure}
  \centering
  \includegraphics[width=.30\linewidth]{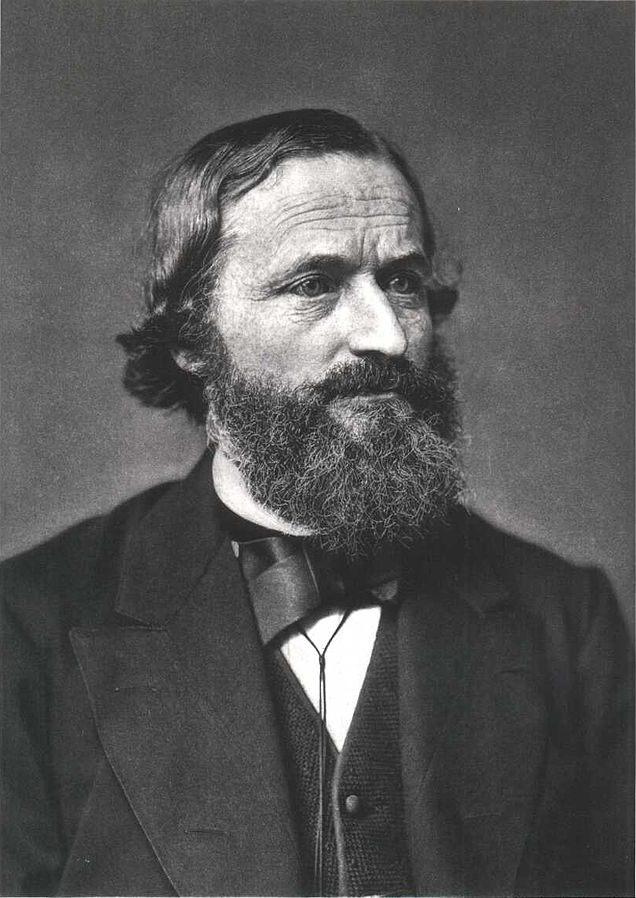}
 \caption{\sl Gustav Kirchhoff
 (Public Domain image from
        Wikipedia, source Smithsonian Libraries).}\label{SIMKIRAoGREENFIDItangeFI}
\end{figure}

For completeness, let us also mention a slightly different model for vibrating strings proposed\footnote{Besides the nonlocal equations for vibrating strings related to~\eqref{leng2tht346y356ef7o8om-ely},
Gustav Kirchhoff provided fundamental contributions
to the theory of electrical circuits, spectroscopy, and black-body radiation.
Which creates quite a confusion, because the names of
Kirchhoff's Law and Kirchhoff equation end up being used in all these topics, with different meanings.

By the way, Kirchhoff's home town was K\"onigsberg, a
historic Prussian city that is now Kaliningrad, in Russia's Kaliningrad Oblast
(a small exclave of the vast Russian state
providing an interesting counterexample to the conjecture that political states are connected
regions, up to additional islands). The reason for which we mention K\"onigsberg is because
one classical and very famous mathematical problem that was settled by Euler is related to the bridges of
K\"onigsberg (see e.g.~\cite{MR2053120} to know more about the K\"{o}nigsberg bridges).} by Gustav Robert Kirchhoff in~\cite{zbMATH02717254}, see Figure~\ref{SIMKIRAoGREENFIDItangeFI}.
The gist of this model is that the elastic coefficient~$\kappa$ in~\eqref{CUSTRJMSONGJMSD-3}
is assumed to be a positive constant but, in reality, the elastic properties of the string may
depend on its elongation (according to the material, the string could either become ``stiff'' 
or ``slack off'' for large
elongations). To account for such a possibility, Kirchhoff proposed that~$\kappa$
instead depends on the length of the deformed string, namely on
\begin{equation}\label{leng2tht346y356ef7o8om-ely}
{\mbox{\footnotesize{\calligra{L}}}}\;\,(t):=\int_0^L \sqrt{1+\big(\partial_x u(x,t)\big)^2}\,dt.\end{equation}
In this model, $\kappa$ in~\eqref{CUSTRJMSONGJMSD-3} must be interpreted as~$\kappa\left(
{\mbox{\footnotesize{\calligra{L}}}}\;\,(t)\right)$, that is a function of the length of the string
(the corresponding equation is thus called Kirchhoff equation). \index{Kirchhoff equation}
Notice that not only in this situation~$\kappa$ varies with time, but also
it depends on a ``nonlocal'' quantity. Indeed, the length of the deformed string in~\eqref{leng2tht346y356ef7o8om-ely} is a ``global'' object: to calculate it, it is not sufficient to
know the shape of the string at a given point, instead full information about its large scale geometry
is required. And, of course, mathematical problems requiring the knowledge of nonlocal
quantities become structurally harder (but often quite interesting, since they aim at
capturing the Big Picture).

\subsection{Elastic membranes}

The description of elastic membranes\index{elastic membrane} can be seen as a multi-dimensional version
of the one of vibrating strings presented in Section~\ref{GUITAR}.
A simple model could be that of considering the membrane as constituted by
a web of infinitesimal strings in the coordinate directions.
In this way, one may think that the vertical force exerted
at a point~$(x,u(x))\in\R^n\times\R$ of the membrane
is the sum of the forces produced by the tension of the infinitesimal strings
in each direction, plus gravity.
Recalling~\eqref{TOTFOSTR}, we may think that, for each~$i\in\{1,\dots,n\}$,
the vertical force exerted by the tension of the infinitesimal string
located in direction~$e_i$ (say, located from~$x-\delta e_i$ and~$x+\delta e_i$) is equal to
$$ 2\kappa \delta T(x)\partial_i\left(
\frac{\partial_i u(x,t)}{\sqrt{1+\big(\partial_i u(x,t)\big)^2}}\right),$$
up to higher order in~$\delta$.

Hence the total vertical force acting on the membrane at~$(x,u(x))\in\R^n\times\R$
is given by
$$ 2\kappa \delta T(x)\sum_{i=1}^n\partial_i\left(
\frac{\partial_i u(x,t)}{\sqrt{1+\big(\partial_i u(x,t)\big)^2}}\right)-m_\delta g,$$
where, in view of~\eqref{MEPSI},
$$m_\delta
=2\rho\delta \sum_{i=1}^n\sqrt{1+\big(\partial_i u(x,t)\big)^2},$$
up to higher orders in~$\delta$, being~$\rho$ the linear densities of the infinitesimal strings in each direction
(that we suppose to be constant and independent of the direction).
In this setting, Newton's Second Law gives that
\begin{equation}\label{CUSTRJMSONGJMSD-287729ou35t4jrg}
\begin{split}
\partial_{tt}u&=\frac{2\kappa \delta T}{m_\delta}\sum_{i=1}^n\partial_i\left(
\frac{\partial_i u}{\sqrt{1+(\partial_i u)^2}}\right)-g\\
&=\frac{\kappa T}{\displaystyle
\rho\sum_{i=1}^n\sqrt{1+(\partial_i u)^2}}
\sum_{i=1}^n\partial_i\left(
\frac{\partial_i u}{\sqrt{1+(\partial_i u)^2}}\right)-g,
\end{split}\end{equation}
which can be considered the higher dimensional\footnote{{F}rom the geometric point
of view, it is interesting to observe that the term~$\displaystyle
\sum_{i=1}^n\partial_i\left(
\frac{\partial_i u}{\sqrt{1+(\partial_i u)^2}}\right)$ in~\eqref{CUSTRJMSONGJMSD-287729ou35t4jrg}
detects the mean curvature \index{mean curvature}
of the membrane
(remarks of this type will be formalized in a general setting
in Theorem~\ref{INCOO}).}
version of~\eqref{CUSTRJMSONGJMSD-2}.
In the small elongation approximation we thus consider~$T:=1$ and disregard the terms that
are quadratic in the displacement~$u$, thus reducing~\eqref{CUSTRJMSONGJMSD-287729ou35t4jrg}
to
\begin{equation}\label{CUSTRJMSONGJMSD-287729ou35t4jrg24qwerty} \partial_{tt}u=\frac{ \kappa }\rho\Delta u-g,\end{equation}
which is the higher dimensional version of~\eqref{CUSTRJMSONGJMSD-3}.

In particular, an elastic membrane in equilibrium is a stationary solution of~\eqref{CUSTRJMSONGJMSD-287729ou35t4jrg24qwerty} and thus satisfies
\begin{equation}\label{CUSTRJMSONGJMSD-287729ou35t4jrg24qwerty44} \frac{ \kappa }\rho\Delta u=g,\end{equation}
which is a form of the so-called Poisson's equation (compare with~\eqref{DAGB-ADkrVoiweLL4re2346ytmngrrUj9}).

In the absence of gravity, \eqref{CUSTRJMSONGJMSD-287729ou35t4jrg24qwerty44} reduces to
the so-called Laplace's equation (compare with~\eqref{DAGB-ADkrVoiweLL4re2346ytmngrrUj10})
\begin{equation*}
\Delta u(x)=0.\end{equation*}
This gives that (at least in the linear approximation) elastic membranes constrained at their boundaries
are described by the graph of harmonic functions.
Once again, this fact brings us a useful hint, suggesting a smoothing effect
for harmonic functions, somewhat inherited from the fact that elastic membranes
do not develop spikes (the elasticity would indeed reduce the peak
to find a balance): this intuition will be mathematically formalized, developed 
and consolidated in Chapter~\ref{CHA-2}.

\subsection{Elasticity theory and torsion of bars}\label{RODSBACKSE}

\begin{figure}
  \centering
  \includegraphics[width=.45\linewidth]{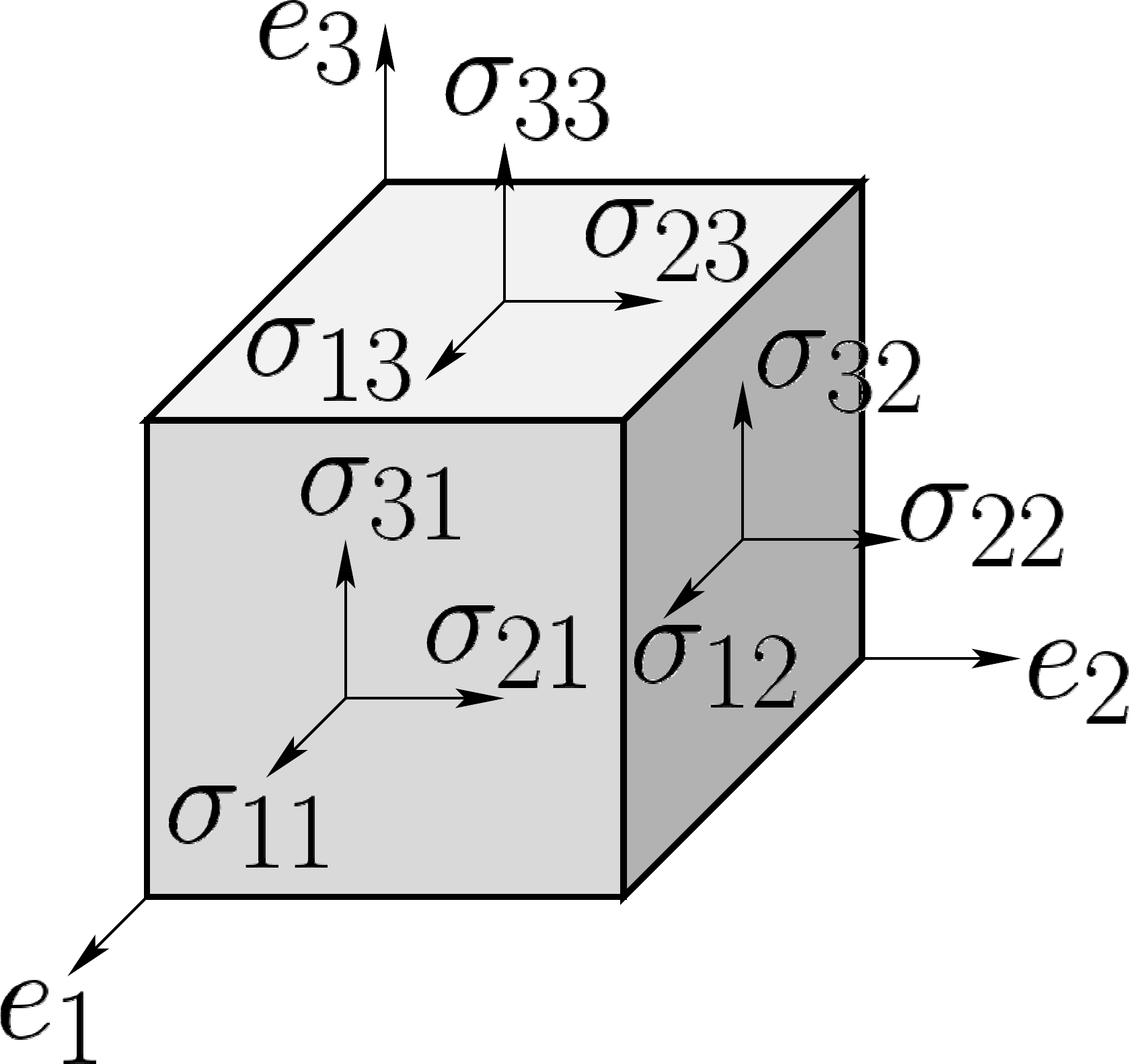}
 \caption{\sl The stress tensor.}\label{STREDItangeFI}
\end{figure}

A topical argument in material sciences deals with the determination of the shape of objects
subject to external forces and to the understanding of the relations between the forces
applied and the resulting deformations. Providing a comprehensive account of this elasticity theory goes well beyond
the scopes of these notes (see e.g.~\cite{MR2779440} and the references therein for further readings):
here we just introduce, in a rather simplified form, some basic concepts which will lead us to the study
of an interesting\footnote{Actually, the final equation that we will obtain in~\eqref{EUMSD-OS-32456i7DNSSE234R-BNMBAR}
is precisely the same as the one arising from fluid dynamics in~\eqref{EUMSD-OS-32456i7DNSSE234R-BNM}. This coincidence shows
once again the extraordinary unifying power of mathematics.}
elliptic equation (which actually will be studied in detail in Sections~\ref{SEC:OVER-SERRINTheorem} and~\ref{BACKSE}).
For concreteness, we focus on the three-dimensional case and we consider an infinitesimal portion of the material under consideration,
modeled as a very tiny cubes with faces oriented along the coordinate axes, see Figure~\ref{STREDItangeFI}.
If the cube has the form~$(0,\delta)^3$, for some small~$\delta>0$, we denote by~$S_1$, $S_2$ and~$S_3$
three facets of the cubes in direction~$e_1$, $e_2$ and~$e_3$ respectively, i.e.
$$ S_1:=\{\delta\}\times(0,\delta)\times(0,\delta),\qquad
S_2:=(0,\delta)\times\{\delta\}\times(0,\delta)\qquad{\mbox{and}}\qquad
S_3:=(0,\delta)\times(0,\delta)\times\{\delta\}.$$
For each~$i$, $j\in\{1,2,3\}$ we denote by~$\sigma_{ij}$ the force per unit of area applied in direction~$e_i$ to the surface~$S_j$.
In jargon, the matrix~$\{\sigma_{ij}\}_{i,j\in\{1,2,3\}}$ is sometimes called ``stress tensor''\index{stress tensor}.

By Newton's Third Law, at equilibrium the forces acting on the above ideal infinitesimal cube
must balance, therefore the forces acting on the facets of the cubes opposite to~$S_1$ (or~$S_2$, or~$S_3$) are the ones discussed above, but with a sign changed.
In particular, given~$i\in\{1,2,3\}$, at equilibrium the total forces per unit of area applied in direction~$e_i$ must balance out
and therefore (noticing that if~$x\in S_j$ then~$x-\delta e_j$ belongs to the facet opposite to~$S_j$) 
$$ 0=\sum_{j=1}^3\Big(\sigma_{ij}(x)-\sigma_{ij}(x-\delta e_j)\Big)=
\delta \sum_{j=1}^3\partial_j\sigma_{ij}(x)+o(\delta),$$
leading to the following balance prescription for the stress tensor:
\begin{equation}\label{356um1ijkialivj}
\sum_{j=1}^3\partial_j\sigma_{ij}=0\qquad{\mbox{for all }}\;i\in\{1,2,3\}.
\end{equation}

Additionally, the angular momentum contribution~$\vartheta_1$ coming from~$S_1$ is proportional to the vector product
between the vectorial force per unit of area on~$S_1$, corresponding to~$\sigma_{11}e_1+\sigma_{21}e_2+\sigma_{31}e_3$,
and the vector joining normally the center of the cube to~$S_1$, namely~$e_1$. Hence, up to normalizing constants,
$$ \vartheta_1=\det\left(\begin{matrix}
e_1 & e_2 & e_3\\
\sigma_{11} & \sigma_{21} & \sigma_{31}\\
1 & 0 &0
\end{matrix}\right)=\sigma_{31}e_2-\sigma_{21}e_3.
$$
Similarly,
the angular momentum contribution~$\vartheta_2$ coming from~$S_2$ is proportional to the vector product
between the vectorial force per unit of area on~$S_2$, corresponding to~$\sigma_{12}e_1+\sigma_{22}e_2+\sigma_{32}e_3$,
and the vector joining normally the center of the cube to~$S_2$, namely~$e_2$, whence
$$ \vartheta_2=\det\left(\begin{matrix}
e_1 & e_2 & e_3\\
\sigma_{12} & \sigma_{22} & \sigma_{32}\\
0 & 1 &0
\end{matrix}\right)=-\sigma_{32}e_1+\sigma_{12}e_3.
$$
And also the angular momentum contribution~$\vartheta_3$ coming from~$S_3$ is proportional to the vector product
between the vectorial force per unit of area on~$S_3$, corresponding to~$\sigma_{13}e_1+\sigma_{23}e_2+\sigma_{33}e_3$,
and the vector joining normally the center of the cube to~$S_2$, namely~$e_3$, whence
$$ \vartheta_3=\det\left(\begin{matrix}
e_1 & e_2 & e_3\\
\sigma_{13} & \sigma_{23} & \sigma_{33}\\
0 & 0 &1
\end{matrix}\right)=\sigma_{23}e_1-\sigma_{13}e_2.
$$
On this account, the
conservation of the total angular momentum entails that, at equilibrium,
$$ 0=\vartheta_1+\vartheta_2+\vartheta_3=(\sigma_{23}-\sigma_{32})e_1
+(\sigma_{31}-\sigma_{13})e_2
+(\sigma_{12}
-\sigma_{21})e_3.$$
This establishes that the stress tensor is symmetric, namely
\begin{equation}\label{SYMME}
\sigma_{ij}=\sigma_{ji}\qquad{\mbox{for all }}\,i,j\in\{1,2,3\}.
\end{equation}
Now, so as to model the distortion of the material caused by these forces, we consider a displacement vector~$v$
describing at any point the change in the configuration of the given body. 
The displacement vector itself
may not be the most interesting object to take into account in an elasticity theory, since
translations and rotations do produce significant displacements without affecting the shape of the body under consideration.
It is therefore common to try to detect the deformation effects that displacements produce
by accounting for relative displacements. Namely, if a point~$x$ is moved to~$S(x):=x+v(x)$, given~$\omega\in\partial B_1$,
the point~$y:=x+\delta\omega$ is moved to~${\mathcal{P}}(y)=y+v(y)=x+\delta\omega+v(x+\delta\omega)$. To compute the corresponding
relative displacement, one checks how~$y-x$ is affected by this transformation, namely we point out that
\begin{eqnarray*}
{\mathcal{P}}(y)-{\mathcal{P}}(x)={\mathcal{P}}(x+\delta \omega)-{\mathcal{P}}(x)=\delta\,D{\mathcal{P}}(x)\omega+O(\delta^2)=D{\mathcal{P}}(x)(y-x)+O(\delta^2).\end{eqnarray*}
Namely, the ``deformation gradient''~$D{\mathcal{P}}$ is a linear approximation measuring relative displacements.

If we want to measure the change of square length, we have that
\begin{equation}\label{UJSvherEDIMAPOMKE}\begin{split}& |{\mathcal{P}}(y)-{\mathcal{P}}(x)|^2-|y-x|^2=|D{\mathcal{P}}(x)(y-x)+O(\delta^2)|^2-|y-x|^2\\&\qquad=
\big(D{\mathcal{P}}(x)(y-x)+O(\delta^2)\big)\cdot\big(D{\mathcal{P}}(x)(y-x)+O(\delta^2)\big)-|y-x|^2\\&\qquad=
\big(D{\mathcal{P}}(x)(y-x)\big)\cdot\big(D{\mathcal{P}}(x)(y-x)\big)-|y-x|^2+O(\delta^3)\\&\qquad=
\big(D{\mathcal{P}}(x)\big)^\top D{\mathcal{P}}(x)(y-x)\cdot(y-x)-|y-x|^2+O(\delta^3)\\&\qquad=\Big(
\big(D{\mathcal{P}}(x)\big)^\top D{\mathcal{P}}(x)-{\rm Id}\Big)(y-x)\cdot(y-x)+O(\delta^3),
\end{split}\end{equation}
where~${\rm Id}$ is the identity matrix and the superscript~``$\top$'' denotes matrix transposition.

It is also useful to observe that
\begin{eqnarray*}&&
D{\mathcal{P}}^\top D{\mathcal{P}}-{\rm Id}=\big({\rm Id}+Dv\big)^\top\big({\rm Id}+Dv)-{\rm Id}=
Dv+Dv^\top+Dv^\top Dv.
\end{eqnarray*}
In particular, if one assumes that the displacements are small, the quadratic term~$Dv^\top Dv$
is negligible with respect to the ``symmetric gradient''~$Dv+Dv^\top$ which is of linear type.
Therefore, it is customary to consider the approximation~$
D{\mathcal{P}}^\top D{\mathcal{P}}-{\rm Id}\simeq Dv+Dv^\top$ and rewrite~\eqref{UJSvherEDIMAPOMKE} as
\begin{equation*}
|{\mathcal{P}}(y)-{\mathcal{P}}(x)|^2-|y-x|^2\simeq \Big(Dv(x)+\big(Dv(x)\big)^\top\Big)
(y-x)\cdot(y-x)
\end{equation*}
or equivalently
\begin{equation}\label{STRTENDE}
|{\mathcal{P}}(y)-{\mathcal{P}}(x)|^2-|y-x|^2\simeq \sum_{i,j=1}^3 \e_{ij}(x)
(y_i-x_i)(y_j-x_j),
\end{equation}
where
\begin{equation}\label{STRTENDE-2}
\e_{ij}(x):=\Big(Dv(x)+\big(Dv(x)\big)^\top\Big)_{ij}
=\partial_i v_j(x)+\partial_j v_i(x).
\end{equation}
In jargon,~$\e_{ij}$ is called the ``strain tensor''\index{strain tensor} (and notice the symmetry feature that~$\e_{ij}=\e_{ji}$)
and~\eqref{STRTENDE} expresses in a quantitative way (though after some simplifying approximation)
the relative change in the position of points under the deformation.

We now relate strain and stress.
To this end, once again,
the simplest possible ansatz, leading to a linear theory of elasticity, is now that
deformations (or, more precisely, strains, to avoid translations and rotations) are directly proportional
to the forces excerpt (this is indeed a possible rephrasing of Hooke's Law for linear elastic materials).
That is, neglecting proportionality constants, we will suppose that
\begin{equation}\label{STRESTAP}
\e_{ij}=\sigma_{ij}.\end{equation}
It is interesting to remark that this identity is also algebraically compatible with the fact that
both the tensors above are symmetric.
\medskip

\begin{figure}
  \centering
  \includegraphics[width=.5\linewidth]{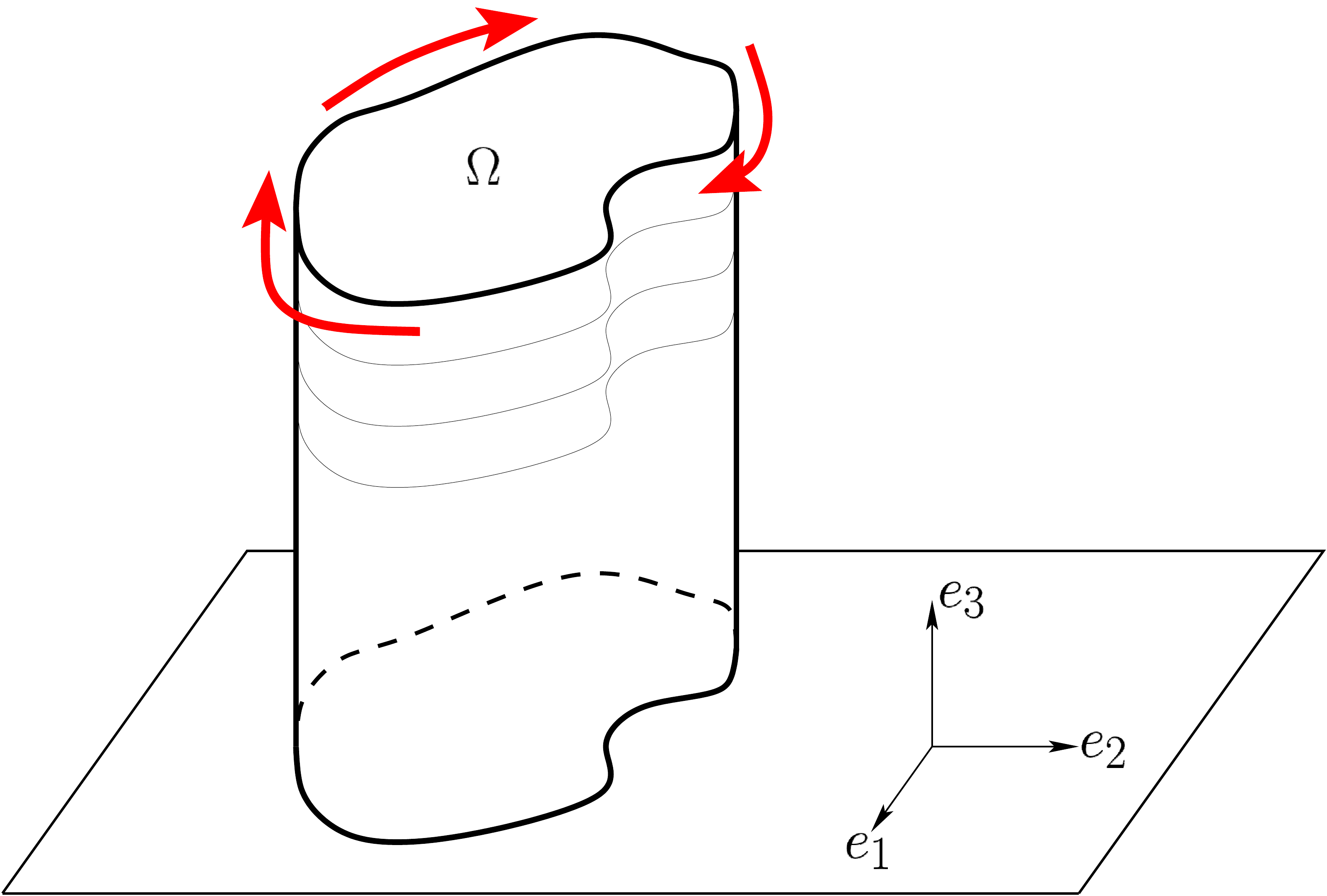}
 \caption{\sl A bar subject to torsion.}\label{BARSTREDItangeFI}
\end{figure}

We now apply the general framework of the linear theory of elasticity showcased so far
to the special situations of straight bars subject to torsion. This will be indeed one of the founding motivations
in Sections~\ref{SEC:OVER-SERRINTheorem} and~\ref{BACKSE}. To model a straight bar, we consider
a body of the form
\begin{equation}\label{HANXTHECRPR4RO4O}
\Omega\times[0,+\infty),\end{equation} with~$\Omega\subset\R^{2}$
smooth and contractible, see Figure~\ref{BARSTREDItangeFI}.

We assume that such a bar is constrained at the level~$\{x_3=0\}$ but it gets ``twisted'' from the top.
The effect of the twist would be to make the cross section rotate by some angle, depending on the level that we take into account.
Assuming again a linear behavior, we suppose that the rotation performed at level~$\{x_3=\zeta\}$
is proportional to~$\zeta$ (i.e., the angle of rotation increases linearly with the height,
which is not so unreasonable, at least for small heights, since the bar is constrained at ground zero).
Accordingly, we model the twist at level~$\{x_3=\zeta\}$ by a rotation of angle~$\Theta_\zeta= c_0\zeta$, for some~$c_0>0$.
The corresponding displacement sends therefore the point~$x=(x_1,x_2,\zeta)$
to a point~$z:=x+v(x)$, with the horizontal component of~$v$ being prescribed by the above rotation by the angle~$\Theta_\zeta$,
namely, for small heights~$\zeta$,
\begin{eqnarray*} \left( \begin{matrix} z_1\\z_2\end{matrix}\right)&=&
\left( \begin{matrix} \cos\Theta_\zeta&\sin\Theta_\zeta\\-\sin\Theta_\zeta&\cos\Theta_\zeta\end{matrix}\right)
\left( \begin{matrix} x_1\\x_2\end{matrix}\right)\\&=&
\left( \begin{matrix} \cos(c_0\zeta)&\sin(c_0\zeta)\\-\sin(c_0\zeta)&\cos(c_0\zeta)\end{matrix}\right)
\left( \begin{matrix} x_1\\x_2\end{matrix}\right)\\&=&
\left( \begin{matrix} 1+O(\zeta^2)& c_0\zeta+O(\zeta^2)\\-c_0\zeta+O(\zeta^2)&1+O(\zeta^2)\end{matrix}\right)
\left( \begin{matrix} x_1\\x_2\end{matrix}\right).
\end{eqnarray*}
This gives that
\begin{equation*} 
\left( \begin{matrix} v_1(x)\\v_2(x)\end{matrix}\right)=
\left( \begin{matrix} z_1\\z_2\end{matrix}\right)-
\left( \begin{matrix} x_1\\x_2\end{matrix}\right)
=
\left( \begin{matrix} O(\zeta^2)& c_0\zeta+O(\zeta^2)\\-c_0\zeta+O(\zeta^2)&O(\zeta^2)\end{matrix}\right)
\left( \begin{matrix} x_1\\x_2\end{matrix}\right)=
\left( \begin{matrix} c_0\zeta x_2\\-c_0\zeta x_1\end{matrix}\right)+O(\zeta^2).
\end{equation*}
Hence, by sloppily\footnote{This just because here we are playing around with motivations.
{F}rom next chapter on, this kind of sloppiness will no longer be tolerated!} neglecting the higher order terms, and using the fact that~$x=(x_1,x_2,x_3)\in\{x_3=\zeta\}$,
we will simply write that
\begin{equation}\label{STRESTAP-3} 
v_1(x)= c_0 x_2 x_3\qquad{\mbox{and}}\qquad v_2(x)=-c_0 x_1 x_3.
\end{equation}

As for the vertical component of the displacement, a common ansatz is to consider it independent of the height:
roughly speaking, the torsion of the bar may exhibit a distortion of the cross section out of its own plane
(known in jargon as ``warping deformation''\index{warping deformation}), but usually this effect is less visible than other types of deformations
and we are assuming here that each layer of the bar (i.e., each horizontal cross section) ``equally pushes
the next one'', regardless of its height (once again, this is not a completely unrealistic assumption
especially if we are restricting our attention to the layers close to ground zero). To model that this warping effect
does not depend on the height, we assume that the vertical component of the displacement vector is independent of~$x_3$, that is
\begin{equation}\label{STRESTAP-4} v_3(x)=w(x_1,x_2),\end{equation}
for some function~$w$.

Now, from~\eqref{STRTENDE-2}, \eqref{STRESTAP}, \eqref{STRESTAP-3} and~\eqref{STRESTAP-4} we arrive at
$$ \sigma_{13}=\e_{13}=\partial_1 v_3+\partial_3 v_1 =\partial_1 w+c_0x_2.$$
Similarly,
$$ \sigma_{23}=\e_{23}=\partial_2 v_3+\partial_3 v_2 =\partial_2 w-c_0x_1$$
and
$$ \sigma_{33}=\e_{33}=2\partial_3 v_3=0.$$
Interestingly, the above expressions for~$\sigma_{i3}$ with~$i\in\{1,2,3\}$ are all independent of the height~$x_3$
(which was possibly not obvious from the beginning). Furthermore, recalling the
balance prescription in~\eqref{356um1ijkialivj} and the symmetry property in~\eqref{SYMME},
\begin{eqnarray*}&& 0=\sum_{j=1}^3\partial_j\sigma_{3j}=\sum_{j=1}^3\partial_j\sigma_{j3}=\partial_1\sigma_{13}+\partial_2\sigma_{23}.
\end{eqnarray*}
That is, the differential form~$\sigma_{23}\,dx_1-\sigma_{13}\,dx_2$ is closed, and therefore exact, thanks to the
Poincar\'e Lemma (see e.g.~\cite[Theorem~8.3.8]{MR1209437}). This yields that there exists a function~$\phi=\phi(x_1,x_2)$,
sometimes refereed to with the name of ``warping potential''\index{warping potential}, such that~$d\phi=\sigma_{23}\,dx_1-\sigma_{13}\,dx_2$.

This gives that~$\partial_1\phi=\sigma_{23}$ and~$\partial_2\phi=-\sigma_{13}$. Consequently,
\begin{equation}\begin{split}\label{TRUNPDOmrfgEPPJO}&
\Delta\phi=\partial_1\sigma_{23}-\partial_2\sigma_{13}=
\partial_1(\partial_2 w-c_0x_1)-\partial_2(\partial_1 w+c_0x_2)\\&\qquad\qquad\qquad=
\partial_{12} w-c_0-\partial_{12} w-c_0=-2c_0.\end{split}
\end{equation}

Another notion of interest for engineering purposes is the ``traction''\index{traction}~${\mathcal{T}}$ along the boundary of the bar,
which is defined as the normal component of the vertical stress. More explicitly, if we define~$\sigma_3:=(\sigma_{31},\sigma_{32},\sigma_{33})$ and we consider the outer unit normal~$\nu^\star$ of the bar, we set
\begin{equation}\label{CA2LT} {\mathcal{T}}:=\sigma_3\cdot\nu^\star.\end{equation}
As a matter of fact, since the straight bar has the form given in~\eqref{HANXTHECRPR4RO4O},
we can identify~$\nu^\star\in\R^3$ with the vector~$(\nu,0)$, being~$\nu\in\R^2$ the 
outer unit normal at~$\partial\Omega$ and consequently
$$ {\mathcal{T}}=\sigma_{31} \nu_1+\sigma_{32}\nu_2=\sigma_{13} \nu_1+\sigma_{23}\nu_2=
-\partial_2\phi\nu_1+\partial_1\phi\nu_2=\nabla\phi\cdot\tau,
$$
where~$\tau:=(-\nu_1,\nu_2)$ can be taken as the unit tangent vector along~$\partial\Omega$ (say, clockwise oriented).
In this spirit, we find that the traction is the tangential derivative of the warping potential.

In the special situation in which the bar is subject to no traction, it follows that~$\nabla\phi\cdot\tau=0$ on~$\partial\Omega$. This gives that~$\phi$ is constant along~$\partial\Omega$, because if~$\partial\Omega$ is locally parameterized by a curve~$\gamma:(-1,1)\to\R^2$ with~$\dot\gamma$ proportional to~$\tau$, we have that
$$ \frac{d}{dt}\phi(\gamma(t))=\nabla\phi(\gamma(t))\cdot\dot\gamma(t)=0.$$
Thus, if the traction of the bar vanishes, we can write that
\begin{equation}\label{TRUNPDOmrfgEPPJO2}
{\mbox{$\phi=c_2$ on~$\partial\Omega$,}}\end{equation}
for some~$c_2\in\R$.

Recalling~\eqref{CA2LT}, we call the magnitude of the traction on the surface of the bar~$\partial\Omega$ the scalar~$|\sigma_3|$.
We remark that
$$ |\nabla\phi|=\sqrt{ (\partial_1\phi)^2+(\partial_2\phi)^2}=\sqrt{\sigma_{23} ^2+\sigma_{13}^2}=
\sqrt{\sigma_{23} ^2+\sigma_{13}^2+\sigma_{33}^2}=|\sigma_3|$$
therefore the magnitude of the traction coincides also with the norm of the gradient of the warping potential.

Thus, a consequence of~\eqref{TRUNPDOmrfgEPPJO2} is that~$|\partial_\nu\phi|=|\nabla\phi|$,
hence in this situation the magnitude of the traction coincides also (possibly up to a sign) with the normal
derivative of the warping potential along the surface of the bar.

As a consequence, in the special situation in which the bar is subject to no traction and
the magnitude of the traction is constant, we have that
\begin{equation}\label{TRUNPDOmrfgEPPJO2b}
{\mbox{$\partial_\nu\phi=c_3$ on~$\partial\Omega$,}}\end{equation}
for some~$c_3\in\R$.

Thus, if we define
$$ u:=\frac{c_2-\phi}{2c_0},$$
we deduce from~\eqref{TRUNPDOmrfgEPPJO},
\eqref{TRUNPDOmrfgEPPJO2} and~\eqref{TRUNPDOmrfgEPPJO2b}
that for a straight bar
subject to torsion and no traction, if the magnitude of the traction along the surface of the bar is constant then there exists a solution of
\begin{equation}\label{EUMSD-OS-32456i7DNSSE234R-BNMBAR}
\begin{dcases}
\Delta u=1&{\mbox{ in }}\Omega,\\
u=0&{\mbox{ on }}\partial\Omega,\\
\partial_\nu u=c&{\mbox{ on }}\partial\Omega.
\end{dcases}\end{equation}
Notice the perfect coincidence\footnote{As noticed in footnote~\ref{OJSLLAFOr5O3OSMEr3r}
on page~\pageref{OJSLLAFOr5O3OSMEr3r},
equation~\eqref{EUMSD-OS-32456i7DNSSE234R-BNMBAR} possesses a solution when~$\Omega$
is a disk, that is when the bar has a circular cross section.
We will discuss in Sections~\ref{SEC:OVER-SERRINTheorem} and~\ref{BACKSE}
whether bars of other shapes maintain the same property of presenting no traction and constant
traction magnitude along their surface.}
of this set of prescriptions with that in~\eqref{EUMSD-OS-32456i7DNSSE234R-BNM},
which also reveals a possibly unexpected connection between the dynamics of viscous fluids
and the elastic reactions of bars subject to torsion.

\begin{figure}
  \centering
  \includegraphics[width=.5\linewidth]{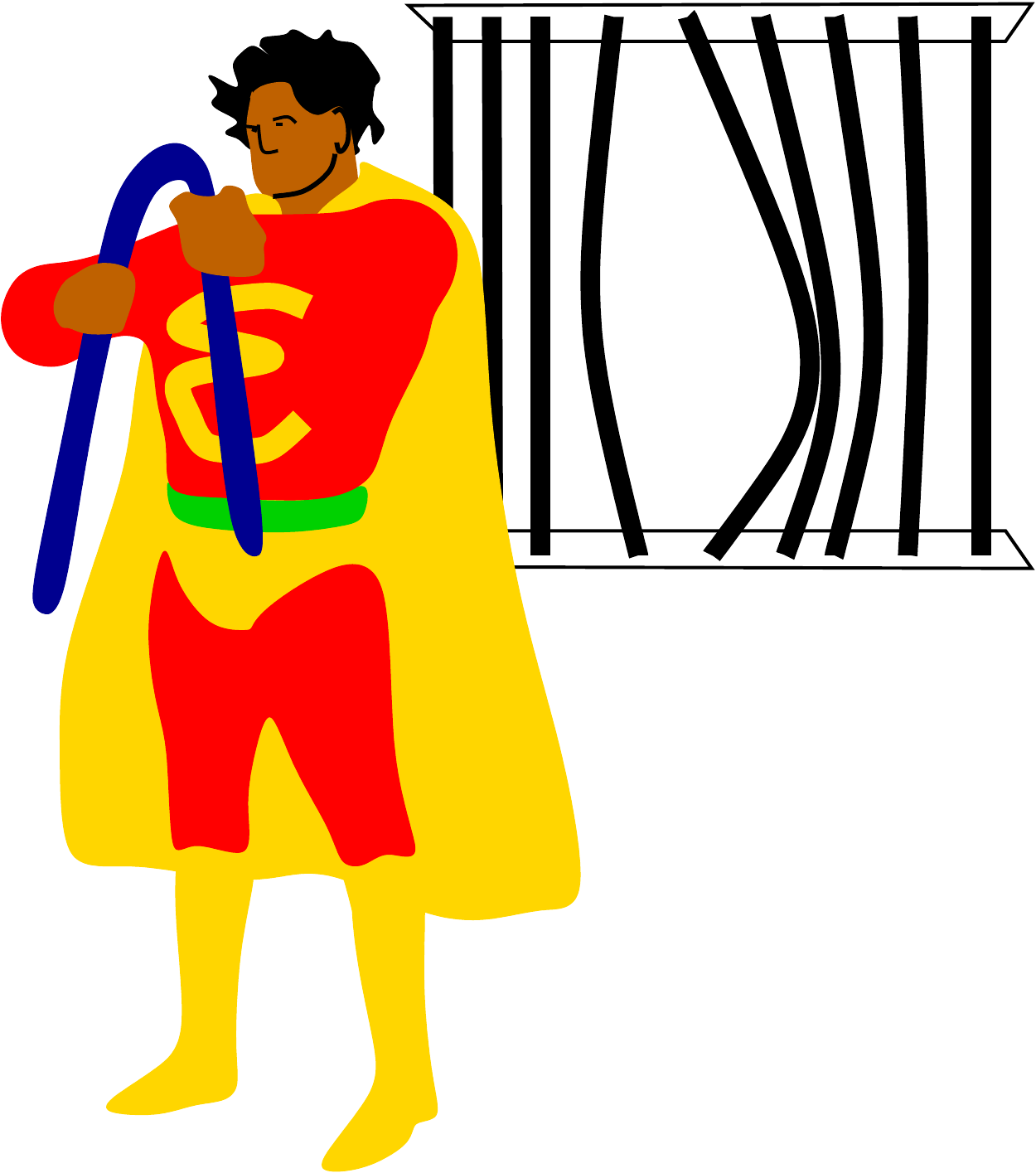}
 \caption{\sl Bending steel beams with bare hands by the power of~$\e$.}\label{T5w36FGHSN0dEcdhPSUPERSdhbM356U2LSMSUJNS-G5O2F4O2R}
\end{figure}

\subsection{Bending beams and plates}

\begin{figure}
                \centering
                \includegraphics[width=.6\linewidth]{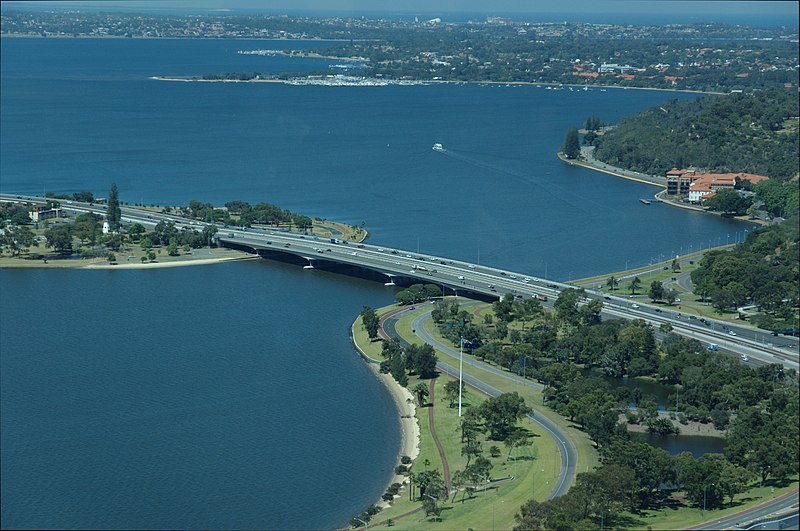}
        \caption{\sl Narrows Bridge, Perth (image by Speddie23 from
        Wikipedia,
        licensed under the Creative Commons Attribution-Share Alike 4.0 International license).}\label{2HPERNANFOJED231}
\end{figure}

\begin{figure}
                \centering
                \includegraphics[width=.55\linewidth]{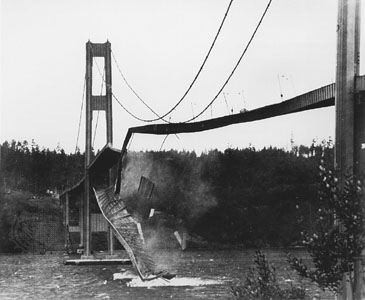}
        \caption{\sl The Tacoma Bridge falling into the strait (Public Domain image from
        Wikipedia).}\label{LOTACOMAOJED}
\end{figure}

Now we deal with deflection of beams, which is a classical topic in all superhero comics, see Figure~\ref{T5w36FGHSN0dEcdhPSUPERSdhbM356U2LSMSUJNS-G5O2F4O2R}.
Besides, the question has obvious applications\footnote{Sometimes mathematics plays an important role also in tragic events, such as the collapse of structures. One of the most famous events of this type was the collapse of the Tacoma Bridge in~1940, only four months after its opening, see Figure~\ref{LOTACOMAOJED}.

Fortunately, there was only one fatality of the Tacoma Bridge disaster, namely Tubby, the dog (a black male Cocker Spaniel) of a journalist who was left in a car on the bridge.
People attempted to rescue Tubby, who was however too terrified to leave the car and bit one of the rescuers.
Tubby died when the bridge fell and neither his body nor the car was ever recovered, probably because
the swift tides quickly moved the car away from the ruins of the bridge.

See {\tt https://en.wikipedia.org/wiki/File:Tacoma\_Narrows\_Bridge\_destruction.ogv} for an impressive video of the Tacoma Bridge collapse, as well as~\cite{MR4468373} and the references therein
for a mathematical investigation of the causes of the disaster.}
in engineering, see e.g. Figure~\ref{2HPERNANFOJED231}.
The beam equation \index{beam equation}
\begin{equation}\label{BARREbi96}
u''''=q
\end{equation}
describes\footnote{The linear model of elasticity accounting for the load-carrying and deflection of beams
is often called ``Euler-Bernoulli beam theory''.
No confusion should arise between~\eqref{BARREbi96} and the model presented
in Section~\ref{RODSBACKSE} where instead the torsion theory
of bars, and not the deformation under a load, was taken into account.}
the small deflection~$u$ of a homogeneous beam at equilibrium in terms of the applied load~$q$.

In~\eqref{BARREbi96}, $u$ and~$q$ are functions of one variable~$x\in\R$ (the beam is assumed to be infinitely long).
In a nutshell, the derivation of~\eqref{BARREbi96} from prime principles relies on balance of forces and moments.
Roughly speaking, the bending of the beam produces a compression force on one side
of the beam and a tensile stress on the other side: these deformation
stresses will be modeled via an elastic Hooke's Law and produce
turning forces which, at equilibrium, together with the corresponding moments,
need to be balanced with the external load on the beam.

The details go as follows. We assume that the beam is displayed along the horizontal axis and slightly deformed
into a graph of the form~$z=u(x)$ (the beam may well be three-dimensional, in which case we assume
that its
transversal sections possess some given shape that remains essentially invariant under small bending, and the three-dimensional
coordinates are denoted, as usual, by~$(x,y,z)$).

\begin{figure}
                \centering
                \includegraphics[width=.55\linewidth]{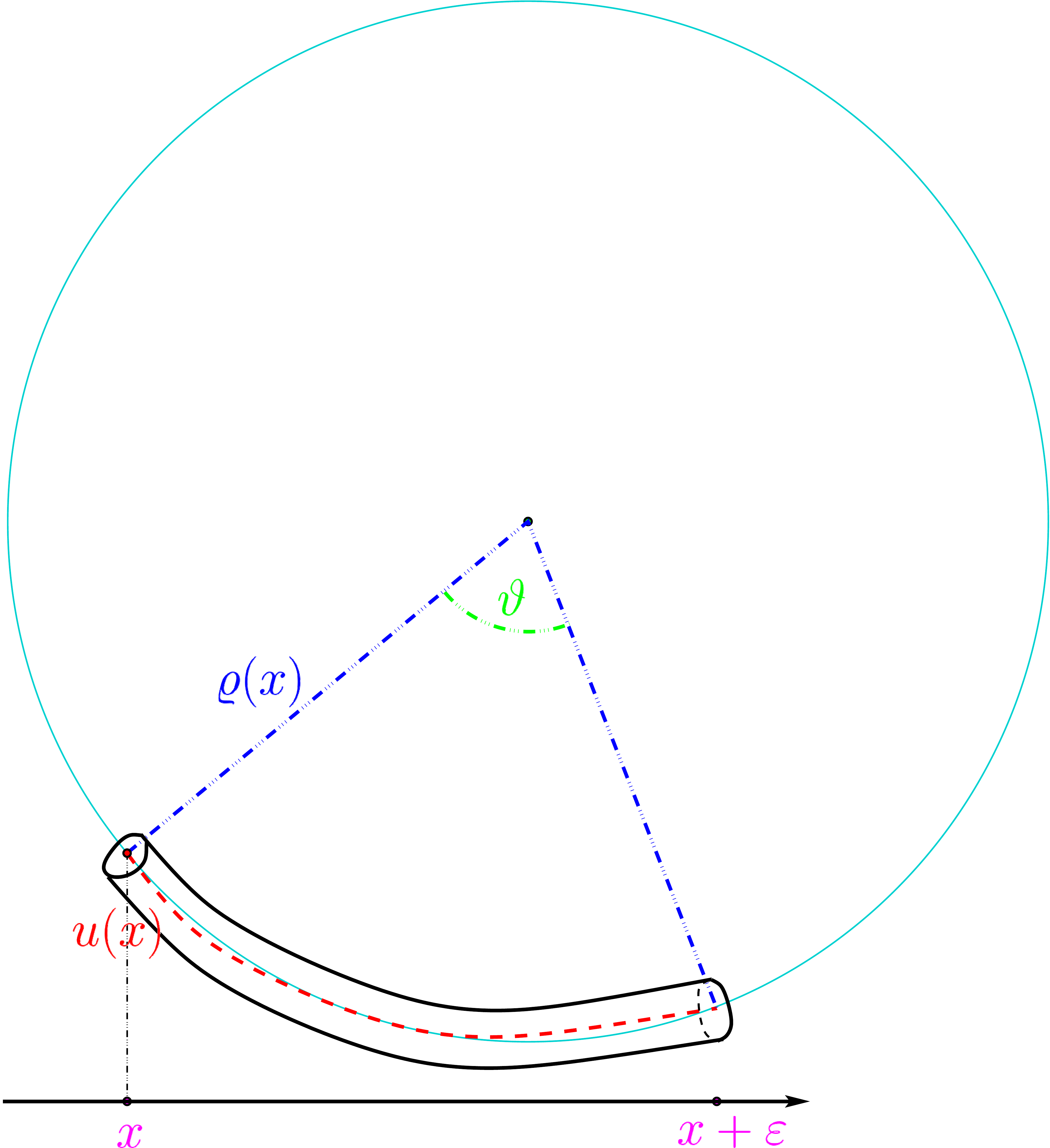}
        \caption{\sl A deflected beam in the $(x,z)$-plane with an osculating circle.}\label{IBEMDENS345UMHDNOJHNFOJED}
\end{figure}

More precisely, the graph~$z=u(x)$ describes the ``neutral fiber'' of the beam,
that is a fiber which maintains infinitesimally the original length that it possessed at rest, see Figure~\ref{IBEMDENS345UMHDNOJHNFOJED}.
The transversal sections of the beam that were perpendicular to the neutral fiber before the beam deforms are assumed to remain perpendicular to the neutral fiber after the bending.
We stress that this neutral fiber may not be located at the mid-height of the beam. In any case, we assume that on one side
of this neutral fiber compression takes place and on the other side tension occurs
(specifically, in Figure~\ref{IBEMDENS345UMHDNOJHNFOJED}, the upper fibers are in compression and the lower fibers are under tension).

To efficiently describe the elastic forces produced by this interplay of compression and tension, it is profitable to consider
the osculating circle at a given point of the neutral fiber, see Figure~\ref{IBEMDENS345UMHDNOJHNFOJED}.
Considering an infinitesimal quantity~$\e>0$ and an infinitesimal element of the beam between~$x$ and~$x+\e$,
we thereby replace the graph of~$u$ with a small portion of the osculating circle of the graph of~$u$ tangent to~$(x,u(x))$.
We denote by~$\varrho(x)$ the radius of this osculating circle and by~$\vartheta$ the corresponding polar angle infinitesimally
joining~$(x,u(x))$ to~$(x+\e,u(x+\e))$.

We consider an upper fiber located at distance~$z$ above the neutral fiber and we
try to quantify the force produced by compression.
To this end, we consider the ratio~$\upsilon(z)$
between the variation of infinitesimal
length of the fiber and its original length before bending
and we suppose that the corresponding compression force is linear with respect to~$\upsilon(z)$
(linearity with respect to length variations is indeed the main ingredient of
elasticity according to Hooke's Law, and note that the bigger the distance~$z$ from the neutral fiber the bigger
the compression). We also observe that the infinitesimal length of this fiber after bending (in the osculating circle approximation)
is given by~$(\varrho(x)-z)\vartheta$ while its original length before bending was equal to that of the neutral fiber,
which is~$\varrho(x)\vartheta$.
Therefore,
$$\upsilon(z)=\frac{\varrho(x)\vartheta-(\varrho(x)-z)\vartheta}{\varrho(x)\vartheta}
=\frac{z\vartheta}{\varrho(x)\vartheta}=\frac{z}{\varrho(x)}$$
and accordingly the compression force related to an upper fiber at distance~$z$
from the neutral fiber is proportional to~$\frac{z}{\varrho(x)}$.

Similarly, the lower fiber located at distance~$z$ below the neutral
one would produce by tension a force which is
proportional to~$\frac{z}{\varrho(x)}$
(notice that the forces produced by tension are opposed to the ones produced by compression).
Thus, if we suppose that the elastic modulus of the material is the same for compression and tension
(which is another classical ansatz in the elasticity theory according to Hooke's Law) and we use the convention that a positive~$z$
corresponds to the upper fibers and a negative~$z$
corresponds to the lower fibers, the resulting elastic force~$F$ is approximately oriented along the tangential axis and equals to~$\frac{Ez}{\varrho(x)}$, where~$E>0$ denotes the elastic modulus of the material (that we suppose to be constant).

We stress that the fact that the force~$F$ changes sign with~$z$ (accordingly, that compression and tension produce forces
in opposite directions above and below the neutral fiber) produces a total torque~$M(x)$.
Indeed, to compute this torque at the point~$(x,u(x))$, for each~$z$ we take into consideration
the vector product between the position vector and the corresponding force, namely
if~$(y,z)$ belongs to the cross section~$A$ of the beam corresponding to~$(x,u(x))$,
the magnitude of the torque at~$(y,z)$ is approximately given by~$zF=\frac{Ez^2}{\varrho(x)}$.

The magnitude of the total torque~$M(x)$ is thereby the corresponding integral for~$(y,z)\in A$, that is
\begin{equation}\label{uJNSEMM0olSdfg2} M(x)=\iint_A \frac{Ez^2}{\varrho(x)}\,dy\,dz.\end{equation}
It is also customary\footnote{Of course, the axial second moment of area depends on the shape of the cross section of the beam.
For instance, for a circular cross section we have that~$A=\{(y,z)\in\R^2$ s.t. $y^2+z^2<r^2\}$, for some~$r>0$, and therefore
in this case
$$I=\iint_{\{ y^2+z^2<r^2\}} z^2\,dy\,dz=2\int_{-r}^{r}z^2\sqrt{r^2-z^2} \,dz=\frac{\pi r^4}{4}.
$$
Instead, if the cross section is a rectangle~$(-a,a)\times(-b,b)$ for some~$a$, $b>0$,
$$I=\iint_{(-a,a)\times(-b,b) } z^2\,dy\,dz=\frac{4ab^3}{3}.
$$}
to define the axial second moment of area \index{axial second moment of area} as
$$ I:=\iint_A z^2\,dy\,dz.$$
With this notation, we can write~\eqref{uJNSEMM0olSdfg2} as
\begin{equation}\label{uJNSEMM0olSdfg223} M(x)=\frac{EI}{\varrho(x)}.\end{equation}
We also recall that the osculating radius~$\varrho(x)$ can be written as~$\frac{(1+(u'(x))^2)^{3/2}}{u''(x)}$,
see e.g.~\cite{MR1807240}.
Hence, for small deformations, the osculating radius~$\varrho(x)$ can be approximated by~$\frac{1}{u''(x)}$
and then, with this approximation, equation~\eqref{uJNSEMM0olSdfg223} reduces to
\begin{equation}\label{uJNSEMM0olSdfg224} M(x)= EI\,u''(x).\end{equation}

\begin{figure}
                \centering
                \includegraphics[width=.65\linewidth]{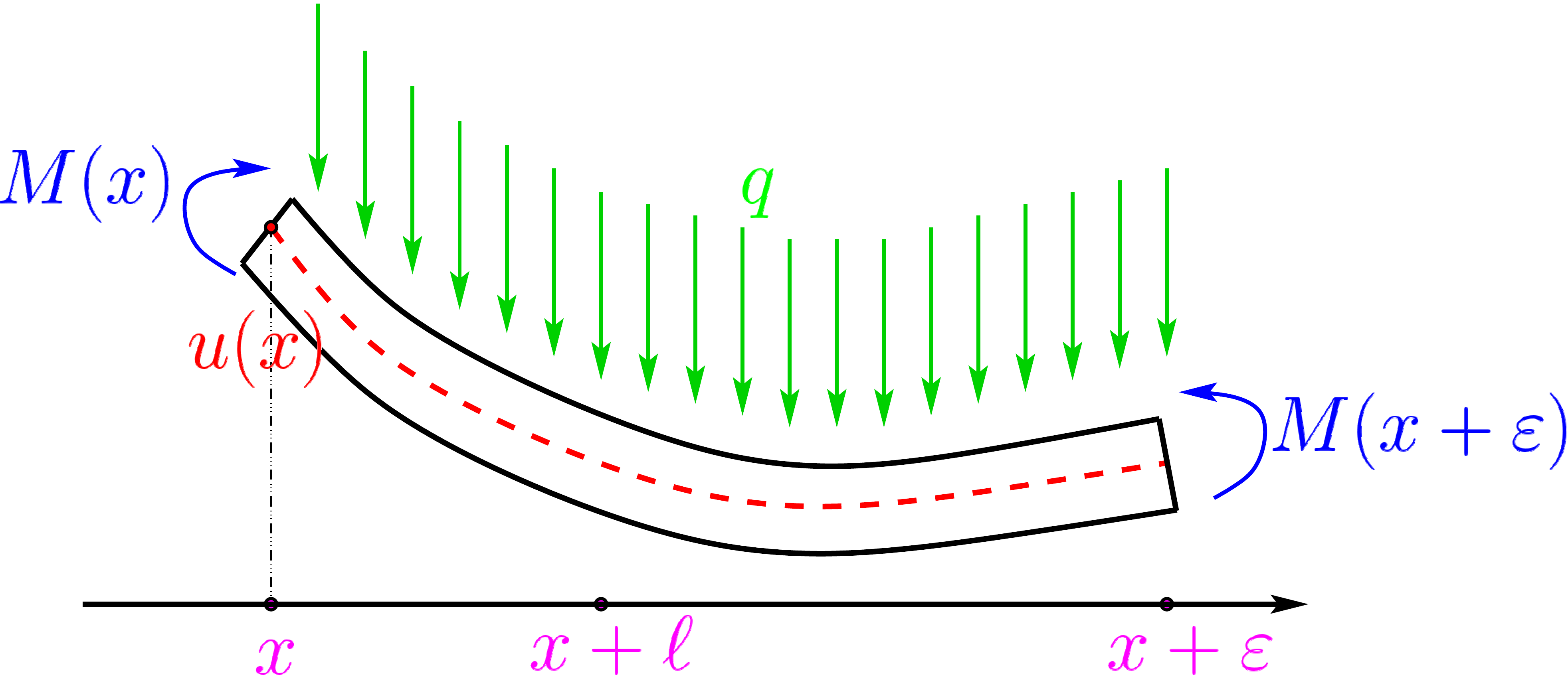}
        \caption{\sl Loads and moments of a beam.}\label{INS345UMHDNBENDIOJHNFOJED}
\end{figure}

To carry on with the derivation of~\eqref{BARREbi96},
we now take into consideration the balance of moments, see Figure~\ref{INS345UMHDNBENDIOJHNFOJED}.
For this, we assume that the beam is subject to a vertical distributed load of magnitude~$q=q(x)$
and we compute its torque at the reference point corresponding to~$x+\e$. For this,
given~$\ell\in(0,\e)$, we calculate the magnitude of the vector product between the distributed load at the point~$x+\ell$,
which is~$(0,-q(x+\ell))$ and the vector joining the reference point to the point corresponding to~$x+\ell$, which is~$\big(x+\ell,u(x+\ell)\big)-\big(x+\e,u(x+\e)\big)=\big( \ell-\e, u(x+\ell)-u(x+\e)\big)$, thus finding the quantity~$ q(x+\ell)(\ell-\e)$.
Correspondingly, the magnitude of the total torque produced by the external load on the beam between~$x$ and~$x+\e$ corresponds to
\begin{equation}\label{BAMO-MAHNdfo-21} \int_0^\e q(x+\ell)(\ell-\e)\,d\ell=
\int_0^\e \big(q(x)+O(\e)\big)(\ell-\e)\,d\ell=
-\frac{\e^2q(x)}2+O(\e^3).\end{equation}
Moreover,
\begin{equation}\label{BAMO-MAHNdfo-22} M(x+\e)-M(x)=\e \partial_xM (x)+\frac{\e^2}2 \partial_x^2 M(x)+O(\e^3).\end{equation}
The balance of moments at the second order in~$\e$ between~\eqref{BAMO-MAHNdfo-21} and~\eqref{BAMO-MAHNdfo-22}
leads\footnote{The first order in~$\e$ in~\eqref{BAMO-MAHNdfo-22} reveals the existence of a shear force
acting on the faces of the beam: namely, the first order in~$\e$ moments balance would produce that this shear force equates
the term~$\partial_xM$ (possibly up to a sign convention).}
to
$$ \partial_x^2 M(x)= q(x).$$
By combining this with~\eqref{uJNSEMM0olSdfg224} we obtain
$$ q(x)=\partial_x^2( EI\,u''(x))=EI\,u''''(x),$$
which corresponds to the beam equation in~\eqref{BARREbi96} (up to normalizing constants).
See e.g.~\cite{MR4298495} for further details on the beam equation and on related topics.\medskip

\begin{figure}
                \centering
                \includegraphics[width=.35\linewidth]{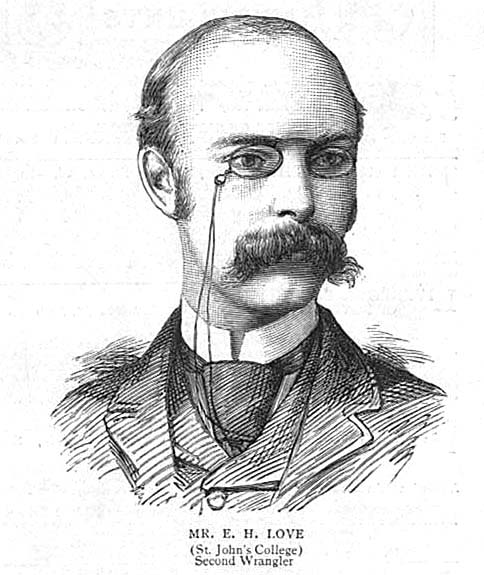}
        \caption{\sl Augustus Edward Hough Love (Public Domain image from
        Wikipedia).}\label{LOVUMHDNOJHNFOJED}
\end{figure}

It is interesting to generalize the model discussed so far to comprise the case of\footnote{The mathematical model commonly used to describe stresses and deformations in thin plates 
is sometimes called ``Kirchhoff-Love theory'',
after Gustav Robert Kirchhoff and Augustus Edward Hough Love.
We have already met Kirchhoff on page~\pageref{SIMKIRAoGREENFIDItangeFI}.
Besides his work in elasticity theory, Love dedicated himself to the study of
the spin-orbit locking, that is the phenomenon occurring when an astronomical body
always has the same face toward the object it is orbiting
(for example, up to some minor variability,
the same side of the Moon always faces the Earth). 
See Figure~\ref{LOVUMHDNOJHNFOJED} for a wood engraved print representing Love taken from the British Newspaper The Graphic.}
bending plates. For instance, one could model a thin plate as a graph~$x_{n+1}=u(x_1,\dots,x_{n})$, with~$x=(x_1,\dots,x_n)\in\R^n$
(of course, bearing in mind that concretely~$n=2$ and the corresponding plate is two-dimensional).
In this setting, the beam equation in~\eqref{BARREbi96} leads to the plate equation \index{plate equation}
\begin{equation}\label{BARREbi96-PLA}
\Delta^2 u=q,
\end{equation}
where~$\Delta^2$ is the Laplace operator applied twice.

To deduce the plate equation in~\eqref{BARREbi96-PLA} from the beam equation in~\eqref{BARREbi96} one can argue as follows.
Let us suppose that the vertical bending force of a horizontal plate is described, for small deformations, by a linear differential operator~${\mathcal{L}}$ of the form
\begin{equation}\label{BARREbi96-PLA-DAD2} {\mathcal{L}}:=\sum_{{\alpha\in\N^n}\atop{|\alpha|\le m}}c_\alpha\partial^\alpha,\end{equation}
for some~$m\in\N$ and~$c_\alpha\in\R$. Under this assumption, the force balance
would lead to the equation~${\mathcal{L}}u=q$, being~$q$ the external load.

Then, to check~\eqref{BARREbi96-PLA}, one needs to show that
\begin{equation}\label{BARREbi96-PLA-DAD}
{\mathcal{L}}=\Delta^2 .
\end{equation}
For this, we assume that, consistently with the model of beam deformation in~\eqref{BARREbi96},
the operator~${\mathcal{L}}$ acts as a fourth derivative on any one-dimensional function.
Namely, let us suppose that if~$u_\omega(x):=u_0(\omega\cdot x)$ for some~$u_0:\R\to\R$ and~$\omega\in\partial B_1$
then
\begin{equation}\label{BARREbi96-PLA-DAD3} {\mathcal{L}}u_\omega(x)=u_0''''(\omega\cdot x).\end{equation}
We observe that, for all~$j\in\{1,\dots,n\}$,
$$ \partial_j u_\omega(x)=\omega_j \,u_0'(\omega\cdot x)$$
and therefore, if~$\alpha=(\alpha_1,\dots,\alpha_n)\in\N^n$,
$$ \partial^\alpha u_\omega(x)=\partial_1^{\alpha_1}\dots\partial_n^{\alpha_n} u_\omega(x)
=\omega_1^{\alpha_1}\dots\omega_n^{\alpha_n} \,u_0^{(|\alpha|)}(\omega\cdot x)=\omega^\alpha\,u_0^{(|\alpha|)}(\omega\cdot x).
$$
Hence, by~\eqref{BARREbi96-PLA-DAD2} and~\eqref{BARREbi96-PLA-DAD3},
$$ u_0''''(\omega\cdot x)={\mathcal{L}}u_\omega(x)=\sum_{{\alpha\in\N^n}\atop{|\alpha|\le m}}c_\alpha\, \omega^{\alpha} \,u_0^{(|\alpha|)}(\omega\cdot x)
.$$
Since~$u_0$ is arbitrary, we infer that~$c_\alpha=0$ unless~$|\alpha|=4$ and
$$ 1=\sum_{{\alpha\in\N^n}\atop{|\alpha|=4}}c_\alpha\, \omega^{\alpha}
.$$
Now, given~$X\in\R^n\setminus\{0\}$ we let~$\rho:=|X|$ and~$\omega:=\frac{X}{|X|}$ and we deduce that
\begin{eqnarray*}&& \sum_{i,j=1}^n X_i^2 X_j^2=\left(\sum_{i=1}^n X_i^2\right)^2
=\big(|X|^2\big)^2=|X|^4=|X|^4\sum_{{\alpha\in\N^n}\atop{|\alpha|=4}}c_\alpha\, \omega^{\alpha}\\&&\qquad=
\sum_{{\alpha\in\N^n}\atop{|\alpha|=4}}c_\alpha\, \rho^{|\alpha|}\omega^{\alpha}=\sum_{{\alpha\in\N^n}\atop{|\alpha|=4}}c_\alpha\,\big( \rho \omega\big)^{\alpha}=\sum_{{\alpha\in\N^n}\atop{|\alpha|=4}}c_\alpha\,X^{\alpha}
.\end{eqnarray*}
Consequently, by the Identity Principle for polynomials,
$$ c_\alpha=\begin{dcases}
1 & {\mbox{ if $\alpha=2e_i+2e_j$ for some }}i,j\in\{1,\dots,n\},\\
0&{\mbox{ otherwise.}}
\end{dcases}$$
{F}rom this and~\eqref{BARREbi96-PLA-DAD2} we arrive at
$$ {\mathcal{L}}=\sum_{i,j=1}^n \partial^2_i\partial_j^2=\sum_{i=1}^n \partial^2_i\Delta=\Delta^2$$
and this establishes~\eqref{BARREbi96-PLA-DAD}, as desired.

\subsection{Gravitation and electrostatics}\label{SECT:Gravitation and electrostatics}

One of the reasons for which the Laplace operator became very popular, especially in the eighteenth and\index{gravitation} 
nineteenth centuries, is its pivotal role in the description of the gravity and electromagnetic\index{electrostatics} potentials.
The physical intuition of these phenomena was actually very helpful for the beautiful minds of those ages
in understanding a number of very deep concepts which paved the way to the modern theory of
partial differential equations and which are now widely used in many branches of mathematics, physics
and engineering.

We give here a couple of motivations\footnote{Also, a different but related
perspective on electrostatics and magnetism will be given in the forthcoming
Section~\ref{MAXW}.} relating the gravity field and harmonic functions
(similar arguments would link the electrostatic field and harmonic functions as well,
just changing the notion of mass with that of electric charge, and possibly allowing
a sign change, given the fact that unlike charges attract each other, but
like charges repel).\medskip

One observation is that the gravity force (in~$\R^3$) is inversely proportional to the square of the distance,
namely a point mass located at~$x$ is attracted by a Permanent Center of Gravity located at the origin
via a force~$f(x)$ equal (up to dimensional constant) to~$-\frac{x}{|x|^3}$, with the minus
sign to stress that this force is attractive, hence tends to reduce the distance of
the mass located in~$x$ to the origin. 

It is readily seen that~$f$ is originated by a gravity potential~$u$: more precisely, if we set~$u(x):=\frac{1}{|x|}$,
we have that~$\nabla u(x)=-\frac{x}{|x|^3}=f(x)$. Remarkably, the gravity potential is harmonic away\footnote{We will learn in equation~\eqref{WRBSY0perIKSMaismoBBeJ}
a more handy way to perform the calculation in~\eqref{WRBSY0perIKSMaismoBBeJPREB}.}
from
the origin (here it is important to work in dimension~$3$) since, if~$x=(x_1,x_2,x_3)\in\R^3\setminus\{0\}$,
\begin{equation}\label{WRBSY0perIKSMaismoBBeJPREB}
\Delta u(x)=-\sum_{i=1}^3 \partial_i \frac{x_i}{|x|^3}
=-\sum_{i=1}^3 \left(\frac{1}{|x|^3}
-\frac{3 x_i^2}{|x|^5}\right)=-\frac{3}{|x|^3}
+\frac{3 |x|^2}{|x|^5}=0.
\end{equation}
\medskip

These considerations can be recast in a possibly more general, and more ``geometric'',
framework, also allowing a higher dimensional presentation.
For this, we first employ the fact that the
work done by the gravitational force depends only on initial and final positions, and not on the path between them: this gives that the gravitational vector field~$f$ is conservative and thus can be written as the gradient of some function~$u$, which is usually called the gravitational potential.

Furthermore, Gau{\ss}'s Flux Law states that the flux of the gravity
field through an arbitrary closed surface enclosing no masses (or no electric) charges is necessarily zero (roughly speaking,
the field lines going into the region enclosed by the surface balance exactly the ones coming out).
Thus, from the Divergence Theorem, for every bounded region~$\Omega\subset\R^n$ (say, with smooth boundary)
$$ 0=\int_{\partial\Omega} f(x)\cdot\nu(x)\,d{\mathcal{H}}_x^{n-1}=
\int_{\Omega} \div f(x)\,dx=\int_{\Omega} \div (\nabla u(x))\,dx=\int_{\Omega} \Delta u(x)\,dx.$$
Since~$\Omega$ is arbitrary we thus conclude that the gravity potential (as well as the electrostatic potential) is harmonic away from the masses that are present in the environment.

The study of functions that are harmonic away from a point singularity will be
the main topic of the forthcoming Section~\ref{lfundsP-S}. There we will also appreciate
how different dimensions affect the precise expression of the above potentials.

\subsection{Classical electromagnetism}\label{MAXW}
The theory of classical electromagnetism\index{electromagnetism} has its roots in the so-called Maxwell's equations\index{Maxwell's equations}
that describe the behavior and mutual interaction of electric and magnetic fields.
These equations are somewhat a unified version of several specific information arising from
Gau{\ss}'s Flux Law, Faraday's Law and
Amp\`ere's Law.
In the vacuum (hence, in a region of~$\R^3$ with no electric charges and no electric currents)
Maxwell's equations take the form
\begin{equation}\label{KMSJUNDNSMDEJDN} \div E =0,\qquad \curl E=-{\frac{\partial B }{\partial t}},\qquad
\div B=0
\qquad{\mbox{and}}\qquad\curl B=\frac1{c^2}{\frac{\partial E }{\partial t}}.\end{equation}
Here above, $E$ is the electric field,
$B$ is the magnetic field and~$c$ is a physical constant (corresponding to the speed of light in vacuum).

We also observe that, for all smooth vector fields~$v:\R^3\to\R^3$,
\begin{equation}\label{KMSJUNDNSMDEJDN2}
\nabla (\div v)-\curl(\curl v)=\Delta v.
\end{equation}
To check this classical identity in an elementary way, given~$v=(v_1,v_2,v_3)$,
we write~$v=v^{(1)}+v^{(2)}+v^{(3)}$,
where~$v^{(1)}:=(v_1,0,0)$, $v^{(2)}:=(0,v_2,0)$ and~$v^{(3)}:=(0,0,v_3)$,
and we remark that, thanks to the linear structure of~\eqref{KMSJUNDNSMDEJDN2},
it suffices to check it for the vector fields~$v^{(1)}$, $v^{(2)}$ and~$v^{(3)}$.
We can focus our attention on~$v^{(1)}$, up to reordering coordinates.
In this situation, since, for every smooth vector field~$V:\R^3\to\R^3$,
$$ \curl V=\left(
{\frac{\partial V_3}{\partial x_2}}-{\frac{\partial V_2}{\partial x_3}} ,\;
{\frac{\partial V_1}{\partial x_3}}-{\frac{\partial V_3}{\partial x_1}},\;
{\frac{\partial V_2}{\partial x_1}}-{\frac{\partial V_1}{\partial x_2}}\right),$$
we have that
$$ W:=\curl v^{(1)}=\big(0 ,\partial_3 v_1, -\partial_2 v_1\big)$$
and subsequently
\begin{eqnarray*}&&
\nabla (\div v^{(1)})-\curl(\curl v^{(1)})\\&=&
\nabla(\partial_1 v_1)
-\curl W\\
&=&\left(\partial_{11} v_1,\,
\partial_{12} v_1,\,
\partial_{13} v_1
\right)
-
\left(
{\frac{\partial W_3}{\partial x_2}}-{\frac{\partial W_2}{\partial x_3}} ,\;-
{\frac{\partial W_3}{\partial x_1}},\;
{\frac{\partial W_2}{\partial x_1}}\right)\\
&=&\left(\partial_{11} v_1,\,
\partial_{12} v_1,\,
\partial_{13} v_1
\right)-
\left(
-\partial_{22} v_1-\partial_{33} v_1 ,\,
\partial_{12} v_1,\,\partial_{13} v_1\right)\\&=&
\left(\partial_{11} v_1+
\partial_{22} v_1+\partial_{33} v_1 ,\,0,\,0\right)\\&=&\Delta v^{(1)},
\end{eqnarray*}
thus establishing~\eqref{KMSJUNDNSMDEJDN2}.

Now, in the light of~\eqref{KMSJUNDNSMDEJDN} and~\eqref{KMSJUNDNSMDEJDN2},
\begin{equation}\label{EFEDLFIemfdL}
\partial_{tt}E=c^2\partial_t(\curl B)=c^2\curl(\partial_t B)=-c^2\curl(\curl E)
=c^2\big(\Delta E-\nabla (\div E)\big)=c^2\Delta E
\end{equation}
and similarly
\begin{equation*}
\partial_{tt}B=-\partial_t(\curl E)=-\curl(\partial_t E)=-c^2\curl(\curl B)=
c^2\big(\Delta B-\nabla (\div B)\big)=c^2\Delta B.
\end{equation*}
These observations show that
both the electric and the magnetic fields propagate as waves do, since they satisfy
the wave equation in~\eqref{WAYEBJJD121t416JH098327uyrhdhc832nbcM}.

\subsection{Quantum-mechanical systems}

A popular equation governing the wave function of a quantum-mechanical system
was proposed by  Erwin Schr\"odinger
(this equation formed the basis for the work that resulted in 
Schr\"odinger's 1933 Nobel Prize).

In a nutshell, the Schr\"odinger equation reads
\begin{equation}\label{SCHREQOR}
-i\hbar\partial_t\psi=\frac{\hbar^2 \Delta \psi}{2m}-V\psi,
\end{equation}
where~$i=\sqrt{-1}$, $\psi=\psi(x,t)$ is the wavefunction\footnote{Roughly speaking, a wavefunction is
a function that assigns a complex number to each point~$x\in\R^3$ at each time~$t\in\R$.
The magnitude squared of this function
represents the probability density of measuring the particle as being at the point~$x$ at a given time~$t$.}
of a particle of mass~$m$ subject to a potential~$V$
and~$\hbar$ is a suitable positive physical constant (called the ``reduced Planck constant''\index{reduced Planck constant}).

Particularly interesting solutions of~\eqref{SCHREQOR}
are the so-called standing waves, namely solutions which oscillate in time but whose amplitude profile does not change in space.
{F}rom the mathematical point of view, these solutions are of the form
$$ \psi(x,t)=u(x)\,e^{i\phi t},$$
where~$u$ is a real valued function and~$\phi>0$ is the time frequency of oscillation.

Since for standing waves we have that
$$ \partial_t\psi(x,t)=i\phi\,u(x)\,e^{i\phi t}
\qquad{\mbox{and}}\qquad
\Delta \psi(x,t)=\Delta u(x)\,e^{i\phi t},$$
these specia5l solutions of~\eqref{SCHREQOR} are actually solutions\footnote{According to the classification
presented in footnote~\ref{CLASSIFICATIONFOOTN}
on page~\pageref{CLASSIFICATIONFOOTN}, equation~\eqref{CLASSIFICATIONFOOTNSV}
is of elliptic type. When~$V=V(x)$ equation~\eqref{CLASSIFICATIONFOOTNSV}
is linear in~$u$. More complex situations arise when~$V$ also depends on the magnitude of the wavefunction,
say~$V=V(x,u(x))$ or even~$V=V(u(x))$, since in this situation equation~\eqref{CLASSIFICATIONFOOTNSV}
is not anymore linear in~$u$ (though it is linear with respect to the Hessian of~$u$). A study of these ``semilinear'' equations \index{semilinear equation}
will be started
on page~\pageref{Pohozaev Identity}. \label{BIS Pohozaev Identity}
Further comments on semilinear equations will be given in footnote~\ref{TRIS Pohozaev Identity}
on page~\pageref{TRIS Pohozaev Identity}.}
of
\begin{equation}\label{CLASSIFICATIONFOOTNSV}
\hbar \phi u 
=\frac{\hbar^2 \Delta u }{2m}-V u.
\end{equation}

Going back to the initial discussion of this section,
one may wonder why and how Schr\"odinger introduced the equation in~\eqref{SCHREQOR} above.
This would be a rather long and complicated story and a short answer
was given in~\cite[Chapter~16]{MR0213079}:
``Where did we get that from? Nowhere. It's not possible to derive it from anything you know. It came out of the mind of Schr\"odinger, invented in his struggle to find an understanding of the experimental observations of the real world''.
However, several quite convincing motivations of the Schr\"odinger equation\index{Schr\"odinger equation} are possible (see e.g.~\cite{2006physics2, 2006physics, MR3393677, YAN2021} and the references therein)
and we present one of them here below, based on the notion of canonical quantization\index{canonical quantization}.

To this end, we recall that the energy~${\mathcal{E}}$ of a photon is proportional to its temporal frequency~$\phi$, according to the so-called Planck-Einstein energy-frequency relation
$$ {\mathcal{E}}=\hbar\phi.$$
Also, Einstein's relativistic energy formula reads
\begin{equation}\label{ENEJHND8iujMNSID-77} {\mathcal{E}}=\sqrt{ (|p|\, c)^2 + (mc^2)^2},\end{equation}
being~$p$ the relativistic momentum, $m$ the mass at rest and~$c$ the speed of light. For massless particles, such as photons, this reduces to \begin{equation}\label{ENEJHND8iujMNSID} {\mathcal{E}}= |p|\,c,\end{equation} from which we obtain
\begin{equation}\label{PPho} |p|=\frac{\mathcal{E}}c=\frac{\hbar\phi}{c}.\end{equation}

Now, we relate a photon to its electric field~$E$ traveling through space and manifesting itself as a solution of the wave equation (recall~\eqref{EFEDLFIemfdL}). That is, we consider a direction of propagation~$\varpi\in\partial B_1$, a spatial frequency~$\kappa$ and a temporal frequency~$\phi$ and we suppose that,
for large distances from the photon's source, the field is modeled by a simple plane wave, say
$$ E=E_0 e^{i(\kappa\varpi\cdot x-\phi t)},$$ where~$E_0\in(0,+\infty)$.
More precisely, by~\eqref{EFEDLFIemfdL},
$$ -\phi^2 E_0 e^{i(\kappa\varpi\cdot x-\phi t)}=
\partial_{tt} E=c^2\Delta E=-c^2 \kappa^2 E_0 e^{i(\kappa\varpi\cdot x-\phi t)}$$
therefore
\begin{equation}\label{PPho22536-0i2rkjfmMS}
\phi=c \kappa,\end{equation} thus relating temporal and spatial frequencies of the photon.

We thus obtain from~\eqref{PPho}
that
\begin{equation}\label{PPho22536}
|p|=\hbar\kappa\end{equation} and therefore
\begin{equation}\label{PPho2}|p|E=\hbar\kappa E=-i\hbar\nabla E\cdot\varpi.\end{equation}
Taking the direction of the momentum to coincide with the spatial direction~$\varpi$ of the wave, i.e. taking~$p=|p|\varpi$,
we rewrite~\eqref{PPho2} in the form
\begin{equation}\label{PPho​3}pE=-i\hbar\nabla E.\end{equation}

We can also reconsider~\eqref{ENEJHND8iujMNSID} in view of~\eqref{PPho22536-0i2rkjfmMS} and~\eqref{PPho22536} and write that
$${\mathcal{E}}= |p|\,c=\hbar c\kappa=\hbar\phi$$ and consequently
\begin{equation}\label{ENEJHND8iujMNSID-2krwjeg}
{\mathcal{E}}E=\hbar \phi E=i\hbar\partial_t E.
\end{equation}

With this, without the aim of being exhaustive, we can recall some of the main ideas leading to canonical quantization.
This procedure, originally introduced by Paul Dirac in his 1926 PhD thesis, wishes to recast a classical theory into a quantum
one by preserving (as much as possible) its symmetries
and formal structures. Since the Hamiltonian formalism is one of the key tools to understand symmetries of classical mechanical systems, a natural idea in this framework is to try to (at least partially) preserve the Hamiltonian structure in the quantum description of nature.

To this end, a simple observation is that the Hamiltonian formalism relies on conjugated variables~$q$ and~$p$, with the physical meaning of position and momentum. The first goal of the canonical quantization is therefore to rephrase position and momentum in a way that is compatible with the quantum description of the observables.

In particular, in quantum mechanics,
all significant features
of a particle are contained in a certain state~$\psi$. 
The observables are represented by operators acting on states. For instance, the eigenvalues of an operator represent the values
of the measurements of the corresponding fundamental states of a particle (namely, the eigenfunctions, or eigenstates);
since any state is represented as a linear combination of eigenstates, the application of an operator to the state~$\psi$
corresponds to the determination of a measurable parameter and the physical act of measuring corresponds to make the
values of the state ``collapse'' from a superposition of eigenstates to a single eigenstate due to the
interaction with the external world.

With these basic principia in mind, we therefore aim at replacing the classical position and momentum variables~$q$ and~$p$ with two operators, say~$Q$ and~$P$, by preserving the original structure as much as possible.

The most natural setting is therefore to consider the position operator~$Q$ as the operator that corresponds to the position evaluation of a state (simply, applying~$Q$ to a state~$\psi$ being the evaluation of the function~$\psi$ at a given point in space).

As for the corresponding momentum operator, in light of~\eqref{PPho​3}, a natural choice is to take~$P$ as the differential operator~$-i \hbar \nabla$, thus obtaining a quantum analogue of~\eqref{PPho​3} via the relation
\begin{equation}\label{ENEJHND8iujMNSID-2krwjeg-21-ktg22} P\psi=-i\hbar\nabla\psi.\end{equation}

Pushing this analogy a bit further, we can also obtain a quantum counterpart~${\mathcal{H}}$ of
the total energy~${\mathcal{E}}$ (say, the Hamiltonian) in~\eqref{ENEJHND8iujMNSID-2krwjeg}.
In this setting, the operator analogue of~\eqref{ENEJHND8iujMNSID-2krwjeg} would then be
\begin{equation}\label{ENEJHND8iujMNSID-2krwjeg-21-ktg}
{\mathcal{H}}\psi=-i\hbar\partial_t \psi.
\end{equation}

While this canonical formalism for quantum mechanics was obtained by analyzing the special case of the wave produced by a photon,
we can believe that the same protocol governs the evolution of the wave function of a particle (not necessarily a photon)
subject to a given potential~$V$. For instance, for a particle
of mass~$m>0$ one can push the previous construction by considering the quantum analogue of the mechanical energy
(perhaps neglecting for the moment some relativistic effects,
which we will briefly discuss in the forthcoming footnote~\ref{BEJMNDSIMinBArhgCLAinrT}).
Indeed, in classical mechanics the total energy would be given by the sum of the kinetic energy
of a particle (equal to~$\frac{m|v|^2}{2}$, being~$v$ its velocity) and the potential energy~$V$. That is,
recalling the classical momentum definition~$p=mv$,
$$ {\mathcal{E}}=\frac{m|v|^2}{2}+V=\frac{|p|^2}{2m}+V.$$
By formally applying the quantization in~\eqref{ENEJHND8iujMNSID-2krwjeg-21-ktg22} and~\eqref{ENEJHND8iujMNSID-2krwjeg-21-ktg}, we thus obtain its quantum counterpart as
$$ i\hbar\partial_t={\mathcal{H}}=\frac{|P|^2}{2m}+V=
\frac{P\cdot P}{2m}+V=-
\frac{\hbar^2 \nabla\cdot\nabla}{2m}+V=-\frac{\hbar^2 \Delta}{2m}+V,
$$
that is the Schr\"odinger equation\footnote{It is interesting to point out that variations
in the previous arguments lead to other equations of interest: for instance, 
thinking over the full relativistic energy formula in~\eqref{ENEJHND8iujMNSID-77} and formally inserting the quantum mechanical operator in~\eqref{PPho​3}, the quantum analogue of the kinetic energy becomes
$$ \sqrt{ -\hbar^2\, c^2\Delta + (mc^2)^2},$$
that is the ``square root of the Laplacian''\index{square root of the Laplacian} (see e.g.~\cite{MR3967804}
for a basic introduction to this very interesting object). Thus,
in the presence of an external potential~$V$, the conservation of the full energy and the quantization
in~\eqref{ENEJHND8iujMNSID-2krwjeg-21-ktg} lead to the balance
$$-i\hbar\partial_t =\sqrt{ -\hbar^2\, c^2\Delta + (mc^2)^2}-V,$$
which produces the equation
$$ -i\hbar\partial_t \psi=\sqrt{ -\hbar^2\, c^2\Delta \psi+ (mc^2)^2}-V\psi.$$
This equation and variants of it are usually \label{BEJMNDSIMinBArhgCLAinrT}
called ``relativistic Schr\"odinger equations''\index{relativistic Schr\"odinger equation}, see e.g. equation~(1.4) in~\cite{MR1054115}.

Another approach for taking into account relativistic effects consists in taking the square of the
relativistic energy formula in~\eqref{ENEJHND8iujMNSID-77}, thus writing that
$$ {\mathcal{E}}^2= (|p|\, c)^2 + (mc^2)^2.$$
In the absence of external potentials (i.e., if this represents the square of the total energy of the system), one can think about exploiting the quantization procedure in~\eqref{PPho​3}
and~\eqref{ENEJHND8iujMNSID-2krwjeg-21-ktg}, formally finding that
$$ -\hbar^2\partial_{tt}= -\hbar^2 c^2\Delta+ (mc^2)^2.$$
This leads to the equation
$$ \partial_{tt}\psi=  c^2\Delta\psi-\frac{ m^2c^4\psi}{\hbar^2},$$
which is often referred to with the name of ``Klein-Gordon equation''\index{Klein-Gordon equation},
see e.g. equation~(34) in~\cite{MR734851}.} in~\eqref{SCHREQOR}.

\subsection{Bouncing balls and whispers}
\begin{figure}
                \centering
                \includegraphics[width=.65\linewidth]{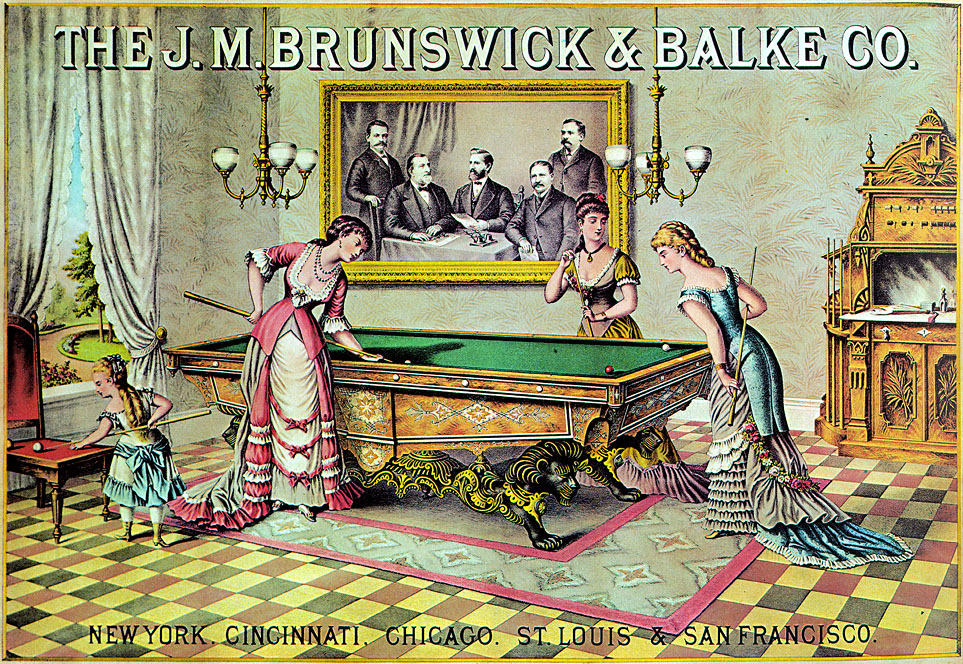}
        \caption{\sl Advertising poster, early 1880s (Public Domain image from
        Wikipedia).}\label{HAFOUMBILLI923NOJHNFOJEDdw4I}
\end{figure}

\begin{figure}
                \centering
                \includegraphics[width=.2\linewidth]{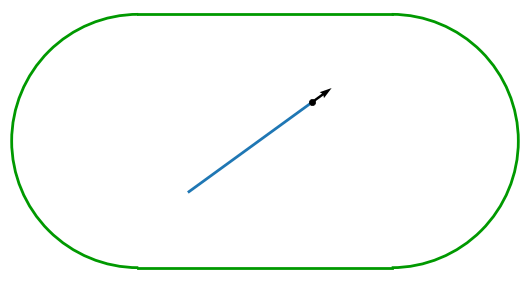}$\quad$
                \includegraphics[width=.2\linewidth]{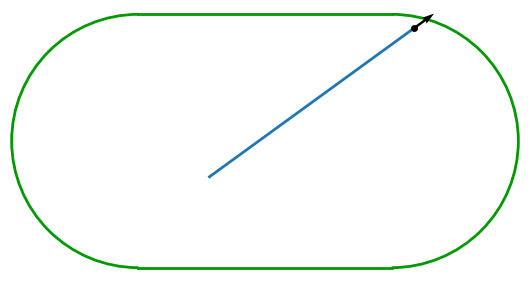}$\quad$
                \includegraphics[width=.2\linewidth]{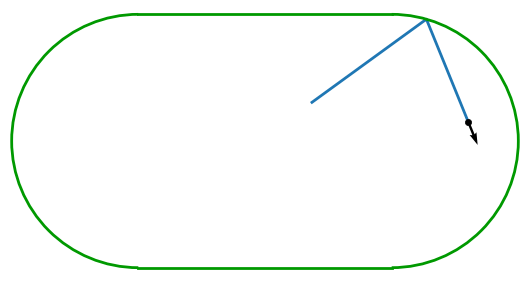}$\quad$
                \includegraphics[width=.2\linewidth]{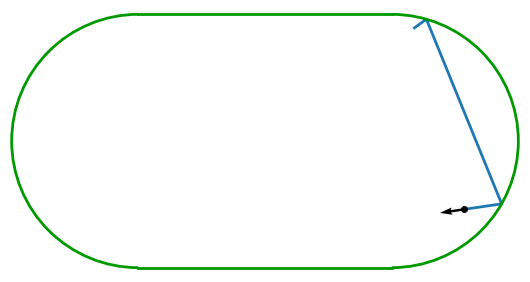}\\
                \includegraphics[width=.2\linewidth]{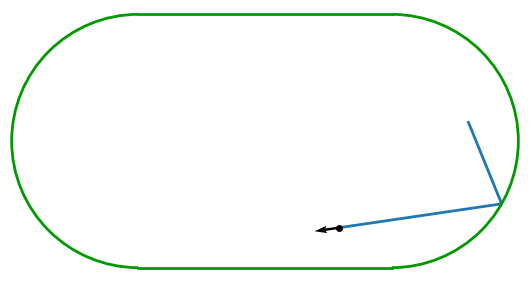}$\quad$
                \includegraphics[width=.2\linewidth]{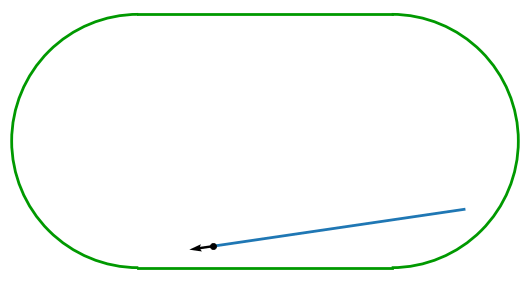}$\quad$
                \includegraphics[width=.2\linewidth]{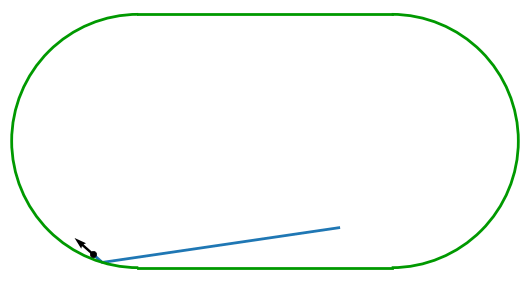}$\quad$
                \includegraphics[width=.2\linewidth]{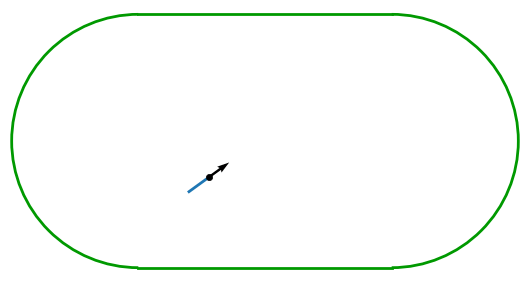}
        \caption{\sl A dynamical billiard (animation by George Datseris, image from
        Wikipedia,
        licensed under the Creative Commons Attribution-Share Alike 4.0 International license).}\label{HSTAGde0ioJHNFOJEDdw4I}
\end{figure}

The popular game of billiard, see Figure~\ref{HAFOUMBILLI923NOJHNFOJEDdw4I},
has a mathematical counterpart in
the study of dynamical billiards, namely a dynamical system of a material particle confined in a (typically, bounded and sufficiently smooth) region~$\Omega$ of~$\R^n$ (or of a more general manifold)
which alternates between the free motion in the interior of~$\Omega$
(given by a straight line travelled at unit speed, or, for general manifolds
by a geodesic path) and specular reflections from the boundary.
See Figure~\ref{HSTAGde0ioJHNFOJEDdw4I} for a visual animation\footnote{Actually,
Figure~\ref{HSTAGde0ioJHNFOJEDdw4I} depicts the motion of a particle
in a famous example of billiard introduced in~\cite{MR530154}.} of
a dynamic billiard.

More precisely,
the billiard flow in the region~$\Omega$ corresponds to the motion of a material point which is moving in~$\Omega$ with constant (e.g., unit)
velocity in the interior of~$\Omega$ making reflections at the boundary of~$\Omega$
according to the usual law of geometric optics, namely prescribing that
the angle of incidence is equal to
the angle of reflection.\medskip

It is interesting to notice that this motion can be obtained, at least in the limit and at least
exploiting some formal expansions and approximations,
from a smooth Hamiltonian system. Namely, one can consider a large parameter~$M>0$ and a potential~$V_M:\R^n\to\R$ defined as
$$ V_M(q):=\begin{dcases}
0 & {\mbox{ if }}q\in\Omega,\\
\displaystyle
\frac{M}2\big({\rm dist}(q,\partial\Omega)\big)^2& {\mbox{ if }}q\in\R^n\setminus\Omega.
\end{dcases}$$

Then, at least at a formal level,
\begin{equation}\label{FOHAMM}\begin{split}&
{\mbox{the Hamiltonian }}H_M(q,p):=\frac{|p|^2}2+V_M(q)\\&
{\mbox{formally recovers the billiard motion as }}M\to+\infty.\end{split}
\end{equation}
To check this, we observe that the equations of motion induced by~$H_M$ are
$$ \begin{dcases}
\dot q=p,\\
\dot p=-\nabla V_M(q).
\end{dcases}$$
Consequently, as long as the trajectory lies in~$\Omega$ we have that~$\dot p=0$
and thus~$q$ evolves linearly (precisely as in the billiard motion, and this property
is stable as~$M\to+\infty$).

We thereby investigate what happens when a particle hits the boundary of~$\Omega$
and we will show that the motion induced by~$H_M$ in this case reduces to the billiard reflection law when~$M\to+\infty$. Let us assume that the particle hits the boundary of~$\Omega$ at
some time~$t_0$, which we can assume to be~$0$,
some point~$p_0\in\partial\Omega$, which, up to a translation can be taken to be the origin.
Also, up to a rotation, we can assume that \begin{equation}\label{012ofek029rfeiwohg3uiwfegv-ytrgfFe7ywfhenNARTFGVIvr67002}{\mbox{the exterior unit
normal of~$\Omega$ at the origin equals~$e_n$.}}\end{equation}
We also suppose, for simplicity, that the velocity of the particle when reaching~$\Omega$ is not tangent to~$\partial\Omega$, hence~$p_n(0)>0$. In particular, if~$p(0)=0$,
we can assume that~$p(t)\in\R^n\setminus\Omega$ for all~$t\in(0,t_1)$, for some~$t_1>0$.

Actually, by conservation of energy, we have that~$|{\rm dist}(q(t),\partial\Omega)|\le\frac{ C}{\sqrt{M}}$
(that is, the particle cannot leave a small neighborhood of~$\Omega$)
and therefore
$$ \nabla V_M(q(t))=M{\rm dist}(q(t),\partial\Omega)\nabla {\rm dist}(q(t),\partial\Omega)=
M{\rm dist}(q(t),\partial\Omega)\nu(\pi(q(t)),$$
being~$\pi$ the projection along~$\partial\Omega$.

Additionally,
$$ \nu(\pi(q(t))=\nu(\pi(0))+O(|\pi(q(t))-\pi(0)|)=
\nu(0)+O(|q(t)|)=e_n+O\left(\frac1{\sqrt{M}}\right),$$
due to~\eqref{012ofek029rfeiwohg3uiwfegv-ytrgfFe7ywfhenNARTFGVIvr67002},
hence for the formal justification of~\eqref{FOHAMM}
we disregard this higher order correction and write simply
$$ \nabla V_M(q(t))\simeq M{\rm dist}(q(t),\partial\Omega)e_n.$$

Again by energy conservation, we know that~$|\dot q(t)|=|p(t)|\le C$, therefore
we can suppose that~$q(t)$ is in a small neighborhood of the origin as long as~$t_1$ is small enough:
hence, using~\eqref{012ofek029rfeiwohg3uiwfegv-ytrgfFe7ywfhenNARTFGVIvr67002},
at a formal level we can use the approximation
$$ {\rm dist}(q(t),\partial\Omega)\simeq q_n(t),$$
where we use the notation~$q=(q',q_n)\in\R^{n-1}\times\R$.

As a result, for all~$t\in(0,t_1)$, if~$t_1$ is small enough,
\begin{equation*}\begin{split}& \ddot q(t)=\dot p(t)=-\nabla V_M(q(t))\simeq -M{\rm dist}(q(t),\partial\Omega)e_n\simeq-
M q_n(t)e_n.
\end{split}\end{equation*}
Therefore,
$$\begin{dcases}
\ddot q'(t)\simeq0,\\ \ddot q_n(t)\simeq-Mq_n(t)
\end{dcases}$$
and we observe that
\begin{equation*}
q(0)=0\qquad{\mbox{and}}\qquad \dot q(0)=p(0).
\end{equation*}
In this way, solving the above ordinary differential equations, we obtain that formally
\begin{equation}\label{lkjb89io-oijohicsbFedf3vdJcbsdR4Hs2L5O3ks4d2ss23oer-1t} q'(t)\simeq p'(0)t\qquad{\mbox{and}}\qquad
q_n(t)\simeq \frac{p_n(0)}{\sqrt{M}}\sin\big(\sqrt{M}t\big).\end{equation}

Now we detect the ``reentering'' of the particle in the domain~$\Omega$.
To this end,
assuming~$M$ appropriately large, we claim that
\begin{equation}\label{iujhbdNArryujf0RR92i3rjtgCahsdn9B7yuhdfGTGJA893-1}
{\mbox{$q(t)\in\R^n\setminus\Omega$
for all~$t\in\left(0,\frac{\pi}{\sqrt{M}}-\frac{1}{\sqrt{M}\,\ln M}\right]$}}\end{equation}
and that
\begin{equation}\label{iujhbdNArryujf0RR92i3rjtgCahsdn9B7yuhdfGTGJA893-2}
{\mbox{there exists }}t_M\in \left(\frac{\pi}{\sqrt{M}}-\frac{1}{\sqrt{M}\,\ln M},\,
\frac{\pi}{\sqrt{M}}+\frac{1}{\sqrt{M}\,\ln M}\right] {\mbox{such that~$q(t_M)\in\Omega$.}}\end{equation}
To prove this, consistently with~\eqref{012ofek029rfeiwohg3uiwfegv-ytrgfFe7ywfhenNARTFGVIvr67002},
we parameterize~$\Omega$ near the origin as a subgraph of a function~$f:\R^{n-1}\to\R$
with~$f(0)=0$ and~$\partial_kf(0)=0$ for each~$k\in\{1,\dots,n-1\}$.

\begin{figure}
  \centering
  \includegraphics[height=.34\textheight]{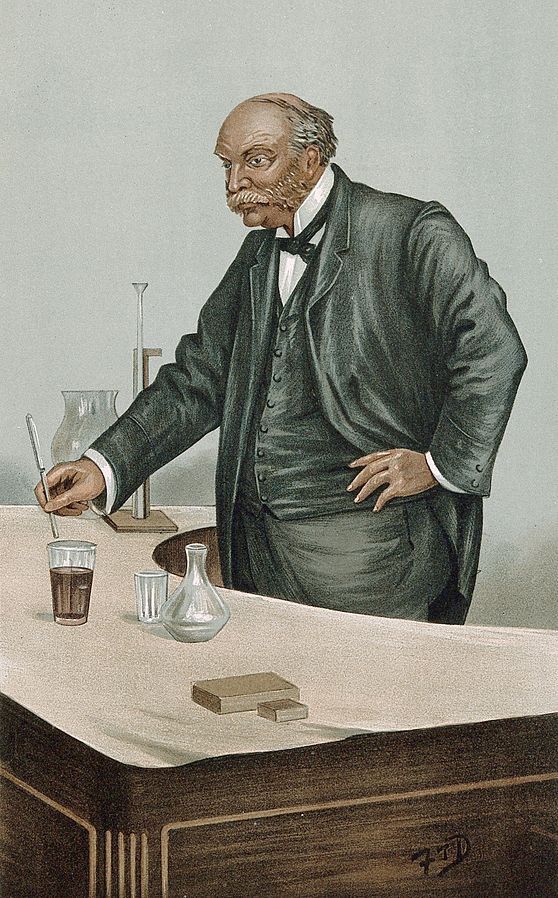}
 \caption{\sl Caricature of Lord Rayleigh
 (Public Domain image from
 Wikipedia).}\label{LORFASAaLOTGHSdNI8imCAttgeFI}
\end{figure}

In this setting, if~$t\in\left(0,\frac{\pi}{\sqrt{M}}-\frac{1}{\sqrt{M}\,\ln M}\right)$,
using the approximation in~\eqref{lkjb89io-oijohicsbFedf3vdJcbsdR4Hs2L5O3ks4d2ss23oer-1t},
\begin{equation}\label{lkjb89io-oijohicsbFedf3vdJcbsdR4Hs2L5O3ks4d2ss23oer-1t2}
\begin{split}&
q_n(t)-f(q'(t))\ge q_n(t)-C|q'(t)|^2\simeq
\frac{p_n(0)}{\sqrt{M}}\sin\big(\sqrt{M}t\big)- C|p'(0)|^2t^2\\&\qquad=
\frac{p_n(0)}{\sqrt{M}}\left( \sin\big(\sqrt{M}t\big)- \frac{C\big(\sqrt{M}t\big)^2}{\sqrt{M}}\right),
\end{split}\end{equation}
up to renaming~$C>0$ line after line.

Thus, if~$t\in\left(0,\frac{1}{\sqrt{M}}\right)$ we use that~$\sin\big(\sqrt{M}t\big)\ge\frac{\sqrt{M}t}2$,
whence, in the formal approximation above,
\begin{eqnarray*}&&
q_n(t)-f(q'(t))\ge
\frac{p_n(0)}{2\sqrt{M}}\left( \sqrt{M}t- \frac{C\big(\sqrt{M}t\big)^2}{\sqrt{M}}\right)=
\frac{p_n(0)\,t}{2}\left( 1- \frac{C \sqrt{M}t}{\sqrt{M}}\right)\ge
\frac{p_n(0)\,t}{2}\left( 1- \frac{C }{\sqrt{M}}\right)>0.
\end{eqnarray*}
If instead~$t\in\left(\frac1{\sqrt{M}},\frac{\pi}{\sqrt{M}}-\frac{1}{\sqrt{M}\,\ln M}\right]$,
we infer from~\eqref{lkjb89io-oijohicsbFedf3vdJcbsdR4Hs2L5O3ks4d2ss23oer-1t2} that
\begin{eqnarray*}&&
q_n(t)-f(q'(t))\ge
\frac{p_n(0)}{\sqrt{M}}\left( \sin\left(\pi-\frac1{\ln M}\right)- \frac{C}{\sqrt{M}}\right)\\&&\qquad=
\frac{p_n(0)}{\sqrt{M}}\left( \sin\left(\frac1{\ln M}\right)- \frac{C}{\sqrt{M}}\right)\ge
\frac{p_n(0)}{2\sqrt{M}}\left(  \frac1{\ln M}- \frac{C}{\sqrt{M}}\right)>0
\end{eqnarray*}
and so~\eqref{iujhbdNArryujf0RR92i3rjtgCahsdn9B7yuhdfGTGJA893-1} follows.

\begin{figure}
  \centering
  \includegraphics[height=.24\textheight]{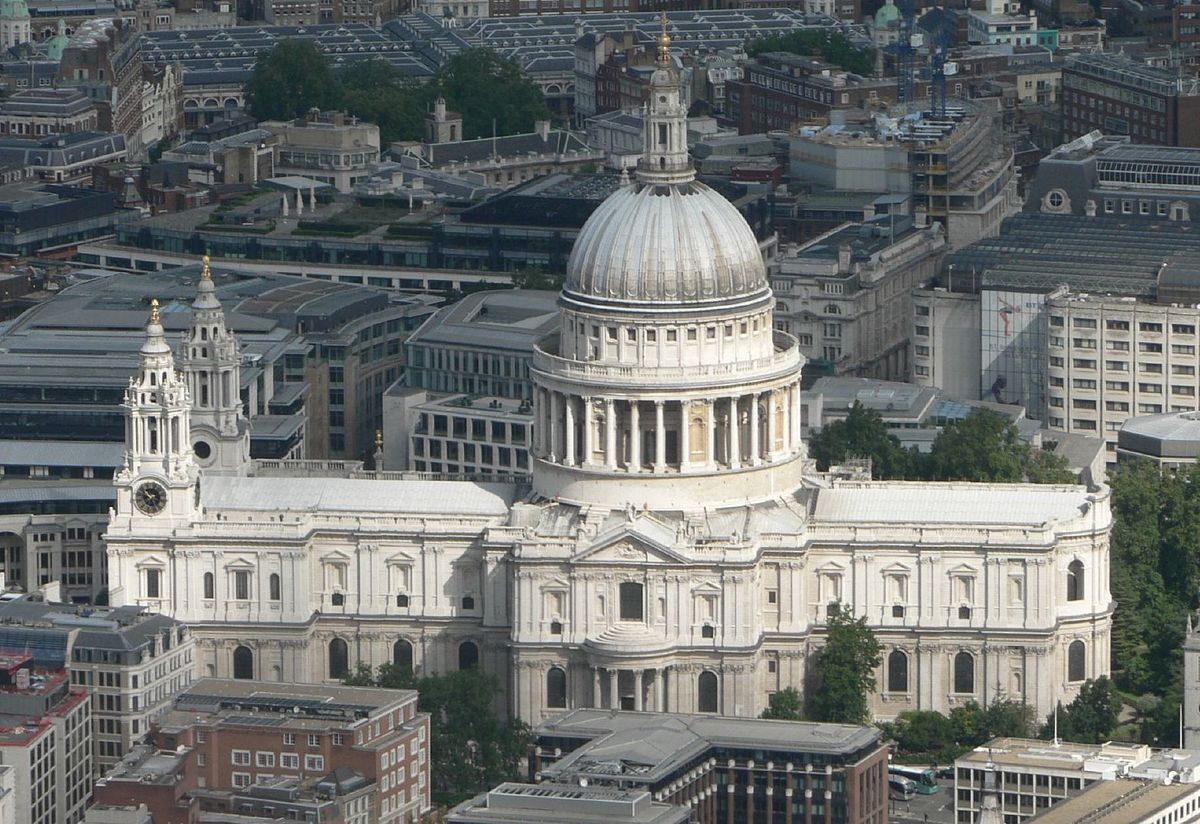} $\quad$
  \includegraphics[height=.24\textheight]{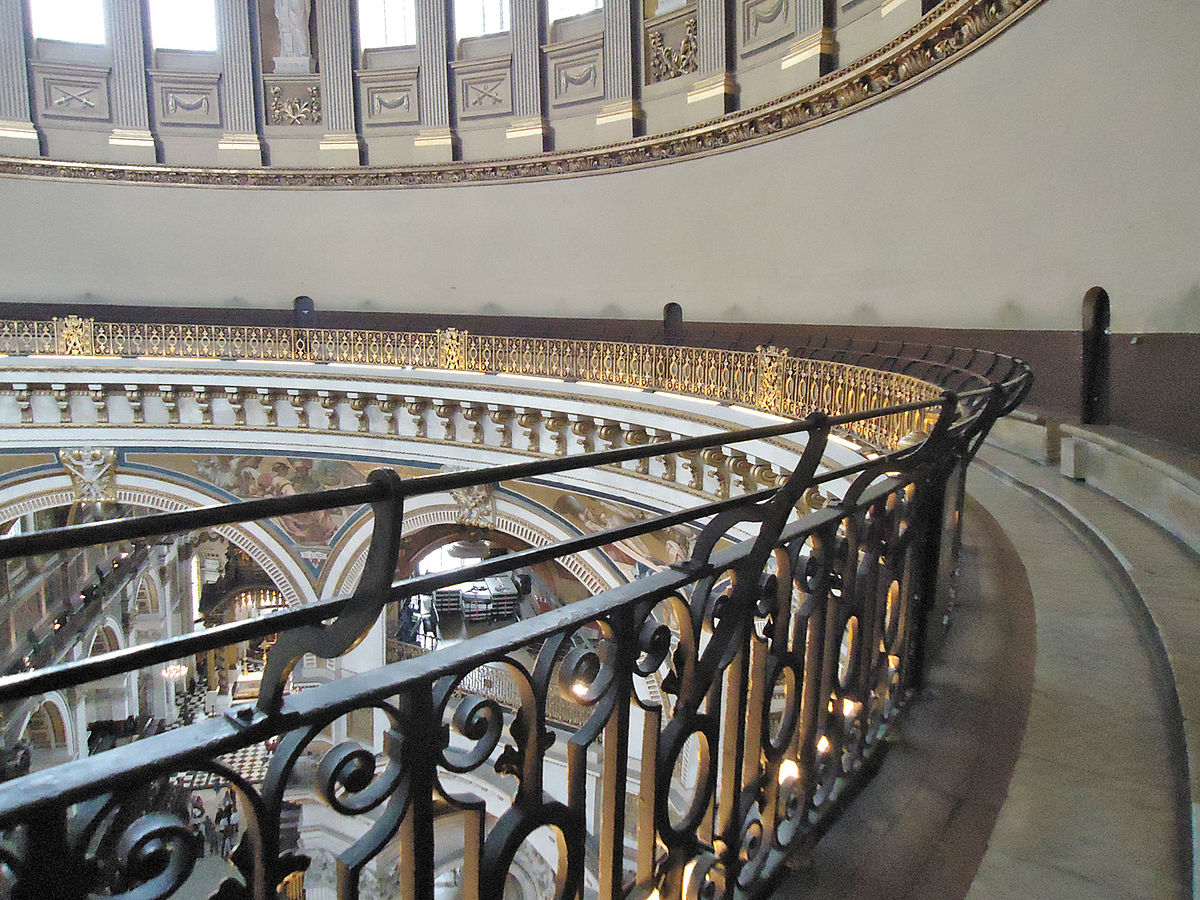}
 \caption{\sl Left: aerial view of Cathedral Church of St Paul the Apostle, London;
 right: its whispering gallery 
 (images from
 Wikipedia,
 photo by Mark Fosh,
 licensed under the Creative Commons Attribution 2.0 Generic license for the first,
 photo by Femtoquake, licensed under the Creative Commons Attribution-Share Alike 3.0 Unported license
 for the second).}\label{NI8imCAttgeFI}
\end{figure}

Furthermore,
\begin{eqnarray*}&&
q_n\left(\frac\pi{\sqrt{M}}+\frac1{{\sqrt{M}}\,\ln M}\right)-f\left(q'\left(\frac\pi{\sqrt{M}}+\frac1{{\sqrt{M}}\,\ln M}\right)\right)\\&\le&
q_n\left(\frac\pi{\sqrt{M}}+\frac1{{\sqrt{M}}\,\ln M}\right)+C\left|q'\left(\frac\pi{\sqrt{M}}+\frac1{{\sqrt{M}}\,\ln M}\right)\right|^2\\&\simeq&
\frac{p_n(0)}{\sqrt{M}}\sin\left(\pi+\frac1{\ln M}\right)+\frac{C}{M}\\&=&
-\frac{p_n(0)}{\sqrt{M}}\sin\left(\frac1{\ln M}\right)+\frac{C}{M}\\&<&0,
\end{eqnarray*}
that gives~\eqref{iujhbdNArryujf0RR92i3rjtgCahsdn9B7yuhdfGTGJA893-2}.

Therefore, the velocity of the particle reentering~$\Omega$ is of the form
$$ p(t_M)=\dot q(t_M)\simeq\left(
p'(0),\,p_n(0)\cos\big(\sqrt{M}t_M\big)
\right)=\Big(
p'(0),\,p_n(0)\cos\big(\pi+\theta_M\big)
\Big),
$$
with~$\theta_M\in\left( -\frac{1}{\ln M},\,\frac{1}{\ln M}\right]$.

In the limit as~$M\to+\infty$ this bouncing velocity becomes~$
\big(
p'(0),\,p_n(0)\cos\pi
\big)=\big(
p'(0),\,-p_n(0)\big)$, which is precisely the elastic reflection of geometric optics,
due to~\eqref{012ofek029rfeiwohg3uiwfegv-ytrgfFe7ywfhenNARTFGVIvr67002}.

This completes the formal proof of~\eqref{FOHAMM}.\medskip

The formal Hamiltonian presentation of billiards
as an infinite potential well in~\eqref{FOHAMM} is sometimes written as a description
of the billiard motion as a Hamiltonian system run by the Hamiltonian
\begin{equation}\label{HIKS-wkemfMIt4rgy5Ujhd-234L34S}
H_\infty(q,p):=\frac{|p|^2}2+V_\infty(q),\end{equation} with
$$ V_\infty(q):=\begin{dcases}
0 & {\mbox{ if }}q\in\Omega,\\
+\infty& {\mbox{ if }}q\in\R^n\setminus\Omega,
\end{dcases}$$
whatever this expression means.

\begin{figure}
  \centering
  \includegraphics[height=.18\textheight]{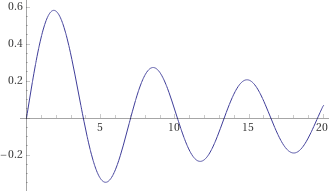} $\qquad$
  \includegraphics[height=.18\textheight]{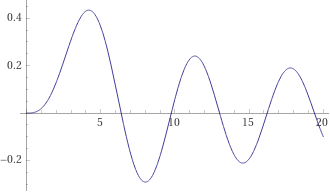} \\
  $\,$\\
    \includegraphics[height=.18\textheight]{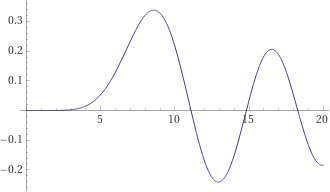} $\qquad$
  \includegraphics[height=.18\textheight]{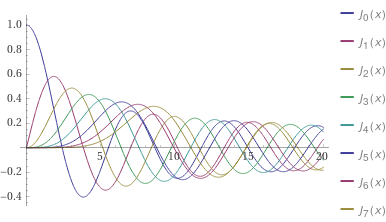}
 \caption{\sl Bessel functions of the first kind~$J_m$ (corresponding respectively
 to~$m=1$, $m=3$, $m=7$ and~$m\in\{0,\dots,7\}$).}\label{NI8imCAt5ujk883tgeFIBESSEL}
\end{figure}

One of the advantages of such a (perhaps teetering, but quite practical) description
is that it readily opens the possibility to study a quantum analogue of the classical
billiards, simply by considering
the quantization proposed in~\eqref{ENEJHND8iujMNSID-2krwjeg-21-ktg22}.
This method suggests to replace the classical Hamiltonian in~\eqref{HIKS-wkemfMIt4rgy5Ujhd-234L34S}
with an operator of the form
$$ \frac{|-i\hbar\nabla|^2}2+V_\infty(x)
=\frac{\hbar^2}2\sum_{k=1}^n \partial_k^2+V_\infty(x)=
\frac{\hbar^2}2\Delta+V_\infty(x).$$
The corresponding partial differential equation would thereby become
$$ Eu(x)=\frac{\hbar^2}2\Delta u(x)+V_\infty(x)u(x),$$
for some scalar~$E$.

Now, since~$V_\infty=+\infty$ outside~$\Omega$, to make sense of the right hand side above
one typically assume that~$u$ vanishes outside~$\Omega$ (that is,
the condition that the billiard particle is confined in~$\Omega$ is reflected into
a restriction on the support of the corresponding wave function in its quantum analogue).
This leads to the problem~$ Eu=\frac{\hbar^2}2\Delta u$, with~$u=0$ outside~$\Omega$,
or, up to a change of notation,
\begin{equation}\label{owjef924-1934i5tj1} \begin{dcases}
\Delta u=-\lambda u & {\mbox{ in }}\Omega,\\
u=0 &{\mbox{ on }}\partial\Omega,
\end{dcases}\end{equation}
where~$\lambda:=-\frac{2E}{\hbar^2}$ corresponds now to a large scalar
(given the smallness of the
reduced Planck constant~$\hbar$).
Equation~\eqref{owjef924-1934i5tj1} is sometimes called
\index{Helmholtz equation}
the Helmholtz equation.\medskip

\begin{figure}
  \centering
  \includegraphics[height=.18\textheight]{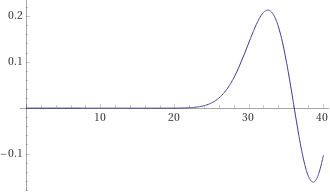} $\qquad$
  \includegraphics[height=.18\textheight]{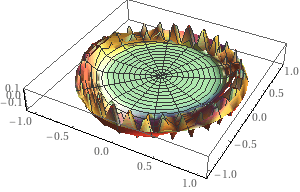}
 \caption{\sl Bessel function of the first kind~$J_{30}$ and the corresponding
 eigenfunction according to~\eqref{9oikm-0ok0infttimes}.}\label{COEIGENI8imCAt5ujk883tgeFIBESSEL}
\end{figure}

In analogy with the finite dimensional case of matrices, \eqref{owjef924-1934i5tj1} is considered as
\index{eigenvalue problem} an eigenvalue problem for the Laplace operator (in this case,
with homogeneous Dirichlet conditions).
The study of quantum billiards is thus often quite related to the analysis of an eigenvalue
problem, and often specifically focused on the case of large eigenvalues.\medskip

Interestingly, problems of this type also surface in the theory of sound\footnote{And
of course similar phenomena can take place for light and other electromagnetic
radiations, making whispering phenomena very attractive also for technological reasons.}
and provide an essential ingredient for the understanding of
whispering gallery waves. \index{whispering gallery waves}
This phenomenon was deeply investigated by 
John William Strutt, 3rd Baron Rayleigh (most commonly addressed as Lord\footnote{Here
is a funny anecdote about Lord Rayleigh.
Given his Anglican faith, at some point Lord Rayleigh wanted to include a
religious quotation from the Bible at the opening of a collection of papers of him.
He was discouraged from doing so
by the staff of the Cambridge University Press, as he reported:
``When I was bringing out my Scientific Papers I proposed a motto from the Psalms, {\em The Works of the Lord are great, sought out of all them that have pleasure therein}. The Secretary to the Press suggested with many apologies that the reader might suppose that I was the {\em Lord}''.

The moral of the story, if any, could also be that if one happens to be a member of the aristocracy
it could be advisable not to include the noble title in the author's name.}
Rayleigh,
see Figure~\ref{LORFASAaLOTGHSdNI8imCAttgeFI} for a caricature of him published in the London magazine Vanity Fair in~1899). 

To study the whispering phenomenon, Lord Rayleigh took inspiration from the
concrete case of the whispering gallery in St Paul's Cathedral, London,
see Figure~\ref{NI8imCAttgeFI}, in which whispers can be heard clearly from one side
to the other of the gallery.

Actually, Lord Rayleigh proposed several approaches to the question. In~\cite{14786441008636993}
he related the whispering phenomenon with problem~\eqref{owjef924-1934i5tj1}.
Specifically, we can model the whispering gallery in Figure~\ref{NI8imCAttgeFI} as
a planar disk, say~$\Omega=B_1\subset\R^2$.
We can take the model of propagation of sound waves put forth in Section~\ref{PROSAWA0qiodjwfl0owfej23t98y4823rty}, that is, in the light of~\eqref{WAYEBJJD121t416JH098327uyrhdhc832nbcM}, we consider solutions of
\begin{equation}\label{MR1625845-EQ1} \partial_{tt}\rho=c^2 \Delta\rho,\end{equation}
where~$\rho=\rho(x,t)$ represents the density of the air, as perturbed by the propagating sound,
and~$c>0$ is the speed of this propagation.

We assume that the density is constant, say~$\rho_0$, along~$\partial\Omega$,
we consider the eigenfunctions of the Laplacian with zero boundary datum along~$\partial\Omega$,
see~\cite[Section~6.5.1]{MR1625845} and we look for solutions in the form
\begin{equation}\label{MR1625845-EQ2} \rho(x,t)=\rho_0+\sum_{k=1}^{+\infty} \rho_k(t)\,\eta_k(x),\end{equation}
being~$\eta_k$ the eigenfunctions corresponding to the eigenvalue~$\lambda_k$,
with~$\lambda_k$ nondecreasing and~$\lambda_k\to+\infty$ as~$k\to+\infty$
(and, say~$\|\eta_k\|_{L^2(\Omega)}=1$ as a possible normalization).

\begin{figure}
                \centering
                \includegraphics[width=.34\linewidth]{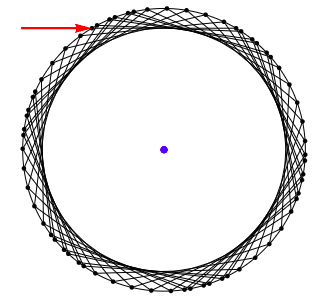}
        \caption{\sl Dynamics of circular billiard (obtained from the free
        resource {\tt https://demonstrations.wolfram.com/DynamicBilliardsInEllipse/}).}\label{HSCAUSttI}
\end{figure}

By plugging~\eqref{MR1625845-EQ2} into~\eqref{MR1625845-EQ1},
we obtain the equation~$\ddot\rho_k(t)=-c^2\lambda_k\rho_k(t)$,
giving e.g.~$\rho_k(t)=b_k\sin\big(c\sqrt{\lambda_k}t\big)$ for some~$b_k\in\R$, whence
$$ \rho(x,t)=\rho_0+\sum_{k=1}^{+\infty} b_k\sin\big(c\sqrt{\lambda_k}t\big)\,\eta_k(x).$$
This suggests that the propagation of sound is modeled by the linear superposition
of harmonic modes of the form~$b_k\sin\big(c\sqrt{\lambda_k}t\big)\,\eta_k(x)$.
Lord Rayleigh's idea is thus to detect modes corresponding to eigenfunctions~$\eta_k$
whose mass is mostly concentrated in the vicinity of~$\partial\Omega$.
{F}rom the technical point of view,
Lord Rayleigh's analysis consists in using polar coordinates~$(r,\vartheta)$ to write
conveniently~$\eta_k(x)=\beta_k(r)\,\alpha_k(\vartheta)$, for suitable functions~$\alpha_k$ and~$\beta_k$,
and obtain the radial component~$\beta_k$ in terms of ``known'' special functions,
in this case the so-called Bessel functions of the first kind. Then, he chooses conveniently the parameters
to exhibit cases in which the mass of these special functions is mostly concentrated near the boundary.

More explicitly, one can check that, for all~$m\in\N$, the function
\begin{equation}\label{9oikm-0ok0infttimes}
[0,+\infty)\times[0,2\pi)\ni
(r,\vartheta)\longmapsto J_m\big(\sqrt{\lambda}r\big)\,\cos(m\vartheta)\end{equation}
satisfies~\eqref{owjef924-1934i5tj1} provided that\footnote{Actually,
the Bessel differential equation~\eqref{OKSDD-1rfBEkdhhseLL76LD}
has two families of solutions,
say~$J_m$ and~$Y_m$.
The solutions in the first family are the ones
considered here: they are called
Bessel functions of the first kind and they
are continuous at the origin.
The solutions in the second family
are singular at the origin, therefore they do not
provide solutions for eigenvalue problems
in a disk (for this reason, they do not appear in the calculations
given in these pages): they are
called Bessel functions of the second kind,
or Neumann functions or Weber functions.
See e.g.~\cite{MR1349110} and the references therein for additional
information.}
the special function~$J_m$ is a solution of
\begin{equation}\label{OKSDD-1rfBEkdhhseLL76LD}\begin{dcases}
s^2 J_m''(s)+s J_m'(s)+\big( s^2 
-m^2\big) J_m(s)=0 & {\mbox{ for all }}s>0,\\
J_m\big(\sqrt{\lambda}\big)=1,&
\end{dcases}
\end{equation}
known as Bessel differential equation \index{Bessel differential equation}
(the boundary condition here being related to the boundary datum in~\eqref{owjef924-1934i5tj1}).
See Figures~\ref{NI8imCAt5ujk883tgeFIBESSEL} and~\ref{COEIGENI8imCAt5ujk883tgeFIBESSEL} for a sketch of these functions.

As we see from Figures~\ref{NI8imCAt5ujk883tgeFIBESSEL} and~\ref{COEIGENI8imCAt5ujk883tgeFIBESSEL},
these functions~$J_m$ for large values of~$m$ have the tendency to become flatter and flatter at the origin
and concentrate proportionally more mass near their first root. This fact can be proven rigorously,
see e.g. page~227 in~\cite{MR1349110}, and it can also be heuristically guessed from the Bessel differential equation,
since if near~$0$ we make the ansatz that~$J_m(s)=s^\ell+o(s^\ell)$ for some~$\ell\ge0$ we formally obtain
from~\eqref{OKSDD-1rfBEkdhhseLL76LD} that
\begin{eqnarray*}0&=&s^2 \big(\ell(\ell-1)s^{\ell-2}+o(s^{\ell-2})\big)
+s \big(\ell s^{\ell-1}+o(s^{\ell-1})\big)+\big( s^2 
-m^2\big) \big(s^\ell+o(s^\ell)\big)\\&=&
\big(\ell(\ell-1)+\ell -m^2\big) s^{\ell}+o(s^\ell)\\&=&
\big(\ell^2 -m^2\big) s^{\ell}+o(s^\ell),
\end{eqnarray*}
hence~$0=\ell^2-m^2+o(1)$ and consequently~$\ell=m$.

The possible localization of eigenfunctions in a small region around the boundary
of the domain is particularly evident in the second image of Figure~\ref{COEIGENI8imCAt5ujk883tgeFIBESSEL}.
\medskip

The whispering phenomena described by the concentration of eigenfunctions
in the proximity of the boundary can also be considered the counterpart of the caustics
in dynamical billiards (namely, curves to which the trajectory remains tangent): in concrete cases,
when these caustics are localized near the boundary, they can provide a confinement for
the trajectory, see Figure~\ref{HSCAUSttI} (to be compared with the second image
in Figure~\ref{COEIGENI8imCAt5ujk883tgeFIBESSEL}).
\medskip

See~\cite{MR3038118} and the references therein for further information about
localization of eigenfunctions, Bessel functions,
whispering gallery waves and their links to bouncing balls.\medskip

In general, dynamical billiards have been widely studied from many mathematical viewpoints,
also as excellent training camps for a number of complex systems, yet many problems 
related to billiards are still completely open
and mathematicians often feel to be behind the eight ball.
See~\cite{MR3617212} and the references therein
for additional information about billiards and their intimate connections
with the Laplace operator.

\subsection{Strange attractions}

Another pioneering work by Lord Rayleigh (see Figure~\ref{LORFASAaLOTGHSdNI8imCAttgeFI})
consisted in the study of thermal convection
\index{thermal convection}
problems, addressed in~\cite{zbMATH02613875}.
The problem can be explained by considering, for simplicity,
a two-dimensional incompressible and viscous fluid in a gravity field, with coordinates~$(x,z)\in\R^2$.
The top and bottom surfaces of the fluid are horizontal, say at level~$z=H$ and~$z=0$ respectively,
supposed to be kept at constant temperature.
Up to normalization, we suppose that the top surface is at temperature~$\overline{T}$ and
the bottom surface at some temperature~$\underline{T}$, which we suppose to be larger than~$\overline{T}$.
 
\begin{figure}
                \centering
                \includegraphics[width=.75\linewidth]{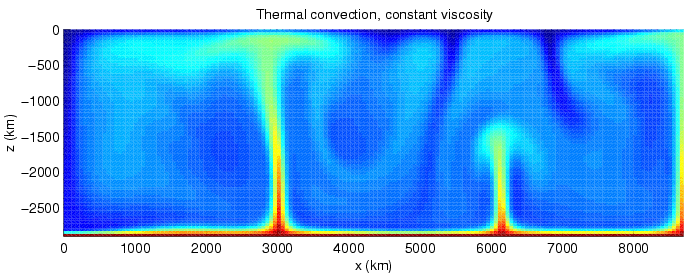}
        \caption{\sl Simulation of thermal convection in the Earth's mantle
        (by Harroschmeling; image from
        Wikipedia,
        licensed under the Creative Commons Attribution-Share Alike 3.0 Unported license).}\label{2HAFODAKEDCONVE7solFUMHDNOJHNFOJED231}
\end{figure}

The idea of this model is that gravity tends to push down the molecules of the fluid,
but heat from the bottom may favor upward motions, thus producing interesting convective patterns,
see e.g. Figure~\ref{2HAFODAKEDCONVE7solFUMHDNOJHNFOJED231}.

To describe this phenomenon, one can utilize the Navier-Stokes equation in~\eqref{EUMSD-OS-DNS}
complemented with the incompressibility condition in~\eqref{INCOMPRE}
and the assumption that the top and bottom surfaces remain unchanged: writing
the velocity of the fluid as~$v=(v_x,v_z)\in\R^2$,
this leads to
\begin{equation}\label{GBS-conMA-EQ-LRGHS-01}
\begin{dcases}
\rho \partial_t v+\rho (v\cdot\nabla) v=\mu\Delta v
-\rho g(0,1)-\nabla p&{\mbox{ in }}\R\times(0,H),\\
\div v=0 &{\mbox{ in }}\R\times(0,H),\\
v_z(x,0)=v_z(x,H)=0.&
\end{dcases}\end{equation}
We recall that~$\rho$ in this framework is the density of the fluid, $t$ is time and~$p$ is pressure.

To simplify this setting, one can perform some convenient approximations.
First of all, the density is supposed to depend only on the temperature~$T$, which in turn
evolves according to the heat equation in~\eqref{DAGB-ADkrVoiweLL4re2346ytmngrrUj7}, with no external sources, 
but with the caveat that the particles are moving with the velocity field. This gives that
$$ \frac{d}{dt} T(x(t),z(t),t)=\kappa\Delta T$$
for some diffusion coefficient~$\kappa>0$, and therefore, since~$(\dot x(t),\dot z(t))=v(x(t),z(t),t)$,
$$\partial_tT+v\cdot\nabla T=\kappa\Delta T .$$
We suppose that the higher the temperature, the lower the density of the fluid:
this can be understood
by considering that when the temperature increases the molecules move faster
and bump into each other more frequently, consequently spreading apart,
taking up more space and reducing density (and this reduction of density at higher temperature
is the reason for which high temperature parcels of fluid happen to be lighter and
have the tendency to move up). The simplest possible model to consider is that
in which the density variation is linear with respect to the temperature, say~$\rho=\rho(T)=\rho_0\big(1-\alpha( T-\overline{T})\big)$,
for some~$\alpha>0$, which we consider as a small parameter (that is, the density does
depend on the temperature, but it is almost constant).

The other simplifying assumption is that the velocity and the velocity variations of the fluid are small.
Therefore, we will disregard quadratic terms in~$v$, as well as terms involving the product
of the two small quantities~$\alpha$ and~$v$ (or~$\alpha$ and the derivatives of~$v$).
With this motivation, we have that
\begin{equation}\label{SIMPLEXRELO}
\rho \partial_t v+\rho (v\cdot\nabla) v=
\rho_0\big(1-\alpha( T-\overline{T})\big)\partial_t v+\rho_0\big(1-\alpha( T-\overline{T})\big) (v\cdot\nabla) v
\simeq \rho_0\partial_t v.
\end{equation}

These observations and~\eqref{GBS-conMA-EQ-LRGHS-01} suggest\footnote{One can
observe that, when passing from~\eqref{GBS-conMA-EQ-LRGHS-01}
to~\eqref{GBS-conMA-EQ-LRGHS-02},
essentially we have just considered the density to be constant
except when it appears in terms multiplied by the gravity acceleration~$g$.
In a sense, this approximation is based on the ansatz that
the difference in density is negligible in itself but the gravity field is supposed to be
sufficiently strong to make it appreciable.
This method is called \index{Boussinesq approximation} Boussinesq approximation.}
to study the problem
\begin{equation}\label{GBS-conMA-EQ-LRGHS-02}
\begin{dcases}
\rho_0 \partial_t v=\mu\Delta v
-\rho_0\big(1-\alpha( T-\overline{T})\big)g(0,1)-\nabla p&{\mbox{ in }}\R\times(0,H),\\
\div v=0 &{\mbox{ in }}\R\times(0,H),\\
v_z(x,0)=v_z(x,H)=0,&\\
\partial_tT+v\cdot\nabla T=\kappa\Delta T& {\mbox{ in }}\R\times(0,H),\\
T(x,0)=\underline{T},\\
T(x,H)=\overline{T}.
\end{dcases}\end{equation}

We observe that this system of equations possesses a simple solution
in which the fluid is at rest, i.e. $v=0$, and the temperature increases linearly, i.e. $T(x,z)=\underline{T}+
\frac{\big(\overline{T}-\underline{T}\big)z}H$. Indeed, substituting these relations in~\eqref{GBS-conMA-EQ-LRGHS-02},
we obtain the pressure prescription
$$\nabla p=-\rho_0\big(1-\alpha( T-\overline{T})\big)g(0,1)=-\rho_0\left[1-\alpha\big(\underline{T}-\overline{T}\big)\left( 1-\frac{z}H\right)\right]g(0,1),$$
which
is satisfied by choosing
$$p(x,z)=p_0-\rho_0g\left[z-\alpha\big(\underline{T}-\overline{T}\big)\left( z-\frac{z^2}{2H}\right)\right],$$
for some~$p_0\in\R$.

\begin{figure}
                \centering
                \includegraphics[width=14cm,height=4.5cm]{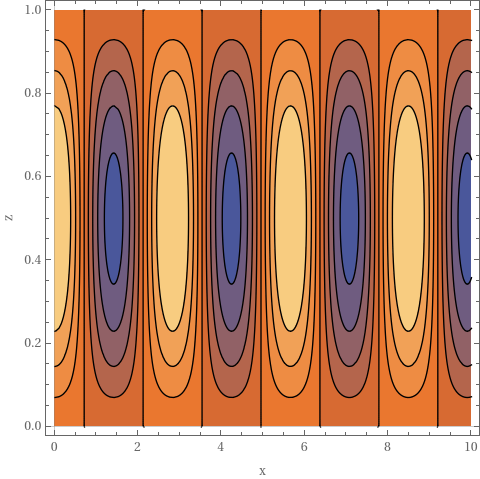}
        \caption{\sl Level sets of the function~$\cos\left(\frac{\pi x}{\sqrt2}\right)\,\sin(\pi z)$.}\label{2HAFLORROJED231MURA2}
\end{figure}

It is therefore interesting to seek solutions of~\eqref{GBS-conMA-EQ-LRGHS-02} which
bifurcate from these basic states: for this we look for solutions of~\eqref{GBS-conMA-EQ-LRGHS-02}
in which
\begin{equation}\label{TSOSTTHET} T=\underline{T}+\frac{\big(\overline{T}-\underline{T}\big)z}H+\Theta\qquad{\mbox{and}}\qquad
p=p_0-\rho_0g\left[z-\alpha\big(\underline{T}-\overline{T}\big)\left( z-\frac{z^2}{2H}\right)\right]+P,\end{equation}
for some unknown~$\Theta=\Theta(x,z)$ and~$P=P(x,z)$.

We insert~\eqref{TSOSTTHET} into~\eqref{GBS-conMA-EQ-LRGHS-02} and we obtain
the following system of equations for the new unknown functions:
\begin{equation}\label{GBS-conMA-EQ-LRGHS-03}
\begin{dcases}
\rho_0 \partial_t v=\mu\Delta v
+\alpha\rho_0g \Theta\,(0,1)
-\nabla P
&{\mbox{ in }}\R\times(0,H),\\
\div v=0 &{\mbox{ in }}\R\times(0,H),\\
v_z(x,0)=v_z(x,H)=0,&\\
\partial_t\Theta-\frac{\big(\underline{T}-\overline{T}\big)v_z}H
+v\cdot\nabla\Theta
=\kappa\Delta \Theta& {\mbox{ in }}\R\times(0,H),\\
\Theta(x,0)=0,\\
\Theta(x,H)=0.
\end{dcases}\end{equation}

\begin{figure}
                \centering
                \includegraphics[width=.6\linewidth]{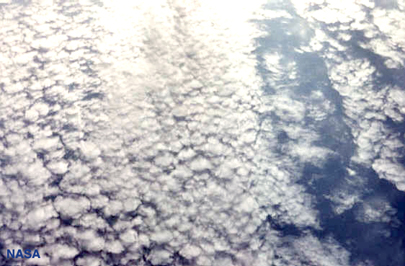}
        \caption{\sl Altocumulus clouds formed by convective activity as seen from shuttle
        (Public Domain image from
 Wikipedia).}\label{2auth29HAFLORROJED231MURA2}
\end{figure}

One can also consider a scaled version of~\eqref{GBS-conMA-EQ-LRGHS-03}
by setting
\begin{equation}\label{scaling00}\begin{split}
&(\bar{x},\bar{z}):=\frac{(x,z)}H, \qquad\bar{t}:=\frac{\kappa t}{H^2},\qquad
\bar{v}(\bar{x},\bar{z},\bar{t}):=\frac{Hv(H\bar{x},H\bar{z},H^2\kappa^{-1}\bar{t})}\kappa\\&
\bar\Theta(\bar{x},\bar{z},\bar{t}):=\frac{\alpha\rho_0gH^3\Theta(H\bar{x},H\bar{z},H^2\kappa^{-1}\bar{t})}{\mu\kappa}
\qquad{\mbox{and}}\qquad
\bar{P}(\bar{x},\bar{z},\bar{t}):=\frac{H^2P(H\bar{x},H\bar{z},H^2\kappa^{-1}\bar{t})}{\mu\kappa},\end{split}
\end{equation}
thus obtaining
\begin{equation}\label{GBS-conMA-EQ-LRGHS-04}
\begin{dcases}
\partial_{\bar t}\bar v=\frac{\mu}{\rho_0\kappa}\left( \Delta_{\bar{x}\bar{z}}\bar v
+\bar\Theta\,(0,1)
-\nabla\bar P\right)
&{\mbox{ in }}\R\times(0,1),\\
\div \bar v=0 &{\mbox{ in }}\R\times(0,1),\\
\bar v_{\bar z}(x,0)=\bar v_{\bar z}(x,1)=0,&\\
\partial_{\bar t}\bar\Theta+\bar v\cdot\nabla\bar\Theta= \Delta\bar \Theta+\frac{\alpha\rho_0gH^3\big(\underline{T}-\overline{T}\big)\,\bar v_{\bar z}}{\mu\kappa}& {\mbox{ in }}\R\times(0,1),\\
\bar\Theta(x,0)=0,\\
\bar\Theta(x,1)=0.
\end{dcases}\end{equation}
The parameter
$${\mathcal{P}}:=\frac{\mu}{\rho_0\kappa}$$ is called \index{Prandtl number}
the Prandtl number.

The parameter
$${\mathcal{R}}:=\frac{\alpha\rho_0gH^3\big(\underline{T}-\overline{T}\big)}{\mu\kappa}$$
is called the Rayleigh number \index{Rayleigh number} and takes into consideration the thermal gradient
(notice indeed that~${\mathcal{R}}$ is large when the temperature gap~$\underline{T}-\overline{T}$
between the bottom and the top of the fluid is large).

{F}rom now on, to ease the notation, we drop the bars in~\eqref{GBS-conMA-EQ-LRGHS-04}.
Also, as observed in~\eqref{0pod-45ia7n7b9e0r} since we are considering a two-dimensional setting,
the divergence free condition in~\eqref{GBS-conMA-EQ-LRGHS-04} allows us to write~$v=(-\partial_z\psi,\partial_x\psi)$,
for some scalar function~$\psi:\R\times(0,1)\to\R$.

\begin{figure}
                \centering
                \includegraphics[width=.21\linewidth]{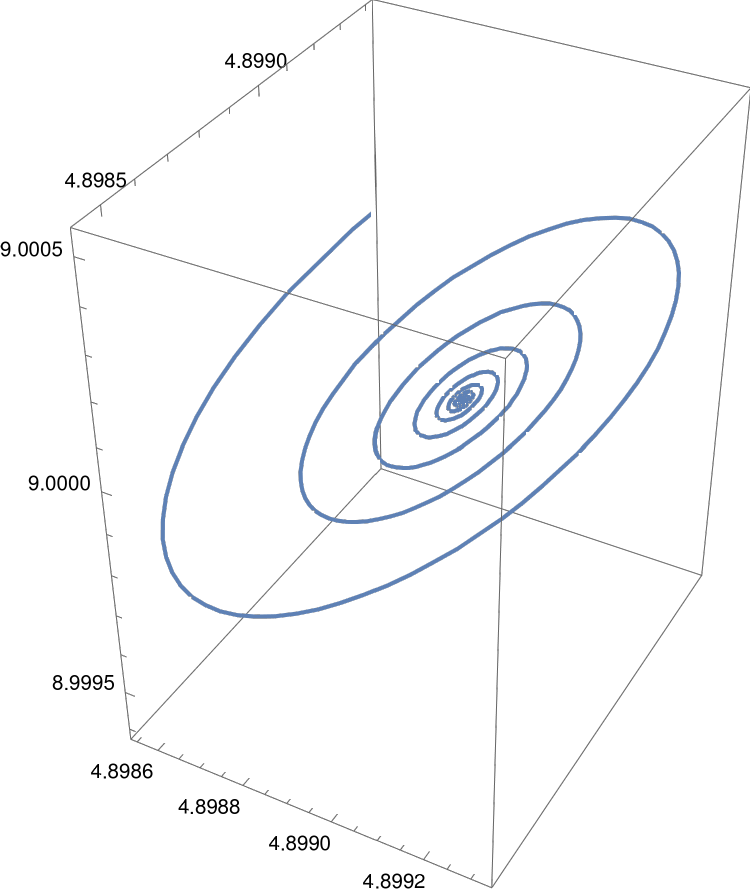}$\quad$
                \includegraphics[width=.21\linewidth]{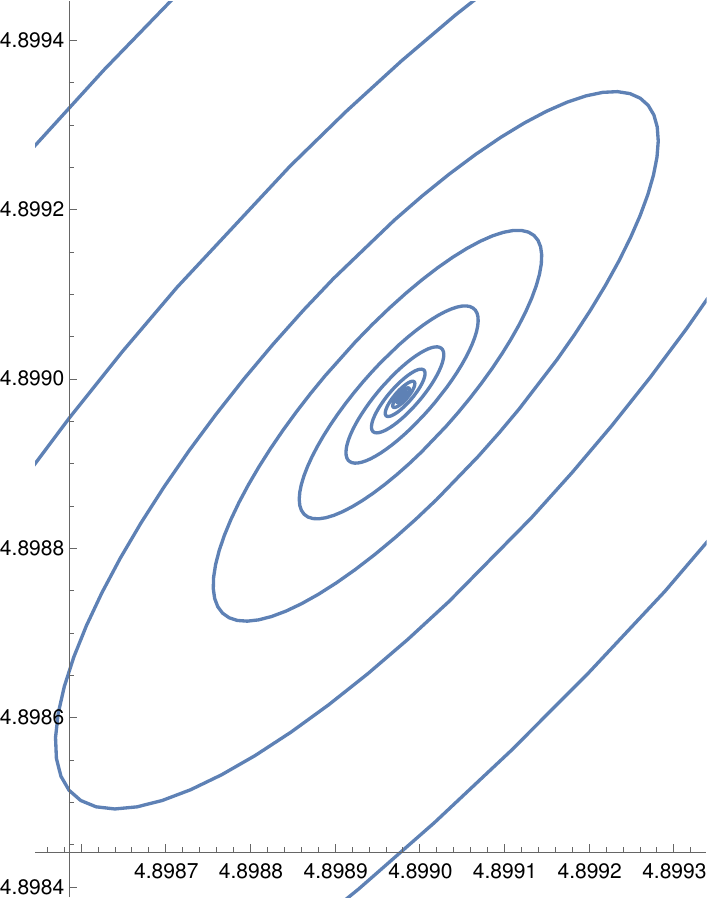}$\quad$
                \includegraphics[width=.21\linewidth]{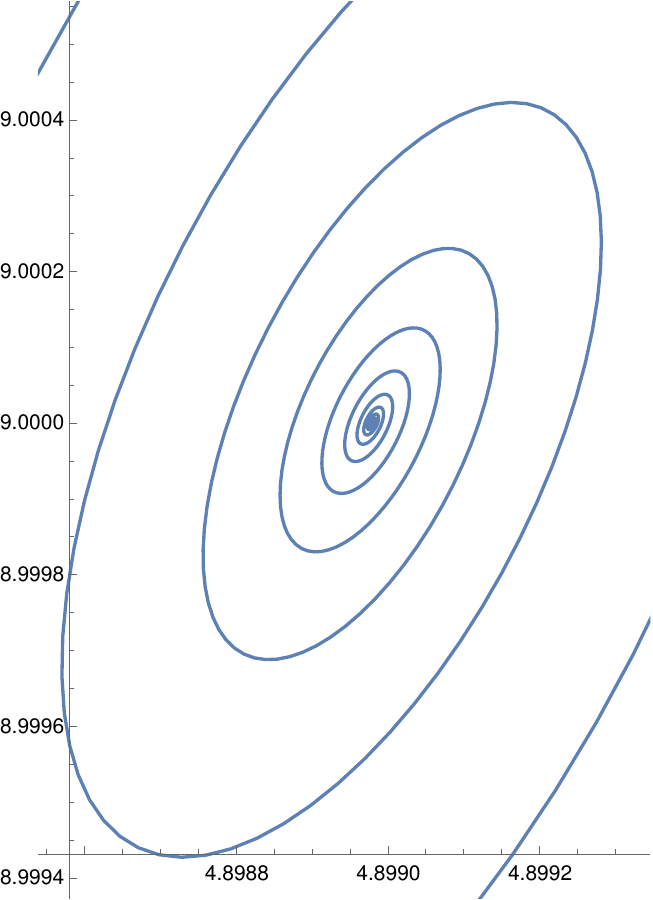}$\quad$
                \includegraphics[width=.21\linewidth]{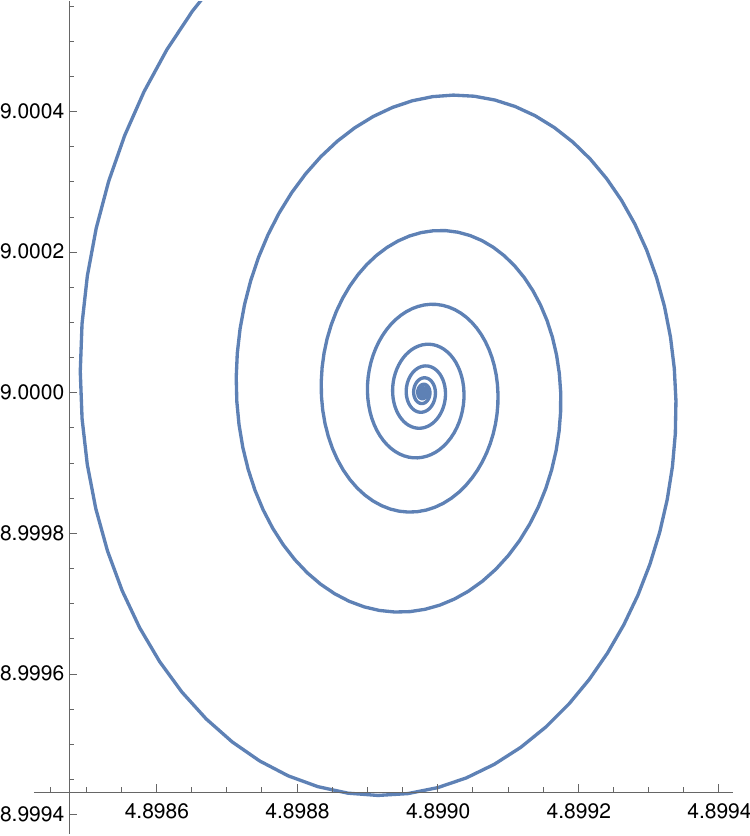}              
        \caption{\sl Mathematica simulation of the trajectory
of the Lorenz system with parameters~$\sigma:= 10$, $\rho:=10$, $\beta:=\frac83$,
starting at the point~$(10,10,10)$
(three-dimensional plot and projections onto the coordinate planes).}\label{2HAFLORRE31MURA2-STAB}
\end{figure}

In this way, we have that
$$ v\cdot\nabla\Theta=(-\partial_z\psi,\partial_x\psi)\cdot(\partial_x\Theta,\partial_z\Theta)=
-\partial_z\psi\partial_x\Theta+\partial_x\psi\partial_z\Theta=
\det\left( \begin{matrix} \partial_x\psi&\partial_z\psi\cr \partial_x\Theta&\partial_z\Theta
\end{matrix}\right)=:{\mathcal{J}}(\psi,\Theta).
$$
Therefore, the temperature equation in~\eqref{GBS-conMA-EQ-LRGHS-04} can be written in the form
$$
\partial_{t}\Theta+{\mathcal{J}}(\psi,\Theta)= \Delta \Theta+{\mathcal{R}}\partial_x\psi.$$

As a result, \eqref{GBS-conMA-EQ-LRGHS-04} becomes
\begin{equation}\label{GBS-conMA-EQ-LRGHS-05}
\begin{dcases}\partial_t \partial_z\psi={\mathcal{P}}\left( \Delta \partial_z\psi
+\partial_x P\right)
&{\mbox{ in }}\R\times(0,1),\\
\partial_t \partial_x\psi={\mathcal{P}}\left( \Delta\partial_x\psi
+\Theta
-\partial_z P\right) &{\mbox{ in }}\R\times(0,1),\\
\partial_x\psi(x,0)=\partial_x\psi(x,1)=0,&\\
\partial_{t}\Theta+{\mathcal{J}}(\psi,\Theta)= \Delta \Theta+{\mathcal{R}}\partial_x\psi& {\mbox{ in }}\R\times(0,1),\\
\Theta(x,0)=\Theta(x,1)=0.
\end{dcases}\end{equation}
We now perform an additional approximation, by considering both~$\psi$ and~$\Theta$ as small
quantities, together with their derivatives, and thus dropping the quadratic term in~\eqref{GBS-conMA-EQ-LRGHS-05}.
With this, we reduce ourselves to
\begin{equation}\label{GBS-conMA-EQ-LRGHS-05b}
\begin{dcases}\partial_t \partial_z\psi={\mathcal{P}}\left( \Delta \partial_z\psi
+\partial_x P\right)
&{\mbox{ in }}\R\times(0,1),\\
\partial_t \partial_x\psi={\mathcal{P}}\left( \Delta\partial_x\psi
+\Theta
-\partial_z P\right) &{\mbox{ in }}\R\times(0,1),\\
\partial_x\psi(x,0)=\partial_x\psi(x,1)=0,&\\
\partial_{t}\Theta= \Delta \Theta+{\mathcal{R}}\partial_x\psi& {\mbox{ in }}\R\times(0,1),\\
\Theta(x,0)=\Theta(x,1)=0.
\end{dcases}\end{equation}
We observe that
\begin{equation}\label{PKM-3kfk3902-2-2}
{\mbox{the first two equations in~\eqref{GBS-conMA-EQ-LRGHS-05b} are equivalent to }}\,
\partial_t \Delta\psi={\mathcal{P}}\big( \Delta^2\psi
+\partial_x \Theta\big),
\end{equation}
where~$\Delta^2$ is the Laplace operator applied twice. 

\begin{figure}
                \centering
                \includegraphics[width=.21\linewidth]{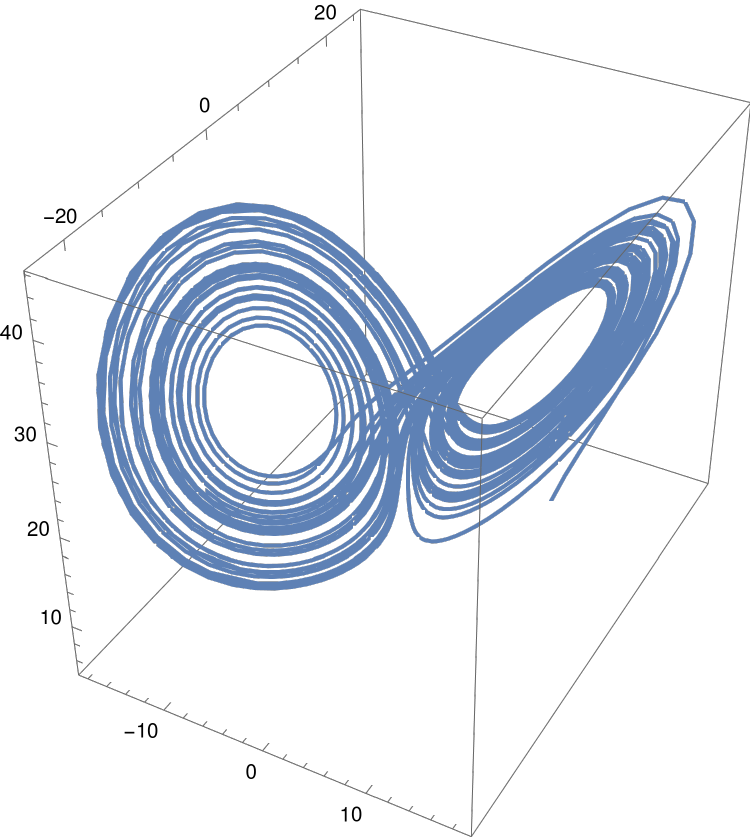}$\quad$
                \includegraphics[width=.21\linewidth]{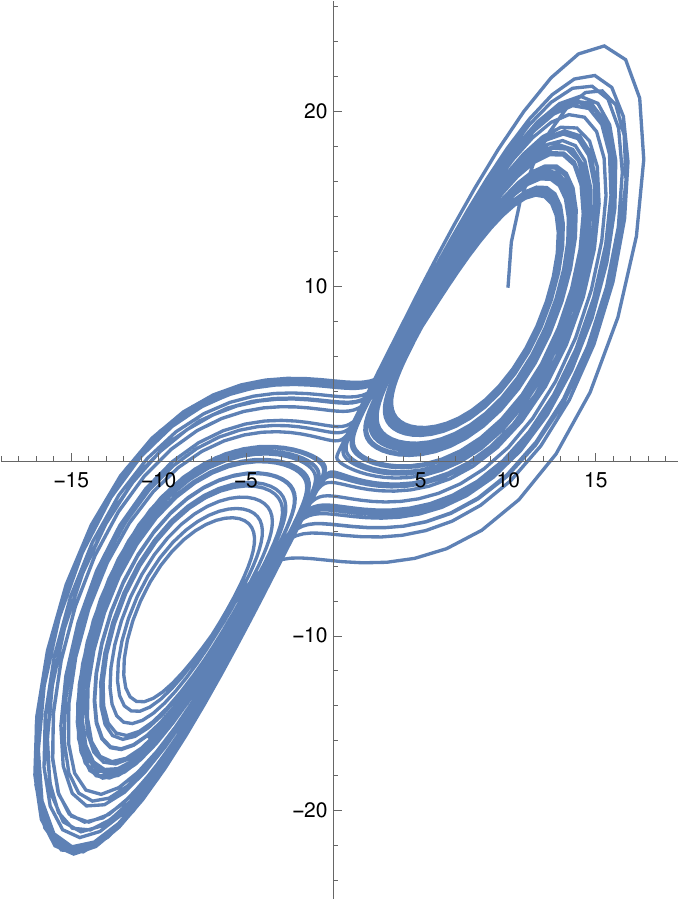}$\quad$
                \includegraphics[width=.21\linewidth]{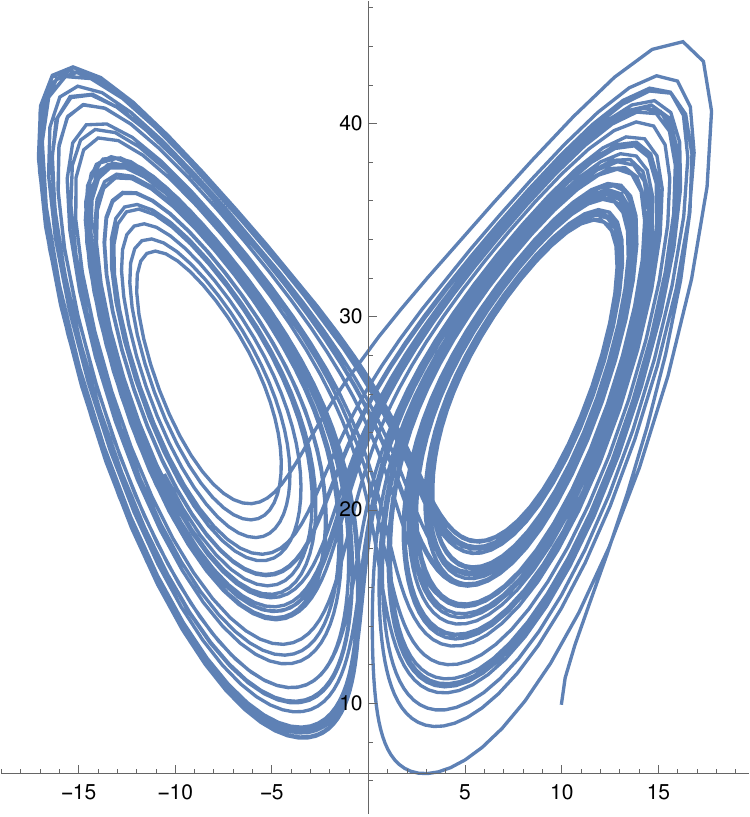}$\quad$
                \includegraphics[width=.21\linewidth]{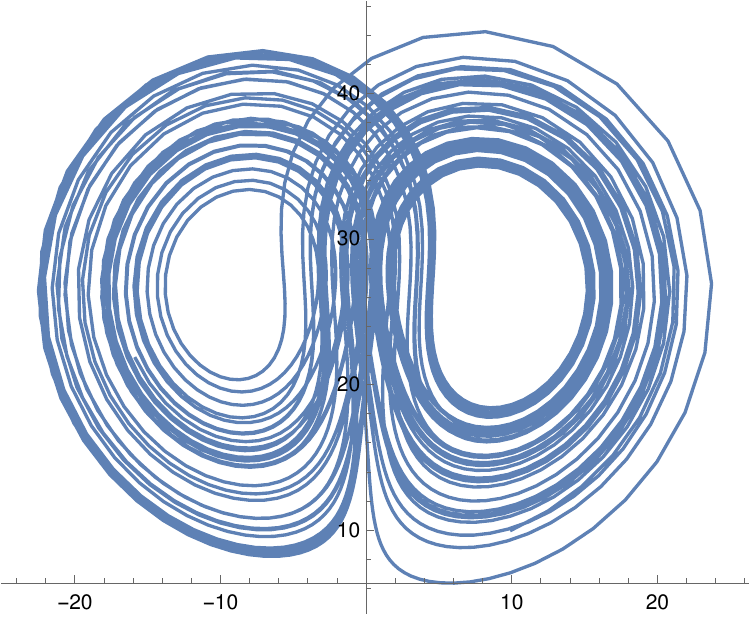}           
        \caption{\sl Mathematica simulation of the trajectory
of the Lorenz system with parameters~$\sigma:= 10$, $\rho:=28$, $\beta:=\frac83$,
starting at the point~$(10,10,10)$
(three-dimensional plot and projections onto the coordinate planes).}\label{2HAFLORRE31MURA2}
\end{figure}

Indeed, if the first two equations in~\eqref{GBS-conMA-EQ-LRGHS-05b} hold true, then
\begin{eqnarray*}
\partial_t \Delta\psi&=&
\partial_x\Big(\partial_t \partial_x\psi\Big)+\partial_z\Big(\partial_t \partial_z\psi\Big)
\\&=&
\partial_x\Big({\mathcal{P}}\big( \Delta\partial_x\psi
+\Theta
-\partial_z P\big)\Big)+\partial_z\Big({\mathcal{P}}\big( \Delta\partial_z\psi
+ \partial_x P\big)\Big)
\\&=&{\mathcal{P}}\big( \Delta^2\psi
+\partial_x \Theta\big).
\end{eqnarray*}
Conversely, if the equation in~\eqref{PKM-3kfk3902-2-2} is satisfied,
we consider the differential form
$$ \omega:=
\left(\frac1{\mathcal{P}}\partial_t \partial_z\psi-\Delta \partial_z\psi\right)\,dx
-\left(\frac1{\mathcal{P}}\partial_t \partial_x\psi- \Delta\partial_x\psi
-\Theta\right)\,dz$$
and we observe that
\begin{eqnarray*}
d\omega&=&
-\left[\partial_z\left(\frac1{\mathcal{P}}\partial_t \partial_z\psi-\Delta \partial_z\psi\right)
+\partial_x\left(\frac1{\mathcal{P}}\partial_t \partial_x\psi- \Delta\partial_x\psi
-\Theta\right)\right]\,dx\wedge dz
\\&=&-\frac1{{\mathcal{P}}}\big[\partial_t\Delta\psi-{\mathcal{P}}\Delta^2\psi
-{\mathcal{P}}\partial_x\Theta\big]\,dx\wedge dz\\&=&0.
\end{eqnarray*}
Consequently, using the
Poincar\'e Lemma (see e.g.~\cite[Theorem~8.3.8]{MR1209437}), we obtain that there exists a scalar function~$P:\R\times(0,1)\to\R$ for which~$\omega=dP$, hence
$$ \frac1{\mathcal{P}}\partial_t \partial_z\psi-\Delta \partial_z\psi=\partial_xP\qquad
{\mbox{ and }}\qquad
-\frac1{\mathcal{P}}\partial_t \partial_x\psi+ \Delta\partial_x\psi
+\Theta=\partial_zP,$$ which are precisely the first two equations in~\eqref{GBS-conMA-EQ-LRGHS-05b},
thereby completing the proof of~\eqref{PKM-3kfk3902-2-2}.

As a result, by~\eqref{PKM-3kfk3902-2-2} we can rewrite~\eqref{GBS-conMA-EQ-LRGHS-05b} in the form
\begin{equation}\label{GBS-conMA-EQ-LRGHS-05c}
\begin{dcases}\partial_t \Delta\psi={\mathcal{P}}\big( \Delta^2\psi
+\partial_x \Theta\big) &{\mbox{ in }}\R\times(0,1),\\
\partial_x\psi(x,0)=\partial_x\psi(x,1)=0,&\\
\partial_{t}\Theta= \Delta \Theta+{\mathcal{R}}\partial_x\psi& {\mbox{ in }}\R\times(0,1),\\
\Theta(x,0)=\Theta(x,1)=0.
\end{dcases}\end{equation}
Note that this formulation is convenient not only because it involves one less equation,
but also because we got rid of the dependence on the pressure, hence we also have one less unknown.

One can now look for solutions in the form of exponential and trigonometric functions,
by taking the pressure to be constant.
For instance, if we seek solutions of~\eqref{GBS-conMA-EQ-LRGHS-05c} in the form
\begin{equation}\label{TSOSTTHET134}
\begin{split}&\psi= \psi_0\,\sin(qx)\,\sin(k\pi z)\,e^{\lambda t},\\ &
\Theta= \Theta_0\,\cos(qx)\,\sin(k\pi z)\,e^{\lambda t},
\end{split}\end{equation}
for some~$q\in\R$, $k\in\N$ and~$\lambda\in\cOMPL$, when we insert these expressions into~\eqref{GBS-conMA-EQ-LRGHS-05c} we obtain that
\begin{equation*}
\begin{dcases}-\lambda(q^2+k^2\pi^2)\psi_0={\mathcal{P}}\big( (q^2+k^2\pi^2)^2\psi_0-q\Theta_0\big),\\
\lambda\Theta_0=-(q^2+k^2\pi^2)\Theta_0+{\mathcal{R}}q\psi_0.
\end{dcases}\end{equation*}
This system can be conveniently written in a matrix form by
$$ MV=\lambda V,\qquad{\mbox{where }}\quad V:=\left(\begin{matrix}\psi_0\\ \Theta_0\end{matrix}\right),\qquad
M:=\left(\begin{matrix}-{\mathcal{P}}\gamma^2 &{\mathcal{P}}q/\gamma^2\cr
{\mathcal{R}}q&-\gamma^2
\end{matrix}\right)\quad
{\mbox{ and }}\quad\gamma:=\sqrt{q^2+k^2\pi^2}.$$
The question is now {\em for what choice of physical parameters this system
allows solutions corresponding to~$\lambda>0$}. This is extremely relevant in practice,
because in such a situation the basic state in~\eqref{TSOSTTHET},
corresponding to a steady fluid with linearly increasing temperature, becomes unstable,
due to the exponential in time divergence in~\eqref{TSOSTTHET134}.

\begin{figure}
  \centering
  \includegraphics[height=5cm]{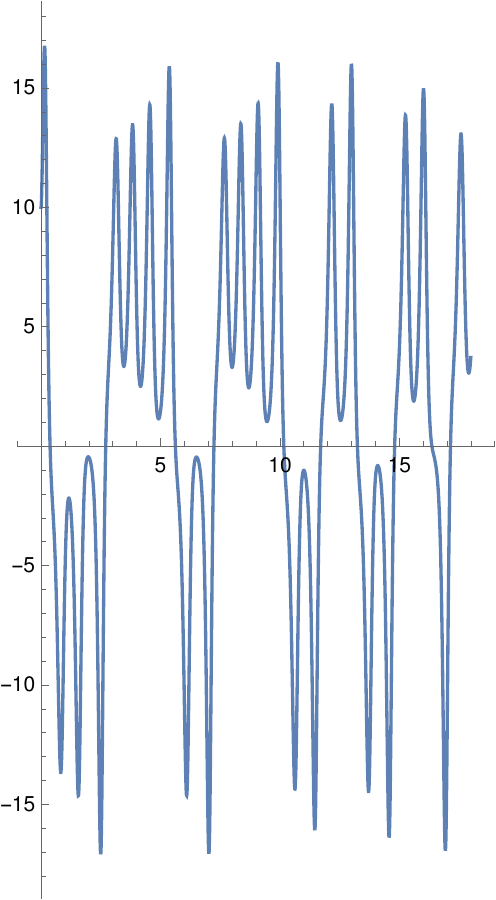} $\quad$ 
  \includegraphics[height=5cm]{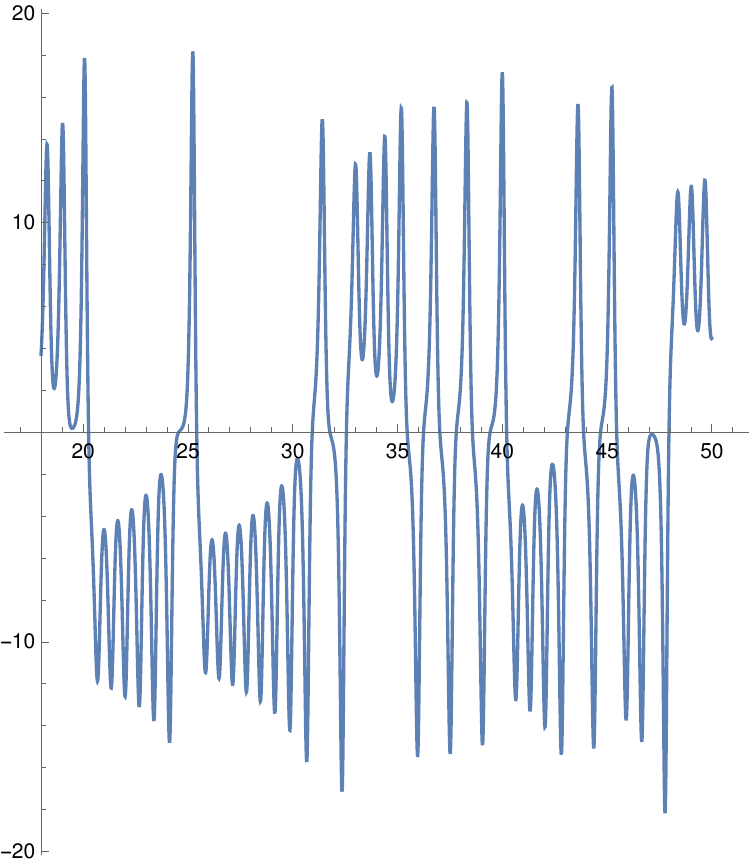}
  \\
  \includegraphics[height=5cm]{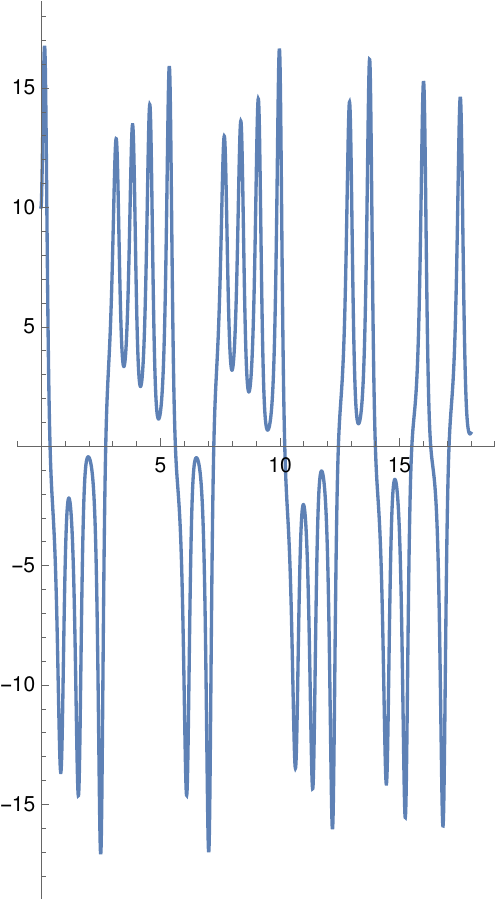}$\quad$
  \includegraphics[height=5cm]{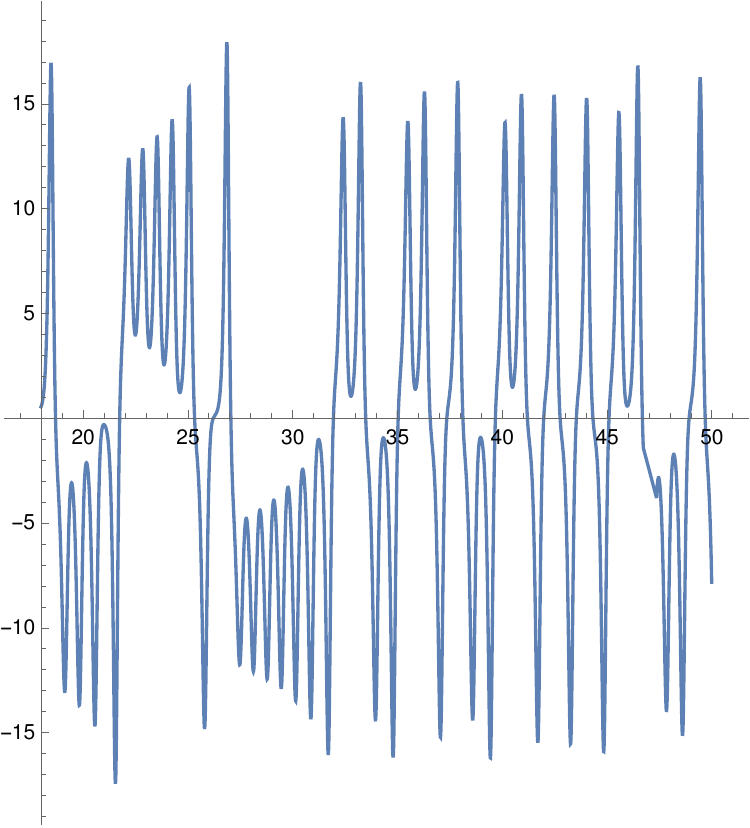}
 \caption{\sl Mathematica simulation of the $X$-trajectories
of the Lorenz system with parameters~$\sigma:= 10$, $\rho:=28$, $\beta:=\frac83$,
starting at the points~$(10,10,10)$ (above) and~$\left(10+\frac1{100},10,10\right)$ (below).}\label{NEc11AahMA0rtXcdtt345pdkfDIRDItangeFI}
\end{figure}

Thus, finding a positive eigenvalue~$\lambda$ for~$M$ is equivalent to look for~$\lambda>0$ solving
$$ 0=(-{\mathcal{P}}\gamma^2-\lambda)(-\gamma^2-\lambda)-\frac{{\mathcal{P}}{\mathcal{R}}q^2}{\gamma^2}=
\lambda^2+({\mathcal{P}}+1)\gamma^2\lambda+\frac{{\mathcal{P}}}{\gamma^2}\left(\gamma^6
- {\mathcal{R}}q^2\right),
$$
which, solving in~$\lambda$, gives
\begin{eqnarray*} 2\lambda&=&
-({\mathcal{P}}+1)\gamma^2\pm\sqrt{ ({\mathcal{P}}+1)^2\gamma^4-
\frac{4{\mathcal{P}}}{\gamma^2}\left(\gamma^6- {\mathcal{R}}q^2\right)}\\&=&
-({\mathcal{P}}+1)\gamma^2\pm\frac{\sqrt{
( {\mathcal{P}}-1)^2\gamma^6 + 4 {\mathcal{P}}{\mathcal{R}}q^2 }}{\gamma}
.\end{eqnarray*}
Notice that the above radicand is nonnegative
hence a real solution does exist: however, to have~$\lambda>0$ we must impose that $$
({\mathcal{P}}+1)^2\gamma^6<
( {\mathcal{P}}-1)^2\gamma^6 + 4 {\mathcal{P}}{\mathcal{R}}q^2,$$
that is
\begin{equation}\label{JMS-TEBFDrey}
{\mathcal{R}}>\frac{\gamma^6}{q^2}.
\end{equation}
This condition thus ensures that 
if the Rayleigh number is sufficiently large, i.e., when the thermal gradient
between the bottom and the top of the fluid is sufficiently large, the (approximated) convection problem
presents new stable solutions which bifurcate from the basic state as prescribed in~\eqref{TSOSTTHET134}.
We also note that~\eqref{JMS-TEBFDrey} can be written in the form
$$ {\mathcal{R}}>\frac{(q^2+k^2\pi^2)^3}{q^2},$$
the minimal threshold of~${\mathcal{R}}>\frac{27\pi^4}{4}$
corresponding to the choice~$k:=1$ and~$q:=\frac\pi{\sqrt2}$.

\begin{figure}
  \centering
  \includegraphics[height=6.9cm]{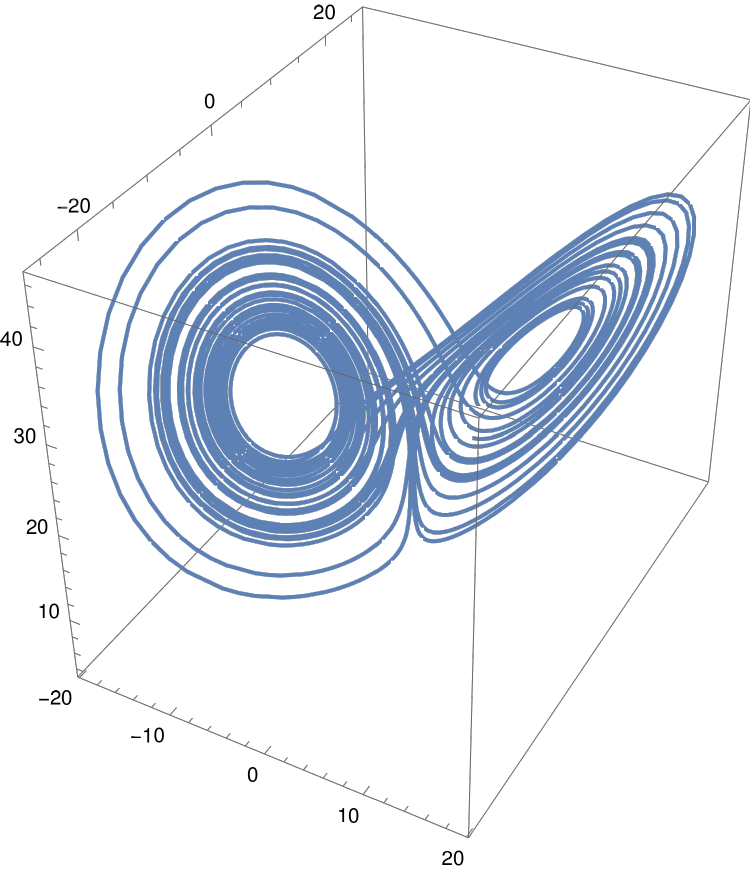}$\quad$
  \includegraphics[height=6.9cm]{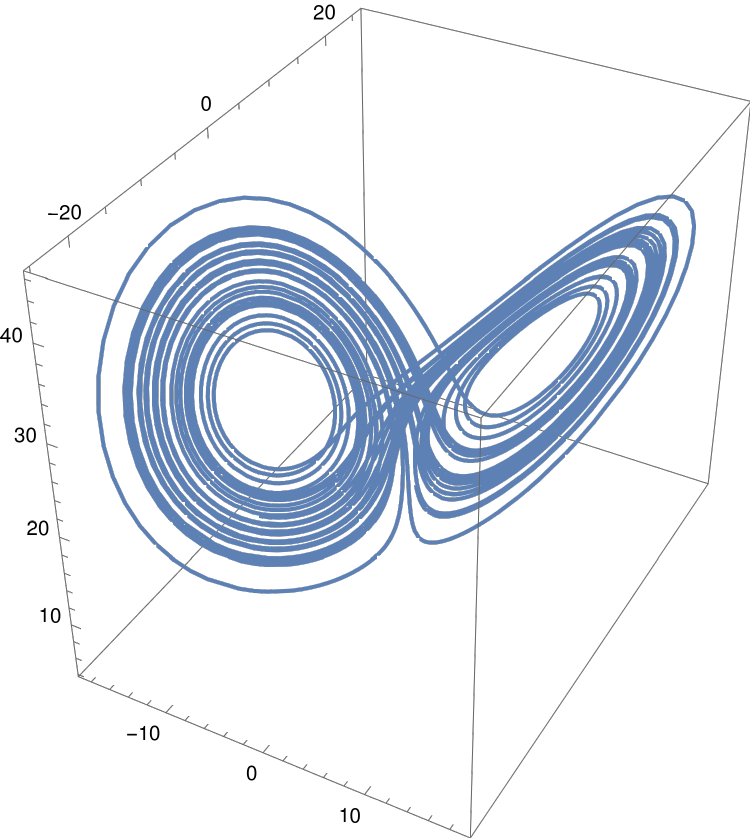}
 \caption{\sl Mathematica simulation of the trajectories
of the Lorenz system with parameters~$\sigma:= 10$, $\rho:=28$, $\beta:=\frac83$,
starting at the points~$(10,10,10)$ (left) and~$\left(10+\frac1{100},10,10\right)$ (right) for~$t\in[18,50]$.}\label{NEc11AahMA0rtXcdtt34LOREatt5pdkfDIRDItangeFI}
\end{figure}

It is thus interesting to plot the level sets of the function~$\cos\left(\frac{\pi x}{\sqrt2}\right)\,\sin(\pi z)$,
since they correspond to the long time dominant asymptotics of the temperature correction in~\eqref{TSOSTTHET134},
see Figure~\ref{2HAFLORROJED231MURA2}.

Actually, Figure~\ref{2HAFLORROJED231MURA2} describes a typical regular pattern arising in convection phenomena,
with the formation of distinctive convection cells known as \index{B\'enard cells} B\'enard cells,
see e.g. Figure~\ref{2auth29HAFLORROJED231MURA2}.
See also~\cite{MR1629996} and the references therein for
additional information about convection patterns.
\medskip

Interestingly, the description of the convection phenomena by Lord Rayleigh
is intimately related with one of the paradigmatic models used in the study of chaotic dynamical
systems. Indeed, mathematician and meteorologist Edward Norton Lorenz in~\cite{DeterministicNonperiodicFlow}
proposed to retake~\eqref{GBS-conMA-EQ-LRGHS-05c} and consider also the nonlinear interactions
of the form~$(v\cdot\nabla) v$ and~${\mathcal{J}}(\psi,\Theta)$ that were dropped
to have a simpler and more manageable expression (recall~\eqref{SIMPLEXRELO}
and~\eqref{GBS-conMA-EQ-LRGHS-05}).

This leads\footnote{We point out that~${\mathcal{J}}(\psi,\Delta\psi)$ in~\eqref{GBS-conMA-EQ-LRGHS-05c-NONLI} is obtained
by keeping the term~$\rho_0(v\cdot\nabla)v$ in~\eqref{SIMPLEXRELO}. In this way,
the first equation in~\eqref{GBS-conMA-EQ-LRGHS-03} becomes
$$ \rho_0 \partial_t v+\rho_0(v\cdot\nabla)v=\mu\Delta v
+\alpha\rho_0g \Theta\,(0,1)-\nabla P\quad
{\mbox{ in }}\R\times(0,H)$$
and therefore, after the scaling in~\eqref{scaling00} (and dropping the bars),
the first equation in~\eqref{GBS-conMA-EQ-LRGHS-04} becomes
$$
\partial_{ t}\bar v+(v\cdot\nabla)v={\mathcal{P}}\left( \Delta  v
+\Theta\,(0,1)-\nabla P\right)\quad
{\mbox{ in }}\R\times(0,1).$$
Hence, using the divergence free condition in~\eqref{GBS-conMA-EQ-LRGHS-04}, we write~$v=(-\partial_z\psi,\partial_x\psi)$,
for some scalar function~$\psi:\R\times(0,1)\to\R$, and we replace the first two equations in~\eqref{GBS-conMA-EQ-LRGHS-05} with
$$ \begin{dcases}\partial_t \partial_z\psi+\partial_x\psi\Delta\psi={\mathcal{P}}\left( \Delta \partial_z\psi
+\partial_x P\right)
&{\mbox{ in }}\R\times(0,1),\\
\partial_t \partial_x\psi-\partial_z\psi\Delta\psi={\mathcal{P}}\left( \Delta\partial_x\psi
+\Theta
-\partial_z P\right) &{\mbox{ in }}\R\times(0,1).
\end{dcases} $$
Assuming that both~$\psi$ and~$\Theta$ are small quantity, we obtain that the first two equations
in~\eqref{GBS-conMA-EQ-LRGHS-05b} become
$$\begin{dcases}\partial_t \partial_z\psi+\partial_x\psi\Delta\psi={\mathcal{P}}\left( \Delta \partial_z\psi
+\partial_x P\right)
&{\mbox{ in }}\R\times(0,1),\\
\partial_t \partial_x\psi-\partial_z\psi\Delta\psi={\mathcal{P}}\left( \Delta\partial_x\psi
+\Theta
-\partial_z P\right) &{\mbox{ in }}\R\times(0,1).
\end{dcases}$$
These equations are equivalent to
$$ \partial_t\Delta\psi+{\mathcal{J}}(\psi,\Delta\psi)={\mathcal{P}}\big(\Delta^2\psi +\partial_x\Theta\big),$$
which is the first equation in~\eqref{GBS-conMA-EQ-LRGHS-05c-NONLI}.}to the study of the system
\begin{equation}\label{GBS-conMA-EQ-LRGHS-05c-NONLI}
\begin{dcases}\partial_t \Delta\psi+{\mathcal{J}}(\psi,\Delta\psi)={\mathcal{P}}\big( \Delta^2\psi
+\partial_x \Theta\big) &{\mbox{ in }}\R\times(0,1),\\
\partial_x\psi(x,0)=\partial_x\psi(x,1)=0,&\\
\partial_{t}\Theta+{\mathcal{J}}(\psi,\Theta)= \Delta \Theta+{\mathcal{R}}\partial_x\psi& {\mbox{ in }}\R\times(0,1),\\
\Theta(x,0)=\Theta(x,1)=0.
\end{dcases}\end{equation}
Inspired by~\eqref{TSOSTTHET134}, Lorenz proposed to look for approximate solutions in the form
\begin{equation}\label{TSOSTTHET134-LORE}
\begin{split}&\psi= {\mathcal{A}}(t)\,\sin(qx)\,\sin(\pi z),\\ &
\Theta= {\mathcal{B}}(t)\,\cos(qx)\,\sin(\pi z)+{\mathcal{C}}(t)\sin(2\pi z),
\end{split}\end{equation}
for some~$q\in\R$ and functions~${\mathcal{A}}$, ${\mathcal{B}}$ and~${\mathcal{C}}$ to be determined.

Notice that, with this choice, the nonlinear term~${\mathcal{J}}(\psi,\Delta\psi)=0$.

\begin{figure}
  \centering
  \includegraphics[height=6.9cm]{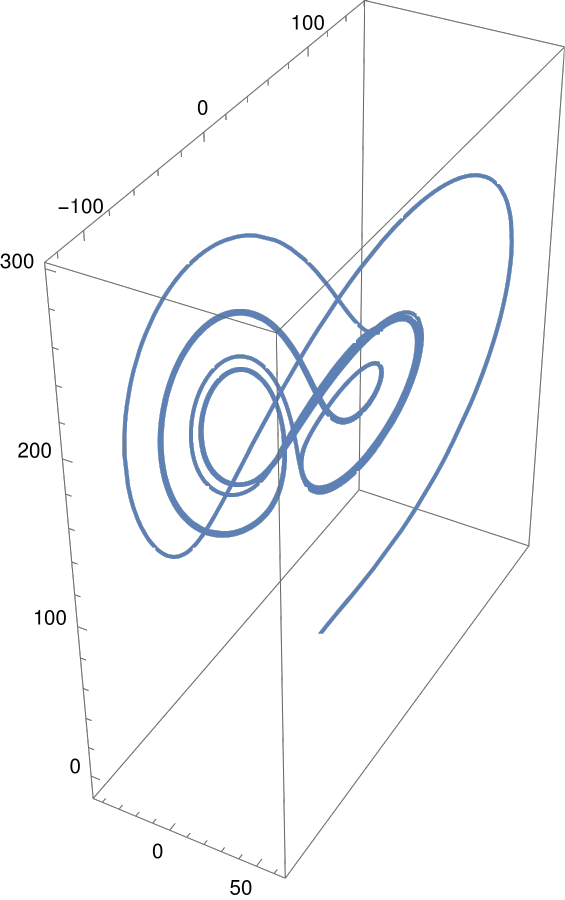}$\quad$
  \includegraphics[height=6.9cm]{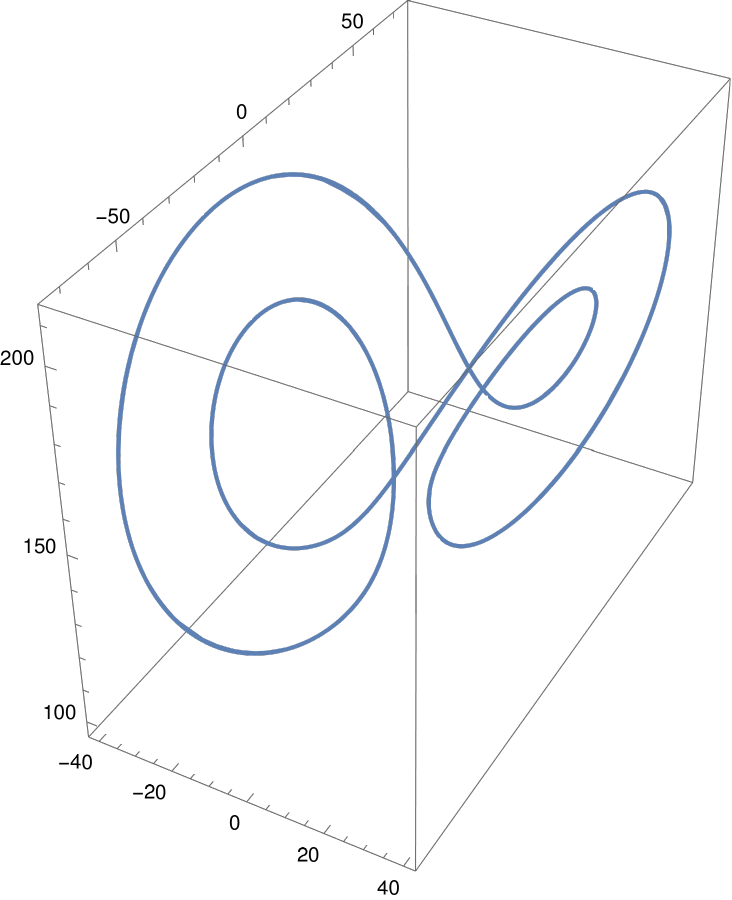}
 \caption{\sl Mathematica simulations of the trajectories
of the Lorenz system with parameters~$\sigma:= 10$, $\rho:=155$, $\beta:=\frac83$,
starting at the point~$(10,10,10)$ with~$t\in[0,100]$ (left) and~$t\in[90,100]$ (right).}\label{NEc11AahMA4srdtfyvhbjnklm9a99345tRDItangeFI}
\end{figure}

By plugging~\eqref{TSOSTTHET134-LORE} into~\eqref{GBS-conMA-EQ-LRGHS-05c-NONLI} one obtains
\begin{equation}\label{TSOSTTHET134-LORE-sinf49}
\begin{dcases}
-(q^2+\pi^2) \dot{\mathcal{A}}
={\mathcal{P}}\big( (q^2+\pi^2)^2{\mathcal{A}}-q{\mathcal{B}}\big),\\
\Big[\dot{\mathcal{B}}+(q^2+\pi^2) {\mathcal{B}}-{\mathcal{R}}q{\mathcal{A}}
\Big]\,\cos(qx)\,\sin(\pi z)+\left[
\frac{\pi q{\mathcal{A}}{\mathcal{B}}}2+
\dot{\mathcal{C}}+4\pi^2{\mathcal{C}}\right]\sin(2\pi z)\\ \qquad\qquad
+2\pi q{\mathcal{A}}{\mathcal{C}}\cos(qx)\,\sin(\pi z)\cos(2\pi z)
=0.\end{dcases}\end{equation}
It is also useful to apply the identity
\begin{eqnarray*} 2\sin(\pi z)\cos(2\pi z)
=\frac{(e^{\pi z}-e^{-\pi z})(e^{2\pi z}+e^{-2\pi z})}{2i}
=
\frac{e^{3\pi z}-e^{-3\pi z}-e^{\pi z}+e^{-\pi z}}{2i}
=\sin(3\pi z)-\sin(\pi z)
\end{eqnarray*}
to manipulate the last term in~\eqref{TSOSTTHET134-LORE-sinf49} and reduce ourselves to
\begin{equation}\label{TSOSTTHET134-LORE-sinf49-2}
\begin{dcases}
-(q^2+\pi^2) \dot{\mathcal{A}}
={\mathcal{P}}\big( (q^2+\pi^2)^2{\mathcal{A}}-q{\mathcal{B}}\big),\\
\Big[\dot{\mathcal{B}}+(q^2+\pi^2) {\mathcal{B}}-{\mathcal{R}}q{\mathcal{A}}
-\pi q{\mathcal{A}}{\mathcal{C}}\Big]\,\cos(qx)\,\sin(\pi z)+\left[
\frac{\pi q{\mathcal{A}}{\mathcal{B}}}2+
\dot{\mathcal{C}}+4\pi^2{\mathcal{C}}\right]\sin(2\pi z)\\ \qquad\qquad
+\pi q{\mathcal{A}}{\mathcal{C}}\cos(qx)\,\sin(3\pi z)
=0.\end{dcases}\end{equation}
Now, comparing with the approximation put forth in~\eqref{TSOSTTHET134-LORE},
the term containing~$\sin(3\pi z)$ has a spatial dependence of higher frequency than
what the previous assumption was accounting for. Therefore, 
the suggestion by Lorenz was to {\em drop this term} and simplify~\eqref{TSOSTTHET134-LORE-sinf49-2}
into
\begin{equation*}
\begin{dcases}
-(q^2+\pi^2) \dot{\mathcal{A}}
={\mathcal{P}}\big( (q^2+\pi^2)^2{\mathcal{A}}-q{\mathcal{B}}\big),\\
\Big[\dot{\mathcal{B}}+(q^2+\pi^2) {\mathcal{B}}-{\mathcal{R}}q{\mathcal{A}}
-\pi q{\mathcal{A}}{\mathcal{C}}\Big]\,\cos(qx)\,\sin(\pi z)+\left[
\frac{\pi q{\mathcal{A}}{\mathcal{B}}}2+
\dot{\mathcal{C}}+4\pi^2{\mathcal{C}}\right]\sin(2\pi z)
=0,\end{dcases}\end{equation*}
leading to
\begin{equation}\label{TSOSTTHET134-LORE-sinf49-34}
\begin{dcases}
-(q^2+\pi^2) \dot{\mathcal{A}}
={\mathcal{P}}\big( (q^2+\pi^2)^2{\mathcal{A}}-q{\mathcal{B}}\big),\\
\dot{\mathcal{B}}={\mathcal{R}}q{\mathcal{A}}-(q^2+\pi^2) {\mathcal{B}}+\pi q{\mathcal{A}}{\mathcal{C}},\\
-\dot{\mathcal{C}}=
\frac{\pi q{\mathcal{A}}{\mathcal{B}}}2+4\pi^2{\mathcal{C}}.
\end{dcases}\end{equation}

\begin{figure}
  \centering
  \includegraphics[height=6.9cm]{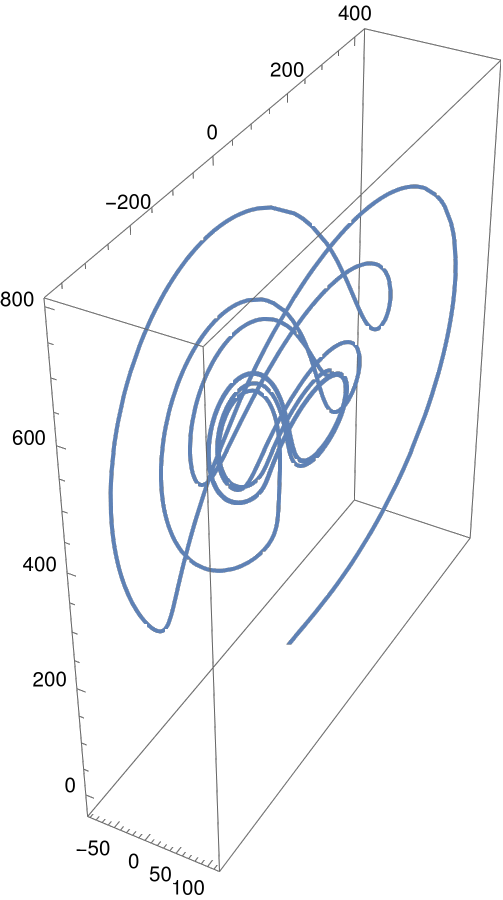}$\quad$
  \includegraphics[height=6.9cm]{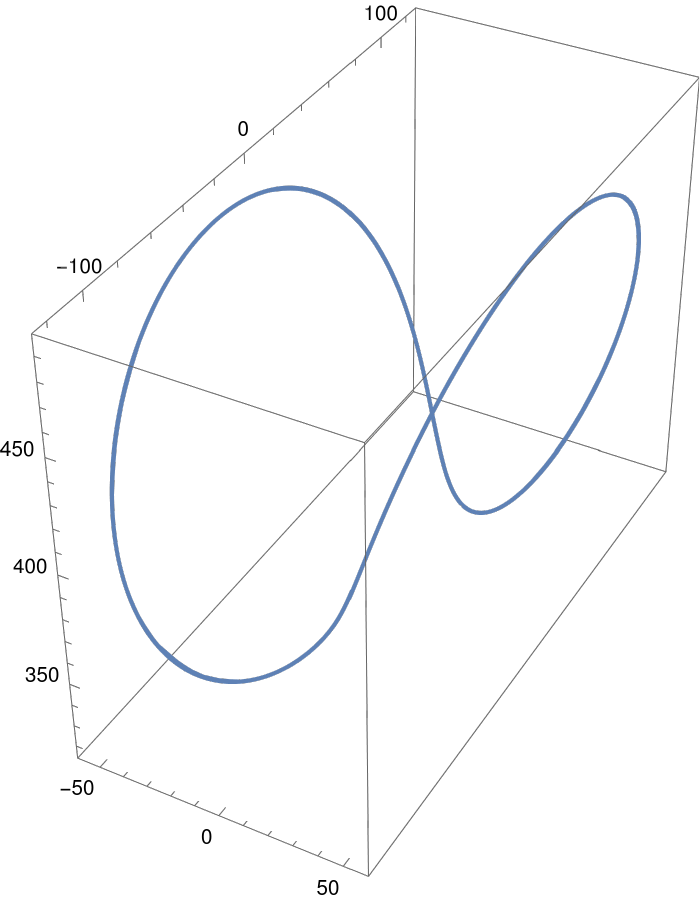}
 \caption{\sl Mathematica simulations of the trajectories
of the Lorenz system with parameters~$\sigma:= 10$, $\rho:=400$, $\beta:=\frac83$,
starting at the point~$(10,10,10)$ with~$t\in[0,100]$ (left) and~$t\in[90,100]$ (right).}\label{NEc11AahMA4srdtfyvhbj400nklm9a99345tRDItangeFI}
\end{figure}

To make these equations look more attractive, it is customary to perform some cosmetic modifications
by defining
\begin{eqnarray*}&&
\omega:=\frac{q^2}{(q^2+\pi^2)^3},\qquad\sigma:={\mathcal{P}},\qquad
\beta:=\frac{4\pi^2}{q^2+\pi^2},\qquad\rho:=\omega {\mathcal{R}},\\ &&
X(t):=\frac{\pi q}{\sqrt2(q^2+\pi^2)}\,{\mathcal{A}}\left(\frac{t}{q^2+\pi^2} \right),\qquad
Y(t):=\frac{\pi \omega}{\sqrt2}\,{\mathcal{B}}\left(\frac{t}{q^2+\pi^2} \right)\\ {\mbox{and}}\quad&&
Z(t):=-{\pi\omega}{\mathcal{C}}\left(\frac{t}{q^2+\pi^2} \right).
\end{eqnarray*}
With this, using the notation~$\tau:=\frac{t}{q^2+\pi^2}$,
one deduces from~\eqref{TSOSTTHET134-LORE-sinf49-34} that
\begin{eqnarray*}&&
\dot X(t)=\frac{\pi q}{\sqrt2(q^2+\pi^2)^2} \dot{\mathcal{A}}(\tau)
=-\frac{\pi{\mathcal{P}}q}{\sqrt2(q^2+\pi^2)^3}\big( (q^2+\pi^2)^2{\mathcal{A}}(\tau)-q{\mathcal{B}}(\tau)\big)\\&&\qquad\qquad
=\frac{\pi{\mathcal{P}}q^2}{\sqrt2(q^2+\pi^2)^3}
{\mathcal{B}}(\tau)
-\frac{\pi{\mathcal{P}}q}{\sqrt2(q^2+\pi^2)}\,{\mathcal{A}}(\tau)
=\frac{{\mathcal{P}}q^2}{\omega(q^2+\pi^2)^3}
Y(t)-{\mathcal{P}}X(t)\\&&\qquad\qquad=\sigma(Y(t)-X(t)),
\end{eqnarray*}
that
\begin{eqnarray*}&&
\dot Y(t)=\frac{\pi \omega}{\sqrt2(q^2+\pi^2)}\,\dot{\mathcal{B}}(\tau)=
\frac{\pi \omega}{\sqrt2(q^2+\pi^2)}
\Big({\mathcal{R}}q{\mathcal{A}}(\tau)-(q^2+\pi^2) {\mathcal{B}}(\tau)+\pi q{\mathcal{A}}(\tau){\mathcal{C}}(\tau)\Big)\\&&\qquad\qquad=
\frac{\pi \omega}{\sqrt2(q^2+\pi^2)}
\left(\frac{\sqrt2{\mathcal{R}(q^2+\pi^2)}}{\pi}\,X(t)
-\frac{\sqrt2(q^2+\pi^2)}{\pi\omega} \,Y(t)-\frac{\sqrt2(q^2+\pi^2)}{\pi\omega}\,X(t)\,Z(t)\right)\\&&\qquad\qquad\;=
\rho X(t)
-Y(t)- X(t)\,Z(t)
\end{eqnarray*}
and that
\begin{eqnarray*}&&
\dot Z(t)=-\frac{\pi\omega}{q^2+\pi^2}\,\dot{\mathcal{C}}(\tau)=
\frac{\pi\omega}{q^2+\pi^2}\left(
\frac{\pi q{\mathcal{A}}(\tau){\mathcal{B}}(\tau)}2+4\pi^2{\mathcal{C}}(\tau)\right)
\\&&\qquad\qquad=
\frac{\pi\omega}{q^2+\pi^2}\left(
\frac{q^2+\pi^2}{\pi\omega}\,X(t)\,Y(t)-\frac{4\pi}\omega\,Z(t)\right)=X(t)\,Y(t)-\beta Z(t).
\end{eqnarray*}

\begin{figure}
  \centering
  \includegraphics[height=2.8cm,width=0.9\textwidth]{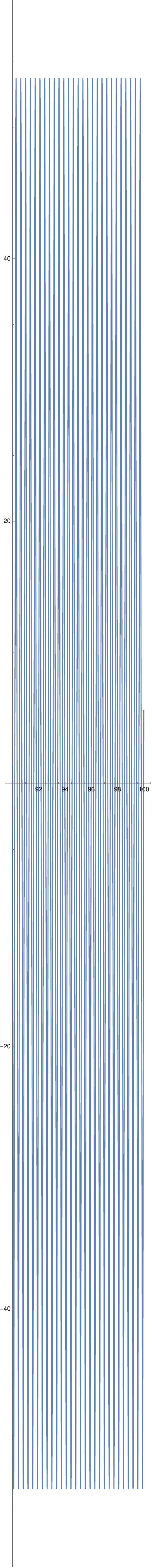}\\
  \includegraphics[height=2.8cm,width=0.9\textwidth]{LOR-Y}\\
  \includegraphics[height=2.8cm,width=0.9\textwidth]{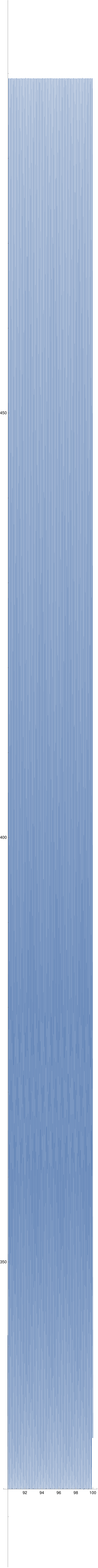}
  \caption{\sl Mathematica simulations of the trajectories $X$, $Y$ and~$Z$
of the Lorenz system with parameters~$\sigma:= 10$, $\rho:=400$, $\beta:=\frac83$,
starting at the point~$(10,10,10)$ with~$t\in[90,100]$.}\label{NEc11AahMA4srdtfyvhbj400nklm9a99345tRDItangeFI-2la2}
\end{figure}

Hence, we can write~\eqref{TSOSTTHET134-LORE-sinf49-34} in the simpler form
\begin{equation}\label{RTYLORR}\begin{dcases}
\dot X&=\sigma (Y-X),\\
\dot Y&=X(\rho -Z)-Y,\\
\dot Z&=XY-\beta Z.\end{dcases}
\end{equation}
The set of these three ordinary differential equations is known as the \index{Lorenz system}
Lorenz system. Notice that~$\rho$ is proportional to the Rayleigh number~${\mathcal{R}}$, hence
recalling~\eqref{JMS-TEBFDrey}, we may suspect that the specific value of~$\rho$
may play an important role in a possible bifurcation diagram for this system.

This is indeed the case, though a full understanding of the complexity of the Lorenz system
overcomes the present knowledge on the subject: as pointed out in~\cite[Section~15.1]{MR3293130},
``we are decades (if not centuries) away from
rigorously understanding all of the fascinating dynamical phenomena that
occur as the parameters change''.

A first glimpse on the importance of the values of the parameters in~\eqref{RTYLORR}
can be understood by studying the equilibria of the system. For example (see e.g.~\cite[Section~14.2]{MR3293130})
one can show that~\eqref{RTYLORR} always presents the origin as an equilibrium,
and two additional equilibria
$$Q_\pm := \left(\pm\sqrt{ \beta(\rho - 1)},\, \pm\sqrt{ \beta(\rho - 1)},\, \rho - 1\right)$$
when~$\rho>1$. When~$\rho<1$ the origin is a sink, but it becomes a saddle when~$\rho>1$. Assuming that~$\sigma>1+\beta$,
the equilibria in~$Q_\pm$ are sinks when~$\rho\in(1,\rho_\star)$, with
\begin{equation}\label{eq34567u4YHSa5tio3qwrt34n}\rho_\star:=\frac{\sigma(\sigma+\beta+3)}{\sigma-\beta-1},\end{equation}
but when~$\rho>\rho_\star$ the eigenvalues of the linearized system
at~$Q_\pm$ become purely imaginary and, very roughly speaking, the local dynamics tend to circulate around these points, with opposite velocities.

\begin{figure}
  \centering
  \includegraphics[height=.09\textheight]{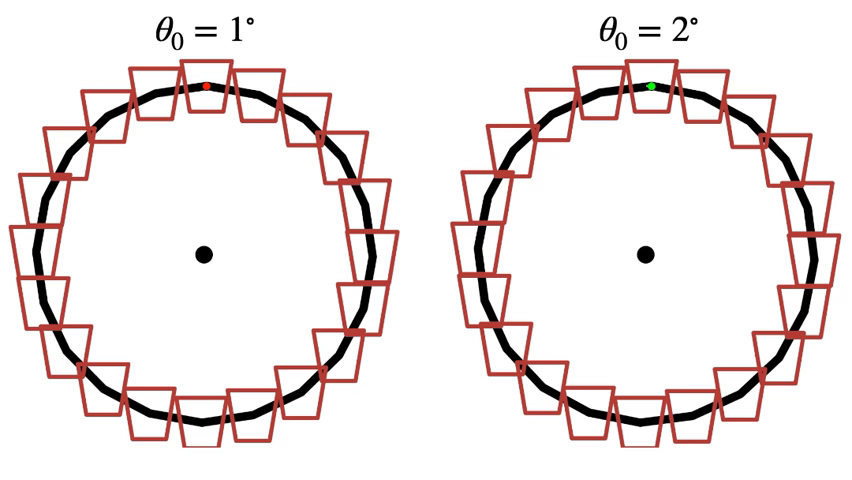}$\qquad\longmapsto\qquad$
 \includegraphics[height=.09\textheight]{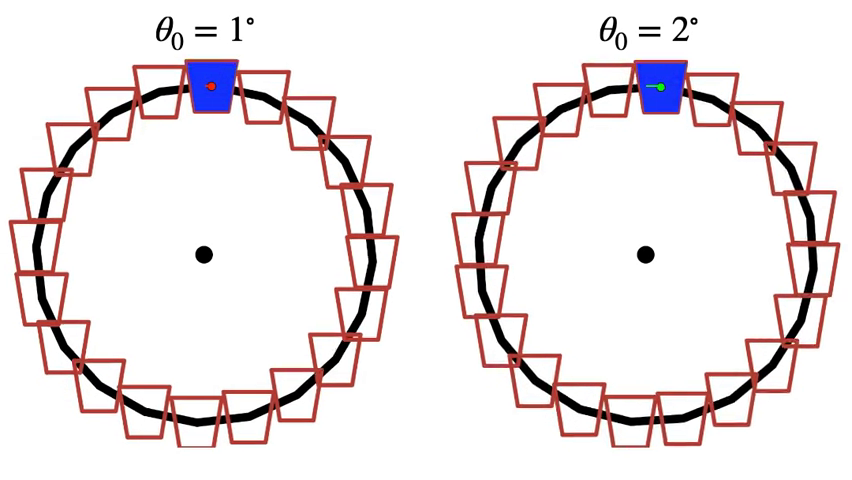}$\qquad\longmapsto\qquad$
 \includegraphics[height=.09\textheight]{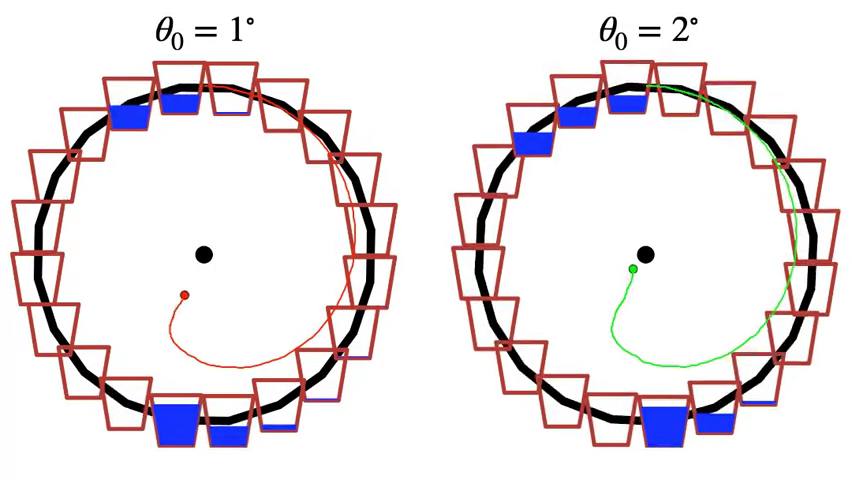}\\ $\,$
 \\
 \includegraphics[height=.09\textheight]{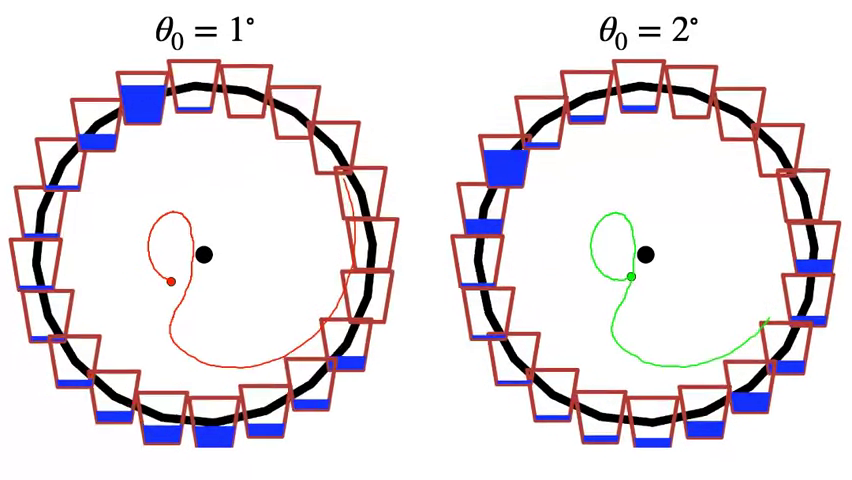}$\qquad\longmapsto\qquad$
 \includegraphics[height=.09\textheight]{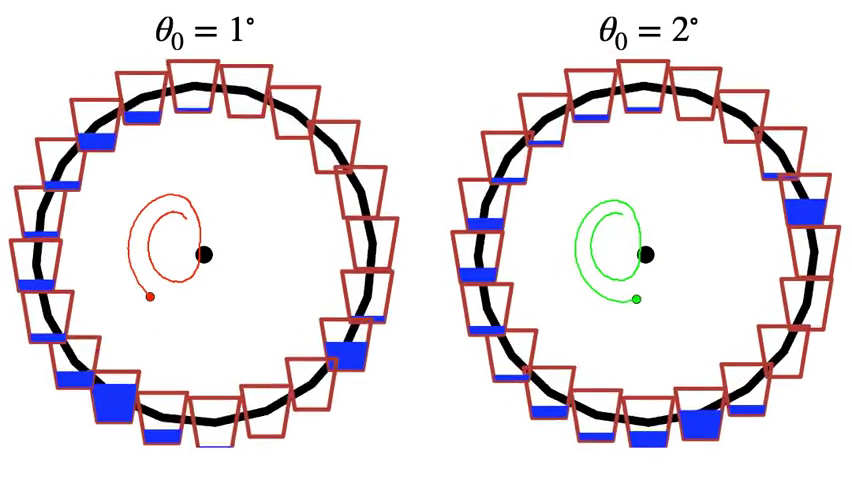}$\qquad\longmapsto\qquad$
 \includegraphics[height=.09\textheight]{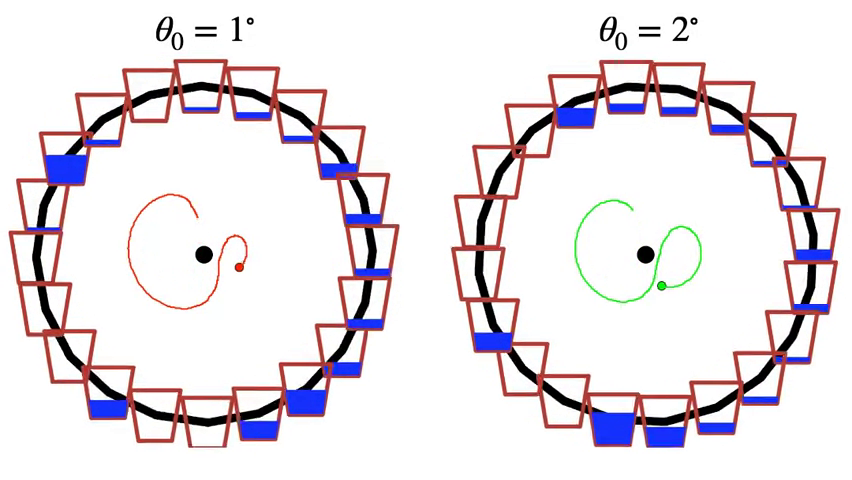}\\ $\,$
 \\
 \includegraphics[height=.09\textheight]{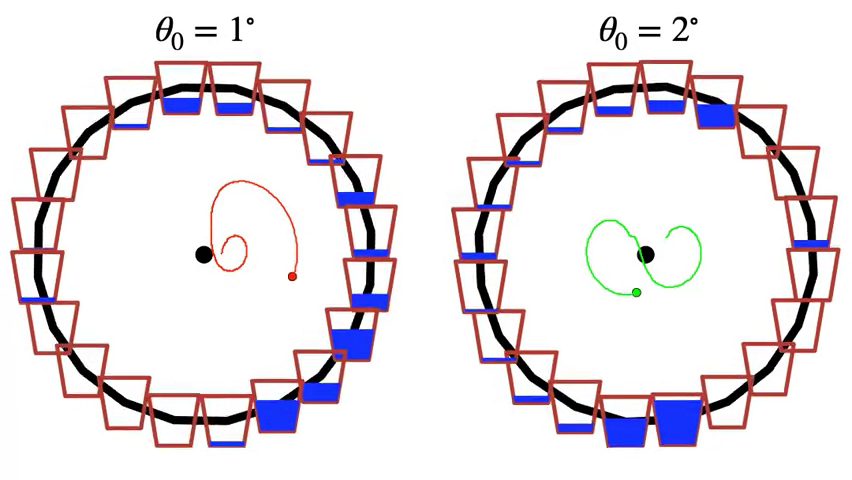}$\qquad\longmapsto\qquad$
 \includegraphics[height=.09\textheight]{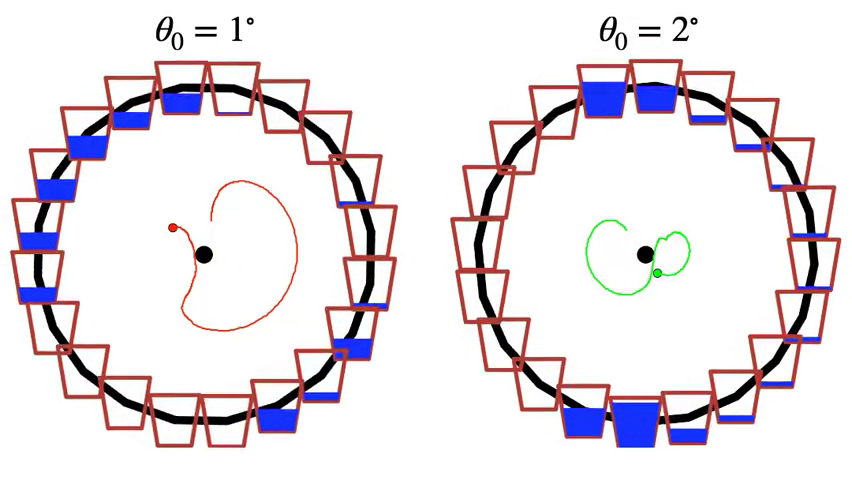}$\qquad\longmapsto\qquad$
 \includegraphics[height=.09\textheight]{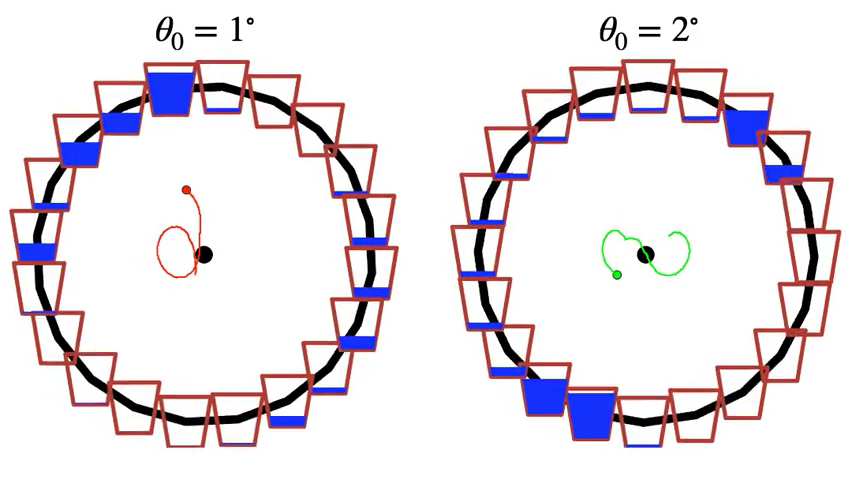}
 \caption{\sl Simulation of two waterwheels with initial angles differing by~$1$ degree
 (extract from a video by Aiyopasta from
        Wikipedia,
        licensed under the Creative Commons Attribution-Share Alike 4.0 International license).}\label{G24WHEE3NFID23ItangeFI243ORNUL}
\end{figure}

Also, the system shrinks volumes exponentially fast
and trajectories are confined within a bounded region (see e.g.~\cite[page~309]{MR3293130}),
thus suggesting that the system presents an \index{attractor} ``attractor''
(roughly speaking, a set to which trajectories
asymptotically approach in the course of dynamic evolution).
Actually, before the introduction of
the Lorenz system, the only types of stable attractors known in differential equations
were equilibria and closed orbits, while the Lorenz system exhibited a very
different structure:
namely, for certain values of the parameters, the Lorenz  system possesses what has come to
be known as \index{strange attractor} a ``strange attractor''
(see e.g.~\cite{MR1870856, MR2043795}), that is, roughly speaking,
an object of fractal dimension which is strongly sensitive to initial conditions.

Without going into
the difficult details of this matter, let us mention that, as a model case,
Lorenz considered~$\sigma:= 10$ and~$\beta:=\frac83$, for which~\eqref{eq34567u4YHSa5tio3qwrt34n}
produces~$\rho_\star=\frac{470}{19}=24.7368421\dots$;
see Figure~\ref{2HAFLORRE31MURA2-STAB} 
for a plot of the relatively regular dynamics occurring\footnote{See also
{\tt https://itp.uni-frankfurt.de/$\sim$gros/Vorlesungen/SO/simulation\_example/},
{\tt http://www.malinc.se/m/Lorenz.php}
and {\tt https://fusion809.github.io/Lorenz/}
for interactive animated simulations of the Lorenz system.}
when~$\rho:=10<\rho_\star$.
Instead, for larger values of~$\rho$ a number of complex phenomena arise,
see Figure~\ref{2HAFLORRE31MURA2}
for a sketch of the case~$\rho:=28>\rho_\star$, which was the one originally presented\footnote{Concerning the traditional
choice~$\sigma:= 10$, $\rho:=28$, $\beta:=\frac83$,
as stressed
in~\cite{SOK2:L}, ``there is 
nothing magical about these values; indeed, any other values within a fairly
wide range would produce qualitatively similar behavior'':
in~\cite{SOK2:L} it is also stressed that however some authors have tried
to infer from these rather arbitrary values
a ``positivity ratio'' which, with apparently ``no theoretical or empirical justification for
the use of differential equations drawn from fluid dynamics'',
should predict ``an individual's degree of
flourishing''. The authors introducing
this questionable positivity ratio accomplished however some kind of academic prestige,
since this ratio ``has had an extensive influence on the field
of positive psychology'', their work ``has been frequently cited'', and popular books
are ``devoted to expounding this
{\em huge discovery}, which has also been enthusiastically brought to a wider audience''
in dissemination books.

Well, sometimes producing an easy and user-friendly recipe
to cope with complex and very topical problems does pay off
(but sometimes it doesn't).}
by Lorenz: notice in this picture the very complicated
structure of the trajectory, which seems repeatedly to approach and then spiral
away from one of the two equilibria~$Q_\pm$ to move close to the other, then
escaping again towards the previous one, jumping back and forward on and on,
producing dense leaves which are interwoven in a rather tangled way.

Additionally, the resulting attractor seems to be very sensitive\footnote{Lorenz himself was astonished by the discovery of the high sensitivity
on the initial data. In his own words (see
{\tt https://eapsweb.mit.edu/sites/default/files/Scientist\_by\_Choice.pdf}):
``At first I suspected trouble with the computer, which occurred fairly often, but,
when I compared the new solution step by step with the older one, I found that at
first the solutions were the same, and then they would differ by one unit in the last
decimal place, and then the differences would become larger and larger''.}
to initial conditions, showing
a great deal of unpredictability which is typically considered one of the main ingredients of ``chaos''.
To understand this sensitivity, we can look at Figure~\ref{NEc11AahMA0rtXcdtt345pdkfDIRDItangeFI},
which plots the $X$-trajectories
of the Lorenz system
starting at the points~$(10,10,10)$ and~$\left(10+\frac1{100},10,10\right)$:
while these two trajectories are almost
indistinguishable for small times (e.g., about~$t=18$, in the left images)
they become completely different and essentially unrelated as time flows (e.g.,
for~$t\in[18,50]$ in the images on the right).

Notice however that while
a tiny change in the initial position may result in
drastic changes in the eventual behavior of the orbits, the trajectories end up approaching
the same attractor, just traveled at different instants of time: this
is shown for instance in Figure~\ref{NEc11AahMA0rtXcdtt34LOREatt5pdkfDIRDItangeFI},
in which the global geometric pattern described by the journey of the two trajectories
happen to be quite similar.

Moreover, the Lorenz system
seems to show alternating patterns of chaos and periodic motion, depending
on the values of~$\rho$ (see e.g.~\cite[Chapter~4]{MR681294}),
with ``windows'' of the values of~$\rho$ allowing for periodic behaviors,
see Figure~\ref{NEc11AahMA4srdtfyvhbjnklm9a99345tRDItangeFI}
for an example of this situation.

Furthermore, for very large values of~$\rho$, perhaps quite surprisingly,
the dynamics of the Lorenz system simplifies again, with globally
attracting limit cycles (see~\cite[page~109]{MR2310642}
and the references therein). This phenomenon is showcased here in Figures~\ref{NEc11AahMA4srdtfyvhbj400nklm9a99345tRDItangeFI}
and~\ref{NEc11AahMA4srdtfyvhbj400nklm9a99345tRDItangeFI-2la2},
in which the periodic limit trajectory appears quite distinctly.

We stress however that these claims are at the moment
mostly based only on numerical experiments
and we are still in need of complete analytic proofs. \medskip

Interestingly, there is also a simple and practical mechanical model of a particular
case of the Lorenz system,
which was realized by Willem Malkus in the 1960s
and consists of a \index{chaotic waterwheel} chaotic waterwheel.
The model consists indeed of a toy waterwheel with
cups. A water source at the top pours
a constant flow in at the top bucket, making the wheel turn.
At low flow rates, the wheel turns in one direction or the other, but
at high flow rates the wheel spins one way, and then the other,
since water filled cups are heavier, oppose the spinning of the wheel and
make it turn the other way. Furthermore, the waterwheel exhibits a sensitive dependence
on the initial data, see e.g. Figure~\ref{G24WHEE3NFID23ItangeFI243ORNUL}
in which the initial angles of the wheels differ by only~$1$ degree, but the evolution of the plotted centers of mass, as well as the disposition of water in the buckets, diverge significantly in time.
The unpredictable manner in which spins reverse orientation in this
model is precisely a manifestation of the chaotic patterns exhibited by
the Lorenz system, see~\cite{LESLIE}, \cite[Section~9.1]{MR3837141}
and the references therein for a detailed mathematical discussion of the
chaotic waterwheel and its link to the Lorenz system.\medskip

All in all, the analysis pioneered by Lord Rayleigh
and Edward Lorenz was not only important towards a partial understanding of
a rather simplified model for convection; instead, it truly opened up new
horizons, showcasing a number of novel and breathtaking phenomena which were later discovered
virtually in any area of science.\medskip

See~\cite{MR681294, MR3293130, MR3837141} for a thorough introduction to
the dynamics of the Lorenz system.
For additional information about the close relationship between convection models in fluid and the Lorenz system
see e.g.~\cite{MR882723, MR1263025, MR2178136, MR4283024} and the references therein.

\subsection{Who wants to be a millionaire?}

One way to become a millionaire is to solve a Millennium Prize Problem (see pages~\pageref{106106}
and~\pageref{106106-bisPE}), but this is probably by far the most difficult way to get rich.
Also, one should take into consideration the possibility that mathematics may not make us rich,
but it makes us happy, which is inestimable.

\begin{figure}
  \centering
  \includegraphics[width=.55\linewidth]{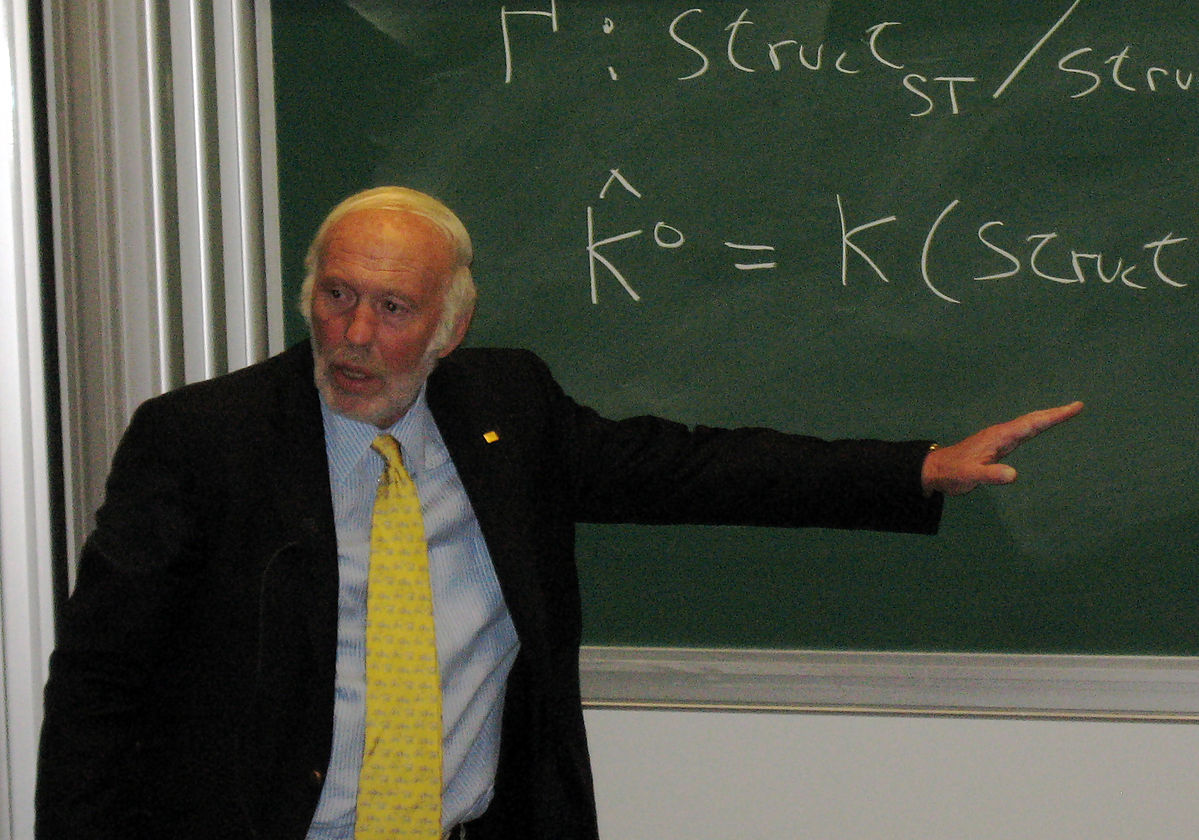}
 \caption{\sl James Harris Simons in 2007
 (image from
        Wikipedia,
        licensed under the 
 Creative Commons Attribution-Share Alike 2.0 Germany license).}\label{SIMoGREENFIDItangeFI}
\end{figure}

There are however truly outstanding mathematicians, such as Jim Simons (see Figure~\ref{SIMoGREENFIDItangeFI})
who, after having revolutionized the theory of minimal surfaces and the topological quantum field theory,
and also after having broken codes for the NSA
during the Cold War,
decided to improve his portfolio:
allegedly, Simons' net worth is estimated to be~$25.2\times 10^9$ US \$,
making him the 66th-richest person in the world, and likely the richest\footnote{Simons has not forgotten
his love for mathematics. Not only because he owns a motor yacht named Archimedes,
but also because he
co-founded the Simons Foundation,
supporting projects related to scientific research, education and health, and
he sponsored the Simons Institute for the Theory of Computing
and the Mathematical Sciences Research Institute at Berkeley.

Here is an interesting story by the way. While Simons was working for the Institute for Defense Analyses by cracking Russian codes, his boss General Maxwell Taylor wrote an article in the New York Times Magazine strongly supporting the Vietnam War. Then, Simons published a counter-editorial in the Times, claiming that not everyone who worked for Taylor subscribed to his views. This resulted in an interview of Simons from Newsweek. When Simons told his superiors about this interview, they fired him right away.

Allegedly, Simons said ``Getting fired once can be a good experience. You just don't want to make a habit of it''.}
among the mathematicians.
His financial success is related to the foundation of
a quantitative hedge fund (see below for the notion of hedging)
which trades using quantitative mathematical models, so mathematics can be sometimes useful after all.

One of the chief uses of partial differential equations
in mathematical finance is related to the Black-Scholes equation\index{Black-Scholes equation}
\begin{equation}\label{B:LA:CK}
\frac{\partial V}{\partial t}+{\frac{1}{2}}\sigma ^{2}S^{2}
\frac{\partial ^{2}V}{\partial S^2}=rV-rS\frac{\partial V}{\partial S},\end{equation}
named after Fischer Black and Myron Scholes.

The meaning of this equation is the following.
Suppose that we have a stock\index{stock} (that is, a security that represents the ownership of a fraction of a corporation)
whose price is denoted by~$S$.
The price of a stock fluctuates based on supply and demand
(if more people want to buy than sell it, the price will rise) hence we take~$S=S(t)$ to be a function of time.

To model the fluctuations of~$S$ in the simplest possible way, we can make two assumptions.

On the one hand, we can assume that the stock refers to some company
with revenues, earnings, dividends, etc., and we can estimate the growth rate
of the company and of the corresponding stocks: for this, we assume that the variation of~$S$
in the unit of time is proportional, by some factor~$\mu$, to the value of~$S$
(the more stocks one possesses, and the higher the coefficient~$\mu$, the more one receives from the growth of the company's revenues; actually positive values of~$\mu$ correspond to a growth
and negative ones to a degrowth).

On the other hand, the market presents a number of uncertainties which are almost
impossible to fully take into account, such as economic crises, political developments,
real estate bubbles, technological innovations, pandemics, landing of aliens from outer space, etc.
If we have no specific information on the matter, we can just assume that this volatility
is modeled by a Brownian motion (pretty much as the random walk presented on page~\pageref{RANDOW}).
This is accounted for by a process~$W(t)$ that ``randomly wiggles up and down'' and
that quantifies the source of uncertainty in the price. This term is modulated by a volatility coefficient~$\sigma\in[0,+\infty)$ (the higher this coefficient, the bigger the impact of the uncertainties in the price oscillation).

\begin{figure}
  \centering
  \includegraphics[width=.65\linewidth]{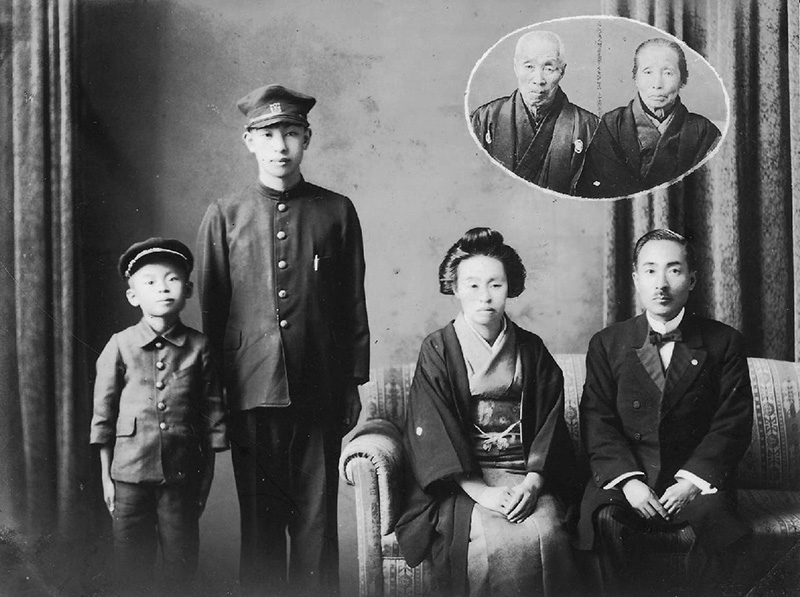}
 \caption{\sl The It\^{o} family (Public Domain image from
 Wikipedia).}\label{ITONFIDItangeFI}
\end{figure}

Assuming that these two effects contribute somewhat consistently to the variation in time of the stock,
we therefore write that\footnote{Equation~\eqref{3456yuioiuhgfdr5678iokmnbvfdr5678ijhgfde45678ijhgfdety3hfe2}
is sometimes called the ``geometric Brownian motion''\index{geometric Brownian motion}
equation. It is a very popular tool in quantitative finance since it relies
on the simple, but often reasonable, assumption that the relative variation~$\frac{dS}S$
of a value follows a diffusive process.}
\begin{equation}\label{3456yuioiuhgfdr5678iokmnbvfdr5678ijhgfde45678ijhgfdety3hfe2}
dS=\mu S\,dt+\sigma S\, dW.\end{equation}
Strictly speaking, this is not really a differential equation,
since we cannot ``divide both terms by~$dt$'', because Brownian motions
are ``not differentiable'': see e.g.~\cite[Section~3.4.2]{MR3154922}
for a precise statement and for an elegant introduction to stochastic calculus.
Indeed, \eqref{3456yuioiuhgfdr5678iokmnbvfdr5678ijhgfde45678ijhgfdety3hfe2}
is a stochastic differential equation\index{stochastic differential equation},
which is a type of equations that we do not treat here,
dealing with it just at an intuitive level. Basically, the only bit of information that
we need from stochastic calculus is the so-called\footnote{Kiyosi It\^{o}
was the pioneer of the stochastic calculus that brings his name.
See Figure~\ref{ITONFIDItangeFI}
for a picture of Kiyosi It\^{o}
at age~22 (Kiyosi is the second from the left; the first from the left is his brother
Seiz\^{o}, who also became a mathematician).

The heuristic explanation of It\^o's Chain Rule~\eqref{ITOPIOTNNRRRESIDEKKA}
is that, formally, ``$(dW)^2=dt$'', in agreement with the fact that, in the random walk,
the unit of time depends ``quadratically'' on the unit of space, according to~\eqref{SCAKLI-UJs903th}.
This (by keeping linear orders in~$dt$, and up to quadratic in~$dW$) leads to
\begin{eqnarray*}&&
H(X+dX,t+dt)-H(X,t)\\&=& \partial_X H\,dX+\partial_t H\,dt+\frac12 \partial_X^2 H\,(dX)^2+\frac12 \partial_t^2 H\,(dt)^2
\\&&\qquad+\partial_{Xt}^2 H\,dX\,dt+o((dX)^2+(dt)^2) \\
&=& \partial_X H\,(F\,dt+G\,dW)+\partial_t H\,dt+\frac12 \partial_X^2 H\,(F\,dt+G\,dW)^2
+\frac12 \partial_t^2 H\,(dt)^2\\&&\qquad+\partial_{Xt}^2 H\,(F\,dt+G\,dW)\,dt+o(dt\,dW+(dW)^2+(dt)^2)
\\&=&\partial_X H\,(F\,dt+G\,dW)+\partial_t H\,dt+\frac{G^2}2 \partial_X^2 H\,(dW)^2
+o(dt+(dW)^2)
\\&=&
\left( \partial_t H +F \partial_X H+\frac{G^2}{2} \partial_X^2 H\right)\,dt+G\partial_X H\,dW+o(dt),\end{eqnarray*}
that gives It\^o's Chain Rule~\eqref{ITOPIOTNNRRRESIDEKKA}.

See e.g.~\cite{MR3154922} and the references therein for a more exhaustive treatment of It\^o's Chain Rule and a more
accurate approach to its proof.}
It\^o's Chain Rule\index{It\^o's Chain Rule} (see e.g. Sections~1.3,  4.3, 4.4 and in particular the theorem on page~80 of~\cite{MR3154922}) according to which if~$dX=F\,dt+G\,dW$,
then the random variable~$H(X(t),t)$ satisfies
\begin{equation} \label{ITOPIOTNNRRRESIDEKKA}
dH(X,t)=\left( \partial_t H +F \partial_X H+\frac{G^2}{2} \partial_X^2 H\right)\,dt+G\partial_X H\,dW.\end{equation}
Let us now come back to our financial setting.
Suppose that we perceive the buying of a given stock as too risky
and we prefer instead to buy only ``the possibility of buying the stock
in the future'', if this turns out to be convenient.
This type of contract is called ``option''\index{option}
and conveys its owner
the right, but not the obligation, to buy the stock at a specified strike price either
prior to or on a specified date.
For simplicity, let us consider the case in which the option allows, but not obliges, the holder
to buy the stock for a specified price~$K$ (called ``strike price''\index{strike price}
in jargon)
at a specific time~$T$, that is on the date of option maturity:
these options are called in jargon ``European options''\index{European option};
the options which instead 
can be exercised any time up to and including the date of expiration
are called ``American options''\index{American option}
(perhaps the name originated
by the alleged fact that one type of option was commonly
traded in London and the other in New York).

We denote by~$V$ the price of this option
and we suppose that~$V$ depends on time and on the value~$S$ of the underlying stock, that is~$V=V(t,S)$, for~$t$ less than the maturity time~$T$.
For instance, at the maturity date the value of this option is either~$S(T)-K$,
namely the difference between the actual value of the stock~$S$ at time~$T$
and the strike price~$K$ that was agreed at the beginning, if~$S(T)-K$ is positive
(in which case the option turned out to be a convenient investment),
or null, if~$S(T)-K$ is negative (in this case the value of the stock is lower than the strike price,
so no need to use the option to pay more, if one is interested in the stock, they can just forget
about the option and buy the stock at the market price). More specifically, these observations yield that\footnote{For American options, one also has that
$$ V(t, S(t))\ge\max\{S(t)-K,0\} $$
for all~$ t\in[0,T]$,
the inequality coming from the fact that one can, but is not obliged to, utilize the option at time~$t$
and this provides an ``obstacle problem''\index{obstacle problem}
for equation~\eqref{B:LA:CK}, since it forces solutions of~\eqref{B:LA:CK}
to stay above a constraint
(or to drop the equation when they meet the constraint;
this situation is significantly more complicated than that of a single equation, and this is the reason for which
we limited ourselves here to the case of European options).

With respect to~\eqref{TERM01-2pwf}, we may consider it as a ``terminal condition''\index{terminal condition}
for equation~\eqref{B:LA:CK}.

Interestingly, with respect to the classification discussed in footnote~\ref{CLASSIFICATIONFOOTN} on page~\pageref{CLASSIFICATIONFOOTN},
we have that equation~\eqref{B:LA:CK} is of parabolic type, but the sign in front of the time derivative
is opposite to the case of the heat equation (compare with~\eqref{DAGB-ADkrVoiweLL4re2346ytmngrrUj7};
a closer relation \label{BSDAGB-ADkrVoiweLL4re2346ytmngrrUj7}
between the Black-Scholes equation and the heat equation will be discussed on page~\pageref{BSDAGB-ADkrVoiweLL4re2346ytmngrrUj72}).
In this sense, equation~\eqref{B:LA:CK} shares some similarity with heat type equations,
but going backwards in time. This makes sense, since while the heat equation
smooths out differences as time flows, the value of an option becomes quite rigid close to its maturity
date.}
\begin{equation}\label{TERM01-2pwf} V(T,S)=\max\{S-K,0\}.\end{equation}

Our objective is thus to understand the time evolution of the value~$V$ of the option:
this is important, since if one possesses an option, may decide to sell it at some time~$t$,
knowing just the strike price~$K$ agreed at the beginning and the present value of the stock~$S(t)$,
hence it is relevant to establish a fair price~$V(t,S(t))$ for the option
under this information (since one would like to answer this question for all possible values
of the stocks, this reduces to understand~$V(t,S)$ just knowing the strike price~$K$,
with this one would substitute~$S=S(t)$ to obtain the fair price of the option).

To model the time evolution of~$V$, and thus provide a convincing derivation of
the Black-Scholes equation in~\eqref{B:LA:CK}, we 
use a method called in finance ``delta hedging''\index{delta hedging},
namely we adopt a trading strategy that reduces, or ``hedges'', the risk related to the volatility
of the process.
For this, we suppose to buy a quantity~$\alpha=\alpha(t,S(t))$ of options and some quantity~$\beta=\beta(t,S(t))$ of the stocks
(here, $\alpha$ and~$\beta$ are in the reals, negative numbers formally corresponding
to options or stock that are actually sold).
The goal is to choose~$\alpha$ and~$\beta $ to fully cover the risk
(or, at least, to cover the risk to the best of our possibilities). More precisely, we consider a portfolio~$P:=\alpha V+\beta S$
and we compare this investment with another one that carries zero risk,
that is some operation which ensure a certain (or, more realistically, almost certain)
future return and (virtually) no possibility of loss: for instance a treasury bill, or treasury bond,
namely a bond in which the face value is repaid at the time of maturity with some interest.
We denote by~$r\ge0$ the risk-free interest rate (this will be precisely the parameter~$r$
appearing in~\eqref{B:LA:CK}). The value of this zero risk portfolio~$P_0$ is thus
described by~${dP_0}=r P_0\,dt$.

Since, via a suitable choice of~$\alpha$ and~$\beta$, the portfolio~$P$ is set up to have zero risk, it must be equivalent
to the portfolio~$P_0$, otherwise there would be a better choice to gain money at zero risk! This says that
also~$ dP=rP\,dt$ and therefore
\begin{equation}\label{93456yuioiuhgfdr5678iokmnbvfdr5678ijhgfde45678ijhgfdety3hfe3CP8}
\big(r\alpha V+r\beta S\big)\,dt=
rP\,dt =dP.\end{equation}
We stress that the portfolio~$P$ is assumed to be self-financing\index{self-financing portfolio},
namely there are no inflows or outflows of money. More specifically, at every instant of time,
the purchase of a new option or stock must be financed by the sale of an old one: 
hence, since there is no external infusion or withdrawal of money,
the variation of the value of the portfolio corresponds precisely to the variation of the values
of the options and stocks possessed, that is
\begin{equation}\label{93456yuioiuhgfdr5678iokmnbvfdr5678ijhgfde45678ijhgfdety3hfe3CP89} dP=\alpha\, dV+\beta \,dS.\end{equation}
We stress that this does not mean necessarily that~$\alpha$ and~$\beta$ are constants,
just that the market value of the portfolio~$P$ at a given time
equals the purchase value of the new portfolio.
Namely, if in the infinitesimal interval of time~$dt$ the value of the options has increased respectively by~$dV$
and the value of the stocks by~$dS$, then the gain obtained would correspond to the number of options
possessed times the options' increment value (that is~$\alpha$ times~$dV$) plus the number of stocks
possessed times the stocks' increment value (that is~$\beta$ times~$dS$): this gives
that the total gain in the infinitesimal interval of time~$dt$ is~$\alpha\, dV+\beta \,dS$, which is precisely
the right hand side in~\eqref{93456yuioiuhgfdr5678iokmnbvfdr5678ijhgfde45678ijhgfdety3hfe3CP89}.
Accordingly, the balance prescribed by~\eqref{93456yuioiuhgfdr5678iokmnbvfdr5678ijhgfde45678ijhgfdety3hfe3CP89}
states that this gain is used precisely to reinvest to enlarge the portfolio itself
(this, in the optimistic scenario that there is a gain; if there is a loss then the right hand side of~\eqref{93456yuioiuhgfdr5678iokmnbvfdr5678ijhgfde45678ijhgfdety3hfe3CP89} is negative and
the portfolio needs to be correspondingly shrunk).

Now we let~$\varpi$ be the ratio\footnote{In most of the literature,
this ratio (actually, minus this ratio) is denoted by either~$\Delta$ or~$\delta$, whence the name
of delta hedging. Here, for typographical convenience, and to avoid confusion with
the Laplace operator and with increments and derivatives, we preferred to call it~$\varpi$.} between~$\beta$ and~$\alpha$, namely we take~$\varpi=\varpi(t,S(t))$ such that~$\beta=\alpha\varpi$.
Then~\eqref{93456yuioiuhgfdr5678iokmnbvfdr5678ijhgfde45678ijhgfdety3hfe3CP8}
can be stated in the form
\begin{equation*}
\big(r\alpha V+r\alpha\varpi S\big)\,dt=
\big(r\alpha V+r\beta S\big)\,dt=dP
=\alpha\,dV+\beta\,dS
=\alpha\,dV+\alpha\varpi\,dS
\end{equation*}
and therefore (assuming for simplicity~$\alpha\ne0$)
\begin{equation}\label{93456yuioiuhgfdr5678iokmnbvfdr5678ijhgfde45678ijhgfdety3hfe38}
\big(r V+r\varpi S\big)\,dt
=dV+\varpi\,dS.
\end{equation}
We also point out that
\begin{equation*}
dV=\left( \partial_t V +\mu S \partial_S V+\frac{\sigma^2 S^2}{2} \partial_S^2 V\right)\,dt+\sigma S \partial_S V\,dW
,\end{equation*}
thanks to~\eqref{ITOPIOTNNRRRESIDEKKA}.

{F}rom this, \eqref{3456yuioiuhgfdr5678iokmnbvfdr5678ijhgfde45678ijhgfdety3hfe2} and~\eqref{93456yuioiuhgfdr5678iokmnbvfdr5678ijhgfde45678ijhgfdety3hfe38} we infer that
\begin{equation}\label{0olLur2345609iuygfadTbSan2sDMrote}
\begin{split}
0&=dV+\varpi\,dS-
\big(rV+r\varpi S\big)\,dt\\
&=
\left( \partial_t V +\mu S\partial_S V+\frac{\sigma^2 S^2}{2} \partial_S^2 V\right)\,dt+\sigma S\partial_S V\,dW
+\varpi\,dS-
\big(rV+r\varpi S\big)\,dt\\&=
\left( \partial_t V +\mu S\partial_S V+\frac{\sigma^2 S^2}{2} \partial_S^2 V
-rV-r\varpi S
\right)\,dt+\sigma S\partial_S V\,dW
+\varpi\,\big(\mu S\,dt+\sigma S\, dW\big)
\\&=
\left( \partial_t V +\mu S\big(\partial_S V+\varpi\big)+\frac{\sigma^2 S^2}{2} \partial_S^2 V
-rV-r\varpi S
\right)\,dt+ \sigma S\big(\partial_S V+\varpi \big)\,dW.\end{split}
\end{equation}
Consequently, since~$\varpi$ is chosen to make the portfolio free from the risk caused by the stochastic term
this gives that one must choose~$\varpi$ such that~$\partial_SV+\varpi=0$
(note that this cancels in~\eqref{0olLur2345609iuygfadTbSan2sDMrote}
the term involving~$W$ and the linear term in~$\sigma$). That is, we take
\begin{equation}\label{KSPDJ0ow98yihg94hgOHdpfoikgrtohyjrp-you45jhn945955-104}
\varpi=-{\partial_SV}
\end{equation}
and we thus plug\footnote{We observe that in financial models, it is customary to assume that~$\partial_SV\in[0,1]$.
Indeed, if the price of a stock increases, we may expect that the price of the corresponding option also increases
(since it would allow us the possibility of buying a valuable stock at a prescribed price), which implies that~$\partial_S V\ge0$.
Moreover, one expect the price of an option to increase slower than the price of a stock, since the change of value
of the option is somehow the consequence subsequent to the change of value of the stock, which suggests that~$\partial_SV\le1$.

In this spirit, the prescription in~\eqref{KSPDJ0ow98yihg94hgOHdpfoikgrtohyjrp-you45jhn945955-104} reads~$\frac\beta\alpha=\varpi\le0$, giving that, in our notation, $\alpha$ and~$\beta$ have opposite sign.
This is in agreement with the intuition that, say, the purchase of stocks (corresponding to a positive~$\beta$)
corresponds to a sell of options (corresponding to a negative~$\alpha$).}
this information into~\eqref{0olLur2345609iuygfadTbSan2sDMrote},
concluding that
\begin{equation*}
0=\left( \partial_t V +\frac{\sigma^2 S^2}{2} \partial_S^2 V
-rV+rS{\partial_SV}\right)\,dt
,\end{equation*}
which leads to the the Black-Scholes equation in~\eqref{B:LA:CK}, as desired.\medskip

It is also interesting to revisit the Black-Scholes equation in~\eqref{B:LA:CK}
in view of the financial model that it describes. We note in particular that the
right hand side of~\eqref{B:LA:CK} is described only in terms
of the risk-free interest rate~$r$: this is therefore the risk-free part of the investment
and it consists of the superposition of a long term strategy dictated by the option~$V$ (which will become effective only at its maturity time) and a short term one embodied by the stock~$S$ (which fluctuates
in a market subject to randomness).

On the left hand side of~\eqref{B:LA:CK} instead we see two terms. The first is the time derivative
of the option. The second term on the left hand side of~\eqref{B:LA:CK}
is more ``geometric'' since it reflects the convexity properties of the dependence of the option
on the underlying stock.\medskip

It is also useful to relate the Black-Scholes equation in~\eqref{B:LA:CK} to \label{BSDAGB-ADkrVoiweLL4re2346ytmngrrUj72}
the classical heat equation. As already observed on page~\pageref{BSDAGB-ADkrVoiweLL4re2346ytmngrrUj7},
to make this connection one needs to revert the arrow of time, hence it is convenient to look at a new time~$\tau:=T-t$ which transforms the terminal condition~\eqref{TERM01-2pwf} into an initial condition.
Furthermore, it comes in handy to introduce the nonlinear transformation~$x=\ln \left({\frac{S}{K}}\right)+\left(r-{\frac{\sigma^{2}}{2}}\right)\tau$ in order to simplify the right hand side of~\eqref{B:LA:CK}. More specifically, one sets
$$ u(x,\tau):=e^{r\tau} \,V\left(T-\tau,Ke^{x-\left(r-\frac{\sigma^2}2\right)\tau}\right)$$
and observes that if~$V$ solves the Black-Scholes equation in~\eqref{B:LA:CK} then,
using the notation~$S:=Ke^{x-\left(r-\frac{\sigma^2}2\right)\tau}$ and~$t:=T-\tau$,
\begin{eqnarray*}&&
e^{-r\tau}\left[
\frac{\sigma^2}2 \partial_{xx} u(x,\tau)-\partial_\tau u(x,\tau)\right]\\&=&
\frac{\sigma^2}2 \partial_{x} 
\left[ Ke^{x-\left(r-\frac{\sigma^2}2\right)\tau}\partial_S
V\left(T-\tau,Ke^{x-\left(r-\frac{\sigma^2}2\right)\tau}\right)\right]
-r  V\left(T-\tau,Ke^{x-\left(r-\frac{\sigma^2}2\right)\tau}\right)\\&&\qquad
+Ke^{x-\left(r-\frac{\sigma^2}2\right)\tau} \left(r-\frac{\sigma^2}2\right)\partial_S
V\left(T-\tau,Ke^{x-\left(r-\frac{\sigma^2}2\right)\tau}\right)\\&&\qquad
+\partial_t V\left(T-\tau,Ke^{x-\left(r-\frac{\sigma^2}2\right)\tau}\right)
\\&=&
\frac{\sigma^2}2 
\left[ 
Ke^{x-\left(r-\frac{\sigma^2}2\right)\tau}\partial_S
V\left(T-\tau,Ke^{x-\left(r-\frac{\sigma^2}2\right)\tau}\right)
+
\left(Ke^{x-\left(r-\frac{\sigma^2}2\right)\tau}\right)^2\partial_{SS}
V\left(T-\tau,Ke^{x-\left(r-\frac{\sigma^2}2\right)\tau}\right)
\right]\\&&\qquad
-r  V\left(T-\tau,Ke^{x-\left(r-\frac{\sigma^2}2\right)\tau}\right)
+Ke^{x-\left(r-\frac{\sigma^2}2\right)\tau} \left(r-\frac{\sigma^2}2\right)\partial_S
V\left(T-\tau,Ke^{x-\left(r-\frac{\sigma^2}2\right)\tau}\right)\\&&\qquad
+\partial_t V\left(T-\tau,Ke^{x-\left(r-\frac{\sigma^2}2\right)\tau}\right)\\
&=&\frac{\sigma^2}2 
\left[ S\partial_S
V(t,S)+S^2\partial_{SS}
V(t,S)
\right]-r  V(t,S)+S \left(r-\frac{\sigma^2}2\right)\partial_S
V(t,S)+\partial_t V(t,S)\\
&=&\frac{\sigma^2 S^2}2\partial_{SS}V(t,S)-r  V(t,S)+rS\partial_S
V(t,S)+\partial_t V(t,S)\\&=&0.
\end{eqnarray*}
Also,
$$ u(x,0)=V\left(T,Ke^{x}\right).$$
Consequently, the Black-Scholes equation in~\eqref{B:LA:CK} and the terminal condition~\eqref{TERM01-2pwf} 
produce the following heat diffusion problem for time~$\tau>0$ with initial datum:
\begin{equation}\label{CALOCA}
\begin{dcases}
&\partial_\tau u(x,\tau)=\displaystyle\frac{\sigma^2}2 \partial_{xx} u(x,\tau),\\
&u(x,0)=K(e^x-1)\,\chi_{(0,+\infty)}(x),
\end{dcases}
\end{equation}
where, as usual, we used the notation for the characteristic function of a set, namely
$$ \chi_A(x):=\begin{dcases}1 & {\mbox{ if }}x\in A,\\
0 &{\mbox{ otherwise.}}\end{dcases}$$

\begin{figure}
  \centering
    \includegraphics[width=.45\linewidth]{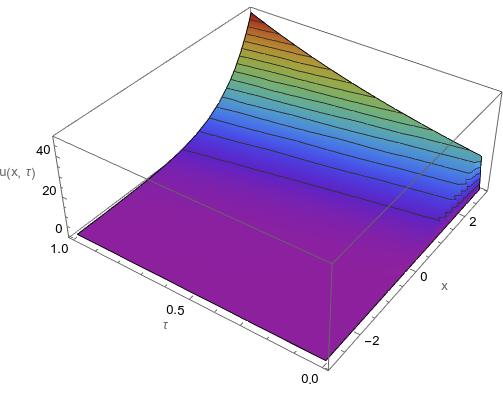}
 \caption{\sl Plot of a solution of~\eqref{CALOCA}
 (with~$K:=1$ and~$\sigma:=\sqrt2$).}\label{CALDRHSo0mNEdvspp-02iX32h90-SPIKERTERRAeFI-b}
\end{figure}

The plot of a solution of~\eqref{CALOCA} is sketched in Figure~\ref{CALDRHSo0mNEdvspp-02iX32h90-SPIKERTERRAeFI-b}.
The advantage of this formulation is that one can focus on the solution of the classical heat equation in~\eqref{CALOCA}
and deduce useful information on the solution of the Black-Scholes equation in~\eqref{B:LA:CK}
by transforming back to the original variables.

See e.g.~\cite{MR2200584, MR3115237} for additional readings on the Black-Scholes equation\footnote{Of course,
as every equation, the Black-Scholes equation
must be taken with a pinch of salt.
For instance, according to Ian Stewart's article appeared in the British newspaper
The Observer
{\tt https://www.theguardian.com/science/2012/feb/12/black-scholes-equation-credit-crunch}
the Black-Scholes equation was ``one ingredient in a rich stew of financial irresponsibility, political ineptitude,
perverse incentives and lax regulation'' that ultimately led to
the 2007 financial crisis.
``The equation itself wasn't the real problem. It was useful, it was precise, and its limitations were clearly stated [...]. The trouble was its potential for abuse''. In particular, we stress that the equation relies on the knowledge of the market volatility~$\sigma$.
``This is a measure of how erratically its market value changes. The equation assumes that the asset's volatility remains the same for the lifetime of the option, which needs not be correct. Volatility can be estimated by statistical analysis of price movements but it can't be measured in a precise, foolproof way, and estimates may not match reality [...].
 The equation also assumes that there are no transaction costs, no limits on short-selling and that money can always be lent and borrowed at a known, fixed, risk-free interest rate. Again, reality is often very different''.} and related models.
 
\subsection{Diffusion of transition probabilities}\label{9-42coMS02urjf0iurjfpnewnblkfxd}

In a stochastic model, one is often interested in the so-called transition probability,\index{transition probability}
namely the likelihood of transitioning from one state to another.

\begin{figure}
                \centering
                \includegraphics[width=.55\linewidth]{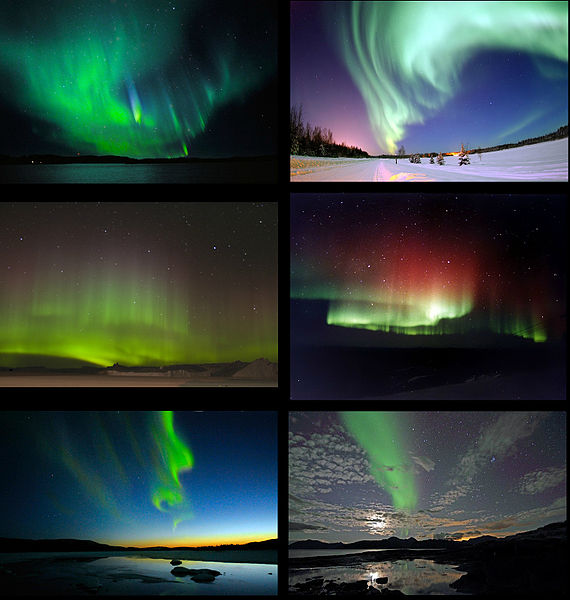}
        \caption{\sl Collection of pictures of aurorae from around the world (photos by
        Mila Zinkova, Samuel Blanc, Joshua Strang, Varjisakka and Jerry Magnum Porsbjer;
        images from
        Wikipedia,
        licensed under the Creative Commons Attribution-Share Alike 1.0 Generic license).}\label{AUR2HAFODAKEDFUMSPierre-SimonLaplacldRGIRA4AXELEUMDJOMNFHARLROA7789GIJ7solFUMHDNOJHNFOJED231}
\end{figure}

\begin{figure}
                \centering
                \includegraphics[width=.35\linewidth]{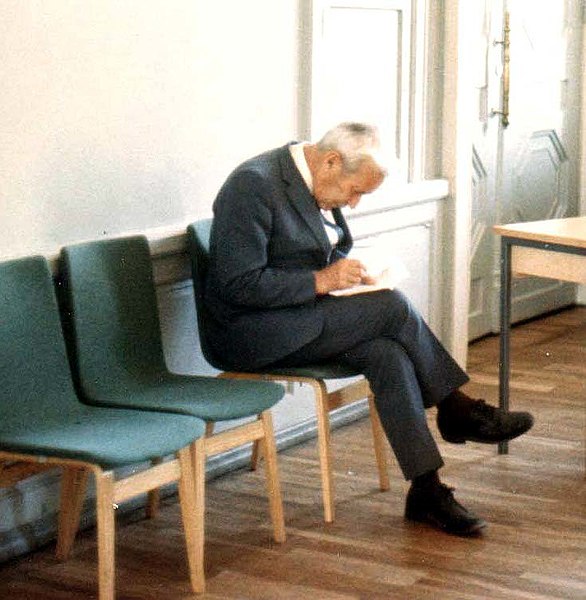}
        \caption{\sl Getting ready for a big presentation (photo by Terrence L. Fine; image from
        Wikipedia,
        licensed under the Creative Commons Attribution-Share Alike 3.0 Unported license).}\label{2HAFODAKEDFUMSPierre-SimonLaplacldRGIRA4AXELEUMDJOMNFHARLROA7KOLM789GIJ7solFUMHDNOJHNFOJED231}
\end{figure}

In particular, Markov processes\index{Markov process} are stochastic models ``without memory'', i.e.
in which predictions about future outcomes can be done based solely on the knowledge of its present state.
In this setting, if~$P(x_0,t_0| x,t)$ denotes the probability of going from the point~$x_0$ at time~$t_0$
to the point~$x$ at time~$ t$ (say, with~$t>t_0$),
we have that such a probability can be expressed as the superposition of all the
probabilities of going first to a point~$y$ at some time~$s\in(t_0,t)$
and from there reaching~$x$ at time~$t$: more explicitly, given~$s\in (t_0,t)$, it holds that
\begin{equation}\label{CHAKOL}
P(x_0,t_0|x,t)=\int_{\R^n} P(x_0,t_0|y,s)\,
P(y,s|x,t)\,dy.
\end{equation}
This relation is often called
Chapman-Kolmogorov equation\index{Chapman-Kolmogorov equation},
after\footnote{Chapman was a mathematician and geophysicist. At age 16 he entered the University of Manchester with a scholarship (he was the last student selected) and graduated with an engineering degree, but his passion for mathematics drove him to study for one further year to take a mathematics degree.

Besides his eminent contributions to stochastic processes, Chapman contributed to the theory of geomagnetism, investigating the beautiful phenomenon of aurorae (see Figure~\ref{AUR2HAFODAKEDFUMSPierre-SimonLaplacldRGIRA4AXELEUMDJOMNFHARLROA7789GIJ7solFUMHDNOJHNFOJED231}) in relation to the Earth's magnetic field and the solar wind, 
contributing to the understanding of the photochemical mechanisms that produce the ozone layer,
and predicting the existence of the magnetosphere (confirmed experimentally 30 years later).

Kolmogorov's scientific contributions are paramount \label{KOLMONOTE}
and comprise basically all fields of mathematics, including probability, topology, logic, mathematical physics, harmonic analysis, mathematical biology and numerical analysis. His mathematical talent appeared quite early in his life: at the age of five, he wrote his first mathematical paper, published in the journal of his school. The content of the paper was the formula
$$ \sum_{k=0}^{N-1} (2k+1)=N^2,$$
which is easy peasy for all of us, but not quite a piece of cake for a five-year-old boy (by the way, the little kid also became editor of the mathematical section of the school journal).

His former student and prominent mathematician Vladimir Igorevich Arnol\cprime d used to place Kolmogorov in his top 5 mathematicians list (with Poincar\'e, Gau{\ss}, Euler and Newton, see Figure~\ref{2HAFODAKEDFUMSPierre-SimonLaplacldRGIRA4AXELEUMDJOMNFHARLROA7KOLM789GIJ7solFUMHDNOJHNFOJED231}
for a picture of Kolmogorov preparing his talk during a conference in Tallin in 1973).

A mark on Kolmogorov's reputation is however produced by his involvement in the so-called "Luzin Affair" during the Great Purge in~1936 \label{PURFO2}
(see also footnote~\ref{PURFO} on page~\pageref{PURFO}: on this occasion, Nikolai Nikolaevich Luzin, prominent mathematical analyst and point-set topologist and Kolmogorov's former doctoral advisor (and author of the classical Luzin's Theorem, see e.g.~\cite[Theorem~4.20]{MR3381284}) became a target of Stalin's regime. Luzin was accused of plagiarism, nepotism, of being an enemy of the Soviet people and a servant to ``fascistoid science'' (whatever it means), as confirmed by the fact that Luzin published some of his results in foreign journals.

In view of this allegation, Luzin lost his academic position (not too bad after all, 
his advisor Dimitri Fyodorovich Egorov, author of Egorov's Theorem, see~\cite[Theorem~4.17]{MR3381284}, was arrested during a previous purge and died after a hunger strike initiated in prison).

Regrettably, Kolmogorov took part in the hearing at the Commission of the Academy of Sciences of the USSR, where the allegations against Luzin were formalized. The question of whether Kolmogorov was coerced by the police into testifying against Luzin, possibly
using as a threat an alledged long-lasting homosexual relationship,
remains a topic of speculation, see e.g.~\cite{MR2526973}. In any case, a sad story of envy, violence, indifference, discrimination and brutality.
Which is good to know, since
those who fail to learn from the mistakes done in the past
are doomed to repeat them.}
Sydney Chapman and Andrey Nikolaevich Kolmogorov.

It is clearly desirable, given an initial point~$x_0$ at an initial time~$t_0$,
to know more about the likelihood of being somewhere else at some future time, that is to have more information on
the function~$p(x,t):=P(x_0,t_0|x,t)$. 

\begin{figure}
                \centering
                \includegraphics[width=.6\linewidth]{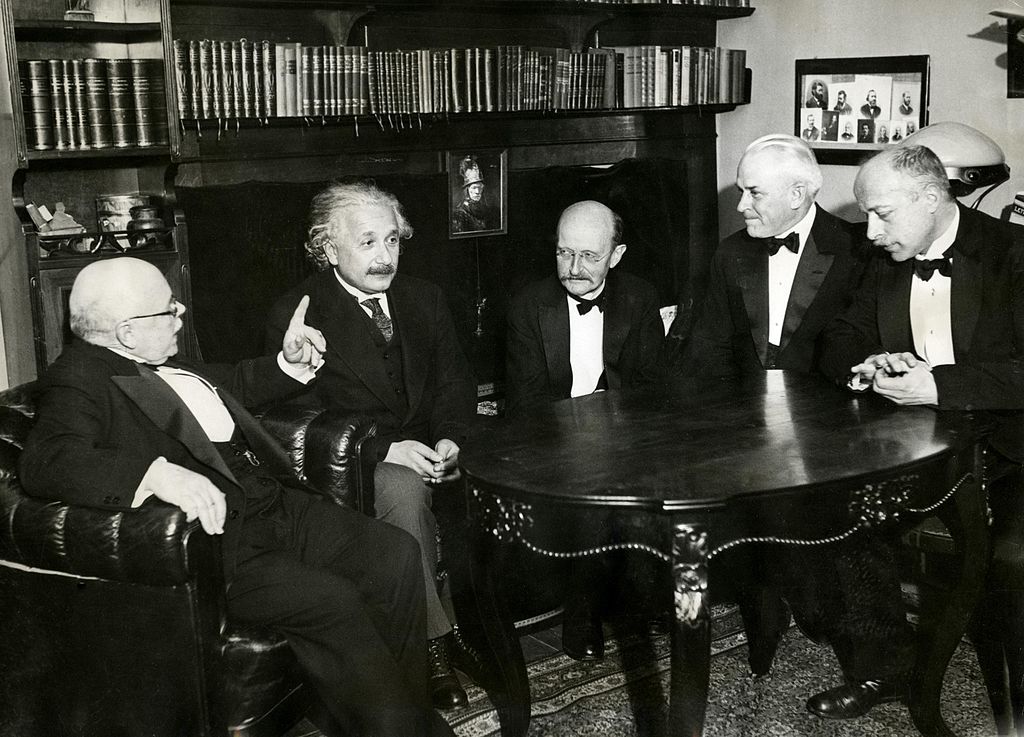}
        \caption{\sl Nobel Laureates only
       (Public Domain image from
 Wikipedia).}\label{2HAFODAKEDFUMSPierre-SimonLaplacldRGIRA4AXELEUMDJOMNFHARLROA7KOLM789GIJ7solFUMHDNOJHNFOJED2312}
\end{figure}

This is obtained
through a partial differential equation of the type\footnote{Equation~\eqref{FOKPLA}
is usually called Fokker-Planck equation\index{Fokker-Planck equation}. 
The name comes after Adriaan Dani\"el Fokker and Max Karl Ernst Ludwig Planck.

Fokker was a physicist and a musician, inventor of
a 31-tone equal-tempered organ
(we recall that
the 12-tone equal temperament
is the most common musical system today, recall the musical digression on page~\pageref{JS-PGAGPA}).
A picture of Fokker playing his organ is available at {\tt https://upload.wikimedia.org/wikipedia/en/a/ae/AdriaanFokker.jpg}

Planck is the famous originator of quantum theory: fortunately, he did not follow the advice of
his professor Philipp von Jolly at the 
Ludwig-Maximilians-Universit\"at M\"unchen, according to whom it was not worth to go into physics
since ``in this field, almost everything is already discovered, and all that remains is to fill a few holes''.

The Planck family lived for many years in a villa in Berlin-Grunewald, located in Wangenheimstra{\ss}e 21,
who became a gathering place for other professors, such as Albert Einstein and theologian Adolf von Harnack 
(brother of the mathematician portrayed on Figure~\ref{HARAXELHARLROAGIJ7soloDItangeFI}
on \label{gHARAXELHARLROAGIJ7soloDItangeFIPAGI} page~\pageref{HARAXELHARLROAGIJ7soloDItangeFIPAGI}).

Planck family suffered hard times due to the deplorable political situation around.
At the onset of World War I, Planck was one of the signatories of the so-called
Manifesto of the Ninety-Three
stating an unequivocal support of German military actions; among the other ninety-two prominent signatories,
Fritz Haber, Adolf von Harnack, Philipp Lenard (Nobel Prize for cathode rays research,
future supporter of the Nazi ideology and future crusader against Einstein's relativity theory)
and Felix Klein (the outstanding mathematician who devise the ``the Klein bottle''). During World War I,
Planck's oldest son Karl was killed in action at the Battle of Verdun
and his second son Erwin was taken prisoner by the French.

When the Nazis came to power, Planck was 74 and, when asked
to gather distinguished professors to issue a public proclamation against the expulsion of Jewish academics,
Planck replied, ``If you are able to gather today 30 such gentlemen, then tomorrow 150 others will come and speak against it, because they are eager to take over the positions of the others''.
However, Planck tried to avoid the expulsion of Fritz Haber
(who had pioneered chemical warfare
by introducing the use of poisonous gases during World War I and 
personally overseeing the first use of chlorine gas during the Second Battle of Ypres,
and who had received the Nobel Prize in Chemistry in 1918 for the invention of a process to synthesize ammonia from nitrogen and hydrogen gases, with broad application in the synthesis of fertilizers and explosives). 
Planck's support in this circumstance was unsuccessful
and Haber died in exile the following year.

Planck's son Erwin was later arrested by the Gestapo following the attempted assassination of Hitler on 20 July 1944,
sentenced to death and executed by hanging.

On a positive note, see Figure~\ref{2HAFODAKEDFUMSPierre-SimonLaplacldRGIRA4AXELEUMDJOMNFHARLROA7KOLM789GIJ7solFUMHDNOJHNFOJED2312}
for a social dinner
among Nobel Laureates: 
from left to right: Walther Nernst (1920 Nobel Prize for Chemistry), Albert Einstein
(1921 Nobel Prize for Physics), Max Planck (1918 Nobel Prize for Physics),
Robert A. Millikan (1923 Nobel Prize for Physics)
and Max von Laue (1914 Nobel Prize for Physics).}
\begin{equation}\label{FOKPLA}
\partial_t p(x,t)\,dx=-\div\big(\alpha(x,t)\,p(x,t)\big)+ 
\sum_{i,j=1}^n \partial_{ij}\big(\beta_{ij}(x,t)\,p(x,t)\big)
,\end{equation}
where~$\alpha$ and~$\beta_{ij}$ are called drift vector and diffusion tensor, respectively.

We now give a motivation for equation~\eqref{FOKPLA}.
Let~$\varphi$ be a smooth and compactly supported function in~$\R^n$.
Using the Chapman-Kolmogorov equation in~\eqref{CHAKOL}, we know that, given~$s\in (t_0,t)$,
\begin{eqnarray*}&&
\int_{\R^n}\varphi(x)\,p(x,t)\,dx=
\int_{\R^n}\varphi(x)\,P(x_0,t_0|x,t)\,dx\\&&\qquad=\iint_{\R^n\times\R^n} P(x_0,t_0|y,s)\,
P(y,s|x,t)\,\varphi(x)\,dx\,dy
\end{eqnarray*}
and, as a result, differentiating in~$t$,
\begin{equation}\label{frasuetartialarphi}
\int_{\R^n}\varphi(x)\,\partial_t p(x,t)\,dx=\iint_{\R^n\times\R^n } P(x_0,t_0|y,s)\,
\partial_t P(y,s|x,t)\,\varphi(x)\,dx\,dy
.\end{equation}
It is also useful to recall that the total probability is normalized to~$1$, namely
\begin{equation*} \int_{\R^n} P(y,s|x,t)\,dx=1,\end{equation*}
whence, taking a derivative in~$t$,
\begin{equation}\label{BASDBAefTTrotlk4453} \int_{\R^n} \partial_t P(y,s|x,t)\,dx=0.\end{equation}
Given~$y\in\R^n$, we also use a Taylor expansion to write that, when~$|x-y|$ is small,
$$ \varphi(x)=\varphi(y)+\nabla\varphi(y)\cdot(x-y)+\frac12 D^2\varphi(y)(x-y)\cdot(x-y)+O(|x-y|^3).$$
This and~\eqref{BASDBAefTTrotlk4453} give that
\begin{equation}\label{poe0586uyrhqwertyuiolkjhgfdszxcvbnm0987654}
\begin{split}
&\int_{\R^n} \partial_t P(y,s|x,t)\,\varphi(x)\,dx\\
=&\int_{\R^n} \partial_t P(y,s|x,t)\,\left(\varphi(y)+\nabla\varphi(y)\cdot(x-y)+\frac12 D^2\varphi(y)(x-y)\cdot(x-y)+O(|x-y|^3)
\right)\,dx\\
=&\int_{\R^n} \partial_t P(y,s|x,t)\,\left(\nabla\varphi(y)\cdot(x-y)+\frac12 D^2\varphi(y)(x-y)\cdot(x-y)+O(|x-y|^3)
\right)\,dx\\
=&\,\alpha(y,s,t)\cdot\nabla\varphi(y)+\frac12 \sum_{i,j=1}^n\beta_{ij}(y,t)\partial_{ij}\varphi(y)+
\int_{\R^n} \partial_t P(y,s|x,t)\,O(|x-y|^3)\,dx,
\end{split}\end{equation}
where\footnote{To be precise, we note that
the latter integral in~\eqref{poe0586uyrhqwertyuiolkjhgfdszxcvbnm0987654}
\label{LATENRINT}
occurs only in the support of~$\varphi$, which is a bounded set.}
\begin{eqnarray*}
&& \alpha(y,s,t):=
\int_{\R^n} \partial_t P(y,s|x,t)\,(x-y)\,dx\\
{\mbox{and }}&&
\beta_{ij}(y,s,t):=\frac12\int_{\R^n} \partial_t P(y,s|x,t)\,(x_i-y_i)(x_j-y_j)\,dx.
\end{eqnarray*}
Combining this and~\eqref{frasuetartialarphi}, we find that
\begin{equation}\label{NEHDF-HSJDJI2}
\begin{split}&
\int_{\R^n}\varphi(x)\,\partial_t p(x,t)\,dx\\
=\,&\int_{\R^n} 
\Bigg[
\alpha(y,s,t)\cdot\nabla\varphi(y)+ \sum_{i,j=1}^n\beta_{ij}(y,s,t)\partial_{ij}\varphi(y)\\&\qquad\qquad+
\int_{\R^n} \partial_t P(y,s|x,t)\,O(|x-y|^3)\,dx\Bigg]\,p(y,s)\,dy \\
=\,&
\int_{\R^n} \alpha(y,s,t)\cdot\nabla\varphi(y)\,p(y,s)\,dy\\&\qquad\qquad+\int_{\R^n} 
\sum_{i,j=1}^n\beta_{ij}(y,s,t)\partial_{ij}\varphi(y)\,p(y,s)\,dy +{\mathcal{R}}(s,t),\end{split}\end{equation}
where
\begin{eqnarray*}{\mathcal{R}}(s,t):=
\iint_{\R^n\times\R^n} \partial_t P(y,s|x,t)\,O(|x-y|^3)\,dx\,dy.
\end{eqnarray*}
We suppose additionally that the stochastic process is homogeneous in time, namely
the transition probability depends only on the points~$x_0$ and~$x$ and on the elapsed time~$t-t_0$, that is
$$ P(x_0,t_0+\tau|x,t+\tau)=P(x_0,t_0|x,t)$$
for all~$\tau\in\R$.

Differentiating in~$\tau$, this gives that
$$ \partial_{t_0} P(x_0,t_0|x,t)+\partial_{t} P(x_0,t_0|x,t)=0.$$
As a consequence, observing that~$
P(y,t|x,t)=0$
when~$x\ne y$ (unless one is provided with the ``gift of ubiquity''),
\begin{eqnarray*}
\alpha(y,t)&:=&\lim_{\e\searrow0}
\frac1\e\int_{t-\e}^t\alpha(y,s,t)\,ds\\&=&\lim_{\e\searrow0}
\frac1\e
\iint_{\R^n\times(t-\e,t)} \partial_t P(y,s|x,t)\,(x-y)\,dx\,ds\\&
=&-\lim_{\e\searrow0}
\frac1\e\iint_{\R^n\times(t-\e,t)} \partial_s P(y,s|x,t)\,(x-y)\,dx\,ds\\&
=&\lim_{\e\searrow0}
\frac1\e\int_{\R^n}\Big( P(y,t-\e|x,t)-P(y,t|x,t)\Big)
\,(x-y)\,dx\\&
=&\lim_{\e\searrow0}
\frac1\e\int_{\R^n}P(y,t-\e|x,t)\,(x-y)\,dx
\\&=&\lim_{\e\searrow0}
\int_{\R^n}\frac{P(y,0|x,\e)\,(x-y)}\e\,dx.
\end{eqnarray*}
Similarly,
\[ \beta_{ij}(y,t):=\lim_{\e\searrow0}
\frac1\e\int_{t-\e}^t\beta_{ij}(y,s,t)\,ds=\lim_{\e\searrow0}
\frac12\int_{\R^n} \frac{ P(y,0|x,\e)}{\e}\,(x_i-y_i)(x_j-y_j)\,dx\]
and
\begin{equation}\label{NEHDF-HSJDJI}
{\mathcal{R}}(t):=\lim_{\e\searrow0}
\frac1\e\int_{t-\e}^t{\mathcal{R}}(s,t)\,ds=\lim_{\e\searrow0}
\iint_{\R^n\times\R^n} \frac{ P(y,0|x,\e)}{\e}\,O(|x-y|^3)\,dx\,dy.
\end{equation}
Formally, we also have that
\begin{eqnarray*}&&
\lim_{\e\searrow0}\frac1\e
\iint_{\R^n\times(t-\e,t)} \alpha(y,s,t)\cdot\nabla\varphi(y)\,p(y,s)\,dy\,ds\\&=&
\lim_{\e\searrow0}\frac1\e\left(
\iint_{\R^n\times(t-\e,t)} \alpha(y,s,t)\cdot\nabla\varphi(y)\,p(y,t)\,dy\,ds\right.\\&&\qquad\quad\left.+
\iint_{\R^n\times(t-\e,t)} \alpha(y,s,t)\cdot\nabla\varphi(y)\,\big(p(y,s)-p(y,t)\big)\,dy\,ds\right)\\&=&
\int_{\R^n} \alpha(y,t)\cdot\nabla\varphi(y)\,p(y,t)\,dy+
\lim_{\e\searrow0}\frac1\e
\iint_{\R^n\times(t-\e,t)} \alpha(y,s,t)\cdot\nabla\varphi(y)\,\big(p(y,s)-p(y,t)\big)\,dy\,ds
\\&=&
\int_{\R^n} \alpha(y,t)\cdot\nabla\varphi(y)\,p(y,t)\,dy+
\lim_{\e\searrow0}\frac1\e
\iint_{\R^n\times(t-\e,t)} |\nabla\varphi(y)|\,O(\e)\,dy\,ds\\&=&
\int_{\R^n} \alpha(y,t)\cdot\nabla\varphi(y)\,p(y,t)\,dy
\end{eqnarray*}
and, in a similar way,
\begin{eqnarray*}
\lim_{\e\searrow0}\frac1\e
\iint_{\R^n\times(t-\e,t)} \beta_{ij}(y,s,t)\partial_{ij}\varphi(y)\,p(y,s)\,dy
=\int_{\R^n} \beta_{ij}(y,t)\partial_{ij}\varphi(y)\,p(y,t)\,dy.\end{eqnarray*}
Therefore, integrating~\eqref{NEHDF-HSJDJI2} in~$s\in(t-\e,\e)$, dividing by~$\e$ and sending~$\e\searrow0$,
we see that
\begin{equation}\label{NEHDF-HSJDJI2COSetyt6} 
\begin{split}&
\int_{\R^n}\varphi(x)\,\partial_t p(x,t)\,dx
=
\int_{\R^n} \alpha(y,t)\cdot\nabla\varphi(y)\,p(y,t)\,dy\\&\qquad\qquad\qquad+\int_{\R^n} 
\sum_{i,j=1}^n\beta_{ij}(y,t)\partial_{ij}\varphi(y)\,p(y,t)\,dy +{\mathcal{R}}(t).\end{split}\end{equation}

\begin{figure}
  \centering
  \includegraphics[height=.18\textheight]{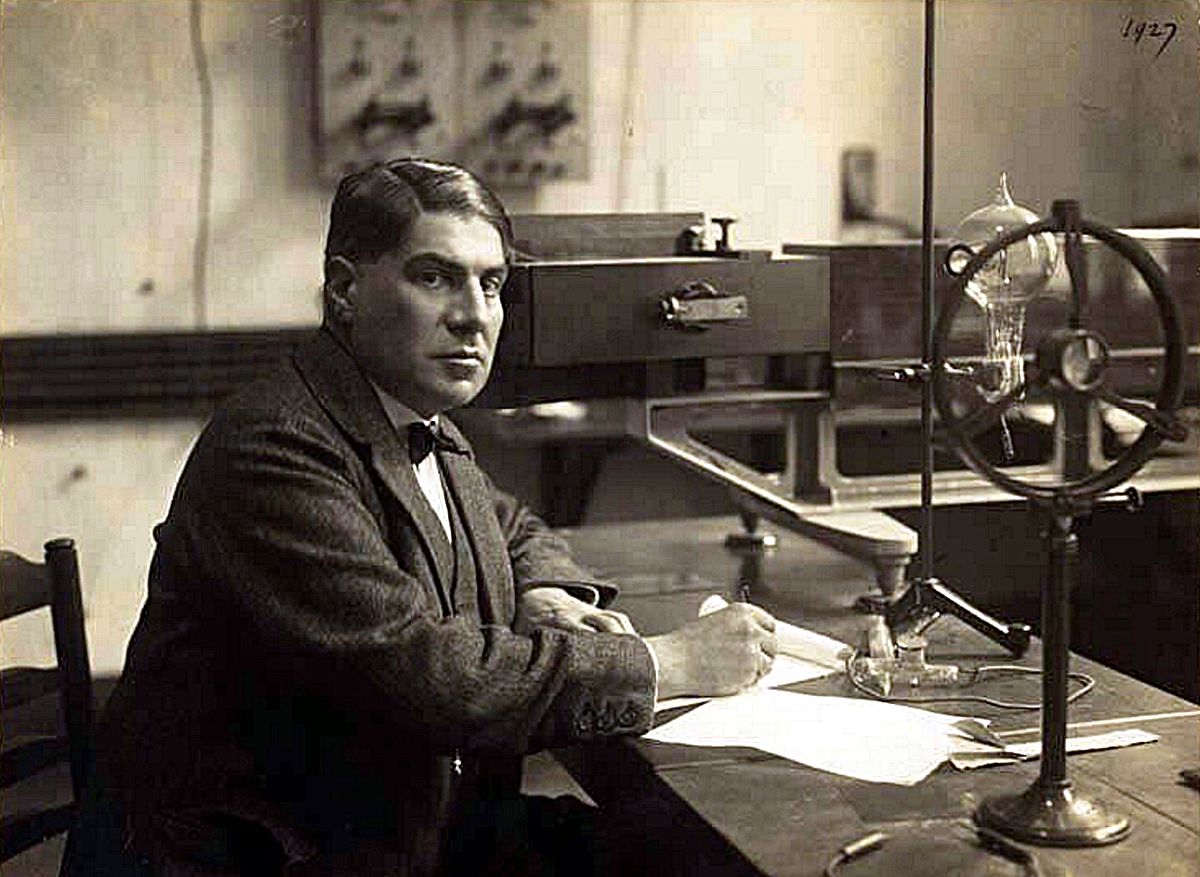}$\quad$
  \includegraphics[height=.18\textheight]{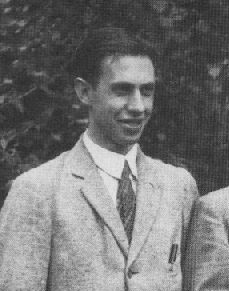}
 \caption{\sl Leonard Ornstein and George Uhlenbeck (Public Domain images from
 Wikipedia).}\label{G243NFIDItangeFI243ORNUL}
\end{figure}

Now, for a sufficiently small time~$\e$, if the stochastic process is continuous
we expect the probability of going from~$y$ to~$x$ in the elapsed time~$\e$
to be ``rather small'' except when~$|x-y|$ is small. It is therefore customary to consider the quantity~${\mathcal{R}}(t)$ in~\eqref{NEHDF-HSJDJI} as a ``negligible term''
(one may also want to recall the observation in footnote~\ref{LATENRINT}
to get rid of the contribution in~${\mathcal{R}}(t)$ coming from infinity).
Hence, taking the liberty of formally dropping it from the equation, we reduce~\eqref{NEHDF-HSJDJI2COSetyt6} to
\begin{equation*}\begin{split}&
\int_{\R^n}\varphi(x)\,\partial_t p(x,t)\,dx=\int_{\R^n} \alpha(y,t)\cdot\nabla\varphi(y)\,p(y,t)\,dy\\&\qquad\qquad+ \int_{\R^n} 
\sum_{i,j=1}^n\beta_{ij}(y,t)\partial_{ij}\varphi(y)\,p(y,t)\,dy.\end{split}\end{equation*}
Integrating by parts (and renaming $y$ into~$x$ in the variable of integration) we obtain
\begin{equation*}\begin{split}&
\int_{\R^n}\varphi(x)\,\partial_t p(x,t)\,dx\\&\qquad=-\int_{\R^n} \varphi(x)\,\div\big(\alpha(x,t)\,p(x,t)\big)\,\,dx+\int_{\R^n} 
\sum_{i,j=1}^n \varphi(x)\,\partial_{ij}\big(\beta_{ij}(x,t)\,p(x,t)\big)\,dx,\end{split}\end{equation*}
which, by the arbitrariness of~$\varphi$, produces the
Fokker-Planck equation in~\eqref{FOKPLA}, as desired.
\medskip

\begin{figure}
                \centering
                \includegraphics[width=.42\linewidth]{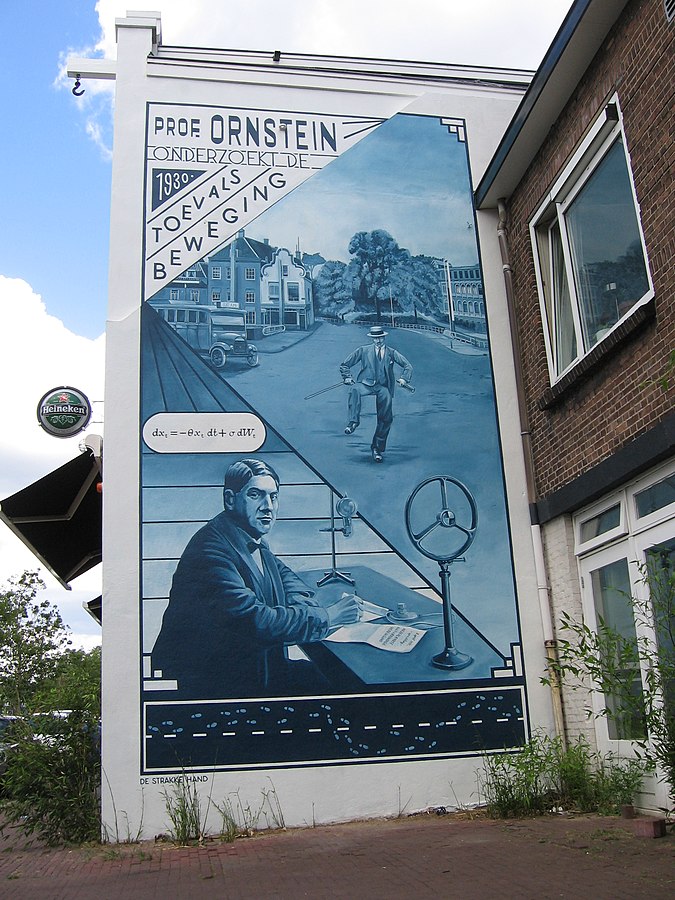}
        \caption{\sl Leonard Ornstein mural (image by Hansmuller from
        Wikipedia,
        licensed under the Creative Commons Attribution-Share Alike 4.0 International license).}\label{2HAFODAKEDFUMSPierre-SimonLaplacldRGIRA4AXELEUMDJOMNFHARLROA7789GIJ7solFUMHDNOJHNFOJED231MURA1}
\end{figure}

Another possible derivation of the Fokker-Planck equation in~\eqref{FOKPLA} comes from
stochastic differential equations (we present this argument in dimension~$1$ for the sake of simplicity, but the same
ideas would carry over in higher dimensions too). Namely, if~$dX=F\,dt+G\,dW$
we can use It\^o's Chain Rule in~\eqref{ITOPIOTNNRRRESIDEKKA}
to every random variable~$H(X(t),t)$, integrate and find that
\begin{equation*} 
H(X(T),T)-H(X(0),0)=
\int_0^T \,dH=
\int_0^T \left[
\left( \partial_t H +F \partial_X H+\frac{G^2}{2} \partial_X^2 H\right)\,dt+G\partial_X H \,dW\right].\end{equation*}
In particular, if~$H=H(X(t))$,
\begin{equation}\label{WIENE-01} 
H(X(T))-H(X(0))=
\int_0^T \left[
\left(F \partial_X H+\frac{G^2}{2} \partial_X^2 H\right)\,dt+G\partial_X H \,dW\right].\end{equation}
Now we take the expected value~${\mathbb{E}}$ of this identity (say, at a given time~$T$). For this, it is useful to recall that,
for a ``reasonable'' $\phi$, we have that
\begin{equation}\label{WIENE-02} {\mathbb{E}}\left(\int_0^T \phi\,dW\right)=0,\end{equation}
see e.g.~\cite[equation~(ii) on page~68]{MR3154922}
(the bottom line of this formula being that~$W$ presents independent increments hence the expectation of a function~$\phi$
that is adapted to the stochastic process will vanish).

Thus, it follows from~\eqref{WIENE-01} and~\eqref{WIENE-02} that
\begin{equation}\label{WIENE-05} {\mathbb{E}}\Big(
H(X(T))-H(X(0))\Big)={\mathbb{E}}\left(
\int_0^T
\left( F \partial_X H+\frac{G^2}{2} \partial_X^2 H\right)\,dt\right).\end{equation}
Now we recall the definition of expected value of a random variable~$Y$ (see e.g.~\cite[equation~(4) on page~14]{MR3154922}), according to which
$$ {\mathbb{E}}(Y)=\int_{\R^n} Y(x) \,p(x,T)\,dx.$$
In this way, equation~\eqref{WIENE-05} becomes
\begin{equation*}\begin{split}& \int_{\R^n} \Big(
H(x)-H(X(0))\Big)\,p(x,T)\,dx\\&\quad=
\iint_{\R^n\times(0,T)}
\left( F(x,t) \partial_x H(x)+\frac{G^2(x,t)}{2} \partial_x^2 H(x)\right)\,p(x,T)\,dx\,dt.\end{split}\end{equation*}
We now write this identity for~$T+\e$, subtract the one for~$T$, divide by~$\e$
and formally send~$\e\searrow0$:
namely, recalling~\eqref{BASDBAefTTrotlk4453}, 
\begin{eqnarray*}&&
\int_{\R^n} H(x)\,\partial_t p(x,T)\,dx\\&=&
\int_{\R^n} \Big(
H(x)-H(X(0))\Big)\,\partial_t p(x,T)\,dx\\
&=&\lim_{\e\searrow0}
\int_{\R^n} \Big(
H(x)-H(X(0))\Big)\,\frac{p(x,T+\e)-p(x,T)}\e\,dx\\
&=&
\lim_{\e\searrow0}\frac1\e
\left[\int_{\R^n} \Big(
H(x)-H(X(0))\Big)\,p(x,T+\e)\,dx-
\int_{\R^n} \Big(
H(x)-H(X(0))\Big)\,p(x,T)\,dx\right]\\&=&\lim_{\e\searrow0}\frac1\e
\Bigg[\iint_{\R^n\times(0,T+\e)}
\left( F(x,t) \partial_x H(x)+\frac{G^2(x,t)}{2} \partial_x^2 H(x)\right)\,p(x,T+\e)\,dx\,dt\\&&\qquad
-\iint_{\R^n\times(0,T)}
\left( F(x,t) \partial_x H(x)+\frac{G^2(x,t)}{2} \partial_x^2 H(x)\right)\,p(x,T)\,dx\,dt\Bigg]\\&=&\lim_{\e\searrow0}
\Bigg[\iint_{\R^n\times(0,T+\e)}
\left( F(x,t) \partial_x H(x)+\frac{G^2(x,t)}{2} \partial_x^2 H(x)\right)\,\frac{p(x,T+\e)-p(x,T)}\e\,dx\,dt\\&&\qquad
+\frac1\e\iint_{\R^n\times(T,T+\e)}
\left( F(x,t) \partial_x H(x)+\frac{G^2(x,t)}{2} \partial_x^2 H(x)\right)\,p(x,T)\,dx\,dt\Bigg]\\&=&
{\mathcal{S}}(T)
+\int_{\R^n}
\left( F(x,T) \partial_x H(x)+\frac{G^2(x,T)}{2} \partial_x^2 H(x)\right)\,p(x,T)\,dx,
\end{eqnarray*}
where
$$ {\mathcal{S}}(T):=\iint_{\R^n\times(0,T)}
\left( F(x,t) \partial_x H(x)+\frac{G^2(x,t)}{2} \partial_x^2 H(x)\right)\,\partial_t p(x,T)\,dx\,dt.$$

\begin{figure}
                \centering
                \includegraphics[width=.4\linewidth]{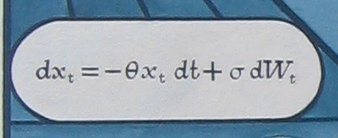}
        \caption{\sl The Ornstein-Uhlenbeck process, detail from the Leonard Ornstein mural (image by Hansmuller from
        Wikipedia,
        licensed under the Creative Commons Attribution-Share Alike 4.0 International license).}\label{2HAFODAKEDFUMSPierre-SimonLaplacldRGIRA4AXELEUMDJOMNFHARLROA7789GIJ7solFUMHDNOJHNFOJED231MURA2}
\end{figure}

Since the term~${\mathcal{S}}(T)$ is small when~$T$ is small, if we neglect it from the previous
computation we reduce the problem to
\[ \int_{\R^n} H(x)\,\partial_t p(x,T)\,dx=
\int_{\R^n}
\left( F(x,T) \partial_x H(x)+\frac{G^2(x,T)}{2} \partial_x^2 H(x)\right)\,p(x,T)\,dx.\]
Integrating by parts in~$x$, we thus conclude that
\[ \int_{\R^n} H(x)\,\partial_t p(x,T)\,dx=
\int_{\R^n}
\left( -\partial_x\Big(F(x,T) \,p(x,T)\Big)+\frac12\partial_x^2\Big(
{G^2(x,T)}\, p(x,T)\Big)\right)\,H(x)\,dx.\]
{F}rom the arbitrariness of~$H$
we thus obtain
\begin{equation}\label{KMS76ygiTSg12rZasdG} \partial_t p=-\partial_x \big(F p\big)+\frac12\partial_x^2\big({G^2} p\big),\end{equation}
which is the Fokker-Planck equation~\eqref{FOKPLA} in this setting:
interestingly,
in the expression~\eqref{KMS76ygiTSg12rZasdG} one can directly related the drift coefficient of the Fokker-Planck equation
with the ``deterministic'' term of the stochastic equation and the diffusion coefficient of the Fokker-Planck equation
with the ``random'' term of the stochastic equation.
\medskip

\begin{figure}
                \centering
                \includegraphics[width=.6\linewidth]{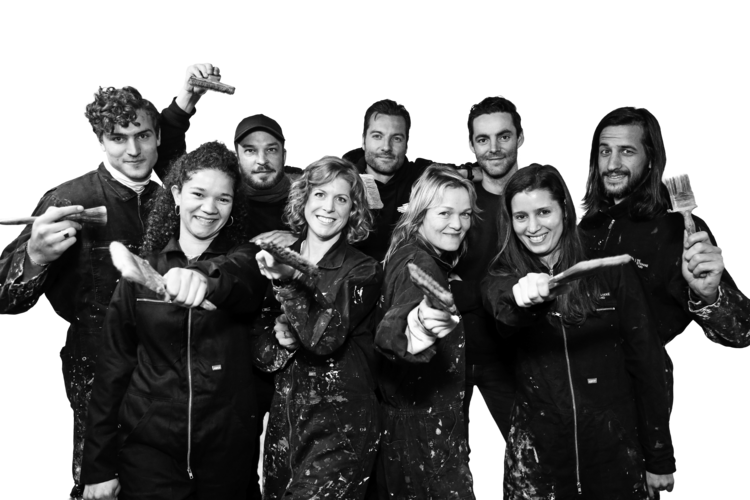}
        \caption{\sl De Strakke Hand Team (image from De Strakke Hand website
{\tt  http://www.destrakkehand.nl/about}).}\label{2HAFODAKEDFUMS-HANDGIRA4AXELEUMDJOMNFHARLROA7789GIJ7solFUMHDNOJHNFOJED231MURA2}
\end{figure}

A classical application for the Fokker-Planck equation in~\eqref{KMS76ygiTSg12rZasdG} is given by the
so-called\footnote{The process in~\eqref{ORST-ULHL-S-04} is named after
Leonard Salomon Ornstein and George Eugene Uhlenbeck, see Figure~\ref{G243NFIDItangeFI243ORNUL}.

Compare also the photo of Ornstein in his lab with
Figure~\ref{2HAFODAKEDFUMSPierre-SimonLaplacldRGIRA4AXELEUMDJOMNFHARLROA7789GIJ7solFUMHDNOJHNFOJED231MURA1}, showing a very elegant
mural in Utrecht, The Netherlands, depicting Ornstein, 
the random walk (embodied by a drunkard, check the bottle of
booze in his left hand) and the Ornstein-Uhlenbeck
process (note the formula next to Ornstein's head,
as enlarged in Figure~\ref{2HAFODAKEDFUMSPierre-SimonLaplacldRGIRA4AXELEUMDJOMNFHARLROA7789GIJ7solFUMHDNOJHNFOJED231MURA2}).

The mural has been painted by the Dutch painting collective
De Strakke Hand, see Figure~\ref{2HAFODAKEDFUMS-HANDGIRA4AXELEUMDJOMNFHARLROA7789GIJ7solFUMHDNOJHNFOJED231MURA2}.
See also {\tt http://www.destrakkehand.nl/} for other beautiful pieces of art in the urban landscape.

Sometimes, models similar to~\eqref{ORST-ULHL-S-04} are also referred to with the name
of Langevin equation, \index{Langevin equation}
after Paul Langevin.} Ornstein-Uhlenbeck\index{Ornstein-Uhlenbeck process} stochastic differential equation.
In this case, one considers in the background a stochastic differential equation of the form
\begin{equation}\label{ORST-ULHL-S-04} dX=-\vartheta \,X\,dt+\sigma \,dW,
\end{equation}
for suitable positive constant parameters~$\vartheta$ and~$\sigma$.

To develop an intuition of the process in~\eqref{ORST-ULHL-S-04}, we can imagine an
elastic spring in the presence of a large dumping and thermal fluctuations. In a nutshell, Hooke's Law would prescribe
an equation of motion of the type~$\ddot X=-\kappa X$, with~$\kappa>0$ being the elastic constant
of the spring. In presence of friction, the previous equation becomes~$\ddot X=-\kappa X-\gamma \dot X$,
being~$\gamma$ the friction coefficient. Rewriting this equation in the form~$\frac{\ddot X}{\gamma}=-\vartheta X-\dot X$,
with~$\vartheta:=\kappa/\gamma$, for large values of elastic and dumping coefficients we can reduce ourselves
to~$0=-\vartheta X-\dot X$, that we can formally write as~$dX=-\vartheta X\,dt$. In this sense,
the Ornstein-Uhlenbeck equation in~\eqref{ORST-ULHL-S-04} consists simply in adding to this model a random fluctuation
(e.g. due to thermal deviations)
and therefore it can be considered as describing a noisy relaxation process of an oscillator.
\medskip

Notice that the setting in~\eqref{ORST-ULHL-S-04} corresponds to that in~\eqref{KMS76ygiTSg12rZasdG}
with the choices~$F:=-\vartheta \,X$ and~$G:=\sigma$, whence the corresponding 
Fokker-Planck equation for the Ornstein-Uhlenbeck process takes the form
\begin{equation}\label{ORST-ULHL-S-05} \partial_t p=\vartheta\partial_x \big(x p\big)+\frac{\sigma^2}2\partial_x^2 p.\end{equation}
It is interesting to observe that the drift term modulated by~$\vartheta$ in~\eqref{ORST-ULHL-S-05}
corresponds to an attraction towards the center~$x=0$ for the transition probability~$p$.
To convince ourselves of this fact one can just compare~\eqref{ORST-ULHL-S-05}
with~\eqref{CHEMLKSDEQ}, which was describing the evolution of a biological population following a chemical attractant:
notice that~$w$ in~\eqref{CHEMLKSDEQ} would correspond, up to constants, to~$-\vartheta x^2$ in~\eqref{ORST-ULHL-S-05},
hence~$p$ is ``attracted'' towards the higher values of~$-\vartheta x^2$ (namely, $x=0$) in the same way
as a biological population is attracted towards the higher density regions of the chemotactic factor.

\begin{figure}
                \centering
                \includegraphics[width=.43\linewidth]{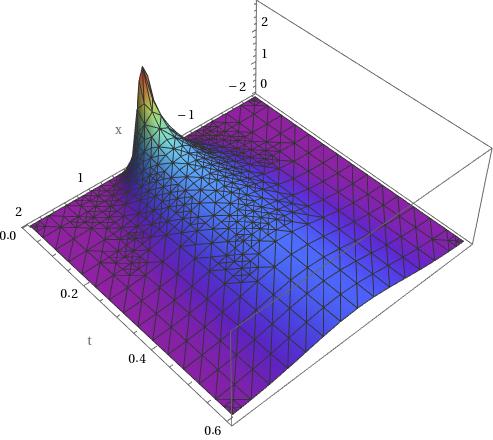}
        \caption{\sl Plotting the evolution of the transition probability for the Ornstein-Uhlenbeck process
        according to the Fokker-Planck equation in~\eqref{ORST-ULHL-S-05}.}\label{RIVEKELVHARGIRA4AXELHARLROA7789GIJ7soFOKPLAULlFUMHDNOJHNFOJED}
\end{figure}

A plot of the solution~$p$ of~\eqref{ORST-ULHL-S-05}, for a suitable choice of the parameters~$\vartheta$ and~$\sigma$, is sketched in Figure~\ref{RIVEKELVHARGIRA4AXELHARLROA7789GIJ7soFOKPLAULlFUMHDNOJHNFOJED}:
the initial datum at~$t=0$ corresponds to the probability of finding the random particle located at~$x=1$
(and for this reason the picture exhibits a blowup at~$(x,t)=(1,0)$) and we can appreciate that, as time flows,
the solution has the tendency to move its mass towards~$x=0$.

\begin{figure}
                \centering
                \includegraphics[width=.21\linewidth]{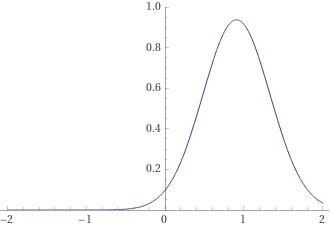}$\quad$
                \includegraphics[width=.21\linewidth]{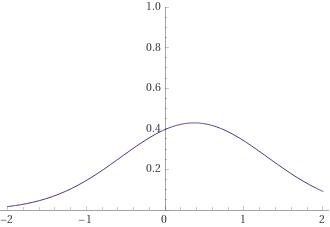}$\quad$
                \includegraphics[width=.21\linewidth]{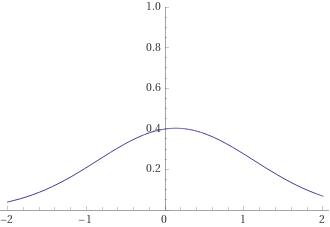}$\quad$
                \includegraphics[width=.21\linewidth]{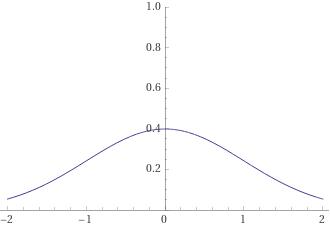}
        \caption{\sl Evolution of the transition probability for the Ornstein-Uhlenbeck process
        at subsequent times.}\label{RIVEKELVHARGIRA4AXELHARLROA7789GIJ7soFOKPLAULlFUMHDNOJHNFOJEDsu}
\end{figure}

See also Figure~\ref{RIVEKELVHARGIRA4AXELHARLROA7789GIJ7soFOKPLAULlFUMHDNOJHNFOJEDsu} for
frames of this solution at times~$t=0.1$, $1$, $2$, $10$.

\begin{figure}
                \centering
                \includegraphics[height=.23\textheight]{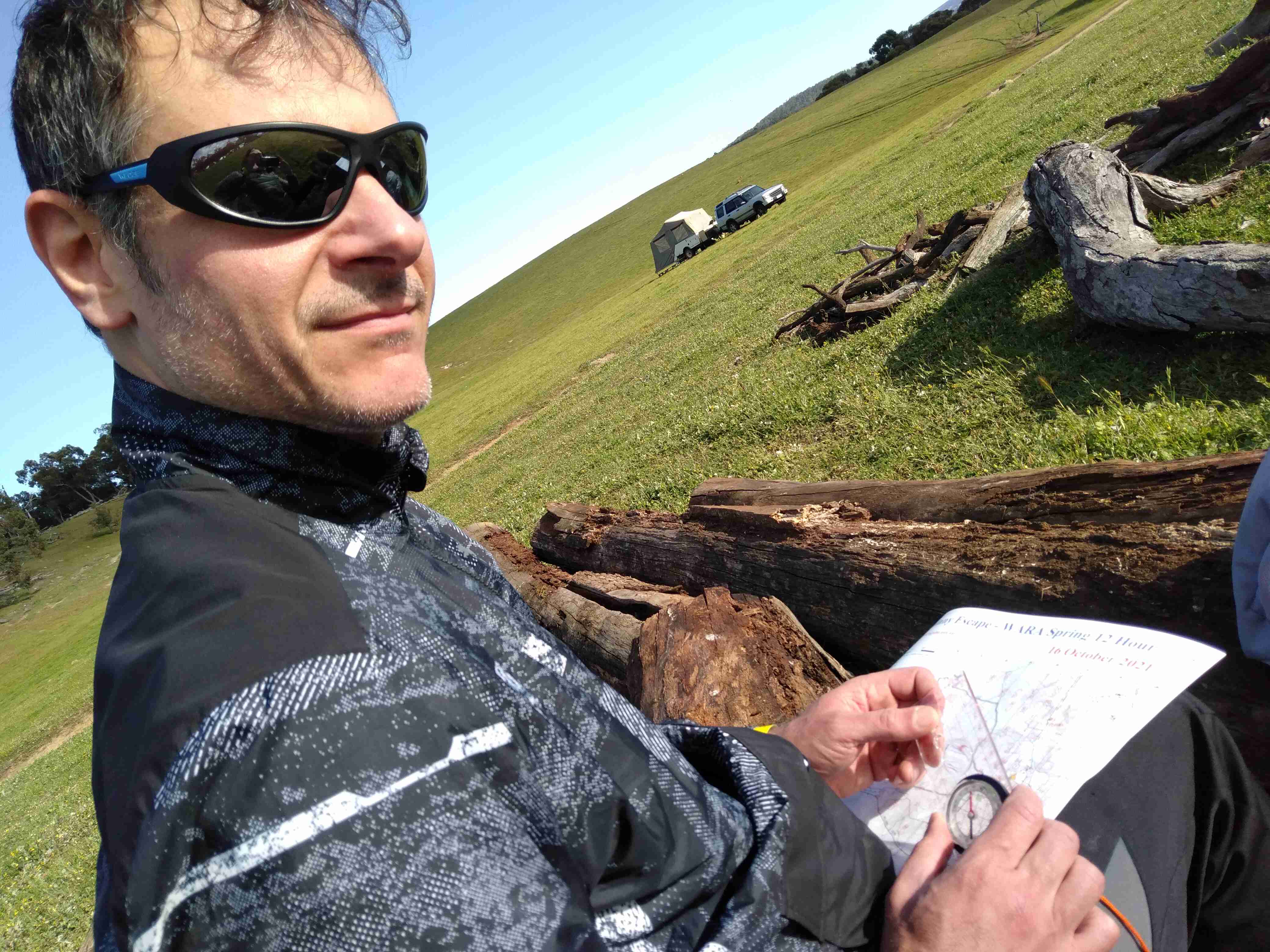}$\,$
                \includegraphics[height=.23\textheight]{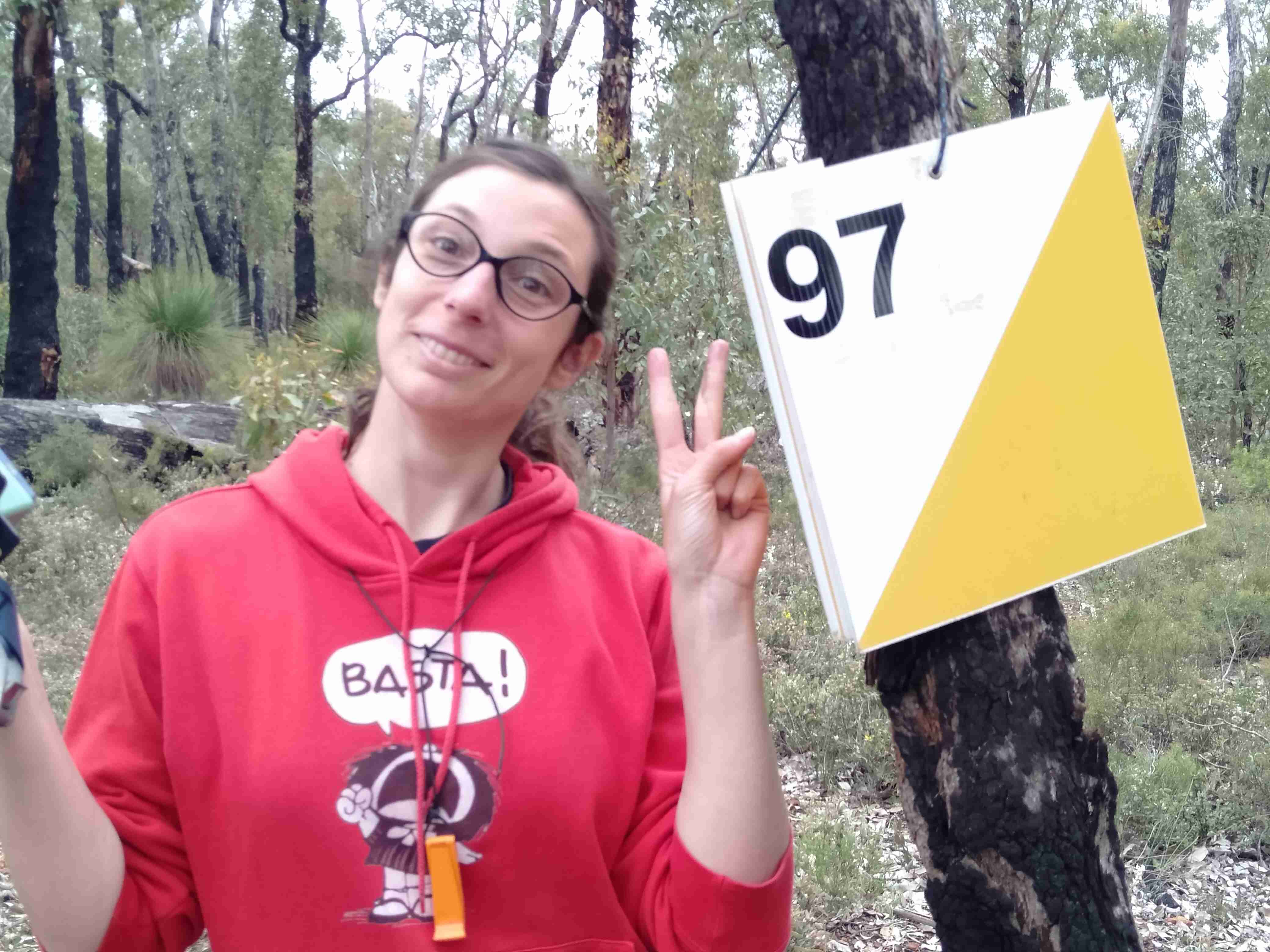}$\,$
                \includegraphics[height=.23\textheight]{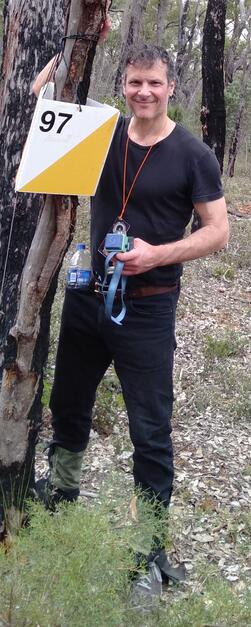}
        \caption{\sl The difference between planning routes on a map and getting there.}\label{RIVEKELVHARGIRA4AXELHARLROA7789GIJ7soFOKPLAULlFUMHDNOJHNFOMAPJEDsu}
\end{figure}

The situation in which the general tendency of reaching an objective (such as the origin
in Figures~\ref{RIVEKELVHARGIRA4AXELHARLROA7789GIJ7soFOKPLAULlFUMHDNOJHNFOJED}
and~\ref{RIVEKELVHARGIRA4AXELHARLROA7789GIJ7soFOKPLAULlFUMHDNOJHNFOJEDsu})
is weighted against the possibility of diffusing around and missing it
is actually quite common in our everyday experience, hence
the plots in Figures~\ref{RIVEKELVHARGIRA4AXELHARLROA7789GIJ7soFOKPLAULlFUMHDNOJHNFOJED}
and~\ref{RIVEKELVHARGIRA4AXELHARLROA7789GIJ7soFOKPLAULlFUMHDNOJHNFOJEDsu}
should be somewhat close to our intuition. One can think for instance
to the case of rogaining, in which, in principle compass bearings point you straight to
the next control, but in reality one has to account for some fluctuation around the goal
(due to some miscalculations of angles and distances for which the teammate has to be blamed, see e.g. Figure~\ref{RIVEKELVHARGIRA4AXELHARLROA7789GIJ7soFOKPLAULlFUMHDNOJHNFOMAPJEDsu},
to the roughness of the terrain, to possible modifications of the landscape, to possible inaccuracies of the map, etc.): accordingly, the evolution of an initial position is realistically not given
by a single point (that is, a Dirac Delta Function located at a moving point)
but rather by a probability distribution describing how likely one is exactly at that point.
In this scenario, the origin in Figures~\ref{RIVEKELVHARGIRA4AXELHARLROA7789GIJ7soFOKPLAULlFUMHDNOJHNFOJED}
and~\ref{RIVEKELVHARGIRA4AXELHARLROA7789GIJ7soFOKPLAULlFUMHDNOJHNFOJEDsu} represents the ideal destination corresponding to the maximal value of the probability distribution, but the tails of this distribution represent the possibility of having missed the control during navigation
due to whatever fluctuation (at that point, the rogainer will have to seek the control patiently
and accurately; knowing the probability density would be perhaps a good indication on how far one should search).\medskip

For additional information about Brownian motion,
random disturbances and the Fokker-Planck equation,
see~\cite{MR2962, MR34975, MR987631, MR2676235, MR3408563} and the references therein.

\subsection{Phase coexistence models}

A topical problem in material sciences is the understanding of the separation
patterns between different phases of a given substance.
The two phases could be related to molecule orientation, magnetization, 
multi-component alloy systems, state of the matter (including superconductive
or superfluid states). In this context, the description of a physical system often relies on the understanding of suitable ``order parameters'', that are suitable quantities capable of distinguishing different phases of a given system.

For instance, the magnetization vector~$M$ expresses the density of magnetic dipole moment in a magnetic material. Roughly speaking, at each point~$x$ of a given magnetic material, the direction of~$M(x)$ reproduces the direction of the magnetic field near~$x$ and the magnitude~$|M(x)|$ represents the amount of magnetization.

Another classical example of order parameter arises in nematic crystals. \index{nematic crystal}These materials consist of long and thin molecules that prefer to align with one another. In this situation, the two ends of the elongated molecules are essentially indistinguishable, therefore an efficient order parameter for these crystals is given by a vector~$N(x)$ modulo a sign (that is, an element of the projective plane, or equivalently a vector on a hemisphere with opposing points along the equator identified). That is, while magnetic materials tend to displace the magnetic vectors in a parallel fashion, the nematic vectors can be equally displayed either in a parallel or antiparallel form (the technological advantage is thus that these crystals can be aligned by an external magnetic or electric field, as it happens for instance in the liquid crystal displays for video games).\medskip

While the above mentioned~$M$ and~$N$ are vectorial order parameters,
there are also ideal cases in which a scalar order parameter is sufficiently efficient to
describe some features of a system. For example, nematic crystals
may exhibit at some sufficiently large scale a ``preferred direction'' along which most
of its elements are aligned, which is called in jargon ``local nematic director''. \index{local nematic director} In this context, scalar state parameters which are often adopted for practical purposes
are the functions of the angle~$\vartheta$ between the liquid-crystal molecular axis and the local nematic director.

\begin{figure}
  \centering
  \includegraphics[width=.8\linewidth]{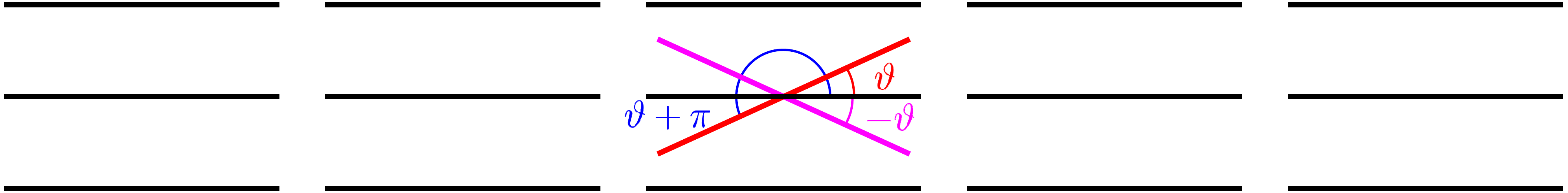}
 \caption{\sl A scalar order parameter using the local nematic director.}\label{NEMADIRDItangeFI}
\end{figure}

The simplest of this function may be constructed as follows. We take, up to a rotation, the 
local nematic director to be the first vector~$e_1$ of the Euclidean basis
and we consider a state parameter~$S(\vartheta)$ such that:
\begin{itemize}
\item[(i).] $S$ is $2\pi$-periodic (consistently with the fact that~$\vartheta$ is an angle),
\item[(ii).] $S(\vartheta)=S(\vartheta+\pi)$ (consistently with the fact that
the two endings of the nematic crystals are indistinguishable),
\item[(iii).] $S(\vartheta)=S(-\vartheta)$ (consistently with the fact that this
order parameter should account for the ``angular distance'' of a molecule to the local nematic director),
\item[(iv).] $S(\vartheta)\le 1=S(0)$ (normalizing the maximal of the state parameter with~$1$,
corresponding to the full alignment case),
\item[(v).] $\displaystyle\fint_{\partial B_1} S\big(\arccos(\omega\cdot e_1)\big)\,d{\mathcal{H}}^{n-1}_\omega=0$
(normalizing at~$0$ the case of random directions, using the notation~$\omega\cdot e_1=\cos\vartheta$),
\end{itemize}
see Figure~\ref{NEMADIRDItangeFI}.

By~(i), one can seek~$S$ in the form of a Fourier Series, say
$$ S(\vartheta)=\frac{a_0}2+\sum_{j=1}^{+\infty} \left(a_j\cos(j\vartheta)+b_j\sin(j\vartheta)\right).$$
By~(iii), we have that~$b_j=0$ for all~$j\in\{1,2,\dots\}$. Hence, using~(ii),
$$ 0=S(\vartheta+\pi)-S(\vartheta)
=\sum_{j=1}^{+\infty} a_j\left( \cos(j(\vartheta+\pi))-\cos(j\vartheta)\right)
=-2\sum_{{j\ge1}\atop{j{\tiny\mbox{ odd}}}} a_j \cos(j\vartheta)
$$
and therefore~$a_j=0$ for all~$j$ odd.

Thus, neglecting higher orders, one can take
$$ S(\vartheta)=\frac{a_0}2+a_2\cos(2\vartheta).$$
By~(iv), we have that~$\frac{a_0}2=1-a_2$, whence, in view of~(v),
\begin{eqnarray*}&&
0=\fint_{\partial B_1} S(\arccos(\omega\cdot e_1))\,d{\mathcal{H}}^{n-1}_\omega
=1-a_2+a_2\fint_{\partial B_1} \cos(2\arccos\omega_1)\,d{\mathcal{H}}^{n-1}_\omega\\&&\qquad\qquad=
1-a_2+a_2\fint_{\partial B_1} \big(2\cos^2(\arccos\omega_1)-1\big)\,d{\mathcal{H}}^{n-1}_\omega\\&&\qquad\qquad
=1-2a_2+2a_2\fint_{\partial B_1} \omega_1^2\,d{\mathcal{H}}^{n-1}_\omega\\&&\qquad\qquad
=1-2a_2+\frac{2a_2}{n},
\end{eqnarray*}
which, in the concrete case~$n=3$, yields that~$a_2=3/4$.

As a result,
$$ S(\vartheta)=\frac14+\frac34\cos(2\vartheta)=\frac12\big(3\cos^2\vartheta-1\big).$$
Interestingly, we will discover that this coincides with the
second Legendre polynomial\index{Legendre polynomial}
(see footnote~\ref{LISTACON} on page~\pageref{LISTACON}).
\medskip

Coming back to the study of phase transitions,
Landau \index{Landau theory}
theory\footnote{This theory is named after Lev Davidovich Landau, 
1962 Nobel Prize in Physics for his mathematical theory of superfluidity. \label{PURFO}
In~1938, Landau wrote a leaflet, joint with Mosey Korets, condemning
Stalin and the People's Commissariat for Internal Affairs (the forerunner of
MVD and KGB) on the occasion of the Great Purge in the Soviet Union
(see also footnote~\ref{PURFO2} on page~\pageref{PURFO2}).
According to the leaflet,
Stalin and his acolytes
``can only beat defenseless prisoners, catch unsuspecting innocent people, plunder national wealth and invent ridiculous trials against nonexistent conspiracies''.
The authors of the leaflet were sure that
``Proletariat [...] will throw off the fascist dictator and his clique''.

Things didn't work out as expected and Landau was arrested and held in the Lubyanka prison for more
than one year, see Figure~\ref{DAUGREENFIDItangeFI}.

There were however other features which possibly compensated Landau's dissident leaflet
in the eyes of Stalin. Landau led a team of mathematicians supporting Soviet atomic and hydrogen bomb development, thus receiving for his work the Stalin Prize in 1949 and 1953
and the title Hero of Socialist Labor in 1954. We should be always suspicious about ourselves
when we receive too many prizes and awards, because,
as all Spider-Man fans know, ``with great power comes great responsibility''.

By the way, according to Wikipedia,
``Landau believed in free love rather than monogamy and encouraged his wife
and his students to practise free love. However, his wife was not enthusiastic''.}
attempts to understand the phase transitions related to continuous order parameters. Its main ansatz is that the equilibria of the system should come from a ``free energy'' produced by the order parameter under consideration. This energy may also be sensitive to the temperature of the system. For instance, one can consider an order parameter~$\eta$ for a given system at temperature~$T$ and describe the energy~$E$ of the system as a function of~$\eta$ and~$T$ (for the sake of simplicity,
we take here~$\eta$ to be a scalar order parameter, but the case
of vector valued order parameters can be treated in a similar manner). 

Up to a normalization, one can suppose that the energy corresponding to~$\eta=0$ is zero. Moreover, in many concrete cases, the energy of the physical system is invariant if we exchange~$\eta$ with~$-\eta$ (this is the case, for instance, of magnetic materials, since exchanging the North with the South pole should not alter the energy state of the system, and it is also the case of nematic crystals, since the ending of their molecules is symmetric, recall e.g. point~(iii) here above).

These considerations lead to an energy~$E(\eta,T)$ which is even in~$\eta$ and such that~$E(0,T)=0$.
The corresponding Taylor expansion
must therefore contain only even terms and takes the form
$$ E(\eta,T)=a(T)\eta^2+b(T)\eta^4+\dots$$
and we will indeed neglect higher order terms, thus reducing simply to the case
\begin{equation}\label{ENEFRE} E(\eta,T)=a(T)\eta^2+b(T)\eta^4.\end{equation}
The energy coefficients~$a(T)$ and~$b(T)$ depend on the temperature~$T$
and their signs play a decisive role in the formation and separation of phases. More specifically,
one usually considers stable solutions in correspondence to minimal energy levels:
therefore, to avoid minimizers corresponding to an energy equal to~$-\infty$, a natural structural assumption
is to suppose that~$b(T)>0$ for every temperature~$T$. Instead, to allow for
possible phase changes, one can assume that~$a(T)$ changes sign above and below some critical temperature~$T_c$,
for instance taking~$a(T)>0$ when~$T<T_c$ and~$a(T)<0$ when~$T>T_c$.

\begin{figure}
  \centering
  \includegraphics[width=.4\linewidth]{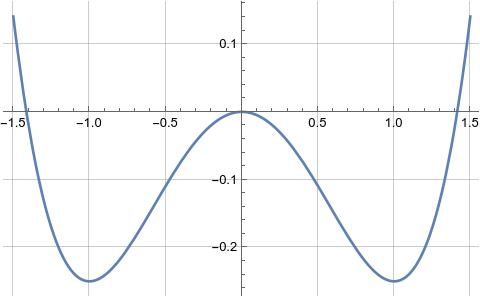}
 \caption{\sl Plot of the function~$\eta\mapsto\frac{(\eta^2-1)^2-1}4$.}\label{DOPPIOPODItangeFI}
\end{figure}

In this setting, below the critical temperature the free energy exhibits the null value of the order parameter~$\eta$
as its only minimizer, but above the critical temperature the stable
phase corresponds to~$\pm\eta_0$ where~$\eta_0$ is the minimizer for~$E$ given by
$$ \eta_0:=\sqrt{-\frac{a}{2b}},$$
see Figure~\ref{DOPPIOPODItangeFI}.

In particular, if given~$T>T_c$, one wishes to normalize the stable phases at the levels~$\pm1$
(as done for instance in point~(iv) here above), up to normalizing factors one can choose~$a(T)=-\frac12$ and~$b(T)=\frac14$,
so that~\eqref{ENEFRE} reduces to
\begin{equation}\label{DWHADEKLTSTMONAKN2}E(\eta)=-\frac12\eta^2+\frac14\eta^4=\frac{(\eta^2-1)^2-1}4\end{equation}
and the stable phases are the zeros of
\begin{equation}\label{DWHADEKLTSTMONAKN}
E'(\eta)=\eta^3-\eta.\end{equation}\medskip

It is interesting to remark that this analysis of the free energy is useful to detect
the stable phases, namely~$\pm1$ in~\eqref{DWHADEKLTSTMONAKN},
but it does not give any information on the separation between the two possible phases:
indeed, all configurations only attaining the stable phases~$\pm1$ would indeed be
zeros of~\eqref{DWHADEKLTSTMONAKN}, as well as minimizers of~\eqref{DWHADEKLTSTMONAKN2}.
However, in many concrete situations, one expects the separation between phases to
be somewhat minimal as well: one often experiences situations in which the two phases
occupy two separate ``bulks'' of the material which are separated by an interface, thus showing
a distinctive coarsening and phase separation.
This phenomenon is possibly the outcome of a ``ferromagnetic''
effect, which tends to align the direction of the molecules of the material,
hence avoiding oscillations of the order parameter.

To include this phenomenon into the phase separation model, one
can modify the energy in~\eqref{DWHADEKLTSTMONAKN2} by adding
a small penalization term which charges the formation of interfaces.
The simplest version of this procedure is typically to add to~\eqref{DWHADEKLTSTMONAKN2}
a small ``gradient term'' (a gradient indeed detects the ``local oscillation'' of phases), i.e.
replace~\eqref{DWHADEKLTSTMONAKN2} by
\begin{equation}\label{JohnWCahnandSam AllenE} {\mathcal{G}}(\eta)=\frac{\e}2\int_\Omega|\nabla\eta(x)|^2\,dx+\int_\Omega
\frac{(\eta^2(x)-1)^2-1}4\,dx,\end{equation}
for a small parameter~$\e>0$ (where~$\Omega$ is the region occupied by the material).

Looking for the corresponding minimal points of the full energy~${\mathcal{G}}$ leads to the equation
\begin{equation}\label{JohnWCahnandSam Allen}\e
\Delta\eta=\eta^3-\eta,
\end{equation}
which is often referred\footnote{Equation~\eqref{JohnWCahnandSam Allen}
is named after John Cahn and Sam Allen, who presented
a theoretical treatment of metal alloys~\cite{ALLEN19751017}.

With respect to the classification
presented in footnote~\ref{CLASSIFICATIONFOOTN} on page~\pageref{CLASSIFICATIONFOOTN},
we see that equation~\eqref{JohnWCahnandSam Allen} is of elliptic type.

We think that it is interesting to relate the ``democratic'' features of the Laplace operator,
according to the discussions on pages~\pageref{DEEFFNVST}
and~\pageref{DEEFFNVSTBIS}, to the ferromagnetic tendency of some materials which avoids
oscillations in spins and molecule alignments.

Furthermore, equation~\eqref{JohnWCahnandSam Allen} provides one of the chief examples
of ``semilinear equations'', within a structure that we will investigate further  \index{semilinear equation}
from page~\pageref{Pohozaev Identity} on. See also footnote~\ref{BIS Pohozaev Identity}
on page~\pageref{BIS Pohozaev Identity}. \label{TRIS Pohozaev Identity}

One can also appreciate the link between minimizers of the energy functional in~\eqref{JohnWCahnandSam AllenE}
and interfaces of minimal area, by using
the Cauchy-Schwarz Inequality and the Coarea Formula (see e.g.~\cite{MR3409135}), namely setting~$W(\eta):=
\frac{(\eta^2-1)^2}4$ and
noticing that
\begin{eqnarray*}&&{\mathcal{G}}(\eta)+\frac{|\Omega|}4=
\frac{\e}2\int_\Omega|\nabla\eta(x)|^2\,dx+\int_\Omega W(\eta(x))\,dx
\ge2\int_\Omega \sqrt{ \frac{\e}2|\nabla\eta(x)|^2W(\eta(x))}\,dx\\&&\qquad=
\sqrt{ 2{\e}}\int_\Omega |\nabla\eta(x)|\sqrt{W(\eta(x))}\,dx=\sqrt{2{\e}}
\int_\R \left[\int_{\{\eta=\tau\}} \sqrt{W(\eta(x))}\,d{\mathcal{H}}^{n-1}_x
\right]\,d\tau.
\end{eqnarray*}
Hence, heuristically, it should be convenient for a minimizer~$\eta_\e$ to sit ``as much as possible''
into the minima of~$W$ (which correspond to the stable phases~$\pm1$): if we assume that the level sets of~$\eta_\e$
are surfaces ``more or less parallel'' to an interface~${\mathcal{S}}$
and~$\eta_\e$ approaches ``quite fast'' the stable phases~$\pm1$, the previous
computation suggests that the minimization of~${\mathcal{G}}$ reduces, as a first approximation,
to the minimization of the area of the interface~${\mathcal{S}}$.

Of course, it is not easy to transform the above heuristic argument into a rigorous proof, since the setting in~\eqref{JohnWCahnandSam Allen} is that of a singular perturbation \index{singular perturbation},
namely the small parameter~$\e$ affects precisely the most significant term in equation~\eqref{JohnWCahnandSam Allen}.
As a matter of fact, the coherent development of arguments of this sort required the introduction
of a novel notion of asymptotics called~$\Gamma$-convergence, \index{$\Gamma$-convergence}
see~\cite{MR448194}.}
to with the name of Allen-Cahn equation. \index{Allen-Cahn equation}
\medskip

See e.g.~\cite{PHASEBOOK} and the references therein for additional details
on phase separation models.

\begin{figure}
  \centering
  \includegraphics[width=.45\linewidth]{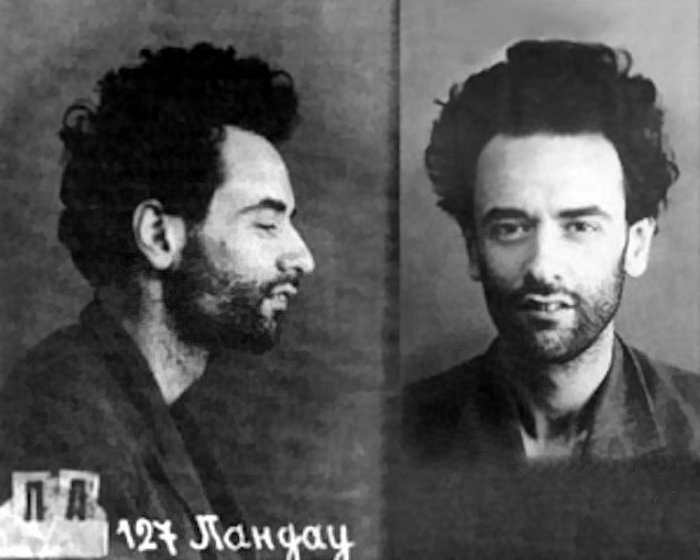}
 \caption{\sl Lev Landau in prison, photo by the People's Commissariat for Internal Affairs (Public Domain image from
 Wikipedia).}\label{DAUGREENFIDItangeFI}
\end{figure}

\subsection{Growth of interfaces}

A number of scientific problems of interest are associated with the growth of the profile of suitable surfaces which describe, for instance, tumors, flame fronts, clusters, etc. 

\begin{figure}
                \centering
                \includegraphics[width=.45\linewidth]{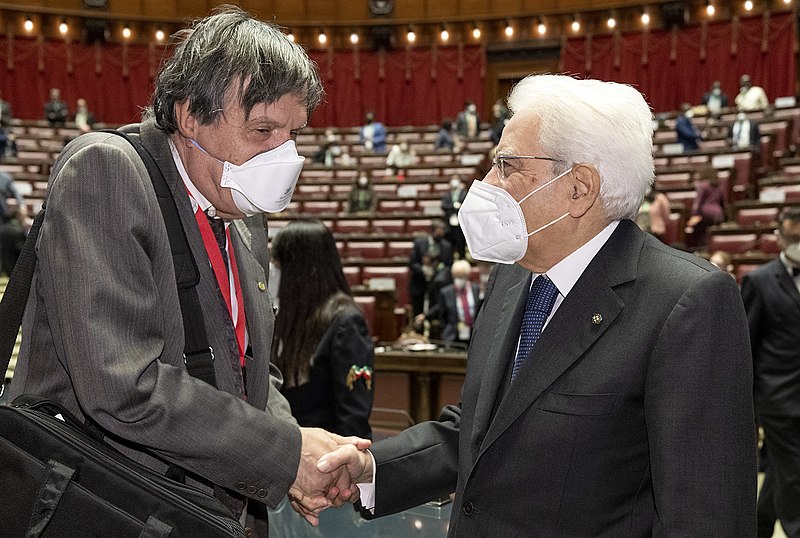}
        \caption{\sl Opening ceremony at the Italian Parliament (October 10, 2021;
        attribution: Presidenza della Repubblica {\tt https://www.quirinale.it/elementi/60152\#\&gid=1\&pid=10}).}\label{RIVEKELVHARGIRA4AXELHARLROA7789GIJ7solFUMHDNOJHNFOPARISJED}
\end{figure}

A typical model describing this phenomenon is given  by the so-called\footnote{Equation~\eqref{KPZEQUA} is named after Mehran Kardar, Giorgio Parisi and Yi-Cheng Zhang, who introduced this model in~\cite{zbMATH05115911}.
Parisi is also known for his contributions in
quantum fluctuations, spin glasses and whirling flocks of birds. See Figure~\ref{RIVEKELVHARGIRA4AXELHARLROA7789GIJ7solFUMHDNOJHNFOPARISJED}
in which the Italian President Sergio Mattarella (on the right) congratulates Parisi (on the left)
for having received the 2021 Nobel Prize in Physics (they are wearing face masks due to
COVID regulations).}
Kardar-Parisi-Zhang equation\index{Kardar-Parisi-Zhang equation}
\begin{equation}\label{KPZEQUA}
\partial_t h=\mu \Delta h+\lambda |\nabla h|^2+f,
\end{equation}
where~$\mu>0$ is a diffusion parameter, $\lambda$ is a parameter related to the
surface growth and~$f$ is a forcing term (which could be either a deterministic function
or a stochastic term, such as a white noise). Of course, when~$\lambda=0$ the setting in~\eqref{DAGB-ADkrVoiweLL4re2346ytmngrrUj7}
essentially boils down to that of the heat equation in~\eqref{DAGB-ADkrVoiweLL4re2346ytmngrrUj7}, but when~$\lambda\ne0$
the equation presents a nonlinear term of geometric importance.

Indeed, the motivation underpinning~\eqref{KPZEQUA} goes as follows. Suppose
that we have some aggregate with an active zone of growth on its surface
which is described, at a given time~$t$, by the graph of a certain function~$y=g(x,t)$.
We consider an ``infinitesimal'' time step~$\tau$ in which the surface growth takes
place and suppose that such a growth is regulated by three ingredients:
first, in time~$\left(t,t+\frac\tau3\right)$, some particle is added to the surface,
then, in time~$\left(t+\frac\tau3,t+\frac{2\tau}3\right)$ a forcing term kicks in,
finally, in time~$\left(t+\frac{2\tau}3,t+\tau\right)$ some diffusion takes place (the order
in which these phenomena occurs is not really important, but we fix this order just to keep a concrete case in mind).

\begin{figure}
                \centering
                \includegraphics[width=.70\linewidth]{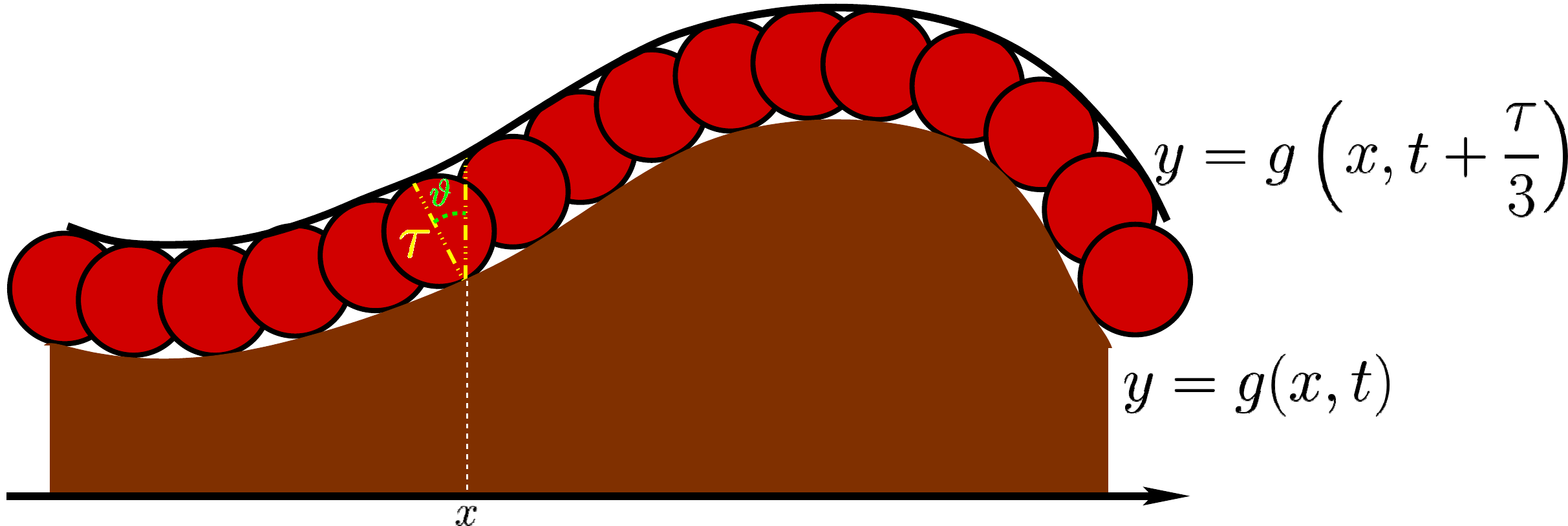}
        \caption{\sl Time evolution of a growing interface.}\label{AaCCRIVEKELVHARGIRA4AXELHARLROA7789GIJ7solFUMHDNOJHNFOPARISJED}
\end{figure}

To describe the first step, we assume that spherical particles of diameter~$\tau$
are added uniformly along the surface of the cluster, as described in Figure~\ref{AaCCRIVEKELVHARGIRA4AXELHARLROA7789GIJ7solFUMHDNOJHNFOPARISJED}.
In this setting, if $\vartheta$ denotes the angle between the normal
of the spherical particle and the vertical direction, we have that the infinitesimal vertical surface growth
corresponds to~$\frac\tau{\cos\vartheta}$. Hence we write that
\begin{equation}\label{KPZ-Eq-1} g\left(x,t+\frac\tau3\right)=g(x,t)+\frac\tau{\cos\vartheta}.\end{equation}
Actually, the normal of the spherical particle agrees with the normal of the original surface~$\frac{\left(-\nabla g(x,t),1\right)}{\sqrt{1+|\nabla g(x,t)|^2}}$ and accordingly
$$\cos\vartheta=\frac{\left(-\nabla g(x,t),1\right)}{\sqrt{1+|\nabla g(x,t)|^2}}\cdot e_n=
\frac{1}{\sqrt{1+|\nabla g(x,t)|^2}}.$$
We substitute this information into~\eqref{KPZ-Eq-1} and we find that
\begin{equation}\label{KPZ-Eq-2} g\left(x,t+\frac\tau3\right)=g(x,t)+ \tau \sqrt{1+|\nabla g(x,t)|^2}.\end{equation}

Now we describe the second step, namely the action of a forcing term: this effect is encoded in the equation
\begin{equation}\label{KPZ-Eq-3} g\left(x,t+\frac{2\tau}3\right)= g\left(x,t+\frac\tau3\right)+ \tau f\left(x,t+\frac\tau3\right).\end{equation}

As for the third step, we assume that diffusion takes place in the form
\begin{equation}\label{KPZ-Eq-4} g(x,t+\tau)= g\left(x,t+\frac{2\tau}3\right)+ \tau \Delta g\left(x,t+\frac{2\tau}3\right)
.\end{equation}

Thus, by collecting the observations in~\eqref{KPZ-Eq-2}, \eqref{KPZ-Eq-3} and~\eqref{KPZ-Eq-4},
\begin{eqnarray*}
&&\frac{g(x,t+\tau)-g(x,t)}\tau\\&=&\frac{\displaystyle g(x,t+\tau)-g\left(x,t+\frac{2\tau}3\right)}\tau
+\frac{ \displaystyle g\left(x,t+\frac{2\tau}3\right)-g\left(x,t+\frac{\tau}3\right)}\tau+
\frac{\displaystyle  g\left(x,t+\frac{\tau}3\right)-g(x,t)}\tau\\&=& \Delta g\left(x,t+\frac{2\tau}3\right)
+ f\left(x,t+\frac\tau3\right)+\sqrt{1+|\nabla g(x,t)|^2},
\end{eqnarray*}
whence, formally taking the limit as~$\tau\searrow0$,
\begin{equation}\label{KPZ-Eq-5} \partial_t g(x,t)= \Delta g(x,t)
+ f(x,t)+\sqrt{1+|\nabla g(x,t)|^2}
.\end{equation}
When the slope of the growing interface is small, one can consider~$|\nabla g|$ as a small perturbation
and employ the approximation~$\sqrt{1+|\nabla g|^2}\simeq 1+\frac{|\nabla g|^2}2$.
In this framework, one reduces~\eqref{KPZ-Eq-5} to
\begin{equation*} \partial_t g= \Delta g
+ f+1+\frac{1}{2}|\nabla g|^2.\end{equation*}
Considering the velocity shift~$h(x,t):=g(x,t)-t$, we thus conclude that
\begin{equation*} \partial_t h=
\partial_t g-1=
\Delta g
+ f+\frac{1}{2}|\nabla g|^2
=\Delta h
+ f+\frac{1}{2}|\nabla h|^2,
\end{equation*}
providing a motivation for the
Kardar-Parisi-Zhang equation in~\eqref{KPZEQUA}.\medskip

We mention that, in spite of its apparent simplicity, equations such as the one in~\eqref{KPZEQUA}
do capture complex phenomena and lead to the study of important
universality classes sharing the same characteristic exponents and scale invariant limits,
see e.g.~\cite{MR2930377}.

\begin{figure}
  \centering
  \includegraphics[width=.3\linewidth]{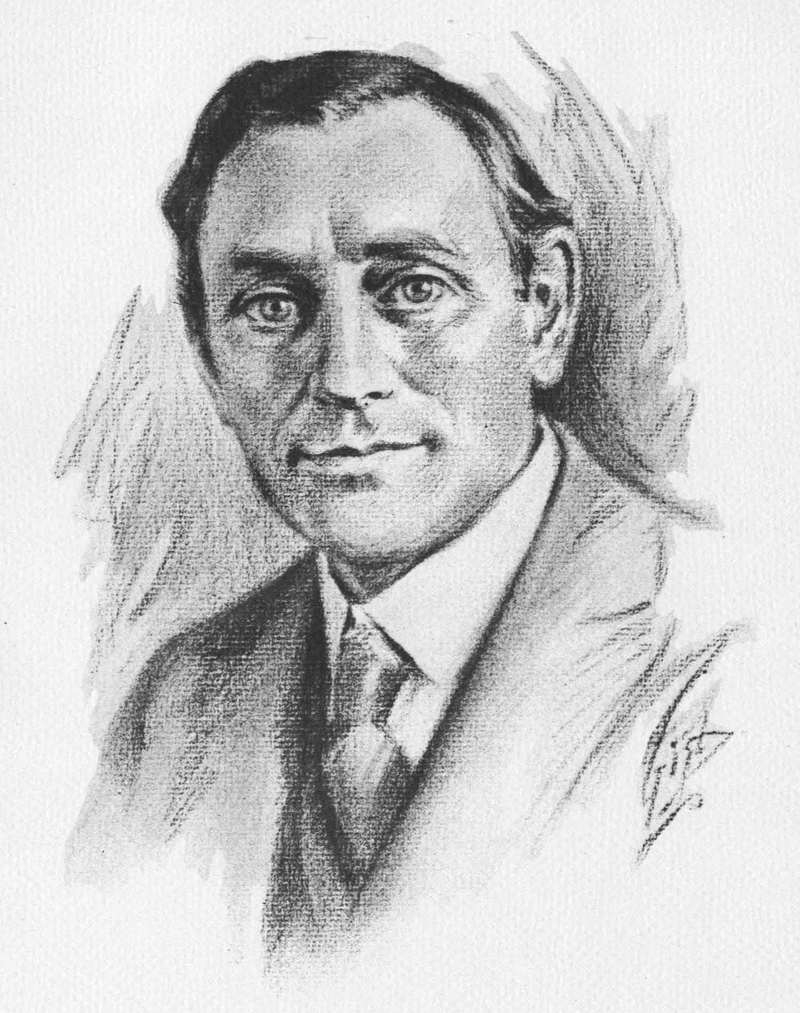}
 \caption{\sl Sketch of Harry Bateman (Public Domain image from
 Wikipedia).}\label{BATMGREENFIDItangeFI}
\end{figure}

Furthermore, equation~\eqref{KPZEQUA} is also directly related to other classical equations.
In particular, defining~$v:=\nabla h$, we have that
\begin{eqnarray*}&&
\partial_t v+v\cdot\nabla v-\mu\Delta v-\nabla f=\nabla\Big(
\partial_t h-\mu\Delta h-f\Big)+\sum_{j=1}^n \partial_jh\, \partial_j \nabla h\\&&\qquad=\nabla\Big(
\partial_t h-\mu\Delta h-f\Big)+\frac12\nabla \left(\sum_{j=1}^n |\partial_jh |^2\right)\\&&\qquad=\nabla\left(
\partial_t h-\mu\Delta h-f+\frac12|\nabla h|^2\right),
\end{eqnarray*}
from which it follows that the
Kardar-Parisi-Zhang equation in~\eqref{KPZEQUA} with~$\lambda:=\frac12$
is equivalent to the vectorial equation
\begin{equation}\label{VERBURG}
\partial_t v+v\cdot\nabla v=\mu\Delta v+\nabla f,
\end{equation}
which is a version of the Navier-Stokes equation in~\eqref{EUMSD-OS-DNS} in the absence of gravity
(in this setting, $v$ corresponds to the speed of a given fluid, $\mu$ is its viscosity coefficient
and~$f$ is minus the pressure).

For constant pressure, in dimension~$n=1$, equation~\eqref{VERBURG} reduces to the scalar equation
\begin{equation}\label{VERBURG2}
\partial_t v+v\partial_x v=\mu\partial_{xx} v,
\end{equation}
which is known in jargon\footnote{Equation~\eqref{VERBURG2}
is named after Harry Bateman and Johannes Martinus Burgers.
See Figure~\ref{BATMGREENFIDItangeFI} for a nice drawing of Harry Bateman
from the 1931 yearbook of the California Institute of Technology.} as the \index{Bateman-Burgers equation} Bateman-Burgers equation,
which models the speed of a fluid in a thin ideal pipe.

\begin{figure}
  \centering
  \includegraphics[height=.13\textheight]{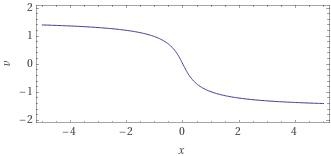} $\quad$
    \includegraphics[height=.13\textheight]{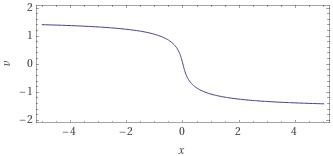} \\
  \includegraphics[height=.13\textheight]{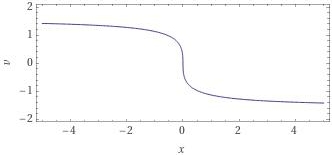} $\quad$
  \includegraphics[height=.13\textheight]{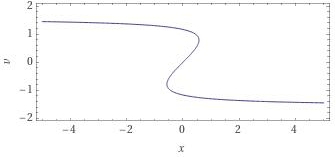}
 \caption{\sl Implicit plot of~$v=\arctan( tv-x)$ for~$t\in\{0.5,\,0.75,\,0.99,\,2\}$.}\label{SINGBATMGREENFIDItangeFI}
\end{figure}

In this context, it is interesting to point out that in the absence of viscosity equation~\eqref{VERBURG2} boils down to
\begin{equation}\label{VERBURG3}
\partial_t v+v\partial_x v=0,
\end{equation}
known as\footnote{Note that the inviscid Bateman-Burgers equation in~\eqref{VERBURG3}
can be seen as a particular case of the Korteweg-de Vries equation
in~\eqref{KDV-or} with~$a:=1$, $b:=0$ and~$c:=1$.}
the inviscid Bateman-Burgers equation. A feature of the solutions of~\eqref{VERBURG3}
is that they may develop ``shock waves'', \index{shock wave}
i.e. singularities at a breaking time. For instance,
it is readily seen that the relation~$v(x,t)=\arctan\big( tv(x,t)-x\big)$
provides, via the Implicit Function Theorem, a solution of~\eqref{VERBURG3} for small times, with~$v(x,0)=-\arctan x$.
This solution exhibits a shock wave since otherwise, for~$t>1$,
$$ |\partial_x v(0,t)|=\left|\frac{t\partial_x v(0,t)}{\cos^2\big(tv(0,t)\big)}\right|\ge | t\partial_x v(0,t)|>|\partial_x v(0,t)|,
$$
which is a contradiction. See Figure~\ref{SINGBATMGREENFIDItangeFI}
for a sketch of the formation of such a shock wave.

\begin{figure}
  \centering
  \includegraphics[height=.37\textheight]{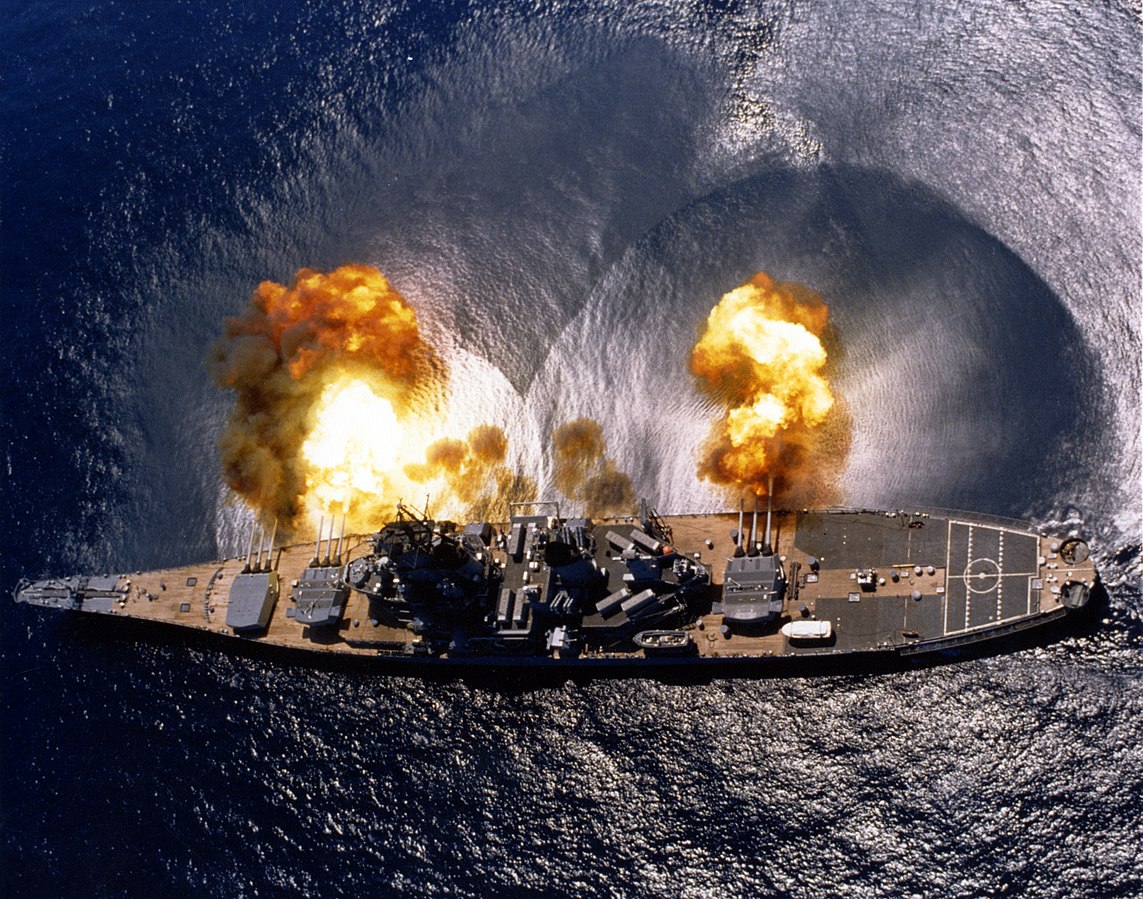}
 \caption{\sl The battleship USS Iowa of the United States Navy
 during a training exercise in Puerto Rico (Public Domain image from
 Wikipedia).}\label{SINGBATMGREENFIDItanSHO4geFI}
\end{figure}

See e.g.~\cite[pages~140, 195, 204]{MR1625845} for more information on shock waves
in the context of the Bateman-Burgers equation.
See also Figure~\ref{SINGBATMGREENFIDItanSHO4geFI} for a rather dramatic view of
expanding spherical atmospheric shock waves from a gun firing on the surface of the water.

\subsection{What to do when you're stuck in a traffic jam}

\begin{figure}
  \centering
  \includegraphics[height=.29\textheight]{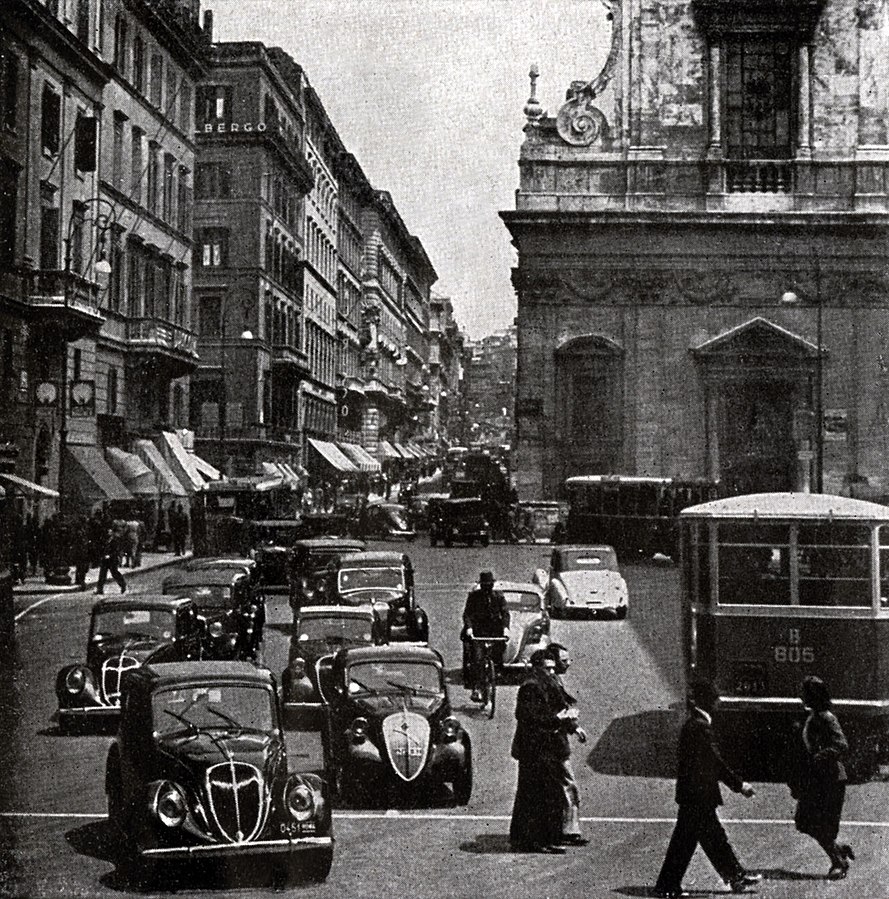}$\quad$
    \includegraphics[height=.29\textheight]{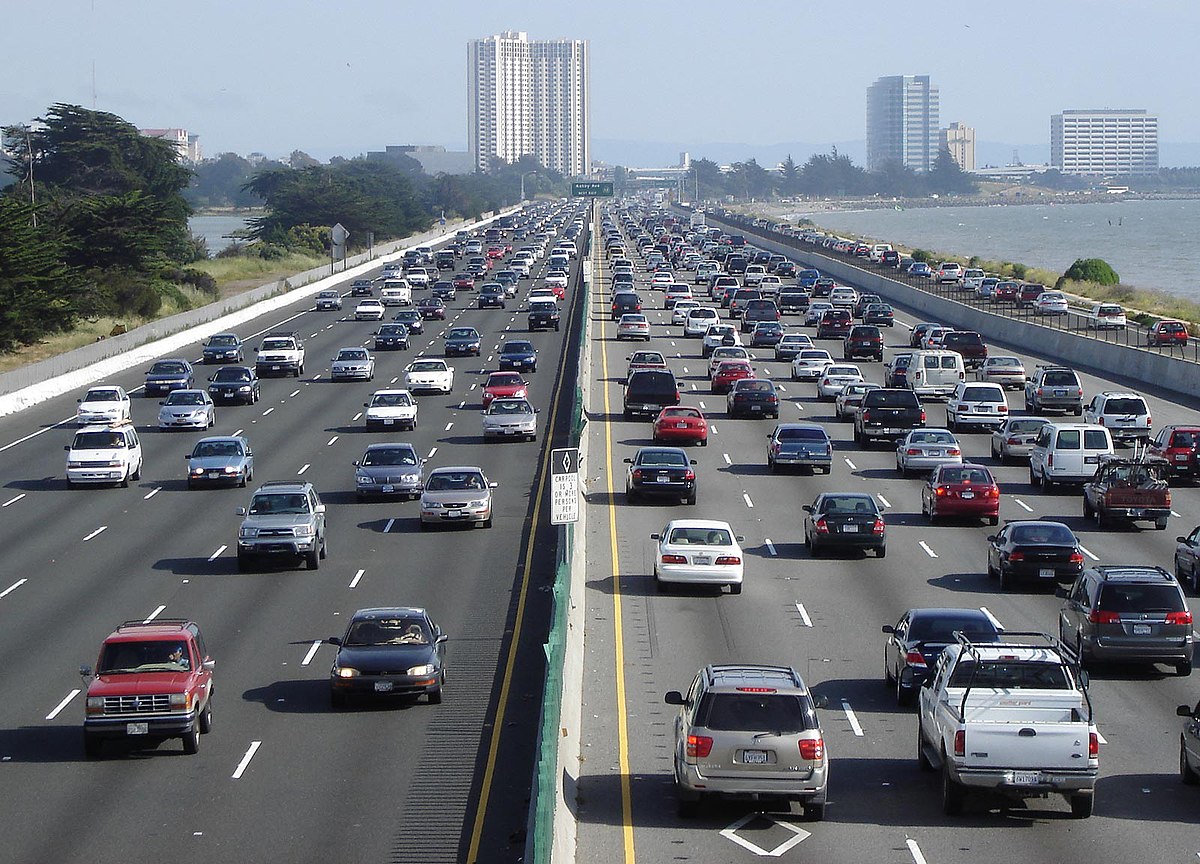}
 \caption{\sl Traffic in Rome in 1939 and in Berkeley in 2005 (left, Public Domain image from Wikipedia; right, photo by Minesweeper, image from Wikipedia,
licensed under the GNU Free Documentation License).}\label{2SINGBATM8ijr3oHNSDOe8iujrf8887uhhhgbdmnGREENFIDItanSHO4geFI}
\end{figure}

Traffic jams (see Figure~\ref{2SINGBATM8ijr3oHNSDOe8iujrf8887uhhhgbdmnGREENFIDItanSHO4geFI}) are maddening, so let's try to understand a bit about traffic flow problems and what to do in the unlucky circumstance in which we get stuck in a blockage.

For this, we consider a highway modeled by the real line, with material points (the cars) moving on it. We denote by~$\rho=\rho(x,t)$ the density of cars at a given point~$x$ at time~$t$. Since (in the absence of accidents) the total number of cars on the highway is preserved, the continuity equation~\eqref{OJS-PJDN-0IHGDOIUGDBV02ujrfMTE} prescribes that
\begin{equation}\label{BUEHRGHSNDISJ09238rujf9ihSIFTCR7UT} \partial_t\rho+\partial_x(\rho v)=0,\end{equation}
where~$v=v(x,t)$ denotes the corresponding velocity of the traffic flow at the point~$x$ at time~$t$.

To model in a somewhat realistic way how drivers react to the traffic, one may suppose that the velocity depends on the density of cars. Typically, if the highway is empty (i.e., if~$\rho = 0$), the drivers go as quick as they are allowed by the speed limit (i.e. reaching the maximal speed allowed, denoted by~$V$), but in heavy traffic situations the drivers slow down, and possibly stop in a tailback where the cars are bumper to bumper (i.e., the velocity is zero when the density reaches a maximal threshold, that we denote by~$R$).

As an example of velocity function with this property, one can consider the linear function
$$ v(x,t):= \frac{\big(R-\rho(x,t)\big) V}{R}.$$

By inserting this into~\eqref{BUEHRGHSNDISJ09238rujf9ihSIFTCR7UT}, one finds
$$ \partial_t\rho+\partial_x\left( \frac{(R-\rho) V\rho }{R}\right)=0.$$
Substituting for
$$u(x,t):=\frac{R-2\rho\left(x,\frac{t}{2V}\right)}{2R},$$ 
which corresponds to
$$ \rho(x,t)=\frac{R}{2}\left(1- 2u\left(x,2Vt\right)\right),$$ 
we find
$$\partial_t u+\partial_x\left(\frac{u^2}2\right)=0,$$ \index{Bateman-Burgers equation}
which corresponds to the inviscid Bateman-Burgers equation in~\eqref{VERBURG3}.

The fact that this equation produces shock waves (as depicted in Figure~\ref{SINGBATMGREENFIDItangeFI}) is certainly no good news for the drivers (but likely very good news for auto body repairs and insurance companies!).

The bottom line is that it's always your choice, what you do with the moment: if you're stuck in a traffic jam, you can get angry and honk your horn, or work on your mathematics.

\subsection{Reaching equilibrium. Or maybe not?}

\begin{figure}
                \centering
                \includegraphics[width=.48\linewidth]{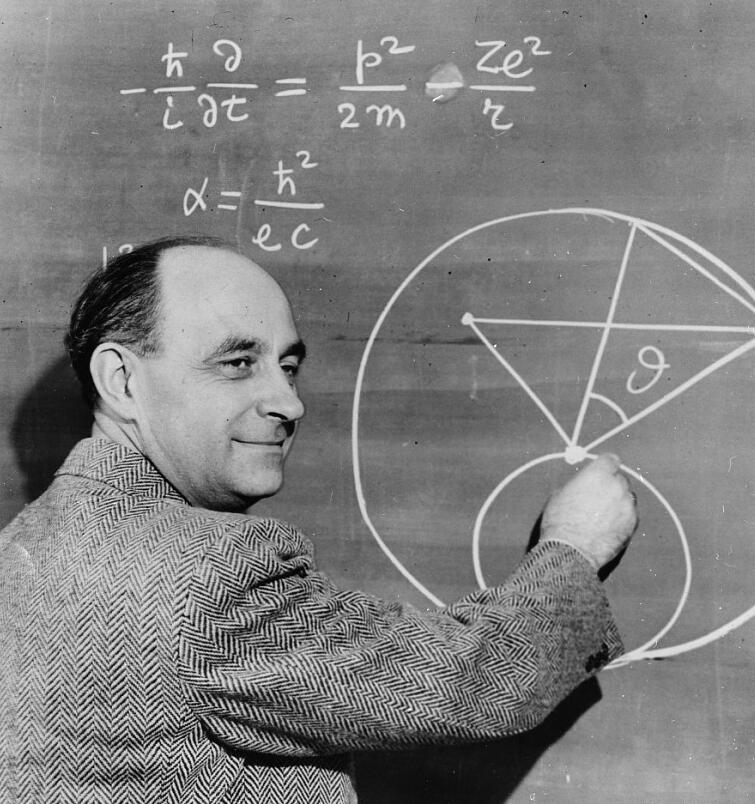}
        \caption{\sl Enrico Fermi at the blackboard (image from the Smithsonian Institution Archives).}\label{2HAFODAKEDFUMFERmNFOJED231MURA1} 
        %https://www.si.edu/object/enrico-fermi-1901-1954%3Asiris_arc_289414
\end{figure}

One can believe that, over long periods of time, a system tends to evolve naturally towards a disordered state in which
energy is shared equally among all of its various forms.
This was certainly the belief of Enrico Fermi \index{Fermi-Pasta-Ulam-Tsingou problem|(}
(according to~\cite{WRONG} and~\cite[page~128]{zbMATH06559660}), see Figure~\ref{2HAFODAKEDFUMFERmNFOJED231MURA1}.

But life is full of surprises.
Fermi left his lab located in Via Panisperna, Rome, in 1938 to escape racial laws that affected his Jewish wife, Laura Capon,
and moved to the United States. After some years, Fermi joined the secret laboratory in Los Alamos
to design and build the first atomic bomb, to which he contributed\footnote{In particular,
Fermi constructed the world's first artificial nuclear reactor
and, on the first detonation of a nuclear weapon, he calculated the energy that was going to be
released as blast
(with that, he allegedly ended up
taking bets on whether or not the whole atmosphere would ignite due to the bomb's explosion, and if so whether it would incinerate just the state or the entire planet).

Before moving to Los Alamos, Fermi was already one of the world leaders
in both theoretical physics and experimental physics. He had been
awarded the 1938 Nobel Prize in Physics for his discovery of induced radioactivity by neutron bombardment and of chemical elements with atomic numbers greater than that of uranium.
Fermi's prominent contributions range from statistical mechanics to particle physics,
pioneering the theory of weak interaction.

In~\cite{1922RendL21F}, he introduced the moving system of reference adapted to lines and surfaces
that is nowadays called ``Fermi Coordinates'', which are broadly used \label{PAGFermiCoordinates}
in geometry, analysis and mathematical physics \index{Fermi Coordinates}
(we will also utilize them on page~\pageref{FECOORD}).} very substantially. 

\begin{figure}
                \centering
                \includegraphics[width=.29\linewidth]{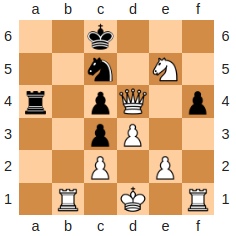}
        \caption{\sl MANIAC I beating a human player at Los Alamos chess.}\label{LOSA2HAFOULAM:URA18} 
\end{figure}

\begin{figure}
                \centering
                \includegraphics[width=.3\linewidth]{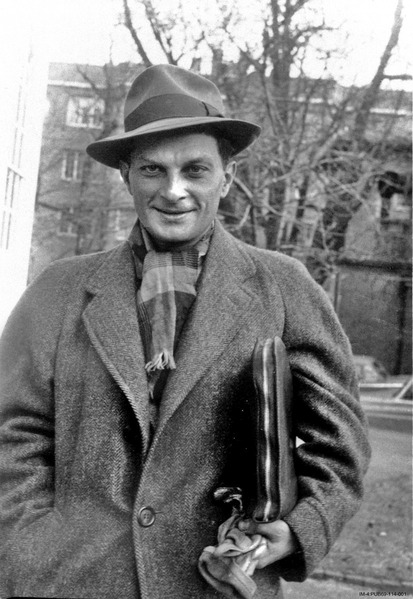}
        \caption{\sl Stanislaw Ulam (image from the Los Alamos National Laboratory).}\label{2HAFOULAM:URA18} 
\end{figure}

\begin{figure}
                \centering
                \includegraphics[width=.21\linewidth]{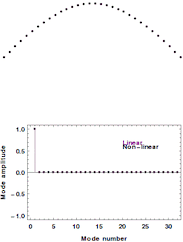}$\quad$
                                \includegraphics[width=.21\linewidth]{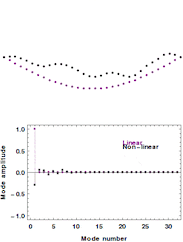}$\quad$
                                                                \includegraphics[width=.21\linewidth]{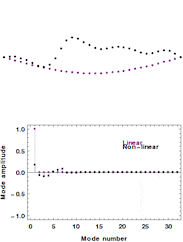}$\quad$
                                \includegraphics[width=.21\linewidth]{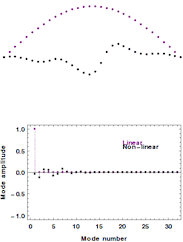}\\
                                                                                                \includegraphics[width=.21\linewidth]{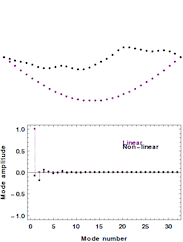}$\quad$
                                                                                                                                \includegraphics[width=.21\linewidth]{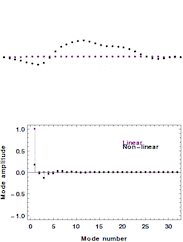}$\quad$
                                                                              \includegraphics[width=.21\linewidth]{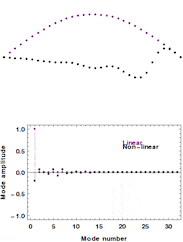}                                                                                                                  $\quad$\includegraphics[width=.22\linewidth]{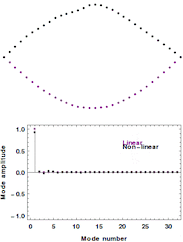}
        \caption{\sl An elastic chain with a small perturbation that initially distributes the energy among its modes, but later on shifts back almost all the energy to the original mode
        (extract from a simulation by Jacopo Bertolotti, image from
        Wikipedia, licensed under the Creative Commons CC0 1.0 Universal Public Domain Dedication).}\label{LOEXP3LAM:URA18} 
\end{figure}

Fermi remained in close contact with the lab in Los Alamos even after the end of the war
and was very interested in the development of the vacuum-tube computer\footnote{In retrospect, the computational powers of these
computers were perhaps quite limited. In 1956, MANIAC I would have become the first computer to defeat a human being in chess
(the human player was a beginner, whose name is unknown to us,
who had been taught the rules of chess only in the preceding week,
the final position after the 23rd move Ne5\# being depicted in Figure~\ref{LOSA2HAFOULAM:URA18}).
Well, not precisely chess, but a chess variant on a $6\times6$ chessboard in which bishops
were removed due to the limited amount of memory and computing power.
This variant of chess is nowadays called Los Alamos chess or, due to
the missing bishops, anti-clerical chess.
The rules are as in chess except that
there is no pawn initial double-step move, no en passant capture
and no castling (and, of course, pawns may not promote to bishops).}
MANIAC I (Mathematical Analyzer Numerical Integrator and Automatic Computer Model I),
which Fermi regarded not merely as a tool for calculating,
but rather as a useful device to explore and better understand the complicated dynamics of physical systems.
For this, Fermi had extensive discussions with
John Pasta and\footnote{Ulam was also a prominent scientific figure.
Among the other accomplishments of his career, he 
invented cellular automaton and
the Monte Carlo method. See~\cite{MR0485098}
for Ulam's autobiography.}
Stanislaw Ulam (see Figure~\ref{2HAFOULAM:URA18})
to formulate a physical problem which was simple to state
but whose solution would require a computation lengthy enough to be unfeasible by
pencil and paper (even\footnote{Fermi was famous for his ability to make exceptionally good approximate calculations.
His renowned reputation on this topic is the reason for which nowadays
problems that require brilliant applications of dimensional analysis and approximation are named ``Fermi quizzes''.} for Fermi's remarkable computing skills).
The problem finally chosen by Fermi, Pasta and Ulam was that
of some masses connected by elastic springs: this was a problem
easy to visualize and related to everyday experience (e.g., a string instrument, such as a
guitar). Fermi's expectation was that, after a short time, a nonlinear perturbation of this elastic system
would have reached an energy equipartition.
As Ulam reported (see~\cite[page~256]{zbMATH06559660})
``the original objective had been to see at what rate the energy of the string, initially put into
a single sine wave (the note was struck as one tone), would gradually develop higher tones
with the harmonics, and how the shape would finally become a mess both in the form
of the string and in the way the energy was distributed among higher and higher modes.
Nothing of the sort happened. To our surprise\footnote{The fact that Fermi's original belief about systems (almost always) evolving naturally to energy equipartition was incorrect does not diminish his genius: on the contrary, the experiment run with MANIAC I was one of the cornerstone of the modern theory of dynamical systems, in light of ergodic theory, KAM theory, chaos, etc.}
the string started playing a game of musical
chairs, only between several low notes, and perhaps even more amazingly, after what would
have been several hundred ordinary up and down vibrations, it came back almost exactly to
its original sinusoidal shape''. See Figure~\ref{LOEXP3LAM:URA18} for a sketch of such a
dumbfounding phenomenon.
\medskip

In further detail, the mathematical framework chosen by
Fermi, Pasta and Ulam was to describe the solution of the chain of oscillators in Fourier Series
and to consider the ``harmonic energy'' associated to each mode of such a decomposition.
Their original expectation was to see the system exploring the available energy surface
and relax towards an equipartition state in which, for instance, all the harmonic energies (or at least
their time average) would approach the same value. Instead, the harmonic energies seemed
to exhibit a recurrent behavior, coming back almost to the initial configuration after some excursions,
and the time average of the harmonic energies seemed to relax not towards a constant distribution
but to a distribution that looked exponentially decreasing with the modes (in particular,
the system seemed to exhibit the
formation of a ``packet of modes'' with geometrically decaying energies). These behaviors
were perceived by Fermi, Pasta and Ulam as paradoxical with respect to the equilibrium theory
of statistical mechanics.
\medskip

\begin{figure}
                \centering
                \includegraphics[width=.35\linewidth]{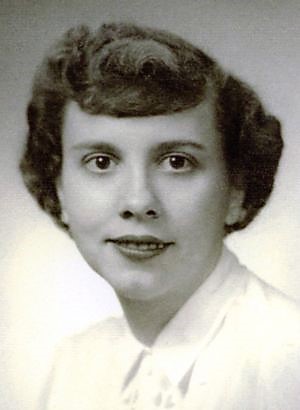}
        \caption{\sl Mary Tsingou (image from the Los Alamos National Laboratory).}\label{2HAFRtgbstsOULAM:URA18}
\end{figure}

As a hystorical remark, Fermi, Pasta and Ulam were experts in mathematics as well as theoretical, nuclear and computational physics, but not specialized programmers,
therefore they contacted one of the programmers of the MANIAC I, Mary Tsingou (see Figure~\ref{2HAFRtgbstsOULAM:URA18}) to write the code and implement the simulation: according to Tsingou,
interestingly, Fermi, Pasta and Ulam
``didn't know anything\footnote{But perhaps Tsingou's statement was also a bit ungenerous:
according to Ulam (see the preface to~\cite{FEPAU}),
``during one summer Fermi learned very rapidly how to program problems for the electronic computers and he not only could plan the general outline and construct the so-called flow diagram but would work out himself the actual coding of the whole problem in detail''.} about programming. They set up the equations, and I did all the programming''.
Even the action of programming was at the time quite different from what programmers do today:
according to Pasta (as reported on~\cite[page~359]{MR1462745})
``the program was of course punched on
cards. A DO loop was executed by the operator feeding in the deck of cards over
and over''.
\medskip

Short after this numerical experiment took place, Fermi was invited to give a
honorary lecture at the annual American Mathematical Society meeting and he
intended to talk about it. Unfortunately,
he became ill before the meeting, so his
lecture never took place. Fermi died of inoperable stomach cancer and the results of this work
were published in~\cite{FEPAU} after his death
(Enrico Fermi, John Pasta and Stanislaw Ulam being the authors of the article;
the footnote on page~979 of~\cite{FEPAU} reporting the acknowledgement
``We thank Miss Mary Tsingou for efficient coding of the problems and for running the computations on the Los Alamos MANIAC machine'').
The problem has been established in the scientific literature as the Fermi-Pasta-Ulam problem, but was renamed in 2008
as the Fermi-Pasta-Ulam-Tsingou problem
to grant attribution not only to the ideators of the mathematical physics problem but to the computer programmer as well.\medskip

Let us now see more specifically what this problem was and how (perhaps quite surprisingly)
it relates to the theory of partial differential equations (see also~\cite{MR1462745} for further details).
Given~$h>0$ and~$N\in\N$, let us put some material points of equal mass~$m$ along a line
at some rest positions located at~$jh$, with~$j\in\{0,\dots,N\}$.
Let us connect\footnote{In the original numerical experiment~\cite{FEPAU}, \label{FEconstrained}
the model focused on a case in which the first and last material points were constrained
at their rest position, to have a situation with fixed endpoints, but later on
periodic arrays of such a configuration were taken into account as well:
so, to keep things as simple as possible, we will skate over
the detail of the precise boundary conditions chosen
(a general perception in the literature is that the specific boundary condition taken
is possibly not the most relevant feature in this particular setting, but see~\cite{MR2398159}
for a precise discussion of similarities and differences).}
these material points by identical elastic springs. 

\begin{figure}
                \centering
                \includegraphics[width=.95\linewidth]{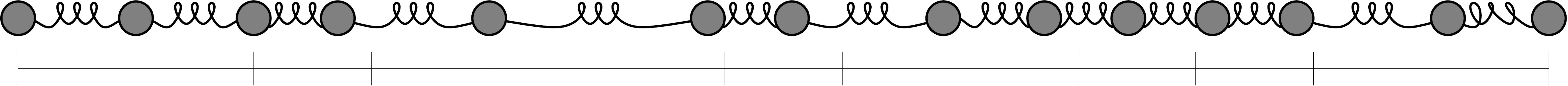}
        \caption{\sl A chain of oscillators.}\label{C24tHasdA3O2SLOSA2HAFOULAM:URA18} 
\end{figure}

This system mimics a string with restoring forces (we will see more on that in Section~\ref{GUITAR}).
If we move horizontally some of the material points from their rest position, the springs exert an elastic
force producing oscillations. We denote by~$X_j(t)$ the position of the $j$th material point at time~$t$.
The force acting on the $j$th material point is coming from the springs to its left and right sides:
assuming perfect elasticity,
since the left spring produces a force proportional to~$-(X_j-X_{j-1})$
and the right spring a force proportional to~$(X_{j+1}-X_{j})$, the total force acting on 
the $j$th material point has the form
\begin{equation}\label{EFFEj9o4PujFE} f_j:=-\kappa (X_j-X_{j-1})+\kappa(X_{j+1}-X_{j}),\end{equation}
for some elastic coefficient~$\kappa>0$, see Figure~\ref{C24tHasdA3O2SLOSA2HAFOULAM:URA18}.

Thus, we can also introduce a potential energy
\begin{equation}\label{EFFEj9o4PujFE-3}
P_j(X):=\frac{\kappa}2\left( (X_j-X_{j-1})^2+(X_{j+1}-X_{j})^2\right),\end{equation}
where~$X:=(X_0,\dots,X_N)$ and rewrite the force in~\eqref{EFFEj9o4PujFE} as
\begin{equation}\label{EFFEj9o4PujFE-2} f_j=-\partial_{X_j} P_j(X).\end{equation}

The idea proposed by Fermi is to consider a small anelastic perturbation of this ideal problem.
For simplicity, one can just add a small additional potential energy term of cubic (instead of quadratic) type,
thus replacing~\eqref{EFFEj9o4PujFE-3} by
\[
\widetilde P_j(X):=
\frac{\kappa}2\left( (X_j-X_{j-1})^2+(X_{j+1}-X_{j})^2\right)+
\frac{\e}3\left( (X_j-X_{j-1})^3+(X_{j+1}-X_{j})^3\right),
\]
where~$\e>0$ is a small parameter.

Correspondingly, using the notation~$Y_j:=X_j-X_{j-1}$,
the force in~\eqref{EFFEj9o4PujFE-2} would become
\begin{eqnarray*}
f_j&=&-\partial_{X_j} \widetilde P_j(X)\\
&=&-
\kappa\left( (X_j-X_{j-1})-(X_{j+1}-X_{j})\right)-{\e}\left( (X_j-X_{j-1})^2-(X_{j+1}-X_{j})^2\right)
\\&=& -\kappa\left( Y_j-Y_{j+1}\right)-{\e}\left( Y_j^2-Y_{j+1}^2\right)
\\&=&-\kappa\left( Y_j-Y_{j+1}\right)\left(1+\frac\e\kappa
\left( Y_j+Y_{j+1}\right)\right)
\\&=& \kappa\left( X_{j+1}+X_{j-1}-2X_j\right)\left(1+\frac\e\kappa
\left( X_{j+1}-X_{j-1}\right)\right).
\end{eqnarray*}
By Newton's Law, this gives that
$$ m\ddot{X}_j=\kappa\left( X_{j+1}+X_{j-1}-2X_j\right)\left(1+\frac\e\kappa
\left( X_{j+1}-X_{j-1}\right)\right).$$\medskip

To understand the link between these equations of motion and the theory
of partial differential equations, it is useful to choose~$\e:=\kappa h$
and define~$c:=\sqrt{\frac\kappa{m}}\,h$. In this way, we have that
\begin{equation}\label{DDhbF0oPajmd7tMijyhdpefr570-9-h} \ddot{X}_j(t)=\frac{c^2}{h^2}\left( X_{j+1}(t)+X_{j-1}(t)-2X_j(t)\right)\left(1+h
\left( X_{j+1}(t)-X_{j-1}(t)\right)\right).\end{equation}
It is also useful to introduce a smooth function~$u(x, t)$ to try to measures (in some appropriate limit sense) the
displacement at time~$t$ of the particle corresponding to the equilibrium position~$x$. 
In principle, the equilibria are discrete and of the form~$x=jh$ with~$j\in\{0,\dots,N\}$,
but the continuous limit would correspond to the asymptotics as~$h\searrow0$. Hence, the above function~$u$
is modeled such that~$u(jh,t)=X_j(t)$ (and, say, interpolated to be defined as~$u(x,t)$ for a continuum of~$x$).

We thereby consider the formal Taylor expansions as~$h\searrow0$
\begin{eqnarray*}&&
\frac{X_{j+1}(t)+X_{j-1}(t)-2X_j(t)}{h^2}\\&=&\frac{u(jh+h,t)+u(jh-h,t)-2u(jh,t)}{h^2}\\& =&
\frac{1}{h^2}\Bigg( \left(u(jh,t)+h\partial_xu(jh,t)+\frac{h^2}2\partial_x^2u(jh,t)
+\frac{h^3}6\partial_x^3u(jh,t)+\frac{h^4}{24}\partial_x^4u(jh,t)
\right)\\&&\quad\quad
+\left(u(jh,t)-h\partial_xu(jh,t)+\frac{h^2}2\partial_x^2u(jh,t)
-\frac{h^3}6\partial_x^3u(jh,t)+\frac{h^4}{24}\partial_x^4u(jh,t)
\right)\\&&\quad\quad-2u(jh,t)+o(h^4)
\Bigg)\\
&=&\partial_x^2u(jh,t)+\frac{h^2}{12}\partial_x^4u(jh,t)
+o(h^2)
\end{eqnarray*}
and
\begin{eqnarray*}&&
X_{j+1}(t)-X_{j-1}(t)\\&=&u(jh+h,t)-u(jh-h,t)\\&=&
\big(u(jh,t)+h\partial_xu(jh,t)\big)
-
\big(u(jh,t)-h\partial_xu(jh,t)\big)+o(h)\\&=&
2h\partial_xu(jh,t)+o(h).
\end{eqnarray*}
Hence, since also~$ \ddot X_j(t)=\partial_t^2 u(jh,t)$,
we formally rewrite~\eqref{DDhbF0oPajmd7tMijyhdpefr570-9-h} as
\begin{equation*} \begin{split}\partial_t^2 u(jh,t)\,&= c^2\left( 
\partial_x^2u(jh,t)+\frac{h^2}{12}\partial_x^4u(jh,t)
+o(h^2)\right)\left(1+h
\Big( 2h\partial_xu(jh,t)+o(h)\Big)\right)\\
&=
c^2\left( \partial_x^2u(jh,t)+2h^2\partial_xu(jh,t)\partial_x^2u(jh,t)+\frac{h^2}{12}\partial_x^4u(jh,t)
+o(h^2)\right).
\end{split}\end{equation*}
Formally writing~$x$ in the place of~$jh$ we thus obtain\footnote{Interestingly, if we sent~$h\searrow0$ in~\eqref{SUERF:PAHNDFGE-2},
we would just obtain the vibrating string equation~$\partial_t^2 u=c^2\partial_x^2 u$,
to be compared with~\eqref{CUSTRJMSONGJMSD-3}.
This is, in a sense, not too surprising, since sending~$h\searrow0$ ``too early''
would correspond to ``neglect the anelastic perturbation'', since~$\e$ is proportional to~$h$.
The finer analysis that we perform to get~\eqref{4f-u9oikhe002d}
is therefore necessary to account for the anelastic term in the continuous limit.}
\begin{equation}\label{SUERF:PAHNDFGE-2}
\frac1{c^2}\partial_t^2 u(x,t)-\partial_x^2u(x,t)=
2h^2\partial_xu(x,t)\partial_x^2u(x,t)+\frac{h^2}{12}\partial_x^4u(x,t)+o(h^2).
\end{equation}

Defining~$U(\xi,\tau):=u\left(\xi+\frac{\tau}{h^2},\frac{\tau}{ch^2}\right)$, it follows that\footnote{One can recall that this substitution
is similar to the one utilized to ``surf the wave'' on page~\pageref{SUERF:PAHNDFGE}.}
$$ u(x,t)=U(\xi,\tau),\quad{\mbox{ where $\xi=x-ct$ and $\tau=ch^2 t$}}$$
and therefore, omitting the variables for short, for each~$m\in\N$,
\begin{eqnarray*}
&&\partial_x^{m}u=\partial_\xi^m U,\\
&&\partial_{t} u=-c\partial_\xi U+ch^2\partial_\tau u\\
{\mbox{and }}&&
\partial_{t}^2 u=c^2\partial_\xi^2 U+c^2h^4\partial_{\tau}^2 U-2c^2h^2\partial_{\xi\tau}U.
\end{eqnarray*}
In particular, 
$$ \frac1{c^2}\partial_t^2 u-\partial_x^2u=h^4\partial_{\tau}^2 U-2h^2\partial_{\xi\tau}U
$$
and~\eqref{SUERF:PAHNDFGE-2} becomes
\begin{equation*}
h^4\partial_{\tau}^2 U-2h^2\partial_{\xi\tau}U=
2h^2\partial_\xi U\partial_\xi^2 U+\frac{h^2}{12}\partial_\xi^4 U+o(h^2)
\end{equation*}
and thus, dividing by~$2h^2$,
\begin{equation}\label{4f-u9oikhe002d}
-\partial_{\xi\tau}U=
\partial_\xi U\partial_\xi^2 U+\frac{1}{24}\partial_\xi^4 U+o(1).
\end{equation}
Setting~$v:=\partial_\xi U$ and formally taking the limit as~$h\searrow0$ we find
\begin{equation}\label{INCarFVU023iekrjmfsuJMS}
-\partial_{\tau}v=
v\partial_\xi v+\frac{1}{24}\partial_\xi^3v,
\end{equation} \index{Korteweg-de Vries equation|(}
which corresponds to the Korteweg-de Vries equation in~\eqref{KDV-or}
(say, with~$a:=1$, $b:=1/24$ and~$c:=1$, the numerology being unimportant here).

The link between the Fermi-Pasta-Ulam-Tsingou problem and the Korteweg-de Vries equation
is not a mere curiosity based on formal (though a bit sophisticated) Taylor expansions: instead, it has a deep impact on the understanding
of the surprising phenomenon highlighted by the oscillators ending up relocating almost all the energy
on the initial mode, as discovered\footnote{The name ``soliton'' \index{soliton}
was indeed coniated in~\cite{zbMATH05826808}
to describe the ``solitary waves'' which persist after
nondestructive interactions. Likely, the name ``soliton''
is meant to stress the similarities of the behavior of such waves
with that of interacting ``particles''
(since elementary particles have names which end in ``-on'', such as
``electron'', ``proton'', ``neutron'', ``lepton'', ``boson'', etc.,
why not naming ``solitons'' these
wave packets that maintain their shape while propagating and interacting?).}
in~\cite{zbMATH05826808}.
Roughly speaking, the existence of solitons in the Korteweg-de Vries equation
which essentially recover their initial shape after interactions (as observed by~\cite{zbMATH05826808})
is considered the counterpart in the framework of partial differential
equations of the near recurrence to the initial state in discretized weakly-nonlinear strings.
The topic is very complex and still under an intense scientific investigation,
hence we do not aim at being exhaustive here.
Yet, in a sense, the description of the solutions of the Korteweg-de Vries equation~\eqref{INCarFVU023iekrjmfsuJMS}
explored in~\cite{zbMATH05826808}
goes as follows. At the beginning the term~$\partial_{\tau}v+
v\partial_\xi v$
in~\eqref{INCarFVU023iekrjmfsuJMS} ``dominates'' and the solution steepens in regions where it has a negative slope (this is indeed the typical behavior of solutions of 
the Bateman-Burgers equation~\eqref{VERBURG3}, \index{Bateman-Burgers equation}
recall Figure~\ref{SINGBATMGREENFIDItangeFI}). But, as time flows, the third derivative term~$\partial_\xi^3v$ in equation~\eqref{INCarFVU023iekrjmfsuJMS} takes over\footnote{The importance
of the third order term in the Korteweg-de Vries equation and its smoothing \label{GITARFKMAOSTGABafOKS02o3t-2}
effect also appeared in the discussion on page~\pageref{GITARFKMAOSTGABafOKS02o3t}.}
and helps prevent the formation of a discontinuities, producing instead some oscillations on the left of the front. Later on, the amplitudes of these oscillations grow till each oscillation achieves an almost steady amplitude (increasing from left to right).
These bumps begin to move with a speed proportional to its amplitude (hence the bumps on the right move faster than the ones on the left, causing these bumps to spread apart).
Since the domain is bounded\footnote{Strictly speaking, periodic boundary conditions were considered in~\cite{zbMATH05826808}, rather than fixed endpoint in~\cite{FEPAU},
but we are glossing over this detail, see footnote~\ref{FEconstrained}.} these bumps will eventually meet again and interact. But quite surprisingly, shortly after the interaction, the bumps reappear essentially with the same shape.

This is a remarkable property of these solitons which are able to ``pass through
each other'' by maintaining their silhouette.
The analogy between the behavior of these solitons and the almost return to initial condition
detected in the Fermi-Pasta-Ulam-Tsingou setting is that,
from time to time, the solitons come back to the positions they had initially, hence
almost restoring the initial condition.

The dynamics of the solitons is clearly explained in~\cite{zbMATH05826808} with the following wording:
``When solitons of very different amplitude approach, their trajectories deviate from straight lines (accelerate) as they "pass through" one another. During the overlap time interval the joint amplitude of the interacting solitons decreases (in contradistinction to what would happen if two pulses overlapped linearly). [...] When the amplitudes of approaching solitons are comparable
they seem to exchange amplitudes and therefore velocities. [...] At [a certain] time~$T_R$ all the solitons arrive {\em almost} in the same phase and almost reconstruct the initial state through the nonlinear interaction. This process
proceeds onwards, and at~$2T_R$ one again has a ``near recurrence'' which is not as good as the first recurrence. [...] Because the solitons are remarkably stable entities, preserving their identity through numerous interactions, one would expect this system to exhibit thermalization (complete energy sharing among the corresponding linear normal modes) only after extremely long times, if ever''.

For a clear picture of these solitons see Figure~1 in~\cite{zbMATH05826808}.\medskip

The investigation of the Fermi-Pasta-Ulam-Tsingou problem, as well as its continuous
counterparts such as the Korteweg-de Vries equation, is still very active and, in spite of remarkable
progresses, a definite conclusion has not been reached till now.
A mostly unanimous belief among experts is that the Fermi-Pasta-Ulam-Tsingou problem
develops chaotic features for essentially all the initial data, but the time required 
to observe these phenomena may overcome the lifetime of a physical system
and significantly challenge the reliability of the numerical simulations.
In any case, it was experimentally detected in~\cite{zbMATH03241199}
that the recurrent phenomena appeared in the original Fermi-Pasta-Ulam-Tsingou problem
somewhat disappear if initial data with sufficiently high energy are considered: in this case,
it seems that the system quickly relaxes to energy equipartition
and the time averages of all the harmonic energies approach the same value
(compare e.g. Figures~3 and~4 on~\cite[page~240]{MR3445504}
to appreciate the effect of increasing the initial energy per oscillator).
Also, some numerical data seem to indicate that, even
when energy is initially given to a small subset of modes of low frequency, for very long times the system does tend towards a
statistical equilibrium, identified by energy equipartition, see e.g.~\cite{MR1192748, BERCHIALLA2004167, MR2826618, MR4146721}.

In this perspective, the state detected in the original experiment~\cite{FEPAU},
in which, at low energy, a packet of modes, in addition to the initially excited one, enters the game,
possibly returning close to the initial condition,
is only apparently stationary and should be instead considered only
a ``metastable state'', since on a much longer time scale\footnote{For completeness, we also mention
that the features of the chains of oscillators in dimension higher
than one are quite different from those of the
Fermi-Pasta-Ulam-Tsingou problem
and, for instance, already in dimension two
the times necessary to numerically observe states of equilibrium are substantially shorter than in dimension one, see e.g.~\cite{MR2133459}.

The dependence of chaotic behaviors upon the number of dimensions
is also a general feature of the modern theory of dynamical systems, since in low dimension
some systems possess a ``topological obstruction'' to chaos, see e.g.~\cite{MR2238867}.}
the system should evolve towards statistical equilibrium
and energy equipartition. Yet, complete mathematical proofs for these phenomena seem to be missing.

For perturbative approaches to the
Fermi-Pasta-Ulam-Tsingou problem (with periodic conditions and a small energy assumption
depending on the number of oscillators) see~\cite{MR1831098, MR2367202}:
in this setting, for small energies, one rigorously obtains many recurrent solutions of the problem of
quasi-periodic type.

Another important question related to the Fermi-Pasta-Ulam-Tsingou problem
deals with the ``thermodynamic limit'' \index{thermodynamic limit}
in which
the number of oscillators~$N$ in the chain tends to infinity
while keeping the total energy proportional to~$N$. This would be perhaps the most significant
question on the topic in relation to statistical mechanics, whose 
foundations heavily depend on the understanding of systems with an arbitrarily large number of interacting particles.
This problem is also largely open: for questions of this kind,
as stated in~\cite[Section~3]{MR3445504},
``of course numerics can
just give some indications, while a definitive result can only come from a theoretical
result, which is the only one able to attain the limit~$N=\infty$''.
In this context,
the results related to formation of the packet and their persistence in chains of~$N$ oscillators
are typically confined to the regime in which the total energy is of order~$N^{-3}$.
Though this limitation is unsatisfactory for the physical applications
(confining the study to that of large oscillators with extremely low energy),
it is somewhat remarkable that different theoretical approaches lead to an assumption
about the same energy regime (on the other hand, ``numerics do not provide any evidence of
changes in the dynamics when energy is increased beyond this limit'',
see~\cite[page~251]{MR3445504}).

The rigorous link between the Fermi-Pasta-Ulam-Tsingou problem
and the Korteweg-de Vries equation is also a delicate issue.
One of the structural problems in this approximation is that the Fermi-Pasta-Ulam-Tsingou system
turns out to be a ``singular perturbation'' of the Korteweg-de Vries equation, in the sense that
the approximation error contains higher order derivatives. See~\cite[Section~4.1]{MR3445504}
for rigorous results about this asymptotic problem. 

The influence of the theory of solitons on the rigorous study of the Fermi-Pasta-Ulam-Tsingou problem
is showcased in the series of papers~\cite{MR1726667, MR1912298, MR2023440, MR2023441}.
See also~\cite{MR2402016} for an account of the various facets of the
Fermi-Pasta-Ulam-Tsingou problem and~\cite{MR2133458, MR4294860}
for rigorous links between the Fermi-Pasta-Ulam-Tsingou problem
and the Korteweg-de Vries equation.

But let us also stress a significant
difference between the dynamics of the Fermi-Pasta-Ulam-Tsingou system and the one of 
the Korteweg-de Vries equation: while the recurrent states are nowadays expected to be a rare or transient phenomenon
in chain of oscillators,
the Korteweg-de Vries equation is a ``completely integrable'' system (see e.g.~\cite[pages~277--284]{MR3470076}) and
all its solutions (under appropriate boundary conditions)
are periodic, quasi-periodic, or almost-periodic in time, see~\cite[pages~ix, 1--5]{MR1997070}.
Hence, differently from what we should expect for the Fermi-Pasta-Ulam-Tsingou problem,
solutions of the Korteweg-de Vries equation return infinitely often arbitrarily close to their initial state. \index{Fermi-Pasta-Ulam-Tsingou problem|)}\index{Korteweg-de Vries equation|)}

\subsection{Definitions come later on}

\begin{figure}
  \centering
  \includegraphics[width=.3\linewidth]{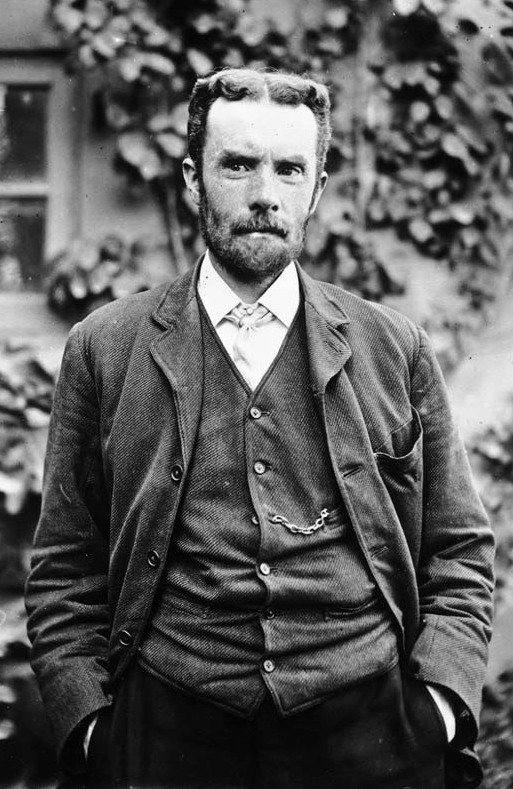}
 \caption{\sl Oliver Heaviside (Public Domain image from
 Wikipedia).}\label{NIRDG0-HEA02-43X1qwrt5t43REENFIDItangeFI}
\end{figure}

Given~$a$, $b$, $c\ge0$, the equation
\begin{equation}\label{HEATELEEQ0}
\frac{\partial ^{2} u}{\partial x^{2}}-a\frac{\partial ^{2} u}{\ \partial t^{2}}=
b \frac{\partial u }{\partial t} +c u
\end{equation}
is called\footnote{The inventor of equation~\eqref{HEATELEEQ0}
was Oliver Heaviside, see Figure~\ref{NIRDG0-HEA02-43X1qwrt5t43REENFIDItangeFI}.

Heaviside's uncle was Sir Charles Wheatstone, co-inventor of the first commercially successful telegraph. Following his uncle's advice, Heaviside became a telegraph operator and an electrician, continuing to study and do science while working.

His contributions to mathematics, physics and engineering were deep and remained as classical contributions: they include the use of complex numbers to solve circuit analysis differential equations, the calculus of the deformations of an electromagnetic field surrounding a moving charge and the derivation of the magnetic force on a moving charged particle (the first contribution towards the understanding of the Lorentz Force) and the prediction of the existence of a reflective layer of the ionosphere (which allowed radio waves radiated into the sky to return to Earth beyond the horizon).

Several works by Heaviside
were of great practical use, including the possible exploitation of loading coils in telephone and telegraph lines to increase their self-induction and correct the distortion which they suffered.

Several years later, some American telecommunications companies hired their own scientists to extend Heaviside's work and adapt the use of coils previously introduced by Heaviside. Some corporations later offered Heaviside money in exchange for his rights, but he declined to accept any money unless the company were to give him full recognition for his discoveries and inventions (it turns out that Heaviside remained for all his life chronically poor).

Heaviside was perhaps a bit eccentric too. He was a firm opponent of Einstein's theory of relativity, possibly beyond reasonable scientific arguments, and at the end of his life he developed a very strong aversion to meeting people and became a recluse. He died at age 74 after falling from a ladder.

See~\cite{MR939170} for a very captivating and detailed biography of Oliver Heaviside.

Heaviside said: ``Mathematics is an experimental science, and definitions do not come first, but later on. They make themselves, when the nature of the subject has developed itself''.}
the telegrapher's equation. \index{telegrapher's equation}

In equation~\eqref{HEATELEEQ0}, we have that~$u=u(x,t)$ with~$x\in\R$, $t\in[0,+\infty)$.
When~$a:=0$, $c:=0$ and~$b>0$, equation~\eqref{HEATELEEQ0} boils down to
the heat equation in~\eqref{DAGB-ADkrVoiweLL4re2346ytmngrrUj7}.
When~$b:=0$, $c:=0$ and~$a>0$, it reduces\footnote{With respect to the terminology
in footnote~\ref{CLASSIFICATIONFOOTN} on page~\pageref{CLASSIFICATIONFOOTN},
we note that equation~\eqref{HEATELEEQ0} is hyperbolic when~$a>0$, parabolic when~$a:=0<b$
and elliptic when~$a:=0$ and~$b:=0$.}
to the wave equation
in~\eqref{WAYEBJJD121t416JH098327uyrhdhc832nbcM}.

\begin{figure}
  \centering
  \includegraphics[width=.18\linewidth]{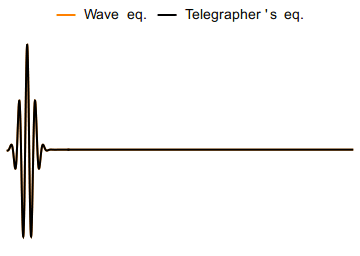}$\;$
    \includegraphics[width=.18\linewidth]{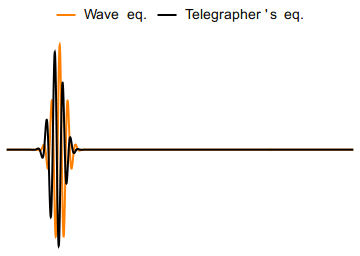}$\;$
  \includegraphics[width=.18\linewidth]{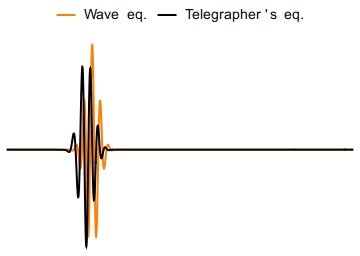}$\;$
  \includegraphics[width=.18\linewidth]{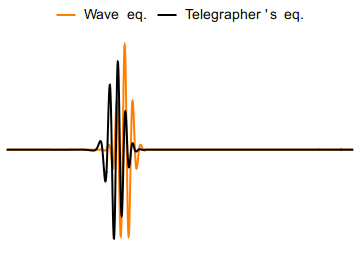}$\;$
  \includegraphics[width=.18\linewidth]{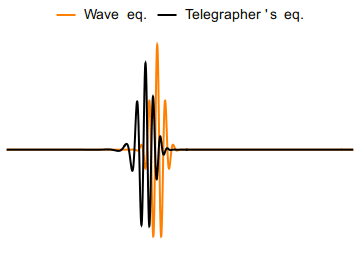}\\
  \includegraphics[width=.18\linewidth]{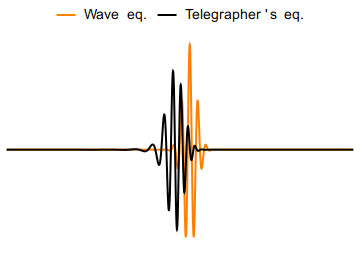}$\;$
  \includegraphics[width=.18\linewidth]{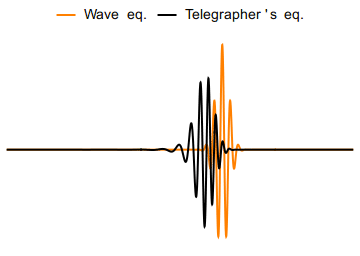}$\;$
  \includegraphics[width=.18\linewidth]{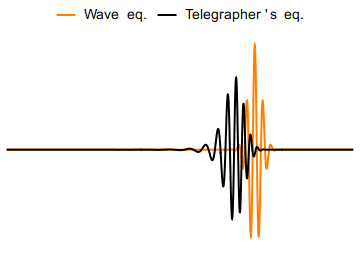}$\;$
  \includegraphics[width=.18\linewidth]{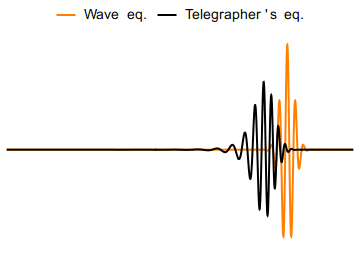}$\;$
  \includegraphics[width=.18\linewidth]{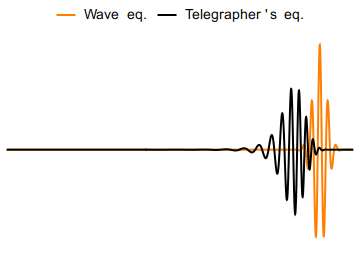}
 \caption{\sl Comparison between the evolution in time of the solution of the wave equation and of the
 telegrapher's equation (author
 Jacopo Bertolotti, {\tt https://twitter.com/j\_bertolotti/status/1172517281374572551},
 images from
 Wikipedia, licensed under the Creative Commons CC0 1.0 Universal Public Domain Dedication).}\label{wP2DAILAPIJAM STELEHGREEMIItange3FIlAsAL13458gANSE}
\end{figure}

When~$a$, $b$, $c>0$ however the behavior of the solutions of~\eqref{HEATELEEQ0} are interestingly different from
those of the wave equation,
since solutions of the telegrapher's equation present damping and dispersion effects
that are not present in the wave equation.
See e.g. Figure~\ref{wP2DAILAPIJAM STELEHGREEMIItange3FIlAsAL13458gANSE}
in which the evolution of the solution of the telegrapher's equation 
exhibits a visible \index{dispersion}
dispersion, in which the velocity of travel depends on the frequency thus
enlarging the support of the solution, and losses of intensity, causing\footnote{These dispersive effects, in which the speed of a signal is not constant but depends on \label{ihkbdDAvbsdewI089JHGSFDCokwfIUyb-X3-DISm} frequency, are typically not good for applications: for instance, when some frequencies travel along long transmission lines at a higher velocity than others, dispersion may cause the signals to be unrecognizable. More about this on page~\pageref{KSM0wodje-3ijgHHANSVV123rtgiowiejGBSna0oerL0}.}
the peaks of the traveling front to reduce over time.

\begin{figure}
  \centering
  \includegraphics[height=.22\textheight]{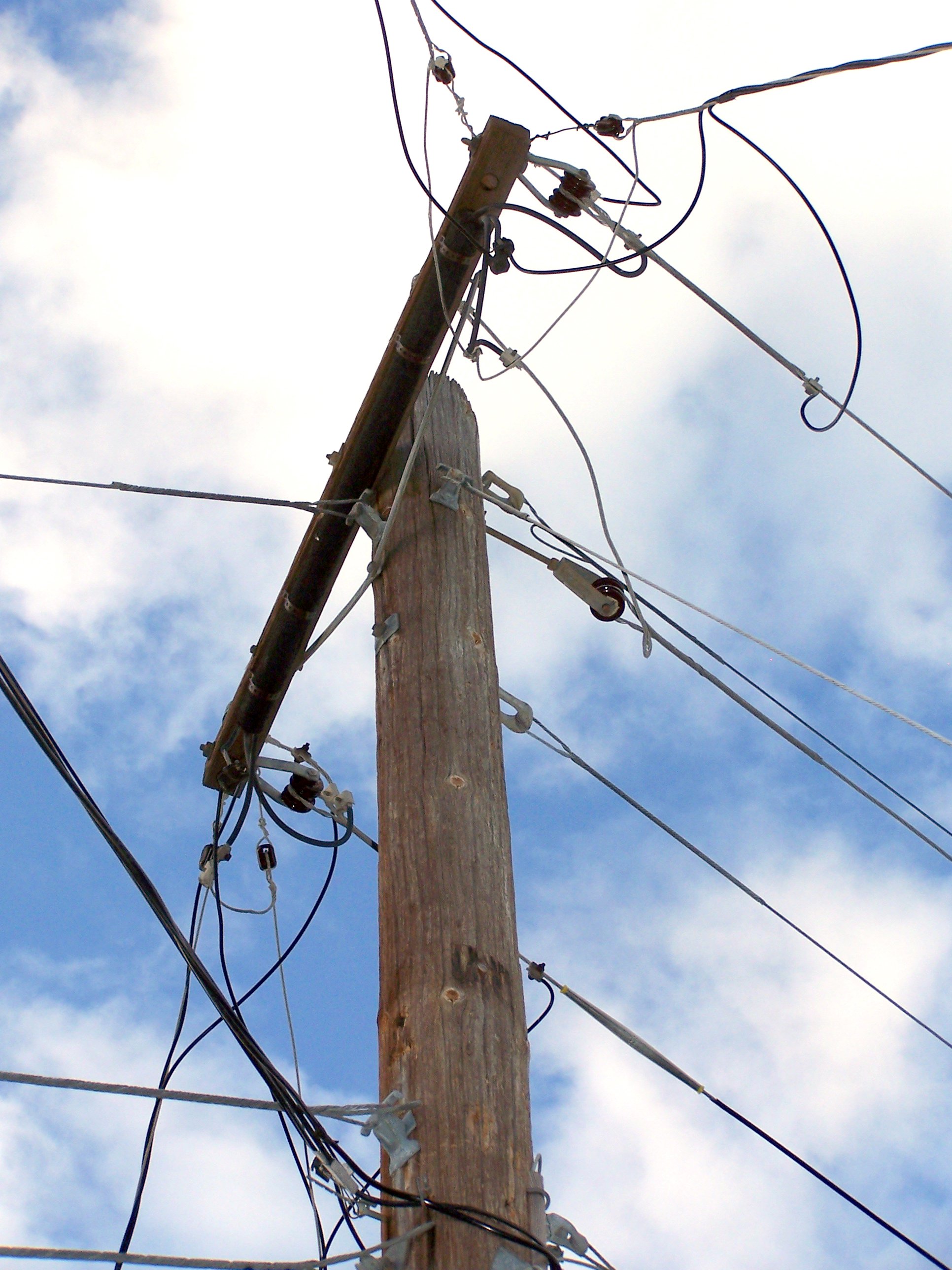} $\quad$
  \includegraphics[height=.22\textheight]{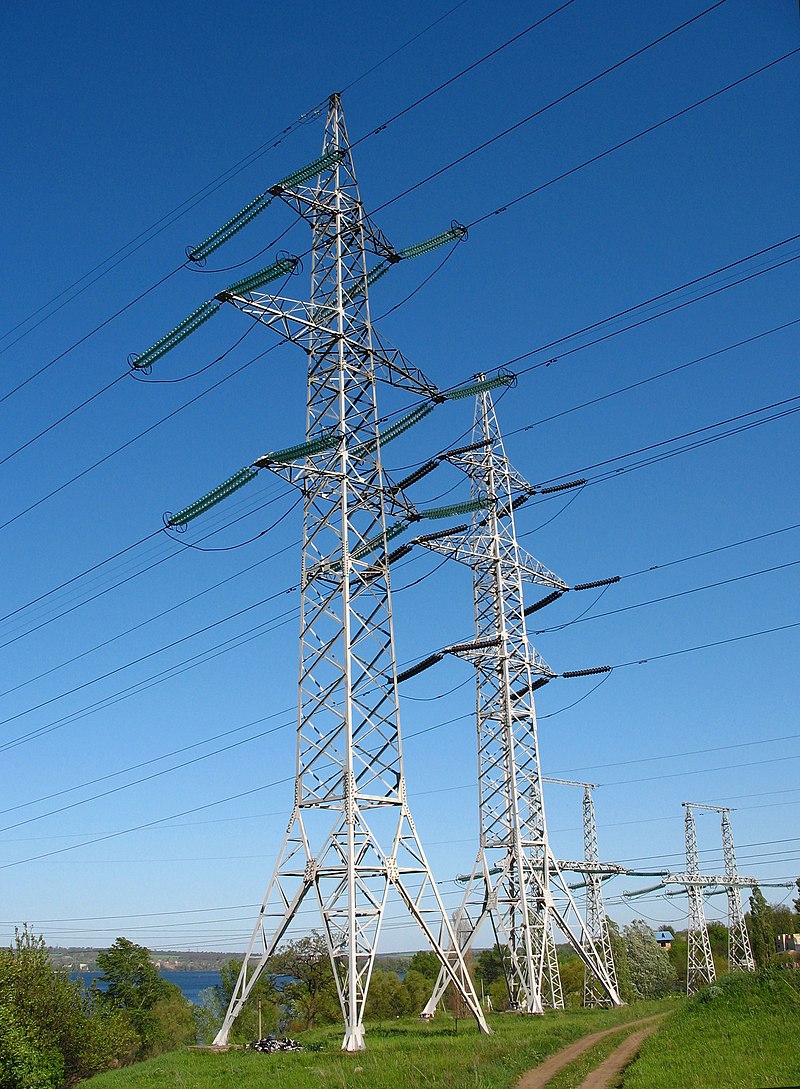} $\quad$
  \includegraphics[height=.22\textheight]{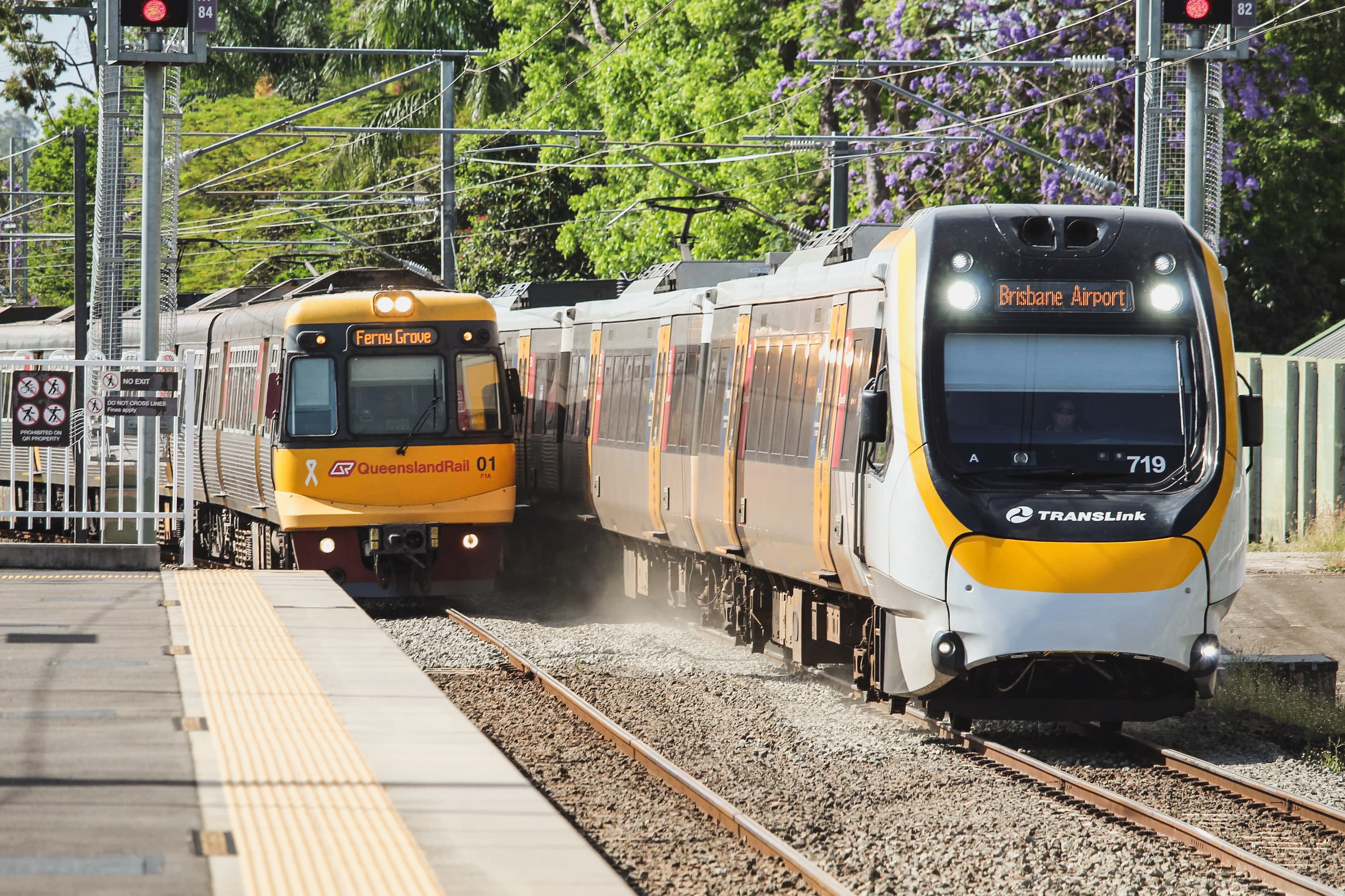}
 \caption{\sl Left: utility pole for a telephone line; center: electrical wires; 
 right: overhead lines in Queensland
 (images from
 Wikipedia, Public Domain for the first, photo by Novoklimov, licensed under the Creative Commons Attribution-Share Alike 4.0 International, licensed under the Unported license
 for the second, photo by Pytomelon87 licensed under the Creative Commons Attribution-Share Alike 4.0 International license for the third).}\label{NI8imWIRERDGREENFIDItangeFI}
\end{figure}

\begin{figure}
  \centering
  \includegraphics[height=.42\textheight]{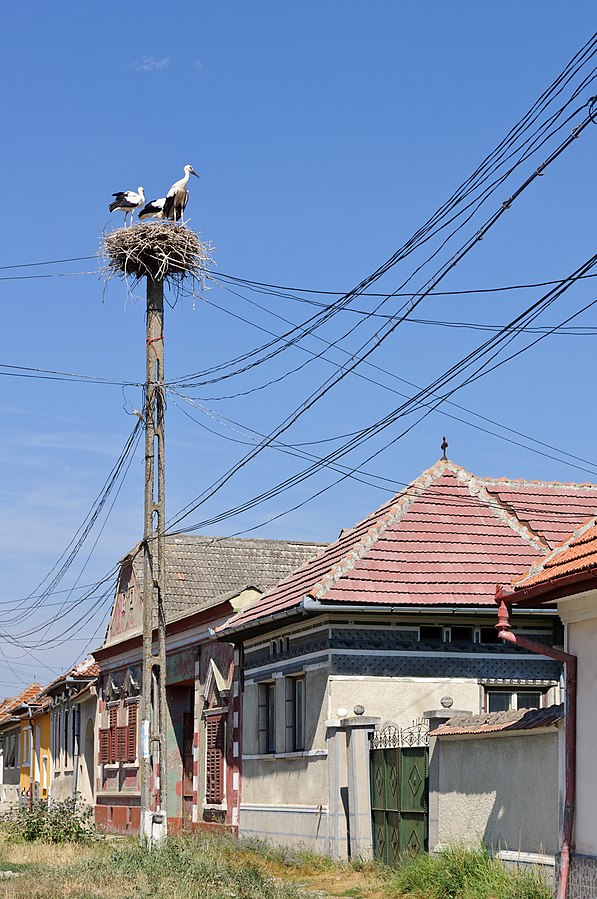}
 \caption{\sl White storks nesting on an utility pole
 (photo by Myrabella; image from
 Wikipedia, licensed under the 
 Creative Commons Attribution-Share Alike 3.0 Unported license).}\label{CICONI8imWIRERDGREENFIDItangeFI}
\end{figure}

The solution~$u$ in equation~\eqref{HEATELEEQ0} describes the electric current intensity (or, in a similar manner,
the voltage) along a transmission line.
This transmission line can be a telegraph wire, a overhead electrical conductor, a telephone line, etc.,
see Figure~\ref{NI8imWIRERDGREENFIDItangeFI}. See also Figure~\ref{CICONI8imWIRERDGREENFIDItangeFI}
for a peculiar use of transmission lines.

To get to the bottom of the telegrapher's equation~\eqref{HEATELEEQ0}, we consider
an electric line made of two parallel electrical wires.
There is a voltage difference between the wires (e.g., produced by an electric generator
``at the end of the wires'': since the wires in our idealized model are straight lines, the generator is essentially
``at infinity'').

Thinking at the two wires just as parallel straight lines would be however reductive,
since:
\begin{itemize}
\item each straight line is actually made by a conductor which presents some
distributed electric resistance (say, modeled by a series resistor) and also some
distributed inductance (e.g., due to the magnetic field around the wires and to self-inductance, modeled by a series inductor),
\item furthermore the dielectric material separating the two conductors is also not completely neutral from
the electric point of view, since it can carry
capacitance and conductance (modeled by
a shunt capacitor and resistor located between each ``infinitesimal'' portion~$dx$ of the conductors).
\end{itemize}

To describe these phenomena,
as customary in the physical description of electric phenomena, we 
reserve the name~$I$ for the current, $V$ for the voltage, $R$ for the 
resistance, $L$ for the inductance and~$C$ for the capacity. We also use the letter~$G$ to 
denote the conductance, which is the reciprocal of a resistance.
The parameters~$R$, $L$, $C$ and~$G$ will be treated as structural constants
and we will be interested in the description of the functions~$I=I(x,t)$
and~$V=V(x,t)$ with respect to the position~$x\in\R$ on the transmission line
and the time~$t\ge0$.

\begin{figure}
  \centering
  \includegraphics[width=.42\linewidth]{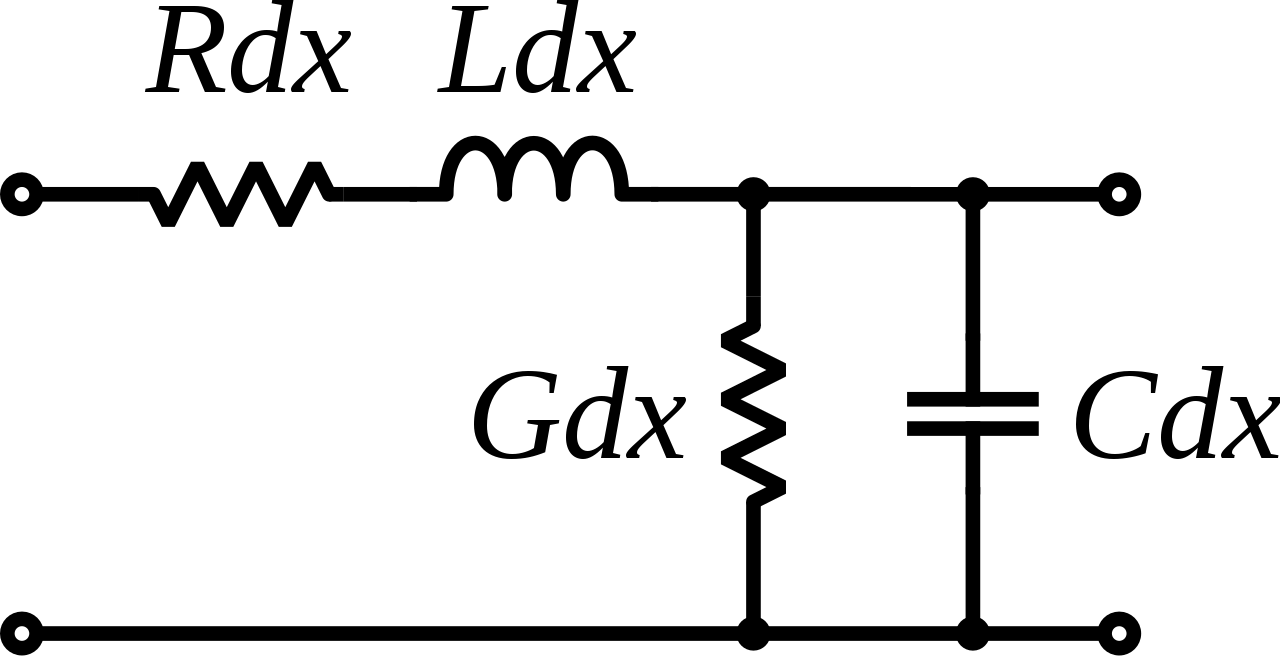}
 \caption{\sl Schematic for an elemental length~$dx$ of transmission line (image by Omegatron from
 Wikipedia, licensed under the Creative Commons Attribution-Share Alike 3.0 Unported license).}\label{P2DAILAGREEMIItangeLINE-3pe3FIlAALANSE}
\end{figure}

More specifically, for concreteness we focus on the current in the upper conductor, see Figure~\ref{P2DAILAGREEMIItangeLINE-3pe3FIlAALANSE}: in this setting,
the distributed resistance in the infinitesimal elemental length~$dx$ of the conductor
is denoted by~$R\,dx$
and the distributed inductance is denoted by~$L\,dx$.

The conductance of the dielectric material separating the two conductors
(accounting for bulk conductivity of the dielectric and dielectric loss) is denoted by~$G\,dx$
and the capacity by~$C\,dx$ (notice that, in ``real life'' there is no wire connecting the top and the bottom
cable, but Figure~\ref{P2DAILAGREEMIItangeLINE-3pe3FIlAALANSE} translates the behavior of the dielectric
between the two conductors into the language of electric circuits).

\begin{figure}
  \centering
  \includegraphics[width=.75\linewidth]{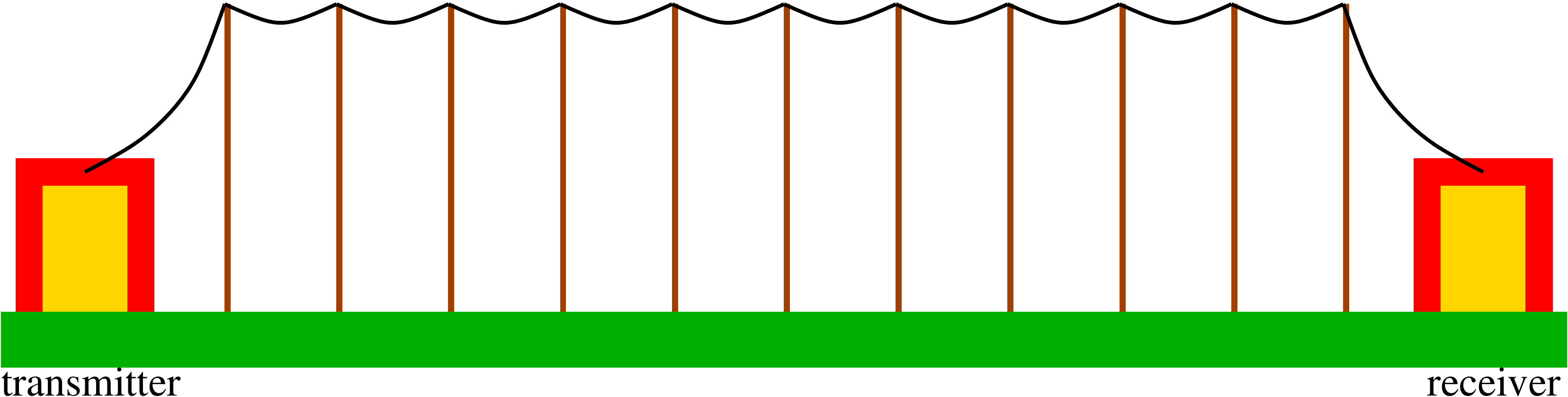}
 \caption{\sl Transmission line for a telegraph, using Earth as return path.}\label{P2DAILAGREEMIItangeLINE-3pe3FIlAALANSE-2}
\end{figure}

Now, to describe the current flowing through the upper wire, we denote by~$I(x,t)$ the current intensity on the left end
of the elementary upper conductor in Figure~\ref{P2DAILAGREEMIItangeLINE-3pe3FIlAALANSE}
and by~$I(x+dx,t)$ the one on the right end (say, with the convention that
the current is traveling from left to right).
Similarly, we denote by~$V(x,t)$ the voltage on the left end
of the elementary upper conductor and by~$V(x+dx,t)$ the one on the right end
(actually, $V$ would stand for the difference of voltage from the upper and the lower conductors;
for simplicity\footnote{Actually, telegraphs originally used two wires, utilizing a forward and a return paths
in a closed circuit to move energy along the transmission line.
But it was then noticed that Earth itself can be used as the return path. In this setting, the bottom wire is just the ground, which is normalized to be approximatively at constant zero volt (from practical purposes however, Earth is an adequate, but certainly not optimal, return path). Thus, in the situation in which the transmission line is modeled by one forward cable using Earth as a return path, the poles are also considered as dielectric and provide some capacitance and conductance, see Figure~\ref{P2DAILAGREEMIItangeLINE-3pe3FIlAALANSE-2}. In this case,
the vertical elements related to~$G\,dx$ and~$C\,dx$ in Figure~\ref{P2DAILAGREEMIItangeLINE-3pe3FIlAALANSE}
can be thought as ``concrete objects'', such as the telegraph poles, located at an ``infinitesimal'' distance~$dx$
at a large scale (in which the transmitter and the receiver are located ``at infinity'').

See e.g. {\tt https://youtu.be/ySuUZEjARPY} for a very thorough video about electric transmission.} one can think that the conductor at the bottom is at zero voltage).
Let also~$\iota_1$ be the current between the nodes of the upper conductor (oriented left to right),
and~$\iota_2$ and~$\iota_3$ be the current through the 
shunt capacitor~$C\,dx$ and resistor~$G\,dx$ respectively (oriented downwards).

With this notation, by Kirchhoff's Junction Rule (or simply by conservation of charge), we have that
$$ I(x,t)=\iota_1+\iota_3\qquad{\mbox{and}}\qquad
\iota_1=I(x+dx,t)+\iota_2.$$
Also, by the Laws of Ohm and Faraday,
$$ V(x,t)=\frac{\iota_3}{G\,dx}\qquad{\mbox{and}}\qquad
V(x+dx,t)-V(x,t)=-RI(x,t)\,dx-L\partial_t I(x,t)\,dx.$$
Additionally, by the definition of capacity,
$$\partial_tV(x,t)=\frac{\iota_2}{C\,dx}.$$

These considerations lead to
\begin{eqnarray*}
\partial_x V(x,t)\,dx\simeq V(x+dx,t)-V(x,t)=-RI(x,t)\,dx-L\partial_t I(x,t)\,dx
\end{eqnarray*}
and
\begin{eqnarray*}&&
\partial_x I(x,t)\,dx\simeq
I(x+dx,t)-I(x,t)=I(x+dx,t)-\iota_1-\iota_3\\&&\qquad
=-\iota_2-\iota_3=-C\partial_t V(x,t)\,dx-GV(x,t)\,dx.
\end{eqnarray*}
Consequently,
\begin{equation}\label{ihkbdDAvbsdewI089JHGSFDCokwfIUyb-X1}
\begin{dcases}
\partial_x V=-L \partial_t I-R I,\\
\partial_x I=-C\partial_t V-GV.
\end{dcases}
\end{equation}
We can take the derivative with respect to~$x$ of the first equation
and compare with the derivative with respect to~$t$ of the second equation
in~\eqref{ihkbdDAvbsdewI089JHGSFDCokwfIUyb-X1}, finding that
\begin{equation}\label{ihkbdDAvbsdewI089JHGSFDCokwfIUyb-X2}
\partial_x^2 V-LC \partial_t^{2} V=(RC+GL) \partial_t V+GRV.
\end{equation}

Similarly, we can take the derivative with respect to~$t$ of the first equation
and compare with the derivative with respect to~$x$ of the second equation
in~\eqref{ihkbdDAvbsdewI089JHGSFDCokwfIUyb-X1}, which yields that
\begin{equation}\label{ihkbdDAvbsdewI089JHGSFDCokwfIUyb-X3}
\partial_x^{2} I-LC\partial_t^{2}I=(RC+GL) \partial_t I+GRI.
\end{equation}
We stress that, up to replacing voltage and current, the structure of~\eqref{ihkbdDAvbsdewI089JHGSFDCokwfIUyb-X2}
and~\eqref{ihkbdDAvbsdewI089JHGSFDCokwfIUyb-X3} is the same.
Also, \eqref{ihkbdDAvbsdewI089JHGSFDCokwfIUyb-X2} and~\eqref{ihkbdDAvbsdewI089JHGSFDCokwfIUyb-X3}
give that both voltage and current are solutions of the telegrapher's equation in~\eqref{HEATELEEQ0}.

It is interesting to note that, in this situation, the structural parameters
in the telegrapher's equation in~\eqref{HEATELEEQ0}
correspond to~$a:=LC$, $b:=RC+GL$ and~$c:=GR$.
Thus, a lossless transmission in which the roles of resistances can be neglected corresponds to~$R:=0$ and~$G:=0$
(this is a perfect conductor offering no resistance and a perfect electrical insulator as a dielectric allowing for no
conductance), which in turn produce~$b:=0$ and~$c:=0$ in~\eqref{HEATELEEQ0}, thus
reducing to the wave equation.
With this observation in mind, one can also have a second look at Figure~\ref{wP2DAILAPIJAM STELEHGREEMIItange3FIlAsAL13458gANSE} to appreciate how the damping and dispersion effects \index{dispersion}
are the outcome of the imperfect conductors and dielectric in terms of resistance and electrical insulation.
\medskip

Let us finally go back to the dispersion problem \index{dispersion} (see footnote~\ref{ihkbdDAvbsdewI089JHGSFDCokwfIUyb-X3-DISm}) which annoyingly affects the transmission of signals along extended transmission lines. To this end, given~$\kappa>0$, let us set\begin{equation}\label{90iojf023ijgHHANSVV123rtgiowiejGBSna0oerL0} v_\kappa:=\sqrt{\frac{4(GR+\kappa^2)}{(GL+RC)^2+4LC\kappa^2}}\qquad{\mbox{and}}\qquad \lambda_\kappa:=\frac{(GL+RC)\, v_\kappa}{2},\end{equation}and observe that the function\begin{equation}\label{90iojf023ijgHHANSVV123rtgiowiejGBSna0oerL} I(x,t):=e^{-\lambda_\kappa x}\cos\big(\kappa(x-v_\kappa t)\big)\end{equation}is a solution of~\eqref{ihkbdDAvbsdewI089JHGSFDCokwfIUyb-X3}.

The function in~\eqref{90iojf023ijgHHANSVV123rtgiowiejGBSna0oerL} presents an interesting shape dictated by the coefficients involved in this construction. Indeed, on the one hand, the exponential term~$e^{-\lambda_\kappa x}$ in~\eqref{90iojf023ijgHHANSVV123rtgiowiejGBSna0oerL} describes the attenuation of the signal along a transmission line. On the other hand the cosine term in~\eqref{90iojf023ijgHHANSVV123rtgiowiejGBSna0oerL} describes a travelling wave, modulated by~$\kappa$. The speed at which this way moves is~$v_\kappa$.
The fact that this speed depends on~$\kappa$ is, in principle, a major obstacle to useful signal transmission, since different frequencies of a given signal would travel along the line at different paces, producing significant distortion phenomena and quickly making the original signal almost unrecognizable.

This problem was brilliantly solved by Oliver Heaviside who figured out that a special tuning of the parameters would produce distortionless transmission lines! Indeed, no dispersion would occur if the line parameters exhibited the following ratio:\begin{equation}\label{KSM0wodje-3ijgHHANSVV123rtgiowiejGBSna0oerL0}\frac{R}L=\frac{G}C.\end{equation}
As a matter of fact, in this case~$GL=RC$ and we see that~\eqref{90iojf023ijgHHANSVV123rtgiowiejGBSna0oerL0} boils down to $$v_\kappa=\frac1{\sqrt{LC}}\qquad{\mbox{and}}\qquad \lambda_\kappa=\sqrt{GR}.$$
In particular, all frequencies travel with the same speed along the transmission line and they all present the same attenuation feature. Though condition~\eqref{KSM0wodje-3ijgHHANSVV123rtgiowiejGBSna0oerL0} is rather ideal (a bit of dispersion and distortion is unavoidable after all), the mathematical approach by Heaviside and the theoretical solution proposed in~\eqref{KSM0wodje-3ijgHHANSVV123rtgiowiejGBSna0oerL0} worked, and long distance transmission lines finally became feasible.

\subsection{Nerd sniping}\label{NESNISE}

\begin{figure}
                \centering
                \includegraphics[width=.95\linewidth]{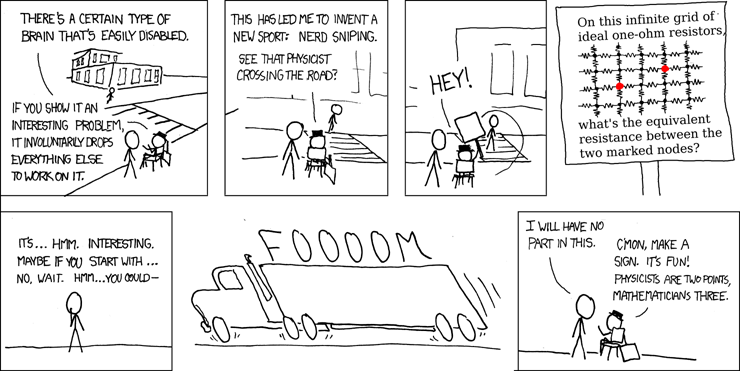}
        \caption{\sl Nerd sniping {\tt https://xkcd.com/356/},
        licensed under the Creative Commons Attribution-NonCommercial 2.5 License.}
        \label{NEMFDRDSIMNIPAI3}
\end{figure}

\begin{figure}
                \centering
                                \includegraphics[width=.4\linewidth]{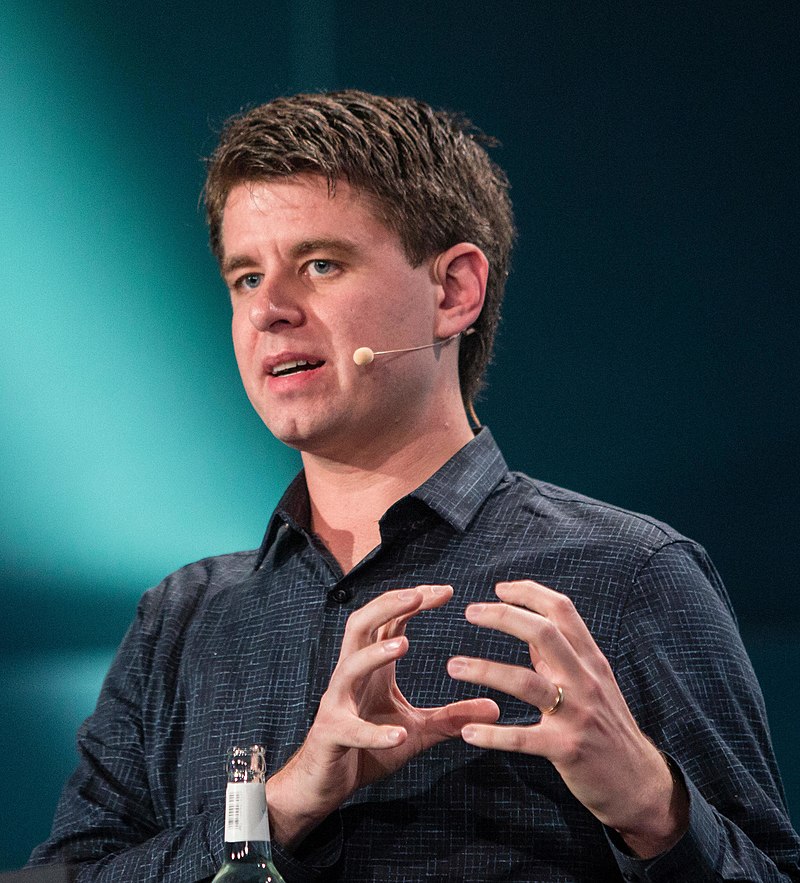}\\
                                                \includegraphics[width=.4\linewidth]{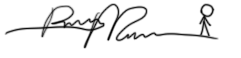}
        \caption{\sl Randall Munroe and his signature (images from
 Wikipedia, the first by re:publica/Jan Zappner,
 licensed under the Creative Commons Attribution 2.0 Generic license, the second in the Public Domain).}
        \label{NEMFDRDSIMNIPAI31-B}
\end{figure}

Nerd sniping is a new sport invented in {\tt https://xkcd.com/356/}
in which the value of mathematicians is frightfully appreciated,
see Figure~\ref{NEMFDRDSIMNIPAI3}. Here, we discuss the solution\footnote{The infinite resistor problem
happened to be studied quite in detail, see in particular~\cite{Cserti},\\
{\tt https://www.explainxkcd.com/wiki/index.php/356:\_Nerd\_Sniping}\\
{\tt https://www.mathpages.com/home/kmath668/kmath668.htm}\\
{\tt https://physics.stackexchange.com/questions/2072/}\\ $\phantom{123456789}$ {\tt on-this-infinite-grid-of-resistors-whats-the-equivalent-resistance}\\
and the references therein.

Strictly speaking, without additional information, the problem is possibly not uniquely posed, but in these pages we confine ourselves
to the case in which the solution is assumed to decay ``fast enough'' at infinity to validate the formal computations
showcased here, though we gloss over any technical aspect related to uniqueness, decay, convergence
and infinite cancellations (but let us mention, for instance, that if~$\Gamma(k)$ is a solution of~\eqref{FU:NE:023}, then
so is~$\Gamma(k)+k_1$ and the equivalent resistance~$R$ in~\eqref{12wrt13kdelta7}  would thus
be affected by an additional term~$-4$, showing a uniqueness issue unless we add some decay assumption on~$\Gamma$).

Also, the webcomic xkcd was created by Randall Patrick Munroe, see Figure~\ref{NEMFDRDSIMNIPAI31-B}.}
of the infinite resistor problem presented in the fourth cartoon of Figure~\ref{NEMFDRDSIMNIPAI3}
and its strong relation\footnote{See also~\cite{MR920811}
for a comprehensive treatment of the deep links between elliptic partial differential equations,
random walks and electric networks. In addition, see~\cite[page~231]{PURCE},
\cite[Chapter~3]{MR2535945}
and the references therein for further information about chains of infinitely many resistors.

Apropos,
according to the
Google Labs Aptitude Test\\
{\tt http://googlesystem.blogspot.com/2005/12/google-labs-aptitude-test.html},\\
Google used the nerd sniping problem as a recruitment question.}
with elliptic partial differential equations (we mean, not the relation between nerds and
elliptic partial differential equations, but the one between the infinite resistor problem
and a discrete variant of elliptic partial differential equations). Actually, the answer to the problem is
\begin{equation}\label{NERD:ANS}
\frac4\pi-\frac12,\end{equation}
but you'll see, it's about the journey, not the destination.

We follow here an approach based on the fundamental solution
of a suitable operator reminiscent of the Laplacian. The methods related to fundamental solutions will be refined in
Section~\ref{lfundsP-S}. For simplicity, we work here in dimension~$2$, but many of the arguments that we present are applicable in higher dimensions and further generality, see~\cite{Cserti} and the references therein. Given~$\e>0$ and a function~$u:\Z^2\to\R$, for all~$x\in\Z^2$ we define
\begin{equation}\label{DISCRELEA} {\mathcal{L}} u(x):=\sum_{j=1}^2\Big( u(x+e_j)+u(x-e_j)-2u(x)\Big).\end{equation}
Notice that the operator~${\mathcal{L}}$ can be seen as a discrete version\footnote{For more information
on the discretization of
the Laplace operator (and for the opportunities and dangers entailed by discretization), see e.g.~\cite[Chapter~8]{MR2398759}.}
of the Laplacian. Our main goal will be to determine the fundamental solution of this operator, namely a function~$\Gamma:\Z^2\to\R$ such that
\begin{equation}\label{FU:NE:023} {\mathcal{L}}\,\Gamma =\delta_0,\end{equation}
where~$\delta_0$ is the Dirac Delta Function at the origin in the ``discrete sense'', meaning that
$$ \sum_{k\in\Z^2} \varphi(k)\,\delta_0(k)=\varphi(0)$$
for all~$\varphi:\Z^2\to\R$.

The convenience of using this fundamental solution is due to the following observation. We put coordinates in~$\R^2$ such that the nodes of the infinite resistor problem presented in the fourth cartoon of Figure~\ref{NEMFDRDSIMNIPAI3} correspond to the lattice~$\Z^2$.
We can also suppose that the two red points in the lattice
of Figure~\ref{NEMFDRDSIMNIPAI3} correspond to~$(0,0)$ and~$(2,1)$.
Then, to test the resistors, we aim at constructing a distribution of voltage~$V:\Z^2\to\R$ which produces a current with unit intensity from~$(0,0)$ to~$(2,1)$ and no further net current in the circuit.
Hence, we apply Ohm's Law at each node~$k$ of the lattice, with respect to its first neighbors~$k\pm e_j$, for~$j\in\{1,2\}$.
Namely, by Ohm's Law, since all resistors of Figure~\ref{NEMFDRDSIMNIPAI3} are unitary, the current flowing into the node~$k$ from its neighbors~$k\pm e_j$ is equal to
$$ \sum_{j=1}^2 (V({k+ e_j})-V(k))+\sum_{j=1}^2 (V({k- e_j})-V(k)),$$
which in turn equals to~${\mathcal{L}}V(k)$.
Hence, the voltage distribution flowing one unit of current from~$(0,0)$ to~$(2,1)$ is such that
\begin{equation}\label{12wrt13kdelta0}
{\mathcal{L}}V(k)=
\delta_{(2,1)}(k)-\delta_{(0,0)}(k)
=\delta_{0}(k-(2,1))-\delta_{0}(k).
\end{equation}
If we found such a voltage distribution~$V$, we can apply Ohm's Law again
and determine the requested equivalent resistance~$R$ between
the two red nodes in Figure~\ref{NEMFDRDSIMNIPAI3} via the relation
\begin{equation}\label{12wrt13kdelta1} R=V(0,0)-V(2,1).\end{equation}
But to find such a~$V$, it comes in handy to use the fundamental solution~$\Gamma$:
indeed, if we find~$\Gamma$ as in~\eqref{FU:NE:023} it suffices to set
\begin{equation}\label{12wrt13kdelta2} V(k):=\Gamma(k-(2,1))-\Gamma(k)\end{equation}
and observe that this is a solution of~\eqref{12wrt13kdelta0}.

Summarizing, in view of~\eqref{12wrt13kdelta1} and~\eqref{12wrt13kdelta2}, once we determine the fundamental solution~$\Gamma$ we can also find the desired equivalent resistance~$R$ by the formal relation
\begin{equation} \label{12wrt13kdelta7} R=V(0,0)-V(2,1)=
\Big( \Gamma(-2,-1)-\Gamma(0,0)\Big)
- \Big( \Gamma(0,0)-\Gamma(2,1)\Big)
.\end{equation}

Hence, we focus now on determining the fundamental solution~$\Gamma$. To this end, we deduce from~\eqref{FU:NE:023} that, for every~$\xi\in\R^2$,
\begin{equation}\label{PMSTABSNdavsRAFSDs}\begin{split}&
1=\sum_{k\in\Z^2} \delta_0(k)\,e^{2\pi ik\cdot\xi} =\sum_{k\in\Z^2}{\mathcal{L}}\,\Gamma (k)\,e^{2\pi ik\cdot\xi}=
\sum_{{j\in\{1,2\}}\atop{k\in\Z^2}}\Big(\Gamma(k+e_j)+\Gamma(k-e_j)-2\Gamma(k)\Big)\,e^{2\pi ik\cdot\xi}\\&\qquad\qquad=
\sum_{{j\in\{1,2\}}\atop{m\in\Z^2}} \Gamma(m)\,e^{2\pi i(m-e_j)\cdot \xi}+
\sum_{{j\in\{1,2\}}\atop{m\in\Z^2}}\Gamma(m)\,e^{2\pi i(m+e_j)\cdot \xi}-2\sum_{{j\in\{1,2\}}\atop{m\in\Z^2}}\Gamma(m)\,e^{2\pi im\cdot \xi}\\&\qquad\qquad=
\sum_{{j\in\{1,2\}}\atop{m\in\Z^2}} \Big(e^{-2\pi i\xi_j}+e^{2\pi i\xi_j}-2\Big)\,\Gamma(m)\,e^{2\pi im\cdot \xi} =-\mu(\xi)
\sum_{{m\in\Z^2}} \Gamma(m)\,e^{2\pi im\cdot \xi}
,\end{split}\end{equation}
where
$$ \mu(\xi):=2\sum_{j=1}^2\big( 1-\cos(2\pi \xi_j)\big).$$
We remark that~$\mu$ is a periodic function, namely~$\mu(\xi+\ell)=\mu(\xi)$ for every~$\ell\in\Z^2$ and~$\xi\in\R^2$.
Thus, also~$\frac1\mu$ is periodic and then we can write this function as a Fourier Series.
For this, to be precise, since~$\mu$ is nonnegative but may vanish,
it is convenient to pick~$\delta>0$ and define~$\mu_\delta(\xi):=\max\{\delta,\mu(\xi)\}$,
expand the function~$\frac1{\mu_\delta}$ in Fourier Series and then pass~$\delta\searrow0$ one way or another.
Namely, we have that
$$ \frac1{\mu_\delta(\xi)}=\sum_{m\in\Z^2} c_{m,\delta} \,e^{2\pi i m\cdot\xi},$$
where
$$ c_{m,\delta}:=\int_{Q} \frac{e^{-2\pi i m\cdot\eta}}{\mu_\delta(\eta)}\,d\eta\;\qquad{\mbox{and}}\qquad \; Q:=\left(-\frac12,\frac12\right)^2.$$
Thus, if~$\phi:\R^2\to\R$ is a smooth and periodic function
and $$ \widehat\phi_m:=\int_{Q} \phi(\xi)\,e^{-2\pi im\cdot \xi}\,d\xi$$
is the Fourier coefficient of~$\phi$, we formally deduce from~\eqref{PMSTABSNdavsRAFSDs} that
\begin{equation}\label{N2A3S3DwedabN023VhwVaq}\begin{split}&
\lim_{\delta\searrow0}
\sum_{m\in\Z^2} c_{m,\delta} \,\widehat\phi_{-m}=
\lim_{\delta\searrow0}\int_{Q}
\sum_{m\in\Z^2} c_{m,\delta} \,\phi(\xi)\,e^{2\pi i m\cdot\xi}\,d\xi=
\lim_{\delta\searrow0}\int_{Q} \frac{\phi(\xi)}{\mu_\delta(\xi)}\,d\xi=
\int_{Q} \frac{\phi(\xi)}{\mu(\xi)}\,d\xi\\&\qquad\qquad=
-\int_{Q}\sum_{{m\in\Z^2}} \Gamma(m)\,e^{2\pi im\cdot \xi}\,\phi(\xi)\,d\xi=-\sum_{{m\in\Z^2}} \Gamma(m)\,
\widehat\phi_{-m}.
\end{split}\end{equation}
We can thereby apply~\eqref{N2A3S3DwedabN023VhwVaq}
to the function~$\psi(\xi):=\phi(\xi)-\phi(0)$ and find that
\begin{equation}\label{76t54wedMSKAMS0987654DIJMD2r3teg}
\lim_{\delta\searrow0}
\sum_{m\in\Z^2}\int_{Q} \widehat\psi_{-m} \frac{e^{-2\pi i m\cdot\eta}\,d\eta}{\mu_\delta(\eta)}=
\lim_{\delta\searrow0}
\sum_{m\in\Z^2} c_{m,\delta} \,\widehat\psi_{-m}=
-\sum_{{m\in\Z^2}} \Gamma(m)\,
\widehat\psi_{-m}.
\end{equation}
Also, since
$$ \sum_{m\in\Z^2} \widehat\psi_{-m}=\sum_{m\in\Z^2} \widehat\psi_m=
\sum_{m\in\Z^2} \widehat\psi_m\,e^{2\pi im\cdot 0}=\psi(0)=0,$$
we have that
\begin{eqnarray*}
\int_{Q}\sum_{m\in\Z^2} \widehat\psi_{-m}\,\frac{d\eta}{\mu_\delta(\eta)}=0.
\end{eqnarray*}
This and~\eqref{76t54wedMSKAMS0987654DIJMD2r3teg} lead to
\begin{equation}\label{76t54wedMSKAMS0987654DIJMD2r3teg3}
\lim_{\delta\searrow0}
\sum_{m\in\Z^2}\int_{Q} \widehat\psi_{-m} \frac{(e^{-2\pi i m\cdot\eta}-1)\,d\eta}{\mu_\delta(\eta)}=
-\sum_{{m\in\Z^2}} \Gamma(m)\,
\widehat\psi_{-m}.
\end{equation}
We also point out that~$\mu_\delta$ is an even function and therefore
$$ \int_{Q} \frac{\sin(2\pi m\cdot\eta)\,d\eta}{\mu_\delta(\eta)}=0.$$
{F}rom this and~\eqref{76t54wedMSKAMS0987654DIJMD2r3teg3} it follows that
\begin{equation}\label{76t54wedMSKAMS0987654DIJMD2r3teg2}
\lim_{\delta\searrow0}
\sum_{m\in\Z^2}\int_{Q} \widehat\psi_{-m} \frac{\big(\cos(2\pi m\cdot\eta)-1\big)\,d\eta}{\mu_\delta(\eta)}=
-\sum_{{m\in\Z^2}} \Gamma(m)\,
\widehat\psi_{-m}.
\end{equation}

The advantage of~\eqref{76t54wedMSKAMS0987654DIJMD2r3teg2} with respect to~\eqref{76t54wedMSKAMS0987654DIJMD2r3teg} is that the singularity at the denominator
in the first integral is compensated by the vanishing of the corresponding numerator, thus allowing us
to pass to the limit inside the integral and write that
\begin{equation}\label{76t54wedMSKAMS0987654DIJMD2r3teg5}
\sum_{m\in\Z^2}\int_{Q} \widehat\psi_{-m} \frac{\big(\cos(2\pi m\cdot\eta)-1\big)\,d\eta}{\mu(\eta)}=
-\sum_{{m\in\Z^2}} \Gamma(m)\,
\widehat\psi_{-m}.
\end{equation}
We also observe that
$$\widehat\psi_0=\widehat\phi_0-\phi(0)=
\widehat\phi_0-\sum_{m\in\Z^2}\widehat\phi_m e^{2\pi i m\cdot0}=\widehat\phi_0-\sum_{m\in\Z^2}\widehat\phi_{-m}
$$ and~$
\widehat\psi_m=\widehat\phi_m$ for all~$m\in\Z^2\setminus\{0\}$. Plugging this information into~\eqref{76t54wedMSKAMS0987654DIJMD2r3teg5} we find that
\begin{equation*}
\sum_{m\in\Z^2}\int_{Q} \widehat\phi_{-m} \frac{\big(\cos(2\pi m\cdot\eta)-1\big)\,d\eta}{\mu(\eta)}=
\sum_{{m\in\Z^2}} \big(\Gamma(0)-\Gamma(m)\big)\,
\widehat\phi_{-m}.
\end{equation*}
The arbitrariness of~$\phi$ thus gives that
\begin{equation*}
\begin{split}
\Gamma(m)-\Gamma(0)&=-\int_{Q} \frac{\big(\cos(2\pi m\cdot\eta)-1\big)\,d\eta}{\mu(\eta)}\\&=\frac12
\iint_{(-1/2,1/2)\times(-1/2,1/2)} \frac{ 1-\cos(2\pi m_1\eta_1+2\pi m_2\eta_2) }{2-\cos(2\pi \eta_1)-\cos(2\pi \eta_2)}\,d\eta_1\,d\eta_2.
\end{split}\end{equation*}
We insert this into the equivalent resistance relation~\eqref{12wrt13kdelta7} and we arrive at
\begin{equation} \label{76t54wedMSKAMS0987654DIJMD2r3teg5b90} R=
\iint_{(-1/2,1/2)\times(-1/2,1/2)} \frac{ 1-\cos(4\pi \eta_1+2\pi \eta_2) }{2-\cos(2\pi \eta_1)-\cos(2\pi \eta_2)}\,d\eta_1\,d\eta_2
.\end{equation}
The good news is that this is an integral involving simple trigonometric functions, so we should feel confident to solve it explicitly thanks to our consolidated calculus skills. For instance, we can proceed as follows. We note that
\begin{equation} \label{76t54wedMSKAMS0987654DIJMD2r3teg5b}\begin{split}\cos(4\pi \eta_1+2\pi \eta_2)\,&=
\cos(4\pi \eta_1)\cos(2\pi \eta_2)
-\sin(4\pi \eta_1)\sin(2\pi \eta_2)\\&=
\Big( 1-2\sin^2(2\pi \eta_1)\Big)\cos(2\pi \eta_2)
-2\sin(2\pi \eta_1)\cos(2\pi \eta_1)\sin(2\pi \eta_2).
\end{split}\end{equation}
Hence, if~$\tau_1:=\tan(\pi \eta_1)$ and~$\tau_2:=\tan(\pi \eta_2)$, using in~\eqref{76t54wedMSKAMS0987654DIJMD2r3teg5b} the ``tangent half-angle formulas''
$$ \sin(2\pi \eta_j)=\frac{2\tau_j}{1+\tau_j^2}\qquad{\mbox{and}}\qquad
\cos(2\pi \eta_j)=\frac{1-\tau_j^2}{1+\tau_j^2},$$
we find that
%%% \begin{equation*}\begin{split}\cos(4\pi \eta_1+2\pi \eta_2)\,&=
%%% \left( 1-2\left( \frac{2\tau_1}{1+\tau_1^2} \right)^2\right)\frac{1-\tau_2^2}{1+\tau_2^2}
%%% - \frac{8\tau_1(1-\tau_1^2)\tau_2}{(1+\tau_1^2)^2\,(1+\tau_2^2)}\\&
%%% =\frac{-8 \tau_1 (1-\tau_1^2) \tau_2 + ( \tau_1^4-6 \tau_1^2 + 1) (1 - \tau_2^2)}{(1+\tau_1^2)^2 \,(1+\tau_2^2)}
%%% \end{split}\end{equation*}
%%% and therefore
\begin{equation}\label{76t54wedMSKAMS0987654DIJMD2r3teg5c}1-\cos(4\pi \eta_1+2\pi \eta_2)=\frac{2 (2 \tau_1 + \tau_2 - \tau_1^2 \tau_2)^2}{(1 + \tau_1^2)^2 \,(1 + \tau_2^2)}.
\end{equation}
Similarly,
\begin{equation}\label{76t54wedMSKAMS0987654DIJMD2r3teg5d}
2-\cos(2\pi \eta_1)-\cos(2\pi \eta_2)=
2-\frac{1-\tau_1^2}{1+\tau_1^2}-\frac{1-\tau_2^2}{1+\tau_2^2}=
\frac{2 (\tau_1^2 + \tau_2^2 + 2 \tau_1^2 \tau_2^2)}{(1+\tau_1^2) (1+\tau_2^2)}.
\end{equation}
Thus, since
$$ d\eta_j=\frac{d\tau_j}{\pi (1+\tau_j^2)},$$
we infer from~\eqref{76t54wedMSKAMS0987654DIJMD2r3teg5c} and~\eqref{76t54wedMSKAMS0987654DIJMD2r3teg5d} that
\begin{eqnarray*}
\frac{ 1-\cos(4\pi \eta_1+2\pi \eta_2) }{2-\cos(2\pi \eta_1)-\cos(2\pi \eta_2)}\,d\eta_1\,d\eta_2=
\frac{ (\tau_1^2 \tau_2 - 2 \tau_1 - \tau_2)^2}{\pi^2 (1+\tau_1^2 )^2 \,(1+\tau_2^2 ) \,(2 \tau_1^2 \tau_2^2 + \tau_1^2 + \tau_2^2)}\,d\tau_1\,d\tau_2
\end{eqnarray*}
This and~\eqref{76t54wedMSKAMS0987654DIJMD2r3teg5b90} lead to
\begin{equation*} \begin{split}
R&=
\iint_{\R\times\R} 
\frac{(\tau_1^2 \tau_2 - 2 \tau_1 - \tau_2)^2}{\pi^2
(1+\tau_1^2 )^2 \,(1+\tau_2^2 ) \,(2 \tau_1^2 \tau_2^2 + \tau_1^2 + \tau_2^2)}\,d\tau_1\,d\tau_2\\&=\iint_{\R\times\R} 
\frac{(x y^2 -x - 2 y)^2}{\pi^2 (1+x^2 ) \,(1+y^2 )^2 \,(2 x^2 y^2 + x^2 + y^2)}\,dx\,dy
,\end{split}\end{equation*}
where we simplified notation using the variables~$(x,y)$ in place of~$(\tau_2,\tau_1)$.
Since the latter integral remains invariant under the map~$(x,y)\mapsto(-x,-y)$ we can reduce the previous equation to
\begin{equation}\label{76t54wedMSKAMS0987654DIJMD2r3teg5dwqr326tygh}
R=\iint_{\R\times(0,+\infty)} 
\frac{2 (x y^2 -x - 2 y)^2}{\pi^2
(1+x^2 ) \,(1+y^2 )^2 \,(2 x^2 y^2 + x^2 + y^2)}\,dx\,dy
.\end{equation}

Now one can check that a primitive of the function~$\frac{(1+y^2)(x y^2 -x - 2 y)^2}{(1+x^2 ) (2 x^2 y^2 + x^2 + y^2)} $ in the variable~$x$ is
\begin{eqnarray*} \Phi(x,y)&:=&2 y (y^2 - 1) \ln\frac{1+x^2}{2 x^2 y^2 + x^2 + y^2} + (y^4 - 6 y^2 + 1) \arctan x \\&&\qquad+ \frac{y ( 10 y^2-y^4 + 3) \arctan\frac{x \sqrt{2 y^2 + 1}}{y}}{\sqrt{2 y^2 + 1}}. \end{eqnarray*}
Hence, for every~$y\in(0,+\infty)$,
\begin{eqnarray*} \Phi(\pm\infty,y)&=&2 y (y^2 - 1) \ln\frac{1}{2 y^2 + 1} \pm\frac{\pi}{2}\left((y^4 - 6 y^2 + 1) +\frac{y ( 10 y^2-y^4 + 3) }{\sqrt{2 y^2 + 1}}\right). \end{eqnarray*}
As a result, for all~$y\in(0,+\infty)$,
\begin{eqnarray*}
\int_\R \frac{(1+y^2)(x y^2 -x - 2 y)^2}{(1+x^2 ) (2 x^2 y^2 + x^2 + y^2)} \,dx&=&
\Phi(+\infty,y)-\Phi(-\infty,y)\\&=&
\pi  (y^4 - 6 y^2 + 1) +\frac{\pi y ( 10 y^2-y^4 + 3) }{\sqrt{2 y^2 + 1}}\end{eqnarray*}
and consequently
\begin{eqnarray*}
\int_\R \frac{2 (x y^2 -x - 2 y)^2}{\pi^2
(1+x^2 ) \,(1+y^2 )^2 \,(2 x^2 y^2 + x^2 + y^2)}\,dx &=&
\frac{2(y^4 - 6 y^2 + 1)}{\pi (1+y^2)^3} +
\frac{2 y ( 10 y^2-y^4 + 3)}{\pi (1+y^2)^3\sqrt{2 y^2 + 1}}.
\end{eqnarray*}
Combining this with~\eqref{76t54wedMSKAMS0987654DIJMD2r3teg5dwqr326tygh} we find that
\begin{eqnarray*}
R&=&\int_{0}^{+\infty} 
\frac{2(y^4 - 6 y^2 + 1)}{\pi (1+y^2)^3} \,dy+\int_{0}^{+\infty}
\frac{2 y ( 10 y^2-y^4 + 3)}{\pi (1+y^2)^3\sqrt{2 y^2 + 1}}\,dy
\\&=& \left. \frac{2(y - y^3)}{\pi (1 + y^2)^2}\right|_{y=0}^{y=+\infty}
-\left. \frac2{\pi}\left(\frac{2 \sqrt{2 y^2 + 1}}{(1+y^2)^2} +\arctan\sqrt{2 y^2 + 1}\right)\right|_{y=0}^{y=+\infty}\\
&=&-\frac2{\pi}\left(\frac{\pi}2-2
-\arctan 1\right)\\&=&\frac{\pi}4-\frac12,
\end{eqnarray*}
in accordance with~\eqref{NERD:ANS}, as desired.

\subsection{The streetsweeper problem}\label{I2K243S3M:0oikrfg-4-2rtghj}

\begin{figure}
                \centering
                \includegraphics[height=.25\textheight]{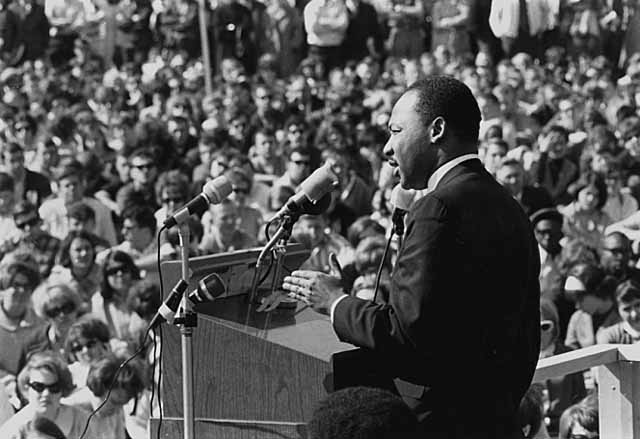}
        \caption{\sl Martin Luther King Jr. speaking against the Vietnam War
at the University of Minnesota
(image of the Minnesota Historical Society from
Wikipedia, licensed under the 
Creative Commons Attribution-Share Alike 2.0
Generic license).}\label{HMLKlFUMHDNOJHNFOJED}
\end{figure}

\begin{figure}
                \centering
                \includegraphics[height=.3\textheight]{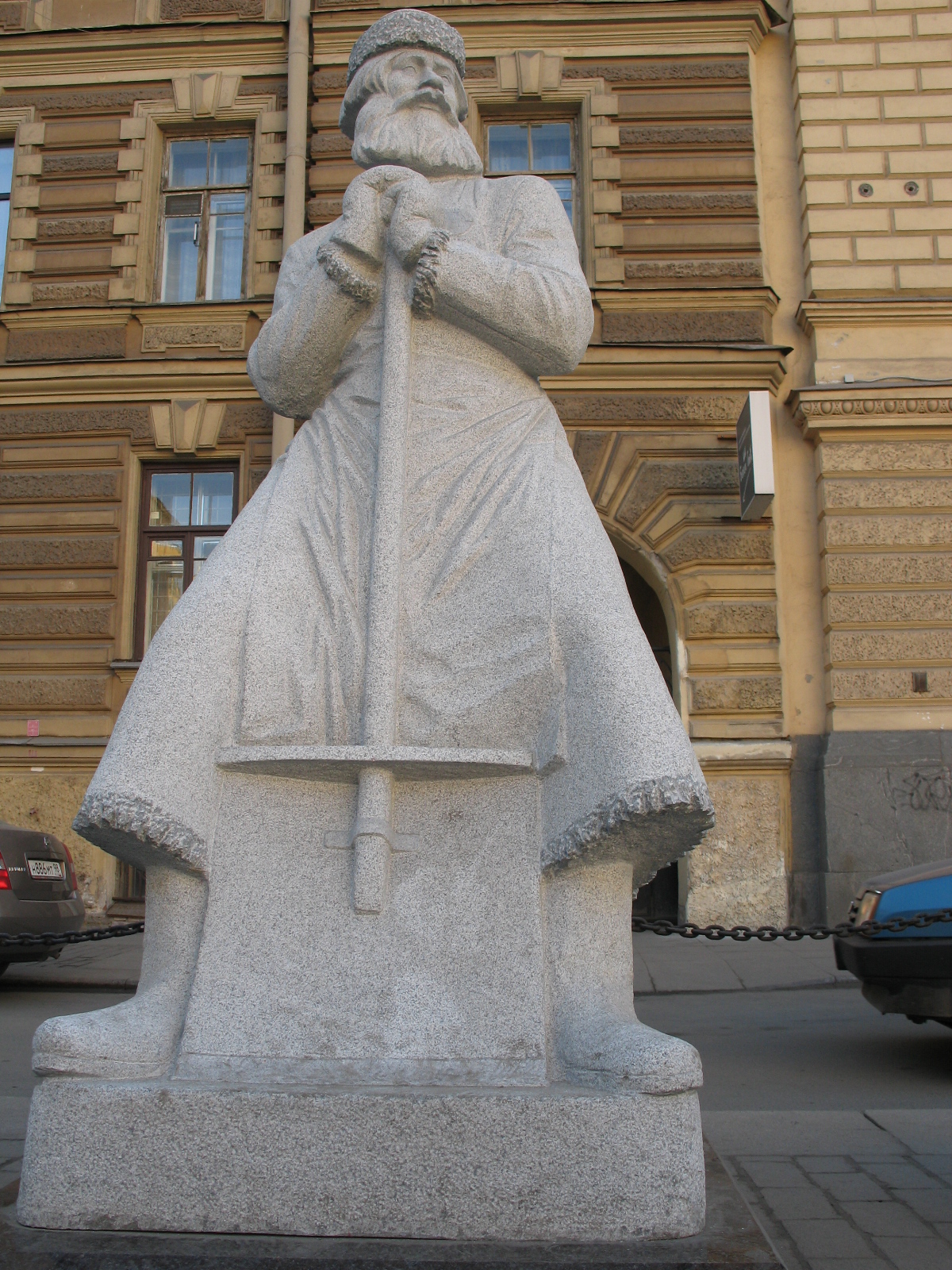} $\,$
                \includegraphics[height=.3\textheight]{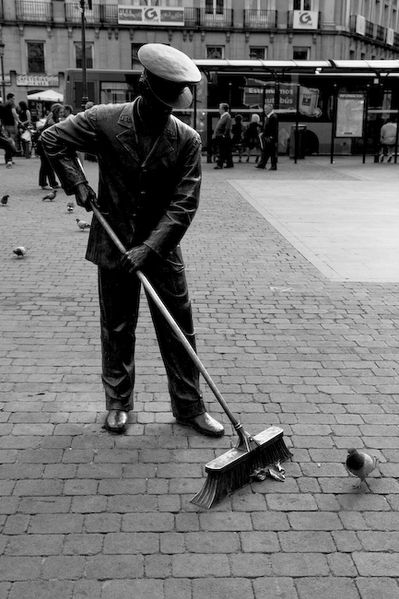}
        \caption{\sl Left: Monument of streetsweeper in Saint Petersburg (Public Domain image from
        Wikipedia).
Right: Monument of streetsweeper in Madrid
(photo by Alex Pascual Guardia, image from
Wikipedia, licensed under the Creative Commons Attribution 2.0 Generic license).}\label{HASSWE7solFUMHDNOJHNFOJED}
\end{figure}

Martin Luther King Jr. (see Figure~\ref{HMLKlFUMHDNOJHNFOJED})
once said ``What I'm saying to you this morning, my friends, even if it falls your lot to be a streetsweeper,
go on out and sweep streets like Michelangelo painted pictures;
sweep streets like Handel and Beethoven composed music;
sweep streets like Shakespeare wrote poetry;
sweep streets so well that all the host of heaven and earth will have to
pause and say: Here lived a great streetsweeper who swept his job well''.
At the end of the day, whatever we do and whatever the reason,
the competence and dedication that we put in what we do makes 
all the difference in the world,
and a large number of tasks and heavy duties are essential for the prosperity of the whole community. For this reason, and many others, the streetsweeper
is represented in several monuments worldwide, see e.g. Figure~\ref{HASSWE7solFUMHDNOJHNFOJED}.

\begin{figure}
                \centering
                \includegraphics[height=.3\textheight]{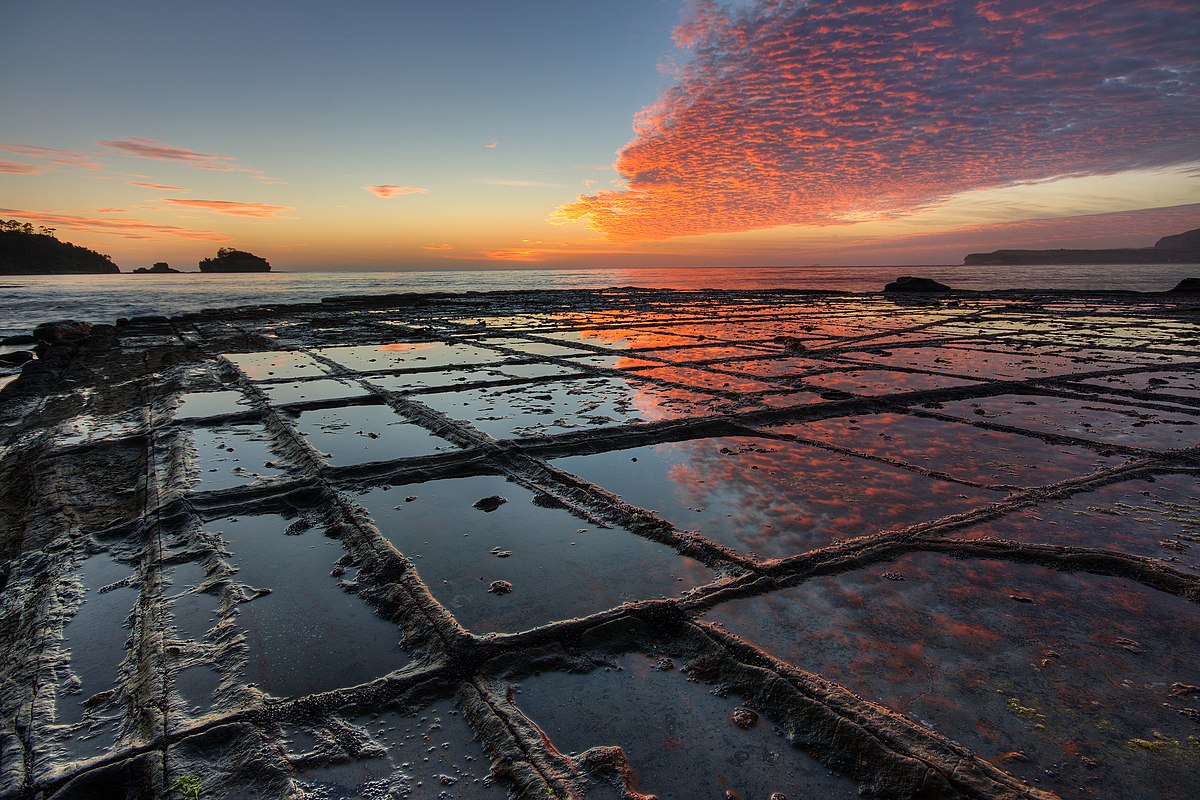} $\,$
                \includegraphics[height=.3\textheight]{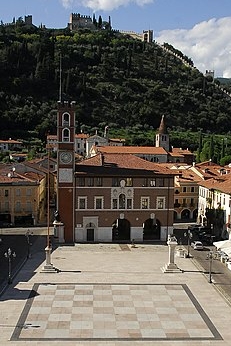}
        \caption{\sl Left: natural tessellated pavement at Eaglehawk Neck, Tasmania (photo by J. J. Harrison, image from
        Wikipedia, licensed under the Creative Commons Attribution-Share Alike 2.5 Generic license).
Right: Piazza degli Scacchi, Marostica, Italy
(detail of a photo by Rino Porrovecchio, image from
        Wikipedia, licensed under the  Creative Commons Attribution-Share Alike 2.0 Generic license).}\label{TAS89EEMIItange3FIlAALANSE}
\end{figure}

\begin{figure}
                \centering
                \includegraphics[width=0.9\textwidth]{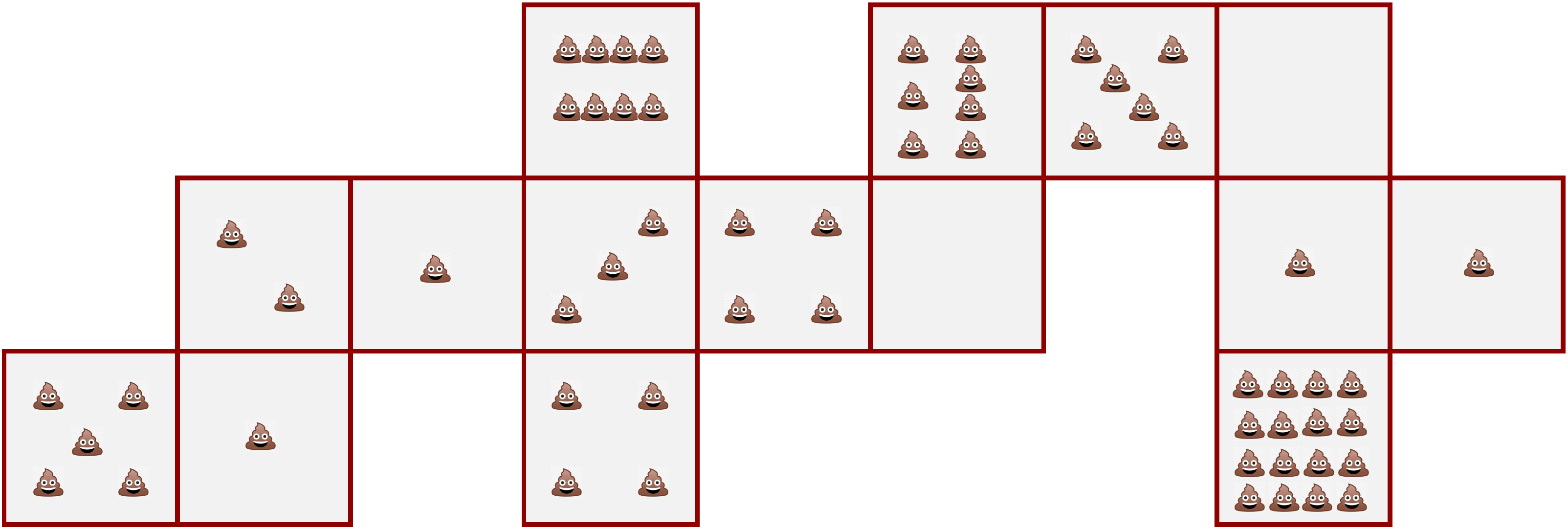}
        \caption{\sl A tessellated floor to be cleaned.}\label{TBvcDNOJHNFOJED}
\end{figure}

The problem\footnote{This problem was kindly suggested to us by
H\'{e}ctor Chang-Lara, see also the beautiful video {\tt https://www.youtube.com/watch?v=--iUT5IYxZs}
for further inspiration.}
that we present now is also inspired by the noble figure of the streetsweeper.
Suppose that the task of a streetsweeper is to clean a
finite collection of tiles in the floor (see Figure~\ref{TAS89EEMIItange3FIlAALANSE}
for concrete example of tesselated floors). To make a mathematical model describing this job,
we can consider a tiling of the plane~$\R^2$ made of square tiles of the form~$\left(k_1-\frac12,k_1+\frac12\right)\times\left(k_2-\frac12,k_2+\frac12\right)$ with~$k=(k_1,k_2)\in\Z^2$. In this way, every tile can be indexed by
the integer coordinates~$k$ of its center.
We suppose that initially every tile contains a certain amount of rubbish.
For this, we write that the tile centered at~$k\in\Z^2$ contains initially a quantity~$\mu_0(k)$ of rubbish.
In this way, we can think that~$\mu_0:\Omega\to[0,+\infty)$
where~$\Omega\subset\Z^2$ is the finite collection of the centers of the tiles, see Figure~\ref{TBvcDNOJHNFOJED}.
For simplicity, we suppose that
\begin{equation}\label{TASKIN-0}{\mbox{each tile of~$\Omega$ possesses at least one adjacent tile in~$\Omega$.}}\end{equation}
This is to exclude the somewhat trivial case of having single and isolated tiles.

Then, the cleaning process works as follows: the streetsweeper selects a tile of~$\Omega$ (say, the one centered at some~$k\in\Z^2$)
and distributes the rubbish present there into the four adjacent tiles (namely,
the ones centered at~$k+(0,1)$, $k-(0,1)$, $k+(1,0)$ and~$k-(1,0)$), regardless if they belong to~$\Omega$ or not,
and then repeats the process by selecting another tile, distributing the rubbish into the adjacent tiles, and so on.
The rubbish that ends up outside~$\Omega$ is left there (in particular,
tiles outside~$\Omega$ but adjacent to tiles in~$\Omega$ get overtime filled in by the rubbish that gets
removed from~$\Omega$).

The goal of the streetsweeper is thus to clean~$\Omega$ from all the rubbish.
Notice that \begin{equation}\label{TASKIN}\begin{split}&{\mbox{this cleaning task
can only be accomplished ``in the limit'',}}\\ &{\mbox{that is only after ``infinitely many iterations of the process''.}}\end{split}\end{equation}
Indeed suppose, by contradiction, that after a finite number of iterations the streetsweeper
has removed all the rubbish from~$\Omega$. Consider the latest tile that the streetsweeper has cleaned,
and let~$\mu_\star>0$ be the amount of rubbish located there right before being cleaned.
Then, the latest operation of the streetsweeper was to remove~$\mu_\star$ from the tile and place~$\frac{\mu_\star}4$
on each of the adjacent tiles. But, by~\eqref{TASKIN-0}, we know that at least one of this adjacent tiles belong to~$\Omega$,
therefore the whole of~$\Omega$ has not been cleaned. This contradiction proves~\eqref{TASKIN}, as desired.

\begin{figure}
                \centering
                \includegraphics[width=.9\textwidth]{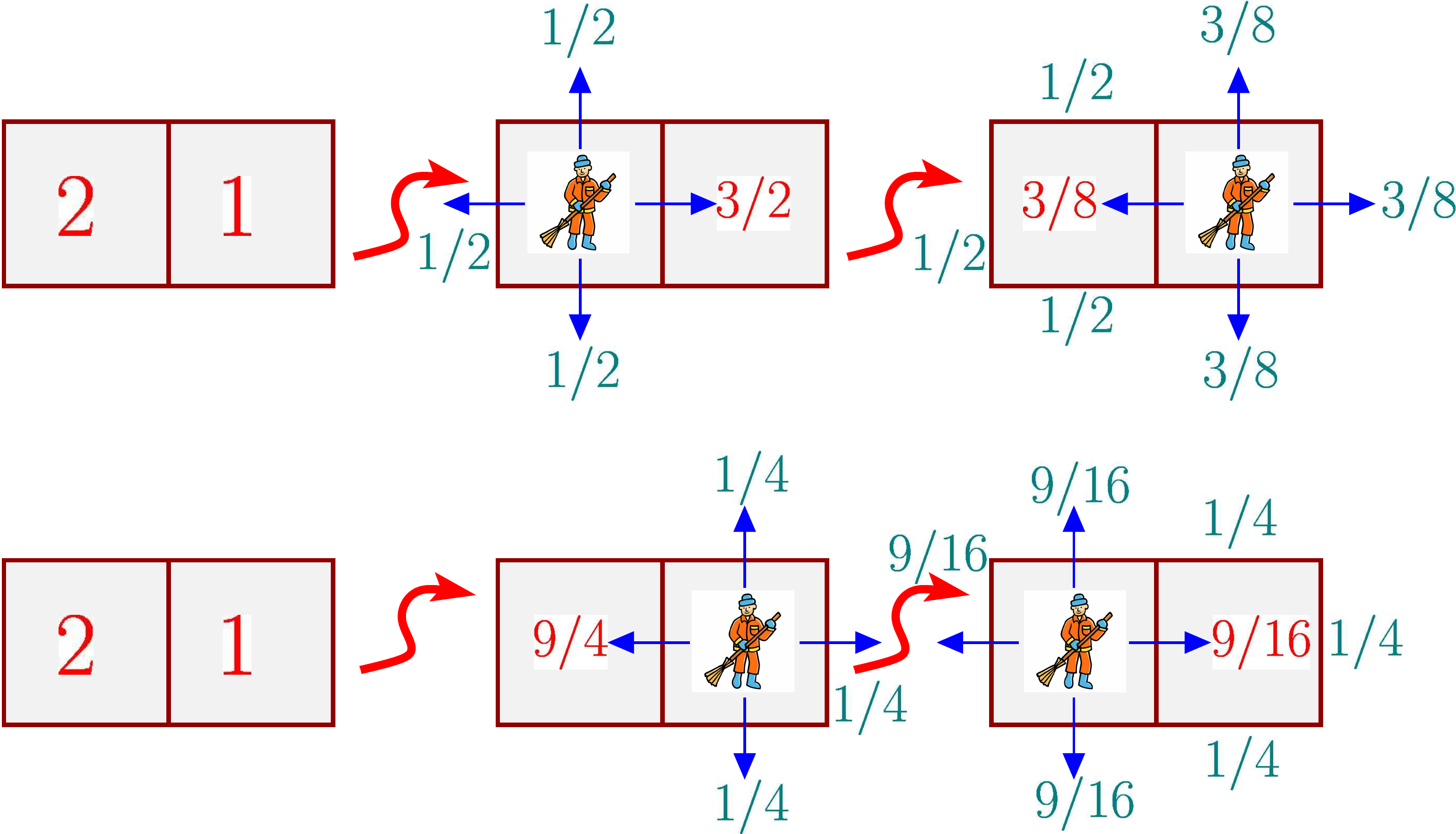}
        \caption{\sl Two different cleaning strategies in a very simple example.}\label{TAS89EEMIItange3FIlAALANSESwerEqXATWO}
\end{figure}

To better model the cleaning process, we denote by~$\mu_j:\Z^2\to[0,+\infty)$
the amount of rubbish located at each tile after~$j$ iterations of the clean process
(in this way, having denoted by~$\mu_0$ the initial rubbish distribution,
we have that~$\mu_1$ stands for the rubbish distribution after one step, and so on).
We stress that~$\mu_j$ does depend on the strategy chosen by the 
streetsweeper to clean the floor, and this strategy is not ``commutative''
(that is, cleaning one tile and then another produces, in general, a different result from
the same operation performed in the reverse order).
As an example for this lack of commutativity, we can consider the simple case in which~$\Omega:=\{(0,0),$ $(0,1)\}$,
$\mu_0(0,0)=2$ and~$\mu_0(0,1)=1$.
If the streetsweeper decided to clean first the tile centered in~$(0,0)$ and then the one centered in~$(0,1)$,
the corresponding result would be~$\mu_1(0,0)=0$,
$\mu_1(0,1)=\frac32$, $\mu_2(0,0)=\frac38$ and~$\mu_2(0,1)=0$.
Instead,
if the streetsweeper decided to clean first the tile centered in~$(0,1)$ and then the one centered in~$(0,0)$,
the corresponding result would be~$\mu_1(0,0)=\frac94$,
$\mu_1(0,1)=0$, $\mu_2(0,0)=0$ and~$\mu_2(0,1)=\frac9{16}$, thus showing that the cleaning process is
not commutative and the cleaning strategy plays a role for the subsequent rubbish distributions
(see Figure~\ref{TAS89EEMIItange3FIlAALANSESwerEqXATWO}
for a sketch of the two different cleaning strategies).

However, quite surprisingly, \begin{equation}\label{SWEE}\begin{split}&{\mbox{if the streetsweeper manages to clean~$\Omega$ in the limit,}}\\ &{\mbox{then the limit distribution of rubbish outside~$\Omega$}}\\
&{\mbox{is independent of the successful cleaning strategy chosen.}}\end{split}\end{equation}
This is a rather deep result, and to prove it we will need to introduce and understand some bits
of elliptic partial differential equations (at least in the version of discrete mathematics that we have
already encountered in Section~\ref{NESNISE}).
To this end, it is convenient to denote by~$u_j:\Z^2\to[0,+\infty)$ the function that keeps a record
of ``how much rubbish went out in total from a given tile after~$j$ iterations of the cleaning process''
(with~$u_0:=0$).
Notice that
\begin{equation}\label{PERmuinftyk} u_j(k)=\begin{dcases}
u_{j-1}(k)+\mu_{j-1}(k) & \begin{matrix}{\mbox{if the tile centered at~$k$ is the one selected}}&\\
{\mbox{for the $j$th iteration of the cleaning process,}}&\end{matrix}\\
\\
u_{j-1}(k)&{\mbox{otherwise.}}
\end{dcases}\end{equation}
In particular, for a given~$k\in\Z^2$, the sequence~$u_{j}(k)$ is monotone and therefore we can define
\begin{equation}\label{ULIbar} u_\infty(k):=\lim_{j\to+\infty}u_j(k).\end{equation}
Furthermore, since the cleaning process is only performed at the tiles in~$\Omega$,
no rubbish is taken outside the tiles centered in~$\Z^2\setminus\Omega$, meaning that~$u_j(k)=0$,
whence
\begin{equation}\label{ULIbar-ZER}{\mbox{$u_\infty(k)=0$ for every~$k\in\Z^2\setminus\Omega$.}}\end{equation}

We can also consider~$\widetilde u_j:=u_j-u_{j-1}$ to be the rubbish taken out at the $j$th step of the cleaning process.
Note that~$\widetilde u_j(k)\ge0$ for all~$k\in\Z^2$.
We stress that, at the~$j$th iteration, at the tile centered at~$k$ 
the rubbish can increase or decrease: specifically, it decreases by~$\widetilde u_j(k)$
(which is strictly positive, coinciding with~$\mu_{j-1}(k)$, only if~$k$ is the tile selected to be cleaned at the $j$th step of the cleaning process) and increases by
$$\frac14\sum_{{i\in\{1,2\}}\atop{\sigma\in\{-1,+1\}}}\widetilde u_j(k+\sigma e_i)$$
(which is the amount of rubbish possibly received by the adjacent tiles).
This observation can be formalized through the formula
$$ \mu_j(k)=\mu_{j-1}(k)-\widetilde u_j(k)+\frac14\sum_{{i\in\{1,2\}}\atop{\sigma\in\{-1,+1\}}}\widetilde u_j(k+\sigma e_i).$$
This leads to
\begin{equation*}\begin{split}&
\mu_j(k)+ u_j(k)-\frac14\sum_{{i\in\{1,2\}}\atop{\sigma\in\{-1,+1\}}} u_j(k+\sigma e_i)\\
=\;&
\mu_j(k)+\widetilde u_j(k)-\frac14\sum_{{i\in\{1,2\}}\atop{\sigma\in\{-1,+1\}}} \widetilde u_j(k+\sigma e_i)
+u_{j-1}(k)-\frac14\sum_{{i\in\{1,2\}}\atop{\sigma\in\{-1,+1\}}} u_{j-1}(k+\sigma e_i)\\
=\;&\mu_{j-1}(k)+ u_{j-1}(k)-\frac14\sum_{{i\in\{1,2\}}\atop{\sigma\in\{-1,+1\}}} u_{j-1}(k+\sigma e_i)
\end{split}\end{equation*}
and therefore, by iteration, for each~$\ell\in\{0,\dots,j\}$,
\begin{equation*}\begin{split}&
\mu_j(k)+ u_j(k)-\frac14\sum_{{i\in\{1,2\}}\atop{\sigma\in\{-1,+1\}}} u_j(k+\sigma e_i)\\
=\;&\mu_{j-\ell}(k)+ u_{j-\ell}(k)-\frac14\sum_{{i\in\{1,2\}}\atop{\sigma\in\{-1,+1\}}} u_{j-\ell}(k+\sigma e_i).
\end{split}\end{equation*}
{F}rom this, we arrive at
\begin{equation*}\begin{split}&
\mu_j(k)+ u_j(k)-\frac14\sum_{{i\in\{1,2\}}\atop{\sigma\in\{-1,+1\}}} u_j(k+\sigma e_i)\\
=\;&\mu_{0}(k)+ u_{0}(k)-\frac14\sum_{{i\in\{1,2\}}\atop{\sigma\in\{-1,+1\}}} u_{0}(k+\sigma e_i)\\=\;&\mu_0(k),
\end{split}\end{equation*}
which can be seen as a ``conservation law for the rubbish''
(that is, the initial rubbish~$\mu_0$ cannot disappear).

Hence, using the notation for the discrete Laplacian introduced in~\eqref{DISCRELEA}, we can write that
\begin{equation}\label{ULIbar-ZER-67}
\mu_j(k)-\frac14{\mathcal{L}} u_j(k)
=\mu_{0}(k).\end{equation}
This and~\eqref{ULIbar} also allow us to write that
\begin{equation}\label{EXILIMU}
\mu_{0}(k)+\frac14{\mathcal{L}} u_\infty(k)
=\mu_{0}(k)+\lim_{j\to+\infty}\frac14{\mathcal{L}} u_j(k)
=\lim_{j\to+\infty}\mu_j(k)=:\mu_\infty(k),
\end{equation}
thus ensuring that the final distribution of rubbish is indeed well-defined.

This observation also entails that
\begin{equation}\label{CLEPOSSI}
{\mbox{the streetsweeper can always find a successful strategy to clean any given floor,}}
\end{equation}
namely it is always possible to select the sequence of tiles to clean such that~$\mu_\infty(k)=0$ for all~$k\in\Omega$.
As a matter of fact, it suffices for the streetsweeper to
\begin{equation}\label{GOOD:STRAT}
{\mbox{pick a strategy that cleans each tile infinitely many times.}}\end{equation}
This means that for every~$k\in\Omega$ there exists a sequence~$j^{(k)}_m$
such that~$j^{(k)}_m\to+\infty$ as~$m\to+\infty$ and~$\mu_{j^{(k)}_m}(k)=0$
(note indeed that this condition follows from cleaning the tile centered at~$k$ in the~$j^{(k)}_m$th step of the process).
It thus follows from the existence of the limit in~\eqref{EXILIMU} that
$$\mu_\infty(k)=\lim_{m\to+\infty}\mu_{j^{(k)}_m}(k)=0,$$
which establishes~\eqref{CLEPOSSI}, as desired.

Let us now go back to the discrete Laplacian and make the following\footnote{Observations as in~\eqref{MAXPLEDISCERET}
will naturally lead to the development of the theory of Maximum Principles \index{Maximum Principle}
for elliptic partial differential equations. See the forthcoming Section~\ref{MAXPLEDISCERET:SE} for more details
on this topic.} observation:
\begin{equation}\label{MAXPLEDISCERET}
\begin{split}&
{\mbox{if~$w:\Z^2\to\R$ is such that~${\mathcal{L}}w(k)\le0$ for all~$k\in\Omega$
and~$w(k)\ge0$ for all~$k\in\Z^2\setminus\Omega$,}}\\ &{\mbox{then~$w(k)\ge0$ in~$\Omega$.}}\end{split}
\end{equation}
Indeed, suppose not and pick~$k_\star\in\Omega$ such that
$$ \min_{k\in\Omega}w(k)=w(k_\star)<0.$$
It thus follow that
\begin{equation} \label{MAXPLEDISCERET:SE2}
0\ge{\mathcal{L}} w(k_\star)=\sum_{i=1}^2\Big( w(k_\star+e_i)+w(k_\star-e_i)-2w(k_\star)\Big).
\end{equation}
Denoting by~$\ell$ the number of tiles that are adjacent to the one centered at~$k_\star$ and belong to~$\Omega$,
this leads to
\begin{eqnarray*}
4w(k_\star)\ge \sum_{{{i\in\{1,2\}}\atop{\sigma\in\{-1,+1\}}}\atop{k_\star+\sigma e_i\in\Omega}} w(k_\star+\sigma e_i)
+ \sum_{{{i\in\{1,2\}}\atop{\sigma\in\{-1,+1\}}}\atop{k_\star+\sigma e_i\in\Z^2\setminus\Omega}} w(k_\star+\sigma e_i)
\ge
\ell w(k_\star)+0,
\end{eqnarray*}
and therefore
\begin{equation}\label{sowivr496b897546jilkiujhygtfd}
(4-\ell)w(k_\star)\ge0.\end{equation}
Notice that~$\ell\le4$. If~$\ell<4$ then we would deduce from~\eqref{sowivr496b897546jilkiujhygtfd}
that~$w(k_\star)\ge0$, thus obtaining a contradiction. 

As a consequence, we conclude that~$\ell=4$, that is the tile centered at~$k_\star$ is surrounded by tiles belonging
to~$\Omega$. As a result, we infer from~\eqref{MAXPLEDISCERET:SE2}
that~$w(k_\star+e_i)=w(k_\star-e_i)=w(k_\star)$ for all~$i\in\{1,2\}$.
Iterating this argument, we obtain that~$k_\star\pm m e_i\in\Omega$ for every~$i\in\{1,2\}$ and every~$m\in\N$,
thus contradicting the finiteness of~$\Omega$, and this contradiction completes the proof of~\eqref{MAXPLEDISCERET}.

Among the many interesting consequences of~\eqref{MAXPLEDISCERET},
by creatively exploiting an argument\footnote{Other approaches
are possible as well, including a finite dimensional minimization argument.
As a matter of fact, to establish~\eqref{SWEE}, only the uniqueness
result in~\eqref{EXUNIDISCREESLA} will be exploited, but, given the elegance
of the argument provided, we decided to state and prove~\eqref{EXUNIDISCREESLA}
in its complete form. Also, we will use the existence result in~\eqref{EXUNIDISCREESLA}
in the forthcoming equation~\eqref{POISSODISC845}.

Comparing the terminology with that of footnote~\ref{FOO:POIEQUATYHN}
on page~\pageref{FOO:POIEQUATYHN}, one can consider~\eqref{EXUNIDISCREESLA}
as a discrete version of the \index{Poisson equation}
Poisson equation.} in linear algebra, we deduce an interesting
existence and uniqueness result\footnote{For a corresponding 
existence and uniqueness result about the classical Laplace operator, we will have to introduce more sophisticated
mathematical tools and patiently wait for
Corollary~\ref{S-coroEXIS-M023}.}
for the discrete Laplacian. Namely, for every~$f:\Omega\to\R$ and every~$g:\Z^2\setminus\Omega\to\R$
\begin{equation}\label{EXUNIDISCREESLA}\begin{split}
{\mbox{there exists a unique~$U:\Z^2\to\R$ such that}}\\
\begin{dcases}
{\mathcal{L}} U(k)=f(k)&{\mbox{ for all }}k\in\Omega,\\
U(k)=g(k)&{\mbox{ for all }}k\in\Z^2\setminus\Omega.
\end{dcases}\end{split}
\end{equation}
To check this, we observe that, without loss of generality, we can suppose that
\begin{equation}\label{EXUNIDISCREESLA-HOM-90}
{\mbox{$g(k)=0$
for all~$k\in\Z^2\setminus\Omega$,}}\end{equation} since solutions of~\eqref{EXUNIDISCREESLA}
correspond to solutions of
\begin{equation}\label{EXUNIDISCREESLA-HOM}
\begin{dcases}
{\mathcal{L}} \widehat U(k)=\widehat f(k)&{\mbox{ for all }}k\in\Omega,\\
\widehat U(k)=0&{\mbox{ for all }}k\in\Z^2\setminus\Omega
\end{dcases}
\end{equation}
under the transformations defined, for all~$k\in\Z^2$,
by~$\widehat U(k):=U(k)-G(k)$
and~$\widehat f(k):=f(k)-{\mathcal{L}}G(k)$,
where
$$
G(k):=\begin{dcases} g(k) & {\mbox{ if }}k\in\Z^2\setminus\Omega,\\
0&{\mbox{ otherwise.}}\end{dcases}$$
Thus, we focus on the establishment of an existence and uniqueness theory of~\eqref{EXUNIDISCREESLA}
under the additional information in~\eqref{EXUNIDISCREESLA-HOM-90}.
To achieve this, 
the main observation is that, after setting~$U(k):=g(k)=0$ for all~$k\in\Z^2\setminus\Omega$,
to solve~\eqref{EXUNIDISCREESLA} one only needs to determine the values of~$U$ at the finitely many points~$k_1,\dots,k_N$
of~$\Omega$ (here we are denoting by~$N$ the finite cardinality of~$\Omega$). Therefore, since the operator~${\mathcal{L}} $ is linear,
we can consider the vector~$X:=(U(k_1),\dots,U(k_N))\in\R^N$
and write the first equation in~\eqref{EXUNIDISCREESLA} as
$$AX=Y,$$ for a suitable
$N\times N$ matrix~$A$ and a suitable vector~$Y\in\R^N$.
We notice that~$Y$ depends only the values of~$f$, which are given once and for all:
more explicitly, in light of~\eqref{EXUNIDISCREESLA-HOM-90}, for all~$k\in\Omega$, one can
recast the first equation in~\eqref{EXUNIDISCREESLA} in the form
\begin{equation*}
\begin{split}
\sum_{i=1}^2\Big( U(k+e_i)\chi_\Omega(k+e_i)
+U(k-e_i)\chi_\Omega(k-e_i)-2U(k)\Big)=
{\mathcal{L}} U(k)=f(k)
\end{split}
\end{equation*}
and we stress that the latter term is given in dependence of~$f$
and can be considered as an $N$-dimensional array~$Y$ (with one entry for each~$k\in\Omega$),
while the first term is a linear map applied to the vector~$X$.

Consequently, to prove that the first equation in~\eqref{EXUNIDISCREESLA} admits one and only one solution,
we need to check that the matrix~$A$ is invertible.
To this end, by the Fundamental Theorem of Linear Algebra (see e.g.~\cite[Proposition~1.3.6]{zbMATH07269795}), it suffices
to check that 
\begin{equation}\label{MAXPLEDISCERET-LININ}
{\mbox{the linear map associated to~$A$ is injective.}}\end{equation} For this, suppose that there exists~$X_\star=(X_{\star,1},\dots,X_{\star,N})\in\R^N$ such that
\begin{equation}
AX_\star=0.\end{equation}
By letting~$U_\star(k_j):=X_{\star,j}$ for every~$j\in\{1,\dots,N\}$, 
as well as~$U_\star(k):=0$ for every~$k\in\Z^2\setminus\Omega$,
we have constructed a function~$U_\star:\Z^2\to\R$ satisfying
\begin{equation*}
\begin{dcases}
{\mathcal{L}}U_\star(k)=0&{\mbox{ for all }}k\in\Omega,\\
U_\star(k)=0&{\mbox{ for all }}k\in\Z^2\setminus\Omega.
\end{dcases}
\end{equation*}
By applying~\eqref{MAXPLEDISCERET} to both~$U_\star$ and~$-U_\star$, we thus
conclude that~$U_\star(k)=0$ for each~$k\in\Z^2$, hence~$X_\star=0$,
which finishes the proof of~\eqref{MAXPLEDISCERET-LININ}.

Now, with this preliminary work, we can complete the proof of the surprising result in~\eqref{SWEE}
by arguing as follows. We consider two different
successful cleaning strategies (the existence of a cleaning strategy being warranted by~\eqref{CLEPOSSI})
and we denote by~$\mu_j$ and~$\mu_j'$ the corresponding rubbish densities after~$j$ steps
of each strategy. 
Note that~$\mu_0'=\mu_0$ since the initial rubbish distribution is given.
Since the strategies are successful, we know that~$\mu_\infty(k)=0=\mu'_\infty(k)$
for every~$k\in\Omega$ (recall that these limit densities are well-defined owing to~\eqref{EXILIMU}).

Thus, to prove~\eqref{SWEE}, we need to show that
\begin{equation}\label{SWEEpre}
{\mbox{$\mu_\infty(k)=\mu'_\infty(k)$ for every $k\in\Z^2$.}}
\end{equation}
For this, we let~$u_j$ and~$u'_j$ be the functions accounting for the total rubbish taken out
after $j$ steps of the cleaning process,
corresponding to the two cleaning strategies (as formalized in~\eqref{PERmuinftyk}).
Thanks to~\eqref{ULIbar}, we can also consider their limit configurations~$ u_\infty$ and~$u_\infty'$.
We observe that, as a consequence of~\eqref{ULIbar-ZER-67}, for all~$k\in\Omega$,
\begin{equation}\label{eforalinZetmimega}
-\frac14{\mathcal{L}} u_\infty(k)=
\mu_\infty(k)-\frac14{\mathcal{L}} u_\infty(k)
=\lim_{j\to+\infty}\left(
\mu_j(k)-\frac14{\mathcal{L}} u_j(k)\right)
=\mu_{0}(k)\end{equation}
and similarly
\[ -\frac14{\mathcal{L}} u_\infty'(k)=\mu_{0}(k).\]
These observations and~\eqref{ULIbar-ZER} yield that both~$u_\infty$ and~$u_\infty'$ are solutions
of
\begin{equation*}
\begin{dcases}
{\mathcal{L}}U(k)=-4\mu_0(k)&{\mbox{ for all }}k\in\Omega,\\
U(k)=0&{\mbox{ for all }}k\in\Z^2\setminus\Omega.
\end{dcases}
\end{equation*}
Hence, by the uniqueness
result in~\eqref{EXUNIDISCREESLA}, it follows that~$u_\infty=u_\infty'$. {F}rom this and~\eqref{EXILIMU}
we deduce that
\[ \mu_\infty(k)=
\mu_{0}(k)+\frac14{\mathcal{L}} u_\infty(k)=\mu_{0}(k)+\frac14{\mathcal{L}} u_\infty'(k)=\mu_\infty'(k).
\] This completes the proof of~\eqref{SWEEpre},
and therefore of the surprising claim in~\eqref{SWEE}.\medskip

\begin{figure}
                \centering
                \includegraphics[width=.65\textwidth]{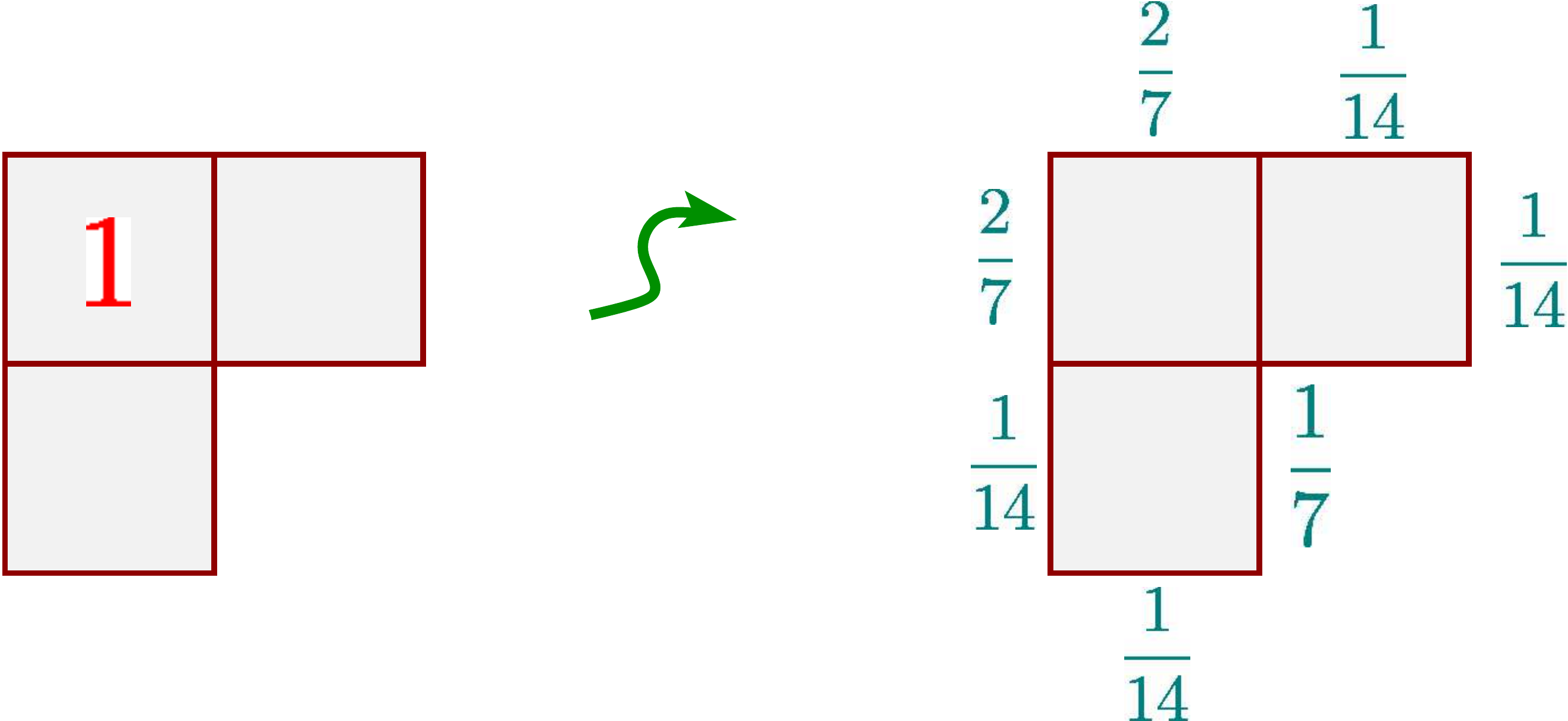}
        \caption{\sl A simple example workable by pencil and paper.}\label{1435-CLE1ANSESwerEqXATWO}
\end{figure}

It is now very tempting to switch on our personal computer and write down a short code
and see some simulations about such a sweeping process, nicely checking that different successful cleaning strategies
do produce the same final distribution, so to have a ``numerical confirmation''
for the rigorous result obtained in~\eqref{SWEE}. It is however also useful to
work out at least a simple example just by pencil and paper:
after all, according to Albert Einstein,
``computers are incredibly fast, accurate, and stupid. Human beings are incredibly slow, inaccurate, and brilliant. Together they are powerful beyond imagination''.

\begin{figure}
                \centering
                \includegraphics[width=.95\textwidth]{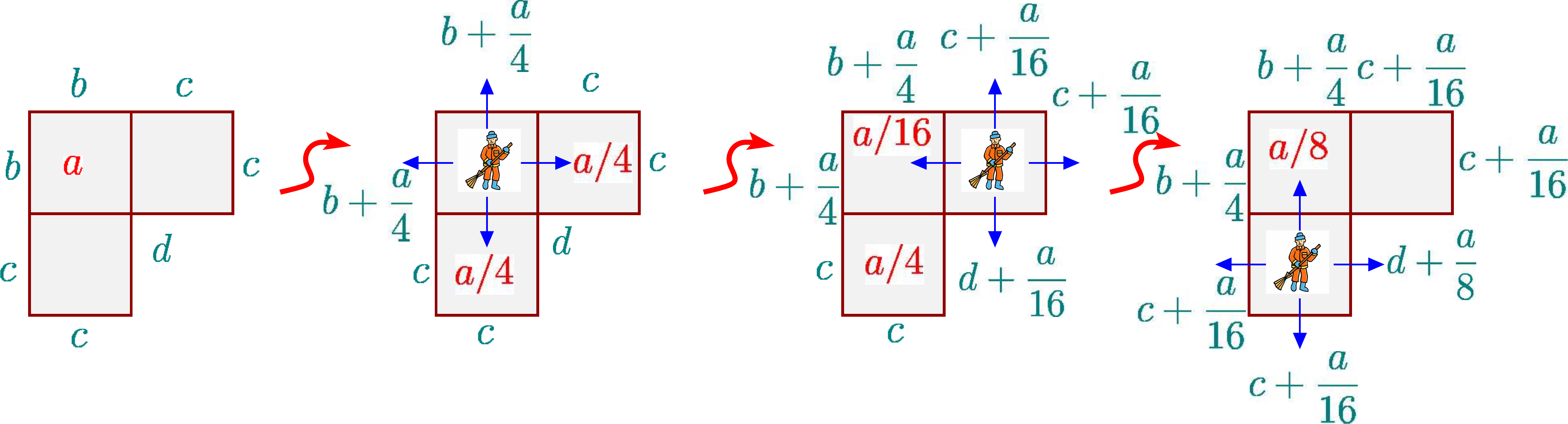}
        \caption{\sl A cleaning strategy in a simple example.}\label{CLE1ANSESwerEqXATWO}
\end{figure}

\begin{figure}
                \centering
                \includegraphics[width=.65\textwidth]{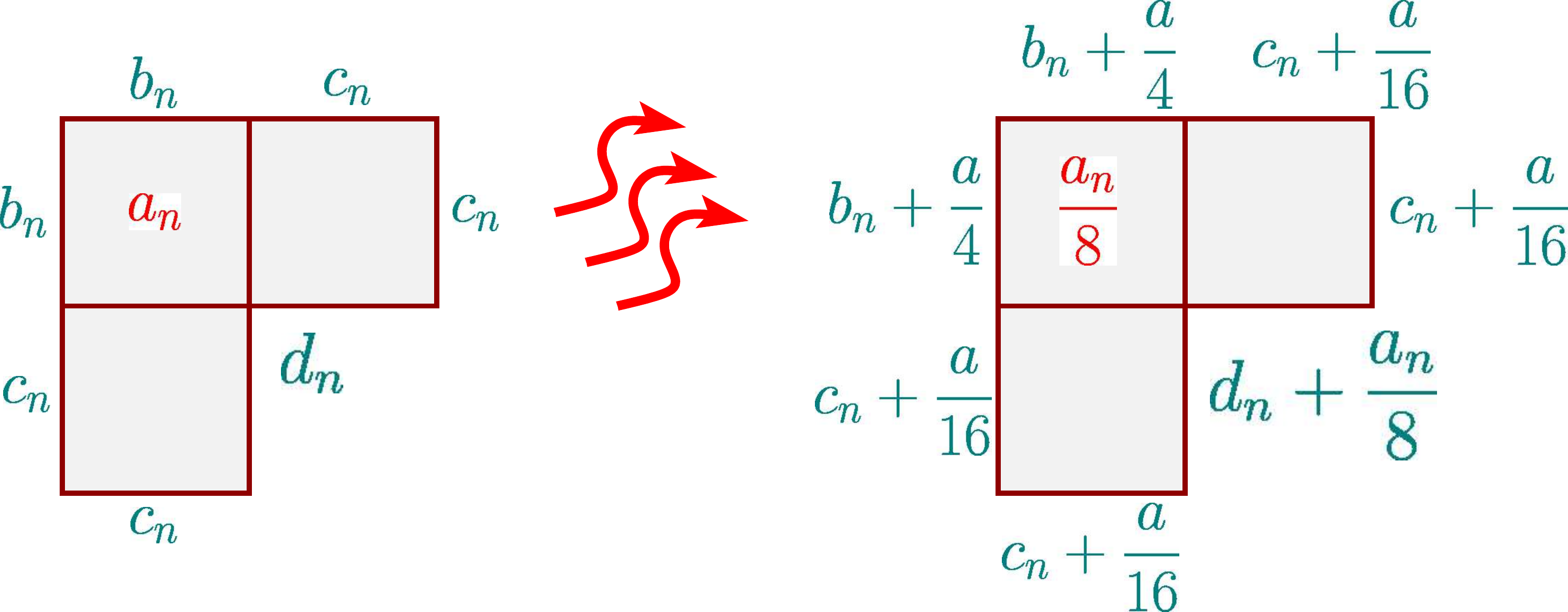}
        \caption{\sl Summary of the cleaning strategy of Figure~\ref{CLE1ANSESwerEqXATWO}.}\label{2CLE1ANSESwerEqXATWO}
\end{figure}

A simple and explicit example
presents the additional boon to clearly relate the final distribution
to the geometry of the domain. For instance, leaving no room for doubt, one can consider the case in which~$\Omega:=\{(0,0),$
$(0,1)$, $(1,1)\}$ and
$$\mu_0(k):=\chi_{\{(0,1)\}}(k)=\begin{dcases}
1&{\mbox{ if }}k=(0,1),\\
0& {\mbox{ otherwise.}}
\end{dcases}$$
In this case, just by pencil and paper, we see that
\begin{equation}\label{POI:AMAIKJM0}
\mu_\infty(k)=\begin{dcases}
2/7 & {\mbox{ if }}k\in\{(-1,1),\,(0,2)\},\\
1/7 & {\mbox{ if }}k=(1,0),\\
1/14 & {\mbox{ if }}k\in\{(-1,0),\,(0,-1),\,(2,1),\,(1,2)\},
\end{dcases}
\end{equation}
see Figure~\ref{1435-CLE1ANSESwerEqXATWO}. To check~\eqref{POI:AMAIKJM0},
one can follow the cleaning strategy depicted in Figure~\ref{CLE1ANSESwerEqXATWO}
(for peace of mind, one can recall that this is indeed a successful cleaning strategy,
thanks to~\eqref{GOOD:STRAT}, and any other successful cleaning strategy would produce, in view of~\eqref{SWEE}, the same
final distribution).
Combining together three steps of the strategy in Figure~\ref{CLE1ANSESwerEqXATWO}, one obtains the iterative situation
skecthed in Figure~\ref{2CLE1ANSESwerEqXATWO}:
namely, the corner square in~$\Omega$ is filled by an amount~$a_n$ of rubbish, with~$a_{n+1}=\frac{a_n}8$
and~$a_0=1$, the two squares outside~$\Omega$ adjacent to it
are filled by an amount~$b_n$ of rubbish, with~$b_{n+1}=b_n+\frac{a_n}4$ and~$b_0=0$, then we see
four squares outside~$\Omega$ filled by an amount~$c_n$ of rubbish, with~$c_{n+1}=c_n+\frac{a_n}{16}$
and~$c_0=0$, and finally one square
filled by an amount~$d_n$ of rubbish, with~$d_{n+1}=d_n+\frac{a_n}8$
and~$d_0=0$. 

As a consequence, we find that~$a_n=\frac{1}{8^n}$ (confirming that~$a_n\to0$ as~$n\to+\infty$,
hence providing a successful cleaning strategy), that
$$ b_n=b_0+\frac14\sum_{i=0}^{n-1}a_i
=\frac14\sum_{i=0}^{n-1}\frac1{8^i}
\longrightarrow\frac14\sum_{i=0}^{+\infty}\frac1{8^i}=\frac27,
$$
that
$$ c_n=c_0+\frac1{16}\sum_{i=0}^{n-1}a_i
=\frac1{16}\sum_{i=0}^{n-1}\frac1{8^i}
\longrightarrow\frac1{16}\sum_{i=0}^{+\infty}\frac1{8^i}=\frac1{14}
$$
and that
$$ d_n=d_0+\frac1{8}\sum_{i=0}^{n-1}a_i
=\frac1{8}\sum_{i=0}^{n-1}\frac1{8^i}
\longrightarrow\frac1{8}\sum_{i=0}^{+\infty}\frac1{8^i}=\frac1{7}.
$$
These observations give~\eqref{POI:AMAIKJM0}, as desired.

\medskip

This sweeping problem provides additional intuition for many other methodologies typical of elliptic partial differential equations. We give indeed another example of a remarkable construction that can be performed in this (apparently innocent) discrete setting.
Given~$\ell\in\Omega$, we let~$G_\ell:\Z^2\to\R$ be the function~$u_\infty$
obtained in~\eqref{ULIbar} when the initial rubbish distribution has the form~$\mu_0(k):=\frac14\,\chi_{\{\ell\}}(k)$.
Then, we claim that the function
\begin{equation}\label{CONTECONT09} u(k):=\sum_{\ell\in\Omega} f(\ell)\,G_\ell(k)\end{equation}
is the unique solution\footnote{The function~$G_\ell(k)$ can be seen
as a discrete counterpart of an object which is called ``Green Function'' \index{Green Function}
in the continuous setting and which will be addressed in further detail in Section~\ref{Green Function:SECT}.
In particular, one can compare~\eqref{CONTECONT09} and~\eqref{CONTECONT10} here
with the continuous framework described in footnote~\ref{CONTECONT09F}
on page~\pageref{CONTECONT09F}.}
of
\begin{equation}\label{CONTECONT10}
\begin{dcases}
-{\mathcal{L}}u(k)=f(k)&{\mbox{ for all }}k\in\Omega,\\
u(k)=0&{\mbox{ for all }}k\in\Z^2\setminus\Omega.
\end{dcases}
\end{equation}
To check this, it suffices to prove that~$u$ in~\eqref{CONTECONT09} is a solution of~\eqref{CONTECONT10},
since the uniqueness claim follows from~\eqref{EXUNIDISCREESLA}.

Also, $G_\ell(k)=0$ for all~$k\in\Z^2\setminus\Omega$, owing to~\eqref{ULIbar-ZER}, and consequently~$u(k)=0$ for all~$k\in\Z^2\setminus\Omega$.

Furthermore, by~\eqref{eforalinZetmimega}, if~$k\in\Omega$ then
$$ -{\mathcal{L}} G_\ell(k)=\chi_{\{\ell\}}(k)$$
and therefore, by~\eqref{CONTECONT09},
$$ {\mathcal{L}} u(k)=\sum_{\ell\in\Omega} f(\ell)\,{\mathcal{L}} G_\ell(k)=
-\sum_{\ell\in\Omega} f(\ell)\,\chi_{\{\ell\}}(k)=-f(k).$$
These considerations establish the validity of~\eqref{CONTECONT10}, as desired.\medskip

Along these lines, exploiting the existence and uniqueness theory of~\eqref{EXUNIDISCREESLA},
given~$\ell\in\Z^2\setminus\Omega$ one can consider the unique function~$P_\ell:\Z^2\to\R$ such that
\begin{equation}\label{POISSODISC845}
\begin{dcases}
{\mathcal{L}}P_\ell (k)=0&{\mbox{ for all }}k\in\Omega,\\
P_\ell(k)=\chi_{\{\ell\}}(k)&{\mbox{ for all }}k\in\Z^2\setminus\Omega.
\end{dcases}
\end{equation}
This function also provides a useful representation formula, somewhat related to~\eqref{CONTECONT09}
and~\eqref{CONTECONT10}: more specifically, we have that
the function
\begin{equation}\label{CONTECONT09POI} v(k):=\sum_{\ell\in\Z^2\setminus\Omega} g(\ell)\,P_\ell(k)\end{equation}
is the unique solution\footnote{Comparing with the framework in Section~\ref{SEC:POIKER-S},
one can consider~$P_\ell(k)$ as a discrete version of the Poisson Kernel. \index{Poisson Kernel}
In particular, equations~\eqref{CONTECONT09POI} and~\eqref{CONTECONT10POI}
here can be seen as a discrete counterpart of Theorem~\ref{POIBALL1}
in the continuous setting.} of
\begin{equation}\label{CONTECONT10POI}
\begin{dcases}
{\mathcal{L}}v(k)=0&{\mbox{ for all }}k\in\Omega,\\
v(k)=g(k)&{\mbox{ for all }}k\in\Z^2\setminus\Omega.
\end{dcases}
\end{equation}
For this, once again, in light of the uniqueness result in~\eqref{EXUNIDISCREESLA}, it is sufficient to check
that~$v$ is a solution of~\eqref{CONTECONT10POI}. But this follows by linearity, since, combining~\eqref{POISSODISC845}
and~\eqref{CONTECONT09POI} we see that when~$k\in\Omega$
$$ {\mathcal{L}}v(k)=\sum_{\ell\in\Z^2\setminus\Omega} g(\ell)\,{\mathcal{L}} P_\ell(k)=0$$
and when~$k\in\Z^2\setminus\Omega$
$$ v(k)=\sum_{\ell\in\Z^2\setminus\Omega} g(\ell)\,\chi_{\{\ell\}}(k)=g(k).$$

\subsection{Image processing}\label{I2K243S3M:0oikrfg-4}

A common operation nowadays is to process a given (say, for simplicity, grayscale) image in order to remove noise or enhance some visual effect. 

To describe this procedure, one can denote by~$u$ the brightness of the image at a certain pixel of the screen and we may think that~$u=u(x,t)$, where~$x$ lies in some domain of~$\R^2$ (or even~$\R^n$, for the sake of generality) that represents the position on the screen (which we consider, in the limit, as a continuum of pixels) and~$t$ is time. With this notation, the initial image~$u(x,0)$ undergoes some evolution, with the aim of improving its aspects, according to some standards that we can choose.

The simplest possibility would be to consider a heat equation for the brightness function~$u$. This would have the advantage of possibly removing the imperfection from the image manifested, e.g., in abrupt variations of its brightness caused by impurities. Also, due to the smoothing effect of the heat equation, such a process does not introduce additional spurious details.

\begin{figure}
  \centering
      \includegraphics[height=.2\textheight]{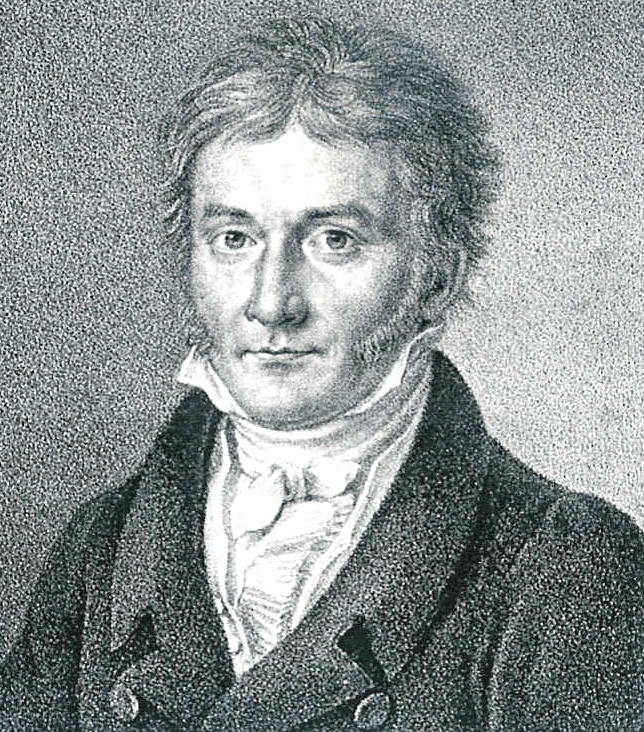}$\quad$
    \includegraphics[height=.2\textheight]{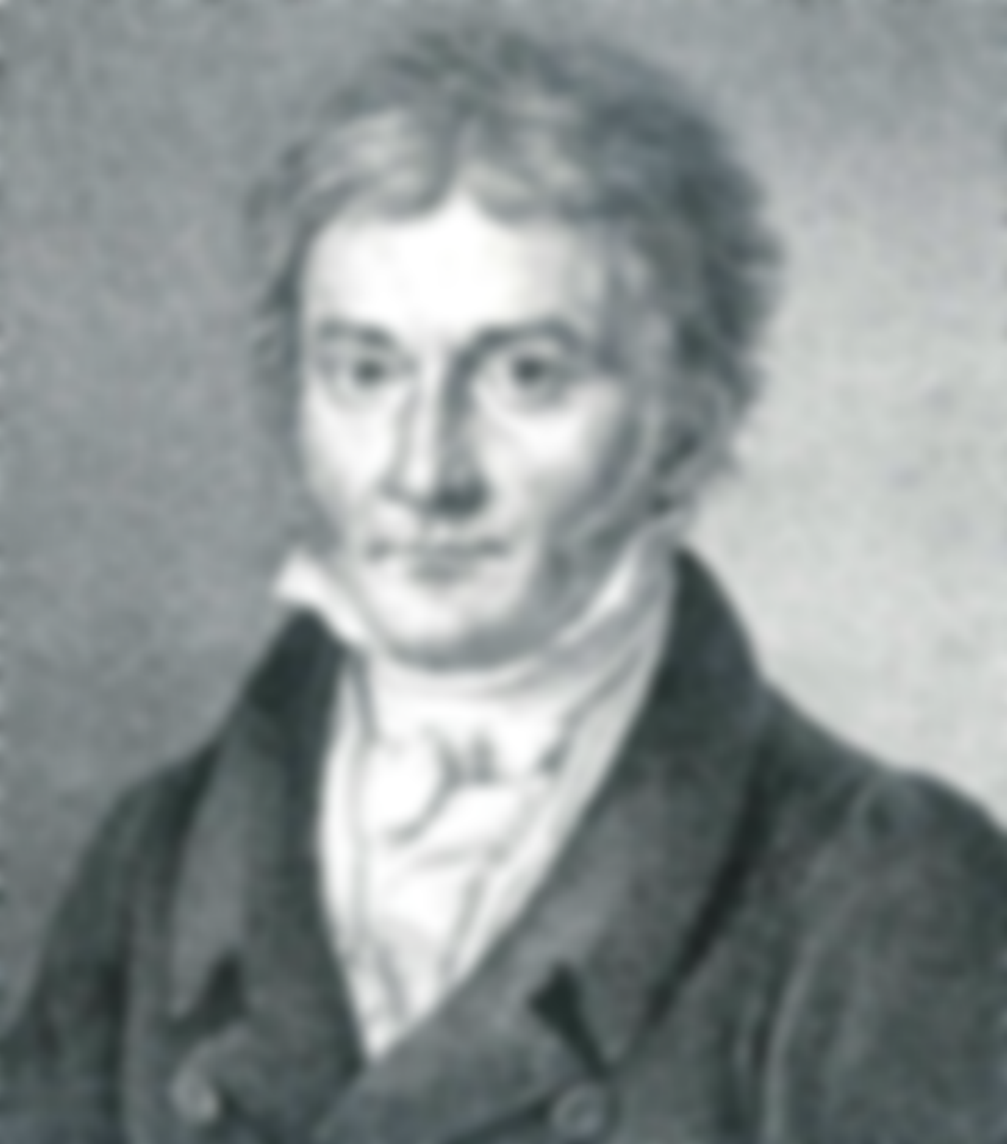}
 \caption{\sl Portrait of Gau{\ss} 
by artist Siegfried Detlev Bendixen published in Astronomische Nachrichten in 1828
(Public Domain images from
 Wikipedia) and a Gau{\ss}ian filtered image of it (obtained via Mathematica).}\label{VILORi-1nsFI0231}
\end{figure}

These benefits in the use of the heat equation, however, present an inconvenient drawback caused by the coarsening of the image resolutions. While this downside is perhaps unavoidable after all, since some kind of averaging effect is necessary precisely to get rid of noises and impurities, the most significant hiccup is that this blurring does not necessarily respect the natural boundaries of the original image, see
Figure~\ref{VILORi-1nsFI0231} for a clear sketch\footnote{The choice of using a portrait of Gau{\ss} here to check the different image processing procedures is completely arbitrary. For many years, research articles used for this scope a 512$\times$512px standard test version of model Lena Fors\'en (stage name Lenna, previously Lena Soderberg, born Sj\"o\"oblom). This standard image was a closeup detail of Lenna wearing a hat (the complete original image being a full length portrait of Lenna wearing a hat, published as the centerfold of the November 1972 issue of Playboy magazine). 

Allegedly, the Lenna standard image has become one of the most used images in computer history. After an initial concern about copyright infringement, it seems that over time Playboy has implicitly decided to overlook
the wide and free use of the Lenna image for scientific purposes --
however, at the beginning,
Playboy tried to restrain the use of its copyrighted material
e.g. by writing a letter to the journal Optical Engineering, see {\tt https://doi.org/10.1117/12.60707},
whose content was: ``It has come to our attention that you have used a
portion of the centerfold photograph of our November 1972 PLAYBOY PLAYMATE OF THE
MONTH Lenna Sj\"o\"oblom, in your July 1991 issue of
Optical Engineering magazine... Playboy Enterprises, Inc., the publisher of PLAYBOY magazine,
owns the copyright in and to this photograph. As fellow publishers, we're sure you understand
the need for us to protect our proprietary rights...''; later on, anyway,
Playboy seemed to be pleased about this phenomenon, see e.g.
{\tt http://www.lenna.org/playboy\_backups/index.html}
and {\tt http://www.lenna.org/playboy\_backups/lena.html},
since the story was included in ``The World History of Playboy''
and the magazine seemed proud of the perspective that ``the image of this Playboy Playmate can remain the standard reference image for comparing compression technologies into the 21st century''.

Due to the popularity of this image, Lena Fors\'en was invited as a guest at the 50th annual Conference of the Society for Imaging Science and Technology in 1997 and as a guest of honor at the banquet of IEEE ICIP in 2015: on the latter occasion, she also delivered a speech and chaired the best paper award ceremony.

The use of the Lenna picture has however become controversial, often considered as being degrading to women and of detrimental impact on aspiring female students in computer science. Several scientific journals nowadays explicitly and strongly discourage the use of the Lenna image, others do not consider new submissions containing the image.

In January 2019, in an interview to the monthly American magazine Wired, Lena Fors\'en seemed to have declared to be ``really proud of that picture''
(``the only note of regret she expressed was that she wasn't better compensated'',
according to {\tt https://www.wired.com/story/finding-lena-the-patron-saint-of-jpegs/}), but
when asked ``if she had heard anything about the recent controversy around her image, she seemed alarmed at the thought that she could have a part in hurting or discouraging young women''.
In November of the same year Lena Fors\'en also took part to a short documentary titled Losing Lena,
by Australian cinematographer Anna Howard,
which aimed at galvanizing efforts to end the use of her image in technology research. On that occasion, Lena Fors\'en explicitly stated {\em``I retired from modeling a long time ago. It's time I retired from tech, too... Let's commit to losing me''}.

For the sake of completeness, we also mention that the Lenna picture was
shot by photographer Dwight Hooker
and
the ``unknown researcher''
who first scanned the Lenna image is sometimes (see e.g.
{\tt https://www.ee.cityu.edu.hk/$\sim$lmpo/len\-na/Len\-na97.html})
reported to be
William K. Pratt.

Also, it could be just a coincidence, but allegedly
the Lenna issue (i.e., November 1972) is said to be Playboy's best selling issue ever,
with more than~$7\times 10^6$ copies sold.} of this pitfall.\medskip

It would be instead more desirable to accomplish a reduction of image noise without removing significant content from the image, such as edges and details which are relevant for the interpretation of the image.\medskip

\begin{figure}
  \centering
      \includegraphics[height=.2\textheight]{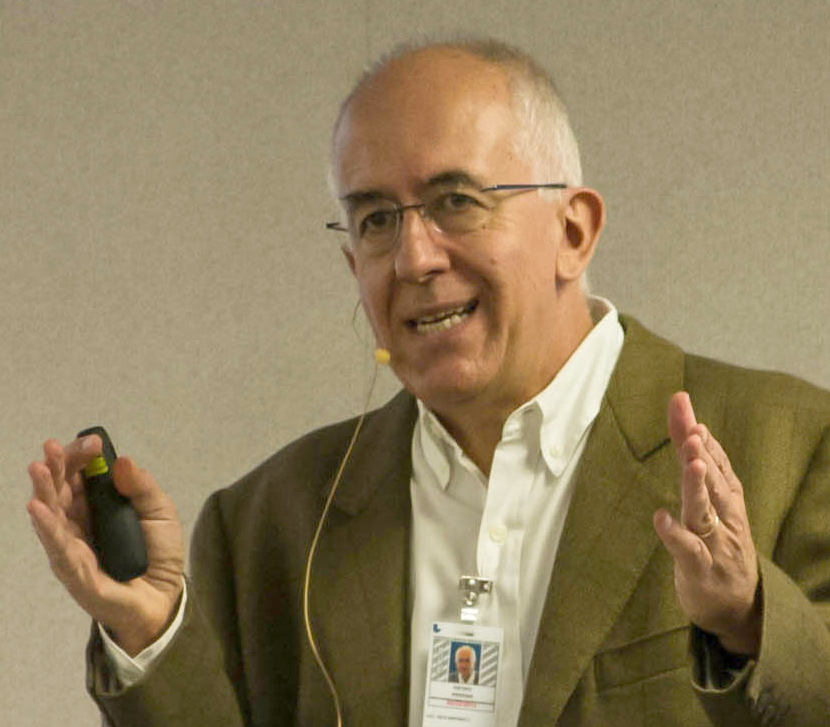}$\quad$
    \includegraphics[height=.2\textheight]{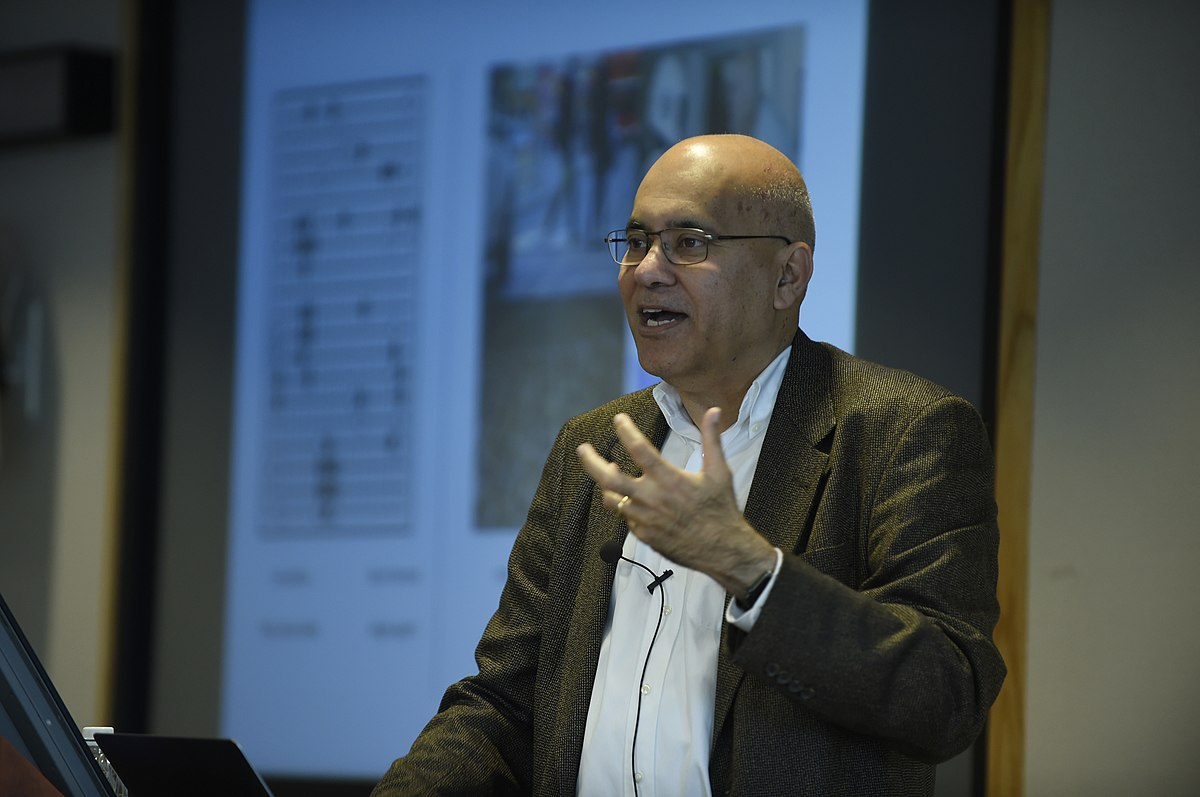}
 \caption{\sl Pietro Perona and Jitendra Malik in action
 (Public Domain images from
 Wikipedia).}\label{VIL:24EMALihnsFI}
\end{figure}

Though several methods have been proposed for this purpose, we do not intend here to give an exhaustive presentation of them, but rather to recall a classical method proposed by Pietro Perona and Jitendra Malik~\cite{56205PERMAL}
(see Figure~\ref{VIL:24EMALihnsFI}).
The gist of their technique consists in replacing the classical heat equation
with a ``shape-adapted'' smoothing process in which, roughly speaking,
the diffusion coefficient, instead of being a constant, is a function induced by the image itself.
This allows one to remove noise from digital images without blurring edges, since
one can somewhat encourage the smoothing effects within regions in which
the original image is already sufficiently smooth,
but rather suppress these diffusion effects across strong edges, along which
the brightness function changes abruptly.\medskip

To make this strategy clear, we remark that in a perfectly black and white image (without any gray tone),
we can model~$u$ as a step function taking values in~$\{0,1\}$. In this ideal situation,
the gradient of~$u$ vanishes both in the white regions (corresponding, say, to~$u=1$)
and in the black regions (where~$u=0$). The gradient is instead ``concentrated''
along the sharp separations between black and white areas, along which it formally attains
an infinite value.\medskip

For grayscale images, we can therefore interpret edges and details as thin lines along which
the gradient of the brightness function reaches its maximal norm.
The idea of the Perona-Malik method is thus to force the diffusion coefficient to be small
for large values of~$|\nabla u|$. If possible, it could be also convenient to have an evolution equation
with divergence structure (e.g., to interpret it as a gradient descent). The \index{Perona-Malik diffusion}
Perona-Malik equation is thereby a nonlinear diffusion equation of the form
\begin{equation}\label{MAL:equi10-1}
\partial_t u(x,t)=\div\Big( g\big(|\nabla u(x,t)|\big)\,\nabla u(x,t)\Big).
\end{equation}
While the choice~$g:=1$ gives back the standard heat equation, to preserve the
region boundaries of the images by suppressing diffusion at the sharp edges, the Perona-Malik method
consists in choosing~$g$ to be a nonnegative
monotonically decreasing function approaching zero at infinity. As models
for this type of functions, it is proposed in~\cite{56205PERMAL} to take, for instance, among the others,
\begin{equation}\label{PS:MOD23we} g(r) := \frac{1}{ 1 +(r/K)^2},\end{equation}
for some~$K>0$.

The Perona-Malik method manages to maintain some of the advantage of the classical heat equation
(for example, a suitable maximum principle still avoids the creation of
new features in the image when passing from fine to coarse scale, see~\cite[Appendix~A]{56205PERMAL}), however some instabilities
can arise for large gradients (though we skate around these difficulties at this level).\medskip

To convince ourselves of the fact that
the Perona-Malik method tends to preserve, and possibly sharpen, the brightness
edges, we can assume for simplicity that one of these high gradient edges passes through the origin,
with gradient vector in the direction of~$e_n$.
At a very small scale, we can suppose that the image presents an edge orthogonal to~$e_n$,
whence we approximate the brightness function~$u$
near the origin with a one-dimensional function depending only on~$x_n$ and with Taylor expansion
$$ u(x)=c+A x_n+\frac{ Mx_n^2}2+\frac{L x_n^3}{6}+o(x_n^3),$$
for some~$c$, $M$, $L\in\R$ and a ``very large'' slope~$A>0$.

We observe that the condition that~$|\nabla u|$ is maximal (hence critical) at the origin entails that
$$ 0=\partial_n\left.\frac{|\nabla u(x)|^2}2\right|_{x=0}
=\nabla u(0)\cdot\nabla \partial_n u(0)=AM,
$$
therefore~$M=0$ and accordingly, near the origin,
$$ u(x)=c+Ax_n+\frac{L x_n^3}{6}+o(x_n^3).$$
Furthermore, since~$|\nabla u|$ is actually maximal at the origin we have that
$$ 0\ge\partial^2_n\left.\frac{|\nabla u(x)|^2}2\right|_{x=0}
=|\nabla\partial_nu(0)|^2+\nabla u(0)\cdot\nabla\partial^2_nu(0)=AL,
$$
hence~$L\le0$, and suppose to be in a ``nondegenerate'' situation, hence
\begin{equation}\label{PS:MOD23we2}L<0.\end{equation}
Thus, setting~$\phi(r):=g(r)r$, we deduce from~\eqref{MAL:equi10-1} that
$$ \partial_t u=\div\left( \phi\big(|\nabla u|\big)\,\frac{\nabla u}{|\nabla u|}\right)
=\partial_n\Big( \phi\big(\partial_n u\big)\Big).
$$
As a result, setting~$U:=\partial_nu=A+\frac{L x_n^2}{2}+o(x_n^2)$, we see that
$$ \partial_t U=\partial_n^2 \big( \phi(U)\big)=\phi''(U)(\partial_nU)^2+\phi'(U)\partial_n^2U=\phi'(U)L+o(1).
$$
In the model case~\eqref{PS:MOD23we}, we have that~$\phi(r) = \frac{r}{ 1 +(r/K)^2}$ and consequently, in this situation,
$$ \partial_t U=\frac{K^2 L (K^2 - U^2)}{(K^2 + U^2)^2}+o(1)>0$$
when~$U(0)=A>K$, due to~\eqref{PS:MOD23we2}: this shows that if the slope of the edge
is large enough, then
the slope itself will increase with
time and the edge becomes sharper (for small slope, instead, the smoothing
effect of diffusion may prevail).\medskip

\begin{figure}
  \centering
      \includegraphics[height=.2\textheight]{ORIFI1.jpg}$\quad$
    \includegraphics[height=.2\textheight]{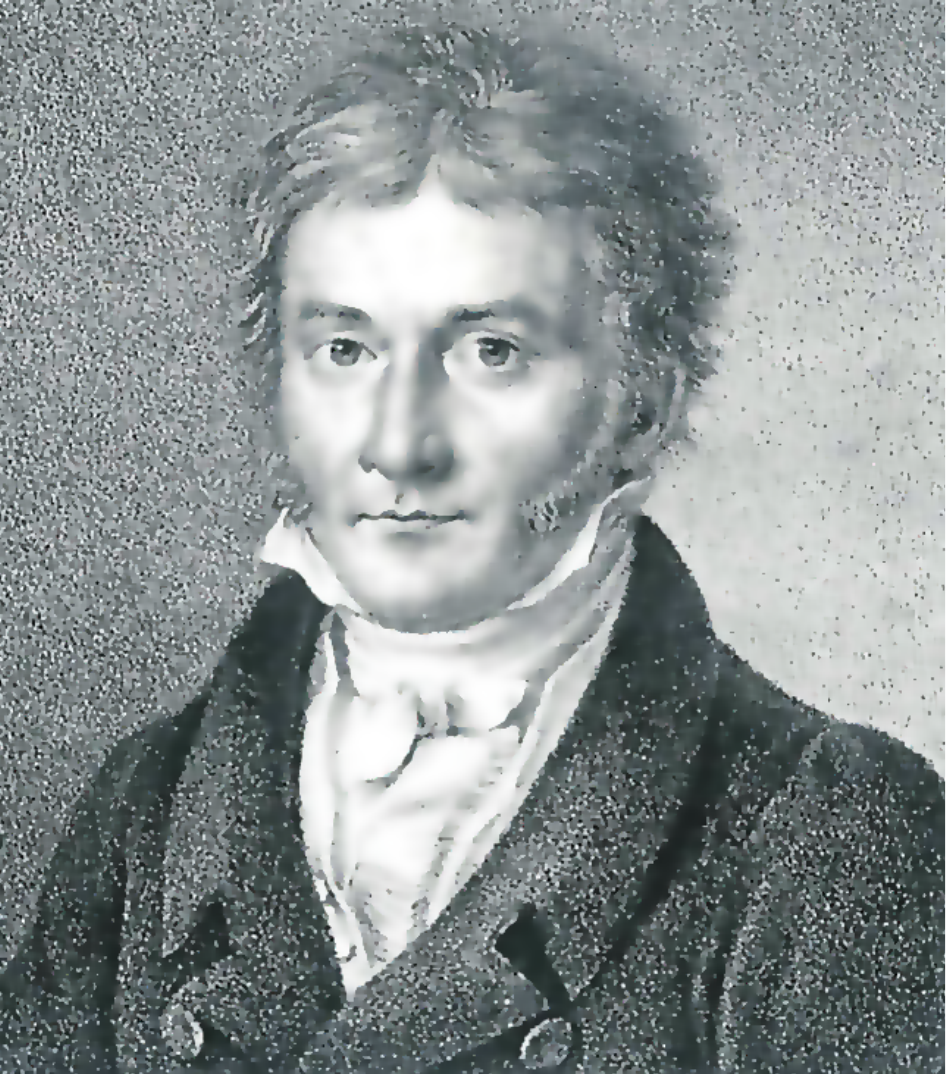}
 \caption{\sl Portrait of Gau{\ss} 
by artist Siegfried Detlev Bendixen published in Astronomische Nachrichten in 1828
(Public Domain images from
 Wikipedia) and a Perona-Malik filtered image of it (obtained via Mathematica).}\label{VILORi-1nsFI0232}
\end{figure}

The remarkable effect of sharpening edges while denoising smooth regions is depicted in
Figure~\ref{VILORi-1nsFI0232} (to be compared with Figure~\ref{VILORi-1nsFI0231}).\medskip

We remark that the Perona-Malik equation~\eqref{MAL:equi10-1} is also related to
(but structurally different from) the level set description of the mean curvature flow \index{mean curvature flow}
\begin{equation}\label{1-123OKS-234}
\partial_t u(x,t)=|\nabla u|\div\left( \frac{\nabla u(x,t)}{|\nabla u(x,t)|}\right),\end{equation}
see e.g.~\cite{MR1100206}.\medskip

To show, at least formally, that the level sets of a smooth and nondegenerate solution~$u$ of~\eqref{1-123OKS-234} evolve with normal velocity equal to their mean curvature, one can proceed as follows. We suppose that, say, the level set~$\{u(\cdot,t)=0\}$ corresponds to
a smooth and bounded hypersurface~$M_t$ given by an embedding~$x_t$. In this framework, 
we aim at showing that the normal velocity (pointing inward near convex portions) of~$M_t$
is equal to its mean curvature~$H_t$, namely that \index{mean curvature}
\begin{equation}\label{ROLLSTO-MECA}\nu_t\cdot
\partial_t x_t=-H_t,
\end{equation}
where~$\nu_t$ is the exterior normal to~$M_t=\partial\{u(\cdot,t)>0\}$ at the point identified by~$x_t$.

In order to check this, we pick a time~$t$ and a point, say the origin, on~$M_t$ and we set normal coordinates around it.
In this way, $-\frac{\nabla u(0,t)}{|\nabla u(0,t)|}=\nu_t(0)=e_n$ and, for every~$j\in\{1,\dots,n-1\}$, we have that~$\partial_j x_t\cdot e_n=0$.

As a result, denoting by~$\eta=(\eta_1,\dots,\eta_{n-1})$ the coordinates of the embedding~$x_t$
(hence~$\partial_jx_t=\frac{\partial x_t}{\partial\eta_j}$, and we assume that~$x_t(0)=0$), we obtain that
if~$\varpi$ is any smooth vector field such that~$|\varpi(x)|=1$ for all~$x\in\R^n$ and~$\varpi(0)=\nu_t(0)$,
then
$$0=\left.\partial_n \frac{|\varpi(x)|^2}{2}\right|_{x=0}= \varpi(0)\cdot\partial_n\varpi(0)=\partial_n\varpi_n(0)$$
and thus
\begin{equation}\label{djeiotyeuirtuerghreuig}\begin{split}
&\left.
\nabla\left(\frac{\nabla u(x,t)}{|\nabla u(x,t)|}\cdot\varpi(x)\right)\cdot\varpi(x)\right|_{x=0}=\sum_{k=1}^n\left.
\partial_n\left(\frac{\partial_k u(x,t)\,\varpi_k(x)}{|\nabla u(x,t)|}\right)\right|_{x=0}\\ &\qquad=
\frac{\partial_{nn} u(0,t)}{|\nabla u(0,t)|}-\frac{(\partial_n u(0,t))^2\partial_{nn}u(0,t)
}{|\nabla u(0,t)|^3}=
0
\end{split}\end{equation}
and as a byproduct\footnote{Notice indeed that the quantity in~\eqref{djeiotyeuirtuerghreuig}
corresponds to the normal component of the divergence of the normal vector field

A more precise description of the mean curvature will be given in~\eqref{MC}.
Here, we are essentially taking for granted that the mean curvature is the tangential divergence of the normal field,
a concept that will be better clarified in~\eqref{VB-d1}.
} we arrive at
$$ H_t(0)=-\left.\div\left( \frac{\nabla u(x,t)}{|\nabla u(x,t)|}\right)\right|_{x=0}.$$

\begin{figure}
  \centering
      \includegraphics[height=.2\textheight]{ORIFI1.jpg}$\quad$
    \includegraphics[height=.2\textheight]{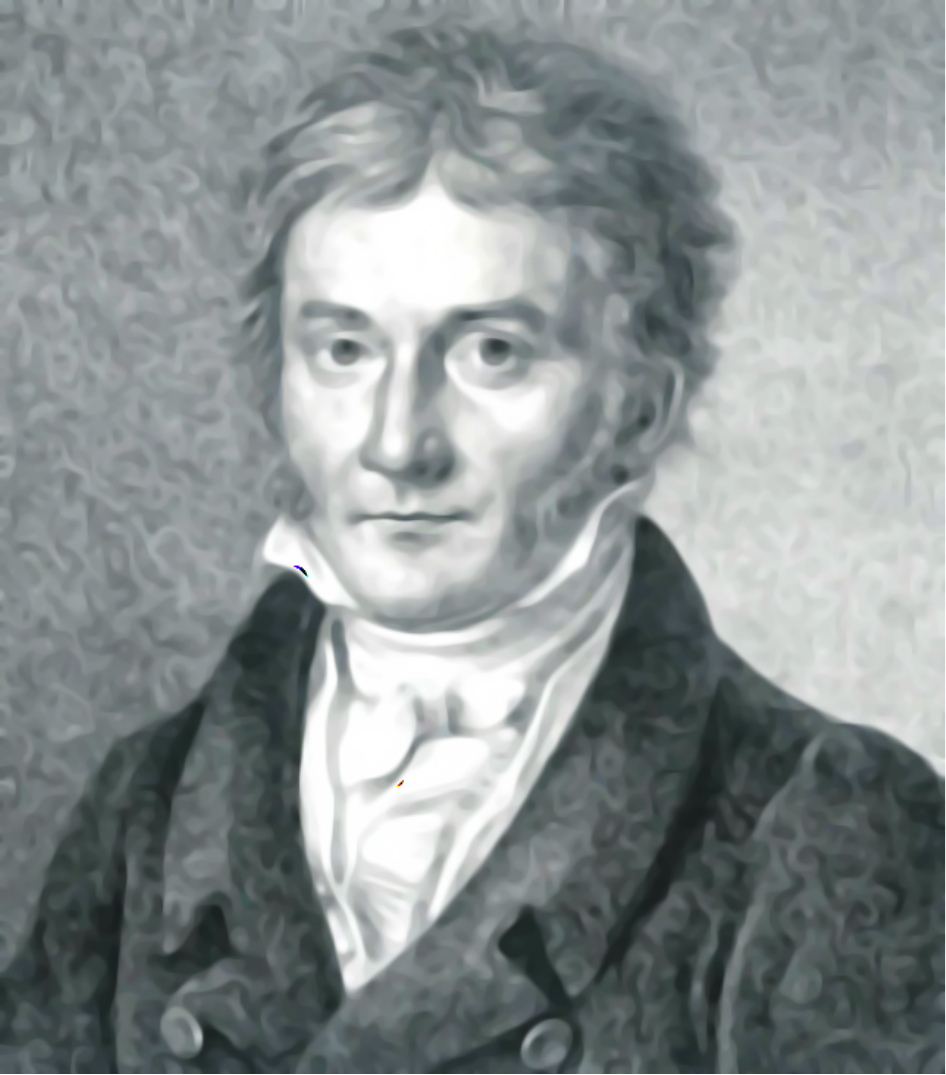}
 \caption{\sl Portrait of Gau{\ss} 
by artist Siegfried Detlev Bendixen published in Astronomische Nachrichten in 1828
(Public Domain images from
 Wikipedia) and a mean curvature flow filtered image of it (obtained via Mathematica).}\label{VILORi-1nsFI0232BIS}
\end{figure}

Hence, since
\begin{eqnarray*}
&&0=\partial_t \big(u(x_t(\eta),t)\big)=\nabla u(x_t(\eta),t)\cdot\partial_t x_t(\eta)+\partial_tu(x_t(\eta),t),
\end{eqnarray*}
we have that
\begin{eqnarray*}
\nu_t(0)\cdot\partial_t x_t(0)=-\frac{\nabla u(0,t)}{|\nabla u(0,t)|}\cdot\partial_t x_t(0)
=\frac{\partial_tu(0,t)}{|\nabla u(0,t)|}=\left.\div\left( \frac{\nabla u(x,t)}{|\nabla u(x,t)|}\right)\right|_{x=0}=-H_t(0).
\end{eqnarray*}
This gives~\eqref{ROLLSTO-MECA}, as desired.

Other geometric evolution equations will appear in Section~\ref{ROLSTP}.\medskip

See Figure~\ref{VILORi-1nsFI0232BIS} for an application of the (mean) curvature flow
filter. Once again, by comparing Figures~\ref{VILORi-1nsFI0231},
\ref{VILORi-1nsFI0232} and~\ref{VILORi-1nsFI0232BIS}, we have a practical confirmation
on how sharply the Perona-Malik method outperforms both classical diffusion and curvature flows
in maintaining definite meaningful edges while smoothing intermediate brightness regions.

\subsection{Artificial intelligence and machine learning}\label{ITMLDL}

We have already encountered in Sections~\ref{NESNISE} and~\ref{I2K243S3M:0oikrfg-4-2rtghj}
some fascinating problems related to discrete versions of the Laplace
operator. Here we provide a concrete situation in which these discrete Laplacians
naturally surface.

\begin{figure}
  \centering
  \includegraphics[height=6cm]{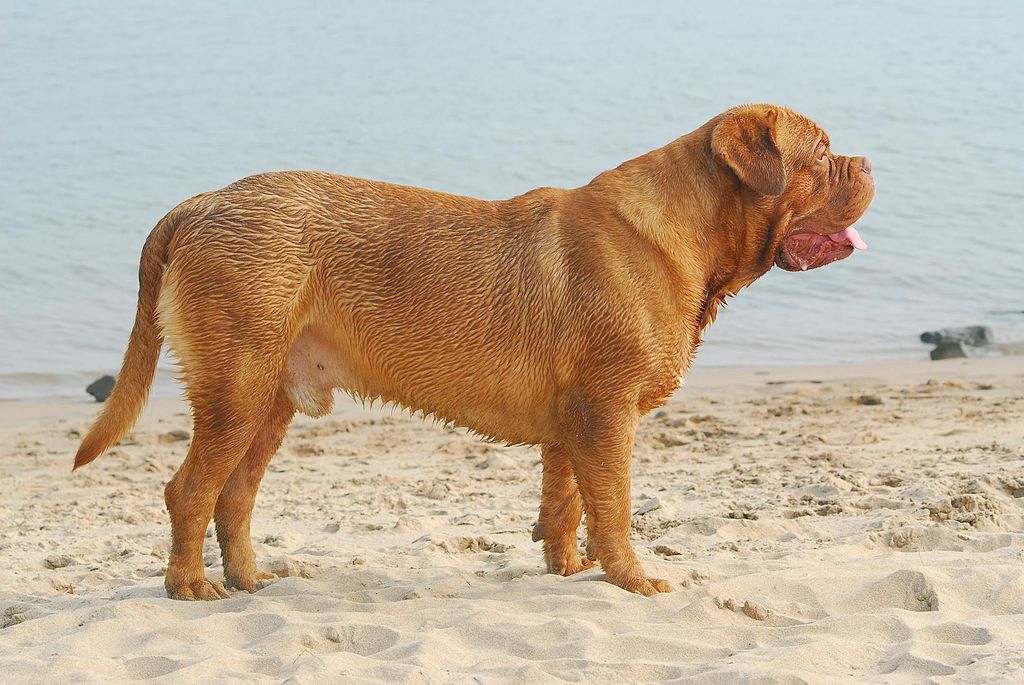}$\quad$
  \includegraphics[height=6cm]{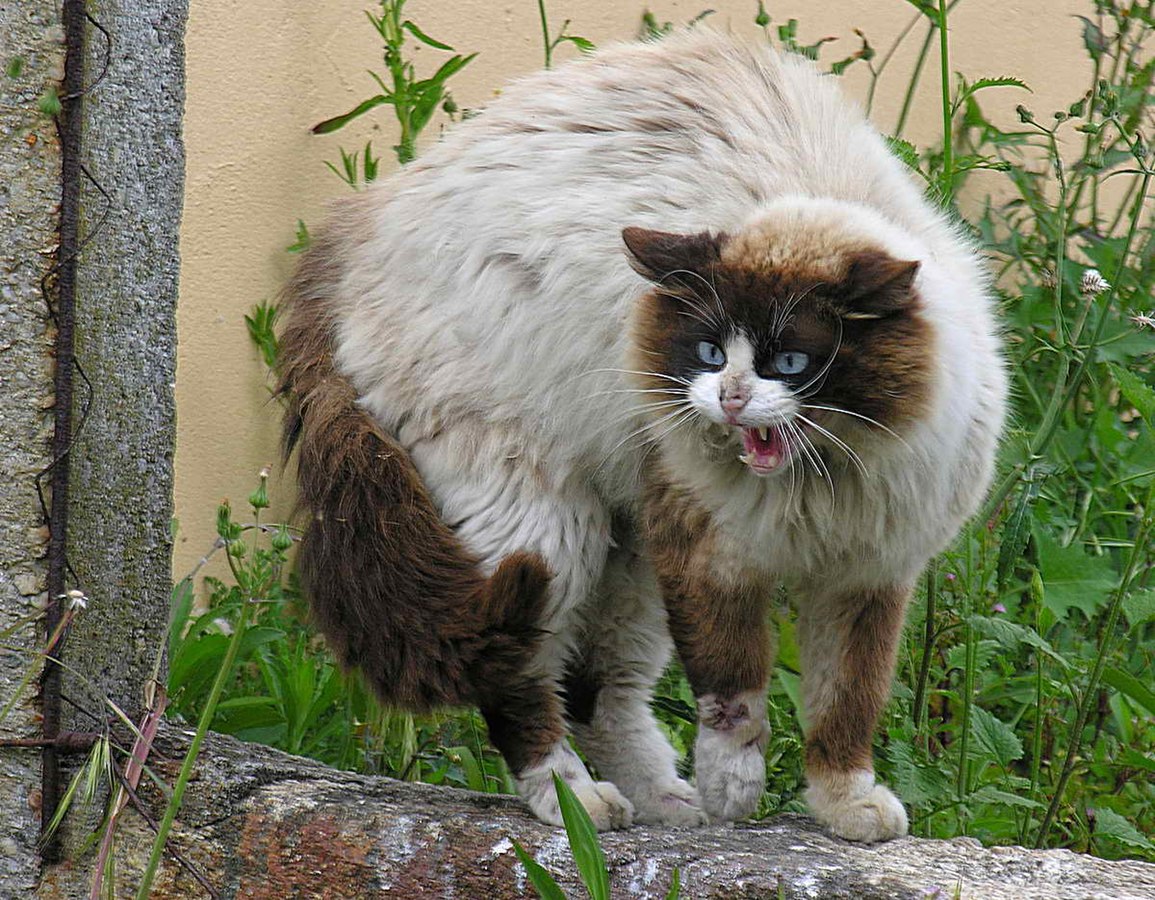}
 \caption{\sl Distinguishing a dog from a cat
 (images from
 Wikipedia; left image by DJT, licensed under the Creative Commons Attribution-Share Alike 2.5 Generic license,
 right image by Luis Miguel Bugallo S\'anchez, licensed under the Creative Commons Attribution 2.0 Generic
 license).}\label{CATTOGAN}
\end{figure}

The setting is that of image classification (not to be confused with
the image processing presented in Section~\ref{I2K243S3M:0oikrfg-4}):
the idea is to train a computer to identify specific objects in photos.
For this, the simplest model to keep in mind is to have a bunch of (say,
grayscale, for simplicity) images which represent some classes of targets to be identified
(again, for simplicity, we can think of two targets, ``cats'' and ``dogs'', see Figure~\ref{CATTOGAN}). The goal is thus to train
the machine to {\em distinguish cats from dogs} by looking at the photos.

Without aiming to be exhaustive, we mention some hints on how to implement such
a type of learning. First of all, let us transform the photos of cats and dogs into well-defined
mathematical concepts. Suppose that all the images have the same size, say each image consisting of~$d$ pixels.
Then, each image can be considered as an element of~$X:=[0,1]^d$,
in which, for every~$m\in\{1,\dots,d\}$, the $m$th coordinate of an element~$x\in X$
would correspond to the brightness level of the $m$th pixel (scaled to be a number between~$0$ and~$1$).
Also, each image has a label attached to it which specifies its content: we can think such a label
as an element of the set~$Y:=[-1,+1]$ (for example, the value~$-1$ corresponding to ``dog''
and~$+1$ to ``cat'', as an example of binary class).

The gist is now to show to the machine some data sample~$S$ of these labeled images,
say
\begin{equation}\label{01oeuf0iefhgwiqelfyiwt743tyghbrhbhu555ig335823asA}
(x_1,y_1),\dots(x_M,y_M)\in X\times Y,\end{equation} which can be used as a useful training
to learn how to independently label images.

For instance, the machine could have access to a certain class~${\mathcal{F}}$ of
functions~$f:X\to Y$ that could use to label new photos: in this way, when seeing a new image~$x\in X$,
the machine would produce a value~$f(x)$ which is the label ``guessed''
as appropriate: in our example, $f(x)=-1$ (or~$f(x)$ very close to~$-1$) would say that photo~$x$ is identified to be a dog,
$f(x)=+1$ (or~$f(x)$ very close to~$+1$) that photo~$x$ is identified to be a cat,
and values such as~$f(x)=0$ (or~$f(x)$ quite far from both~$-1$ and~$+1$)
would denote that the machine is uncertain about the attribution.

In order to make the machine select correctly the function needed for this identification process,
we can think to use the sample~$S$ to build some \index{loss function}
empirical loss function
which penalizes wrong attributions. Consider for example the case of a quadratic loss function
$$ L_f(x_1,\dots,x_M,y_1,\dots,y_M):=\sum_{j=1}^M \big(y_j-f(x_j)\big)^2,$$
which penalizes any possible choices of~$f\in{\mathcal{F}}$ unless it produces the right identification
on the sample~$S$ (which is~$f(x_j)=y_j$).

To let the previous example look more independent on the number of elements of the sample, one can also consider
$$ L_f(x_1,\dots,x_M,y_1,\dots,y_M):=\frac1M\sum_{j=1}^M \big(y_j-f(x_j)\big)^2.$$

Another example to keep in mind can be that of a ``cross entropy loss function'' of the type
$$ L_f(x_1,\dots,x_M,y_1,\dots,y_M):=\sum_{j=1}^M \left(
(1+y_j) \ln\frac{2}{1+f(x_j)}+(1-y_j)\ln\frac{2}{1-f(x_j)}
\right).$$
To develop an intuition of this loss function, one can notice that if, say, the $j$th image is a dog, then~$y_j=-1$
and the corresponding addend boils down to~$\ln\frac{2}{1-f(x_j)}$, which
is minimized indeed by the correct identification~$f(x_j)=-1$, producing a positive (possibly infinite) term whenever~$f(x_j)\in(-1,+1]$; differently from the previous cases, however,
this loss function would produce an infinite value when~$f(x_j)=+1$,
thus penalizing heavily the predictions that are confident but wrong.
\medskip

\begin{figure}
  \centering
  \includegraphics[height=6cm]{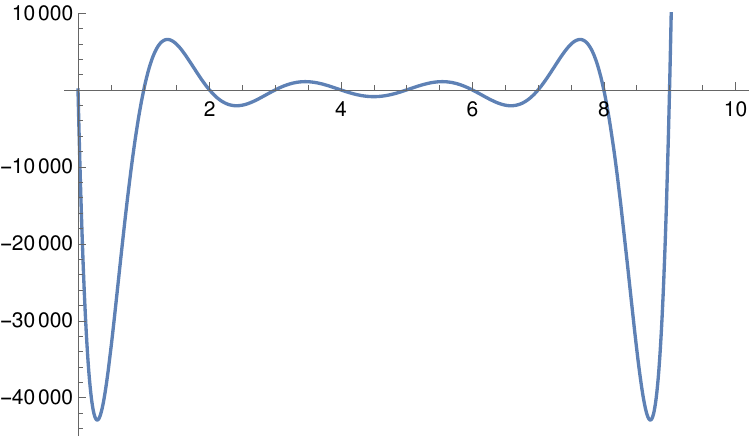}
 \caption{\sl Plot of the fuction~$t\mapsto t^{10} - 45 t^9 + 870 t^8 - 9450 t^7 + 63273 t^6 - 269325 t^5 + 723680 t^4 - 1172700 t^3 + 1026576 t^2 - 362880 t$.}\label{C22CAM-23OGAN}
\end{figure}

Now, in light of these discussions, one could argue that to efficiently train a machine it suffices
to find the minimizer~$f$ of some empirical loss function~$L_f$ (e.g., one of those above, or a similar one).
But there is a catch. Working in this way, our minimizer could just be any function~$f_\star$ such that
\begin{equation}\label{9ielnacw23r4t2f-xxj}
{\mbox{$f_\star(x_j)=y_j$
for all~$j\in\{1,\dots,M\}$.}}\end{equation} In this way, our machine has just ``memorized'' all the given pictures,
but it has discovered no pattern at all which can be used for future predictions.
Actually, it is quite likely that in front of a new image~$\widetilde x\not\in S$
the value~$f_\star(\widetilde x)$ may have little to do with the object represented in the image~$\widetilde x$: indeed, without putting any additional restrictions, there are simply ``too many functions''
satisfying a finite number of pointwise constraints, as in~\eqref{9ielnacw23r4t2f-xxj}.

This situation, in which the machine fits excellently against its training data without\footnote{As an evident example of overfitting, one can also consider a sample set given by the points 
in the plane
$$\big\{(0,0),\;(0,1),\;(0,2),\;(0,3),\;(0,4),\;(0,5),\;(0,6),\;(0,7),\;(0,8),\;(0,9)
\big\}.$$
Then, both the graph of the function~$f_0$ constantly equal to zero and that of the function
$$ \overline{f}(t):=t^{10} - 45 t^9 + 870 t^8 - 9450 t^7 + 63273 t^6 - 269325 t^5 + 723680 t^4 - 1172700 t^3 + 1026576 t^2 - 362880 t$$
would fit the data (see Figure~\ref{C22CAM-23OGAN}).

However, the function~$f_0$ constantly equal to zero would probably ``generalize'': indeed,
in our everyday experience, we may expect the next element to be taken into account to be the point~$(0,10)$, which actually would lie on the graph of~$f_0$. Instead, the prediction of the unnecessarily
and oscillatory function~$\overline{f}$ would produce~$\overline{f}(10)=3628800$, that would be quite off from
the correct value.

This is quite of an indication that simpler functions should be preferred to avoid overfitting,
also because oscillatory functions tend to produce values outside the interpolation points which are almost completely unrelated from the prescribed interpolation.

First and foremost, machine learning is an excellent playground to understand human learning too:
for machines, as well as for humans,
any sort of at least partially supervised learning may stem out of some training on a finite number of examples,
but, after that, to provide solid outcomes, it needs to
generalize to a broader class from which the examples are drawn, otherwise it is not real learning,
it is mere case-by-case memorization.

In a sense, for machines as well as for humans, the ultimate goal of learning is not to learn what to learn, but rather
to learn how to learn.}
being able to
generalize to new sets of data, is called ``overfitting''. \index{overfitting}
\medskip

\begin{figure}
  \centering
  \includegraphics[height=7cm]{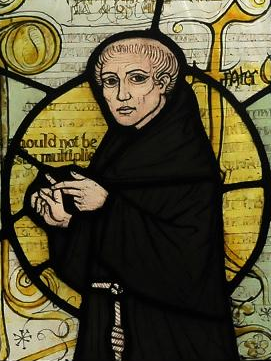}
 \caption{\sl A stained glass window at an English church allegedly representing
William of Ockham (image by Moscarlop from
 Wikipedia, licensed under the Creative Commons Attribution-Share Alike 3.0 Unported license).}\label{CCAM-23OGAN}
\end{figure}

To overcome this difficulty it is therefore convenient to appropriately select our functional space~${\mathcal{F}}$ (called in jargon ``hypothesis space'') \index{hypothesis space}
and to somehow
invite the machine to select, among the possible optimizing functions, the ``simplest'' one
(taking the ansatz that nature, after all, complies with Ockham's razor
and, among possible explanations, the simpler one is to be preferred, see Figure~\ref{CCAM-23OGAN}
to meet William of Ockham). Though dealing with the fine problems posed by the suitable selection
of hypothesis spaces goes well beyond the scopes of these notes (and relies on clever uses
of functional analysis), 
we can still try to get our hand in some aspects of this theory. 

\begin{figure}
  \centering
  \includegraphics[height=6cm]{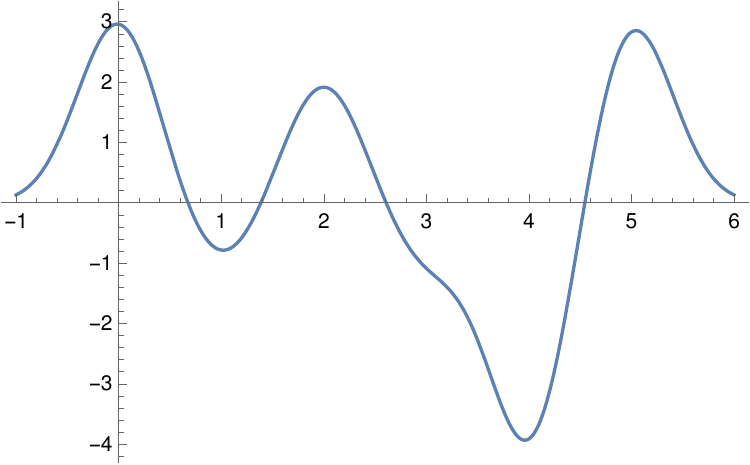}
 \caption{\sl Plot of the fuction~$x\mapsto 
 3e^{-\pi x^2} - e^{-\pi (x-1)^2} + 2 e^{-\pi (x-2)^2} - 
e^{-\pi (x-3)^2} - 4 e^{-\pi (x-4)^2} + 3 e^{-\pi (x-5)^2}$.}\label{C24522CAM123456789-23OGAN}
\end{figure}

To make things as simple as possible, though the method is very general,
we can just focus on a concrete case in which we try to approximate a given location by a Gau{\ss}ian,
e.g. a point~$(p,q)$ in the plane by the function~$t\mapsto q e^{-\pi(t-p)^2}$. In analogy with this idea,
one can define\footnote{See e.g.~\cite{MR2677883} for more information about Gau{\ss}ian kernels.
The choice of~$\pi$ in the exponent of the kernel is not relevant, we chose this normalization here
just because it would produce a total integral equal to~$1$.} a kernel, given~$\xi_1$, $\xi_2\in X$, as \label{IDJ:KETSPAE}
$$ K(\xi_1,\xi_2):= e^{-\pi(\xi_1-\xi_2)^2}
$$
and pick a large family of elements of~$X$, say~$\{\xi_1,\dots,\xi_D\}$ and define, for all~$i\in\{1,\dots,D\}$,
\begin{equation}\label{C24522CAM123456789-23OGAN:EQ2} \phi_i(x):= K(x,\xi_i).\end{equation}
Then, a natural space of functions capable of extracting information from this family
is that of the functions~$f$ that can be written as a linear combination of~$\phi_1,\dots,\phi_D$, say the functions of the form
\begin{equation}\label{C24522CAM123456789-23OGAN:EQ} X\ni x\longmapsto f(x;\theta):=\sum_{i=1}^D \theta_i\,\phi_i(x),\end{equation}
for some~$\theta=(\theta_1,\dots,\theta_D)\in\R^D$, see e.g. Figure~\ref{C24522CAM123456789-23OGAN}
for a simplified visual sketch of such ``linear expansion\footnote{After all, expressions as those in~\eqref{C24522CAM123456789-23OGAN:EQ} are not too different from a truncation of a Fourier series:
at hearth, an approximating Fourier sum is nothing but a linear superposition of nonlinear functions,
with the aspiration of reconstructing ``essentially'' every interesting function.

Instead of sines and cosines, maybe peaked functions as in~\eqref{C24522CAM123456789-23OGAN:EQ2}
have better chances for reconstructing functions out of pointwise constraints, as in~\eqref{9ielnacw23r4t2f-xxj}, hence can better serve to the scope of machine learning.

To build a serious theory out of the nonsense presented here, one has to leverage the analysis of
\index{reproducing kernel Hilbert spaces}
reproducing kernel Hilbert spaces, see e.g.~\cite{MR3526117}.

Also, here, we surf over the important problem of which type of functions can be reconstructed by the superposition of given functions (i.e., whether the sums of linear combinations of suitable translations are dense in interesting functional spaces). This type of questions often rely on functional analysis methods (such as the Stone-Weierstra{\ss} Theorem, the Hahn-Banach Theorem, the Riesz Representation Theorem, the Wiener Tauberian Theorem, the Kolmogorov-Arnol\cprime d Superposition Theorem, etc.). See e.g.~\cite{MR1015670, FUNAHASHI1989183, HORNIK1989359, HORNIK1991251, MR1178852, HORNIK19931069, LESHNO1993861, MR1819645, MR1883281, MR2677883, BA-2109-09710} for more information about density and approximation results (for instance, see Theorems~3.1 and~5.1 in~\cite{MR1819645} for specific approximation \label{OLKSD-APPOTJSD}
and interpolation results).}
of a nonlinear function''.

\begin{figure}
  \centering
  \includegraphics[width=.7\linewidth]{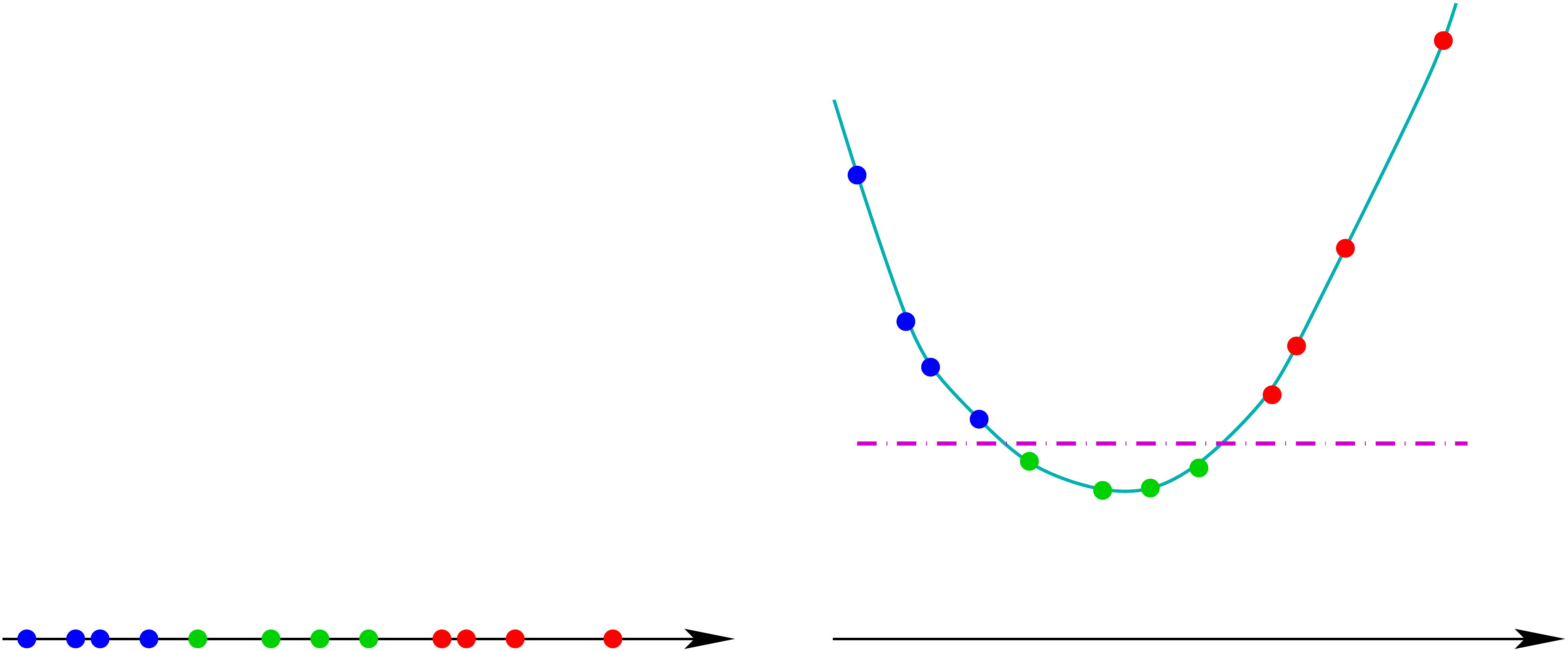}
 \caption{\sl A machine trying to learn how to pick good kiwis by embedding methods and separation by linear functions: the kiwis with the appropriate level of ripening are depicted here in green color, the ones in blue are still unripe and the ones in red are gone off.}\label{SIRKIW1fs0I}
\end{figure}

With this, a natural choice of hypothesis space~${\mathcal{F}}$ is that of functions~$f$
as in~\eqref{C24522CAM123456789-23OGAN:EQ} and the problem of minimizing
a loss functional becomes a finite dimensional\footnote{In general, the dimension~$D$ here is larger
than the dimension~$d+1$ of the space in which the sample~$S$ is taken (recall~\eqref{01oeuf0iefhgwiqelfyiwt743tyghbrhbhu555ig335823asA}), in order to allow the machine
to compare with the available data in the training sample and have sufficient parameters to optimize.

The convenience of embedding the data into a larger dimensional space is indeed very common in machine
learning and it is sometimes explained via a simple geometric analogy. For example,
let us suppose that we want to train our machine to recognize the best kiwis in our garden
and, for simplicity, suppose that the choice is made only in virtue of the ripening of the fruit. Of course,
kiwis that are either too unripe or gone off should be avoided, the ones that are to be selected have
a level of ripening lying in some appropriate interval. A structured way for a machine to select
such an interval therefore consists in embedding these ripening data into a higher dimensional
space and then separating the convenient ones by using a simple linear function, see Figure~\ref{SIRKIW1fs0I}.} minimization problem over the set
of parameters~$\theta\in\R^D$.\medskip

With this setting of hypothesis space~${\mathcal{F}}$ we can also come back to the overfitting problem,
trying now to avoid it by adding to the empirical loss function~$L_f$ a term which penalizes unnecessary
complicated oscillations. For example, one can consider a penalized loss function given by
\begin{equation}\label{0powjfeg:economical} L_f(x_1,\dots,x_M,y_1,\dots,y_M)+\lambda |\theta |^2,\end{equation}
for some~$\lambda>0$, where, e.g., $|\theta|$ is just the Euclidean norm in~$\R^D$: in this way, functions with ``more economical'' expansions would be preferred, and the larger the value of~$\lambda$ the larger the
penalty on the ``complexity'' of~$f$ (other choices of penalization are certainly possible
and methods of this sort are sometimes referred to as ``regularizations'').

Minimizing methods of empirical functionals as in~\eqref{0powjfeg:economical} are
indeed the core of the so-called {\em supervised machine learning}, in which
only labeled data sets are utilized for the training. In spite of the fact that this procedure
is often quite satisfactory, the method presents margin of improvements. As a matter of fact,
in many natural processes, learning occurs in a rather {\em semi-supervised} regime: \label{SSMLR65JS}
for example, a child is exposed to many new phenomena
and only few of them are directly connected to a specific label, yet a relatively small amount of feedback is sufficient to allow the child to make consistent learning progress.

This suggests that unlabeled data can also be usefully processed and used to
extract valuable information. To build a semi-supervised learning algorithm, one can modify the setting in~\eqref{C24522CAM123456789-23OGAN:EQ} by comprising in the minimization problem also
some unlabeled examples~$x_{M+1},\dots,x_{M+P}\in X$
(note that these examples are not labeled, namely, differently from the sample in~\eqref{01oeuf0iefhgwiqelfyiwt743tyghbrhbhu555ig335823asA}, they are not attached to any element of~$Y$).
To use these unlabeled examples, one must somewhat quantify how ``close
any element is to any other''.
To this end, given~$i$, $j\in\{1,\dots,M+P\}$, it is customary to endow the pair~$(x_i,x_j)\in X\times X$
with a weight~$w_{ij}\ge0$ representing their ``similarity'' (the case~$w_{ij}=0$ corresponding to
dissimilar objects and\footnote{Of course, extracting sufficient features to
assign automatically a precise degree of similarity to
unstructured data like images is a rather sophisticated task in itself, ultimately linked to pattern recognition,
see e.g.~\cite{MIMAGECALELF}.
For the purposes of this simplified explanation, we are just assuming that
our data set is endowed with these similarity weights
(regardless on whether these weights were assigned by humans
or by machines): what matters at this level is that this
unlabeled data set can now be used for the semi-supervised learning algorithm that
we are presenting.} higher values of~$w_{ij}$ indicating a higher degree of kinship).
The similarity relation is supposed to be symmetric, hence~$w_{ji}=w_{ij}$.

In this setting, one can impose to the empirical loss function an additional cost
if the assessment of~$x_i$ is far from that of~$x_j$ anytime the objects~$x_i$ and~$x_j$ are strongly related
(e.g., they have a high value of the similarity weight~$w_{ij}$): for example, one can consider an additional cost of the form
\begin{equation}\label{20346-PKJS-214ikt} \frac{\gamma}{(M+P)^2}\sum_{i,j=1}^{M+P} w_{ij}\big(f(x_i)-f(x_j)\big)^2,\end{equation}
for some~$\gamma>0$.

{F}rom this and~\eqref{0powjfeg:economical} we may conclude that an interesting loss function,
to be minimized within a suitable set of parameters that describes a convenient hypothesis space
(recall~\eqref{C24522CAM123456789-23OGAN:EQ}), in the framework of semi-supervised machine learning
takes the form
\begin{equation*} L_f(x_1,\dots,x_M,y_1,\dots,y_M)+\lambda |\theta |^2+\frac{\gamma}{(M+P)^2}\sum_{i,j=1}^{M+P} w_{ij}\big(f(x_i)-f(x_j)\big)^2.\end{equation*}

{F}rom the perspective of partial differential equations, the term in~\eqref{20346-PKJS-214ikt} is of particular interest,
since the minimization of loss functions containing this term produces discrete equations
involving the operator
$$ \sum_{j=1}^{M+P} w_{ij}\big(f(x_i)-f(x_j)\big).$$
In particular, when
$$w_{ij}=\begin{dcases} 1 & {\mbox{if }}|i-j|=1,\\
0&{\mbox{otherwise}},\end{dcases}$$ this operator reduces to the discrete Laplacian
studied in Sections~\ref{NESNISE} and~\ref{I2K243S3M:0oikrfg-4-2rtghj}
(compare with~\eqref{DISCRELEA}).\medskip

Once we understand the training process of a machine as the minimization over a set of parameters~$\theta$ of a suitable loss function, it is interesting to mention as well a fashionable direction of research related to the use of parameters and nonlinear functions which is somewhat modeled on the way natural organisms ``learn'', hence having interesting connections with neurology and cognitive sciences, and which paved the way\footnote{On the one hand, it is truly fascinating to compare machine and human learning, and this can certainly enrich cross-disciplinary research involving mathematics, engineering, computer sciences, biology and neuroscience. On the other hand, one has also to be aware that, as pointed out in~\cite[page~164]{MR3617773}, ``modern neural network research, however, is guided by many mathematical and engineering disciplines, and the goal of neural networks is not to perfectly model the brain. It is best to think of feedforward networks as function approximation machines that are designed to
achieve statistical generalization, occasionally drawing some insights from what we know about the brain, rather than as models of brain function''.}
to the theory of \index{deep learning} ``deep learning'', where
the adjective ``deep'' refers to a performance improvement due to the presence of additional
layers in the artificial neural architecture.

Namely, rather than restricting our attention only to the kernel methods described on page~\pageref{IDJ:KETSPAE},
a convenient way to obtain a simple, but possibly effective, hypothesis space is to produce a ``rich family'' of functions by composing repeatedly an affine function and a nonlinear one. After all, this is the way in which natural organisms develop their perceptions: for instance, if a neuron receives a sufficiently intense stimulus (corresponding to a significant voltage change over a short interval), it generates an electrochemical pulse which can reach another neuron; these traveling signals may be excitatory or inhibitory, increasing or reducing the voltage. That is, very roughly speaking, when a neuron receives an information can either suppress it or transmit it (more or less in an ``all-or-nothing'' fashion, according, say, to the intensity of the stimulus) and, in the reception and transmission process to another neuron the signal can be ``adjusted'' (that is increased, decreased, or modulated).

Thanks to mathematics, one can reproduce this method for a silicon based intelligence. In this situation, a ``neuron'' is replaced by a function, sometimes called \index{nonlinear activation function}
nonlinear activation function (after all, a neuron is just a function which gives different values
to model different reactions, or no reaction at all, to a given stimulus).

The literature related to machine learning presents plenty of nonlinear activation functions related to neural networks. In some sense, the ones closest to the intuition coming from biological neurons are activators which ``switch from~$0$ to~$1$'', such as functions~$\sigma:\R\to\R$ such that~$\sigma(t)\to0$ as~$t\to-\infty$ and~$\sigma(t)\to1$ as~$t\to+\infty$. This type of functions are sometimes called \index{sigmoids} ``sigmoids''. \label{SIGRE}

However, for computational purposes it is often convenient to pick simpler activation functions: with this respect, a very popular choice is to select as activation function the ``positive part'' of a real number, namely the function
$$ \R\ni t\longmapsto t_+:=\max\{t,0\}.$$
This function is sometimes dubbed with fancy names, such as ``rectified linear unit'', or \index{ReLU} ``ReLU'' for short.

Roughly speaking, a sigmoid mimics a neural response of the type ``all-or-nothing'', up to a possible small transition between the regime in which ``nothing happens'' to that in which ``something happens'', while the ReLU reproduces a ``ramp type'' reaction in which ``nothing happens below a threshold, after which the reaction is linear''.

The technical advantage of considering a ReLU, besides its simplicity, is that it is a convex and piecewise linear function. Additionally, a sigmoid can be easily constructed as a linear superposition of two ReLUs, since
$$ t_+ - (t-1)_+ =\begin{dcases} 1 & {\mbox{ if }} t\ge1,\\
t & {\mbox{ if }}t\in(0,1),\\ 0 &{\mbox{ if }}t\le0, 
\end{dcases}$$
hence using the ReLUs may end up being technically simpler, but roughly ``as good as'', using the (perhaps conceptually more intuitive) sigmoids, thus obtaining a performance boost\footnote{For higher
accuracy, the ReLU is also sometimes replaced by the so-called ``leaky ReLU''
$$ \begin{dcases} t & {\mbox{ if }} t\ge0,\\
-ct & {\mbox{ if }}t<0, 
\end{dcases}$$
for~$c\in(0,+\infty)$.} of the method.

The specific choice of nonlinear activation functions is somehow related to the choice
of loss function and these two design elements are connected since
the configuration of the neural layer frames the prediction problem of the machine, while the choice of the loss function calculates the predictive error for a given frame of the problem.

\begin{figure}
                \centering
                \includegraphics[width=0.5\linewidth]{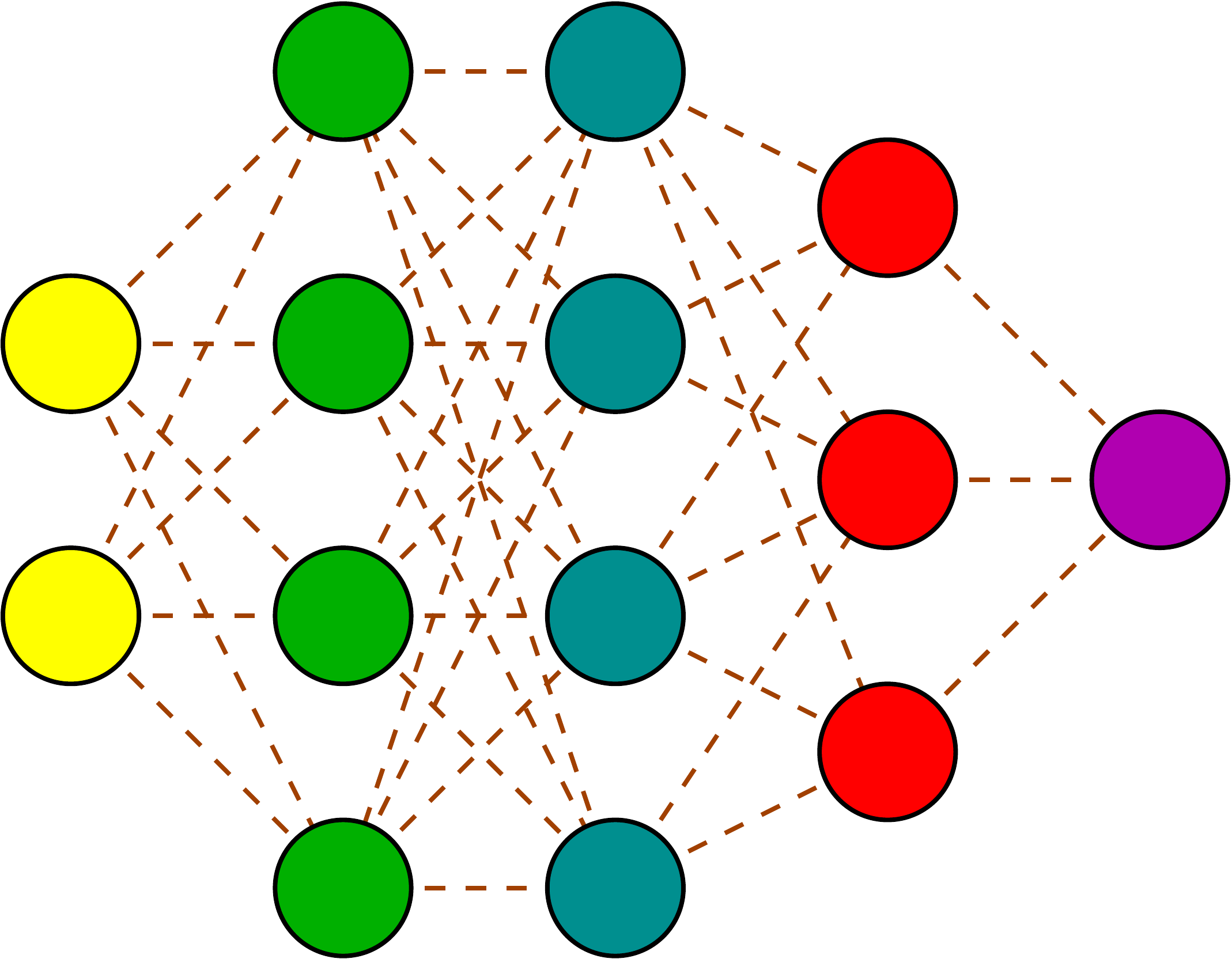}
        \caption{\sl Scheme of a neural network.}\label{23HAFCh-12X231}
\end{figure}

One of the strengths of deep learning consists of building
multiple processing layers, by placing a number of ``hidden'' computational levels
between the initial input and the final predictive output.
In this way, maybe,
an artificial network can ``imitate'' the way humans acquire sophisticated types of knowledge,
since these hidden layers in the artificial neural structure may somewhat reproduce
multiple levels\footnote{Roughly speaking, one can imagine that the hidden layers in a neural network correspond to different level of perceptions: for instance, if, given a collection of pictures, the task of the machine is to distinguish a kangaroo from a shark, one can imagine that the first hidden layer focuses on the identification of the edges of the picture, the second on the detection of the corners (recognizable as collection of edges), the third on the perception of object parts (detectable as collection of corners and contours), the fourth on the mutual interaction between object parts, and so on. But of course, this approach to describe the layers of a neural network may well be too simplistic, and too much ``human oriented''. See however the forthcoming description
of convolutional neural networks on page~\pageref{GENFIKDt576KJSllIKDSqasay7T1yHL}.}
of ``abstraction'' in which the system represents and analyzes the data.

Roughly speaking, the construction of the additional layers added to enhance the possibility for a 
machine to learn from the data is sketched in
Figure~\ref{23HAFCh-12X231}, depicting a simple\footnote{We strongly recommend to try to understand computer science, and in particular deep learning, first with your computer switched off, using only pencil and paper instead. In this spirit, for instance, a simple learning algorithm is presented in details in Section~6.1 of~\cite{MR3617773}.}
neural network, with an input, consisting of a vector~$(x_1,x_2)$, represented by the yellow dots, three ``hidden'' layers, portrayed by the green, cyan and red dots, and an output, given by a scalar and illustrated by the magenta dot.

The brown dotted lines in Figure~\ref{23HAFCh-12X231} represent the actions of matrices of weights, and the composition with a nonlinear activator function. For instance, suppose that the brown dotted lines entering the top green dot correspond to the weights~$\frac12$, coming from the upper yellow dot, and~$\sqrt\pi$, coming from the bottom yellow dot, then the value assigned to the top green dot corresponds to~$g\left(\frac{i_1}{2}+\sqrt\pi i_2\right)$.
In general, if~$k\in\{1,\dots,4\}$, the~$k$th green dot is assigned\footnote{More generally, one
can consider affine transformations like~$\theta_{1,1,k}x_1+\theta_{1,2,k}x_2+\theta_{1,3,k}$ instead of the linear one~$\theta_{1,1,k}x_1+\theta_{1,2,k}x_2$.
But this boils down mostly to a notational convention, since it would suffice to replace~$(x_1,x_2)$ with~$(x_1,x_2,x_3):=(x_1,x_2,1)$
to reduce the affine setting to the linear one.}
to some value~$h_{1,k}:=g(\theta_{1,1,k}x_1+\theta_{1,2,k}x_2)$, for some weights~$\theta_{1,1,k}$ and~$\theta_{1,2,k}$.

And so on, looking at the subsequent hidden layer of cyan dots, the~$k$th cyan dot is assigned to some value
$$h_{2,k}:=g(\theta_{2,1,k}h_{1,1}+\theta_{2,2,k}h_{1,2}+\theta_{2,3,k}h_{1,3}+\theta_{2,4,k}h_{1,4}),$$
for some weights~$\theta_{2,1,k}$, $\theta_{2,2,k}$, $\theta_{2,3,k}$ and~$\theta_{2,4,k}$, or, in a short vectorial notation,
$$ h_{2,k}(x)=g(\theta_{2,k}\cdot h_1(x))=g\left( \theta_{2,k}\cdot g(\theta_{1}\cdot x) \right),$$
with the slight notational abuse in which~$g$ denotes both a scalar and a vectorial function.
\medskip

This network structure\footnote{Coming back to the approximation problems described in footnote~\ref{OLKSD-APPOTJSD} on page~\pageref{OLKSD-APPOTJSD},
we mention that the possibility of using additional
hidden layers in a neural network enhances the capability of approximating a given function. That is, two
basic techniques can be adopted to approximate
a given function by the outputs of a neural network.
The first is to build a network with many neurons but
without hidden layers. In this case,
the approximation relies on the linear superposition
of a large number of nonlinear activation functions, and this situation,
often denoted the ``arbitrary-width case'', corresponds
to having only one nonlinear step between the input and
the output of the form
$$\sum_{i=1}^N C_{i}\sigma\left(\sum_{{1\le j\le n}}\theta_{ij} x_j+b_i\right),$$
see e.g.~\cite{MR1015670, MR1178852, MR1819645, FUNAHASHI1989183, HORNIK1989359, HORNIK1991251, HORNIK19931069}. In this case, the ``width'' is given, for instance, by the number~$N$ of superpositions
of neurons considered here above.

Note that the introduction of hidden layers amplifies in principle the possibility of approximating a given function, without spoiling the possibility of using a network with no hidden layers, because
one can simply approximate the identity function with additional later layers, for which it suffices
to pick a point~$t_0$ for which~$\sigma'(t_0)\ne0$ and note that
$$ \frac1\e\left[\sigma\left( \frac{\e t}{\sigma'(t_0)}+t_0\right)-\sigma(t_0)\right]=t+O(\e).$$

In this spirit, the other basic approximation method, 
called ``arbitrary-depth case'', consists in approximating a given function by using a neural network with
an arbitrary number of hidden layers, each with at most some prescribed number of neurons,
see e.g.~\cite{NIPS201732cbf687, math7100992, pmlr-v125-kidger20a, ZHOU2020787}.
In this situation, the ``depth'' is given by the number of hidden layers. 

All these ``universal approximation'' results put perhaps some pressure 
on the scientists utilizing machine learning:
as stated in~\cite[Section~3]{HORNIK1989359},
``standard multilayer feedforward networks are capable of approximating any measurable function to any desired degree of accuracy, in a very specific and satisfying
sense [...] This implies
that any lack of success in applications must arise
from inadequate learning, insufficient numbers of
hidden units or the lack of a deterministic relationship between input and target''.}
is sometimes referred to as ``feedforward'', since the information flows\footnote{Enhanced algorithms
allowing the information to go back from the loss functional to the network itself are instead called ``backpropagation''.}
from the input through the intermediate layers and finally to the output: interestingly, the training procedure typically occurs only at the level of the output, which needs to be as close as possible to the real data, as prescribed by the loss function (in this sense, the learning algorithm uses the hidden layers to produce an optimal output, but the training data do not specify directly what individual layer should do).\medskip
\begin{figure}
  \centering
  \includegraphics[height=7.1cm]{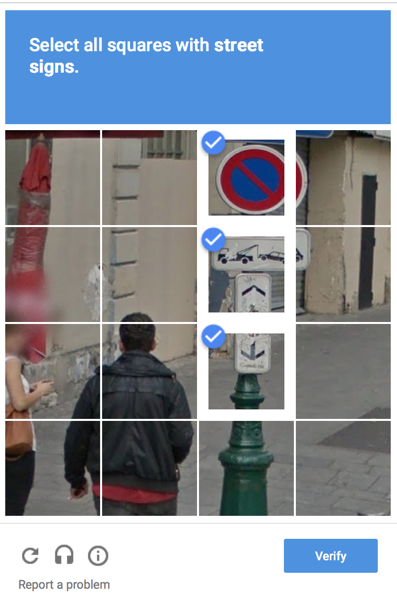}$\quad$
   \includegraphics[height=7.1cm]{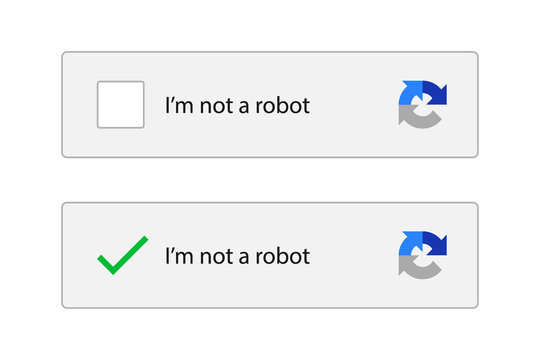}
 \caption{\sl CAPTCHA (Completely Automated Public Turing test to tell Computers and Humans Apart):
 I'm not robot (am I?).}\label{SIRCGCAP0AU0oedfs0I}
\end{figure}

We also remark that, while the machine learning process typically relies on optimized sequences
of linear and nonlinear operations, and the precise choice of nonlinearity may play
a role in the efficiency of the algorithm (recall the discussion about sigmoids and ReLUs
on page~\pageref{SIGRE}), also the choice of the linear operation can contribute to
the achievement of performant learning outcomes.
Without going into details, for example, we recall that the standard linear superposition
gets sometimes replaced, or complemented, by a convolution operation (see~\cite[Chapter~9]{MR3617773}). Roughly speaking, the idea for such convolutional neural networks
can be understood by considering for simplicity the pixels of a black and white image
as our input (corresponding, according to the color of the pixel, to either the value~$0$ or~$1$).
A feature detector is applied to this input by applying a discrete convolution: say that the pixels
of the input image are labeled as~$x_{ij}\in\{0,1\}$ with~$i\in\{1,\dots,N_1\}$ and~$j\in\{1,\dots,N_2\}$,
the action of a ``feature detector'' produces a value~$\widetilde x_{\ell m}\in\R$
with~$\ell\in\{1,\dots,\widetilde N_1\}$ and~$m
\in\{1,\dots,\widetilde N_2\}$ via the relation
$$ \widetilde x_{\ell m}=\sum_{{1\le i\le N_1}\atop{1\le j\le N_2}} K_\theta(\ell-i,m-j) \,x_{ij}.$$
Namely, kernels~$K_\theta$ are used to perform linear
convolutional operations, the result of which
is called ``feature map'' (note that in our example
the input image is~$(N_1\times N_2)$-dimensional
while the feature map is~$(\widetilde N_1\times\widetilde N_2)$-dimensional):
therefore, a convolution operation may be convenient not only
to relate the information encoded into a pixel with that of the neighboring pixels,
but also to reduce the size of the input image whenever convenient, thus
accelerating the algorithm (yet, possibly at the cost of
losing some information, but sometimes to recognize a pattern
one uses some broadly seen features, without needing to dig down to their minutiae).
The linear operation of convolution is then combined with some nonlinear activation
function, possibly within an architecture of hidden neural layers
dealing with collections of feature maps, and, as above,
the parameters~$\theta$ describing the kernels are optimized by the minimization
of a suitable loss functional,
with the goal of training the machine against a given sample.

In terms of actions performed on the image, different convolution kernels
can correspond to \label{GENFIKDt576KJSllIKDSqasay7T1yHL} blurring, 
sharpening contours, detection of edges, thus, in some sense, the training process of
the kernel parameters
allows the machine to detect and interpret features such as lines,
corners and contours, by selecting the optimal
feature detector for the appropriate task at the convenient layer of
the network, with the nonlinear activation function
maximizing the effect of the feature determinant for the required task. In this respect, convolutional networks are
intimately related to their neuroscientific counterpart\footnote{And nowadays
convolutional networks can outperform humans in several visual recognition tasks.
For instance, machine
learning outperforms humans by a factor of two on a traffic sign recognition benchmark, see~\cite[Section~3.5]{CIRE6248110}.

So, when on internet some websites waste our time asking us to prove that we are not a robot
by clicking on the images of traffic signs, we'd better ask the help of a robot, see Figure~\ref{SIRCGCAP0AU0oedfs0I}.}
in the primary visual cortex of the brain (see~\cite[Section~9.10]{MR3617773}
for further details on the relation between neuroscience and convolutional networks). 
\medskip

Let us mention that in the numerical simulations of partial differential equations, neural networks have been widely exploited in view of their efficiency in terms of computational time. Typically, in these simulations, neural networks do not replace the classical methods of numerical analysis (such as finite elements or finite difference algorithms): indeed, not only these classical methods still produce more accurate results than those obtained through machine learning, but also they are used precisely to produce the numerical solutions used to train the machine via the neural network. In this, the machine learning approach to simulations of partial differential equations sometimes aims not at producing more precise results, but rather quick answers, even at the risk of being slightly suboptimal.

Very interestingly, neural networks can also be used to suggest new perspectives expanding our mathematical knowledge of partial differential equations and helping mathematicians draw a road map towards the proof of very difficult theorems, see e.g. {\tt
https://www.quantamagazine.org/deep-learning-poised-to-blow-up-famed-fluid-equations-20220412/} and the references therein.\medskip

For additional and more comprehensive information\footnote{We do not address in these pages the latest developments of artificial intelligence and machine learning. It suffices to mention that some typical issues of the classical machine learning methods are related to vanishing gradient problems (affecting the gradient descents used e.g. for minimizing the loss function), slow learning algorithms (due to the many weights to be optimized) and high learning costs (due to the power consumption employed by a machine while training large datasets). Also, an additional challenge for classical networks is to deal with data which vary with time.

A fashionable answer to these problems is given by \index{reservoir computing} reservoir computing. Roughly speaking, a reservoir is a high dimensional nonlinear dynamical system, playing the role of a recurrent network of interconnected neurons. The reservoir typically receives an input, e.g. through a weighted linear combination of time-series data, which are processed by the nonlinear recursive features of the reservoir, thus producing an output, e.g. again via a  weighted linear combination. 

In this architecture, typically only the weights related to the linear combinations entering and exiting the reservoir are trained against the data sample, while the reservoir itself is not subject to a specific training (though some of the parameters of the reservoir typically require a bespoke choice for high performances).

The learning algorithms can be speeded up in this way, since the size of the training set is often downsized, and the necessary computing resources can be reduced, since the number of adjustable weights to be optimized is lower than in the classical neural architectures.

Interestingly, a physical system can be used as a reservoir, thus replacing a step usually performed by silicon chips by possibly cheaper physical systems and devices (also, a link with quantum computing may arise, since the nonlinear system utilized can be borrowed from either classical or quantum mechanics).

{F}rom the perspective of cognitive neurosciences, one can also observe (differences and) similarities between the use of reservoirs in artificial intelligence and some traits of the prefrontal cortex in the human brain, see e.g. Figure~2 in~\cite{HANDSON}. See also~\cite{GAUTHIER} and the references therein for further information about reservoir computing.}
about machine learning, seen from various perspectives, see for instance
Adam Oberman's lec\-ture {\tt https://www.youtube.com/watch?v=s-rMfJsf35Y},
Dejan Slepcev's lec\-ture
{\tt https://www.youtube.com/watch?v=bLlgdzg094g},
the text\-books~\cite{MR3617773, MR3931734}, the articles~\cite{MR2274444, MR4070203}, the set of notes by Leonardo Ferreira Guilhoto
{\tt https://math.uchicago.edu/$\sim$may/REU2018/REUPapers/Guilhoto.pdf},
and the references therein.

\subsection{Cutting networks}

Other interesting applications of the discrete Laplacians presented in Sections~\ref{NESNISE}, \ref{I2K243S3M:0oikrfg-4-2rtghj} and~\ref{ITMLDL} surface when studying networks.
Roughly speaking, a network (also known as graph) is a (finite) collection of objects in some mutual relation. To better visualize this structure, these objects are depicted as vertices (also called nodes) and vertices that are related to each other are connected by an edge (also called link). Formally, a network thus consists of a pair~$(V,E)$, where~$V$ is the collection of vertices and~$E$ is the collection of edges.

We note that an edge in the network is identified by its endpoints, which are vertices of the network. Namely, given~$v$, $w\in V$, an edge connecting~$v$ and~$w$ is denoted by~$(v,w)$; in our setting, the link between two vertices has no preferred direction, therefore the edge~$(v,w)$ coincides with~$(w,v)$.

If the vertices~$v$ and~$w$ are connected\footnote{For simplicity, here we do not allow edges from a vertex to itself (loops), therefore~$v\sim v$ never holds.} by an edge, i.e., if~$(v,w)\in E$,
we will write~$v\sim w$.

We will tacitly assume that the network under consideration is connected, i.e. for each pair of vertices~$v$ and~$w$, there is a sequence of edges~$(v_0,v_1),(v_1,v_2),\dots,(v_{k-1},v_k)\in E$ such that~$v_0=v$, $v_k=w$, $v_i\ne v_j$ unless~$i=j\in\{0,\dots,k\}$, and~$(v_i,v_{i+1})\ne(v_j,v_{j+1})$ unless~$i=j\in\{0,\dots,k-1\}$ (that is, we can connect any two vertices in our connected network).

An edge~$(v,w)\in E$ is said to be incident to a vertex~$z\in V$ if either~$ v = z$ or~$w = z$.

\begin{figure}
  \centering
  \includegraphics[width=.65\linewidth]{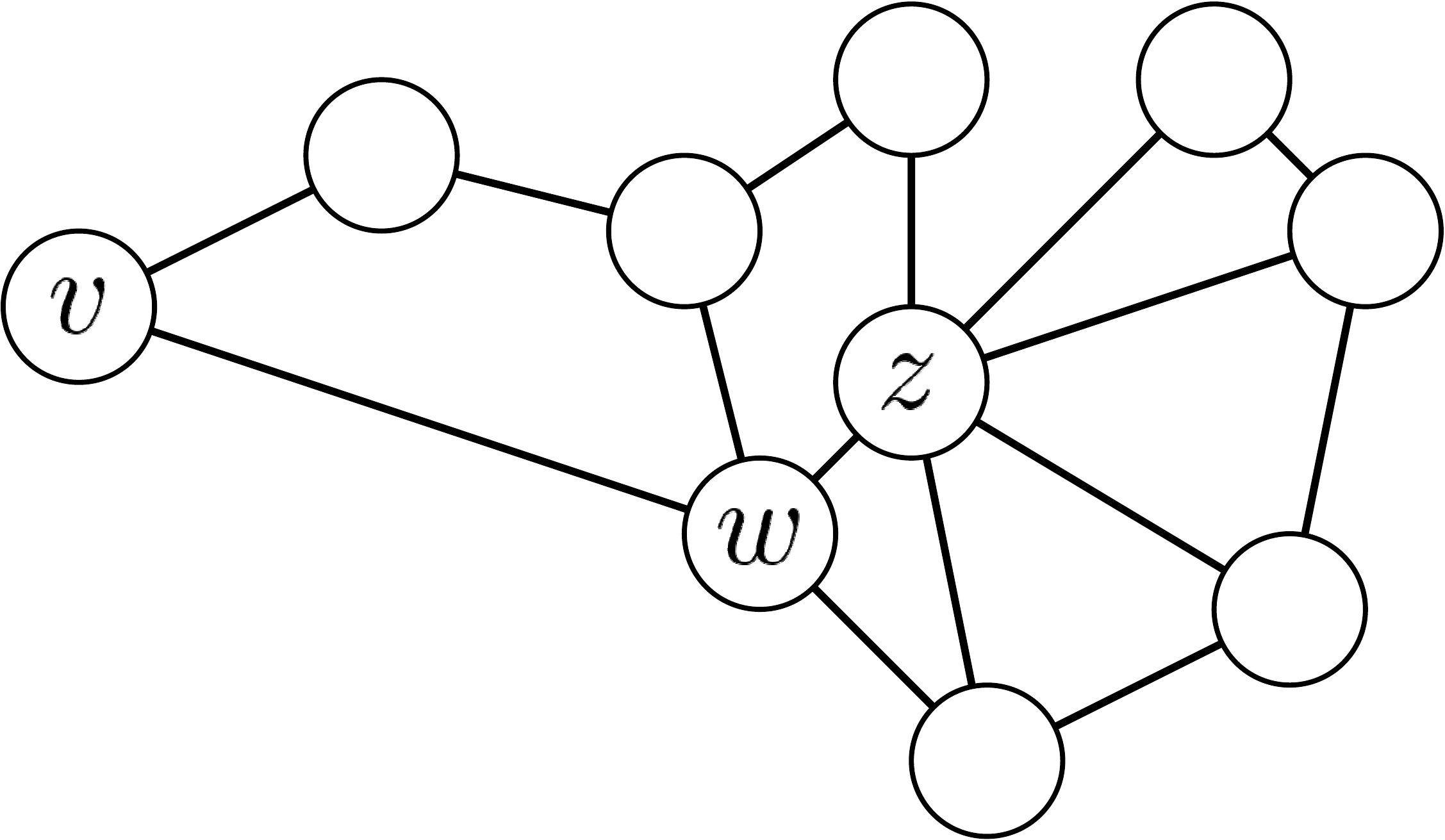}
 \caption{\sl A network.}\label{1P2NeT2e3FIlAALANSE-2}
\end{figure}

The degree of a vertex of a network is the number of edges that are incident to the vertex (for instance, in Figure~\ref{1P2NeT2e3FIlAALANSE-2}, the vertices~$v$ and~$w$ are connected by an edge,
and so are the vertices~$w$ and~$z$,
the vertex~$v$ has degree~$2$,
the vertex~$w$ has degree~$4$, and the vertex~$z$ has degree~$6$). If~$v\in V$ is a vertex, the degree of~$v$ is often denoted by~$\deg(v)$.

The volume of a collection~$S\subseteq V$ of vertices is the sum of the degrees of vertices in~$S$, namely one sets
$$ |S|:=\sum_{v\in S} \deg(v).$$
In a sense, the volume of~$S$ describes ``how big'' the collection~$S$ of vertices is, by taking into account not simply the number of vertices in~$S$ but rather weighting this number by the ``importance'' of these vertices, as encoded by their degrees.

A cut of a network is a partition of its vertices into two subsets, say~$S$ and~$T:=V\setminus S$. The boundary (also called cut-set) of such a cut is the set
$$ \partial S:=\big\{(v,w) \in E {\mbox{ s.t. }} v\in S {\mbox{ and }} w\in T\big\} ,$$
that is the collection of edges that have one endpoint in~$S$ and the other endpoint in~$T$. We observe that, by construction, $\partial S=\partial T$.

The size of such a cut is the number of edges crossing the cut, that is the cardinality of~$\partial S$,
which will be denoted by~$|\partial S|$.
For instance, if we take~$S:=\{v,w,z\}$ in Figure~\ref{1P2NeT2e3FIlAALANSE-2},
we have that~$|S|=12$ and~$|\partial S|=8$.\medskip

In a social network, \index{social network}
every vertex represents a person and an edge represents some link between people
(e.g., the edge~$(v,w)$ can represent the fact that two people know each other).
In this sense, one may be interested in using the notion of cut for a social network
to split the community~$V$ into two subgroups~$S$ and~$T$ such that
people in the same subgroup are ``likely to know each other'' and people in different
subgroups are ``less likely to know each other''.

Another example of network consists in a collection of handwritten digits:
in this situation, for instance, we can consider vertices corresponding to
handwritten ``threes'' or ``eights'' and two vertices are linked by an edge if they ``look sufficiently similar to each other''.
In this situation, we may want to implement an efficient algorithm to
tell handwritten threes and eights apart, i.e. to cut the network consisting of these digits
into two subsets, one collecting the digits which are ``likely to be three''
and the other collecting the digits which\footnote{To make things more accurate,
one could even consider networks with weights, assigning a higher importance
to a link between people who know each other well,
or to handwritten digits which are particularly similar to each other.
For simplicity, we do not consider this case here.
A procedure like this was already mentioned in the context of semi-supervised
machine learning, see page~\pageref{SSMLR65JS}.}
are ``likely to be eight''.\medskip

\begin{figure}
  \centering
  \includegraphics[width=.75\linewidth]{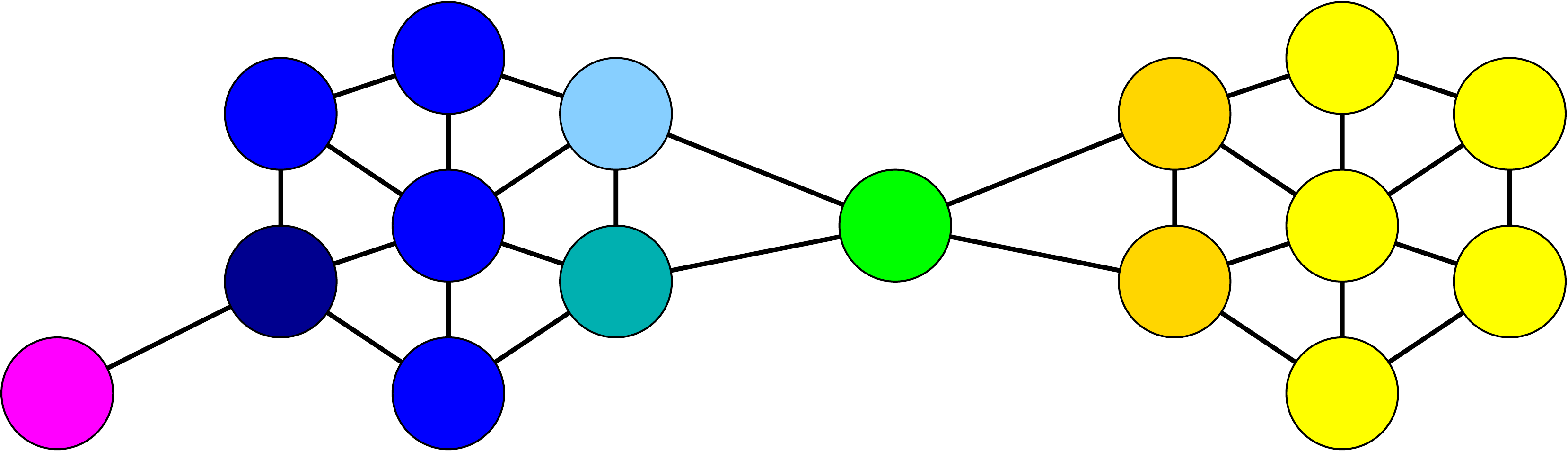}
 \caption{\sl A colorful network.}\label{COLOT2e3FIlAALANSE-2}
\end{figure}

Of course, in choosing a cut in these examples, one would like to reduce the size of the cut as much as possible, but a greedy choice in this can be highly counterproductive, in view of possible outliers.
To get into the swing of this phenomenon, one can look at the case depicted in Figure~\ref{COLOT2e3FIlAALANSE-2}.
In this situation, one would like to distinguish ``blue'' and ``yellow'' colors, and vertices are linked
if the colors are sufficiently similar. Notice the presence of a ``green'' vertex (which looks
sufficiently similar to both a couple of pale blue vertices and a couple of intense yellow vertices)
and of a ``magenta'' vertex (which looks
sufficiently similar a dark blue vertex but different from all the yellow vertices).

In this case, the ``greedy'' choice of just minimizing the size of the cut would
lead to consider the magenta vertex as separate from all the others, since this cut has size~$1$,
but this operation would fail badly in telling blues and yellows apart. It would be instead more appropriate
to select a cut which groups together, for instance, all the blues (maintaining the magenta vertex
in this group as a negligible error) and separately all the yellows
(keeping perhaps the green vertex
in this group as a negligible error): this cut has indeed a slightly larger size, equal to~$2$, but
has the advantage of splitting the network into two groups of equal numerousness.

To implement this idea, it is convenient to introduce the Cheeger ratio \index{Cheeger ratio} of a cut, given by
\begin{equation}\label{EHOSS} h(S):=\frac{|\partial S|}{\min\big\{|S|,|T|\big\}},\end{equation}
and we observe that~$h(S)=h(T)$, since~$T=V\setminus S$.

The Cheeger constant \index{Cheeger constant}
of the network (also known as the conductance of the network) is then defined by
\begin{equation}\label{EHOSS0} h:=\min h(S),\end{equation}
where the minimum above is taken among all possible choices of subsets~$S$ of~$V$.

We observe that finding a cut that attains the Cheeger constant may actually be much more successful
in the task of separating blues and yellows, especially when one knows to start with that
the network collects two classes of individuals of approximately equal sizes,
since the minimization method in~\eqref{EHOSS0} finds more ``balanced'' separations: for instance, in Figure~\ref{COLOT2e3FIlAALANSE-2},
the cut corresponding to the isolation of the magenta vertex produces a Cheeger ratio of~$1$,
while the cut dividing the blues (magenta included) and the yellows (green included)
provides a Cheeger ratio strictly smaller than~$1$.

For practical purposes, the task of finding a cut attaining the Cheeger constant
may be very hard to accomplish, since, in principle, for a network consisting of~$n$ vertices
one has~$2^n$ possible cuts to take into account (or about~$2^{n-1}$ if one exploits the symmetry between
the cut identified by~$S$ and the one identified by~$T=V\setminus S$), making the procedure unfeasible
for large (and interesting) networks,
such as the ones of prohibitively large size arising
in social sciences, biology or internet data. It could be therefore convenient
to reduce the problem to a smaller and more manageable one, possibly producing an explicit algorithm,
in such a way that this new problem maintains similar features to that proposed by
the Cheeger constant. In particular, it would be ideal to replace the Cheeger constant
by another more tractable quantity which is close to zero
if and only if the Cheeger constant is close to zero, so to somewhat maintain the efficiency of the cut
in the spirit of the Cheeger constant.

An interesting solution for this question is provided by the spectral
gap of a network, given by\footnote{The quantity~$R(g)$ is sometimes called
Rayleigh quotient. \index{Rayleigh quotient}
{F}rom the computational viewpoint, a good approximation for the optimizers of the
Rayleigh quotient can be obtained through power iteration algorithms, see e.g.~\cite{zbMATH01049347}.

As a notational remark, in the sum in the numerator of~\eqref{ORTOGRAFI0},
we have that~$v\sim w$ if and only if~$w\sim v$, but these vertices $v$ and $w$ are counted
only once in the sum.}
\begin{equation}\label{ORTOGRAFI0} \lambda:=\inf R(g),\qquad{\mbox{with}}\qquad
R(g):=\frac{\displaystyle\sum_{v\sim w}(g(v)-g(w))^2}{\displaystyle\sum_{v\in V} \deg(v)\,g^2(v)},\end{equation}
where the above infimum is taken over functions~$g:V\to\R$ (not identically equal to zero) satisfying
\begin{equation}\label{ORTOGRAFI} \sum_{v\in V} \deg(v)\,g(v)=0.\end{equation}

\begin{figure}
  \centering
  \includegraphics[width=.35\linewidth]{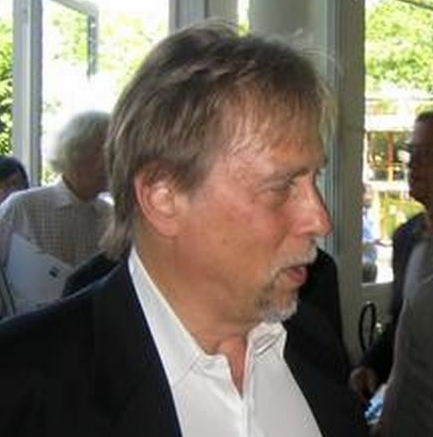}
 \caption{\sl Jeff Cheeger in 2007
 (closeup of an image from
        Wikipedia,
        licensed under the 
 Creative Commons Attribution-Share Alike 2.0 Germany license).}\label{CHEEeFIDItangeFI}
\end{figure}

For our purposes here, the interest of the spectral gap~$\lambda$ is that
it detects the Cheeger constant within a quadratic factor of the optimum, since we have that
\begin{equation}\label{CHEE:INE}
2h\ge\lambda\ge\frac{h^2}{2},
\end{equation}
which is called the \index{Cheeger Inequality}
(discrete) Cheeger Inequality\footnote{See below
for a refinement of this inequality and for its proof.
This type of estimates first appeared in~\cite{MR0402831}
in the context of geometric analysis. Interestingly, mathematical topics happen to be intensively intertwined: in~\cite{MR0402831}, Cheeger ``wishes  to  thank  J.  Simons  for  helpful  conversations, etc.''. See Figure~\ref{SIMoGREENFIDItangeFI} for a picture of
Simons and Figure~\ref{CHEEeFIDItangeFI} for a picture Cheeger.}
for networks.

Also, the understanding of this inequality entails the detection of efficient methods in the search of ``good cuts'',
according to the following algorithm.
One first seeks an optimizer~$g_\star$ for~$\lambda$ in~\eqref{ORTOGRAFI0}.
Then, if the network consists of~$n$ vertices, one reorders them in the form~$v_1,\dots,v_n$ such that
\begin{equation}\label{ORDIgg9} g_\star(v_1)\ge g_\star(v_2)\ge\dots\ge g_\star(v_n).\end{equation}
For every~$j\in\{1,\dots,n\}$, one considers the cut induced by
\begin{equation}\label{02ol34-CHEE:INE:1}
S_j:=\{v_1,\dots,v_j\}\end{equation} and, in the notation of~\eqref{EHOSS}, looks at the optimizer
\begin{equation}\label{02ol34-CHEE:INE:2} \alpha:=\min_{j\in\{1,\dots,n\}}h(S_j).\end{equation}
Now, it turns out that the quantity~$\alpha$, as detected by this explicit algorithm,
is already a ``sufficiently good'' estimate for an optimal cut in terms of the Cheeger constant,
since we have that
\begin{equation}\label{CHEE:INE:2}
2h\ge\lambda\ge\frac{\alpha^2}{2}\ge\frac{h^2}{2},
\end{equation}
which can be considered as a refined version of~\eqref{CHEE:INE}.

We also remark that the link between the spectral gap and the discrete Laplacian
discussed in Sections~\ref{NESNISE}, \ref{I2K243S3M:0oikrfg-4-2rtghj} and~\ref{ITMLDL} becomes clear by
assuming, up to multiplying~$g_\star$ by a scalar factor, that
$$ \sum_{v\in V} \deg(v)\,g_\star^2(v)=1$$
and by
noticing that, by the minimality of~$g_\star$, for every~$\phi:V\to\R$ satisfying~\eqref{ORTOGRAFI}
and every small~$\e$ we have that
\begin{eqnarray*}
0&\ge& R(g_\star)-R(g_\star+\e\phi)\\&=&
R(g_\star)-\frac{\displaystyle\sum_{v\sim w}\Big(
(g_\star(v)-g_\star(w))+\e(\phi(v)-\phi(w))\Big)^2}{\displaystyle\sum_{v\in V} \deg(v)\,(g_\star(v)+\e\phi(v))^2}\\
&=&
R(g_\star)-\frac{\displaystyle\sum_{v\sim w}
(g_\star(v)-g_\star(w))^2+2\e\displaystyle\sum_{v\sim w}
(g_\star(v)-g_\star(w))(\phi(v)-\phi(w)) +o(\e)}{\displaystyle\sum_{v\in V} \deg(v)\,g_\star^2(v)+2\e
\displaystyle\sum_{v\in V} \deg(v)\,g_\star(v)\phi(v)+o(\e)}\\&=&
R(g_\star)-\frac{R(g_\star)+2\e\displaystyle\sum_{v\sim w}
(g_\star(v)-g_\star(w))(\phi(v)-\phi(w)) +o(\e)}{1+2\e
\displaystyle\sum_{v\in V} \deg(v)\,g_\star(v)\phi(v)+o(\e)}\\
&=&-2\e\left(\sum_{v\sim w}
(g_\star(v)-g_\star(w))(\phi(v)-\phi(w))-\lambda\sum_{v\in V} \deg(v)\,g_\star(v)\phi(v)
\right)+o(\e)
\end{eqnarray*}
and therefore
\begin{equation}\label{ORTOGRAFI2} \sum_{v\sim w}
(g_\star(v)-g_\star(w))(\phi(v)-\phi(w))=\lambda\sum_{v\in V} \deg(v)\,g_\star(v)\phi(v).\end{equation}
In particular, for a given~$\overline v\in V$, choosing
$$ \phi(v):=\begin{dcases}
\delta(\overline v) & {\mbox{ if }}v=\overline v,
\\ -\deg(\overline v) &{\mbox{ if }}v\ne\overline v,
\end{dcases}$$
with
$$ \delta(\overline v):=\displaystyle\sum_{v\ne\overline v}\deg(v),$$
we see that~$\phi$ fulfills~\eqref{ORTOGRAFI},
and thus~\eqref{ORTOGRAFI2} reduces to
\begin{eqnarray*}&&
\big(\delta(\overline v)+ \deg(\overline v)\big)\sum_{w\sim \overline v}
(g_\star(\overline v)-g_\star(w)) =\lambda
\left(\deg(\overline v)\,g_\star(\overline v)\,\delta(\overline v)-\deg(\overline v)
\sum_{v\in V\setminus\{\overline v\}} \deg(v)\,g_\star(v)
\right)\\&&\qquad=
\lambda\deg(\overline v)
\big(\delta(\overline v)+\deg(\overline v)
\big)\,g_\star(\overline v),
\end{eqnarray*}
that is
\begin{equation*} {\mathcal{L}}g_\star(\overline v)=-\lambda g_\star(\overline v),\end{equation*}
where\footnote{Up to a possible sign convention, the operator
$$ \sum_{{w\in V}\atop{w\sim v}}g(w)-{\deg (v)} g(v)$$
is sometimes called the Laplacian of a network, see e.g.~\cite[page~64]{MR3821579}.
In this sense, the operator in~\eqref{equNOLAn} plays the role of
a normalized Laplacian on a network, see e.g. equation~(1) in~\cite{MR2416605}.
We also observe that if~$g$ has a minimum at~$v$, then~${\mathcal{L}}g(v)\ge0$.}
\begin{equation}\label{equNOLAn} {\mathcal{L}}g(v):=\frac1{\deg (v)}\sum_{{w\in V}\atop{w\sim v}}
(g (w)-g(v))=
\frac1{\deg (v)}\sum_{w\sim v} g(w)-g(v),\end{equation}
to be compared with~\eqref{DISCRELEA}.\medskip

It is now time to prove~\eqref{CHEE:INE:2}
(which in turn entails~\eqref{CHEE:INE}). To this end, we first consider a partition of vertices~$S_\star$
that achieves the Cheeger constant in~\eqref{EHOSS0}. Up to replacing~$S_\star$ with~$V\setminus S_\star$,
we can assume that~$|S_\star|\le|V\setminus S_\star|$, therefore
\begin{equation}\label{PIUPISASSDA}
|S_\star|\le\frac{|V|}2
\end{equation} and
$$ h=h(S_\star)=\frac{|\partial S_\star|}{|S_\star|}.$$
Let also
$$ \widetilde g:=\chi_{S_\star}-\frac{|S_\star|}{|V|}.$$
We remark that
\begin{eqnarray*}
\sum_{v\in V} \deg(v)\,\widetilde g(v)=
\sum_{v\in S_\star}\deg (v)-\frac{|S_\star|}{|V|}\sum_{v\in V} \deg(v)=
|S_\star|-\frac{|S_\star|}{|V|}\,|V|=0,
\end{eqnarray*}
and therefore~$\widetilde g$ fulfills~\eqref{ORTOGRAFI}.

As a result, by~\eqref{ORTOGRAFI0},
\begin{equation}\label{SOCHEEFA}
\lambda\le\frac{\displaystyle\sum_{v\sim w}(\widetilde g(v)-\widetilde g(w))^2}{\displaystyle\sum_{v\in V} \deg(v)\,\widetilde g^2(v)}
.\end{equation}
Moreover,
\begin{eqnarray*}&&
\sum_{v\sim w}(\widetilde g(v)-\widetilde g(w))^2=\sum_{v\sim w}(\chi_{S_\star}(v)-\chi_{S_\star}(w))^2
=\sum_{S_\star\ni v\sim w\in V\setminus S_\star} 1=|\partial S_\star|
\end{eqnarray*}
and, by~\eqref{PIUPISASSDA},
\begin{eqnarray*}&&
\sum_{v\in V} \deg(v)\,\widetilde g^2(v)=\sum_{v\in V} \deg(v)\,\left[
\left(1-\frac{2|S_\star|}{|V|}\right)\chi_{S_\star}+\frac{|S_\star|^2}{|V|^2}\right]\\&&\qquad=
\left(1-\frac{2|S_\star|}{|V|}\right)\sum_{v\in S_\star} \deg(v)+\frac{|S_\star|^2}{|V|^2}\sum_{v\in V} \deg(v)=
\left(1-\frac{2|S_\star|}{|V|}\right)|S_\star|+\frac{|S_\star|^2}{|V|^2}\,|V|\\&&\qquad=|S_\star|-
\frac{|S_\star|^2}{|V|}\ge\frac{|S_\star|}{2}.\end{eqnarray*}
By plugging this information into~\eqref{SOCHEEFA} we arrive at
\begin{equation}\label{90o24-3iatar:01}
\lambda\le \frac{2|\partial S_\star|}{|S_\star|}=2h.
\end{equation}

Now we recall the notations in~\eqref{02ol34-CHEE:INE:1} and~\eqref{02ol34-CHEE:INE:2}
and we pick the largest~$r\in\{1,\dots,n\}$ for which~$|S_r|\le\frac{|V|}{2}$
(this is possible since~$S_n=\{v_1,\dots,v_n\}=V$). This gives that
\begin{equation}\label{RDEFMIS}
\begin{split}
&|S_j|\le|V\setminus S_j|\qquad{\mbox{for all }}\,j\in\{1,\dots,r\}\\
{\mbox{and }}\qquad\;&|S_j|>|V\setminus S_j|\qquad{\mbox{for all }}\,j\in\{r+1,\dots,n\}.
\end{split}\end{equation}

We also observe that, for every~$c\in\R$ and any function~$g:V\to\R$
satisfying~\eqref{ORTOGRAFI},
$$ \sum_{v\in V} \deg(v) (g(v)-c)^2= \sum_{v\in V} \deg(v) g^2(v)+c^2|V|\ge\sum_{v\in V} \deg(v) g^2(v),$$
therefore
we see that
\begin{equation}\label{fvgb:olpPg64r-g00}
\lambda=
\frac{\displaystyle\sum_{v\sim w}(g_\star(v)-g_\star(w))^2}{\displaystyle\sum_{v\in V} \deg(v)\,g_\star^2(v)}
\ge\frac{\displaystyle\sum_{v\sim w}(g_\star(v)-g_\star(w))^2}{\displaystyle\sum_{v\in V} \deg(v)\,(g_\star(v)-g_\star(v_r))^2}.
\end{equation}
Now we set
\begin{equation}\label{988iuf9wieug938fweuigrfbvBSTwAVFk}\begin{split} &g_1(v):=\begin{dcases}
g_\star(v)-g_\star(v_r) & {\mbox{ if }}g_\star(v)\ge g_\star(v_r),\\
0&{\mbox{otherwise}},
\end{dcases}\\
{\mbox{and }}\qquad&
g_2(v):=\begin{dcases}
g_\star(v_r)-g_\star(v) & {\mbox{ if }}g_\star(v)\le g_\star(v_r),\\
0&{\mbox{otherwise}}.
\end{dcases}
\end{split}\end{equation}
We point out that~$g_\star(v)-g_\star(v_r)=g_1(v)-g_2(v)$ and~$g_1(v)g_2(v)=0$, hence
\begin{equation}\label{fvgb:olpPg64r-g01} (g_\star(v)-g_\star(v_r))^2=(g_1(v)-g_2(v))^2=g_1^2(v)+g_2^2(v).\end{equation}
Furthermore,
\begin{eqnarray*}
&&(g_\star(v)-g_\star(w))^2=(g_\star(v)-g_\star(v_r)+g_\star(v_r)-g_\star(w))^2\\&&\qquad=
(g_1(v)-g_2(v)-g_1(w)+g_2(w))^2\\&&\qquad=
(g_1(v)-g_1(w))^2+(g_2(v)-g_2(w))^2-2(g_1(v)-g_1(w))(g_2(v)-g_2(w)).
\end{eqnarray*}
We also observe that~$g_1(v)\ge0$ and~$g_2(v)\ge0$, therefore
\begin{equation*}
(g_1(v)-g_1(w))(g_2(v)-g_2(w))
=-g_1(v)g_2(w)-g_1(w)g_2(v)\le0.
\end{equation*}
As a consequence,
$$ (g_\star(v)-g_\star(w))^2\ge
(g_1(v)-g_1(w))^2+(g_2(v)-g_2(w))^2.$$
{F}rom this observation, \eqref{fvgb:olpPg64r-g00} and~\eqref{fvgb:olpPg64r-g01} we conclude that
\begin{equation}\label{MediAS-0}
\lambda
\ge\frac{\displaystyle\sum_{v\sim w}
\Big((g_1(v)-g_1(w))^2+(g_2(v)-g_2(w))^2\Big)
}{\displaystyle\sum_{v\in V} \deg(v)\,(g_1^2(v)+g_2^2(v))}.
\end{equation}

We now recall the Mediant Inequality\footnote{To prove~\eqref{MediAS}, one can
suppose that~$\frac{a}{c}\le\frac{b}{d}$ and notice that
$$ \frac{a+b}{c+d}-\frac{a}{c}=\frac{bc-ad}{c(c+d)}=\frac{d}{c+d}\left(\frac{b}{d}-\frac{a}{c}\right)\ge0.$$} for all~$a$, $b$, $c$, $d\ge0$
\begin{equation}\label{MediAS}
\frac{a+b}{c+d}\ge\min\left\{\frac{a}c,\frac{b}d\right\}\end{equation}
and we thereby deduce from~\eqref{MediAS-0} that
\begin{equation}\label{YtEL3c4r2KJ:901oi2ejr3hfg} \lambda
\ge\min\left\{\frac{\displaystyle\sum_{v\sim w}(g_1(v)-g_1(w))^2}{\displaystyle\sum_{v\in V} \deg(v)\,g_1^2(v)},
\; \frac{\displaystyle\sum_{v\sim w}(g_2(v)-g_2(w))^2
}{\displaystyle\sum_{v\in V} \deg(v)\, g_2^2(v)}\right\}=
\min\Big\{ R(g_1),\,R(g_2)\Big\}.
\end{equation}
Now we assume that~$R(g_1)\le R(g_2)$, the other case being similar, and rewrite~\eqref{YtEL3c4r2KJ:901oi2ejr3hfg} in the form
\begin{equation} \label{0OIPOJHKB456789swr4roij-098uij9-098uhnmytsdfvghj}\lambda
\ge R(g_1).\end{equation}

Now we point out that, by the Cauchy-Schwarz Inequality, for every~$g:V\to\R$,
\begin{equation*}
\sum_{v\sim w}|g^2(v)-g^2(w)|=
\sum_{v\sim w}\Big( |g(v)+g(w)|\,|g(v)-g(w)|\Big)\le
\sqrt{\sum_{v\sim w} (g(v)+g(w))^2}\,\sqrt{\sum_{v\sim w}(g(v)-g(w))^2}
\end{equation*}
and therefore
\begin{equation}\label{MBIyhK86gtCtcdeyBCS9i4rtg}\begin{split}&
R(g)=\frac{\displaystyle\sum_{v\sim w}(g(v)-g(w))^2}{\displaystyle\sum_{v\in V} \deg(v)\,g^2(v)}=
\frac{\displaystyle\sum_{v\sim w}(g(v)+g(w))^2\displaystyle\sum_{v\sim w}(g(v)-g(w))^2}{
\displaystyle\sum_{v\sim w}(g(v)+g(w))^2
\displaystyle\sum_{v\in V} \deg(v)\,g^2(v)}\\
&\qquad\qquad\qquad\ge\frac{\displaystyle\left(\sum_{v\sim w}|g^2(v)-g^2(w)|\right)^2}{
\displaystyle\sum_{v\sim w}(g(v)+g(w))^2
\displaystyle\sum_{v\in V} \deg(v)\,g^2(v)}.
\end{split}\end{equation}
Moreover, if~$i<j\in\{1,\dots,n\}$, using a telescoping sum we have that
$$ g^2(v_i)-g^2(v_j)=\sum_{k=i}^{j-1}(g^2(v_k)-g^2(v_{k+1}))$$
and thus, exchanging the order of the summations,
\begin{equation}\label{MCON6tgbGTI81NccKq}\begin{split}&
\sum_{v\sim w}|g^2(v)-g^2(w)|=
\sum_{{v\sim w}\atop{g(v)\ge g(w)}}(g^2(v)-g^2(w))\\&\qquad
=\sum_{{1\le i<j\le n} \atop{g(v_i)\ge g(v_j)}}\chi_E\big((v_i,v_j)\big)(g^2(v_i)-g^2(v_j))\\&\qquad=
\sum_{{1\le i<j\le n}\atop{g(v_i)\ge g(v_j)}}\sum_{i\le k\le j-1} \chi_E\big((v_i,v_j)\big)(g^2(v_k)-g^2(v_{k+1}))
\\&\qquad=
\sum_{1\le k\le n-1}\sum_{{{1\le i\le k}\atop{k+1\le j\le n}}\atop{g(v_i)\ge g(v_j)}} \chi_E\big((v_i,v_j)\big)(g^2(v_k)-g^2(v_{k+1})).
\end{split}\end{equation}

Besides, it is interesting to observe that, by~\eqref{ORDIgg9} and~\eqref{988iuf9wieug938fweuigrfbvBSTwAVFk},
\begin{equation*}
\begin{split} &g_1(v_j)=\begin{dcases}
g_\star(v_j)-g_\star(v_r) & {\mbox{ if }}j\in\{1,\dots,r-1\},\\
0&{\mbox{otherwise}},
\end{dcases}
\end{split}\end{equation*}
and therefore~$g_1(v_i)\ge g_1(v_j)$ if and only if~$i\le j$.

Hence, we can specialize~\eqref{MCON6tgbGTI81NccKq} in the case~$g:=g_1$, obtaining the simpler expression
\begin{equation}\label{MCON6tgbGTI81NccK}\begin{split}&
\sum_{v\sim w}|g_1^2(v)-g_1^2(w)|=
\sum_{1\le k\le n-1}\sum_{{{1\le i\le k}\atop{k+1\le j\le n}}} \chi_E\big((v_i,v_j)\big)(g_1^2(v_k)-g_1^2(v_{k+1})).
\end{split}\end{equation}

Notice also that the expression~$\sum_{{1\le i\le k}\atop{k+1\le j\le n}} \chi_E\big((v_i,v_j)\big)$
counts all the edges joining vertices in~$\{v_1,\dots,v_k\}$ and in~$\{v_{k+1},\dots,v_n\}$, whence,
recalling the notation in~\eqref{02ol34-CHEE:INE:1},
$$ \sum_{{1\le i\le k}\atop{k+1\le j\le n}} \chi_E\big((v_i,v_j)\big)=|\partial S_k|.$$
This and~\eqref{MCON6tgbGTI81NccK} yield that
\begin{equation}\label{MCON6tgbGTI81NccK2}
\sum_{v\sim w}|g^2_1(v)-g^2_1(w)|=
\sum_{1\le k\le n-1}|\partial S_k| (g^2_1(v_k)-g^2_1(v_{k+1})).
\end{equation}

In addition, recalling~\eqref{02ol34-CHEE:INE:2},
$$ \alpha\le h(S_k)=\frac{|\partial S_k|}{\mu_k},\qquad{\mbox{where}}\qquad\mu_k:=\min\big\{|S_k|,|V\setminus S_k|\big\}.$$

As a result,
\begin{eqnarray*}
&& \sum_{1\le k\le n-1}|\partial S_k| (g_1^2(v_k)-g_1^2(v_{k+1}))\ge
\alpha\sum_{1\le k\le n-1}\mu_{k} (g_1^2(v_k)-g_1^2(v_{k+1}))\\&&\qquad=
\alpha\Bigg[ \sum_{1\le k\le n-1}\mu_{k} g_1^2(v_k)
-\sum_{2\le k\le n}\mu_{k-1} g_1^2(v_{k})\Bigg]
\\&&\qquad=
\alpha\Bigg[ \sum_{2\le k\le n-1}\big(\mu_{k}
-\mu_{k-1}\big) g_1^2(v_{k})
+\mu_{1} g_1^2(v_1)-
\mu_{n-1}g_1^2(v_n)
\Bigg].
\end{eqnarray*}
We also observe that~$V\setminus S_n=\varnothing$, whence~$\mu_n=0$. Hence, setting~$S_0:=\varnothing$ we also have that~$\mu_0=0$, thus  we obtain the compact expression
\begin{equation}\label{MCON6tgbGTI81NccK5}\begin{split}&
\sum_{1\le k\le n-1}|\partial S_k| (g_1^2(v_k)-g_1^2(v_{k+1}))\ge
\alpha \sum_{1\le k\le n}\big(\mu_{k} -\mu_{k-1}\big) g_1^2(v_{k})\\&\qquad\qquad=
\alpha \sum_{1\le k\le r-1}\big(\mu_{k} -\mu_{k-1}\big) g_1^2(v_{k}).\end{split}
\end{equation}

It is now useful to observe that, by~\eqref{RDEFMIS},
$ \mu_j=|S_j|$ for all~$j\in\{1,\dots,r\}$. Therefore, if~$k\le r-1$ then
\begin{eqnarray*}
\mu_{k} -\mu_{k-1}=|S_k|-|S_{k-1}|=\sum_{j=1}^k\deg(v_j)-\sum_{j=1}^{k-1}\deg(v_j)=\deg(v_k).
\end{eqnarray*}

Thus, we deduce from~\eqref{MCON6tgbGTI81NccK5} that
$$ \sum_{1\le k\le n-1}|\partial S_k| (g_1^2(v_k)-g_1^2(v_{k+1}))\ge
\alpha \sum_{1\le k\le r-1}\deg(v_k) g_1^2(v_{k})=
\alpha \sum_{v\in V}\deg(v) g_1^2(v).$$
Combining this inequality with~\eqref{MCON6tgbGTI81NccK2}, we have that
\begin{equation*}
\sum_{v\sim w}|g^2_1(v)-g^2_1(w)|\ge\alpha \sum_{v\in V}\deg(v) g_1^2(v).
\end{equation*}
As a consequence, recalling~\eqref{MBIyhK86gtCtcdeyBCS9i4rtg},
\begin{equation}\label{MBIyhK86gtCtcdeyBCS9i4rtg22}
R(g_1)\ge\frac{\displaystyle \alpha^2 \sum_{v\sim w}\deg(v )g_1^2(v)}{\displaystyle\sum_{v\sim w}(g_1(v)+g_1(w))^2}.
\end{equation}

Also, for every~$h:V\times V\to\R$ such that~$h(v,w)=h(w,v)$, we have that
$$ \sum_{v\sim w} h(v,w)=\frac12\sum_{1\le i\le n}\sum_{{1\le j\le n}\atop{v_i\sim v_j}} h(v_i,v_j)$$
where the factor~$\frac12$ is used to avoid the double counting of the edge~$(v_i,v_j)=(v_j,v_i)$, therefore,
for all~$g:V\to\R$,
\begin{eqnarray*}&&\sum_{v\sim w}(g(v)+g(w))^2\le2\sum_{v\sim w}(g^2(v)+g^2(w))
=\sum_{1\le i\le n}\sum_{{1\le j\le n}\atop{v_i\sim v_j}} (g^2(v_i)+g^2(v_j))\\&&\qquad=2\sum_{1\le i\le n}\sum_{{1\le j\le n}\atop{v_i\sim v_j}}g^2(v_i)=2\sum_{1\le i\le n }\deg(v_i)g^2(v_i)=2\sum_{v\in V}\deg(v )g^2(v).
\end{eqnarray*}
This and~\eqref{MBIyhK86gtCtcdeyBCS9i4rtg22} give
\begin{equation*}
R(g_1)\ge\frac{\displaystyle \alpha^2 }{2}.
\end{equation*}

Hence, in light of~\eqref{0OIPOJHKB456789swr4roij-098uij9-098uhnmytsdfvghj},
\begin{equation} \label{0OIPOJHKB456789swr4roij-098uij9-098uhnmytsdfvghj-2}\lambda\ge \frac{\alpha^2 }{2}.\end{equation}
Also, comparing~\eqref{EHOSS0} and~\eqref{02ol34-CHEE:INE:2}, we see that~$\alpha\ge\lambda$.
Combining this information with~\eqref{90o24-3iatar:01}, and~\eqref{0OIPOJHKB456789swr4roij-098uij9-098uhnmytsdfvghj-2}, we obtain~\eqref{CHEE:INE:2}, as desired.
\medskip

See~\cite{MR2744229} for further information about the Cheeger Inequality, its connection
to network partitions, random walks and
algorithm used to rank web pages by search engines
(such as the PageRank algorithm emploited by Google Search).
For more information about the problem of detecting communities in a network, see e.g.~\cite{MR2282139}.
See also~\cite{MR1421568, MR3821579} and the references therein for further readings on related topics.

\subsection{The wear of the rolling stones}\label{ROLSTP}
\begin{figure}
  \centering
  \includegraphics[width=.39\linewidth]{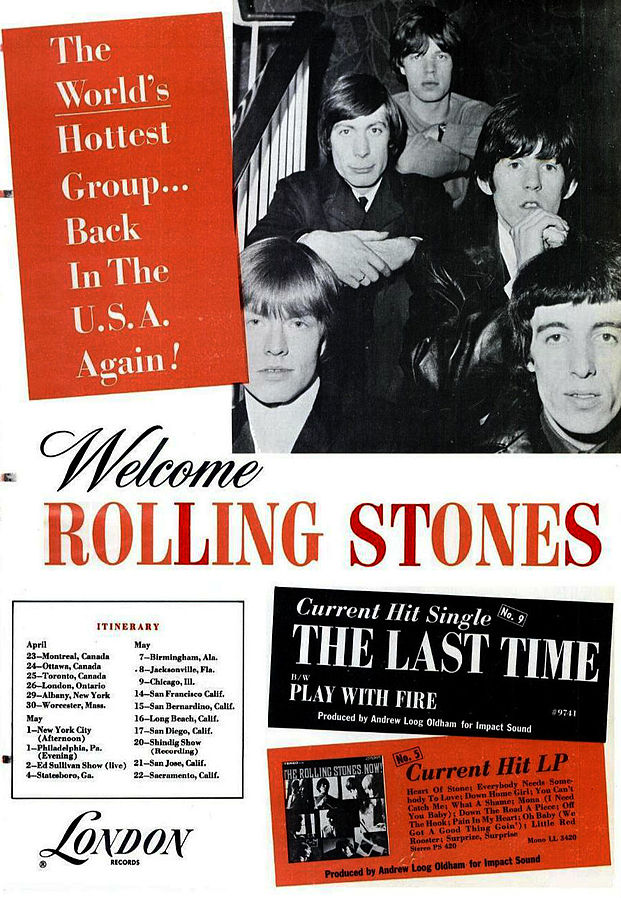}
 \caption{\sl A 1965 trade ad for a tour of a popular English rock band
 (Public Domain image from
 Wikipedia).}\label{BATROLSngeFI}
\end{figure}

\begin{figure}
  \centering
  \includegraphics[width=.62\linewidth]{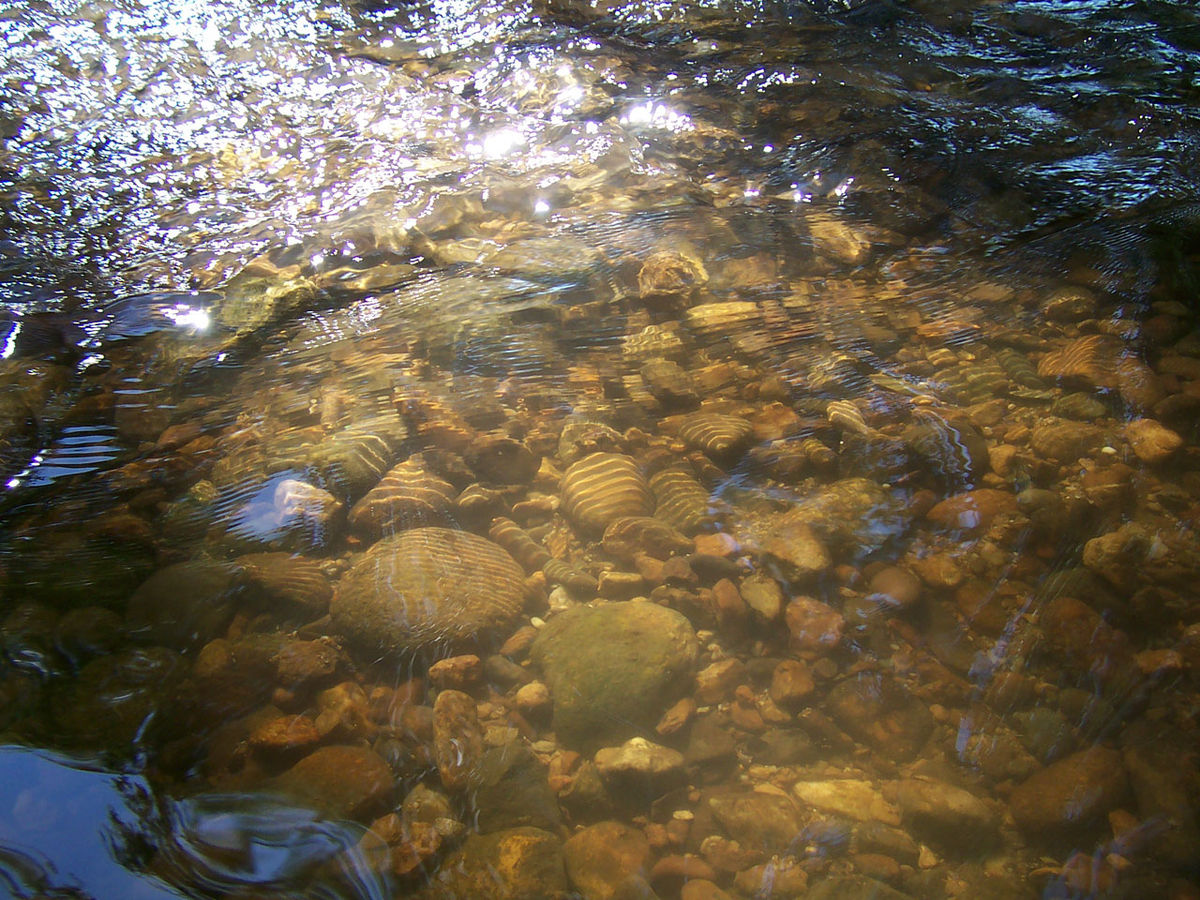}
 \caption{\sl Stones in a riverbed (photo by Fir0002, image from
 Wikipedia, licensed under the Creative Commons Attribution-Share Alike 3.0 Unported license).}\label{P2Dq23we4rtyuikIILANSE}
\end{figure}

Let us now deal with the rolling stones. No, not the Rolling Stones (see Figure~\ref{BATROLSngeFI}), our objective here is to understand the wearing process of stones on beaches and riverbeds which are pounded by waves beautifully shaping their forms (see Figure~\ref{P2Dq23we4rtyuikIILANSE}).

This model was introduced by William J. Firey in~\cite{MR362045} and describes the boundary of the stone at time~$t$ as a two-dimensional, smooth, bounded, strictly convex surface~$M_t$ corresponding to an embedding~$x_t$. In this framework, Firey's equation for the rolling stones reads
\begin{equation}\label{ROLLSTO}
\partial_t x_t=-K_t\,\nu_t,
\end{equation}
where~$\nu_t$ is the exterior normal to~$M_t$ at the point identified by~$x_t$ and~$K_t$ is the Gau{\ss} curvature. \index{Gau{\ss} curvature}

Equation~\eqref{ROLLSTO} is known in jargon as the Gau{\ss} curvature\footnote{Notice that~\eqref{ROLLSTO} prescribes that the normal velocity of the evolving surface coincides with its Gau{\ss} curvature (i.e., the product of its principal curvatures). For planar curves, both the Gau{\ss} curvature and the mean curvature flows coincide with the curve shortening flow. See~\cite{MR4249616} and the references therein for thorough presentations
of geometric flows.

It is also interesting to point out that mean curvature flows, Gau{\ss} curvature flows,
and more general types of geometric flows also play a pivotal role in the description of
bushfire spread, see~\cite{WWMCS}.

A remarkable feature of stones is that, as observed in~\cite{MR362045},
``often stones on beaches pounded by waves wear into quite smooth, regular shapes, sometimes apparently ellipsoidal and even spherical''. It was indeed established in~\cite{MR1714339} that surfaces moving by Gau{\ss} curvature become spherical as they contract to a point.} flow.

To understand the rationale of~\eqref{ROLLSTO}, one can imagine that the wear of a stone comes from its tumbling over an abrasive plane. Suppose that each direction has equal likelihood of
specifying a contact position between the stone and this ideal abrasive plane.
Given a point~$p$ on the surface of the stone, the wear is then proportional to 
\begin{equation}\label{LVuhnGL10ijd9O}\begin{split}&
{\mbox{the ratio between the (infinitesimal) measure of the set of directions for which}}
\\&{\mbox{the abrasive plane
touches the stone in an (infinitesimal) neighborhood of~$p$ and}}\\&{\mbox{the (infinitesimal) measure of the stone's surface of this
(infinitesimal) neighborhood of~$p$.}}\end{split}\end{equation} For the sake of simplicity,
in this process, one disregards any possible
dynamic effect of the global shape of the stone on the tumbling process.

Specifically, up to a rigid motion, we can suppose that~$p$ is the origin of
our reference frame and 
the tangent plane to the stone in~$0$ is horizontal.
In this way, the ($n$-dimensional) stone is parameterized in a small neighborhood of~$0$
by the sublevel sets of some function~$\varphi:\R^{n-1}\to\R$, with~$\varphi(0)=0$ and~$\nabla\varphi(0)=0$.

We observe that, when~$x\in \R^{n-1}$ lies near~$0$,
$$ \varphi(x)=\frac{D^2\varphi(0)x\cdot x}2+ o(|x|^2).$$
Also, the exterior normal to the stone in the vicinity of~$0$ takes the form
\begin{equation}\label{bab02erwR2R23A1}\begin{split}& \nu(x)=\frac{(-\nabla\varphi(x),1)}{\sqrt{1+|\nabla \varphi(x)|^2}}
=\frac{\Big(-D^2\varphi(0)x+o(|x|),1\Big)}{\sqrt{1+\Big|D^2\varphi(0)x+o(|x|)\Big|^2}}\\&\qquad\qquad\qquad=
\frac{\big(-D^2\varphi(0)x,1\big)}{\sqrt{1+\big|D^2\varphi(0)x\big|^2}}+o(|x|)=
-D^2\varphi(0)x+o(|x|).\end{split}
\end{equation}

\begin{figure}
  \centering
  \includegraphics[width=.7\linewidth]{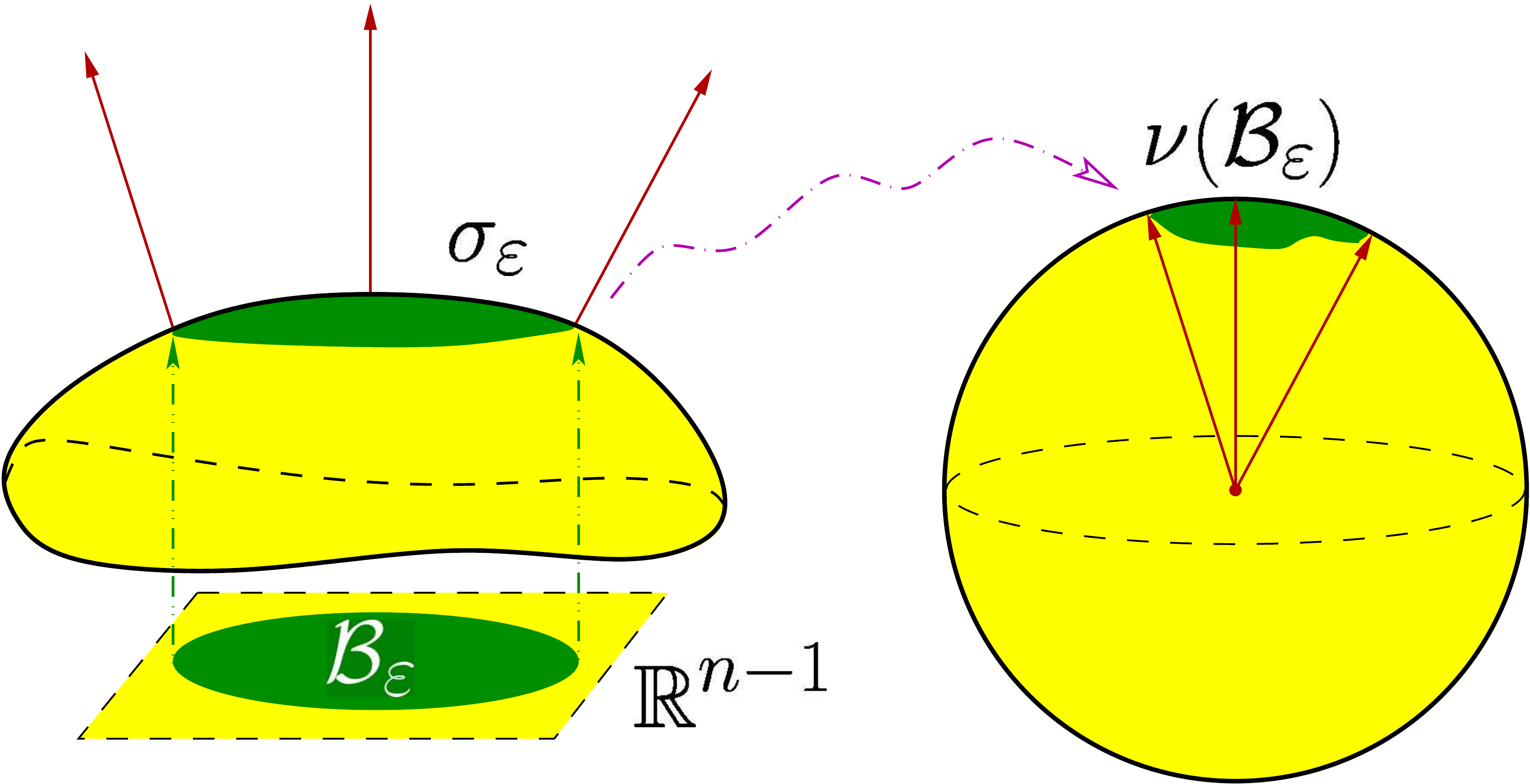}
 \caption{\sl The normal map of a surface (sometimes called Gau{\ss} map).}\label{SIRCGAU0oedfs0I}
\end{figure}

Accordingly, the measure of the set of directions for which the abrasive plane
touches the stone near the origin (i.e., the
measure of the possible tangent planes in the vicinity of the origin) is the $(n-1)$-dimensional
measure of~$\nu({\mathcal{B}}_\e)$ in~$\partial B_1$, where
$$ {\mathcal{B}}_\e:=\{x\in\R^{n-1}{\mbox{ s.t. }}|x|<\e\},$$
see Figure~\ref{SIRCGAU0oedfs0I}.

It is thus convenient to write the surface~$\nu({\mathcal{B}}_\e)$ in a graphical form: namely,
since~$|\nu|=1$, the points~$(y,y_n)\in\R^{n-1}\times\R$ on the surface~$\nu({\mathcal{B}}_\e)$ can be written in the form
$$y_n=\sqrt{1-|y|^2}=:\psi(y),$$
with, owing to~\eqref{bab02erwR2R23A1},
$$y=-D^2\varphi(0)x+o(|x|),$$ for~$x\in {\mathcal{B}}_\e$.

Thus, denoting by~$\widetilde{{\mathcal{B}}_\e}$ the above set of~$y$'s,
and assuming for definiteness that the Hessian matrix~$D^2\varphi(0)$ is invertible,
it follows that
\begin{eqnarray*}&&| \nu({\mathcal{B}}_\e )|=\int_{\widetilde{{\mathcal{B}}_\e}}
\sqrt{1+|\nabla \psi(y)|^2}\,dy=\int_{\widetilde{{\mathcal{B}}_\e}}
\sqrt{1+\left|\frac{y}{\sqrt{1-|y|^2}}\right|^2}\,dy=\int_{\widetilde{{\mathcal{B}}_\e}}
{ \frac{1}{\sqrt{1-|y|^2}} }\,dy\\&&\qquad\qquad\qquad
=\int_{{{\mathcal{B}}_\e}}
\frac1{{\sqrt{1-\left|{ -D^2\varphi(0)x}+o(|x|)\right|^2}}}\,
\left|\det\left( -D^2\varphi(0)+o(1)\right)\right|\,dx\\&&\qquad\qquad\qquad=
(1+o(1))\,\left|\det D^2\varphi(0)\right|\,|{\mathcal{B}}_\e|.
\end{eqnarray*}
If we assume that the stone is a convex body, then~$D^2\varphi(0)\ge0$ and accordingly
\begin{equation*} \big|\det D^2\varphi(0)\big|=\det D^2\varphi(0)=K(0),\end{equation*}
being~$K$ the Gau{\ss} curvature. {F}rom this we arrive at
\begin{eqnarray*} |\nu({\mathcal{B}}_\e )|&=&
(1+o(1))\,K(0)\,|{\mathcal{B}}_\e|\\&=&
K(0)\,\int_{{\mathcal{B}}_\e}\sqrt{1+|\nabla\varphi(x)|^2}\,dx+o(\e^{n-1})
\\&=&K(0)
\,|\sigma_\e|+o(\e^{n-1})
\end{eqnarray*}
being~$\sigma_\e$ the $(n-1)$-dimensional surface describing the stone in the vicinity of the origin
(that is, the $(n-1)$-dimensional surface described by~$\{(x,\varphi(x))$ with
$x\in{\mathcal{B}}_\e\}$).

Thus, by~\eqref{LVuhnGL10ijd9O}, the wear inducing the normal velocity is proportional to
$$\lim_{\e\searrow0}\frac{|\nu({\mathcal{B}}_\e )|}{|\sigma_\e|}=K(0).$$
This justifies equation~\eqref{ROLLSTO}.\medskip

Of course, equation~\eqref{ROLLSTO} presents a natural generalization to the case in which
an additional normal velocity takes part, for instance due to a function~$g$, leading to
\begin{equation*}
\partial_t x_t=(g-K_t)\,\nu_t.
\end{equation*}
Stationary solutions correspond in this case to solutions of
\begin{equation}\label{ROLLSTO2} K=g,\end{equation}
that is a form of prescribed Gau{\ss} curvature equation. \index{prescribed Gau{\ss} curvature equation}

In particular, if we are interested in a graphical framework in which the $n$-dimensional surface under
consideration is (locally) the graph of some function~$u:\R^n\to\R$ (and thus, up to a slight change of notation, also~$g:\R^n\to\R$), we have that
$$ K=\frac{\det D^2 u}{(1+|\nabla u|^2)^{\frac{n+2}2}},$$
see e.g.~\cite[equation~(4.15)]{MR787227}, and thus~\eqref{ROLLSTO2} can be recast in the form
\begin{equation*}
\det D^2u(x) =g(x)(1+|\nabla u(x)|^2)^{\frac{n+2}2}.
\end{equation*}
This type of prescribed Gau{\ss} curvature equation is in fact a particular case of the so-called
Monge-Amp\`ere equation \index{Monge-Amp\`ere equation}
\begin{equation*}
\det D^2u =f,
\end{equation*}
which arises in several aspects of science, including
antenna design, image processing, optimal mass transport, etc.

We refer to~\cite{MR1829162, MR2483373, MR3237759, MR3617963} for further information
about the prescribed Gau{\ss} curvature equation and the
Monge-Amp\`ere equation.\medskip

It is interesting to observe similarities and differences between the Monge-Amp\`ere equation
and the Poisson equation (see footnote~\ref{FOO:POIEQUATYHN} on page~\pageref{FOO:POIEQUATYHN}).
First of all, while the Laplacian operator appearing in the Poisson equation can be interpreted
as the {\em trace} of the Hessian matrix of the solution, the leading operator in the Monge-Amp\`ere equation is the {\em determinant} of the Hessian matrix.

With respect to this algebraic structure, the Monge-Amp\`ere equation
and the Poisson equation are special cases of the so-called $k$-Hessian equations, \index{$k$-Hessian equation}
focusing on $k$th elementary symmetric polynomials\footnote{Specifically,
the $k$th elementary symmetric polynomial of~$n$ variables~$\lambda=(\lambda_1,\dots,\lambda_n)$ is defined by
$$ S_k(\lambda) := \sum_{1\le i_1<\dots<i_k\le n} \lambda_{i_1}\dots
\lambda_{i_k}.$$
The $k$-Hessian equation thus deals with~$S_k(\lambda)$ when~$\lambda$ is the collection
of the eigenvalues (with multiplicity) of the Hessian matrix~$D^2u$.

In particular, $S_1(\lambda)=\lambda_{1}+\dots+\lambda_n$, which produces the Laplace
operator, while~$S_n(\lambda)=\lambda_{1}\dots\lambda_n$, which produces the determinant
operator and leads to the Monge-Amp\`ere equation.}
of the Hessian matrix, see~\cite{MR2500526}.\medskip

Differently from the Poisson equation,
the Monge-Amp\`ere equation shows a severe nonlinear dependence on the Hessian of the solution
and it is indeed an interesting case of
fully nonlinear equation. \index{fully nonlinear equation}
It is however instructive to look at its linearization (see~\cite{MR1439555}): namely, 
if~$A$ is an invertible matrix, $B$ is a matrix, and~$\e$ a small parameter,
$$\det(A+\e B)=\det A\,\big(1+\e{\rm tr} (A^{-1}B)\big)+o(\e)$$
and therefore, if one considers a function~$u=u_0+\e v$ and assumes that~$D^2u_0$ is invertible, then
\begin{equation}\label{MA:009-01}
\begin{split}& \det D^2u=\det(D^2u_0+\e D^2v)=
\det(D^2u_0)\Big(1+\e{\rm tr}\big((D^2u_0)^{-1} D^2v\big)\Big)+o(\e)
\\&\qquad=
\det(D^2u_0)+\e \det(D^2u_0)\sum_{i,j=1}^n a_{ij} \partial_{ij} v+o(\e),\end{split}
\end{equation}
where
\begin{equation}\label{MA:009-02}
{\mbox{$a_{ij}$ is the $(i,j)$th entry of the inverse matrix of~$D^2u_0$.}}
\end{equation} Considering the first order in~$\e$ in~\eqref{MA:009-01}
and dividing by~$\det(D^2 u_0)$ (which is nonzero, since we assumed~$D^2u_0$
to be invertible), the linearized Monge-Amp\`ere equation \index{linearized Monge-Amp\`ere equation}
is therefore written as \label{LINEMONGEALI}
\begin{equation}\label{C2ALPHACHAPET77}
\sum_{i,j=1}^n a_{ij} \partial_{ij} v=h,\end{equation}
for a suitable~$h$ and with the notation in~\eqref{MA:009-02}.\medskip

A particularly interesting structure in equation~\eqref{C2ALPHACHAPET77} arises
under a suitable ``uniform convexity'' assumption on~$u_0$: namely, if we assume that~$D^2 u_0$
is bounded and strictly
positive definite, then also its inverse is bounded and strictly
positive definite. This observation and~\eqref{MA:009-02}
yield that, in this situation, equation~\eqref{C2ALPHACHAPET77} becomes elliptic
(according to the
classification in footnote~\ref{CLASSIFICATIONFOOTN}
on page~\pageref{CLASSIFICATIONFOOTN}).
\medskip

We will get a glimpse on the linear elliptic equations as in~\eqref{C2ALPHACHAPET77}
in the forthcoming Chapters~\ref{C2ALPHACHAPET}
and~\ref{C2ALPHACHAPETil2}.

\subsection{Bushrangers and outlaws}

\begin{figure}
  \centering
  \includegraphics[height=.29\textheight]{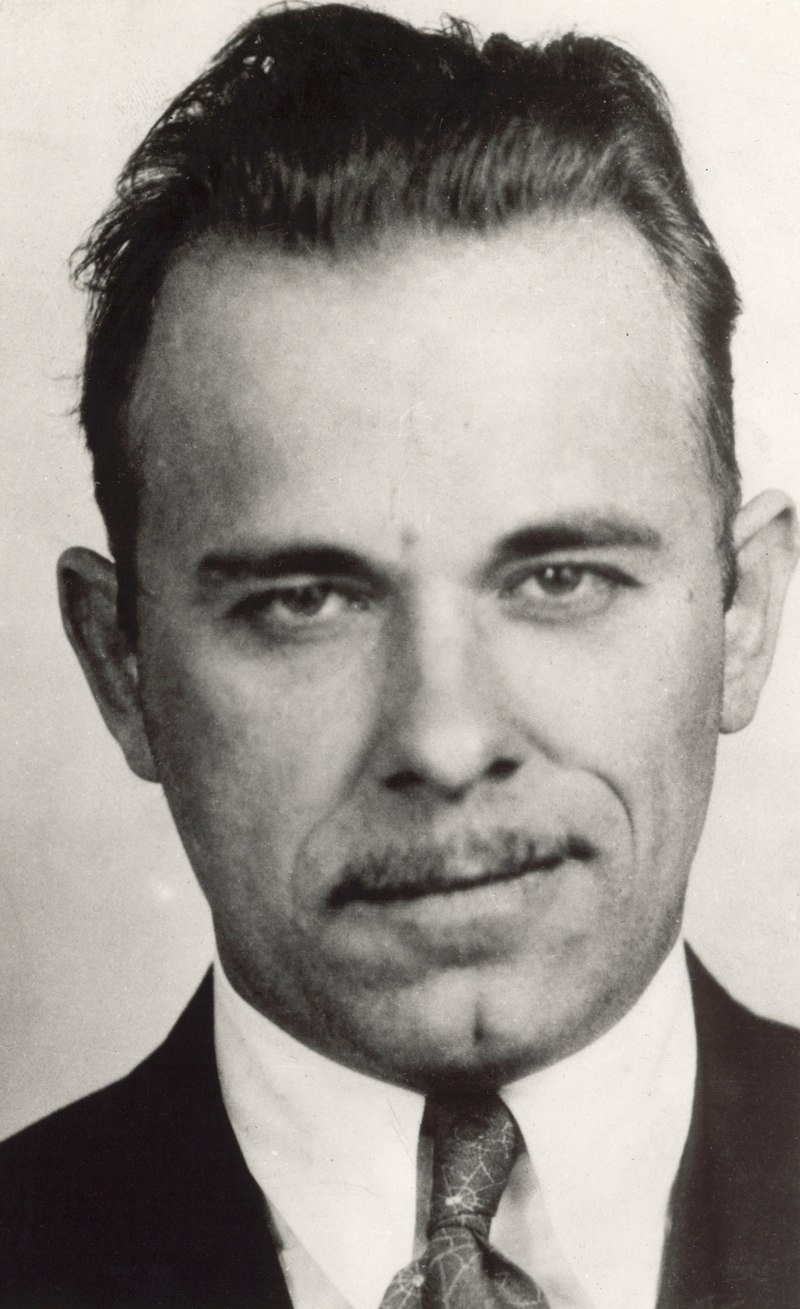}$\quad$
      \includegraphics[height=.29\textheight]{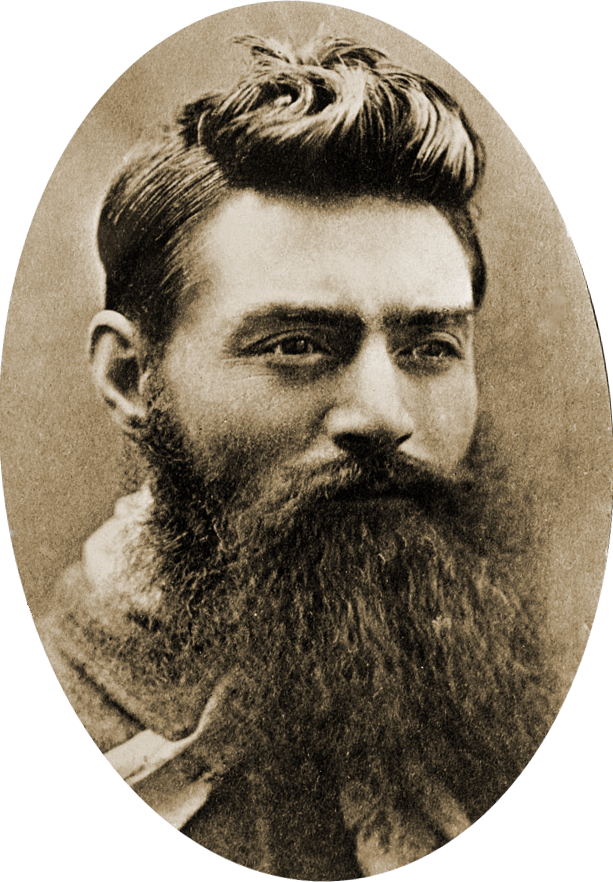}$\quad$
    \includegraphics[height=.29\textheight]{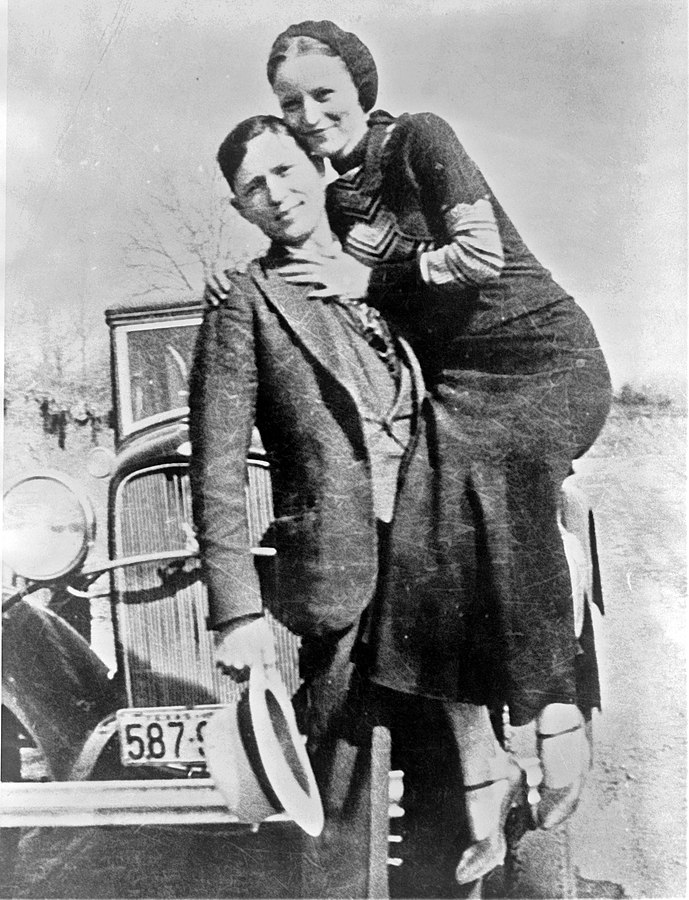}
 \caption{\sl John Dillinger,
 Ned Kelly, Bonnie and Clyde
 (Public Domain images from
 Wikipedia).}\label{VIL:24ENFIDItangeFI}
\end{figure}

\begin{figure}
  \centering
  \includegraphics[width=.53\linewidth]{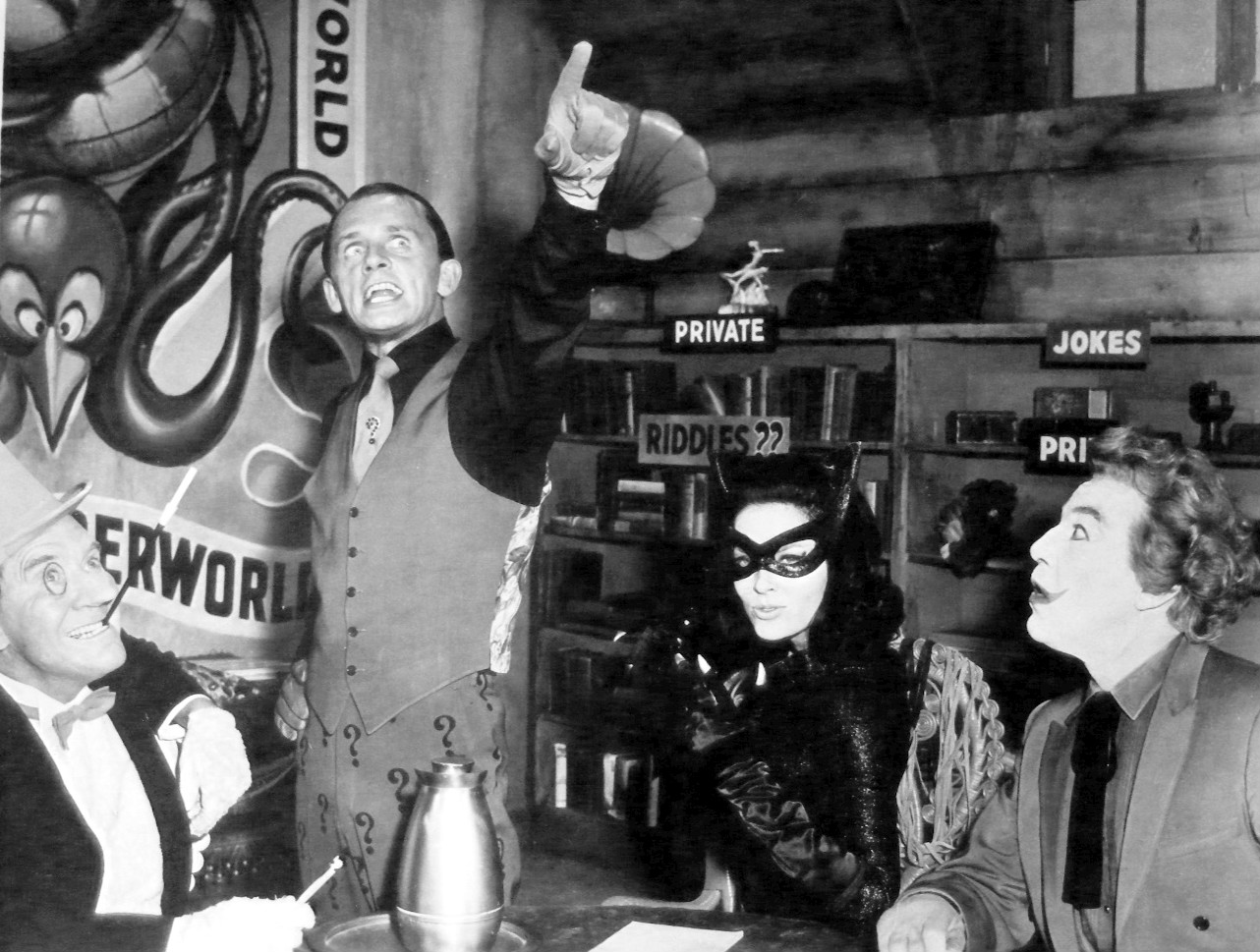}
 \caption{\sl Batman villains: the Penguin, the Riddler, the Catwoman and the Joker in the 1966 film Batman 
 (Public Domain image from
 Wikipedia).}\label{BATMENFIDItangeFI}
\end{figure}

Criminal activities are a common concern for states and citizens
(though criminals, bushrangers and heroic outlaws are sometimes glamorized with a romantic aura, see Figure~\ref{VIL:24ENFIDItangeFI},
and supervillains are popular stock characters in many comic books, see Figure~\ref{BATMENFIDItangeFI}).
Mathematics can certainly play a role in the understanding of criminal
activities,
partial differential equations lend themselves well to a realistic description of several situations
related to transgression of law,
and several research articles have been recently devoted to this topic (see e.g.~\cite{MR2438215, MR2671615}).

These models are actually closely intertwined with the mathematical descriptions of epidemics (see Section~\ref{EPISECT:d})
since they aim at describing ``collective\footnote{The mathematical analysis of collective behaviors has an interesting predecessor in science fiction literature. Namely, the eminent writer (and biochemist) Isaac Asimov in some of his novels discussed about {\em psychohistory}, i.e. a fictional science combining history, sociology, statistics and mathematics in order to make predictions, in average, about the behavior of large groups of people.

As all disciplines based on mathematics, psychohistory possessed its own axioms (1. that the population whose behavior was modeled should be sufficiently large, 2. that the population should remain in ignorance of the results of the application of psychohistorical analyses because if it is aware, the group changes its behavior, etc.). And as all disciplines based on mathematics, psychohistory possessed its own limitations (but disciplines possessing no limitation are a bunch of crooked claptrap). In Asimov's novels, the limitations of 
psychohistory are showcased by its incapability of predicting the appearance of the Mule (mutant warlord endowed with immensely powerful psychic skills, suffering from megalomania and psychopathic paranoia, and aiming at conquering the galaxy).

Perhaps the catch in psychohistory, which was supposed to foresee even thousands of years to come, happens to be the role of the single individuals (as a statistical science, it deals with large groups of people, while the Mule is a ``statistical outlier'') and its inability to quickly evolve to comprise situations in which its own axioms are overturned (e.g., not taking into account the probability of the Mule's mutation). All in all, real, as well as fictitious, predictive sciences are of great importance, but they must be taken with a pinch of salt.}
behaviors'' \index{collective behaviors}
and ``social phenomena'' which
present several treats reminiscent of the spread of an infection. Remarkably, similar models have also
been introduced to describe, for instance, the diffusion of new ideas and scientific knowledge~\cite{cp513-536, MR3471879, MR3519855},
the diffusion of a product in the market~\cite{BASS:BI}, the emergence and establishment of a political party~\cite{MR3412972},
the propagation of riots~\cite{MR4231145}, etc.\medskip

Here, we recall a model for criminal activity put forth in~\cite{MR2671615}. Such a model leverages two basic ans\"atze
from criminology: namely, the {\em repeated victimization theory}, postulating
that places at which certain crimes have been
committed incur a higher risk of this crime being repeated, leading
to a situation in which crime in an area produces more crime,
and the {\em routine activity theory}, according to which
motivating agents, potential benefits and opportunities are key factors for a crime to take place.

To translate these ingredients into a mathematical framework, one can consider the evolution of two functions,
namely the function~$v=v(x,t)$, that describes the ``propensity for people to commit a crime'' at position~$x\in\R^n$ at time~$t>0$, and the function~$u=u(x,t)$, accounting for the number of crimes really occurring at position~$x\in\R^n$ at time~$t>0$. Notice that while~$u$ can be concretely observed and measured in
some explicit, empiric or objective way (e.g., by the number of incidents reported by the police), 
the function~$v$ is mostly a lumped quantity arising from many complex social interactions.
Nonetheless, the susceptibility function~$v$ plays a decisive role to activate and modulate the growth of~$u$,
and conversely the activity level~$u$ induces a feedback on~$v$.

More specifically, to model the number of crimes one can 
assume that the positive values of~$v$ correspond to a tendency to induce criminal activities,
postulate that the number of crimes~$u$, as time flows,
would approach zero exponentially fast (say, like~$\exp(-\lambda_u\, t)$ for some~$\lambda_u>0$)
in the absence of any propensity to commit a crime, but would increase in the presence
of positive values of~$v$ (say, by a number~$\Lambda(v)$, being~$\Lambda$ an increasing function in~$(0,+\infty)$
such that~$\Lambda=0$ in~$(-\infty,0]$). One can also suppose that the presence of a crime in a given location
may trigger unlawful enterprises in neighboring areas and it is suggestive to model this diffusion process via a heat equation,
say with diffusion coefficient~$d_u\ge0$. These considerations lead to an evolution equation for the criminal activities of the form
\begin{equation}\label{0-0-B-soc-1}
\partial_t u(x,t)=d_u\,\Delta u(x,t)+\Lambda\big(v(x,t)\big)-\lambda_u\, u(x,t).
\end{equation}

As for the tendency~$v$ of committing crimes, one can assume that there is an ``innate propension'' of the population
to felonies and quantify it by some function~$\sigma$: this~$\sigma$ may indeed depend on space, to allow for
spatial inhomogeneities, such as bad and dangerous neighborhoods, and even on time, since for instance improved life
conditions, access to healthcare and education systems,
and promotion of cultural activities may contribute to reducing the value of~$\sigma$. Note that while positive values of~$\sigma$ correspond to a natural inclination towards criminal activities,
negative values indicate a natural anticrime tendency.
Thus, in the absence of feedback from~$u$, one could suppose that, as time goes,~$v$
would approach~$\sigma$ exponentially fast (say, like~$\exp(-\lambda_v\, t)$ for some~$\lambda_v>0$).
One can also suppose that the crime propension in a given location
may trigger criminal tendencies in neighboring areas as well,
and again one could model this diffusion process via a heat equation,
say with diffusion coefficient~$d_v\ge0$.
 
Also, in light of the repeated victimization and
routine activity theories,
the feedback of criminal activities on criminal inclination can be modeled by supposing a growth
of criminal intention in direct proportion with the number of critical events, say by a proportionality coefficient~$\varrho$,
possibly varying in time and space (concretely, this~$\varrho$ could measure the benefit of committing the crime,
positive values corresponding to an advantageous payoff of the criminal activity and negative values accounting
for negative consequences of unlawful actions). These observations lead to an evolution equation
for the criminal tendency~$v$ of the form
\begin{equation}\label{0-0-B-soc-2}
\partial_t v(x,t)=d_v\,\Delta v(x,t)-\lambda_v\, \big(v(x,t)-\sigma(x,t)\big)+\varrho(x,t)\,u(x,t).
\end{equation}
The system of equations~\eqref{0-0-B-soc-1} and~\eqref{0-0-B-soc-2} can be seen as an activity/susceptibility model
linking criminal activities with social tendencies and dispositions, and it could also be complemented
with an evolution equation for the benefit function~$\varrho$.

\subsection{Fighting cancer using differential equations}

While body's cells in a normal condition grow and multiply in a controlled way, 
some cells can become abnormal and keep growing, possibly ending up forming a lump called tumor. 
The name of
cancer thus comprises this group of diseases involving abnormal cell growth with the potential to spread to other parts of the body.

There are more than 100 types of cancers which affect humans and
cancer is a leading cause of death worldwide, usually reported to be the second most common cause of
death, after heart diseases. Also,
as the world population is growing and aging, the global number of cancer deaths appears to be increasing.

\begin{figure}
  \centering
  \includegraphics[width=.9\linewidth]{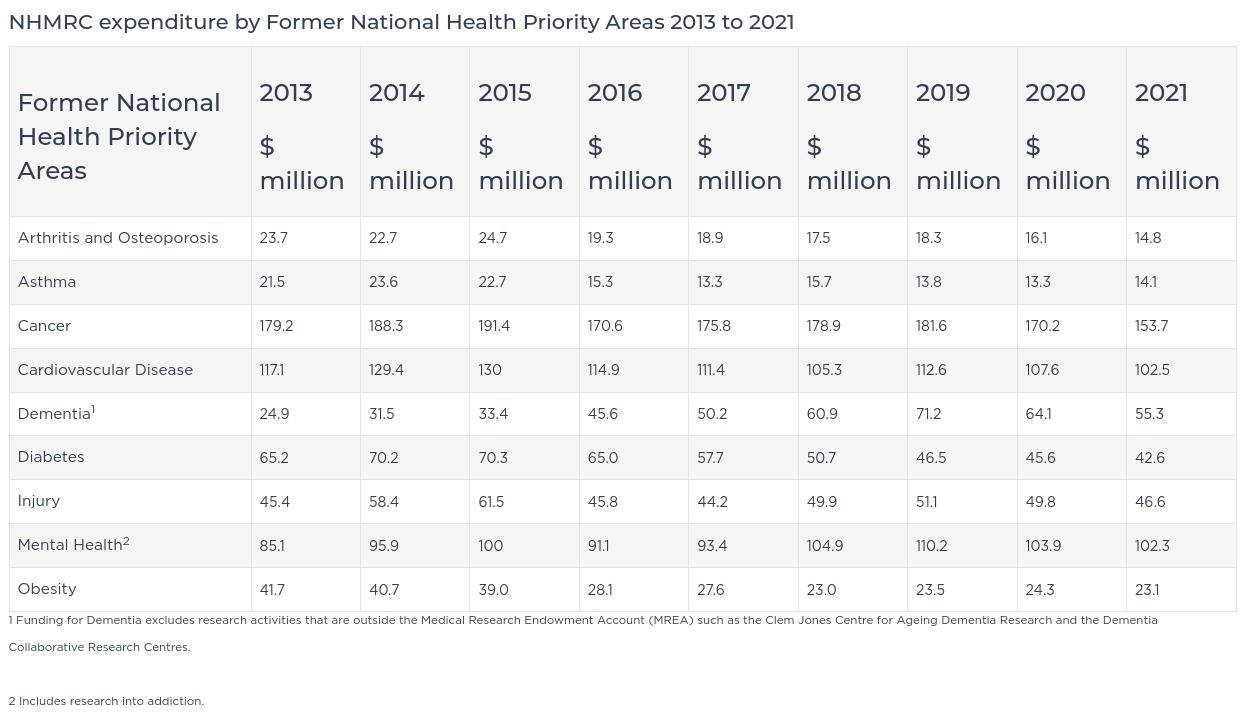}
 \caption{\sl National Health and Medical Research Council funding for major diseases,
table from
{\tt https://www.nhmrc.gov.au/funding/data-research/research-funding-statistics-and-data}.}\label{CAn9ikjd0300293ierjfjTVnTHSnyhJUiFI01}
\end{figure}

\begin{figure}
  \centering
  \includegraphics[width=.9\linewidth]{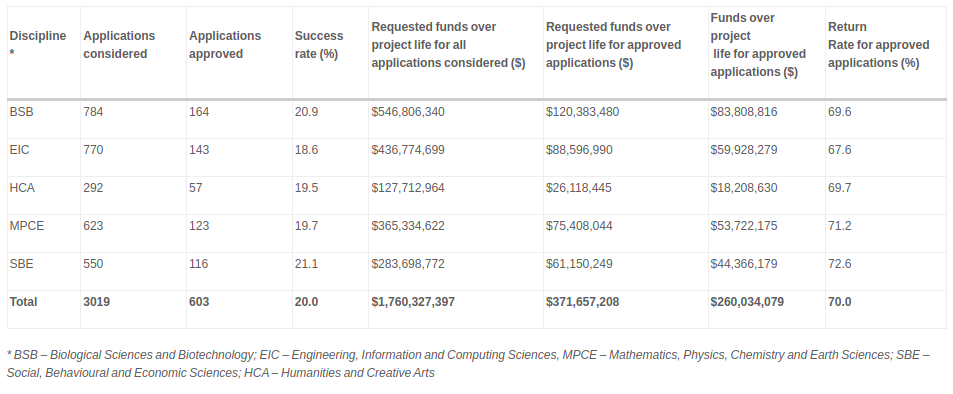}
 \caption{\sl 
Australian Research Council funding for Discovery Projects,
table from {\tt https://www.arc.gov.au/grants/grant-outcomes/selection-outcome-reports/selection-report-discovery-projects-2021}.}\label{CAn9ikjd0300293ierjfjTVnTHSnyhJUiFI02}
\end{figure}

In view of its tragic impact on our lives, and of the obvious importance of
improving early detection methods and treatments,
the support for research to cure cancer is usually prioritized. For instance, in~2021
the expenditure specifically targeting cancer research by the
Australian National Health and Medical Research Council
amounted to~$153.7\times 10^6$ Australian dollars,
i.e., more than~$27\%$ of the total investment\footnote{These figures can be compared, for instance, with those
coming from the Australian Research Council Discovery Projects:
e.g. the total funds allocated for this type of schemes to approved (less than~$20\%$
of the total) applications in all disciplines
of mathematics, physics, chemistry, and earth sciences, amounted in~2021 to
about~$53,7\times 10^6$ Australian dollars, see Figure~\ref{CAn9ikjd0300293ierjfjTVnTHSnyhJUiFI02}.}, see Figure~\ref{CAn9ikjd0300293ierjfjTVnTHSnyhJUiFI01}.
\medskip

Though it is not realistic that mathematics by itself could provide the keystone to defeat cancer, mathematics, and in particular partial differential equations, have provided several tools useful to a deeper and quantitative understanding of many issues related to cancer. Among them, we present here a direction related to cell adhesion as a form of fighting the malignant progression of cancers.

More specifically, in most of the cases, the fatal outcome of cancer is caused not the by the
primary tumor in itself,
but to the invasion by malignant cells to a different body part from where the tumor started
(this phenomenon is named metastasis,
the newly pathological body locations being called metastasis sites, or
metastases). Hence, to expand our knowledge about cancer it could be of great importance
to understand how malignant cells can migrate, and how fast,
and ideally it would be exceptionally useful to detect conditions which can reduce or suppress
tumor invasiveness.

Among these conditions, a pivotal role could be played by \index{cell adhesion}
cell adhesion. This process allows cells to interact and attaches each cell to neighboring ones through specialized molecules of the cell surface. In normal conditions, cell adhesion
maintains cells together and strengthens contact between cells. However,
cancer metastasis can leverage dysfunction of cell adhesion which allows pathological cells
to escape their site of origin and spread through organism, see e.g.~\cite{MTP11}.

Therefore,  one of the directions undertaken to fight cancer is to investigate the role
of cell adhesion as metastasis suppressors, by following the idea that
\begin{equation}\label{L45y5I6T1U44C5E6MI}
{\mbox{cell adhesion is involved in limiting tumor cell migration,}}\end{equation}
see e.g.~\cite{MOH12}.

In relation to this, here we describe
a mathematical model, first introduced in~\cite{MR2279324} for the one-dimensional case,
which analyzes cell spread as induced by standard diffusion
and moderated by cell adhesion, and based on the equation
\begin{equation}\label{MODECA}\begin{split}&
\partial_t u(x,t)=\kappa 
\Delta u(x,t)-\div\Big(u(x,t)\,{\mathcal{K}}_u(x,t)\Big),\\
{\mbox{with }}\quad&
{\mathcal{K}}_u(x,t):= \int_{ B_R}  u(x+y,t)\,\frac{w(|y|)y}{|y|}\,dy
.\end{split}
\end{equation}
Here above, $u=u(x,t)$ denotes the density of cells at some point~$x\in\R^n$ (or more generally at some point in a domain of~$\R^n$) at time~$t>0$. Also, $\kappa\in(0,+\infty)$
and~$w$ is a nonnegative function. 
\medskip

The idea to establish the model in~\eqref{MODECA} is to revisit the derivation of the heat equation put forth in Section~\ref{HEATSEC:02w} by including in the setting the effect of cell adhesion too. To this end, we assume that there is no cell birth or death in our system. The conservation of cell mass thus allows us to follow the strategy described in~\eqref{DAGB-ADkrVoiweLL4re2346ytmngrrUj}: namely, we suppose that the variation of density in a given region of space~$\Omega$ is the flux of some heat flux vector~$B(x,t)$, that is
\begin{equation}\label{1KJSDNc9wov03ifgu9430eirqh89h32Xyr8tgyh4btr2} \partial_t \int_{\Omega} u(x,t)\,dt =-\int_{\partial\Omega} B(x,t)\cdot\nu(x)\,d{\mathcal{H}}^{n-1}_x.\end{equation}
To keep things as simple as possible, one can assume that this flux vector arises as the superposition of two independent effects: classical diffusion (which can be modeled via the Fick's Law in~\eqref{DAGB-ADkrVoiweLL4re2346ytmngrrUj3}) and an additional adhesive flux vector~$A(x,t)$, that is
\begin{equation}\label{1KJSDNc9wov03ifgu9430eirqh89h32Xyr8tgyh4btr} B(x,t)=-\kappa\,\nabla u(x,t)+A(x,t),\end{equation}
for some (say, constant) $\kappa>0$.

To describe the adhesive flux vector~$A$, we can assume that it arises due to molecule interactions within a given distance from a reference point~$x$, say, occurring in a ball~$B_R(x)$. Thus, we suppose that~$A$ is proportional to the adhesive force acting on the cells.
This force comes from the interaction of the biological material located at the point~$x$ with the one located somewhere else in the space (and, say, it is proportional to such biological materials
and directed along the vector joining their positions). Hence, if we assume, for simplicity, that this interaction takes place within a finite distance, this leads to
\begin{equation}\label{LATER67890owdjfc} \begin{split}A(x,t)\,&=
\int_{B_R(x)} u(x,t)\,u(y,t)\,w(|x-y|)\,\frac{(x-y)\,dy}{|x-y|}\\&=
\iint_{(0,R)\times(\partial B_1)} u(x,t)\,u(x+re,t)\,r^{n-1}w(r)\,e\,dr\,d{\mathcal{H}}^{n-1}_e,
\end{split}\end{equation}
where~$w$ is a nonnegative function weighting the interaction force with respect to the distance:
for instance, one could consider~$w$ to be a decreasing function of the distance, or, alternatively, to keep the model as simple as possible, one could also take~$w$ to be equal to~$\frac{1}{|B_R|}$,
in which case the adhesion force contribution would be simply the product between the density
at the point~$x$ and the average density in the ball~$B_R(x)$.

{F}rom~\eqref{1KJSDNc9wov03ifgu9430eirqh89h32Xyr8tgyh4btr} and~\eqref{LATER67890owdjfc}
we arrive at
\begin{equation*} B(x,t)=-\kappa\,\nabla u(x,t)+
\iint_{(0,R)\times(\partial B_1)} u(x,t)\,u(x+re,t)\,r^{n-1}w(r)\,e\,dr\,d{\mathcal{H}}^{n-1}_e.\end{equation*}
Hence, in light of~\eqref{1KJSDNc9wov03ifgu9430eirqh89h32Xyr8tgyh4btr2},
\begin{equation}\label{1KJSDNc9wov03ifgu9430eirqh89h32Xyr8tgyh4btrbi} \begin{split}&
\int_{\Omega}\partial_t u(x,t)\,dt =\kappa\int_{\partial\Omega} 
\nabla u(x,t)\cdot\nu(x)\,d{\mathcal{H}}^{n-1}_x\\&\qquad-
\iiint_{(\partial\Omega)\times(0,R)\times(\partial B_1)} u(x,t)\,u(x+re,t)\,r^{n-1}w(r)\,e
\cdot\nu(x)\,d{\mathcal{H}}^{n-1}_x
\,dr\,d{\mathcal{H}}^{n-1}_e.\end{split}\end{equation}

As a side comment, we notice that the sign difference between the two terms in the right hand side of~\eqref{1KJSDNc9wov03ifgu9430eirqh89h32Xyr8tgyh4btrbi} has an important practical outcome.
To appreciate this difference, one can focus, for instance, on the case of a region~$\Omega$
with high density of cells.
Indeed, recalling that~$\nu$ is the external unit normal (by the notation introduced on page~\pageref{71KJSDNc9wov03ifgu9430eirqh89h32Xyr8tgyh4btr2}), 
the first term in the right hand side of~\eqref{1KJSDNc9wov03ifgu9430eirqh89h32Xyr8tgyh4btrbi} is the one
produced by the Fick's Law and describes the tendency, induced by the gradient, for {\em high density regions to push cells out},
because high density regions correspond to gradients pointing inwards, 
yielding a negative contribution 
to the first term in the right hand side of~\eqref{1KJSDNc9wov03ifgu9430eirqh89h32Xyr8tgyh4btrbi},
which in turn tries to make the density decrease (being the time derivative of the density
on the left hand side of~\eqref{1KJSDNc9wov03ifgu9430eirqh89h32Xyr8tgyh4btrbi}).

Instead, for the second term in the right hand side of~\eqref{1KJSDNc9wov03ifgu9430eirqh89h32Xyr8tgyh4btrbi}, coming from cell adhesion, we can make these considerations.
When~$\Omega$ corresponds to a high density region, given~$x\in\partial\Omega$,
we can imagine~$u(x+re,t)$ to be large when~$e$ points inwards and
small when~$e$ points outwards. Hence, the integrand in the second term of the right hand side of~\eqref{1KJSDNc9wov03ifgu9430eirqh89h32Xyr8tgyh4btrbi} is large when~$e\cdot\nu(x)\le0$
and small when~$e\cdot\nu(x)\ge0$. That is, in case of a high density region~$\Omega$,
the last integral in~\eqref{1KJSDNc9wov03ifgu9430eirqh89h32Xyr8tgyh4btrbi}
provides a negative contribution.
But given the minus sign in front of it, this says that the cell adhesion is
trying to {\em increase the density of a high density region and pull cells in},
which is consistent with the idea that adhesion forces try to make cells to cluster together.

This comment clarifies that Fick's Law and adhesion forces end up having
opposite tendencies
and can be considered as a first indication that the model presented is consistent with the fact highlighted in~\eqref{L45y5I6T1U44C5E6MI} (see page~\pageref{apa45L45y5I6T1U44C5E6MI} below for additional
confirming evidence).\medskip

Now, to obtain~\eqref{MODECA}, one can write~\eqref{1KJSDNc9wov03ifgu9430eirqh89h32Xyr8tgyh4btrbi}
in the form \begin{eqnarray*}&&
\int_{\Omega}\partial_t u(x,t)\,dt =\kappa\int_{\partial\Omega} 
\nabla u(x,t)\cdot\nu(x)\,d{\mathcal{H}}^{n-1}_x\\&&\qquad-
\iint_{(\partial\Omega)\times B_R} u(x,t)\,u(x+y,t)\,w(|y|)\,\frac{y}{|y|}
\cdot\nu(x)\,d{\mathcal{H}}^{n-1}_x
\,dy\end{eqnarray*}
and then
combine this with the Divergence Theorem, finding that
\begin{eqnarray*}&&
\int_{\Omega}\partial_t u(x,t)\,dt =\kappa\int_{\Omega} 
\Delta u(x,t)\,dx\\&&\qquad-\iint_{\Omega\times B_R} \div_x\left( u(x,t)\,u(x+y,t)
\,\frac{w(|y|)\,y}{|y|}\right)
\,dx
\,dy.\end{eqnarray*}
Since~$\Omega$ is arbitrary, we thereby find that
\begin{eqnarray*}
\partial_t u(x,t)&=&\kappa 
\Delta u(x,t)-\int_{ B_R} \div_x\left( u(x,t)\,u(x+y,t)\,\frac{w(|y|)\,y}{|y|}\right)
\,dy,
\end{eqnarray*}
which corresponds to~\eqref{MODECA}.
\medskip

Some further comments on~\eqref{MODECA} are in order.
First of all, the one-dimensional case of~\eqref{MODECA} reduces to
\begin{equation}\label{MODECA-1D}\begin{split}&
\partial_t u=\kappa 
\partial_{xx} u-\partial_x\Big(u\,{\mathcal{K}}_u\Big),\\
{\mbox{with }}\quad&
{\mathcal{K}}_u(x,t):= \int_{-R}^R u(x+y,t)\,\varpi(y)\,dy
,\end{split}
\end{equation}
for an odd symmetric function~$\varpi:\R\to\R$ with~$\varpi\ge0$ in~$(0,+\infty)$.

Furthermore, equation~\eqref{MODECA} (as well as~\eqref{MODECA-1D}) has a nonlocal character,
because of the last term. Indeed, differently from the other terms involved,
to compute~${\mathcal{K}}_u$ (and therefore its derivatives) at a given~$(x,t)$,
it does not suffice to know~$u$ in an arbitrarily small neighborhood of~$(x,t)$,
since a full knowledge of the values of~$u$ in the whole of~$B_R(x)$ (or in a neighborhood of it,
when we want to compute derivatives) is needed.

There are however convenient approximations of~\eqref{MODECA} that attempt to reduce
the model to a more standard, and local, partial differential equation. For example,
assuming small oscillation densities, one could formally look at the expansion
$$ u(x+y,t)\simeq u(x,t)+\nabla u(x,t)\cdot y+\frac12\,D^2u(x,t)y\cdot y+\dots,$$
where higher orders are dismissed.

This approximation, via the definition of~${\mathcal{K}}_u$ in~\eqref{MODECA}, translates into
\begin{eqnarray*}
{\mathcal{K}}_u(x,t)\cdot e_j&\simeq& \int_{ B_R}  \left(u(x,t)+\nabla u(x,t)\cdot y+\frac12\,D^2u(x,t)y\cdot y\right)\,\frac{w(|y|)y_j}{|y|}\,dy+\dots\\ &=&
\int_{ B_R}  \big(\nabla u(x,t)\cdot y\big)\,\frac{w(|y|)y_j}{|y|}\,dy+\dots\\
&=& \xi_j\cdot\nabla u(x,t)+\dots,
\end{eqnarray*}
where
$$ \xi_j:=\int_{ B_R} y\,\frac{w(|y|)y_j}{|y|}\,dy\in\R^n.$$
Correspondingly, we have that
\begin{eqnarray*}
\div\Big(u\,{\mathcal{K}}_u\Big)&=&
\sum_{j=1}^n \partial_j\Big(u\,{\mathcal{K}}_u\cdot e_j\Big)+\dots
\\&=&\sum_{j=1}^n\xi_j \cdot\partial_j\Big(u\,\nabla u\Big)+\dots\\
\\&=&\sum_{j,k=1}^n\xi_{j,k} \,\partial_j\Big(u\,\partial_k u\Big)+\dots,
\end{eqnarray*}
where
$$ \xi_{j,k}:=\xi_j\cdot e_k
=\int_{ B_R} y_k\,\frac{w(|y|)y_j}{|y|}\,dy
=\delta_{jk}\int_{ B_R} \frac{w(|y|)y_j^2}{|y|}\,dy.$$
Therefore
\begin{eqnarray*}
\div\Big(u\,{\mathcal{K}}_u\Big)&=&
\sum_{j=1}^n
c_j\partial_j\Big(u\,\partial_j u\Big)+\dots,
\end{eqnarray*}
with
\begin{equation}\label{GRBDDEGSOsugdjc9uytgfiu3grfu-8} c_j:=\int_{ B_R} \frac{w(|y|)y_j^2}{|y|}\,dy\ge0.\end{equation}
As a result, if one is willing to take such an approximation, the nonlocal equation~\eqref{MODECA}
gets replaced by the local partial differential equation
\begin{equation}\label{GRBDDEGSOsugdjc9uytgfiu3grfu}
\partial_t u=\kappa 
\Delta u-\sum_{j=1}^n
c_j\partial_j\Big(u\,\partial_j u\Big).\end{equation}

It is also interesting to observe that this approximate equation also confirms the agreement between the model presented here and
the fact showcased in~\eqref{L45y5I6T1U44C5E6MI}. To see this, let us consider a point~$x_0$
and a time~$t_0$ where the density achieves a local maximum. We are going to check that
the effect of the first term in the right hand side of~\eqref{GRBDDEGSOsugdjc9uytgfiu3grfu} is to
make the density at~$x_0$ decrease, thus sending out cells, but conversely
the tendency of the second term in the right hand side of~\eqref{GRBDDEGSOsugdjc9uytgfiu3grfu}
is to make such a density increase, hence trying to keep the cells in their original high density location,
and thus providing another confirmation of~\eqref{L45y5I6T1U44C5E6MI}. \label{apa45L45y5I6T1U44C5E6MI}

For this, we suppose that~$x_0=0$, and thus,
up to higher orders, for~$x$ in the vicinity of~$0$, using our local maximality assumption we write
$$ u(x,t_0)=u(0,t_0)-\frac{Mx\cdot x}{2}+\dots,$$
for a symmetric and nonnegative definite matrix~$M$, and therefore~\eqref{GRBDDEGSOsugdjc9uytgfiu3grfu}
gives that
\begin{equation}\label{GRBDDEGSOsugdjc9uytgfiu3grfu-9}
\partial_t u(0,t_0)=-\kappa\,{\rm Tr}M 
+\sum_{j=1}^n
c_j \,u(0,t_0)\,M_{jj}.\end{equation}
By our assumptions on~$M$ and~\eqref{GRBDDEGSOsugdjc9uytgfiu3grfu-8} (and the fact that~$u(0,t_0)\ge0$, being a density),
it follows that the first term in the right hand side of~\eqref{GRBDDEGSOsugdjc9uytgfiu3grfu-9}
is indeed less than or equal to zero, while the second term is greater than or equal to zero,
as claimed.

\subsection{When you are a mathematician}\label{Whenyouareamathematician}

Don Pedro Calder\'{o}n was a pragmatic person. Urologist by profession, he had five children, three daughters, Margarita Isabel (known as Nenacha), Matilde Jr, and María Teresa, and two sons, Alberto Pedro and Calixto Pedro.

\begin{figure}
                \centering
                \includegraphics[width=.43\linewidth]{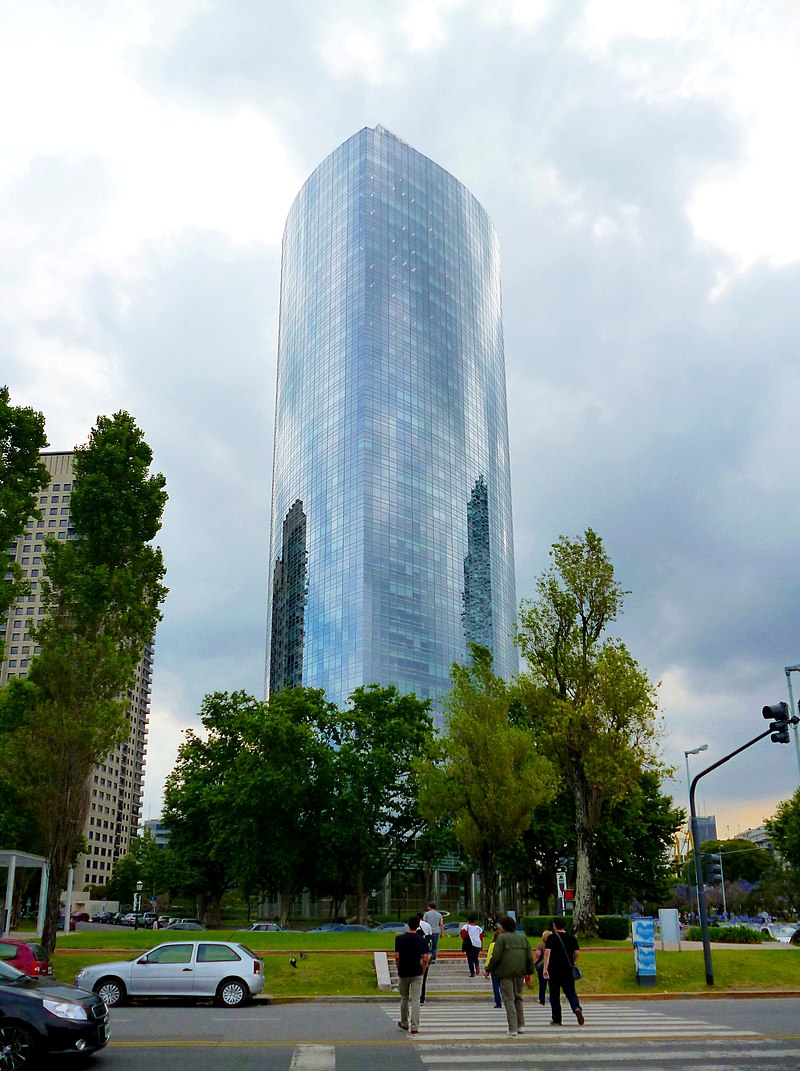}
        \caption{\sl YPF Tower (photo by Allan Aguilar, image from
        Wikipedia,
        licensed under the Creative Commons Attribution 3.0 Unported license).}\label{HA000lFUMHD09876567890NFOJED}
\end{figure}

Don Pedro had a natural affinity for mathematics that he tried to transmit (together with its own name!) to his sons at an early age, by challenging them with rapid mental calculations at the dinner table. Nonetheless, he was a pragmatic person: as such, he firmly believed that a decent person could not make a living as a mathematician and to have a good career one needed to study practical subjects.

\begin{figure}
                \centering
                \includegraphics[height=.53\textwidth]{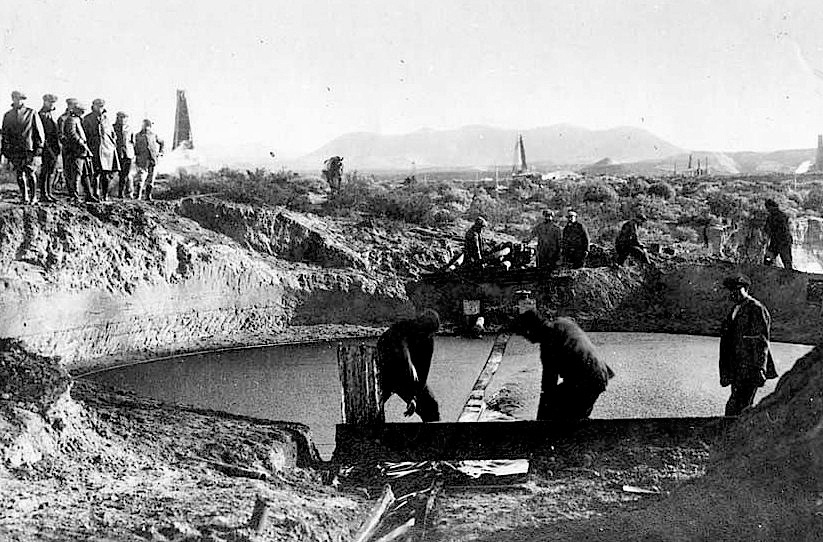}
        \caption{\sl Group of YPF workers at an oil well (Public Domain images from Wikipedia).}\label{YXowrert602394859409-lHAFNEWYDH-3-D-SOFKO}
\end{figure}

Don Pedro happened to be quite a convincing person too, since, in spite of his true love for mathematics, Alberto ended up graduating in civil engineering and, soon after his graduation, working as an engineer for the geophysical division of the Yacimentos Petrolíferos Fiscales, or YPF for short, the Argentinian state oil company (see Figure~\ref{HA000lFUMHD09876567890NFOJED} for the high-rise building located in the Puerto Madero barrio of Buenos Aires hosting nowadays the headquarters of this company, and Figure~\ref{YXowrert602394859409-lHAFNEWYDH-3-D-SOFKO} for a group of YPF workers in the Patagonian Chubut Province).

Unfortunately, or fortunately for Alberto's real passion and for the History of Mathematics, Alberto's work at YPF, according to his recollection, was very interesting, but he was not well treated there, therefore he changed his mind and decided to be a professional mathematician, in spite of his father's resolute advice.

Don Pedro's forecasts turned out to be inaccurate: Alberto soon became a prominent mathematician, with a good salary paid by the University of Chicago.

Don Pedro's expectations on his sons remained unfulfilled, since also Calixto, who firstly enrolled in an engineering degree, rapidly switched to a mathematics degree (also thanks to the strong moral and financial support by his elder brother, who at that point had already become a professor). Thus, Calixto became a professional mathematician (and happened to become the thesis advisor of Luis Caffarelli).

The moral of the story is perhaps that parents don't always know what's best for their sons, as well as for their daughters. 

Or perhaps that mathematics can very often offer rewarding careers and a pleasant lifestyle.

Or perhaps that when you are a mathematician (and you can well be a mathematician even with a degree in civil engineering, or with no degree at all, since being a mathematician has to do with mind and heart, and little to do with paperwork) you end up creating beautiful mathematics even while you are working for a national oil company that treats you badly. This is precisely what happened to Alberto Calder\'{o}n during his employment period with the Yacimentos Petrolíferos Fiscales, let's see why
(and see also footnote~\ref{Whenyouareamathematician2} on page~\pageref{Whenyouareamathematician2} for more about Alberto Calder\'{o}n).

The electrical conductivity of an inhomogeneous body varies significantly, due to changes of temperature, density and ionic concentration. Therefore, variations of the electrical conductivity of the earth may reveal the presence of oil or gas, helping with the identification
of productive zones to drill. No wonder this was an interesting topic for YPF!

The issue is therefore to measure the electrical conductivity of some regions in the earth in a way which is as simple, affordable and efficient as possible. The idea of drilling the ground to examine samples, for instance, does not really fulfill these requirements, being rather expensive, risky, invasive, as well as environmentally problematic. Instead, how about performing a geophysical prospection for oil and gas exploration\footnote{Not only this method can be useful for the detection of oil and gas, it can also
provide valuable information about the water contaminated with hydrocarbons which is produced during the extraction. This is called in jargon ``produced water'' and needs to be managed correctly to avoid serious
environmental consequences. Therefore, it is important, before the extraction, to have an estimate on the proportion of water that will be present in the so-called ``multiphase flow''
(i.e. the mixed combination of different substances, such as oil, gas and water) which will be taken out of the ground.

Actually, besides its environmental importance, the methods presented here have also an economical impact, since
water is typically separated from the oil and gas by a company specialising in this process, not by the extractor company itself, therefore a careful estimate of these additional costs plays a role in the extractor company's financial plan.

Another important application of the determination of the electrical conductivity, is its use, along with other measurements, to enhance drilling safety by preventing blowout. 
In this framework, a high variation of the values of the conductivity may indicate that the properties of the rock at a certain
depth have changed and actions have to be taken to prevent a retaining wall from collapsing. 

Once again, it would be ideal to detect relative proportions of oil, gas and water without going through complicated, expensive or intrusive procedures, such as sampling methods or gamma rays produced by radioactive sources. The method that Calder\'{o}n studied and appropriate modifications of it have the potential to provide helpful strategies.

To appreciate how electrical conductivity can be efficiently utilized to tell different substances apart, let us compare the different orders of magnitude (measured in siemens per meter) of several media, such as gold~$4.11\times10^7$, sea water~$4.8$, drinking water from~$5\times10^{-4}$ to~$5\times10^{-2}$,
air from~$10^{-15}$ to~$10^{-9}$.
In general, gas on its own has very small ability of conducting electricity, but its conductivity increases in the presence of ions;
also electric charge carriers are absent in oil, hence oil has a very small conductivity too.
The electric conductivity of the soil change considerably according to its salinity, but, to have a rough idea, sand has a relatively low conductivity (from~$1\times10^{-3}$ to~$10m\times10^{-3}$), silt has a medium conductivity (from~$8\times10^{-3}$ to~$8\times10^{-2}$) and clay has a high conductivity (from~$2\times10^{-2}$ to~$8\times10^{-1}$).}
only through simple measurements on the surface of a region? What if one could use some cheap electrical methods for this goal, such as voltage and current measurements on the surface of the earth?\medskip

To provide a mathematical model for the problem under consideration, we recall that the electrical conductivity measures how easily electric current can flow through a given material. More specifically, the resistivity~$\rho$ at a particular point of a given material is defined\footnote{To clarify the jargon, we recall that both the terms ``resistance'' and ``resistivity'' describe how difficult it is to make electrical current flow through a material. The only difference is that resistivity is an intrinsic property of the material, while resistance depends on the macroscopical object the material is shaped into.
For example, the electrical resistance~$R$ of a conductor of uniform cross section with sectional area~$A$ and length~$\ell$ is related to the electrical resistivity~$\rho$ through the formula$$R= \frac{\rho\ell }{A}.$$
That is, the resistivity is a specific property of the material, the resistance depends on the area and length of the conducting wire.

Strictly speaking, for alternating currents (say, in sinusoidal form) it is also appropriate to replace the notion of resistance with that of ``electrical impedance'', which is the ratio of voltage between the terminals of a circuit and the current flowing through it (both in complex notation, so that the impedance is a complex number, the real part of it corresponding to the resistance).

To be totally accurate, here we should take into account this notion of impedance (or of its reciprocal, which is called ``admittance''), since in practice the electrical impedance tomography that we describe uses either alternating currents at a single frequency or multiple frequencies. Since the basic ideas remain unchanged, we surf over this detail.}
as the ratio of the electric field~$E$ to the density of the current~$J$ that~$E$ creates at that point.
The conductivity~$\sigma$ is thus defined as the reciprocal of the resistivity~$\rho$. 

In this way,
$$ \rho=\frac{E}{J}\qquad{\mbox{and}}\qquad
J=\frac{E}{\rho}=\sigma E.$$
Accordingly, for a given electric field~$E$, a higher conductivity material produces more current flow~$J$ than a low conductivity material (hence the name of ``conductivity'').

Since the electrostatic force is conservative, we can also write~$E=\nabla u$, for some potential function~$u$. Also, by Amp\`ere's Law, we know that the curl of the magnetic field~$H$ is equal to the electric current density~$J$ and consequently
$$ \curl H=J=\sigma E=\sigma\nabla u.$$
As a result, recalling the vector calculus identity in~\eqref{YY08ojew9oik903hyto3298ythgfjb8ihyri-1urX},
\begin{equation} \label{CALDE1}0= \div(\curl H)=\div(\sigma\nabla u).\end{equation}
The geological problem considered by Calder\'{o}n was then to apply a boundary voltage potential~$\phi$ to some given region~$\Omega$. In this framework, equation~\eqref{CALDE1} is taking place in the domain~$\Omega$ and is complemented by the boundary condition
\begin{equation}\label{CALDE2}
u=\phi\,{\mbox{ on }}\,\partial\Omega.
\end{equation}
In this setting, one can measure the energy needed to maintain the potential~$\phi$ along~$\partial\Omega$ as
$$ Q_\sigma(\phi)\,:=\,\int_\Omega\sigma(x)|\nabla u(x)|^2\,dx.$$
Note that, by the Divergence Theorem, \eqref{CALDE1} and~\eqref{CALDE2},
\begin{eqnarray*} 
Q_\sigma(\phi)&=&\int_\Omega\Big( \div\big( \sigma(x)  u(x) \nabla u(x)\big)
-u(x)\div\big( \sigma(x) \nabla u(x)\big) \Big)\,dx\\&=&
\int_\Omega \div\big( \sigma(x)  u(x) \nabla u(x)\big)\,dx\\&=&\int_{\partial\Omega} u(x)\sigma(x) \partial_\nu u(x) \,d{\mathcal{H}}^{n-1}_x\\&=&\int_{\partial\Omega} \phi(x) \sigma(x) \partial_\nu u(x) \,d{\mathcal{H}}^{n-1}_x.\end{eqnarray*}
We stress that, for all~$x\in\partial\Omega$, the quantity~$\sigma(x) \partial_\nu u(x)$ accounts for the normal component of the electric current at the boundary.

\begin{figure}
  \centering
  \includegraphics[width=.47\linewidth]{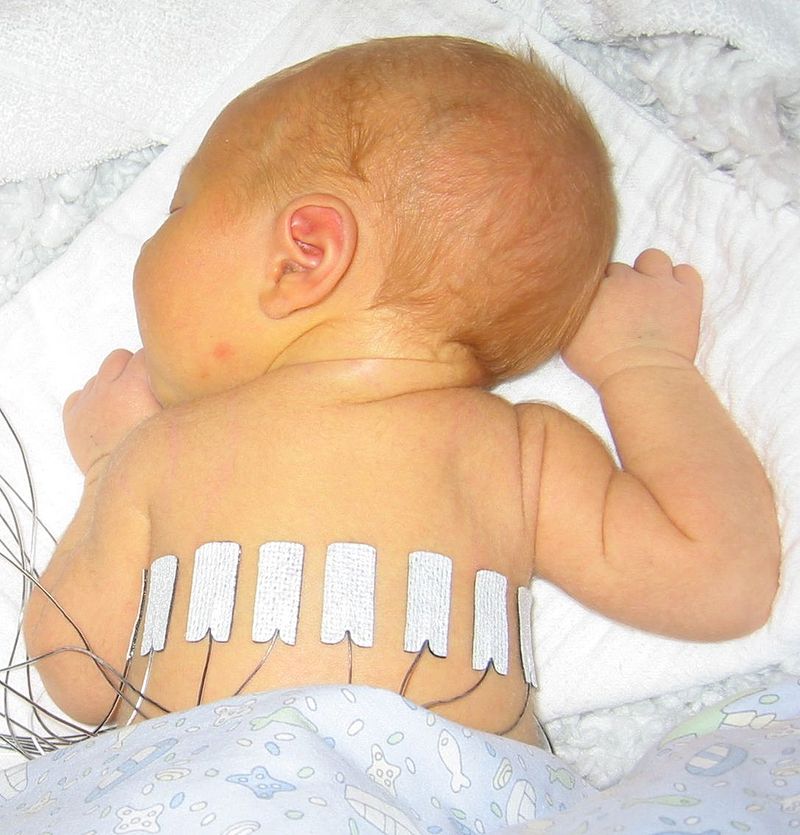}
 \caption{\sl Electrodes attached for electrical impedance tomography lung monitoring
 on a ten-day old spontaneously breathing neonate
 (photo by S. Heinrich, H. Schiffmann, A. Frerichs, A. Klockgether-Radke and I. Frerichs; image from
 Wikipedia, licensed under the Creative Commons Attribution-Share Alike 3.0 Unported license).}\label{22nbgt5bfr4uh8u0uIUHB9203ruofeh-0987689iokmnbg7tug}
\end{figure}

Calder\'{o}n's problem \index{Calder\'{o}n's problem} thus consisted in whether or not one can recover, at least approximately, the electrical conductivity~$\sigma$ in~$\Omega$ by the knowledge of~$Q_\sigma(\phi)$ which can be measured along~$\partial\Omega$: notice the ``noninvasivity" of these measurements in order to determine the inner structure of the material, since the necessary detections can be carried out along the boundary of~$\Omega$, i.e. without penetrating inside the domain.

Some questions posed explicitly by Calder\'{o}n in his seminal paper~\cite{MR590275} were: is~$\sigma$ uniquely determined by~$Q_\sigma(\phi)$? Can one calculate~$\sigma$ in terms of~$Q_\sigma(\phi)$?

Though the problem remains largely open in its generality, Calder\'{o}n already provided several interesting results, for instance about the possibility of determining uniquely the electrical conductivity~$\sigma$ from the knowledge of~$Q_\sigma(\phi)$, up to a small error, under the assumption that~$\sigma$ is sufficiently close to a constant.

Quite remarkably, for two-dimensional domains~$\Omega$, the electrical conductivity~$\sigma$ is completely determined by~$Q_\sigma(\phi)$, see~\cite{MR2195135}. See also~\cite{MR2435340, MR2472875, MR3460047} for more information\footnote{We stress that Calder\'{o}n's problem is perhaps the prototypical ``inverse problem'', \index{inverse problem} in which the objective is not to derive information on the solution of a certain problem given a set of data and of physical laws, but rather to obtain information about a physical object from observable measurements.

In other words, typical mathematical problems aim to predict physical observations given the values of the parameters defining the model and the underlying physical laws. Instead, an inverse problem utilizes the results of the observations to infer the values of the parameters characterizing the system under investigation.

Curiously, the groundbreaking paper~\cite{MR590275} has remained the only article by Calder\'{o}n on the topic of inverse problems.}
about Calder\'{o}n's problem.
\medskip

The idea of determining the internal structure of a body through electrical measurements on its surface is called in jargon ``electrical impedance tomography'' \index{electrical impedance tomography} and it is useful not only to provide an educated guess on whether or not it is worth drilling the soil with the expectation of finding oil or gas. Indeed, given that it is a noninvasive,
nondestructive, portable, and inexpensive methodology, one of the promising applications of this idea lies in the field of medical diagnostics, as pioneered in~\cite{WT89-03Lip}.

{F}rom the clinical standpoint, the key fact is that electrical conductivity highly varies among different tissues and it is deeply influenced by the movement of fluids and gases within a given tissue (for instance, muscle and blood have higher electric conductivity than fat and bones). The use of the electrical impedance tomography in this situation is typically realized by attaching electrodes to the skin around the body part to be examined, see Figure~\ref{22nbgt5bfr4uh8u0uIUHB9203ruofeh-0987689iokmnbg7tug}. 

\begin{figure}
  \centering
  \includegraphics[width=.44\linewidth]{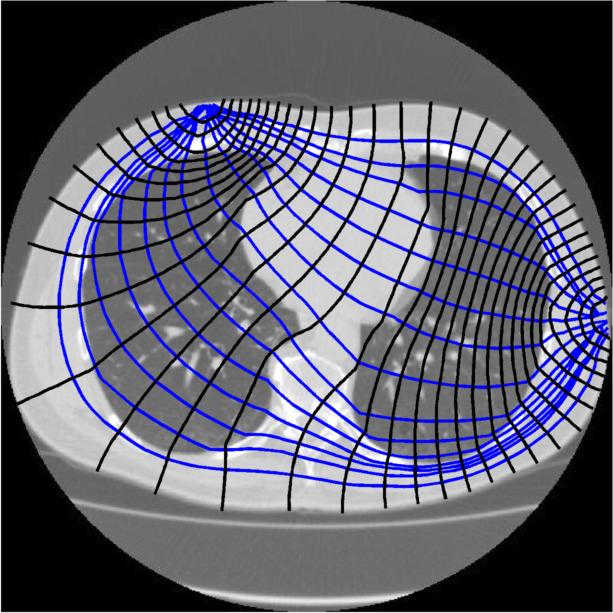}
 \caption{\sl A cross section of a human thorax with electrical current streamlines and equipotentials
 (imahe by Andy Adler from
 Wikipedia, licensed under the Creative Commons Attribution-Share Alike 3.0 Unported license).}\label{nbgt5bfr4uh8u0uIUHB9203ruofeh-0987689iokmnbg7tug}
\end{figure}

The streamlines and equipotential surfaces are detected and, after several repetitions in different configurations, the results are elaborated and reconstructed into a two-dimensional image which shows the bending of streamlines and equipotentials caused by the change in conductivity between different tissues, see Figure~\ref{nbgt5bfr4uh8u0uIUHB9203ruofeh-0987689iokmnbg7tug}.

A topical clinical application consists for instance in the monitoring of lung functions. The electrical impedance tomography possesses in this case very promising potentials, since the conductivity of lung tissues is about five times lower than the soft tissues in the thorax, thus allowing for a significant contrast. Furthermore, the conductivity highly oscillates depending on the presence of the air in the lungs, since air plays the role of an insulating factor, thus detecting vital activities such as inspiration and expiration.

The use of electrical impedance tomography for lungs monitoring is particularly important to detect possible injuries caused by mechanical ventilation, thus helping physicians in a bespoke adjusting of ventilators' settings to provide lung ventilations to patients needing it without increasing too much the risk of collateral damages caused by this artificial process.\medskip

Another important medical application of electrical impedance tomography relates to breast imaging for cancer detection, to provide a possible alternative to mammography, which uses more invasive ionizing radiations, and magnetic resonance imaging, which employs as contrast media chemical substances which can cause problems to patients with kidney failure (also, due to their low specificity, these classical techniques are subject to a relatively high rate of very distressive false positives). In this setting, the electrical impedance tomography has promising potentials, due to the different electrical conductivities of normal and malignant breast tissues.
\medskip 

A further application is related to the monitoring of brain function and neuronal activity. In this context, the effectiveness of the method relies on the fact that cerebral hypoglycemia and hypoxemia, ischemia and haemorrhage are typically associated with significant changes of the electrical properties of the brain.\medskip

In some cases, the clinical applications require a refinement of the mathematical techniques. For instance, muscle tissues typically present great differences in the electrical conductivities in different directions (e.g. the conductivity of the cardiac muscle changes by a factor 2.7 circa between the transverse and the longitudinal directions). This requires mathematicians to take into account also the case of ``anisotropic conductivities'': in this situation,
rather than considering a scalar function~$\sigma$ in~\eqref{CALDE1}, it is customary to account for conductivity variations in different directions by taking into consideration a symmetric matrix valued function~$\sigma(x)=\{\sigma_{ij}(x)\}_{i,j\in\{1,\dots,n\}}$.
This casts~\eqref{CALDE1} into the more general form\footnote{For the reader interested in geometric analysis
and differential geometry, let us mention that in dimension~$n\ge3$ equation~\eqref{K0987678908yhM-olk4DIF} corresponds to a geometric equation on a Riemannian manifold with boundary, see~\cite{MR1029119}. The development of this geometric intuition was very fruitful to construct counterexamples to the uniqueness of
the electric conductivity~$\sigma$ in~\eqref{K0987678908yhM-olk4DIF} (this, in the geometric setting, corresponds to the uniqueness
of a Riemann metric related to
harmonic functions; the differential operator used for this is the Laplace-Beltrami operator which will be presented in some detail in Section~\ref{BELS}).}
\begin{equation}\label{K0987678908yhM-olk4DIF} \sum_{i,j=1}^n\partial_i(\sigma_{ij}\partial_j u)=0.\end{equation}
\medskip

The mathematical difficulties of dealing with a nonlinear and ill-posed inverse problem translate into severe practical problems that the adoption of the electrical impedance tomography for medical purposes has still to overcome. In particular, the image reconstruction methods involved are hard since there are typically more than one possible solution corresponding to a three-dimensional region projected onto a two-dimensional surface (see e.g.~\cite{MR2024725}
for a mathematical description of the lack of uniqueness problem). Often, to have a reliable reconstruction, a detailed knowledge of the
exact dimensions and shape of the organ under consideration is required, as well as a very precise electrode location. In addition, sometimes the fluctuations of electrical quantities due to anomalies are still not sufficiently broad to compensate for reconstruction errors.

For all these reasons, the clinical use of electrical impedance tomography is still considered mostly experimental. However, some commercial devices relying on electrical impedance tomography
have been already introduced to monitor lung function in intensive care patients, and hopefully future coordinated progress in the mathematical understanding of the problem, in the precision of the technology used for the measurements and in the algorithms for the image reconstruction will lead in the future to a secure use of electrical impedance tomography in several areas of great social and medical importance.

\subsection{If you are shy, do as the Romulans do}

The power of invisibility is a very cool superpower to possess, especially if one is shy. This power is actually an evergreen classic. It appeared already around year 1100 in the Welsh mythology (a ``mantle of invisibility'' being one of King Arthur's most valuable possessions) and it surfaced also in Norse mythologies, Celtic stories, Grimm's fairy tales and Japanese traditions (in which the ``kakuremino'' is a useful raincoat providing also invisibility).

\begin{figure}
  \centering
  \includegraphics[height=6.9cm]{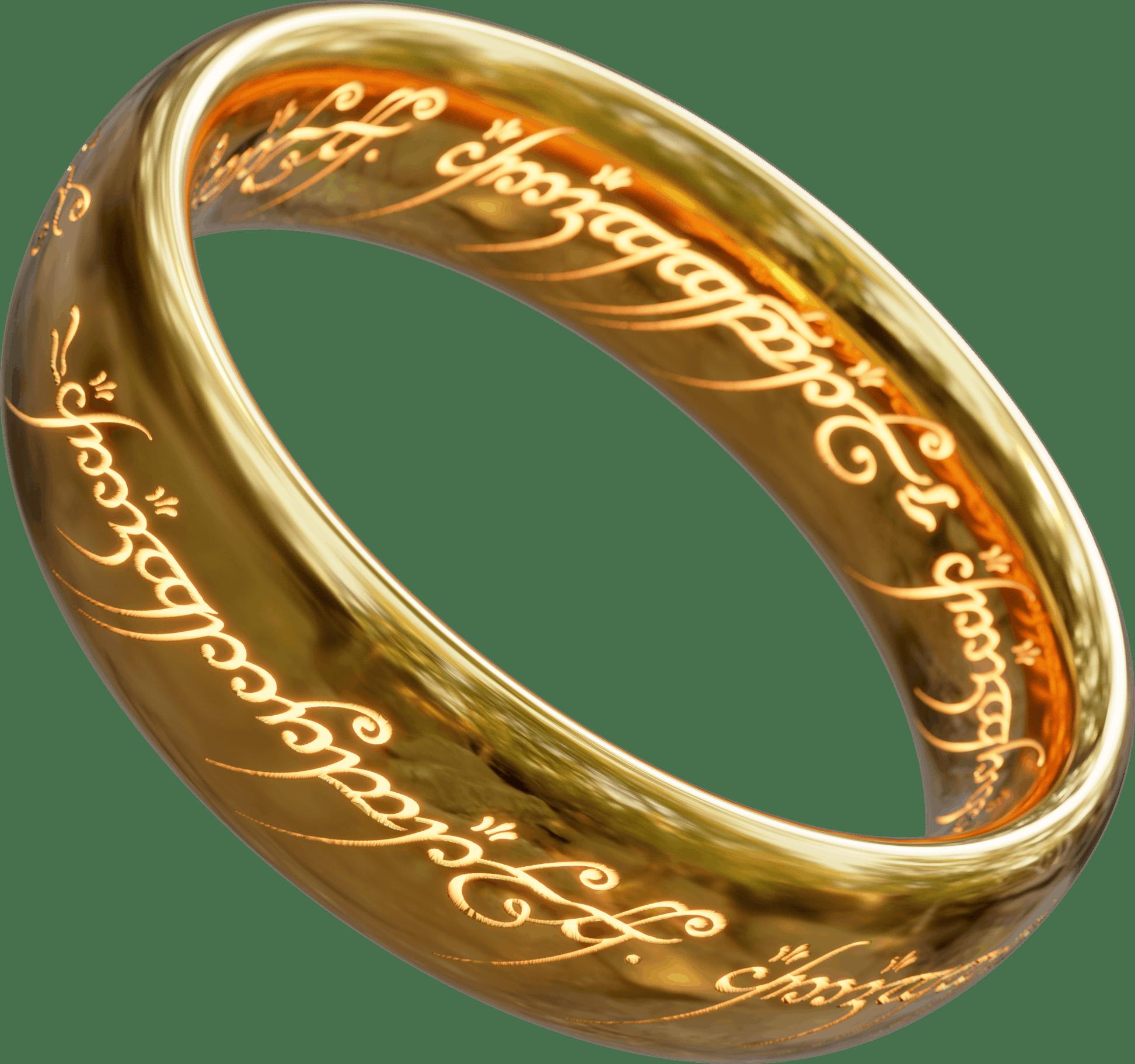}
 \caption{\sl The One Ring
 (image from Wikipedia, own work by Peter J. Yost,
        licensed under the Creative Commons Attribution-Share Alike 4.0 International license).}\label{1RinGG0i3ujehNIRDGREENFIDItangeFI}
\end{figure}

Invisibility plays also a pivotal role in J. R. R. Tolkien's novels ``The Hobbit'' and ``The Lord of the Rings'', in which the One Ring is presented as a magic device forged by the Dark Lord Sauron granting the wearer invisibility (see Figure~\ref{1RinGG0i3ujehNIRDGREENFIDItangeFI}; as additional benefits, any inherent power of its owner would be automatically amplified, the life of a mortal possessor of this device would extend indefinitely, aging would be stopped, and domination over the wills of others would be guaranteed, at the sole side effect of a malevolent and manipulative agency on the possessor's ego).

\begin{figure}
  \centering
  \includegraphics[height=8.5cm]{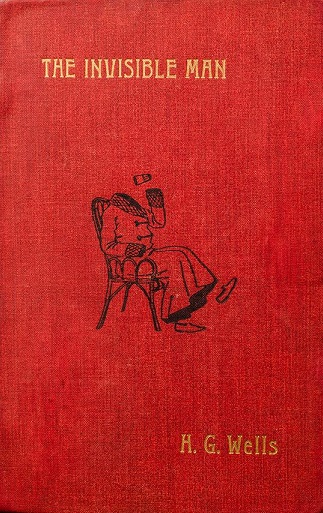}  $\quad$ \includegraphics[height=8.5cm]{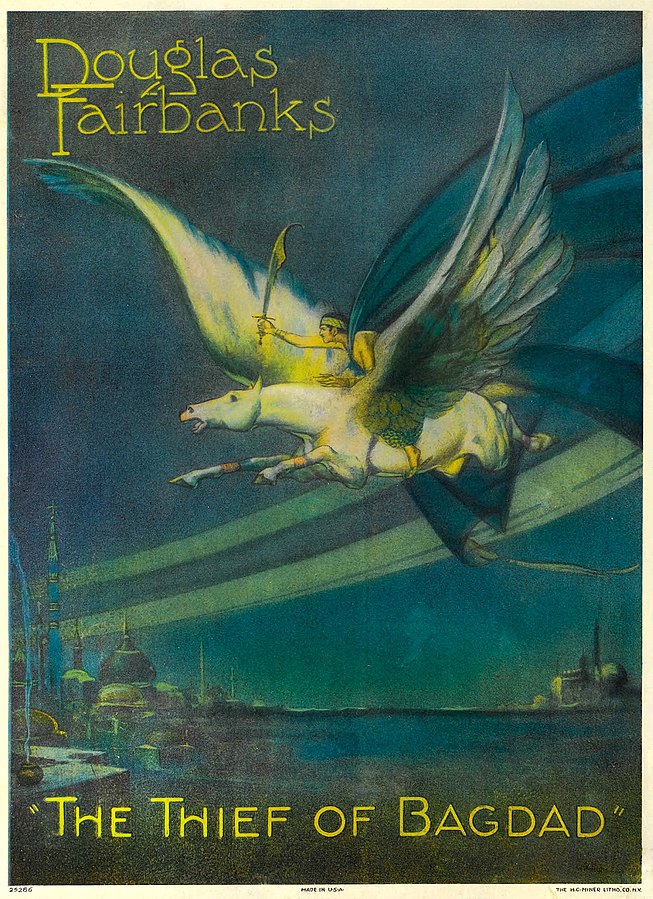}
 \caption{\sl First edition cover of ``The Invisible Man'' (left, 1897) and 
 poster for ``The Thief of Bagdad'' (right, 1924; 
Public Domain images from
 Wikipedia).}\label{09i090ohbdRFVFn0987666561qazdcfyhjk78Evdfu8VVbsdic-F2I4}
\end{figure}

Also in cinema, invisibility dates back at least to 1924 (see Figure~\ref{09i090ohbdRFVFn0987666561qazdcfyhjk78Evdfu8VVbsdic-F2I4}), since the main character Ahmed of the silent film ``The thief of Bagdad'' gains at some point a cloak of invisibility (but no big deal, the movie also presents a crystal ball from the eye of a giant idol showing whatever one wants to see, a magic powder turning people into whatever they wish, a flying carpet, a magic rope, and a magic apple which cures any disease including\footnote{Not sure how to eat an apple {\em after} death, though.} death).

Invisibility is also one of the cornerstones of science fiction. One of the first and most influential science fiction novels, ``The Invisible Man'' by H. G. Wells (see Figure~\ref{09i090ohbdRFVFn0987666561qazdcfyhjk78Evdfu8VVbsdic-F2I4}) is about a mysterious man, Griffin, referred to as ``the stranger'', who happens to be a scientist who manages to render himself invisible (sadly, failing in his attempt to reverse this condition).

In recent novels and movies, let us mention that cloaks of invisibility are abundantly employed in the Harry Potter series.

With comics, the power of invisibility has maybe come to an extreme. For instance, Susan (Sue) Storm, the Invisible Woman of the Fantastic Four (and the first female hero in the Marvel Universe) possesses extended invisibility powers allowing her to do more than just turn invisible (she can also make the objects and the people around her invisible, as well as creating invisible\footnote{Not sure what a {\em visible} force field is, though.} force fields).\medskip

All good, but the main question is: how would you manage to turn yourself invisible? 

The idea of Griffin in the novel by H. G. Wells is very ingenuous: using his knowledge of optics, he invents a way to modify the refractive index of his own body making it precisely equal to that of air (in this way, essentially, he neither absorbs or reflects light, appearing transparent to everybody). Very ingenious indeed, but maybe presenting some side effects, since, as pointed out by \label{Yakov Perelman}
Yakov Perelman in his book ``Physics can be fun'', this method should have man made Griffin not only invisible but blind as well (because the eye needs to absorb incoming light to transmit visual signals to the brain). On top of that, the fact that Griffin could not restore its refractive index to its normal value suggests that this technique may be too intrusive and not easily reversible, so perhaps not recommendable.

As for Sue Storm, her superhuman powers are the outcome of the exposure to massive cosmic radiation
(which happens to be a clich\'e in superheroic plots): again, this doesn't sound like a viable method, due to possible side effects.

As for Harry Potter, his cloak of invisibility was passed down to him directly by  his father (but the cloak was previously owned by Death), so, once again, unless you come from a wealthy family of sorcerers, this strategy may not be the most practicable one.\medskip

\begin{figure}
  \centering
  \includegraphics[height=9.7cm]{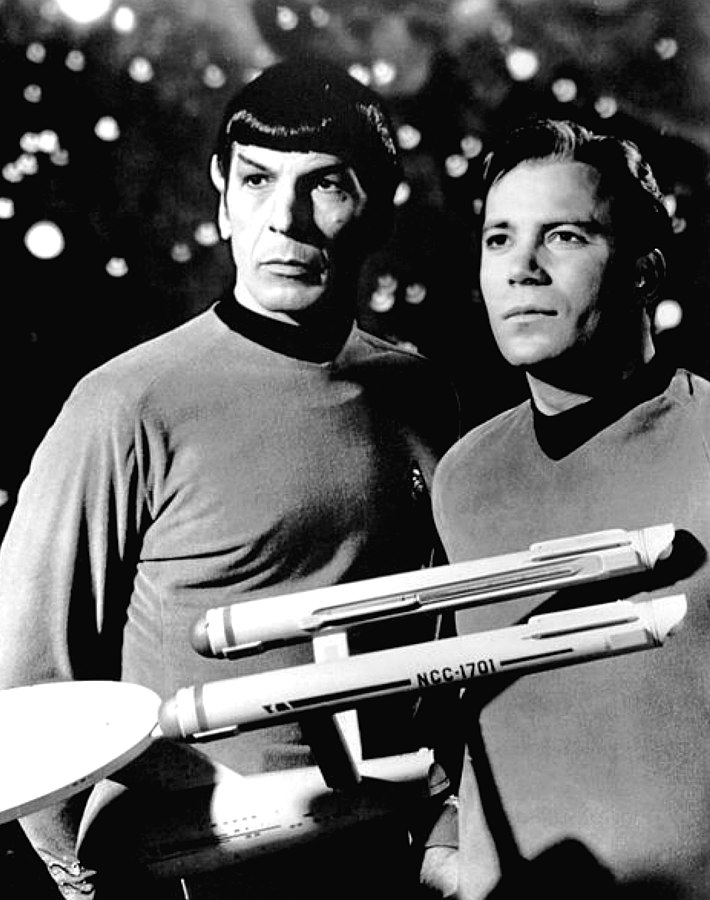}
 \caption{\sl Mr. Spock and Captain Kirk from the television series Star Trek
 (with the USS Enterprise starship passing by on the bottom; Public Domain
 image from Wikipedia).}\label{STGHBS1Riek89nGG0i3ujpoijrvtd4aehNIRDGREENFIDItangeFs2I56}
\end{figure}

Fortunately however, there is a rather promising strategy that was put forth in the 1966 episode ``Balance of Terror'' of the \label{1s2ddc23425qwedfgboprenc0813ti0a4g5gda}
science fiction saga of ``Star Trek'', see Figure~\ref{STGHBS1Riek89nGG0i3ujpoijrvtd4aehNIRDGREENFIDItangeFs2I56}. Specifically, in that episode, our heroes of the USS Enterprise starship had to battle their archenemies Romulans who employed a new type of space vessel which was capable of remaining invisible. The Romulan technology was based\footnote{But after all the idea of the Romulans of guiding light rays so that they are reflected past their spaceship (so that the viewer would have the impression that the spaceship is not there) is a natural one, since it is at the basis of various optical illusions, such as the mirages in the desert.

Anyway, while the Romulans were created specifically for the episode ``Balance of Terror'',
they have remained in the Star Trek saga to embody an aggressive, warlike, alien species having founded an interstellar empire. Not to be confused with the Vulcans, or Vulcanians, who are a peaceful, methodical, logic-oriented species, and perhaps quite unemotional, kind of alien people (Mr. Spock being a prominent representative). Both Romulans and Vulcans are however characterized by pointed ears and
slanted eyebrows, so the two populations have very similar physiological traits.

Yet, perhaps the primary antagonists in the classical Star Trek saga were the Klingons,
an alien species of brutal moods practicing feudalism and authoritarianism.
Likely, the use of the bloodthirsty Klingons over the militaristic but technologically sophisticated Romulans
in the original series was motivated by politics and budget reasons.
Indeed, the makeup necessary to make the Romulans was at that time rather costly, while the one
for the Klingons was cheap and quick (the application of shoe polish being enough to
provide Klingons with the bronze skin which was basically their sole physical difference with respect to humans). Also, in the middle of the Cold War, the physical appearance of Klingons was probably meant to
suggest orientalism, at a time when memories of Japanese actions during World War II were still fresh,
and their feudalistic society was perhaps alluding to Soviet Union, with the 
``United Federation of Planets'' to which the Enterprise crew belongs playing the role of the United States.} on a {\em selective bending of light} (which was dubbed ``cloaking device'' in the 1968 episode ``The Enterprise Incident'').\medskip

\begin{figure}
  \centering
  \includegraphics[width=.56\linewidth]{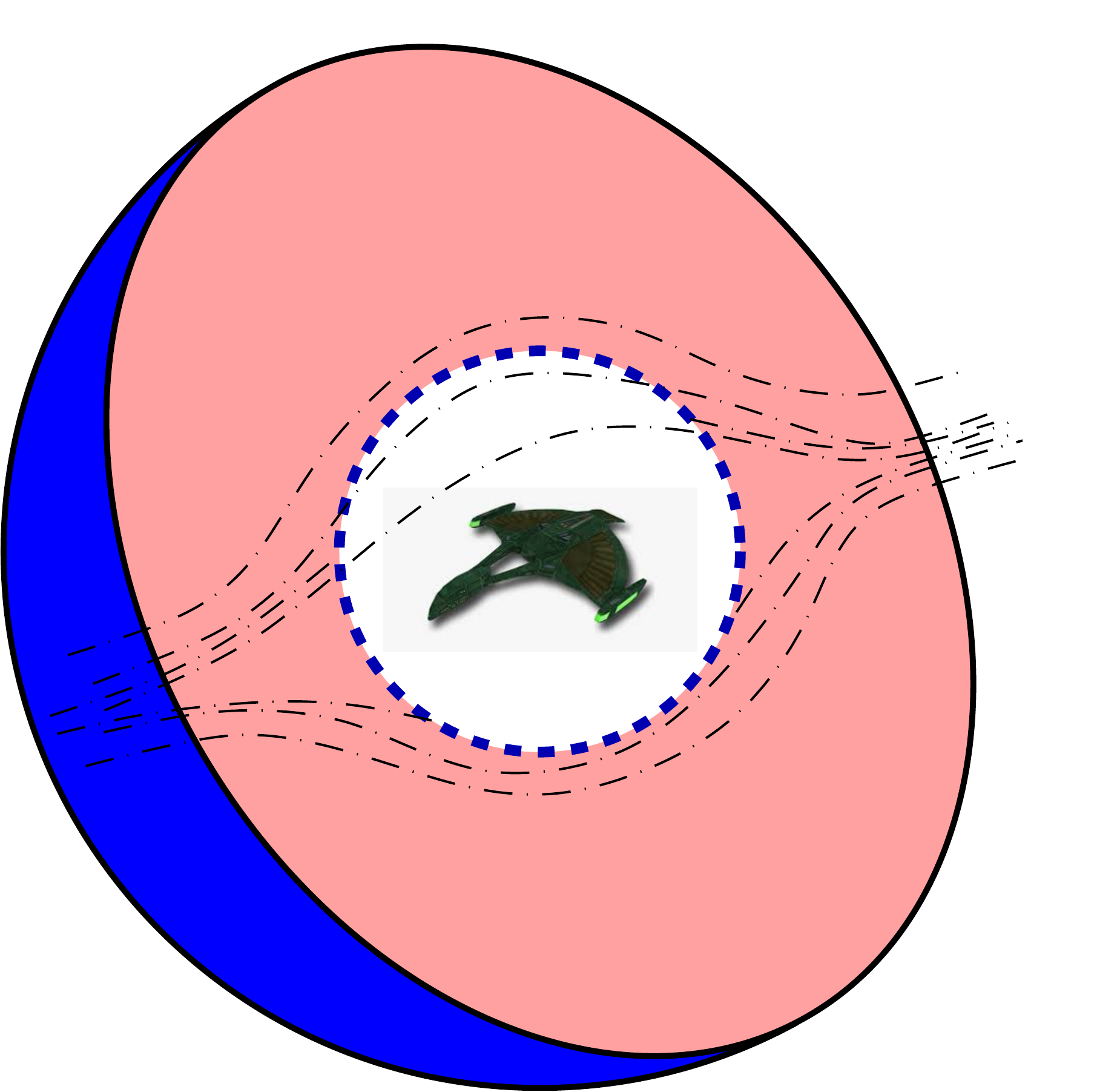}
 \caption{\sl Making a Romulan spacecraft invisible
 by guiding the light rays away via Maxwell's equations.
 Deflecting light rays around the vessel and restoring them on the other side
 make the observer believe that light
 had passed through empty space.}\label{0okjnh78uigfi34bqtb76b67c3b2b376b5vbhfuy75O5MT8R9AD-7ItngeFI}
\end{figure}

Since then, the dream of many mathematicians consisted in making this idea work in practice, just using a bit of theory of partial differential equations. The gist is that light rays, and in general electromagnetic fields, are described by Maxwell's equations. These equations depend on the specific transmission medium via the permittivity and permeability of the substance under consideration. Therefore, to build an invisibility device according to Romulans, one should find, for instance, a medium with suitable permittivity and permeability forcing the light rays to bend around a ``hole'': in this way, a space vessel located inside such a hole would be invisible from the outside, see Figure~\ref{0okjnh78uigfi34bqtb76b67c3b2b376b5vbhfuy75O5MT8R9AD-7ItngeFI}.

Some practical questions still remain: how can one determine, at least at a mathematical level, the pointwise values of the permittivity and permeability of the medium in order to guide the lines around the desired hole? And can one do this for real?

To address these questions, let us start with the mathematical aspect of the problem (see the two adjacent
articles~\cite{OCMULA0ipj3} and~\cite{PHBDTP2wesd67uPKMRlm3} for the pioneering aspects of the theory, as well as for links with complex analysis and Riemannian geometry).
The calculations to implement these strategies rely on the mathematical property of Maxwell's equations to maintain their structure\footnote{ The ``form invariance'' of Maxwell's equations
is one of the main building blocks of the so-called
``transformation optics'', \index{transformation optics} which aims at using changing of coordinate system as a tool to reduce the wave propagation in inhomogeneous materials
to that in the standard space.}
under all spatial transformations. In a nutshell, this invariance property allows one, say, to pick a spatial \label{6CHR55} transformation, for example, which maps a ball into an annulus\footnote{{F}rom a technical point of view, such a transformation is singular, since the topologies of the initial and final domains are different, and this causes some mathematical difficulty on the notion of (weak) solution to consider, but here we surf over this detail (the interested reader may look at Section~3.5
of~\cite{MR2481110}).} and which coincides with the identity map outside the ball. Therefore, a solution of Maxwell's equations in the usual space would be mapped into a solution of Maxwell's equations in a space with inhomogeneous permittivity and permeability in which the light rays cannot enter the hole that has been produced inside the annulus. At that point, the mathematical task is completed and it ``suffices'' to build, at least approximately, a material with those specific values of permittivity and permeability to obtain a nice cloaking of any object placed inside the hole.

Let us see in some further detail how to perform these mathematical transformations and how to construct the desired material for the cloaking.

To start with, let us check the invariance of Maxwell's equations with respect to spatial transformations. To this end, consider Maxwell's equations in~$\R^3$ for the electric field~$E$ and the magnetic field~$B$, without sources or currents, with permittivity~$\e$ and permeability~$\mu$, namely, up to normalizations, \index{Maxwell's equations}
\begin{equation}\label{MAXLWDJOHK-oiwjfeknv}
\begin{dcases}
\div (\e E) = 0,\\
\curl E = -\mu\partial_t B ,\\ 
\div (\mu B) = 0,\\
\curl B = \e\partial_t E.
\end{dcases}\end{equation}
In this setting, the space variable will be denoted by~$x\in\R^3$, and~$\mu=\mu(x)$ and~$\e=\e(x)$ are matrix valued functions, to account for possible inhomogeneity and anisotropy of the medium.

A spatial change of coordinates~$\bar{x}=\bar{x}(x)$ maintains the structure in~\eqref{MAXLWDJOHK-oiwjfeknv}. Indeed, considering the Jacobian matrix
\begin{equation} \label{MAXLWDJOHK-oiwjfeknv-098765456789oiujhygfdcvbhtre32345t6y7u890oiujhgvcxsw3wsxdrtR}A_{ij}:=\frac{\partial\bar x_i}{\partial x_j},\end{equation}
with~$i$, $j\in\{1,2,3\}$, one can define the transformed quantities as
\begin{eqnarray*}&&
\bar{E}(\bar{x},t):= A^{-T}( x) \,E(x,t), \qquad
\bar{B}(\bar{x},t):= A^{-T}(x) \,B(x,t),\\&&
\bar\mu(\bar{x}):=\frac{1}{{\mathcal{J}}(x)}\,A(x)\,\mu(x)\, A^T(x)\qquad{\mbox{and}}\qquad
\bar\e(\bar{x}):=\frac{1}{{\mathcal{J}}(x)}\,A(x)\,\e(x)\, A^T(x),
\end{eqnarray*}
where~$x=x(\bar{x})$, $A^{-T}$ is the inverse of the transpose (which coincides with the transpose of the inverse) of~$A$, and~${\mathcal{J}}:=\det A$.

In this way, for each~$i$, $j\in\{1,2,3\}$,
\begin{eqnarray*}
\delta_{ij}=\partial_{\bar{x}_i} \bar{x}_j=
\sum_{k=1}^3 \partial_{x_k} \bar{x}_j\partial_{\bar{x}_i}x_k
=\sum_{k=1}^3 A_{jk} \partial_{\bar{x}_i}x_k,
\end{eqnarray*}
showing that the inverse of~$A$ is the matrix with entries
\begin{equation}\label{KLO023oi4krjt-okPht7GGubu-000}A^{-1}_{ij}=\partial_{\bar{x}_j}x_i.\end{equation}

We also calculate\footnote{Some of the tedious computations in these pages could be elegantly streamlined by the use of the Levi-Civita symbols for curls and determinants (see e.g. Chapter~9 in~\cite{MR3136419}
for an introduction of this notation and the appendix in~\cite{PhysRevE77036611}
for a compact proof of the form invariance of Maxwell's equations). We avoid here this quicker, but more sophisticated, approach in order to maintain the calculations in an elementary form.} that
\begin{equation}\label{KLO023oi4krjt-okPht7GGubu-00}
\partial_{x_k} A_{ij}=\partial_{x_k} \partial_{x_j}\bar x_i=\partial_{x_j}A_{ik}
,\end{equation}
that
\begin{equation}\label{KLO023oi4krjt-okPht7GGubu-01}
\e E=\Big({{\mathcal{J}}}{A^{-1}\bar\e A^{-T}}\Big)\Big(A^T\bar E\Big)={{\mathcal{J}}}{A^{-1}\bar\e \bar E}\end{equation}
and that
\begin{equation}\label{KLO023oi4krjt-okPht7GGubu-02} \partial_{x_j}=\sum_{i=1}^3 \partial_{x_j}\bar x_i\,\partial_{\bar x_i}=\sum_{i=1}^3 A_{ij}\partial_{\bar x_i}.\end{equation}

The latter identity also entails that
\begin{equation}\label{KLO023oi4krjt-okPht7GGubu-02-b} 
\sum_{j=1}^3A_{jk}^{-1}\partial_{x_j}=\sum_{i,j=1}^3 A_{jk}^{-1}A_{ij}\partial_{\bar x_i}=
\sum_{i=1}^3 \delta_{ik}\partial_{\bar x_i}=\partial_{\bar x_k}
.\end{equation}

We also recall that
$$ A_{ij}^{-1}=\frac{C_{ji}}{\mathcal{J}},$$
where the cofactor matrix has entries
$$ C_{\ell m}:=(-1)^{\ell+m}M_{\ell m},$$
being~$M_{\ell m}$ the minor of the entry in the~$\ell$th row and~$m$th column (that is, the determinant of the submatrix of~$A$ formed by deleting the~$\ell$th row and~$m$th column).

In this way,
\begin{equation}\label{KLO023oi4krjt-okPht7GGubu-03} {{\mathcal{J}}}A^{-1}_{jk}=C_{kj}=(-1)^{k+j}M_{kj}.
\end{equation}

We claim that, for all~$k\in\{1,2,3\}$,
\begin{equation}\label{KLO023oi4krjt-okPht7GGubu-04} \sum_{j=1}^3\partial_{x_j}\big({{\mathcal{J}}}A^{-1}_{jk}\big)=0.
\end{equation}
We choose~$k=1$ (the cases~$k=2$ and~$k=3$ being similar) and we employ~\eqref{KLO023oi4krjt-okPht7GGubu-00} and~\eqref{KLO023oi4krjt-okPht7GGubu-03} to see that
\begin{eqnarray*}
\sum_{j=1}^3 \partial_{x_j}\big({{\mathcal{J}}}A^{-1}_{j1}\big)&=&-\sum_{j=1}^3 (-1)^{j}\partial_{x_j}M_{1j}\\
&=& \partial_{x_1}M_{11}-\partial_{x_2}M_{12}+\partial_{x_3}M_{13}\\&=&
\partial_{x_1}(A_{22}A_{33}-A_{23}A_{32})-\partial_{x_2}(A_{21}A_{33}-A_{23}A_{31})+\partial_{x_3}(A_{21}A_{32}-A_{22}A_{31})\\&=&
\partial_{x_1}A_{22} A_{33}+A_{22}\partial_{x_1}A_{33}
-\partial_{x_1}A_{23}A_{32}-A_{23}\partial_{x_1}A_{32}\\&&\qquad
-\partial_{x_2}A_{21}A_{33}-A_{21}\partial_{x_2}A_{33}
+\partial_{x_2}A_{23}A_{31}+A_{23}\partial_{x_2}A_{31}\\&&\qquad
+\partial_{x_3}A_{21}A_{32}+A_{21}\partial_{x_3}A_{32}
-\partial_{x_3}A_{22}A_{31}-A_{22}\partial_{x_3}A_{31}\\&=&
\partial_{x_1}A_{22} A_{33}+A_{22}\partial_{x_1}A_{33}
-\partial_{x_1}A_{23}A_{32}-A_{23}\partial_{x_1}A_{32}\\&&\qquad
-\partial_{x_1}A_{22}A_{33}-A_{21}\partial_{x_2}A_{33}
+\partial_{x_2}A_{23}A_{31}+A_{23}\partial_{x_1}A_{32}\\&&\qquad
+\partial_{x_1}A_{23}A_{32}+A_{21}\partial_{x_2}A_{33}
-\partial_{x_2}A_{23}A_{31}-A_{22}\partial_{x_1}A_{33}\\&=&0.
\end{eqnarray*}
The claim in~\eqref{KLO023oi4krjt-okPht7GGubu-04} is thereby established.

As a result, using~\eqref{KLO023oi4krjt-okPht7GGubu-01},
\eqref{KLO023oi4krjt-okPht7GGubu-02}, \eqref{KLO023oi4krjt-okPht7GGubu-02-b} and~\eqref{KLO023oi4krjt-okPht7GGubu-04},
\begin{equation*}
\begin{split}
& 0=\div_x (\e E)=\div_x\left({{\mathcal{J}}}{A^{-1}\bar\e \bar E}\right)=
\sum_{j=1}^3 \partial_{x_j} \left({{\mathcal{J}}}{A^{-1}\bar\e \bar E}\right)_j=
\sum_{j,k=1}^3 \partial_{x_j} \left({{\mathcal{J}}}{A^{-1}_{jk}(\bar\e \bar E)_k}\right)\\&\qquad=
\sum_{j,k=1}^3 \partial_{x_j}\big({\mathcal{J}} A^{-1}_{jk}\big)(\bar\e \bar E)_k
+\sum_{j,k=1}^3 {\mathcal{J}} A^{-1}_{jk} \partial_{x_j} \big(\bar\e \bar E\big)_k=0
+\sum_{k=1}^3 {\mathcal{J}} \partial_{\bar x_k} \big(\bar\e \bar E\big)_k\\&\qquad={\mathcal{J}}\div_{\bar x}(\bar\e \bar E\big).
\end{split}\end{equation*}

Accordingly,
\begin{equation}\label{KLO023oi4krjt-okPht7GGubu-02-100}
\div_{\bar x}(\bar\e \bar E\big)=0
\end{equation}
and similarly
\begin{equation}\label{KLO023oi4krjt-okPht7GGubu-02-101}
\div_{\bar x}(\bar\mu \bar B\big)=0.
\end{equation}

Furthermore, for all~$i$, $j\in\{1,2,3\}$ it follows from~\eqref{KLO023oi4krjt-okPht7GGubu-02-b} that
\begin{eqnarray*}&& \partial_{\bar{x}_i}\bar E_j=
\partial_{\bar{x}_i}(A^{-T}E)_j=\sum_{m=1}^3 A_{mi}^{-1}\partial_{x_m}(A^{-T}E)_j
=\sum_{m,\ell=1}^3 A_{mi}^{-1}\,\partial_{x_m}(A^{-1}_{\ell j}E_\ell)\\
&&\qquad=\sum_{m,\ell=1}^3 A_{mi}^{-1}\,\partial_{x_m}A^{-1}_{\ell j}\,E_\ell+\sum_{m,\ell=1}^3 A_{mi}^{-1}\,A^{-1}_{\ell j}\,\partial_{x_m}E_\ell=\sum_{\ell=1}^3 \partial_{\bar x_i}A^{-1}_{\ell j}\,E_\ell+\sum_{m,\ell=1}^3 A_{mi}^{-1}\,A^{-1}_{\ell j}\,\partial_{x_m}E_\ell
\end{eqnarray*}
and therefore
\begin{equation}\label{KLO023oi4krjt-okPht7GGubu-02-108}\begin{split}&
\partial_{\bar{x}_i}\bar E_j-\partial_{\bar{x}_j}\bar E_i =
\sum_{\ell=1}^3 \Big(\partial_{\bar x_i}A^{-1}_{\ell j}-\partial_{\bar x_j}A^{-1}_{\ell i}\Big)E_\ell+\sum_{m,\ell=1}^3
\Big(A_{mi}^{-1}\,A^{-1}_{\ell j}-A_{mj}^{-1}\,A^{-1}_{\ell i}\Big)\partial_{x_m}E_\ell. \end{split}
\end{equation}

In addition, by~\eqref{KLO023oi4krjt-okPht7GGubu-000},
$$ \partial_{\bar x_i}A^{-1}_{\ell j}=\partial_{\bar x_i}\partial_{\bar{x}_j}x_\ell=\partial_{\bar x_j}A^{-1}_{\ell i}
,$$ which reduces~\eqref{KLO023oi4krjt-okPht7GGubu-02-108} to
\begin{equation}\label{KLO023oi4krjt-okPht7GGubu-02-108F8ujh-09G}\begin{split}
\partial_{\bar{x}_i}\bar E_j-\partial_{\bar{x}_j}\bar E_i &=0+\sum_{m,\ell=1}^3
\Big(A_{mi}^{-1}\,A^{-1}_{\ell j}-A_{mj}^{-1}\,A^{-1}_{\ell i}\Big)\partial_{x_m}E_\ell
\\&=\sum_{{1\le m,\ell\le3}\atop{m\ne \ell}}
\Big(A_{mi}^{-1}\,A^{-1}_{\ell j}-A_{mj}^{-1}\,A^{-1}_{\ell i}\Big)\partial_{x_m}E_\ell\\&=\sum_{{1\le m,\ell\le3}\atop{m\ne \ell}}A_{mi}^{-1}\,A^{-1}_{\ell j} \,\partial_{x_m}E_\ell
-\sum_{{1\le m,\ell\le3}\atop{m\ne \ell}} A_{\ell j}^{-1}\,A^{-1}_{m i}\,\partial_{x_\ell}E_m\\&=
\sum_{{1\le m,\ell\le3}\atop{m\ne \ell}}A_{mi}^{-1}\,A^{-1}_{\ell j} \Big(\partial_{x_m}E_\ell-\partial_{x_\ell}E_m\Big)\\&=
\sum_{{1\le m,\ell\le3}\atop{m< \ell}}A_{mi}^{-1}\,A^{-1}_{\ell j} \Big(\partial_{x_m}E_\ell-\partial_{x_\ell}E_m\Big)+
\sum_{{1\le m,\ell\le3}\atop{m> \ell}}A_{mi}^{-1}\,A^{-1}_{\ell j} \Big(\partial_{x_m}E_\ell-\partial_{x_\ell}E_m\Big)\\&=
\sum_{{1\le m,\ell\le3}\atop{m< \ell}}A_{mi}^{-1}\,A^{-1}_{\ell j} \Big(\partial_{x_m}E_\ell-\partial_{x_\ell}E_m\Big)-
\sum_{{1\le m,\ell\le3}\atop{m< \ell}}A_{\ell i}^{-1}\,A^{-1}_{m j} \Big(\partial_{x_m}E_\ell-\partial_{x_\ell}E_m\Big)\\&=
\sum_{{1\le m,\ell\le3}\atop{m< \ell}}\Big(A_{mi}^{-1}\,A^{-1}_{\ell j} -A_{\ell i}^{-1}\,A^{-1}_{m j}\Big)\Big(\partial_{x_m}E_\ell-\partial_{x_\ell}E_m\Big).\end{split}
\end{equation}

Now, given~$q\in\{1,2,3\}$ we let~$\alpha_q$, $\beta_q\in\{1,2,3\}$ be such\footnote{More explicitly,
$$ (\alpha_1,\beta_1):=(2,3),\qquad (\alpha_2,\beta_2):=(1,3)\qquad{\mbox{and}}\qquad
(\alpha_3,\beta_3):=(1,2).$$}
that~$\alpha_q<\beta_q$
and~$\{q,\alpha_q,\beta_q\}=\{1,2,3\}$.

Let also~$R$ and~$\bar R$ stand for a short notation of~$\curl_x E$ and~$\curl_{\bar x}\bar E$, respectively.

In this way,
$$ R_q=(-1)^{q+1} \Big(\partial_{x_{\alpha_q}}E_{\beta_q}-\partial_{x_{\beta_q}}E_{\alpha_q}\Big)$$
and, for all~$p\in\{1,2,3\}$, we can rewrite~\eqref{KLO023oi4krjt-okPht7GGubu-02-108F8ujh-09G} in the form
\begin{equation*}
\begin{split}
\bar R_{p}&=(-1)^{p+1}\Big(\partial_{\bar{x}_{\alpha_p}}\bar E_{\beta_p}-\partial_{\bar{x}_{\beta_p}}\bar E_{\alpha_p} \Big)\\&=(-1)^{p+1}
\sum_{{1\le m,\ell\le3}\atop{m< \ell}}\Big(A_{m \alpha_p}^{-1}\,A^{-1}_{\ell \beta_p} -A_{\ell \alpha_p}^{-1}\,A^{-1}_{m \beta_p}\Big)\Big(\partial_{x_m}E_\ell-\partial_{x_\ell}E_m\Big)\\&=(-1)^{p+1}
\sum_{q=1}^3\Big(A_{\alpha_q \alpha_p}^{-1}\,A^{-1}_{{\beta_q} \beta_p} -A_{{\beta_q} \alpha_p}^{-1}\,A^{-1}_{\alpha_q \beta_p}\Big)\Big(\partial_{x_{\alpha_q}}E_{\beta_q}-\partial_{x_{\beta_q}}E_{\alpha_q}\Big)\\&=
\sum_{q=1}^3(-1)^{p+q}\Big(A_{\alpha_q \alpha_p}^{-1}\,A^{-1}_{{\beta_q} \beta_p} -A_{{\beta_q} \alpha_p}^{-1}\,A^{-1}_{\alpha_q \beta_p}\Big)R_q.
\end{split}
\end{equation*}

On this account, for all~$d\in\{1,2,3\}$,
\begin{equation}\label{KLO023oi4krjt-okPht7GGubu-02-108F8ujh-09G-1}
\begin{split}
(A^{-1}\bar R)_d=\sum_{p=1}^3 A^{-1}_{dp}\bar R_p=
\sum_{p,q=1}^3(-1)^{p+q}\Big(A^{-1}_{dp}\,A_{\alpha_q \alpha_p}^{-1}\,A^{-1}_{{\beta_q} \beta_p} -A^{-1}_{dp}\,A_{{\beta_q} \alpha_p}^{-1}\,A^{-1}_{\alpha_q \beta_p}\Big)R_q.
\end{split}
\end{equation}

We now claim that if~$q\ne d$ then
\begin{equation}\label{KLO023oi4krjt-okPht7GGubu-02-108F8ujh-09G-2}
\sum_{p=1}^3(-1)^{p}\Big(A^{-1}_{dp}\,A_{\alpha_q \alpha_p}^{-1}\,A^{-1}_{{\beta_q} \beta_p} -A^{-1}_{dp}\,A_{{\beta_q} \alpha_p}^{-1}\,A^{-1}_{\alpha_q \beta_p}\Big)=0.
\end{equation}
Indeed, if~$q\ne d$ then~$d\in\{\alpha_q,\beta_q\}$. Let us assume that~$d=\alpha_q$ (the case~$d=\beta_q$ being similar). In this situation,
\begin{eqnarray*}
&&\sum_{p=1}^3(-1)^{p}\Big(A^{-1}_{dp}\,A_{\alpha_q \alpha_p}^{-1}\,A^{-1}_{{\beta_q} \beta_p} -A^{-1}_{dp}\,A_{{\beta_q} \alpha_p}^{-1}\,A^{-1}_{\alpha_q \beta_p}\Big)\\&=&
\sum_{p=1}^3(-1)^{p}\Big(A^{-1}_{\alpha_q p}\,A_{\alpha_q \alpha_p}^{-1}\,A^{-1}_{{\beta_q} \beta_p} -A^{-1}_{\alpha_q p}\,A_{{\beta_q} \alpha_p}^{-1}\,A^{-1}_{\alpha_q \beta_p}\Big)\\&=&
-\Big(A^{-1}_{\alpha_q 1}\,A_{\alpha_q 2}^{-1}\,A^{-1}_{{\beta_q} 3} -A^{-1}_{\alpha_q 1}\,A_{{\beta_q} 2}^{-1}\,A^{-1}_{\alpha_q 3}\Big)
+\Big(A^{-1}_{\alpha_q 2}\,A_{\alpha_q 1}^{-1}\,A^{-1}_{{\beta_q}3} -A^{-1}_{\alpha_q 2}\,A_{{\beta_q}1}^{-1}\,A^{-1}_{\alpha_q 3}\Big)\\&&\qquad
-\Big(A^{-1}_{\alpha_q 3}\,A_{\alpha_q1}^{-1}\,A^{-1}_{{\beta_q}2} -A^{-1}_{\alpha_q 3}\,A_{{\beta_q} 1}^{-1}\,A^{-1}_{\alpha_q2}\Big),
\end{eqnarray*}
which vanishes (canceling the first term with the third, the second with the fifth, and the fourth with the sixth).
This calculation proves~\eqref{KLO023oi4krjt-okPht7GGubu-02-108F8ujh-09G-2}.

Hence, plugging~\eqref{KLO023oi4krjt-okPht7GGubu-02-108F8ujh-09G-2} into~\eqref{KLO023oi4krjt-okPht7GGubu-02-108F8ujh-09G-1}, we arrive at
\begin{equation}\label{DETEKLO023oi4krjt-okPht7GGubu-02-108F8ujh-09G-2D}
\begin{split}
(A^{-1}\bar R)_d=
\sum_{p=1}^3(-1)^{p+d}\Big(A^{-1}_{dp}\,A_{\alpha_d \alpha_p}^{-1}\,A^{-1}_{{\beta_d} \beta_p} -A^{-1}_{dp}\,A_{{\beta_d} \alpha_p}^{-1}\,A^{-1}_{\alpha_d \beta_p}\Big)R_d.
\end{split}
\end{equation}

Now we claim that for every matrix~$M\in{\rm Mat}(3\times3)$ and every~$d\in\{1,2,3\}$
\begin{equation}\label{DETEKLO023oi4krjt-okPht7GGubu-02-108F8ujh-09G-2}
\sum_{p=1}^3(-1)^{p+d}\Big(M_{dp}\,M_{\alpha_d \alpha_p}\,M_{{\beta_d} \beta_p} -M_{dp}\,M_{{\beta_d} \alpha_p}\,M_{\alpha_d \beta_p}\Big)=\det M.\end{equation}
To check this identity one can rely on the ``axiomatic approach'' to determinant functions (see e.g.~\cite[Theorem 8.14]{SCHAUMM}), stating that the determinant is the unique multilinear and alternating function on matrices which returns~$1$ when applied to the identity. More explicitly, to prove~\eqref{DETEKLO023oi4krjt-okPht7GGubu-02-108F8ujh-09G-2} it suffices to check three properties: namely, if we call~$M_i$ the $i$th row of a matrix~$M$ and~$
{\mathcal{D}}(M_1,M_2,M_3)$ the function in the left hand side of~\eqref{DETEKLO023oi4krjt-okPht7GGubu-02-108F8ujh-09G-2} (thinking about it as a function acting on rows), our aim is to check that
\begin{eqnarray}
\nonumber &&{\mbox{if~$i\in\{1,2,3\}$ and~$M_i=av+bw$ for some~$a$, $b\in\R$ and~$v$, $w\in\R^3$}},\\&&\qquad\qquad{\mbox{then~${\mathcal{D}}(\dots,M_i,\dots)=a{\mathcal{D}}(\dots,v,\dots)+b{\mathcal{D}}(\dots,w\dots)$,}}\label{DETEKLO023oi4krjt-okPht7GGubu-02-108F8ujh-09G-2i} \\
\label{DETEKLO023oi4krjt-okPht7GGubu-02-108F8ujh-09G-2ii} && {\mbox{if~$M_i=M_j$ for some~$i\ne j\in\{1,2,3\}$, then~${\mathcal{D}}(M_1,M_2,M_3)=0$}}\\
\label{DETEKLO023oi4krjt-okPht7GGubu-02-108F8ujh-09G-2iii} {\mbox{and }}&& {\mathcal{D}}(e_1,e_2,e_3)=1.
\end{eqnarray}
To check~\eqref{DETEKLO023oi4krjt-okPht7GGubu-02-108F8ujh-09G-2i}, let us suppose that~$i=d$
(the cases in which~$i=\alpha_d$ and~$i=\beta_d$ being completely analogous). In this situation, we have that
\begin{eqnarray*}
&&M_{dp}\,M_{\alpha_d \alpha_p}\,M_{{\beta_d} \beta_p} -M_{dp}\,M_{{\beta_d} \alpha_p}\,M_{\alpha_d \beta_p}\\
&&\qquad=(av+bw)_p\,M_{\alpha_d \alpha_p}\,M_{{\beta_d} \beta_p} -(av+bw)_p\,M_{{\beta_d} \alpha_p}\,M_{\alpha_d \beta_p}\\
&&\qquad=av_p\,M_{\alpha_d \alpha_p}\,M_{{\beta_d} \beta_p} 
+bw_p\,M_{\alpha_d \alpha_p}\,M_{{\beta_d} \beta_p} 
\\&&\qquad\qquad -av_p\,M_{{\beta_d} \alpha_p}\,M_{\alpha_d \beta_p}
-bw_p\,M_{{\beta_d} \alpha_p}\,M_{\alpha_d \beta_p},
\end{eqnarray*}
from which~\eqref{DETEKLO023oi4krjt-okPht7GGubu-02-108F8ujh-09G-2i} plainly follows.

To establish~\eqref{DETEKLO023oi4krjt-okPht7GGubu-02-108F8ujh-09G-2ii} one should consider three cases,
$M_d=M_{\alpha_d}$, $M_d=M_{\beta_d}$ and~$M_{\alpha_d}=M_{\beta_d}$.
The first and the second case however are similar, so we focus on the proof of the second and third cases.

Thus, let us suppose that~$M_d=M_{\beta_d}$. Then,
\begin{eqnarray*}
(-1)^d{\mathcal{D}}(M_1,M_2,M_3)&=&
\sum_{p=1}^3(-1)^{p}\Big(M_{dp}\,M_{\alpha_d \alpha_p}\,M_{d \beta_p} -M_{dp}\,M_{d \alpha_p}\,M_{\alpha_d \beta_p}\Big)\\&=&
-\Big(M_{d1}\,M_{\alpha_d 2}\,M_{d 3} -M_{d1}\,M_{d 2}\,M_{\alpha_d3}\Big)
\\&&\qquad
+\Big(M_{d2}\,M_{\alpha_d1}\,M_{d 3} -M_{d2}\,M_{d1}\,M_{\alpha_d 3}\Big)
\\&&\qquad-\Big(M_{d3}\,M_{\alpha_d1}\,M_{d 2} -M_{d3}\,M_{d 1}\,M_{\alpha_d 2}\Big)
\end{eqnarray*}
which vanishes after simplifying the first term with the sixth, the second with the fourth, and the third with the fifth.

This establishes~\eqref{DETEKLO023oi4krjt-okPht7GGubu-02-108F8ujh-09G-2ii} in this case and thus we take now into consideration the case~$M_{\alpha_d}=M_{\beta_d}$. In this situation,
$$ {\mathcal{D}}(M_1,M_2,M_3)=
\sum_{p=1}^3(-1)^{p+d}\Big(M_{dp}\,M_{\alpha_d \alpha_p}\,M_{{\alpha_d} \beta_p} -M_{dp}\,M_{{\alpha_d} \alpha_p}\,M_{\alpha_d \beta_p}\Big),$$
which clearly vanishes and the proof of~\eqref{DETEKLO023oi4krjt-okPht7GGubu-02-108F8ujh-09G-2ii} is thereby complete.

Moreover,
\begin{eqnarray*}
{\mathcal{D}}(e_1,e_2,e_3)&=&
\sum_{p=1}^3(-1)^{p+d}\Big(\delta_{dp}\,\delta_{\alpha_d \alpha_p}\,\delta_{{\beta_d} \beta_p} -\delta_{dp}\,\delta_{{\beta_d} \alpha_p}\,\delta_{\alpha_d \beta_p}\Big)\\&=&(-1)^{d+d}\Big(\delta_{\alpha_d \alpha_d}\,\delta_{{\beta_d} \beta_d} -\delta_{{\beta_d} \alpha_d}\,\delta_{\alpha_d \beta_d}\Big)\\&=&1-0,
\end{eqnarray*}
giving~\eqref{DETEKLO023oi4krjt-okPht7GGubu-02-108F8ujh-09G-2iii}.

This completes the proof of~\eqref{DETEKLO023oi4krjt-okPht7GGubu-02-108F8ujh-09G-2}.

Thus, we plug~\eqref{DETEKLO023oi4krjt-okPht7GGubu-02-108F8ujh-09G-2} into~\eqref{DETEKLO023oi4krjt-okPht7GGubu-02-108F8ujh-09G-2D} and we obtain that
\begin{eqnarray*}
(A^{-1}\bar R)_d=
\det(A^{-1})R_d=\frac{R_d}{\mathcal{J}},
\end{eqnarray*}
that is~$A^{-1}\bar R=\frac{R}{\mathcal{J}}$.

As a result,
\begin{eqnarray*}
\curl_{\bar x}\bar E=\bar R=\frac{AR}{\mathcal{J}}=\frac{A\curl_xE}{\mathcal{J}}=-\frac{A\mu \partial_t B}{\mathcal{J}}=-\frac{A({\mathcal{J}}A^{-1}\bar\mu A^{-T}) \partial_t (A^T\bar B)}{\mathcal{J}}=-\bar\mu\partial_t\bar B\end{eqnarray*}
and similarly
$$\curl_{\bar x} \bar B-\bar\e\partial_t \bar E=0 .$$

These observations, together with~\eqref{KLO023oi4krjt-okPht7GGubu-02-100} and~\eqref{KLO023oi4krjt-okPht7GGubu-02-101}, establish the form invariance of Maxwell's equations, i.e. the property
that Maxwell's equations~\eqref{MAXLWDJOHK-oiwjfeknv} maintain their structure under every spatial transformation.\medskip

\begin{figure}
  \centering
  \includegraphics[width=.67\linewidth]{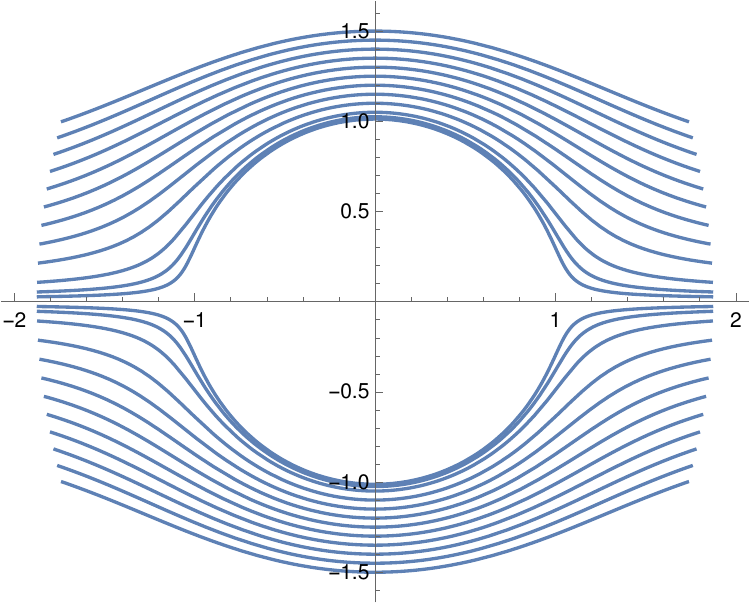}
 \caption{\sl Steering light around a hole, as analytically described in~\eqref{MAXLWDJOHK-oiwjfeknv-098765456789oiujhygfdcvbhtre32345t6y7u890oiujhgvcxsw3wsxdrtR2} (with parameters~$h\in
\{\pm0.025,\pm0.05,\pm0.1,\pm0.2,\pm0.3,\pm0.4,\pm0.5,\pm0.6,\pm0.7,\pm0.8,\pm0.9,\pm1\}$).
To be compared with Figure~\ref{0okjnh78uigfi34bqtb76b67c3b2b376b5vbhfuy75O5MT8R9AD-7ItngeFI}.}\label{SIRCOMTRADItangeFI-BENL09}
\end{figure}

As discussed on page~\pageref{6CHR55}, this feature of Maxwell's equations of
``looking the same'' in every coordinate system
can be conveniently interpreted, rather than just as the representation of the same electromagnetic properties of a material written in different coordinates, as the representation of different electromagnetic properties in the usual space. 

Namely,
we consider a spatial transformation to map a ball into an annulus, i.e. a ring-shaped object with a hole in the middle.
The form invariance of Maxwell's equations entails that light cannot pass through the hole (hence objects
can be made invisible by placing them inside the hole).
Since the transformation is designed to coincide with the identity map outside the ball, the light rays will behave ``as usual'' far from the invisible object:
the transformed space corresponds to a medium with inhomogeneous permittivity and permeability, and we will discuss below how to construct in practice a material with the appropriate structure to
respond to electric and magnetic fields in such a prescribed way.\medskip

To map a ball, say~$B_2$ (or better to say, a ball minus a point, say~$B_2\setminus\{0\}$), into an annulus, say~$B_2\setminus B_1$, one can consider the transformation
$$ \bar{x}:=\begin{dcases}\displaystyle\frac{|x|+2}{2|x|}\,x &{\mbox{ if }}x\in B_2\setminus \{0\},\\
x&{\mbox{ if }}x\in\R^3\setminus B_2.\end{dcases}$$
Since~$|\bar x|=\frac{|x|+2}{2}\in[1,2]$ for all~$x\in B_2\setminus \{0\}$ we have that~$\bar x(B_2\setminus \{0\})=B_2\setminus B_1$. Also, $\bar x=x$ outside~$B_2$, hence this change of coordinates possesses the features described above.

In this case, the Jacobian matrix in~\eqref{MAXLWDJOHK-oiwjfeknv-098765456789oiujhygfdcvbhtre32345t6y7u890oiujhgvcxsw3wsxdrtR} in the relevant region~$B_2\setminus\{0\}$ becomes
$$ A_{ij}=\partial_{x_j} \left(\frac{|x|+2}{2|x|}\,x_i\right)=\frac{|x|+2}{2|x|}\,\delta_{ij}-\frac{x_i\,x_j}{|x|^3}.$$
%%%That is,
%%%$$ A=\left(\begin{matrix}
%%%\frac{|x|+2}{2|x|}-\frac{x_1^2}{|x|^3} & -\frac{x_1\,x_2}{|x|^3} & -\frac{x_1\,x_3}{|x|^3}\cr
%%%-\frac{x_1\,x_2}{|x|^3} & \frac{|x|+2}{2|x|}-\frac{x_2^2}{|x|^3} & -\frac{x_2\,x_3}{|x|^3}\cr
%%%-\frac{x_1\,x_3}{|x|^3} & -\frac{x_2\,x_3}{|x|^3} & \frac{|x|+2}{2|x|}-\frac{x_3^2}{|x|^3}
%%%\end{matrix}\right).$$
Therefore, in this situation, the matrix~$A$ is symmetric,
$$ {\mathcal{J}}=\det A=\frac{(|x| + 2)^2 |x|^3 }{8 |x|^5}$$
and accordingly the new permittivity~$\bar\e$ and permeability~$\bar\mu$ obtained from the homogeneous case~$\e=\mu=1$ are equal to
$$ \frac{1}{{\mathcal{J}}}\,A\, A^T=\frac{1}{{\mathcal{J}}}\,A^2,$$
with~$A$ as above.

To appreciate the bending of a light ray induced by this transformation, for a given parameter~$h\in(-1,1)$,
we can consider the parallel segments
$$\ell_h:=\left\{ (t,h,0), {\mbox{ with }} |t|<\sqrt{3}\right\}.$$
Then,
\begin{equation}\label{MAXLWDJOHK-oiwjfeknv-098765456789oiujhygfdcvbhtre32345t6y7u890oiujhgvcxsw3wsxdrtR2} \bar{x}(\ell_h)=\frac{\sqrt{t^2+h^2}+2}{2\sqrt{t^2+h^2}}\,(t,h,0),\end{equation}
see Figure~\ref{SIRCOMTRADItangeFI-BENL09} for a visual representation of such a hole which cannot be entered by light rays.\medskip

Now let us briefly discuss how, in practice, there are potentials to build efficient electromagnetic cloaking devices. The idea is that one can use special media, called ``metamaterials'', \index{metamaterials}
whose physical parameters can be accurately designed to steer light around a hidden region of space.
These media are
artificial composite, usually assembled from multiple elements, such as metals and special plastics, often arranged in repeating patterns and presenting small structures (such as tiny needles).
The metamaterial design can get quite sophisticated, involving thousands of elements with complex geometry,
see e.g.~\cite{u8htbgnSC67uj11iE4567890paleijert73uYHityb58c9e-09, WOOD2009379}.

Interestingly, not only mathematics is pivotal in the development of metamaterials, but also metamaterials
have been recently used to
perform mathematical operations, such as differentiation, integration or convolution, so we can expect a nice interplay between mathematics and technology in this field, see~\cite{MR3156078}.\medskip

\begin{figure}
  \centering
  \includegraphics[width=.3\linewidth]{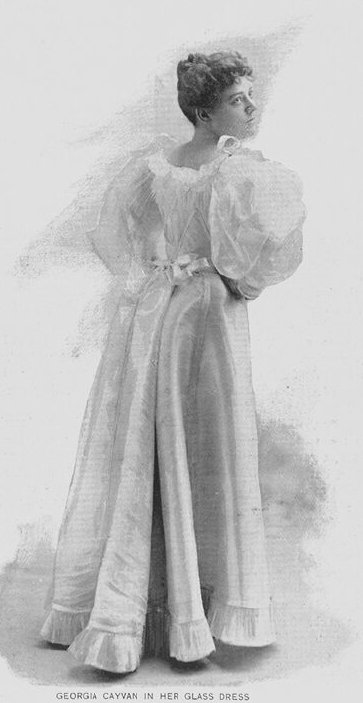}
 \caption{\sl The famous actress Georgia Cayvan posing in her glass dress
 (Public Domain image from
 Wikipedia).}\label{1NIRglassD-12435eFI}
\end{figure}

As a matter of curiosity, let us mention that the first patent for artificial materials of this type probably dates back to~1880, when Prussian-American inventor Hermann Hammesfahr registered a new material, known as
fiberglass. The outmost application of this material at that time was the design of a new dress that was very successfully presented at the 1893 World's Exposition in Chicago (the electric light bulb and the Ferris wheel
also debuted at the same fair, but allegedly the ``glass dress'' resulted far more sensational, see Figure~\ref{1NIRglassD-12435eFI}).
\medskip

However, going back to the topic of invisibility and cloaking devices,
at this point a natural remark arises. If the light bends around a shell,
any object enclosed within this region must be necessarily blind, since it cannot communicate electromagnetically\footnote{Interestingly, for the same reason, cloaked objects have no shadow.
All in all, not only would observers be incapable of seeing a cloaked object, they would also be unaware that something has been hidden!} with the outside world. That is, the observation made by Yakov Perelman
about Wells' Invisible Man (recall the discussion on page~\pageref{Yakov Perelman}) also applies in the case of
electromagnetic cloaking (though under different physical circumstances).
Several scientists have tried to get around this somewhat unavoidable issue, allowing some kind of information exchange between the inside and outside of the cloak, possibly at the cost of losing perfect invisibility, or by
placing the invisible outside the cloak utilizing so-called ``complementary media'' which can optically ``cancel'' a certain region of space at a certain
frequency to create an ``antiobject''
see e.g.~\cite{OPENCLO, PhysRevLett102093901, 184195-12p4e4ow3ejfhkgrb92eyrihf} 
(one could also try to rely on the use of different frequencies, e.g. cloaking for visible lights and using
infrared vision).
\medskip

It is also worth retaking inspiration from the Star Trek series mentioned on page~\pageref{1s2ddc23425qwedfgboprenc0813ti0a4g5gda}. Probably to make the plot more intriguing,
in the Star Trek saga while a space vessel is using its cloaking device, it is not allowed to fire weapons, use defensive shields, or operate transporters. To do any of these things, the vessel must ``decloak''.
Quite interestingly, this has a counterpart in the mathematical theory of cloaking.
Namely, what we described here is the cloaking of ``passive'' objects, i.e. objects without
internal currents: in this case, cloaking for Maxwell's equations is indeed mathematically possible.
Instead, if an object is ``active'' (e.g., it presents internal currents, as supposedly should occur for
a space vessel firing weapons, employing defensive shields, or operating transporters), even the mathematical theory of cloaking poses additional restrictions and
the complete invisibility seems more problematic\footnote{{F}rom the mathematical point of view,
this has to deal with the fact that ``finite energy'' solutions to Maxwell's equations do not
exist in the presence of active sources inside the cloaked
region.}
and a ``double coating'' may be necessary, see~\cite{MR2336363}, and one may need to
line the inner surface of the cloak with a perfectly reflecting material,
or induce extraordinary electric and magnetic surface voltages
at the inner surface of the cloak, see~\cite{PhysRevLett100063904}.

Roughly speaking,
active electromagnetic objects are more difficult to cloak
because they are constantly emitting or interacting with electromagnetic waves
(similar difficulties arise when cloaking moving objects);
to circumvent issues of this type one may also consider approximate cloaks, see e.g.~\cite{MR3192429}.
\medskip

Moreover, we have discussed here some basic features of ``passive'' ways of cloaking an object, but there are also ``active'' cloaking methods, corresponding to the recent technology of noise-canceling headphones: roughly speaking, one could try to ``cancel'' signals by actively producing ``opposite'' signals
(see~\cite{NELSON, MR2349998} for the case of acoustic waves and~\cite{PhysRevX3041011} for electromagnetic cloaking). A drawback of these active cloaking devices is that they may be subject to field instabilities and the design of the cloak also become even more difficult and complicated,
see e.g.~\cite{ALITALO200922}, thus posing severe challenges in practical applications.
\medskip

We stress that the use of the technical word ``active'' in the context of cloaking could be sometimes confusing,
since this may refer to either the object to be cloaked, in which some electromagnetic activity may take place,
or to the cloaking device,  which may exploit active sources specifically designed
to cancel suitable signals (or both, see~\cite[Remark~2]{MR2811880}).
\medskip

Let us also mention that here we considered the cloaking problem in a rather broad generality: additional assumptions on the object to be cloaked or on the angle of incidence of the light rays may allow more effective cloaking strategies that are specifically tailored to a specific situation, see e.g.~\cite{u8htbgnSC67uj11iE4567890paleijert73uYHityb58c9e-08,
u8htbgnSC67uj11iE4567890paleijert73uYHityb58c9e-09, u8htbgnSC67uj11iE4567890paleijert73uYHityb58c9e-12, u8htbgnSC67uj11iE4567890paleijert73uYHityb58c9e-13}.
\medskip

Notice that the ideas related to invisibility and cloaking, though at the moment
still in their early stages of development,
can have fruitful applications also beyond their original goals (and hopefully beyond the obvious military scope). For instance, if one improved the practical possibility of manipulating light rays, this would allow the creation of powerful microscopes and telescopes and maybe even faster computers.

Certainly, the mathematical difficulties for it are enormous, and so are the technological hurdles to be overcome. For instance, to effectively obtain invisibility, a cloak should efficiently interact at once with all frequencies of light rays, or at least the ones in the visible spectrum, while the parameters in the equation seem to be sensible to the change of frequency due to dispersive phenomena. The steps towards a possible improvement in this context may include the combination of metamaterials and dielectric materials, see e.g.~\cite{VALENTINE}.

On the one hand, at the moment, metamaterials seem to be employed fruitfully mostly at microwave\footnote{What happens is that typically, for cloaking devices,
one has to engineer the metamaterial at a scale smaller than the wavelength of radiation. 
For the optical radiation of visible light, the wavelength is less than a micron, therefore the structure of a metamaterial effective in this situation has to rely on nanotechnology.

For this reason, cloaking devices often work well for a single wavelength, or for a small wavelength range.

On a positive, or maybe negative, note, let us mention that
for radars the wavelength is about 3 cm, so one can relatively easily design metamaterials at the scale of millimeters to hide objects from radars, with military-specific projects.

By the way, possible military uses of cloaking devices are not limited to making vehicles and weapons invisible. For example, cloaking may be employed to deactivate
structures that interrupt signals, since one could cloak these structures to allow signals to pass by freely (in a sense, this leverages the blindness of the cloaked object
turning it into an advantage).

Also, steering light rays may make objects bigger or smaller (or different).} frequencies, see~\cite{0ojkhniatip1342a2n421d78536},
with additional difficulties to assemble and design them above this threshold; on the other hand, metamaterials engineering at optical frequencies has been also tested, see~\cite{oknGABS9oikmIJSODKJOLuiO0orlgP0002}.

Another possible hindrance to the development of cloaking devices may be the cost of the materials employed. A promising line of research tries to build cloaks out of ordinary materials, such as calcite, see e.g.~\cite{0987890-poijhgfvu87ygfetyhgfdsx20987654567890987}.\medskip

The cloaking problem for Maxwell's equations is also related with the electrostatic problem
presented in Section~\ref{Whenyouareamathematician}: indeed, a natural analogue in that case
is to construct examples of objects that cannot be detected by electrostatic measurements performed along the 
boundary (for examples, electrostatically cloaked complex objects that look the same
under electrical impedance tomography as an isotropic material).
See~\cite{MR2481110} for additional information on this type of problems, as well as for related
questions about electromagnetic wormholes.

\subsection{Fighting a pandemic using differential equations}\label{EPISECT:d}

\begin{figure}
  \centering
  \includegraphics[width=.79\linewidth]{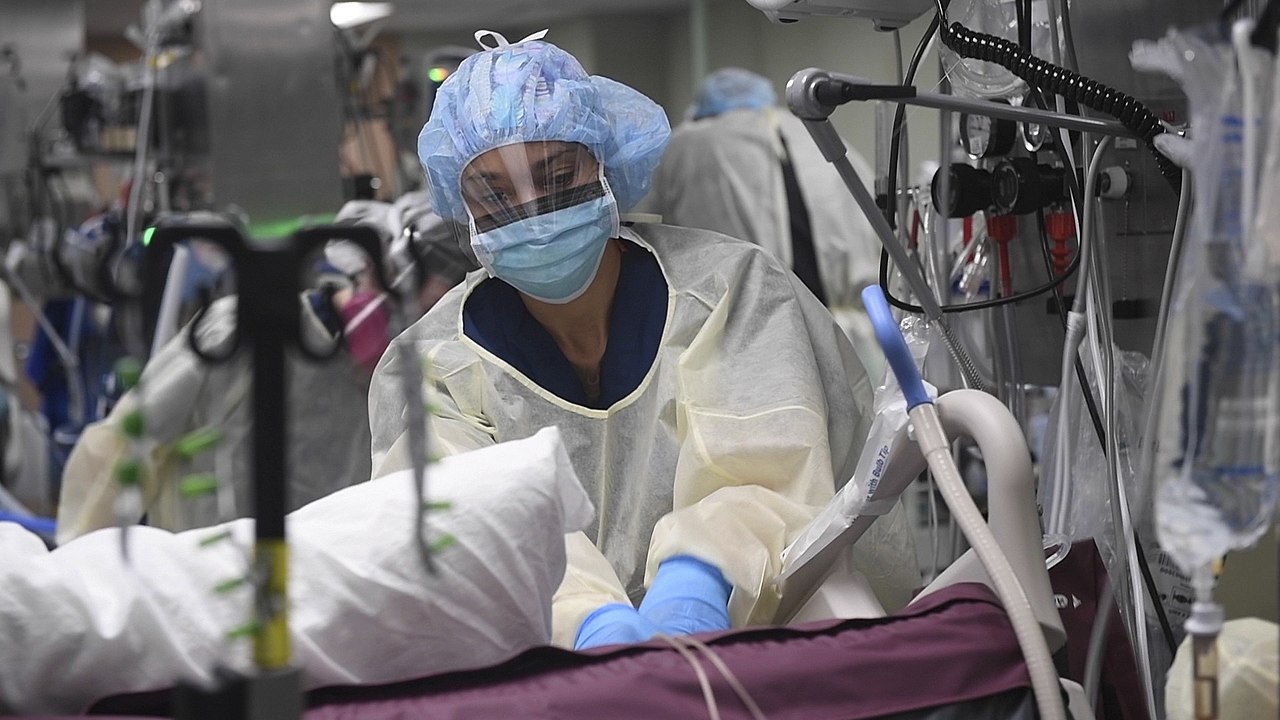}
 \caption{\sl Natasha McClinton, a surgical nurse, prepares a patient for a procedure in a COVID-19 intensive care unit
 (Public Domain image from
 Wikipedia).}\label{NIRDGREENFIDItangeFI}
\end{figure}

A topical subject nowadays consists in the spread of a virus among a given population, see Figure~\ref{NIRDGREENFIDItangeFI}.
Mathematics has established a solid reputation in trying to understand the fundamental properties of
epidemic disease diffusion and has often helped taking farsighted decisions in the critical moments
of global pandemics.

In this set of notes, we certainly do not aim at presenting all possible mathematical models
that can be effectively used to describe infectious diseases, nor at giving a fully ``realistic'' description of
an epidemic in the real world. Just, we briefly recall a very classical approach which, in spite of its simplifications
and limitations, can already highlight the great potential of mathematics in dealing with epidemiology
and can obviously open the possibility of presenting more sophisticated and bespoke models to address more
realistically concrete cases of epidemics.

The setting that we recall here 
builds on the research (among the others) of Sir Ronald Ross, Hilda Phoebe Hudson, William Ogilvy Kermack and Anderson Gray McKendrick and is called the \index{SIR model} ``SIR model'', not because of Ross' formal honorific address, but
because it divides the whole population into three compartments:
\begin{itemize}
\item {\em Susceptible:} individuals who might become infected if exposed to the infectious agent (e.g., a virus),
\item {\em Infected (and Infectious):} individuals who are currently infected and can transmit the
infection to susceptible individuals (e.g., by contact),
\item {\em Recovered (or, better to say, Removed):} individuals who, after being
infected, become immune to the infection
(in this category, it is also common to place the people who died for the infection,
since, like the ones that recovered, they cannot contribute to the spread of the disease).\end{itemize}
This is an example of compartmental model, \index{compartmental model}
since it splits a given population into ``compartments'' (three compartments, in this case,
that are labeled with the letters~$S$, for susceptible, $I$, for infected, and~$R$, for recovered),
and several types of related models are indeed broadly utilized in epidemiology.

For simplicity, we can consider a large population and denote by~$S$, $I$ and~$R$ the proportion
of susceptible, infected and recovered individuals, respectively: in this way, $S$, $I$, $R\in[0,1]$
(the value~$1$ corresponding to~$100\%$ of the population,
the value~$\frac12$ to~$50\%$, and so on). Also, possibly adopting the above mentioned convention of counting casualties
in~$R$, we can suppose that the total population remains constant, hence
\begin{equation} \label{SIR:NCO}S+I+R=1.\end{equation}
The SIR model then specifies the evolution~$S(t)$, $I(t)$ and~$R(t)$
of the different individuals over time according 
transition rates between the different compartments.
First of all, in a unit of time, one assumes that some susceptible individual may become infected.
The ansatz in this situation is that the number of new infected people in the unit of time is proportional
to the number of susceptible individuals~$S$ (the higher this number, the easier is that someone catches the disease)
and to the number of infected~$I$ (the more the infected people, the easier for the epidemic to spread).
Accordingly, if we denote by~$\alpha\ge0$ this proportionality coefficient,
in a unit of time, we have that~$\alpha SI$ individuals transit from the compartment
of susceptible individuals to that of infected (the parameter~$\alpha$ can be seen as
a ``transmission rate'').

Concurrently, some infected people can recover. This number is taken to be proportional
to the number of infected:
e.g., if the medicine is effective in a given proportion of treatments,
the higher the number of people receiving the medical treatment, the higher the number of recovered patients
(also in the pessimistic scenario of counting casualties in this compartment,
the higher the number of infected individuals, the higher the number of possible deaths).
Accordingly, if we denote by~$\beta\in[0,1]$ the proportionality coefficient involved in
this transition,
in a unit of time, we have that~$\beta I$ individuals move from the compartment
of infected to that of recovered (the parameter~$\beta$ thus plays the role of a ``recovery rate'').

\begin{figure}
  \centering
  \includegraphics[width=.4\linewidth]{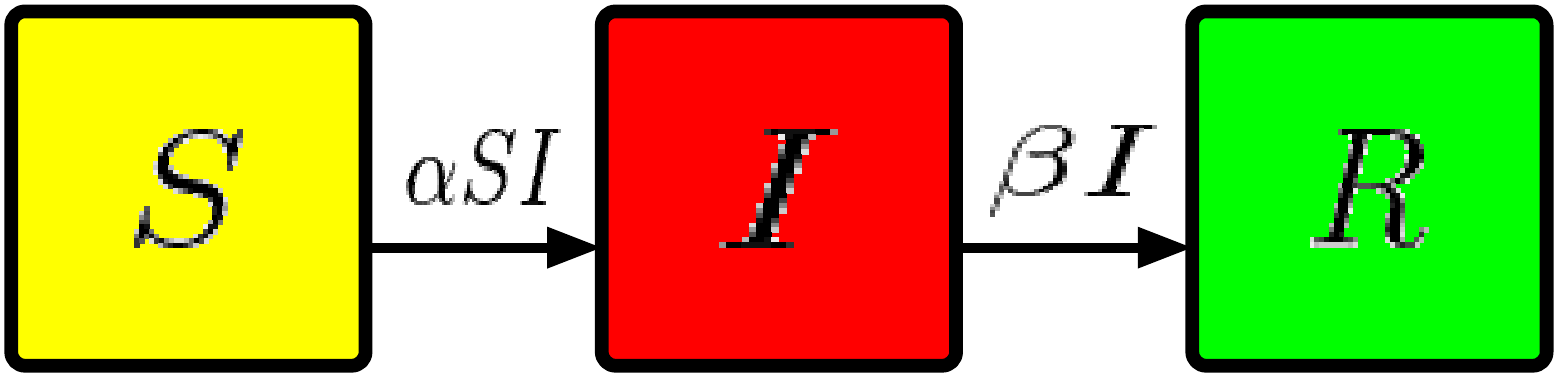}
 \caption{\sl States in a SIR epidemic model and transition between compartments.}\label{SIRCOMTRADItangeFI}
\end{figure}

The transition between compartments described in this way is summarized in Figure~\ref{SIRCOMTRADItangeFI}.

It is certainly convenient to translate Figure~\ref{SIRCOMTRADItangeFI} into a mathematical formulation:
this is done by writing the system of ordinary differential equations corresponding to this compartmental transit,
namely
\begin{equation}\label{SIR-MOD}
\begin{dcases}
\dot{S}=-\alpha SI,\\
\dot{I}=\alpha SI-\beta I,\\
\dot{R}=\beta I,
\end{dcases}
\end{equation}
where the dot represents the derivative with respect to time,
$S=S(t)$, $I=I(t)$ and~$R=R(t)$. The system in~\eqref{SIR-MOD}
is usually taken as the mathematical description of the SIR model.
It is interesting to observe that, by~\eqref{SIR-MOD},
$$ \frac{d}{dt}(S+I+R)=
\dot{S}+\dot{I}+\dot{R}=-\alpha SI+(\alpha SI-\beta I)+\beta I=0,$$
consistently with~\eqref{SIR:NCO}.

Similarly, one can also focus only on the first two equations in~\eqref{SIR-MOD}
to determine the time evolution of the number of susceptible and infected individuals
and then obtain as a byproduct the number of recovered individuals by using~\eqref{SIR:NCO}.
In this way, one can reduce~\eqref{SIR-MOD} to the system of two ordinary differential equations
\begin{equation}\label{SIRRIDO}
\begin{dcases}
\dot{S}=-\alpha SI,\\
\dot{I}=\alpha SI-\beta I.
\end{dcases}
\end{equation}

In spite of its exceptional conceptual simplicity (after all, we are just saying that
susceptible individuals may become infected and
infected individuals may recover, or possibly die), the SIR model already showcases some
very interesting information about the spread of a disease.

First of all, the first equation in~\eqref{SIRRIDO} already suggests that
the infection occurs by the contact between infectious and susceptible
people (since the quantity~$SI$ is a good model for ``random encounters'' between~$S$ and~$I$).
As a consequence, the parameter~$\alpha$ in~\eqref{SIRRIDO} can be considered as
a ``contact rate''.
For this reason, measures of
social distancing and lockdowns aim at reducing contacts between people, hence at reducing the value of~$\alpha$.

The parameter~$\beta$ in the second equation in~\eqref{SIRRIDO} instead reduces the number of infected
people, correspondingly \footnote{Beware that since casualties are also counted among recoveries,
a tragic way to increase the number of recoveries consists in killing infectious people.
This is a classical topic for science fiction horror films, such as
``The Crazies'' by George A. Romero.}
raising the number of recovers.
To increase the parameter~$ \beta$, one effective way is clearly to have better medicines,
since more effective treatments facilitate healing and accelerate recovery rates.
Another useful measure relies on the early detection of infected people and their
isolation in hospitals or quarantine areas: in this way, infected people are de facto
removed from the dynamics of~\eqref{SIRRIDO}
and, from the mathematical point of view, transit to the~$R$ compartment (even if they are not yet recovered
from the medical perspective).

\begin{figure}
                \centering
                \includegraphics[height=.12\textheight]{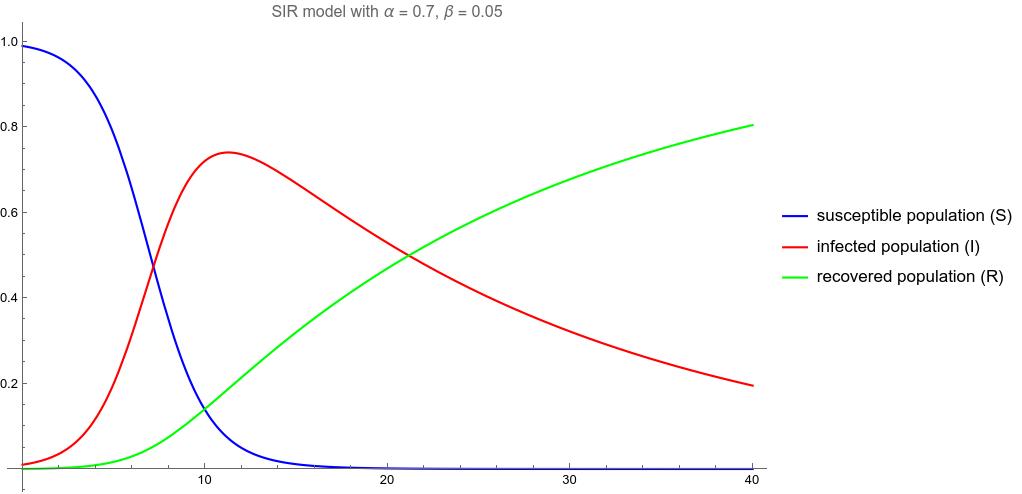}$\quad$
                                \includegraphics[height=.12\textheight]{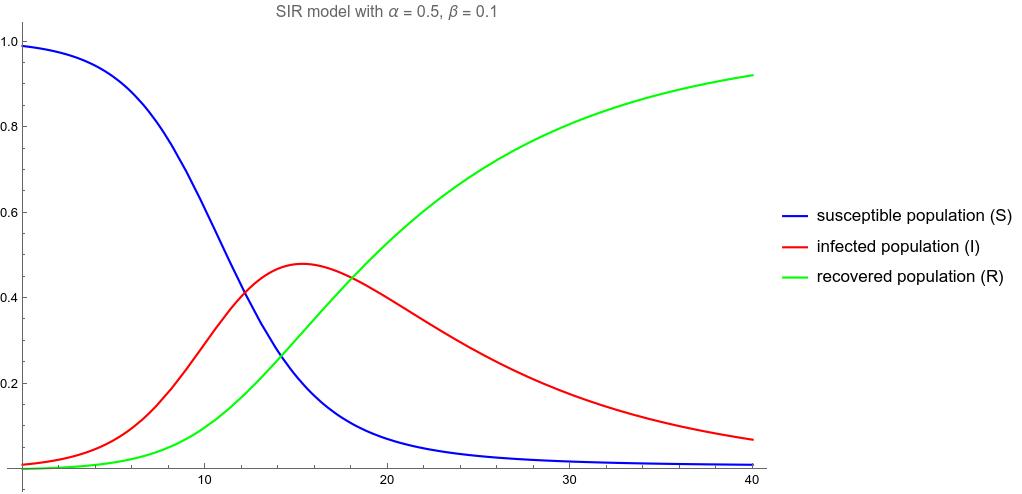}$\quad$
                \includegraphics[height=.12\textheight]{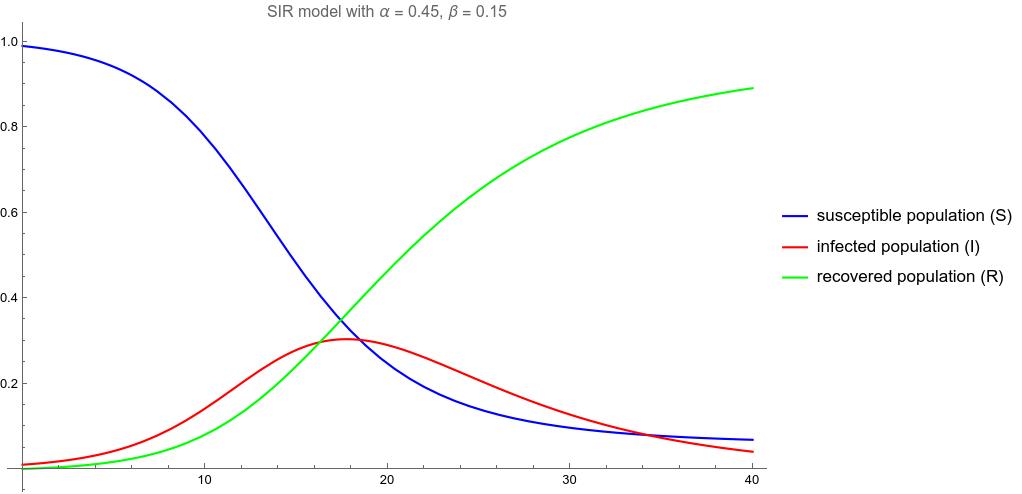}
        \caption{\sl Flattening the curve of infection: the parameters~$(\alpha,\beta)$
are chosen here in~$\{(0.7,0.05),\;(0.5,0.1),\;(0.45,0.15)\}$.}\label{RIVEKELVHASPREPINOJHNFOJED}
\end{figure}

See Figure~\ref{RIVEKELVHASPREPINOJHNFOJED} to appreciate how decreasing the transmission rate~$\alpha$
and increasing the recovery rate~$\beta$
can help slowing an epidemic's spread.\medskip

In the description of epidemic diseases, it is also interesting to consider the basic reproduction number
(often called ``R naught'' in jargon), e.g. defined by \index{R naught}
\begin{equation}\label{miMNDfgkdaPKSMo0liRnou53} {\mathcal{R}}_o:=\frac\alpha\beta .\end{equation}
Its importance lies in the fact that, by~\eqref{SIRRIDO},
\begin{equation}\label{miMNDfgkdaPKSMo0li}
\dot{I}=\beta I\left(\frac\alpha\beta S-1\right)=\beta I\left({\mathcal{R}}_o\, S-1\right)
\end{equation}
and thus, since~$S\in[0,1]$,
$$ \dot{I}\le\beta I\left({\mathcal{R}}_o-1\right).$$
Hence, when
\begin{equation}\label{2rterOPT14534yhSHE}
{\mathcal{R}}_o<1,\end{equation} we find that
$$ \frac{d}{dt} \ln I(t)=\frac{\dot{I}(t)}{I(t)}\le \beta \left({\mathcal{R}}_o-1\right)=-\beta \left(1-{\mathcal{R}}_o\right),$$
leading to
$$ I(t)\le I(0)\, e^{-\beta \left(1-{\mathcal{R}}_o\right)t} ,$$
which means that the number of infected people decreases exponentially fast and the epidemic is
likely to be overcome sufficiently quickly (conversely, when~${\mathcal{R}}_o>1$
one can expect that the number of infected people will grow exponentially, with obvious tragic consequences).

The computation in~\eqref{miMNDfgkdaPKSMo0li} also reveals an interesting
information about the possibility of embanking the spread of the epidemic by reducing the
number of susceptible individuals: indeed, when
\begin{equation}\label{HERDI}{\mathcal{R}}_o\, S<1\end{equation} a similar computation
as above would lead to an exponential decreasing of infected individuals.
This is related to the notion of ``herd immunity'' \index{herd immunity},
namely in obtaining a community protection due to a critical proportion of the population having become immune to the disease:
notice indeed that condition~\eqref{HERDI} prescribes that, to reach this situation,
at least a proportion~${\mathcal{H}}_o:=1-\frac{1}{{\mathcal{R}}_o}$ must have become immune.

If we wish to play with some numbers, an initial estimate of
the World Health Organization on~${\mathcal{R}}_o$ for COVID-19
took into account a possibility of~$ {\mathcal{R}}_o\simeq 2.4$,
thus largely violating the optimistic threshold in~\eqref{2rterOPT14534yhSHE}.
Also, with this value of~${\mathcal{R}}_o$, one would obtain a herd immunity threshold~${\mathcal{H}}_o\simeq
1-\frac1{2.4}=0.58\overline{3}$, giving that about~$60\%$ of the population
must be immune to protect the community from an exponential spread of epidemic
(but not from the disease\footnote{We stress that reaching the herd immunity
does not stop the epidemic overnight.
It only makes the curve of newly infected
people bend down and approach to zero exponentially fast. But, according to
the SIR model, other individuals will get infected even after having reached
the status of herd immunity and the curve representing the total number
of infected people will keep increasing towards a horizontal asymptote.

We also remark that we are not accounting here for important realistic
features such as the possible diminished effectiveness of a vaccine on a rapidly mutating virus, the possibility
of repeated infections within the population, the different personal responses of individuals to the virus exposition, etc.}
in itself).

Reaching such a threshold of immunization without a vaccination campaign is likely extremely dangerous
(since it entails that a vast majority of the population must catch the disease). Since, realistically, vaccines are also not~$ 100\%$
effective in avoiding the dissemination of the virus, if one aims at reaching a herd immunity via a vaccination campaign,
the efficacy of the vaccine is an important parameter to take into account.
For instance, if~$80\%$ of the population is vaccinated with a vaccine which is effective in~$70\%$
of the cases in avoiding new infections, we have that~$56\%$ of the population
is immune to the disease (and this number is still below the~$60\%$ threshold
discussed above for herd immunity).

Of course, we are not proposing here concrete measures to deal with COVID-19,
and the numbers plugged in for the naive computations above must not be taken seriously,
simply we wanted to stress how important is to deal with epidemics using a scientific approach
and how useful mathematics can be for such a life-or-death situation.\medskip

Another striking application of the SIR model consists in the
prediction of the epidemic peak, namely of the maximal number~$I_{\max}$ of people
that get simultaneously infected during the epidemic.

Indeed, one can deduce from~\eqref{SIRRIDO} that
\begin{equation}\label{IMAXESYT} I_{\max} \simeq\frac{{\mathcal{R}}_o- 1-\ln {\mathcal{R}}_o}{{\mathcal{R}}_o}.\end{equation}
Notice that as~${\mathcal{R}}_o\to+\infty$ we have that~$ I_{\max}\to1$, which corresponds to the total population
ending up being infected at the same time.
We also remark that estimating~$I_{\max}$ is of great strategic importance,
since to accommodate patients, one always wants to have a number of available beds in the hospitals
commensurate to the number of people in need of hospitalization.

To prove~\eqref{IMAXESYT}, one can pick a time~$t_\star$ at which we expect the number of infected
individuals to be maximal, i.e. such that $I_{\max}=I(t_\star)$. Then, by~\eqref{SIRRIDO},
$$ 0=\frac{dI}{dt}(t_\star)=\alpha \,S(t_\star)\,I(t_\star)-\beta\, I(t_\star),$$
from which we infer that
$$ S(t_\star)=\frac\beta\alpha=\frac1{{\mathcal{R}}_o}.$$
Notice that this is not quite what~\eqref{IMAXESYT} is looking for, since~$ S(t_\star)$ represents only
the number of susceptible individuals at the time maximizing the number of infected people.
Nevertheless, the structure of~\eqref{SIRRIDO} comes in handy now and provides that
\begin{eqnarray*} \frac{d(I+S)}{dt}&=& (\alpha SI-\beta I)-\alpha SI\\
&=&-\beta I\\
&=&-\frac{\alpha\beta SI}{\alpha S}\\
&=&\frac{\beta}{\alpha S}\frac{dS}{dt}\\
&=&\frac{1}{{\mathcal{R}}_o}\frac{d}{dt}(\ln S).
\end{eqnarray*}
Therefore,
$$ I(t)+S(t)=I(0)+S(0)+\frac{\ln(S(t))-\ln(S(0))}{{\mathcal{R}}_o}$$
and, as a result,
\begin{eqnarray*}
I_{\max}&=&I(t_\star)\\
&=&I(0)+S(0)+\frac{\ln(S(t_\star))-\ln(S(0))}{{\mathcal{R}}_o}-S(t_\star)\\
&=&I(0)+S(0)+\frac{-\ln {\mathcal{R}}_o-\ln(S(0))}{{\mathcal{R}}_o}-\frac1{{\mathcal{R}}_o}\\&=&
I(0)+S(0)-
\frac{1+\ln {\mathcal{R}}_o+\ln(S(0))}{{\mathcal{R}}_o}.\end{eqnarray*}
If we assume that
at the beginning of the epidemic there is nobody that is recovering, namely\footnote{No confusion should
arise between the initial number of recovered~$R(0)$ and the structural parameter~${\mathcal{R}}_o$
in~\eqref{miMNDfgkdaPKSMo0liRnou53}.} that~$R(0)=0$, we
deduce from~\eqref{SIR:NCO} that~$I(0)+S(0)=1$,
and thus
$$ I_{\max}=1-
\frac{1+\ln {\mathcal{R}}_o+\ln(S(0))}{{\mathcal{R}}_o},$$
that completes the proof of~\eqref{IMAXESYT}.
\medskip

We also mention that explicit solutions of the SIR model are available by using
an implicit ``time reparameterization''. Specifically, one considers
the integral of~$I$ as a ``new time''~$\tau$ by setting
\begin{equation}\label{EXPLISIR2} \tau(t):=\int_0^t I(s)\,ds.\end{equation}
Note that~$\frac{\tau}{t}$ represents the average number of
infected people at time~$t$. In this new variable, the solution of~\eqref{SIR-MOD} becomes explicit
and takes the form
\begin{equation}\label{EXPLISIR}\begin{dcases}
\tilde S(\tau)=S(0)\,e^{-\alpha \tau},\\
\tilde I(\tau)=I(0)+S(0)(1-e^{-\alpha \tau})-\beta \tau,\\
\tilde R(\tau)=R(0)+\beta \tau.
\end{dcases}
\end{equation}
Of course, one of the limitations of this explicit formulation is that the new solutions~$\tilde S$,
$\tilde I$ and~$\tilde R$ are functions of an implicitly defined time~$\tau$ through the relations~$\tilde S(\tau(t))=S(t)$,
$\tilde I(\tau(t))=I(t)$ and~$\tilde R(\tau(t))=R(t)$, whence the explicit representation in~\eqref{EXPLISIR}
is not always particularly pleasant for concrete calculations (see however the forthcoming equation~\eqref{34EXPLISIR}
for a reformulation of the time parameter problem).

To check the validity of~\eqref{EXPLISIR} one can proceed by implicit differentiation,
observing that from~\eqref{SIR-MOD} it follows that
\begin{equation*}\begin{dcases}\displaystyle
\frac{d\tilde R}{d\tau}=\beta,\\ \displaystyle
\frac{d\tilde S}{d\tau}=-\alpha \tilde{S},
\end{dcases}
\end{equation*}
which are explicitly solvable leading to~$\tilde S(\tau)$ and~$\tilde R(\tau)$ in~\eqref{EXPLISIR}.
Then, using
$$\frac{d\tilde I}{d\tau}=\alpha \tilde S-\beta ,$$
and exploiting the previous explicit solution for~$\tilde S(\tau)$ one
obtains the expression of~$\tilde I(\tau)$ in~\eqref{EXPLISIR}.

Having completed the proof of~\eqref{EXPLISIR},
we also observe that, by~\eqref{EXPLISIR2},
$$ \frac{d\tau}{dt}=\tilde I(\tau),$$
from which one obtains an explicit integral formula relating the old and the new times of the type
\begin{equation} \label{34EXPLISIR}
t=\int_0^\tau\frac{ds}{I(0)+S(0)(1-e^{-\alpha s})-\beta s}.\end{equation}
Once again, this time reparameterization is explicit, but the integral cannot be usually
written in terms of standard elementary mathematical functions.
\medskip

We stress that the classical SIR model in~\eqref{SIRRIDO} does not take into account the possible
mobility of individuals (roughly speaking, the whole population in~\eqref{SIRRIDO}
is concentrated at the same place). There are of course a number of models available in the literature
which account for possible spatial displacements of the population: for instance, if one models
the space configuration by a variable~$x$ and assumes that the population moves according to a random walk
as in Section~\ref{CHEMOTX}, one may think about the possibility of replacing the ordinary differential
equations in~\eqref{SIRRIDO} with a system of partial differential equations of the type
\begin{equation}\label{9-42coMS02urjf}
\begin{dcases}
\partial_t{S}=\mu\Delta S-\alpha SI,\\
\partial_t{I}=\nu\Delta I+\alpha SI-\beta I.
\end{dcases}
\end{equation}
In this case, $S=S(x,t)$, $I=I(x,t)$ and the diffusion coefficients~$\mu$, $\nu\in[0,+\infty)$
may also be different to account for the possibility of a different speed between susceptible and infected people
(e.g., in a situation in which the diffusivity of infected individuals is limited by the illness itself,
or by the responsibility of self-isolating individuals with symptoms).

Restrictions on travel are of course a consequential measure to avoid or limit the spread of an epidemic
due to the diffusive nature of the system in~\eqref{9-42coMS02urjf}.
\medskip

See for instance~\cite{zbMATH05276332} for additional information on the SIR models
and for many related models utilized in mathematical epidemiology.\medskip

Now, for completeness, also in view of the discussion presented in
Section~\ref{9-42coMS02urjf0iurjfpnewnblkfxd},
we recall here a variation of the SIR model based on the evolution of
the probability density~$p(x,y,t)$ corresponding to having, at time~$t$,
a proportion~$x\in[0,1]$ of susceptible individuals and~$y\in[0,1]$ of infected ones,
as proposed in~\cite{MR2782833}.
The core of this model is to obtain a partial differential equation of the form
\begin{equation}\label{MR2782833:EQUA}
\partial_t p+\div(p F)=0,\end{equation}
where the above divergence is taken with respect to the variables~$(x,y)$
and the vector field~$F$ corresponds to the right hand side of
the SIR model in~\eqref{SIRRIDO}, that is
\begin{equation}\label{MR2782833:EQUA2}
F(x,y):=\left( -\alpha xy,\; \alpha xy-\beta y\right).
\end{equation}
We remark that~\eqref{MR2782833:EQUA} has the form of a
transport equation \index{transport equation} as in~\eqref{OJS-PJDN-0IHGDOIUGDBV02ujrfMTE}.

To deduce~\eqref{MR2782833:EQUA} from the principles of the SIR model one can argue as follows.
Given~$n$, $m\in\N$ with~$n+m\le N$, we consider the probability~$P(n,m,t)$ of having~$n$ susceptible and~$m$ infected at time~$t$, being~$N$  the total number of individuals
(notice that, in view of the constancy of the population in~\eqref{SIR:NCO}, the number of recovered is thus necessarily~$N-n-m$). We argue that, given an ``infinitesimal'' time increment~$\tau$, the probability~$P(n,m,t+\tau)$ of having~$n$ susceptible and~$m$ infected at time~$t+\tau$ is built by the superposition of three occurrences, namely:\begin{itemize}\item either, during the elapsed time~$\tau$, an infected individual may have recovered: this corresponds to the probability~$P(n,m+1,t)$ of having one more infected at time~$t$, times the probability of recovering, which, inspired by Figure~\ref{SIRCOMTRADItangeFI}, we take to be equal to~$\frac{\beta(m+1)}{N}$,\item or, during the elapsed time~$\tau$, a susceptible individual may have got infected: this corresponds to the probability~$P(n+1,m-1,t)$ of having one more susceptible and one less infected at time~$t$, times the probability\footnote{In all this discussion, we are
implicitly assuming that~$\frac{\beta(m+1)}{N}$, $\frac{\alpha (n+1)(m-1)}{N^2}\in[0,1]$,
consistently with the standard notion of probability. This hypothesis is justified, for instance,
if the number of infected~$m$ is sufficiently small with respect to the total population~$N$.}
of getting infected, which, inspired by Figure~\ref{SIRCOMTRADItangeFI}, we take to be equal to~$\frac{\alpha (n+1)(m-1)}{N^2}$,\item or, during the elapsed time~$\tau$, no susceptible individual got infected
and no infected got recovered: this corresponds to the probability~$P(n,m,t)$ of having the same number of susceptible and infected individuals at time~$t$, times the probability of having no infections or recovery, which is equal to the remaining probability~$1-
\frac{\beta m}{N}-\frac{\alpha nm}{N^2}$.\end{itemize}
Due to these considerations, we write that\begin{equation}\label{SIRMOPDEMD}\begin{split}&P(n,m,t+\tau)=\frac{\beta(m+1)}{N}P(n,m+1,t)+\frac{\alpha (n+1)(m-1)}{N^2}P(n+1,m-1,t)\\&\qquad\qquad\qquad\qquad\qquad+
\left(1-\frac{\beta m}N-\frac{\alpha nm}{N^2}\right)P(n,m,t).\end{split}\end{equation}Now we try to achieve a continuum model by taking the limit as~$N\to+\infty$, corresponding to a large population. For this, for all~$x$, $y\in[0,1]$
we  define $$p(x,y,t):=NP(Nx,Ny,t).$$
In this way, we deduce from~\eqref{SIRMOPDEMD} that\begin{eqnarray*}&&\partial_t p(x,y,t)\\&\simeq&\frac{p(x,y,t+\tau)-p(x,y,t)}{\tau}\\&=&\frac{N}{\tau}\Big[ P(Nx,Ny,t+\tau)-P(Nx,Ny,t)\Big]\\&=&\frac{N}{\tau}\Bigg[ \frac{\beta(Ny+1)}{N}P(Nx,Ny+1,t)+
\frac{\alpha (Nx+1)(Ny-1)}{N^2}P(Nx+1,Ny-1,t)\\&&\qquad\qquad+
\left(1- \beta y-\alpha xy\right)P(Nx,Ny,t)-P(Nx,Ny,t)\Bigg]\\
&=&\frac{1}{\tau}\Bigg[ 
\beta\left(y+\frac1{N}\right)p\left(x,y+\frac1N,t\right)+
\alpha \left(x+\frac1N\right)\left(y-\frac1N\right)p\left(x+\frac1N,y-\frac1N,t\right)\\&&\qquad\qquad-
\left( \beta y+ \alpha xy\right)p(x,y,t)
\Bigg]\\&=&
\frac{1}{\tau}\Bigg[B\left(x,y+\frac1N,t\right)+A\left(x+\frac1N,y-\frac1N,t\right)-B(x,y,t)-A(x,y,t)
\Bigg]\\&\simeq&
\frac{1}{\tau}\Bigg[\frac1N\partial_y B(x,y,t)+
\frac1N\left(\partial_x A(x,y,t)-\partial_yA(x,y,t)\right) \Bigg]
,\end{eqnarray*}
where
$$ A(x,y,t):=\alpha xy\, p(x,y,t)
\qquad{\mbox{and}}\qquad
B(x,y,t):=\beta y \,p(x,y,t).$$
Hence, by choosing a scaling between the time step and the total population such that~$\tau N\simeq1$,
we reduce the infectious disease spreading problem to 
$$ \partial_t p=\partial_x A+\partial_y (B-A)=
\partial_x  (\alpha xy p)+\partial_y\big((\beta y-\alpha xy)p\big),
$$
which, due to the choice of the vector field in~\eqref{MR2782833:EQUA2}, corresponds to the desired
transport equation in~\eqref{MR2782833:EQUA}.\medskip

For further motivations about elliptic partial differential equations, see e.g. page~\pageref{0uojf29249-45kpkfdSmd11493839429efv}
here, Chapter~2 in~\cite{MR0029463},
Chapter~4 in~\cite{MR1149383},
Chapter~1 in~\cite{MR1160355},
pages~160--168 in~\cite{MR1398883},
the introduction of~\cite{MR1625845},
Chapter~1 in~\cite{MR2007714},
Section~1.4 in~\cite{MR2164768}, Chapter~1 in~\cite{MR2398759},
Chapters~4 and~5 in~\cite{MR2488750},
Section~1.2 in~\cite{MR3362185},
the appendix of~\cite{ugjbfvCSnxcAT35Mowjfgoe3AGSHDKS}
and the references therein.\medskip

Having established
solid heuristic motivations to study partial differential equations
and especially the Laplace operator, we now enter the mathematical theory
of these objects.
{F}rom now on, we are not allow any longer to ``waive hands'',
provide formal expansions or take additional simplifying assumptions,
instead we will have to fully comply to the rigorous mathematical precision.

\section{What is the Laplacian, after all?}

We now dive into some mathematical setting
which turns out to be handy when dealing with partial differential equations.
Given an open set~$\Omega\subseteq\R^n$,
a function~$u\in C^2(\Omega)$ and a point~$x\in\Omega$, we consider the ``Laplace operator''~$\Delta$
(which we have already informally met and used, from~\eqref{DAGB-ADkrVoiweLL4re2346ytmngrrUj7} on)
defined by
\begin{equation}\label{LA:DIFF}
\Delta u(x):=\sum_{j=1}^n \frac{\partial^2u}{\partial x_j^2}(x).\end{equation}
We point out that the Laplacian is invariant under translation,
namely
\begin{equation}\label{TRAL} \Delta \big(u(x+y)\big)=\Delta u(x+y)\qquad{\mbox{for every~$x\in\Omega$ and~$y\in\R^n$ such that~$x+y\in\Omega$,}}\end{equation}
and it
possesses a ``divergence form''\index{divergence form} structure, since
\begin{equation}\label{DITRO} \Delta u={\rm div}(\nabla u).\end{equation}
As a consequence of this and of the Divergence Theorem, if~$\Omega$
has a boundary of class~$C^1$,
the average of the Laplacian
can be reconstructed by the normal flow through the boundary of the domain, namely
\begin{equation}\label{1111DIV0934} \int_\Omega \Delta u(x)\,dx=\int_{\partial\Omega } \nabla u(x)\cdot \nu(x)\,d{\mathcal{H}}^{n-1}_x,\end{equation}
where~$\nu$ denotes the external unit normal of~$\Omega$ and~${\mathcal{H}}^{n-1}$
is the $(n-1)$-dimensional Hausdorff measure.\medskip

It is also useful to remark that the Laplacian measures the ``infinitesimal distance''
between the value of a function at a given point and the average nearby, namely:

\begin{theorem}\label{KAHAR2}
Let~$x_0\in\R^n$ and~$r>0$. Suppose that~$u\in C^2(B_r(x_0))$.
Then,
\begin{equation}\label{LI} \lim_{\rho\searrow0} \frac1{\rho^2}\left(
\fint_{B_\rho(x_0)} u(x)\,dx-u(x_0)\right)=\frac{1}{2(n+2)}\,\Delta u(x_0).\end{equation}
\end{theorem}

\begin{proof}
We start by recalling a standard relation between the Lebesgue measure of
the unit ball and the $(n-1)$-dimensional Hausdorff measure of the unit sphere. Namely,
using polar coordinates, we see that
\begin{equation}\label{B1}
|B_1|={\mathcal{H}}^{n-1}(\partial B_1)\,\int_0^1 t^{n-1}\,dt=\frac{ {\mathcal{H}}^{n-1}(\partial B_1)}{n}.\end{equation}

We now point out a useful cancellation property. If~$T_j$ is the reflection across the~$j$th coordinate, i.e.
$$ T_j (y_1,\dots,y_n)=(y_1,\dots,y_{j-1},-y_j,y_{j+1},\dots,y_n),$$
and~$g: B_\rho\to \R$ is such that~$g(T_j(y))=-g(y)$ for all~$y\in B_\rho$, using the change of variable~$Y:=T_j(y)$
it follows that
\begin{equation}\label{22345XXse5q6u5at3i332o2n} \int_{B_\rho} g(y)\,dy=-\int_{B_\rho} g(T_j(y))\,dy=-\int_{B_\rho} g(Y)\,dY,\end{equation}
and therefore
$$ \int_{B_\rho} g(y)\,dy=0.$$
Applying this observation to the function~$y_j$ for every~$j\in\{1,\dots,n\}$,
we find that
\begin{equation}\label{CAN1} \int_{B_\rho} y_j\,dy=0.\end{equation}
Similarly, applying the previous observation to the function~$y_j y_k$
for every~$j\ne k\in\{1,\dots,n\}$,
\begin{equation}\label{CAN2} \int_{B_\rho} y_j y_k\,dy=0.\end{equation}

Now we perform an explicit computation in the ball by taking advantage
of the symmetry between coordinate exchanges.
Namely, for every~$j\in\{1,\dots,n\}$, using polar coordinates we have that
\begin{equation}\label{H66}
n\int_{B_\rho} y_{j}^2\,dy=
\sum_{k=1}^n\int_{B_\rho} y_k^2\,dy=
\int_{B_\rho} |y|^2\,dy={\mathcal{H}}^{n-1}(\partial B_1)\,\int_0^\rho t^{n+1}\,dt=
\frac{ {\mathcal{H}}^{n-1}(\partial B_1)\, \rho^{n+2} }{n+2}.
\end{equation}

Now, we focus on the proof
of the desired result in~\eqref{LI}. For this,
given~$\rho>0$ and~$x\in B_\rho(x_0)$, we use the Taylor expansion
\begin{equation}\label{TAYTAY}
u(x)=u(x_0)+\nabla u(x_0)\cdot (x-x_0)+\frac12 D^2 u(x_\star)(x-x_0)\cdot(x-x_0),\end{equation}
for a suitable~$x_\star \in B_\rho(x_0)$ and we find that
\begin{eqnarray*}
\fint_{B_\rho(x_0)} u(x)\,dx-u(x_0)&=&\fint_{B_\rho(x_0)}\big( u(x)-u(x_0)\big)\,dx\\&=&
\fint_{B_\rho(x_0)} \left(
\nabla u(x_0)\cdot (x-x_0)+\frac12 D^2 u(x_\star)(x-x_0)\cdot(x-x_0)\right)
\,dx\\&=&
\sum_{j=1}^n
\fint_{B_\rho} \frac{\partial u}{\partial x_j}(x_0)\,y_j\,dy+\frac12\sum_{j,k=1}^n
\fint_{B_\rho} \frac{\partial^2 u}{\partial x_j\partial x_k}(x_\star)\,y_j y_k\,dy\\&=&\frac12
\sum_{j,k=1}^n
\fint_{B_\rho} \frac{\partial^2 u}{\partial x_j\partial x_k}(x_\star)\,y_j y_k\,dy,
\end{eqnarray*}
where the latter identity is a consequence of~\eqref{CAN1}
and~\eqref{CAN2}.

Consequently, by~\eqref{H66},
\begin{equation}\label{ETA0}\begin{split} \fint_{B_\rho(x_0)} u(x)\,dx-u(x_0)\,&=\,\frac12
\sum_{j,k=1}^n
\fint_{B_\rho} \frac{\partial^2 u}{\partial x_j\partial x_k}(x_\star)\,y_j y_k\,dy\\&=\,\frac12
\sum_{j=1}^n
\fint_{B_\rho} \frac{\partial^2 u}{\partial x_j^2 }(x_0)\,y_j^2\,dy+\eta(\rho)
,\end{split}\end{equation}
where
\begin{equation}\label{ETA} \eta(\rho):=
\frac12 \sum_{j,k=1}^n\fint_{B_\rho}\left(\frac{\partial^2 u}{\partial x_j\partial x_k}(x_\star)-
\frac{\partial^2 u}{\partial x_j\partial x_k}(x_0)\right)\,y_j y_k\,dy.\end{equation}

Now we observe that
\begin{equation}\label{H9}
\lim_{\rho\searrow0}\frac{\eta(\rho)}{\rho^{2}}=0.
\end{equation}
Indeed, since~$x_\star$ approaches~$x_0$ as~$\rho\searrow0$,
given any~$\e>0$, if~$\rho$ is small enough we have that
$$ \left|\frac{\partial^2 u}{\partial x_j\partial x_k}(x_\star)-
\frac{\partial^2 u}{\partial x_j\partial x_k}(x_0)\right|\le\e$$
and consequently, using again~\eqref{H66},
$$ |\eta(\rho)|\le
\frac{\e}2 \sum_{j=1}^n\fint_{B_\rho} y_j^2\,dy=
\frac{ {\mathcal{H}}^{n-1}(\partial B_1)\, \rho^{2}\,\e }{2(n+2)\,|B_1|}
.$$
{F}rom this, the claim in~\eqref{H9} plainly follows.

The desired result in~\eqref{LI} is now a direct consequence of~\eqref{B1},
\eqref{H66}, \eqref{ETA0}
and~\eqref{H9}, since
\begin{eqnarray*}&&
\lim_{\rho\searrow0} \frac1{\rho^2}\left(
\fint_{B_\rho(x_0)} u(x)\,dx-u(x_0)\right)=
\lim_{\rho\searrow0} \left(
\frac1{2\rho^2}
\sum_{j=1}^n
\fint_{B_\rho} \frac{\partial^2 u}{\partial x_j^2}(x_0)\,y_j^2\,dy+\frac{\eta(\rho)}{\rho^2}\right)
\\&&\qquad=
\frac{ {\mathcal{H}}^{n-1}(\partial B_1) }{2n(n+2)\,|B_1|}
\sum_{j=1}^n \frac{\partial^2 u}{\partial x_j^2}(x_0)=\frac{1}{2(n+2)}\,\Delta u(x_0),
\end{eqnarray*}
as desired.
\end{proof}

For completeness, as a variant of Theorem~\ref{KAHAR2},
we also point out a similar result for the limit of the spherical
averages:

\begin{theorem}\label{KAHAR2-SPHE}
Let~$x_0\in\R^n$ and~$r>0$. Suppose that~$u\in C^2(B_r(x_0))$.
Then,
\begin{equation*} \lim_{\rho\searrow0} \frac1{\rho^2}\left(
\fint_{\partial
B_\rho(x_0)} u(x)\,d{\mathcal{H}}^{n-1}_x-u(x_0)\right)=
\frac{1}{2n}\,\Delta u(x_0).\end{equation*}
\end{theorem}

\begin{proof} By the
Taylor expansion in~\eqref{TAYTAY} and two odd cancellation
arguments, as~$\rho\searrow0$ we have
\begin{eqnarray*}
&& \fint_{\partial
B_\rho(x_0)} u(x)\,d{\mathcal{H}}^{n-1}_x-u(x_0)\\&=&
\fint_{\partial
B_\rho(x_0)} \left(
\nabla u(x_0)\cdot (x-x_0)+\frac12 D^2 u(x_0)(x-x_0)\cdot(x-x_0)\right)
\,d{\mathcal{H}}^{n-1}_x
+o(\rho^2)\\
&=&
\frac12\,\fint_{\partial
B_\rho(x_0)} D^2 u(x_0)(x-x_0)\cdot(x-x_0)
\,d{\mathcal{H}}^{n-1}_x
+o(\rho^2)\\&=&
\frac12\,\sum_{i,j=1}^n\fint_{\partial
B_\rho} \partial_{ij} u(x_0)y_iy_j
\,d{\mathcal{H}}^{n-1}_y
+o(\rho^2)\\&=&
\frac12\,\sum_{i=1}^n\fint_{\partial
B_\rho} \partial_{ii} u(x_0)y_i^2
\,d{\mathcal{H}}^{n-1}_y
+o(\rho^{2}).
\end{eqnarray*}
Now, for every~$i\in\{1,\dots,n\}$, we have that
$$ n\,\int_{\partial B_\rho} y_i^2
\,d{\mathcal{H}}^{n-1}_y=\sum_{k=1}^n
\int_{\partial B_\rho} y_k^2
\,d{\mathcal{H}}^{n-1}_y=\int_{\partial B_\rho} |y|^2
\,d{\mathcal{H}}^{n-1}_y=
\int_{\partial B_\rho} \rho^2
\,d{\mathcal{H}}^{n-1}_y=\rho^{2}\,{\mathcal{H}}^{n-1}(\partial B_\rho)
$$
and consequently
\begin{eqnarray*}
&& \fint_{\partial
B_\rho(x_0)} u(x)\,d{\mathcal{H}}^{n-1}_x-u(x_0)=
\frac1{2n}\,\sum_{i=1}^n \partial_{ii} u(x_0)\rho^2
+o(\rho^2)=\frac{\rho^2}{2n}\,\Delta u(x_0)+o(\rho^2),
\end{eqnarray*}
that plainly leads to the desired result.
\end{proof}

We point out that
Theorems~\ref{KAHAR2} and~\ref{KAHAR2-SPHE}
are special cases of a more general type of result, known in the literature
as Pizzetti's Formula\index{Pizzetti's Formula}, involving the
higher order Laplace operator
$$ \Delta^k:=\underbrace{\Delta\dots\Delta}_{\footnotesize{\mbox{$k$ times}}},$$
with~$k\in\N$ and the convention that~$\Delta^0$ is the identity operator.
We recall this setting for the sake of completeness:

\begin{proposition}
Let~$x_0\in\R^n$ and~$r>0$. Suppose that~$u\in C^{2N}(B_r(x_0))$. Then, as~$\rho\searrow0$,
\begin{eqnarray}&&
\label{PALLAPIZZE1} \fint_{B_\rho(x_0)}u(x)\,dx
=
n\, \Gamma\left(\frac{n}2\right)\sum_{k=0}^N
\frac{\rho^{2k}}{2^{2k+1}\,
k!\,\Gamma\left(\frac{n}2+k+1\right)}\Delta^k u(x_0)+o(\rho^{2N})
\\
{\mbox{and }}&&
\label{PALLAPIZZE2} \fint_{\partial B_\rho(x_0)}u(x)\,d{\mathcal{H}}^{n-1}_x=
\Gamma\left(\frac{n}2\right)\sum_{k=0}^N
\frac{\rho^{2k}}{2^{2k}\,k!\,\Gamma\left(\frac{n}2+k\right)}\Delta^k u(x_0)+o(\rho^{2N}).\end{eqnarray}
Here above, $\Gamma$ stands for the Euler Gamma Function, defined,
for every~$z\in\cOMPL$ with~$\Re z>0$, by
\begin{equation}\label{EULEROGA}
\Gamma (z):=\int _{0}^{+\infty }t^{z-1}e^{-t}\,dt.\end{equation}
\end{proposition}

\begin{proof} We start with an observation that can be considered
as a version of the Multinomial Theorem for the Laplace operator. We claim that,
for every~$k\in\N$,
\begin{equation}\label{MULTOLPALPA}
\Delta^k u=\sum_{{\alpha\in\N^n}\atop{|\alpha|=k}}\frac{k!}{\alpha!}\,D^{2\alpha}u.
\end{equation}
Here, we are using the multi-index notation for which~$\alpha=(\alpha_1,\dots,\alpha_n)$,
$|\alpha|=\alpha_1+\dots+\alpha_n$, $\alpha!=\alpha_1!\dots\alpha_n!$
and~$2\alpha=(2\alpha_1,\dots,2\alpha_n)$.
We prove~\eqref{MULTOLPALPA} by induction over~$k$.
Indeed, when~$k=0$, 
the claim in~\eqref{MULTOLPALPA} reduces to the true identity~$u=u$.
Then, we suppose that~\eqref{MULTOLPALPA} holds true for the index~$k$
and observe that, being~$e_i$ the $i$th element of the Euclidean basis,
\begin{eqnarray*}&&
\Delta^{k+1} u=\Delta\left(
\sum_{{\alpha\in\N^n}\atop{|\alpha|=k}}\frac{k!}{\alpha!}\,D^{2\alpha}u\right)=
\sum_{i=1}^n\frac{\partial^2}{\partial x_i^2}
\left(
\sum_{{\alpha\in\N^n}\atop{|\alpha|=k}}\frac{k!}{\alpha!}\,D^{2\alpha}u\right)
=
\sum_{{{\alpha\in\N^n}\atop{|\alpha|=k}}\atop{1\le i\le n}}
\frac{k!}{\alpha!}\,D^{2(\alpha+e_i)}u\\&&\qquad\qquad\qquad=
\sum_{{{\alpha\in\N^n}\atop{|\alpha|=k}}\atop{1\le i\le n}}
\frac{k!\,(\alpha_i+1)}{(\alpha+e_i)!}\,D^{2(\alpha+e_i)}u=
\sum_{{{\beta\in\N^n}\atop{|\beta|=k+1}}\atop{1\le i\le n}}
\frac{k!\,\beta_i}{\beta !}\,D^{2\beta}u=
\sum_{{{\beta\in\N^n}\atop{|\beta|=k+1}}}
\frac{k!\,|\beta|}{\beta!}\,D^{2\beta}u\\&&\qquad\qquad\qquad=
\sum_{{{\beta\in\N^n}\atop{|\beta|=k+1}}}
\frac{k!\,(k+1)}{\beta!}\,D^{2\beta}u=
\sum_{{{\beta\in\N^n}\atop{|\beta|=k+1}}}
\frac{(k+1)!}{\beta!}\,D^{2\beta}u,
\end{eqnarray*}
which completes the inductive step and provides the proof of~\eqref{MULTOLPALPA}.

Now we point out that, by odd symmetry
of the function~$x_i\mapsto x_i^{\alpha_i}$,
if~$\alpha_i$ is an odd integer then~$ \int_{B_\rho} x^{\alpha}\,dx=0$ 
and consequently
\begin{equation}\label{MULTOLPALPA2} \int_{B_\rho} x^{\alpha}\,dx=0\qquad{\mbox{unless $\alpha=2\beta$ for some~$\beta\in\N^n$.}}
\end{equation}
In the same way,
\begin{equation}\label{MULTOLPALPA3} \int_{\partial B_\rho} x^{\alpha}\,d{\mathcal{H}}^{n-1}_x=0\qquad{\mbox{unless $\alpha=2\beta$ for some~$\beta\in\N^n$.}}
\end{equation}

It is now convenient to perform some specific calculation
related to spherical integrals. 
So as to achieve this goal,
recalling the definition~\eqref{EULEROGA} of
Euler Gamma Function
and
using the substitution~$\tau:=t^2$,
we point out that, for every~$b>-1/2$,
\begin{equation}\label{KSM-ONSN Dcakgau934} \int_0^{+\infty} t^{2b} e^{-t^2}\,dt=
\frac12\,\int_0^{+\infty} \tau^{\frac{2b-1}2} e^{-\tau} \,d\tau
=\frac{1}2\,\Gamma\left({\frac{2b+1}2}\right).\end{equation}
Consequently, for every~$\alpha\in\N^n$,
\begin{equation}\label{KSM-ONSN Dcakgau9342}
\int_{\R^n} x^{2\alpha}e^{-|x|^2}\,dx=\prod_{i=1}^n
\int_{\R} x_i^{2\alpha_i}e^{-x_i^2}\,dx_i=2^n\,
\prod_{i=1}^n
\int_{0}^{+\infty} x_i^{2\alpha_i}e^{-x_i^2}\,dx_i=
\prod_{i=1}^n\Gamma\left({\frac{2\alpha_i+1}2}\right).
\end{equation}
On the other hand, using polar coordinates and once again~\eqref{KSM-ONSN Dcakgau934},
\begin{equation*}
\int_{\R^n} x^{2\alpha}e^{-|x|^2}\,dx=\int_{\partial B_1}
\left(\int_0^{+\infty} r^{2|\alpha|+n-1}
\omega^{2\alpha}e^{-r^2}\,dr
\right)\,d{\mathcal{H}}^{n-1}_\omega=
\frac{1}2\,\Gamma\left({\frac{2|\alpha|+n}2}\right)\,
\int_{\partial B_1}
\omega^{2\alpha}\,d{\mathcal{H}}^{n-1}_\omega.
\end{equation*}
Comparing this with~\eqref{KSM-ONSN Dcakgau9342},
we deduce that
\begin{equation}\label{KSMX:att89i45tgsbad2glifdf0athbnergr}
\int_{\partial B_1}
\omega^{2\alpha}\,d{\mathcal{H}}^{n-1}_\omega
=\frac{\displaystyle2\prod_{i=1}^n\Gamma\left({\frac{2\alpha_i+1}2}\right)}{
\displaystyle\Gamma\left({\frac{2|\alpha|+n}2}\right)
}.
\end{equation}

It is also instructive to recall that, integrating by parts,
\begin{equation}\label{GAMMAPIU1} \Gamma (z+1)=\int _{0}^{+\infty }t^{z}e^{-t}\,dt=-
\int _{0}^{+\infty }t^{z}\,\frac{d}{dt}e^{-t}\,dt
=\int _{0}^{+\infty }zt^{z-1}e^{-t}\,dt=z\,\Gamma(z),\end{equation}
from which it follows by induction that
$$ \Gamma(j+1)=j!\qquad{\mbox{for every }}j\in\N.$$
Furthermore, making use of~\eqref{KSM-ONSN Dcakgau934}  with~$b:=0$,
\begin{equation}\label{ZKMdgamzo34}
\Gamma\left(\frac12\right)=2\int_0^{+\infty}e^{-t^2}\,dt=
\int_{-\infty}^{+\infty}e^{-t^2}\,dt=\sqrt{\pi}.
\end{equation}

Now we claim that, for every~$j\in\N$,
\begin{equation}\label{9ijn8uh9ijGSDNlsodgue0olg:SL467A}
\Gamma\left(j+1+\frac12\right)=\frac{\sqrt\pi\,(2j+1)!}{2^{2j+1} \,j!}
\end{equation}
This statement is actually a particular case
of the Legendre Duplication Formula,
or of the 
Gau{\ss} Multiplication Formula, but we provide a direct proof of~\eqref{9ijn8uh9ijGSDNlsodgue0olg:SL467A}
for the facility of the reader. For this, we argue by induction over~$j$.
We remark that
\begin{equation}\label{GAMMA3mez} \Gamma\left(1+\frac12\right)=\frac12\,\Gamma\left(\frac12\right)
=
\frac{\sqrt\pi}{2}
\end{equation}
thanks to~\eqref{GAMMAPIU1} and~\eqref{ZKMdgamzo34},
and this gives the claim in~\eqref{9ijn8uh9ijGSDNlsodgue0olg:SL467A}
when~$j=0$.

Suppose now that~\eqref{9ijn8uh9ijGSDNlsodgue0olg:SL467A}
holds true for~$j$. Then, we use again~\eqref{GAMMAPIU1} to obtain that
\begin{eqnarray*}
&&\Gamma\left(j+2+\frac12\right)=\left(j+1+\frac12\right)\,
\Gamma\left(j+1+\frac12\right)=
\left(j+1+\frac12\right)\,\frac{\sqrt\pi\,(2j+1)!}{2^{2j+1} \,j!}\\&&\qquad=
\left(2j+3\right)\,\frac{\sqrt\pi\,(2j+1)!}{2^{2j+2} \,j!}=
\left(2j+3\right)\left(2j+2\right)\,\frac{\sqrt\pi\,(2j+1)!}{2^{2j+3}\,(j+1) \,j!}
=\frac{\sqrt\pi\,(2j+3)!}{2^{2j+3} \,(j+1)!},
\end{eqnarray*}
thus completing the inductive step and establishing~\eqref{9ijn8uh9ijGSDNlsodgue0olg:SL467A}.

Now we claim that
\begin{equation}\label{GRfasvdDPqu}
\Gamma\left({\frac{2\alpha_i+1}2}\right)= 2^{1-2\alpha_i}\,\varsigma(\alpha_i)\,\sqrt\pi\prod_{j=0}^{\alpha_i-1}(\alpha_i
+j)
,\end{equation}
where $$ \varsigma(\alpha_i):=\begin{dcases}
1 & {\mbox{ if }}\alpha_i\ne0,\\ \\
\frac12& {\mbox{ if }}\alpha_i=0.
\end{dcases}$$
Indeed, by~\eqref{9ijn8uh9ijGSDNlsodgue0olg:SL467A},
if~$\alpha_i\ne0$ then
$$ \Gamma\left({\frac{2\alpha_i+1}2}\right)= 
\frac{\sqrt\pi\,(2\alpha_i-1)!}{2^{2\alpha_i-1} \,(\alpha_i-1)!}=2^{1-2\alpha_i}\,\sqrt\pi\prod_{j=0}^{\alpha_i-1}(\alpha_i
+j).$$
This proves~\eqref{GRfasvdDPqu} when~$\alpha_i\ne0$.

If instead~$\alpha_i=0$, we recall that, as usual,
the empty product is~$1$ by convention, and we use~\eqref{ZKMdgamzo34} to see that in this case
\begin{eqnarray*}
&&\Gamma\left({\frac{2\alpha_i+1}2}\right)-2^{1-2\alpha_i}\, \varsigma(\alpha_i)\,\sqrt\pi\prod_{j=0}^{\alpha_i-1}(\alpha_i
+j)=\Gamma\left({\frac{1}2}\right)- 2^{0}\,\sqrt\pi=0.
\end{eqnarray*}
The proof of~\eqref{GRfasvdDPqu} is thereby complete.

We also observe that
\begin{equation}\label{Sm90KmSweovj-mkthc9klps9okiduih1}
\prod_{i=1}^n\varsigma(\alpha_i)=
\frac1{2^{\theta(\alpha)}},
\end{equation}
where~$\theta(\alpha)$ is the number of vanishing components of the multi-index~$\alpha$.

Thus, by~\eqref{GRfasvdDPqu},
$$ 2\prod_{i=1}^n\Gamma\left({\frac{2\alpha_i+1}2}\right)=
2^{n+1-2|\alpha|-\theta(\alpha)}\,\pi^{\frac{n}{2}}
\displaystyle\prod_{i=1}^n\prod_{j=0}^{\alpha_i-1}(\alpha_i
+j)$$
and therefore,
in light of~\eqref{KSMX:att89i45tgsbad2glifdf0athbnergr},
\begin{equation}\label{KSMX:att89i45tgsbad2glifdf0athbnergr-2}
\int_{\partial B_1}
\omega^{2\alpha}\,d{\mathcal{H}}^{n-1}_\omega
=\frac{2^{n+1-2|\alpha| -\theta(\alpha)}\,\pi^{\frac{n}{2}}
\displaystyle\prod_{i=1}^n\prod_{j=0}^{\alpha_i-1}(\alpha_i
+j)
}{\displaystyle\Gamma\left({\frac{2|\alpha|+n}2}\right)}
.\end{equation}

Consequently,
\begin{equation}\label{8-9-46-249-2394-84rjf-238rikf223-4yci-239rij7n-12era}
\begin{split} &\int_{B_1} x^{2\alpha}\,dx=
\int_{\partial B_1}\left(\int_0^1 r^{2|\alpha|+n-1}\omega^{2\alpha}\,dr\right)\,
d{\mathcal{H}}^{n-1}_\omega=
\frac{1}{2|\alpha|+n}\int_{\partial B_1}
\omega^{2\alpha}\,
d{\mathcal{H}}^{n-1}_\omega \\&\qquad\qquad\qquad\qquad
=
\frac{2^{n+1-2|\alpha| -\theta(\alpha)}\,\pi^{\frac{n}{2}}
\displaystyle\prod_{i=1}^n\prod_{j=0}^{\alpha_i-1}(\alpha_i
+j)}{\displaystyle(2|\alpha|+n)\,\Gamma\left({\frac{2|\alpha|+n}2}\right)}
.\end{split}
\end{equation}

We also note that
\begin{equation}\label{Sm90KmSweovj-mkthc9klps9okiduih2}
(2\beta_i)!=
2\varsigma(\beta_i)\,\beta_i!\;\prod_{j=0}^{\beta_i-1}(\beta_i+j).
\end{equation}
Indeed, if~$\beta_i\ne0$ then
\begin{eqnarray*}&& (2\beta_i)!=
2\beta_i\;(2\beta_i-1)!=2\beta_i
\prod_{j=1}^{2\beta_i-1} j=2\beta_i
\left(\prod_{k=1}^{\beta_i-1}k\right)
\left(\prod_{k=\beta_i}^{2\beta_i-1}k\right)\\&&\qquad\qquad\quad=2\beta_i
(\beta_i-1)!\;
\left(\prod_{j=0}^{\beta_i-1}(\beta_i+ j)\right)=2\varsigma(\beta_i)\,\beta_i!\;\prod_{j=0}^{\beta_i-1}(\beta_i+j)
,\end{eqnarray*}
while if~$\beta_i=0$ then
\begin{eqnarray*}&& (2\beta_i)!-2\varsigma(\beta_i)\,\beta_i!\;\prod_{j=0}^{\beta_i-1}(\beta_i+j)=
1-2\varsigma(\beta_i)=0.
\end{eqnarray*}
These observations establish~\eqref{Sm90KmSweovj-mkthc9klps9okiduih2}.

{F}rom~\eqref{Sm90KmSweovj-mkthc9klps9okiduih1} and~\eqref{Sm90KmSweovj-mkthc9klps9okiduih2} we arrive at
\begin{equation}\label{Sm90KmSweovj-mkthc9klps9okiduih3}\begin{split}&
\frac{2^{n-\theta(\beta)}\,
\displaystyle\prod_{i=1}^n\prod_{j=0}^{\beta_i-1}(\beta_i
+j)
}{(2\beta)!}=\frac{
\displaystyle\prod_{i=1}^n\left(2\varsigma(\beta_i)\prod_{j=0}^{\beta_i-1}(\beta_i
+j)\right)
}{(2\beta)!}\\&\qquad=
\prod_{i=1}^n
\frac{
\displaystyle2\varsigma(\beta_i)\prod_{j=0}^{\beta_i-1}(\beta_i
+j)
}{(2\beta_i)!}
=\prod_{i=1}^n\frac1{\beta_i!}=
\frac{1}{ \beta!}.\end{split}
\end{equation}

We now prove~\eqref{PALLAPIZZE1}. For this, up to a translation,
we can reduce to the case~$x_0=0$. We exploit the Taylor expansion
\begin{equation}\label{MULTOLPALPA4} u(x)=\sum_{{\alpha\in\N^n}\atop{|\alpha|\le 2N}}
\frac{D^\alpha u(0)}{\alpha!}\,x^\alpha+o(x^{2N})\end{equation}
and we average over~$B_\rho$, thus finding that, for small~$\rho$,
\begin{eqnarray*}
\fint_{B_\rho}u(x)\,dx&=&\sum_{{\alpha\in\N^n}\atop{|\alpha|\le 2N}}
\frac{D^\alpha u(0)}{\alpha!}\,\fint_{B_\rho}x^\alpha\,dx+o(\rho^{2N})\\
&=&\sum_{{\beta\in\N^n}\atop{|\beta|\le N}}
\frac{D^{2\beta} u(0)}{(2\beta)!}\,\fint_{B_\rho}x^{2\beta}\,dx+o(\rho^{2N})\\
&=&\sum_{{\beta\in\N^n}\atop{|\beta|\le N}}
\frac{D^{2\beta} u(0)\,\rho^{2|\beta|}}{(2\beta)!}\,\fint_{B_1}y^{2\beta}\,dy
+o(\rho^{2N})
\\&=&\sum_{k=0}^N\sum_{{\beta\in\N^n}\atop{|\beta|=k}}
\frac{D^{2\beta} u(0)\,\rho^{2k}}{(2\beta)!}\,\fint_{B_1}y^{2\beta}\,dy+o(\rho^{2N})
\\&=&
\sum_{k=0}^N\sum_{{\beta\in\N^n}\atop{|\beta|=k}}
\frac{D^{2\beta} u(0)\,\rho^{2k}}{(2\beta)!\,|B_1|}\,\frac{2^{n+1-2k -\theta(\beta)}\,
\pi^{\frac{n}{2}}\displaystyle\prod_{i=1}^n\prod_{j=0}^{\beta_i-1}(\beta_i
+j)
}{\displaystyle(2k+n)\,\Gamma\left({\frac{2k+n}2}\right)}
+o(\rho^{2N})\\&=&
\sum_{k=0}^N\sum_{{\beta\in\N^n}\atop{|\beta|=k}}
\frac{D^{2\beta} u(0)\,\rho^{2k}}{ \beta!\,|B_1|}\,
\frac{2^{1-2k}\,\pi^{\frac{n}{2}}
}{\displaystyle(2k+n)\,\Gamma\left({\frac{2k+n}2}\right)}
+o(\rho^{2N})
\\&=&
\sum_{k=0}^N
\frac{\Delta^k u(0)\,\rho^{2k}}{ |B_1|}\,
\frac{2^{1-2k}\,\pi^{\frac{n }{2}}
}{k!\,\displaystyle(2k+n)\,\Gamma\left({\frac{2k+n}2}\right)}
+o(\rho^{2N}),\end{eqnarray*}
thanks to~\eqref{MULTOLPALPA}, \eqref{MULTOLPALPA2},
\eqref{8-9-46-249-2394-84rjf-238rikf223-4yci-239rij7n-12era}
and~\eqref{Sm90KmSweovj-mkthc9klps9okiduih2}.

Hence, since, owing to~\eqref{B1},
\eqref{KSMX:att89i45tgsbad2glifdf0athbnergr}
(used here with~$\alpha:=0$) and~\eqref{ZKMdgamzo34},
\begin{equation}\label{MULTOLPALPA280kdsmc8724yrth0923}
\Gamma\left({\frac{n}2}\right)
=\frac{2 }{{\mathcal{H}}^{n-1}(\partial B_1)}\,\left(
\Gamma\left({\frac{1}2}\right)\right)^n=\frac{2\pi^{\frac{n}2} }{{\mathcal{H}}^{n-1}(\partial B_1)}
=\frac{2\pi^{\frac{n}2}}{n\,| B_1|},\end{equation}
we obtain that
\begin{equation*}
\begin{split}
\fint_{B_\rho}u(x)\,dx\,&=\,n\,
\Gamma\left({\frac{n}2}\right)
\sum_{k=0}^N
\frac{ 
\Delta^k u(0)\,\rho^{2k}}{2^{2k}\,k!\,\displaystyle(2k+n)\,\Gamma\left({\frac{2k+n}2}\right)}
+o(\rho^{2N})\\
&=\,n\,
\Gamma\left({\frac{n}2}\right)
\sum_{k=0}^N
\frac{\Delta^k u(0)\,\rho^{2k}
}{2^{2k+1}\,k!\,\displaystyle\Gamma\left({\frac{n}2}+k+1\right)}
+o(\rho^{2N}),
\end{split}
\end{equation*}
where~\eqref{GAMMAPIU1} has been used again in the last line.
This completes the proof of~\eqref{PALLAPIZZE1}.

Now, to establish~\eqref{PALLAPIZZE2},
we exploit~\eqref{MULTOLPALPA}, \eqref{MULTOLPALPA3}, \eqref{KSMX:att89i45tgsbad2glifdf0athbnergr-2},
\eqref{Sm90KmSweovj-mkthc9klps9okiduih2} and~\eqref{MULTOLPALPA4}, and we see that
\begin{eqnarray*}
\fint_{\partial B_\rho}u(x)\,d{\mathcal{H}}^{n-1}_x
&=&\sum_{{\alpha\in\N^n}\atop{|\alpha|\le 2N}}
\frac{D^\alpha u(0)}{\alpha!}\,\fint_{\partial B_\rho}x^\alpha\,d{\mathcal{H}}^{n-1}_x+o(\rho^{2N})\\
&=&\sum_{{\beta\in\N^n}\atop{|\beta|\le N}}
\frac{D^{2\beta} u(0)}{(2\beta)!}\,\fint_{\partial B_\rho}x^{2\beta}\,
d{\mathcal{H}}^{n-1}_x+o(\rho^{2N})\\
&=&\sum_{{\beta\in\N^n}\atop{|\beta|\le N}}
\frac{D^{2\beta} u(0)\,\rho^{2|\beta|}}{(2\beta)!}\,\fint_{\partial
B_1}y^{2\beta}\,d{\mathcal{H}}^{n-1}_y
+o(\rho^{2N})
\\&=&\sum_{k=0}^N
\sum_{{\beta\in\N^n}\atop{|\beta|=k}}
\frac{D^{2\beta} u(0)\,\rho^{2k}}{(2\beta)!}\,\fint_{\partial
B_1}y^{2\beta}\,d{\mathcal{H}}^{n-1}_y
+o(\rho^{2N})\\
&=&\sum_{k=0}^N
\sum_{{\beta\in\N^n}\atop{|\beta|=k}}
\frac{D^{2\beta} u(0)\,\rho^{2k}}{(2\beta)!\,{\mathcal{H}}^{n-1}(\partial B_1)}\,
\frac{2^{n+1-2k-\theta(\beta)}\,\pi^{\frac{n}{2}}
\displaystyle\prod_{i=1}^n\prod_{j=0}^{\beta_i-1}(\beta_i
+j)
}{\displaystyle\Gamma\left({\frac{2k+n}2}\right)}
+o(\rho^{2N})\\&=&\sum_{k=0}^N
\sum_{{\beta\in\N^n}\atop{|\beta|=k}}
\frac{D^{2\beta} u(0)\,\rho^{2k}}{ \beta!\,{\mathcal{H}}^{n-1}(\partial B_1)}\,
\frac{2^{1-2k}\,\pi^{\frac{n}{2}}
}{\displaystyle\Gamma\left({\frac{2k+n}2}\right)}
+o(\rho^{2N})\\&=&
\sum_{k=0}^N
\frac{\Delta^k u(0)\,\rho^{2k}}{k!\,{\mathcal{H}}^{n-1}(\partial B_1)}\,
\frac{2^{1-2k}\,\pi^{\frac{n}{2}}
}{\displaystyle\Gamma\left({\frac{2k+n}2}\right)}
+o(\rho^{2N})
.\end{eqnarray*}
This, together with~\eqref{MULTOLPALPA280kdsmc8724yrth0923}, proves~\eqref{PALLAPIZZE2}, as desired
(alternatively, one could have proved one between~\eqref{PALLAPIZZE1}
and~\eqref{PALLAPIZZE2} with an explicit expression for the remainder
and obtained the other by either integration or differentiation
in~$\rho$).
\end{proof}

We stress that Theorems~\ref{KAHAR2} and~\ref{KAHAR2-SPHE}
are particular cases of~\eqref{PALLAPIZZE1}
and~\eqref{PALLAPIZZE2}
respectively, just corresponding to the case~$N:=1$.
\medskip

In spite of its simple flavor, Theorem~\ref{KAHAR2} reveals one of the fundamental
features of the Laplacian, which will also play an important role in the Mean Value Formula
that will be described in Theorem~\ref{KAHAR}. Also, it immediately leads to
the fact that the Laplace operator is invariant under rotation, namely:

\begin{corollary}\label{KSMD:ROTAGSZOKA}
Let~${\mathcal{R}}:\R^n\to\R^n$ be a rotation
and~$u\in C^2(\R^n)$. Let~$u_{\mathcal{R}}(x):=u({\mathcal{R}}x)$.

Then,
$$ \Delta u_{\mathcal{R}}(x)=\Delta u({\mathcal{R}}x).$$
\end{corollary}

\begin{proof} Let~$x_0\in\R^n$.
By Theorem~\ref{KAHAR2}, using the change of variable~$y:={\mathcal{R}}(x-x_0)$
and the notation~$v(x):=u(x+{\mathcal{R}}x_0)$, and exploiting the translation invariance in~\eqref{TRAL},
we have
\begin{eqnarray*}
&& \frac{1}{2(n+2)}\,\Delta u_{\mathcal{R}}(x_0)=
\lim_{\rho\searrow0} \frac1{\rho^2}\left(
\fint_{B_\rho(x_0)} u_{\mathcal{R}}(x)\,dx-u_{\mathcal{R}}(x_0)\right)=
\lim_{\rho\searrow0} \frac1{\rho^2}\left(
\fint_{B_\rho(x_0)} u({\mathcal{R}}x)\,dx-u({\mathcal{R}}x_0)\right)\\&&\qquad=
\lim_{\rho\searrow0} \frac1{\rho^2}\left(
\fint_{B_\rho} u(y+{\mathcal{R}}x_0)\,dy-u({\mathcal{R}}x_0)\right)=
\lim_{\rho\searrow0} \frac1{\rho^2}\left(
\fint_{B_\rho} v(y)\,dy-v(0)\right)\\&&\qquad=\frac{1}{2(n+2)}\,\Delta v(0)=
\frac{1}{2(n+2)}\,\Delta u({\mathcal{R}}x_0)
.\end{eqnarray*}
This proves the desired result (for a proof not relying on
Theorem~\ref{KAHAR2} but rather on a direct computation in matrix form
see e.g.~\cite[page 172]{MR3938717}).
\end{proof}

Of course, while the integral characterization of the Laplacian presented in Theorem~\ref{KAHAR2}
is conceptually very useful and reveals a deep geometric structure of the operator,
the explicit differential structure in~\eqref{LA:DIFF} is often
simpler to exploit for explicit calculations.
As an example, we recall the following Bochner Identity\index{Bochner Identity}:

\begin{lemma}\label{BOCHIFRENTYHAS}
For a given~$C^3$ function~$u$,
$$ \Delta \left( \frac{|\nabla u|^{2}}{2}\right)=
\nabla( \Delta u)\cdot\nabla u
+|D^2u|^{2}
,$$
where
\begin{equation}\label{D22H}
|D^2u|^2:=\sum_{i,j=1}^n (\partial_{ij} u)^2
.\end{equation}
\end{lemma}

\begin{proof} By a direct computation,
\begin{equation*}\begin{split}&
\Delta \left( \frac{|\nabla u|^{2}}{2} \right)=\frac12\,\sum_{i,j=1}^n
\partial_{ii} (\partial_j u)^2=
\sum_{i,j=1}^n
\partial_i (\partial_j u\,\partial_{ij}u)\\&\qquad
=\sum_{i,j=1}^n\Big(
(\partial_{ij}u)^2
+\partial_j u\,\partial_{iij}u)
\Big)=|D^2u|^2+\sum_{j=1}^n\partial_ju\,\partial_j(\Delta u)
.\qedhere\end{split}\end{equation*}
\end{proof}

\chapter{The Laplace operator and harmonic functions}\label{CHA-2}

This chapter is devoted to the analysis of the Laplace operator
and of the functions that lie in its kernel.

\section{The Laplacian and the Mean Value Formula}

Among the several properties that a given function may possess,\index{Mean Value Formula}
a very relevant one is ``harmonicity'', corresponding to the vanishing of
the trace of the Hessian matrix (in particular, these functions
are ``saddle-looking'' with respect to their tangent planes at every point).
The precise setting that we consider is the following:

\begin{definition}\label{KSMD:PSDLKNIE0sdf}
Given an open set~$\Omega\subseteq\R^n$ and
a function~$u\in C^2(\Omega)$, we say that~$u$ is harmonic in~$\Omega$
if~$\Delta u(x)=0$ for every~$x\in\Omega$.
\end{definition}

For example, constant and linear functions are harmonic in all~$\R^n$.
Also, the functions~$u:\R^2\to\R$ given by~$u(x_1,x_2)=x_1 x_2$,
$u(x_1,x_2)=x_1^2- x_2^2$ and~$u(x_1,x_2)=e^{x_1}\sin x_2$ are harmonic.

Other examples of harmonic functions in domains of~$\R^2$ can be obtained via complex analysis,
identifying~$(x,y)\in\R^2$ with~$z=x+iy\in \cOMPL$,
since
\begin{equation}\label{KS:RUDVA0o1keLLetyh439}
{\mbox{the real and imaginary parts of holomorphic functions are harmonic,}}\end{equation}
see~\cite[Chapter 11]{MR0210528}.
In particular, for every~$j\in\N$, using the notation~$r=|z|=\sqrt{x^2+y^2}$
and~$z=|z| e^{i\vartheta}=re^{i\vartheta}$,
the functions
\begin{equation}\label{REj} \Re z^j =r^j \cos(j\vartheta)\qquad{\mbox{and}}\qquad
\Im z^j =r^j \sin(j\vartheta)\end{equation}
are harmonic in all~$\R^2$, and so are the functions
$$ \Re e^z =e^x \cos y\qquad{\mbox{and}}\qquad
\Im e^z =e^x \sin y.$$
In addition, given~$\alpha>0$, the function~$z\mapsto z^\alpha:=r^\alpha e^{i\alpha\vartheta}$ is 
well-defined and holomorphic in~$\vartheta\in(-\pi,\pi)$, hence
$$ \Re z^\alpha =r^\alpha \cos(\alpha\vartheta)\qquad{\mbox{and}}\qquad
\Im z^\alpha =r^\alpha \sin(\alpha\vartheta)$$
are harmonic in~$\R^2\setminus\ell$, being~$\ell:=(-\infty,0]\times\{0\}$.\medskip

One of the most striking properties of harmonic functions
is that their value at any point is precisely equal to the average
of the values around such point.
In this sense,
the values attained by harmonic functions happen to be ``perfectly balanced'',
according to the following result:

\begin{theorem}\label{KAHAR}
Given an open set~$\Omega\subseteq\R^n$ and
a function~$u\in L^1_{\rm loc}(\Omega)$,
the following conditions are equivalent:
\begin{itemize}
\item[(i).] The function~$u$
belongs to~$C^2(\Omega)$ and
is harmonic in $\Omega$.
\item[(ii).] For almost
every~$x_0\in\Omega$ and almost
every~$r>0$ such that~$B_r(x_0)\Subset\Omega$,
we have that
$$ u(x_0)=\fint_{\partial B_r(x_0)} u(x)\,d{\mathcal{H}}^{n-1}_x.$$
\item[(iii).] For almost
every~$x_0\in\Omega$ and almost every~$r>0$ such that~$B_r(x_0)\Subset\Omega$,
we have that
$$ u(x_0)=\fint_{B_r(x_0)} u(x)\,dx.$$
\end{itemize}

Additionally, if~$u$ satisfies any of the equivalent
conditions~(i), (ii) or~(iii), then\footnote{The equivalence
between conditions~(i), (ii) and~(iii)
highlights an interesting regularizing property of the Laplace
operator, since locally integrable functions
satisfying either~(ii) or~(iii) turn out to belong to~$C^2(\Omega)$
hence their Laplacian can be computed pointwise and it is equal to zero.

A similar regularizing effect will be highlighted by the forthcoming
Lemma~\ref{WEYL}.

Also, the last statement
in Theorem~\ref{KAHAR} concerning the smoothness of~$u$
gives that harmonic functions are automatically~$C^\infty(\Omega)$:
this statement will be strengthened in Theorem~\ref{KMS:HAN0-0},
where we will show in fact that harmonic functions
are real analytic.}
it belongs to~$C^\infty(\Omega)$.
\end{theorem}

\begin{proof} We start by showing that
\begin{equation}\label{MEG9ifjliio0oO1}
{\mbox{if~$u$ satisfies either~(ii) or~(iii),
then~$u\in C^\infty(\Omega)$,}}
\end{equation}
up to redefining~$u$ in a set of null Lebesgue measure.
To this end, we use a mollification argument.
We\footnote{As customary, here and in the following,
the subscript~$0$ in~$C^\infty_0$ means ``with compact support in''.} take~$\tau\in C^\infty_0(B_1,\,[0,+\infty))$ to
be radially symmetric and such that~$\int_{B_1}\tau(x)\,dx=1$.
Given~$\eta>0$, we let~$\tau_\eta(x):=\frac1{\eta^n}
\tau\left(\frac{x}{\eta}\right)$
and define~$u_\eta:=u*\tau_\eta$.
We pick a point~$\overline{x}\in\Omega$
and~$R>0$
such that~$B_{2R}(\overline{x})\Subset\Omega$, and we show
that, when~$\eta\in(0,R)$,
\begin{equation}\label{MEG9ifjliio0oO2}
{\mbox{if~$u$ satisfies either~(ii) or~(iii),
then~$u=u_\eta$ a.e. in~$B_{R}(\overline{x})$.}}
\end{equation}
Indeed, if~$u$ satisfies~(ii), 
for each~$x\in B_{R}(\overline{x})$ and~$r\in(0,\eta]$
we have that~$B_r(x)\subseteq
B_\eta(x)\subseteq B_{R+\eta}(\overline{x})\Subset\Omega$.
Hence,
using polar coordinates (see e.g.~\cite[Theorem~3.12]{MR3409135}),
for almost any~$x\in B_{R}(\overline{x})$,
\begin{eqnarray*}
u_\eta(x)&=&\int_{B_\eta(x)} \tau_\eta(x-y)\,u(y)\,dy
\\&=&\int_0^\eta \left[\int_{\partial B_r(x)}
\tau_\eta(x-\omega)\,u(\omega)\,d{\mathcal{H}}^{n-1}_\omega
\right]\,dr\\&=&\int_0^\eta \left[\int_{\partial B_r(x)}
\tau_\eta(re_1)\,u(\omega)\,d{\mathcal{H}}^{n-1}_\omega
\right]\,dr\\&=&\int_0^\eta \left[{\tau_\eta(re_1)}\,{{\mathcal{H}}^{n-1}(\partial B_r)}\,
\fint_{\partial B_r(x)}
u(\omega)\,d{\mathcal{H}}^{n-1}_\omega
\right]\,dr\\&=&u(x)\,
\int_0^\eta {\tau_\eta(re_1)}\,{{\mathcal{H}}^{n-1}(\partial B_r)}\,\,dr.
\end{eqnarray*}
Accordingly, since
\begin{eqnarray*}&& 1=\int_{B_\eta}\tau_\eta(y)\,dy
=\int_0^\eta \left[\int_{\partial B_r}
\tau_\eta(\omega)\,d{\mathcal{H}}^{n-1}_\omega
\right]\,dr\\&&\qquad\qquad=
\int_0^\eta \left[\int_{\partial B_r}
\tau_\eta(re_1)\,d{\mathcal{H}}^{n-1}_\omega
\right]\,dr=\int_0^\eta {\tau_\eta(re_1)}\,{{\mathcal{H}}^{n-1}(\partial B_r)}\,\,dr
,\end{eqnarray*}
we gather that~$u_\eta(x)=u(x)$ for almost all~$x\in
B_{R}(\overline{x})$
and, as a result,~\eqref{MEG9ifjliio0oO2}
holds true when condition~(ii) is satisfied.

Also, if~$u$ fulfills condition~(iii), then by polar coordinates
and~\eqref{B1},
for almost every~$x_0\in\Omega$ and almost
every~$r>0$ such that~$B_r(x_0)\Subset\Omega$,
we have that
\begin{eqnarray*}&&
\fint_{\partial B_r(x_0)} u(y)\,dy=\frac1{r^{n-1}\mathcal{H}^{n-1}(\partial B_1)}
\,\frac{d}{dr}\int_{ B_r(x_0)} u(y)\,dy
=\frac1{n r^{n-1}}
\,\frac{d}{dr}\left( r^n\fint_{ B_r(x_0)} u(y)\,dy\right)\\&&\qquad=
\frac{u(x_0)}{n r^{n-1}}
\,\frac{d}{dr}\left( r^n\right)=u(x_0),
\end{eqnarray*}
which gives~(ii).
This reduces us to the previous case, and therefore the proof
of~\eqref{MEG9ifjliio0oO2} is complete.

In turn, we have that~\eqref{MEG9ifjliio0oO2}
entails~\eqref{MEG9ifjliio0oO1}, as desired.
\medskip

Now we will show that~(i) implies~(ii) which implies~(iii) which implies~(i).\medskip

Let us assume that~(i) holds true and let~$B_r(x_0)\Subset\Omega$.
Then, $u$ is harmonic in~$B_\rho(x_0)$ for every~$\rho\in(0,r]$. Accordingly, by
the application of the Divergence Theorem given in equation~\eqref{1111DIV0934},
$$ 0=
\int_{B_\rho(x_0)} \Delta u(x)\,dx=\int_{\partial B_\rho(x_0) } \nabla u(x)\cdot \nu(x)\,d{\mathcal{H}}^{n-1}_x
=\frac1\rho\,\int_{\partial B_\rho(x_0) } \nabla u(x)\cdot (x-x_0)\,d{\mathcal{H}}^{n-1}_x.$$
On the other hand,
\begin{eqnarray*}&& \frac{d}{d\rho} \left(\fint_{\partial B_\rho(x_0)} u(x)\,d{\mathcal{H}}^{n-1}_x\right)
=\frac{d}{d\rho} \left(\fint_{\partial B_1} u(x_0+\rho\omega)\,d{\mathcal{H}}^{n-1}_\omega\right)\\&&\qquad=
\fint_{\partial B_1} \nabla u(x_0+\rho\omega)\cdot\omega\,d{\mathcal{H}}^{n-1}_\omega=\frac1\rho\,\fint_{\partial B_\rho(x_0) } \nabla u(x)\cdot (x-x_0)\,d{\mathcal{H}}^{n-1}_x.\end{eqnarray*}
These observations entail that, for every~$\rho\in(0,r]$,
$$ \frac{d}{d\rho} \left(\fint_{\partial B_\rho(x_0)} u(x)\,d{\mathcal{H}}^{n-1}_x\right)=0,$$
whence
$$ \fint_{\partial B_\rho(x_0)} u(x)\,d{\mathcal{H}}^{n-1}_x {\mbox{ is constant for all }}\rho\in(0,r].$$
In particular,
$$ \fint_{\partial B_r(x_0)} u(x)\,d{\mathcal{H}}^{n-1}_x=\lim_{\rho\searrow0}
\fint_{\partial B_\rho(x_0)} u(x)\,d{\mathcal{H}}^{n-1}_x=u(x_0)$$
and this shows that~(ii) holds true.\medskip

The fact that~(ii) implies~(iii) is a consequence of polar coordinates. Indeed, if~(ii) is satisfied, then
\begin{eqnarray*}&&
\fint_{B_r(x_0)} u(x)\,dx=\frac{1}{|B_r|} \int_0^r \left(
\int_{\partial B_\rho(x_0)} u(x)\,d{\mathcal{H}}^{n-1}_x\right)\,d\rho
\\&&\qquad\qquad=\frac{{\mathcal{H}}^{n-1}(\partial B_1)}{|B_1| \,r^n} \int_0^r \left( 
\rho^{n-1}
\fint_{\partial B_\rho(x_0)} u(x)\,d{\mathcal{H}}^{n-1}_x\right)\,d\rho\\&&
\qquad\qquad=
\frac{{\mathcal{H}}^{n-1}(\partial B_1)\,u(x_0)}{|B_1| \,r^n} \int_0^r
\rho^{n-1}\,d\rho=\frac{{\mathcal{H}}^{n-1}(\partial B_1)\,u(x_0)}{|B_1| \,n}.
\end{eqnarray*}
This and~\eqref{B1} yield that
$$ \fint_{B_r(x_0)} u(x)\,dx=u(x_0),$$
that is~(iii).
\medskip

Let us now suppose that~(iii) holds true. Then, by~\eqref{MEG9ifjliio0oO1}
and
Theorem~\ref{KAHAR2}, for every~$x_0\in\Omega$,
$$0=
\lim_{r\searrow0} \frac1{r^2}\left(
\fint_{B_r(x_0)} u(x)\,dx-u(x_0)\right)=\frac{1}{2(n+2)}\,\Delta u(x_0),
$$
thus showing the validity of~(i).
\end{proof}

For a comprehensive survey on
mean value properties of harmonic and closely related functions, see~\cite{MR1321628}.
A simple byproduct of the Mean Value Formula in Theorem~\ref{KAHAR}
is the following interesting geometric observation:

\begin{corollary}\label{e5vbhfueweuwiuet0000} Let~$n\ge2$.
A harmonic function does not possess isolated zeroes.
\end{corollary}

\begin{proof} Let~$u$ be harmonic in some open set~$\Omega\subseteq\R^n$
and suppose that~$u(x_0)=0$. Arguing for a contradiction,
we suppose that there exists~$r>0$ such that~$B_r(x_0)\Subset\Omega$
and~$u\ne0$ in~$B_r(x_0)\setminus\{x_0\}$. Thus, by continuity and the fact that~$B_r(x_0)\setminus \{x_0\}$
is a connected set when~$n\ge2$,
by possibly replacing~$u$ with~$-u$,
we can suppose that~$u>0$ in~$B_r(x_0)\setminus\{x_0\}$. This and
Theorem~\ref{KAHAR}(iii) yield that
$$ 0=u(x_0)=\fint_{B_r(x_0)} u(x)\,dx>0,$$
which is a contradiction.
\end{proof}

We remark that Corollary~\ref{e5vbhfueweuwiuet0000} does not hold true
when~$n=1$, since the function~$u(x)=x$ for all~$x\in\R$ is harmonic but
possesses an isolated zero.
\medskip

It can be useful to stress that the
``almost every~$x_0\in\Omega$'' and ``almost
every~$r>0$'' in Theorem~\ref{KAHAR}(ii)-(iii)
can be replaced by the simpler
``every~$x_0\in\Omega$'' and ``every~$r>0$'' thanks to a continuity argument:
the details go as follows.

\begin{corollary}\label{KAHAR-OVUN}
Given an open set~$\Omega\subseteq\R^n$ and
a function~$u\in C(\Omega)$,
the following conditions are equivalent:
\begin{itemize}
\item[(i).] The function~$u$
belongs to~$C^2(\Omega)$ and
is harmonic in $\Omega$.
\item[(ii).] For every~$x_0\in\Omega$ and 
every~$r>0$ such that~$B_r(x_0)\Subset\Omega$,
we have that
$$ u(x_0)=\fint_{\partial B_r(x_0)} u(x)\,d{\mathcal{H}}^{n-1}_x.$$
\item[(iii).] For
every~$x_0\in\Omega$ and every~$r>0$ such that~$B_r(x_0)\Subset\Omega$,
we have that
$$ u(x_0)=\fint_{B_r(x_0)} u(x)\,dx.$$
\end{itemize}
\end{corollary}

\begin{proof} We prove that condition~(i) is equivalent to condition~(ii)
(similarly, one can prove the equivalence between conditions~(i) and~(iii)).
Assume~(i) here. Then, condition~(i) in Theorem~\ref{KAHAR} holds true,
which entails condition~(ii) in Theorem~\ref{KAHAR}.
This gives that condition~(ii) here is satisfied for almost every~$x_0\in\Omega$ and almost
every~$r>0$. Let now~$\bar x\in\Omega$ and~$\bar r>0$ such that~$B_{\bar r}(\bar x)\Subset\Omega$.
Let~$x^{(j)}$ be a sequence converging to~$\bar{x}$ as~$j\to+\infty$,
with~$x^{(j)}$ belonging to the above mentioned set of full measure for which~(ii) holds true.
Let also~$r\in (\bar{r}-3|\bar{x}-x^{(j)}|,\bar{r}-2|\bar{x}-x^{(j)}|)$
and notice that~$B_r(x^{(j)})\subseteq B_{\bar r}(\bar{x})\Subset\Omega$. Accordingly,
we can find~$r_j\in(\bar{r}-3|\bar{x}-x^{(j)}|,\bar{r}-2|\bar{x}-x^{(j)}|)$ such that
$$ u(x^{(j)})=\fint_{\partial B_{r_j}(x^{(j)})} u(x)\,d{\mathcal{H}}^{n-1}_x.$$
Passing to the limit as~$j\to+\infty$ and using the continuity of~$u$, we obtain~(ii) here, as desired.

Suppose now that~(ii) here holds true.
Then, condition~(ii) in Theorem~\ref{KAHAR} holds true,
which entails condition~(i) in Theorem~\ref{KAHAR}, that is condition~(i) here.
\end{proof}

It is instructive to emphasize that the notion of harmonicity is ``local'': since it relies
on the value of the derivative of a function at points, if~${\mathcal{I}}$
is a sets of indices, $\Omega_i\subseteq\R^n$
are open sets for every~$i\in {\mathcal{I}}$
and~$u$ is harmonic in each of the~$\Omega_i$, then
\begin{equation}\label{2432rf4g-2rfjv-1k2ermfAKSd56yu-124rktgmMqe22rtS-23t4yhtj0}{\mbox{$u$ is harmonic in }}\bigcup_{i\in{\mathcal{I}}}\Omega_i.\end{equation} As a result, the Mean Value Formula in Theorem~\ref{KAHAR}
(or its modification in Corollary~\ref{KAHAR-OVUN})
can be also localized, according to this observation:

\begin{corollary}\label{2432rf4g-2rfjv-1k2ermfAKSd56yu-124rktgmMqe22rtS-23t4yhtj}
Given an open set~$\Omega\subseteq\R^n$ and
a function~$u\in C(\Omega)$,
the following conditions are equivalent:
\begin{itemize}
\item[(i).] The function~$u$
belongs to~$C^2(\Omega)$ and
is harmonic in $\Omega$.
\item[(ii).] For
every~$x_0\in\Omega$ there exists~$r_0>0$ such that~$B_{r_0}(x_0)\Subset\Omega$ and for
every~$r\in(0,r_0)$
we have that
$$ u(x_0)=\fint_{\partial B_r(x_0)} u(x)\,d{\mathcal{H}}^{n-1}_x.$$
\item[(iii).] For
every~$x_0\in\Omega$ there exists~$r_0>0$ such that~$B_{r_0}(x_0)\Subset\Omega$ and for
every~$r\in(0,r_0)$
we have that
$$ u(x_0)=\fint_{B_r(x_0)} u(x)\,dx.$$
\end{itemize}
\end{corollary}

The complete proof of Corollary~\ref{2432rf4g-2rfjv-1k2ermfAKSd56yu-124rktgmMqe22rtS-23t4yhtj}
will be postponed to page~\pageref{2432rf4g-2rfjv-1k2ermfAKSd56yu-124rktgmMqe22rtS-23t4yhtj.LAPA}
since it will rely on the existence of a ``harmonic replacement'',
namely of a harmonic function for given boundary data in a ball
(though perhaps intuitively obvious, in our line of reasoning this will be ensured by
the construction of an explicit solution kernel in Theorems~\ref{POIBALL1}
and~\ref{POIBALL}).\medskip
 
A useful consequence of Corollary~\ref{2432rf4g-2rfjv-1k2ermfAKSd56yu-124rktgmMqe22rtS-23t4yhtj} is the so-called
Schwarz reflection principle, which allows to extend a harmonic function vanishing
on a hyperplane by odd reflection:

\begin{lemma}\label{RIFLESCHZ}
Let~$\Omega\subseteq\R^n$ be an open set.
Let
\begin{eqnarray*}&& \Omega_+:=\Big\{x=(x',x_n)\in\Omega{\mbox{ s.t. }}x_n>0\Big\},\\
&& \Omega_0:=\Big\{x=(x',x_n)\in\Omega{\mbox{ s.t. }}x_n=0\Big\},\\
&& \Omega_-:=
\Big\{x=(x',x_n){\mbox{ s.t. $x_n<0$ and }}(x',-x_n)\in\Omega_+\Big\}
\\{\mbox{and }}&&
\Omega_\star:=\Omega_+\cup\Omega_0\cup\Omega_-.\end{eqnarray*}
Let~$u\in C^2(\Omega_+)\cap C(\Omega_+\cup \Omega_0)$ be harmonic
and such that~$u(x)=0$ along~$\Omega_0$.

Then, the function
$$ \Omega_\star\ni x\mapsto u_\star(x)=u_\star(x',x_n)
:=\begin{dcases} u(x)& {\mbox{ if }}x\in\Omega_+\cup\Omega_0,\\-
u(x',-x_n)& {\mbox{ if }}x\in\Omega_-
\end{dcases}$$
is harmonic in~$\Omega_\star$.
\end{lemma}

\begin{proof}
We observe that
\begin{equation}\label{09io29203yoiewfgh7328ryfu237o8trgfoy2fr73twefg83246tyg89ergfwu-1203ury}
{\mbox{$u_\star$ is continuous and~$u_\star(x',0)=0$
for all~$x=(x',0)\in\Omega_0$.}}\end{equation}
Furthermore, if~$x\in\Omega_-$ then~$\Delta u_\star(x)=-\Delta u(x',-x_n)=0$
and accordingly~$u_\star$ is harmonic in~$\Omega_+\cup\Omega_-$, due
to~\eqref{2432rf4g-2rfjv-1k2ermfAKSd56yu-124rktgmMqe22rtS-23t4yhtj}.
This and Corollary~\ref{2432rf4g-2rfjv-1k2ermfAKSd56yu-124rktgmMqe22rtS-23t4yhtj}(iii) give that
\begin{equation}\label{09io29203yoiewfgh7328ryfu237o8trgfoy2fr73twefg83246tyg89ergfwu-1203ury2}
\begin{split}&{\mbox{for every~$\bar x\in\Omega_+\cup\Omega_-$
there exists~$\bar{r}>0$ such that for all~$r\in(0,\bar{r})$ we have that
}} \\&\fint_{B_r(\bar x)}u_\star(x)\,dx=u_\star(\bar x).\end{split}\end{equation}
Now we take~$x_0\in\Omega_0$ and~$r>0$ such that~$B_r(x_0)\Subset\Omega_\star$
and we observe that
\begin{equation*}
\begin{split}&
\fint_{B_r(x_0)}u_\star(x)\,dx=
\frac{1}{|B_r|}\left(
\int_{B_r(x_0)\cap \Omega_+}u_\star(x)\,dx+
\int_{B_r(x_0)\cap \Omega_-}u_\star(x)\,dx
\right)
\\&\qquad=
\frac{1}{|B_r|}\left(
\int_{B_r(x_0)\cap \Omega_+}u(x)\,dx-
\int_{B_r(x_0)\cap \Omega_-}u(x',-x_n)\,dx
\right)\\&\qquad=
\frac{1}{|B_r|}\left(
\int_{B_r(x_0)\cap \Omega_+}u(x)\,dx-
\int_{B_r(x_0)\cap \Omega_+}u(x',x_n)\,dx
\right)=0=u_\star(x_0),
\end{split}
\end{equation*}
thanks to~\eqref{09io29203yoiewfgh7328ryfu237o8trgfoy2fr73twefg83246tyg89ergfwu-1203ury}.

{F}rom this and~\eqref{09io29203yoiewfgh7328ryfu237o8trgfoy2fr73twefg83246tyg89ergfwu-1203ury2},
it follows that for every~$\bar x\in\Omega_\star$
there exists~$\bar{r}>0$ such that for all~$r\in(0,\bar{r})$ we have that~$\fint_{B_r(\bar x)}u_\star(x)\,dx=u_\star(\bar x)$.
This and Corollary~\ref{2432rf4g-2rfjv-1k2ermfAKSd56yu-124rktgmMqe22rtS-23t4yhtj} give that~$u_\star$ is harmonic in~$\Omega_\star$,
as desired.\end{proof}

The following classical identities, which are useful variations
of the Divergence Theorem, are often very helpful
and they are named\footnote{One of the cornerstones
of the foundation of potential theory
was indeed George Green's 1828 essay~\cite{GRESSAY},
see Figure~\ref{GREENFIDItangeFI},
which founded the mathematical theory of electricity and magnetism. The first edition of the essay
was printed for the author and sold on a subscription basis to only 51 people.
In 1850--1854, the essay was transcribed in~\cite{MR1578654, MR1578805, MR1578862}, with several
typographical corrections and a reference section added. Interestingly,
George Green was almost entirely self-taught, having
received only about one year of formal schooling, between the ages of~8 and~9.

By the way, the story of George Green is truly amazing. George Green's
father, also called George Green, was a baker in Nottingham, the UK town linked to the legend of Robin Hood.
Actually this legend does reflect long-lasting social problems
in the area which also affected the Green family.
In particular, at some point bakers were blamed for the incessant rise of the price of bread
and crowds of people broke into bakers to steal food and, in this circumstance, 
the Green family's bakery was also attacked.

The bakery was however probably doing well from the financial point of view, since, the year after these riots,
the little George was sent to allegedly the best and most expensive school in Nottingham, where he was taught
for four terms, which, as mentioned above, amount to his whole formal training:
at nine, the boy starts working in his father's bakery business.

And this business kept being profitable, allowing the Green family to buy a land
and build a brick wind corn-mill
(the mill was renovated in 1986 and is now a science center, see Figure~\ref{P2DAILAGREEMIItange3FIlAALANSE}).
George Green Jr. then fell in love with the daughter
of the manager of the mill, named Jane Smith. They never married but they had together seven children.

Green also joined the Nottingham Subscription Library, thus finding access to a few scientific books and articles,
and also some works published in other countries.

Green used to study and do mathematics on the top floor of the mill and it is probably here
that the famous essay by him was conceived and written.
The~51 subscribers who bought the first edition of the book (at the price of~7 pounds and
6~pence) where probably for a vast majority members of the Nottingham Subscription Library and likely
their mathematical proficiency was insufficient to fully appreciate the content of the essay.
However, among these subscribers there was also Sir Edward Thomas Bromhead, 2nd Baronet, wealthy landowner, mathematician and founder of the Analytical Society, a precursor of the Cambridge Philosophical Society. Bromhead
realized that the essay was the production of a brilliant scientist and
invited Green to send any further papers to the Royal Society of London, the Royal Society of Edinburgh and the Cambridge Philosophical Society. Green took Bromhead's offer as mere politeness and did not respond for two years.

Then, following Bromhead's encouragement, Green wrote three further papers and, having
accumulated considerable wealth and land owned, 
was able to abandon his miller duties, pursue mathematical studies
and, aged nearly forty, enroll as an undergraduate at the University of Cambridge.

There are rumors that, at Cambridge, Green 
may have succumbed to alcohol, possibly losing the endorsement of
his earlier supporters (we forgot to mention that Bromhead was also approached by local ministers
to help to establish a Temperance Society).

George Green died age~48 and
the Nottingham Review published the following short obituary:
``we believe he was the son of a miller, residing near Nottingham, but having a taste for study, he applied his gifted mind to the science of mathematics [...]. Had his life been prolonged, he might have stood eminently high as a mathematician''.
They obviously did not understand that he had already stood most eminently among his contemporaries
and left an essay that would have revolutionized the history of science.

Historians of science are unsure how Green managed to acquire his formidable mathematical knowledge
with so little formal training. One possibility however, is that the Leibniz-Newton calculus controversy
(the silly dispute about who had first invented calculus) resounded 
to a great disadvantage of the English school, which remained locked into the Newtonian notation of calculus,
stubbornly rejecting the notation introduced by Leibnitz and adopted
by continental mathematicians, which ultimately proved to be more flexible and effective.
It is possible that, being self-taught, Green had the possibility of getting in contact with Leibnitz's notation
(or possibly develop his own approach to calculus and analysis) without a rigid bias of an academia
influenced by dummy politics and sterile nationalism (well, anyway, a good notation
is always helpful, but in Green's case personal talent and inventiveness certainly made the difference).

For further readings on the figure of George Green see~\cite{MR1016105, zbMATH00495972, MR1829410, CHA:L2003}.} ``Green's Identities''\index{Green's Identities}:

\begin{figure}
  \centering
  \includegraphics[width=.45\linewidth]{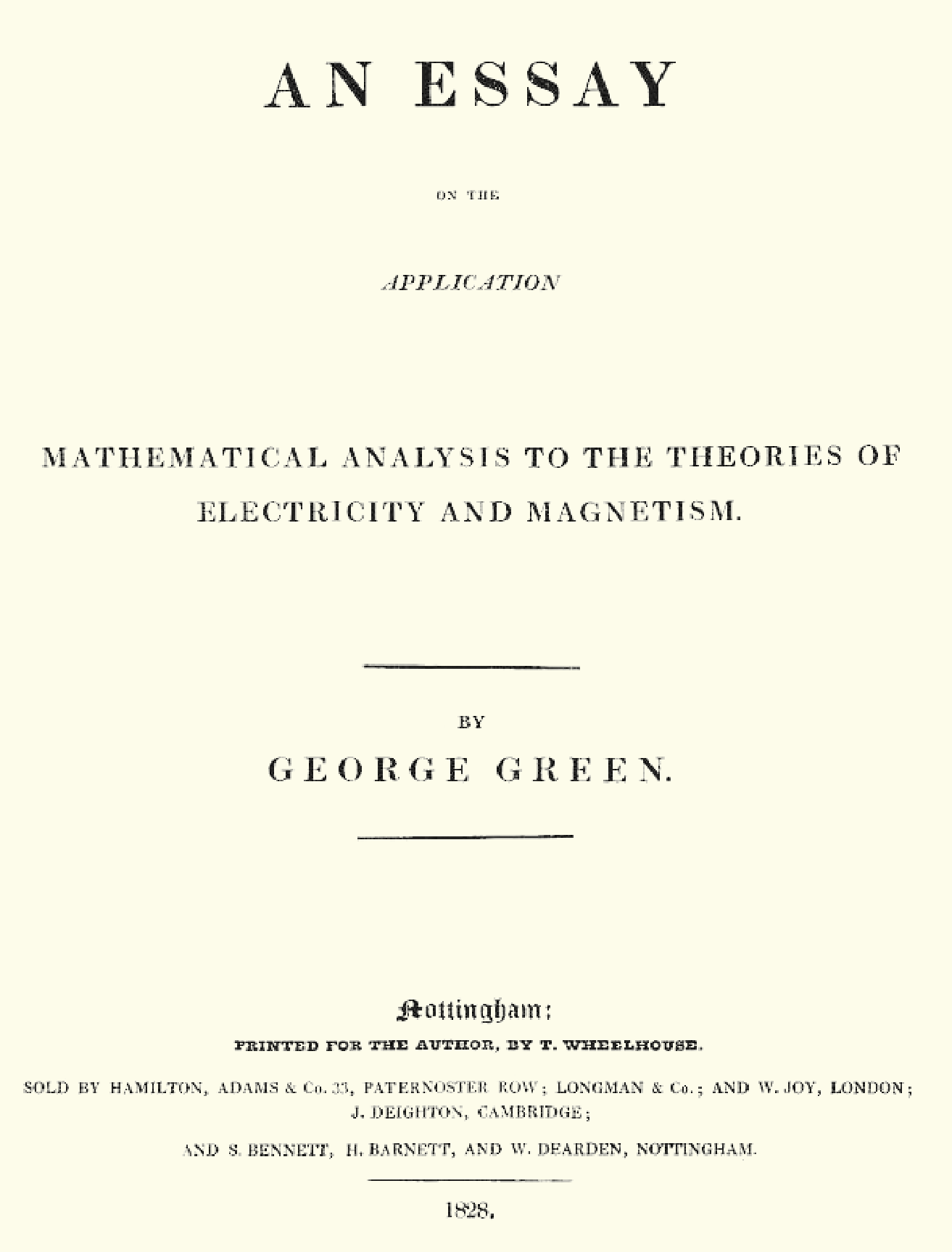}
 \caption{\sl The title page to Green's original essay (Public Domain image from
 Wikipedia).}\label{GREENFIDItangeFI}
\end{figure}

\begin{lemma}\label{GREEN}
Let~$\Omega\subseteq\R^n$ be a bounded open set of class~$C^1$, with exterior normal~$\nu$,
and let~$\varphi$, $\psi\in C^2(\Omega)\cap C^1(\overline\Omega)$.

Then,
\begin{equation}\label{GRr1}
\int_\Omega \Big( \varphi(x)\Delta \psi(x)+\nabla \varphi(x)\cdot\nabla \psi(x)\Big)\,dx=
\int_{\partial\Omega} \varphi(x)\frac{\partial \psi}{\partial\nu}(x)\,d{\mathcal{H}}^{n-1}_x
\end{equation}
and
\begin{equation}\label{GRr2}
\int_\Omega \Big( \varphi(x)\Delta \psi(x)-\psi(x)\Delta \varphi(x)\Big)\,dx=
\int_{\partial\Omega}\left(\varphi(x)\frac{\partial \psi}{\partial\nu}(x)-
\psi(x)\frac{\partial \varphi}{\partial\nu}(x)\right)\,d{\mathcal{H}}^{n-1}_x
.\end{equation}
\end{lemma}

\begin{proof} By the Divergence Theorem,
\begin{eqnarray*}&&
\int_\Omega \Big( \varphi(x)\Delta \psi(x)+\nabla \varphi(x)\cdot\nabla \psi(x)\Big)\,dx=
\int_\Omega \div\big( \varphi(x)\nabla \psi(x)\big)\,dx
=
\int_{\partial\Omega} \varphi(x)\frac{\partial \psi}{\partial\nu}(x)\,d{\mathcal{H}}^{n-1}_x
,\end{eqnarray*}
that is~\eqref{GRr1}.

Also, exchanging the roles of~$\varphi$ and~$\psi$ in~\eqref{GRr1},
\[
\int_\Omega \Big( \psi(x)\Delta \varphi(x)+\nabla \psi(x)\cdot\nabla \varphi(x)\Big)\,dx=
\int_{\partial\Omega} \psi(x)\frac{\partial \varphi}{\partial\nu}(x)\,d{\mathcal{H}}^{n-1}_x
.\]
Subtracting this from~\eqref{GRr1} we obtain~\eqref{GRr2}.
\end{proof}

We observe that identity~\eqref{1111DIV0934} can now be considered as a special
case of~\eqref{GRr1}. If creatively exploited,
Green's Identities are very useful to deduce important integral formulas,
which in turn entail structural information on several relevant equations.
As a prototype of this idea, we recall the classical Poho\v{z}aev Identity\index{Poho\v{z}aev Identity}
(see~\cite{MR0192184}; actually, an early occurrence of the Poho\v{z}aev Identity was obtained already in~\cite[equation~(2)]{MR2456}):

\begin{theorem}\label{Pohozaev Identity}
Let~$\Omega$ be a bounded open set in~$\R^n$ with~$C^1$ boundary
and~$u\in C^2(\Omega)\cap C^1(\overline\Omega)$ be a solution of
\begin{equation}\label{EQU-peru} \begin{dcases}
\Delta u=f (u)& {\mbox{ in }}\Omega,\\
u=0&{\mbox{ on }}\partial\Omega,
\end{dcases}\end{equation}
for some~$f\in L^\infty_{\rm loc}(\R)$.

\begin{figure}
  \centering
  \includegraphics[width=.52\linewidth]{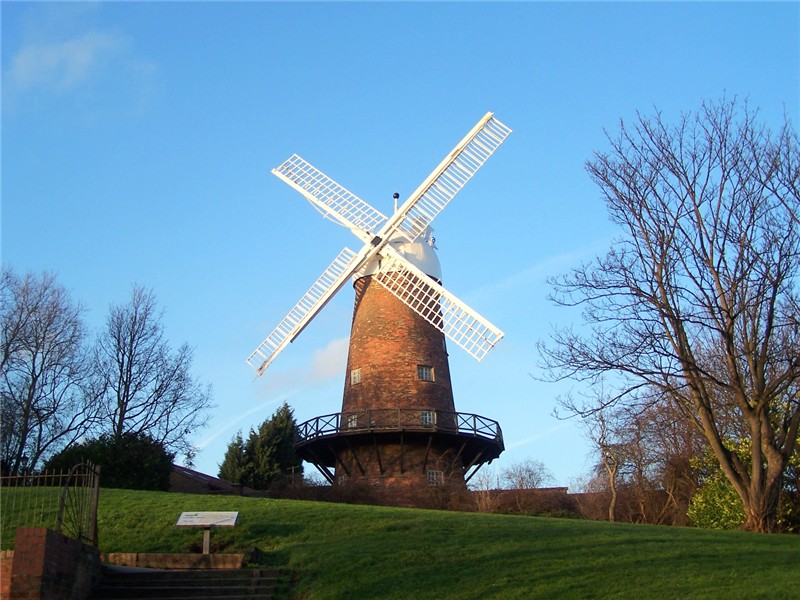}
 \caption{\sl Green's mill (photo by Kev747, image from
 Wikipedia, licensed under the Creative Commons Attribution-Share Alike 3.0 Unported license).}\label{P2DAILAGREEMIItange3FIlAALANSE}
\end{figure}

Let also
$$ F(r):=\int_0^r f(t)\,dt.$$
Then,
$$ \frac12\,\int_{\partial\Omega} (\partial_\nu u(x))^2\,(x\cdot\nu(x))\,d{\mathcal{H}}^{n-1}_x=
\frac{n-2}2\,\int_\Omega u(x) \,f(u(x))\,dx-n\int_\Omega F(u(x))\,dx.
$$
\end{theorem}

\begin{proof}
The idea of the proof is to test the equation against
the radial derivative~$\nabla u(x)\cdot x$ using suitable integration by parts.
Namely, from~\eqref{EQU-peru},
\begin{equation*}
\begin{split}&
\int_\Omega \Delta u(x)\,(\nabla u(x)\cdot x)\,dx
=\int_\Omega f(u(x))(\nabla u(x)\cdot x)\,dx\\&\qquad
=\int_\Omega \nabla \big( F(u(x))\big)\cdot x\,dx
=\int_\Omega\Big( \div \big( F(u(x))\, x\big) -nF(u(x))\Big)\,dx.\end{split}
\end{equation*}
This and the Divergence Theorem,
recalling the boundary condition in~\eqref{EQU-peru},
give that
\begin{equation}\label{SPL8ygbT2Ter}
\begin{split}&
\int_\Omega \Delta u(x)\,(\nabla u(x)\cdot x)\,dx+
n\int_\Omega F(u(x))\,dx=
\int_\Omega \div \big( F(u(x))\, x\big)\,dx\\&\qquad
=\int_{\partial\Omega} F(u(x))\, (x\cdot\nu(x))\,d{\mathcal{H}}^{n-1}_x
=\int_{\partial\Omega} F(0)\, (x\cdot\nu(x))\,d{\mathcal{H}}^{n-1}_x
=0
.\end{split}\end{equation}
Furthermore,
\begin{eqnarray*}&&
\Delta u\,(\nabla u\cdot x)=\div\Big((\nabla u\cdot x)\nabla u\Big)-
\nabla(\nabla u\cdot x)\cdot\nabla u
=\div\Big((\nabla u\cdot x)\nabla u\Big)-\sum_{i,j=1}^n
\partial_{ij}u\,\partial_i u\,x_j-|\nabla u|^2\\&&\qquad
=\div\Big((\nabla u\cdot x)\nabla u\Big)-\frac12\,\sum_{j=1}^n
\partial_{j}|\nabla u|^2\,x_j-|\nabla u|^2
=\div\Big((\nabla u\cdot x)\nabla u\Big)-\frac12\,\nabla (|\nabla u|^2)\cdot
x-|\nabla u|^2\\&&\qquad=
\div\Big((\nabla u\cdot x)\nabla u\Big)-\frac12\,\left[
\div\Big(|\nabla u|^2\,x\Big)
-n|\nabla u|^2
\right]-|\nabla u|^2\\&&\qquad=
\div\Big((\nabla u\cdot x)\nabla u\Big)-\frac12\,
\div\Big(|\nabla u|^2\,x\Big)
+\frac{n-2}2|\nabla u|^2.
\end{eqnarray*}Thus, making use again
of the Divergence Theorem and of the first Green's Identity~\eqref{GRr1}
we find that
\begin{equation}\begin{split}\label{dtv654574v6} &\int_\Omega \Delta u(x)\,(\nabla u(x)\cdot x)\,dx\\
=\;&\int_{\partial\Omega}
(\nabla u(x)\cdot x)\partial_\nu u(x)\,d{\mathcal{H}}^{n-1}_x
-\frac12\,\int_{\partial\Omega}|\nabla u(x)|^2\,x\cdot\nu(x)\,d{\mathcal{H}}^{n-1}_x
+\frac{n-2}2\,\int_\Omega|\nabla u(x)|^2\,dx\\
=\;&\int_{\partial\Omega}
(\nabla u(x)\cdot x)\partial_\nu u(x)\,d{\mathcal{H}}^{n-1}_x
-\frac12\,\int_{\partial\Omega}|\nabla u(x)|^2\,x\cdot\nu(x)\,d{\mathcal{H}}^{n-1}_x
-\frac{n-2}2\,\int_\Omega\Delta u(x)\,u(x)\,dx
.\end{split}\end{equation}
Now we observe that~$\nabla u=\pm|\nabla u|\nu$ on~$\partial\Omega$, and therefore, for every~$x\in\Omega$,
$$(\nabla u(x)\cdot x)\partial_\nu u(x)=|\nabla u(x)|^2(\nu(x)\cdot x)(\nu(x)\cdot\nu(x))=|\nabla u(x)|^2(\nu(x)\cdot x).
$$
Plugging this information into~\eqref{dtv654574v6}, we find that
$$\int_\Omega \Delta u(x)\,(\nabla u(x)\cdot x)\,dx=
\frac12\,\int_{\partial\Omega}|\nabla u(x)|^2\,x\cdot\nu(x)\,d{\mathcal{H}}^{n-1}_x
-\frac{n-2}2\,\int_\Omega\Delta u(x)\,u(x)\,dx.$$
Combining this and~\eqref{SPL8ygbT2Ter}, we obtain the desired result.
\end{proof}

Equations as in~\eqref{EQU-peru}
are often called ``semilinear''  \index{semilinear equation}
since they are not linear in~$u$ (unless the source term~$f$
is linear) but they are linear in the second derivative of~$u$.
Interestingly,
solutions of semilinear equations\footnote{We refer to the footnotes
on pages~\pageref{BIS Pohozaev Identity}
and~\pageref{TRIS Pohozaev Identity} for motivational comments about semilinear equations.} enjoy the special
feature of having constant Laplacian along their level sets,
namely if~$u$ solves~\eqref{EQU-peru},
given any~$c\in\R$, we have that~$\Delta u=f(c)$
on~$\{u=c\}$.
\medskip

As a consequence of the Poho\v{z}aev Identity
in Theorem~\ref{Pohozaev Identity}, one obtains
nonexistence results, as the one in the forthcoming Corollary~\ref{corpo}.
For this, we give the following definition:

\begin{definition} Let~$\Omega\subseteq \R^n$. Given~$x_0\in\Omega$, we say that~$\Omega$ is starshaped
with respect to~$x_0$ if for every~$x\in\Omega$ we have that~$tx+(1-t)x_0\in\Omega$ for all~$t\in[0,1]$.

Furthermore, we say that~$\Omega$ is starshaped if there exists~$x_0\in\Omega$
such that~$\Omega$ is starshaped with respect to~$x_0$.
\end{definition}

With this, we give the following nonexistence result:

\begin{corollary}\label{corpo}
Let~$n\ge 3$ and~$p > \frac{n+ 2}{n - 2}$.
Let~$\Omega$ be a bounded starshaped
open set in~$\R^n$ with~$C^1$ boundary.
Let~$u\in C^2(\Omega)\cap C^1(\overline\Omega)$ be a solution of
\begin{equation*} \begin{dcases}
\Delta u=-|u|^{p-1}u& {\mbox{ in }}\Omega,\\
u=0&{\mbox{ on }}\partial\Omega.
\end{dcases}\end{equation*}
Then, $u$ vanishes identically.
\end{corollary}

\begin{proof} Up to a translation, we suppose that
\begin{equation}\label{ISpapkrfeppa}
{\mbox{$\Omega$ is
starshaped with respect to the origin.}}\end{equation}
 We claim that
\begin{equation}\label{PDMASNDc}
{\mbox{$x\cdot\nu(x)\geq0$ for every~$x\in\partial\Omega$.}}\end{equation}
To check this, given~$x_0\in\partial\Omega$, we write~$\Omega$
in the vicinity of~$x_0$ as the superlevels of some function~$\Phi\in C^1(\R^n)$
with~$\nabla\Phi(x_0)\ne0$,
that is we take~$\rho>0$ such that~$\Omega\cap B_\rho(x_0)=\{\Phi>0\}\cap
B_\rho(x_0)$. In this way, we have that~$\nu=-\frac{\nabla\Phi}{|\nabla\Phi|}$
on~$\partial\Omega$.
Also, by~\eqref{ISpapkrfeppa},
we have that~$t x_0\in\overline\Omega$ for every~$t\in[0,1]$.
As a result, for~$t\in[0,1]$ sufficiently close to~$1$,
we have~$\Phi(tx_0)\ge0$. Therefore,
$$ 0\ge
\lim_{t\nearrow1}\frac{\Phi(tx_0)}{t-1}=
\lim_{t\nearrow1}\frac{\Phi(tx_0)-\Phi(x_0)}{t-1}=\nabla\Phi(x_0)\cdot x_0=
-|\nabla\Phi(x_0)|\,\nu(x_0)\cdot x_0.
$$
This proves~\eqref{PDMASNDc}.

We now exploit the
Poho\v{z}aev Identity in Theorem~\ref{Pohozaev Identity} with~$f(u):=-|u|^{p-1}u$,
and hence~$F(r):=-\frac{|r|^{p+1}}{p+1}$. In this way, using~\eqref{PDMASNDc},
we find that
\begin{eqnarray*}&& 0\le \frac12\,\int_{\partial\Omega} 
(\partial_\nu u(x))^2\,(x\cdot\nu(x))\,d{\mathcal{H}}^{n-1}_x=-
\frac{n-2}2\,\int_\Omega |u(x)|^{p+1} \,dx+
\frac{n}{p+1}\,\int_\Omega |u(x)|^{p+1}\,dx\\&&\qquad\qquad=
\frac{p(2-n)+n+2}{2(p+1)}\,
\int_\Omega |u(x)|^{p+1} \,dx\le0.
\end{eqnarray*}
In particular,
$$ \frac{p(2-n)+n+2}{2(p+1)}\,\int_\Omega |u(x)|^{p+1} \,dx=0,$$
from which the desired result follows.\end{proof}

A natural question is whether or not
the average over balls and spheres in the Mean Value Formulas of Theorem~\ref{KAHAR}
can be substituted with averages on different sets.
As we will see in Section~\ref{KURAN-SEC},
this is not the case and in fact the geometry of the balls and spheres
play a decisive
role in the Mean Value Formula (this classical problem
was pioneered in~\cite{MR140700, MR177124, MR279320, MR320348, MR1021402}).

\section{Weak solutions}\label{WERGSB-EAFGBS-cojwnedfDICewew4N8wedibfn}

We present here a classical result\index{weak solution}
often referred to with the name of Weyl's Lemma\index{Weyl's Lemma}:

\begin{lemma}\label{WEYL}
Let~$\Omega\subseteq\R^n$ be an open set
and let~$u\in L^1_{\rm loc}(\Omega)$. Assume that
\begin{equation}\label{GRAVSSEFORC2}
\int_\Omega u(x)\,\Delta\varphi(x)\,dx=0\qquad{\mbox{for every }}\varphi\in
C^\infty_0(\Omega).\end{equation}
Then, $u$ is harmonic in~$\Omega$.
\end{lemma}

\begin{proof} We stress that the desired claim
follows directly from~\eqref{GRAVSSEFORC2} and the
second Green's Identity~\eqref{GRr2} when~$u\in C^2(\Omega)$.

If instead~$u$ is merely locally integrable in~$\Omega$,
we use a mollification argument.
To this end,
we take~$\tau\in C^\infty_0(B_1,\,[0,+\infty))$
with~$\int_{B_1}\tau(x)\,dx=1$.
Given~$\eta>0$, we let~$\tau_\eta(x):=\frac1{\eta^n}
\tau\left(\frac{x}{\eta}\right)$
and define~$u_\eta:=u*\tau_\eta$.
Then, given~$x_0\in\Omega$ and~$\rho>0$ such that~$B_{2\rho}(x_0)\Subset\Omega$,
for all~$\varphi\in C^\infty_0(B_\rho(x_0))$ and all~$\eta\in(0,\rho)$ we have that
\begin{equation*}\begin{split}&
\int_\Omega u_\eta(x)\,\Delta\varphi(x)\,dx=
\iint_{B_{2\rho}(x_0)\times\Omega} u(y)\,\tau_\eta(x-y)\Delta\varphi(x)\,dx\,dy\\&\qquad=
\int_\Omega u(x)\,(\tau_\eta*\Delta\varphi)(x)\,dx
=\int_\Omega u(x)\,\Delta\varphi_\eta(x)\,dx,\end{split}
\end{equation*}
where~$\varphi_\eta:=\varphi*\tau_\eta\in C^\infty_0(B_{2\rho}(x_0))\subseteq
C^\infty_0(\Omega)$.
As a result, from~\eqref{GRAVSSEFORC2}, we deduce that
$$ \int_\Omega u_\eta(x)\,\Delta\varphi(x)\,dx=0\qquad{\mbox{ for every }}
\varphi\in C^\infty_0(B_\rho(x_0)),$$
as long as~$\eta\in(0,\rho)$.
Since~$u_\eta\in C^2(B_\rho(x_0))$,
this gives that~$u_\eta$ is harmonic in~$B_\rho(x_0)$.
Owing to this and to
Theorem~\ref{KAHAR}(iii),
for every ball~$B_r(\overline{x})\Subset B_\rho(x_0)$,
we have that
$$ u_\eta(\overline{x})=\fint_{B_r(\overline{x})} u_\eta(x)\,dx.$$
We now send~$\eta\searrow0$ and (see e.g. Theorems~9.6 and~9.13
in~\cite{MR3381284}) we conclude that,
whenever~$B_r(\overline{x})\Subset B_\rho(x_0)$
and~$\overline{x}$ is a Lebesgue
density point for~$u$,
\begin{equation} \label{3LS2345Dscdfoj74586vvyjf}
u(\overline{x})=\fint_{B_r(\overline{x})} u(x)\,dx.\end{equation}
Furthermore, by the Dominated Convergence Theorem,
for every~$\widetilde{x}\in\Omega$ and~$r>0$ such that~$B_r(
\widetilde{x})\subset\Omega$,
$$ \lim_{p\to\widetilde{x}}\int_{B_r(p)} u(x)\,dx=
\int_{B_r(\widetilde{x})} u(x)\,dx.$$
This and~\eqref{3LS2345Dscdfoj74586vvyjf} give that, up
to continuously extend~$u$ in a set of null Lebesgue measure
in~$B_\rho(x_0)$, we have that
\[ u(\overline{x})=\fint_{B_r(\overline{x})} u(x)\,dx\qquad{\mbox{for every }}
\overline{x}\in B_\rho(x_0),\]
as long as~$B_r(\overline{x})\Subset B_\rho(x_0)$.
Using again Theorem~\ref{KAHAR}, we thereby conclude that~$u$
is harmonic in~$B_\rho(x_0)$. This gives that~$\Delta u(x_0)=0$
for every~$x_0\in\Omega$, as desired.
\end{proof}

\begin{corollary}\label{S-OS-OSSKKHEPARHA}
Let~$\Omega\subseteq\R^n$ be an open set
and~$u_k$ be a sequence of harmonic functions in~$\Omega$.
Suppose that~$u_k\to u$ in~$L^1_{\rm loc}(\Omega)$. Then, $u$
is harmonic in~$\Omega$.
\end{corollary}

\begin{proof}
Let~$\Omega'\Subset\Omega$ and~$\varphi\in
C^\infty_0(\Omega')$.
By the
second Green's Identity~\eqref{GRr2} we know that
\begin{eqnarray*}&&
\left|\int_\Omega u(x)\,\Delta\varphi(x)\,dx\right|
\le
\int_{\Omega'} |u(x)-u_k(x)|\,|\Delta\varphi(x)|\,dx
+\left|\int_\Omega u_k(x)\,\Delta\varphi(x)\,dx\right|\\&&\qquad
\le \|\varphi\|_{C^2(\Omega')}\,\|u-u_k\|_{L^1(\Omega')}
+\left|\int_\Omega \Delta u_k(x)\,\varphi(x)\,dx\right|=\|\varphi\|_{C^2(\Omega')}\,\|u-u_k\|_{L^1(\Omega')}.
\end{eqnarray*}
Hence, sending~$k\to+\infty$,
$$ \int_\Omega u(x)\,\Delta\varphi(x)\,dx=0.$$
The desired result thus follows from
Lemma~\ref{WEYL}.
\end{proof}

An alternative proof of Corollary~\ref{S-OS-OSSKKHEPARHA}
can be obtained also using directly the 
Mean Value Formula in Theorem~\ref{KAHAR}(iii). \medskip

See also~\cite{MR3362185, MR3380662}
and the references therein for a more complete discussion
on the role played by ``weak'' or ``distributional''
formulations of partial differential equations and a careful discussion
of the functional analysis methods involved in such a theory.
\medskip

A classical application of the weak setting of partial differential
equations is provided by Kato's Inequality\index{Kato's Inequality}, see~\cite{MR333833},
as presented in the following result.
For this, we use the standard notation, for every~$r\in\R$,
$$\sign(r):=\begin{dcases} \frac{r}{|r|}
&{\mbox{ if }}r\ne0,\\0&{\mbox{ if }}r=0.\end{dcases}$$ 

\begin{theorem} \label{KATOTH} Let
\begin{equation}\label{L1co} u\in L^1_{\rm loc}(\R^n)\end{equation}
be such that there exists~$f\in L^1_{\rm loc}(\R^n)$ satisfying
\begin{equation}\label{L1co2} \int_{\R^n} u(x)\,\Delta\psi(x)\,dx
=\int_{\R^n} f(x)\,\psi(x)\,dx\qquad{\mbox{ for all }}\psi\in C^\infty_0(\R^n). \end{equation}

Then, for\footnote{In jargon,
condition~\ref{L1co2} can be rewritten by stating that the Laplacian
of~$u$, as defined in the weak sense, is actually a locally integrable function
that is denoted by~$f$ (and
this is of course the case for smooth functions~$u$).

Similarly, equation~\eqref{MSN:LASsmmS2} can be written as
$$\Delta|u|\ge\sign(u)\,\Delta u$$ in the weak sense.

This can be also considered as a ``limit case'' of the following
observation: if~$\Phi\in C^1(\R)$ is a convex function, then
$$\Phi
(u(x\pm he_i))-\Phi(u(x))\ge\Phi'(u(x))
(u(x\pm he_i)-u(x))$$ and accordingly
$$ \Phi(u(x+ he_i))+\Phi(u(x- he_i))-2\Phi(u(x))\ge
\Phi'(u(x))
(u(x+he_i)+u(x- he_i)-2u(x)),$$
which leads to
$$ \Delta\Big(\Phi(u(x))\Big)\ge\Phi'(u(x))\,\Delta u(x).$$
With this respect, Kato's Inequality~\eqref{MSN:LASsmmS2}
corresponds, formally, to the limit case in which~$\Phi(t)=|t|$,
that produces~$\Phi'(t)=\sign(t)$ (at least when~$t\ne0$,
and the corresponding inequality holding true in the weak sense).} every~$\varphi\in C^\infty_0(\R^n,\,[0,+\infty))$,
\begin{equation}\label{MSN:LASsmmS2}\int_{\R^n} |u(x)|\,\Delta\varphi(x)\,dx\geq\int_{\R^n}\sign(u(x))\,\varphi(x)\, f(x)\,dx.\end{equation}\end{theorem} 

\begin{proof}
Given~$\eta>0$, we let~$\tau_\eta(x):=\frac1{\eta^n}
\tau\left(\frac{x}{\eta}\right)$
and define~$u_\eta:=u*\tau_\eta$. In this setting, by~\eqref{L1co},
possibly up to subsequences, we have that~$u_\eta$ converges
to~$u$ in~$L^1_{\rm loc}(\R^n)$
(see~\cite[Theorem 9.6]{MR3381284}) and almost everywhere,
with additionally~$|u_\eta|\le h$ for a suitable~$h\in L^1_{\rm loc}(\R^n)$
(see~\cite[Theorem~4.9]{MR2759829}).
For this reason,
\begin{equation}\label{NOSMDBERPPLSBE} \int_{\R^n} |u(x)|\,
\Delta\varphi(x)\,dx
=\lim_{\eta\searrow0}\int_{\R^n} |u_\eta(x)|\,\Delta\varphi(x)\,dx.\end{equation}
Furthermore, 
for every~$\psi\in C^\infty_0(\R^n)$,
\begin{eqnarray*}&& \int_{\R^n} u_\eta(x)\,\Delta\psi(x)\,dx
=\int_{\R^n}\left[
\int_{\R^n} u(y)\,\tau_\eta(x-y)\,\Delta\psi(x)\,dx\right]\,dy
=\int_{\R^n} u(y)\,\Delta\psi_\eta(y)\,dy\\&&\qquad
=\int_{\R^n} f(y)\,\psi_\eta(y)\,dy
=\int_{\R^n}\left[
\int_{\R^n} f(y)\,\psi(x)\,\tau_\eta(x-y)\,dy
\right]\,dx=\int_{\R^n} f_\eta(x)\,\psi(x)\,dx,
\end{eqnarray*}
and accordingly~$\Delta u_\eta=f_\eta$.

In this way, possibly extracting a subsequence,
we deduce that~$\Delta u_\eta$ converges
to~$f$ in~$L^1_{\rm loc}(\R^n)$
(see~\cite[Theorem 9.6]{MR3381284}) and almost everywhere,
with additionally~$|\Delta u_\eta|\le H$ for a suitable~$H\in L^1_{\rm loc}(\R^n)$
(see~\cite[Theorem~4.9]{MR2759829}).

Now we let~$\e>0$ and set~$v_{\e,\eta}(x):=\sqrt{(u_\eta(x))^2+\e^2}$.
In this way, the function~$v_{\e,\eta}$ belongs to~$C^\infty(\R^n)$
and~$\sign(v_{\e,\eta}(x))=1$ for all~$x\in\R^n$. Furthermore, $$ 2v_{\e,\eta}
\nabla v_{\e,\eta}=\nabla v_{\e,\eta}^2=\nabla(u_\eta^2+\e^2)=2u_\eta
\nabla u_\eta.$$ Consequently,
$$ |u_\eta|\,|\nabla v_{\e,\eta}|\le v_{\e,\eta}\,|\nabla v_{\e,\eta}|= |u_\eta
|\,|\nabla u_\eta|$$
and
$$ |\nabla v_{\e,\eta}|^2+v_{\e,\eta}\Delta v_{\e,\eta}
= \div(v_{\e,\eta}\nabla v_{\e,\eta})=\div(u_\eta
\nabla u_\eta)=|\nabla u_\eta|^2+u_\eta
\Delta u_\eta.$$ As a result, $$ v_{\e,\eta}\Delta v_{\e,\eta}
=|\nabla u_\eta|^2-|\nabla v_{\e,\eta}|^2+u_\eta
\Delta u_\eta\geq u_\eta\Delta u_\eta.$$
We thus define~$\sigma_{\e,\eta}(x):=\frac{u_\eta(x)}{v_{\e,\eta}(x)}$
and find that~$\Delta v_{\e,\eta}\geq\sigma_{\e,\eta}\Delta u_\eta
$, and then
\begin{equation}\label{KA78:01} \int_{\R^n} v_{\e,\eta}(x)\,\Delta\varphi(x)
\,dx\geq\int_{\R^n}\sigma_{\e,\eta}
(x)\,\varphi(x)\, \Delta u_\eta(x)\,dx,\end{equation}
for every~$\varphi\in C^\infty_0(\R^n,\,[0,+\infty))$.

It is also helpful to observe that~$|v_{\e,\eta}|=v_{\e,\eta}\le
|u_\eta|+\e\le h+1$ and that, for a.e.~$x\in\R^n$,
$$ \lim_{\eta\searrow0}v_{\e,\eta}(x)=\sqrt{(u(x))^2+\e}.$$
We can therefore exploit the Dominated
Convergence Theorem to find that
\begin{equation}\label{CO6789SEDJSMBIMserD9ofkvRASHJD}
\lim_{\eta\searrow0}\int_{\R^n} v_{\e,\eta}(x)\,\Delta\varphi(x)\,dx
=\int_{\R^n} \sqrt{(u(x))^2+\e}
\;\Delta\varphi(x)\,dx.
\end{equation}
Additionally, for a.e.~$x\in\R^n$,
$$ \lim_{\eta\searrow0}\sigma_{\e,\eta}(x)=
\frac{u(x)}{\sqrt{(u(x))^2+\e}} ,$$
and~$|\sigma_{\e,\eta}|\le1$. Hence, by the Dominated Convergence Theorem,
\begin{equation*}
\lim_{\eta\searrow0}
\int_{\R^n}\sigma_{\e,\eta}
(x)\,\varphi(x)\, \Delta u_\eta(x)\,dx=
\int_{\R^n}\frac{u(x)}{\sqrt{(u(x))^2+\e}}\,\varphi(x)\, f(x)\,dx.\end{equation*}
Combining this fact with~\eqref{CO6789SEDJSMBIMserD9ofkvRASHJD},
we can pass~\eqref{KA78:01} to the limit as~$\eta\searrow0$
and see that,
for every~$\varphi\in C^\infty_0(\R^n,\,[0,+\infty))$,
\begin{equation*} \int_{\R^n} \sqrt{(u(x))^2+\e}
\;\Delta\varphi(x)\,dx\geq\int_{\R^n}\frac{u(x)}{\sqrt{(u(x))^2+\e}}
\,\varphi(x)\, f(x)\,dx.\end{equation*}
By sending~$\e\searrow0$, we thereby obtain
the desired result in~\eqref{MSN:LASsmmS2}.\end{proof}

As a consequence of 
Kato's Inequality in Theorem~\ref{KATOTH}, we present a classification
result for global weak solutions of the equation~$\Delta u=Vu+cu$,
see~\cite{MR0493420} for additional details.

\begin{corollary} Let~$V\in L^2_{\rm loc}(\R^n,\,[0,+\infty))$.
Let~$u\in L^2(\R^n)$ and assume that, for
every~$\varphi\in C^\infty_0(\R^n)$,
\begin{equation}\label{LSKMD-SKD}\int_{\R^n} u(x)\,\Delta\varphi(x)\,dx
=\int_{\R^n} V(x)\,u(x)\,\varphi(x)\,dx.\end{equation}
Then, $u$ vanishes identically.\end{corollary}

\begin{proof} To exploit Theorem~\ref{KATOTH},
we notice that condition~\eqref{L1co2} is fulfilled
with~$f(x):= V(x)\,u(x)\in  L^1_{\rm loc}(\R^n)$, thanks to~\eqref{LSKMD-SKD}.
Therefore, in light of~\eqref{MSN:LASsmmS2},
for every~$\varphi\in C^\infty_0(\R^n,\,[0,+\infty))$,
\begin{equation}\label{GDBFNIqrtJSMFFIDIRENCEXNCKFNEJFINFNIFG}
\int_{\R^n} |u(x)|\,\Delta\varphi(x)\,dx
\geq\int_{\R^n}\sign(u(x))\,\varphi(x)\,
V(x)\,u(x)\,dx.\end{equation}
Now
we take~$\tau\in C^\infty_0(B_1,\,[0,+\infty))$
with~$\int_{B_1}\tau(x)\,dx=1$.
Given~$\eta>0$, we let~$\tau_\eta(x):=\frac1{\eta^n}
\tau\left(\frac{x}{\eta}\right)$
and define~$w_\eta:=|u|*\tau_\eta$. Notice that~$|u|\in L^2(\R^n)$
and therefore~$w_\eta\to|u|$
in~$L^2(\R^n)$
(see e.g.~\cite[Theorem 9.6]{MR3381284}).
As a result, possibly up to a subsequence,
there exists~$h\in L^2(\R^n)$ such that
\begin{equation}\label{KS-DUEMD-L2ANSAVSAILA}
|w_\eta(x)|\le h(x)\end{equation} 
a.e.~$x\in\R^n$ (see e.g.~\cite[Theorem~4.9]{MR2759829}).

Additionally,
for every~$\varphi\in C^\infty_0(\R^n,\,[0,+\infty))$,
\begin{eqnarray*}
&& \int_{\R^n} w_\eta(x)\,\Delta\varphi(x)\,dx
=\int_{\R^n} \left[\int_{\R^n}
|u(y)|\,\tau_\eta(x-y)\,\Delta\varphi(x)\,dx\right]\,dy
\\&&\qquad
=\int_{\R^n} |u(y)|\,\Delta\varphi_\eta(y)\,dy
\geq\int_{\R^n}\sign(u(y))\,\varphi_\eta(y)\,
V(y)\,u(y)\,dy\\&&\qquad
=\int_{\R^n}\varphi_\eta(y)\,
V(y) \,|u(y)|\,dy\ge0
,\end{eqnarray*}
thanks to~\eqref{GDBFNIqrtJSMFFIDIRENCEXNCKFNEJFINFNIFG}.

As a result, we find that~$\Delta w_\eta\ge0$.
Thus, if~$R>1$ and~$\xi_R
\in C^\infty_0(B_R,\,[0,1])$ with~$\xi_R=1$ in~$B_{R-1}$
and~$|\nabla\xi_R|\le 2$, letting~$\zeta_R:=\xi_R^2$
we find that
\begin{eqnarray*}&&
\int_{\R^n} \zeta_R(x)\,|\nabla w_\eta(x)|^2\,dx
=-\int_{\R^n} \div\big(\zeta_R(x)\,\nabla w_\eta(x)\big)\,w_\eta(x)\,dx\\
&&\qquad=
-\int_{\R^n} \nabla\zeta_R(x)\cdot\nabla w_\eta(x) \,w_\eta(x)\,dx
-\int_{\R^n} \zeta_R(x)\,\Delta w_\eta(x)\,w_\eta(x)\,dx\\&&\qquad\le
-\int_{\R^n} \nabla\zeta_R(x)\cdot\nabla w_\eta(x) \,w_\eta(x)\,dx=
-2\int_{\R^n} \xi_R(x)\,\nabla\xi_R(x)\cdot\nabla w_\eta(x) \,w_\eta(x)\,dx
\\&&\qquad\le\frac12\int_{\R^n} \zeta_R(x)\,|\nabla w_\eta(x)|^2\,dx+
2\int_{\R^n}|\nabla\xi_R(x)|^2 \,(w_\eta(x))^2\,dx. \end{eqnarray*}
For this reason, and recalling~\eqref{KS-DUEMD-L2ANSAVSAILA},
\begin{eqnarray*}&&\frac12
\int_{\R^n} \zeta_R(x)\,|\nabla w_\eta(x)|^2\,dx\le
8\int_{\R^n\setminus B_{R-1}}(w_\eta(x))^2\,dx\le 8\int_{\R^n\setminus B_{R-1}}(h(x))^2\,dx
\end{eqnarray*}
and therefore
\begin{eqnarray*}&&\int_{\R^n} |\nabla w_\eta(x)|^2\,dx=\lim_{R\to+\infty}
\int_{B_{R-1}} |\nabla w_\eta(x)|^2\,dx
\le\lim_{R\to+\infty}
\int_{\R^n} \zeta_R(x)\,|\nabla w_\eta(x)|^2\,dx\\&&\qquad\le16
\lim_{R\to+\infty}
\int_{\R^n\setminus B_{R-1}}(h(x))^2\,dx=0.
\end{eqnarray*}
This leads to~$w_\eta$ being constant, and thus constantly equal to
zero, due to~\eqref{KS-DUEMD-L2ANSAVSAILA}. {F}rom this,
taking the limit as~$\eta\searrow0$, we find that~$u$ is
constantly equal to
zero as well.\end{proof}

\section{The Laplace-Beltrami operator}\label{BELS}

\begin{figure}
  \centering
  \includegraphics[width=.7\linewidth]{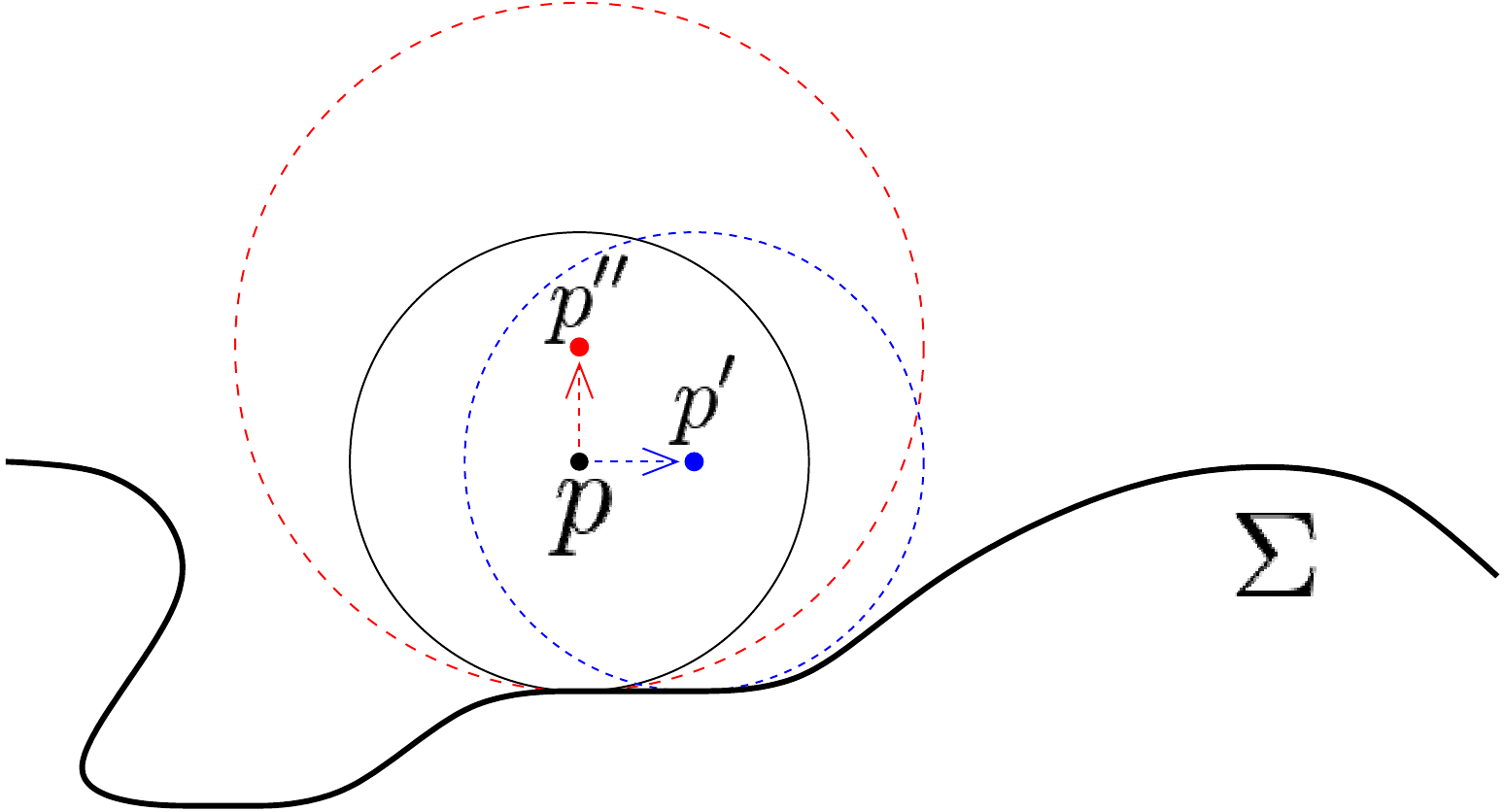}
 \caption{\sl Taking derivatives of the distance function.}\label{DItangeFI}
\end{figure}

The Laplace operator in~$\R^n$ is actually a ``special case''
of a more general operator acting
on functions defined on manifolds embedded in the
Euclidean space (or, even more generally, on Riemannian and pseudo-Riemannian manifolds). 
For concreteness, though more general settings
can be taken into account (see also the comments on page~\pageref{DAC3aC22}), we consider here the case of
a hypersurface~$\Sigma=\partial E$ of class~$C^3$, \label{DAC3aC2}
for a bounded and open set~$E\subseteq\R^n$,
and we denote by~$\nu$ its unit exterior normal and by~$d_\Sigma$ the \label{KJSMDISTANCEKDMF023945YT}
signed distance function\index{signed distance function}
to~$\Sigma$ (say, with the convention that~$d_\Sigma\ge0$ in~$E$ and~$d_\Sigma\le0$ in~$\R^n\setminus E$). We point out that~$d_\Sigma$ is also of class~$C^3$ \label{LAKSJHDKSDJFH0193287493utjgQUEMD}
in a suitably small neighborhood~${\mathcal{N}}$ of~$\Sigma$,
and for every~$x\in{\mathcal{N}}$ there exists a unique point~$\pi_\Sigma(x)\in\Sigma$ (often called the ``projection\index{projection} of~$x$
onto~$\Sigma$'') such that
\begin{equation}\label{PROIE}
x=\pi_\Sigma(x)-d_\Sigma(x)\,\nu(\pi_\Sigma(x)),\end{equation}
moreover $\pi_\Sigma$ is of class~$C^2({\mathcal{N}})$ and
\begin{equation}\label{PROIE2} \nabla d_\Sigma(x)=-\nu(\pi_\Sigma(x)),\end{equation}
see\footnote{The intuition behind~\eqref{PROIE2}
is sketched in Figure~\ref{DItangeFI}. Roughly speaking,
one can consider a point~$p$
and measure its distance from~$\Sigma$ by considering the ball
centered at~$p$ and tangent to~$\Sigma$. Taking derivatives
of the distance function with respect to ``tangential directions''
corresponds to moving~$p$ infinitesimally
to the point~$p'$ and considering the ball
centered at~$p'$ and tangent to~$\Sigma$: since~$\Sigma$
detaches ``quadratically'' from its tangent hyperplane at~$p$,
this new ball is a small perturbation of the translation of the original ball
(thus producing a zero tangential derivative).

Instead,
taking normal derivatives
of the distance function
corresponds to moving~$p$ infinitesimally
to the point~$p''$ and considering the ball
centered at~$p''$ and tangent to~$\Sigma$:  
in this case, the new ball has a radius
equal to the one of the old ball, plus the distance between~$p$ and~$p''$,
up to small perturbations (and this produces a unit normal derivative).

The minus sign in~\eqref{PROIE2} is due to the fact that the outer normal of~$E$
points towards the region in which the sign distance is negative.}
e.g. Lemma~14.16 in~\cite{MR1814364} or Appendix~B
in~\cite{MR775682}. For our purposes~${\mathcal{N}}$ will be always supposed to be a conveniently
small neighborhood of~$\Sigma$.\medskip

Given a function~$u:\Sigma\to\R$, this framework allows us to
define the
``normal extension\index{normal extension} of~$u$ outside~$\Sigma$'' for all~$x\in{\mathcal{N}}$ as
\begin{equation}\label{DEEST} u_{\mbox{\scriptsize{ext}}}(x):= u(\pi_\Sigma(x)).\end{equation}
Notice that if~$p\in\Sigma$ and~$|t|$ is sufficiently small such that~$p+t\nu(p)\in{\mathcal{N}}$,
then~$
u_{\mbox{\scriptsize{ext}}}(p+t\nu(p))=u(p)$.
As a result, 
\begin{equation}\label{CE}\begin{split}&
0=\frac{d}{dt} u(p)=\frac{d}{dt} u_{\mbox{\scriptsize{ext}}}(p+t\nu(p))=\nabla
u_{\mbox{\scriptsize{ext}}}(p+t\nu(p))\cdot\nu(p)\\&\qquad\qquad=
\nabla
u_{\mbox{\scriptsize{ext}}}(p+t\nu(p))\cdot\nu_{\mbox{\scriptsize{ext}}}(p+t\nu(p)).\end{split}
\end{equation}
\medskip

The Laplace-Beltrami operator\index{Laplace-Beltrami operator} of a function~$u\in C^2(\Sigma)$ is then defined, for each~$p\in\Sigma$, by
\begin{equation}\label{BEL} \Delta_\Sigma u(p):=\Delta u_{\mbox{\scriptsize{ext}}}(p),\end{equation}
where~$\Delta$ represents here the standard Laplacian acting on functions in~$C^2({\mathcal{N}})$.
\medskip

Interestingly, the
Laplace-Beltrami operator is compatible with the
gradient structure intrinsic to~$\Sigma$. To this end, one defines
the ``tangential gradient''\index{tangential gradient} as the projection onto the tangent plane, namely,
for every~$f\in C^1({\mathcal{N}})$ and any~$p\in{\mathcal{N}}$,
\begin{equation}\label{IGRA}
\nabla_\Sigma f(p):=\nabla f(p)-\big(\nabla f(p)\cdot\nu_{\mbox{\scriptsize{ext}}}(p)\big)\,\nu_{\mbox{\scriptsize{ext}}}(p).\end{equation}
Also, if~$f\in C^1(\Sigma)$, we define its tangential gradient via the tangential
gradient of the normal extension, namely
\begin{equation}\label{INvjkds-3}
\nabla_\Sigma f:=\nabla_\Sigma f_{\mbox{\scriptsize{ext}}}.\end{equation}
We observe that, in view of~\eqref{CE}, on~$\Sigma$ we have that
\begin{equation}\label{INvjkds-1}
\nabla f_{\mbox{\scriptsize{ext}}}\cdot\nu=0
\end{equation}
and thus, by~\eqref{IGRA},
\begin{equation}\label{INvjkds-2} \nabla_\Sigma f_{\mbox{\scriptsize{ext}}}
=\nabla f_{\mbox{\scriptsize{ext}}}.\end{equation}

In analogy with the tangential gradient defined in~\eqref{IGRA},
one can introduce the ``tangential divergence''\index{tangential divergence}
of a vector field~$F\in C^1({\mathcal{N}},\R^n)$ at points of~${\mathcal{N}}$ as
\begin{equation}\label{VB-d1} \div_\Sigma F:=\div F-\nabla(F\cdot\nu_{\mbox{\scriptsize{ext}}})\cdot\nu_{\mbox{\scriptsize{ext}}}.\end{equation}
Roughly speaking, one is ``removing'' here the normal contribution
of the full divergence.
Also, if~$F\in C^1(\Sigma,\R^n)$, one defines its tangential 
divergence as the tangential divergence of its
normal extension, namely
\begin{equation}\label{VB-d2} \div_\Sigma F:=\div_\Sigma F_{\mbox{\scriptsize{ext}}},\end{equation}
where~$F_{\mbox{\scriptsize{ext}}}$ is the vector field obtained
by the normal extension of all the components of~$F$. In this
situation, we point out that, for each~$p\in{\mathcal{N}}$,
$$ F_{\mbox{\scriptsize{ext}}}(p)\cdot\nu_{\mbox{\scriptsize{ext}}}(p)=
F(\pi_\Sigma(p))\cdot\nu(\pi_\Sigma(p))=
(F\cdot\nu)_{\mbox{\scriptsize{ext}}}(p),$$
whence, in light of~\eqref{CE}, $\nabla
\big( F_{\mbox{\scriptsize{ext}}}\cdot\nu_{\mbox{\scriptsize{ext}}}\big)
\cdot\nu=0$ on~$\Sigma$.
Combining this with~\eqref{VB-d1} and~\eqref{VB-d2}, it follows that
\begin{equation}\label{0909} \div_\Sigma F=\div F_{\mbox{\scriptsize{ext}}}\qquad{\mbox{on }}\;\Sigma.\end{equation}

Concerning the definition of tangential gradient and divergence,
a caveat should be taken into account: namely,
given a function~$u\in C^1({\mathcal{N}})$
(or a vector field~$F\in C^1({\mathcal{N}},\R^n)$)
one can consider the restriction~$u\big|_\Sigma
\in C^1(\Sigma)$
(or~$F\big|_\Sigma\in C^1(\Sigma,\R^n)$)
and then compute the tangential gradient
of~$u\big|_\Sigma$ (or the tangential divergence of~$F\big|_\Sigma$) on~$\Sigma$, according to
definition~\eqref{INvjkds-3} (or definition~\eqref{VB-d2}),
that is using the normal extension
defined in~\eqref{DEEST}. The value obtained
in this way coincides with the tangential gradient of~$u$
computed via definition~\eqref{IGRA} (or the tangential divergence of~$F$ computed via definition~\eqref{VB-d1}) evaluated at~$\Sigma$. As a matter of fact, the values of a function
on~$\Sigma$ suffice to compute its tangential gradient
(as well as the values of a vector field
on~$\Sigma$ suffice to compute its tangential divergence), according to the following observation:

\begin{lemma}\label{PRIKD}
Let~$u$, $\widetilde u\in C^1({\mathcal{N}})$ be
such that~$u=\widetilde u$ on~$\Sigma$
and let~$\nabla_\Sigma u$ and~$\nabla_\Sigma \widetilde u$ be computed as in~\eqref{IGRA}. Then, on~$\Sigma$
we have that
\begin{equation}\label{FEU1}
\nabla_\Sigma u=
\nabla_\Sigma \widetilde u.\end{equation}

Furthermore,
let~$F$, $\widetilde F\in C^1({\mathcal{N}},\R^n)$ be
such that~$F=\widetilde F$ on~$\Sigma$
and let~$\div_\Sigma F$ and~$\div_\Sigma \widetilde F$ be computed as in~\eqref{VB-d1}. Then, on~$\Sigma$
we have that
\begin{equation}\label{FEU2}
\div_\Sigma F=
\div_\Sigma \widetilde F.\end{equation}
\end{lemma}

\begin{proof}
First of all, we observe that, if~$F=(F_1,\dots,F_n)$, then
\begin{equation}\label{FEU}
\div_\Sigma F=\sum_{j=1}^n\nabla_\Sigma F_j\cdot e_j.
\end{equation}
Indeed, using~\eqref{IGRA} and~\eqref{VB-d1}, and then also~\eqref{CE},
we see that
\begin{eqnarray*}&&
\div_\Sigma F-\sum_{j=1}^n\nabla_\Sigma F_j\cdot e_j\\&=&
\sum_{j=1}^n\partial_j F_j
-\nabla\left(\sum_{j=1}^nF_je_j\cdot\nu_{\mbox{\scriptsize{ext}}}\right)\cdot\nu_{\mbox{\scriptsize{ext}}}
-\sum_{j=1}^n
\left(\nabla F_j-\big(\nabla F_j\cdot\nu_{\mbox{\scriptsize{ext}}}\big)\,\nu_{\mbox{\scriptsize{ext}}}
\right)\cdot e_j\\&=&
-\nabla\left(\sum_{j=1}^nF_je_j\cdot\nu_{\mbox{\scriptsize{ext}}}\right)\cdot\nu_{\mbox{\scriptsize{ext}}}
+\sum_{j=1}^n
\big(\nabla F_j\cdot\nu_{\mbox{\scriptsize{ext}}}\big)\big(\nu_{\mbox{\scriptsize{ext}}}
\cdot e_j\big)\\&=&
-\sum_{j=1}^nF_j\nabla\big(e_j\cdot\nu_{\mbox{\scriptsize{ext}}}\big)\cdot\nu_{\mbox{\scriptsize{ext}}}\\&=&0,
\end{eqnarray*}
thus proving~\eqref{FEU}.

Now we prove~\eqref{FEU1}. For this, we set~$w:=u-\widetilde u$
and we remark that~$w=0$ on~$\Sigma$. Given a point~$p$ of~$\Sigma$, we suppose
that in a neighborhood of~$p$ the hypersurface~$\Sigma$ is parameterized by the graph
of a function~$\psi:\R^{n-1}\to\R$ (up to renumbering the variables, we also assume that this
graph occurs in the $n$th coordinate direction, with the set~$E$
lying above the graph), namely there exists~$r>0$ such that
\begin{equation} \label{uojw-S29-32jfewnb}
B_r(p)\cap E=\{x_n>\psi(x')\}\cap B_r(p).\end{equation}
Notice that, on~$\Sigma$,
\begin{equation}\label{NBORM} \nu=\frac{\left( \nabla'\psi,-1\right)}{\sqrt{1+|\nabla' \psi|^2}},\end{equation}
where the notation
\begin{equation}\label{NABLEP}
\nabla':=(\partial_1,\dots,\partial_{n-1})\end{equation}
has been used.

In this way, in the vicinity of~$p$ we can write that~$w(x',\psi(x'))=0$ and
$$ 0=\nabla' \Big(w(x',\psi(x'))\Big)=\nabla'w(x',\psi(x'))+\partial_nw(x',\psi(x'))\,\nabla'\psi(x').
$$
Consequently, using~\eqref{IGRA}, we find that, on~$\Sigma$,
in the vicinity of~$p$,
\begin{eqnarray*}
\nabla_\Sigma u-\nabla_\Sigma\widetilde u&=&
\nabla w-\big(\nabla w\cdot\nu\big)\,\nu\\&=&
\Big( -\partial_nw\,\nabla'\psi,\partial_nw\Big)-\left(
\Big( -\partial_nw\,\nabla'\psi,\partial_nw\Big)\cdot
\frac{\left( \nabla'\psi,-1\right)}{\sqrt{1+|\nabla' \psi|^2}}\right)\frac{\left( \nabla'\psi,-1\right)}{\sqrt{1+|\nabla' \psi|^2}}\\
&=&\Big( -\partial_nw\,\nabla'\psi,\partial_nw\Big)+
\frac{\partial_nw\left( |\nabla'\psi|^2+1\right)}{\sqrt{1+|\nabla' \psi|^2}}\frac{\left( \nabla'\psi,-1\right)}{\sqrt{1+|\nabla' \psi|^2}}\\
&=&\Big( -\partial_nw\,\nabla'\psi,\partial_nw\Big)+
\Big(\partial_nw\, \nabla'\psi,-\partial_nw\Big)\\
&=&0,
\end{eqnarray*}
which establishes~\eqref{FEU1}.

To prove~\eqref{FEU2}, we exploit~\eqref{FEU1} (applied to the scalar component
functions~$F_j$ and~$\widetilde F_j$) and~\eqref{FEU}
to compute that
\begin{eqnarray*}
&&\div_\Sigma F-
\div_\Sigma \widetilde F
=\sum_{j=1}^n\nabla_\Sigma F_j\cdot e_j-
\sum_{j=1}^n\nabla_\Sigma \widetilde F_j\cdot e_j
=0.
\end{eqnarray*}
This completes the proof of~\eqref{FEU2}.
\end{proof}

In a nutshell, the content of Lemma~\ref{PRIKD}
is that different extensions of a smooth object defined only on~$\Sigma$ do not alter the tangential ``first order''
operators, since the tangent hyperplane of~$\Sigma$ ``detaches quadratically''
from~$\Sigma$ (we will find however that
``second order'' operators are sensitive to different
types of extensions, see~\eqref{MN:ikmc} below).
\medskip

We observe that the Laplace-Beltrami operator possesses
a ``tangential divergence form structure''\index{divergence form}, to be compared with
the classical one in~\eqref{DITRO}, namely:

\begin{lemma}\label{GB}
For every~$u\in C^2(\Sigma)$, on~$\Sigma$ we have that
$$ \Delta_\Sigma u=\div_\Sigma(\nabla_\Sigma u).$$
\end{lemma}

\begin{proof} Using in order~\eqref{BEL}, \eqref{INvjkds-3}, \eqref{INvjkds-2},
\eqref{VB-d1} and~\eqref{CE},
we see that, on~$\Sigma$,
\begin{equation*}\begin{split}&
\Delta_\Sigma u-\div_\Sigma(\nabla_\Sigma u)=
\Delta u_{\mbox{\scriptsize{ext}}}-\div_\Sigma
(\nabla_\Sigma u_{\mbox{\scriptsize{ext}}})=\Delta u_{\mbox{\scriptsize{ext}}}-\div_\Sigma
(\nabla u_{\mbox{\scriptsize{ext}}})\\&\qquad=
\Delta u_{\mbox{\scriptsize{ext}}}-
\div (\nabla u_{\mbox{\scriptsize{ext}}})-\nabla(\nabla u_{\mbox{\scriptsize{ext}}}
\cdot\nu_{\mbox{\scriptsize{ext}}})\cdot\nu
=0.
\qedhere\end{split}\end{equation*}
\end{proof}

For further use, it is now useful to recall an asymptotic
result about the sets obtained by ``thickening''~$\Sigma$.
For general and precise formulas
computing tubular neighborhoods of hypersurfaces
see e.g. \cite[Theorem~1]{MR2984315}.

\begin{lemma}\label{243}
Let~$\alpha$ be a continuous function on~$\Sigma$
and~$\beta$ be a continuous function on~${\mathcal{N}}$.
Let~$\e>0$ and
\begin{equation}\label{KMSD-fpblf} \Sigma_\e(\alpha):=\big\{
p+t\nu(p),\;p\in\Sigma,\;0\le t\le \e\alpha(p)
\big\}.\end{equation}
Then, as~$\e\searrow0$,
\begin{equation}\label{dxs} \int_{\Sigma_\e(\alpha)} \beta(y)\,dy= \e
\int_\Sigma \alpha_+(p)\,\beta(p)\,d{\mathcal{H}}^{n-1}_p+o(\e),
\end{equation}
where~$\alpha_+(p):=\max\{\alpha(p),\,0\}$.
\end{lemma}

\begin{figure}
  \centering
  \includegraphics[width=.6\linewidth]{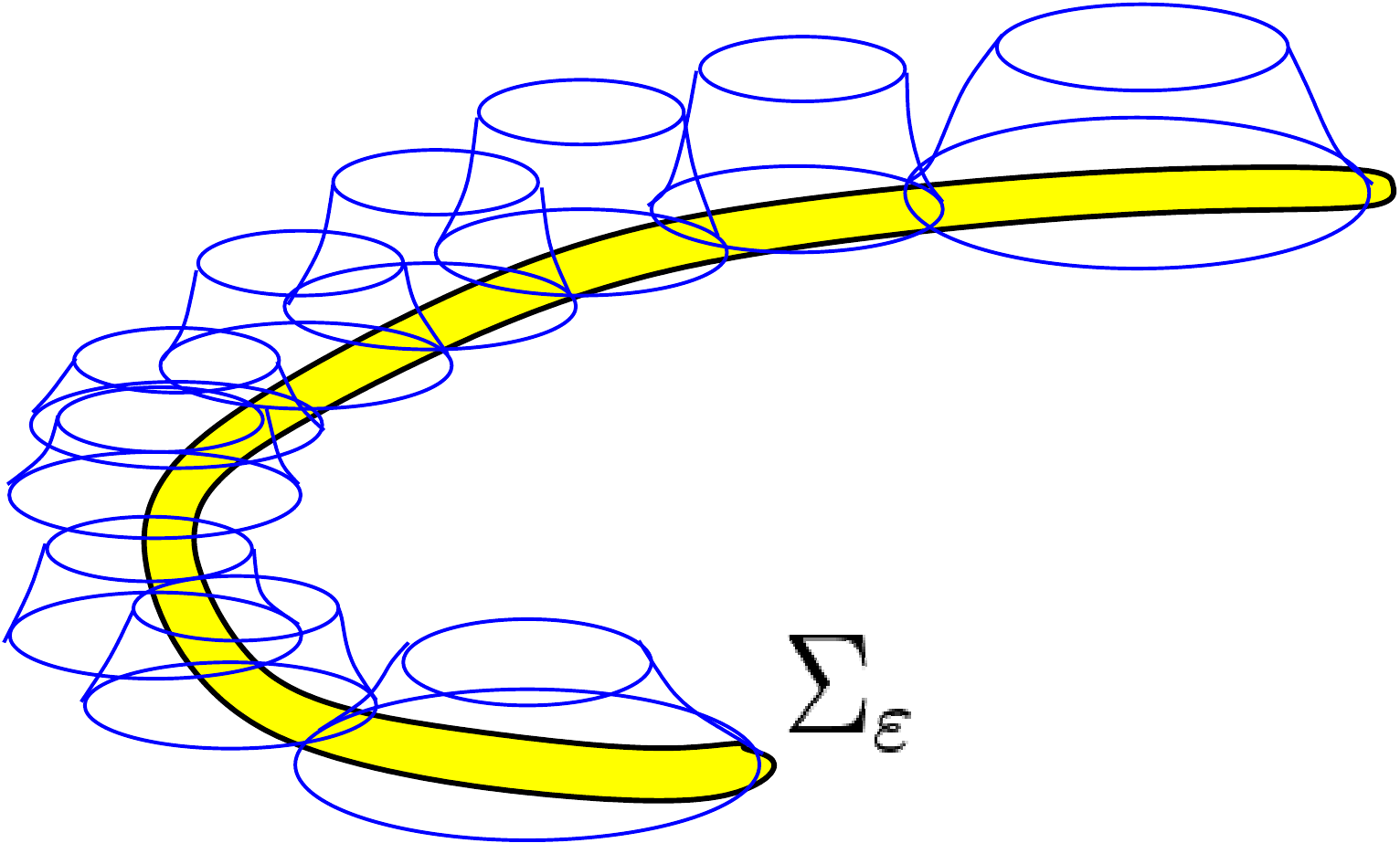}
 \caption{\sl Local charts for~$\Sigma$ and a partition of unity.}\label{PPPP}
\end{figure}

\begin{proof} In local coordinates, we write a surface element of~$\Sigma$
as a graph of a function
\begin{equation}\label{psi9itg}\psi:U\subseteq\R^{n-1}\to\R,\end{equation}
say in the $n$th direction, with normal as in~\eqref{NBORM}.
In this way, points~$y$ in this element of~$\Sigma_\e(\alpha)$ are of the form
\begin{equation}\label{psi9itg2} y=(x',\psi(x'))+
\frac{t\,\left( \nabla'\psi(x'),-1\right)}{\sqrt{1+|\nabla' \psi(x')|^2}}
,\qquad{\mbox{
with~$x'\in U$ and~$0<t<\e\alpha(x',\psi(x'))$.}}\end{equation}
That is, one can consider a partition of unity
(see e.g.~\cite[page~192]{MR0426007}) made of functions~$\phi_i\in C^\infty_0({\mathcal{N}},\,[0,1])$
with~${i\in\N}$ and finite overlapping supports, each compactly contained
in a local chart of~$\Sigma$,
such that~$\sum_{i\in\N}\phi_i=1$ in a given neighborhood~${\mathcal{N}}'$
of~$\Sigma$ (with~$\Sigma\subseteq{\mathcal{N}}'\Subset{\mathcal{N}}$, see Figure~\ref{PPPP}).
Then, letting~$\beta_i:=\beta\phi_i$, it suffices to prove~\eqref{dxs} with~$\beta$ replaced by~$\beta_i$,
since
\begin{equation*} \int_{\Sigma_\e(\alpha)} \beta(y)\,dy=
\sum_{i\in\N}\int_{\Sigma_\e(\alpha)} \beta_i(y)\,dy
\qquad{\mbox{and}}\qquad
\int_\Sigma \alpha_+(p)\,\beta(p)\,d{\mathcal{H}}^{n-1}_p=
\sum_{i\in\N}\int_\Sigma \alpha_+(p)\,\beta_i(p)\,d{\mathcal{H}}^{n-1}_p.
\end{equation*}
Therefore, from now on, to prove~\eqref{dxs}, up to replacing~$\beta$
with~$\beta_i$, we can suppose that, in the support of~$\beta$,
the hypersurface~$\Sigma$ is a graph of a function~$\psi$ as in~\eqref{psi9itg2},
say in the $n$th direction, and, for small~$\e$, the tubular neighborhood~$ {\Sigma_\e(\alpha)}$
can be written as the set of points~$y$ in~\eqref{psi9itg2}.

For convenience, one can denote~$x_n:=t$ and~$x:=(x',x_n)$
in~\eqref{psi9itg2}, and thus describe~$ {\Sigma_\e(\alpha)}$ in the support of~$\beta$
as the collection of points
$$ y=(x',\psi(x'))+
\frac{x_n\,\left( \nabla'\psi(x'),-1\right)}{\sqrt{1+|\nabla' \psi(x')|^2}}
,$$
with~$x\in U\times[0,\e\alpha(x',\psi(x'))]$
(if the latter quantity is well-defined, i.e. if~$\alpha(x',\psi(x'))\ge0$).
Notice in particular that~$x_n=O(\e)$ and we thus consider,
for small~$\e$, the change of variable relating~$y$
and~$x$, with
$$ \frac{\partial y}{\partial x}=\left(
\begin{matrix}
\partial_{x_1} y_1&\;\dots\;&\partial_{x_{n-1}} y_1&\;\,\partial_{x_n} y_1
\\ &\ddots\\
\partial_{x_1} y_{n-1}&
\;\dots\;&\partial_{x_{n-1}} y_{n-1}&\;\,\partial_{x_n} y_{n-1}\\
\partial_{x_1} y_{n}&\;\dots\;&\partial_{x_{n-1}} y_n&\;\,\partial_{x_n} y_n
\end{matrix}
\right)
=\left(
\begin{matrix}
1&\;\dots\;0&\;\,\partial_{x_1} \psi/R
\\ &\ddots\\
0&\;\dots\;1&\;\,\partial_{x_{n-1}} \psi/R\\
\partial_{x_1} \psi&\;\dots\;\partial_{x_{n-1}} \psi&\;\,-1/R
\end{matrix}
\right)+O(\e),$$
with~$R:=\sqrt{1+|\nabla' \psi(x')|^2}$.

As a result,
$$ \left|\det\frac{\partial y}{\partial x}
\right|=\frac{
(\partial_{x_1} \psi)^2+\dots(\partial_{x_{n-1}} \psi)^2+1}{R}
+O(\e)=
R+O(\e).$$
We remark that~$R\,dx'$ is the surface element on~$\Sigma$
(see e.g.~\cite[page 125]{MR3409135}), hence we write
the volume element of~\eqref{dxs} in the form
$$ dy=\left|\det\frac{\partial y}{\partial x}
\right|\,dx
=(R(x')+O(\e))\,dx=
d{\mathcal{H}}^{n-1}_p\,dx_n+
O(\e)\,dx.$$
That is,
\begin{eqnarray*}
&&\int_{\Sigma_\e(\alpha)} \beta(y)\,dy=
\int_{\Sigma_\e(\alpha)} \big(\beta(\pi_\Sigma(y))+o(1)\big)\,dy\\&&\qquad=
\int_\Sigma\left[ \int_0^{\e\alpha_+(p)}
\big(\beta(p)+o(1)\big)\,dx_n\right]\,d{\mathcal{H}}^{n-1}_p
+o(\e)\\&&\qquad=\e
\int_\Sigma \alpha_+(p)\,\beta(p)\,d{\mathcal{H}}^{n-1}_p+o(\e),
\end{eqnarray*}
giving~\eqref{dxs}
as desired.\end{proof}

The tangential differential setting provides another useful form of ``integration by parts formula''
according to Theorem~\ref{GBNAmn} below. Differently from the Euclidean case,
this result takes into account an additional term coming from the geometry of~$\Sigma$.
For this,
it is useful to introduce the mean curvature\index{mean curvature}
at a point~$x\in\Sigma$, defined as
\begin{equation}\label{MC}
H(x):=\div_\Sigma \nu(x).
\end{equation}
See e.g. Section~1.2 in~\cite{MR3230079} for a geometric
description of the mean curvature.
Then, we have the following result, sometimes called the ``Tangential Divergence Theorem''\index{Tangential Divergence Theorem}:

\begin{theorem}\label{GBNAmn}
For every~$F\in C^1(\Sigma,\R^n)$ and~$\varphi\in C^1(\Sigma)$,
$$ \int_\Sigma \div_\Sigma F(x)\;\varphi(x)\,d{\mathcal{H}}^{n-1}_x=
\int_\Sigma F(x)\cdot\Big(H(x)\nu(x)\varphi(x)-\nabla_\Sigma\varphi(x)\Big)\,d{\mathcal{H}}^{n-1}_x.$$
\end{theorem}

\begin{proof} Given~$\e>0$, to be taken conveniently small,
we consider a tubular neighborhood of~$\Sigma$ of radius~$\e$, namely we set~$
\Sigma_\e$ as in~\eqref{KMSD-fpblf} with~$\alpha:=1$.

\begin{figure}
  \centering
  \includegraphics[width=.7\linewidth]{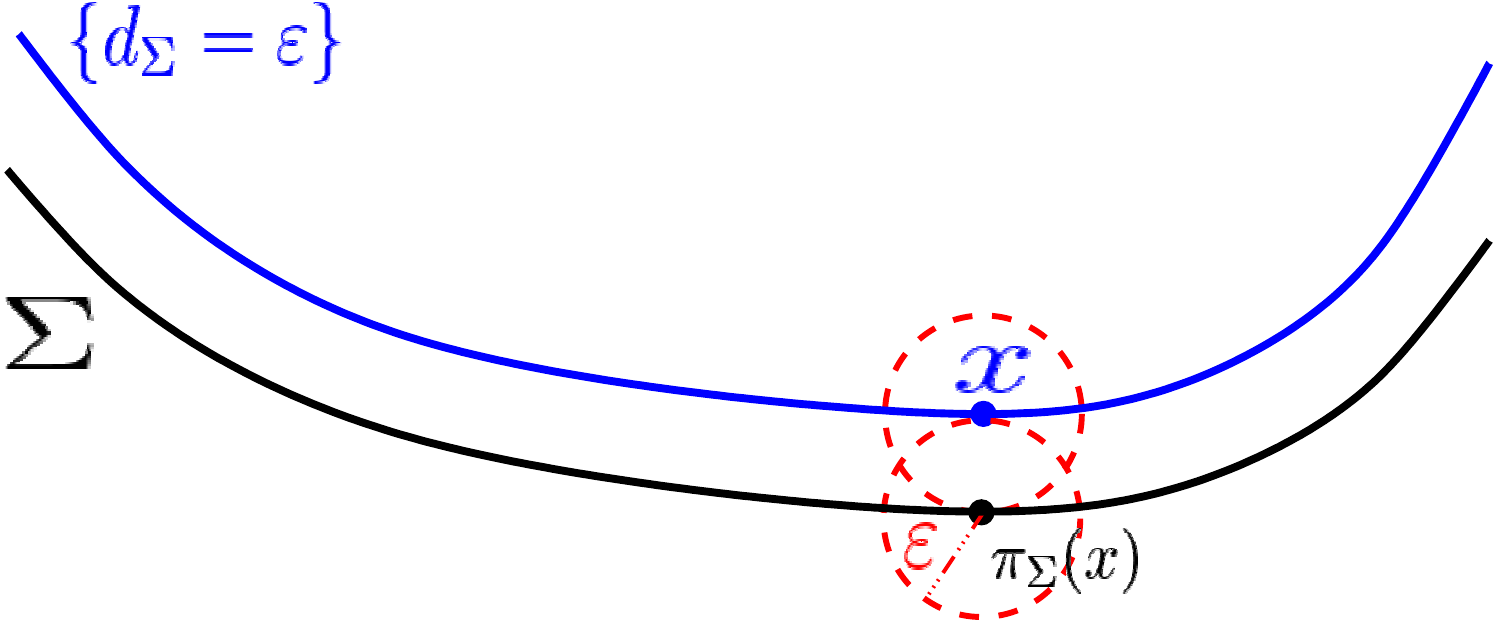}
 \caption{\sl The hypersurface~$\Sigma$ and the ``parallel hypersurface''
 at distance~$\varepsilon$.}\label{DIFI}
\end{figure}

Let us now analyze the exterior normal~$\nu_{\Sigma_\e}$
along~$\partial \Sigma_\e$.
We stress that~$\partial \Sigma_\e=\{d_\Sigma=\e\}\cup\{d_\Sigma=-\e\}$
and thus we denote by~$\nu^{(\pm)}_\e$ the exterior normal
of~$\Sigma_\e$ along~$\{d_\Sigma=\pm\e\}$, and, for clarity,
by~$\nu_\Sigma$ the exterior normal of~$E$ along~$\Sigma$, see Figure~\ref{DIFI}.
In this way, if~$x\in\{d_\Sigma=\e\}$ the exterior normal~$\nu^{(+)}_\e$
at~$x\in\{d_\Sigma=\e\}$ is minus
the exterior normal~$\nu_\Sigma$ of~$E$ at~$\pi_\Sigma(x)$,
while the
exterior normal~$\nu^{(-)}_\e$ at~$x\in\{d_\Sigma=-\e\}$ is plus
the exterior normal~$\nu_\Sigma$ of~$E$ at~$\pi_\Sigma(x)$, that is
\begin{equation}\label{221}
\begin{split}&
\nu^{(+)}_\e(x)=-\nu_\Sigma(\pi_\Sigma(x))=-\nu_{\mbox{\scriptsize{ext}}}(x)\qquad{\mbox{if }}x\in\{d_\Sigma=\e\}\\
{\mbox{and }}\qquad&
\nu^{(-)}_\e(x)=\nu_\Sigma(\pi_\Sigma(x))=\nu_{\mbox{\scriptsize{ext}}}(x)\qquad{\mbox{if }}x\in\{d_\Sigma=-\e\}.
\end{split}\end{equation}
Moreover, by~\eqref{dxs},
used here with~$\alpha:=1$
and~$\beta:=
\div F_{\mbox{\scriptsize{ext}}}\;\varphi
+ F\cdot\nabla\varphi_{\mbox{\scriptsize{ext}}}$,
\begin{equation}\label{BGio}
\begin{split}
&\int_\Sigma \div_\Sigma F(x)\;\varphi(x)\,d{\mathcal{H}}^{n-1}_x+\int_\Sigma F(x)\cdot\nabla_\Sigma\varphi(x)\,d{\mathcal{H}}^{n-1}_x\\=\,&
\int_\Sigma \Big( \div F_{\mbox{\scriptsize{ext}}}(x)\;\varphi(x)
+ F(x)\cdot\nabla\varphi_{\mbox{\scriptsize{ext}}}(x)\Big)\,d{\mathcal{H}}^{n-1}_x
\\=\,&
\lim_{\e\searrow0}\frac1{\e}\,
\int_{\Sigma_\e} \Big( \div F_{\mbox{\scriptsize{ext}}}(x)\;\varphi_{\mbox{\scriptsize{ext}}}
(x)+ F_{\mbox{\scriptsize{ext}}}(x)\cdot
\nabla\varphi_{\mbox{\scriptsize{ext}}}(x)\Big)\,dx\\=\,&
\lim_{\e\searrow0}\frac1{\e}\,
\int_{\Sigma_\e} \div \Big(
\varphi_{\mbox{\scriptsize{ext}}}
(x)\;F_{\mbox{\scriptsize{ext}}}(x)\Big)\,dx\\=\,&\lim_{\e\searrow0}\frac1{\e}\,
\int_{\Sigma_\e} \div G_{\mbox{\scriptsize{ext}}}(x)\,dx
,\end{split}\end{equation}
where~$G:=\varphi F$.

We also point out that, if~$\widetilde G:=G\cdot\nu$,
\begin{eqnarray*}&&
\div\Big( (G_{\mbox{\scriptsize{ext}}}\cdot\nu_{\mbox{\scriptsize{ext}}})\nu_{\mbox{\scriptsize{ext}}}\Big)-
(G_{\mbox{\scriptsize{ext}}}\cdot\nu_{\mbox{\scriptsize{ext}}})\div\nu_{\mbox{\scriptsize{ext}}}=
\nabla(G_{\mbox{\scriptsize{ext}}}\cdot\nu_{\mbox{\scriptsize{ext}}})\cdot\nu_{\mbox{\scriptsize{ext}}}=\nabla \widetilde G_{\mbox{\scriptsize{ext}}}\cdot\nu_{\mbox{\scriptsize{ext}}}=0,
\end{eqnarray*}
thanks to~\eqref{CE}. Consequently,
\begin{eqnarray*}&&
\int_\Sigma H(x) F(x)\cdot\nu(x)\varphi(x)\,d{\mathcal{H}}^{n-1}_x=
\int_\Sigma (G(x)\cdot\nu(x)) \div\nu_{\mbox{\scriptsize{ext}}}(x) \,d{\mathcal{H}}^{n-1}_x\\&&\qquad
=\int_\Sigma(G_{\mbox{\scriptsize{ext}}}(x)\cdot\nu_{\mbox{\scriptsize{ext}}}(x)) \div\nu_{\mbox{\scriptsize{ext}}}(x) \,d{\mathcal{H}}^{n-1}_x=\lim_{\e\searrow0}\frac1{
\e}\,
\int_{\Sigma_\e}(G_{\mbox{\scriptsize{ext}}}(x)\cdot\nu_{\mbox{\scriptsize{ext}}}(x))\div\nu_{\mbox{\scriptsize{ext}}}(x)\,dx\\&&\qquad=\lim_{\e\searrow0}\frac1{
\e}\,
\int_{\Sigma_\e}
\div\Big( (G_{\mbox{\scriptsize{ext}}}(x)\cdot\nu_{\mbox{\scriptsize{ext}}}(x)
)\nu_{\mbox{\scriptsize{ext}}}(x)\Big)\,dx.
\end{eqnarray*}
This and~\eqref{BGio}, together with~\eqref{221}, lead to
\begin{equation*}\begin{split}&
\int_\Sigma \div_\Sigma F(x)\;\varphi(x)\,d{\mathcal{H}}^{n-1}_x+
\int_\Sigma F(x)\cdot\nabla_\Sigma\varphi(x)\,d{\mathcal{H}}^{n-1}_x
-\int_\Sigma H(x) F(x)\cdot\nu(x)\varphi(x)\,d{\mathcal{H}}^{n-1}_x\\=\,&
\lim_{\e\searrow0}\frac1{\e}\,
\int_{\Sigma_\e} \div\Big( G_{\mbox{\scriptsize{ext}}}(x)-
(G_{\mbox{\scriptsize{ext}}}(x)
\cdot\nu_{\mbox{\scriptsize{ext}}}(x))\nu_{\mbox{\scriptsize{ext}}}(x)\Big)\,dx\\=\,&
\lim_{\e\searrow0}\frac1{\e}\,
\int_{\partial\Sigma_\e}
\Big( G_{\mbox{\scriptsize{ext}}}(x)-
(G_{\mbox{\scriptsize{ext}}}(x)
\cdot\nu_{\mbox{\scriptsize{ext}}}(x))\nu_{\mbox{\scriptsize{ext}}}(x)\Big)
\cdot\nu_{\Sigma_\e}(x)
\,d{\mathcal{H}}^{n-1}_x\\=\,&
\lim_{\e\searrow0}\frac1{\e}\,
\Bigg(-\int_{\{d_\Sigma=\e\}}
\Big( G_{\mbox{\scriptsize{ext}}}(x)-
(G_{\mbox{\scriptsize{ext}}}(x)
\cdot\nu_{\mbox{\scriptsize{ext}}}(x))\nu_{\mbox{\scriptsize{ext}}}(x)\Big)
\cdot\nu_{\mbox{\scriptsize{ext}}}(x)
\,d{\mathcal{H}}^{n-1}_x\\&\qquad+
\int_{\{d_\Sigma=-\e\}}
\Big( G_{\mbox{\scriptsize{ext}}}(x)-
(G_{\mbox{\scriptsize{ext}}}(x)
\cdot\nu_{\mbox{\scriptsize{ext}}}(x))\nu_{\mbox{\scriptsize{ext}}}(x)\Big)
\cdot\nu_{\mbox{\scriptsize{ext}}}(x)
\,d{\mathcal{H}}^{n-1}_x
\Bigg)
\\=\,&
\lim_{\e\searrow0}\frac1{\e}\,
\Bigg(-\int_{\{d_\Sigma=\e\}}
\Big( (G_{\mbox{\scriptsize{ext}}}(x)
\cdot\nu_{\mbox{\scriptsize{ext}}}(x))-(G_{\mbox{\scriptsize{ext}}}(x)
\cdot\nu_{\mbox{\scriptsize{ext}}}(x))\Big)
\,d{\mathcal{H}}^{n-1}_x\\&\qquad+
\int_{\{d_\Sigma=-\e\}}
\Big( (G_{\mbox{\scriptsize{ext}}}(x)
\cdot\nu_{\mbox{\scriptsize{ext}}}(x))-(G_{\mbox{\scriptsize{ext}}}(x)
\cdot\nu_{\mbox{\scriptsize{ext}}}(x))\Big)
\,d{\mathcal{H}}^{n-1}_x
\Bigg)
\\=\,&\lim_{\e\searrow0}\frac1{\e}\,\big(0-0\big)
\\=\,&0.\qedhere\end{split}
\end{equation*}
\end{proof}

The content of Theorem~\ref{GBNAmn} happens to be a special case
of some ``Minkowski Integral Formulas''\index{Minkowski Integral Formula}, see e.g.~\cite{MR2522595, MR3382197}.
For instance, a classical identity is the following one:

\begin{corollary}\label{MIK}
We have that
$$ {\mathcal{H}}^{n-1}(\Sigma)=\frac1{n-1}
\int_\Sigma H(x)\,x\cdot\nu(x)\,d{\mathcal{H}}^{n-1}_x.$$
\end{corollary}

\begin{proof} For every~$x\in{\mathcal{N}}$,
let~$F(x):=x$. 
Then, by Lemma~\ref{PRIKD} and~\eqref{VB-d1}, recalling also~\eqref{CE}, for~$x\in\Sigma$ we have that
\begin{eqnarray*}&&
\div_\Sigma F(x)=\div F(x)-\nabla(F(x)\cdot\nu_{\mbox{\scriptsize{ext}}}(x))\cdot\nu(x)=
n-\nabla(x\cdot\nu_{\mbox{\scriptsize{ext}}}(x))\cdot\nu(x)\\&&\qquad=n-1-\sum_{j=1}^n
x_j\nabla(\nu_{\mbox{\scriptsize{ext}}}(x)\cdot e_j)\cdot\nu(x)=n-1.
\end{eqnarray*}
The desired result follows by exploiting the Tangential Divergence Theorem
(see Theorem~\ref{GBNAmn})
with~$\varphi:=1$.
\end{proof}

{F}rom the Tangential Divergence Theorem~\ref{GBNAmn}, one also obtains:

\begin{corollary}
For every~$u\in C^2(\Sigma)$ and~$\varphi\in C^1(\Sigma)$,
\begin{equation}\label{PAR-091} \int_\Sigma \Delta_\Sigma u(x)\,\varphi(x)\,d{\mathcal{H}}^{n-1}_x=
-\int_\Sigma \nabla_\Sigma u(x)\cdot
\nabla_\Sigma\varphi(x)\,d{\mathcal{H}}^{n-1}_x.\end{equation}
Moreover, for every~$u\in C^2(\Sigma)$ and~$\varphi\in C^2(\Sigma)$,
\begin{equation}\label{PAR-092}
\int_\Sigma \Delta_\Sigma u(x)\,\varphi(x)\,d{\mathcal{H}}^{n-1}_x=
\int_\Sigma u(x)\,\Delta_\Sigma\varphi(x)\,d{\mathcal{H}}^{n-1}_x.\end{equation}
\end{corollary}

\begin{proof} Let~$u\in C^2(\Sigma)$
and~$F:=\nabla_\Sigma u$. Then, by~\eqref{INvjkds-3}, \eqref{INvjkds-1}
and~\eqref{INvjkds-2}, on~$\Sigma$ we have that
$$ F\cdot\nu=\nabla_\Sigma u_{\mbox{\scriptsize{ext}}}\cdot\nu=\nabla
u_{\mbox{\scriptsize{ext}}}\cdot\nu=0.$$
Hence, 
if~$\varphi\in C^1(\Sigma)$,
by Lemma~\ref{GB} and Theorem~\ref{GBNAmn},
\begin{eqnarray*}&&
\int_\Sigma \Delta_\Sigma u(x)\,\varphi(x)\,d{\mathcal{H}}^{n-1}_x+
\int_\Sigma \nabla_\Sigma u(x)\cdot
\nabla_\Sigma\varphi(x)\,d{\mathcal{H}}^{n-1}_x\\&=&
\int_\Sigma \div_\Sigma F(x)\,\varphi(x)\,d{\mathcal{H}}^{n-1}_x+
\int_\Sigma F(x)\cdot
\nabla_\Sigma\varphi(x)\,d{\mathcal{H}}^{n-1}_x\\&=&
\int_\Sigma H(x)\varphi(x)F(x)\cdot\nu(x)
\,d{\mathcal{H}}^{n-1}_x\\&=&0.
\end{eqnarray*}
This establishes~\eqref{PAR-091}.

Now we suppose that also~$\varphi\in C^2(\Sigma)$,
and we write~\eqref{PAR-091} exchanging the roles of~$u$ and~$\varphi$,
namely
$$ \int_\Sigma \Delta_\Sigma \varphi(x)\,u(x)\,d{\mathcal{H}}^{n-1}_x=
-\int_\Sigma \nabla_\Sigma \varphi(x)\cdot
\nabla_\Sigma u(x)\,d{\mathcal{H}}^{n-1}_x.$$
{F}rom this and~\eqref{PAR-091}, the desired result in~\eqref{PAR-092}
plainly follows.
\end{proof}

Now we give an explicit formula for the mean curvature\index{mean curvature}
with respect to the function that describes locally the hypersurface~$\Sigma$.
We use the notation in~\eqref{NABLEP}
and~$\div':=\partial_1+\dots+\partial_{n-1}$. With this, we have:

\begin{theorem}\label{INCOO}
Let~$p=(p',p_n)\in \Sigma$.
Assume that there exists~$\rho>0$
such that~$E$ in~$B_\rho(p)$ can be written as a supergraph
of a function~$\psi:\R^{n-1}\to\R$, namely
$$ E\cap B_\rho(p)=\{ x_n>\psi(x')\}\cap B_\rho(p).$$
Then, at~$p$,
\begin{equation}\label{LAMEANCURVA} H=\div'\left(
\frac{\nabla' \psi}{\sqrt{1+|\nabla' \psi|^2}}
\right).\end{equation}
\end{theorem}

\begin{proof} Given~$x=(x',x_n)\in\R^n$, we define
$$ \nu_\star(x):=\frac{\left( \nabla'\psi(x'),-1\right)}{\sqrt{1+|\nabla' \psi(x')|^2}}.$$
By~\eqref{NBORM}, we know that~$\nu_\star$ coincides with~$\nu$
(and thus with~$ \nu_{\mbox{\scriptsize{ext}}}$) on~$\Sigma$ in the vicinity of~$p$
and therefore,
using~\eqref{MC}
and successively~\eqref{VB-d1}, \eqref{VB-d2}
and~\eqref{FEU2}, we have that,
for all~$x\in\Sigma$
in the vicinity of~$p$,
\begin{equation}\label{kPmsi}
H(x)=\div_\Sigma \nu(x)
=\div_\Sigma \nu_{\mbox{\scriptsize{ext}}}(x)=\div_\Sigma \nu_\star=
\div \nu_\star-\nabla(\nu_\star\cdot\nu_{\mbox{\scriptsize{ext}}})\cdot\nu.
\end{equation}
We also observe that~$\nu_\star(x)$ does not depend on~$x_n$ and therefore
\begin{equation}\label{LKJNnbvdceij84} \div \nu_\star=
\div'\left(
\frac{\nabla' \psi}{\sqrt{1+|\nabla' \psi|^2}}
\right).\end{equation}
Furthermore, for every~$j\in\{1,\dots,n\}$,
$$ 0=\partial_j \frac12=\partial_j\frac{ |\nu_\star|^2}2=
\frac12\sum_{m=1}^n \partial_j (\nu_\star\cdot e_m)^2=
\sum_{m=1}^n (\nu_\star\cdot e_m)\,
\partial_j(\nu_\star\cdot e_m).
$$
That is, on~$\Sigma$,
$$ \sum_{m=1}^n (\nu\cdot e_m)\,
\partial_j(\nu_\star\cdot e_m)=0.$$
Similarly, for every~$j\in\{1,\dots,n\}$, on~$\Sigma$,
$$ \sum_{m=1}^n (\nu\cdot e_m)\,
\partial_j(\nu_{\mbox{\scriptsize{ext}}}\cdot e_m)=0.$$
Hence, on~$\Sigma$,
\begin{eqnarray*}
\nabla(\nu_\star\cdot\nu_{\mbox{\scriptsize{ext}}})\cdot\nu&=&\sum_{j=1}^n
\partial_j(\nu_\star\cdot\nu_{\mbox{\scriptsize{ext}}})\,(\nu\cdot e_j)
\\&=&
\sum_{j,m=1}^n
\partial_j\Big(
(\nu_\star\cdot e_m)(\nu_{\mbox{\scriptsize{ext}}}\cdot e_m)\Big)\,(\nu\cdot e_j)
\\&=&
\sum_{j,m=1}^n\left[
(\nu_{\mbox{\scriptsize{ext}}}\cdot e_m)(\nu\cdot e_j)\,
\partial_j
(\nu_\star\cdot e_m)
+(\nu_\star\cdot e_m)(\nu\cdot e_j)\,\partial_j(\nu_{\mbox{\scriptsize{ext}}}\cdot e_m)\right]\\&=&
\sum_{j,m=1}^n\left[
(\nu\cdot e_m)(\nu\cdot e_j)\,
\partial_j
(\nu_\star\cdot e_m)
+(\nu\cdot e_m)(\nu\cdot e_j)\,\partial_j(\nu_{\mbox{\scriptsize{ext}}}\cdot e_m)\right]\\&=&0.
\end{eqnarray*}
This, \eqref{kPmsi}
and~\eqref{LKJNnbvdceij84} yield the desired result.
\end{proof}

As a consequence of Theorem~\ref{INCOO}, we also point out that
the mean curvature is monotone with respect to set inclusions:

\begin{corollary}\label{ILCOMS}
Let~$E$, $F\subseteq\R^n$ be bounded and open sets
with~$C^3$ boundary, with mean curvature~$H_{\partial E}$ and~$H_{\partial F}$, respectively.
Assume that~$E\subseteq F$ and that~$p\in (\partial E)\cap(\partial F)$. Then,~$H_{\partial E}(p)\ge
H_{\partial F}(p)$.
\end{corollary}

\begin{proof} We remark that the external normal~$\nu_{\partial E}$ of~$E$
coincides with the external normal~$\nu_{\partial F}$ of~$F$ at~$p$,
hence we can write both~$E$ and~$F$ as graphs with respect to a common coordinate direction.
That is, without loss of generality, we may assume that there exist~$\rho>0$ and two functions~$\psi_E$,
$\psi_F:\R^{n-1}\to\R$ such that~$p=(p',p_n)=(p',\psi_E(p'))=(p',\psi_F(p'))$,
$$ E\cap B_\rho(p)=\{ x_n>\psi_E(x')\}\cap B_\rho(p)\qquad{\mbox{ and }}\qquad
F\cap B_\rho(p)=\{ x_n>\psi_F(x')\}\cap B_\rho(p).$$
Also, since~$E\subseteq F$, we have that~$\psi_E\ge\psi_F$ in the vicinity of~$p'$.
Therefore, the function~$\phi:=\psi_E-\psi_F$ has a local minimum at~$p'$, entailing that~$\nabla'\phi(p')=0$
and the Hessian of~$\phi$ at~$p$, that we denote by~$D^2_{x'}\phi(p)$,
is a nonnegative definite matrix. {F}rom this and Theorem~\ref{INCOO} we deduce that
at~$p'$
\begin{equation}\label{TBO-PLA-S0}\begin{split}H_{\partial E}
-H_{\partial F}\,&=\div'\left(
\frac{\nabla' \psi_E}{\sqrt{1+|\nabla' \psi_E|^2}}
\right)-\div'\left(
\frac{\nabla' \psi_F}{\sqrt{1+|\nabla' \psi_F|^2}}
\right)\\&=\frac{\Delta' \psi_E}{\sqrt{1+|\nabla' \psi_E|^2}}-\frac{\Delta' \psi_F}{\sqrt{1+|\nabla' \psi_F|^2}}
-\sum_{i,j=1}^{n-1}\left(\frac{\partial_{ij}\psi_E\partial_i\psi_E\,\partial_j\,\psi_E}{(1+|\nabla' \psi_E|^2)^{3/2}}
-\frac{\partial_{ij}\psi_F\,\partial_i\psi_F\,\partial_j\psi_F}{(1+|\nabla' \psi_F|^2)^{3/2}}\right)\\&=
\frac{1}{({1+|\nabla' \psi_E|^2})^{3/2}}\left(
\big({1+|\nabla' \psi_E|^2}\big)\sum_{i=1}^n(D^2_{x'}\phi(p) \,e_i)\cdot e_i-{
(D^2_{x'}\phi(p) \,\nabla'\psi_E)\cdot \nabla'\psi_E}
\right).\end{split}\end{equation}

We now observe that if~$M\in{\rm Mat}(m\times m)$
is symmetric and nonnegative definite, and~$v\in\R^m$, then
\begin{equation}\label{TBO-PLA-S}
|v|^2\sum_{i=1}^m Me_i\cdot e_i\ge Mv\cdot v.
\end{equation}
To prove this,
we use the Spectral Theorem to find
an orthonormal basis~$\{\eta_1,\dots,\eta_m\}$ of~$\R^m$
consisting of eigenvectors of~$M$. Namely,
$$ M\eta_i=\lambda \eta_i\qquad{\mbox{for all }}i\in\{1,\dots,m\},$$
with~$\lambda_1,\dots,\lambda_m\ge0$.
Thus, we write~$v$ in terms of this new basis as
$$ v=\sum_{i=1}^m (v\cdot \eta_i)\eta_i$$
and we have that
\begin{eqnarray*}
Mv\cdot v=\sum_{i,j=1}^m (v\cdot \eta_i)
(v\cdot \eta_j)M\eta_i\cdot\eta_j=
\sum_{i,j=1}^m \lambda_i (v\cdot \eta_i)
(v\cdot \eta_j)(\eta_i\cdot\eta_j)
=
\sum_{i=1}^m \lambda_i (v\cdot \eta_i)^2.\end{eqnarray*}
Applying this identity with~$e_k$ instead of~$v$,
we also obtain that, for each~$k\in\{1,\dots,m\}$,
\begin{eqnarray*}
Me_k\cdot e_k=
\sum_{i=1}^m \lambda_i (e_k\cdot \eta_i)^2.\end{eqnarray*}
As a result,
\begin{equation}\label{KND st2T2}\begin{split}&
|v|^2\sum_{k=1}^m Me_k\cdot e_k-Mv\cdot v
=
\sum_{i=1}^m \lambda_i \left(\sum_{k=1}^m|v|^2(e_k\cdot \eta_i)^2
-
(v\cdot \eta_i)^2
\right)\\&\qquad\qquad\ge|v|^2\sum_{i=1}^m \lambda_i \left(\sum_{k=1}^m(e_k\cdot \eta_i)^2
-
1
\right).\end{split}
\end{equation}
Hence, since
$$ 1=|\eta_i|^2=\left|
\sum_{k=1}^m(\eta_i\cdot e_k)e_k
\right|^2=\sum_{k=1}^m(\eta_i\cdot e_k)^2,$$
we obtain~\eqref{TBO-PLA-S} as a consequence of~\eqref{KND st2T2}.

{F}rom~\eqref{TBO-PLA-S0}
and~\eqref{TBO-PLA-S}
we conclude that~$H_{\partial E}-H_{\partial F}\ge0$, as desired.
\end{proof}

\begin{corollary}\label{COsdRPOH}
There exists~$p\in\Sigma$ such that~$H(p)>0$.
\end{corollary}

\begin{proof} Since~$E$ is bounded, we have that~$E\subseteq B_R$ for some~$R>0$.
By decreasing~$R$ till the external ball touches~$\partial E$, we can thus suppose that~$E\subseteq B_R$
and there exists~$p\in(\partial E)\cap(\partial B_R)$. Then, by
Corollary~\ref{ILCOMS}, $H(p)\ge\frac{n-1}{R}>0$.\end{proof}

To further analyze the interplay between the Laplace-Beltrami operator and the mean curvature of~$\Sigma$,
one can also introduce the ``tangential Laplacian'', defined, for every~$p\in\Sigma$
and every~$u\in C^2({\mathcal{N}})$, as
\begin{equation}\label{TANLAP}
\Delta_T u(p):=\Delta u(p)-\big( D^2 u(p)\,\nu(p)\big)\cdot\nu(p).\end{equation}
That is, the tangential Laplacian is the ``$(n-1)$-dimensional Laplacian''
on the tangent plane to~$\Sigma$ at~$p$.
We stress that the tangential Laplacian does not coincide with the Laplace-Beltrami operator
and in fact their difference is related to the mean curvature of~$\Sigma$, as pointed out by the following result:

\begin{theorem} \label{225}
Let~$u\in C^2({\mathcal{N}})$. Then, on~$\Sigma$,
$$ \Delta_\Sigma u = \Delta_T u-H\nabla u\cdot \nu.$$
\end{theorem}

\begin{proof} Using Lemma~\ref{GB}, \eqref{IGRA}
and~\eqref{0909}, we see that, on~$\Sigma$,
\begin{eqnarray*}
\Delta_\Sigma u&=&\div_\Sigma(\nabla_\Sigma u)\\
&=&\div_\Sigma\Big(\nabla u-\big(\nabla u\cdot\nu\big)\,\nu\Big)\\
&=&\div\Big(\nabla u-\big(\nabla u\cdot\nu\big)\,\nu\Big)_{\mbox{\scriptsize{ext}}}\\
&=&\div\Big((\nabla u)_{\mbox{\scriptsize{ext}}}\Big)-\nabla
\Big( \big(\nabla u\cdot\nu\big)_{\mbox{\scriptsize{ext}}}\Big)\cdot\nu -\big(\nabla u\cdot\nu\big)\div\nu_{\mbox{\scriptsize{ext}}}.
\end{eqnarray*}
Thus, by~\eqref{CE}, \eqref{0909} and~\eqref{MC},
\begin{eqnarray*}
\Delta_\Sigma u&=&
\div_\Sigma(\nabla u)-0-\big(\nabla u\cdot\nu\big)H.
\end{eqnarray*}
As a consequence, by~\eqref{VB-d1}, and using again~\eqref{CE}, on~$\Sigma$,
\begin{eqnarray*}
\Delta_\Sigma u&=&
\div( \nabla u)-\nabla\big( \nabla u\cdot\nu_{\mbox{\scriptsize{ext}}}\big)\cdot\nu
-\big(\nabla u\cdot\nu\big)H\\&=&
\Delta u-(D^2 u\,\nu)\cdot\nu-\sum_{j=1}^n\nabla (\nu_{\mbox{\scriptsize{ext}}}
\cdot e_j)\cdot\nu\, \partial_j u
-\big(\nabla u\cdot\nu\big)H\\&=&
\Delta_T u-0-\big(\nabla u\cdot\nu\big)H,
\end{eqnarray*}
as desired.\end{proof}

\begin{figure}
  \centering
  \includegraphics[width=.4\linewidth]{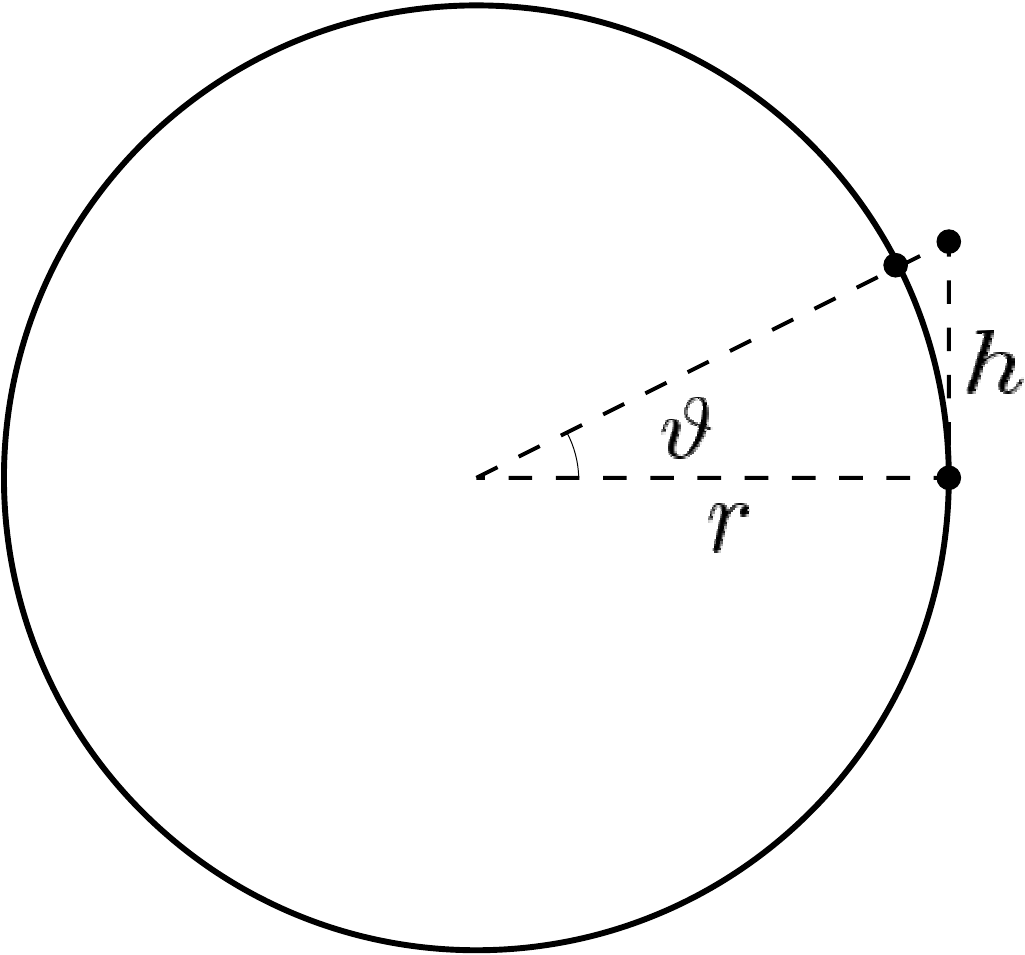}
 \caption{\sl The role of the curvature in second order operators.} \label{DIFI22}
\end{figure}

An intuitive explanation for the role played by the curvature of~$\Sigma$ in the computation
of second order operator naturally arises from Figure~\ref{DIFI22}. Namely, near a given point~$p\in\Sigma$,
the second order behavior of the hypersurface is well-approximated by its osculating sphere,
hence, for simplicity, we presented in Figure~\ref{DIFI22} the case~$n=2$, $\Sigma=\partial B_r$
and, up to a rotation, $p=re_1$. We consider a ``tangential infinitesimal increment''~$he_2$.
Notice from Figure~\ref{DIFI22} that~$h=r\tan\vartheta$, and accordingly~$\cos\vartheta=\frac{r}{\sqrt{r^2+h^2}}=1-\frac{h^2}{2r^2}+o(h^2)$ and~$\sin\vartheta=\frac{h\cos\vartheta}{r}=\frac{h}{r}+o(h^2)$.

Then, the second order incremental quotient of~$u_{\mbox{\scriptsize{ext}}}$ in the tangential direction~$e_2$ is therefore
\begin{eqnarray*}&& \frac{u_{\mbox{\scriptsize{ext}}}(p+he_2)+
u_{\mbox{\scriptsize{ext}}}(p-he_2)-2u_{\mbox{\scriptsize{ext}}}(p)}{h^2}=
\frac{u(r\cos\vartheta,r\sin\vartheta)+
u(r\cos\vartheta,-r\sin\vartheta)-2u(r,0)}{h^2}\\&&\quad=
\frac{u\left(r-\frac{h^2}{2r}+o(h^2),h+o(h^2)\right)+
u\left(r-\frac{h^2}{2r}+o(h^2),-h+o(h^2)\right)-2u(r,0)}{h^2}\\&&\quad=
\frac{1}{h^2}\Bigg(\left(-\partial_1 u(r,0)\frac{h^2}{2r}
+\partial_2 u(r,0)h+\frac12\partial_{22} u(r,0)h^2\right)\\&&\qquad\qquad\qquad\qquad+
\left(-\partial_1 u(r,0)\frac{h^2}{2r}
-\partial_2 u(r,0)h+\frac12\partial_{22} u(r,0)h^2\right)
+o(h^2)\Bigg)\\&&\quad=
-\frac{1}{r}\,\partial_1 u(r,0)+\partial_{22} u(r,0)+o(1),
\end{eqnarray*}
which, as~$h\searrow0$, recalling that the curvature of the circle is~$\frac1r$,
converges to~$-H\partial_\nu u(p)+\partial_{22} u(p)$, and this corresponds to the identity found in
Theorem~\ref{225}.\medskip

A simple consequence of
Theorem~\ref{225} is also the expression of the Laplacian along level sets\index{Laplacian along level sets}
of the function~$u$, according to the following observation:

\begin{corollary}\label{C228}
Let~$u\in C^2({\mathcal{N}})$. Suppose that~$u$ is constant along~$\Sigma$. Then, on~$\Sigma$,
$$ \Delta u= H\,\partial_\nu u+\partial_{\nu\nu} u
.$$
\end{corollary}

\begin{proof} Let~$c\in\R$ be such that~$u=c$ on~$\Sigma$.
Then, $u_{\mbox{\scriptsize{ext}}}(x)=u(\pi_\Sigma(x))=c$ for
all~$x\in{\mathcal{N}}$, whence~$\Delta_\Sigma u=0$.
Consequently,
one combines Theorem~\ref{225} and the definition
of tangential Laplacian in~\eqref{TANLAP} to see that, on~$\Sigma$,
\begin{equation*}
H\nabla u\cdot \nu=
\Delta_\Sigma u +H\nabla u\cdot \nu=
\Delta_T u=\Delta u-( D^2 u\,\nu\big)\cdot\nu
.\qedhere\end{equation*}
\end{proof}

For completeness, we now give another formula to compute the mean curvature
\index{mean curvature} of
the boundary of a smooth set. This formula is perhaps not very handy for explicit calculations, but it has the
conceptual advantage of expressing the mean curvature as a local average of a set against its complement
(thus entailing that a set whose boundary has zero mean curvature has the property that the volume density of
the set is well compensated at each boundary point by the volume density of its complement). The precise result
goes as follows:

\begin{lemma}
For every~$p\in\Sigma=\partial E$,
\begin{equation}\label{IniZibnepbu} H(p)=\lim_{r\searrow0} \frac{1}{c r^{n+1}}\int_{B_r(p)}\big(\chi_{\R^n\setminus E}(x)-\chi_E(x)\big)\,dx,\end{equation}
for a suitable positive constant~$c$ depending only on~$n$.
\end{lemma}

\begin{proof}
Up to a translation and a rotation, we assume that the point~$p$ is the origin and that~$E$ in the vicinity
of the origin is the superlevel set of a function~$\psi$, as in~\eqref{uojw-S29-32jfewnb}, with~$\psi(0)=0$ and~$\nabla\psi(0)=0$.
Consequently, given~$\e>0$ there exists~$r_\e>0$ such that if~$r\in(0,r_\e)$ then
\begin{equation}\label{0okS2890io3rhnfwe-4PA} \begin{split}&U:=B_r\cap\left\{
x_n\ge \frac12 D^2\psi(0)x'\cdot x'+\e|x'|^2\right\}\subseteq E\cap B_r\\
&\qquad\subseteq B_r\cap\left\{
x_n\ge \frac12 D^2\psi(0)x'\cdot x'-\e|x'|^2\right\}=:V,\end{split}\end{equation}
see Figure~\ref{S2890io3rhnfwe-4PA}.

\begin{figure}
  \centering
  \includegraphics[width=.4\linewidth]{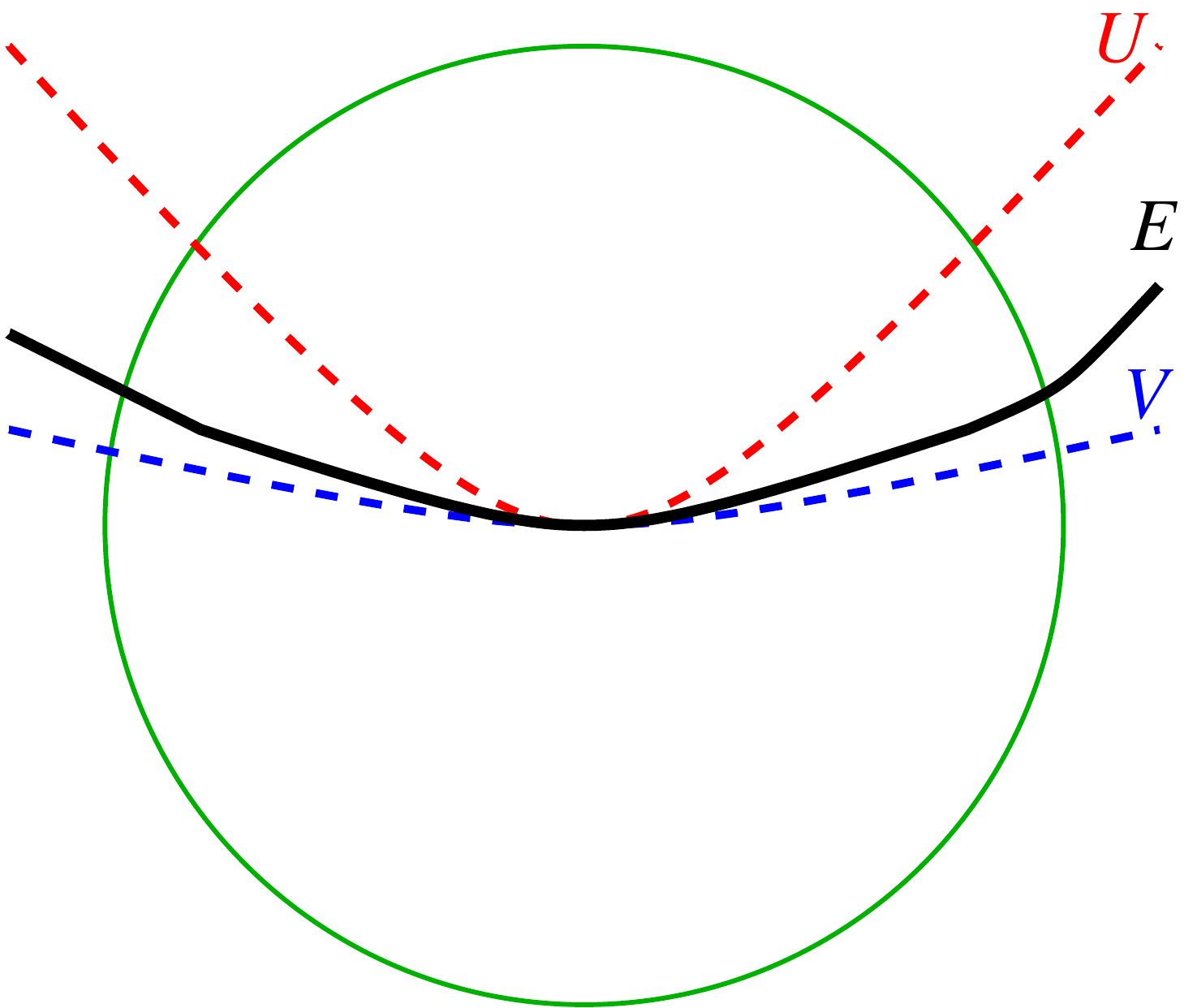}
 \caption{\sl The sets in~\eqref{0okS2890io3rhnfwe-4PA}.}  \label{S2890io3rhnfwe-4PA}
\end{figure}

Accordingly, if~$y\in B_r\setminus V$ then~$\chi_{\R^n\setminus E}(y)-\chi_E(y)=1$
and if~$y\in U$ then~$\chi_{\R^n\setminus E}(y)-\chi_E(y)=-1$. It is therefore convenient to define
$$ S:=B_r\cap\left\{
x_n- \frac12 D^2\psi(0)x'\cdot x'\in\big(-\e|x'|^2,\,\e|x'|^2\big)\right\}$$
and note that
$$ B_r= U\cup S\cup( B_r\setminus V).$$
These observations give that
\begin{equation}\label{APCOSMFe}
\int_{B_r}\big(\chi_{\R^n\setminus E}(x)-\chi_E(x)\big)\,dx
=|B_r\setminus V|-|U|+\int_S\big(\chi_{\R^n\setminus E}(x)-\chi_E(x)\big)\,dx.
\end{equation}
Moreover,
\begin{eqnarray*}
|S|\le2\e\int_{|x'|< r}|x'|^2\,dx' =O(\e r^{n+1}),
\end{eqnarray*}
whence we infer from~\eqref{APCOSMFe} that
\begin{equation}\label{APCOSMFe-32}
\int_{B_r}\big(\chi_{\R^n\setminus E}(x)-\chi_E(x)\big)\,dx
=|B_r\setminus V|-|U|+O(\e r^{n+1}).
\end{equation}

Similarly,
\begin{eqnarray*}
|U|=\left| B_r\cap\left\{x_n\ge \frac12 D^2\psi(0)x'\cdot x'\right\}\right|+O(\e r^{n+1})
\end{eqnarray*}
and
\begin{eqnarray*}
|B_r\setminus V|&=&\left| B_r\cap\left\{x_n< \frac12 D^2\psi(0)x'\cdot x'\right\}\right|+O(\e r^{n+1})\\&=&
\left| B_r\cap\left\{x_n> -\frac12 D^2\psi(0)x'\cdot x'\right\}\right|+O(\e r^{n+1}),
\end{eqnarray*}
where the change of variable~$x_n\mapsto -x_n$ has been used in the last equality.

{F}rom these observations, using the short notation~$M(x'):=\frac12 D^2\psi(0)x'\cdot x'$ and~\eqref{APCOSMFe-32} we find that
\begin{equation}\label{APCOSMFe-32090sd}
\begin{split}&
\int_{B_r}\big(\chi_{\R^n\setminus E}(x)-\chi_E(x)\big)\,dx\\
&\qquad\quad=\,
\left| B_r\cap\big\{x_n> -M(x')\big\}\right|-\left| B_r\cap\big\{x_n\ge M( x')\big\}\right|+O(\e r^{n+1}).
\end{split}\end{equation}
Now, letting~$C:=\frac12 |D^2\psi(0)|$,
we have that
\begin{equation}\label{GAVBuTGBStnao0os-1}\begin{split}
\left| B_r\cap\big\{x_n\ge M( x')\big\}\right|&=\left| B_r\cap\big\{M(x')\le x_n\le C|x'|^2\big\}\right|+
\left| B_r\cap\big\{x_n\ge C|x'|^2\big\}\right|\\&=
\left| B_r\cap\big\{M(x')\le x_n\le C|x'|^2\big\}\right|+A_r\\
&=E_r-
F_r
+A_r,
\end{split}\end{equation}
where
\begin{eqnarray*}&&A_r:=\left| B_r\cap\big\{x_n\ge C|x'|^2\big\}\right|,\\&&
E_r:=\left| \big\{|x'|<r\big\}\cap\big\{M(x')\le x_n\le C|x'|^2\big\}\right|\\{\mbox{and}}\qquad&&
F_r:=\left| \big\{|x'|<r<|x|\big\}\cap\big\{M(x')\le x_n\le C|x'|^2\big\}\right|.
\end{eqnarray*}
Additionally, if~$x\in F_r$ then
$$r^2\le|x|^2= |x'|^2+x_n^2\le|x'|^2+(C|x'|^2)^2=|x'|^2+C^2|x'|^4\le (1+C^2r^2)|x'|^2
\le (1+\e)|x'|^2,$$
as long as~$r$ is small enough. This entails that
$$ F_r\subseteq \left\{ |x'|\in\left[ \frac{r}{\sqrt{1+\e}},r\right]\right\}
\cap\big\{| x_n|\le C|x'|^2\big\},$$
whence~$|F_r|=O(\e r^{n+1})$.

{F}rom this observation and~\eqref{GAVBuTGBStnao0os-1} we arrive at
\begin{eqnarray*}
\left| B_r\cap\big\{x_n\ge M( x')\big\}\right|&=&E_r+A_r+O(\e r^{n+1})\\
&=& \left| \big\{|x'|<r\big\}\cap\big\{M(x')\le x_n\le C|x'|^2\big\}\right|+A_r+O(\e r^{n+1}).\end{eqnarray*}
Thus, replacing~$M$ with~$-M$ and noticing that~$A_r$ remains unchanged under this modification,
\begin{eqnarray*}
\left| B_r\cap\big\{x_n> -M( x')\big\}\right|&=&
\left| \big\{|x'|<r\big\}\cap\big\{-M(x')\le x_n\le C|x'|^2\big\}\right|+A_r+O(\e r^{n+1}).\end{eqnarray*}

Using these remarks and~\eqref{APCOSMFe-32090sd} we thereby conclude that
\begin{eqnarray*}
&&
\int_{B_r}\big(\chi_{\R^n\setminus E}(x)-\chi_E(x)\big)\,dx\\
&&\qquad\quad=\,
\left(
\left| \big\{|x'|<r\big\}\cap\big\{-M(x')\le x_n\le C|x'|^2\big\}\right|+A_r\right)\\&&\qquad\qquad\quad
-\left( \left| \big\{|x'|<r\big\}\cap\big\{M(x')\le x_n\le C|x'|^2\big\}\right|+A_r\right)+O(\e r^{n+1})\\
&&\qquad\quad=\,
\left| \big\{|x'|<r\big\}\cap\big\{-M(x')\le x_n\le C|x'|^2\big\}\right|\\&&\qquad\qquad\quad
- \left| \big\{|x'|<r\big\}\cap\big\{M(x')\le x_n\le C|x'|^2\big\}\right|+O(\e r^{n+1})\\
&&\qquad\quad=\,\int_{\{|x'|<r\}} \left(\int_{-M(x')}^{M(x')} dx_n\right)\,dx'
+O(\e r^{n+1})\\&&\qquad\quad=\,
D^2\psi(0)\int_{\{|x'|<r\}} x'\cdot x'\,dx'
+O(\e r^{n+1})\\&&\qquad\quad=\,
\sum_{i=1}^{n-1} \partial^2_i \psi(0)\int_{\{|x'|<r\}} x_i^2\,dx'
+O(\e r^{n+1})\\&&\qquad\quad=\,
\frac1{n-1}\sum_{i=1}^{n-1} \partial^2_i \psi(0)\int_{\{|x'|<r\}} |x'|^2\,dx'
+O(\e r^{n+1})\\&&\qquad\quad=\,
cr^{n+1} \Delta\psi(0)
+O(\e r^{n+1}),
\end{eqnarray*}
for a suitable~$c>0$, from which~\eqref{IniZibnepbu} plainly follows in view of~\eqref{LAMEANCURVA} .
\end{proof}

\section{The Laplace-Beltrami operator in local coordinates}

\begin{figure}
  \centering
  \includegraphics[width=.4\linewidth]{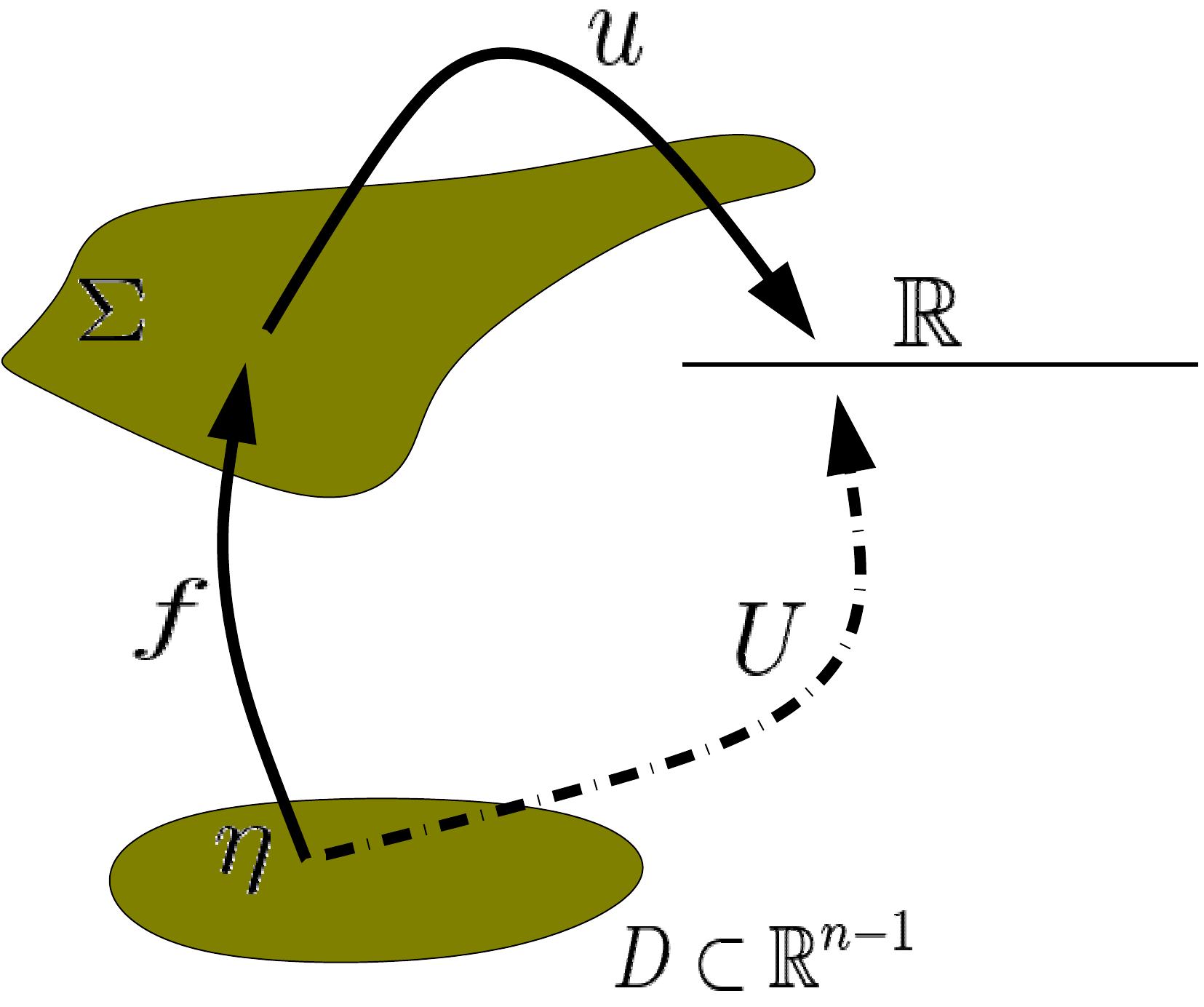}
 \caption{\sl Local coordinates for~$\Sigma$.}  \label{SPA}
\end{figure}

Regarding the setting in Section~\ref{BELS},
it is sometimes useful to express geometric
differential operators (such as tangential gradients and divergences,
as well as the 
Laplace-Beltrami operator) in ``local coordinates''
with respect to the manifold~$\Sigma$.
That is, one supposes that an element of~$\Sigma$
is locally parameterized by a diffeomorphism~$f:D\to\Sigma$,
for some domain~$D\subset\R^{n-1}$, see Figure~\ref{SPA}.
One considers the metrics~${\mbox{\Large{\calligra{g}}}}$,
which can be identified with a symmetric $(n-1)\times(n-1)$ matrix with elements
\begin{equation}\label{METRIX}
g_{ij}:=\frac{\partial f}{\partial\eta_i}\cdot\frac{\partial f}{\partial\eta_j}.\end{equation}
We observe that the matrix~${\mbox{\Large{\calligra{g}}}}\,$
is the product between the transpose of~$D_\eta f$ and~$D_\eta f$ itself, being~$(D_\eta f)_{kj}=
\frac{\partial f}{\partial\eta_j}\cdot e_k$, for~$j\in\{1,\dots,n-1\}$ and~$k\in\{1,\dots,n\}$. Therefore,
by the Binet-Cauchy Formula (see e.g.~\cite[Theorem~3.7]{MR3409135}), the square of
the determinant of~${\mbox{\Large{\calligra{g}}}}\,$
is the sum of the squares
of the determinants of each~$ (n-1)\times(n-1)$-submatrix of~$D_\eta f$, namely
\begin{equation}\label{DETTE} \det {\mbox{\Large{\calligra{g}}}}=\sqrt{\sum_{k=1}^n\left[\;
\det\left( \begin{matrix}\frac{\partial f}{\partial\eta}\cdot e_1
\\
\vdots\\
\frac{\partial f}{\partial\eta}\cdot e_{k-1}\\
\frac{\partial f}{\partial\eta}\cdot e_{k+1}\\ \vdots\\
\frac{\partial f}{\partial\eta}\cdot e_n\end{matrix}
\right)\;\right]^2
}.\end{equation}
For instance, if~$\Sigma$ is locally the graph (say, in the $n$th direction)
of a function~$\psi$, one can take~$f(\eta):=(\eta,\psi(\eta))$,
notice that in this case~$\frac{\partial f}{\partial\eta_i}=\left(e_i',\,\frac{\partial \psi}{\partial\eta_i}\right)$,
which are linearly independent,
being~$\{e'_1,\dots,e'_{n-1}\}$ the Euclidean basis of~$\R^{n-1}$,
obtain that
$$g_{ij}=\delta_{ij}+\frac{\partial\psi}{\partial\eta_i}\,\frac{\partial \psi}{\partial\eta_j}$$
and deduce from~\eqref{DETTE} that
\begin{equation*} \det {\mbox{\Large{\calligra{g}}}}=\sqrt{1+|\partial_\eta\psi|^2
}.\end{equation*}
We also remark that, for every~$v=(v_1,\dots,v_{n-1})\in\R^{n-1}$,
$$ ({\mbox{\Large{\calligra{g}}}}\, v)\cdot v=\sum_{i,j=1}^{n-1}\left(
\delta_{ij}+\frac{\partial\psi}{\partial\eta_i}\,\frac{\partial \psi}{\partial\eta_j}\right)v_i\,v_j=
|v|^2+\left(\frac{\partial\psi}{\partial\eta}\cdot v\right)^2\ge|v|^2,
$$
whence~${\mbox{\Large{\calligra{g}}}}\,$ is positive definite.

In particular, we have that~${\mbox{\Large{\calligra{g}}}}\,$ is invertible
and we denote the corresponding inverse matrix by~$g^{ij}$.
We also set~$g:=\det{\mbox{\Large{\calligra{g}}}}\,$.

In this setting, given~$u:\Sigma\to\R$, we can consider~$
U:D\to\R$ defined by~$U(\eta):=u(f(\eta))$.
Moreover,
it is convenient to identify vectors~$\tau\in\R^n$ which are tangent to~$\Sigma$
at a given point with suitable vectors in~$\R^{n-1}$. This identification
should also reconstruct the Euclidean
scalar products of tangent vectors in~$\R^n$
from a suitable product between the associated vectors in~$\R^{n-1}$. To this end,
given~$p\in\Sigma$, we can assume that~$0\in D$ and~$p=f(0)$ and,
for all~$i\in\{1,\dots,n\}$, we consider the tangent vectors
$$ \tau_i:=\frac{d}{dt} f(te_i)\Big|_{t=0}
= \frac{\partial f}{\partial\eta_i}(0)
.$$
Given
a tangent vector~$V\in\R^n$ to~$\Sigma$ at~$p$
we write it in the form
\begin{equation}\label{0oLSPpSLkmdc} V=\sum_{i=1}^{n-1} V^i\,\tau_i
.\end{equation}
Thus, we can associate to the tangent vector~$V$
its $(n-1)$-dimensional representation
\begin{equation}\label{THETAVE}\Theta(V)=(\Theta^1(V),\dots,\Theta^{n-1}(V))=(V^1,\dots,V^{n-1})\in\R^{n-1}.\end{equation}
In this way, considering the scalar product induced by
the metrics between~$a=(a^1,\dots,a^{n-1})$, $b=(b^1,\dots,b^{n-1})\in\R^{n-1}$ defined as
$$ g(a,b):=\sum_{i,j=1}^{n-1} g_{ij} \,a^i\, b^j,$$
it follows that, for every tangent vector~$V$, $W\in\R^n$,
\begin{equation}\label{MANMGBmhsl04599} V\cdot W=
\sum_{i,j=1}^n V^i\,W^j\,{\tau_i\cdot\tau_j}
=\sum_{i,j=1}^n \Theta^i(V)\,\Theta^j(W)\,g_{ij}=
g\big(\Theta(V),\Theta(W)\big).
\end{equation}
We also recall a handy way to use the metric to reconstruct the coordinates of a tangent
vector with respect to the (not necessarily
orthonomal) basis of the tangent space given by~$\{\tau_1,\dots,\tau_{n-1}\}$. 

\begin{lemma}
Let~$V$ be a tangent vector field on~$\Sigma$ as in~\eqref{0oLSPpSLkmdc}, 
and recall the notation in~\eqref{THETAVE}.
Then, for every~$j\in\{1,\dots,n-1\}$,
\begin{equation}\label{LAMNjjma834RSIi2} \Theta^j(V) = \sum_{i=1}^{n-1} g^{ij}\,V\cdot\frac{\partial f}{\partial\eta_i}.\end{equation}
Moreover,
\begin{equation}\label{LAMNjjma834RSIi} V=
\sum_{i,j=1}^{n-1}
g^{ij}\,V\cdot\frac{\partial f}{\partial\eta_i}\,\frac{\partial f}{\partial\eta_j}
.\end{equation}\end{lemma}

\begin{proof} It follows from~\eqref{0oLSPpSLkmdc}, that, for every~$j\in\{1,\dots,n-1\}$,
$$ V\cdot\frac{\partial f}{\partial\eta_j}=
\sum_{i=1}^{n-1} V^i\,\tau_i\cdot
\frac{\partial f}{\partial\eta_j}
=\sum_{i=1}^{n-1} V^i\frac{\partial f}{\partial\eta_i}\cdot
\frac{\partial f}{\partial\eta_j}=
\sum_{i=1}^{n-1} g_{ij}V^i.$$
That is, if we consider the vectors~$W:=\left(V\cdot\frac{\partial f}{\partial\eta_1},\dots,
V\cdot\frac{\partial f}{\partial\eta_{n-1}}
\right)$ and~$\Theta:=\Theta(V)=(V^1,\dots,V^{n-1})$, in matrix notation we have that~$W={{\mbox{\Large{\calligra{g}}}}}\,\Theta $
and thus~$\Theta ={{\mbox{\Large{\calligra{g}}}}}^{-1}W$, which gives, for each~$j\in\{1,\dots,n-1\}$,
$$\Theta^j=\sum_{i=1}^{n-1}g^{ij} W_i,$$
from which~\eqref{LAMNjjma834RSIi2} plainly follows. 

In view of~\eqref{LAMNjjma834RSIi2}, we see that~\eqref{0oLSPpSLkmdc} reduces to
\begin{equation*} V=\sum_{j=1}^{n-1} V^j\frac{\partial f}{\partial\eta_j}=
\sum_{i,j=1}^{n-1}
g^{ij}\,V\cdot\frac{\partial f}{\partial\eta_i}\,\frac{\partial f}{\partial\eta_j}
,\end{equation*}
which is~\eqref{LAMNjjma834RSIi}.\end{proof}

Formula~\eqref{LAMNjjma834RSIi} can be straightforwardly extended to all (not necessarily tangent) vector fields:

\begin{corollary}\label{NAMSblamnd872BAP01}
If~$V$ is a vector field on~$\Sigma$,
\begin{equation*} V=(V\cdot\nu)\nu+
\sum_{i,j=1}^{n-1}
g^{ij}\,V\cdot\frac{\partial f}{\partial\eta_i}\,\frac{\partial f}{\partial\eta_j}
.\end{equation*}
\end{corollary}

\begin{proof}
{F}rom~\eqref{LAMNjjma834RSIi},
considering the tangent vector field~$\widetilde{V}:=V-(V\cdot\nu)\nu$, we see that
\begin{equation*}\begin{split}& V-(V\cdot\nu)\nu=\widetilde{V}=
\sum_{i,j=1}^{n-1}
g^{ij}\,\widetilde{V}\cdot\frac{\partial f}{\partial\eta_i}\,\frac{\partial f}{\partial\eta_j}\\&\qquad\qquad=
\sum_{i,j=1}^{n-1}
g^{ij}\,\Big(V-(V\cdot\nu)\nu\Big)\cdot\frac{\partial f}{\partial\eta_i}\,\frac{\partial f}{\partial\eta_j}=\sum_{i,j=1}^{n-1}
g^{ij}\,{V}\cdot\frac{\partial f}{\partial\eta_i}\,\frac{\partial f}{\partial\eta_j}
.\qedhere\end{split}\end{equation*}\end{proof}

The above setting allows us to write the tangential gradient in local coordinates, according
to the following result:

\begin{lemma} \label{0905t0987654828232-548c8x90c0c0}
Let~$u\in C^1({\mathcal{N}})$.
For every~$i\in\{1,\dots,n-1\}$, on~$\Sigma$,
$$ \Theta^i(\nabla_\Sigma u) = \sum_{j=1}^{n-1} g^{ij}\frac{\partial U}{\partial\eta_j}.$$
\end{lemma}

\begin{proof} We stress that~$\nabla_\Sigma u$ is a tangent vector since,
by~\eqref{IGRA},
$$ \nabla_\Sigma u\cdot\nu=
\Big(\nabla u-\big(\nabla u\cdot\nu\big)\,\nu\Big)\cdot\nu=\big(\nabla u\cdot\nu\big)-
\big(\nabla u\cdot\nu\big)\nu\cdot\nu=0.$$
Also, by the Chain Rule, for all~$i\in\{1,\dots,n-1\}$,
$$\frac{\partial U}{\partial\eta_i}(\eta)=\frac{\partial }{\partial\eta_i}(u(f(\eta))=\nabla u(f(\eta))\cdot\frac{\partial U}{\partial\eta_i}(f(\eta))$$
and thus, in light of~\eqref{LAMNjjma834RSIi2},
\begin{equation*}
\begin{split}& \Theta^j(\nabla_\Sigma u) = \sum_{i=1}^{n-1} g^{ij}\,\nabla_\Sigma u\cdot\frac{\partial f}{\partial\eta_i}=
\sum_{i=1}^{n-1} g^{ij}\,\Big(\nabla u-\big(\nabla u\cdot\nu\big)\,\nu\Big)\cdot\frac{\partial f}{\partial\eta_i}\\&\qquad\quad\qquad=
\sum_{i=1}^{n-1} g^{ij}\, \nabla u\cdot\frac{\partial f}{\partial\eta_i}=\sum_{i=1}^{n-1} g^{ij}\frac{\partial U}{\partial\eta_i}.
\qedhere\end{split}\end{equation*}
\end{proof}

Correspondingly to Lemma~\ref{0905t0987654828232-548c8x90c0c0}, one can write the tangential divergence
in local coordinates as follows: 

\begin{lemma}\label{SDMDOF-IR}
Let~$F\in C^1(\Sigma,\R^n)$ be a tangent vector field. Then, on~$\Sigma$,
$$ \div_\Sigma F=
\sum_{i=1}^{n-1}\frac1{\sqrt{g}}\frac{\partial}{\partial\eta_i}\Big(
\sqrt{g}\,\Theta^i(F)
\Big).$$
\end{lemma}

\begin{proof} We consider a smooth function~$\varphi$ supported in a local chart and
let~$\Phi(\eta):=\varphi(f(\eta))$.
We recall (see~\cite[125-126]{MR3409135})
that the surface element of~$\Sigma$ in local coordinates
can be written as~$\sqrt{g}\,d\eta$, therefore, by Theorem~\ref{GBNAmn},
and using that~$F$ is tangential,
\begin{eqnarray*}&& \div_\Sigma F\,\Phi\,\sqrt{g}\,d\eta=\div_\Sigma F\,\varphi\,d{\mathcal{H}}^{n-1}=
H F\cdot \nu\,\varphi\,d{\mathcal{H}}^{n-1}-F\cdot
\nabla_\Sigma\varphi\,d{\mathcal{H}}^{n-1}=
-F\cdot
\nabla_\Sigma\varphi\,d{\mathcal{H}}^{n-1}
.\end{eqnarray*}
This and~\eqref{MANMGBmhsl04599} lead to
$$ \div_\Sigma F\,\Phi\,\sqrt{g}\,d\eta=
-\sum_{i,j=1}^n g_{ij}\, \Theta^i(F)\,\Theta^j(\nabla_\Sigma\varphi)\,d{\mathcal{H}}^{n-1}.$$
Since, in light of Lemma~\ref{0905t0987654828232-548c8x90c0c0}, 
$$ \Theta^j(\nabla_\Sigma \varphi) = \sum_{k=1}^{n-1} g^{jk}\frac{\partial \Phi}{\partial\eta_k},$$
we thereby obtain that
$$ \div_\Sigma F\,\Phi\,\sqrt{g}\,d\eta=
-\sum_{i,j,k=1}^n g_{ij}\, g^{jk}\,\Theta^i(F)\,
\frac{\partial \Phi}{\partial\eta_k}\,\sqrt{g}\,d\eta
=-\sum_{i=1}^n \Theta^i(F)\,
\frac{\partial \Phi}{\partial\eta_i}\,\sqrt{g}\,d\eta.$$
The latter term can be integrated by parts, giving
$$ \div_\Sigma F\,\Phi\,\sqrt{g}\,d\eta=
\sum_{i=1}^n 
\frac{\partial }{\partial\eta_i}\Big(
\Theta^i(F)\,
\,\sqrt{g}\Big)\Phi\,d\eta.$$
Consequently, since~$\varphi$ (and hence~$\Phi$) is an arbitrary 
test function, we conclude that
\begin{equation*}
\div_\Sigma F\, \sqrt{g}=
\sum_{i=1}^n 
\frac{\partial }{\partial\eta_i}\Big(
\Theta^i(F)\,
\,\sqrt{g}\Big).
\qedhere\end{equation*}
\end{proof}

\begin{corollary}\label{GVBSKD:VAYSD:SI}
Let~$F\in C^1(\Sigma,\R^n)$. Then, on~$\Sigma$,
$$ \div_\Sigma F=
\sum_{i=1}^{n-1}\frac1{\sqrt{g}}\frac{\partial}{\partial\eta_i}\Big(
\sqrt{g}\,\Theta^i\big(\widetilde F\big)
\Big)+H(F\cdot\nu),$$
where~$\widetilde{F}:=F-(F\cdot\nu)\nu$.\end{corollary}

\begin{proof}
By~\eqref{VB-d1}, we know that the tangential divergence is
linear with respect to the vector fields. 
Therefore, 
\begin{equation}\label{0OKMS8138OMS}
\div_\Sigma F=
\div_\Sigma \widetilde{F}+\div_\Sigma \big((F\cdot\nu)\nu\big)
.\end{equation}
We also remark that if~$\varphi\in C^1(\Sigma)$, on~$\Sigma$ we have that
$$ \div_\Sigma(\varphi\nu)=\div((\varphi\nu)_{\mbox{\scriptsize{ext}}})
-\nabla\big( (\varphi\nu)_{\mbox{\scriptsize{ext}}}\cdot\nu_{\mbox{\scriptsize{ext}}}\big)\cdot\nu_{\mbox{\scriptsize{ext}}}
=\nabla\varphi_{\mbox{\scriptsize{ext}}}\cdot\nu_{\mbox{\scriptsize{ext}}}+
\varphi\div\nu_{\mbox{\scriptsize{ext}}}
-\nabla\varphi_{\mbox{\scriptsize{ext}}}\cdot\nu_{\mbox{\scriptsize{ext}}}
=H\varphi,$$
thanks to~\eqref{CE}
and~\eqref{MC}. This and~\eqref{0OKMS8138OMS} entail that
\begin{equation*}
\div_\Sigma F=
\div_\Sigma \widetilde{F}+H(F\cdot\nu)
.\end{equation*}
Since~$\widetilde{F}$ is a tangential vector field, we can employ Lemma~\ref{SDMDOF-IR}
and obtain the desired result.\end{proof}

As a variant of
Corollary~\ref{GVBSKD:VAYSD:SI}, we also have the
following useful expression in coordinates
of the tangential divergence of a vector field:

\begin{corollary}\label{LSDcmevandtimsflVVSND}
Let~$F\in C^1({\mathcal{N}},\R^n)$. Then, on~$\Sigma$,
$$ \div_\Sigma F=
\sum_{i,j=1}^{n-1}g^{ij}\left(DF\,\frac{\partial f}{\partial\eta_i}\right)\cdot
\frac{\partial f}{\partial\eta_j}.
$$\end{corollary}

\begin{proof} 
For every~$k\in\{1,\dots,n\}$,
by Corollary~\ref{NAMSblamnd872BAP01}, applied here to~$V:=\frac{\partial F}{\partial x_k}$, we see that, 
\begin{equation*} \frac{\partial F}{\partial x_k}=
\left(\frac{\partial F}{\partial x_k}\cdot\nu\right)\nu+
\sum_{i,j=1}^{n-1}
g^{ij}\,\frac{\partial F}{\partial x_k}\cdot\frac{\partial f}{\partial\eta_i}\,\frac{\partial f}{\partial\eta_j}
.\end{equation*}
As a result,
\begin{equation*} \div F=\sum_{k=1}^n
\frac{\partial F}{\partial x_k}\cdot e_k=\sum_{k=1}^n
\left(\frac{\partial F}{\partial x_k}\cdot\nu\right)\nu\cdot e_k+
\sum_{{1\le i, j\le n-1}\atop{1\le k\le n}}
g^{ij}\,\frac{\partial F}{\partial x_k}\cdot\frac{\partial f}{\partial\eta_i}\,\frac{\partial f}{\partial\eta_j}\cdot e_k.
\end{equation*}
{F}rom this, we deduce that
\begin{eqnarray*}
&&\sum_{i,j=1}^{n-1}g^{ij}\left(DF\,\frac{\partial f}{\partial\eta_j}\right)\cdot
\frac{\partial f}{\partial\eta_i}
=
\sum_{{1\le i, j\le n-1}\atop{1\le k\le n}}g^{ij} \frac{\partial F}{\partial x_k}\cdot
\frac{\partial{f}}{\partial\eta_i}\;\frac{\partial f}{\partial\eta_j}\cdot e_k\\&&\qquad=\div F-\sum_{k=1}^n
\left(\frac{\partial F}{\partial x_k}\cdot\nu\right)\nu\cdot e_k=
\div F-(DF\,\nu)\cdot\nu.
\end{eqnarray*}
That is, recalling the definition of tangential divergence in~\eqref{VB-d1}, and exploiting also~\eqref{CE},
\begin{equation*}\begin{split}
&\sum_{i,j=1}^{n-1}g^{ij}\left(DF\,\frac{\partial f}{\partial\eta_j}\right)\cdot
\frac{\partial f}{\partial\eta_i}-\div_\Sigma F=\sum_{i,j=1}^{n-1}g^{ij}\left(DF\,\frac{\partial f}{\partial\eta_j}\right)\cdot
\frac{\partial f}{\partial\eta_i}-\div F+\nabla(F\cdot\nu_{\mbox{\scriptsize{ext}}})\cdot\nu_{\mbox{\scriptsize{ext}}}
\\&\qquad
=\nabla(F\cdot\nu_{\mbox{\scriptsize{ext}}})\cdot\nu_{\mbox{\scriptsize{ext}}}-(DF\,\nu)\cdot\nu=\sum_{k,m=1}^n (F\cdot e_m)\,\left(\frac{\partial\nu_{\mbox{\scriptsize{ext}}}}{\partial x_k}\cdot e_m\right)
\,(\nu_{\mbox{\scriptsize{ext}}}\cdot e_k)\\&\qquad
=\sum_{m=1}^n (F\cdot e_m)\,\nabla(\nu_{\mbox{\scriptsize{ext}}}\cdot e_m)\cdot\nu_{\mbox{\scriptsize{ext}}}=0
.\qedhere
\end{split}\end{equation*}
\end{proof}

As a direct consequence of
Lemmata~\ref{GB}, \ref{0905t0987654828232-548c8x90c0c0} and~\ref{SDMDOF-IR},
we obtain that
the Laplace-Beltrami operator can be written in terms of these local coordinates, according to the following result:

\begin{corollary}\label{Jns89ijgm5jjgfjgjg}
Let~$u\in C^2({\mathcal{N}})$. Then,
\begin{equation}\label{Jns89ijgm5jjgfjgjg2}
\Delta_\Sigma u =\sum_{i,j=1}^{n-1}\frac1{\sqrt{g}}\frac{\partial}{\partial\eta_i}\left(
\sqrt{g}\,g^{ij}
\frac{\partial U}{\partial\eta_j}
\right).
\end{equation}
\end{corollary}

As a technical remark, we observe that, for the purpose of defining the Laplace-Beltrami operator, the regularity of~$\Sigma$
can be reduced from~$C^3$ to~$C^2$. \label{DAC3aC22}
Indeed, on page~\pageref{DAC3aC2} we assumed~$\Sigma$ of class~$C^3$, which was used
to deduce that~$\nu$ was of class~$C^2$, and accordingly, by~\eqref{PROIE}, that~$\pi_\Sigma$
was of class~$C^2$. {F}rom this, we deduced that~$u_{\mbox{\scriptsize{ext}}}$ was of class~$C^2$
if so was~$u:\Sigma\to\R$, thanks to~\eqref{DEEST}, and the $C^2$ regularity of~$u_{\mbox{\scriptsize{ext}}}$
was utilized to give the definition of the
Laplace-Beltrami operator in~\eqref{BEL}.
On the other hand, thanks to Corollary~\ref{Jns89ijgm5jjgfjgjg},
we can now point out that a $C^2$ regularity assumption on~$\Sigma$ would suffice
for formula~\eqref{Jns89ijgm5jjgfjgjg2}. Hence, taking~\eqref{Jns89ijgm5jjgfjgjg2} instead of~\eqref{BEL}
as
definition of the
Laplace-Beltrami operator would allow us to work with~$\Sigma$ of class~$C^2$
(alternatively, one can define the Laplace-Beltrami operator as in~\eqref{BEL} for~$\Sigma$ of class~$C^3$
and then extend it by approximation when~$\Sigma$ is of class~$C^2$, using~\eqref{Jns89ijgm5jjgfjgjg2}
to pass to the limit).

For additional information on the Laplace-Beltrami operator see e.g.~\cite{MR1481970, MR1462892, MR3726907}
and the references therein.

\section{The Laplacian in spherical coordinates}

Let~$u:\R^n\setminus\{0\}\to\R$.
For every~$r>0$ and~$\vartheta\in\partial B_1$, we define
\begin{equation*}
u_0(r,\vartheta):= u(r\vartheta).
\end{equation*}
In this setting, one can compute the Laplace operator in spherical coordinates according\index{Laplacian in spherical coordinates}
to the following formula:

\begin{theorem} \label{SPH}
Let~$u\in C^2(\R^n\setminus\{0\})$. Then,
$$ \Delta u(x)=\partial_{rr}u_0(r,\vartheta)+\frac{n-1}{r}\,\partial_{r}u_0(r,\vartheta)+\Delta_{\partial B_r}u(r\vartheta),$$
for every~$x\in\R^n\setminus\{0\}$, where~$r:=|x|$ and~$\vartheta:=\frac{x}{|x|}$.
\end{theorem}

\begin{proof} We take~$E:=B_r$ and~$\Sigma:=\partial B_r$. In this way, we have that~$\nu(x)=\frac{x}{|x|}$
for every~$x\in\partial B_r$.
Then, by~\eqref{DEEST}, $\nu_{\mbox{\scriptsize{ext}}}(x)= \frac{x}{|x|}$.
Hence, we deduce from~\eqref{0909} and~\eqref{MC}
that, for every~$x\in\partial B_r$,
\begin{equation} \label{8yujbggd}
H(x)=\div_{\partial B_r} \nu(x)=\div \nu_{\mbox{\scriptsize{ext}}}(x)=\div\frac{x}{|x|}=
\frac{n-1}{|x|}=\frac{n-1}r.\end{equation}
In addition, if~$x\in\partial B_r$,
\begin{eqnarray*}
\nabla u(x)\cdot \nu(x)=\lim_{h\to0}\frac{u(x+h\nu(x))-u(x)}{h}=\lim_{h\to0}\frac{u_0(r+h,\vartheta)-u_0(r,\vartheta)}{h}=\partial_ru_0(r,\vartheta)
\end{eqnarray*}
and
\begin{eqnarray*}&& (D^2 u(x)\,\nu(x))\cdot\nu(x)=\lim_{h\to0}\frac{u(x+h\nu(x))+u(x-h\nu(x))-2u(x)}{h^2}\\&&\qquad=
\lim_{h\to0}\frac{u_0(r+h,\vartheta)+u_0(r-h,\vartheta)-2u_0(r,\vartheta)}{h^2}=\partial_{rr}u_0(r,\vartheta).
\end{eqnarray*}
Using these identities and~\eqref{8yujbggd} in combination with
Theorem~\ref{225}, we conclude that, at~$x\in\partial B_r$,
\begin{equation*} \Delta_{\partial B_r} u = \Delta_T u-H\nabla u\cdot \nu=
\Delta u-(D^2 u\,\nu)\cdot\nu-\frac{n-1}r \,\partial_r u_0=
\Delta u-\partial_{rr}u_0-\frac{n-1}r \,\partial_r u_0
.\qedhere\end{equation*}
\end{proof}

Given~$r>0$, an equivalent formulation
of Theorem~\ref{SPH} can be given by replacing
the Laplace-Beltrami of~$u$ along~$\partial B_r$
with the Laplace-Beltrami
of the map~$\partial B_1\ni\vartheta\mapsto u_0(r,\vartheta)$
along~$\partial B_1$, as follows:

\begin{theorem} \label{SPHECOO}
Let~$u\in C^2(\R^n\setminus\{0\})$. Then,
$$ \Delta u(x)=\partial_{rr}u_0(r,\vartheta)+
\frac{n-1}{r}\,\partial_{r}u_0(r,\vartheta)+\frac1{r^2}\,
\Delta_{\partial B_1}u_0(r,\vartheta),$$
for every~$x\in\R^n\setminus\{0\}$, where~$r:=|x|$ and~$\vartheta:=\frac{x}{|x|}$.
\end{theorem}

\begin{proof} In light of Theorem~\ref{SPH}, it suffices to show that,
given~$r>0$,
\begin{equation}\label{TRGSATTHNBELBRTYUTT9}
\Delta_{\partial B_r}u(r\vartheta)=\frac1{r^2}\Delta_{\partial B_1}u_0(r,\vartheta)
\qquad{\mbox{for all }}\vartheta\in\partial B_1.
\end{equation}
To this end, we perform a careful scaling argument.
Given~$r>0$, for every~$\vartheta\in\partial B_1$
we let~$v^{(r)}(\vartheta):=u_0(r,\vartheta)$.
Thus, the 
normal extension of~$v^{(r)}$ outside~$\partial B_1$,
as introduced in~\eqref{DEEST}, is
the function
\begin{equation}\label{WRBSY0perIKSMaismoBBeJ}
w^{(r)}(x):=v^{(r)}\left(\frac{x}{|x|}\right)=
u_0\left(r,\frac{x}{|x|}\right)=
u\left(\frac{rx}{|x|}\right).\end{equation} Notice that~$w^{(r)}(\lambda x)=
u_0\left(\frac{r\lambda x}{|\lambda x|}\right)=u_0\left(\frac{rx}{|x|}\right)=w^{(r)}(x)$,
whence~$w^{(r)}$ is positively homogeneous of degree zero.
Therefore, the function~$\phi^{(r)}:=\Delta w^{(r)}$ is homogeneous
of degree~$-2$.
Thus, recalling the Laplace-Beltrami definition
in~\eqref{BEL}, for every~$y\in\partial B_1$,
\begin{equation}\label{CQWSXRF3O5M4P4AS5S6784P234}
\Delta_{\partial B_1} u_0(r,y)=\Delta_{\partial B_1} v^{(r)}(y)=
\Delta w^{(r)}(y)=\phi^{(r)}(y)=r^2\,\phi^{(r)}(ry).\end{equation}
Additionally, we have that
the 
normal extension of~$u$ outside~$\partial B_r$
is the function~$\R^n\setminus\{0\}\ni x\mapsto u\left(\frac{rx}{|x|}\right)$,
which coincides with~$w^{(r)}(x)$, thanks to~\eqref{WRBSY0perIKSMaismoBBeJ}.
Therefore, by the Laplace-Beltrami definition
in~\eqref{BEL} we have that for every~$x\in\partial B_r$
$$ \Delta_{\partial B_r} u(x)=\Delta w^{(r)}(x)=\phi^{(r)}(x)
=\phi^{(r)}(ry),$$
where~$y:=\frac{x}r$. Hence, comparing with~\eqref{CQWSXRF3O5M4P4AS5S6784P234},
we conclude that~$\Delta_{\partial B_r} u(ry)=r^{-2}
\Delta_{\partial B_1} u_0(r,y)$ for every~$y\in\partial B_1$,
and this completes the proof of~\eqref{TRGSATTHNBELBRTYUTT9}.
\end{proof}

As a special case of
Theorem~\ref{SPH}, one obtains that if~$u\in C^2(\R^n\setminus\{0\})$ is radially symmetric,
i.e.~$u(x)=u_0(|x|)$ for some~$u_0:(0,+\infty)\to\R$,
then
\begin{equation}\label{ROTSE} \Delta u(x)=u_0''(r)+\frac{n-1}{r}\,u_0'(r),\end{equation}
with~$r=|x|$. Of course, this formula can also be obtained
by direct computations, see e.g.~\cite[page 21]{MR1625845}.\medskip

We also point out that if~$\alpha\in\R$ and~$u_\alpha(x):=|x|^\alpha$, then, taking~$\Sigma:=\partial B_1$,
we have that~$u_\alpha=1$ on~$\partial B_1$,
\begin{equation}\label{MN:ikmc}
\Delta u_\alpha=\alpha\,(n+\alpha-2),\qquad\Delta_T u_\alpha=\alpha\,(n-1)\qquad{\mbox{and}}\qquad
\Delta_{\partial B_1} u_\alpha=0\qquad{\mbox{on }}\,\partial B_1.
\end{equation}
In particular, notice that the normal
extension of the function identically
equal to~$1$ on~$\partial B_1$ according to definition~\eqref{DEEST}
is the function identically
equal to~$1$ on~$\R^n$, corresponding to~$u_\alpha$ with~$\alpha:=0$.
The other values of~$\alpha$ provide different extensions
of the function identically
equal to~$1$ on~$\partial B_1$ and, in this case,
the Laplace-Beltrami operator is not equal to the full Laplacian
of these extensions. Namely, geometric second order
operators can be reconstructed by full operators by extension,
but they are sensitive to the different type of extension chosen
(differently from first order operators,
as discussed in Lemma~\ref{PRIKD}).

\section{The Kelvin Transform}\label{LLOO2}

\begin{figure}
                \centering
                \includegraphics[width=.65\linewidth]{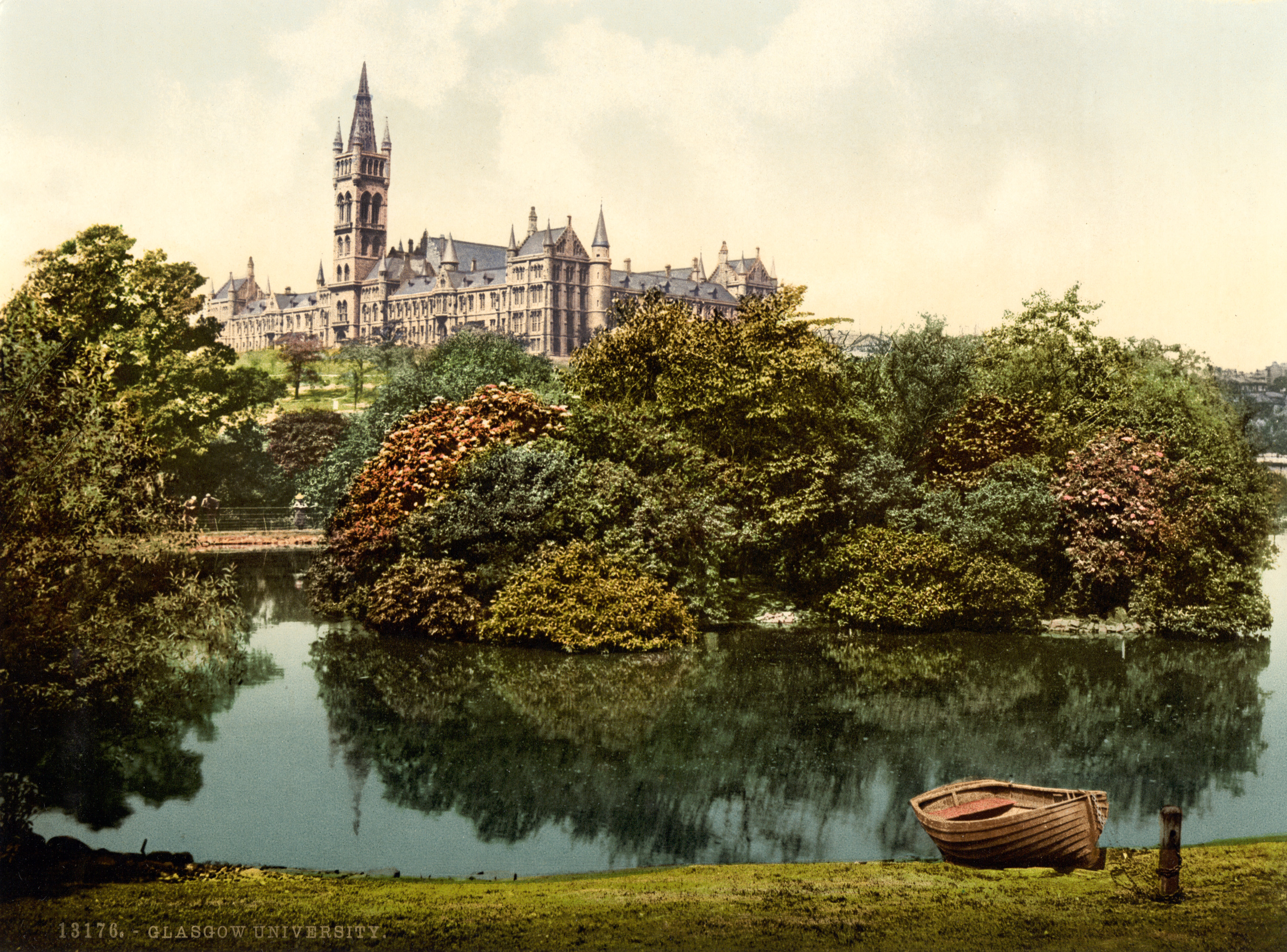}
        \caption{\sl Meander of the River Kelvin with the Gilmorehill campus of the University of Glasgow (Public Domain image from
        Wikipedia).}\label{RIVEKELVHARGIRA4AXELHARLROA7789GIJ7solFUMHDNOJHNFOJED}
\end{figure}

The Kelvin Transform\index{Kelvin Transform}
is a useful tool\footnote{The Kelvin Transform is named after 
William Thomson, 1st Baron Kelvin. Actually, Thomson himself
was named Baron Kelvin after the River Kelvin that flows past 
the University of Glasgow where Thomson used to work.
See Figure~\ref{RIVEKELVHARGIRA4AXELHARLROA7789GIJ7solFUMHDNOJHNFOJED}
for a photochrom print of the river and the university dating back to the end of nineteenth century.
See also Figure~\ref{KELVHARGIRA4AXELHARLROA7789GIJ7solFUMHDNOJHNFOJED}
for a caricature of Lord Kelvin (by caricaturist 
Sir Leslie Matthew Ward, a.k.a. Spy),
published in the magazine Vanity Fair in 1897.}
to reduce the analysis of partial differential
equations in exterior domains to that of interior ones
(and vice versa). 

Moreover, this transformation
possesses a number of geometric and analytic properties
that make it handy in several occasions.
To introduce\footnote{For the sake of clarity,
we presented the Kelvin Transform with the aim
of highlighting its remarkable analytic, algebraic and geometric
properties. On the other hand,
the Kelvin Transform is naturally motivated by some
important physical considerations inspired by the method of image
charges: the reader who wishes to go straight to this motivation
can look at Section~\ref{LOO}.

Also, for simplicity, we take~$\partial B_1$ as the reference set which
remains invariant for the Kelvin Transform, but we observe that a similar theory
can be carried out by leaving invariant~$\partial B_R$ (for this, one can either proceed by
scaling or change~\eqref{KSM:J:KELVIN1}
into~${\mathcal{K}}(x):=\frac{R^2x}{|x|^2}$.} it, we define, for all~$x\in\R^n\setminus\{0\}$,
\begin{equation}\label{KSM:J:KELVIN1}
{\mathcal{K}}(x):=\frac{x}{|x|^2}.\end{equation}
We list here some interesting properties of the Kelvin Transform:

\begin{lemma}\label{KSMLD0eroglh0345ptyhl0-bkfk4h578hrh}
The Kelvin Transform is an
involution, meaning that~${\mathcal{K}}({\mathcal{K}}(x))=x$
for all~$x\in\R^n\setminus\{0\}$.

Also, for all~$x$, $y\in\R^n\setminus\{0\}$,
\begin{eqnarray}
\label{KSM:J:KELVIN2}&& |{\mathcal{K}}(x)|\,|x|=1,\\
\label{KSM:J:KELVIN3}&& \frac{|x|^2}{1-|x|^2}=\frac{1}{|{\mathcal{K}}(x)|^2-1},\\
\label{KSM:J:KELVIN4}
&& \frac{{\mathcal{K}}(x)\cdot {\mathcal{K}}(y)}{|{\mathcal{K}}(x)|\,|{\mathcal{K}}(y)|}=\frac{x\cdot y}{|x|\,|y|},\\
\label{KSM:J:KELVIN5}&& |{\mathcal{K}}(x)-{\mathcal{K}}(y)|=\frac{|x-y|}{|x|\,|y|}.
\end{eqnarray}
Moreover\footnote{We point out that~\eqref{KSM:J:KELVIN4}
states that the Kelvin Transform is
angle preserving. Furthermore~\eqref{KSM:J:KELVIN8}
gives that the Kelvin Transform is conformal, since its Jacobian matrix
is a scalar function times an orthogonal matrix: indeed,\label{VEDHNDOCOSIJNDFSUCJFM}
$$ \sum_{k=1}^n
\left(
\delta_{ik}-\frac{2x_ix_k}{|x|^2}
\right)\left(
\delta_{kj}-\frac{2x_kx_j}{|x|^2}
\right)=\delta_{ij}.$$}
for each~$e\in\partial B_1$,
\begin{equation}\label{KSM:J:KELVIN9}
|x|\,|{\mathcal{K}}(x)-e|=|x-e|
,\end{equation}
and, for each~$i$, $j\in\{1,\dots,n\}$,
\begin{equation}\label{KSM:J:KELVIN8}
\partial_{x_j}{\mathcal{K}}(x)\cdot e_i=\frac{1}{|x|^2}\left(
\delta_{ij}-\frac{2x_ix_j}{|x|^2}
\right).
\end{equation}
\end{lemma}

\begin{proof} We remark that~\eqref{KSM:J:KELVIN2} is a direct consequence of the definition
in~\eqref{KSM:J:KELVIN1}.
Thus,
$$ {\mathcal{K}}({\mathcal{K}}(x))=\frac{{\mathcal{K}}(x)}{|{\mathcal{K}}(x)|^2}=
|x|^2{\mathcal{K}}(x)=x,$$
that proves the involution property.

{F}rom~\eqref{KSM:J:KELVIN2} we also deduce that
$$ \frac{1}{|{\mathcal{K}}(x)|^2-1}=
\frac{1}{|x|^{-2}-1}=\frac{|x|^2}{1-|x|^{2}},$$
that is~\eqref{KSM:J:KELVIN3}.

\begin{figure}
                \centering
                \includegraphics[width=.35\linewidth]{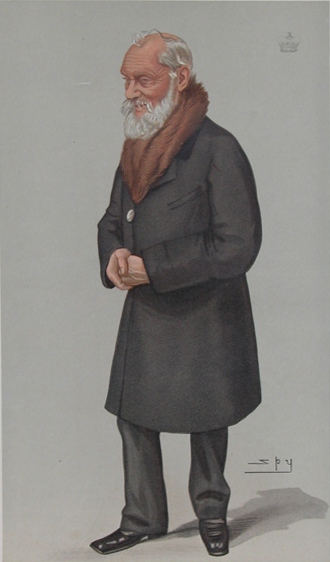}
        \caption{\sl Caricature of Lord Kelvin (Public Domain image from
        Wikipedia).}\label{KELVHARGIRA4AXELHARLROA7789GIJ7solFUMHDNOJHNFOJED}
\end{figure}

The claim in~\eqref{KSM:J:KELVIN4}
is also a straightforward byproduct of the definition
in~\eqref{KSM:J:KELVIN1} and~\eqref{KSM:J:KELVIN2}. In addition,
using~\eqref{KSM:J:KELVIN2} in combination with~\eqref{KSM:J:KELVIN4},
\begin{eqnarray*}&&
|x|^2|y|^2
|{\mathcal{K}}(x)-{\mathcal{K}}(y)|^2=
|x|^2|y|^2\big( |{\mathcal{K}}(x)|^2+|{\mathcal{K}}(y)|^2-2{\mathcal{K}}(x)\cdot {\mathcal{K}}(y)\big)\\&&\qquad
=|x|^2|y|^2\left( |{\mathcal{K}}(x)|^2+|{\mathcal{K}}(y)|^2-2|{\mathcal{K}}(x)|\,| {\mathcal{K}}(y)|
\frac{x\cdot y}{|x|\,|y|}
\right)=
|x|^2|y|^2\left( \frac{1}{|x|^2}+\frac1{|y|^2}-
\frac{2x\cdot y}{|x|^2\,|y|^2}
\right)\\&&\qquad
=\big( |y|^2+|x|^2-
2x\cdot y
\big)=
|x-y|^2,
\end{eqnarray*}
which establishes~\eqref{KSM:J:KELVIN5}.

Additionally, using~\eqref{KSM:J:KELVIN2} once again, if~$e\in\partial B_1$ then
\begin{eqnarray*}&& |x|^2\,|{\mathcal{K}}(x)-e|^2-|x-e|^2=
|x|^2\Big(|{\mathcal{K}}(x)|^2+1-2{\mathcal{K}}(x)\cdot e\Big)-
\Big(|x|^2+1-2x\cdot e\Big)\\&&\qquad=
1+|x|^2-2|x|^2 {\mathcal{K}}(x)\cdot e-\Big(|x|^2+1-2x\cdot e\Big)=
0,\end{eqnarray*}
and thus\footnote{Of course, many of the identities
in Lemma~\ref{KSMLD0eroglh0345ptyhl0-bkfk4h578hrh} can be proved
using different strategies. For instance, one
could also deduce~\eqref{KSM:J:KELVIN9}
directly from~\eqref{KSM:J:KELVIN5}
by noticing that, when~$|e|=1$,
$$ |{\mathcal{K}}(x)-e|=
|{\mathcal{K}}(x)-{\mathcal{K}}(e)|
=\frac{|x-e|}{|x|\,|e|}=\frac{|x-e|}{|x|}.$$} we have proved~\eqref{KSM:J:KELVIN9}.

Finally,
$$ \partial_{x_j}{\mathcal{K}}(x)\cdot e_i=\partial_{x_j}
\frac{x\cdot e_i}{|x|^2}
=\frac{|x|^2\delta_{ij}-2x_ix_j}{|x|^4},
$$
that is~\eqref{KSM:J:KELVIN8}.
\end{proof}

We remark that the Kelvin Transform
is an ``inversion of the sphere'', namely, in view of~\eqref{KSM:J:KELVIN2},
we have that~${\mathcal{K}}(B_1\setminus\{0\})=\R^n\setminus B_1$,
${\mathcal{K}}(\R^n\setminus B_1)=B_1\setminus\{0\}$
and~${\mathcal{K}}(\partial B_1)=\partial B_1$.

Also, the Kelvin Transform acts naturally on functions,
in a nicely compatible way with respect to the Laplace operator. Namely,
setting
\begin{equation}\label{SOLLE:0} u_{\mathcal{K}}(x):= |x|^{2-n}u({\mathcal{K}}(x)),\end{equation}
we have:

\begin{theorem}\label{ATTn12}
If~$v:=u_{\mathcal{K}}$, then~$v_{\mathcal{K}}=u$.

Moreover\footnote{In the second part of Theorem~\ref{ATTn12}, it is convenient to exclude
the case~$n=1$ to avoid integrability issues (from the technical point of view,
this condition is used in~\eqref{ATTn1}
to avoid contributions coming from infinity). As an illustrative example
of the loss of integrability in~\eqref{KS:09876543209876543lkjhgfduerhfnv3ue-0}
when~$n=1$, one can consider a function~$v\in C^\infty_0((-1,1))$ such that~$v(x)=1$ for each~$x\in\left[-\frac12,\frac12\right]$
and notice that, in view of~\eqref{SOLLE:0}, $v_{\mathcal{K}}(x)= |x| v\left(\frac{x}{|x|^2}\right)=|x|=x$
for each~$x\ge2$. Accordingly, in this case~$\nabla v_{\mathcal{K}}=1$ in~$(2,+\infty)$, whence
$$ \int_{\R^n} |\nabla
v_{\mathcal{K}}(x)|^2\,dx\ge\int_2^{+\infty} dx=+\infty.$$
However, the claim in~\eqref{KS:09876543209876543lkjhgfduerhfnv3ue}
is valid in dimension~$1$ as well, as it can be checked
by differentiation (though it is arguably not very useful in this case).}
if~$n\ge2$ and~$u\in C^1(\overline{B_1})$ (respectively, if~$u\in C^1(\R^n\setminus B_1)$),
then, for every~$v\in C^\infty_0(B_1)$
(respectively, $v\in C^\infty_0(\R^n\setminus B_1)$,
\begin{equation}\label{KS:09876543209876543lkjhgfduerhfnv3ue-0}
\int_{\R^n}\nabla u_{\mathcal{K}}(x)\cdot\nabla
v_{\mathcal{K}}(x)\,dx=
\int_{\R^n}\nabla u(x)\cdot\nabla
v(x)\,dx
.\end{equation}

In addition to that,
if~$u\in C^2(\overline{B_1})$ (respectively, if~$u\in C^2(\R^n\setminus B_1)$), for every~$x\in\R^n\setminus B_1$ (respectively,
for every~$x\in B_1\setminus\{0\}$) it holds that
\begin{equation}\label{KS:09876543209876543lkjhgfduerhfnv3ue}
\Delta u_{\mathcal{K}}(x)=
\frac{1}{|x|^{n+2}}\,
\Delta u\big({\mathcal{K}}(x)\big).\end{equation}
\end{theorem}

\begin{proof} We observe that, if~$v:=u_{\mathcal{K}}$, then
$$ v({\mathcal{K}}(x))=u_{\mathcal{K}}({\mathcal{K}}(x))=
|{\mathcal{K}}(x)|^{2-n}u\big({\mathcal{K}}({\mathcal{K}}(x))\big)=|x|^{n-2}
u(x),
$$
thanks to~\eqref{SOLLE:0} and Lemma~\ref{KSMLD0eroglh0345ptyhl0-bkfk4h578hrh},
and, for this reason,
$$ v_{\mathcal{K}}(x)=
|x|^{2-n}v({\mathcal{K}}(x))=
u(x).$$
This proves the desired involution property.

Now we establish~\eqref{KS:09876543209876543lkjhgfduerhfnv3ue-0}
and~\eqref{KS:09876543209876543lkjhgfduerhfnv3ue}.
For this, we consider~$v\in C^\infty_0(B_1)$
(respectively, $v\in C^\infty_0(\R^n\setminus B_1)$.
We let~$M(x)$ be the matrix with entries~$
\delta_{ij}-\frac{2x_ix_j}{|x|^2}$ for all~$i$, $j\in\{1,\dots,n\}$
and we recall that~$M(x)$ is orthogonal (see the footnote
on page~\pageref{VEDHNDOCOSIJNDFSUCJFM}).
As a result, we have that~$\det M(x) =1$.
Using~\eqref{KSM:J:KELVIN8}, we see that the change
of variable~$y={\mathcal{K}}(x)$ leads to~$dy=\left|
\det D{\mathcal{K}}(x)\right|\,dx=\left|
\det\frac{M(x)}{|x|^2}\right|\,dx
=\frac{dx}{|x|^{2n}}$.

Moreover, by~\eqref{KSM:J:KELVIN8}, 
$$ \nabla u_{\mathcal{K}}(x)=\nabla\Big( |x|^{2-n}u({\mathcal{K}}(x))\Big)=
(2-n)|x|^{-n}xu({\mathcal{K}}(x))+|x|^{-n}M(x)\nabla u({\mathcal{K}}(x)),$$
and a similar formula holds for~$v$ replacing~$u$.

Therefore,
\begin{equation}\label{JHSNDsikcmsdSfgaisdwedsagfdsgf39450ksdxc}\begin{split}&
\int_{\R^n}\nabla u_{\mathcal{K}}(x)\cdot\nabla
v_{\mathcal{K}}(x)\,dx=
\int_{\R^n}
\Big[
(2-n)^2|x|^{2-2n}u({\mathcal{K}}(x))v({\mathcal{K}}(x))\\&\qquad
+(2-n) |x|^{-2n}u({\mathcal{K}}(x))\big(M(x)\nabla v({\mathcal{K}}(x))\big)\cdot x
\\&\qquad+(2-n) |x|^{-2n}v({\mathcal{K}}(x))\big(M(x)\nabla u({\mathcal{K}}(x))\big)\cdot x
\\&\qquad+|x|^{-2n}\big(M(x)\nabla u({\mathcal{K}}(x))\big)\cdot
\big(M(x)\nabla v({\mathcal{K}}(x))\big)
\Big]\,dx.
\end{split}\end{equation}
We also remark that
\begin{eqnarray*}&&
\div\Big(|x|^{2-2n}x u({\mathcal{K}}(x))v({\mathcal{K}}(x))\Big)\\&=&
\div\Big(|x|^{2-2n}x\Big) u({\mathcal{K}}(x))v({\mathcal{K}}(x))\\&&\qquad
+|x|^{2-2n}x\cdot\nabla\Big( u({\mathcal{K}}(x))\Big)
v({\mathcal{K}}(x))+
|x|^{2-2n}x u({\mathcal{K}}(x))\cdot\nabla\Big(v({\mathcal{K}}(x))\Big)\\&=&
(2-n)|x|^{2-2n}u({\mathcal{K}}(x))v({\mathcal{K}}(x))\\&&\qquad+
|x|^{-2n}x\cdot\big(M(x)\nabla u({\mathcal{K}}(x))\big)
v({\mathcal{K}}(x))
+|x|^{-2n}x\cdot\big(M(x)\nabla v({\mathcal{K}}(x))\big)
u({\mathcal{K}}(x)).
\end{eqnarray*}
This and the Divergence Theorem give that
\begin{equation}\label{ATTn1}
\begin{split}&
\int_{\R^n}
\Big[
(2-n)^2|x|^{2-2n}u({\mathcal{K}}(x))v({\mathcal{K}}(x))
+(2-n) |x|^{-2n}u({\mathcal{K}}(x))\big(M(x)\nabla v({\mathcal{K}}(x))\big)\cdot x
\\&\qquad+(2-n) |x|^{-2n}v({\mathcal{K}}(x))\big(M(x)\nabla u({\mathcal{K}}(x))\big)\cdot x
\Big]\,dx=0.\end{split}
\end{equation}
Plugging this information into~\eqref{JHSNDsikcmsdSfgaisdwedsagfdsgf39450ksdxc},
we thus conclude that
$$ \int_{\R^n}\nabla u_{\mathcal{K}}(x)\cdot\nabla
v_{\mathcal{K}}(x)\,dx=
\int_{\R^n}|x|^{-2n}\big(M(x)\nabla u({\mathcal{K}}(x))\big)\cdot
\big(M(x)\nabla v({\mathcal{K}}(x))\big)\,dx.$$
Therefore, using the orthogonality of~$M(x)$
and changing variable,
\begin{eqnarray*}
\int_{\R^n}\nabla u_{\mathcal{K}}(x)\cdot\nabla
v_{\mathcal{K}}(x)\,dx&=&\int_{\R^n}|x|^{-2n}\nabla u({\mathcal{K}}(x))\cdot\nabla
v({\mathcal{K}}(x))\,dx\\
&=&\int_{\R^n}\nabla u(y)\cdot\nabla
v(y)\,dy.
\end{eqnarray*}
This proves~\eqref{KS:09876543209876543lkjhgfduerhfnv3ue-0}.

Now we prove~\eqref{KS:09876543209876543lkjhgfduerhfnv3ue}.
For this, we utilize~\eqref{KS:09876543209876543lkjhgfduerhfnv3ue-0}
to see that
\begin{eqnarray*}
&&
-\int_{\R^n}|{\mathcal{K}}(y)|^{n+2}\,\Delta u_{\mathcal{K}}({\mathcal{K}}(y))\,
v(y)\,dy=
-\int_{\R^n}|y|^{-n-2}\,\Delta u_{\mathcal{K}}({\mathcal{K}}(y))\,
v(y)\,dy\\&&\qquad=
-\int_{\R^n}|x|^{2-n}\,\Delta u_{\mathcal{K}}(x)\,
v({\mathcal{K}}(x))\,dx=
-\int_{\R^n}\Delta u_{\mathcal{K}}(x)\,
v_{\mathcal{K}}(x)\,dx\\&&\qquad
=
\int_{\R^n}\nabla u_{\mathcal{K}}(x)\cdot\nabla
v_{\mathcal{K}}(x)\,dx=\int_{\R^n}\nabla u(y)\cdot\nabla
v(y)\,dy
=-
\int_{\R^n}\Delta u(y)\,
v(y)\,dy.
\end{eqnarray*}
This yields that~$|{\mathcal{K}}(y)|^{n+2}\,\Delta u_{\mathcal{K}}({\mathcal{K}}(y))=
\Delta u(y)$, from which we obtain~\eqref{KS:09876543209876543lkjhgfduerhfnv3ue}.
\end{proof}

Another interesting geometric property of the Kelvin Transform is that it
carries spheres and hyperplanes into spheres and hyperplanes:

\begin{proposition}\label{h1ds3rfdsK3456790-dskpodnglkbd3565EL-DIFfffI}
The Kelvin Transform carries
\begin{itemize}
\item spheres not passing through the origin
into spheres not passing through the origin,
\item spheres passing through the origin into hyperplanes not passing through the origin,
\item and hyperplanes not passing through the origin into
spheres passing through the origin.
\end{itemize}
Also, the Kelvin Transform leaves invariant all the
hyperplanes passing through the origin.\medskip

More explicitly, if~$p\in\R^n$, $r\in(0,+\infty)$, $\omega\in\partial B_1$ and~$c\in\R$ we have that
\begin{equation}\label{4567k-0olk-45rt-HSMgAolfLIkdgnN89kAm-01}
\begin{split}
&{\mathcal{K}}\Big( \{ x\in\R^n\setminus\{0\} {\mbox{ s.t. }} |x-p|^2=r^2\}\Big)\\
&\qquad\qquad=\begin{dcases}
\displaystyle\left\{ x\in\R^n\setminus\{0\} {\mbox{ s.t. }} \left|x-\frac{p}{|p|^2-r^2}\right|^2=\frac{r^2}{(|p|^2-r^2)^2}\right\}&{\mbox{ if }}r\ne|p|
,\\
\\
\displaystyle\left\{ x\in\R^n\setminus\{0\} {\mbox{ s.t. }} p\cdot x=\frac12\right\}&{\mbox{ if }}r=|p|,
\end{dcases}\end{split}\end{equation}
and
\begin{equation}\label{4567k-0olk-45rt-HSMgAolfLIkdgnN89kAm-02}
\begin{split}
&{\mathcal{K}}\Big( \{ x\in\R^n\setminus\{0\} {\mbox{ s.t. }} \omega\cdot x=c\}\Big)\\&\qquad\qquad=\begin{dcases}
\displaystyle\left\{ x\in\R^n\setminus\{0\} {\mbox{ s.t. }} \left|x-\frac{\omega}{2c}\right|^2=\frac{1}{4c^2}\right\}&{\mbox{ if }}c\ne0
,\\
\\
\displaystyle\left\{ x\in\R^n\setminus\{0\} {\mbox{ s.t. }} \omega\cdot x=0\right\}&{\mbox{ if }}c=0.
\end{dcases}\end{split}\end{equation}
\end{proposition}

\begin{proof} Let~$x\ne0$ and~$y:={\mathcal{K}}(x)$.
Since, by Lemma~\ref{KSMLD0eroglh0345ptyhl0-bkfk4h578hrh},
the Kelvin Transform is an involution, we have that~$x={\mathcal{K}}(y)=\frac{y}{|y|^2}$. Therefore, if~$|x-p|^2=r^2$ then
\begin{equation*}
\begin{split}&
0=|y|^2\left( |x-p|^2-r^2\right)=|y|^2\left( \left|\frac{y}{|y|^2}-p\right|^2-r^2\right)\\&\qquad\qquad
=\left|\frac{y}{|y|}-|y|p\right|^2-r^2|y|^2=
1-2 p\cdot y+(|p|^2-r^2)|y|^2.\end{split}
\end{equation*}
Now, if~$r=|p|$ this boils down to~$p\cdot y=\frac12$. If instead~$r\ne|p|$, then
\begin{eqnarray*}&&
0=\frac{1}{|p|^2-r^2}-\frac{2 p\cdot y}{|p|^2-r^2}+|y|^2=
\frac{1}{|p|^2-r^2}-\frac{|p|^2}{(|p|^2-r^2)^2}+\left|y-\frac{p}{|p|^2-r^2}\right|^2\\&&\qquad\qquad=-
\frac{r^2}{(|p|^2-r^2)^2}+\left|y-\frac{p}{|p|^2-r^2}\right|^2.
\end{eqnarray*}

Furthermore, we notice that, since~$y:={\mathcal{K}}(x)=\frac{x}{|x|^2}$,
\begin{eqnarray*}
|y-p|^2=|y|^2-2p\cdot y +|p|^2=\frac1{|x|^2}-\frac{2p\cdot x}{|x|^2}+|p|^2.
\end{eqnarray*}
Thus, if~$r=|p|$ and~$p\cdot x=\frac12$, then
$$|y-p|^2=
\frac1{|x|^2}-\frac{1}{|x|^2}+|p|^2=|p|^2=r^2.$$
If instead~$r\ne|p|$ and~$\left|x-\frac{p}{|p|^2-r^2}\right|^2=\frac{r^2}{(|p|^2-r^2)^2}$, then
\begin{eqnarray*}
|y-p|^2=\frac1{|x|^2}-\frac{1}{|x|^2}\big(|x|^2(|p|^2-r^2)+1\big)+|p|^2
=\frac1{|x|^2}-|p|^2+r^2-\frac1{|x|^2}
+|p|^2=r^2.
\end{eqnarray*}
These considerations establish~\eqref{4567k-0olk-45rt-HSMgAolfLIkdgnN89kAm-01}.

Now, if~$\omega\cdot x=c$, we find that~$\omega\cdot y=c|y|^2$.
This reduces to~$\omega\cdot y=0$ if~$c=0$. Instead, if~$c\ne0$,
$$ 0=|y|^2-\frac{\omega\cdot y}{c}=
\left| y-\frac{\omega}{2c}\right|^2-\frac{1}{4c^2}.$$
Viceversa, if~$\omega\cdot x=0$, then~$\omega\cdot y=0=c$.
If instead~$\left| x-\frac{\omega}{2c}\right|^2=\frac{1}{4c^2}$, then
$$\omega\cdot y=\frac{\omega\cdot x}{|x|^2}=\frac{c|x|^2}{|x|^2}=c,
$$
thus completing the proof of~\eqref{4567k-0olk-45rt-HSMgAolfLIkdgnN89kAm-02}.
\end{proof}

\begin{figure}
  \centering
  \includegraphics[width=.6\linewidth]{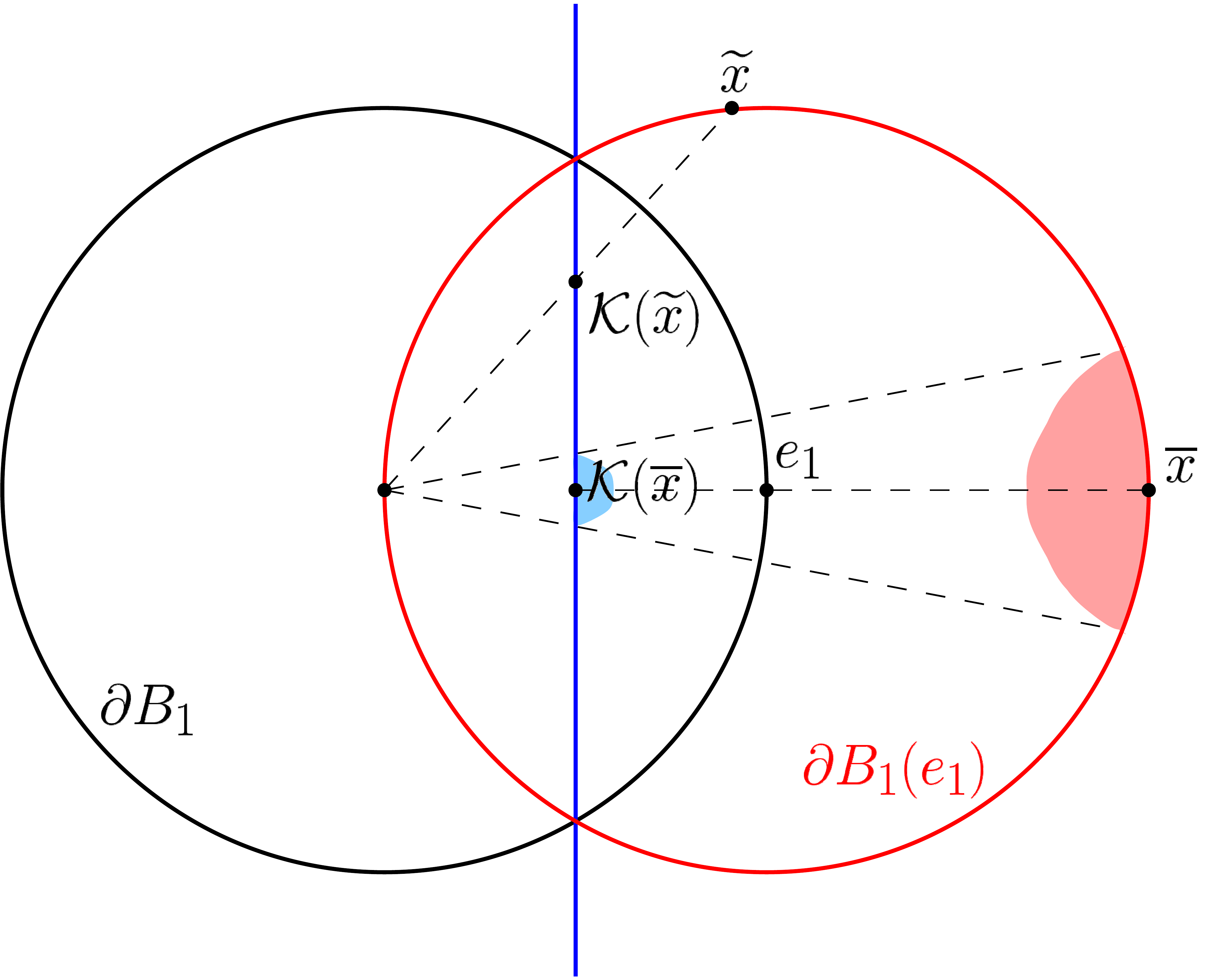}
 \caption{\sl The Kelvin Transform carries the sphere~$\partial B_1(e_1)$ into the hyperplane~$\left\{x_1=\frac12\right\}$.
 A small neighborhood of~$2e_1$ in~$B_1(e_1)$ is carried into a small neighborhood of~$\frac{e_1}2$
 in~$\left\{x_1>\frac12\right\}$.}  \label{K6EL-DIFfffI}
\end{figure}

We refer to Figure~\ref{K6EL-DIFfffI} for a graphical representation
of the situation discussed in detail in Proposition~\ref{h1ds3rfdsK3456790-dskpodnglkbd3565EL-DIFfffI}.

\section{The fundamental solution}\label{lfundsP-S}

Here we will describe the notion of fundamental solution\index{fundamental solution of the Laplace equation} of the Laplace\footnote{Usually, \label{FOO:POIEQUATYHN}
though the notation is certainly not uniform across the literature, the equation~$\Delta u(x)=0$
for all~$x$ in~$\Omega$ is often referred to with the name of ``Laplace equation'', \index{Laplace equation}
and correspondingly~$\Delta u(x)=f(x)$ with the name of ``Poisson equation''. \index{Poisson equation}
When the
Laplace equation (or, sometimes, the Poisson equation) is complemented
with a boundary datum~$u=g$ along~$\partial\Omega$ people speak about the ``Dirichlet problem''.}
equation,
which is physically motivated by the electrostatic (or gravitational) potential produced
by a point charge (or a point mass). For this, the isotropy and homogeneity of the ambient space
play a crucial role by inducing rotational and translational symmetries.

To start with, in view of the polar coordinates representation of the Laplacian in~\eqref{ROTSE}, 
one can explicitly find all the rotationally symmetric harmonic functions in~$\R^n\setminus\{0\}$,
according to the following observation:

\begin{lemma}\label{MHNvKSD}
If~$v\in C^2(\R^n\setminus\{0\})$ (respectively, if~$v\in C^2(B_R\setminus\{0\})$
for some ball~$B_R\subset\R^n$)
is rotationally symmetric and harmonic in~$\R^n\setminus\{0\}$
(respectively, in~$B_R\setminus\{0\}$),
then there exist~$a$, $b\in\R$ such that, for every~$x\in\R^n\setminus\{0\}$
(respectively, for every~$x\in B_R\setminus\{0\}$), we have that
$$ v(x)=\begin{dcases}
\frac{a}{|x|^{n-2}}+b&{\mbox{ if }}n\ne2,\\
a\ln |x|+b&{\mbox{ if }}n=2.
\end{dcases}$$
\end{lemma}

\begin{proof} We argue in~$\R^n\setminus\{0\}$, the case in~$B_R\setminus\{0\}$ being similar.
We use the notation~$r:=|x|$ and let~$v_0:(0,+\infty)\to\R$ be such that~$v(x)=v_0(|x|)$.
Then, by~\eqref{ROTSE}, in~$\R^n\setminus\{0\}$,
\begin{equation*}
0=\Delta v=v_0''+\frac{n-1}r\,v_0'=r^{1-n}\,\frac{d}{dr} \big( r^{n-1} v_0'\big),
\end{equation*}
therefore there exists~$a\in\R$ such that~$r^{n-1}v_0'(r)=a$ for every~$r>0$. Integrating this expression,
we obtain that there exists~$b\in\R$ such that
$$ v_0(r)=\begin{dcases}
\frac{a r^{2-n}}{2-n}+b&{\mbox{ if }}n\ne2,\\
a\ln r+b&{\mbox{ if }}n=2.
\end{dcases}$$
Hence, by renaming the constant~$a$,
\begin{equation*} v_0(r)=\begin{dcases}
\frac{a }{r^{n-2}}+b&{\mbox{ if }}n\ne2,\\
a\ln r+b&{\mbox{ if }}n=2.
\end{dcases}\end{equation*}
and the desired result plainly follows.
\end{proof}

\begin{figure}
  \centering
  \includegraphics[width=.4\linewidth]{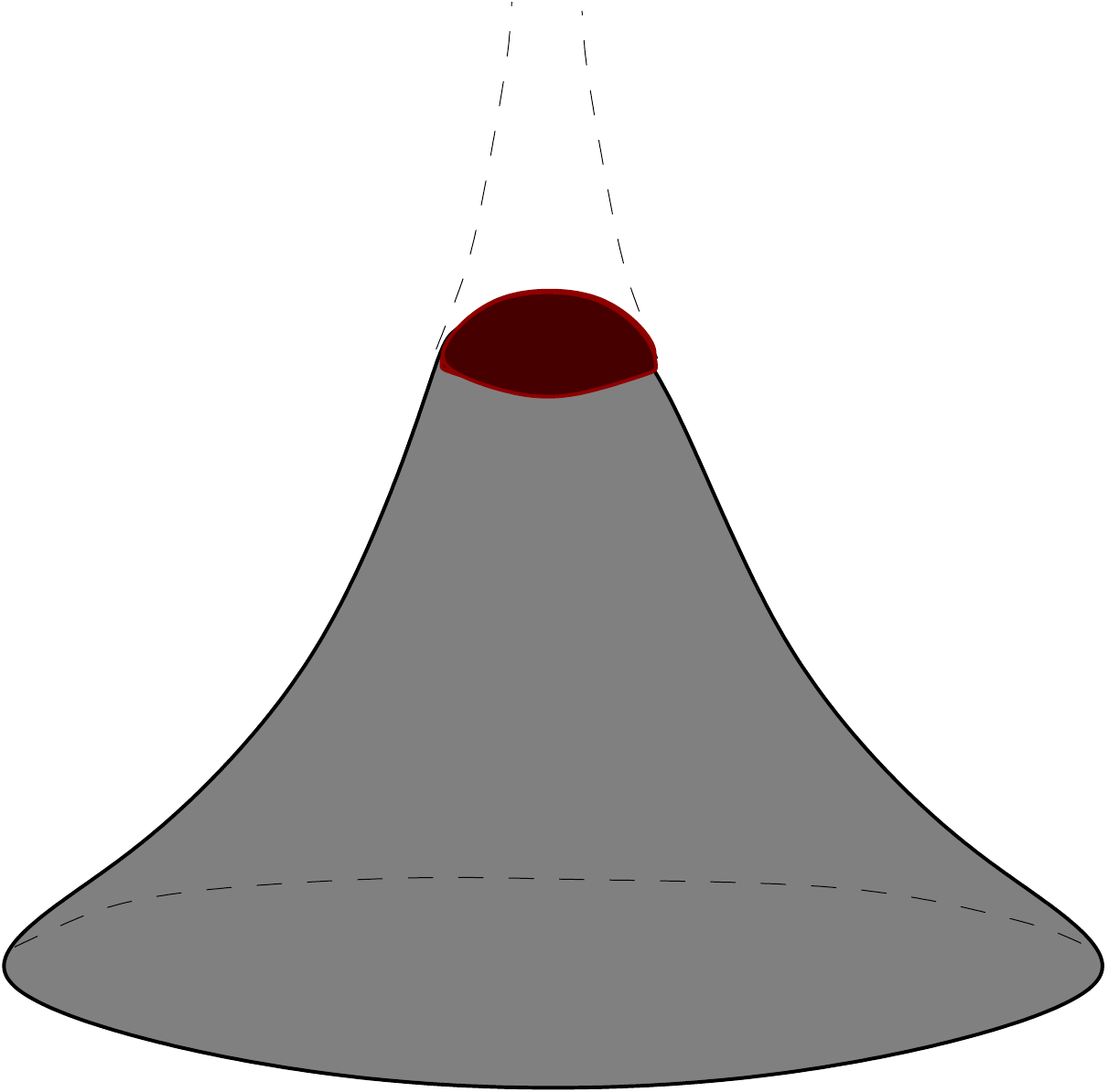}
 \caption{\sl The function~$\Gamma_\rho$ (if~$n\ge3$).}  \label{DIFfffI}
\end{figure}

The functions introduced in Lemma~\ref{MHNvKSD}
play a pivotal role in the development of the theory of elliptic
partial differential equations, since not only, as observed in Lemma~\ref{MHNvKSD},
they satisfy pointwise outside the origin the equation~$\Delta v=0$,
but also, as we will see now, up to normalization constants,
they provide the ``fundamental solutions'' of the Laplace operators,
that is their Laplacian (in a suitable distributional sense,
as we will precise in Theorem~\ref{KScvbMS:0okrmt4ht}) is (minus)
the Dirac Delta Function
at the origin.
To clarify this concept, in light of Lemma~\ref{MHNvKSD},
it is useful to define
\begin{equation}\label{OVERLINEV} \overline v(r)=\begin{dcases}
\frac{1}{r^{n-2}}&{\mbox{ if }}n\ne2,\\
-\ln r&{\mbox{ if }}n=2
\end{dcases}\end{equation}
and choose the normalizing constant
\begin{equation}\label{COIENEE} c_n:=\begin{dcases}\frac1{n(n-2)\,| B_1|}
& {\mbox{ if }} n\ne2,\\
\frac1{2\pi} & {\mbox{ if }} n=2.\end{dcases}\end{equation}
The reason for choosing~$c_n$ in this way is that, for every~$\rho>0$, the function
\begin{equation}\label{COIENEE-0986ytufgkbv-0-rjfeonvnb2GDB} \overline\Gamma_\rho(r):=\begin{dcases}
\frac{\rho^2-r^2}{2n |B_\rho|}+c_n \,\overline v(\rho) & {\mbox{ if }}r\in(0,\rho),\\ c_n\,
\overline v(r) & {\mbox{ if }}r\in[\rho,+\infty)
\end{dcases}\end{equation}
is such that~$\overline\Gamma_\rho'\in C^{0,1}((0,+\infty))$. Let also
\begin{equation}\label{GAROGRA}
\Gamma_\rho(x):=\overline\Gamma_\rho(|x|).\end{equation}
We remark that~$\Gamma_\rho$ is obtained, roughly speaking, by ``gluing a paraboloid near the origin''
to the harmonic function outside the origin that was introduced in
Lemma~\ref{MHNvKSD}, see Figure~\ref{DIFfffI}. Moreover, this paraboloid is normalized to
have Laplacian that integrates to~$-1$ in~$B_\rho$, namely~$\Gamma_\rho$
can be extended to a function in~$C^{1,1}(\R^n)$ satisfying
\begin{equation}\label{DEGAM} \Delta \Gamma_\rho=\begin{dcases}
-\frac{1}{|B_\rho|}&{\mbox{ in }}B_\rho,\\
0&{\mbox{ in }}\R^n\setminus\overline{B_\rho}.
\end{dcases}\end{equation}
We then consider the formal limit as~$\rho\searrow0$ of~$\Gamma_\rho$, namely we define
\begin{equation}\label{GAMMAFU} \Gamma(x):=c_n\overline v(|x|)=
\begin{dcases}
\frac{c_n}{|x|^{n-2}}&{\mbox{ if }}n\ne2,\\
-c_n\ln |x|&{\mbox{ if }}n=2
.\end{dcases}\end{equation}
We stress that when~$n=3$ this function represents, up to a
normalizing constant, the electrostatic potential
generated by a point charge (as well as the gravitational potential generated by a point mass,
or the equilibrium temperature produced by a concentrated heat source),
and we can exploit this physical motivation in any dimension as well,
at least\footnote{It is
usually a challenging task to understand higher dimensions, and the change of the fundamental
solution when~$n=2$ is a deep feature to keep in mind. To ease the intuition and possibly recover, at least
at a heuristic level, lower dimensional cases from higher dimensional ones, we propose some
reflections about how electrostatic potentials of a line of uniformly distributed charges
in~$\R^n$ produces the fundamental solutions in~$\R^{n-1}$ (up to
some renormalization that plays a role in low dimension). For this,
we consider a uniform distribution of charges
on the line~$L:=\{x=(x',x_n)\in\R^{n-1}\times\R$
s.t. $x'=0\}$, which (up to physical constant) produces an electrostatic
potential at the point~$p=(p',0)\in\R^{n-1}\times\{0\}$ of the type
$$ U(r):=\frac1{c_n}\int_L \Gamma_n(p-x)\,d{\mathcal{H}}^{n-1}_x
=\int_L \overline{v}(|p-x|)\,d{\mathcal{H}}^{n-1}_x
=\int_{-\infty}^{+\infty}\overline{v}_n\left(\sqrt{|p'|^2+x_n^2}\right)\,dx_n
=\int_{-\infty}^{+\infty}\overline{v}_n\left(\sqrt{r^2+x_n^2}\right)\,dx_n,$$
where~$r:=|p'|$, and~$\overline{v}_n$ and~$\Gamma_n$
are as in~\eqref{OVERLINEV} and~\eqref{GAMMAFU},
with the subscript~$n$ exploited to underline the dependence of these functions upon the dimension.

Strangely enough, the easiest case to understand is the high dimensional one: namely,
when~$n\ge4$ we obtain
$$ U(r)=\int_{-\infty}^{+\infty}\left( r^2+x_n^2\right)^{\frac{2-n}2}\,dx_n=
r^{1-n}
\int_{-\infty}^{+\infty}\left( 1+t^2\right)^{\frac{2-n}2}\,dt=cr^{3-n}=c\overline{v}_{n-1}(r),
$$
where~$c>0$, showing that the potential of the charged line produces
the fundamental solution in one dimension less, up to a normalizing constant.

The lower dimensional cases~$n\in\{2,3\}$ are instead more tricky,
since the corresponding constant~$c$ would diverge. To make the argument
work (at least at a heuristic level) one needs to perform a ``renormalization''
procedure, formally ``subtracting infinity'' to the potential (as a matter of fact,
potentials are always defined ``up to additive constants'', since the physical forces
come from their derivatives). Concretely, when~$n\in\{2,3\}$ one has to replace the previous definition of~$U$
by the following renormalized one:
$$ U(r):=\lim_{R\to+\infty}\left(
\int_{-R}^{+R}\overline{v}_n\left(\sqrt{r^2+x_n^2}\right)\,dx_n-\phi_n(R)\right),$$
for a suitable renormalization~$\phi_n(R)$. The previous computations
give that one can choose~$\phi_n$ to be identically zero when~$n\ge3$,
but, as we will see now, the cases~$n\in\{2,3\}$ do require a more specific choice.

If~$n=3$, we choose~$\phi(R):=2\ln(2R)$, and we thereby obtain
\begin{eqnarray*}
U(r)&=&\lim_{R\to+\infty}\left(
\int_{-R}^{+R}\frac1{\sqrt{r^2+x_n^2}}\,dx_n-2\ln(2R)\right)\\&=&
\lim_{R\to+\infty}\left(\ln\left(1+\frac{2 R (\sqrt{r^2 + R^2} + R)}{r^2} \right)
-\ln(4R^2)\right)
\\&=&\lim_{R\to+\infty}
\ln\frac{r^2+2 R (\sqrt{r^2 + R^2} + R)}{4r^2 R^2} \\&=&
\ln\frac{1}{r^2}\\&=&2\overline{v}_2 (r).
\end{eqnarray*}
If instead~$n=2$, we choose~$\phi(R):=-2R(\ln R+1)$, and we have
\begin{eqnarray*}
U(r)&=&\lim_{R\to+\infty}\left(
-\frac12\int_{-R}^{+R}\ln\big(r^2+x_n^2\big)\,dx_n+2R(\ln R+1)\right)\\
&=&\lim_{R\to+\infty}\left(
-2 r \arctan\frac{R}r -R \ln (r^2 + R^2)+2R\ln R\right)\\
&=&\lim_{R\to+\infty}\left(
-2 r \arctan\frac{R}r -R \ln \frac{r^2 + R^2}{R^2}\right)\\
&=&\lim_{R\to+\infty}\left(
-2 r \arctan\frac{R}r -R \ln \left(1+\frac{r^2}{R^2}\right)\right)\\
&=&-\pi r\\&=&-\pi \overline{v}_1(r).
\end{eqnarray*}
The negative sign here above is consistent with the fact that the constant
in~\eqref{COIENEE} is negative when~$n=1$.

Related approaches for constructing lower dimensional effective potentials
deal with the interaction of a charged line with a ``test line'' (rather than a ``test point'').

Another approach to recover dimension~$1$ directly from dimension~$3$ is to consider the renormalized
electrostatic potential of a charged plate (say,~$\R^2\times\{0\}$) at the point~$(0,0,r)$: in this case
the computation, up to dimensional constants, would involve a corrector of the type~$2\pi R$ and go as follows:
\begin{eqnarray*}
U(r)&=&\lim_{R\to+\infty} \left[\int_0^R \left(\int_{\{|x'|=\rho,\;x_3=0\}}\frac{d{\mathcal{H}}^{2}_x}{|x-(0,0,r)|}\right)\,d\rho-2\pi R\right]\\
&=&\lim_{R\to+\infty} \left[\int_0^R \frac{2\pi\rho\,d\rho}{\sqrt{\rho^2+r^2}}-2\pi R\right]\\&=&
\lim_{R\to+\infty}2\pi\sqrt{R^2 + r^2}-2\pi r-2\pi R\\
&=&-2\pi r\\&=&-2\pi \overline{v}_1(r).
\end{eqnarray*}}
to facilitate our mathematical intuition.
In our setting, the function~$\Gamma$ is the ``fundamental solution''
of the Laplace operator, in the sense made precise by the following result:

\begin{theorem}\label{KScvbMS:0okrmt4ht}
For every~$\varphi\in C^\infty_0(\R^n)$,
$$ \int_{\R^n}\Gamma(x)\, \Delta\varphi(x) \,dx=-\varphi(0).$$
\end{theorem}

\begin{proof} Since~$\overline\Gamma_\rho'\in C^{0,1}((0,+\infty))$,
there exists~$C>0$ such that~$|\overline\Gamma_\rho'|+|\overline\Gamma_\rho''|\le C$ a.e. in~$\R^n$.
As a result, for every~$i$, $j\in\{1,\dots,n\}$, and a.e.~$x\in\R^n\setminus B_{\rho/2}$,
$$ |\partial_{ij}\Gamma_\rho(x)|=\left|
\overline\Gamma_\rho''(|x|)\,\frac{x_i\,x_j}{|x|^2}
+
\overline\Gamma_\rho'(|x|)\,\frac{|x|^2\delta_{ij}-x_i\,x_j}{|x|^3}
\right|\le 4C\left(1+\frac1\rho\right).
$$
As a result, letting~$A_{\rho,\e}:=B_{\rho+\e}\setminus B_{\rho-\e}$,
\begin{equation}\label{OK-plkdm-1}\lim_{\e\searrow0}
\left|
\int_{A_{\rho,\e}}\Delta\Gamma_\rho(x)\,\varphi(x) \,dx
\right|\le \lim_{\e\searrow0}
4Cn\,\|\varphi\|_{L^\infty(\R^n)}\left(1+\frac1\rho\right)\,\big|A_{\rho,\e}\big|=0.
\end{equation}
Also, 
\begin{equation}\label{OK-plkdm-2}\lim_{\e\searrow0}
\left|
\int_{A_{\rho,\e}}\Gamma_\rho(x)\,\Delta\varphi(x) \,dx
\right|\le \lim_{\e\searrow0}n\,\|\Gamma_\rho\|_{L^\infty(S)}
\,\|D^2\varphi\|_{L^\infty(\R^n)}\,\big|A_{\rho,\e}\big|=0,
\end{equation}
where~$S$ is the support of~$\varphi$.

Now we point out that if~$F:\R^n\to\R$ is a continuous vector field, then
\begin{equation}\label{CVF}\lim_{\e\searrow0}
\int_{\partial A_{\rho,\e}} F(x)\cdot\nu(x)\,d{\mathcal{H}}^{n-1}_x=0.
\end{equation}
For this, given~$\eta>0$, we consider a smooth vector field~$F_\eta$ such that~$\|F-F_\eta\|_{L^\infty(B_{2\rho},\R^n)}\le\eta$ and we use the Divergence Theorem to see that
\begin{eqnarray*}&&
\lim_{\e\searrow0}\left|
\int_{\partial A_{\rho,\e}} F(x)\cdot\nu(x)\,d{\mathcal{H}}^{n-1}_x
\right|\\&\le&
\lim_{\e\searrow0}
\|F-F_\eta\|_{L^\infty(\partial A_{\rho,\e},\R^n)}\,{\mathcal{H}}^{n-1}(\partial A_{\rho,\e})+
\left|
\int_{\partial A_{\rho,\e}} F_\eta(x)\cdot\nu(x)\,d{\mathcal{H}}^{n-1}_x
\right|
\\&\le&2\eta{\mathcal{H}}^{n-1}(\partial B_\rho)+\lim_{\e\searrow0}
\left|
\int_{ A_{\rho,\e}}\div F_\eta(x)\,dx
\right|\\&\le&
2\eta{\mathcal{H}}^{n-1}(\partial B_\rho)+\lim_{\e\searrow0}
n\,\|F_\eta\|_{C^1(B_{2\rho},\R^n)}\,| A_{\rho,\e}|
\\&=&2\eta{\mathcal{H}}^{n-1}(\partial B_\rho).
\end{eqnarray*}
We now send~$\eta\searrow0$ and we conclude the proof of~\eqref{CVF}.

As a consequence of~\eqref{CVF}, we find that
$$ \lim_{\e\searrow0}
\int_{\partial A_{\rho,\e}}
\Gamma_\rho(x)\frac{\partial \varphi}{\partial\nu}(x) \,d{\mathcal{H}}^{n-1}_x
=0\qquad{\mbox{ and }}\qquad
\lim_{\e\searrow0}
\int_{\partial A_{\rho,\e}}
\varphi(x)\frac{\partial \Gamma_\rho}{\partial\nu}(x)\,d{\mathcal{H}}^{n-1}_x
=0.$$
Hence, by~\eqref{OK-plkdm-1}, \eqref{OK-plkdm-2}
and
the second Green's Identity (recall~\eqref{GRr2}),
\begin{equation}\label{IJSfabnxcstukcIDjmcopiougf087}
\begin{split}& \int_{\R^n}\Gamma_\rho(x)\, \Delta\varphi(x) \,dx=\lim_{\e\searrow0}
\int_{\R^n\setminus A_{\rho,\e}}\Gamma_\rho(x)\, \Delta\varphi(x) \,dx
\\&\qquad=\lim_{\e\searrow0}\left[
\int_{\R^n\setminus A_{\rho,\e}}\Delta\Gamma_\rho(x)\, \varphi(x) \,dx-
\int_{\partial A_{\rho,\e}}\left(\Gamma_\rho(x)\frac{\partial \varphi}{\partial\nu}(x)-
\varphi(x)\frac{\partial \Gamma_\rho}{\partial\nu}(x)\right)\,d{\mathcal{H}}^{n-1}_x\right]\\&\qquad=
\int_{\R^n }\Delta\Gamma_\rho(x)\, \varphi(x) \,dx.
\end{split}\end{equation}
Thus, using~\eqref{DEGAM},
$$  \int_{\R^n}\Gamma_\rho(x)\, \Delta\varphi(x) \,dx=
-\fint_{B_\rho} \varphi(x) \,dx.$$
The desired result now follows be sending~$\rho\searrow0$
(notice that~$|\Gamma_\rho|\le\Gamma$ which is locally integrable,
hence the Dominated Convergence Theorem can be utilized here).
\end{proof}

A variant\footnote{For completeness, we observe that Theorem~\ref{KScvbMS:0okrmt4ht-VAR}
provides an alternative approach towards the Mean Value Formula in Theorem~\ref{KAHAR}.
For instance, by using~\eqref{300G} with~$\varphi:=-1$ and~$\Omega:=B_r(x_0)$ we find that
\begin{eqnarray*}
1=-\int_{\partial B_r(x_0)} \frac{\partial \Gamma}{\partial\nu}(x-x_0)\,d{\mathcal{H}}^{n-1}_x
=-\left.\frac{\partial \Gamma}{\partial\nu}\right|_{\partial B_r}\,{\mathcal{H}}^{n-1}(\partial B_r).
\end{eqnarray*}
Therefore, if~$\varphi$ is harmonic in~$B_R(x_0)$ and~$r\in(0,R)$ it follows from~\eqref{300G} that
\begin{eqnarray*}
\varphi(x_0)&=&\int_{\partial B_r(x_0)}
\left(\Gamma (x-x_0)\frac{\partial \varphi}{\partial\nu}(x)-
\varphi(x)\frac{\partial \Gamma}{\partial\nu}(x-x_0)\right)\,d{\mathcal{H}}^{n-1}_x\\&=&
\Gamma\big|_{\partial B_r}
\int_{\partial B_r(x_0)}\frac{\partial \varphi}{\partial\nu}(x)\,d{\mathcal{H}}^{n-1}_x-\left.\frac{\partial \Gamma}{\partial\nu}\right|_{\partial B_r}\int_{\partial B_r(x_0)}\varphi(x)\,d{\mathcal{H}}^{n-1}_x\\&=&
\Gamma\big|_{\partial B_r}
\int_{ B_r(x_0)}\Delta\varphi(x)\,dx+\frac1{{\mathcal{H}}^{n-1}(\partial B_r)}
\int_{\partial B_r(x_0)}\varphi(x)\,d{\mathcal{H}}^{n-1}_x\\
&=&\fint_{\partial B_r(x_0)}\varphi(x)\,d{\mathcal{H}}^{n-1}_x,\end{eqnarray*}
which corresponds to Theorem~\ref{KAHAR}(ii).} of Theorem~\ref{KScvbMS:0okrmt4ht} goes as follows:

\begin{theorem}\label{KScvbMS:0okrmt4ht-VAR}
Let~$\Omega$ be a bounded open set in~$\R^n$
with~$C^1$ boundary. Let~$x_0\in\Omega$.
Then, for every~$\varphi\in C^2(\Omega)\cap C^1(\overline\Omega)$,
\begin{equation}\label{300G} \int_{\Omega}
\Gamma(x-x_0)\, \Delta\varphi(x) \,dx-
\int_{\partial\Omega}
\left(\Gamma (x-x_0)\frac{\partial \varphi}{\partial\nu}(x)-
\varphi(x)\frac{\partial \Gamma}{\partial\nu}(x-x_0)\right)\,d{\mathcal{H}}^{n-1}_x
=-\varphi(x_0).\end{equation}
\end{theorem}

\begin{proof} Up to replacing~$\varphi(x)$ with~$\widetilde\varphi(x):=\varphi(x+x_0)$,
we can assume that~$x_0=0$.
We take~$\rho>0$ so small such that~$B_\rho\Subset\Omega$,
and thus~$\Gamma_\rho=\Gamma$ in a neighborhood of~$\partial\Omega$.
Then, we replace~\eqref{IJSfabnxcstukcIDjmcopiougf087}
in this framework by
\begin{equation*}
\begin{split}& \int_{\Omega}
\Gamma_\rho(x)\, \Delta\varphi(x) \,dx=\lim_{\e\searrow0}
\int_{\Omega\setminus A_{\rho,\e}}\Gamma_\rho(x)\, \Delta\varphi(x) \,dx
\\&\qquad=\lim_{\e\searrow0}\left[
\int_{\Omega\setminus A_{\rho,\e}}\Delta\Gamma_\rho(x)\, \varphi(x) \,dx+
\int_{\partial\Omega}\left(\Gamma_\rho(x)\frac{\partial \varphi}{\partial\nu}(x)-
\varphi(x)\frac{\partial \Gamma_\rho}{\partial\nu}(x)\right)\,d{\mathcal{H}}^{n-1}_x
\right.\\&\qquad\qquad\qquad\left.-
\int_{\partial A_{\rho,\e}}\left(\Gamma_\rho(x)\frac{\partial \varphi}{\partial\nu}(x)-
\varphi(x)\frac{\partial \Gamma_\rho}{\partial\nu}(x)\right)\,d{\mathcal{H}}^{n-1}_x\right]\\&\qquad=
\int_{\Omega }\Delta\Gamma_\rho(x)\, \varphi(x) \,dx+
\int_{\partial \Omega}\left(\Gamma(x)\frac{\partial \varphi}{\partial\nu}(x)-
\varphi(x)\frac{\partial \Gamma }{\partial\nu}(x)\right)\,d{\mathcal{H}}^{n-1}_x
\end{split}\end{equation*}
and we conclude as in the proof of Theorem~\ref{KScvbMS:0okrmt4ht}
by sending~$\rho\searrow0$.
\end{proof}

The identity~\eqref{300G}
is sometimes called ``Green's Representation Formula''\index{Green's Representation Formula}: interestingly,
it allows one to reconstruct the pointwise value of a function from
its Laplacian in a given domain and its value and the values
of its normal derivative at the boundary of the domain.\medskip

The fundamental solution
can be used to construct regular solutions of the Poisson equation
with sufficiently regular right hand side. As an example, we provide the
following result (a more precise version will follow from the 
Schauder estimates in Proposition~\ref{SI:IN:SE:DR}).

\begin{proposition}\label{ESI:NICE:FB}
Let~$f\in C^1_0(\R^n)$. Then, the function~$v:=-\Gamma*f$ belongs to~$C^2(\R^n)$
and satisfies~$\Delta v=f$ in~$\R^n$.
\end{proposition}

\begin{proof}
We notice that~$v$ is well-defined since~$\Gamma$ is locally integrable
and~${f}$ is bounded and with bounded support.
In fact, $v\in L^\infty(\R^n)$.
Also, by Theorem~\ref{KScvbMS:0okrmt4ht}, for every~$\varphi\in C^\infty_0(\R^n)$
and every~$x\in\R^n$ we know that
$$ 
\int_{\R^n}\Gamma(Y-x)\, \Delta\varphi(Y) \,dY=
\int_{\R^n}\Gamma(y)\, \Delta\varphi(x+y) \,dy=-\varphi(x)$$
and therefore
\begin{equation}\label{KA-OS-ikjfmf-S-kCONd}
\begin{split}& 
\int_{\R^{n}}v(Y)\,
\Delta\varphi(Y) \,dY=-
\iint_{\R^{2n}}\Gamma(Y-x)\,{f}(x)\,
\Delta\varphi(Y) \,dx\,dY=\int_{\R^n}\varphi(x)\,
{f}(x)\,dx.
\end{split}\end{equation}
We observe that if~$\phi\in L^\infty(\R^n)\cap L^1(\R^n)$ then, for all~$i\in\{1,\dots,n\}$,
\begin{equation}\label{PRIMA98DER}
\partial_i(\Gamma*\phi)=(\partial_i\Gamma)*\phi.
\end{equation}
To check this, given~$\rho>0$ we use the regularization~$\Gamma_\rho$ of the fundamental solution
introduced in~\eqref{COIENEE-0986ytufgkbv-0-rjfeonvnb2GDB} and~\eqref{GAROGRA}.
We define~$\psi:=\Gamma*\phi$ and~$\psi_\rho:=\Gamma_\rho*\phi$ and we observe that, for every~$x\in\R^n$,
\begin{equation}\label{IHFSHNTHSYDHMSWMDTKFDOFLL}
\begin{split}&
|\psi_\rho(x)-\psi(x)|\le
\frac{c_n}{2n|B_\rho|}\int_{B_\rho} (\rho^2-|y|^2)|\phi(y)|\,dy\\
&\qquad\le
\frac{C\,\|\phi\|_{L^\infty(\R^n)}}{\rho^n}
\int_0^\rho(\rho^2-r^2)\,r^{n-1}\,dr\le C\,\|\phi\|_{L^\infty(\R^n)}\,\rho^2
\end{split}\end{equation}
for some constant~$C>0$ depending only on~$n$ and possibly varying from line to line.

Additionally,
\begin{eqnarray*}
&&\left|\partial_i \psi_\rho(x)-(\partial_i\Gamma)*\phi(x)\right|=
\left|\partial_i\Gamma_\rho*\phi(x)-(\partial_i\Gamma)*\phi(x)\right|\\&&\qquad\le
\frac{C}{\rho^n}\int_{B_\rho}|y|\,|\phi(x-y)|\,dy+C\int_{B_\rho}|y|^{1-n}\,|\phi(x-y)|\,dy
\le C\,\|\phi\|_{L^\infty(\R^n)}\,\rho.
\end{eqnarray*}
{F}rom this and~\eqref{IHFSHNTHSYDHMSWMDTKFDOFLL} we infer that~$\psi_\rho$ converges uniformly to~$\psi$
and its derivative to~$(\partial_i\Gamma)*\phi$, whence it follows that~$\partial_i\psi=(\partial_i\Gamma)*\phi$,
thus establishing~\eqref{PRIMA98DER}.

We also remark that~$\partial_j(\Gamma*f)=\Gamma*(\partial_jf)$
(see e.g.~\cite[Theorem~9.3]{MR3381284}) and therefore it follows from~\eqref{PRIMA98DER} that~$
\partial_{ij}^2 v=-(\partial_i\Gamma)*(\partial_j f)$.
Since~$\partial_j f$ is continuous and compactly supported
and~$\partial_i\Gamma$ is locally integrable, we have that~$\partial_{ij}^2 v\in L^\infty(\R^n)$.
Moreover, we suppose that the support of~$f$ is contained in a bounded set~$\Omega$ and, for every~$x\in\R^n$,
we denote by~$\Omega_x$ the bounded set containing all the points~$y$ such that~$x-y$
belong to~$\displaystyle\bigcup_{p\in\Omega} B_1(p)$.
Then, if~$x\in\R^n$ and~$x_k\to x$ as~$k\to+\infty$, we have that
$$ |\partial_i\Gamma(y)\partial_j f(x_k-y)|\le\|{f}\|_{C^1(\R^n)}\,|\nabla\Gamma(y)|\,\chi_{\Omega_x}(y),$$
and the latter is an integrable function of~$y\in\R^n$.
As a result, by the Dominated Convergence Theorem,
$$ \lim_{k\to+\infty}\partial_{ij}^2v(x_k)=
-\lim_{k\to+\infty}\int_{\R^n}\partial_i\Gamma(y)\partial_j f(x_k-y)\,dy
=\int_{\R^n}\partial_i\Gamma(y)\partial_j f(x-y)\,dy=\partial_{ij}^2v(x),$$
whence~$
v\in C^2(\R^n)$.
This and~\eqref{KA-OS-ikjfmf-S-kCONd} yield that~$\Delta v=f$.
\end{proof}

The function~$-\Gamma*f$ in Proposition~\ref{ESI:NICE:FB}
(or sometimes, up to a sign convention, the function~$\Gamma*f$)
is often referred to with the name of Newtonian potential\index{Newtonian potential}.\medskip

A neat intuition for the result in Proposition~\ref{ESI:NICE:FB} comes from
physical motivations: namely, since the fundamental solution~$-\Gamma(x-y)$ (disregarding physical constants and possible sign conventions) corresponds to the
gravity\index{gravitation} field at the point~$x$ generated by a pointwise mass located at the point~$y$, it follows that
the quantity~$-\Gamma*f(x)$ (which agrees with
the the Newtonian potential\index{Newtonian potential})
corresponds to the gravitational field at the point~$x$ generated by a distribution~$f$ of masses, since it can be
seen as the superposition of~$\Gamma(x-y)f(y)\,dy$. By
Gau{\ss}' Law, the flux of the field through a given surface~$\partial\Omega$ (seen as the boundary
of a bounded set~$\Omega$) corresponds to the
total mass comprised inside~$\partial\Omega$, namely
$$ \int_\Omega f(x)\,dx=-\int_{\partial\Omega} \nabla(\Gamma*f)(x)\cdot\nu(x)\,d{\mathcal{H}}^{n-1}_x.$$
Using the Divergence Theorem on the latest integral we thereby find that
$$ \int_\Omega f(x)\,dx=-\int_{\Omega}\div\Big( \nabla(\Gamma*f)(x)\Big)\,dx
=-\int_\Omega \Delta(\Gamma*f(x))\,dx.$$
Hence, since~$\Omega$ is an arbitrary domain, we obtain that~$f=-\Delta(\Gamma*f)$,
which is precisely the content of Proposition~\ref{ESI:NICE:FB}.
In this sense, Proposition~\ref{ESI:NICE:FB} is a mathematically structured statement
of Gau{\ss}' Law (with some care devoted on the regularity of the distribution of mass~$f$
required for the result to hold true).
\medskip

As a byproduct of our analysis, one can also deal with the gravitational potential of a homogeneous ball and show that at all external points this potential is equal
to the potential of the material point of the same mass placed in its center (see~\eqref{OJHSN-IOU0}
for an application to physical geodesy):

\begin{lemma}\label{QUESTOL}
Let~$R>0$ and~$x\in\R^n\setminus \overline{B_R}$. Then,$$ \int_{B_R} \Gamma(x-y)\,dy=|B_R|\,\Gamma(x).$$\end{lemma}

\begin{proof} Given~$y\in B_R$, we have that~$\Delta\Gamma(x-y)=0$ for every~$x\in\R^n\setminus B_R$. In particular,
by the Green's Identity~\eqref{GRr2},
if~$r\in(0,R)$ and~$\psi(y):=\frac{|y|^2-r^2}{2n}$,\begin{eqnarray*}&&\int_{B_r}
\Gamma(x-y)\,dy=\int_{B_r} \Gamma(x-y)\Delta \psi(y)\,dy\\&&
\qquad=\int_{B_r}\psi(y)\Delta \Gamma(x-y) \,dy+
\int_{\partial B_r}\left(\Gamma(x-y)\frac{\partial \psi}{\partial\nu}(y)-
\psi(y)\frac{\partial \Gamma}{\partial\nu}(x-y)\right)\,d{\mathcal{H}}^{n-1}_y
\\&&\qquad=0+\int_{\partial B_r}\left(\frac{r}n\,\Gamma(x-y)-0\right)\,d{\mathcal{H}}^{n-1}_y=\frac{r}n\,\int_{\partial B_r} \Gamma(x-y)\,d{\mathcal{H}}^{n-1}_y
.\end{eqnarray*}As a result, defining$$ \Psi(r):=r^{-n}\int_{B_r}\Gamma(x-y)\,dy,$$and using
polar coordinates (see e.g.~\cite[Theorem~3.12]{MR3409135}), for~$r\in(0,R)$ we have
\begin{eqnarray*}\Psi'(r)=-nr^{-n-1}\int_{B_r}\Gamma(x-y)\,dy
+
r^{-n}\int_{\partial B_r}\Gamma(x-y)\,d{\mathcal{H}}^{n-1}_y
=0.\end{eqnarray*}
This gives that~$\Psi$ is constant in~$(0,R]$ and thus\begin{equation*} \Gamma(x)=\lim_{r\searrow0}\fint_{B_r}\Gamma(x-y)\,dy=\lim_{r\searrow0}\frac{\Psi(r)}{|B_1|}=\lim_{r\nearrow R}\frac{\Psi(r)}{|B_1|}=
\fint_{B_R}\Gamma(x-y)\,dy.\qedhere\end{equation*}
\end{proof}

\section{Back to the Kelvin Transform:
the method of image
charges}\label{LOO}

We discuss here a simple, but very influential, technique from electrostatics\index{electrostatics}
that naturally leads to the Kelvin Transform introduced in Section~\ref{LLOO2}.\index{method of image charges}
This method uses the fundamental solution
showcased in Section~\ref{lfundsP-S}, interpreting it as the electrostatic
potential generated by a point charge (the method can also be
considered as an inspiration for the construction
of the Green Function
of the ball in the forthcoming Theorem~\ref{MS:SKMD344567yjghS-1}).
\medskip

The details of this motivation go as follows.
Suppose that we have a positive unit point charge located at some point
in the ball~$B_1$. Is it possible to place a negative point charge (not necessarily a
unit charge) somewhere outside of the ball in order to make~$\partial B_1$
a surface with constant, say zero, potential? If so, can we determine the position and
intensity of this
auxiliary charge?\medskip

For the mathematical setting for this situation we suppose\footnote{The case~$n=2$
requires some conceptual modification: \label{CONSDcenrfPTHDNULAMDALSIKLASKGIMDOVAS}
when~$n=2$, rather than imposing that~$\partial B_1$
is a surface with zero potential, one requires the weaker condition
that it is an equipotential surface, and the value of the potential
may depend on~$x_0$. To compensate
this weaker potential condition however, in dimension~$2$
one can additionally impose that the image charge has
also unit intensity (but it is negatively charged), that
is, up to a sign, both the original charge and the mirror charge have
the same intensity.
In this way, equation~\eqref{hifbdnHDNDNbdsfidbiwtg49673fggMS}
when~$n=2$ is replaced by
$$ \ln|x-x_0|-
\ln |x-T(x_0)|=\beta(x_0)\qquad{\mbox{for every }}x\in\partial B_1,$$
for a suitable~$\beta(x_0)$.
That is, setting~$\gamma(x_0):=e^{\beta(x_0)}$,
$$ 1+|x_0|^2-2x\cdot x_0=|x-x_0|^2=(\gamma(x_0))^2 |x-T(x_0)|^2=
(\gamma(x_0))^2(1+|T(x_0)|^2-2x\cdot T(x_0)).
$$
This identity is the counterpart of~\eqref{2345tyKS:KJHGFDOIUYTROIHGFD0987654}
when~$n=2$ (with the notation~$\mu(x_0)=\frac{1}{(\gamma(x_0))^2}$),
hence the computation provided in these pages (and leading to~\eqref{2345tyKS:KJHGFDOIUYTROIHGFD0987654BIS})
would give~$\gamma(x_0)=|x_0|$ and hence~$\beta(x_0)=\ln|x_0|$.}
that~$n\ne2$. In this way, up to a normalizing constant,
we can suppose that the electrostatic potential at some point~$x\in\R^n$
generated by
a positive unit charge located at~$x_0\in B_1$ is equal to~$\frac1{|x-x_0|^{n-2}}$,
recall~\eqref{GAMMAFU}.
So, suppose that we place a charge of some intensity~$-\alpha(x_0)\in(-\infty,0)$
at some point~$T(x_0)\in\R^n\setminus B_1$. The potential
generated by this auxiliary charge is equal to~$-\frac{\alpha(x_0)}{|x-T(x_0)|^{n-2}}$.

Therefore, the condition that~$\partial B_1$ is a surface with zero potential
boils down to the relation
\begin{equation}\label{hifbdnHDNDNbdsfidbiwtg49673fggMS}
\frac1{|x-x_0|^{n-2}}-
\frac{\alpha(x_0)}{|x-T(x_0)|^{n-2}}=0\qquad{\mbox{for every }}x\in\partial B_1.
\end{equation}
Setting~$\mu(x_0):=(\alpha(x_0))^{\frac2{n-2}}$, we therefore obtain that,
for every~$x\in\partial B_1$,
\begin{equation}\label{2345tyKS:KJHGFDOIUYTROIHGFD0987654}\begin{split}&
1+|T(x_0)|^2-2x\cdot T(x_0)
=|x|^2+|T(x_0)|^2-2x\cdot T(x_0)=
|x-T(x_0)|^2=(\alpha(x_0))^{\frac2{n-2}}|x-x_0|^2\\&\qquad=
\mu(x_0)\,|x-x_0|^2=\mu(x)\,\big(|x|+|x_0|^2-2x\cdot x_0\big)=
\mu(x_0)\,\big(1+|x_0|^2-2x\cdot x_0\big),\end{split}
\end{equation}
that is
\begin{equation}\label{LAN DcindfgbfdizPM-10293tu4y}
1+|T(x_0)|^2-\mu(x_0)\,\big(1+|x_0|^2\big)
=2x\cdot \big( T(x_0)-
\mu(x_0)\, x_0\big)
\qquad{\mbox{for every }}x\in\partial B_1.
\end{equation}
Now, given~$x\in\partial B_1$, we exploit~\eqref{LAN DcindfgbfdizPM-10293tu4y}
both for~$x$ and for~$-x\in\partial B_1$, finding that
$$ 2x\cdot \big( T(x_0)-
\mu(x_0)\, x_0\big)=1+|T(x_0)|^2-\mu(x_0)\,\big(1+|x_0|^2\big)=
-2x\cdot \big( T(x_0)-
\mu(x_0)\, x_0\big).$$
Consequently,
\begin{equation}\label{IDKMsubiadkfgbK9idfsaj9irtkgki0d-1}
x\cdot \big( T(x_0)-
\mu(x_0)\, x_0\big)=0
\qquad{\mbox{for every }}x\in\partial B_1
\end{equation}
and thus~\eqref{LAN DcindfgbfdizPM-10293tu4y} reduces to
\begin{equation}\label{IDKMsubiadkfgbK9idfsaj9irtkgki0d-2}
1+|T(x_0)|^2-\mu(x_0)\,\big(1+|x_0|^2\big)
=0.
\end{equation}
As a matter of fact, by choosing~$x$ as
an element of the Euclidean basis in~\eqref{IDKMsubiadkfgbK9idfsaj9irtkgki0d-1}
we obtain that
\begin{equation}\label{KS:rujmaldmcgab dtedyfwiosxicaqs3d}
T(x_0)\cdot e_i=
\mu(x_0)\, x_0\cdot e_i\qquad{\mbox{for every }}i\in\{1,\dots,n\}
\end{equation}
and therefore
$$ |T(x_0)|^2=\sum_{i=1}^n(T(x_0)\cdot e_i)^2
=(\mu(x_0))^2\,\sum_{i=1}^n( x_0\cdot e_i)^2
=(\mu(x_0))^2\,|x_0|^2.$$
Substituting this information into~\eqref{IDKMsubiadkfgbK9idfsaj9irtkgki0d-2}
we find that
$$ 0=1+(\mu(x_0))^2\,|x_0|^2-\mu(x_0)\,\big(1+|x_0|^2\big)$$
and accordingly either~$\mu(x_0)=1$ or~$\mu(x_0)=\frac1{|x_0|^2}$.
However, the possibility~$\mu(x_0)=1$ must be ruled out,
otherwise, by~\eqref{KS:rujmaldmcgab dtedyfwiosxicaqs3d}, we would
have that~$\R^n\setminus B_1\ni T(x_0)=x_0\in B_1$ which is
a contradiction.

In this way, we have establishes that
\begin{equation}\label{2345tyKS:KJHGFDOIUYTROIHGFD0987654BIS}
\mu(x_0)=\frac1{|x_0|^2},\end{equation}
and thus necessarily~$x_0\ne0$ if we want our problem to have a solution,
whence
\begin{equation}\label{SOLLE:1}
\alpha(x_0)=\frac1{|x_0|^{n-2}}.\end{equation}
Furthermore, we obtain from~\eqref{KS:rujmaldmcgab dtedyfwiosxicaqs3d} that
\begin{equation}\label{SOLLE:2}
T(x_0)=\frac{x_0}{|x_0|^2}.
\end{equation}
Thus, conditions~\eqref{SOLLE:1}
and~\eqref{SOLLE:2} provide the solution of our
image electrostatic charge problem, specifying, respectively,
the intensity of the image charge and its spatial location.

Remarkably, the spatial location of the image charge in~\eqref{SOLLE:2}
coincides with the Kelvin Transform in~\eqref{KSM:J:KELVIN1}
and the charge intensity in~\eqref{SOLLE:1} coincides with
the multiplicative factor in the functional action of the 
Kelvin Transform as defined in~\eqref{SOLLE:0}.

\section{Maximum Principles}\label{MAXPLEDISCERET:SE}

Maximum Principles are one of the cornerstones\index{Maximum Principle}
of the theory of elliptic partial differential equations
and they can also
be seen as the real analysis counterpart of the Maximum Modulus Principle
which is in turn one of the backbones of complex analysis (see e.g.~\cite{MR0210528}).
Roughly speaking, the main idea is that if~$\Delta u=0$ then
each point of the graph of~$u$ is necessarily a ``saddle point''
and therefore~$u$ cannot have local maxima or minima.
In these pages, we will obtain the simplest
possible versions of the Maximum Principle for the Laplace operator.
First of all, we present
a statement, often called
Strong Maximum Principle, that goes as follows:

\begin{theorem}\label{STRONGMAXPLE1}
Let~$\Omega\subseteq\R^n$ be open and connected,
and let~$u\in C^2(\Omega)$.

\begin{itemize}
\item[(i).] If~$\Delta u\ge0$ in~$\Omega$
and there exists~$\overline{x}\in\Omega$
such that~$u(\overline{x})=\sup_\Omega u$, then u is constant.
\item[(ii).] If~$\Delta u\le0$ 
in~$\Omega$ and there exists~$\underline{x}\in\Omega$
such that~$u(\underline{x})=\inf_\Omega u$, then u is constant.
\item[(iii).] If~$u$ is harmonic
in~$\Omega$, then it cannot
attain an interior maximum or minimum value unless it is constant.\end{itemize}
\end{theorem}

\begin{proof} We observe that the claim in~(ii) follows from~(i)
by changing~$u$ into~$-u$. Also, the claim in~(iii) follows from~(i) and~(ii).
Therefore, it suffices to prove~(i). To this end,
we define
\begin{equation}\label{S:osl4995ffdssSloamed} {\mathcal{U}}:=\left\{x\in\Omega{\mbox{ s.t. }}u(x)=\sup_\Omega u\right\}
.\end{equation}
We remark that~$\overline{x}\in{\mathcal{U}}$
and thus~${\mathcal{U}}\ne\varnothing$. The continuity of~$u$
also gives that~${\mathcal{U}}$ is closed in~$\Omega$.

\begin{figure}
  \centering
  \includegraphics[width=.3\linewidth]{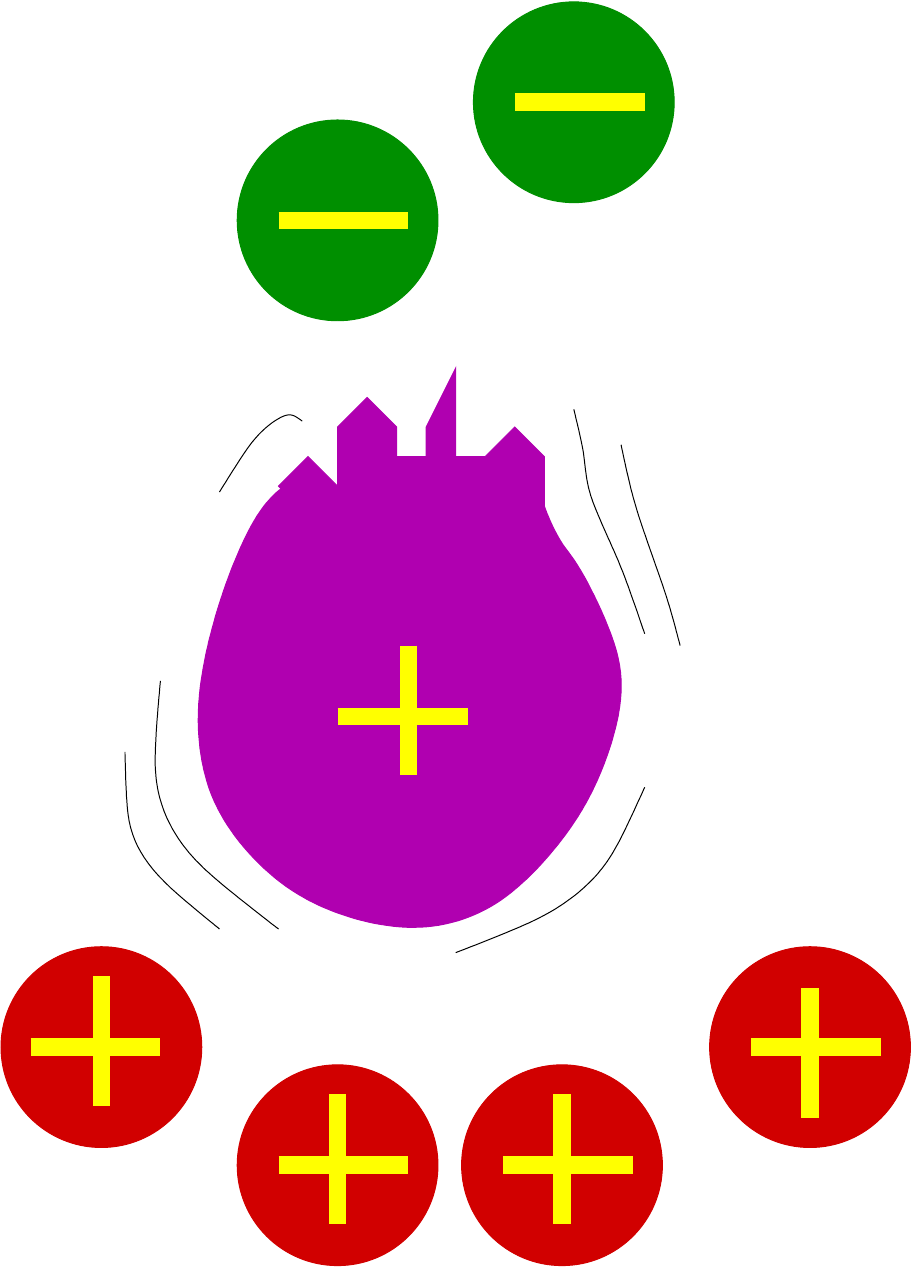}
 \caption{\sl Electromagnetic levitation, after 
 Ren\'e Magritte's The Castle of the Pyrenees.}\label{2354HA4MAGR66RH788AR9C0O0U1N2TERELAGIJ7soloDItangeFI}
\end{figure}

Moreover, if~$r>0$ is such that~$B_r(\overline{x})\Subset\Omega$,
we let~$H(x):=\Gamma(x-\overline{x})-\Gamma(re_1)$
and we deduce from Green's Representation Formula~\eqref{300G}
and the Divergence Theorem (recall also the setting of the fundamental
solution in~\eqref{COIENEE}
and~\eqref{GAMMAFU}, as well as the measure theoretic
identity in~\eqref{B1})
that
\begin{equation}\label{0987654hdhdnknn098765hjoiuytre-0987654adsdoodo}
\begin{split}&
0\le\int_{B_r(\overline{x})}
H(x)\,\Delta u(x)\,dx=
\int_{B_r(\overline{x})}
\Gamma(x-\overline{x})\,\Delta u(x)\,dx-\Gamma(re_1)
\int_{B_r(\overline{x})}\Delta u(x)\,dx
\\&\qquad
=
\int_{\partial B_r(\overline{x})}
\left(\Gamma (x-\overline{x})\frac{\partial u}{\partial\nu}(x)-
u(x)\frac{\partial \Gamma}{\partial\nu}(x-\overline{x})\right)\,d{\mathcal{H}}^{n-1}_x
\\&\qquad\qquad-u(\overline{x})-\Gamma(re_1)
\int_{\partial B_r(\overline{x})}\frac{\partial u}{\partial\nu}(x)\,d{\mathcal{H}}^{n-1}_x
\\&\qquad
=\frac{r}{n\,|B_1|}\,\int_{\partial B_r(\overline{x})}
\frac{u(x)}{|
x-\overline{x}|^n}\,d{\mathcal{H}}^{n-1}_x
-\sup_\Omega u\\&\qquad= \fint_{\partial B_r(\overline{x})}
u(x)\,d{\mathcal{H}}^{n-1}_x
-\sup_\Omega u.
\end{split}
\end{equation}
Now we take~$\rho>0$ such that~$B_\rho(\overline{x})\Subset\Omega$
and we claim that
\begin{equation}\label{TGASBDCryuiTTskdfroppnuo}
{\mbox{for every~$p\in B_\rho(\overline{x})$ it holds that~$u(p)=\,$}}\sup_\Omega u.
\end{equation}
Indeed, suppose not. Then, we have that
$$ \fint_{B_\rho(\overline{x})}u(x)\,dx<\sup_\Omega u.$$
Thus, using polar coordinates and~\eqref{0987654hdhdnknn098765hjoiuytre-0987654adsdoodo},
\begin{eqnarray*}
&&\sup_\Omega u>\frac{1}{| B_\rho|}\,
\int_0^\rho\left(
\int_{\partial B_r(\overline{x})}u(x)\,d{\mathcal{H}}^{n-1}_x
\right)\,dr\ge
\frac{{\mathcal{H}}^{n-1}(\partial B_1)}{|B_1|\rho^{n}}\,
\int_0^\rho\left(r^{n-1}
\sup_\Omega u
\right)\,dr\\&&\qquad\qquad=
\frac{{\mathcal{H}}^{n-1}(\partial B_1)}{n\,|B_1|}\,\sup_\Omega u
=\sup_\Omega u.
\end{eqnarray*}
This contradiction proves~\eqref{TGASBDCryuiTTskdfroppnuo}.

As a consequence
of~\eqref{TGASBDCryuiTTskdfroppnuo}, we have that~${\mathcal{U}}$
is open. Then, by the connectedness of~$\Omega$, it follows that~${\mathcal{U}}=\Omega$
and consequently~$u(x)=\sup_\Omega u$ for all~$x\in\Omega$.
\end{proof}

A more general version of this result, without assuming the function~$u$
to be smooth, will be given in Lemma~\ref{LplrfeELLtahsmedschdfb4camntM}.\medskip

\begin{figure}
                \centering
                \includegraphics[width=.48\linewidth]{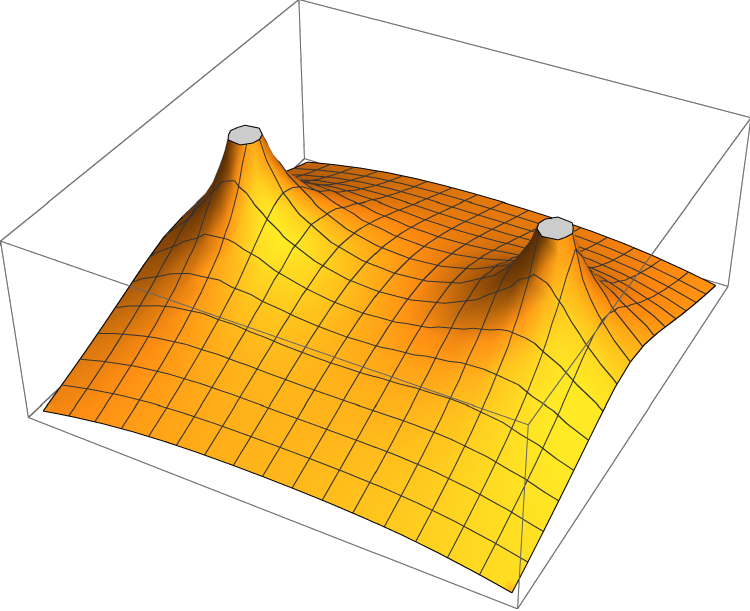}
        \caption{\sl Electrostatic potential generated by two point charges in the plane located at~$e_1$ and~$-e_1$,
        corresponding to~$-\frac{1}{2\pi}\left( \ln|x-e_1|+\ln|x+e_1|\right)$. Note that the origin is a saddle point.}\label{2ed2HAFOTkAKA123JE2qasD231SEL}
\end{figure}

\begin{figure}
                \centering
                \includegraphics[width=.65\linewidth]{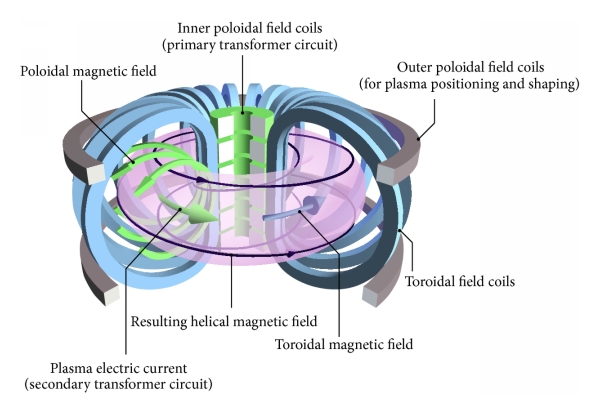}
        \caption{\sl Schematic of a tokamak chamber
        (image from Wikipedia, taken from~\cite{MR3224327},
        licensed under the Creative Commons Attribution-Share Alike 4.0 International license).}\label{2ed2HAFOTkAKA123JED231}
\end{figure}

An interesting consequence of the Strong Maximum Principle in Theorem~\ref{STRONGMAXPLE1}
is the observation, originally due to \index{Earnshaw's Theorem}
Samuel Earnshaw, according to which\footnote{Of course, if we put a positive charge at~$e_1$ and another at~$-e_1$,
a positive charge at the origin would be in equilibrium (by symmetry), but this would be an unstable equilibrium,
since a small perturbation of the position of the charge at the origin in the~$e_2$ direction would make it drift away.

Given a collection of point charges in~$\R^n$, with~$n\ge2$, the corresponding harmonic electrostatic potential~$u$
would generate a force field~$F=-\nabla u$
(recall the discussion on electromagnetic fields in Section~\ref{SECT:Gravitation and electrostatics}). The equilibria correspond to the zeros of~$F$,
and the notion of stability adopted here
means that if we move the particle a little bit away from the equilibrium the electrostatic force will tend to bring it back to the equilibrium.

Hence, if, say, the origin were a stable equilibrium
in the above setting, then~$F(0)=0$ and, if~$x$ is close enough to the origin,
then~$F(x)\cdot x\le0$; that is, $\nabla u(0)=0$ and, if~$x$ is close enough to the origin,~$\nabla u(x)\cdot x\ge0$.
In particular, there exists~$\rho>0$ such that if~$x\in B_\rho\setminus\{0\}$ and~$\omega:=\frac{x}{|x|}$ then
$$ u(x)-u(0)=u(|x|\omega)-u(0)=\int_0^{|x|} \frac{d}{dr} u(r\omega)\,dr=
\int_0^{|x|} \nabla u(r\omega)\cdot\omega\,dr=\frac{1}{r}\,
\int_0^{|x|} \nabla u(r\omega)\cdot(r\omega)\,dr\ge0.$$
In particular, a stable equilibrium at the origin would produce a local minimum for~$u$,
which would violate Theorem~\ref{STRONGMAXPLE1}(iii).
This proves Earnshaw's Theorem in~\eqref{1Earnshaw}.

Actually, Theorem~\ref{STRONGMAXPLE1}(iii) gives, more generally, that
the electrostatic potential cannot
have a maximum or a minimum value at any point in space not occupied
by an electric charge
(saddle points are instead possible and would correspond to unstable
equilibria, as discussed above, see Figure~\ref{2ed2HAFOTkAKA123JE2qasD231SEL}).

We stress that the notion of stability used here is the one formulated
by Maxwell in~\cite{MR0063293},
according to which an equilibrium is stable if the field at points close to the equilibrium pushes back
toward the equilibrium itself. This notion of stability is not the same as
the Lyapunov stability (but a version of Earnshaw's Theorem also holds
in terms of Lyapunov stability, see~\cite[Theorem~1.2]{MR3445477}).

An interesting consequence of Earnshaw's Theorem is that
it is not possible to achieve a stable static levitation by only using a combination of fixed magnets and electric charges, see Figure~\ref{2354HA4MAGR66RH788AR9C0O0U1N2TERELAGIJ7soloDItangeFI}
(and note that including a constant gravitation field would just modify the potential by an additional
harmonic term~$gx_n$, being~$g$ the constant gravity acceleration).
Also, if one wants to trap charged particles by only using
an electromagnetic field, necessarily the particles have
to be dynamically kept in motion (no stable rest equilibrium being allowed),
which is one of the inspiring principles for the construction
of a tokamak, see Figure~\ref{2ed2HAFOTkAKA123JED231}
(and recall Section~\ref{PLASECT-0uorwehfg} for some basic notions
of plasma models).

See~\cite{FALLA}
and the references therein for more information on Earnshaw's Theorem.}
\begin{equation}\label{1Earnshaw}\begin{split}&
{\mbox{a collection of point charges cannot
be maintained}}\\ &{\mbox{in a stable stationary electrostatic equilibrium.}}\end{split}\end{equation}

\medskip

As a byproduct of
the Strong Maximum Principle in Theorem~\ref{STRONGMAXPLE1},
we have a Weak Maximum Principle, as follows:\index{Maximum Principle}

\begin{corollary}\label{WEAKMAXPLE}
Let~$\Omega\subseteq\R^n$ be open, bounded and connected,
and let~$u\in C^2(\Omega)\cap C(\overline\Omega)$.

\begin{itemize}
\item[(i).] If~$\Delta u\ge0$ in~$\Omega$,
then~$\sup_\Omega u=\sup_{\partial\Omega} u$.
\item[(ii).] If~$\Delta u\le0$ in~$\Omega$,
then~$\inf_\Omega u=\inf_{\partial\Omega} u$.
\item[(iii).] If~$u$ is harmonic
in~$\Omega$, then, for every~$x\in\Omega$,
$$\inf_{\partial\Omega}u\le u(x)\le\sup_{\partial\Omega} u.$$
\end{itemize}
\end{corollary}

\begin{proof} Since claim~(ii) follows from~(i) up to changing a sign in~$u$
and claim~(iii) is a consequence of~(i) and~(ii), we focus on the proof of~(i).
To this end, we consider a point~$p\in\overline{\Omega}$ such that~$u(p)=\max_{\overline\Omega}u$.
If~$p\in\partial\Omega$, then the claim in~(i) plainly follows.
Hence we can assume that~$p\in\Omega$. In this way, we have
that~$\sup_{ \Omega} u=u(p)$
and accordingly, in light of
Theorem~\ref{STRONGMAXPLE1}, we have that~$u$ is constant in~$\Omega$,
and this gives claim~(i) as desired.
\end{proof}

A sharper version of Corollary~\ref{WEAKMAXPLE}
not relying on smoothness assumptions on~$u$
will be presented later on in Corollary~\ref{PERSUBAWEAK}.

A useful uniqueness result for follows straightforwardly
from 
the Weak Maximum Principle in Corollary~\ref{WEAKMAXPLE}:\index{Maximum Principle}

\begin{corollary}\label{UNIQUENESSTHEOREM}
Let~$\Omega\subseteq\R^n$ be open, bounded and connected,
and let~$u$, $v\in C^2(\Omega)\cap C(\overline\Omega)$.

Assume that~$\Delta u=\Delta v$ in~$\Omega$ and that~$u=v$
on~$\partial\Omega$. Then, $u(x)=v(x)$ for every~$x\in\Omega$.
\end{corollary}

\begin{proof} The function~$w:=u-v$ is harmonic in~$\Omega$
and vanishes along~$\partial\Omega$. The desired result
then plainly follows from
Corollary~\ref{WEAKMAXPLE}(iii).
\end{proof}

In general, the Maximum Principle fails
in the case of unbounded domains. For instance,
for all~$c_1,c_2,c_3\in\R$,
the function~$u(x):=
c_1x_n+c_2 x_{n-1}x_n+c_3(x_n^3-3x_{n-1}^2x_n )$
is harmonic in the halfspace~$\{x_n>0\}$ and vanishes along~$\{x_n=0\}$,
without being identically zero (unless~$c_1=c_2=c_3=0$). Similarly, for all~$k\in\N$, $k\ge1$,
the function~$u(x):=e^{kx_{n-1}}\sin(kx_n)$
is harmonic in the slab~$\{x_n\in(0,\pi)\}$
and vanishes on
the boundary of the slab, without being identically zero.
In addition, in the external domain~$\R^n\setminus B_1$,
the function~$u(x):=\Gamma(e_1)-\Gamma(x)$ is harmonic
with~$u=0$ on~$\partial(\R^n\setminus B_1)$, but~$u$ is not
identically zero (here~$\Gamma$ is the
fundamental solution presented in Section~\ref{lfundsP-S}).

However, there are suitable assumptions
on the solution and on the domain that guarantee the validity
of a suitable Maximum Principle also for unbounded regions,
see e.g.~\cite{MR1258192, MR1329831, MR1395408, MR1892934}
and the references therein
(also, some of the Maximum Principles
for unbounded domains can be related to
the so-called {P}hragm\'{e}n-{L}indel\"{o}f Principle,
see e.g.~\cite{MR1501888}). Here we just recall a prototypical
situation:\index{Maximum Principle}

\begin{proposition} \label{UNDIFKFCOKDMTFDEBSF}
Let~$\Omega$ be a halfspace.
Let~$u\in C^2(\Omega)$ be uniformly continuous, bounded from above
and harmonic
in~$\Omega$, with~$u\le0$ on~$\partial\Omega$. 

Then, $u\le0$ in~$\Omega$.
\end{proposition}

\begin{proof} Up to a rigid motion, we can suppose that
\begin{equation}\label{WEUISDNindvDTifbIMPKMS}
\Omega=
\{x_n>0\}.\end{equation} Given~$\e>0$, we define~$v(x):=u(x)-\e x_n$.
We observe that~$v\le u$ in~$\overline\Omega$, thanks to~\eqref{WEUISDNindvDTifbIMPKMS},
and therefore
\begin{equation}\label{WEUISDNindvDTifbIMPKMS2}
v\le0 {\mbox{ on }}\partial\Omega=\{x_n=0\}.
\end{equation}
Now, we take~$\sigma$ to be the supremum of~$v$.
Let also~$p_k\in\Omega$ be a sequence such that~$v(p_k)\to\sigma$
as~$k\to+\infty$. We write~$p_k=(p_k',p_{k,n})\in\R^{n-1}\times[0,+\infty)$
and, up to disregarding a finite number of indices, we can assume that~$
v(p_k)\ge \sigma-1$.

Then, we have that~$\sigma-1\le u(p_k)-\e p_{k,n}\le \sup_\Omega u-\e p_{k,n}$,
from which it follows that~$p_{k,n}$ is a bounded sequence.
As a result, up to a subsequence, we assume that~$p_{k,n}\to q$
for some~$q\in[0,+\infty)$ as~$k\to+\infty$.

Thus, we define~$v_k(x):=v(x'+p_k', x_n)$.
Additionally, $v_k$ is uniformly equicontinuous, thanks to our uniform continuity assumption on~$u$.
Moreover, 
\begin{eqnarray*}
&&|v_k(0,q)|\le |v_k(0,p_{k,n})-v_k(0,q)|+ |v_k(0,p_{k,n})|\\&&\qquad\le|u(p_k', p_{k,n})-u(p_k',q)|+\e|p_{k,n}-q|+|v(p_k)|
\le 1+\varepsilon +\sigma+1,
\end{eqnarray*}
for large~$k$, which, together with the equicontinuity, gives that~$v_k$ is
locally equibounded. Consequently, by the Arzel\`a-Ascoli Theorem,
up to a subsequence we can assume 
that~$v_k$ locally uniformly
converges to some~$v_\infty$ as~$k\to+\infty$.

As a consequence,
\begin{eqnarray*}&& \sigma\ge v_\infty(0,q)=\lim_{k\to+\infty} v_k(0,q)\ge
\lim_{k\to+\infty}\big( v_k(0,p_{k,n}) - | v_k(0,q)-v_k(0,p_{k,n})|\big)\\&&\qquad\qquad\qquad\qquad
= \lim_{k\to+\infty} v(p_k',p_{k,n}) -0
=\sigma,\end{eqnarray*}
that is
\begin{equation}\label{eokfmgthasmdaccKMD:5678}
v_\infty(0,q)=\max_{\overline\Omega} v_\infty=\sigma.\end{equation}
Furthermore, by~\eqref{WEUISDNindvDTifbIMPKMS2},
for all~$x'\in\R^{n-1}$,
\begin{equation}\label{dfgheokfmgthasmdaccKMD:5678BIS}
v_\infty(x',0)=\lim_{k\to+\infty}v(x'+p_k', 0)
\le0.\end{equation}
Moreover, for every~$\varphi\in C^\infty_0(\Omega)$,
$$ \int_\Omega v_\infty(x)\,\Delta\varphi(x)\,dx
=\lim_{k\to+\infty}\int_\Omega v_k(x)\,\Delta\varphi(x)\,dx
=0,$$
and thus~$v_\infty$ is harmonic in~$\Omega$,
due to the Weyl's Lemma in~\eqref{WEYL}.

{F}rom this, \eqref{eokfmgthasmdaccKMD:5678}
and the Maximum Principle in
Theorem~\ref{STRONGMAXPLE1}(iii),
it follows that
either
\begin{equation}\label{MS:POOSK-BILASNSTPAMLOS3-01}
v_\infty(x)=\sigma {\mbox{ for all }}x\in\Omega
\end{equation}
or~$(0,q)\in\partial\Omega=\{x_n=0\}$,
and therefore
\begin{equation}\label{MS:POOSK-BILASNSTPAMLOS3-02}
q=0.
\end{equation}
In any case, from either~\eqref{MS:POOSK-BILASNSTPAMLOS3-01}
or~\eqref{MS:POOSK-BILASNSTPAMLOS3-02}, it follows that~$v_\infty(0,q)=
v_\infty(0,0)$.
Consequently, by~\eqref{eokfmgthasmdaccKMD:5678}
and~\eqref{dfgheokfmgthasmdaccKMD:5678BIS},
we deduce that~$\sigma=v_\infty(0,q)=v_\infty(0,0)\le0$.

Therefore, we find that~$u(x)\le\e x_n$ for every~$x\in\Omega$,
and thus the desired result follows by sending~$\e\searrow0$.
\end{proof}

For further reference, we also point out a Maximum Principle in the Sobolev space setting
(for the sake of simplicity,
we focus on a concrete result, though the method
of proof is more general, see e.g.~\cite[Theorem~8.1]{MR1814364}). As customary, here and in the following,
we denote by~${\mathcal{D}}^{1,2}(\Omega)$ the closure
of~$C^\infty_0(\Omega)$ with respect to the seminorm~$[u]_{H^1(\Omega)}:=\sqrt{\int_\Omega|\nabla u(x)|^2\,dx}$:

\begin{lemma}\label{SOBOMAXPLEH10}
Let~$n\ge3$ and~${\mathcal{U}}\subset\R^n$ be open and with uniformly Lipschitz boundary.
Let~$u\in C^2({\mathcal{U}})\cap {\mathcal{D}}^{1,2}(\R^n)$.
Suppose that~$u$ is harmonic in~${\mathcal{U}}$ and that~$u\ge0$ on~$\partial{\mathcal{U}}$ (in the trace sense).
Suppose also that~$|{\mathcal{U}}|=+\infty$.

Then, $u\ge0$ a.e. in~${\mathcal{U}}$.\index{Maximum Principle}
\end{lemma}

\begin{proof} Let~$v:=\max\{-u,0\}$. Since~$v$ belongs to~${\mathcal{D}}^{1,2}(\R^n)$,
we take a sequence of functions~$v_j\in C^\infty_0(\R^n)$ such that~$\nabla v_j\to\nabla v$
in~$L^2(\R^n)$. Thus, using the fact that~$u$ is harmonic and the first 
Green's Identity in~\eqref{GRr1},
we see that
$$ \int_{\mathcal{U}}\nabla u(x)\cdot \nabla v_j(x)\,dx=0$$
and therefore, sending~$j\to+\infty$,
$$ \int_{\mathcal{U}}\nabla u(x)\cdot \nabla v(x)\,dx=0.$$
As a consequence, since~$\nabla v=-\nabla u\chi_{\{u<0\}}$ (see e.g.~\cite[Corollary~6.18]{MR1817225}),
$$ 0=-\int_{{\mathcal{U}}\cap\{u<0\}}|\nabla u(x)|^2\,dx
=-\int_{{\mathcal{U}}}|\nabla v(x)|^2\,dx.$$
This gives that~$v$ is necessarily constant in~${\mathcal{U}}$ (see e.g.~\cite[Theorem~6.11]{MR1817225}),
say~$v(x)=c$ a.e.~$x\in{\mathcal{U}}$, for some~$c\in\R$.

Hence, using that~$u\in L^{\frac{2n}{n-2}}(\R^n)$
(thanks to
the Gagliardo-Nirenberg-Sobolev Inequality, see e.g.~\cite[Theorem~12.4]{MR2527916}),
$$ +\infty>\int_{{\mathcal{U}}} |u(x)|^{\frac{2n}{n-2}}\,dx\ge\int_{{\mathcal{U}}} |v(x)|^{\frac{2n}{n-2}}\,dx=
\int_{{\mathcal{U}}} |c|^{\frac{2n}{n-2}}\,dx.$$
Since~${\mathcal{U}}$ has infinite measure, we thus infer that~$c=0$.

Recapitulating, we have proven that~$v=\max\{-u,0\}=0$ a.e. in~${\mathcal{U}}$,
from which we obtain the desired result.
\end{proof}

\section{The Green Function}\label{Green Function:SECT}

Let~$n\ge2$. Given
an open set~$\Omega\subset\R^n$
with~$C^1$ boundary and a point~$x_0\in\Omega$,
we consider the fundamental solution~$\Gamma$
presented
in Theorem~\ref{KScvbMS:0okrmt4ht-VAR}
and set~$\Gamma^{(x_0)}(x):=\Gamma(x-x_0)$.
We also consider a function~$\Psi^{(x_0)}\in C^2(\Omega)\cap C(\overline\Omega)$ that satisfies
\begin{equation}\label{L:S:ASM S}\begin{dcases}
\Delta\Psi^{(x_0)}=0&{\mbox{ in }}\Omega,\\
\Psi^{(x_0)}=\Gamma^{(x_0)}&{\mbox{ on }}\partial\Omega.
\end{dcases}\end{equation}
In the literature, the function~$\Psi^{(x_0)}$ is sometimes called Robin Function.\index{Robin Function}

We stress that~$\Gamma^{(x_0)}$ is always finite on~$\partial\Omega$,
hence the right hand side of the boundary condition in~\eqref{L:S:ASM S}
is well-defined. Also, if~$\Omega$ is bounded, then
the solution to~\eqref{L:S:ASM S} is unique, due to the
Maximum Principle in Corollary~\ref{UNIQUENESSTHEOREM}.

When~$\Omega$ is unbounded, the solution
to~\eqref{L:S:ASM S} is not unique: for instance if~$\Omega=\{x_n>0\}$
and~$\Psi^{(x_0)}$ is a solution of~\eqref{L:S:ASM S} then, for example,
so is~$\Psi^{(x_0)}+c_1x_n+c_2 x_{n-1}x_n+c_3(x_n^3-3x_{n-1}^2x_n )$,
for all~$c_1,c_2,c_3\in\R$.

In our setting, the function~$\Psi^{(x_0)}$
will act as a ``corrector for the fundamental solution'': namely,
for all~$x\in\R^n\setminus\{x_0\}$
we define
\begin{equation}\label{FLGREEN} G(x,x_0):=\Gamma(x-x_0)-\Psi^{(x_0)}(x).\end{equation}
The function~$G$ is called the Green Function of~$\Omega$.\index{Green Function}

In the next lemma\footnote{It is worth observing \label{POSI:GREEN:Ut}
that the Green Function is positive in a bounded set~$\Omega$. Indeed, by the singular behavior of the fundamental solution at the pole, we know that if~$\delta>0$ is sufficiently small, the harmonic function~$G(\cdot,x_0)$ is positive along~$\partial B_\delta(x_0)\subset\Omega$. Since this function vanishes along~$\partial\Omega$, we infer from the Maximum Principle that~$G(\cdot,x_0)>0$ in~$\Omega\setminus B_\delta(x_0)$. Thus, since~$\delta$ can be taken arbitrarily small, it follows that~$G(\cdot,x_0)>0$ in~$\Omega$.} we provide a symmetry result for the Green Function:

\begin{lemma}\label{SYJMDMMAD}
For every~$x$, $y\in\overline\Omega$ with~$x\ne y$,
\begin{equation*}
G(x,y)=G(y,x).
\end{equation*}
\end{lemma}

\begin{proof} Let~$a$, $b\in\Omega$ with~$a\ne b$.
Let~$\rho>0$ sufficiently small such that~$B_\rho(a)\cup B_\rho(b)\Subset\Omega$
and~$B_\rho(a)\cap B_\rho(b)=\varnothing$. Let~$\Omega_\rho:=
\Omega\setminus(B_\rho(a)\cup B_\rho(b))$, $\alpha(x):=G(x,a)$, $\beta(x):=G(x,b)$. Notice that, if~$x\in\partial\Omega$,
then~$\alpha(x)=\Gamma(x-a)-\Psi^{(a)}(x)=\Gamma(x-a)-\Gamma^{(a)}(x)=0$,
and similarly~$\beta(x)=0$. Thus, making use of
the second Green's Identity~\eqref{GRr2}, we have
\begin{equation*}
\begin{split}&\int_{\Omega_\rho} \Big( \alpha(x)\Delta \beta(x)-\beta(x)\Delta \alpha(x)\Big)\,dx=
\int_{\partial{\Omega_\rho}}\left(\alpha(x)\frac{\partial \beta}{\partial\nu}(x)-
\beta(x)\frac{\partial\alpha}{\partial\nu}(x)\right)\,d{\mathcal{H}}^{n-1}_x\\
&\qquad=
\int_{\partial{B_\rho}(a)}\left(
\beta(x)\frac{\partial\alpha}{\partial\nu}(x)-
\alpha(x)\frac{\partial \beta}{\partial\nu}(x)\right)\,d{\mathcal{H}}^{n-1}_x
+
\int_{\partial{B_\rho}(b)}\left(
\beta(x)\frac{\partial\alpha}{\partial\nu}(x)-
\alpha(x)\frac{\partial \beta}{\partial\nu}(x)\right)\,d{\mathcal{H}}^{n-1}_x.
\end{split}
\end{equation*}
Since
in~$\Omega_\rho$ it holds
that~$\Delta\alpha=\Delta\Gamma^{(a)}-\Delta\Psi^{(a)}(x)=0$
and similarly~$\Delta\beta=0$, we thereby conclude that
\begin{equation}\label{134rty-9DM-2}
0=
\int_{\partial{B_\rho}(a)}\left(
\beta(x)\frac{\partial\alpha}{\partial\nu}(x)-
\alpha(x)\frac{\partial \beta}{\partial\nu}(x)\right)\,d{\mathcal{H}}^{n-1}_x
+
\int_{\partial{B_\rho}(b)}\left(
\beta(x)\frac{\partial\alpha}{\partial\nu}(x)-
\alpha(x)\frac{\partial \beta}{\partial\nu}(x)\right)\,d{\mathcal{H}}^{n-1}_x.
\end{equation}
In addition,
$$\lim_{\rho\searrow0}
\fint_{\partial{B_\rho}(a)}
\beta(x)\frac{\partial}{\partial\nu} \Psi^{(a)}(x)\,d{\mathcal{H}}^{n-1}_x=
\beta(a)\frac{\partial}{\partial\nu} \Psi^{(a)}(a)=
\big(\Gamma(a-b)-\Psi^{(b)}(a)\big)\Psi^{(a)}(a)$$
and
$$\lim_{\rho\searrow0}
\fint_{\partial{B_\rho}(a)}
\Psi^{(a)}(x)\frac{\partial \beta}{\partial\nu}(x)\,d{\mathcal{H}}^{n-1}_x
=\Psi^{(a)}(a)\frac{\partial \beta}{\partial\nu}(a)
=\Psi^{(a)}(a) \frac{\partial \Gamma}{\partial\nu}(a-b)-\Psi^{(a)}(a) \frac{\partial \Psi^{(b)}}{\partial\nu}(a),
$$
whence
$$\lim_{\rho\searrow0}
\int_{\partial{B_\rho}(a)}
\beta(x)\frac{\partial}{\partial\nu} \Psi^{(a)}(x)\,d{\mathcal{H}}^{n-1}_x=0\qquad{\mbox{and}}\qquad
\lim_{\rho\searrow0}
\int_{\partial{B_\rho}(a)}
\Psi^{(a)}(x)\frac{\partial \beta}{\partial\nu}(x)\,d{\mathcal{H}}^{n-1}_x
=0.$$
For this reason,
\begin{equation}\label{ascikiuads9ddwfvdbSND2}
\begin{split}
&\lim_{\rho\searrow0}
\int_{\partial{B_\rho}(a)}\left(
\beta(x)\frac{\partial\alpha}{\partial\nu}(x)-
\alpha(x)\frac{\partial \beta}{\partial\nu}(x)\right)\,d{\mathcal{H}}^{n-1}_x\\
=\,&
\lim_{\rho\searrow0}
\int_{\partial{B_\rho}(a)}\left(
\beta(x)\frac{\partial}{\partial\nu}\big(\Gamma(x-a)-\Psi^{(a)}(x)\big)-
\big(\Gamma(x-a)-\Psi^{(a)}(x)\big)\frac{\partial \beta}{\partial\nu}(x)\right)\,d{\mathcal{H}}^{n-1}_x\\
=\,&
\lim_{\rho\searrow0}
\int_{\partial{B_\rho}(a)}\left(\beta(x)\frac{\partial\Gamma}{\partial\nu}(x-a)-
\Gamma(x-a)
\frac{\partial \beta}{\partial\nu}(x)
\right)\,d{\mathcal{H}}^{n-1}_x.
\end{split}\end{equation}
Now, since, by~\eqref{GAMMAFU},
$$ \int_{\partial{B_\rho}(a)}
\Gamma(x-a)
\frac{\partial \beta}{\partial\nu}(x)\,d{\mathcal{H}}^{n-1}_x=
\begin{dcases}
\frac{c_n}{\rho^{n-2}}\int_{\partial{B_\rho}(a)}
\frac{\partial \beta}{\partial\nu}(x)\,d{\mathcal{H}}^{n-1}_x & {\mbox{ if $n\ne2$,}}\cr
-c_n\ln\rho\int_{\partial{B_\rho}(a)}
\frac{\partial \beta}{\partial\nu}(x)\,d{\mathcal{H}}^{n-1}_x & {\mbox{ if $n=2$,}}\cr
\end{dcases}$$
we find that
\begin{equation}\label{ascikiuads9ddwfvdbSND}
\lim_{\rho\searrow0}
\int_{\partial{B_\rho}(a)}
\Gamma(x-a)
\frac{\partial \beta}{\partial\nu}(x)\,d{\mathcal{H}}^{n-1}_x=0.
\end{equation}
Besides,
\begin{eqnarray*}&&\lim_{\rho\searrow0}
\int_{\partial{B_\rho}(a)}
\beta(x)\frac{\partial\Gamma}{\partial\nu}(x-a)
\,d{\mathcal{H}}^{n-1}_x=
\lim_{\rho\searrow0}
\frac{\widetilde{c}_n}{\rho^{n-1}}\int_{\partial{B_\rho}(a)}
\beta(x)
\,d{\mathcal{H}}^{n-1}_x=C_n\,\beta(a)
\end{eqnarray*}
for some~$\widetilde{c}_n$, $C_n\in\R$.

Plugging this information and~\eqref{ascikiuads9ddwfvdbSND}
into~\eqref{ascikiuads9ddwfvdbSND2}, we find that
\begin{equation}\label{134rty-9DM}
\lim_{\rho\searrow0}
\int_{\partial{B_\rho}(a)}\left(
\beta(x)\frac{\partial\alpha}{\partial\nu}(x)-
\alpha(x)\frac{\partial \beta}{\partial\nu}(x)\right)\,d{\mathcal{H}}^{n-1}_x=
C_n\,\beta(a).
\end{equation}
Exchanging the roles of~$a$ and~$b$ (and of~$\alpha$ and~$\beta$),
we also have that
\begin{equation*}
\lim_{\rho\searrow0}
\int_{\partial{B_\rho}(b)}\left(
\alpha(x)\frac{\partial\beta}{\partial\nu}(x)-
\beta(x)\frac{\partial \alpha}{\partial\nu}(x)\right)\,d{\mathcal{H}}^{n-1}_x=
C_n\,\alpha(b).
\end{equation*}
{F}rom this and~\eqref{134rty-9DM}, we deduce from~\eqref{134rty-9DM-2}
that~$\alpha(b)=\beta(a)$. That is,
\begin{equation*}G(b,a)=
\alpha(b)=\beta(a)=G(a,b).\qedhere
\end{equation*}
\end{proof}

The importance of the Green Function in the theory of partial differential
equations lies in the possibility of reconstructing a solution
from its boundary value. Specifically,
we already know from Green's Representation Formula
in~\eqref{300G} how to reconstruct a function from its Laplacian
and the boundary values of the function {\em and
of its normal derivative}:
the Green Function allows us to ``minimize the necessary information''
reconstructing the function
merely from its Laplacian
and the boundary values (without the necessity
of knowing additionally the values of the normal derivative at the boundary).
Indeed, we have that:

\begin{theorem}\label{MS:SKMD344567yjghS}
Let~$\Omega$ be a bounded open set in~$\R^n$
with~$C^1$ boundary.
Then, for every~$u\in C^2(\Omega)\cap C^1(\overline\Omega)$,
$$ u(x)=-\int_\Omega G(x,y)\,\Delta u(y)\,dy-
\int_{\partial \Omega}u(y)\,\frac{\partial G}{\partial\nu}(x,y)\,d{\mathcal{H}}^{n-1}_y,$$
where, for every~$x\in\Omega$ and~$y\in\partial\Omega$,
\begin{equation}\label{OKJN-NOAYHNSTAZ} \frac{\partial G}{\partial\nu}(x,y):=
\nabla_yG(x,y)\cdot\nu(y)=
\sum_{i=1}^n \frac{\partial G}{\partial y_i}(x,y)\,\nu(y)\cdot e_i.\end{equation}
\end{theorem}

\begin{proof} We note that~$\Psi^{(x)}(y)=\Psi^{(y)}(x)$, due to the symmetry
property of the Green Function in Lemma~\ref{SYJMDMMAD},
the one of the fundamental solution in~\eqref{GAMMAFU}
and the setting in~\eqref{FLGREEN}.
As a consequence, in the notation of~\eqref{OKJN-NOAYHNSTAZ},
for all~$x\in\Omega$ and~$y\in\partial\Omega$,
\begin{eqnarray*}&&
\frac{\partial \Psi^{(x)}}{\partial\nu}(y)=
\sum_{i=1}^n \frac{\partial }{\partial y_i}\big(\Psi^{(x)}(y)\big)\,\nu(y)\cdot e_i=
\sum_{i=1}^n \frac{\partial }{\partial y_i}\big(\Psi^{(y)}(x)\big)\,\nu(y)\cdot e_i
\\
&&\quad
=
\sum_{i=1}^n \frac{\partial }{\partial y_i}\big(
\Gamma(x-y)-G(x,y)
\big)\,\nu(y)\cdot e_i=-\nabla \Gamma(x-y)\cdot\nu(y)
-\sum_{i=1}^n \frac{\partial }{\partial y_i}G(x,y)\,\nu(y)\cdot e_i
\\&&\quad=
-\frac{\partial \Gamma}{\partial\nu}(x-y)
-\frac{\partial G}{\partial\nu}(x,y).
\end{eqnarray*}
Thus, we use Green's Representation Formula~\eqref{300G}
combined with~\eqref{L:S:ASM S}, recalling also
the second Green's Identity~\eqref{GRr2}, to see that
\begin{equation*}
\begin{split}
&\int_\Omega G(x,y)\,\Delta u(y)\,dy+
\int_{\partial \Omega}u(y)\,\frac{\partial G}{\partial\nu}(x,y)\,d{\mathcal{H}}^{n-1}_y
\\=\,&\int_\Omega \Gamma(x-y)\,\Delta u(y)\,dy
-\int_\Omega \Psi^{(y)}(x)\,\Delta u(y)\,dy-
\int_{\partial \Omega}u(y)\,\frac{\partial \Psi^{(x)}}{\partial\nu}(y)\,d{\mathcal{H}}^{n-1}_y\\&\qquad\qquad
-\int_{\partial \Omega}u(y)\,\frac{\partial \Gamma}{\partial\nu}(x-y)\,d{\mathcal{H}}^{n-1}_y
\\ =\,& \int_\Omega \Gamma(x-y)\,\Delta u(y)\,dy-
\int_{\partial \Omega}\Psi^{(x)}(y)\,\frac{\partial u}{\partial\nu}(y)\,d{\mathcal{H}}^{n-1}_y
-\int_{\partial \Omega}u(y)\,\frac{\partial \Gamma}{\partial\nu}(x-y)\,d{\mathcal{H}}^{n-1}_y\\=\,&
\int_\Omega \Gamma(y-x)\,\Delta u(y)\,dy-
\int_{\partial \Omega}\Gamma(y-x)\,\frac{\partial u}{\partial\nu}(y)\,d{\mathcal{H}}^{n-1}_y
+\int_{\partial \Omega}u(y)\,\frac{\partial \Gamma}{\partial\nu}(y-x)\,d{\mathcal{H}}^{n-1}_y\\=\,&
-u(x). 
\qedhere\end{split}\end{equation*}
\end{proof}

A straightforward consequence of Theorem~\ref{MS:SKMD344567yjghS}
is a representation result for solutions of the boundary value problem
\begin{equation}\label{MSIKMPOOKSssom x}\begin{dcases}
-\Delta u =f& {\mbox{ in }}\Omega,\cr
u=g& {\mbox{ on }}\partial\Omega,\end{dcases}
\end{equation}
as detailed\footnote{In particular, by Corollary~\ref{SKMD:DOK0oerlfgdit},
solutions of
\begin{equation*}\begin{dcases}
-\Delta u =f& {\mbox{ in }}\Omega,\cr
u=0& {\mbox{ on }}\partial\Omega,\end{dcases}
\end{equation*}
take the form \label{CONTECONT09F}
$$ u(x)=\int_\Omega f(y)\,G(x,y)\,dy.$$
This can be seen as the continuous counterpart of formula~\eqref{CONTECONT09}
which was obtained in the discrete setting.}
in the following result:

\begin{corollary}\label{SKMD:DOK0oerlfgdit}
Let~$\Omega$ be a bounded open set in~$\R^n$
with~$C^1$ boundary.
If~$u\in C^2(\Omega)\cap C^1(\overline\Omega)$
solves~\eqref{MSIKMPOOKSssom x} for some~$f:\Omega\to\R$
and~$g:\partial\Omega\to\R$,
then, for every~$x\in\Omega$,
$$ u(x)=\int_\Omega f(y)\,G(x,y)\,dy-\int_{\partial \Omega}
g(y)\,\frac{\partial G}{\partial\nu}(x,y)\,d{\mathcal{H}}^{n-1}_y,$$
where the notation in~\eqref{OKJN-NOAYHNSTAZ} has been used.
\end{corollary}

Given the result in Corollary~\ref{SKMD:DOK0oerlfgdit}, to have an explicit
expression of the solution
of the boundary value problem~\eqref{MSIKMPOOKSssom x},
it would be desirable to know an explicit expression of the Green Function~$G$.
Unfortunately, this is rarely possible in concrete cases,
but there are at least two interesting situations
in which the Green Function can be written in terms
of elementary functions. These two examples
are given by the cases in which the domain is either a ball or a halfspace,
as discussed in the forthcoming Theorems~\ref{MS:SKMD344567yjghS-1}
and~\ref{MS:SKMD344567yjghS-2}.

\begin{figure}
  \centering
  \includegraphics[width=.4\linewidth]{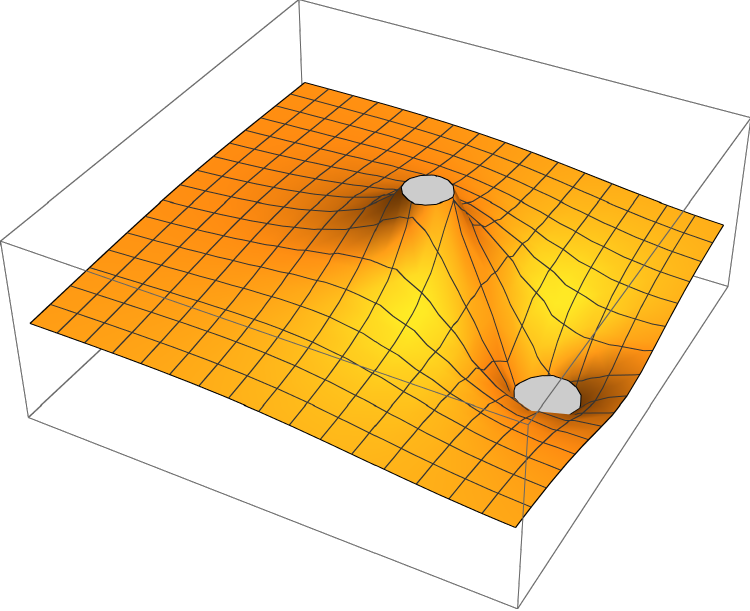}$\qquad$
  \includegraphics[width=.4\linewidth]{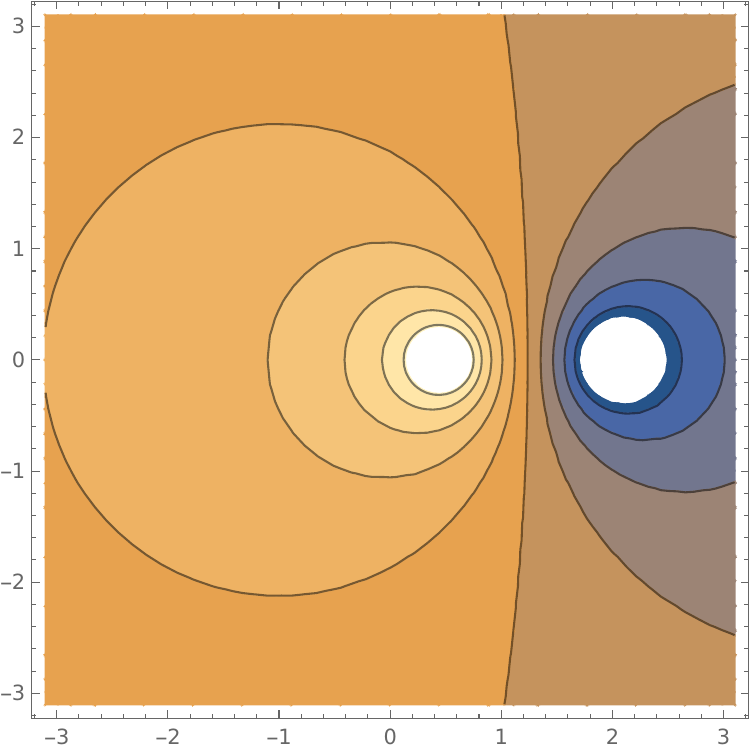}
 \caption{\sl Plot of the Green Function of the unit
ball in dimension~$2$
when~$x_0=\left(\frac12,0\right)$,
and its level sets.
Notice that~$\partial B_1$ is a zero level set.}\label{GREEBsoloDItangeFI}
\end{figure}

\begin{figure}
  \centering
  \includegraphics[width=.4\linewidth]{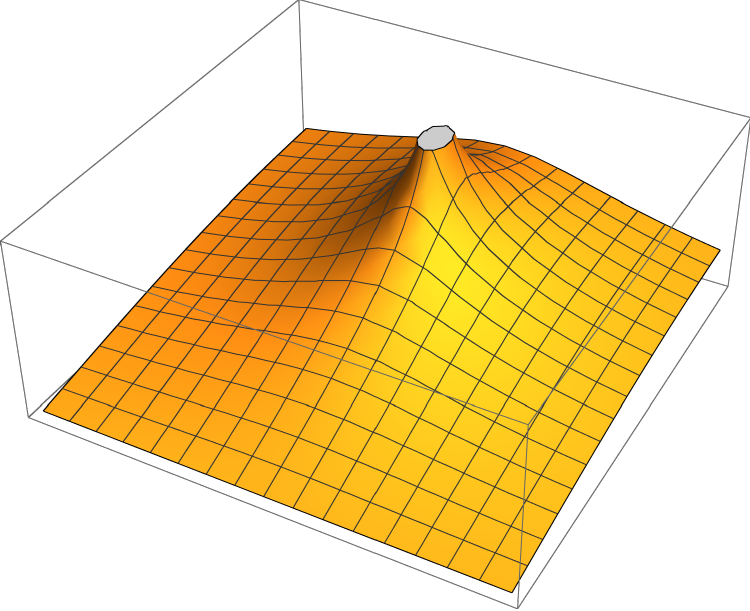}$\qquad$
  \includegraphics[width=.4\linewidth]{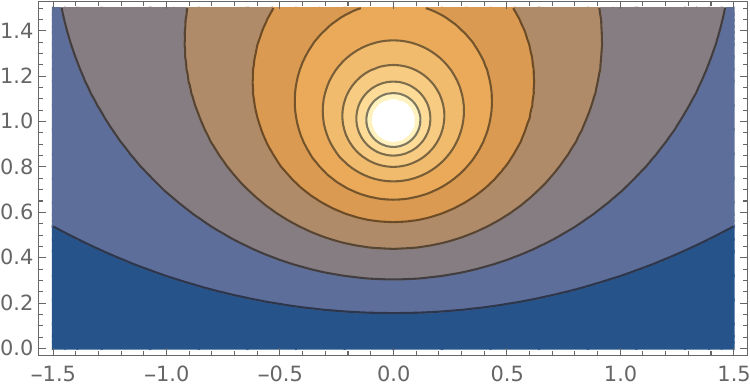}
 \caption{\sl Plot of the Green Function of the halfplane
when~$x_0=(0,1)$,
and its level sets.
Notice that~$\R\times\{0\}$ is a zero level set.}\label{SEMIGREEBsoloDItangeFI}
\end{figure}

\begin{theorem}\label{MS:SKMD344567yjghS-1}
The Green Function\footnote{The physical intuition underneath 
Theorem~\ref{MS:SKMD344567yjghS-1} can be grasped by taking,
up to scaling, $R:=1$ and by recalling the 
method of electrostatic image charges as discussed in Section~\ref{LOO}.
Namely, by construction, we have that in the setting of
Theorem~\ref{MS:SKMD344567yjghS-1},
$-\Delta G(x,x_0)=\Delta\Psi^{(x_0)}(x)-\Delta\Gamma(x-x_0)=\delta_{x_0}(x)$
for all~$x\in B_1$. In this sense~$G(x,x_0)$ corresponds, up to normalizing
constants, to
the potential of a unit electrostatic charge located at~$x_0$.
On the other hand, we know that, if~$x\in\partial B_1$,
then~$G(x,x_0)= \Gamma(x-x_0)-\Psi^{(x_0)}(x)=0$,
whence~$\partial B_1$ is an equipotential surface. The determination
of the Green Function of~$B_1$ is therefore reduced to
that of an electrostatic potential generated by a unit charge in~$x_0$
and an auxiliary charge outside~$B_1$ that makes~$\partial B_1$
an equipotential surface.
We know from Section~\ref{LOO} (see in particular~\eqref{SOLLE:1}
and~\eqref{SOLLE:2},
and compare with the Kelvin Transform in~\eqref{KSM:J:KELVIN1})
that to produce such an electrostatic potential (say, when~$n\ne2$),
it suffices
to place a charge of intensity~$-\frac1{|x_0|^{n-2}}$
at the position~${\mathcal{K}}(x_0)$. Thus, the corresponding
electrostatic potential takes the form~$\frac{1}{|x-x_0|^{n-2}}
-\frac1{|x_0|^{n-2}\,|x-{\mathcal{K}}(x_0)|^{n-2}}$,
which agrees with~\eqref{GERTBAL} up to normalizing constants,
and this gives a concrete justification for the result
in Theorem~\ref{MS:SKMD344567yjghS-1} (say, when~$n\ne2$, and the case~$n=2$
is discussed in the footnote of page~\pageref{CONSDcenrfPTHDNULAMDALSIKLASKGIMDOVAS}).
See also Figure~\ref{GREEBsoloDItangeFI} for graphical images
of the Green Function of a ball (to be compared also with Figure~\ref{SEMIGREEBsoloDItangeFI}
that represents a Green Function for a halfspace).}
of the ball~$B_R$ is
\begin{equation}\label{GERTBAL}
G(x,x_0)=\begin{dcases}\Gamma(x-x_0)-\Gamma\left(\frac{|x_0|}R\,\left(
x-\frac{R^2x_0}{|x_0|^2}
\right)\right)&{\mbox{ if }}x_0\in B_R\setminus\{0\},\\
\Gamma(x)-\Gamma(Re_1)&{\mbox{ if }}x_0=0
.\end{dcases}\end{equation}
\end{theorem}

\begin{proof} We remark that the function defined in~\eqref{GERTBAL}
when~$x_0\in B_R\setminus\{0\}$ is continuously extended to~$x_0=0$ since
\begin{eqnarray*}&& 
\frac{|x_0|}R\,\left|
x-\frac{R^2x_0}{|x_0|^2}
\right|=\frac{|x_0|}R\,
\frac{\big| |x_0|^2 x-R^2x_0\big|}{|x_0|^2}=R\,
\left| \frac{|x_0| x}{R^2}-\frac{x_0}{|x_0|}\right|\\&&\qquad=R\sqrt{1-
\frac{2 x\cdot x_0}{R^2}
+\frac{|x_0|^2 |x|^2}{R^4}
}=R(1+O(|x_0|)).
\end{eqnarray*}
Thus, to prove the claim in Theorem~\ref{MS:SKMD344567yjghS-1},
we consider
$$ \Psi^{(x_0)}(x):
=\begin{dcases}\Gamma\left(\frac{|x_0|}R\,\left(
x-\frac{R^2x_0}{|x_0|^2}
\right)\right)&{\mbox{ if }}x_0\in B_R\setminus\{0\},\\
\Gamma(Re_1)&{\mbox{ if }}x_0=0
\end{dcases}$$
and we need to show that
\begin{equation}\label{cdsecnche dcdknontinmdcecjehia SN}
{\mbox{$\Psi^{(x_0)}$ satisfies~\eqref{L:S:ASM S}
with~$\Omega:=B_R$.}}\end{equation}
As a matter of fact, since~$\left| \frac{R^2x_0}{|x_0|^2}\right|=
\frac{R^2}{|x_0|}>R$, we have that~$
\Delta\Psi^{(x_0)}(x)=0$ for every~$x\in B_R$.
Moreover, if~$x\in\partial B_R$, recalling~\eqref{KSM:J:KELVIN1}
and~\eqref{KSM:J:KELVIN9}, we see that
\begin{equation}\label{UJdgvbf fficpqlwdfdvbia0comPucjv}
\left|\frac{|x_0|}R\,\left(
x-\frac{R^2x_0}{|x_0|^2}\right)\right|=|x_0|\,
\left|\frac{x}{R}-\frac{Rx_0}{|x_0|^2}\right|=
|x_0|\,
\left|\frac{x}{R}-{\mathcal{K}}\left(\frac{x_0}{R}\right)\right|=R\,
\left|\frac{x}{R}- \frac{x_0}{R}\right|=|x-x_0|.\end{equation}
Therefore, for every~$x\in\partial B_R$,
we have that~$ \Psi^{(x_0)}(x)=\Gamma(x-x_0)=
\Gamma^{(x_0)}(x)$,
and these observations establish~\eqref{cdsecnche dcdknontinmdcecjehia SN}.
\end{proof}

\begin{theorem}\label{MS:SKMD344567yjghS-2}
A Green Function\footnote{The physical intuition underneath 
Theorem~\ref{MS:SKMD344567yjghS-2} is that the
method of electrostatic image charges, as described in Section~\ref{LOO}
when the domain is a ball, dramatically simplifies
when the domain is a halfspace,
since in this case it suffices to place
an image charge with the same intensity
in a symmetrical position with respect to the boundary of the halfspace.}
of the halfspace~$\{x_n>0\}$ is $$
G(x,x_0)=\Gamma(x-x_0)-\Gamma({\mathcal{R}}(x)-x_0),$$
where~${\mathcal{R}}(x_1,\dots,x_{n-1},x_n)=(x_1,\dots,x_{n-1},-x_n)$.
\end{theorem}

\begin{proof} We let~$\Psi^{(x_0)}(x):=\Gamma({\mathcal{R}}(x)-x_0)$
and we need to show that
\begin{equation}\label{cdsecnche dcdknontinmdcecjehia SN-2}
{\mbox{$\Psi^{(x_0)}$ satisfies~\eqref{L:S:ASM S}
with~$\Omega:=\{x_n>0\}$.}}\end{equation}
For this, we observe that~$|{\mathcal{R}}(x)-x_0|=|x-{\mathcal{R}}(x_0)|$,
hence~$\Psi^{(x_0)}(x)=\Gamma(x-{\mathcal{R}}(x_0))$
and consequently~$\Delta\Psi^{(x_0)}=0$ in~$\{x_n>0\}$.
Furthermore, if~$x_n=0$ then~${\mathcal{R}}(x)=x$
and~$\Psi^{(x_0)}(x)=\Gamma(x-x_0)$. These remarks
prove the validity of~\eqref{cdsecnche dcdknontinmdcecjehia SN-2}.
\end{proof}

For an extensive treatment of Green Functions see~\cite{MR1888091} and the references therein.
We now  recall an interesting relation
about the surface area of the level sets of the Green Function
(we focus on the case of the Laplace operator, but the approach
to the problem is rather general, see~\cite{MR0145191} for full details):

\begin{lemma}\label{VDEFGRE-LE}
Let~$\Omega$ be a bounded and open subset of~$\R^n$ and let~$x_0\in\Omega$.
Let~$G$ be the Green Function of~$\Omega$ and, for every~$t\ge0$, set
\begin{equation}\label{VDEFGRE} v(t):=\Big|
\big\{ x\in\Omega{\mbox{ s.t. }}G(x,x_0)>t\big\}
\Big|.\end{equation}
Then,
$$ \big(v(t)\big)^{\frac{2-n}{n}}-|\Omega|^{\frac{2-n}{n}}\ge
n(n-2)|B_1|^{\frac2{n}}t.$$
\end{lemma}

\begin{proof} Let~$S(t):={\mathcal{H}}^{n-1}\big(\big\{ x\in\Omega{\mbox{ s.t. }}
G(x,x_0)=t\big\}\big)$.
By the Cauchy-Schwarz Inequality,
\begin{equation}\label{PTRFVSOMGSHDNBODMS9}\begin{split}
S(t)\,&=\,\int_{\{G(\cdot,x_0)=t\}}d{\mathcal{H}}^{n-1}_x
\\&=\,\int_{\{G(\cdot,x_0)=t\}}
\frac{\sqrt{|\nabla G(x,x_0)|}\,d{\mathcal{H}}^{n-1}_x }{\sqrt{|\nabla G(x,x_0)|}}
\\ &\le\,\sqrt{
\int_{\{G(\cdot,x_0)=t\}}
|\nabla G(x,x_0)|\,d{\mathcal{H}}^{n-1}_x \;
\int_{\{G(\cdot,x_0)=t\}}
\frac{d{\mathcal{H}}^{n-1}_x }{|\nabla G(x,x_0)|}}.\end{split}
\end{equation}
Additionally, by the Coarea Formula (see e.g.~\cite{MR3409135}),
\begin{eqnarray*}
v(t)&=&\int_{\Omega}\frac{\chi_{(t,+\infty)}(G(x,x_0))\;|\nabla G(x,x_0)|}{|\nabla G(x,x_0)|}\,dx
\\&=&\int_\R\left[
\int_{\{G(\cdot,x_0)=s\}}
\frac{\chi_{(t,+\infty)}(G(x,x_0)) }{|\nabla G(x,x_0)|}
\,d{\mathcal{H}}^{n-1}_x\right]\,ds\\&=&\int^{+\infty}_t\left[
\int_{\{G(\cdot,x_0)=s\}}
\frac{d{\mathcal{H}}^{n-1}_x }{|\nabla G(x,x_0)|}
\right]\,ds
\end{eqnarray*}
and therefore, for a.e.~$t\ge0$,
\begin{equation*}
-v'(t)=\int_{\{G(\cdot,x_0)=t\}}
\frac{d{\mathcal{H}}^{n-1}_x }{|\nabla G(x,x_0)|}.
\end{equation*}
By combining this and~\eqref{PTRFVSOMGSHDNBODMS9}, we obtain
\begin{equation}\label{TYHSNOTTUSGALSXYJMD}
-v'(t)\ge
\left(\int_{\{G(\cdot,x_0)=t\}}
|\nabla G(x,x_0)|\,d{\mathcal{H}}^{n-1}_x\right)^{-1}\,\big(S(t)\big)^2.
\end{equation}
By the Isoperimetric Inequality (see e.g.~\cite{MR775682}),
$$ 
\frac{S(t)}{\big(v(t)\big)^{\frac{n-1}{n}}}\ge
\frac{{\mathcal{H}}^{n-1}(\partial B_1)}{|B_1|^{\frac{n-1}{n}}},$$
which, together with~\eqref{TYHSNOTTUSGALSXYJMD},
gives that
\begin{equation}\label{USC-10-0}
-v'(t)\ge
\left(\int_{\{G(\cdot,x_0)=t\}}
|\nabla G(x,x_0)|\,d{\mathcal{H}}^{n-1}_x\right)^{-1}\,
\frac{\big({\mathcal{H}}^{n-1}(\partial B_1)\big)^2}{|B_1|^{\frac{2(n-1)}{n}}}\,
\big(v(t)\big)^{\frac{2(n-1)}{n}}.
\end{equation}
We stress that, a.e.~$t\ge0$, the set~$\{G(\cdot,x_0)=t\}$
is a manifold of class~$C^1$ which corresponds to the boundary of~$\{G(\cdot,x_0)>t\}$,
because~$\nabla G(\cdot,x_0)\ne0$ along such set (as a consequence
of the Morse-Sard Theorem, see e.g. Theorem~4.3 in~\cite{MR1503449}).
As a consequence, a.e.~$t\ge0$, along~$\{G(\cdot,x_0)=t\}$ we have that
$$ |\nabla G(\cdot,x_0)|=- \nabla G(\cdot,x_0)\cdot\nu(\cdot).$$
Notice that in this setting~$\nu(\cdot)$ denotes the exterior unit normal to the set~$\{G(\cdot,x_0)>t\}$.
 
Also, since~$G(x_0,x_0)=+\infty$, we have that~$x_0\in\{G(\cdot,x_0)>t\}$.
Consequently, using that~$G$ is harmonic in~$\Omega\setminus\{x_0\}\supseteq
\{G(\cdot,x_0)< t\}$,
\begin{eqnarray*}&&
\int_{\{G(\cdot,x_0)=t\}}
|\nabla G(x,x_0)|\,d{\mathcal{H}}^{n-1}_x=
-\int_{\{G(\cdot,x_0)=t\}}
\nabla G(x,x_0)\cdot\nu(x)\,d{\mathcal{H}}^{n-1}_x\\
&&\qquad=
\int_{\partial\Omega}
\nabla G(x,x_0)\cdot\nu(x)\,d{\mathcal{H}}^{n-1}_x
-\int_{\{G(\cdot,x_0)=t\}}
\nabla G(x,x_0)\cdot\nu(x)\,d{\mathcal{H}}^{n-1}_x\\&&\qquad\qquad\qquad
-\int_{\partial\Omega}
\nabla G(x,x_0)\cdot\nu(x)\,d{\mathcal{H}}^{n-1}_x
\\&&\qquad=
\int_{\{G(\cdot,x_0)<t\}}
\Delta G(x,x_0)\,d{\mathcal{H}}^{n-1}_x-\int_{\partial\Omega}
\nabla G(x,x_0)\cdot\nu(x)\,d{\mathcal{H}}^{n-1}_x
\\&&\qquad=-
\int_{\partial\Omega}
\nabla G(x,x_0)\cdot\nu(x)\,d{\mathcal{H}}^{n-1}_x.
\end{eqnarray*}
Combining this with Theorem~\ref{MS:SKMD344567yjghS} (used here with~$u:=1$)
we conclude that
\[
\int_{\{G(\cdot,x_0)=t\}}
|\nabla G(x,x_0)|\,d{\mathcal{H}}^{n-1}_x=1.
\]
This and~\eqref{USC-10-0} lead to
\begin{equation}\label{MAXLOS4567C4RO5gN0987AS}
-v'(t)\ge
\frac{\big({\mathcal{H}}^{n-1}(\partial B_1)\big)^2}{|B_1|^{\frac{2(n-1)}{n}}}\,
\big(v(t)\big)^{\frac{2(n-1)}{n}}.
\end{equation}
Since~$G$ is nonnegative by the Maximum Principle, we also know that
$$ v(0)=
\Big|
\big\{ x\in\Omega{\mbox{ s.t. }}G(x,x_0)>0\big\}
\Big|=|\Omega|.$$
{F}rom this and~\eqref{MAXLOS4567C4RO5gN0987AS} it follows that, for every~$T\ge0$,
\begin{eqnarray*}&&
\big(v(T)\big)^{\frac{2-n}{n}}-|\Omega|^{\frac{2-n}{n}}=
\big(v(T)\big)^{\frac{2-n}{n}}-\big(v(0)\big)^{\frac{2-n}{n}}=
\int_0^T \frac{d}{dt}\big(v(t)\big)^{\frac{2-n}{n}}\,dt\\
&&\qquad=
-\frac{n-2}{n}
\int_0^T\big(v(t)\big)^{\frac{2(1-n)}{n}}\,v'(t)\,dt\ge
\frac{(n-2)\big({\mathcal{H}}^{n-1}(\partial B_1)\big)^2\,T}{n\,|B_1|^{\frac{2(n-1)}{n}}}
,\end{eqnarray*}
yielding the desired result in view of~\eqref{B1}.
\end{proof}

In relation to Lemma~\ref{VDEFGRE-LE},
we also recall an interesting bound on the Robin Function
evaluated at the central point~$x_0$ (as pointed out in equation~(2.14)
of~\cite{MR1021402}):

\begin{lemma}\label{WEI-SHA}
Let~$n\ge3$. 
Let~$\Omega$ be a bounded and open subset of~$\R^n$ and let~$x_0\in\Omega$.
Then, the Robin Function of~$\Omega$ satisfies\footnote{We observe that the inequality in Lemma~\ref{WEI-SHA}
is sharp, since equality is attained when~$\Omega=B_1$ and~$x_0=0$: indeed,
in this case, by~\eqref{COIENEE}, \eqref{GAMMAFU}, \eqref{GERTBAL},
\begin{eqnarray*}&&
(n-2)\Psi^{(x_0)}(x_0)-
\frac{ 1}{n |B_1|^{\frac2n} |\Omega|^{\frac{n-2}{n}} }
=(n-2)\Psi^{(0)}(0)-\frac{1}{n|B_1|}\\
&&\qquad\qquad
=(n-2)\Gamma(e_1)-\frac{1}{n|B_1|}=c_n(n-2)-\frac{1}{n|B_1|}=0.
\end{eqnarray*}}
$$ (n-2)\Psi^{(x_0)}(x_0)\geq
\frac{ 1}{n|B_1|^{\frac2n} |\Omega|^{\frac{n-2}{n}} }
.$$
\end{lemma}

\begin{proof} Without loss of generality, up to a translation, we suppose that~$x_0=0$.
By~\eqref{COIENEE} and~\eqref{GAMMAFU}, we have that~$\Gamma(x)=
\frac{1}{{n(n-2)\,| B_1|}\,|x|^{n-2}}$ and therefore, recalling~\eqref{FLGREEN},
\begin{equation}\label{INVEFUSER} |x|
=\frac{1}{\big(n(n-2)\,| B_1|\big)^{\frac1{n-2}}\,\big(\Gamma(x)\big)^{\frac1{n-2}}}
=\frac{1}{\big(n(n-2)\,| B_1|\big)^{\frac1{n-2}}\,\big(G(x,0)+\Psi^{(0)}(x)\big)^{\frac1{n-2}}}.
\end{equation}
Let~$v$ be as in~\eqref{VDEFGRE}. We have, for a.e.~$t>0$,
\begin{eqnarray*}
v(t)&=&\frac1n\,\int_{\{G(\cdot,0)>t\}} \div(x)\,dx
\\&=&\frac1n\,\int_{\{G(\cdot,0)=t\}} x\cdot\nu(x)\,d{\mathcal{H}}^{n-1}_x\\
&\le&\frac1n\,\int_{\{G(\cdot,0)=t\}}|x|\,d{\mathcal{H}}^{n-1}_x\\
&=&\frac1n\,\int_{\{G(\cdot,0)=t\}}
\frac{d{\mathcal{H}}^{n-1}_x}{\big(n(n-2)\,| B_1|\big)^{\frac1{n-2}}\,\big(G(x,0)+\Psi^{(0)}(x)\big)^{\frac1{n-2}}}
\\&=&\frac1n\,\int_{\{G(\cdot,0)=t\}}
\frac{d{\mathcal{H}}^{n-1}_x}{\big(n(n-2)\,| B_1|\big)^{\frac1{n-2}}\,\big(t+\Psi^{(0)}(x)\big)^{\frac1{n-2}}}.
\end{eqnarray*}
We point out that, as~$t\to+\infty$, the set~${\{G(\cdot,0)\ge t\}}$ shrinks towards
the origin and thus, if~$x\in\{G(\cdot,0)=t\}$,
$$ \big(t+\Psi^{(0)}(x)\big)^{\frac1{2-n}}=
t^{\frac1{2-n}}\left(1+\frac{\Psi^{(0)}(x)}{t}\right)^{\frac1{2-n}}
=t^{\frac1{2-n}}\left(1+\frac1{2-n}\frac{\Psi^{(0)}(x)}{t}+O\left(\frac1{t^2}\right)\right).
$$
Therefore
\begin{equation}\label{STIMAHG-2}
\begin{split}
v(t)\,\leq\,\frac{t^{\frac1{2-n}}}{n\big(n(n-2)\,| B_1|\big)^{\frac1{n-2}}}\,\int_{\{G(\cdot,0)=t\}}
\left(1+\frac{\Psi^{(0)}(x)}{(2-n)t}+O\left(\frac1{t^2}\right)\right)
\,d{\mathcal{H}}^{n-1}_x.\end{split}\end{equation}
Now we let
$$ \rho_0:=\frac{1}{\big(n(n-2)\,| B_1|\big)^{\frac1{n-2}}}$$
and use the notation~$\e:=1/t$. For every~$(\rho,\omega,\e)\in[0,+\infty)\times(\partial B_1)\times[0,1]$,
we consider the function
\begin{equation}\label{INVEFUSER000}
F(\rho,\omega,\e):=\rho-\frac{1}{\big(n(n-2)\,| B_1|\big)^{\frac1{n-2}}\,\Big(1+\e\Psi^{(0)}
\big(\e^{\frac1{n-2}}\rho\omega\big)\Big)^{\frac1{n-2}}}
\end{equation}
and we note that~$F(\rho_0,\omega,0)=0$.
Since~$\partial_\rho F(\rho,\omega,0)=1$, by the Implicit Function Theorem, we find a
smooth function~$\rho=\rho(\omega,\e)$ such that~$\rho(\omega,\e)=\rho_0+\rho_1(\omega)\e+o(\e)$,
for some~$\rho_1:\partial B_1\to\R$,
and~$F\big(\rho(\omega,\e),\omega,\e\big)=0$, that is
\begin{equation}\label{qqwe3565477898899}
\rho(\omega,\e)=\frac{1}{\big(n(n-2)\,| B_1|\big)^{\frac1{n-2}}\,\Big(1+\e\Psi^{(0)}
\big(\e^{\frac1{n-2}}\rho(\omega,\e)\omega\big)\Big)^{\frac1{n-2}}}.\end{equation}
We also observe that
\begin{equation}\label{sw243536b657ipoiuytre}
\rho_1(\omega)=\rho_1:=-\frac{\rho_0\,\Psi^{(0)}(0)}{n-2}.
\end{equation}
Indeed, by~\eqref{qqwe3565477898899},
\begin{eqnarray*}
\rho_0+\rho_1(\omega)\e+o(\e)&=&\frac{1}{\big(n(n-2)\,| B_1|\big)^{\frac1{n-2}}\,\Big(1+\e\Psi^{(0)}
\big(\e^{\frac1{n-2}}(\rho_0+\rho_1(\omega)\e+o(\e))\omega\big)\Big)^{\frac1{n-2}}}\\&
=&\frac{\rho_0}{\Big(1+\e\Psi^{(0)}
\big(\e^{\frac1{n-2}}(\rho_0+\rho_1(\omega)\e+o(\e))\omega\big)\Big)^{\frac1{n-2}}}\\&
=&\frac{\rho_0}{\Big(1+\e\Psi^{(0)}(0)+o(\e)\Big)^{\frac1{n-2}}}\\&=&
\rho_0 \left(1-\frac{\e\,\Psi^{(0)}(0)}{n-2}+o(\e)\right),
\end{eqnarray*}
from which the desired result in~\eqref{sw243536b657ipoiuytre} follows.

Thus, we consider the map~$\partial B_1\ni\omega
\mapsto \zeta(\omega,\e):=\e^{\frac1{n-2}}\rho(\omega,\e)\omega =
\e^{\frac1{n-2}}(\rho_0+\rho_1\e+o(\e))\omega$ and we observe that, by~\eqref{INVEFUSER},
\begin{eqnarray*}
\Big(G(\zeta(\omega,\e),0)+\Psi^{(0)}
\big(\zeta(\omega,\e)\big)\Big)^{\frac1{n-2}} &=&
\frac{1}{\big(n(n-2)\,| B_1|\big)^{\frac1{n-2}}\,|\zeta(\omega,\e)|}
\\&=&\frac{1}{\big(n(n-2)\,| B_1|\big)^{\frac1{n-2}}\,
\e^{\frac1{n-2}}\rho(\omega,\e)}\\&=&
\frac{\Big(1+\e\Psi^{(0)}
\big(\e^{\frac1{n-2}}\rho(\omega,\e)\omega\big)\Big)^{\frac1{n-2}}
}{\e^{\frac1{n-2}}}\\&=&
\left(\frac1\e+\Psi^{(0)}
\big(\e^{\frac1{n-2}}\rho(\omega,\e)\omega\big)\right)^{\frac1{n-2}}
\\&=&
\left(\frac1\e+\Psi^{(0)}
\big(\zeta(\omega,\e)\big)\right)^{\frac1{n-2}}
\end{eqnarray*}
and accordingly~$G(\zeta(\omega,\e),0)=1/\e$.
This says that
\begin{equation}\label{diu4t8vb4t8tvt87t76vb8}
\{G(\cdot,0)=1/\e\}\supseteq\zeta(\partial B_1,\e).\end{equation}

We also claim that
\begin{equation}\label{pior485745867856897649}
{\mbox{if~$x\in\{G(\cdot,0)=1/\e\}$ then there exists~$\omega\in\partial B_1$
such that~$x=\e^{\frac1{n-2}}\rho(\omega,\e)\omega$.}}
\end{equation}
To prove this, we let~$x\in\{G(\cdot,0)=1/\e\}$ and we define~$\omega:=x/|x|$ and~$\rho:=|x|/\e^{\frac1{n-2}}$.
In this way, $x=\e^{\frac1{n-2}}\rho\omega$.
Hence, recalling~\eqref{INVEFUSER} and~\eqref{INVEFUSER000}, we have that
\begin{eqnarray*}&&
 F(\rho,\omega,\e)=\rho-\frac{1}{\big(n(n-2)\,| B_1|\big)^{\frac1{n-2}}\,\Big(1+\e\Psi^{(0)}
\big(\e^{\frac1{n-2}}\rho\omega\big)\Big)^{\frac1{n-2}}}\\
&&\quad=\frac{|x|}{\e^{\frac1{n-2}}}-\frac{1}{\big(n(n-2)\,| B_1|\big)^{\frac1{n-2}}\,\Big(1+\e\Psi^{(0)}
\big(\e^{\frac1{n-2}}\rho\omega\big)\Big)^{\frac1{n-2}}}\\&&\quad=
\frac{1}{\e^{\frac1{n-2}}\big(n(n-2)\,| B_1|\big)^{\frac1{n-2}}\,\big(G(x,0)+\Psi^{(0)}(x)\big)^{\frac1{n-2}}}-
\frac{1}{\big(n(n-2)\,| B_1|\big)^{\frac1{n-2}}\,\Big(1+\e\Psi^{(0)}
\big(\e^{\frac1{n-2}}\rho\omega\big)\Big)^{\frac1{n-2}}}\\&&\quad=
\frac{1}{\e^{\frac1{n-2}}\big(n(n-2)\,| B_1|\big)^{\frac1{n-2}}\,\big(1/\e+\Psi^{(0)}(x)\big)^{\frac1{n-2}}}-
\frac{1}{\big(n(n-2)\,| B_1|\big)^{\frac1{n-2}}\,\Big(1+\e\Psi^{(0)}(x)\Big)^{\frac1{n-2}}}\\&&\quad=0.
\end{eqnarray*}
As a consequence of this and the uniqueness statement in the Implicit Function Theorem, we obtain
that~$\rho=\rho(\omega,\e)$, and therefore~$x=\e^{\frac1{n-2}}\rho(\omega,\e)\omega$.
This establishes~\eqref{pior485745867856897649}.

{F}rom~\eqref{diu4t8vb4t8tvt87t76vb8} and~\eqref{pior485745867856897649},
we deduce that~$\{G(\cdot,0)=t\}=\{G(\cdot,0)=1/\e\}=\zeta(\partial B_1,\e)$,
which is a dilation of~$\partial B_1$ by a factor~$\e^{\frac1{n-2}}\big(\rho_0+\rho_1\e+o(\e)\big)=t^{\frac1{2-n}}\left(\rho_0+\frac{\rho_1}t+o\left(\frac1t\right)\right)=:\rho(t)$. Notice that~$\rho(t)\to0$ as~$t\to+\infty$.
%%small perturbation of the dilation of~$\partial B_1$ by a factor~$\e^{\frac1{n-2}}\rho_0$
%%(that is, a small perturbation of~$\partial B_{\e^{\frac1{n-2}}\rho_0}=
%%\partial B_{t^{\frac1{2-n}}\rho_0}$). Also, the size of the perturbation is of order~$\e^{\frac1{n-2}+1}$, which in turn is~$o\left(\frac1t\right)$.

We insert this information into~\eqref{STIMAHG-2} and we see that
\begin{equation}\begin{split}\label{sdwecrvbygrgy75657694597597oik}
v(t)\leq\,&\frac{t^{\frac1{2-n}}}{n\big(n(n-2)\,| B_1|\big)^{\frac1{n-2}}}\,\int_{\partial B_{\rho(t)}}
\left(1+\frac{\Psi^{(0)}(x)}{(2-n)t}+o\left(\frac1{t}\right)\right)
\,d{\mathcal{H}}^{n-1}_x\\
=\,&\frac{t^{\frac1{2-n}}\left(1+o\left(\frac1t\right)\right)}{n\big(n(n-2)\,| B_1|\big)^{\frac1{n-2}}}\left[
(\rho(t))^{n-1}{\mathcal{H}}^{n-1}(\partial B_1)+
\frac{1}{(2-n)t}\int_{\partial B_{\rho(t)}}
\Psi^{(0)}(x)\,d{\mathcal{H}}^{n-1}_x\right]\\=\,&
\frac{t^{\frac1{2-n}}(\rho(t))^{n-1}{\mathcal{H}}^{n-1}(\partial B_1)
\left(1+o\left(\frac1t\right)\right)}{n\big(n(n-2)\,| B_1|\big)^{\frac1{n-2}}}\left[
1+
\frac{1
}{(2-n)t}\fint_{\partial B_{\rho(t)}}
\Psi^{(0)}(x)\,d{\mathcal{H}}^{n-1}_x\right]\\=\,&
\frac{t^{\frac1{2-n}}(\rho(t))^{n-1}{\mathcal{H}}^{n-1}(\partial B_1)
\left(1+o\left(\frac1t\right)\right)}{n\big(n(n-2)\,| B_1|\big)^{\frac1{n-2}}}\left[
1+
\frac{1
}{(2-n)t}\big(\Psi^{(0)}(0)+o(1)\big)\right].
\end{split}\end{equation}
Now we observe that
\begin{eqnarray*}
&&t^{\frac1{2-n}}(\rho(t))^{n-1}=t^{\frac1{2-n}}\left(t^{\frac{1}{2-n}}\left(\rho_0+\frac{\rho_1}t+o\left(\frac1t\right)\right)\right)^{n-1}
=t^{\frac{n}{2-n}}\rho_0^{n-1}\left(1+\frac{\rho_1}{\rho_0\,t}+o\left(\frac1t\right)\right)^{n-1}\\&&\qquad
=t^{\frac{n}{2-n}}
\rho_0^{n-1}\left(1+\frac{(n-1)\rho_1}{\rho_0\,t}+o\left(\frac1t\right)\right).
\end{eqnarray*}
{F}rom this, \eqref{B1} and~\eqref{sdwecrvbygrgy75657694597597oik} we deduce that
\begin{equation*}\begin{split}
v(t)\leq\,&
\frac{t^{\frac{n}{2-n}}
\rho_0^{n-1}| B_1|\left(1+\frac{(n-1)\rho_1}{\rho_0\,t}+o\left(\frac1t\right)\right)
\left(1+o\left(\frac1t\right)\right)}{\big(n(n-2)\,| B_1|\big)^{\frac1{n-2}}}\left[
1+
\frac{1
}{(2-n)t}\big(\Psi^{(0)}(0)+o(1)\big)\right]\\ =\,&
t^{\frac{n}{2-n}}
\rho_0^{n}| B_1|\left(1+\frac{(n-1)\rho_1}{\rho_0\,t}+o\left(\frac1t\right)\right)
\left[1+\frac{1
}{(2-n)t}\big(\Psi^{(0)}(0)+o(1)\big)\right]\\=\,&
t^{\frac{n}{2-n}}
\rho_0^{n}| B_1|\left[1+\frac1t\left(\frac{(n-1)\rho_1}{\rho_0}+\frac{\Psi^{(0)}(0)}{2-n}+o(1)\right)\right]
.\end{split}\end{equation*}
Recalling~\eqref{sw243536b657ipoiuytre}, we observe that
$$\frac{(n-1)\rho_1}{\rho_0}+\frac{\Psi^{(0)}(0)}{2-n}=\frac{(n-1)\Psi^{(0)}(0)}{2-n}+\frac{\Psi^{(0)}(0)}{2-n}=
\frac{n\Psi^{(0)}(0)}{2-n},
$$
and therefore
\begin{equation*}\begin{split}
v(t)\leq
t^{\frac{n}{2-n}}
\rho_0^{n}| B_1|\left[1+\frac1t\left(\frac{n\Psi^{(0)}(0)}{2-n}+o(1)\right)\right].
\end{split}\end{equation*}
Accordingly,
\begin{eqnarray*}
\big(v(t)\big)^{\frac{2-n}{n}}&\le&
t\,\rho_0^{2-n}| B_1|^{\frac{2-n}{n}}\left[1+\frac1t\left(\frac{n\Psi^{(0)}(0)}{2-n}+o(1)\right)\right]^{\frac{2-n}{n}}\\&=&
t\,\rho_0^{2-n}| B_1|^{\frac{2-n}{n}}\left[1+\frac1t\left(\Psi^{(0)}(0)+o(1)\right)\right]\\&=&
t\,n(n-2) | B_1|^{\frac{2}{n}}\left[1+\frac1t\left(\Psi^{(0)}(0)+o(1)\right)\right].
\end{eqnarray*}
{F}rom this and Lemma~\ref{VDEFGRE-LE} we deduce that
\begin{equation*}
|\Omega|^{\frac{2-n}{n}}+
n(n-2)|B_1|^{\frac2{n}}t
\le t\,n(n-2) | B_1|^{\frac{2}{n}}\left[1+\frac1t\left(\Psi^{(0)}(0)+o(1)\right)\right]
.\end{equation*}
Simplifying one term, 
we conclude that
$$ |\Omega|^{\frac{2-n}{n}}\le  n(n-2) | B_1|^{\frac{2}{n}}\left(\Psi^{(0)}(0)+o(1)\right).
$$
Sending~$t\to+\infty$, we obtain the desired result.
\end{proof}

\section{The Poisson Kernel}\label{SEC:POIKER-S}

The Poisson Kernel is defined\index{Poisson Kernel}
as (minus) the normal derivative of the Green Function.
More specifically, 
given
an open set~$\Omega\subset\R^n$
with~$C^1$ boundary and a point~$x_0\in\partial\Omega$,
we let~$\Omega\ni x\mapsto G(x,x_0)$ be the Green Function of~$\Omega$,
as presented in~\eqref{FLGREEN}.
Then, for every~$x\in\Omega$
we define the Poisson Kernel of~$\Omega$ by
\begin{equation}\label{DELMDFsyhcTRhdjfmveobne6tYUStt}
P(x,x_0):=-\frac{\partial G}{\partial\nu}(x,x_0),
\end{equation}
where the notation of~\eqref{OKJN-NOAYHNSTAZ} has been used.

We observe that if~$x\ne x_0$,
\begin{equation}\label{PfghjklX99ERHILD}
\begin{split}&
-\Delta_x P(x,x_0)=-\Delta_x\left(
\nabla_yG(x,y)\cdot\nu(y)
\right)
=-\Delta_x\left(
\nabla_y\Gamma(x-y)\cdot\nu(y)-
\nabla_y\Psi^{(y)}(x)\cdot\nu(y)
\right)\\&\qquad=-
\nabla_y-\Delta_x\Gamma(x-y)\cdot\nu(y)+
\nabla_y-\Delta_x\Psi^{(y)}(x)\cdot\nu(y)
=0,\end{split}
\end{equation}
thanks to~\eqref{L:S:ASM S}
and~\eqref{FLGREEN}.

The importance of this kernel\footnote{As a byproduct of the positivity of the Green function (see footnote~\ref{POSI:GREEN:Ut} on page~\pageref{POSI:GREEN:Ut}) and of the fact that it vanishes along~$\partial\Omega$, we have that the Poisson Kernel is nonnegative (and, in fact, strictly positive, by the Strong Maximum Principle) in a bounded set~$\Omega$.} lies in the following
representation result for harmonic functions with prescribed boundary values:

\begin{theorem}\label{POIBALL1}
If~$\Omega\subseteq\R^n$ is an open and bounded set with~$C^1$
boundary, and~$u\in C^2(\Omega)\cap C(\overline\Omega)$
solves
\begin{equation}\label{kjhgfdKMSD:OSKMIDEGBD1}\begin{dcases}
\Delta u =0& {\mbox{ in }}\Omega,\cr
u=g& {\mbox{ on }}\partial\Omega,\end{dcases}
\end{equation}
for some~$g:\partial\Omega\to\R$,
then, for every~$x\in\Omega$,
\begin{equation}\label{kjhgfdKMSD:OSKMIDEGBD2} u(x)=\int_{\partial \Omega}
g(y)\,P(x,y)\,d{\mathcal{H}}^{n-1}_y.\end{equation}
Moreover,
\begin{equation}\label{L903:9455:Sw44}
\int_{\partial\Omega} P(x,y)\,d{\mathcal{H}}_y^{n-1}=1.\end{equation}

Furthermore, if~$u$ is defined as in~\eqref{kjhgfdKMSD:OSKMIDEGBD2} and~$g$ is continuous,
then~$u\in C^2(\Omega)\cap C(\overline{\Omega})$ and~$u$ satisfies~\eqref{kjhgfdKMSD:OSKMIDEGBD1}.
\end{theorem}

\begin{proof} When~$u\in C^1(\overline{\Omega})$ this is in fact a particular case of
Corollary~\ref{SKMD:DOK0oerlfgdit}, used here with~$f:=0$
(and, en passant, this observation proves also~\eqref{L903:9455:Sw44},
since the function~$u$ identically equal to~$1$ solves~\eqref{kjhgfdKMSD:OSKMIDEGBD1}
with~$g:=1$,
hence~\eqref{L903:9455:Sw44}
follows from~\eqref{kjhgfdKMSD:OSKMIDEGBD2} applied to a function as regular as we wish).

If instead one wants to prove~\eqref{kjhgfdKMSD:OSKMIDEGBD2}
by only assuming that~$u\in C(\overline\Omega)$, one can reduce to
the case~$u\in C^1(\overline{\Omega})$ via a mollification argument. For this,
we take~$\tau\in C^\infty_0(B_1,\,[0,+\infty))$
with~$\int_{B_1}\tau(x)\,dx=1$.
Given~$\eta>0$, we let~$\tau_\eta(x):=\frac1{\eta^n}
\tau\left(\frac{x}{\eta}\right)$
and define~$u_\eta:=u*\tau_\eta$ (recall e.g.~\cite[Theorem 9.8]{MR3381284}
for the uniform approximation properties of~$u_\eta$ to~$u$ when~$\eta\searrow0$). We also consider~$\Omega_\eta$
to be the points~$p$ in~$\Omega$ for which~$B_{2\eta}(p)\Subset\Omega$.
Then, $u_\eta\in C^2(\overline{\Omega_\eta})$ and therefore,
setting~$g_\eta(x):=u_\eta(x)$ for all~$x\in\partial\Omega_\eta$,
we can write that, for all~$x\in\Omega_\eta$,
\begin{equation*} u_\eta(x)=\int_{\partial \Omega_\eta}
g_\eta(y)\,P(x,y)\,d{\mathcal{H}}^{n-1}_y.\end{equation*}
Thus, for every~$x\in\Omega$,
\begin{equation*} \begin{split}&
u(x)
-\int_{\partial \Omega}
g(y)\,P(x,y)\,d{\mathcal{H}}^{n-1}_y
=\lim_{\eta\searrow0}
u_\eta(x)
-\int_{\partial \Omega}
u(y)\,P(x,y)\,d{\mathcal{H}}^{n-1}_y
\\&\qquad=\lim_{\eta\searrow0}\int_{\partial \Omega_\eta}
g_\eta(y)\,P(x,y)\,d{\mathcal{H}}^{n-1}_y
-\int_{\partial \Omega_\eta}
u(y)\,P(x,y)\,d{\mathcal{H}}^{n-1}_y\\&\qquad
=
\lim_{\eta\searrow0}\int_{\partial \Omega_\eta}
\big(u_\eta(y)-u(y)\big)\,P(x,y)\,d{\mathcal{H}}^{n-1}_y
=0,\end{split}
\end{equation*}
thanks to~\eqref{L903:9455:Sw44},
and this establishes~\eqref{kjhgfdKMSD:OSKMIDEGBD2}
in the claimed generality.

We now prove that if~$u$ is defined as in~\eqref{kjhgfdKMSD:OSKMIDEGBD2} and~$g$ is continuous,
then~$u\in C^2(\Omega)\cap C(\overline{\Omega})$ and~$u$ satisfies~\eqref{kjhgfdKMSD:OSKMIDEGBD1}.
We start by noticing that~$u$ is harmonic
since so is the Poisson Kernel (since so are the Green Function
and the fundamental solution outside its singularity, recall~\eqref{FLGREEN}
and~\eqref{DELMDFsyhcTRhdjfmveobne6tYUStt}). Thus, we focus
on proving that~$u$ attains the boundary datum~$g$ continuously.
For this, the gist is to exploit the fact that the Poisson Kernel has a single singularity and integrates to~$1$. The details of the proof go as follows. Let~$x_0\in\partial\Omega$ and~$\e>0$. We use the continuity of~$g$ to pick a~$\delta>0$ such that~$|g(x)-g(x_0)|\le\e$ for each~$x\in(\partial\Omega)\cap B_\delta(x_0)$. We also take a sequence~$x_k\in\Omega$
converging to~$x_0$ as~$k\to+\infty$. Without loss of generality, we can suppose that~$x_k\in B_\delta(x_0)$.
Then, by~\eqref{kjhgfdKMSD:OSKMIDEGBD2} and~\eqref{L903:9455:Sw44},
\begin{eqnarray*}
|u(x_k)-g(x_0)|&=&\left|\int_{\partial \Omega}
g(y)\,P(x_k,y)\,d{\mathcal{H}}^{n-1}_y-g(x_0)\right|\\
&=&\left|\int_{\partial \Omega}
\big(g(y)-g(x_0)\big)\,P(x_k,y)\,d{\mathcal{H}}^{n-1}_y\right|
\\&\le&\e\int_{(\partial \Omega)\cap B_\delta(x_0)}
P(x_k,y)\,d{\mathcal{H}}^{n-1}_y+
\int_{(\partial \Omega)\setminus B_\delta(x_0)}
\big|g(y)-g(x_0)\big|\,P(x_k,y)\,d{\mathcal{H}}^{n-1}_y\\&\le&\e+
\int_{(\partial \Omega)\setminus B_\delta(x_0)}
\big|g(y)-g(x_0)\big|\,P(x_k,y)\,d{\mathcal{H}}^{n-1}_y.
\end{eqnarray*}
Sending~$k\to+\infty$ we thereby deduce that
$$\limsup_{k\to+\infty}|u(x_k)-g(x_0)|\le\e.$$
Hence, taking now~$\e$ as small as we wish,
$$\lim_{k\to+\infty}|u(x_k)-g(x_0)|=0,$$
as desired.
\end{proof}

Since, in light of~\eqref{COIENEE} and~\eqref{GAMMAFU},
we know that~$\nabla\Gamma(x)=-\frac{x}{n\, | B_1|\,|x|^n}$,
we straightforwardly deduce from Theorems~\ref{MS:SKMD344567yjghS-1}
and~\ref{MS:SKMD344567yjghS-2} the explicit expression
of the Poisson Kernels of the ball and of the halfspace (see Figure~\ref{POS34S45S345EEBsoloDItangeFI}
for a plot of the case of the unit ball and Figure~\ref{POS34S45S345EEBsoloDItangeF67894g95696569I}
for the case of the halfspace):

\begin{figure}
  \centering
  \includegraphics[width=.4\linewidth]{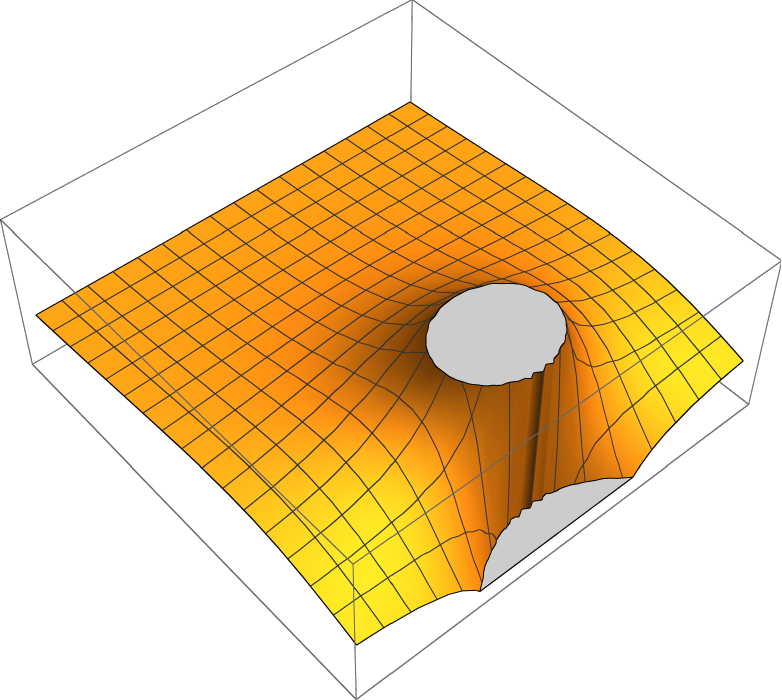}$\qquad$
  \includegraphics[width=.4\linewidth]{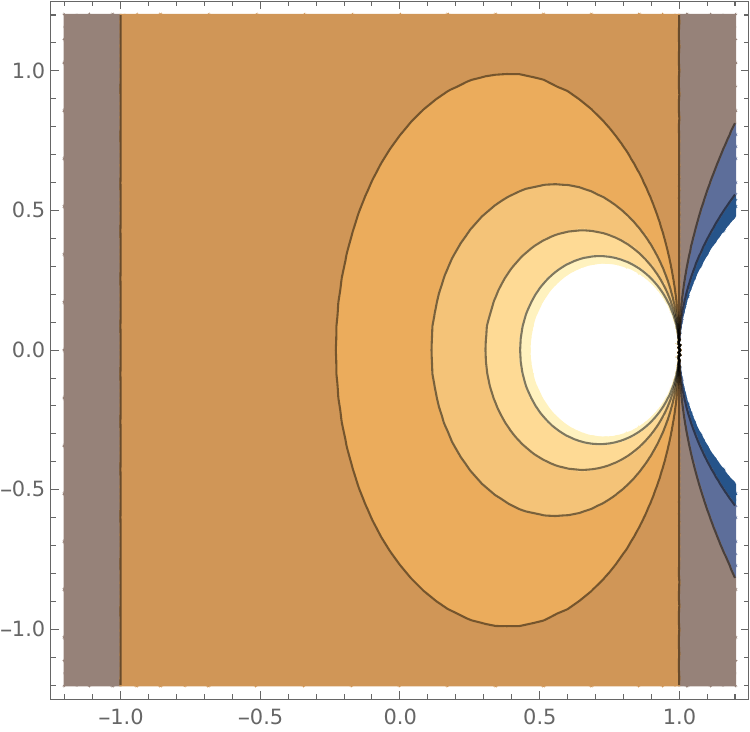}
 \caption{\sl The Poisson Kernel for the unit ball and its level sets.}\label{POS34S45S345EEBsoloDItangeFI}
\end{figure}

\begin{theorem} \label{POIBALL}
The Poisson Kernel
of the ball~$B_R$ is\index{Poisson Kernel}
\begin{equation*}P(x,x_0)=\frac{R^2-|x|^2}{n\,|B_1|\,R\,|x-x_0|^n},
\end{equation*}
for every~$x\in B_R$ and~$x_0\in\partial B_R$. 
\end{theorem}

\begin{proof}
We recall the notation in~\eqref{DELMDFsyhcTRhdjfmveobne6tYUStt}
and~\eqref{OKJN-NOAYHNSTAZ}.
Thus, from Theorem~\ref{MS:SKMD344567yjghS-1}, \eqref{DELMDFsyhcTRhdjfmveobne6tYUStt}
and~\eqref{UJdgvbf fficpqlwdfdvbia0comPucjv}, we have that,
if~$x\in B_R\setminus\{0\}$,
\begin{eqnarray*}
&&-P(x,x_0)=\nabla_{x_0} G(x,x_0)\cdot\frac{x_0}R=
\nabla_{x_0} \big(G(x_0,x)\big)\cdot\frac{x_0}R\\&&\qquad
=
\nabla_{x_0}\left[
\Gamma(x_0-x)-\Gamma\left(\frac{|x|}R\,\left(
x_0-\frac{R^2x}{|x|^2}
\right)\right)
\right]\cdot\frac{x_0}R\\&&\qquad=-
\frac{x_0}{n\,|B_1|\,R}\cdot\left(\frac{x_0-x}{|x_0-x|^n}-
\frac{\frac{|x|^2}{R^2}\,\left(
x_0-\frac{R^2x }{|x |^2}
\right)}{\left|\frac{|x |}R\,\left(
x_0-\frac{R^2x }{|x |^2}
\right)\right|^n}
\right)\\&&\qquad=-
\frac{x_0}{n\,|B_1|\,R\,|x-x_0|^n}\cdot\left(x_0-x-\frac{|x|^2}{R^2}\,\left(
x_0-\frac{R^2x }{|x|^2}\right)\right)\\&&\qquad=
-
\frac{x_0}{n\,|B_1|\,R\,|x-x_0|^n}\cdot\left(x_0-\frac{|x|^2x_0}{R^2}\right)
\\&&\qquad=
-
\frac{|x_0|^2}{n\,|B_1|\,R\,|x-x_0|^n}\left(1-\frac{|x|^2}{R^2}\right)\\&&\qquad=
-
\frac{R^2-|x|^2}{n\,|B_1|\,R\,|x-x_0|^n},
\end{eqnarray*}
which is the desired result when~$x\in B_R\setminus~\{0\}$. 

If instead~$x=0$, then Theorem~\ref{MS:SKMD344567yjghS-1}
and~\eqref{DELMDFsyhcTRhdjfmveobne6tYUStt} lead to
\begin{eqnarray*}&&
-P(0,x_0)=\nabla_{x_0} G(0,x_0)\cdot\frac{x_0}R=
\nabla_{x_0} \big(G(x_0,0)\big)\cdot\frac{x_0}R
=\nabla_{x_0}
\big( \Gamma(x_0)-\Gamma(R)\big)\cdot\frac{x_0}R\\&&\qquad=
-\frac{x_0}{n\, | B_1|\,|x_0|^n}\cdot\frac{x_0}R=-\frac{R^2}{n\, | B_1|\,R\,|x_0|^n},
\end{eqnarray*}
which completes the proof of Theorem~\ref{POIBALL}.
\end{proof}

\begin{figure}
  \centering
  \includegraphics[width=.4\linewidth]{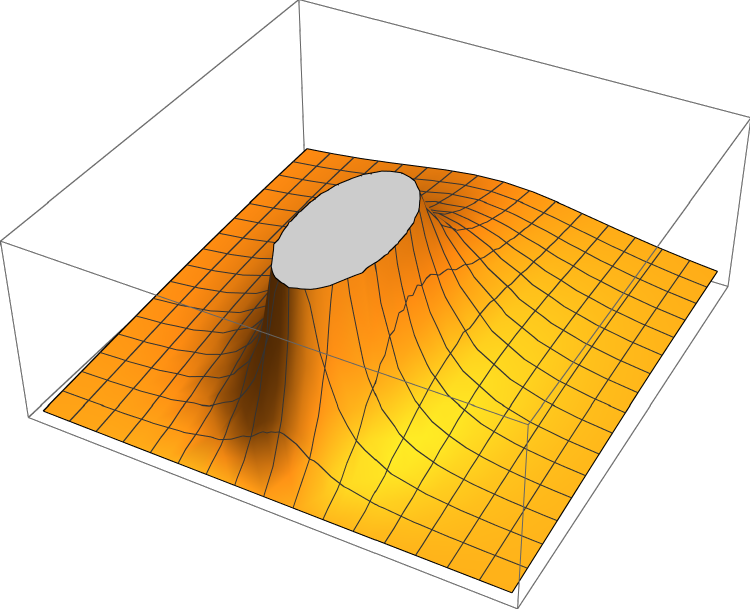}$\qquad$
  \includegraphics[width=.4\linewidth]{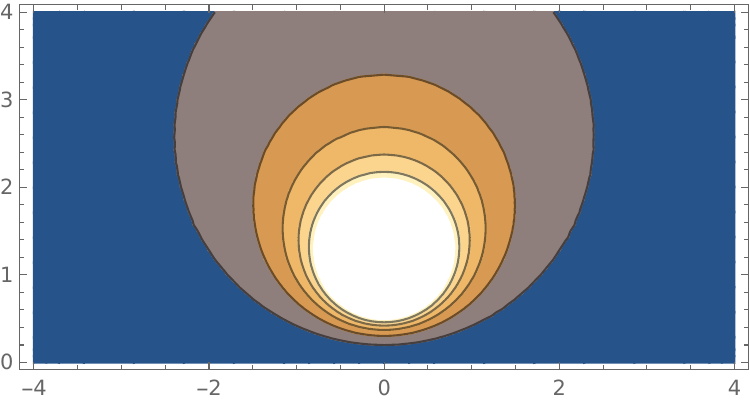}
 \caption{\sl The Poisson Kernel for the halfplane and its level sets.}\label{POS34S45S345EEBsoloDItangeF67894g95696569I}
\end{figure}

\begin{theorem} A Poisson Kernel
of the halfspace~$\{x_n>0\}$ is \begin{equation}\label{PO3456M4567Pt67LAMMNO-1}
P(x,x_0)=\frac2{n|B_1|}\,\frac{x \cdot e_n}{|x-x_0|^n}
,\end{equation}
for every~$x \in \{x_n>0\}$ and~$x_0\in \{x_n=0\}$. 

Also, if~$u\in C^2(\{x_n>0\})$ is a bounded\index{Poisson Kernel}
and uniformly continuous
solution
of
\begin{equation}\label{PO3456M4567Pt67LAMMNO-100}
\begin{dcases}
\Delta u =0& {\mbox{ in }}\R^{n-1}\times(0,+\infty),\cr
u=g& {\mbox{ on }}\R^{n-1}\times\{0\},\end{dcases}
\end{equation}
for some~$g:\partial\Omega\to\R$,
then, for every~$x=(x',x_n)\in\R^{n-1}\times(0,+\infty)$,
\begin{equation}\label{PO3456M4567Pt67LAMMNO-2} u(x)=\frac2{n|B_1|}\,
\int_{\R^{n-1}}\frac{x_n\,
g(y)}{\big(|x'-y|^2+x_n^2\big)^{\frac{n}2}}\,\,dy.\end{equation}
\end{theorem}

\begin{proof} 
{F}rom Theorem~\ref{MS:SKMD344567yjghS-2},
for every~$x_0$ with~$x_0\cdot e_n=0$,
\begin{equation*}\begin{split}&
P(x,x_0)=\nabla_{x_0}\big(\Gamma(x-x_0)-\Gamma({\mathcal{R}}(x)-x_0)\big)
\cdot e_n=\big(
-\nabla\Gamma(x-x_0)+\nabla\Gamma({\mathcal{R}}(x)-x_0)\big)\cdot e_n
\\&\qquad= 
\frac1{n\, | B_1|\,|x-x_0|^n}
\big(
x-x_0-({\mathcal{R}}(x)-x_0)\big)\cdot e_n
=\frac{2x\cdot e_n}{n\, | B_1|\,|x-x_0|^n}\end{split}
\end{equation*}
and this proves~\eqref{PO3456M4567Pt67LAMMNO-1}.

Now we prove~\eqref{PO3456M4567Pt67LAMMNO-2}.
To this end, we first check that the function introduced in~\eqref{PO3456M4567Pt67LAMMNO-2}
is indeed a solution of~\eqref{PO3456M4567Pt67LAMMNO-100}.
For this,
we
take~$x=(x',x_n)\in\R^{n-1}\times\R$
and we define
$$\Theta(x):=\frac{x_n}{\big(|x'|^2+x_n^2\big)^{\frac{n}2}}=\frac{x_n}{|x|^n}.
$$ We observe that, by the expression
of the Laplace operator in spherical coordinates~\eqref{ROTSE},
\begin{equation}\label{LISNtGABSanisoldvISS4PKM}
\begin{split}&
\Delta\Theta(x)=\Delta (x_n\,|x|^{-n})
=x_n\,\Delta|x|^{-n}+2\nabla x_n\cdot\nabla |x|^{-n}
\\&\qquad= x_n\,\big(
n(n+1)|x|^{-n-2}
-n(n-1)|x|^{-n-2}
\big)-2n|x|^{-n-2} x_n=0.\end{split}
\end{equation}
Furthermore, for each~$i\in\{1,\dots,n\}$ and~$h\in\left(0, \frac{x_n}2\right)$,
\begin{eqnarray*}&& |\Theta(x+he_i)+\Theta(x-he_i)-2\Theta(x)|\le
\sup_{\xi\in B_h(x)} |D^2\Theta(\xi)|\le
C\,\sup_{\xi=(\xi',\xi_n)\in B_h(x)} \frac{\xi_n}{|\xi|^{n+2}}\\&&\quad\quad
\le
\frac{C\,x_n}{\big((|x'|-h)^2+x_n^2\big)^{\frac{n+2}2}}\le
\begin{dcases}
\frac{C}{ x_n^{n+1}} & {\mbox{ if }} |x'|\le 8x_n,\cr
\frac{C\,x_n}{ 
\big(|x'|^2+x_n^2\big)^{\frac{n+2}2}} & {\mbox{ if }} |x'|\ge 8x_n,
\end{dcases}\\&&\quad\quad
\le\frac{C\,x_n}{ 
\big(|x'|^2+x_n^2\big)^{\frac{n+2}2}}
\end{eqnarray*}
for some constant~$C>0$, possibly varying from line to line.

Consequently, since,
for a given~$x_n>0$, the function~$\R^{n-1}\ni x'\mapsto\frac{x_n}{ 
\big(|x'|^2+x_n^2\big)^{\frac{n+2}2}}$ belongs to~$L^1(\R^{n-1})$,
recalling that~$g$ is bounded, we deduce from the Dominated Convergence Theorem and~\eqref{LISNtGABSanisoldvISS4PKM} that
\begin{eqnarray*}&& \Delta\left(
\int_{\R^{n-1}}
g(y) \Theta(x'-y,x_n)\,\,dy\right)\\&=&
\lim_{h\searrow0}\frac1{h^2}\,\sum_{i=1}^{n}\int_{\R^{n-1}}
g(y) \Big(\Theta\big( (x'-y,x_n)+he_i\big)
+\Theta\big( (x'-y,x_n)-he_i\big)-2\Theta (x'-y,x_n)
\Big)\,\,dy\\&=&
\lim_{h\searrow0}\frac1{h^2}\,\sum_{i=1}^{n}\int_{\R^{n-1}}
g(x'-\zeta) \Big(\Theta\big( (\zeta,x_n)+he_i\big)
+\Theta\big( (\zeta,x_n)-he_i\big)-2\Theta (\zeta,x_n)
\Big)\,\,d\zeta\\&=&
\sum_{i=1}^{n}\int_{\R^{n-1}}
g(x'-\zeta) \frac{\partial^2 \Theta}{\partial x_i^2}(\zeta,x_n)\,\,d\zeta\\&=&
\int_{\R^{n-1}}
g(x'-\zeta) \Delta \Theta(\zeta,x_n)\,\,d\zeta\\&=&0.
\end{eqnarray*}
This shows that~$u$ satisfies the differential equation in~\eqref{PO3456M4567Pt67LAMMNO-100}.

Besides, using the notation of the Euler Gamma Function,
recalling~\eqref{MULTOLPALPA280kdsmc8724yrth0923}, 
$$
\Gamma\left({\frac{n}2}\right)
=\frac{2\pi^{\frac{n}2} }{{\mathcal{H}}^{n-1}(\partial B_1)}
=\frac{2\pi^{\frac{n}2}}{n\,| B_1|},$$
hence, if~$B'_1$ is the unit ball in dimension~$n-1$,
$$ \frac{{\mathcal{H}}^{n-2}(\partial B_1')}{n\,| B_1|}=
\frac{{\mathcal{H}}^{n-2}(\partial B_1')\,\Gamma\left({\frac{n}2}\right)}{2\pi^{\frac{n}2}}=
\frac{2\pi^{\frac{n-1}2}\,\Gamma\left({\frac{n}2}\right)}{2\pi^{\frac{n}2}\,\Gamma\left({\frac{n-1}2}\right)}=
\frac{\Gamma\left({\frac{n}2}\right)}{\sqrt\pi\,\Gamma\left({\frac{n-1}2}\right)}.
$$
Therefore, 
using the substitution~$t:=\frac{x'-y}{\eta}$
and the Dominated Convergence Theorem, for each~$x'\in\R^{n-1}$,
\begin{equation}\label{bYNUUTFn7nujv5tv6g9g3}\begin{split}&
\lim_{\eta\searrow0} u(x',\eta)=\lim_{\eta\searrow0}\frac2{n|B_1|}\,
\int_{\R^{n-1}}\frac{\eta\,
g(y)}{\big(|x'-y|^2+\eta^2\big)^{\frac{n}2}}\,\,dy\\
&\qquad=\lim_{\eta\searrow0} \frac2{n|B_1|}\,
\int_{\R^{n-1}}\frac{
g(x'-\eta t)}{\big( |t|^2+1\big)^{\frac{n}2}}\,\,dt=\frac2{n|B_1|}\,
g(x')\,\int_{\R^{n-1}}\frac{dt}{\big( |t|^2+1\big)^{\frac{n}2}}\\&\qquad
=\frac{2{\mathcal{H}}^{n-2}(\partial B_1')}{n|B_1|}\,
g(x')\,\int_{0}^{+\infty}\frac{\rho^{n-2}\,d\rho}{\big( \rho^2+1\big)^{\frac{n}2}}=
\frac{2\Gamma\left({\frac{n}2}\right)}{\sqrt\pi\,\Gamma\left({\frac{n-1}2}\right)}\,
g(x')\,\int_{0}^{+\infty}\frac{\rho^{n-2}\,d\rho}{\big( \rho^2+1\big)^{\frac{n}2}}
.
\end{split}\end{equation}

Now we claim that, for all~$n\in\N$, $n\ge2$, it holds that
\begin{equation}\label{OL-RICOLRHO}
\int_{0}^{+\infty}\frac{\rho^{n-2}\,d\rho}{\big( \rho^2+1\big)^{\frac{n}2}}=
\frac{\sqrt\pi\,\Gamma\left({\frac{n-1}2}\right)}{2\Gamma\left({\frac{n}2}\right)}.
\end{equation}
To prove this, we denote by~${\mathcal{J}}_n$ the left hand side of~\eqref{OL-RICOLRHO}
and we observe that
$$ \frac{d}{d\rho}\frac1{\big( \rho^2+1\big)^{\frac{n}2}}
=-\frac{n\rho}{\big( \rho^2+1\big)^{\frac{n+2}2}}
$$
and consequently, integrating by parts,
\begin{equation}\label{0OKMAMM0okm0okm00kKJSqeweIJ}
\begin{split}& {\mathcal{J}}_{n+2}=
\int_{0}^{+\infty}\frac{\rho^{n}\,d\rho}{\big( \rho^2+1\big)^{\frac{n+2}2}}=
-\frac1n\,\int_{0}^{+\infty}
\rho^{n-1}\,\left(\frac{d}{d\rho}\frac1{\big( \rho^2+1\big)^{\frac{n}2}}\right)
\,d\rho\\&\qquad\qquad=
\frac{n-1}n\,\int_{0}^{+\infty}
\frac{\rho^{n-2}}{\big( \rho^2+1\big)^{\frac{n}2}}
\,d\rho=
\frac{n-1}n\,
{\mathcal{J}}_n.
\end{split}\end{equation}
In addition,
\begin{equation}\label{ZKMdgamzo34:1}
{\mathcal{J}}_2=
\int_{0}^{+\infty}\frac{ d\rho}{ \rho^2+1 }=\frac\pi2,
\end{equation}
and, using the substitution~$r:=\rho^2+1$,
\begin{equation}\label{ZKMdgamzo34:2}
{\mathcal{J}}_3=
\int_{0}^{+\infty}\frac{\rho\,d\rho}{\big( \rho^2+1\big)^{\frac{3}2}}=\frac12\,\int_{1}^{+\infty}\frac{dr}{r^{\frac{3}2}}=1.
\end{equation}
We stress that~\eqref{ZKMdgamzo34:1} is precisely~\eqref{OL-RICOLRHO}
when~$n=2$, in light of~\eqref{ZKMdgamzo34}.
Similarly we have that~\eqref{ZKMdgamzo34:2} is precisely~\eqref{OL-RICOLRHO}
when~$n=3$, owing to~\eqref{GAMMA3mez}.

Therefore, to complete the proof of~\eqref{OL-RICOLRHO}, we argue by induction.
Having established the basis of
the induction for~$n\in\{2,3\}$ in~\eqref{ZKMdgamzo34:1} and~\eqref{ZKMdgamzo34:2},
we can suppose that~\eqref{OL-RICOLRHO} holds true for all~$n\le N$ with~$N\in\N\cap[3,+\infty)$
and we prove it for~$N+1$.
To this end, we make use of~\eqref{0OKMAMM0okm0okm00kKJSqeweIJ}, we recall the basic
property of the Euler Gamma Function
in~\eqref{GAMMAPIU1} and we gather that
$$ {\mathcal{J}}_{N+1}=
\frac{N-2}{N-1}\,
{\mathcal{J}}_{N-1}
=\frac{N-2}{N-1}\,
\frac{\sqrt\pi\,\Gamma\left({\frac{N-2}2}\right)}{2\Gamma\left({\frac{N-1}2}\right)}
=\frac{\frac{N-2}2}{\frac{N-1}2}\,
\frac{\sqrt\pi\,\Gamma\left({\frac{N-2}2}\right)}{2\Gamma\left({\frac{N-1}2}\right)}=
\frac{\sqrt\pi\,\Gamma\left({\frac{N}2}\right)}{2\Gamma\left({\frac{N+1}2}\right)}
$$
and this completes the proof of~\eqref{OL-RICOLRHO}.

Now, from~\eqref{bYNUUTFn7nujv5tv6g9g3}
and~\eqref{OL-RICOLRHO}, we conclude that
$$\lim_{\eta\searrow0} u(x',\eta)=g(x'),$$
which checks that~$u$ satisfies the boundary condition in~\eqref{PO3456M4567Pt67LAMMNO-100}.

Additionally, we remark that the solution to~\eqref{PO3456M4567Pt67LAMMNO-100}
in the class of bounded 
and uniformly continuous
functions is unique,
thanks to the Maximum Principle in Proposition~\ref{UNDIFKFCOKDMTFDEBSF},
and this establishes~\eqref{PO3456M4567Pt67LAMMNO-2},
as desired.
\end{proof}

We refer to~\cite{MR1272823} for several beautiful geometric observations
on the Poisson Kernels of the ball and of the halfplane, also in connection with conformal and hyperbolic geometries.
See also~\cite{MR1158146, MR1446490} for interesting geometric connections with complex analysis.\medskip

An additional
interesting physical motivation of the Poisson Kernel can be provided in terms of the ``double layer'' potential\index{double layer potential}, that is the electrostatic\index{electrostatics} potential associated with a dipole\index{dipole} distribution around the smooth
surface~$\Sigma:=\partial\Omega$.
Roughly speaking, if we place some positive (respectively, negative) charges in the vicinity of the external (respectively, internal) side of~$\Sigma$, these charges will almost ``cancel'' each other but the first nonnegligible order can be related directly to the Poisson Kernel, up to a harmonic correction
(a precise statement will be given in equation~\eqref{DOUPOSHAR} below).
To make this idea more explicit, 
let~$g$ be a smooth function defined on~$\Sigma$ and~$\nu$ be the external normal along~$\Sigma$. Given~$\e>0$,
we endow the surface~$\Sigma_+:=\{y+\e\nu(y),\;y\in\Sigma\}$ with a charge density modeled on~$g$,
and the surface~$\Sigma_-:=\{y-\e\nu(y),\;y\in\Sigma\}$ with a charge density modeled on~$-g$.
In view of the comment right after~\eqref{GAMMAFU}, for each~$x\in\Omega$, the electrostatic potential obtained in this way can be written as
\begin{equation}\label{KSM-LAPSOOTENXKDZIASD03}
u_\e(x)\,:=\,\int_{\Sigma} g(y)\,\Gamma\big(x-(y+\e\nu(y))\big)\,d{\mathcal{H}}^{n-1}_y
-\int_{\Sigma} g(y)\,\Gamma\big(x-(y-\e\nu(y))\big)\,d{\mathcal{H}}^{n-1}_y
.\end{equation}
We observe that
\begin{eqnarray*}
\int_{\Sigma} g(y)\,\Gamma\big(x-(y\pm\e\nu(y))\big)\,d{\mathcal{H}}^{n-1}_y
=
\int_{\Sigma} g(y)\,\Big[\Gamma(x-y)
\mp\e\nabla\Gamma (x-y)\cdot\nu(y)\Big]\,d{\mathcal{H}}^{n-1}_y
+o(\e).
\end{eqnarray*}
This and~\eqref{KSM-LAPSOOTENXKDZIASD03} lead to
$$ u_\e(x)=-
2\e\int_{\Sigma} g(y)\nabla\Gamma (x-y)\cdot\nu(y)\,d{\mathcal{H}}^{n-1}_y
+o(\e).$$
Therefore, recalling the definition
of the Robin Function in~\eqref{L:S:ASM S}, that of the Green Function in~\eqref{FLGREEN}, and that of
the Poisson Kernel in~\eqref{DELMDFsyhcTRhdjfmveobne6tYUStt}, and making use of
the notation in~\eqref{OKJN-NOAYHNSTAZ} and of the symmetry of the fundamental solution,
we see that (up to the factor~$2$ that we omit) the kernel related to the double layer potential has the form
\begin{equation}\label{DOUPOSHAR}\begin{split}
-\nabla \Gamma(x-y)\cdot\nu(y)\,&=\,\nabla_y \Gamma(x-y)\cdot\nu(y)
\\&=\,\nabla_y G(x,y)\cdot\nu(y)-\nabla_y\Psi^{(y)}(x)\cdot\nu(y)
\\&=\,\frac{\partial G}{\partial\nu}(x,y)-\nabla_y\Psi^{(y)}(x)\cdot\nu(y)
\\&=\,-P(x,y)-\nabla_y\Psi^{(y)}(x)\cdot\nu(y).\end{split}
\end{equation}
Since the Robin Function is harmonic, this equation expresses the Poisson Kernel
as a suitable harmonic correction of the double layer potential.
It is indeed a classical approach to use the double layer potential to solve boundary value problems for elliptic equations, see e.g.~\cite[Chapter~III]{MR1306729}.\medskip

As an additional remark, we point out that the Poisson Kernel
provides a useful existence theory (with somewhat explicit solutions in some cases)
for harmonic functions with given boundary data. This allows us to complete
the proof of Corollary~\ref{2432rf4g-2rfjv-1k2ermfAKSd56yu-124rktgmMqe22rtS-23t4yhtj}
by arguing as follows.

\begin{proof}[Proof of Corollary~\ref{2432rf4g-2rfjv-1k2ermfAKSd56yu-124rktgmMqe22rtS-23t4yhtj}] We prove that condition~(i) is equivalent to condition~(ii) \label{2432rf4g-2rfjv-1k2ermfAKSd56yu-124rktgmMqe22rtS-23t4yhtj.LAPA}
(similarly, one can prove the equivalence between conditions~(i) and~(iii)).
On the one hand, condition~(i) here coincides with condition~(i) in Corollary~\ref{KAHAR-OVUN}
which entails condition~(ii) in Corollary~\ref{KAHAR-OVUN}, which in turn entails
condition~(ii) here.

On the other hand, if condition~(ii) here holds,
then also condition~(iii) holds (by polar coordinates), so it suffices to check that condition~(iii) implies condition~(i). To this purpose, we pick~$x_0\in\Omega$ and~$r_0>0$ such that~$B_{r_0}(x_0)\Subset\Omega$, we consider
(owing to Theorems~\ref{POIBALL1} and~\ref{POIBALL})
the harmonic function~$v\in C^2(B_{r_0}(x_0))\cap C(\overline{B_{r_0}(x_0)})$ such that~$v=u$ along~$\partial B_{r_0}(x_0)$ and we claim that
\begin{equation}\label{RFSD:REPLACEMENT}
{\mbox{$v=u$ in~$B_{r_0}(x_0)$.}}\end{equation}
Indeed, if not, up to a translation we can assume that~$x_0=0$, we can define~$w:=u-v$ and suppose, without loss of generality, that the set~$S:=\{x\in B_{r_0}$ s.t. $w>0\}$ is not empty.

As a consequence, since~$w$ is continuous and vanishes along~$(\partial B_{r_0})\cap(\partial S)$, we can find~$x_\star\in S\subseteq B_{r_0}$ such that
$$ \max_{\overline{B_r}}w=\max_{\overline S} w= w(x_\star)>0.$$
Among all the possible points with this maximization property, again by the continuity of~$w$, we can pick the one with biggest norm, namely, up to renaming~$x_\star$, we can assume that
\begin{equation}\label{9ojdfvd9fe:9qwyfihgb0293NHTR4t}
{\mbox{$w(x_\star)\ge w(x)$ for every~$x\in B_r$ and that~$w(x_\star)> w(x)$ for every~$x\in B_r$ with~$|x|>|x_\star|$.}}\end{equation}

Moreover, by~(iii) and the Mean Value Formula for~$v$, we can find~$r_\star>0$ small enough such that
$$ w(x_\star)=\fint_{B_{r_\star}(x_\star)} w(x)\,dx.$$
But, in light of~\eqref{9ojdfvd9fe:9qwyfihgb0293NHTR4t},
$$ \fint_{B_{r_\star}(x_\star)} w(x)\,dx>w(x_\star),$$
providing a contradiction.

This proves~\eqref{RFSD:REPLACEMENT}, which in turn, since~$v$ is harmonic, yields~(i).
\end{proof}

\section{Point singularities for harmonic functions}

A natural question is what types of blowup's harmonic functions can
exhibit at a given point. For instance, the fundamental solution in~\eqref{GAMMAFU} provides
an example of a harmonic function in~$B_1\setminus\{0\}$ which,
when~$n\ge2$, diverges at the origin. It turns out that milder blowup's
than the one offered by the fundamental solution are not allowed,
as pointed out in the following result:

\begin{theorem}\label{0oiujhgv098j-23rtgrJSPsdfgdh56uh}
Let~$n\ge2$, $r>0$ and suppose\footnote{As a matter of fact, Theorem~\ref{0oiujhgv098j-23rtgrJSPsdfgdh56uh} remains valid
when~$n=1$, since in this case~$u$ vanishes identically. Indeed, when~$n=1$, one
deduces from~\eqref{GAMMAFU} and~\eqref{MDIVINZ} that
$$ \lim_{x\to0}\frac{u(x)}{|x|}=0$$
and thus, since in this situation, $u$ is affine in~$(-r,0)$ and~$(0,r)$,
necessarily~$u(x)=0$ for all~$x\in(-r,0)\cup(0,r)$.} that~$u$ is harmonic in~$B_r\setminus\{0\}$.
Assume also that
\begin{equation}\label{MDIVINZ}
\lim_{x\to0}\frac{u(x)}{\Gamma(x)}=0,
\end{equation}
where~$\Gamma$ is the fundamental solution in~\eqref{GAMMAFU}.
Then, $u$ can be extended to a harmonic function in the whole of~$B_r$.
\end{theorem}

\begin{proof}
Let~$\rho\in(0,r)$. Using the Poisson Kernel
in Theorems~\ref{POIBALL1} and~\ref{POIBALL}, we can construct a harmonic
function~$v$ in~$B_\rho$ such that~$v=u$ along~$\partial B_\rho$.
Let~$\e>0$. Define~$\alpha:=u-v+\e\Gamma+\e$
and~$\beta:=u-v-\e\Gamma-\e$.
By~\eqref{MDIVINZ}, we can find~$\delta\in(0,\rho)\cap(0,\e)$ such that
if~$|x|\in(0,\delta]$ then~$|u(x)|\le\frac\e2\,\Gamma(x)$.
Since~$\Gamma$ diverges at the origin,
by possibly taken~$\delta$ smaller, we can also suppose that
if~$|x|\in(0,\delta]$ then~$\|v\|_{L^\infty(B_\rho)}\le\frac\e2\,\Gamma(x)$.

As a result, if~$x\in\partial B_\delta$,
$$ \alpha(x)\ge\e\qquad{\mbox{and}}\qquad\beta(x)\le-\e.$$
Moreover, if~$x\in\partial B_\rho$,
$$ \alpha(x)=\e\Gamma(x)+\e\ge\e\qquad{\mbox{and}}\qquad\beta(x)
=-\e\Gamma(x)-\e
\le-\e.$$
Since~$\alpha$ and~$\beta$ are harmonic in~$B_\rho\setminus B_\delta$
we thereby deduce from the
Maximum Principle in Corollary~\ref{WEAKMAXPLE} that~$\alpha\ge\e$
and~$\beta\le-\e$ in~$B_\rho\setminus B_\delta$.
Accordingly, in~$B_\rho\setminus B_\delta$,
$$ u-v+\e\Gamma+\e\ge\e \qquad{\mbox{and}}\qquad
u-v-\e\Gamma-\e\le-\e.$$
By taking~$\e$ as small as we wish, we thus conclude that
in~$B_\rho$
$$ u-v\ge0 \qquad{\mbox{and}}\qquad
u-v\le0.$$
That is, $u$ coincides with the harmonic function~$v$ in~$B_\rho$ and the desired
result is established. \end{proof}

\begin{figure}
  \centering
  \includegraphics[width=.4\linewidth]{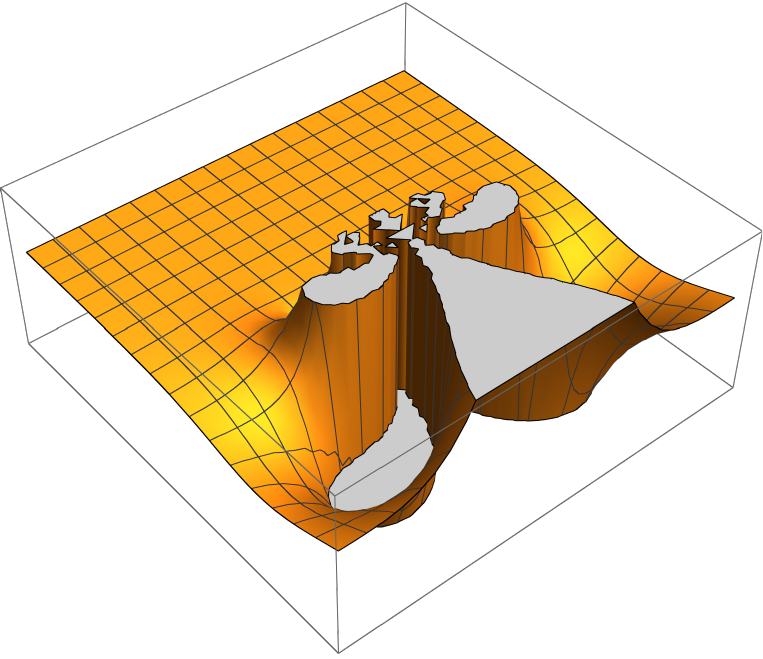}$\qquad$
  \includegraphics[width=.4\linewidth]{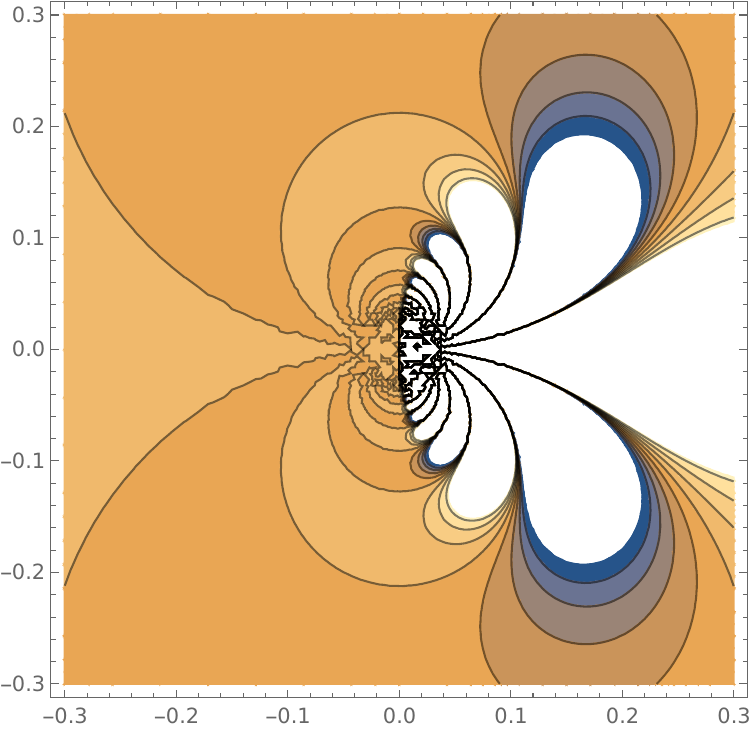}
 \caption{\sl Plot of the function in~\eqref{9203iohgr9490E542uu56Z0264446h}
and of its level sets.}\label{G02iur32-8u5utgjh93uythg9000458utjy477556n346geFI}
\end{figure}

We stress that assumption~\eqref{MDIVINZ}
cannot be removed from the statement of Theorem~\ref{0oiujhgv098j-23rtgrJSPsdfgdh56uh}: for instance, if we consider the function
\begin{equation}\label{9203iohgr9490E542uu56Z0264446h} \R^2\setminus\{0\}\ni(x,y)\longmapsto u(x,y):=\Re \left( \exp\left(\frac{1}{x+iy} \right)\right)=
e^{\frac{x}{x^2+y^2}}\cos\left({\frac{y}{x^2+y^2}}\right)
\end{equation}
we see that~$u$ is harmonic in~$\R^2\setminus\{0\}$ (being the real part of a holomorphic function outside the origin), but it presents an essential singularity at the origin (see Figure~\ref{G02iur32-8u5utgjh93uythg9000458utjy477556n346geFI} for a diagram of this function).\medskip

The study of removable singularities is a very broad and classical topic
of investigation in partial differential equations, see e.g.~\cite{MR859333}
for further information.

\section{Further geometric insights on harmonic functions in balls}

Here we present a number of fascinating results
focused on some geometric properties
of harmonic functions. The core of the discussion will
start only on page~\pageref{MALM} (and the expert reader is
welcome to jump directly there): before that, we start
by discussing a change of variable formula for a
spherical integral in~$\R^n$:

\begin{lemma}\label{COV9INSPJERE}
Let~$n\ge2$, $R>1\ge r>0$, and~$Q$ be a diffeomorphism
of~$B_R\setminus B_r$ such that~$Q(\partial B_\rho)=\partial B_\rho$
for each~$\rho\in[r,R]$. Let~$g:\R^n\to\R$ be
a continuous function.

Then,
\begin{equation}\label{DHNETERGD01324u54}
\int_{\partial B_1} g(\omega)\,d{\mathcal{H}}^{n-1}_\omega=
\int_{\partial B_1}
g\left( Q(\omega)\right)\,
|\det DQ(\omega)|\,d{\mathcal{H}}^{n-1}_\omega.\end{equation}
\end{lemma}

\begin{proof}
The change of variable~$x:=Q(y)$ leads to
\begin{equation*}
\begin{split}&
\int_{\partial B_1} g(\omega)\,d{\mathcal{H}}^{n-1}_\omega
=\frac{n}{R^n-r^n}
\int_r^R \left[ \int_{\partial B_1} \rho^{n-1} g(\omega)\,d{\mathcal{H}}^{n-1}_\omega\right]\,d\rho\\&\qquad=\frac{n}{R^n-r^n}
\int_{B_R\setminus B_r} g\left( \frac{x}{|x|}\right)\,dx=
\frac{n}{R^n-r^n}
\int_{B_R\setminus B_r} g\left( \frac{Q(y)}{|Q(y)|}\right)\,|\det DQ(y)|\,dy\\
&\qquad=
\frac{n}{R^n-r^n}\int_r^R\left[
\int_{\partial B_1}\rho^{n-1}
g\left( \frac{Q(\rho\omega)}{|Q(\rho\omega)|}\right)\,|\det DQ(\rho\omega)|\,d{\mathcal{H}}^{n-1}_\omega\right]\,d\rho
.\end{split}\end{equation*}
In particular, picking~$R:=1+\e$ and~$r:=1$,
\begin{equation*}
\begin{split}&
\int_{\partial B_1} g(\omega)\,d{\mathcal{H}}^{n-1}_\omega=
\frac{n}{(1+\e)^n-1}
\int_1^{1+\e}\left[
\int_{\partial B_1}\rho^{n-1}
g\left( \frac{Q(\rho\omega)}{|Q(\rho\omega)|}\right)\,|\det DQ(\rho\omega)|\,d{\mathcal{H}}^{n-1}_\omega\right]\,d\rho.
\end{split}\end{equation*}
Taking the limit as~$\e\searrow0$, this equation and
L'H\^opital's Rule lead to
\begin{equation*}
\begin{split}&
\int_{\partial B_1} g(\omega)\,d{\mathcal{H}}^{n-1}_\omega\\&\qquad=\lim_{\e\searrow0}
\frac{1}{(1+\e)^{n-1}}
\int_{\partial B_1}({1+\e})^{n-1}
g\left( \frac{Q(({1+\e})\omega)}{|Q(({1+\e})\omega)|}\right)\,
|\det DQ(({1+\e})\omega)|\,d{\mathcal{H}}^{n-1}_\omega
\\&\qquad=
\int_{\partial B_1}
g\left( \frac{Q(\omega)}{|Q(\omega)|}\right)\,
|\det DQ(\omega)|\,d{\mathcal{H}}^{n-1}_\omega\\&\qquad=
\int_{\partial B_1}
g\left( Q(\omega)\right)\,
|\det DQ(\omega)|\,d{\mathcal{H}}^{n-1}_\omega.\qedhere
\end{split}\end{equation*}
\end{proof}

\begin{figure}
  \centering
  \includegraphics[width=.55\linewidth]{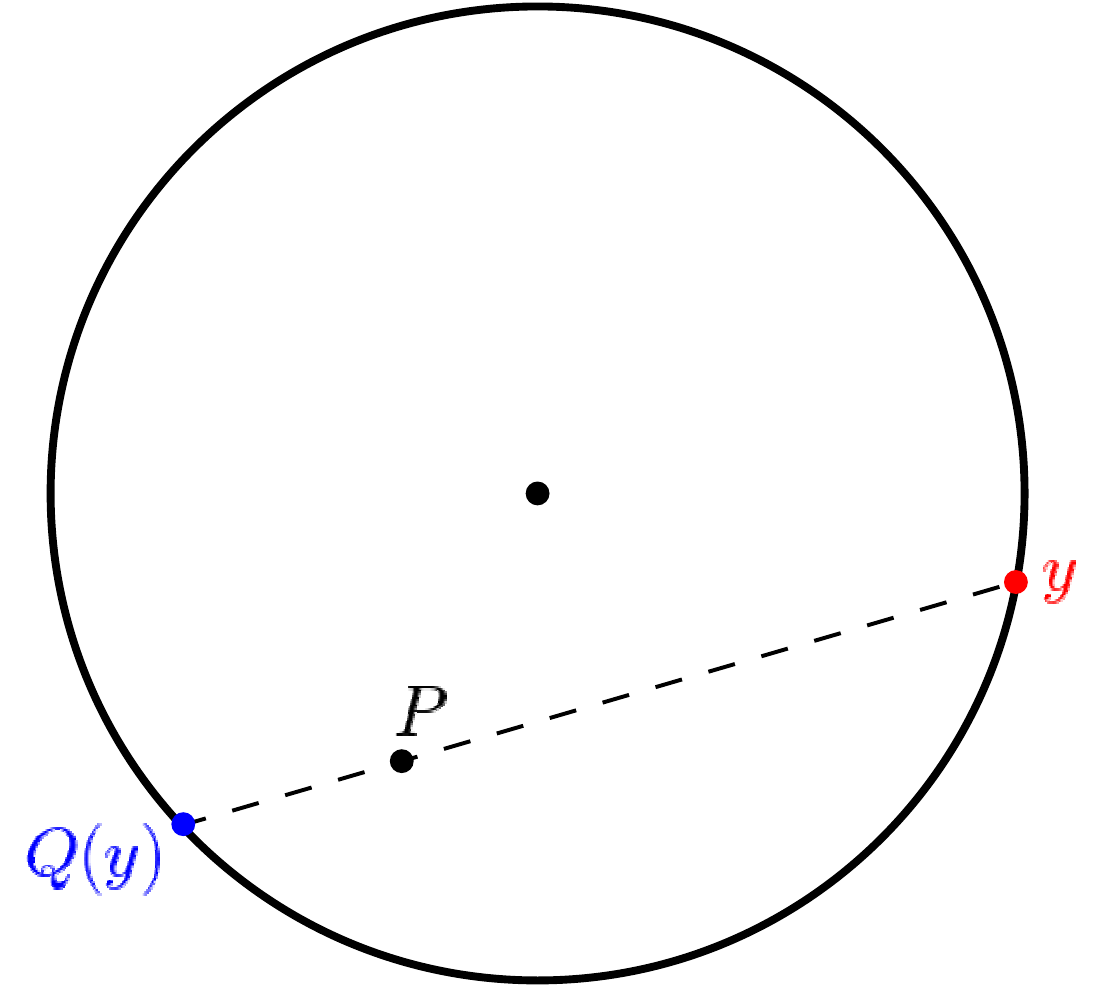}
 \caption{\sl The map~$Q:\partial B_1\to\partial B_1$ in~\eqref{QDEFI2}.}\label{2354HA466RH788AR9C0O0U1N2TERELAGIJ7soloDItangeFI}
\end{figure}

We now specify Lemma~\ref{COV9INSPJERE} for a particular choice of
the map~$Q$ (for neat geometric interpretations in the plane,
see~\cite{MR86885, MR1272823}). Namely,
given~$P\in B_1$, for all~$y\in B_2\setminus B_{(1+|P|)/2}$ we define
\begin{equation}\label{QDEFI2} Q(y):=y-\frac{2(P-y)\cdot y}{|P-y|^2}\,(P-y).
\end{equation}
We observe that~$Q(y)$ belongs to the straight line~$\left\{ 
y+t(P-y),\;t\in\R\right\}$ passing through~$y$ and~$P$.
Furthermore, 
\begin{eqnarray*} |Q(y)|^2=
|y|^2+
\frac{4((P-y)\cdot y)^2}{|P-y|^2}-
\frac{4(P-y)\cdot y}{|P-y|^2}\,(P-y)\cdot y
=|y|^2
\end{eqnarray*}
and therefore~$Q$ maps~$\partial B_\rho$ into itself, 
for each~$\rho\in\left[\frac{1+|P|}{2},2\right]$,
see Figure~\ref{2354HA466RH788AR9C0O0U1N2TERELAGIJ7soloDItangeFI} for a representation of the action of~$Q$
on the unit circumference in the plane.
In this setting, we have:

\begin{lemma}\label{DSFHUDSIPANFADIKDMMAEDA098765UJ}
Let~$n\ge2$, $P\in B_1$ and~$Q$ be as in~\eqref{QDEFI2}.

Then, for every~$\omega\in\partial B_1$ we have that
\begin{equation}\label{FGHJdsfgvbh0987654ikTYHSJ}
|\det DQ(\omega)|=
\left(\frac{ 1-|P|^2 }{|P-\omega|^2}\right)^{n-1}.
\end{equation}
Moreover, for every continuous function~$g:\R^n\to\R$,
\begin{equation}\label{DHNETERGD01324u54-2}
\int_{\partial B_1} g(\omega)\,d{\mathcal{H}}^{n-1}_\omega=
\int_{\partial B_1}
g\left( Q(\omega)\right)\,\left(\frac{ 1-|P|^2 }{|P-\omega|^2}\right)^{n-1}\,d{\mathcal{H}}^{n-1}_\omega.\end{equation}
\end{lemma}

\begin{proof}
In the light of~\eqref{QDEFI2}, we point out that
\begin{equation}\label{TOR-POISS-PRE}
Q(y)-P=-\left(1+\frac{2(P-y)\cdot y}{|P-y|^2}\right)\,(P-y)=
-\frac{|P|^2-|y|^2}{|P-y|^2}\,(P-y)
\end{equation}
and therefore, for every~$\omega\in\partial B_1$,
\begin{equation*}
\frac{|P-Q(\omega)|}{|P-\omega|}=
\frac{1-|P|^2}{|P-\omega|^2}.\end{equation*}
We now pick a point~$\omega\in\partial B_1$
and we consider the plane containing
the origin, $P$ and~$\omega$.
Up to a rotation, we suppose that
this plane is the one given by the first coordinates,
namely~$\Pi:=\{x_3=\dots=x_n=0\}$ (of course, if~$n=2$, this part simplifies
since~$\Pi$ would coincide with the whole plane of interest~$\R^2$).
Consequently, by~\eqref{QDEFI2},
we know that
if~$y\in\Pi$ then~$Q(y)\in\Pi$. Since~$\omega$, $e_1$, $e_2\in\Pi$,
this gives that~$Q(\omega+\e e_1)$
and~$Q(\omega+\e e_2)$ belong to~$\Pi$ for all~$\e\in\R$.
Accordingly, $\partial_1 Q(\omega)$ and~$\partial_2 Q(\omega)$
belong to~$\Pi$, that is
$$ \partial_1 Q(\omega)=\big(\partial_1 Q_1(\omega),
\partial_1 Q_2(\omega),0,\dots,0\big)\qquad{\mbox{and}}\qquad
\partial_2 Q(\omega)=\big(\partial_2 Q_1(\omega),
\partial_2 Q_2(\omega),0,\dots,0\big).$$
For this reason,
$$ DQ(\omega)=
\left(\begin{matrix}
\partial_1 Q_1(\omega)&\partial_1 Q_2(\omega)&0&\dots&0\\
\partial_2 Q_1(\omega)&\partial_2 Q_2(\omega)&0&\dots&0\\
\partial_3 Q_1(\omega)&\partial_3 Q_2(\omega)&\partial_3 Q_3(\omega)&\dots&\partial_3 Q_n(\omega)\\
& &\ddots&\\
\partial_n Q_1(\omega)&\partial_n Q_2(\omega)&
\partial_n Q_3(\omega)&\dots&\partial_n Q_n(\omega)
\end{matrix}
\right)$$
which entails that
\begin{equation}\label{LEO0}
|\det DQ(\omega)|=\left|
\det\left(
\begin{matrix}
\partial_1 Q_1(\omega)&\partial_1 Q_2(\omega)\\
\partial_2 Q_1(\omega)&\partial_2 Q_2(\omega)
\end{matrix}
\right)
\right|\,\left|
\det\left(\begin{matrix}
\partial_3 Q_3(\omega)&\dots&\partial_3 Q_n(\omega)\\
&\ddots&\\
\partial_n Q_3(\omega)&\dots&\partial_n Q_n(\omega)
\end{matrix}
\right)
\right|.
\end{equation}
Also, by~\eqref{TOR-POISS-PRE},
for every~$j\in\{1,\dots,n\}$,
\begin{equation}\label{09-008-POS-OPLSMSKSM-023943586455}
\partial_j Q(\omega)=
\frac{2\omega_j}{|P-\omega|^2}\,(P-\omega)
+\frac{2(1-|P|^2)(P_j-\omega_j)}{|P-\omega|^4}\,(P-\omega)
-\frac{1-|P|^2}{|P-\omega|^2}\,e_j
\end{equation}
Therefore, using again that~$
\omega_3=\dots=\omega_n=P_3=\dots=P_n=0$,
we deduce that, for every~$j\in\{3,\dots,n\}$,
\begin{eqnarray*}
\partial_j Q(\omega)=-
\frac{1-|P|^2}{|P-\omega|^2}\,e_j.
\end{eqnarray*}
As a result, the matrix~$\left(\begin{matrix}
\partial_3 Q_3(\omega)&\dots&\partial_3 Q_n(\omega)\\
&\ddots&\\
\partial_n Q_3(\omega)&\dots&\partial_n Q_n(\omega)
\end{matrix}
\right)$ is~$ -\frac{1-|P|^2}{|P-\omega|^2} $ times the identity
matrix (intended here as a square $(n-2)\times(n-2)$ matrix).
This observation
and~\eqref{LEO0} lead to
\begin{equation}\label{KMDCNFMONDUYJNDFFJNVEHDNIDI}
|\det DQ(\omega)|=
\left(
\frac{1-|P|^2}{|P-\omega|^2}
\right)^{n-2}\;
\left|
\det\left(
\begin{matrix}
\partial_1 Q_1(\omega)&\partial_1 Q_2(\omega)\\
\partial_2 Q_1(\omega)&\partial_2 Q_2(\omega)
\end{matrix}
\right)
\right|.
\end{equation}
Now, we rewrite~\eqref{09-008-POS-OPLSMSKSM-023943586455}
in the form
\begin{equation*}
\partial_j Q(\omega)=-
\frac{1-|P|^2}{|P-\omega|^2}\,e_j
+\alpha_j(\omega)\,(P-\omega),\qquad{\mbox{where}}\qquad
\alpha_j(\omega):=
\frac{2\omega_j}{|P-\omega|^2}
+\frac{2(1-|P|^2)(P_j-\omega_j)}{|P-\omega|^4}.
\end{equation*}
Consequently, letting~$v(\omega):=(P_1-\omega_1,P_2-\omega_2)$
\begin{equation}\label{KHNS:DEBNDNALINTNSATA}
\begin{split}
&\det\left(
\begin{matrix}
\partial_1 Q_1(\omega)&\partial_1 Q_2(\omega)\\
\partial_2 Q_1(\omega)&\partial_2 Q_2(\omega)
\end{matrix}
\right)
=\det\left(
\begin{matrix}-
\frac{1-|P|^2}{|P-\omega|^2}\,(1,0)
+\alpha_1(\omega)\,v(\omega)
\\-
\frac{1-|P|^2}{|P-\omega|^2}\,(0,1)
+\alpha_2(\omega)\,v(\omega)
\end{matrix}
\right)\\&\qquad=\left(
\frac{1-|P|^2}{|P-\omega|^2}\right)^2\det\left(
\begin{matrix}
1&0\\0&1\end{matrix}
\right)-
\alpha_2(\omega)\frac{1-|P|^2}{|P-\omega|^2}\det\left(
\begin{matrix}
(1,0)
\\
v(\omega)
\end{matrix}
\right)
-\alpha_1(\omega)\frac{1-|P|^2}{|P-\omega|^2}
\det\left(
\begin{matrix}
v(\omega)
\\
(0,1)\end{matrix}
\right)\\&\qquad=\left(
\frac{1-|P|^2}{|P-\omega|^2}\right)^2-
\alpha_2(\omega)\,(P_2-\omega_2)\frac{1-|P|^2}{|P-\omega|^2}
-\alpha_1(\omega)\,(P_1-\omega_1)\frac{1-|P|^2}{|P-\omega|^2}
\\&\qquad=\left(
\frac{1-|P|^2}{|P-\omega|^2}\right)\left[\left(
\frac{1-|P|^2}{|P-\omega|^2}\right)-
\alpha_2(\omega)\,(P_2-\omega_2)
-\alpha_1(\omega)\,(P_1-\omega_1)\right].
\end{split}\end{equation}
We also point out that
\begin{eqnarray*}&&
\alpha_1(\omega)\,(P_1-\omega_1)+\alpha_2(\omega)\,(P_2-\omega_2)\\
&=&\sum_{j=1}^2\left(
\frac{2\omega_j(P_j-\omega_j)}{|P-\omega|^2}
+\frac{2(1-|P|^2)(P_j-\omega_j)^2}{|P-\omega|^4}\right)\\&=&
\frac{2(P\cdot\omega-1)}{|P-\omega|^2}
+\frac{2(1-|P|^2)}{|P-\omega|^2}\\
&=&\frac{2(P\cdot\omega-|P|^2)}{|P-\omega|^2}
\end{eqnarray*}
and therefore
\begin{eqnarray*}&&
\frac{1-|P|^2}{|P-\omega|^2}-
\alpha_1(\omega)\,(P_1-\omega_1)-\alpha_2(\omega)\,(P_2-\omega_2)
=
\frac{1-2P\cdot\omega+|P|^2}{|P-\omega|^2}=1.\end{eqnarray*}
Plugging this information into~\eqref{KHNS:DEBNDNALINTNSATA},
we find that
$$
\det\left(
\begin{matrix}
\partial_1 Q_1(\omega)&\partial_1 Q_2(\omega)\\
\partial_2 Q_1(\omega)&\partial_2 Q_2(\omega)
\end{matrix}
\right)=
\frac{ 1-|P|^2 }{|P-\omega|^2}.$$
Hence, recalling~\eqref{KMDCNFMONDUYJNDFFJNVEHDNIDI}
we obtain~\eqref{FGHJdsfgvbh0987654ikTYHSJ}.
This and~\eqref{DHNETERGD01324u54} give~\eqref{DHNETERGD01324u54-2},
as desired.
\end{proof}

\begin{figure}
  \centering
  \includegraphics[width=.55\linewidth]{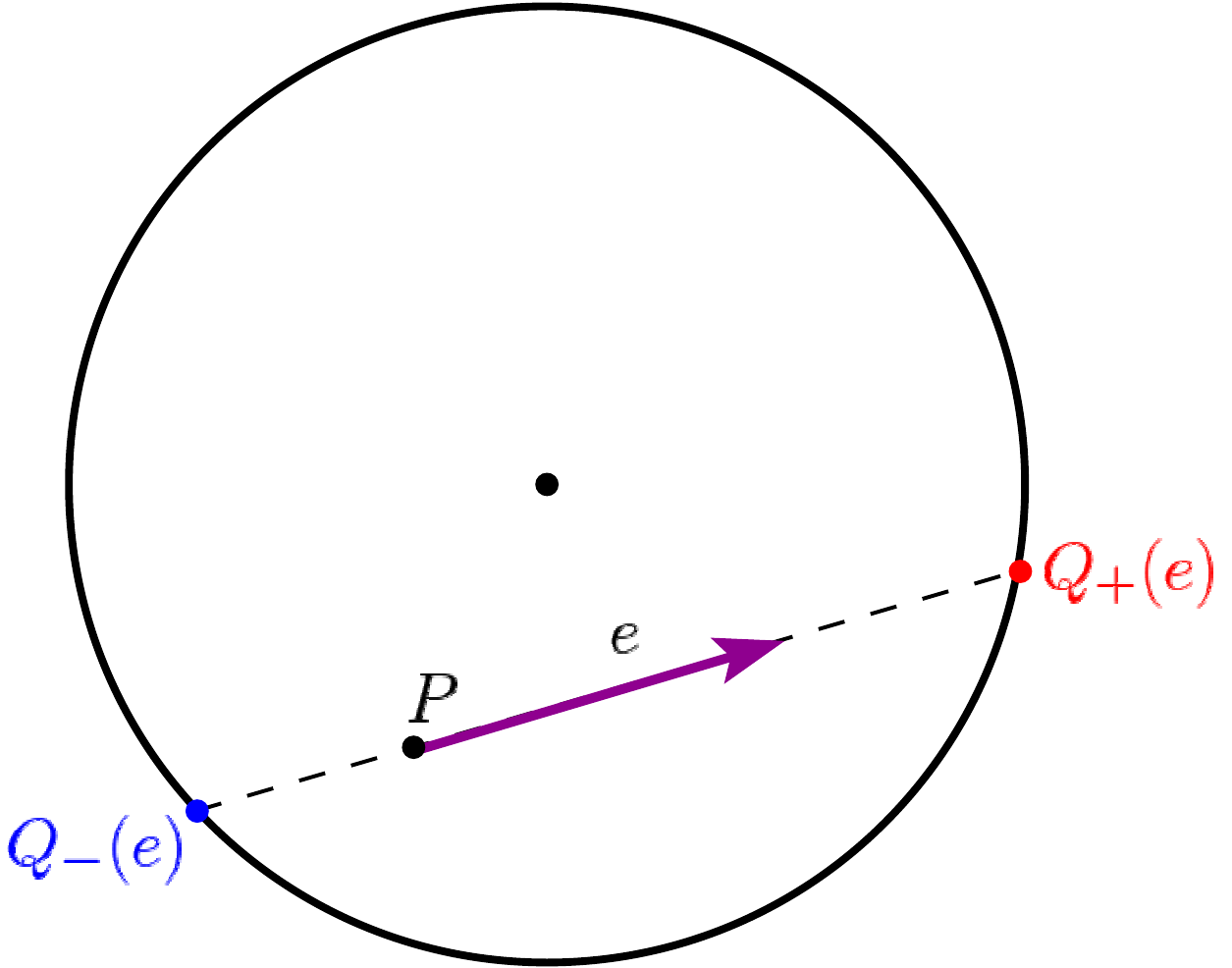}
 \caption{\sl The maps~$Q_\pm:\partial B_1\to\partial B_1$ in~\eqref{PIUJSMENP}.}\label{2354HA466RH7dss8AR9C2O0U1N3TERELSGIJ7so8oDI3angeFI}
\end{figure}

Another interesting change of variable, which can be seen
as a slight modification of the one discussed in Lemma~\ref{DSFHUDSIPANFADIKDMMAEDA098765UJ}
is the following.
For every~$P\in B_1$
and every~$e\in\partial B_1$, we consider the points~$Q_+(e)$
and~$Q_-(e)$ lying in the intersection between the straight line
of direction~$e$ passing through~$P$ and~$\partial B_1$, namely
let~$r_+(e)>0>r_-(e)$ be such that
\begin{equation}\label{PIUJSMENP}
Q_+(e)=P+r_+(e) \,e\in\partial B_1\qquad{\mbox{and}}
\qquad Q_-(e)=P+r_-(e) \,e\in\partial B_1.
\end{equation} 
See Figure~\ref{2354HA466RH7dss8AR9C2O0U1N3TERELSGIJ7so8oDI3angeFI} for a sketch of this configuration.
Also, by~\eqref{PIUJSMENP},
\begin{equation}\label{EUIND0OSP92345yt} 1=|Q_\pm(e)|^2=|P+r_\pm(e) \,e|^2=
|P|^2+(r_\pm(e))^2+ 2 P\cdot e\,r_\pm(e).\end{equation}
Therefore, 
\begin{equation}\label{CHISMD-0193248OKE}
r_\pm(e)=-P\cdot e\pm\sqrt{D(e)},\qquad{\mbox{where}}\qquad
D(e):=(P\cdot e)^2-|P|^2+1.
\end{equation}
The maps~$Q_\pm$ can also be extended in~$\R^n\setminus\{0\}$
by considering the homogeneous extension of degree~$1$:
in this way, for every~$x\in\R^n\setminus\{0\}$ we have that
\begin{equation}\label{OMOGENEIT7}
Q_\pm(x)=|x|\,Q_\pm\left(\frac{x}{|x|}\right).\end{equation}
This and~\eqref{EUIND0OSP92345yt} give that, for each~$\rho>0$,
if~$x\in \partial B_\rho$ then~$|Q_\pm(x)|=\rho\,\left|
Q_\pm\left(\frac{x}{|x|}\right)\right|=\rho$, hence~$Q_\pm$ maps~$\partial B_\rho$
into itself. Thus, in analogy
with Lemma~\ref{DSFHUDSIPANFADIKDMMAEDA098765UJ}, we have the following result:

\begin{lemma}\label{7s9ik9z-lAmco-ASpOAKS}
Let~$n\ge2$, $P\in B_1$ and~$Q_\pm$ be as in~\eqref{PIUJSMENP}.

Then, for every~$\omega\in\partial B_1$ we have that
\begin{equation}\label{VEBDNFKRGZINDMKBVINSKDSA-21}
|\det DQ_\pm(\omega)|= \frac{(\pm r_\pm(\omega))^n}{1-|P|^2-r_\pm(\omega)P\cdot \omega}.
\end{equation}
Moreover, for every continuous function~$g:\R^n\to\R$,
\begin{equation}\label{VEBDNFKRGZINDMKBVINSKDSA-22}
\int_{\partial B_1} g(\omega)\,d{\mathcal{H}}^{n-1}_\omega=
\int_{\partial B_1}
g\left( Q_\pm(\omega)\right)\, \frac{(\pm r_\pm(\omega))^n}{1-|P|^2-r_\pm(\omega)P\cdot \omega}
\,d{\mathcal{H}}^{n-1}_\omega.\end{equation}
\end{lemma}

\begin{proof} Let~$\omega\in\partial B_1$ and let~$O$ denote the origin.
Up to a rotation, we can suppose that the points~$O$, $P$ and~$P+\omega$ lie
in the plane~$\{x_3=\dots=x_n=0\}$.
Also, up to a further rotation in this plane, we can suppose that~$\omega=e_1$. Thus, by~\eqref{OMOGENEIT7},
$$ Q_\pm(\omega+\e e_1)=
Q_\pm((1+\e) e_1)=
(1+\e)\,Q_\pm(e_1)=(1+\e)\big(Q_{\pm,1}(e_1),Q_{\pm,2}(e_1),0,\dots,0\big).$$
Consequently,
$$ \partial_1Q_\pm(\omega )=\big(Q_{\pm,1}(e_1),Q_{\pm,2}(e_1),0,\dots,0\big).$$
Let now~$j\in\{2,\dots,n\}$. We have that~$|\omega+\e e_j|=|e_1+\e e_j|=1+o(\e)$, therefore,
recalling~\eqref{PIUJSMENP} and~\eqref{CHISMD-0193248OKE},
\begin{eqnarray*}&&
Q_\pm(\omega+\e e_j)\\&=&(1+o(\e))
\,Q_\pm\left(\frac{\omega+\e e_j}{|\omega+\e e_j|}\right)\\
&=& P+r_\pm\left(\frac{\omega+\e e_j}{|\omega+\e e_j|}\right)\frac{\omega+\e e_j}{|\omega+\e e_j|}+o(\e)
\\&=&
P+\left(-
P\cdot \frac{\omega+\e e_j}{|\omega+\e e_j|}\pm\sqrt{
\left(P\cdot \frac{\omega+\e e_j}{|\omega+\e e_j|}\right)^2-|P|^2+1}\right) (\omega+\e e_j)+o(\e)
\\&=&
P+\left(-
P\cdot (\omega+\e e_j)\pm\sqrt{
\left(P\cdot (\omega+\e e_j)\right)^2-|P|^2+1}\right) (\omega+\e e_j)+o(\e)
\\&=&
P+\left(-
P\cdot \omega-\e P\cdot e_j\pm\sqrt{(P\cdot\omega)^2+2\e(P\cdot\omega)(P\cdot e_j)-|P|^2+1}\right) (\omega+\e e_j)+o(\e)
\\&=&
P+\left(-
P\cdot e_1-\e P\cdot e_j\pm\sqrt{(P\cdot e_1)^2-|P|^2+1}
\pm\frac{\e(P\cdot e_1)(P\cdot e_j)}{\sqrt{(P\cdot e_1)^2-|P|^2+1}}
\right) (e_1+\e e_j)+o(\e).
\end{eqnarray*}
Taking the first order in~$\e$, we find that
\begin{eqnarray*}
\partial_jQ_\pm(\omega )&=&
\left( -P\cdot e_j
\pm\frac{(P\cdot e_1)(P\cdot e_j)}{\sqrt{(P\cdot e_1)^2-|P|^2+1}}
\right) e_1
+
\left(-P\cdot e_1\pm\sqrt{(P\cdot e_1)^2-|P|^2+1}\right)  e_j\\
&=& \alpha_j e_1+r_\pm(e_1) e_j,
\end{eqnarray*}
where
$$ \alpha_j:=-P\cdot e_j
\pm\frac{(P\cdot e_1)(P\cdot e_j)}{\sqrt{(P\cdot e_1)^2-|P|^2+1}}=-\frac{P\cdot e_j\big(\sqrt{D(e_1)}\mp P\cdot e_1\big)}{\sqrt{D(e_1)}}=\mp\frac{P\cdot e_j \,r_\pm(e_1)}{\sqrt{D(e_1)}}.$$
These observations lead to
$$ DQ_\pm(\omega)=\left(\begin{matrix}
Q_{\pm,1}(e_1)&Q_{\pm,2}(e_1)&0&0&0&\dots&0\\
\alpha_2&r_\pm(e_1)&0&0&0&\dots&0\\
\alpha_3&0&r_\pm(e_1)&0&0&\dots&0\\
\alpha_4&0&0&r_\pm(e_1)&0&\dots&0\\
\alpha_5&0&0&0&r_\pm(e_1)&\dots&0\\
& & &\ddots& & &\\
\alpha_n&0&0&0&0&\dots&r_\pm(e_1)
\end{matrix}\right)$$
and therefore
\begin{equation}\label{TGBSIKMwqerrSDFGHJKTYUJMSUGJSBYGHBuasjbd}
|\det DQ_\pm(\omega)|=
\left|r_\pm(e_1)^{n-2}\det
\left(\begin{matrix}
Q_{\pm,1}(e_1)&Q_{\pm,2}(e_1)\\
\alpha_2&r_\pm(e_1)
\end{matrix}\right)
\right|.
\end{equation}
We also note that
\begin{eqnarray*}
Q_\pm(e_1)&=&(P\cdot e_1,P\cdot e_2,0,\dots,0)+
\left(-P\cdot e_1\pm\sqrt{(P\cdot e_1)^2-|P|^2+1}
\right)e_1\\&=&\left(\pm\sqrt{(P\cdot e_1)^2-|P|^2+1},
P\cdot e_2,0,\dots,0
\right)\\&=&\left(\pm\sqrt{D(e_1)},
P\cdot e_2,0,\dots,0
\right)
\end{eqnarray*}
and consequently
\begin{eqnarray*}&&
\det
\left(\begin{matrix}
Q_{\pm,1}(e_1)&Q_{\pm,2}(e_1)\\
\alpha_2&r_\pm(e_1)
\end{matrix}\right)=Q_{\pm,1}(e_1)r_\pm(e_1)-
Q_{\pm,2}(e_1)\alpha_2=
\pm\sqrt{D(e_1)}r_\pm(e_1)\pm 
\frac{(P\cdot e_2)^2 \,r_\pm(e_1)}{\sqrt{D(e_1)}}\\
&&\qquad=\pm\frac{r_\pm(e_1)}{\sqrt{D(e_1)}}\Big(D(e_1)
+ (P\cdot e_2)^2
\Big)=\pm\frac{r_\pm(e_1)}{\sqrt{D(e_1)}}\Big((P\cdot e_1)^2-|P|^2+1
+ (P\cdot e_2)^2
\Big)\\&&\qquad=\pm\frac{r_\pm(e_1)}{\sqrt{D(e_1)}}=\frac{r_\pm(e_1)}{r_\pm(e_1)+P\cdot e_1}
=\frac{(r_\pm(e_1))^2}{(r_\pm(e_1))^2+P\cdot e_1\,r_\pm(e_1)}.
\end{eqnarray*}
Since
\begin{equation}\begin{split} \label{2134.1317}
&(r_\pm(e_1))^2+P\cdot e_1\,r_\pm(e_1)\\ =\,&
(P\cdot e_1)^2+D(e_1)\mp2P\cdot e_1\sqrt{D(e_1)}+P\cdot e_1\,r_\pm(e_1)\\
=\,&(P\cdot e_1)^2+D(e_1)\mp2P\cdot e_1\sqrt{D(e_1)}-(P\cdot e_1)^2
\pm(P\cdot e_1)\sqrt{D(e_1)}\\
=\,&D(e_1)\mp P\cdot e_1\sqrt{D(e_1)}\\=\,&
1-|P|^2+(P\cdot e_1)^2\mp P\cdot e_1\sqrt{D(e_1)}\\=\,&
1-|P|^2-(P\cdot e_1)\Big(-P\cdot e_1\pm \sqrt{D(e_1)}\Big)\\=\,&
1-|P|^2-(P\cdot e_1)r_\pm(e_1),
\end{split}\end{equation}
we arrive at
$$ \det
\left(\begin{matrix}
Q_{\pm,1}(e_1)&Q_{\pm,2}(e_1)\\
\alpha_2&r_\pm(e_1)
\end{matrix}\right)=\frac{(r_\pm(e_1))^2}{1-|P|^2-(P\cdot e_1)r_\pm(e_1)}.$$
Thus, retaking~\eqref{TGBSIKMwqerrSDFGHJKTYUJMSUGJSBYGHBuasjbd}, \begin{equation}\label{JS-0OKS56yhn7tgufwe9ufger}
|\det DQ_\pm(\omega)|=\left|
\frac{(r_\pm(e_1))^n}{1-|P|^2-(P\cdot e_1)r_\pm(e_1)}\right|.\end{equation}
It is also useful to notice that
\begin{eqnarray*}&&
1-|P|^2-(P\cdot e_1)r_\pm(e_1)=
1-(P\cdot e_1)^2-(P\cdot e_2)^2-(P\cdot e_1)\Big(-(P\cdot e_1)\pm\sqrt{D(e_1)}\Big)\\
&&\qquad=1-(P\cdot e_1)^2-(P\cdot e_2)^2-(P\cdot e_1)\big(-(P\cdot e_1)\pm\sqrt{1-(P\cdot e_2)^2}\big)\\&&\qquad=1-(P\cdot e_2)^2\mp(P\cdot e_1)\sqrt{1-(P\cdot e_2)^2}\\&&\qquad=
\sqrt{1-(P\cdot e_2)^2}\,\Big(\sqrt{1-(P\cdot e_2)^2}\mp(P\cdot e_1)\Big)\\&&\qquad\ge
\sqrt{1-(P\cdot e_2)^2}\,\Big(\sqrt{1-(P\cdot e_2)^2}-|P\cdot e_1|\Big)\\&&\qquad=
\frac{\sqrt{1-(P\cdot e_2)^2}}{\sqrt{1-(P\cdot e_2)^2}+|P\cdot e_1|}
\,\Big( 1-(P\cdot e_2)^2-(P\cdot e_1)^2\Big)
\\&&\qquad=
\frac{\sqrt{1-(P\cdot e_2)^2}}{\sqrt{1-(P\cdot e_2)^2}+|P\cdot e_1|}
\,\Big( 1-|P|^2\Big)\ge0,
\end{eqnarray*}
which, together with~\eqref{JS-0OKS56yhn7tgufwe9ufger}, establishes~\eqref{VEBDNFKRGZINDMKBVINSKDSA-21}.
This and~\eqref{DHNETERGD01324u54} lead to~\eqref{VEBDNFKRGZINDMKBVINSKDSA-22}.
\end{proof}

Now we are ready to dive into the analysis of some geometric features of harmonic
functions. In this framework,
a classical result, \label{MALM}
sometimes referred to with the name of Malmheden's Theorem\index{Malmheden's Theorem},
provides an elegant geometric algorithm for solving the Dirichlet problem in a ball,
see~\cite{NEUMANN, zbMATH03014632, MR86885, MR1272823, MR2734448}.
In a nutshell, the main idea for this algorithm is:
\begin{itemize}
\item to
consider a point~$P$ in a given ball, whose boundary is endowed by some
continuous datum~$f$,
\item then,
take an arbitrary chord passing through the point~$P$
and calculate the value at~$P$ of the linear function
that interpolates the values of~$f$ at the endpoints of the chord~$L$,
\item finally, compute the average of these values over all possible chords through~$P$.\end{itemize}
Quite remarkably, this procedure produces the harmonic function in the ball\footnote{We stress that
this procedure essentially works only with balls and cannot be generalized
to other domains, see~\cite{MR162953}
and Section~4 in~\cite{MR2734448}.
Instead, a generalization of this procedure to
cross sections of higher dimension has been provided in
Section~3 of~\cite{MR2734448}.}
with datum~$f$ on the boundary.

The details of this construction go as follows.
Let~$f:\partial B_1\to\R$ be a continuous function.
Let~$L_e$ be the segment joining~$Q_-(e)$ to~$Q_+(e)$, as defined in~\eqref{PIUJSMENP},
and 
\begin{equation}\label{ELLEeMA}\begin{split}&
{\mbox{let~$\ell_e$ be the linear (or, better to say, affine) function on~$L_e$}}\\&{\mbox{such that~$\ell_e(Q_-(e))=f(Q_-(e))$ and~$\ell_e(Q_+(e))=f(Q_+(e))$.}}\end{split}\end{equation}
Let also
\begin{equation}\label{DEUGRGBD} u(P):=\fint_{\partial B_1} \ell_e(P)\,d{\mathcal{H}}^{n-1}_e.\end{equation}
Notice that the function~$u$ is constructed by averaging linear interpolations of the boundary datum
in each direction.

\begin{theorem}\label{HARMOu3}
The function~$u$ in~\eqref{DEUGRGBD}
is harmonic in~$B_1$, continuous in~$\overline{B_1}$,
and~$u=f$ on~$\partial B_1$.
\end{theorem}

We observe that when~$P$ is the origin Theorem~\ref{HARMOu3} reduces to the Mean Value Formula in Theorem~\ref{KAHAR}(ii).
Also, when~$n=1$ Theorem~\ref{HARMOu3} boils down to the fact that
harmonic functions are linear (hence, Theorem~\ref{HARMOu3} is interesting
only when~$n\ge2$).

Following~\cite{MR2734448}, we will give two proofs of Theorem~\ref{HARMOu3}, 
one now,
exploiting Lemma~\ref{7s9ik9z-lAmco-ASpOAKS}, and one 
starting on page~\pageref{5ddg-5115566},
relying on harmonic polynomials.

\begin{proof}[Proof of Theorem~\ref{HARMOu3}] 
First, we show that, for every~$\zeta\in\partial B_1$,
\begin{equation}\label{HARMOu2}
\lim_{P\to\zeta} u(P)=f(\zeta).
\end{equation}
To this end, we 
point out that the linear function~$\ell_e$ on~$L_e:=\big\{
P+se,\;s\in[r_-(e),r_+(e)]\big\}$
such that~$\ell_e(Q_-(e))=f(Q_-(e))$ has the form
\begin{equation}\label{EELS}
\ell_e(P+se)= \frac{\big(f(Q_+(e))-f(Q_-(e))\big)s +
r_+(e)f(Q_-(e))-r_-(e)f(Q_+(e))
}{r_+(e)-r_-(e)}.\end{equation}
Furthermore, 
in view of~\eqref{CHISMD-0193248OKE}, the quantities~$Q_\pm(e)$ and~$r_\pm(e)$
depend continuously on~$P$,
with either~$Q_-(e)$ or~$Q_+(e)$ approaching~$\zeta$
and correspondingly either to~$r_-(e)$ or~$r_+(e)$ approaching~$0$
as~$P\to\zeta$.
As a result, by~\eqref{EELS},
\begin{eqnarray*}
\lim_{P\to\zeta}\ell_e(P)= \lim_{P\to\zeta}
\frac{
r_+(e)f(Q_-(e))-r_-(e)f(Q_+(e))
}{r_+(e)-r_-(e)}
=f(\zeta).
\end{eqnarray*}
{F}rom this
and~\eqref{DEUGRGBD} we obtain~\eqref{HARMOu2}, as desired.

Furthermore, in light of~\eqref{DEUGRGBD} and~\eqref{EELS},
\begin{equation}\label{AM-0okr8t-12ersfZZg}\begin{split} u(P)\,&=\,\fint_{\partial B_1} 
\frac{
r_+(e)f(Q_-(e))-r_-(e)f(Q_+(e))
}{r_+(e)-r_-(e)}\,d{\mathcal{H}}^{n-1}_e
\\&=\,\fint_{\partial B_1} 
\frac{
r_+(e)f(Q_-(e))
}{r_+(e)-r_-(e)}\,d{\mathcal{H}}^{n-1}_e
-\fint_{\partial B_1} 
\frac{r_-(e)f(Q_+(e))
}{r_+(e)-r_-(e)}\,d{\mathcal{H}}^{n-1}_e
\\&=\,2\fint_{\partial B_1} 
\frac{
r_+(e)f(Q_-(e))
}{r_+(e)-r_-(e)}\,d{\mathcal{H}}^{n-1}_e.
\end{split}
\end{equation}
Moreover, recalling~\eqref{CHISMD-0193248OKE},
\begin{equation}\label{ANT6}r_+(e) r_-(e)
=(P\cdot e)^2-D(e)=
|P|^2-1.\end{equation}
Consequently, owing also to~\eqref{CHISMD-0193248OKE} and~\eqref{2134.1317},
\begin{equation}\label{90KS-098765RYUISOoscuSS0jsdmc}
\begin{split}&
\frac{2r_+(e)}{r_+(e)-r_-(e)}=\frac{r_+(e)}{\sqrt{D(e)}}=-\frac{r_+(e)}{r_-(e)+P\cdot e}
=-\frac{r_+(e)r_-(e)}{(r_-(e))^2+P\cdot e\,r_-(e)}\\
&\qquad\qquad
=-\frac{r_+(e)r_-(e)}{1-|P|^2-(P\cdot e)r_-(e)}=\frac{1-|P|^2}{1-|P|^2-(P\cdot e)r_-(e)}.
\end{split}\end{equation}
Thus, since~$|P-Q_-(e)|=|r_-(e)|=-r_-(e)$, we deduce from~\eqref{VEBDNFKRGZINDMKBVINSKDSA-22}, applied here with~$g(e):=
\frac{f(e)\,(1-|P|^2)}{|P-e|^n}$,
and~\eqref{AM-0okr8t-12ersfZZg} that
\begin{eqnarray*}
u(P)&=&\fint_{\partial B_1} \frac{f(Q_-(e))\,(1-|P|^2)}{1-|P|^2-(P\cdot e)r_-(e)}
\,d{\mathcal{H}}^{n-1}_e\\
&=&\fint_{\partial B_1}
\frac{f(Q_-(e))\,(1-|P|^2)}{|P-Q_-(e)|^n}\, \frac{(-r_-(e))^n}{{1-|P|^2-(P\cdot e)r_-(e)}}
\,d{\mathcal{H}}^{n-1}_e\\
&=&\fint_{\partial B_1}
g(Q_-(e))\, \frac{(-r_-(e))^n}{{1-|P|^2-(P\cdot e)r_-(e)}}
\,d{\mathcal{H}}^{n-1}_e\\&=&\fint_{\partial B_1}
g(e)\,d{\mathcal{H}}^{n-1}_e
\\&=&\fint_{\partial B_1}
\frac{f(e)\,(1-|P|^2)}{|P-e|^n}\,d{\mathcal{H}}^{n-1}_e.
\end{eqnarray*}
That is, recalling~\eqref{B1}
$$ u(P)=\frac1{n\,|B_1|}
\int_{\partial B_1}
\frac{f(e)\,(1-|P|^2)}{|P-e|^n}\,d{\mathcal{H}}^{n-1}_e.$$
Recognizing the Poisson Kernel of the ball, in view of Theorems~\ref{POIBALL1}
and~\ref{POIBALL}, we obtain that~$u$ is harmonic in~$B_1$, as desired.
\end{proof}

Harmonic functions in the plane
enjoy special geometric features (see e.g.~\cite{MR1272823}). In particular,
the following result holds true:

\begin{theorem}\label{SCH}
Let~$n=2$, $f:\partial B_1\to\R$ be a continuous function and~$P\in B_1$. For every~$\omega\in\partial B_1$, let
$$Q(\omega):=\omega-\frac{2(P-\omega)\cdot \omega}{|P-\omega|^2}\,(P-\omega).$$
Let
\begin{equation}\label{DEFUMCD} u(P):=\fint_{\partial B_1} f(Q(\omega))\,d{\mathcal{H}}^{1}_\omega.\end{equation}
Then, $u$
is harmonic in~$B_1$, continuous in~$\overline{B_1}$,
and~$u=f$ on~$\partial B_1$.
\end{theorem}

\begin{figure}
  \centering
  \includegraphics[width=.55\linewidth]{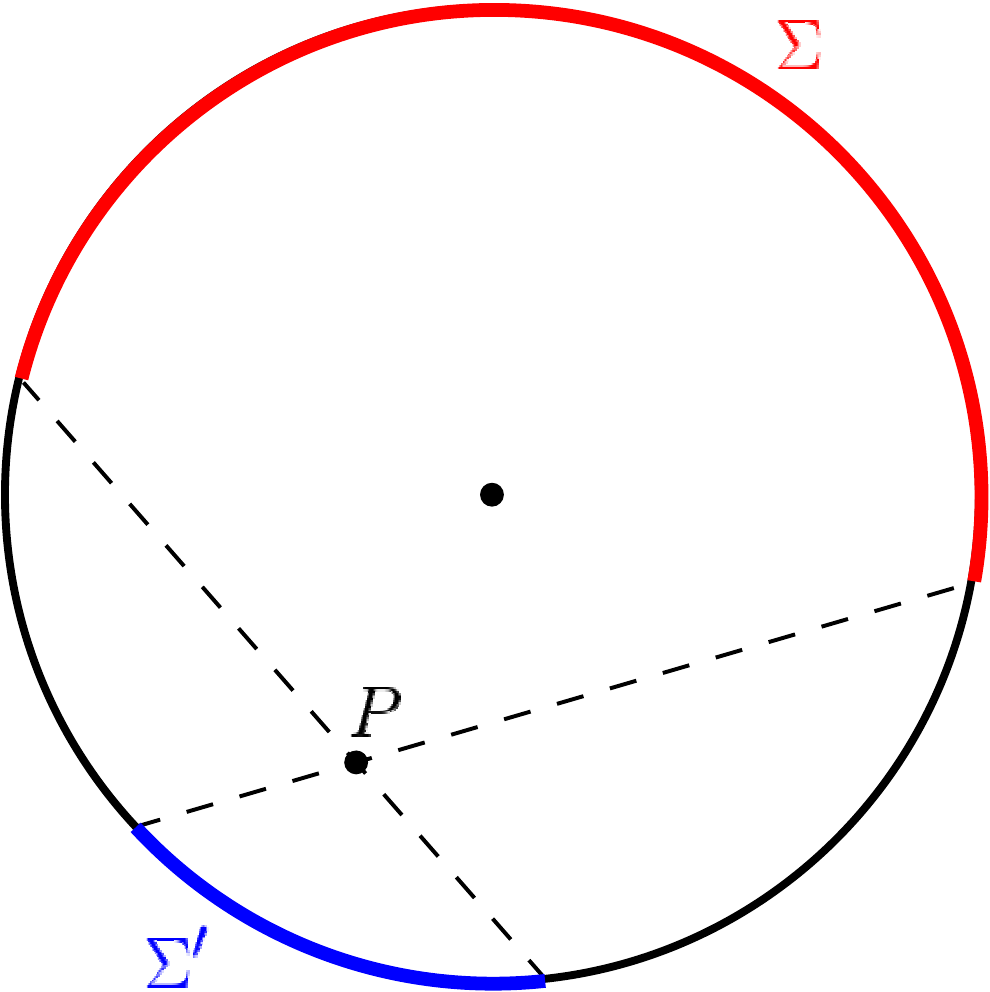}
 \caption{\sl Finding the temperature in a disk just by using a ruler.}\label{2555Hq466Rd7dss8xR9C2L0U1N3zEREkSGkJ0so8oDl3an0eFL}
\end{figure}

We observe that when~$P$ is the origin Theorem~\ref{SCH} reduces to the Mean Value Formula in Theorem~\ref{KAHAR}(ii).
Also, Theorem~\ref{SCH} has a neat geometric interpretation, which was probably of inspiration
for Hermann Schwarz's classic work~\cite{MR0392470}. Namely, suppose that the temperature on the boundary
of the planar unit disk is equal to~$1$ along a given arc~$\Sigma$ and~$0$ elsewhere.
Then, given a point~$P$ in the disk, the temperature at~$P$ can be calculated precisely via
a simple, and purely geometric, protocol: namely, one projects the arc~$\Sigma$ through the focal point~$P$
obtaining a ``conjugated arc''~$\Sigma'$, and the temperature at~$P$ is exactly equal to the length of~$\Sigma'$,
divided by~$2\pi$ (notice indeed that this is the content of~\eqref{DEFUMCD}
in the limit case in which~$f=\chi_\Sigma$), see Figure~\ref{2555Hq466Rd7dss8xR9C2L0U1N3zEREkSGkJ0so8oDl3an0eFL}.

\begin{proof}[Proof of Theorem~\ref{SCH}] Setting~$\varpi:=Q(\omega)$, we exploit the spherical change of variable in~\eqref{DHNETERGD01324u54-2} (recall that here~$n=2$). In this way, we deduce from~\eqref{DEFUMCD} that
$$ u(P)=\fint_{\partial B_1} 
\frac{f(\varpi)( 1-|P|^2) }{|P-\varpi|^2}
\,d{\mathcal{H}}^{1}_\varpi.$$
This expression coincides with the one obtained by the Poisson Kernel representation in
Theorems~\ref{POIBALL1} and~\ref{POIBALL}, therefore~$u$ is harmonic in~$B_1$
and attains the datum~$f$ along~$\partial B_1$.
\end{proof}

\begin{proof}[Another proof of Theorem~\ref{SCH}] The idea of the proof is to reduce to
Theorem~\ref{HARMOu3}. We use the notation in~\eqref{PIUJSMENP} and~\eqref{ELLEeMA}: hence, recalling~\eqref{EELS},
$$ \ell_e(P)= \frac{r_+(e)f(Q_-(e))-r_-(e)f(Q_+(e))
}{r_+(e)-r_-(e)}.$$
Notice also that~$Q(Q_\pm(e))=Q_\mp(e)$.
By~\eqref{DEUGRGBD}  and Theorem~\ref{HARMOu3},
we know that the harmonic function in~$B_1$
with boundary datum equal to~$f$ can be written as
\begin{eqnarray*}&&
v(P):=\fint_{\partial B_1} \ell_e(P)\,d{\mathcal{H}}^{1}_e=
\fint_{\partial B_1} \frac{r_+(e)f(Q_-(e))}{r_+(e)-r_-(e)}\,d{\mathcal{H}}^{1}_e
-\fint_{\partial B_1} \frac{r_-(e)f(Q_+(e))
}{r_+(e)-r_-(e)}\,d{\mathcal{H}}^{1}_e.
\end{eqnarray*}
Our goal is thus to show that~$v(P)=u(P)$, with which the desired result would be proved.
To this end, we exploit~\eqref{VEBDNFKRGZINDMKBVINSKDSA-22}, applied with~$g(\omega):=f(Q(\omega))$
and~\eqref{90KS-098765RYUISOoscuSS0jsdmc} to see that
\begin{eqnarray*}&&
u(P)=\fint_{\partial B_1} f(Q(\omega))\,d{\mathcal{H}}^{1}_\omega=\fint_{\partial B_1} g(\omega)\,d{\mathcal{H}}^{1}_\omega=
\fint_{\partial B_1}
g\left( Q_-(\omega)\right)\, \frac{( r_-(\omega))^2}{1-|P|^2-r_-(\omega)P\cdot \omega}
\,d{\mathcal{H}}^{1}_\omega\\
&&\qquad\qquad=
\fint_{\partial B_1}
f\left( Q(Q_-(\omega))\right)\, \frac{( r_-(\omega))^2}{1-|P|^2-r_-(\omega)P\cdot \omega}
\,d{\mathcal{H}}^{1}_\omega\\&&\qquad\qquad=
\fint_{\partial B_1}
f(Q_+(\omega))\, \frac{( r_-(\omega))^2}{1-|P|^2-r_-(\omega)P\cdot \omega}
\,d{\mathcal{H}}^{1}_\omega\\&&\qquad\qquad=
\fint_{\partial B_1}
f(Q_+(\omega))\, \frac{2r_+(\omega)( r_-(\omega))^2}{(r_+(\omega)-r_-(\omega))(1-|P|^2)}
\,d{\mathcal{H}}^{1}_\omega.
\end{eqnarray*}
This and~\eqref{ANT6} lead to
$$ u(P)=
-\fint_{\partial B_1}
f(Q_+(\omega))\, \frac{2r_-(\omega)}{r_+(\omega)-r_-(\omega)}
\,d{\mathcal{H}}^{1}_\omega=v(P),$$
as desired.\end{proof}

\begin{figure}
  \centering
  \includegraphics[width=.55\linewidth]{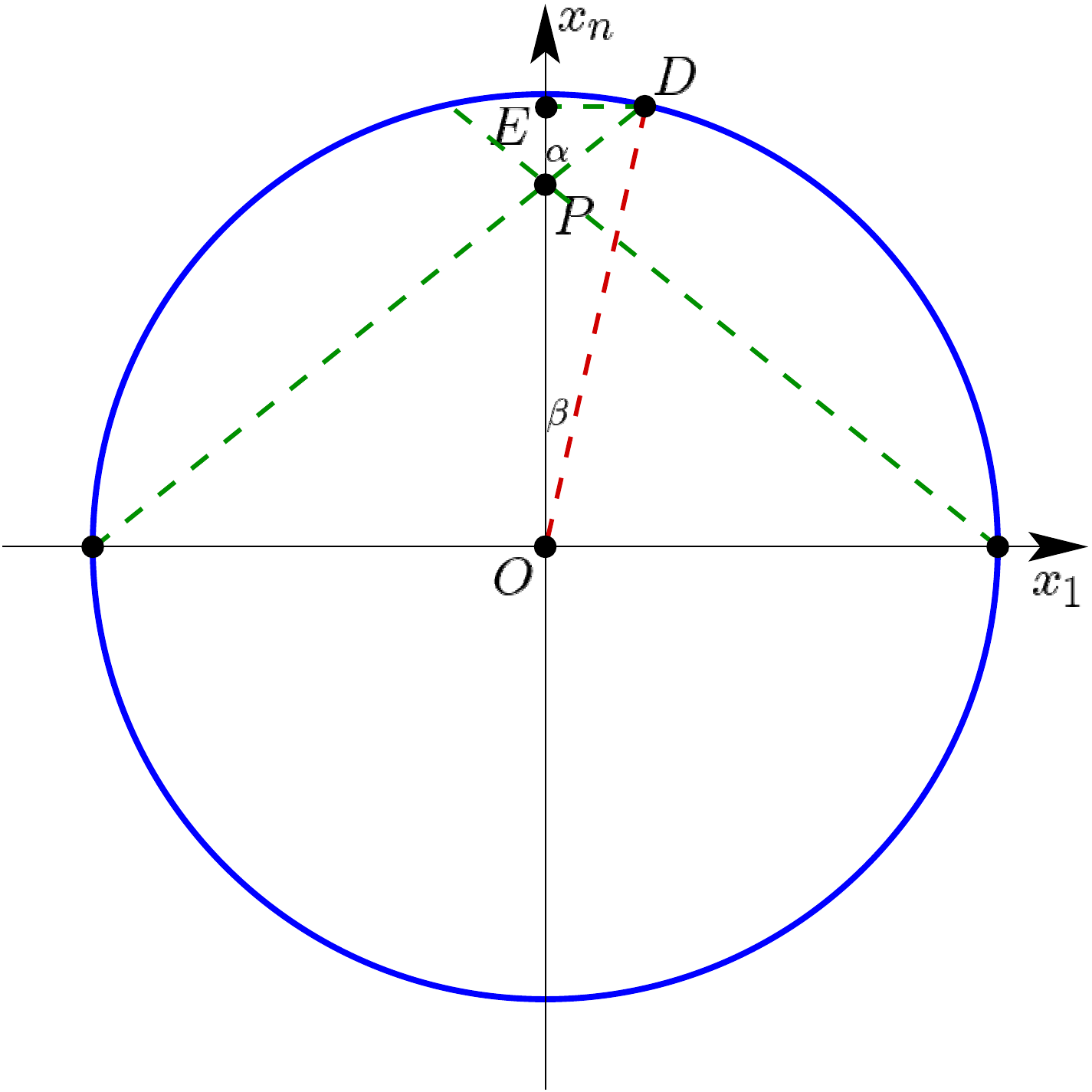}
 \caption{\sl A counterexample to Theorem~\ref{SCH} when~$n\ge3$.}\label{HARHARCOUNTERELAGIJ7soloDItangeFI}
\end{figure}

We stress that Theorem~\ref{SCH} only holds in the plane. To construct a counterexample in dimension~$n\ge3$,
one can consider the function~$f(x_1,\dots,x_n):=\chi_{(-\infty,0)}(x_n)$ and the corresponding
function~$u$ as in~\eqref{DEFUMCD} (technically, since~$f$ is discontinuous,
one should first mollify~$f$ and then pass to the limit). We observe that~$u$ is not harmonic when~$n\ge3$.
Indeed, let~$P:=(0,\dots,0,1-\e)=(1-\e) e_n$, for some small~$\e>0$. Then,
up to a multiplicative constant,
\begin{equation}\label{SPHECAP}
\begin{split}&
{\mbox{$u(P)$ coincides with the surface area
of the spherical cap obtained}}\\&{\mbox{by projecting the lower halfsphere~$\{x_n<0\}\cap\partial B_1$
through the point~$P$.}}\end{split}\end{equation} One can estimate this surface area by elementary geometry, see
Figure~\ref{HARHARCOUNTERELAGIJ7soloDItangeFI}. Indeed, with respect to this picture, if~$\alpha$
denotes the angle~$\widehat{EPD}$, we see that this angle is also equal to the angle~$\widehat{OPe_1}$
and therefore, looking at the triangle with vertices~$O$, $P$ and~$e_1$,
$$ \tan\alpha=\frac{1}{\overline{OP}}=\frac1{1-\e}=1+\e+o(\e).
$$
This leads to
$$\alpha=\arctan(1+\e+o(\e))=\frac\pi4+\frac{\e}2+o(\e).$$
As a consequence,
\begin{eqnarray*}&&
\cos\alpha=\frac{\sqrt{2}}{2}-\frac{\sqrt{2}\,\e}{4}+o(\e)\\
{\mbox{and }}&&\sin\alpha=\frac{\sqrt{2}}{2}+\frac{\sqrt{2}\,\e}{4}+o(\e).
\end{eqnarray*}
Also, applying the Law of cosines to the triangle with vertices~$O$, $P$ and~$D$,
\begin{eqnarray*} &&1=\overline{OD}^2=\overline{OP}^2+\overline{PD}^2
-2\,\overline{OP}\;\overline{PD}\cos(\pi-\alpha)=
(1-\e)^2+\overline{PD}^2
+2(1-\e)\overline{PD}\cos \alpha\\&&\qquad\qquad=
1-2\e+\overline{PD}^2
+\left(\sqrt2-\frac{3\sqrt{2}\;\e}2\right)
\overline{PD}+o(\e).
\end{eqnarray*}
As a result,
$$ \overline{PD}^2
+\left(\sqrt2-\frac{3\sqrt{2}\;\e}2\right)
\overline{PD}-2\e+o(\e)=0.$$
By solving this quadratic equation, we find that
$$ \overline{PD}=
\frac{ \sqrt{4+4\e} + 3 \e - 2 }{2\sqrt2}+o(\e)
=\sqrt{2}\,\e+o(\e).$$
Consequently, if~$\beta$ denotes the angle~$\widehat{POD}$, by applying the Law of sines to
the triangle with vertices~$O$, $P$ and~$D$, we deduce that
$$ \sin\beta=\frac{\overline{PD}\;\sin(\pi-\alpha)}{\overline{OD}}
=\overline{PD}\;\sin\alpha =
\sqrt{2}\,\e\,\left(
\frac{\sqrt{2}}{2}+\frac{\sqrt{2}\,\e}{4}\right)+o(\e)=\e+o(\e).
$$
Accordingly, we have that~$\beta=\e+o(\e)$ and thus, by~\eqref{SPHECAP},
\begin{equation}\label{7ygvBTSYSTUNDFETsm09876-0mdsTTA}
u(P)\le C\e^{n-1},
\end{equation}
for some~$C>0$ depending only on~$n$.

On the other hand, a lower bound for~$u(P)$ can be obtained directly from the Poisson Kernel
of the ball. Namely, in light of Theorems~\ref{POIBALL1} and~\ref{POIBALL},
\begin{eqnarray*}&& u(P)=
\int_{\{y_n<0\}\cap\partial B_1}
\frac{1-|P|^2}{n\,|B_1|\,|P-y|^n}\,d{\mathcal{H}}^{n-1}_y
=\int_{\{y_n<0\}\cap\partial B_1}
\frac{1-(1-\e)^2}{n\,|B_1|\,\big( |y'|^2+(1-\e-y_n)^2\big)^{\frac{n}2}}\,d{\mathcal{H}}^{n-1}_y
\\&&\qquad\qquad\qquad\qquad\ge
\int_{\{y_n<0\}\cap\partial B_1}
\frac{\e}{n\,|B_1|\,4^{\frac{n}2}}\,d{\mathcal{H}}^{n-1}_y=c\e,
\end{eqnarray*}
for a suitable~$c>0$ depending only on~$n$. Comparing this with~\eqref{7ygvBTSYSTUNDFETsm09876-0mdsTTA},
we find that~$c\le C\e^{n-2}$, which, for small~$\e$,
cannot hold true when~$n\ge3$.

\section{Analyticity of harmonic functions}

A direct consequence
of the
Poisson Kernel representation obtained in Theorem~\ref{POIBALL}
is the interior regularity of harmonic functions, as pointed out in the following result:

\begin{theorem} \label{KMS:HAN0-0}
Harmonic functions are real analytic.
\end{theorem}

\begin{proof}
Let~$\Omega\subseteq\R^n$ be an open set and let~$u$ be harmonic in~$\Omega$.
Let~$p\in\Omega$ and pick~$R>0$ such that~$B_R(p)\Subset\Omega$.
Up to a translation, we suppose that~$p=0$ and thus we exploit
Theorems~\ref{POIBALL1}
and~\ref{POIBALL} to write, for every~$x\in B_R$,
\begin{equation}\label{SDfhSUImdxfctaOIKSBLodke} u(x)=\int_{\partial B_R}\frac{(R^2-|x|^2)\,u(y)}{n\,|B_1|\,R\,|x-y|^n}\,d{\mathcal{H}}_y^{n-1}.\end{equation}

Now, if~$y\in\partial B_R$, we 
have that \begin{equation}\label{LMEXJABINTAYUY2}
{\mbox{the function~$B_{R/2}\ni x\mapsto |x-y|^2=|x|^2+|y|^2+2x\cdot y$ is real analytic,}}\end{equation}
since it is a polynomial.

Also, for every~$x\in B_{R/2}$ and~$y\in\partial B_R$,
\begin{equation}\label{LMEXJABINTAYUY3}
|x-y|^2\ge (|y|-|x|)^2\ge \frac{R^2}4\qquad{\mbox{and}}\qquad|x-y|^2\le(|x|+|y|)^2\le \frac{9R^2}4.\end{equation}

Additionally,
\begin{equation}\label{LMEXJABINTAYUY}
{\mbox{the function~$\left(\frac{R^2}8,\frac{9R^2}{2}\right)\ni t\mapsto t^{-\frac{n}2}$
is real analytic.}}\end{equation}
Indeed, we can write~$\sigma:=\frac{16}{R^2}\left(
t-\frac{R^2}{16}\right)
\in(1,72)$, and, applying the Binomial Series, find that
the function~$\left(1+\sigma
\right)^{-\frac{n}2}$ is real analytic, namely
$$ \left(1+\sigma
\right)^{-\frac{n}2}=
\sum_{k=0}^{+\infty} \frac1{k!}\,\prod_{j=0}^{k-1}\left(-\frac{n}2-j\right)\,\sigma^k
.$$
Since~$\sigma^k$ can be written as a Binomial Expansion in powers of~$t$, we have just established~\eqref{LMEXJABINTAYUY}.

Thus, by~\eqref{LMEXJABINTAYUY2}, \eqref{LMEXJABINTAYUY3} and~\eqref{LMEXJABINTAYUY},
we deduce that the function~$B_{R/2}\ni x\mapsto (|x-y|^2)^{-\frac{n}2}=|x-y|^{-n}$
is real analytic (see e.g. the composition result
in~\cite[Proposition~1.6.7]{MR1916029}). In this way, we write
$$ |x-Re_1|^{-n}=\sum_{\alpha\in\N^n} c_\alpha \,x^\alpha$$
and this series is uniformly convergent for~$x\in B_{\rho}$, for some~$\rho\in\left(0,\frac{R}4\right)$.

As a matter of fact, considering a rotation~${\mathcal{R}}_y$ such that~$\frac{y}{|y|}={\mathcal{R}}_y \,e_1$,
and setting~$X:={\mathcal{R}}_y^{-1} x$, for all~$x\in B_{\rho}$ we have that~$X\in B_{\rho}$
and so we obtain the uniformly convergent series expression
$$ |x-y|^{-n}=|{\mathcal{R}}_yX-R\,{\mathcal{R}}_y \,e_1|^{-n}=|X-Re_1|^{-n}=
\sum_{\alpha\in\N^n} c_\alpha\,X^\alpha
=\sum_{\alpha\in\N^n} c_\alpha\,({\mathcal{R}}_y^{-1} x)^\alpha=\sum_{\alpha\in\N^n} C_\alpha(y) x^\alpha
$$
for suitable coefficients~$C_\alpha(y)$.

This and~\eqref{SDfhSUImdxfctaOIKSBLodke} give that, if~$x\in B_\rho$,
$$ u(x)=\int_{\partial B_R}
\sum_{\alpha\in\N^n}
\frac{C_\alpha(y)\,(R^2-|x|^2)\,u(y)\,x^\alpha}{n\,|B_1|\,R}\,d{\mathcal{H}}_y^{n-1}
=
\sum_{\alpha\in\N^n} C_\alpha\,(R^2-|x|^2)\,x^\alpha,$$
where
$$ C_\alpha:=
\int_{\partial B_R}
\frac{C_\alpha(y)\,u(y)}{n\,|B_1|\,R}\,d{\mathcal{H}}_y^{n-1}
.$$
This shows that~$u$
can be written as a Taylor series in a neighborhood the origin and thus it establishes the desired result.
\end{proof}

For a different proof of
Theorem~\ref{KMS:HAN0-0}
see e.g.~\cite[pages 31-32]{MR1625845}.

\section{The Harnack Inequality}\label{HARNAFOR}

As a consequence of the
Poisson Kernel obtained in Theorem~\ref{POIBALL}, we have the 
following
result that provides a sharp and explicit
control on the oscillation of harmonic 
functions:\index{Harnack Inequality}

\begin{theorem}\label{PRE:HA}
Let~$x_0\in\R^n$ and~$R>r>0$.
Assume that~$u$ is nonnegative and harmonic in~$B_R(x_0)$.
Then, for every~$x\in B_r(x_0)$,
\begin{equation}\label{0-PREJA} \left(\frac{R}{R+r}\right)^{n-2}\,
\frac{R-r}{R+r}\,u(x_0)\le
u(x)\le\left(\frac{R}{R-r}\right)^{n-2}\,\frac{R+r}{R-r}\,u(x_0).\end{equation}\end{theorem}

\begin{proof} Up to a translation,
we suppose~$x_0=0$. We also take~$\rho\in(r,R)$.
Thus, by Theorems~\ref{POIBALL1} and~\ref{POIBALL},
and recalling~\eqref{B1},
for every~$x\in B_\rho$ we have that
\begin{equation}\label{KS MIJMNSINSTHNMEGJSUISK90OS}
\begin{split}&
u(x)=\int_{\partial B_\rho}\frac{(\rho^2-|x|^2)\,u(y)}{n\,|B_1|\,\rho\,|x-y|^n}\,
d{\mathcal{H}}^{n-1}_y=
\int_{\partial B_\rho}\frac{(\rho^2-|x|^2)\,u(y)}{{\mathcal{H}}^{n-1}(\partial B_1)\,\rho\,|x-y|^n}\,
d{\mathcal{H}}^{n-1}_y\\&\qquad\qquad=\rho^{n-2}
\fint_{\partial B_\rho}\frac{(\rho^2-|x|^2)\,u(y)}{
|x-y|^n}\,
d{\mathcal{H}}^{n-1}_y.\end{split}
\end{equation}
Additionally, if~$x\in B_r$ and~$y\in\partial B_\rho$,
\begin{eqnarray*}&&\frac{\rho^2-r^2}{ (\rho+r)^n}\le
\frac{\rho^2-|x|^2}{ (|y|+|x|)^n}\le\frac{\rho^2-|x|^2}{ |x-y|^n}\\
{\mbox{and }}&&\frac{\rho^2-|x|^2}{ |x-y|^n}\le
\frac{\rho^2-|x|^2}{ (|y|-|x|)^n}=\frac{\rho^2-|x|^2}{ (\rho-|x|)^n}=
\frac{\rho+|x|}{ (\rho-|x|)^{n-1}}\le\frac{\rho+r}{ (\rho-r)^{n-1}}.
\end{eqnarray*}
These inequalities, \eqref{KS MIJMNSINSTHNMEGJSUISK90OS} and
the 
Mean Value Formula in Theorem~\ref{KAHAR}(ii) give that,
if~$x\in B_r$,
\begin{eqnarray*}&& \left(\frac{\rho}{\rho+r}\right)^{n-2}\,
\frac{\rho-r}{\rho+r}\,u(0)
=
\frac{\rho^{n-2}(\rho^2-r^2)}{ (\rho+r)^n}
\,u(0)=
\frac{\rho^{n-2}(\rho^2-r^2)}{ (\rho+r)^n}
\fint_{\partial B_\rho} u(y)\,
d{\mathcal{H}}^{n-1}_y\\&&\qquad
\le
\rho^{n-2}
\fint_{\partial B_\rho}\frac{(\rho^2-|x|^2)\,u(y)}{
|x-y|^n}\,
d{\mathcal{H}}^{n-1}_y=u(x)
\end{eqnarray*}
and also
\begin{eqnarray*}
u(x)\le
\frac{\rho^{n-2}(\rho+r)}{
(\rho-r)^{n-1}}\,
\fint_{\partial B_\rho}u(y)\,
d{\mathcal{H}}^{n-1}_y=\frac{\rho^{n-2}(\rho+r)}{
(\rho-r)^{n-1}}\,u(0)=
\left(\frac{\rho}{\rho-r}\right)^{n-2}\,\frac{\rho+r}{\rho-r}\,u(0).
\end{eqnarray*}
Sending~$\rho\nearrow R$ we obtain the desired result.\end{proof}

We point out that the estimate in Theorem~\ref{PRE:HA}
is optimal.
Indeed, when~$n=1$, the function~$u(x):=x-x_0+R$
satisfies$$ \frac{R-r}R\,u(x_0)=R-r=u(x_0-r)$$
and$$ u(x_0+r)=R+r=\frac{R+r}R\,u(x_0),$$
attaining the bounds presented in~\eqref{0-PREJA}
when~$n=1$. Moreover, when~$n\ge2$,
in light of Theorem~\ref{POIBALL},
given~$\e>0$, we can consider the Poisson Kernel~$P_\e(\cdot,\cdot)$
for~$B_{R+\e}$ and define
\begin{equation}\label{HESEM} u(x):=P_\e(x,(R+\e)e_1)=
\frac{(R+\e)^2-|x|^2}{n\,|B_1|\,(R+\e)\,|x-(R+\e)e_1|^n}.\end{equation}
We observe that~$u$ is harmonic in~$B_R$, thanks to~\eqref{PfghjklX99ERHILD}.
Moreover, if~$x:=re_1$,
$$ \frac{u(x)}{u(0)}=
\frac{\big((R+\e)^2-r^2\big)(R+\e)^n}{ (R+\e)^2
|r-(R+\e)|^n},$$
which, as~$\e\searrow0$, attains the right hand side in~\eqref{0-PREJA}.

Similarly, if~$x:=-re_1$,
$$ \frac{u(x)}{u(0)}=
\frac{\big((R+\e)^2-r^2\big)(R+\e)^n}{ (R+\e)^2
|r+(R+\e)|^n},$$
which, as~$\e\searrow0$, attains the left hand side in~\eqref{0-PREJA}.
These observations show that the bounds in~\eqref{0-PREJA}
are sharp and cannot be improved.\medskip

\begin{figure}
  \centering
  \includegraphics[width=.2\linewidth]{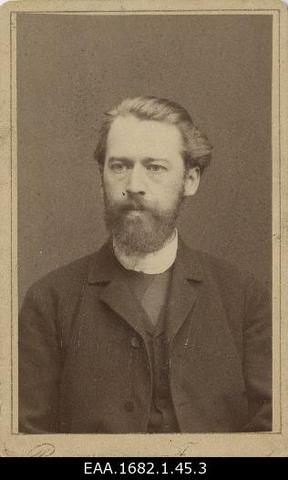}
 \caption{\sl Carl Gustav Axel Harnack (Public Domain image from \label{HARAXELHARLROAGIJ7soloDItangeFIPAGI}
 Wikipedia).}\label{HARAXELHARLROAGIJ7soloDItangeFI}
\end{figure}

As a straightforward consequence\footnote{We mention that an alternative
proof of Corollary~\ref{HAR:BAL} that does not rely
on Theorem~\ref{PRE:HA} can be obtained directly from the Mean Value Formula
in Theorem~\ref{KAHAR}(iii): namely, given two points~$p$, $q\in B_{R/8}(x_0)$,
we have that~$B_{R/8}(q)\subseteq B_{R/2}(p)\subseteq B_{3R/4}(x_0)$ and therefore, for a nonnegative and harmonic~$u$
in~$B_R(x_0)$,
$$ u(p)=\fint_{B_{R/2}(p)}u(x)\,dx\ge\frac1{|B_{R/2}|}\int_{B_{R/8}(q)}u(x)\,dx
=\frac{|B_{R/8}|}{|B_{R/2}|}\fint_{B_{R/8}(q)}u(x)\,dx=\frac{|B_{1/8}|}{|B_{1/2}|}\,u(q)$$
and consequently
$$ \sup_{B_{R/8}(x_0)}u\le\frac{|B_{1/2}|}{|B_{1/8}|}\,\inf_{B_{R/8}(x_0)}u,$$
from which Corollary~\ref{HAR:BAL} would follow from a covering argument.

In our presentation, the relevance of Theorem~\ref{PRE:HA}
is however to produce sharp constants
in~\eqref{0-PREJA}.}
of Theorem~\ref{PRE:HA},
one obtains the following result, which is often
referred to with the name of\footnote{Results
of this type are named after the author of the pioneer work~\cite{GRUND}. The gruesome fate
of Carl Gustav Axel Harnack (mathematician) seems to be that of being too often mistaken
for his twin brother Carl Gustav Adolf von Harnack (theologian, historian, author of
many religious publications, signatory of a public statement 
in support of German military actions in World War I and ennobled with the addition of von to his name;
see also footnote~\ref{gHARAXELHARLROAGIJ7soloDItangeFIPAGI} on page~\pageref{gHARAXELHARLROAGIJ7soloDItangeFIPAGI}).
For a while, Wikipedia erroneously displayed a picture of Adolf in place of that of Axel:
the error was spotted since Adolf's picture was representing a rather senior person,
while Axel prematurely died at the age of~37 after suffering from health problems
(Adolf  long outlived him and died at the age of~79, confirming that personal interests and physical conditions
are not always akin between twin brothers).

Quite surprisingly for us, the prestigious
Harnack Medal of the Max Planck Society is named after Adolf, not after Axel.
This is possibly due to the fact that Adolf was one of the founders, as well as the
first president, of the
Kaiser-Wilhelm-Gesellschaft zur F\"orderung der Wissenschaften
(Society for the Advancement of Science named in honor of the German Emperor
Wilhelm~II)
whose functions were later taken over by the Max Planck Society.

To partially compensate
his sorry fate,
a picture of Axel is given here in Figure~\ref{HARAXELHARLROAGIJ7soloDItangeFI}.}
Harnack Inequality:\index{Harnack Inequality}

\begin{corollary}\label{HAR:BAL}
There exists a constant~$C>1$, depending only on~$n$,
such that for every~$x_0\in\R^n$, every~$R>0$ and every
nonnegative, harmonic function~$u$ in~$B_R(x_0)$, we have that
$$\sup_{B_{R/2}(x_0)}u\le C\,\inf_{B_{R/2}(x_0)}u.$$
\end{corollary}

\begin{proof} Take~$r:=R/2$ in Theorem~\ref{PRE:HA}.
\end{proof}

A more general version of the Harnack Inequality in Corollary~\ref{HAR:BAL}
can be obtained by a covering argument for arbitrary connected
domains (with a constant dependent on the domain):

\begin{figure}
  \centering
  \includegraphics[width=.55\linewidth]{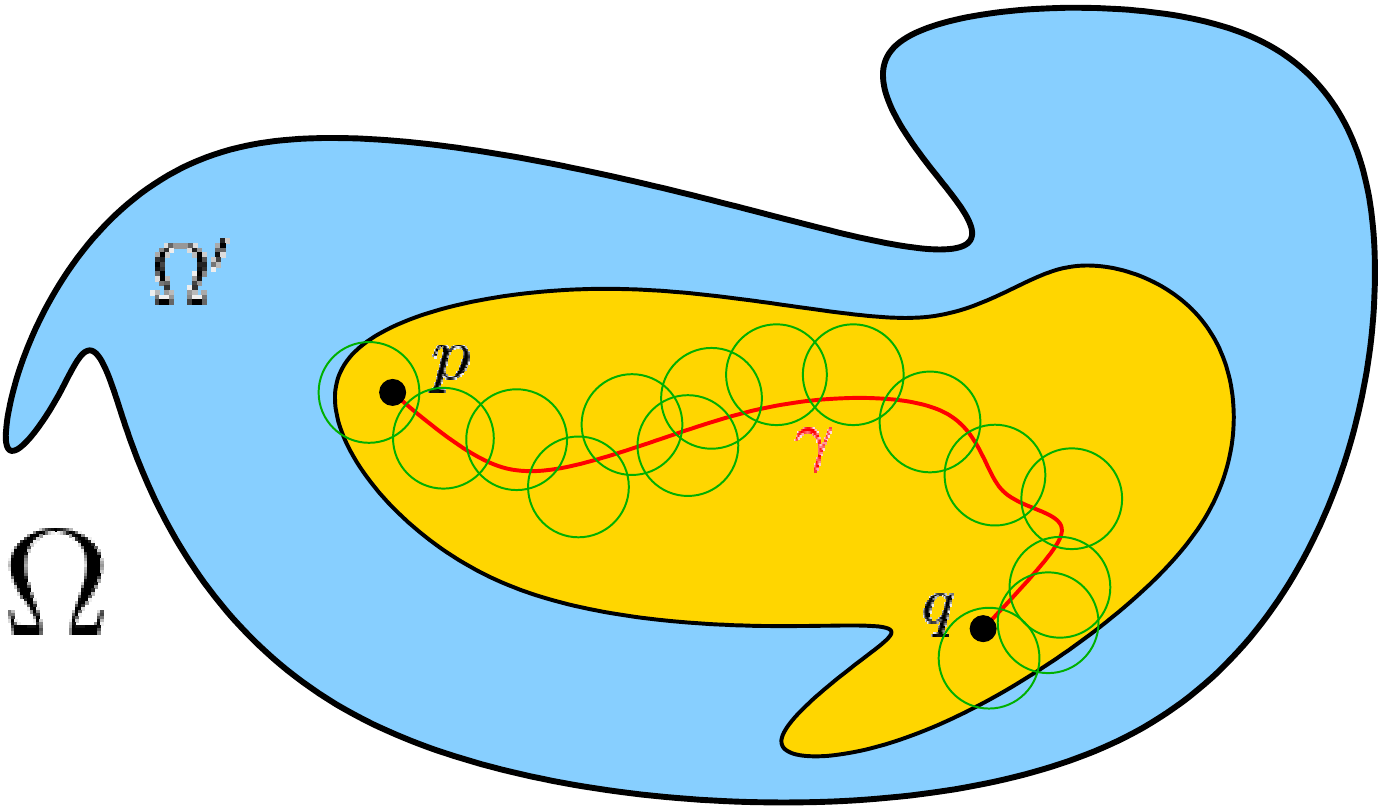}
 \caption{\sl Proof of Corollary~\ref{NSDio}.}\label{HARHARLAGIJ7soloDItangeFI}
\end{figure}

\begin{corollary}\label{NSDio}
Let~$\Omega\subseteq\R^n$ be open.
Assume that~$u$ is nonnegative\index{Harnack Inequality}
and harmonic in~$\Omega$ and let~$\Omega'\Subset\Omega$
be open, connected and bounded.

Then, there exists a constant~$C>1$ depending only on~$n$,
$\Omega'$ and~$\Omega$ such that
$$\sup_{\Omega'}u\le C\,\inf_{\Omega'}u.$$
\end{corollary}

\begin{proof} 
We let~$D$ be the distance between~$\Omega'$ and~$\partial \Omega$
and set~$R:=D/4$ (if~$D$ is finite; if~$D=+\infty$ just take~$R:=1$).
Let~$p$, $q\in\overline{\Omega'}$ be such
that~$\max_{\overline{\Omega'}}u=u(p)$ and~$\min_{\overline{\Omega'}}u=u(q)$.
By the connectedness of~$\Omega'$,
we can take a closed arc~$\gamma$
joining~$p$ and~$q$ that lies in~$\Omega'$. 
Since
$$ \overline{\Omega'}\subseteq\bigcup_{x\in\Omega'} B_{R/2}(x),$$
by compactness we can find~$x^{(1)},\dots,x^{(N)}\in\Omega'$ such that
$$ \Omega'\subseteq\bigcup_{i=1}^N B_{R/2}(x^{(i)}),$$
with~$N$ depending only on~$\Omega$ and~$\Omega'$.
We let~$i_1,\dots,i_{N'}\in\{1,\dots,N\}$ be such that~$\gamma\subseteq B_{R/2}(x^{(i_1)})\cup\dots\cup B_{R/2}(x^{(i_{N'})})$.
Without loss of generality, we can think that all the~$x^{(i_j)}$ are different.
Additionally, we can conveniently reorder these balls
(and possibly disregard balls that do not contribute to the covering of~$\gamma$).
Namely, we define~$B^{(1)}:=B_{R/2}(x^{(i_1)})$.
If~$\gamma$ is all contained in~$B^{(1)}$ we stop, otherwise, by the continuity of the arc~$\gamma$,
there must be a point~$P_1$ of~$\gamma$ on~$\partial B^{(1)}$.
Hence there must exist a ball of the collection~$\{B_{R/2}(x^{(i_2)}),\dots,B_{R/2}(x^{(i_M)})$
that covers~$P_1$ and we call~$B^{(2)}$ this ball.
Again, if~$\gamma$ is contained in~$B^{(1)}\cup B^{(2)}$ we stop,
otherwise we consider a point~$P_2\in\partial(B^{(1)}\cup B^{(2)})$
and a ball of the remaining collection, that we denote by~$B^{(3)}$,
which contains~$P_3$. Recursively, we have found balls~$B^{(1)},\dots,B^{(M)}$
such that~$M\le N'\le N$, with
$$ \gamma\subseteq B^{(1)}\cup \dots\cup B^{(M)}$$
and with the additional property that, for all~$j\in\{1,\dots,M-1\}$,
$$ \big( B^{(1)}\cup \dots\cup B^{(j)}\big)\cap B^{(j+1)}\ne\varnothing,$$
see Figure~\ref{HARHARLAGIJ7soloDItangeFI}.

In particular, taking~$\zeta_j$ in the intersection above
and applying Corollary~\ref{HAR:BAL}, we see that
\begin{equation}\label{78:MA:SPE:MA:LIS:T56ter:090704020-01}
\sup_{B^{(j+1)}}u\le C \inf_{B^{(j+1)}}u \le C u(\zeta_j)\le C\sup_{B^{(1)}\cup \dots\cup B^{(j)}}u
\end{equation}
and
\begin{equation}\label{78:MA:SPE:MA:LIS:T56ter:090704020-02}
\inf_{B^{(1)}\cup \dots\cup B^{(j)}}u\le u(\zeta_j)\le\sup_{B^{(j+1)}}u\le C \inf_{B^{(j+1)}}u,
\end{equation}
for some~$C>1$ depending only on~$n$.

Now we claim that, for every~$j\in\{1,\dots,M\}$,
\begin{equation}\label{78:MA:SPE:MA:LIS:T56ter:090704020-03}
\sup_{B^{(1)}\cup \dots\cup B^{(j)}}u\le C^{j-1}\sup_{B^{(1)}}u.
\end{equation}
To prove this, we argue by induction. When~$j=1$ the claim is obvious, hence we suppose the claim
true for some index~$j\in\{1,\dots,M-1\}$ and we aim at establishing it for the index~$j+1$.
To this end, we make use of~\eqref{78:MA:SPE:MA:LIS:T56ter:090704020-01} and of the inductive assumption
to find that
\begin{eqnarray*}
\sup_{B^{(1)}\cup \dots\cup B^{(j+1)}}u&=&\max\left\{
\sup_{B^{(1)}\cup \dots\cup B^{(j)}}u,\;
\sup_{B^{(j+1)}}u
\right\}\\&\le&C\sup_{B^{(1)}\cup \dots\cup B^{(j)}}u\\&\le&C^{j}\sup_{B^{(1)}}u.
\end{eqnarray*}
This completes the inductive step and establishes~\eqref{78:MA:SPE:MA:LIS:T56ter:090704020-03}.

In analogy with~\eqref{78:MA:SPE:MA:LIS:T56ter:090704020-03},
we claim that
\begin{equation}\label{78:MA:SPE:MA:LIS:T56ter:090704020-04}
\inf_{B^{(1)}\cup \dots\cup B^{(j)}}u\ge C^{1-j}\inf_{B^{(1)}}u.
\end{equation}
Once again, for this we argue by induction. Since the desired result is obvious when~$j=1$,
we suppose that the claim in~\eqref{78:MA:SPE:MA:LIS:T56ter:090704020-04} holds
true for some index~$j\in\{1,\dots,M-1\}$ and we aim at establishing it for the index~$j+1$.
For this, using~\eqref{78:MA:SPE:MA:LIS:T56ter:090704020-02} and the inductive assumption
we see that
\begin{eqnarray*}
\inf_{B^{(1)}\cup \dots\cup B^{(j+1)}}u&=&\min\left\{
\inf_{B^{(1)}\cup \dots\cup B^{(j)}}u,\;
\inf_{B^{(j+1)}}u
\right\}\\&\ge&C^{-1}\inf_{B^{(1)}\cup \dots\cup B^{(j)}}u\\&\ge&C^{-j}\inf_{B^{(1)}}u,
\end{eqnarray*}
which finishes the proof of~\eqref{78:MA:SPE:MA:LIS:T56ter:090704020-04}.

Now we take~$j_q\in\{1,\dots,M\}$ such that~$q\in B^{(j_q)}$
and~$j_p\in\{1,\dots,M\}$ such that~$p\in B^{(j_p)}$. Combining~\eqref{78:MA:SPE:MA:LIS:T56ter:090704020-03}
and~\eqref{78:MA:SPE:MA:LIS:T56ter:090704020-04}, and applying Corollary~\ref{HAR:BAL}
once more, we thus deduce that
\begin{equation*}\begin{split}& \max_{\overline{\Omega'}}u=u(p)\le
\sup_{B^{(j_p)}}u\le\sup_{B^{(1)}\cup\dots\cup B^{(j_p)}}u
\le C^{j_p-1}\sup_{B^{(1)}}u
\le C^{j_p}\inf_{B^{(1)}}u\\&\qquad
\le C^{j_p+j_q-1}\inf_{B^{(1)}\cup\dots\cup B^{(j_q)}}u
\le C^{j_p+j_q-1}\inf_{B^{(j_q)}}u
\le C^{2N-1} u(q)
=C^{2N-1}\,\min_{\overline{\Omega'}}u
.\qedhere\end{split}\end{equation*}
\end{proof}

Harnack Inequality, in its various forms, is certainly of paramount
importance, since, roughly speaking, it states that, in the interior, all the values of harmonic functions
are essentially comparable (and this, in particular,
provides a quantitative version of the Maximum Principle in Theorem~\ref{STRONGMAXPLE1}(iii)). \medskip

Here we present a variant of Corollary~\ref{S-OS-OSSKKHEPARHA}
that exploits the Harnack Inequality and it is sometimes
called the Harnack Convergence Theorem:\index{Harnack Convergence Theorem}

\begin{corollary}
Let~$\Omega\subseteq\R^n$ be open and connected.
Let~$u_k$ be a monotone increasing sequence of harmonic functions
in~$\Omega$ and suppose that there exists~$y\in\Omega$
such that the sequence~$u_k(y)$ is bounded from above.

Then, the sequence~$u_k$
converges locally uniformly in~$\Omega$ to a
harmonic function.
\end{corollary}

\begin{proof} Let~$\Omega'\Subset\Omega$ be open and bounded.
Let~$\Omega''\Subset\Omega$ be connected, open, bounded
and such that~$\Omega'\cup\{y\}\subseteq\Omega''$.

By construction, the sequence~$u_k(y)$
is convergent. In particular, given~$\e>0$, there exists~$k_\e\in\N$
such that for all~$k$, $k'\ge k_\e$, say with~$k\ge k'$, we have that~$|u_k(y)-u_{k'}(y)|\le\e$.

We observe in addition that the function~$u_k-u_{k'}$ is
harmonic and nonnegative.
Thus, for every~$x\in\Omega''$,
$$ |u_k(x)-u_{k'}(x)|=u_k(x)-u_{k'}(x)\le \sup_{\Omega''}(u_k-u_{k'})\le
C\,\inf_{\Omega''}(u_k-u_{k'})\le C\,(u_k(y)-u_{k'}(y))
\le C\e,$$
as long as~$k\ge k'\ge k_\e$,
thanks to Corollary~\ref{NSDio}.

For this reason, we have that~$u_k$ converges uniformly in~$\Omega''$
(and therefore in~$\Omega'$) to some function~$u$.
Hence, recalling 
Corollary~\ref{S-OS-OSSKKHEPARHA}, we obtain that~$u$
is harmonic.\end{proof}

A boundary counterpart of the Harnack Inequality will be briefly
discussed in Section~\ref{BOUNAHSRGA9876543}.

\section{The Hopf Lemma}\label{HOPFSECTH}

Concerning the boundary behavior of harmonic functions,
we recall here a (simple version of a) classical result,
often referred to with the name of Hopf Lemma:\index{Hopf Lemma}

\begin{lemma}\label{JJS:PA}
Let~$\Omega$ be a bounded, connected, open set with~$C^1$ boundary.
Let~$u\in C^2(\Omega)\cap C^1(\overline\Omega)$
be harmonic in~$\Omega$.
Assume that~$x_0\in\partial\Omega$ is such that
$$ \max_{\partial\Omega} u=u(x_0).$$
Suppose also that there exists a ball~$B\subseteq\Omega$
such that~$x_0\in\partial B$.

Then, either~$ \partial_\nu u(x_0)>0$ or~$u$ is constant in~$\Omega$.\end{lemma}

\begin{figure}
  \centering
  \includegraphics[width=.45\linewidth]{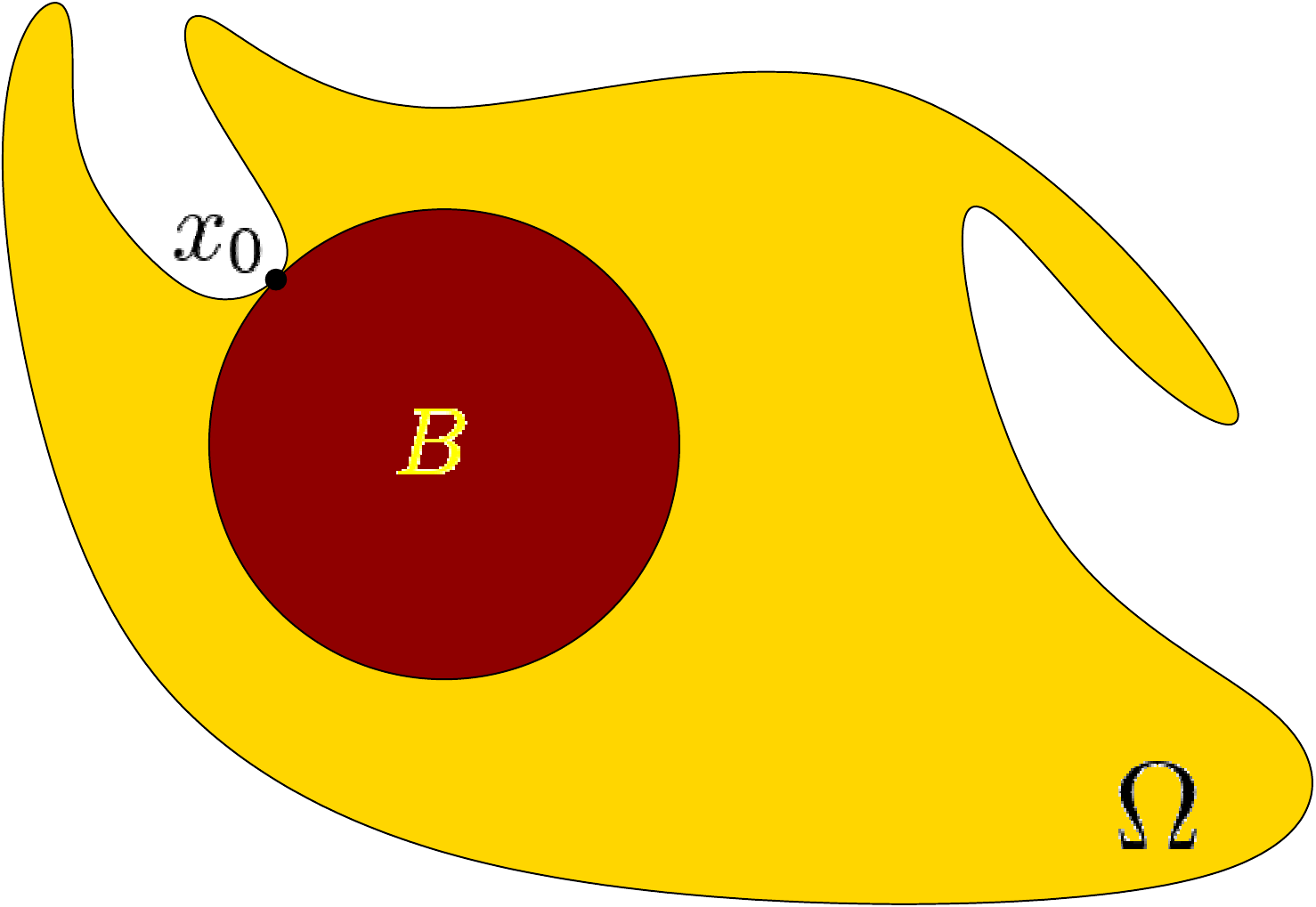}
 \caption{\sl The interior ball condition\index{interior ball condition} in Lemma~\ref{JJS:PA}.}\label{HOPFLAGIJ7soloDItangeFIHOP}
\end{figure}

\begin{proof} By the
Weak Maximum Principle in Corollary~\ref{WEAKMAXPLE}(iii), we already know that~$u(x_0)=\max_{\partial\Omega} u
=\max_{\overline\Omega} u$.
Hence~$u(x_0)=\max_{\partial B} u$.

Consequently, if~$\max_{\partial B} u
=\min_{\partial B} u$, it follows that~$u$ is constantly equal to~$u(x_0)$ along~$\partial B$.
In particular, $u$ attains an interior maximum in~$\Omega$
and accordingly, by the Strong Maximum Principle in Theorem~\ref{STRONGMAXPLE1}(iii),
we obtain that~$u$ is constant and we are done.

For this reason, we can suppose that~$\max_{\partial B} u>\min_{\partial B} u$.
As a result, by the continuity of~$u$, we can take~$x_\star\in\partial B$
such that~$u(x_\star)=\min_{\partial B} u$ and also~$\rho\in\left(0,\frac{|x_0-x_\star|}{2}\right)$ and~$a>0$
such that for all~$y\in (\partial B)\cap B_\rho(x_\star)$
we have
\begin{equation}\label{JASweLLASMMD03}
u(y)+a\le\max_{\partial B}u=
\max_{\partial\Omega} u.
\end{equation}

Now,
up to a translation, we can assume that the ball~$B$
is centered at the origin, hence we have that~$B=B_R$ for some~$R>0$.
Furthermore, up to a rotation, since~$x_0\in\partial B=\partial B_R$,
we can suppose that~$x_0=Re_1$. Then, the 
exterior normal of~$B_R$
at~$x_0$ coincides with the exterior
normal of~$\Omega$ at~$x_0$ and thus, without any ambiguity,
we can write~$\nu(x_0)=e_1$.
Let also~$t>0$, to be taken as small as we wish in the following,
and~$x_t:=x_0-t\nu(x_0)=(R-t)e_1$.
Thus, we exploit the Poisson Kernel representation in Theorems~\ref{POIBALL1} and~\ref{POIBALL},
as well as~\eqref{JASweLLASMMD03}
and obtain that
\begin{eqnarray*}&&
u(x_0)-u(x_t)=\max_{\partial\Omega} u-\int_{\partial B_R} u(y)\,P(x_t,y)\,d{\mathcal{H}}^{n-1}_y
=\int_{\partial B_R} \left(\max_{\partial\Omega} u-u(y)\right)\,P(x_t,y)\,d{\mathcal{H}}^{n-1}_y\\
&&\qquad=
\frac{R^2-|x_t|^2}{n\,|B_1|\,R}
\int_{\partial B_R} 
\left(\max_{\partial\Omega} u-u(y)\right)\,\frac{d{\mathcal{H}}^{n-1}_y}{|x_t-y|^n}\ge
\frac{a\,(R^2-|x_t|^2)}{n\,|B_1|\,R}
\int_{(\partial B_R)\cap B_\rho(x_\star)} 
\frac{d{\mathcal{H}}^{n-1}_y}{|x_t-y|^n}.
\end{eqnarray*}
Hence, since, if~$y\in B_\rho(x_\star)$,
\begin{eqnarray*}&& |x_t-y|=|x_0-t\nu(x_0)-y|\ge |x_0-t\nu(x_0)-x_\star|-|x_\star-y|\ge
|x_0-x_\star|-t-\rho\\&&\qquad\qquad\qquad\ge\frac{|x_0-x_\star|}2-t\ge\frac{|x_0-x_\star|}4\end{eqnarray*}
and
$$ R^2-|x_t|^2=R^2-(R-t)^2=2Rt-t^2\ge Rt$$
as long as~$t$ is small enough, we conclude that
\begin{eqnarray*}&&
u(x_0)-u(x_t)\ge
\frac{4^n at}{n\,|B_1|}
\int_{(\partial B_R)\cap B_\rho(x_\star)} 
\frac{d{\mathcal{H}}^{n-1}_y}{|x_0-x_\star|^n}=\frac{4^n a\,
{\mathcal{H}}^{n-1}((\partial B_R)\cap B_\rho(x_\star))\,
t}{n\,|B_1|\,|x_0-x_\star|^n}.
\end{eqnarray*}
Dividing by~$t$ and sending~$t\searrow0$ we obtain the desired result.
\end{proof}

Different proofs
of Lemma~\ref{JJS:PA}
can be obtained by constructing suitable barriers
and utilizing the Maximum Principle, see e.g.~\cite{MR1814364}
(as a matter of fact, we will use a barrier method here to prove a useful generalization
of Lemma~\ref{JJS:PA} as given in Lemma~\ref{MS-LEHOPF2-M}).

We stress that the interior ball condition at~$x_0$ in Lemma~\ref{JJS:PA}
(that is, the existence of the ball~$B$ tangent from inside to~$\partial\Omega$
and containing~$x_0$ on its boundary, see Figure~\ref{HOPFLAGIJ7soloDItangeFIHOP}) can be relaxed
but\footnote{Also, the assumption that the boundary of~$\Omega$ is of class~$C^1$ cannot be removed from Lemma~\ref{JJS:PA}.
Not only this assumption is needed to properly define the exterior normal~$\nu$, but it is also necessary
to avoid cases in which~$u$ could be smooth in the vicinity of the point~$x_0$, but~$\nabla u(x_0)=0$.
As an example, one can consider~$\Omega:=(0,1)^2\subset\R^2$, $x_0:=(0,0)$ and~$u(x,y):=-xy$, since in this situation one has that~$x_0\in\partial\Omega$,
$u(x_0)=0=\max_{\overline\Omega}u$, but~$\nabla u(x_0)=0$.}
not completely dropped, see e.g.~\cite{MR3513140}
and the references therein (as well as Appendix~\ref{MAHH}
and Lemma~\ref{SerrinsCornerLemma} in this set of notes) for further details.
We also remark that this interior ball condition is satisfied if the boundary
of~$\Omega$ is of class~$C^{1,1}$, see e.g.~\cite[Lemma~A.1]{MR3436398}.
See also~\cite{MR3012036, MR2356201} and the references
therein
for additional details on the
Hopf Lemma and on related topics.\medskip

\begin{figure}
  \centering
  \includegraphics[width=.45\linewidth]{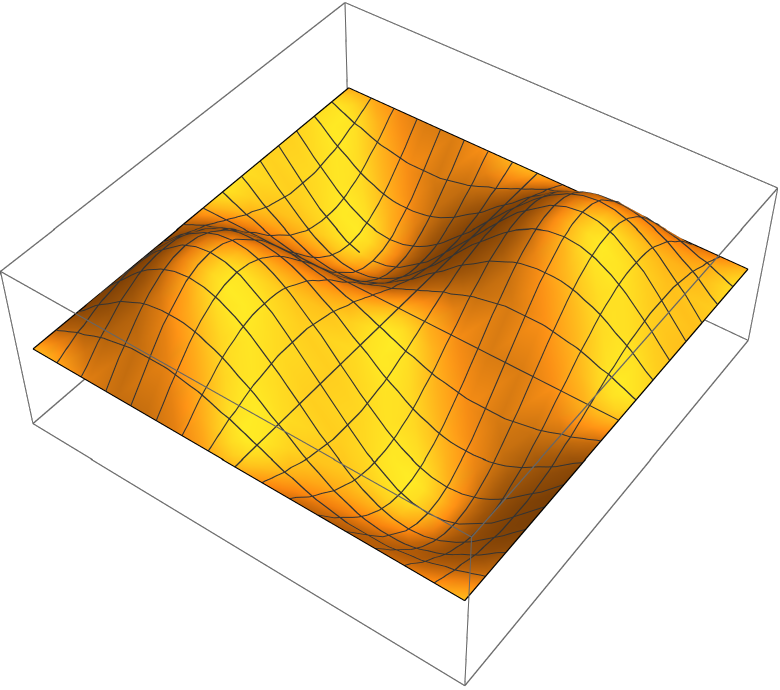}$\qquad$
  \includegraphics[width=.3\linewidth]{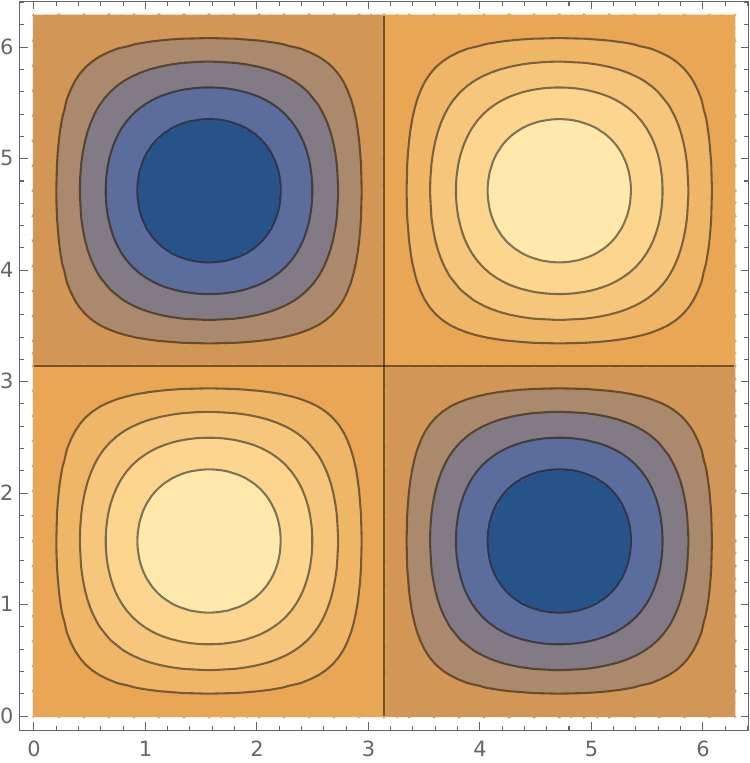}
 \caption{\sl Graph and level sets of the function~$u(x_1,x_2)=\sin x_1\sin x_2$.}\label{LwAS1IewN24XGIJ7soloD354hItangeFIrfdg}
\end{figure}

In terms of generality of the methods
used here, it is worth pointing out that the Strong Maximum Principle
in~\eqref{STRONGMAXPLE1} fails if the Laplacian is replaced by more general operators
with zero order terms, that is solutions of equations of the form~$\Delta u(x)+c(x)u(x)=0$
may develop interior minima and maxima, even for ``nice'' coefficients~$c$
(as an example, we observe that the function~$\R\ni x\mapsto u(x):=\cos x$
satisfies~$\Delta u+u=0$ in~$\R$ and it possesses interior minima and maxima; another example
is provided by the function~$\R^2\ni x=(x_1,x_2)\mapsto u(x)=\sin x_1\sin x_2$
which satisfies~$\Delta u+2u=0$ in~$\R^2$, see Figure~\ref{LwAS1IewN24XGIJ7soloD354hItangeFIrfdg}).
However, a variation of Theorem~\ref{STRONGMAXPLE1} holds true under suitable sign assumptions
(e.g. on the coefficient~$c$ or on the solution).
In this spirit, we point out here a result that will be employed in Section~\ref{DEDXOJD-GIDANSNI}.

\begin{theorem}\label{MAXPLECON45678CG-PM}
Let~$\Omega\subseteq\R^n$ be open and connected.

Let~$a_{ij}$, $b_i$, $c \in C(\overline\Omega)$. Assume that\index{Maximum Principle}
for every~$\xi=(\xi_1,\dots,\xi_n)\in\partial B_1$
\begin{equation}\label{ELLIPTIC-PERFVDJ}
\sum_{i,j=1}^n a_{ij}(x)\xi_i\xi_j\in [\lambda,\Lambda],
\end{equation}
for some~$\Lambda\ge\lambda>0$.

Let~$u\in C^2(\Omega)$ be such that
\begin{equation}\label{CIHNSINSikpomf98S} \begin{dcases}
\sum_{i,j=1}^na_{ij}(x)\partial_{ij}u(x)+\sum_{i=1}^n b_i(x)\partial_i u(x) +c(x)u(x)\ge0 &{\mbox{ for all }}x\in \Omega,\\
u(x)\le0&{\mbox{ for all }}x\in\Omega.
\end{dcases}\end{equation}

Then, $u<0$ in~$\Omega$, unless it is constantly equal to zero.
\end{theorem}

The reason for which we consider quite general operators in~\eqref{CIHNSINSikpomf98S}
(and not just, say, the Laplacian plus possibly a zero order term)
is that these operators will pop up naturally when we will deal with soap bubbles in Section~\ref{MPPALEB}:
compare, in particular, with equation~\eqref{orollarHOPF12ELLIPTIC-PERFVDJ-LEM-32-CORo};
fortunately this additional generality does not really provide major complications
in the proofs, whose structure would remain pretty much the same even in the case of
a Laplacian plus a zero order term.

As we will see more specifically in Section~\ref{KPJMDarJMSiaCLiMMSciAmoenbvD8Sijf-23},
the assumption in~\eqref{ELLIPTIC-PERFVDJ} can be seen as an ``ellipticity requirement''
on the leading coefficients of the equation (recall also footnote~\ref{CLASSIFICATIONFOOTN}
on page~\pageref{CLASSIFICATIONFOOTN}): roughly speaking, \label{CLASSIFICATIONFOOTN2}
this type of conditions are usually taken to ensure some kinship to the Laplace operator 
and, for instance, guarantee\footnote{To see why condition~\eqref{ELLIPTIC-PERFVDJ}
is advantageous when dealing with a Maximum Principle,
we point out that such a requirement ``behaves nicely'' with respect
to the maxima of a function. For instance,
if~$u$ has a local maximum at some point~$\overline{x}$, then
$$\sum_{i,j=1}^na_{ij}(\overline x)\partial_{ij}u(\overline x)\le0.$$
Since we will repeatedly exploit this observation, let us give a quick proof of it.
By construction,
the Hessian matrix of~$u$ at~$\overline{x}$
is nonpositive and therefore we can denote by~$M=\{M_{ij}\}_{i,j\in\{1,\dots,n\}}$ the square root (in the matrix sense, see e.g.~\cite[page~125]{MR2337395})
of minus the Hessian matrix of~$u$ at~$\overline{x}$. In this way, we have that
$$ -\partial_{ij}u(\overline{x})=\sum_{k=1}^n M_{ik}M_{jk}.$$
We also let~$\zeta_k:=(M_{1k},\dots,M_{nk})$, so that
$$ |\zeta_k|^2=\sum_{\ell=1}^n M_{\ell k}^2.$$
As a result, from the ellipticity assumption in~\eqref{ELLIPTIC-PERFVDJ} we arrive at
\begin{eqnarray*}
&&-\sum_{i,j=1}^na_{ij}(\overline x)\partial_{ij}u(\overline x)\\
&=&\sum_{i,j,k=1}^na_{ij}(\overline x) M_{ik}M_{jk}
\\&\ge&\lambda\sum_{k=1}^n |\zeta_k|^2\\
&=&\lambda\sum_{k,\ell=1}^n |M_{\ell k}|^2\\
&\ge&0.
\end{eqnarray*}}
the validity of a suitable Maximum Principle.

The proof of Theorem~\ref{MAXPLECON45678CG-PM} makes use of a generalization of the Hopf Lemma
presented in Lemma~\ref{JJS:PA}, which goes as follows:

\begin{figure}
  \centering
  \includegraphics[width=.45\linewidth]{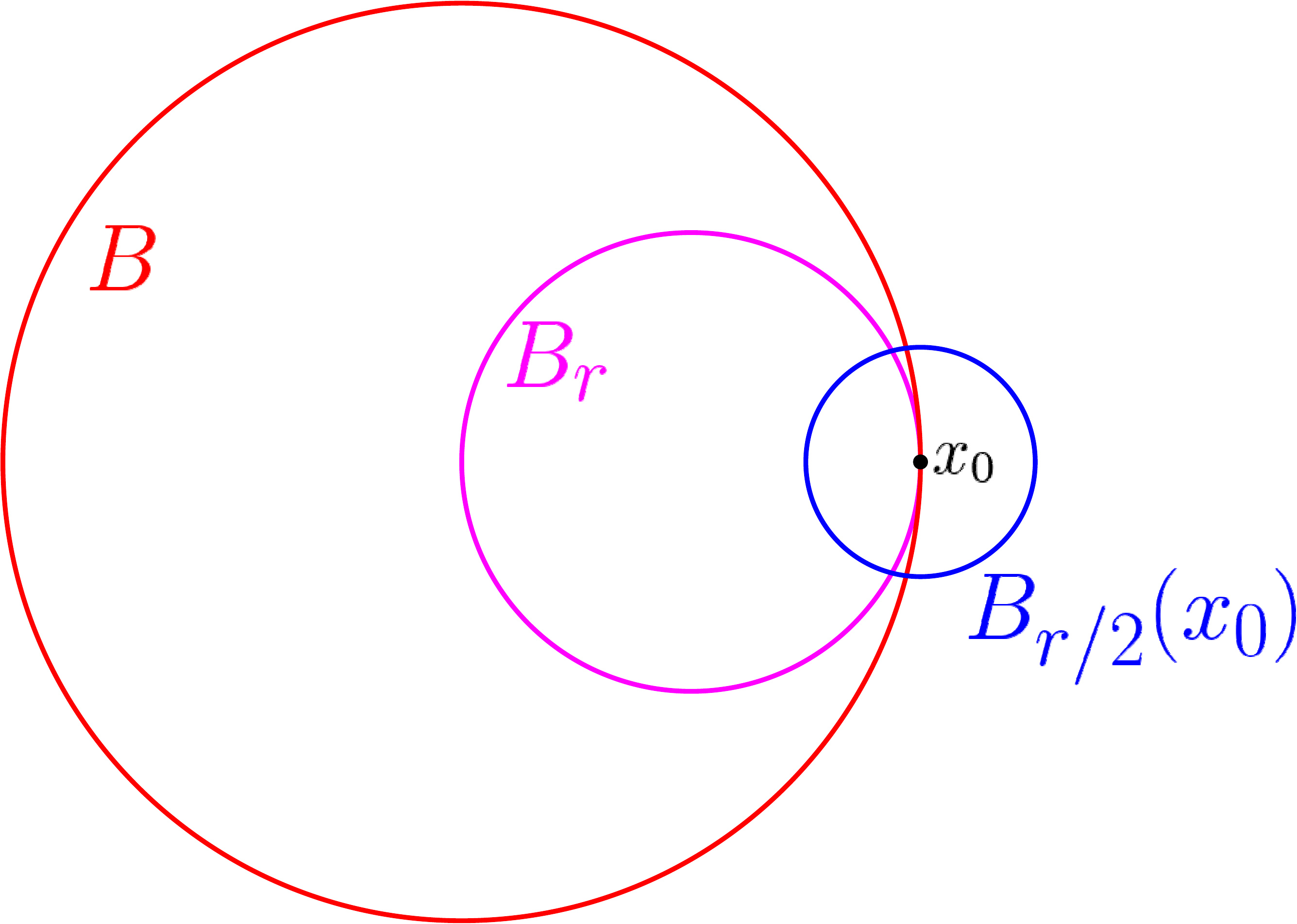}
 \caption{\sl The geometry involved in the proof of Lemma~\ref{MS-LEHOPF2-M}.}\label{LwAS1IewN24XGIJ7soloD354hItangeFIrfd72iqwhrjfngpoewjnikdhgebidkkwdbbXVSdng}
\end{figure}

\begin{lemma}\label{MS-LEHOPF2-M}
Let~$B$ be an open ball in~$\R^n$ with outward normal~$\nu$.

Let~$a_{ij}$, $b_i$, $c\in C(\overline{B})$.
Assume the ellipticity condition in~\eqref{ELLIPTIC-PERFVDJ}
and suppose also that~$c\le0$ in~$B$.

Let~$u\in C^2(B)\cap C^1(\overline{B})$ be a solution of
$$ 
\sum_{i,j=1}^na_{ij}(x)\partial_{ij}u(x)+\sum_{i=1}^n b_i(x)\partial_i u(x) +c(x)u(x)\ge0\qquad{\mbox{for all }}\,x\in B.$$
Assume that~$x_0\in\partial B$ is such that~$u(x_0)\ge0$ and~$u(x)< u(x_0)$ for every~$x\in B$.\index{Hopf Lemma}

Then,~$ \partial_\nu u(x_0)>0$.\end{lemma}

\begin{proof}
Let~$\widetilde B$ be a ball contained in~$B$, with radius half of the radius of~$B$ and such that~$x_0\in\partial\widetilde B$.
Up to a translation, we suppose that~$\widetilde B$ is centered at the origin, namely~$\widetilde B=B_r$ for some~$r>0$.
Let~$S:=B_r\cap B_{r/2}(x_0)$ and notice that~$B_r\cap (\partial B_{r/2}(x_0))\Subset B$, see Figure~\ref{LwAS1IewN24XGIJ7soloD354hItangeFIrfd72iqwhrjfngpoewjnikdhgebidkkwdbbXVSdng}, and therefore
$$ \sup_{B_r\cap (\partial B_{r/2}(x_0))} u<u(x_0).$$
Hence, there exists~$\delta_0>0$ such that
\begin{equation}\label{9kmSDHC578055FliSBA42M212MA}
\sup_{B_r\cap (\partial B_{r/2}(x_0))} u+\delta_0\le u(x_0).\end{equation}
Let also
\begin{equation}\label{LwAS1IewN24XGIJ7soloD354hItangeFIrfd72iqwhrjfngpoewjnikdhgebidkkwdbbXVSdng9028032utDFVBmeKTRA920irkjtd}
h(x):=e^{-\alpha|x|^2}-e^{-\alpha r^2},\end{equation} with~$\alpha>1$ to be conveniently chosen.
We point out that, for all~$x\in B_r$, $h(x)\ge0$ and, using the sign of~$c$
and the ellipticity condition in~\eqref{ELLIPTIC-PERFVDJ}, we find that, for all~$x\in B_r\setminus B_{r/4}$,
\begin{equation*}
\begin{split}&
\sum_{i,j=1}^na_{ij}(x)\partial_{ij}h(x)+\sum_{i=1}^n b_i(x)\partial_i h(x) +c(x)h(x)\\&\qquad=
2\alpha e^{-\alpha|x|^2}\sum_{i,j=1}^na_{ij}(x) \big(2\alpha x_i x_j - \delta_{ij} \big)-2\alpha e^{-\alpha|x|^2}\sum_{i=1}^n b_i(x)x_i +c(x)\big(e^{-\alpha|x|^2}-e^{-\alpha r^2}\big)
\\&\qquad\ge2\alpha e^{-\alpha|x|^2}\sum_{i,j=1}^na_{ij}(x) \big(2\alpha x_i x_j- \delta_{ij} \big)-2\alpha e^{-\alpha|x|^2}\sum_{i=1}^n b_i(x)x_i +c(x) e^{-\alpha|x|^2}\\&\qquad=
2\alpha e^{-\alpha|x|^2}\left[\sum_{i,j=1}^na_{ij}(x) \big(2\alpha x_i x_j - \delta_{ij} \big)-\sum_{i=1}^n b_i(x)x_i +\frac{c(x)}{2\alpha}\right]
\\&\qquad\ge
2\alpha e^{-\alpha|x|^2}\left[ 2\alpha \lambda|x|^2- n^2\sup_{i,j\in\{1,\dots,n\}}\|a_{ij}\|_{L^\infty(B)}
-nr\sup_{i\in\{1,\dots,n\}}\|b_i\|_{L^\infty(B)}
-\frac{\|c\|_{L^\infty(B)}}{2}\right]\\&\qquad\ge
2\alpha e^{-\alpha|x|^2}\left[ 2\alpha \lambda|x|^2- C\,(1+r)\right]\\&\qquad\ge
2\alpha e^{-\alpha|x|^2}\left[ \frac{\alpha \lambda r^2}8- C\,(1+r)\right],\end{split}\end{equation*}
where~$C>0$ depends only on~$n$ and on the structural bounds on the coefficients of the equation.

As a result, for all~$x\in B_r\setminus B_{r/4}$,
\begin{equation}\label{SLc870xvexe}
\sum_{i,j=1}^na_{ij}(x)\partial_{ij}h(x)+\sum_{i=1}^n b_i(x)\partial_i h(x) +c(x)h(x)\ge 2\alpha e^{-\alpha|x|^2},
\end{equation}
as long as~$\alpha$ is sufficiently large (and this determines~$\alpha$ once and for all).

Now, let~$\e\in(0,\delta_0)$, to be chosen conveniently small in what follows,
and~$v_\e:=u+\e h$. We point out that
\begin{equation}\label{8ihkjnNSX782fgAIK-3rk3}
v_\e(x_0)= u(x_0)\ge0.
\end{equation}
We claim that
\begin{equation}\label{YTgMAHNanOMCCRMIAHYKLEMTNAMIUOJHN}
\sup_S v_\e=\sup_{\partial S} v_\e.
\end{equation}
Indeed, suppose not. Then~$v_\e$ possesses an interior maximum~$x_\e$ in~$S$.
For this reason, $\nabla v_\e(x_\e)=0$ and~$D^2v_\e(x_\e)\le0$, which, using again 
the ellipticity condition in~\eqref{ELLIPTIC-PERFVDJ},
yields that
$$\sum_{i,j=1}^na_{ij}(x_\e)\partial_{ij}v_\e(x_\e)+\sum_{i=1}^n b_i(x_\e)\partial_i v_\e(x_\e)+c(x_\e)v_\e(x_\e)\le c(x_\e)v_\e(x_\e).
$$
Moreover, $v_\e(x_\e)\ge v_\e(x_0)\ge0$, thanks to~\eqref{8ihkjnNSX782fgAIK-3rk3}.
Using these observations, the sign of~$c$ and~\eqref{SLc870xvexe}, we infer that
\begin{eqnarray*}
0 &\ge&c(x_\e)v_\e(x_\e)
\\&\ge&\sum_{i,j=1}^na_{ij}(x_\e)\partial_{ij}v_\e(x_\e)+\sum_{i=1}^n b_i(x_\e)\partial_i v_\e(x_\e)+c(x_\e)v_\e(x_\e)
\\&\ge&\e \left[\sum_{i,j=1}^na_{ij}(x_\e)\partial_{ij}h(x_\e)+\sum_{i=1}^n b_i(x_\e)\partial_i h(x_\e)+c(x_\e)h(x_\e)\right]\\
&\ge&2\alpha\e e^{-\alpha|x_\e|^2}.
\end{eqnarray*}
But the latter term is strictly positive and this contradiction establishes~\eqref{YTgMAHNanOMCCRMIAHYKLEMTNAMIUOJHN}.

We also observe that if~$x\in\partial B_r$ then~$h(x)=0$ and therefore~$v_\e(x)=u(x)\le u(x_0)\le v_\e(x_0)$.
This gives that
\begin{equation}\label{YTgMAHNanOMCCRMIAHYKLEMTNAMIUOJHN-2}
\sup_{(\partial B_r)\cap B_{r/2}(x_0)} v_\e= v_\e(x_0).
\end{equation}
Furthermore, by~\eqref{9kmSDHC578055FliSBA42M212MA},
if~$x\in B_r\cap (\partial B_{r/2}(x_0))$ then
$$ v_\e(x)\le u(x)+\e \|h\|_{L^\infty(B_r)}\le u(x)+\delta_0\le u(x_0)\le v_\e(x_0),$$
provided that~$\e$ is chosen suitably small (and this determines~$\e$ once and for all).
Accordingly,
\begin{equation*}
\sup_{B_r\cap (\partial B_{r/2}(x_0))} v_\e\le v_\e(x_0).
\end{equation*}
With this, \eqref{YTgMAHNanOMCCRMIAHYKLEMTNAMIUOJHN} and~\eqref{YTgMAHNanOMCCRMIAHYKLEMTNAMIUOJHN-2} we see that
\begin{eqnarray*}
\sup_S v_\e=\max\left\{
\sup_{(\partial B_r)\cap B_{r/2}(x_0)} v_\e,\;
\sup_{B_r\cap (\partial B_{r/2}(x_0))} v_\e
\right\}=v_\e(x_0).
\end{eqnarray*}
As a consequence,
$$ 0\le\partial_\nu v_\e(x_0)=\partial_\nu u(x_0)+\e \partial_\nu h(x_0)
=\partial_\nu u(x_0)-2\alpha\e |x_0|e^{-\alpha|x_0|^2}<\partial_\nu u(x_0),$$
as desired.
\end{proof}

\begin{figure}
  \centering
  \includegraphics[width=.45\linewidth]{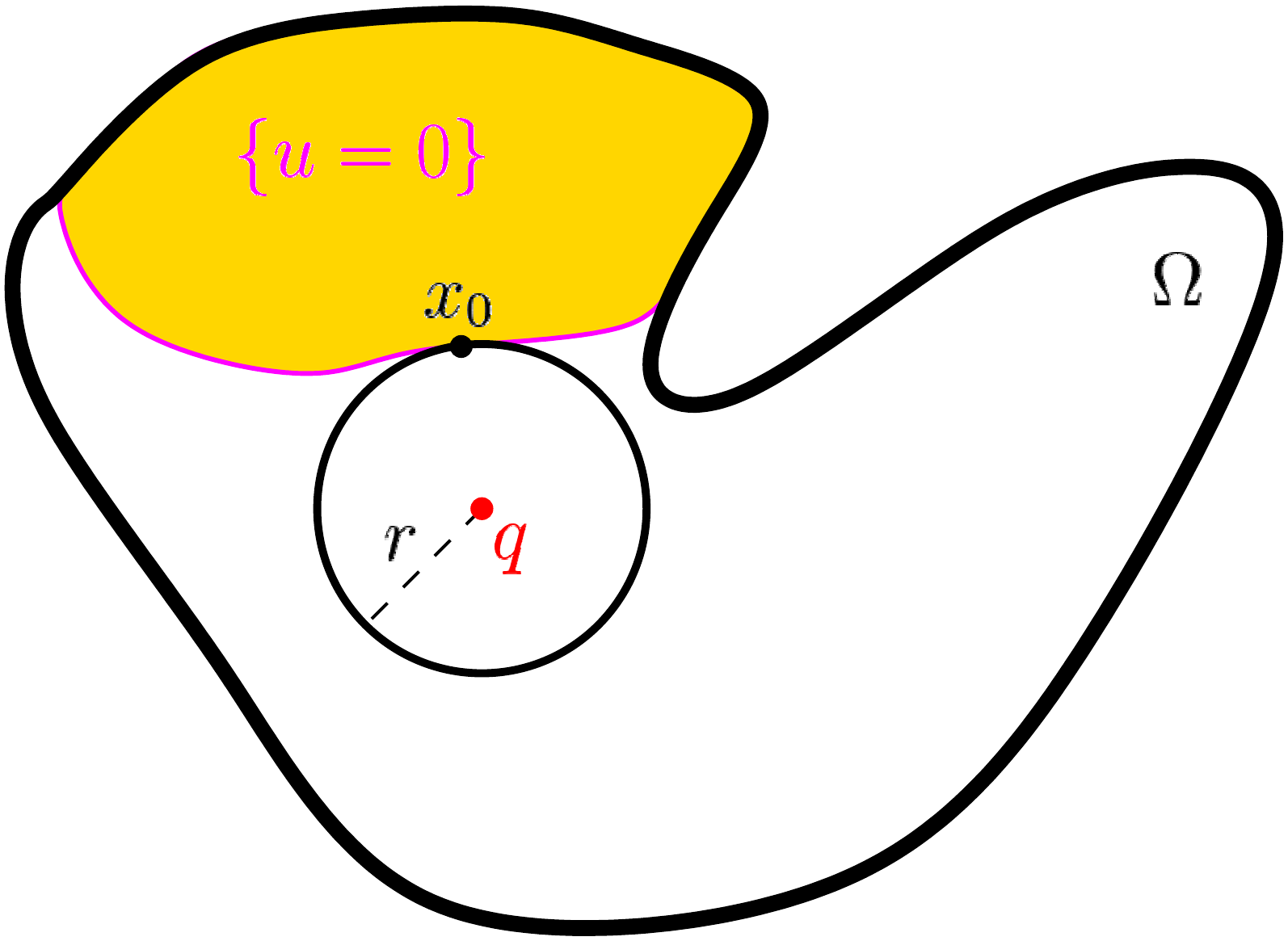}
 \caption{\sl The geometry involved in the proof of Theorem~\ref{MAXPLECON45678CG-PM}.}\label{Lw5678u2teighfie80ikmf0KShisfvnodo-9loD354hItangeFIrfd72iqwhrjfngpoewjnikdhgebidkkwdbbXVSdng}
\end{figure}

A consequence of Lemma~\ref{MS-LEHOPF2-M} is that, if the solution~$u$ has a sign, then the sign
assumption on~$c$ can be dropped,
according to
the following statement:

\begin{corollary}\label{MS-LEHOPF2-M-NEWCORO}
Let~$B$ be an open ball in~$\R^n$ with outward normal~$\nu$.

Let~$a_{ij}$, $b_i$, $c\in C(\overline{B})$.
Assume the ellipticity condition in~\eqref{ELLIPTIC-PERFVDJ}.

Let~$u\in C^2(B)\cap C^1(\overline{B})$ be a solution of
$$ 
\sum_{i,j=1}^na_{ij}(x)\partial_{ij}u(x)+\sum_{i=1}^n b_i(x)\partial_i u(x) +c(x)u(x)\ge0\qquad{\mbox{for all }}\,x\in B.$$
Assume that~$x_0\in\partial B$ is such that~$u(x_0)=0$ and~$u(x)< 0$ for every~$x\in B$.

Then,~$ \partial_\nu u(x_0)>0$.
\end{corollary}

\begin{proof}
We notice that, for every~$x\in B$,
\begin{eqnarray*}
0&\le& \sum_{i,j=1}^na_{ij}(x)\partial_{ij}u(x)+\sum_{i=1}^n b_i(x)\partial_i u(x) +c(x)u(x)\\&=&
\sum_{i,j=1}^na_{ij}(x)\partial_{ij}u(x)+\sum_{i=1}^n b_i(x)\partial_i u(x) +c^+(x)u(x)-c^-(x)u(x)\\&\le&
\sum_{i,j=1}^na_{ij}(x)\partial_{ij}u(x)+\sum_{i=1}^n b_i(x)\partial_i u(x) -c^-(x)u(x).
\end{eqnarray*}
Accordingly, the assumptions of Lemma~\ref{MS-LEHOPF2-M} are satisfied taking~$c:=-c^-\le0$,
and therefore~$ \partial_\nu u(x_0)>0$,
as desired.
\end{proof}

\begin{proof}[Proof of Theorem~\ref{MAXPLECON45678CG-PM}] Let~${\mathcal{U}}:=\{x\in\Omega$ s.t. $u(x)=0\}$ and note that~$
{\mathcal{U}}$ is a closed set in~$\Omega$.
If~${\mathcal{U}}=\varnothing$ then~$u<0$ in~$\Omega$ and we are done,
hence we can suppose that~${\mathcal{U}}\ne\varnothing$.
Our goal is to show that~${\mathcal{U}}=\Omega$, since in this way we would
get that~$u$ vanishes identically, as claimed in Theorem~\ref{MAXPLECON45678CG-PM}.
For this, we argue by contradiction and assume that~${\mathcal{U}}$
is strictly contained in~$\Omega$. In particular, there exists~$p\in{\mathcal{U}}$
such that for all~$j\in\N$ we have that~$B_{1/j}(p)\setminus{\mathcal{U}}\ne\varnothing$
(otherwise~${\mathcal{U}}$ would be also open, and then coincide with~$\Omega$ by connectedness).

As a consequence, we can find a point~$q\in \Omega\setminus{\mathcal{U}}$ with~${\rm dist}(q,\partial\Omega)>|q-p|$.
That is,
$$ {\rm dist}(q,\partial\Omega)>|q-p|
\ge{\rm dist}(q,{\mathcal{U}})=:r,$$
yielding that
\begin{equation*}
{\mbox{$B_r(q)$ is contained in~$\Omega\setminus{\mathcal{U}}$ and
there exists~$x_0\in(\partial B_r(q))\cap{\mathcal{U}}$.}}\end{equation*}
See Figure~\ref{Lw5678u2teighfie80ikmf0KShisfvnodo-9loD354hItangeFIrfd72iqwhrjfngpoewjnikdhgebidkkwdbbXVSdng}
for a sketch of this configuration.

Now we let
\begin{equation}\label{CPIUOMENODEFI}
c^+:=\max\{c,0\}\qquad{\mbox{ and }}\qquad c^-:=\max\{-c,0\}.\end{equation} In this way, $c^+\ge0$, $c^-\ge0$
and~$c=c^+-c^-$. Therefore, using the sign of~$u$, for every~$x\in\Omega$
\begin{eqnarray*}
&&\sum_{i,j=1}^na_{ij}(x)\partial_{ij}u(x)+\sum_{i=1}^n b_i(x)\partial_i u(x)-c^-(x)u(x)
\\&=&
\sum_{i,j=1}^na_{ij}(x)\partial_{ij}u(x)+\sum_{i=1}^n b_i(x)\partial_i u(x)+c(x)u(x)
-c^+(x)u(x)\\
&\ge&-c^+(x)u(x)\\&\ge&0.
\end{eqnarray*}
Given the sign of~$c^-$, we are therefore in position of using Lemma~\ref{MS-LEHOPF2-M}
in this context. With the aid of Lemma~\ref{MS-LEHOPF2-M} (used here with~$c$ replaced by~$-c^-$, and~$B$ replaced by~$B_r(q)$) we thus deduce that
either~$u$ is constant (and thus constantly equal to zero, since~${\mathcal{U}}$ is nonvoid),
or~$\nabla u(x_0)\ne0$. But the latter condition cannot hold (since~$x_0$ is an interior maximum,
hence~$\nabla u(x_0)=0$), whence the proof of the desired result is complete. \end{proof}

Putting together Theorem~\ref{MAXPLECON45678CG-PM} and Lemma~\ref{MS-LEHOPF2-M}
we have:

\begin{corollary}\label{HOPF12ELLIPTIC-PERFVDJ-LEM-32-CORo}
Let~$B$ be an open ball in~$\R^n$ with outward normal~$\nu$.
Let~$a_{ij}$, $b_i \in C(\overline{B})$ and assume the ellipticity condition in~\eqref{ELLIPTIC-PERFVDJ}.

Let~$u\in C^2(B)\cap C^1(\overline{B})$.
Suppose that
$$ \begin{dcases}
\sum_{i,j=1}^na_{ij}(x)\partial_{ij}u(x)+\sum_{i=1}^n b_i(x)\partial_i u(x) \ge0 &{\mbox{ for all }}x\in B,\\
u(x)\le0&{\mbox{ for all }}x\in B.
\end{dcases}$$
Let~$x_0\in\partial B$ be such that~$u(x_0)=0$.\index{Hopf Lemma}
Then, either~$ \partial_\nu u(x_0)>0$ or~$u$ vanishes identically.
\end{corollary}

\begin{proof} Suppose that~$u$ does not vanish identically in~$B$. Then, by
Theorem~\ref{MAXPLECON45678CG-PM}, we have that~$u<0$ in~$B$.
Consequently, we can use Lemma~\ref{MS-LEHOPF2-M}
and infer that~$\partial_\nu u(x_0)>0$.
\end{proof}

A related result is the so-called Maximum Principle for narrow domains:\index{Maximum Principle for narrow domains}
this result focuses on domains which are constrained in a sufficiently thin strip and it goes as follows.

\begin{theorem}\label{M:PLE:NARROW}
Let~$\Omega$ be an open and bounded subset of~$\R^n$ with boundary
of class~$C^1$. Let~$c\in C(\overline\Omega)$
and assume that
\begin{equation}\label{MCONSNARR}\begin{split}&
\Omega\Subset\big\{
x=(x_1,\dots,x_n)\in\R^n {\mbox{ s.t. }}x_1\in(0,d)
\big\}\\&\qquad{\mbox{for some~$d>0$ with $
{\|c^+\|_{L^{\infty}(\Omega)}^{1/2}}\,d\in[0,
\pi]$.}}\end{split}\end{equation}
Let~$u\in C^2(\Omega)\cap C^1(\overline\Omega)$ be a solution of
$$\begin{dcases}
\Delta u(x)+c(x)u(x)\ge0 &{\mbox{ for all }}x\in\Omega,\\
u(x)\le0&{\mbox{ for all }}x\in\partial\Omega.
\end{dcases}$$
Then, $u\le0$ in~$\Omega$.

Also, if~$\Omega$ is connected,
then~$u<0$ in~$\Omega$, unless it is constantly equal to zero.
\end{theorem}

{F}rom the geometric point of view, condition~\eqref{MCONSNARR} requires the domain~$\Omega$
to lie within a slab of sufficiently small width, which justifies the name
of Maximum Principle for narrow domains
(the position of the slab is irrelevant,
up to a translation and a rotation,
recall~\eqref{TRAL} and Corollary~\ref{KSMD:ROTAGSZOKA}).

\begin{proof}[Proof of Theorem~\ref{M:PLE:NARROW}]
First of all, we show that
\begin{equation}\label{14KPLC2rluiaLP00}
{\mbox{$u(x)\le0$ for every~$x\in\Omega$.}}
\end{equation}
For this, up to reducing our analysis to each connected component of~$\Omega$,
we can assume that
\begin{equation}\label{14KPLC2rluiaLP}
{\mbox{$\Omega$ is connected.}}
\end{equation}
Let
$$ w(x):=\sin\frac{\pi x_1}{d}$$
and observe that~$w>0$ in~$\overline\Omega$. Thus, if
$$ \Lambda:=
\frac{\|u\|_{L^\infty(\Omega)}}{ \displaystyle\min_{\overline\Omega}w }+1
$$
we see that, for all~$x\in\overline\Omega$,
$$ \Lambda w(x)-u(x)\ge\Lambda\min_{\overline\Omega}w-\|u\|_{L^\infty(\Omega)}>0.$$
We thereby define
\begin{equation}\label{f245ggor2al21lgq2a} \Lambda_0:=\inf\big\{\Lambda {\mbox{ s.t. $\Lambda w(x)-u(x)>0$ for every $x\in\Omega$}}
\big\}.\end{equation}

Now, we claim that
\begin{equation}\label{KMSDFGHJSIJD9iry93uwgfoiqwt7893uowig9fuwegfiGBSIMfjgmVkfRNKu4gl3a2teft6tejgs2bb}
\Lambda_0\le0.
\end{equation}
Indeed, suppose not, namely~$\Lambda_0>0$. Let~$v:=u-\Lambda_0 w$ and observe that, for all~$x\in\Omega$,
$$ \Delta v(x)+c(x)v(x)\ge
-\Lambda_0\,\Big(
\Delta w(x)+c(x)w(x)
\Big)=-\Lambda_0\,\left(
-\frac{\pi^2}{d^2}+c(x)\right)\,w(x).
$$
{F}rom this and the fact (recall~\eqref{MCONSNARR}) that
$$ c(x)\le c^+(x)\le \|c^+\|_{L^\infty(\Omega)}\le\frac{\pi^2}{d^2},$$
we conclude that~$\Delta v(x)+c(x)v(x)\ge0$ for all~$x\in\Omega$.
Also~$v\le0$ in~$\Omega$, due to~\eqref{f245ggor2al21lgq2a},
and therefore, in the setting of~\eqref{14KPLC2rluiaLP}, we are in position of using Theorem~\ref{MAXPLECON45678CG-PM} and deduce that
\begin{equation}\label{1294825akm19lTe}
{\mbox{either~$v<0$ in~$\Omega$ or~$v$ vanishes identically.}}\end{equation}
We also remark that, on~$\partial\Omega$, it holds that~$v<u\le0$. Combining this information with~\eqref{1294825akm19lTe}
we infer that~$v<0$ in~$\overline\Omega$. As a result,
$$ \e_0:=-\max_{\overline\Omega} v>0$$
and therefore, for all~$x\in\Omega$ and~$\Lambda\in(\Lambda_0-\e_0,\Lambda_0]$,
$$ \Lambda w(x)-u(x)=(\Lambda-\Lambda_0)w(x)-v(x)\ge\e_0
-|\Lambda-\Lambda_0|=\e_0-\Lambda_0+\Lambda>0.$$
This is a contradiction with the infimum property
in~\eqref{f245ggor2al21lgq2a} and so it completes the proof of~\eqref{KMSDFGHJSIJD9iry93uwgfoiqwt7893uowig9fuwegfiGBSIMfjgmVkfRNKu4gl3a2teft6tejgs2bb}. 

{F}rom~\eqref{KMSDFGHJSIJD9iry93uwgfoiqwt7893uowig9fuwegfiGBSIMfjgmVkfRNKu4gl3a2teft6tejgs2bb} 
we deduce~\eqref{14KPLC2rluiaLP00}.

This and Theorem~\ref{MAXPLECON45678CG-PM} give the desired result.
\end{proof}

A variant of the Maximum Principle for narrow domains
in Theorem~\ref{M:PLE:NARROW} is the so-called Maximum Principle for small volume domains\index{Maximum Principle for small volume domains},
which goes as follows.

\begin{theorem}\label{M:PLE:SMALLV}
There exists~$c_0\in(0,1)$, depending only on~$n$, such that the following statement holds true.

Let~$\Omega$ be an open and bounded subset of~$\R^n$ with boundary
of class~$C^1$. Let~$c\in C(\overline\Omega)$
and assume that
\begin{equation}\label{22qwrf2KMSDFGHJSIJD9iry93uwgfoiqwt7893uowig9fuwegfiGBSIMfjgmVkfRNKu4gl3a2teft6tejgs}
\|c^+\|_{L^{\infty}(\Omega)}\,|\Omega|^{\frac2n}<c_0.\end{equation}
Let~$u\in C^2(\Omega)\cap C^1(\overline\Omega)$ be a solution of
$$\begin{dcases}
\Delta u(x)+c(x)u(x)\ge0 &{\mbox{ for all }}x\in\Omega,\\
u(x)\le0&{\mbox{ for all }}x\in\partial\Omega.
\end{dcases}$$
Then, $u\le0$ in~$\Omega$.

Also, if~$\Omega$ is connected,
then~$u<0$ in~$\Omega$, unless it is constantly equal to zero.
\end{theorem}

This result is actually a particular case of a more general statement.
The core of the idea, dating back to Guido Stampacchia~\cite{MR0251373}, 
is that suitable Maximum Principles continue to remain valid
for ``small'' linear perturbations of the Laplacian. Note that the main structural condition
of Theorem~\ref{M:PLE:SMALLV} (namely
hypothesis~\eqref{22qwrf2KMSDFGHJSIJD9iry93uwgfoiqwt7893uowig9fuwegfiGBSIMfjgmVkfRNKu4gl3a2teft6tejgs})
is satisfied provided that the Lebesgue measure of the domain~$\Omega$
is sufficiently small (possibly in dependence of the size of~$c^+$) and this justifies
the name of Maximum Principle for small volume domains.
The reader can appreciate similarities and differences
between the narrow domain condition~\eqref{MCONSNARR}
and the small volume condition~\eqref{22qwrf2KMSDFGHJSIJD9iry93uwgfoiqwt7893uowig9fuwegfiGBSIMfjgmVkfRNKu4gl3a2teft6tejgs}.

\begin{proof} [Proof of Theorem~\ref{M:PLE:SMALLV}]
We will employ\footnote{When~$n\ge3$ one can use instead
the Sobolev-Gagliardo-Nirenberg's Inequality (see e.g.~\cite[Theorem~9.9]{MR2759829}),
as done in~\cite[Lemma~1]{MR1662746}. Using Nash Inequality here we can treat
all the dimensions~$n$ at the same time.} a classical inequality\index{Nash Inequality}
introduced by John Nash~\cite{MR100158},
stating that for all~$v\in L^1(\R^n)\cap W^{1,2}(\R^n)$ we have that
\begin{equation}\label{NASH:IN}
\|v\|_{L^{2}(\R^{n})}^{\frac{n+2}n}\leq C\,\|v\|_{L^{1}(\R^{n})}^{\frac2n}\,\|\nabla v\|_{L^{2}(\R^{n})},
\end{equation}
for some~$C>0$ depending only on~$n$.
For completeness, we recall the elegant proof of this inequality
(for another elegant proof based on rearrangements see~\cite[Theorem~8.13]{MR1817225}). We can assume that~$v$ does not
vanish identically (otherwise~\eqref{NASH:IN} is obvious) and we
take the Fourier Transform~$\widehat v(\xi)$
of~$v(x)$ and we notice that~$\partial_j v(x)$ coincides with the inverse Fourier Transform
of~$2\pi i\xi_j \widehat v(\xi)$, for all~$j\in\{1,\dots,n\}$.
Consequently, Plancherel Theorem gives that, for all~$\rho>0$,
\begin{eqnarray*}&& \int_{\R^n\setminus B_\rho}|\widehat v(\xi)|^2\,d\xi
\le\frac{1}{\rho^2}
\int_{\R^n\setminus B_\rho}|\xi|^2\,|\widehat v(\xi)|^2\,d\xi\le\frac{1}{\rho^2}
\int_{\R^n}|\xi|^2\,|\widehat v(\xi)|^2\,d\xi\\&&\qquad\qquad
=\frac{1}{4\pi^2\rho^2}
\int_{\R^n}|\nabla v(x)|^2\,dx\le
\frac{\|\nabla v\|_{L^{2}(\R^{n})}^2}{4\pi^2\rho^2}.
\end{eqnarray*}
Additionally,
$$ |\widehat v(\xi)|=\left|\int_{\R^n} v(x)\,e^{2\pi i x\cdot\xi}\,dx
\right|\le\int_{\R^n} |v(x)|\,dx=\|v\|_{L^1(\R^n)}$$
and therefore
$$ \int_{B_\rho}|\widehat v(\xi)|^2\,d\xi
\le |B_1|\,\rho^n\,\|v\|_{L^1(\R^n)}^2.$$
{F}rom these observations, and using Plancherel Theorem again, it follows that
$$ \| v\|_{L^{2}(\R^{n})}^2=
\int_{\R^n}|\widehat v(\xi)|^2\,d\xi\le
\frac{\|\nabla v\|_{L^{2}(\R^{n})}^2}{4\pi^2\rho^2}+
|B_1|\,\rho^n\,\|v\|_{L^1(\R^n)}^2.$$
With this, the claim in~\eqref{NASH:IN} follows by picking~$\rho:=\left(\frac{\|\nabla v\|_{L^{2}(\R^{n})}}{\| v\|_{L^{1}(\R^{n})}}\right)^{\frac2{n+2}}$.

Now, if~$v\in W^{1,2}_0(\Omega)$, we deduce from~\eqref{NASH:IN}
and the Cauchy-Schwarz Inequality that
$$ \|v\|_{L^{2}(\R^{n})}^{\frac{n+2}n}\leq C\,
\left(|\Omega| \int_\Omega|v(x)|^2\,dx\right)^{\frac1n}\,\|\nabla v\|_{L^{2}(\R^{n})}=
C\,|\Omega|^{\frac1n}\,\|v\|_{L^{2}(\R^{n})}^{\frac{2}n}\,\|\nabla v\|_{L^{2}(\R^{n})}
$$
and therefore\footnote{Equation~\eqref{PONASH}
can also be seen as a version of the Poincar\'e Inequality with an explicit
dependence of the constant on the volume of the domain. Compare e.g. with~\cite[Theorem~12.17]{MR2527916}.}
\begin{equation}\label{PONASH} \|v\|_{L^{2}(\R^{n})}\leq C\,
|\Omega|^{\frac1n}\,\|\nabla v\|_{L^{2}(\R^{n})}.
\end{equation}

Now, we show that
\begin{equation}\label{KMSDFGHJSIJD9iry93uwgfoiqwt7893uowig9fuwegfiGBSIMfjgmVkfRNKu4gl3a2teft6tejgs2}
u(x)\le0\qquad{\mbox{for every }}\,x\in\Omega.
\end{equation}
To this end, we note that
\begin{equation}\label{MDABUHSOMENAG6ALSPBPIN8aay}
{\mbox{$u^+=0$ along~$\partial\Omega$}},\end{equation}
thus~$u^+\in W^{1,2}(\R^n)$, once we extend it by zero outside~$\Omega$,
see e.g.~\cite[Proposition 9.18]{MR2759829}, and
\begin{equation}\label{IHK-pjlfeiwhfgewiukboighboiurae6yhHOSIM}
\begin{split} 0\,&\le\,\int_\Omega\big(\Delta u(x)+c(x)u(x)\big)\,u^+(x)\,dx\\&=\,
\int_\Omega \Delta u(x)\,u^+(x)\,dx
+\int_\Omega c(x)\,(u^+(x))^2\,dx
\\&\le\,
\int_\Omega \Big( \div\big( u^+(x)\,\nabla u(x)\big)-\nabla u^+(x)\cdot\nabla u(x)
\Big)\,dx
+\int_\Omega c^+(x)\,(u^+(x))^2\,dx\\&=\,
\int_{\partial\Omega} u^+(x)\,\partial_\nu u(x)\,d{\mathcal{H}}^{n-1}_x
-\int_\Omega |\nabla u^+(x)|^2\,dx
+\int_\Omega c^+(x)\,(u^+(x))^2\,dx\\&=\,
-\int_\Omega |\nabla u^+(x)|^2\,dx
+\int_\Omega c^+(x)\,(u^+(x))^2\,dx.
\end{split}\end{equation}
Besides, by~\eqref{PONASH},
$$ \int_\Omega c^+(x)\,(u^+(x))^2\,dx\le
\|c^+\|_{L^\infty(\Omega)}\,\|u^+\|_{L^2(\Omega)}^2\le
C\, \|c^+\|_{L^\infty(\Omega)}\,|\Omega|^{\frac2n}\,\|\nabla u^+\|_{L^{2}(\R^{n})}^2.$$
{F}rom this and~\eqref{IHK-pjlfeiwhfgewiukboighboiurae6yhHOSIM} we arrive at
$$ \|\nabla u^+\|_{L^{2}(\R^{n})}^2=\int_\Omega |\nabla u^+(x)|^2\,dx\le
\int_\Omega c^+(x)\,(u^+(x))^2\,dx\le C\, \|c^+\|_{L^\infty(\Omega)}\,|\Omega|^{\frac2n}\,\|\nabla u^+\|_{L^{2}(\R^{n})}^2.$$
This and~\eqref{22qwrf2KMSDFGHJSIJD9iry93uwgfoiqwt7893uowig9fuwegfiGBSIMfjgmVkfRNKu4gl3a2teft6tejgs}
give that~$\|\nabla u^+\|_{L^{2}(\Omega)}=0$ and therefore~$u^+$ is constant in each connected
component of~$\Omega$.
{F}rom this and~\eqref{MDABUHSOMENAG6ALSPBPIN8aay} we obtain~\eqref{KMSDFGHJSIJD9iry93uwgfoiqwt7893uowig9fuwegfiGBSIMfjgmVkfRNKu4gl3a2teft6tejgs2}. 

The desired result then follows from~\eqref{KMSDFGHJSIJD9iry93uwgfoiqwt7893uowig9fuwegfiGBSIMfjgmVkfRNKu4gl3a2teft6tejgs2} and
Theorem~\ref{MAXPLECON45678CG-PM}.
\end{proof}

Several of the results presented here actually hold true in further generality,
see e.g.~\cite{MR333220, MR544879, MR1662746, MR2356201} and the references therein.
See also~\cite[Theorem~2.32]{MR2777537}
for a formulation of the Maximum Principle for small volume domains due to Sathamangalam Ranga Iyengar Srinivasa Varadhan.

\section{Cauchy's Estimates}

Now we deduce from the regularity result
in Theorem~\ref{KMS:HAN0-0} and the
Mean Value Formula in Theorem~\ref{KAHAR}(iii)
a useful set of explicit bounds on the derivative of harmonic functions,
often referred to with the name\footnote{For related Cauchy's Estimates
in the complex variable setting, see e.g.~\cite[Lemma~10.7]{MR3839273}.} of
Cauchy's Estimates:\index{Cauchy's Estimates}

\begin{theorem}\label{CAUESTIMTH}
Let~$\Omega\subseteq\R^n$ be open and let~$u$ be harmonic in~$\Omega$.
Assume that~$B_r(x_0)\Subset\Omega$. Then, for every~$k\in\N$ and every~$\alpha\in\N^n$
with~$|\alpha|=k$ we have that
$$ \left|
\frac{\partial^\alpha u}{\partial x^\alpha}(x_0)
\right|\le \frac{C_k\,\|u\|_{L^1(B_r(x_0))}}{r^{n+k}}.$$
The constant~$C_k$ above can be taken of the form
\begin{equation}\label{CIKAPP} C_0:=\frac{1}{|B_1|}\qquad{\mbox{and}}\qquad
C_k:=\frac{ (2^{n+1} n k)^k}{|B_1|}\qquad{\mbox{for all }}k\ge1.\end{equation}
\end{theorem}

\begin{proof} We stress
that, by
the regularity result
in Theorem~\ref{KMS:HAN0-0}, we know that~$u$ is differentiable
as many times as we want and all its derivatives are harmonic.
We suppose, up to a translation, that~$x_0=0$ and we argue
by induction over~$k$.

Thus, we use the Mean Value Formula in Theorem~\ref{KAHAR}(iii)
to see that
$$ |u(0)|=\left|\fint_{B_r} u(x)\,dx\right|\le
\fint_{B_r}| u(x)|\,dx=
\frac{\|u\|_{L^1(B_r)}}{|B_1|\,r^n}
$$
and this is the desired estimate for~$k=0$.

When~$k=1$, we define~$v_i:=\partial_{x_i} u$ for every~$i\in\{1,\dots,n\}$ and
we observe that~$v_i
$ is harmonic in~$\Omega$. Hence, we can exploit the
Mean Value Formula in Theorem~\ref{KAHAR}(iii) (applied here to~$v_i$) 
and the Divergence Theorem, to see that
\begin{equation}\label{43vh9hfue75843769}
\begin{split}
&\partial_{x_i}u(0)=v_i(0)=
\fint_{B_{r/2}} v_i(x)\,dx=
\fint_{B_{r/2}} \partial_{x_i} u(x)\,dx\\&\qquad=
\frac{2^n}{|B_1|\,r^n}\int_{B_{r/2}} \div(u(x)e_i)\,dx=
\frac{2^{n+1}}{|B_1|\,r^{n+1}}\int_{\partial B_{r/2}} u(x)x_i\,d{\mathcal{H}}^{n-1}_x.
\end{split}\end{equation}
Now, we observe that
for every~$q\in\partial B_{r/2}$ we have that~$B_{r/2}(q)\subseteq B_r$.
Consequently, using Cauchy's Estimate for~$k=0$,
$$ |u(q)|
\le \frac{C_0\, 2^n\,\|u\|_{L^1(B_{r/2}(q))}}{r^{n}}\le \frac{ 2^n\,\|u\|_{L^1(B_{r})}}{|B_1|\,r^{n}}
.$$
Plugging this information into~\eqref{43vh9hfue75843769}, we thereby obtain that
$$
|\partial_{x_i} u(0)|\le \frac{2^{n+1}n\, \|u\|_{L^1(B_{r})}}{|B_1|\,r^{n+1}},
$$
which proves the desired estimate when~$k=1$.

Now, we argue recursively and we suppose that the desired estimate holds true for the index~$k\in\N$.
We pick~$\alpha=(\alpha_1,\dots,\alpha_n)\in\N^n$ with~$|\alpha|=k+1$
and, up to exchanging the order of the variables, we assume that~$\alpha_1\ne0$.
Hence, we can write~$\alpha=\beta+e_1$, with~$\beta:=
(\alpha_1-1,\alpha_2,\dots,\alpha_n)$. Notice that~$|\beta|=k$.
Let also~$w:=\partial^\beta_x u$
and~$v:=\partial^\alpha_x u=\partial_{x_1}w$.
Since~$v
$ is harmonic in~$\Omega$, applying on~$v$ the
Mean Value Formula in Theorem~\ref{KAHAR}(iii)
and the Divergence Theorem, for all~$\rho\in(0,r]$,
\begin{equation}\label{KMMARGBnoKKSMMLDSshciCLAALMSifUFFPOP}
\begin{split}
&\partial_{x}^\alpha u(0)=v(0)=
\fint_{B_\rho} v(x)\,dx=
\fint_{B_\rho} \partial_{x_1} w(x)\,dx\\&\qquad=
\frac1{|B_1|\,\rho^n}\int_{B_\rho} \div(w(x)e_1)\,dx=
\frac1{|B_1|\,\rho^{n+1}}\int_{\partial B_\rho} w(x)x_1\,d{\mathcal{H}}^{n-1}_x.
\end{split}\end{equation}
Now, we choose~$\rho:=\frac{r}{k+1}$
and~$R:=k\rho$: in this way,
for every~$q\in\partial B_{\rho}$
we have that~$B_{R}(q)\subseteq B_r$. Thus, the inductive assumption
gives that
$$ |w(q)|=\left|
\frac{\partial^\beta u}{\partial x^\beta}(q)
\right|\le \frac{C_k\,\|u\|_{L^1(B_{R}(q))}}{R^{n+k}}\le\frac{C_k\,\|u\|_{L^1(B_{r})}}{R^{n+k}}
.$$
As a result,
$$ \int_{\partial B_{\rho} }|w(x)|\,|x_1|\,d{\mathcal{H}}^{n-1}_x\le
\rho
\int_{\partial B_{\rho}} |w(x)|\,d{\mathcal{H}}^{n-1}_x\le
\frac{C_k\,\rho^n\,{\mathcal{H}}^{n-1}(\partial B_1)\,\|u\|_{L^1(B_{r})}}{R^{n+k}}.
$$
Combined with~\eqref{KMMARGBnoKKSMMLDSshciCLAALMSifUFFPOP},
and recalling~\eqref{B1},
this information leads to
\begin{eqnarray*}
|\partial_{x}^\alpha u(0)|&\le&
\frac{C_k\,{\mathcal{H}}^{n-1}(\partial B_1)\,\|u\|_{L^1(B_{r})}}{\rho\,|B_1|\,R^{n+k}}\\&=&
\frac{C_k\,n\,\|u\|_{L^1(B_{r})}}{\rho\,R^{n+k}}\\&=&
\frac{C_k\,(k+1)^{n+k+1}\,n\,\|u\|_{L^1(B_{r})}}{k^{n+k} r^{n+k+1}}\\&=&
\frac{(2^{n+1} n k)^{k}\,(k+1)^{n+k+1}\,n\,\|u\|_{L^1(B_{r})}}{|B_1|\,k^{n+k} r^{n+k+1}}
\\&\leq&
\frac{(2^{n+1} n (k+1))^{k+1}\,\|u\|_{L^1(B_{r})}}{|B_1|\,r^{n+k+1}}\\&
=& \frac{C_{k+1}\,\|u\|_{L^1(B_r )}}{r^{n+k+1}}
\end{eqnarray*}
and this completes the inductive step.
\end{proof}

For completeness, we observe that Cauchy's Estimates
(and in particular the explicit value of~$C_k$
determined in~\eqref{CIKAPP})
can be exploited to provide a
proof of the real analyticity of harmonic functions
alternative to the one presented for Theorem~\ref{KMS:HAN0-0} here (see e.g.~\cite{MR1625845}).

For further reference, we also use Cauchy's Estimates in Theorem~\ref{CAUESTIMTH} to deduce a Maximum Principle in Lebesgue spaces for unbounded domains
which can be seen as a counterpart of Corollary~\ref{WEAKMAXPLE}. 

\begin{lemma}\label{VodmvVikfAndfdfNAJowqrq0aqjf}
Let~$p\ge1$ and~$M\ge0$. Let~$\Omega\subseteq\R^n$ be open and connected
and~$u\in C^2(\Omega)\cap
L^p(\Omega)$ be harmonic in~$\Omega$. 

Suppose that
\begin{equation}\label{1-0-0-0234}
{\mbox{$\R^n\setminus\Omega$ is a bounded set}}\end{equation}
and
\begin{equation}\label{0987654--plkmnS-iuh}
\limsup_{\Omega\ni y\to z}u(y)\le M\quad{\mbox{ for all~$z\in\partial\Omega$.}}
\end{equation}
Then, $$\sup_{\Omega} u\le M.$$\end{lemma} 

\begin{proof} We argue for a contradiction and suppose that \begin{equation}\label{0987654--plkmnS-iuh2}\sup_\Omega u\ge M+a,\end{equation} for some~$a>0$. We also take a sequence~$x_k\in\Omega$ such that $$\lim_{k\to+\infty}u(x_k)=\sup_\Omega u.$$ We distinguish two cases, either~$x_k$ is bounded or not.

If~$x_k$ is bounded, up to a subsequence we can assume that~$x_k\to\overline{x}$, for some~$\overline{x}\in\overline\Omega$.
By~\eqref{0987654--plkmnS-iuh}, we infer that~$\overline{x}\in \Omega$. Thus, if~$\rho>0$ is so small that~$B_\rho(\overline{x})\subset\Omega$ we have that~$u\in C(B_\rho(\overline{x}))$ and thus
$$ u(\overline{x})=\lim_{k\to+\infty}u(x_k)=\sup_\Omega u.$$
This and the Strong Maximum Principle in Theorem~\ref{STRONGMAXPLE1}(iii) yield that~$u$ is constant in~$\Omega$,
say~$u=c$ in~$\Omega$ for some~$c\in\R$.

On the one hand, it follows from~\eqref{0987654--plkmnS-iuh} that~$c\le M$.
On the other hand, by~\eqref{0987654--plkmnS-iuh2}, we have that~$c\ge M+a>M$, which is
a contradiction.

We can therefore focus on the case in which~$x_k$ is unbounded. By~\eqref{1-0-0-0234},
we can suppose that~$\R^n\setminus\Omega\Subset B_R$, for some~$R>0$, hence~$B_{r_k}(x_k)\Subset\Omega$
for~$r_k:=|x_k|-(R+1)$ (notice that~$r_k>0$ if~$k$ is large enough).
As a result, using Cauchy's Estimates in Theorem~\ref{CAUESTIMTH}
and the H\"older Inequality, we see that, for every~$x\in B_{r_k/2}(x_k)$,
$$ |\nabla u(x)|\le \frac{C\,\|u\|_{L^1(B_{r_k}(x_k))}}{r_k^{n+1}}\le
\frac{C\,\|u\|_{L^p(\Omega)}\,r_k^{\frac{n(p-1)}{p}}}{r_k^{n+1}}=
\frac{C\,\|u\|_{L^p(\Omega)}}{r_k^{1+\frac{n}p}}
$$
for some constant~$C>0$, depending only on~$n$ and possibly varying from step to step.

Consequently, taking~$k$ sufficiently large such that~$r_k\ge2$ and~$u(x_k)\ge\sup_\Omega u-\frac{a}2$,
for every~$x\in B_1(x_k)$ we have that
$$ u(x)\ge u(x_k)-\|\nabla u\|_{L^\infty(B_1(x_k))}\ge
\sup_\Omega u-\frac{a}2-\frac{C\,\|u\|_{L^p(\Omega)}}{r_k^{1+\frac{n}p}}\ge
M+\frac{a}2-\frac{C\,\|u\|_{L^p(\Omega)}}{r_k^{1+\frac{n}p}}
.$$
That is, if~$k$ is sufficiently large, for every~$x\in B_1(x_k)$
$$ u(x)\ge
M+\frac{a}4\ge\frac{a}4,$$
from which we obtain a contradiction since~$u\in L^p(\Omega)$.
The desired result is thereby proved.
\end{proof}

\section{The Weak Harnack Inequality}

In this section, we present a series of results
that are often called in the literature, possibly under different forms,
the ``Weak Harnack Inequality''.\index{Weak Harnack Inequality}
In this jargon, the word ``weak'' should not be considered as reductive,
but rather highlighting that the 
estimates obtained do not rely only on ``strong'' notions
such as the pointwise values of a solution (as it happens,
for instance, in the Harnack Inequality presented in Section~\ref{HARNAFOR})
but rather on ``weak'' notions such as measures of level sets
and integrals of the solution. 
We do not address here the most general statement
of this theory: for a comprehensive vision, see~\cite{MR0093649,
MR100158,
MR170091,
MR159138,
MR525227,
MR1351007,
MR1814364,
MR2291922,
MR2777537,
MR3380557}
and the references therein.

To start this topic we show,
roughly speaking, that
a harmonic function
bounded from above and with a positive measure set
of values below its maximum is necessarily poitwise well
separated from its maximum inside the domain.
The quantitative formulation of this statement goes as follows:

\begin{theorem} \label{WEHA}Let~$R>r>0$
and~$\eta>0$.
Let~$\Omega$ be an open subset of~$\R^n$ with~$B_R\Subset\Omega$.
Let~$u\in C^2(\Omega)$ be harmonic in~$\Omega$ and such that
\begin{equation}\label{ETA:01}\sup_{B_R}u\le 1.\end{equation}
Assume that
\begin{equation}\label{ETA:02} \big|\{x\in B_R {\mbox{ s.t. }}u(x)
\le0\}\big|\ge\eta.\end{equation}

Then, there exists~$c>0$, depending only on~$n$, $\eta$, $r$ and~$R$,
such that
\begin{equation}\label{DEGDB:DESID} u\le 1-c\qquad{\mbox{in }}\,B_{r}.\end{equation}\end{theorem}

\begin{proof} By the Mean Value Formula in Theorem~\ref{KAHAR}(iii),
\eqref{ETA:01} and~\eqref{ETA:02},
\begin{eqnarray*}&&u(0)=\fint_{B_R} u(x)\,dx
=\frac{1}{|B_R|}
\left(\int_{B_R\cap\{u\le 0\}}u(x)\,dx
+\int_{B_R\cap\{u>0\}}u(x)\,dx
\right)\\&&\qquad\qquad\qquad\le
\frac{1}{|B_R|}\int_{B_R\cap\{u>0\}}u(x)\,dx\le
\frac{|B_R|-\eta}{|B_R|}=1-c_\star,
\end{eqnarray*}
where~$c_\star:=\frac\eta{|B_R|}>0$.

Now we use\footnote{For the sake of generality, we are presenting here a
proof which is independent of the Harnack Inequality that we have discussed so far. However,
if one is willing to use the Harnack Inequality, one can obtain the desired result by considering the harmonic function~$v:=1-u\ge0$ in~$B_R$, with the observation that~$v(0)=1-u(0)\ge c_\star$,
yielding in~$B_r$ that
$$ c_\star\le v(0)\le\sup_{B_r}v \le C\inf_{B_r} v=
C-C\sup_{B_r} u,$$
and thus$$\sup_{B_r} u\le 1-\frac{c_\star}C,$$ as desired.}
Cauchy's Estimate in Theorem~\ref{CAUESTIMTH}
with~$k=1$.
In particular, we take~$C_1$ to be the constant introduced in~\eqref{CIKAPP}
with~$k:=1$, we set~$r_0:=\min\left\{\frac{R}{2},
\frac{c_\star R}{2^{n+2} C_1\,|B_1|}\right\}$ and
we deduce from Theorem~\ref{CAUESTIMTH}
that, for every~$x\in B_{r_0}$,
\begin{equation}\label{WE:OBVsdsdbfS0}
u(x)\le u(0)+r_0\sup_{B_{r_0}} |\nabla u|\le
1-c_\star+\frac{C_1 r_0\|u\|_{L^1(B_R)}}{(R/2)^{n+1}}\le
1-c_\star+\frac{2^{n+1} C_1 r_0\,|B_1| }{R}\le1-\frac{c_\star}2.
\end{equation}
If~$r_0\ge r$, this proves the desired result.
If instead~$r_0<r$, we can iterate the previous
inequality.
To this end,
we define recursively
\begin{equation}\label{RKRECDF}
r_{k+1}:=\frac{(R-r_k)r_0+Rr_k}{R}.
\end{equation} 
Notice that
\begin{equation}\label{POMIA} r_{k+1}\le R.\end{equation}
Indeed, the claim is true for~$r_0$ and then, proceeding inductively,
$$ r_{k+1}=\frac{Rr_0+(R-r_0)r_k}{R}\le
\frac{Rr_0+(R-r_0)R}{R}=R,
$$
which establishes~\eqref{POMIA}.

Now we claim that, for every~$k\in\N$,
there exists~$c_k>0$ such that
\begin{equation}\label{WE:OBVsdsdbfS}
\sup_{B_{r_k}}u\le 1-c_k.
\end{equation}
We prove this by induction. Indeed, when~$k=0$,
the claim in~\eqref{WE:OBVsdsdbfS} follows from~\eqref{WE:OBVsdsdbfS0}
by taking~$c_0:=\frac{c_\star}2$.
Suppose now that~\eqref{WE:OBVsdsdbfS} holds true for some index~$k$
and let~$p\in B_{r_k}$.
We define
$$v(x):=\frac1{c_k}\,\left[
u\left(p+\frac{(R-r_k)x}{R}\right)-1+c_k\right]$$
and we observe that if~$x\in B_R$ then
$$ \left|p+\frac{(R-r_k)x}{R}\right|\le r_k+\frac{(R-r_k)|x|}{R}<R,$$
and thus~$v$ is harmonic in~$B_R$.
Additionally, we have that~$v$ is bounded from above by~$1$.
Furthermore, we observe that
$$ \lim_{\e\searrow0}
\frac{r_k+(R-r_k)\e}{R}\,r_k+
\frac{(R-r_k)\e}{R}=
\frac{r_k}{R}\,r_k<r_k,$$
whence we can
take~$\e_k\in(0,1)$ so small that
$$ \frac{r_k+(R-r_k)\e_k}{R}\,r_k+
\frac{(R-r_k)\e_k}{R}<r_k.$$
With this choice, it follows that, if~$q_k:=(\e_k-1)p$
and~$x\in B_{\e_k}(q_k)$,
we have that
\begin{eqnarray*}&& \left| p+\frac{(R-r_k)x}{R}\right|=
\left| p+\frac{(R-r_k)(x-q_k)}{R}+\frac{(R-r_k)q_k}{R}\right|\\&&\qquad
\le \left| p+\frac{(R-r_k)(\e_k-1)p}{R}\right|+\left|\frac{(R-r_k)(x-q_k)}{R}\right|\\
&&\qquad
< \frac{r_k+(R-r_k)\e_k}{R}\,r_k+
\frac{(R-r_k)\e_k}{R}\\&&\qquad<r_k
\end{eqnarray*}
and therefore~$p+\frac{(R-r_k)x}{R}\in B_{r_k}$.

Consequently,
by the inductive
hypothesis,
\begin{eqnarray*}&& \{x\in B_R {\mbox{ s.t. }}v(x)
\le0\}=\left\{x\in B_R {\mbox{ s.t. }}u\left(p+\frac{(R-r_k)x}{R}\right)\le
1-c_k\right\}\\&&\qquad\qquad\supseteq
\left\{x\in B_R {\mbox{ s.t. }}
p+\frac{(R-r_k)x}{R}\in B_{r_k}\right\}\supseteq
B_{\e_k}(q_k)\end{eqnarray*}
and, as a result,
$$ \eta_k:=\big|\{x\in B_R {\mbox{ s.t. }}v(x)
\le0\}\big|>0.$$
We can therefore apply inequality~\eqref{WE:OBVsdsdbfS0}
to~$v$ and find that
$$ \sup_{B_{r_0}} v\le1-\widetilde{c}_k,$$
for some~$\widetilde{c}_k>0$.
This gives that, for every~$x\in B_{r_0}$,
$$
u\left(p+\frac{(R-r_k)x}{R}\right)
\le c_k(1-\widetilde{c}_k)+1-c_k
=1-c_k\,\widetilde{c}_k$$
and therefore
$$ \sup_{B_{(R-r_k)r_0/R}(p)}u\le1-c_k\,\widetilde{c}_k.$$
Since~$p$ is an arbitrary point in~$ B_{r_k}$, we thus infer that
$$ \sup_{B_{((R-r_k)r_0+Rr_k)/R}}u\le1-c_k\,\widetilde{c}_k,$$
which, in light of the definition of~$r_{k+1}$
in~\eqref{RKRECDF}, concludes the inductive step and establishes~\eqref{WE:OBVsdsdbfS}.

We also stress that, by~\eqref{RKRECDF} and~\eqref{POMIA},
$$ r_{k+1}-r_k=\frac{(R-r_k)r_0}{R}\ge0$$
hence~$r_{k+1}$ is a bounded and monotone increasing
sequence. We can therefore define~$R_\star$ to be its limit
and we deduce from~\eqref{RKRECDF} that
$$ R_\star=\frac{(R-R_\star)r_0+R R_\star}{R},$$
leading to~$ (R-R_\star)r_0=0$ and thus~$R=R_\star$.

This observation gives that
$$ \lim_{k\to+\infty} r_k=R>r,$$
hence we can choose~$k_\star\in\N$ such that~$r_{k_\star}>r$.
In this way, the desired claim in~\eqref{DEGDB:DESID}
follows from~\eqref{WE:OBVsdsdbfS} with~$k:=k_\star$.
\end{proof}

As a variant of Theorem~\ref{WEHA},
one can control the supremum of a harmonic function by its
$L^p$-norm
(interestingly, the case~$p\in(0,1]$ is also allowed).
This result will be stated precisely in Theorem~\ref{CORODASOLTFJ}
and its proof will also rely on the following very useful
observation which allows us to ``reabsorb a term'' on the left
hand side of a scaled inequality:

\begin{lemma}\label{S-O-GIAHJAQ}
Let~$a>b$ and~$\varphi:[a,b]\to\R$
be a bounded function. Assume that there exist~$\vartheta\in[0,1)$,
$A\ge0$, $B\ge0$ and~$\gamma>0$ such that,
for every~$r$, $R\in[a,b]$ with~$r<R$
we have that
\begin{equation}\label{BA7u6y5t4SE}
\varphi(r)\le\vartheta\varphi(R)+\frac{A}{(R-r)^\gamma}+B.\end{equation}

Then, there exists~$C>0$, depending only on~$\vartheta$
and~$\gamma$ such that for every~$r$, $R\in[a,b]$ with~$r<R$
we have that
$$ \varphi(r)\le C\,\left(\frac{A}{(R-r)^\gamma}+B\right).$$
\end{lemma}

\begin{proof} Fix~$r$, $R\in[a,b]$ with~$r<R$
and let~$\tau:=\frac{1+\vartheta^{\frac1\gamma}}2$.
Notice that
\begin{equation}\label{MORAUTAU}
\tau\in\big(\vartheta^{\frac1\gamma},1\big)\end{equation}
since~$\vartheta<1$. We also define~$r_0:=r$ and recursively,
for each~$k\in\N$,
$$ r_{k+1}:=r_k+ (1-\tau)\,\tau^k\,(R-r).$$
We point out that the sequence~$r_k$ is monotone increasing,
thanks to~\eqref{MORAUTAU} and can be written in closed form as
\begin{eqnarray*}r_{k+1}&=&
r_{k-1}+ (1-\tau)\,\tau^{k-1}\,(R-r)
+ (1-\tau)\,\tau^k\,(R-r)\\&=&r_{k-2}+ (1-\tau)\,\tau^{k-2}\,(R-r)
r_{k-1}+ (1-\tau)\,\tau^{k-1}\,(R-r)
+ (1-\tau)\,\tau^k\,(R-r)\\&=&\dots\\&=&
r_0+\sum_{j=0}^k
(1-\tau)\,\tau^j\,(R-r)\\
&=&r+\frac{(1-\tau)(R-r) (1-\tau^{k+1})}{1-\tau}
.\end{eqnarray*}
As a result,
\begin{equation*} \lim_{k\to+\infty}r_k=r+\frac{(1-\tau)(R-r) }{1-\tau}
=R.\end{equation*}
Now we claim that, for every~$k\in\N$,
\begin{equation}\label{INDUSTRAWUIDMFWHIANFANDD}
\varphi(r)\le\vartheta^k\varphi(r_k)+\left(
\frac{A}{(1-\tau)^\gamma(R-r)^\gamma}+B
\right)\sum_{i=0}^{k-1}\left(\frac{\vartheta}{\tau^\gamma}\right)^i.
\end{equation}
Indeed, when~$k=0$ the above formula is obvious.
When~$k=1$, we apply~\eqref{BA7u6y5t4SE} (used here with~$r_1$
in place of~$R$) and see that
$$ \varphi(r)\le\vartheta\varphi(r_1)+\frac{A}{(r_1-r)^\gamma}+B
=\vartheta\varphi(r_1)+\frac{A}{(1-\tau)^\gamma(R-r)^\gamma}+B,$$
and this is precisely~\eqref{INDUSTRAWUIDMFWHIANFANDD}
with~$k=1$.

Thus, to complete the proof of~\eqref{INDUSTRAWUIDMFWHIANFANDD},
we argue recursively, assuming that this formula
holds true for the index~$k$ and we prove it for the index~$k+1$.
To this end, we use the inductive assumption
and~\eqref{BA7u6y5t4SE} (utilized here with~$r_k$ in place of~$r$
and~$r_{k+1}$
in place of~$R$) and see that
\begin{eqnarray*}
\varphi(r)&\le&\vartheta^k\varphi(r_k)+\left(
\frac{A}{(1-\tau)^\gamma(R-r)^\gamma}+B
\right)\sum_{i=0}^{k-1}\left(\frac{\vartheta}{\tau^\gamma}\right)^i
\\&\le&\vartheta^k\left(
\vartheta\varphi(r_{k+1})+\frac{A}{(r_{k+1}-r_k)^\gamma}+B
\right)+\left(
\frac{A}{(1-\tau)^\gamma(R-r)^\gamma}+B
\right)\sum_{i=0}^{k-1}\left(\frac{\vartheta}{\tau^\gamma}\right)^i\\
&=&\vartheta^{k+1}\varphi(r_{k+1})+\vartheta^k\left(
\frac{A}{
(1-\tau)^\gamma\,\tau^{\gamma k}\,(R-r)^\gamma}+B\right)+\left(
\frac{A}{(1-\tau)^\gamma(R-r)^\gamma}+B
\right)\sum_{i=0}^{k-1}\left(\frac{\vartheta}{\tau^\gamma}\right)^i\\
&\le&\vartheta^{k+1}\varphi(r_{k+1})+\vartheta^k\left(
\frac{A}{
(1-\tau)^\gamma\,\tau^{\gamma k}\,(R-r)^\gamma}+\frac{B}{\tau^{\gamma k}}\right)+\left(
\frac{A}{(1-\tau)^\gamma(R-r)^\gamma}+B
\right)\sum_{i=0}^{k-1}\left(\frac{\vartheta}{\tau^\gamma}\right)^i\\&=&
\vartheta^{k+1}\varphi(r_{k+1})+\left(
\frac{A}{(1-\tau)^\gamma(R-r)^\gamma}+B
\right)\sum_{i=0}^{k}\left(\frac{\vartheta}{\tau^\gamma}\right)^i,
\end{eqnarray*}
and this proves~\eqref{INDUSTRAWUIDMFWHIANFANDD}.

We also observe that~$\frac{\vartheta}{\tau^\gamma}<1$,
thanks to~\eqref{MORAUTAU},
and we can therefore send~$k\to+\infty$ in~\eqref{INDUSTRAWUIDMFWHIANFANDD},
concluding that
\begin{equation*} \varphi(r)\le \left(
\frac{A}{(1-\tau)^\gamma(R-r)^\gamma}+B
\right)\sum_{i=0}^{+\infty}\left(\frac{\vartheta}{\tau^\gamma}\right)^i=
\left(
\frac{A}{(1-\tau)^\gamma(R-r)^\gamma}+B
\right)\,\frac{1}{1-\frac{\vartheta}{\tau^\gamma}}.\qedhere\end{equation*}
\end{proof}

\begin{theorem}\label{CORODASOLTFJ}
Let~$R>r>0$ and~$p>0$.
Let~$\Omega$ be an open subset of~$\R^n$ with~$B_R\Subset\Omega$.
Let~$u\in C^2(\Omega)$ be harmonic in~$\Omega$.

Then, there exists~$C(n,p)>0$, depending only on~$n$ and~$p$,
such that
\begin{equation}\label{CORODASOLTFJ0}
\sup_{B_r}u_+\le\frac{ C(n,p)}{(R-r)^{\frac{n}p}}\left( \int_{B_R} (u_+(x))^p\,dx\right)^{\frac1p}.\end{equation}
\end{theorem}

\begin{proof} We first establish~\eqref{CORODASOLTFJ0} when~$p\ge1$.

For this, let~$\zeta\in B_R$. We use
the Mean Value Formula in Theorem~\ref{KAHAR}(iii) to see that,
if~$\rho\in(0,R-|\zeta|]$ then
\begin{equation*} u(\zeta)=\fint_{B_\rho(\zeta)} u(x)\,dx
\le
\fint_{B_\rho(\zeta)} u_+(x)\,dx.
\end{equation*}
Thus, if~$p\ge1$, one can use the H\"{o}lder Inequality and find that,
for all~$\zeta\in B_r$,
$$ u(\zeta)\leq\fint_{B_{R-r}(\zeta)} u_+(x)\,dx\le
\left( \fint_{B_{R-r}(\zeta)}( u_+(x))^p\,dx\right)^{\frac1p}\le
\left( \frac1{|B_{R-r}|}\int_{B_R}( u_+(x))^p\,dx\right)^{\frac1p}.
$$
Since the above estimate is obviously also true when~$u(\zeta)\le0$,
we thus deduce that
$$ u_+(\zeta)\le
\left( \frac1{|B_{R-r}|}\int_{B_R}( u_+(x))^p\,dx\right)^{\frac1p},$$
from which we obtain~\eqref{CORODASOLTFJ0} when~$p\ge1$.

Suppose now that~$p\in(0,1)$. 
In this case, we have that
$$ \int_{B_R} u_+(x)\,dx=
\int_{B_R} (u_+(x))^{p+(1-p)}\,dx\le\left(\sup_{B_R} u_+\right)^{1-p}
\int_{B_R} (u_+(x))^p\,dx.$$
Therefore, using~\eqref{CORODASOLTFJ0} with~$p=1$ and setting
$$ \varphi(t):=\sup_{B_t} u_+,$$
we find that
$$ \varphi(r)=\sup_{B_r} u_+\le
\frac{ C(n,1)}{(R-r)^n}\int_{B_R} u_+(x)\,dx
\le \frac{ C(n,1)}{(R-r)^n}\left(\varphi(R)\right)^{1-p}
\int_{B_R} (u_+(x))^p\,dx.$$
We thus exploit Young's Inequality with exponents~$\frac1{1-p}$
and~$\frac1p$, obtaining that
\begin{equation*}
\varphi(r)\le\frac{\varphi(R)}2+
\frac{C}{(R-r)^{\frac{n}p}}\left(
\int_{B_R} (u_+(x))^p\,dx\right)^{\frac1p}.
\end{equation*}
We can therefore utilize
Lemma~\ref{S-O-GIAHJAQ} and conclude, up to renaming~$C$, that
\begin{equation*} \varphi(r)\le \frac{C}{(R-r)^{\frac{n}p}}\left(
\int_{B_R} (u_+(x))^p\,dx\right)^{\frac1p}.\qedhere\end{equation*}
\end{proof}

A counterpart of Theorem~\ref{CORODASOLTFJ}
(actually relying merely on the Harnack Inequality presented
in Section~\ref{HARNAFOR})
is given by the following result, in which
the $L^p$-norm of a nonnegative harmonic
function is bounded from above by the infimum
in a smaller ball.

\begin{theorem}\label{COROD8765ASOLTFJ-WEAKSOR}
Let~$R_0>R>r>0$ and~$p>0$.
Let~$\Omega$ be an open subset of~$\R^n$
with~$B_{R_0}\Subset\Omega$.
Let~$u\in C^2(\Omega)$ be nonnegative and harmonic in~$\Omega$.

Then, there exists~$c>0$, depending only on~$n$, $p$, $r$, $R$ and~$R_0$
such that
\begin{equation}\label{CORODA09876547899SOLTFJ0}
\inf_{B_r}u \ge c\,\left( \int_{B_R} (u(x))^p\,dx
\right)^{\frac1p}.\end{equation}
\end{theorem}

\begin{proof} 
By the
Harnack Inequality in Theorem~\ref{PRE:HA},
for every~$x\in B_R$,
$$ u(x)\le\left(\frac{R_0}{R_0-R}\right)^{n-2}\,\frac{R_0+R}{R_0-R}\,u(0)
.$$
As a result,
\begin{equation}\label{OKSMIBJSGIUJSFBGRGEJAOKD}
\int_{B_R} (u(x))^p\,dx\le
\left(\frac{R_0}{R_0-R}\right)^{(n-2)p}\,\frac{(R_0+R)^p}{(R_0-R)^p}
\,(u(0))^p\,|B_R|.
\end{equation}
Also, using again the
Harnack Inequality in Theorem~\ref{PRE:HA},
$$
\left(\frac{R_0}{R_0+r}\right)^{n-2}\,
\frac{R_0-r}{R_0+r}\,u(0)\le\inf_{B_r}
u.$$
Combining this and~\eqref{OKSMIBJSGIUJSFBGRGEJAOKD},
we conclude that
\begin{equation*}
\int_{B_R} (u(x))^p\,dx\le\frac{(R_0+R)^p\,(R_0+r)^{(n-1)p}}{
(R_0-R)^{(n-1)p}\,(R_0-r)^p}\, |B_R|\,
\left(\inf_{B_r}u\right)^p.\qedhere
\end{equation*}
\end{proof}

We observe that, when~$p\in(0,1]$,
Theorem~\ref{COROD8765ASOLTFJ-WEAKSOR} holds true also\footnote{Indeed,
if~$p\in(0,1]$, one considers~$P:=\frac1p\ge1$
and uses the H\"{o}lder Inequality and the
Mean Value Formula in Theorem~\ref{KAHAR}(iii) to see that
\begin{equation*}
\fint_{B_{R}} (u(x))^p\,dx\le
\left(\fint_{B_{R}} u(x)\,dx\right)^{\frac1P}=(u(0))^p.\end{equation*}
Thus, since, by the
Harnack Inequality in Theorem~\ref{PRE:HA},
\begin{equation*} \left(\frac{R}{R+r}\right)^{n-2}\,
\frac{R-r}{R+r}\,u(0)\le\inf_{B_r}u,
\end{equation*}
it follows that, say when~$n\ge2$,
$$\left(\fint_{B_{R}} (u(x))^p\,dx\right)^{\frac1p}\le u(0)\le
\left(\frac{R+r}{R}\right)^{n-2}\,
\frac{R+r}{R-r}\,\inf_{B_r}u\le
\left(\frac{2R}{R}\right)^{n-2}\,
\frac{2R}{R-r}\,\inf_{B_r}u=
\frac{2^{n-1}\,R}{R-r}\,\inf_{B_r}u
$$
and this gives~\eqref{CORODA09876547899SOLTFJ0}
when~$p\in(0,1]$ with~$R_0=R$.}
with~$R_0=R$. Instead, for large~$p$,
it is necessary to assume~$R_0>R$ since, if~$u$ is as in~\eqref{HESEM},
then, for small~$\e$,
\begin{eqnarray*}
\int_{B_R} (u(x))^p\,dx&\ge&
\int_{B_\e((R-\e)e_1)}
\frac{((R+\e)^2-|x|^2)^p}{n^p\,|B_1|^p\,(R+\e)^p\,|x-(R+\e)e_1|^{np}}\,dx\\
&\ge&
\int_{B_\e((R-\e)e_1)}
\frac{((R+\e)^2-R^2)^p}{n^p\,|B_1|^p\,(2R)^p\,(3\e)^{np}}\,dx
\\&\ge&\frac{c}{\e^{np-p-n}}
\end{eqnarray*}
for some~$c>0$ depending only on~$n$, $p$ and~$R$,
and the latter term in the above chain of inequalities
diverges as~$\e\searrow0$ when~$n\ge2$ and~$p$ is sufficiently large
(for further discussions on the range of~$p$ also in more general
contexts
see e.g.~\cite[Theorem~8.18 and Problem~9.12]{MR1814364}
and~\cite[Theorem~4.15]{MR2777537}).

\section{The Boundary Harnack Inequality}\label{BOUNAHSRGA9876543}

For completeness, and without aiming at being exhaustive,
we recall here a boundary version of the Harnack Inequality. For this, we use the notation\index{Boundary Harnack Inequality}
$$ B_r^+:=\{x=(x_1,\dots,x_n){\mbox{ s.t. }}x_n>0\}.$$

\begin{theorem} Let~$u$, $v\in C^2(B_2^+)\cap C^1(\overline{B_2^+})$ be harmonic functions in~$B_2^+$.
Assume that~$u$, $v>0$ in~$B_2^+$
and~$u=v=0$ along~$\{x_n=0\}$. Then, there exists a constant~$C>1$ only depending on~$n$ such that, for every~$x\in B_{1/8}^+$,
\begin{equation}\label{BOUNGHARNBAS} \frac{1}C\,\frac{u(e_n/8)}{v(e_n/8)}\le\frac{u(x)}{v(x)}\le C\,\frac{u(e_n/8)}{v(e_n/8)}.\end{equation}
\end{theorem}

\begin{proof} We observe that, by possibly replacing~$u$ with~$\frac{u}{u(e_n/8)}$
and~$v$ with~$\frac{v}{v(e_n/8)}$, we can assume that
\begin{equation}\label{NORMAEN8}
u\left(\frac{e_n}8\right)=1=v\left(\frac{e_n}8\right).
\end{equation}
Let
$$ B_r^-:=\{x=(x_1,\dots,x_n){\mbox{ s.t. }}x_n<0\}.$$
Given~$x=(x_1,\dots,x_n)$ we let~$x_\star:=(x_1,\dots,-x_n)$. For every~$x\in B_2$, let also
$$ u_\star(x):=\begin{dcases}
u(x) & {\mbox{ if }}x\in \overline{B_2^+},\\
-u(x_\star)& {\mbox{ if }}x\in \overline{B_2^-}
.\end{dcases}
$$
Notice that this is a good definition since~$u=0$ and~$x_\star=x$ along~$\{x_n=0\}$.
We observe that
\begin{equation}\label{HARMEXTDISP}
{\mbox{$u_\star$ is harmonic in $B_2$.}}
\end{equation}
Indeed, if~$\varphi\in C^\infty_0(B_2)$, we let~$\psi(x):=\varphi(x_\star)$.
Thus, the use of the second
Green Identity~\eqref{GRr2}
leads to
\begin{eqnarray*}
\int_{B_2} u_\star(x)\,\Delta\varphi(x)\,dx&=&
\lim_{\e\searrow0} \int_{B_2\cap\{x_n>\e\}} u_\star(x)\,\Delta\varphi(x)\,dx+
\int_{B_2\cap\{x_n<-\e\}} u_\star(x)\,\Delta\varphi(x)\,dx\\&=&
\lim_{\e\searrow0} \int_{B_2\cap\{x_n>\e\}} u(x)\,\Delta\varphi(x)\,dx-
\int_{B_2\cap\{x_n>\e\}} u(x)\,\Delta\psi(x)\,dx\\&=&
\lim_{\e\searrow0} \int_{B_2\cap\{x_n=\e\}} \big(
u(x)\,\partial_\nu\varphi(x) -\varphi(x)\,\partial_\nu u(x)\big)\,dx\\&&\qquad-
\lim_{\e\searrow0} \int_{B_2\cap\{x_n=\e\}} \big(
u(x)\,\partial_\nu\psi(x) -\psi(x)\,\partial_\nu u(x)\big)\,dx\\&=&
\lim_{\e\searrow0} \int_{B_2\cap\{x_n=\e\}} \big( \psi(x)-\varphi(x)\big)\,\partial_\nu u(x)\,dx\\&=&0.\end{eqnarray*}
In view of
Weyl's Lemma (i.e., recall Lemma~\ref{WEYL}) 
this proves~\eqref{HARMEXTDISP}.

As a result, we can exploit the Poisson Kernel representation
(recall Theorems~\ref{POIBALL1} and~\ref{POIBALL})
and obtain that, for all~$x\in B_1^+$,
\begin{equation}\label{LUNSNNDFDGAN MDJTAYHHAHNDJIHASLN}
\begin{split}
n\,|B_1|\, u(x)\,=\,&
\int_{\partial B_1} 
\frac{(1-|x|^2)\,u_\star(y)}{|x-y|^n}\,d{\mathcal{H}}^{n-1}_y\\=\,&
\int_{\partial B_1^+} 
\frac{(1-|x|^2)\,u(y)}{|x-y|^n}\,d{\mathcal{H}}^{n-1}_y
-
\int_{\partial B_1^-} 
\frac{(1-|x|^2)\,u(y_\star)}{|x-y|^n}\,d{\mathcal{H}}^{n-1}_y
\\=\,&
\int_{\partial B_1^+} 
\frac{(1-|x|^2)\,u(y)}{|x-y|^n}\,d{\mathcal{H}}^{n-1}_y
-
\int_{\partial B_1^+} 
\frac{(1-|x|^2)\,u(y)}{|x-y_\star|^n}\,d{\mathcal{H}}^{n-1}_y\\=\,&\int_{\partial B_1^+}\left(\frac{1}{|x-y|^n}-\frac{1}{|x-y_\star|^n}
\right)
(1-|x|^2)\,u(y)\,
d{\mathcal{H}}^{n-1}_y.
\end{split}
\end{equation}

Now we point out that, if~$a$, $b\in B_2\setminus B_{1/2}$ then
\begin{equation}\label{QUNESIEMSPVRISDMOJDNC}\begin{split}&
|a|^{-n}-|b|^{-n}=|b+(a-b)|^{-n}-|b|^{-n}\\&\qquad
=-n\int_0^1 |b+t(a-b)|^{-n-2}(b+t(a-b))\cdot (a-b)\,dt .\end{split}
\end{equation}
Also, if~$|a-b|\le1/4$, for each~$t\in[0,1]$ we have that
$$ |b+t(a-b)|\ge |b|-|a-b|\ge\frac12-\frac14=\frac14.$$
Consequently,
\begin{equation}\label{VGRMSDCHKIAP098766kJA6}
|a|^{-n}-|b|^{-n}\le n\int_0^1 |b+t(a-b)|^{-n-1}|a-b|\,dt \le 4^{n+1}n\,|a-b|
.\end{equation}
Also, if~$y\in\partial B_1$ with~$y_n\ge 3/4$
and~$x\in B_{1/4}$, for every~$t\in[0,1]$
we have that
$$x_n+y_n-2tx_n\ge \frac34+(1-2t)x_n\ge\frac12,$$
whence
we can use~\eqref{QUNESIEMSPVRISDMOJDNC}
with~$a:=x-y$ and~$b:=x_\star-y$ and find that, for every~$x\in B_{1/4}^+$,
\begin{equation}\label{QUNESIEMSPVRISDMOJDNC-2}\begin{split}&
|x-y|^{-n}-|x-y_\star|^{-n}=|x-y|^{-n}-|x_\star-y|^{-n}\\&\qquad
=-2nx_n\int_0^1 |x_\star-y+2tx_n e_n|^{-n-2}(x_\star-y+2tx_n e_n)\cdot e_n\,dt \\&\qquad
=2nx_n\int_0^1 |x_\star-y+2tx_n e_n|^{-n-2}(x_n+y_n-2tx_n )\,dt 
\\&\qquad\ge
\frac{nx_n}{4^{n+2}}
.\end{split}
\end{equation}
Thus, for every~$x\in B_{1/8}^+$, we can use~\eqref{VGRMSDCHKIAP098766kJA6}
with~$a:=x-y$ and~$b:=x_\star-y$ to obtain that
\begin{equation}\label{0okmBYHNSDPJIAMSK}
\int_{\partial B_1^+}\left(\frac{1}{|x-y|^n}-\frac{1}{|x-y_\star|^n}
\right)
(1-|x|^2)\,u(y)\,
d{\mathcal{H}}^{n-1}_y\le
Cx_n\,\int_{\partial B_1^+}u(y)\,
d{\mathcal{H}}^{n-1}_y.
\end{equation}
Besides, for every~$x\in B_{1/8}^+$, we notice that, for every~$y\in B_1^+$, it holds that~$
|x_n-y_n|\le |x_n|+|y_n|= x_n+y_n$,
whence
$$ |x-y_\star|^2=\sum_{i=1}^{n-1} (x_i-y_i)^2+(x_n+y_n)^2
\ge\sum_{i=1}^{n-1} (x_i-y_i)^2+(x_n-y_n)^2=|x-y|^2.$$
Accordingly,
we can utilize~\eqref{QUNESIEMSPVRISDMOJDNC-2} to obtain that
\begin{equation}\label{gfsUJNKYJHJJA MVAJMDD}\begin{split}&
\int_{\partial B_1^+}\left(\frac{1}{|x-y|^n}-\frac{1}{|x-y_\star|^n}
\right)
(1-|x|^2)\,u(y)\,
d{\mathcal{H}}^{n-1}_y\\&\qquad\ge
\int_{(\partial B_1^+)\cap\{y_n\ge3/4\}}\left(\frac{1}{|x-y|^n}-\frac{1}{|x-y_\star|^n}
\right)
(1-|x|^2)\,u(y)\,
d{\mathcal{H}}^{n-1}_y\\&\qquad
\ge
\frac{x_n}C\,\int_{(\partial B_1^+)\cap\{y_n\ge3/4\}}u(y)\,
d{\mathcal{H}}^{n-1}_y.\end{split}
\end{equation}
for a suitable~$C>1$.

Now we claim that
\begin{equation}\label{DESFVDVKJV7AIJM}
\sup_{[-1,1]^{n-1}\times\left[0,\frac{1}{32}\right]} u\le C,
\end{equation}
for some~$C>1$ (possibly larger than the previous~$C$).
To check this, given~$p=(p',p_n)\in \left[-\frac{33}{32},\frac{33}{32}\right]^{n-1}\times\left[0,\frac1{16}\right]$
we apply the Harnack Inequality of Corollary~\ref{HAR:BAL} in~$B_{p_n/2}(p)$
to see that
$$ u(p',p_n)\le \sup_{B_{p_n/2}(p)} u\le C\inf_{B_{p_n/2}(p)} u\le C u\left(p', \frac{3}2p_n\right)$$
for some~$C>1$.
We can therefore iterate this inequality, letting~$\eta(p_n):=\frac{\ln(1/(8p_n))}{\ln(3/2)}$
and choosing~$k\in\N\cap [\eta(p_n),\eta(p_n)+1)$, thus finding that
$$ u(p)\le C^k u\left(p', \left(\frac{3}2\right)^k p_n\right).
$$
Since
\begin{eqnarray*}&&
\frac{1}{8p_n}=e^{\ln(1/(8p_n))}=\left(\frac{3}2\right)^{\eta(p_n)} \le\left(\frac{3}2\right)^k\\
{\mbox{and }}&& \left(\frac{3}2\right)^k\le \left(\frac{3}2\right)^{\eta(p_n)+1}=
\frac{3}{16p_n},
\end{eqnarray*}
we thereby obtain that
\begin{equation}\label{zxfghuiZXDFGUI908765U9KM}
u(p)\le C^{\eta(p_n)+1}\sup_{\left[-\frac{33}{32},\frac{33}{32}\right]^{n-1}\times\left[\frac18,\frac{3}{16}\right]} u\le
\frac{C}{p_n^{C_\star}}\;
\sup_{\left[-\frac{33}{32},\frac{33}{32}\right]^{n-1}\times\left[\frac18,\frac{3}{16}\right]} u,
\end{equation}
for some~$C_\star>0$,
up to renaming~$C>1$ at each step.

Now, let~$q\in [-1,1]^{n-1}\times\left[0,\frac{1}{32}\right]$. We exploit the
local boundedness result in Theorem~\ref{CORODASOLTFJ}
(specifically, the estimate in~\eqref{CORODASOLTFJ0},
used here with~$p:=1/(2C_\star)$) and, making use of~\eqref{zxfghuiZXDFGUI908765U9KM}, we find that, possibly renaming~$C$,
\begin{equation}\label{DESFVDVKJV7AIJM2}\begin{split}
&u(q)=u_\star(q)\le C\left( \int_{
\left[-\frac{33}{32},\frac{33}{32}\right]^{n-1}\times\left[-\frac1{16},\frac1{16}\right]
} |u_\star(p)|^{\frac1{2C_\star}}\,dp\right)^{2C_\star}\\&\qquad\qquad\qquad\qquad
=
C\left(2 \int_{
\left[-\frac{33}{32},\frac{33}{32}\right]^{n-1}\times\left[0,\frac1{16}\right]
} |u(p)|^{\frac1{2C_\star}}\,dp\right)^{2C_\star}\\&\qquad\qquad\qquad\qquad
\le
C\left(\int_{
\left[-\frac{33}{32},\frac{33}{32}\right]^{n-1}\times\left[0,\frac1{16}\right]
} \frac1{\sqrt{p_n}}
\,dp\right)^{2C_\star}\,
\sup_{\left[-\frac{33}{32},\frac{33}{32}\right]^{n-1}\times\left[\frac18,\frac{3}{16}\right]} u\\
&
\qquad\qquad\qquad\qquad\le C\,
\sup_{\left[-\frac{33}{32},\frac{33}{32}\right]^{n-1}\times\left[\frac18,\frac{3}{16}\right]} u
.\end{split}\end{equation}
Additionally, recalling the normalization in~\eqref{NORMAEN8},
using the Harnack Inequality in Corollary~\ref{NSDio} with~$\Omega:=\left(-\frac{35}{32},\frac{35}{32}\right)^{n-1}\times\left(\frac1{16},\frac{1}{4}\right)$
and~$\Omega':=\left(-\frac{33}{32},\frac{33}{32}\right)^{n-1}\times\left(\frac18,\frac{3}{16}\right)$, and possibly renaming~$C$,
$$ \sup_{\left[-\frac{33}{32},\frac{33}{32}\right]^{n-1}\times\left[\frac18,\frac{3}{16}\right]} u
\le C\inf_{\left[-\frac{33}{32},\frac{33}{32}\right]^{n-1}\times\left[\frac18,\frac{3}{16}\right]} u\le
Cu\left(\frac{e_n}8\right)=C.$$
This and~\eqref{DESFVDVKJV7AIJM2} yield~\eqref{DESFVDVKJV7AIJM},
as desired.

In addition, using the Harnack Inequality in Corollary~\ref{NSDio} with~$\Omega:=
B_{3/2}^+\cap \left\{ x_n >\frac{1}{64}\right\}$
and~$\Omega':=B_1^+\cap\left\{x_n>\frac{1}{32}\right\}$,
\begin{equation*}
\sup_{ B_1^+\cap\left\{x_n>\frac{1}{32}\right\}} u\le C\inf_{B_1^+\cap\left\{x_n>\frac{1}{32}\right\}}u.
\end{equation*}
This information and
the normalization
in~\eqref{NORMAEN8} lead to
\begin{equation}\label{DESFVDVKJV7AIJM-A}
1=u\left(\frac{e_n}8\right)
\le
\sup_{ B_1^+\cap\left\{x_n>\frac{1}{32}\right\}} u\le C\inf_{B_1^+\cap\left\{x_n>\frac{1}{32}\right\}}u
\end{equation}
and
\begin{equation}\label{DESFVDVKJV7AIJM-B}
\sup_{ B_1^+\cap\left\{x_n>\frac{1}{32}\right\}} u\le C\inf_{B_1^+\cap\left\{x_n>\frac{1}{32}\right\}}
u\le Cu\left(\frac{e_n}8\right)=C.
\end{equation}
As a byproduct of~\eqref{DESFVDVKJV7AIJM}
and~\eqref{DESFVDVKJV7AIJM-B}, we have that
$$ \int_{\partial B_1^+}u(y)\,
d{\mathcal{H}}^{n-1}_y\le C$$
and thus, in light of~\eqref{0okmBYHNSDPJIAMSK},
if~$x\in B_{1/8}^+$,
\begin{equation}\label{CIMDDS-OSKMDSMIAMO-a}
\int_{\partial B_1^+}\left(\frac{1}{|x-y|^n}-\frac{1}{|x-y_\star|^n}
\right)
(1-|x|^2)\,u(y)\,
d{\mathcal{H}}^{n-1}_y\le
Cx_n,
\end{equation}
up to renaming~$C$.

Similarly, from~\eqref{gfsUJNKYJHJJA MVAJMDD} and~\eqref{DESFVDVKJV7AIJM-A},
and possibly renaming~$C$ once again,
if~$x\in B_{1/8}^+$,
\begin{equation}\label{CIMDDS-OSKMDSMIAMO-b}
\int_{\partial B_1^+}\left(\frac{1}{|x-y|^n}-\frac{1}{|x-y_\star|^n}
\right)
(1-|x|^2)\,u(y)\,
d{\mathcal{H}}^{n-1}_y\ge
\frac{x_n}C.\end{equation}

Therefore, gathering~\eqref{LUNSNNDFDGAN MDJTAYHHAHNDJIHASLN},
\eqref{CIMDDS-OSKMDSMIAMO-a}
and~\eqref{CIMDDS-OSKMDSMIAMO-b},
and possibly renaming~$C$, we deduce that, for all~$x\in B_{1/8}^+$, 
$$ \frac{x_n}C\le u(x)\le C{x_n},$$
and a similar estimate holds for~$v$ as well.

Consequently, for all~$x\in B_{1/8}^+$,
$$ \frac{u(x)}{v(x)}\le\frac{u(x)}{x_n/C}\le\frac{Cx_n}{x_n/C}=C^2$$
and 
$$ \frac{u(x)}{v(x)}\ge\frac{u(x)}{Cx_n}\ge\frac{x_n/C}{Cx_n}=\frac1{C^2}.$$
These observations complete the proof of~\eqref{BOUNGHARNBAS}.
\end{proof}

For a thorough presentation of Boundary Harnack Inequalities (also known
as Carleson Estimates) see~\cite{MR293114,
MR466593,
MR513884,
MR513885,
MR620271,
MR676988,
MR803243,
MR2145284,
MR2464701,
MR3393271,
MR3556055,
MR4093736}
and the references therein.

\begin{figure}
                \centering
                \includegraphics[width=.4\linewidth]{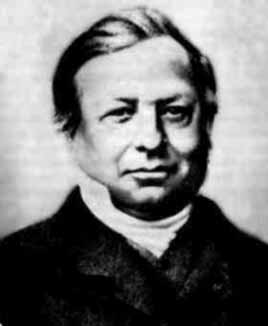}
        \caption{\sl Portrait of Joseph Liouville
        (Public Domain image from
        Wikipedia).}\label{LiPORLROA7789GIJ7sLIOU1polFUMHDNOJHNFOJEDSOL1}
\end{figure}

\section{Liouville's Theorem}

\begin{figure}
                \centering
                \includegraphics[width=.34\linewidth]{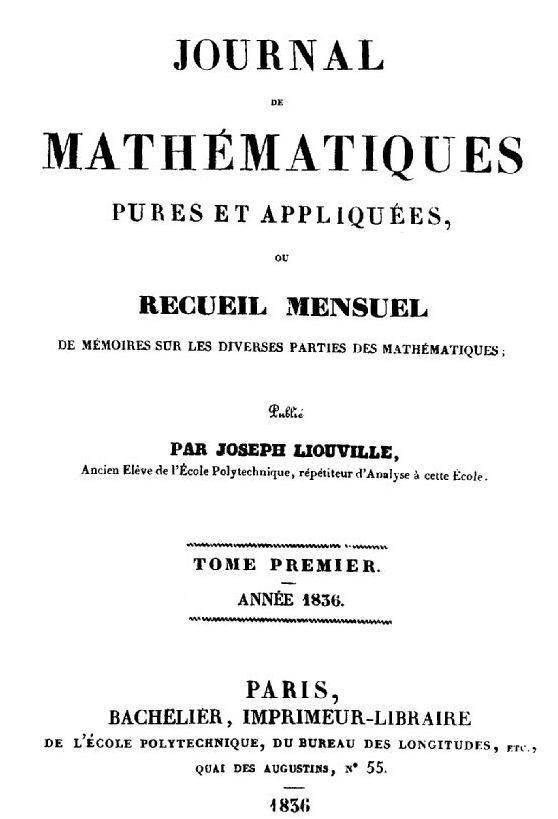}
        \caption{\sl Title page of the first volume of Liouville's Journal
        (Public Domain image from
        Wikipedia).}\label{LiHAFOUMSldRGIRA4AXELHARLROA7789GIJ7sLIOU1polFUMHDNOJHNFOJEDSOL1}
\end{figure}

A classical result known with the name\footnote{Liouville's contributions to science are variegate and cover a large number of topics such as electrodynamics, potential theory, spectral theory, theory of heat, fractional calculus, number theory, rational mechanics and differential geometry, and different classical results go nowadays under the name of ``Liouville's Theorem''.
The one related to harmonic functions (and complex analysis) was probably found around 1844 by Liouville while working on doubly periodic meromorphic functions (nowadays called ``elliptic functions'' in the complex analysis jargon) and its initial proof relied on trigonometric series (later on, he also proposed several different approaches).

The priority about Liouville's Theorem is however controversial. In particular, Cauchy, eager to secure his priority in the field of complex analysis, called attention to some of his results from 1843-1844 stating that ``if a function~$f(z)$ of the real or imaginary variable~$z$ always remain continuous, and consequently always finite, it is reduced to a simple constant''. In hindsight, Cauchy's statement, though perhaps foregoing or more general than Liouville's one, would be considered today as poorly stated, both because of the ambiguity between ``finite'' and ``bounded'' and, more importantly, because at this stage Cauchy seems to erroneously believe that continuity is equivalent to holomorphy ({\em aliquando bonus dormitat Homerus}: this mistake seems to have been repeated by Cauchy on other occasions as well).

In any case, the priority debate on Liouville's Theorem is still unsettled (but, after all, who cares?). See~\cite{10230723632205, MR1066463} for more information about Liouville's biography and mathematics, and about the history of Liouville's Theorem.
See also Figure~\ref{LiPORLROA7789GIJ7sLIOU1polFUMHDNOJHNFOJEDSOL1}
for a portrait of Liouville and Figure~\ref{HAFOUMSLi:bblFUMHDNOJHNFOJED2423686900-4} for a portrait of
Cauchy
(next to a rather obscure figure, Libri;
see footnote~\ref{PALISOJA0987654AIKJ}
on page~\pageref{PALISOJA0987654AIKJ} for some crazy stories also about this guy, Libri; here it suffices to say that Liouville discovered some plagiarism, and errors, by Libri, who at some point
became Liouville's arch-enemy, see~\cite{MR1066463} for more details about these events).

Interestingly, the first application Liouville made of his/Cauchy's result was an elegant proof of the Fundamental Theorem of Algebra which is still standard in today's textbooks; yet, Liouville had left this proof unpublished (this proof can be read from Liouville's manuscripts available at the Biblioth\`eque de l'Institut de France, Ms 3617(5), page 85).

We also recall that in 1836 Liouville founded a prominent scientific journal of mathematics named ``Journal de Math{\'e}matiques Pures et Appliqu{\'e}es''.
The journal is still published on a monthly basis and it is often nicknamed ``Liouville's Journal'' after its creator.
See Figure~\ref{LiHAFOUMSldRGIRA4AXELHARLROA7789GIJ7sLIOU1polFUMHDNOJHNFOJEDSOL1}
for the title page of the journal's first volume.}
of Liouville's Theorem states that
bounded harmonic functions in the whole of~$\R^n$ are necessarily constant:
as a matter of fact, to obtain such a result, a one-sided bound (i.e., from either 
above or below) is sufficient, as we discuss in the next result.\index{Liouville's Theorem}

\begin{theorem}\label{LIOVE}
If~$u$ is harmonic in~$\R^n$ and either bounded from above or from below, then~$u$ is constant.
\end{theorem}

\begin{proof}
Up to replacing~$u$ with~$-u$, we can assume that~$u$ is bounded below.
Then, up to replacing~$u$ with~$u-\inf_{\R^n}u$, we can suppose that
\begin{equation}\label{SUjhncbtrafr4a5cdfgtrP}
{\mbox{$u\ge0$ in $\R^n$.}}
\end{equation}
We now exploit the
first derivative Cauchy's Estimate in Theorem~\ref{CAUESTIMTH}. Namely,
given~$x_0\in\R^n$ and~$r>0$, we have that
$$ \left|\nabla u(x_0)
\right|\le \frac{C\,\|u\|_{L^1(B_r(x_0))}}{r^{n+1}},$$
for some~$C>0$ depending only on~$n$.
This, \eqref{SUjhncbtrafr4a5cdfgtrP} and the
Mean Value Formula in Theorem~\ref{KAHAR}(iii)
lead to
$$ \left|\nabla u(x_0)
\right|\le \frac{C}{r^{n+1}}\,
\int_{B_r(x_0)}u(x)\,dx=
\frac{C\,|B_1|}{r}\,
\fint_{B_r(x_0)}u(x)\,dx=\frac{C\,|B_1|\,u(x_0)}{r}.
$$
Sending~$r\to+\infty$, we thereby find that~$\nabla u(x_0)=0$.
Since this is valid for all~$x_0\in\R^n$, we conclude that~$u$ is constant, as desired.
\end{proof}

Many different proofs of Theorem~\ref{LIOVE} can be performed.
For example, here is one using
the Harnack Inequality instead of Cauchy's Estimate:

\begin{proof}[Another proof of Theorem~\ref{LIOVE}] Thanks to~\eqref{SUjhncbtrafr4a5cdfgtrP},
one can employ Theorem~\ref{PRE:HA} and take the limit as~$R\to+\infty$ in~\eqref{0-PREJA}. This,
for every~$x\in\R^n$, leads to

\begin{eqnarray*} u(0)=\lim_{R\to+\infty}\left(\frac{R}{R+r}\right)^{n-2}\,
\frac{R-r}{R+r}\,u(0)\le
u(x)\le\lim_{R\to+\infty}\left(\frac{R}{R-r}\right)^{n-2}\,\frac{R+r}{R-r}\,u(0)=u(0)\end{eqnarray*}
and therefore~$u(x)=u(0)$.\end{proof}

A variant of Theorem~\ref{LIOVE}
goes as follows:

\begin{theorem}\label{LIOCOT6TROL}
Let~$u$ be harmonic in~$\R^n$ and suppose that
$$ |u(x)| \le M\,(|x|+1)^k$$
for some~$M\ge0$ and~$k\in\N$. 

Then, $u$ is a polynomial of degree at most~$k$.
\end{theorem}

\begin{proof} Let~$x_0\in\R^n$ and pick~$r>|x_0|+1$. We notice that
$$ \|u\|_{L^1(B_r(x_0))}\le
M\int_{B_r(x_0)} (|x|+1)^k\,dx
\le M\int_{B_{2r}} (|x|+1)^k\,dx
\le C r^{k+n},$$
for some~$C>0$ depending only on~$M$, $k$ and~$n$.

Therefore, by the
Cauchy's Estimate in Theorem~\ref{CAUESTIMTH}, if~$\alpha\in \N^n$ with~$|\alpha|=k+1$,
$$ \left|
\frac{\partial^\alpha u}{\partial x^\alpha}(x_0)
\right|\le \frac{C_{k+1}\,\|u\|_{L^1(B_r(x_0))}}{r^{n+k+1}}\le\frac{C_{k+1}\,C}{r}.
$$
Hence, we can send~$r\to+\infty$
and deduce that~$\frac{\partial^\alpha u}{\partial x^\alpha}(x_0)=0$.

Since this is valid for all~$x_0\in\R^n$, we conclude that the derivatives of order~$k+1$
of~$u$ are identically zero. As a result, we use a Taylor expansion near the origin to see that
$$ u(x)=\sum_{|\alpha|\le k}\frac{\partial^\alpha u(0)}{k!}\,x^\alpha.$$
In particular, $u$ coincides with a polynomial of degree at most~$k$
near the origin. Since~$u$ is real analytic (recall Theorem~\ref{KMS:HAN0-0}),
making use of the Unique Continuation Principle
for analytic functions (see e.g.~\cite[page~67]{MR1916029})
we conclude that~$u$
coincides with this
polynomial in the whole of~$\R^n$.
\end{proof}

\section{Harmonic polynomials and spherical harmonics}\label{SLLD:SPJEDHHA}

Harmonic polynomials are polynomials~$p:\R^n\to\R$ which are harmonic,\index{harmonic polynomial}
i.e.~$\Delta p=0$ in~$\R^n$. 
The restrictions
of homogeneous\footnote{As customary,
a polynomial is called homogeneous of degree~$k$
if all its nonzero terms have the same degree~$k$, namely if it is of the form
$$ P(x_1,\dots,x_n)=\sum_{{a_1,\dots,a_n\in\N}\atop{a_1+\dots+a_n=k}}
c_{a_1,\dots,a_n}\,x_1^{a_1}\dots x_n^{a_n},$$
for some coefficients~$c_{a_1,\dots,a_n}\in\R$. We stress that, with
this notation, zero is a homogeneous polynomial of any degree (it suffices
to take all the coefficients to be zero). This setting \label{Le3AGOZUE0035}
is convenient since it makes
the family of homogeneous polynomials of degree $k$ a vector space.}
harmonic polynomials to the sphere~$\partial B_1$
are called spherical harmonics\index{spherical harmonic} (from now on,
in this section, we will suppose~$n\ge2$ to avoid the
trivial case of~$n=1$ in which~$\partial B_1$ consists of two points).
More specifically,
the restrictions of homogeneous harmonic polynomials of degree~$k$
to the sphere~$\partial B_1$
are called spherical harmonics of degree~$k$.

We start with an approximation result by polynomials on the sphere:

\begin{lemma}\label{OLISEDE}
The set of polynomials is dense in~$L^2(\partial B_1)$.
\end{lemma}

\begin{proof} Let~$f\in L^2(\partial B_1)$ and~$\e>0$.
We can find a continuous function~$g$ such that~$\|f-g\|_{L^2(\partial B_1)}<\e$
(indeed, up to a partition of unity, we can work in a local chart,
and reduce the approximation problem to that in a domain of~$\R^{n-1}$,
where standard convolution methods, see e.g.~\cite[Theorem 9.6]{MR3381284},
can be applied). 
Then, let~$\psi\in C^\infty_0( B_4\setminus B_{1/4})$ with~$\psi=1$ in~$B_2\setminus B_{1/2}$
and~$h(x):=\psi(x)\,g\left(\frac{x}{|x|}\right)$. We observe that~$h$ is a continuous function on~$\R^n$.
Hence, by the
Stone-Weierstra{\ss} Theorem
(see e.g.~\cite[Lemma~2.1]{MR3626547}), there exists a polynomial~$p$ such that~$\|h-p\|_{L^\infty(B_4)}<\e$.
As a result,
\begin{equation*}\begin{split}&
\|f-p\|_{L^2(\partial B_1)}\le\|f-g\|_{L^2(\partial B_1)}+\|g-p\|_{L^2(\partial B_1)}=
\|f-g\|_{L^2(\partial B_1)}+\|h-p\|_{L^2(\partial B_1)}\\&\qquad<\e+\sqrt{\int_{\partial B_1} |h(\omega)-p(\omega)|^2\,d{\mathcal{H}}^{n-1}_\omega}<\left(1+\sqrt{{\mathcal{H}}^{n-1}(\partial B_1) }\right)\e.
\qedhere\end{split}
\end{equation*}\end{proof}

Now we show that polynomials on the sphere can be harmonically
and polinomially extended to the whole of~$\R^n$.

\begin{lemma}\label{L2efsdfpolti6}
Let~$p$ be a polynomial of degree at most~$k$. Then, there exists a unique harmonic function~$h$
in~$\R^n$ such that~$h=p$ on~$\partial B_1$. Moreover, $h$ is a polynomial of degree at most~$k$.

More specifically, $h=p-(1-|x|^2)\widetilde{p}$,
for a suitable polynomial~$\widetilde{p}$ of degree at most~$k-2$.
\end{lemma}

\begin{proof} For every polynomial~$q$ of degree at most~$k$, we define~$Tq(x)=\Delta\big((1-|x|^2)q(x)\big)$
and we remark that~$Tq$ is also a polynomial of degree at most~$k$. Notice that~$T$ is linear.

Furthermore, $T$ is injective: to check this, suppose that~$Tq(x)=0$ for every~$x\in\R^n$,
and set~$\widetilde{q}(x):=(1-|x|^2)q(x)$. Note that~$\widetilde{q}$ is harmonic in~$B_1$
and vanishes along~$\partial B_1$, whence, by the Maximum Principle in
Corollary~\ref{WEAKMAXPLE}(iii),
necessarily~$\widetilde{q}$ vanishes identically in~$B_1$. This entails that~$q$
vanishes identically in~$B_1$, and also in~$\R^n$ by the
Unique Continuation Principle
for analytic functions (see e.g.~\cite[page~67]{MR1916029}).

Since the space of polynomials of degree at most~$k$ is finite dimensional, we thus
deduce that~$T$ is also surjective. Hence, taking~$h:=p-(1-|x|^2)T^{-1}(\Delta p)$,
we have that~$h$ is a polynomial of degree at most~$k$, that~$h=p$ along~$\partial B_1$,
and that
$$\Delta h=\Delta p-\Delta\big((1-|x|^2)T^{-1}(\Delta p)\big)
=\Delta p-T\big(T^{-1}(\Delta p)\big)=0.$$
Accordingly, $h$ satisfies the desired claim.

Suppose now that there exists another harmonic function~$f$ such that~$f=p$ on~$\partial B_1$
and let~$u:=f-h$. We have that~$\Delta u=0$ in~$B_1$ and~$u=0$
along~$\partial B_1$.
Accordingly, by the Maximum Principle in Corollary~\ref{WEAKMAXPLE}(iii)
we find that~$u$ vanishes identically and thus~$f=h$, which shows the uniqueness of~$h$.
\end{proof}

In the following Corollaries~\ref{DENSPO2} and~\ref{DENSPO}
we deduce some useful consequence of
Lemma~\ref{L2efsdfpolti6}:

\begin{corollary}\label{DENSPO2}
Let~$p$ be a homogeneous polynomial of degree~$k\ge2$.
Then, there exist a harmonic
homogeneous polynomial~$h$ of degree~$k$
and a homogeneous polynomial~$q$ of degree~$k-2$
such that
\begin{equation}\label{FAVscoocc8qiwfhgia2}p =h +|x|^2q .\end{equation}
Also, this representation is unique.
\end{corollary}

\begin{proof}
By Lemma~\ref{L2efsdfpolti6}, we can write~$p =H -(1-|x|^2)\widetilde{p} $,
where~$H$ is a harmonic polynomial of degree at most~$k$,
$H=p$ on~$\partial B_1$ and~$\widetilde{p} $ is a
polynomial of degree at most~$k-2$.
We let~$h$ to be the homogeneous part of~$H$ of degree~$k$, namely we write~$H=h+\widetilde{h}$,
where~$h$ is a harmonic
homogeneous polynomial of degree~$k$
and~$\widetilde{h}$ is a polynomial of degree at most~$k-1$.

We also take~$q$ to be
the homogeneous part of~$\widetilde{p}$ of degree~$k-2$, namely we write~$
\widetilde{p}=q+\widetilde{q}$,
where~$q$ is a homogeneous polynomial of degree~$k-2$
and~$\widetilde{q}$ is a polynomial of degree at most~$k-3$.
As a side remark, we stress that~$p$, $h$ and~$q$ are also allowed
to be zero in our setting (recall the footnote about
homogeneous polynomials on page~\pageref{Le3AGOZUE0035}).

As a result,
\begin{equation} \label{FAVscoocc8qiwfhgia} p=
h+\widetilde{h}- (1-|x|^2)(q+\widetilde{q})=
h+|x|^2q+P,
\end{equation}
with~$P:=\widetilde{h}-q-\widetilde{q}+|x|^2\widetilde{q}$.
Notice that~$P$ is a polynomial of degree at most~$k-1$.
On the other hand, by construction, $p-h-|x|^2q$
is 
a homogeneous polynomial of degree~$k$.
These observations and~\eqref{FAVscoocc8qiwfhgia} yield that~$P$
is necessarily zero. {F}rom this and~\eqref{FAVscoocc8qiwfhgia}
we obtain the desired claim in~\eqref{FAVscoocc8qiwfhgia2}.

To complete the desired result, it remains to check
that the representation in~\eqref{FAVscoocc8qiwfhgia2}
is unique. For this, by taking the difference of two representations,
it suffices to prove that
\begin{equation}\label{SMDC:9eorfkjfjf}
\begin{split}&
{\mbox{if a harmonic function~$\overline{h}$
and a homogeneous polynomial~$\overline{q}$ of degree~$k-2$
are such that}}
\\&{\mbox{$\overline{h}(x)+|x|^2\overline{q}(x)=0$ for all~$x\in\R^n$,
then both~$\overline{h}$ and~$\overline{q}$ are identically zero.}}
\end{split}\end{equation}
As a matter of fact, to prove this, 
it is enough to show that
\begin{equation}\label{SMdsdDC:9eorfkjfjf-2}
\begin{split}&
{\mbox{if a harmonic function~$\overline{h}$
and a homogeneous polynomial~$\overline{q}$ of degree~$k-2$
are such that}}
\\&{\mbox{$\overline{h}(x)+|x|^2\overline{q}(x)=0$ for all~$x\in\R^n$,
then for every~$j\in\N\cap\left[1,\frac{k}2\right]$}}\\&{\mbox{there exists
a homogeneous polynomial~$q_j$ of degree~$k-2j$
such that}}\\&{\mbox{$\overline{h}(x)+|x|^{2j}q_j(x)=0$ for all~$x\in\R^n$.}}
\end{split}\end{equation}
Indeed, if~\eqref{SMdsdDC:9eorfkjfjf-2} holds true,
then we can choose
$$j_\star:=\begin{dcases}
\frac{k}{2} & {\mbox{ if $k$ is even,}}\\
\frac{k-1}{2} & {\mbox{ if $k$ is odd,}}
\end{dcases}$$
and notice that~$k-2j_\star\in\{ 0,1\}$. Accordingly, the degree of
the polynomial~$q_{j_\star}$ is either zero or one.
Then, we have that~$q_{j_\star}(x)=a\cdot x+b$ for some~$a\in\R^n$ and~$b\in\R$
and, as a result,
\begin{equation}\label{2456KSJND:0eorifgk9568} \overline{h}(x)=-|x|^{2j_\star}
q_{j_\star}(x)=-|x|^{2j_\star}(a\cdot x+b).\end{equation}
{F}rom this we find that
\begin{eqnarray*}
&& 0=\Delta\overline{h}=-\Delta\Big( -|x|^{2j_\star}
(a\cdot x+b)\Big)=
2j_\star(2j_\star+n-2)|x|^{2j_\star-2}(a\cdot x+b) +
2j_\star |x|^{2j_\star-2}x\cdot a\\
&&\qquad\qquad=2j_\star\big(2j_\star+n-1\big)|x|^{2j_\star-2}
a\cdot x +2j_\star(2j_\star+n-2)b|x|^{2j_\star-2}
.
\end{eqnarray*}
This gives that~$a=0$ and~$b=0$, which in turn, combined with~\eqref{2456KSJND:0eorifgk9568}, yields that~$\overline{h}$
vanishes identically. {F}rom this and the fact that~$|x|^2\overline{q}=
-\overline{h}$, we conclude that~$\overline{q}$
vanishes identically as well, thus showing the validity
of~\eqref{SMDC:9eorfkjfjf}.

In view of these observations, to show the validity of~\eqref{SMDC:9eorfkjfjf},
and thus complete the proof of Corollary~\ref{DENSPO2},
it only remains to check~\eqref{SMdsdDC:9eorfkjfjf-2}.

The proof of~\eqref{SMdsdDC:9eorfkjfjf-2} is by induction over~$j$.
When~$j=1$, one takes~$q_j:=\overline{q}$.
Then, suppose that
there exists
a homogeneous polynomial~$q_j$ of degree~$k-2j$
with~$j\le\frac{k-2}2$ and
such that~$\overline{h}(x)+|x|^{2j}q_j(x)=0$ for all~$x\in\R^n$.
We observe that
\begin{equation}\label{KM:okdf928jjfjfHSHDN}0=
\Delta\big( \overline{h}+|x|^{2j}q_j\big)=\Delta\big(|x|^{2j}q_j\big)=
2j(2j+n-2)|x|^{2j-2}q_j +|x|^{2j}\Delta q_j+
4j |x|^{2j-2}x\cdot\nabla q_j.\end{equation}
Also, by homogeneity,
$$ x\cdot\nabla q_j(x)=\left.\frac{d}{dt} q_j(tx)\right|_{t=1}=\left.
\frac{d}{dt} \Big( t^{k-2j}q_j(x)\Big)\right|_{t=1}
=(k-2j) q_j(x),$$
and thus, from~\eqref{KM:okdf928jjfjfHSHDN},
\begin{eqnarray*}0&=&
2j(2j+n-2)|x|^{2j-2}q_j +|x|^{2j}\Delta q_j+
4j |x|^{2j-2}(k-2j) q_j\\&=&
2j(n-2+2k-2j)|x|^{2j-2}q_j +|x|^{2j}\Delta q_j.
\end{eqnarray*}
Therefore, since~$n-2+2k-2j\ge n-2+2k-(k-2)
>0$,
we can define
$$ q_{j+1}:=-\frac{\Delta q_j}{2j(n-2+2k-2j)}$$
and find that
$$ 0=q_j +\frac{|x|^{2}\Delta q_j}{2j(n-2+2k-2j)}=
q_j-|x|^2q_{j+1}.$$
In this way,
$$ \overline{h}+|x|^{2(j+1)}q_{j+1}=
\overline{h}+|x|^{2j}q_j+|x|^{2j}(|x|^2q_{j+1}-q_j)=\overline{h}+|x|^{2j}q_j=0,$$
which shows the validity of~\eqref{SMdsdDC:9eorfkjfjf-2}, as desired.
\end{proof}

\begin{corollary}\label{DENSPO}
Harmonic polynomials are dense in~$L^2(\partial B_1)$.
\end{corollary}

\begin{proof} Let~$f\in L^2(\partial B_1)$ and~$\e>0$.
By Lemma~\ref{OLISEDE}, we can find a polynomial~$p$ such that~$\|f-p\|_{L^2(\partial B_1)}<\e$.
Also, by Lemma~\ref{L2efsdfpolti6}, there exists a harmonic polynomial~$h$ such that~$\|p-h\|_{L^2(\partial B_1)}<\e$.
As a byproduct, $\|f-h\|_{L^2(\partial B_1)}<2\e$.
\end{proof}

It is also useful to remark that spherical harmonics of different degrees are orthogonal.
For this, we first observe that
there is a strong relation between harmonic polynomials and the eigenfunctions
of the Laplace-Beltrami operator\index{eigenfunction
of the Laplace-Beltrami operator}:

\begin{lemma}\label{8yiuhgdshh:LE}
Spherical harmonics of degree~$k$
are precisely the eigenfunctions of the Laplace-Beltrami operator on~$\partial B_1$
with\footnote{We will see in Theorem~\ref{8yiuhgdshh:CO}
that the eigenvalues listed in Lemma~\ref{8yiuhgdshh:LE}
are indeed all the eigenvalues of the Laplace-Beltrami operator on the sphere.} eigenvalue~$-k( n+k-2 )$.

That is, let~$h\in C^\infty(\partial B_1)$. Then, the following conditions are equivalent:
\begin{itemize}
\item[(i).] $h$ is a spherical harmonic of degree~$k$.
\item[(ii).] $h$
is an eigenfunction of the Laplace-Beltrami operator on the sphere with eigenvalue~$-k( n+k-2 )$.
\end{itemize}
\end{lemma}

\begin{proof} Let~$h:\partial B_1\to\R$ be a spherical harmonic
of degree~$k$. Then, there exists
a homogeneous harmonic polynomial~$p:\R^n\to\R$ of degree~$k$
such that~$h=p$ on~$\partial B_1$.
In particular, by~\eqref{BEL}, we have that
\begin{equation}\label{98jhnerLMpapakdrapa}
\Delta_{\partial B_1}h=\Delta_{\partial B_1} p.\end{equation}
This and Theorem~\ref{SPH} yield that, on~$\partial B_1$,
$$ 0=\Delta p=\partial_{rr}p_0+({n-1})\,\partial_{r}p_0+\Delta_{\partial B_1}h.$$
where
\begin{equation}\label{pioskd-2}
{\mbox{$p_0(r,\vartheta):=p(r\vartheta)$ for each~$r>0$ and~$\vartheta\in\partial B_1$.}}\end{equation}
By homogeneity, we know that~$p_0(r,\vartheta)=r^k p(\vartheta)$,
therefore
\begin{equation}\label{pioskd-3}\begin{split}&
\partial_{r}p_0(r,\vartheta)=kr^{k-1}p(\vartheta)=kr^{k-1}h(\vartheta)\\
{\mbox{and}}\qquad&\partial_{rr}p_0(r,\vartheta)=k(k-1)r^{k-2}p(\vartheta)=k(k-1)r^{k-2}h(\vartheta).\end{split}\end{equation}
Consequently, on~$\partial B_1$,
$$ 0=k(k-1)h+(n-1)kh+\Delta_{\partial B_1}h=k( n+k-2 )h+\Delta_{\partial B_1}h
.$$
Accordingly, we have that~$h$ is an
eigenfunction of the Laplace-Beltrami operator with eigenvalue~$-k( n+k-2 )$.

Viceversa, suppose now that~$h$ is an
eigenfunction of the Laplace-Beltrami operator on~$\partial B_1$ \index{eigenvalue problem}
with eigenvalue~$-k( n+k-2 )$. Let~$p(x):=|x|^k h\left(\frac{x}{|x|}\right)$.
Let~$p_0$ be as in~\eqref{pioskd-2} and notice that~$p$ is homogeneous of degree~$k$. Thus,
exploiting Theorem~\ref{SPH}, \eqref{98jhnerLMpapakdrapa} and~\eqref{pioskd-3}, we see that, on~$\partial B_1$,
$$ \Delta p =\partial_{rr}p_0+(n-1)\,\partial_{r}p_0+\,\Delta_{\partial B_1}h
=k(k-1)r^{k-2}h+(n-2)kr^{k-2}h- k(n+k-2)h=0.$$
Since, by homogeneity, for every~$\lambda>0$ and~$x\in\R^n$,
$$ \frac{\partial^2p}{\partial x_i\partial x_j}(\lambda x)=\frac1{\lambda^2}\,
\frac{\partial^2}{\partial x_i\partial x_j}( p(\lambda x) )=\frac1{\lambda^2}\,\frac{\partial^2}{\partial x_i\partial x_j}( \lambda^k p( x) )
 =\lambda^{k-2}\frac{\partial^2p}{\partial x_i\partial x_j}(x),
$$ we thus conclude that, for each~$x\in\R^n\setminus\{0\}$,
$$ \Delta p(x)=\Delta p\left( |x|\,\frac{x}{|x|}\right)=|x|^{k-2}\Delta p\left( \frac{x}{|x|}\right)=0,
$$
and so, by continuity, $\Delta p=0$ in~$\R^n$.

In addition, for every~$x\in\R^n$,
$$ |p(x)|\le |x|^k\max_{\partial B_1}|h|,$$
hence, by the Liouville Theorem (recall Theorem~\ref{LIOCOT6TROL}), it follows that~$p$ is a polynomial of degree~$k$.

Since~$p$ coincides with~$h$ on~$\partial B_1$, this gives that~$h$
is a spherical harmonics of degree~$k$.
\end{proof}

With this, we are in the position of proving that spherical harmonics
are orthogonal:

\begin{lemma}\label{DIGGDJKDL}
Let~$h$ and~$g$ be spherical harmonics
of different degree. Then,
$$ \int_{\partial B_1} h(x)\,g(x)\,d{\mathcal{H}}^{n-1}_x=0.$$
\end{lemma}

\begin{proof} Let~$k$ and~$m$ be the degrees of~$h$ and~$g$, respectively.
Let also~$\phi(t):=t( t+n-2 )$. Since~$\phi'(t)=2t+n-2\ge n>0$ for all~$t\ge1$,
it follows that~$\phi$ is strictly increasing in~$[1,+\infty)$ and in particular, if~$k\ne m$,
then
\begin{equation}\label{664f3298456-3495}
k( n+k-2 )\ne m( m+n-2 ).\end{equation}
Moreover, by Lemma~\ref{8yiuhgdshh:LE} and the self-adjointness property pointed out in~\eqref{PAR-092},
\begin{eqnarray*}&&-k( n+k-2 )\int_{\partial B_1} h(x)\,g(x)\,d{\mathcal{H}}^{n-1}_x=
\int_{\partial B_1} \Delta_{\partial B_1}h(x)\,g(x)\,d{\mathcal{H}}^{n-1}_x\\&&\qquad=
\int_{\partial B_1}h(x)\, \Delta_{\partial B_1}g(x)\,d{\mathcal{H}}^{n-1}_x=-m( m+n-2 )\int_{\partial B_1} h(x)\,g(x)\,d{\mathcal{H}}^{n-1}_x
.\end{eqnarray*}
This and~\eqref{664f3298456-3495} entail the desired result.
\end{proof}

With this, we can actually strengthen Corollary~\ref{DENSPO}
by approximating any function in~$L^2(\partial B_1)$
using harmonic polynomials of increasing degree
and modifying each step of the approximation by a polynomial
of higher degree:

\begin{corollary} For every~$f\in L^2(\partial B_1)$
there exists a sequence~$h_k$ of spherical harmonics of degree
equal to~$k$ such that
$$ \lim_{N\to+\infty}\left\|f-\sum_{k=0}^N h_k 
\right\|_{L^2(\partial B_1)}=0.$$
\end{corollary}

\begin{proof} The proof exploits some classical
methods from functional analysis and Fourier Series. We give
a self-contained argument for the facility of the reader.
For every~$k\in\N\setminus\{0\}$ we consider the vector
space of spherical harmonics of degree
precisely equal to~$k$. This space
is finite dimensional (since so is the space of polynomials)
and therefore we can consider an orthonormal basis (with respect to the scalar product in~$L^2(\partial B_1)$)
that we will denoted by~$\{\eta_{1,k},\dots,\eta_{n_k,k}\}$.
By Lemma~\ref{DIGGDJKDL}, the set~$\{\eta_{j,k}\}_{{k\in\N}\atop{1\le j\le n_k}}$
is an orthonormal family, that is
\begin{equation}\label{iksdmcdelta} \langle
\eta_{j,k} 
,\eta_{i,m}
\rangle_{L^2(\partial\Omega)}= \delta_{km}\delta_{ij}
.\end{equation}
Let
$$ h_k:=\sum_{1\le j\le n_k} \langle
f,\eta_{j,k}
\rangle_{L^2(\partial\Omega)}\,\eta_{j,k}$$
and notice that~$h_k$ is a spherical harmonic of degree
equal to~$k$. 

Now we take~$\e>0$ and employ Corollary~\ref{DENSPO}
to find a linear combination of spherical harmonics~$h^{(\e)}$ such that~$\|f
-h^{(\e)}\|_{L^2(\partial\Omega)}^2<\e$. We write
$$ h^{(\e)}=\sum_{k=0}^{N^{(\e)}}
h^{(\e)}_k,$$
for suitable spherical harmonics~$h^{(\e)}_k$ of degree precisely equal to~$k$
and some~${N^{(\e)}}\in\N$.
We also express each~$h^{(\e)}_k$ in terms of the above orthonormal
basis by writing
$$ h^{(\e)}_k=\sum_{1\le j\le n_k} a_{jk}^{(\e)}\,\eta_{j,k},$$
for suitable~$a_{jk}^{(\e)}\in\N$.
All in all, we find that
$$ h^{(\e)}=\sum_{{0\le k\le N^{(\e)}}\atop{1\le j\le n_k}}
a_{jk}^{(\e)}\,\eta_{j,k}$$
and accordingly
$$ h^{(\e)}\pm \sum_{k=0}^{N^{(\e)} } h_k
=\sum_{{0\le k\le N^{(\e)}}\atop{1\le j\le n_k}}\Big(
a_{jk}^{(\e)}\pm\langle
f,\eta_{j,k}
\rangle_{L^2(\partial\Omega)}\Big)\,\eta_{j,k}.$$
In light of~\eqref{iksdmcdelta}, we thereby have that
\begin{equation*}
\begin{split}&
\left\|f-\sum_{k=0}^{N^{(\e)} } h_k \right\|_{L^2(\partial\Omega)}^2-\e
\\&\qquad<
\left\|f-\sum_{k=0}^{N^{(\e)} }h_k \right\|_{L^2(\partial\Omega)}^2-
\|f-h^{(\e)}\|_{L^2(\partial\Omega)}^2\\&\qquad=
\|f\|_{L^2(\partial\Omega)}^2+
\left\langle \sum_{m=0}^{N^{(\e)} }h_m,\sum_{k=0}^{N^{(\e)} }h_k
\right\rangle_{L^2(\partial\Omega)}
-2\left\langle f,\sum_{k=0}^{N^{(\e)} }h_k
\right\rangle_{L^2(\partial\Omega)}\\&\qquad\qquad
-\|f\|_{L^2(\partial\Omega)}^2-\|h^{(\e)}\|_{L^2(\partial\Omega)}^2
+2\left\langle f,h^{(\e)}
\right\rangle_{L^2(\partial\Omega)}\\&\qquad=
\left\langle 
\sum_{{0\le m\le N^{(\e)}}\atop{1\le i\le n_m}} \langle
f,\eta_{i,m}
\rangle_{L^2(\partial\Omega)}\,\eta_{i,m}
,\,\sum_{{0\le k\le N^{(\e)}}\atop{1\le j\le n_k}} \langle
f,\eta_{j,k}
\rangle_{L^2(\partial\Omega)}\,\eta_{j,k}
\right\rangle_{L^2(\partial\Omega)}\\&\qquad\qquad
-2\left\langle f,\,\sum_{{0\le k\le N^{(\e)}}\atop{1\le j\le n_k}} \langle
f,\eta_{j,k}
\rangle_{L^2(\partial\Omega)}\,\eta_{j,k}
\right\rangle_{L^2(\partial\Omega)}\\&\qquad\qquad
-
\left\langle
\sum_{{0\le m\le N^{(\e)}}\atop{1\le i\le n_k}}
a_{im}^{(\e)}\,\eta_{i,m},\,
\sum_{{0\le k\le N^{(\e)}}\atop{1\le j\le n_k}}
a_{jk}^{(\e)}\,\eta_{j,k}
\right\rangle_{L^2(\partial\Omega)}
+2\left\langle f,\,\sum_{{0\le k\le N^{(\e)}}\atop{1\le j\le n_k}}
a_{jk}^{(\e)}\,\eta_{j,k}
\right\rangle_{L^2(\partial\Omega)}\\&\qquad=
-\sum_{{0\le k\le N^{(\e)}}\atop{1\le j\le n_k}} \langle
f,\eta_{j,k}
\rangle_{L^2(\partial\Omega)}^2-
\sum_{{0\le k\le N^{(\e)}}\atop{1\le j\le n_k}}
(a_{jk}^{(\e)})^2
+2\sum_{{0\le k\le N^{(\e)}}\atop{1\le j\le n_k}}
a_{jk}^{(\e)}\,\langle
f,\eta_{j,k}
\rangle_{L^2(\partial\Omega)}\\&\qquad=-
\sum_{{0\le k\le N^{(\e)}}\atop{1\le j\le n_k}}\Big(
\langle
f,\eta_{j,k}
\rangle_{L^2(\partial\Omega)}-a_{jk}^{(\e)}
\Big)^2\\&\qquad\le0.\qedhere
\end{split}\end{equation*}
\end{proof}

We will now establish that the eigenvalues found
in Lemma~\ref{8yiuhgdshh:LE}
are indeed all the possible eigenvalues of the Laplace-Beltrami operator\index{eigenvalue of the Laplace-Beltrami operator}
on the sphere:

\begin{theorem}\label{8yiuhgdshh:CO}
The eigenvalues of the Laplace-Beltrami operator
on~$\partial B_1$ are~$-k( n+k-2 )$ with~$k\in\N\setminus\{0\}$.
\end{theorem}

\begin{proof} We recall that we are assuming~$n\ge2$ and we use the
notation~$x=(x_1,\dots,x_n)$ and~$z=x_1+ix_2\in\cOMPL$.
By~\eqref{REj},
we know that, for every~$k\in\N$, the following function is harmonic
in~$\R^2$ (hence in~$\R^n$ by trivial extension in the variables~$(x_3,\dots,x_n)$):
$$ p_k(x_1,x_2):= \Re z^k=\Re(x_1+ix_2)^k=\Re\left(
\sum_{m=0}^k\left({{k}\atop{m}}\right) i^m x_1^{k-m} x_2^{m} 
\right)=
\sum_{{0\le m\le k}\atop{m\in2\Z}}\left({{k}\atop{m}}\right) (-1)^{\frac{m}2} x_1^{k-m} x_2^{m}
.$$
Notice that~$p_k$ is a homogeneous polynomial of degree~$k$
and thus, in light of Lemma~\ref{8yiuhgdshh:LE}, the corresponding spherical harmonic
defined by~$h_k(x):=p_k(x)$ for all~$x\in\partial B_1$
is an eigenfunction of the Laplace-Beltrami operator on~$\partial B_1$
with eigenvalue~$-k( n+k-2 )$.

This shows that the set of eigenvalues of the Laplace-Beltrami operator on~$\partial B_1$
contains~$\big\{-k( n+k-2 ),$ $k\in\N\setminus\{0\}\big\}$.

Suppose now that~$\lambda$ is an
eigenvalue of the Laplace-Beltrami operator on~$\partial B_1$.
We aim at showing that~$\lambda\in
\big\{-k( n+k-2 ),$ $k\in\N\setminus\{0\}\big\}$.
To this end, we argue by contradiction  and suppose that
\begin{equation}\label{dicj9qhdiwbiamo} 
\lambda\not\in
\big\{-k( n+k-2 ),\;k\in\N\setminus\{0\}\big\}.
\end{equation}
For all~$k\in\N\setminus\{0\}$ we observe that the set
of spherical harmonics of degree~$k$ is a vector space,
and it is finite dimensional (since so is the space of polynomials).
Thus we consider an orthonormal basis for such a space
(with respect to the scalar product in~$L^2(\partial B_1)$),
denoted by~$\{E_{1,k},\dots,E_{n_k,k}\}$
and the projection
\begin{equation}\label{9ijn avadindcvv} \Pi_k \varphi(x):=\sum_{j=1}^{n_k}\int_{\partial B_1}\varphi(y)\,
E_{j,k}(y)\,d{\mathcal{H}}^{n-1}_y\,E_{j,k}(x)=
\sum_{j=1}^{n_k}\langle\varphi,
E_{j,k}\rangle_{L^2(\partial B_1)}\,E_{j,k}(x).\end{equation}
Furthermore, by Lemma~\ref{8yiuhgdshh:LE}
and the self-adjointness property in~\eqref{PAR-092},
if~$\varphi$ is an eigenfunction corresponding to~$\lambda$, say normalized such that
\begin{equation}\label{AJNClitgsze}
\|\varphi\|_{L^2(\partial B_1)}=1
,\end{equation}
then
\begin{eqnarray*} &&-k( n+k-2 )\langle\varphi,
E_{j,k}\rangle_{L^2(\partial B_1)}=-
\int_{\partial B_1}\varphi(y)\,k( n+k-2 )
E_{j,k}(y)\,d{\mathcal{H}}^{n-1}_y\\&&\qquad=
\int_{\partial B_1}\varphi(y)\,\Delta_{\partial B_1}
E_{j,k}(y)\,d{\mathcal{H}}^{n-1}_y=
\int_{\partial B_1}\Delta_{\partial B_1}\varphi(y)\,
E_{j,k}(y)\,d{\mathcal{H}}^{n-1}_y\\&&\qquad=\lambda
\int_{\partial B_1}\varphi(y)\,
E_{j,k}(y)\,d{\mathcal{H}}^{n-1}_y=\lambda\langle\varphi,
E_{j,k}\rangle_{L^2(\partial B_1)}
\end{eqnarray*}
and accordingly
$$ \left(\frac{\lambda}{k( n+k-2 )} +1\right)\,
\langle\varphi,
E_{j,k}\rangle_{L^2(\partial B_1)}=0.
$$
With this information and~\eqref{dicj9qhdiwbiamo}, we deduce that~$\langle\varphi,
E_{j,k}\rangle_{L^2(\partial B_1)}=0$ for each~$j\in\{1,\dots,n_k\}$.
This and~\eqref{9ijn avadindcvv} yield that~$\Pi_k \varphi=0$
for every~$k\in\N\setminus\{0\}$.

Now, for each~$\e>0$, we use
Corollary~\ref{DENSPO} to pick a spherical harmonic~$h_\e$
such that~$\|\varphi-h_\e\|_{L^2(\partial B_1)}<\e$.
We let~$k_\e$ be the degree of~$h_\e$ and we observe that
the projection of~$h_\e$ on the space of
spherical harmonics with degree~$k_\e$ is~$h_\e$ itself.
As a result, we have that
\begin{eqnarray*}&&\| h_\e\|_{L^2(\partial B_1)}=
\| h_\e-\Pi_{k_\e} \varphi\|_{L^2(\partial B_1)}
=\| \Pi_{k_\e}h_\e-\Pi_{k_\e} \varphi\|_{L^2(\partial B_1)}\\&&\qquad
=\| \Pi_{k_\e}(h_\e-\varphi)\|_{L^2(\partial B_1)}\le\| h_\e-\varphi\|_{L^2(\partial B_1)}<\e.
\end{eqnarray*}
Consequently,
$$ \|\varphi\|_{L^2(\partial B_1)}\le\|\varphi-h_\e\|_{L^2(\partial B_1)}+
\| h_\e\|_{L^2(\partial B_1)}<2\e.$$
Since~$\e$ is arbitrary, this entails that~$\|\varphi\|_{L^2(\partial B_1)}=0$,
in contradiction with the normalization in~\eqref{AJNClitgsze}.
\end{proof}

We can now compute the dimension of the vector spaces
of spherical harmonics:

\begin{theorem} \label{DIMSPOM}
The space of harmonic homogeneous polynomials
of degree~$k$ has dimension\footnote{When~$k\ge1$, after a simple algebraic computation,
one can alternatively write the quantity in~\eqref{ALEHD:DJFJFJ} as
$$ \frac{n+2k-2}{k}\left({k+n-3}\atop{k-1}\right) .$$
}
\begin{equation}\label{ALEHD:DJFJFJ} \begin{dcases}
1 & {\mbox{ if }}k=0,\\
n& {\mbox{ if }}k=1,\\
\left({k+n-1}\atop{n-1}\right)-\left({k+n-3}\atop{n-1}\right) & {\mbox{ if }}k\ge2.
\end{dcases}\end{equation}
\end{theorem}

\begin{proof} We consider separately the cases~$k=0$,
$k=1$ and~$k\ge2$.

When~$k=0$,
the only homogeneous
polynomials are the constants, thus producing a one-dimensional
vector space. 

When~$k=1$,
we only have linear functions (which are the homogeneous
polynomials of degree one, and are also harmonic),
thus producing an $n$-dimensional
vector space.

Thus, we can focus now on the proof of~\eqref{ALEHD:DJFJFJ}
in the case~$k\ge2$.
We denote by~$D_k$ the dimension of the
space of the homogeneous polynomials in~$\R^n$ of degree~$k$.
Let also~$D_k^\star$ the dimension of the
space of the
harmonic homogeneous polynomials in~$\R^n$ of degree~$k$.

In view of~\eqref{FAVscoocc8qiwfhgia2}, we know that~$D_k=D_k^\star+D_{k-2}$,
and therefore
\begin{equation}\label{56:09erfikSFGFDke6ay}\begin{split}&
{\mbox{the dimension~$D_k^\star$
of the space of harmonic homogeneous polynomials of degree~$k$}}\\&{\mbox{is equal to~$D_k-D_{k-2}$.}}\end{split}
\end{equation}

Now we prove that
\begin{equation}\label{56:09erfikSFGFDke6ay2}
D_k=\left({k+n-1}\atop{n-1}\right).
\end{equation}
This is an algebraic (or, actually, combinatorial)
separate statement on polynomials (nothing to do with
harmonic functions).
To prove it, we observe that
for every~$a_1,\dots,a_n\in\N$ with~$a_1+\dots+a_n=k$,
the homogeneous polynomials, or actually monomials, of the form~$
x_1^{a_1}\dots x_n^{a_n}$ are all linearly independent, due to the Identity Principle for polynomials.
Since every homogeneous polynomial is written as the sum of such monomials,
we have that
\begin{equation}D_k=\#\Big\{
a_1,\dots,a_n\in\N {\mbox{ with }}a_1+\dots+a_n=k
\Big\}.\end{equation}
Now we claim that
\begin{equation}\label{ILVCLKVVLA89CLA908ytCKAla089-434-1}
\#\Big\{
a_1,\dots,a_n\in\N {\mbox{ with }}a_1+\dots+a_n=k
\Big\}=
\left({k+n-1}\atop{n-1}\right).
\end{equation}
To prove this, 
we define
\begin{equation}\label{ILVCLKVVLA89CLA908ytCKAla089-434-2} {\mathcal{A}}_{n,k}:=
\Big\{
(a_1,\dots,a_n)\in\N^n {\mbox{ with }}a_1+\dots+a_n=k
\Big\}.\end{equation}
We let~${\mathcal{T}}:{\mathcal{A}}_{n,k}\to {\mathcal{A}}_{n-1,k}\cup {\mathcal{A}}_{n,k-1}$ defined as
$$ {\mathcal{T}}(a_1,\dots,a_n):=\begin{dcases}
(a_1,\dots,a_{n-1}) & {\mbox{ if }}a_n=0,\\
(a_1,\dots,a_{n-1},a_n-1) & {\mbox{ if }}a_n\ne0.
\end{dcases}$$
Notice indeed that, for each~$(a_1,\dots,a_n)\in{\mathcal{A}}_{n,k}$,
if~$a_n=0$ then~${\mathcal{T}}(a_1,\dots,a_n)=(a_1,\dots,a_{n-1}) \in{\mathcal{A}}_{n-1,k}$,
while if~$a_n\ne0$ then~${\mathcal{T}}(a_1,\dots,a_n)=(a_1,\dots,a_{n-1},a_n-1) \in{\mathcal{A}}_{n,k-1}$.

Furthermore, 
\begin{equation}\label{T-0o-0o-oiuygfds9uyfd7ycf7trfg01}
{\mbox{${\mathcal{T}}$ is injective.}}\end{equation}
Indeed, let~$(a_1,\dots,a_n)$, $(b_1,\dots,b_n)\in{\mathcal{A}}_{n,k}$ and suppose that~$
{\mathcal{T}}(a_1,\dots,a_n)={\mathcal{T}}(b_1,\dots,b_n)=:\tau$. If~$\tau\in{\mathcal{A}}_{n-1,k}$,
then necessarily~$a_n=b_n=0$ and~$\tau=(a_1,\dots,a_{n-1})=(b_1,\dots,b_{n-1})$,
leading to~$(a_1,\dots,a_n)=(b_1,\dots,b_n)$.

If instead~$\tau\in{\mathcal{A}}_{n,k-1}$, then necessarily~$a_n\ne0$ and~$b_n\ne0$,
and moreover~$\tau=(a_1,\dots,a_{n-1},a_n-1)=(b_1,\dots,b_{n-1},b_n-1)$,
showing again that~$(a_1,\dots,a_n)=(b_1,\dots,b_n)$.
These observations complete the proof of~\eqref{T-0o-0o-oiuygfds9uyfd7ycf7trfg01}.

Additionally, we have that
\begin{equation}\label{T-0o-0o-oiuygfds9uyfd7ycf7trfg02}
{\mbox{${\mathcal{T}}$ is surjective.}}\end{equation}
Indeed, if~$(\vartheta_1,\dots,\vartheta_{n-1})\in{\mathcal{A}}_{n-1,k}$,
we let~$a_i:=\vartheta_i$ for each~$i\in\{1,\dots,n-1\}$ and~$a_n:=0$;
in this way, we see that~$(a_1,\dots,a_n)\in{\mathcal{A}}_{n,k}$
and~${\mathcal{T}}(a_1,\dots,a_n)=
(\vartheta_1,\dots,\vartheta_{n-1})$.
If instead~$(\vartheta_1,\dots,\vartheta_{n})\in{\mathcal{A}}_{n,k-1}$,
then we set~$a_i:=\vartheta_i$ for each~$i\in\{1,\dots,n-1\}$ and~$a_n:=\vartheta_n+1$;
in this way, we see that~$(a_1,\dots,a_n)\in{\mathcal{A}}_{n,k}$
and~${\mathcal{T}}(a_1,\dots,a_n)=
(\vartheta_1,\dots,\vartheta_{n})$.
These considerations prove~\eqref{T-0o-0o-oiuygfds9uyfd7ycf7trfg02}.

As a result, we deduce from~\eqref{T-0o-0o-oiuygfds9uyfd7ycf7trfg01}
and~\eqref{T-0o-0o-oiuygfds9uyfd7ycf7trfg02}
that
\begin{equation}\label{LIN798DGYTreMNAHER-dfV9ikmdKInmndG9oidf00}
\#{\mathcal{A}}_{n,k}=\#{\mathcal{A}}_{n-1,k}+\# {\mathcal{A}}_{n,k-1}.
\end{equation}

Now we claim that
\begin{equation}\label{LIN798DGYTreMNAHER-dfV9ikmdKInmndG9oidf0}
\#{\mathcal{A}}_{n,k}=
\left({k+n-1}\atop{n-1}\right).
\end{equation}
We prove this by induction over~$N:=k+n$.
When~$N=1$, we have that necessarily~$n=1$ and~$k=0$,
thus~${\mathcal{A}}_{n,k}={\mathcal{A}}_{1,0}$ reduces to the zero element, whence
$$ \#{\mathcal{A}}_{1,0}
=1=\left({0}\atop{0}\right).$$
As for the inductive step, we consider~$N=k+n\ge2$
and assume~\eqref{LIN798DGYTreMNAHER-dfV9ikmdKInmndG9oidf0}
to be true for the index~$N-1\ge1$. Thus, exploiting~\eqref{LIN798DGYTreMNAHER-dfV9ikmdKInmndG9oidf00},
\begin{eqnarray*}&&
\#{\mathcal{A}}_{n,k}=
\#{\mathcal{A}}_{n-1,k}+\# {\mathcal{A}}_{n,k-1}=
\left({k+n-2}\atop{n-2}\right)
+
\left({k+n-2}\atop{n-1}\right)\\&&\qquad\qquad\quad=
\frac{(k+n-2)!}{k!\,(n-2)!}+\frac{(k+n-2)!}{(k-1)!\,(n-1)!}=
\frac{(k+n-2)!}{k!\,(n-1)!}\big((n-1)+k
\big)\\&&\qquad\qquad\quad
=\frac{(k+n-1)!}{k!\,(n-1)!}
=\left({k+n-1}\atop{n-1}\right).
\end{eqnarray*}
This completes the inductive step and establishes~\eqref{LIN798DGYTreMNAHER-dfV9ikmdKInmndG9oidf0},
as desired.

The claim in~\eqref{ILVCLKVVLA89CLA908ytCKAla089-434-1}
now plainly follows from~\eqref{ILVCLKVVLA89CLA908ytCKAla089-434-2}
and~\eqref{LIN798DGYTreMNAHER-dfV9ikmdKInmndG9oidf0}.

Accordingly,
the desired claim in~\eqref{ALEHD:DJFJFJ} is a straightforward
consequence of~\eqref{56:09erfikSFGFDke6ay}
and~\eqref{ILVCLKVVLA89CLA908ytCKAla089-434-1}.\end{proof}

In view of
Theorem~\ref{DIMSPOM},
we now aim at constructing explicit bases
of the space of harmonic homogeneous polynomials
of degree~$k$ (this is of clear practical importance,
since one can reduce a number of general computations
to those related to a specific basis). This is indeed the content of
the forthcoming Section~\ref{sec:leg}.

\section{Spherical harmonics and Legendre polynomials}\label{sec:leg}

We now turn our attention to some special class of spherical harmonics with explicit and convenient
algebraic properties. To this end, 
we consider a basis for the space of harmonic homogeneous polynomials
of degree~$k$
and the corresponding spherical harmonics of degree~$k$.
After a Gram-Schmidt process, we find a maximal linearly independent set of spherical harmonics~$\{Y_{k,1},\dots,Y_{k,N_k}\}$
that are orthonormal with respect to the scalar product in~$L^2(\partial B_1)$
(here~$N_k$ is as in~\eqref{ALEHD:DJFJFJ},
thanks to Theorem~\ref{DIMSPOM}).
We stress that, in light of Lemma~\ref{DIGGDJKDL},
\begin{equation}\label{LORZRPO} \int_{\partial B_1} Y_{k,j}(x)\,Y_{m,i}(x)\,d{\mathcal{H}}^{n-1}_x= \delta_{km}\,\delta_{ji}.\end{equation}
For each~$x$, $y\in\partial B_1$, we also define
\begin{equation}\label{ASSDM} F_k(x,y):=\sum_{j=1}^{N_k} Y_{k,j}(x)\,Y_{k,j}(y)\end{equation}
and we point out a symmetry property of this function:

\begin{lemma}\label{LE21}
The function~$F_k$ is invariant under rotations. That is, if~${\mathcal{R}}$
is a rotation matrix, then
$$ F_k({\mathcal{R}}x,{\mathcal{R}}y)=F_k(x,y)
\qquad{\mbox{for every~$x$, $y\in\partial B_1$.}}$$
\end{lemma}

\begin{proof} If~$p$ is a harmonic homogeneous polynomial of degree~$k$, then
so is~$\widetilde{p}(x):=p({\mathcal{R}}x)$, thanks to the invariance under rotation of the Laplace equation
(recall Corollary~\ref{KSMD:ROTAGSZOKA}).
Accordingly, rotations preserve spherical harmonics and~$\widetilde{Y}_{k,j}(x):={Y}_{k,j}({\mathcal{R}}x)$
is also a spherical harmonic
of degree~$k$. Accordingly, we can write~$\widetilde{Y}_{k,j}$ in terms of the
orthonormal basis~$\{Y_{k,1},\dots,Y_{k,N_k}\}$ as
$$ \widetilde{Y}_{k,j}=\sum_{i=1}^{N_k} M_{ij}\, Y_{k,i},$$
for some~$M_{ij}\in\R$. We let~$M$ be
the $(N_k\times N_k)$ matrix produced by~$M_{ij}$ and we note that
\begin{equation}\label{LORZRPO2}
{\mbox{the matrix~$M$ is orthogonal,}}
\end{equation}
since, by~\eqref{LORZRPO},
\begin{eqnarray*}
&&\sum_{i=1}^{N_k} M_{ij}M_{im}=
\sum_{i,\ell=1}^{N_k} M_{ij}M_{\ell m}\delta_{\ell i}
=\sum_{i,\ell=1}^{N_k} M_{ij}M_{\ell m}\int_{\partial B_1} Y_{k,\ell}(x)\,Y_{k,i}(x)\,d{\mathcal{H}}^{n-1}_x\\&&\qquad
=\int_{\partial B_1} \widetilde{Y}_{k,j}(x)\,\widetilde{Y}_{k,m}(x)\,d{\mathcal{H}}^{n-1}_x
=\int_{\partial B_1} {Y}_{k,j}({\mathcal{R}}x)\,{Y}_{k,m}({\mathcal{R}}x)\,d{\mathcal{H}}^{n-1}_x\\&&\qquad=\int_{\partial B_1} {Y}_{k,j}(y)\,{Y}_{k,m}(y)\,d{\mathcal{H}}^{n-1}_y=
\delta_{j m}.
\end{eqnarray*}
Hence, using~\eqref{LORZRPO2}, we have that the transpose of~$M$ is orthogonal too and therefore
\begin{equation*}\begin{split}&
F_k({\mathcal{R}}x,{\mathcal{R}}y)=\sum_{j=1}^{N_k} Y_{k,j}({\mathcal{R}}x)\,Y_{k,j}({\mathcal{R}}y)=\sum_{j=1}^{N_k} \widetilde{Y}_{k,j}(x)\,\widetilde{Y}_{k,j}(y)\\&\qquad=
\sum_{j,i,\ell=1}^{N_k} M_{ij}\,M_{\ell j} \,{Y}_{k,i}(x)\,{Y}_{k,\ell}(y)=\sum_{i,\ell=1}^{N_k} \delta_{i \ell} \,{Y}_{k,i}(x)\,{Y}_{k,\ell}(y)\\&\qquad=\sum_{j=1}^{N_k} {Y}_{k,j}(x)\,{Y}_{k,j}(y)=F_k(x,y).\qedhere
\end{split}\end{equation*}
\end{proof}

\begin{corollary}\label{KSCODROISTOAZO}
There exists\footnote{Up to multiplicative constants,
this polynomial will be precisely
the Legendre polynomial that will be introduced in Theorem~\ref{cq810ua}.
This point will be addressed specifically in the forthcoming Theorem~\ref{CARRLegendre}.} a polynomial~$P_k$ such that, for every~$x$, $y\in\partial B_1$,
$$F_k(x,y)=P_k(x\cdot y).$$
\end{corollary}

\begin{figure}
  \centering
  \includegraphics[width=.35\linewidth]{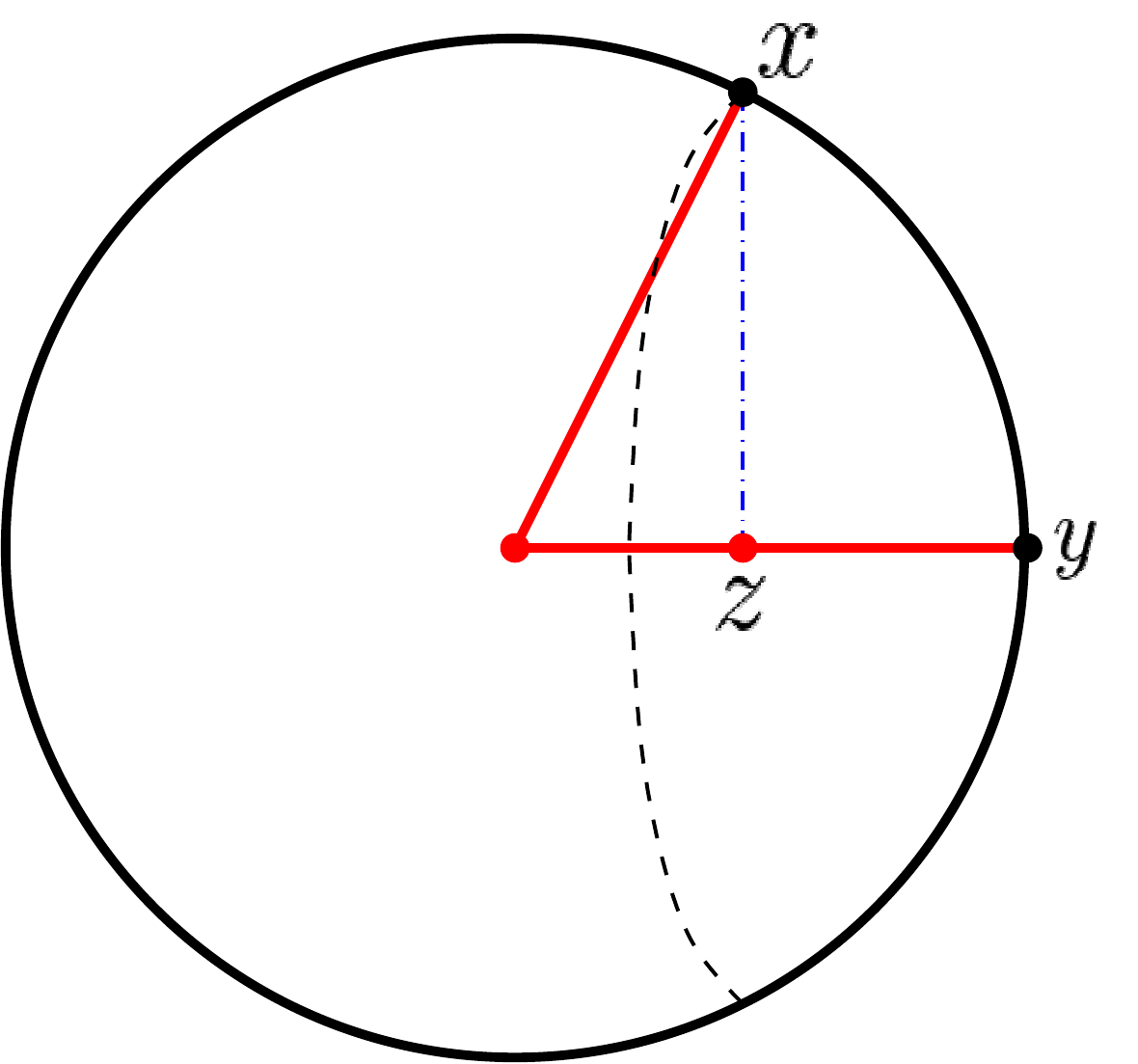}
 \caption{\sl A rotation used in the proof of Corollary~\ref{KSCODROISTOAZO}:
 notice that~$|z|=x\cdot y=t$ and~$|y-z|=\sqrt{1-t^2}$.}\label{LAGIJ7soloDItangeFI}
\end{figure}

\begin{proof} First of all, we note that, being~$Y_{k,j}$ a spherical harmonic,
we can consider it as the restriction to~$\partial B_1$ of a homogeneous polynomial
of degree~$k$, say
$$ Y_{k,j}(x)=\sum_{{\alpha^{(j)}\in\N^n}\atop{|\alpha^{(j)}|=k}} c_{[\alpha^{(j)},j,k]} \;x^{\alpha^{(j)}}$$
for suitable~$ c_{[\alpha^{(j)},j,k]} \in\R$.

This and~\eqref{ASSDM}
give that also~$F_k$ can be considered as the
restriction to~$(\partial B_1)\times(\partial B_1)$
of a suitable polynomial in~$\R^n\times\R^n$
of the form
\begin{equation}\label{YMSIdkan23456oMSMD0odkfgmSyjmta}
F_k(x,y)=\sum_{{1\le j\le N_k}\atop{{
{{{\alpha^{(j)}\in\N^n}\atop{|\alpha^{(j)}|=k}}}
\atop{{\beta^{(j)}\in\N^n}\atop{|\beta^{(j)}|=k}} }}}
c_{[\alpha^{(j)},j,k]} \;c_{[\beta^{(j)},j,k]} \;x^{\alpha^{(j)}}\;x^{\beta^{(j)}}
\end{equation}

Now let~$t:=x\cdot y$. We consider a rotation~${\mathcal{R}}_{(x,y)}$ such that~${\mathcal{R}}_{(x,y)}y=e_1$ and~${\mathcal{R}}_{(x,y)}x=(t,\sqrt{1-t^2},0,\dots,0)=te_1+\sqrt{1-t^2}\,e_2$,
see Figure~\ref{LAGIJ7soloDItangeFI}. By Lemma~\ref{LE21}, we know that
\begin{equation}\label{tgb HTRb8uwejdfby83}
F_k(x,y)=F_k({\mathcal{R}}_{(x,y)}x,\,{\mathcal{R}}_{(x,y)}y)=F_k\Big(te_1+\sqrt{1-t^2}\,e_2,e_1\Big).
\end{equation}
In particular, by~\eqref{YMSIdkan23456oMSMD0odkfgmSyjmta},
one can write~$F_k(x,y)$ as a polynomial in~$t$ and~$\tau:=\sqrt{1-t^2}$, say
\begin{equation}\label{tgb HTRb8uwejdfby84}\begin{split}
F_k(x,y)\,&=\,
\sum_{{1\le j\le N_k}\atop{{
{{{\alpha^{(j)}\in\N^n}\atop{|\alpha^{(j)}|=k}}}
\atop{{\beta^{(j)}\in\N^n}\atop{|\beta^{(j)}|=k}} }}}
c_{[\alpha^{(j)},j,k]} \;c_{[\beta^{(j)},j,k]} \;
(te_1+\tau e_2)^{\alpha^{(j)}}\;e_1^{\beta^{(j)}}
\\&=\,\sum_{{1\le j\le N_k}\atop{{
{{{\alpha^{(j)}\in\N^2}\atop{|\alpha^{(j)}|=k}}} }}}
c_{[\alpha^{(j)},j,k]} \;c_{[ke_1,j,k]} \;
t^{\alpha^{(j)}_1} \,\tau^{\alpha^{(j)}_2}
\\&=\,\sum_{{0\le i,j\le k}\atop{i+j=k}} c_{ij} \,t^i\,\tau^{j},
\end{split}
\end{equation}
for suitable~$c_{ij}\in\R$.

Now we compose~${\mathcal{R}}_{(x,y)}$ with a further rotation~${\mathcal{R}}$ of~$\pi$ radiants around the first axis
and let~${\mathcal{R}}_{(x,y)}':={\mathcal{R}}{\mathcal{R}}_{(x,y)}$. In this way~${\mathcal{R}}_{(x,y)}'y={\mathcal{R}}e_1=e_1$ and~${\mathcal{R}}_{(x,y)}'x={\mathcal{R}}(t,\sqrt{1-t^2},0,\dots,0)=(t,-\sqrt{1-t^2},0,\dots,0)=
te_1-\sqrt{1-t^2}\,e_2$, see again Figure~\ref{LAGIJ7soloDItangeFI}. Arguing as in~\eqref{tgb HTRb8uwejdfby83}
we thus obtain that
\begin{equation*}
F_k(x,y)=F_k\Big(te_1-\sqrt{1-t^2}\,e_2,e_1\Big).
\end{equation*}
That is, in the notation of~\eqref{tgb HTRb8uwejdfby84}
$$ F_k(x,y)=\sum_{{0\le i,j\le k}\atop{i+j=k}} c_{ij} \,t^i\,(-\tau)^{j}.$$
Comparing this and~\eqref{tgb HTRb8uwejdfby84}, we obtain that
$$ \sum_{{0\le i,j\le k}\atop{i+j=k}} c_{ij} \,t^i\,\tau^{j}=\sum_{{0\le i,j\le k}\atop{i+j=k}}
(-1)^j c_{ij} \,t^i\,\tau^{j},$$
and thus
\begin{equation}\label{0-LSSS-PSLspet9YTYUNAPS}
\sum_{{{0\le i,j\le k}\atop{i+j=k}}} d_{ij} \,t^i\,\tau^{j}=0,\end{equation}
where
\begin{equation}\label{THNSPOCOPROGUWAROLDEL-1} d_{ij}:=(1-(-1)^j)c_{ij}=\begin{dcases}
0 &{\mbox{ if $j$ is even,}}\\
2c_{ij} &{\mbox{ if $j$ is odd.}}
\end{dcases}\end{equation}
We claim that
\begin{equation}\label{THNSPOCOPROGUWAROLDEL-2}
{\mbox{$
d_{ij}=0$ for each $i,j\in\{1,\dots,k\}$ with~$i+j=k$.}}
\end{equation}
Indeed, suppose not and take~$d_{IJ}\ne0$ with~$J$ as small as possible
and~$I\in\{1,\dots,k\}$ with~$I+J=k$.
That is, from~\eqref{0-LSSS-PSLspet9YTYUNAPS},
$$0=\sum_{{{1\le j\le k}\atop{i=k-j}}} d_{ij} \,t^i\,\tau^{j}
=\sum_{{{J\le j\le k}\atop{i=k-j}}} d_{ij} \,t^i\,\tau^{j}
=
d_{IJ} \,t^{k-J}\,\tau^J+
\sum_{{{J+1\le j\le k}\atop{i=k-j}}} d_{ij} \,t^i\,\tau^{j}.$$
Dividing by~$\tau^J$, we conclude that
$$0=
d_{IJ} \,t^{k-J}+
\sum_{{{J+1\le j\le k}\atop{i=k-j}}} d_{ij} \,t^i\,\tau^{j-J}.$$
Hence, picking~$\tau:=0$ (which corresponds to~$t:=1$), we find that~$0=d_{IJ}$,
against our original assumption.

This contradiction proves~\eqref{THNSPOCOPROGUWAROLDEL-2}
and thus, recalling~\eqref{THNSPOCOPROGUWAROLDEL-1},
we obtain that~$c_{ij}=0$ whenever~$j$ is odd.

As a result, from~\eqref{tgb HTRb8uwejdfby84},
$$ F_k(x,y)=\sum_{{{0\le i,j\le k}\atop{i+j=k}}\atop{j\in2\N}}
c_{ij} \,t^i\,\tau^{j}=
\sum_{{{0\le m\le k/2}\atop{i=k-2m}}} c_{ij} \,t^i\,(1-t^2)^m.$$
Since the latter is a polynomial in~$t$, the desired result is established.
\end{proof}

\begin{corollary}\label{KMSDCsdpoddpmensvhisbhni6789fkc00PA}
For each~$\omega\in\partial B_1$,
$$ F_k(\omega,\omega)=\frac{N_k}{{\mathcal{H}}^{n-1}(\partial B_1)}.$$
\end{corollary}

\begin{proof} Using~\eqref{LORZRPO} and~\eqref{ASSDM},
\begin{eqnarray*}
\int_{\partial B_1}F_k(\omega,\omega)
\,d{\mathcal{H}}^{n-1}_\omega
=\sum_{j=1}^{N_k} \int_{\partial B_1}\big(Y_{k,j}(\omega)\big)^2
\,d{\mathcal{H}}^{n-1}_\omega=N_k.\end{eqnarray*}
In addition, owing to Corollary~\ref{KSCODROISTOAZO},
we see that~$F_k(\omega,\omega)=P_k(\omega\cdot \omega)=P_k(1)$, which is
a constant. On this account, we conclude that
\begin{equation*}
{\mathcal{H}}^{n-1}(\partial B_1)\,F_k(\omega,\omega)=
\int_{\partial B_1}F_k(\omega,\omega)
\,d{\mathcal{H}}^{n-1}_\omega=N_k.\qedhere
\end{equation*}
\end{proof}

In view of Lemma~\ref{LE21}, and, more generally, in view of the obvious symmetries of the sphere,
it is natural to investigate
the spherical harmonics that remain
invariant under rotations\footnote{To investigate this
type of symmetry, it is convenient to consider~$n\ge3$
since in the plane the only rotation which fixes an axis is obviously
the identity.} about a given (say, the first) coordinate axis. These are known in jargon as Legendre polynomials\index{Legendre polynomial}
(or sometimes, for~$n\ge4$, as Gegenbauer or ultraspherical polynomials) and\index{Gegenbauer polynomial}
are\index{ultraspherical polynomial} characterized by the following result:

\begin{theorem}\label{cq810ua}
There exists a unique spherical harmonic~$L_k$ of degree~$k$ such that
\begin{equation}\label{CARA}
\begin{split}&
L_k(e_1)=1\\ {\mbox{and }}\qquad&
L_k({\mathcal{R}}x)=L_k(x)\qquad{\mbox{for all $x\in\partial B_1$ and all rotations~${\mathcal{R}}$ such that~${\mathcal{R}}e_1=e_1$.}}
\end{split}\end{equation}

Moreover, \begin{equation}\label{CARAB}{\mbox{$L_k$ is a polynomial in~$x_1$,}}\end{equation}
namely there exists a polynomial~$P_k:\R\to\R$
such that
$L_k(x)=P_k(x_1)$ for
all~$x\in\partial B_1$.

If~$j\ne k$, we also have that
\begin{equation}\label{LEGGETEM}
\int_{-1}^1 P_j(t)\,P_k(t)\,(1-t^2)^{\frac{n-3}2}\,dt=0\end{equation}
and
\begin{equation}\label{LEGGETEM2}\begin{split}&{\mbox{$
\{P_0,\dots,P_k\}$ is a basis for the space of
polynomials}}\\&{\mbox{in one variable with degree less than or equal to~$k$.}}\end{split}\end{equation}

Additionally,
the following recurrence relation holds true:
\begin{equation}\label{LEGGEB} (n+k-2)P_{k+1}(t)
=(n+2k-2)tP_{k}(t)-kP_{k-1}(t).\end{equation}

Furthermore, 
\begin{equation}\label{YHNDinho4ryng0657dmgu223344}
{\mbox{$P_k$ has degree equal to~$k$}}\end{equation} and it satisfies
the differential equation
\begin{equation}\label{LEGGE3}
\frac{d}{dt}\left[ (1-t^2)^{\frac{n-1}2}\,{\frac{dP_{k} }{dt}}(t)\right]+
k(n+k-2)(1-t^2)^{\frac{n-3}2}\,P_k(t)
=0.\end{equation}

In addition\footnote{The identity in~\eqref{LEGGE1}
is often referred to with the name of Rodrigues' Formula.
Moreover~\eqref{LEGGE3} is sometimes called
Legendre's Differential Equation
and~\eqref{LEGGEB}
Bonnet's Recursion Formula.}
\begin{equation}\label{LEGGE1} P_k(t)=\frac{(-1)^k}{
\,(1-t^2)^{\frac{n-3}2}\,\displaystyle\prod_{j=1}^k (n+2j-3)}\,\frac{d^k}{dt^k}(1-t^2)^{\frac{n+2k-3}2}.\end{equation}

Moreover, $P_k(t)$ corresponds
to the coefficients in a formal expansion in powers of~$\tau$
of the generating function
\begin{equation}\label{TRVSBlatsds} \big(1-2t\tau+\tau^{2}\big)^{\frac{2-n}2}
=\sum _{k=0}^{+\infty }\mu_k\,P_k(t)\,\tau^k
,\end{equation}
where\footnote{We notice that, when~$n=3$, the quantity~$\mu_k$
in~\eqref{TRVSBlatsds2} is equal to~$1$.}
\begin{equation}\label{TRVSBlatsds2}
\mu_k:=\prod_{j=1}^{k}\frac{j+n-3}{j}
.
\end{equation}

Besides,
\begin{equation}\label{YHNDinho4ryng0657dmgu2233}
P_{k}(-t)= (-1)^kP_k(t).\end{equation}
\end{theorem}

\begin{proof} First off, we construct an example
of spherical harmonic of degree~$k$ satisfying~\eqref{CARA}
(just to be sure\footnote{Strictly speaking, this
step is arguably
not even necessary, since the proof of the uniqueness statement
will also produce a natural candidate
in the forthcoming formulas~\eqref{Cbavmoiu7taso}, \eqref{RECU1}
and~\eqref{RECU2}, but we believe that it is extremely instructive
to rely the Legendre polynomial directly to the fundamental solution, to develop
a direct intuition of the role they play in the Laplace equation.
Furthermore, the uniqueness result and the generating formula
identity in~\eqref{TRVSBlatsds} will be used to deduce other identities
in a rather short and handy fashion.}
that we are not speaking about the empty set!). 
This will also highlight the relation with the fundamental solution as presented in~\eqref{TRVSBlatsds}.
We take~$x$, $y\in\R^n$ with~$0<r:=|x|<|y|=:\varrho$. We also define~$\vartheta$ to be the angle between~$x$ and~$y$,
that is~$\vartheta:=\arccos\frac{x\cdot y}{r\varrho}$. Let also~$\tau:=\frac{r}{\varrho}$ and~$t:=\cos\vartheta$. Then,
\begin{equation}\label{KMS:098765432kjgfas}
|x-y|^{2-n}=\big({|x|^2-2x\cdot y+|y|^2}
\big)^{\frac{2-n}2}
=\frac{1}{
\varrho^{n-2} |\tau^2-2t\tau +1|^{\frac{n-2}2}}
.\end{equation}
Now, if~$z\in\cOMPL$, then
$$z^2-2tz+1=z^2-2\cos\vartheta\,z+1=(z-e^{i\vartheta})(z-e^{-i\vartheta})$$
and thus, if~$|z|<1$,
$$|z^2-2tz+1|=|z-e^{i\vartheta}|\,|z-e^{-i\vartheta}|
\ge (|e^{i\vartheta}|-|z|)(|e^{-i\vartheta}|-|z|)=(1-|z|)^2>0.
$$
As a consequence, we can expand~\eqref{KMS:098765432kjgfas}
in power series of~$\tau$ (which converges
as~$\tau\in[0,1)$,
or even as~$\tau\in\cOMPL$ with~$|\tau|<1$), thus obtaining that
\begin{equation}\label{KM Sdans0otkyeile}
|x-y|^{2-n}
=\frac{1}{
\varrho^{n-2} }\sum_{k=0}^{+\infty} C_k(t)\,\tau^k\end{equation}
for suitable coefficients~$C_k(t)$.

More precisely,
\begin{equation}\label{JS:PAYDSRI} k!\,C_k(t)=
\left.\frac{d^k}{d\tau^k}
\frac{1}{ |\tau^2-2t\tau +1|^{\frac{n-2}2}
}\right|_{\tau=0}
\end{equation}
which shows that~$C_k(t)$ depends polynomially on~$t$, and actually~$C_k(t)$
has degree less than or equal to~$k$. Thus, in the setting of~\eqref{TRVSBlatsds2},
we can define~$P_k:=\frac{C_k}{\mu_k}$ and obtain that\footnote{As a matter
of fact, the statement in~\eqref{YHNDinho4ryng0657dmgu223344BIS}
will be sharpened once we prove~\eqref{YHNDinho4ryng0657dmgu223344}.}
\begin{equation}\label{YHNDinho4ryng0657dmgu223344BIS}
{\mbox{$P_k(t)$ is polynomial of~$t$
of degree less than or equal to~$k$.}}\end{equation}
Additionally, by~\eqref{KM Sdans0otkyeile}, if~$y:=\varrho e_1$,
letting~$Q_k(x):=|x|^kP_k\left(\frac{x_1}{|x|}\right)$,
\begin{equation}\label{ODncocumlydoem} 
|x-y|^{2-n} =
\frac{1}{
\varrho^{n-2} }\sum_{k=0}^{+\infty} \mu_k\,P_k(t)\,\left(\frac{r}{\varrho}\right)^k=
\frac{1}{
\varrho^{n-2} }\sum_{k=0}^{+\infty} \frac{\mu_k}{\varrho^k}\,Q_k(x).
\end{equation}

We now claim that
\begin{equation}\label{YHNDinho4ryng0657dmgu}
{\mbox{$Q_k$ is a homogeneous polynomial of degree~$k$.}}
\end{equation}
For this, since the homogeneity is obvious from the definition of~$Q_k$,
it remains to show that~$Q_k$ is a polynomial. To this end,
we observe that, by~\eqref{JS:PAYDSRI}, one can write~$C_k(t)=
\left.\frac{d^k}{d\tau^k}
\phi(t,\tau)\right|_{\tau=0}$, for a suitable function~$\phi$ such that~$\phi(-t,\tau)=\phi(t,-\tau)$.
For this reason,
$$ C_k(-t)=\left.\frac{d^k}{d\tau^k}
\phi(-t,\tau)\right|_{\tau=0}=\left.\frac{d^k}{d\tau^k}
\phi(t,-\tau)\right|_{\tau=0}=(-1)^k\,C_k(t),$$
and correspondingly
\begin{equation}\label{YHNDinho4ryng0657dmgu22}P_k(-t)=(-1)^k\,P_k(t).
\end{equation}
Hence, each~$P_k$ is the linear combination of monomials
of the type~$t^j$ with~$j\le k$ and~$j$ even if~$k$ is even,
and~$j$ odd if~$k$ is odd.
Consequently, $Q_k(x)$
is the linear combination of terms
of the type~$|x|^k\left(\frac{x_1}{|x|}\right)^j=|x|^{k-j} x_1^j=|x|^{m} x_1^j$,
where~$m:=k-j$ is even, and this completes the proof of~\eqref{YHNDinho4ryng0657dmgu}.

Now, exploiting~\eqref{ODncocumlydoem},
for every~$\varphi\in C^\infty_0(B_\varrho)$ we have that
\begin{eqnarray*}
&& 0=
\int_{\R^n}\varphi(x)\,\Delta|x-y|^{2-n}\,dx =
\int_{\R^n}\Delta\varphi(x)\,|x-y|^{2-n}\,dx\\&&\qquad =
\int_{\R^n}\Delta\varphi(x)\,
\frac{1}{
\varrho^{n-2} }\sum_{k=0}^{+\infty} \frac{\mu_k}{\varrho^k}\,Q_k(x)\,dx=
\frac{1}{
\varrho^{n-2} }\sum_{k=0}^{+\infty} \frac{\mu_k}{\varrho^k}
\int_{\R^n}\Delta\varphi(x)\,Q_k(x)\,dx\\
&&\qquad=
\frac{1}{
\varrho^{n-2} }\sum_{k=0}^{+\infty} \frac{\mu_k}{\varrho^k}
\int_{\R^n}\varphi(x)\,\Delta Q_k(x)\,dx.
\end{eqnarray*}
As a consequence, by
the Identity Principle for power series we deduce that
$$ \mu_k\, \int_{\R^n}\varphi(x)\,\Delta Q_k(x)\,dx=0.$$
Thus, recalling~\eqref{TRVSBlatsds2},
$$ \int_{\R^n}\varphi(x)\,\Delta Q_k(x)\,dx=0,$$
which gives that~$\Delta Q_k(x)=0$ for every~$x\in B_\varrho$.
Since this identity is now valid for all~$\varrho>0$, we thus conclude that~$\Delta Q_k(x)=0$
for every~$x\in\R^n$. This and~\eqref{YHNDinho4ryng0657dmgu}
give that the function
\begin{equation}\label{givenby}
\partial B_1\ni x\mapsto L_k(x):=Q_k(x)=P_k(x_1)\end{equation}
is a spherical harmonic.
We are left to prove~\eqref{CARA}. For this, we exploit~\eqref{ODncocumlydoem} with~$x:=re_1$ and~$y:=\varrho e_1$, with~$0<r<\varrho=1$,
that produce~$\vartheta=0$ and thus~$t=1$, to see that
\begin{equation*}
(1-r)^{2-n} = \sum_{k=0}^{+\infty} \mu_k\,P_k(1)\,r^k.
\end{equation*}
Comparing every order in~$r$ with~\eqref{TRVSBlatsds2}, we thus conclude that~$1=P_k(1)=Q_k(e_1)=L_k(e_1)$.
Moreover, if~${\mathcal{R}}$ is a rotation that fixes~$e_1$, we have
that~$L_k({\mathcal{R}}x)=P_k({\mathcal{R}}x\cdot e_1)=P_k(x_1)=L_k(x)$, completing the proof of~\eqref{CARA}.

We also observe that~\eqref{TRVSBlatsds} holds true in this setting, thanks
to~\eqref{KMS:098765432kjgfas}
and~\eqref{ODncocumlydoem}.
\medskip

Having established one example of spherical harmonics with the desired properties
(as well as the validity of~\eqref{TRVSBlatsds} for this example)
we now turn our attention to the uniqueness claim in the statement of Theorem~\ref{cq810ua}
and to the proof of~\eqref{CARAB}, \eqref{LEGGETEM}, \eqref{LEGGETEM2}, \eqref{LEGGEB},
\eqref{YHNDinho4ryng0657dmgu223344}, \eqref{LEGGE3},
\eqref{LEGGE1} and~\eqref{YHNDinho4ryng0657dmgu2233}.\medskip

More specifically, for the uniqueness claim we argue as follows.
We take a spherical harmonic~$L_k$ of degree~$k$ satisfying~\eqref{CARA}
and we consider the corresponding homogeneous and harmonic polynomial~$Q_k$ which coincides with~$L_k$
along~$\partial B_1$. We extrapolate from every monomial of~$Q_k$ the factor depending on~$x_1$,
namely we write
\begin{equation}\label{IMSPDssacuJD} Q_k(x)=\sum_{j=0}^k x_1^j\,h_{k-j,k}(x_2,\dots,x_n),\end{equation}
with~$h_{i,k}$ homogeneous polynomials in~$x_2,\dots,x_n$ of degree~$i$. Thus, if~${\mathcal{R}}$ is a rotation that fixes the first axis
and we denote by~$\widetilde{{\mathcal{R}}}$ its action\footnote{To be precise,
one can define~$\overline{{\mathcal{R}}}:\R^{n-1}\to\R^n$ as~$\overline{{\mathcal{R}}}X:=
{{\mathcal{R}}}(0,X)$ and notice that~$\overline{{\mathcal{R}}} X\cdot e_1=
(0,X)\cdot e_1=0$. Accordingly~$\overline{{\mathcal{R}}}:\R^{n-1}\to\{0\}\times\R^{n-1}$
and therefore, defining
$$ \widetilde{{\mathcal{R}}}X:=\sum_{i=2}^n \big((\overline{{\mathcal{R}}}X)\cdot e_2\big)e_2,$$
we have that~$\widetilde{{\mathcal{R}}}:\R^{n-1}\to\R^{n-1}$
and~$\overline{{\mathcal{R}}}X=(0,\widetilde{{\mathcal{R}}}X)$.
Furthermore, if~$x=(x_1,X)\in\R\times\R^{n-1}$,
$$ {{\mathcal{R}}}x=
{{\mathcal{R}}}\big(x_1e_1+(0,X)\big)=x_1{\mathcal{R}}e_1+\overline{{\mathcal{R}}}(0,X)=
x_1e_1+(0,\widetilde{{\mathcal{R}}}X)=(x_1,\widetilde{{\mathcal{R}}}X)
.$$

We also remark that, conversely, every rotation~$\widetilde{{\mathcal{R}}}$
on~$\R^{n-1}$ produces a rotation~${\mathcal{R}}$ on~$\R^n$ that fixes the first
coordinate axis, simply by defining~${{\mathcal{R}}}x:=(x_1,\widetilde{{\mathcal{R}}}X)$.}
on the variables~$X:= (x_2,\dots,x_n)\in\R^{n-1}$,
\begin{eqnarray*}&&
0=|x|^k\big(L_k({\mathcal{R}}x)-L_k(x)\big)=|{\mathcal{R}}x|^{k}L_k({\mathcal{R}}x)-|x|^k L_k(x)=Q_k({\mathcal{R}}x)-Q_k(x)\\&&\qquad\qquad=
\sum_{j=0}^k x_1^j\,\Big(h_{k-j,k}(\widetilde{{\mathcal{R}}}X)-h_{k-j,k}(X)\Big)
\end{eqnarray*}
and consequently~$h_{k-j,k}(\widetilde{{\mathcal{R}}}X)=h_{k-j,k}(X)$
for all~$X\in\R^{n-1}$.
That is, each~$h_{i,k}$ only depends on~$|X|$, hence
we write~$h_{i,k}(X)=\widetilde{h}_{i,k}(|X|)$ for
some~$\widetilde{h}_{i,k}:\R\to\R$. As a matter of fact, by the homogeneity of~$
h_{i,k}$, for every~$r\ge0$,
$$ \widetilde{h}_{i,k}(r)=h_{i,k}(re_n)=r^i h_{i,k}(e_n)=c_{i,k} r^i,$$
where~$c_{i,k}:=h_{i,k}(e_n)\in\R$.

Thus, since~$h_{i,k}$ is a polynomial and~$h_{i,k}(X)=\widetilde{h}_{i,k}(|X|)=c_{i,k}\,|X|^i$,
necessarily~$c_{i,k}=0$ whenever~$i$ is odd. As a result,
\begin{equation}\label{Cbavmoiu7taso} Q_k(x)=\sum_{j=0}^k c_{k-j,k}\,x_1^j\,|X|^{k-j}=
\sum_{i=1}^k c_{i,k}\,x_1^{k-i}\,|X|^{i}=\sum_{{\ell\in\N}\atop{0\le\ell\le k/2}}
c_{2\ell,k}\,x_1^{k-2\ell}\,|X|^{2\ell}.\end{equation}
To show the desired uniqueness claim, it thus suffices to show that these~$
c_{2\ell,k}$ are
uniquely determined. First of all, by~\eqref{CARA} and~\eqref{Cbavmoiu7taso},
\begin{equation}\label{RECU1} 1=L_k(e_1)=Q_k(e_1)=c_{0,k}.
\end{equation}
Additionally, from the harmonicity of~$Q_k$,
\begin{eqnarray*}
0&=&\Delta Q_k(x)\\
&=&\sum_{{\ell\in\N}\atop{0\le\ell\le k/2}}
c_{2\ell,k}\,\Delta(x_1^{k-2\ell}\,|X|^{2\ell})\\&=&
\sum_{{\ell\in\N}\atop{0\le\ell\le k/2}}
c_{2\ell,k}\,\partial^2_{x_1}(x_1^{k-2\ell}\,|X|^{2\ell})
+
\sum_{{{\ell\in\N}\atop{0\le\ell\le k/2}}\atop{2\le i\le n}}
c_{2\ell,k}\,\partial^2_{x_i}(x_1^{k-2\ell}\,|X|^{2\ell})\\&=&
\sum_{{\ell\in\N}\atop{0\le\ell\le (k-2)/2}}(k-2\ell)(k-2\ell-1)\,
c_{2\ell,k}\,x_1^{k-2\ell-2}\,|X|^{2\ell}\\&&\qquad
+
\sum_{{{\ell\in\N}\atop{1\le\ell\le k/2}}\atop{2\le i\le n}}
2\ell(2\ell+n-3)\,c_{2\ell,k}\,x_1^{k-2\ell}\,|X|^{2\ell-2}.
\end{eqnarray*}
That is, changing indices' names,
\begin{eqnarray*}
0=\sum_{{m\in\N}\atop{1\le m\le k/2}}\Big((k-2m+2)(k-2m+1)\,
c_{2m-2,k}+
2m(2m+n-3)\,c_{2m,k}
\Big)\,x_1^{k-2m}\,|X|^{2m-2}
.\end{eqnarray*}
This produces the recurrence relation
\begin{equation}\label{RECU2}
c_{2m,k}=
\frac{(k-2m+2)(k-2m+1)\,
c_{2m-2,k}}{
2m(2m+n-3)},
\end{equation}
that, joined with~\eqref{RECU1},
determines uniquely all the coefficients.
This completes the proof of the uniqueness claim in the statement of Theorem~\ref{cq810ua}.
\medskip

We observe that in light of the example constructed and of the uniqueness
statement, the unique spherical harmonic~$L_k$ of degree~$k$
satisfying~\eqref{CARA} is given by~\eqref{givenby}. Hence, \eqref{CARAB}
is also satisfied.

Moreover, we remark that~\eqref{TRVSBlatsds}
and~\eqref{YHNDinho4ryng0657dmgu2233}
have been already proved in this special case
 (recall also~\eqref{YHNDinho4ryng0657dmgu22}),
hence they follow from the uniqueness statement as well.

\medskip

Furthermore, if~$f\in L^\infty([-1,1])$, writing
\[{\mbox{$\partial B_1\cap\{x_1>0\}=\big\{x_1=\sqrt{1-|X|^2}$ with~$X\in\R^{n-1}$
and~$|X|<1\big\}$,}} \]
and letting~$\omega_{n-2}$ the $(n-2)$-dimensional Hausdorff measure
measure of the unit sphere in~$\R^{n-1}$,
we see that
\begin{equation*}
\int_{\partial B_1\cap\{x_1>0\}} f(x_1)\,d{\mathcal{H}}^{n-1}_x=
\int_{{X\in\R^{n-1}}\atop{|X|<1}} \frac{f\big(\sqrt{1-|X|^2}\big)}{\sqrt{1-|X|^2}}
\,dX=\omega_{n-2}\int_{0}^1 \frac{f\big(\sqrt{1-\rho^2}\big)}{\sqrt{1-\rho^2}}
\,\rho^{n-2}\,d\rho.
\end{equation*}
Thus, the substitution~$t:=\sqrt{1-\rho^2}$ leads to
\begin{equation*}
\int_{\partial B_1\cap\{x_1>0\}} f(x_1)\,d{\mathcal{H}}^{n-1}_x=
\omega_{n-2}\int_{0}^1 f(t)
\,(1-t^2)^{\frac{n-3}2}\,dt.
\end{equation*}
Applying this also to~$\widetilde{f}(t):=f(-t)$, we obtain the following useful integral formula
on the sphere for functions of one variable:
\begin{equation*}
\begin{split}
&\int_{\partial B_1} f(x_1)\,d{\mathcal{H}}^{n-1}_x
=\int_{\partial B_1\cap\{x_1>0\}} f(x_1)\,d{\mathcal{H}}^{n-1}_x+
\int_{\partial B_1\cap\{x_1<0\}} f(x_1)\,d{\mathcal{H}}^{n-1}_x\\&\qquad=
\omega_{n-2}\int_{0}^1 f(t)
\,(1-t^2)^{\frac{n-3}2}\,dt+
\int_{\partial B_1\cap\{x_1>0\}} \widetilde{f}(x_1)\,d{\mathcal{H}}^{n-1}_x\\&\qquad
=\omega_{n-2}\left[\int_{0}^1 f(t)
\,(1-t^2)^{\frac{n-3}2}\,dt+\int_{0}^1 \widetilde{f}(t)
\,(1-t^2)^{\frac{n-3}2}\,dt\right]\\&\qquad=
\omega_{n-2}\int_{-1}^1 f(t)
\,(1-t^2)^{\frac{n-3}2}\,dt.\end{split}
\end{equation*}
As a result,
\begin{eqnarray*}\omega_{n-2}
\int_{-1}^1 P_j(t)\,P_k(t)\,(1-t^2)^{\frac{n-3}2}\,dt=\int_{\partial B_1} P_j(x_1)\,P_k(x_1)\,d{\mathcal{H}}^{n-1}_x
=\int_{\partial B_1} L_j(x_1)\,L_k(x_1)\,d{\mathcal{H}}^{n-1}_x.
\end{eqnarray*}
The latter quantity is zero, thanks to the orthogonality condition for spherical harmonics of different degree in Lemma~\ref{DIGGDJKDL}, and this proves~\eqref{LEGGETEM}.\medskip

Also, the set~$\{P_0,\dots,P_k\}$
contains~$k+1$ elements and they are
necessarily linearly independent, thanks to~\eqref{LEGGETEM},
whence~\eqref{LEGGETEM2} plainly follows.\medskip

Now we focus on the proof of~\eqref{LEGGEB}. 
For this, we differentiate~\eqref{TRVSBlatsds}
with respect to~$\tau$, finding that
\[ \frac{(2-n)(\tau-t)}{({1-2t\tau+\tau^{2}})^{\frac{n}2}}
=\sum _{k=1}^{+\infty }k\mu_k\,P_k(t)\,\tau^{k-1}.\]
Comparing with~\eqref{TRVSBlatsds} we thereby gather that
\begin{equation*}
\frac{(2-n)(\tau-t)}{{1-2t\tau+\tau^{2}}}\sum _{k=0}^{+\infty }\mu_k\,P_k(t)\,\tau^k
=\sum _{k=1}^{+\infty }k\mu_k\,P_k(t)\,\tau^{k-1}.\end{equation*} 
As a result,
\begin{eqnarray*}&&
(2-n)\sum _{k=1}^{+\infty }\mu_{k-1}\,P_{k-1}(t)\,\tau^k-
(2-n)\sum _{k=0}^{+\infty }\mu_k\,t\,P_k(t)\,\tau^k\\&&\qquad
=
(2-n)\sum _{k=0}^{+\infty }\mu_k\,P_k(t)\,\tau^{k+1}-
(2-n)\sum _{k=0}^{+\infty }\mu_k\,t\,P_k(t)\,\tau^k\\&&\qquad=
(2-n)(\tau-t)\sum _{k=0}^{+\infty }\mu_k\,P_k(t)\,\tau^k
=({1-2t\tau+\tau^{2}})\sum _{k=1}^{+\infty }k\mu_k\,P_k(t)\,\tau^{k-1}\\&&\qquad=
\sum _{k=1}^{+\infty }k\mu_k\,P_k(t)\,\tau^{k-1}
-2\sum _{k=1}^{+\infty }k\mu_k\,t\,P_k(t)\,\tau^{k}
+\sum _{k=1}^{+\infty }k\mu_k\,P_k(t)\,\tau^{k+1}\\&&\qquad=
\sum _{k=0}^{+\infty }(k+1)\mu_{k+1}\,P_{k+1}(t)\,\tau^{k}
-2\sum _{k=1}^{+\infty }k\mu_k\,t\,P_k(t)\,\tau^{k}
+\sum _{k=2}^{+\infty }(k-1)\mu_{k-1}\,P_{k-1}(t)\,\tau^{k}.
\end{eqnarray*}
Taking the~$k$th order of this identity, it follows that
\begin{equation}\label{3085vb9trghjgrbk}
(2-n)\mu_{k-1}\,P_{k-1}-(2-n)\mu_k \,tP_k=
(k+1)\mu_{k+1}\,P_{k+1}-2k\mu_k\,tP_k+(k-1)\mu_{k-1}\,P_{k-1}.\end{equation}
Since from~\eqref{TRVSBlatsds2} we have that
$$ \mu_{k+1}=\frac{n+k-2}{k+1}\,\mu_k\qquad
{\mbox{and}}\qquad \mu_{k-1}=\frac{n}{n+k-3}\,\mu_k,$$
the identity in~\eqref{3085vb9trghjgrbk}
leads to~\eqref{LEGGEB}, as desired.
\medskip

To prove~\eqref{YHNDinho4ryng0657dmgu223344} we argue by induction.
First of all, using~\eqref{JS:PAYDSRI}
we see that the desired claim holds true for~$k\in\{0,1\}$.
Indeed,
\begin{equation}\begin{split}\label{0e1PL}
& P_0(t)=\frac{C_0(t)}{\mu_0}=1
\\{\mbox{and }} \qquad &
P_1(t)=\frac{C_1(t)}{\mu_1}=\frac{(n-2)t}{n-2}=t.
\end{split}\end{equation}
To perform the inductive step, we suppose that~$P_j$ has degree~$j$ for all~$j\in\N\cap[0,k]$.
Then, we exploit~\eqref{LEGGEB} to deduce that~$P_{k+1}$ has degree~$k+1$
and this establishes~\eqref{YHNDinho4ryng0657dmgu223344}.
\medskip

Now we prove~\eqref{LEGGE3}. To this end, we notice that
\begin{equation}\label{GIOmnsde895col34yleokmdfg4}
\begin{split}
\frac{d}{dt}\left[ (1-t^2)^{\frac{n-1}2}\,{\frac{dP_{k} }{dt}}(t)\right]\,&=\,
-(n-1)t(1-t^2)^{\frac{n-3}2}\,{\frac{dP_{k} }{dt}}(t)
+
(1-t^2)^{\frac{n-1}2}\,{\frac{d^2P_{k} }{dt^2}}(t)\\&
=\,(1-t^2)^{\frac{n-3}2}\,T_k(t),\end{split}\end{equation}
where \begin{equation*}T_k(t):=
(1-n)t \,{\frac{dP_{k} }{dt}}(t)+
(1-t^2)\,{\frac{d^2P_{k} }{dt^2}}(t).
\end{equation*}
We also observe that~$T_k$ is a polynomial
of degree less than or equal to~$k$,
owing to~\eqref{YHNDinho4ryng0657dmgu223344BIS}.
Consequently, in light of~\eqref{LEGGETEM2},
we can write~$T_k$ as a linear combination of~$\{P_0,\dots,P_k\}$,
namely there exist~$s_{j,k}\in\R$ such that
\begin{equation}\label{SCRIUTGSt} T_k=\sum_{j=0}^k s_{j,k}\,P_j.\end{equation}
As a result, for each~$i\in\{0,\dots,k\}$,
recalling~\eqref{LEGGETEM} we get that
\begin{equation}\label{GIOmnsde895col34yleokmdfg4B}\begin{split}
\int_{-1}^1 P_i(t)\,T_k(t)\,(1-t^2)^{\frac{n-3}2}\,dt\,&=\,
\sum_{j=0}^k s_{j,k}
\int_{-1}^1 P_i(t)\,P_j(t)\,(1-t^2)^{\frac{n-3}2}\,dt\\&=\,
s_{i,k}
\int_{-1}^1 (P_i(t))^2\,(1-t^2)^{\frac{n-3}2}\,dt.\end{split}
\end{equation}
In addition, 
$$ \lim_{t\to\pm1} (1-t^2)^{\frac{n-1}2}=0.$$
Therefore, we can exploit~\eqref{GIOmnsde895col34yleokmdfg4}
and~\eqref{GIOmnsde895col34yleokmdfg4B} and integrate by parts twice, to find that
\begin{equation}\label{SCRIUTGSt2}
\begin{split}&
s_{i,k}
\int_{-1}^1 (P_i(t))^2\,(1-t^2)^{\frac{n-3}2}\,dt=\int_{-1}^1
(1-t^2)^{\frac{n-3}2}\,T_k(t)\,P_i(t)\,dt\\&\qquad=
\int_{-1}^1
\frac{d}{dt}\left[ (1-t^2)^{\frac{n-1}2}\,{\frac{dP_{k} }{dt}}(t)\right]\,P_i(t)\,dt
=-\int_{-1}^1 (1-t^2)^{\frac{n-1}2}\,{\frac{dP_{k} }{dt}}(t)\,{\frac{dP_{i} }{dt}}(t)\,dt\\&\qquad=
\int_{-1}^1
\frac{d}{dt}\left[ (1-t^2)^{\frac{n-1}2}\,{\frac{dP_i }{dt}}(t)\right]\,P_k(t)\,dt=
\int_{-1}^1(1-t^2)^{\frac{n-3}2}\,T_i(t)\,P_k(t)\,dt.
\end{split}\end{equation}
We stress that~\eqref{GIOmnsde895col34yleokmdfg4} was exploited again in the last identity, but used
with the index~$i$ instead of~$k$. Thus, on account of~\eqref{LEGGETEM}, \eqref{SCRIUTGSt}
and~\eqref{SCRIUTGSt2}, for every~$i\in\{0,\dots,k-1\}$,
\begin{equation}\label{BUnvdt9tto0}
s_{i,k}
\int_{-1}^1 (P_i(t))^2\,(1-t^2)^{\frac{n-3}2}\,dt=s_{j,i}\sum_{j=1}^i
\int_{-1}^1(1-t^2)^{\frac{n-3}2}\,P_j(t)\,P_k(t)\,dt=0.
\end{equation}
Since, by~\eqref{CARA}, we know that~$P_i(1)=L_i(e_1)=1$,
we deduce from~\eqref{BUnvdt9tto0} that~$s_{i,k}=0$ whenever~$i\in\{0,\dots,k-1\}$.
As a byproduct, we infer from~\eqref{SCRIUTGSt} that~$T_k= s_{k,k}\,P_k$
and consequently~\eqref{GIOmnsde895col34yleokmdfg4} boils down to
\begin{equation}\label{MSonydwihTHDiubu45G}
\frac{d}{dt}\left[ (1-t^2)^{\frac{n-1}2}\,{\frac{dP_{k} }{dt}}(t)\right]\,=\,s_{k,k}\,(1-t^2)^{\frac{n-3}2}P_k(t)
.\end{equation}
Now we denote by~$d_k$ the highest power of~$t$
in the monomial expansion of~$P_k$, namely, owing to~\eqref{YHNDinho4ryng0657dmgu223344}, we write~$P_k(t)=d_k t^{k}+\widetilde{P}_k(t)$,
where~$d_k\ne0$, and~$\widetilde{P}_k$ is a polynomial of degree strictly lower that~$k$
(and possibly zero). 
Therefore, as~$t\to+\infty$,
\begin{eqnarray*}
&&\frac{d}{dt}\left[ (1-t^2)^{\frac{n-1}2}\,{\frac{dP_{k} }{dt}}(t)\right]=
\frac{d}{dt}\left[ (1-t^2)^{\frac{n-1}2}\,
\left(kd_k t^{k-1}+{\frac{d\widetilde{P}_k}{dt}}(t)
\right)\right]\\&&\qquad=
-(n-1)t(1-t^2)^{\frac{n-3}2}\,\left(kd_k t^{k-1}+{\frac{d\widetilde{P}_k}{dt}}(t)
\right)\\&&\qquad\qquad\qquad+
(1-t^2)^{\frac{n-1}2}\,
\left(k(k-1)d_k t^{k-2}+
{\frac{d^2\widetilde{P}_k}{dt^2}}(t)
\right)\\&&\qquad=\Big(
(-1)^{\frac{n-1}2}(n-1)t^{n-2}+O(t^{n-3})\,\Big(kd_k t^{k-1}+O(t^{k-2}\Big)\\&&\qquad\qquad\qquad+\Big(
(-1)^{\frac{n-1}2} t^{n-1}+O(t^{n-2})\Big)\Big(k(k-1)d_k t^{k-2}+
O(t^{k-3})\Big)\\&&\qquad=
(-1)^{\frac{n-1}2}
kd_k(n+k-2)\,
t^{n+k-3}
+O(t^{n+k-4})
\end{eqnarray*}
and
\begin{eqnarray*}
(1-t^2)^{\frac{n-3}2}P_k(t)&=&\Big((-1)^{\frac{n-3}2}t^{n-3}+O(t^{n-4})\Big)\Big(d_kt^{k}+O(t^{k-1})\Big)\\&=&
(-1)^{\frac{n-3}2}d_kt^{n+k-3}+O(t^{n+k-4}).
\end{eqnarray*}
{F}rom these observations and~\eqref{MSonydwihTHDiubu45G}
we deduce that, for large~$t$,
$$ (-1)^{\frac{n-1}2}
kd_k(n+k-2)\,
t^{n+k-3}=(-1)^{\frac{n-3}2}s_{k,k}d_kt^{n+k-3}+O(t^{n+k-4})$$
and therefore~$s_{k,k}=-k (n+k-2)$.
Combining this and~\eqref{MSonydwihTHDiubu45G}
we obtain the desired result in~\eqref{LEGGE3}.\medskip

Now we prove~\eqref{LEGGE1}.
For this, we note that, for all~$\alpha\in\R$ and~$\ell\in\N$,
\begin{equation}\label{KSMCLSinbdupolii}
\frac{d^\ell}{dt^\ell}(1-t^2)^{\frac{\alpha}2}=(1-t^2)^{\frac{\alpha-2\ell}2} p_{\alpha,\ell}(t),
\end{equation}
for a suitable polynomial~$p_{\alpha,\ell}$ of degree less than or equal to~$\ell$.
To check this, we argue by induction. When~$\ell=0$
the claim in~\eqref{KSMCLSinbdupolii} is obviously true with~$p_{\alpha,0}:=1$.
Suppose now that~\eqref{KSMCLSinbdupolii} is true for the index~$\ell$. Then, taking one more derivative,
we have that
\begin{eqnarray*}&& \frac{d^{\ell+1}}{dt^{\ell+1}}(1-t^2)^{\frac{\alpha}2}=
\frac{d}{dt}\Big( (1-t^2)^{\frac{\alpha-2\ell}2} p_{\alpha,\ell}(t)\Big)\\&&\qquad
=(2\ell-\alpha)t(1-t^2)^{\frac{\alpha-2(\ell+1)}2} p_{\alpha,\ell}(t)
+
(1-t^2)^{\frac{\alpha-2\ell}2} \frac{dp_{\alpha,\ell}}{dt}(t)\\&&\qquad
=(1-t^2)^{\frac{\alpha-2(\ell+1)}2} p_{\alpha,\ell+1}(t),\end{eqnarray*}
with~$p_{\alpha,\ell+1}:=
(2\ell-\alpha)t p_{\alpha,\ell}
+(1-t^2)\frac{dp_{\alpha,\ell}}{dt}$ that is a polynomial of degree up to~$\ell+1$.
This completes the inductive step and establishes~\eqref{KSMCLSinbdupolii} as desired.

Now we define
$$ \zeta_k:=
\frac{(-1)^k}{\displaystyle\prod_{j=1}^k (n+2j-3)}\qquad{\mbox{and}}\qquad
P^\star_k(t):=\zeta_k\,(1-t^2)^{\frac{3-n}2}
\frac{d^k}{dt^k}(1-t^2)^{\frac{n+2k-3}2}.$$
Observe that~$P^\star_k$ is precisely the right hand side of~\eqref{LEGGE1}.

Now we claim that, if~$j<k$ and~$i\in\{0,\dots,j\}$, then
\begin{equation}\label{L:IND:INT:ASMDdefer}
\int_{-1}^1 
\frac{d^j}{dt^j}(1-t^2)^{\frac{n+2j-3}2}\,
\frac{d^k}{dt^k}(1-t^2)^{\frac{n+2k-3}2}
\,(1-t^2)^{\frac{3-n}2}\,dt
=
\int_{-1}^1 q_{i,j}(t)\,\frac{d^{k-i}}{dt^{k-i}}(1-t^2)^{\frac{n+2k-3}2}\,dt.
\end{equation}
Here~$q_{i,j}$ is a suitable polynomial of degree less than or equal to~$j-i$.
To prove this, we argue by induction over~$i$.
When~$i=0$, we define~$q_{0,j}(t):=
\frac{d^j}{dt^j}(1-t^2)^{\frac{n+2j-3}2}
\,(1-t^2)^{\frac{3-n}2}$ and we have that the identity in~\eqref{L:IND:INT:ASMDdefer}
holds true. What is more~$q_{0,j}$ is a
polynomial of degree less than or equal to~$j$, thanks to~\eqref{KSMCLSinbdupolii}.

To perform the inductive step, we suppose that the desired claim is true for all indices up to some~$i$, with~$i\le j-1$,
and we aim at proving it for the index~$i+1$. To this end, we let~$\psi_{i,k}(t):=
\frac{d^{k-i-1}}{dt^{k-i-1}}(1-t^2)^{\frac{n+2k-3}2}$. 
Using~\eqref{KSMCLSinbdupolii}, we have that~$\psi_{i,k}(t)=
(1-t^2)^{\frac{n+2i-1}2}\widetilde{q}_{i,k}(t)$, for a suitable polynomial~$\widetilde{q}_{i,k}$. In particular,
$$ \psi_{i,k}(\pm1)=
(1-(\pm1)^2)^{\frac{n+2i-1}2}\widetilde{q}_{i,k}(\pm1)
=0.$$
Thus, by the inductive assumption
and an integration by parts,
\begin{eqnarray*}
&&
\int_{-1}^1 
\frac{d^j}{dt^j}(1-t^2)^{\frac{n+2j-3}2}\,
\frac{d^k}{dt^k}(1-t^2)^{\frac{n+2k-3}2}
\,(1-t^2)^{\frac{3-n}2}\,dt
=
\int_{-1}^1 q_{i,j}(t)\,\frac{d\psi_{i,k}}{dt}(t)\,dt\\
&&\qquad=q_{i,j}(1)\psi_{i,k}(1)-q_{i,j}(-1)\psi_{i,k}(-1)
-\int_{-1}^1 \frac{d q_{i,j}}{dt}(t)\,\psi_{i,k}(t)\,dt\\
&&\qquad=0+\int_{-1}^1 q_{i+1,j}(t)\,\frac{d^{k-i-1}}{dt^{k-i-1}}(1-t^2)^{\frac{n+2k-3}2}\,dt
\end{eqnarray*}
with~$q_{i+1,j}:=-\frac{d q_{i,j}}{dt}$ and the latter is a polynomial with degree at most~$j-i-1$.
The proof of~\eqref{L:IND:INT:ASMDdefer} is thereby complete.

We also observe that, if~$j<k$,
\begin{equation}\label{96787wPIUJSssMMS}
\left. \frac{d^{k-j-1}}{dt^{k-j-1}}(1-t^2)^{\frac{n+2k-3}2}\right|_{t=\pm1}=0.
\end{equation}
Indeed, by~\eqref{KSMCLSinbdupolii}, we know that~$\frac{d^{k-j-1}}{dt^{k-j-1}}(1-t^2)^{\frac{n+2k-3}2}$
is equal to some polynomial times~$(1-t^2)^{\frac{n+2j-1}2}$, and the latter function vanishes at~$t=\pm1$,
thus proving~\eqref{96787wPIUJSssMMS}.

Now we exploit~\eqref{L:IND:INT:ASMDdefer} with~$i:=j$.
Notice that the polynomial~$q_{j,j}(t)$ has degree zero, hence it is a constant (still denoted by~$q_{j,j}$ for simplicity).
Therefore, if~$j<k$, using also~\eqref{96787wPIUJSssMMS}, we have that
\begin{eqnarray*}&&
\int_{-1}^1 P_j^\star(t)\,P_k^\star(t)\,(1-t^2)^{\frac{n-3}2}\,dt=
\zeta_j\,\zeta_k\,\int_{-1}^1
\frac{d^j}{dt^j}(1-t^2)^{\frac{n+2j-3}2}
\frac{d^k}{dt^k}(1-t^2)^{\frac{n+2k-3}2}
(1-t^2)^{\frac{3-n}2}\,dt
\\&&\qquad=
\zeta_j\,\zeta_k\,q_{j,j}\,
\int_{-1}^1 \frac{d^{k-j}}{dt^{k-j}}(1-t^2)^{\frac{n+2k-3}2}\,dt\\&&\qquad=
\zeta_j\,\zeta_k\,q_{j,j}\,\left(
\left. \frac{d^{k-j-1}}{dt^{k-j-1}}(1-t^2)^{\frac{n+2k-3}2}\right|_{t=1}
-\left. \frac{d^{k-j-1}}{dt^{k-j-1}}(1-t^2)^{\frac{n+2k-3}2}\right|_{t=-1}
\right)=0,
\end{eqnarray*}
and therefore, up to exchanging~$j$ and~$k$,
\begin{equation}\label{ByhdclkdERnicsdes967y2INT}
\int_{-1}^1 P_j^\star(t)\,P_k^\star(t)\,(1-t^2)^{\frac{n-3}2}\,dt=0\quad{\mbox{ whenever~$j\ne k$.}}\end{equation}

Now we claim that
\begin{equation}\label{ByhdclkdERnicsdes967y2}
{\mbox{$P_k^\star(t)=P_k(t)$ for all~$k\in\N$.}}
\end{equation}
To this end, we argue by induction.
Recalling~\eqref{0e1PL}, we remark that~$P_0^\star(t)=1=P_0(t)$,
hence~\eqref{ByhdclkdERnicsdes967y2} is valid for~$k=0$.
Let us now suppose that~\eqref{ByhdclkdERnicsdes967y2} holds true for all indices~$\{0,\dots,k-1\}$
with~$k\ge1$ and we prove it for the index~$k$. For this, we remark that~$P_k^\star$ is a polynomial of degree at most~$k$, owing to~\eqref{KSMCLSinbdupolii}. Hence, by~\eqref{LEGGETEM2},
we can write~$P_k^\star$ as a linear combination of~$\{P_0,\dots,P_k\}$,
namely there exist~$\sigma_{j,k}\in\R$ such that
\begin{equation*}P^\star_k=\sum_{j=0}^k \sigma_{j,k}\,P_j.\end{equation*}
This and the inductive assumption lead to
\begin{equation*}P^\star_k=\sigma_{k,k}\,P_k+\sum_{j=0}^{k-1} \sigma_{j,k}\,P_j^\star.\end{equation*}
Therefore, for each~$i\in\{1,\dots,k-1\}$, exploiting~\eqref{LEGGETEM}
and~\eqref{ByhdclkdERnicsdes967y2INT}
we have that
\begin{eqnarray*}
0&=&
\int_{-1}^1 P_i^\star(t)\,P_k^\star(t)\,(1-t^2)^{\frac{n-3}2}\,dt
\\&=&\sigma_{k,k}\,
\int_{-1}^1 P_i^\star(t)\,P_k(t)\,(1-t^2)^{\frac{n-3}2}\,dt
+\sum_{j=0}^{k-1} \sigma_{j,k}\,\int_{-1}^1 P_i^\star(t)
P_j^\star(t)\,(1-t^2)^{\frac{n-3}2}\,dt\\&=&\sigma_{k,k}\,
\int_{-1}^1 P_i (t)\,P_k(t)\,(1-t^2)^{\frac{n-3}2}\,dt
+ \sigma_{i,k}\,\int_{-1}^1 (P_i(t))^2\,(1-t^2)^{\frac{n-3}2}\,dt\\
&=&\sigma_{i,k}\,\int_{-1}^1 (P_i(t))^2\,(1-t^2)^{\frac{n-3}2}\,dt.
\end{eqnarray*}
Thus, since~$P_i(1)=1$ due to~\eqref{CARA}, we gather that~$\sigma_{i,k}=0$
for all~$i\in\{1,\dots,k-1\}$, and accordingly
\begin{equation}\label{NSNienfTEFNNFasd}
P^\star_k=\sigma_{k,k}\,P_k.\end{equation}
Consequently, to complete the proof of~\eqref{ByhdclkdERnicsdes967y2}, it remains to show that
\begin{equation}\label{LKSncngikfdloe}
\sigma_{k,k}=1.\end{equation}
For this, we note that if~$j\in\{1,\dots,k\}$,
$$\left|\lim_{t\nearrow1}
(1-t)^{\frac{3-n}2} \frac{d^{k-j}}{dt^{k-j}}(1-t)^{\frac{n+2k-3}2}\right|=\left|\lim_{y\searrow0}
y^{\frac{3-n}2} \frac{d^{k-j}}{dy^{k-j}}y^{\frac{n+2k-3}2}\right|=0.
$$
Thus, by the General Leibniz Rule for higher order derivative,
\begin{eqnarray*}&&
\lim_{t\nearrow1}
(1-t^2)^{\frac{3-n}2}
\frac{d^k}{dt^k}(1-t^2)^{\frac{n+2k-3}2}=
\lim_{t\nearrow1}(1+t)^{\frac{3-n}2}
(1-t)^{\frac{3-n}2}
\frac{d^k}{dt^k}\Big((1+t)^{\frac{n+2k-3}2}(1-t)^{\frac{n+2k-3}2}\Big)\\
&&\qquad=
2^{\frac{3-n}2}\lim_{t\nearrow1}
(1-t)^{\frac{3-n}2}\sum_{j=0}^k{{k}\choose{j}}
\frac{d^j}{dt^j}(1+t)^{\frac{n+2k-3}2} \frac{d^{k-j}}{dt^{k-j}}(1-t)^{\frac{n+2k-3}2}\\&&\qquad
=2^{\frac{3-n}2}\lim_{t\nearrow1}
(1-t)^{\frac{3-n}2}
(1+t)^{\frac{n+2k-3}2} \frac{d^{k}}{dt^{k}}(1-t)^{\frac{n+2k-3}2}\\&&\qquad
=2^{k}\lim_{t\nearrow1}
(1-t)^{\frac{3-n}2} \frac{d^{k}}{dt^{k}}(1-t)^{\frac{n+2k-3}2}=
(-1)^k\,2^{k}\,\prod_{i=0}^{k-1} {\frac{n+2(k-i)-3}2}\\&&\qquad=
(-1)^k\,\prod_{i=0}^{k-1} (n+2(k-i)-3)=(-1)^k\,\prod_{j=1}^{k} (n+2j-3)=\frac1{\zeta_k}.
\end{eqnarray*}
This shows that~$P^\star_k(1)=1$ and thus, comparing with~\eqref{NSNienfTEFNNFasd},
we find that~$\sigma_{k,k}=1$.
This proves~\eqref{LKSncngikfdloe}, and thus~\eqref{ByhdclkdERnicsdes967y2},
from which we obtain~\eqref{LEGGE1}, as desired.
\end{proof}

For a list of the first ten Legendre polynomials in dimension~$3$, see e.g.
the footnote on page~\pageref{LISTACON}
(at this level, we will not exploit much the explicit values of the 
Legendre polynomials but rather their
fundamental algebraic structure).\medskip

Now we point out an additional identity for Legendre polynomials
in terms of the function introduced in~\eqref{ASSDM}
(as a byproduct, this characterizes the Legendre polynomials
precisely as the polynomials for which Corollary~\ref{KSCODROISTOAZO}
holds true, up to normalization constants).

\begin{theorem}\label{CARRLegendre} Let~$P_k$ be the Legendre polynomials in Theorem~\ref{cq810ua}
and~$F_k$ be as in~\eqref{ASSDM}, with~$N_k$ as in~\eqref{ALEHD:DJFJFJ}. 

Then, for every~$x$, $y\in\partial B_1$,
$$ P_k( x\cdot y)=\frac{{\mathcal{H}}^{n-1}(\partial B_1)}{N_k} \,F_k(x,y).$$
\end{theorem}

\begin{proof} For every~$\omega\in\partial B_1$, we define
$$ L^\star_k(\omega):=\frac{{\mathcal{H}}^{n-1}(\partial B_1)}{N_k} \,F_k(\omega,e_1).$$
By~\eqref{ASSDM}, we know that~$L^\star_\omega$
is a spherical harmonic. Moreover, by
Corollary~\ref{KMSDCsdpoddpmensvhisbhni6789fkc00PA},
$$ L^\star_k(e_1)=\frac{{\mathcal{H}}^{n-1}(\partial B_1)}{N_k} \,F_k(e_1,e_1)=1.$$
Additionally, if~${\mathcal{R}}$ is a rotation on~$\R^n$ such that~${\mathcal{R}}e_1=e_1$,
we infer from Lemma~\ref{LE21} that
$$ L^\star_k({\mathcal{R}}\omega)
=\frac{{\mathcal{H}}^{n-1}(\partial B_1)}{N_k} \,F_k({\mathcal{R}}\omega,e_1)=
\frac{{\mathcal{H}}^{n-1}(\partial B_1)}{N_k} \,F_k({\mathcal{R}}\omega,{\mathcal{R}}e_1)=
\frac{{\mathcal{H}}^{n-1}(\partial B_1)}{N_k} \,F_k(\omega,e_1)
= L^\star_k(\omega).$$
These observations and the
uniqueness claim in Theorem~\ref{cq810ua} entail that~$L^\star_k=L_k$.

Therefore, if~$x$, $y\in\partial B_1$
and~${\mathcal{R}}_y$ is a rotation on~$\R^n$
such that~${\mathcal{R}}_y y=e_1$,
we let~$\omega:={\mathcal{R}}_yx$. In this way,
noticing that
$$ \omega\cdot e_1={\mathcal{R}}_yx\cdot e_1=
x\cdot {\mathcal{R}}_y^{-1}e_1=x\cdot y,$$
and using Lemma~\ref{LE21} once again,
we conclude that
\begin{equation*}
\begin{split}&\frac{{\mathcal{H}}^{n-1}(\partial B_1)}{N_k} \,F_k(x,y)=
\frac{{\mathcal{H}}^{n-1}(\partial B_1)}{N_k} \,F_k({\mathcal{R}}_yx,{\mathcal{R}}_yy)
=
\frac{{\mathcal{H}}^{n-1}(\partial B_1)}{N_k} \,F_k(\omega,e_1)
\\&\qquad\qquad=L^\star_k(\omega)=L_k(\omega)=P_k(\omega_1)
=P_k(x\cdot y).
\qedhere\end{split}\end{equation*}
\end{proof}

\begin{corollary} We have that
\begin{eqnarray}\label{SCOMDCh9wiecickyhw9edivn}&&
\max_{[-1,1]}|P_k|\le1\\{\mbox{and, for every~$x\in\partial B_1$, }}&&\label{SCOMDCh9wiecickyhw9edivn2}
\int_{\partial B_1}\big(
P_k(x\cdot y)\big)^2\,d{\mathcal{H}}^{n-1}_y=\frac{ {\mathcal{H}}^{n-1}(\partial B_1) }{N_k} \,.
\end{eqnarray}
\end{corollary}

\begin{proof}
Given~$t\in[-1,1]$, we use Theorem~\ref{CARRLegendre} 
with~$x:=(t,\sqrt{1-t^2},0,\dots,0)$ and~$y:=e_1$
and we see that
\begin{eqnarray*}&&\big(
P_k(t)\big)^2=\big(
P_k( x\cdot y)\big)^2=\frac{\big({\mathcal{H}}^{n-1}(\partial B_1)\big)^2}{N_k^2} \,
\big(F_k(x,y)\big)^2
=
\frac{\big({\mathcal{H}}^{n-1}(\partial B_1)\big)^2}{N_k^2} \,\left(
\sum_{j=1}^{N_k} Y_{k,j}(x)\,Y_{k,j}(y)\right)^2
.\end{eqnarray*}
This and the Cauchy-Schwarz Inequality, together with a further application of Theorem~\ref{CARRLegendre},
lead to
\begin{eqnarray*}&& \big(
P_k(t)\big)^2\le \left(
\frac{ {\mathcal{H}}^{n-1}(\partial B_1) }{N_k}
\sum_{j=1}^{N_k} \big(Y_{k,j}(x)\big)^2
\right)\left(
\frac{ {\mathcal{H}}^{n-1}(\partial B_1) }{N_k}
\sum_{j=1}^{N_k} \big(Y_{k,j}(y)\big)^2
\right)\\&&\qquad=P_k(x\cdot x)\,P_k(y\cdot y)=P_k(1)\,P_k(1).
\end{eqnarray*}
This and the fact that~$P_k(1)=1$
(recall the normalization condition in~\eqref{CARA}) give~\eqref{SCOMDCh9wiecickyhw9edivn},
as desired.

Furthermore, making use again Theorem~\ref{CARRLegendre},
and recalling the normalization in~\eqref{CARA}
and the orthogonality condition in~\eqref{LORZRPO},
\begin{eqnarray*}&&
\int_{\partial B_1}\big(
P_k(x\cdot y)\big)^2\,d{\mathcal{H}}^{n-1}_y
=
\frac{\big({\mathcal{H}}^{n-1}(\partial B_1)\big)^2}{N_k^2} \,
\sum_{i,j=1}^{N_k} 
\int_{\partial B_1}
Y_{k,i}(x)\,Y_{k,i}(y)\,
Y_{k,j}(x)\,Y_{k,j}(y)
\,d{\mathcal{H}}^{n-1}_y\\&&\qquad=
\frac{\big({\mathcal{H}}^{n-1}(\partial B_1)\big)^2}{N_k^2} \,
\sum_{j=1}^{N_k} 
\big(Y_{k,j}(x)\big)^2=\frac{{\mathcal{H}}^{n-1}(\partial B_1)}{N_k}\,P_k(x\cdot x)=
\frac{{\mathcal{H}}^{n-1}(\partial B_1)}{N_k}\,P_k(1)=
\frac{{\mathcal{H}}^{n-1}(\partial B_1)}{N_k},
\end{eqnarray*}
thus establishing~\eqref{SCOMDCh9wiecickyhw9edivn2}.
\end{proof}

Now we provide a technical observation related to the algebraic nondegeneracy
of spherical harmonics:

\begin{lemma}\label{sNOKMDCsand02or3tgREFGDBammdpo}
Let~$m\in\N\cap[1, N_k]$ and~$\{j_1,\dots,j_m\}\subseteq\{1,\dots,N_k\}$. Then, there
exist~$x^{(1)},\dots,x^{(m)}\in\partial B_1$
such that
\begin{equation}\label{DE:MA:AL}
\det\left( \begin{matrix}
Y_{k,j_1}(x^{(1)}) &\dots& Y_{k,j_1}(x^{(m)})\cr
\vdots\cr
Y_{k,j_m}(x^{(1)}) &\dots& Y_{k,j_m}(x^{(m)})
\end{matrix}
\right)\ne0.
\end{equation}
\end{lemma}

\begin{proof}
We argue by induction over~$m$. First of all, by~\eqref{LORZRPO},
we know that~$\int_{\partial B_1} (Y_{k,j}(x))^2\,d{\mathcal{H}}^{n-1}_x=1$
and this gives~\eqref{DE:MA:AL} when~$m=1$.

Suppose now that~\eqref{DE:MA:AL} is satisfied for some~$m\in\N\cap[1,N_k-1]$.
To establish it for~$m+1$, for every~$x\in\partial B_1$ we consider the function
\begin{equation}\label{PHIdet} \phi(x):=\det\left( \begin{matrix}
Y_{k,j_1}(x^{(1)}) &\dots& Y_{k,j_1}(x^{(m)})&Y_{k,j_1}(x)\cr
\vdots\cr
Y_{k,j_m}(x^{(1)}) &\dots& Y_{k,j_m}(x^{(m)})&Y_{k,j_m}(x)\cr
Y_{k,j_{m+1}}(x^{(1)})&\dots&Y_{k,j_{m+1}}(x^{(m)})&Y_{k,j_{m+1}}(x)
\end{matrix}
\right).\end{equation}
Notice that the desired result is proved once there exists~$x^{(m+1)}\in\partial B_1$
for which~$\phi(x^{(m+1)})\ne0$. Suppose, by contradiction, that~$\phi$ vanishes
identically on~$\partial B_1$.
Then, expanding the determinant in~\eqref{PHIdet} down the last column,
for every~$x\in\partial B_1$ we have that
$$ 0=\phi(x)=\sum_{i=1}^{m+1} {\mathcal{D}}_i(x^{(1)},\dots,x^{(m)})\,Y_{k,j_i}(x),$$
where~${\mathcal{D}}_i$ denotes the corresponding cofactor minor determinant.
Hence, multiplying by~$Y_{k,j_\ell}(x)$, integrating over~$x\in\partial B_1$,
and recalling~\eqref{LORZRPO},
$$ 0= {\mathcal{D}}_\ell(x^{(1)},\dots,x^{(m)})\qquad{\mbox{ for all }}\ell\in\{1,\dots,m+1\}.$$
In particular,
$$0={\mathcal{D}}_{m+1}(x^{(1)},\dots,x^{(m)})=\det\left( \begin{matrix}
Y_{k,j_1}(x^{(1)}) &\dots& Y_{k,j_1}(x^{(m)})\cr
\vdots\cr
Y_{k,j_m}(x^{(1)}) &\dots& Y_{k,j_m}(x^{(m)})
\end{matrix}
\right).$$
But the latter term is different from zero, in view of the inductive hypothesis,
and thus we have reached a contradiction and completed the proof of~\eqref{DE:MA:AL}.
\end{proof}

With this, we are in the position of writing every spherical harmonic
as the superposition of Legendre polynomials in suitable directions
(or as a suitable average):

\begin{theorem}
Let~$Y_k$ be a spherical harmonic of degree~$k$.
Then, there exist~$a_1,\dots,a_{N_k}\in\R$
and~$x^{(1)},\dots,x^{(N_k)}\in\partial B_1$ such that
\begin{equation}\label{KMS-tegrmrmi-he-da0palde893jfgevsta023rtra1}
Y_k(x)=\sum_{j=1}^{N_k}a_j\,P_k(x\cdot x^{(j)})\qquad{\mbox{for all }}x\in\partial B_1.
\end{equation}

Furthermore,
\begin{equation}\label{KMS-tegrmrmi-he-da0palde893jfgevsta023rtra2}
Y_k(x)=\frac{N_k}{{\mathcal{H}}^{n-1}(\partial B_1)}
\int_{\partial B_1} Y_k(\omega)\,P_k(x\cdot\omega)\,d{\mathcal{H}}^{n-1}_\omega
\qquad{\mbox{for all }}x\in\partial B_1.
\end{equation}
\end{theorem}

\begin{proof} Using Theorem~\ref{CARRLegendre}
and Lemma~\ref{sNOKMDCsand02or3tgREFGDBammdpo},
for every~$x\in\partial B_1$, for all~$j\in\{1,\dots,N_k\}$,
$$ P_k( x\cdot x^{(j)})=\frac{{\mathcal{H}}^{n-1}(\partial B_1)}{N_k} \,F_k(x,x^{(j)})=
\frac{{\mathcal{H}}^{n-1}(\partial B_1)}{N_k}\,\sum_{j=1}^{N_k} Y_{k,j}(x)\,Y_{k,j}(x^{(j)}),
$$
that is
\begin{equation}\label{egatsdhBRIDMgnaca}\left(
\begin{matrix}
P_k( x\cdot x^{(1)})\\ \vdots \\P_k( x\cdot x^{(N_k)})
\end{matrix}\right)=\frac{{\mathcal{H}}^{n-1}(\partial B_1)}{N_k}\,
\left( \begin{matrix}
Y_{k,1}(x^{(1)}) &\dots& Y_{k,1}(x^{(N_k)})\cr
\vdots\cr
Y_{k,1}(x^{(1)}) &\dots& Y_{k,N_k}(x^{(N_k)})
\end{matrix}
\right)\left(
\begin{matrix}
Y_{k,1}( x)\\ \vdots \\Y_{k,N_k}( x)
\end{matrix}\right),
\end{equation}
and the matrix in the right hand side of~\eqref{egatsdhBRIDMgnaca} is invertible.
Hence, inverting this matrix in~\eqref{egatsdhBRIDMgnaca}, we obtain that,
for each~$m\in\{1,\dots,N_k\}$,
$$ Y_{k,m}(x)=\sum_{j=1}^{N_k} c_{j,k,m}\,P_k(x\cdot x^{(j)}),$$
for suitable~$c_{j,k,m}\in\R$.
Also, since~$\{Y_{k,1},\dots,Y_{k,N_k}\}$ is a basis,
\begin{equation}\label{KMS Dkrtgf64hfngcheka-sd}
Y_{k}(x)=\sum_{m=1}^{N_k} a_{k,m}\,Y_{k,m}(x),\end{equation}
for suitable~$a_{k,m}\in\R$, and we thereby conclude that
\begin{equation*}
Y_{k}(x)=\sum_{j=1}^{N_k} \left(
\sum_{m=1}^{N_k} a_{k,m}\,c_{j,k,m}\right)\,P_k(x\cdot x^{(j)}),
\end{equation*}
which gives the desired result in~\eqref{KMS-tegrmrmi-he-da0palde893jfgevsta023rtra1}
with~$a_j:=
\sum_{m=1}^{N_k} a_{k,m}\,c_{j,k,m}$.\medskip

Now we prove~\eqref{KMS-tegrmrmi-he-da0palde893jfgevsta023rtra2}.
To this end, we use again~\eqref{KMS Dkrtgf64hfngcheka-sd},
joined with the additional identity in
Theorem~\ref{CARRLegendre} and the orthogonality condition
in~\eqref{LORZRPO}, to find that
\begin{eqnarray*}
&&\frac{N_k}{{\mathcal{H}}^{n-1}(\partial B_1)}
\int_{\partial B_1} Y_k(\omega)\,P_k(x\cdot\omega)\,d{\mathcal{H}}^{n-1}_\omega\\
&=&\int_{\partial B_1} Y_k(\omega)\,F_k(x,\omega)\,d{\mathcal{H}}^{n-1}_\omega\\
&=&\sum_{j=1}^{N_k}\int_{\partial B_1} Y_k(\omega)\, Y_{k,j}(x)\,Y_{k,j}(\omega)
\,d{\mathcal{H}}^{n-1}_\omega\\
&=&\sum_{j,m=1}^{N_k}\int_{\partial B_1} a_{k,m}\,Y_{k,m}(\omega)\, Y_{k,j}(x)\,Y_{k,j}(\omega)
\,d{\mathcal{H}}^{n-1}_\omega\\
&=&\sum_{m=1}^{N_k} a_{k,m}\,Y_{k,m}(x)
\\&=&Y_k(x),
\end{eqnarray*}
and this completes the proof of~\eqref{KMS-tegrmrmi-he-da0palde893jfgevsta023rtra2}.
\end{proof}

For a concrete application of spherical harmonics and Legendre Polynomials
to the evaluation of the gravitational field, see Appendix~\ref{GOEAPL}.
See also~\cite{MR1355826, MR1805196, MR2254749, MR3290046} for an extensive treatment
of spherical harmonics.\medskip

Following~\cite{MR2734448}, we now
give another proof of Theorem~\ref{HARMOu3}. This proof is interesting
since it uses harmonic polynomials (and moreover, as pointed out
in~\cite{MR2734448}, it is flexible enough to address also the case
of higher order operators).

\begin{proof}[Another proof of Theorem~\ref{HARMOu3}]\label{5ddg-5115566}
We know that the function~$f$ is attained by~$u$ at the boundary, owing to~\eqref{HARMOu2}.

We show that
\begin{equation}\label{HARMOu}
{\mbox{$u$ is harmonic in~$B_1$.}}
\end{equation}
For this, we consider a homogeneous harmonic polynomial~$F$, say of degree~$m\ge1$,
and we define~$f(y):=F(y-P)$ for every~$y\in\partial B_1$.
Furthermore, we let
\begin{equation} \label{KJSMSduIJD}
u_F(P):=\fint_{\partial B_1} \ell_e(P)\,d{\mathcal{H}}^{n-1}_e.\end{equation}
Notice that this definition coincides with that in~\eqref{DEUGRGBD},
but we emphasize here the dependence of the function~$u$ on~$F$.
We also observe that
\begin{equation}\label{KJSMSduIJD2}
F(0)=0
\end{equation}
and, recalling the notation in~\eqref{PIUJSMENP},
\begin{eqnarray*}
&&f(Q_\pm(e))=F(Q_\pm(e)-P)
=F\big(r_\pm(e)\,e\big)=
(r_\pm(e))^m F(e)
.\end{eqnarray*}
As a result, recalling~\eqref{EELS},
\begin{equation}\label{KMSD-okjnbfe8u5pskd}
\begin{split}
\ell_e(P)\,&=\, \frac{
r_+(e)f(Q_-(e))-r_-(e)f(Q_+(e))
}{r_+(e)-r_-(e)}
\\&=\, \frac{
r_+(e)(r_-(e))^m -r_-(e) (r_+(e))^m
}{r_+(e)-r_-(e)}\,F(e)
\\&=\, 
r_+(e) r_-(e)\,
\frac{(r_-(e))^{m-1} -(r_+(e))^{m-1}
}{r_+(e)-r_-(e)}\,F(e).
\end{split}
\end{equation}
Now we use that
$$ (r_-(e)-r_+(e))\sum_{j=0}^{m-2} (r_-(e))^j(r_+(e))^{m-2-j}
=(r_-(e))^{m-1}-(r_+(e))^{m-1}.$$
With this and~\eqref{ANT6}, we write~\eqref{KMSD-okjnbfe8u5pskd}
in the form
\begin{equation}\label{SL:SKM4D}
\ell_e(P) =
(1-|P|^2)\,F(e)
\sum_{j=0}^{m-2} (r_-(e))^j (r_+(e))^{m-2-j} .
\end{equation}
Now we remark that, if~$j\in\N\cap\left[0,\frac{m-2}2\right]$,
\begin{eqnarray*}&&
(r_-(e))^j (r_+(e))^{m-2-j}+
(r_-(e))^{m-2-j} (r_+(e))^{j}=
(r_-(e)\,r_+(e))^{j}
\Big(
(r_+(e))^{m-2-2j}
+
(r_-(e))^{m-2-2j} 
\Big)\\&&\qquad=(|P|^2-1)^j
\Big(
(r_+(e))^{m-2-2j}
+
(r_-(e))^{m-2-2j} 
\Big).
\end{eqnarray*}
As a consequence, using the change of index~$J:=m-2-j$,
\begin{eqnarray*}&&
\sum_{j=0}^{m-2} (r_-(e))^j (r_+(e))^{m-2-j}\\& =&
\sum_{0\le j<(m-2)/2} (r_-(e))^j (r_+(e))^{m-2-j} 
+\sum_{(m-2)/2<j\le m-2} (r_-(e))^j (r_+(e))^{m-2-j}
\\&&\qquad \qquad +\chi_{2\N}(m) (r_-(e))^{\frac{m-2}2} (r_+(e))^{\frac{m-2}2}\\&=&
\sum_{0\le j<(m-2)/2} (r_-(e))^j (r_+(e))^{m-2-j} 
+\sum_{0\le J<(m-2)/2} (r_-(e))^{m-2-J} (r_+(e))^{j}
\\&&\qquad \qquad+\chi_{2\N}(m) (r_-(e))^{\frac{m-2}2} (r_+(e))^{\frac{m-2}2}
\\&=&
\sum_{0\le j<(m-2)/2}\Big( (r_-(e))^j (r_+(e))^{m-2-j} +
(r_-(e))^{m-2-J} (r_+(e))^{j}\Big)
+\chi_{2\N}(m) (r_-(e))^{\frac{m-2}2} (r_+(e))^{\frac{m-2}2}\\&=&
\sum_{0\le j<(m-2)/2}(|P|^2-1)^j
\Big((r_+(e))^{m-2-2j}+(r_-(e))^{m-2-2j} \Big)\\&&\qquad\qquad
+\chi_{2\N}(m) 
(|P|^2-1)^{\frac{m-2}2}
\Big((r_+(e))^{\frac{m-2}2}+
(r_-(e))^{\frac{m-2}2} \Big)
\\&=&\sum_{
j\in\N\cap\left[0,\frac{m-2}2\right]
} C_{j,P}\,\Big(
(r_+(e))^{m-2-2j}
+
(r_-(e))^{m-2-2j} 
\Big),
\end{eqnarray*}
for suitable coefficients~$C_{j,P}$ not depending on~$e$.
Plugging this information into~\eqref{SL:SKM4D}, we find that
\begin{equation}\label{KS-O-SM04t46}
\ell_e(P)=
F(e)\,
\sum_{
j\in\N\cap\left[0,\frac{m-2}2\right]
} C_{j,P}\,\Big(
(r_+(e))^{m-2-2j}
+
(r_-(e))^{m-2-2j} 
\Big),
\end{equation}
up to renaming the coefficients~$C_{j,P}$.

We also point out that, exploiting~\eqref{CHISMD-0193248OKE},
for each~$j\in\N\cap\left[0,\frac{m-2}2\right]$
\begin{eqnarray*}&&
(r_+(e))^{m-2-2j}
+
(r_-(e))^{m-2-2j} \\&=&
\Big(-P\cdot e+\sqrt{D(e)}\Big)^{m-2-2j}+\Big(-P\cdot e-\sqrt{D(e)}\Big)^{m-2-2j}
\\
&=&\sum_{k=0}^{m-2-2j}
\binom{m-2-2j}{k}
\left(
(-P\cdot e)^{m-2-2j-k} \big(\sqrt{D(e)}\big)^{k}
+(-P\cdot e)^{m-2-2j-k} \big(-\sqrt{D(e)}\big)^k
\right)\\
&=&
\sum_{{0\le k\le m-2-2j}\atop{k \in2\N}}
\binom{m-2-2j}{k}
\left(
(-P\cdot e)^{m-2-2j-k} \big(\sqrt{D(e)}\big)^{k}
+(-P\cdot e)^{m-2-2j-k} \big(-\sqrt{D(e)}\big)^k
\right)
\\&=&
2\sum_{0\le i\le (m-2-2j)/2}
\binom{m-2-2j}{2i}
(-P\cdot e)^{m-2-2j-2i} (D(e) )^{i}
\\&=&
2\sum_{0\le i\le (m-2-2j)/2}
\binom{m-2-2j}{2i}
(-P\cdot e)^{m-2-2j-2i} \big(
(P\cdot e)^2-|P|^2+1
\big)^{i},
\end{eqnarray*}
which is a polynomial of degree at most~$m-2$
in the variable~$P\cdot e$ (and therefore in the variable~$e$), that we denote by~$\Pi_{j,P}$.
We infer from this and~\eqref{KS-O-SM04t46} that
\begin{equation}\label{S5MJDO0-5-eirkeggj}
\ell_e(P)=\Pi_P(e)\,F(e),\end{equation}
for a suitable polynomial~$\Pi_P$ of degree at most~$m-2$.

Now we exploit Lemma~\ref{L2efsdfpolti6}
to find a harmonic polynomial~$H_P$ of degree at most~$m-2$
such that~$H_P=\Pi_P$ along~$\partial B_1$.
Hence, since~$F$ was supposed to be a
homogeneous harmonic polynomial of degree~$m$, we deduce from Lemma~\ref{DIGGDJKDL} that
$$ \int_{\partial B_1} F(e)\,\Pi_P(e)\,d{\mathcal{H}}^{n-1}_e=
\int_{\partial B_1} F(e)\,H_P(e)\,d{\mathcal{H}}^{n-1}_e=0.$$
This, \eqref{KJSMSduIJD}, \eqref{KJSMSduIJD2}
and~\eqref{S5MJDO0-5-eirkeggj} lead to
\begin{equation}\label{S334S}
u_F(P)=\fint_{\partial B_1} \ell_e(P)\,d{\mathcal{H}}^{n-1}_e=
\fint_{\partial B_1} 
\Pi_P(e)\,F(e)
\,d{\mathcal{H}}^{n-1}_e=0=F(0).\end{equation}
We note that~$F$ is the unique (by Corollary~\ref{UNIQUENESSTHEOREM})
harmonic function in~$B_1(-P)$ with its datum on~$\partial B_1(-P)$.
We can thus reconstruct its values from the Poisson Kernel
in Theorems~\ref{POIBALL1} and~\ref{POIBALL}, namely we write,
for all~$x\in B_1(-P)$,
$$ F(x)=\int_{\partial B_1(-P)}
F(y)\,P_{B_1(-P)}(x,y)\,d{\mathcal{H}}^{n-1}_y,$$
where~$P_{B_1(-P)} $ denotes the Poisson Kernel of~$B_1(-P)$.

This and~\eqref{S334S} give that
\begin{equation}\label{KSDM-oKD3DD}
u_F(P)=\int_{\partial B_1(-P)}
F(y)\,P_{B_1(-P)}(0,y)\,d{\mathcal{H}}^{n-1}_y.\end{equation}
We stress that this result was obtained by assuming that~$F$
is a homogeneous harmonic polynomial, but we are now
going to extend~\eqref{KSDM-oKD3DD} to any continuous function~$F$ on~$\partial B_1(-P)$. 
To this end, we first remark that, by linearity, we see that~\eqref{KSDM-oKD3DD} holds true when~$F$
is a harmonic polynomial (not necessarily homogeneous). Accordingly,
we make use of Corollary~\ref{DENSPO} and we find
a sequence of harmonic polynomials~$F_j$ that converge to~$F$ in~$L^2(\partial B_1(-P))$:
we thus deduce from~\eqref{KSDM-oKD3DD} that
\begin{equation}\begin{split}\label{05943h5v866767467}
\lim_{j\to+\infty}u_{F_j}(P)=\,&\lim_{j\to+\infty}
\int_{\partial B_1(-P)}
F_j(y)\,P_{B_1(-P)}(0,y)\,d{\mathcal{H}}^{n-1}_y\\
=\,&\int_{\partial B_1(-P)}
F(y)\,P_{B_1(-P)}(0,y)\,d{\mathcal{H}}^{n-1}_y
\\=\,&
\int_{\partial B_1}
F(y+P)\,P_{B_1(-P)}(0,y+P)\,d{\mathcal{H}}^{n-1}_y
\\=\,&
\int_{\partial B_1}
f(y)\,P_{B_1}(P,y)\,d{\mathcal{H}}^{n-1}_y
,\end{split}\end{equation}
where~$P_{B_1}$ denotes the Poisson Kernel of~$B_1$.

On the other hand, we recall~\eqref{EELS} and we write
\begin{eqnarray*}
\ell_e(P)=\frac{r_+(e)f(Q_-(e))-r_-(e)f(Q_+(e))}{r_+(e)-r_-(e)}=
\frac{r_+(e)F(Q_-(e)-P)-r_-(e)F(Q_+(e)-P)}{r_+(e)-r_-(e)}.
\end{eqnarray*}
This and~\eqref{KJSMSduIJD} give that
\begin{eqnarray*}&&
\big|u_{F}(P)-u_{F_j}(P)\big|\,{\mathcal{H}}^{n-1}(\partial B_1)\\&\le&
\int_{\partial B_1} 
\left(\frac{r_+(e)}{r_+(e)-r_-(e)}\big|F(Q_-(e)-P)-F_j(Q_-(e)-P)\big|
\right. \\&&\qquad\qquad \left.+
\frac{-r_-(e)}{r_+(e)-r_-(e)}\big|F(Q_+(e)-P)-F_j(Q_+(e)-P)\big|
\right)
\,d{\mathcal{H}}^{n-1}_e\\ &\le&
\sqrt{\int_{\partial B_1} 
\left(\frac{r_+(e)}{r_+(e)-r_-(e)}\right)^2\,d{\mathcal{H}}^{n-1}_e}\cdot\sqrt{\int_{\partial B_1} 
\big|F(Q_-(e)-P)-F_j(Q_-(e)-P)\big|^2\,d{\mathcal{H}}^{n-1}_e}\\&&\qquad+
\sqrt{\int_{\partial B_1} 
\left(\frac{-r_-(e)}{r_+(e)-r_-(e)}\right)^2\,d{\mathcal{H}}^{n-1}_e}\cdot\sqrt{\int_{\partial B_1} 
\big|F(Q_+(e)-P)-F_j(Q_+(e)-P)\big|^2
\,d{\mathcal{H}}^{n-1}_e}.
\end{eqnarray*}
We notice that~$r_+(e)$, $-r_-(e)\le2$ and that
$$r_+(e)-r_-(e)=2\sqrt{D(e)}=2\sqrt{(P\cdot e)^2-|P|^2+1}\ge 2\sqrt{1-|P|^2}.$$
As a consequence,
\begin{equation}\begin{split}\label{ud8435bv798t7tbv8y9}
\big|u_{F}(P)-u_{F_j}(P)\big|\le\,& C_{P}\left(\sqrt{\int_{\partial B_1} 
\big|F(Q_-(e)-P)-F_j(Q_-(e)-P)\big|^2\,d{\mathcal{H}}^{n-1}_e}\right. \\&\qquad\left.+\sqrt{\int_{\partial B_1} 
\big|F(Q_+(e)-P)-F_j(Q_+(e)-P)\big|^2
\,d{\mathcal{H}}^{n-1}_e}\right),\end{split}\end{equation}
for some~$C_{P}>0$ depending on~$P$ and~$n$.

Also, exploiting~\eqref{VEBDNFKRGZINDMKBVINSKDSA-22}
with
$$g(\omega):=\big|F(\omega-P)-F_j(\omega-P)\big|^2 \frac{1-|P|^2-r_-\left( \frac{P-Q_-(\omega)}{|P-Q_-(\omega)|}
\right)P\cdot \frac{P-Q_-(\omega)}{|P-Q_-(\omega)|}}{\left|r_-\left( \frac{P-Q_-(\omega)}{|P-Q_-(\omega)|}
\right)\right|^n}
,$$ we obtain that
\begin{eqnarray*}&&
\int_{\partial B_1} 
\big|F(Q_-(e)-P)-F_j(Q_-(e)-P)\big|^2\,d{\mathcal{H}}^{n-1}_e\\&&\qquad\qquad=
\int_{\partial B_1} 
g(Q_-(e))\, \frac{|r_-(e)|^n}{1-|P|^2-r_-(e)P\cdot e}
\,d{\mathcal{H}}^{n-1}_e=\int_{\partial B_1} 
g(e)\,d{\mathcal{H}}^{n-1}_e.
\end{eqnarray*}
Now we observe that
$$ \frac{1-|P|^2-r_-\left( \frac{P-Q_-(\omega)}{|P-Q_-(\omega)|}
\right)P\cdot \frac{P-Q_-(\omega)}{|P-Q_-(\omega)|}}{\left|r_-\left( \frac{P-Q_-(\omega)}{|P-Q_-(\omega)|}
\right)\right|^n}\le \frac{4}{(1-|P|)^n},
$$
which gives that
\begin{eqnarray*}
&&\int_{\partial B_1} 
\big|F(Q_-(e)-P)-F_j(Q_-(e)-P)\big|^2\,d{\mathcal{H}}^{n-1}_e\le C_P\int_{\partial B_1} 
\big|F(e-P)-F_j(e-P)\big|^2\,d{\mathcal{H}}^{n-1}_e\\&&\qquad\qquad= C_P\int_{\partial B_1(-P)} 
\big|F(\zeta)-F_j(\zeta)\big|^2\,d{\mathcal{H}}^{n-1}_\zeta,
\end{eqnarray*}
up to renaming~$C_P>0$. Similarly, 
\begin{eqnarray*}
&&\int_{\partial B_1} 
\big|F(Q_+(e)-P)-F_j(Q_+(e)-P)\big|^2\,d{\mathcal{H}}^{n-1}_e\le  C_P\int_{\partial B_1(-P)} 
\big|F(\zeta)-F_j(\zeta)\big|^2\,d{\mathcal{H}}^{n-1}_\zeta.
\end{eqnarray*}
Therefore, recalling~\eqref{ud8435bv798t7tbv8y9}, we obtain that
\begin{equation*}\begin{split}
\big|u_{F}(P)-u_{F_j}(P)\big|\le\,& C_P
\int_{\partial B_1(-P)} 
\big|F(\zeta)-F_j(\zeta)\big|^2\,d{\mathcal{H}}^{n-1}_\zeta,
\end{split}\end{equation*}
up to renaming~$C_P$ once again. Consequently,
$$ \lim_{j\to+\infty} u_{F_j}(P)=u_{F}(P).$$
{F}rom this and~\eqref{05943h5v866767467} we deduce that
$$ u_F(P)=\int_{\partial B_1}
f(y)\,P_{B_1}(P,y)\,d{\mathcal{H}}^{n-1}_y.$$
Consequently, by the Poisson Kernel representation in Theorems~\ref{POIBALL1} and~\ref{POIBALL},
we have that~$u_F$ is harmonic
and this completes the proof of~\eqref{HARMOu}.

The proof of Theorem~\ref{HARMOu3} is thus completed in view
of~\eqref{HARMOu2} and~\eqref{HARMOu}.
\end{proof}

\section{Almansi's Formula}

As a short digression, we now point out that
the representation method introduced in Lemma~\ref{L2efsdfpolti6}
has a natural counterpart in the classical Almansi's Formula\index{Almansi's Formula}, according to which
any function satisfying~$\Delta^ku=0$ can be represented
by functions that are harmonic (that is, polyharmonic functions
can be algebraically reconstructed from harmonic ones).
The details go as follows:

\begin{lemma}\label{353ALMAMAL}
Let~$k\in\N$, $k\ge2$.
Let~$\Omega$ be a starshaped open subset of~$\R^n$.
Then, for every~$u\in C^{2k}(\overline\Omega)$ such that
\begin{equation}\label{ALMA:1}
{\mbox{$\Delta^ku=0$
in~$\Omega$}}\end{equation}
there exist~$v$, $w\in C^{2(k-1)}(\overline\Omega)$
such that
\begin{equation}\label{ALMA:2}
{\mbox{$\Delta^{k-1}v=\Delta^{k-1}w=0$
in~$\Omega$}}\end{equation} and
\begin{equation}\label{ALMA:3} u(x)=|x|^2\,v(x)+w(x)\qquad{\mbox{for every }}x\in\Omega.\end{equation}
\end{lemma}

\begin{proof} First of all, we claim that, for every~$i\in\N$, $i\ge1$,
and any smooth function~$f$,
\begin{equation}\label{iONBDFverpoverufgvvRFDV}
\Delta^i(|x|^2f)=c_i\Delta^{i-1}f+d_ix\cdot\nabla\Delta^{i-1}f+|x|^2\Delta^i f,
\end{equation}
with
$$ c_i:=2ni+4i(i-1)\qquad{\mbox{and}}\qquad
d_i:=4i .$$
To check this, we argue by induction.
When~$i=1$, the claim in~\eqref{iONBDFverpoverufgvvRFDV}
reduces to~$\Delta(|x|^2f)=2n f+4 x\cdot\nabla f+|x|^2\Delta f$,
which is indeed true. Let us now suppose that~\eqref{iONBDFverpoverufgvvRFDV}
holds true for some index~$i$. Then, we find that
\begin{eqnarray*}
\Delta^{i+1}(|x|^2f)&=&\Delta\Big(
c_i\Delta^{i-1}f+d_ix\cdot\nabla\Delta^{i-1}f+|x|^2\Delta^i f\Big)\\
&=&c_i\Delta^{i}f+d_ix\cdot\nabla\Delta^{i}f
+2d_i\Delta^{i}f
+2n\Delta^i f
+4 x\cdot\nabla\Delta^i f+
|x|^2\Delta^{i+1} f\\&=&
(c_i+2d_i+2n)\Delta^{i}f
+(d_i+4)x\cdot\nabla\Delta^{i}f
+|x|^2\Delta^{i+1} f
\\&=&(2ni+4i(i-1)+8i+2n)\Delta^{i}f
+(4i+4)x\cdot\nabla\Delta^{i}f
+|x|^2\Delta^{i+1} f\\&=&c_{i+1}\Delta^{i}f
+d_{i+1} x\cdot\nabla\Delta^{i}f
+|x|^2\Delta^{i+1} f,
\end{eqnarray*}
which completes the inductive step and establishes the validity of~\eqref{iONBDFverpoverufgvvRFDV}.

Now, we define
\begin{equation}\label{ALMALPA}
\alpha:=\frac{c_{k-1}}{d_{k-1}}\in(0,+\infty).\end{equation}
Up to a translation, we can assume that~$\Omega$
is starshaped with respect to the origin and then, for every~$x\in\Omega$,
we can define
$$ T(x):=\frac1{d_{k-1}}\int_0^1 t^{\alpha-1}\Delta^{k-1}u(tx)\,dt.$$
We remark that, for all~$x\in\Omega$,
\begin{equation}\label{ALMA:5}
\Delta T(x)=\frac1{d_{k-1}}\int_0^1 t^{\alpha+1}\Delta^{k}u(tx)\,dt=0,
\end{equation}
thanks to~\eqref{ALMA:1}, and furthermore
\begin{equation}\label{1ALMALPA}
\begin{split}
x\cdot\nabla T(x)\,&=\,\frac1{d_{k-1}}\int_0^1 t^{\alpha}x\cdot\nabla\Delta^{k-1}u(tx)\,dt
\\& =\,\frac1{d_{k-1}
}\int_0^1\left( \frac{d}{dt}\Big(t^{\alpha} \Delta^{k-1}u(tx)\Big)-\alpha
t^{\alpha-1} \Delta^{k-1}u(tx)
\right)\,dt\\&=\,
\frac1{d_{k-1}
} \Delta^{k-1}u(x)-\alpha T(x).
\end{split}
\end{equation}
Now we take a function~$\zeta_1$ such
that~$\Delta\zeta_1=T$ in~$\Omega$: for this, for instance,
we can define
$$ \widetilde{T}:=\begin{dcases} T &{\mbox{ in }}\Omega,\\
0&{\mbox{ in }}\R^n\setminus\Omega,\end{dcases}$$
and exploit the fundamental solution in Theorem~\ref{KScvbMS:0okrmt4ht-VAR}
and set~$\zeta_1:=-\Gamma*\widetilde{T}$.
Similarly, we take a function~$\zeta_2$ such
that~$\Delta\zeta_2=\zeta_1$ in~$\Omega$,
and recursively for all~$j\in \{2,\dots,k-2\}$
find functions~$\zeta_j$ such that~$\Delta\zeta_j=\zeta_{j-1}$ in~$\Omega$.
In particular, choosing~$v:=\zeta_{k-2}$, we find that, in~$\Omega$,
\begin{equation}\label{ALMA:4}
\Delta^{k-2}v=\Delta^{k-2}\zeta_{k-2}=
\Delta^{k-3}(\Delta\zeta_{k-2})=\Delta^{k-3}\zeta_{k-3}=\dots=
\Delta \zeta_{1}=T.
\end{equation}
Let also~$w(x):=u(x)-|x|^2v(x)$. In this way,
the claim in~\eqref{ALMA:3} is satisfied. Moreover,
by~\eqref{ALMA:5} and~\eqref{ALMA:4},
\begin{equation}\label{3ALMALPA}
\Delta^{k-1}v=\Delta T=0,
\end{equation}
and also,
recalling in addition~\eqref{iONBDFverpoverufgvvRFDV}, \eqref{ALMALPA}
and~\eqref{1ALMALPA},
\begin{equation}\label{2ALMALPA}\begin{split}
\Delta^{k-1}w\,&=\,
\Delta^{k-1} u -\Delta^{k-1}\big(|x|^2v\big)\\&=\,
\Delta^{k-1} u
-
c_{k-1}\Delta^{k-2}v-d_{k-1} x\cdot\nabla\Delta^{k-2}v-|x|^2\Delta^{k-1}v
\\&=\,
\Delta^{k-1} u
-
c_{k-1} T-d_{k-1} x\cdot\nabla T\\&=\,d_{k-1}\,\left(\frac1{d_{k-1}}
\Delta^{k-1} u
-
\alpha T-  x\cdot\nabla T\right)\\&=\,0.
\end{split}\end{equation}
Thus, the claim in~\eqref{ALMA:2} follows from~\eqref{3ALMALPA}
and~\eqref{2ALMALPA}.
\end{proof}

\begin{corollary}
Let~$k\in\N$, $k\ge1$.
Let~$\Omega$ be a starshaped open subset of~$\R^n$.
Then, every solution~$u\in C^{2k}(\overline\Omega)$ of
\[
{\mbox{$\Delta^ku=0$
in~$\Omega$}}\]
can be written in the form
\begin{equation}\label{13MDpijwefgurqgehwd85867igtru3iefmmxvabndf}
u(x)=\sum_{j=0}^{k-1} |x|^{2j} v^{(j)}(x),
\end{equation}
with~$v^{(0)},\dots,v^{(k-1)}\in C^2(\overline\Omega)$
and harmonic in~$\Omega$.\end{corollary}

\begin{proof}
We argue by induction over~$k$. When~$k=1$ it suffices to define~$v^{(0)}:=u$.
Suppose now that the desired claim holds true for some index~$k-1$, with~$k\ge2$,
and take~$u\in C^{2k}(\overline\Omega)$
such that~$\Delta^ku=0$
in~$\Omega$. We invoke Lemma~\ref{353ALMAMAL}
and we find~$v$, $w\in C^{2(k-1)}(\overline\Omega)$
such that~$\Delta^{k-1}v=\Delta^{k-1}w=0$
in~$\Omega$ and~$u(x)=|x|^2\,v(x)+w(x)$.
We can apply the inductive assumption on~$v$ and~$w$
and obtain
functions~$\alpha^{(0)},\dots,\alpha^{(k-2)},
\beta^{(0)},\dots,\beta^{(k-2)}\in C^2(\overline\Omega)$
and harmonic in~$\Omega$
such that
$$ v(x)=\sum_{j=0}^{k-2} |x|^{2j} \alpha^{(j)}(x)\qquad{\mbox{ and }}\qquad
w(x)=\sum_{j=0}^{k-2} |x|^{2j} \beta^{(j)}(x).$$
In this way, we have that
\begin{eqnarray*}
u(x)&=&
|x|^2\,\sum_{j=0}^{k-2} |x|^{2j} \alpha^{(j)}(x)+\sum_{j=0}^{k-2} |x|^{2j} \beta^{(j)}(x)\\
&=&
\sum_{j=1}^{k-1} |x|^{2j} \alpha^{(j-1)}(x)+\sum_{j=0}^{k-2} |x|^{2j} \beta^{(j)}(x).
\end{eqnarray*}
Hence, the claim in~\eqref{13MDpijwefgurqgehwd85867igtru3iefmmxvabndf} holds true with
\begin{equation*} v^{(j)}:=
\begin{dcases}
\beta^{(0)}&{\mbox{ if }}j=0,\\
\alpha^{(j-1)}+\beta^{(j)}&{\mbox{ if }}j\in\{1,\dots,k-2\},\\
\alpha^{(k-2)}
&{\mbox{ if }}j=k-1,
\end{dcases}\end{equation*}
which completes the proof of the desired result.
\end{proof}

We recall here one of the the main results, motivated by elasticity
theory, that inspired the above discussions 
(see~\cite[pages 15-16]{zbMATH02668224}
and also~\cite{MR1103852} and the references therein):

\begin{corollary}
Let~$g$, $h\in C(\partial B_1)$.
Then, the solution~$u\in C^4(\overline{B_1})$ of
$$ \begin{dcases}\Delta^2 u=0 & {\mbox{ in }}B_1,\\
u=g & {\mbox{ on }}\partial B_1,\\
\partial_\nu u=h& {\mbox{ on }}\partial B_1
\end{dcases}$$
can be explicitly written as
$$ u(x)=-\frac{(1-|x|^2)^2}{4n\,|B_1|}
\int_{\partial B_1}\left[
\frac{g(y)}{|x-y|^{n+2}}\left(
(n-4)|x-y|^2-n(1-|x|^2)
\right)+\frac{2h(y)}{|x-y|^{n}}
\right]\,d{\mathcal{H}}^{n-1}_y.$$
\end{corollary}

\begin{proof} We exploit 
Lemma~\ref{353ALMAMAL} with~$k:=2$
and set~$W:=v+w$. Notice that~$W$ is harmonic and,
in light of~\eqref{ALMA:3},
\begin{equation}\label{sjwi478v5bb5767nin8n8}
u(x)=(|x|^2-1)v(x)+W(x).\end{equation}
In particular, we have that~$g=u=W$ along~$\partial B_1$, whence we can exploit the
Poisson Kernel 
representation in Theorems~\ref{POIBALL1}
and~\ref{POIBALL}
and conclude that, for every~$x\in B_1$,
\begin{equation}\label{PW93-09}
W(x)=\int_{\partial B_1}\frac{(1-|x|^2)\,g(y)}{n\,|B_1|\,|x-y|^n}
\,d{\mathcal{H}}^{n-1}_y.\end{equation}
Additionally, setting~$V(x):=2v(x)+x\cdot\nabla W(x)$,
we have that, on~$\partial B_1$,
$$ h=\partial_\nu u
=\partial_\nu\big((|x|^2-1)v+W\big)
= 2v+\partial_\nu W=2v+x\cdot\nabla W=V.$$
Thus, since~$\Delta V=2\Delta v+
x\cdot\nabla \Delta W+ 2\Delta W=0$,
we can use again the Poisson Kernel 
representation in Theorems~\ref{POIBALL1}
and~\ref{POIBALL} and find that, for every~$x\in B_1$,
$$ V(x)=
\int_{\partial B_1} \frac{(1-|x|^2)\,h(y)}{n\,|B_1|\,|x-y|^n}\,d{\mathcal{H}}_y^{n-1}
$$
As a result, recalling~\eqref{PW93-09},
\begin{eqnarray*}
v(x)&=&\frac{V(x)-x\cdot\nabla W(x)}2\\
&=&\frac12\,\left(\int_{\partial B_1} \frac{(1-|x|^2)\,h(y)}{n\,|B_1|\,|x-y|^n}\,d{\mathcal{H}}_y^{n-1}
-x\cdot\nabla \left(\int_{\partial B_1} \frac{(1-|x|^2)\,g(y)}{n\,|B_1|\,|x-y|^n}\,d{\mathcal{H}}_y^{n-1}\right)
\right)\\&=&\frac1{2\,n\,|B_1|}\left(\int_{\partial B_1} \frac{(1-|x|^2)\,h(y)}{|x-y|^n}\,d{\mathcal{H}}_y^{n-1}
+\int_{\partial B_1} \frac{2|x|^2\,g(y)}{ |x-y|^n}\,d{\mathcal{H}}_y^{n-1}\right.\\&&\qquad\left.
+n\int_{\partial B_1} \frac{x\cdot(x-y)\,(1-|x|^2)\,g(y)}{ |x-y|^{n+2}}\,d{\mathcal{H}}_y^{n-1}
\right)
.\end{eqnarray*}
This and~\eqref{sjwi478v5bb5767nin8n8}, exploiting~\eqref{PW93-09} once more, give that
\begin{eqnarray*}
u(x)&=&
\frac{|x|^2-1}{2\,n\,|B_1|}\left(\int_{\partial B_1} \frac{(1-|x|^2)\,h(y)}{|x-y|^n}\,d{\mathcal{H}}_y^{n-1}
+\int_{\partial B_1} \frac{2|x|^2\,g(y)}{ |x-y|^n}\,d{\mathcal{H}}_y^{n-1}\right.\\&&\qquad\left.
+n\int_{\partial B_1} \frac{x\cdot(x-y)\,(1-|x|^2)\,g(y)}{ |x-y|^{n+2}}\,d{\mathcal{H}}_y^{n-1}
\right)+\int_{\partial B_1}\frac{(1-|x|^2)\,g(y)}{n\,|B_1|\,|x-y|^n}
\,d{\mathcal{H}}^{n-1}_y\\&=&\frac{|x|^2-1}{2\,n\,|B_1|}\left(\int_{\partial B_1} \frac{(1-|x|^2)\,h(y)}{|x-y|^n}\,d{\mathcal{H}}_y^{n-1}
+\int_{\partial B_1} \frac{2|x|^2\,g(y)}{ |x-y|^n}\,d{\mathcal{H}}_y^{n-1}\right.\\&&\qquad\left.
+n\int_{\partial B_1} \frac{x\cdot(x-y)\,(1-|x|^2)\,g(y)}{ |x-y|^{n+2}}\,d{\mathcal{H}}_y^{n-1}
-\int_{\partial B_1}\frac{2g(y)}{|x-y|^n}
\,d{\mathcal{H}}^{n-1}_y\right)\\&=&\frac{|x|^2-1}{2\,n\,|B_1|}\left(\int_{\partial B_1} \frac{(1-|x|^2)\,h(y)}{|x-y|^n}\,d{\mathcal{H}}_y^{n-1}
\right.\\&&\qquad\left.+\int_{\partial B_1} \frac{(|x|^2-1)(2|x-y|^2-nx\cdot(x-y))\,g(y)}{ |x-y|^{n+2}}\,d{\mathcal{H}}_y^{n-1}\right)\\&=&
-\frac{(1-|x|^2)^2}{2\,n\,|B_1|}\left(\int_{\partial B_1} \frac{h(y)}{|x-y|^n}\,d{\mathcal{H}}_y^{n-1}
-\int_{\partial B_1} \frac{(2|x-y|^2-nx\cdot(x-y))\,g(y)}{ |x-y|^{n+2}}\,d{\mathcal{H}}_y^{n-1}\right)
.\end{eqnarray*}
Thus, since, for all~$y\in\partial B_1$,
\begin{eqnarray*} &&-2\big(2|x-y|^2-nx\cdot(x-y)\big)=-4|x-y|^2+2nx\cdot(x-y)\\&&\qquad=
-4|x-y|^2+n(2|x|^2-2x\cdot y)=
-4|x-y|^2+n(|x|^2+|x|^2+|y|^2-2x\cdot y-1)\\&&\qquad=
-4|x-y|^2+n(|x|^2-1+|x-y|^2)=(n-4)|x-y|^2-n(1-|x|^2),
\end{eqnarray*}
we obtain the desired result.
\end{proof}

\begin{figure}
  \centering
  \includegraphics[width=.55\linewidth]{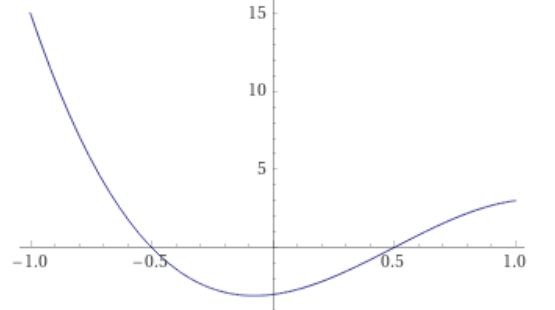}
 \caption{\sl Plot of the function in~\eqref{BILLAL}.}\label{BILLAsoloDItangeFI}
\end{figure}

We stress that the treatment of higher order operators that we give in these notes is far from being
exhaustive and the results presented here should not leave the false impression that
higher order operators can be generally reduced to the case of the Laplacian.
As a matter of fact, in spite of some similarities (such as the existence of
a convenient Green Function for the ball, see~\cite{BOGGIO}),
even for the biharmonic operator 
important structural differences arise.
Just to mention one, we stress that, differently from harmonic functions, biharmonic functions
do not satisfy in general the Maximum Principle: for instance (see~\cite{MR1267051}), as an instructive
example, one can consider, for a suitably small~$\e>0$,
the function
\begin{equation}\label{MASdnhcbor7uia0-23e} u(x,y):=\Big((4 - 3x)(1-x)^2-\e\Big)(
x^2+25 y^2-1)^2\end{equation}
on the
ellipse~$\Omega:=\{ (x,y)\in\R^2 {\mbox{ s.t. }} x^2+25y^2<1\}$.
One can compute explicitly that~$ \Delta^2 u>0$
in~$\Omega$.
Since also~$u=0$ on~$\partial\Omega$, the validity of a Maximum Principle (compare e.g.
with Corollary~\ref{WEAKMAXPLE})
would entail that~$u$
should not change sign in~$\Omega$, while, by inspection~$u(0,0)>0>u(1-
\e,0)$,
thus showing that the Maximum Principle is violated\footnote{It is interesting to remark
that the example in~\eqref{MASdnhcbor7uia0-23e} also satisfies~$\nabla u=0$ along~$\partial\Omega$,
hence even adding such a ``natural'' boundary prescription the Maximum Principle is still violated.

Instead, under the boundary condition~$\Delta u=0$ on~$\partial\Omega$, the Maximum Principle
for the biharmonic operator holds true since it boils down to the case of the Laplacian: namely if
$$ \begin{dcases}
\Delta^2 u\ge0 & {\mbox{ in }}\Omega,\\
u=\Delta u=0 &{\mbox{ on }}\partial\Omega,
\end{dcases}$$
it suffices to define~$v:=\Delta u$, observe that
$$ \begin{dcases}
\Delta v\ge0 & {\mbox{ in }}\Omega,\\
v=0 &{\mbox{ on }}\partial\Omega,
\end{dcases}$$
deduce from
Theorem~\ref{STRONGMAXPLE1} and Corollary~\ref{WEAKMAXPLE} that~$v\le0$ in~$\Omega$
(and in fact either~$v(x)=0$ for all~$x\in\Omega$, or~$v(x)<0$ for all~$x\in\Omega$), infer that
$$ \begin{dcases}
\Delta u\le0 & {\mbox{ in }}\Omega,\\
u=0 &{\mbox{ on }}\partial\Omega,
\end{dcases}$$
apply again Theorem~\ref{STRONGMAXPLE1} and Corollary~\ref{WEAKMAXPLE}
and conclude that~$u\ge0$ in~$\Omega$
(and in fact either~$u(x)=0$ for all~$x\in\Omega$, or~$u(x)>0$ for all~$x\in\Omega$).}
in this case.

As a one-dimensional example, one can also take into account the cubic function
\begin{equation}\label{BILLAL}(-1,1)\ni x\mapsto u(x)=(1 - 2 x)  (1+2 x)(2 x - 3)=
-8 x^3 + 12 x^2 + 2 x - 3\end{equation}
which satisfies~$u(-1)=15>0$, $u(1)=3>0$ and~$\Delta^2 u=u''''=0$;
since~$u(0)=-3<0$,
the Maximum Principle is violated
in this case as well (see Figure~\ref{BILLAsoloDItangeFI}
for a plot of this function).

See~\cite{MR46440, MR1925907, MR2928398} for
other counterexamples
to Maximum Principles for the biLaplace operator.
See~\cite{MR3514727, MR3908840} for counterexamples
to Maximum Principles for the triLaplacian and the quadriLaplacian.
See also~\cite{MR3554425, MR4181195} for several insight on the lack of Maximum Principles
for higher order operators, as well as~\cite{MR1406091, MR2667016} and the references therein for a thorough
discussion of polyharmonic problems.

\chapter{Upper semicontinuity and subharmonicity}

In this chapter we introduce and discuss the
notion of subharmonic functions, that is, roughly speaking
of functions that ``stay below'' the corresponding
harmonic functions with the same boundary data.
The essence of the theory of subharmonic functions is
extremely beautiful, and their structure is pivotal to
understand the behavior of the solutions of elliptic
equations. Nevertheless, if one aims to a theory of subharmonic functions that remains ``closed''
under useful and natural operations, then a number of subtleties arise and a convenient functional framework becomes necessary.
For this, in the forthcoming section,
we recall the notion of upper semicontinuous functions, which will provide
a convenient scenario for dealing with
subharmonic functions in Section~\ref{SEMICSEC}.

\section{Semicontinuous functions}

An often useful ancillary tool of mathematical analysis
is provided by
semicontinuous functions. This is a weaker notion
than continuity, since it allows ``discontinuities in one direction''.
Roughly speaking, a function is upper (respectively, lower) semi-continuous at a point
if the values of the function nearby the point are not much higher (respectively, lower) than the value attained precisely
at the point. Differently than continuous functions, 
upper (respectively, lower) semi-continuous functions may possess values nearby the point that
are much lower (respectively, higher) than the value
at the point. The notion of semicontinuity will be exploited in Section~\ref{SEMICSEC}
to deal with subharmonic and superharmonic functions. Before that,
we will introduce the formal definition of semicontinuity, discuss
the main properties of semicontinuous functions and develop
some intuition of this setting
(for this, we will mainly focus on upper semicontinuity, since
the case of lower semicontinuity is then easily obtained by a sign change).

\begin{definition}
Given an open set~$\Omega\subseteq\R^n$ and
a function~$u:\Omega\to\R\cup\{-\infty\}$, we say that~$u$ is upper semicontinuous\index{upper semicontinuous} in~$\Omega$
if
\begin{equation}\label{likm8ijk9-0eofkjjgjjn} \limsup_{x\to x_0} u(x)\le u(x_0),\end{equation}
for each~$x_0\in\Omega$. Also, we say that~$u:\Omega\to\R\cup\{+\infty\}$ is lower semicontinuous\index{lower semicontinuous} in~$\Omega$
if the function~$-u$ is upper semicontinuous.\end{definition}

As an example, we have that the function
$$ u(x):=\begin{dcases}
1 & {\mbox{ if }}x\ge0,\\
0 & {\mbox{ if }}x<0
\end{dcases}$$
is lower semicontinuous.
As it is easily seen, a function is continuous in~$\Omega$ if and only if it is both
upper and lower semicontinuous.
Also, upper and lower semicontinuous functions are measurable
(see e.g.~\cite[Corollary~4.16]{MR3381284}). Moreover,
semicontinuity is preserved under suitably monotone limits,
as remarked in the following result:

\begin{lemma}\label{65co99865bSdsaJ-04}
The pointwise
limit of a decreasing sequence
of upper semicontinuous functions is upper semicontinuous.
\end{lemma}

\begin{proof} Let~$u_k$ be a decreasing sequence
of upper semicontinuous functions in some domain~$\Omega$
and let~$u$ be its pointwise limit.
Then, for every~$j\in\N$ and~$x\in\Omega$, we have that~$u_{k+j}(x)\le u_k(x)$, therefore
$$ u(x)=\lim_{j\to+\infty} u_{k+j}(x)\le u_k(x).$$
As a result, using the upper semicontinuity of~$u_k$, for each~$x_0\in\Omega$,
\begin{equation*}
\limsup_{x\to x_0}u(x)\le\limsup_{x\to x_0} u_k(x)\le u_k(x_0).
\end{equation*}Thus, taking the limit as~$k\to+\infty$,
\begin{equation*}
\limsup_{x\to x_0}u(x)\le\lim_{k\to+\infty} u_k(x_0)=u(x_0).\qedhere
\end{equation*}
\end{proof}

We stress that the decresing monotonicity assumption on the sequence of
functions in Lemma~\ref{65co99865bSdsaJ-04} cannot be removed: for instance,
one can consider the sequence of function~$u_k(x):=-e^{-k|x|^2}$
and observe that it converges pointwise to~$u(x):=-\chi_{\{0\}}(x)$ which is
not upper semicontinuous.

Decreasing sequences, as the ones in Lemma~\ref{65co99865bSdsaJ-04}
can indeed be a useful tool to approximate upper semicontinuous functions,
as explicitly observed in the next result:

\begin{lemma}\label{ANS-BOSDS9id2344we}
If~$u$ is upper semicontinuous in~$\Omega$ and bounded above,
then there exists a decreasing sequence of Lipschitz
continuous functions
that converges to~$u$ at every point of~$\Omega$.
\end{lemma}

\begin{proof} 
%%%%% By restricting to connected components, we can suppose that~$\Omega$ is connected. 
We can assume that~$u$ is not identically~$-\infty$
(otherwise we consider the decreasing sequence of constant functions~$u_j:=-j$
and we are done).
We also consider the extension of~$u$ to the whole of~$\R^n$ given, for every~$x\in\R^n$, by
$$ \overline{u}(x):=\begin{dcases}
u(x) & {\mbox{ if }}x\in\Omega,\\
\displaystyle\limsup_{\Omega\ni y\to x}u(y)& {\mbox{ if }}x\in\partial\Omega,\\
-\infty& {\mbox{ if }}x\in\R^n\setminus\overline{\Omega}.
\end{dcases}$$
We remark that
\begin{equation}\label{poisemi}
{\mbox{$\overline{u}$ is upper semicontinuous in~$\R^n$.}}\end{equation}
Indeed, $\overline{u}$ is obviously upper semicontinuous in~$\Omega$
and in~$\R^n\setminus\overline{\Omega}$.
Furthermore, we consider a point~$x_0\in\partial\Omega$ and
we take a sequence~$x_k\to x_0$ as~$k\to+\infty$.
Up to a subsequence, we can suppose that the limit of~$\overline u(x_k)$ exists and
realizes the $\limsup$ at the point~$x_0$, that is
$$ \lim_{k\to+\infty}\overline u(x_k)=\limsup_{x\to x_0}\overline u(x).$$
Hence, to prove~\eqref{poisemi},
we aim at showing that
\begin{equation}\label{poisemi-PRE}
\lim_{k\to+\infty}\overline u(x_k)\le\overline u(x_0).
\end{equation}
To accomplish this goal,
we distinguish two cases: either~$x_k\in\partial\Omega$ for finitely many indices~$k$,
or~$x_k\in\partial\Omega$ for infinitely many indices~$k$.

Let us first deal with the case in which~$x_k\in\partial\Omega$ for finitely many indices~$k$.
In this situation, if~$x_k$ belonged to~$\Omega$ only for finitely many indices~$k$ we would have that
for all~$k$ large enough~$x_k\in\R^n\setminus\overline{\Omega}$ and thus~$u(x_k)=-\infty$,
making the claim in~\eqref{poisemi-PRE} obvious. This allows us to suppose that~$x_k$ belongs to~$\Omega$ for infinitely many indices~$k$, that we denote by~$k_j$. In this setting, we conclude that
$$ \limsup_{x\to x_0}\overline u(x)=\lim_{k\to+\infty}\overline u(x_k)=\lim_{j\to+\infty}\overline u(x_{k_j})
=\lim_{j\to+\infty}u(x_{k_j})\le\limsup_{\Omega\ni y\to x_0}u(y)=\overline u(x_0),$$
which is~\eqref{poisemi-PRE}.

Hence it remains to consider the case
in which~$x_k\in\partial\Omega$ for infinitely many indices~$k$.
In this case, by the definition of~$\overline{u}$,
we have that there exists~$y_k\in B_{1/k}(x_k)\cap\Omega$
such that~$\overline{u}(x_k)\le u(y_k)+\frac1k$. Also,
$$ |y_k-x_0|\le |y_k-x_k|+|x_k-x_0|\le \frac1k +|x_k-x_0|,$$
which gives that~$y_k\to x_0$ as~$k\to+\infty$.
As a consequence,
$$ \overline{u}(x_0)=\limsup_{\Omega\ni y\to x}u(y)
\ge \lim_{k\to+\infty}u(y_k)\ge
\lim_{k\to+\infty}\left( \overline{u}(x_k)-\frac1k\right)
=\lim_{k\to+\infty}\overline{u}(x_k).
$$
These considerations 
complete the proof of~\eqref{poisemi-PRE}, and thus of~\eqref{poisemi}.

Moreover, $\overline{u}$ is bounded above.
In particular,
\begin{equation}\label{EMmoKS3}
\overline{u}\le M,\end{equation}
for some~$M\ge0$.

For every~$x\in\R^n$, we utilize a ``sup-convolution method'' by setting
\begin{equation}\label{KS:p34rid} u_j(x):=\sup_{y\in\R^n} \overline{u}(y)-j\,|x-y|.\end{equation}
We observe that, for every~$x$, $y\in\R^n$,
$$ \overline{u}(y)-(j+1)\,|x-y|
=-|x-y|+\big(\overline{u}(y)-j\,|x-y|\big)
\le -|x-y|+u_j(x)\le u_j(x),$$
hence~$u_{j+1}(x)\le u_j(x)$ and thus~$u_j$ is a decreasing sequence.

Moreover, picking~$y:=x$ in~\eqref{KS:p34rid},
\begin{equation}\label{RINgindoogie4-2}
u_j(x)\ge\overline{u}(x)\qquad{\mbox{ for every }}x\in\R^n.
\end{equation}

In addition, given two points~$x$, $x_0\in\R^n$, for every~$y\in\R^n$,
$$ u_j(x_0)\ge\overline{u}(y)-j\,|y-x_0|\ge\overline{u}(y)-j\,|y-x|-j\,|x-x_0|.$$
Taking the supremum in~$y$, we thus see that
$$ u_j(x_0)\ge u_j(x)-j\,|x-x_0|,$$
whence, exchanging the roles of~$x$ and~$x_0$,
$$ |u_j(x)-u_j(x_0)|\le j\,|x-x_0|.$$
This proves that~$u_j$ is Lipschitz continuous.

Now we show that
\begin{equation}\label{RINgindoogie4-1}
\lim_{j\to+\infty}u_j(x)=\overline{u}(x)\qquad{\mbox{ for every }}x\in\R^n.
\end{equation}
By~\eqref{RINgindoogie4-2}, it suffices to show that
\begin{equation*}
\limsup_{j\to+\infty}u_j(x)\le\overline{u}(x)\qquad{\mbox{ for every }}x\in\R^n.
\end{equation*}
Suppose not. Then, there exist~$x\in\R^n$ and~$a>0$ such that
\begin{equation}\label{T0oTHSismfg}
\limsup_{j\to+\infty}u_j(x)>-\infty
\end{equation}
and
\begin{equation}\label{T0oTHSismf}
\limsup_{j\to+\infty}u_j(x)\ge a+\overline{u}(x).
\end{equation}
For each~$j\in\N$ we use the definition
in~\eqref{KS:p34rid} and take~$y_j\in\R^n$ such that
\begin{equation}\label{KSMD:034ik3rjwjejdd} u_j(x)\le \overline{u}(y_j)-j\,|x-y_j|+\frac{a}{2}.\end{equation}
This and~\eqref{T0oTHSismf} lead to
\begin{equation}\label{BANlbokfonrsc5}
\limsup_{j\to+\infty}\Big(
\overline{u}(y_j)-j\,|x-y_j|\Big)\ge \frac{a}2+\overline{u}(x).
\end{equation}
Now we exploit
the upper semicontinuity of~$\overline{u}$ and we take~$\rho>0$
sufficiently small such that if~$y\in B_\rho(x)$ then~$\overline{u}(y)\le\overline{u}(x)+\frac{a}{4}$.
We stress that~$y_j\not\in B_\rho(x)$, otherwise
$$ \overline{u}(y_j)-j\,|x-y_j|\le 
\overline{u}(y_j)\le
\overline{u}(x)+\frac{a}{4},$$
in contradiction with~\eqref{BANlbokfonrsc5}.

But then, recalling~\eqref{EMmoKS3},
$$ \overline{u}(y_j)-j\,|x-y_j|\le M-j\rho\to-\infty\qquad{\mbox{as }}\,j\to+\infty.$$
This and~\eqref{KSMD:034ik3rjwjejdd} give that
$$ \lim_{j\to+\infty} u_j(x)=-\infty,$$
which is in contradiction with~\eqref{T0oTHSismfg}.
This establishes~\eqref{RINgindoogie4-1},
which in turn yields the desired result.
\end{proof}

Also, upper semicontinuous functions attain their maximum in compact sets:

\begin{lemma}\label{ABOABP}
If~$u$ is upper semicontinuous in~$\Omega$
and~$K\Subset\Omega$ is compact, then the maximum of~$u$
in~$K$ is attained and it is finite.
\end{lemma}

\begin{proof}
One takes a maximizing sequence and proceeds as in the classical
proof of Weierstra{\ss}' Extreme Value Theorem, exploiting~\eqref{likm8ijk9-0eofkjjgjjn}.
\end{proof}

By combining Lemmata~\ref{ANS-BOSDS9id2344we} and~\ref{ABOABP}, we obtain:

\begin{corollary}\label{ANS-BOSDS9id2344we-BIS}
Let~$u$ be upper semicontinuous in~$\Omega$. Let~$\Omega'\Subset\Omega$ be an open set.
Then there exists a decreasing sequence of Lipschitz
continuous functions
that converges to~$u$ at every point of~$\Omega'$.
\end{corollary}

\begin{figure}
  \centering
  \includegraphics[width=.55\linewidth]{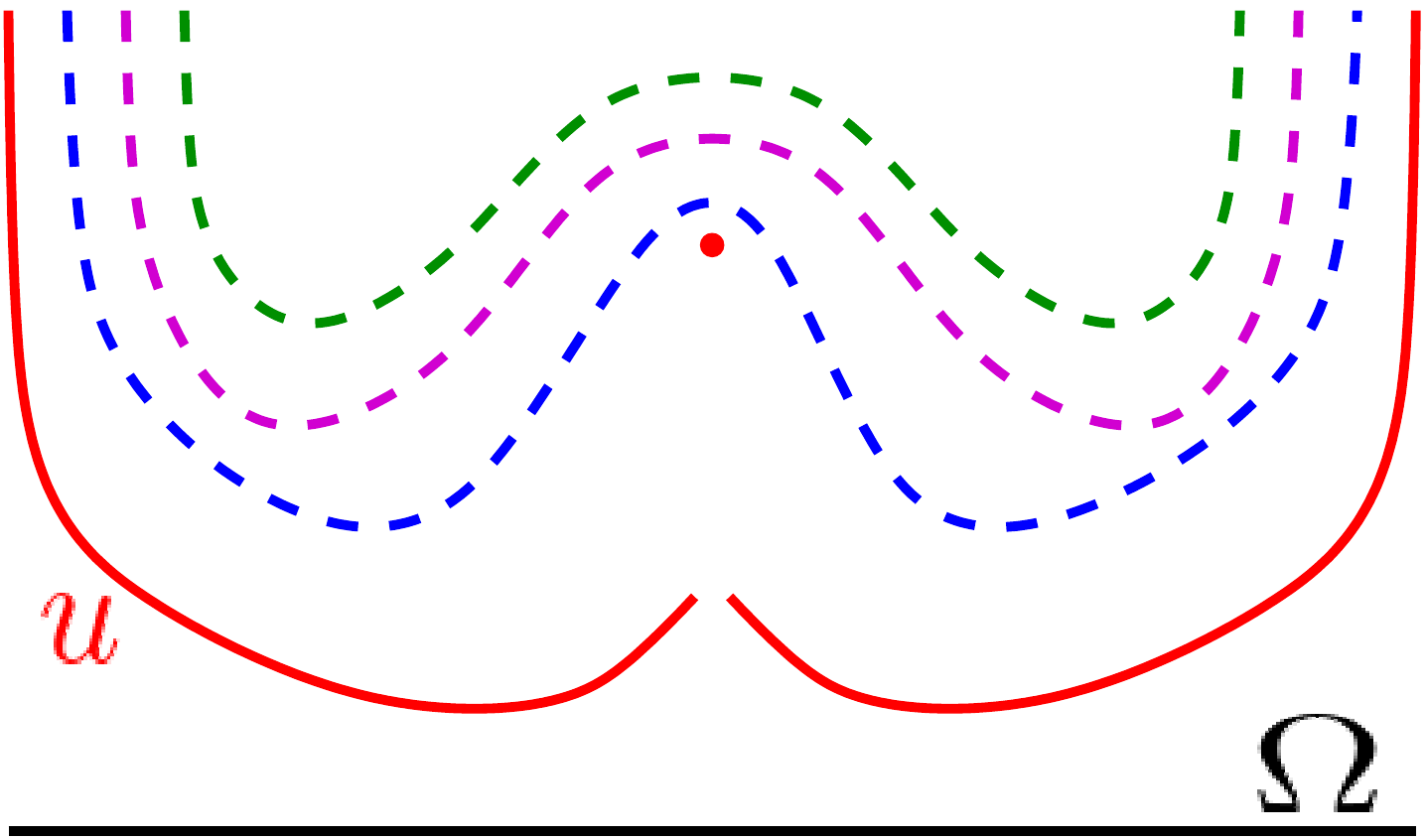}
 \caption{\sl Regular approximation from above of an upper semicontinuous function.}\label{soloDItangeFI}
\end{figure}

See Figure~\ref{soloDItangeFI} for a visualization of a regular and decreasing approximation
of an upper semicontinuous function (notice that one cannot obtain a similar approximation from below).

We can thus characterize upper semicontinuous functions as follows:

\begin{corollary}\label{754PCHAYPP}
The following conditions are equivalent:
\begin{itemize}
\item[(i).] $u$ is upper semicontinuous in~$\Omega$.
\item[(ii).]  For each~$\alpha\in[-\infty,+\infty)$ the supergraph~$\{u\ge\alpha\}$ is closed in~$\Omega$.
\item[(iii).] For every open set~$\Omega'\Subset\Omega$, there exists a decreasing sequence of functions~$u_j\in C^\infty(\Omega')$
that converges to~$u$ at every point of~$\Omega'$.\end{itemize}
\end{corollary}

\begin{proof} We check that~(i) and~(iii) are equivalent.
Indeed, if~(i) holds true, given~$\Omega'\Subset\Omega$, we consider~$\Omega''$ such that~$\Omega'\Subset\Omega''\Subset\Omega$
and we
utilize
Corollary~\ref{ANS-BOSDS9id2344we-BIS}
to find a decreasing sequence of Lipschitz
continuous functions~$v_j$
that converges to~$u$ at every point of~$\Omega''$.
We define~$w_j:=v_j+\frac1{2^j}$ and we consider~$\tau\in C^\infty_0
(B_1,\,[0,+\infty))$
with~$\int_{B_1}\tau(x)\,dx=1$. Given~$\e>0$, we let~$\tau_\e(x):=\frac1{\e^n}\tau\left(\frac{x}{\e}\right)$
and~$w_{j,\e}:=w_j*\tau_\e$.
We stress that~$w_{j,\e}\in C^\infty(\Omega')$ and~$w_{j,\e}$ converges to~$w_j$
uniformly in~$\Omega'$ as~$\e\searrow0$
(see e.g.~\cite[Theorem~9.8]{MR3381284}). Consequently, we can find~$\e_j$ such that~$\|
w_{j,\e_j}-w_j\|_{L^\infty(\Omega')}\le\frac1{2^{j+2}}$.
In this way, setting~$u_j:=w_{j,\e_j}$,
we have that, for every~$x\in\Omega'$,
$$ \lim_{j\to+\infty} u_j(x)=\lim_{j\to+\infty} \Big(\big(w_{j,\e_j}-w_j(x)\big)+w_j(x)\Big)
=\lim_{j\to+\infty}w_j(x)=\lim_{j\to+\infty}\left(
v_j(x)+\frac1{2^j}\right)=\lim_{j\to+\infty}v_j(x)=u(x),$$
that is~$u_j$ converges to~$u$ at each point of~$\Omega'$. It only remains to prove that the convergence
occurs in a monotonically decreasing way. To this end, we point out that
\begin{equation}\label{POmbwfyhintpel}
u_j=w_{j,\e_j}\le w_j+\frac1{2^{j+2}}=v_j+\frac1{2^j}+\frac1{2^{j+2}}=v_j+\frac{5}{2^{j+2}}
\end{equation}
and
$$ u_j=w_{j,\e_j}\ge
w_j-\frac1{2^{j+2}}
=v_j+\frac1{2^j}-\frac1{2^{j+2}}=
v_j+\frac{3}{2^{j+2}}
\ge
v_{j+1}+\frac3{2^{j+2}}.$$ 
We rephrase the latter inequality as
$$ u_{j-1}\ge v_{j}+\frac3{2^{j+1}},$$
that is
$$ v_j\le u_{j-1}-\frac{3}{2^{j+1}}.$$
As a byproduct of this and~\eqref{POmbwfyhintpel},
$$ u_j\le v_j+\frac{5}{2^{j+2}}\le u_{j-1}-\frac{3}{2^{j+1}}+\frac{5}{2^{j+2}}=
u_{j-1}-\frac{1}{2^{j+2}}<u_{j-1},
$$
confirming thereby that~$u_j$ is a decreasing sequence. This shows that~(i) implies~(iii).\medskip

Moreover, if~(iii) holds true, then we can apply Lemma~\ref{65co99865bSdsaJ-04}
and conclude that~$u$ is upper semicontinuous.
This argument shows that~(iii) implies~(i).\medskip

Let us now check that~(i) implies~(ii). 
To this end, given~$\alpha\in[-\infty,+\infty)$,
we distinguish two cases. If~$\alpha=-\infty$, then~$\{u\ge\alpha\}=\Omega$,
which is closed in~$\Omega$, and then~(ii) is fulfilled in this case.

If instead~$\alpha\in\R$,
we
take a sequence~$x_k\in\Omega$ such that~$u(x_k)\ge\alpha$
and suppose that~$x_k\to x_\star\in\Omega$ as~$k\to+\infty$. Then, by~\eqref{likm8ijk9-0eofkjjgjjn},
$$ \alpha\le\lim_{k\to+\infty} u(x_k)\le \limsup_{x\to x_\star} u(x)\le u(x_\star),$$
hence~$x_\star\in\{u\ge\alpha\}$, thereby establishing the validity of~(ii).\medskip

Thus, to complete the proof of Corollary~\ref{754PCHAYPP},
we now check that~(ii) implies~(i). For this, let~$x_\star\in\Omega$ and~$\alpha:=u(x_\star)\in[-\infty,+\infty)$.
Pick also a sequence~$x_k\in\Omega$ such that~$x_k\to x_\star$ as~$k\to+\infty$ and
$$ \lim_{k\to+\infty} u(x_k)=\limsup_{x\to x_\star} u(x).$$
Hence, to check that~(i) holds true, we need to show that
\begin{equation}\label{INVIO}
\lim_{k\to+\infty} u(x_k)\le u(x_\star).
\end{equation}
Suppose not, and distinguish whether~$\alpha\in\R$
or~$\alpha=-\infty$. Suppose first that~$\alpha\in\R$: then, if~\eqref{INVIO} is violated,
we have that~$u(x_k)\ge u(x_\star)+b=\alpha+b$, for some~$b>0$ and infinitely many indices~$k$.
In particular~$x_k\in\{u\ge\alpha+b\}$ and thus, by~(ii), also~$x_\star\in\{u\ge\alpha+b\}$. This gives that~$
\alpha=u(x_\star)\ge\alpha+b$, which is a contradiction.

Suppose now that~$\alpha=-\infty$. If~\eqref{INVIO} did not hold, then~$u(x_k)\ge d$, for some~$d\in\R$
and infinitely many indices~$k$.
In particular~$x_k\in\{u\ge d\}$ and thus, by~(ii), also~$x_\star\in\{u\ge d\}$. This gives that~$
-\infty=\alpha=u(x_\star)\ge d$, which is a contradiction.
The proof of~\eqref{INVIO}, and thus of
Corollary~\ref{754PCHAYPP}, is thereby complete.
\end{proof}

We emphasize that the convergence in Corollary~\ref{754PCHAYPP}(iii)
is in the everywhere sense, not only almost everywhere: as an example,
one can consider the decreasing sequence of functions~$u_j(x):=\frac1j$
which converge almost everywhere to~$ u(x):=-\chi_{\{0\}}(x)$, but~$u$ is not upper semicontinuous. 

We also point out
that upper semicontinuity is preserved by taking infima (as well as finitely many maxima), as stated in the next two
results:

\begin{lemma}\label{KSMD:0orkegm}
The maximum of finitely many upper semicontinuous functions
is upper semicontinuous.\end{lemma}

\begin{proof} Let~$u_1,\dots,u_N$ be upper semicontinuous
in a given~$\Omega\subseteq\R^n$.
We aim at showing that~$u:=\max\{u_1,\dots,u_N\}$ is also upper semicontinuous.
Suppose not. Then there exist~$x_\star\in\Omega$,
$a>0$ and a sequence of points~$x_k$ converging to~$x_\star$ as~$k\to+\infty$
such that~$u(x_k)\ge u(x_\star)+a$. Notice in particular that~$x_k\ne x_\star$,
hence~$\{x_k\}_{k\in\N}$ is an infinite set. As a result, being the set~$\{1,\dots,N\}$
finite, there must exist~$i\in\{1,\dots,N\}$ such that~$u_i(x_k)=
u(x_k)\ge u(x_\star)+a$ for infinitely many indices~$k$. Consequently~$u_i(x_k)
\ge u_i(x_\star)+a$ for infinitely many indices~$k$ and, for this reason,
$$ \limsup_{x\to x_\star}u_i(x)\ge\limsup_{k\to+\infty}u_i(x_k)\ge
u_i(x_\star)+a,$$
and this is in contradiction with the upper semicontinuity of~$u_i$.
\end{proof}

\begin{lemma}\label{MS:MINSINDFF}
The infimum of a family of upper semicontinuous functions is upper semicontinuous.
\end{lemma}

\begin{proof} Let~$\Omega$ be an open set.
We consider a set of indices~${\mathcal{J}}$ and, for each~$j\in{\mathcal{J}}$, an upper semicontinuous function~$u_j:\Omega\to[-\infty,+\infty)$. We let~$u(x):=\inf_{j\in{\mathcal{J}}}u_j(x)$ and we check that~$u$ is upper semicontinuous. For this,
we observe that, for every~$\alpha\in[-\infty,+\infty)$,
\begin{equation*}
\{u\geq\alpha\}=\bigcap_{j\in{\mathcal{J}}}\{ u_j\geq\alpha\}.
\end{equation*}
Since, by Corollary~\ref{754PCHAYPP}(ii), the set~$\{ u_j\geq\alpha\}$
is closed, we thus infer that also~$\{u\geq\alpha\}$ is closed.
Thus, using again Corollary~\ref{754PCHAYPP}(ii), we obtain
that~$u$ is upper semicontinuous.
\end{proof}

Regarding Lemmata~\ref{KSMD:0orkegm}
and~\ref{MS:MINSINDFF}, we notice that
it is not true in general that the 
supremum of an infinite family of upper semicontinuous functions is upper semicontinuous:
as a counterexample, one can take~$\Omega:=(-1,1)$
and~$u_j(x):=-e^{-\frac{|x|^2}{j}}$: indeed, in this case
$$ \sup_j u_j(x)= -\chi_{\{0\}}(x),$$
which is not upper semicontinuous.

We also mention that Lemma~\ref{MS:MINSINDFF} has a natural counterpart,
leading to another characterization of upper semicontinuous functions
(in addition to the ones described in 
Corollary~\ref{754PCHAYPP}):

\begin{corollary}
The following conditions are equivalent:
\begin{itemize}
\item[(i).] $u$ is upper semicontinuous in~$\Omega$.
\item[(ii).] $u$ is locally the infimum of a family of upper
semicontinuous functions.
\item[(iii).] $u$ is locally the infimum of a family of $C^\infty$ functions.\end{itemize}
\end{corollary}

\begin{proof} We show that~(iii) implies~(ii) which implies~(i)
which implies~(iii). Indeed: since~$C^\infty$ functions
are upper semicontinuous, obviously~(iii) implies~(ii).
Moreover, as a consequence of Lemma~\ref{MS:MINSINDFF},
we have that~(ii) implies~(i). Finally, in view of Corollary~\ref{754PCHAYPP}(iii),
we know that~(i) implies~(iii).
\end{proof}

It is also useful
to define the upper semicontinuous envelope of
a function~$u:\Omega\to[-\infty,+\infty)$ that is locally bounded above\footnote{The
assumption that~$u$ is locally bounded above \label{SKM:SIJKDM023o3rtetergsrib}
is convenient to ensure that the set in the right hand side of~\eqref{S-asle0o4}
is nonempty, namely that there exists at least one upper semicontinuous function
above~$u$ (for instance, the constant function equal to the supremum of~$u$).} by
\begin{equation}\label{S-asle0o4} U_u (x):=\inf\big\{ v(x),
{\mbox{ among all $v$ which are upper semicontinuous on
$\Omega$ and $v\ge u$}} \big\}.
\end{equation}
We point out that, by construction,
\begin{equation}\label{Semxeeq}
U_u \ge u.\end{equation}
Then, we have:

\begin{lemma}\label{ENCEBENFD}
The upper semicontinuous envelope of
a locally bounded above function~$u$ is an upper semicontinuous function.

Moreover, $u$ coincides with its
upper semicontinuous envelope if and only if~$u$
is upper semicontinuous.
\end{lemma}

\begin{proof} Let~$x_\star\in\Omega$. Let~$x_k\in\Omega$ be such that~$x_k\to x_\star$ as~$k\to+\infty$. Given~$\e>0$, we take~$v_\e$ which is upper semicontinuous on $\Omega$, $v_\e\ge u$
and~$v_\e(x_\star)\le U_u(x_\star)+\e$.
In this way, since~$U_u\le v_\e$, we have that
$$ \lim_{k\to+\infty}U_u(x_k)\le
\lim_{k\to+\infty} v_\e(x_k)\le
v_\e(x_\star)\le U_u(x_\star)+\e.$$
Sending~$\e\searrow0$, we get that
$$ \lim_{k\to+\infty}U_u(x_k)\le U_u(x_\star)$$
and this shows that~$U_u$ is upper semicontinuous.

As a byproduct, one also has that if~$u$ coincides with~$U_u$ then~$u$ is upper semicontinuous.

Now suppose that~$u$ is upper semicontinuous.
Then, one can choose~$v:=u$ in the definition~\eqref{S-asle0o4} 
of upper semicontinuous envelope, thus concluding that
\begin{equation*} U_u(x)\le
\inf\big\{ v(x), {\mbox{ among all $v$ which are upper
semicontinuous on $\Omega$ and $v\ge u$}}\big\}\le u(x)
\end{equation*}
This and~\eqref{Semxeeq} entail that~$U_u=u$, as desired.
\end{proof}

Alternative definitions of the upper semicontinuous envelope are as follows:

\begin{lemma}\label{Snrdzampqoraur03rt}
Let~$u:\Omega\to[-\infty,+\infty)$ be locally bounded
and let~$U_u$ be its upper semicontinuous envolope as in~\eqref{S-asle0o4}.
Then, for each~$x\in\Omega$,
$$ U_u(x)=\inf_{r>0}\sup_{\Omega\cap B_r(x)}u=
\max\left\{ \limsup_{y\to x}u(y),\,u(x)\right\}.$$
\end{lemma}

\begin{proof} First, we claim that
\begin{equation}\label{dvabvansyoti}
U_u(x)\le \inf_{r>0}\sup_{\Omega\cap B_r(x)}u.
\end{equation}
For this, let~$\e>0$ and take~$r_\e>0$ such that
\begin{equation}\label{comedy8di234P}\e+\inf_{r>0}\sup_{\Omega\cap B_r(x)}u
\ge\sup_{\Omega\cap B_{r_\e}(x)}u.\end{equation}
For all~$y\in\Omega$, let also
$$ v(y):=\begin{dcases}
\sup_{\Omega\cap B_{r_\e}(x)}u & {\mbox{ if }}y\in\Omega\cap B_{r_\e}(x),\\
\sup_\Omega u& {\mbox{ if }}y\in\Omega\setminus B_{r_\e}(x).
\end{dcases}$$
We observe that~$v$ is upper semicontinuous and it is larger than
or equal to~$u$,
therefore, by~\eqref{S-asle0o4}, we have that
$$U_u(x)\le v(x)=\sup_{\Omega\cap B_{r_\e}(x)}u.$$
{F}rom this and~\eqref{comedy8di234P} we conclude that
$$ U_u(x)\le\e+\inf_{r>0}\sup_{\Omega\cap B_r(x)}u,$$
whence~\eqref{dvabvansyoti} follows by sending~$\e\searrow0$.

Now we show that
\begin{equation}\label{LSnalassednwoe94}
\inf_{r>0}\sup_{\Omega\cap B_r(x)}u\le
\max\left\{ \limsup_{y\to x}u(y),\,u(x)\right\}
\end{equation}
For this, we argue by contradiction and suppose that
there exists~$a>0$ such that, for every~$k\in\N$,
$$ \sup_{\Omega\cap B_{1/k}(x)}u\ge a+
\max\left\{ \limsup_{y\to x}u(y),\,u(x)\right\}.$$
Let also~$x_k\in\Omega\cap B_{1/k}(x)$ be such that~$
u(x_k)\ge\sup_{\Omega\cap B_{1/k}(x)}u-\frac{a}2$. In this way,
\begin{equation}\label{ukx03} u(x_k)\ge \frac{a}2+
\max\left\{ \limsup_{y\to x}u(y),\,u(x)\right\}.\end{equation}
In particular, we have that~$u(x_k)\ge \frac{a}2+u(x)$, which shows that~$x_k\ne x$,
and consequently, recalling~\cite[Definition~2.6]{MR2016192},
since~$|x_k-x|\le\frac1k$,
$$ \limsup_{k\to+\infty}u(x_k)\le\limsup_{y\to x}u(y).$$
This provides a contradiction with~\eqref{ukx03},
thus proving~\eqref{LSnalassednwoe94}.

Now let~$v$ be upper semicontinuous and $v\ge u$. Then,
$$ v(x)\ge\limsup_{y\to x}v(y)\ge\limsup_{y\to x}u(y).$$
Taking the infimum over~$v$, we thereby deduce that~$U_u(x)\ge \limsup_{y\to x}u(y)$.
This and~\eqref{Semxeeq} yield that
$$ \max\left\{ \limsup_{y\to x}u(y),\,u(x)\right\}\le U_u(x).$$
Combining this inequality with~\eqref{dvabvansyoti} and~\eqref{LSnalassednwoe94}
the desired result plainly follows.
\end{proof}

As a side remark, we point out that if two functions are locally bounded
and coincide almost everywhere,
it is not necessarily
For instance let~$u$ be the function identically zero on~$\R$,
and~$\widetilde u:=\chi_{\Q}$. Then,
by Lemma~\ref{Snrdzampqoraur03rt} and the density of the rationals, for every~$x\in\R$,
$$ U_{\widetilde u}(x)\ge
\limsup_{y\to x}\chi_{\Q}(y)=1>0=U_u(x).$$

\section{Subharmonic and superharmonic functions}\label{SEMICSEC}

We now consider a class of functions
that extends the one of harmonic functions introduced in Definition~\ref{KSMD:PSDLKNIE0sdf}.
Roughly speaking, harmonic
functions are exceptionally special since
they lie precisely in the kernel of the Laplacian and
they satisfy an equality with their means.
Of course, an equality is made by the common validity
of two inequalities with different signs: hence 
we will consider here functions that satisfy one
of the opposite inequalities characterizing harmonic functions
(but maybe not both of them). On the one hand,
these functions carry ``half of the information''
with respect to harmonic functions. On the other hand,
we will see that these functions are still quite special
and enjoy a number of important properties (and moreover
they play a crucial role in the development of the theory
of elliptic partial differential equations, especially
in view of a one-sided Maximum Principle that these functions
satisfy).

We start with an equivalence result that
will characterize the main properties enjoyed by this
pivotal class of functions.

\begin{theorem}\label{SUBHS-EQ}
Let~$\Omega\subseteq\R^n$ be an open set and~$u: \Omega\to\R\cup\{-\infty\}$. The following conditions are equivalent:
\begin{itemize}
\item[(i).] 
The function~$u$ is upper semicontinuous and for every~$x\in\Omega$
and~$r>0$ with~$B_r(x)\Subset\Omega$ we have that
\begin{equation*} u(x)\le\fint_{B_r(x)} u(y)\,dy.\end{equation*}
\item[(ii).] The function~$u$ is upper semicontinuous and for every~$x\in\Omega$
there exists~$r_0>0$ with~$
B_{r_0}(x)\Subset\Omega$ such that, for each~$r\in(0,r_0]$,
\begin{equation*} u(x)\le\fint_{B_{r}(x)} u(y)\,dy.\end{equation*}
\item[(iii).] The function~$u$ is upper semicontinuous and for every~$x\in\Omega$
and~$r>0$ with~$B_r(x)\Subset\Omega$ we have that
\begin{equation*} u(x)\le\fint_{\partial B_r(x)} u(y)\,d{\mathcal{H}}^{n-1}_y.\end{equation*}
\item[(iv).] The function~$u$ is upper semicontinuous and for every~$x\in\Omega$
there exists~$r_0>0$ with~$
B_{r_0}(x)\Subset\Omega$ such that, for each~$r\in(0,r_0]$,
\begin{equation*} u(x)\le\fint_{\partial B_{r}(x)} u(y)\,d{\mathcal{H}}^{n-1}_y.\end{equation*}
\item[(v).] 
Either~$u$ is constantly\footnote{Notice that we are allowing~$u$
to be possibly constantly~$-\infty$ in~$\Omega$, though of course
this provides a trivial case in
Theorem~\ref{SUBHS-EQ}. Other settings
are possible as well, see e.g.~\cite[Definition~3.1.2]{MR1801253}.
A slight abuse of notation is used in~(v) and~(vi) since, in the usual
sense of Lebesgue spaces, the function~$u$ in~(v) and~(vi)
is defined only up to sets
of null measure. More precisely,
when needed, we will be able to define~$u$
as the limit of its volume average.
See in particular~\eqref{AKSM:PSDKMMD:034r-1}
and~\eqref{pi9ltar} below. The setting in~(v) is sometimes
denoted by ``weakly subharmonic functions'', see in particular~\cite{MR0107746, MR177186, MR2906766}.
See also~\cite{MR123096} for a parabolic analogue.} equal to~$-\infty$ in a connected
component of~$\Omega$, or~$u\in L^1_{\rm loc}(\Omega)$ and for every~$\varphi\in
C^\infty_0(\Omega,\,[0,+\infty))$,
\begin{equation}\label{vuuvrdcsbaP0lqefo} \int_\Omega u(x)\,\Delta\varphi(x)\,dx\ge0.\end{equation}
\item[(vi).] Either~$u$ is constantly equal to~$-\infty$ in a connected
component of~$\Omega$, or~$u\in L^1_{\rm loc}(\Omega)$ and 
for every~$x\in\Omega$
there exists~$r_0>0$ such that
for every~$\varphi\in
C^\infty_0(B_{r_0}(x),\,[0,+\infty))$,
$$ \int_{B_{r_0}(x)} u(y)\,\Delta\varphi(y)\,dy\ge0.$$
\item[(vii).] The function~$u$ is upper semicontinuous and for every
bounded open sets~$\Omega'\Subset\Omega$ and every harmonic
function~$h\in C^2(\Omega')\cap
C(\overline{\Omega'})$ such that~$u\le h$ on~$\partial\Omega'$
we have that~$u\le h$ in~$\Omega'$.
\item[(viii).] The function~$u$ is upper semicontinuous and for every~$x\in\Omega$
there exists~$r_0>0$ such that~$B_{r_0}(x)\Subset\Omega$ and
for every harmonic function~$h\in C^2(B_{r_0}(x))\cap
C(\overline{B_{r_0}(x)})$ such that~$u\le h$ on~$\partial B_{r_0}(x)$
we have that~$u\le h$ in~$B_{r_0}(x)$.
\end{itemize}
\end{theorem}

\begin{proof} To start with, we prove that
\begin{equation}\label{PA-001}
{\mbox{the statements~(i), (iii), (v) and (vii) are all equivalent.}}
\end{equation}
The strategy
that we adopt for this purpose is
to show that~(i) implies~(v),
(v) implies~(vii), (vii) implies~(iii) and~(iii) implies~(i).\medskip 

First, we show that statement~(i)
implies statement~(v).
To this end,
we can suppose that
\begin{equation}\label{KSffggMD2}{\mbox{$u$ is not
constantly equal to~$-\infty$ in any connected component,}}\end{equation}
otherwise we are done. We claim that
\begin{equation}\label{okmdOffSODLKMF}
u\in L^1_{\rm loc}(\Omega).\end{equation}
For this, we consider a connected component~$\widetilde\Omega$
of~$\Omega$ and we define~$\Omega_0$
to be the collection of points~$y\in\widetilde\Omega$
such that there exists~$\rho>0$ such that~$B_\rho(y)\Subset\widetilde\Omega$
and~$u\in L^1(B_\rho(y))$. Our goal is to show that
\begin{equation}\label{okmdOffSODLKMF2}
\Omega_0=\widetilde\Omega.\end{equation}
The proof is by contradiction: suppose not and
let~$y_0\in\widetilde\Omega\setminus\Omega_0$.
Choose~$\rho_0>0$ sufficiently small such that~$B_{3\rho_0}(y_0)\Subset
\widetilde\Omega$. Let now~$\xi\in B_{\rho_0}(y_0)$.
We observe that
\begin{equation}\label{SOjnfdli346lli34udf}
u\not\in L^1(B_{2\rho_0}(\xi)),\end{equation}
otherwise, since~$B_{\rho_0}(y_0)\subseteq B_{2\rho_0}(\xi)$,
we would have that~$u\in L^1(B_{\rho_0}(y_0))$,
contradicting the assumption that~$y_0\not\in\Omega_0$.

Using~\eqref{SOjnfdli346lli34udf} and the fact that~$u$ is bounded
from above (recall Lemma~\ref{ABOABP}), we find that
$$ \int_{B_{2\rho_0}(\xi)} u(y)\,dy=-\infty.$$
This and the setting in~(i) give that~$u(\xi)=-\infty$.
Since~$\xi$ was an arbitrary point in~$B_{\rho_0}(y_0)$, we have thereby
established that~$u=-\infty$ in~$B_{\rho_0}(y_0)$.
In particular, we have that~$B_{\rho_0}(y_0)\subseteq\widetilde\Omega\setminus\Omega_0$.
As a result, since~$y_0$ is taken to be any
point of~$\widetilde\Omega\setminus\Omega_0$, we have
that~$\widetilde\Omega\setminus\Omega_0$ is open.
Since~$\widetilde\Omega\setminus\Omega_0$ is also closed
in the topology of~$\widetilde\Omega$, we have thus
shown that~$\widetilde\Omega\setminus\Omega_0=\widetilde\Omega$,
that is
\begin{equation}\label{KSffggMD}
\Omega_0=\varnothing.\end{equation}
We now exploit~\eqref{KSffggMD2} and we take~$\zeta_0\in\widetilde\Omega$
such that~$u(\zeta_0)\not=-\infty$. This and~(i) give that,
if~$r_0>0$ is sufficiently small such that~$B_{r_0}(\zeta_0)\Subset\widetilde\Omega$,
then
$$ \int_{ B_{r_0}(\zeta_0)}u(y)\,dy>-\infty.$$
This and the boundedness from above of~$u$
(recall again Lemma~\ref{ABOABP}) yield that~$u\in L^1(B_{r_0}(\zeta_0))$,
which is in contradiction with~\eqref{KSffggMD}.
This contradiction establishes~\eqref{okmdOffSODLKMF2},
from which~\eqref{okmdOffSODLKMF} plainly follows.

Hence, by~\eqref{okmdOffSODLKMF},
to complete the proof of~(v), we take~$\varphi\in
C^\infty_0(\Omega,\,[0,+\infty))$. We let
$$ f_\rho(x):=\frac1{\rho^2}\left(
\fint_{B_\rho(x)} \varphi(y)\,dy-\varphi(x)\right)$$
and we
employ
Theorem~\ref{KAHAR2} to see that, for every~$x\in\R^n$,
$$\lim_{\rho\searrow0} f_\rho(x)=
\frac{1}{2(n+2)}\,\Delta \varphi(x).$$
In addition, by a Taylor expansion and an odd symmetry cancellation,
for small~$\rho>0$ we have that
\begin{eqnarray*}&& \left|
\int_{B_\rho(x)} \big(\varphi(y)-\varphi(x)\big)\,dy
\right|=
\left|
\int_{B_\rho(x)} \Big(
\nabla\varphi(x)\cdot(y-x)+O(|y-x|^2)\Big)\,dy
\right|\\&&\qquad=
\left|
\int_{B_\rho(x)} O(|y-x|^2)\,dy
\right|=O(\rho^{n+2}),
\end{eqnarray*}
with a bound depending only on~$n$ and~$\|\varphi\|_{C^2(\R^n)}$.
This gives that~$f_\rho\in L^\infty(\R^n)$.
As a consequence, recalling~\eqref{okmdOffSODLKMF},
we can utilize the Dominated Convergence Theorem and
conclude that, for every~$\Omega'\Subset\Omega$,
\begin{equation}\label{POSIKpisdfjvarr}
\lim_{\rho\searrow0}\int_{\Omega'} u(x)\,f_\rho(x)\,dx=
\int_{\Omega'}\left(
\lim_{\rho\searrow0} u(x)\,f_\rho(x)\right)\,dx=\frac{1}{2(n+2)}\,
\int_{\Omega'} u(x)\,\Delta\varphi(x)\,dx.
\end{equation}
Furthermore, for all~$\rho>0$ sufficiently small, and choosing~$\Omega'$
such that the support of~$\varphi$ is contained in~$\Omega'$,
\begin{eqnarray*}
&&\int_{\Omega'} u(x)\,f_\rho(x)\,dx=
\frac1{\rho^2}\int_{\Omega'} u(x)\,
\left(
\fint_{B_\rho(x)}\varphi(y)\,dy-\varphi(x)\right)\,dx
\\&&\qquad=\frac1{\rho^{2}\,|B_\rho|}
\int_{\R^n}
\left(
\int_{B_\rho(y)}u(x)\,\varphi(y)\,dx\right)\,dy
-\frac{1}{\rho^2}
\int_{\Omega'} u(x)\,\varphi(x)\,dx\\&&\qquad\ge
\frac1{\rho^{2}}
\int_{\R^n}
u(y)\,\varphi(y)\,dy
-\frac{1}{\rho^2}
\int_{\Omega'} u(x)\,\varphi(x)\,dx=0,
\end{eqnarray*}
thanks to~(i). Combining this with~\eqref{POSIKpisdfjvarr},
we find that
$$ \int_\Omega u(x)\,\Delta\varphi(x)\,dx
=\int_{\Omega'} u(x)\,\Delta\varphi(x)\,dx
\ge0,$$
which is the desired claim in~(v).\medskip

Now we show that~(v) implies~(vii). To this end, given~$\Omega'\Subset\Omega$,
we can restrict our analysis to the connected component of~$\Omega$ containing~$\Omega'$
and assume that~$u$ is not identically equal to~$-\infty$ there, otherwise~(vii) is
obvious.
Hence we can suppose that~$u\in L^1_{\rm loc}(\Omega)$. We also remark that, if~$B_r(x)\Subset\Omega$, then
\begin{equation}\label{SMSDconti} \lim_{p\to x}\fint_{B_r(p)}u(y)\,dy=\fint_{B_r(x)}u(y)\,dy.\end{equation}
To check this, let~$p_j$ be a sequence of points converging to~$x$ as~$j\to+\infty$
and let~$f_j:=u \chi_{B_r(p_j)}$. Notice that~$|f_j|\le u\in L^1_{\rm loc}(\Omega)$.
Consequently, by the Dominated Convergence Theorem,
$$ \lim_{j\to+\infty}\int_{B_r(p_j)}u(y)\,dy
=\int_{B_r(x)}u(y)\,dy,$$
and this proves~\eqref{SMSDconti}.

The gist now is to consider~$x\in\Omega$
and, for every~$r>0$ such that~$B_r(x)\Subset\Omega$, the
volume average
$$ {\mathcal{V}}(r):=\fint_{B_r(x)}u(y)\,dy.$$
We claim that
\begin{equation}\label{AKSM:PSDKMMD:034r-1}
{\mbox{if~\eqref{vuuvrdcsbaP0lqefo}
holds true, the the function~${\mathcal{V}}$
is monotone increasing in~$r$.}}\end{equation}
To check this,
we let~$\varrho\ge\rho>0$
such that~$B_\varrho(x)\Subset\Omega$.
We consider the regularizations~$\Gamma_\varrho$ and~$\Gamma_\rho$
of the fundamental solution that was introduced in~\eqref{GAROGRA}
and we observe that~$\Gamma_\varrho\le\Gamma_\rho$.
Thus, setting, for each~$y\in\R^n$,
$$\varphi_{\varrho,\rho}(y):=\Gamma_\rho(y-x)-\Gamma_\varrho(y-x),$$
we have that~$\varphi_{\varrho,\rho}\in C^{1,1}_0(B_\varrho(x))$.

Now we utilize the statement in~(v).
For this, we remark that, by density, the inequality~\eqref{vuuvrdcsbaP0lqefo}
holds true for all~$\varphi\in C^{1,1}_0(\Omega)$, therefore,
recalling
the basic property of the regularized fundamental solution~\eqref{DEGAM},
\begin{equation*}\begin{split}&0\le
\int_\Omega u(y)\,\Delta\varphi_{\varrho,\rho}(y)\,dy=
\int_{B_\rho(x)} u(y)\,\Delta\Gamma_{\rho}(y-x)\,dy
-\int_{B_\varrho(x)} u(y)\,\Delta\Gamma_{\varrho}(y-x)\,dy\\&\qquad\qquad\qquad=
-
\fint_{B_\rho(x)} u(y)\,dy+\fint_{B_\varrho(x)} u(y)\,dy,\end{split}
\end{equation*}
and this establishes~\eqref{AKSM:PSDKMMD:034r-1}.

As a consequence of~\eqref{AKSM:PSDKMMD:034r-1} the limit
as~$r\searrow0$ of~${\mathcal{V}}(r)$ exists (and this,
for every given~$x\in\Omega$).
Hence
(recall the Lebesgue's Differentiation Theorem,
see e.g.~\cite[Theorem 7.15]{MR3381284}) we can redefine~$u$
in a set of null measure and suppose that, for every~$x\in\Omega$,
\begin{equation}\label{pi9ltar}
u(x)=\lim_{r\searrow0}\fint_{B_r(x)} u(y)\,dy=\lim_{r\searrow0}{\mathcal{V}}(r).
\end{equation}
With this setting, we have that~$u$
is the pointwise limit of the decreasing sequence of continuous functions~$u_j(x):=\fint_{B_{1/j}(x)}u(y)\,dy$
(the continuity being a consequence of~\eqref{SMSDconti}),
hence
\begin{equation}\label{RECU}
{\mbox{$u$ is upper semicontinuous,}}\end{equation} thanks to
Lemma~\ref{65co99865bSdsaJ-04}.

An interesting byproduct of the monotonicity of the volume average~${\mathcal{V}}$ with respect to~$r$
and the integral in polar coordinates (see e.g.~\cite[Theorem 3.12]{MR3409135}) is that, possibly excluding
a set of zero Lebesgue measure of values of~$r$,
\begin{eqnarray*} &&0\le\frac{d}{dr}\left(\fint_{B_r(x)}u(y)\,dy
\right)=
\frac{1}{|B_r|}
\frac{d}{dr}\left(\int_{B_r(x)}u(y)\,dy
\right)
-\frac{n}{r}\fint_{B_r(x)}u(y)\,dy\\&&\qquad\qquad=
\frac{1}{|B_r|}
\int_{\partial B_r(x)}u(y)\,d{\mathcal{H}}^{n-1}_y
-\frac{n}{r}\fint_{B_r(x)}u(y)\,dy.
\end{eqnarray*}
As a result,
$$ \fint_{B_r(x)}u(y)\,dy\le
\frac{r}{n\,|B_r|}
\int_{\partial B_r(x)}u(y)\,d{\mathcal{H}}^{n-1}_y
=\fint_{\partial B_r(x)}u(y)\,d{\mathcal{H}}^{n-1}_y
.$$
Consequently,
using again~\eqref{pi9ltar} and the increasing monotonicity of~${\mathcal{V}}$, for almost every~$r$,
\begin{equation}\label{09iKlsmS-03045kscsd}
u(x)\le\fint_{B_r(x)} u(y)\,dy\le\fint_{\partial B_r(x)}u(y)\,d{\mathcal{H}}^{n-1}_y.
\end{equation}
Now, let~$h$ be as requested in~(vii). Given~$\e>0$, we let~$h_\e:=h+\e$.
Notice that~$h_\e\ge u+\e$ on~$\partial\Omega'$.
We claim that there
exists~$\delta_\e>0$ such that
\begin{equation}\label{INTOS}
h_\e\ge u+\frac\e2\,{\mbox{ in }}\,\Omega_\e':=\bigcup_{p
\in\partial\Omega'} \big(B_{\delta_\e}(p)\cap\overline{\Omega'}\big).
\end{equation}
Indeed, suppose not. Then, for every~$j$ there exists~$q_j\in\overline{\Omega'}$
at distance less than~$\frac1j$ from~$\partial\Omega'$ and such that~$
h_\e(q_j)< u(q_j)+\frac\e2$. Since~$\overline{\Omega'}$
is compact, up to a subsequence, we can suppose that~$q_j\to
q\in\overline{\Omega'}$ as~$j\to+\infty$. By construction~$q\in\partial\Omega'$.
Hence, by the upper semicontinuity of~$u$,
$$ h(q)+\e=h_\e(q)=\lim_{j\to+\infty}h_\e(q_j)
\le\limsup_{j\to+\infty} u(q_j)+\frac\e2\le u(q)+\frac\e2,
$$
hence~$h(q)<u(q)$, in contradiction with the assumptions in~(vii).
This completes the proof of~\eqref{INTOS}.

Now, we consider a radial function~$\tau\in C^\infty_0(B_1,\,[0,+\infty))$
with~$\int_{B_1}\tau(x)\,dx=1$.
Given~$\eta>0$, we let~$\tau_\eta(x):=\frac1{\eta^n}
\tau\left(\frac{x}{\eta}\right)$
and define
\begin{equation}\label{CONVETA}
u_\eta:=u*\tau_\eta.\end{equation} We claim that, if~$\eta>0$ is sufficiently small,
\begin{equation}\label{SM:0oedKSKDd35D}
{\mbox{$u_\eta\ge u$ in~$\overline{\Omega'}$.}}
\end{equation}
Indeed, integrating in polar coordinates (see e.g.~\cite[Theorem 3.12]{MR3409135}),
\begin{equation}\label{conwsdfisces093mat}
\begin{split}
&u_\eta(x)=\int_{B_\eta} u(x-y)\,\tau_\eta(y)\,dy=
\int_0^\eta \left(
\int_{\partial B_\rho}
u(x-\zeta)\,\tau_\eta(\zeta)\,d{\mathcal{H}}^{n-1}_\zeta
\right)\,d\rho\\&\qquad\qquad=
\int_0^\eta \tau_\eta(\rho e_1)\left(
\int_{\partial B_\rho}
u(x-\zeta)\,d{\mathcal{H}}^{n-1}_\zeta
\right)\,d\rho
\end{split}\end{equation}
Consequently, recalling~\eqref{09iKlsmS-03045kscsd},
\begin{eqnarray*}&&
u_\eta(x)\ge
u(x)\,
\int_0^\eta \tau_\eta(\rho e_1) {\mathcal{H}}^{n-1}(\partial B_\rho)\,d\rho
=u(x)\,\int_{B_\eta} \tau_\eta(y)\,dy=u(x),
\end{eqnarray*}
which proves~\eqref{SM:0oedKSKDd35D}.

Now we observe that, for small~$\eta>0$,
\begin{equation}\label{KSM0942tgfdvgdsyte}
{\mbox{$\Delta u_\eta\ge0$ in~$\overline{\Omega'}$.}}
\end{equation}
Indeed, for every~$\varphi\in C^\infty_0(\Omega,\,[0,+\infty))$, we let~$\varphi_\eta:=\varphi*\tau_\eta$
and we remark that~$\varphi_\eta\in C^\infty_0(\Omega,\,[0,+\infty))$
as long as~$\eta$ is sufficiently small. Then, using~(v),
\begin{eqnarray*}
0\le\int_\Omega u(y)\,\Delta\varphi_\eta(y)\,dy=
\int_{\R^n}\left(\int_{\R^n} u(y)\,\varphi(z)\,\Delta\tau_\eta(y-z)\,dz\right)\,dy=\int_{\R^n}\Delta u_\eta(z)\varphi(z)
\,dz\end{eqnarray*}
and thus~\eqref{KSM0942tgfdvgdsyte} since~$\Delta u_\eta$ is continuous in~$\overline{\Omega'}$.

We also set~$h_{\e,\eta}:=h_\e*\tau_\eta$ and we observe that, if~$\eta>0$ is conveniently small,
then~$\Delta h_{\e,\eta}=(\Delta h_\e)*\tau_\eta=(\Delta h)*\tau_\eta=0$ in~$\overline{\Omega'}$.
{F}rom this and~\eqref{KSM0942tgfdvgdsyte}, setting~$v_{\e,\eta}:=u_\eta-h_{\e,\eta}$, we deduce that
\begin{equation}\label{KDM:99495OEKMF}
{\mbox{$\Delta v_{\e,\eta}\ge0$ in~$\overline{\Omega'}$.}}\end{equation}
Also, if~$\eta$ is sufficiently small, possibly in dependence of~$\e$, we know
from~\eqref{INTOS} that, for every~$x\in\partial\Omega'$,
$$ v_{\e,\eta}(x)
=(u-h_\e)*\tau_\eta(x)=
\int_{B_\eta(x)} \big(u(y)-h_\e(y)\big)\,\tau_\eta(x-y)\,dy<0.
$$
Consequently, by~\eqref{KDM:99495OEKMF} and the Maximum Principle in
Corollary~\ref{WEAKMAXPLE}(i),
we conclude that~$v_{\e,\eta}(x)\le0$ for every~$x\in\Omega'$.
For this reason, and recalling~\eqref{SM:0oedKSKDd35D},
by taking limits, for every~$x\in\Omega'$,
\begin{eqnarray*}&&
0\ge\lim_{\e\searrow0}\left(\lim_{\eta\searrow0}
v_{\e,\eta}(x)
\right)=
\lim_{\e\searrow0}\left(\lim_{\eta\searrow0}(
u_\eta(x)-h_{\e,\eta}(x))
\right)\ge\lim_{\e\searrow0}\left(\lim_{\eta\searrow0}(
u(x)-h_{\e,\eta}(x))
\right)\\&&\qquad\qquad=
\lim_{\e\searrow0}\left(
u(x)-h_{\e}(x)
\right)=u(x)-h(x).
\end{eqnarray*}
This proves that~(v) implies~(vii).\medskip

Now we prove that~(vii) implies~(iii).
For this, pick a ball~$B_r(x)\Subset\Omega$.
Let also~$R>r$ be such that~$B_R(x)\Subset\Omega$.
Given~$j\in\N$, we exploit
Corollary~\ref{754PCHAYPP}(iii)
to find a decreasing sequence of functions~$u_j\in C^\infty(B_R(x))$
that converges to~$u$ at every point of~$B_R(x)$.
Let~$h_j$ be the harmonic function in~$B_r(x)$ that
coincides with~$u_j$ on~$\partial B_r(x)$.
The existence of~$h_j$ is warranted, for instance, by 
the Poisson Kernel representation in
Theorems~\ref{POIBALL1}
and~\ref{POIBALL}.
Also, we know that~$h_j \in C^2(B_r(x))\cap
C(\overline{B_r(x)})$. Since~$h_j=u_j\ge u$ on~$\partial B_r(x)$,
we deduce from~(vii) and the Mean Value Formula
for harmonic functions 
(recall Theorem~\ref{KAHAR}(ii)) that
\begin{equation}\label{uhjux3x} u(x)\le h_j(x)=
\lim_{\rho\searrow r}
\fint_{\partial B_\rho(x)} h_j(y)\,d{\mathcal{H}}^{n-1}_y=
\fint_{\partial B_r(x)} h_j(y)\,d{\mathcal{H}}^{n-1}_y=
\fint_{\partial B_r(x)} u_j(y)\,d{\mathcal{H}}^{n-1}_y.\end{equation}
Now we pass~$j\to+\infty$. For this, we stress that~$u_j\le u_1$
and that~$u_j$ converges to~$u$ everywhere in~$B_r(x)$
in a monotone decreasing way: therefore
we can employ the Monotone Convergence Theorem (see e.g.~\cite[Theorem 10.27(ii)]{MR3381284}, and notice
that a convergence of~$u_j$ in the almost everywhere sense
with respect to Lebesgue measure would not suffice to exploit this result).
In this way, we infer that
$$\lim_{j\to+\infty}\fint_{\partial B_r(x)} u_j(y)\,d{\mathcal{H}}^{n-1}_y=\fint_{\partial B_r(x)} u(y)\,d{\mathcal{H}}^{n-1}_y.$$
Combining this with~\eqref{uhjux3x} we have established that~(iii) holds true,
thus proveing that~(vii) implies~(iii).\medskip

Now we prove that~(iii) implies~(i).
For this, we use polar coordinates (see e.g.~\cite[Theorem 3.12]{MR3409135}) to find that
\begin{eqnarray*}&& \fint_{B_r(x)} \big(u(y)-u(x)\big)\,dy=\frac{1}{|B_r|}\int_0^r\left(
\int_{\partial B_\rho(x)} \big(u(\zeta)-u(x)\big)\,d{\mathcal{H}}^{n-1}_\zeta\right)\,d\rho\\&&\qquad\qquad=
\frac{{\mathcal{H}}^{n-1}(\partial B_1)}{|B_r|}\int_0^r\left(\rho^{n-1}
\fint_{\partial B_\rho(x)} \big(u(\zeta)-u(x)\big)\,d{\mathcal{H}}^{n-1}_\zeta\right)\,d\rho.
\end{eqnarray*}
Thus, if~(iii) holds true, then the latter integrand is positive, and one obtains~(i), as desired.
This completes the proof of~\eqref{PA-001}.\medskip

As a matter of fact, given a ball~$B_{r_0}(x)\Subset\Omega$,
we can now
apply~\eqref{PA-001} to this ball instead of~$\Omega$:
in this way, we obtain that
\begin{equation}\label{PA-002}
{\mbox{the statements~(ii), (iv), (vi) and (viii) are all equivalent.}}
\end{equation}

In light of~\eqref{PA-001} and~\eqref{PA-002}, to complete the proof
of Theorem~\ref{SUBHS-EQ}, it is enough to show that
the statements~(v) and~(vi) are equivalent. As a matter of fact,
since~(v) obviously implies~(vi), it suffices to prove that
\begin{equation}\label{PA-003}
{\mbox{statement~(vi) implies statement~(v).}}
\end{equation}
To check this, we can suppose that~$u\in L^1_{\rm loc}(\Omega)$,
otherwise we are done,
and we
let~$\varphi\in
C^\infty_0(\Omega,\,[0,+\infty))$, as requested in~(v).
We can therefore denote by~$K$ a compact subset of~$\Omega$
such that~$\varphi=0$ outside~$K$. For each~$x\in\Omega$
we let~$r_0=r_0(x)>0$ as in~(vi) and we notice that
$$ K\subseteq\bigcup_{x\in K} B_{r_0(x)}(x).$$
In view of the compactness of~$K$, we can take a finite subcover,
thus finding~$x_1,\dots,x_N\in\Omega$ such that
$$ K\subseteq\bigcup_{j=1}^N B_{{r_0(x_j)}}(x_j).$$
Then we consider a partition of unity
made of functions~$\phi_1\in
C^\infty(B_{r_0(x_1)}(x_1),\,[0,1])$, $\dots$, $\phi_N\in C^\infty(B_{r_0(x_N)}(x_N),\,[0,1])$,
with finite overlapping supports, such that
$$ \sum_{j=1}^N\phi_j=1\qquad{\mbox{in }}\,\bigcup_{j=1}^N B_{{r_0(x_j)}/2}(x_j).$$
Notice in particular that
$$ \sum_{j=1}^N\nabla\phi_j=0\qquad{\mbox{and}}\qquad
\sum_{j=1}^N\Delta\phi_j=0
\qquad{\mbox{in }}\,\bigcup_{j=1}^N B_{{r_0(x_j)}/2}(x_j).$$
Thus, if we define~$\varphi_j:=\varphi\phi_j$,
we see that
\begin{equation}\label{9-0LSMD:alscm} \Delta\varphi=\Delta\varphi\chi_K=
\sum_{j=1}^N\Delta\varphi\phi_j=
\sum_{j=1}^N
\Big(\Delta(\varphi\phi_j)-\varphi\Delta\phi_j
-2\nabla\varphi\cdot\nabla\phi_j\Big)=\sum_{j=1}^N\Delta\varphi_j.
\end{equation}
Also, since~$\varphi_j\in C^\infty(B_{r_0(x_j)}(x_j),\,[0,1])$,
by~(vi) we know that
$$ \int_{\Omega} u(y)\,\Delta\varphi_j(y)\,dy=\int_{B_{r_0(x_j)}(x_j)}u(y)\,\Delta\varphi_j(y)\,dy\ge0.$$
This and~\eqref{9-0LSMD:alscm} yield that
$$\int_\Omega u(y)\,\Delta\varphi(y)\,dy=\sum_{j=1}^N\int_\Omega u(y)\,\Delta\varphi_j(y)\,dy\ge0,$$
that is statement~(v).
This establishes~\eqref{PA-003} and thus completes the proof
of Theorem~\ref{SUBHS-EQ}.
\end{proof}

\begin{definition}\label{SHAITY}
Given an open set~$\Omega\subseteq\R^n$,
an upper semicontinuous
function~$u: \Omega\to\R\cup\{-\infty\}$
is said to be subharmonic\index{subharmonic function} if any of the equivalent properties
listed in Theorem~\ref{SUBHS-EQ} holds true.

A function~$u:\Omega\to\R\cup\{+\infty\}$ is said to be superharmonic\index{superharmonic function}
if~$-u$ is subharmonic.
\end{definition}

The name ``subharmonic'' is clearly
a legacy of Theorem~\ref{SUBHS-EQ}(vii).
It is also interesting to notice that
the statements~(ii), (iv), (vi) and (viii) in
Theorem~\ref{SUBHS-EQ} are simply the localized versions
of the statements~(i), (iii), (v) and (vii): in this sense,
a function is subharmonic in~$\Omega$
if and only if it is subharmonic in any subdomain of~$\Omega$.

As a consequence of Theorem~\ref{KScvbMS:0okrmt4ht}
and Theorem~\ref{SUBHS-EQ}(v), we also have that
\begin{equation}\label{Coamscfunconftahf}
{\mbox{$-\Gamma$
is subharmonic,}}\end{equation} being~$\Gamma$ the fundamental solution
in~\eqref{GAMMAFU}
(extended as~$\Gamma(0):=+\infty$ when~$n\ge2$, which also clarifies
a good reason for allowing upper semicontinuous
functions with values in~$[-\infty,+\infty)$).

Moreover, we observe that, in virtue of the local integrability property in
Theorem~\ref{SUBHS-EQ}(v), we have that if~$u$
is subharmonic in~$\Omega$ then~$\{|u|=+\infty\}$
has null measure in~$\Omega$, unless~$u=-\infty$ in the whole of a connected
component of~$\Omega$.
\medskip

When the function is smooth, the notion of subharmonicity boils down to an inequality
in the corresponding Laplace equation, as clarified by the following observation:

\begin{lemma}\label{FACI098}
Given an open set~$\Omega\subseteq\R^n$,
a function~$u\in C^2(\Omega,\R\cup\{-\infty\})$
is subharmonic if and only if~$\Delta u\ge0$ in~$\Omega$.
\end{lemma}

\begin{proof} Using Theorem~\ref{SUBHS-EQ}(v) and the second Green's Identity~\eqref{GRr2}, we have
that a function~$u: \in C^2(\Omega,\R\cup\{-\infty\})$
is subharmonic if and only if
$$ \int_\Omega \Delta u(x)\,\varphi(x)\,dx\ge0$$
for every~$\varphi\in
C^\infty_0(\Omega,\,[0,+\infty))$,
from which the desired result follows.
\end{proof}

As a byproduct of the above observations, we also have that
the averages of a subharmonic function increase with radius:

\begin{corollary}
Let~$u$ be subharmonic in~$\Omega$. Let~$x\in\Omega$
and~$R>r>0$ be such that~$B_R(x)\Subset\Omega$. Then,
\begin{equation}\label{TGbickdtollibstaiks}
\begin{split} &\fint_{ B_r(x)} u(y)\,dy\le
\fint_{ B_R(x)} u(y)\,dy\\{\mbox{and }}\qquad
& \fint_{\partial B_r(x)} u(y)\,d{\mathcal{H}}^{n-1}_y\le
\fint_{\partial B_R(x)} u(y)\,d{\mathcal{H}}^{n-1}_y.
\end{split}\end{equation}
\end{corollary}

\begin{proof} Without loss of generality, we can reduce to the case of a connected domain~$\Omega$.
Also, we can suppose that~$u$ is not identically~$-\infty$ in~$\Omega$, otherwise the desired
claims are obviously satisfied.
Hence, we exploit Theorem~\ref{SUBHS-EQ}(v)
and we conclude that~$u
\in L^1_{\rm loc}(\Omega)$
and~\eqref{vuuvrdcsbaP0lqefo} holds true.
This and~\eqref{AKSM:PSDKMMD:034r-1} give the first inequality in~\eqref{TGbickdtollibstaiks}.

Now, let~$R'>R$ be such that~$B_{R'}(x)\Subset\Omega$.
We consider the convolution~$u_\eta$ introduced 
in~\eqref{CONVETA}. Given~$\varrho\in(r,R')$,
we take~$h_{\eta,\varrho}$ to be the harmonic function in~$B_\varrho(x)$
that coincides with~$u_\eta$ along~$\partial B_\varrho(x)$.
Notice that~$h_{\eta,\varrho}\in C^2(B_\varrho(x))\cap C(\overline{B_\varrho(x)})$.
Since~$u_\eta$ is subharmonic in~$B_\varrho(x)$, 
thanks to~\eqref{KSM0942tgfdvgdsyte} and Lemma~\ref{FACI098}, we thereby deduce from Theorem~\ref{SUBHS-EQ}(vii)
that
\begin{equation}\label{KStd643misk}
{\mbox{$u_\eta\le h_{\eta,\varrho}$ in $B_\varrho(x)$.}}
\end{equation}
Now we take~$h_{\eta,r}$ to be the harmonic function in~$B_r(x)$
that coincides with~$u_\eta$ along~$\partial B_r(x)$.
Notice that~$h_{\eta,r}=u_\eta\le h_{\eta,\varrho}$
along~$\partial B_r(x)$, due to~\eqref{KStd643misk}.
{F}rom this and Theorem~\ref{SUBHS-EQ}(vii),
we conclude that~$h_{\eta,r}\le h_{\eta,\varrho}$
in~$B_r(x)$.
As a consequence, using the Mean Value Formula in Theorem~\ref{KAHAR}(ii),
\begin{eqnarray*}&&
\fint_{\partial B_r(x)} u_{\eta}(y)\,d{\mathcal{H}}^{n-1}_y=
\fint_{\partial B_r(x)} h_{\eta,r}(y)\,d{\mathcal{H}}^{n-1}_y=
h_{\eta,r}(x)\le h_{\eta,\varrho}(x)\\&&\qquad\qquad=\fint_{\partial B_\varrho(x)} h_{\eta,\varrho}(y)\,d{\mathcal{H}}^{n-1}_y
=\fint_{\partial B_\varrho(x)} u_{\eta}(y)\,d{\mathcal{H}}^{n-1}_y.\end{eqnarray*}
This and~\eqref{SM:0oedKSKDd35D} yield that
$$ 
\fint_{\partial B_r(x)} u(y)\,d{\mathcal{H}}^{n-1}_y\le\fint_{\partial B_\varrho(x)} u_{\eta}(y)\,d{\mathcal{H}}^{n-1}_y,$$
that is
$$ 
\varrho^{n-1}\,{\mathcal{H}}^{n-1}(\partial B_1)
\fint_{\partial B_r(x)} u(y)\,d{\mathcal{H}}^{n-1}_y\le\int_{\partial B_\varrho(x)} u_{\eta}(y)\,d{\mathcal{H}}^{n-1}_y.$$
We now integrate this inequality
in~$\varrho\in(R-\delta,R+\delta)$, with~$\delta\in(0,R'-R)$, using
polar coordinates in the right hand side (see e.g.~\cite[Theorem 3.12]{MR3409135})
thus finding that
$$ \frac{(R+\delta)^{n}-(R-\delta)^{n}}{n}\,
\,{\mathcal{H}}^{n-1}(\partial B_1)\fint_{\partial B_r(x)} u(y)\,d{\mathcal{H}}^{n-1}_y\le
\int_{ B_{R+\delta}(x)\setminus B_{R-\delta}(x)} u_{\eta}(y)\,dy.$$
Now we exploit that~$u_\eta\to u$ in~$L^1(B_{R'}(x)$
(see e.g.~\cite[Theorem~9.6]{MR3381284}) and we thus
conclude that
$$ \frac{(R+\delta)^{n}-(R-\delta)^{n}}{n}\,
\,{\mathcal{H}}^{n-1}(\partial B_1)\fint_{\partial B_r(x)} u(y)\,d{\mathcal{H}}^{n-1}_y\le
\int_{ B_{R+\delta}(x)\setminus B_{R-\delta}(x)} u(y)\,dy.$$
Hence, we combine
polar coordinates (see e.g.~\cite[Theorem 3.12]{MR3409135})
and L'H\^opital's Rule and we find that
\begin{eqnarray*}&&2\int_{\partial B_{R}(x)} u(y)\,d{\mathcal{H}}^{n-1}_y=
\lim_{\delta\searrow0} \frac{d}{d\delta}\left(
\int_{ B_{R+\delta}(x)} u(y)\,dy
-\int_{ B_{R-\delta}(x)} u(y)\,dy
\right)\\&&\qquad=
\lim_{\delta\searrow0}\frac{\frac{d}{d\delta}
\int_{ B_{R+\delta}(x)\setminus B_{R-\delta}(x)} u(y)\,dy
}{\frac{d}{d\delta}\delta}
=
\lim_{\delta\searrow0}\frac{
\int_{ B_{R+\delta}(x)\setminus B_{R-\delta}(x)} u(y)\,dy
}{\delta}\\&&\qquad\ge\lim_{\delta\searrow0}
\frac{(R+\delta)^{n}-(R-\delta)^{n}}{n\,\delta}\,
\,{\mathcal{H}}^{n-1}(\partial B_1)
\fint_{\partial B_r(x)} u(y)\,d{\mathcal{H}}^{n-1}_y\\&&\qquad=
2R^{n-1}
\,{\mathcal{H}}^{n-1}(\partial B_1)
\fint_{\partial B_r(x)} u(y)\,d{\mathcal{H}}^{n-1}_y=
2
\,{\mathcal{H}}^{n-1}(\partial B_R)\fint_{\partial B_r(x)} u(y)\,d{\mathcal{H}}^{n-1}_y,
\end{eqnarray*}
which is the second inequality in~\eqref{TGbickdtollibstaiks}.
\end{proof}

Using Theorem~\ref{SUBHS-EQ} we can also
strengthen the Strong Maximum Principle given in
Theorem~\ref{STRONGMAXPLE1} by removing the smoothness
assumption on~$u$:

\begin{lemma}\label{LplrfeELLtahsmedschdfb4camntM}
Let~$\Omega\subseteq\R^n$ be open and connected,
and let~$u$ be subharmonic in~$\Omega$.
If there exists~$\overline{x}\in\Omega$
such that~$u(\overline{x})=\sup_\Omega u$, then u is constant.
\end{lemma}

\begin{proof} Suppose that~$u(\overline{x})=\sup_\Omega u$
for some~$\overline{x}\in\Omega$.
We argue as in the proof of Theorem~\ref{SUBHS-EQ}, considering
the set~${\mathcal{U}}$ in~\eqref{S:osl4995ffdssSloamed},
which is nonvoid
since~$\overline{x}\in{\mathcal{U}}$.

Also, 
\begin{equation}\label{BRFSVBcrgaaVOPAAVIKSAVTDHJ}
{\mbox{the set~${\mathcal{U}}$
is closed in~$\Omega$,}}\end{equation} since if~$x_k\in{\mathcal{U}}$
and~$x_k\to\widetilde{x}$ as~$k\to+\infty$,
we deduce from the upper semicontinuity condition~\eqref{likm8ijk9-0eofkjjgjjn} that
$$ \sup_\Omega u=
\limsup_{k\to+\infty} u(x_k)\le u(\widetilde{x})\le\sup_\Omega u,$$
showing that~$\widetilde{x}\in{\mathcal{U}}$
and establishing~\eqref{BRFSVBcrgaaVOPAAVIKSAVTDHJ}.

Hence, to prove the claim in Lemma~\ref{LplrfeELLtahsmedschdfb4camntM},
it suffices to show that
\begin{equation}\label{BRFSVBcrgaaVOPAAVIKSAVTDHJ2}
{\mbox{the set~${\mathcal{U}}$
is also open.}}\end{equation}
To this end, we consider~$r>0$ such that~$B_r(
\overline{x})\Subset\Omega$
and we utilize Theorem~\ref{SUBHS-EQ} to see that
$$ \sup_\Omega u=u(\overline{x})\le\fint_{B_r(\overline{x})} u(y)\,dy\le
\sup_\Omega u,$$
and accordingly~$u=\sup_\Omega u$ a.e. in~$B_r(\overline{x})$.

For this reason, for every~$\widetilde{p}\in B_r(\overline{x})$,
we can take a sequence~$p_k\in B_r(\overline{x})$
such that~$u(p_k)=\sup_\Omega u$
and~$p_k\to \widetilde{p}$ as~$k\to+\infty$. We thereby deduce
from the upper semicontinuity of~$u$ (recall~\eqref{likm8ijk9-0eofkjjgjjn})
that
$$ \sup_\Omega u=\limsup_{k\to+\infty} u(p_k)\le u(\widetilde{p})\leq\sup_\Omega u,$$
and accordingly~$u(\widetilde{p})=\sup_\Omega u$.
This shows that~$u=\sup_\Omega u$ everywhere in~$B_r(\overline{x})$
and it completes the proof of~\eqref{BRFSVBcrgaaVOPAAVIKSAVTDHJ2},
as desired.
\end{proof}

\begin{corollary}\label{PERSUBAWEAK}
Let~$\widetilde\Omega\subseteq\R^n$ be open and connected,
and let~$u$ be subharmonic in~$\widetilde\Omega$. Let also~$\Omega\Subset\widetilde\Omega$ be bounded.
Then,
$$\sup_\Omega u=\sup_{\partial\Omega} u.$$
\end{corollary}

\begin{proof} The proof of Corollary~\ref{PERSUBAWEAK} goes through
by applying
Lemma~\ref{LplrfeELLtahsmedschdfb4camntM}
in lieu of
Theorem~\ref{STRONGMAXPLE1}. Indeed, one uses Lemma~\ref{ABOABP}
to ensure that~$u$ attains a maximum in~$\overline\Omega$, that is
there exists~$p\in\overline\Omega$ such that~$u(p)=\max_{\overline\Omega}u$.
If~$p\in\partial\Omega$, then the desired claim follows.

If instead~$p\in\Omega$, then we use
Lemma~\ref{LplrfeELLtahsmedschdfb4camntM} to say that~$u$ is constant in~$\Omega$.
Now, we take~$\widetilde x\in\partial\Omega$ and a sequence of points~$x_k\in\Omega$
such that~$x_k\to\widetilde x$ as~$k\to+\infty$. Accordingly, we exploit the upper semicontinuity of~$u$
to obtain that
$$ \max_{\overline\Omega}u=u(p)=\lim_{k\to+\infty}u(x_k)\le \limsup_{x\to\widetilde x}u(x)\le u(\widetilde x)
\le \max_{\overline\Omega}u,
$$
which gives the desired result also in this case.
\end{proof}

In addition, it follows from Theorem~\ref{SUBHS-EQ}
that the set of subharmonic functions is endowed of a conical structure,
namely if~$u$ and~$v$ are subharmonic in~$\Omega$
and~$a$, $b\in[0,+\infty)$, then the function~$au+bv$ is
subharmonic in~$\Omega$ as well. \medskip

An additional interesting feature is that the spherical mean of subharmonic
functions always converges to the value at a point
(and this holds everywhere, not only almost everywhere
as warranted by Lebesgue's Differentiation Theorem,
see e.g.~\cite[Theorem 7.15]{MR3381284}).
Indeed, we have the following observation:

\begin{lemma}\label{SCELTA}
If~$u$ is subharmonic in~$\Omega$ and~$x\in\Omega$, then
$$ \lim_{r\searrow0}\fint_{B_{r}(x)} u(y)\,dy=\lim_{r\searrow0}
\fint_{\partial B_r(x)} u(y)\,d{\mathcal{H}}^{n-1}_y=\limsup_{y\to x} u(y)=u(x).$$
\end{lemma}

\begin{proof} Given~$\e>0$, in view of~\eqref{likm8ijk9-0eofkjjgjjn},
we know that there exists~$\delta>0$ such that, for every~$y\in
B_\delta(x)$, $u(y)\le u(x)+\e$. Thus,
if~$r\in(0,\delta)$,
\begin{equation*}
\fint_{B_{r}(x)} u(y)\,dy\le u(x)+\e\qquad{\mbox{ and }}\qquad
\fint_{\partial B_r(x)} u(y)\,d{\mathcal{H}}^{n-1}_y\le u(x)+\e.
\end{equation*}
As a result,
\begin{equation*}\limsup_{r\searrow0}
\fint_{B_{r}(x)} u(y)\,dy\le u(x)+\e\qquad{\mbox{ and }}\qquad
\limsup_{r\searrow0}
\fint_{\partial B_r(x)} u(y)\,d{\mathcal{H}}^{n-1}_y\le u(x)+\e.
\end{equation*}
Therefore, since~$\e$ is arbitrary,
\begin{equation}\label{98i0us9dos:02eo3fkejnvojfenvd}
\limsup_{r\searrow0}
\fint_{B_{r}(x)} u(y)\,dy\le u(x) \qquad{\mbox{ and }}\qquad
\limsup_{r\searrow0}
\fint_{\partial B_r(x)} u(y)\,d{\mathcal{H}}^{n-1}_y\le u(x).
\end{equation}
On the other hand,
by Theorem~\ref{SUBHS-EQ}(i) and~(iii),
\begin{equation*}
\liminf_{r\searrow0}
\fint_{B_{r}(x)} u(y)\,dy\ge u(x) \qquad{\mbox{ and }}\qquad
\liminf_{r\searrow0}
\fint_{\partial B_r(x)} u(y)\,d{\mathcal{H}}^{n-1}_y\ge u(x).
\end{equation*}
{F}rom this and~\eqref{98i0us9dos:02eo3fkejnvojfenvd},
we conclude that
\begin{equation} \label{56:098re34u29r3ihfejb77569}
\lim_{r\searrow0}\fint_{B_{r}(x)} u(y)\,dy=\lim_{r\searrow0}
\fint_{\partial B_r(x)} u(y)\,d{\mathcal{H}}^{n-1}_y=u(x).\end{equation}

In addition, by Theorem~\ref{SUBHS-EQ}(i),
for every~$k\in\N$ (say,
sufficiently large such that~$B_{1/k}(x)\Subset\Omega$)
there exists~$x_k\in B_{1/k}(x)$ such that~$u(x_k)\ge u(x)$.
As a consequence,
$$ \limsup_{y\to x} u(y)\ge \lim_{k\to+\infty}u(x_k)\ge u(x),$$
which, combined with~\eqref{likm8ijk9-0eofkjjgjjn} leads to
$$ \limsup_{y\to x} u(y)= u(x).$$
{F}rom this and~\eqref{56:098re34u29r3ihfejb77569}
the desired result plainly follows.
\end{proof}

\begin{corollary}
Two subharmonic functions that coincide almost everywhere
necessarily coincide everywhere.
\end{corollary}

\begin{proof} Let~$u$ and~$v$ be subharmonic and equal almost everywhere.
Then, for every~$x$ and every~$r>0$ such that~$B_r(x)$ lies in the subharmonicity domain of~$u$
and~$v$,
$$ \int_{B_{r}(x)} u(y)\,dy=\int_{B_{r}(x)} v(y)\,dy.$$
Thus, recalling Lemma~\ref{SCELTA}, for every point~$x$ we have that
\begin{equation*} u(x)=\lim_{r\searrow0}\fint_{B_{r}(x)} u(y)\,dy
=\lim_{r\searrow0}\fint_{B_{r}(x)} v(y)\,dy=v(x).\qedhere\end{equation*}
\end{proof}

The one-dimensional subharmonic functions
reduce to the convex ones, as remarked in the following observation:

\begin{lemma}\label{DO:cD00}
If~$n=1$ and~$u$ is subharmonic in
an open interval~$I$, then~$u$ is convex in~$I$.
\end{lemma}

\begin{proof} Let~$a$, $b\in I$ with~$a<b$ and, for every~$x\in[a,b]$,
let~$h(x):=\frac{(u(b)-u(a))(x-a)}{b-a}+u(a)$.
Since~$h''=0$, we have that~$h$ is harmonic
and thus we can use Theorem~\ref{SUBHS-EQ}(vii), deducing that, for
every~$t\in[0,1]$, if~$x:=tb+(1-t)a$,
\begin{equation*} u(tb+(1-t)a)=u(x)\le h(x)=
t\,(u(b)-u(a))+u(a)
=t\,u(b)+(1-t)\,u(a).\qedhere\end{equation*}
\end{proof}

The converse of Lemma~\ref{DO:cD00} holds true in every dimension:

\begin{lemma}\label{DO:cD}
If~$u$ is convex in~$\Omega$, then it is subharmonic
in~$\Omega$.
\end{lemma}

\begin{proof}
We recall that a convex function in~$\R^n$ is continuous,
see e.g.~\cite[Theorem~6.7(i)]{MR3409135}.
Moreover,
let~$x\in \Omega$ and let~$r_0>0$ be such that~$B_{r_0}(x)\Subset\Omega$.
Take the supporting plane for~$u$ at~$x$, namely consider~$\omega\in\R^n$
such that~$u(y)\ge\omega\cdot(y-x)+u(x)$ for all~$y$ in~$\Omega$.
Then, for all~$r\in(0,r_0)$,
\begin{equation*} \fint_{B_{r}(x)} u(y)\,dy\ge\fint_{B_{r}(x)}\big(
\omega\cdot(y-x)+u(x)\big)\,dy=
\omega\cdot
\fint_{B_{r}}\eta\,d\eta+u(x)=0+u(x),
\end{equation*}
thanks to an odd symmetry cancellation. This and
Theorem~\ref{SUBHS-EQ}(ii) yields the desired result.
\end{proof}

We now turn our attention to some operations that preserve subharmonicity:

\begin{lemma}\label{POMAX}
The pointwise maximum of two subharmonic functions is subharmonic.
\end{lemma}

\begin{proof} If~$u$ and~$v$ are subharmonic in~$\Omega\Supset B_r(x)$
and~$w:=\max\{u,v\}$,
then we know that~$w$ is upper semicontinuous, due to Lemma~\ref{KSMD:0orkegm}, and
\begin{equation*}
\fint_{ B_r(x)}w(y)\,dy\ge\max\left\{
\fint_{ B_r(x)}u(y)\,dy,\,\fint_{ B_r(x)}v(y)\,dy
\right\}\ge\max\{u(x),v(x)\}=w(x).
\qedhere\end{equation*}
\end{proof}

Of course, the previous result can be extended inductively to the
pointwise maximum of finitely many subharmonic functions.
The case of infinitely many subharmonic functions requires
the extra assumption of upper semicontinuity:

\begin{lemma}\label{UPSEUB}
Let~${\mathcal{I}}$ be a set of indices and let~$u_i$ be subharmonic
in~$\Omega$ for all~$i\in{\mathcal{I}}$.
Let~$u:=\sup_{i\in{\mathcal{I}}}u_i$ and suppose that~$u$ is bounded above.
Then, the upper semicontinuous envelope of~$u$ is subharmonic in~$\Omega$.

Also, if~$u$ is
upper semicontinuous, then~$u$ is also subharmonic in~$\Omega$.
\end{lemma}

\begin{proof} We recall that the
upper semicontinuous envelope~$U_u$ of~$u$
was introduced in~\eqref{S-asle0o4}.
Up to reducing to a connected component, we can suppose that~$\Omega$
is connected.
We can assume that at least one~$u_i$ is not identically equal to~$-\infty$
in~$\Omega$, otherwise~$u$ would also be~$-\infty$
in~$\Omega$, thus trivially leading to the desired result.
Hence, by
Theorem~\ref{SUBHS-EQ}(vi), we can assume that
\begin{equation}\label{SKMcpksl1}
u\in L^1_{\rm loc}(\Omega).
\end{equation}
We take
a connected open set~$\Omega'\Subset\Omega$
and~$r\in(0,{\rm dist}(\Omega',\partial\Omega))$.
Let~$x\in\Omega'$, then, for all~$i\in{\mathcal{I}}$,
$$ u_{i}(x)\le\fint_{ B_r(x)}u_{i}(y)\,dy\le\fint_{ B_r(x)}u(y)\,dy,$$
and therefore, taking the supremum,
\begin{equation}\label{OZkd9oio3c3} u(x)\le\fint_{ B_r(x)}u(y)\,dy.\end{equation}
This is not sufficient to state that~$u$ is subharmonic (unless~$u$
is upper semicontinuous), nevertheless, 
for all~$x\in\Omega'$, we deduce from~\eqref{OZkd9oio3c3} that~$u(x)\le v_r(x)$,
where
$$ v_r(x):=\fint_{ B_r(x)}u(y)\,dy.$$
Also, by~\eqref{SKMcpksl1} and the Dominated Convergence Theorem,
we infer that~$v_r$ is continuous.
Consequently, by the definition of
upper semicontinuous envelope~$U_u$ in~\eqref{S-asle0o4},
we deduce that~$U_u\le v_r$ in~$\Omega'$.
{F}rom this and~\eqref{Semxeeq}
we obtain that, for all~$x\in\Omega'$,
\begin{equation}\label{PALlospo627fa7212erfre} U_u(x)\le v_r(x)=\fint_{ B_r(x)}u(y)\,dy\le
\fint_{ B_r(x)}U_u(y)\,dy.\end{equation}
This and the fact that~$U_u$ is upper semicontinuous
(recall Lemma~\ref{ENCEBENFD}) show that~$U_u$ is subharmonic.

Moreover, if~$u$ is
upper semicontinuous, then~$u$ coincides with~$U_u$, thanks to Lemma~\ref{ENCEBENFD},
and therefore~$u$ is subharmonic as well.
\end{proof}

The assumption that~$u$ is bounded above
in Lemma~\ref{UPSEUB} is taken to avoid problematic examples
such as the following one:
consider the fundamental solution~$\Gamma$
in~\eqref{GAMMAFU}, extended as~$\Gamma(0):=+\infty$,
and let~$u_j(x):=-\Gamma(x)+j$.
Then, by~\eqref{Coamscfunconftahf}, we know that~$u_j$ is subharmonic,
but
$$ \sup_j u_j(x)=\begin{dcases}
-\infty&{\mbox{ if }}x=0,\\+\infty&{\mbox{ otherwise,}}
\end{dcases}
$$
and for such a function one cannot define an upper semicontinuous envelope
(see the
footnote on page~\pageref{SKM:SIJKDM023o3rtetergsrib}).

Further, the assumption that~$u$ is upper semicontinuous cannot be removed
from Lemma~\ref{UPSEUB}. Indeed, 
when~$n\ge2$, we can\footnote{The case~$n=1$
does not lead to similar
counterexamples, see Lemma~\ref{DO:cD00}
and recall that the supremum of convex functions is convex, hence continuous.
In a sense, the semicontinuity setting is superfluous in dimension~$1$,
since subharmonicity in that case entails automatically continuity.}
define~$u_j(x):=-\frac{\Gamma(x)}{j}$. Recalling~\eqref{Coamscfunconftahf},
we have that~$u_j$
is subharmonic in~$B_1$ and
$$ \sup_j u_j(x)=\begin{dcases}
-\infty&{\mbox{ if }}x=0,\\0&{\mbox{ otherwise,}}
\end{dcases}
$$
which is not upper semicontinuous (however, its
upper semicontinuous envelope is the function identically
zero, which is subharmonic, in agreement with
Lemma~\ref{UPSEUB}).
In spite of this type of examples, it is worth observing
that the ``pathological'' points
in which the supremum of subharmonic functions
does not agree with its upper semicontinuous envelope is negligible in the measure
theoretic sense:

\begin{corollary}\label{MS-098765roikjhgf-9ijHJSKorhmEEMI4366K7ASBDD}
Let~${\mathcal{I}}$ be a set of indices and let~$u_i$ be subharmonic
in~$\Omega$ for all~$i\in{\mathcal{I}}$.
Let~$u:=\sup_{i\in{\mathcal{I}}}u_i$ and
suppose that~$u$ is bounded above.
Then, $u$ coincides almost everywhere in~$\Omega$
with its upper semicontinuous envelope.

In particular, $u$ coincides almost everywhere in~$\Omega$
with a subharmonic function.
\end{corollary}

\begin{proof} As discussed in~\eqref{SKMcpksl1},
we can focus on the case in which~$
u\in L^1_{\rm loc}(\Omega)$.
We let~$U_u$ be the upper semicontinuous envelope of~$u$.
Thus, using~\eqref{PALlospo627fa7212erfre},
$$ U_u(x)\le \fint_{ B_r(x)}u(y)\,dy,$$ hence,
if~$x$ is a Lebesgue density point for~$u$,
\begin{equation*} U_u(x)\le \lim_{r\searrow0}\fint_{ B_r(x)}u(y)\,dy=u(x).
\end{equation*}
{F}rom this and~\eqref{Semxeeq} it follows that~$U_u=u$ at all 
Lebesgue density points for~$u$,
hence almost everywhere
in~$\Omega$ (see e.g.~\cite[Theorem~7.13]{MR3381284}).
Since, in light of Lemma~\ref{UPSEUB}, we already know that~$U_u$
is subharmonic, we obtain as a byproduct that~$u$ coincides almost everywhere in~$\Omega$
with a subharmonic function, as desired.
\end{proof}

For further details about the set at which
the supremum of superharmonic functions may differ from
its upper semicontinuous envelope, see e.g.~\cite[Theorems~3.7.5 and~5.7.1]{MR1801253}.\medskip

Now we analyze the pointwise limits of subharmonic functions:

\begin{theorem}\label{SUBseq}
The pointwise
limit of a decreasing sequence
of subharmonic functions is subharmonic.
\end{theorem}

\begin{proof}
Let~$u_k$ be a decreasing sequence
of subharmonic functions in some domain~$\Omega$
and let~$u$ be its pointwise limit.
By Lemma~\ref{65co99865bSdsaJ-04}, we already know that~$u$ is
upper semicontinuous. 

Now suppose that~$B_r(x)\Subset\Omega$.
Then, by Lemma~\ref{ABOABP}, we have that~$u_1$
is bounded from above in~$B_r(x)$, say~$u_1(y)\le M$ for every~$y\in B_r(x)$.
As a result, also~$u_k(y)\le M$ for every~$y\in B_r(x)$.
We can therefore employ the Monotone Convergence Theorem
(see e.g.~\cite[Theorem~5.32(ii)]{MR3381284}) and conclude that
\begin{equation} \label{PuhnSGGlfus}\lim_{k\to+\infty}
\fint_{B_r(x)} u_k(y)\,dy=\fint_{B_r(x)} u(y)\,dy.\end{equation}
Also, $u_k\ge u_{k+j}$ for all~$j\in\N$, whence~$u_k\ge u$.
This, Theorem~\ref{SUBHS-EQ}(i) and~\eqref{PuhnSGGlfus} yield that
\begin{equation*}
u(x)\le
\lim_{k\to+\infty}u_k(x)
\le\lim_{k\to+\infty}
\fint_{B_r(x)} u_k(y)\,dy=\fint_{B_r(x)} u(y)\,dy.\qedhere\end{equation*}
\end{proof}

\begin{figure}
  \centering
  \includegraphics[width=.5\linewidth]{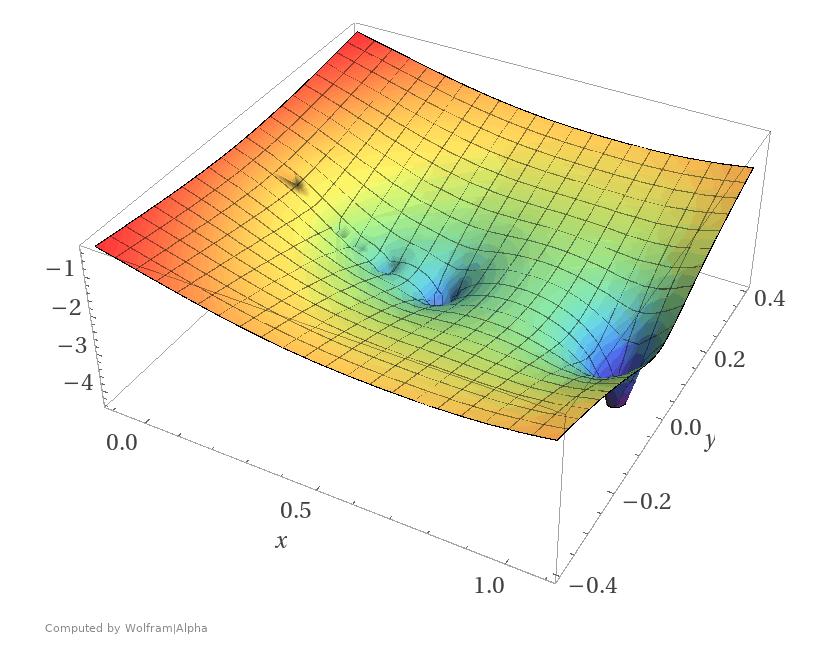}$\qquad$\includegraphics[width=.4\linewidth]{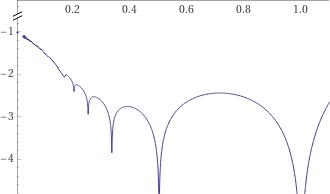}
 \caption{\sl A computer plot of the function~$u_k$ when~$n=2$
 for some ``large''~$k$, and a section of the graph in the $e_2$-direction.}\label{ukkDItangg9845ghjgeFI}
\end{figure}

As a confirming example for Theorem~\ref{SUBseq}, one
can consider the sequence~$u_k:\R^n\to\R$ given by
\begin{equation}\label{ESEMPIETTP} u_k(x):=-\sum_{j=1}^k \frac1{2^j}\,{\Gamma\left(x-\frac{e_n}j\right)},
\end{equation}
where~$\Gamma$ is the fundamental solution in~\eqref{GAMMAFU}
(extended as~$\Gamma(0):=+\infty$ when~$n\ge2$), see Figure~\ref{ukkDItangg9845ghjgeFI}.
Notice that~$u_k$ is subharmonic, thanks to~\eqref{Coamscfunconftahf},
the sequence~$u_k$ is decreasing, and
\begin{equation}\label{OSK0seues}\lim_{k\to+\infty}u_k(x)=
u(x):=
-\sum_{j=1}^{+\infty} \frac1{2^j}\,{\Gamma\left(x-\frac{e_n}j\right)}. \end{equation}
We observe that
$$ |u(0)|\le\sum_{j=1}^{+\infty} \frac{j^{n-2}+\ln j}{2^j}<+\infty,$$
while, when~$n\ge2$,
$$ {\mbox{$u(x)=-\infty$ if~$x\in$}} \bigcup_{j=1}^{+\infty}\frac{e_n}j.$$
We observe that~$u$
is upper semicontinuous, due to Lemma~\ref{65co99865bSdsaJ-04},
but not continuous when~$n\ge2$, since
\begin{equation}\label{N333O234Nsdscdjngiqwufvh}
\liminf_{x\to 0}u(x)\le
\lim_{j\to+\infty}
u\left(\frac{e_n}j\right)=-\infty<u(0).\end{equation}
This example\footnote{The example constructed in~\eqref{ESEMPIETTP}
is that of a discontinuous, but unbounded, subharmonic function.
If instead one wants to construct an example
of a locally bounded and discontinuous
subharmonic function, it suffices to consider~$u$
as in~\eqref{OSK0seues} and define~$\widetilde{u}(x):=\max\{u(x),\,
u(0)-1\}$.
In this way, we have that~$\widetilde{u}$ is subharmonic
by Lemma~\ref{POMAX}, and, when~$n\ge2$, 
formula~\eqref{N333O234Nsdscdjngiqwufvh} gets replaced by
$$\liminf_{x\to 0}\widetilde{u}(x)\le
\lim_{j\to+\infty}
\widetilde{u}\left(\frac{e_n}j\right)=\max\{-\infty,\,u(0)-1\}=u(0)-1<u(0)=\widetilde{u}(0).$$
We also observe that when~$n\ge3$ we have that~$u(0)-1\le\widetilde{u}\le0$,
hence~$\widetilde{u}$ is globally bounded.
On the other hand, in dimension~$2$,
one cannot construct global, discontinuous,
and globally bounded subharmonic functions, due to the forthcoming Theorem~\ref{LIOSU}.

The example of a discontinuous 
subharmonic function in dimension~$n\ge2$
does not have a counterpart in dimension~$1$
since one-dimensional subharmonic functions are necessarily convex
(hence continuous), according to Lemma~\ref{DO:cD00}.}
also clarifies a structural reason
to include the upper semicontinuity assumption
in the definition of subharmonic functions
(however, an approach to
subharmonicity without upper semicontinuity
is also possible, see~\cite{MR1454488}; see also~\cite[Chapter~7]{MR1801253}
for further details on possible discontinuity points
of subharmonic functions
and~\cite{MR0222317} for a theory of subharmonic functions
under continuity assumptions).
\medskip

Interestingly, 
the example in~\eqref{OSK0seues} can be generalized to produce
rather surprising and highly discontinuous subharmonic functions.
For instance, if we enumerate the rationals in~$B_1$
as~$\Q\cap B_1=\{q_j\}_{j\in\N}$,
one can define
$$ u_k(x):=-\sum_{j=1}^{+\infty} \frac1{2^j}\,{\Gamma(x-q_j)},$$
and observe that~$u_k$ is subharmonic, thanks to~\eqref{Coamscfunconftahf},
and accordingly
\begin{equation}\label{OSK0seues2} u(x):=
-\sum_{j=1}^{+\infty} \frac1{2^j}\,{\Gamma(x-q_j)}\end{equation}
is subharmonic, due to
Theorem~\ref{SUBseq}. 
We observe that
\begin{equation}\label{ijndiwenowp}
u\in L^1(B_1),\end{equation} since
$$ \int_{B_1}u(x)\,dx\le
\sum_{j=1}^{+\infty} \frac1{2^j}\,\int_{B_1}|\Gamma(x-q_j)|\,dx
=\sum_{j=1}^{+\infty} \frac1{2^j}\,\int_{B_1(q_j)}|\Gamma(y)|\,dy\le
\int_{B_2}|\Gamma(y)|\,dy\,
\sum_{j=1}^{+\infty} \frac1{2^j}<+\infty.
$$
The special feature of~$u$ is that it is
finite almost everywhere in~$B_1$, thanks to~\eqref{ijndiwenowp},
but, when~$n\geq2$, infinite on~$\Q\cap B_1$ (hence in a dense set).
\medskip

\begin{figure}
  \centering
  \includegraphics[width=.7\linewidth]{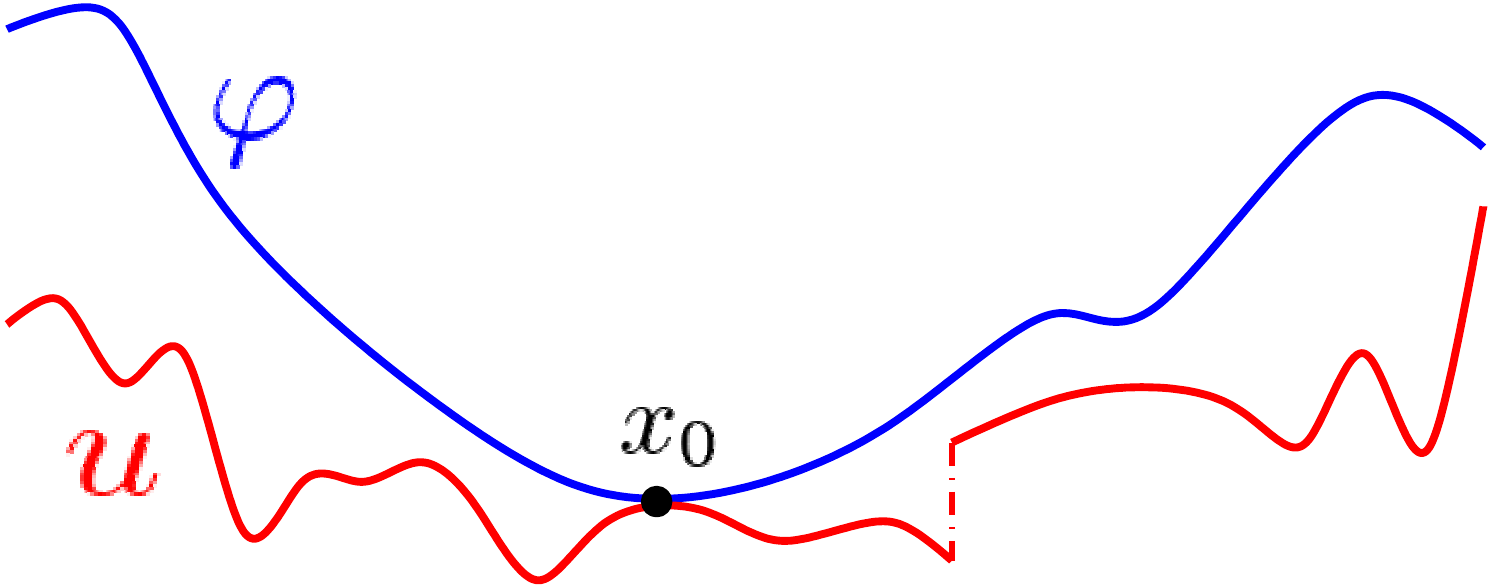}
 \caption{\sl A function~$u-\varphi$ having a local maximum
at~$x_0$ (up to a vertical translation, $\varphi$ is touching~$u$ from above
at~$x_0$).}\label{DviscoItangeFI}
\end{figure}

We now point out a further characterization 
of subharmonic functions motivated by the theory of viscosity 
solutions (see e.g.~\cite{MR1351007}
and Figure~\ref{DviscoItangeFI}).

\begin{theorem}\label{VISCSOL0}
Let~$u:\Omega\to[-\infty,+\infty)$ be upper semicontinuous.
Then, $u$ is subharmonic in~$\Omega$
if and only if, for every~$x_0 \in \Omega$,
every~$\rho>0$ such that~$B_\rho(x_0)\subseteq\Omega$
and every~$\varphi\in C^2(B_\rho(x_0))$
such that the function~$u-\varphi$ has a local maximum
at~$x_0$, it holds that~$\Delta\varphi(x_0)\ge0$.
\end{theorem}

\begin{proof} Assume first that~$u$ is subharmonic
and let~$x_0$, $\rho$ and~$\varphi$
as in the statement of Theorem~\ref{VISCSOL0}.
Let~$r\in(0,\rho)$ and~$h$ be the harmonic
function coinciding with~$\varphi$ along~$\partial B_r(x_0)$
(the existence and continuity of~$h$ being warranted by the Poisson Kernel representation theory
discussed in Theorems~\ref{POIBALL1}
and~\ref{POIBALL}).
Let also~$\widetilde h:=h+(u-\varphi)(x_0)$.
Notice that, if~$y\in\partial B_r(x_0)$,
$$ u(y)=(u-\varphi)(y)+\varphi(y)\le
(u-\varphi)(x_0)+\varphi(y)=
(u-\varphi)(x_0)+h(y)=
\widetilde h(y).$$
Consequently, by
Theorem~\ref{SUBHS-EQ}(vii), we know that~$u\le\widetilde h$
in~$B_r(x_0)$. Therefore, using the Mean Value Property
in Theorem~\ref{KAHAR}(ii),
$$ \varphi(x_0)=u(x_0)+(\varphi-u)(x_0)\le
\widetilde h(x_0)+(\varphi-u)(x_0)=h(x_0)=\fint_{\partial B_r(x_0)}
h(y)\,d{\mathcal{H}}^{n-1}_y
=\fint_{\partial B_r(x_0)}\varphi(y)\,d{\mathcal{H}}^{n-1}_y.$$
{F}rom this and Theorem~\ref{KAHAR2-SPHE}
we deduce that
\begin{equation*}0\le \lim_{r\searrow0} \frac1{r^2}\left(
\fint_{\partial
B_r(x_0)} \varphi(x)\,d{\mathcal{H}}^{n-1}_x-\varphi(x_0)\right)=
\frac{1}{2(n+2)}\,\Delta \varphi(x_0).\end{equation*}

This proves that the first claim in Theorem~\ref{VISCSOL0} implies
the second one.
We now show that, conversely, the second
claim in Theorem~\ref{VISCSOL0} implies
the first one. To this end, 
we take~$x\in\Omega$
and~$r_0>0$ such that~$B_{r_0}(x)\Subset\Omega$.
Let also~$h\in C^2(B_{r_0}(x))\cap
C(\overline{B_{r_0}(x)})$ be a harmonic function
such that~$u\le h$ on~$\partial B_{r_0}(x)$.
In light of Theorem~\ref{SUBHS-EQ}(viii), to prove that~$u$
is subharmonic, and thus complete the proof of
Theorem~\ref{VISCSOL0}, it suffices to check that~$
u\le h$ in~$B_{r_0}(x)$. For this, we argue by contradiction
and suppose that there exist~$x_\star\in B_{r_0}(x)$ and~$a>0$
such that~$h(x_\star)+a\le u(x_\star)$. 
Given~$\e>0$, to be taken sufficiently small, we define
$$\varphi_\e(y):=h(y)-\frac{\e\,|y-x_\star|^2}{n}.$$
We stress that the function~$u-\varphi_\e$ is upper semicontinuous
and therefore, in view of Lemma~\ref{ABOABP},
it attains a maximum in~$\overline{B_{r_0}(x)}$.
Namely, there exists~$X\in\overline{B_{r_0}(x)}$
such that~$(u-\varphi_\e)(X)\ge(u-\varphi_\e)(y)$
for all~$y\in\overline{B_{r_0}(x)}$
and, in particular,
$$ (u-\varphi_\e)(X)\ge(u-\varphi_\e)(x_\star)
=(u-h)(x_\star)
\ge a,$$
while, if~$y\in\partial B_{r_0}(x)$,
$$ (u-\varphi_\e)(y)
=(u-h)(y)+\frac{\e\,|y-x_\star|^2}{n}\le
\frac{\e\,|y-x_\star|^2}{n}\le
\frac{\e\,r_0^2}{n}<a
,$$
as long as~$\e>0$ is sufficiently small.
These observations yield that~$X\not\in\partial B_{r_0}(x)$,
hence~$X$ is a local maximum for~$u-\varphi_\e$.
Since we are assuming that the second statement in
Theorem~\ref{VISCSOL0} is satisfied,
we infer that~$\Delta\varphi_\e(X)\ge0$.
That is,
$$ 0\le\Delta\varphi_\e(X)=\Delta
h(X)-2\e=-2\e.$$
This is a contradiction and the proof of
Theorem~\ref{VISCSOL0} is thereby complete.
\end{proof}

Since subharmonic functions
can exhibit quite complicated behaviors (recall the examples
in~\eqref{OSK0seues} and~\eqref{OSK0seues2}) it is often technically
convenient to reduce the analysis to
smooth subharmonic functions that approach ``as well as possible''
a given subharmonic function. For instance, we have the following result:

\begin{theorem}
Every subharmonic function can be locally approximated by a decreasing
sequence of $C^\infty$ subharmonic functions.
\end{theorem}

\begin{proof} Let~$u$ be subharmonic in~$\Omega$.
By possibly reducing to connecting components, we can suppose that~$\Omega$
is connected. We can also assume that~$u$ is not identically equal to~$-\infty$,
otherwise constant functions will be sufficient for the desired approximation.
As a result,
recalling Theorem~\ref{SUBHS-EQ}(v),
we have that~$u\in L^1_{\rm loc}(\Omega)$,
and we consider the mollification~$u_\eta$ as in~\eqref{CONVETA}.
Given~$\Omega'\Subset\Omega$, if~$\eta>0$ is sufficiently small,
we recall that
\begin{equation}\label{MIbvekdisdcmosn}
{\mbox{$u_\eta\ge u$ in~$\Omega'$,}}\end{equation}
in light of~\eqref{SM:0oedKSKDd35D}. Furthermore,
by~\eqref{KSM0942tgfdvgdsyte}, we know that~$u_\eta$ is subharmonic
in~$\Omega'$.

Additionally, recalling~\eqref{conwsdfisces093mat}
\begin{equation}\label{JSOND-oewdkfm} u_\eta(x)=
\int_0^\eta \tau_\eta(\rho e_1)\left(
\int_{\partial B_\rho}
u(x-\zeta)\,d{\mathcal{H}}^{n-1}_\zeta
\right)\,d\rho,\end{equation}
and therefore,
making use of~\eqref{TGbickdtollibstaiks}
\begin{eqnarray*} u_\eta(x)&\le&\int_0^\eta \tau_\eta(\rho e_1)
\,{\mathcal{H}}^{n-1}(\partial B_r)\,d\rho \;
\fint_{\partial B_\eta(x)} u(\zeta
)\,d{\mathcal{H}}^{n-1}_\zeta
\\&=&\int_0^\eta\left(\int_{\partial B_\rho }
\tau_\eta(y)\,d{\mathcal{H}}^{n-1}_y\right)\,d\rho \;
\fint_{\partial B_\eta(x)} u(\zeta
)\,d{\mathcal{H}}^{n-1}_\zeta\\&=&\int_{B_\eta}\tau_\eta(y)\,d
y\;
\fint_{\partial B_\eta(x)} u(\zeta
)\,d{\mathcal{H}}^{n-1}_\zeta\\&=&\fint_{\partial B_\eta(x)} u(\zeta
)\,d{\mathcal{H}}^{n-1}_\zeta.
\end{eqnarray*}
Consequently, recalling Lemma~\ref{SCELTA},
$$ \limsup_{\eta\searrow0}u_\eta(x)\le
\limsup_{\eta\searrow0}\fint_{\partial B_\eta(x)} u(\zeta
)\,d{\mathcal{H}}^{n-1}_\zeta=u(x).$$
This and~\eqref{MIbvekdisdcmosn} give that~$u_\eta$ converges pointwise
to~$u$ everywhere in~$\Omega'$.

Therefore, to complete the proof of the desired result,
it only remains to prove that~$u_\eta$ is monotone,
namely that~$u_\eta\le u_{\eta'}$ whenever~$\eta\le\eta'$.
For this, we let~$\kappa:=\frac{\eta'}{\eta}\ge1$
and we remark that
$$ \frac1{\kappa^n}\tau_\eta\left( \frac{x}{\kappa}\right)
=\frac1{\kappa^n\eta^n}\tau\left( \frac{x}{\kappa\eta}\right)
=\tau_{\kappa\eta}(x),$$
and consequently, using again~\eqref{JSOND-oewdkfm}
and changing variable~$r:=\kappa \rho$,
\begin{eqnarray*}
u_\eta(x)&=&\frac1\kappa\,\int_0^{\kappa\eta} \tau_\eta\left(\frac{r e_1}\kappa\right)\left(
\int_{\partial B_{r/\kappa}}
u(x-\zeta)\,d{\mathcal{H}}^{n-1}_\zeta
\right)\,dr\\
&=&\frac{1}{\kappa}\,\int_0^{\kappa\eta} \tau_\eta\left(\frac{r e_1}\kappa\right)\,
{\mathcal{H}}^{n-1}(\partial B_{r/\kappa})\,\left(
\fint_{\partial B_{r/\kappa}}
u(x-\zeta)\,d{\mathcal{H}}^{n-1}_\zeta
\right)\,dr
\\&=&\int_0^{\kappa\eta} \tau_{\kappa\eta}(r e_1)\,{\mathcal{H}}^{n-1}(
\partial B_{r})\,\left(
\fint_{\partial B_{r/\kappa}}
u(x-\zeta)\,d{\mathcal{H}}^{n-1}_\zeta
\right)\,dr.
\end{eqnarray*}
This and~\eqref{TGbickdtollibstaiks}, with a further use of~\eqref{JSOND-oewdkfm}, lead to
\[ u_\eta(x)\leq
\int_0^{\kappa\eta} \tau_{\kappa\eta}(r e_1)\,{\mathcal{H}}^{n-1}(
\partial B_{r})\,\left(
\fint_{\partial B_r}
u(x-\zeta)\,d{\mathcal{H}}^{n-1}_\zeta
\right)\,dr=u_{\kappa\eta}(x),\]
which establishes the desired monotonicity property.
\end{proof}

Subharmonic functions in two-dimensional domains enjoy some
special properties. For instance,
the classical Hadamard Three-Circle-Theorem\index{Hadamard Three-Circle-Theorem}
says that the maximum on a circle of radius~$\rho$
of a function that is subharmonic in a planar annulus
is a convex function of~$\ln \rho$. More explicitly, we have that:

\begin{theorem}
Let~$R>r>0$. Let~$\Omega$ be an open subset of~$\R^2$ such that
$$ \big\{x\in\R^2{\mbox{ s.t. }}|x|\in (r,R)\big\}\Subset\Omega.$$
Let~$u$ be subharmonic in~$\Omega$ and, for every~$t\in[\ln r,\ln R]$ let
$$ {\mathcal{M}}(t):=\sup_{\partial B_{e^t}} u.$$
Then, ${\mathcal{M}}$ is a convex function.

Furthermore, if, for~$\rho\in[r,R]$,
\begin{equation}\label{MRHO} M(\rho):= \sup_{\partial B_\rho}u,\end{equation}
we have that, for every~$r_2>\rho>r_1$,
with~$R>r_2>r_1>r$,
\begin{equation}\label{S:3CIT}
M(\rho)\le \frac{
M(r_1)\,(\ln r_2-\ln \rho)
+M(r_2)\,(\ln \rho-\ln r_1)
}{\ln r_2-\ln r_1}.
\end{equation}
\end{theorem}

\begin{proof} For every~$x\in\R^2\setminus\{0\}$, we define
$$ v(x):=u(x)-
\frac{
M(r_1)\,(\ln r_2-\ln |x|)
+M(r_2)\,(\ln |x|-\ln r_1)
}{\ln r_2-\ln r_1}.$$
Since~$\Delta(\ln|x|)=0$ if~$x\in\R^2\setminus\{0\}$, we see that~$v$
is subharmonic in~$\Omega\setminus\{0\}$. Accordingly, by the
Maximum Principle in Corollary~\ref{PERSUBAWEAK},
for every~$x\in A:=\big\{x\in\R^2{\mbox{ s.t. }}|x|\in (r_1,r_2)\big\}$,
\begin{eqnarray*}
&&v(x)\le\sup_{y\in\partial A}v(y)=\max\left\{
\sup_{y\in\partial B_{r_1}}v(y)
,\,\sup_{y\in\partial B_{r_2}}v(y)
\right\}\\&&\qquad\qquad=
\max\left\{
\sup_{y\in\partial B_{r_1}} u(y)-M(r_1),\,
\sup_{y\in\partial B_{r_2}} u(y)-M(r_2)\right\}=0.\end{eqnarray*}
This leads to
\begin{eqnarray*}&& M(\rho)-
\frac{
M(r_1)\,(\ln r_2-\ln\rho)
+M(r_2)\,(\ln \rho-\ln r_1)
}{\ln r_2-\ln r_1}
\\&&\qquad=\sup_{x\in\partial B_\rho}
u(x)-
\frac{
M(r_1)\,(\ln r_2-\ln |x|)
+M(r_2)\,(\ln |x|-\ln r_1)
}{\ln r_2-\ln r_1}
=\sup_{x\in\partial B_\rho}v(x)\le0,\end{eqnarray*}
which is~\eqref{S:3CIT}.

Now, if~$t_1$, $t_2\in[\ln r,\ln R]$, with~$t_2>t_1$,
we let~$r_1:=e^{t_1}$
and~$r_2:=e^{t_2}$. Let also~$\vartheta\in[0,1]$ and~$\rho:=e^{(1-\vartheta)t_1+\vartheta t_2}$. 
We also notice that~${\mathcal{M}}(t)=M(e^t)$.
Then, from~\eqref{S:3CIT},
\begin{eqnarray*}
0&\le& \frac{
M(r_1)\,(\ln r_2-\ln \rho)
+M(r_2)\,(\ln \rho-\ln r_1)
}{\ln r_2-\ln r_1}-M(\rho)
\\&=&(1-\vartheta)\,M(e^{t_1})+\vartheta\,M(e^{t_2})
-M(e^{(1-\vartheta)t_1+\vartheta t_2})
\\&=&
(1-\vartheta)\,{\mathcal{M}}(t_1)+\vartheta\,{\mathcal{M}}(t_2)
-{\mathcal{M}}({(1-\vartheta)t_1+\vartheta t_2})
,\end{eqnarray*}
thus establishing the convexity of~$ {\mathcal{M}}$.
\end{proof}

Likely, the name of Three-Circle-Theorem
is due to inequality~\eqref{S:3CIT} that involves the circles of radius~$r_2$, $\rho$ and~$r_1$.\medskip

An interesting consequence of
the Hadamard Three-Circle-Theorem
is a reinforcement of Liouville's Theorem (see Theorem~\ref{LIOVE})
in dimension~2, as follows:

\begin{theorem}\label{LIOSU}
A function which is subharmonic in the whole of~$\R^2$
and bounded from above is necessarily constant.\index{Liouville's Theorem}
\end{theorem}

\begin{proof}
Let~$u$ be subharmonic in the whole of~$\R^2$
with~$u\le C$ for some~$C>0$, and let~$M$ be as in~\eqref{MRHO}.
On the one hand,
by the Maximum Principle in Corollary~\ref{PERSUBAWEAK}, if~$\rho_2>\rho_1>0$,
$$ M(\rho_2)= \sup_{\partial B_{\rho_2}}u
=\sup_{\overline{B_{\rho_2}}}u\ge
\sup_{\overline{B_{\rho_1}}}u=\sup_{{\partial B_{\rho_1}}}u
=M(\rho_1),$$
hence
\begin{equation}\label{NONDECRM}
{\mbox{$M$ is nondecreasing}}.
\end{equation}
On the other hand, in light of~\eqref{S:3CIT}
(applied here by sending~$r_2\to+\infty$), we see that, if~$\rho>r_1$,
\begin{eqnarray*}&&
M(\rho)\le\lim_{r_2\to+\infty} \frac{
M(r_1)\,(\ln r_2-\ln \rho)
+M(r_2)\,(\ln \rho-\ln r_1)
}{\ln r_2-\ln r_1}\\&&\qquad\qquad
\le \lim_{r_2\to+\infty} \frac{
M(r_1)\,(\ln r_2-\ln \rho)
+C\,(\ln \rho-\ln r_1)
}{\ln r_2-\ln r_1}= M(r_1).\end{eqnarray*}
This and~\eqref{NONDECRM} yield that, for all~$\rho>r_1$,
$$ \sup_{\overline{B_{\rho}}}u=M(\rho)=M(r_1)=
\sup_{\overline{B_{r_1}}}u$$
and therefore, by the Strong Maximum Principle in Lemma~\ref{LplrfeELLtahsmedschdfb4camntM},
we have that~$u$ is constant.
\end{proof}

We observe that the assumption on the boundedness from above
in Theorem~\ref{LIOSU} cannot be replaced by a boundedness from
below (for example, the function~$u(x)=|x|^2$ is subharmonic
in~$\R^2$ and bounded below, without being constant). Notice that
this is an interesting structural difference with respect to
the Liouville's Theorem (Theorem~\ref{LIOVE}).

Furthermore, Theorem~\ref{LIOSU} 
holds true also in dimension~1, thanks to the convexity equivalence in Lemma~\ref{DO:cD00},
but it does not carry over to dimension~$n\ge3$
(not even assuming both a bound from above and from below). For instance, for~$n\ge3$,
the function~$u(x):=-\frac1{\left(1+|x|^2\right)^{\frac{n-2}2}}$ satisfies
$$ \Delta u(x)=\frac{n(n-2)}{\left(1+|x|^2\right)^{\frac{n-2}2}}\ge0,$$
hence~$u$ is subharmonic everywhere, bounded, but not constant.\medskip

For further readings on subharmonic functions,
see also~\cite{MR0344979, MR0460672, MR1301332, MR2569331, MR3675703} and the references therein.

For the interesting connections between the Three-Circle-Theorem and the Riemann Conjecture see~\cite[page~188]{MR1854455},
\cite[page~332]{MR882550} and~\cite[pages~70--72]{MR933558}.

\begin{figure}
                \centering
                \includegraphics[width=.7\linewidth]{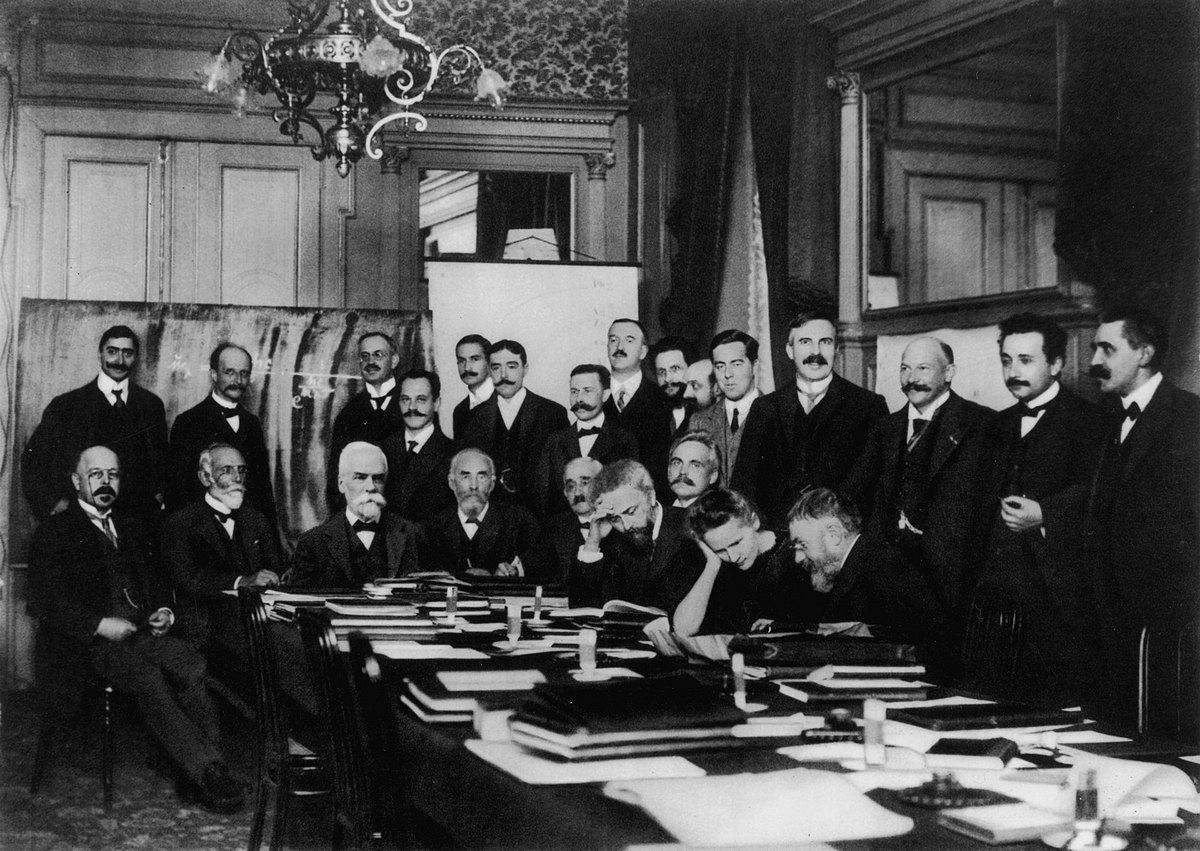}
        \caption{\sl 
        Photograph of the 
        1911 Solvay Conference on Physics.
         Seated: W.~Nernst, M.~Brillouin, E.~Solvay, H.~Lorentz, E.~Warburg, J.~Perrin, W.~Wien, M.~{S}k\l odowska Curie
         and H.~Poincar\'e. Standing: R.~Goldschmidt, M.~Planck, H.~Rubens, A.~Sommerfeld, F.~Lindemann, M.~de Broglie, M.~Knudsen, F.~Hasen\"ohrl, G.~Hostelet, E.~Herzen, J.~H.~Jeans, E.~Rutherford, H.~Kamerlingh Onnes, A.~Einstein and P.~Langevin
        (Public Domain image from
        Wikipedia).}\label{HAFOUMSldRGIRA4AXELHARLROA7789GIJ7solFUMHDNOJHNFOJEDSOL1}
\end{figure}

\medskip

In connection with Maximum Principles, we point out some further properties of the
Green Function:

\begin{lemma}
Let~$n\ge2$ and consider an open set~$\Omega\subset\R^n$
with~$C^1$ boundary.

Let~$G$ be the Green Function of~$\Omega$, as introduced in Section~\ref{Green Function:SECT}. Let also~$\Gamma$ be the fundamental solution, as introduced in Section~\ref{lfundsP-S}.

Then, for all~$x$, $x_0\in\Omega$,
\begin{equation}\label{GREENPOS01we}G(x,x_0)\ge0\end{equation}
and
\begin{equation}\label{GREENPOS01we-2}G(x,x_0)\le\Gamma(|x-x_0|).\end{equation}
\end{lemma}

\begin{proof}Using the notation in~\eqref{COIENEE-0986ytufgkbv-0-rjfeonvnb2GDB}
and~\eqref{GAROGRA} together with~\eqref{FLGREEN} and Lemma~\ref{SYJMDMMAD} we have that
$$ G(x_0,x)=G(x,x_0)=\lim_{\rho\searrow0}G_\rho(x,x_0),$$
where
\begin{equation}\label{GREENPOS01we-20} G_\rho(x,x_0):=\Gamma_\rho(x-x_0)-\Psi^{(x_0)}(x).\end{equation}

Also, when~$\rho$ is smaller than the distance of~$x_0$ from~$\partial\Omega$,
for all~$x\in\partial\Omega$ we have that~$\Gamma_\rho(x-x_0)=\Gamma(x-x_0)$
and consequently~$G_\rho(\cdot,x_0)=0$ along~$\partial\Omega$.

Since~$G_\rho$ is superharmonic, due to~\eqref{DEGAM} and~\eqref{L:S:ASM S}, we infer from the Maximum Principle in Corollary~\ref{PERSUBAWEAK} that, for every~$x\in\Omega$,
$$ G_\rho(x,x_0)\ge\inf_{y\in\partial\Omega}G_\rho(y,x_0)=0$$
and therefore~\eqref{GREENPOS01we} plainly follows.

Also, from~\ref{L:S:ASM S} and the Maximum Principle we infer that~$\Psi^{(x_0)}(x)\ge0$ and accordingly, by~\eqref{GREENPOS01we-20}, we find that~$G_\rho(x,x_0)\le\Gamma_\rho(x-x_0)$. Passing to the limit in~$\rho$, this gives~\eqref{GREENPOS01we-2}, as desired.
\end{proof}

In view of the above observations, we now point out a bound in Lebesgue spaces.
For more general results, see~\cite[Theorems~8.15, 8.16, 8.17, 8.18, 8.25 and~8.26]{MR1814364} (here, the situation is much simpler than the general case,
since we are focusing on functions vanishing along the boundary of a given domain).

\begin{proposition}\label{MALEBEAPP}
Let~$n\ge2$ and~$\Omega$ be a bounded and open subset of~$\R^n$ with boundary of class~$C^1$.
Let~$p>\frac{n}2$ and~$f:\Omega\to \R$ with~$f^-\in L^p(\Omega)$.

Let~$u\in C^2(\Omega)\cap C^1(\overline\Omega)$ be a solution of
$$ \begin{dcases}\Delta u\ge f&{\mbox{ in }}\Omega,\\
u=0&{\mbox{ on }}\partial\Omega.\end{dcases}
$$
Then,
$$ \sup_{\Omega }u\le C\|f^-\|_{L^p(\Omega)},$$
for some~$C>0$ depending only on~$n$, $p$ and $\Omega$.
\end{proposition}

\begin{proof} Up to considering connected components, we may suppose that~$\Omega$ is connected.
Let~$x_0\in \Omega'$. By 
Theorem~\ref{MS:SKMD344567yjghS},
\begin{equation} \label{01-23iorhg02i3ruohg4152ty3ruj981yt29wtty987t0987678iekjf}
-u(x_0)=\int_\Omega G(x_0,y)\,\Delta u(y)\,dy,\end{equation}
where~$G$ is the Green Function of~$\Omega$.

On that account, we infer from~\eqref{GREENPOS01we}
and~\eqref{01-23iorhg02i3ruohg4152ty3ruj981yt29wtty987t0987678iekjf} that
\begin{equation*}
-u(x_0)\ge\int_\Omega G(x_0,y)\,f(y)\,dy\ge-\int_\Omega G(x_0,y)\,|f^-(y)|\,dy.\end{equation*}

%We also observe that
%$$ \sup_{x\in\partial\Omega}|\Gamma(x-x_0)|\le C,$$
%for some~$C>0$ depending only on~$n$, $\Omega$ and~$\Omega'$.
%
%Hence, by the Maximum Principle, also~$|\Psi^{(x_0)}|\le C$ in~$\Omega$, leading to
%$$ |G_\rho(x,x_0)|\le |\Gamma(x-x_0)|+|\Psi^{(x_0)}(x)|\le
%|\Gamma(x-x_0)|+C\le C|\Gamma(x-x_0)|,$$
%up to renaming~$C$ at each step of our calculations.

Additionally, we observe that~$\Gamma$ is locally in~$L^q$, where~$q\in[1,+\infty)$ is the dual exponent of~$p$, thanks to the assumption on the range of~$p$, and consequently, by~\eqref{GREENPOS01we-2} and the H\"older Inequality,
\begin{equation*}
u(x_0)\le\int_\Omega G(x_0,y)\,|f^-(y)|\,dy\le
\int_\Omega \Gamma(x_0-y)\,|f^-(y)|\,dy\le C\,\|f^-\|_{L^p(\Omega)},
\end{equation*}
for some~$C>0$ depending only on~$n$, $p$ and $\Omega$.
\end{proof}

One can compare Proposition~\ref{MALEBEAPP} with Theorem~\ref{CORODASOLTFJ}
which dealt with harmonic functions, but without any assumption on the boundary data.

It is interesting to observe that the exponent~$p$ in Proposition~\ref{MALEBEAPP}
is\footnote{For completeness, we mention that in the case of operators in nondivergence form
the threshold for the exponent~$p$ somehow switches from~$\frac{n}2$ to~$n$,
see~\cite[Theorems~9.20 and~9.22]{MR1814364}
and~\cite[Theorem~3.2]{MR1351007}.
The search for the optimal exponent for general elliptic equations
is an active field of research related to a conjecture by Carlo Pucci~\cite{MR208150}, see e.g.~\cite{MR2472875, MR4077156}
for further details.} essentially optimal. As an example, one can consider~$\rho\in\left(0,\frac14\right)$ and the approximation of the fundamental solution~$\Gamma_\rho$ as in~\eqref{GAROGRA}, and thus define~$u_\rho(x):=\Gamma_\rho(x)-\Gamma_\rho(e_1)$. In this way, we have that~$u_\rho\in C^{1,1}(\overline{B_1})$ with~$u_\rho=0$ on~$\partial B_1$ and, recalling~\eqref{DEGAM},
$$ \Delta u_\rho=\begin{dcases}
-\frac{1}{|B_\rho|}&{\mbox{ in }}B_\rho,\\
0&{\mbox{ in }}\R^n\setminus\overline{B_\rho}.
\end{dcases}$$
In this situation, recalling the notation in~\eqref{OVERLINEV}, we see that
\begin{eqnarray*}&&
\frac{\sup_{B_{1/2}} u_\rho}{\|(\Delta u_\rho)^-\|_{L^{n/2}(B_1)}}=
\frac{u_\rho(0)}{|B_\rho|^{\frac{2-n}n}}=
\frac{1}{|B_\rho|^{\frac{2-n}n}}\left(
\frac{\rho^2}{2n |B_\rho|}+c_n \big(\overline v(\rho)-\overline v(1)\big)
\right).
\end{eqnarray*}
In particular, when~$n=2$,
\begin{eqnarray*}&&
\frac{\sup_{B_{1/2}} u_\rho}{\|(\Delta u_\rho)^-\|_{L^{n/2}(B_1)}}=\frac{1}{4\pi}-c_2 \ln\rho
,\end{eqnarray*}
which diverges as~$\rho\searrow0$, and this implies the optimality of the exponent~$p$ in Proposition~\ref{MALEBEAPP}.
\medskip

Also when~$n\ge3$, this exponent cannot be improved, as showcased by the following example.
For~$\e\in(0,1)$ we consider
$$ f_\e(x):=\begin{dcases}-\displaystyle\frac{1}{|x|^{2}|\ln (|x|/2)|}&{\mbox{ if }}x\in B_1\setminus B_\e,\\ 0 &{\mbox{ if }}x\in B_\e\end{dcases}$$
and we define
$$ B_1\ni x\longmapsto u_\e(x)=-\int_{B_1} f_\e(y)\,
\left[\Gamma(x-y)-\Gamma\left( |y|\,\left(
x-\frac{y}{|y|^2}\right)\right)\right]\,dy.$$
By Corollary~\ref{SKMD:DOK0oerlfgdit} and Theorem~\ref{MS:SKMD344567yjghS-1}, we have that
$$\begin{dcases}
\Delta u_\e=f_\e&{\mbox{ in }}B_1,\\
u_\e=0&{\mbox{ on }}\partial B_1.
\end{dcases}$$
Notice that, on the one hand, if~$n\ge3$,
\begin{eqnarray*}
\|f_\e^-\|_{L^{n/2}(B_1)}&=&
\left( 
\int_{B_1\setminus B_\e} \frac{1}{|x|^{n}|\ln (|x|/2)|^{\frac{2}n}}\,dx
\right)^{\frac{2}n}\\&\le& C
\left( 
\int_\e^1 \frac{1}{r |\ln (r/2)|^{\frac{n}2}}\,dr \right)^{\frac{2}n}\\&=& C
\left( 
\int_{\ln2}^{-\ln(\e/2)} \frac{1}{t^{\frac{n}2}}\,dt\right)^{\frac{2}n}\\&\le& C
\left( 
\int_{\ln2}^{+\infty} \frac{1}{t^{\frac{n}2}}\,dt\right)^{\frac{2}n}\\&\le& C
,\end{eqnarray*}
up to renaming~$C$ at each line.

On the other hand, if~$n\ge3$ and~$\rho>0$ is small enough (and~$\e$ is appropriately small possibly in dependence of~$\rho$),
\begin{eqnarray*}u_\e(0)&=&-\int_{B_1} f_\e(y)\,
\big[\Gamma(y)-\Gamma(e_1)\big]\,dy\\&\ge&
\int_{B_{\rho}\setminus B_\e}
\frac{\Gamma(y)-\Gamma(e_1)}{|y|^{2}|\ln (|y|/2)|}\,dy\\&\ge&
\frac1{C}\int_{B_\rho\setminus B_\e}
\frac{dy}{|y|^{n}|\ln (|y|/2)|}
\\&=&\frac1{C}\int_{\e}^\rho
\frac{dr}{r|\ln (r/2)|}\\&=&\frac1{C}\int_{-\ln(\rho/2)}^{-\ln(\e/2)}
\frac{dt}{t}\\&\ge& \frac{\ln|\ln(\e/2)|}C.
\end{eqnarray*}
{F}rom these observations we arrive at\begin{eqnarray*}&&
\frac{\sup_{B_{1/2}} u_\e}{\|(\Delta u_\e)^-\|_{L^{n/2}(B_1)}}\ge\frac{\ln|\ln(\e/2)|}C
,\end{eqnarray*}
which diverges as~$\e\searrow0$, and this implies the optimality of the exponent~$p$ in Proposition~\ref{MALEBEAPP} also when~$n\ge3$.
\medskip

An extension of Proposition~\ref{MALEBEAPP} not assuming vanishing data along the boundary will be considered in Proposition~\ref{MALEBEAPP-BIS}, exploiting
(differently from the classical literature) the forthcoming Calder\'{o}n-Zygmund estimates
presented in Section~\ref{SEC:CaldZygmund estimates} (actually, only the conceptually
simpler portion of these estimates in class~$W^{1,p}$).

\section{The Perron method}

\begin{figure}
                \centering
                \includegraphics[width=.47\linewidth]{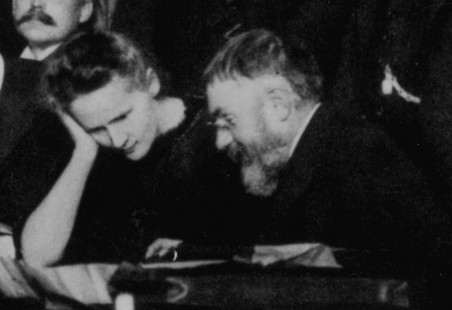}
        \caption{\sl 
        Detail from Figure~\ref{HAFOUMSldRGIRA4AXELHARLROA7789GIJ7solFUMHDNOJHNFOJEDSOL1}
        (Public Domain image from
        Wikipedia).}\label{HAFOUMSldRGIRA4AXELHARLROA7789GIJ7solFUMHDNOJHNFOJEDSOL2}
\end{figure}

We address here the problem of existence of harmonic functions
with prescribed boundary conditions. Several approaches can be taken for this problem
and for more general ones,
leveraging variational methods and functional analysis
(see e.g. Sections~2.2.5(b) and~6.2 in~\cite{MR1625845}).
Here, we follow instead a classical intuition by\footnote{Poincar\'e
is often considered the last universalist,
i.e. the last scientist who was capable of
deeply understanding and
revolutionizing the set of contemporary knowledge seen as a whole.
He was a mathematician, theoretical physicist, engineer and philosopher of science,
often relying on his outstanding intuition and capability of visual representation.
He laid the foundations of chaos theory and topology.
His work contributed to the birth of special relativity, having detected relativistic velocity transformations 
and a perfect invariance of all of Maxwell's equations\index{Maxwell's equations} and having proposed the existence
of gravitational waves.

One of his conjectures became one of the \label{106106-bisPE}
Millennium Prize Problems (see footnote~\ref{106106} on page~\pageref{106106}),
and actually the first (and, up to today, only) question to receive a complete answer
(by Grigori Yakovlevich, a.k.a. Grisha, Perelman~\cite{2003math7245P}, who
rejected the prize of ~$10^6$ US \$).

We think that a nice example of enthusiasm and scientific fervour
is provided in Figures~\ref{HAFOUMSldRGIRA4AXELHARLROA7789GIJ7solFUMHDNOJHNFOJEDSOL1} and~\ref{HAFOUMSldRGIRA4AXELHARLROA7789GIJ7solFUMHDNOJHNFOJEDSOL2}: rather than posing
for the official photograph, Marie {S}k\l odowska (discoverer of the elements polonium and radium,
1903 Nobel Prize in Physics, 1911
Nobel Prize in Chemistry)
and Henri Poincar\'e
keep intensively discussing Science.

By the way, Poincar\'e's multifaceted scientific talents were for his peers allegedly as proverbial as
his clumsiness and his inability of drawing:
it seems that
his lack of drawing skill might even have jeopardized Poincar\'e's admittance to the \'Ecole Polytechnique,
due to a perfect 0 in the exam of wash drawing
(fortunately, his examiners were sufficiently magnanimous or impressed by the marks obtained
in the mathematical tests to change the 0 into a 1, which was sufficient for an overall pass).

In spite of some poor drawings by Poincar\'e, which the reader can find e.g. in~\cite{BART},
Poincar\'e heavily relied on figures to favor intuition and claimed that
``figures first of all make up for the infirmity of our intellect by calling on the aid of our senses; but not only this. It is worth repeating that geometry is the art of reasoning well from badly drawn figures'', see~\cite{zbMATH02678889}.

Poincar\'e was certainly conscious of his limits as an illustrator: when, during his studies of celestial mechanics
in~\cite{MR0087814},
he discovered the complex web arising from
the intersections of stable and unstable manifolds (ultimately leading to chaotic patterns in dynamical systems), Poincar\'e
wrote ``The complexity of this figure, which I will not even attempt to draw, is striking''.} Henri Poincar\'e~\cite{MR1505534},
as developed by Oskar Perron~\cite{MR1544619}
and Robert Remak~\cite{MR1544666}, formalized in what
is nowadays called the ``Perron method''. In a nutshell,
the idea is to obtain a harmonic function by an increasing
sequence of subharmonic ones, in a way which is similar
to that of creating a stable electrostatic potential by sweeping out charges from inside the domain.
The bottom line of the method is therefore to consider the ``largest
subharmonic function'', check that it is indeed harmonic and
that, under suitable conditions, it meets the boundary datum.

{F}rom the technical point of view, an important step in this construction
is the fact that we can already construct harmonic functions
given a boundary datum in balls, thanks to the Poisson Kernel method
developed in Theorems~\ref{POIBALL1} and~\ref{POIBALL} (this will allow us
to ``lift'' a subharmonic function in a ball by replacing it with the harmonic
function with the same boundary data). The construction also relies on
the Maximum Principle, which, among the other useful information,
allows one to conclude that the above mentioned harmonic lifting,
when glued together to a subharmonic function,
produces a subharmonic function as well.

Here are the technical details related to\index{harmonic lift}
the above mentioned ``harmonic lift'' technique:

\begin{lemma}\label{L:LIFT}
Let~$\Omega$ be an open and bounded
subset of~$\R^n$. Let~$R>0$ and assume that~$B_R\Subset\Omega$.
Let~$w\in C(\Omega)$ be subharmonic in~$\Omega$.

Then, there exists a unique function~$W\in C(\Omega)$ 
which is harmonic in~$B_R$, subharmonic in~$\Omega$ and such that
\begin{equation}\label{E:LIFT}
\begin{dcases}{\mbox{$W=w$ \; in \;$\Omega\setminus B_R$}},\\
{\mbox{$W\geq w$\; in \;$\Omega$.}}
\end{dcases}\end{equation}
\end{lemma}

\begin{proof}
The uniqueness claim follows by
applying Corollary~\ref{UNIQUENESSTHEOREM} in~$ B_R$.

We now focus on the existence claim. To this end, we utilize the Poisson Kernel of~$B_R$,
in the light of Theorems~\ref{POIBALL1} and~\ref{POIBALL}. In this way,
we can find a harmonic function~$h$ in~$B_R$
with datum~$w$ along~$\partial B_R$. We observe that
\begin{equation}\label{SMD-OSKMDM-O-kmdfgfg}
{\mbox{$h$ is continuous in~$\overline{B_R}$}}\end{equation}
since if~$\zeta\in\partial B_R$ and~$\eta\in B_R$ (with~$\rho:=|\eta-\zeta|$ to be taken conveniently small), then, given~$\delta\in\left[2\rho,\frac{R}8\right]$,
\begin{eqnarray*}&&
|h(\eta)-w(\zeta)|\\&=&\left|\int_{\partial B_R} w(y)\frac{R^2-|\eta|^2}{n\,|B_1|\,R\,|\eta-y|^n}\,d{\mathcal{H}}^{n-1}_y
-w(\zeta)\right|\\&\le&
\int_{\partial B_R} |w(y)-w(\zeta)|\frac{R^2-|\eta|^2}{n\,|B_1|\,R\,|\eta-y|^n}\,d{\mathcal{H}}^{n-1}_y\\
\\&\le& C\,(R^2-|\eta|^2)\left[
\int_{(\partial B_R)\cap B_\delta(\eta)} |w(y)-w(\zeta)|\frac{d{\mathcal{H}}^{n-1}_y}{|\eta-y|^n}+\frac1{\delta^n}
\int_{\partial B_R\setminus B_\delta(\eta)} |w(y)-w(\zeta)|\,d{\mathcal{H}}^{n-1}_y
\right]\\&\le& C\,(R-|\eta|)\left[\sup_{x\in B_{2\delta}(\zeta)}|w(x)-w(\zeta)|
\int_{(\partial B_R)\cap B_\delta(\eta)} \frac{d{\mathcal{H}}^{n-1}_y}{|\eta-y|^n}+\frac1{\delta^n}
\right]\\&\le& C\sigma\left[
\sup_{x\in B_{2\delta}(\zeta)}|w(x)-w(\zeta)|
\sum_{j=0}^{\ln(\delta/\sigma)}
\int_{(\partial B_R)\cap (B_{\delta/e^j}(\eta)\setminus B_{\delta/e^{j+1}}(\eta))} \frac{d{\mathcal{H}}^{n-1}_y}{|\eta-y|^n}
+\frac1{\delta^n}
\right]\\&\le& C\sigma\left[
\sup_{x\in B_{2\delta}(\zeta)}|w(x)-w(\zeta)|
\sum_{j=0}^{\ln(\delta/\sigma)}\frac{(\delta/e^j)^{n-1}}{(\delta/e^{j+1})^n}
+\frac1{\delta^n}
\right]\\&\le& C\sigma\left[
\sup_{x\in B_{2\delta}(\zeta)}|w(x)-w(\zeta)|
\sum_{j=0}^{\ln(\delta/\sigma)}\frac{ e^j }{\delta}
+\frac1{\delta^n}
\right]\\&\le& C\sigma\left[\frac1\sigma
\sup_{x\in B_{2\delta}(\zeta)}|w(x)-w(\zeta)|
+\frac1{\delta^n}
\right]
\end{eqnarray*}
for some~$C>0$ depending only on~$n$, $R$ and~$\|w\|_{L^\infty(\partial B_R)}$,
and possibly varying from line to line, where~$\sigma:=R-|\eta|$ denotes the distance between~$\eta$ and~$\partial B_R$.

That is, taking a sequence~$\eta_m\in B_R$ such that~$\eta_m\to\zeta\in\partial B_R$ as~$m\to+\infty$,
letting~$\sigma_m:=R-|\eta_m|$ and noticing that~$\sigma_m\to0$,
\begin{eqnarray*}
\lim_{m\to+\infty}|h(\eta_m)-w(\zeta)|
&\le&\lim_{m\to+\infty}
C\,\left[\sup_{x\in B_{2\delta}(\zeta)}|w(x)-w(\zeta)|
+\frac{\sigma_m}{\delta^n}
\right]\\
&=&C\,\sup_{x\in B_{2\delta}(\zeta)}|w(x)-w(\zeta)|,
\end{eqnarray*}
for all~$\delta\in\left(0,\frac{R}8\right]$. We can thus send~$\delta\searrow0$ and conclude that
\begin{eqnarray*}
\lim_{m\to+\infty}|h(\eta_m)-w(\zeta)|=0,\end{eqnarray*}
which proves~\eqref{SMD-OSKMDM-O-kmdfgfg}.

Now we define
\begin{equation}\label{equ8itrh547hrhsfwetn} W(x):=\begin{dcases}
h(x)&{\mbox{ if }}x\in B_R,\\
w(x)&{\mbox{ if }}x\in \Omega\setminus B_R
\end{dcases}\end{equation}
and we remark that
\begin{equation}\label{JE4CUSNSPENODSLSEND-1}W\in C(\overline\Omega),\end{equation}
owing to~\eqref{SMD-OSKMDM-O-kmdfgfg}.

Also, the function~$w-h$ is subharmonic in~$B_R$. Hence, by Corollary~\ref{PERSUBAWEAK},
$$\sup_{B_R}(w-h)=\lim_{r\nearrow R}\sup_{B_r} (w-h)=
\lim_{r\nearrow R}\sup_{\partial B_r} (w-h)=\sup_{\partial B_R} (w-h)=0.$$
As a consequence, we have that~$w\le h$ in~$B_R$, whence
\begin{equation}\label{JE4CUSNSPENODSLSEND}
{\mbox{$W\ge w$ in~$\Omega$.}}\end{equation}
We also remark that
\begin{equation}\label{JE4CUSNSPENODSLSEND-2}
{\mbox{$W$ is subharmonic in~$\Omega$.}}\end{equation}
To check this, let~$\widetilde{x}\in\Omega$
and~$\widetilde{r}>0$ such that~$\widetilde{B}:=B_{\widetilde{r}}(\widetilde{x})\Subset\Omega$.
Consider a harmonic function~$\widetilde{h}$ which is continuous in the closure of~${\widetilde{B}}$
such that~$W\le \widetilde{h}$ on~$\partial \widetilde{B}$.
By~\eqref{JE4CUSNSPENODSLSEND}, we know that~$w\le W\le \widetilde{h}$ on~$\partial \widetilde{B}$.
Hence, recalling Theorem~\ref{SUBHS-EQ}(viii)
(as well as Definition~\ref{SHAITY}), we deduce that~$w\le \widetilde{h}$ in the whole of~$ \widetilde{B}$.
As a result,
\begin{equation}\label{JE4CUSNSPENODSLSEND-19093}
{\mbox{in~$\widetilde{B}\setminus B_R$ we have that~$W=w\le \widetilde{h}$.}}\end{equation}
Additionally, since~$\partial(\widetilde{B}\cap B_R)$
is contained in
the closure of~$\widetilde{B}$, we know that~$W\le \widetilde{h}$
on~$\partial(\widetilde{B}\cap B_R)$ and therefore,
using again the Maximum Principle of Corollary~\ref{PERSUBAWEAK} (up to the boundary of~$\widetilde{B}\cap B_R$),
we find that
$$ \sup_{\widetilde{B}\cap B_R}(W- \widetilde{h})=\sup_{\partial(\widetilde{B}\cap B_R)}(
W- \widetilde{h})\le0.$$
Combining this and~\eqref{JE4CUSNSPENODSLSEND-19093}
we find that~$W\le \widetilde{h}$
in~$\widetilde{B}$. This and
Theorem~\ref{SUBHS-EQ}(viii) lead to~\eqref{JE4CUSNSPENODSLSEND-2}, as desired.

The existence claim
is thus a consequence of~\eqref{equ8itrh547hrhsfwetn},
\eqref{JE4CUSNSPENODSLSEND-1},
\eqref{JE4CUSNSPENODSLSEND} and~\eqref{JE4CUSNSPENODSLSEND-2}.
\end{proof}

With this, we can address the Perron method\index{Perron method} for constructing harmonic functions,
as follows:
 
\begin{theorem}\label{PERRO}
Let~$\Omega\subset\R^n$ be open and bounded.
Let~$g$ be a bounded function on~$\partial\Omega$ and define
\begin{equation}\label{KL:PER} u(x):=\sup v(x),\end{equation}
where the above supremum is taken over all functions~$v\in C(\Omega)$
which are subharmonic in~$\Omega$ and such that 
\begin{equation}\label{VSOTTO} \limsup_{\Omega\ni x\to p}v(x)\le g(p)\qquad{\mbox{for every~$p\in\partial\Omega$}}.\end{equation}

Then, $u$ is harmonic in~$\Omega$.
\end{theorem}

\begin{proof} Without loss of generality, up to restricting ourselves to connected
components, we may suppose that~$\Omega$ is connected. We show that
if~$v\in C(\Omega)$
is subharmonic in~$\Omega$ and satisfies~\eqref{VSOTTO}, then
\begin{equation}\label{KSM-0987rpolkjbg9uIGVS56789UJHDD} \sup_{\Omega} v\le \sup_{\partial\Omega}g.\end{equation}
For this, we consider a sequence of points~$p_j\in\Omega$ such that
$$ \lim_{j\to+\infty} v(p_j)= \sup_{\Omega} v.$$
Up to a subsequence, we can suppose that~$p_j\to p_\star$ as~$j\to+\infty$,
for some~$p_\star\in\overline\Omega$.

Now, if~$p_\star\in\partial\Omega$ we deduce from~\eqref{VSOTTO}
that~$\sup_{\Omega} v\le g(p_\star)$, from which~\eqref{KSM-0987rpolkjbg9uIGVS56789UJHDD} follows
at once, therefore we can focus on the case in which~$p_\star\in\Omega$.

In this case,
\begin{equation}\label{swr48v54b8y7589yt54y6754754754u6ii}
v(p_\star)=\lim_{j\to+\infty} v(p_j)= \sup_{\Omega} v.
\end{equation}
Furthermore, for every~$j\in\N$, we define
$$ \Omega_j:=\big\{ x\in \Omega\; {\mbox{ s.t. }}\; B_{1/j}(x)\Subset\Omega
\big\}
$$
and we notice that there exists~$j_0\in\N$ such that~$p_\star\in\Omega_j$
for all~$j\ge j_0$.
Consequently, using the Maximum Principle in
Corollary~\ref{PERSUBAWEAK} we find that
\begin{equation}\label{d3i9v3buy478vyb4}
v(p_\star)\le \sup_{\partial \Omega_j} v.\end{equation}
Furthermore, for every~$\e>0$, we have that there exists~$p_\star^{(j)}\in\partial\Omega_j$ such that
\begin{equation}\label{d3i9v3buy478vyb422}
\sup_{\partial \Omega_j} v\le v\big(p_\star^{(j)}\big)+\e.\end{equation}
Notice that~$p_\star^{(j)}$ converges to some~$\overline{p}_\star\in\partial\Omega$ as~$j\to+\infty$.
{F}rom this, \eqref{VSOTTO}, \eqref{d3i9v3buy478vyb4} and~\eqref{d3i9v3buy478vyb422}, we deduce that
$$ v(p_\star)\le \limsup_{j\to+\infty} v\big(p_\star^{(j)}\big)+\e\le \sup_{\partial\Omega} g+\e.
$$
Hence, sending~$\e$ to zero and recalling~\eqref{swr48v54b8y7589yt54y6754754754u6ii},
we obtain~\eqref{KSM-0987rpolkjbg9uIGVS56789UJHDD}, as desired.

As a consequence of~\eqref{KSM-0987rpolkjbg9uIGVS56789UJHDD}
we also have that
\begin{equation}\label{KL:PER2}
{\mbox{the supremum in~\eqref{KL:PER} is finite, hence~$u$ is well-defined as a function in the reals.}}\end{equation}

Now we pick~$x_0\in\Omega$ 
and~$R>0$ such that~$B_{R}(x_0)\Subset\Omega$
and we show that
\begin{equation}\label{EC5436uCdOLQMNS}
\begin{split}
&{\mbox{there exist a sequence
of functions~$W_k\in C(\Omega)$}}\\&{\mbox{which are subharmonic in~$\Omega$,
harmonic in~$B_{R}(x_0)$,}}\\&{\mbox{with }}
\limsup_{\Omega\ni x\to p}W_k(x)\le g(p)
{\mbox{ for every~$p\in\partial\Omega$
and such that~$W_k\le u$ in~$\Omega$,}}\\&{\mbox{and a function~$H$ which is harmonic in~$B_{R}(x_0)$}}\\&{\mbox{such that~$W_k$ converges
to~$H$ locally uniformly in~$B_{R}(x_0)$}}\\&{\mbox{and~$u(x_0)=H(x_0)$.}}
\end{split}\end{equation}
To prove this, up to a translation, we can suppose that~$x_0=0$.
Then, by~\eqref{KL:PER}, we can take a sequence of functions~$v_k\in C(\Omega)$
which are subharmonic in~$\Omega$, 
such that~\eqref{VSOTTO} is satisfied and
\begin{equation}\label{0ihg09iuhbgv09iuhgv98uygv8yg87yg8gwfb} v_k(0)\le u(0)\le v_k(0)+\frac1k.\end{equation}
Let~$w_k(x):=\max\left\{ v_k(x),\inf_{\partial\Omega}g\right\}$.
We observe that~$w_k$ is also continuous in~$\Omega$,
that
$$\limsup_{\Omega\ni x\to p}w_k(x)\le g(p)$$
for every~$p\in\partial\Omega$, and also
that~$w_k$ is subharmonic in~$\Omega$, thanks to Lemma~\ref{POMAX}.
In particular, we have that~$u\ge w_k$. Since in addition~$w_k\ge v_k$, we deduce from~\eqref{0ihg09iuhbgv09iuhgv98uygv8yg87yg8gwfb}
that
\begin{equation}\label{0ihg09iuhbgv09iuhgv98uygv8yg87yg8gwfb-2} w_k(0)\le u(0)\le w_k(0)+\frac1k.\end{equation}
Now we take~$R>0$ sufficiently small such that~$B_R\Subset\Omega$. We take~$W_k$ to be the harmonic
lift of~$w_k$ in~$B_R$, according to Lemma~\ref{L:LIFT}.
Using~\eqref{E:LIFT}, 
we see that~$W_k$ satisfies~\eqref{VSOTTO}, and so
we have that
\begin{equation}\label{ULANS}
{\mbox{$u\ge W_k$ in~$\Omega$.}}
\end{equation}
and in particular~$u(0)\ge W_k(0)$.
This, \eqref{E:LIFT} and~\eqref{0ihg09iuhbgv09iuhgv98uygv8yg87yg8gwfb-2}
lead to
\begin{equation}\label{0oik-98uy-9uhg0987654tguvJ} W_k(0)\le u(0)\le W_k(0)+\frac1k.\end{equation}

Now we exploit Cauchy's Estimates in Theorem~\ref{CAUESTIMTH} to see that,
for every~$\alpha$, $\beta\in(0,1)$ with~$\alpha<\beta$,
\begin{equation}\label{KMSSCOILSIDD}\sup_{B_{\alpha R}} |\nabla W_k|\le C\|W_k\|_{L^\infty(B_{\beta R}(0))},\end{equation}
for some~$C>0$ depending only on~$n$, $R$, $\alpha$ and~$\beta$.

By construction, we know that
$$ \inf_{\partial\Omega}g\le
\max\left\{ v_k(x),\inf_{\partial\Omega}g\right\}=w_k(x).
$$
Moreover, recalling~\eqref{KSM-0987rpolkjbg9uIGVS56789UJHDD},
$$ w_k(x)=\max\left\{ v_k(x),\inf_{\partial\Omega}g\right\}\le \sup_{\partial\Omega}g.$$
These observations and the Maximum Principle in
Corollary~\ref{WEAKMAXPLE} yield that
\begin{equation}\label{MS-OLS-OKMS6568765}
\sup_{B_R}|W_k|=\sup_{\partial B_R}|W_k|=\sup_{\partial B_R}|w_k|\le\|g\|_{L^\infty(\partial\Omega)}.\end{equation}
Hence, utilizing~\eqref{KMSSCOILSIDD}, we deduce that~$W_k$
is uniformly equicontinuous and equibounded in~$B_{\alpha R}$. As a consequence,
by the Arzel\`a-Ascoli Theorem,
up to a subsequence we can assume 
that~$W_k$ converges locally uniformly in~$B_{\alpha R}$ to some function~$H$.
We also know that~$u(0)=
H(0)$, due to~\eqref{0oik-98uy-9uhg0987654tguvJ},
and that~$H$ is harmonic in~$B_{\alpha R}$, thanks to
Corollary~\ref{S-OS-OSSKKHEPARHA}.
The proof of~\eqref{EC5436uCdOLQMNS} is thereby completed by collecting all the previous information.

Now, given~$x_0\in\Omega$, we claim that
\begin{equation}\label{FIPEERMOCN}
{\mbox{$u$ is harmonic
in a neighborhood of~$x_0$.}}\end{equation}
Notice that this would
complete the proof of
Theorem~\ref{PERRO}, since being harmonic is a local property, thanks to
Theorem~\ref{KAHAR}(iii).

Hence, to complete the desired result in~\eqref{FIPEERMOCN}, it suffices to take~$R$ and~$H$ as in~\eqref{EC5436uCdOLQMNS} and
show that
\begin{equation}\label{KS-0okMNSNS-okf9253659843679043ythgbge}
{\mbox{$u=H$ \; in \; $B_{R}(x_0)$.}}\end{equation} By~\eqref{EC5436uCdOLQMNS}, we already know that
$$u\ge\lim_{k\to+\infty}W_k(x)=H(x)$$
for every~$x\in B_{R}(x_0)$, hence, to prove~\eqref{KS-0okMNSNS-okf9253659843679043ythgbge},
we can argue by contradiction and suppose that there exists a point~$x_0^\sharp\in B_{R}(x_0)$ such that~$u(x_0^\sharp)\ge H(x_0^\sharp)+2a$, for some~$a>0$.

In particular, by~\eqref{KL:PER}, we can find
a function~$v^\sharp\in C(\Omega)$
which is subharmonic in~$\Omega$, such that
$$ \limsup_{\Omega\ni x\to p}v^\sharp(x)\le g(p)\qquad{\mbox{for every~$p\in\partial\Omega$}}
$$ and satisfying
\begin{equation}\label{099-013495958j3j} v^\sharp(x_0^\sharp)\ge H(x_0^\sharp)+a.\end{equation}
Let~$w^\sharp_k(x):=\max\{ v^\sharp(x),W_k(x)\}$, with~$W_k$ as in~\eqref{EC5436uCdOLQMNS}.
We observe that~$w^\sharp_k
\in C(\Omega)$, that
$$ \limsup_{\Omega\ni x\to p}w_k^\sharp(x)\le g(p)\qquad{\mbox{for every~$p\in\partial\Omega$}}
$$ and also
that~$w_k^\sharp$ is subharmonic in~$\Omega$, thanks to Lemma~\ref{POMAX}.

Now we take~$W_k^\sharp$ as the harmonic lift of~$w_k^\sharp$ in the ball~$B_R(x_0)$,
according to Lemma~\ref{L:LIFT}.
As in~\eqref{KMSSCOILSIDD} and~\eqref{MS-OLS-OKMS6568765}, we see that~$W_k^\sharp$
is locally
uniformly equicontinuous and equibounded in~$B_R(x_0)$ and therefore,
up to a subsequence, it converges locally uniformly in~$B_R(x_0)$ to some function~$H^\sharp$
which is harmonic in~$B_R(x_0)$.

We stress that, for every~$x\in B_R(x_0)$,
$$ W_k(x)\le w^\sharp_k(x)\le W^\sharp_k(x)
$$
and therefore
\begin{equation}\label{MSLD-opKSMv5609986pm8sidpm} H(x)\le H^\sharp(x)\qquad{\mbox{for every \, $x\in B_R(x_0)$.}}
\end{equation}
Moreover, by~\eqref{099-013495958j3j},
$$ H(x_0^\sharp)+a\le v^\sharp(x_0^\sharp)
\le w^\sharp_k(x_0^\sharp)\le W^\sharp_k(x_0^\sharp),$$
from which we arrive at
$$ H(x_0^\sharp)<H(x_0^\sharp)+a\le H^\sharp(x_0^\sharp).
$$
Combining this information with~\eqref{MSLD-opKSMv5609986pm8sidpm} and the Strong Maximum Principle
in Theorem~\ref{STRONGMAXPLE1}(iii), we deduce that
\begin{equation}\label{EC5436uCdOLQMNSx0}
{\mbox{$H(x)<H^\sharp(x)$ for every~$x\in B_R(x_0)$.}}\end{equation}
Furthermore, by the properties of the harmonic lift in Lemma~\ref{L:LIFT} and~\eqref{KL:PER}, we know that~$u\ge W_k^\sharp$ in~$\Omega$.
Consequently, we have that~$u\ge H^\sharp$ in~$B_R(x_0)$.
{F}rom this, \eqref{EC5436uCdOLQMNS} and~\eqref{EC5436uCdOLQMNSx0} we deduce that
$$
u(x_0)=H(x_0)<H^\sharp(x_0)\le u(x_0),
$$
which provides the desired contradiction.
\end{proof}

We want to emphasize that Theorem~\ref{PERRO} in itself does not
ensure that the harmonic function constructed there attains continuously the boundary datum~$g$ along~$\partial\Omega$.
This however is guaranteed if~$\Omega$ is sufficiently regular. We introduce some terminology
to clarify this point.

\begin{figure}
  \centering
  \includegraphics[width=.35\linewidth]{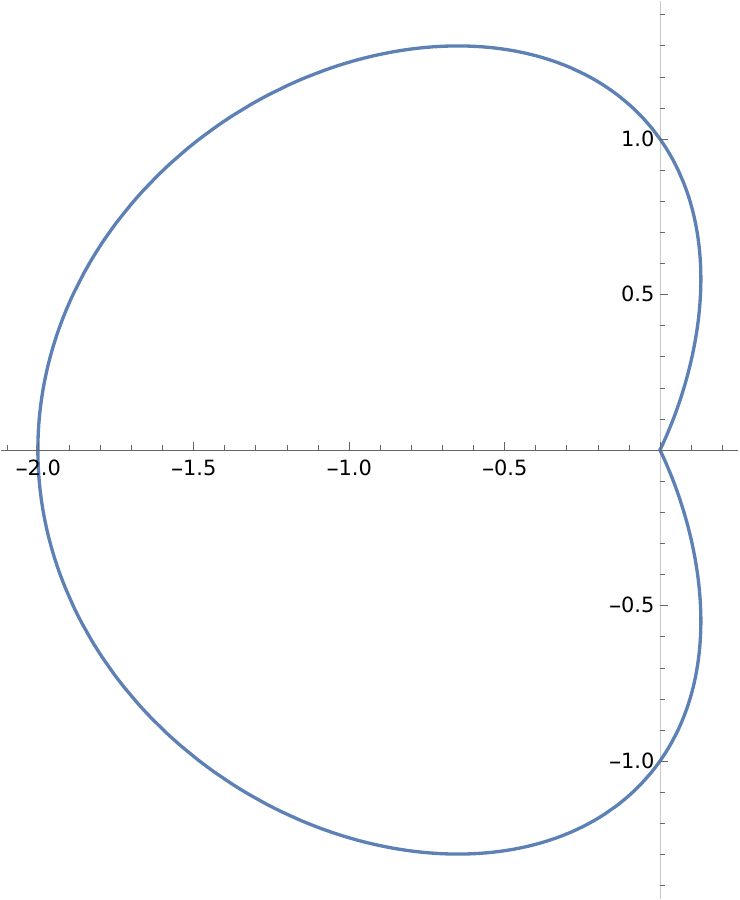}$\quad$
  \includegraphics[width=.5\linewidth]{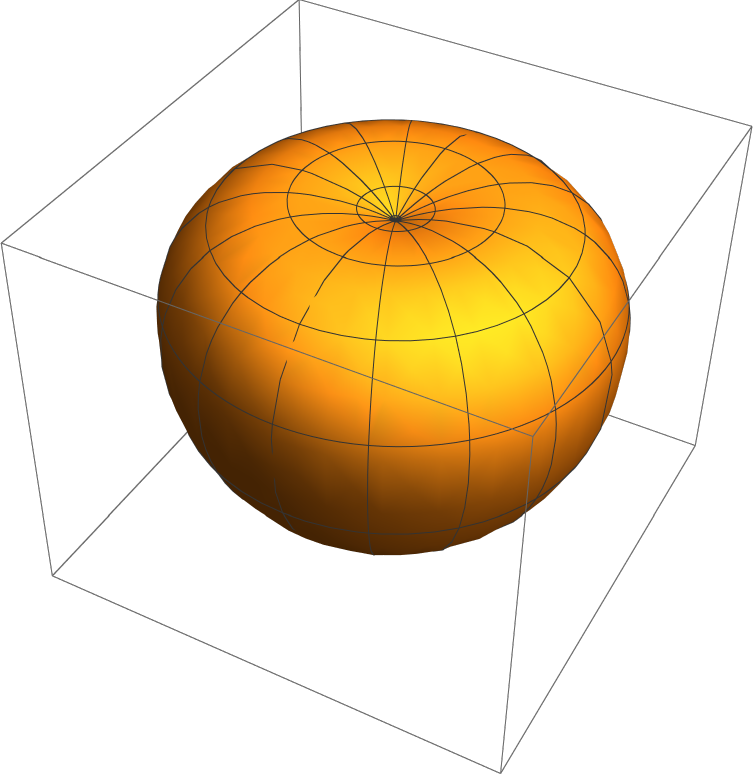}$\quad$
 \caption{\sl Domains in~$\R^2$ and~$\R^3$ satisfying the exterior cone condition.
 These are parametric plots of~$\big(x(t),y(t)\big)=\big((1-\cos t)\sin |t/2|\cos t,(1-\cos t)\sin t\big)\,$
 and~$\big(x(s,t),y(s,t),z(s,t)\big)=\big(
 (1-\cos t) \sin t\sin s,\,
 (1-\cos t)\sin t\cos s,\,
 (1-\cos t)\sin |t/2|\cos t \big)$, respectively, with~$s\in(-\pi/2,\pi/2]$, $t\in(-\pi,\pi]$.
 }\label{T536EcdhPOaIJodidpo3tju24ylkiliR2543654yh2432RAe87658FI6MI4M356U2LSMSUJNS-G5O2F4O2R-08jC1ON}
\end{figure}

We say that an open set~$\Omega\subset\R^n$
satisfies the exterior cone condition\index{exterior cone condition} if
for every point~$x\in\partial\Omega$
there exists a finite right circular\footnote{As usual, given~$b\in(0,\pi)$ and~$e\in\partial B_1$, a right circular cone of opening~$b$ in direction~$e$ is an object of the form
$$ \left\{x\in\R^n{\mbox{ s.t. }}\frac{x\cdot e}{|x|}>\cos b\right\}.$$
A finite right circular cone then takes the form
$$ \left\{x\in B_\rho{\mbox{ s.t. }}\frac{x\cdot e}{|x|}>\cos b\right\},$$
for some~$\rho>0$.} cone~$K$ with vertex at~$x$ such that
\begin{equation}\label{CONOK}\overline{K}\cap\overline\Omega =\{x\}.\end{equation}
See Figure~\ref{T536EcdhPOaIJodidpo3tju24ylkiliR2543654yh2432RAe87658FI6MI4M356U2LSMSUJNS-G5O2F4O2R-08jC1ON}
for examples in~$\R^2$ and~$\R^3$ of domains satisfying the exterior cone condition.

We also say that an open set~$\Omega\subset\R^n$
is accessible by simple arcs from the exterior\index{accessible by simple arcs from the exterior} if for every~$x\in\partial\Omega$
there exists an injective continuous map~$\gamma:[0,1]\to\R^n$ such that~$\gamma(t)\in\R^n\setminus\overline\Omega$
for all~$t\in[0,1)$ and~$\gamma(1)=x$.

As a side remark, we point out that if~$\Omega$ satisfies the exterior ball condition
then it is also accessible by simple arcs from the exterior.

With this notation, we have:

\begin{theorem}\label{PERRO-2}
Let~$\Omega\subset\R^n$ be open,
bounded and connected. Assume that either~$n=2$ and~$\Omega$
is accessible by simple arcs from the exterior, or~$n\ge3$ and~$\Omega$
satisfies the exterior cone\footnote{For completeness, we observe that
the proof of Theorem~\ref{PERRO-2} becomes technically easier \label{DEPHIBARRPERPRER8-F}
for domains with additional regularity, see e.g.~\cite[equation~(2.34)]{MR1814364}
for a barrier in the case of domains
satisfying the exterior ball condition. With this respect,
one says that an open set~$\Omega\subset\R^n$
satisfies the exterior ball condition\index{exterior ball condition} if for every~$x\in\partial\Omega$ there exists~$R>0$
and~$y\in\R^n$ such that~$\overline{B_R(y)}\cap\overline{\Omega}=\{x\}$.
We observe that domains with boundary of class~$C^{1,1}$ satisfy
the exterior ball condition,
see e.g.~\cite[Lemma~A.1]{MR3436398}.
Notice also that this condition is somewhat the counterpart of the interior
ball condition that was assumed in the
Hopf Lemma (namely, Lemma~\ref{JJS:PA}). We also stress that domains satisfying the exterior
ball condition obviously satisfy the exterior
cone condition as well.

For additional results on the Perron method for domains satisfying the
exterior cone condition see~\cite{MR171990, MR221087, MR316884, MR425346, MR603385, MR1696765}.}
condition.

Let~$g$ be a
continuous function on~$\partial\Omega$ and let~$u$ be as in~\eqref{KL:PER}.
Then, $u\in C(\overline\Omega)$ and
\begin{equation}\label{IMEND-PMA2456hgfdhAA}\lim_{\Omega\ni x\to p}u(x)=g(p)\end{equation}
for all~$p\in\partial\Omega$.
\end{theorem}

\begin{proof} We start with a general observation
related to superharmonic functions. For this, suppose
that~$\Omega\subset\R^n$ is open
and bounded, and let~$p\in\partial\Omega$. Let~$\rho>0$ and~$\beta:\overline{\Omega\cap B_\rho(p)}\to\R$
such that~$\beta\in C({\Omega\cap B_\rho(p)})$.
Assume that~$\beta$ is superharmonic in~$\Omega\cap B_\rho(p)$
and let
\begin{equation}\label{SETOMEG}
m:=\inf_{\Omega\cap(B_\rho(p)\setminus B_{\rho/2}(p))}\beta\qquad{\mbox{and}}\qquad
\omega(x):=\begin{dcases}
\min\{m,\beta(x)\}& {\mbox{ if }}x\in\overline\Omega\cap B_{\rho/2}(p),\\
m& {\mbox{ if }}x\in\overline\Omega\setminus B_{\rho/2}(p).
\end{dcases}\end{equation}
We point out that
\begin{equation}\label{2MA:0sod03044d4}
\omega\in C(\Omega).\end{equation}
We claim that
\begin{equation}\label{KSM-0987y6t5rghn-oytgghNS-SUP-2}
{\mbox{$\omega$ is superharmonic in~$\Omega$.}}\end{equation}
For this, 
we first employ Lemma~\ref{POMAX}
to see that
\begin{equation}\label{KSM-0987y6t5rghn-oytgghNS-SUP-1}
{\mbox{$\omega$ is superharmonic in~$\Omega\cap B_{\rho/2}(p)$.}}\end{equation}
Furthermore,
we observe that~$\omega$ is constant, hence
superharmonic in~$\Omega\setminus\overline{B_{\rho/2}(p)}$.
This, \eqref{KSM-0987y6t5rghn-oytgghNS-SUP-1} and the fact
that superharmonicity is a local property (recall
Theorem~\ref{SUBHS-EQ}(viii)) gives that, to prove~\eqref{KSM-0987y6t5rghn-oytgghNS-SUP-2}
it suffices to consider an arbitrary point~$q\in
\Omega\cap\partial B_{\rho/2}(p)$ and an arbitrary small radius~$r_0$
(in particular, $B_{r_0}(q)\Subset\Omega$)
and check that, for each~$r\in(0,r_0]$,
\begin{equation}\label{KSM-0987y6t5rghn-oytgghNS-SUP-3}
\omega(q)\ge\fint_{B_{r}(q)} \omega(x)\,dx.
\end{equation}
To establish this, we observe that~$\omega\le m$ and thus
\begin{eqnarray*}
\fint_{B_{r}(q)} \omega(x)\,dx
\le m=\omega(q),
\end{eqnarray*}
which gives~\eqref{KSM-0987y6t5rghn-oytgghNS-SUP-3},
and thus~\eqref{KSM-0987y6t5rghn-oytgghNS-SUP-2}, as desired.

Now we point out some general remarks about the existence of a barrier
(for this, we will specify the functions~$\beta$ and~$\omega$ introduced above). Namely, suppose again that~$\Omega\subset\R^n$ is open
and bounded, and let~$g$ be a bounded and continuous function on~$\partial\Omega$
and~$u$ be as in~\eqref{KL:PER}.
Let~$p\in\partial\Omega$. We have that
\begin{equation}\label{JKA:BARRLPAER}
\begin{split}&
{\mbox{if there exist~$\rho>0$ and a function~$\beta: \overline{\Omega\cap B_\rho(p)}\to\R$,}}\\
&{\mbox{which is continuous in~$\Omega\cap B_\rho(p)$,}}\\&{\mbox{superharmonic in~$\Omega\cap B_\rho(p)$, strictly positive in~$\overline{\Omega\cap B_\rho(p)}\setminus\{p\}$}}\\&{\mbox{and such that }}
\lim_{\Omega\ni x\to p}\beta(x)=0,\\&
{\mbox{then }}\lim_{\Omega\ni x\to p}u(x)=g(p)
.\end{split}
\end{equation}
To prove this, we adopt the setting in~\eqref{SETOMEG}.
We have that
$$m>0= \lim_{\Omega\ni x\to p}\beta(x)$$ and therefore
\begin{equation}\label{MA:0sod03044d4}
{\mbox{if~$q\in\overline{\Omega}\setminus\{p\}$
then }}
\liminf_{\Omega\ni x\to q}\omega(x)> \lim_{\Omega\ni x\to p} \beta(x)=\lim_{\Omega\ni x\to p} \omega(x)=0.\end{equation}
Now, given~$\e>0$, we exploit the continuity of~$g$ to pick a~$\delta>0$
such that
\begin{equation}\label{INESP-EPS-2}
{\mbox{if $x\in (\partial\Omega)\cap B_\delta(p)$ we have that~$|g(x)-g(p)|\le\e$.}}\end{equation}
Also, since~$\iota:=\inf_{\Omega\setminus B_\delta(p)}\omega>0$,
due to~\eqref{MA:0sod03044d4}, by taking~$K:=\frac{2\|g\|_{L^\infty(\partial\Omega)}}\iota$
we have that,
for all~$x\in\Omega\setminus B_\delta(p)$
\begin{equation}\label{INESP-EPS-23} K \omega(x)\ge K \iota= 2\|g\|_{L^\infty(\partial\Omega)}.\end{equation}
We now define
$$ \omega_\pm(x):=g(p)\pm\e\pm K\omega(x).$$
We claim that
\begin{equation}\label{SOPRASOTTOMEGA}
\limsup_{\Omega\ni x\to q}\omega_-(x)\le g\le\liminf_{\Omega\ni x\to q}\omega_+(x)\quad{\mbox{ for every }}\,q\in\partial\Omega.
\end{equation}
For this, we take~$q\in\partial\Omega$ and we distinguish two cases.
If~$q\in B_\delta(p)$, we have that
$$ \limsup_{\Omega\ni x\to q}\omega_-(x)-g(q)\le\limsup_{\Omega\ni x\to q}\big(g(p)-\e-K\omega(x)\big)
-g(p)+\e=-K\liminf_{\Omega\ni x\to q}\omega(x)\le0
,$$
thanks to~\eqref{MA:0sod03044d4} and~\eqref{INESP-EPS-2}, and similarly
$$ \liminf_{\Omega\ni x\to q}\omega_+(x)-g(q)\ge
\liminf_{\Omega\ni x\to q}\big(g(p)+\e+ K\omega(x)\big)-g(p)-\e=K\liminf_{\Omega\ni x\to q}
\omega(x)\ge0.
$$
These considerations give that~\eqref{SOPRASOTTOMEGA} holds true
in~$(\partial\Omega)\cap B_\delta(p)$, hence we consider now the case in which~$q\in
(\partial\Omega)\setminus B_\delta(p)$. In this situation,
we utilize~\eqref{INESP-EPS-23} to see that $$\limsup_{\Omega\ni x\to q}
\omega_-(x)-g(q)=\limsup_{\Omega\ni x\to q}\big(
g(p)-\e- K\omega(x)\big)-g(q)\le
g(p)-\e- 2\|g\|_{L^\infty(\partial\Omega)}-g(q)\le-\e<0
$$
and in the same way
$$ \liminf_{\Omega\ni x\to q}\omega_+(x)-g(q)=\limsup_{\Omega\ni x\to q}\big(
g(p)+\e+ K\omega(x)\big)-g(q)\ge
g(p)+\e+ 2\|g\|_{L^\infty(\partial\Omega)}-g(q)\ge\e>0.
$$
These observations complete the proof of~\eqref{SOPRASOTTOMEGA}.

Furthermore, $\omega_+$ is superharmonic in~$\Omega$
and~$\omega_-$ is subharmonic in~$\Omega$, thanks to~\eqref{KSM-0987y6t5rghn-oytgghNS-SUP-2}. Also, these functions
are continuous in~$\overline\Omega$, due to~\eqref{2MA:0sod03044d4}.
By~\eqref{KL:PER} and~\eqref{SOPRASOTTOMEGA}, we infer that~$u\ge\omega_-$,
and as a result
\begin{equation*}
\lim_{\Omega\ni x\to p}u(x)\ge
\lim_{\Omega\ni x\to p}\omega_-(x) =
g(p)-\e.
\end{equation*}
By taking~$\e$ arbitrarily small we arrive at
\begin{equation}\label{JSN-COOMNSalBOSRFSD0}
\lim_{\Omega\ni x\to p}u(x)\ge g(p).
\end{equation}
Additionally,
for every function~$v\in C(\Omega)$
which is subharmonic in~$\Omega$ and satisfies~\eqref{VSOTTO}, we have
that~$v-\omega_+$ is subharmonic in~$\Omega$.

Also, we claim that
\begin{equation}\label{132367548484dewv9t8bv8uy8}
\sup_{\Omega}(v-\omega_+)\le0.
\end{equation}
To prove it, we argue by contradiction and suppose that there exists a sequence of points~$p_k\in\Omega$
such that
$$\lim_{k\to+\infty}(v-\omega_+)(p_k)=\sup_{\Omega}(v-\omega_+)\ge a>0.$$
Since~$\Omega$ is bounded, we have that~$p_k$ converges, up to a subsequence, to some~$\overline{p}\in\overline\Omega$.
We observe that~$\overline{p}\in\Omega$, otherwise, if~$\overline{p}\in\partial\Omega$, in light of~\eqref{VSOTTO}
and~\eqref{SOPRASOTTOMEGA} we would have
$$ 0<a\le \lim_{k\to+\infty}(v-\omega_+)(p_k)\le \limsup_{\Omega\ni x\to\overline{p}}(v-\omega_+)(x)\le g(\overline{p})-
g(\overline{p})=0,
$$
which is a contradiction. Hence, exploiting Lemma~\ref{LplrfeELLtahsmedschdfb4camntM}
we conclude that~$v-\omega_+$ is constant in~$\Omega$, and in particular~$(v-\omega_+)(x)\ge a$ for every~$x\in\Omega$.

Accordingly, for every~$q\in\partial\Omega$,
$$ a\le \limsup_{\Omega\ni x\to q}(v-\omega_+)(x)\le g(q)-g(q)=0,
$$
which is a contradiction, and therefore the proof of~\eqref{132367548484dewv9t8bv8uy8} is complete.

\begin{figure}
  \centering
  \includegraphics[width=.5\linewidth]{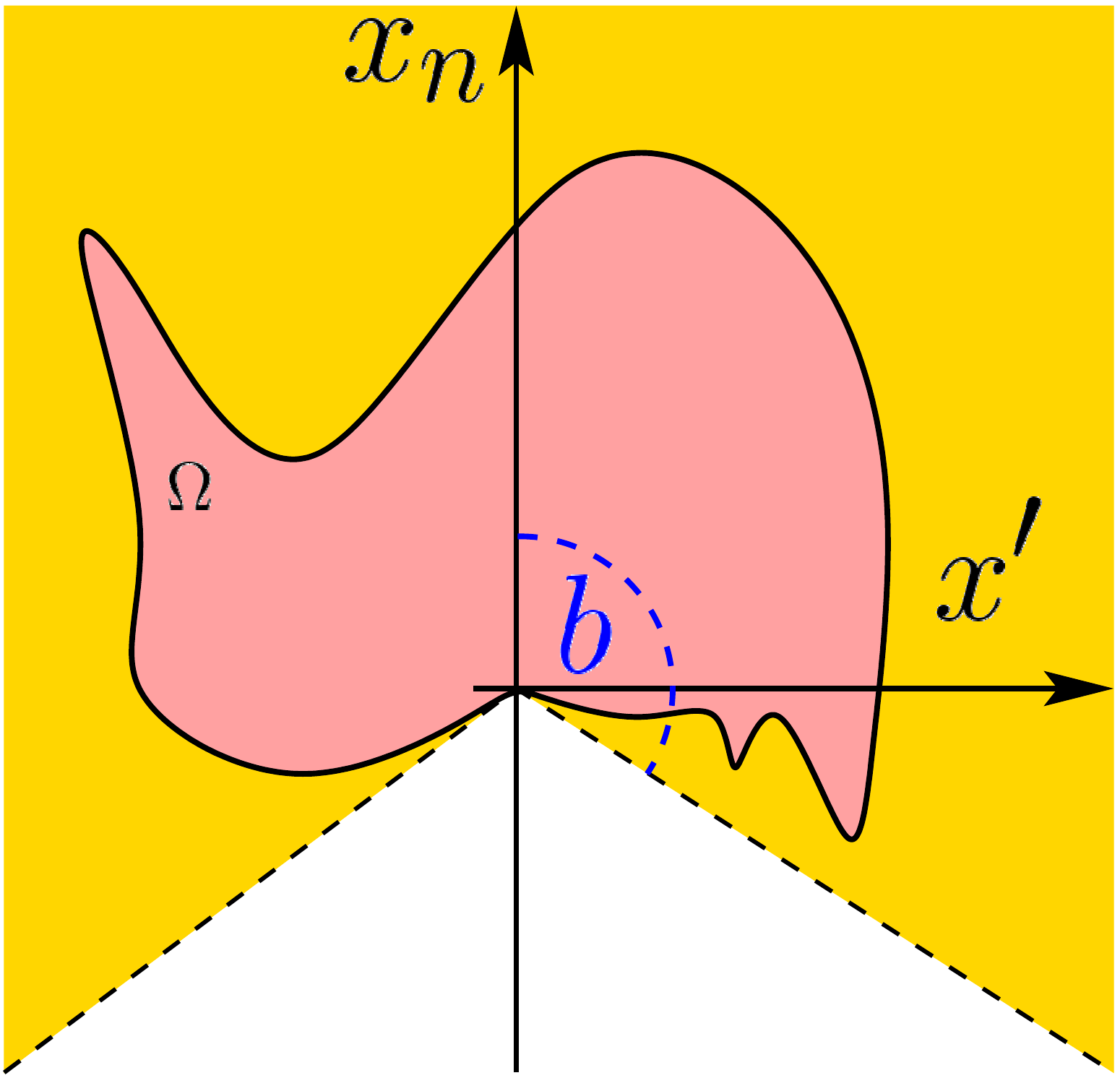}
 \caption{\sl A domain satisfying the exterior cone condition
 at the origin.}\label{T536FGHSN0dEcdhPOaIJodidpo3tju24ylkiliR2543654yh2432RAe87658FI6MI4M356U2LSMSUJNS-G5O2F4O2R}
\end{figure}

Therefore, thanks to~\eqref{132367548484dewv9t8bv8uy8}, we have that~$v\le\omega_+$ in~$\Omega$.
Taking the supremum over such functions~$v$ and recalling the definition of~$u$ in~\eqref{KL:PER},
we thereby find that~$u\le\omega_+$ in~$\Omega$.

Accordingly,
\begin{equation*}
\lim_{\Omega\ni x\to p}u(x)\le \lim_{\Omega\ni x\to p}\omega_+(x)
= g(p)+\e.
\end{equation*}
By taking~$\e$ as small as we like we conclude that
\begin{equation*}
\lim_{\Omega\ni x\to p}u(x)\le g(p).
\end{equation*}
Combining this with~\eqref{JSN-COOMNSalBOSRFSD0}, we obtain~\eqref{JKA:BARRLPAER},
as desired.

Now, in light of~\eqref{JKA:BARRLPAER}, in order to complete the proof
of Theorem~\ref{PERRO-2}, it suffices to check that the domains
in the hypothesis of Theorem~\ref{PERRO-2} allow the construction of a barrier~$\beta$
as requested in~\eqref{JKA:BARRLPAER}. To this end, we treat separately the cases~$n=2$ and~$n\ge3$.

When~$n=2$, let~$p\in\Omega$ and suppose, up to a translation, that~$p$ is the origin.
Since~$\Omega$ is accessible by simple arcs from the exterior,
there exists an injective continuous map~$\gamma:[0,1]\to\R^n$ such that~$\gamma(t)\in\R^n\setminus\overline\Omega$
for all~$t\in[0,1)$ and~$\gamma(1)=0$. In particular, we have that~$\gamma(0)\ne\gamma(1)=0$, hence we
can take~$\rho\in(0,\min\{1,|\gamma(0)|\})$.

We use the polar representation
in complex coordinates~$z=r e^{i\vartheta}$ for points~$z\in\Omega\cap B_\rho$,
with~$r>0$ and~$\vartheta\in\R$.
We notice that the curve~$\gamma$ provides a branch cut
for the definition of the polar angle~$\vartheta$ and thus for the one of the
complex logarithm in~$\Omega\cap B_\rho$.
Notice also that if~$z=r e^{i\vartheta}\in\Omega\cap B_\rho$ then~$|z|<\rho<1$, whence~$\ln z\ne0$. As a result, we can define
the holomorphic function~$-\frac1{\ln z}$. This gives rise (see e.g.~\cite[Chapter 11]{MR0210528}) to the harmonic
(hence superharmonic)
function
$$ \Omega\cap B_\rho\ni z=re^{i\vartheta}\mapsto\beta(z):=\Re\left(-\frac1{\ln z}\right)=-\frac{\ln r}{\ln^2 r+\vartheta^2}.$$
Notice that
$$ \lim_{z\to0}\beta(z)=0<
\frac{-\ln \widetilde{r}}{\ln^2  \widetilde{r}+ \widetilde\vartheta^2}
=\lim_{z\to\widetilde{z}}\beta(z),$$
for every~$\Omega\cap B_\rho\ni \widetilde{z}=\widetilde{r}e^{i\widetilde\vartheta}\ne0$.
These considerations show that~$\beta$ is an appropriate barrier in the sense of~\eqref{JKA:BARRLPAER},
whence~\eqref{IMEND-PMA2456hgfdhAA} follows in this case.

\begin{figure}
  \centering
  \includegraphics[width=.4\linewidth]{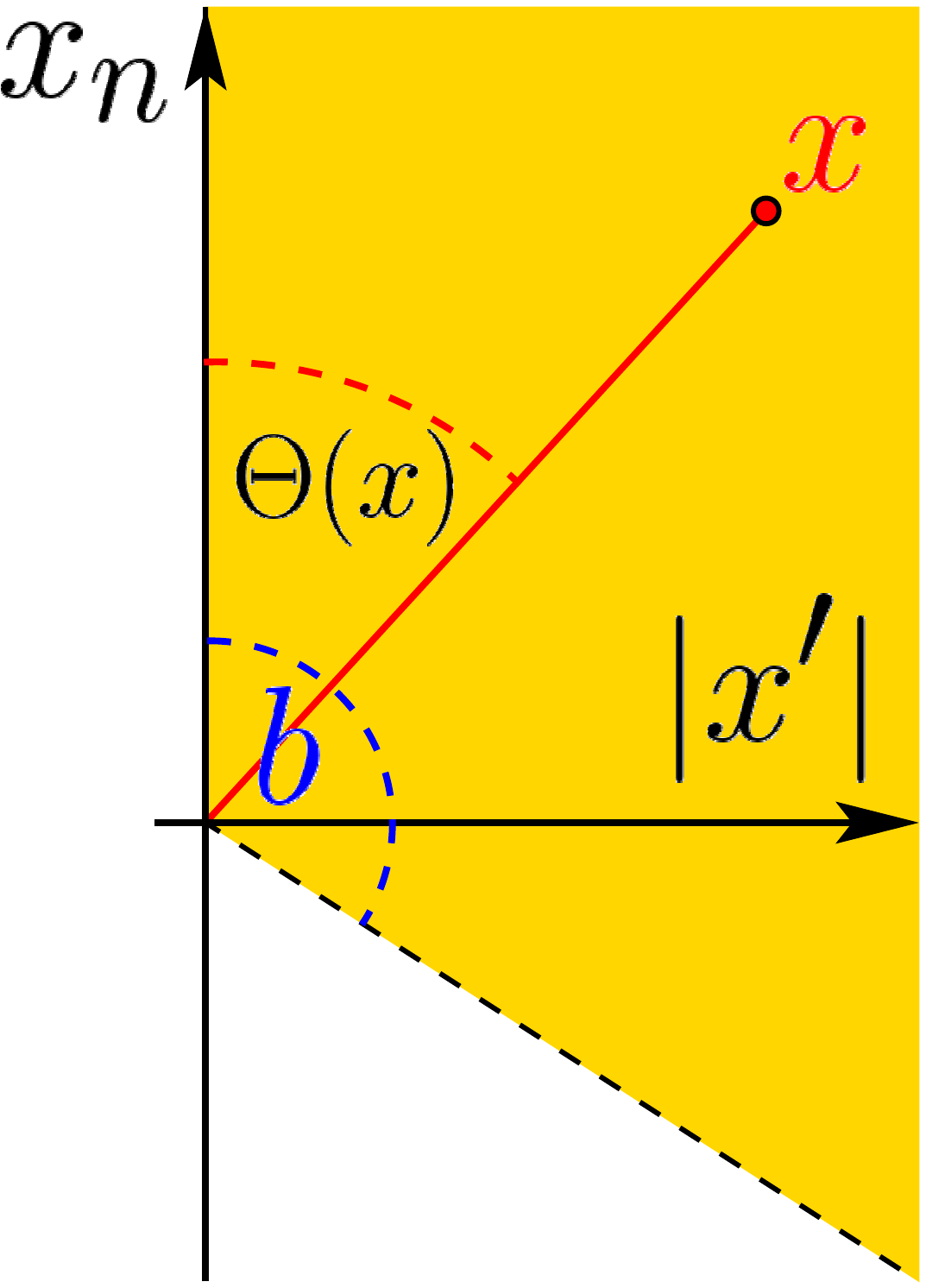}
 \caption{\sl The cone in~\eqref{MI4M356U2LSMSUJNS-G5O2F4O2R}.}\label{T536EcdhPOaIJodidpo3tju24ylkiliR2543654yh2432RAe87658FI6MI4M356U2LSMSUJNS-G5O2F4O2R}
\end{figure}

Let us now deal with the case~$n\ge3$. For this, 
given~$b\in(0,\pi)$, we consider\footnote{For the sake of completeness,
we point out that another approach
to the proof of Theorem~\ref{PERRO-2} when~$n\ge3$
consists in constructing barriers of the form~$|x|^\lambda f\left(\frac{x}{|x|}\right)$
for a suitable~$\lambda>0$, being~$f$ the first eigenfunction for the Laplace-Beltrami operator on the sphere
with vanishing datum outside a given cone.
See also~\cite[Theorem~8.27 on pages~173--174]{MR0261018}
for a different barrier in the exterior cone condition setting.

We also mention that dealing with the exterior cone condition in dimension~$n=2$
is technically easier, since one can rely on complex analysis and seek barriers of the form~$\Im z^\lambda=
r^\lambda\sin(\lambda\vartheta)$ with~$\lambda>0$ large enough.}
the cone
\begin{equation} \label{MI4M356U2LSMSUJNS-G5O2F4O2R}
{\mathcal{C}}_b:=\left\{
x=(x_1,\dots,x_n)\in\R^n {\mbox{ s.t. }}\frac{x_n}{|x|}>\cos b
\right\},\end{equation}
see Figures~\ref{T536FGHSN0dEcdhPOaIJodidpo3tju24ylkiliR2543654yh2432RAe87658FI6MI4M356U2LSMSUJNS-G5O2F4O2R}
and~\ref{T536EcdhPOaIJodidpo3tju24ylkiliR2543654yh2432RAe87658FI6MI4M356U2LSMSUJNS-G5O2F4O2R},
and we deal with a preliminary
computation about the Laplacian of a given function~$w$ of the form
$${\mathcal{C}}_b\ni x\mapsto
w(x)=w_0(|x|,\Theta(x)),$$
where~$\Theta(x)\in[0,b)$ is such that
\begin{equation}\label{0OKMSMS-2245ne}
x_n=|x|\cos(\Theta({x}))\qquad {\mbox{ and }}\qquad |x'|=|x|\sin(\Theta({x})).
\end{equation}
 
We pick a point~$\overline{x}\ne0$.
Up to a rotation (recall
the rotational invariance of the Laplacian according to Corollary~\ref{KSMD:ROTAGSZOKA})
we suppose that~$\overline{x}$ lies in the plane spanned by the first and the last elements
of the Euclidean basis, namely~$\overline{x}=(\overline{x}_1,0,\dots,0,\overline{x}_n)$.
We define~$\overline{x}^\perp:=(\overline{x}_n,0,\dots,0,-\overline{x}_1)$
and take an orthonormal frame~$\{\eta_1,\dots,\eta_n\}$ with
$$\eta_1:=\frac{\overline{x}}{|\overline{x}|}\qquad{\mbox{and}}\qquad\eta_2:=\frac{\overline{x}^\perp}{|\overline{x}|}.$$

We remark that if~$j\in\{3,\dots,n\}$ then~$\overline{x}\cdot\eta_j=|\overline{x}|\eta_1\cdot\eta_j=0$, whence, for~$\e$ small,
$$ |\overline{x}+\e\eta_j|=
\sqrt{|\overline{x}|^2+2\e\overline{x}\cdot\eta_j+\e^2}=
|\overline{x}|
\sqrt{1+\frac{\e^2}{|\overline{x}|^2}}=
|\overline{x}|+\frac{\e^2}{2|\overline{x}|}+o(\e^2).$$
Also,
$$ \overline{x}_n\overline{x}-\overline{x}_1\overline{x}^\perp=
(\overline{x}_1\overline{x}_n,0,\dots,0,\overline{x}_n^2)
-(\overline{x}_1\overline{x}_n,0,\dots,0,-\overline{x}_1^2)=
(0,0,\dots,0,|\overline{x}|^2)=|\overline{x}|^2 e_n,$$
and accordingly
$$ e_n=\frac{
\overline{x}_n\eta_1-\overline{x}_1\eta_2
}{|\overline{x}|}.$$
This gives that~$\eta_j\cdot e_n=0$ for all~$j\in\{3,\dots,n\}$, leading to
$$ (\overline{x}+\e\eta_j)\cdot e_n=
\overline{x}_n.$$
As a result,
\begin{eqnarray*}&& \Theta(\overline x+\e\eta_j)=\arccos\left(\frac{\overline x_n}{|\overline{x}|+\frac{\e^2}{2|\overline{x}|}+o(\e^2)}\right)=\arccos\left(
\frac{\overline x_n}{|\overline{x}|}-\frac{\e^2 \overline x_n}{2|\overline{x}|^3}+o(\e^2)\right)\\&&\quad\quad=
\arccos\left(
\frac{\overline x_n}{|\overline{x}|}\right)+\frac{\e^2 \overline x_n}{2|\overline{x}|^2\sqrt{|\overline{x}|^2-\overline x_n^2}}
+o(\e^2)=\Theta(\overline{x})
+\frac{\e^2 \overline x_n}{2|\overline{x}|^2|\overline x_1|}
+o(\e^2).
\end{eqnarray*}
{F}rom these observations we deduce that, for each~$j\in\{3,\dots,n\}$,
\begin{eqnarray*}
\partial^2_{\eta_j\eta_j} w(\overline{x})&=&\lim_{\e\to0}\frac{
w_0(|\overline x+\e\eta_j|,\Theta(\overline x+\e\eta_j))
+w_0(|\overline x-\e\eta_j|,\Theta(\overline x-\e\eta_j))
-2w_0(|\overline x|,\Theta(\overline x))
}{\e^2}
\\&=&2\lim_{\e\to0}\frac{
w_0\left(
|\overline{x}|+\frac{\e^2}{2|\overline{x}|}+o(\e^2)
,\Theta(\overline{x})
+\frac{\e^2 \overline x_n}{2|\overline{x}|^2|\overline x_1|}
+o(\e^2)\right)
-w_0(|\overline x|,\Theta(\overline x))
}{\e^2}
\\&=&
\frac{\partial_1 w_0 (|\overline{x}|,\Theta(\overline{x}) )}{|\overline{x}|}
+
\frac{\overline x_n\,\partial_2 w_0 (|\overline{x}|,\Theta(\overline{x}) )
}{|\overline{x}|^2|\overline x_1|}.
\end{eqnarray*}

Moreover,
$$ |\overline{x}+\e\eta_1|=
\left|
\frac{|\overline{x}|\overline{x}+\e\overline{x}}{|\overline{x}|}\right|
=|\overline{x}|+\e.$$
Also,~$(\overline{x}+\e\eta_1)\cdot e_n=\overline{x}_n+
\frac{\e\overline{x}_n}{|\overline{x}|}$, whence
\begin{eqnarray*}&&
\Theta(\overline x+\e\eta_1)=
\arccos\left(\frac{(\overline{x}+\e\eta_1)\cdot e_n}{|\overline{x}+\e\eta_1|}\right)
\\&&\qquad=
\arccos\left(\frac{\overline{x}_n+
\frac{\e\overline{x}_n}{|\overline{x}|}}{|\overline{x}|+\e}\right)=
\arccos\left( \frac{\overline{x}_n}{|\overline{x}|}\right)
=\Theta(\overline x).
\end{eqnarray*}
Consequently,
\begin{eqnarray*}
\partial^2_{\eta_1\eta_1} w(\overline{x})&=&\lim_{\e\to0}\frac{
w_0(|\overline x+\e\eta_1|,\Theta(\overline x+\e\eta_1))
+w_0(|\overline x-\e\eta_1|,\Theta(\overline x-\e\eta_1))
-2w_0(|\overline x|,\Theta(\overline x))
}{\e^2}
\\&=&
\lim_{\e\to0}\frac{
w_0(|\overline x|+\e,\Theta(\overline x))
+w_0(|\overline x|-\e,\Theta(\overline x))
-2w_0(|\overline x|,\Theta(\overline x))
}{\e^2}
\\&=&\partial^2_{11}w_0(|\overline x|,\Theta(\overline x)).
\end{eqnarray*}

We also note that
$$ |\overline{x}+\e\eta_2|=
\left|\overline{x}+
\frac{\e\overline{x}^\perp}{|\overline{x}|}\right|=
\sqrt{|\overline{x}|^2+
\frac{2\e\overline x\cdot\overline{x}^\perp}{|\overline{x}|}
+\frac{\e^2|\overline{x}^\perp|^2}{|\overline{x}|^2}
}=\sqrt{|\overline{x}|^2+\e^2}=|\overline{x}|+\frac{\e^2}{2|\overline{x}|}+o(\e^2).
$$
Furthermore,~$(\overline{x}+\e\eta_2)\cdot e_n=\overline{x}_n-
\frac{\e\overline{x}_1}{|\overline{x}|}$, whence
\begin{eqnarray*}&&
\Theta(\overline x+\e\eta_2)=
\arccos\left(\frac{\overline{x}_n-
\frac{\e\overline{x}_1}{|\overline{x}|}}{
|\overline{x}|+\frac{\e^2}{2|\overline{x}|}+o(\e^2)
}\right)=
\arccos\left(
\left(\frac{\overline{x}_n}{|\overline x|}-
\frac{\e\overline{x}_1}{|\overline{x}|^2}\right)
\left(1-\frac{\e^2}{2|\overline{x}|^2}+o(\e^2)
\right)
\right)\\&&\qquad=\arccos\left(
\frac{\overline{x}_n}{|\overline x|}-
\frac{\e\overline{x}_1}{|\overline{x}|^2}
-\frac{\e^2\overline x_n}{2|\overline{x}|^3}+o(\e^2)
\right)=\Theta(\overline{x})
+\frac{\e\overline{x}_1}{|\overline{x}| |\overline{x}_1|}+\frac{\e^2\overline{x}_n}{2|\overline{x}|^2|\overline{x}_1|^2}
%%%%  -\frac{\e^2 \left(\overline{x}_n (\overline{x}_n^2 
%%%%  -|\overline{x}|^2 + \overline{x}_1^2)\right)}{
%%%%  2 \left(|\overline{x}|^3 (|\overline{x}|^2 - \overline{x}_n^2) \sqrt{1 - 
%%%%  \frac{\overline{x}_n^2}{|\overline{x}|^2}}\right)}
+o(\e^2).\end{eqnarray*}
This yields that
\begin{eqnarray*}&&
\partial^2_{\eta_2\eta_2} w(\overline{x})\\&=&\lim_{\e\to0}\frac{
w_0(|\overline x+\e\eta_2|,\Theta(\overline x+\e\eta_2))
+w_0(|\overline x-\e\eta_2|,\Theta(\overline x-\e\eta_2))
-2w_0(|\overline x|,\Theta(\overline x))
}{\e^2}
\\&=&\lim_{\e\to0}\frac{1}{\e^2}\Bigg(
w_0\left(|\overline{x}|+
\frac{\e^2}{2|\overline{x}|}+o(\e^2),\Theta(\overline{x})
+\frac{\e\overline{x}_1}{|\overline{x}| |\overline{x}_1|}+\frac{\e^2\overline{x}_n}{2|\overline{x}|^2|\overline{x}_1|^2}
+o(\e^2)\right)\\&&\qquad\qquad
+w_0\left(|\overline{x}|+
\frac{\e^2}{2|\overline{x}|}+o(\e^2),\Theta(\overline{x})
-\frac{\e\overline{x}_1}{|\overline{x}| |\overline{x}_1|}+\frac{\e^2\overline{x}_n}{2|\overline{x}|^2|\overline{x}_1|^2}
+o(\e^2)\right)
-2w_0(|\overline x|,\Theta(\overline x))\Bigg)
\\&=&
\frac{
\partial_{1}w_0(\overline{x},\Theta(\overline{x}))}{|\overline{x}|}
+\frac{\partial_{2}w_0(\overline{x},\Theta(\overline{x}))\,\overline{x}_n
}{|\overline{x}|^2|\overline{x}_1|}
+\frac{\partial^2_{22}w_0(\overline{x},\Theta(\overline{x}))}{|\overline{x}|^2}.
\end{eqnarray*}
Gathering these observations, we conclude that
\begin{eqnarray*}
\Delta w(\overline{x})&=&\sum_{j=1}^n
\partial^2_{\eta_j\eta_j} w(\overline{x})\\
&=&
\partial^2_{11}w_0(\overline{x},\Theta(\overline{x}))
+\frac{\partial^2_{22}w_0(\overline{x},\Theta(\overline{x}))
}{|\overline{x}|^2}
+\sum_{j=2}^n\left[\frac{\partial_1 w_0 (|\overline{x}|,\Theta(\overline{x}) )}{|\overline{x}|}
+
\frac{\overline x_n\,\partial_2 w_0 (|\overline{x}|,\Theta(\overline{x}) )
}{|\overline{x}|^2|\overline x_1|}\right]\\
&=&
\partial^2_{11}w_0(\overline{x},\Theta(\overline{x}))
+\frac{\partial^2_{22}w_0(\overline{x},\Theta(\overline{x}))
}{|\overline{x}|^2}
+(n-1)\left[\frac{\partial_1 w_0 (|\overline{x}|,\Theta(\overline{x}) )}{|\overline{x}|}
+
\frac{\overline x_n\,\partial_2 w_0 (|\overline{x}|,\Theta(\overline{x}) )
}{|\overline{x}|^2|\overline x_1|}\right]
.
\end{eqnarray*}
In particular, if~$w(x)=|x|^\lambda f(\Theta(x))$ for some~$\lambda\in\R$
(i.e., $w_0(r,\theta)=r^\lambda f(\theta)$) we get that
\begin{eqnarray*}
\Delta w(\overline{x})&=&
\lambda(\lambda-1)|\overline{x}|^{\lambda-2}f(\Theta(\overline{x}))
+
|\overline{x}|^{\lambda- 2}f''(\Theta(\overline{x}))
\\&&\quad+(n-1)
\left[ \lambda|\overline{x}|^{\lambda-2}f(\Theta(\overline{x}) )
+\frac{\overline x_n\, |\overline{x}|^{\lambda -2}
f'(\Theta(\overline{x}) )}{|\overline x_1|}
\right].
\end{eqnarray*}
{F}rom this and~\eqref{0OKMSMS-2245ne} we obtain that 
\begin{equation}\label{sdk3i9v58b9648694689677}
|\overline{x}|^{2-\lambda}\Delta w(\overline{x})=
\lambda(\lambda+n-2)f(\Theta(\overline{x}))
+f''(\Theta(\overline{x}))+(n-1)\cot(\Theta(\overline{x}) )
f'(\Theta(\overline{x}) ).
\end{equation}
We remark that~$\cot(\Theta(\overline{x}) )>\cot b\ge-|\cot b|$ for all~$\overline{x}\in{\mathcal{C}}_b$.

Now we set~$K:=1+2(n-1)|\cot b|$ and we consider a small parameter~$\mu\in(0,1)$.
We choose
$$ f(\theta):=1-\frac{2\mu}{K^2} \left(e^{K\theta}+1\right)+\frac{2\mu\theta}K $$
and we observe that, for all~$\theta\in[0,\pi)$,
\begin{eqnarray*}
|f(\theta)|\le2,\qquad
f'(\theta)=\frac{2\mu}{K}\big(1-e^{K\theta}\big)
\qquad{\mbox{ and }}\qquad
f''(\theta)=-{2\mu}e^{K\theta},
\end{eqnarray*}
provided that~$\mu$ is sufficiently small.

Plugging this information into~\eqref{sdk3i9v58b9648694689677} and assuming~$\lambda\in(0,\mu^2]$ and~$\mu$ small enough, we obtain that
\begin{eqnarray*}
|\overline{x}|^{2-\lambda}\Delta w(\overline{x})&\le&
2\lambda n
-{2\mu}e^{K\Theta(\overline{x})}
+\frac{2\mu (n-1)\cot(\Theta(\overline{x}) )}{K}
\big(1-e^{K\Theta(\overline{x})}\big)\\
&\le&2\mu^2 n
-{2\mu}e^{K\Theta(\overline{x})}
+\frac{2\mu (n-1)\,|\cot b|}{K}
\big(e^{K\Theta(\overline{x})}-1\big)\\
&\le&2\mu^2 n
-{2\mu}e^{K\Theta(\overline{x})}
+\mu 
\big(e^{K\Theta(\overline{x})}-1\big)\\&=&
2\mu^2 n
-{\mu}e^{K\Theta(\overline{x})}
-\mu\\&\le&-{\mu}e^{K\Theta(\overline{x})},
\end{eqnarray*}
which is negative.

Also, $w>0$ in~$\overline{{\mathcal{C}}_b}\setminus\{0\}$
since
$$ f(\theta)\ge 1- 2\mu \left(e^{K\pi}+1\right)-
2\mu\pi\ge\frac12$$
as long as~$\mu$ is small enough.
This gives that~$w$ satisfies the requirements in~\eqref{JKA:BARRLPAER}, as desired.
\end{proof}

{F}rom Theorems~\ref{PERRO} and~\ref{PERRO-2} we deduce the following\footnote{A more precise version of
Corollary~\ref{S-coroEXIS-M023} will follow from the 
Schauder estimates in Proposition~\ref{SI:IN:SE:DR}.} existence result
for the Dirichlet problem in sufficiently regular domains:

\begin{corollary}\label{S-coroEXIS-M023}
Let~$\Omega\subset\R^n$ be open,
bounded and connected. Assume that either~$n=2$ and~$\Omega$
is accessible by arcs from the exterior, or that~$n\ge3$ and~$\Omega$
satisfies the exterior cone condition.

Let~$f\in C^1(\Omega')$ for some open set~$\Omega'\Supset\Omega$
and~$g\in C(\partial\Omega)$.

Then, there exists a unique solution~$u\in C^2(\Omega)\cap C(\overline\Omega)$
of the Dirichlet problem
$$ \begin{dcases}
\Delta u=f &{\mbox{ in }}\Omega,\\
u=g &{\mbox{ on }}\partial\Omega.
\end{dcases}$$
\end{corollary}

\begin{proof} The uniqueness claim is a consequence of
Corollary~\ref{UNIQUENESSTHEOREM}.

As for the existence claim, without loss of generality we can suppose that~$\Omega'$ is bounded
and we consider another bounded open set~$\Omega''$ such that~$\Omega\Subset\Omega''\Subset\Omega'$.
We take~$\tau\in C^\infty_0(\Omega')$ with~$\tau=1$ in~$\Omega''$ and extend~$f$ to the whole of~$\R^n$ by setting
$$ \widetilde{f}(x):=
\begin{dcases}
\tau(x)\,f(x) &{\mbox{ if }}x\in\Omega',\\
0 &{\mbox{ if }}x\in\R^n\setminus\Omega'.
\end{dcases}$$
We apply Proposition~\ref{ESI:NICE:FB} (with~$f$ there replaced by~$\widetilde{f}$ here)
and we construct a function~$v\in C^2(\R^n)$ such that~$\Delta v=\widetilde{f}$ in~$\R^n$.
In particular, we have that~$\Delta v=f$ in~$\Omega$.

Now, for each~$x\in\partial\Omega$, we set~$\widetilde g(x):=g(x)-v(x)$.
We remark that~$\widetilde{g}\in C(\partial\Omega)$, hence we can employ
Theorems~\ref{PERRO} and~\ref{PERRO-2} and find a function~$\widetilde{u}\in C(\overline\Omega)$
which is harmonic in~$\Omega$ and equal to~$\widetilde g$ on~$\partial\Omega$.
Thus, we set~$u:=v+\widetilde{u}$. We see that~$u\in C^2(\Omega)\cap C(\overline\Omega)$,
since both~$v$ and~$\widetilde{u}$ belong to these spaces, that~$\Delta
u=\Delta v+\Delta\widetilde{u}=f$ in~$\Omega$, and that~$u=v+\widetilde{g}=g$ along~$\partial\Omega$.
\end{proof}

Concerning the result in
Corollary~\ref{S-coroEXIS-M023},
we stress that
some geometric conditions on the domain are
needed to ensure the continuity of the solution
at the boundary. Indeed, as an example, 
we remark that 
\begin{equation}\label{SALT6}\begin{split}&
{\mbox{there exists no harmonic function~$u$
in~$B_1\setminus\{0\}$ that is continuous in~$\overline{B_1\setminus\{0\}}=
\overline{B_1}$}}\\&{\mbox{and such that~$u=0$ on~$\partial B_1$ and~$u(0)=1$,}}\end{split}\end{equation}
otherwise, by Theorem~\ref{0oiujhgv098j-23rtgrJSPsdfgdh56uh},
such a function could be extended to a harmonic function in~$B_1$,
which should necessarily vanish identically in view of the
Maximum Principle in Corollary~\ref{WEAKMAXPLE}.

For completeness, we point out that another proof of~\eqref{SALT6}
can be obtained by exploiting spherical averages
and the classification result in Lemma~\ref{MHNvKSD}. Indeed,
given~$r\in(0,1)$ and~$\vartheta\in\partial B_1$, let~$u_0(r,\vartheta):=u(r\vartheta)$.
For all~$x\in \overline{B_1}$, let also
$$ v(x):=\fint_{\partial B_1} u(|x|\vartheta)\,d{\mathcal{H}}^{n-1}_\vartheta=
\fint_{\partial B_1} u_0(|x|,\vartheta)\,d{\mathcal{H}}^{n-1}_\vartheta.$$
By the spherical representation of the Laplacian in Theorem~\ref{SPHECOO}, we know that,
for all~$x\in B_1\setminus\{0\}$, \begin{equation}\begin{split}\label{Byphericalpresentatitlaciorem}
\Delta v(x)\,&=\,
\fint_{\partial B_1} \left(
\partial_{rr}u_0(|x|,\vartheta)+
\frac{n-1}{|x|}\,\partial_{r}u_0(|x|,\vartheta)+\frac1{|x|^2}\,
\Delta_{\partial B_1}u_0(|x|,\vartheta)
\right)\,d{\mathcal{H}}^{n-1}_\vartheta\\&
=\,\fint_{\partial B_1} \Delta u(|x|\vartheta)\,d{\mathcal{H}}^{n-1}_\vartheta\end{split}\end{equation}
and accordingly $\Delta v(x)=0$.
Thus, since~$v$ is rotationally symmetric, we can apply
Lemma~\ref{MHNvKSD} and gather that, for every~$x\in B_1\setminus\{0\}$,
\begin{equation}\label{89-jdj4955t-1} v(x)=\begin{dcases}
\frac{a}{|x|^{n-2}}+b&{\mbox{ if }}n\ne2,\\
a\ln |x|+b&{\mbox{ if }}n=2
\end{dcases}\end{equation}
for some~$a$, $b\in\R$. Now, if~$u$ were continuously attaining the datum~$1$
at the origin, for every~$\e>0$ there would exist~$\delta>0$ such that if~$x\in B_\delta$ then~$u(x)\in [1-\e,1+\e]$, hence~$v(x)\in[1-\e,1+\e]$ too, namely
$$ \lim_{x\to0}v(x)=1.$$
This and~\eqref{89-jdj4955t-1} yield that necessarily~$a=0$, and also~$b=1$.
As a result, we have that~$v$ is identically equal to~$1$ in~$B_1$.
On the other hand, if~$u$ were continuous in~$\overline B_1$
and vanishing along~$\partial B_1$, there would exist~$\rho\in(0,1)$ such that~$u(x)\le\frac12$
for all~$|x|\in[\rho,1]$, leading to
$$ 1=v(\rho e_n)=\fint_{\partial B_1} u(\rho\vartheta)\,d{\mathcal{H}}^{n-1}_\vartheta
\le\frac12.$$
This contradiction provides an alternative proof of~\eqref{SALT6}.

Note that~\eqref{SALT6} highlights that
that the assumption
of being accessible by arcs from the exterior for a planar\footnote{It is interesting
to ponder on the physical counterpart of this example: one might
translate this construction into the phenomenon of a perfectly elastic planar membrane
that is constrained along a circle and with a sharp needle placed at the center of the circle.
The mathematical example suggests that the needle would brake the membrane and
pass through it.}
domain in Corollary~\ref{S-coroEXIS-M023}
cannot be removed.

\begin{figure}
  \centering
  \includegraphics[width=.44\linewidth]{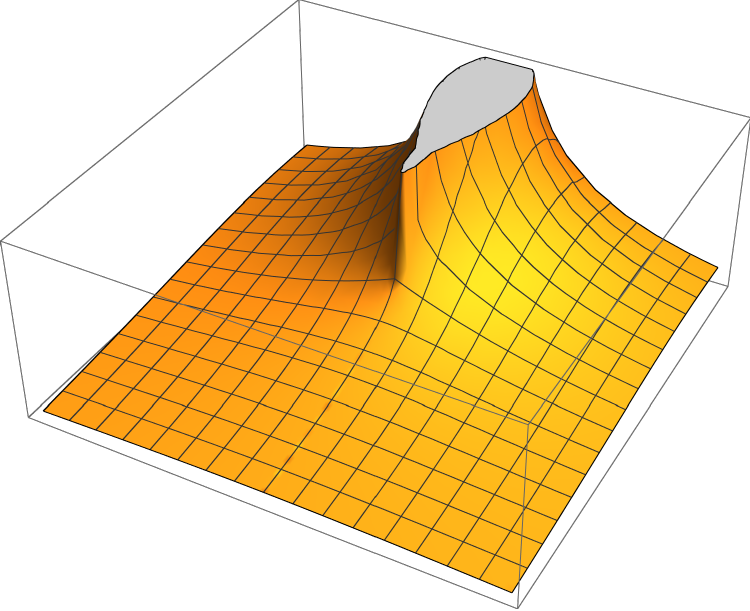}$\quad$\includegraphics[width=.44\linewidth]{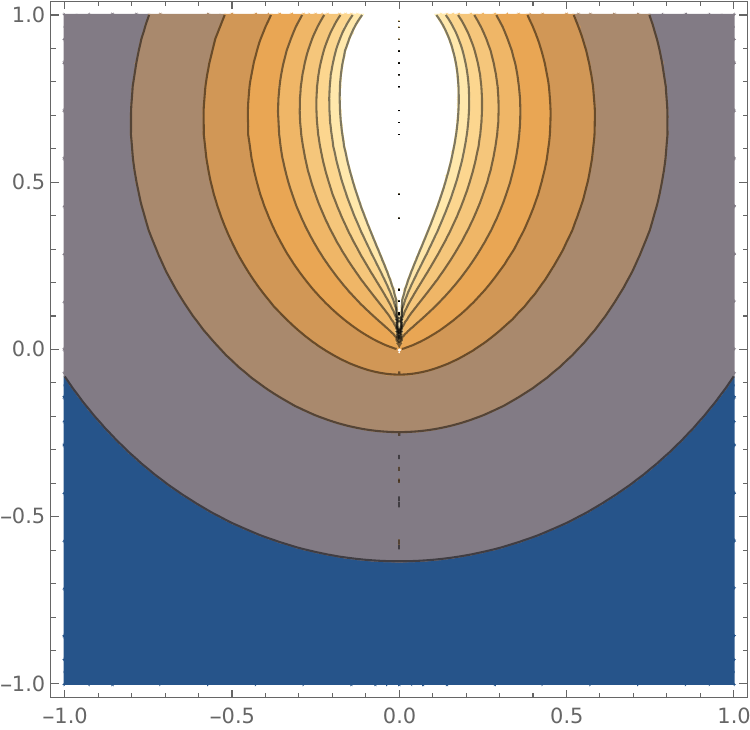}
 \caption{\sl The graph of~$v$ in~\eqref{Abdr5oguNSt5h6eJSDGLSobe}
 and of its level sets as a function of~$(r,x_3)$.}\label{67390-38494hf29T536EcdhPOaIJodidpo3tju24ylkiliR2543654yh2432RAe87658FI6tfh564}
\end{figure}

In the same spirit, we stress that
when~$n\ge3$,
classical examples due to Henri Lebesgue
(and sometimes
referred to with the name of ``Lebesgue spine''\index{Lebesgue spine}, \label{SILEBE678NE}
see e.g.~\cite[page 304, figure 19]{MR1013360},
\cite[Section~7.2]{MR1306729} and~\cite[pages~175--176]{MR0261018})
exhibit cuspidal
domains in~$\R^3$ with a prescribed
continuous datum at the boundary
in which the solution
cannot meet the boundary datum in a continuous way.
Let us briefly recall here one of these examples. \label{CUSPI-01}
The main idea underpinning this example is to look at the electrostatic potential of a charged segment, with a suitable choice of the charge distribution, then a cuspidal domain having this segment in its complement
will do the job.
To make this idea concrete, we observe that, given~$a$, $b\in\R$ with~$a^2 + (b + s)^2>0$,
if
$$ \Phi(s,a,b):=
\sqrt{a^2 + (b + s)^2} +\frac{b}{2}\ln\left(1-\frac{b + s}{\sqrt{a^2 + (b + s)^2}}\right)
-\frac{b}{2}\ln\left(1+\frac{b + s}{\sqrt{a^2 + (b + s)^2}}\right),$$
then
\begin{equation*}\begin{split}& \frac{\partial\Phi}{\partial s}(s,a,b)=
\frac{s}{\sqrt{a^2 + (s + b)^2}}.\end{split}\end{equation*}
Therefore, if we define~$r:=\sqrt{x_1^2+x_2^2}$ and
\begin{equation}\label{DEV-v-LEBE} v(x)=v(x_1,\dots,x_n):=
\int_0^1\frac{s}{|(x_1,x_2,x_3)-(0,0,s)|}\,ds=
\int_0^1\frac{s}{\sqrt{r^2 + (s -x_3)^2}}\,ds,\end{equation}
we have that
\begin{equation}\label{Abdr5oguNSt5h6eJSDGLSobe}
\begin{split} v(x)\,=\,&\int_0^1\frac{\partial\Phi}{\partial s}(s,r,-x_3)\,ds\\
=\,&\Phi(1,r,-x_3)-\Phi(0,r,-x_3)
\\=\,&
\sqrt{r^2 + (1 -x_3)^2} - \frac{x_3}2 \ln\left(1 - \frac{1 -x_3}{\sqrt{r^2 + (1 -x_3)^2}}\right) +
\frac{x_3}{2} \ln\left(1+ \frac{1 -x_3}{\sqrt{r^2 + (1 -x_3)^2}}\right)\\&\qquad\qquad\quad\qquad
-\sqrt{r^2 +x_3^2} + \frac{x_3}2 \ln\left(1 +\frac{x_3}{\sqrt{r^2 + x_3^2}}\right)
-\frac{x_3}2 \ln\left(1 +-\frac{x_3}{\sqrt{r^2 + x_3^2}}\right)\\=\,&
\sqrt{r^2 + (1 -x_3)^2} -\sqrt{r^2 +x_3^2} \\&\qquad\qquad\quad\qquad+ \frac{x_3}2 
\ln\frac{\left({\sqrt{r^2 + (1 -x_3)^2}+1 -x_3}\right)
\left(\sqrt{r^2 + x_3^2}+x_3\right)}{\left(\sqrt{r^2 + (1 -x_3)^2}
-1+x_3\right)\left(\sqrt{r^2 + x_3^2}-x_3\right)} .
\end{split}\end{equation}
See Figure~\ref{67390-38494hf29T536EcdhPOaIJodidpo3tju24ylkiliR2543654yh2432RAe87658FI6tfh564} for the graph of~$v$ and of its level sets as a function of~$(r,x_3)$.

\begin{figure}
  \centering
  \includegraphics[width=.5\linewidth]{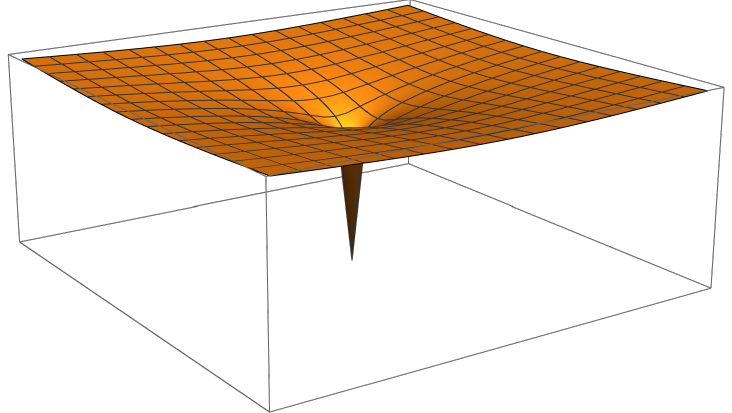}
 \caption{\sl The boundary of the domain in~\eqref{LEBVEasfdghRAeFI2}
 near the origin (when~$n=3$).}\label{LEBVEasfdghRAeFI}
\end{figure}

Now we consider a bounded domain~$\Omega$ whose boundary is~$C^\infty$ outside the origin and such that
\begin{equation}\label{LEBVEasfdghRAeFI2} \Omega\cap B_{1/2}=\left\{x=(x_1,\dots,x_n)\in\R^n {\mbox{ s.t. }} x_3<-\frac{1}{2\ln|(x_1,x_2)|}
\right\}.\end{equation}
See Figure~\ref{LEBVEasfdghRAeFI} for a sketch of the cuspidal boundary of this domain
in the proximity of the origin.

We point out that~$v$ is continuous on~$\partial\Omega$, because,
by~\eqref{Abdr5oguNSt5h6eJSDGLSobe},
\begin{equation} \label{LOGLIM}\begin{split}&\lim_{{x\to0}\atop{x_3=-1/(2\ln|(x_1,x_2)|)}}v(x)
=1- \lim_{r\to0}\frac{1}{4\ln r} \,\ln\frac{\left({\sqrt{r^2 + \left(1 +\frac{1}{2\ln r}\right)^2}+1
+\frac{1}{2\ln r}}\right)
\left(\sqrt{r^2 + 
\frac{1}{4\ln^2 r}}-\frac{1}{2\ln r}\right)}{\left(
\sqrt{r^2 + \left(1 +\frac{1}{2\ln r}\right)^2}-1-\frac{1}{2\ln r}
\right)\left(\sqrt{r^2 + \frac{1}{4\ln^2 r}} +
\frac{1}{2\ln r}
\right)} \\&\qquad=1- \lim_{r\to0}\frac{1}{4\ln r} \,\ln\frac{\left({\sqrt{r^2 + \left(1 +\frac{1}{2\ln r}\right)^2}+1
+\frac{1}{2\ln r}}\right)^2
\left(\sqrt{r^2 + 
\frac{1}{4\ln^2 r}}-\frac{1}{2\ln r}\right)^2}{r^4} 
\\
&\qquad=
1- \lim_{r\to0}\frac{1}{4\ln r} \,\ln\frac{\left(2+O\left(\frac1{\ln r}\right)\right)^2\frac1{4\ln^2 r }
\left(\sqrt{4r^2\ln^2r + 1}+1\right)^2}{r^4} \\
&\qquad=1- \lim_{r\to0}\frac{1}{4\ln r} \,\ln\frac{\left(2+O\left(\frac1{\ln r}\right)\right)^2\frac1{4\ln^2 r }
\left(2+O(r^2\ln^2r )\right)^2}{r^4}\\&\qquad=
1- \lim_{r\to0}\frac{1}{4\ln r} \,\left\{\ln\left[\left(2+O\left(\frac1{\ln r}\right)\right)^2\frac1{4\ln^2 r }
\left(2+O(r^2\ln^2r )\right)^2\right]-\ln r^4\right\}\\&\qquad=
2-\lim_{r\to0}\frac{1}{4\ln r} \,\ln\left[\left(2+O\left(\frac1{\ln r}\right)\right)^2\frac1{4\ln^2 r }
\left(2+O(r^2\ln^2r )\right)^2\right]\\&\qquad=
2-\lim_{r\to0}\frac{1}{4\ln r} \,\ln
\frac{4+o(1)}{4\ln^2 r }\\&\qquad=2+\lim_{r\to0}\frac{\ln
(4\ln^2 r )}{4\ln r}\\&\qquad=2
.\end{split}\end{equation}

Now, let~$I:=\{te_n,$ $t\in[0,1]\}$. By~\eqref{DEV-v-LEBE}, noticing that~$v$ depends only on~$X:=(x_1,x_2,x_3)$, and recognizing the fundamental solution in~\eqref{GAMMAFU}
in dimension~$n=3$, we see that~$v\in C^\infty(\R^n\setminus I)\subset C^\infty(\Omega)$ and,
for every~$x\in\R^n\setminus I$,
\begin{equation*}
\Delta v(x)=\Delta_{X}v(x)=
\int_0^1\Delta_X\frac{s}{|X-(0,0,s)|}\,ds=0.
\end{equation*}
These observations and the uniqueness result in Corollary~\ref{UNIQUENESSTHEOREM} give
that~$v$ is the unique harmonic function in~$\Omega$ coinciding with its own boundary
values along~$\partial\Omega$.

However, recalling~\eqref{DEV-v-LEBE}, if~$x_3<0$ then
$$ v(x)=
\int_0^1\frac{s}{\sqrt{r^2 + (s +|x_3|)^2}}\,ds
\le\int_0^1\frac{s}{\sqrt{s^2}}\,ds=1,$$
hence, comparing with~\eqref{LOGLIM}, we see that~$v$ is not continuous at the origin.
This example shows that the exterior cone condition in Corollary~\ref{S-coroEXIS-M023}
cannot be removed.\label{CUSPI-02}
\medskip

A natural question arising from Corollary~\ref{S-coroEXIS-M023}
is whether one can characterize precisely the domains for which
a solution of the Dirichlet problem that is continuous up to the boundary
always exists. We will address this question via a criterion due to Norbert Wiener
in the forthcoming Theorem~\ref{NOWI}.
\medskip

As a technical remark, it is interesting to observe that the existence results obtained
by the Perron method are conceptually different from the ones obtained
by energy or variational methods (compare e.g. Sections~2.2.5(b) and~6.2 in~\cite{MR1625845}).
In particular, solutions constructed via the Perron method do not necessarily
possess finite energy. To highlight this phenomenon, we recall a classical example
due to Jacques Hadamard~\cite{MR1504545} (see also~\cite{MR1611073}
for a comprehensive historical perspective). For this, we use polar coordinates~$(r,\vartheta)$
in the plane and consider the function on~$B_1\subset\R^2$ defined by
$$ u(r,\vartheta):=\sum_{k=1}^{+\infty}
\frac{r^{k!}}{k^2}\sin(k!\vartheta).$$
Recalling the representation of the Laplacian in polar coordinates discussed
in Theorem~\ref{SPHECOO}, we know that, for every~$k\in\N$,
\begin{eqnarray*}&& \Delta\left( r^{k!}\sin(k!\vartheta)\right)
=\left(\partial_{rr}+
\frac{1}{r}\partial_r+\frac{1}{r^2}\partial_{\vartheta\vartheta}\right)
\left( r^{k!}\sin(k!\vartheta)\right)\\&&\qquad=\Big(
k!(k!-1) 
+k!
- (k!)^2 \Big)r^{k!-2}\sin(k!\vartheta)=
0\end{eqnarray*}
and consequently the function
$$ u_N(r,\vartheta):=\sum_{k=1}^{N}
\frac{r^{k!}}{k^2}\sin(k!\vartheta)$$
is harmonic in~$B_1$ for all~$N\in\N\cap[1,+\infty)$.

Moreover,
$$ |u(r,\vartheta)-u_N(r,\vartheta)|\le\sum_{k=N+1}^{+\infty}
\frac{r^{k!}}{k^2}\le\sum_{k=N+1}^{+\infty}
\frac{1}{k^2},$$
which is infinitesimal as~$N\to+\infty$. This shows that~$u_N$ converges uniformly to~$u$
in~$\overline{B_1}$ and therefore~$u$ is continuous in~$\overline{B_1}$ and harmonic in~$B_1$ (recall
Corollary~\ref{S-OS-OSSKKHEPARHA}, and in fact it is the only function with these properties, due to the
uniqueness result in Corollary~\ref{UNIQUENESSTHEOREM}).

But~$u$ has infinite energy, that is
\begin{equation}\label{ADFVSRYUANSOKMASCLAPAKS2}
\int_{B_1}|\nabla u|^2=+\infty.
\end{equation}
To check this, we argue by contradiction, assuming the converse, say~$\int_{B_1}|\nabla u|^2\le M$
for some~$M\in[0,+\infty)$. Then, since~$|\nabla u|\ge\frac{1}{r}|\partial_\vartheta u|$,
$$ M\ge\iint_{(0,1)\times(0,2\pi)} \frac{1}{r}|\partial_\vartheta u(r,\vartheta)|^2\,dr\,d\vartheta$$
and thus, a.e.~$r\in(0,1)$,
$$ \int_{(0,2\pi)} |\partial_\vartheta u(r,\vartheta)|^2\,d\vartheta<+\infty.$$
This says that, a.e.~$r\in(0,1)$, the map~$(0,2\pi)\ni\vartheta\mapsto f(\vartheta):=\partial_\vartheta u(r,\vartheta)$ belongs to~$L^2((0,2\pi))$, and we can therefore consider the corresponding
Fourier coefficients (up to normalization)
$$ f_\ell:=\int_{(0,2\pi)} f(\vartheta)\,e^{- i\vartheta \ell} \,d\vartheta$$
and by Plancherel Theorem (see e.g.~\cite{MR3381284}) we find that
\begin{equation}\label{ADFVSRYUANSOKMASCLAPAKS} \sum_{\ell\in\Z} |f_\ell|^2<+\infty.\end{equation}
Note that, via an integration by parts,
\begin{eqnarray*}f_\ell&=&\int_{(0,2\pi)} \partial_\vartheta u(r,\vartheta)\,e^{- i\vartheta \ell} \,d\vartheta\\&
=& i\ell\int_{(0,2\pi)} u(r,\vartheta)\,e^{- i\vartheta \ell} \,d\vartheta \\&=& i\ell\lim_{N\to+\infty}\int_{(0,2\pi)} u_N(r,\vartheta)\,e^{- i\vartheta \ell} \,d\vartheta
\\&=& i\ell\lim_{N\to+\infty}\sum_{k=1}^{N}\int_{(0,2\pi)}
\frac{r^{k!}}{k^2}\sin(k!\vartheta)\,e^{- i\vartheta \ell} \,d\vartheta\\&=&
\pi\lim_{N\to+\infty} \sum_{k=1}^{N} 
\frac{\ell\,\delta_{\ell,k!}\,r^{k!} }{k^2}.
\end{eqnarray*}
In particular, for every~$m\in\N$
$$ |f_{m!}|=\left|
\pi\lim_{N\to+\infty} \sum_{k=1}^{N} 
\frac{m!\,\delta_{m!,k!}\,r^{k!}}{k^2}
\right|=\pi 
\frac{m!\,r^{m!}}{m^2}
,$$
which entails that
$$ \sum_{\ell\in\Z} |f_\ell|^2\ge\sum_{m\in\N} |f_{m!}|^2=+\infty.$$
This is in contradiction with~\eqref{ADFVSRYUANSOKMASCLAPAKS}
and the proof of~\eqref{ADFVSRYUANSOKMASCLAPAKS2} is thus complete.
\medskip

For a thorough comparison between energy techniques and the Perron method see~\cite{MR2409177}. 
\medskip

For additional observations on the Perron method
and its relation with the ``obstacle problem'' see e.g.~\cite{MR0222317, MR0261018, MR1814364, MR3675703}.
A different method for solving the Dirichlet problem
was also proposed by Poincar\'e and relied on subsequent harmonic replacements
in chains of balls covering the domain, see e.g.~\cite[Section~1.4.3]{MR3099262}
for further details on this technique.

\chapter{Equations in nondivergence form: $C^{2,\alpha}$-regularity theory}\label{C2ALPHACHAPET}

\section{Hints and limitations for a smooth regularity theory}

We discuss here some bits of the regularity theory
for elliptic equations (for a thorough presentation, see e.g.~\cite{MR1625845, MR1814364, MR2777537}
and the references therein).

As a first observation, we remark that the situation for partial differential equations
is structurally very different from that of ordinary differential equations.
Indeed, for (say, linear) ordinary differential equations the regularity theory is mostly straightforward:
for instance, if~$I$ is an interval in~$\R$ and~$u:I\to\R$ solves~$\ddot{u}(t)=f(t)$ for all~$t\in I$
then ``$u$ is two derivatives better than~$f$'', as it can be seen by integrating twice the equation.

For partial differential equations, the situation can be significantly more complicated.
For general equations, one cannot expect any kind of regularity at all: for example,
given any function~$u_0:\R\to\R$, the function~$u(x_1,x_2):=u_0(x_1-x_2)$ is
a solution of
\begin{equation}\label{mpVpEjsodvnMSMDomebcvbNThbndoAnsdjGenfLO-KSM-323ty}
\partial_{11}u-\partial_{22} u=0.
\end{equation}
While in principle to write such an equation one needs to assume that~$u_0$
is twice differentiable, no quantitative estimate on the second derivative of~$u_0$ comes into
play, therefore smooth solutions of~\eqref{mpVpEjsodvnMSMDomebcvbNThbndoAnsdjGenfLO-KSM-323ty}
can exhibit arbitrarily large derivatives. This suggests
that, in a suitable jargon, equation~\eqref{mpVpEjsodvnMSMDomebcvbNThbndoAnsdjGenfLO-KSM-323ty} can
be interpreted in a ``weak form'' which does not require the solution~$u$
to possess basically any regularity at all,
and no regularity can be inferred for
this type of weak solutions of~\eqref{mpVpEjsodvnMSMDomebcvbNThbndoAnsdjGenfLO-KSM-323ty}.

The situation for the Laplace operator is on the other hand radically different.
For instance, the weak solutions of~$\Delta u=0$,
as discussed in Section~\ref{WERGSB-EAFGBS-cojwnedfDICewew4N8wedibfn},
happen to be as smooth as we wish.
This suggests that the Laplace operator possesses some kind of ``regularizing effect''.
Indeed, even if the ``algebraic difference'' between the Laplacian
and the operator in~\eqref{mpVpEjsodvnMSMDomebcvbNThbndoAnsdjGenfLO-KSM-323ty}
seems negligible (just a minus sign!), the ``geometric structure''
underpinning the Laplace operator tends to ``average out'' differences and oscillations
(recall e.g. Theorem~\ref{KAHAR2} to appreciate how the Laplacian tends
to revert the pointwise values of a function
to the average nearby).

In dimension~$2$, another strong indication that solutions of equations driven by the Laplace operator
enjoy additional rigidity and regularity properties comes from the fact that harmonic functions
in planar domains can be seen as real (or imaginary) parts of holomorphic functions
in complex domains (with the identification of~$(x,y)\in\R^2$ with~$x+iy\in\cOMPL$): with respect
to this, we recall that the existence of a complex derivative in a neighbourhood entails
that the complex function is actually infinitely differentiable and, in fact, real analytic.

For global solutions (i.e., solutions in the whole of~$\R^n$), another quantitative hint
of the regularity enjoyed by solutions of~$\Delta u(x)=f(x)$ for all~$x\in\R^n$
and for~$f$ in the Schwartz space of smooth functions whose derivatives are rapidly decreasing at infinity
comes from Fourier analysis: indeed, after a Fourier Transform, the equation becomes~$-4\pi^2|\xi|^2\widehat u(\xi)=\widehat f(\xi)$.
Thus, recalling that, for each~$j$, $m\in\{1,\dots,n\}$, ${\mathcal {F}}\left(\partial^2_{jm} u\right)(\xi)=
-4\pi^2 \xi_j\xi_m\widehat u(\xi)$, it follows from the Plancherel Theorem
that
\begin{eqnarray*}&&\| D^2u\|_{L^2(\R^n)}=\sqrt{\sum_{j,m=1}^n \int_{\R^n} (\partial^2_{jm} u(x))^2\,dx}
=\sqrt{\sum_{j,m=1}^n \int_{\R^n} \left|
{\mathcal {F}}\left(\partial^2_{jm} u\right)(\xi)\right|^2\,d\xi}\\&&\qquad
=4\pi^2\sqrt{\sum_{j,m=1}^n \int_{\R^n} |\xi_j|^2|\xi_m|^2|\widehat u(\xi)|^2\,d\xi}
\le
2\sqrt{2}\pi^2\sqrt{\sum_{j,m=1}^n \int_{\R^n} (|\xi_j|^4+|\xi_m|^4)|\widehat u(\xi)|^2\,d\xi}
\\&&\qquad=
4\pi^2\sqrt{\sum_{j,m=1}^n \int_{\R^n} |\xi_j|^4|\widehat u(\xi)|^2\,d\xi}
\le
4\pi^2 {n}\,\sqrt{\int_{\R^n} |\xi|^4|\widehat u(\xi)|^2\,d\xi}\\&&\qquad
={n}\,\sqrt{\int_{\R^n} |\widehat f(\xi)|^2\,d\xi}={n}\,\sqrt{\int_{\R^n} |f(x)|^2\,dx}={n}\,
\|f\|_{L^2(\R^n)},
\end{eqnarray*}
which can be seen as a global bound in~$L^2(\R^n)$ of the second derivative of the solution~$u$
with respect to the bound~$L^2(\R^n)$ of the source term~$f$.\medskip

While these heuristic (or very specific)
considerations suggest that solutions of Laplace-like (in jargon, elliptic)
equations should possess additional regularity, detecting this regularity is typically
an extremely hard task, it relies on a number of technical (albeit beautiful) arguments
and, in general, it is very difficult to ``guess the right results'' only based\footnote{For instance,
the example pointed out in Theorem~\ref{CONSDTYROE-Pse3m4pa235sfbscd3jv} reveals
that the classical notion of continuity alone is not the appropriate one to deal effectively
with partial differential equations (even in the more ``amenable'' case of the Laplace
operator): roughly speaking, differently from what happens for ordinary differential equations,
derivatives in different directions may end up being singular, but their singularities cancel out
after summation. However, this in itself does not fully explain the complexity
of the situation, since cancellations of singularities are not allowed when the source term is zero
(since harmonic functions, no matter how weak their definition is,
end up being real analytic, recall Theorems~\ref{KAHAR} and~\ref{KMS:HAN0-0},
as well as Lemma~\ref{WEYL}). 

The fact that will be apparent from the forthcoming pages
is that the appropriate functional spaces
to account for the continuous regularity theory of equations driven by the Laplacian,
and by more general elliptic operators, are the H\"older spaces~$C^{k,\alpha}$ with~$\alpha\in(0,1)$.
In this sense, not only one needs to quantify precisely the notion of continuity,
but also to take into account a sort of ``fractional'' regularity encoded in the power~$\alpha$
(for harmonic functions, or for solutions of equations with a $C^\infty$ source term,
one does not really see this issue since the equation can be differentiated
infinitely many times making any fractional type regularity be incorporated into the
subsequent continuous space with integer regularity: say, when the source term
can be differentiated infinitely many times the~$C^{k,\alpha}$ structure ends up being merged
into the stronger~$C^{k+1}$ framework, making the fractional exponent scenario much less visible).

The appearance and significance of this intrinsically fractional structure within the framework of classical
partial differential equations are perhaps surprising and even mysterious.
To be honest, we do not have a compelling and exhaustive explanation
on why ``intuitively'' the appearance of fractional type regularity should have been
expected just from our ``basic perception of the universe'', i.e. without relying on the
technical mathematical aspects of the equations. However, it may be suggestive to think
that the universe is perhaps much more fractal and fractional than what it appears to the naked eye:
many fundamental patterns of the world arise from fractional equations
(see e.g. the introduction of~\cite{CDV-MDG-19})
and the Laplace operator is, after all, one operator of a broader family
of objects which share a number of geometric properties. In this sense, the classical Laplacian
can be seen as a limit of fractional operators for which regularity in fractional spaces appears
naturally. Accepting this kind of reasoning, one can suggestively think that
the fractional regularity in H\"older spaces for the Laplacian is just the ``vestigial feature''
of this limit process.}
on ``physical intuition of what's going on''.
As an example (see e.g.~\cite[Problem~4.9]{MR1814364} and~\cite[Remark~3.14]{MR2777537}), we remark that solutions of~$\Delta u(x)=f(x)$ for a continuous function~$f$
do not need\footnote{For further reference, \label{FAV-0kmSRTS9rtboTEFSftS}
we point out that the examples provided in the two proofs of
Theorem~\ref{CONSDTYROE-Pse3m4pa235sfbscd3jv} that we give here
also show that there exist~$f\in L^\infty(B_1\setminus\{0\})$ and~$u\in C^2(B_1\setminus\{0\})$
for which~$\Delta u=f$ in~$B_1\setminus\{0\}$ but~$D^2u\not\in L^\infty({B_1}\setminus\{0\})$.

The analogue of this observation in~$L^1(B_1\setminus\{0\})$ will be discussed in
the forthcoming Theorem~\ref{FAV-0kmSRTS9rtboTEFSftS1}.} to be continuously twice differentiable (hence
$u$ is not ``always two derivatives better than~$f$''!):

\begin{theorem}\label{CONSDTYROE-Pse3m4pa235sfbscd3jv} 
There exists~$f\in C(\overline{B_1})$ for which the Dirichlet problem
\begin{equation}\label{DIRINON} \begin{dcases}
\Delta u=f & {\mbox{ in }}B_1,\\
u=0 & {\mbox{ on }}\partial B_1
\end{dcases}\end{equation}
does not admit any solution~$u\in C^2(B_1)\cap C(\overline{B_1})$.
\end{theorem}

\begin{proof} Up to a scaling argument, we prove the desired
result for~$B_{1/2}$ instead of~$B_1$. We let
\begin{equation}\label{MMS-NUS-OuohenNunSQidmqi4U5A3Lsdnvd} f(x):=
\frac{(x_1^2 -x_2^2) (2(n+2) \ln|x| - 1)}{4\,|x|^2 (-\ln|x|)^{3/2}},\end{equation}
for every~$(x_1,x_2,\dots,x_n)\in B_{1/2}$,
and we observe that~$f\in C(\overline{B_{1/2}})$ since
$$ \lim_{x\to0}|f(x)|\le
\lim_{x\to0}
\frac{2(n+2) |\ln|x||+ 1}{4(-\ln|x|)^{3/2}}=0.
$$
We claim that the function~$f$
in~\eqref{MMS-NUS-OuohenNunSQidmqi4U5A3Lsdnvd}
does not admit any solution~$u\in C^2(B_{1/2})\cap C(\overline{B_{1/2}})$
of the Dirichlet problem in~\eqref{DIRINON}. To prove this,
we argue by contradiction and we suppose that a solution exists. We consider the function
\begin{equation}\label{MMS-NUS-OuohenNunSQidmqi4U5A3Lsdnvd2} B_{1/2}\ni x=(x_1,x_2,\dots,x_n)\longmapsto v(x):=\sqrt{-\ln|x|}\,(x_1^2-x_2^2)\end{equation}
and we point out that~$v\in C^\infty(B_{1/2}\setminus\{0\})\cap C(\overline{B_{1/2}})$.
Also, by a direct computation we see that~$\Delta v=f$ in~$B_{1/2}\setminus\{0\}$.
Consequently, the function~$w:=v-u$ belongs to~$C^2(B_{1/2}\setminus\{0\})\cap C(\overline{B_{1/2}})$
and satisfies
\begin{equation*} \begin{dcases}
\Delta w=0 & {\mbox{ in }}B_{1/2}\setminus\{0\},\\
w=v & {\mbox{ on }}\partial B_{1/2}.\end{dcases}\end{equation*}
By Theorem~\ref{0oiujhgv098j-23rtgrJSPsdfgdh56uh}, we can extend~$w$
to a harmonic function in the whole of~$B_{1/2}$: in particular,
there exists~$w_\star\in C^2(B_{1/2})$ such that~$w_\star=w$ in~$B_{1/2}\setminus\{0\}$.
As a result, 
\begin{equation}\label{78jpikjhyFUYSV8yeoiuew82yrCOmnsdf}{\mbox{the function~$v_\star:=u+w_\star$ belongs to~$C^2(B_{1/2})$}}\end{equation}
and it is such that~$v_\star=v$ in~$B_{1/2}\setminus\{0\}$.

Therefore
\begin{eqnarray*}&& \lim_{x_2\searrow0}\partial_{11} v_\star(0,x_2, 0,\dots,0)=
\lim_{x_2\searrow0}\partial_{11} v(0,x_2, 0,\dots,0)=\lim_{{x_2\searrow0}}
\frac{ 1-4|\ln x_2|}{2\sqrt{- \ln  x_2 }}=+\infty,
\end{eqnarray*}
in contradiction with~\eqref{78jpikjhyFUYSV8yeoiuew82yrCOmnsdf}.
\end{proof}

\begin{proof}[Another proof of Theorem~\ref{CONSDTYROE-Pse3m4pa235sfbscd3jv}]
There are other possibilities for the functions~$f$ and~$v$ in~\eqref{MMS-NUS-OuohenNunSQidmqi4U5A3Lsdnvd} and~\eqref{MMS-NUS-OuohenNunSQidmqi4U5A3Lsdnvd2}.
An instructive alternative is provided by a superposition method and goes as follows.
One can replace~\eqref{MMS-NUS-OuohenNunSQidmqi4U5A3Lsdnvd2} with
$$ v(x):=\sum_{k=1}^{+\infty}\frac{Q(e^k x)}{e^{2k}\,k},$$
where~$Q(x):=\eta(x)h(x)$, where~$\eta\in C^\infty_0(B_1,[0,1])$ with~$\eta=1$ in~$B_{1/2}$
and~$h(x):=x_1^2-x_2^2$. Also, one replaces~\eqref{MMS-NUS-OuohenNunSQidmqi4U5A3Lsdnvd}
with
$$ f(x):=\Delta v(x).$$
To establish Theorem~\ref{CONSDTYROE-Pse3m4pa235sfbscd3jv},
according to the previous proof, one has to check that
\begin{eqnarray}
\label{CSUdijcINdfyhTUndRgbdEb-1}&& {\mbox{the definitions of~$v$ and~$f$ are well-posed,}}\\
\label{CSUdijcINdfyhTUndRgbdEb-2}&& f\in C(\overline{B_1}),\\
\label{CSUdijcINdfyhTUndRgbdEb-3}&& {\mbox{$v\not\in C^2(B_1)$.}}
\end{eqnarray}
For this, we observe that the definition of~$v$ is well-posed since~$\displaystyle\sum_{k=1}^{+\infty}\frac{1}{e^{2k}\,k}<+\infty$.
Furthermore, we have that~$\Delta Q=h \Delta\eta +\eta\Delta h+2\nabla\eta\cdot\nabla h=
h\Delta\eta+2\nabla\eta\cdot\nabla h$. As a result, if
$$ v_m(x):=\sum_{k=1}^{m}\frac{Q(e^k x)}{e^{2k}\,k},$$
we have that
$$ \Delta v_m(x)=\sum_{k=1}^{m}\frac{\Delta Q(e^k x)}{ k}
=\sum_{k=1}^{m}\frac{h(e^k x)\Delta\eta(e^k x)+2\nabla\eta(e^k x)\cdot\nabla h(e^k x)}{k}
.$$
Noticing that~$\eta(e^kx)=0$ for all~$x\in\R^n\setminus B_{e^{-k}}$, we
have that, for all~$x\in B_2$,
\begin{eqnarray*}
\big|h(e^k x)\Delta\eta(e^k x)+2\nabla\eta(e^k x)\cdot\nabla h(e^k x)\big|
\le C\big(\| h\|_{L^\infty(B_{e^{-k}})}+\| \nabla h\|_{L^\infty(B_{e^{-k}})}\big)\le Ce^{-k},
\end{eqnarray*}
for some~$C>0$,
entailing that~$\Delta v_m$ converges uniformly in~$B_2$
to a continuous function~$\widetilde{f}$.

As a result, for every~$\varphi\in C^\infty_0(B_1)$,
\begin{eqnarray*}
\int_{B_1} v(x)\,\Delta \varphi(x)\,dx=
\lim_{m\to+\infty}\int_{B_1} v_m(x)\,\Delta \varphi(x)\,dx=
\lim_{m\to+\infty}\int_{B_1} \Delta v_m(x)\,\varphi(x)\,dx
=\int_{B_1} \widetilde f(x)\,\varphi(x)\,dx.
\end{eqnarray*}
This gives that~$f$ is indeed well-defined and coincides with~$\widetilde f$,
thus completing the proof of~\eqref{CSUdijcINdfyhTUndRgbdEb-1}
and~\eqref{CSUdijcINdfyhTUndRgbdEb-2}.

In addition, since~$\displaystyle\sum_{k=1}^{+\infty}\frac{1}{e^{k}\,k}<+\infty$, we have that
$$ \partial_1 v(x)=\sum_{k=1}^{+\infty}\frac{\partial_1Q(e^k x)}{e^{k}\,k}.$$
In particular, since~$h(0)=0$ and~$\partial_1h(0)=0$, whence~$\partial_1 Q(0)=0$, we have that~$\partial_1v(0)=0$.
As a consequence, if~$\e\in(0,1)$, to be taken as small as we wish in what follows,
we see that
\begin{eqnarray*}&&
\big|\partial_1 v(\e e_1)-\partial_1 v(0)\big|=
\big|\partial_1 v(\e e_1)\big|=\left|
\sum_{k=1}^{+\infty}\frac{\partial_1Q(e^k \e e_1)}{e^{k}\,k}
\right|\\&&\quad\qquad=\left|
\sum_{k=1}^{+\infty}\frac{\partial_1\eta(e^k \e e_1)\,h(e^k \e e_1)+
\partial_1 h(e^k \e e_1)\,\eta(e^k \e e_1)
}{e^{k}\,k}
\right|\\&&\quad\qquad=\left|
\sum_{1\le k\le|\ln\e|}\frac{e^{2k} \e^2\,\partial_1\eta(e^k \e e_1)+2
e^k \e\,\eta(e^k \e e_1)
}{e^{k}\,k}
\right|\\&&\quad\qquad=\e\,\left|
\sum_{1\le k\le|\ln\e|}\frac{e^{k} \e\,\partial_1\eta(e^k \e e_1)+2
\eta(e^k \e e_1)
}{k}
\right|.
\end{eqnarray*}
Now, let~$C_0:=2+\|\eta\|_{C^1(\R^n)}$ and note that
$$ \sum_{|\ln\e|-\ln C_0\le k\le|\ln\e|}\frac{|e^{k} \e\,\partial_1\eta(e^k \e e_1)+2
\eta(e^k \e e_1)|
}{k}\le\sum_{|\ln\e|-\ln C_0\le k\le|\ln\e|}\frac{C_0}{k}\le\frac{C}{|\ln\e|}
,$$
for some~$C>0$.

Also, if~$k\in[1,\,|\ln\e|-\ln C_0)$ then~$e^k\e\le\frac1{C_0}\le\frac12$ and therefore~$\eta(e^k \e e_1)=1$,
leading to
\begin{eqnarray*}
&&e^{k} \e\,\partial_1\eta(e^k \e e_1)+2
\eta(e^k \e e_1)\ge 2-\frac{|\partial_1\eta(e^k \e e_1)|}{C_0}\ge1.
\end{eqnarray*}
Thanks to these observations, we find that
\begin{eqnarray*}&&
\frac{\big|\partial_1 v(\e e_1)-\partial_1 v(0)\big|}\e\ge
\left|
\sum_{1\le k<|\ln\e|-\ln C_0}\frac{e^{k} \e\,\partial_1\eta(e^k \e e_1)+2
\eta(e^k \e e_1)
}{k}
\right|-\frac{C}{|\ln\e|}\\&&\qquad\qquad\ge
\sum_{1\le k<|\ln\e|-\ln C_0}\frac{1
}{k}-\frac{C}{|\ln\e|}.
\end{eqnarray*}
This proves~\eqref{CSUdijcINdfyhTUndRgbdEb-3}: otherwise, if~$v$ belonged to~$ C^2(B_1)$, we would have that
$$\partial_{11}v(0)=\lim_{\e\searrow0}\frac{\big|\partial_1 v(\e e_1)-\partial_1 v(0)\big|}\e
\ge\lim_{\e\searrow0}\sum_{1\le k<|\ln\e|-\ln C_0}\frac{1
}{k}-\frac{C}{|\ln\e|}=\sum_{k=1}^{+\infty}\frac{1
}{k}=+\infty,$$
which is a contradiction.
\end{proof}

In view of Theorem~\ref{CONSDTYROE-Pse3m4pa235sfbscd3jv}, the regularity theory
for elliptic equations (and even for the ``simplest possible'' case given by solutions
of the Poisson equation~$\Delta u=f$) relies first on the detection  of the ``appropriate
functional spaces'' to take into account. The regularity theory that we present here below
is the so-called Schauder theory, named after its inventor\footnote{Being
a Polish mathematician of Jewish origin,
after the invasion of German troops in 1941, Juliusz Schauder
wrote to the German mathematician Ludwig Bieberbach pleading for his support. 
Instead, Bieberbach passed his letter to the Gestapo which arrested and executed Schauder.
The regularity theory of elliptic partial differential \label{JNDtsXVSBNS-KASrSDa3c5h4245e2r2}
equations happens to present a number of tragic events related to World War II,
in which scientific creativity, sense of responsibility, heroism and nobility of spirit
had to confront against pusillanimity, treacheries, greediness and pettiness.
See also footnote~\ref{JNDtsXVSBNS-KASrSDa3c5h4245e2r}
on page~\pageref{JNDtsXVSBNS-KASrSDa3c5h4245e2r}.}
Juliusz Schauder~\cite{MR1545448}).
In a nutshell, the bottom line of this theory is that while Theorem~\ref{CONSDTYROE-Pse3m4pa235sfbscd3jv}
says that if~$\Delta u=f\in C$ does not yield~$u\in C^2$,
this statement becomes essentially correct if one works instead in H\"older spaces,
namely~$\Delta u=f\in C^\alpha$ yields~$u\in C^{2,\alpha}$, provided that~$\alpha\in(0,1)$. Below are
the technical details of this fascinating method.

\section{Potential theory and Schauder estimates for the Laplace operator}\label{KPJMDarJMSiaCLiMMSciAmoenbvD8Sijf-21}

To develop the regularity theory in H\"older spaces for elliptic equations,\index{Schauder estimates|(}
we start by considering the case of equations driven by the Laplace operator.
The main results\footnote{As customary, for every~$k\in\N$
and~$\alpha\in(0,1]$ we will use the norm notation \label{u9oj-NNp-OjtmRM-NosdTHSAMSATInfOdeN-N336OS235yS}
\begin{eqnarray*}
&&[u]_{C^\alpha(\Omega)}:=\sup_{{x,y\in\Omega}\atop{x\ne y}}\frac{|u(x)-u(y)|}{|x-y|^\alpha}\\
{\mbox{and }}&&\|u\|_{C^{k,\alpha}(\Omega)}:=\sum_{{\beta\in\N^n}\atop{|\beta|\le k}}
\|D^\beta u\|_{L^\infty(\Omega)}+\sum_{{\beta\in\N^n}\atop{|\beta|= k}}[D^\beta u]_{C^\alpha(\Omega)}.
\end{eqnarray*}
We also use the notation
$$[u]_{C^{k,\alpha}(\Omega)}:=\sum_{{\beta\in\N^n}\atop{|\beta|= k}}[D^\beta u]_{C^\alpha(\Omega)}.$$
A function is said to belong to~$C^{k,\alpha}(\Omega)$
when~$\|u\|_{C^{k,\alpha}(\Omega)}<+\infty$.
Also, the space~$C^{0,\alpha}(\Omega)$ will be often denoted as~$C^\alpha(\Omega)$ for short.
We stress that, in our notation, these H\"older spaces are meant in a ``global'' \label{LAS0ThgdRujfghlyhn2UPksdssi984n7tE62ps}
sense, and instead we say that~$u\in C^{k,\alpha}_{\rm loc}(\Omega)$
when~$u\in C^{k,\alpha}(\Omega')$ for every open and bounded set~$\Omega'$ such that~$\overline{\Omega'}\subset\Omega$.

See e.g.~\cite[Chapter~3]{MR1406091}, \cite{MR3588287},
\cite{ugjbfvCSnxcAT35Mowjfgoe3AGSHDKS} and the references therein
for a thoroughgoing introduction
to H\"older spaces (keeping in mind that the notation used
is not necessarily uniform across the existing literature).} that we want to address are the following ones, dealing
with\footnote{For completeness, we mention that \label{qupodjh9320jd-29rtjgjjjnn3-x1023N}
the choices of the balls, say~$B_1$ and~$B_{1/2}$, in the regularity results
is somewhat arbitrary and done for the sake of simplicity; one can
state regularity results for instance in balls~$B_R$ and~$B_r$ with~$R>r>0$,
see e.g. the forthcoming Corollary~\ref{qupodjh9320jd-29rtjgjjjnn3-x1023} for future details.} interior estimates (Theorem~\ref{SCHAUDER-INTE})
and local estimates\footnote{For boundary regularity results under various
boundary conditions, see e.g.~\cite{MR601389}.}
at the boundary (Theorem~\ref{SCHAUDER-BOUTH}):

\begin{theorem}\label{SCHAUDER-INTE} Let~$f\in C^\alpha(B_1)$ for some~$\alpha\in(0,1)$. Let~$u\in C^2(B_1)$ be a solution of
\begin{equation}\label{IGSBNHJjusjdvnICtsghdfnEJHACANosdat2R4A9rfX20} \Delta u=f\quad{\mbox{in }}\,B_1.\end{equation}
Then, there exists~$C>0$, depending only on~$n$ and~$\alpha$, such that
$$ \|u\|_{C^{2,\alpha}(B_{1/2})}\le C\,\Big(\|u\|_{L^\infty(B_1)}+\|f\|_{C^\alpha(B_1)}\Big).$$
\end{theorem}

\begin{theorem}\label{SCHAUDER-BOUTH}
Let
\begin{equation}\label{9oj3falphaB10}\begin{split}&
B_r^+:=\{x=(x',x_n)\in B_r {\mbox{ s.t. }}x_n>0\}\\
and\quad&B_r^0:=\{x=(x',x_n)\in B_r {\mbox{ s.t. }}x_n=0\}.\end{split}\end{equation}
Let~$f\in C^\alpha(B_1^+)$ and~$g\in C^{2,\alpha}(B_1^0)$
for some~$\alpha\in(0,1)$. 
Let~$u\in C^2(B_1^+)\cap C(B_1^+\cup B_1^0)$ be a solution of
\begin{equation}\label{poknISOdu-oM-EN0d-f0123-s}\begin{dcases} \Delta u=f\quad{\mbox{in }}\,B_1^+,\\
u=g \quad{\mbox{on }}\,B_1^0.\end{dcases}\end{equation}
Then, there exists~$C>0$, depending only on~$n$ and~$\alpha$, such that
$$ \|u\|_{C^{2,\alpha}(B_{1/2}^+)}\le C\,\Big(\|u\|_{L^\infty(B_1^+)}+\|g\|_{C^{2,\alpha}(B_1^0)}+\|f\|_{C^\alpha(B_1^+)}\Big).$$
\end{theorem}

Results such as the ones in Theorems~\ref{SCHAUDER-INTE} and~\ref{SCHAUDER-BOUTH}
are important, especially considered that in general it is difficult
to express solutions of partial differential equations in a simple and explicit form:
nonetheless, concrete estimates as the ones above are, say, ``almost as good''
as explicit formulas for solutions, in the sense that they are suitable for continuity methods (as
developed in the forthcoming Theorem~\ref{Theorem6.14GT}) and for stable numerical schemes.

The proofs of Theorems~\ref{SCHAUDER-INTE} and~\ref{SCHAUDER-BOUTH}
rely on several deep observations of independent interest.
To address the proof of Theorem~\ref{SCHAUDER-INTE} the strategy is
to try to reduce the problem to the case of harmonic functions
(and this is advantageous because we already know a lot on the regularity theory
of harmonic functions, by exploiting for instance
Cauchy's Estimates in Theorem~\ref{CAUESTIMTH}).

A similar concept for reducing to harmonic functions was exploited
in Corollary~\ref{S-coroEXIS-M023} and it used the idea of subtracting to the given solution
a ``special'' one for which some explicit computations come handy. \label{S-PANEW-TS0jNAIJSkJA9kBBEkS}
Namely, the special solution that one can consider is the Newtonian potential~$-\Gamma*f$ in\index{Newtonian potential} Proposition~\ref{ESI:NICE:FB}.

On the one hand, we have already encountered in the proof of Proposition~\ref{ESI:NICE:FB}
one of the main advantages of working with the Newtonian potential: namely one can utilize
its convolution structure to take derivatives. On the other hand, the proof of Proposition~\ref{ESI:NICE:FB}
also highlighted a conceptual limitation of this approach, due to the singular behavior of
the fundamental solution~$\Gamma$: indeed, derivatives of order~$k$ of~$\Gamma$
behave at the origin like~$|x|^{2-n-k}$, which is integrable only when~$k<2$. This hurdle suggests
that a similar procedure becomes problematic if one is interested (as in Theorem~\ref{SCHAUDER-INTE})
in derivatives of order~$2$: as a matter of fact, a new ingredient
is needed to deal with this significant hindrance, which will be bypassed thanks
to a careful analysis of the singular integral induced by the second derivatives of the fundamental
solution (as we will see, this new analysis will leverage the symmetry
properties of the singular integral kernel associated with second derivatives
in order to detect suitable cancellations).\medskip

Here are the technical details to make this strategy work and prove Theorem~\ref{SCHAUDER-INTE}.
First off, we explicitly express second derivatives
of the Newtonian potential\index{Newtonian potential} in terms of a singular integral kernel (this is
not obvious, since, due to the discussion above, we cannot just place two
derivatives on the fundamental solution, due to the lack of local integrability,
and indeed the final result makes it also appear an additional Dirac Delta Function):

\begin{proposition}\label{SI:IN:SE:DR}
Let~$\Gamma$ be the fundamental solution in~\eqref{GAMMAFU} and~$f\in C^\alpha_0(\R^n)$,
for some~$\alpha\in(0,1]$. Set~$v:=-\Gamma*f$. Then, for all~$i$, $j\in\{1,\dots,n\}$,
\begin{equation}\label{POTNDEWSECDERIVA} \partial_{ij}v(x)=\int_{\R^n}\partial_{ij}\Gamma(x-y)\big(f(x)-f(y)\big)\,dy+\frac{\delta_{ij}}n\,f(x).\end{equation}
\end{proposition}

Here and in what follows, improper and singular integrals are taken in the principal value sense
by averaging outside the singularity: for instance
\begin{equation}\label{12456y34-2rtgkK123343524Sqew53M-THD-qdofkCVNXM}
\int_{\R^n}\partial_{ij}\Gamma(x-y)\big(f(x)-f(y)\big)\,dy:=\lim_{R\to+\infty}
\int_{B_R(x)}\partial_{ij}\Gamma(x-y)\big(f(x)-f(y)\big)\,dy.\end{equation}
Notice indeed that~$\partial_{ij}\Gamma$ is not integrable, hence the integral
on the left hand side of the previous equation is not a standard Lebesgue integral
and instead requires the limit definition introduced on the right hand side.
Also, of course, no confusion should arise in equation~\eqref{POTNDEWSECDERIVA}
between the second order partial derivative~$\partial_{ij}$
and the Kronecker notation~$\delta_{ij}$
(recall the notation introduced in footnote~\ref{Kronecker-notation}
on page~\pageref{Kronecker-notation}). Interestingly, it follows from the harmonicity of~$\Gamma$ outside
its singular point and the assumption that~$f\in C^\alpha_0(\R^n)$ that
\begin{eqnarray*}&&
\left|\int_{\R^n}\Delta\Gamma(x-y)\big(f(x)-f(y)\big)\,dy\right|=
\lim_{r\searrow0}\left|\int_{B_r(x)}\Delta\Gamma(x-y)\big(f(x)-f(y)\big)\,dy\right|\\&&\qquad\qquad\le
\lim_{r\searrow0}C\left|\int_{B_r(x)}|x-y|^{\alpha-n}\,dy\right|=\lim_{r\searrow0}Cr^\alpha=0,
\end{eqnarray*}
therefore~\eqref{POTNDEWSECDERIVA} also entails that
\begin{equation}\label{SOLPONE-0}
{\mbox{if~$f\in C^\alpha_0(\R^n)$ then~$-\Delta(\Gamma*f)=f$.}}\end{equation}
This observation shows that the result in Proposition~\ref{ESI:NICE:FB} can be sharpened
by replacing the assumption~$f\in C^1_0(\R^n)$ with~$f\in C^\alpha_0(\R^n)$ for some~$\alpha\in(0,1)$.
Correspondingly, 
\begin{equation}\label{SqerLPONE-0-9ieyhdfweohgiowegb74uthghhs1334756ujhgf}\begin{split}&
{\mbox{the assumption~$f\in C^1(\Omega)$ in Corollary~\ref{S-coroEXIS-M023}}}\\&{\mbox{can be relaxed to~$f\in C^\alpha(\Omega)$ for some~$\alpha\in(0,1)$.}}\end{split}\end{equation}
When needed, the observation in~\eqref{SOLPONE-0} can be even sharpened, since it actually holds that
\begin{equation}\label{SOLPONE-0-LOCALE}
{\mbox{if~$x_0\in\R^n$, $r>0$ and~$f\in C^\alpha(B_r(x_0))$ then~$-\Delta(\Gamma*f)(x_0)=f(x_0)$.}}\end{equation}
To prove this,
it is enough to consider a cutoff function~$\varphi\in C^\infty_0(B_{3r/4}(x_0))$
with~$\varphi=1$ in~$B_{r/2}(x_0)$,
define~$f_1:=f\varphi$ and~$f_2:=f(1-\varphi)$ and notice that~$f_1\in C^\alpha_0(\R^n)$. 
We can also extend~$f_2$ as a bounded function on the whole of~$\R^n$, say, by setting~$f_2:=0$
outside~$B_r(x_0)$.
Accordingly, by~\eqref{SOLPONE-0}
and the fact that~$f_2=0$ in~$B_{r/2}(x_0)$, for all~$x\in B_{r/4}(x_0)$ we have that
\begin{equation}\label{2uj4NSIC5ChjdEjmdvLL99PiPMMSN09f0MsdUnsdCCkfE03}\begin{split}&
-\Delta(\Gamma*f)(x)=-\Delta(\Gamma*f_1)(x)-\Delta(\Gamma*f_2)(x)
=f_1(x)-\Delta\int_{\R^n}\Gamma(x-y)\,f_2(y)\,dy\\&\qquad\qquad\qquad\qquad\qquad
=f(x)-\Delta\int_{\R^n\setminus B_{r/2}(x_0)}\Gamma(x-y)\,f_2(y)\,dy.
\end{split}\end{equation}
Now, if~$x\in B_{r/4}(x_0)$ and~$y\in\R^n\setminus B_{r/2}(x_0)$ then~$|x-y|\ge|y|-|x|\ge\frac{r}4$.
Consequently, by~\eqref{2uj4NSIC5ChjdEjmdvLL99PiPMMSN09f0MsdUnsdCCkfE03} and
the harmonicity of~$\Gamma$ away from its singularity, for all~$x\in B_{r/4}(x_0)$ we have that
$$ -\Delta(\Gamma*f)(x)=f(x)-\int_{\R^n\setminus B_{r/2}(x_0)}\Delta\Gamma(x-y)\,f_2(y)\,dy
=f(x),$$
from which one obtains~\eqref{SOLPONE-0-LOCALE}.

\begin{proof}[Proof of Proposition~\ref{SI:IN:SE:DR}] This argument is a refinement
of that produced in the proof of~\eqref{PRIMA98DER}: namely, given~$\rho>0$, we use the regularization~$\Gamma_\rho$ of the fundamental solution introduced in~\eqref{COIENEE-0986ytufgkbv-0-rjfeonvnb2GDB} and~\eqref{GAROGRA}
and use a limit argument. The technical details go as follows.
Let~$V:=\partial_i v=-\partial_i(\Gamma*f)$ and~$V_\rho:=-\partial_i(\Gamma_\rho*f)$.
We know from~\eqref{PRIMA98DER} that~$V=-(\partial_i\Gamma)*f$.

Moreover, for all~$x\in\R^n$,
\begin{equation}\label{6787iudsfeVAskdncfgFwasfdgfFASBdcNidkHKihfvbnLI3Sd0}
\begin{split}&
|V_\rho(x)-V(x)|=|(\partial_i\Gamma)*f(x)-(\partial_i\Gamma_\rho)*f(x)|\\&\qquad\le
\frac{C}{\rho^n}\int_{B_\rho} |y|\,|f(x-y)|\,dy+C\int_{B_\rho}|y|^{1-n}\,|f(x-y)|\,dy
\le C\,\|f\|_{L^\infty(\R^n)}\,\rho,\end{split}
\end{equation}
for some constant~$C>0$ depending only on~$n$.

Additionally, for all~$j\in\{1,\dots,n\}$,
\begin{equation}\label{NS-0iJMS-92gfiuwgflvsuy3frvhchfDAhvsdbMiwoion9hgnnO8ik}\begin{split}&
\left| \partial_j V_\rho(x)-\int_{\R^n}\partial_{ij}\Gamma(x-y)\big(f(x)-f(y)\big)\,dy-\frac{\delta_{ij}}n\,f(x)
\right|\\=\;&
\left| (\partial_{ij} \Gamma_\rho)*f(x)+\int_{\R^n}\partial_{ij}\Gamma(x-y)\big(f(x)-f(y)\big)\,dy+\frac{\delta_{ij}}n\,f(x)
\right|\\
=\;&\left|
\int_{\R^n}\partial_{ij}\Gamma_\rho(x-y)\,f(y)\,dy+
\int_{\R^n}\partial_{ij}\Gamma(x-y)\big(f(x)-f(y)\big)\,dy+\frac{\delta_{ij}}n\,f(x)\right|\\
=\;&\left|
\int_{\R^n}\partial_{ij}\Gamma_\rho(x-y)\,f(x)\,dy
-\int_{\R^n}\partial_{ij}\Gamma_\rho(x-y)\,\big(f(x)-f(y)\big)\,dy\right.\\&\qquad\left.\qquad
+\int_{\R^n}\partial_{ij}\Gamma(x-y)\big(f(x)-f(y)\big)\,dy+\frac{\delta_{ij}}n\,f(x)\right|\\
=\;&\left|\int_{\R^n}\partial_{ij}\Gamma_\rho(x-y)\,f(x)\,dy
\right.\\&\qquad\left.\qquad
+\int_{B_\rho(x)}\big(\partial_{ij}\Gamma(x-y)-\partial_{ij}\Gamma_\rho(x-y)\big)\big(f(x)-f(y)\big)\,dy+\frac{\delta_{ij}}n\,f(x)
\right|\\
\leq\;&\left|\int_{\R^n}\partial_{ij}\Gamma_\rho(x-y)\,f(x)\,dy
+\frac{\delta_{ij}}n\,f(x)
\right|+C\,\|f\|_{C^\alpha(\R^n)}\int_{B_\rho(x)}|x-y|^{\alpha-n}\,dy\\
\leq\;&\left|-\int_{B_\rho(x)}
\frac{\delta_{ij} f(x)}{n |B_\rho|}\,\,dy+\int_{\R^n\setminus B_\rho(x)}
\partial_{ij}\Gamma (x-y)\,f(x)\,\,dy
+\frac{\delta_{ij}}n\,f(x)
\right|+C\,\|f\|_{C^\alpha(\R^n)}\,\rho^\alpha\\
=\;&\left| \int_{\R^n\setminus B_\rho(x)}
\partial_{ij}\Gamma (x-y)\,f(x)\,\,dy
\right|+C\,\|f\|_{C^\alpha(\R^n)}\,\rho^\alpha,
\end{split}\end{equation}
up to renaming~$C>0$ from line to line.

We also let
$$ c'_n:=\begin{dcases}
2c_n & {\mbox{ if }}n=2,\\
c_n\,n(n-2)& {\mbox{ if }}n\ne2,
\end{dcases}$$
where~$c_n$ is the constant in~\eqref{COIENEE},
and we
point out that if~$i\ne j$ and~$R>\rho$ then
\begin{eqnarray*}&&
\int_{B_R(x)\setminus B_\rho(x)}\partial_{ij}\Gamma(x-y) \,dy=
c'_n\int_{B_R(x)\setminus B_\rho(x)}\frac{(x_i-y_i)(x_j-y_j)}{|x-y|^{n+2}}\,dy=0,
\end{eqnarray*}
due to odd symmetry.
Moreover,
$$ \int_{B_R(x)\setminus B_\rho(x)}\partial_{ii}\Gamma(x-y) \,dy=
\frac1n\int_{B_R(x)\setminus B_\rho(x)}
\Delta\Gamma(x-y) \,dy=0,$$
thanks to the harmonicity of~$\Gamma$ outside its singularity.

It follows from these observations that, for all~$i$, $j\in\{1,\dots,n\}$,
\begin{equation*} \int_{B_R(x)\setminus B_\rho(x)}\partial_{ij}\Gamma(x-y) \,dy=0\end{equation*}
and thus, in the notation of~\eqref{12456y34-2rtgkK123343524Sqew53M-THD-qdofkCVNXM},
\begin{equation*} \int_{\R^n\setminus B_\rho(x)}\partial_{ij}\Gamma(x-y) \,dy=0.\end{equation*}
{F}rom this
and~\eqref{NS-0iJMS-92gfiuwgflvsuy3frvhchfDAhvsdbMiwoion9hgnnO8ik} we obtain that
\begin{equation*}\begin{split}&
\left| \partial_j V_\rho(x)-\int_{\R^n}\partial_{ij}\Gamma(x-y)\big(f(x)-f(y)\big)\,dy-\frac{\delta_{ij}}n\,f(x)
\right|\le C\,\|f\|_{C^\alpha(\R^n)}\,\rho^\alpha.
\end{split}\end{equation*}
It follows from this estimate and~\eqref{6787iudsfeVAskdncfgFwasfdgfFASBdcNidkHKihfvbnLI3Sd0} that~$V_\rho$
converges uniformly to~$V$ and~$\partial_j V_\rho$ converges uniformly to the right hand side of~\eqref{POTNDEWSECDERIVA}. That is, recalling~\eqref{IHFSHNTHSYDHMSWMDTKFDOFLL},
the function~$-\Gamma_\rho*f$
converges uniformly together with its derivatives up to the second order
(the function converging to~$-\Gamma*f$ and the second derivatives
converging to the right hand side of~\eqref{POTNDEWSECDERIVA})
and this establishes~\eqref{POTNDEWSECDERIVA}, as desired.
\end{proof}

In view of Proposition~\ref{SI:IN:SE:DR}, it is now convenient to study
the property of a singular integral kernel~$K:\R^n\setminus\{0\}\to\R$ satisfying the following properties:
\begin{eqnarray}
\label{KERNE-SINGO-001}&&{\mbox{there exists~$g\in L^\infty(\partial B_1)$ such that for every~$x\in \R^n\setminus\{0\}$}}\nonumber\\
&&{\mbox{we have that }} 
K(x)=\frac{1}{|x|^n}\,g\left(\frac{x}{|x|}\right), \;{\mbox{and}}\\
%%%%%% &&\label{KERNE-SINGO-002}\int_{\partial B_1} g(\omega)\,d{\mathcal{H}}_\omega^{n-1}=0,\\
\label{KERNE-SINGO-003}&& {\mbox{there exists~$C>0$ such that for every~$x\in \R^n\setminus\{0\}$ we have that }}
|\nabla K(x)|\le \frac{C}{|x|^{n+1}}.
\end{eqnarray}
For short, given an open set~${\mathcal{V}}\subseteq\R^n$ and a function~$\Phi\in L^1_{\rm loc}(\R^n\times\R^n)$
we will also use the principal value notation
\begin{equation}\label{KS-TSGDBfojeISLmadne-ojdfnLS-KMDMDND9uehn}
\int_{{\mathcal{V}}}K(x-y)\,\Phi(x,y)\,dy:=
\lim_{{R\to+\infty}\atop{\e\searrow0}}\int_{({\mathcal{V}}\setminus B_\e(x))\cap B_R(x)}K(x-y)\,\Phi(x,y)\,dy,
\end{equation}
whenever the above limits exist.

We will also set
\begin{equation}\label{4.1.14BIS}
Tf(x):=\int_{{\mathcal{V}}}K(x-y)\,\big(f(x)-f(y)\big)\,dy.
\end{equation}

The main step to prove the interior regularity result in Theorem~\ref{SCHAUDER-INTE} and the local result at the boundary in Theorem~\ref{SCHAUDER-BOUTH} is now the following
singular integral regularity estimate:

\begin{lemma}\label{NKSD5678yihnqehr837tCOvnlIDJMcINBSno3Na7syshdrr9ISDcvaHF23}
Assume that~$K$ satisfies~\eqref{KERNE-SINGO-001} and~\eqref{KERNE-SINGO-003}.
Suppose also that, for each~$z_0\in\R^n$ and~$\rho>0$,
\begin{equation}\label{ABBE9}
\left|
\int_{{\mathcal{V}}\setminus B_{\rho}(z_0)}K(z_0-y)\,dy\right|\le
C_0,
\end{equation}
for a suitable~$C_0>0$.

If~$f\in C^\alpha({\mathcal{V}})$ for some~$\alpha\in(0,1)$, then
$$ [Tf]_{C^\alpha({\mathcal{V}})}\le C\,[f]_{C^\alpha({\mathcal{V}})},$$
for a suitable constant~$C>0$ depending only on~$n$, $\alpha$, $K$ and~$C_0$.
\end{lemma}

\begin{proof} Let~$z_0$, $z_1\in{\mathcal{V}}$. Let~$\delta:=|z_0-z_1|$ and~$z_2:=\frac{z_0+z_1}2$.
The gist of the following estimate is that in regions of order~$\delta$ around~$z_2$
all contributions are of order~$\delta^\alpha$ and for the complement of these regions the computation becomes more delicate
and further cancellations must be taken into account. More precisely, we have
\begin{equation}\label{7.02irj-517}
\begin{split}&
|Tf(z_0)-Tf(z_1)|\\  =\,&\left|
\int_{{\mathcal{V}}}K(z_0-y)\big(f(z_0)-f(y)\big)\,dy-\int_{{\mathcal{V}}}K(z_1-y)\big(f(z_1)-f(y)\big)\,dy
\right|\\
\le\,&\left|
\int_{{\mathcal{V}}\setminus B_{2\delta}(z_2)}K(z_0-y)\big(f(z_0)-f(y)\big)\,dy-\int_{{\mathcal{V}}\setminus B_{2\delta}(z_2)}K(z_1-y)\big(f(z_1)-f(y)\big)\,dy
\right|\\&\qquad +
\left|
\int_{{\mathcal{V}}\cap B_{2\delta}(z_2)}K(z_0-y)\big(f(z_0)-f(y)\big)\,dy\right|+\left|\int_{{\mathcal{V}}\cap B_{2\delta}(z_2)}K(z_1-y)\big(f(z_1)-f(y)\big)\,dy
\right|.\end{split}
\end{equation}
Now, for~$j\in\{0,1\}$, if~$y\in {\mathcal{V}}\cap B_{2\delta}(z_2)$ then~$|z_j-y|\le|z_j-z_2|+|z_2-y|<3\delta$, whence
\begin{equation}\label{7.02irj-518}
\begin{split}&
\left|
\int_{{\mathcal{V}}\cap B_{2\delta}(z_2)}K(z_j-y)\big(f(z_j)-f(y)\big)\,dy\right|\le
[f]_{C^\alpha({\mathcal{V}})}\int_{B_{2\delta}(z_2)}|K(z_j-y)|\,|z_j-y|^\alpha\,dy\\&\qquad\le
[f]_{C^\alpha({\mathcal{V}})}\|g\|_{L^\infty(\partial B_1)}\int_{B_{2\delta}(z_2)}|z_j-y|^{\alpha-n}\,dy\le
[f]_{C^\alpha({\mathcal{V}})}\|g\|_{L^\infty(\partial B_1)}\int_{B_{3\delta}}|x|^{\alpha-n}\,dx\\&\qquad\le C\,[f]_{C^\alpha({\mathcal{V}})}\|g\|_{L^\infty(\partial B_1)}\,\delta^\alpha,
\end{split}\end{equation}
for some~$C>0$,
thanks to~\eqref{KERNE-SINGO-001}.

Furthermore, making use of~\eqref{KERNE-SINGO-003}, if~$y\in {\mathcal{V}}\setminus B_{2\delta}(z_2)$ and~$t\in[0,1]$ then
$$ |z_1+t(z_0-z_1)-y|\ge |y-z_2|-|z_2-z_1|-|z_0-z_1|\ge
|y-z_2|-\frac{3\delta}2\ge\frac{3|y-z_2|}4,$$
from which, up to renaming~$C$, we deduce that
\begin{eqnarray*}&&
|K(z_0-y)-K(z_1-y)|\le \int_0^1 |\nabla K(z_1+t(z_0-z_1)-y)|\,dt\,|z_0-z_1|\\
&&\qquad\le\sup_{t\in[0,1]}\frac{C\,|z_0-z_1|}{|z_1+t(z_0-z_1)-y|^{n+1}}\le
\frac{C\,\delta}{|y-z_2|^{n+1}}.
\end{eqnarray*}
For this reason,
\begin{equation}\label{418-OIhgcufs0-2674-9023941-94u8-01} \begin{split}\Xi\,:=\,&\left|
\int_{{\mathcal{V}}\setminus B_{2\delta}(z_2)}K(z_0-y)\big(f(z_0)-f(y)\big)\,dy-\int_{{\mathcal{V}}\setminus B_{2\delta}(z_2)}K(z_1-y)\big(f(z_1)-f(y)\big)\,dy\right|
\\ \le\,&\left|
\int_{{\mathcal{V}}\setminus B_{2\delta}(z_2)}K(z_0-y)\big(f(z_0)-f(y)\big)\,dy-\int_{{\mathcal{V}}\setminus B_{2\delta}(z_2)}K(z_0-y)\big(f(z_1)-f(y)\big)\,dy\right|\\&\qquad+
\left|
\int_{{\mathcal{V}}\setminus B_{2\delta}(z_2)}K(z_0-y)\big(f(z_1)-f(y)\big)\,dy-\int_{{\mathcal{V}}\setminus B_{2\delta}(z_2)}K(z_1-y)\big(f(z_1)-f(y)\big)\,dy\right|\\ \le\,&\left|
\int_{{\mathcal{V}}\setminus B_{2\delta}(z_2)}K(z_0-y)\big(f(z_0)-f(z_1)\big)\,dy
\right| +C\,\delta
\int_{{\mathcal{V}}\setminus B_{2\delta}(z_2)}\frac{\big|f(z_1)-f(y)\big|}{|y-z_2|^{n+1}}\,dy\\ \le\,&
\big|f(z_0)-f(z_1)\big|\,\left(
\left|
\int_{B_{8\delta}(z_0)\setminus B_{2\delta}(z_2)}K(z_0-y)\,dy
\right| +\left|
\int_{{\mathcal{V}}\setminus B_{8\delta}(z_0)}K(z_0-y)\,dy
\right|\right)\\&\qquad+C\,[f]_{C^\alpha({\mathcal{V}})}\,\delta
\int_{\R^n\setminus B_{2\delta}(z_2)}\frac{|z_1-y|^\alpha}{|y-z_2|^{n+1}}\,dy.
\end{split}\end{equation}
We also point out that, in view of~\eqref{KERNE-SINGO-001},
\begin{equation}\label{418-OIhgcufs0-2674-9023941-94u8-02}\int_{B_{8\delta}(z_0)\setminus B_{2\delta}(z_2)}|K(z_0-y)|\,dy\le
\|g\|_{L^\infty(\partial B_1)}
\int_{B_{8\delta}(z_0)\setminus B_{\delta/8}(z_0)}\frac{dy}{|z_0-y|^n}\le C\,\|g\|_{L^\infty(\partial B_1)}.
\end{equation}
Recalling~\eqref{ABBE9}, from the observations in~\eqref{418-OIhgcufs0-2674-9023941-94u8-01} and~\eqref{418-OIhgcufs0-2674-9023941-94u8-02}, we arrive at
\begin{eqnarray*}
\Xi&\le&\big(C\,\|g\|_{L^\infty(\partial B_1)}+C_0\big)
\big|f(z_0)-f(z_1)\big|+C\,[f]_{C^\alpha({\mathcal{V}})}\,\delta
\int_{\R^n\setminus B_{2\delta}(z_2)}\frac{|z_1-y|^\alpha}{|y-z_2|^{n+1}}\,dy\\&\le&
\big(C\|g\|_{L^\infty(\partial B_1)}+C_0\big)\,[f]_{C^\alpha({\mathcal{V}})}
\delta^\alpha+C\,[f]_{C^\alpha({\mathcal{V}})}\,\delta
\int_{\R^n\setminus B_{2\delta}(z_2)}\frac{(|z_1-z_2|+|y-z_2|)^\alpha}{|y-z_2|^{n+1}}\,dy\\
&\le&\big(C\|g\|_{L^\infty(\partial B_1)}+C_0\big)\,[f]_{C^\alpha({\mathcal{V}})}\,
\delta^\alpha+C\,[f]_{C^\alpha({\mathcal{V}})}\,\delta
\int_{\R^n\setminus B_{2\delta}(z_2)}|y-z_2|^{\alpha-n-1}\,dy\\
&\le&\big(C\|g\|_{L^\infty(\partial B_1)}+C_0\big)\,[f]_{C^\alpha({\mathcal{V}})}\,
\delta^\alpha+C\,[f]_{C^\alpha({\mathcal{V}})}\,\delta^\alpha.
\end{eqnarray*}
Combined with~\eqref{7.02irj-517} and~\eqref{7.02irj-518} this yields that
$$ |Tf(z_0)-Tf(z_1)|\le\Xi+C\,[f]_{C^\alpha({\mathcal{V}})}\|g\|_{L^\infty(\partial B_1)}\,\delta^\alpha\le
C\,[f]_{C^\alpha({\mathcal{V}})}\big(\|g\|_{L^\infty(\partial B_1)}+1\big)\,\delta^\alpha,$$
providing the desired result.
\end{proof}

With this, we can prove the following estimate:

\begin{corollary}\label{QUA4}
Let~$f\in C^\alpha_0(\R^n)$ for some~$\alpha\in(0,1)$ and~$v:=-\Gamma*{f}$.

Then, for all~$i$, $j\in\{1,\dots,n\}$,
\begin{equation*}
[\partial_{ij}v]_{C^\alpha(\R^n)}\le
C\,[f]_{C^\alpha(\R^n)}
,\end{equation*}
for some positive constant~$C$ depending only on~$n$ and~$\alpha$.
\end{corollary}

\begin{proof} By Proposition~\ref{SI:IN:SE:DR},
\begin{equation*}
\partial_{ij}v(x)=\int_{\R^n}K(x-y)\big(f(x)-f(y)\big)\,dy+\frac{\delta_{ij}}n\,f(x),\end{equation*}
with
\begin{equation}\label{AV1B}
K(x):=\partial_{ij}\Gamma(x).\end{equation}
That is, in the notation of~\eqref{4.1.14BIS},
\begin{equation}\label{7wuyfh0KSM-owjfh375yth-0okdmHSND-1-2rijfnbvij38gh94uyrgbv-294k5} \partial_{ij}v(x)=T{f}(x)+\frac{\delta_{ij}}n\, f(x)
.\end{equation}
In light of~\eqref{GAMMAFU}, we know that~$K(x)$ agrees, up to a dimensional constant, with~$
\frac{x_i x_j}{|x|^{n+2}}$ if~$i\ne j$, and with~$\frac{|x|^2-nx_i^2}{|x|^{n+2}}$ if~$i=j$.
Consequently assumption~\eqref{KERNE-SINGO-001}
is satisfied with
\begin{equation}\label{9UOHNS-01euowfhvbogdblwiqfgyoerufg8veduifgedirIPSEtbcTIZnowevd99e}
\partial B_1\ni\omega\mapsto g(\omega):=\begin{dcases}
\omega_i\omega_j & {\mbox{ if }}i\ne j,\\
1-n\omega_i^2& {\mbox{ if }}i= j.
\end{dcases}\end{equation}
Also, assumption~\eqref{KERNE-SINGO-003} holds true.

As for~\eqref{ABBE9}, we obtain it as a byproduct of a general strategy based
on cancellations. Namely, in the notation of~\eqref{KERNE-SINGO-001} and~\eqref{AV1B},
we have that
\begin{equation}\label{KERNE-SINGO-002}
\int_{\partial B_1} g(\omega)\,d{\mathcal{H}}_\omega^{n-1}=0.\end{equation}
As a matter of fact, we see
from~\eqref{9UOHNS-01euowfhvbogdblwiqfgyoerufg8veduifgedirIPSEtbcTIZnowevd99e} 
that~\eqref{KERNE-SINGO-002} holds true
when~$i\ne j$ due to an odd symmetry argument, as well as when~$i=j$ since in this case
$$ \fint_{\partial B_1}(1-n\omega_i^2)\,d{\mathcal{H}}_\omega^{n-1}=
1-n\fint_{\partial B_1} \omega_i^2\,d{\mathcal{H}}_\omega^{n-1}=
1-\sum_{k=1}^n\fint_{\partial B_1} \omega_k^2\,d{\mathcal{H}}_\omega^{n-1}=
1- \fint_{\partial B_1}|\omega|^2\,d{\mathcal{H}}_\omega^{n-1}=0.$$
Having established~\eqref{KERNE-SINGO-002}, we exploit it to check that, for each~$z_0\in\R^n$ and~$\rho>0$,
\begin{equation}\label{418-OIhgcufs0-2674-9023941-94u8-03}\int_{\R^n\setminus B_{\rho}(z_0)}K(z_0-y)\,dy=0.\end{equation}
Indeed:
\begin{eqnarray*}&&\int_{\R^n\setminus B_{\rho}(z_0)}K(z_0-y)\,dy=
\int_{\R^n\setminus B_{\rho}(z_0)}g\left( \frac{z_0-y}{|z_0-y|}\right)\frac{dy}{|z_0-y|^n}
=\iint_{[\rho,+\infty)\times(\partial B_1)}g(\omega)\frac{dr\,d{\mathcal{H}}^{n-1}_\omega}{r}=0,\end{eqnarray*}
which gives~\eqref{418-OIhgcufs0-2674-9023941-94u8-03}.

In turn, \eqref{418-OIhgcufs0-2674-9023941-94u8-03} 
entails~\eqref{ABBE9}.
We can thereby exploit Lemma~\ref{NKSD5678yihnqehr837tCOvnlIDJMcINBSno3Na7syshdrr9ISDcvaHF23} and find that
$$ [T f]_{C^\alpha(\R^n)}\le C\,[ f]_{C^\alpha(\R^n)}.$$
{F}rom this estimate and~\eqref{7wuyfh0KSM-owjfh375yth-0okdmHSND-1-2rijfnbvij38gh94uyrgbv-294k5} the desired result follows.\end{proof}

In order to prove Theorem~\ref{SCHAUDER-INTE}, we will also need the following~$L^\infty$-estimate:

\begin{lemma}\label{lemmaspe0}
Let~$f\in C^\alpha_0(\R^n)$ for some~$\alpha\in(0,1)$ and~$v:=-\Gamma*{f}$.

Then, for all~$i$, $j\in\{1,\dots,n\}$,
\begin{equation*}
\|\partial_{ij}v\|_{L^\infty(\R^n)}\le
C\,\|f\|_{C^\alpha(\R^n)}
,\end{equation*}
for some positive constant~$C$ depending only on~$n$, $\alpha$ and the support of~$f$.
\end{lemma}

\begin{proof}
By Proposition~\ref{SI:IN:SE:DR},
\begin{equation}\label{d4t4ut4yu485yu58}
\partial_{ij}v(x)=\int_{\R^n}\partial_{ij}\Gamma(x-y)\big(f(x)-f(y)\big)\,dy+\frac{\delta_{ij}}n\,f(x).\end{equation}
We suppose that the support of~$f$ is contained in~$B_R$ for some~$R>0$, and therefore we have that
there exists~$x_0\in B_{R+1}$ such that~$f(x_0)=0$.

We also set
$$\phi(x):=\int_{\R^n}\partial_{ij}\Gamma(x-y)\big(f(x)-f(y)\big)\,dy$$
and we see that
\begin{eqnarray*} &&
|\phi(x_0)|\le C\int_{\R^n}\frac{\big|f(y)\big|}{|x_0-y|^{n}}\,dy=
C\int_{B_R}\frac{\big|f(y)\big|}{|x_0-y|^{n}}\,dy=
C\int_{B_R}\frac{\big|f(x_0)-f(y)\big|}{|x_0-y|^{n}}\,dy\\&&\qquad\le
C[f]_{C^\alpha(\R^n)}\int_{B_R}|x_0-y|^{\alpha-n}\,dy\le C[f]_{C^\alpha(\R^n)},
\end{eqnarray*}
up to renaming~$C>0$.
As a consequence, for every~$x\in B_{R+2}$,
\begin{eqnarray*}
|\phi(x)|\le |\phi(x_0)|+|\phi(x)-\phi(x_0)|\le  C[f]_{C^\alpha(\R^n)}
+[\phi]_{C^\alpha(B_{R+2})}|x-x_0|^\alpha \le C[f]_{C^\alpha(\R^n)},
\end{eqnarray*}
where we have used Lemma~\ref{NKSD5678yihnqehr837tCOvnlIDJMcINBSno3Na7syshdrr9ISDcvaHF23}
(recall~\eqref{AV1B} and~\eqref{418-OIhgcufs0-2674-9023941-94u8-03} to conclude that~\eqref{ABBE9}
is satisfied).

Furthermore, we notice that if~$x\in \R^n\setminus B_{R+2}$ and~$y\in B_R$, then~$|x-y|\ge|x|-|y|\ge R+2-R=2$, and thus
\begin{eqnarray*}
|\phi(x)|\le C\int_{B_R}\frac{\big|f(y)\big|}{|x-y|^{n}}\,dy\le C\|f\|_{L^\infty(\R^n)}\int_{B_R} \frac{dy}{2^{n}}\le
C\|f\|_{L^\infty(\R^n)}.
\end{eqnarray*}
These considerations show that
$$\|\phi\|_{L^\infty(\R^n)}\le C\|f\|_{C^\alpha(\R^n)}.$$
{F}rom this and~\eqref{d4t4ut4yu485yu58} we obtain the desired estimate.
\end{proof}

With this preliminary work, we can now address\footnote{A conceptually different proof of
Theorem~\ref{SCHAUDER-INTE}, not directly relying on singular integral
estimates for the Newtonian potential, will be presented on page~\pageref{023iurjw94jfJAMAnJHOSo3o24f32tofXZwf}.}
the proof of the desired interior regularity result:

\begin{proof}[Proof of Theorem~\ref{SCHAUDER-INTE}] Let~$\varphi\in C^\infty_0(B_{8/9},\,[0,1])$,
with~$\varphi=1$ in~$B_{3/4}$ and~$|\nabla \varphi|\le 10$. Let also~$\widetilde f:=f\varphi$.
In this way, we have that~$\widetilde f\in C^\alpha_0(B_1)$. We define
\begin{equation}\label{IGSBNHJjusjdvnICtsghdfnEJHACANosdat2R4A9rfX2-9}
v:=-\Gamma*\widetilde{f}\end{equation}
and we deduce from~\eqref{SOLPONE-0} that
\begin{equation}\label{IGSBNHJjusjdvnICtsghdfnEJHACANosdat2R4A9rfX2}{\mbox{$
\Delta v=\widetilde{f}$ in~$\R^n$.}}\end{equation}
Also, by Corollary~\ref{QUA4},
for all~$i$, $j\in\{1,\dots,n\}$,
\begin{equation} \label{LSDNlwnfeECAocOSndfeDEFSGBD-UJNSDfgkrhVEUjfIMSdfdp01}
[\partial_{ij}v]_{C^\alpha(\R^n)}\le C\,[\widetilde f]_{C^\alpha(\R^n)}\le
C\,\|f\|_{C^\alpha(B_1)}
,\end{equation}
up to renaming~$C$.

Moreover, for every~$x\in B_1$,
$$ |v(x)|\le\|\widetilde f\|_{L^\infty(\R^n)}\int_{B_1}\Gamma(x-y)\,dy\le
\|f\|_{L^\infty(B_1)}\int_{B_2}\Gamma(z)\,dz\le C\,\|f\|_{L^\infty(B_1)}$$
and
$$ |\nabla v(x)|\le
\int_{B_1} |\nabla\Gamma(x-y)|\,|\widetilde f(y)|\,dy\le
\|f\|_{L^\infty(B_1)}\int_{B_2}|\nabla\Gamma(z)|\,dz\le C\,\|f\|_{L^\infty(B_1)},
$$
thanks to~\eqref{PRIMA98DER} and~\eqref{IGSBNHJjusjdvnICtsghdfnEJHACANosdat2R4A9rfX2-9}.

Also, by Lemma~\ref{lemmaspe0}, for all~$i$, $j\in\{1,\dots,n\}$,
\begin{eqnarray*}| \partial_{ij}v(x)|
\le C\|\widetilde f\|_{C^\alpha(\R^n)}
\le C\|f\|_{C^\alpha(B_1)}.
\end{eqnarray*}
Accordingly, we deduce from~\eqref{LSDNlwnfeECAocOSndfeDEFSGBD-UJNSDfgkrhVEUjfIMSdfdp01} that
\begin{equation} \label{LSDNlwnfeECAocOSndfeDEFSGBD-UJNSDfgkrhVEUjfIMSdfdp01-2-09uyh-0ks89}
\|v\|_{C^{2,\alpha}(B_1)}\le C\,\|f\|_{C^\alpha(B_1)}
,\end{equation}
up to renaming~$C$ once again.

We now define~$w:=u-v$ and we exploit~\eqref{IGSBNHJjusjdvnICtsghdfnEJHACANosdat2R4A9rfX20}
and~\eqref{IGSBNHJjusjdvnICtsghdfnEJHACANosdat2R4A9rfX2} to conclude that~$\Delta w=0$ in~$B_{3/4}$.
Thus, it follows from Cauchy's Estimates (see Theorem~\ref{CAUESTIMTH}) that
$$ \|w\|_{C^3(B_{1/2})}\le C\,\|w\|_{L^\infty(B_{3/4})}\le
C\,\Big( \|u\|_{L^\infty(B_{3/4})}+\|v\|_{L^\infty(B_{3/4})}\Big)
.$$
Since, for every~$x$, $y\in B_{1/2}$,
$$ |\partial_{ij}w(x)-\partial_{ij}w(y)|\le C\|w\|_{C^3(B_{1/2})}|x-y|\le C\|w\|_{C^3(B_{1/2})}|x-y|^\alpha,$$
up to renaming~$C$, we thus conclude that
\begin{equation*} \|w\|_{C^{2,\alpha}(B_{1/2})}\le C\,\Big( \|u\|_{L^\infty(B_{3/4})}+\|v\|_{L^\infty(B_{3/4})}\Big)
.\end{equation*}
In light of this estimate and~\eqref{LSDNlwnfeECAocOSndfeDEFSGBD-UJNSDfgkrhVEUjfIMSdfdp01-2-09uyh-0ks89},
we have that
\begin{equation*} \begin{split}&\|u\|_{C^{2,\alpha}(B_{1/2})}\le
\|v\|_{C^{2,\alpha}(B_{1/2})}+
\|w\|_{C^{2,\alpha}(B_{1/2})}\le C\,\Big( \|u\|_{L^\infty(B_{3/4})}+
\|v\|_{C^{2,\alpha}(B_{3/4})}
\Big)\\&\qquad\qquad\qquad\qquad\le C\,\Big( \|u\|_{L^\infty(B_{3/4})}+
\|f\|_{C^{\alpha}(B_{1})}
\Big)
.\qedhere\end{split}\end{equation*}
\end{proof}

Now we turn our attention to the local estimates at the boundary to prove Theorem~\ref{SCHAUDER-BOUTH}.
To this end, we point out a variation of the result in Proposition~\ref{SI:IN:SE:DR}
(full details of this modification are provided for the facility of the reader):

\begin{proposition}\label{SI:IN:SE:DR-MEZPA}
Let~$\Gamma$ be the fundamental solution in~\eqref{GAMMAFU}.
Let~$\Omega_0\subset\R^n$ be an open and bounded set with Lipschitz boundary.
Let~$\Omega\subseteq\Omega_0$ be open and bounded,
and~$f\in L^\infty(\R^n)$. Suppose that~$f\in C^\alpha(\Omega)$
and~$f=0$ in~$\R^n\setminus\Omega$.

Set
$$ v(x):=-\int_\Omega \Gamma(x-y)\,f(y)\,dy.$$ Then, for all~$i$, $j\in\{1,\dots,n\}$ and all~$x\in\Omega$,
\begin{equation}\label{po346ro5l8et-0iSTHndfOjmfRhfgMNHGrgbd92ies} \partial_{ij}v(x)=\int_{\Omega_0}\partial_{ij}\Gamma(x-y)\big(f(x)-f(y)\big)\,dy+f(x)\int_{\partial\Omega_0}\partial_i\Gamma(x-y)\,\nu(y)\cdot e_j\,d{\mathcal{H}}^{n-1}_y.\end{equation}
\end{proposition}

\begin{proof}
As in the proof of Proposition~\ref{SI:IN:SE:DR}, given~$\rho>0$
to be taken as small as we wish, we use the regularization~$\Gamma_\rho$ of the fundamental solution introduced in~\eqref{COIENEE-0986ytufgkbv-0-rjfeonvnb2GDB} and~\eqref{GAROGRA}.
We define~$V(x):=\partial_i v(x)$ and
$$V_\rho(x):=-\partial_i\left(
\int_\Omega \Gamma_\rho(x-y)\,f(y)\,dy\right).$$
We know from~\eqref{6787iudsfeVAskdncfgFwasfdgfFASBdcNidkHKihfvbnLI3Sd0}
(applied here with~$f\chi_\Omega$ instead of~$f$) that
\begin{equation}\label{CINdfNjroe876robl3a-93riMS-cnr02-1080}
{\mbox{$V_\rho$ converges uniformly to~$V$.}}\end{equation}
Additionally, if~$x\in\Omega'\Subset\Omega$,
\begin{eqnarray*}
&&\left|\partial_jV_\rho(x)-\int_{\Omega_0}\partial_{ij}\Gamma(x-y)\big(f(x)-f(y)\big)\,dy-f(x)\int_{\partial\Omega_0}\partial_i\Gamma(x-y)\,\nu(y)\cdot e_j\,d{\mathcal{H}}^{n-1}_y\right|\\&=&
\left|
\int_\Omega\partial_{ij} \Gamma_\rho(x-y)\,f(y)\,dy+\int_{\Omega_0}\partial_{ij}\Gamma(x-y)\big(f(x)-f(y)\big)\,dy\right.
\\&&\qquad\left.
+f(x)\int_{\partial\Omega_0}\partial_i\Gamma(x-y)\,\nu(y)
\cdot e_j\,d{\mathcal{H}}^{n-1}_y\right|\\
&=&\left|
\int_{B_\rho(x)}\partial_{ij} \Gamma_\rho(x-y)\,f(y)\,dy+
\int_{\Omega\setminus B_\rho(x)}\partial_{ij} \Gamma (x-y)\,f(y)\,dy
+
\int_{B_\rho(x)}\partial_{ij}\Gamma(x-y)\big(f(x)-f(y)\big)\,dy\right.
\\&&\qquad\left.+
\int_{\Omega_0\setminus B_\rho(x)}\partial_{ij}\Gamma(x-y)\big(f(x)-f(y)\big)\,dy
+f(x)\int_{\partial\Omega_0}\partial_i\Gamma(x-y)\,\nu(y)
\cdot e_j\,d{\mathcal{H}}^{n-1}_y\right|\\
&=&\left|
-\delta_{ij}
\int_{B_\rho(x)}\frac{f(y)}{n\,|B_\rho|}\,dy
+f(x)
\int_{\Omega_0\setminus B_\rho(x)}\partial_{ij}\Gamma(x-y)\,dy
+f(x)\int_{\partial\Omega_0}\partial_i\Gamma(x-y)\,\nu(y)
\cdot e_j\,d{\mathcal{H}}^{n-1}_y\right|\\&&\qquad
+C[f]_{C^\alpha(B_\rho(x))}\int_{B_\rho(x)}|x-y|^{\alpha-n}\,dy.
\end{eqnarray*}
We also notice that
\begin{eqnarray*}&&
\int_{\Omega_0\setminus B_\rho(x)}\partial_{ij}\Gamma(x-y)\,dy
=-\int_{\Omega_0\setminus B_\rho(x)}\div_y \big(\partial_{i}\Gamma(x-y)\,e_j\big)\,dy\\
&&\qquad=-\int_{\partial(\Omega_0\setminus B_\rho(x))}
\partial_{i}\Gamma(x-y)\,\nu(y)\cdot e_j \,d{\mathcal{H}}^{n-1}_y\\
&&\qquad=-\int_{\partial \Omega_0}
\partial_{i}\Gamma(x-y)\,\nu(y)\cdot e_j \,d{\mathcal{H}}^{n-1}_y
+\bar c_n\int_{\partial B_\rho(x)}
\frac{(x_i-y_i)(x_j-y_j)}{|x-y|^{n+1}} \,d{\mathcal{H}}^{n-1}_y,
\end{eqnarray*}
where, recalling~\eqref{COIENEE},
\begin{equation} \label{098iuffdwqLso0djhsRhnPOkbfiEOjmovL90P0-986yr3}\begin{split}
\bar c_n\,&:=\,\begin{dcases}
c_n & {\mbox{ if }}n=2,\\
c_n(n-2)& {\mbox{ if }}n\ne2
\end{dcases}\\ &=\,\frac1{n\,|B_1|}.\end{split}\end{equation}
These observations lead to
\begin{eqnarray*}
&&\left|\partial_jV_\rho(x)-\int_{\Omega_0}\partial_{ij}\Gamma(x-y)\big(f(x)-f(y)\big)\,dy-f(x)\int_{\partial\Omega_0}\partial_i\Gamma(x-y)\,\nu(y)\cdot e_j\,d{\mathcal{H}}^{n-1}_y\right|\\
&\le&
\left|
-\delta_{ij}
\int_{B_\rho(x)}\frac{f(y)}{n\,|B_\rho|}\,dy
+\bar c_n\,f(x)\int_{\partial B_\rho(x)}
\frac{(x_i-y_i)(x_j-y_j)}{|x-y|^{n+1}} \,d{\mathcal{H}}^{n-1}_y
\right|
+C\,[f]_{C^\alpha(B_\rho(x))}\,\rho^\alpha.
\end{eqnarray*}
Now we note that
\begin{eqnarray*}&&
\bar c_n\int_{\partial B_\rho(x)}
\frac{(x_i-y_i)(x_j-y_j)}{|x-y|^{n+1}} \,d{\mathcal{H}}^{n-1}_y=
\bar c_n\delta_{ij}\int_{\partial B_\rho(x)}\frac{(x_i-y_i)^2}{|x-y|^{n+1}} \,d{\mathcal{H}}^{n-1}_y\\&&\qquad=
\frac{\bar c_n\delta_{ij}}{n} \int_{\partial B_\rho(x)}\frac{|x-y|^2}{|x-y|^{n+1}} \,d{\mathcal{H}}^{n-1}_y=
\frac{\bar c_n\delta_{ij}\,{\mathcal{H}}^{n-1}(\partial B_\rho)}{n\,\rho^{n-1}}=
\frac{\delta_{ij}}{n},
\end{eqnarray*}
thanks to~\eqref{B1}.

{F}rom this, we arrive at
\begin{eqnarray*}
&&\left|\partial_jV_\rho(x)-\int_{\Omega_0}\partial_{ij}\Gamma(x-y)\big(f(x)-f(y)\big)\,dy-f(x)\int_{\partial\Omega_0}\partial_i\Gamma(x-y)\,\nu(y)\cdot e_j\,d{\mathcal{H}}^{n-1}_y\right|\\
&\le&
\left|
-\frac{\delta_{ij}}n
\fint_{B_\rho(x)}f(y)\,dy
+\frac{\delta_{ij}\,f(x)}{n}
\right|
+C\,[f]_{C^\alpha(B_\rho(x))}\,\rho^\alpha
\\&\le&C\,[f]_{C^\alpha(B_\rho(x))}\,\rho^\alpha
\end{eqnarray*}
up to renaming~$C>0$.

This entails that~$\partial_jV_\rho$ converges locally uniformly
to the right hand side of~\eqref{po346ro5l8et-0iSTHndfOjmfRhfgMNHGrgbd92ies}.
The desired result thus follows by recalling~\eqref{CINdfNjroe876robl3a-93riMS-cnr02-1080}.
\end{proof}

Now we point out two useful regularity results
for the Newtonian potential\index{Newtonian potential} set in a halfspace. For this, we will use the notation
$$ \R^n_+:=\{x=(x',x_n)\in \R^{n-1}\times\R {\mbox{ s.t. }}x_n>0\}.$$

\begin{proposition}\label{QUA67809}
Let~$\Gamma$ be the fundamental solution in~\eqref{GAMMAFU}
and~$f\in C^\alpha(\R^{n}_+)$. Suppose that~$f(x)=0$ for all~$x\in\R^n_+ $ with~$|x|>\frac{99}{100}$.

Set
$$ v(x):=-\int_{\R^n_+} \Gamma(x-y)\,f(y)\,dy.$$ Then, for all~$i$, $j\in\{1,\dots,n\}$,
\begin{equation}\label{LA-0uojnSHGJSMEnodclsdCOSKdflrTAaXOCasNThdnfESG2}
[\partial_{ij}v]_{C^\alpha(\R^n_+)}\le
C\,[f]_{C^\alpha(\R^n_+)},
\end{equation}
for a suitable constant~$C>0$ depending only on~$n$ and~$\alpha$.
\end{proposition}

\begin{proof} We point out that
\begin{equation}\label{EFVA-LEQIAUTI}
\Delta v(x)=f(x)\qquad{\mbox{for all }}\,x\in\R^n_+
,\end{equation}
owing to~\eqref{SOLPONE-0-LOCALE} (applied here to the function~$f\chi_{\R^n_+}$).

We also note that
\begin{equation}\label{LA-0uojnSHGJSMEnodclsdCOSKdflrTAaXOCasNThdnfESG3}
{\mbox{it suffices to establish~\eqref{LA-0uojnSHGJSMEnodclsdCOSKdflrTAaXOCasNThdnfESG2}
under the additional assumption that either~$i\ne n$ or~$j\ne n$.}}
\end{equation}
Indeed, if~\eqref{LA-0uojnSHGJSMEnodclsdCOSKdflrTAaXOCasNThdnfESG2} 
holds true with~$(i,j)\ne(n,n)$ then one can exploit~\eqref{EFVA-LEQIAUTI}
to write that, in~$\R^n_+$,
$$\partial_{nn}v=f-\sum_{i=1}^{n-1}\partial_{ii}v$$
and thus
$$ [\partial_{nn}v]_{C^\alpha(\R^n_+)}\le
[f]_{C^\alpha(\R^n_+)}+\sum_{i=1}^{n-1}[\partial_{ii}v]_{C^\alpha(\R^n_+)}
\le C\,[f]_{C^\alpha(\R^n_+)},$$
up to renaming~$C$.
This observation establishes~\eqref{LA-0uojnSHGJSMEnodclsdCOSKdflrTAaXOCasNThdnfESG3}.

Thanks to~\eqref{LA-0uojnSHGJSMEnodclsdCOSKdflrTAaXOCasNThdnfESG3},
up to exchanging~$i$ and~$j$,
we can now focus on the case
\begin{equation}\label{HSND-TAGHSNimABL-LASMBA9lNBA-23}
{\mbox{$i\in\{1,\dots,n\}$\;and\;$j\in\{1,\dots,n-1\}$.}}\end{equation}
In this setting, given~$r>0$, we exploit Proposition~\ref{SI:IN:SE:DR-MEZPA}
with~$\Omega:=B^+_{r}$ and~$\Omega_0:=B^+_R$ with~$R>r$.
Noticing that along~$\{x_n=0\}$ we have that~$\nu\cdot e_j=0$, due to~\eqref{HSND-TAGHSNimABL-LASMBA9lNBA-23}, for every~$x\in B_r^+$
we have that
\begin{equation}\label{andrachiamata}\begin{split}
\partial_{ij}v(x)=\;&\int_{B_R^+}\partial_{ij}\Gamma(x-y)\big(f(x)-f(y)\big)\,dy+f(x)\int_{\partial B_R^+}\partial_i\Gamma(x-y)\,\nu(y)\cdot e_j\,d{\mathcal{H}}^{n-1}_y\\
=\;&\int_{B_R^+}\partial_{ij}\Gamma(x-y)\big(f(x)-f(y)\big)\,dy
+f(x)\int_{\partial B_R\cap\{y_n>0\}}\partial_i\Gamma(x-y)\,\nu(y)\cdot e_j\,d{\mathcal{H}}^{n-1}_y.
\end{split}\end{equation}
Now we observe that, if~$x\in B_r^+$ and~$y\in\partial B_R$, as~$R\to+\infty$,
\begin{eqnarray*}
\frac{1}{|x-y|^n}=\frac1{R^n\left| \frac{x}{R}-\frac{y}R\right|^n}=\frac{1+o(1)}{R^n},
\end{eqnarray*}
and therefore, using also that~$j\neq n$,
\begin{eqnarray*}&&
\int_{\partial B_R\cap\{y_n>0\}}\partial_i\Gamma(x-y)\,\nu(y)\cdot e_j\,d{\mathcal{H}}^{n-1}_y
=\frac{c_n(n-2)}R \int_{\partial B_R\cap\{y_n>0\}}\frac{(x_i-y_i)y_j}{|x-y|^n}\,d{\mathcal{H}}^{n-1}_y\\&&\qquad
\qquad
=\frac{c_n(n-2)(1+o(1))}{R^{n+1}} \int_{\partial B_R\cap\{y_n>0\}}(x_i-y_i)y_j\,d{\mathcal{H}}^{n-1}_y\\
&&\qquad\qquad= -\frac{c_n(n-2)(1+o(1))}{R^{n+1}}
\int_{\partial B_R\cap\{y_n>0\}}y_i y_j\,d{\mathcal{H}}^{n-1}_y.
\end{eqnarray*}
Thus, if~$i\neq j$, by odd symmetry, we have that
$$\int_{\partial B_R\cap\{y_n>0\}}\partial_i\Gamma(x-y)\,\nu(y)\cdot e_j\,d{\mathcal{H}}^{n-1}_y=0$$
while if~$i=j\in\{1,\dots,n-1\}$
\begin{eqnarray*}&&
\int_{\partial B_R\cap\{y_n>0\}}\partial_i\Gamma(x-y)\,\nu(y)\cdot e_j\,d{\mathcal{H}}^{n-1}_y
=-\frac{c_n(n-2)(1+o(1))}{R^{n+1}}
\int_{\partial B_R\cap\{y_n>0\}}y_i^2\,d{\mathcal{H}}^{n-1}_y\\&&\qquad\qquad
=-\frac{c_n(n-2)(1+o(1))}{2R^{n+1}}
\int_{\partial B_R}y_i^2\,d{\mathcal{H}}^{n-1}_y=
-\frac{c_n(n-2)(1+o(1))}{2nR^{n+1}}
\int_{\partial B_R}|y|^2\,d{\mathcal{H}}^{n-1}_y\\&&\qquad\qquad=
-\frac{c_n(n-2){\mathcal{H}}^{n-1}(\partial B_1)(1+o(1))}{2n},
\end{eqnarray*}
as~$R\to+\infty$.

Hence, taking the limit as~$R\to+\infty$ in~\eqref{andrachiamata} and using these observations,
we have that, for every~$x\in B_r^+$,
\begin{equation}\label{s20e3wt584376uyidAsdfhtrjtyujtyututuPPP}
\partial_{ij}v(x)=\int_{\R^n_+}\partial_{ij}\Gamma(x-y)\big(f(x)-f(y)\big)\,dy 
-f(x)\delta_{ij}\frac{c_n(n-2){\mathcal{H}}^{n-1}(\partial B_1)}{2n}.\end{equation}
The aim is now to exploit
Lemma~\ref{NKSD5678yihnqehr837tCOvnlIDJMcINBSno3Na7syshdrr9ISDcvaHF23},
with~$K:=\partial_{ij}\Gamma$. Since~$K$ satisfies~\eqref{KERNE-SINGO-001} and~\eqref{KERNE-SINGO-003}, in order to use Lemma~\ref{NKSD5678yihnqehr837tCOvnlIDJMcINBSno3Na7syshdrr9ISDcvaHF23} we have
to check that the condition in~\eqref{ABBE9} is satisfied. To this end,
recalling the principal value notation in~\eqref{KS-TSGDBfojeISLmadne-ojdfnLS-KMDMDND9uehn},
for every~$z_0\in\R^n$ and every~$\rho>0$, \begin{eqnarray*}&&-\int_{\R^n_+\setminus B_\rho(z_0)}\partial_{ij}\Gamma(z_0-y)\,dy=
-\lim_{R\to+\infty}\int_{\R^n_+\cap (B_R(z_0)\setminus B_\rho(z_0))}\partial_{ij}\Gamma(z_0-y)\,dy\\&&\quad
=-\lim_{R\to+\infty}
\int_{\R^n_+\cap (B_R(z_0)\setminus B_\rho(z_0))}\div_y\big(\partial_{i}\Gamma(z_0-y)\,e_j\big)\,dy.
\end{eqnarray*}
As a result, by the Divergence Theorem,
\begin{equation}\label{KeqM-0ol-ekrfm29ifhweoghfougoehvKybMIf}
\left|\int_{\R^n_+\setminus B_\rho(z_0)}\partial_{ij}\Gamma(z_0-y)\,dy\right|\le
\int_{\partial B_R(z_0)} |\nabla\Gamma(z_0-y)| \,d{\mathcal{H}}_y^{n-1}+
\int_{\partial B_\rho(z_0)} |\nabla\Gamma(z_0-y)| \,d{\mathcal{H}}_y^{n-1}.
\end{equation}
We also observe that, for every~$r>0$,
$$ \int_{\partial B_r(z_0)} |\nabla\Gamma(z_0-y)| \,d{\mathcal{H}}_y^{n-1}\le
C \int_{\partial B_\rho(z_0)} \frac{d{\mathcal{H}}_y^{n-1}}{ |z_0-y|^{n-1}}\le C,$$
up to renaming~$C>0$.

Combining this with~\eqref{KeqM-0ol-ekrfm29ifhweoghfougoehvKybMIf}
we see that
\begin{equation}\label{s2u834738vfhgduashsw23y3tyJJJJ}
{\mbox{\eqref{ABBE9} holds true with~$K=\partial_{ij}\Gamma$.}}\end{equation}
Thus, recalling~\eqref{s20e3wt584376uyidAsdfhtrjtyujtyututuPPP}
and exploiting Lemma~\ref{NKSD5678yihnqehr837tCOvnlIDJMcINBSno3Na7syshdrr9ISDcvaHF23}
with~${\mathcal{V}}:=\R^n_+$,
we obtain that
\begin{eqnarray*}
[\partial_{ij} v]_{C^\alpha(B_r^+)}\le [Tf]_{C^\alpha(B_r^+)}+C[f]_{C^\alpha(B_r^+)}\le
[Tf]_{C^\alpha(\R^n_+)}+C[f]_{C^\alpha(\R^n_+)}\le C[f]_{C^\alpha(\R^n_+)},
\end{eqnarray*}
up to renaming~$C$.
This and the arbitrariness of~$r$ give the desired result.
\end{proof}

\begin{proposition}\label{QUA67809BIS}
Let~$\Gamma$ be the fundamental solution in~\eqref{GAMMAFU}
and~$f\in C^\alpha(\R^{n}_+)$. Suppose that~$f(x)=0$ for all~$x\in\R^n_+ $ with~$|x|>\frac{99}{100}$.

Set
$$ v(x):=-\int_{\R^n_+} \Gamma(x-y)\,f(y)\,dy.$$ Then, for all~$i$, $j\in\{1,\dots,n\}$,
\begin{equation}\label{LA-0uojnSHGJSMEnodclsdCOSKdflrTAaXOCasNThdnfESG2BIS}
\|\partial_{ij}v\|_{L^\infty(\R^n_+)}\le
C\,\|f\|_{C^\alpha(\R^n_+)},
\end{equation}
for a suitable constant~$C>0$ depending only on~$n$ and~$\alpha$.
\end{proposition}

\begin{proof} We argue as in the proof of Proposition~\ref{QUA67809}, exploiting~\eqref{EFVA-LEQIAUTI}
to conclude that
it suffices to establish~\eqref{LA-0uojnSHGJSMEnodclsdCOSKdflrTAaXOCasNThdnfESG2BIS}
under the additional assumption that either~$i\ne n$ or~$j\ne n$.
In this setting, one can employ~\eqref{s20e3wt584376uyidAsdfhtrjtyujtyututuPPP} to find that,
for every~$r>0$, in~$B_r^+$,
\begin{equation}\label{sweryr5ty5ye34765437565839}
\partial_{ij}v(x)=\int_{\R^n_+}\partial_{ij}\Gamma(x-y)\big(f(x)-f(y)\big)\,dy 
-f(x)\delta_{ij}\frac{c_n(n-2){\mathcal{H}}^{n-1}(\partial B_1)}{2n}.\end{equation}
Now we set
$$\phi(x):=\int_{\R^n_+}\partial_{ij}\Gamma(x-y)\big(f(x)-f(y)\big)\,dy .$$
We observe that there exists~$x_0\in B_1^+$ such that~$f(x_0)=0$, and therefore
\begin{eqnarray*} &&
|\phi(x_0)|\le C\int_{\R^n_+}\frac{\big|f(y)\big|}{|x_0-y|^{n}}\,dy=
C\int_{B_1^+}\frac{\big|f(y)\big|}{|x_0-y|^{n}}\,dy=
C\int_{B_1^+}\frac{\big|f(x_0)-f(y)\big|}{|x_0-y|^{n}}\,dy\\&&\qquad\le
C[f]_{C^\alpha(\R^n_+)}\int_{B_1^+}|x_0-y|^{\alpha-n}\,dy\le C[f]_{C^\alpha(\R^n_+)},
\end{eqnarray*}
up to renaming~$C>0$.
Consequently, thanks to Lemma~\ref{NKSD5678yihnqehr837tCOvnlIDJMcINBSno3Na7syshdrr9ISDcvaHF23},
applied here with~${\mathcal{V}}:=\R^n_+$
and~$K:=\partial_{ij}\Gamma$ (recall~\eqref{s2u834738vfhgduashsw23y3tyJJJJ}), we have that
for every~$x\in B_4^+$,
\begin{eqnarray*}
|\phi(x)|\le |\phi(x_0)|+|\phi(x)-\phi(x_0)|\le  C[f]_{C^\alpha(\R^n_+)}
+[\phi]_{C^\alpha(B_4^+)}|x-x_0|^\alpha \le C[f]_{C^\alpha(\R^n_+)},
\end{eqnarray*}
up to relabeling~$C>0$.

Furthermore, we notice that if~$x\in \R^n_+\setminus B_4$ and~$y\in B_1^+$, then~$|x-y|\ge|x|-|y|\ge 3$, and thus
\begin{eqnarray*}
|\phi(x)|\le C\int_{B_1^+}\frac{\big|f(y)\big|}{|x-y|^{n}}\,dy\le C\|f\|_{L^\infty(B_1^+)}\int_{B_1^+} \frac{dy}{3^n}\le
C\|f\|_{L^\infty(B_1^+)}.
\end{eqnarray*}
These considerations show that
$$\|\phi\|_{L^\infty(\R^n_+)}\le C\|f\|_{C^\alpha(\R^n_+)}.$$
{F}rom this and~\eqref{sweryr5ty5ye34765437565839} we obtain  that
$$\|\partial_{ij}v\|_{L^\infty(B^+_r)}\le C\|f\|_{C^\alpha(\R^n_+)}.$$
Since~$r$ is any arbitrary positive real number, this gives
the desired estimate.
\end{proof}

We now complete the proof of the desired
local estimates\footnote{A conceptually different proof of
Theorem~\ref{SCHAUDER-BOUTH} will be presented on page~\pageref{234556df0323iurjw924324jfJAMAnJHOSo3o24f32tofXZwf}.}
at the boundary:

\begin{proof}[Proof of Theorem~\ref{SCHAUDER-BOUTH}]
We first observe that
\begin{equation}\label{omMAEwgNGLt4IONkjgfnbdsdtyuijk9024i5SZERSD}
{\mbox{it suffices to prove Theorem~\ref{SCHAUDER-BOUTH} when~$g$ vanishes identically.}}
\end{equation}
Indeed, let~$\widetilde u(x)=\widetilde u(x',x_n):=u(x)-g(x')$
and~$\widetilde f(x)=\widetilde f(x',x_n):=f(x)-\Delta g(x')$.
By~\eqref{poknISOdu-oM-EN0d-f0123-s} we have that $$
\begin{dcases} \Delta \widetilde u=\widetilde f\quad{\mbox{in }}\,B_1^+,\\
\widetilde u= 0\quad{\mbox{on }}\,B_1^0.\end{dcases}$$
Consequently, if Theorem~\ref{SCHAUDER-BOUTH} holds true when~$g$ vanishes identically, we have
\begin{eqnarray*}&& 
\|u\|_{C^{2,\alpha}(B_{1/2}^+)}\le
\|\widetilde u\|_{C^{2,\alpha}(B_{1/2}^+)}+
\|g\|_{C^{2,\alpha}(B_{1}^0)}\le C\,\Big(\|\widetilde u\|_{L^\infty(B_1^+)}+\|\widetilde f\|_{C^\alpha(B_1^+)}\Big)+\|g\|_{C^{2,\alpha}(B_{1}^0)}\\&&\qquad\qquad\qquad\qquad\le
C\,\Big(\|u\|_{L^\infty(B_1^+)}+\|g\|_{C^{2,\alpha}(B_1^0)}+\|f\|_{C^\alpha(B_1^+)}\Big),
\end{eqnarray*}
up to renaming~$C$ at each step of the calculation, and this observation establishes~\eqref{omMAEwgNGLt4IONkjgfnbdsdtyuijk9024i5SZERSD}.

Thus, in the light of~\eqref{omMAEwgNGLt4IONkjgfnbdsdtyuijk9024i5SZERSD}, we focus now on the proof
of Theorem~\ref{SCHAUDER-BOUTH} by taking the additional assumption that
\begin{equation}\label{JNSptgbTSODIAN0n53sdfAC78ZiZCmUAM}
{\mbox{$g$ vanishes identically.}}\end{equation}
Furthermore, by multiplying by a cutoff function~$\varphi\in C^\infty_0(B_{8/9},\,[0,1])$,
with~$\varphi=1$ in~$B_{3/4}$ and~$|\nabla \varphi|\le 10$ as done at the beginning of the proof of Theorem~\ref{SCHAUDER-INTE}
we can suppose that
\begin{equation}\label{1J13545y245M54utjg35645h43wtgS-01293ecyrgfuywiedjkasgcv}
{\mbox{$f=0$ \, in \, $\R^n_+\setminus B_{8/9}^+$. }}\end{equation}
We use the following handy notation to deal with reflections across the horizontal hyperplane:
given~$x=(x',x_n)\in\R^{n-1}\times\R$, we let~$x^\star:=(x',-x_n)$.
Also, since~$f$ is uniformly continuous in~$B_1^+$, it can be extended continuously to~$B_1^0$.
With this notation in mind, for all~$x=(x',x_n)\in \R^n$ we define
\begin{equation} \label{fstarywueut878788}
f^\star(x):=\begin{dcases}
f(x)&{\mbox{ if }}x_n\ge0,
\\f(x^\star)&{\mbox{ if }}x_n<0
\end{dcases}\end{equation}
and we claim that
\begin{equation}\label{HIUASC2A3SHdnofUj9234rRaksmd}
[f^\star]_{C^\alpha(\R^n)}\le [f]_{C^\alpha(\R^n_+)}.
\end{equation}
To check this, let~$x$, $y\in\R^n$.
If both~$x_n$ and~$y_n$ are nonnegative,
then~$|f^\star(x)-f^\star(y)|=|f(x)-f(y)|\le[f]_{C^\alpha(\R^n_+)}|x-y|^\alpha$,
and~\eqref{HIUASC2A3SHdnofUj9234rRaksmd} plainly follows.
Similarly, if both~$x_n$ and~$y_n$ are negative,
then~$|f^\star(x)-f^\star(y)|=|f(x^\star)-f(y^\star)|\le[f]_{C^\alpha(\R^n_+)}|x^\star-y^\star|^\alpha=[f]_{C^\alpha(\R^n_+)}|x-y|^\alpha$,
from which we obtain~\eqref{HIUASC2A3SHdnofUj9234rRaksmd}.
Therefore, up to exchanging~$x$ and~$y$, we may assume that~$x_n\ge0>y_n$.
In this case, we have that~$|x_n+y_n|\le|x_n|+|y_n|=x_n-y_n\le|x_n-y_n|$ and consequently
\begin{eqnarray*}&&|f^\star(x)-f^\star(y)|=|f(x)-f(y^\star)|
\le[f]_{C^\alpha(\R^n_+)}|x-y^\star|^\alpha
=[f]_{C^\alpha(\R^n_+)}\big(|x'-y'|^2+|x_n+y_n|^2\big)^{\frac\alpha2}\\&&\qquad\le
[f]_{C^\alpha(\R^n_+)}\big(|x'-y'|^2+|x_n-y_n|^2\big)^{\frac\alpha2}=[f]_{C^\alpha(\R^n_+)}|x-y|^\alpha,
\end{eqnarray*} which completes the proof of~\eqref{HIUASC2A3SHdnofUj9234rRaksmd}.

Let now
\begin{equation*}
w(x):=\int_{\R^n_+}\big(\Gamma(x-y)-\Gamma(x-y^\star)\big)\,f(y)\,dy.\end{equation*}
We claim that
\begin{equation}\label{44y5uF92UjIMS-097y8tfy-9yiwgjefbHSNA9Nma243565n4n3a9}
\lim_{x_n\searrow0} w(x',x_n)=0.
\end{equation}
Indeed, 
$$ |x-y^\star|=\sqrt{|x'-y'|^2+|x_n+y_n|^2}=|x^\star-y|$$
and therefore,
for every~$y\in B_1$,
$$ \big|\Gamma(x-y)-\Gamma(x-y^\star)\big|
=\big|\Gamma(x-y)-\Gamma(x^\star-y)\big|
\le C\,|x_n|\int_{-1}^1 \frac{dt}{|(x'-y',tx_n-y_n)|^{n-1}}
,$$
for some~$C>0$.
As a result, by~\eqref{1J13545y245M54utjg35645h43wtgS-01293ecyrgfuywiedjkasgcv}
and exploiting Fubini's Theorem and the change of variable~$z:=(x'-y',tx_n-y_n)$,
\begin{eqnarray*}&&
|w(x',x_n)|\le\|f\|_{L^\infty(\R^n_+)}\int_{B_1}\big|\Gamma(x-y)-\Gamma(x-y^\star)\big|\,dy\\&&\qquad
\le C\,|x_n|\int_{-1}^1\left(\int_{B_1}\frac{dy}{|(x'-y',tx_n-y_n)|^{n-1}}\right)\,dt\le 
C\,|x_n| \int_{-1}^1\left(\int_{B_3}\frac{dz}{|z|^{n-1}}\right)\,dt
\le C\,|x_n|,
\end{eqnarray*}
up to renaming~$C>0$,
from which~\eqref{44y5uF92UjIMS-097y8tfy-9yiwgjefbHSNA9Nma243565n4n3a9}
follows.

Let also
$$ v(x):=\int_{\R^n_+}\Gamma(x-y)\,f(y)\,dy\qquad{\mbox{and}}\qquad
W(x):=\int_{\R^n}\Gamma(x-y)\,f^\star(y)\,dy.$$
Using the notation~$\R^n_-:=\{x=(x',x_n)\in\R^n {\mbox{ s.t. }}x_n<0\}$
and the substitution~$z:=y^\star$, we remark that
\begin{equation}\begin{split}\label{4.1.41BIS}
w(x)=\,&
\int_{\R^n_+} \Gamma(x-y)\,f(y)\,dy-\int_{\R^n_+}\Gamma(x-y^\star)\,f(y)\,dy\\
=\,&\int_{\R^n_+} \Gamma(x-y)\,f(y)\,dy-\int_{\R^n_-}\Gamma(x-z)\,f(z^\star)\,dz
\\=\,&\int_{\R^n_+} \Gamma(x-y)\,f(y)\,dy-\int_{\R^n_-}\Gamma(x-z)\,f^\star(z)\,dz
\\=\,&2\int_{\R^n_+} \Gamma(x-y)\,f(y)\,dy-\int_{\R^n}\Gamma(x-z)\,f^\star(z)\,dz\\
=\,& 2v(x)-W(x).
\end{split}\end{equation}
Thus, since, by Corollary~\ref{QUA4} and~\eqref{HIUASC2A3SHdnofUj9234rRaksmd},
$$ [\partial_{ij}W]_{C^\alpha(\R^n)}\le
C\,[f^\star]_{C^\alpha(\R^n)}\le C\,[f]_{C^\alpha(\R^n_+)},$$
and, by Proposition~\ref{QUA67809},
$$ [\partial_{ij}v]_{C^\alpha(\R^n_+)}\le
C\,[f]_{C^\alpha(\R^n_+)},$$
we conclude that, for all~$i$, $j\in\{1,\dots,n\}$,
\begin{equation}\label{0oiqjwqefv9-0okd-203jfj3j342jPM}
[\partial_{ij}w]_{C^\alpha(\R^n_+)}\le2 [\partial_{ij}v]_{C^\alpha(\R^n_+)}+
[\partial_{ij}W]_{C^\alpha(\R^n)}\le C\,[f]_{C^\alpha(\R^n_+)}.\end{equation}
Furthermore, owing to~\eqref{SOLPONE-0-LOCALE}, for all~$x\in\R^n_+$,
\begin{equation}\label{23EGbstt3eIcHi78UmPmfK2S-2093e}
-\Delta w(x)=-2\Delta v(x)+\Delta W(x)=2f(x)-f^\star(x)=f(x).
\end{equation}
We also point out that, for all~$x\in\R^n_+$,
\begin{equation}\label{K334yA2M-lasFI2fevnKrefAcC}
|w(x)|+|\nabla w(x)|\le C\| f\|_{L^\infty(\R^n_+)},\end{equation}
thanks to~\eqref{PRIMA98DER} and~\eqref{IGSBNHJjusjdvnICtsghdfnEJHACANosdat2R4A9rfX2-9}.

Moreover, by~\eqref{HIUASC2A3SHdnofUj9234rRaksmd},
\eqref{4.1.41BIS}, Lemma~\ref{lemmaspe0} and Proposition~\ref{QUA67809BIS}, we have that,
for all~$x\in\R^n_+$,
\begin{equation}\label{K334yA2M-lasFI2fevnKrefAcC22}
|\partial_{ij}w(x)|\le 2|v(x)|+|W(x)|\le C\left(\|f\|_{C^\alpha(\R^n_+)}+\|f^\star\|_{C^\alpha(\R^n)}\right)
\le C\|f\|_{C^\alpha(\R^n_+)}.
.\end{equation}

Now we define~$V:=u+w$. In light of~\eqref{44y5uF92UjIMS-097y8tfy-9yiwgjefbHSNA9Nma243565n4n3a9} and~\eqref{23EGbstt3eIcHi78UmPmfK2S-2093e}, we see that~$V\in C^2(B_1^+)\cap C(B_1^+\cap B_1^0)$,
and, recalling~\eqref{poknISOdu-oM-EN0d-f0123-s} and~\eqref{JNSptgbTSODIAN0n53sdfAC78ZiZCmUAM},
\begin{equation*}\begin{dcases} \Delta V=0\quad{\mbox{in }}\,B_1^+,\\
V=0 \quad{\mbox{on }}\,B_1^0.\end{dcases}\end{equation*}
For this reason and Lemma~\ref{RIFLESCHZ}, we have that
the odd reflection of~$V$ across~$B_1^0$, that we denote by~$V^\sharp$,
is harmonic in~$B_1$. As a result, exploting the 
Cauchy's Estimates in Theorem~\ref{CAUESTIMTH}, and employing~\eqref{K334yA2M-lasFI2fevnKrefAcC},
\begin{eqnarray*}&& \|V\|_{C^3(B_{1/2}^+)}\le\|V^\sharp\|_{C^3(B_{1/2})}\le C\|V^\sharp\|_{L^1(B_{7/8})}\le
C\,\Big(\|u\|_{L^\infty(B_{7/8}^+)}+ \|w\|_{L^\infty(B_{7/8}^+)}\Big)\\&&\qquad\qquad\qquad\qquad\le
C\,\Big(\|u\|_{L^\infty(B_{7/8}^+)}+ \|f\|_{L^\infty(\R^n_+)}\Big).\end{eqnarray*}
The desired result now follows from this estimate and~\eqref{0oiqjwqefv9-0okd-203jfj3j342jPM}
(recall also~\eqref{K334yA2M-lasFI2fevnKrefAcC} and~\eqref{K334yA2M-lasFI2fevnKrefAcC22}).
\end{proof}

By combining the
interior estimates of Theorem~\ref{SCHAUDER-INTE}
and the local estimates at the boundary of Theorem~\ref{SCHAUDER-BOUTH} one can
promptly obtain global estimates up to the boundary when\footnote{We point out that
global regularity results on balls are technically easier than those in general domains
(compare e.g. with the forthcoming Theorem~\ref{ThneGiunfBB032t4jNKS-s3i4}).
Indeed, for general domains a natural strategy is to ``straighten the boundary''
via a diffeomorphism (as it will be sketched
in Figure~\ref{fLOCAPAere8MPyItangeFqI} on page~\pageref{fLOCAPAere8MPyItangeFqI})
and then exploit the local estimates at the boundary obtained in halfballs:
however this procedure changes the operator (in particular, the corresponding equation obtained
after straightening the boundary is not anymore driven by the Laplace operator). This will be for us one
of the main motivations to study operators that are ``slightly more general than the Laplacian'':
the complications arisen from this further generality will be well compensated by the
fact that this class of operator is stable for diffeomorphism, thus allowing to straighten the boundary
(compare with the forthcoming Lemma~\ref{92o3jwelg02jfj-qiewjdfnoewfh9ewohgveioqNONDIFO234RM}).

The reason for which balls are special domains in this setting is that
one can use in this case the Kelvin Transform which
locally straightens the boundary of the ball and also preserves the Laplace operator: this will be
in fact the core of the proof of Theorem~\ref{KS-2rkjoewuj5-2prlfgk-a}.} the domain is a ball,
as pointed out in the following result:

\begin{theorem}\label{KS-2rkjoewuj5-2prlfgk-a}
Let~$f\in C^\alpha(B_1)$ and~$\varphi\in C^{2,\alpha}(\partial B_1)$
for some~$\alpha\in(0,1)$. Let~$u\in C^2(B_1)$ be a solution of
\begin{equation*} \begin{dcases}
\Delta u=f\quad{\mbox{in }}\,B_1,\\
u=\varphi\quad{\mbox{on }}\,\partial B_1.\end{dcases}
\end{equation*}
Then, there exists~$C>0$, depending only on~$n$ and~$\alpha$, such that
$$ \|u\|_{C^{2,\alpha}(B_1)}\le C\,\Big(\|f\|_{C^\alpha(B_1)}+\|\varphi\|_{C^{2,\alpha}(\partial B_1)}\Big).$$
\end{theorem}

\begin{proof} Consider a point~$\overline{x}\in\partial B_1$.
By a (translation of a) Kelvin Transform~${\mathcal{K}}$,
a neighborhood~${\mathcal{B}}(\overline{x})$ of~$\overline{x}$ in~$B_1$ is sent to a half ball~$B^+$
centered at~${\mathcal{K}}(\overline{x})$, see Figure~\ref{K6EL-DIFfffI}
and Proposition~\ref{h1ds3rfdsK3456790-dskpodnglkbd3565EL-DIFfffI}.
Let us denote~$u_\star(x):=u({\mathcal{K}}(x))$ and exploit~\eqref{KS:09876543209876543lkjhgfduerhfnv3ue}
to see that~$\Delta u_\star=f_\star$ in~$B^+$ for a suitable~$f_\star$
with~$\|f_\star\|_{C^\alpha(B^+)}\le
C\,\|f\|_{C^\alpha(B_1)}$. Also, if~$B^0$ denotes the flat part of~$\partial B^+$,
we have that~$u_\star(x)=\varphi({\mathcal{K}}(x))=:\varphi_\star(x)$ for all~$x\in B^0$.
Note also that~$\|\varphi_\star\|_{C^{2,\alpha}(B^0)}\le C\,\|\varphi\|_{C^{2,\alpha}(\partial B_1)}$.
As a result, making use of the local estimates at the boundary
of Theorem~\ref{SCHAUDER-BOUTH},
\begin{eqnarray*}&&
\|u\|_{C^{2,\alpha}(\widetilde{\mathcal{B}}(\overline{x}))}\le C\,
\|u_\star\|_{C^{2,\alpha}(\widetilde B^+)}\le C\,\Big(\|u_\star\|_{L^\infty(B^+)}+\|\varphi_\star\|_{C^{2,\alpha}(B^0)}+\|f_\star\|_{C^\alpha(B^+)}\Big)\\&&\qquad\qquad\qquad\qquad\quad\le
C\,\Big(\|u \|_{L^\infty(B_1)}+\|\varphi \|_{C^{2,\alpha}(\partial B_1)}+\|f \|_{C^\alpha(B_1)}\Big),
\end{eqnarray*}
where~$\widetilde{\mathcal{B}}(\overline{x})\Subset {\mathcal{B}}(\overline{x})$ and~$\widetilde {\mathcal{B}}(\overline{x})$
is sent to~$\widetilde B^+\Subset B^+$.

By combining this with the interior estimates of Theorem~\ref{SCHAUDER-INTE}, we arrive at
\begin{equation} \label{SMEmfenZlkvck4fdMAimfp0lSeHdu}\|u\|_{C^{2,\alpha}(B_1)}\le C\,\Big(
\|u\|_{L^\infty(B_1)}+\|f\|_{C^\alpha(B_1)}+\|\varphi\|_{C^{2,\alpha}(\partial B_1)}\Big).\end{equation}

Now we consider the function
$$ w(x):=u(x)-\|f\|_{L^\infty(B_1)}\,|x|^2.$$
Noticing that~$\Delta w=f-2n\|f\|_{L^\infty(B_1)}\le0$, we deduce from the Weak Maximum Principle in
Corollary~\ref{WEAKMAXPLE}(ii) that, for all~$x\in B_1$,
$$ w(x)\ge\inf_{\partial B_1} w\ge-\|\varphi \|_{L^\infty(\partial B_1)}
-\|f\|_{L^\infty(B_1)}$$
and therefore
$$ u(x)\ge -\|\varphi \|_{L^\infty(\partial B_1)}
-2\|f\|_{L^\infty(B_1)}.$$
Similarly, by considering the function
$$ \widetilde w(x):=u(x)+\|f\|_{L^\infty(B_1)}\,|x|^2$$
and using the Weak Maximum Principle in
Corollary~\ref{WEAKMAXPLE}(i), one sees that,
for all~$x\in B_1$,
$$ u(x)\leq\|\varphi \|_{L^\infty(\partial B_1)}
+2\|f\|_{L^\infty(B_1)}.$$
These observations give that
$$ \|u \|_{L^\infty(B_1)}\le\|\varphi \|_{L^\infty(\partial B_1)}
+2\|f\|_{L^\infty(B_1)}\le
\|\varphi \|_{C^{2,\alpha}(\partial B_1)}
+2\|f\|_{C^{\alpha}(B_1)}.$$
The desired result follows from this inequality and~\eqref{SMEmfenZlkvck4fdMAimfp0lSeHdu}.
\end{proof}

We now mention that the interior estimates in Theorem~\ref{SCHAUDER-INTE} can also be
recast to obtain global estimates in terms of ``scaled norms''.\index{scaled norm}
For this, given an open set~$\Omega\subset\R^n$, it is convenient to consider the distance functions
$$ d_x:={\rm dist}(x,\partial\Omega)=\inf_{p\in\partial\Omega}|x-p|\qquad{\mbox{and}}\qquad
d_{xy}:=\min\{d_x,d_y\},$$
and, for~$k\in\N$ and~$\alpha\in(0,1]$, take into account the weighted quantities
\begin{equation*}
\begin{split}
&\|u\|_{C^{k}_\star(\Omega)}:=\sum_{{\beta\in\N^n}\atop{|\beta|\le k}}
\sup_{x\in\Omega} d_x^{|\beta|}\,|D^\beta u(x)|,\\&
[u]_{C^{k,\alpha}_\star(\Omega)}:=
\sum_{{\beta\in\N^n}\atop{|\beta|= k}}\sup_{{x,y\in\Omega}\atop{x\ne y}}d_{xy}^{k+\alpha}\,\frac{|D^\beta u(x)-D^\beta u(y)|}{|x-y|^\alpha} \\
{\mbox{and }}\;&
\|u\|_{C^{k,\alpha}_\star(\Omega)}:=\|u\|_{C^{k}_\star(\Omega)}+[u]_{C^{k,\alpha}_\star(\Omega)}.\end{split}
\end{equation*}
Given~$m\in\N$, it is also convenient\footnote{One sees that~$\|u\|_{C^{\alpha}_{0,\star}(\Omega)}=\|u\|_{C^{0,\alpha}_\star(\Omega)}$. However, the norms~$\|\cdot\|_{C^{k,\alpha}_\star(\Omega)}$
and~$\|\cdot\|_{C^{\alpha}_{m,\star}(\Omega)}$ are conceptually different in general,
since~$\|\cdot\|_{C^{k,\alpha}_\star(\Omega)}$ involves derivatives up to order~$k$,
scaled with their own order of differentiation, while~$\|\cdot\|_{C^{\alpha}_{m,\star}(\Omega)}$
only involves the function, not its derivatives, but an additional scaling factor of order~$m$ takes place in this case.}
to define
$$ \|u\|_{C^{\alpha}_{m,\star}(\Omega)}:=\sup_{x\in\Omega} d_x^m\,|u(x)|
+\sup_{{x,y\in\Omega}\atop{x\ne y}}d_{xy}^{m+\alpha}\,\frac{|u(x)-u(y)|}{|x-y|^\alpha}.$$
These weighted norms should be compared with the classical ones in footnote~\ref{u9oj-NNp-OjtmRM-NosdTHSAMSATInfOdeN-N336OS235yS} on page~\pageref{u9oj-NNp-OjtmRM-NosdTHSAMSATInfOdeN-N336OS235yS}. With this notation,
we have the following weighted and global regularity result:

\begin{theorem} \label{243557754weighted91274and2i-3-3-520804global928485regularity}
Let~$\Omega\subset\R^n$ be an open set.
Let~$f\in C^\alpha(\Omega)$ for some~$\alpha\in(0,1)$. Let~$u\in C^2(\Omega)$ be a solution of
\begin{equation*} \Delta u=f\quad{\mbox{in }}\,\Omega.\end{equation*}
Then, there exists~$C>0$, depending only on~$n$, $\alpha$ and~$\Omega$, such that
$$ \|u\|_{C^{2,\alpha}_\star(\Omega)}\le C\,\Big(\|u\|_{L^\infty(\Omega)}+\|f\|_{C^\alpha_{2,\star}(\Omega)}\Big).$$
\end{theorem}

We remark that, in spite of their ``global'' flavor, regularity results
as the one in Theorem~\ref{243557754weighted91274and2i-3-3-520804global928485regularity} 
are not quite ``regularity results up to the boundary''
since the weighted norm~$\|\cdot\|_{C^{2,\alpha}_\star(\Omega)}$ allows a possible
degeneracy in the vicinity of the boundary of~$\Omega$;
still, Theorem~\ref{243557754weighted91274and2i-3-3-520804global928485regularity} does provide an interesting global information by giving
a uniform bound on the ``most singular behavior that can possibly happen
at the boundary'', regardless of the regularity of the boundary datum itself
(notice indeed that no regularity for the datum of the solution
along~$\partial\Omega$ is assumed in Theorem~\ref{243557754weighted91274and2i-3-3-520804global928485regularity}).

\begin{proof}[Proof of Theorem~\ref{243557754weighted91274and2i-3-3-520804global928485regularity}] Let~$x_0\in\Omega$ and~$R_0:=\frac{d_{x_0}}{8}$.
Let also~$u_0(x):=u\left(x_0+R_0 x\right)$. Then,
$\partial_i u_0(x)=R_0\partial_iu\left(x_0+R_0 x\right)$
and~$\partial_{ij}u_0(x)=R_0^2\partial_{ij}u\left(x_0+R_0 x\right)$.

As a result,
\begin{equation}\label{9jm-83utjhgbhiahfpegrhg-2ugoiewhg9f32gvberbto43hgubXifwjhnrVSg7tkjfiu43qgt}
d_{x_0}\,|\nabla u(x_0)|+
d_{x_0}^2\,|D^2 u(x_0)|=
\frac{d_{x_0}}{R_0}\,|\nabla u_0(0)|+
\frac{d_{x_0}^2}{R_0^2}\,|D^2 u_0(0)|\le
C\,\Big(|\nabla u_0(0)|+|D^2 u_0(0)|\Big).
\end{equation}
We also remark that~$\Delta u_0(x)=R_0^2\Delta u(x_0+R_0x)=R_0^2 f(x_0+R_0x)=:f_0(x)$
for every~$x\in B_1$. Therefore, in view of
Theorem~\ref{SCHAUDER-INTE},
\begin{equation}\label{LM-PARZIMNWEI} \|u_0\|_{C^{2,\alpha}(B_{1/2})}\le C\,\Big(\|u_0\|_{L^\infty(B_1)}+\|f_0\|_{C^\alpha(B_1)}\Big).\end{equation}
We also observe that if~$x\ne y\in B_1$, then
\begin{equation}\label{J32ledO0MS-NSdf0oA0r9igb32gbfTgdCXk-rKhdTTuikm0}\begin{split} &|f_0(x)-f_0(y)|\\=\;&R_0^2 |f(x_0+R_0x)-f(x_0+R_0y)|\\ \le\;&
[f]_{C^\alpha(B_{R_0}(x_0))}\,R_0^2\,|(x_0+R_0x)-(x_0+R_0y)|^\alpha\\ =\;&
[f]_{C^\alpha(B_{R_0}(x_0))}\,R_0^{2+\alpha}\,|x-y|^\alpha\\ =\;&
\sup_{{X,Y\in B_{R_0}(x_0)}\atop{X\ne Y}}\frac{|f(X)-f(Y)|\,R_0^{2+\alpha}\,|x-y|^\alpha}{|X-Y|^\alpha}.
\end{split}\end{equation}
Also, if~$X\in B_{R_0}(x_0)$ and~$p\in\partial\Omega$ we have that~$|X-p|\ge|x_0-p|-R_0\ge d_{x_0}-R_0\ge\frac{7d_{x_0}}8$
and therefore~$d_X\ge\frac{7d_{x_0}}8$.
Accordingly, if~$X$, $Y\in B_{R_0}(x_0)$, we have that~$d_{XY}\ge\frac{7d_{x_0}}8=7R_0$ and therefore
we infer from~\eqref{J32ledO0MS-NSdf0oA0r9igb32gbfTgdCXk-rKhdTTuikm0} that
\begin{equation*}\begin{split}& |f_0(x)-f_0(y)|\le C\,
\sup_{{X,Y\in B_{R_0}(x_0)}\atop{X\ne Y}}\frac{|f(X)-f(Y)|\,d_{XY}^{2+\alpha}\,|x-y|^\alpha}{|X-Y|^\alpha}\le C\,\|f\|_{C^\alpha_{2,\star}(\Omega)}\,|x-y|^\alpha.
\end{split}\end{equation*}
This gives that~$\|f_0\|_{C^\alpha(B_1)}\le C\,\|f\|_{C^\alpha_{2,\star}(\Omega)}$ and thus, recalling~\eqref{LM-PARZIMNWEI},
we arrive at
\begin{equation}\label{LM-PARZIMNWEI-2}
\|u_0\|_{C^{2,\alpha}(B_{1/2})}\le C\,\Big(\|u\|_{L^\infty(\Omega)}+\|f\|_{C^\alpha_{2,\star}(\Omega)}\Big),\end{equation}
up to renaming~$C$ repeatedly.

{F}rom~\eqref{9jm-83utjhgbhiahfpegrhg-2ugoiewhg9f32gvberbto43hgubXifwjhnrVSg7tkjfiu43qgt} and~\eqref{LM-PARZIMNWEI-2}
we deduce that
\begin{equation*}
d_{x_0}\,|\nabla u(x_0)|+
d_{x_0}^2\,|D^2 u(x_0)|\le C\,\Big(\|u\|_{L^\infty(\Omega)}+\|f\|_{C^\alpha_{2,\star}(\Omega)}\Big)
\end{equation*}
and therefore, taking the supremum over~$x_0\in\Omega$,
\begin{equation}\label{KS:M-263787987kAMScAKFo8jnds}
\|u\|_{C^{2}_\star(\Omega)}\le C\,\Big(\|u\|_{L^\infty(\Omega)}+\|f\|_{C^\alpha_{2,\star}(\Omega)}\Big)
.\end{equation}
Let now~$x_0\ne y_0\in\Omega$. 
We claim that
\begin{equation}\label{FGSBoNNSoEkAMScAKFo8jnds}
d_{x_0 y_0}^{2+\alpha}\,\frac{|D^2 u(x_0)-D^2 u(y_0)|}{|x_0-y_0|^\alpha}\le
C\,\Big(\|u\|_{L^\infty(\Omega)}+\|f\|_{C^\alpha_{2,\star}(\Omega)}\Big)
.\end{equation}
For this, up to exchanging these points, we may suppose that~$d_{y_0}\ge d_{x_0}=
d_{x_0y_0}$. Using the above notation about~$R_0$ and~$u_0$
we distinguish two cases. If~$|x_0-y_0|\ge\frac{R_0}4$,
we have that
\begin{eqnarray*}&& d_{x_0 y_0}^{2+\alpha}\,\frac{|D^2 u(x_0)-D^2 u(y_0)|}{|x_0-y_0|^\alpha}\le
d_{x_0 }^{2+\alpha}\,\frac{|D^2 u(x_0)|+|D^2 u(y_0)|}{(R_0/4)^\alpha}
\le\,C d_{x_0}^2\Big(|D^2 u(x_0)|+|D^2 u(y_0)|\Big)\\&&\qquad
\le\,C \Big(d_{x_0}^2|D^2 u(x_0)|+d_{y_0}^2|D^2 u(y_0)|\Big)\le C\,\|u\|_{C^{2}_\star(\Omega)}.
\end{eqnarray*}
Combining this with~\eqref{KS:M-263787987kAMScAKFo8jnds} we obtain~\eqref{FGSBoNNSoEkAMScAKFo8jnds}
in this case.

Thus, to complete the proof of~\eqref{FGSBoNNSoEkAMScAKFo8jnds} we can
now focus on the case in which~$|x_0-y_0|<\frac{R_0}4$.
In this situation, we set~$z_0:=\frac{y_0-x_0}{R_0}$
and note that~$|z_0|\le\frac{|y_0-x_0|}{R_0}<\frac14$. Hence we find that
\begin{eqnarray*}&&
d_{x_0 y_0}^{2+\alpha}\,\frac{|D^2 u(x_0)-D^2 u(y_0)|}{|x_0-y_0|^\alpha} 
=\frac{d_{x_0}^{2+\alpha}}{R_0^2}\,\frac{|D^2 u_0(0)-D^2 u_0(z_0)|}{|x_0-y_0|^\alpha} \le
C\,d_{x_0}^{\alpha}\,\frac{\|D^2u_0\|_{C^\alpha(B_{1/4})}|z_0|^\alpha}{|x_0-y_0|^\alpha}\\&&\qquad\qquad\le
C\,d_{x_0}^{\alpha}\,\frac{\|D^2u_0\|_{C^\alpha(B_{1/4})}}{R_0^\alpha}\le
C\,\|D^2u_0\|_{C^\alpha(B_{1/4})}.
\end{eqnarray*}
{F}rom this estimate and~\eqref{LM-PARZIMNWEI-2} we obtain that~\eqref{FGSBoNNSoEkAMScAKFo8jnds}
holds true in this case as well.

Hence, taking the supremum in~\eqref{FGSBoNNSoEkAMScAKFo8jnds} we conclude that
\[ [u]_{C^{2,\alpha}_\star(\Omega)}\le C\,\Big(\|u\|_{L^\infty(\Omega)}+\|f\|_{C^\alpha_{2,\star}(\Omega)}\Big).\]
The desired result then follows by combining this inequality with~\eqref{KS:M-263787987kAMScAKFo8jnds}.
\end{proof}

For the sake of completeness, as already remarked in footnote~\ref{qupodjh9320jd-29rtjgjjjnn3-x1023N}
on page~\pageref{qupodjh9320jd-29rtjgjjjnn3-x1023N}, we point out that
typical regularity results can also be stated with balls of different radii (but the different statements
turn out to be typically equivalent, possibly up to an appropriate covering argument).
Though we do not linger
too much on this technical detail, to highlight this situation we formulate
a variant of the interior estimates in Theorem~\ref{SCHAUDER-INTE} 
in balls of general radii:

\begin{corollary}\label{qupodjh9320jd-29rtjgjjjnn3-x1023} Let~$R>r>0$.
Let~$f\in C^\alpha(B_R)$ for some~$\alpha\in(0,1)$. Let~$u\in C^2(B_R)$ be a solution of
\begin{equation*} \Delta u=f\quad{\mbox{in }}\,B_R.\end{equation*}
Then, there exists~$C>0$, depending only on~$n$, $\alpha$, $R$ and~$r$, such that
$$ \|u\|_{C^{2,\alpha}(B_{r})}\le C\,\Big(\|u\|_{L^\infty(B_R)}+\|f\|_{C^\alpha(B_R)}\Big).$$
\end{corollary}

The proof of Corollary~\ref{qupodjh9320jd-29rtjgjjjnn3-x1023} leverages a covering
argument of general use which goes as follows:

\begin{figure}
  \centering
  \includegraphics[width=.4\linewidth]{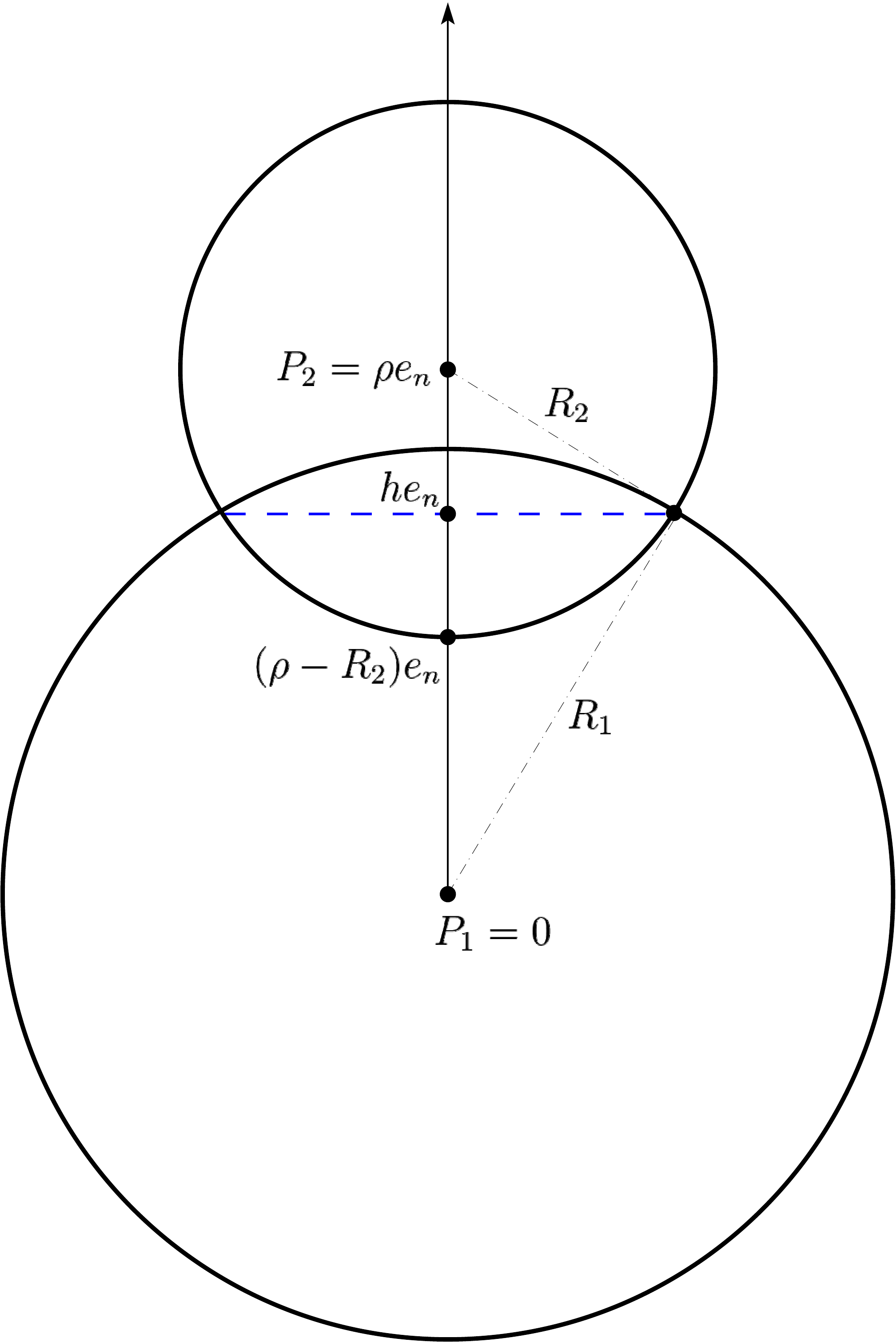}
 \caption{\sl The geometry appearing in the proof of
 Lemma~\ref{COVERINGARG} (1/3).}\label{T536EcdhPOaIJodidpo3tju24ylkiliR2543654yh2432RAe87658FI6tfh564}
\end{figure}

\begin{lemma}\label{COVERINGARG} Let~$k\in\N$ and~$\alpha\in(0,1]$. Let~$R_1$, $R_2>0$, $P_1$, $P_2\in\R^n$, with
\begin{equation}\label{2M32t343S9021hfdhghf9}
\rho:=|P_1-P_2|< R_1+R_2.\end{equation}
Assume that~$u\in C^{k,\alpha}(B_{R_1}(P_1))\cap C^{k,\alpha}(B_{R_2}(P_2))$. Then, $u\in C^{k,\alpha}(B_{R_1}(P_1)\cup B_{R_2}(P_2))$ and $$\|u\|_{C^{k,\alpha}(B_{R_1}(P_1)\cup B_{R_2}(P_2))}\le C\,\Big(
\|u\|_{C^{k,\alpha}(B_{R_1}(P_1))}+\|u\|_{C^{k,\alpha}( B_{R_2}(P_2))}\Big),$$
for some~$C>0$ depending only on~$n$, $R_1$, $R_2$ and~$\rho$.\end{lemma} 

\begin{proof} Certainly, $$\|u\|_{C^{k}(B_{R_1}(P_1)\cup B_{R_2}(P_2))}=
\max\big\{ \|u\|_{C^{k}(B_{R_1}(P_1))},\,\|u\|_{C^{k}(B_{R_2}(P_2))}\big\}.$$ Also, if~$j\in\{1,2\}$ and both~$x$ and~$y\in B_{R_j}(P_j)$ then~$|D^k u(x)-D^k u(y)|\le
C\|u\|_{C^{k,\alpha}(B_{R_j}(P_j))}|x-y|^\alpha$.

In view of these observations, to complete the proof of the desired result, we can focus on the case
\begin{equation}\label{1foeojuatio0UR4JHG048JT-4MD2}
x\in B_{R_1}(P_1)\setminus B_{R_2}(P_2)\qquad{\mbox{and}}\qquad y\in B_{R_2}(P_2)\setminus B_{R_1}(P_1).\end{equation} In this situation, we claim that there exist~$C>0$, depending only on~$n$, $R_1$, $R_2$ and~$\rho$, and~$z\in 
B_{R_1}(P_1)\cap B_{R_2}(P_2)$ such that
\begin{equation}\label{1foeojuatio0UR4JHG048JT-4MD}
|x-z|+|z-y|\le C|x-y|.
\end{equation}
To check this, 
up to swapping indices we can suppose that~$R_1\ge R_2$.
Besides, up to a translation we can suppose that~$P_1$ is the origin.
Also, up to a rotation, we can assume that~$P_2=|P_2|e_n$.
Therefore, $P_2=|P_1-P_2|e_n=\rho e_n$, see Figure~\ref{T536EcdhPOaIJodidpo3tju24ylkiliR2543654yh2432RAe87658FI6tfh564} (in particular,
the dependence on~$P_1$ and~$P_2$ can be dropped, and the constants~$C$
will depend naturally only on~$n$, $R_1$, $R_2$ and~$\rho$).
We also take~$h>0$ such that~$\partial B_{R_1}\cap \partial B_{R_2}(P_2)$ is contained in the
hyperplane~$\{x_n=h\}$ and we denote by
$${\mathcal{B}}_h:=\overline{B_{R_1}\cap B_{R_2}(P_2)} \cap\{ x_n=h\}.$$

\begin{figure}
  \centering
  \includegraphics[width=.4\linewidth]{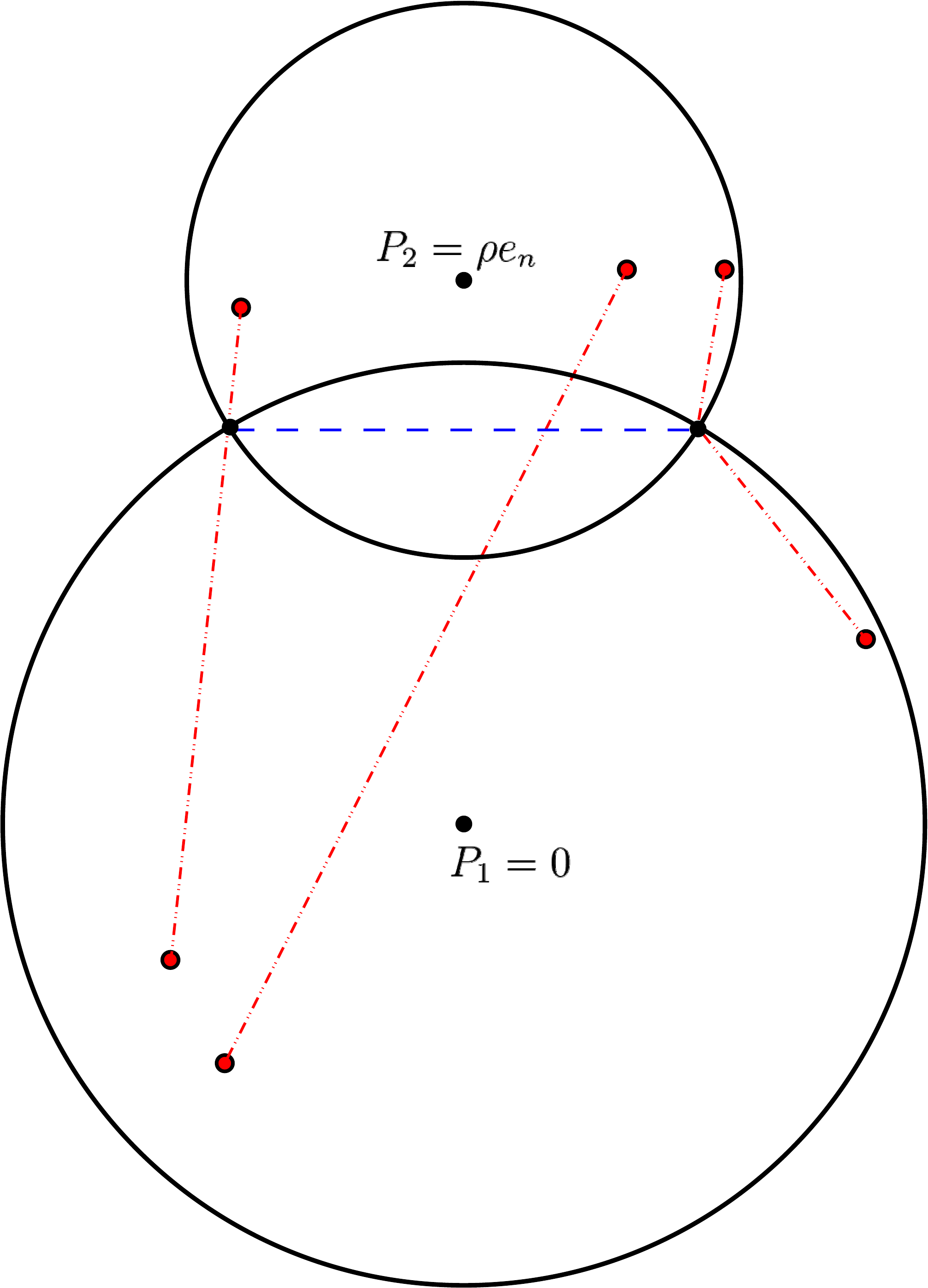}
 \caption{\sl The geometry appearing in the proof of
 Lemma~\ref{COVERINGARG} (2/3).}\label{T536EcdhPOaIJodidpo3tju24ylkiliR2543654yh2432RAe87658FI6tfh5642}
\end{figure}

In this setting, given~$x$ and~$y$
as in~\eqref{1foeojuatio0UR4JHG048JT-4MD2}, we consider the geodesics connecting~$x$ and~$y$
contained in~$\overline{B_{R_1}\cup B_{R_2}(P_2)}$. We remark that, in light
of~\eqref{1foeojuatio0UR4JHG048JT-4MD2}, the geodesics passes through~${\mathcal{B}}_h$
and we let~$\widetilde z$ be a point of the geodesics lying on~${\mathcal{B}}_h$.

We point out that
\begin{equation}\label{os34983764896yhkguhdkjb}\begin{split}&
{\mbox{the geodesics is either a segment passing through~$x$, $y$ and~$\widetilde z$,}}\\ &
{\mbox{or the union of two segments joining~$x$ to~$\widetilde z$ and~$\widetilde z$ to~$y$,}}\end{split}\end{equation}
see Figure~\ref{T536EcdhPOaIJodidpo3tju24ylkiliR2543654yh2432RAe87658FI6tfh5642}. To prove this, suppose that the geodesics is not a segment joining~$x$ and~$y$.
Then, we consider the geodesics joining~$x$ to~$\widetilde z$ and we observe that it has to be a segment due to
the convexity of the ball (recall that~$\widetilde z\in \overline{B_{R_1}}$). Similarly, 
the geodesics joining~$\widetilde z$ to~$y$ has to be a segment, since~$\widetilde z$, $y\in\overline{B_{R_2}(P_2)}$.
These considerations establish~\eqref{os34983764896yhkguhdkjb}.

In light of~\eqref{os34983764896yhkguhdkjb}, we consider the two situations where either
the geodesics is a segment or the union of two segments. 
In the first case, we have that~$|x-\widetilde z|+|y-\widetilde z|=|x-y|$. Hence, if~$\widetilde z\in \big(
B_{R_1}\cap B_{R_2}(P_2)\big) \cap\{ x_n=h\}$, then~\eqref{1foeojuatio0UR4JHG048JT-4MD} follows
by taking~$z:=\widetilde z$. If instead~$\widetilde z\in \partial B_{R_1}\cap \partial B_{R_2}(P_2)$,
then we observe that~$|x-\widetilde z|+|y-\widetilde z|<2|x-y|$, and
then, by continuity, we can find~$z\in \big(B_{R_1}\cap B_{R_2}(P_2)\big) \cap\{ x_n=h\}$
such that~$|x-z|+|z-y|\le\frac32|x-y|$, which gives~\eqref{1foeojuatio0UR4JHG048JT-4MD} in this case as well.

Hence we are left with the case in which
the geodesics is the union of two segments. 
In this case, we have that~$\widetilde z\in \partial B_{R_1}\cap \partial B_{R_2}(P_2)$ (indeed, if not, we could find
a shorter path connecting~$x$ and~$y$, thus contradicting the minimizing property of the geodesics, see Figure~\ref{T536EcdhPOaIJodidpo3tju24ylkiliR2543654yh2432RAe87658FI6tfh5643}).
We consider the angle~$\gamma$
between the vectors~$x-\widetilde z$ and~$y-\widetilde z$ and we observe that~$\gamma\in(\gamma_0,\pi]$,
for some~$\gamma_0\in(0,\pi)$. Then,
\begin{eqnarray*}
|x-y|^2&=&|x-\widetilde z|^2+|\widetilde z-y|^2 -2|x-\widetilde z|\,|\widetilde z-y|\,\cos\gamma\\
&\ge&
|x-\widetilde z|^2+|\widetilde z-y|^2 -2|x-\widetilde z|\,|\widetilde z-y|\,\cos\gamma_0\\
&\ge&
|x-\widetilde z|^2+|\widetilde z-y|^2 -\big(|x-\widetilde z|^2+|\widetilde z-y|^2\big)\cos\gamma_0
\\ &=& \big(1-\cos\gamma_0\big)\,\big(|x-\widetilde z|^2+|\widetilde z-y|^2 \big).
\end{eqnarray*}
As a consequence,
$$\frac{|x-y|}{\sqrt{1-\cos\gamma_0}}\ge \sqrt{|x-\widetilde z|^2+|\widetilde z-y|^2}\ge
\frac1{\sqrt 2} \big(|x-\widetilde z|+|\widetilde z-y|\big).
$$
Thus, by continuity, we can find~$z\in \big(B_{R_1}\cap B_{R_2}(P_2)\big) \cap\{ x_n=h\}$ such
that~\eqref{1foeojuatio0UR4JHG048JT-4MD} holds true for a suitable positive constant~$C$.
This completes the proof of~\eqref{1foeojuatio0UR4JHG048JT-4MD}.

\begin{figure}
  \centering
  \includegraphics[width=.4\linewidth]{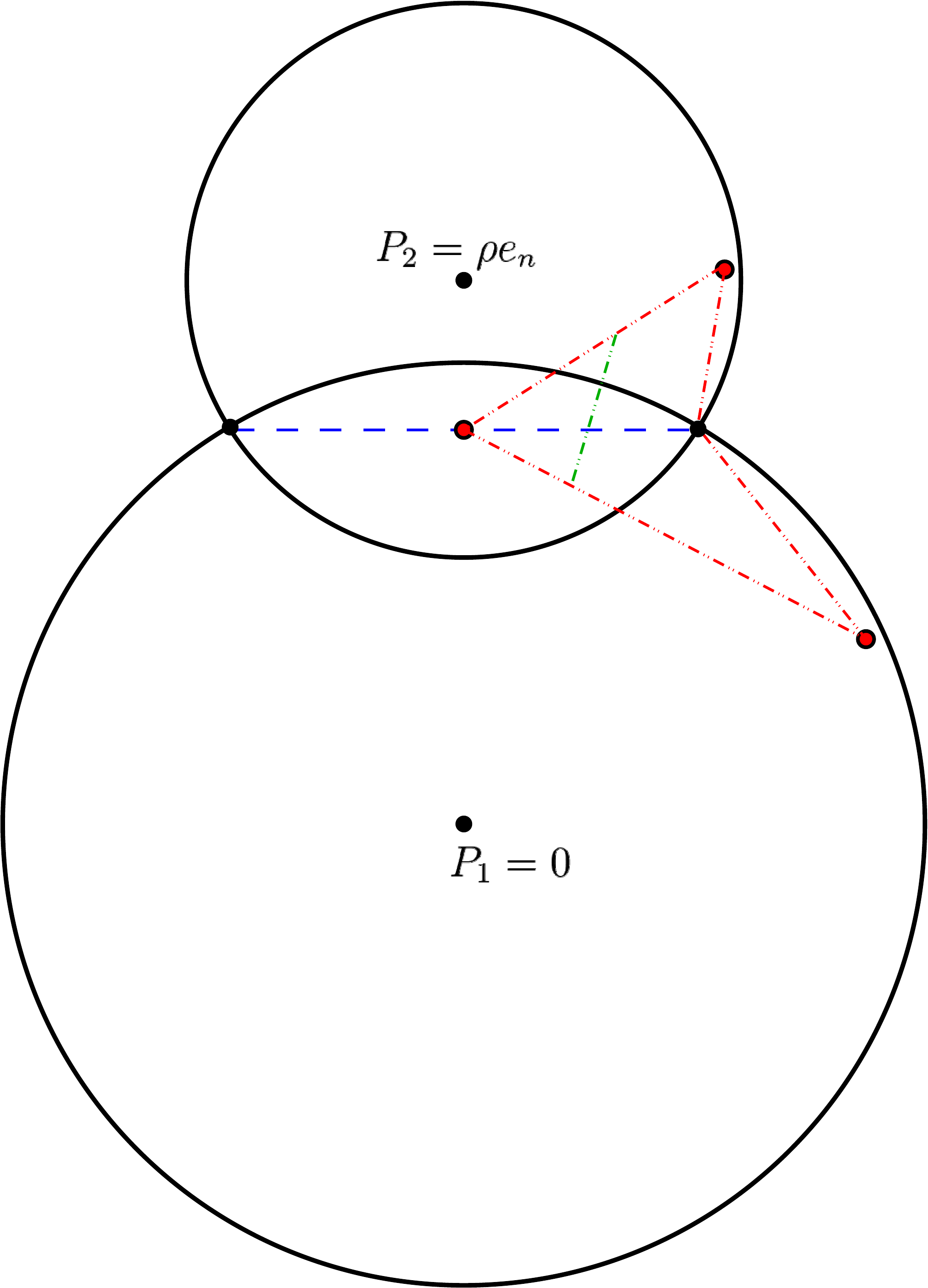}
 \caption{\sl The geometry appearing in the proof of
 Lemma~\ref{COVERINGARG} (3/3).}\label{T536EcdhPOaIJodidpo3tju24ylkiliR2543654yh2432RAe87658FI6tfh5643}
\end{figure}

{F}rom~\eqref{1foeojuatio0UR4JHG048JT-4MD} we deduce that
\begin{eqnarray*}&& |D^ku(x)-D^ku(y)|\le
|D^ku(x)-D^ku(z)|+|D^ku(z)-D^ku(y)|\\&&\qquad\le
C\|u\|_{C^{k,\alpha}(B_{R_1}(P_1))}|x-z|^\alpha
+
C\|u\|_{C^{k,\alpha}(B_{R_2}(P_2))}|z-y|^\alpha\\&&\qquad\le
C\,\Big(\|u\|_{C^{k,\alpha}(B_{R_1}(P_1))}
+
C\|u\|_{C^{k,\alpha}(B_{R_2}(P_2))}\Big)|x-y|^\alpha,\end{eqnarray*}
thus completing the proof of Lemma~\ref{COVERINGARG}.\end{proof}

With this, we can now proceed with the proof of Corollary~\ref{qupodjh9320jd-29rtjgjjjnn3-x1023}.

\begin{proof}[Proof of Corollary~\ref{qupodjh9320jd-29rtjgjjjnn3-x1023}]
We first observe that
\begin{equation}\label{TGB-Corollaryqupodjh9320jd-29rtjgjjjnn3-x1023}
{\mbox{it suffices to establish Corollary~\ref{qupodjh9320jd-29rtjgjjjnn3-x1023}
with~$R:=1$.}}
\end{equation}
To check this, for all~$x\in B_1$ let~$\widetilde u(x):=u(Rx)$ and~$\widetilde f(x):=R^2 f(Rx)$.
Let also~$\widetilde r:=\frac{r}{R}$.
Noticing that~$\Delta \widetilde u(x)=R^2 \Delta u(Rx)=R^2 f(Rx)=\widetilde f(x)$
for all~$x\in B_1$, if we know that Corollary~\ref{qupodjh9320jd-29rtjgjjjnn3-x1023} holds true with~$R=1$ we find that
$$ \|u\|_{C^{2,\alpha}(B_{r})}\le C\|\widetilde u\|_{C^{2,\alpha}(B_{\widetilde r})}
\le C\,\Big(\|\widetilde u\|_{L^\infty(B_1)}+\|\widetilde f\|_{C^\alpha(B_1)}\Big)
\le C\,\Big(\|u\|_{L^\infty(B_R)}+\|f\|_{C^\alpha(B_R)}\Big)
,$$
where~$C$ is a positive constant, possibly varying from step to step and
depending only on~$n$, $\alpha$, $R$ and~$r$.
This proves~\eqref{TGB-Corollaryqupodjh9320jd-29rtjgjjjnn3-x1023}.

Hence, by~\eqref{TGB-Corollaryqupodjh9320jd-29rtjgjjjnn3-x1023}, we can suppose that~$R:=1$.
Also, if~$r\le\frac12$ the desired result follows directly from Theorem~\ref{SCHAUDER-INTE},
therefore we can suppose that~$r\in\left(\frac12,1\right)$. In this setting, we observe that
Corollary~\ref{qupodjh9320jd-29rtjgjjjnn3-x1023} with~$R=1$
is proved once we show that for all~$j\in\N$
there exists~$C>0$, depending only on~$n$, $\alpha$ and~$j$, such that
\begin{equation}\label{TGB-Corollaryqupodjh9320jd-29rtjgjjjnn3-x1023MA2D}
\|u\|_{C^{2,\alpha}(B_{r_j})}\le C\,\Big(\|u\|_{L^\infty(B_1)}+\|f\|_{C^\alpha(B_1)}\Big),\qquad{\mbox{where }}\;r_j:=1-\frac1{2}\left(\frac34\right)^j .\end{equation}
Indeed, once~\eqref{TGB-Corollaryqupodjh9320jd-29rtjgjjjnn3-x1023MA2D}
is established, it is enough to notice that~$r_0=\frac12<r$ and~$r_j\to1>r$ as~$r\to+\infty$,
and then pick~$j\in\N\cap[1,+\infty)$ such that~$r_{j-1}<r\le r_{j}$: in this way, 
Corollary~\ref{qupodjh9320jd-29rtjgjjjnn3-x1023} with~$R=1$ would follow directly from~\eqref{TGB-Corollaryqupodjh9320jd-29rtjgjjjnn3-x1023MA2D}.

Thanks to this observation, we focus on the proof of~\eqref{TGB-Corollaryqupodjh9320jd-29rtjgjjjnn3-x1023MA2D}. To this end, we argue by induction over~$j$. When~$j=0$,
the claim in~\eqref{TGB-Corollaryqupodjh9320jd-29rtjgjjjnn3-x1023MA2D}
follows from Theorem~\ref{SCHAUDER-INTE}. To perform the inductive step,
we now suppose that the statement in~\eqref{TGB-Corollaryqupodjh9320jd-29rtjgjjjnn3-x1023MA2D}
holds true for an index~$j$ and we aim at showing its validity for the index~$j+1$.
For this, we pick a point~$p\in\partial B_{r_j}$
and for every~$x\in B_1$ we define~$\rho_j:=\frac{1-r_j}2$, $
u_j(x):=u(p+\rho_j x)$ and~$f_j(x):=\rho_j^2f(p+\rho_j x)$.
We remark that if~$x\in B_1$ then~$|p+\rho_j x|\le|p|+\rho_j\le r_j+\rho_j=\frac{1+r_j}2<1$ and therefore, for every~$x\in B_1$, $\Delta u_j(x)=\rho_j^2f(p+\rho_j x)=f_j(x)$. We can thereby employ
Theorem~\ref{SCHAUDER-INTE} and infer that
\begin{equation}\label{JSw82bigcup}\begin{split}&
 \|u\|_{C^{2,\alpha}(B_{\rho_j/2}(p))}\le C\| u_j\|_{C^{2,\alpha}(B_{1/2})}
\le C\,\Big(\| u_j\|_{L^\infty(B_1)}+\| f_j\|_{C^\alpha(B_1)}\Big)\\&\qquad\qquad
\le C\,\Big(\|u\|_{L^\infty(B_1)}+\|f\|_{C^\alpha(B_1)}\Big)
,\end{split}\end{equation}
where~$C$ is a positive constant, possibly varying from step to step and
depending only on~$n$, $\alpha$ and~$j$.

Now we note that
$$ \bigcup_{p\in\partial B_{r_j}} B_{\rho_j/2}(p)=
B_{r_j+\rho_j/2}\setminus B_{r_j-\rho_j/2}
$$
and thus we can find~$N\in\N$ depending only on~$n$ and~$j$ and~$p_1,\dots,p_N\in\partial B_{r_j}$
such that
$$ B_{r_j+\frac{\rho_j}2}\subseteq \left(
\bigcup_{\ell=1}^N B_{\rho_j/2}(p_\ell)\right)\cup B_{r_j}.$$
Thus, since
$$ r_j+\frac{\rho_j}2=r_j+\frac{1-r_j}4=\frac14+\frac34\,r_j=
\frac14+\frac34-\frac1{2}\left(\frac34\right)^{j+1}=r_{j+1},
$$
we have that
$$ B_{r_{j+1}}\subseteq \left(
\bigcup_{\ell=1}^N B_{\rho_j/2}(p_\ell)\right)\cup B_{r_j}.$$
We can therefore exploit~\eqref{JSw82bigcup} at the points~$p_1,\dots,p_N\in\partial B_{r_j}$
and combine it with the estimate in~$B_{r_j}$:
in this way, utilizing Lemma~\ref{COVERINGARG}, we conclude that
$$ \|u\|_{C^{2,\alpha}(B_{r_{j+1}})}\le C\,\Big(\|u\|_{L^\infty(B_1)}+\|f\|_{C^\alpha(B_1)}\Big),$$
up to renaming~$C$,
which concludes the inductive step.
This completes the proof of~\eqref{TGB-Corollaryqupodjh9320jd-29rtjgjjjnn3-x1023MA2D}.\end{proof}

For further reference we recall that a useful technical variation of Theorem~\ref{SCHAUDER-BOUTH} considers the case in which
the halfball~$B_1$ is replaced by the intersection of a ball and a halfspace (whose boundary does not
necessarily pass through the center of the ball). For this, we use the notation
\begin{equation}\label{9oj3falphaB10GENE} \begin{split}&
B_r^+(p):=\{x=(x',x_n)\in B_r(p) {\mbox{ s.t. }}x_n>0\}
\\{\mbox{and }}\;&B_r^0(p):=\{x=(x',x_n)\in B_r(p) {\mbox{ s.t. }}x_n=0\}\end{split}\end{equation}
and we have:

\begin{theorem}\label{SCHAUDER-BOUTH-VARIA}
Let~$p=(p',p_n)\in\R^{n-1}\times\R$ with~$p'=0$.
Let~$f\in C^\alpha(B_1^+(p))$ and~$g\in C^{2,\alpha}(B_1^0(p))$
for some~$\alpha\in(0,1)$. 
Let~$u\in C^2(B_1^+(p))\cap C(B_1^+(p)\cup B_1^0(p))$ be a solution of
\begin{equation*}\begin{dcases} \Delta u=f\quad{\mbox{in }}\,B_1^+(p),\\
u=g \quad{\mbox{on }}\,B_1^0(p).\end{dcases}\end{equation*}
Then, there exists~$C>0$, depending only on~$n$ and~$\alpha$, such that
$$ \|u\|_{C^{2,\alpha}(B_{1/2}^+(p))}\le C\,\Big(\|u\|_{L^\infty(B_1^+(p))}+\|g\|_{C^{2,\alpha}(B_1^0(p))}+\|f\|_{C^\alpha(B_1^+(p))}\Big).$$
\end{theorem}

\begin{proof} If~$p_n<-\frac12$ then~$B_{1/2}^+(p)=\varnothing$ and there is nothing to prove,
hence we can assume that
\begin{equation}\label{LS-21ihf0alphaB2p}
p_n\ge-\frac12.\end{equation}
Also, if~$p_n\ge\frac{3}4$, then for all~$x=(x',x_n)\in B_{3/4}(p)$ we have that~$x_n\ge p_n-|x_n-p_n|\ge\frac34-|x-p|>0$
and therefore~$B_{3/4}(p)\subseteq B_1^+(p)$. Thus, in this case we can utilize Corollary~\ref{qupodjh9320jd-29rtjgjjjnn3-x1023}
(centered at~$p$ instead of the origin) and deduce that
$$ \|u\|_{C^{2,\alpha}(B_{1/2}^+(p))}=\|u\|_{C^{2,\alpha}(B_{1/2}(p))}\le C\,\Big(\|u\|_{L^\infty(B_{3/4}(p))}+\|f\|_{C^\alpha(B_{3/4}(p))}\Big).$$
In view of this observation and~\eqref{LS-21ihf0alphaB2p},
we can focus on the case
\begin{equation*}
|p_n|\le\frac34.\end{equation*}
In this situation, we observe that
$$ \lim_{\varrho\searrow0} \left(\frac{3}{4}\right)^2+\left(\varrho+\frac14\right)^2
+\frac32\varrho=\frac{10}{16}<1,$$
hence we can fix a constant~$\varrho\in\left(0,\frac1{10}\right)$ such that
\[\left(\frac{3}{4}\right)^2+\left(\varrho+\frac14\right)^2
+\frac32\varrho<1.\]
We have that, for every~$q'\in\R^{n-1}$ with~$|q'|\le\frac{1}4$,
\begin{equation}\label{KS09-ry3begiatiobelLS-21ihffrac1n}
B_{\varrho}^+(q',0)\subseteq B_1^+(p).
\end{equation}
Indeed, if~$x\in B_{\varrho}^+(q',0)$ then~$x_n>0$ and
$$ |x|\le |x-(q',0)|+|(q',0)|<\varrho+|q'|\le\varrho+\frac14.$$
This gives that
$$ p_n^2+|x|^2<\left(\frac{3}{4}\right)^2+\left(\varrho+\frac14\right)^2$$
and therefore
\begin{eqnarray*}&& |x-p|^2=|x'|^2+(x_n-p_n)^2=|x'|^2+x_n^2+p_n^2-2x_np_n\le
|x|^2+p_n^2+\frac32|x_n|\\
&&\qquad\le|x|^2+p_n^2+\frac32|x-(q',0)|<\left(\frac{3}{4}\right)^2+\left(\varrho+\frac14\right)^2
+\frac32\varrho<1.
\end{eqnarray*}
This proves~\eqref{KS09-ry3begiatiobelLS-21ihffrac1n}.

As a consequence of~\eqref{KS09-ry3begiatiobelLS-21ihffrac1n}, we consider a finite family of
halfballs~$B_{\varrho}^+(q'_1,0),\dots,B_{\varrho}^+(q'_N,0)$, with~$|q'_j|\le\frac14$
for every~$j\in\{1,\dots,N\}$, with~$N$ universal, such that
$$ {\mathcal{S}}:=\left(-\frac14,\frac14\right)^{n-1}\times\left(0,\frac\varrho4\right)\subseteq
\bigcup_{i=1}^N B_{\varrho/2}^+(q'_i,0)$$
and we utilize Theorem~\ref{SCHAUDER-BOUTH} in each of these halfballs, finding that
$$ \|u\|_{C^{2,\alpha}(B_{\varrho/2}^+(q'_i,0))}
\le C\,\Big(\|u\|_{L^\infty(B_1^+(p))}+\|g\|_{C^{2,\alpha}(B_1^0(p))}+\|f\|_{C^\alpha(B_1^+(p))}\Big),$$
and thus, up to renaming constants,
\begin{equation} \label{OSJ-028r932hyfewgf87tgfbvcvdcew7u43trgfubewjTAgbsALSDVACAk}
\|u\|_{C^{2,\alpha}({\mathcal{S}})}
\le C\,\Big(\|u\|_{L^\infty(B_1^+(p))}+\|g\|_{C^{2,\alpha}(B_1^0(p))}+\|f\|_{C^\alpha(B_1^+(p))}\Big).\end{equation}
Similarly, one considers a finite family of balls~$B_{\varrho/9}(Q_1),\dots,B_{\varrho/8}(Q_M)$,
with~$M$ universal, such that
$$ {\mathcal{T}}:=\left\{x=(x',x_n)\in
B_{1-2\varrho}^+(p) {\mbox{ s.t. }}x_n>
\frac\varrho8 \right\}\subseteq
\bigcup_{i=1}^M B_{\varrho/9}(Q_i).$$
One can employ Corollary~\ref{qupodjh9320jd-29rtjgjjjnn3-x1023}
(centered at~$Q_i$ instead of the origin) in each of these balls and conclude that
$$ \|u\|_{C^{2,\alpha}(B_{\varrho/9}(Q_i))}
\le C\,\Big(\|u\|_{L^\infty(B_{\varrho/8}(Q_i))}+\|f\|_{C^\alpha(B_{\varrho/8}(Q_i))}\Big)$$
and therefore, up to renaming~$C$,
\begin{equation*} \|u\|_{C^{2,\alpha}({\mathcal{T}})}
\le C\,\Big(\|u\|_{L^\infty(B_1^+(p))}+\|g\|_{C^{2,\alpha}(B_1^0(p))}+\|f\|_{C^\alpha(B_1^+(p))}\Big).\end{equation*}
{F}rom this and~\eqref{OSJ-028r932hyfewgf87tgfbvcvdcew7u43trgfubewjTAgbsALSDVACAk}, up to renaming~$C$,
we arrive at
\begin{equation}\label{0o92oirjgmathcalT} \|u\|_{C^{2,\alpha}({\mathcal{S}}\cup{\mathcal{T}})}
\le C\,\Big(\|u\|_{L^\infty(B_1^+(p))}+\|g\|_{C^{2,\alpha}(B_1^0(p))}+\|f\|_{C^\alpha(B_1^+(p))}\Big).\end{equation}
We also have that
\begin{equation}\label{0o92oirjgmathcalT-1}
{\mathcal{S}}\cup{\mathcal{T}}\supseteq B_{1/4}^+(p).\end{equation}
Indeed, take~$x=(x',x_n)\in B_{1/4}^+(p)$. If~$x_n>\frac\varrho8$, we use that~$|x-p|<\frac14<1-2\varrho$
and thus~$x\in{\mathcal{T}}$. If instead~$x_n\in\left(0,\frac\varrho8\right]$ we have that, for every~$m\in\{1,\dots,n-1\}$,
$|x_m|\le|x'|=|x'-p'|<\frac14$, whence~$x\in{\mathcal{S}}$, which completes the proof of~\eqref{0o92oirjgmathcalT-1}.

By~\eqref{0o92oirjgmathcalT} and~\eqref{0o92oirjgmathcalT-1} we deduce that
$$ \|u\|_{C^{2,\alpha}(B_{1/4}^+(p)}
\le C\,\Big(\|u\|_{L^\infty(B_1^+(p))}+\|g\|_{C^{2,\alpha}(B_1^0(p))}+\|f\|_{C^\alpha(B_1^+(p))}\Big),$$
leading to the desired result up to scaling and covering.
\end{proof}

\section{Pointwise H\"older spaces}\label{KPJMDarJMSiaCLiMMSciAmoenbvD8Sijf-22}

It is often convenient to investigate the H\"older regularity theory\index{pointwise H\"older space}
in the light of the so-called ``pointwise H\"older spaces''
(see e.g.~\cite{MR1474098}, \cite[Section~9.2]{MR2590673}
and the references therein), that is, roughly speaking, to detect the H\"older regularity of a function
via the pointwise approximation of its Taylor polynomials.

\begin{definition}\label{JOSLi4o245nfIJSU-184ur032hv9qythgf92y98trpgJA4NBf}
Let~$\Omega\subseteq\R^n$ be open.
Let~$k\in\N$ and~$\alpha\in(0,1]$.
Let~$u:\Omega\to\R$ and~$x_0\in\Omega$. We say that~$u$ is~$C^{k,\alpha}$ at the point~$x_0$
in~$\Omega$
(and we write~$u\in C^{k,\alpha}_\Omega(x_0)$) if there exist~$M>0$ and a polynomial~$P$ of degree at most~$k$ such that
\begin{equation*} |u(x)-P(x)|\le M|x-x_0|^{k+\alpha}\qquad{\mbox{for all}}\,x\in\Omega.\end{equation*}
\end{definition}

\begin{figure}
  \centering
  \includegraphics[width=.6\linewidth]{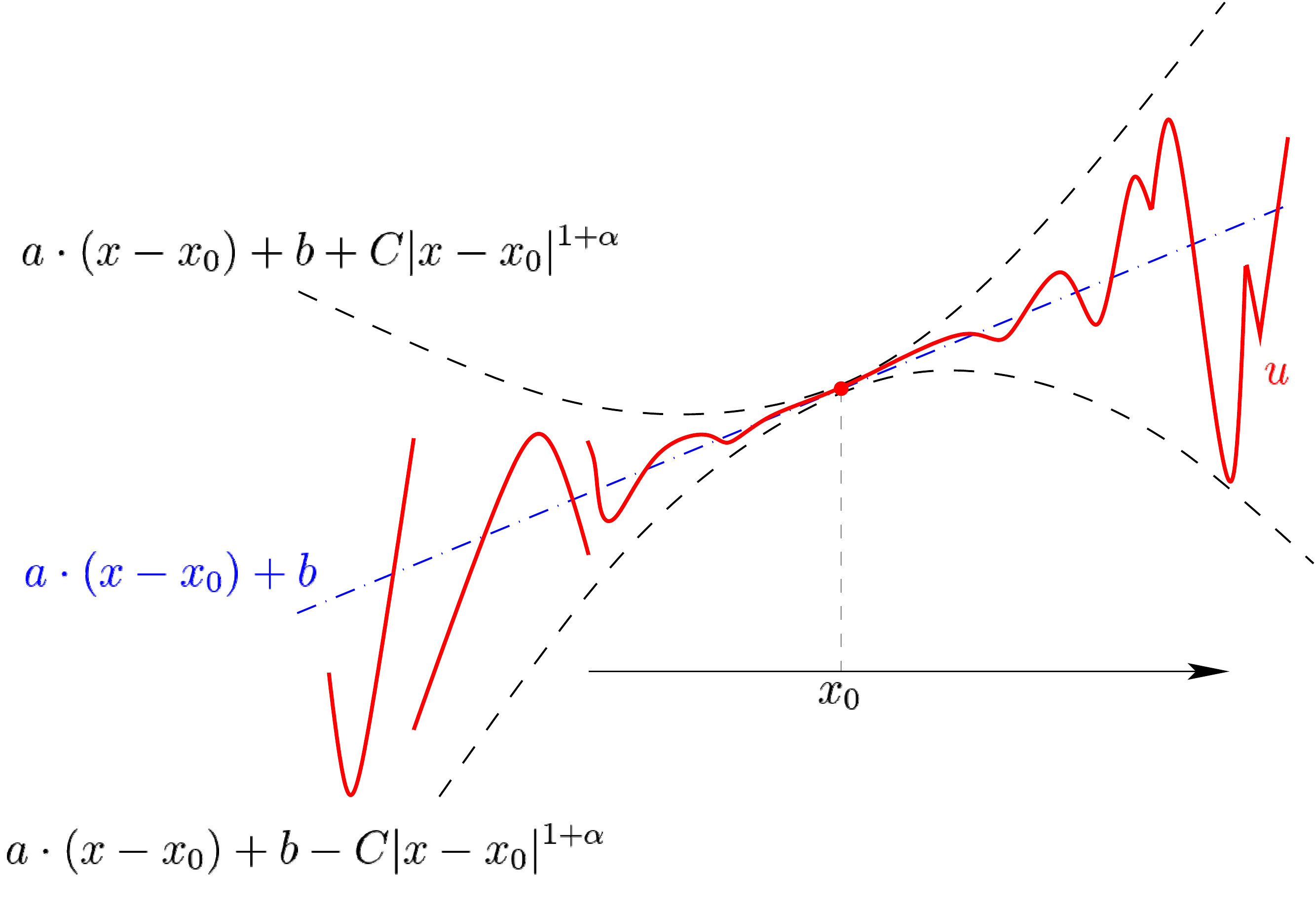}
 \caption{\sl A function in~$C^{1,\alpha}_\R(x_0)$.}\label{T536EcdhtjuylkiliR2543654yh2432RAe87658FI6tfh564}
\end{figure}

See Figure~\ref{T536EcdhtjuylkiliR2543654yh2432RAe87658FI6tfh564} for a pictorial representation of Definition~\ref{JOSLi4o245nfIJSU-184ur032hv9qythgf92y98trpgJA4NBf} when~$k=1$. When the constant~$M$ becomes uniform
within a Lipschitz domain, one recovers from Definition~\ref{JOSLi4o245nfIJSU-184ur032hv9qythgf92y98trpgJA4NBf}
the usual notion of H\"older spaces, according to the following result:

\begin{proposition}\label{XBOUN-AGG-P}
Let~$k\in\N$ and~$\alpha\in(0,1]$.
Let~$\Omega\subseteq\R^n$ be open, bounded and convex.
Let~$u:\Omega\to\R$.
Then, the following conditions are equivalent:
\begin{itemize}
\item[(i).] $u\in C^{k,\alpha}(\Omega)$,
\item[(ii).] There exists~$M\ge0$ such that for every~$x_0\in \Omega$ there exists
a polynomial~$P_{x_0}$ of degree at most~$k$
such that
\begin{equation}\label{BOUN-AGG-P}
|D^\gamma P_{x_0}(x_0)|\le M\qquad{\mbox{for all }}\,\gamma\in\N^n {\mbox{ with }} |\gamma|=m {\mbox{ and }}
m\in\{0,\dots,k\},\end{equation}
and
\begin{equation}\label{BOUN-AGG-P-2}
|u(x)-P_{x_0}(x)|\le M|x-x_0|^{k+\alpha}\qquad{\mbox{for all }}\,x\in\Omega.\end{equation}\end{itemize}
Additionally, if~(i) holds true then~$M$ in~(ii) can be taken of the form~$M:=C\,\|u\|_{C^{k,\alpha}(\Omega)}$
for some~$C>0$ depending only on~$n$, $k$, $\alpha$ and~$\Omega$.

Similarly, if~(ii) holds true, then~$\|u\|_{C^{k,\alpha}(\Omega)}\le CM$,
for some~$C>0$ depending only on~$n$, $k$, $\alpha$ and~$\Omega$.
\end{proposition}

\begin{proof} 
Assume that~(i) holds true. Given~$x_0\in\Omega$ we let~$P_{x_0}$ be the Taylor
polynomial of degree~$k$ of~$u$, namely
$$ P_{x_0}(x):=\sum_{{\beta\in\N^n}\atop{|\beta|\le k}}\frac{D^\beta u(x_0)}{\beta!}\,(x-x_0)^\beta.$$
In particular, for all~$m\in\N\cap[0,k]$ and~$\beta\in\N^n$ with~$|\beta|=m$,
\begin{equation}\label{OKSxgJSF0P890jRknsDFGHSOIdxYSHD} |D^\beta P_{x_0}(x_0)|=|D^\beta u(x_0)|\le 
\|u\|_{C^{k,\alpha}(\Omega)}.\end{equation}
Also, since~$\Omega$ is convex,
for every~$x\in\Omega$ there exists~$\xi(x)\in\Omega$ lying on the segment joining~$x$ and~$x_0$ such that
\begin{eqnarray*}
\left| u(x)-P_{x_0}(x)\right|&=&\left|
\sum_{{\beta\in\N^n}\atop{|\beta|= k}}\frac{D^\beta u(\xi(x))-D^\beta u(x_0)}{\beta!}\,(x-x_0)^\beta\right|\\
&\le&\sum_{{\beta\in\N^n}\atop{|\beta|= k}}\frac{|D^\beta u(\xi(x))-D^\beta u(x_0)|}{\beta!}\,|x-x_0|^k\\
&\le& \|u\|_{C^{k,\alpha}(\Omega)}\sum_{{\beta\in\N^n}\atop{|\beta|= k}}\frac{|\xi(x)-x_0|^\alpha}{\beta!}\,|x-x_0|^k\\
&\le& C\,\|u\|_{C^{k,\alpha}(\Omega)}\,|x-x_0|^{k+\alpha},
\end{eqnarray*}
which, together with~\eqref{OKSxgJSF0P890jRknsDFGHSOIdxYSHD}, gives~(ii).

Now, conversely, suppose that~(ii) is satisfied. We observe that
\begin{equation}\label{13432rfw2trfew5terghtry54yhgfjyhtrh254t5trw547h3hdhdrer-0-P09S}
\begin{split}&
{\mbox{for every polynomial~$P$ in~$\R^n$ of degree at most~$k$, we have that }}
\|P\|_{C^k(\Omega)}\le C\,\|P\|_{L^\infty(\Omega)},
\end{split}
\end{equation}
for some~$C>0$ depending only on~$n$, $k$ and~$\Omega$.
To check this, we let~$N$ be the number of
multi-indices~$\beta\in\N^n$ such that~$|\beta|\le k$. We consider the family~${\mathcal{P}}_k$
of polynomials of~$n$ variables with degree at most~$k$ and the map
$$ \R^N\ni a=\{a_\beta\}_{{\beta\in\N^n}\atop{|\beta|\le k}}\longmapsto
P_a(x):=\sum_{{\beta\in\N^n}\atop{|\beta|\le k}}a_\beta x^\beta\in {\mathcal{P}}_k.$$
We let~$\|a\|_1:=\|P_a\|_{C^{k}(\Omega)}$ and~$\|a\|_2:=\|P_a\|_{L^\infty(\Omega)}$.
We stress that both~$\|\cdot\|_1$ and~$\|\cdot\|_2$ are norms in~$\R^N$.
Since all norms in~$\R^N$ are equivalent (see e.g.~\cite[Theorem~15.26]{MR606198}), we find that
there exists~$C$ as above such that, for every~$a\in\R^N$, it holds that~$\|P_a\|_{C^{k}(\Omega)}=
\|a\|_1\le C\|a\|_2=C\|P_a\|_{L^\infty(\Omega)}$, from which we obtain~\eqref{13432rfw2trfew5terghtry54yhgfjyhtrh254t5trw547h3hdhdrer-0-P09S}.

We also claim that, for all~$p\in\R^n$ and~$r>0$, if~$P$ is a polynomial of degree at most~$k$ and~$j\in\N\cap[0,k]$,
then, for all~$\beta\in\N^n$ with~$|\beta|=j$,
\begin{equation}\label{13432rfw2trfew5terghtry54yhgfjyhtrh254t5trw547h3hdhdrer-0-P09SBIS}
\|D^\beta P\|_{L^\infty(B_r(p))}\le\frac{C}{r^j}\,\|{P}\|_{L^\infty(B_r(p))},
\end{equation}
for some~$C>0$ depending only on~$n$ and~$k$.
To check this, we let~$\widetilde{P}(x):=P(p+rx)$ and we make use of~\eqref{13432rfw2trfew5terghtry54yhgfjyhtrh254t5trw547h3hdhdrer-0-P09S} to find that
\begin{eqnarray*}&&
r^j\sup_{y\in B_r(p)}|D^\beta P(y)|=
r^j\sup_{x\in B_1}|D^\beta P(p+rx)|=
\sup_{x\in B_1}|D^\beta \widetilde{P}(x)|\\&&\qquad\le\|\widetilde{P}\|_{C^k(B_1)}\le
C\|\widetilde{P}\|_{L^\infty(B_1)}=C\|{P}\|_{L^\infty(B_r(p))},
\end{eqnarray*}
proving~\eqref{13432rfw2trfew5terghtry54yhgfjyhtrh254t5trw547h3hdhdrer-0-P09SBIS}.

Now we observe that for all~$x_0$, $\widetilde x_0\in\Omega$, setting~$p:=\frac{x_0+\widetilde x_0}{2}$
and~$r:=|x_0-\widetilde x_0|\in(0,1)$, if~$x\in B_r(p)$ we have that~$|x-x_0|\le|x-p|+|p-x_0|\le 2r$,
and similarly~$|x-\widetilde x_0|\le 2r$, giving that
$$ |P_{x_0}(x)-P_{\widetilde x_0}(x)|\le|P_{x_0}(x)-u(x)|+|u(x)-P_{\widetilde x_0}(x)|\le
M|x-x_0|^{k+\alpha}+M|x-\widetilde x_0|^{k+\alpha}\le
CMr^{k+\alpha}.$$
{F}rom this and~\eqref{13432rfw2trfew5terghtry54yhgfjyhtrh254t5trw547h3hdhdrer-0-P09SBIS} (used
here with~$P:=P_{x_0}-P_{\widetilde x_0}$) it follows that
for all~$x_0$, $\widetilde x_0\in\Omega$, if~$p:=\frac{x_0+\widetilde x_0}{2}$
and~$r:=|x_0-\widetilde x_0|\in(0,1)$, then, for all~$\beta\in\N^n$ with~$|\beta|=j$,
\begin{equation}\label{9027yr98ytr832gv87gbv7t37t3rug874bv875hgfu4eytccn9843n8yhhoerhghgiohi8y43y80}
\|D^\beta P_{x_0}-D^\beta P_{\widetilde x_0}\|_{L^\infty(B_r(p))}\le CMr^{k-j+\alpha}.
\end{equation}
Now we prove that
\begin{equation}\label{CLAIM-POL-0ok4r3e-jX}
\begin{split}&
{\mbox{$u\in C^{k,\alpha}(\Omega)$ and, for every~$x_0\in\Omega$
and~$\beta\in\N^n$ with~$|\beta|\le k$,}} \\&{\mbox{we have that }}D^\beta u(x_0)= D^\beta P_{x_0}(x_0).\end{split}
\end{equation}
This will be accomplished by checking that, for all~$j\in\N\cap[0,k]$,
\begin{equation}\label{CLAIM-POL-0ok4r3e-j}
\begin{split}&
{\mbox{$u\in C^{j,\alpha}(\Omega)$ and, for every~$x_0\in\Omega$
and~$\beta\in\N^n$ with~$|\beta|\le j$,}} \\&{\mbox{we have that }}D^\beta u(x_0)= D^\beta P_{x_0}(x_0).\end{split}
\end{equation}
We prove this claim by induction over~$j$.
For this, we first observe that, for every~$x_0\in\Omega$, condition~(ii) entails that
\begin{eqnarray*}&& \lim_{x\to x_0} |u(x)-P_{x_0}(x_0)|\le
\lim_{x\to x_0}\Big( |u(x)-P_{x_0}(x)|+|P_{x_0}(x)-P_{x_0}(x_0)|\Big)\\&&\qquad
=\lim_{x\to x_0} |u(x)-P_{x_0}(x)|\le\lim_{x\to x_0}
M|x-x_0|^{k+\alpha}=0.\end{eqnarray*}
This says that~$u$ is continuous at~$x_0$ and~$u(x_0)=P_{x_0}(x_0)$.
Furthermore, writing
\begin{equation}\label{POLIGRA018urhfi} P_{x_0}(x)=\sum_{{\beta\in\N^n}\atop{|\beta|\le k}}\frac{D^\beta P_{x_0}(x_0)}{\beta!}\,(x-x_0)^\beta,\end{equation}
if~$\delta\in(0,1)$ and~$e\in\partial B_1$ we have that
$$ |P_{x_0}(x_0+\delta e)-P_{x_0}(x_0)|=
\left|\sum_{{\beta\in\N^n}\atop{|\beta|\le k}}\frac{D^\beta P_{x_0}(x_0)}{\beta!}\,(\delta e)^\beta-P_{x_0}(x_0)
\right|=
\left|
\sum_{{\beta\in\N^n}\atop{1\le|\beta|\le k}}\frac{\delta^{|\beta|}D^\beta P_{x_0}(x_0)}{\beta!}\,e^\beta
\right|\le CM\delta,$$
thanks to~\eqref{BOUN-AGG-P}.

For this reason,
\begin{eqnarray*}&& |u(x_0+\delta e)-u(x_0)|=|u(x_0+\delta e)-P_{x_0}(x_0)|\\&&\qquad\le
|u(x_0+\delta e)-P_{x_0}(x_0+\delta e)|+
|P_{x_0}(x_0+\delta e)-P_{x_0}(x_0)|\le M\delta^{k+\alpha}+CM\delta\le CM\delta^\alpha
\end{eqnarray*}
and therefore~$u\in C^\alpha(\Omega)$.
This is the desired claim in~\eqref{CLAIM-POL-0ok4r3e-j} when~$j=0$.

Now, if~$k\ge1$ and~$i\in\{1,\dots,n\}$
\begin{eqnarray*}&&
\lim_{\e\to0} \left|
\frac{u(x_0+\e e_i)-u(x_0)}{\e}-\partial_i P_{x_0}(x_0)
\right|\\& \le&\lim_{\e\to0}\left( \left|
\frac{u(x_0+\e e_i)-u(x_0)}{\e}-\frac{P_{x_0}(x_0+\e e_i)-P_{x_0}(x_0)}{\e}
\right|
+
\left|
\frac{P_{x_0}(x_0+\e e_i)-P_{x_0}(x_0)}{\e}-\partial_i P_{x_0}(x_0)
\right|
\right)\\&=&\lim_{\e\to0}\left|
\frac{u(x_0+\e e_i)-u(x_0)}{\e}-\frac{P_{x_0}(x_0+\e e_i)-P_{x_0}(x_0)}{\e}
\right|\\&=&\lim_{\e\to0}\left|
\frac{u(x_0+\e e_i)- P_{x_0}(x_0+\e e_i) }{\e}
\right|
\\&\le&\lim_{\e\to0}\frac{M|\e|^{k+\alpha}}{|\e|}\\
& =&0.
\end{eqnarray*}
This gives that~$u$ can be differentiated at least once at~$x_0$ and~$\nabla u(x_0)=\nabla P_{x_0}(x_0)$.
Additionally, differentiating~\eqref{POLIGRA018urhfi}, for all~$i\in\{1,\dots,n\}$,
\begin{equation*}\partial_i P_{x_0}(x)=
\sum_{{\beta\in\N^n}\atop{1\le|\beta|\le k}}\frac{\beta_i\,D^\beta P_{x_0}(x_0)}{\beta!}\,(x-x_0)^{\beta-e_i}.\end{equation*}
Hence, recalling~\eqref{BOUN-AGG-P}, if~$\delta\in(0,1)$ and~$e\in\partial B_1$,
\begin{equation}\begin{split}\label{4.2.6BIS}&
|\partial_i P_{x_0}(x_0+\delta e)-\partial_i P_{x_0}(x_0)|=
\left|\sum_{{\beta\in\N^n}\atop{1\le|\beta|\le k}}\frac{\beta_i\,D^\beta P_{x_0}(x_0)}{\beta!}\,(\delta e)^{\beta-e_i}
-\partial_i P_{x_0}(x_0)
\right|\\&\qquad
=\left| \sum_{{\beta\in\N^n}\atop{2\le|\beta|\le k}} \frac{\beta_i\,D^\beta P_{x_0}(x_0)}{\beta!}\,(\delta e)^{\beta-e_i}
+\sum_{j=1}^{n} \delta_{ij}\partial_j P_{x_0}(x_0)\,(\delta e)^{e_j-e_i}
-\partial_i P_{x_0}(x_0)
\right|\\&\qquad=
\left|
\sum_{2\le|\beta|\le k}\frac{\delta^{|\beta|-1}\beta_i\,D^\beta P_{x_0}(x_0)}{\beta!}\,e^{\beta-e_i}
\right|\le CM\delta.
\end{split}
\end{equation}
Moreover,
\begin{eqnarray*}
&&\left| \frac{u(x_0+\delta e+\e e_i)-u(x_0+\delta e)}{\e}-\frac{u(x_0+\e e_i)-u(x_0)}{\e}\right|
\\&\le&
\left| \frac{P_{x_0+\delta e}(x_0+\delta e+\e e_i)-P_{x_0+\delta e}(x_0+\delta e)}{\e}-\frac{P_{x_0}(x_0+\e e_i)-P_{x_0}(x_0)}{\e}\right|
\\&&\quad+\frac1{|\e|}\Big(|u(x_0+\delta e+\e e_i)-P_{x_0+\delta e}(x_0+\delta e+\e e_i)|+
|u(x_0+\delta e)-P_{x_0+\delta e}(x_0+\delta e)|\\&&\qquad\qquad+
|u(x_0+\e e_i)-P_{x_0}(x_0+\e e_i)|+
|u(x_0)-P_{x_0}(x_0)|\Big)
\\&\le&
\left| \frac{P_{x_0+\delta e}(x_0+\delta e+\e e_i)-P_{x_0+\delta e}(x_0+\delta e)}{\e}-\frac{P_{x_0}(x_0+\e e_i)-P_{x_0}(x_0)}{\e}\right|
\\&&\quad+\frac1{|\e|}\Big(M|\e|^{k+\alpha}+0+M|\e|^{k+\alpha}+0\Big)
,\end{eqnarray*}
leading to, as~$\e\to0$,
$$ |\partial_i u(x_0+\delta e)-\partial_i u(x_0)|\le
|\partial_i P_{x_0+\delta}(x_0+\delta e)-\partial_i P_{x_0}(x_0)|.$$
This, together with~\eqref{9027yr98ytr832gv87gbv7t37t3rug874bv875hgfu4eytccn9843n8yhhoerhghgiohi8y43y80}
and~\eqref{4.2.6BIS},
gives that, if~$x_0\in\Omega$, $e\in\partial B_1$ and~$\delta\in(0,1)$ is sufficiently small such that~$x_0+\delta e\in\Omega$,
then
\begin{eqnarray*}&&
|\partial_i u(x_0+\delta e)-\partial_i u(x_0)|\le
|\partial_i P_{x_0+\delta}(x_0+\delta e)-\partial_i P_{x_0}(x_0+\delta e)|+
|\partial_i P_{x_0}(x_0+\delta e)-\partial_i P_{x_0}(x_0)|\\&&\qquad\qquad\le
CM\delta^{k-1+\alpha}+CM\delta\le CM\delta^{\alpha}.\end{eqnarray*}
This gives that~$\partial_i u\in C^\alpha(\Omega)$ and completes the proof
of~\eqref{CLAIM-POL-0ok4r3e-j} when~$j=1$.

Suppose now recursively that~$j\in\N\cap[0,k-1]$
and~$u\in C^{j,\alpha}(\Omega)$ is such that at
each point~$x_0\in\Omega$ we have~$D^\beta u(x_0)= D^\beta P_{x_0}(x_0) $
as long as~$|\beta|\le j$. Then, recalling~\eqref{9027yr98ytr832gv87gbv7t37t3rug874bv875hgfu4eytccn9843n8yhhoerhghgiohi8y43y80},
for all~$i\in\{1,\dots,n\}$,
\begin{eqnarray*}
&&\lim_{\e\to0} \left|
\frac{D^\beta u(x_0+\e e_i)-D^\beta u(x_0)}{\e}- D^{\beta+e_i} P_{x_0}(x_0)
\right|\\&\le&
\lim_{\e\to0}\left( \left|
\frac{D^\beta u(x_0+\e e_i)-D^\beta u(x_0)}{\e}-\frac{D^\beta P_{x_0}(x_0+\e e_i)-D^\beta P_{x_0}(x_0)}{\e}
\right|\right.\\&&\qquad\qquad\left.+ 
\left|\frac{D^\beta P_{x_0}(x_0+\e e_i)-D^\beta P_{x_0}(x_0)}{\e}-{D^{\beta+e_i} P_{x_0}(x_0)}\right|
\right)\\&=&\lim_{\e\to0} \left|
\frac{D^\beta u(x_0+\e e_i)-D^\beta u(x_0)}{\e}-\frac{D^\beta P_{x_0}(x_0+\e e_i)-D^\beta P_{x_0}(x_0)}{\e}
\right|\\&=&\lim_{\e\to0} \left|
\frac{D^\beta P_{x_0+\e e_i}(x_0+\e e_i)-D^\beta P_{x_0}(x_0)}{\e }-\frac{D^\beta P_{x_0}(x_0+\e e_i)-D^\beta P_{x_0}(x_0)}{\e}
\right|\\&=&\lim_{\e\to0} \left|
\frac{D^\beta P_{x_0+\e e_i}(x_0+\e e_i)- D^\beta P_{x_0}(x_0+\e e_i) }{\e}
\right|\\&\le&\lim_{\e\to0}CM|\e|^{k-j-1+\alpha}\\&\le&\lim_{\e\to0}CM|\e|^{\alpha}\\&=&0.
\end{eqnarray*}
This gives that~$D^{\beta+e_i}u(x_0)=D^{\beta+e_i}P_{x_0}(x_0)$.

Furthermore, from~\eqref{POLIGRA018urhfi} we see that, for every~$\gamma\in\N^n$,
\begin{eqnarray*}
D^\gamma P_{x_0}(x)=\sum_{{{\beta\in\N^n}\atop{|\beta|\le k}}\atop{\gamma\le\beta}}
\prod_{{1\le j\le n}\atop{\gamma_j\ge1}}\beta_j(\beta_j-1)\dots(\beta_j-\gamma_j+1)
\frac{D^\beta P_{x_0}(x_0)}{\beta!}\,(x-x_0)^{\beta-\gamma}.
\end{eqnarray*}
Here above we use the notation that~$\gamma=(\gamma_i)_{1\le i\le n}$, $\beta=(\beta_i)_{1\le i\le n}$
and we say that~$\gamma\le\beta$ if~$\gamma_i\le\beta_i$ for every~$i\in\{1,\dots,n\}$.

Accordingly, if~$x_0\in\Omega$, $e\in\partial B_1$ and~$\delta\in(0,1)$,
\begin{equation}\label{qw25454r6QQQQQQQQ}\begin{split}&
\left| D^\gamma P_{x_0}(x_0+\delta e)-D^\gamma P_{x_0}(x_0)
\right|\\=\;&\left|\sum_{{{\beta\in\N^n}\atop{|\beta|\le k}}\atop{\gamma\le\beta}}
\prod_{{1\le j\le n}\atop{\gamma_j\ge1}}\beta_j(\beta_j-1)\dots(\beta_j-\gamma_j+1)
\frac{D^\beta P_{x_0}(x_0)}{\beta!}\,(\delta e)^{\beta-\gamma}-D^\gamma P_{x_0}(x_0)
\right|\\=\;&
\left|\sum_{{{\beta\in\N^n}\atop{|\gamma|+1\le|\beta|\le k}}\atop{\gamma\le\beta} }
\prod_{{1\le j\le n}\atop{\gamma_j\ge1}}\beta_j(\beta_j-1)\dots(\beta_j-\gamma_j+1)
\frac{D^\beta P_{x_0}(x_0)}{\beta!}\,(\delta e)^{\beta-\gamma}\right. \\&\qquad\left.
+\sum_{{{\beta\in\N^n}\atop{|\gamma|=|\beta|\le k}}\atop{\gamma\le\beta} }
\prod_{{1\le j\le n}\atop{\gamma_j\ge1}}\beta_j(\beta_j-1)\dots(\beta_j-\gamma_j+1)
\frac{D^\beta P_{x_0}(x_0)}{\beta!}\,(\delta e)^{\beta-\gamma}
-D^\gamma P_{x_0}(x_0)
\right|.
\end{split}\end{equation}
Now we notice that if~$|\gamma|=|\beta|$ and~$\gamma\le\beta$ then~$\beta_i=\gamma_i$ for all~$i\in\{1,\dots,n\}$,
and therefore~$\beta=\gamma$ and
$$\prod_{{1\le j\le n}\atop{\gamma_j\ge1}}\beta_j(\beta_j-1)\dots(\beta_j-\gamma_j+1)=\beta!  .$$
This gives that
$$ \sum_{{{\beta\in\N^n}\atop{|\gamma|=|\beta|\le k}}\atop{\gamma\le\beta} }
\prod_{{1\le j\le n}\atop{\gamma_j\ge1}}\beta_j(\beta_j-1)\dots(\beta_j-\gamma_j+1)
\frac{D^\beta P_{x_0}(x_0)}{\beta!}\,(\delta e)^{\beta-\gamma}= D^\gamma P_{x_0}(x_0).$$
Plugging this information into~\eqref{qw25454r6QQQQQQQQ} and
recalling~\eqref{BOUN-AGG-P}, we conclude that
\begin{equation}\label{aggh234poiutsgfkhljohutyryrtg}\begin{split}&
\left| D^\gamma P_{x_0}(x_0+\delta e)-D^\gamma P_{x_0}(x_0)
\right|\\&\qquad
= \left|\sum_{{{\beta\in\N^n}\atop{|\gamma|+1\le|\beta|\le k}}\atop{\gamma\le\beta} }
\prod_{{1\le j\le n}\atop{\gamma_j\ge1}}\beta_j(\beta_j-1)\dots(\beta_j-\gamma_j+1)
\frac{|\delta|^{|\beta|-|\gamma|}\,D^\beta P_{x_0}(x_0)}{\beta!}\, e^{\beta-\gamma}\right|
\le CM\delta.\end{split}
\end{equation}

Consequently,
using~\eqref{aggh234poiutsgfkhljohutyryrtg} with~$\gamma:=\beta+e_i$
and~\eqref{9027yr98ytr832gv87gbv7t37t3rug874bv875hgfu4eytccn9843n8yhhoerhghgiohi8y43y80},
if~$x_0\in\Omega$, $e\in\partial B_1$ and~$\delta\in(0,1)$ is sufficiently small such that~$x_0+\delta e\in\Omega$,
then
\begin{equation}\label{OKMSXJHS33mFHSO70RdgyO0x}\begin{split}
&|D^{\beta+e_i}u(x_0)-D^{\beta+e_i}u(x_0+\delta e)|=
|D^{\beta+e_i}P_{x_0}(x_0)
-D^{\beta+e_i}P_{x_0+\delta e}(x_0+\delta e)|\\&\qquad\le
|D^{\beta+e_i}P_{x_0}(x_0)
-D^{\beta+e_i}P_{x_0}(x_0+\delta e)|+
|D^{\beta+e_i}P_{x_0}(x_0+\delta e)
-D^{\beta+e_i}P_{x_0+\delta e}(x_0+\delta e)|\\&\qquad\le CM\delta+CM\delta^{k-j-1+\alpha}\le CM\delta^\alpha.
\end{split}\end{equation}
Accordingly, we have that~$D^{\beta+e_i}u\in C^\alpha(\Omega)$, which
completes the inductive step. The claim in~\eqref{CLAIM-POL-0ok4r3e-j}
is thereby established, and from this we obtain also~\eqref{CLAIM-POL-0ok4r3e-jX}.

{F}rom~\eqref{BOUN-AGG-P}
and~\eqref{CLAIM-POL-0ok4r3e-jX} it also follows that
$$ \sup_{{x_0\in\Omega}\atop{|\gamma|\le k}}|D^\gamma u(x_0)|= \sup_{{x_0\in\Omega}\atop{|\gamma|\le k}}
|D^\gamma P_{x_0}(x_0)|\le M.$$
Additionally, by~\eqref{OKMSXJHS33mFHSO70RdgyO0x},
we see that~$|D^{\gamma}u(x_0)-D^{\gamma}u(x_0+\delta e)|\le CM\delta^\alpha$ for all~$\gamma\in\N^n$
such that~$|\gamma|\le k$,
and therefore~$\|u\|_{C^{k,\alpha}(\Omega)}\le CM$. This gives that~(i)
holds true, as desired.
\end{proof}

It is worth pointing out that, in view of Definition~\ref{JOSLi4o245nfIJSU-184ur032hv9qythgf92y98trpgJA4NBf}, condition~\eqref{BOUN-AGG-P-2}
gives\footnote{Condition~\eqref{BOUN-AGG-P} is instead mostly technical
and can also be obtained from~\eqref{BOUN-AGG-P-2} (if valid at every point) if one knows~$P_{\widetilde x}$
for a point~$\widetilde x\in\Omega$: indeed,
it follows from~\eqref{BOUN-AGG-P-2} and the boundedness of~$\Omega$
that
$$\|u\|_{L^\infty(\Omega)}=\sup_{x\in\Omega}|u(x)|\le
\sup_{x\in\Omega}\Big( |P_{\widetilde x}(x)|+M|x-\widetilde x|^{k+\alpha}\Big)\le
\|P_{\widetilde x}\|_{L^\infty(\Omega)}+CM,$$
whence
$$ \|P_{x_0}\|_{L^\infty(\Omega)}=\sup_{x\in\Omega}|P_{x_0}(x)|\le
\sup_{x\in\Omega}\Big( |u(x)|+M|x-x_0|^{k+\alpha}\Big)
\le
\|u\|_{L^\infty(\Omega)}+CM\le
C\Big(
\|P_{\widetilde x}\|_{L^\infty(\Omega)}+M
\Big).$$
{F}rom this and the finite dimensional structure of the polynomials
(see e.g.~\eqref{13432rfw2trfew5terghtry54yhgfjyhtrh254t5trw547h3hdhdrer-0-P09S})
one would obtain that~$|D^mP_{x_0}(x_0)|\le C\Big(
\|P_{\widetilde x}\|_{L^\infty(\Omega)}+M
\Big)$ and therefore a uniform information
on~\eqref{BOUN-AGG-P-2} would entail a uniform information on~\eqref{BOUN-AGG-P}.

The small catch of the above argument is that~$M$ would be replaced, up to constants,
by~$\|P_{\widetilde x}\|_{L^\infty(\Omega)}+M$, which is a quantity not immediately
related to the original function~$u$. This small inconvenience can be fixed,
for instance, if one is willing to work with bounded functions:
this idea will be clarified in the forthcoming
Lemma~\ref{POLIAPP-J-GEO}.}
that~$u\in C^{k,\alpha}_\Omega(x_0)$ for every~$x_0\in\Omega$, with a uniform~$M$.

Without uniformity in~$M$, one cannot precisely reconstruct~$C^{k,\alpha}(\Omega)$
out of~$C^{k,\alpha}_\Omega(x_0)$ for~$x_0\in\Omega$, as detailed in the following result:

\begin{lemma}
Let~$k\in\N$ and~$\alpha\in(0,1]$.
Let~$\Omega\subseteq\R^n$ be open, bounded and convex.
We have that
\begin{equation}\label{9087r98732ytfiu874d6h743657bLHJSFTMAtoIYo02}
C^{k,\alpha}(\Omega)\subseteq \bigcup_{x_0\in\Omega}C^{k,\alpha}_\Omega(x_0) \end{equation}
and the above inclusion is strict.
\end{lemma}

\begin{proof}
If~$u\in C^{k,\alpha}(\Omega)$, then condition~(i) in Proposition~\ref{XBOUN-AGG-P}
holds true, leading to~\eqref{BOUN-AGG-P-2}, whence~$u\in C^{k,\alpha}_\Omega(x_0)$ for all~$x_0\in\Omega$,
in light of Definition~\ref{JOSLi4o245nfIJSU-184ur032hv9qythgf92y98trpgJA4NBf}.

This shows the desired inclusion. To check that the inclusion is strict, we 
consider a sequence of mutually disjoint intervals~$I_j\subset(0,+\infty)$
of the form~$I_j=[p_j-\ell_j,p_j+\ell_j]$, with~$p_j\in(0,1)$ and~$p_j\searrow0$ as~$j\to+\infty$.
Possibly shrinking~$I_j$ we can additionally suppose that~$\ell_j\le
\left(\frac{p_j}2\right)^{\frac{k+\alpha}k}$ and in particular~$\ell_j\le\frac{p_j}2$.
Therefore if~$x_1\in I_j$ then
\begin{equation}\label{w3tg-inoiqdjQunAM34S-a}
|x_1|\ge x_1\ge p_j-\ell_j\ge \frac{p_j}2\ge\ell_j^{\frac{k}{k+\alpha}}=\frac{\ell_j}{{\ell_j^{\frac\alpha{k+\alpha}}}}
\ge\frac{|x_1-p_j|}{\ell_j^{\frac\alpha{k+\alpha}}}.
\end{equation}
Let also~$s_j:=\frac{1}{\ell_j^\alpha}$ and let~$\eta_j\in C^\infty_0(I_j,\,[0,1])$
be such that~$\eta=1$ in~$\left[p_j-\frac{\ell_j}2,p_j+\frac{\ell_j}2\right]$.

\begin{figure}
  \centering
  \includegraphics[width=.6\linewidth]{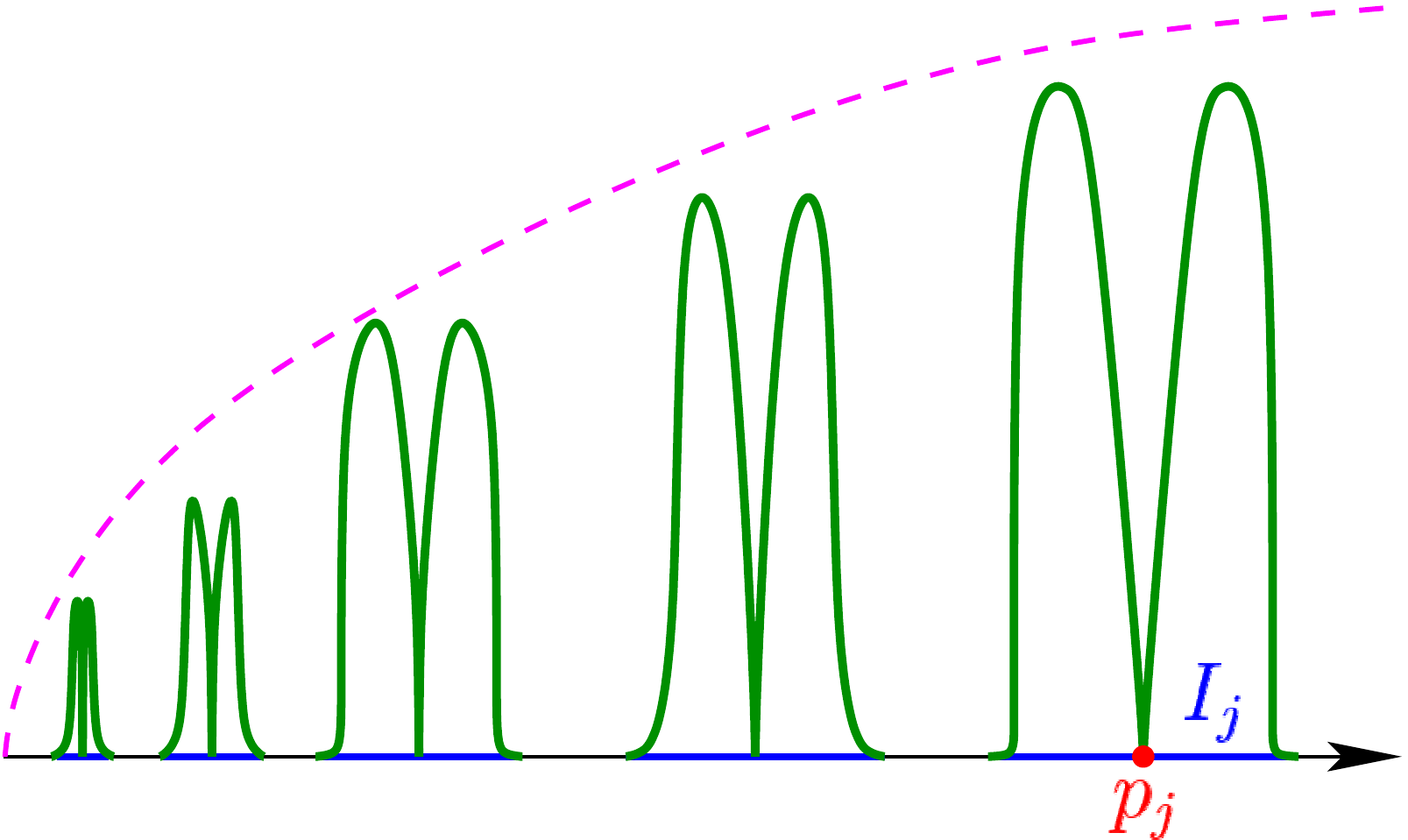}
 \caption{\sl The function~$u_\star$ in~\eqref{T12536Ec32674658769h65876ij658dhtjuylkiliR2543654yh2432RAe87658FI6tfh564EQ} with~$k=0$.}\label{T12536Ec32674658769h65876ij658dhtjuylkiliR2543654yh2432RAe87658FI6tfh564}
\end{figure}

For each~$j\in\N$, we define the function
$$ \R\ni t\mapsto u_j(t):=s_j |t-p_j|^{k+\alpha}\eta_j(t).$$
Notice that the supports of~$u_j$ do not intersect, hence we can define
\begin{equation}\begin{split}\label{T12536Ec32674658769h65876ij658dhtjuylkiliR2543654yh2432RAe87658FI6tfh564EQ}& \R\ni t\longmapsto u_\star(t):=\sum_{j\in\N} u_j(t)
\\
{\mbox{and }}\;\;&\R^n\ni x\longmapsto u(x)=u(x_1,\dots,x_n):=u_\star(x_1).\end{split}\end{equation}
See Figure~\ref{T12536Ec32674658769h65876ij658dhtjuylkiliR2543654yh2432RAe87658FI6tfh564} for a sketch of this situation.
We stress that~$u(x)\ge0$ for all~$x\in\R^n$ and~$u(x)=0$ unless~$x_1$
belongs to one of the disjoint intervals~$I_j$. Accordingly, for every~$x$ such that~$u(x)\ne0$
we denote by~$j_x\in\N$ the unique index for which~$x_1\in I_{j_x}$.

In this way, if~$u(x)\ne0$ then, utilizing~\eqref{w3tg-inoiqdjQunAM34S-a},
\begin{equation}\label{rgf4y544654yh54gJ3265M3256S-ijsa}\begin{split}&
|u(x)|=u(x)=u_{j_x}(x_1)\le s_{j_x} |x_1-p_{j_x}|^{k+\alpha}=\frac{|x_1-p_{j_x}|^{k+\alpha}}{\ell_{j_x}^\alpha}\\&\qquad\qquad=
\left(\frac{|x_1-p_{j_x}|}{\ell_{j_x}^{\frac\alpha{k+\alpha}}}\right)^{{k+\alpha}}\le|x_1|^{{k+\alpha}}\le|x|^{k+\alpha}.\end{split}
\end{equation}
In addition, for every~$x_0\in\R^n\setminus\{0\}$, there exists~$\rho_{x_0}>0$ such that~$u\in C^{k,\alpha}(B_{\rho_{x_0}}(x_0))$
and therefore, in the light of Proposition~\ref{XBOUN-AGG-P},
there exist~$M_{x_0}\ge0$ and polynomial~$P_{x_0}$ of degree at most~$k$
such that~$|u(x)-P_{x_0}(x)|\le M_{x_0}|x-x_0|^{k+\alpha}$ for all~$x\in B_{\rho_{x_0}}(x_0)$
(and actually for all~$x\in\R^n$, up to changing~$M_{x_0}$, since~$u$ is bounded).

It follows from this and~\eqref{rgf4y544654yh54gJ3265M3256S-ijsa} that
Definition~\ref{JOSLi4o245nfIJSU-184ur032hv9qythgf92y98trpgJA4NBf} is fulfilled for every~$x_0\in\R^n$,
taking~$P:=0$ and~$M:=1$ when~$x_0=0$ and~$P:=P_{x_0}$ and~$M:=M_{x_0}$
when~$x_0\ne0$. As a result, 
\begin{equation}\label{9087r98732ytfiu874d6h743657bLHJSFTMAtoIYo01}
{\mbox{$u\in C^{k,\alpha}_{\R^n}(x_0)$ for every~$x_0\in \R^n$.}}\end{equation}

But we claim that, for every~$\rho>0$,
\begin{equation}\label{9087r98732ytfiu874d6h743657bLHJSFTMAtoIYo04}
u\not\in C^{k,\alpha}(B_\rho).
\end{equation}
To check this, we argue by contradiction and we suppose that~$u\in C^{k,\alpha}(B_\rho)$.
Using a Taylor expansion for~$u$ and~\eqref{rgf4y544654yh54gJ3265M3256S-ijsa}, it follows that, for small~$|x|$,
$$ \sum_{{\beta\in\N^n}\atop{|\beta|\le k}}\frac{D^\beta u(0)}{\beta!}x^\beta+o(|x|^k)=o(|x|^k),$$
and therefore~$D^\beta u(0)=0$ for all~$\beta\in\N^n$ with~$|\beta|\le k$.
Thus, picking~$J\in\N$, to be taken as large as we wish,
\begin{equation}\label{838-bIJSUjmmdFNFiGIA45df-53637}
\left|\partial^k_1 u\left(p_J+\frac{\ell_J}8,0,\dots,0\right)\right|=
\left|\partial^k_1 u\left(p_J+\frac{\ell_J}8,0,\dots,0\right)-\partial^k_1 u(0)\right|
\le C\left(p_J+\frac{\ell_J}8\right)^\alpha.\end{equation}
But if~$x_1\in\left(p_J, p_J+\frac{\ell_J}4\right)$
we have that
$$ \partial^k_1 u(x)=c
\,s_J \,|x_1-p_J|^{\alpha}, \qquad{\mbox{with }}\;c:=\prod_{m=1}^k(m+\alpha)$$
and therefore
$$ \partial^k_1 u\left(p_J+\frac{\ell_J}8,0,\dots,0\right)= \frac{c\,s_J\,\ell_J^{\alpha}}{8^{\alpha}}=\frac{c}{8^{\alpha}}.
$$
This and~\eqref{838-bIJSUjmmdFNFiGIA45df-53637} entail that
$$ \frac{c}{8^{\alpha}}\le\lim_{J\to+\infty}C\left(p_J+\frac{\ell_J}8\right)^\alpha=0,$$
which is a contradiction.
This proves~\eqref{9087r98732ytfiu874d6h743657bLHJSFTMAtoIYo04}.

In view of~\eqref{9087r98732ytfiu874d6h743657bLHJSFTMAtoIYo01}
and~\eqref{9087r98732ytfiu874d6h743657bLHJSFTMAtoIYo04},
we infer that the inclusion in~\eqref{9087r98732ytfiu874d6h743657bLHJSFTMAtoIYo02}
is strict.
\end{proof}

One of the technical advantages of
the pointwise H\"older spaces introduced in
Definition~\ref{JOSLi4o245nfIJSU-184ur032hv9qythgf92y98trpgJA4NBf}
is that it suffices to check the validity of the polynomial approximation
along a geometrically decaying sequence of radii, according to the following result
(in which we use the notation of distance function introduced on page~\pageref{LAKSJHDKSDJFH0193287493utjgQUEMD}):

\begin{figure}
  \centering
  \includegraphics[width=.7\linewidth]{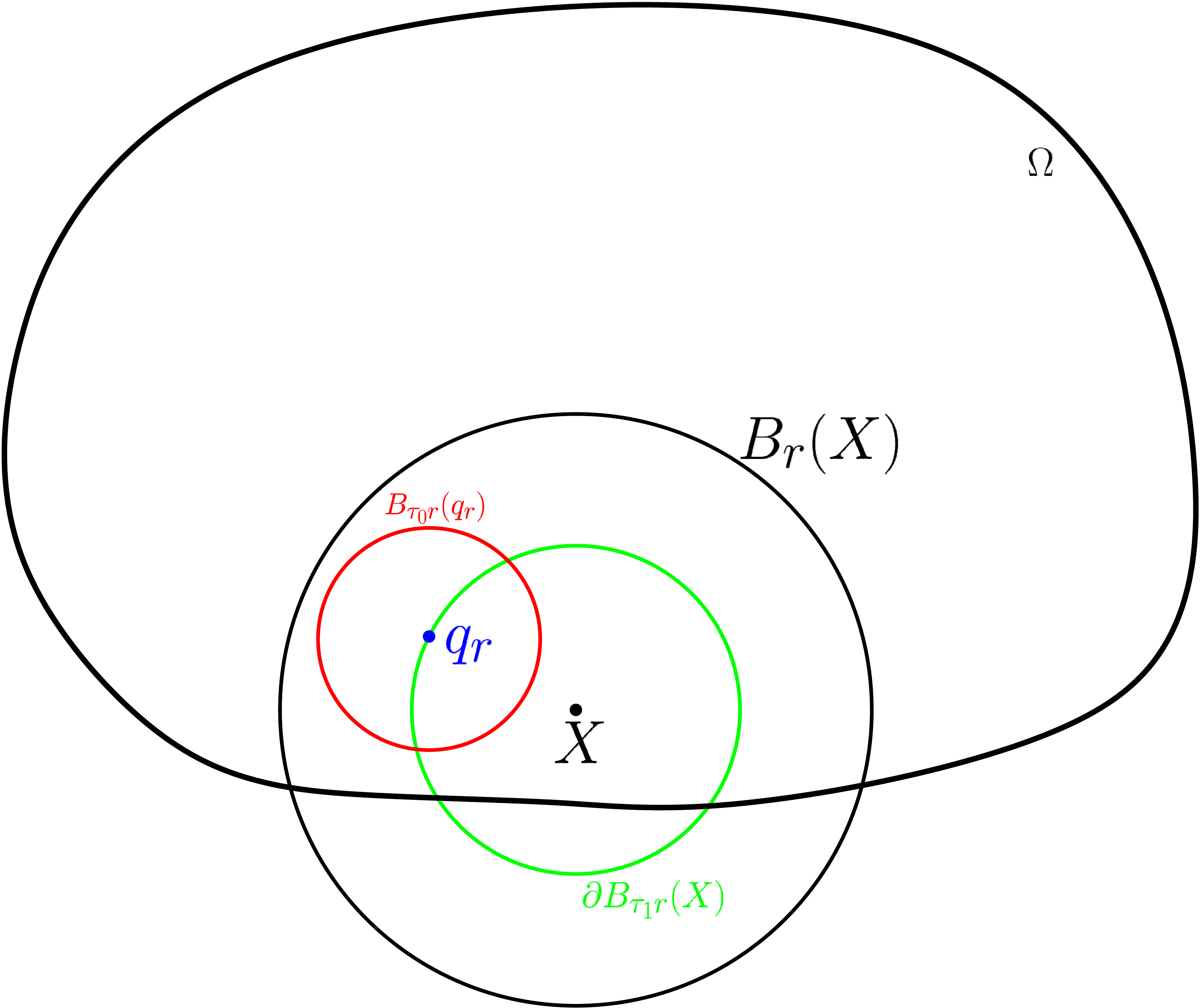}
 \caption{\sl The geometric condition in~\eqref{CLEANBALL}.}\label{23546g11ItangeFI}
\end{figure}

\begin{lemma}\label{POLIAPP-J-GEO}
Let~$k\in\N$, $\alpha\in(0,1]$ and~$\rho\in(0,1)$.
Let~$\Omega\subseteq\R^n$ be open and bounded.

Assume that there exist~$\tau_\star$, $\tau_\sharp$, $\tau_0$, $\tau_1$, $\in(0,1)$ such that for all~$X\in\Omega$
and all~$r>0$ such that~$\frac{ d_{\partial\Omega}(X)}{\tau_\sharp}\le r\le\tau_\star$
there exists~$q_r\in\partial B_{\tau_1 r}(X)$ such that\footnote{See Figure~\ref{23546g11ItangeFI}
for a sketch of condition~\eqref{CLEANBALL}.}
\begin{equation}\label{CLEANBALL}
B_{\tau_0 r}(q_r)\subseteq\Omega\cap B_r(X).
\end{equation}

Let~$x_0\in\Omega$ and~$u\in L^\infty(\Omega)$.
Then, the following conditions are equivalent:
\begin{itemize}
\item[(i).] $u\in C^{k,\alpha}_\Omega(x_0)$ and there exists~$M\ge0$ such that
\begin{equation*} |u(x)-P(x)|\le M|x-x_0|^{k+\alpha}\qquad{\mbox{for all}}\,x\in\Omega.\end{equation*}
\item[(ii).] There exists~$L\ge0$ such that for every~$j\in\N$
there exists a polynomial~$P_j$ of degree at most~$k$ such that
\begin{equation*} |u(x)-P_j(x)|\le L\rho^{(k+\alpha)j}\qquad{\mbox{for all}}\,x\in \Omega\cap B_{\rho^j}(x_0).\end{equation*}\end{itemize}

Also, if~(i) holds true, one can choose~$L:=M$ in~(ii).

If instead~(ii) holds true, one can choose~$M$ in~(i) of the form~$C \Big( \|u\|_{L^\infty(\Omega)}+L\Big)$, for some~$C>0$ depending only on~$n$, $k$, $\alpha$, $\rho$ and~$\Omega$.

Furthermore, the polynomial~$P$ in~(i) satisfies~$|D^\beta P(x_0)|\le C\Big( \|u\|_{L^\infty(\Omega)}+M\Big)$,
for all~$\beta\in\N^n$ with~$|\beta|\le k$,
for some~$C>0$ depending only on~$n$, $k$, $\alpha$, $\rho$ and~$\Omega$.
\end{lemma}

\begin{proof}
Up to a translation, we suppose~$x_0=0$.
If~(i) holds true, then~(ii) also holds true with~$P_j:=P$ and~$L:=M$.

Suppose instead that~(ii) holds true. The gist to prove~(i) is that the polynomials~$P_j$
may vary from scale to scale, but they must change ``very little'' at each iteration, due to
the geometric decay of the sequence of radii, and this fact will entail the existence of a ``limit polynomial''
with uniform estimates. To make this idea work, let us exploit~(ii) to notice that
\begin{equation}\label{BISqualcosaloyuythgb}
\|P_{j+1}-P_j\|_{L^\infty(\Omega\cap B_{\rho^{j+1}})}\le
\|P_{j+1}-u\|_{L^\infty(\Omega\cap B_{\rho^{j+1}})}+\|u-P_j\|_{L^\infty(\Omega\cap B_{\rho^{j }})}\le2L\rho^{(k+\alpha)j}.
\end{equation}

We now point out that~$\|\cdot\|_{L^\infty(B_1)}$ and~$\|\cdot\|_{L^\infty(B_{\tau_0}(\tau_1e_1))}$
are norms in the finite dimensional linear space of polynomials of degree at most~$k$:
by the equivalence of norms in finite dimensional spaces (see e.g.~\cite[Theorem~15.26]{MR606198})
it follows that, for every polynomial~$P$ of degree at most~$k$,
$$ \|P\|_{L^\infty(B_1)}\le C\| P\|_{L^\infty(B_{\tau_0}(\tau_1 e_1))},$$
for some~$C$ depending only on~$n$, $k$, $\tau_0$ and~$\tau_1$.

Now, if~$p\in\partial B_{\tau_1}$ we take a rotation~${\mathcal{R}}$ centered at the origin sending~$\tau_1e_1$ into~$p$.
Given a polynomial~$P$ of degree at most~$k$, we let~$P_{\mathcal{R}}(x):=P({\mathcal{R}}x)$
and we observe that~$P_{\mathcal{R}}$ is also a polynomial of degree at most~$k$.
Consequently, we have that
\begin{equation}\label{lkjhgfdpoiuygtfJSNQTG-ioqwfm34rtgh6yuj-1}
\begin{split}&
\|P\|_{L^\infty(B_1)}=\sup_{y\in B_1}|P(y)|=\sup_{x\in B_1}|P({\mathcal{R}}x)|=
 \|P_{\mathcal{R}}\|_{L^\infty(B_1)}\le C\| P_{\mathcal{R}}\|_{L^\infty(B_{\tau_0}(\tau_1 e_1))}\\&\qquad
=C\sup_{x\in B_{\tau_0}(\tau_1 e_1)} |P({\mathcal{R}}x)|
=C\sup_{y\in B_{\tau_0}(\tau_1 e_1)} |P(y)|
= C\| P\|_{L^\infty(B_{\tau_0}(p))}.\end{split}
\end{equation}
Furthermore, using again the equivalence of the norms in finite dimensional spaces,
\begin{equation}\label{lkjhgfdpoiuygtfJSNQTG-ioqwfm34rtgh6yuj-8}
\|P\|_{L^\infty(B_1)}\le C\| P\|_{L^\infty(B_{\tau_\sharp})}.\end{equation}

Now we claim that,
for every~$i\in\N\cap[0,k]$, every~$\beta\in\N^n$ with~$|\beta|=i$ and every polynomial~$P$ of degree at most~$k$,
\begin{equation}\label{PSFVDJODHBDEBDJRDFGHISNA} \|D^\beta P\|_{L^\infty(\Omega\cap B_{\rho^{j+1}})}\le\frac{C}{\rho^{i(j+1)}}\,\|{P}\|_{L^\infty(\Omega\cap B_{\rho^{j+1}})}.\end{equation}
To check this, we define~$\widetilde P(x):=P(\rho^{j+1}x)$ and observe that
\begin{equation}\label{lkjhgfdpoiuygtfJSNQTG-ioqwfm34rtgh6yuj-2}
\|D^\beta P\|_{L^\infty(\Omega\cap B_{\rho^{j+1}})}
\le
\|D^\beta P\|_{L^\infty(B_{\rho^{j+1}})}\le
\frac{C}{\rho^{i(j+1)}}\,\|{P}\|_{L^\infty(B_{\rho^{j+1}})}=
\frac{C}{\rho^{i(j+1)}}\,\|\widetilde{P}\|_{L^\infty(B_{1})},\end{equation}
thanks to~\eqref{13432rfw2trfew5terghtry54yhgfjyhtrh254t5trw547h3hdhdrer-0-P09SBIS}.

Now we distinguish two cases, either~$ \frac{d_{\partial\Omega}(0)}{{\tau_\sharp}}\le{\rho^{j+1}}$
or~$\frac{d_{\partial\Omega}(0)}{\tau_\sharp}>{\rho^{j+1}}$.

Let us first
consider the case in which~$d_{\partial\Omega}(0)\le\frac{\rho^{j+1}}{\tau_\sharp}$.
In this situation, we can exploit~\eqref{CLEANBALL} and choose~$p\in \partial B_{\tau_1}$ such that~$B_{\tau_0\rho^{j+1}}(\rho^{j+1}p)\subseteq
\Omega\cap B_{\rho^{j+1}}$. Using this, \eqref{lkjhgfdpoiuygtfJSNQTG-ioqwfm34rtgh6yuj-1}
and~\eqref{lkjhgfdpoiuygtfJSNQTG-ioqwfm34rtgh6yuj-2} we obtain that
\begin{eqnarray*}&& \|D^\beta P\|_{L^\infty(\Omega\cap B_{\rho^{j+1}})}\le
\frac{C}{\rho^{i(j+1)}}\,\|\widetilde{P}\|_{L^\infty(B_{1})}\le
\frac{C}{\rho^{i(j+1)}}\,\|\widetilde{P}\|_{L^\infty(B_{\tau_0}(p))}\\&&\qquad
=\frac{C}{\rho^{i(j+1)}}\,\|{P}\|_{L^\infty(B_{\tau_0\rho^{j+1}}(\rho^{j+1}p))}\le\frac{C}{\rho^{i(j+1)}}\,\|{P}\|_{L^\infty(\Omega\cap B_{\rho^{j+1}})},
\end{eqnarray*}
showing that~\eqref{PSFVDJODHBDEBDJRDFGHISNA} holds true
in this case.

Now we focus on the case~$\frac{d_{\partial\Omega}(0)}{\tau_\sharp}>{\rho^{j+1}}$.
In this situation, we have that~$B_{\tau_\sharp\rho^{j+1}}\subseteq B_{d_{\partial\Omega}(0)}\subseteq\Omega$. Hence, 
$$ B_{\tau_\sharp\rho^{j+1}}=B_{\tau_\sharp\rho^{j+1}}\cap B_{\rho^{j+1}}
\subseteq
\Omega\cap B_{\rho^{j+1}}.$$
As a consequence,
by~\eqref{lkjhgfdpoiuygtfJSNQTG-ioqwfm34rtgh6yuj-8} and~\eqref{lkjhgfdpoiuygtfJSNQTG-ioqwfm34rtgh6yuj-2},
\begin{eqnarray*}&&
\|D^\beta P\|_{L^\infty(\Omega\cap B_{\rho^{j+1}})}
\le
\frac{C}{\rho^{i(j+1)}}\,\|\widetilde{P}\|_{L^\infty(B_{1})}\le
\frac{C}{\rho^{i(j+1)}}\,\|\widetilde{P}\|_{L^\infty(B_{\tau_\sharp})}\\&&\qquad=
\frac{C}{\rho^{i(j+1)}}\,\|{P}\|_{L^\infty(B_{\tau_\sharp \rho^{j+1}})}\le
\frac{C}{\rho^{i(j+1)}}\,\|{P}\|_{L^\infty(\Omega\cap B_{\rho^{j+1}})}
,\end{eqnarray*}
which completes the proof of~\eqref{PSFVDJODHBDEBDJRDFGHISNA}.

As a result, from~\eqref{BISqualcosaloyuythgb} and~\eqref{PSFVDJODHBDEBDJRDFGHISNA},
for every~$i\in\N\cap[0,k]$ and every~$\beta\in\N^n$ with~$|\beta|=i$,
\begin{eqnarray*}&& \|D^\beta(P_{j+1}-P_j)\|_{L^\infty(\Omega\cap B_{\rho^{j+1}})}
\le\frac{C}{\rho^{i(j+1)}}\,\|{P_{j+1}-P_j}\|_{L^\infty(  B_{\rho^{j+1}})}\le
CL\rho^{(k+\alpha)j-i(j+1)}
\end{eqnarray*}
for some~$C>0$ depending only on~$n$, $k$ and~$\Omega$.

Consequently, if we write
$$ P_j(x)=\sum_{{\beta\in\N^n}\atop{|\beta|\le k}} a_{\beta,j} x^\beta$$
we have that
$$ |a_{\beta,j+1}-a_{\beta,j}|=\left|\frac{D^\beta P_{j+1}(0)}{\beta!}-
\frac{D^\beta P_{j}(0)}{\beta!}
\right|\le C\|D^\beta(P_{j+1}-P_j)\|_{L^\infty(\Omega\cap B_{\rho^{j+1}})}\le CL\rho^{(k+\alpha)j-|\beta|(j+1)},$$
which in turn entails that for all~$m\in\N\cap[1,+\infty)$
\begin{equation}\label{4254512378tgbM35hA24O2S0jnfg0ujhv} |a_{\beta,j+m}-a_{\beta,j}|\le
\sum_{i=0}^{m-1} |a_{\beta,j+i+1}-a_{\beta,j+i}|\le
CL\sum_{i=0}^{m-1}\rho^{(k+\alpha)(j+i)-|\beta|(j+i+1)}
\le CL\rho^{(k-|\beta|+\alpha)j-|\beta|}.\end{equation}
In particular,
\[ |a_{\beta,j+m}-a_{\beta,j}|\le CL\rho^{ \alpha j-k},\]
which is infinitesimal as~$j\to+\infty$.

This gives that~$a_{\beta,j}$ is a converging sequence, say~$a_{\beta,j}\to a_{\beta}$ as~$j\to+\infty$.
Moreover, passing to the limit as~$m\to+\infty$ in~\eqref{4254512378tgbM35hA24O2S0jnfg0ujhv},
\begin{equation}\label{345yinsetminus} |a_{\beta}-a_{\beta,j}|\le CL\rho^{(k-|\beta|+\alpha)j-|\beta|}.\end{equation}
As a result, letting
$$P(x):=\sum_{{\beta\in\N^n}\atop{|\beta|\le k}} a_\beta x^\beta$$
and
using also~(ii) we see that, for all~$x\in \Omega\cap B_{\rho^j}$,
\begin{equation}\label{9ojer4betalpha}
\begin{split}&
|u(x)-P(x)|\le|u(x)-P_j(x)|+|P_j(x)-P(x)|\le L\rho^{(k+\alpha)j}+\left|
\sum_{{\beta\in\N^n}\atop{|\beta|\le k}} (a_{\beta,j}-a_\beta) x^\beta
\right|\\&\quad\qquad \le L\rho^{(k+\alpha)j}+
CL\sum_{{\beta\in\N^n}\atop{|\beta|\le k}}\rho^{(k+\alpha)j-|\beta|}\le L\rho^{(k+\alpha)j} +CL\rho^{(k+\alpha)j-k}
\le CL\rho^{(k+\alpha)j},
\end{split}
\end{equation}
up to renaming~$C$ in dependence of~$n$, $k$, $\alpha$, $\rho$ and~$\Omega$.

Now we take~$x\in\Omega$ and we distinguish two cases.
If~$x\in B_1$, we pick~$j\in\N$ such that~$\rho^{j+1}\le|x|<\rho^j$
and we exploit~\eqref{9ojer4betalpha} to conclude that
\begin{equation}\label{9ojer4betalpha2}
|u(x)-P(x)|\le CL\rho^{(k+\alpha)j}\le
CL|x|^{k+\alpha},
\end{equation}
up to renaming~$C$.

If instead~$x\in\Omega\setminus B_1$, then
$$ |u(x)-P(x)|\le \|u\|_{L^\infty(\Omega\setminus B_1)}+
\sum_{{\beta\in\N^n}\atop{|\beta|\le k}}|a_\beta|\,|x|^{|\beta|}\le
\Big( \|u\|_{L^\infty(\Omega\setminus B_1)}+\sup_{{\beta\in\N^n}\atop{|\beta|\le k}}|a_\beta|\Big)|x|^{k+\alpha}.
$$
Since, by~\eqref{13432rfw2trfew5terghtry54yhgfjyhtrh254t5trw547h3hdhdrer-0-P09S} and~\eqref{345yinsetminus},
$$ |a_{\beta}|\le |a_{\beta,0}|+ CL\rho^{-|\beta|}=\left|\frac{D^\beta P_0(0)}{\beta!}\right|+ CL\rho^{-|\beta|}
\le
\|P_0\|_{C^k(B_1)}+CL\le C
\|P_0\|_{L^\infty(B_1)}+CL,$$
we thus obtain that, for all~$x\in\Omega\setminus B_1$,
\begin{equation}\label{092yehgfggTGS3inftyOmegaL} |u(x)-P(x)|\le
C \Big( \|u\|_{L^\infty(\Omega\setminus B_1)}+\|P_0\|_{L^\infty(B_1)}+L\Big)|x|^{k+\alpha}
.\end{equation}
Notice also that, for all~$x\in\Omega\cap B_1$,
$$ |P_0(x)|\le|u(x)|+L\le\|u\|_{L^\infty(\Omega)}+L,$$
thanks to~(ii),
and therefore
\begin{equation}\label{3QU2JENRIGBFGBNDA3G32G24A}
\|P_0\|_{L^\infty(\Omega\cap B_1)}\le\|u\|_{L^\infty(\Omega)}+L.
\end{equation}
Moreover, by the equivalence of the norms on finite dimensional spaces, we know that~$\|P_0\|_{L^\infty( B_1)}\le C\|P_0\|_{L^\infty(\Omega\cap B_1)}$, hence we infer from~\eqref{3QU2JENRIGBFGBNDA3G32G24A} that
$$ \|P_0\|_{L^\infty( B_1)}\le C\big(\|u\|_{L^\infty(\Omega)}+L\big).$$
This and~\eqref{092yehgfggTGS3inftyOmegaL} yield that, for all~$x\in\Omega\setminus B_1$,
\begin{equation*} |u(x)-P(x)|\le
C \Big( \|u\|_{L^\infty(\Omega)}+L\Big)|x|^{k+\alpha}
.\end{equation*}
{F}rom this inequality and~\eqref{9ojer4betalpha2}, we see that~(i) holds true with~$M:=
C \Big( \|u\|_{L^\infty(\Omega)}+L\Big)$, as desired.

Additionally, a byproduct of~(i) is that, if~$x\in \Omega$,
$$ |P(x)|\le|u(x)|+CM\le\|u\|_{L^\infty(\Omega)}+CM$$
and thus, by~\eqref{13432rfw2trfew5terghtry54yhgfjyhtrh254t5trw547h3hdhdrer-0-P09S},
for all~$\beta\in\N^n$ with~$|\beta|\le k$,
$$ |D^\beta P(0)|\le\|P\|_{C^k(\Omega)}\le C\,\|P\|_{L^\infty(\Omega)}\le C\big(
\|u\|_{L^\infty(\Omega)}+M\big),$$
thus concluding the proof of Lemma~\ref{POLIAPP-J-GEO}.
\end{proof}

\begin{corollary}\label{i0fU824aiugfepZAkkZfPXBOUN-AGG-P}
Let~$k\in\N$, $\alpha\in(0,1]$, $L\ge0$ and~$\rho\in(0,1)$.
Let~$\Omega\subseteq\R^n$ be open, bounded and convex.

Let~$u\in L^\infty(\R^n)$ and suppose that for every~$x_0\in\Omega$
and every~$j\in\N$
there exists a polynomial~$P_{j,x_0}$ of degree at most~$k$ such that
\begin{equation*}
|u(x)-P_{j,x_0}(x)|\le L\rho^{(k+\alpha)j}\qquad{\mbox{for all}}\,x\in \Omega\cap B_{\rho^j}(x_0)
.\end{equation*}
Then, $u\in C^{k,\alpha}(\Omega)$ and~$\|u\|_{C^{k,\alpha}(\Omega)}\le C(L+\|u\|_{L^\infty(\R^n)})$,
for some~$C>0$ depending only on~$n$, $k$, $\alpha$ and~$\Omega$.
\end{corollary}

\begin{proof} We observe that
\begin{equation}\label{43rSefd5ATDIFSQAWSA}{\mbox{$\Omega$ satisfies condition~\eqref{CLEANBALL}.}}\end{equation}
For this, to start with, we show that~$\partial\Omega$ can be written locally as the graph of a Lipschitz
function, namely
\begin{equation}\label{43rSefd5ATDIFSQAWSAL}\begin{split}&{\mbox{for every~$p\in\Omega$ 
there exists~$\rho_p>0$ such that~$(\partial\Omega)\cap B_{\rho_p}(p)$
is the graph,}}\\&{\mbox{in some direction, of a convex and Lipschitz function.}}\end{split}\end{equation}
To establish this, by the boundedness of~$\Omega$ we can suppose that~$\Omega\Subset B_R$ for some~$R>0$.
Hence, given any~$\omega\in\partial B_1$, the hyperplane~$\Pi_\omega:=\{\omega\cdot x=-R\}$
lies outside~$\Omega$ (otherwise there would be a point~$x\in\Omega\Subset B_R$ with~$|x|\ge|\omega\cdot x|=|-R|=R$, which is a
contradiction). We can thus consider the projection~$P_\omega$ of~$\Omega$ onto~$\Pi_\omega$,
which is open (and convex) in~$\Pi_\omega$ since~$\Omega$ is open
(and convex). For all~$x\in P_\omega$
we denote by~$f_\omega(x)$ the infimal~$t$ for which~$x+t\omega$ belongs to~$\overline\Omega$, see Figure~\ref{conpro}.

\begin{figure}
  \centering
  \includegraphics[width=.70\linewidth]{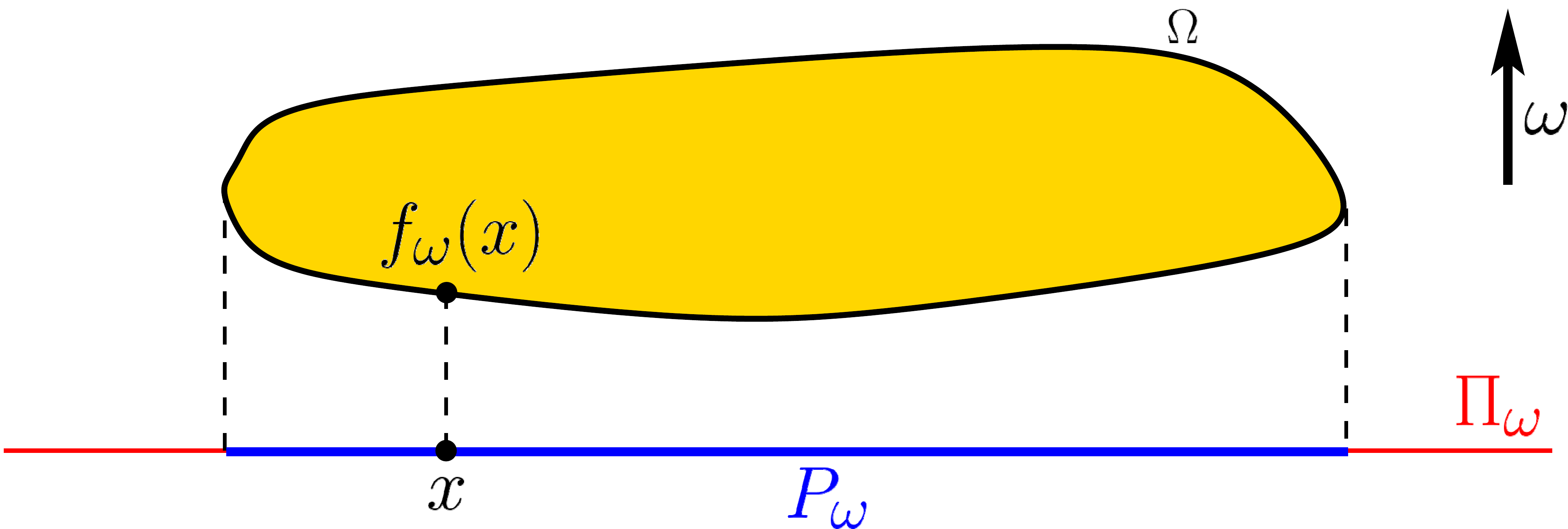}
 \caption{\sl The geometry involved in~\eqref{43rSefd5ATDIFSQAWSAB} and~\eqref{43rSefd5ATDIFSQAWSABC}.}\label{conpro}
\end{figure}

Notice that, by construction,
\begin{equation}\label{43rSefd5ATDIFSQAWSABbb}
f_\omega(x)\in[0,2R].\end{equation}
Also, 
\begin{equation}\label{43rSefd5ATDIFSQAWSAB}
{\mbox{$f_\omega$ is convex (when we look at it as a function in direction~$\omega$).}}\end{equation}
To check this, up to a rotation, we can assume that~$\omega=e_n$. We pick~$(x',-R)$, $(y',-R)\in P_{e_n}$
and we note that~$(x',f_{e_n}(x'))$, $(y',f_{e_n}(y'))\in\partial\Omega$.
{F}rom the convexity of~$\Omega$ we infer that
$$ \Big(tx'+(1-t)y',\,tf_{e_n}(x')+(1-t) f_{e_n}(y')\Big)=
t(x',f_{e_n}(x'))+(1-t)(y',f_{e_n}(y'))\in\overline\Omega$$
for all~$t\in[0,1]$. Hence, by the definition of~$f_{e_n}$, we conclude that
$$ f_{e_n}\big(tx'+(1-t)y'\big)\le tf_{e_n}(x')+(1-t) f_{e_n}(y'),$$
which proves~\eqref{43rSefd5ATDIFSQAWSAB}.

Additionally, we have that
\begin{equation}\label{43rSefd5ATDIFSQAWSABC}
{\mbox{$f_\omega$ is locally Lipschitz.}}\end{equation}
This is a general argument about convexity (see e.g.~\cite[Theorem~6.7]{MR3409135}) and goes as follows.
Again, we suppose~$\omega=e_n$ and we take~$D\Subset P_{e_n}$.
Let~$d>0$ be the distance of~$D$ to~$\partial P_{e_n}$
and let also~$(x',-R)$ and~$(y',-R)$ be two different points in~$D$ with~$|x'-y'|<d$.
We define~$e:=\frac{y'-x'}{|y'-x'|}$ and, for all~$\mu\ge0$,
$$ x'_\mu:=x'+\mu e.$$
We pick~$\mu$ such that~$(x'_\mu,-R)\in\partial P_{e_n}$
and we stress that~$\mu\ge d$.

Besides, if~$t:=\frac{|x'-y'|}{\mu}$ we have that~$t\in\left[0,\frac{d}\mu\right]\subseteq[0,1]$.
Furthermore,
$$ tx'_\mu+(1-t)x'= t(x'_\mu-x')+x'=\frac{|x'-y'|\mu e}{\mu}+x'=(y'-x')+x'=y'$$
and accordingly, by~\eqref{43rSefd5ATDIFSQAWSAB},
$$ f_{e_n}(y')=
f_{e_n}\big(tx'_\mu+(1-t)x'\big)\le tf_{e_n}(x'_\mu)+(1-t) f_{e_n}(x')=f_{e_n}(x')+
t\big(f_{e_n}(x'_\mu)- f_{e_n}(x')\big)
.$$
This and~\eqref{43rSefd5ATDIFSQAWSABbb} give that
$$ f_{e_n}(y')-f_{e_n}(x')\le t
\big(|f_{e_n}(x'_\mu)|+| f_{e_n}(x')|\big)
\le \frac{4R |x'-y'|}{\mu}\le \frac{4R |x'-y'|}{d}.$$
By exchanging~$x'$ and~$y'$ we thereby conclude that, for all~$(x',-R)$, $(y',-R)\in D$ with~$|x'-y'|<d$,
\begin{equation}\label{43rSefd5ATDIFSQAWSABC27hfrj-3}
|f_{e_n}(x')-f_{e_n}(y')|\le\frac{4R |x'-y'|}{d}.
\end{equation}
Additionally, if~$(x',-R)$, $(y',-R)\in D$ with~$|x'-y'|\ge d$, we deduce from~\eqref{43rSefd5ATDIFSQAWSABbb}
that
$$ |f_{e_n}(x')-f_{e_n}(y')|\le 4R\le\frac{4R |x'-y'|}{d}.$$
Hence we combine this inequality with~\eqref{43rSefd5ATDIFSQAWSABC27hfrj-3}
and we see that the proof of~\eqref{43rSefd5ATDIFSQAWSABC} is complete.

{F}rom~\eqref{43rSefd5ATDIFSQAWSAB} and~\eqref{43rSefd5ATDIFSQAWSABC} (and reducing~$P_\omega$ to a slightly smaller subdomain
to obtain a global Lipschitz property)
we thus obtain~\eqref{43rSefd5ATDIFSQAWSAL}, as desired.

Thus, since~$\partial\Omega$ is compact, we use~\eqref{43rSefd5ATDIFSQAWSAL} to cover~$\partial\Omega$
with a finite number of balls~$B_{\rho_1/2}(p_1),$ $\dots$, $B_{\rho_N/2}(p_N)$, with~$p_1,\dots,p_N\in\partial\Omega$,
such that for all~$i\in\{1,\dots,N\}$
we have that~$(\partial\Omega)\cap B_{\rho_i}(p_i)$ can be written as a graph of a convex and Lipschitz function~$f_i$.

With these preliminary observations, we can now address the core of the proof of~\eqref{43rSefd5ATDIFSQAWSA}.
For this we let
\begin{equation}\label{ITAUSP} \tau_\star:=\min\left\{ 1,\frac{\rho_1}{16},\dots,\frac{\rho_N}{16}\right\},\qquad\tau_\sharp:=\frac18\qquad{\mbox{and}}\qquad\tau_1:=\frac12
,\end{equation}
and we pick~$X\in\Omega$. We let~$d(X)$ be the distance from~$X$ to~$\partial\Omega$ and we also
pick~$r\in\left[\frac{d(X)}{\tau_\sharp},\tau_\star\right]$.

In this way, we have that~$B_{d(X)}(X)\subseteq\Omega$
and there exists~$Y\in(\partial B_{d(X)}(X))\cap(\partial\Omega)$.
We observe that if~$i\in\{1,\dots, N\}$ such that~$Y\in B_{\rho_i/2}(p_i)$ then
\begin{equation}\label{displaystoejfbngbghfhfrgh}
B_r(X)\subseteq B_{\rho_i}(p_i), \end{equation}
because if~$z\in B_r(X)$ then
$$|z-p_i|\le|z-X|+|X-Y|+|Y-p_i|<
r+d(X)+\frac{\rho_i}2\le 2r+\frac{\rho_i}2\le2\tau_\star+\frac{\rho_i}2<\rho_i,$$
owing to~\eqref{ITAUSP}.

\begin{figure}
  \centering
  \includegraphics[width=.60\linewidth]{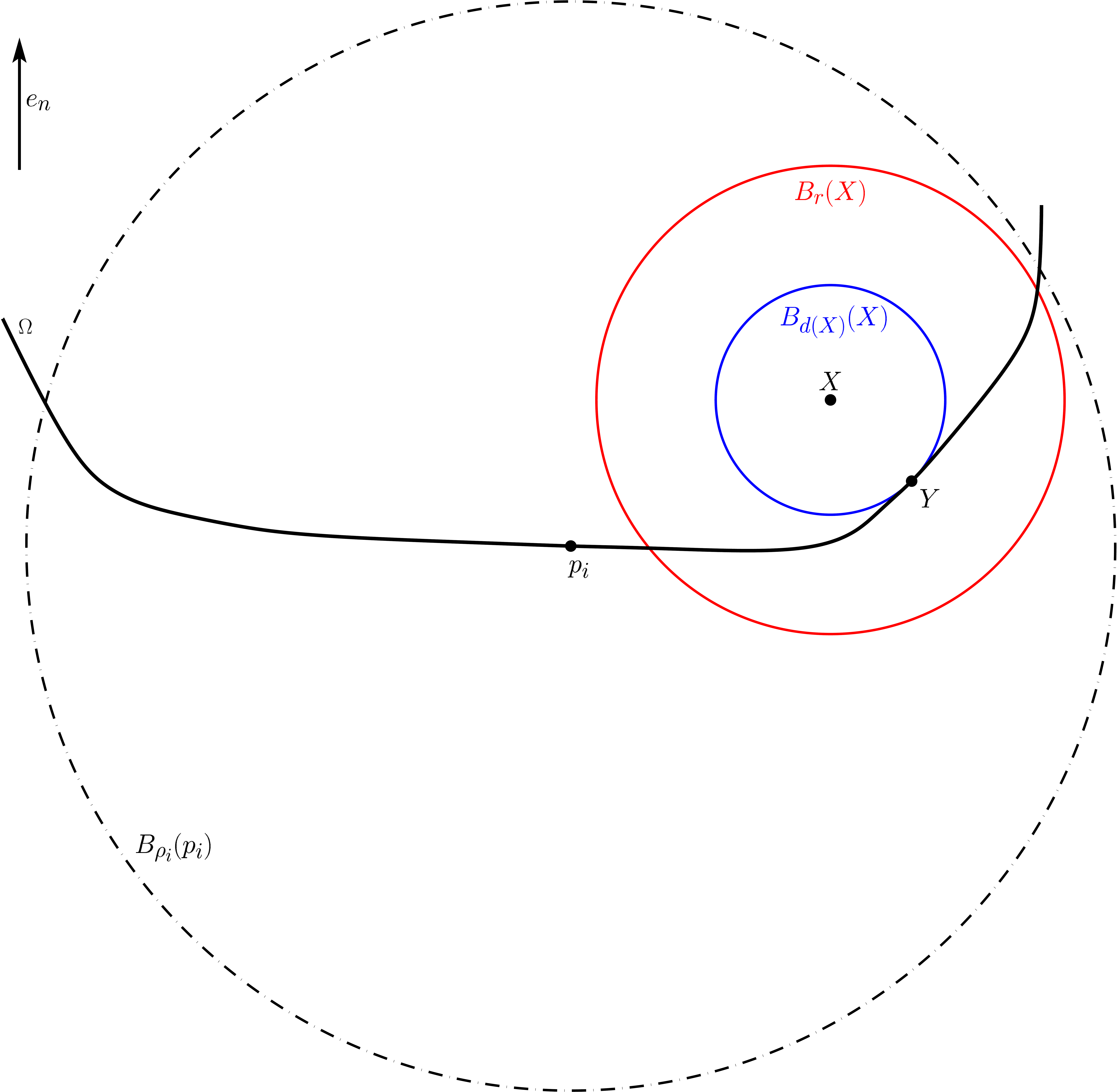}
 \caption{\sl The geometry involved in~\eqref{A0mnRnSIllDmAIJDd9lGgb}.}\label{puota}
\end{figure}

Up to a rotation, we may assume that~$f_i$ is a convex function in the direction of~$e_n$, see Figure~\ref{puota}.
Accordingly, for some~$M>0$ (related to the maximum of the Lipschitz moduli of continuity of
the functions~$f_1,\dots,f_N$), we have that
\begin{equation}\label{A0mnRnSIllDmAIJDd9lGgb} \Big\{ x_n-Y_n\ge M|x'-Y'|\Big\}\cap B_{\rho_i}(p_i)\subseteq\Omega.\end{equation}
Furthermore, since~$|Y-X|=d(X)\le\frac{r}8<\frac{r}2$, hence~$Y\in B_{r/2}(X)$ we can take~$\vartheta\ge0$
such that~$q_r:=Y+\vartheta e_n\in \partial B_{r/2}(X)$. We stress that~$q_r\in\partial_{\tau_1 r}(X)$, due to~\eqref{ITAUSP}.
Also,
\begin{equation}\label{A0mnRnSIllDmAIJDd9} \vartheta=|\vartheta e_n|=|q_r-Y|\ge|q_r-X|-|X-Y|=\frac{r}2-d(X)\ge\frac{r}4.\end{equation}
Similarly,
\begin{equation}\label{A0mnRnSIllDmAIJDd9u}
\vartheta=|\vartheta e_n|=|q_r-Y|\le|q_r-X|+|X-Y|=\frac{r}2+d(X)\le\frac{5r}8.\end{equation}
We also set
\begin{equation}\label{ITAUSP2} \tau_0:=\frac1{16(M+1)}
\end{equation}
and we claim that
\begin{equation}\label{A0mnRnSIllDmAIJDd9l}
B_{\tau_0 r}(q_r)\subseteq\Big\{ x_n-Y_n\ge M|x'-Y'|\Big\}.  
\end{equation}
Indeed, if~$\zeta\in B_{\tau_0 r}(q_r)$ then
\begin{eqnarray*}
&& \zeta_n-Y_n- M|\zeta'-Y'| \ge q_{r,n}-Y_n- M|q_r'-Y'|-(M+1)\tau_0 r\\&&\qquad=
\vartheta- 0-(M+1)\tau_0 r\ge\frac{r}4-(M+1)\tau_0 r>0,
\end{eqnarray*}
thanks to~\eqref{A0mnRnSIllDmAIJDd9}, and this proves~\eqref{A0mnRnSIllDmAIJDd9l}.

We also have that
\begin{equation}\label{A0mnRnSIllDmAIJDd9l6}
B_{\tau_0 r}(q_r)\subseteq B_r(X)
\end{equation}
because if~$\eta\in B_{\tau_0 r}(q_r)$ then, by~\eqref{A0mnRnSIllDmAIJDd9u},
$$ |\eta-X|\le|\eta-q_r|+|q_r-Y|+|Y-X|<\tau_0 r+\vartheta+d(X)
\le \tau_0 r+\frac{5r}8+\frac{r}8<r,
$$
which proves~\eqref{A0mnRnSIllDmAIJDd9l6}.

Now, gathering the results in~\eqref{displaystoejfbngbghfhfrgh}, \eqref{A0mnRnSIllDmAIJDd9lGgb}, \eqref{A0mnRnSIllDmAIJDd9l}
and~\eqref{A0mnRnSIllDmAIJDd9l6}, we find that
$$ B_{\tau_0 r}(q_r)\subseteq\{ x_n-Y_n\ge M|x'-Y'|\}\cap B_r(X)
\subseteq \Omega\cap B_r(X),$$
which is the desired claim in~\eqref{CLEANBALL}.
The proof of~\eqref{43rSefd5ATDIFSQAWSA} is thereby complete.

Thus, in light of~\eqref{43rSefd5ATDIFSQAWSA}, we can use Lemma~\ref{POLIAPP-J-GEO}
and find that, for each~$x_0\in\Omega$ we can find
a polynomial~$P_{x_0}$ of degree at most~$k$ such that
\begin{equation*} \begin{split}
&|u(x)-P_{x_0}(x)|\le M|x-x_0|^{k+\alpha}\qquad{\mbox{for all}}\,x\in \Omega
\\{\mbox{and }}\;\;&|D^\beta P_{x_0}(x_0)|\le M \qquad{\mbox{ for all~$\beta\in\N^n$ with~$|\beta|\le k$}},\end{split}\end{equation*}
with~$M:=C(L+\|u\|_{L^\infty(\R^n)})$.

{F}rom this and Proposition~\ref{XBOUN-AGG-P} the desired result follows.
\end{proof}

We stress that the regularity theory developed by potential analysis, for instance, in the proof of 
Theorem~\ref{SCHAUDER-INTE} relied crucially on a very delicate derivative computation
(as performed in Proposition~\ref{SI:IN:SE:DR})
and in the fine analysis of the singular integral
produced by this computation (as carried out in
Lemma~\ref{NKSD5678yihnqehr837tCOvnlIDJMcINBSno3Na7syshdrr9ISDcvaHF23}).
Instead, the pointwise H\"older spaces do not require at all to take derivatives
and Corollary~\ref{i0fU824aiugfepZAkkZfPXBOUN-AGG-P} allows one to check the pointwise polynomial approximation
only at suitable scales. To highlight the conceptual innovation
of the methods leveraging the pointwise H\"older spaces, we provide here a different proof 
of Theorem~\ref{SCHAUDER-INTE}.

\begin{proof}[Another proof of Theorem~\ref{SCHAUDER-INTE}] \label{023iurjw94jfJAMAnJHOSo3o24f32tofXZwf}
Up to a covering argument (see Lemma~\ref{COVERINGARG}), we can prove Theorem~\ref{SCHAUDER-INTE} with the ball~$B_1$
replaced by the ball~$B_4$ (changing the initial size of the ball will be convenient
in the forthcoming equation~\eqref{OIHSsot1inv2doPAEer20pp},
where we will be then allowed to work in the unit ball for the function that
will occupy most of the computations needed for this argument).

Also, we observe that it suffices to prove Theorem~\ref{SCHAUDER-INTE}
under the additional assumptions that
\begin{equation}\label{HYPETheoremSCHAUDER-INTE}
\|u\|_{L^\infty(B_4)}\le1 {\mbox{\qquad{and}\qquad}}\|f\|_{C^\alpha(B_4)}\le1.
\end{equation}
Indeed, if we know that Theorem~\ref{SCHAUDER-INTE} holds true under this additional assumption,
we define~$\widetilde u:=\frac{u}{\|u\|_{L^\infty(B_4)}+\|f\|_{C^\alpha(B_4)}}$
and~$\widetilde f:=\frac{f}{\|u\|_{L^\infty(B_4)}+\|f\|_{C^\alpha(B_4)}}$,
we observe that~$\|\widetilde u\|_{L^\infty(B_4)}\le1$ and~$\|\widetilde f\|_{C^\alpha(B_4)}\le1$
and we are thereby in the position of concluding that
$$ \frac{\|u\|_{C^{2,\alpha}(B_{1/2})}}{\|u\|_{L^\infty(B_4)}+\|f\|_{C^\alpha(B_4)}}
=\|\widetilde u\|_{C^{2,\alpha}(B_{1/2})}
\le C,$$
thus providing Theorem~\ref{SCHAUDER-INTE} in its full generality.

For this reason, we now assume that the additional hypotheses in~\eqref{HYPETheoremSCHAUDER-INTE}
is satisfied.

In light of Corollary~\ref{i0fU824aiugfepZAkkZfPXBOUN-AGG-P}, we also point out that to establish Theorem~\ref{SCHAUDER-INTE} it is enough to check
that there exists~$\rho\in\left(0,\frac12\right)$ such that for every~$x_0\in B_{1/2}$
and every~$j\in\N$
there exists a polynomial~$P_{j,x_0}$ of degree at most~$2$ such that
\begin{equation}\label{90oihfe8yhfh24inftyB1}
|u(x)-P_{j,x_0}(x)|\le M\rho^{(2+\alpha)j}\qquad{\mbox{for all}}\,x\in B_{\rho^j}(x_0),\end{equation}
where~$M>0$ depends only on~$n$ and~$\alpha$.

For this, we can consider the translation~$v(x):=u(x_0+x)$
and notice that, by~\eqref{HYPETheoremSCHAUDER-INTE},
\begin{equation}\label{OIHSsot1inv2doPAEer20pp}
\begin{split}&\Delta v(x)=f(x_0+x)=:\phi(x) {\mbox{ for all }}x\in B_{1},\\
&\|v\|_{L^\infty(B_{1})}\le\|u\|_{L^\infty(B_4)}\le1\\{\mbox{and }}\;\;&
\|\phi\|_{C^\alpha(B_{1})}\le \|f\|_{C^\alpha(B_4)}\le1.
\end{split}\end{equation}
In this setting, to establish~\eqref{90oihfe8yhfh24inftyB1}, we aim at proving that
\begin{equation}\label{90oihfe8yhfh24inftyB2}
 \begin{split}
&|v(x)-P_{j}(x)|\le M\rho^{(2+\alpha)j}\qquad{\mbox{for all}}\,x\in B_{\rho^j} ,\end{split}\end{equation}
for a suitable polynomial~$P_j$ of degree at most~$2$ and with~$C>0$ depending only
on~$n$ and~$\alpha$.
As a matter of fact, we will prove additionally that we can choose~$P_j$ such that
\begin{equation}\label{90oihfe8yhfh24inftyBT}
\Delta P_j(x)=\phi(0)\quad{\mbox{ for every }}x\in\R^n.
\end{equation}
Quite surprisingly, imposing such an additional constraint will turn out to simplify, rather than complicate,
the proof of~\eqref{90oihfe8yhfh24inftyB2}. In a sense, this additional constraint suggests that
an appropriate strategy to detect regularity is ``to replace the source by a constant one near the origin''
thus reducing to the case of harmonic functions. The technical details are the following.

Our goal is to prove~\eqref{90oihfe8yhfh24inftyB2} and~\eqref{90oihfe8yhfh24inftyBT}
by induction over~$j\in\N$.
For this, we are free to choose~$\rho\in\left(0,\frac12\right)$ conveniently small and~$M\ge1$ conveniently large (in a universal fashion
depending only on~$n$ and~$\alpha$)
and we argue as follows. When~$j=0$ we choose~$P_0(x):=\frac{\phi(0)}{2n}|x|^2$; in this way, equation~\eqref{90oihfe8yhfh24inftyBT}
holds true for~$j=0$ and, for all~$x\in B_{1}$,
$$ |v(x)-P_{0}(x)|\le \|v\|_{L^\infty(B_{1})}+\|\phi\|_{L^\infty(B_{1})}\le2
,$$
giving that~\eqref{90oihfe8yhfh24inftyB2} also holds true for~$j=0$ (as long as we choose~$M\ge2$).

Thus, to perform the inductive step, we now suppose that~\eqref{90oihfe8yhfh24inftyB2} and~\eqref{90oihfe8yhfh24inftyBT}
are satisfied for an index~$j$ and we wish to establish them for the index~$j+1$.
To this end, it is convenient to consider the rescaled solution
$$ B_{1}\ni x\longmapsto V(x):=\frac{v(\rho^j x)-P_j(\rho^j x)}{M\rho^{(2+\alpha)j}}.$$
Using the inductive assumption on~\eqref{90oihfe8yhfh24inftyB2}, we see that
\begin{equation}\label{9u02ojrwr25suppartialB1w}
\sup_{x\in B_1}|V(x)|=\sup_{y\in B_{\rho^j}}\frac{|v(y)-P_j(y)|}{M\rho^{(2+\alpha)j}}\le1.\end{equation}
Moreover, using the inductive assumption on~\eqref{90oihfe8yhfh24inftyBT},
for every~$x\in B_1$,
$$ \Delta V(x)=
\frac{\phi(\rho^j x)-\Delta P_j(\rho^j x)}{M\rho^{\alpha j}}
=\frac{\phi(\rho^j x)-\phi(0)}{M\rho^{\alpha j}}$$
and accordingly
\begin{equation}\label{9u02ojrwr25suppartialB1w2} 
|\Delta V(x)|\le\frac{\|\phi\|_{C^\alpha(B_{1})}\rho^{\alpha j}}{M\rho^{\alpha j}}\le\frac1M.\end{equation}
We also remark that~$B_1(x_0)\subseteq B_2\Subset B_4$,
hence~$V\in C(\overline{B_1})$. 
We thus consider the solution~$w\in C^2(B_1)\cap C(\overline{B_1})$ of
$$ \begin{dcases}
\Delta w=0&{\mbox{ in }}B_1,\\
w=V&{\mbox{ on }}\partial B_1.
\end{dcases}$$
We point out that~$w$ can be constructed, for instance,
by using the Poisson Kernel of the ball, as discussed in Theorems~\ref{POIBALL1} and~\ref{POIBALL}.
Using the Weak Maximum Principle in Corollary~\ref{WEAKMAXPLE}(iii), combined with~\eqref{9u02ojrwr25suppartialB1w}, we know that
$$ \|w\|_{L^\infty(B_1)}\le\sup_{\partial B_1} |w|=\sup_{\partial B_1} |V|\le1.$$
As a consequence, we infer from
Cauchy's Estimates (recall Theorem~\ref{CAUESTIMTH}) that
\begin{equation}\label{s254a24con24352s32eq6427ue} \| w\|_{C^3(B_{1/2})}\le C\|w\|_{L^\infty(B_{1})}\le C.\end{equation}
Let now~$W_\pm(x):=w(x)-V(x)\pm\frac{|x|^2}{2nM}$. Recalling~\eqref{9u02ojrwr25suppartialB1w2},
in~$B_1$ we have that~$ \Delta W_+=-\Delta V+\frac{1}M\ge0$ and similarly~$\Delta W_-\le0$.
As a result, the Weak Maximum Principle in Corollary~\ref{WEAKMAXPLE} gives that, for every~$x\in B_1$,
\begin{eqnarray*}&&
w(x)-V(x)+\frac{|x|^2}{2nM}=W_+(x)\le\sup_{y\in\partial B_1} W_+(y)=\sup_{y\in\partial B_1}\frac{|y|^2}{2nM}\le\frac1M\\
{\mbox{and }}&&
w(x)-V(x)-\frac{|x|^2}{2nM}=W_-(x)\ge\inf_{y\in\partial B_1} W_-(y)=\inf_{y\in\partial B_1}\left(-\frac{|y|^2}{2nM}\right)\ge-\frac1M,
\end{eqnarray*}
and accordingly
\begin{equation}\label{912jjfn4f5t865y2B325r24ho} \|w-V\|_{L^\infty(B_1)}\le\frac2M.\end{equation}
Now we define the quadratic polynomial
$$\widetilde{P}(x):=w(0)+\nabla w(0)\cdot x+\frac12 D^2w(0)x\cdot x$$
and we observe that~$\Delta\widetilde{P}(x)=\Delta w(0)=0$ for all~$x\in\R^n$.

In addition, by~\eqref{s254a24con24352s32eq6427ue} and a Taylor expansion, for all~$x\in B_{1/2}$,
\begin{eqnarray*}
|w(x)-\widetilde{P}(x)|\le \|D^3 w\|_{L^\infty(B_{1/2})}|x|^3\le C|x|^3.
\end{eqnarray*}
We stress that the above~$C$ depends only on~$n$ and~$\alpha$.
In particular, since~$\alpha\in(0,1)$, we can choose~$\rho$ universally small such that~$C\rho^3\le\frac{\rho^{2+\alpha}}2$,
with~$C$ as above,
and conclude that
\begin{equation}\label{034lefr3ac324r3232h32o2}
\sup_{x\in B_\rho}|w(x)-\widetilde{P}(x)|\le  C\rho^3\le\frac{\rho^{2+\alpha}}2.
\end{equation}
We can now choose~$M$ sufficiently large such that~$\frac2M\le\frac{\rho^{2+\alpha}}2$ and deduce from~\eqref{912jjfn4f5t865y2B325r24ho}
and~\eqref{034lefr3ac324r3232h32o2} that
\begin{equation}\label{034lefr3ac324r3232h32o28900-0} \|V-\widetilde{P}\|_{L^\infty(B_\rho)}\le\|V-w\|_{L^\infty(B_\rho)}+\|w-\widetilde{P}\|_{L^\infty(B_\rho)}\le
\frac2M+\frac{\rho^{2+\alpha}}2\le \rho^{2+\alpha}.
\end{equation}
We now scale back by defining
$$P_{j+1}(x):=P_j(x)+M\rho^{(2+\alpha)j}\widetilde{P}(\rho^{-j}x).$$ In this way we have that
\begin{equation}\label{034lefr3ac324r3232h32o2-o02-1} \Delta P_{j+1}(x)=
\Delta P_j(x)+M\rho^{ \alpha j}\Delta\widetilde{P}(\rho^{-j}x)=\phi(0)+0=\phi(0).\end{equation}
Furthermore, owing to~\eqref{034lefr3ac324r3232h32o28900-0},
\begin{equation}\label{034lefr3ac324r3232h32o2-o02-2}\begin{split}&
\sup_{x\in B_{\rho^{j+1}}} |v(x)-P_{j+1}(x)|=\sup_{x\in B_{\rho^{j+1}}}|M\rho^{(2+\alpha)j} V(\rho^{-j} x)+P_j(x)-P_{j+1}(x)|\\&\qquad=
M\rho^{(2+\alpha)j}
\sup_{y\in B_{\rho}}\left| V(y)+\frac{P_j(\rho^j y)-P_{j+1}(\rho^j y)}{M\rho^{(2+\alpha)j}}\right|\\&\qquad=
M\rho^{(2+\alpha)j}
\sup_{y\in B_{\rho}}\left| V(y)-\widetilde{P}(y)\right|
\\&\qquad\le
M\rho^{(2+\alpha)(j+1)}.
\end{split}
\end{equation}
The inductive step is thereby complete:
the observations in~\eqref{034lefr3ac324r3232h32o2-o02-1} and~\eqref{034lefr3ac324r3232h32o2-o02-2}
give~\eqref{90oihfe8yhfh24inftyB2} and~\eqref{90oihfe8yhfh24inftyBT} for the index~$j+1$, as desired.
\end{proof}

The method of polynomial approximation can also be used to
obtain local estimates at the boundary: to show this, we give here another proof of Theorem~\ref{SCHAUDER-BOUTH}.

\begin{proof}[Another proof of Theorem~\ref{SCHAUDER-BOUTH}]
For computational convenience, \label{234556df0323iurjw924324jfJAMAnJHOSo3o24f32tofXZwf}
up to a covering argument (recall Lemma~\ref{COVERINGARG}), we will prove Theorem~\ref{SCHAUDER-BOUTH} with~$B_1^+$
replaced by~$B_4^+$.
Moreover, as observed in~\eqref{omMAEwgNGLt4IONkjgfnbdsdtyuijk9024i5SZERSD},
it suffices to prove Theorem~\ref{SCHAUDER-BOUTH} when~$g$ vanishes identically.
Also, as done in~\eqref{HYPETheoremSCHAUDER-INTE}, it is not restrictive to assume that
\begin{equation}\label{DSikoljfoi3y98tyiuw98743ytgfeBVUjsbdfnU} \|u\|_{L^\infty(B_4^+)}\le1\qquad{\mbox{and}}\qquad \|f\|_{C^\alpha(B_4^+)}\le1.\end{equation}
In addition, by Corollary~\ref{i0fU824aiugfepZAkkZfPXBOUN-AGG-P}, we can obtain
Theorem~\ref{SCHAUDER-BOUTH} by proving that
that there exists~$\rho\in\left(0,\frac12\right)$ such that for every~$x_0\in B_{1/2}^+$
and every~$j\in\N$
there exists a polynomial~$P_{j,x_0}$ of degree at most~$2$ such that
\begin{equation}\label{90oihfe8yhfh24inftyB1-NUOVA}\begin{split}&
|u(x)-P_{j,x_0}(x)|\le M\rho^{(2+\alpha)j}\\&{\mbox{for all}}\,x\in B_{\rho^j}^+(x_0):=\{x=(x',x_n)\in\R^n{\mbox{ s.t. $
|x-x_0|<\rho^j$ and~$x_n>0$\}}},\end{split}\end{equation}
where~$M>0$ depends only on~$n$ and~$\alpha$.

To prove~\eqref{90oihfe8yhfh24inftyB1-NUOVA}, we will argue by induction over~$j\in\N$
and we will prove additionally that we can choose~$P_{j,x_0}$ such that
\begin{equation}\label{90oihfe8yhfh24inftyBT-NUOVA}
\begin{split}&\Delta P_{j,x_0}(x)=f(x_0)\quad{\mbox{ for every }}x\in\R^n\\
{\mbox{and, when~$x_{0,n}\in(0,\rho^j)$, }}\quad\;
&P_{j,x_0}(x',0)=0\quad{\mbox{ for every }}x'\in\R^{n-1}.
\end{split}
\end{equation}
To perform the inductive argument, we choose~$\rho\in\left(0,\frac12\right)$ conveniently small and~$M\ge1$ conveniently large in what follows. When~$j=0$ we choose~$P_{0,x_0}(x):=\frac{f(x_0)\,x_n^2}{2}$; in this way, equation~\eqref{90oihfe8yhfh24inftyBT-NUOVA}
holds true for~$j=0$ and, for all~$x\in B_{1}^+(x_0)$,
$$ |u(x)-P_{0,x_0}(x)|\le \|u\|_{L^\infty(B_{1}^+(x_0))}+\|f\|_{L^\infty(B_{1}^+(x_0))}\le2
,$$
thanks to~\eqref{DSikoljfoi3y98tyiuw98743ytgfeBVUjsbdfnU},
giving that~\eqref{90oihfe8yhfh24inftyB1-NUOVA} also holds true for~$j=0$ provided that~$M\ge2$.

We now perform the inductive step towards the proof of~\eqref{90oihfe8yhfh24inftyB1-NUOVA}
and~\eqref{90oihfe8yhfh24inftyBT-NUOVA},
assuming that these claims hold true for an index~$j$ and aiming at establishing them for the index~$j+1$.

To accomplish this goal,
we let
$$ \overline{B_1}\cap\{x_{0,n}+\rho^j x_n\ge0\}\ni x\longmapsto V(x):=\frac{u(x_0+\rho^j x)-P_{j,x_0}(x_0+\rho^j x)}{M\rho^{(2+\alpha)j}}.$$
By the inductive assumption,
\begin{equation}\label{9u02ojrwr25suppartialB1w-NUOVAt}
\sup_{x\in B_1\cap\{x_{0,n}+\rho^j x_n>0\}}|V(x)|=\sup_{y\in B^+_{\rho^j}(x_0)}\frac{|u(y)-P_{j,x_0}(y)|}{M\rho^{(2+\alpha)j}}\le1\end{equation}
and, for every~$x\in B_1\cap\{x_{0,n}+\rho^j x_n>0\}$,
$$ \Delta V(x)=
\frac{f(x_0+\rho^j x)-\Delta P_j(x_0+\rho^j x)}{M\rho^{\alpha j}}
=\frac{f(x_0+\rho^j x)-f(x_0)}{M\rho^{\alpha j}},$$
leading to
\begin{equation}\label{9u02ojrwr25suppartialB1w2-NUOVA} 
|\Delta V(x)|\le\frac{\|f\|_{C^\alpha(B_{4}^+)}\,\rho^{\alpha j}}{M\rho^{\alpha j}}\le\frac1M.\end{equation}

Now we distinguish two cases,
namely Case~1 corresponding to~$x_{0,n}\ge\rho^{j}$
and Case~2 corresponding to~$x_{0,n}\in(0,\rho^{j})$.
In Case~2 it is convenient to extend~$V$ by an odd reflection across
the hyperplane~$\{x_{0,n}+\rho^j x_n=0\}$
by exploiting the second line of~\eqref{90oihfe8yhfh24inftyBT-NUOVA}
(notice that the second line of~\eqref{90oihfe8yhfh24inftyBT-NUOVA}
becomes void in Case~1). Namely, in Case~2 we observe that, by inductive assumption,
if~$x=(x',x_n)\in \overline{B_1}\cap\{x_{0,n}+\rho^j x_n=0\}$ then
$$ V(x)=\frac{u(x_0'+\rho^j x',0)-P_{j,x_0}(x_0'+\rho^j x',0)}{M\rho^{(2+\alpha)j}}=0.$$
As a result, in Case~2, if~$x=(x',x_n)\in\{x_{0,n}+\rho^j x_n<0\}$ with~$\left(x',-\frac{2x_{0,n}}{\rho^j}-x_n\right)\in\overline{B_1}$ we define
\begin{equation}\label{S254A3L4O23JDILBCHNEDDOJCH} V(x',x_n):=-V\left(x',-\frac{2x_{0,n}}{\rho^j}-x_n\right)\end{equation}
and we have that~$V$ is continuous in~$\overline{\mathcal{B}}$, where~${\mathcal{B}}$
is symmetrization of~${B_1}\cap\{x_{0,n}+\rho^j x_n>0\}$ through
the hyperplane~$\{x_{0,n}+\rho^j x_n=0\}$, i.e.
\begin{equation}\label{BCA2-9oihkf} {\mathcal{B}}:=\Big(
{B_1}\cap\{x_{0,n}+\rho^j x_n\ge0\}\Big)\cup\left\{
x=(x',x_n)\in\{x_{0,n}+\rho^j x_n<0\} {\mbox{ s.t. }}\left(x',-\frac{2x_{0,n}}{\rho^j}-x_n\right)\in B_1
\right\},\end{equation}
see Figure~\ref{SCHPIPEDItangeFI}.

We observe that, by taking~$x_n:=-\frac{x_{0,n}}{\rho^j}$ in~\eqref{S254A3L4O23JDILBCHNEDDOJCH},
we find that, in Case~2, 
$$V\left(x',-\frac{x_{0,n}}{\rho^j}\right)=-V\left(x',-\frac{x_{0,n}}{\rho^j}\right)$$ and thus
\begin{equation}\label{2ioKLiugfwj8OIJNDkfwEFGD9-XCVV}
V\left(x',-\frac{x_{0,n}}{\rho^j}\right)=0.
\end{equation}

\begin{figure}
  \centering
  \includegraphics[width=.95\linewidth]{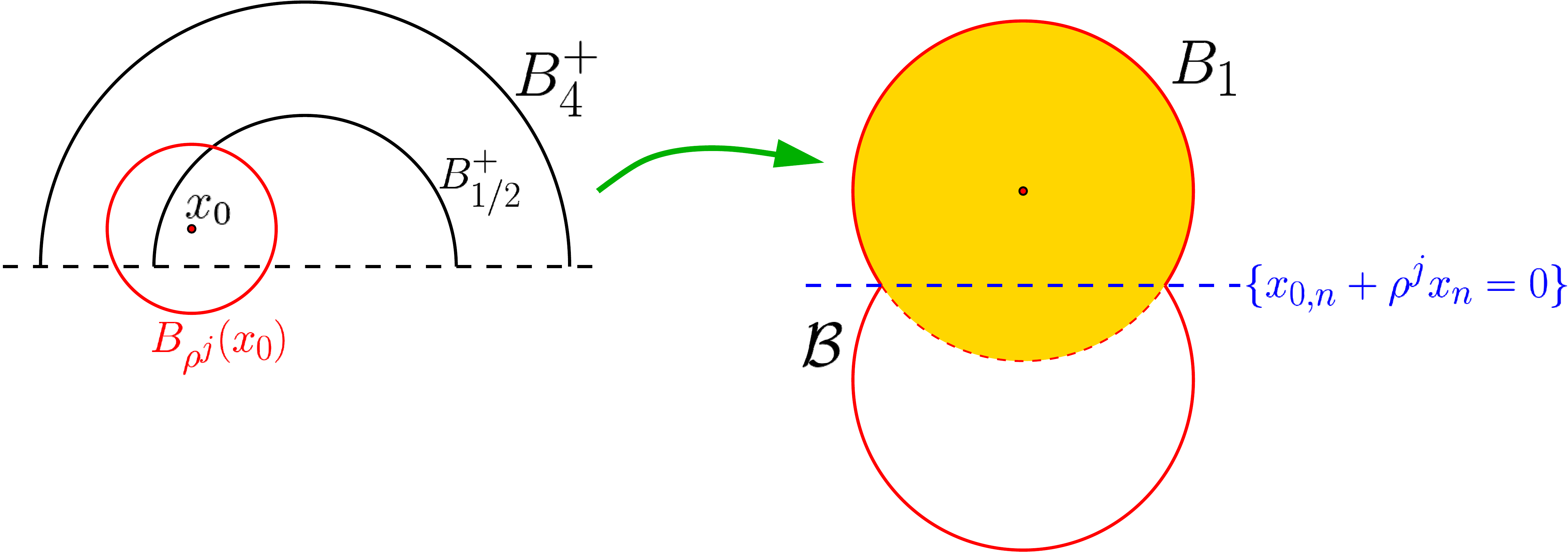}
 \caption{\sl The geometry involved in the alternative
 proof of Theorem~\ref{SCHAUDER-BOUTH} when~$x_{0,n}\in(0,\rho^j)$ (Case~2).}\label{SCHPIPEDItangeFI}
\end{figure}

\begin{figure}
  \centering
  \includegraphics[width=.95\linewidth]{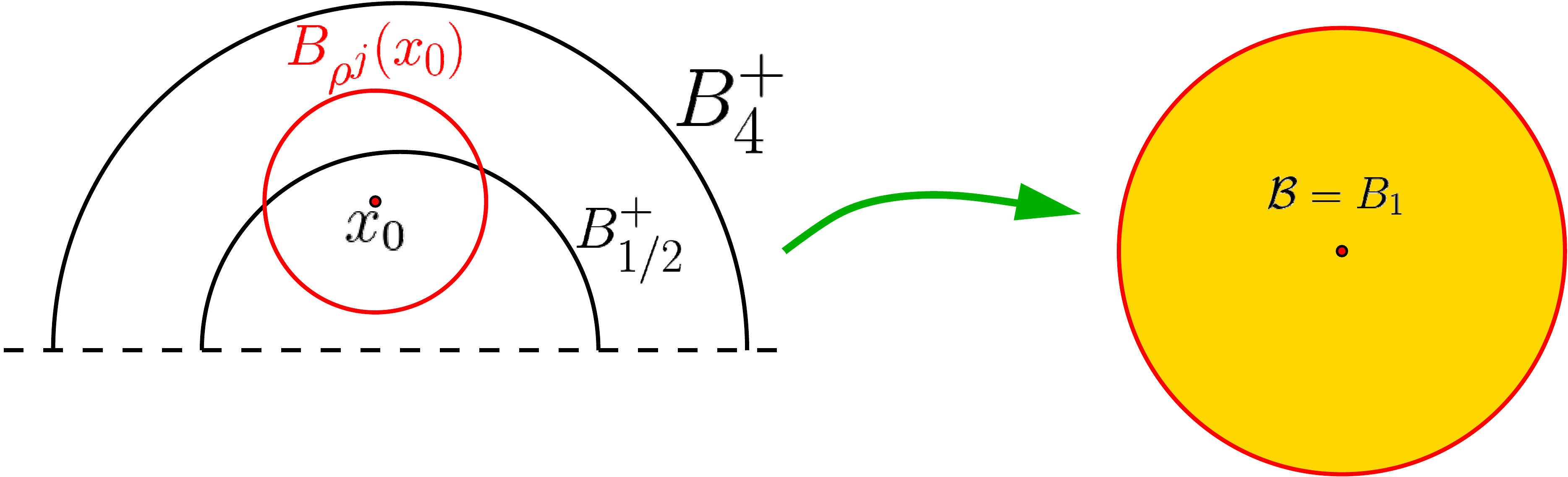}
 \caption{\sl The geometry involved in the alternative
 proof of Theorem~\ref{SCHAUDER-BOUTH} when~$x_{0,n}\ge\rho^j$ (Case~1).}\label{CA1SCHPIPEDItangeFI}
\end{figure}

When we are in Case~1, this odd reflection method is not needed:
to make the notation uniform between the two cases, we can replace~\eqref{BCA2-9oihkf}
in Case~1 with~${\mathcal{B}}:=B_1$
(and no need to extend~$V$), see Figure~\ref{CA1SCHPIPEDItangeFI}.

In both cases, by~\eqref{9u02ojrwr25suppartialB1w-NUOVAt},
\begin{equation}\label{9u02ojrwr25suppartialB1w-NUOVA}
\sup_{x\in {\mathcal{B}}}|V(x)|\le1.\end{equation}

Now, making use of the Perron method (see Corollary~\ref{S-coroEXIS-M023}
and observe that~${\mathcal{B}}$ satisfies the exterior cone condition both in Case~1 and in Case~2),
we consider the solution~$w\in C^2(B_1)\cap C(\overline{B_1})$ of
$$ \begin{dcases}
\Delta w=0&{\mbox{ in }}{\mathcal{B}},\\
w=V&{\mbox{ on }}\partial {\mathcal{B}}.
\end{dcases}$$
It follows from the Weak Maximum Principle in Corollary~\ref{WEAKMAXPLE}(iii) and~\eqref{9u02ojrwr25suppartialB1w-NUOVA} that
$$ \|w\|_{L^\infty({\mathcal{B}})}\le\sup_{\partial {\mathcal{B}}} |w|=\sup_{\partial {\mathcal{B}}} |V|\le1$$
and accordingly, by Cauchy's Estimates (see Theorem~\ref{CAUESTIMTH}),
\begin{equation}\label{s254a24con24352s32eq6427ue-NUOVA} \| w\|_{C^3(
{B_{1/2}}\cap\{x_{0,n}+\rho^j x_n\ge0\}
)}\le C\|w\|_{L^\infty({\mathcal{B}})}\le C.\end{equation}
We remark that, if we are in Case~2, then~$w$ is odd symmetric with respect to the hyperplane~$\{x_{0,n}+\rho^j x_n=0\}$, namely, for all~$x\in{\mathcal{B}}$,
\begin{equation}\label{2ioKLiugfwj8OIJNDkfwEFGD9}
w(x',x_n)=-w\left(x',-\frac{2x_{0,n}}{\rho^j}-x_n\right)
\end{equation}
Indeed, in Case~2, if
$$ W(x',x_n):=w(x',x_n)+w\left(x',-\frac{2x_{0,n}}{\rho^j}-x_n\right),$$
we have that~$\Delta W=0+0=0$ in~${\mathcal{B}}$ and, in light of~\eqref{S254A3L4O23JDILBCHNEDDOJCH}, we have that
$W(x',x_n)=V(x',x_n)+V\left(x',-\frac{2x_{0,n}}{\rho^j}-x_n\right)=0$ for all~$(x',x_n)\in\partial{\mathcal{B}}$.
{F}rom these observations and the uniqueness result in Corollary~\ref{UNIQUENESSTHEOREM} we deduce that~$W$
vanishes identically in~${\mathcal{B}}$, from which we obtain~\eqref{2ioKLiugfwj8OIJNDkfwEFGD9},
as desired.

As a consequence of~\eqref{2ioKLiugfwj8OIJNDkfwEFGD9}, by taking~$x_n:=-\frac{x_{0,n}}{\rho^j}$ there,
we also find that, in Case~2, 
$$w\left(x',-\frac{x_{0,n}}{\rho^j}\right)=-w\left(x',-\frac{x_{0,n}}{\rho^j}\right)$$ and therefore
\begin{equation}\label{2ioKLiugfwj8OIJNDkfwEFGD9-XC}
w\left(x',-\frac{x_{0,n}}{\rho^j}\right)=0.
\end{equation}

Now we define~$W_\pm(x):=w(x)-V(x)\pm\frac{|x|^2}{2nM}$ and we utilize~\eqref{9u02ojrwr25suppartialB1w2-NUOVA},
finding that~$ \Delta W_+=-\Delta V+\frac{1}M\ge0$, and similarly~$\Delta W_-\le0$, in~$B_1\cap\{ x_{0,n}+\rho^j x_n>0\}$
(the careful reader will notice that this set coincides with~$B_1$ in Case~1).
These observations and the Weak Maximum Principle in Corollary~\ref{WEAKMAXPLE} give that, for every~$x\in B_1\cap\{ x_{0,n}+\rho^j x_n>0\}$,
$$
w(x)-V(x)+\frac{|x|^2}{2nM}=W_+(x)\le\sup_{y\in\partial (B_1\cap\{ x_{0,n}+\rho^j x_n>0\})} W_+(y)
=\sup_{y\in\partial( B_1\cap\{ x_{0,n}+\rho^j x_n>0\})}\frac{|y|^2}{2nM}\le\frac1M$$ and $$
w(x)-V(x)-\frac{|x|^2}{2nM}=W_-(x)\ge\inf_{y\in\partial (B_1\cap\{ x_{0,n}+\rho^j x_n>0\})} W_-(y)=\inf_{y\in\partial( B_1\cap\{ x_{0,n}+\rho^j x_n>0\})}\left(-\frac{|y|^2}{2nM}\right)\ge-\frac1M.$$
Notice that here above we used the fact that~$w=V$ on~$\partial B_1$ in Case 1, while in Case 2 we have exploited
the facts that~$w=V$ on~$\partial B_1\cap \{ x_{0,n}+\rho^j x_n\ge 0\}$ and that~$w=V=0$
on~$B_1\cap\{ x_{0,n}+\rho^j x_n=0\}$ thanks to~\eqref{2ioKLiugfwj8OIJNDkfwEFGD9-XCVV}
and~\eqref{2ioKLiugfwj8OIJNDkfwEFGD9-XC}.

As a result,
\begin{equation}\label{912jjfn4f5t865y2B325r24ho-NUOVA} \|w-V\|_{L^\infty(B_1\cap\{ x_{0,n}+\rho^j x_n>0\})}\le\frac2M.\end{equation}
Now we set
$$ \R^n\supset{\mathcal{B}}\ni\zeta_j:=\begin{dcases}
\displaystyle\left(0,\dots,0,-\frac{x_{0,n}}{\rho^j}\right)
&{\mbox{if~$x_{0,n}\in(0,\rho^{j+1})$,}}\\
(0,\dots,0) &{\mbox{otherwise}}
\end{dcases}$$
and define the quadratic polynomial
$$
\widetilde{P}(x):=w(\zeta_j)+\nabla w(\zeta_j)\cdot (x-\zeta_j)+\frac12 D^2w(\zeta_j)(x-\zeta_j)\cdot( x-\zeta_j).$$
Note that~$\Delta\widetilde{P}(x)=\Delta w(\zeta_j)=0$ for all~$x\in\R^n$.

Furthermore, by~\eqref{2ioKLiugfwj8OIJNDkfwEFGD9-XC},
$$ w\left(x',-\frac{x_{0,n}}{\rho^j}\right)=\partial_i w\left(x',-\frac{x_{0,n}}{\rho^j}\right)=
\partial_{ik}w\left(x',-\frac{x_{0,n}}{\rho^j}\right)=0$$
for all~$i\in\{1,\dots,n-1\}$ and~$k\in\{1,\dots,n\}$.

Consequently, if~$x_{0,n}\in(0,\rho^{j+1})$ then
$$ w(\zeta_j)=\partial_i w(\zeta_j)=
\partial_{ik}w(\zeta_j)=0.$$
For this reason, if~$x_{0,n}\in(0,\rho^{j+1})$ then
\begin{equation}\label{0oikj9iyhk3roiyo32ihgf43frtew0nj1}
\widetilde{P}\left(x',-\frac{x_{0,n}}{\rho^j}\right)=
w(\zeta_j)+\nabla w(\zeta_j)\cdot (x',0)+\frac12 D^2w(\zeta_j)(x',0)\cdot( x',0)
=0
\end{equation}
for all~$x'\in\R^{n-1}$.

Besides, in view of~\eqref{s254a24con24352s32eq6427ue-NUOVA} and a Taylor expansion around~$\zeta_j$,
we see that, for all~$x\in B_{1/2}\cap \{x_{0,n}+\rho^jx_n>0 \}$,
\begin{equation}\label{IHKNBSWQHENLWMDLQWMDOxJDA2}
|w(x)-\widetilde{P}(x)|\le \|D^3 w\|_{L^\infty(B_{1/2} \cap \{x_{0,n}+\rho^jx_n>0 \})}|x-\zeta_j|^3\le C|x-\zeta_j|^3.
\end{equation}
We stress that the above~$C$ depends only on~$n$ and~$\alpha$.

We also notice that for all~$x\in B_\rho$ we have that
\begin{equation}\label{IHKNBSWQHENLWMDLQWMDOxJDA1}
|x-\zeta_j|\le 2\rho
\end{equation}
Indeed, if~$x_{0,n}\in(0,\rho^{j+1})$ it holds that
$$ |x-\zeta_j|^2=\left|x+\frac{x_{0,n}}{\rho^j}e_n\right|^2=
|x'|^2+\left(x_n+\frac{x_{0,n}}{\rho^j}\right)^2=
|x|^2+\frac{x_{0,n}^2}{\rho^{2j}}+\frac{2x_n x_{0,n}}{\rho^j}\le
\rho^2+\frac{\rho^{2(j+2)}}{\rho^{2j}}+\frac{2\rho^{j+2}}{\rho^j}=4\rho^2,
$$
otherwise~$|x-\zeta_j|=|x|\le\rho$, which proves~\eqref{IHKNBSWQHENLWMDLQWMDOxJDA1}.

Combining~\eqref{IHKNBSWQHENLWMDLQWMDOxJDA2} and~\eqref{IHKNBSWQHENLWMDLQWMDOxJDA1}, we see that, for all~$x\in B_\rho \cap \{x_{0,n}+\rho^jx_n>0 \}$,
$$ |w(x)-\widetilde{P}(x)|\le C\rho^3,$$
up to renaming~$C$.
In particular, since~$\alpha\in(0,1)$, we can choose~$\rho$ universally small such that~$C\rho^3\le\frac{\rho^{2+\alpha}}2$,
with~$C$ as above, therefore
\begin{equation}\label{034lefr3ac324r3232h32o2-NUOVA}
\sup_{x\in B_\rho \cap \{x_{0,n}+\rho^jx_n>0 \}}|w(x)-\widetilde{P}(x)|\le C\rho^3\le\frac{\rho^{2+\alpha}}2.
\end{equation}
Now we choose~$M$ sufficiently large such that~$\frac2M\le\frac{\rho^{2+\alpha}}2$ and deduce from~\eqref{912jjfn4f5t865y2B325r24ho-NUOVA}
and~\eqref{034lefr3ac324r3232h32o2-NUOVA} that
\begin{equation}\label{034lefr3ac324r3232h32o28900-0-NUOVA} \begin{split}
&\|V-\widetilde{P}\|_{L^\infty(B_\rho
\cap\{ x_{0,n}+\rho^j x_n>0\})}\le\|V-w\|_{L^\infty(B_\rho\cap\{ x_{0,n}+\rho^j x_n>0\})}+\|w-\widetilde{P}\|_{L^\infty(B_\rho\cap\{ x_{0,n}+\rho^j x_n>0\})}\\&\qquad\qquad\qquad\le
\frac2M+\frac{\rho^{2+\alpha}}2\le \rho^{2+\alpha}.\end{split}
\end{equation}
We now scale back by defining
$$P_{j+1,x_0}(x):=P_{j,x_0}(x)+M\rho^{(2+\alpha)j}\widetilde{P}\left(\frac{x-x_0}{\rho^j}\right).$$ Hence, we have that
\begin{equation}\label{034lefr3ac324r3232h32o2-o02-1-NUOVA} \Delta P_{j+1,x_0}(x)=
\Delta P_{j,x_0}(x)+M\rho^{ \alpha j}\Delta\widetilde{P}\left(\frac{x-x_0}{\rho^j}\right)=f(x_0)+0=f(x_0).\end{equation}
Furthermore, if~$x_{0,n}\in(0,\rho^{j+1})$ we can utilize~\eqref{0oikj9iyhk3roiyo32ihgf43frtew0nj1}
and conclude that
\begin{equation}\label{SCQUNQKAMDLKAMMSOHD}
P_{j+1,x_0}(x',0)=P_{j,x_0}(x',0)+M\rho^{(2+\alpha)j}\widetilde{P}
\left(\frac{x'-x_0'}{\rho^j},-\frac{x_{0,n}}{\rho^j}\right)
=0+0=0.
\end{equation}
Furthermore, by~\eqref{034lefr3ac324r3232h32o28900-0-NUOVA},
\begin{equation*}\begin{split}&
\sup_{x\in B_{\rho^{j+1}}^+(x_0)} |u(x)-P_{j+1,x_0}(x)|=\sup_{x\in B_{\rho^{j+1}}^+(x_0)}\left|M\rho^{(2+\alpha)j} V\left(\frac{x-x_0}{\rho^j}\right)+P_{j,x_0}(x)-P_{j+1,x_0}(x)\right|\\&\qquad=
M\rho^{(2+\alpha)j}
\sup_{y\in B_{\rho}\cap\{ x_{0,n}+\rho^j y_n>0\}}\left| V(y)+\frac{P_{j,x_0}(\rho^j y+x_0)-P_{j+1,x_0}(\rho^j y+x_0)}{M\rho^{(2+\alpha)j}}\right|\\&\qquad=
M\rho^{(2+\alpha)j}
\sup_{y\in B_{\rho}\cap\{ x_{0,n}+\rho^j y_n>0\}}\left| V(y)-\widetilde{P}(y)\right|
\\&\qquad\le
M\rho^{(2+\alpha)(j+1)}.
\end{split}
\end{equation*}
This and~\eqref{034lefr3ac324r3232h32o2-o02-1-NUOVA}, together
with~\eqref{SCQUNQKAMDLKAMMSOHD}
when~$x_{0,n}\in(0,\rho^{j+1})$,
establish~\eqref{90oihfe8yhfh24inftyB1-NUOVA}
and~\eqref{90oihfe8yhfh24inftyBT-NUOVA} for the index~$j+1$, as desired, thus completing the proof of the inductive step.
\end{proof}

\section{Schauder estimates for equations in nondivergence form}\label{KPJMDarJMSiaCLiMMSciAmoenbvD8Sijf-23}

Now we turn our attention to a more general class of elliptic equations.
Specifically, rather than looking only at the Laplace operator,
we will consider equations in the form
\begin{equation}\label{NONDIFORM}
\sum_{i,j=1}^n a_{ij}(x)\partial_{ij}u(x)+\sum_{i=1}^n b_i(x)\partial_i u(x)+c(x)u(x)=f(x).\end{equation}
We always suppose that the coefficients~$a_{ij}$ give rise to a symmetric matrix, namely
\begin{equation}\label{NONDIFORM2}
a_{ij}=a_{ji}.\end{equation}
This assumption can be taken essentially without loss of generality, since
if~\eqref{NONDIFORM2} did not hold, it would suffice to define~$\widetilde a_{ij}:=\frac{a_{ij}+a_{ji}}{2}$
and exploit the symmetry of the second derivatives of~$u$ to recast~\eqref{NONDIFORM}
into a similar equation in which the coefficients~$a_{ij}$ are replaced by the symmetric ones~$\widetilde a_{ij}$.

The structural assumption that we take on equation~\eqref{NONDIFORM} is the ``ellipticity of the coefficients~$a_{ij}$'':
more precisely, we will always assume that there exist~$\lambda\in(0,+\infty)$ and~$\Lambda\in[\lambda,+\infty)$ such that
for every~$\xi=(\xi_1,\dots,\xi_n)\in\partial B_1$
\begin{equation}\label{ELLIPTIC}
\sum_{i,j=1}^n a_{ij}(x)\xi_i\xi_j\in [\lambda ,\Lambda ].
\end{equation}
{F}rom the algebraic viewpoint this condition states that all the eigenvalues of the symmetric matrix~$a_{ij}$
are positive, bounded away from zero and bounded (compare with the notion of ellipticity
presented in footnote~\ref{CLASSIFICATIONFOOTN} on page~\pageref{CLASSIFICATIONFOOTN}).

\begin{figure}
  \centering
  \includegraphics[width=.3\linewidth]{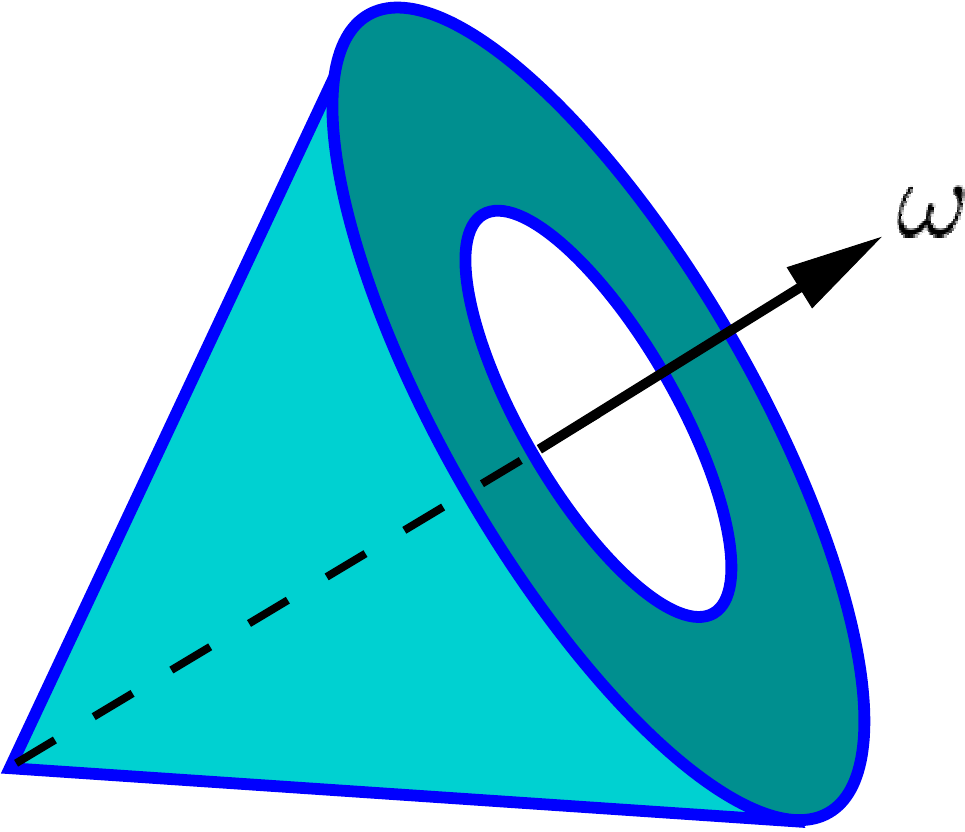}
 \caption{\sl Graphical intuition for the ellipticity condition in~\eqref{ELLIPTIC}: the cone~$K_\omega$
 in~\eqref{ELLIPTICONO}.}\label{T536EcdhPOaIJodidpo3tju24ylkiliRELLIPTICONO2543uj8890j09092po-okf-3edk-3edhf-2-fh564}
\end{figure}

A geometric ``visualization'' of~\eqref{ELLIPTIC} can be obtained, at least when~$n=2$,
considering the following map~$E$ from the space of symmetric matrices into~$\R^3$: \label{CLASSIFICATIONFOOTN3}
$$ E:\left( \begin{matrix} a_{11} & a_{12}\\a_{12}&a_{22}\end{matrix}\right)\longmapsto
\left( \sqrt{2}\,a_{12},\,\frac{a_{11}-a_{22}}{\sqrt2},\,\frac{a_{11}+a_{22}}{\sqrt2}\right).$$
We observe that~$E$ is an isometry, in the sense that if~$A$ is the matrix corresponding to~$\left( \begin{matrix} a_{11} & a_{12}\\a_{12}&a_{22}\end{matrix}\right)$ then
$$ |E(A)|^2=2a_{12}^2+\frac{(a_{11}-a_{22})^2+(a_{11}+a_{22})^2}{2}=2a_{12}^2+a_{11}^2+a_{22}^2=|A|^2.$$
Furthermore, $E$ reverts the trace of the matrix product into the scalar product, since
\begin{eqnarray*} E(A)\cdot E(\widetilde A)&=&\left( \sqrt{2}\,a_{12},\,\frac{a_{11}-a_{22}}{\sqrt2},\,\frac{a_{11}+a_{22}}{\sqrt2}\right)\cdot
\left( \sqrt{2}\,\widetilde a_{12},\,\frac{\widetilde a_{11}-\widetilde a_{22}}{\sqrt2},\,\frac{\widetilde a_{11}+\widetilde a_{22}}{\sqrt2}\right)\\&
=&2a_{12}\widetilde a_{12}+\frac{(a_{11}-a_{22})(\widetilde a_{11}-\widetilde a_{22})
+(a_{11}+a_{22})(\widetilde a_{11}+\widetilde a_{22})}2
\\&=&2a_{12}\widetilde a_{12}+a_{11}\widetilde a_{11}+a_{22}\widetilde a_{22}\\&=&{\rm tr}\, (A\widetilde A)
.\end{eqnarray*}
In particular, if~$\omega:=\frac{E(A)}{|E(A)|}$ and~$
{\mathcal{R}}$ is the space of the matrices~$\Xi$ of the form~$\left( \begin{matrix} \xi_1^2 & \xi_1\xi_2\\
\xi_1\xi_2&\xi_{2}^2\end{matrix}\right)$, then condition~\eqref{ELLIPTIC} gives that, for all~$\Xi\in{\mathcal{R}}\setminus\{0\}$,
\begin{eqnarray*}&&
\omega\cdot \frac{E(\Xi)}{|E(\Xi)|}=\frac{E(A)\cdot E(\Xi)}{|E(A)|\,|E(\Xi)|}=
\frac{{\rm tr}\, (A\Xi)}{|A|\,|\Xi|}=
\frac{a_{11}\,\xi_{1}^2+a_{22}\,\xi_2^2+2a_{12}\,\xi_1\xi_2}{|A|\,\sqrt{\xi_1^4+\xi_2^4+2\xi_1^2\xi_2^2}}\\
&&\qquad=\frac{a_{11}\,\xi_{1}^2+a_{22}\,\xi_2^2+2a_{12}\,\xi_1\xi_2}{|A|\,(\xi_1^2+\xi_2^2)}=
\frac{1}{|A|\,|\xi|^2}\,\sum_{i,j=1}^2 a_{ij}\xi_i\xi_j\in\left[\frac{\lambda}{|A|},\,\frac{\Lambda}{|A|}\right].
\end{eqnarray*}
Accordingly, the ellipticity condition gives in this situation that~$E$ maps~${\mathcal{R}}$ into the cone
\begin{equation}\label{ELLIPTICONO}K_\omega:=\left\{p\in\R^3\setminus\{0\} {\mbox{ s.t. }} \omega\cdot\frac{p}{|p|}\in\left[\frac{\lambda}{|A|},\,\frac{\Lambda}{|A|}\right]
\right\}\cup\{0\},\end{equation}
see Figure~\ref{T536EcdhPOaIJodidpo3tju24ylkiliRELLIPTICONO2543uj8890j09092po-okf-3edk-3edhf-2-fh564}.\medskip

Equations such as the one in~\eqref{NONDIFORM} are often referred to with the name
of ``equations in nondivergence form''\index{nondivergence form} (to distinguish them from the ones surfacing
from a variational structure, compare e.g. with equation~\eqref{DITRO}).
There are sundry reasons for addressing equation~\eqref{NONDIFORM} in such a generality.
One of the most important, in our personal opinion, \label{TOWLemma92o3jwelg02jfj-qiewjdfnoewfh9ewohgveioqNONDIFO234RM}
is that the Laplace operator is not maintained under ``nice transformations of the space''
(except for translations and rotations, in view of~\eqref{TRAL} and Corollary~\ref{KSMD:ROTAGSZOKA})
and it would be instead be desirable to consider a class of equations which is preserved
under natural families of transformations.

With this respect, we point out the following result:

\begin{lemma}\label{92o3jwelg02jfj-qiewjdfnoewfh9ewohgveioqNONDIFO234RM}
Let~$\Omega\subseteq\R^n$ be open and~$u\in C^2(\Omega)$.
Let~$a_{ij}$, $b_i$, $c$, $f\in C^\alpha(\Omega)$ for some~$\alpha\in[0,1]$ and assume that
\begin{equation}\label{oj-LEM-EQUJSN-comSIJMCAi-mbiJAQ-ledm-1}
\sum_{i,j=1}^n a_{ij}(x)\partial_{ij}u(x)+\sum_{i=1}^n b_i(x)\partial_i u(x)+c(x)u(x)=f(x)\qquad{\mbox{for every }}\;x\in\Omega.\end{equation}
Let~$T:\R^n\to\R^n$ be an invertible diffeomorphism of class~$C^{2,\alpha}$ with inverse of class~$C^{2,\alpha}$.
For each~$y\in T(\Omega)$ let~$\widetilde u(y):=u(T^{-1}(y))$.
Then,
\begin{equation}\label{oj-LEM-EQUJSN-comSIJMCAi-mbiJAQ-ledm-2}
\sum_{i,j=1}^n \widetilde a_{ij}(y)\partial_{ij}
\widetilde u(y)+\sum_{i=1}^n \widetilde b_i(y)\partial_i \widetilde u(y)+
\widetilde c(y)\widetilde u(y)=\widetilde f(y)\qquad{\mbox{for every }}\;y\in T(\Omega)\end{equation}
for suitable~$\widetilde a_{ij}$, $\widetilde b_i$, $\widetilde c$, $\widetilde f\in C^\alpha(T(\Omega))$.

Additionally, if~$c\ge0$ (respectively, $c\le0$) in~$\Omega$, then~$\widetilde c\ge0$
(respectively, $\widetilde c\le0$) in~$T(\Omega)$.

Furthermore, if there exist~$\Lambda\ge\lambda>0$ such that
\begin{equation}\label{92o3rlambdaxi2}
\lambda|\xi|^2\le\sum_{i,j=1}^n a_{ij}\xi_i\xi_j\le\Lambda|\xi|^2\qquad{\mbox{for all }}\,\xi=(\xi_1,\dots,\xi_n)\in\R^n,
\end{equation}
then
\begin{equation}\label{92o3rlambdaxi2b}
\widetilde\lambda|\xi|^2\le\sum_{i,j=1}^n \widetilde a_{ij}\xi_i\xi_j\le
\widetilde \Lambda|\xi|^2\qquad{\mbox{for all }}\,\xi=(\xi_1,\dots,\xi_n)\in\R^n,
\end{equation}
for suitable~$\widetilde\Lambda\ge\widetilde\lambda>0$.
\end{lemma}

\begin{proof} We implicitly assume that~$u$ is evaluated at~$x$
and~$\widetilde u$ is evaluated at~$y=T(x)=(T_1(x),\dots,T_n(x))$. Also, $T$ will be implicitly
assumed to be evaluated at~$x$. With this notation,
\begin{eqnarray*}&&
\partial_i u=\sum_{k=1}^n \partial_k \widetilde u \,\partial_i T_k\\
{\mbox{and }}&&\partial_{ij} u=\sum_{k,m=1}^n\partial_{km} \widetilde u\, \partial_i T_k\,\partial_j T_m
+\sum_{k=1}^n\partial_k \widetilde u \,\partial_{ij} T_k.
\end{eqnarray*}
As a result, setting
\begin{equation}\label{PSSKTROFNMJCGLQMGCSDVC}\begin{split}
&\widetilde a_{km}(y):=\sum_{i,j=1}^n a_{ij}(T^{-1}(y))\,\partial_i T_k\,\partial_j T_m ,\\
&\widetilde b_k(y):=
\sum_{i,j=1}^n a_{ij}(T^{-1}(y))\, \partial_{ij} T_k+\sum_{i=1}^n b_i (T^{-1}(y))\,\partial_i T_k,\\
&\widetilde c(y):=c(T^{-1}(y))\\
{\mbox{and }}\quad &\widetilde f(y):=f(T^{-1}(y)),
\end{split}\end{equation}
in view of~\eqref{oj-LEM-EQUJSN-comSIJMCAi-mbiJAQ-ledm-1} (and omitting the variable by following the
previous convention) we have that
\begin{eqnarray*}
\widetilde f-\widetilde c \widetilde u&=&f-cu\\
&=&\sum_{i,j=1}^n a_{ij} \partial_{ij}u +\sum_{i=1}^n b_i\partial_i u\\
&=&\sum_{i,j,k,m=1}^n a_{ij} \partial_i T_k\,\partial_j T_m \,\partial_{km} \widetilde u +
\sum_{i,j,k=1}^n a_{ij}\partial_{ij} T_k \,\partial_k \widetilde u
+\sum_{i,k=1}^n b_i \partial_i T_k \,\partial_k \widetilde u
\\ &=& 
\sum_{k,m=1}^n \widetilde a_{ij} \partial_{km} \widetilde u +
\sum_{k=1}^n \widetilde b_{k}  \partial_k \widetilde u,\end{eqnarray*}
which is~\eqref{oj-LEM-EQUJSN-comSIJMCAi-mbiJAQ-ledm-2}.

Also, the sign of~$\widetilde c$ is clearly the same as the one (if any) of~$c$.
Furthermore, if~$\xi\in\R^n$ we
let~$\eta=(\eta_1,\dots,\eta_n)\in\R^n$ be defined as
$$\eta_i:=\sum_{k=1}^n \partial_i T_k \xi_k.$$
We thus find that
\begin{eqnarray*}
\sum_{k,m=1}^n\widetilde a_{km}\xi_k\xi_m
=\sum_{i,j,k,m=1}^n a_{ij} \partial_i T_k\,\partial_j T_m\,\xi_k\,\xi_m=\sum_{i,j=1}^n a_{ij}\eta_i\eta_j.
\end{eqnarray*}
Hence, if~\eqref{92o3rlambdaxi2} holds true,
\begin{equation}\label{192o3inftyOmegaxi}
\sum_{k,m=1}^n\widetilde a_{km}\xi_k\xi_m
\in[\lambda|\eta|^2,\Lambda|\eta|^2].
\end{equation}
We also observe that, if~$(DT)_{ij}=\partial_i T_j$, then~$\eta=DT\,\xi$ and therefore~$|\eta|\le (\|T\|_{C^1(\Omega)}+1)|\xi|$.
Moreover, since~$\xi=(DT^{-1})\eta$, we have that~$|\xi|\le(\|T^{-1}\|_{C^1(\Omega)}+1)|\eta|$. {F}rom these observations
and~\eqref{192o3inftyOmegaxi} we deduce that
$$ \sum_{k,m=1}^n\widetilde a_{km}\xi_k\xi_m
\in\left[\frac{\lambda}{(\|T^{-1}\|_{C^1(\Omega)}+1)^2}|\xi|^2,\,(\|T\|_{C^1(\Omega)}+1)^2
\Lambda|\xi|^2\right],
$$
which is~\eqref{92o3rlambdaxi2b}.
\end{proof}

The fact that the structure of equation~\eqref{NONDIFORM} is preserved
under suitable transformations (as guaranteed by Lemma~\ref{92o3jwelg02jfj-qiewjdfnoewfh9ewohgveioqNONDIFO234RM})
will play a crucial role in the proof of Theorem~\ref{ThneGiunfBB032t4jNKS-s3i4}.
\medskip

Besides this important property of preservation of elliptic structures,
as warrated by Lemma~\ref{92o3jwelg02jfj-qiewjdfnoewfh9ewohgveioqNONDIFO234RM},
we recall that equations in nondivergence form also naturally arise
from the linearization of more complicated, and nonlinear, equations,
see e.g. the discussion on page~\pageref{LINEMONGEALI}.\medskip

Another natural motivation for equation~\eqref{NONDIFORM}, which is also helpful to develop some familiarity
with this equation, arises in the random movement of a biological population with density~$u$ living
in some environment~$\Sigma$ (this can be seen as a revisitation of the
model presented in Section~\ref{CHEMOTX}). Suppose that~$\Sigma$ presents hills and valleys, \label{0uojf29249-45kpkfdSmd11493839429efv}
for instance it can be modeled like a graph over~$\R^n$ of the form
$$ \Sigma=\{ X=(x,y)\in\R^n\times\R {\mbox{ s.t. }}y=\gamma(x)\},$$
for some (nice) function~$\gamma:\R^n\to\R$ (the biological case
would correspond to~$n=2$). In this notation, we think that~$u$ is a function of
time~$t\in\R$ and position~$X=(x,\gamma(x))\in\Sigma$, hence of~$t\in \R$ and~$x\in\R^n$.

We assume that the biological population reproduces itself,
hence at a given unit of time~$\tau$ the population increases proportionally to the existing
population: at each time step, the proliferation effects would thus add to the population~$u$ an additional term~$c\tau u$:
the case $c=0$ would correspond to no reproduction at all, and we may think that~$c$ is also a function of~$x$
if the environment presents regions that are either more (higher~$c$) or less (smaller~$c$)
favorable for reproduction.

We can also suppose that the population is bred, therefore at each unit of time~$\tau$ some additional population~$\tau f$
is added into the environment. Again, $f$ could depend on~$x$, since the breeders may identify strategic zones
where to insert the newcomers. Also, in some regions~$f$ could be positive,
corresponding to breeders inserting new individuals, or negative, corresponding to breeders killing part of the existing population e.g. for food production (and of course~$f$ can change its sign if both these breeding
activities take place).
The case of~$f$ vanishing identically is also allowed and corresponds to a free biological population
with no breeders.

\begin{figure}
  \centering
  \includegraphics[width=.4\linewidth]{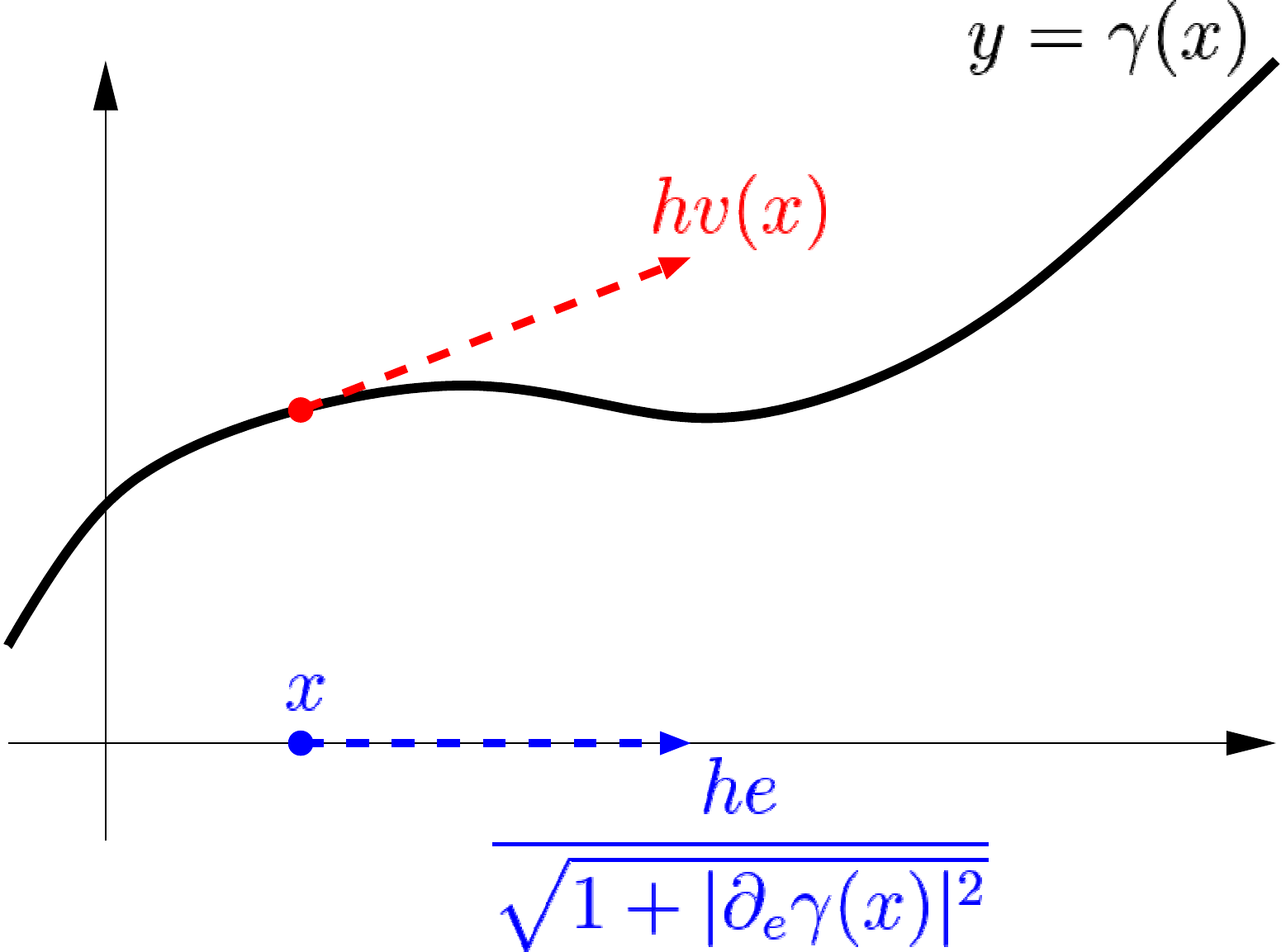}
 \caption{\sl Random displacement of the population leaving on~$\Sigma$ (with~$h$ to be sent to~$0$
 in the limit).}\label{T536EcdhPOaIJodidpo3tju24ylkiliR2543uj8890j09092po-okf-3edk-3edhf-2-fh564}
\end{figure}

We suppose that the biological population moves randomly across the environment.
At each time step, the individuals pick randomly a direction~$e\in \partial B_1$
and move of a unit step~$h$ along~$\Sigma$ in the direction indicated by~$e$.
We stress that, for small~$h$, when an individual moves on the surface~$\Sigma$ by a space length~$h$ in direction
indicated by~$e$,
its projection on~$\R^n$ is moved by a space length~$\frac{h}{\sqrt{1+|\partial_e\gamma|^2}}$
in the horizontal direction~$e$, up to~$o(h)$, 
and the spatial displacement of the population in~$\R^{n+1}$ is~$
hv$, up to~$o(h)$, where
$$ v(x):=\left(\frac{e}{ \sqrt{1+|\partial_e\gamma(x)|^2}},\;
\frac{\partial_e\gamma(x)}{\sqrt{1+|\partial_e\gamma(x)|^2}}
\right),$$
since,
using the notation~$X(x):=(x,\gamma(x))$,
\begin{eqnarray*} X\left(x+\frac{he}{\sqrt{1+|\partial_e\gamma(x)|^2}}\right)-X(x)
&=&\left(\frac{he}{\sqrt{1+|\partial_e\gamma(x)|^2}},\;
\gamma\left( x+\frac{he}{\sqrt{1+|\partial_e\gamma(x)|^2}}\right)-\gamma(x)
\right)\\&=& \left(\frac{he}{\sqrt{1+|\partial_e\gamma(x)|^2}},\;
\frac{h\partial_e\gamma(x)}{\sqrt{1+|\partial_e\gamma(x)|^2}}
\right)+o(h)
\\&=&
h v(x)+o(h),
\end{eqnarray*}
see Figure~\ref{T536EcdhPOaIJodidpo3tju24ylkiliR2543uj8890j09092po-okf-3edk-3edhf-2-fh564}.

With respect to this observation, up to higher orders that we neglect,
the individuals corresponding to the population density~$u(x,t)$
are moved at time~$t+\tau$ by the random process into~$u\left(x+\frac{he}{\sqrt{1+|\partial_e\gamma(x)|^2}},t\right)$,
averaged over~$e\in\partial B_1$, namely
$$ u(x,t+\tau)=\fint_{\partial B_1}u\left(x+\frac{he}{\sqrt{1+|\partial_e\gamma(x)|^2}},t\right)\, d{\mathcal{H}}^{n-1}_e.$$
For flat surfaces in which~$\nabla\gamma=0$ the above would correspond to the classical (isotropic) random walk.

We can also assume that the individuals are mildly subject to gravity, in the sense
that the gravity force tends to let them slightly slide down in the direction of maximal slope
of the environment~$\Sigma$. To model this effect, for instance, we can suppose that at each unit
of time the  individuals corresponding to the population density~$u(x,t)$
are moved by gravity into~$u(x-h^2 d(x)\nabla\gamma(x),t+\tau)$.
The parameter~$d$ models a ``grip'', which could also vary from region to region,
namely~$d=0$ is a perfect adherent situation in which no sliding occurs, and instead a large value of~$d$
would provide a slippy surface.

In this complex (but relatively simple, after all) model, the density at time~$t+\tau$
is thus the superposition of a number of effects,
such as proliferation, breeding, random movement and drifting (one can also remove
some of this effects to concentrate on a simpler model) and we can therefore write that
\begin{equation}\label{OJS0w5pjrqu32adt4gfdd0843k-014kjdd}
\begin{split}
u(x,t+\tau)=& c\tau u(x,t)+\tau f(x,t)+\fint_{\partial B_1} u\left(x+\frac{he}{\sqrt{1+|\partial_e\gamma(x)|^2}},t\right)\,d{\mathcal{H}}^{n-1}_e
\\&\quad+u(x+h^2 d(x)\nabla\gamma(x),t)-u(x,t).\end{split}\end{equation}
By an odd symmetry cancellation, we see that
\begin{eqnarray*}&&
\fint_{\partial B_1} u\left(x+\frac{he}{\sqrt{1+|\partial_e\gamma(x)|^2}},t\right)\,d{\mathcal{H}}^{n-1}_e\\&=&
\fint_{\partial B_1} \Bigg[
u(x,t)+\nabla u(x,t)\cdot\frac{he}{\sqrt{1+|\partial_e\gamma(x)|^2}}\\&&\qquad\qquad
+\frac12\,D^2u(x,t)
\frac{he}{\sqrt{1+|\partial_e\gamma(x)|^2}}\cdot\frac{he}{\sqrt{1+|\partial_e\gamma(x)|^2}}
\Bigg]\,d{\mathcal{H}}^{n-1}_e+o(h^2)\\
&=& u(x,t)+\frac{h^2}2\sum_{i,j=1}^n\fint_{\partial B_1} \frac{\partial_{ij} u(x,t)\,e_i\,e_j}{1+|\partial_e\gamma(x)|^2}
\,d{\mathcal{H}}^{n-1}_e+o(h^2)\\&=&
u(x,t)+ h^2 \sum_{i,j=1}^n a_{ij}(x)\partial_{ij} u(x,t)+o(h^2),
\end{eqnarray*}
where
\begin{equation}\label{4215e21q2uwera32t335i332o5rqn} a_{ij}(x):=\frac12\,
\fint_{\partial B_1} \frac{e_i\,e_j}{1+|\partial_e\gamma(x)|^2}
\,d{\mathcal{H}}^{n-1}_e.\end{equation}

Moreover,
$$ u(x+h^2 d(x)\nabla\gamma(x),t)-u(x,t)=h^2d(x)\nabla u(x,t)\cdot\nabla\gamma(x)+o(h^2)=
h^2\sum_{i=1}^n b_i(x)\partial_i u(x,t)+o(h^2)
,$$
where
$$ b(x):= d(x)\nabla \gamma(x).$$
By these observations,
we can recast~\eqref{OJS0w5pjrqu32adt4gfdd0843k-014kjdd} into
\begin{equation*}
\begin{split}&
u(x,t+\tau)=c\tau u(x,t)+\tau f(x,t)+
u(x,t)+ h^2 \sum_{i,j=1}^n a_{ij}(x)\partial_{ij} u(x,t)+h^2\sum_{i=1}^n b_i(x)\partial_i u(x,t)+o(h^2)
.\end{split}\end{equation*}
Hence, taking the term~$u(x,t)$ to the left hand side, dividing by~$\tau$ and choosing
a quadratic space-time scaling unit, i.e.~$\tau:=h^2$, sending~$h\searrow0$ (or equivalently~$\tau\searrow0$),
we obtain
$$ \partial_tu(x,t)=cu(x,t)+f(x,t)+
\sum_{i,j=1}^n a_{ij}(x)\partial_{ij} u(x,t)+\sum_{i=1}^n b_i(x)\partial_i u(x,t).$$
With respect to this, we see that equation~\eqref{NONDIFORM}
describes the stationary states of this model.
\medskip

We have also encountered a general form of elliptic equation for biological species possibly subject to
drift and chemotaxis, recall~\eqref{OJS-PJDN-0IHGDOIUGDBV02ujrf22NDF}.\medskip

Another example of general elliptic equation was encountered
in the stationary Fokker-Planck setting, recalling equation~\eqref{FOKPLA}.
\medskip

Other motivations for equation~\eqref{NONDIFORM} arise in differential geometry,
when one studies partial differential equations on manifolds (compare e.g. with the Laplace-Beltrami
operator in local coordinates as presented in~\eqref{Jns89ijgm5jjgfjgjg2}).\medskip

Now we discuss some interior estimates 
for equation~\eqref{NONDIFORM} which can be seen as the natural extension\footnote{The assumption~$u\in C^{2,\alpha}(B_1)$
does play a role in the proof of Theorem~\ref{THM:0-T536EcdhPOaIJodidpo3tju24ylkiliR2543uj8890j09092po-okf-3edk-3edhf-2-fh564},
but it can be weakened to the more natural hypothesis~$u\in C^2(B_1)$. Weakening this assumption at this stage
is not completely trivial and we will accomplish this further generality in the forthcoming
Corollary~\ref{THM:0-T536EcdhPOaIJodidpo3tju24ylkiliR2543uj8890j09092po-okf-3edk-3edhf-2-fh564COR}.}
of Theorem~\ref{SCHAUDER-INTE}:

\begin{theorem}\label{THM:0-T536EcdhPOaIJodidpo3tju24ylkiliR2543uj8890j09092po-okf-3edk-3edhf-2-fh564}
Let~$a_{ij}$, $b_i$, $c$, $f\in C^\alpha(B_1)$ for some~$\alpha\in(0,1)$. Assume the ellipticity 
condition in~\eqref{ELLIPTIC}.
Let~$u\in C^{2,\alpha}(B_1)$ be a solution of
\begin{equation*}
\sum_{i,j=1}^n a_{ij} \partial_{ij}u+\sum_{i=1}^n b_i\partial_i u+cu=f
\quad{\mbox{in }}\,B_1.\end{equation*}
Then, there exists~$C>0$, depending only on~$n$,
$\alpha$, $a_{ij}$, $b_i$ and~$c$, such that
\begin{equation}\label{92i45kk5emmaabeCOVERINGARG} \|u\|_{C^{2,\alpha}(B_{1/2})}\le C\,\Big(\|u\|_{L^\infty(B_1)}+\|f\|_{C^\alpha(B_1)}\Big).\end{equation}
\end{theorem}

We also have a counterpart of the local estimates at the boundary given in Theorem~\ref{SCHAUDER-BOUTH},
in the halfball notation stated in~\eqref{9oj3falphaB10},
according to the following result:

\begin{theorem}\label{THM:0-T536EcdhPOaIJodidpo3tju24ylkiliR2543uj8890j09092po-okf-3edk-3edhf-2-fh564-LEATB}
Let~$a_{ij}$, $b_i$, $c$, $f\in C^\alpha(B_1^+)$ and~$g\in
C^{2,\alpha}(B_1^0)$ for some~$\alpha\in(0,1)$. Assume the ellipticity 
condition in~\eqref{ELLIPTIC}.
Let~$u\in C^{2,\alpha}(B_1^+)\cap C(B_1^+\cup B_1^0)$ be a solution of
\begin{equation*}
\begin{dcases} \sum_{i,j=1}^n a_{ij} \partial_{ij}u+\sum_{i=1}^n b_i\partial_i u+cu=f\quad{\mbox{in }}\,B_1^+,\\
u=g \quad{\mbox{on }}\,B_1^0.\end{dcases}\end{equation*}
Then, there exists~$C>0$, depending only on~$n$,
$\alpha$, $a_{ij}$, $b_i$ and~$c$, such that
\begin{equation*} \|u\|_{C^{2,\alpha}(B_{1/2}^+)}\le C\,\Big(\|u\|_{L^\infty(B_1^+)}+
\|g\|_{C^{2,\alpha}(B_1^0)}+\|f\|_{C^\alpha(B_1^+)}\Big).\end{equation*}
\end{theorem}

To deal with the proof of regularity results such as
Theorems~\ref{THM:0-T536EcdhPOaIJodidpo3tju24ylkiliR2543uj8890j09092po-okf-3edk-3edhf-2-fh564}
and~\ref{THM:0-T536EcdhPOaIJodidpo3tju24ylkiliR2543uj8890j09092po-okf-3edk-3edhf-2-fh564-LEATB}, it is often convenient to 
deal with intermediate estimates in which weighted norms
in a ball are bounded by a small constant times
suitable weighted norms
in larger balls. While in principle the weighted norms
in larger balls do not seem helpful to conclude a meaningful bound,
the smallness of the constant in front may allow one
to ``reabsorb'' the inconvenient terms.
This methodology relies on some fine analysis
of scaling and weighted norms and it is of general use (see~\cite[pages~398--399]{MR1459795}):
for our scopes, we limit ourselves to the following result:

\begin{proposition}\label{SVA-LS-LEOSKDNDO}
Let~$k\in\N$, $\alpha\in(0,1]$, $\vartheta\in[0,1]$, $\sigma\ge0$ and $R>0$.
There exists~$\e_0\in(0,1)$ depending only on~$n$ and $k$ such that the following holds true.

Let~$u\in C^{k,\alpha}(B_R)$.
Assume that
for every~$x_0\in B_R$ and every~$\rho\in\left( 0,\frac{R-|x_0|}{2}\right)$
\begin{equation}\label{Rstar8b-1}
{\mathcal{N}}_\vartheta(x_0,\rho)
\le \sigma+\e_0 \,{\mathcal{N}}_\vartheta(x_0,2\rho),
\end{equation}
where we used the notation
$$ {\mathcal{N}}_\vartheta(x_0,\rho):=
\sum_{j=0}^k\rho^j\|D^j u\|_{L^\infty(B_\rho(x_0))}
+\vartheta\rho^{k+\alpha} [u]_{C^{k,\alpha}(B_\rho(x_0))}.$$
Then, there exists~$C>0$, depending only on~$n$, $k$ and $R$,
such that
\begin{equation}\label{08uojbslj94tqquadqquad-1}
\sum_{j=0}^k \|D^j u\|_{L^\infty(B_{R/2})}
+\vartheta [u]_{C^{k,\alpha}(B_{R/2})}
\le C\sigma.\end{equation}
\end{proposition}

A counterpart of Proposition~\ref{SVA-LS-LEOSKDNDO} holds true in the
halfball notation stated in~\eqref{9oj3falphaB10} and generalized in~\eqref{9oj3falphaB10GENE}.
The technical details go as follows:

\begin{proposition}\label{SVA-LS-LEOSKDNDO-HB-2}
Let~$k\in\N$, $\alpha\in(0,1]$, $\vartheta\in[0,1]$, $\sigma\ge0$ and~$R>0$. 
There exists~$\e_0\in(0,1)$ depending only on~$n$ and~$k$
such that the following holds true.

Let~$u\in C^{k,\alpha}(B^+_R)$.
Assume that
for every~$x_0\in B^+_{R}$ and every~$\rho\in\left( 0,\frac{R-|x_0|}{2}\right)$
\begin{equation}\label{Rstar8b-2}
{\mathcal{N}}_\vartheta(x_0,\rho)
\le \sigma+\e_0 \,{\mathcal{N}}_\vartheta(x_0,2\rho),
\end{equation}
where we used the notation
$$ {\mathcal{N}}_\vartheta(x_0,\rho):=
\sum_{j=0}^k\rho^j\|D^j u\|_{L^\infty(B^+_\rho(x_0))}
+\vartheta\rho^{k+\alpha} [u]_{C^{k,\alpha}(B^+_\rho(x_0))}.$$
Then, there exists~$C>0$, depending only on~$n$, $k$
and~$R$, such that
\begin{equation}\label{08uojbslj94tqquadqquad-2}
\sum_{j=0}^k \|D^j u\|_{L^\infty(B^+_{R/2})}
+\vartheta [u]_{C^{k,\alpha}(B^+_{R/2})}
\le C\sigma.\end{equation}
\end{proposition}

We will prove Propositions~\ref{SVA-LS-LEOSKDNDO} and~\ref{SVA-LS-LEOSKDNDO-HB-2} at the same time,
by letting~${\mathcal{B}}_r(x)$ to be equal to~$B_r(x)$ in the setting of
Proposition~\ref{SVA-LS-LEOSKDNDO} and to be equal to~$B_r^+(x)$
in the setting of Proposition~\ref{SVA-LS-LEOSKDNDO-HB-2}:

\begin{proof}[Proof of Propositions~\ref{SVA-LS-LEOSKDNDO} and~\ref{SVA-LS-LEOSKDNDO-HB-2}]
Let
\begin{equation}\label{Qequ232-at33iowan} Q:=\sup_{{x_0\in {\mathcal{B}}_R}\atop{
\rho\in\left( 0,\frac{R-|x_0|}2\right)} }{\mathcal{N}}_\vartheta (x_0, \rho )\end{equation}
and notice that
$$Q\le \sum_{j=0}^k\rho^j\|D^j u\|_{L^\infty({\mathcal{B}}_R)}
+\vartheta\rho^{k+\alpha} [u]_{C^{k,\alpha}({\mathcal{B}}_R)}\le
\sum_{j=0}^kR^j\|D^j u\|_{L^\infty({\mathcal{B}}_R)}
+\vartheta R^{k+\alpha} [u]_{C^{k,\alpha}({\mathcal{B}}_R)}
<+\infty$$ since~$u\in C^{k,\alpha}({\mathcal{B}}_R)$.

Now we take~$x_0\in {\mathcal{B}}_R$ and~$\rho\in \left( 0,\frac{R-|x_0|}2\right)$.
We cover~${\mathcal{B}}_\rho(x_0)$ by a family of balls~${\mathcal{B}}_{\frac{\rho}{8}}(p_i)$
with~$p_i\in {\mathcal{B}}_\rho(x_0)$ and~$i\in\{1,\dots,N\}$; we observe that we can bound~$N$ from above
in dependence only of~$n$. 

We notice also that
$$|p_i|\le|x_0|+\rho<R-2\rho+\rho=R-\rho.$$
As a consequence, we have that~$p_i\in{\mathcal{B}}_R$ and~$\frac{\rho}2\in\left(0,\frac{R-|p_i|}2\right)$.
Hence, we can use either~\eqref{Rstar8b-1} or~\eqref{Rstar8b-2} and find that,
for all~$i\in\{1,\dots,N\}$,
\begin{equation*}
{\mathcal{N}}_\vartheta\left(p_i,\frac\rho{4}\right)
\le \sigma+\e_0{\mathcal{N}}_\vartheta\left(p_i, \frac\rho2\right)
\le \sigma+\e_0 Q
\end{equation*}
and therefore
\begin{equation}\label{0ojw23-ef2bR-rrightk}
\sum_{i=1}^N{\mathcal{N}}_\vartheta\left(p_i,\frac\rho{4}\right)
\le N\sigma+\e_0 NQ
.\end{equation}
Now we observe that, for all~$j\in\{0,\dots,k\}$,
\begin{equation}\label{04efeftrho0KSM-9yrhfb1}\begin{split}& \rho^j\|D^j u\|_{L^\infty({\mathcal{B}}_\rho(x_0))}\le
\rho^j\sum_{i=1}^N \|D^j u\|_{L^\infty({\mathcal{B}}_{\rho/8}(p_i))}=
8^j\sum_{i=1}^N \left(\frac{\rho}{8}\right)^j\|D^j u\|_{L^\infty({\mathcal{B}}_{\rho/8}(p_i))}\\
&\qquad\le 8^j\sum_{i=1}^N{\mathcal{N}}_\vartheta\left(p_i,\frac\rho{8}\right)\le
8^k\sum_{i=1}^N{\mathcal{N}}_\vartheta\left(p_i,\frac\rho{8}\right).
\end{split}\end{equation}

We also claim that
\begin{equation}\label{04efeftrho0KSM-9yrhfb2}
\vartheta\rho^{k+\alpha}[ u]_{C^{k,\alpha}({\mathcal{B}}_{\rho}(x_0))}\le C 8^{k+\alpha}\sum_{i=1}^N{\mathcal{N}}_\vartheta\left(p_{i},\frac\rho{4}\right),
\end{equation}
where~$C$ is a positive constant depending only on~$n$ and~$k$.
 
To check~\eqref{04efeftrho0KSM-9yrhfb2},
let~$x\in {\mathcal{B}}_\rho(x_0)$ and~$i_x\in\{1,\dots,N\}$ be such that~$x\in {\mathcal{B}}_{\frac\rho{8}}(p_{i_x})$.
We pick~$y\in {\mathcal{B}}_\rho(x_0)$ and distinguish two cases: if~$y\in {\mathcal{B}}_{\frac\rho{8}}(x)$, then~$x$, $y\in {\mathcal{B}}_{\frac\rho{4}}(p_{i_x})$
and therefore, if~$\beta\in\N^n$ with~$|\beta|=k$,
\begin{equation}\label{357Be8921urj-ZefFD}
\begin{split}&
\vartheta\rho^{k+\alpha}\frac{|D^\beta u(x)-D^\beta u(y)|}{|x-y|^\alpha}\le\vartheta\rho^{k+\alpha}
[u]_{C^{k,\alpha}( {\mathcal{B}}_{\rho/{4}}(p_{i_x}))}=\vartheta
4^{k+\alpha}
\left(\frac\rho{4}\right)^{k+\alpha}
[u]_{C^{k,\alpha}( {\mathcal{B}}_{\rho/{4}}(p_{i_x}))}\\&\qquad\qquad\le 4^{k+\alpha}
{\mathcal{N}}_\vartheta\left(p_{i_x},\frac\rho{4}\right)\le
4^{k+\alpha}\sum_{i=1}^N{\mathcal{N}}_\vartheta\left(p_{i},\frac\rho{4}\right).
\end{split}
\end{equation}
If instead~$y\not\in {\mathcal{B}}_{\frac\rho{8}}(x)$
we let~$i_y\in\{1,\dots,N\}$ be such that~$y\in {\mathcal{B}}_{\frac\rho{8}}(p_{i_y})$ and we find that,
if~$\beta\in\N^n$ with~$|\beta|=k$,
\begin{eqnarray*}&&
\vartheta\rho^{k+\alpha}\frac{|D^\beta u(x)-D^\beta u(y)|}{|x-y|^\alpha}\le
\rho^{k+\alpha}\frac{|D^\beta u(x)|+|D^\beta u(y)|}{(\rho/8)^\alpha}\\&&\qquad\le 8^\alpha\rho^k
\Big( \|D^\beta u\|_{L^\infty({\mathcal{B}}_{\rho/8}(p_{i_x}))}+\|D^\beta u\|_{L^\infty({\mathcal{B}}_{\rho/8}(p_{i_y}))}\Big)\\
&&\qquad\le
8^\alpha
\sum_{i=1}^N \rho^k\|D^\beta u\|_{L^\infty({\mathcal{B}}_{\rho/8}(p_i))}\le8^{k+\alpha}\sum_{i=1}^N{\mathcal{N}}_\vartheta\left(p_{i},\frac\rho{8}\right).
\end{eqnarray*}
By combining this and~\eqref{357Be8921urj-ZefFD}, we conclude that, for all~$x$, $y\in {\mathcal{B}}_\rho(x_0)$,
$$ \vartheta\rho^{k+\alpha}\frac{|D^\beta u(x)- D^\beta u(y)|}{|x-y|^\alpha}\le
8^{k+\alpha}\sum_{i=1}^N{\mathcal{N}}_\vartheta\left(p_{i},\frac\rho{4}\right).$$
Hence, summing up over all~$\beta\in\N^n$ such that~$|\beta|=k$ (recall the notation in
footnote~\ref{u9oj-NNp-OjtmRM-NosdTHSAMSATInfOdeN-N336OS235yS}), we obtain the desired result in~\eqref{04efeftrho0KSM-9yrhfb2}.

Thus, in light of~\eqref{04efeftrho0KSM-9yrhfb1} and~\eqref{04efeftrho0KSM-9yrhfb2},
$$ {\mathcal{N}}_\vartheta(x_0,\rho)\le C
\left( (k+1)8^k+8^{k+\alpha}\right)\sum_{i=1}^N{\mathcal{N}}_\vartheta\left(p_i,\frac\rho{4}\right)\le C
(k+2)8^{k+1}\sum_{i=1}^N{\mathcal{N}}_\vartheta\left(p_i,\frac\rho{4}\right).
$$
Thus, recalling~\eqref{0ojw23-ef2bR-rrightk},
$$  {\mathcal{N}}_\vartheta(x_0,\rho)\le C
(k+2)8^{k+1}\Big( N\sigma+\e_0 NQ\Big).$$
Taking the supremum over such~$x_0$ and~$\rho$, in view of~\eqref{Qequ232-at33iowan}
we find that
\begin{equation*}
Q\le C(k+2)8^{k+1}\Big( N\sigma+\e_0 NQ\Big).
\end{equation*}
We now choose~$\e_0:=\frac{1}{2NC(k+2)8^{k+1}}$ and we infer that
\begin{equation*}
Q\le C(k+2)8^{k+1}\left( N\sigma+\frac{Q}{2C(k+2)8^{k+1}}\right)=
(k+2)8^{k+1} N\sigma+\frac{Q}2
\end{equation*}
and therefore, reabsorbing one term into the left hand side,
\begin{equation}\label{098-iekjdj90uf3rjh794939} Q\le 2C(k+2)8^{k+1} N\sigma.\end{equation}
Hence, since, for all~$\rho\in\left(0,\frac{R}2\right)$,
$$ Q\ge{\mathcal{N}}_\vartheta (0, \rho )=
\sum_{j=0}^k\rho^j\|D^j u\|_{L^\infty({\mathcal{B}}_\rho)}
+\vartheta\rho^{k+\alpha} [u]_{C^{k,\alpha}({\mathcal{B}}_\rho)},
$$
we deduce from~\eqref{098-iekjdj90uf3rjh794939} that
\begin{eqnarray*} 2C(k+2)8^{k+1} N\sigma&\ge&\lim_{\rho\nearrow \frac{R}2}
\sum_{j=0}^k\rho^j\|D^j u\|_{L^\infty({\mathcal{B}}_\rho)}
+\vartheta\rho^{k+\alpha} [u]_{C^{k,\alpha}({\mathcal{B}}_\rho)}\\&=&
\sum_{j=0}^k \left(\frac{R}2\right)^j\|D^j u\|_{L^\infty({\mathcal{B}}_{R/2})}
+\vartheta \left(\frac{R}2\right)^{k+\alpha} [u]_{C^{k,\alpha}({\mathcal{B}}_{R/2})}.\end{eqnarray*}
{F}rom this, the claims in~\eqref{08uojbslj94tqquadqquad-1} and~\eqref{08uojbslj94tqquadqquad-2}
follow, as desired.
\end{proof}

With this auxiliary result on scaled estimates, we can proceed with
the proof of Theorem~\ref{THM:0-T536EcdhPOaIJodidpo3tju24ylkiliR2543uj8890j09092po-okf-3edk-3edhf-2-fh564}
by arguing as follows.

\begin{proof}[Proof of Theorem~\ref{THM:0-T536EcdhPOaIJodidpo3tju24ylkiliR2543uj8890j09092po-okf-3edk-3edhf-2-fh564}] The gist of the proof is to exploit the result for the Laplace operator
in Theorem~\ref{SCHAUDER-INTE} in a perturbative fashion, since, roughly speaking,
in tiny balls the coefficients~$a_{ij}$ vary very little from constant ones,
and the additional lower order terms induced by~$b_i$ and~$c$ produce, by scaling,
negligible quantities.

To make this argument work, we first observe that
\begin{equation}\label{021irkfe2m4m53a0294OVERINGAR3G}
\begin{split}&
{\mbox{Theorem~\ref{SCHAUDER-INTE} holds true if the Laplacian is replaced by an operator}}\\&
{\mbox{with constant
coefficients of the form}} \sum_{i,j=1}^n \overline{a}_{ij} \partial_{ij}u, \\&{\mbox{where~$\overline{a}_{ij}$ is constant
and fulfills the ellipticity condition in~\eqref{ELLIPTIC},}}\end{split}
\end{equation}
with constants depending also on the elliptic parameters~$\lambda$ and~$\Lambda$
in~\eqref{ELLIPTIC}.
Indeed, if~$u$ is a solution of~$\sum_{i,j=1}^n \overline{a}_{ij} \partial_{ij}u=f$ in~$B_1$
and~$\overline{a}_{ij}$ is as above, then we can define~$\overline{A}$ to be the matrix~$\{\overline{a}_{ij}\}_{i,j\in\{1,\dots,n\}}$ and~$M$ to be the square root of~$\overline{A}$ in the matrix sense, see e.g.~\cite[page~125]{MR2337395}.
We observe that
\begin{equation}\label{9uaois93oifeSDNDjnsOINSverWimdfKAly}
\left|\frac{Mx}{\sqrt{\Lambda}}\right|=
\sqrt{\frac{(Mx)\cdot(Mx)}{\Lambda}}=\sqrt{\frac{\overline{A}x\cdot x}{\Lambda}}.\end{equation}
Thus, if we consider the function~$\overline u(x):=u\left(\frac{Mx}{\sqrt\Lambda}\right)$ we have that, for all~$x\in B_1$,
$$ \left|\frac{Mx}{\sqrt{\Lambda}}\right|<
\sqrt{\frac{\Lambda}{\Lambda}}=1$$
and
$$ \Delta\overline u(x)=\frac{1}\Lambda\sum_{i,j,m=1}^n M_{ mi}M_{mj}\partial_{ij} u\left(\frac{Mx}{\sqrt\Lambda}\right)=
\frac{1}\Lambda \sum_{i,j=1}^n \overline{a}_{ij}\partial_{ij} u\left(\frac{Mx}{\sqrt\Lambda}\right)
=\frac{1}\Lambda\,f\left(\frac{Mx}{\sqrt\Lambda}\right)=:\overline{f}(x).$$
We can therefore employ Theorem~\ref{SCHAUDER-INTE} and conclude that
\begin{equation}\label{9uaois93oifeSDNDjnsOINSverWimdfK24Aly2} \|\overline u\|_{C^{2,\alpha}(B_{1/2})}\le C\,\Big(\|\overline u\|_{L^\infty(B_1)}+\|\overline f\|_{C^\alpha(B_1)}\Big).\end{equation}
Since, by~\eqref{9uaois93oifeSDNDjnsOINSverWimdfKAly},
$$ \left|\frac{Mx}{\sqrt{\Lambda}}\right|\ge\sqrt{\frac{\lambda| x|^2}{\Lambda}}
= \sqrt{\frac{\lambda}{\Lambda}}|x|,$$
up to renaming~$C$
we deduce from~\eqref{9uaois93oifeSDNDjnsOINSverWimdfK24Aly2} that
\begin{equation*} \| u\|_{C^{2,\alpha}\left(B_{
\sqrt{\frac{\lambda}{4\Lambda}}
}\right)}\le C\,\Big(\| u\|_{L^\infty(B_1)}+\| f\|_{C^\alpha(B_1)}\Big).\end{equation*}
{F}rom this and a covering argument (see Lemma~\ref{COVERINGARG}) we obtain~\eqref{021irkfe2m4m53a0294OVERINGAR3G},
as desired.

Now, given~$x_0\in B_{1/2}$ and~$\rho\in\left(0,\frac14\right)$, to be taken suitably small in the following,
we set \begin{eqnarray*}&& \widetilde{a}_{ij}(x):=
a_{ij}(x_0)-a_{ij}(x)\\{\mbox{and }}&& \widetilde{u}(x):=\frac{u(x_0+\rho x)}{\rho^2}.\end{eqnarray*}
For all~$x\in B_1$ we have that
\begin{equation}\label{hinHNDfOIKHDNinEDFDDIAbgerIGGjrtNVMVFkqTESKMFDFRT}\begin{split}
&\sum_{i,j=1}^n a_{ij}(x_0)\partial_{ij}\widetilde{u}(x)=
\sum_{i,j=1}^n a_{ij}(x_0)\partial_{ij}u(x_0+\rho x)\\
&\qquad=\sum_{i,j=1}^n a_{ij}(x_0+\rho x)\partial_{ij}u(x_0+\rho x)+
\sum_{i,j=1}^n \widetilde{a}_{ij}(x_0+\rho x)\partial_{ij}u(x_0+\rho x)\\
&\qquad=f(x_0+\rho x)-\sum_{i=1}^n b_i(x_0+\rho x)\partial_iu(x_0+\rho x)\\&\qquad\qquad\qquad-c(x_0+\rho x)u(x_0+\rho x)+
\sum_{i,j=1}^n \widetilde{a}_{ij}(x_0+\rho x)\partial_{ij}u(x_0+\rho x)\\&\qquad=:\widetilde f(x).
\end{split}\end{equation}
As a result, we deduce from~\eqref{021irkfe2m4m53a0294OVERINGAR3G} that
\begin{equation}\label{SM9ihw320D0-2i5rkwfracrhoight}
\|\widetilde u\|_{C^{2,\alpha}(B_{1/2})}\le C\,\Big(\|\widetilde u\|_{L^\infty(B_1)}+\|\widetilde f\|_{C^\alpha(B_1)}\Big).\end{equation}
We remark that if~$x\in B_1$ then~$|\widetilde u(x)|=\frac{|u(x_0+\rho x)|}{\rho^2}
\le \frac{\|u\|_{L^\infty(B_\rho(x_0))}}{\rho^2}$.
In addition, for every~$k\in\N$ we have that~$D^k u(y)=\rho^{2-k} D^k \widetilde{u}\left(\frac{y-x_0}\rho\right)$
hence, for all~$y$, $z\in B_{\rho/2}(x_0)$,
\begin{eqnarray*}
&&|u(y)|\le \rho^2\|\widetilde{u}\|_{L^\infty(B_{1/2})},\\
&&|\nabla u(y)|\le\rho\|\nabla\widetilde{u}\|_{L^\infty(B_{1/2})},\\
&&|D^2 u(y)|\le \|D^2\widetilde{u}\|_{L^\infty(B_{1/2})}\\
{\mbox{and }}&&
|D^2 u(y)-D^2u(z)|=\left|
D^2 \widetilde{u}\left(\frac{y-x_0}\rho\right)
-D^2 \widetilde{u}\left(\frac{z-x_0}\rho\right)
\right|\le C\rho^{-\alpha} [\widetilde{u}]_{C^{2,\alpha}(B_{1/2})}|x-y|^\alpha.
\end{eqnarray*}
It follows from these observations and~\eqref{SM9ihw320D0-2i5rkwfracrhoight} that
\begin{equation}\label{SM9ihw320D0-2i5rkwfracrhoight-BSjd}\begin{split}&
\sum_{k=0}^2 \rho^{k}\|D^ku\|_{L^\infty(B_{\rho/2}(x_0))}+
\rho^{2+\alpha} [u]_{C^{2,\alpha}(B_{\rho/2}(x_0))}\le
\rho^2\|\widetilde u\|_{C^{2,\alpha}(B_{1/2})}\\&\qquad
\le C\rho^2\Big(\| \widetilde u\|_{L^\infty(B_1)}+\|\widetilde f\|_{C^\alpha(B_1)}\Big)
\le C\,\Big(\| u\|_{L^\infty(B_1)}+\rho^{2}\|\widetilde f\|_{C^\alpha(B_1)}\Big)
.\end{split}\end{equation}
We observe that, for all~$x$, $y\in B_1$ and any function~$g$, 
\[|g(x_0+\rho x)-g(x_0+\rho y)|\le\rho^\alpha [g]_{C^\alpha(B_\rho(x_0))}|x-y|^\alpha.\]
As a consequence, letting~$b:=(b_1,\dots,b_n)$ and~$\widetilde A:=\{\widetilde a_{ij}\}_{i,j\in\{1,\dots,n\}}$,
\begin{equation}\label{SALK9aue298032ytrhfSHDIojgfo-2ot}\begin{split}
[\widetilde f]_{C^\alpha(B_1)}\le\;&
C \rho^\alpha\Big([f]_{C^\alpha(B_\rho(x_0))}+[b]_{C^\alpha(B_\rho(x_0))}\|\nabla u\|_{L^\infty(B_\rho(x_0))}+
\|b\|_{L^\infty(B_\rho(x_0))}[\nabla u]_{C^\alpha(B_\rho(x_0))}\\&\qquad
+[c]_{C^\alpha(B_\rho(x_0))}\|u\|_{L^\infty(B_\rho(x_0))}+
\|c\|_{L^\infty(B_\rho(x_0))}[u]_{C^\alpha(B_\rho(x_0))} \\&\qquad
+[\widetilde A]_{C^\alpha(B_\rho(x_0))}\| D^2u\|_{L^\infty(B_\rho(x_0))}
+\|\widetilde A\|_{L^\infty(B_\rho(x_0))}[D^2u]_{C^\alpha(B_\rho(x_0))}\Big)\\
\le\;& C\rho^\alpha \Big( \|f\|_{C^\alpha(B_\rho(x_0))}+\|u\|_{C^{2}(B_\rho(x_0))}+\|\widetilde A\|_{L^\infty(B_\rho(x_0))}[D^2u]_{C^\alpha(B_\rho(x_0))}
\Big)\\
\le\;& C\rho^\alpha \Big( \|f\|_{C^\alpha(B_\rho(x_0))}+\|u\|_{C^{2}(B_\rho(x_0))}+\rho^\alpha [D^2u]_{C^\alpha(B_\rho(x_0))}
\Big).\end{split}\end{equation}
Additionally,
\begin{eqnarray*}
\|\widetilde f\|_{L^\infty(B_1)}&\le&C\Big(
\|f\|_{L^\infty(B_\rho(x_0))}
+\|u\|_{C^1(B_\rho(x_0))}+\|\widetilde A\|_{L^\infty(B_\rho(x_0))}\|D^2u\|_{L^\infty(B_\rho(x_0))}\Big)
\\&\le& C\Big(
\|f\|_{L^\infty(B_\rho(x_0))}
+\|u\|_{C^1(B_\rho(x_0))}+\rho^\alpha\|D^2u\|_{L^\infty(B_\rho(x_0))}\Big).\end{eqnarray*}
This and~\eqref{SALK9aue298032ytrhfSHDIojgfo-2ot} entail that
$$ \|\widetilde f\|_{C^\alpha(B_1)}\le C \Big( \|f\|_{C^\alpha(B_\rho(x_0))}+\|u\|_{C^{1}(B_\rho(x_0))}+\rho^\alpha\|D^2u\|_{L^\infty(B_\rho(x_0))}+\rho^{2\alpha} [D^2u]_{C^\alpha(B_\rho(x_0))}\Big).$$
{F}rom this and~\eqref{SM9ihw320D0-2i5rkwfracrhoight-BSjd} we arrive at
\begin{equation}\label{KS0poqi038yrue983ythgfhGHSmnlkqfWgQEFb}
\begin{split}&
\sum_{k=0}^2 \rho^{k}\|D^ku\|_{L^\infty(B_{\rho/2}(x_0))}+
\rho^{2+\alpha} [u]_{C^{2,\alpha}(B_{\rho/2}(x_0))}\\
\le\;& C\,\Big(
\|f\|_{C^\alpha(B_1)}+
\| u\|_{L^\infty(B_1)}+
\rho^2\|\nabla u\|_{L^\infty(B_\rho(x_0))}+\rho^{2+\alpha}\|D^2u\|_{L^\infty(B_\rho(x_0))}+\rho^{2+2\alpha} [D^2u]_{C^\alpha(B_\rho(x_0))}
\Big)
\\ \le\;&
C\,\Big(
\|f\|_{C^\alpha(B_1)}+
\| u\|_{L^\infty(B_1)}\Big)+C\rho^\alpha\left(
\sum_{k=0}^2 \rho^{k}\|D^ku\|_{L^\infty(B_{\rho}(x_0))}+
\rho^{2+\alpha} [u]_{C^{2,\alpha}(B_{\rho}(x_0))}\right)
.\end{split}\end{equation}
We are thus in the setting of Proposition~\ref{SVA-LS-LEOSKDNDO}
with~$k:=2$, $\vartheta:=1$ and~$\sigma:=C\,\Big(
\|f\|_{C^\alpha(B_1)}+
\| u\|_{L^\infty(B_1)}\Big)$. If~$\e_0$ is
as in Proposition~\ref{SVA-LS-LEOSKDNDO},
we have that, if~$\rho$ is sufficiently small such that~$C\rho^\alpha\le \e_0$, then
the framework in~\eqref{Rstar8b-1} is fulfilled.
We thereby deduce from~\eqref{08uojbslj94tqquadqquad-1} that
\begin{equation}\label{xse2d32IJSMS-1iur93ufjhr23I-POKS}
\|u\|_{C^{2,\alpha}(B_{1/4})}
\le C\,\Big(
\|f\|_{C^\alpha(B_1)}+
\| u\|_{L^\infty(B_1)}\Big).\end{equation}
Accordingly, we can combine~\eqref{xse2d32IJSMS-1iur93ufjhr23I-POKS}
and a covering argument (see Lemma~\ref{COVERINGARG})
and obtain~\eqref{92i45kk5emmaabeCOVERINGARG}, as desired.
\end{proof}

Now we focus on the local estimates at the boundary:

\begin{proof}[Proof of Theorem~\ref{THM:0-T536EcdhPOaIJodidpo3tju24ylkiliR2543uj8890j09092po-okf-3edk-3edhf-2-fh564-LEATB}] As in~\eqref{omMAEwgNGLt4IONkjgfnbdsdtyuijk9024i5SZERSD},
it suffices to prove Theorem~\ref{THM:0-T536EcdhPOaIJodidpo3tju24ylkiliR2543uj8890j09092po-okf-3edk-3edhf-2-fh564-LEATB} when~$g$ vanishes identically. Then, one can
repeat verbatim the proof of Theorem~\ref{THM:0-T536EcdhPOaIJodidpo3tju24ylkiliR2543uj8890j09092po-okf-3edk-3edhf-2-fh564}
with the following modifications:
\begin{itemize}
\item balls~$B_r$ are replaced by halfballs~$B_r^+$,
\item the use of Theorem~\ref{SCHAUDER-INTE} is replaced by that of Theorem~\ref{SCHAUDER-BOUTH-VARIA},
\item the use of Proposition~\ref{SVA-LS-LEOSKDNDO} is replaced by that of Proposition~\ref{SVA-LS-LEOSKDNDO-HB-2}.\qedhere
\end{itemize}
\end{proof}

\section{Global Schauder estimates and existence theory}

We now employ the previously developed regularity results to obtain a global, up to the boundary, estimate
(in the spirit of the pioneer work by Oliver D. Kellogg on harmonic functions, see~\cite{MR1501602}),
which will be of pivotal use in what follows:

\begin{theorem}\label{ThneGiunfBB032t4jNKS-s3i4}
Let~$\Omega\subset\R^n$ be a bounded open set with boundary of class~$C^{2,\alpha}$
for some~$\alpha\in(0,1)$.
Let~$a_{ij}$, $b_i$, $c$, $f\in C^\alpha(\Omega)$. Let also~$g\in C^{2,\alpha}(\partial\Omega)$.
Assume that the ellipticity 
condition in~\eqref{ELLIPTIC} holds true.

Let~$u\in C^2(\Omega)\cap C(\overline\Omega)$ be a solution of
\begin{equation*}
\begin{dcases}
\sum_{i,j=1}^n a_{ij} \partial_{ij}u+\sum_{i=1}^n b_i\partial_i u+cu=f
&\quad{\mbox{in }}\,\Omega,\\
u=g&\quad{\mbox{on }}\,\partial\Omega.\end{dcases}
\end{equation*}
Then, there exists~$C>0$, depending only on~$n$,
$\alpha$, $\Omega$, $a_{ij}$, $b_i$ and~$c$, such that\footnote{When~$\Omega$ has boundary of class~$C^{2,\alpha}$
and~$g:\partial\Omega\to\R$ the  notion of~$\|g\|_{C^{2,\alpha}(\partial\Omega)}$ can be equivalently
defined in various ways. Rather than working with norms for functions defined on curved boundaries,
it is usually more handy to consider a boundary function as the restriction of a globally defined
function with its appropriate norm. Namely, one can extend~$g$ to a function~$\widetilde g:\overline{\Omega}\to\R$
with~$\widetilde g=g$ along~$\partial\Omega$ and~$\|\widetilde g\|_{C^{2,\alpha}(\Omega)}<+\infty$, see e.g.~\cite[Lemma~6.38]{MR1814364}. Then, one can define~$\|g\|_{C^{2,\alpha}(\partial\Omega)}$
as the infimum of~$\|\widetilde g\|_{C^{2,\alpha}(\Omega)}$ over all possible extensions~$\widetilde g$
as above. See~\cite[Section~6.2]{MR1814364} for further details about this procedure.\label{ESTE678i34rihf}}
\begin{equation}\label{LAMds-STudf-APEImriJEr1-1} \|u\|_{C^{2,\alpha}(\Omega)}\le C\,\Big(\|g\|_{C^{2,\alpha}(\partial\Omega)}+\|f\|_{C^\alpha(\Omega)}+\|u\|_{L^\infty(\Omega)}\Big).\end{equation}
Additionally, if~$c(x)\le0$ for all~$x\in\Omega$ then
\begin{equation} \label{LAMds-STudf-APEImriJEr1-2}
\|u\|_{C^{2,\alpha}(\Omega)}\le C\,\Big(\|g\|_{C^{2,\alpha}(\partial\Omega)}+\|f\|_{C^\alpha(\Omega)}\Big).\end{equation}
\end{theorem}

We stress that, in spite of their similarity, the estimates in~\eqref{LAMds-STudf-APEImriJEr1-1}
and~\eqref{LAMds-STudf-APEImriJEr1-2} present a fundamental difference since the latter does
not involve the solution~$u$ in its right hand side: this ``tiny detail'' will provide an essential
ingredient in the proof of the forthcoming existence result in Theorem~\ref{Theorem6.14GT}.

The proof of~\eqref{LAMds-STudf-APEImriJEr1-2} will rely on a general form of the Weak Maximum Principle
which has also independent interest (compare with Corollary~\ref{WEAKMAXPLE} and see~\cite[Chapters~3, 8 and~9]{MR1814364}
for more general and sharper forms of this result):

\begin{lemma}\label{iojsc6442-46-433-453-542-12345}
Let~$\Omega\subseteq\R^n$ be open and bounded.
Let~$a_{ij}$, $b_i\in L^\infty(\Omega)$.
Assume that the ellipticity 
condition in~\eqref{ELLIPTIC} holds true.
Let~$u\in C^2(\Omega)\cap C(\overline\Omega)$.
Then, we have:
\begin{itemize}
\item[(i).] If
\begin{equation}\label{3e5qu245456a2t24i42o2n} \sum_{i,j=1}^n a_{ij} \partial_{ij}u+\sum_{i=1}^n b_i\partial_i u\ge0
\quad{\mbox{ in }}\,\Omega\end{equation}
then
$$\sup_{\overline\Omega} u=\sup_{\partial\Omega} u.$$
\item[(ii).] If
$$ \sum_{i,j=1}^n a_{ij} \partial_{ij}u+\sum_{i=1}^n b_i\partial_i u\le0
\quad{\mbox{ in }}\,\Omega$$
then
$$\inf_{\overline\Omega} u=\inf_{\partial\Omega} u.$$
\end{itemize}
Here above, $C$ is a positive constant depending only on~$n$, $\lambda$, $\|b_1\|_{L^\infty(\Omega)},\dots,\|b_n\|_{L^\infty(\Omega)}$
and~$\Omega$.
\end{lemma}

\begin{proof} We start by proving~(i). To this end, 
we observe that, as a byproduct of the ellipticity condition in~\eqref{ELLIPTIC}
exploited with~$\xi:=e_1$, we have that
\begin{equation}\label{ci0pwrOOldOrfhndin904pol3a5rMS-2ri32phord}
a_{11}\ge\lambda.\end{equation}
Now, let
$$\gamma:=\frac{1+\sup_{i\in\{1,\dots,n\}}\|b_i\|_{L^\infty(\Omega)} }{\lambda}$$ and~$\eta(x):=e^{\gamma x_1}$.
Let also~$\e>0$ and~$u_\e(x):=u(x)+\e \eta(x)$. We have that~$\partial_i\eta=\delta_{i1}\gamma\eta$
and~$\partial_{ij}\eta=\delta_{i1}\delta_{j1}\gamma^2\eta$.  Consequently by~\eqref{3e5qu245456a2t24i42o2n}
and~\eqref{ci0pwrOOldOrfhndin904pol3a5rMS-2ri32phord},
\begin{equation}\label{kd-NSPSchkerRAkdDenftyOmega}
\begin{split}&
\sum_{i,j=1}^n a_{ij} \partial_{ij}u_\e+\sum_{i=1}^n b_i\partial_i u_\e\\ \ge \;&
\e\left[
\sum_{i,j=1}^n a_{ij} \partial_{ij}\eta+\sum_{i=1}^n b_i\partial_i \eta
\right]\\
=\;&\e\left[\gamma^2 a_{11}+b_1\gamma
\right]\eta\\
\ge\;& \e\left[\gamma^2 \lambda-\sup_{i\in\{1,\dots,n\}}\|b_i\|_{L^\infty(\Omega)}\gamma
\right]\eta\\ =\;&\e\gamma\eta.
\end{split}\end{equation}

We claim that
\begin{equation}\label{an4anis3dentically}
\sup_{\overline\Omega} u_\e=\sup_{\partial\Omega} u_\e.\end{equation}
Indeed, suppose not, then~$u_\e$ presents an interior local maximum
at some point~$\overline{x}_\e$. Therefore, the Hessian matrix of~$u_\e$ at~$\overline{x}_\e$
is nonpositive and the gradient of~$u_\e$ at~$\overline{x}_\e$ vanishes.
Thus, we denote by~$M=\{M_{ij}\}_{i,j\in\{1,\dots,n\}}$ the square root (in the matrix sense, see e.g.~\cite[page~125]{MR2337395})
of minus the Hessian matrix of~$u_\e$ at~$\overline{x}_\e$. In this way, we have that
$$ -\partial_{ij}u_\e(\overline{x}_\e)=\sum_{k=1}^n M_{ik}M_{jk}.$$
Now we use the ellipticity 
condition in~\eqref{ELLIPTIC} to see that
$$ \sum_{i,j=1}^na_{ij}(\overline{x}_\e) M_{ik}M_{jk}\ge\lambda\sum_{m=1}^n |M_{m k}|^2.$$
As a result,
$$ -\sum_{i,j=1}^na_{ij}(\overline{x}_\e) \partial_{ij}u_\e(\overline{x}_\e)=
\sum_{i,j,k=1}^na_{ij}(\overline{x}_\e) M_{ik}M_{jk}\ge\lambda \sum_{m,k=1}^n |M_{m k}|^2\ge0.$$
Recalling~\eqref{kd-NSPSchkerRAkdDenftyOmega}, we thereby obtain a contradiction.
This completes the proof of~\eqref{an4anis3dentically}:
thus, by sending~$\e\searrow0$ in~\eqref{an4anis3dentically},
we establish the desired claim in~(i).

To prove~(ii), one can employ~(i) on the function~$-u$.
\end{proof}

\begin{corollary}\label{9oikr3eg9u2iwhefk9qryfhweiwwrfgfLftyOmega-kdf}
Let~$\Omega\subseteq\R^n$ be open and bounded.
Let~$a_{ij}$, $b_i$, $c$, $f\in L^\infty(\Omega)$.
Assume that~$c(x)\le0$ for all~$x\in\Omega$ and that the ellipticity 
condition in~\eqref{ELLIPTIC} holds true.
Let~$u\in C^2(\Omega)\cap C(\overline\Omega)$ be a solution of
\begin{equation*}
\sum_{i,j=1}^n a_{ij} \partial_{ij}u+\sum_{i=1}^n b_i\partial_i u+cu\ge f
\quad{\mbox{ in }}\,\Omega.\end{equation*}
Then,
\begin{equation}\label{9u2iwhefk9qryfhweiwwrfgfLftyOmega}
\sup_{\overline\Omega} u\le\sup_{\partial\Omega} u^++C\|f\|_{L^\infty(\Omega)}.\end{equation}
Here above, $C$ is a positive constant depending only on~$n$, $\lambda$, $\|b_1\|_{L^\infty(\Omega)},\dots,\|b_n\|_{L^\infty(\Omega)}$
and~$\Omega$.
\end{corollary}

\begin{proof} We observe that
\begin{equation}\label{9u2iwhefk9qryfhweiwwrfgfLftyOmega2}
{\mbox{it suffices to prove the desired result when~$f$ vanishes identically.}}
\end{equation}
Indeed, suppose to know the desired result when~$f$ vanishes identically. Let~$R>0$ be sufficiently large such that~$\Omega\subseteq B_{R}$
and let
$$v(x):=u(x)+ \|f\|_{L^\infty(\Omega)}\eta(x),$$ with~$\eta(x):=K(e^{\gamma x_1}-e^{\gamma R})$,
$K:=\frac{e^{\gamma R}}\gamma$ and~$\gamma :=\frac{1+\|b_1\|_{L^\infty(\Omega)}}\lambda$.
We have that~$\partial_i\eta(x)=\delta_{i1} K\gamma e^{\gamma x_1}$ and~$\partial_{ij}\eta(x)=\delta_{i1}\delta_{j1}K \gamma^2 e^{\gamma x_1}$. As a consequence, noticing that~$e^{\gamma x_1}-e^{\gamma R}\le0$ in~$\Omega$
and recalling~\eqref{ci0pwrOOldOrfhndin904pol3a5rMS-2ri32phord},
we find that, in~$\Omega$,
\begin{eqnarray*}
&&\sum_{i,j=1}^n a_{ij} \partial_{ij}v+\sum_{i=1}^n b_i\partial_i v+cv\\&\ge&
f+\|f\|_{L^\infty(\Omega)}\left[
K\gamma^2 a_{11} e^{\gamma x_1}+K\gamma b_1 e^{\gamma x_1} +cK(e^{\gamma x_1}-e^{\gamma R})\right]\\
\\&\ge&-\|f\|_{L^\infty(\Omega)}+\|f\|_{L^\infty(\Omega)}\left[
\lambda K\gamma^2 e^{\gamma x_1}-K\gamma\| b_1\|_{L^\infty(\Omega)} e^{\gamma x_1} +cK(e^{\gamma x_1}-e^{\gamma R})\right]\\
&\ge&\|f\|_{L^\infty(\Omega)}\left[K \gamma e^{\gamma x_1}(\lambda \gamma-\| b_1\|_{L^\infty(\Omega)}) -1\right]
\\&=&\|f\|_{L^\infty(\Omega)}\left(K \gamma e^{\gamma x_1}-1\right) \\
&\ge&\|f\|_{L^\infty(\Omega)}\left(K \gamma e^{-\gamma R}-1\right)
\\&\geq&0.
\end{eqnarray*}
Thus, if the desired result holds true when~$f$ vanishes identically, we deduce that
$$\sup_{\overline\Omega} v\le\sup_{\partial\Omega} v^+.$$
Hence, since~$v^+\le u^+$ and
$$ v\ge u-\|f\|_{L^\infty(\Omega)} Ke^{\gamma R}=u-\frac{ \|f\|_{L^\infty(\Omega)}\ e^{2\gamma R}}\gamma,$$
we obtain that~\eqref{9u2iwhefk9qryfhweiwwrfgfLftyOmega} holds true. This completes the proof of~\eqref{9u2iwhefk9qryfhweiwwrfgfLftyOmega2}.

Therefore, in the light of~\eqref{9u2iwhefk9qryfhweiwwrfgfLftyOmega2},
we now suppose that~$f$ vanishes identically.
We let
\begin{eqnarray*}
&&\Omega_\star:=\{x\in\Omega {\mbox{ s.t. }}u(x)>0\}
\end{eqnarray*}
and we remark that, if~$x\in\Omega_\star$,
$$ \sum_{i,j=1}^n a_{ij} \partial_{ij}u+\sum_{i=1}^n b_i\partial_i u\ge-cu\ge 0.$$
{F}rom this and Lemma~\ref{iojsc6442-46-433-453-542-12345}(i), applied here on~$\Omega_\star$ in the place of~$\Omega$,
we deduce that
\begin{equation}\label{9u2iwhefk9qryfhweiwwrfgfLftyOmega-3}\sup_{\overline{\Omega_\star}} u=\sup_{\partial\Omega_\star} u.\end{equation}

Now we prove that~\eqref{9u2iwhefk9qryfhweiwwrfgfLftyOmega} holds true with~$f=0$.
Indeed, if not, there exists~$X\in\overline{\Omega}$ and~$a>0$ such that~$u(X)\ge a+\sup_{\partial\Omega} u^+$.
In particular, $u(X)>0$ and thus~$X\in\Omega_\star$.
This and~\eqref{9u2iwhefk9qryfhweiwwrfgfLftyOmega-3} yield that
\begin{equation}\label{9u2iwhefk9qryfhweiwwrfgfLftyOmegaSUPP-0}
a+\sup_{\partial\Omega} u^+\le u(X)\le\sup_{\overline{\Omega_\star}} u=\sup_{\partial\Omega_\star} u.\end{equation}
This also entails that
\begin{equation}\label{9u2iwhefk9qryfhweiwwrfgfLftyOmegaSUPP}
\sup_{\partial\Omega_\star} u>0.\end{equation}
Now, since~$\partial\Omega_\star$ is compact and~$u$ is continuous along~$\partial\Omega_\star\subseteq\overline\Omega$,
the above supremum is attained at some point~$Y\in\partial\Omega_\star$.
We deduce from~\eqref{9u2iwhefk9qryfhweiwwrfgfLftyOmegaSUPP} that~$u(Y)>0$,
hence~$u(Y)=u^+(Y)$.

We claim that
\begin{equation}\label{9u2iwhefk9qryfhweiwwrfgfLftyOmegaSUPP2}
Y\in\partial\Omega.\end{equation}
Indeed, if not, $Y$ would be an interior point of~$\Omega$
and therefore, by continuity, there would exist~$\varrho>0$ such that~$B_\varrho(Y)\subset\Omega$
and~$u(x)\ge\frac{u(Y)}2>0$ for all~$x\in B_\varrho(Y)$. But this would give that~$Y$ is also an interior point
for~$\Omega_\star$, in contradiction with the assumption that~$Y\in\partial\Omega_\star$.
This establishes~\eqref{9u2iwhefk9qryfhweiwwrfgfLftyOmegaSUPP2}.

Then, from~\eqref{9u2iwhefk9qryfhweiwwrfgfLftyOmegaSUPP-0} and~\eqref{9u2iwhefk9qryfhweiwwrfgfLftyOmegaSUPP2},
$$ a+\sup_{\partial\Omega} u^+\le \sup_{\partial\Omega_\star} u=u(Y)=u^+(Y)\le
\sup_{\partial\Omega} u^+.$$
This is a contradiction, which shows that~\eqref{9u2iwhefk9qryfhweiwwrfgfLftyOmega} holds true with~$f=0$,
as desired.
\end{proof}

We can now focus on the proof of the gloabal estimates up to the boundary:

\begin{figure}
  \centering
  \includegraphics[width=.5\linewidth]{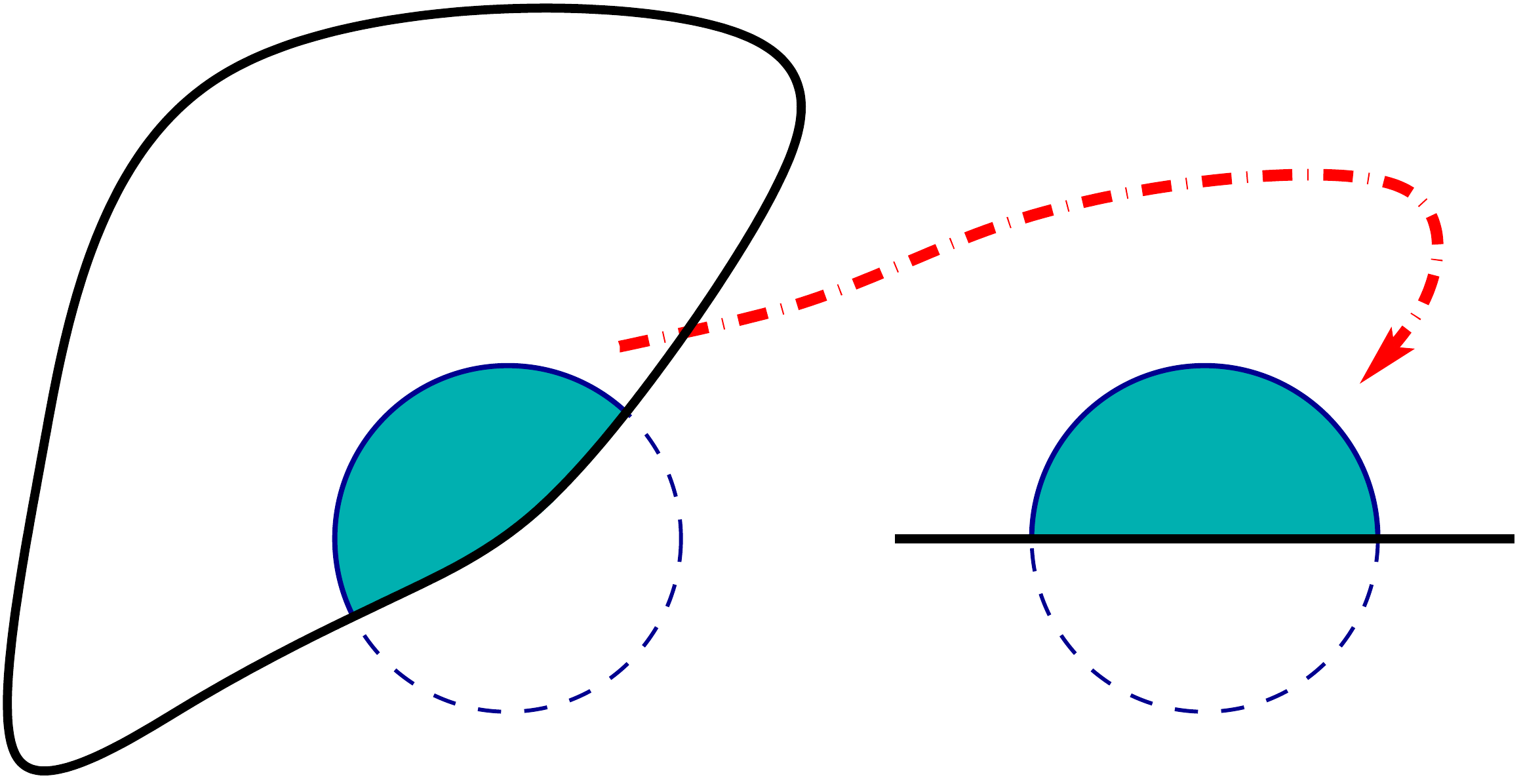}
 \caption{\sl A domain deformation\index{domain deformation} to locally straighten the boundary.}\label{fLOCAPAere8MPyItangeFqI}
\end{figure}

\begin{proof}[Proof of Theorem~\ref{ThneGiunfBB032t4jNKS-s3i4}] The core of the proof
consists in straightening the boundary (see Figure~\ref{fLOCAPAere8MPyItangeFqI})
in order to use the local estimates at the boundary given in Theorem~\ref{THM:0-T536EcdhPOaIJodidpo3tju24ylkiliR2543uj8890j09092po-okf-3edk-3edhf-2-fh564-LEATB} (the remaining part
of the domain lies away from the boundary and one can use there the
interior estimates in Theorem~\ref{THM:0-T536EcdhPOaIJodidpo3tju24ylkiliR2543uj8890j09092po-okf-3edk-3edhf-2-fh564}).
The structural reason for which this argument does the business is that domain deformations
preserve the class of operators that we are considering, according to Lemma~\ref{92o3jwelg02jfj-qiewjdfnoewfh9ewohgveioqNONDIFO234RM} (as a matter of fact, preserving a suitable class of equations
under diffeomorphisms was precisely one of the main motivations
for us to study elliptic equations in such a generality, as discussed on page~\pageref{TOWLemma92o3jwelg02jfj-qiewjdfnoewfh9ewohgveioqNONDIFO234RM}).

The technical details of the proof go as follows.
Since~$\Omega$ is bounded and with boundary of class~$C^{2,\alpha}$, we find~$\rho>0$ and a finite family
of balls~$\{B_\rho(p_i)\}_{i\in\{1,\dots,N\}}$ with~$p_i\in\partial\Omega$ and
\begin{equation}\label{beuationdetildealphauation}
 {\mathcal{N}}:=\{ x\in\Omega {\mbox{ s.t. }} B_{\rho/8}(x)\cap(\partial\Omega)\ne\varnothing\}
\subseteq \bigcup_{i=1}^N B_{\rho/4}(p_i)\end{equation}
such that~$B_\rho(p_i)\cap \Omega$
is equivalent to~$B_\rho^+$ via a diffeomorphism~$T_i$ of class~$C^{2,\alpha}$ with
inverse of class~$C^{2,\alpha}$.

We can thus employ Lemma~\ref{92o3jwelg02jfj-qiewjdfnoewfh9ewohgveioqNONDIFO234RM} in this setting:
that is, setting~$\widetilde u_i(y):=u(T_i^{-1}(y))$ we find from~\eqref{oj-LEM-EQUJSN-comSIJMCAi-mbiJAQ-ledm-2}
that~$\widetilde u_i$ satisfies in~$B_\rho^+$ the same type of equation that~$u$ satisfies in~$B_\rho(p_i)\cap \Omega$,
with ellipticity constants preserved up to the multiplication by structural constants,
and with sign of the zero order coefficient maintained nonpositive.

{F}rom these facts, we are in the position of using the local estimates at the boundary in Theorem~\ref{THM:0-T536EcdhPOaIJodidpo3tju24ylkiliR2543uj8890j09092po-okf-3edk-3edhf-2-fh564-LEATB} on the function~$\widetilde u_i$, arriving at
\begin{equation*} \|\widetilde u_i\|_{C^{2,\alpha}(B_{\rho/2}^+)}\le C\,\Big(\|u\|_{L^\infty(\Omega)}+
\|g\|_{C^{2,\alpha}(\partial\Omega)}+\|f\|_{C^\alpha(\Omega)}\Big).\end{equation*}
As a consequence,
\begin{equation*} \|u\|_{C^{2,\alpha}(\Omega\cap B_{\rho/2}(p_i))}\le C\,\Big(\|u\|_{L^\infty(\Omega)}+
\|g\|_{C^{2,\alpha}(\partial\Omega)}+\|f\|_{C^\alpha(\Omega)}\Big)\end{equation*}
and therefore, by~\eqref{beuationdetildealphauation},
\begin{equation}\label{LAMds-STudf-APEImriJErqr2vbgjqxb4rqtqrw3ww1-1}
\|u\|_{C^{2,\alpha}({\mathcal{N}})}\le C\,\Big(\|u\|_{L^\infty(\Omega)}+
\|g\|_{C^{2,\alpha}(\partial\Omega)}+\|f\|_{C^\alpha(\Omega)}\Big),\end{equation}
up to renaming constants from line to line.

Now we consider a covering of the remaining part of the domain. Namely, we take a finite family of balls such that
\begin{equation}\label{THM:0-T536EcdhPOaIJodidpo3tju24ylkiliR2543uj8890j09092po-okf-3edk-3edhf-2-fh564-09-er}
 {\mathcal{J}}:=\{ x\in\Omega {\mbox{ s.t. }} B_{\rho/9}(x)\cap(\partial\Omega)=\varnothing\}
\subseteq \bigcup_{i=1}^M B_{\rho/4}(q_i),\end{equation}
for suitable~$q_i\in\Omega$. We utilize the
interior estimates of Theorem~\ref{THM:0-T536EcdhPOaIJodidpo3tju24ylkiliR2543uj8890j09092po-okf-3edk-3edhf-2-fh564}
in each ball~$B_{7\rho/20}(q_i)\Supset B_{\rho/4}(q_i)$ and we conclude that
\begin{equation*} \|u\|_{C^{2,\alpha}(B_{\rho/4}(q_i))}\le C\,\Big(\|u\|_{L^\infty(\Omega)}+\|f\|_{C^\alpha(\Omega)}\Big).\end{equation*}
This and~\eqref{THM:0-T536EcdhPOaIJodidpo3tju24ylkiliR2543uj8890j09092po-okf-3edk-3edhf-2-fh564-09-er}
yield that
\begin{equation*} \|u\|_{C^{2,\alpha}({\mathcal{J}})}\le C\,\Big(\|u\|_{L^\infty(\Omega)}+\|f\|_{C^\alpha(\Omega)}\Big).\end{equation*}
This, together with~\eqref{LAMds-STudf-APEImriJErqr2vbgjqxb4rqtqrw3ww1-1} and a covering argument as in
Lemma~\ref{COVERINGARG},
proves~\eqref{LAMds-STudf-APEImriJEr1-1}.

Then, the claim in~\eqref{LAMds-STudf-APEImriJEr1-2} follows from~\eqref{LAMds-STudf-APEImriJEr1-1}
and Corollary~\ref{9oikr3eg9u2iwhefk9qryfhweiwwrfgfLftyOmega-kdf}.
\end{proof}

\begin{figure}
  \centering
  \includegraphics[width=.6\linewidth]{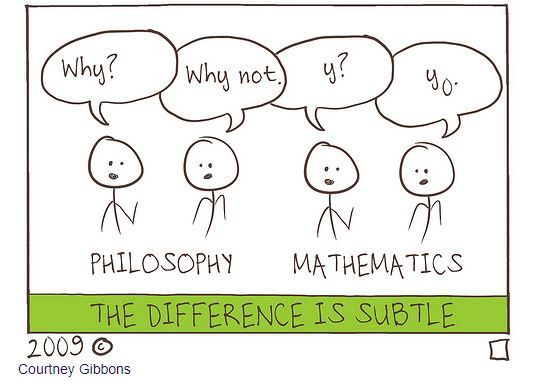}
 \caption{\sl Courtesy of Courtney Gibbons.}\label{fere8MPyItangeFI}
 \end{figure}

We point out that the regularity theory developed so far in this section
(as well as in Sections~\ref{KPJMDarJMSiaCLiMMSciAmoenbvD8Sijf-21}, \ref{KPJMDarJMSiaCLiMMSciAmoenbvD8Sijf-22} and~\ref{KPJMDarJMSiaCLiMMSciAmoenbvD8Sijf-23}) dealt with ``a priori estimates'':
namely, we estimated\footnote{The Latin wording ``a priori'' means ``from the one before'' and usually \index{a priori estimate} \label{FOOAMPRIUoRAnABBA-GRA}
refers to knowledge that is established or developed
without being based on previous experience.
The use of this wording in the context of partial differential equations has become
a standard in the technical jargon after the work of Sergei Natanovich (a.k.a. Serge) Bernstein~\cite{MR1511375}.

As a matter of fact, this business of a priori estimates does pose some philosophical puzzles,
since not only we have developed a whole regularity theory so far for
objects which {\em could not} exist, but actually at a level of generality for which
these objects {\em do not} exist at all, as it will be apparent in the forthcoming
examples~\eqref{MpksmREFVSpaoMSRHsfo0s} and~\eqref{MpksmREFVSpaoMSRHsfo0s-2}
for which, for instance, the application of Theorem~\ref{ThneGiunfBB032t4jNKS-s3i4}
becomes void.

Fortunately, mathematics is no philosophy (see Figure~\ref{fere8MPyItangeFI} for a 
systemic difference)
and the hard work done to prove
Theorem~\ref{ThneGiunfBB032t4jNKS-s3i4} will pay off: specifically, we will obtain in
Theorem~\ref{Theorem6.14GT} that solutions always exist at least under a sign assumption
on the coefficient~$c$. The ``magic'' of Theorem~\ref{ThneGiunfBB032t4jNKS-s3i4}
is that it does prove uniform estimates in a case in which we know that the solution
exists (namely for the Laplace operator) and it thus allows one to use the Implicit Function Theorem (in some form)
to find solutions for operators ``nearby the Laplacian''. But then one can repeat the argument,
precisely because the estimates in
Theorem~\ref{ThneGiunfBB032t4jNKS-s3i4} are uniform, and keep finding solutions for a large
class of operators which are ``connected'' to the Laplacian.

This technique is sometimes called the ``continuity method''. \index{continuity method}
Making this heuristic idea work will be precisely the content of the proof of Theorem~\ref{Theorem6.14GT}.} suitable norms (of the derivatives) of a solution
of a partial differential equation, without knowing that the solution exists at all!

As a matter of fact, there are concrete situations in which the Dirichlet problem
\begin{equation*}
\begin{dcases}
\sum_{i,j=1}^n a_{ij} \partial_{ij}u+\sum_{i=1}^n b_i\partial_i u+cu=f
&\quad{\mbox{in }}\,\Omega,\\
u=g&\quad{\mbox{on }}\,\partial\Omega\end{dcases}
\end{equation*}
has no solution~$u\in C^2(\Omega)\cap C(\overline\Omega)$.
This possibly unpleasant situation is not an artifact of the generality considered for the equation,
and actually explicit examples of nonexistence arise even in low dimensions, even for very smooth domains,
even when the operator~$\sum_{i,j=1}^n a_{ij} \partial_{ij}$ is just the Laplacian.
For example,
the Dirichlet problem \begin{equation}\label{MpksmREFVSpaoMSRHsfo0s}\begin{dcases} u''(x)+u(x)=1 &{\mbox{ for all }}x\in(\pi,2\pi),\\u(\pi)=0,\\ u(2\pi)=0\end{dcases}\end{equation} possesses no solutions. Indeed, the general solution of the above ordinary differential equation has the form~$u(x)=c_1\sin x+c_2\cos x+1$ for constants~$c_1$ and~$c_2$. Imposing~$u(\pi)=0$ we find that~$c_2=1$.
But then~$u(2\pi)=\cos(2\pi)+1=2$, showing that~\eqref{MpksmREFVSpaoMSRHsfo0s} does not admit any solution.\medskip

As for higher dimensional cases, we observe that for~$n\ge2$
the Dirichlet problem \begin{equation}\label{MpksmREFVSpaoMSRHsfo0s-2}\begin{dcases} \Delta u(x)-\frac{n-1}{|x|^2}\,
\nabla u\cdot x+u(x)=1 &{\mbox{ for all }}x\in B_{2\pi}\setminus\overline{B_\pi},\\u(x)=0
&{\mbox{ for all }}x\in(\partial B_{2\pi})\cup(\partial{B_\pi})
\end{dcases}\end{equation} possesses no solutions. To check this,
we suppose by contradiction that a solution~$u$ exists and for all~$x\in B_{2\pi}\setminus\overline{B_\pi}$ we define
$$ v(x):=\fint_{\partial B_1} u(|x|\vartheta)\,d{\mathcal{H}}^{n-1}_\vartheta.$$
We note that, for every~$i\in\{1,\dots,n\}$,
$$ \partial_i v(x)=\sum_{k=1}^n\fint_{\partial B_1} \partial_k u(|x|\vartheta)\,\frac{\vartheta_k x_i}{|x|}\,d{\mathcal{H}}^{n-1}_\vartheta.$$
Also, by~\eqref{Byphericalpresentatitlaciorem},
$$ \Delta v(x)=\fint_{\partial B_1} \Delta u(|x|\vartheta)\,d{\mathcal{H}}^{n-1}_\vartheta.$$
Consequently, for every~$x\in B_{2\pi}\setminus\overline{B_\pi}$,
\begin{eqnarray*}
&&\Delta v(x)-\frac{n-1}{|x|^2}\,
\nabla v(x)\cdot x+v(x)
= \fint_{\partial B_1} \left(\Delta u(|x|\vartheta)
-\frac{n-1}{|x|^2}\,\nabla u(|x|\vartheta)\cdot(|x|\vartheta)+u(|x|\vartheta)\right)\,d{\mathcal{H}}^{n-1}_\vartheta=1,
\end{eqnarray*}
meaning that~$v$ is also a solution of the equation in~\eqref{MpksmREFVSpaoMSRHsfo0s-2}.

Additionally, $v$ is invariant under rotation and therefore, writing~$v(x)=v_0(|x|)$ with~$v_0:\R\to\R$ and recalling
the spherical representation of the Laplacian in Theorem~\ref{SPHECOO}, for all~$r:=|x|\in(\pi,2\pi)$,
$$ \Delta v(x)-\frac{n-1}{|x|^2}\,
\nabla v(x)\cdot x+v(x)=\left(v_0''(r)+
\frac{n-1}{r}\,v_0'(r)\right)-\frac{n-1}{r^2}\,
v_0'(r)\,r+v_0(r)=v_0''(r)+v_0(r)
$$
and therefore~$v_0''(r)+v_0(r)=1$. Since in addition
$$v_0(\pi)=v(\pi e_1)=\fint_{\partial B_1} u(\pi\vartheta)\,d{\mathcal{H}}^{n-1}_\vartheta=0$$
and similarly~$v_0(2\pi)=0$, we find that~$v_0$ is a solution of~\eqref{MpksmREFVSpaoMSRHsfo0s}.
But we already know that~\eqref{MpksmREFVSpaoMSRHsfo0s} does not possess any solution,
hence this is a contradiction that shows that~\eqref{MpksmREFVSpaoMSRHsfo0s-2}
does not admit any solution.\medskip

Rather than being discouraged by the examples in~\eqref{MpksmREFVSpaoMSRHsfo0s}
and~\eqref{MpksmREFVSpaoMSRHsfo0s-2}, we hold onto what we already positively know
(such as the existence result obtained via the Perron method in Corollary~\ref{S-coroEXIS-M023})
and investigate the structure of the previous examples.
In particular, in both these examples, the coefficient in front of the zero order term of the equation is
positive. This turns out to be the major obstruction towards the existence of solutions,
and in fact when the coefficient in front of the zero order term is nonnegative the corresponding
Dirichlet problem does admit a (unique) solution. The precise result goes as follows:

\begin{theorem}\label{Theorem6.14GT} Let~$\Omega\subset\R^n$ be a bounded open set with boundary of class~$C^{2,\alpha}$
for some~$\alpha\in(0,1)$.
Let~$a_{ij}$, $b_i$, $c$, $f\in C^\alpha(\Omega)$. Let also~$g\in C^{2,\alpha}(\partial\Omega)$.
Assume that~$c(x)\le0$ for all~$x\in\Omega$ and that the ellipticity 
condition in~\eqref{ELLIPTIC} holds true.

Then, the Dirichlet problem
\begin{equation}\label{134-24rf-pkwPK-TAGBco6YHSRVrfik-RiodSS02o3ef1}
\begin{dcases}
\sum_{i,j=1}^n a_{ij} \partial_{ij}u+\sum_{i=1}^n b_i\partial_i u+cu=f
&\quad{\mbox{in }}\,\Omega,\\
u=g&\quad{\mbox{on }}\,\partial\Omega\end{dcases}
\end{equation}
has a unique solution in~$C^2(\Omega)\cap C(\overline\Omega)$.

Additionally, such a solution~$u$ belongs to~$C^{2,\alpha}(\Omega)$
and there exists~$C>0$, depending only on~$n$,
$\alpha$, $\Omega$, $a_{ij}$, $b_i$ and~$c$, such that
\begin{equation}\label{134-24rf-pkwPK-TAGBco6YHSRVrfik-RiodSS02o3ef2} \|u\|_{C^{2,\alpha}(\Omega)}\le C\,\Big(\|g\|_{C^{2,\alpha}(\partial\Omega)}+\|f\|_{C^\alpha(\Omega)}\Big).\end{equation}
\end{theorem}

The proof of Theorem~\ref{Theorem6.14GT} will use the a priori estimate
in Theorem~\ref{ThneGiunfBB032t4jNKS-s3i4} and the following general inversion result
for linear operators stating that small perturbations of invertible operators in Banach spaces remain invertible:

\begin{lemma}\label{LEMroofofheoreheorem6.14GT}
Let~${\mathcal{X}}$ and~${\mathcal{Y}}$ be Banach spaces and let~${\mathcal{R}}$ and~${\mathcal{S}}$
be bounded linear operators from~${\mathcal{X}}$ to~${\mathcal{Y}}$.

Assume that~${\mathcal{R}}$ is invertible with bounded inverse.
Assume also that
\begin{equation*}
\|{\mathcal{S}}\|<\frac1{\|{\mathcal{R}}^{-1}\|}.
\end{equation*}
Then, the operator~${\mathcal{R}}+{\mathcal{S}}$ is invertible and
$$ \big\|({\mathcal{R}}+{\mathcal{S}})^{-1}\big\|\le \frac{\|{\mathcal{R}}^{-1}\|}{1-
\|{\mathcal{R}}^{-1} \|\,\|{\mathcal{S}}\|}. $$
\end{lemma}

\begin{proof} We denote by~${\mathcal{I}}$ the identity operator from~${\mathcal{X}}$ to~${\mathcal{X}}$.
Let~${\mathcal{T}}:=-{\mathcal{R}}^{-1} {\mathcal{S}}$ and note that
\begin{equation}\label{summathcalTk}
\|{\mathcal{T}}\|=\| {\mathcal{R}}^{-1} {\mathcal{S}}\|=\sup_{x\in{\mathcal{X}}\setminus\{0\}}\frac{| {\mathcal{R}}^{-1} {\mathcal{S}}x|}{|x|}\le
\sup_{x\in{\mathcal{X}}\setminus\{0\}}\frac{\|{\mathcal{R}}^{-1} \|\,|{\mathcal{S}}x|}{|x|}\le
\|{\mathcal{R}}^{-1} \|\,\|{\mathcal{S}}\|<1.
\end{equation}
We define
$$ {\mathcal{L}}_j:={\mathcal{I}}+{\mathcal{T}}+{\mathcal{T}}^2+\dots+{\mathcal{T}}^j=\sum_{k=0}^j{\mathcal{T}}^k$$
and we observe that for all~$i\ge1$
$$ \|{\mathcal{L}}_{j+i}-{\mathcal{L}}_j\|=\left\|\sum_{k=j+1}^{j+i}{\mathcal{T}}^k\right\|\le
\sum_{k=j+1}^{j+i}\|{\mathcal{T}}\|^k.$$
This and~\eqref{summathcalTk} yield that~${\mathcal{L}}_j$ is a Cauchy sequence of bounded linear operators from~${\mathcal{X}}$ to~${\mathcal{X}}$,
hence there exists a bounded linear operator~${\mathcal{L}}$ such that~$\|{\mathcal{L}}_j-{\mathcal{L}}\|\to0$
as~$j\to+\infty$. In this way, we see that, as~$j\to+\infty$,
\begin{eqnarray*}
&&\big\| {\mathcal{L}}({\mathcal{I}}-{\mathcal{T}})-{\mathcal{I}}\|\le
\big\|{\mathcal{L}}_j ({\mathcal{I}}-{\mathcal{T}})-{\mathcal{I}}\|
+\big\| ({\mathcal{L}}-{\mathcal{L}}_j)({\mathcal{I}}-{\mathcal{T}})\|\\&&\qquad
=\left\| \sum_{k=0}^j{\mathcal{T}}^k({\mathcal{I}}-{\mathcal{T}})-{\mathcal{I}}\right\|+o(1)
=\left\| \sum_{k=0}^j{\mathcal{T}}^k-\sum_{k=0}^j{\mathcal{T}}^{k+1}-{\mathcal{I}}\right\|+o(1)\\&&\qquad=\|{\mathcal{T}}^{j+1}\|+o(1)=o(1),
\end{eqnarray*}
thanks again to~\eqref{summathcalTk}, and accordingly
$$ {\mathcal{L}}({\mathcal{I}}-{\mathcal{T}})-{\mathcal{I}}=0.$$
This entails that
$$ {\mathcal{L}}{\mathcal{R}}^{-1} ({\mathcal{R}}+ {\mathcal{S}})=
{\mathcal{L}}({\mathcal{I}}+{\mathcal{R}}^{-1} {\mathcal{S}})=
{\mathcal{L}}({\mathcal{I}}-{\mathcal{T}})={\mathcal{I}}
$$
and as a consequence the operator~${\mathcal{R}}+ {\mathcal{S}}$ is invertible, with~$ ({\mathcal{R}}+ {\mathcal{S}})^{-1}=
{\mathcal{L}}{\mathcal{R}}^{-1}$.

With this,
\begin{eqnarray*} &&\|({\mathcal{R}}+ {\mathcal{S}})^{-1}\| =
\|{\mathcal{L}}{\mathcal{R}}^{-1}\|=\lim_{j\to+\infty}\|{\mathcal{L}}_j{\mathcal{R}}^{-1}\|
\le\lim_{j\to+\infty}\|{\mathcal{L}}_j\|\,\|{\mathcal{R}}^{-1}\|\\&&\qquad
\le \lim_{j\to+\infty}\sum_{k=0}^j\|{\mathcal{T}}\|^k\,\|{\mathcal{R}}^{-1}\|
=\frac{\|{\mathcal{R}}^{-1}\|}{1-\|{\mathcal{T}}\|}=\frac{\|{\mathcal{R}}^{-1}\|}{1-
\|{\mathcal{R}}^{-1} {\mathcal{S}}\|}\le\frac{\|{\mathcal{R}}^{-1}\|}{1-
\|{\mathcal{R}}^{-1} \|\,\|{\mathcal{S}}\|}
,\end{eqnarray*}
as desired.
\end{proof}

\begin{corollary}\label{i0ojendelope934etting876in979fu}
Let~${\mathcal{X}}$ and~${\mathcal{Y}}$ be Banach spaces and let~${\mathcal{L}}_0$ and~${\mathcal{L}}_1$
be bounded linear operators from~${\mathcal{X}}$ to~${\mathcal{Y}}$. For every~$t\in[0,1]$, let
$$ {\mathcal{L}}_t:=t{\mathcal{L}}_1+(1-t){\mathcal{L}}_0.$$
Assume that~${\mathcal{L}}_0$ is surjective and that there exists~$K>0$ such that
\begin{equation}\label{0248mathcalmathcal2} |x |\le K\,|{\mathcal{L}}_tx|\qquad{\mbox{for all $x\in{\mathcal{X}}$ and~$t\in[0,1]$.}}\end{equation}
Then, 
\begin{equation}\label{0248mathcalmathcal2-B2wefr} {\mbox{for every~$y\in{\mathcal{Y}}$ there exists a unique~$x_y\in{\mathcal{X}}$ such that~${\mathcal{L}}_1 x_y=y$.}}\end{equation}
In addition,
\begin{equation}\label{0248mathcalmathcal2-B2wefr2} |x_y|\le K\,|y|.\end{equation}
\end{corollary}

\begin{proof}
Up to replacing~$K$ in~\eqref{0248mathcalmathcal2} with~$K+\|{\mathcal{L}}_0\|+\|{\mathcal{L}}_1\|$,
we may assume that
\begin{equation}\label{shapeioryhfbfb15436586960dfkiPPPP}
\|{\mathcal{L}}_0\|+\|{\mathcal{L}}_1\|\le K.\end{equation}
Let~$N\in\N\cap\left[ 2K^2,+\infty\right)$. We claim that
\begin{equation}\label{0248mathcalmathcal}
{\mbox{for all~$j\in\{0,\dots,N-1\}$ and all~$t\in$}}\left[\frac{j}N,\frac{j+1}N\right],
{\mbox{the operator~${\mathcal{L}}_t:{\mathcal{X}}\to{\mathcal{Y}}$ is invertible.}}
\end{equation}
We prove this by induction over~$j$. 
If~$j=0$, we let~${\mathcal{R}}:={\mathcal{L}}_0$ and~${\mathcal{S}}:=t({\mathcal{L}}_1-{\mathcal{L}}_0)$.
We point out that~${\mathcal{R}}$ is invertible
(we know that it is surjective, and the injectivity follows from~\eqref{0248mathcalmathcal2} with~$t:=0$).
Also, the norm of~${\mathcal{R}}^{-1}$ is bounded from above by~$K$ (thanks again to~\eqref{0248mathcalmathcal2}
with~$t:=0$). Moreover, for all~$t\in\left[0,\frac{1}N\right]$, in light of~\eqref{shapeioryhfbfb15436586960dfkiPPPP} we have that
$$ \| {\mathcal{S}}\|=\sup_{x\in{\mathcal{X}}\setminus\{0\}}\frac{\big| t({\mathcal{L}}_1-{\mathcal{L}}_0)x\big|}{|x|}
\le\frac1N\,\sup_{x\in{\mathcal{X}}\setminus\{0\}}\frac{| {\mathcal{L}}_1x|+|{\mathcal{L}}_0x|}{|x|}\le
\frac{K}N\le\frac{K^2}{N\,\|{\mathcal{R}}^{-1}\|}\le\frac{1}{2\,\|{\mathcal{R}}^{-1}\|}.
$$
We are therefore in the position of using Lemma~\ref{LEMroofofheoreheorem6.14GT},
thus finding that~${\mathcal{R}}+{\mathcal{S}}={\mathcal{L}}_t$ is invertible.

This establishes~\eqref{0248mathcalmathcal} when~$j=0$ and we now focus on the inductive step.
To this end, we suppose that~\eqref{0248mathcalmathcal} holds true for some index~$j\in\{0,\dots,N-2\}$
and we aim at proving it for the index~$j+1$. For this, we let~${\mathcal{R}}:={\mathcal{L}}_{(j+1)/N}$
and~${\mathcal{S}}:=\left(t-\frac{j+1}N\right)({\mathcal{L}}_1-{\mathcal{L}}_0)$.
We know by inductive assumption that~${\mathcal{R}}$ is invertible
and by~\eqref{0248mathcalmathcal2} with~$t:=\frac{j+1}N$ that~$\|{\mathcal{R}}^{-1}\|\le K$.
Additionally, for every~$t\in\left[\frac{j+1}N,\frac{j+2}N\right]$ we have that
$$ \| {\mathcal{S}}\|=\sup_{x\in{\mathcal{X}}\setminus\{0\}}\frac{\left|
\left(t-\frac{j+1}N\right)({\mathcal{L}}_1-{\mathcal{L}}_0)x\right|}{|x|}\le
\frac1N\,\sup_{x\in{\mathcal{X}}\setminus\{0\}}\frac{|{\mathcal{L}}_1x|+|{\mathcal{L}}_0 x|}{|x|}
\le\frac{K}N\le\frac{K^2}{N\,\|{\mathcal{R}}^{-1}\|}\le\frac{1}{2\,\|{\mathcal{R}}^{-1}\|}.
$$
Thus, we are again in the position of using Lemma~\ref{LEMroofofheoreheorem6.14GT},
and we thereby conclude that~${\mathcal{R}}+{\mathcal{S}}={\mathcal{L}}_{(j+1)/N}+
\left(t-\frac{j+1}N\right)({\mathcal{L}}_1-{\mathcal{L}}_0)={\mathcal{L}}_t$
is invertible for every~$t\in\left[\frac{j+1}N,\frac{j+2}N\right]$, which completes the inductive step
and proves~\eqref{0248mathcalmathcal}.

Now, it follows from~\eqref{0248mathcalmathcal} that the operator~${\mathcal{L}}_t$ is invertible for all~$t\in[0,1]$.
In particular, the operator~${\mathcal{L}}_1$ is invertible, and thus surjective, which entails~\eqref{0248mathcalmathcal2-B2wefr}.

The claim in~\eqref{0248mathcalmathcal2-B2wefr2} then follows by taking~$t:=1$ and~$x:=x_y$ in~\eqref{0248mathcalmathcal2}.
\end{proof}

With the help of the previously developed setting in functional analysis, we can now
establish the existence result for the Dirichlet problem stated in Theorem~\ref{Theorem6.14GT}.
 
\begin{proof}[Proof of Theorem~\ref{Theorem6.14GT}] 
To start with, we prove that
the solution of the Dirichlet problem in~\eqref{134-24rf-pkwPK-TAGBco6YHSRVrfik-RiodSS02o3ef1}
is unique in~$C^2(\Omega)\cap C(\overline\Omega)$.
To check this, suppose that there are two solutions in~$C^2(\Omega)\cap C(\overline\Omega)$
and consider their difference~$v$. Then, $v$ is in turn
a solution in~$C^2(\Omega)\cap C(\overline\Omega)$ of the Dirichlet problem
\begin{equation*}
\begin{dcases}
\sum_{i,j=1}^n a_{ij} \partial_{ij}v+\sum_{i=1}^n b_i\partial_i v+cv=0
&\quad{\mbox{in }}\,\Omega,\\
v=0&\quad{\mbox{on }}\,\partial\Omega.\end{dcases}
\end{equation*}
We can apply Corollary~\ref{9oikr3eg9u2iwhefk9qryfhweiwwrfgfLftyOmega-kdf} to both~$v$ and~$-v$, thus finding that
$$ \sup_{\overline\Omega} v\le\sup_{\partial\Omega} v^+=0\qquad{\mbox{and}}\qquad
\sup_{\overline\Omega} (-v)\le\sup_{\partial\Omega} (-v)^+=0.$$
This tells us that~$v$ vanishes identically, thus establishing the desired uniqueness claim.

Now we focus on the
proof of the solvability of the Dirichlet problem in~\eqref{134-24rf-pkwPK-TAGBco6YHSRVrfik-RiodSS02o3ef1}
and on the regularity estimate in~\eqref{134-24rf-pkwPK-TAGBco6YHSRVrfik-RiodSS02o3ef2}.
We will first prove these claims under the additional assumption that
\begin{equation}\label{NgfdbMb6hODesmMAuyrF4r32a6ctSAmritkjFAvdimPAdunbtfDRAgra}\begin{split}&
{\mbox{the Dirichlet problem}}\\&
\begin{dcases}
\Delta u=f
&\quad{\mbox{in }}\,\Omega,\\
u=g&\quad{\mbox{on }}\,\partial\Omega\end{dcases}
\\&{\mbox{admits a solution in~$C^{2,\alpha}(\Omega)$.}}\end{split}\end{equation}
In this spirit, the gist is that one reduces the solvability of the Dirichlet problem for general elliptic operators
in~\eqref{134-24rf-pkwPK-TAGBco6YHSRVrfik-RiodSS02o3ef1}
(as well as the regularity estimate in~\eqref{134-24rf-pkwPK-TAGBco6YHSRVrfik-RiodSS02o3ef2})
to the case of the Laplace operator,
and this strategy somewhat underlies the fact that the general structure of the elliptic
operator (with the ``right sign of the zero order coefficient'')
does not provide structural obstructions for the solvability of the Dirichlet problem
since, roughly speaking, a domain which is ``good for the Laplacian'' is also ``equally good
for general operators''.
Thus, for the moment we assume the additional hypothesis~\eqref{NgfdbMb6hODesmMAuyrF4r32a6ctSAmritkjFAvdimPAdunbtfDRAgra}
and we will establish the existence of a solution in~$C^{2,\alpha}(\Omega)$ (which is stronger
than~$C^2(\Omega)\cap C(\overline\Omega)$, due to the global character of~$C^{2,\alpha}(\Omega)$,
see the notation in footnote~\ref{LAS0ThgdRujfghlyhn2UPksdssi984n7tE62ps} on page~\pageref{LAS0ThgdRujfghlyhn2UPksdssi984n7tE62ps}) for the Dirichlet problem in~\eqref{134-24rf-pkwPK-TAGBco6YHSRVrfik-RiodSS02o3ef1}
and at the same time we will prove that~\eqref{134-24rf-pkwPK-TAGBco6YHSRVrfik-RiodSS02o3ef2} holds true;
after that, we will then check that all domains with boundary of~$C^{2,\alpha}$ satisfy~\eqref{NgfdbMb6hODesmMAuyrF4r32a6ctSAmritkjFAvdimPAdunbtfDRAgra}
(though a bit convoluted, this will turn out to be an effective strategy,
since the previously built knowledge obtained in this way will play a pivotal role
in the subsequent stages of the proof, as one can appreciate how~\eqref{35636-42gftger3teqrt43ger8y43gfh9uvg98t43hgfeubgerSynfD4ma43ex5Wedrveve3tKMS-2re43hg}
will be exploited to check the validity of~\eqref{NgfdbMb6hODesmMAuyrF4r32a6ctSAmritkjFAvdimPAdunbtfDRAgra}
by using balls as a further intermediate step).

To complete this plan,
as observed in~\eqref{omMAEwgNGLt4IONkjgfnbdsdtyuijk9024i5SZERSD},
and keeping in mind the notation of boundary norms discussed in footnote~\ref{ESTE678i34rihf}
on page~\pageref{ESTE678i34rihf},
we recall that
it suffices to prove the solvability of the Dirichlet problem in~\eqref{134-24rf-pkwPK-TAGBco6YHSRVrfik-RiodSS02o3ef1}
and the regularity estimate in~\eqref{134-24rf-pkwPK-TAGBco6YHSRVrfik-RiodSS02o3ef2} when
\begin{equation}\label{0u9wfehVABSniSHih30-01}
{\mbox{$g$ vanishes identically.}}\end{equation}
We let~${\mathcal{X}}$ be the space of functions~$v$
in~$C^{2,\alpha}(\Omega)$ such that~$v=0$ on~$\partial\Omega$. Let also~${\mathcal{Y}}:=C^{\alpha}(\Omega)$ and
$$ Lu:=\sum_{i,j=1}^n a_{ij} \partial_{ij}u+\sum_{i=1}^n b_i\partial_i u+cu.$$
The idea is to ``connect~$L$ with the Laplace operator'', namely, for each~$t\in[0,1]$, we define
$$ L_tu:=tLu+(1-t)\Delta u.$$
We stress that~$L_t$ is a bounded linear operator from~${\mathcal{X}}$ to~${\mathcal{Y}}$ since~$\|L_tu\|_{C^{\alpha}(\Omega)}\le C\,\|u\|_{C^{2,\alpha}(\Omega)}$. Moreover, for all~$\xi=(\xi_1,\dots,\xi_n)\in\partial B_1$,
\begin{eqnarray*}
\sum_{i,j=1}^n \big(ta_{ij}+(1-t)\delta_{ij}\big)\xi_i\xi_j
\in\big[ t\lambda+(1-t),\,t\Lambda+(1-t)\big]\subseteq \big[\min\{1,\lambda\},\,\max\{1,\Lambda\}\big]
\end{eqnarray*}
and therefore~$L_t$ satisfies the ellipticity condition in~\eqref{ELLIPTIC}, up to changing the structural constants.
As a result, by the a priori estimates in Theorem~\ref{ThneGiunfBB032t4jNKS-s3i4} (recall especially~\eqref{LAMds-STudf-APEImriJEr1-2} and~\eqref{0u9wfehVABSniSHih30-01}),
\begin{equation*} \|u\|_{C^{2,\alpha}(\Omega)}\le C\,\|L_tu\|_{C^\alpha(\Omega)}.\end{equation*}
We are therefore in the setting required by~\eqref{0248mathcalmathcal2}: hence, to exploit
Corollary~\ref{i0ojendelope934etting876in979fu} we only need to check that~$L_0$ is surjective.
For this, since~$L_0$ is the Laplacian, for every~$f\in {\mathcal{Y}}= C^\alpha(\Omega)$ we
can apply~\eqref{NgfdbMb6hODesmMAuyrF4r32a6ctSAmritkjFAvdimPAdunbtfDRAgra}
to find a function~$u_f\in C^{2,\alpha}(\Omega)$
such that~$\Delta u_f=f$ in~$\Omega$ and~$u_f=0$ along~$\partial\Omega$.

{F}rom this, it follows that~$u_f=L_0^{-1}f$ and therefore~$L_0$ is surjective.
We can consequently employ Corollary~\ref{i0ojendelope934etting876in979fu} and deduce that
for every~$f\in C^\alpha(\Omega)$ there exists a unique~$u_f\in C^{2,\alpha}(\Omega)$ with~$u_f=0$ on~$\partial\Omega$
such that~${\mathcal{L}}_1 u_f=f$ in~$\Omega$ and additionally~$\|u_f\|_{C^{2,\alpha}(\Omega)}\le C\,\|f\|_{C^\alpha(\Omega)}$. These observations provide a solution for the Dirichlet problem in~\eqref{134-24rf-pkwPK-TAGBco6YHSRVrfik-RiodSS02o3ef1} which also satisfies the desired estimate in~\eqref{134-24rf-pkwPK-TAGBco6YHSRVrfik-RiodSS02o3ef2} (recall~\eqref{0u9wfehVABSniSHih30-01}).

Hence, 
\begin{equation}\label{35636-42gftger3teqrt43ger8y43gfh9uvg98t43hgfeubgerSynfD4ma43ex5Wedrveve3tKMS-2re43hg}
\begin{split}&
{\mbox{we have completed the proof of Theorem~\ref{Theorem6.14GT}}}\\&{\mbox{under the additional assumption in~\eqref{NgfdbMb6hODesmMAuyrF4r32a6ctSAmritkjFAvdimPAdunbtfDRAgra}.}}\end{split}\end{equation}
To finish, we now need to establish that
\begin{equation}\label{35636-42gftger3teqrt43ger8y43gfh9uvg98t43hgfeubgerSynfD4ma43ex5Wedrveve3tKMS-2re43hgTT}
{\mbox{all domains with boundary of class~$C^{2,\alpha}$
satisfy~\eqref{NgfdbMb6hODesmMAuyrF4r32a6ctSAmritkjFAvdimPAdunbtfDRAgra}. }}\end{equation}
To this end, 
we observe that~\eqref{NgfdbMb6hODesmMAuyrF4r32a6ctSAmritkjFAvdimPAdunbtfDRAgra}
holds true when~$\Omega$ is a ball, thanks to Theorem~\ref{KS-2rkjoewuj5-2prlfgk-a}.
{F}rom this and~\eqref{35636-42gftger3teqrt43ger8y43gfh9uvg98t43hgfeubgerSynfD4ma43ex5Wedrveve3tKMS-2re43hg},
we deduce that
\begin{equation}\label{35636-42gftger3teqrt43ger8y43gfh9uvg98t43hgfeubgerSynfD4ma43ex5Wedrveve3tKMS-2re43hg-2}
{\mbox{Theorem~\ref{Theorem6.14GT} holds true when~$\Omega$ is a ball.}}\end{equation}
Now, to deal with the proof of~\eqref{35636-42gftger3teqrt43ger8y43gfh9uvg98t43hgfeubgerSynfD4ma43ex5Wedrveve3tKMS-2re43hgTT},
we recall that, in light of
Corollary~\ref{S-coroEXIS-M023} (recall also~\eqref{SqerLPONE-0-9ieyhdfweohgiowegb74uthghhs1334756ujhgf}),
we can find a solution~$u\in C^2(\Omega)\cap C(\overline\Omega)$ of the Dirichlet problem for the Laplace operator in~\eqref{NgfdbMb6hODesmMAuyrF4r32a6ctSAmritkjFAvdimPAdunbtfDRAgra}. Therefore,
it only remains to show that
\begin{equation}\label{MAbdyStrtso43ceuvnvmoprojgleantuworqyi7596mnf8thatdfwnedidenga}
{\mbox{this solution belongs to~$C^{2,\alpha}(\Omega)$.}}
\end{equation}
As a matter of fact, we already know from the interior estimates for the Laplace operator given in Theorem~\ref{SCHAUDER-INTE}
that
\begin{equation}\label{M2423AbdyStrtso43ceuvnvmoprojgleantuworqyi7596mnf8thatdfwnedidenga-BJS0-djfd}
u\in C^{2,\alpha}_{\rm loc}(\Omega),\end{equation} therefore we only need to obtain uniform estimates
of~$C^{2,\alpha}$ type in a neighborhood of~$\partial\Omega$. To this end, we take~$\rho>0$ sufficiently small,
a finite family of points~$p_1,\dots,p_N\in\partial\Omega$, and sets with $C^{2,\alpha}$ boundaries~$U_1,\dots,U_N$
with~$B_{\rho/3}(p_j)\cap\Omega\Subset U_j\Subset B_{\rho/2}(p_j)$, such that
\begin{equation}\label{ONS-9uyhdwCAncoo0kdRAYTmdndvp} {\mathcal{N}}:=\{ x\in\Omega {\mbox{ s.t. }} B_{\rho/8}(x)\cap(\partial\Omega)\ne\varnothing\}
\subseteq \bigcup_{j=1}^N B_{\rho/4}(p_j)\end{equation}
and with~$U_j$ equivalent via a~$C^{2,\alpha}$ diffeomorphism~$T_j$
with~$C^{2,\alpha}$ inverse to the ball~$B_1$, see Figure~\ref{TERRAeF5636-42gftger3teqrt43ger8y43gfh9uvg98t43hgfeubgerSynfD4ma43ex5Wedrveve3tKMS-2re43hgTT}.

\begin{figure}
  \centering
  \includegraphics[width=.9\linewidth]{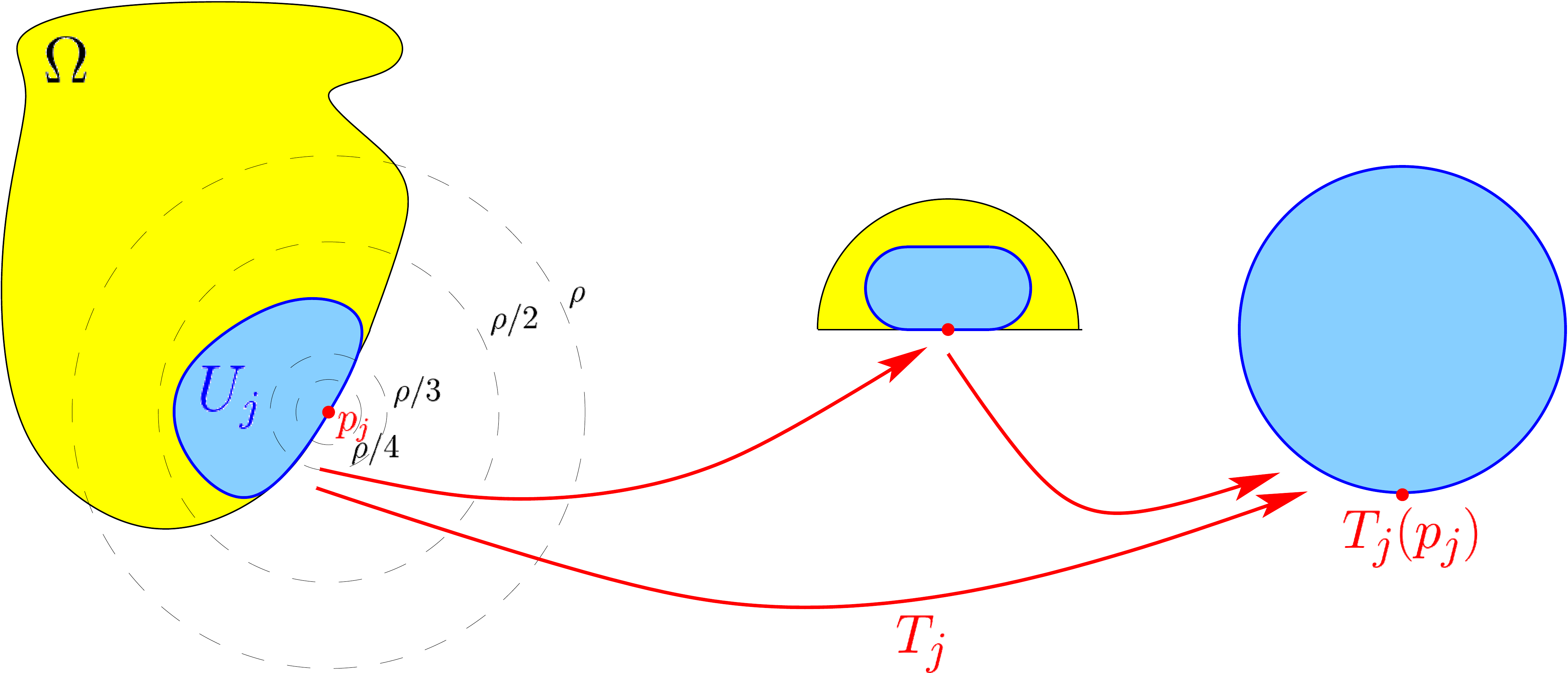}
 \caption{\sl Reducing to a ball to prove~\eqref{35636-42gftger3teqrt43ger8y43gfh9uvg98t43hgfeubgerSynfD4ma43ex5Wedrveve3tKMS-2re43hgTT}.}\label{TERRAeF5636-42gftger3teqrt43ger8y43gfh9uvg98t43hgfeubgerSynfD4ma43ex5Wedrveve3tKMS-2re43hgTT}
\end{figure}

Thus, for each~$y\in T_j(U_j)=B_1$ we let~$\widetilde u(y):=u(T_j^{-1}(y))$ and we employ Lemma~\ref{92o3jwelg02jfj-qiewjdfnoewfh9ewohgveioqNONDIFO234RM} to see that~$\widetilde u$ satisfies an equation
of the type
$$ \sum_{i,j=1}^n \widetilde a_{ij}(y)\partial_{ij}
\widetilde u(y)+\sum_{i=1}^n \widetilde b_i(y)\partial_i \widetilde u(y)=\widetilde f(y)$$
for every~$y\in B_1$, with~$\widetilde a_{ij}$, $\widetilde b_i$, $\widetilde f\in C^\alpha(\Omega)$.
Additionally, for all~$y\in T_j((\partial U_j)\cap(\partial\Omega))$ we have that~$\widetilde u(y)=g(T_j^{-1}(y))=:\widetilde{g}(y)$
and we stress that~$\widetilde{g}$ is also a~$C^{2,\alpha}$ function.

Now we take a sequence of functions~$\widetilde u_k\in C^{2,\alpha}(B_1)$ such that~$\widetilde u_k=\widetilde u=\widetilde g$
in~$T_j(B_{\rho/3}(p_j)\cap\Omega)$ and~$\widetilde u_k\to \widetilde u$ uniformly in~$\overline{B_1}$.
By~\eqref{35636-42gftger3teqrt43ger8y43gfh9uvg98t43hgfeubgerSynfD4ma43ex5Wedrveve3tKMS-2re43hg-2}
we can solve the Dirichlet problem
\begin{equation}\label{STj0s6tOJNS0qpodjwnwhfiuosgfigweVmiE} \begin{dcases}
\sum_{i,j=1}^n \widetilde a_{ij}\partial_{ij}
\widetilde v_k +\sum_{i=1}^n \widetilde b_i \partial_i \widetilde v_k=\widetilde f &{\mbox{ in }}B_1,\\
\widetilde v_k=\widetilde u_k & {\mbox{ on }}\partial B_1,
\end{dcases}\end{equation}
with~$\widetilde v_k\in C^{2,\alpha}(B_1)$.

Also, the function~$\widetilde w_{k\ell}:=\widetilde v_k-\widetilde v_\ell$ satisfies
$$ \begin{dcases}
\sum_{i,j=1}^n \widetilde a_{ij}\partial_{ij}
\widetilde w_{k\ell} +\sum_{i=1}^n \widetilde b_i \partial_i \widetilde w_{k\ell}=0 &{\mbox{ in }}B_1,\\
\widetilde w_{k\ell}=\widetilde u_k-\widetilde u_\ell & {\mbox{ on }}\partial B_1.
\end{dcases}$$
Thus, by Corollary~\ref{9oikr3eg9u2iwhefk9qryfhweiwwrfgfLftyOmega-kdf}
(applied here to both~$\widetilde w_{k\ell}$ and~$-\widetilde w_{k\ell}$ with~$f:=0$), we see that
$$ \|\widetilde w_{k\ell}\|_{L^\infty(B_1)}\le\|\widetilde w_{k\ell}\|_{L^\infty(\partial B_1)}
=\|\widetilde u_{k}-\widetilde u_\ell\|_{L^\infty(\partial B_1)}.$$
This gives that~$\widetilde v_{k}$ converges uniformly in~$B_1$
to some function~$\widetilde v_\star$.
In particular, we have that~$\widetilde v_\star\in C(\overline{B_1})$.

Moreover, using Cauchy's Estimates (see Theorem~\ref{CAUESTIMTH}), for all~$q\in B_1$
and~$r\in(0,1)$ such that~$B_r(q)\Subset B_1$, we see that
$$ \| \widetilde w_{k\ell}\|_{C^3(B_{r/2}(q))}\le \frac{C}{r^{3}}\|
\widetilde w_{k\ell}\|_{L^\infty(B_1)}.$$
This entails that~$\widetilde v_k$ is a Cauchy sequence in~$C^3(B_{r/2}(q))$
hence, in~$B_{r/2}(q)$,
$$ \widetilde f=\lim_{k\to+\infty}
\sum_{k=1}^n \widetilde a_{ij} \partial_{ij} \widetilde v_k +
\sum_{k=1}^n \widetilde b_{i}  \partial_i \widetilde v_k
=\sum_{k=1}^n \widetilde a_{ij} \partial_{ij} \widetilde v_\star +
\sum_{k=1}^n \widetilde b_{i}  \partial_i \widetilde v_\star.$$
This gives that~$\widetilde v_\star\in C^2(B_1)\cap C(\overline{B_1})$
solves the Dirichlet problem
$$ \begin{dcases}
\sum_{k=1}^n \widetilde a_{ij} \partial_{ij} \widetilde v_\star +
\sum_{k=1}^n \widetilde b_{i}  \partial_i \widetilde v_\star=\widetilde f & {\mbox{ in }}B_1,\\
\widetilde v_\star=\widetilde u & {\mbox{ on }}\partial B_1.
\end{dcases}$$
By the uniqueness claim for the Dirichlet problem
(e.g., using again~\eqref{35636-42gftger3teqrt43ger8y43gfh9uvg98t43hgfeubgerSynfD4ma43ex5Wedrveve3tKMS-2re43hg-2}), we conclude that
\begin{equation}\label{ssaerdelmidclsaaKSMche21irf0-whfgewuvnie}
\widetilde v_\star=\widetilde u.
\end{equation}
In addition, we can apply the
local estimates at the boundary given in Theorem~\ref{THM:0-T536EcdhPOaIJodidpo3tju24ylkiliR2543uj8890j09092po-okf-3edk-3edhf-2-fh564-LEATB}
to the function~$\widetilde v_k$
(applied here to a straightening of~$\partial B_1$ in the vicinity of~$T(p_j)$).
In this way, we find that
$$ \|\widetilde v_k\|_{C^{2,\alpha}(T_j( B_{\rho/4}(p_j)\cap\Omega ))}\le C\,\Big(\|\widetilde v_k\|_{L^\infty(B_1)}+
\|\widetilde g\|_{C^{2,\alpha}((T_j((\partial\Omega)\cap B_{\rho/2}(p_j)))}+\|\widetilde f\|_{C^\alpha(B_1)}\Big).$$
Consequently, passing to the limit as~$k\to+\infty$,
$$ \|\widetilde v_\star\|_{C^{2,\alpha}(T_j(B_{\rho/4}(p_j)\cap\Omega))}\le C\,\Big(\|\widetilde v_\star\|_{L^\infty(B_1)}+
\|\widetilde g\|_{C^{2,\alpha}((T_j((\partial\Omega)\cap B_{\rho/2}(p_j)))}+\|\widetilde f\|_{C^\alpha(B_1)}\Big).$$
Hence, by~\eqref{ssaerdelmidclsaaKSMche21irf0-whfgewuvnie},
$$ \|\widetilde u\|_{C^{2,\alpha}(T_j(B_{\rho/4}(p_j)\cap\Omega))}\le C\,\Big(\|\widetilde u\|_{L^\infty(B_1)}+
\|\widetilde g\|_{C^{2,\alpha}((T_j((\partial\Omega)\cap B_{\rho/2}(p_j)))}+\|\widetilde f\|_{C^\alpha(B_1)}\Big).$$
Transforming back, and possibly renaming constants, we obtain that
$$ \| u\|_{C^{2,\alpha}( B_{\rho/4}(p_j) \cap\Omega)}\le C\,\Big(\|u\|_{L^\infty(B_1)}+
\| g\|_{C^{2,\alpha}( \partial\Omega)}+\|f\|_{C^\alpha(\Omega)}\Big).$$
This and~\eqref{ONS-9uyhdwCAncoo0kdRAYTmdndvp} yield that
$$ \| u\|_{C^{2,\alpha}( {\mathcal{N}}) }\le C\,\Big(\|u\|_{L^\infty(B_1)}+
\| g\|_{C^{2,\alpha}( \partial\Omega)}+\|f\|_{C^\alpha(\Omega)}\Big),$$
up to renaming constants once again.

The latter estimate, together with~\eqref{M2423AbdyStrtso43ceuvnvmoprojgleantuworqyi7596mnf8thatdfwnedidenga-BJS0-djfd} (and a covering argument as in Lemma~\ref{COVERINGARG}), gives that 
$$ \| u\|_{C^{2,\alpha}( \Omega) }\le C\,\Big(\|u\|_{L^\infty(B_1)}+
\| g\|_{C^{2,\alpha}( \partial\Omega)}+\|f\|_{C^\alpha(\Omega)}\Big).$$
Using this and Corollary~\ref{9oikr3eg9u2iwhefk9qryfhweiwwrfgfLftyOmega-kdf}, we obtain that
$$ \| u\|_{C^{2,\alpha}( \Omega) }\le C\,\Big(
\| g\|_{C^{2,\alpha}( \partial\Omega)}+\|f\|_{C^\alpha(\Omega)}\Big),$$
which finishes the proof of~\eqref{MAbdyStrtso43ceuvnvmoprojgleantuworqyi7596mnf8thatdfwnedidenga}.
The proof of~\eqref{35636-42gftger3teqrt43ger8y43gfh9uvg98t43hgfeubgerSynfD4ma43ex5Wedrveve3tKMS-2re43hgTT}
is thereby complete as well.
\end{proof}

We observe that the ellipticity assumption 
in~\eqref{ELLIPTIC} is an essential ingredient for the solvability of the Dirichlet problem
obtained in Theorem~\ref{Theorem6.14GT}. For instance,
in this setting, assumption~\eqref{ELLIPTIC} cannot be weakened by assuming
that for every~$\xi=(\xi_1,\dots,\xi_n)\in\partial B_1$
\begin{equation}\label{ELLIPTIC-DEGENERATE}
\sum_{i,j=1}^n a_{ij}(x)\xi_i\xi_j\in [0,\Lambda]
\end{equation}
for some~$\Lambda\ge0$ (that is, one cannot take~$\lambda=0$ in~\eqref{ELLIPTIC}).
As an instructive example of this degeneracy,
we observe that when~$n=2$ the Dirichlet problem \begin{equation}\label{FVbdImRmfKe03}\begin{dcases} x_2^2\partial_{11}u+x_1^2\partial_{22}u-2x_1x_2\partial_{12}u-x_1\partial_1u-x_2\partial_2u=1&{\mbox{ in }}B_1\subset\R^2,\\ u=0 &{\mbox{ on }}\partial B_1\end{dcases}\end{equation} does not possess any solution in~$C^{2,\alpha}(B_1)$ for any~$\alpha\in(0,1)$. Indeed, if such a solution existed, we deduce from the boundary value prescription that~$u(\cos\vartheta,\sin\vartheta)=0$ for all~$\vartheta\in\R$. As a result, taking one derivative in~$\vartheta$, 
\[ 0=-\sin\vartheta\,\partial_1 u(\cos\vartheta,\sin\vartheta)+\cos\vartheta\,\partial_2 u(\cos\vartheta,\sin\vartheta).\] By taking another derivative we find that \begin{eqnarray*}&&0=\sin^2\vartheta\,\partial_{11} u(\cos\vartheta,\sin\vartheta)+\cos^2\vartheta\,\partial_{22} u(\cos\vartheta,\sin\vartheta)-2\sin\vartheta\cos\vartheta\,\partial_{12} u(\cos\vartheta,\sin\vartheta)\\&&\qquad\qquad\qquad\qquad
-\cos\vartheta\,
\partial_1 u(\cos\vartheta,\sin\vartheta)-\sin\vartheta\,\partial_2 u(\cos\vartheta,\sin\vartheta)
.\end{eqnarray*}
But the equation in~\eqref{FVbdImRmfKe03} gives that this quantity is equal to~$1$, providing a contradiction. This shows that no solution in~$C^{2,\alpha}(B_1)$ is admissible for the Dirichlet problem in~\eqref{FVbdImRmfKe03}.
We remark that this example fulfills the degenerate ellipticity condition in~\eqref{ELLIPTIC-DEGENERATE}
(but not the strong version in~\eqref{ELLIPTIC}). Indeed, in this case~$a_{11}=x_2^2$, $a_{22}=x_1^2$
and~$a_{12}=a_{21}=-x_1x_2$ and thus,
for every~$x\in{B_1}$ and~$\xi=(\xi_1,\dots,\xi_n)\in\partial B_1$,
\begin{eqnarray*}
\sum_{i,j=1}^n a_{ij}(x)\xi_i\xi_j
=x_2^2 \xi_1^2+x_1^2\xi_2^2-2x_1x_2\xi_1\xi_2=
(x_2\xi_1-x_1\xi_2)^2
\in [0,4]
\end{eqnarray*}
with the vanishing threshold attained for~$\xi:=\frac{x}{|x|}$ when~$x\in B_1\setminus\{0\}$
(and for all~$\xi\in\partial B_1$ when~$x=0$).\medskip

It is also interesting to observe that the approximation method employed at the end of the proof
of Theorem~\ref{Theorem6.14GT} (namely, the technique adopted for~\eqref{STj0s6tOJNS0qpodjwnwhfiuosgfigweVmiE}) can be of general use
and leads to a number of extensions of Theorem~\ref{Theorem6.14GT} in case of
``less regular data''. As an example of this feature, we show that for continuous boundary
datum one can still solve the Dirichlet problem continuously up to the boundary, as made precise in
the following result:

\begin{theorem}\label{oind-94984u5f-14UNSpas88saleohi4a2T}
Let~$\Omega\subset\R^n$ be a bounded open set with boundary of class~$C^{2,\alpha}$
for some~$\alpha\in(0,1)$.
Let~$a_{ij}$, $b_i$, $c$, $f\in C^\alpha(\Omega)$. Let also~$g\in C(\partial\Omega)$.
Assume that~$c(x)\le0$ for all~$x\in\Omega$ and that the ellipticity 
condition in~\eqref{ELLIPTIC} holds true.

Then, the Dirichlet problem
\begin{equation}\label{oind-94984u5f-14UNSpas88saleohi4a2}
\begin{dcases}
\sum_{i,j=1}^n a_{ij} \partial_{ij}u+\sum_{i=1}^n b_i\partial_i u+cu=f
&\quad{\mbox{in }}\,\Omega,\\
u=g&\quad{\mbox{on }}\,\partial\Omega\end{dcases}
\end{equation}
has a unique solution in~$C^2(\Omega)\cap C(\overline\Omega)$.

Additionally, such a solution~$u$ belongs to~$C^{2,\alpha}_{\rm loc}(\Omega)$ and for every~$\Omega'\Subset\Omega$
there exists~$C>0$, depending only on~$n$,
$\alpha$, $\Omega'$, $\Omega$, $a_{ij}$, $b_i$ and~$c$, such that
\begin{equation}\label{oind-94984u5f-14UNSpas88saleohi4a3}
\|u\|_{C^{2,\alpha}(\Omega')}\le C\,\Big(\|g\|_{L^\infty(\partial\Omega)}+\|f\|_{C^\alpha(\Omega)}\Big).\end{equation}
\end{theorem}

\begin{proof} The uniqueness claim follows from Corollary~\ref{9oikr3eg9u2iwhefk9qryfhweiwwrfgfLftyOmega-kdf},
hence we focus on the existence claim.

Let~$g_k\in C^{2,\alpha}(\partial\Omega)$ be such that~$g_k\to g$ in~$L^\infty(\partial\Omega)$.
We can use Theorem~\ref{Theorem6.14GT} and find~$u_k\in C^{2,\alpha}(\Omega)$
that solves the Dirichlet problem
\begin{equation}\label{oind-94984u5f-14UNSpas88saleohi4a1}
\begin{dcases}
\sum_{i,j=1}^n a_{ij} \partial_{ij}u_k+\sum_{i=1}^n b_i\partial_i u_k+cu_k=f
&\quad{\mbox{in }}\,\Omega,\\
u_k=g_k&\quad{\mbox{on }}\,\partial\Omega.\end{dcases}
\end{equation}
We also let~$w_{k\ell}:=u_k-u_\ell$ and we observe that
\begin{equation*}
\begin{dcases}
\sum_{i,j=1}^n a_{ij} \partial_{ij}w_{k\ell}+\sum_{i=1}^n b_i\partial_i w_{k\ell}+cw_{k\ell}=0
&\quad{\mbox{in }}\,\Omega,\\
w_{k\ell}=g_k-g_\ell&\quad{\mbox{on }}\,\partial\Omega.\end{dcases}
\end{equation*}
{F}rom Corollary~\ref{9oikr3eg9u2iwhefk9qryfhweiwwrfgfLftyOmega-kdf} we deduce that
$$ \| u_k-u_\ell\|_{L^\infty(\Omega)}=\| w_{k\ell}\|_{L^\infty(\Omega)}
\le\| w_{k\ell}\|_{L^\infty(\partial\Omega)}=
\| g_k-g_\ell\|_{L^\infty(\partial\Omega)}
.$$
As a result, $u_k$ is a Cauchy sequence in~$L^\infty(\Omega)$,
hence it converges uniformly in~$\Omega$ to some function~$u\in C(\overline\Omega)$.

Furthermore, for all~$x_0\in\Omega$ and all~$r>0$ such that~$B_{r}(x_0)\Subset\Omega$,
we can use the interior estimates in Theorem~\ref{SCHAUDER-INTE} and deduce that
$$ \|u_k\|_{C^{2,\alpha}(B_{r/2}(x_0))}\le C_r\,\Big(\|u_k\|_{L^\infty(B_r(x_0))}+\|f\|_{C^\alpha(B_r(x_0))}\Big)\le
C_r\,\Big(\|u_k\|_{L^\infty(\Omega)}+\|f\|_{C^\alpha(\Omega)}\Big).$$
with~$C_r>0$, depending only on~$r$, $n$ and~$\alpha$.

Since, using again Corollary~\ref{9oikr3eg9u2iwhefk9qryfhweiwwrfgfLftyOmega-kdf}, it holds that
$$ \|u_k\|_{L^\infty(\Omega)}\le\|u_k\|_{L^\infty(\partial\Omega)} +C\|f\|_{L^\infty(\Omega)}=
\|g_k\|_{L^\infty(\partial\Omega)} +C\|f\|_{L^\infty(\Omega)},$$
we find that
\begin{equation}\label{oind-94984u5f-14UNSpas88saleohi4a} \|u_k\|_{C^{2,\alpha}(B_{r/2}(x_0))}\le C_r\,\Big(\|g_k\|_{L^\infty(\partial\Omega)}+\|f\|_{C^\alpha(\Omega)}\Big).\end{equation}
In particular, for large~$k$,
$$ \|u_k\|_{C^{2,\alpha}(B_{r/2}(x_0))}\le 
C_r\,\Big(1+\|g\|_{L^\infty(\partial\Omega)}+\|f\|_{C^\alpha(\Omega)}\Big).$$
This gives a uniform bound for~$\|u_k\|_{C^{2,\alpha}(B_{r/2}(x_0))}$
and therefore, up to a subsequence, we have that~$u_k\to u$ in~$C^2(\Omega')$ for all~$\Omega'\Subset\Omega$.
We can therefore pass to the limit as~$k\to+\infty$ in~\eqref{oind-94984u5f-14UNSpas88saleohi4a1}
and obtain~\eqref{oind-94984u5f-14UNSpas88saleohi4a2}.
Similarly,
passing to the limit as~$k\to+\infty$ in~\eqref{oind-94984u5f-14UNSpas88saleohi4a},
we obtain~\eqref{oind-94984u5f-14UNSpas88saleohi4a3}.
\end{proof}

As a side remark, we observe that Theorem~\ref{oind-94984u5f-14UNSpas88saleohi4a2T}
allows one to say that~$C^2$ solutions are automatically~$C^{2,\alpha}_{\rm loc}$
if the data permit and according the ``spurious'' assumption~$u\in C^{2,\alpha}(B_1)$
in Theorem~\ref{THM:0-T536EcdhPOaIJodidpo3tju24ylkiliR2543uj8890j09092po-okf-3edk-3edhf-2-fh564}
can be weakened in favor of~$u\in C^2(B_1)$. The details of these observations are contained in the following result:

\begin{corollary}\label{THM:0-T536EcdhPOaIJodidpo3tju24ylkiliR2543uj8890j09092po-okf-3edk-3edhf-2-fh564COR}
Let~$\Omega\subset\R^n$ be open and bounded.
Let~$a_{ij}$, $b_i$, $c$, $f\in C^\alpha(\Omega)$ for some~$\alpha\in(0,1)$. Assume the ellipticity 
condition in~\eqref{ELLIPTIC}.
Let~$u\in C^{2}(\Omega)$ be a solution of
\begin{equation*}
\sum_{i,j=1}^n a_{ij} \partial_{ij}u+\sum_{i=1}^n b_i\partial_i u+cu=f
\quad{\mbox{in }}\,\Omega.\end{equation*}
Then, $u\in C^{2,\alpha}_{\rm loc}(\Omega)$ and for all~$\Omega'\Subset\Omega$
there exists~$C>0$, depending only on~$n$, $\Omega$, $\Omega'$,
$\alpha$, $a_{ij}$, $b_i$ and~$c$, such that
\begin{equation}\label{hgFYl0W7HJBo945DnAGN-1} \|u\|_{C^{2,\alpha}(\Omega')}\le C\,\Big(\|u\|_{L^\infty(\Omega)}+\|f\|_{C^\alpha(\Omega)}\Big).\end{equation}
\end{corollary}

\begin{proof} We take~$\Omega''$ and~$\Omega'''$ to be open, with $C^{2,\alpha}$ boundary and such that~$\Omega'\Subset\Omega''\Subset\Omega'''\Subset\Omega$. Notice that~$u$ is continuous along~$\partial\Omega'''$
and~$\widetilde f:=f-cu\in C^\alpha(\Omega''')$,
hence we can employ\footnote{Notice that we can apply 
Theorem~\ref{oind-94984u5f-14UNSpas88saleohi4a2T} for problem~\eqref{wevafkgrgicanapply}
since~$c=0$ in this setting.}
Theorem~\ref{oind-94984u5f-14UNSpas88saleohi4a2T} and find a solution~$v\in C^2(\Omega''')\cap C(\overline{\Omega'''})$
of the Dirichlet problem
\begin{equation}\label{wevafkgrgicanapply}
\begin{dcases}
\sum_{i,j=1}^n a_{ij} \partial_{ij}v+\sum_{i=1}^n b_i\partial_i v=\widetilde f
&\quad{\mbox{in }}\,\Omega''',\\
v=u&\quad{\mbox{on }}\,\partial\Omega'''\end{dcases}
\end{equation}
with~$v\in C^{2,\alpha}_{\rm loc}(\Omega''')$.
%% \begin{equation}\label{hgFYl0W7HJBo945DnAGN-2X}
%% \|v\|_{C^{2,\alpha}(\Omega'')}\le C\,\Big(\|u\|_{L^\infty(\partial\Omega''')}+\|\widetilde f\|_{C^\alpha(\Omega''')}\Big)\le
%% C\,\Big(\|u\|_{C^\alpha(\Omega''')}+\|f\|_{C^\alpha(\Omega''')}\Big).\end{equation}

Furthermore,
we observe that the function~$w:=u-v$ satisfies
\begin{equation*}
\begin{dcases}
\sum_{i,j=1}^n a_{ij} \partial_{ij}w+\sum_{i=1}^n b_i\partial_i w=0
&\quad{\mbox{in }}\,\Omega''',\\
w=0&\quad{\mbox{on }}\,\partial\Omega'''\end{dcases}
\end{equation*}
and accordingly
Corollary~\ref{9oikr3eg9u2iwhefk9qryfhweiwwrfgfLftyOmega-kdf} entails that
$$ \|u-v\|_{L^\infty(\Omega''')}=\|w\|_{L^\infty(\Omega''')}\le
\|w\|_{L^\infty(\partial\Omega''')}=0.$$
Therefore, $u=v$ in~$\Omega'''$, and as a consequence~$u\in C^{2,\alpha}(\Omega'')$.

Hence, we can use Theorem~\ref{THM:0-T536EcdhPOaIJodidpo3tju24ylkiliR2543uj8890j09092po-okf-3edk-3edhf-2-fh564}
and a covering argument (recall Lemma~\ref{COVERINGARG})
to conclude that
$$  \|u\|_{C^{2,\alpha}(\Omega')}\le C\,\Big(\|u\|_{L^\infty(\Omega'')}+\|f\|_{C^\alpha(\Omega'')}\Big),$$
which implies the desired estimate in~\eqref{hgFYl0W7HJBo945DnAGN-1}.
\end{proof}

We refer to~\cite{MR1814364, SHASU, MR3012036, MR3931747} for further information
about Schauder estimates and for solvability results for Dirichlet problems.
It is worth pointing out that Schauder estimates can be also used
to establish the existence theory of general  type of Dirichlet problems by
revisiting the Perron method, see e.g.~\cite[Theorem~6.11]{MR1814364}.\medskip

It is also useful to notice that Schauder estimates may be ``iterated'', or ``bootstrapped'', in order to
achieve higher regularity, whenever the data allow one to do so. In particular, we have the following result:

\begin{theorem}\label{9827hfCXVTHM:0-T536EcdhPOaIJodidfpo3tju24ylkiliR2543uj8890j09092po-okf-3edk-3edhf-2-fh564COR}
Let~$k\in\N$.
Let~$\Omega\subset\R^n$ be open and bounded.
Let~$a_{ij}$, $b_i$, $c$, $f\in C^{k,\alpha}(\Omega)$ for some~$\alpha\in(0,1)$. Assume the ellipticity 
condition in~\eqref{ELLIPTIC}.
Let~$u\in C^{2}(\Omega)$ be a solution of
\begin{equation}\label{K-edqtbdjkti904MEbvfSjmng10-VjhdKl-PSj}
\sum_{i,j=1}^n a_{ij} \partial_{ij}u+\sum_{i=1}^n b_i\partial_i u+cu=f
\quad{\mbox{in }}\,\Omega.\end{equation}
Then, $u\in C^{k+2,\alpha}_{\rm loc}(\Omega)$ and for all~$\Omega'\Subset\Omega$
there exists~$C>0$, depending only on~$n$, $k$, $\Omega$, $\Omega'$,
$\alpha$, $a_{ij}$, $b_i$ and~$c$, such that
\begin{equation*} \|u\|_{C^{k+2,\alpha}(\Omega')}\le C\,\Big(\|u\|_{L^\infty(\Omega)}+\|f\|_{C^{k,\alpha}(\Omega)}\Big).\end{equation*}
\end{theorem}

\begin{proof} When~$k=0$, the result is contained in
Corollary~\ref{THM:0-T536EcdhPOaIJodidpo3tju24ylkiliR2543uj8890j09092po-okf-3edk-3edhf-2-fh564COR}.
To prove it for every~$k\in\N$, we thus proceed by induction,
assuming that the desired result holds true for some~$\ell\in\N\cap[0,k-1]$ and proving its validity for~$\ell+1$.

The gist of the argument is that one can ``differentiate the equation''
and observe that the derivative of~$u$ satisfies a similar equation, thus reducing
the regularity of the derivative of~$u$ to the previously known step of the induction.
There is only one small catch in this argument: namely, the inductive assumption only gives
that~$u\in C^{\ell+2,\alpha}_{\rm loc}$
but to apply the previously known step to the derivative of~$u$ we would need to know that this derivative is in~$C^2$.
That is, this strategy would only work when~$\ell+2\ge3$, that is~$\ell\ge1$. This would create an issue precisely
for going to~$\ell=0$ to~$\ell=1$ in the inductive step.
To avoid this caveat, it is convenient to work (at least when~$\ell=0$, but also for all~$\ell\in\N\cap[0,k-1]$ for the sake of uniformity)
with discrete increments rather then derivatives. The technical adjustments needed are as follows.

Let~$\beta\in\N^n$ with~$|\beta|=\ell\in\N$ and~$v:=\partial^\beta u$.
We stress that this is a good definition, since the inductive assumption gives us that~$u\in C^{\ell+2,\alpha}_{\rm loc}(\Omega)$
and in fact~$v\in C^{2,\alpha}_{\rm loc}(\Omega)$. Differentiating~\eqref{K-edqtbdjkti904MEbvfSjmng10-VjhdKl-PSj}
$\ell$ times, we find that, in~$\Omega$,
\begin{equation} \sum_{i,j=1}^n a_{ij} \partial_{ij} v+\sum_{i=1}^n \widetilde b_i\partial_i v+\widetilde cv=\widetilde f,\end{equation}
for suitable~$\widetilde b_i$, $\widetilde c$, $\widetilde f\in C^{k-\ell,\alpha}(\Omega)$.

Given~$h>0$ and~$e\in \partial B_1$, we define
$$ w_h(x):=\frac{v(x+he)-v(x)}{h}.
$$
Hence we take~$\Omega''$ and~$\Omega'''$
such that~$\Omega'\Subset\Omega''\Subset\Omega'''\Subset\Omega$ and we deduce that if~$h$
is conveniently small then in~$\Omega'''$ we have that
\begin{equation*} \sum_{i,j=1}^n a_{ij} \partial_{ij} w_h+\sum_{i=1}^n \widetilde b_{i,h}\partial_i w_h+\widetilde c_h w_h=\widetilde f_h,\end{equation*}
for suitable~$\widetilde b_{i,h}$, $\widetilde c_h$, $\widetilde f_h\in C^{k-\ell-1,\alpha}(\Omega)$.
Since~$k-\ell-1\ge0$, we can utilize
Corollary~\ref{THM:0-T536EcdhPOaIJodidpo3tju24ylkiliR2543uj8890j09092po-okf-3edk-3edhf-2-fh564COR}
and deduce that
\begin{equation*} \|w_h\|_{C^{2,\alpha}(\Omega')}\le C\,\Big(\|w_h\|_{L^\infty(\Omega'')}
+\|\widetilde f_h\|_{C^\alpha(\Omega'')}\Big)\le
C\,\Big(\|v\|_{C^1(\Omega''')}+\|f\|_{C^{k,\alpha}(\Omega)}\Big).\end{equation*}
This and the inductive assumption yield that
\begin{equation}\label{nGsaASmnZAjfZekf8LA0-Pjf} \|w_h\|_{C^{2,\alpha}(\Omega')}\le
C\,\Big(\|u\|_{L^\infty(\Omega)}+\|f\|_{C^{k,\alpha}(\Omega)}\Big),\end{equation}
up to renaming~$C$.

We can therefore apply the the Arzel\`a-Ascoli Theorem and obtain, up to a subsequence, that~$w_h$
converges as~$h\searrow0$ in~$C^2(\Omega')$. As a result, passing to the limit as~$h\searrow0$ in~\eqref{nGsaASmnZAjfZekf8LA0-Pjf},
$$ \|\partial_e v\|_{C^{2,\alpha}(\Omega')}\le
C\,\Big(\|u\|_{L^\infty(\Omega)}+\|f\|_{C^{k,\alpha}(\Omega)}\Big).$$
{F}rom this we infer that
$$ \|D^{\ell+3} u\|_{L^\infty(\Omega')}+[D^{\ell+3} u]_{C^\alpha(\Omega')}\le
C\,\Big(\|u\|_{L^\infty(\Omega)}+\|f\|_{C^{k,\alpha}(\Omega)}\Big).$$
This and the inductive assumption lead to
$$ \|u\|_{C^{\ell+3,\alpha}(\Omega')}\le
C\,\Big(\|u\|_{L^\infty(\Omega)}+\|f\|_{C^{k,\alpha}(\Omega)}\Big),$$
which completes the inductive step.
\end{proof}

It follows from Theorem~\ref{9827hfCXVTHM:0-T536EcdhPOaIJodidfpo3tju24ylkiliR2543uj8890j09092po-okf-3edk-3edhf-2-fh564COR} that solutions are~$C^\infty$ when so are~$a_{ij}$, $b_i$, $c$, $f$.
A regularity theory in the class of real analytic functions is also possible, see e.g.~\cite[Section~6.6]{MR2492985}
and~\cite[pages 207--210]{MR598466}
(or~\cite[Theorem~9.5.1]{MR1996773} for a more general setting).
See also~\cite{MR75415, MR89334, MR106336, MR107081}
for additional information on real analytic regularity theory.

A counterpart of Theorem~\ref{9827hfCXVTHM:0-T536EcdhPOaIJodidfpo3tju24ylkiliR2543uj8890j09092po-okf-3edk-3edhf-2-fh564COR} in terms of boundary regularity holds true as well, see~\cite[Theorem~6.19]{MR1814364} for full\index{Schauder estimates|)} details.

\chapter{Equations in nondivergence form: $W^{2,p}$-regularity theory}\label{C2ALPHACHAPETil2}

\section{Hints and limitations for a regularity theory in Lebesgue spaces}

The main motivation for Chapter~\ref{C2ALPHACHAPET}
has been to understand if, and in which sense, solutions of elliptic equations are ``two derivatives better
than the source term''. While in Chapter~\ref{C2ALPHACHAPET} we focus our
attention to H\"older spaces (which are the closest possible replacement of
classical spaces of continuous functions, in view of
the pathological examples presented in Theorem~\ref{CONSDTYROE-Pse3m4pa235sfbscd3jv}),
our aim is now to consider Lebesgue spaces~$L^p$:
after all, a function is ``nice'' if either it is continuous, possibly together with its derivatives,
or if it has some integrability properties, possibly together with its derivatives, therefore
the topic presented in this chapter can certainly be seen
as complementary to that of Chapter~\ref{C2ALPHACHAPET}. Of course,
in terms of applications, it is absolutely crucial to possess several forms
of a regularity theory, since different occasions (either related to continuity or
integrability properties) naturally occur in several problems of interest.\medskip

The first objective of this chapter is thus to understand
whether or not solutions of~$\Delta u=f$ with~$f\in L^p$
happen to be in the Sobolev space~$W^{2,p}$ (see e.g.~\cite{MR1625845, MR1814364, MR2527916, MR2759829, MR2895178}
and the references therein for a complete introduction to Sobolev spaces).

We already know (recall footnote~\ref{FAV-0kmSRTS9rtboTEFSftS} on page~\pageref{FAV-0kmSRTS9rtboTEFSftS})
that when~$p=\infty$ such a regularity theory does not hold.
The case~$p=1$ is also out of reach, as pointed out by the following result:

\begin{theorem}\label{FAV-0kmSRTS9rtboTEFSftS1} Let~$n\ge2$.
There exist a set of null measure~${\mathcal{Z}}$, $f\in L^1(B_1\setminus{\mathcal{Z}})$ and~$u\in C^2(B_1\setminus{\mathcal{Z}})$
such that~$\Delta u=f$ in~$B_1\setminus{\mathcal{Z}}$ but~$D^2u\not\in L^1({B_1}\setminus{\mathcal{Z}})$.
\end{theorem}

\begin{proof} Given~$x=(x_1,\dots,x_n)\in\R^n$ we use the notation~$\widehat x:=(x_1,x_2)\in\R^2$
and~$\widetilde x:=(x_3,\dots,x_n)\in\R^{n-2}$ (of course, when~$n=2$ we have that~$\widehat x=x$
and the definition of~$\widetilde x$ can be dismissed).
Let
\begin{eqnarray*}&&
{\mathcal{Z}}:=\big\{x=(\widehat x,\widetilde x)\in\R^2\times\R^{n-2} {\mbox{ s.t. }}\widehat x=0\big\},\\
&& u_0(r):=\ln\left(\ln\left(\frac{e}{r}\right)\right),\\&&
u(x):=u_0(|\widehat x|)=\ln\left(\ln\left(\frac{e}{|\widehat x|}\right)\right)\\
{\mbox{and }}&&f(x):=\Delta u(x)=u_0''(|\widehat x|)+\frac{1}{|\widehat x|}u_0'(|\widehat x|)=
-\frac{ 1}{|\widehat x|^2 \,\ln^2\left(\frac{e}{|\widehat x|} \right)}
.
\end{eqnarray*}
Notice that~${\mathcal{Z}}$ has null Lebesgue measure.
Thus, using polar coordinates in~$\R^2$ and the substitution~$t:=\ln\left(\frac{e}{r}\right)$,
and possibly replacing~$C$ from line to line,
\begin{eqnarray*}
\|f\|_{L^1(B_1\setminus{\mathcal{Z}})}&\le&\int_{ 
\{|\widehat x|\in(0,1)\}\times\{|\widetilde x|<1\}
 }
\frac{1}{|\widehat x|^2 \,\ln^2\left(\frac{e}{|\widehat x|} \right)}
\,dx\\
&=&C\,\int_{ \{|\widehat x|\in(0,1)\}
 }
\frac{1}{|\widehat x|^2 \,\ln^2\left(\frac{e}{|\widehat x|} \right)}
\,d\widehat{x}
\\&=&
C\,\int_0^1 
\frac{dr}{r \,\ln^2\left(\frac{e}{r} \right)}
\\&=&
C\,\int_1^{+\infty} 
\frac{dt}{t^2} ,
\end{eqnarray*}
which is finite.

We also remark that
\begin{equation}\label{MsnbdTGDfZjngfGJsJAqKomfAL865}
{\mathcal{W}}:=
\left\{x=(\widehat x,\widetilde x)\in\R^2\times\R^{n-2} {\mbox{ s.t. }}|\widehat x|\in\left(0,\frac12\right)
{\mbox{ and }} |\widetilde x|\in\left[0,\frac12\right)\right\}
\;\subseteq\;
B_1\setminus{\mathcal{Z}}.
\end{equation}
Indeed, if~$x=(\widehat x,\widetilde x)\in{\mathcal{W}}$ then~$\widehat x\ne0$ and
$$ |x|^2=|\widehat x|^2+|\widetilde x|^2<\frac14+\frac14<1,$$
proving~\eqref{MsnbdTGDfZjngfGJsJAqKomfAL865}.

Additionally,
$$ \partial_{12}u(x)=\frac{ x_1x_2}{|\widehat x|^4\ln\left(\frac{e}{|\widehat x|}\right)}\left(2-\frac1{\ln\left(\frac{e}{|\widehat x|}\right)}\right).$$
Since, if~$x\in B_1\setminus{\mathcal{Z}}$ then~$ \ln\left(\frac{e}{|\widehat x|}\right)\ge\ln e=1$, we infer that
$$ |\partial_{12}u(x)|\ge\frac{ |x_1x_2|}{|\widehat x|^4\ln\left(\frac{e}{|\widehat x|}\right)}.$$
For this reason, making use of~\eqref{MsnbdTGDfZjngfGJsJAqKomfAL865} we deduce that
\begin{equation*}\begin{split}
\|D^2 u\|_{L^1(B_1\setminus{\mathcal{Z}})}\,&\ge\,
\int_{ {\mathcal{W}} }\frac{ |x_1x_2|}{|\widehat x|^4\ln\left(\frac{e}{|\widehat x|}\right)}\,dx\\
&=\,C
\int_{ \{|\widehat x|\in(0,1/2)\} }\frac{ |x_1x_2|}{|\widehat x|^4\ln\left(\frac{e}{|\widehat x|}\right)}\,d\widehat{x}\\&=\,C\,
\iint_{ (0,2\pi)\times(0,1/2)}
\frac{ |\sin\vartheta\,\cos\vartheta|}{r\,\ln\left(\frac{e}{r}\right)}\,d\vartheta\,dr\\
&=\,C\,
\int_{0}^{1/2}
\frac{ dr}{r\,\ln\left(\frac{e}{r}\right)}\\&=\,C\,
\int_{\ln(2e)}^{+\infty}
\frac{ dt}{t}  ,
\end{split}\end{equation*}
which is infinite.
\end{proof}

\begin{figure}
  \centering
  \includegraphics[height=5.9cm]{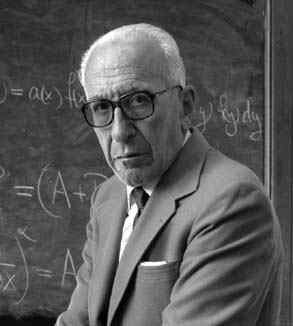}$\quad$  \includegraphics[height=5.9cm]{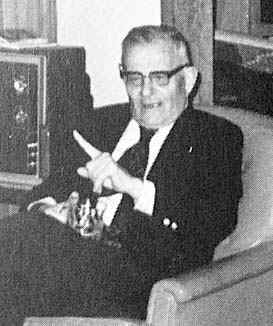}
 \caption{\sl Alberto Calder\'on and Antoni Zygmund (photos by Paul Halmos from MacTutor History of Mathematics Archive,
 licensed under the Creative Commons Attribution-Share Alike 4.0 International license).}\label{SPCAZIGM}
\end{figure}

In the forthcoming pages we will see that the cases~$p=1$ and~$p=\infty$ are the only
exceptional cases for a regularity theory in Lebesgue spaces. That is, roughly
speaking, solutions of~$\Delta u=f$ with~$f\in L^p$ do
happen to be in the Sobolev space~$W^{2,p}$ when~$p\in(1,+\infty)$.
The techniques used for this regularity theory were chiefly introduced by
Alberto Pedro Calder\'on and Antoni Zygmund~\cite{MR52553, MR84633}.
The concepts developed in this setting are very deep and fascinating, but, as it often happens
with truly beautiful ideas, not completely easy to digest
at a first glance (not by chance Calder\'on and Zygmund are actually the founders\footnote{Here is a nice anecdote about the initial meeting between Calder\'{o}n and Zygmund, see~\cite{MR2435330}. In 1948,
during a scientific visit to Buenos Aires, Zygmund delivered a two-month course based on one of his books. The lectures were attended by many young Argentine mathematicians, including, of course, Alberto Calder\'{o}n (who had
graduated in~1947 with a degree in civil engineering, pursued, instead of a degree in mathematics, to follow his father's wishes).

Each of the attendees had to present a topic from the text. Zygmund appeared to be increasingly agitated by Calder\'{o}n's exposition, till he interrupted the speaker abruptly to ask
where he had read the material that he was presenting.

Calder\'{o}n replied that this material was certainly coming from Zygmund's book, but Zygmund vehemently informed the audience that this material was not coming from his book.

After the lecture,  Zygmund took Calder\'{o}n aside to further investigate the matter. Calder\'{o}n finally confessed that he did try to read the material from the book, but after the first couple of lines, instead of turning the page, he had figured out how to develop the arguments by himself in a new and original way.

Obviously, Zygmund immediately recognized Calder\'{o}n's uncommon mathematical skills and invited him to Chicago to study with him. And this was the birth of one of the most successful collaborations in the recent history of mathematics.

See Figure~\ref{SPCAZIGM} for pictures of Calder\'{o}n and Zygmund. See also Section~\ref{Whenyouareamathematician} for Calder\'{o}n's prominent contribution \label{Whenyouareamathematician2} to inverse problems.}
of the so-called
``Chicago School of hard Analysis''). 
We now dive into some details of this construction;
for additional readings on the Calder\'{o}n-Zygmund theory
see also~\cite{MR0290095, MR1104656, MR1640159, MR1814364, MR1996773, MR2435520, MR2508404, MR3099262} and the references therein.

\section{Potential theory and Calder\'{o}n-Zygmund estimates for the Laplace operator}

With respect to the H\"older spaces setting of Chapter~\ref{C2ALPHACHAPET},
the regularity theory in Lebesgue spaces presents an obvious 
initial difficulty, since functions in~$W^{2,p}$ are not necessarily twice differentiable\index{Calder\'{o}n-Zygmund estimates|(}
and therefore one cannot give a pointwise meaning to the equation involved
and should instead rely on weak formulations.
However, to circumvent this technical hurdle, our strategy will be to start working with smooth
objects, obtaining a priori estimates (recall footnote~\ref{FOOAMPRIUoRAnABBA-GRA}
on page~\pageref{FOOAMPRIUoRAnABBA-GRA}), with bounds in Sobolev and Lebesgue spaces
that are independent on the smoothness assumed on the solution.

As already pointed out on page~\pageref{S-PANEW-TS0jNAIJSkJA9kBBEkS},
a cunning manoeuvre to deal with regularity theory is
to try to reduce the problem to the case of harmonic functions
by subtracting to the given solution
the Newtonian potential~$-\Gamma*f$ studied in Proposition~\ref{ESI:NICE:FB}.
The core of the theory then consists in obtaining suitable bounds on
the derivatives of the Newtonian potential
(specifically, Section~\ref{KPJMDarJMSiaCLiMMSciAmoenbvD8Sijf-21} focused
on such an analysis in H\"older spaces and we now focus on the case of Lebesgue spaces).
\medskip

As a preliminary observation,
we point out that the regularity theory in~$L^2$ is somewhat special with respect to the
general case of~$L^p$ with~$p\in(1,+\infty)$, since when~$p=2$
one can simply perform an integration by parts, according to the following remarks. First of all,
if~$R>0$, $x\in \R^n\setminus B_{2R}$ and~$f\in C^{0,1}_0(B_R)$ then
\begin{equation}\label{MSrruj0GAnTBGSnak-1}\begin{split}&
\left|\int_{\R^n}\partial_{ij}\Gamma(x-y)\big(f(x)-f(y)\big)\,dy\right|=
\left|\int_{\R^n}\partial_{ij}\Gamma(x-y)\,f(y)\,dy\right|\\ &\qquad\le
C\|f\|_{L^\infty(B_R)}\int_{B_R}\frac{dy}{|x-y|^n}\le\frac{CR^n\|f\|_{L^\infty(B_R)}}{|x|^n},\end{split}
\end{equation}
up to renaming~$C$ at each stage of the computation.

Hence, by~\eqref{POTNDEWSECDERIVA}, if~$v:=-\Gamma*f$ is the Newtonian potential of~$f$ then, for all~$x\in \R^n\setminus B_{2R}$,
$$ |\partial_{ij}v(x)|\le \left|\int_{\R^n}\partial_{ij}\Gamma(x-y)\big(f(x)-f(y)\big)\,dy\right|+\frac{1}n\,|f(x)|\le
\frac{CR^n\|f\|_{L^\infty(B_R)}}{|x|^n}$$
and therefore
\begin{equation*}
\begin{split}&
\int_{\R^n}|D^2v(x)|^2\,dx=\lim_{M\to+\infty}\sum_{i,j=1}^n\int_{B_M}|\partial_{ij}v(x)|^2\,dx\\&\quad=
\lim_{M\to+\infty}\sum_{i,j=1}^n\int_{B_M}\Big[ \div\big(\partial_{ij}v(x) \partial_jv(x)\,e_i\big)-
\partial_{iij}v(x) \partial_jv(x)\Big]\,dx\\&\quad=
\lim_{M\to+\infty}\sum_{i,j=1}^n\left[ \int_{\partial B_M} \partial_{ij}v(x) \partial_jv(x)\,e_i\cdot\nu(x)\,d{\mathcal{H}}^{n-1}_x-\int_{B_M}\partial_{iij}v(x) \partial_jv(x)\,dx
\right]\\&\quad=
\lim_{M\to+\infty}\sum_{i,j=1}^n\left[ {\mathcal{H}}^{n-1}(\partial B_M)\;O\left(\frac{1}{M^n}\right)
-\int_{B_M}\big[\div\big(\partial_{ii}v(x) \partial_jv(x)\,e_j\big)-\partial_{ii}v(x) \partial_{jj}v(x)
\big]\,dx
\right]\\&\quad=
\lim_{M\to+\infty}\sum_{i,j=1}^n\left[ O\left(\frac{1}{M}\right)
-\int_{\partial B_M} \partial_{ii}v(x) \partial_jv(x)\,e_j\cdot\nu(x)\,d{\mathcal{H}}^{n-1}_x
+\int_{B_M} \partial_{ii}v(x) \partial_{jj}v(x)\,dx
\right]\\&\quad=
\lim_{M\to+\infty}\sum_{i,j=1}^n\left[ O\left(\frac{1}{M}\right)
+\int_{B_M} \partial_{ii}v(x) \partial_{jj}v(x)\,dx
\right]\\&\quad=
\sum_{i,j=1}^n\int_{\R^n} \partial_{ii}v(x) \partial_{jj}v(x)\,dx
\\&\quad= \int_{\R^n} |\Delta v(x)|^2\,dx.
\end{split}
\end{equation*}
This and Proposition~\ref{ESI:NICE:FB} give that, for all~$f\in C^{0,1}_0(\R^n)$,
\begin{equation}\label{o2ur9032hie-L2theorMS-yYNS-1ik842}
\|D^2(\Gamma*f)\|_{L^2(\R^n)}^2=
\int_{\R^n}|D^2(\Gamma*f)(x)|^2\,dx=\int_{\R^n} |f(x)|^2\,dx
=\|f\|_{L^2(\R^n)}^2,\end{equation}
that is the Newtonian potential\index{Newtonian potential} is ``two derivative better than~$f$
in the sense of~$L^2$''.
\medskip

To develop instead the more complex regularity theory in~$L^p$ for all~$p\in(1,+\infty)$,
in view of Proposition~\ref{SI:IN:SE:DR}, it is convenient to study
the property of a singular integral kernel~$K:\R^n\setminus\{0\}\to\R$
satisfying~\eqref{KERNE-SINGO-001},
\eqref{KERNE-SINGO-003} and (recalling~\eqref{KERNE-SINGO-002})
\begin{equation}\label{ABBE9-NOBA}
\int_{B_\rho(x)\setminus B_{\delta}(x)}K(x-y)\,dy=0\qquad{\mbox{for all~$x\in\R^n$ and~$\rho>\delta>0$,}}
\end{equation}
by defining
\begin{equation}\label{ojndPP08uj3rHOmaT7u0-n0epsjf} Tf(x):=\int_{\R^n} K(x-y)\,\big(f(x)-f(y)\big)\,dy,\end{equation}
where the principal value notation in~\eqref{KS-TSGDBfojeISLmadne-ojdfnLS-KMDMDND9uehn}
is implicitly understood.

In this framework, one of the cornerstones of the theory is to establish the following result
for singular integrals:

\begin{theorem}\label{T:CAZY-PIVO}
Let~$p\in(1,+\infty)$ and~$f\in C^{0,1}_0(\R^n)$. Then, there exists a positive constant~$C$, depending only on~$n$, $p$ and~$K$,
such that
$$ \|Tf\|_{L^p(\R^n)}\le C\,\|f\|_{L^p(\R^n)}.$$
\end{theorem}

To address this result, it is opportune to underline the fact that the operator in~\eqref{ojndPP08uj3rHOmaT7u0-n0epsjf}
is ``of order zero''. As a matter of fact, in the literature there are several possible
definitions of ``order of an operator'', but what we mean here is just that
the operator in~\eqref{ojndPP08uj3rHOmaT7u0-n0epsjf} behaves ``neutrally'' with respect to scaling.
More explicitly, if~$f_r(x):=f\left(\frac{x}r\right)$, we (informally) say that an operator~$S$ has order~$\gamma\in\R$ if
$$ Sf_r(x)=r^{-\gamma}\, Sf\left(\frac{x}r\right).$$
That is, in this setting, operators of order~$\gamma$ pick up a factor~$r^{-\gamma}$ under a dilation.
Typically, the sign of~$\gamma$ has a significant impact on the possible regularizing
effects of a given operator. Roughly speaking, when~$\gamma>0$,
the picture becomes ``more irregular at small scales'', suggesting a worsening of the regularity
by the operator. 

For instance, if~$Sf(x):=\partial_1 f(x)$, we have that~$Sf_r(x)=\frac{1}r\partial_1f\left(\frac{x}r\right)=
r^{-1}Sf\left(\frac{x}r\right)$, hence in this case~$\gamma=1$. In general,
taking~$m$ derivatives of a function produces an operator of order~$m\in\N$. Differential operators
confirm the idea that operators with order~$\gamma>0$ worsen the regularity of the function to
which they apply (by definition, the derivative of a function is less regular than the function itself,
having ``one derivative less'').

Another instructive example is given by integral operators: for instance, if~$f:\R\to\R$ and 
$$Sf(x):=\int_0^x f(y)\,dy,$$
then $$Sf_r(x)=\int_0^x f\left(\frac{y}r\right)\,dy=r\int_0^{x/r} f(t)\,dt=rSf\left(\frac{y}r\right),$$
showing that in this case~$\gamma=-1$. Thus, integral operators confirm the intuition
that their action ``improves the regularity'' of the function to which the operators apply 
(by definition, the primitive of a function is more regular than the function itself,
having ``one derivative more'').

As an example of operators of order~$\gamma<0$, one can consider also
$$Sf(x):=\int_{\R^n}|y|^{-\gamma-n}f(x-y)\,dy $$
(the Newtonian potential\index{Newtonian potential} is of this form when~$n\ge3$, with~$\gamma=-2$): in this situation,
$$ Sf_r(x)=\int_{\R^n}|y|^{-\gamma-n}f\left(\frac{x-y}r\right)\,dy=
r^{-\gamma}\int_{\R^n}|z|^{-\gamma-n}f\left(\frac{x}r-z\right)\,dz
= r^{-\gamma} Sf\left(\frac{x}r\right).$$
With respect to the above terminology, the operator in~\eqref{ojndPP08uj3rHOmaT7u0-n0epsjf} is of order zero,
since it commutes with dilations thanks to~\eqref{KERNE-SINGO-001}: indeed,
\begin{equation}
\label{98yiFUIS8GSC6lvEVwds-INVAJSDIAL}\begin{split}
Tf_r(x)\,&=\,\int_{\R^n} K(x-y)\,\left(f\left(\frac{x}r\right)-f\left(\frac{y}r\right)\right)\,dy\\&=\,
\int_{\R^n} g\left(\frac{x-y}{|x-y|}\right)\,\left(f\left(\frac{x}r\right)-f\left(\frac{y}r\right)\right)\,\frac{dy}{|x-y|^n}\\&=\,
r^n\int_{\R^n} g\left(\frac{x-rz}{|x-rz|}\right)\,\left(f\left(\frac{x}r\right)-f(z)\right)\,\frac{dz}{|x-rz|^n}\\&=\,
\int_{\R^n} g\left(\frac{\frac{x}r-z}{\left|\frac{x}r-z\right|}\right)\,\left(f\left(\frac{x}r\right)-f(z)\right)\,\frac{dz}{\left|\frac{x}r-z\right|^n}
\\& =\,
\int_{\R^n} K\left(\frac{x}r-z\right)\,\left(f\left(\frac{x}r\right)-f(z)\right)\,dz\\&=\,Tf\left(\frac{x}r\right).
\end{split}\end{equation}
Another classical example of operator of order zero is given by
the Hilbert Transform\index{Hilbert Transform} of a function~$f :\R\to\R$ defined by
$$ Hf(x):=\frac1\pi\int_\R\frac{f(x)-f(y)}{x-y}\,dy.$$
The Hilbert Transform naturally appears in
fractional calculus and complex analysis
(since, on the boundary of the complex halfplane, two harmonic conjugated functions of a holomorphic function
are related via the Hilbert Transform, see e.g.~\cite[Example~1.10]{CDV-MDG-19} for further details on this).
To confirm that the Hilbert Transform has order zero, one computes that
$$ Hf_r(x)=\frac1\pi\int_\R\frac{f\left(\frac{x}r\right)-f\left(\frac{y}r\right)}{x-y}\,dy
=\frac{r}\pi\int_\R\frac{f\left(\frac{x}r\right)-f(z)}{x-rz}\,dz
=\frac{1}\pi\int_\R\frac{f\left(\frac{x}r\right)-f(z)}{\frac{x}r-z}\,dz=Hf\left(\frac{x}r\right).$$
Interestingly, the Hilbert Transform is a ``combination of a discrete differential and an integral operator'',
in the sense that the ``incremental quotient'' in its integrand acts, roughly speaking, with the order of one derivative
which is compensated by the integration operation.\medskip

Having developed some familiarity with the concept of order of an operator and with its
important connection with the regularity issues, we now aim at introducing a notion of Fourier Transform
for the operator~$T$ in~\eqref{ojndPP08uj3rHOmaT7u0-n0epsjf}. This is not completely straightforward,
since~$K\not\in L^1(\R^n)$ and therefore some care is needed in applying Fourier Transforms to the convolutions
with~$K$. The result that we need in this context is the following:

\begin{lemma} \label{MIJNDdefvPOSujmTDZIOPA3}
Let~$f\in C^{0,1}_0(\R^n)$. Then,
\begin{equation}\label{MSrruj0GAnTBGSnak-12}
Tf\in L^2(\R^n).\end{equation}
Moreover, there exists a function~$\widetilde K\in L^\infty(\R^n)$ such that
\begin{eqnarray}
\label{MIJNDdefvPOSujmTDZIOPA1}&&\widetilde K(t\xi)=\widetilde K(\xi)\qquad{\mbox{for all~$\xi\in\R^n$ and~$t\in(0,+\infty)$}}\\
{\mbox{and }}&&\label{MIJNDdefvPOSujmTDZIOPA2}
\widehat{Tf}(\xi)=\widetilde K(\xi)\,\widehat f(\xi).
\end{eqnarray}
\end{lemma}

Here above, we have denoted by~$\widehat f$ the Fourier Transform of~$f$.
In this spirit, the claim in~\eqref{MIJNDdefvPOSujmTDZIOPA2} identifies a Fourier multiplier
for the operator~$T$ (interestingly, Lemma~\ref{MIJNDdefvPOSujmTDZIOPA3}
also says that this multiplier is bounded and positively homogeneous of degree zero).

In our framework, the interest of Lemma~\ref{MIJNDdefvPOSujmTDZIOPA3}
is that it will allow us to exploit Fourier methods to
establish a theory in~$L^2(\R^n)$ for singular integral operators
(this will be accomplished in the forthcoming equation~\eqref{L2THESINGOPE},
which can be seen as a general counterpart of the estimate on the Newtonian
potential obtained in~\eqref{o2ur9032hie-L2theorMS-yYNS-1ik842} via the Divergence Theorem).

\begin{proof}[Proof of Lemma~\ref{MIJNDdefvPOSujmTDZIOPA3}]
To prove~\eqref{MSrruj0GAnTBGSnak-12} we proceed as in~\eqref{MSrruj0GAnTBGSnak-1}.
That is, we take~$R>0$ sufficiently large such that the support of~$f$ is contained in~$B_R$.
Thus, if~$x\in \R^n\setminus B_{2R}$ then
\begin{eqnarray*}&& \left|\int_{\R^n} K(x-y)\big(f(x)-f(y)\big)\,dy\right|=
\left|\int_{\R^n}K(x-y)\,f(y)\,dy\right|\\ &&\qquad\le
C\|f\|_{L^\infty(B_R)}\int_{B_R}\frac{dy}{|x-y|^n}\le\frac{CR^n\,\|f\|_{L^\infty(B_R)}}{|x|^n},\end{eqnarray*}
thanks to~\eqref{KERNE-SINGO-001}.

If instead~$x\in B_{2R}$ we use~\eqref{KERNE-SINGO-001} and~\eqref{ABBE9-NOBA} to see that
\begin{eqnarray*} &&\left|\int_{\R^n} K(x-y)\big(f(x)-f(y)\big)\,dy\right|\\&\le&
\left|\int_{B_{4R}(x)} K(x-y)\big(f(x)-f(y)\big)\,dy\right|+
\left|\int_{\R^n\setminus B_{4R}(x)} K(x-y)\,f(x)\,dy\right|
\\& \le&
\|f\|_{C^{0,1}(\R^n)} \int_{B_{4R}(x)}|K(x-y)|\,|x-y|\,dy+0
\\ &\le&
C \|f\|_{C^{0,1}(\R^n)} \int_{B_{4R}(x)}\frac{dy}{|x-y|^{n-1}}\,dy\\ &\le&
C R\,\|f\|_{C^{0,1}(\R^n)}.
\end{eqnarray*}
{F}rom these observations we infer that
\begin{eqnarray*}
\|Tf\|^2_{L^2(\R^n)}\le C\left[ R^2\,\|f\|_{C^{0,1}(\R^n)}^2\,|B_{2R}|+R^{2n}\,\|f\|_{L^\infty(B_R)}^2
\int_{\R^n\setminus B_{2R}} \frac{dx}{|x|^{2n}}\right]<+\infty,
\end{eqnarray*}
which establishes~\eqref{MSrruj0GAnTBGSnak-12}.

Owing to~\eqref{MSrruj0GAnTBGSnak-12}, we can therefore consider the Fourier Transform of~$Tf$
in~$L^2(\R^n)$. To this end, we let~$\e\in(0,1)$, to be taken arbitrarily small in what follows, and define
\begin{equation}\label{TAGLIOKETR}
K_\e(x):=\chi_{B_{1/\e}\setminus B_\e}(x)\,K(x).\end{equation}
Let also
\begin{equation}\label{TAGLIOKETR0} T_\e f(x):=\int_{\R^n} K_\e(x-y)\big(f(x)-f(y)\big)\,dy.\end{equation}
We claim that
\begin{equation}\label{MSrruj0GAnTBGSnak-123}
\lim_{\e\searrow0}\|T f-T_\e f\|_{L^2(\R^n)}=0.
\end{equation}
To prove this, we recall~\eqref{ABBE9-NOBA} and observe that
\begin{equation*}
\begin{split}&
|Tf(x)-T_\e f(x)|\\&\qquad=\left|\int_{\R^n} \big(K(x-y)-K_\e(x-y)\big)\big(f(x)-f(y)\big)\,dy\right|\\
&\qquad= \left|\int_{B_{\e}(x)\cup(\R^n\setminus B_{1/\e}(x))}  K(x-y)\big(f(x)-f(y)\big)\,dy\right|\\&\qquad\le 
\left|\int_{B_{\e}(x)}  K(x-y)\big(f(x)-f(y)\big)\,dy\right|+
\left|\int_{\R^n\setminus B_{1/\e}(x)}  K(x-y)\,f(y) \,dy\right|+0\\&\qquad\le 
C\,\|f\|_{C^{0,1}(\R^n)}\,\chi_{B_{R+1}}(x)\int_{B_\e(x)}  \frac{dy}{|x-y|^{n-1}}+C\,\|f\|_{L^\infty(\R^n)}\int_{B_R\setminus B_{1/\e}(x)}  \frac{dy}{|x-y|^n}\\&\qquad\le 
C\e\,\|f\|_{C^{0,1}(\R^n)}\,\chi_{B_{R+1}}(x)+C\,\|f\|_{L^\infty(\R^n)}\int_{B_R\setminus B_{1/\e}(x)}  \frac{dy}{|x-y|^n}.
\end{split}\end{equation*}
We also observe that if~$B_R\setminus B_{1/\e}(x)\ne\varnothing$ then there exists~$z\in B_R$ with~$|x-z|\ge\frac1\e$,
and therefore~$|x|\ge|x-z|-|z|\ge\frac1\e-R\ge\frac{1}{2\e}$, as long as~$\e$ is small enough. As a result,
\begin{equation*} \int_{B_R\setminus B_{1/\e}(x)}  \frac{dy}{|x-y|^n}\le
\chi_{\R^n\setminus B_{1/(2\e)}}(x)\,\int_{B_R}  \frac{dy}{|x-y|^n}\le
\frac{C R^n}{|x|^n}\,\chi_{\R^n\setminus B_{1/(2\e)}}(x).
\end{equation*}
Using these remarks, we arrive at
\begin{eqnarray*}
\|T f-T_\e f\|_{L^2(\R^n)}^2
&\le& C\,\|f\|_{C^{0,1}(\R^n)}^2 \left[ \int_{\R^n}
\left(\e^2\,\chi_{B_{R+1}}(x)+
\frac{R^{2n}}{|x|^{2n}}\,\chi_{\R^n\setminus B_{1/(2\e)}}(x)\right)\,dx\right]\\
&\le& C\|f\|_{C^{0,1}(\R^n)}^2\left[\e^2\,(R+1)^n+\e^nR^{2n}\right]
\end{eqnarray*}
and this proves~\eqref{MSrruj0GAnTBGSnak-123}.

Now we take into account the Fourier Transform of~$T_\e f$ (that we denote by either~$\widehat{T_\e f}$
or~${\mathcal{F}}( {T_\e f})$.
For this, we note that
\begin{eqnarray*}&& T_\e f(x)=\int_{B_{1/\e}(x)\setminus B_\e(x)} K(x-y)\big(f(x)-f(y)\big)\,dy
=-\int_{B_{1/\e}(x)\setminus B_\e(x)} K(x-y)\,f(y)\,dy\\&&\qquad\qquad\qquad=-\int_{\R^n} K_\e(x-y)\,f(y)\,dy=-K_\e*f(x)
,\end{eqnarray*}
thanks to~\eqref{ABBE9-NOBA}, whence
\begin{equation}\label{MSrruj0GAnTBGSnak-1234}
\widehat{T_\e f}(\xi)\,=\,-\widehat{K_\e}(\xi)\,\widehat f(\xi).
\end{equation}
We let~$\xi=|\xi| \omega$ for some~$\omega\in\partial B_1$ and we have that
\begin{equation}\label{ABBE9-N2r432-42OBA-2}\begin{split}
\widehat{K_\e}(\xi)\,&=\,\int_{\R^n} K_\e(x)\,e^{-2\pi i x\cdot\xi}\,dx\\&
=\,\int_{B_{1/\e}\setminus B_\e} K(x)\,e^{-2\pi i |\xi| x\cdot \omega}\,dx\\&
=\,\frac1{|\xi|^n}\int_{B_{|\xi|/\e}\setminus B_{\e|\xi|}} K\left(\frac{y}{|\xi|}\right)\,e^{-2\pi i y\cdot \omega}\,dy\\&
=\, \int_{B_{|\xi|/\e}\setminus B_{\e|\xi|}} K(y)\,e^{-2\pi i y\cdot \omega}\,dy.
\end{split}\end{equation}
Now, for all~$\xi\in\R^n\setminus\{0\}$ we define
\begin{equation}\label{VFKSMSMLOMSSTENDMDSTR8JS}\begin{split}& \widetilde K(\xi)\,:=\, -
\int_{\R^n\setminus B_{1}} K(y)\,e^{-2\pi i y\cdot \omega}\,dy
-
\int_{B_{1}} K(y)\,(e^{-2\pi i y\cdot \omega}-1)\,dy
,\end{split}\end{equation}
where~$\omega=\frac\xi{|\xi|}$.

We observe that~$\widetilde K$ is a bounded function. Indeed,
using the fact that
\begin{equation}\label{OIMTSICONDPVEM08i0u3jt8jPmnoh8CFr2j9e}
\omega\cdot\nabla_y (e^{-2\pi i y\cdot \omega})=-2\pi i e^{-2\pi i y\cdot \omega}\end{equation}
and recalling~\eqref{KERNE-SINGO-001} and~\eqref{KERNE-SINGO-003},
we see that
\begin{eqnarray*}&&
\left|\int_{\R^n\setminus B_{1}} K(y)\,e^{-2\pi i y\cdot \omega}\,dy\right|
= C \left|\int_{\R^n\setminus B_{1}} K(y)\omega\cdot\nabla_y (e^{-2\pi i y\cdot \omega})
\,dy\right|\\&&\qquad=
C\left|\int_{\R^n\setminus B_{1}} \Big[\div\Big(K(y)\,e^{-2\pi i y\cdot \omega}\,\omega\Big)
-\nabla K(y)\cdot\omega \,e^{-2\pi i y\cdot \omega}
\Big]\,dy\right| \\&&\qquad \le
\lim_{\varrho\to+\infty} C\left|
\int_{B_{\varrho}\setminus B_{1}} \Big[\div\Big(K(y)\,e^{-2\pi i y\cdot \omega}\,\omega\Big)\,dy\right|
+C\int_{\R^n\setminus B_{1}} \frac{dy}{|y|^{n+1}}\\&&\qquad\le
\lim_{\varrho\to+\infty} C\left|
\int_{\partial(B_{\varrho}\setminus B_{1})} K(y)\,e^{-2\pi i y\cdot \omega}\,\omega\cdot\nu(y)
\,d{\mathcal{H}}^{n-1}_y\right| +C\\&&\qquad \le
\lim_{\varrho\to+\infty} C\left|\int_{\partial B_\varrho}\frac{d{\mathcal{H}}^{n-1}_y}{|y|^{n}}+
\int_{\partial B_{1}} \frac{d{\mathcal{H}}^{n-1}_y}{|y|^{n}}
\right| +C\\&&\qquad \le
\lim_{\varrho\to+\infty} \frac{C}{\varrho} +C\\&&\qquad \le C,
\end{eqnarray*}
up to relabeling~$C>0$ line after line. Consequently, using again~\eqref{KERNE-SINGO-001},
\begin{equation}\label{OJSMNAMNAJDJAGIAGGIAMDMADKJM9DMASbd}
|\widetilde K(\xi)|\le \left|\int_{\R^n\setminus B_{1}} K(y)\,e^{-2\pi i y\cdot \omega}\,dy\right|+
C\int_{B_{1}} \frac{|e^{-2\pi i y\cdot \omega}-1|}{|y|^n}\,dy\le
C+
C\int_{B_{1}} \frac{dy}{|y|^{n-1}}\le C,
\end{equation}
up to renaming~$C$ at each stage of the calculation, as desired.

Also, the homogeneity claimed in~\eqref{MIJNDdefvPOSujmTDZIOPA1} holds true, since~$\widetilde K(\xi)$
only depends on~$\omega=\frac\xi{|\xi|}$.

We claim that for all~$\xi\in\R^n\setminus\{0\}$
\begin{equation}\label{OJSMNAMNAJDJAGIAGGIAMDMADKJM9DMASb}
\lim_{\e\searrow0}\widehat{K_\e}(\xi) =-\widetilde K(\xi).
\end{equation}
To check this, we pick~$\xi\in\R^n\setminus\{0\}$ and we take~$\e$ small enough such that~$1\in\left(\e|\xi|,
\,\frac{|\xi|}\e\right)$. Hence, we exploit~\eqref{ABBE9-NOBA}
and~\eqref{ABBE9-N2r432-42OBA-2} to see that
\begin{eqnarray*}
\widehat{K_\e}(\xi)&=& \int_{B_{|\xi|/\e}\setminus B_{\e|\xi|}} K(y)\,e^{-2\pi i y\cdot \omega}\,dy\\
&=&\int_{B_{|\xi|/\e}\setminus B_{1}} K(y)\,e^{-2\pi i y\cdot \omega}\,dy
+
\int_{B_{1}\setminus B_{\e|\xi|}} K(y)\,e^{-2\pi i y\cdot \omega}\,dy\\
&=&\int_{B_{|\xi|/\e}\setminus B_1} K(y)\,e^{-2\pi i y\cdot \omega}\,dy
+
\int_{B_{1}\setminus B_{\e|\xi|}} K(y)\,(e^{-2\pi i y\cdot \omega}-1)\,dy.
\end{eqnarray*}
{F}rom this, \eqref{VFKSMSMLOMSSTENDMDSTR8JS} and~\eqref{OIMTSICONDPVEM08i0u3jt8jPmnoh8CFr2j9e},
we arrive at
\begin{eqnarray*}
|\widehat{K_\e}(\xi) +\widetilde K(\xi)|&\le&\left|
\int_{\R^n\setminus B_{|\xi|/\e}} K(y)\,e^{-2\pi i y\cdot \omega}\,dy
\right|+\left|
\int_{B_{\e|\xi|}} K(y)\,(e^{-2\pi i y\cdot \omega}-1)\,dy\right|
\\&\le&
C\,\left(
\left|
\int_{\R^n\setminus B_{|\xi|/\e}} K(y)\,\omega\cdot\nabla_y (e^{-2\pi i y\cdot \omega})\,dy
\right|
+\int_{B_{\e|\xi|}} \frac{dy}{|y|^{n-1}}\right)\\&\le&
C\,\left(
\left|
\int_{\R^n\setminus B_{|\xi|/\e}} \Big[ \div\Big(K(y)\,e^{-2\pi i y\cdot \omega}\,\omega\Big)
-\nabla K(y)\cdot\omega
\,e^{-2\pi i y\cdot \omega}\Big]\,dy
\right|
+ \e|\xi| \right)\\&\le&\lim_{\varrho\to+\infty}
C\,\left(
\left|
\int_{B_\varrho\setminus B_{|\xi|/\e}} \div\Big(K(y)\,e^{-2\pi i y\cdot \omega}\,\omega\Big)\,dy
\right|+
\int_{\R^n\setminus B_{|\xi|/\e}}\frac{dy}{|y|^{n+1}}
+ \e|\xi| \right)\\&\le&
\lim_{\varrho\to+\infty}
C\,\left(
\left|
\int_{\partial(B_\varrho\setminus B_{|\xi|/\e})} K(y)\,e^{-2\pi i y\cdot \omega}\,\omega\cdot\nu(y)\,d{\mathcal{H}}_y^{n-1}
\right|+
\int_{\R^n\setminus B_{|\xi|/\e}}\frac{dy}{|y|^{n+1}}
+ \e|\xi| \right)\\&\le&
\lim_{\varrho\to+\infty}
C\,\left(
\int_{\partial B_\varrho } \frac{d{\mathcal{H}}_y^{n-1}}{|y|^n}
+\int_{\partial B_{|\xi|/\e}} \frac{d{\mathcal{H}}_y^{n-1}}{|y|^n}
+\frac{\e}{|\xi|}
+ \e|\xi| \right)\\&\le&
\lim_{\varrho\to+\infty}
C\,\left(
\frac{1}\varrho
+\frac{\e}{|\xi|}
+ \e|\xi| \right)\\&=&C\,\left(
\frac{\e}{|\xi|}
+ \e|\xi| \right),
\end{eqnarray*}
whence~\eqref{OJSMNAMNAJDJAGIAGGIAMDMADKJM9DMASb} plainly follows, as desired.

Now we claim that
\begin{equation}\label{OJSMNAMNAJDJAGIAGGIAMDMADKJM9DMASbc}
\sup_{\e\in(0,1)}\|\widehat{K_\e}\|_{L^\infty(\R^n)}<+\infty.
\end{equation}
To this end, we repeatedly exploit~\eqref{ABBE9-NOBA} and~\eqref{ABBE9-N2r432-42OBA-2}
(together with~\eqref{KERNE-SINGO-001} and~\eqref{KERNE-SINGO-003})
by distinguishing three cases, depending on whether~$1>\frac{|\xi|}\e$,
or~$\frac{|\xi|}\e\ge1\ge\e|\xi|$, or~$\e|\xi|>1$.

Let us first assume that~$1>\frac{|\xi|}\e$. Then,
\begin{eqnarray*}
|\widehat{K_\e}(\xi)|=\left|
\int_{B_{|\xi|/\e}\setminus B_{\e|\xi|}} K(y)\,(e^{-2\pi i y\cdot \omega}-1)\,dy\right|
\le C\,\int_{B_{|\xi|/\e}} \frac{dy}{|y|^{n-1}}
\le \frac{C\,|\xi|}\e \le C,
\end{eqnarray*}
giving~\eqref{OJSMNAMNAJDJAGIAGGIAMDMADKJM9DMASbc} in this case.

Let us now suppose that~$\frac{|\xi|}\e\ge1\ge\e|\xi|$.
In this situation, it is appropriate to use again~\eqref{OIMTSICONDPVEM08i0u3jt8jPmnoh8CFr2j9e}
to see that
\begin{eqnarray*}
|\widehat{K_\e}(\xi)|&=&\left|
\int_{B_{|\xi|/\e}\setminus B_{1}} K(y)\,e^{-2\pi i y\cdot \omega}\,dy+
\int_{B_{1}\setminus B_{\e|\xi|}} K(y)\,e^{-2\pi i y\cdot \omega}\,dy\right|\\&\le&
C\left|\int_{B_{|\xi|/\e}\setminus B_{1}} K(y)\,\omega\cdot\nabla_y (e^{-2\pi i y\cdot \omega})\,dy\right|+\left|
\int_{B_{1}\setminus B_{\e|\xi|}} K(y)\,(e^{-2\pi i y\cdot \omega}-1)\,dy\right|\\&\le&
C\left|\int_{B_{|\xi|/\e}\setminus B_{1}} \Big[\div\Big(K(y)\,e^{-2\pi i y\cdot \omega}\omega\Big)
-\nabla K(y)\cdot\omega\,e^{-2\pi i y\cdot \omega}\Big]\,dy\right|+C
\int_{B_{1}} \frac{dy}{|y|^{n-1}}
\\&\le&
C\left|\int_{\partial(B_{|\xi|/\e}\setminus B_{1})} K(y)\,e^{-2\pi i y\cdot \omega}\omega\cdot\nu(y)\,d{\mathcal{H}}^{n-1}_y\right|+
C\int_{\R^n\setminus B_1}\frac{dy}{|y|^{n+1}}+C
\\&\le&
C \int_{\partial B_{|\xi|/\e}} \frac{d{\mathcal{H}}^{n-1}_y}{|y|^n}+
C \int_{\partial B_{1}} \frac{d{\mathcal{H}}^{n-1}_y}{|y|^n}+C\\&\le& \frac{C\e}{|\xi|}+C\\&\le&C,
\end{eqnarray*}
which establishes~\eqref{OJSMNAMNAJDJAGIAGGIAMDMADKJM9DMASbc} in this case.

It remains to consider the case~$\e|\xi|>1$: in this situation, we make use again of~\eqref{OIMTSICONDPVEM08i0u3jt8jPmnoh8CFr2j9e}
and we have that
\begin{eqnarray*}
|\widehat{K_\e}(\xi)|&=&\left|\int_{B_{|\xi|/\e}\setminus B_{\e|\xi|}} K(y)\,e^{-2\pi i y\cdot \omega}\,dy\right|\\&\le&
C\left|\int_{B_{|\xi|/\e}\setminus B_{\e|\xi|}} K(y)\,\omega\cdot\nabla_y (e^{-2\pi i y\cdot \omega})\,dy\right|\\&=&
C\left|\int_{B_{|\xi|/\e}\setminus B_{\e|\xi|}} \Big[ \div\Big( K(y)\,e^{-2\pi i y\cdot\omega}\,\omega
\Big)-\nabla K(y)\cdot\omega\,e^{-2\pi i y\cdot \omega}\Big]\,dy\right|
\\&\le&
C\left|\int_{\partial(B_{|\xi|/\e}\setminus B_{\e|\xi|})} K(y)\,e^{-2\pi i y\cdot\omega}\,\omega\cdot\nu(y)\,d{\mathcal{H}}^{n-1}_y\right|+
C\int_{\R^n\setminus B_{\e|\xi|}} \frac{dy}{|y|^{n+1}}\\&\le&
C \int_{\partial B_{|\xi|/\e}} \frac{d{\mathcal{H}}^{n-1}_y}{|y|^n}+
C \int_{\partial B_{\e|\xi|}} \frac{d{\mathcal{H}}^{n-1}_y}{|y|^n}+\frac{C}{\e|\xi|}\\&\le&
\frac{C\e}{|\xi|}
+\frac{C}{\e|\xi|}\\&\le&\frac{C}{\e|\xi|}\\&\le&C,
\end{eqnarray*}
which completes the proof of~\eqref{OJSMNAMNAJDJAGIAGGIAMDMADKJM9DMASbc}.

Now we show that
\begin{equation}\label{OJSMNAMNAJDJAGIAGGIAMDMADKJM9DMAS}
\lim_{\e\searrow0}\|
\widehat{K_\e}\,\widehat f
+\widetilde K\,\widehat f\|_{L^2(\R^n)}=0.
\end{equation}
For this, we utilize~\eqref{OJSMNAMNAJDJAGIAGGIAMDMADKJM9DMASbd}
and~\eqref{OJSMNAMNAJDJAGIAGGIAMDMADKJM9DMASbc}
and note that
$$ |\widehat{K_\e}(\xi)\,\widehat f(\xi)
+\widetilde K(\xi)\,\widehat f(\xi)|^2=|\widehat{K_\e}(\xi)+\widetilde K(\xi)|^2\,|\widehat f(\xi)|^2\le
C\,|\widehat f(\xi)|^2
$$
and the latter function belongs to~$L^1(\R^n)$, since~$f\in C^{0,1}_0(\R^n)$.
{F}rom this observation, \eqref{OJSMNAMNAJDJAGIAGGIAMDMADKJM9DMASb} and the
Dominated Convergence Theorem we obtain~\eqref{OJSMNAMNAJDJAGIAGGIAMDMADKJM9DMAS},
as desired.

As a consequence, by~\eqref{MSrruj0GAnTBGSnak-123}, \eqref{MSrruj0GAnTBGSnak-1234}
and~\eqref{OJSMNAMNAJDJAGIAGGIAMDMADKJM9DMAS},
and making use of the Plancherel Theorem,
\begin{eqnarray*}&&
\|\widehat{Tf}-\widetilde K\,\widehat f\|_{L^2(\R^n)}\le\lim_{\e\searrow0}
\|\widehat{Tf}-\widehat{T_\e f}\|_{L^2(\R^n)}+\|\widehat{T_\e f}-\widetilde K\,\widehat f\|_{L^2(\R^n)}\\&&\qquad
=\lim_{\e\searrow0}
\|{Tf}-{T_\e f}\|_{L^2(\R^n)}+\|
\widehat{K_\e}\,\widehat f
+\widetilde K\,\widehat f\|_{L^2(\R^n)}
=0,
\end{eqnarray*}
which establishes~\eqref{MIJNDdefvPOSujmTDZIOPA2}.
\end{proof}

Another crucial ingredients for the proof of Theorem~\ref{T:CAZY-PIVO} consists in a deeper understanding of
the $L^1$ case. On the one hand, 
the counterexample presented in Theorem~\ref{FAV-0kmSRTS9rtboTEFSftS1} prevents
us from a ``completely satisfactory theory in~$L^1$''. On the other hand,
it still leaves the possibility open for a ``weak theory in~$L^1$'', in the following sense.\index{weak theory in~$L^1$}

If a function~$g$ belongs to~$L^1(\R^n)$ then, for all~$\lambda>0$,
$$ \|g\|_{L^1(\R^n)}\ge\int_{\{ |g|\ge\lambda\}}|g(x)|\,dx\ge\int_{\{ |g|\ge\lambda\}}\lambda\,dx=
\lambda\,\big| \{ |g|\ge\lambda\}\big|$$
and therefore
\begin{equation}\label{35MS9uyhdfiowghf93wutygfuywegf7wftqfwguweghSERdfgbLFVBdfNOMDUFJNTInNOA}
\big| \{ |g|\ge\lambda\}\big|\le\frac{\|g\|_{L^1(\R^n)}}{\lambda}.
\end{equation}
This information is often referred to with the name of Chebyshev's Inequality.
When~$g:=Tf$, we know from the
counterexample in Theorem~\ref{FAV-0kmSRTS9rtboTEFSftS1} that~\eqref{35MS9uyhdfiowghf93wutygfuywegf7wftqfwguweghSERdfgbLFVBdfNOMDUFJNTInNOA}
may well become void, simply because~$Tf$ may not belong to~$L^1(\R^n)$
and the right hand side of~\eqref{35MS9uyhdfiowghf93wutygfuywegf7wftqfwguweghSERdfgbLFVBdfNOMDUFJNTInNOA}
may accordingly become infinite. However, fortunately this is not the end of the story,
since when~$g:=Tf$ a variant of~\eqref{35MS9uyhdfiowghf93wutygfuywegf7wftqfwguweghSERdfgbLFVBdfNOMDUFJNTInNOA} holds true,
simply by replacing the norm on the right hand side with~$\|f\|_{L^1(\R^n)}$, up to a constant.
The precise result in this setting goes as follows:

\begin{theorem}\label{hoiFC-AJS-OJMS-OKSMBDkldmrmt-lfhqwhig}
Let~$f\in C^{0,1}_0(\R^n)$. Then, for every~$\lambda>0$,
$$ \big| \{ |Tf|\ge\lambda\}\big|\le\frac{C\,\|f\|_{L^1(\R^n)}}{\lambda},$$
for a suitable constant~$C>0$, depending only on~$n$ and~$K$.
\end{theorem}

This result is very deep and to prove it in the most efficient way we can start with
some simple observations which help us to ``normalize the picture'' of
Theorem~\ref{hoiFC-AJS-OJMS-OKSMBDkldmrmt-lfhqwhig}. First of all, by possibly replacing~$f$
with~$\frac{f}{\lambda}$, we see that
\begin{equation}\label{hoiFC-AJS-OJMS-OKSMBDkldmrmt-lfhqwhig-NORM1}
{\mbox{it is enough to prove Theorem~\ref{hoiFC-AJS-OJMS-OKSMBDkldmrmt-lfhqwhig} when~$\lambda=1$.}}
\end{equation}
In addition,
\begin{equation}\label{hoiFC-AJS-OJMS-OKSMBDkldmrmt-lfhqwhig-NORM2}
{\mbox{it is enough to prove Theorem~\ref{hoiFC-AJS-OJMS-OKSMBDkldmrmt-lfhqwhig} when~$\|f\|_{L^1(\R^n)}=1$.}}
\end{equation}
To check this, note that if~$\|f\|_{L^1(\R^n)}=0$ then~$f$ and~$Tf$ vanish identically,
hence the result in Theorem~\ref{hoiFC-AJS-OJMS-OKSMBDkldmrmt-lfhqwhig} is obvious.
Thus, we can suppose that~$\|f\|_{L^1(\R^n)}>0$ and let~$\widetilde f(x):=f\big( \|f\|_{L^1(\R^n)}^{1/n}x\big)$.
Then, using that~$T$ commutes with dilations (recall~\eqref{98yiFUIS8GSC6lvEVwds-INVAJSDIAL}),
$$ T\widetilde f(x)= Tf\big( \|f\|_{L^1(\R^n)}^{1/n}x\big).$$
Also,
$$ \|\widetilde f\|_{L^1(\R^n)}=\int_{\R^n}
\big| f\big( \|f\|_{L^1(\R^n)}^{1/n}x\big)\big| \,dx=\frac{1}{\|f\|_{L^1(\R^n)}}
\int_{\R^n}|f(y)|\,dy=1.$$
Hence, if we knew that Theorem~\ref{hoiFC-AJS-OJMS-OKSMBDkldmrmt-lfhqwhig}
held true for functions with unit norm in~$L^1(\R^n)$, we would conclude that
\begin{eqnarray*}
\frac{\big| \{ |Tf|\ge\lambda\}\big|\}}{\| f\|_{L^1(\R^n)}}
=\big| \{ |T\widetilde f|\ge\lambda\}\big|
\le\frac{C\,\|\widetilde f\|_{L^1(\R^n)}}{\lambda}=\frac{C}{\lambda},
\end{eqnarray*}
that is, Theorem~\ref{hoiFC-AJS-OJMS-OKSMBDkldmrmt-lfhqwhig} would hold for~$f$ as well,
proving~\eqref{hoiFC-AJS-OJMS-OKSMBDkldmrmt-lfhqwhig-NORM2}.

Owing to~\eqref{hoiFC-AJS-OJMS-OKSMBDkldmrmt-lfhqwhig-NORM1} and~\eqref{hoiFC-AJS-OJMS-OKSMBDkldmrmt-lfhqwhig-NORM2}, to prove Theorem~\ref{hoiFC-AJS-OJMS-OKSMBDkldmrmt-lfhqwhig} it suffices to
focus on the renormalized situation in which
\begin{equation}\label{hoiFC-AJS-OJMS-OKSMBDkldmrmt-lfhqwhig-NORM3}\begin{split}&
{\mbox{we assume that~$f\in C^{0,1}_0(\R^n)$ with~$\|f\|_{L^1(\R^n)}=1$}}\\&{\mbox{and we aim at showing that
}}\big| \{ |Tf|\ge1\}\big|\le C.\end{split}
\end{equation}
In this setting, another useful observation is that for ``nice functions''
one can rely on the~$L^2$ theory.
That is, the observation in~\eqref{o2ur9032hie-L2theorMS-yYNS-1ik842} can be recast into
a general framework of Fourier Transforms as follows: first of all,
by Lemma~\ref{MIJNDdefvPOSujmTDZIOPA3} and Plancherel Theorem,
\begin{equation}\label{L2THESINGOPE} \| {Tf}\|_{L^2(\R^n)}=\| \widehat{Tf}\|_{L^2(\R^n)}=\|\widetilde K\,\widehat f\|_{L^2(\R^n)}
\le\|\widetilde K\|_{L^\infty(\R^n)}\,\|\widehat f\|_{L^2(\R^n)}\le C\| f\|_{L^2(\R^n)}.\end{equation}
As a result,
$$ C\| f\|_{L^2(\R^n)}^2\ge\int_{\{ |Tf|\ge1\}}|Tf(x)|^2\,dx\ge\big| \{ |Tf|\ge1\}\big|.$$
In particular, if~$\|f\|_{L^\infty(\R^n)}\le1$ and the setting of~\eqref{hoiFC-AJS-OJMS-OKSMBDkldmrmt-lfhqwhig-NORM3} holds true, then
$$ \| f\|_{L^2(\R^n)}^2=\int_{\R^n}|f(x)|^2\,dx\le\int_{\R^n}|f(x)|\,dx=\|f\|_{L^1(\R^n)}=1$$
and consequently~$\big| \{ |Tf|\ge1\}\big|\le C$.

Summarizing, in the setting of~\eqref{hoiFC-AJS-OJMS-OKSMBDkldmrmt-lfhqwhig-NORM3},
\begin{equation}\label{hoiFC-AJS-OJMS-OKSMBDkldmrmt-lfhqwhig-NORM4}\begin{split}&
{\mbox{if additionally~$\|f\|_{L^\infty(\R^n)}\le1$ then }}\big| \{ |Tf|\ge1\}\big|\le C\end{split}
\end{equation}
and we are done.\medskip

\begin{figure}
  \centering
  \includegraphics[width=.6\linewidth]{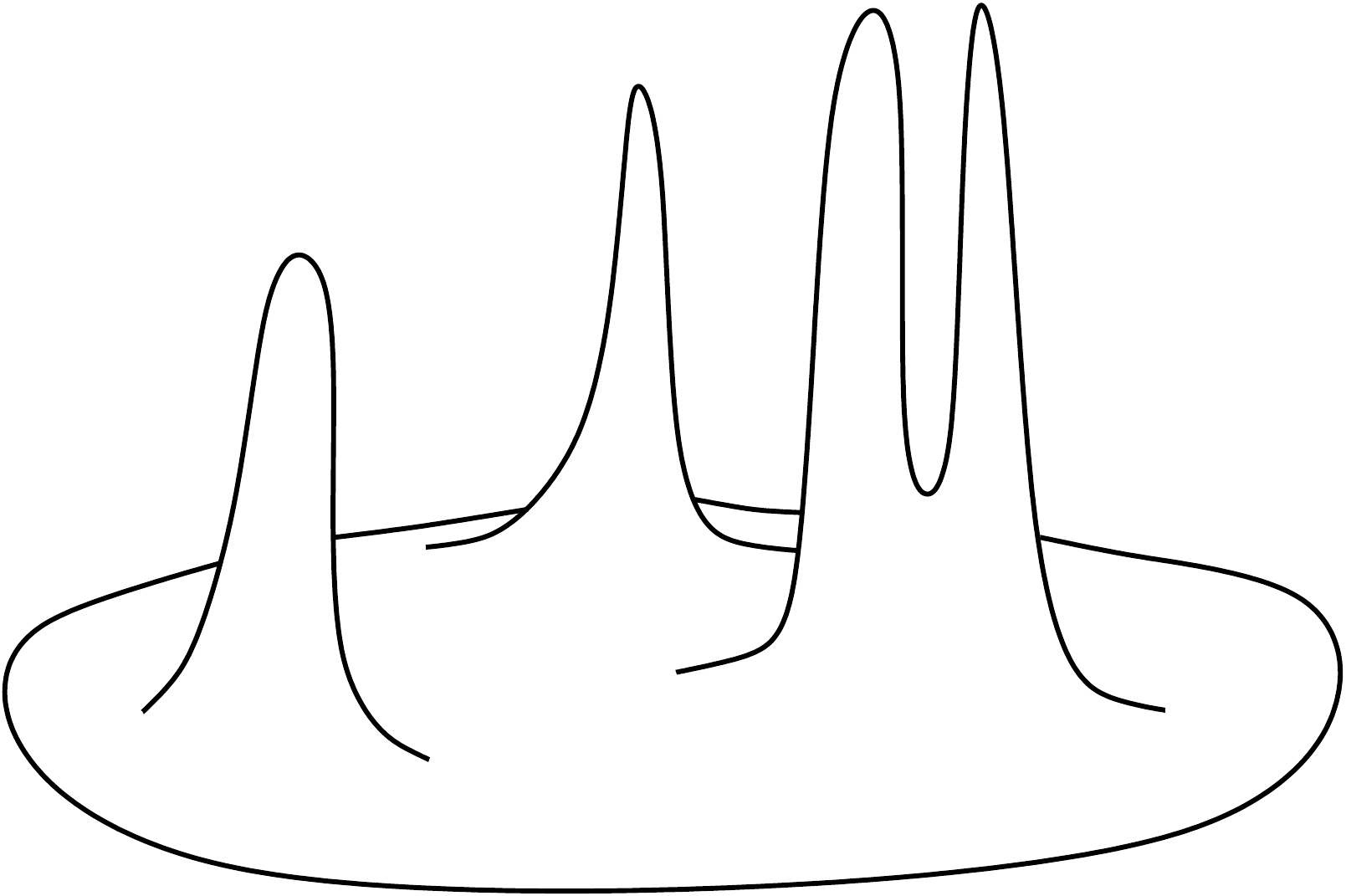}
 \caption{\sl A function developing spikes.}\label{SPIKERTERRAeFI}
\end{figure}

Thus, in light of~\eqref{hoiFC-AJS-OJMS-OKSMBDkldmrmt-lfhqwhig-NORM3}
and~\eqref{hoiFC-AJS-OJMS-OKSMBDkldmrmt-lfhqwhig-NORM4},
to prove Theorem~\ref{hoiFC-AJS-OJMS-OKSMBDkldmrmt-lfhqwhig} it only remains to understand the case \label{CI-CUSNDBEDOCOMDDOSStni0PMS-oop}
of smooth and compactly supported functions with unit mass in~$L^1(\R^n)$ whose absolute value is not bounded by~$1$:
roughly speaking, this is the case in which the function ``develops spikes'', see Figure~\ref{SPIKERTERRAeFI}.
To develop some intuition handy to address this case, let us consider the ``worse case scenario'' of a function developing spikes,
that is the one of the sum of Dirac Delta Functions (well, well, strictly speaking
the Dirac Delta Function is not a function! but so what, at this point we are just trying
to understand what additional difficulties spikes introduce to the setting in~\eqref{hoiFC-AJS-OJMS-OKSMBDkldmrmt-lfhqwhig-NORM3} and to unveil a useful device to cope with them). For instance, if~$f$ were~$\delta_{x_0}$, then,
at a formal level, recalling~\eqref{KERNE-SINGO-001} and~\eqref{ABBE9-NOBA} we reduce~\eqref{ojndPP08uj3rHOmaT7u0-n0epsjf} to
\begin{equation}\label{RAWAMSFA8BSPEKragasNZO04SIMA9ktJA3AMD5}
\begin{split}& Tf(x)=\int_{\R^n} K(x-y)\,\big(\delta_{x_0}(x)-\delta_{x_0}(y)\big)\,dy=
0-\int_{\R^n} K(x-y)\,\delta_{x_0}(y)\,dy\\&\qquad\qquad\qquad\qquad=-K(x-x_0)=-\frac{1}{|x-x_0|^n}\,g\left(\frac{x-x_0}{|x-x_0|}\right).\end{split}\end{equation}
As a consequence,
\begin{eqnarray*} \{ |Tf|\ge1\}&=&\left\{
x\in\R^n {\mbox{ s.t. }} \left|g\left(\frac{x-x_0}{|x-x_0|}\right)\right|\ge|x-x_0|^n
\right\}\\&=&\left\{
y+x_0\in\R^n {\mbox{ s.t. }} \left|g\left(\frac{y}{|y|}\right)\right|\ge|y|^n
\right\}\\&=&\left\{
te+x_0\in\R^n {\mbox{ s.t. $e\in\partial B_1$, $t\ge0$ and }} |g(e)|\ge t^n
\right\}\\&=&\left\{
x+x_0\in\R^n {\mbox{ s.t. }} x\in{\mathcal{B}},
\right\}
\end{eqnarray*}
where
$$ {\mathcal{B}}:=\Big\{
x=te {\mbox{ s.t. $e\in\partial B_1$, $t\ge0$ and }} |g(e)|\ge t^n
\Big\}.$$
For instance, when~$g$ is identically equal to~$1$ we have that~${\mathcal{B}}$ is the closure of the unit ball
and~$\{ |Tf|\ge1\}$ is the closure of the unit ball centered at~$x_0$. So, with an abuse of notation
(that we allow to ourselves, being this discussion already rather heuristic), we consider~${\mathcal{B}}$
to be a ``ball induced by the kernel~$K$'' and we realize that when~$f$ is the Dirac Delta Function centered at~$x_0$
then~$\{ |Tf|\ge1\}$ is the ball induced by the kernel~$K$ centered at~$x_0$.
Since, in polar coordinates, this cases produces
$$ |\{ |Tf|\ge1\}|=|{\mathcal{B}}|=\int_{\partial B_1}\left[
\int_{0}^{|g(e)|^{1/n}} t^{n-1}\,dt 
\right]\,d{\mathcal{H}}_e^{n-1}\le
\int_{\partial B_1}\left[
\int_{0}^{\|g\|_{L^\infty}^{1/n}} t^{n-1}\,dt 
\right]\,d{\mathcal{H}}_e^{n-1}
\le C,
$$
we see that the desired result in~\eqref{hoiFC-AJS-OJMS-OKSMBDkldmrmt-lfhqwhig-NORM3}
would be accomplished in this model case.

\begin{figure}
  \centering
    \includegraphics[width=.55\linewidth]{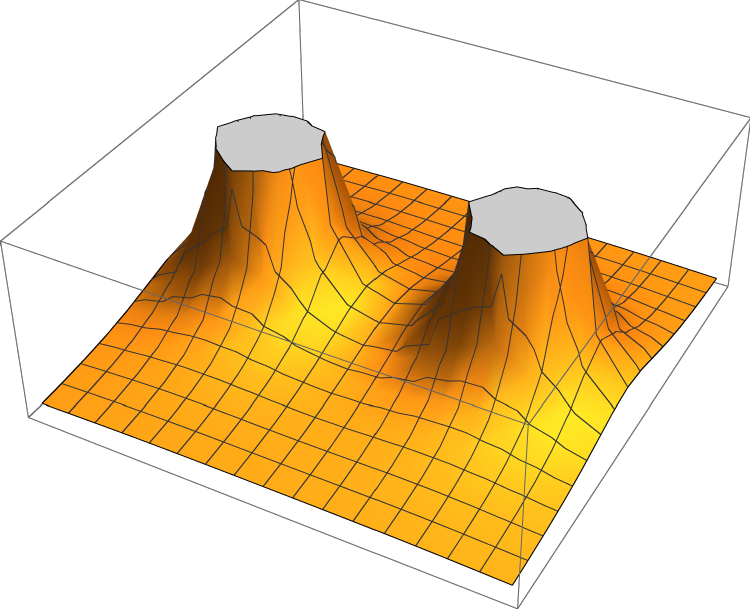}$\quad$
  \includegraphics[width=.3\linewidth]{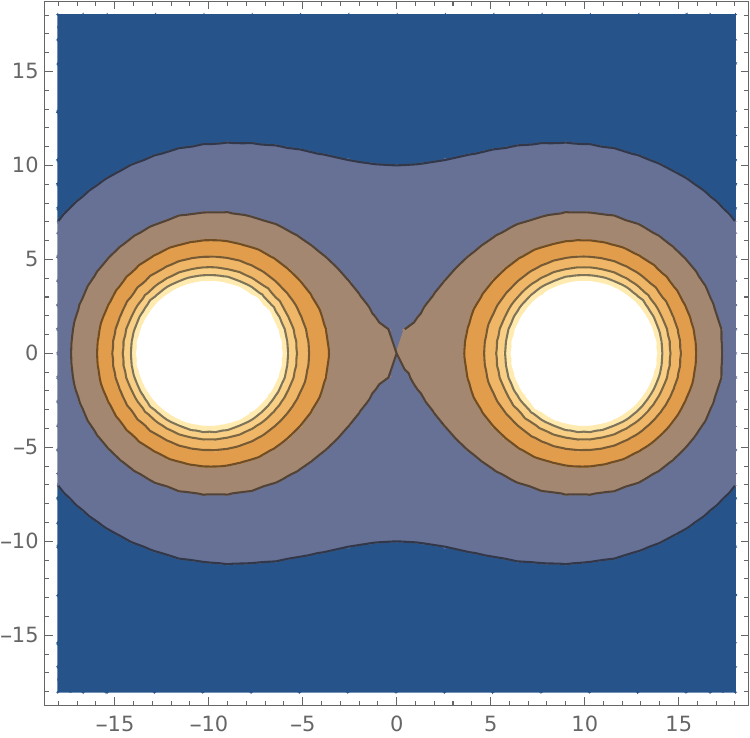}
 \caption{\sl The function in~\eqref{RAWAMSFA8BSPEKragasNZO04SIMA9ktJA3AMD8} and its level sets
 when~$g=1$, $n=2$, $x_0=10e_1$, $y_0=-10e_1$ and the factor~$-\frac12$ is dropped.}\label{DRHSo0mNEdvspp-02iX32h90-SPIKERTERRAeFI-b}
\end{figure}

We have therefore learned that when~$f$ is a single Dirac Delta Function
the effect on~$\{ |Tf|\ge1\}$ is that of producing some kind of ball~${\mathcal{B}}$
and we now consider the next model situation in which~$f$ is the sum of two Dirac Delta Functions
(normalized to maintain total unit mass), say~$\frac{\delta_{x_0}+\delta_{y_0}}2$.
In this case we can imagine that when~$x_0$ and~$y_0$ are extremely close,
the situation pretty much reduces to that of a single Dirac Delta Function. Instead,
when~$x_0$ and~$y_0$ are ``far away'' then one can use~\eqref{RAWAMSFA8BSPEKragasNZO04SIMA9ktJA3AMD5}
and find that
\begin{equation}\label{RAWAMSFA8BSPEKragasNZO04SIMA9ktJA3AMD8} Tf(x)=-\frac12\,\left[
\frac{1}{|x-x_0|^n}\,g\left(\frac{x-x_0}{|x-x_0|}\right)
+\frac{1}{|x-y_0|^n}\,g\left(\frac{x-y_0}{|x-y_0|}\right)\right],\end{equation}
which is basically the sum of two singularities located at a country mile from each other,
see e.g. Figure~\ref{DRHSo0mNEdvspp-02iX32h90-SPIKERTERRAeFI-b} for the case~$g$ identcally equal to~$1$. Hence, for very remote locations of~$x_0$ and~$y_0$,
if one neglects, at least\footnote{We kindly ask the scrupulous reader to hold their fire: we will come back
to the issue created by the tails in a second.}
for the moment, the ``overlapping tails'' of the two functions in the right hand side of~\eqref{RAWAMSFA8BSPEKragasNZO04SIMA9ktJA3AMD8}, the set~$\{ |Tf|\ge1\}$
is essentially made by two disjoint balls induced by the kernels (one centered at~$x_0$ and one centered at~$y_0$)
of the form
$$ {\mathcal{B}}':=\Big\{
x=te {\mbox{ s.t. $e\in\partial B_1$, $t\ge0$ and }} |g(e)|\ge 2t^n
\Big\}.$$
Since, substituting~$\tau:=2^{\frac1n}t$,
$$ |{\mathcal{B}}'|=\int_{\partial B_1}\left[
\int_{0}^{(|g(e)|/2)^{1/n}} t^{n-1}\,dt 
\right]\,d{\mathcal{H}}_e^{n-1}=\frac12\int_{\partial B_1}\left[
\int_{0}^{|g(e)|^{1/n}} \tau^{n-1}\,dt 
\right]\,d{\mathcal{H}}_e^{n-1}=\frac{|{\mathcal{B}}|}2
$$
we find that for the sum of two Dirac Delta Functions with singularities located far away,
the measure of~$\{ |Tf|\ge1\}$ is, as a first approximation, essentially
twice that of the ball~${\mathcal{B}}'$, that is precisely that of~${\mathcal{B}}$.
This observation suggests that also the case of two Dirac Delta Functions
can be reduced to that of a single Dirac Delta Function and is therefore well under control:
however, this heuristic also points out that the way to cope with the case of two (or more)
Dirac Delta Functions is to distinguish whether the singularities are ``close'' (in which case
one has to gather all the singularities together and somewhat replace this cluster with a single
Dirac Delta Function, possibly with an appropriate weight) or ``far away'' (in which case
each produces a measure comparable, give or take, to that produced by a weighted 
Dirac Delta Function).\medskip

At this point, however, the scrupulous reader has certainly spotted a noteworthy weakness
in the above argument: namely, neglecting the interactions of the tails
of remote spikes is clumsily too crude, since the decay of the functions in~\eqref{RAWAMSFA8BSPEKragasNZO04SIMA9ktJA3AMD8} is too slow, and not integrable at infinity.
This is a considerable drawback, since it entails that, in principle, the more spikes the function
develops, the bigger the outcome produced by the tail interactions, which would, again in principle,
prevent us from obtaining the desired result (recall that the constant in~\eqref{hoiFC-AJS-OJMS-OKSMBDkldmrmt-lfhqwhig-NORM3} is allowed to depend only on~$n$ and~$K$, and surely not on the number of spikes
of the given function, nor on their precise location).

But we can still correct the previous, somewhat sloppy, approximation argument
by exploiting again the cancellations of the kernel provided by~\eqref{ABBE9-NOBA}:
to this end, roughly speaking, one has to locally subtract to the spike its own average
and reabsorb this additional average term into a ``universally bounded'' part of the function.
In this way, on the one hand, the spikes minus their averages provide a bunch of localized functions
with zero average, for which the kernel provides an additional cancellation that makes
the tail interactions uniformly integrable; on the other hand, for the universally bounded part of the function
one can exploit the observation in~\eqref{hoiFC-AJS-OJMS-OKSMBDkldmrmt-lfhqwhig-NORM4} (possibly, with~$1$ replaced by a universal constant).
The combination of these facts will establish the desired result in~\eqref{hoiFC-AJS-OJMS-OKSMBDkldmrmt-lfhqwhig-NORM3} after a careful inspection, as we will see in the forthcoming pages.

\begin{figure}
  \centering
  \includegraphics[width=.6\linewidth]{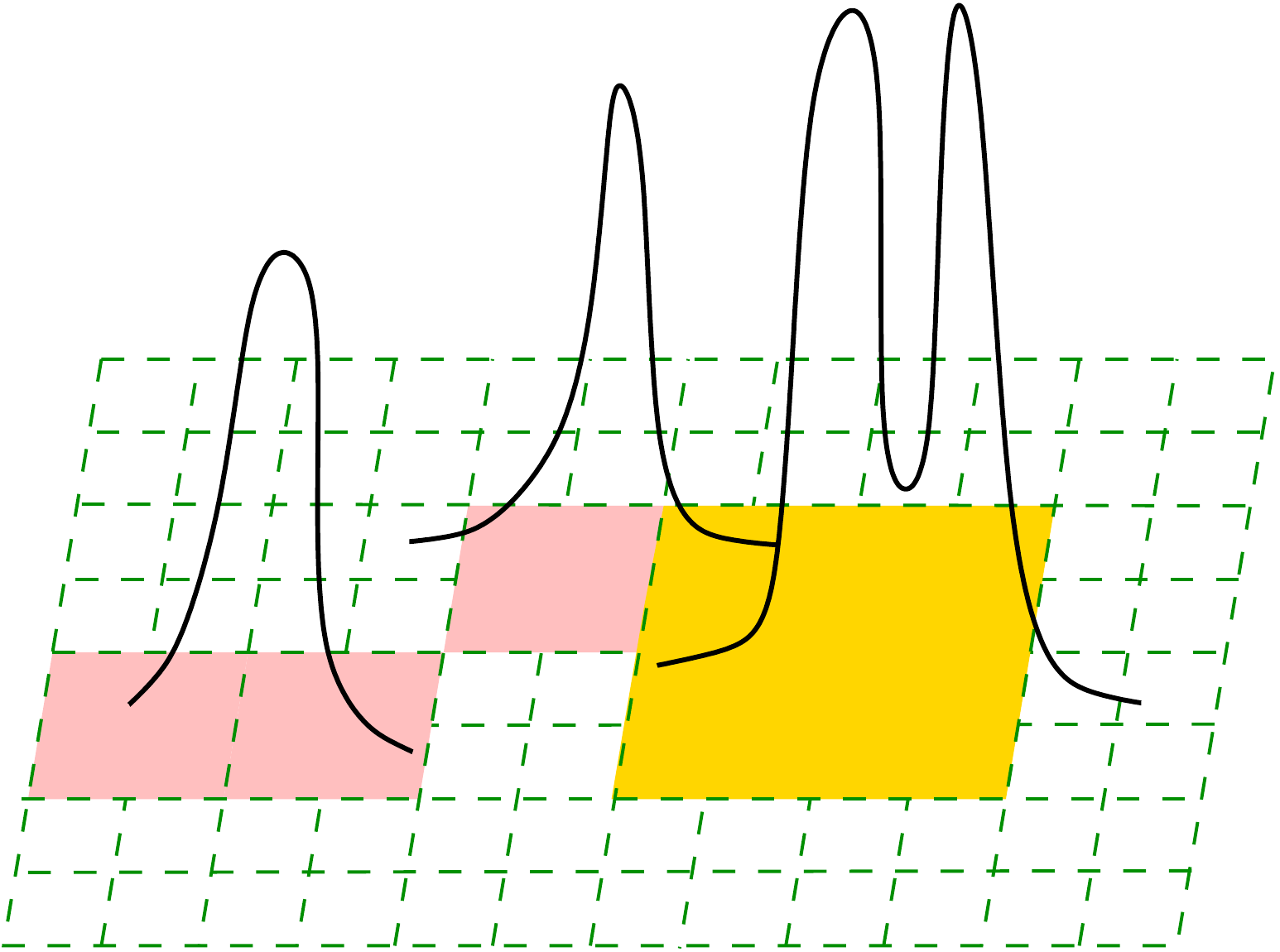}
 \caption{\sl Cube decomposition to detect and gather together the spikes of Figure~\ref{SPIKERTERRAeFI}.}\label{90-SPIKERTERRAeFI-b}
\end{figure}

The Ultimate Question of Life, the Universe, and Everything
is therefore how to make this heuristic argument work and how to make it quantitative
in a sufficiently general setting which includes that in~\eqref{hoiFC-AJS-OJMS-OKSMBDkldmrmt-lfhqwhig-NORM3}.
The Answer to this Question comes from the so-called Calder\'{o}n-Zygmund
cube decomposition\index{Calder\'{o}n-Zygmund
cube decomposition|(}, which is a brilliant way to detect spikes
for a given function with unit mass in~$L^1(\R^n)$.
The idea of this procedure is that one can split the whole space~$\R^n$ into cubes, say of unit side.
In each of these cubes, the average of the function is less than or equal to~$1$ (because the total
mass of the function is equal to~$1$). However, in very small cubes located in the vicinity of the spikes
the average of the function can be very large. Therefore, to locate the position of the possible spikes,
one can split dyadically the original cubes and inspect the average in these new cubes.
If the average of the function is stricly larger than~$1$, we place the cube into a new family (roughly speaking, in this cube,
we have ``almost detected a spike'', or perhaps more than one spikes that are sufficiently close together and that
can thereby be clustered in such a cube); if instead the average in the cube is still less than or equal to~$1$
we keep splitting the cube dyadically, see Figure~\ref{90-SPIKERTERRAeFI-b}.
This protocol allows us to partition~$\R^n$ into two subsets: a subset collecting all the cubes
in the new family (in which the average of the function is stricly larger than~$1$,
though such average cannot be too large, given the fact that the average in the preceding cube
was less than or equal to~$1$) and a subset containing all the other points (and one can check that
the value of the function at these points cannot overcome~$1$, otherwise the point would
fall, sooner or later after dyadic divisions, into one of the cubes of the family).\medskip

To make this argument more transparent and precise, we state and prove the following result:

\begin{theorem}\label{CI-CUSNDBEDOCOMDDOSStni0PMS}
Let~$f\in C^{0,1}_0(\R^n)$ with~$\|f\|_{L^1(\R^n)}=1$. Then, there exists a countable
family of cubes~$\{Q_j\}_{j\in\N}$ with sides parallel to the Cartesian axes such that
\begin{eqnarray}&&\label{kjytr4M542SMIHSRFGBDFRTYUDuMAjdfPERghvLEKESEMFCOS1k-2}
\fint_{Q_j} |f(x)|\,dx\in(1,2^n]\qquad{\mbox{for all }}\,j\in\N,\\&&
\label{kjytr4M542SMIHSRFGBDFRTYUDuMAjdfPERghvLEKESEMFCOS1k-2bis}
\sum_{j\in\N}|Q_j|\le1
\\ {\mbox{and }}&&\label{kjytr4M542SMIHSRFGBDFRTYUDuMAjdfPERghvLEKESEMFCOS1k-3}
|f(x)|\le1\qquad{\mbox{for almost all }}\,x\in\R^n\setminus \left( \bigcup_{j\in\N} Q_j\right)
.\end{eqnarray}
\end{theorem}

\begin{proof}
We start by considering cubes of side~$1$ of the form~$(m_1,m_1+1)\times\dots\times(m_n,m_n+1)$ with~$m_1,\dots,m_n\in\Z$.
The family of all these cubes will be denoted by~${\mathcal{Q}}_0$.
If~$Q$ is any of these cubes, we have that~$|Q|=1$ and
$$ \fint_{Q} |f(x)|\,dx=\int_Q |f(x)|\,dx\le\|f\|_{L^1(\R^n)}=1.$$
We now denote by~${\mathcal{Q}}_1$ the cubes of side~$\frac12$ obtained by dividing into two equal parts
each side of the cubes in the family~${\mathcal{Q}}_0$
(i.e., each cube in~${\mathcal{Q}}_0$ is split into~$2^n$ cubes to produce the family~$
{\mathcal{Q}}_1$; the cubes in~${\mathcal{Q}}_1$ are also congruent, open and disjoint).

We collect from~${\mathcal{Q}}_1$ the cubes~$Q$ (if any) for which
$$ \fint_{Q} |f(x)|\,dx>1$$
and we call~${\mathcal{F}}_1$ the family containing all the cubes of this type, i.e.
$$ {\mathcal{F}}_1:=\left\{
Q\in{\mathcal{Q}}_1{\mbox{ s.t. }}\fint_{Q} |f(x)|\,dx>1
\right\}.$$
Then we keep splitting the remaining cubes in~${\mathcal{Q}}_1\setminus {\mathcal{F}}_1$,
collecting them in a family called~${\mathcal{Q}}_2$,
we identify those cubes~$Q$ in~${\mathcal{Q}}_2$ (if any) for which~$\fint_{Q} |f(x)|\,dx>1$
and we place them into a family~${\mathcal{F}}_2$, then we keep splitting the remaining cubes in~${\mathcal{Q}}_2\setminus {\mathcal{F}}_2$, collecting them in a family called~${\mathcal{Q}}_3$, and so on.

In this way, the cubes in~${\mathcal{Q}}_j$ have side~$\frac1{2^j}$
and the cubes~$Q$ in~${\mathcal{F}}_j$ have the property that
\begin{equation}\label{kjytr4M542SMIHSRFGBDFRTYUDuMAjdfPERghvLEKESEMFCOS1k-1}
\fint_{Q} |f(x)|\,dx>1.
\end{equation}
We define
$${\mathcal{F}}:=\bigcup_{j=1}^{+\infty}{\mathcal{F}}_j$$
and we write~${\mathcal{F}}=\{Q_j\}_{j\in\N}$
for suitable cubes~$Q_j$.

We also point out that if~$Q$ is in~${\mathcal{F}}_j$ then its predecessor was in~${\mathcal{Q}}_{j-1}\setminus
{\mathcal{F}}_{j-1}$ (with the notation that~${\mathcal{F}}_{0}$ is void): that is, 
if~$Q$ is in~${\mathcal{F}}_j$ then there exists~$\widetilde Q\in{\mathcal{Q}}_{j-1}\setminus
{\mathcal{F}}_{j-1}$ such that~$|\widetilde Q|=2^n|Q|$ and
$$ \fint_{\widetilde{Q}} |f(x)|\,dx\le1.$$
Consequently,
$$ \fint_{Q} |f(x)|\,dx=\frac{2^n}{|\widetilde Q|}\int_{Q} |f(x)|\,dx\le 2^n\fint_{\widetilde{Q}} |f(x)|\,dx\le2^n.$$
Combining this and~\eqref{kjytr4M542SMIHSRFGBDFRTYUDuMAjdfPERghvLEKESEMFCOS1k-1}
we obtain~\eqref{kjytr4M542SMIHSRFGBDFRTYUDuMAjdfPERghvLEKESEMFCOS1k-2}.

As a byproduct of~\eqref{kjytr4M542SMIHSRFGBDFRTYUDuMAjdfPERghvLEKESEMFCOS1k-2} we also have that
$$\sum_{j\in\N}|Q_j|\le\sum_{j\in\N}\int_{Q_j}|f(x)|\,dx\le\int_{\R^n}|f(x)|\,dx=\|f\|_{L^1(\R^n)}=1,$$
which proves~\eqref{kjytr4M542SMIHSRFGBDFRTYUDuMAjdfPERghvLEKESEMFCOS1k-2bis}, as desired.

Now we prove~\eqref{kjytr4M542SMIHSRFGBDFRTYUDuMAjdfPERghvLEKESEMFCOS1k-3}.
For this, we pick~$x\in\R^n\setminus \left( \displaystyle\bigcup_{j\in\N} Q_j\right)$
to be a Lebesgue point for~$|f|$: we recall that these points have full measure,
thanks to Lebesgue's Differentiation Theorem,
see e.g.~\cite[Theorem~7.16]{MR3381284}. We can also suppose that the coordinates~$x_1,\dots,x_n$ of~$x$
do not belong to the grid induced by the dyadic cubes (since this is a set of measure zero).
In this situation, since~$x$ lies outside the cubes of the family~${\mathcal{F}}$, we have that for every~$m$
there exists a cube~$K_m$ of side~$\frac{1}{2^m}$ such that~$x\in K_m$ and
$$ \fint_{K_m}|f(y)|\,dy\le1.$$
Note that if~$K^\star_m:=\left(x_1-\frac{2}{2^m},x_1+\frac{2}{2^m}\right)\times\dots\times
\left(x_n-\frac{2}{2^m},x_n+\frac{2}{2^m}\right)$ we have that~$K^\star_m\supseteq K_m$
and~$|K^\star_m|=\frac{4^n}{2^{mn}}={4^n}|K_m|$, therefore~$K_m$ shrinks regularly to~$x$
(in the terminology of~\cite[Theorem 7.16]{MR3381284}).
Accordingly,
$$ |f(x)|=\lim_{m\to+\infty}\fint_{K_m}|f(y)|\,dy\le1,$$
which proves~\eqref{kjytr4M542SMIHSRFGBDFRTYUDuMAjdfPERghvLEKESEMFCOS1k-3}.
\end{proof}

We stress once again that~\eqref{kjytr4M542SMIHSRFGBDFRTYUDuMAjdfPERghvLEKESEMFCOS1k-2} reveals
a special feature of the Calder\'{o}n-Zygmund
cube decomposition since it detects the ``correct scale'' of the cubes to take into account:
if in lieu of this accurate decomposition one considered the simple
decomposition obtained by replacing every single spike with a constant function (say, equal to~$1$)
in a cube,
with the size of the cubes depending on the height of the given spike, it is very likely that these
cubes would end up overlapping between themselves and the values of the function
obtained in this way would run completely out of control. The core of the Calder\'{o}n-Zygmund
cube decomposition, as showcased in~\eqref{kjytr4M542SMIHSRFGBDFRTYUDuMAjdfPERghvLEKESEMFCOS1k-2},
is instead to choose the cubes in such a way that the ``measure of all the spikes inside the cube''
(as quantified by~$\int_{Q_j}|f(x)|\,dx$) is precisely comparable to the measure of the cube itself (that is~$|Q_j|$).\medskip

With the aid of Theorem~\ref{CI-CUSNDBEDOCOMDDOSStni0PMS}, we can thus retake the heuristic
strategy introduced on page~\pageref{CI-CUSNDBEDOCOMDDOSStni0PMS-oop}, and make it work!
The master plan is thus to use a cube decomposition to detect
spikes, clustering together spikes that are sufficiently close to lie in the same cube~$Q_j$.
One then splits the function~$f$ into two functions: on the one side,
a ``bad part'', that will be denoted by~${\mbox{\footnotesize{\calligra{B}}}}\;$, which collects all the spikes~$b_j$ located in all the cubes~$Q_j$,
being careful in subtracting in each of these cubes the average of the corresponding spikes
(so to make it average zero in its support, which is important to exploit kernel cancellations in the tail interaction and deduce that away from~$Q_j$
a favorable estimate holds, see the forthcoming equation~\eqref{XVBSFVAKTHNS955KA2M56I7Z8A-094-6}
and e.g. Figure~2 in~\cite{MR1640159}
for a classical picture about this method); on the other side,
a ``good part'', that will be denoted by~${\mbox{\footnotesize{\calligra{G}}}}\;$, which accounts for the remaining terms
of~$f$ coming from~\eqref{kjytr4M542SMIHSRFGBDFRTYUDuMAjdfPERghvLEKESEMFCOS1k-3},
as well as from the average of the spikes (these contributions will be bounded
uniformly allowing us to employ~\eqref{hoiFC-AJS-OJMS-OKSMBDkldmrmt-lfhqwhig-NORM4} on~${\mbox{\footnotesize{\calligra{G}}}}\;$).\medskip

The mathematical details of this blueprint are the following. 
In the setting of Theorem~\ref{CI-CUSNDBEDOCOMDDOSStni0PMS}, we define
\begin{equation}\label{Tbfrac12}
\begin{split}
&b_j(x):=\begin{dcases}
\displaystyle f(x)-\fint_{Q_j}f(y)\,dy & {\mbox{ if }}x\in Q_j,\\
0&{\mbox{ otherwise,}}
\end{dcases}\\
&{\mbox{\footnotesize{\calligra{B}}}}\;(x):=\sum_{j\in\N}b_j(x)\\{\mbox{and }}\;&
{\mbox{\footnotesize{\calligra{G}}}}\;(x):=f(x)-{\mbox{\footnotesize{\calligra{B}}}}\;(x).
\end{split}\end{equation}
Notice that
\begin{equation}\label{kjytr4M542SMIHSRFGBDFRTYUDuMAjdfPERghvLEKESEMFCOS1k-3987}\fint_{Q_j} b_j(x)\,dx=0\qquad{\mbox{ for all }}\,j\in\N
\end{equation}
and
\begin{equation}\label{kjytr4M542SMIHSRFGBDFRTYUDuMAjdfPERghvLEKESEMFCOS1k-39} {\mbox{\footnotesize{\calligra{G}}}}\;(x)=\begin{dcases}
\displaystyle f(x) & {\mbox{ if }}x\in\R^n\setminus\bigcup_{j\in\N} Q_j,\\
\displaystyle \fint_{Q_j}f(y)\,dy & {\mbox{ if }}x\in Q_j.
\end{dcases}\end{equation}\index{Calder\'{o}n-Zygmund
cube decomposition|)}
With this notation, we have:

\begin{lemma} \label{eto776aggiunto}
We have that
\begin{eqnarray}
\label{XVBSFVAKTHNS955KA2M56I7Z8A-094-1}&&\big|{\mbox{\footnotesize{\calligra{G}}}}\;(x)\big|
\le2^n \qquad{\mbox{ for all }}\,x\in\R^n,\\
\label{XVBSFVAKTHNS955KA2M56I7Z8A-094-2}&&\big\|{\mbox{\footnotesize{\calligra{G}}}}\;\;\big\|_{L^1(\R^n)}\le1\\{\mbox{and }}
\label{XVBSFVAKTHNS955KA2M56I7Z8A-094-3}&&\fint_{Q_j}|b_j(x)|\,dx\le2^{n+1}\qquad{\mbox{ for all }}\,j\in\N.
\end{eqnarray}
Furthermore, if~$Q'_j$ denotes the cube with the same center of~$Q_j$ and side $(3\sqrt{n}+1)$ times the side of~$Q_j$,
\begin{equation}\label{XVBSFVAKTHNS955KA2M56I7Z8A-094-6}
\int_{\R^n\setminus Q_j'}|Tb_j(x)|\,dx\le C\,|Q_j|.
\end{equation}
Additionally,
\begin{equation}\int_{\R^n\setminus \bigcup_{j\in\N}Q_j'}\big|T{\mbox{\footnotesize{\calligra{B}}}}\;(x)\big|\,dx\le C
.\label{XVBSFVAKTHNS955KA2M56I7Z8A-094-5}\end{equation}
Here,~$C$ is a positive constant depending only on~$n$ and~$K$.
\end{lemma}

\begin{proof} If~$x\in Q_j$ we deduce from~\eqref{kjytr4M542SMIHSRFGBDFRTYUDuMAjdfPERghvLEKESEMFCOS1k-2}
and~\eqref{kjytr4M542SMIHSRFGBDFRTYUDuMAjdfPERghvLEKESEMFCOS1k-39} that
$$ \big|{\mbox{\footnotesize{\calligra{G}}}}\;(x)\big|\le\fint_{Q_j} |f(x)|\,dx\le2^n.$$
{F}rom this and~\eqref{kjytr4M542SMIHSRFGBDFRTYUDuMAjdfPERghvLEKESEMFCOS1k-3} we obtain~\eqref{XVBSFVAKTHNS955KA2M56I7Z8A-094-1}.

Moreover, using again~\eqref{kjytr4M542SMIHSRFGBDFRTYUDuMAjdfPERghvLEKESEMFCOS1k-39},
\begin{eqnarray*}
\int_{\R^n}\big|{\mbox{\footnotesize{\calligra{G}}}}\;
(x)\big|\,dx&\le&\int_{\R^n\setminus\bigcup_{j\in\N} Q_j}|f(x)|\,dx+\sum_{j\in\N}
\int_{Q_j}\left|\fint_{Q_j}f(y)\,dy\right|\,dx\\
&\le&
\int_{\R^n\setminus\bigcup_{j\in\N} Q_j}|f(x)|\,dx+
\sum_{j\in\N} \int_{Q_j}|f(y)|\,dy\\
&=&\|f\|_{L^1(\R^n)},
\end{eqnarray*}
which gives~\eqref{XVBSFVAKTHNS955KA2M56I7Z8A-094-2}.

Furthermore, if~$x\in Q_j$,
$$|b_j(x)|\le|f(x)|+\fint_{Q_j}|f(y)|\,dy$$
and therefore, recalling~\eqref{kjytr4M542SMIHSRFGBDFRTYUDuMAjdfPERghvLEKESEMFCOS1k-2},
$$ \fint_{Q_j}|b_j(x)|\,dx
\le
\fint_{Q_j}|f(x)|\,dx
+\fint_{Q_j}\left(
\fint_{Q_j}|f(y)|\,dy\right)\,dx=
2\fint_{Q_j}|f(x)|\,dx\le2^{n+1},$$
establishing~\eqref{XVBSFVAKTHNS955KA2M56I7Z8A-094-3}.

Now we denote by~$p_j$ the center of~$Q_j$ and employ~\eqref{ABBE9-NOBA},
\eqref{ojndPP08uj3rHOmaT7u0-n0epsjf}
and~\eqref{kjytr4M542SMIHSRFGBDFRTYUDuMAjdfPERghvLEKESEMFCOS1k-3987} to see that,
in the principal value sense,
\begin{eqnarray*}&&
Tb_j(x)=\int_{\R^n} K(x-y)\,\big(b_j(x)-b_j(y)\big)\,dy=
-\int_{\R^n} K(x-y)\,b_j(y)\,dy\\&&\qquad\qquad=-\int_{Q_j} K(x-y)\,b_j(y)\,dy=
\int_{Q_j}\big(K(x-p_j)- K(x-y)\big)\,b_j(y)\,dy.
\end{eqnarray*}
Notice that here we have exploited in a crucial way the fact that~$b_j$ averages to zero in~$Q_j$.
This allows us to improve the nonintegrable decay of the kernel by using~\eqref{KERNE-SINGO-003},
according to the following calculation. 
If~$\ell_j $ denotes the length of the side of~$Q_j$, $y\in Q_j$ and~$x\in
\R^n\setminus Q_j'$, then~$|x-p_j|\ge \frac{3\sqrt{n} \ell_j}2$ and consequently,
for all~$t\in(0,1)$,
$$ |x-(1-t)p_j-ty|=|x-p_j+t(p_j-y)|\ge|x-p_j|-|p_j-y|\ge|x-p_j|-\sqrt{n}\ell_j\ge\frac{|x-p_j|}3.
$$
As a result, if~$y\in Q_j$ and~$x\in
\R^n\setminus Q_j'$,
\begin{eqnarray*}&& |K(x-p_j)-K(x-y)|\le\int_0^1 |\nabla K(x-(1-t)p_j-ty)|\,|y-p_j|\,dt
\\&&\qquad\le C\int_0^1 \frac{|y-p_j|}{|x-(1-t)p_j-ty|^{n+1}}\,dt
\le C\ell_j\int_0^1 \frac{dt}{|x-p_j|^{n+1}}\,dt\le\frac{C\ell_j}{|x-p_j|^{n+1}}.
\end{eqnarray*}
We stress that this decay is one power better than the one of the kernel~$K$ and,
in particular, provides an integrable tail. In this way, we find that, for all~$x\in
\R^n\setminus Q_j'$,
\begin{eqnarray*}
|Tb_j(x)|\le \int_{Q_j}\big|K(x-p_j)- K(x-y)\big|\,|b_j(y)|\,dy\le
\frac{C\ell_j}{|x-p_j|^{n+1}} \int_{Q_j}|b_j(y)|\,dy.
\end{eqnarray*}
This and~\eqref{XVBSFVAKTHNS955KA2M56I7Z8A-094-3} give that, for all~$x\in
\R^n\setminus Q_j'$,
\begin{eqnarray*}
|Tb_j(x)|\le \frac{C\ell_j\,|Q_j|}{|x-p_j|^{n+1}}.\end{eqnarray*}
We can now integrate in~$x$ and find that
\begin{eqnarray*}
&&\int_{\R^n\setminus Q_j'}|Tb_j(x)|\,dx\le \int_{\R^n\setminus B_{\ell_j/2}(p_j)}
\frac{C\ell_j\,|Q_j|\,dx}{|x-p_j|^{n+1}}\le
C\ell_j\,|Q_j|\int_{\ell_j/2}^{+\infty}\frac1{\rho^2}\,d\rho\le
C\,|Q_j|
\end{eqnarray*}
and this proves~\eqref{XVBSFVAKTHNS955KA2M56I7Z8A-094-6}.

{F}rom this and~\eqref{kjytr4M542SMIHSRFGBDFRTYUDuMAjdfPERghvLEKESEMFCOS1k-2bis}, we also arrive at
$$ 
\int_{\R^n\setminus \bigcup_{k\in\N}Q_k'}\big|T{\mbox{\footnotesize{\calligra{B}}}}\;(x)\big|\,dx\le
\sum_{j\in\N}\int_{\R^n\setminus \bigcup_{k\in\N}Q_k'}|Tb_j(x)|\,dx\le
\sum_{j\in\N}\int_{\R^n\setminus Q_j'}|Tb_j(x)|\,dx
\le C\,\sum_{j\in\N} |Q_j|\le C,$$
thus establishing~\eqref{XVBSFVAKTHNS955KA2M56I7Z8A-094-5}, as desired.
\end{proof}

With this, we can now complete the proof of Theorem~\ref{hoiFC-AJS-OJMS-OKSMBDkldmrmt-lfhqwhig} in the following way.

\begin{proof}[Proof of Theorem~\ref{hoiFC-AJS-OJMS-OKSMBDkldmrmt-lfhqwhig}]
We work in the renormalized setting provided by~\eqref{hoiFC-AJS-OJMS-OKSMBDkldmrmt-lfhqwhig-NORM3} and we use the notation of Lemma~\ref{eto776aggiunto}.
In light of~\eqref{XVBSFVAKTHNS955KA2M56I7Z8A-094-1} and~\eqref{XVBSFVAKTHNS955KA2M56I7Z8A-094-2},
up to renaming constants, we can exploit~\eqref{hoiFC-AJS-OJMS-OKSMBDkldmrmt-lfhqwhig-NORM4} for the function~${\mbox{\footnotesize{\calligra{G}}}}\;$
and obtain that
\begin{equation}\label{kjytr4M542SMIHSRFGBDFRTYUDuMAjdfPERghvLEKESEMFCOS1k-2bis-0945-XVBSNDICAVBSSOC}
\left| \left\{\big|T{\mbox{\footnotesize{\calligra{G}}}}\;\;\big|\ge\frac12\right\}\right|\le C.
\end{equation}
Now we let
$$ {\mathcal{U}}:=\bigcup_{j\in\N}Q_j'.$$
By~\eqref{kjytr4M542SMIHSRFGBDFRTYUDuMAjdfPERghvLEKESEMFCOS1k-2bis},
\begin{equation}\label{kjytr4M542SMIHSRFGBDFRTYUDuMAjdfPERghvLEKESEMFCOS1k-2bis-0945}
|{\mathcal{U}}|\le\sum_{j\in\N}|Q_j'|\le C
\sum_{j\in\N}|Q_j|\le C.\end{equation}
Also, by~\eqref{XVBSFVAKTHNS955KA2M56I7Z8A-094-5},
$$
C\ge\int_{\R^n\setminus {\mathcal{U}}}\big|T{\mbox{\footnotesize{\calligra{B}}}}\;(x)\big|\,dx\ge
\int_{(\R^n\setminus {\mathcal{U}})\cap\{ \big|T{\mbox{\footnotesize{\calligra{B}}}}\;\;\big|\ge1/2\} }\big|T{\mbox{\footnotesize{\calligra{B}}}}\;(x)\big|\,dx\ge
\frac12\left|
\big(\R^n\setminus {\mathcal{U}}\big)\cap\left\{\big|T{\mbox{\footnotesize{\calligra{B}}}}\;\;\big|\ge\frac12\right\}
\right|.
$$
This and~\eqref{kjytr4M542SMIHSRFGBDFRTYUDuMAjdfPERghvLEKESEMFCOS1k-2bis-0945}
yield that
\begin{equation}\label{kjytr4M542SMIHSRFGBDFRTYUDuMAjdfPERghvLEKESEMFCOS1k-2bis-0945-XVBSNDICAVBSSO}
\left|
\left\{\big|T{\mbox{\footnotesize{\calligra{B}}}}\;\;\big|\ge\frac12\right\}
\right|\le
\left|
\big(\R^n\setminus {\mathcal{U}}\big)\cap\left\{\big|T{\mbox{\footnotesize{\calligra{B}}}}\;\;\big|\ge\frac12\right\}
\right|
+
|{\mathcal{U}}|
\le C.
\end{equation}
We also observe that
\begin{equation}\label{kjytr4M542SMIHSRFGBDFRTYUDuMAjdfPERghvLEKESEMFCOS1k-2bis-0945-XVBSNDICAVBSSO2}
\left\{|Tf|\ge 1\right\}\,
\subseteq\,\left\{\big|T{\mbox{\footnotesize{\calligra{B}}}}\;\;\big|\ge\frac12\right\}
\cup\left\{\big|T{\mbox{\footnotesize{\calligra{G}}}}\;\;\big|\ge\frac12\right\}.
\end{equation}
Indeed, suppose not. Then, there exists~$x\in\left\{|Tf|\ge 1\right\}$
with~$\big|T{\mbox{\footnotesize{\calligra{B}}}}\;(x)\big|<\frac12$ and~$\big|T{\mbox{\footnotesize{\calligra{G}}}}\;(x)\big|<\frac12$.
Consequently, by the definition of~${\mbox{\footnotesize{\calligra{G}}}}\;$ in~\eqref{Tbfrac12},
$$ 1\le|Tf(x)|=\big|T{\mbox{\footnotesize{\calligra{B}}}}\;(x)+T{\mbox{\footnotesize{\calligra{G}}}}\;(x)\big|\le\big|T{\mbox{\footnotesize{\calligra{B}}}}\;(x)\big|+\big|T{\mbox{\footnotesize{\calligra{G}}}}\;(x)\big|<\frac12+\frac12=1.$$
This is a contradiction and therefore the claim in~\eqref{kjytr4M542SMIHSRFGBDFRTYUDuMAjdfPERghvLEKESEMFCOS1k-2bis-0945-XVBSNDICAVBSSO2} holds true.

Now, by~\eqref{kjytr4M542SMIHSRFGBDFRTYUDuMAjdfPERghvLEKESEMFCOS1k-2bis-0945-XVBSNDICAVBSSOC},
\eqref{kjytr4M542SMIHSRFGBDFRTYUDuMAjdfPERghvLEKESEMFCOS1k-2bis-0945-XVBSNDICAVBSSO} and~\eqref{kjytr4M542SMIHSRFGBDFRTYUDuMAjdfPERghvLEKESEMFCOS1k-2bis-0945-XVBSNDICAVBSSO2},
\begin{eqnarray*}
\left|\left\{|Tf|\ge 1\right\}\right|\le\left|\left\{\big|T{\mbox{\footnotesize{\calligra{B}}}}\;\;\big|\ge\frac12\right\}\right|+
\left|\left\{\big|T{\mbox{\footnotesize{\calligra{G}}}}\;\;\big|\ge\frac12\right\}\right|\le C,
\end{eqnarray*}
up to renaming~$C$ as usual. This shows that the claim in~\eqref{hoiFC-AJS-OJMS-OKSMBDkldmrmt-lfhqwhig-NORM3}
holds true and the proof of Theorem~\ref{hoiFC-AJS-OJMS-OKSMBDkldmrmt-lfhqwhig}
is thereby complete.
\end{proof}

It is instructive to observe that the estimate in Theorem~\ref{hoiFC-AJS-OJMS-OKSMBDkldmrmt-lfhqwhig} is
essentially optimal. For example, we can take~$\phi\in C^\infty_0(B_1,[0,+\infty))$ with~$\int_{B_1}\phi(x)\,dx=1$
and define~$f_j(x):=j^n\phi(jx)$. We also assume that the function~$g$ in~\eqref{KERNE-SINGO-001}
satisfies~$|g|\ge a$ for some~$a>0$ in some region~$\Sigma\subseteq\partial B_1$ with~${\mathcal{H}}^{n-1}(\Sigma)>0$. We also fix~$\lambda_0>0$ and take~$\lambda\in(0,\lambda_0)$.
In this setting, we have that
$$ \|f_j\|_{L^1(\R^n)}=j^n\int_{\R^n}\phi(jx)\,dx=\int_{B_1}\phi(y)\,dy=1.$$
Moreover, by~\eqref{ABBE9-NOBA} and~\eqref{ojndPP08uj3rHOmaT7u0-n0epsjf},
\begin{eqnarray*}&& |Tf_j(x)|=\left|\int_{\R^n} K(x-y)\,f_j(y)\,dy\right|=
\left|j^n\int_{\R^n} K(x-y)\,\phi(j y)\,dy\right|=
\left|\int_{B_1} K\left(x-\frac{z}j\right)\,\phi(z)\,dz\right|\\&&\qquad\qquad
\ge
\left|K(x)\int_{B_1} \phi(z)\,dz\right|-
\left|\int_{B_1} \left(K\left(x-\frac{z}j\right)-K(x)\right)\,\phi(z)\,dz\right|\\&&\qquad\qquad
\ge |K(x)|-\int_{B_1} \left|K\left(x-\frac{z}j\right)-K(x)\right|\,\phi(z)\,dz
\\&&\qquad\qquad
= \frac{1}{|x|^n}\left|g\left( \frac{x}{|x|}\right)\right|-\int_{B_1} \left|K\left(x-\frac{z}j\right)-K(x)\right|\,\phi(z)\,dz.
\end{eqnarray*}
Thus, if
$$ \Sigma_\lambda:=\left\{ r\omega, {\mbox{ with $\omega\in\Sigma $ and }} r\in \left[
\left( \frac{a}{4\lambda}\right)^{\frac1n},\left( \frac{a}{2\lambda}\right)^{\frac1n}
\right]\right\}$$
and~$x\in\Sigma_\lambda$, we have that~$|x|^n\le\frac{a}{2\lambda}$ and~$\left|g\left( \frac{x}{|x|}\right)\right|\ge a$.
{F}rom these observations we arrive at
\begin{equation}\label{QCNSPVNSPCS} |Tf_j(x)|\ge 2\lambda-\int_{B_1} \left|K\left(x-\frac{z}j\right)-K(x)\right|\,\phi(z)\,dz.\end{equation}
Also, if~$z\in B_1$ and~$t\in[0,1]$,
$$ \left|x-\frac{z}j\right|\ge|x|-\frac1j\ge\frac{|x|}2+\frac12\,\left( \frac{a}{4\lambda}\right)^{\frac1n}-\frac1j\ge
\frac{|x|}2$$
as long as~$j$ is large enough. In this setting, using~\eqref{KERNE-SINGO-003},
\begin{eqnarray*}&&
\int_{B_1} \left|K\left(x-\frac{z}j\right)-K(x)\right|\,\phi(z)\,dz\le
\frac{\|\phi\|_{L^\infty(B_1)}}j\iint_{B_1\times[0,1]}\left|\nabla K\left(x-\frac{z}j\right)\right|\,dz\,dt
\\&&\qquad\qquad
\le\frac{C}j\iint_{B_1\times[0,1]}\frac{dz\,dt}{
\left|x-\frac{z}j\right|^{n+1}}\le
\frac{C}j\iint_{B_1\times[0,1]}\frac{dz\,dt}{
|x|^{n+1}}\le\frac{C\lambda^{\frac{n+1}{n}} }{a^{\frac{n+1}{n}}j}\le\lambda
\end{eqnarray*}
for sufficiently large~$j$ (depending also on~$\lambda_0$)
and consequently, by~\eqref{QCNSPVNSPCS}, $|Tf_j(x)|\ge\lambda$.

This gives that~$\{|Tf_j|\ge\lambda\}\supseteq\Sigma_\lambda$ for~$j$ large enough, whence
\begin{eqnarray*}
\frac{\big| \{ |Tf_j|\ge\lambda\}\big|}{\|f_j\|_{L^1(\R^n)}}=\big| \{ |Tf_j|\ge\lambda\}\big|\ge
|\Sigma_\lambda|=\iint_{\left[
\left( \frac{a}{4\lambda}\right)^{\frac1n},\left( \frac{a}{2\lambda}\right)^{\frac1n}
\right]\times\Sigma} r^{n-1}\,dr\,d{\mathcal{H}}^{n-1}_\omega
=\frac{ a\,{\mathcal{H}}^{n-1}(\Sigma) }{4n\lambda} , 
\end{eqnarray*}
showing the optimality of Theorem~\ref{hoiFC-AJS-OJMS-OKSMBDkldmrmt-lfhqwhig}.\medskip

Another important ingredient towards the proof of
Theorem~\ref{T:CAZY-PIVO}
consists in a classical interpolation\index{interpolation} result\footnote{Marcinkiewicz has been one of the brilliant students and collaborators of \label{JNDtsXVSBNS-KASrSDa3c5h4245e2r}
Antoni
Zygmund, as well as a collaborator of Juliusz Schauder (this biographical
fact provides also a conceptual bridge between the results
presented here and those in Chapter~\ref{C2ALPHACHAPET}). In~1939
Marcinkiewicz was
taken as a Polish prisoner of war to a Soviet camp.
Shortly before dying at the age of 30, he
gave his manuscripts to his parents, who unfortunately were also imprisoned and died of hunger in a camp.
Unsurprisingly, no trace of these manuscripts was ever found,
adding one more tragic story to the uncountable atrocities of World War II.
See~\cite{MR0168434} for a brief description by Zygmund of related events
and of the lost of Marcinkiewicz's latest mathematical works.
Compare also with footnote~\ref{JNDtsXVSBNS-KASrSDa3c5h4245e2r2} on page~\pageref{JNDtsXVSBNS-KASrSDa3c5h4245e2r2}.}
due to
J\'{o}zef Marcinkiewicz~\cite{zbMATH03033497}
(see also~\cite{MR80887}). The general idea behind interpolation,
roughly speaking, is that if one knows something about a functional space~$X$
and something about a functional space~$Z$ ideally it is possible to infer some information
about a functional space~$Y$ which ``lies between~$X$ and~$Z$''.
An illustrative example can be taken into account when~$X:=L^P(\R^n)$ and~$Z:=L^R(\R^n)$,
with~$R> P\ge1$, and~$Y:=L^Q(\R^n)$ with~$Q\in(P,R)$. In this case,
the H\"older Inequality with exponents~$\frac{R-P}{R-Q}$ and~$\frac{R-P}{Q-P}$
gives that if~$v\in L^P(\R^n)\cap L^R(\R^n)$, then~$v\in L^Q(\R^n)$, with
\begin{eqnarray*}&&
\|v\|_{L^Q(\R^n)}=\left(
\int_{\R^n}|v(x)|^Q\,dx
\right)^{\frac1Q}=\left(
\int_{\R^n}|v(x)|^{\frac{(R-Q)P}{R-P}+\frac{(Q-P)R}{R-P}}\,dx
\right)^{\frac1Q}\\&&\qquad\qquad\qquad\le
\left(
\int_{\R^n}|v(x)|^{P}\,dx
\right)^{ \frac{R-Q}{(R-P)Q} }
\left(
\int_{\R^n}|v(x)|^{{R}}\,dx
\right)^{\frac{Q-P}{(R-P)Q}}
=\|v\|_{L^P(\R^n)}^{ \frac{(R-Q)P}{(R-P)Q} }\,
\|v\|_{L^R(\R^n)}^{\frac{(Q-P)R}{(R-P)Q}}.
\end{eqnarray*}
The interpolation result that we need in our framework is much more sophisticated than that
and, in light of Theorem~\ref{hoiFC-AJS-OJMS-OKSMBDkldmrmt-lfhqwhig},
addresses a ``weak theory in~$L^p$''.
For this, we say that an operator~$S$ acting on the space of functions is {\em sublinear} if for every functions~$f$ and~$g$
it holds that
\begin{equation}\label{35MS9uyhdfiowghf93wutygfuywegf7wftqfwguweghSERdfgbLFVBdfNOMDUFJNTInNOA41}
|S(f+g)(x)|\le|Sf(x)|+|Sg(x)|\qquad{\mbox{ for all }}\,x\in\R^n.
\end{equation}
For all~$p\in[1,+\infty]$ we say that~$S$ is of {\em weak type~$(p,p)$} if there exists~$C>0$
such that, for all~$f\in L^p(\R^n)$ and~$\lambda>0$,
\begin{equation}\label{KS03irjg9043hgv0jhf0h9ugvoewtr8732gDCAdncavbikgKAmtR3EW}
\big| \{ |Sf|\ge\lambda\}\big|\le\frac{C\,\|f\|^p_{L^p(\R^n)}}{\lambda^p}.
\end{equation}
Comparing with the Chebyshev's Inequality in~\eqref{35MS9uyhdfiowghf93wutygfuywegf7wftqfwguweghSERdfgbLFVBdfNOMDUFJNTInNOA},
one may think, heuristically, that when~$S$ is of weak type~$(p,p)$
one has that ``$Sf$ is almost in~$L^p(\R^n)$''.

One also says that~$S$ is of {\em strong type~$(p,p)$} if there exists~$C>0$
such that, for all~$f\in L^p(\R^n)$,
\begin{equation}\label{KS03irjg9043hgv0jhf0h9ugvoewtr8732gDCAdncavbikgKAmtR3E}
\|Sf\|_{L^p(\R^n)}\le C\,\|f\|_{L^p(\R^n)},\end{equation}
that is~$S$ is a bounded operator on~$L^p(\R^n)$.

Using the Chebyshev's Inequality in~\eqref{35MS9uyhdfiowghf93wutygfuywegf7wftqfwguweghSERdfgbLFVBdfNOMDUFJNTInNOA},
one sees that
\begin{equation}\label{KS03irjg9043hgv0jhf0h9ugvoewtr8732gDCAdncavbikgKAmtR3E-090-090o}
{\mbox{if~$S$ is of strong type~$(p,p)$ then it is also of weak type~$(p,p)$.}}\end{equation}

The interpolation result that we use in this framework is the following one:

\begin{figure}
  \centering
  \includegraphics[width=.65\linewidth]{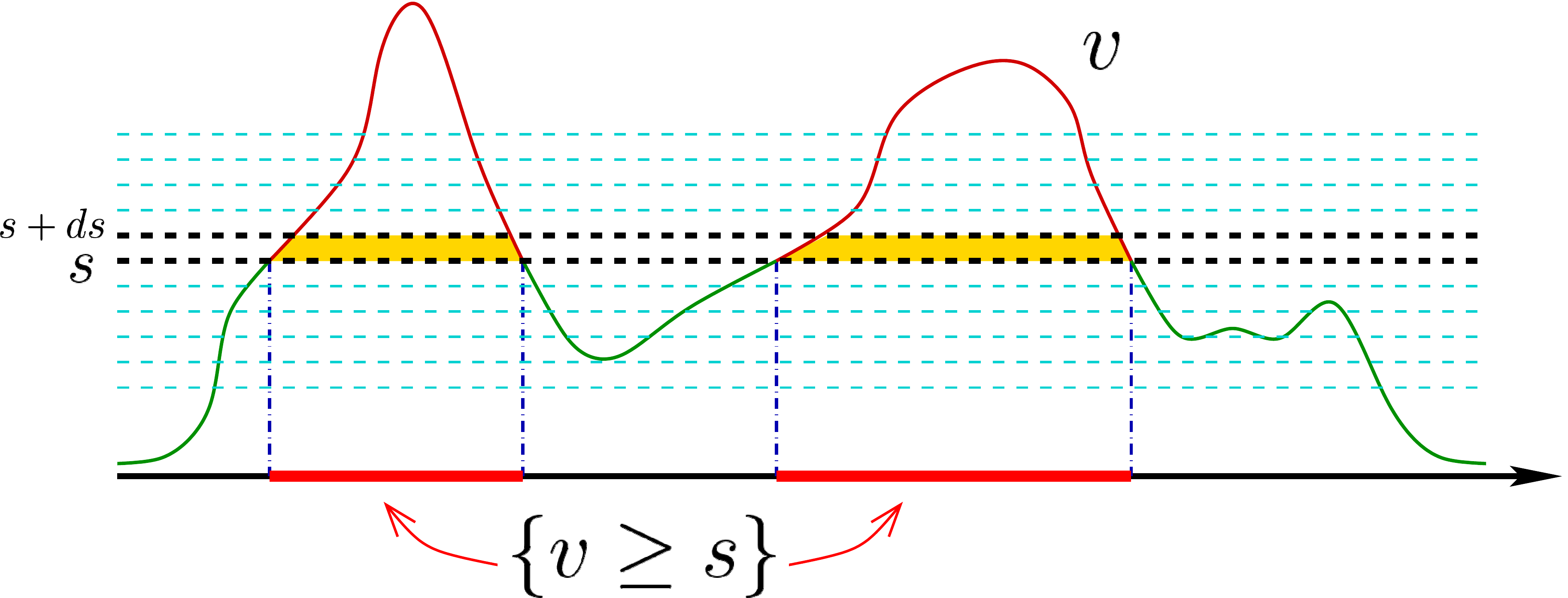}
 \caption{\sl The Layer Cake Representation Formula when~$v\ge0$ and~$m=1$.}\label{2356457SUEGGIUTERRAeFICAKE}
\end{figure}

\begin{theorem}\label{UIMSINTEPAMANDR}
Let~$p$, $r\in[1,+\infty)$ with~$p<r$.
Let~$S$ be a sublinear operator acting on the space of functions.
Assume that~$S$ is of weak type~$(p,p)$ and of weak type~$(r,r)$.

Then, $S$ is of strong type~$(q,q)$ for every~$q\in(p,r)$.\medskip

Also, the corresponding constant in~\eqref{KS03irjg9043hgv0jhf0h9ugvoewtr8732gDCAdncavbikgKAmtR3E}
for~$L^q(\R^n)$ depends only on~$n$, $p$, $q$, $r$ and the corresponding constants in~\eqref{KS03irjg9043hgv0jhf0h9ugvoewtr8732gDCAdncavbikgKAmtR3EW} for~$L^p(\R^n)$ and~$L^r(\R^n)$.
\end{theorem}

\begin{figure}
  \centering
  \includegraphics[width=.65\linewidth]{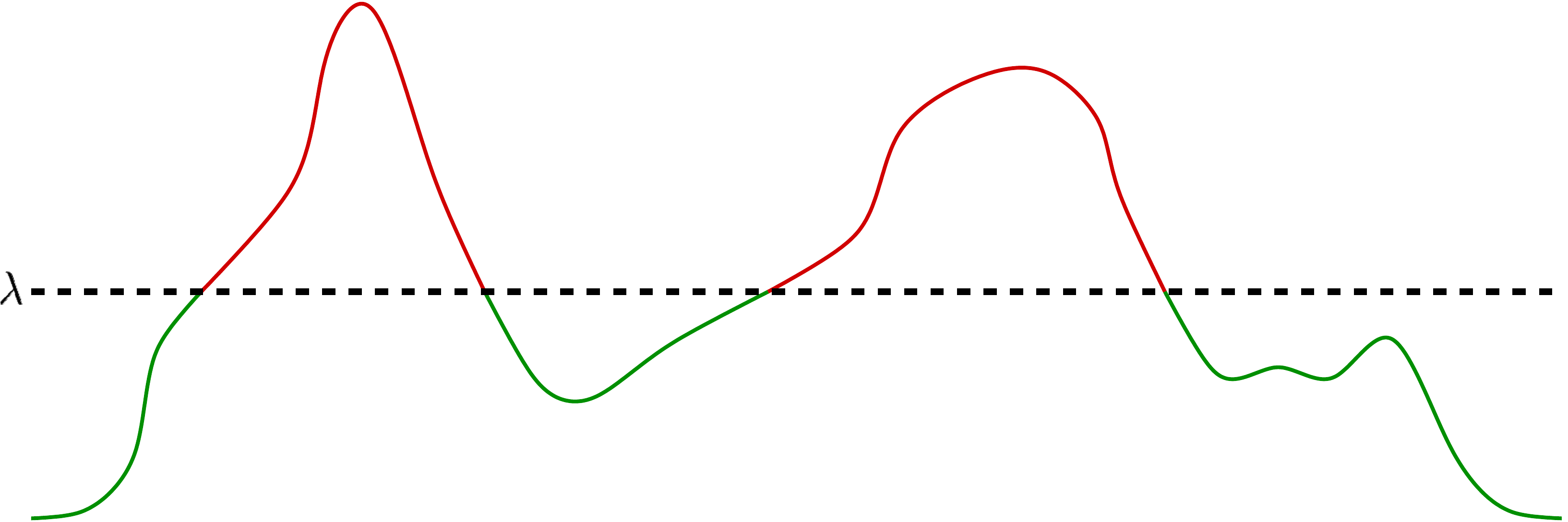}
 \caption{\sl Splitting a function into its ``top'' and ``bottom'' parts.}\label{2356457SUEGGIUTERRAeFI}
\end{figure}

\begin{proof} We will repeatedly exploit the Layer Cake Representation Formula (valid for
all~$m\in(0,+\infty)$, see e.g.~\cite[Theorem~5.51]{MR3381284}, see Figure~\ref{2356457SUEGGIUTERRAeFICAKE}
for a schematic representation when~$v\ge0$ and~$m=1$)
\begin{equation} \label{LCKRPWZMALLIOMPSII}
\int_{\R^n} |v(x)|^m\,dx=m\int_0^{+\infty} s^{m-1}\big|\{|v|\ge s\}\big|\,ds.\end{equation}
Let~$\lambda>0$. We split a function~$f$ into its ``top and bottom parts'', that we denote by~$f_t$ and~$f_b$ respectively, see Figure~\ref{2356457SUEGGIUTERRAeFI}, that is we define
\begin{eqnarray*}
f_t:=f\chi_{\{|f|\ge\lambda\}}\qquad{\mbox{and}}\qquad
f_b:=f\chi_{\{|f|<\lambda\}}.
\end{eqnarray*}
Notice that~$f=f_t+f_b$
Since~$S$ is sublinear, by~\eqref{35MS9uyhdfiowghf93wutygfuywegf7wftqfwguweghSERdfgbLFVBdfNOMDUFJNTInNOA41} we know that~$|Sf|\le|S f_t|+|Sf_b|$ and consequently
\begin{equation*}
\big\{ |Sf|\ge\lambda\big\}\,\subseteq\,\left\{|Sf_t|\ge\frac\lambda2\right\}\cup\left\{|Sf_b|\ge\frac\lambda2\right\}.
\end{equation*}
As a result,
\begin{equation}\label{2356457SUEGGIUTERRAeFI-CHSJDLUNDRB}
\big|\big\{ |Sf|\ge\lambda\big\}\big|\,\leq\,
\left|\left\{|Sf_t|\ge\frac\lambda2\right\}\right|+\left|\left\{|Sf_b|\ge\frac\lambda2\right\}\right|.
\end{equation}
Our goal is thus to estimate the two quantities on the right hand side of~\eqref{2356457SUEGGIUTERRAeFI-CHSJDLUNDRB}
by using the assumptions that~$S$ is of weak type~$(p,p)$ and of weak type~$(r,r)$.
Specifically, when dealing with the bottom function~$f_b$,
we will account for the ``tail'' of the function~$f$ and therefore,
to enhance the integrability property as much as possible, it will be convenient
to take into consideration the largest exponent~$r$ (because a small number elevated to a large
power becomes even smaller). Instead, when dealing with the top function~$f_t$,
we will consider the peaks of the function and therefore it
will be appropriate to utilize the smallest exponent~$p$ (because a large number elevated to a small
power becomes less dangerous).

The technical details go as follows. For the bottom function, using that~$S$ is of weak type~$(r,r)$
and~\eqref{LCKRPWZMALLIOMPSII} we see that
\begin{equation}
\begin{split}\label{LCKRPWZMALLIOMPSIIEM}&
\left|\left\{|Sf_b|\ge\frac\lambda2\right\}\right|\le
\frac{C\,\|f_b\|^r_{L^r(\R^n)}}{\lambda^r}=\frac{Cr}{\lambda^r}
\int_0^{+\infty} s^{r-1}\big|\{|f_b|\ge s\}\big|\,ds\\&\qquad\qquad\quad\qquad\qquad=\frac{Cr}{\lambda^r}
\int_0^{\lambda} s^{r-1}\big|\{|f|\ge s\}\big|\,ds.\end{split}
\end{equation}
Similarly, for the top function,
using that~$S$ is of weak type~$(p,p)$
and~\eqref{LCKRPWZMALLIOMPSII} we see that
\begin{equation*}
\begin{split}&
\left|\left\{|Sf_t|\ge\frac\lambda2\right\}\right|\le
\frac{C\,\|f_t\|^p_{L^p(\R^n)}}{\lambda^p}=\frac{Cp}{\lambda^p}
\int_0^{+\infty} s^{p-1}\big|\{|f_t|\ge s\}\big|\,ds\\&\quad\qquad\qquad\qquad\qquad=\frac{Cp}{\lambda^p}
\int_{\lambda}^{+\infty} s^{p-1}\big|\{|f|\ge s\}\big|\,ds.\end{split}
\end{equation*}
{F}rom this, \eqref{LCKRPWZMALLIOMPSII},
\eqref{2356457SUEGGIUTERRAeFI-CHSJDLUNDRB},
and~\eqref{LCKRPWZMALLIOMPSIIEM} we find that
\begin{eqnarray*}&&
\int_{\R^n} |Sf(x)|^q\,dx\\&=& q\int_0^{+\infty} \lambda^{q-1}\big|\{|Sf|\ge \lambda\}\big|\,d\lambda
\\&\le&q\int_0^{+\infty} \lambda^{q-1}\left|\left\{|Sf_t|\ge \frac\lambda2\right\}\right|\,d\lambda+
q\int_0^{+\infty} \lambda^{q-1}\left|\left\{|Sf_b|\ge \frac\lambda2\right\}\right|\,d\lambda\\
&\le&C\left[\int_0^{+\infty} \left(\int_\lambda^{+\infty}
\lambda^{q-p-1}s^{p-1}\big|\{|f|\ge s\}\big|\,ds\right)\,d\lambda+
\int_0^{+\infty} \left(\int_0^\lambda
\lambda^{q-r-1} s^{r-1}\big|\{|f|\ge s\}\big|\,ds\right)\,d\lambda\right]
,\end{eqnarray*}
up to renaming constants.

Thus, switching the order of integration and exploiting that~$p<q<r$,
\begin{eqnarray*}&&
\int_{\R^n} |Sf(x)|^q\,dx\\
&\le&C\left[\int_0^{+\infty} \left(\int_0^s
\lambda^{q-p-1}s^{p-1}\big|\{|f|\ge s\}\big|\,d\lambda\right)\,ds+
\int_0^{+\infty} \left(\int_s^{+\infty}
\lambda^{q-r-1} s^{r-1}\big|\{|f|\ge s\}\big|\,d\lambda\right)\,ds\right]\\
&\le&C\left[\int_0^{+\infty} s^{q-p+p-1}\big|\{|f|\ge s\}\big|\,ds+
\int_0^{+\infty} s^{q-r+r-1}\big|\{|f|\ge s\}\big|\,ds\right]\\
&=& C\int_0^{+\infty} s^{q-1}\big|\{|f|\ge s\}\big|\,ds.
\end{eqnarray*}
Hence, making use of~\eqref{LCKRPWZMALLIOMPSII} once again,
$$ \int_{\R^n} |Sf(x)|^q\,dx\le C\int_{\R^n}|f(x)|^q\,dx,$$
which is the desired result.
\end{proof}

With the previous work and a duality argument, we are now in the position of completing the
proof of Theorem~\ref{T:CAZY-PIVO}.

\begin{proof}[Proof of Theorem~\ref{T:CAZY-PIVO}] When~$p=2$, the desired result follows from~\eqref{L2THESINGOPE}.

When~$p\in(1,2)$ we use the interpolation theory: namely,
we know from~\eqref{L2THESINGOPE} that~$T$ is of strong type~$(2,2)$,
and therefore of weak type~$(2,2)$, thanks to~\eqref{KS03irjg9043hgv0jhf0h9ugvoewtr8732gDCAdncavbikgKAmtR3E-090-090o}.
Also, we know from Theorem~\ref{hoiFC-AJS-OJMS-OKSMBDkldmrmt-lfhqwhig} that~$T$ is of weak type~$(1,1)$.
Hence, it follows from Theorem~\ref{UIMSINTEPAMANDR} that~$T$ is of strong type~$(p,p)$ for all~$p\in(1,2)$
and this proves Theorem~\ref{T:CAZY-PIVO} when~$p\in(1,2)$.

It remains to prove Theorem~\ref{T:CAZY-PIVO} when~$p>2$. For this we use a duality argument.
Let~$q:=\frac{p}{p-1}\in(1,2)$ be the conjugated exponent of~$p$ and let~$g\in C^{0,1}_0(\R^n)$.
Let also~$\widetilde{K}(x):=K(-x)$ and~$\widetilde T$ be as in~\eqref{ojndPP08uj3rHOmaT7u0-n0epsjf}
but with~$\widetilde{K}$ replacing~$K$, that is
\begin{equation*} \widetilde Tf(x):=\int_{\R^n} \widetilde K(x-y)\,\big(f(x)-f(y)\big)\,dy.\end{equation*}
We stress that~$\widetilde{K}$ also satisfies the structural assumptions~\eqref{KERNE-SINGO-001},
\eqref{KERNE-SINGO-003} and~\eqref{ABBE9-NOBA}. Consequently, since we have already
established Theorem~\ref{T:CAZY-PIVO} when~$p\in(1,2)$, we have that
\begin{equation}\label{FUBYTSJMNELL0} 
 \|\widetilde T g\|_{L^q(\R^n)}\le C\,\|g\|_{L^q(\R^n)}.\end{equation}
Now we claim that
\begin{equation}\label{FUBYTSJMNELL}
\int_{\R^n} Tf(x)\,g(x)\,dx=\int_{\R^n} \widetilde Tg(x)\,f(x)\,dx.
\end{equation}
For this, let~$K_\e$ be as in~\eqref{TAGLIOKETR} and~$T_\e$ be as in~\eqref{TAGLIOKETR0}.
Recalling~\eqref{MSrruj0GAnTBGSnak-123} we have that
\begin{eqnarray*}
\lim_{\e\searrow0} \int_{\R^n} |Tf(x)-T_\e f(x)|\,|g(x)|\,dx\le
\lim_{\e\searrow0}\|T f-T_\e f\|_{L^2(\R^n)}\,\|g\|_{L^2(\R^n)}=0.
\end{eqnarray*}
Similarly, if~$\widetilde K_\e$ and~$\widetilde T_\e$ are as in~\eqref{TAGLIOKETR} and~\eqref{TAGLIOKETR0}
but with~$\widetilde K$ in place of~$K$, then
$$ \lim_{\e\searrow0} \int_{\R^n} |\widetilde T g(x)-\widetilde T_\e g(x)|\,|f(x)|\,dx=0.$$
As a result, using~\eqref{ABBE9-NOBA},
\begin{eqnarray*}&&
\int_{\R^n} Tf(x)\,g(x)\,dx-\int_{\R^n} \widetilde Tg(x)\,f(x)\,dx\\&=&
\lim_{\e\searrow0}
\int_{\R^n} T_\e f(x)\,g(x)\,dx-\int_{\R^n} \widetilde T_\e g(x)\,f(x)\,dx\\
\\&=&
\lim_{\e\searrow0}
\int_{\R^n} \left( \int_{\R^n} K_\e(x-y)\big(f(x)-f(y)\big)\,dy\right)\,g(x)\,dx\\&&\qquad
-\int_{\R^n} \left( \int_{\R^n} \widetilde K_\e(x-y)\big(g(x)-g(y)\big)\,dy\right)\,f(x)\,dx\\&=&
\lim_{\e\searrow0}\left[
-\int_{\R^n} \left( \int_{\R^n} K_\e(x-y)\,f(y)\,dy\right)\,g(x)\,dx
+\int_{\R^n} \left( \int_{\R^n} K_\e(y-x)g(y)\,dy\right)\,f(x)\,dx\right]
\\&=&0,
\end{eqnarray*}
owing to Fubini's Theorem. The proof of~\eqref{FUBYTSJMNELL} is thereby complete.

Thus, by~\eqref{FUBYTSJMNELL0} and~\eqref{FUBYTSJMNELL},
\begin{eqnarray*}&&
\left|\int_{\R^n} Tf(x)\,g(x)\,dx\right|=\left|\int_{\R^n} \widetilde Tg(x)\,f(x)\,dx\right|
\le\int_{\R^n} |\widetilde Tg(x)|\,|f(x)|\,dx\\&&\quad\qquad\le
\|\widetilde Tg\|_{L^q(\R^n)}\,\|f\|_{L^p(\R^n)}\le C\|g\|_{L^q(\R^n)}\,\|f\|_{L^p(\R^n)}.
\end{eqnarray*}
By the density of~$C^{1,0}_0(\R^n)$ in~$L^q(\R^n)$
(see e.g.~\cite[Theorem~1.91]{MR2895178}), this estimate remains valid for all~$g\in L^q(\R^n)$.
Accordingly, the linear functional~${\mathcal{L}}:L^q(\R^n)\to\R$ defined by
$$ {\mathcal{L}}(g):=\int_{\R^n} Tf(x)\,g(x)\,dx$$
is bounded and therefore, by duality (see e.g.~\cite[Theorem~5.12.1]{MR3753707}),
there exists a unique function~$\phi\in L^p(\Omega)$ such that~${\mathcal{L}}(g)=\int_{\R^n} \phi(x)\,g(x)\,dx$
for all~$g\in L^q(\R^n)$ and
$$ \|\phi\|_{L^p(\Omega)}=\|{\mathcal{L}}\|=\sup_{g\in L^q(\Omega)\setminus\{0\}}\frac{|{\mathcal{L}}(g)|}{\|g\|_{L^q(\Omega)}}.$$
Since necessarily~$\phi=Tf$, we obtain that
$$ \|Tf\|_{L^p(\Omega)}=\|\phi\|_{L^p(\Omega)}=
\sup_{g\in L^q(\Omega)\setminus\{0\}}\frac{1}{\|g\|_{L^q(\Omega)}}
\left|\int_{\R^n} Tf(x)\,g(x)\,dx\right|\le C\|f\|_{L^p(\R^n)}.$$
This finishes the proof of Theorem~\ref{T:CAZY-PIVO}.
\end{proof}

We observe that~Theorem~\ref{T:CAZY-PIVO} carries over to all functions~$f\in L^p(\R^n)$,
by the density of~$C^{1,0}_0(\R^n)$ in~$L^p(\R^n)$
(see e.g.~\cite[Theorem~1.91]{MR2895178}), that is~$T$ extends by continuity to a bounded linear operator from~$L^p(\R^n)$ to~$L^p(\R^n)$:
we stress however that for a general function~$f$ in~$L^p(\R^n)$ the pointwise
definition of~$Tf$ in~\eqref{ojndPP08uj3rHOmaT7u0-n0epsjf}
does not necessarily makes sense and one has to intend~$Tf$ in this setting as the limit in~$L^p(\R^n)$
of~$Tf_j$, where~$f_j\in C^{0,1}_0(\R^n)$ and~$\|f_j-f\|_{L^p(\R^n)}\to0$ as~$j\to+\infty$.
%%%% With this extended definition of~$Tf$ in mind, we state explicitly the above observation for the sake of completeness:
%%%% 
%%%% \begin{corollary}
%%%% Let~$p\in(1,+\infty)$. Then, there exists a positive constant~$C$, depending only on~$n$, $p$ and~$K$,
%%%% such that
%%%% $$ \|Tf\|_{L^p(\R^n)}\le C\,\|f\|_{L^p(\R^n)}$$
%%%% for all~$f\in L^p(\R^n)$.
%%%% \end{corollary}

As a consequence of Theorem~\ref{T:CAZY-PIVO}, we have:

\begin{corollary}\label{CXf43Unfsa0-430u4noerngTheoremT:CAZY-PIVO}
Let~$p\in(1,+\infty)$.
Let~$\Gamma$ be the fundamental solution in~\eqref{GAMMAFU} and~$f\in C^{0,1}_0(\R^n)$.
Let~$v:=-\Gamma*f$. Then, $v\in W^{2,p}(\R^n)$ and
$$\|D^2 v\|_{L^p(\R^n)}\le C\,\|f\|_{L^p(\R^n)},$$
where~$C>0$ depends only on~$n$ and~$p$.
\end{corollary}

\begin{proof} By~\eqref{POTNDEWSECDERIVA},
\eqref{9UOHNS-01euowfhvbogdblwiqfgyoerufg8veduifgedirIPSEtbcTIZnowevd99e}
and~\eqref{KERNE-SINGO-002}, we can write~$\partial_{ij}v$ as~$Tf$ satisfying conditions~\eqref{KERNE-SINGO-001},
\eqref{KERNE-SINGO-003} and~\eqref{ABBE9-NOBA}. Hence the desired result is now a direct consequence of
Theorem~\ref{T:CAZY-PIVO}.
\end{proof}

We are now in the position of addressing the core of the regularity theory in $L^p$ spaces.
We start with the following interior estimates:

\begin{theorem}\label{982nksjdvnICtsghdfnEJHACANosd90hyri32t8043hger-3506o5y}
Let~$p\in(1,+\infty)$ and~$f\in C^{0,1}(B_1)$. Assume that~$u\in C^2(B_1)$ is a solution of~$\Delta u=f$ in~$B_1$.

Then,
\begin{equation*}
\|u\|_{W^{2,p}(B_{1/2})}\le C\,\Big( \|u\|_{L^p(B_1)}+\|f\|_{L^p(B_1)}\Big),\end{equation*}
where~$C>0$ depends only on~$n$ and~$p$.
\end{theorem}

\begin{proof} The attentive reader will find close connections between
this proof and that of Theorem~\ref{SCHAUDER-INTE}. The main additional ingredient here
is of course the Newtonian potential estimate in~$L^p$ spaces provided by Corollary~\ref{CXf43Unfsa0-430u4noerngTheoremT:CAZY-PIVO}. We will also employ the following useful inequality of general use, valid for all measurable
functions~$g\in W^{2,p}((0,1)^n)$ and all~$\e\in\left(0,\frac1{4}\right)$:
\begin{equation}\label{KS-0ujwfn-092urjnf-i2jrfoneh9843hgb}
\|\nabla g\|_{L^p((0,1)^n)}\le \e\|D^2 g\|_{L^p((0,1)^n)}+\frac{C}{\e}\,\|g\|_{L^p((0,1)^n)},
\end{equation}
with~$C>0$ depending only on~$n$ and~$p$.
The above is a very particular case
of the Gagliardo-Nirenberg Interpolation Inequality (see~\cite[Theorem~1.2]{MR4237368}, used here with~$k:=2$, $j:=1$, $r:=q:=p$, see also~\cite{MR109295, MR109940} for the original articles on this topic and~\cite[Sections~12.5, 13.3, 16.1 and~16.2]{MR2527916}
for a thorough presentation of interpolation methods). 
We give a detailed proof of~\eqref{KS-0ujwfn-092urjnf-i2jrfoneh9843hgb} for completeness. First of all,
we observe that
\begin{equation}\label{KS-0ujwfn-092urjnf-i2jrfoneh9843hgb2ZM}
{\mbox{it is enough to establish~\eqref{KS-0ujwfn-092urjnf-i2jrfoneh9843hgb} when }}\frac1{3\e}\in\N .\end{equation}
Indeed, if~$\e\in\left(0,\frac14\right)$ we let~$N\in\N$ such that~$\frac1{3\e}\in \left[N,N+1 \right)$.
Hence~$\frac{1}{3(N+1)}\in\left(0,\e\right]\subseteq
\left(0,\frac14\right)$ and thus, if~\eqref{KS-0ujwfn-092urjnf-i2jrfoneh9843hgb}
is valid for~$\e$ replaced by~$\frac{1}{3(N+1)}$ we have that
\begin{eqnarray*}&& \|\nabla g\|_{L^p((0,1)^n)}\le \frac{1}{3(N+1)}\,\|D^2 g\|_{L^p((0,1)^n)}+3C(N+1)\,\|g\|_{L^p((0,1)^n)}\\&&\qquad\le
\e\|D^2 g\|_{L^p((0,1)^n)}+6CN\,\|g\|_{L^p((0,1)^n)}\le
\e\|D^2 g\|_{L^p((0,1)^n)}+\frac{2C}{\e}\,\|g\|_{L^p((0,1)^n)}
,\end{eqnarray*}
giving~\eqref{KS-0ujwfn-092urjnf-i2jrfoneh9843hgb} for~$\e$ as well. This establishes~\eqref{KS-0ujwfn-092urjnf-i2jrfoneh9843hgb2ZM}.

In view of~\eqref{KS-0ujwfn-092urjnf-i2jrfoneh9843hgb2ZM},
from now on we will suppose that~$\frac1{3\e}\in\N$.

Additionally,
we point out that we can reduce to dimension~$n=1$, namely
\begin{equation}\label{KS-0ujwfn-092urjnf-i2jrfoneh9843hgb2}
{\mbox{it is enough to establish~\eqref{KS-0ujwfn-092urjnf-i2jrfoneh9843hgb} when~$n=1$
and~$g\in W^{2,p}((0,1))$.}}
\end{equation}
Indeed, if~$g\in W^{2,p}((0,1)^n)$, given~$x'=(x_1,\dots,x_{n-1})\in(0,1)^{n-1}$
we can define~$\widetilde g(x_n):=g(x',x_n)$. If~\eqref{KS-0ujwfn-092urjnf-i2jrfoneh9843hgb} holds true when~$n=1$ then,
up to rename constants,
\begin{eqnarray*}&& 
\int_{0}^1 |\partial_n g(x',x_n)|^p\,dx_n=
\int_{0}^1 |\widetilde g'(x_n)|^p\,dx_n\le
\e^p\int_{0}^1 |\widetilde g''(x_n)|^p\,dx_n+\frac{C}{\e^p}\int_{0}^1 |\widetilde g(x_n)|^p\,dx_n
\\&&\qquad=\e^p\int_{0}^1 |\partial_{nn} g(x',x_n)|^p\,dx_n+\frac{C}{\e^p}\int_{0}^1 |g(x',x_n)|^p\,dx_n.\end{eqnarray*}
Thus, integrating over~$x'\in(0,1)^{n-1}$,
\begin{eqnarray*}&& 
\int_{(0,1)^n} |\partial_n g(x)|^p\,dx\le\e^p\int_{(0,1)^n} |\partial_{nn} g(x)|^p\,dx+\frac{C}{\e^p}\int_{(0,1)^n} |g(x)|^p\,dx\\&&\qquad\qquad
\le \e^p\|D^2g\|_{L^p((0,1)^n)}^p+\frac{C}{\e^p}\,\|g\|_{L^p((0,1)^n)}^p.
\end{eqnarray*}
By exchanging the order of the variables, this entails that
\begin{eqnarray*}&& 
\int_{(0,1)^n} |\partial_j g(x)|^p\,dx\le\e^p\|D^2g\|_{L^p((0,1)^n)}^p+\frac{C}{\e^p}\,\|g\|_{L^p((0,1)^n)}^p\qquad{\mbox{for all }}\;j\in\{1,\dots,n\},
\end{eqnarray*}
which yields~\eqref{KS-0ujwfn-092urjnf-i2jrfoneh9843hgb} in its full generality.
This proves~\eqref{KS-0ujwfn-092urjnf-i2jrfoneh9843hgb2}.

Hence, in view of~\eqref{KS-0ujwfn-092urjnf-i2jrfoneh9843hgb2}, we now focus on the proof
of~\eqref{KS-0ujwfn-092urjnf-i2jrfoneh9843hgb} when~$n=1$. Additionally, by the density of~$C^\infty((0,1))$
in~$W^{2,p}((0,1))$ (see e.g.~\cite[Theorem~2, page~251]{MR1625845}), we can assume that~$g\in C^\infty((0,1))$.
By the Calculus Mean Value Theorem,
given~$x> y\in(0,1)$, there exists~$z_{x,y}\in[y,x]$ such that
$$ g(x)-g(y)=g'(z_{x,y})(x-y).$$
Therefore, for all~$w\in(0,1)$,
$$ g'(w)=g'(w)-g'(z_{x,y}) +\frac{g(x)-g(y)}{x-y}
=\int^w_{z_{x,y}}g''(t)\,dt+\frac{g(x)-g(y)}{x-y}.$$
In particular, we consider an interval of the form~$[a,a+3\e]$, for some~$a\in(0,1-3\e)$,
we assume that~$y\in[a,a+\e]$ and~$x\in[a+2\e,a+3\e]$. In this way,~$x-y\ge\e$ and thus, for all~$w\in[a,a+3\e]$,
$$ |g'(w)|\le \int_a^{a+3\e}|g''(t)|\,dt+\frac{|g(x)|+|g(y)|}{\e}.$$
As a consequence, possibly renaming~$C$ at each stage of the calculation,
\begin{eqnarray*} |g'(w)|^p&\le& 
\left( \int_a^{a+3\e}|g''(t)|\,dt+\frac{|g(x)|+|g(y)|}{\e}\right)^p\\&
\le& C\left( \left(\int_a^{a+3\e}|g''(t)|\,dt\right)^p+\frac{|g(x)|^p+|g(y)|^p}{\e^p}\right)
\\&\le&
C\left( \e^{p-1}\int_a^{a+3\e}|g''(t)|^p\,dt +\frac{|g(x)|^p+|g(y)|^p}{\e^p}\right).
\end{eqnarray*}
Integrating over~$x\in[a+2\e,a+3\e]$, $y\in[a,a+\e]$ and~$
w\in[a,a+3\e]$, and dividing by~$\e^2$, we obtain that
\begin{eqnarray*}
\int_{a}^{a+3\e} |g'(w)|^p\,dw&\le& 
C\left( \e^{p}\int_a^{a+3\e}|g''(t)|^p\,dt +\frac{1}{\e^{p}}\int_{a+2\e}^{a+3\e}
|g(x)|^p\,dx+\frac{1}{\e^{p}}\int_{a+\e}^{a+\e}|g(y)|^p\,dy
\right)\\&\le& 
C\left( \e^{p}\int_a^{a+3\e}|g''(t)|^p\,dt +\frac{1}{\e^{p}}\int_{a}^{a+3\e}
|g(s)|^p\,ds
\right).
\end{eqnarray*}
By taking~$a:=3\e j$ with~$j\in\left\{0,1,\dots, \frac1{3\e}-1\right\}$ we conclude that
\begin{eqnarray*}
\int_0^1 |g'(w)|^p\,dw&=&\sum_{j\in\left\{0,1,\dots, \frac1{3\e}-1\right\}
} \int_{3\e j}^{3\e(j+1)} |g'(w)|^p\\&\le&C\sum_{j\in\left\{0,1,\dots, \frac1{3\e}-1\right\}} 
\left( \e^{p}\int_{3\e j}^{3\e(j+1) }|g''(t)|^p\,dt +\frac{1}{\e^{p}}\int_{3\e j}^{3\e(j+1)} 
|g(s)|^p\,ds
\right)\\
&=&
C\left( \e^{p}\int_0^1|g''(t)|^p\,dt +\frac{1}{\e^{p}}\int_0^1
|g(s)|^p\,ds
\right)
\end{eqnarray*}
and this completes the proof of~\eqref{KS-0ujwfn-092urjnf-i2jrfoneh9843hgb}. 

Now we deal with the core of the proof of Theorem~\ref{982nksjdvnICtsghdfnEJHACANosd90hyri32t8043hger-3506o5y}.
First of all, conveniently extending~$f$ outside~$B_1$, see e.g.~\cite[Lemma~6.37]{MR1814364},
\begin{equation}\label{034widetildeuvarphiu}
{\mbox{we can suppose that~$f\in C^{0,1}_0(B_2)$, with~$\|f\|_{L^p(\R^n)}\le C\|f\|_{L^p(B_1)}$.}}
\end{equation}
Thus, we consider the\index{Newtonian potential} Newtonian potential~$v:=-\Gamma*f$.
We know from~\eqref{SOLPONE-0}
that~$\Delta v=f$ in~$\R^n$. As a result,
defining~$w:=u-v$, we have that~$w$ is harmonic in~$B_1$.
Accordingly, it follows from Cauchy's Estimates (see Theorem~\ref{CAUESTIMTH}) that
$$ \|D^2 w\|_{L^\infty(B_{1/2})}\le C\,\|w\|_{L^1(B_{1})}\le C\Big(\|u\|_{L^1(B_1)}+\|v\|_{L^1(B_1)}\Big).$$
{F}rom this inequality and Corollary~\ref{CXf43Unfsa0-430u4noerngTheoremT:CAZY-PIVO}, recalling~\eqref{034widetildeuvarphiu},
we infer that
\begin{equation}\label{0u32pjt43gnab4lan}
\begin{split} \|D^2 u\|_{L^p(B_{1/2})}\,&\le\,\|D^2 v\|_{L^p(B_{1/2})}+
\|D^2 w\|_{L^p(B_{1/2})} \\
& \le \,C\Big(\|u\|_{L^1(B_1)}+\|v\|_{L^1(B_1)}+\|f\|_{L^p(\R^n)}\Big)\\&\le\,
C\left(\|u\|_{L^1(B_1)}+\iint_{B_1\times B_2} |\Gamma(x-y)|\,|f(y)|\,dx\,dy+\|f\|_{L^p(\R^n)}\right)\\&\le\,
C\left(\|u\|_{L^1(B_1)}+\iint_{B_3\times B_2} |\Gamma(z)|\,|f(y)|\,dz\,dy+\|f\|_{L^p(\R^n)}\right)\\&\le\,
C\Big(\|u\|_{L^1(B_1)}+\|f\|_{L^1(B_2)}+\|f\|_{L^p(\R^n)}\Big)\\&\le\,
C\Big(\|u\|_{L^1(B_1)}+\|f\|_{L^p(B_1)}\Big).
\end{split}\end{equation}
Now, by~\eqref{KS-0ujwfn-092urjnf-i2jrfoneh9843hgb},
$$ \|\nabla u\|_{L^p(B_{1/2})}\le C\left( \|D^2 u\|_{L^p(B_{1/2})}+\|u\|_{L^p(B_{1/2})}\right).$$
This and~\eqref{0u32pjt43gnab4lan}, up to renaming constants, yield that
\begin{eqnarray*}
\|u\|_{W^{2,p}(B_{1/2})}&\le& C\,\Big( \|D^2 u\|_{L^p(B_{1/2})}+\|\nabla u\|_{L^p(B_{1/2})}+\|u\|_{L^p(B_{1/2})}\Big)\\
&\le& C\,\Big( \|D^2 u\|_{L^p(B_{1/2})}+\|u\|_{L^p(B_{1/2})}\Big)
\\&\le& C\,\Big( \|u\|_{L^1(B_1)}+\|f\|_{L^p(B_1)}+\|u\|_{L^p(B_{1/2})}\Big)
\end{eqnarray*}
from which the desired result follows.
\end{proof}

As a technical observation, we point out that the assumption that~$u$ belongs to~$C^2(B_1)$
in the statement of Theorem~\ref{982nksjdvnICtsghdfnEJHACANosd90hyri32t8043hger-3506o5y}
has been taken mainly to have a pointwise definition of the corresponding equation~$\Delta u=f$.
This setting can be weakened in several forms, for instance we point out that
in the statement of Theorem~\ref{982nksjdvnICtsghdfnEJHACANosd90hyri32t8043hger-3506o5y}
one can relax the conditions~$u\in C^2(B_1)$ and~$f\in C^{0,1}(B_1)$ to~$u\in W^{2,p}(B_1)$
and~$f\in L^p(B_1)$ respectively. In this framework, however, the equation is not necessarily satisfied
at every point of~$B_1$, but only a.e. in~$B_1$ (which, being~$u\in W^{2,p}(B_1)$,
also says that the equation provides an identity between functions in~$L^p(B_1)$;
in jargon, these are called ``strong solutions'',
to distinguish with the ``weak solutions'' obtained by integrating by parts as in~\eqref{GRAVSSEFORC2};
the terminology might be confusing anyway, since ``strong'' solutions are ``weaker''
than classical smooth solutions). The counterpart of Theorem~\ref{982nksjdvnICtsghdfnEJHACANosd90hyri32t8043hger-3506o5y}
in this setting goes as follows:

\begin{corollary}\label{THM:0-T536EcdhPOaIJodidpo3tju24ylkiliR2543uj8890j09092po-okf-3edk-3edhf-2-fh564COROLP}
Let~$p\in(1,+\infty)$ and~$f\in L^p(B_1)$. Assume that~$u\in W^{2,p}(B_1)$ is a solution of~$\Delta u=f$ a.e. in~$B_1$.

Then,
\begin{equation*}
\|u\|_{W^{2,p}(B_{1/2})}\le C\,\Big( \|u\|_{L^p(B_1)}+\|f\|_{L^p(B_1)}\Big),\end{equation*}
where~$C>0$ depends only on~$n$ and~$p$.
\end{corollary}

\begin{proof}
We use a mollification argument.
We take~$\tau\in C^\infty_0(B_1,\,[0,+\infty))$
with~$\int_{B_1}\tau(x)\,dx=1$.
Given~$\eta\in\left(0,\frac{1}{10}\right)$, we let~$\tau_\eta(x):=\frac1{\eta^n}
\tau\left(\frac{x}{\eta}\right)$
and define~$u_\eta:=u*\tau_\eta$ and~$f_\eta:=f*\tau_\eta$. We observe that~$u_\eta\to u$ in~$W^{2,p}(B_{3/4})$
and~$f_\eta\to f$ in~$L^p(B_{3/4})$
(see e.g.~\cite[Theorem~9.6]{MR3381284}). 

Also, for all~$\varphi\in C^\infty_0(B_{3/4})$ and~$i\in\{1,\dots,n\}$,
\begin{eqnarray*}&&
\int_{\R^n}\partial_i u_\eta(x)\,\varphi(x)\,dx=
-\int_{\R^n} u_\eta(x)\,\partial_i\varphi(x)\,dx=
-\iint_{\R^n\times\R^n} u(x-y)\,\tau(y)\,\partial_i\varphi(x)\,dx\,dy\\&&\qquad
=\iint_{\R^n\times\R^n} \partial_i u(x-y)\,\tau(y)\,\varphi(x)\,dx\,dy=
\int_{\R^n}(\partial_i u*\tau_\eta)(x)\,\varphi(x)\,dx,
\end{eqnarray*}
that is~$\partial_i u_\eta=\partial_i u*\tau_\eta$ in~$B_{3/4}$
(notice that this holds pointwise at every point of~$B_{3/4}$, since~$u_\eta\in C^\infty(B_{3/4})$). Iterating this information, we find that
\begin{equation*}
\Delta u_\eta=\Delta u*\tau_\eta=f*\tau_\eta=f_\eta\quad{\mbox{ in }}\,B_{3/4}.\end{equation*}
As a result, we can invoke Theorem~\ref{982nksjdvnICtsghdfnEJHACANosd90hyri32t8043hger-3506o5y}
(applied here with~$B_{3/4}$ instead of~$B_1$ and renaming constants accordingly), finding that
$$ \|u_\eta\|_{W^{2,p}(B_{1/2})}\le C\,\Big( \|u_\eta\|_{L^p(B_{3/4})}+\|f_\eta\|_{L^p(B_{3/4})}\Big).$$
The desired result thus follows by taking the limit as~$\eta\searrow0$.
\end{proof}

\section{Calder\'{o}n-Zygmund estimates for equations
in nondivergence form}\label{SEC:CaldZygmund estimates}

As discussed in Section~\ref{KPJMDarJMSiaCLiMMSciAmoenbvD8Sijf-23},
once we established a nice regularity estimate for the Poisson equation, it becomes very desirable to
extend it to a more general class of elliptic equations belonging to some ``natural class''
which can be treated as a local perturbation of the Laplace operator. In this setting,
we have a natural counterpart of the interior estimates in Corollary~\ref{THM:0-T536EcdhPOaIJodidpo3tju24ylkiliR2543uj8890j09092po-okf-3edk-3edhf-2-fh564COROLP}
for general elliptic equations, in a form similar to that of
Theorem~\ref{THM:0-T536EcdhPOaIJodidpo3tju24ylkiliR2543uj8890j09092po-okf-3edk-3edhf-2-fh564},
designed here to address the~$L^p$ theory.

\begin{theorem}\label{LSmvvikSDnkjazikgAzldDrgmIDMmo5r523}
Let~$a_{ij}\in C(B_1)$ satisfy the ellipticity 
condition in~\eqref{ELLIPTIC}. Let~$b_i$, $c\in L^\infty(B_1)$. Let~$f\in L^p(B_1)$ for some~$p\in(1,+\infty)$.
Let~$u\in W^{2,p}(B_1)$ be a solution of
\begin{equation*}
\sum_{i,j=1}^n a_{ij} \partial_{ij}u+\sum_{i=1}^n b_i\partial_i u+cu=f
\quad{\mbox{a.e. in }}\,B_1.\end{equation*}
Then, there exists~$C>0$, depending only on~$n$,
$p$, $a_{ij}$, $b_i$ and~$c$, such that
\begin{equation*} \|u\|_{W^{2,p}(B_{1/2})}\le C\,\Big(\|u\|_{L^p(B_1)}+\|f\|_{L^p(B_1)}\Big).\end{equation*}
\end{theorem}

The proof of Theorem~\ref{LSmvvikSDnkjazikgAzldDrgmIDMmo5r523}
leverages a result for scaled estimates which can be seen as the counterpart
of Proposition~\ref{SVA-LS-LEOSKDNDO} in the~$L^p$ framework.
This auxiliary result is the following:

\begin{proposition}\label{LPSVA-LS-LEOSKDNDO}
Let~$k\in\N$, $p\in[1,+\infty)$, $\sigma\ge0$ and $R>0$. 
There exists~$\e_0\in(0,1)$ depending only on~$n$ and~$k$
such that the following holds true.

Let~$u\in W^{k,p}(B_R)$.
Assume that
for every~$x_0\in B_{R}$ and every~$\rho\in\left( 0,\frac{R-|x_0|}{2}\right) $
\begin{equation}\label{rEGa}
\sum_{j=0}^k \rho^j
\|D^j u\|_{L^{p}(B_\rho(x_0))}
\le \sigma+\e_0 \,\sum_{j=0}^k (2\rho)^j
\|D^j u\|_{L^{p}(B_{2\rho}(x_0))}.
\end{equation}
Then, there exists~$C>0$, depending only on~$n$, $k$ and~$R$
such that
\begin{equation*}\|u\|_{W^{k,p}(B_{R/2})}
\le C\sigma.\end{equation*}
\end{proposition}

\begin{proof}
Let
\begin{equation*} Q:=\sup_{ {x_0\in {{B}}_R}\atop{
\rho\in\left( 0,\frac{R-|x_0|}2\right)} } \sum_{j=0}^k \rho^j
\|D^j u\|_{L^{p}(B_\rho(x_0))}\end{equation*}
and notice that
$$Q\le (k+1)(1+R)^k\,\|u\|_{W^{k,p}(B_{R})}
<+\infty.$$
Now we take~$x_0\in B_R$ and~$\rho\in \left( 0,\frac{R-|x_0|}2\right)$.
We cover~$B_\rho(x_0)$ by a finite family of balls~$B_{\frac{\rho}{4}}(p_i)$
with~$p_i\in B_\rho(x_0)$ and~$i\in\{1,\dots,N\}$, where~$N$ is bounded from above
in dependence only of~$n$. 

Since~$|p_i|\le|x_0|+\rho<R-2\rho+\rho=R-\rho$, we have that~$p_i\in B_R$
and~$\frac{\rho}2\in\left(0,\frac{R-|p_i|}2\right)$. Thus,
for all~$i\in\{1,\dots,N\}$,
\begin{equation*}
\sum_{j=0}^k \left(\frac{\rho}4\right)^j
\|D^j u\|_{L^{p}(B_{\rho/4}(p_i))}
\le \sigma+\e_0 
\sum_{j=0}^k \left(\frac{\rho}2\right)^j
\|D^j u\|_{L^{p}(B_{\rho/2}(p_i))}
\le \sigma+\e_0 Q,
\end{equation*}
thanks to~\eqref{rEGa}.
Therefore,
\begin{equation*}
\sum_{j=0}^k\rho^j
\|D^j u\|_{L^{p}(B_{\rho}(x_0))}
\le
\sum_{j=0}^k 4^j \sum_{i=1}^N\left(\frac{\rho}4\right)^j
\|D^j u\|_{L^{p}(B_{\rho/4}(p_i))}
\le 4^k N\sigma+\e_0 4^kNQ\big)
.\end{equation*}
Taking the supremum over such~$x_0$ and~$\rho$,
we find that
\begin{equation*}
Q\le 4^k N\sigma+\e_0 4^k NQ.
\end{equation*}
Hence, choosing~$\e_0:=\frac{1}{2^{2k+1}N }$ and reabsorbing one term into the left hand side,
\begin{equation*} Q\le2^{2k+1} N \sigma .\end{equation*}
Since
$$Q\ge\lim_{\rho\nearrow R/2} \sum_{j=0}^k \rho^j
\|D^j u\|_{L^{p}(B_{\rho})}=\sum_{j=0}^k \left(\frac{R}2\right)^j
\|D^j u\|_{L^{p}(B_{R/2})}
,$$
the desired result follows.\end{proof}

With this, we can address the proof of Theorem~\ref{LSmvvikSDnkjazikgAzldDrgmIDMmo5r523}.
The argument is a modification of that presented in the proof of Theorem~\ref{THM:0-T536EcdhPOaIJodidpo3tju24ylkiliR2543uj8890j09092po-okf-3edk-3edhf-2-fh564}
(here, we will invoke Corollary~\ref{THM:0-T536EcdhPOaIJodidpo3tju24ylkiliR2543uj8890j09092po-okf-3edk-3edhf-2-fh564COROLP}
instead of Theorem~\ref{SCHAUDER-INTE}
and Proposition~\ref{LPSVA-LS-LEOSKDNDO}
instead of Proposition~\ref{SVA-LS-LEOSKDNDO})
and it utilizes the fact that, zooming in close to a point, the picture reduces to a small perturbation of the Laplace operator.

\begin{proof}[Proof of Theorem~\ref{LSmvvikSDnkjazikgAzldDrgmIDMmo5r523}]
As observed in~\eqref{021irkfe2m4m53a0294OVERINGAR3G}, up to an affine transformation,
Corollary~\ref{THM:0-T536EcdhPOaIJodidpo3tju24ylkiliR2543uj8890j09092po-okf-3edk-3edhf-2-fh564COROLP}
holds true if the Laplacian is replaced by an operator with constant
coefficients of the form~$\sum_{i,j=1}^n \overline{a}_{ij} \partial_{ij}u$
where~$\overline{a}_{ij}$ is constant
and fulfills the ellipticity condition in~\eqref{ELLIPTIC}.

Given~$x_0\in B_{1/2}$ and~$\rho\in\left(0,\frac14\right)$, to be taken suitably small in the following,
we set \begin{eqnarray*}&& \widetilde{a}_{ij}(x):=
a_{ij}(x_0)-a_{ij}(x),\\&&\widetilde{u}(x):=\frac{u(x_0+\rho x)}{\rho^2}\\{\mbox{and }}&& 
\widetilde{f}(x):=f(x_0+\rho x)-\sum_{i=1}^n b_i(x_0+\rho x)\partial_iu(x_0+\rho x)\\&&\qquad\qquad\qquad-c(x_0+\rho x)u(x_0+\rho x)+
\sum_{i,j=1}^n \widetilde{a}_{ij}(x_0+\rho x)\partial_{ij}u(x_0+\rho x)
.\end{eqnarray*}
By~\eqref{hinHNDfOIKHDNinEDFDDIAbgerIGGjrtNVMVFkqTESKMFDFRT}, a.e.~$x\in B_1$ we have that
\begin{eqnarray*}
&&\sum_{i,j=1}^n a_{ij}(x_0)\partial_{ij}\widetilde{u}(x)
=\widetilde f(x).
\end{eqnarray*}
Thus, by Corollary~\ref{THM:0-T536EcdhPOaIJodidpo3tju24ylkiliR2543uj8890j09092po-okf-3edk-3edhf-2-fh564COROLP},
\begin{equation}\label{RESM9ihw320D0-2i5rkwfracrhoightinLP}
\|\widetilde u\|_{W^{2,p}(B_{1/2})}\le C\,\Big(\|\widetilde u\|_{L^p(B_1)}+\|\widetilde f\|_{L^p(B_1)}\Big).\end{equation}
Now we take~$\delta>0$, to be chosen conveniently small here below, and we point out that
$$ \sup_{x\in B_1}|\widetilde{a}_{ij}(x_0+\rho x)|=\sup_{y\in B_\rho(x_0)}|\widetilde{a}_{ij}(y)|=
\sup_{y\in B_\rho(x_0)}|a_{ij}(x_0)-a_{ij}(y)|\le\delta,$$
as long as~$\rho\in(0,\delta)$ is small enough in dependence of~$\delta$, thanks to the continuity assumption on the coefficients~$a_{ij}$.

As a result, for a.e.~$x\in B_1$,
\begin{eqnarray*}
|\widetilde{f}(x)|&\le&C\,\Big(|f(x_0+\rho x)|+|\nabla u(x_0+\rho x)|+|u(x_0+\rho x)|+
\delta|D^2 u(x_0+\rho x)|\Big),
\end{eqnarray*}
leading to
\begin{eqnarray*}
\|\widetilde f\|_{L^p(B_1)}&\le&C\,\left(
\int_{B_1}\Big(
|f(x_0+\rho x)|^p+|\nabla u(x_0+\rho x)|^p+|u(x_0+\rho x)|^p+
\delta^p|D^2 u(x_0+\rho x)|^p\Big)\,dx
\right)^{\frac1p}\\&=&\frac{C}{\rho^{\frac{n}p}}\left(
\int_{B_\rho(x_0)}\Big(
|f(y)|^p+|\nabla u(y)|^p+|u(y)|^p+
\delta^p|D^2 u(y)|^p\Big)\,dy
\right)^{\frac1p}\\&\le&\frac{C}{\rho^{\frac{n}p}}\Big(\|f\|_{L^p(B_\rho(x_0))}+
\|\nabla u\|_{L^p(B_\rho(x_0))}+\|u\|_{L^p(B_\rho(x_0))}+\delta\|D^2u\|_{L^p(B_\rho(x_0))}
\Big).
\end{eqnarray*}
For this reason and~\eqref{RESM9ihw320D0-2i5rkwfracrhoightinLP},
\begin{equation}\label{SM9ihw320D0-2i5rkwfracrhoightinLP}\begin{split}&
\rho^{\frac{n}p}\|\widetilde u\|_{W^{2,p}(B_{1/2})}\\&\quad\le C\,\left[\rho^{\frac{n}p}\|\widetilde u\|_{L^p(B_1)}+
\|f\|_{L^p(B_\rho(x_0))}+
\|\nabla u\|_{L^p(B_\rho(x_0))}+\|u\|_{L^p(B_\rho(x_0))}+\delta\|D^2u\|_{L^p(B_\rho(x_0))}
\right].\end{split}\end{equation}
We also observe that, for every~$j\in\N$, we have that~$D^j u(y)=\rho^{2-j} D^j \widetilde{u}\left(\frac{y-x_0}\rho\right)$ hence
\begin{eqnarray*}&& \|D^j u\|_{L^p(B_r(x_0))}
=\rho^{2-j}\left(\int_{B_r(x_0)} \left|D^j \widetilde{u}\left(\frac{y-x_0}\rho\right)\right|^p\,dy\right)^{\frac1p}\\&&\qquad\quad=
\rho^{\frac{n}p+2-j}\left(\int_{B_{r/\rho}} \left|D^j \widetilde{u}(x)\right|^p\,dx\right)^{\frac1p}=\rho^{\frac{n}p+2-j}\, 
\|D^j\widetilde u\|_{L^p(B_{r/\rho})}
\end{eqnarray*}
and therefore~\eqref{SM9ihw320D0-2i5rkwfracrhoightinLP} becomes
\begin{equation*}\begin{split}&
\sum_{j=0}^2\rho^{j}\|D^j u\|_{L^p(B_{\rho/2}(x_0))}\\&\quad=\rho^{\frac{n}p+2}\sum_{j=0}^2
\|D^j\widetilde u\|_{L^p(B_{1/2})}\\&\quad\le C\,\Big[
\| u\|_{L^p(B_{\rho}(x_0))}+\rho^2\|f\|_{L^p(B_\rho(x_0))}+\rho^2
\|\nabla u\|_{L^p(B_\rho(x_0))}+\rho^2\|u\|_{L^p(B_\rho(x_0))}+\delta\rho^2\|D^2u\|_{L^p(B_\rho(x_0))}
\Big]\\&\quad\le C\,\Big[
\| u\|_{L^p(B_{1})}+\|f\|_{L^p(B_1)}+\delta\rho
\|\nabla u\|_{L^p(B_\rho(x_0))}+\delta\|u\|_{L^p(B_\rho(x_0))}+\delta\rho^2\|D^2u\|_{L^p(B_\rho(x_0))}
\Big]\\
&\qquad=C\,\Big[
\| u\|_{L^p(B_{1})}+\|f\|_{L^p(B_1)}\Big]+C\delta\sum_{j=0}^2\rho^{j}\|D^j u\|_{L^p(B_{\rho}(x_0))}.
\end{split}\end{equation*}
{F}rom this and Proposition~\ref{LPSVA-LS-LEOSKDNDO} (used here
with~$\sigma:=C\big(\| u\|_{L^p(B_{1})}+\|f\|_{L^p(B_1)}\big)$ and~$\e_0:=C\delta$)
, we obtain the desired result.
\end{proof}

We mention that
the $L^p$ theory exposed here is flexible enough to address also the local estimates at
the boundary, thus providing the $L^p$ counterpart of Theorems~\ref{SCHAUDER-BOUTH}
and~\ref{THM:0-T536EcdhPOaIJodidpo3tju24ylkiliR2543uj8890j09092po-okf-3edk-3edhf-2-fh564-LEATB}.
In this setting, one has:

\begin{theorem}\label{JSP0wir-2gjw0egoewnbh30ehf8uy3wgbBELOCVC}
Let the halfball notation in~\eqref{9oj3falphaB10} hold true.
Let~$a_{ij}\in C(\overline{B_1^+})$ satisfy the ellipticity 
condition in~\eqref{ELLIPTIC}. Let~$b_i$, $c\in L^\infty(B_1^+)$. Let~$f\in L^p(B_1^+)$ and~$g\in
W^{2,p}(B_1^+)$
for some~$p\in(1,+\infty)$.
Let~$u\in W^{2,p}(B_1^+)$ be a solution of
\begin{equation}\label{MD-IUD-Xrkhg9}
\begin{dcases} \sum_{i,j=1}^n a_{ij} \partial_{ij}u+\sum_{i=1}^n b_i\partial_i u+cu=f\quad{\mbox{a.e. in }}\,B_1^+,\\
u=g \quad{\mbox{on }}\,B_1^0.\end{dcases}\end{equation}
Then, there exists~$C>0$, depending only on~$n$,
$p$, $a_{ij}$, $b_i$ and~$c$, such that
\begin{equation}\label{0ojlqd0eo-2fkmv-3fgmvws0fPJON} \|u\|_{W^{2,p}(B_{1/2}^+)}\le C\,\Big(\|u\|_{L^p(B_1^+)}+
\|g\|_{W^{2,p}(B_1^+)}+\|f\|_{L^p(B_1^+)}\Big).\end{equation}
\end{theorem}

As usual in the Sobolev space setting, the condition~$u=g$ in~\eqref{MD-IUD-Xrkhg9}
is intended in the sense of traces for~$W^{1,p}$ functions (see e.g.~\cite[Chapter~18]{MR2527916}
for a thorough discussion of the trace operator).

\begin{proof}[Proof of Theorem~\ref{JSP0wir-2gjw0egoewnbh30ehf8uy3wgbBELOCVC}] By possibly replacing~$u$ by~$u-g$, 
\begin{equation}\label{FGANFGGNRSIYTYTIOTGENAEVFK0}
{\mbox{we can reduce ourselves to the case in which~$g$ vanishes identically.}}\end{equation}
Also, 
\begin{equation}\label{FGANFGGNRSIYTYTIOTGENAEVFK}\begin{split}&
{\mbox{it is enough to consider the case in which~$a_{ij}=\delta_{ij}$, $b_i=0$ and~$c=0$,}}\\
&{\mbox{that is
the case of the Laplace operator.}}\end{split}\end{equation}
Indeed, if~\eqref{0ojlqd0eo-2fkmv-3fgmvws0fPJON} is proved
for the Laplace operator (and therefore for operators
with constant coefficients, up to an affine transformation)
the general case follows by perturbing the coefficients and reabsorbing the
lower order terms, precisely as done in the
proof of Theorem~\ref{LSmvvikSDnkjazikgAzldDrgmIDMmo5r523}
(but using here halfballs instead of balls).

Hence, we work under the simplifying assumptions in~\eqref{FGANFGGNRSIYTYTIOTGENAEVFK0}
and~\eqref{FGANFGGNRSIYTYTIOTGENAEVFK}.
In this setting, we can proceed with a reflection and mollification argument.
Let~$\varphi\in C^\infty_0(B_{7/8},[0,1])$ with~$\varphi=1$ in~$B_{3/4}$. Let also~$U:=u\varphi\in W^{2,p}(B_1^+)$
and notice that~$\Delta U=f\varphi+2\nabla u\cdot\nabla\varphi+u\Delta\varphi=:F$ a.e. in~$B_1^+$.
Moreover, we extend~$u$ and~$U$ to the whole of~$B_1$ by odd reflection, setting~$u(x',x_n):=-u(x',-x_n)$ and~$U(x',x_n):=-U(x',-x_n)$ for all~$x=(x',x_n)\in B_1$ with~$x_n<0$, and also~$f$ and~$F$
as done in~\eqref{fstarywueut878788}.
Given~$\e>0$,
we note that~$U\in W^{2,p}(B_1\cap\{|x_n|>\e\})$ and~$\Delta U=F$ a.e. in~$
B_1\cap\{|x_n|>\e\}$ (hence also in the $L^p$ sense) and therefore
\begin{equation}\label{CNSJPDMONAMSVTKBBAT}
\int_{ B_1\cap\{|x_n|>\e\} }F(x)\Phi(x)\,dx=
-\int_{ B_1\cap\{|x_n|>\e\} }\nabla U(x)\cdot\nabla\Phi(x)\,dx\end{equation}
for every~$\Phi\in C^\infty_0(B_1\cap\{|x_n|>\e\})$.

Now we take~$\psi\in C^\infty_0(B_1)$ and~$\zeta\in C^\infty(\R,[0,1])$ be even with~$\zeta=0$ in~$[-1,1]$
and~$\zeta=1$ in~$\R\setminus(-2,2)$. Let~$\zeta_\e(x):=\zeta\left(\frac{x_n}{\e}\right)$.
Using~\eqref{CNSJPDMONAMSVTKBBAT} with~$\Phi:=\zeta_\e \psi$,
we notice that
\begin{eqnarray*}
&&\left|\int_{\R^n} F(x)\psi(x)\,dx+\int_{\R^n}\nabla U(x)\cdot\nabla \psi(x)\,dx\right|
\\&\le&
\left|\int_{\R^n} F(x)\zeta_\e(x)\psi(x)\,dx+\int_{\R^n}\nabla U(x)\cdot\nabla \big(\zeta_\e(x)\psi(x)\big)\,dx\right|
+
\int_{\R^n} \big|1-\zeta_\e(x)\big|\,|F(x)|\,|\psi(x)|\,dx\\&&\qquad
+
\left|\int_{\R^n}\nabla U(x)\cdot\nabla \Big(\big(1-\zeta_\e(x)\big)\psi(x) \Big)\,dx\right|
\\&\le&0+
\int_{\R^n\cap\{|x_n|<2\e\}} |F(x)|\,|\psi(x)|\,dx
+
\left|\int_{\R^n\cap\{|x_n|<2\e\}}\nabla U(x)\cdot\nabla \Big(\big(1-\zeta_\e(x)\big)\psi(x) \Big)\,dx\right|\\&\le&
C\|F\|_{L^p(\{|x_n|<2\e\})}
+\int_{\R^n\cap\{ |x_n|<2\e\}}|\nabla U(x)|\,|\nabla\psi(x)|\,dx \\&&\qquad+
\left|\int_{\R^n\cap\{|x_n|<2\e\}}\partial_n U(x)\partial_n \zeta_\e(x)\psi(x)\,dx\right|\\&\le&C\|F\|_{L^p(\{|x_n|<2\e\})}
+C\|\nabla U\|_{L^p(\{|x_n|<2\e\})}
\\&&\qquad+
\left|\int_{\R^n\cap\{x_n\in(0,2\e)\}}\partial_n U(x)\partial_n \zeta_\e(x)\psi(x)\,dx+
\int_{\R^n\cap\{x_n\in(-2\e,0)\}}\partial_n U(x',-x_n)\partial_n \zeta_\e(x)\psi(x)\,dx\right|.
\end{eqnarray*}
Taking the limit as~$\e\searrow0$ and noticing that
$$ \partial_n \zeta_\e(x)=\frac1\e\zeta'\left(\frac{x_n}\e\right),$$
recalling also that~$\zeta$ is even, whence~$\zeta'$ is odd,
we find that
\begin{eqnarray*}
&&\left|\int_{\R^n} F(x)\psi(x)\,dx+\int_{\R^n}\nabla U(x)\cdot\nabla \psi(x)\,dx\right|
\\&\le&\lim_{\e\searrow0}\frac1\e
\left|\int_{\R^n\cap\{x_n\in(0,2\e)\}}\partial_n U(x) \zeta'\left(\frac{x_n}\e\right)\psi(x)\,dx\right.\\&&\qquad\left.-
\int_{\R^n\cap\{x_n\in(0,2\e)\}}\partial_n U(x',x_n) \zeta'\left(\frac{x_n}\e\right)\psi(x',-x_n)\,dx\right|\\
&\le&
\lim_{\e\searrow0}\frac{C}\e
\int_{\R^n\cap\{x_n\in(0,2\e)\}}|\partial_n U(x)|\,\big|\psi(x',x_n)-\psi(x',-x_n)\big|\,dx
\\&\le&\lim_{\e\searrow0}\int_{\R^n\cap\{x_n\in(0,2\e)\}}|\partial_n U(x)|\,dx
\\&\le&\lim_{\e\searrow0}\| \nabla U\|_{L^p(\R^n\cap\{x_n\in(0,2\e)\})}\\&=&0.
\end{eqnarray*}
As a consequence,
\begin{equation}\label{ITERSNBS-SMSPS0284ujf131M}
\int_{\R^n} F(x)\psi(x)\,dx=-\int_{\R^n}\nabla U(x)\cdot\nabla \psi(x)\,dx.
\end{equation}
Given~$\eta>0$, we take a mollifier~$\tau_\eta$ and the corresponding
mollification~$U_\eta:=U*\tau_\eta\in C^\infty_0(B_{8/9})$ of~$U$.
Let also~$F_\eta:=F*\tau_\eta$.
By~\eqref{ITERSNBS-SMSPS0284ujf131M}, for all~$\psi\in C^\infty_0(B_{7/8})$,
letting~$\psi_\eta:=\psi*\tau_\eta$ (which is supported in~$B_1$ for~$\eta$ small enough), we have that
\begin{eqnarray*}
&&\int_{\R^n} F_\eta(x)\psi(x)\,dx=
\iint_{\R^n\times\R^n} F(y)\tau_\eta(x-y)\psi(x)\,dx\,dy
=\int_{\R^n} F(y)\psi_\eta(y)\,dx\\&&\qquad=-\int_{\R^n}\nabla U(y)\cdot\nabla \psi_\eta(y)\,dy
=-\int_{\R^n}\nabla U_\eta(x)\cdot\nabla \psi(x)\,dx.
\end{eqnarray*}
{F}rom this and the fact that~$U_\eta$ is smooth we arrive at~$\Delta U_\eta=F_\eta$ in~$B_{7/8}$.
We can therefore invoke the interior regularity estimate
in Theorem~\ref{LSmvvikSDnkjazikgAzldDrgmIDMmo5r523}
(or even the simpler formulation given in 
Corollary~\ref{THM:0-T536EcdhPOaIJodidpo3tju24ylkiliR2543uj8890j09092po-okf-3edk-3edhf-2-fh564COROLP}) and deduce that
\begin{equation}\label{JSONJNEAikvCKmfkCANpuKblY2}
\|U_\eta\|_{W^{2,p}(B_{1/2})}\le C\,\Big( \|U_\eta\|_{L^p(B_{7/8})}+\|F_\eta\|_{L^p(B_{7/8})}\Big).\end{equation}
Also, since~$\Delta(U_{\eta}-U_{\eta'})=F_\eta-F_{\eta'}$, we obtain in the same way that
\begin{equation}\label{JSONJNEAikvCKmfkCANpuKblY}
\|U_\eta-U_{\eta'}\|_{W^{2,p}(B_{1/2})}\le C\,\Big( \|U_\eta-U_{\eta'}\|_{L^p(B_{7/8})}+\|F_\eta-F_{\eta'}\|_{L^p(B_{7/8})}\Big).\end{equation}
Since~$U_\eta\to U$ and~$F_\eta\to F$ in~$L^p(B_{7/8})$
(see e.g.~\cite[Theorem 9.6]{MR3381284}), we deduce from~\eqref{JSONJNEAikvCKmfkCANpuKblY} that~$U_\eta\to U$ in~$W^{2,p}(B_{1/2})$. This and~\eqref{JSONJNEAikvCKmfkCANpuKblY2} yield that
\begin{equation}\label{JSONJNEAikvCKmfkCANpuKblY3}
\|U\|_{W^{2,p}(B_{1/2})}\le C\,\Big( \|U\|_{L^p(B_{7/8})}+\|F\|_{L^p(B_{7/8})}\Big).\end{equation}
Observe in addition that~$U=u$ in~$B_{1/2}$ and~$|U|\le|u|$. Accordingly,
we deduce from~\eqref{JSONJNEAikvCKmfkCANpuKblY3} that
\begin{eqnarray*}
\|u\|_{W^{2,p}(B_{1/2})}&\le& C\,\Big( \|u\|_{L^p(B_{7/8})}+\|F\|_{L^p(B_{7/8})}\Big)\\
&\le&C\,\Big( \|u\|_{L^p(B_{7/8})}+\|f\|_{L^p(B_{7/8})}+\|\nabla u\|_{L^p(B_{7/8})}\Big).
\end{eqnarray*}
Consequently, by interpolation (recall~\eqref{KS-0ujwfn-092urjnf-i2jrfoneh9843hgb}),
$$ \|u\|_{W^{2,p}(B_{1/2})}
\le C\,\Big( \|u\|_{L^p(B_{7/8})}+\|f\|_{L^p(B_{7/8})}\Big)+
\e_0 \|D^2 u\|_{L^p(B_{7/8})},$$
with~$\e_0$ as in Proposition~\ref{LPSVA-LS-LEOSKDNDO}.
By scaling and using Proposition~\ref{LPSVA-LS-LEOSKDNDO},
we conclude that
\[ \|u\|_{W^{2,p}(B_{1/2})}\le C\,\Big(\|u\|_{L^p(B_1)}+
+\|f\|_{L^p(B_1)}\Big),\]
which yields the desired result in~\eqref{0ojlqd0eo-2fkmv-3fgmvws0fPJON}
in this case. 

This completes the proof of Theorem~\ref{JSP0wir-2gjw0egoewnbh30ehf8uy3wgbBELOCVC}
when the operator is the Laplacian and~$g$ vanishes identically
(and therefore 
in its full generality, thanks to~\eqref{FGANFGGNRSIYTYTIOTGENAEVFK0}
and~\eqref{FGANFGGNRSIYTYTIOTGENAEVFK}).
\end{proof}

\begin{figure}
  \centering
  \includegraphics[width=.45\linewidth]{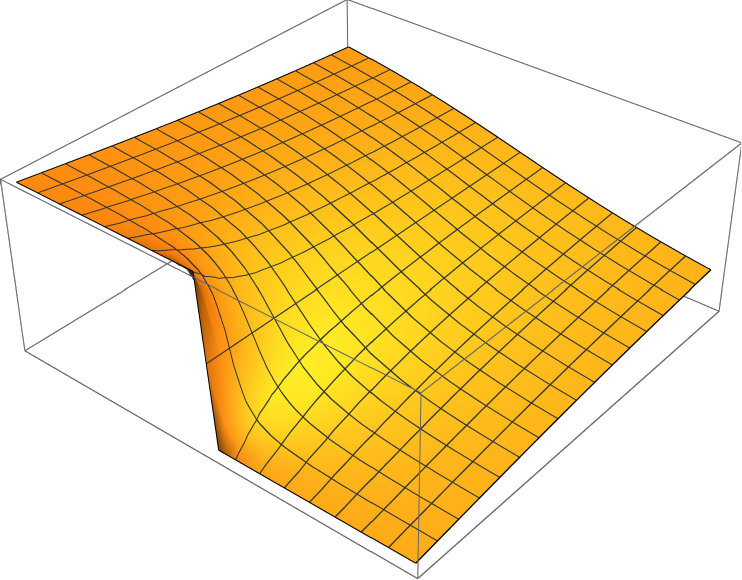}$\qquad$\includegraphics[width=.45\linewidth]{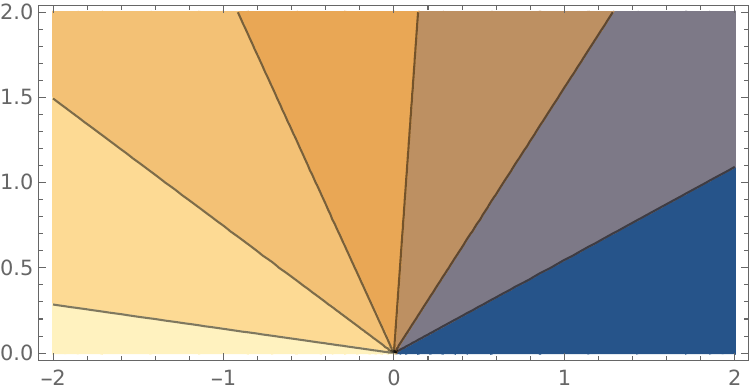}
 \caption{\sl Graph of the function~$u$ in~\eqref{ukAKMS dplwjmfno-kDItangg9845ghjgeFI312} and of its level sets.}\label{ukAKMS dplwjmfno-kDItangg9845ghjgeFI31}
\end{figure}

We remark that
when~$g$ in~\eqref{MD-IUD-Xrkhg9}
is not sufficiently regular, then the local estimate at the boundary in the setting
of Theorem~\ref{JSP0wir-2gjw0egoewnbh30ehf8uy3wgbBELOCVC} may fail.
As an example, one can consider the angular function in the halfplane~$\R^2_+:=\R\times(0,+\infty)$ given by
\begin{equation}\label{ukAKMS dplwjmfno-kDItangg9845ghjgeFI312} u(x)=u(x_1,x_2):=\arccos\frac{x_1}{\sqrt{x_1^2+x_2^2}},
\end{equation}
see Figure~\ref{ukAKMS dplwjmfno-kDItangg9845ghjgeFI31}.

We note that~$u\in C^2(\R^2_+)$ and~$\Delta u=0$ in~$\R^2_+$.
%%%% Moreover,
%%%% $$ \nabla u(x)=\left(-\frac{|x_2|}{|x|^2},\,
%%%% \frac{x_1 x_2}{|x_2|\,|x|^2}\right),$$
%%%% and therefore
%%%% $$ |\nabla u(x)|=\frac{1}{|x|^2}\sqrt{
%%%% x_2^2+x_1^2}=\frac{1}{|x|}.$$
%%%% This says that if~$p\in(1,2)$ then
%%%% $$ \int_{B_1^+}|\nabla u(x)|^p\,dx
%%%% =\pi\int_0^1 \rho^{1-p}\,d\rho=
%%%% \frac{\pi}{2-p}
%%%% <+\infty,$$
%%%% hence~$u\in W^{1,p}(B_1^+)$.

\begin{figure}
  \centering
  \includegraphics[width=.45\linewidth]{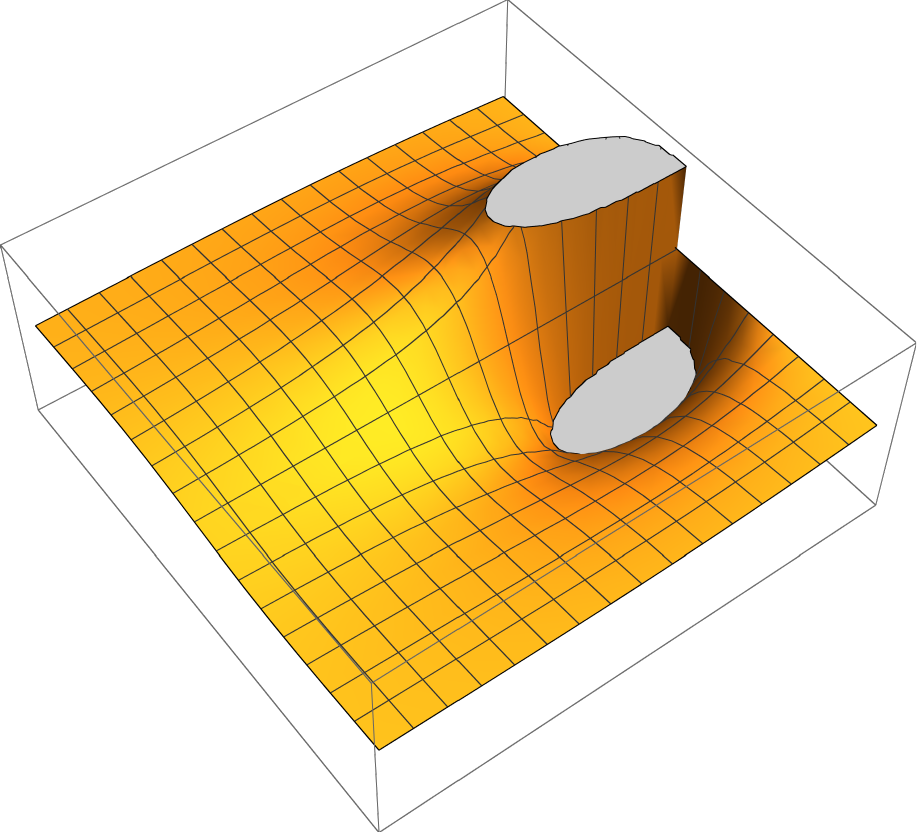}$\qquad$\includegraphics[width=.45\linewidth]{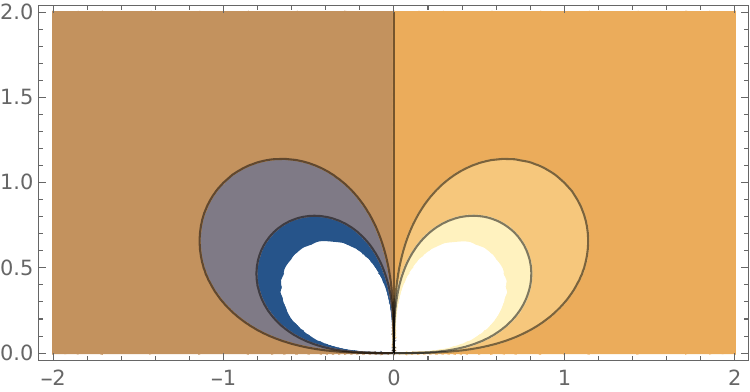}
 \caption{\sl Graph of the function~$\partial_{11}u$
 and of its level sets, with~$u$ as in~\eqref{ukAKMS dplwjmfno-kDItangg9845ghjgeFI312}.}\label{ukAKMS dplwjmfno-kDItangg9845ghjgeFI31TT}
\end{figure}

Additionally,
$$ \lim_{x_2\searrow0}u(x)=
\arccos\frac{x_1}{|x_1|}=\begin{dcases}
0 & {\mbox{ if }}x_1>0,\\
\pi& {\mbox{ if }}x_1<0.
\end{dcases}$$
But the estimate in the thesis of Theorem~\ref{JSP0wir-2gjw0egoewnbh30ehf8uy3wgbBELOCVC} does not hold, since,
for some~$c>0$,
\begin{eqnarray*}&&
\|u\|_{W^{2,p}(B_{1/2}^+)}^p\ge
\int_{B_{1/2}^+}|\partial_{11}u(x)|^p\,dx
=\int_{B_{1/2}^+}\left|
\frac{2 x_1 \,|x_2|}{|x|^4}
\right|^p\,dx=
2^p\int_{B_{1/2}^+}
\frac{|x_1|^p \,|x_2|^p}{|x|^{4p}}
\,dx\\
&&\qquad=2^{p-1}\int_{B_{1/2}}
\frac{|x_1|^p \,|x_2|^p}{|x|^{4p}}
\,dx=2^{p-1}\iint_{(0,1/2)\times(\partial B_1)}\rho^{1-2p}
\,|\omega_1|^p \,|\omega_2|^p
\,d\rho\,d{\mathcal{H}}^1_\omega\\
&&\qquad=c\int_0^{1/2}\rho^{1-2p}\,d\rho
\end{eqnarray*}
which diverges for all~$p\in(1,+\infty)$.
See Figure~\ref{ukAKMS dplwjmfno-kDItangg9845ghjgeFI31TT}
for a plot of~$\partial_{11}u$ and its level sets.\medskip

We now recall a global, up to the boundary, version of the regularity
estimates in~$L^p$ class that follows from interior and local boundary estimates.
This can be seen as the natural counterpart in Lebesgue spaces of
Theorems~\ref{KS-2rkjoewuj5-2prlfgk-a} and~\ref{ThneGiunfBB032t4jNKS-s3i4}.

\begin{theorem}\label{JOSLNGEITHSHLPESTGLOUPTHBUO}
Let~$\Omega\subset\R^n$ be a bounded open set with boundary of class~$C^{1,1}$.
Let~$a_{ij}\in C(\overline{\Omega})$ satisfy the ellipticity 
condition in~\eqref{ELLIPTIC}. Let~$b_i$, $c\in L^\infty(\Omega)$. Let~$f\in L^p(\Omega)$ and~$g\in
W^{2,p}(\Omega)$
for some~$p\in(1,+\infty)$.

Let~$u\in W^{2,p}(\Omega)$ be a solution of
\begin{equation*}
\begin{dcases}
\sum_{i,j=1}^n a_{ij} \partial_{ij}u+\sum_{i=1}^n b_i\partial_i u+cu=f
&\quad{\mbox{a.e. in }}\,\Omega,\\
u=g&\quad{\mbox{on }}\,\partial\Omega.\end{dcases}
\end{equation*}
Then, there exists~$C>0$, depending only\footnote{Actually, browsing through the proof, one can be more
precise by saying that the dependence of~$C$ here is on~$n$,
$p$, the~$C^{1,1}$ regularity parameters of~$\Omega$, the modulus of continuity of~$a_{ij}$
and its ellipticity constants, $\|b_i\|_{L^\infty(\Omega)}$ and~$\|c\|_{L^\infty(\Omega)}$. \label{j0o2r5Omegaand}} on~$n$,
$p$, $\Omega$, $a_{ij}$, $b_i$ and~$c$, such that
\begin{equation}\label{ojhnwfed-2irjgTYHSJDISA2R44G78jS} \|u\|_{W^{2,p}(\Omega)}\le C\,\Big(\|g\|_{W^{2,p}(\Omega)}+\|f\|_{L^p(\Omega)}+\|u\|_{L^p(\Omega)}\Big).\end{equation}
\end{theorem}

\begin{proof} This argument is a variation of that presented in the proof of Theorem~\ref{ThneGiunfBB032t4jNKS-s3i4},
using the regularity results in Lebesgue spaces in place of those in H\"older spaces.
Namely, we straighten the boundary of~$\Omega$ by taking~$\rho>0$ a finite family
of balls~$\{B_\rho(p_i)\}_{i\in\{1,\dots,N\}}$ with~$p_i\in\partial\Omega$ such that
\begin{equation*}
{\mathcal{N}}:=\{ x\in\Omega {\mbox{ s.t. }} B_{\rho/8}(x)\cap(\partial\Omega)\ne\varnothing\}
\subseteq \bigcup_{i=1}^N B_{\rho/4}(p_i)\end{equation*}
such that~$B_\rho(p_i)\cap \Omega$
is equivalent to~$B_\rho^+$ via a diffeomorphism~$T_i$ of class~$C^{1,1}$ with
inverse of class~$C^{1,1}$.

Recalling~\eqref{PSSKTROFNMJCGLQMGCSDVC}, we see that the function~$\widetilde u_i(y):=u(T_i^{-1}(y))$
satisfies in~$B_\rho^+$ the same type of equation that~$u$ satisfies in~$B_\rho(p_i)\cap \Omega$,
(with structural constants of the coefficients possibly multiplied by a different constant as well).
Thus, we employ the local estimates at the boundary in~\eqref{0ojlqd0eo-2fkmv-3fgmvws0fPJON}, obtaining that
\begin{equation*} \|\widetilde u_i\|_{W^{2,p}(B_{\rho/2}^+)}\le C\,\Big(\|u\|_{L^p(\Omega)}+
\|g\|_{W^{2,p}(\Omega)}+\|f\|_{L^p(\Omega)}\Big)\end{equation*}
and consequently
\begin{equation*} \|u\|_{W^{2,p}(\Omega\cap B_{\rho/2}(p_i))}\le C\,\Big(\|u\|_{L^p(\Omega)}+
\|g\|_{W^{2,p}(\Omega)}+\|f\|_{L^p(\Omega)}\Big),\end{equation*}
which in turn leads to
\begin{equation}\label{LAMds-STudf-APEImriJErqr2vbgjqxb4rqtqrw3ww1-1-2}
\|u\|_{W^{2,p}({\mathcal{N}})}\le C\,\Big(\|u\|_{L^p(\Omega)}+
\|g\|_{W^{2,p}(\Omega)}+\|f\|_{L^p(\Omega)}\Big),\end{equation}
up to renaming constants from line to line.

Now we consider a finite family of balls such that
\begin{equation}\label{THM:0-T536EcdhPOaIJodidpo3tju24ylkiliR2543uj8890j09092po-okf-3edk-3edhf-2-fh564-09-e22r}
 {\mathcal{J}}:=\{ x\in\Omega {\mbox{ s.t. }} B_{\rho/9}(x)\cap(\partial\Omega)=\varnothing\}
\subseteq \bigcup_{i=1}^M B_{\rho/4}(q_i),\end{equation}
for suitable~$q_i\in\Omega$. By the
interior estimates in Theorem~\ref{LSmvvikSDnkjazikgAzldDrgmIDMmo5r523},
\begin{equation*} \|u\|_{W^{2,p}(B_{\rho/4}(q_i))}\le C\,\Big(\|u\|_{L^p(\Omega)}+\|f\|_{L^p(\Omega)}\Big).\end{equation*}
This and~\eqref{THM:0-T536EcdhPOaIJodidpo3tju24ylkiliR2543uj8890j09092po-okf-3edk-3edhf-2-fh564-09-e22r}
yield that
\begin{equation*} \|u\|_{W^{2,p}({\mathcal{J}})}\le C\,\Big(\|u\|_{L^p(\Omega)}+\|f\|_{W^{2,p}(\Omega)}\Big).\end{equation*}
This, together with~\eqref{LAMds-STudf-APEImriJErqr2vbgjqxb4rqtqrw3ww1-1-2}, yields the desired result.
\end{proof}

It is worth noticing that~\eqref{ojhnwfed-2irjgTYHSJDISA2R44G78jS} is the counterpart in Lebesgue spaces of the global estimate
in H\"older spaces established in~\eqref{LAMds-STudf-APEImriJEr1-1}. For the sake of completeness,
we mention that a counterpart of~\eqref{LAMds-STudf-APEImriJEr1-2} holds true as well, namely,
in the setting of Theorem~\ref{JOSLNGEITHSHLPESTGLOUPTHBUO},
if~$c(x)\le0$ for a.e.~$x\in\Omega$, then
\begin{equation}\label{JOSLNGEITHSHLPESTGLOUP-jofvnMTTOO}
\|u\|_{W^{2,p}(\Omega)}\le C\,\Big(\|g\|_{W^{2,p}(\Omega)}+\|f\|_{L^p(\Omega)}\Big),\end{equation}
see e.g.~\cite[Lemma~9.17]{MR1814364}.
This improvement of~\eqref{ojhnwfed-2irjgTYHSJDISA2R44G78jS} is very delicate from the technical point of view\footnote{For the case of the Laplace operator, we observe that~\eqref{JOSLNGEITHSHLPESTGLOUP-jofvnMTTOO}
can be obtained \label{JOSLNGEITHSHLPESTGLOUP-jofvnMTTOO-NOT}
by compactness via the standard Maximum Principle by arguing as follows. 
By possibly replacing~$u$ by~$u-g$, we can reduce ourselves to the case in which~$g$ vanishes identically.
Then, if~\eqref{JOSLNGEITHSHLPESTGLOUP-jofvnMTTOO} were violated, there would exist
sequences~$u_j\in W^{2,p}(\Omega)\cap W^{1,p}_0(\Omega)$ and~$f_j\in L^p(\Omega)$ such that~$\Delta u_j=f_j$ in~$\Omega$
and~$\|u_j\|_{L^p(\Omega)}> 2j\|f_j\|_{L^p(\Omega)}$. This and~\eqref{ojhnwfed-2irjgTYHSJDISA2R44G78jS} yield that
$$ 2j\|f_j\|_{L^p(\Omega)}<
\|u_j\|_{W^{2,p}(\Omega)}\le C\,\Big(\|f_j\|_{L^p(\Omega)}+\|u_j\|_{L^p(\Omega)}\Big)
$$
and therefore~$\|u_j\|_{L^{p}(\Omega)}> j\|f_j\|_{L^p(\Omega)}$ for~$j$ large enough.

Thus, setting~$\widetilde u_j:=\frac{u_j}{\|u_j\|_{L^p(\Omega)}}$
and~$\widetilde f_j:=\frac{f_j}{\|u_j\|_{L^p(\Omega)}}$, noticing that
$$ \|\widetilde u_j\|_{W^{2,p}(\Omega)}=
\frac{\|u_j\|_{W^{2,p}(\Omega)}}{\|u_j\|_{L^p(\Omega)}}\le\frac{
C\,\Big(\|f_j\|_{L^p(\Omega)}+\|u_j\|_{L^p(\Omega)}\Big)
}{\|u_j\|_{L^p(\Omega)}}\le C
,$$
we find that~$\widetilde u_j\in W^{2,p}(\Omega)\cap W^{1,p}_0(\Omega)$ solves~$\Delta \widetilde u_j=\widetilde f_j$ in~$\Omega$, with~$\|\widetilde u_j\|_{W^{2,p}(\Omega)}=1$
and~$\|\widetilde f_j\|_{L^p(\Omega)}<\frac1j$. Consequently, up to a subsequence, $\widetilde u_j$ converges as~$j\to+\infty$ to some~$\widetilde u$
strongly in~$W^{1,p}_0(\Omega)$ and weakly in~$W^{2,p}(\Omega)$, and~$\widetilde f_j\to0$ in~$L^p(\Omega)$.

As a result, for all~$\varphi\in C^\infty_0(\Omega)$,
$$ \int_\Omega \Delta \widetilde u(x)\,\varphi(x)\,dx=\lim_{j\to+\infty}\int_\Omega \Delta \widetilde u_j(x)\,\varphi(x)\,dx
=\lim_{j\to+\infty}\int_\Omega \widetilde f_j(x)\,\varphi(x)\,dx=0.$$
This gives that
%%%%% $\widetilde u$ belongs to~$C^2(\Omega)$ and it is harmonic 
%%%%% in~$\Omega$ (recall Lemma~\ref{WEYL}), yielding that 
for all~$\varphi\in C^\infty_0(\Omega)$
$$ \int_\Omega\nabla \widetilde u(x)\cdot\nabla\varphi(x)\,dx=0.$$
Taking a sequence of functions~$\varphi_k\in C^\infty_0(\Omega)$ such that~$\varphi_k\to\widetilde u$ in~$W^{1,p}_0(\Omega)$, we thereby conclude that
$$ \int_\Omega|\nabla \widetilde u(x)|^2\,dx=0$$
and accordingly~$\widetilde u$ vanishes identically.

But this is a contradiction with the fact that
$$ \|\widetilde u\|_{L^p(\Omega)}=\lim_{j\to+\infty}
\|\widetilde u_j\|_{L^p(\Omega)}=1$$
and this establishes~\eqref{JOSLNGEITHSHLPESTGLOUP-jofvnMTTOO} for the Laplace operator.

The same argument would mainly carry over if one possessed a convenient Maximum Principle for the operator under
consideration, see the proof of~Lemma~9.17 in~\cite{MR1814364}.}
since it relies on
a Maximum Principle which is conceptually very different
from that in Corollary~\ref{9oikr3eg9u2iwhefk9qryfhweiwwrfgfLftyOmega-kdf}, which was used to prove~\eqref{LAMds-STudf-APEImriJEr1-2} in the H\"older spaces framework:
indeed, the additional difficulty for solutions in~$W^{2,p}$
is that they do not need to satisfy the equation at every point of the domain,
making the classical Maximum Principles and the barrier methods unsuitable.
To overcome this problem, one needs to develop the so-called
Aleksandrov-Bakel\cprime man-Pucci Maximum Principle~\cite{MR0126604, MR0147776, MR214905},
see e.g.~\cite[Theorems~9.1, 9.5 and~9.15]{MR1814364},
and this is essentially the main technical tool to obtain~\eqref{JOSLNGEITHSHLPESTGLOUP-jofvnMTTOO}.

{F}rom the global estimate in~\eqref{JOSLNGEITHSHLPESTGLOUP-jofvnMTTOO},
one also obtains the following existence (and uniqueness,
thanks again to the Aleksandrov-Bakel\cprime man-Pucci Maximum Principle)
result for
the Dirichlet problem in Lebesgue spaces:

\begin{theorem}\label{21980987y4h2nf9325teuationdetildealphauation-2TH}
Let~$\Omega\subset\R^n$ be a bounded open set with boundary of class~$C^{1,1}$.
Let~$a_{ij}\in C(\overline{\Omega})$ satisfy the ellipticity 
condition in~\eqref{ELLIPTIC}. Let~$b_i$, $c\in L^\infty(\Omega)$, with~$c(x)\le0$ for
a.e.~$x\in\Omega$. Let~$f\in L^p(\Omega)$ and~$g\in
W^{2,p}(\Omega)$
for some~$p\in(1,+\infty)$.

Then, the Dirichlet problem
\begin{equation}\label{24r3t4DI21980987y4h2nf9325teuationdetildealphauation-2TH}
\begin{dcases}
\sum_{i,j=1}^n a_{ij} \partial_{ij}u+\sum_{i=1}^n b_i\partial_i u+cu=f
&\quad{\mbox{a.e. in }}\,\Omega,\\
u=g&\quad{\mbox{on }}\,\partial\Omega\end{dcases}
\end{equation}
admits a unique solution~$u\in W^{2,p}(\Omega)$.

Moreover, there exists~$C>0$, depending only on~$n$,
$p$, $\Omega$, $a_{ij}$, $b_i$ and~$c$, such that
\begin{equation}\label{24r3t4DI21980987y4h2nf9325teuationdetildealphauation-2TH02rwjfngGHSBN}
\|u\|_{W^{2,p}(\Omega)}\le C\,\Big(\|g\|_{W^{2,p}(\Omega)}+\|f\|_{L^p(\Omega)}\Big).\end{equation}
\end{theorem}

See also~\cite[Section~9.6]{MR1814364} for further details on the Dirichlet
problem in Lebesgue spaces.

We also mention that a higher regularity result
that can be considered as the counterpart of
Theorem~\ref{9827hfCXVTHM:0-T536EcdhPOaIJodidfpo3tju24ylkiliR2543uj8890j09092po-okf-3edk-3edhf-2-fh564COR}
in the Lebesgue spaces setting holds true:
see e.g.~\cite[Theorem~9.19]{MR1814364} for a regularity\index{Calder\'{o}n-Zygmund estimates|)}
theory in the class~$W^{k,p}$.\medskip

As an interesting application of the Calder\'{o}n-Zygmund estimates
(in fact, of the technically simpler part dealing just with estimates in class~$W^{1,p}$)
we present now a variant of Proposition~\ref{MALEBEAPP}
which does not require vanishing boundary conditions:

\begin{proposition}\label{MALEBEAPP-BIS}
Let~$n\ge2$ and~$\Omega$ be a bounded and open subset of~$\R^n$.

Let~$p>\frac{n}2$ and~$f:\Omega\to \R$ with~$f\in L^p(\Omega)$.

Let~$u\in C^2(\Omega)\cap L^1(\Omega)$ be a solution of
$$ \Delta u= f\quad{\mbox{ in }}\Omega.
$$
Then, for every open set~$\Omega'\Subset\Omega$,
\begin{equation} \label{MALEBEAPP-BIS-eq}\sup_{\Omega' }|u|\le C\Big( \|u\|_{L^1(\Omega)}+\|f\|_{L^p(\Omega)}\Big),\end{equation}
for some~$C>0$ depending only on~$n$, $p$, $\Omega$ and~$\Omega'$.
\end{proposition}

\begin{proof} Up to reducing to connected components, we will implicitly assume here that the domains under consideration are connected.
We consider open sets~$\Omega''$ and~$\Omega'''$ with smooth boundaries and such that~$ \Omega'\Subset\Omega''\Subset\Omega'''\Subset\Omega$.
We pick~$\phi\in C^\infty_0(\Omega''',[0,1])$ with~$\phi=1$ in~$\Omega''$
and define~$v:=\phi u$.

In this way, $v=0$ on~$\partial\Omega'''$ and, in~$\Omega'''$,
$$ \Delta v=\Delta \phi u+2\nabla \phi \cdot\nabla u+\phi \Delta u=:g.$$
Thus, we can employ Proposition~\ref{MALEBEAPP} and deduce that
\begin{equation}\label{L2SmvvikSDnkjaz34567890-iuygtf9ikgAzldDrgmIDMmo5r523} \sup_{\Omega'' }v\le C\|g^-\|_{L^p(\Omega''')}.\end{equation}
We stress that, up to renaming~$C$ at each level of our calculations,
\begin{eqnarray*}
g^-(x)&=&\max\Big\{0,-\Delta \phi (x)u(x)-2\nabla \phi(x) \cdot\nabla u(x)-\phi(x) \Delta u(x) \Big\}\\&\le&
C\Big(|u(x)|+|\nabla u(x)|+|\Delta u(x)|
\Big)
\end{eqnarray*}
and therefore
$$ \|g^-\|_{L^p(\Omega''')}\le C\Big(\|u\|_{L^p(\Omega''')}+\|\nabla u\|_{L^p(\Omega''')}+\|\Delta u\|_{L^p(\Omega''')}\Big).$$
But taking into consideration Theorem~\ref{LSmvvikSDnkjazikgAzldDrgmIDMmo5r523} (in its simpler formulation for first derivatives)
up to a covering argument we also know that
$$ \|\nabla u\|_{L^p(\Omega''')}%\le \| u\|_{W^{2,p}(\Omega''')}
\le C\Big(\|u\|_{L^p(\Omega''')}+\|\Delta u\|_{L^p(\Omega''')}\Big)$$
and consequently
$$ \|g^-\|_{L^p(\Omega''')}\le C\Big(\|u\|_{L^p(\Omega''')}+\|\Delta u\|_{L^p(\Omega''')}\Big).$$
This inequality, \eqref{L2SmvvikSDnkjaz34567890-iuygtf9ikgAzldDrgmIDMmo5r523} and the fact that~$v=u$ in~$\Omega''$
give that
\begin{equation*} \sup_{\Omega'' }u\le C\Big(\|u\|_{L^p(\Omega''')}+\|\Delta u\|_{L^p(\Omega''')}\Big).\end{equation*}

Applying this to~$-u$ in lieu of~$u$, we conclude that\footnote{A quick
way to obtain~\eqref{L2SmvvikSDnkjaz34567890-iuygtf9ikgAzldDrgmIDMmo5r523-b}
would also be to use Theorem~\ref{LSmvvikSDnkjazikgAzldDrgmIDMmo5r523}
in its full standing (rather than its simpler form for first derivatives), obtaining second derivative estimates such as
\begin{equation*} \|u\|_{W^{2,p}(\Omega'')}\le C\,\Big(\|u\|_{L^p(\Omega''')}+\|\Delta u\|_{L^p(\Omega''')}\Big).\end{equation*}
Since~$p>\frac{n}2$, this and the Sobolev Embedding Theorem (see e.g.~\cite[Theorem~6(ii) on page~270]{MR1625845}) lead to~\eqref{L2SmvvikSDnkjaz34567890-iuygtf9ikgAzldDrgmIDMmo5r523-b}.}
\begin{equation}\label{L2SmvvikSDnkjaz34567890-iuygtf9ikgAzldDrgmIDMmo5r523-b} \sup_{\Omega'' }|u|\le C\Big(\|u\|_{L^p(\Omega''')}+\|\Delta u\|_{L^p(\Omega''')}\Big).\end{equation}
This establishes~\eqref{MALEBEAPP-BIS-eq} with the $L^1$-norm of the solution replaced by its~$L^p$-norm.

To obtain~\eqref{MALEBEAPP-BIS-eq} in its full generality we proceed as follows.
Given~$R>r>0$, we can rescale and apply~\eqref{L2SmvvikSDnkjaz34567890-iuygtf9ikgAzldDrgmIDMmo5r523-b}
to balls of radius~$R-r$ centered at points of~$B_r$, finding that
\begin{equation}\label{L2SmvvikSDnkjaz34567890-iuygtf9ikgAzldDrgmIDMmo5r523-ba} \sup_{B_r }|u|\le \frac{C}{(R-r)^{\frac{n}{p}}}\Big(\|u\|_{L^p(B_R)}+(R-r)^2\|\Delta u\|_{L^p(B_R)}\Big).\end{equation}

Now we claim that for every~$\e>0$ there exists~$C_\e>0$ such that
\begin{equation}\label{L2SmvvikSDnkjaz34567890-iuygtf9ikgAzldDrgmIDMmo5r523-bb}
\|u\|_{L^p(B_R)}\le \e(R-r)^{\frac{n}p}\sup_{B_R}|u|+\frac{C_\e\|u\|_{L^1(B_R)}}{(R-r)^{\frac{n(p-1)}{p}}}.
\end{equation}
Indeed, we have that
\begin{eqnarray*}
\|u\|_{L^p(B_R)}\le \left(\sup_{B_R} |u|^{p-1}\int_{B_R}|u(x)|\,dx\right)^{\frac1p}
\end{eqnarray*}
and the desired claim follows from Young's Inequality with exponents~$\frac{p}{p-1}$ and~${p}$.

As a result, if
$$ \varphi(\rho):=\sup_{ B_\rho}|u|,$$
by choosing~$\e>0$ conveniently small we deduce from~\eqref{L2SmvvikSDnkjaz34567890-iuygtf9ikgAzldDrgmIDMmo5r523-ba}
and~\eqref{L2SmvvikSDnkjaz34567890-iuygtf9ikgAzldDrgmIDMmo5r523-bb} that
$$ \varphi(r)\le \frac{ \varphi(R)}2
+\frac{C}{(R-r)^{{n}}}\|u\|_{L^1(B_R)}+\frac{C}{(R-r)^{\frac{n}{p}-2}} \|f\|_{L^p(B_R)}.
$$
Since we are dealing with local estimates, we can suppose that~$R\le R_0$ for some~$R_0>0$, therefore~$(R-r)^{\frac{n}{p}-2}=(R-r)^{n- \frac{n(p-1)}{p}-2}
\ge R_0^{- \frac{n(p-1)}{p}-2} (R-r)^{n}$, yielding that
$$  \varphi(r)\le \frac{ \varphi(R)}2
+\frac{C\big(\|u\|_{L^1(B_R)}+\|f\|_{L^p(B_R)}\big)}{(R-r)^{n}}.$$

We are thereby in the position of using Lemma~\ref{S-O-GIAHJAQ} and infer that
$$ \sup_{B_r}|u|=\varphi(r)\le \frac{C}{(R-r)^{n}}\Big(\|u\|_{L^1(B_R)}+\|f\|_{L^p(B_R)}\Big).$$
{F}rom this, a covering argument gives~\eqref{MALEBEAPP-BIS-eq}, as desired.
\end{proof}

A consequence of the above result is the following\footnote{Interestingly,
results such as the one in Corollary~\ref{MALEBEAPP-BIS-eqoihknfUNArareDMBE} are usually obtained by energy methods, relying on divergence form operators, while the approach presented here is of nondivergence type.}
observation:

\begin{corollary}\label{MALEBEAPP-BIS-eqoihknfUNArareDMBE}
Let~$n\ge2$.
Let~$\Omega$, $\Omega'$ be bounded open sets of~$\R^n$ with~$\Omega'\Subset\Omega$.
Let~$a\in L^r(\Omega)$ for some~$r>\frac{n}2$.

Let~$u\in C^2(\Omega)\cap L^1(\Omega)$ be a solution of~$ \Delta u=au$ in~$\Omega$.

Then, there exists a positive constant~$C$, depending only on~$n$, $r$, $\Omega$, $\Omega'$ and~$\|a\|_{L^r(\Omega)}$, such that
$$ \|u\|_{L^\infty(\Omega')}\le C \|u\|_{L^1(\Omega)} .$$
\end{corollary}

\begin{proof}Up to a covering argument, we deal with balls, we pick~$p\in\left(\frac{n}2,r\right)$, we revisit~\eqref{L2SmvvikSDnkjaz34567890-iuygtf9ikgAzldDrgmIDMmo5r523-ba} in the present setting
and we use the H\"older Inequality with exponents~$\frac{r}{p}>1$ and~$\frac{r}{r-p}$
to see that, given~$R>r>0$,
\begin{equation} \label{L2Sm78XcvTyvvikSDnkjaz34567890-iuygtf9ikgAzldDrgmIDMmo5r523-ba89-1}\begin{split}\sup_{B_r }|u|&\le \frac{C}{(R-r)^{\frac{n}{p}}}\Big(\|u\|_{L^p(B_R)}+(R-r)^2\|au\|_{L^p(B_R)}\Big)\\
&\le\frac{C}{(R-r)^{\frac{n}{p}}}\Big(\|u\|_{L^p(B_R)}+(R-r)^2\|a\|_{L^r(B_R)}
\|u\|_{L^{\frac{rp}{r-p}}(B_R)}\Big).
\end{split}\end{equation}

We also utilize~\eqref{L2SmvvikSDnkjaz34567890-iuygtf9ikgAzldDrgmIDMmo5r523-bb}
with~$p$ replaced by~${\frac{rp}{r-p}}$ and we find that
\begin{equation}\label{L2Sm78XcvTyvvikSDnkjaz34567890-iuygtf9ikgAzldDrgmIDMmo5r523-ba89-2}
\|u\|_{L^{\frac{rp}{r-p}}(B_R)}\le \e(R-r)^{\frac{n(r-p)}{rp}}\sup_{B_R}|u|+\frac{C_\e\|u\|_{L^1(B_R)}}{(R-r)^{\frac{n( rp-r+p)}{rp}}}.
\end{equation}

We also remark that
$$ 2+\frac{n(r-p)}{rp}-\frac{n}p=\frac{2r-n}r>0.$$
On this account, combining~\eqref{L2SmvvikSDnkjaz34567890-iuygtf9ikgAzldDrgmIDMmo5r523-ba}, \eqref{L2Sm78XcvTyvvikSDnkjaz34567890-iuygtf9ikgAzldDrgmIDMmo5r523-ba89-1} and~\eqref{L2Sm78XcvTyvvikSDnkjaz34567890-iuygtf9ikgAzldDrgmIDMmo5r523-ba89-2},
assuming~$R\le R_0$ for some given~$R_0>0$ and defining
$$ \varphi(\rho):=\sup_{ B_\rho}|u|,$$up to redefining constants
we conclude that
\begin{eqnarray*}
\varphi(r)&=&\sup_{B_r }|u|\\
&\le&\frac{C}{(R-r)^{\frac{n}{p}}}\Big(\|u\|_{L^p(B_R)}+(R-r)^2\|a\|_{L^r(B_R)}
\|u\|_{L^{\frac{rp}{r-p}}(B_R)}\Big)\\&\le&
\frac{C}{(R-r)^{\frac{n}{p}}}\Bigg[
\e(R-r)^{\frac{n}p}\sup_{B_R}|u|+\frac{C_\e\|u\|_{L^1(B_R)}}{(R-r)^{\frac{n(p-1)}{p}}}\\&&\qquad\qquad\quad
+(R-r)^2\Bigg(
\e(R-r)^{\frac{n(r-p)}{rp}}\sup_{B_R}|u|+\frac{C_\e\|u\|_{L^1(B_R)}}{(R-r)^{\frac{n( rp-r+p)}{rp}}}\Bigg)\Bigg]\\&\le&
C\left(\e\sup_{B_R}|u|
+\frac{ \|u\|_{L^1(B_R)} }{(R-r)^{n}}\right)\\&\le&\frac{\varphi(R)}2+\frac{ C\|u\|_{L^1(B_R)} }{(R-r)^{n}},
\end{eqnarray*}
as long as we choose~$\e$ appropriately small.

We are thereby in the position of using Lemma~\ref{S-O-GIAHJAQ} and infer that
$$ \sup_{B_r}|u|=\varphi(r)\le \frac{C\|u\|_{L^1(B_R)}}{(R-r)^{n}}
,$$
from which we obtain the desired result.
\end{proof}

\chapter{The Dirichlet problem in the light of capacity theory}

In this chapter we will recall the notion of capacity and exploit it to determine in which cases
the solution of the Dirichlet problem obtained via Perron methods (recall Corollary~\ref{S-coroEXIS-M023})
is continuous up to the boundary. 
This is indeed a delicate issue in light of counterexamples such as
the Lebesgue spine discussed on page~\pageref{SILEBE678NE},
see~\cite[Exercise~10, page~334]{MR0222317}.

\section{Capacitance and capacity}\label{Capacitance-and-capacity}

The mathematical concept of capacity was initiated\index{capacity|(}
by Gustave Choquet~\cite{MR80760}, see also~\cite{MR867115} for a
historical account on the creation and development of this theory.

Different versions of the capacity theory
are available in the literature: here we present a treatment that we find sufficiently close to the physical
intuition (though a number of technical details, especially to deal with general sets,
appear to be unavoidable).

Given a bounded and open set~$\Omega\subset\R^n$ with boundary
of class~$C^{2,\alpha}$, for some~$\alpha\in(0,1)$, we define
the capacity of~$\Omega$ as
\begin{equation}\label{DEFCAP} \cAPAC(\Omega):=\inf_{{v\in C^\infty_0(\R^n)}\atop{v=1 {\text{ in }}\Omega}}
\int_{\R^n\setminus\overline\Omega}|\nabla v(x)|^2\,dx.\end{equation}
By the density results in Sobolev spaces (see e.g.~\cite[Theorem~10.29]{MR2527916})
one can also rephrase~\eqref{DEFCAP} in the form
\begin{equation}\label{DEFCAP-SOBO} \cAPAC(\Omega)=\inf_{{v\in {\mathcal{D}}^{1,2}(\R^n)}\atop{v=1 {\text{ in }}\Omega}}
\int_{\R^n\setminus\overline\Omega}|\nabla v(x)|^2\,dx.\end{equation}
This definition of capacity that we adopted in~\eqref{DEFCAP-SOBO}
is meaningful only when~$n\ge3$, since:

\begin{lemma}
Let~$n\in\{1,2\}$. Then, $\cAPAC(\Omega)=0$
for every bounded and open set~$\Omega\subset\R^n$ with boundary
of class~$C^{2,\alpha}$ for some~$\alpha\in(0,1)$.
\end{lemma}

\begin{proof} Let~$n=1$ and~$\e>0$. Let~$\Omega\Subset(-R,R)$ for some~$R>0$ and
$$ v_\e(x)=\begin{dcases}
1 & {\mbox{ if }}x\in[-R,R],\\
\e\left(R+\displaystyle\frac{1}\e-|x|\right)& {\mbox{ if }}|x|\in\left(R,R+\frac1\e\right),\\0
& {\mbox{ if }}|x|\in\left[R+\frac1\e,+\infty\right).
\end{dcases}$$
Then, by~\eqref{DEFCAP-SOBO},
$$ \cAPAC(\Omega)\le
\int_{\R\setminus\overline\Omega}|v'_\e(x)|^2\,dx=
2\e.$$
Hence, taking~$\e$ as small as we wish, we conclude that~$\cAPAC(\Omega)=0$.

Let now~$n=2$ and take~$\e>0$ so small that~$\Omega\Subset B_{1/\e}$. Let also
$$ w_\e(x)=\begin{dcases}1&{\mbox{ if }}|x|\in\left[0,\frac1\e\right],\\
\displaystyle\frac{2\ln\e-\ln|x|}{\ln\e} & {\mbox{ if }}|x|\in\left(\frac1\e,\frac1{\e^2}\right),\\0
& {\mbox{ if }}|x|\in\left[\frac1{\e^2},+\infty\right).
\end{dcases}$$
Then, by~\eqref{DEFCAP-SOBO} we have
$$ \cAPAC(\Omega)\le
\int_{\R^2\setminus\overline\Omega}|\nabla w_\e(x)|^2\,dx=\frac1{\ln^2\e}\,
\int_{B_{1/\e^2}\setminus B_{1/\e}}\frac{dx}{|x|^2}=\frac{2\pi}{\ln^2\e}\,
\int_{{1/\e}}^{{1/\e^2}}\frac{d\rho}{\rho}=\frac{2\pi}{|\ln\e|}.
$$
Taking~$\e$ as small as we wish, we conclude that~$\cAPAC(\Omega)=0$ in this case too.
\end{proof}

We now clarify the link between the notion of capacity
and that of harmonic functions vanishing at infinity and with prescribed
value along the boundary of a domain:

\begin{proposition}\label{AM-CAPA}
Let~$n\ge3$
and~$\Omega\subset\R^n$ be a bounded and open set with boundary of class~$C^{2,\alpha}$ for some~$\alpha\in(0,1)$.

Then, there exists a unique function~$u\in C^2(\R^n\setminus\Omega)\cap {\mathcal{D}}^{1,2}(\R^n\setminus\Omega)$ such that
$$ \begin{dcases}
\Delta u=0 \quad {\mbox{ in~$\,\R^n\setminus\overline\Omega$,}}\\
u=1 \quad{\mbox{ on $\,\partial\Omega$.}}
\end{dcases}$$
Additionally, $0\le u(x)\le1$ for all~$x\in\R^n\setminus\overline\Omega$,
\begin{equation}\label{jiwur8b4v784578476768}
\lim_{|x|\to+\infty} u(x)=0
\end{equation} and
\begin{equation}\label{CAB9IJSalskdmnfcx34546fsdQsdfasdAOSL} \cAPAC(\Omega)=\int_{\R^n\setminus\overline\Omega}|\nabla u(x)|^2\,dx.\end{equation}
\end{proposition}

With a slight abuse of notation, we may implicitly assume that the function~$u$ in Proposition~\ref{AM-CAPA}
is defined in all~$\R^n$, just by setting~$u$ to be constantly equal to~$1$ in~$\Omega$.

\begin{proof}[Proof of Proposition~\ref{AM-CAPA}] Let~$k_0\in\N$ be such that~$B_{k_0}\Supset\Omega$.
Let~$v\in C^\infty_0(B_{k_0})$ such that~$v=1$ in~$\Omega$.
For every~$k\in\N\cap[k_0,+\infty)$ we consider the Sobolev space~$X_k$
of the functions~$w\in {\mathcal{D}}^{1,2}(\R^n)$ such that~$w=v$ in~$\Omega\cup (\R^n\setminus B_k)$.
We take~$u_k$ to be the minimizer in~$X_k$ of the functional
$$ X_k\ni w\longmapsto \int_{\R^n\setminus\overline\Omega}|\nabla w(x)|^2\,dx.$$
We stress that such a minimizer exists, by the Direct Method of the Calculus of Variations
see e.g.~\cite{MR772245} and~\cite[Theorems~11.10, 12.15 and~15.29]{MR2527916}.

We point out that, if
$$ u_k^\star(x):=\begin{dcases}
1 & {\mbox{ if }}u_k(x)>1,\\
u_k(x)& {\mbox{ if }}u_k(x)\in[0,1],\\ 0
& {\mbox{ if }}u_k(x)<0,
\end{dcases}$$
then
$$ \int_{\R^n\setminus\overline\Omega}|\nabla u_k^\star(x)|^2\,dx
=\int_{(\R^n\setminus\overline\Omega)\cap\{ 0\le u_k\le1\}}|\nabla u_k(x)|^2\,dx
\le\int_{\R^n\setminus\overline\Omega}|\nabla u_k (x)|^2\,dx.$$ Thus,
up to replacing~$u_k$ by~$u_k^\star$ we can suppose that
\begin{equation}\label{LA-SORROSDTTOSVOLMIANNDKAa3E89SD}
0\le u_k(x)\le1\qquad{\mbox{for all }}\,x\in\R^n\setminus\overline\Omega.\end{equation}
Moreover, since~$X_{k_0}\subseteq X_k$ for all~$k\in\N\cap[k_0,+\infty)$, the minimization property in~$X_k$ yields that
\begin{equation}\label{0oijnSIMSi4a3d34ms00dn-AIsmcAS} \int_{\R^n\setminus\overline\Omega}|\nabla u_k(x)|^2\,dx
\le \int_{\R^n\setminus\overline\Omega}|\nabla u_{k_0}(x)|^2\,dx=:C_0\in[0,+\infty),\end{equation}
and therefore, up to a subsequence, we can suppose that~$u_k$ converges to some function~$u$
in~$L^2_{\rm loc}(\R^n)$
and~$\nabla u_k$ converges to~$\nabla u$
weakly in~$L^2(\R^n)$.

Since the minimizer~$u_k$ is harmonic in~$\R^n\setminus\overline\Omega$,
so is~$u$, thanks to Corollary~\ref{S-OS-OSSKKHEPARHA}.
Also, since~$u_k=v=1$ in~$\Omega$, also~$u=1$ in~$\Omega$.
This gives that~$u$ is an admissible function for the infimum procedure
in the right hand side of~\eqref{DEFCAP-SOBO}. 

We also point out that, by the above mentioned weak convergence,
\begin{equation}\label{WEAKNOSE-22}\begin{split}
0\,&\le\,\liminf_{k\to+\infty} 
\int_{\R^n\setminus\overline\Omega}|\nabla u_{k}(x)-\nabla u(x)|^2\,dx\\&
=\,\liminf_{k\to+\infty}\int_{\R^n\setminus\overline\Omega}\Big(|\nabla u(x)|^2-2
\nabla u_{k}(x)\cdot\nabla u(x)+|\nabla u_k(x)|^2\Big)\,dx\\&
=\,-\int_{\R^n\setminus\overline\Omega}|\nabla u(x)|^2\,dx+
\liminf_{k\to+\infty}\int_{\R^n\setminus\overline\Omega}|\nabla u_k(x)|^2\,dx.
\end{split}\end{equation}
In particular, by~\eqref{0oijnSIMSi4a3d34ms00dn-AIsmcAS},
\begin{equation}\label{WEAKNOSE}
\int_{\R^n\setminus\overline\Omega}|\nabla u(x)|^2\,dx
\le C_0.
\end{equation}
Now we show that
\begin{equation}\label{6yh7uyhHSSN9oik0o2r3}
\int_{\R^n\setminus\overline\Omega}|\nabla u(x)|^2\,dx=
\inf_{{\zeta\in {\mathcal{D}}^{1,2}(\R^n)}\atop{\zeta=1 {\text{ in }}\Omega}}
\int_{\R^n\setminus\overline\Omega}|\nabla \zeta(x)|^2\,dx.
\end{equation}
Indeed, suppose not, then there exist~$b>0$ and~$\zeta\in {\mathcal{D}}^{1,2}(\R^n)$
with~$\zeta=1$ in~$\Omega$ such that
$$ \int_{\R^n\setminus\overline\Omega}|\nabla \zeta(x)|^2\,dx+b\le
\int_{\R^n\setminus\overline\Omega}|\nabla u(x)|^2\,dx.$$
By the density results in Sobolev spaces (see e.g.~\cite[Theorem~10.29]{MR2527916}),
we can find~$w\in C^\infty_0(\R^n)$ with~$w=1$ in~$\Omega$ such that
$$ \int_{\R^n\setminus\overline\Omega}|\nabla w(x)|^2\,dx\le\frac{b}2+
\int_{\R^n\setminus\overline\Omega}|\nabla \zeta(x)|^2\,dx.$$
We take~$k_\star\in\N\cap[k_0,+\infty)$ sufficiently large such that the support of~$w$ is contained in~$B_{k_\star}$.
As a consequence, we have that~$w\in X_k$ for all~$k\ge k_\star$, and therefore
$$ \int_{\R^n\setminus\overline\Omega}|\nabla u_k(x)|^2\,dx\le
\int_{\R^n\setminus\overline\Omega}|\nabla w(x)|^2\,dx.$$
By collecting these items of information and utilizing~\eqref{WEAKNOSE-22}, we see that
\begin{eqnarray*}&&
\frac{b}2=b-\frac{b}2\le\left(
\int_{\R^n\setminus\overline\Omega}|\nabla u(x)|^2\,dx-
\int_{\R^n\setminus\overline\Omega}|\nabla \zeta(x)|^2\,dx\right)+
\left(\int_{\R^n\setminus\overline\Omega}|\nabla \zeta(x)|^2\,dx-
\int_{\R^n\setminus\overline\Omega}|\nabla w(x)|^2\,dx\right)\\&&\quad\qquad\le
\liminf_{k\to+\infty}\int_{\R^n\setminus\overline\Omega}|\nabla u_k(x)|^2\,dx-\int_{\R^n\setminus\overline\Omega}|\nabla w(x)|^2\,dx
\le0.
\end{eqnarray*}
This contradiction proves~\eqref{6yh7uyhHSSN9oik0o2r3}.

Thus, the claim in~\eqref{CAB9IJSalskdmnfcx34546fsdQsdfasdAOSL}
follows from~\eqref{DEFCAP-SOBO}
and~\eqref{6yh7uyhHSSN9oik0o2r3}.

We also observe that
\begin{equation}\label{REGC1BOU}
u\in C(\R^n\setminus\Omega).
\end{equation}
To prove this, we let~$z\in\partial\Omega$
and use the regularity of~$\Omega$ to find a ball~$B_{\varrho_0}(z_0)$ such that~$B_{\varrho_0}(z_0)\subseteq\Omega$ and~$z\in(\partial\Omega)\cap(\partial B_{\varrho_0}(z_0))$.
We define
\begin{equation}\label{DEPHIBARRPERPRER8} \Phi(x):=\begin{dcases}
u(x)-\frac{\varrho_0^{n-2}}{|x-z_0|^{n-2}}&{\mbox{ if }}x\in\R^n\setminus B_{\varrho_0}(z_0),\\
u(x)-1&{\mbox{ if }}x\in B_{\varrho_0}(z_0)\end{dcases}\end{equation}
and we observe that~$\Phi\in {\mathcal{D}}^{1,2}(\R^n\setminus\Omega)$;
moreover, if~$x\in\partial\Omega$ then~$|x-z_0|\ge\varrho_0$, hence,
in the Sobolev trace sense, we have that~$\Phi\ge u-1=0$
along~$\partial\Omega$. Additionally, $\Phi$
is harmonic in~$\R^n\setminus\overline\Omega$.
Therefore, we can exploit Lemma~\ref{SOBOMAXPLEH10}
(used here with~${\mathcal{U}}:=\R^n\setminus\overline\Omega$) and gather that~$\Phi\ge0$.

This leads to
$$ \liminf_{x\to z}u(x)\ge\liminf_{x\to z}\frac{\varrho_0^{n-2}}{|x-z_0|^{n-2}}=
\frac{\varrho_0^{n-2}}{|z-z_0|^{n-2}}=1.
$$
Since by construction~$u\le1$ (recall~\eqref{LA-SORROSDTTOSVOLMIANNDKAa3E89SD}), we thus obtain that
$$ \lim_{x\to z}u(x)=1,$$
thus proving~\eqref{REGC1BOU}.

Now we show that
\begin{equation}\label{ANvibDVISMvaosdIAsNOA}
\lim_{|x|\to+\infty} u(x)=0.
\end{equation}
To this end, equation~\eqref{WEAKNOSE}
and the Gagliardo-Nirenberg-Sobolev Inequality (see e.g.~\cite[Theorem~11.2]{MR2527916}) yield that, since~$n\ge3$,
\begin{equation}\label{nJMSIgbhnGnsdifNGkmdEGSI68LY} \|u\|_{L^{\frac{2n}{n-2}}(\R^n\setminus\overline\Omega)}\le C_1,\end{equation}
for some~$C_1>0$.

Furthermore, by~\eqref{LA-SORROSDTTOSVOLMIANNDKAa3E89SD}, we have that~$\|u\|_{L^\infty({\R^n\setminus\overline\Omega})}$. {F}rom this and 
Cauchy's Estimate in Theorem~\ref{CAUESTIMTH} we infer that,
for every~$x_0\in\R^n$ such that~$B_2(x_0)\cap\Omega=\varnothing$,
\begin{equation}\label{nJMSIgbhnGnsdifNGkmdEGSI68LY2} \left|\nabla u(x_0)
\right|\le C,\end{equation}
for some constant~$C>0$ depending only on~$n$.

We now recall that\footnote{To prove~\eqref{MS-oiuyhgfvb0987654-SPd-THA0SbaSKMlsdl},
one can argue by contradiction, assuming that there exists a sequence~$p_j\in\R^n$
such that~$|p_j|\to+\infty$ as~$j\to+\infty$ and such that~$|\phi(p_j)|\ge a$, for some~$a>0$.
Then one takes~$\rho>0$ such that~$|\phi(x)-\phi(p_j)|\le\frac{a}2$ for every~$x\in B_\rho(p_j)$.
Since, up to a subsequence, we can assume that the balls~$B_\rho(p_j)$ are disjoint, we find that
$$ \|\phi\|_{L^q(\R^n)}^q\ge\sum_{j=0}^{+\infty}\int_{B_\rho(p_j)}
|\phi(x)|^q\,dx
\ge\sum_{j=0}^{+\infty}\int_{B_\rho(p_j)}\big(|\phi(p_j)|-
|\phi(x)-\phi(p_j)|\big)^q\,dx\ge
\sum_{j=0}^{+\infty}\int_{B_\rho(p_j)}
\left(\frac{a}2\right)^q\,dx=+\infty,$$
providing the desired contradiction.}
\begin{equation}\label{MS-oiuyhgfvb0987654-SPd-THA0SbaSKMlsdl}\begin{split}
&{\mbox{if a uniformly continuous function~$\phi$ is in $L^q(\R^n)$ for some~$q\in[1,+\infty)$,}}\\&{\mbox{then }}\,\lim_{|x|\to+\infty}\phi(x)=0.\end{split}
\end{equation}
{F}rom~\eqref{nJMSIgbhnGnsdifNGkmdEGSI68LY},
\eqref{nJMSIgbhnGnsdifNGkmdEGSI68LY2}
and~\eqref{MS-oiuyhgfvb0987654-SPd-THA0SbaSKMlsdl}
we infer~\eqref{ANvibDVISMvaosdIAsNOA}, as desired.

Furthermore, we claim that the function~$u$ that we constructed is unique. For this, suppose that~$u$
and~$v$ satisfy the theses of Proposition~\ref{AM-CAPA} and let~$w:=\frac{u+v}2$.
Since~$w=1$ on~$\partial\Omega$, we infer from~\eqref{DEFCAP-SOBO} that
\begin{eqnarray*}&&
\int_{\R^n\setminus\overline\Omega}|\nabla u(x)-\nabla v(x)|^2\,dx=
\int_{\R^n\setminus\overline\Omega}\Big(|\nabla u(x)|^2+|\nabla v(x)|^2-2\nabla u(x)\cdot\nabla v(x)\Big)\,dx\\&&\qquad=2\cAPAC(\Omega)-2\int_{\R^n\setminus\overline\Omega}\nabla u(x)\cdot\nabla v(x)\,dx\le
2\int_{\R^n\setminus\overline\Omega}|\nabla w(x)|^2\,dx-2
\int_{\R^n\setminus\overline\Omega}\nabla u(x)\cdot\nabla v(x)\,dx\\
&&\qquad=
\int_{\R^n\setminus\overline\Omega}\frac{|\nabla u(x)|^2+|\nabla v(x)|^2+2\nabla u(x)\cdot\nabla v(x)}{2}\,dx-2
\int_{\R^n\setminus\overline\Omega}\nabla u(x)\cdot\nabla v(x)\,dx\\&&\qquad=\frac12\int_{\R^n\setminus\overline\Omega}|\nabla u(x)-\nabla v(x)|^2\,dx.
\end{eqnarray*}
This gives that~$\int_{\R^n\setminus\overline\Omega}|\nabla u(x)-\nabla v(x)|^2\,dx=0$,
and therefore the function~$u-v$ is constant. By looking at the values at infinity,
we deduce that~$u-v$ is constantly equal to zero,
thus establishing the uniqueness claim in Proposition~\ref{AM-CAPA}.

To complete the proof of Proposition~\ref{AM-CAPA}, it remains to show that~$u\in C^2(\R^n\setminus\Omega)$. 
For this, we consider the problem
$$ \begin{dcases}
\Delta v=0 \quad {\mbox{ in~$\,B_{k_0}\setminus\overline\Omega$,}}\\
v=1 \quad{\mbox{ on $\,\partial\Omega$,}}\\
v=u \quad{\mbox{ on $\,\partial B_{k_0}$,}}
\end{dcases}$$
and exploit Theorem~\ref{Theorem6.14GT} to say that it admits a unique solution~$v\in C^2(B_{k_0}\setminus\overline\Omega)
\cap C(\overline {B_{k_0}}\setminus\Omega)$.
Hence, by~\eqref{REGC1BOU} and the uniqueness statement we deduce that~$v=u$. As a consequence of this
and~\eqref{134-24rf-pkwPK-TAGBco6YHSRVrfik-RiodSS02o3ef2} in Theorem~\ref{Theorem6.14GT}
we obtain that~$u\in C^2(\R^n\setminus\Omega)$,
as desired.
\end{proof}

The statement in formula~\eqref{jiwur8b4v784578476768} of Proposition~\ref{AM-CAPA} can be made more precise, as follows:

\begin{proposition}\label{95743jhyrgeyPROPiefergerhe}
Let~$n\ge3$
and~$\Omega\subset\R^n$ be a bounded and open set with boundary of class~$C^{2,\alpha}$ for some~$\alpha\in(0,1)$.
Let~$u$ be the function given by Proposition~\ref{AM-CAPA}.
Assume that~$\Omega\subset B_{R_0}$, for some~$R_0>0$.
Then, for every~$x\in\R^n\setminus B_{R_0}$,
\begin{equation*}
u(x)\le\frac{R_0^{n-2}}{|x|^{n-2}}.
\end{equation*}
\end{proposition}

\begin{proof}
For every~$\e\in(0,1)$, we set
$$ v_\e(x):=\frac{R_0^{n-2}}{|x|^{n-2}} +\e.$$
Also, in light of~\eqref{jiwur8b4v784578476768}, there exists~$R_\e\ge1/\e$ such
that~$|u(x)|\le\frac{\e}2$ for every~$x\in\R^n\setminus B_{R_\e}$.
We notice that
\begin{eqnarray*}
&& v_\e(x)=1+\e>1\ge u(x) \qquad {\mbox{ on }}\partial B_{R_0}\\
{\mbox{and }}&&
v_\e(x)\ge \e> \frac{\e}2\ge u(x)
\qquad {\mbox{ on }}\partial B_{R_\e}.\end{eqnarray*}
Furthermore, both~$u$ and~$v_\e$ are harmonic in~$ B_{R_\e}\setminus B_{R_0}$, and therefore by the Weak Maximum
Principle in Corollary~\ref{WEAKMAXPLE}, we have that
$$ u(x)\le v_\e(x)=\frac{R_0^{n-2}}{|x|^{n-2}} +\e
\qquad{\mbox{ for every }}x\in B_{R_\e}\setminus B_{R_0}.$$
Accordingly, sending~$\e\searrow0$ we obtain the desired result.
\end{proof}

The result in Proposition~\ref{AM-CAPA} highlights the strong connection (actually,
basically, the coincidence, up to physical constants and different terminologies) between
the mathematical definition of capacity in~\eqref{DEFCAP} and
the physical notion of capacitance\index{capacitance} (or, more precisely, of
self-capacitance, see e.g.~\cite{MR3012344})
for an isolated conductor, which is the amount of electric charge that must be added to 
the conductor to raise its electric potential by one unit. \index{electrostatics}
To highlight this physical motivation, we make the following observation:

\begin{proposition}\label{propfoutr65465548789}
Let~$n\ge3$
and~$\Omega\subset\R^n$ be a bounded and open set with boundary of class~$C^{2,\alpha}$ for some~$\alpha\in(0,1)$.
Let~$u$ be the function given by Proposition~\ref{AM-CAPA}.
Then,
\begin{equation}\label{CAB9IJSalskdmnfcx34546fsdQsdfasdAOSL-C1} \cAPAC(\Omega)=\int_{\R^n\setminus\overline\Omega}|\nabla u(x)|^2\,dx=
-\int_{\partial\Omega}\frac{\partial u}{\partial\nu}(x)\,d{\mathcal{H}}^{n-1}_x.\end{equation}
\end{proposition}

\begin{proof} Since~$\Omega$ is bounded, we can consider a radius~$R_0>0$ such that~$\Omega\subset B_{R_0}$.
We take~$R>2R_0$ and we observe that, 
by the first Green's Identity in~\eqref{GRr1},
\begin{equation}\begin{split}\label{s9385v9bsadfaghjfkmkjhtrjhtrjttrj657909877derjv4eolriyu}
&
\int_{B_{R}\setminus\overline\Omega}|\nabla u(x)|^2\,dx
=\int_{\partial(B_{R}\setminus\overline\Omega)}u(x)\frac{\partial u}{\partial\nu}(x)\,d{\mathcal{H}}^{n-1}_x\\&\qquad\qquad=
\int_{\partial B_{R}}u(x)\frac{\partial u}{\partial\nu}(x)\,d{\mathcal{H}}^{n-1}_x-
\int_{\partial\Omega}u(x)\frac{\partial u}{\partial\nu}(x)\,d{\mathcal{H}}^{n-1}_x\\
&\qquad\qquad=
\int_{\partial B_{R}}u(x)\frac{\partial u}{\partial\nu}(x)\,d{\mathcal{H}}^{n-1}_x-
\int_{\partial\Omega}\frac{\partial u}{\partial\nu}(x)\,d{\mathcal{H}}^{n-1}_x
.
\end{split}\end{equation}
We observe that
\begin{equation}\label{s9385v9bsadfaghjfkmkjhtrjhtrjttrj657909877derjv4eolriyu22}
\lim_{R\to+\infty} \int_{B_{R}\setminus\overline\Omega}|\nabla u(x)|^2\,dx=
\int_{\R^n\setminus\overline\Omega}|\nabla u(x)|^2\,dx,
\end{equation}
thanks to the Dominated Convergence Theorem. 

We claim that
\begin{equation}\label{9485bg8vhfdshgfhghh}
\lim_{R\to+\infty}\int_{\partial B_{R}}u(x)\frac{\partial u}{\partial\nu}(x)\,d{\mathcal{H}}^{n-1}_x=0.
\end{equation}
To prove this, we consider a point~$x_0\in \R^n\setminus\overline\Omega$ such that~$|x_0|\ge 2R$.
In light of Proposition~\ref{95743jhyrgeyPROPiefergerhe}, we have that, for every~$x\in\R^n\setminus B_R$,
\begin{equation}\label{948tv49b8b8vy8y85by854hfjdshfgoikjhgfds}
|u(x)|\le\frac{C}{|x|^{n-2}}\le\frac{C}{R^{n-2}}\end{equation}
and accordingly~$\|u\|_{L^1(B_{R/4}(x_0))}\le C R^2$,
up to renaming~$C>0$. Hence, exploiting the Cauchy's Estimates in Theorem~\ref{CAUESTIMTH} we see that
\begin{equation}\label{948tv49b8b8vy8y85by854hfjdshfgoikjhgfds00}
|\nabla u(x_0)|\le \frac{C\|u\|_{L^1(B_{R/4}(x_0))}}{R^{n+1}}
\le\frac{CR^2}{R^{n+1}}=C R^{1-n},\end{equation}
up to relabeling~$C>0$.

As a consequence of this and~\eqref{948tv49b8b8vy8y85by854hfjdshfgoikjhgfds},
$$
\left|\int_{\partial B_{R}}u(x)\frac{\partial u}{\partial\nu}(x)\,d{\mathcal{H}}^{n-1}_x\right|\le CR^{2-n},
$$
from which~\eqref{9485bg8vhfdshgfhghh} follows.

Plugging~\eqref{s9385v9bsadfaghjfkmkjhtrjhtrjttrj657909877derjv4eolriyu22} and~\eqref{9485bg8vhfdshgfhghh}
into~\eqref{s9385v9bsadfaghjfkmkjhtrjhtrjttrj657909877derjv4eolriyu}, we obtain that
$$ \int_{\R^n\setminus\overline\Omega}|\nabla u(x)|^2\,dx
=-\int_{\partial\Omega}\frac{\partial u}{\partial\nu}(x)\,d{\mathcal{H}}^{n-1}_x.$$
Then the desired result follows by recalling~\eqref{CAB9IJSalskdmnfcx34546fsdQsdfasdAOSL}.
\end{proof}

Retaking the physical motivation discussed before Proposition~\ref{propfoutr65465548789}, we observe that 
at the equilibrium, the boundary of the conductor~$\Omega$
is an equipotential energy, thus, denoting by~$u:\R^n\setminus\Omega\to\R$ the pontential
generated by the conductor, up to a normalization we can suppose that~$u=1$ on~$\partial\Omega$ (and~$u$ is a harmonic function,
recalling the comments after~\eqref{GAMMAFU}, and we normalize this potential to
vanish at infinity). The corresponding electric field is~$-\nabla u$ and, by 
Gau{\ss}' Law the total charge of the conductor can be calculated by taking a large ball~$B_R\Supset\Omega$ and, omitting physical constants,
it is equal to the flow of the electric field through~$\partial B_R$, namely
\begin{eqnarray*}&& Q=-\int_{\partial B_R}\frac{\partial u}{\partial\nu}(x)\,d{\mathcal{H}}^{n-1}
=-\int_{\partial (B_R\setminus\Omega)}\frac{\partial u}{\partial\nu}(x)\,d{\mathcal{H}}^{n-1}
-\int_{\partial\Omega}\frac{\partial u}{\partial\nu}(x)\,d{\mathcal{H}}^{n-1}\\&&\qquad\qquad=-
\int_{ B_R\setminus\Omega }\div(\nabla u(x))\,d{\mathcal{H}}^{n-1}
-\int_{\partial\Omega}\frac{\partial u}{\partial\nu}(x)\,d{\mathcal{H}}^{n-1}=-\int_{\partial\Omega}\frac{\partial u}{\partial\nu}(x)\,d{\mathcal{H}}^{n-1}.
\end{eqnarray*}
Since the potential~$V$ along the surface of the conductor is normalized to~$1$, this
total charge corresponds to the self-capacitance~$\frac{Q}{V}$. In turn, this quantity
corresponds to the capacity as defined in~\eqref{DEFCAP}, thanks to~\eqref{CAB9IJSalskdmnfcx34546fsdQsdfasdAOSL-C1}, thus confirming
the essential coincidence between the mathematical notion of capacity and
the physical concept of self-capacitance.
In this spirit, in the physical jargon the function~$u$ in Proposition~\ref{AM-CAPA}
is often referred to with the name\index{conductor potential} of ``conductor potential''.\medskip

Such a conductor potential, when known, can be efficiently exploited to compute the capacity
of a given set, as exemplified by the following result:

\begin{corollary}\label{CAPAPALLA}
Let~$n\ge3$, $R>0$ and~$c_n$ be the positive constant introduced in~\eqref{COIENEE}. 
Then,
$$ \cAPAC(B_R)=\frac{R^{n-2}}{c_n}.$$
\end{corollary}

\begin{proof} Let~$u(x):=\frac{R^{n-2}}{|x|^{n-2}}$. We have that~$u$ is harmonic in the complement of~$B_R$
and attains value~$1$ along~$\partial B_R$. Thus,
\begin{eqnarray*}\int_{\R^n\setminus\overline{B_R}}|\nabla u(x)|^2\,dx=
{\mathcal{H}}^{n-1}(\partial B_1) \int_R^{+\infty}\left(\frac{(n-2)R^{n-2}}{\rho^{n-1}}\right)^2\,\rho^{n-1}\,d\rho
=n\,| B_1| (n-2)R^{n-2} =\frac{R^{n-2}}{c_n},\end{eqnarray*}
thanks to~\eqref{B1}. The desired result thus follows from~\eqref{CAB9IJSalskdmnfcx34546fsdQsdfasdAOSL}.
\end{proof}

The pointwise values of the conductor potential can be also conveniently exploited
to estimate the value of the capacity of the original set, according to the following
observation:

\begin{proposition}\label{ILBOCONLEDIST}
Let~$n\ge3$, $\Omega\subset\R^n$ be a bounded and open set with boundary
of class~$C^{2,\alpha}$ for some~$\alpha\in(0,1)$,
and~$u$ be the conductor potential of~$\Omega$ as given in Proposition~\ref{AM-CAPA}.

Let~$x\in \R^n\setminus\overline\Omega$. Let~$R>r>0$ be such that~$B_r(x)\subseteq\R^n\setminus\Omega$ and~$B_R(x)\supseteq\Omega$.

Then,
$$ \frac{c_n \cAPAC(\Omega)}{R^{n-2}}\le u(x)\le\frac{c_n \cAPAC(\Omega)}{r^{n-2}},$$
where~$c_n$ is the positive constant in~\eqref{COIENEE}.
\end{proposition}

\begin{proof} Up to a translation, we suppose that~$x=0$. Take~$\Omega'\Supset\Omega$
such that~$0\not\in\Omega'$. We let~$\Gamma$ be
the fundamental solution in~\eqref{GAMMAFU} and exploit Theorem~\ref{KScvbMS:0okrmt4ht-VAR}
with~$\varphi:=u$ in the set~$\R^n \setminus\overline\Omega$.
In this way, we find that\footnote{Notice that one should argue as in the proof of Proposition~\ref{propfoutr65465548789},
by exploiting Theorem~\ref{KScvbMS:0okrmt4ht-VAR} in~$B_R\setminus\overline\Omega$ and then sending~$R\to+\infty$. The convergence is then guaranteed by the estimates in~\eqref{948tv49b8b8vy8y85by854hfjdshfgoikjhgfds} 
and~\eqref{948tv49b8b8vy8y85by854hfjdshfgoikjhgfds00}.}
\begin{equation*}
\int_{\partial\Omega}
\left(\Gamma (x)\frac{\partial u}{\partial\nu}(x)-
\frac{\partial \Gamma}{\partial\nu}(x)\right)\,d{\mathcal{H}}^{n-1}_x
=-u(0).
\end{equation*}
Moreover, by the Divergence Theorem
\begin{equation*}
\int_{\partial\Omega}
\frac{\partial \Gamma}{\partial\nu}(x)\,d{\mathcal{H}}^{n-1}_x=\int_\Omega \Delta\Gamma(x)\,dx
=0.
\end{equation*}
As a consequence,
\begin{equation}
\label{9iujhbSGB9jSmmSLP-erpfe} 
u(0)=-\int_{\partial\Omega}\Gamma (x)\frac{\partial u}{\partial\nu}(x)\,d{\mathcal{H}}^{n-1}_x
.\end{equation}

Also, since~$u$ attains its maximum along~$\partial\Omega$, we know that~$\frac{\partial u}{\partial\nu}\le0$ along~$\partial\Omega$. Therefore, for every~$x\in\partial\Omega$
$$
-\Gamma (x)\frac{\partial u}{\partial\nu}(x)
=-\frac{c_n}{|x|^{n-2}} \frac{\partial u}{\partial\nu}(x)\in
\left[
-\frac{c_n}{R^{n-2}}\,\frac{\partial u}{\partial\nu}(x),\,
-\frac{c_n}{r^{n-2}}\,\frac{\partial u}{\partial\nu}(x)
\right].
$$
This, \eqref{CAB9IJSalskdmnfcx34546fsdQsdfasdAOSL-C1} and~\eqref{9iujhbSGB9jSmmSLP-erpfe} give that
\begin{equation*}
\begin{split}
u(0)\,&=\,
-\int_{\partial\Omega}\Gamma (x)\frac{\partial u}{\partial\nu}(x)\,d{\mathcal{H}}^{n-1}_x
\\&\in\,\left[
-\frac{c_n}{R^{n-2}}\,\int_{\partial\Omega}\frac{\partial u}{\partial\nu}(x)\,d{\mathcal{H}}^{n-1}_x,\,
-\frac{c_n}{r^{n-2}}\,\int_{\partial\Omega}\frac{\partial u}{\partial\nu}(x)\,d{\mathcal{H}}^{n-1}_x
\right]\\&=\,
\left[
\frac{c_n \cAPAC(\Omega)}{R^{n-2}},\,
\frac{c_n \cAPAC(\Omega)}{r^{n-2}}
\right].\qedhere
\end{split}\end{equation*}\end{proof}

The following is a useful property of the capacity as a set function:

\begin{lemma}\label{PTEJSN5678dsdvfg9CALsdfLDM-98uyhnSM6789SMMSA}
Let~$\Omega_1$, $\Omega_2\subset\R^n$ be bounded and open sets with boundary of class~$C^{2,\alpha}$
for some~$\alpha\in(0,1)$, and suppose that~$\Omega_1\cup\Omega_2$ has boundary of class~$C^{2,\alpha}$.

Then,
$$ \cAPAC(\Omega_1)\le\cAPAC(\Omega_1\cup\Omega_2)\le
\cAPAC(\Omega_1)+\cAPAC(\Omega_2).$$
\end{lemma}

\begin{proof} Let~$\e>0$. We use the definition in~\eqref{DEFCAP} 
to find a function~$v_\e\in C^\infty_0(\R^n)$ with~$v_\e=1$ in~$\Omega_1\cup\Omega_2$
such that
$$ \cAPAC(\Omega_1\cup\Omega_2)+\e\ge
\int_{\R^n\setminus\overline{\Omega_1\cup\Omega_2}}|\nabla v_\e(x)|^2\,dx
.$$
Also, since~$v_\e=1$ in~$\Omega_2$, whence~$\nabla v_\e=0$ a.e. in~$\Omega_2$ (see
e.g.~\cite[Theorem~6.19]{MR1817225}), we have that
$$ \int_{\R^n\setminus\overline{\Omega_1\cup\Omega_2}}|\nabla v_\e(x)|^2\,dx=
\int_{\R^n\setminus\overline{\Omega_1}}|\nabla v_\e(x)|^2\,dx.$$
{F}rom these observations and using again the definition in~\eqref{DEFCAP} we find that
$$ \cAPAC(\Omega_1\cup\Omega_2)+\e\ge
\int_{\R^n\setminus\overline{\Omega_1}}|\nabla v_\e(x)|^2\,dx\ge\cAPAC(\Omega_1).$$
By taking~$\e$ as small as we wish, we obtain
\begin{equation} \label{AMS-oIUHGVS-SXCVBNCOHSDNMCCIKSMDIAOIJSHB}\cAPAC(\Omega_1\cup\Omega_2)\ge\cAPAC(\Omega_1).\end{equation}

Now, we make use of the definition in~\eqref{DEFCAP} to find~$v_{1,\e}$, $v_{2,\e}\in C^\infty_0(\R^n)$ with~$v_{1,\e}=1$ in~$\Omega_1$ and~$v_{2,\e}=1$ in~$\Omega_2$
such that
$$ \cAPAC(\Omega_1)+\e\ge\int_{\R^n\setminus\overline{\Omega_1}}|\nabla v_{1,\e}(x)|^2\,dx
\qquad{\mbox{and}}\qquad
\cAPAC(\Omega_1)+\e\ge\int_{\R^n\setminus\overline{\Omega_2}}|\nabla v_{2,\e}(x)|^2\,dx.$$
As in~\eqref{LA-SORROSDTTOSVOLMIANNDKAa3E89SD}, without loss of generality
we can suppose that~$v_{1,\e}(x)$, $v_{2,\e}(x)\in[0,1]$ for all~$x\in\R^n$.
Let~$v^\star_\e:=\max\{v_{1,\e},v_{2,\e}\}$. Then~$v^\star_\e=1$ in~$\Omega_1\cup\Omega_2$.
Furthermore,~$v^\star_\e$ is Lipschitz continuous and compactly supported,
which together with~\eqref{DEFCAP-SOBO} leads to
\begin{eqnarray*}&& 
\int_{\R^n\setminus\overline{\Omega_1}}|\nabla v_{1,\e}(x)|^2\,dx
+\int_{\R^n\setminus\overline{\Omega_2}}|\nabla v_{2,\e}(x)|^2\,dx
=
\int_{\R^n}|\nabla v_{1,\e}(x)|^2\,dx
+\int_{\R^n}|\nabla v_{2,\e}(x)|^2\,dx\\&&\qquad
\ge
\int_{\{v_{1,\e}>v_{2,\e}\}}|\nabla v_{1,\e}(x)|^2\,dx
+\int_{\{v_{1,\e}\le v_{2,\e}\}}|\nabla v_{2,\e}(x)|^2\,dx
=
\int_{\R^n}|\nabla v^\star_\e(x)|^2\,dx\\&&\qquad=
\int_{\R^n\setminus\overline{\Omega_1\cup\Omega_2}}|\nabla v^\star_\e(x)|^2\,dx
\ge\cAPAC(\Omega_1\cup\Omega_2).\end{eqnarray*}
Using these observations, we obtain that
$$ \cAPAC(\Omega_1)+\cAPAC(\Omega_2)+2\e\ge
\cAPAC(\Omega_1\cup\Omega_2)$$
and thus, taking~$\e$ arbitrarily small,
$$ \cAPAC(\Omega_1)+\cAPAC(\Omega_2)\ge
\cAPAC(\Omega_1\cup\Omega_2).$$
{F}rom this\footnote{Though we do not need this property here, we observe that by considering also the function~$\min\{v_{1,\e},v_{2,\e}\}$ this calculation shows that
$$ \cAPAC(\Omega_1) + \cAPAC(\Omega_2) \ge \cAPAC(\Omega_1\cup\Omega_2)+\cAPAC(\Omega_1\cap\Omega_2).$$
In jargon, this says that the capacity is a \index{submodular set function} ``submodular'' set function.}
and~\eqref{AMS-oIUHGVS-SXCVBNCOHSDNMCCIKSMDIAOIJSHB},
the desired result plainly follows.
\end{proof}

A counterpart of the result in Lemma~\ref{PTEJSN5678dsdvfg9CALsdfLDM-98uyhnSM6789SMMSA}
deals with the pointwise value of the conductor potentials:

\begin{lemma}\label{KJSM-098765-TFI-IMP245S-PLSIUDND-94}
Let~$n\ge3$, $\Omega_1$, $\Omega_2\subset\R^n$ be bounded and open sets with boundary of class~$C^{2,\alpha}$
for some~$\alpha\in(0,1)$, and suppose that~$\Omega_1\cup\Omega_2$ has boundary of class~$C^{2,\alpha}$.
For~$i\in\{1,2\}$, let~$u_i$ be the conductor potential of~$\Omega_i$ as given in Proposition~\ref{AM-CAPA}.
Let also~$u_3$ be the conductor potential of~$\Omega_1\cup\Omega_2$.

Then, for every~$x\in\R^n$,
$$ u_3(x)\le u_1(x)+u_2(x).$$
\end{lemma}

\begin{proof} In light of Proposition~\ref{AM-CAPA}, the function~$u:=u_1+u_2-u_3$ is harmonic
outside~$\Omega_1\cup\Omega_2$ and it belongs to~${\mathcal{D}}^{1,2}(\R^n)$. Furthermore,
on~$\partial(\Omega_1\cup\Omega_2)\subseteq(\partial\Omega_1)\cup(\partial\Omega_2)$
we have that either~$u_1=1$ or~$u_2=1$, and therefore~$u\ge1-u_3\ge0$.

As a result, we can apply
Lemma~\ref{SOBOMAXPLEH10} and deduce that~$u\ge0$ in~$\R^n$, from which the desired result follows.
\end{proof}

Importantly, one can define
the capacity of any bounded (not necessarily open, nor with regular boundary) set~$\Omega$ by
setting
\begin{equation}\label{CAPGE} \cAPAC(\Omega):=\inf\cAPAC(\Omega'),\end{equation}
where the infimum above is taken over all the open and bounded sets~$\Omega'$
with boundary of class~$C^{2,\alpha}$, for some~$\alpha\in(0,1)$, such that~$\Omega\subset\Omega'$.

\begin{lemma}\label{MS-0okjn8ijn7uhn5tfvAJM456SOISKJMD}
The definition in~\eqref{CAPGE} coincides with that in~\eqref{DEFCAP}
whenever~$\Omega$ is open and bounded with boundary of class~$C^{2,\alpha}$ for some~$\alpha\in(0,1)$.
\end{lemma}

\begin{proof} To avoid confusion, in this proof we denote by~``$\widetilde{\cAPAC}$''
the capacity defined in the left hand side of~\eqref{CAPGE} and reserve the notation of~``$\cAPAC$''
for the one in~\eqref{DEFCAP}.

Suppose that~$\Omega$ is open and bounded with boundary of class~$C^{2,\alpha}$ for some~$\alpha\in(0,1)$.
Let also~$\Omega'$ be open and bounded with boundary of class~$C^{2,\alpha}$, with~$\Omega\subset\Omega'$. 
We know that~$\cAPAC(\Omega)\le\cAPAC(\Omega')$, therefore, by~\eqref{CAPGE},
by taking the infimum over such sets~$\Omega'$, we have
that~$ \cAPAC(\Omega)\le\inf\cAPAC(\Omega')=\widetilde{\cAPAC}(\Omega)$. On the other hand,
by taking~$\Omega':=\Omega$ as a candidate in the infimum in~\eqref{CAPGE},
we have that~$\widetilde{\cAPAC}(\Omega)\le\cAPAC(\Omega)$. These observations show that the
definitions in~\eqref{DEFCAP} and~\eqref{CAPGE} coincide if~$\Omega$
is open and bounded with boundary of class~$C^{2,\alpha}$.
\end{proof}

Additionally, one can show that the results in
Lemma~\ref{PTEJSN5678dsdvfg9CALsdfLDM-98uyhnSM6789SMMSA} hold true
for every bounded sets~$\Omega_1$, $\Omega_2\subset\R^n$, adopting the definition
in~\eqref{CAPGE}. 
For this, the strategy would be to use the definition
in~\eqref{CAPGE} to find approximating sets~$\Omega_{1,j}$ and~$\Omega_{2,j}$
which are bounded and open with boundary of class~$C^{2,\alpha}$ and
contain~$\Omega_1$ and~$\Omega_2$, respectively. One needs to modify
the sets~$\Omega_{1,j}$ and~$\Omega_{2,j}$ in such a way that~$\Omega_{1,j}\cup\Omega_{2,j}$
has boundary of class~$C^{2,\alpha}$ as well, since in this setting one can exploit
Lemma~\ref{PTEJSN5678dsdvfg9CALsdfLDM-98uyhnSM6789SMMSA}
and then obtain the desired result by taking the limit. The difficulty is that one would like to perform
such modifications to the sets altering the capacity in a controlled way.

To prove that one can find such sets, we first
establish the following observation that relates the capacity of the fattening of a
set with the capacity of the set itself. 

\begin{lemma}\label{dewuityghft435t4609876543345589uhgdsfs}
Let~$n\ge3$, $\Omega$ be a bounded and open set with boundary of class~$C^{2,\alpha}$
for some~$\alpha\in(0,1)$. For every~$\rho\in(0,1)$, let~$\Omega_\rho$
be a bounded and open set with boundary of class~$C^{2,\alpha}$ such that
$$\Omega_\rho \subset \bigcup_{p\in\Omega}B_\rho(p).$$

Then, there exists~$\rho_\Omega\in(0,1)$ such that, for every~$\rho\in(0,\rho_\Omega)$,
$$ \cAPAC(\Omega_\rho)\le \left(1+\rho^{1/2}\right) \cAPAC(\Omega) +C\, \rho^{1/2},$$
for some~$C>0$, depending on~$n$ and~$\Omega$.
\end{lemma}

\begin{proof}
We consider the conductor potential~$u$ for~$\Omega$ as given by
Proposition~\ref{AM-CAPA}, and we take a function~$\varphi\in C^\infty(\R,[0,1])$ such that~$\varphi=1$ in~$(-\infty,1]$
and~$\varphi=0$ in~$(2,+\infty)$. We also recall the definition of the signed distance function~$d_{\partial\Omega}$
to~$\partial\Omega$ given at the beginning of Section~\ref{BELS}
(that is, $d_{\partial\Omega}\ge 0$ in~$\Omega$ and~$d_{\partial\Omega}\le0$ in~$\R^n\setminus\Omega$),
and we set~$d_{\R^n\setminus\Omega}:=-d_{\partial\Omega}$ and
$$\varphi_\rho(x):= \varphi\left(\frac{d_{\R^n\setminus\Omega}(x)}{\rho}
\right).
$$
Notice that
\begin{eqnarray*}
&&{\mbox{if~$d_{\R^n\setminus\Omega}(x)<\rho$, then~$\varphi_\rho(x)=1$,}}\\
&&{\mbox{if~$d_{\R^n\setminus\Omega}(x)>2\rho$, then~$\varphi_\rho(x)=0$}}\\
&&{\mbox{and }} \; |\nabla\varphi_\rho|\le \frac{C}{\rho} \chi_{(\rho,2\rho)}(d_{\R^n\setminus\Omega}), 
\;{\mbox{ for some~$C>0$}}.
\end{eqnarray*}
We also set~$u_\rho:=\varphi_\rho+(1-\varphi_\rho)u$ and we observe that~$u_\rho(x)=1$ if~$x\in\Omega_\rho$
and~$u_\rho(x)=u(x)$ if~$|x|$ is sufficiently large, and therefore~$u_\rho\in{\mathcal{D}}^{1,2}(\R^n)$.
As a consequence, exploiting also the Cauchy-Schwarz Inequality, for~$\e\in(0,1)$ (to be specified later
in the proof),
\begin{eqnarray*}&&
\cAPAC(\Omega_\rho)\le\int_{\R^n}|\nabla u_\rho(x)|^2\,dx=
\int_{\R^n}\big|\nabla\varphi_\rho(x)\big(1-u(x)\big)+\big(1-\varphi_\rho(x)\big)
\nabla u(x)\big|^2\,dx\\&&\qquad\le
\left(1+\frac1\e\right)\int_{\{ d_{\R^n\setminus\Omega}(x)\in(\rho,2\rho)\}}|\nabla\varphi_\rho(x)|^2
\big(1-u(x)\big)^2\,dx
+ (1+\e)\int_{\R^n}\big(1-\varphi_\rho(x)\big)^2|\nabla u(x)|^2\,dx \\&&\qquad\le
\frac{C}{\rho^2}\left(1+\frac1\e\right)\int_{\{ d_{\R^n\setminus\Omega}(x)\in(\rho,2\rho)\}}\big(1-u(x)\big)^2\,dx
+ (1+\e)\int_{\R^n}|\nabla u(x)|^2\,dx \\&&\qquad=
\frac{C}{\rho^2}\left(1+\frac1\e\right)\int_{\{ d_{\R^n\setminus\Omega}(x)\in(\rho,2\rho)\}}\big(1-u(x)\big)^2\,dx
+ (1+\e)\cAPAC(\Omega),
\end{eqnarray*}
where we have used~\eqref{CAB9IJSalskdmnfcx34546fsdQsdfasdAOSL} in the last line.
Now we recall that~$u\in C^{0,1}(\R^n\setminus\Omega)$, and therefore, if~$x\in\{ d_{\R^n\setminus\Omega}(x)\in(\rho,2\rho)\}$
then~$(1-u(x))^2\le C \rho^2$, for some~$C>0$. Accordingly,
$$ \cAPAC(\Omega_\rho)\le
C\rho\left(1+\frac1\e\right)+ (1+\e)\cAPAC(\Omega),$$
up to renaming constants.

Now we choose~$\e:=\sqrt{\frac{C\rho}{\cAPAC(\Omega)}}$ and we conclude that
$$ \cAPAC(\Omega_\rho)\le
C\rho+ \frac{C\rho}\e+ (1+\e)\cAPAC(\Omega)=C\rho+ \cAPAC(\Omega)+\sqrt{C\rho\cAPAC(\Omega)},
$$
up to relabeling~$C>0$.
Thus, by exploiting again the Cauchy-Schwarz Inequality,
$$ \cAPAC(\Omega_\rho)\le C\rho+ \cAPAC(\Omega)+ \sqrt{(C\rho^{1/2})(\rho^{1/2}\cAPAC(\Omega))}
\le C\rho+ \cAPAC(\Omega)+C\rho^{1/2}+ \rho^{1/2}\cAPAC(\Omega),
$$
from which the desired result follows.
\end{proof}

We now employ
Lemma~\ref{dewuityghft435t4609876543345589uhgdsfs}
to modify the approximating sets keeping their capacities under control.
See Figure~\ref{2rerMPM124r32t34ey5rthtfgETHTERRAeFI}
for a sketch of this construction.

\begin{lemma}\label{pew4r83657846547836689438y00}
Let~$n\ge3$, $\Omega_1$, $\Omega_2\subset\R^n$ be bounded sets. For every~$\e\in(0,1)$,
let~$\Omega_{1,\e}$, $\Omega_{2,\e}\subset\R^n$
be bounded and open sets with boundary of class~$C^{2,\alpha}$, for some~$\alpha\in(0,1)$, such
that~$\Omega_i\subset\Omega_{i,\e}$ and with the property that
\begin{equation}\label{pew4r83657846547836689438y}
\lim_{\e\searrow0}\cAPAC(\Omega_{i,\e})=\cAPAC(\Omega_i),\end{equation}
for every~$i\in\{1,2\}$.

Then, for every~$\rho\in(0,1)$, there exist bounded and open
sets~$\Omega_{1,\e,\rho}$, $\Omega_{2,\e,\rho}\subset\R^n$
with boundary of class~$C^{2,\alpha}$ such
that
$$ \Omega_{i,\e}\subset\Omega_{i,\e, \rho}\subset \bigcup_{p\in \Omega_{i,\e}} B_\rho(p)$$ and
\begin{equation}\label{pew4r83657846547836689438y22}
\cAPAC(\Omega_i)\le
\cAPAC(\Omega_{i,\e,\rho})\le \left(1+\rho^{1/2}\right) \cAPAC(\Omega_{i}) +C_\e\, \rho^{1/2}+\phi(\e),
\end{equation}
for some~$C_\e>0$ and some function~$\phi$ such that~$\phi(\e)\to0$
as~$\e\searrow0$, for every~$i\in\{1,2\}$, and~$\Omega_{1,\e,\rho}\cup\Omega_{2,\e,\rho}$ has boundary
of class~$C^{2,\alpha}$.
\end{lemma}

\begin{proof}
It follows from~\eqref{pew4r83657846547836689438y} that
\begin{equation}\label{pewtrsaDFSSGFDHG777777}
\cAPAC(\Omega_{i,\e})\le\cAPAC(\Omega_{i})+\phi(\e) \end{equation}
for some function~$\phi$ such that~$\phi(\e)\to0$
as~$\e\searrow0$, for every~$i\in\{1,2\}$.
Furthermore, for every~$\rho\in(0,1)$, we consider sets~$\Omega_{i,\e,\rho}$ with boundary of class~$C^{2,\alpha}$
such that
$$ \Omega_{i,\e,\rho}\subset \bigcup_{p\in\Omega_{i,\e}}B_\rho(p),
$$
for every~$i\in\{1,2\}$,
and~$\Omega_{1, \e, \rho}\cup\Omega_{2,\e, \rho}$ has boundary of class~$C^{2,\alpha}$.
See Figure~\ref{2rerMPM124r32t34ey5rthtfgETHTERRAeFI}
for a sketch of this construction.
We notice that~$\Omega_{i,\e}\subset\Omega_{i, \e, \rho}$ for every~$i\in\{1,2\}$.

\begin{figure}
  \centering
 \includegraphics[width=.27\linewidth]{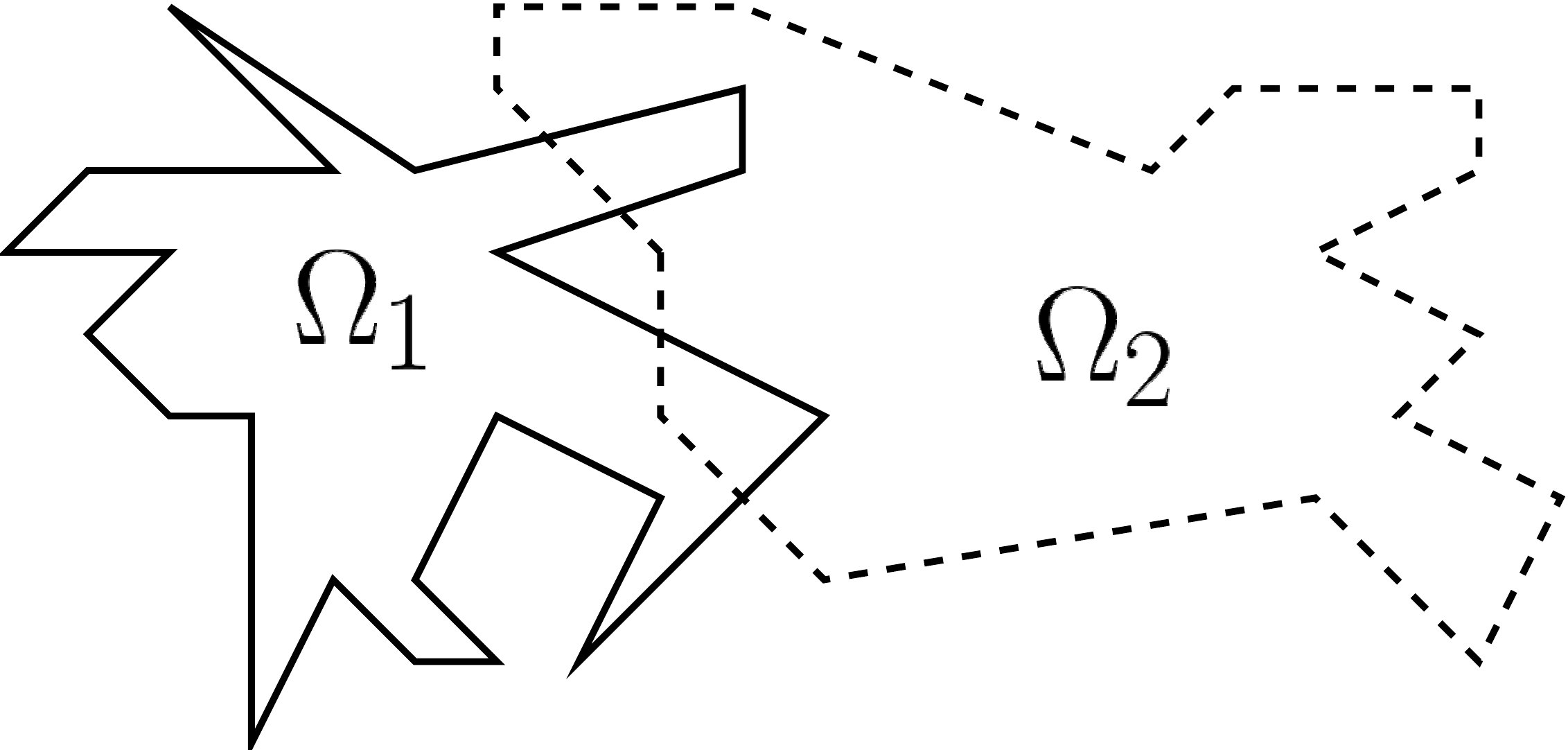}$\qquad$
    \includegraphics[width=.27\linewidth]{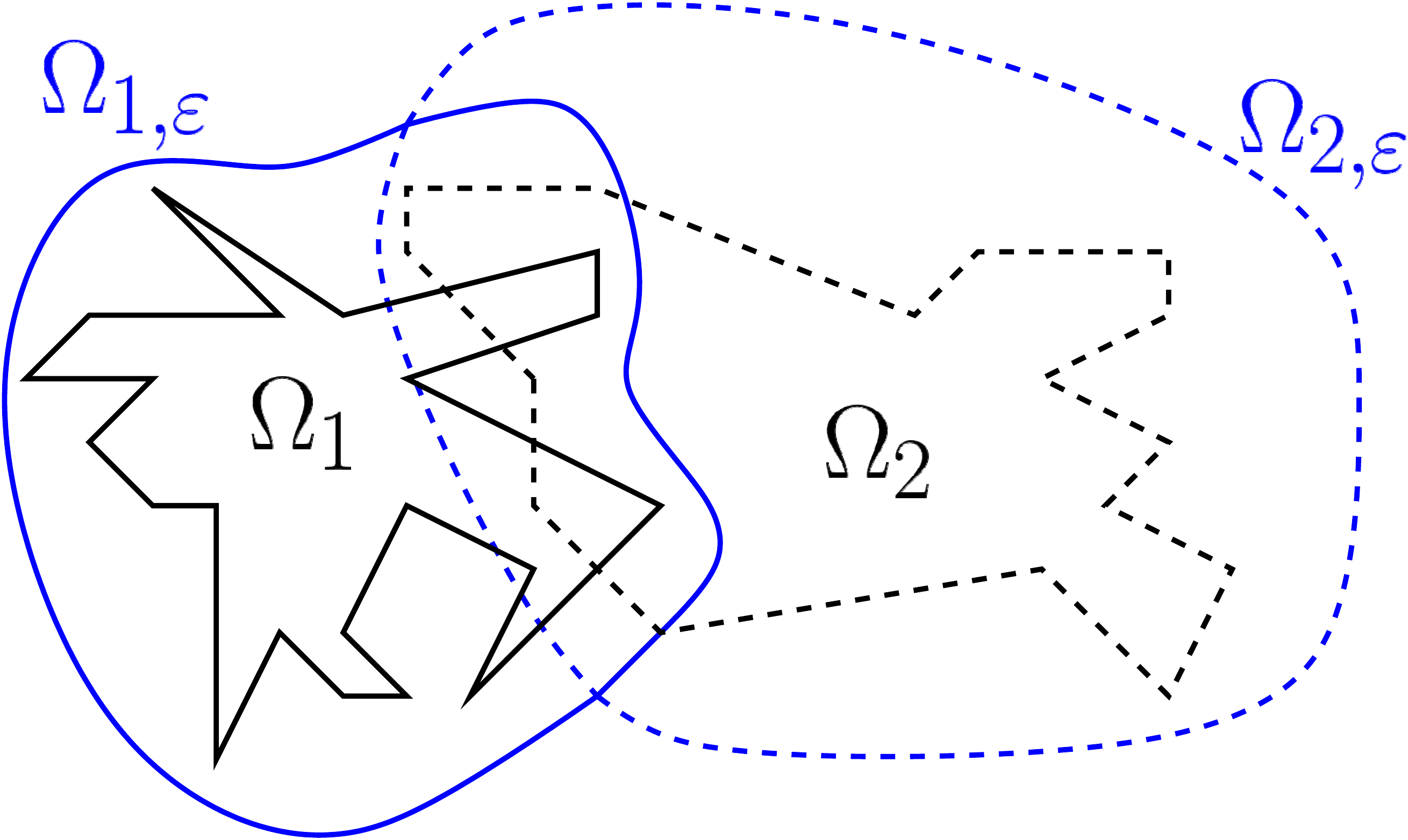}$\qquad$
      \includegraphics[width=.27\linewidth]{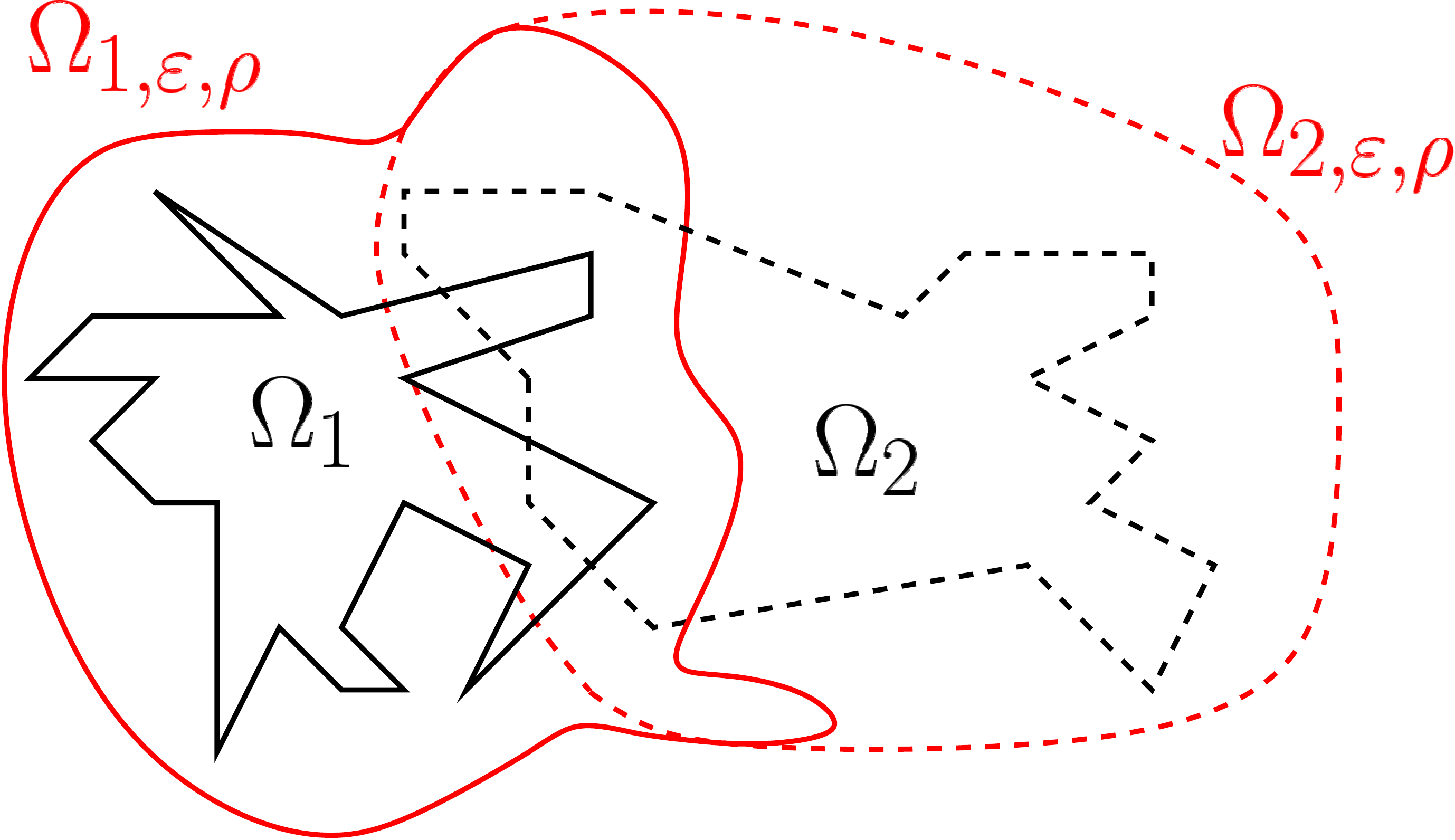}
 \caption{\sl The sets~$\Omega_1$, $\Omega_2$, $\Omega_{1,\e}$, $\Omega_{2,\e}$,
 $\Omega_{1,\e,\rho}$, $\Omega_{2,\e,\rho}$
 involved in the proof of~\eqref{e29rhgreghrereh11111222222}.}\label{2rerMPM124r32t34ey5rthtfgETHTERRAeFI}
\end{figure}

Accordingly, in order to complete the proof of Lemma~\ref{pew4r83657846547836689438y00},
it remains that to prove~\eqref{pew4r83657846547836689438y22}.
For this, we observe that,
in light of Lemma~\ref{dewuityghft435t4609876543345589uhgdsfs}, if~$\rho$ is sufficiently small,
$$ \cAPAC(\Omega_{i,\e,\rho})\le \left(1+\rho^{1/2}\right) \cAPAC(\Omega_{i,\e}) +C_\e\, \rho^{1/2}$$
for some~$C_\e>0$, for every~$i\in\{1,2\}$. This information, together with~\eqref{pewtrsaDFSSGFDHG777777},
gives that
\begin{eqnarray*}
\cAPAC(\Omega_{i,\e,\rho})&\le& \left(1+\rho^{1/2}\right) \big(\cAPAC(\Omega_{i})+\phi(\e)\big) +C_\e\, \rho^{1/2}\\&\le&
\left(1+\rho^{1/2}\right) \cAPAC(\Omega_{i})+C_\e\, \rho^{1/2}+\phi(\e),
\end{eqnarray*}
up to relabeling~$C_\e$, as desired.
\end{proof}

As a consequence of Lemma~\ref{pew4r83657846547836689438y00}, we have the following:

\begin{corollary}\label{coroutile}
Let~$n\ge3$, $\Omega_1$, $\Omega_2\subset\R^n$ be bounded sets. Then, for every~$j\in\N$ there exist
bounded and open sets~$\Omega_{1,j}$ and~$\Omega_{2,j}$ with boundary of class~$C^{2,\alpha}$,
for some~$\alpha\in(0,1)$, such that~$\Omega_i\subset\Omega_{i,j}$
and
\begin{equation*}
\lim_{j\to+\infty}\cAPAC(\Omega_{i,j})=\cAPAC(\Omega_i),\end{equation*}
for every~$i\in\{1,2\}$, and~$\Omega_{1,j}\cup\Omega_{2,j}$ has boundary of class~$C^{2,\alpha}$.
\end{corollary}

\begin{proof}
By the definition in~\eqref{CAPGE}, for every~$\e\in(0,1)$,
there exist bounded and open sets~$\Omega_{1,\e}$, $\Omega_{2,\e}$
with boundary of class~$C^{2,\alpha}$ such
that~$\Omega_i\subset\Omega_{i,\e}$ and with the property that
\begin{equation*}
\lim_{\e\searrow0}\cAPAC(\Omega_{i,\e})=\cAPAC(\Omega_i),\end{equation*}
for every~$i\in\{1,2\}$.

Hence, we are in the position of applying Lemma~\ref{pew4r83657846547836689438y00} to 
find, for every~$\rho\in(0,1)$, bounded and open
sets~$\Omega_{1,\e,\rho}$, $\Omega_{2,\e,\rho}$
with boundary of class~$C^{2,\alpha}$ such
that
$$ \Omega_{i,\e}\subset\Omega_{i,\e, \rho}\subset \bigcup_{p\in \Omega_{i,\e}} B_\rho(p)$$ and
such that~\eqref{pew4r83657846547836689438y22} holds true for every~$i\in\{1,2\}$.
Moreover, $\Omega_{1,\e,\rho}\cup\Omega_{2,\e,\rho}$ has boundary
of class~$C^{2,\alpha}$.

In particular, taking~$\e_j:=1/j$ in~\eqref{pew4r83657846547836689438y22}, we see that, for every~$i\in\{1,2\}$,
$$ \cAPAC(\Omega_i)\le
\cAPAC(\Omega_{i,\e_j,\rho})\le \left(1+\rho^{1/2}\right) \cAPAC(\Omega_{i}) +C_j\, \rho^{1/2} +\phi\left(\frac1j\right),$$
for some~$C_j>0$, with~$\phi\left(\frac1j\right)\to0$ as~$j\to+\infty$. Now, choosing
$$\rho_j:=\min\left\{\frac1j,\,\frac1{jC_j^2 }\right\},$$
we obtain that
$$ \cAPAC(\Omega_i)\le
\cAPAC(\Omega_{i,\e_j,\rho_j})\le \left(1+\frac1{\sqrt{j}}\right) \cAPAC(\Omega_{i}) +\frac1{\sqrt{j}} +\phi\left(\frac1j\right),$$
for every~$i\in\{1,2\}$. Consequently,
$$ \lim_{j\to+\infty}\cAPAC(\Omega_{i,\e_j,\rho_j})=\cAPAC(\Omega_i),$$
for every~$i\in\{1,2\}$. Thus, setting~$\Omega_{i,j}:=\Omega_{i,\e_j,\rho_j}$ for every~$i\in\{1,2\}$, we
obtain the desired result.
\end{proof}

With these observations at hand, we can now prove the following:

\begin{lemma}\label{lemakoe84654656548489}
The results in
Lemma~\ref{PTEJSN5678dsdvfg9CALsdfLDM-98uyhnSM6789SMMSA} hold true
for every bounded sets~$\Omega_1$, $\Omega_2\subset\R^n$, adopting the definition
in~\eqref{CAPGE}.
\end{lemma}

\begin{proof}
Let now~$\Omega_1$, $\Omega_2\subset\R^n$ be bounded.
Since~$\Omega_1\subset\Omega_1\cup\Omega_2$, by the definition
in~\eqref{CAPGE} we have that~$\cAPAC(\Omega_1)\le\cAPAC(\Omega_1\cup\Omega_2)$.

Hence, to complete the proof of Lemma~\ref{lemakoe84654656548489}, we have to show that
\begin{equation}\label{e29rhgreghrereh11111222222}
\cAPAC(\Omega_1\cup\Omega_2)\le \cAPAC(\Omega_1)+\cAPAC(\Omega_2).
\end{equation}
For this, we let~$\e>0$ and we see that there exist bounded and open
sets~$\Omega_{1,\e}$ and~$\Omega_{2,\e}$ with boundary of class~$C^{2,\alpha}$
such that~$\Omega_i\subset\Omega_{i,\e}$ and
\begin{equation*}
\lim_{\e\searrow0}\cAPAC(\Omega_{i,\e})=\cAPAC(\Omega_{i}) \end{equation*} for every~$i\in\{1,2\}$.
Furthermore, from Lemma~\ref{pew4r83657846547836689438y00}, we have that,
for every~$\rho\in(0,1)$, 
there exist bounded and open
sets~$\Omega_{1,\e,\rho}$, $\Omega_{2,\e,\rho}\subset\R^n$
with boundary of class~$C^{2,\alpha}$ such
that~$\Omega_{i,\e}\subset\Omega_{i,\e, \rho}$ and
\begin{equation}\label{pew4r83657846547836689438y22FIN}
\cAPAC(\Omega_{i,\e,\rho})\le \left(1+\rho^{1/2}\right) \cAPAC(\Omega_{i}) +C_\e\, \rho^{1/2}+\phi(\e),
\end{equation}
for some~$C_\e>0$, for every~$i\in\{1,2\}$, and~$\Omega_{1,\e,\rho}\cup\Omega_{2,\e,\rho}$ has boundary
of class~$C^{2,\alpha}$.

{F}rom this, and
exploiting Lemma~\ref{PTEJSN5678dsdvfg9CALsdfLDM-98uyhnSM6789SMMSA},
we obtain that
\begin{eqnarray*}
\cAPAC(\Omega_1)+\cAPAC(\Omega_2)&\ge&
\frac{\cAPAC(\Omega_{1,\e,\rho})+\cAPAC(\Omega_{2,\e,\rho})-2C_\e \rho^{1/2}-2\phi(\e) }{1+\rho^{1/2}}\\&\ge&
\frac{\cAPAC(\Omega_{1,\e,\rho}\cup\Omega_{2,\e,\rho})-2C_\e \rho^{1/2}-2\phi(\e) }{1+\rho^{1/2}}\\&\ge&
\frac{\cAPAC(\Omega_{1}\cup\Omega_{2})-2C_\e \rho^{1/2}-2\phi(\e) }{1+\rho^{1/2}}.
\end{eqnarray*}
Accordingly, sending~$\rho\searrow0$,
$$\cAPAC(\Omega_1)+\cAPAC(\Omega_2) +2\phi(\e )\ge 
\cAPAC(\Omega_{1}\cup\Omega_{2})$$
and sending~$\e\searrow0$ we obtain~\eqref{e29rhgreghrereh11111222222}, as desired.
\end{proof}

It is interesting to observe that the notion of
conductor potential introduced in
Proposition~\ref{AM-CAPA} can also be extended to arbitrary bounded sets
by leveraging the setting in~\eqref{CAPGE}. In this situation,
the conductor potential is not necessarily continuous at the boundary of the set,
yet it maintains useful harmonicity and energy properties:

\begin{corollary}\label{gujbmSIJKc-qpwruofygheibdowi3egufewbuigf3289eufgvosud81ye2rit23456gh}
Let~$n\ge3$
and~$\Omega\subset\R^n$ be a bounded set.

Then, there exists a unique~$u\in C^2(\R^n\setminus\overline\Omega)\cap {\mathcal{D}}^{1,2}(\R^n)$ such that
\begin{equation}\label{LK:0o0o2krndwk35ykhdf} \begin{dcases}
\Delta u=0 \quad {\mbox{ in~$\,\R^n\setminus\overline\Omega$,}}\\
u=1 \quad{\mbox{ a.e. in $\Omega$.}}
\end{dcases}\end{equation}
Additionally, $0\le u(x)\le1$ for a.e.~$x\in\R^n$ and
\begin{equation}\label{CAHBSN123s4345453124135dcHNSdfz2354tyhdq3z4erifjg-1} \cAPAC(\Omega)=\int_{\R^n\setminus\overline\Omega}|\nabla u(x)|^2\,dx.\end{equation}
Furthermore,
if~$x\in \R^n\setminus\overline\Omega$ and~$R>r>0$ are such that~$B_r(x)\subseteq\R^n\setminus\Omega$ and~$B_R(x)\supseteq\Omega$,
then,
\begin{equation}\label{CAHBSN123s4345453124135dcHNSdfz2354tyhdq3z4erifjg-2} \frac{c_n \cAPAC(\Omega)}{R^{n-2}}\le u(x)\le\frac{c_n \cAPAC(\Omega)}{r^{n-2}},\end{equation}
where~$c_n$ is the positive constant in~\eqref{COIENEE}.
\end{corollary}

\begin{proof} For every~$j\in\N$, we use~\eqref{CAPGE} to pick a set~$\Omega_j$
which is open, bounded and
with boundary of class~$C^{2,\alpha}$ such that~$\Omega\Subset\Omega_j$ and
\begin{equation}\label{0oikj8i24ks9ihf94i23203762423}
\lim_{j\to+\infty}
\cAPAC(\Omega_j)=
\cAPAC(\Omega).\end{equation}
Without loss of generality, we can suppose that
\begin{equation}\label{6.1.30BIS}
\Omega_j\subset \bigcup_{p\in\Omega} B_{1/j}(p),
\end{equation}
in light of Lemma~\ref{PTEJSN5678dsdvfg9CALsdfLDM-98uyhnSM6789SMMSA}, since reducing the domain
also reduces the capacity.

Furthermore, we claim that
\begin{equation}\label{e346324563562566757000000kipoioioioioi}
{\mbox{without loss of generality, 
we can suppose that~$\Omega_{j+1}\subset\Omega_j$.}}\end{equation}
To prove it, we take~$\rho_j\in(0,1)$ (to be specified later in the proof)
and we consider bounded and open sets~$\Omega^{\#}_j$ with boundary of class~$C^{2,\alpha}$ and such that
$$ \bigcup_{p\in\Omega_j} B_{\rho_j/2}(p)\subset
\Omega^{\#}_j \subset  \bigcup_{p\in\Omega_j} B_{\rho_j}(p).
$$
Now we define~$\widetilde\Omega_j:=\Omega_1 \cap \cdots\cap \Omega_j$
and we notice that
\begin{equation}\label{qwertyuioplkjhgfdsazxcvbnmmnbvcxzasdfghjk0}
\widetilde\Omega_{j+1}\subset\widetilde\Omega_j.\end{equation} Furthermore,
we consider sets~$\widehat\Omega_j$ with boundary of class~$C^{2,\alpha}$ and such that
\begin{equation}\label{qwertyuioplkjhgfdsazxcvbnmmnbvcxzasdfghjk}
\bigcup_{p\in\widetilde\Omega_j} B_{\rho_j/4}(p)\subset
\widehat\Omega_j \subset  \bigcup_{p\in\widetilde\Omega_j} B_{\rho_j/2}(p).
\end{equation}
We remark that~$\widehat\Omega_j\subset\Omega^{\#}_j$. Using this and
Lemma~\ref{dewuityghft435t4609876543345589uhgdsfs}, we see that
\begin{eqnarray*}
\cAPAC(\widehat\Omega_j)\le \cAPAC(\Omega^{\#}_j)\le
\left(1+\rho_j^{1/2}\right) \cAPAC(\Omega_j) +C_j\, \rho_j^{1/2},
\end{eqnarray*}
for some~$C_j>0$.
Thus, choosing~$\rho_0:=1$ and
\begin{equation}\label{qwertyuioplkjhgfdsazxcvbnmmnbvcxzasdfghjk1}
\rho_j:=\min\left\{\frac1j,\,\frac1{jC_j^2},\,\frac{\rho_{j-1}}{10}
\right\},\end{equation}
we obtain that
$$ \cAPAC(\widehat\Omega_j)\le
\left(1+\frac1{\sqrt{j}}\right) \cAPAC(\Omega_j) +\frac1{\sqrt{j}}.$$
As a consequence of this and~\eqref{0oikj8i24ks9ihf94i23203762423},
$$  \lim_{j\to+\infty}
\cAPAC(\widehat\Omega_j)=
\cAPAC(\Omega).$$
Hence, to complete the proof of~\eqref{e346324563562566757000000kipoioioioioi},
it remains to check that
\begin{equation*}
\widehat\Omega_{j+1}\subset\widehat\Omega_j.
\end{equation*}
Indeed, by~\eqref{qwertyuioplkjhgfdsazxcvbnmmnbvcxzasdfghjk0},
\eqref{qwertyuioplkjhgfdsazxcvbnmmnbvcxzasdfghjk} and~\eqref{qwertyuioplkjhgfdsazxcvbnmmnbvcxzasdfghjk1},
$$\widehat\Omega_{j+1}\subset \bigcup_{p\in\widetilde\Omega_{j+1}} B_{\rho_{j+1}/2}(p)
\subset \bigcup_{p\in\widetilde\Omega_{j}} B_{\rho_{j+1}/2}(p)
\subset \bigcup_{p\in\widetilde\Omega_{j}} B_{\rho_{j}/20}(p)\subset\widehat\Omega_j,
$$
as desired.

Accordingly, from now on we assume that~$\Omega_{j+1}\subset\Omega_j$.
Moreover, we can also suppose that there exists~$R>0$ such that~$\Omega_j\subset B_R$ for every~$j\in\N$,
and therefore, if~$u_j$ denotes the conductor potential of the set~$\Omega_j$, for every~$j\in\N$,
$$ \int_{\R^n\setminus\overline\Omega}|\nabla u_{j}(x)|^2\,dx =\cAPAC(\Omega_j)
\le\cAPAC(B_R).$$
Consequently, up to a subsequence, we can assume that
\begin{equation}\label{CONJCONDENS}\begin{split}&
{\mbox{$\nabla u_j$ converges weakly in~$L^2(\R^n)$,}}\\&{\mbox{and
$u_j$ converges strongly in~$L^2_{\rm loc}(\R^n)$
and a.e. in~$\R^n$ to a function~$u$,}}\end{split}\end{equation} which is harmonic in~$\R^n\setminus\overline\Omega$
thanks to Corollary~\ref{S-OS-OSSKKHEPARHA}.

Furthermore, we have that~$0\le u\le1$, since the same holds for~$u_j$, and~$u=1$ a.e. in~$\Omega$, thanks to the
pointwise convergence in~\eqref{CONJCONDENS}. 

Now, to prove~\eqref{CAHBSN123s4345453124135dcHNSdfz2354tyhdq3z4erifjg-1}, for every~$j$, $k\in\N$,
we define~$v_{j,k}:=\frac{u_{j+k}+u_j}2$ and we notice that~$v_{j,k}=1$ in~$\Omega_{j+k}$. Accordingly,
\begin{eqnarray*}&&
\cAPAC(\Omega_{j+k})\le \int_{\R^n}|\nabla v_{j,k}(x)|^2\,dx=\frac14\int_{\R^n}\big(
|\nabla u_{j+k}(x)|^2 +|\nabla u_j(x)|^2 +2\nabla u_{j+k}(x)\cdot\nabla u_j(x)\big)\,dx\\
&&\qquad\le \frac12\int_{\R^n}\big(|\nabla u_{j+k}(x)|^2 +|\nabla u_j(x)|^2 \big)\,dx=
\frac12\big(\cAPAC(\Omega_{j+k})+\cAPAC(\Omega_j)\big).
\end{eqnarray*}
In light of~\eqref{0oikj8i24ks9ihf94i23203762423}, this gives that
\begin{equation}\label{dwueryt5y345lllllll}
\lim_{j\to+\infty}\int_{\R^n}|\nabla v_{j,k}(x)|^2\,dx=\cAPAC(\Omega).
\end{equation}
Also, we observe that
\begin{eqnarray*}&&\frac14
\int_{\R^n} |\nabla u_{j+k}(x)-\nabla u_j(x)|^2 \,dx =\frac14\int_{\R^n}\big(
|\nabla u_{j+k}(x)|^2 +|\nabla u_j(x)|^2 -2\nabla u_{j+k}(x)\cdot\nabla u_j(x)\big)\,dx\\&&\qquad\qquad=\frac14
\int_{\R^n}\big(2
|\nabla u_{j+k}(x)|^2 +2|\nabla u_j(x)|^2 -4|\nabla v_{j,k}(x)|^2\big)\,dx\\&&\qquad\qquad=\frac12\big(\cAPAC(\Omega_{j+k})
+\cAPAC(\Omega_j)\big)-\int_{\R^n}|\nabla v_{j,k}(x)|^2\,dx.
\end{eqnarray*}
{F}rom this, \eqref{0oikj8i24ks9ihf94i23203762423} and~\eqref{dwueryt5y345lllllll}, we obtain that
\begin{equation}\label{lskworeiteriyter439548795470}
{\mbox{$\nabla u_j$
is a Cauchy sequence in~$L^2(\R^n)$, and therefore it strongly converges.}}\end{equation}
Consequently, recalling also~\eqref{0oikj8i24ks9ihf94i23203762423},
$$ \cAPAC(\Omega)=\lim_{j\to+\infty}
\cAPAC(\Omega_{j})=
\lim_{j\to+\infty}\int_{\R^n\setminus\overline\Omega}|\nabla u_{j}(x)|^2\,dx
=\int_{\R^n\setminus\overline\Omega}|\nabla u(x)|^2\,dx,
$$
which establishes~\eqref{CAHBSN123s4345453124135dcHNSdfz2354tyhdq3z4erifjg-1}, as desired.

We now prove the uniqueness statement. For this, we set
$$ X_0:=\left\{
v\in C^\infty_0(\R^n)\; {\mbox{ such that there exists an open set~$\Omega'$ such that~$\Omega\Subset\Omega'$
and~$v=1$ in~$\Omega'$}}
\right\} $$
and we define~$X$ to be the closure of~$X_0$ in~${\mathcal{D}}^{1,2}(\R^n)$. We claim that
\begin{equation}\label{uspaziogiusto}
u\in X.
\end{equation}
To prove it, we recall~\eqref{CONJCONDENS} and we notice that~$u_j\in C^{2,\alpha}(\R^n\setminus\Omega_j)
\cap{\mathcal{D}}^{1,2}(\R^n\setminus\Omega_j)\cap L^{\frac{2n}{n-2}}(\R^n\setminus\Omega_j)$, thanks to
the Gagliardo-Nirenberg-Sobolev Inequality, see e.g.~\cite[Theorem~12.4]{MR2527916}.

We take~$R>0$ such that~$\overline\Omega\subset\Omega_j\subset B_R$ and we define
a function~$\tau\in C^\infty_0(\R^n,[0,1])$
such that~$\tau=1$ in~$B_{1}$ and~$\tau=0$ in~$\R^n\setminus B_{2}$.
We set~$\tau_{R}(x):=\tau\left(\frac{x}R\right)$ and, for every~$j\in\N$, $a_j:=\tau_{2R} u_j$
and~$b_j:=(1-\tau_{2R})u_j$.
Notice that~$a_j\in C^{2,\alpha}_0(\R^n\setminus\overline\Omega_j)$.
Thus, given~$\e>0$, taking a mollifier~$\eta_\e$ and the corresponding mollification~$a_{j,\e}:=a_j*\eta_\e$,
we have that there exists~$\e_j$ sufficiently small such that
\begin{equation}\label{dwewrw3500099878787877676776756564534}
a_{j,\e_j}\in X_0.\end{equation}
Indeed, $a_{j,\e}\in C^\infty_0(\R^n)$. Also,
\begin{equation}\label{posdewtiovbn5yufghjhjgfdt60}
{\mbox{if $y\in\Omega_j\subset B_R$, then~$a_j(y)=\tau_{2R}(y)u_j(y)=u_j(y)=1$.}}\end{equation}
We denote by~$d_j$ the distance between~$\Omega$ and~$\partial\Omega_j$,
that is, $$d_j:=\inf_{{x\in\Omega}\atop{y\in\partial\Omega_j}} |x-y|
=\min_{{x\in\overline\Omega}\atop{y\in\partial\Omega_j}} |x-y|,$$
and we observe that~$d_j>0$. We
take an open set~$\widehat\Omega_j$ with boundary of class~$C^{2,\alpha}$ such that
$$\Omega\Subset\widehat\Omega_j\subset
\left\{x\in\Omega_j \,{\mbox{ s.t. }}\, d_{\partial\Omega_j}(x)>\frac{d_j}2\right\}.
$$
Here~$d_{\partial\Omega_j}$ is the signed distance function to~$\partial\Omega_j$
given at the beginning of Section~\ref{BELS}
(that is, $d_{\partial\Omega_j}\ge 0$ in~$\Omega_j$ and~$d_{\partial\Omega_j}\le0$ in~$\R^n\setminus\Omega_j$).
In this way, taking~$\e\in(0,d_j/4)$, if~$x\in\widehat\Omega_j$ and~$y\in B_{\e}(x)$, then~$y\in\Omega_j$, and therefore,
using~\eqref{posdewtiovbn5yufghjhjgfdt60},
$$ a_{j,\e}(x)=\int_{B_\e(x)}a_j(y)\eta_\e(x-y)\,dy = \int_{B_\e(x)}\eta_\e(x-y)\,dy=1.$$
This establishes~\eqref{dwewrw3500099878787877676776756564534}.

We also have that
\begin{equation}\label{e3453489v54368b5654n76786}
{\mbox{$a_{j,\e}$ converges to~$a_j$ in~${\mathcal{D}}^{1,2}(\R^n)$ as~$\e\searrow0$,}}\end{equation}
see e.g. Theorem~9.6 in~\cite{MR3381284}.

Furthermore, since~$b_j\in {\mathcal{D}}^{1,2}(\R^n)$, there exists~$b_{j,k}\in C^{\infty}_0(\R^n)$ such that
\begin{equation}\label{e3453489v54368b5654n7678622}
{\mbox{$b_{j,k}$
converges to~$b_j$ in~${\mathcal{D}}^{1,2}(\R^n)$ as~$k\to+\infty$.}}\end{equation}
Without loss of generality,
we can suppose that
\begin{equation}\label{lakshs4343554}
{\mbox{$b_{j,k}=0$ in~$B_R\supset \overline\Omega$.}}\end{equation}
Indeed, if not, define~$\widetilde{b}_{j,k}:=(1-\tau_{R})b_{j,k}$ and observe that~$\widetilde{b}_{j,k}\in C^{\infty}_0(\R^n)$
and~$\widetilde{b}_{j,k}=0$ in~$B_R\supset\overline\Omega$.
Also, $\widetilde{b}_{j,k}-b_j=b_{j,k}-b_j$
in~$\R^n\setminus B_{2R}$ and~$\widetilde{b}_{j,k}-b_j=(1-\tau_{R})b_{j,k}$
in~$B_{2R}$, and so
\begin{eqnarray*}&&
\int_{\R^n} |\nabla (\widetilde{b}_{j,k}-b_j)(x)|^2\,dx= \int_{B_{2R}}|\nabla (\widetilde{b}_{j,k}-b_j)(x)|^2\,dx
+\int_{\R^n\setminus B_{2R}}|\nabla (\widetilde{b}_{j,k}-b_j(x))|^2\,dx\\&&\qquad=
\int_{B_{2R}\setminus B_R}|\nabla \big((1-\tau_{R}(x))b_{j,k}\big)(x)|^2\,dx
+\int_{\R^n\setminus B_{2R}}|\nabla ({b}_{j,k}-b_j)(x)|^2\,dx\\&&\qquad\le 2\left(
\int_{B_{2R}\setminus B_R}|\nabla \big(1-\tau_{R}(x)\big)|^2|b_{j,k}(x)|^2\,dx +
\int_{B_{2R}\setminus B_R}|\nabla b_{j,k}(x)|^2\,dx\right)\\&&\qquad\qquad
+\int_{\R^n\setminus B_{2R}}|\nabla ({b}_{j,k}-b_j)(x)|^2\,dx\\&&\qquad \le
2\left(
\int_{B_{2R}\setminus B_R}|\nabla \big(1-\tau_{R}(x)\big)|^2|b_{j,k}(x)|^2\,dx +
\int_{\R^n}|\nabla ({b}_{j,k}-b_j)(x)|^2\,dx\right)\\&&\qquad\le
C\left(\frac1{R^2}
\int_{B_{2R}\setminus B_R}|b_{j,k}(x)-b_j(x)|^2\,dx +
\int_{\R^n}|\nabla ({b}_{j,k}-b_j)(x)|^2\,dx\right),
\end{eqnarray*}
for some~$C>0$.
Hence, thanks to~\eqref{e3453489v54368b5654n7678622}, we deduce from this computation that~$\widetilde{b}_{j,k}$
converges to~$b_j$ in~${\mathcal{D}}^{1,2}(\R^n)$ as~$k\to+\infty$, and this completes the proof of~\eqref{lakshs4343554}.

Now, in light of~\eqref{e3453489v54368b5654n76786} and~\eqref{e3453489v54368b5654n7678622}, 
we can take~$\e_j$ sufficiently small
(possibly smaller than~$\e_j$ in~\eqref{dwewrw3500099878787877676776756564534}), and~$k_j$
sufficiently large such that
\begin{eqnarray*}&&
\int_{\R^n}\big|\nabla \big(a_{j,\e_j}-a_j\big)(x)\big|^2\,dx\le\frac1j \\
{\mbox{and }} && \int_{\R^n}\big|\nabla \big(b_{j,k_j}-b_j\big)(x)\big|^2\,dx\le\frac1j.
\end{eqnarray*}
As a consequence, noticing that~$a_j+b_j=u_j$, if~$\e_j$ sufficiently small and~$k_j$ is
sufficiently large,
$$\int_{\R^n}\big|\nabla \big(a_{j,\e_j}+b_{j,k_j}-u_j\big)(x)\big|^2\,dx\le\frac4j.$$
{F}rom this and~\eqref{lskworeiteriyter439548795470}, we conclude that
\begin{equation}\label{e4ut854utewoti5497654iuyportiumyb65uoi65ujkgfj}\begin{split}
& \lim_{j\to+\infty}\sqrt{\int_{\R^n}\big|\nabla \big(a_{j,\e_j}+b_{j,k_j}-u\big)(x)\big|^2\,dx}\\&\qquad
\le\lim_{j\to+\infty}\left(
\sqrt{\int_{\R^n}\big|\nabla \big(a_{j,\e_j}+b_{j,k_j}-u_j\big)(x)\big|^2\,dx}+
\sqrt{\int_{\R^n}\big|\nabla \big(u_j-u\big)(x)\big|^2\,dx}\right)=0.\end{split}\end{equation}
Moreover, from~\eqref{dwewrw3500099878787877676776756564534} and~\eqref{lakshs4343554}, it follows
that~$a_{j,\e_j}+b_{j,k_j}\in X_0$. This fact and~\eqref{e4ut854utewoti5497654iuyportiumyb65uoi65ujkgfj}
establish~\eqref{uspaziogiusto}. 

Now we take~$u_1$, $u_2\in C^2(\R^n\setminus\overline\Omega)\cap {\mathcal{D}}^{1,2}(\R^n)$ to be solutions
of~\eqref{LK:0o0o2krndwk35ykhdf}. Thanks to~\eqref{uspaziogiusto}, we have that~$u_1$, $u_2\in X$, and therefore
there exist functions~$\varphi_{1,j}$, $\varphi_{2,j}\in C^\infty_0(\R^n)$ such that~$\varphi_{1,j}$ and~$\varphi_{2,j}$
converge to~$u_1$ and~$u_2$ in~${\mathcal{D}}^{1,2}(\R^n)$ as~$j\to+\infty$, respectively, and open
sets~$\Omega'_{1,j}$ and~$\Omega'_{2,j}$ such that~$\Omega\Subset\Omega'_{1,j}\cap\Omega'_{2,j}$
and~$\varphi_{1,j}=1$ in~$\Omega'_{1,j}$ and~$\varphi_{j,2}=1$ in~$\Omega'_{2,j}$.
We set~$\widetilde{u}:=u_1-u_2$ and~$\varphi_j:=\varphi_{1,j}-\varphi_{2,j}$, and we notice that~$\varphi_j\in C^\infty_0(\R^n)$
and~$\varphi_j$ converges to~$\widetilde u$ in~${\mathcal{D}}^{1,2}(\R^n)$ as~$j\to+\infty$. Also, setting~$\Omega'_j:=
\Omega'_{1,j}\cap\Omega'_{2,j}$, we have that~$\Omega'_j$ is an open set and~$\Omega\Subset\Omega'_j$. Moreover,
$\varphi_j=0$ in~$\Omega'_j$. We also remark that~$\widetilde{u}$ is harmonic in~$\R^n\setminus\overline\Omega$.

Now, we take an open set~$\widetilde\Omega_j$ with boundary of class~$C^{2,\alpha}$ such that~$\Omega\Subset
\widetilde\Omega_j\Subset\Omega'_j$, and 
we exploit the
first Green's Identity in~\eqref{GRr1} to see that
$$
\int_{\R^n\setminus\overline\Omega} \nabla \varphi_j(x)\cdot\nabla \widetilde u(x) \,dx=
\int_{\partial\widetilde\Omega_j} \varphi_j(x)\frac{\partial \widetilde u}{\partial\nu}(x)\,d{\mathcal{H}}^{n-1}_x=0.$$
Consequently, passing to the limit in~$j$,
$$ \int_{\R^n\setminus\overline\Omega} |\nabla \widetilde{u}(x)|^2\,dx=0,$$
from which it follows that~$\widetilde u$ is constant. Since~$\widetilde u\in L^{\frac{2n}{n-2}}( \R^n\setminus\overline\Omega)$
(thanks to
the Gagliardo-Nirenberg-Sobolev Inequality, see e.g.~\cite[Theorem~12.4]{MR2527916}), this gives that~$\widetilde u$ vanishes identically,
whence~$u_1$ coincides with~$u_2$ and the uniqueness statement is thereby proved.

Now we focus on the proof of~\eqref{CAHBSN123s4345453124135dcHNSdfz2354tyhdq3z4erifjg-2}.
To this end, we employ Proposition~\ref{ILBOCONLEDIST} and we see that,
if~$x\in \R^n\setminus\overline{\Omega_j}$ and~$B_{r_j}(x)\subseteq\R^n\setminus\Omega_j$ and~$B_{R_j}(x)\supseteq\Omega_j$, then
\begin{equation}\label{UAHJ-thremprmnfvaidEDSjnta} \frac{c_n \cAPAC(\Omega_j)}{R_j^{n-2}}\le u_j(x)\le\frac{c_n \cAPAC(\Omega_j)}{r_j^{n-2}}.\end{equation}
Now we take a set~${\mathcal{Z}}\subset \R^n$ with null Lebesgue measure and
such that~$u_j\to u$ in~$ \R^n\setminus{\mathcal{Z}}$. We pick~$x\in\R^n\setminus(\overline{\Omega}\cup{\mathcal{Z}})$
and consider~$B_{r}(x)\subseteq\R^n\setminus\Omega$ and~$B_{R}(x)\supseteq\Omega$.
Moreover, in view of~\eqref{6.1.30BIS}, we have that~$B_{r_j}(x)\subseteq\R^n\setminus\Omega_j$ and~$B_{R_j}(x)\supseteq\Omega_j$,
where~$r_j:=r-\frac2j$ and~$R_j:=R+\frac2j$. In this setting we can exploit~\eqref{UAHJ-thremprmnfvaidEDSjnta} and then send~$j\to+\infty$. In this way, we conclude that~\eqref{CAHBSN123s4345453124135dcHNSdfz2354tyhdq3z4erifjg-2} holds true for all~$x\in
\R^n\setminus(\overline{\Omega}\cup{\mathcal{Z}})$, that is a.e. in~$\R^n\setminus\overline{\Omega}$. Since~$u$ is harmonic, hence continuous
in~$\R^n\setminus\overline{\Omega}$, this guarantees that~\eqref{CAHBSN123s4345453124135dcHNSdfz2354tyhdq3z4erifjg-2} holds true for all~$x\in
\R^n\setminus\overline{\Omega}$, as desired.
\end{proof}

In the spirit of Lemma~\ref{lemakoe84654656548489} we also have:

\begin{lemma}\label{KJSM-098765-TFI-IMP245S-PLSIUDND-95}
The result in
Lemma~\ref{KJSM-098765-TFI-IMP245S-PLSIUDND-94} holds true a.e. in~$\R^n$
for every bounded sets~$\Omega_1$, $\Omega_2\subset\R^n$, adopting the definition
in~\eqref{CAPGE}.
\end{lemma}

\begin{proof} 
Let~$u_1$, $u_2$ and~$u_3$ be the conductor potentials relative to the sets~$\Omega_1$, $\Omega_2$
and~$\Omega_1\cup\Omega_2$, respectively, according to Corollary~\ref{gujbmSIJKc-qpwruofygheibdowi3egufewbuigf3289eufgvosud81ye2rit23456gh}.

Moreover, in light of Corollary~\ref{coroutile}, for~$i\in\{1,2\}$,
we can find open and bounded sets~$\Omega_{i,j}$ with boundaries
of class~$C^{2,\alpha}$ such that~$\Omega_{i}\subset\Omega_{i,j}$ and
$$\lim_{j\to+\infty}\cAPAC(\Omega_{i,j})=\cAPAC(\Omega_i)$$
and~$\Omega_{1,j}\cup\Omega_{2,j}$ has boundary of class~$C^{2,\alpha}$.
We also let~$u_{i,j}$ be the conductor potentials relative to~$\Omega_{i,j}$
and we exploit~\eqref{CONJCONDENS} and~\eqref{lskworeiteriyter439548795470}
to say that~$u_{i,j}$ converges to~$u_i$ in~$\R^n$ (up to a set of measure zero that we disregard)
as~$j\to+\infty$.

Let now~$u:=u_{1,j}+u_{2,j}-u_3$. This function is harmonic
outside~$\Omega_{1,j}\cup\Omega_{2,j}$ and it belongs to~${\mathcal{D}}^{1,2}(\R^n)$. Furthermore,
on~$\partial(\Omega_{1,j}\cup\Omega_{2,j})\subseteq(\partial\Omega_{1,j})\cup(\partial\Omega_{2,j})$
we have that either~$u_{1,j}=1$ or~$u_{2,j}=1$, and therefore~$u\ge1-u_3\ge0$.

As a result, we can apply
Lemma~\ref{SOBOMAXPLEH10} and deduce that~$u\ge0$ in~$\R^n$.
This yields that~$u_{1,j}+u_{2,j}\ge u_3$, whence, sending~$j\to+\infty$,
we see that~$u_1+u_2\ge u_3$ a.e. in~$\R^n$.
\end{proof}

The conductor potentials for general bounded sets
(as constructed in Corollary~\ref{gujbmSIJKc-qpwruofygheibdowi3egufewbuigf3289eufgvosud81ye2rit23456gh}) is not necessarily
continuous along the boundary of~$\Omega$ (differently from the case
of sets with boundaries of class~$C^{2,\alpha}$ treated in
Proposition~\ref{AM-CAPA}). Hence, in the general setting of Corollary~\ref{gujbmSIJKc-qpwruofygheibdowi3egufewbuigf3289eufgvosud81ye2rit23456gh},
the conductor potential may not need to attain value~$1$ continuously 
along~$\partial\Omega$, but only in the Sobolev trace sense.
The possible continuity
of the conductor potentials (locally in the complement of a given domain)
is closely related to the solvability of the Dirichlet
problem, since, roughly speaking, the conductor potential can play the role of a useful barrier
at boundary points. In this framework, we have the following result:

\begin{proposition}\label{GOA-prokd-P-13o}
Let~$n\ge3$ and~$\Omega\subset\R^n$ be open, bounded and connected.

Then, the following statements are equivalent:
\begin{itemize}
\item[(i).] For every~$p\in\partial\Omega$ there exists a sequence~$\rho_j>0$ with~$\rho_j\to0$ as~$j\to+\infty$
such that if~$u_{p,j}$ is the conductor potential of~$B_{\rho_j}(p)\setminus\Omega$, as given by Corollary~\ref{gujbmSIJKc-qpwruofygheibdowi3egufewbuigf3289eufgvosud81ye2rit23456gh}, we have that, for all~$j\in\N$,
\begin{equation}\label{89contpotenap0} \lim_{ \Omega\ni x\to p} u_{p,j}(x)=1.\end{equation}
\item[(ii).] For every~$p\in\partial\Omega$ and every sequence~$\rho_j>0$ with~$\rho_j\to0$ as~$j\to+\infty$
such that if~$u_{p,j}$ is the conductor potential of~$B_{\rho_j}(p)\setminus\Omega$, as given by Corollary~\ref{gujbmSIJKc-qpwruofygheibdowi3egufewbuigf3289eufgvosud81ye2rit23456gh}, we have that, for all~$j\in\N$,
\begin{equation*} \lim_{ \Omega\ni x\to p} u_{p,j}(x)=1.\end{equation*}
\item[(iii).] For every~$g\in C(\partial\Omega)$
there exists a unique solution~$u\in C^2(\Omega)\cap C(\overline\Omega)$
of the Dirichlet problem
$$ \begin{dcases}
\Delta u=0 &{\mbox{ in }}\Omega,\\
u=g &{\mbox{ on }}\partial\Omega.
\end{dcases}$$\end{itemize}
\end{proposition}

\begin{proof} We prove that~(i) implies~(iii) which implies~(ii) (this completes
the proof, since obviously~(ii) implies~(i)).
Assume~(i), pick a point~$p\in\partial\Omega$ and let
$$\beta:=1-\sum_{k=1}^{+\infty}\frac{u_{p,k}}{2^k}.$$
We use~\eqref{89contpotenap0} to see that
$$\lim_{\Omega\ni x\to p}\beta(x)=
1-\sum_{k=1}^{+\infty}\frac{u_{p,k}(p)}{2^k}=
1-\sum_{k=1}^{+\infty}\frac{1}{2^k}=
0.$$
We also observe that~$\beta$ is harmonic (and therefore superharmonic) in~$\Omega$.

We claim that
\begin{equation}\label{S-0k3rfuweWrsdgbERkdspdwe}
{\mbox{$\beta>0$\, in\,~$\overline{\Omega\cap B_{\rho_1}(p)}\setminus\{p\}$.}}\end{equation}
Indeed, by construction~$0\le\beta\le1$, hence if~\eqref{S-0k3rfuweWrsdgbERkdspdwe} were not true
there would exist~$q\in\overline{\Omega\cap B_{\rho_1}(p)}\setminus\{p\}$ such that~$\beta(q)=0$.
That is, there would exist a sequence~$q_j\in\Omega\cap B_{\rho_1}(p)$ such that~$q_j\to q\ne p$
and~$\beta(q_j)\to0$, that is
\begin{equation}\label{pIJD653io8grgtVEKovwoifgee} \lim_{j\to+\infty}\sum_{k=1}^{+\infty}\frac{u_{p,k}(q_j)}{2^k}=1.\end{equation}
We let~$\rho':=|p-q|>0$ and~$j'\in\N$ be such that~$|q-q_j|\le\frac{\rho'}4$
for all~$j\ge j'$. We also take~$j''\in\N$ such that~$\rho_j\le\frac{\rho'}8$ for all~$j\ge j''$
and we define~$J:=\max\{j',j''\}$.
We observe that if~$j\ge J$ and~$y\in B_{\rho'/4}(q_j)$ then~$|p-y|\ge|p-q|-|q-q_j|-|q_j-y|\ge\rho'-\frac{\rho'}4-\frac{\rho'}4=\frac{\rho'}2>\rho_J$, whence
$$B_{\rho'/4}(q_j)\,\subseteq\,
\R^n\setminus B_{\rho_J}(p) \,\subseteq\,
\R^n\setminus(B_{\rho_J}(p)\setminus\Omega).$$
Therefore, recalling~\eqref{CAHBSN123s4345453124135dcHNSdfz2354tyhdq3z4erifjg-2}, if~$j\ge J$,
\begin{equation}\label{P3653I3ONSBIKDSPIOSMNDHD8UUVE} u_{p,J}(q_j)\le\frac{c_n \cAPAC(B_{\rho_J}(p)\setminus\Omega)}{(\rho'/4)^{n-2}}.\end{equation}
Also,~$B_{\rho_J}(p)\setminus\Omega\subseteq B_{\rho_J}(p)\subseteq B_{\rho'/8}(p)$, and we thereby
deduce from Lemma~\ref{PTEJSN5678dsdvfg9CALsdfLDM-98uyhnSM6789SMMSA}
(recall Lemma~\ref{lemakoe84654656548489}) that~$\cAPAC(B_{\rho_J}(p)\setminus\Omega)\le
\cAPAC(B_{\rho'/8}(p))$. This and Corollary~\ref{CAPAPALLA} yield that
$$ \cAPAC(B_{\rho_J}(p)\setminus\Omega)\le \frac{(\rho'/8)^{n-2}}{c_n}.$$
As a result, using~\eqref{P3653I3ONSBIKDSPIOSMNDHD8UUVE}, if~$j\ge J$,
$$ u_{p,J}(q_j)\le\frac{(\rho'/8)^{n-2}}{(\rho'/4)^{n-2}}=\frac1{2^{n-2}}.$$
Thus, in the light of~\eqref{pIJD653io8grgtVEKovwoifgee},
\begin{eqnarray*}
1&=& \lim_{j\to+\infty}\sum_{k=1}^{+\infty}\frac{u_{p,k}(q_j)}{2^k}
\\&\le& \lim_{j\to+\infty}\sum_{{k\ge1}\atop{k\ne J}}\frac{u_{p,k}(q_j)}{2^k}+\frac1{2^{n-2+J}}\\
\\&\le& \sum_{{k\ge1}\atop{k\ne J}}\frac{1}{2^k}+\frac1{2^{n-2+J}}\\
\\&=& \sum_{{k\ge1}}\frac{1}{2^k}+\frac1{2^{n-2+J}}-\frac1{2^J}\\&=&1+\frac1{2^{n-2+J}}-\frac1{2^J}
\\&<&1.
\end{eqnarray*}
This is a contradiction,
and therefore we have completed the proof of~\eqref{S-0k3rfuweWrsdgbERkdspdwe}.

These pieces of information, together with~\eqref{JKA:BARRLPAER}, show that
if~$u$ is the harmonic function constructed by the Perron method in
Theorem~\ref{PERRO} then~$u$ belongs to~$ C(\overline\Omega)$, hence it is a solution of the Dirichlet
problem in~(iii) (the uniqueness of such a solution
is entailed by Corollary~\ref{UNIQUENESSTHEOREM}).

Let us now assume that~(iii) holds true. We take~$p\in\partial\Omega$ and~$\rho>0$.
Let~$u_{p,\rho}$ be the conductor potential of~$B_{\rho}(p)\setminus\Omega$, as given by Corollary~\ref{gujbmSIJKc-qpwruofygheibdowi3egufewbuigf3289eufgvosud81ye2rit23456gh}.
To establish~(ii), it suffices to show that
\begin{equation}\label{ijliimhbdQjsdxocsaqkP934Mlaptye}
\lim_{\Omega\ni x\to p}u_{p,\rho}(x)=1.
\end{equation}
To this end, let~$R>\rho$ and~${\mathcal{B}}:=B_{\rho}(p)\setminus\Omega$.
We consider boundary data equal to~$0$ on~$\partial B_{R}(p)$
and equal to~$1$ on~$\partial {\mathcal{B}}$
and we construct the corresponding harmonic function~$v$
in~$B_{R}(p)\setminus\overline{\mathcal{B}}$, as given by
the Perron method in Theorem~\ref{PERRO}.
We also define~$g(x):=|x-p|$ and
use~(iii) to find a solution~$\beta\in C^2({\Omega})\cap C(\overline{\Omega})$
of the Dirichlet problem
$$ \begin{dcases}
\Delta \beta=0 &{\mbox{ in }}\Omega,\\
\beta=g &{\mbox{ on }}\partial\Omega.
\end{dcases}$$
We claim that
\begin{equation}\label{09iuygf06tf57tyd89etyfgeoendash12-q0wrf9eihg}
{\mbox{$\beta(x)>0$ for every~$x\in\overline{\Omega\cap B_\rho(p)}\setminus\{p\}$.}}
\end{equation}
Indeed, by the Weak Maximum Principle in Corollary~\ref{WEAKMAXPLE}(iii) we know that, for all~$x\in\Omega$,
$$\beta(x)\ge\inf_{\partial\Omega}g=g(p)=0.$$
This and the Strong Maximum Principle in Theorem~\ref{STRONGMAXPLE1} give that~$\beta>0$
in~$\Omega$. Combining this with the fact that~$\beta(x)=g(x)>0$ for all~$x\in(\partial\Omega)\setminus\{p\}$,
we obtain~\eqref{09iuygf06tf57tyd89etyfgeoendash12-q0wrf9eihg},
as desired.

Moreover,
$$ \lim_{\Omega\ni x\to p}\beta(x)=g(p)=0.$$
This and~\eqref{09iuygf06tf57tyd89etyfgeoendash12-q0wrf9eihg}
allow us to exploit the barrier technique in~\eqref{JKA:BARRLPAER} for
the function~$v$ and conclude that
\begin{equation}\label{P0oik-9ierjfo-Xmmde-203uohgsmdno34}
\lim_{B_R(p)\setminus\overline{\mathcal{B}}\ni x\to p}v(x)=1.
\end{equation}
We also remark that for every~$y\in\partial B_R(p)$
\begin{equation}\label{KL:PER-perlav-qui934-2S}
\lim_{B_R(p)\setminus\overline{\mathcal{B}}\ni x\to y}v(x)=0,
\end{equation}
since along~$\partial B_R(p)$ we can use the exterior cone condition,
exploit the barrier constructed in the proof of Theorem~\ref{PERRO-2} (see in particular the computations
between formula~\eqref{MI4M356U2LSMSUJNS-G5O2F4O2R} and the end of the proof)
and rely on~\eqref{JKA:BARRLPAER}.

Additionally,
\begin{equation}\label{KL:PER-perlav-qui934-2}
\sup_{{{B_{R}(p)\setminus\overline{\mathcal{B}}}}}v(x)\le1.
\end{equation}
Indeed, let~${\mathcal{F}}$ be the class of functions~$w\in C({{B_{R}(p)\setminus\overline{\mathcal{B}}}})$
which are subharmonic in~${B_{R}(p)\setminus\overline{\mathcal{B}}}$ and such that
\begin{eqnarray*}&& 
\limsup_{{{B_{R}(p)\setminus\overline{\mathcal{B}}}\ni x\to y}}w(x)\le 0\qquad{\mbox{for every~$y\in\partial B_R(p)$}}
\\ {\mbox{and }}&&\limsup_{{B_{R}(p)\setminus\overline{\mathcal{B}}}\ni x\to z}w(x)\le 1\qquad{\mbox{for every~$z\in\partial {\mathcal{B}}$.}}
\end{eqnarray*}
We know by the Perron method construction in
Theorem~\ref{PERRO} that
\begin{equation}\label{KL:PER-perlav-qui934} v(x):=\sup_{w\in{\mathcal{F}} }w(x).\end{equation}
Also, if~$w\in{\mathcal{F}}$, 
we take a sequence of points~$x_k\in B_R(p)\setminus \overline{\mathcal{B}}$ such that
$$ \lim_{k\to+\infty} w(x_k)=\sup_{B_R(p)\setminus \overline{\mathcal{B}}} w.$$
Up to a subsequence, we can suppose that~$x_k$ converges to
some~$\overline{x}\in\overline{B_R(p)}\setminus \mathcal{B}$ as~$k\to+\infty$.
Now, in light of Lemma~\ref{LplrfeELLtahsmedschdfb4camntM}, we have that~$\overline{x}\in
\partial\big(B_R(p)\setminus \overline{\mathcal{B}}\big)$, and therefore
$$ \sup_{B_R(p)\setminus \overline{\mathcal{B}}} w=
\limsup_{k\to+\infty} w(x_k)\le 1.$$
This and~\eqref{KL:PER-perlav-qui934} give~\eqref{KL:PER-perlav-qui934-2}.

Now we recall~\eqref{CONJCONDENS}, namely that the 
conductor potential~$u_{p,\rho}$ 
is obtained as the a.e. limit of functions~$u_j$ taking values in~$[0,1]$, which are harmonic
outside a set~${\mathcal{U}}_j$ that approaches~$B_{\rho}(p)\setminus\Omega$
from outside
as~$j\to+\infty$, with~$u_j=1$ in~$ {\mathcal{U}}_j\supseteq
{\mathcal{B}}$. We can also suppose that~${\mathcal{U}}_j\Subset B_R(p)$.

In particular, the function~$\widetilde{u}_j:=u_j-v$
is harmonic in~$B_{R}(p)\setminus\overline{\mathcal{U}_j}$
and nonnegative on~$\partial({B_{R}(p)\setminus\overline{\mathcal{U}}_j})$, owing to~\eqref{KL:PER-perlav-qui934-2S} and~\eqref{KL:PER-perlav-qui934-2}.
Accordingly, by the Maximum Principle in
Corollary~\ref{WEAKMAXPLE}(iii), we have that~$\widetilde{u}_j\ge0$, i.e.~$u_j\ge v$,
in~$ {B_{R}(p)\setminus\overline{\mathcal{B}}}$.

Also, we can set~$v$ to be equal to~$1$ in~${\mathcal{B}}$, and therefore, since~$u_j=1$ in~${\mathcal{B}}$, 
we infer that~$u_j\ge v$ in~$B_R(p)$.
In particular, recalling~\eqref{P0oik-9ierjfo-Xmmde-203uohgsmdno34},
$$ \lim_{{\mathcal{B}}\ni x\to p}u_{p,\rho}(x)=\lim_{{\mathcal{B}}\ni x\to p}\lim_{j\to+\infty}u_j(x)
\ge\lim_{{\mathcal{B}}\ni x\to p}v(x)=1.$$
{F}rom this we obtain~\eqref{ijliimhbdQjsdxocsaqkP934Mlaptye}, as desired.
\end{proof}

Now we reconsider Lemma~\ref{PTEJSN5678dsdvfg9CALsdfLDM-98uyhnSM6789SMMSA}:
such a result suggests
that the capacity is a way of measuring how large a set is: however, in light of
Proposition~\ref{AM-CAPA} (and of the comments following its proof),
this notion of largeness is not quite about the usual size of an object
but rather on its capability at holding charge. Nonetheless, 
the notion of capacity can be compared to Lebesgue
and Hausdorff measures (see e.g.~\cite{MR3409135} for an introduction
on the concept of Hausdorff measure, as well as for further information
about capacities and measures; see also~\cite[Section~2.1.7]{MR1461542} for additional details and a general
setting):

\begin{theorem}
Let~$n\ge3$ and $\Omega$ be a bounded subset of~$\R^n$. 
Then,
\begin{equation}\label{CAPSOESOT} c|\Omega|^{\frac{n-2}n}\le\cAPAC(\Omega)\le C{\mathcal{H}}^{n-2}(\Omega),\end{equation}
for suitable positive constants~$c$ and~$C$, depending only on~$n$.
\end{theorem}

\begin{proof} 
We start by proving the first inequality in~\eqref{CAPSOESOT}.
To this end, we take~$u_j$ as the conductor potential of a sequence of bounded open
sets~$\Omega_j$ approximating~$\Omega$ from outside, in the setting of~\eqref{0oikj8i24ks9ihf94i23203762423} and~\eqref{CONJCONDENS},
and we use the Sobolev-Gagliardo-Nirenberg's Inequality (see
e.g.~\cite[Theorem~11.2]{MR2527916}) to find that
\begin{eqnarray*}&& |\Omega|^{\frac{n-2}{2n}}\le|\Omega_j|^{\frac{n-2}{2n}}=
\left(\int_{\Omega_j} |u_j(x)|^{\frac{2n}{n-2}}\,dx\right)^{\frac{n-2}{2n}}
\le
\left(\int_{\R^n} |u_j(x)|^{\frac{2n}{n-2}}\,dx\right)^{\frac{n-2}{2n}}\\&&\qquad\le C
\left(\int_{\R^n} |\nabla u_j(x)|^2\,dx\right)^{\frac{1}{2}}= C
\left(\cAPAC(\Omega_j)\right)^{\frac{1}{2}}\le C
\left(\cAPAC(\Omega)+\frac1j\right)^{\frac{1}{2}},\end{eqnarray*}
for some~$C>0$.
This establishes the first inequality in~\eqref{CAPSOESOT} by sending~$j\to+\infty$.

To prove the second inequality in~\eqref{CAPSOESOT}, let~$\delta>0$ and suppose that $$\Omega\subseteq \bigcup_{i\in\N} B_{r_i}(x_i),$$ for suitable~$x_i\in\R^n$ and~$r_i\in(0,\delta]$. Then, given~$\e>0$, $$\overline\Omega\subseteq \bigcup_{j\in\N} B_{r_i+\e}(x_i),$$
and thus, by compactness, there exists a finite set of indices~$i_1,\dots,i_N\in\N$ such that $$\Omega\subseteq\overline\Omega\subseteq \bigcup_{j=0}^N B_{r_{i_j}+\e}(x_{i_j}).$$ Therefore, by 
Lemma~\ref{PTEJSN5678dsdvfg9CALsdfLDM-98uyhnSM6789SMMSA}
(as refined in Lemma~\ref{lemakoe84654656548489}),
$$ \cAPAC(\Omega)\le\sum_{j=0}^N \cAPAC(B_{r_{i_j}+\e}(x_{i_j})).$$
Hence, sending~$\e\searrow0$,
$$ \cAPAC(\Omega)\le\sum_{j=0}^N \cAPAC(B_{r_{i_j}}(x_{i_j}))\le
\sum_{i\in\N} \cAPAC(B_{r_i}(x_i)).$$ This and Corollary~\ref{CAPAPALLA}
give that $$ c_n\cAPAC(\Omega)\le
\sum_{i\in\N} r_i^{n-2}.$$ Taking the infimum over all these possible coverings of~$\Omega$ and then sending~$\delta\searrow0$ we obtain the second inequality in~\eqref{CAPSOESOT}, as desired.\index{capacity|)}
\end{proof}

\section{Wiener's criterion}

Coming back to the existence result for the Dirichlet problem in Corollary~\ref{S-coroEXIS-M023}, we are now in the
position of characterizing, via the notion of capacity, the domains which always allow the continuous solvability
of the Dirichlet problem up to the boundary. This classical result is due to\footnote{Besides his eminent contributions
to partial differential equations and stochastic processes, Wiener is often considered the originator of cybernetics
and a forerunner of artificial intelligence.
This is possibly the reason for which the synthesizers player Luigi Tonet
composed an experimental, and slightly dystopic, piece of electronic music with title ``Dedicated to Norbert Wiener''
(with lyrics ``We shall continue our life and thank you for your gift of eternity, thank you Norbert Wiener, thank you for everything'').

For a short time, Wiener was also a journalist for the Boston Herald,
an American daily newspaper (allegedly he was fired quite soon
for his reluctance to write favorable articles about an influential politician).}
Norbert Wiener~\cite{WIENE} and it is often referred to as ``Wiener's criterion''.
The bottom line of this criterion is that the Dirichlet problem being solvable
is equivalent for the domain to allow ``enough space'' outside each of its
boundary points (compare e.g. with the conditions of external cone or external ball
discussed in Theorem~\ref{PERRO-2}
and keep in mind the cuspidal counterexample presented on pages~\pageref{CUSPI-01}--\pageref{CUSPI-02}) and that to precisely quantify this
notion one needs to look at the capacities of small external neighborhoods\index{Wiener's criterion}
of every boundary point. The formal result can be stated as follows:

\begin{theorem}\label{NOWI}
Let~$n\ge3$ and~$\Omega$ be an open and bounded subset of~$\R^n$.
The following conditions are equivalent:
\begin{itemize}
\item[(i).] For every~$g\in C(\partial\Omega)$ there exists a function~$u\in C^2(\Omega)\cap C(\overline\Omega)$ that solves the Dirichlet problem
$$ \begin{dcases}
\Delta u=0 & {\mbox{ in }}\Omega,\\
u=g & {\mbox{ on }}\partial\Omega.
\end{dcases}$$
\item[(ii).] For every~$p\in\partial\Omega$ and every~$\lambda\in(0,1)$,
\begin{equation}\label{SOMMAIWIE} \sum_{j=0}^{+\infty}\frac{\cAPAC\Big(\big(
B_{\lambda^j}(p)\setminus {B_{\lambda^{j+1}}(p)}\big)\setminus\Omega\Big)}{\lambda^{(n-2)j}}=+\infty.\end{equation}
\end{itemize}\end{theorem}

\begin{proof} The idea of the proof is to consider a point~$p\in\partial\Omega$
and relate the convergence or divergence of the series in~\eqref{SOMMAIWIE}
with the barrier condition in~\eqref{JKA:BARRLPAER}. 

Up to a translation, we suppose that~$p$
is the origin.
We notice that given~$\eta\in(0,1)$, if~$\mu:=\eta^2$,
\begin{equation}\label{CVBNJLSOMDOABSBBIJMSDKIAMSDFA}
\begin{split}&
\sum_{j=0}^{+\infty}\frac{\cAPAC\Big(\big(
B_{\eta^j}\setminus{ B_{\eta^{j+1}}}\big)\setminus\Omega\Big)}{\eta^{(n-2)j}}=+\infty\\&\qquad
\qquad{\mbox{if and only if}}\qquad
\sum_{j=0}^{+\infty}\frac{\cAPAC\Big(\big(
B_{\mu^j}\setminus{ B_{\mu^{j+1}}}\big)\setminus\Omega\Big)}{\mu^{(n-2)j}}=+\infty.\end{split}
\end{equation}
To check this, we observe that
$$ B_{\mu^j}\setminus{ B_{\mu^{j+1}}}=
B_{\eta^{2j}}\setminus {B_{\eta^{2j+2}} }=
\big(B_{\eta^{2j}}\setminus {B_{\eta^{2j+1}} }\big)\cup
\big({B_{\eta^{2j+1}}}\setminus {B_{\eta^{2j+2}} }\big).
$$
This entails that, for each~$i\in\{0,1\}$,
$$ B_{\eta^{2j+i}}\setminus {B_{\eta^{2j+i+1}} }\subseteq
B_{\mu^j}\setminus{ B_{\mu^{j+1}}}\subseteq
\big(B_{\eta^{2j}}\setminus {B_{\eta^{2j+1}} }\big)\cup
\big({B_{\eta^{2j+1}}}\setminus {B_{\eta^{2j+2}} }\big)
,$$
whence, by Lemma~\ref{PTEJSN5678dsdvfg9CALsdfLDM-98uyhnSM6789SMMSA}
(in its general formulation, as warranted by Lemma~\ref{lemakoe84654656548489}),
\begin{equation*}
\begin{split}&\cAPAC\Big(\big(B_{\eta^{2j+i}}\setminus{ B_{\eta^{2j+i+1}}}\big)\setminus\Omega\Big)\le
\cAPAC\Big(\big(B_{\mu^j}\setminus{ B_{\mu^{j+1}}}\big)\setminus\Omega\Big)\\
&\qquad\le\cAPAC\Big(\big(B_{\eta^{2j}}\setminus{ B_{\eta^{2j+1}}}\big)\setminus\Omega\Big)
+\cAPAC\Big(\big(B_{\eta^{2j+1}}\setminus{ B_{\eta^{2j+2}}}\big)\setminus\Omega\Big).
\end{split}\end{equation*}
Dividing by~$\mu^{(n-2)j}=\eta^{2(n-2)j}$ and summing up, we obtain
\begin{equation}\label{STEOAMS876543PIBV08AM}
\begin{split}&\eta^{(n-2)i}\sum_{j=0}^{+\infty}
\frac{\cAPAC\Big(\big(B_{\eta^{2j+i}}\setminus{ B_{\eta^{2j+i+1}}}\big)\setminus\Omega\Big)}{\eta^{(n-2)(2j+i)}}\le\sum_{j=0}^{+\infty}\frac{
\cAPAC\Big(\big(B_{\mu^j}\setminus{ B_{\mu^{j+1}}}\big)\setminus\Omega\Big)}{\mu^{(n-2)j}}\\
&\qquad\le\sum_{j=0}^{+\infty}
\frac{\cAPAC\Big(\big(B_{\eta^{2j}}\setminus{ B_{\eta^{2j+1}}}\big)\setminus\Omega\Big)}{\eta^{2(n-2)j}}
+\eta^{n-2}\sum_{j=0}^{+\infty}\frac{\cAPAC\Big(\big(B_{\eta^{2j+1}}\setminus{ B_{\eta^{2j+2}}}\big)\setminus\Omega\Big)}{\eta^{(n-2)(2j+1)}}.
\end{split}\end{equation}
Now, if the first series in~\eqref{CVBNJLSOMDOABSBBIJMSDKIAMSDFA} converges,
we use~\eqref{STEOAMS876543PIBV08AM} and the fact that~$\mu\le1$ to find that
\begin{eqnarray*}&&
\sum_{j=0}^{+\infty}\frac{
\cAPAC\Big(\big(B_{\mu^j}\setminus{ B_{\mu^{j+1}}}\big)\setminus\Omega\Big)}{\mu^{(n-2)j}}
\le\sum_{j=0}^{+\infty}
\frac{\cAPAC\Big(\big(B_{\eta^{2j}}\setminus{ B_{\eta^{2j+1}}}\big)\setminus\Omega\Big)}{\eta^{2(n-2)j}}
+\sum_{j=0}^{+\infty}\frac{\cAPAC\Big(\big(B_{\eta^{2j+1}}\setminus{ B_{\eta^{2j+2}}}\big)\setminus\Omega\Big)}{\eta^{(n-2)(2j+1)}}\\&&\qquad\qquad=
\sum_{j=0}^{+\infty}\frac{
\cAPAC\Big(\big(B_{\eta^j}\setminus{ B_{\eta^{j+1}}}\big)\setminus\Omega\Big)}{\eta^{(n-2)j}}<+\infty,
\end{eqnarray*}
showing that the second series in~\eqref{CVBNJLSOMDOABSBBIJMSDKIAMSDFA} converges as well.

If instead the first series in~\eqref{CVBNJLSOMDOABSBBIJMSDKIAMSDFA} diverges,
we exploit~\eqref{STEOAMS876543PIBV08AM} by dividing the first inequality there by~$\eta^{(n-2)i}$
and summing over~$i\in\{0,1\}$: in this way, we obtain that
\begin{eqnarray*}&&+\infty=
\sum_{j=0}^{+\infty}\frac{
\cAPAC\Big(\big(B_{\eta^j}\setminus{ B_{\eta^{j+1}}}\big)\setminus\Omega\Big)}{\eta^{(n-2)j}}=
\sum_{i=0}^1\left[\sum_{j=0}^{+\infty}
\frac{\cAPAC\Big(\big(B_{\eta^{2j+i}}\setminus{ B_{\eta^{2j+i+1}}}\big)\setminus\Omega\Big)}{\eta^{(n-2)(2j+i)}}\right]\\&&\qquad\qquad\le
\left(1+\frac1{\eta^{n-2}}\right)
\sum_{j=0}^{+\infty}\frac{
\cAPAC\Big(\big(B_{\mu^j}\setminus{ B_{\mu^{j+1}}}\big)\setminus\Omega\Big)}{\mu^{(n-2)j}}
\end{eqnarray*}
and this proves that the second series in~\eqref{CVBNJLSOMDOABSBBIJMSDKIAMSDFA} diverges as well, thus completing the proof of~\eqref{CVBNJLSOMDOABSBBIJMSDKIAMSDFA}.

We also point out that, given~$\lambda\in(0,1)$ and~$\e\in\left(0,\min\{\lambda,1-\lambda\}\right)$, there exist~$\underline\lambda\in(0,\e)$ and~$\overline\lambda\in(1-\e,1)$ such that
\begin{equation}\label{CVBNJLSOMDOABSBBIJMSDKIAMSDFA-EPS}
\begin{split}&
\sum_{j=0}^{+\infty}\frac{\cAPAC\Big(\big(
B_{\lambda^j}\setminus{ B_{\lambda^{j+1}}}\big)\setminus\Omega\Big)}{\lambda^{(n-2)j}}=+\infty\\&\qquad
\qquad{\mbox{if and only if}}\qquad
\sum_{j=0}^{+\infty}\frac{\cAPAC\Big(\big(
B_{\underline\lambda^j}\setminus{ B_{\underline\lambda^{j+1}}}\big)\setminus\Omega\Big)}{\underline\lambda^{(n-2)j}}=+\infty\\&\qquad
\qquad{\mbox{and}}\qquad
\sum_{j=0}^{+\infty}\frac{\cAPAC\Big(\big(
B_{\overline\lambda^j}\setminus{ B_{\overline\lambda^{j+1}}}\big)\setminus\Omega\Big)}{\overline\lambda^{(n-2)j}}=+\infty,
\end{split}
\end{equation}
that is the convergence or divergence of the series in~\eqref{SOMMAIWIE}
for a given~$\lambda\in(0,1)$ can be equivalently shifted to that
corresponding to another parameter as close as we wish to~$0$, as well as to~$1$.
Notice that the claim in~\eqref{CVBNJLSOMDOABSBBIJMSDKIAMSDFA-EPS}
is a direct consequence of that in~\eqref{CVBNJLSOMDOABSBBIJMSDKIAMSDFA}, up to an iteration, by taking~$\lambda:=\eta$ and~$\ell\in\N$ sufficiently large so that~$\underline\lambda:=\eta^{2\ell}\in(0,\e)$ (or by taking~$\lambda:=\mu$
and~$\ell\in\N$ sufficiently large so that~$\overline\lambda:=\eta^{1/(2\ell)}\in(1-\e,1)$).

We also point out that the capacity behaves nicely under scaling.
Namely, if~$\kappa>0$ and~$\kappa E:=\{\lambda x$, $x\in E\}$,
if a function~$v$ is equal to~$1$ in~$E$, then the function~$v_\kappa(x):=v\left(\frac{x}\kappa\right)$ is equal to~$1$ in~$\kappa E$, hence it follows from~\eqref{DEFCAP} that
\begin{equation}\label{SCALA-COSCAPAXCP}
\cAPAC(\kappa E):=\kappa^{n-2}\cAPAC(E).\end{equation}
Now, let us
suppose that condition~(ii) in the statement of Theorem~\ref{NOWI} holds true,
namely that the series in~\eqref{SOMMAIWIE} diverges for some~$\lambda\in(0,1)$.
Let also~$\e>0$, to be taken arbitrarily small in what follows. We exploit~\eqref{CVBNJLSOMDOABSBBIJMSDKIAMSDFA-EPS}
that allows us to replace the previous $\lambda$ with a new one, still called~$\lambda$
for simplicity, with the additional property that~$\lambda\in(1-\e,1)$.
We also take~$K\in\N$ large enough (possibly in dependence of~$\lambda$ and~$\e$)
such that~$\lambda^K\in\left(0,\frac\e4\right)$.
For~$\ell\in\{1,\dots,K\}$, we consider the series
$$ {\mathcal{S}}_\ell(\Omega):=
\sum_{j=0}^{+\infty}\frac{\cAPAC\Big(\big(
B_{\lambda^{Kj+\ell-1}}\setminus {B_{\lambda^{Kj+\ell}}}\big)\setminus\Omega\Big)}{\lambda^{(n-2)(Kj+\ell-1)}}.
$$
By~\eqref{SCALA-COSCAPAXCP}, we know that
\begin{eqnarray*} &&\cAPAC\Big(\big(
B_{\lambda^{Kj+\ell-1}}\setminus {B_{\lambda^{Kj+\ell}}}\big)\setminus\Omega\Big)=
\cAPAC\left(\lambda^{\ell-1}\left(\big(
B_{\lambda^{Kj}}\setminus {B_{\lambda^{Kj+1}}}\big)\setminus\frac{\Omega}{\lambda^{\ell-1}}\right)\right)\\&&\qquad=\lambda^{(n-2)(\ell-1)}
\cAPAC\left( \big(
B_{\lambda^{Kj}}\setminus {B_{\lambda^{Kj+1}}}\big)\setminus\frac{\Omega}{\lambda^{\ell-1}}\right)
\end{eqnarray*}
and therefore
$$ {\mathcal{S}}_\ell(\Omega):=
\sum_{j=0}^{+\infty}\frac{\lambda^{(n-2)(\ell-1)}\cAPAC\left( \big(
B_{\lambda^{Kj}}\setminus {B_{\lambda^{Kj+1}}}\big)\setminus\frac{\Omega}{\lambda^{\ell-1}}\right)}{\lambda^{(n-2)(Kj+\ell-1)}}
=
{\mathcal{S}}_1\left(\frac{\Omega}{\lambda^{\ell-1}}\right).
$$
{F}rom this, we arrive at
\begin{eqnarray*}&&+\infty=
\sum_{j=0}^{+\infty}\frac{\cAPAC\Big(\big(
B_{\lambda^j}\setminus {B_{\lambda^{j+1}}}\big)\setminus\Omega\Big)}{\lambda^{(n-2)j}}=\sum_{\ell=1}^{K}
{\mathcal{S}}_\ell(\Omega)=\sum_{\ell=1}^{K}{\mathcal{S}}_1\left(\frac{\Omega}{\lambda^{\ell-1}}\right)\end{eqnarray*}
and, as a consequence, we can find~$\ell_\star\in\{1,\dots,K\}$ such that if
\begin{equation}\label{MS-olSSTANNDRRND-po}
\Omega_\star:=\frac{\Omega}{\lambda^{\ell_\star-1}}\end{equation} it holds that
\begin{equation}\label{098765-0987654-ojnbX0013xsa}
{\mathcal{S}}_1(\Omega_\star)=+\infty.
\end{equation}
We now remark that
\begin{equation}\label{FUNcal201-1P-c1}
\lambda-\lambda^K\le \lambda^{1-K}-1.
\end{equation}
To check this, for all~$t\in\R$, we consider the function~$\phi(t):=\lambda^{1+t}-\lambda^{K+t}$
and we note that
$$ \phi'(t)=\ln\lambda\,(\lambda^{1+t}-\lambda^{K+t})=
-\lambda^{1+t}|\ln\lambda|\,(1-\lambda^{K-1})\le0.
$$
As a result~$
\lambda^{1-K}-1=
\phi(-K)\ge\phi(0)=\lambda-\lambda^{K}$
and this gives~\eqref{FUNcal201-1P-c1}.

Now, we consider the
conductor potential~$v_j$ of the set~$E_j:=\big(
B_{\lambda^j}\setminus {B_{\lambda^{j+1}}}\big)\setminus\Omega_\star$,
as constructed in Corollary~\ref{gujbmSIJKc-qpwruofygheibdowi3egufewbuigf3289eufgvosud81ye2rit23456gh}.
Given~$m'\ge m\in\N$, and taking~$K$ as above, we define
$$ V_{m,m'}:=\sum_{i=m}^{m'} v_{Ki}.$$
Notice that
\begin{equation}\label{1-VodmvVikfAndfdfNAJowqrq0aqjf}
{\mbox{$V_{m,m'}$ is harmonic
outside~$E_{m,m'}:=$}}\bigcup_{i=m}^{m'} E_{Ki}.\end{equation}
Thus, the strategy is to estimate~$V_{m,m'}$ in a small neighborhood of~$E_{m,m'}$ and
then use the Maximum Principle in Lemma~\ref{VodmvVikfAndfdfNAJowqrq0aqjf}
to obtain a bound for~$V_{m,m'}$ outside~$E_{m,m'}$.

To employ this strategy, given~$\delta\in(0,1)$, we consider a~$\delta$-neighborhood of~$E_{m,m'}$ and
we observe that if~$x$ lies in this~$\delta$-neighborhood then there exists~$i_0\in\{m,\dots,m'\}$
such that~$x$ lies in a~$\delta$-neighborhood of~$E_{Ki_0}$. We claim that, if~$\delta$ is
sufficiently small, for every~$i\neq i_0$,
\begin{equation}\label{mcdhfje2134565768970987654}
B_{(\lambda-\lambda^K)\lambda^{Ki}-\delta}(x)\cap E_{Ki}=\varnothing.
\end{equation}
Indeed, suppose by contradiction that there exists~$y\in B_{(\lambda-\lambda^K)\lambda^{Ki}-\delta
}(x)\cap E_{Ki}$. Then, 
$$ |x-y|<(\lambda-\lambda^K)\lambda^{Ki}-\delta \qquad {\mbox{and}}\qquad
\lambda^{Ki+1}\le |y|<\lambda^{Ki}.$$
Using these inequalities, we see that
\begin{eqnarray*}
\lambda^{Ki_0}+\delta>|x|\ge |y|-|x-y|\ge \lambda^{Ki+1}-(\lambda-\lambda^K)\lambda^{Ki}+\delta
=\lambda^{K(i+1)}+\delta,
\end{eqnarray*}
which gives that~$i_0<i+1$. On the other hand, exploiting also~\eqref{FUNcal201-1P-c1},
\begin{eqnarray*}
\lambda^{Ki_0+1}-\delta<|x|\le |x-y|+|y|<(\lambda-\lambda^K)\lambda^{Ki}-\delta+ \lambda^{Ki}=
(\lambda-\lambda^K+1)\lambda^{Ki}-\delta\le \lambda^{1-K}\lambda^{Ki}-\delta,
\end{eqnarray*}
which gives instead that~$i_0>i-1$. Hence, we have that~$i=i_0$ and this is a contradiction, which
establishes~\eqref{mcdhfje2134565768970987654}.

In light of~\eqref{mcdhfje2134565768970987654},
we can employ~\eqref{CAHBSN123s4345453124135dcHNSdfz2354tyhdq3z4erifjg-2}
(with~$\Omega:=E_{Ki}$, $u:=v_{Ki}$ and~$r:=(\lambda-\lambda^K)\lambda^{Ki}-\delta$)
to find that
if~$x$ lies in a~$\delta$-neighborhood of~$E_{Ki_0}$ then, as~$\delta\searrow0$,
$$ v_{Ki}(x)\le\frac{c_n (1+o(1))\cAPAC(E_{Ki})}{(\lambda-\lambda^K)^{n-2}\;\lambda^{K(n-2)i}}.
$$
As a result, since~$v_{Ki_0}\le1$, we obtain that, if~$x$ lies in a~$\delta$-neighborhood of~$E_{Ki_0}$ then
\begin{eqnarray*}
V_{m,m'}(x)&\le& 1+\sum_{{m\le i\le m'}\atop{i\ne i_0}} v_{Ki}(x)\\
&\le&1+\sum_{{m\le i\le m'}\atop{i\ne i_0}}
\frac{c_n (1+o(1))\cAPAC(E_{Ki})}{(\lambda-\lambda^K)^{n-2}\;\lambda^{K(n-2)i}}.
\end{eqnarray*}
That is, if~$z\in\partial E_{m,m'}$,
\begin{eqnarray*}
\lim_{\R^n\setminus\overline{E_{m,m'}}\ni x\to z}V_{m,m'}(x)&\le&
1+\sum_{{m\le i\le m'}\atop{i\ne i_0}}
\frac{c_n \cAPAC(E_{Ki})}{(\lambda-\lambda^K)^{n-2}\;\lambda^{K(n-2)i}}\\
&\le&\frac1{(\lambda-\lambda^K)^{n-2}}
\left(1+\sum_{{i=m}}^{m'}
\frac{c_n \cAPAC(E_{Ki})}{\lambda^{K(n-2)i}}\right).
\end{eqnarray*}
As a consequence, by~\eqref{1-VodmvVikfAndfdfNAJowqrq0aqjf} and
Lemma~\ref{VodmvVikfAndfdfNAJowqrq0aqjf} we infer that, for all~$x\in\R^n\setminus E_{m,m'}$,
$$ V_{m,m'}(x)\le\frac1{(\lambda-\lambda^K)^{n-2}}
\left(1+\sum_{{i=m}}^{m'}
\frac{c_n \cAPAC(E_{Ki})}{\lambda^{K(n-2)i}}\right).$$
Therefore
\begin{equation}\label{110-0-33r-0139-134-35d-w4rfgkn2jdpKsTYHS} W_{m,m'}(x):=
\left(1+\sum_{{i=m}}^{m'}
\frac{c_n \cAPAC(E_{Ki})}{\lambda^{K(n-2)i}}\right)^{-1}{(\lambda-\lambda^K)^{n-2}}\,
V_{m,m'}(x)\le1.\end{equation}

Now, our objective is to exploit
Proposition~\ref{GOA-prokd-P-13o}. To this end, given~$\rho>0$ we consider
the conductor potential~$u_{p,\rho}^\star$ of~$B_{\rho/(\lambda^{\ell_\star-1})} \setminus\Omega_\star$, as given by Corollary~\ref{gujbmSIJKc-qpwruofygheibdowi3egufewbuigf3289eufgvosud81ye2rit23456gh}. We take~$m'$ so large that
\begin{equation}\label{0-0-33r-0139-134-35d-w4rfgkn2jdpKsTYHS}
\lambda^{Km'}<\frac\rho2\end{equation} and we claim that, for every~$x\in\R^n$,
\begin{equation}\label{C0NSIKDOSNESERBLMDVEunSaDXCAJDi0}
W_{m,m'}(x)\le u_{p,\rho}^\star(x).
\end{equation}
To check this, we use an approximating sequence~$u_j$ of conductor potentials
as in~\eqref{CONJCONDENS} which are harmonic outside sets~$Z_j$ with
boundaries of class~$C^{2,\alpha}$ which approach~$B_{\rho/(\lambda^{\ell_\star-1})} \setminus\Omega_\star$
from the outside.
We notice that the function~$\phi:=u_{j}-W_{m,m'}$ is harmonic outside~$Z_j$, thanks to~\eqref{0-0-33r-0139-134-35d-w4rfgkn2jdpKsTYHS}. Moreover, $\phi\ge0$ along~$\partial Z_j$
(as well as inside~$Z_j$), owing to~\eqref{110-0-33r-0139-134-35d-w4rfgkn2jdpKsTYHS}.
Since~$\phi\in {\mathcal{D}}^{1,2}(\R^n)$, we thus obtain from
Lemma~\ref{SOBOMAXPLEH10}
that~$\phi\ge0$ in~$\R^n\setminus Z_j$. {F}rom this, we find that~$u_{j}\ge W_{m,m'}$ and thus,
sending~$j\to+\infty$, we obtain~\eqref{C0NSIKDOSNESERBLMDVEunSaDXCAJDi0}, as desired.

As a consequence of~\eqref{C0NSIKDOSNESERBLMDVEunSaDXCAJDi0}, we have that
\begin{equation}\label{CHVSdomonbmbp-de-oishfvo-ca9ikaq}
\begin{split}&
\lim_{ \Omega_\star\ni x\to 0} u_{p,\rho}^\star(x)\ge\lim_{ \Omega_\star\ni x\to 0}W_{m,m'}(x)
=\left(1+\sum_{{i=m}}^{m'}
\frac{c_n \cAPAC(E_{Ki})}{\lambda^{K(n-2)i}}\right)^{-1}{(\lambda-\lambda^K)^{n-2}}\,
\lim_{ \Omega_\star\ni x\to 0}V_{m,m'}(x)\\&\qquad=
\left(1+\sum_{{i=m}}^{m'}
\frac{c_n \cAPAC(E_{Ki})}{\lambda^{K(n-2)i}}\right)^{-1}{(\lambda-\lambda^K)^{n-2}}\,
\sum_{i=m}^{m'}\lim_{ \Omega_\star\ni x\to 0} v_{Ki}(x).
\end{split}
\end{equation}
Now we observe that
$$ B_{\lambda^{Ki}+|x|}(x)\supseteq B_{\lambda^{Ki}}\supseteq E_{Ki},$$
and therefore, by~\eqref{CAHBSN123s4345453124135dcHNSdfz2354tyhdq3z4erifjg-2}
(used here with~$\Omega:=E_{Ki}$, $u:=v_{Ki}$ and~$R:=\lambda^{Ki}+|x|$),
$$  v_{Ki}(x)\ge\frac{c_n \cAPAC(E_{Ki})}{(\lambda^{Ki}+|x|)^{n-2}}.$$
{F}rom this and~\eqref{CHVSdomonbmbp-de-oishfvo-ca9ikaq} we arrive at
\begin{equation}\label{NSDciKDDNFsoiaAMS-OKSQUA0dkf}
\lim_{ \Omega_\star\ni x\to 0} u_{p,\rho}^\star(x)\ge
\left(1+\sum_{{i=m}}^{m'}
\frac{c_n \cAPAC(E_{Ki})}{\lambda^{K(n-2)i}}\right)^{-1}{(\lambda-\lambda^K)^{n-2}}\,
\sum_{i=m}^{m'}\frac{c_n \cAPAC(E_{Ki})}{\lambda^{K(n-2)i}}.
\end{equation}
Now we recall~\eqref{098765-0987654-ojnbX0013xsa}, according to which
$$ +\infty={\mathcal{S}}_1(\Omega_\star)=
\sum_{i=0}^{+\infty}\frac{ \cAPAC\left( \big(
B_{\lambda^{Ki}}\setminus {B_{\lambda^{Ki+1}}}\big)\setminus \Omega_\star\right)}{\lambda^{K(n-2)i}}=
\sum_{i=0}^{+\infty}\frac{ \cAPAC(E_{Ki})}{\lambda^{K(n-2)i}}.
$$
This gives that
$$\sum_{i=m}^{+\infty}\frac{ \cAPAC(E_{Ki})}{\lambda^{K(n-2)i}}=+\infty$$
and therefore, by sending~$m'\to+\infty$ in~\eqref{NSDciKDDNFsoiaAMS-OKSQUA0dkf}, we deduce that
$$ \lim_{ \Omega_\star\ni x\to 0} u_{p,\rho}^\star(x)\ge(\lambda-\lambda^K)^{n-2}\ge\left(1-\e-\frac{\e}4\right)^{n-2}.$$
Hence, sending~$\e\searrow0$,
$$ \lim_{ \Omega_\star\ni x\to 0} u_{p,\rho}^\star(x)\ge1.$$
Since conductor potentials are bounded by~$1$, we thus conclude that
$$ \lim_{ \Omega_\star\ni x\to 0} u_{p,\rho}^\star(x)=1.$$
By scaling back and recalling~\eqref{MS-olSSTANNDRRND-po}, if~$u_{p,\rho}$
is the conductor potential of~$B_{\rho} \setminus\Omega$,
as given in Corollary~\ref{gujbmSIJKc-qpwruofygheibdowi3egufewbuigf3289eufgvosud81ye2rit23456gh},
we find that
$$ \lim_{ \Omega\ni x\to 0} u_{p,\rho}(x)=1.$$
We can accordingly employ Proposition~\ref{GOA-prokd-P-13o}
and obtain the continuous solvability of the Dirichlet problem,
thus establishing condition~(i) in Theorem~\ref{NOWI}.

To complete the proof of Theorem~\ref{NOWI} we now establish the converse statement,
namely that if the series in~\eqref{SOMMAIWIE} converges the statement in~(i) does not hold.
To this end, we argue for a contradiction. Suppose that the statement in~(i) of Theorem~\ref{NOWI}
holds true, take~$p\in\partial\Omega$ (say, $p=0$ up to a translation)
and recall Proposition~\ref{GOA-prokd-P-13o} to consider the infinitesimal sequence of radii~$\rho_j:=\lambda^j$
and have that if~$u_{j}$ is the conductor potential of~$B_{\rho_j}\setminus\Omega$, as given by Corollary~\ref{gujbmSIJKc-qpwruofygheibdowi3egufewbuigf3289eufgvosud81ye2rit23456gh}, we have that, for each~$j\in\N$,
\begin{equation}\label{89contpo4568790jh-tenap0} \lim_{ \Omega\ni x\to 0} u_{j}(x)=1.\end{equation}
In particular, we can use the convergence of the series in~\eqref{SOMMAIWIE} to pick~$m\in\N$ large enough such that
\begin{equation}\label{0ok0o90uytfcikasmcolalijfnggkkaal99304} \sum_{j=m}^{+\infty}\frac{\cAPAC\Big(\big(
B_{\lambda^j}\setminus {B_{\lambda^{j+1}}}\big)\setminus\Omega\Big)}{\lambda^{(n-2)j}}\le\frac{\lambda^{n-2}}{10\,c_n}\end{equation}
(where~$c_n$ is the positive constant in~\eqref{COIENEE})
and then utilize~\eqref{89contpo4568790jh-tenap0}
to find~$r_0>0$ sufficiently small such that if~$x\in\Omega\cap\overline{ B_{r_0}}$ then
\begin{equation}\label{89contpo4568790jh-tenap0-097} u_{m}(x)\ge\frac{9}{10}.\end{equation}
Now we consider~$m'\in\N\cap[m+1,+\infty)$ sufficiently large such that~$\frac{\rho_{m'}^{n-2}}{(r_0-2\rho_{m'})^{n-2}}
\le\frac1{10}$.
In this way, if~$x\in\partial B_{r_0}$ then~$B_{r_0-2\rho_{m'}}(x)\subseteq\R^n\setminus B_{\rho_{m'}}$ 
and therefore, by~\eqref{CAHBSN123s4345453124135dcHNSdfz2354tyhdq3z4erifjg-2}
(applied here with~$B_{\rho_{m'}}\setminus\Omega$ in place of~$\Omega$, $u:=u_{m'}$
and~$r:=r_0-2\rho_{m'}$),
$$ u_{m'}(x)\le\frac{c_n \cAPAC(B_{\rho_{m'}}\setminus\Omega)}{(r_0-2\rho_{m'})^{n-2}}.$$
Thus, recalling Lemma~\ref{PTEJSN5678dsdvfg9CALsdfLDM-98uyhnSM6789SMMSA} 
(in the light of Lemma~\ref{lemakoe84654656548489})
and Corollary~\ref{CAPAPALLA}, if~$x\in\partial B_{r_0}$,
\begin{equation}\label{09876-2345678-83uejdn3874} u_{m'}(x)\le\frac{c_n \cAPAC(B_{\rho_{m'}})}{(r_0-2\rho_{m'})^{n-2}}=
\frac{\rho_{m'}^{n-2}}{(r_0-2\rho_{m'})^{n-2}}\le\frac1{10}.
\end{equation}

Now we denote by~$U_{m,m'} $ the conductor potential, as given by Corollary~\ref{gujbmSIJKc-qpwruofygheibdowi3egufewbuigf3289eufgvosud81ye2rit23456gh}, with respect to the set~$(B_{\rho_m}\setminus B_{\rho_{m'}})\setminus\Omega$.
We make use of Lemma~\ref{KJSM-098765-TFI-IMP245S-PLSIUDND-94} (see also Lemma~\ref{KJSM-098765-TFI-IMP245S-PLSIUDND-95}) and we find that
$$ u_m\le u_{m'}+U_{m,m'}.$$
Hence, using~\eqref{89contpo4568790jh-tenap0-097} and~\eqref{09876-2345678-83uejdn3874},
if~$x\in\partial B_{r_0}$,
\begin{equation} \label{89-0987654-ilSM-Djmsd-rumc-8iedma023}
U_{m,m'}(x)\ge u_m(x)-u_{m'}(x)\ge\frac9{10}-\frac1{10}=\frac45.\end{equation}
Now we claim that
\begin{equation}\label{MOSND1IAPKSLEEUMDFRONBDOERPENMDO}
U_{m,m'}(x)\ge\frac45\qquad{\mbox{a.e. }}x\in B_{r_0}.\end{equation}
To establish this, we use the
approximation method in~\eqref{CONJCONDENS}, recalling that~$U_{m,m'}$ can be seen
as the a.e. limit as~$j\to+\infty$ of the sequence of conductor potentials~$U_{m,m',j}$
corresponding to open and bounded sets~${\mathcal{U}}_{m,m',j}$ with
boundary of class~$C^{2,\alpha}$ which approximate~$(B_{\rho_m}\setminus B_{\rho_{m'}})\setminus\Omega$ from the exterior.
In this way, the function~$\phi_{m,m',j}:=U_{m,m',j}-U_{m,m'}
$ is harmonic outside~${\mathcal{U}}_{m,m',j}$, it belongs to~${\mathcal{D}}^{1,2}(\R^n)$, and along~$\partial{\mathcal{U}}_{m,m',j}$ we have that~$\phi_{m,m',j}:=1-U_{m,m'}\ge0$.
Consequently, recalling Lemma~\ref{SOBOMAXPLEH10}, we have that~$\phi_{m,m',j}\ge0$
in~$\R^n$. This and~\eqref{89-0987654-ilSM-Djmsd-rumc-8iedma023} give that for every~$x\in
\partial B_{r_0}$,
\begin{equation} \label{89-0987654-ilSM-Djmsd-rumc-8iedma023-90}
U_{m,m',j}(x)\ge U_{m,m'}(x)\geq\frac45.\end{equation}
Let now~$\e\in\left(0,\frac12\right)$ and
$$U_{m,m',j}^\star:=\begin{dcases}
U_{m,m',j} & {\mbox{ in }}\R^n\setminus B_{r_0},\\
\max\left\{
U_{m,m',j} ,
\frac45-\e\right\}& {\mbox{ in }} B_{r_0}.
\end{dcases}$$
We stress that~$U_{m,m',j}^\star\in {\mathcal{D}}^{1,2}(\R^n)$, thanks to~\eqref{89-0987654-ilSM-Djmsd-rumc-8iedma023-90}. Furthermore,
$$ \int_{\R^n\setminus {\mathcal{U}}_{m,m',j}}|\nabla U_{m,m',j}^\star(x)|^2\,dx
\le\int_{\R^n\setminus {\mathcal{U}}_{m,m',j}}|\nabla U_{m,m',j}(x)|^2\,dx.$$
This and the minimality property of~$U_{m,m',j}$ (as in~\eqref{DEFCAP-SOBO} 
and~\eqref{CAB9IJSalskdmnfcx34546fsdQsdfasdAOSL}) yield that
\begin{eqnarray*} 0&=&\int_{\R^n\setminus {\mathcal{U}}_{m,m',j}}|\nabla U_{m,m',j}(x)|^2\,dx-
\int_{\R^n\setminus {\mathcal{U}}_{m,m',j}}|\nabla U_{m,m',j}^\star(x)|^2\,dx\\&
=&\int_{B_{r_0}\cap\left\{ U_{m,m',j}<\frac45-\e\right\}}|\nabla U_{m,m',j}(x)|^2\,dx.\end{eqnarray*}
{F}rom this, it follows that~$U_{m,m',j}\ge\frac45-\e$ in~$B_{r_0}$.
Thus, sending~$\e\searrow0$, we conclude that~$U_{m,m',j}\ge\frac45$ in~$B_{r_0}$.
Now, sending~$j\to+\infty$, we obtain~\eqref{MOSND1IAPKSLEEUMDFRONBDOERPENMDO}.

Now we write
$$ (B_{\rho_m}\setminus B_{\rho_{m'}})\setminus\Omega
=\bigcup_{i=m}^{m'-1}
(B_{\rho_i}\setminus B_{\rho_{i+1}})\setminus\Omega$$
and we denote by~$w_i$ the conductor potential of~$(B_{\rho_i}\setminus B_{\rho_{i+1}})\setminus\Omega$, as given in Corollary~\ref{gujbmSIJKc-qpwruofygheibdowi3egufewbuigf3289eufgvosud81ye2rit23456gh}. We then exploit
Lemma~\ref{KJSM-098765-TFI-IMP245S-PLSIUDND-94} (as refined in Lemma~\ref{KJSM-098765-TFI-IMP245S-PLSIUDND-95}) to see that
\begin{equation}\label{654jdsj289uwgfweibfi4bgeahsjdlamfovlla-0} U_{m,m'}\le\sum_{i=m}^{m'-1} w_i.\end{equation}
We also observe that~$B_{\rho_{i+1}}$
lies in the complement of~$(B_{\rho_i}\setminus B_{\rho_{i+1}})\setminus\Omega$
and therefore, in view of~\eqref{CAHBSN123s4345453124135dcHNSdfz2354tyhdq3z4erifjg-2}
(used here with~$\big(B_{\rho_i}\setminus B_{\rho_{i+1}}\big)\setminus\Omega$
in place of~$\Omega$, $u:=w_i$ and~$r:=\rho_{i+1}$),
$$ w_i(0)\le\frac{c_n \cAPAC\Big(\big(B_{\rho_i}\setminus B_{\rho_{i+1}}\big)\setminus\Omega\Big)}{\rho_{i+1}^{n-2}}=\frac{c_n \cAPAC\Big(\big(B_{\lambda^i}\setminus B_{\lambda^{i+1}}\big)\setminus\Omega\Big)}{\lambda^{n-2}\,\lambda^{(n-2)i}}.$$
This, \eqref{0ok0o90uytfcikasmcolalijfnggkkaal99304} and~\eqref{654jdsj289uwgfweibfi4bgeahsjdlamfovlla-0} lead to
$$ U_{m,m'}(0)\le\frac1{\lambda^{n-2}}
\sum_{i=m}^{m'-1}
\frac{c_n \cAPAC\Big(\big(B_{\lambda^i}\setminus B_{\lambda^{i+1}}\big)\setminus\Omega\Big)}{\lambda^{(n-2)i}}\le\frac1{10}.$$
This is in contradiction with~\eqref{MOSND1IAPKSLEEUMDFRONBDOERPENMDO}
and the proof of Theorem~\ref{NOWI} is complete.
\end{proof}

By inspecting its proof, one also sees that
the Wiener's criterion in Theorem~\ref{NOWI} is actually a ``pointwise''
statement, namely it gives a sharp condition on the solvability
of the Dirichlet problem by a function that attains a prescribed boundary datum
at a given point, depending on condition~\eqref{SOMMAIWIE}
(this type of points are sometimes referred to with the name of
``regular points'' for the Dirichlet problem).\medskip

We also observe that the Wiener's criterion in Theorem~\ref{NOWI} 
contains the previous existence result for the Dirichlet problem in Theorem~\ref{PERRO-2}
as a byproduct. Indeed, take a bounded open set satisfying the exterior cone condition.
Up to a rigid motion, we can assume that the origin lies on the boundary of~$\Omega$
and we use the exterior cone condition to write that
\begin{equation}\label{0-9402ruyfewh-0-31-04-72498753jhfh}
\big\{x=(x',x_n)\in B_{r_0} {\mbox{ s.t. }} x_n\ge b|x'|\big\}\subseteq\R^n\setminus\Omega,
\end{equation}
for some~$r_0\in(0,1)$ and $b>0$.

We let~$i\in\N\cap\left( \frac{|\ln r_0|}{\ln2},+\infty\right)$
and take~$\lambda:=\frac12$.
We observe that~$\lambda^i<r_0$.

We define
$$ p_i:=\frac{3}{2^{i+2}}\,e_n\qquad{\mbox{and}}\qquad
r_i:=\frac{1}{(1+b)\,2^{i+3}}$$
and we claim that
\begin{equation}\label{0-0-verisjfun-unattjnfe-mo03}
\big(B_{\lambda^i}\setminus {B_{\lambda^{i+1}}}\big)\setminus\Omega\;\supseteq\;
B_{r_i}(p_i).
\end{equation}
To check this, let~$x=(x',x_n)\in B_{r_i}(p_i)$. We have that
\begin{equation}\label{0-0-verisjfun-unattjnfe-mo01}
|x|\le|p_i|+r_i\le
\frac{3}{2^{i+2}}+\frac{1}{2^{i+3}}=\frac{7}{2^{i+3}}<\frac{8}{2^{i+3}}=\frac{1}{2^{i}}
\end{equation}
and
\begin{equation}\label{0-0-verisjfun-unattjnfe-mo02}
|x|\ge|p_i|-r_i\ge
\frac{3}{2^{i+2}}-\frac{1}{2^{i+3}}=\frac{5}{2^{i+3}}>\frac{4}{2^{i+3}}=\frac{1}{2^{i+1}}.
\end{equation}
Furthermore,
\begin{eqnarray*}
x\cdot e_n-b|x'|\ge
p_i\cdot e_n-r_i-b r_i=
\frac{3}{2^{i+2}}-(1+b)r_i
=\frac{3}{2^{i+2}}-\frac{1}{2^{i+3}}>0.
\end{eqnarray*}
{F}rom this, \eqref{0-9402ruyfewh-0-31-04-72498753jhfh},
\eqref{0-0-verisjfun-unattjnfe-mo01} and~\eqref{0-0-verisjfun-unattjnfe-mo02}
we obtain~\eqref{0-0-verisjfun-unattjnfe-mo03}, as desired.

Now we combine
Lemma~\ref{PTEJSN5678dsdvfg9CALsdfLDM-98uyhnSM6789SMMSA} (recall also
Lemma~\ref{lemakoe84654656548489}),
Corollary~\ref{CAPAPALLA} and~\eqref{0-0-verisjfun-unattjnfe-mo03}: in this way, we infer that
$$ \cAPAC\Big(
\big(B_{\lambda^i}\setminus {B_{\lambda^{i+1}}}\big)\setminus\Omega
\Big)\ge\cAPAC\big(B_{r_i}(p_i)\big)=\frac{r_i^{n-2}}{c_n}
=
\frac1{c_n\,(1+b)^{n-2}\,2^{(i+3)(n-2)}}=\frac{c}{\lambda^{(n-2)i}},$$
for some~$c>0$ depending only on~$n$ and~$b$.

For this reason,
$$ \sum_{i=0}^{+\infty}\frac{\cAPAC\Big(\big(
B_{\lambda^i}\setminus {B_{\lambda^{i+1}}}\big)\setminus\Omega\Big)}{\lambda^{(n-2)i}}=+\infty.$$
Hence, we can apply Theorem~\ref{NOWI}
and obtain that the Dirichlet problem can be solved by a harmonic function
up to the boundary when~$\Omega$ is an open bounded set satisfying the exterior cone condition.\medskip

The Wiener's criterion in Theorem~\ref{NOWI} 
also contains as a byproduct the Lebesgue spine discussed on page~\pageref{SILEBE678NE},
see~\cite[Exercise~10, page~334]{MR0222317}
(see also~\cite{Cai2007PotentialFO} for related details on ellipsoidal potentials).\medskip

For thoroughgoing explanations on the notions
of capacity see e.g.~\cite{MR0222317, MR0225986, MR0492152, MR1411441, MR1440907, MR1447439, MR1814344, MR2305115, MR2867756, MR3409135, MR3675703, MR3791463}.

\chapter{Some interesting problems coming from the Poisson equation}

We collect in this chapter some
classical and fascinating mathematical problems. The kinship between these problems
is the use of elliptic partial differential equations and
a geometric flavor of the problem and its solution.
More specifically, 
in Section~\ref{SEC:TheSoapBubbleTheorem} we will present the
Soap Bubble Theorem (according to which soap bubbles are round).
This result is quite related to the isoperimetric problem, discussed in
Section~\ref{ISOPE:SECTI} (giving that balls minimize perimeter for a given volume).
Then, in Section~\ref{SEC:OVER-SERRINTheorem}
we will focus on a classical overdetermined problem
(roughly speaking, addressing the
type of domains which allow overabundant boundary prescriptions)
and in Section~\ref{KURAN-SEC}
we deal with the converse of the Mean Value Formula
(determining the type of domains for which a claim such as the one in
Theorem~\ref{KAHAR}
holds true).

The proofs will use a number of astute integration by parts and integral formulas.
A different, and even more geometric approach, to soap bubbles
and overdetermined problems will be presented in Chapter~\ref{MO:PLAN}.

\section{The Soap Bubble Theorem}\label{SEC:TheSoapBubbleTheorem}

A classical result\index{Soap Bubble Theorem|(}
of geometric analysis due to
Aleksandr Danilovi\v{c} Aleksandrov
(see~\cite{MR0150710})
is that a bounded domain whose boundary has constant mean curvature
\index{mean curvature} is
necessarily a ball. Since soap bubbles
are obtained by the balance between surface tension and 
air pressure, they provide nice concrete examples of constant mean curvature surfaces,
and this result is thereby often referred to as the ``Soap Bubble Theorem''
(further details will be given in Theorem~\ref{SOAPTH} and Corollaries~\ref{CO-D-SC}
and~\ref{CO-D} below).
Following~\cite{MR645791, MR996826, MR4124125}, we give here a proof of this result based only on the existence
result for the Dirichlet problem provided by Theorem~\ref{PERRO-2}
and on the geometric identities presented in Section~\ref{BELS}
(the original approach based on the reflection principle method
will be also discussed in Section~\ref{MPPALEB}).
\medskip

We start by introducing some notation.
Given a $C^2$ function~$u$, recalling the setting in~\eqref{D22H}
we define
\begin{equation}\label{DEFICDU} {\mathcal{D}}_u:=n\,|D^2u|^2-(\Delta u)^2.\end{equation}
Exploiting the Cauchy-Schwarz Inequality, we see that
\begin{equation*} (\Delta u)^2=
\left(\sum_{i=1}^n \partial_{ii} u\right)^2\le
\sum_{i=1}^n (\partial_{ii} u)^2\; \sum_{i=1}^n 1\le n\,|D^2u|^2,
\end{equation*}
therefore
\begin{equation}\label{CHA -04}
{\mathcal{D}}_u\ge0
\end{equation}
and
\begin{equation}\label{HOLDUJN}
{\mbox{equality holds if and only if the Hessian of~$u$ is a
scalar function times the identity matrix.}}\end{equation}

It is also convenient to
consider an auxiliary equation
set on a bounded, connected and open set~$\Omega\subseteq\R^n$
with~$C^{2,\alpha}$ boundary, with~$\alpha\in(0,1)$, whose mean curvature will be denoted by~$H$.
Namely, in this setting, we consider solutions of
\begin{equation}\label{EQ:AUX} \begin{dcases}
\Delta u =1 & {\mbox{ in }}\Omega,\\
u=0&{\mbox { on }}\partial\Omega.
\end{dcases}\end{equation}
In this setting, it is useful to
observe that, by~\eqref{1111DIV0934},
\begin{equation} \label{HBMUL0}\int_{\partial\Omega}\partial_\nu u(x)
\,d{\mathcal{H}}^{n-1}_x=\int_\Omega \Delta u(x)\,dx=|\Omega|,\end{equation}
and that, by Corollary~\ref{C228}, on~$\partial\Omega$,
\begin{equation}\label{HBMUL} 1=\Delta u=H\,\partial_\nu u+\partial_{\nu\nu} u
.\end{equation}
Also, we have the following cornerstone result
(see e.g.~\cite{MR645791}):

\begin{lemma}\label{LA:MNSBB}
Consider a bounded, connected and open set~$\Omega\subseteq\R^n$
with~$C^{2,\alpha}$ boundary for some~$\alpha\in(0,1)$.
Let~$u\in C^3(\Omega)\cap
C^2(\overline{\Omega})$ be a solution of~\eqref{EQ:AUX}.

Then,
\begin{equation}\label{ANd-3}
\int_{\partial\Omega}\partial_{\nu\nu} u(x)\,\partial_\nu u(x)\,d{\mathcal{H}}^{n-1}_x\ge
\frac{|\Omega|}{n},\end{equation}
\begin{equation}\label{ANd-3A:SK0trp}
\int_\Omega |D^2u(x)|^2\,dx\ge
\frac{|\Omega|}{n}\end{equation}
and
\begin{equation}\label{LANUAN}
\int_{\partial\Omega}H(x)\,(\partial_\nu u(x))^2\,d{\mathcal{H}}^{n-1}_x\le
\frac{n-1}{n}\,|\Omega|.
\end{equation}
Also, the following conditions are equivalent:
\begin{equation}\label{EQUIVABP-00}
{\mbox{equality holds instead of inequality
in~\eqref{ANd-3},}}\end{equation}
\begin{equation}\label{EQUIVABP-006}
{\mbox{equality holds instead of inequality
in~\eqref{ANd-3A:SK0trp},}}\end{equation}
\begin{equation}\label{EQUIVABP-0}
{\mbox{equality holds instead of inequality
in~\eqref{LANUAN},}}\end{equation}
\begin{equation}\label{EQUIVABP-1}
{\mbox{${\mathcal{D}}_u$ is identically zero in~$\Omega$,}}\end{equation}
\begin{equation}\label{EQUIVABP-2}
{\mbox{there exist~$x_0\in\R^n$ and~$c\in\R$ such that }}\,
u(x)=\frac{|x-x_0|^2}{2n}-c\,{\mbox{ for every~$x\in\Omega$}}.
\end{equation}

Moreover, if one of the equivalent conditions~\eqref{EQUIVABP-00}, \eqref{EQUIVABP-006},
\eqref{EQUIVABP-0},
\eqref{EQUIVABP-1}
and~\eqref{EQUIVABP-2} is satisfied, then~$\Omega$ is necessarily a ball.
\end{lemma}

\begin{proof}
Multiplying~\eqref{HBMUL} by~$\partial_\nu u$ and integrating over~$\partial\Omega$,
we obtain that
\begin{equation*}\int_{\partial\Omega}\partial_\nu u(x)
\,d{\mathcal{H}}^{n-1}_x=
\int_{\partial\Omega}H(x)\,(\partial_\nu u(x))^2\,d{\mathcal{H}}^{n-1}_x
+\int_{\partial\Omega}
\partial_{\nu\nu} u(x)\,\partial_\nu u(x)\,d{\mathcal{H}}^{n-1}_x.
\end{equation*}
This and~\eqref{HBMUL0} give that
\begin{equation}\label{SM-p5ncu}
\int_{\partial\Omega}H(x)\,(\partial_\nu u(x))^2\,d{\mathcal{H}}^{n-1}_x
+\int_{\partial\Omega}
\partial_{\nu\nu} u(x)\,\partial_\nu u(x)\,d{\mathcal{H}}^{n-1}_x=|\Omega|.
\end{equation}
Also, since~$u=0$ along~$\partial\Omega$, we have
that~$\nabla u =|\nabla u|\nu$ on~$\partial\Omega$.
As a result, on~$\partial\Omega$,
$$ \frac12\partial_\nu |\nabla u|^2=\nabla u\cdot\partial_\nu\nabla u
=|\nabla u|\nu\cdot\nabla\partial_\nu u=|\nabla u|\,\partial_{\nu\nu}u=
\nabla u \cdot\nu\,\partial_{\nu\nu}u=\partial_{\nu}u\,\partial_{\nu\nu}u
.$$
Therefore,
using the first
Green Identity~\eqref{GRr1}
and the Bochner Identity in Lemma~\ref{BOCHIFRENTYHAS},
\begin{equation}\label{CHA -03}
\begin{split}&
\int_{\partial\Omega}\partial_{\nu\nu} u(x)\,\partial_\nu u(x)\,d{\mathcal{H}}^{n-1}_x=
\frac12\,
\int_{\partial\Omega}\partial_\nu|\nabla u(x)|^2\,d{\mathcal{H}}^{n-1}_x
=\frac12\,
\int_\Omega \Delta|\nabla u(x)|^2\,dx\\&\qquad\qquad=
\int_\Omega \Big(|D^2u(x)|^2+\nabla u(x)\cdot\nabla\Delta u(x)\Big)\,dx=
\int_\Omega |D^2u(x)|^2\,dx.\end{split}
\end{equation}
{F}rom~\eqref{CHA -03} and~\eqref{CHA -04}, it follows that
\begin{equation*}
\int_{\partial\Omega}\partial_{\nu\nu} u(x)\,\partial_\nu u(x)\,d{\mathcal{H}}^{n-1}_x=
\frac1n
\int_{\Omega}\Big( {\mathcal{D}}_u(x)+(\Delta u(x))^2
\Big)\,dx
=\frac1n
\int_{\Omega}\Big( {\mathcal{D}}_u(x)+1
\Big)\,dx
\ge
\frac{|\Omega|}n,\end{equation*}
with equality holding if and only if~${\mathcal{D}}_u$ vanishes identically in~$\Omega$,
that is if and only if~\eqref{HOLDUJN} holds true at any point of~$\Omega$.
Namely,
we have established~\eqref{ANd-3} and, recalling~\eqref{CHA -03},
also~\eqref{ANd-3A:SK0trp}, and moreover
we have proved
the fact that equality
in~\eqref{ANd-3} is equivalent to equality
in~\eqref{ANd-3A:SK0trp} and to~\eqref{EQUIVABP-1}.
That is, \eqref{EQUIVABP-00},
\eqref{EQUIVABP-006}
and~\eqref{EQUIVABP-1} are equivalent.

Furthermore, in light of~\eqref{SM-p5ncu} and~\eqref{CHA -03}, we have that~\eqref{ANd-3},
\eqref{ANd-3A:SK0trp} and~\eqref{LANUAN} are equivalent.
We also point out
that~\eqref{EQUIVABP-00} (as well as~\eqref{EQUIVABP-006})
is equivalent to~\eqref{EQUIVABP-0},
thanks to~\eqref{SM-p5ncu} (and~\eqref{CHA -03}).

It remains to show that~\eqref{EQUIVABP-1}
is equivalent to~\eqref{EQUIVABP-2}. On the one hand, if~\eqref{EQUIVABP-2}
holds true, then one can perform a direct computation and show that~\eqref{HOLDUJN}
is satisfied, hence~\eqref{EQUIVABP-1} holds true as well.
Conversely, if~\eqref{EQUIVABP-1} is satisfied,
we have that~$D^2 u=\varphi\,{\rm Id}$,
where~${\rm Id}$ is the identity matrix and~$\varphi$ is some function.
As a result, we find that~$ 1=\Delta u=n\varphi$ and thus~$\varphi$ must be constantly equal to~$\frac1n$.
Consequently,
$$ \partial_{ij} u(x)=\frac{\delta_{ij}}n\qquad{\mbox{for all }}x\in \Omega.$$
This identity can be integrated, producing that
$$ u(x)=\frac{x_1^2+\dots+x_n^2}{2n}+\omega\cdot x\qquad{\mbox{for all }}x\in \Omega,$$
for some~$\omega\in\R^n$. We can therefore ``complete the square'' and obtain
$$ u(x)=\left|\frac{x}{\sqrt{2n}}+\frac{\sqrt{n}\,\omega}{\sqrt2}\right|^2-
\frac{n\,|\omega|^2}{2}
\qquad{\mbox{for all }}x\in \Omega.$$
This gives~\eqref{EQUIVABP-2} and completes the proof of the equivalence between~\eqref{EQUIVABP-00},
\eqref{EQUIVABP-006},
\eqref{EQUIVABP-0},
\eqref{EQUIVABP-1} and~\eqref{EQUIVABP-2}.

Suppose now that
one of the equivalent conditions~\eqref{EQUIVABP-00}, \eqref{EQUIVABP-006},
\eqref{EQUIVABP-0}, \eqref{EQUIVABP-1}
and~\eqref{EQUIVABP-2} is satisfied. Since these conditions are equivalent,
we may assume that~\eqref{EQUIVABP-2} holds true, and therefore
the level set~$ \{u=0\}$ is the sphere~$\partial B_{\sqrt{2nc}}(x_0)$, hence~$\partial \Omega$ is a sphere
and~$\Omega$ is a ball.
\end{proof}

The result in Lemma~\ref{LA:MNSBB} is pivotal for many applications.
A classical one is the so-called Heintze-Karcher Inequality\index{Heintze-Karcher Inequality}, which was established in~\cite{MR533065}
and, as pointed out in~\cite{MR996826}, can be proved directly
by using the tools developed so far:

\begin{lemma}\label{HBDD}
Let~$\alpha\in(0,1)$.
Consider a bounded, connected and open set~$\Omega\subseteq\R^n$
with~$C^{2,\alpha}$ boundary for some~$\alpha\in(0,1)$, with strictly positive mean curvature~$H$.

Then,
\begin{equation}\label{HEI}
\int_{\partial\Omega}\frac{d{\mathcal{H}}^{n-1}_x}{H(x)}\ge \frac{n}{n-1}\,|\Omega|.
\end{equation}
Also, equality holds true if and only if~$\Omega$ is a ball.
\end{lemma}

\begin{proof}
By Theorem~\ref{Theorem6.14GT}, there exists a solution~$u\in C^3(\Omega)\cap
C^2(\overline{\Omega})$ of~\eqref{EQ:AUX}.

By~\eqref{HBMUL0} and the Cauchy-Schwarz Inequality,
\begin{equation*}
|\Omega|^2=\left(
\int_{\partial\Omega}\partial_\nu u(x)\,d{\mathcal{H}}^{n-1}_x\right)^2\le
\int_{\partial\Omega}H(x)\,(\partial_\nu u(x))^2\,d{\mathcal{H}}^{n-1}_x\;
\int_{\partial\Omega}\frac{d{\mathcal{H}}^{n-1}_x}{H(x)}.\end{equation*}
This and~\eqref{LANUAN} give~\eqref{HEI}.
Also, if equality holds in~\eqref{HEI}, then equality must hold in~\eqref{LANUAN},
that is condition~\eqref{EQUIVABP-0} is fulfilled.
Accordingly, by
Lemma~\ref{LA:MNSBB}, $\Omega$ is necessarily a ball.

Viceversa, if~$\Omega$ is a ball of radius~$R$, then~$H=\frac{n-1}{R}$, and therefore
\begin{equation*} \int_{\partial\Omega}\frac{d{\mathcal{H}}^{n-1}_x}{H(x)}-\frac{ n\,|\Omega|}{n-1}
=\frac{R\,{\mathcal{H}}^{n-1}(\partial B_R)}{n-1}-\frac{n\,|B_R|}{n-1}
=0.\qedhere\end{equation*}
\end{proof}

As pointed out in~\cite{MR996826},
combining the
Heintze-Karcher Inequality of Lemma~\ref{HBDD}
and the Minkowski Identity in Corollary~\ref{MIK},
a straightforward proof of the Soap Bubble Theorem can be obtained.
The details go as follows:

\begin{theorem}\label{SOAPTH} Let~$\alpha\in(0,1)$.
Consider a bounded, connected
and open set~$\Omega\subseteq\R^n$,
and assume that~$\Sigma:=\partial \Omega$ is of class~$C^{2,\alpha}$.
Assume that the mean curvature of~$\Sigma$ is constant.
Then, $\Omega$ is necessarily a ball.
\end{theorem}

\begin{proof} 
By Corollary~\ref{COsdRPOH} we deduce that the mean curvature~$H$ of~$\Sigma$
is necessarily a positive constant. Hence, in light of the
Minkowski Identity in Corollary~\ref{MIK},
$$ 
\int_{\Sigma}\frac{d{\mathcal{H}}^{n-1}_x}{H(x)}=
\frac{{\mathcal{H}}^{n-1}(\Sigma)}{H}=\frac1{(n-1)H}
\int_\Sigma H(x)\,x\cdot\nu(x)\,d{\mathcal{H}}^{n-1}_x=
\frac{1}{n-1}
\int_\Sigma x\cdot\nu(x)\,d{\mathcal{H}}^{n-1}_x.$$
This and the Divergence Theorem yield that
$$ 
\int_{\Sigma}\frac{d{\mathcal{H}}^{n-1}_x}{H(x)}=\frac{1}{n-1}
\int_\Omega\div( x)\,dx=\frac{n}{n-1}\,|\Omega|.$$
That is, equality holds in the
Heintze-Karcher Inequality~\eqref{HEI}
and accordingly Lemma~\ref{HBDD}
guarantees that~$\Omega$ is a ball.\index{Soap Bubble Theorem|)}
\end{proof}

\section{The isoperimetric problem}\label{ISOPE:SECTI}

The Soap Bubble Theorem has strong connections\index{isoperimetric
problem|(} with the so-called ``isoperimetric
problem'', that is the problem of finding the set of prescribed volume
with smallest possible surface area. We will discuss this connection in
Corollaries~\ref{CO-D-SC} and~\ref{CO-D} below. To this end,
we first point out the following ``first variation formula''
for the area functional:

\begin{lemma}\label{334}
Let~$\alpha\in(0,1)$.
Let~$\Omega\subseteq\R^n$ be a bounded and
connected set of class~$C^{2,\alpha}$. Let~$\rho>0$ and~$p\in\partial\Omega$.
Let~$T>0$
and~$\Phi\in C^{2,\alpha}(\R^n\times[-T,T],\R^n)$
be such that~$\Phi(x,0)=x$ for every~$x\in \R^n$.
Let~$\vartheta(x):=\partial_t\Phi(x,0)$ and~$\Omega_t:=\Phi(\Omega,t)$.
Then,
$$ \frac{d}{dt} {\mathcal{H}}^{n-1}(\partial\Omega_t)\Big|_{t=0}=
\int_{\partial\Omega} H(x)\,\vartheta(x)\cdot\nu(x)\,d{\mathcal{H}}^{n-1}_x.$$
\end{lemma}

\begin{proof}  In local coordinates, we can assume that~$\partial\Omega$
is locally parameterized by some
diffeomorphism~$f:D\to\Sigma$,
for some domain~$D\subset\R^{n-1}$, see Figure~\ref{SPA}.
We consider the metrics~$g_{ij}$ as in~\eqref{METRIX},
corresponding to the matrix~${\mbox{\Large{\calligra{g}}}}$.
We let~$\Sigma$ be the portion of~$\partial\Omega$ parameterized by this set
of coordinates and~$\widetilde\Sigma:=\Phi(\Sigma,t)$ for some~$t$ (with~$|t|$ to be taken as small
as we wish in what follows). The induced metrics on~$\widetilde\Sigma$ will be
denoted by~$\widetilde{g}_{ij}$ and the corresponding matrix by~$\widetilde{ {\mbox{\Large{\calligra{g}}}} }$.
Notice that~$\widetilde\Sigma$ is parameterized by~$\widetilde{f}(\eta):=\Phi(f(\eta),t)=f(\eta)+t\vartheta(f(\eta))+o(t)$, with~$\eta\in D$.

We also let \begin{equation}\label{UIS0thisca9j93u5tzspr}
m_{ij}:=
\sum_{k=1}^{n}\left[ \frac{\partial \vartheta}{\partial x_k}\cdot
\frac{\partial{f}}{\partial\eta_i}\,\frac{\partial f_k}{\partial\eta_j}+
\frac{\partial \vartheta}{\partial x_k}\cdot
\frac{\partial{f}}{\partial\eta_j}\,\frac{\partial f_k}{\partial\eta_i}\right]\end{equation} and denote by~${\mbox{\Large{\calligra{m}}}}$ the associated matrix.
Since, for every~$j\in\{1,\dots,n-1\}$,
$$ \frac{\partial \widetilde{f}}{\partial\eta_j}(\eta)=\frac{\partial{f}}{\partial\eta_j}(\eta)+t
\sum_{k=1}^{n}\frac{\partial \vartheta}{\partial x_k}(f(\eta))\,\frac{\partial f}{\partial\eta_j}(\eta)\cdot e_k+o(t),
$$
we infer from the definition of metrics in~\eqref{METRIX} that
$$ \widetilde{g}_{ij}:=\frac{\partial \widetilde{f}}{\partial\eta_i}\cdot\frac{\partial\widetilde{ f}}{\partial\eta_j}=
g_{ij}+t
\sum_{k=1}^{n}\left[ \frac{\partial \vartheta}{\partial x_k}\cdot
\frac{\partial{f}}{\partial\eta_i}\,\frac{\partial f_k}{\partial\eta_j}+
\frac{\partial \vartheta}{\partial x_k}\cdot
\frac{\partial{f}}{\partial\eta_j}\,\frac{\partial f_k}{\partial\eta_i}\right]
+o(t)=g_{ij}+t m_{ij}+o(t),
$$
where~$f_k:=f\cdot e_k$ (also, from now on, $f$ and its derivatives are evaluated at~$\eta\in D$
and~$\vartheta$ and its derivatives are evaluated at~$f(\eta)$, but we omit this dependence for
the sake of brevity).
 As a result, $$\widetilde{\mbox{\Large{\calligra{g}}}}=
{\mbox{\Large{\calligra{g}}}}+t{\mbox{\Large{\calligra{m}}}}+o(t)=
{\mbox{\Large{\calligra{g}}}}\,\Big({\rm Id}+t{\mbox{\Large{\calligra{g}}}}^{-1}{\mbox{\Large{\calligra{m}}}}+o(t)
\Big),$$ whence
\begin{eqnarray*}&&
\det \widetilde{\mbox{\Large{\calligra{g}}}}=
\det {\mbox{\Large{\calligra{g}}}}\;\det\Big({\rm Id}+t{\mbox{\Large{\calligra{g}}}}^{-1}{\mbox{\Large{\calligra{m}}}}+o(t)
\Big)=\det {\mbox{\Large{\calligra{g}}}}\; \left(1+t\sum_{i,j=1}^{n-1} g^{j i}m_{ij}+o(t)
\right).
\end{eqnarray*}
Consequently,
the surface element of~$\widetilde\Sigma$ in local coordinates (see e.g.~\cite[pages 125-126]{MR3409135}), that we denote with a slight abuse
of notation as~$d\widetilde\Sigma$,
can be written as
\begin{equation}\label{Chade24DS0D} d\widetilde\Sigma=
\sqrt{\det \widetilde{\mbox{\Large{\calligra{g}}}}\,}\;d\eta=
\sqrt{\det {\mbox{\Large{\calligra{g}}}}\,}\;\left(1+\frac{t}2\sum_{j=1}^{n-1} g^{ji}m_{ij}+o(t)
\right)
\;d\eta=\left(1+\frac{t}2\sum_{i,j=1}^{n-1} g^{ij}m_{ij}+o(t)
\right)\,d\Sigma.
\end{equation}
Using this and the definition of~$m_{ij}$ in~\eqref{UIS0thisca9j93u5tzspr},
\begin{eqnarray*} &&\frac{1}2\sum_{i,j=1}^{n-1} g^{ij}m_{ij}=
\sum_{{1\le i, j\le n-1}\atop{1\le k\le n}}g^{ij} \frac{\partial \vartheta}{\partial x_k}\cdot
\frac{\partial{f}}{\partial\eta_i}\,\frac{\partial f_k}{\partial\eta_j}\\&&\qquad=
\sum_{{1\le i, j\le n-1}\atop{1\le k\le n}}g^{ij} \frac{\partial \vartheta}{\partial x_k}\cdot
\frac{\partial{f}}{\partial\eta_i}\;\frac{\partial f}{\partial\eta_j}\cdot e_k
=\sum_{i, j=1}^{n-1}g^{ij} \left(D \vartheta
\frac{\partial{f}}{\partial\eta_i}\right)\cdot
\frac{\partial f}{\partial\eta_j}.
\end{eqnarray*}
Therefore, in light of Corollary~\ref{LSDcmevandtimsflVVSND},
$$ \frac{1}2\sum_{j=1}^{n-1} g^{ij}m_{ij}=
\div_\Sigma \vartheta.$$
This and~\eqref{Chade24DS0D} lead to
\begin{equation*} d\widetilde\Sigma=
\big(1+t\div_\Sigma\vartheta+o(t)
\big)\,d\Sigma.
\end{equation*}
As a result,
$$ {\mathcal{H}}^{n-1}(\partial\Omega_t)=
{\mathcal{H}}^{n-1}(\partial\Omega)+t\int_{\partial\Omega}\div_\Sigma\vartheta\,d{\mathcal{H}}^{n-1}+o(t).$$
Hence, recalling the Tangential Divergence Theorem~\ref{GBNAmn} (used here with~$\varphi:=1$),
$$ {\mathcal{H}}^{n-1}(\partial\Omega_t)={\mathcal{H}}^{n-1}(\partial\Omega)+
t\int_{\partial\Omega} H(x)\,\vartheta(x)\cdot\nu(x)\,d{\mathcal{H}}^{n-1}_x+o(t),$$
that produces the desired result.
\end{proof}

Now, we provide some general observations
on vector fields and one-parameter families of diffeomorphisms.
Firstly, we discuss\footnote{On the one hand,
we believe that
Lemma~\ref{FLUSSETTO} is handy and can be conveniently
used in several occasion. On the other hand,
we also point out that the use of Lemma~\ref{FLUSSETTO}
in these notes can be replaced by the flow of an ordinary differential
equation with a prescribed velocity combined with
a space dilation to preserve a given volume,
see e.g. the alternative proof of Corollary~\ref{CO-D-SC}
presented on page~\pageref{KMSmamx15meuntTTYrumi}
and the alternative argument for the proof
of Lemma~\ref{HADA}
given in the footnote of page~\pageref{KMSmamx15meuntTTYrumi2}.}
how to construct a one-parameter family of diffeomorphisms
that preserves the volume of a set and having as initial velocity
a zero-flux vector field:

\begin{lemma}\label{FLUSSETTO}
Let~$\Omega$ be a bounded subset of~$\R^n$
with~$C^2$ boundary and exterior unit normal~$\nu$.
Let~$\vartheta\in C^1(\partial\Omega,\R^n)$ such that
\begin{equation}\label{exnonc}
\int_{\partial\Omega} \vartheta(x)\cdot\nu(x)\,d{\mathcal{H}}^{n-1}_x=0.
\end{equation}
Then, there exist~$T>0$, a neighborhood~${\mathcal{N}}$
of~$\partial\Omega$ and~$\Phi\in C^1((\Omega\cup{\mathcal{N}})\times [-T,T],\R^n)$
such that
\begin{equation}\label{KMD-ef}
{\mbox{$\Phi(\cdot,t)$ is a diffeomorphism between~$\Omega\cup{\mathcal{N}}$
and~$\Phi(\Omega\cup{\mathcal{N}},t)$ for all~$t\in[-T,T]$,}}\end{equation}
\begin{equation}\label{7ds8u8fjjd}
\frac{\partial}{\partial t}\Phi(x,t)\Big|_{t=0}=\vartheta(x)\qquad{\mbox{ for all }}x\in\partial\Omega
\end{equation}
and
\begin{equation}\label{MEAS-r}
{\mbox{the Lebesgue measure of
the set~$\Phi(\Omega,t)$
is equal to the Lebesgue measure of~$\Omega$ for all~$t\in[-T,T]$.}}\end{equation}
\end{lemma}

\begin{proof} We consider
the normal extension of~$\vartheta$ in a suitably small
neighborhood~${\mathcal{N}}$
of~$\partial\Omega$, as defined
in~\eqref{DEEST}. Moreover, we take
a smaller neighborhood~${\mathcal{N}}'\Subset{\mathcal{N}}$
of~$\partial\Omega$ and~$\vartheta^\star\in C^1(\Omega\cup{\mathcal{N}},\R^n)$ such that~$\vartheta^\star=\vartheta_{\mbox{\scriptsize{ext}}}$
in~${\mathcal{N}}'$ and~$\vartheta^\star=0$ in~$\Omega\setminus{\mathcal{N}}$.

Hence, we let~$X^t$ be the flow associated with
the vector field~$\vartheta^\star$, i.e., for~$|t|$ sufficiently small,
the solution of the ordinary differential equation
$$ \begin{dcases}
\frac{d}{dt} X^t(x)=\vartheta^\star(X^t(x)),\\
X^0(x)=x,
\end{dcases}$$
for all~$x\in\Omega\cup {\mathcal{N}}'$.

We also take~$\nu^\star\in C^1(\Omega\cup{\mathcal{N}},\R^n)$ such that~$\nu^\star=\nu_{\mbox{\scriptsize{ext}}}$
in~${\mathcal{N}}'$ and~$\nu^\star=0$ in~$\Omega\setminus{\mathcal{N}}$.
With this notation, 
for~$x\in\Omega\cup{\mathcal{N}}$,
$s\in\R$ and~$|t|$ small, we also define
$$ f(x,t,s):=X^t(x)+s\nu^\star(x).$$
Let also~$ F(t,s)$ be the Lebesgue measure of the set~$f(\Omega,t,s)$,
that is~$F(t,s):=|f(\Omega,t,s)|$. We stress that
$$ f(\Omega,0,s)=\{
G(x,s),\;\,x\in\Omega
\},$$
where~$G(x,s):=x+s\nu^\star(x)$.
We also remark that, for small~$|s|$, we have~$D_x 
G(x,s)={\rm Id}+sD_x\nu^\star(x)$,
where~${\rm Id}$ is the identity matrix, and consequently
$$ \det D_xG(x,s)=1+s\div\nu^\star(x)+o(s).$$
Using the change of variable~$y:=G(x,s)$,
this entails that, for small~$|s|$,
$$ F(0,s)=\int_{f(\Omega,0,s)} dy=\int_{G(\Omega,s)}dy=
\int_\Omega |\det D_xG(x,s)|\,dx=
\int_\Omega \big(1+s\div\nu^\star(x)\big)\,dx+o(s).$$
Hence, from the Divergence Theorem,
\begin{eqnarray*}&& F(0,s)=|\Omega|+s\int_\Omega
\div\nu^\star(x)\,dx+o(s)
=F(0,0)+s\int_{\partial\Omega}\nu^\star(x)\cdot\nu(x)\,
d{\mathcal{H}}^{n-1}_x+o(s)\\&&\qquad\qquad\qquad=F(0,0)+s
{\mathcal{H}}^{n-1}(\partial\Omega)+o(s).\end{eqnarray*}
As a result, we have that
\begin{equation}\label{JSODNND0}
\partial_s F(0,0)=
{\mathcal{H}}^{n-1}(\partial\Omega)\ne0\end{equation} and thus, by the Implicit Function Theorem, there exists a $C^1$ function~$\gamma$ such that~$\gamma(0)=0$ and,
when~$|t|$ is small enough,
\begin{equation}\label{JSODNND} F(t,\gamma(t))=F(0,0),\end{equation}
that is, setting~$\Phi(x,t):=f(x,t,\gamma(t))$,
$$ |\Phi(\Omega,t)|=|f(\Omega,t,\gamma(t))|=
F(t,\gamma(t))=F(0,0) =|f(\Omega,0,0)|=|\Omega|,$$
thus establishing~\eqref{MEAS-r}.

We also point out that~$ D_x \Phi(\Omega\cup{\mathcal{N}},0)={\rm Id}$,
hence~\eqref{KMD-ef} is a consequence of the Inverse Function Theorem.

In addition, for small~$|t|$,
$$ D_x X^t(x)={\rm Id}+t D_x \vartheta^\star(x)+o(t),$$
hence~$\det D_x X^t(x)=1+t \div\vartheta^\star(x)+o(t)$, and, for this reason,
\begin{eqnarray*}&&
F(t,0)=|f(\Omega,t,0)|=\int_{f(\Omega,t,0)}dy=\int_{X^t(\Omega)}dy=
\int_\Omega
\big(1+t \div\vartheta^\star(x)\big)\,dx+o(t)\\&&\qquad
=|\Omega|+t\int_\Omega \div\vartheta^\star(x)\,dx+o(t)=
F(0,0)+t
\int_{\partial\Omega} \vartheta(x)\cdot\nu(x)\,d{\mathcal{H}}^{n-1}_x+o(t)=
F(0,0)+o(t),
\end{eqnarray*}
where~\eqref{exnonc}
has been used in the last step.

Therefore, $\partial_t F(0,0)=0$ and accordingly,
in light of~\eqref{JSODNND0} and~\eqref{JSODNND},
$$ 0=\frac{d}{dt}
F(t,\gamma(t))\Big|_{t=0}
=\partial_t F(0,0)+\partial_s F(0,0)\,\gamma'(0)=
0+{\mathcal{H}}^{n-1}(\Omega)\,\gamma'(0),$$
from which we deduce that~$\gamma'(0)=0$.

Making use of this, we conclude that, for all~$x\in\partial\Omega$,
\begin{eqnarray*}&& \frac{\partial}{\partial t}\Phi(x,t)\Big|_{t=0}=
\frac{\partial}{\partial t}
f(x,t,\gamma(t))\Big|_{t=0}=
\partial_t f(x,t,\gamma(t))+
\partial_s f(x,t,\gamma(t))\,\gamma'(t)\Big|_{t=0}\\&&\qquad
=
\frac{d}{dt} X^t(x)
\Big|_{t=0} +\nu_{\mbox{\scriptsize{ext}}}(x)\,\gamma'(0)=
\vartheta(X^t(x))
\Big|_{t=0}
=\vartheta(x),\end{eqnarray*}
and this proves~\eqref{7ds8u8fjjd}.\end{proof}

Now we point out that a function
which averages to zero on the boundary of a domain against
all normal fluxes of zero-flux vector fields must necessarily be constant:

\begin{lemma}\label{3349056}
Let~$\Omega$ be a bounded subset of~$\R^n$
with~$C^1$ boundary and exterior unit normal~$\nu$.
Let~$f$ be a continuous function on~$\partial\Omega$ such that
$$ \int_{\partial\Omega} f(x)\vartheta(x)\cdot\nu(x)\,d{\mathcal{H}}^{n-1}_x
=0$$
for every~$\vartheta\in C^1(\partial\Omega,\R^n)$ such that~$
\int_{\partial\Omega} \vartheta(x)\cdot\nu(x)\,d{\mathcal{H}}^{n-1}_x=0$.

Then, $f$ is constant on~$\partial\Omega$.
\end{lemma}

\begin{proof}
Let~$f_\e\in C^\infty(\partial\Omega)$ such that~$f_\e\to f$ uniformly
on~$\partial\Omega$ as~$\e\to0$. Let also
$$ \vartheta_\e(x):=\left(f_\e(x)-\fint_{\partial\Omega} f_\e(y)\,d{\mathcal{H}}^{n-1}_y\right)\nu(x).$$
Then, we have that
$$ \int_{\partial\Omega} \vartheta_\e(x)\cdot\nu(x)\,d{\mathcal{H}}^{n-1}_x=
\int_{\partial\Omega} \left(
f_\e(x)-\fint_{\partial\Omega} f_\e(y)\,d{\mathcal{H}}^{n-1}_y
\right)\,d{\mathcal{H}}^{n-1}_x
=0$$
and accordingly
\begin{eqnarray*}&& 0=
\fint_{\partial\Omega} f(x)\vartheta_\e(x)\cdot\nu(x)\,d{\mathcal{H}}^{n-1}_x
=
\fint_{\partial\Omega} f(x)
\left(
f_\e(x)-\fint_{\partial\Omega} f_\e(y)\,d{\mathcal{H}}^{n-1}_y
\right)\,d{\mathcal{H}}^{n-1}_x\\&&\qquad=
\fint_{\partial\Omega} f(x)\,f_\e(x)\,dx
-\left(\fint_{\partial\Omega} f(y)\,d{\mathcal{H}}^{n-1}_y\right)
\left(\fint_{\partial\Omega} f_\e(y)\,d{\mathcal{H}}^{n-1}_y\right).\end{eqnarray*}
Therefore, passing to the limit as~$\e\to0$,
\begin{eqnarray*}&&
\fint_{\partial\Omega}
\left(
f(x)-\fint_{\partial\Omega} f(y)\,d{\mathcal{H}}^{n-1}_y
\right)^2
\,d{\mathcal{H}}^{n-1}_x\\&=&
\fint_{\partial\Omega}
\left[
f^2(x)+\left(\fint_{\partial\Omega} f(y)\,d{\mathcal{H}}^{n-1}_y\right)^2
-2f(x)\,\fint_{\partial\Omega} f(y)\,d{\mathcal{H}}^{n-1}_y
\right]
\,d{\mathcal{H}}^{n-1}_x\\&=&
\fint_{\partial\Omega}f^2(x)\,dx-
\left(\fint_{\partial\Omega} f(y)\,d{\mathcal{H}}^{n-1}_y\right)^2\\&=&\lim_{\e\to0}
\left[\fint_{\partial\Omega} f(x)\,f_\e(x)\,dx
-\left(\fint_{\partial\Omega} f(y)\,d{\mathcal{H}}^{n-1}_y\right)
\left(\fint_{\partial\Omega} f_\e(y)\,d{\mathcal{H}}^{n-1}_y\right)
\right]
\\&=&0,
\end{eqnarray*}
and this shows that~$f(x)=\fint_{\partial\Omega} f(y)\,d{\mathcal{H}}^{n-1}_y
$ for every~$x\in\partial\Omega$.\end{proof}

We can now address the isoperimetric
problem of determining the shape minimizing surface
for a prescribed volume.
Notice that we do not address here
the delicate problem of the existence (and smoothness)
of the surface minimizers for a fixed volume
(this is one of the classical topics in the calculus of variations,
see e.g.~\cite{MR775682, MR2976521, MR3497381}).

\begin{corollary}\label{CO-D-SC}
Let~$\alpha\in(0,1)$.
Let~$\Omega\subseteq\R^n$ be a bounded, open and
connected set of class~$C^{2,\alpha}$ such that~${\mathcal{H}}^{n-1}(\partial\Omega)\le
{\mathcal{H}}^{n-1}(\partial\widetilde\Omega)$ for all~$\widetilde\Omega$
among the bounded, open and
connected sets of class~$C^{2,\alpha}$
with a given
volume. Then, the mean curvature of~$\partial\Omega$
is necessarily constant.
\end{corollary}

\begin{proof} We take a zero-flux vector field~$\vartheta$
and the corresponding one-parameter family of diffeomorphisms~$\Phi(\cdot,t)$
constructed in Lemma~\ref{FLUSSETTO}.
In light of~\eqref{MEAS-r} and of the minimality of~$\Omega$
we have that~${\mathcal{H}}^{n-1}(\partial\Omega)\le
{\mathcal{H}}^{n-1}(\partial\Omega^{(t)})$, for all~$t\in(-T,T)$,
for a suitable~$T>0$,
where~$\Omega_t:=\Phi(\Omega,t)$. In particular,
recalling Lemma~\ref{334},
$$ 0=\frac{d}{dt} {\mathcal{H}}^{n-1}(\partial\Omega_t)\Big|_{t=0}=
\int_{\partial\Omega} H(x)\,\vartheta(x)\cdot\nu(x)\,d{\mathcal{H}}^{n-1}_x.$$
Since this is valid for all
zero-flux vector fields~$\vartheta$,
we can apply Lemma~\ref{3349056} (used here with~$f:=H$)
and deduce that~$H$ is constant along~$\partial\Omega$.
\end{proof}

For completeness, we give also a proof of Corollary~\ref{CO-D-SC}
that does not rely on the construction of the
family of diffeomorphisms in Lemma~\ref{FLUSSETTO}
that preserves the volume of the domain,
but rather exploits the special properties of scaling of the problem
under consideration:

\begin{proof}[Another proof of Corollary~\ref{CO-D-SC}] \label{KMSmamx15meuntTTYrumi}
We consider any
vector field~$\vartheta\in C^1(\R^n,\R^n)$ and the corresponding solution~$\Phi^t$
of the associated Cauchy problem
$$ \begin{dcases}
\frac{d}{dt}\Phi^t(x)=\vartheta(\Phi^t(x)),\\
\Phi^0(x)=x.
\end{dcases}$$
We point out that the solution of the above problem
is well-defined for~$t\in[-T,T]$, for a suitable~$T>0$,
due to the
Existence and Uniqueness Theorem 
for ordinary differential equations.
We set~$\Omega_t:=\Phi^t(\Omega)$. In this situation,
the volume of~$\Omega_t$ is not necessarily the same of the volume
of~$\Omega$, therefore it is convenient to define~$\Omega_t^\star$ as the dilation of~$\Omega_t$ with the same
volume as~$\Omega$.
To this end, we point out that
\begin{equation}\label{TAYikdc} |\Omega_t|=\int_\Omega (1+t\div\vartheta(x))\,dx+o(t)
=|\Omega|+
t\int_{\partial\Omega}\vartheta(x)\cdot\nu(x)\,d{\mathcal{H}}^{n-1}_x
+o(t).\end{equation}
In particular, we have that~$|\Omega_t|>0$
if~$|t|$ is sufficiently small, whence we can define~$
\mu(t):=\frac{|\Omega|}{|\Omega_t|}$.
Let also
$$ \Omega_t^\star:=\big(\mu(t)\big)^{\frac1n}\,\Omega_t.$$
Notice that~$|\Omega_t^\star|=\mu(t)\,|\Omega_t|=|\Omega|$,
thanks to the scaling properties of the Lebesgue measure,
hence, by the minimality property of~$\Omega$,
we infer that
\begin{equation}\label{9090987787654864}
\frac{d}{dt} {\mathcal{H}}^{n-1}(\partial\Omega_t^\star)\Big|_{t=0}=0.
\end{equation}
We also remark that
\begin{equation}\label{6y6y02}
{\mathcal{H}}^{n-1}(\partial\Omega_t^\star)=
\big(\mu(t)\big)^{\frac{n-1}n}\,
{\mathcal{H}}^{n-1}(\partial\Omega_t),
\end{equation}
thanks to the scaling invariance of the Hausdorff measure.
Moreover, by~\eqref{TAYikdc},
\begin{eqnarray*}&&
\big(\mu(t)\big)^{\frac{n-1}n}=\left(
\frac{|\Omega_t|}{|\Omega|}
\right)^{\frac{1-n}n}
=
\left(1+\frac{t}{|\Omega|}
\int_{\partial\Omega}\vartheta(x)\cdot\nu(x)\,d{\mathcal{H}}^{n-1}_x
+o(t)\right)^{\frac{1-n}n}\\&&\qquad\qquad\qquad
=1+\frac{(1-n)\,t}{n\,|\Omega|}
\int_{\partial\Omega}\vartheta(x)\cdot\nu(x)\,d{\mathcal{H}}^{n-1}_x
+o(t).\end{eqnarray*}
As a result, we have that
$$\frac{d}{dt} 
\big(\mu(t)\big)^{\frac{n-1}n}
\Big|_{t=0}=
\frac{(1-n)}{n\,|\Omega|}
\int_{\partial\Omega}\vartheta(x)\cdot\nu(x)\,d{\mathcal{H}}^{n-1}_x.$$
Besides, we know that~$\mu(0)=1$ and thus, by~\eqref{6y6y02},
\begin{eqnarray*}&&\frac{d}{dt} {\mathcal{H}}^{n-1}(\partial\Omega_t^\star)\Big|_{t=0}
=
\frac{d}{dt}\big(\mu(t)\big)^{\frac{n-1}n}\Big|_{t=0}\,
{\mathcal{H}}^{n-1}(\partial\Omega)+
\big(\mu(0)\big)^{\frac{n-1}n}\,\frac{d}{dt}
{\mathcal{H}}^{n-1}(\partial\Omega_t)\Big|_{t=0}\\&&\qquad=
\frac{(1-n)\,{\mathcal{H}}^{n-1}(\partial\Omega)}{n\,|\Omega|}
\int_{\partial\Omega}\vartheta(x)\cdot\nu(x)\,d{\mathcal{H}}^{n-1}_x
+\frac{d}{dt}
{\mathcal{H}}^{n-1}(\partial\Omega_t)\Big|_{t=0}.\end{eqnarray*}
Combining this information with~\eqref{9090987787654864}
and Lemma~\ref{334}, we gather that
\begin{eqnarray*} 0&=&
\frac{(1-n)\,{\mathcal{H}}^{n-1}(\partial\Omega)}{n\,|\Omega|}
\int_{\partial\Omega}\vartheta(x)\cdot\nu(x)\,d{\mathcal{H}}^{n-1}_x
+\int_{\partial\Omega} H(x)\,\vartheta(x)\cdot\nu(x)\,d{\mathcal{H}}^{n-1}_x\\&=&
\int_{\partial\Omega} \left(
\frac{(1-n)\,{\mathcal{H}}^{n-1}(\partial\Omega)}{n\,|\Omega|}+
H(x)\right)\,\vartheta(x)\cdot\nu(x)\,d{\mathcal{H}}^{n-1}_x.
\end{eqnarray*}
This and the fact that~$\vartheta$ is an arbitrary vector field yield
that~$H$ is constantly equal to~$\frac{(n-1)\,{\mathcal{H}}^{n-1}(\partial\Omega)}{n\,|\Omega|}$
along~$\partial\Omega$.
\end{proof}

Combining Corollary~\ref{CO-D-SC} with the
Soap Bubble Theorem in~\ref{SOAPTH} we immediately
obtain the desired classification result for the isoperimetric problem:

\begin{corollary}\label{CO-D}
Let~$\alpha\in(0,1)$
and~$\Omega\subseteq\R^n$ be a bounded, open and
connected set of class~$C^{2,\alpha}$ such that~${\mathcal{H}}^{n-1}(\partial\Omega)\le
{\mathcal{H}}^{n-1}(\partial\widetilde\Omega)$
among all the 
bounded, open and connected sets~$\widetilde\Omega$
of class~$C^{2,\alpha}$
with a given
volume. Then, $\Omega$ is a ball.
\end{corollary}

\begin{sloppypar}
See Andrejs Treibergs' review
{\tt http://www.math.utah.edu/~treiberg/isoperim/isop.pdf}
for several proofs\footnote{For completeness, we mention that the smoothness assumption on the set for the isoperimetric problem can be relaxed, see e.g.~\cite{MR98331, MR0493768, MR2147710}.} of the isoperimetric problem. 

For further readings about the Soap Bubble Theorem and related topics see also
Nicola Garofalo's notes {\tt https://www.math.purdue.edu/$\sim$garofalo/Soap\_bubble.pdf}
and the references therein.

\index{isoperimetric
problem|)}
\end{sloppypar}

\section{Serrin's overdetermined problem}\label{SEC:OVER-SERRINTheorem}

In the celebrated article~\cite{MR333220}, James Serrin\index{Serrin's overdetermined problem|(}
established an important result about the radial symmetry
of solutions of an ``overdetermined'' Poisson equation.
Here, the overdetermination comes from the fact
that two boundary conditions are prescribed, instead of only one.
As clearly discussed in~\cite{MR333220},
the problem arose from considerations in physics (especially
fluid dynamics and mechanics) and in fact the work came
in response to a question posed by Roger Fosdick, one of Serrin's colleagues at
the University of Minnesota.
As remarked in~\cite{MR333220},
from the physical viewpoint, Serrin's result entails that:
\begin{itemize}
\item given
``a viscous incompressible fluid moving in straight parallel streamlines
through a straight pipe of given cross sectional form [...]
the tangential stress on the
pipe wall is the same at all points of the wall if and only if the pipe has a circular
cross section'';
\item
given ``a solid straight bar subject to torsion, the magnitude of the resulting traction
which occurs at the surface of the bar is independent of position if and only if the
bar has a circular cross section''.\end{itemize}

We recall that the physical description of the viscous flow in a straight pipe was given in Section~\ref{EUMSD-OS-32456i7DNSSE234R:SEC}
and the derivation of the equation for straight bars under torsion was presented in
Section~\ref{RODSBACKSE}.
\medskip

In Corollary~\ref{SVP} below, we also recall another motivation for Serrin's overdetermined equation dealing with
the shape of a prismatic bar that
maximizes the torsional rigidity.\medskip

Here, we provide a proof of a specific case of Serrin's overdetermined problem using
a strategy developed by Hans Weinberger in~\cite{MR333221}:
interestingly,
this proof came out shortly after
the original argument by Serrin and it is contained
in the article following Serrin's one and published in the same journal.
These early works resulted in an incredibly fruitful research field which still possesses a distinguished role in the contemporary research.
See also~\cite{MR0043486, MR1021402, MR1626395, MR3802818, MR4124125} and the references
therein for further details
about Serrin's problem and alternative points of view about it;
the original approach based on the moving plane method will be also presented in Section~\ref{BACKSE}.
See also~\cite{MR4230553} for detailed
physical motivations for Serrin's overdetermined problem.

\begin{theorem}\label{SERR}
Let~$\Omega\subseteq\R^n$ be a bounded, open and connected set
of class~$C^2$
and let~$u\in C^2(\overline{\Omega})$ be a solution of
\begin{equation}\label{HESE} \begin{dcases}
\Delta u=1&{\mbox{ in }}\Omega,\\
u=0&{\mbox{ on }}\partial\Omega,\\
\partial_\nu u=c&{\mbox{ on }}\partial\Omega,
\end{dcases}\end{equation}
for some~$c\in\R$. Then, $\Omega$ is necessarily a ball.
\end{theorem}

We stress that the equation considered in~\cite{MR333220}
is actually more general than the one presented here, but we stick
to the model case in~\eqref{HESE} since it allows one to develop different approaches
and perspectives to the problem (more general equations
instead can only be treated by some of the methods presented).

\begin{proof}[Proof of Theorem~\ref{SERR}] By~\eqref{HBMUL0},
\begin{equation}\label{3418bis}
|\Omega|=\int_{\partial\Omega}
\partial_\nu u(x)\,d{\mathcal{H}}^{n-1}_x=c\,{\mathcal{H}}^{n-1}(\partial\Omega).\end{equation}
Accordingly, the Poho\v{z}aev Identity
in Theorem~\ref{Pohozaev Identity} and the Divergence Theorem
give that
\begin{equation}\label{POGREGBD}
\begin{split}&
\frac{|\Omega|^2 \,n}{2\,({\mathcal{H}}^{n-1}(\partial\Omega))^2}
=
\frac{c^2 \,n}{2}=
\frac{c^2}{2|\Omega|}
\int_{\Omega} \div(x)\,dx
=
\frac{c^2}{2|\Omega|}
\int_{\partial\Omega} x\cdot\nu(x)\,d{\mathcal{H}}^{n-1}_x\\&\qquad=
\frac1{2|\Omega|}\,\int_{\partial\Omega} (\partial_\nu u(x))^2\,(x\cdot\nu(x))
\,d{\mathcal{H}}^{n-1}_x=
\frac{n-2}{2|\Omega|}\,\int_\Omega u(x)\,dx-\frac{n}{|\Omega|}\int_\Omega u(x)\,dx\\&\qquad=-
\frac{n+2}{2|\Omega|}\,\int_\Omega u(x)\,dx
.\end{split}
\end{equation}
Now we define
$$ {\mathcal{P}}:=\frac{n\,|\nabla u|^2}{2}-u.$$
Notice that
\begin{equation*}
{\mathcal{P}}=\frac{n\,c^2}{2}\qquad{\mbox{on }}\,\partial\Omega.
\end{equation*}
Furthermore, by
the Bochner Identity in Lemma~\ref{BOCHIFRENTYHAS} and
the setting in~\eqref{DEFICDU}, in~$\Omega$
\begin{equation}\label{DEPP} \Delta{\mathcal{P}}=n\,\nabla( \Delta u)\cdot\nabla u
+n\,|D^2u|^{2}
-\Delta u=n|D^2u|^{2}
-(\Delta u)^2={\mathcal{D}}_u.\end{equation}
Consequently, by the inequality in~\eqref{CHA -04},
we have that~$\Delta{\mathcal{P}}
\ge0$ in~$\Omega$.
{F}rom this and the Maximum Principle in Corollary~\ref{WEAKMAXPLE}(i),
we conclude that, for every~$x\in\Omega$,
$$ {\mathcal{P}}(x)\le \sup_{\partial\Omega} {\mathcal{P}}=\frac{n\,c^2}{2},$$
and, more precisely, by
the Strong Maximum Principle in Theorem~\ref{STRONGMAXPLE1}(i),
either
\begin{equation}\label{CON-WEI1}
{\mathcal{P}}(x)<\frac{n\,c^2}{2}
\qquad{\mbox{for every }}\,x\in\Omega
\end{equation}
or
\begin{equation}\label{CON-WEI2}
{\mathcal{P}}(x)=\frac{n\,c^2}{2}
\qquad{\mbox{for every }}\,x\in\Omega.
\end{equation}
We claim that
\begin{equation}\label{89:hold}
{\mbox{condition~\eqref{CON-WEI1} cannot hold true.}}
\end{equation}
To prove this statement, suppose, by contradiction, that~\eqref{CON-WEI1} is satisfied.
Then, using
the first Green's Identity~\eqref{GRr1} (with~$\varphi:=\psi:=u$)
and formula~\eqref{POGREGBD},
$$ \frac{n\,c^2\,|\Omega|}{2}
>\int_\Omega
{\mathcal{P}}(x)\,dx=
\frac{n}{2}\,\int_\Omega
|\nabla u(x)|^2\,dx
-\int_\Omega u(x)\,dx
=
-\frac{n+2}{2}\,\int_\Omega u(x)\,dx
=\frac{n|\Omega|^3}{2\,({\mathcal{H}}^{n-1}(\partial\Omega))^2},$$
which is in contradiction with~\eqref{3418bis},
and this proves~\eqref{89:hold}.

Accordingly, by~\eqref{89:hold},
we deduce that necessarily~\eqref{CON-WEI2}
holds true, and therefore~$\Delta{\mathcal{P}}=0$
in~$\Omega$.

This and~\eqref{DEPP} entail that~${\mathcal{D}}_u$ vanishes identically
in~$\Omega$. That is, condition~\eqref{EQUIVABP-1}
is fulfilled, whence~$\Omega$ is necessarily a ball,
thanks to Lemma~\ref{LA:MNSBB}.
\end{proof}

As an additional motivation for Serrin's overdetermined equation,
we recall the so-called
Saint Venant problem\index{Saint Venant problem} (named after
Adh\'emar Jean Claude Barr\'e de Saint Venant, one of the pioneers of the modern theory of
elasticity and hydraulic engineering).
The problem that we take into account is that of a bar of cross section~$\Omega$,
that is a rigid beam with shape given by~$\Omega\times\R$ (in concrete cases,
one takes~$\Omega\subseteq\R^2$, but we deal with the case~$\Omega\subseteq\R^n$
for the sake of generality). The Saint Venant problem
is focused on the cross sections~$\Omega$ that
maximize, among competitors with the same volume,
the ``torsional rigidity''\index{torsional rigidity}
\begin{equation}\label{TAUP} \tau(\Omega):=\max_{{v\in H^1_0(\Omega)}\atop{v\not\equiv0}}
\frac{\left(\int_\Omega v(x)\,dx\right)^2}{\int_\Omega |\nabla v(x)|^2\,dx}.\end{equation}
The above maximum is attained by the Sobolev Embedding,
see e.g.~\cite{MR1625845}.
We now exploit the fact that,
to maximize~$\tau$,
the cross section~$\Omega$ is necessarily stationary under domain variations:
this method is widely used in shape optimization theory and leads to the
very interesting concept of\index{domain derivative}
``domain derivative'', see e.g.~\cite{MR3791463}. Instead, we do not address here
the problem of the existence (and smoothness)
of the maximizing domain for the torsional rigidity,
see Example~2.8 in~\cite{MR1217590};
see also~\cite{MR2371551} for a general approach towards the existence of extremals
of shape optimization problems (as a matter of fact, the forthcoming Lemma~\ref{HADA}
can be seen as a special case of the
Hadamard Variational Formula presented in~\cite[Chapter~5]{MR3791463}
to which we refer for full details on a very delicate type of computations).

\begin{lemma}\label{HADA}
Let~$\alpha\in(0,1)$.
Let~$\Omega\subseteq\R^n$ be a connected open set of class~$C^{2,\alpha}$
that maximizes
the torsional rigidity~$\tau$ among the 
bounded, open and connected sets
of class~$C^{2,\alpha}$
with a given
volume. Then, there exists a function~$u\in C^2(\overline\Omega)$
such that
$$ \begin{dcases}
\Delta u=1&{\mbox{ in }}\Omega,\\
u=0&{\mbox{ on }}\partial\Omega,\\
\partial_\nu u=c&{\mbox{ on }}\partial\Omega,
\end{dcases}$$
for some~$c\in\R$.
\end{lemma}

\begin{proof}
We let
$$ J(v):=
\int_\Omega |\nabla v(x)|^2\,dx$$
and consider the minimizer~$u_\star$ of~$J$ among every candidate~$v\in H^1_0(\Omega)$
with~$\int_\Omega v(x)\,dx=1$. Then, by the theory of Lagrange multipliers (see e.g.~\cite{MR1625845}),
there exists~$\lambda\in\R$ such that, for every~$\varphi\in C^\infty_0(\Omega)$,
\begin{equation*}
\int_\Omega\nabla u_\star(x)\cdot\nabla\varphi(x)\,dx=\lambda\int_\Omega \varphi(x)\,dx.
\end{equation*}
Hence, setting~$\widetilde u(x):=u_\star(x)+\frac{\lambda\,|x|^2}{2n}$, it follows that~$\widetilde u$
is a weak solution of~$\Delta \widetilde u=0$ in~$\Omega$, since
\begin{eqnarray*}
\int_{\Omega}\nabla \widetilde u(x)\cdot\nabla\varphi(x)\,dx&=&
\int_{\Omega}\nabla u_\star(x)\cdot\nabla\varphi(x)\,dx-
\int_\Omega \nabla\left(\frac{\lambda\,|x|^2}{2n}\right)\cdot\nabla\varphi(x)\,dx\\
&=&\lambda \int_{\Omega}\varphi(x)\,dx-
\int_\Omega \Delta\left(\frac{\lambda\,|x|^2}{2n}\right)\varphi(x)\,dx\\&=&
\lambda \int_{\Omega}\varphi(x)\,dx-\lambda
\int_\Omega\varphi(x)\,dx\\&=&0,
\end{eqnarray*}
for every~$\varphi\in C^\infty_0(\Omega)$.

This and Weyl's Lemma (recall Lemma~\ref{WEYL}) give that~$\widetilde u$ is actually~$C^2(\Omega)$
and harmonic. As a result~$u_\star\in C^2(\Omega)$ and solves~$\Delta u_\star=-\lambda$ in~$\Omega$.
As a matter of fact, we have that~$u_\star\in C^2(\overline{\Omega})$, due to Theorem~\ref{Theorem6.14GT}, and thus~$u_\star=0$ along~$\partial\Omega$.

We also remark that, for every~$v\in H^1_0(\Omega)$, with~$v\not\equiv0$, setting~$w:=\frac{v}{\int_\Omega v(x)\,dx}$,
we have that~$\int_\Omega w(x)\,dx=1$ and thus
$$
\frac{\int_\Omega |\nabla v(x)|^2\,dx}{\left(\int_\Omega v(x)\,dx\right)^2}=\int_\Omega |\nabla w(x)|^2\,dx
=J(w)\ge J(u_\star)=\frac{\int_\Omega |\nabla u_\star(x)|^2\,dx}{\left(\int_\Omega u_\star(x)\,dx\right)^2}.$$
Consequently,
by the definition of torsional rigidity in~\eqref{TAUP}
and the first Green's Identity in~\eqref{GRr1}, \begin{equation}\label{TAUUsU2} \tau(\Omega)=
\frac{\left(\int_\Omega u_\star(x)\,dx\right)^2}{\int_\Omega |\nabla u_\star(x)|^2\,dx}
=\frac{\int_\Omega u_\star(x)\,dx}{\lambda}.
\end{equation}
Then, defining~$u:=-\frac{u_\star}\lambda$,
we find that
\begin{equation} \label{SOL00POIT}\begin{dcases}
\Delta u=1&{\mbox{ in }}\Omega,\\
u=0&{\mbox{ on }}\partial\Omega,
\end{dcases}\end{equation}
and, by~\eqref{TAUUsU2},
\begin{equation}\label{TAUUsU} \tau(\Omega)=-
\int_\Omega u(x)\,dx.
\end{equation}

We now consider a family\footnote{Correspondingly\label{KMSmamx15meuntTTYrumi2}
to the alternative proof of Corollary~\ref{CO-D-SC}
presented on page~\pageref{KMSmamx15meuntTTYrumi},
we mention that it is also possible in the proof of
Lemma~\ref{HADA}
to bypass the construction of the diffeorphisms of
Lemma~\ref{FLUSSETTO} by an appropriate scaling argument.
Namely, one can instead
consider any
vector field~$\vartheta\in C^1(\R^n,\R^n)$ and the
corresponding solution~$\Phi^t$
of the associated Cauchy problem,
and then look at the competitor set~$
\left( \frac{|\Omega|}{|\Phi^t(\Omega)|}\right)^{\frac1n}\,\Phi^t(\Omega)$.
Notice that this set comes from a dilation
of the flow of~$\Omega$ with velocity~$\vartheta$,
hence it can be written as~$\widetilde\Phi^t(\Omega)$,
with~$\widetilde\Phi^t(x):=\left( \frac{|\Omega|}{|\Phi^t(\Omega)|}\right)^{\frac1n}\Phi^t(x)$.
We also point out that
\begin{eqnarray*}\widetilde\Phi^t(x)&=&
\left( \frac{|\Omega|}{| \Omega|+t\int_\Omega\div\vartheta(x)\,dx}\right)^{\frac1n}
\big(x+t \vartheta(x)\big)+o(t)\\&=&
\left( 1-\frac{t}{n\,|\Omega|}\int_\Omega\div\vartheta(x)\,dx \right)
\big(x+t \vartheta(x)\big)+o(t)\\&=&
x+t\widetilde\vartheta(x)+o(t),
\end{eqnarray*}
with
$$\widetilde\vartheta(x):=\vartheta(x)-\frac{x}{n\,|\Omega|}\int_\Omega\div\vartheta(x)\,dx,$$
and we observe that~$\widetilde\vartheta$ has zero-flux
along~$\partial\Omega$.
Then, the computations presented in these pages
remain almost unchanged, just by replacing~$\vartheta$
with~$\widetilde\vartheta$.}
of smooth diffeorphisms~$\Phi$
and a smooth
vector field~$\vartheta$
as in 
Lemma~\ref{FLUSSETTO}.
Notice that
for every~$x\in\R^n$, we have that~$\Phi(x,t)=x+t\vartheta(x)+o(t)$ as~$t\to0$,
and~$|\Phi(\Omega,t)|=|\Omega|$.
We remark that~$D_x \Phi(x,t)={\rm Id}+
tD_x\vartheta(x)+o(t)$, where~${\rm Id}$ is the identity matrix,
and consequently~$\det D_x \Phi(x,t)=1+t \div\vartheta(x)+o(t)$.

We also consider~$u^{(t)}$ to be the corresponding solution of~\eqref{SOL00POIT}
when~$\Omega$ is replaced by~$\Phi (\Omega,t)$, namely
\begin{equation*}\begin{dcases}
\Delta u^{(t)}=1&{\mbox{ in }}\Phi (\Omega,t),\\
u^{(t)}=0&{\mbox{ on }}\partial(\Phi(\Omega,t)).
\end{dcases}\end{equation*}
{F}rom the above setting and the regularity theory for the Poisson equation (recall 
Theorem~\ref{Theorem6.14GT}),
we have that, for all~$t\in[-1,1]$,
\begin{equation}\label{REGIUNSD}
\big\|u^{(t)}\big\|_{C^{2,\alpha}(\Phi(\Omega,t))}\le C,
\end{equation}
for some~$C>0$ and~$\alpha\in(0,1)$
independent of~$t$.
Also, by~\eqref{TAUUsU}, we know that
\begin{equation*} \tau(\Phi(\Omega,t))=-
\int_{\Phi (\Omega,t)} u^{(t)}(x)\,dx,
\end{equation*}
whence if~$\Omega$ is a maximizing set for a given volume constraint it follows that
\begin{equation}\label{OCHCJ-0}
\int_{\Phi (\Omega,t)} u^{(t)}(x)\,dx-\int_{\Omega} u(x)\,dx
=\tau(\Omega)
-\tau(\Phi (\Omega,t))=
o(t).
\end{equation}

Now, we point out that
\begin{equation}\label{Idhige-101}
\begin{split}
\int_{\Phi (\Omega,t)} u^{(t)}(y)\,dy\,&=
\int_{\Omega} u^{(t)}(\Phi (x,t))\,(1+t\div\vartheta(x))\,dx+o(t)\\&=
\int_{\Omega} u^{(t)}(\Phi (x,t))\,\Delta u(x)\,dx+
t\int_{\Omega} u^{(t)}(\Phi (x,t))\, \div\vartheta(x)\,dx
+o(t)\\&=-
\int_{\Omega} \nabla\Big(u^{(t)}(\Phi (x,t))\Big)\cdot\nabla u(x)\,dx+
t I_1
+o(t),
\end{split}\end{equation}
where
$$ I_1:=\int_{\Omega} u^{(t)}(\Phi (x,t))\, \div\vartheta(x)\,dx.$$
Furthermore, for every~$i\in\{1,\dots,n\}$ and~$x\in\Omega$,
\begin{equation*}
\partial_i\Big(u^{(t)}(\Phi (x,t))\Big)=\nabla u^{(t)}(\Phi (x,t))\cdot\partial_i
\Phi (x,t)=
\partial_iu^{(t)}(\Phi(x,t))+
t\nabla u^{(t)}(\Phi (x,t))\cdot\partial_i\vartheta(x)+o(t)
\end{equation*}
and
\begin{equation}\label{9i6ytdgfffshP}
\begin{split}
\Delta\Big(u^{(t)}(\Phi(x,t))\Big)\,&=
\sum_{i=1}^n\partial_i\nabla
u^{(t)}(\Phi (x,t))\cdot\partial_i\Phi (x,t)
+t\sum_{i=1}^n\big(D^2 u^{(t)}(\Phi (x,t))\partial_i\Phi (x,t)\big)\cdot\partial_i\vartheta(x)\\&\qquad
+t\nabla u^{(t)}(\Phi (x,t))\cdot\Delta\vartheta(x)
+o(t)\\&=
1+2t\sum_{i=1}^n\partial_i\nabla
u^{(t)}(\Phi (x,t))\cdot\partial_i\vartheta(x)
+t\nabla u^{(t)}(\Phi (x,t))\cdot\Delta\vartheta(x)
+o(t)
.\end{split}\end{equation}
As a result,
\begin{equation}\label{Idhige-102}
\begin{split}
\int_{\Omega} u(x)\,dx\,&=
\int_{\Omega} u(x)\,\Delta\Big(u^{(t)}(\Phi (x,t))\Big)\,dx
-2tI_2-tI_3+o(t)\\&=-
\int_{\Omega} \nabla u(x)\cdot\nabla\Big(u^{(t)}(\Phi (x,t))\Big)\,dx
-2tI_2-tI_3+o(t),\end{split}
\end{equation}
where
\begin{eqnarray*}
&&I_2:=\sum_{i=1}^n\int_{\Omega} u(x)\,
\partial_i\nabla
u^{(t)}(\Phi (x,t))\cdot\partial_i\vartheta(x)
\,dx
\\{\mbox{and }}&&I_3:=\int_{\Omega} u(x)\,\nabla u^{(t)}(\Phi (x,t))\cdot\Delta\vartheta(x)\,dx.
\end{eqnarray*}
Combining~\eqref{Idhige-101} and~\eqref{Idhige-102},
after simplifying one term we find that
$$ \int_{\Phi (\Omega,t)} u^{(t)}(x)\,dx-\int_{\Omega} u(x)\,dx
=t(I_1+2I_2+I_3)+o(t).$$
This and~\eqref{OCHCJ-0} lead to
\begin{equation}\label{9uyfdx09uygP-ds}
I_1+2I_2+I_3=
o(1)\qquad{\mbox{ as~$t\to0$.}}
\end{equation}

Now we consider the function~$w^{(t)}(x):=
u^{(t)}(\Phi (x,t))-u(x)$
and we claim that, as~$t\to0$,
\begin{equation}\label{OZzerbv}
{\mbox{$w^{(t)}\to0$ in~$W^{2,p}(\Omega)$,}}
\end{equation}
for every~$p>1$.
To prove this,
we utilize~\eqref{SOL00POIT}, \eqref{REGIUNSD}
and~\eqref{9i6ytdgfffshP} to see that~$|\Delta w^{(t)}|\le \bar{C}t$, for some~$\bar{C}>0$ independent of~$t$.
Hence, by 
the Calder\'{o}n-Zygmund regularity theory for the Poisson equation
(see~\eqref{JOSLNGEITHSHLPESTGLOUP-jofvnMTTOO}
and footnote~\ref{JOSLNGEITHSHLPESTGLOUP-jofvnMTTOO-NOT}
on page~\pageref{JOSLNGEITHSHLPESTGLOUP-jofvnMTTOO-NOT}),
we obtain~\eqref{OZzerbv}, as desired.

By using~\eqref{9uyfdx09uygP-ds} and~\eqref{OZzerbv}, we thus conclude that
\begin{equation}\label{MSka789okvJ90odfkna0pfk09iuhdfv-s9qd}
\begin{split}
0\,&=\,\lim_{t\to0}I_1+2I_2+I_3\\&=\,
\int_{\Omega} \left( u(x)\, \div\vartheta(x)
+2u(x)\sum_{i=1}^n
\partial_i\nabla
u(x)\cdot\partial_i\vartheta(x)
+u(x)\,\nabla u(x)\cdot\Delta\vartheta(x)\right)\,dx.\end{split}
\end{equation}
We also remark that
\begin{eqnarray*}
\div\big( u (D\vartheta\nabla u)\big)&=&
\sum_{i=1}^n\partial_i\big( u\,\partial_i\vartheta\cdot\nabla u
\big)\\&=&(D\vartheta\nabla u)\cdot\nabla u+u\nabla u\cdot\Delta\vartheta+u
\partial_i\nabla u\cdot\partial_i\vartheta.
\end{eqnarray*}
Comparing with~\eqref{MSka789okvJ90odfkna0pfk09iuhdfv-s9qd}
and using the Divergence Theorem, we thereby conclude that
\begin{equation}\label{KSMD0193expmusjvusirfmgvfuj}
\begin{split}
0\,&=
\int_{\Omega} \left( u(x)\, \div\vartheta(x)
+u(x)\sum_{i=1}^n
\partial_i\nabla
u(x)\cdot\partial_i\vartheta(x)\right.\\&\qquad\qquad\left.{\phantom{\sum_{i=1}^n}}
+\div\Big( u(x) \big(D\vartheta(x)\nabla u(x)\big)\Big)
-\big(D\vartheta(x)\nabla u(x)\big)\cdot\nabla u(x)
\right)\,dx\\
&=\int_{\Omega} \left( u(x)\, \div\vartheta(x)
+u(x)\sum_{i=1}^n
\partial_i\nabla
u(x)\cdot\partial_i\vartheta(x)
-\big(D\vartheta(x)\nabla u(x)\big)\cdot\nabla u(x)
\right)\,dx
.\end{split}
\end{equation}
Additionally,
$$ \div(u\vartheta)=\nabla u\cdot\vartheta+u\div\vartheta.$$
Using this information and~\eqref{KSMD0193expmusjvusirfmgvfuj},
and exploiting once more the Divergence Theorem, we see that
\begin{equation}\label{89idjf-9253it04j-05697OMS}
\begin{split}
0\,&=\int_{\Omega} \left( \div\big(u(x)\,\vartheta(x)\big)-
\nabla u(x)\cdot\vartheta(x){\phantom{\sum_{i=1}^n}}\right.\\&\qquad
\qquad\left.
+u(x)\sum_{i=1}^n
\partial_i\nabla
u(x)\cdot\partial_i\vartheta(x)
-\big(D\vartheta(x)\nabla u(x)\big)\cdot\nabla u(x)
\right)\,dx\\&=\int_{\Omega} \left( -\nabla u(x)\cdot\vartheta(x)
+u(x)\sum_{i=1}^n
\partial_i\nabla
u(x)\cdot\partial_i\vartheta(x)
-\big(D\vartheta(x)\nabla u(x)\big)\cdot\nabla u(x)
\right)\,dx
.\end{split}
\end{equation}
Besides, recalling~\eqref{SOL00POIT},
we observe that
\begin{eqnarray*}
\div\big( (\nabla u\cdot \vartheta)\nabla u\big)&=&
\sum_{i=1}^n\partial_i(\nabla u\cdot \vartheta)\partial_iu+(\nabla u\cdot \vartheta)\Delta u
\\&=&\sum_{i=1}^n\partial_i\nabla u\cdot \vartheta\partial_iu+
(D\vartheta \nabla u)\cdot\nabla u+\nabla u\cdot \vartheta,
\end{eqnarray*}
and also
\begin{eqnarray*}
\sum_{i=1}^n\partial_i\nabla u\cdot \vartheta\partial_iu+u\sum_{i=1}^n
\partial_i\nabla u\cdot \partial_i\vartheta 
&=&\sum_{i=1}^n\partial_i\nabla u\cdot \vartheta\partial_iu+\sum_{i=1}^n
\partial_i\nabla u\cdot \partial_i\vartheta u
-\nabla \Delta u\cdot \vartheta u
\\&=&
\sum_{i=1}^n \partial_i\big(\partial_i\nabla u\cdot \vartheta u\big)
\\&=&\div V,
\end{eqnarray*}
where~$V$ is the vector field defined as~$V\cdot e_i:=
\partial_i\nabla u\cdot \vartheta u$, for each~$i\in\{1,\dots,n\}$.
We also remark that~$V=0$ along~$\partial\Omega$, thanks to
the boundary condition in~\eqref{SOL00POIT}.
Consequently, we can rewrite~\eqref{89idjf-9253it04j-05697OMS}
in the form
\begin{equation*}
\begin{split}
0\,&=\int_{\Omega} \left( \sum_{i=1}^n\partial_i\nabla u(x)
\cdot \vartheta(x)\partial_iu(x)
-\div\Big( \big(\nabla u(x)\cdot \vartheta(x)\big)\nabla u(x)\Big)
+u(x)\sum_{i=1}^n
\partial_i\nabla
u(x)\cdot\partial_i\vartheta(x)
\right)\,dx\\&=\int_{\Omega} \left( \div V(x)
-\div\Big( \big(\nabla u(x)\cdot \vartheta(x)\big)\nabla u(x)\Big)
\right)\,dx\\&=\int_{\partial\Omega}\big(\nabla u(x)\cdot \vartheta(x)\big)\big(
\nabla u(x)\cdot\nu(x)\big)\,d{\mathcal{H}}^{n-1}_x
.\end{split}
\end{equation*}
Since~$u=0$ along~$\partial\Omega$,
hence~$\nabla u=|\nabla u|\,\nu$
along~$\partial\Omega$, we thereby infer that
\begin{equation*}
\int_{\partial\Omega}\big( \partial_\nu u(x)\big)^2\,
\vartheta(x)\cdot\nu(x)\,d{\mathcal{H}}^{n-1}_x=0.\end{equation*}
This and Lemma~\ref{3349056} gives that\footnote{It is interesting to point out that
Lemma~\ref{HADA} also entails a rather explicit motivation
for Serrin's overdetermined problem~\eqref{HESE}
in terms of optimal heating.
Indeed, suppose that a region~$\Omega$ is given.
The external environment lies at constant, say zero, temperature,
and~$\Omega$ is heated in a uniform and homogeneous way.
Up to normalizing constants, the equilibrium configuration is therefore described by the equation
\begin{equation*} \begin{dcases}
-\Delta U_\Omega=C&{\mbox{ in }}\Omega,\\
U_\Omega=0&{\mbox{ on }}\partial\Omega,
\end{dcases}\end{equation*}
where~$U_\Omega$ represents the temperature in~$\Omega$
and~$C\in(0,+\infty)$ is the heat source term described
by the uniform heating of~$\Omega$.
In this setting, a natural question is to determine the domain~$\Omega$ which maximizes
(among the smooth, bounded, connected sets of a given volume)
the average temperature
$$ T(\Omega):=\int_\Omega U_\Omega(x)\,dx.$$
Setting~$u_\Omega:=-\frac{U_\Omega}{C}$ we have that
\begin{equation*} \begin{dcases}
\Delta u_\Omega=1&{\mbox{ in }}\Omega,\\
u_\Omega=0&{\mbox{ on }}\partial\Omega.
\end{dcases}\end{equation*}
Moreover, in view of~\eqref{TAUUsU},
$$ T(\Omega)=-C\int_\Omega u_\Omega(x)\,dx
=C\tau(\Omega).$$
Therefore, the determination of the optimal shape for~$\Omega$
aiming at maximizing the average temperature is equivalent
to the Saint Venant problem (in particular, the optimally heated domain is a ball, thanks to
Corollary~\ref{SVP}). Furthermore, the optimal temperature~$U_\Omega$
for such a domain satisfies Serrin's problem~\eqref{HESE},
since~$\partial_\nu U_\Omega=-C\partial_\nu u_\Omega$ is constant along~$\partial\Omega$,
owing to~\eqref{ONG}.} \begin{equation}\label{ONG}
{\mbox{$\partial_\nu u$
is constant along~$\partial\Omega$.}}\end{equation}
{F}rom this fact and~\eqref{SOL00POIT} the desired result follows.
\end{proof}

Combining Lemma~\ref{HADA} with
Theorem~\ref{SERR}, we obtain the desired classification result
for the Saint Venant problem:\index{Saint Venant problem}

\begin{corollary}\label{SVP} Let~$\alpha\in(0,1)$.
Let~$\Omega\subseteq\R^n$ be a bounded and connected
set of class~$C^{2,\alpha}$
that maximizes
the torsional rigidity~$\tau$ among the
bounded, open and connected sets
of class~$C^{2,\alpha}$
with a given
volume. Then, $\Omega$ is necessarily a ball.\end{corollary}

Other proofs of Corollary~\ref{SVP}
can be obtained from~\eqref{TAUUsU} and rearrangements methods (see~\cite{MR0043486, MR601601}
and also~\cite{MR3412286} and the references therein).\index{Serrin's overdetermined problem|)}

\section{The converse of the Mean Value Formula}\label{KURAN-SEC}

We address here a natural question arising
from the Mean Value Formula\index{Mean Value Formula}
for harmonic functions 
(recall Theorem~\ref{KAHAR}). Namely, suppose that
we are given a set~$\Omega$ and a point~$x_0\in\Omega$
and we know that the value at~$x_0$ of
every harmonic function in~$\Omega$
coincides with the volume average of the function in~$\Omega$ (or, alternatively,
with the surface average along~$\partial\Omega$): then, what can we say
about the set~$\Omega$?

Remarkably, the answer is that~$\Omega$ is necessarily
a ball, and the point~$x_0$ is necessarily the center of such ball:
as a result, balls are the only sets for which
the Mean Value Formula
in Theorem~\ref{KAHAR} holds true.

This beautiful result (obtained in~\cite{MR320348, MR1021402})
follows from the following general statement (see in particular~(i), (ii) and~(v) below):

\begin{theorem}
Let~$n\ge2$ and~$\Omega$ be a bounded and open set with boundary of class~$C^1$. Let also~$x_0\in \Omega$
and~$G$ be the Green Function of~$\Omega$, as defined in~\eqref{FLGREEN}. Then,
the following conditions are equivalent:\index{converse of the Mean Value Formula}
\begin{itemize}
\item[(i).] For every harmonic function~$u\in L^1(\Omega)$, we have that
$$ u(x_0)=\fint_\Omega u(x)\,dx.$$
\item[(ii).] For every harmonic function~$u\in C(\overline\Omega)$, we have that
$$ u(x_0)=\fint_{\partial\Omega} u(x)\,d{\mathcal{H}}^{n-1}_x.$$
\item[(iii).] There exist~$c_1\in\R\setminus\{0\}$, $c_2\in\R$ and~$w\in C^2(\Omega)\cap C^1(\overline\Omega)$ such that
$$ \begin{dcases}
\Delta w(x)= c_1 & {\mbox{ for all }}x\in\Omega,\\w(x)=0
& {\mbox{ for all }}x\in\partial\Omega,\\
\nabla w(y)\cdot\nu(y)=c_2\nabla_y G(x_0,y)\cdot\nu(y)& {\mbox{ for all }}x\in\partial\Omega.
\end{dcases}$$
\item[(iv).] The function~$\partial \Omega\ni y\mapsto \nabla_y G(x_0,y)\cdot\nu(y)$ is constant.
\item[(v).] $\Omega$ is a ball centered at~$x_0$.
\end{itemize}
\end{theorem}

\begin{proof} We show that~(iv) implies~(ii), which implies~(v), which implies~(iv), and that~(iii) implies~(i),
which implies~(v),
which implies~(iii).

For this, let us assume that~(iv) holds true, namely~$\nabla_y G(x_0,y)\cdot\nu(y)=c$
for every~$y\in\partial\Omega$, for some~$c\in\R$, and let~$u$ be harmonic in~$\Omega$ and continuous in~$\overline\Omega$.
Let also~$\Omega'\Subset\Omega$, with~$\Omega'\ni x_0$.
By Theorem~\ref{MS:SKMD344567yjghS} (and recalling also Lemma~\ref{WEYL}),
\begin{eqnarray*}&&-u(x_0)=
\int_{\Omega'} \Delta u(y)G(x_0,y)\,dy+
\int_{\partial\Omega'} u(y)\frac{\partial G}{\partial\nu}(x_0,y)\,d{\mathcal{H}}^{n-1}_y=\int_{\partial\Omega'} u(y)\frac{\partial G}{\partial\nu}(x_0,y)\,d{\mathcal{H}}^{n-1}_y.
\end{eqnarray*}
Hence, taking~$\Omega'$ as close as we wish to~$\Omega$,
\begin{equation}\label{0987-09876k09876-cx-pkm}
u(x_0)=-
\int_{\partial\Omega} u(y)\frac{\partial G}{\partial\nu}(x_0,y)\,d{\mathcal{H}}^{n-1}_y=
-c\int_{\partial{\Omega}} u(y) \,d{\mathcal{H}}^{n-1}_y.
\end{equation}
In particular, choosing~$u:=1$, we find that~$1=-c{\mathcal{H}}^{n-1}(\partial\Omega)$.
Plugging this information back into~\eqref{0987-09876k09876-cx-pkm}, we obtain~(ii),
as desired.

Let us now suppose that~(ii) holds true and let~$u$ be harmonic in~$\Omega$
and such that~$u(y)=\frac1{{\mathcal{H}}^{n-1}(\partial\Omega)}+\nabla_yG(x_0,y)\cdot\nu(y)$ for every~$y\in\partial\Omega$
(notice that such a function exists thanks to Theorem~\ref{PERRO-2}
since the boundary of class~$C^1$ provides the required external cone).
By~(ii) and Theorem~\ref{MS:SKMD344567yjghS}, we have that, for every~$\Omega'\Subset\Omega$,
with~$\Omega'\ni x_0$,
\begin{equation*}\begin{split}&
0=u(x_0)+
\int_{\partial\Omega'}u(y)\nabla_yG(x_0,y)\cdot\nu(y)\,d{\mathcal{H}}^{n-1}_y\\
&\qquad\qquad= \fint_{\partial\Omega'} u(y)\,dy+
\int_{\partial\Omega'}u(y)\nabla_yG(x_0,y)\cdot\nu(y)\,d{\mathcal{H}}^{n-1}_y\\
&\qquad\qquad=
\int_{\partial\Omega'}u(y)\left( \frac1{{\mathcal{H}}^{n-1}(\partial\Omega')}+\nabla_yG(x_0,y)\cdot\nu(y)\right)\,d{\mathcal{H}}^{n-1}_y
.\end{split}
\end{equation*}
Therefore, taking~$\Omega'$ as close as we wish to~$\Omega$,
\begin{eqnarray*} 0&=&
\int_{\partial\Omega}u(y)\left( \frac1{{\mathcal{H}}^{n-1}(\partial\Omega')}+\nabla_yG(x_0,y)\cdot\nu(y)\right)\,d{\mathcal{H}}^{n-1}_y\\&=&\int_{\partial\Omega}
\left( \frac1{{\mathcal{H}}^{n-1}(\partial\Omega')}+\nabla_yG(x_0,y)\cdot\nu(y)\right)^2\,d{\mathcal{H}}^{n-1}_y
,\end{eqnarray*}
leading to
\begin{equation}\label{KA:WEIND}
\nabla_yG(x_0,y)\cdot\nu(y)=-\frac1{{\mathcal{H}}^{n-1}(\partial\Omega)}\quad{\mbox{ for each }}y\in\partial\Omega.
\end{equation}
%%  In particular, by the Divergence Theorem,
%%  \begin{equation*}
%%  \begin{split}&
%%  \int_{\partial\Omega} (y-x_0)\cdot\nu(y)\,
%%  \nabla_yG(x_0,y)\cdot\nu(y)\,d{\mathcal{H}}^{n-1}_y
%%  =-\frac1{{\mathcal{H}}^{n-1}(\partial\Omega)}
%%  \,\int_{\partial\Omega} (y-x_0)\cdot\nu(y)\,d{\mathcal{H}}^{n-1}_y\\&\qquad
%%  =-\frac1{{\mathcal{H}}^{n-1}(\partial\Omega)}\,\int_\Omega \div_y(y-x_0)\,dy
%%       =-\frac{n|\Omega|}{{\mathcal{H}}^{n-1}(\partial\Omega)}.
%%  \end{split}\end{equation*}
Moreover, since~$G(x_0,\cdot)=0$ along~$\partial\Omega$, we know that~$\nabla_y G(x_0,y)=\pm|\nabla_y G(x_0,y)|\nu(y)$
for all~$y\in\partial\Omega$, therefore we infer from~\eqref{KA:WEIND} that~$\frac1{{\mathcal{H}}^{n-1}(\partial\Omega)}=
|\nabla_yG(x_0,y)|$ for all~$y\in\partial\Omega$. For this reason, we have that
\begin{equation}\label{NONSDPODIA}
\begin{split}&
\left|\int_{\partial\Omega} (y-x_0)\cdot\nabla_yG(x_0,y)\,d{\mathcal{H}}^{n-1}_y\right|=
\left|\int_{\partial\Omega}|\nabla_yG(x_0,y)|\, (y-x_0)\cdot  \nu(y)\,d{\mathcal{H}}^{n-1}_y\right|
\\&\qquad\qquad=
\frac1{{\mathcal{H}}^{n-1}(\partial\Omega)}
\left|\int_{\partial\Omega}(y-x_0)\cdot  \nu(y)\,d{\mathcal{H}}^{n-1}_y\right|
\\&\qquad\qquad=
\frac1{{\mathcal{H}}^{n-1}(\partial\Omega)}
\left|\int_{\Omega}{\rm{div}} (y-x_0)\,dy\right|
=\frac{n|\Omega|}{{\mathcal{H}}^{n-1}(\partial\Omega)}.
\end{split}\end{equation}
We now define
$$ h(y):=(y-x_0)\cdot\nabla_y G(x_0,y)+(n-2)G(x_0,y).$$
By~\eqref{GAMMAFU},
\begin{eqnarray*}
(y-x_0)\cdot\nabla_y \Gamma(x_0-y)&=&
\begin{dcases}
-\frac{c_n(n-2)}{|x_0-y|^{n-2}}&{\mbox{ if }}n\ne2,\\
- c_n&{\mbox{ if }}n=2,
\end{dcases}\\&=&
\begin{dcases}
-(n-2)\Gamma(x_0-y)&{\mbox{ if }}n\ne2,\\
-c_n&{\mbox{ if }}n=2
.\end{dcases}
\end{eqnarray*}
Thus, recalling the Robin Function in~\eqref{L:S:ASM S} and~\eqref{FLGREEN}, as well as the value of~$c_n$ when~$n=2$
given in~\eqref{COIENEE}, we see that, for every~$y\in\Omega$,
\begin{eqnarray*}
h(y)&=& (y-x_0)\cdot\nabla_y \Gamma(x_0-y)+(n-2)\Gamma(x_0-y)
-(y-x_0)\cdot\nabla_y \Psi^{(x_0)}(y)-(n-2)\Psi^{(x_0)}(y)\\
&=& c_\star-(y-x_0)\cdot\nabla_y \Psi^{(x_0)}(y)-(n-2)\Psi^{(x_0)}(y),
\end{eqnarray*}
with
$$c_\star=\begin{dcases}
\displaystyle-\frac{1}{2\pi} & {\mbox{ if }}n=2,\\
0&{\mbox{ otherwise }},
\end{dcases}$$ which in particular
gives that~$h$ is harmonic in~$\Omega$ since so is the Robin Function.

Consequently, we can apply~(ii) to the function~$h$ and infer that
\begin{equation*}\begin{split}&c_\star-(n-2)\Psi^{(x_0)}(x_0)
=
h(x_0)=\fint_{\partial\Omega} h(y)\,d{\mathcal{H}}^{n-1}_y\\&\qquad=
\fint_{\partial\Omega} \Big((y-x_0)\cdot\nabla_y G(x_0,y)+(n-2)G(x_0,y)\Big)\,d{\mathcal{H}}^{n-1}_y\\&\qquad
=\fint_{\partial\Omega} (y-x_0)\cdot\nabla_y G(x_0,y)\,d{\mathcal{H}}^{n-1}_y.
\end{split}\end{equation*}
This and~\eqref{NONSDPODIA} lead to
\begin{equation}\label{DIS890T00IUD097jUEXNC6CfA}\begin{split}&\big| c_\star-(n-2)\Psi^{(x_0)}(x_0)\big|
=\frac{n|\Omega|}{\big({\mathcal{H}}^{n-1}(\partial\Omega)\big)^2}.
\end{split}\end{equation}
We now distinguish two cases: when~$n=2$, equation~\eqref{DIS890T00IUD097jUEXNC6CfA}
boils down to
$$ \frac1{2\pi}=
\frac{2|\Omega|}{\big({\mathcal{H}}^{1}(\partial\Omega)\big)^2}$$
that is~$\big({\mathcal{H}}^{1}(\partial\Omega)\big)^2=4\pi|\Omega|$.
Hence, the set~$\Omega$ is a minimizer for the isoperimetric problem in the plane (see Corollary~\ref{CO-D})
and accordingly it is necessarily a ball.

If instead~$n\ne2$, we exploit the lower bound on the Robin Function pointed out in Lemma~\ref{WEI-SHA}
and we deduce from~\eqref{DIS890T00IUD097jUEXNC6CfA} that
$$
\frac{ 1}{n|B_1|^{\frac2n}|\Omega|^{\frac{n-2}{n}} }
\le\big| (n-2)\Psi^{(x_0)}(x_0)\big|=
\frac{n|\Omega|}{\big({\mathcal{H}}^{n-1}(\partial\Omega)\big)^2}
$$
or equivalently, in light of~\eqref{B1},
$$
\frac{ |B_1|^{n-1}}{\big({\mathcal{H}}^{n-1}(\partial B_1)\big)^{n}}\le
\frac{ |\Omega|^{n-1}}{\big({\mathcal{H}}^{n-1}(\partial\Omega)\big)^{n}}.$$
Hence, also in this case~$\Omega$ is a minimizer for the isoperimetric problem in the plane (see Corollary~\ref{CO-D})
and accordingly it is necessarily a ball.

In view of these observations, we can write that~$\Omega=B_R(p)$ for some~$p\in\R^n$ and~$R>0$.
Plugging this information back into~(ii), given~$\varpi\in\R^n$ and choosing~$u(x):=\varpi\cdot (x-p)$, by odd symmetry we find that
$$ \varpi\cdot(x_0-p)=u(x_0)=\fint_{\partial B_R(p)} u(x)\,d{\mathcal{H}}^{n-1}_x=0,$$
from which we conclude that~$p=x_0$. The proof of~(v) is thereby complete.

If instead~(v) holds true, then~(iv) plainly follows from Theorem~\ref{MS:SKMD344567yjghS-1}.

Let us now assume that~(iii) holds true and let~$u$ be harmonic in~$\Omega$.
Then, by the Green's Identity in~\eqref{GRr2},
\begin{equation}\label{CXO-VACXNMCOCISNOVPS}
\begin{split}&
c_1\int_\Omega u(y)\,dy
=
\int_\Omega \Big( u(y)\Delta w(y)-w(y)\Delta u(y)\Big)\,dy=
\int_{\partial\Omega}\left(u(y)\frac{\partial w}{\partial\nu}(y)-
w(y)\frac{\partial u}{\partial\nu}(y)\right)\,d{\mathcal{H}}^{n-1}_y\\&\qquad=c_2
\int_{\partial\Omega} u(y)\nabla_yG(x_0,y)\cdot \nu(y)\,d{\mathcal{H}}^{n-1}_y
.\end{split}\end{equation}
Also, by Theorem~\ref{MS:SKMD344567yjghS},
$$ u(x_0)=-
\int_{\partial \Omega}u(y)\nabla_yG(x_0,y)\cdot\nu(y)\,d{\mathcal{H}}^{n-1}_y,$$
which combined with~\eqref{CXO-VACXNMCOCISNOVPS} gives that
$$ -\frac{c_2}{c_1} u(x_0)=\int_\Omega u(y)\,dy.$$
In particular, choosing~$u$ to be constant, we find that~$-\frac{c_2}{c_1}=|\Omega|$,
from which we obtain~(i).

Now, let us suppose that~(i) holds true. Up to a translation, we assume that~$x_0=0$.
Given~$i$, $j\in\{1,\dots,n\}$,
we take~$\widetilde{u}$ to be the harmonic function in~$\Omega$
such that
\begin{equation}\label{SPLBsUCA5L7L}
\widetilde{u}(x)=x_i\nu_j(x)-x_j\nu_i(x)\quad{\mbox{ for all~$x\in\partial\Omega$}}.\end{equation}
We observe that, for every~$i$, $j\in\{1,\dots,n\}$,
\begin{equation}\label{SPLBsUCA5L7LBIS}
\delta_{ij}\widetilde{u}(x)=0\quad{\mbox{ for all~$x\in\partial\Omega$}}.\end{equation}
We also define~$u(x):=x_i\partial_j\widetilde{u}(x)-x_j\partial_i\widetilde{u}(x)$.
We stress that, for each~$x\in\Omega$,
$$ \Delta u(x)=
2\nabla x_i\cdot\partial_j\nabla\widetilde{u}(x)-2\nabla x_j\cdot\partial_i\nabla\widetilde{u}(x)=2\partial_{ij}\widetilde{u}(x)-2\partial_{ij}\widetilde{u}(x)=0.
$$
In this way, by~(i),
\begin{eqnarray*}&& 0=|\Omega|\,u(0)=\int_\Omega u(x)\,dx=
\int_\Omega\Big( x_i\partial_j\widetilde{u}(x)- x_j\partial_i\widetilde{u}(x)\Big)\,dx
.\end{eqnarray*}
Since, by the Divergence Theorem and~\eqref{SPLBsUCA5L7LBIS}
$$ \int_\Omega x_i\partial_j\widetilde{u}(x)\,dx=
\int_\Omega\Big( \div\big( x_i \widetilde{u}(x)\,e_j \big)-\delta_{ij}\widetilde{u}(x)
\Big)\,dx=
\int_{\partial\Omega} x_i \,\widetilde{u}(x)\,\nu_j(x)\,d{\mathcal{H}}^{n-1}_x,
$$
we thereby find that
\begin{equation*} 0=
\int_{\partial\Omega} \Big(
x_i \,\widetilde{u}(x)\,\nu_j(x)
-x_j \,\widetilde{u}(x)\,\nu_i(x)\Big)
\,d{\mathcal{H}}^{n-1}_x=\int_{\partial\Omega}\widetilde{u}(x)\, \Big(
x_i \,\nu_j(x)
-x_j \,\nu_i(x)\Big)
\,d{\mathcal{H}}^{n-1}_x
.\end{equation*}
Hence, recalling~\eqref{SPLBsUCA5L7L},
$$ 0=\int_{\partial\Omega} \Big(
x_i \,\nu_j(x)
-x_j \,\nu_i(x)\Big)^2
\,d{\mathcal{H}}^{n-1}_x,
$$
which entails that~$x_i \,\nu_j(x)
-x_j \,\nu_i(x)=0$ for every~$x\in\partial\Omega$.

Thus, for every~$x\in\partial\Omega$,
$$ x\cdot\nu(x)\,\nu_j(x)=\sum_{i=1}^nx_i \,\nu_i(x)\nu_j(x)
=\sum_{i=1}^n x_j \,\nu_i^2(x)=x_j
$$
and, as a result, $x$ and~$\nu(x)$ are parallel, for every~$x\in\partial\Omega$.

Consequently, we write~$\nu(x)=g(x)x$, for a suitable scalar function~$g$
and we stress that~$1=|g(x)|\,|x|$, from which we see that~$g$ never vanishes.
We claim that
\begin{equation}\label{CONNESP}
{\mbox{each connected component of~$\partial\Omega$ is a sphere centered at the origin.}}
\end{equation}
To check this, let~$\Sigma$ be a connected component of~$\partial\Omega$,
pick a point~$q\in\Sigma$ and let~$R:=|q|$. Take also another point~$\widetilde{q}$
and a curve~$\gamma:[0,1]\to\Sigma$ such that~$\gamma(0)=q$ and~$\gamma(1)=\widetilde{q}$.
Using that~$\dot\gamma(t)$ is a tangent vector, we have that
\begin{eqnarray*} && |\widetilde{q}|^2 -R^2=
|\gamma(1)|^2-|\gamma(0)|^2=\int_0^1
\frac{d}{dt}|\gamma(t)|^2\,dt\\
&&\qquad=
2\int_0^1
\gamma(t)\cdot\dot\gamma(t)\,dt=2\int_0^1\frac{1}{g(\gamma(t))}\,
\nu(\gamma(t))\cdot\dot\gamma(t)\,dt=0,
\end{eqnarray*}
which proves~\eqref{CONNESP}.

{F}rom~\eqref{CONNESP} and the fact that~$0\in\Omega$, we find that either~$\Omega$ is a ball centered at the origin, or
\begin{equation}\label{ADANN} \Omega=B_{R_0}\cup\left(\bigcup_{k=1}^M  B_{R_k}\setminus \overline{B_{r_k}}\right)\end{equation}
for some~$M\in\N$ and suitable radii~$R_M>r_M>R_{M-1}>r_{m-1}\dots>R_1>r_1>R_0>0$. Thus,
to establish~(v), we need to rule out the possibility of the additional annuli in~\eqref{ADANN}. To this end, we remark
that~$\Omega$ is connected: otherwise, if~$\Omega=\Omega_1\cup\Omega_2$ for suitable open sets~$\Omega_1$ and~$\Omega_2$
whose closure have empty intersection, we can suppose that~$0\in\Omega_1$ and let~$u:=\chi_{\Omega_2}$.
In this setting, we have that~$u$ is harmonic in~$\Omega$ and thus, by~(i),
$$ 0=u(0)=\fint_\Omega u(x)\,dx=\frac{|\Omega_2|}{|\Omega|}>0.$$
This is a contradiction and we have therefore\footnote{An alternative proof that~(i) implies~(v) exploits the Poisson Kernel and goes
as follows. Up to a translation, one can suppose that~$x_0=0$.
Take~$r>0$ such that~$B_r\subseteq\Omega$ with~$p\in(\partial B_r)\cap(\partial\Omega)$.
Let
$$ u(x):=\frac{|x|^2-r^2}{|x-p|^n}+r^{2-n}.$$
Notice that~$u$ is harmonic in~$\R^n\setminus\{p\}\supseteq\Omega$. Also, if~$x\in\R^n\setminus B_r$
then~$u(x)\ge r^{2-n}$. Furthermore, we have that~$u(0)=-\frac{r^2}{|p|^n}+r^{2-n}=0$.
As a result, using~(i) and the Mean Value Formula in Theorem~\ref{KAHAR}(iii),
\begin{eqnarray*}&&
0=u(0)=\fint_\Omega u(y)\,dy=
\frac1{|\Omega|}\left( \int_{B_r} u(y)\,dy+\int_{\Omega\setminus B_r} u(y)\,dy \right)
=\frac1{|\Omega|}\left( |B_r|\,u(0)+\int_{\Omega\setminus B_r} u(y)\,dy \right)\\&&\qquad\qquad\qquad
=\frac1{|\Omega|}\int_{\Omega\setminus B_r} u(y)\,dy
\ge\frac{|\Omega\setminus B_r|}{|\Omega|}\,r^{2-n}.
\end{eqnarray*}
This gives that~$|\Omega\setminus B_r|=0$ and therefore~$\Omega=B_r$.}
established~(v), as desired.

Suppose now that~(v) holds true.
We take~$w$ to be 
the solution of
$$ \begin{dcases}
\Delta w=-1 &{\mbox{ in }}\Omega,
\\ w=0&{\mbox{ on }}\partial\Omega
\end{dcases}$$
and~$u$ be
harmonic in~$\Omega$ and such that~$u(y)=\partial_\nu w(y)-|\Omega|\nabla_yG(x_0,y)\cdot\nu(y)$
for all~$y\in\partial\Omega$ (we recall that such solutions exist thanks to Theorem~\ref{PERRO-2}).
Then, by the Mean Value Formula in Theorem~\ref{KAHAR}(iii) and the Green's Identity in~\eqref{GRr2},
\begin{equation}\label{CXO-VACX-2-NMCOCISNOVPS}
\begin{split}&
-|\Omega|\,u(x_0)=
-\int_\Omega u(y)\,dy=\int_\Omega u(y)\,\Delta w(y)\,dy
=
\int_\Omega \Big( u(y)\Delta w(y)-w(y)\Delta u(y)\Big)\,dy\\&\qquad=
\int_{\partial\Omega}\left(u(y)\frac{\partial w}{\partial\nu}(y)-
w(y)\frac{\partial u}{\partial\nu}(y)\right)\,d{\mathcal{H}}^{n-1}_y\\&\qquad=
\int_{\partial\Omega} \left(
\frac{\partial w}{\partial\nu}(y)-|\Omega|\nabla_yG(x_0,y)\cdot\nu(y)
\right)\frac{\partial w}{\partial\nu}(y)\,d{\mathcal{H}}^{n-1}_y
.\end{split}\end{equation}
Additionally, by Theorem~\ref{MS:SKMD344567yjghS},
\begin{eqnarray*} u(x_0)&=&-
\int_{\partial \Omega}u(y)\nabla_yG(x_0,y)\cdot\nu(y)\,d{\mathcal{H}}^{n-1}_y\\&=&-
\int_{\partial \Omega}\left(
\frac{\partial w}{\partial\nu}(y)-|\Omega|\nabla_yG(x_0,y)\cdot\nu(y)
\right)\nabla_yG(x_0,y)\cdot\nu(y)\,d{\mathcal{H}}^{n-1}_y.\end{eqnarray*}
Comparing this identity with~\eqref{CXO-VACX-2-NMCOCISNOVPS}, we conclude that
\begin{eqnarray*}
0&=&\int_{\partial\Omega} \left(
\frac{\partial w}{\partial\nu}(y)-|\Omega|\nabla_yG(x_0,y)\cdot\nu(y)
\right)\frac{\partial w}{\partial\nu}(y)\,d{\mathcal{H}}^{n-1}_y
+|\Omega|\,u(x_0)\\&=&\int_{\partial\Omega} \left(
\frac{\partial w}{\partial\nu}(y)-|\Omega|\nabla_yG(x_0,y)\cdot\nu(y)
\right)\frac{\partial w}{\partial\nu}(y)\,d{\mathcal{H}}^{n-1}_y\\&&\quad-
\int_{\partial \Omega}\left(
\frac{\partial w}{\partial\nu}(y)-|\Omega|\nabla_yG(x_0,y)\cdot\nu(y)
\right)|\Omega|\nabla_yG(x_0,y)\cdot\nu(y)\,d{\mathcal{H}}^{n-1}_y\\&=&\int_{\partial\Omega} \left(
\frac{\partial w}{\partial\nu}(y)-|\Omega|\nabla_yG(x_0,y)\cdot\nu(y)
\right)^2\,d{\mathcal{H}}^{n-1}_y
\end{eqnarray*}
and therefore
$$ \frac{\partial w}{\partial\nu}(y)-|\Omega|\nabla_yG(x_0,y)\cdot\nu(y)
=0\qquad{\mbox{for every }}y\in\partial\Omega.$$
This says that~$w$ solves the equation in~(iii) with~$c_1:=-1$ and~$c_2:=|\Omega|$, as desired.
\end{proof}

\chapter{The moving plane method}\label{MO:PLAN}

In this chapter we will revisit some of the previously showcased results,\index{moving plane method|(}
such as Alexandrov's Soap Bubble Theorem presented in Section~\ref{SEC:TheSoapBubbleTheorem}
and Serrin's overdetermined problem presented in Section~\ref{SEC:OVER-SERRINTheorem}:
differently from the previous exposition, we exhibit here somewhat the earliest, and possibly most ``visual'',
approach to these questions, namely the so-called moving plane method. In doing so, we will
also discuss a classical result by Basilis Gidas, Wei-Ming Ni and Louis Nirenberg~\cite{MR544879}.

In outline, the quintessence of the method consists in detecting the symmetry
with respect to a given hyperplane by comparing the original picture with the reflected one: to do so,
one begins by placing the reflection hyperplane ``very far'' from the region of interest,
so that the reflection is void, and then slowly moves the hyperplane closer and closer
to the domain under investigation, starts looking at initial reflections, which typically (in the best case scenarios) produce
an easily visualizable inclusion (thus allowing the method ``to start'') and then detects the ``critical position''
at which the picture becomes symmetric.

\begin{figure}
  \centering
  \includegraphics[width=.18\linewidth]{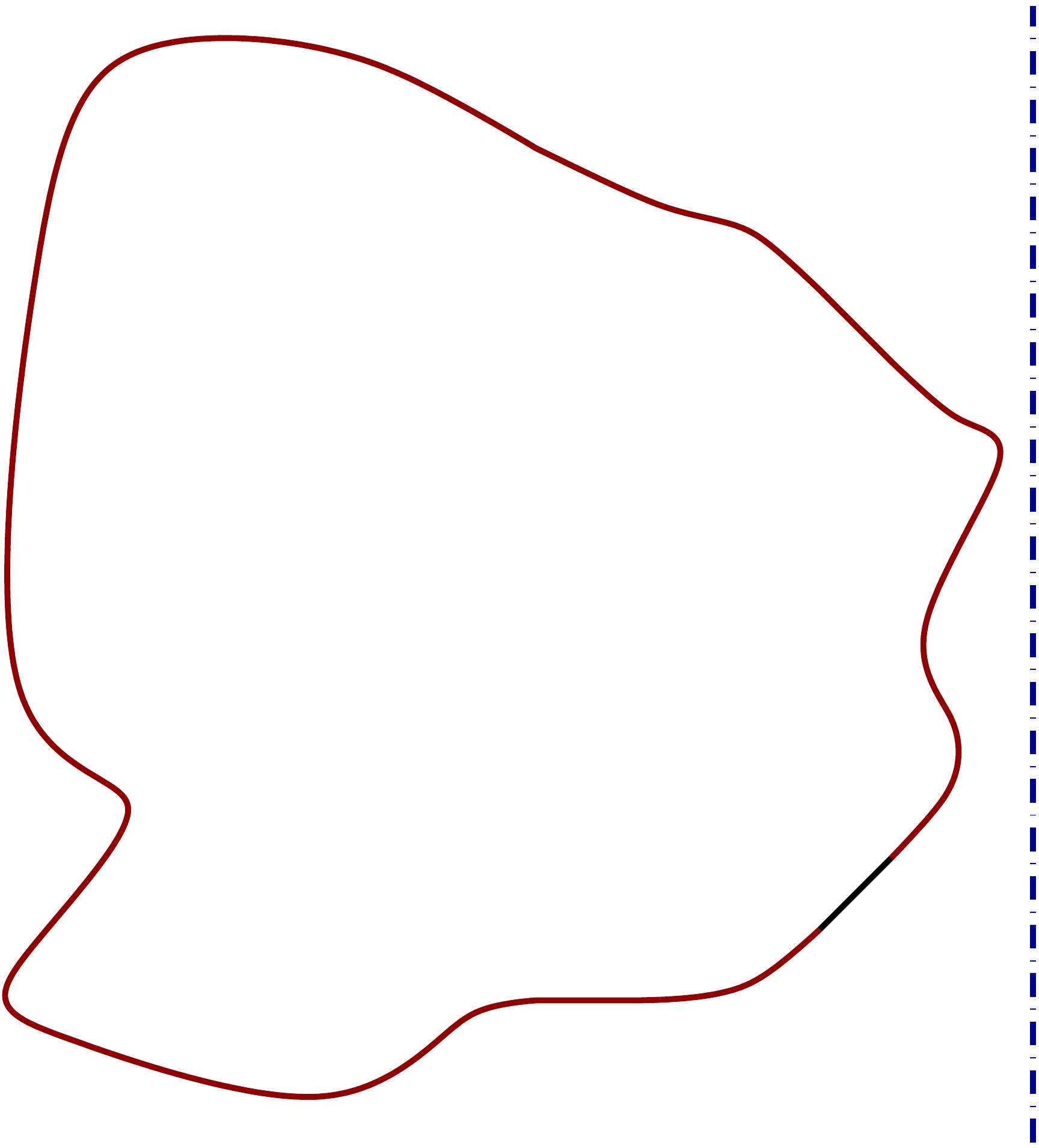}$\qquad$
    \includegraphics[width=.18\linewidth]{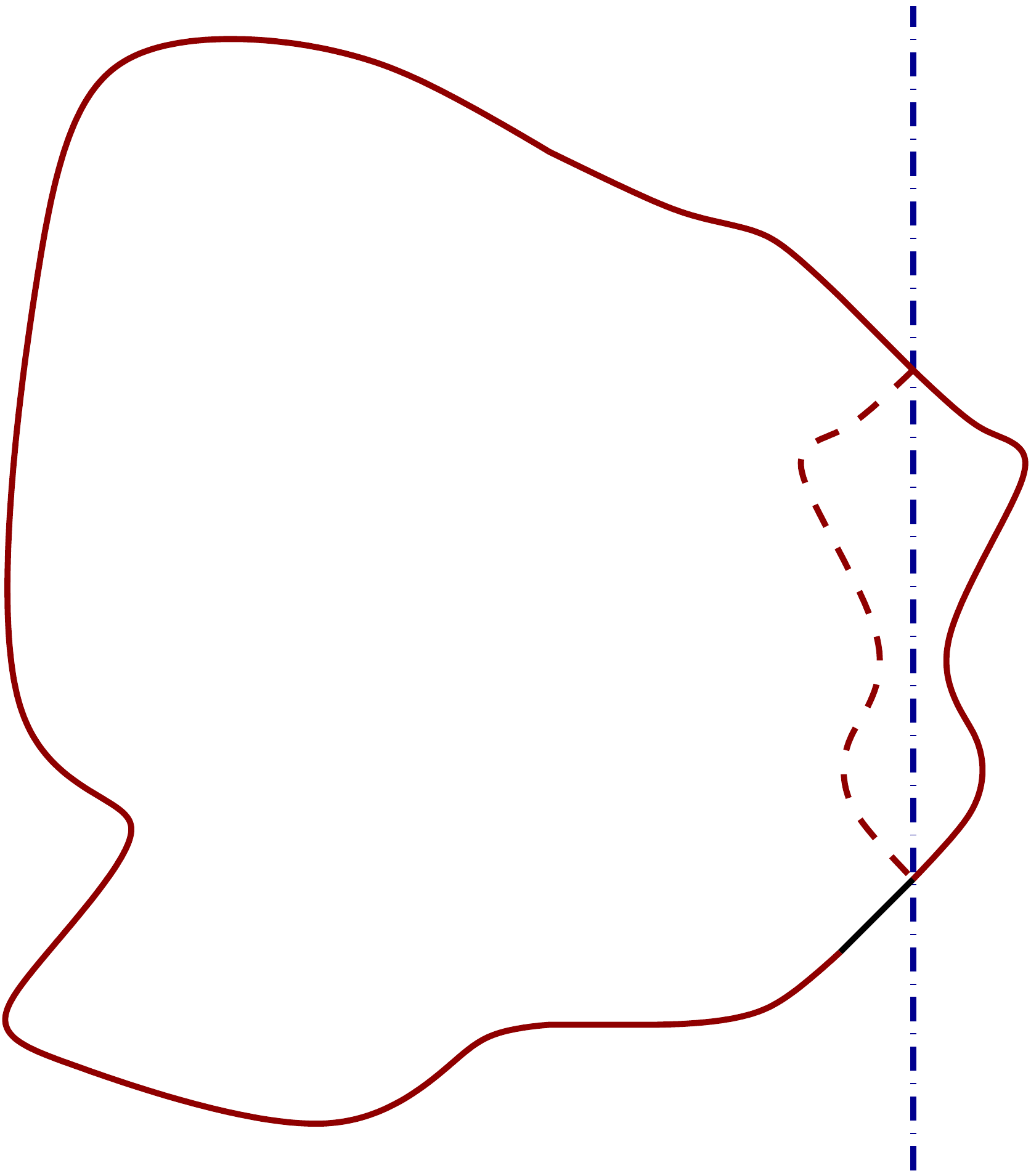}$\qquad$
      \includegraphics[width=.18\linewidth]{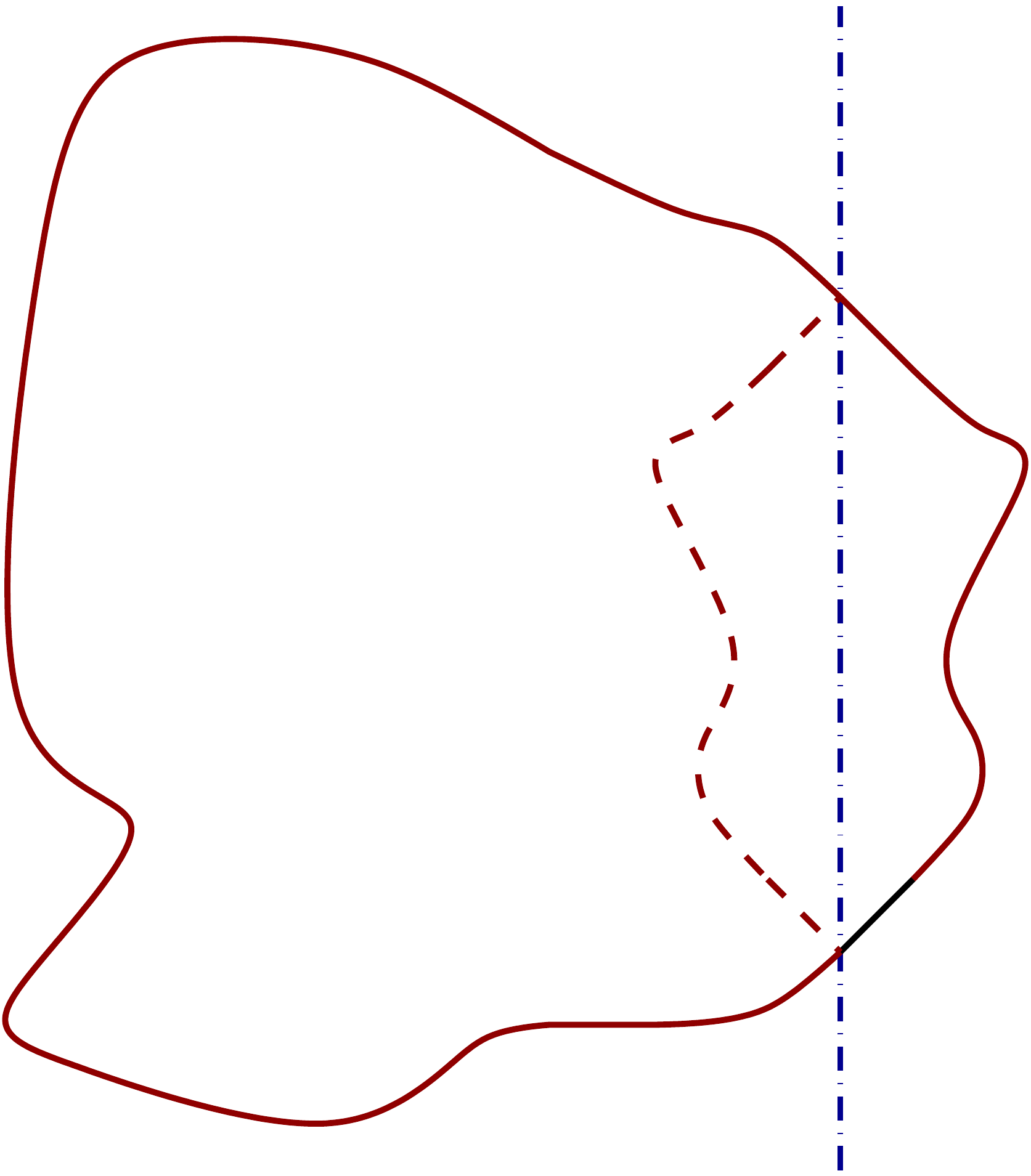}$\qquad$
        \includegraphics[width=.18\linewidth]{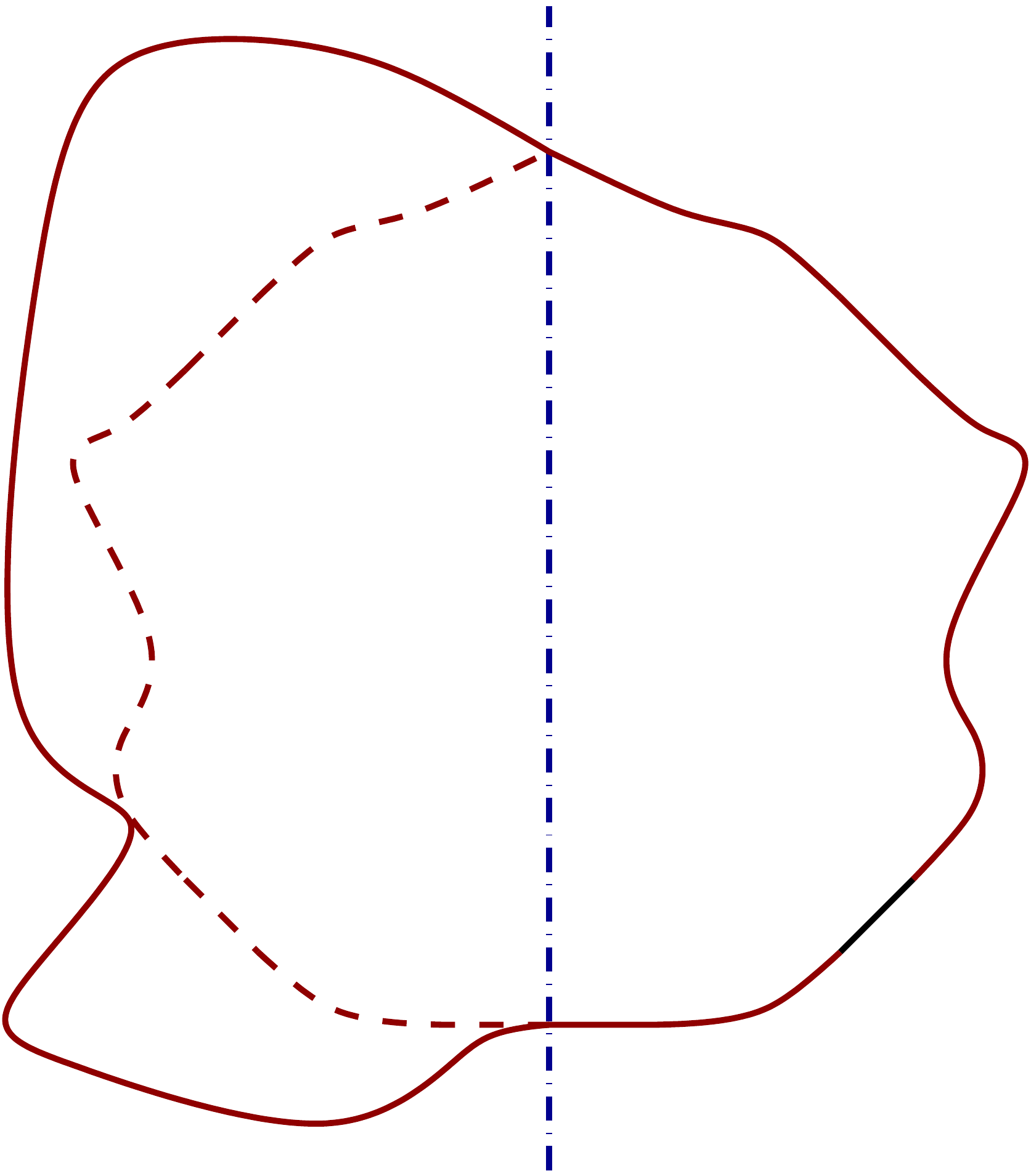}
 \caption{\sl The moving plane method up and running.}\label{MPMETHTERRAeFI}
\end{figure}

For instance, in the Soap Bubble Theorem, the critical position is the
``first contact'' occurrence between the reflected domain and the original one
(or, better to say, stops the method to the latest position in which the reflected domain is included into the original one),
see Figure~\ref{MPMETHTERRAeFI} for a sequence of frames trying to dynamically
depict the moving plane method in this situation (of course, for different problems,
such as the problem of Gidas, Ni and Nirenberg,
the critical position may have a different notion, also related to the equation into consideration).
At the contact point(s), some clever use of the Maximum Principle typically allows one to conclude that the reflected and original
domains necessarily coincide, thus providing a symmetry result with respect to this special hyperplane (by repeating
the procedure using hyperplanes in different directions, one can also deduce rotational symmetry).

\begin{figure}
  \centering
  \includegraphics[width=.44\linewidth]{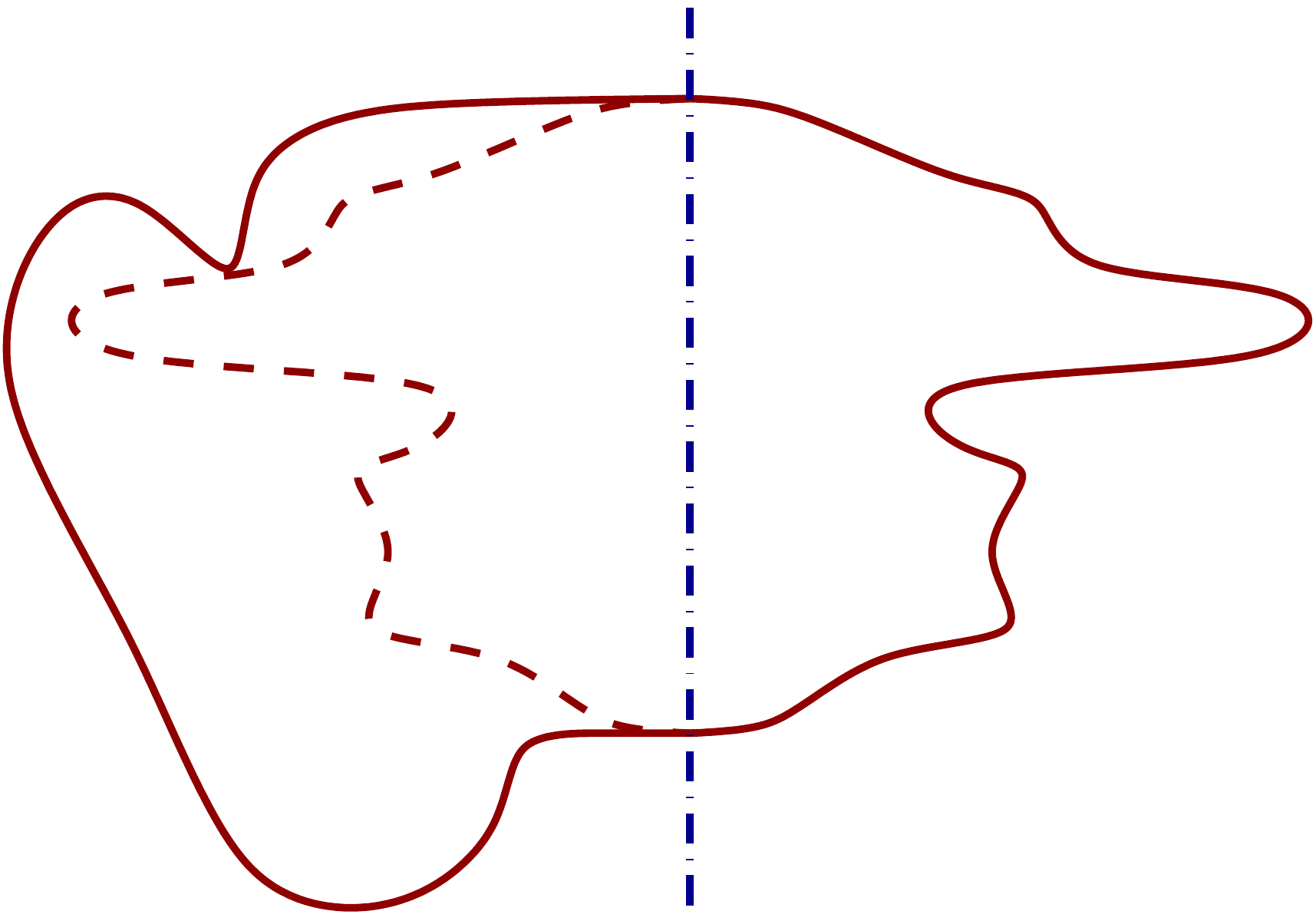}$\qquad\quad$
    \includegraphics[width=.44\linewidth]{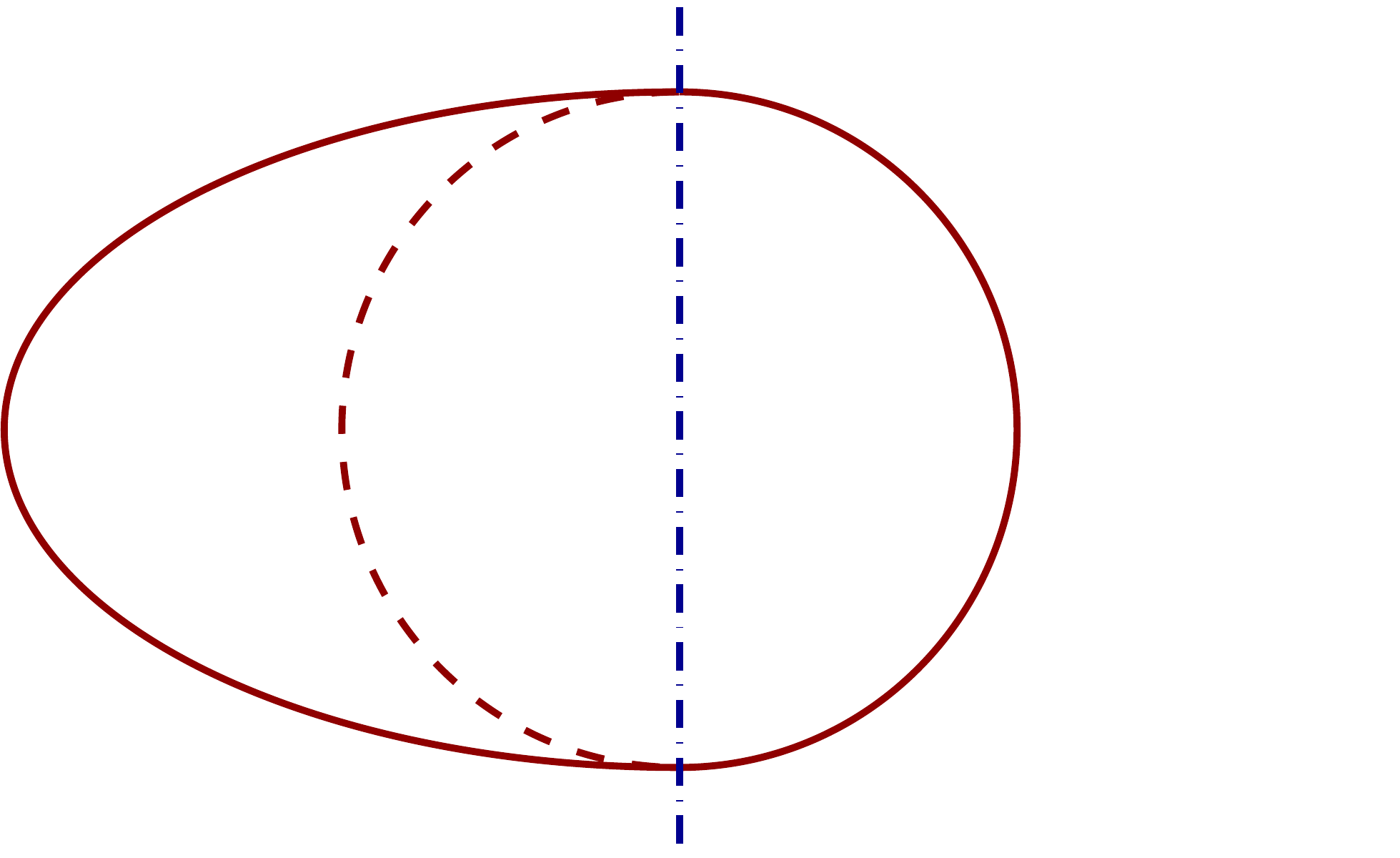}
 \caption{\sl Interior and boundary touching points.}\label{MPME242566THTERRAeFI}
\end{figure}

Of course, a number of conceptual and technical hindrances arise when trying to concretely implement this
elegant idea. A common treat is that one needs to distinguish the situation in which an ``interior'' touching point
between the reflected and the original domains arises and the case in which a ``boundary'' obstruction takes place,
see the last frame in Figure~\ref{MPMETHTERRAeFI}
and Figure~\ref{MPME242566THTERRAeFI} for possible pictures of these different scenarios.

\begin{figure}
  \centering
  \includegraphics[width=.39\linewidth]{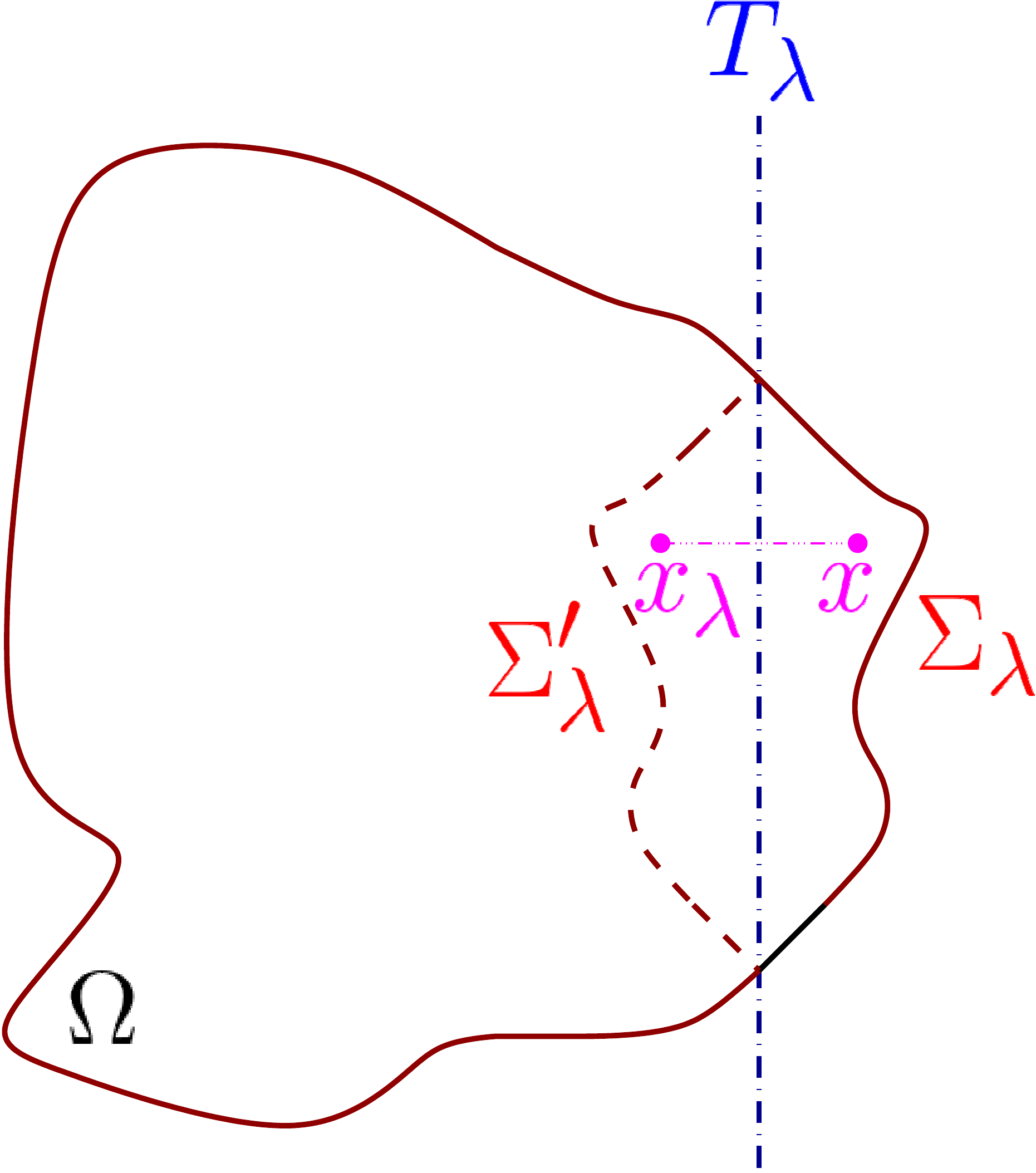}
 \caption{\sl Notation adopted for the moving plane method.}\label{MPM3675869E242566THTERRAeFI}
\end{figure}

Throughout this chapter, we will employ the notation sketched in Figure~\ref{MPM3675869E242566THTERRAeFI}.
Namely, given a set~$\Omega$ and~$\lambda\in\R$, we let
\begin{equation}\label{RIFLECONVE-eMPM3675869E242566THTERRAeFI2}
\begin{split}&
T_\lambda:=\{ x=(x_1,\dots,x_n)\in\R^n {\mbox{ s.t. }}x_1=\lambda\},\\&
\Sigma_\lambda:=\{ x=(x_1,\dots,x_n)\in\Omega {\mbox{ s.t. }}x_1>\lambda\},\\&
x_\lambda:=(2\lambda-x_1,x_2,\dots,x_n)\\ {\mbox{and }}\;&
\Sigma_\lambda':=\{ x_\lambda {\mbox{ s.t. }}x\in\Sigma_\lambda\}.
\end{split}\end{equation}
Notice that the point~$x_\lambda$ is the reflection of the point~$x=(x_1,\dots,x_n)$ with respect to the hyperplane~$T_\lambda$. Also, $\Sigma_\lambda'$ is the reflection of the region~$\Sigma_\lambda$
of~$\Omega$ lying to the right of~$T_\lambda$.
\medskip

\begin{figure}
  \centering
  \includegraphics[width=.44\linewidth]{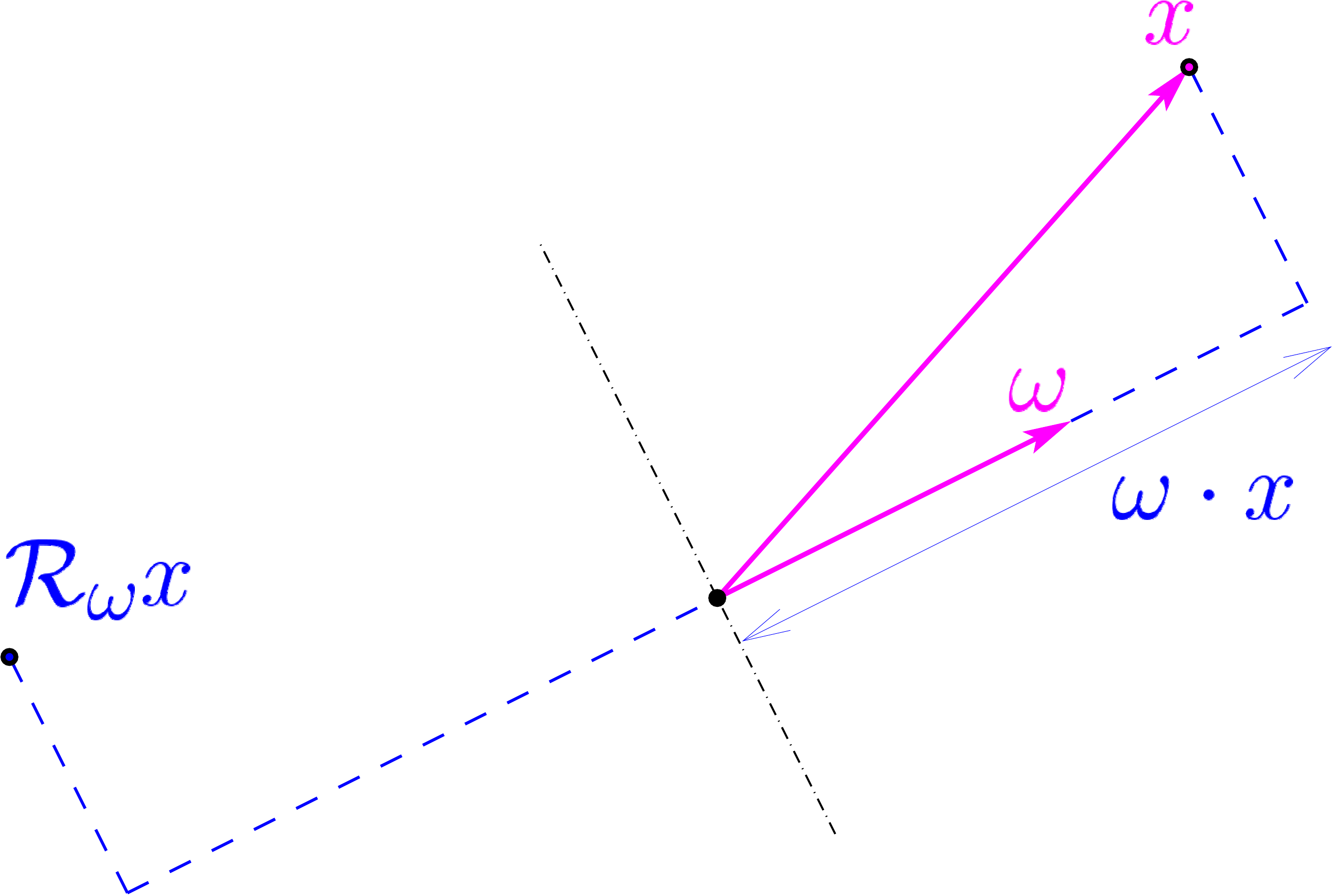}
 \caption{\sl Reflection of a point~$x$ in direction~$\omega$.}\label{RIFLECONVE-eMPM3675869E242566THTERRAeFI}
\end{figure}

Moreover, given~$\omega\in \partial B_1$, we define the reflection of a point~$p$ with respect to the plane~$\{\omega\cdot x=0\}$ by
\begin{equation}\label{CHIQSCECTQCSNLCD123rAS} {\mathcal{R}}_\omega p:=p-2(\omega\cdot p)\omega,\end{equation}
see Figure~\ref{RIFLECONVE-eMPM3675869E242566THTERRAeFI}
(notice also that when~$\omega=e_1$ the definition above boils down to the one of~$x_\lambda$ in~\eqref{RIFLECONVE-eMPM3675869E242566THTERRAeFI2} with~$\lambda=0$, consistently with the reflection with respect to the hyperplane~$\{x_1=0\}$).
\medskip

We refer to~\cite{MR1159383, MR1190345, MR3932952}
and the references therein
for further information about the moving plane method
and its variants in a number of different situations. Now we start watching the moving plane method
in action in some concrete situations.

\section{Back to the Soap Bubble Theorem}\label{MPPALEB}

We now retake the discussion on the Soap Bubble Theorem\index{Soap Bubble Theorem|(}
(presented in Theorem~\ref{SOAPTH}) and we outline here a different proof based on the reflection technique
depicted in Figure~\ref{MPMETHTERRAeFI}.

For this proof, we first establish the additional (perhaps not completely obvious,
see e.g.~\cite[Lemma~2.2 on page~147]{MR1013786})
result that reflections in every directions are sufficient to
prove rotational symmetry, namely:

\begin{lemma}\label{H3O7M8A5C2O21M1EF5K7A90654}
Let~$\Omega$ be as in Theorem~\ref{SOAPTH}
and assume that it admits
a hyperplane of symmetry in every direction (that is, for each~$\omega\in\partial B_1$ there exists~$x_\omega\in\R^n$ such that~${\mathcal{R}}_\omega (p-x_\omega)+x_\omega\in \Omega$
for all~$p\in \Omega$).

Then, $\Omega$ is spherically symmetric (that is, there exists~$x_0$ and a rotation~$R_0$ such that~$R_0(p-x_0)+x_0\in \Omega$ for every~$p\in \Omega$).
\end{lemma}

For this, we point out a general observation\footnote{Lemma~\ref{JNSMPERADVBDDOBYHGGJUINMRARYS3123-4}
is in fact a particular case of the general fact from linear algebra that
the reflections generate the orthogonal group, according to the Cartan-Dieudonn\'e Theorem.
For simplicity, here we present an elementary approach to the problem under consideration.}
relating reflections and rotations.

\begin{lemma}\label{JNSMPERADVBDDOBYHGGJUINMRARYS3123-4}
Let~$u:\R^n\to\R$. Suppose that for every~$\omega\in \partial B_1$ and~$x\in\R^n$ we have that~$u({\mathcal{R}}_\omega x)=u(x)$. 

Then, $u$ is rotationally invariant, namely there exists~$u_0:\R\to\R$ such that~$u(x)=u_0(|x|)$ for all~$x\in\R^n$.\end{lemma}

\begin{proof} Let~$r>0$ and~$p$, $q\in\partial B_r$. We claim that
\begin{equation}\label{90uojges3et43mi54n87u9s}
u(p)=u(q).
\end{equation}
To check this, we can suppose that~$p\ne q$, otherwise we are done
and we let~$\omega:=\frac{p-q}{|p-q|}$. Thus, we have that
\begin{eqnarray*}
&&u(p)=u\big({\mathcal{R}}_\omega p\big)=u\left( p-2(\omega\cdot p)\omega\right)
=u\left( p-2\left(\frac{p-q}{|p-q|}\cdot p\right)\frac{p-q}{|p-q|}\right)
\\&&\qquad=u\left( p-2\frac{(r^2-p\cdot q)(p-q)}{|p-q|^2}\right)=u\left( p-2\frac{(r^2-p\cdot q)(p-q)}{r^2+r^2-2p\cdot q}\right)
=u\left( p-(p-q)\right)=u(q),
\end{eqnarray*}
which proves~\eqref{90uojges3et43mi54n87u9s}.

Thus, we define~$u_0(r):=u(re_1)$ and~\eqref{90uojges3et43mi54n87u9s} comes in handy to see that,
for all~$x\in\R^n$, $u(x)=u(|x|e_1)=u_0(|x|)$, as desired.
\end{proof}

\begin{proof}[Proof of Lemma~\ref{H3O7M8A5C2O21M1EF5K7A90654}]
Up to a translation, we can suppose that the baricenter of~$\Omega$ is the origin, namely that
$$ \int_\Omega x\,dx=0.$$
Then, using the isometric change of variable induced by the reflection
that leaves~$\Omega$ invariant,
that is, looking at the transformation~$x={\mathcal{R}}_\omega (y-x_\omega)+x_\omega$,
we see that
\begin{eqnarray*}&& 0=\int_\Omega\big({\mathcal{R}}_\omega (y-x_\omega)+x_\omega\big)\,dy=
\int_\Omega\big( (y-x_\omega)
-2\omega\cdot(y-x_\omega)\omega
+x_\omega\big)\,dy\\&&\qquad\qquad\qquad=-|\Omega|\,x_\omega+2|\Omega|(\omega\cdot x_\omega)\omega+|\Omega|\,x_\omega=
2|\Omega|(\omega\cdot x_\omega)\omega
\end{eqnarray*}
and consequently
\begin{equation*}
\omega\cdot x_\omega=0.
\end{equation*}
As a result, 
$$ {\mathcal{R}}_\omega (p-x_\omega)+x_\omega=
(p-x_\omega)-2\omega\cdot (p-x_\omega)\omega+x_\omega
=p-2(\omega\cdot p)\omega={\mathcal{R}}_\omega p$$
and thus
the invariance of~$\Omega$ under hyperplane reflections reads\footnote{Geometrically, \eqref{MRMKMJP9io3rn374yth4g} states the interesting property that the reflections hyperplanes
all pass through the baricenter, which was not completely obvious at the beginning.}
\begin{equation}\label{MRMKMJP9io3rn374yth4g}
{\mbox{${\mathcal{R}}_\omega p\in \Omega$ for all~$p\in \Omega$.}}\end{equation}
This gives that~$\chi_\Omega({\mathcal{R}}_\omega x)=\chi_\Omega(x)$ for all~$x\in\R^n$.
We can thereby employ Lemma~\ref{H3O7M8A5C2O21M1EF5K7A90654}
with~$u:=\chi_\Omega$ and deduce that~$\chi_\Omega$ is radial, hence so is~$\Omega$.
\end{proof}

Thanks to these preliminary observations, we can now give a proof of the Soap Bubble Theorem
by utilizing the moving plane method, seeking for symmetry under hyperplane reflections.
For this, we follow the procedure described in Figure~\ref{MPMETHTERRAeFI}
(notice that since~$\partial \Omega$ is of class~$C^{2,\alpha}$ at least at the beginning of the reflections,
the flipped region~$\Sigma_\lambda'$ lies inside~$\Omega$, hence we can start the method).
Thus, we can define
\begin{equation}\label{800ntWAdw3G}
\lambda_0:=\inf\Big\{\lambda\in\R {\mbox{ s.t. $\Sigma'_t\subset\Omega$ 
for all $t>\lambda$}}
\Big\}.
\end{equation}
In this critical configuration, we have that~$\Sigma_{\lambda_0}'$ becomes tangent to~$\partial\Omega$
at some point (otherwise the distance between~$\partial\Sigma_{\lambda_0}'$ and~$\partial\Omega$
would be strictly positive, thus allowing a small further push towards the left of the hyperplane~$T_\lambda$,
in contradiction with the infimum property in~\eqref{800ntWAdw3G}).

Hence, as highlighted in the last frame in Figure~\ref{MPMETHTERRAeFI}
and Figure~\ref{MPME242566THTERRAeFI}, it is convenient to subdivide the critical case arising in~\eqref{800ntWAdw3G}
into two possible\footnote{As a matter of fact, 
both cases~(i) and~(ii) in~\eqref{IDUECASIMPS}
can also occur simultaneously, see e.g. the first picture in
Figure~\ref{MPME242566THTERRAeFI}. In this scenario, all the better,
one is free to pick the point~$p_0\in(\partial \Sigma_{\lambda_0}')\cap(\partial\Omega)$ either on the hyperplane~$T_{\lambda_0}$ or outside it, at their own preference.}
subcases:
\begin{equation}\label{IDUECASIMPS}
\begin{split}
&{\mbox{(i). $\qquad \Sigma_{\lambda_0}'$ is tangent to~$\partial\Omega$ at a point~$p_0$ {\em not lying} on the
hyperplane~$T_{\lambda_0}$,}}
\\&{\mbox{(ii). $\qquad\Sigma_{\lambda_0}'$ is tangent to~$\partial\Omega$ at a point~$p_0$ {\em lying} on the
hyperplane~$T_{\lambda_0}$.}}
\end{split}
\end{equation}
One could argue that case~(i) is somewhat ``generic'' and one has to be quite unlucky to bump into case~(ii),
yet one cannot rule out this perhaps uncommon and less pleasant possibility. For this, one relies
on suitable versions of the Maximum Principle,
possibly exploiting different versions for case~(i) and case~(ii). 
The strategy is to compare the original set and the reflected one at the tangent point
described in~\eqref{IDUECASIMPS}. For this, one locally parameterizes~$\partial\Omega$
over its tangent space by a function that satisfies an elliptic equation (specifically, in this case,
the mean curvature equation) and then applies the Maximum Principle
to the function obtained by subtracting the original solution and the reflected one.

To make the argument work, one also needs to investigate a bit more the geometry
of~$\Sigma'_{\lambda_0}$: on the one hand, in general~$\Sigma'_{\lambda_0}$ is not necessarily convex, and, even worse, not necessarily connected, see Figure~\ref{JMDSJORGFDFN8SIGMA0h5fdKSLE-21ItangeFIrfd72iqwhrjfngpoewjnikdhgebidkkwdbbXVSdng}.
On the other hand, we have that~$\Sigma'_{\lambda_0}$
is always ``convex in the~$e_1$ direction'', as made precise in the following result:

\begin{figure}
  \centering
  \includegraphics[width=.4\linewidth]{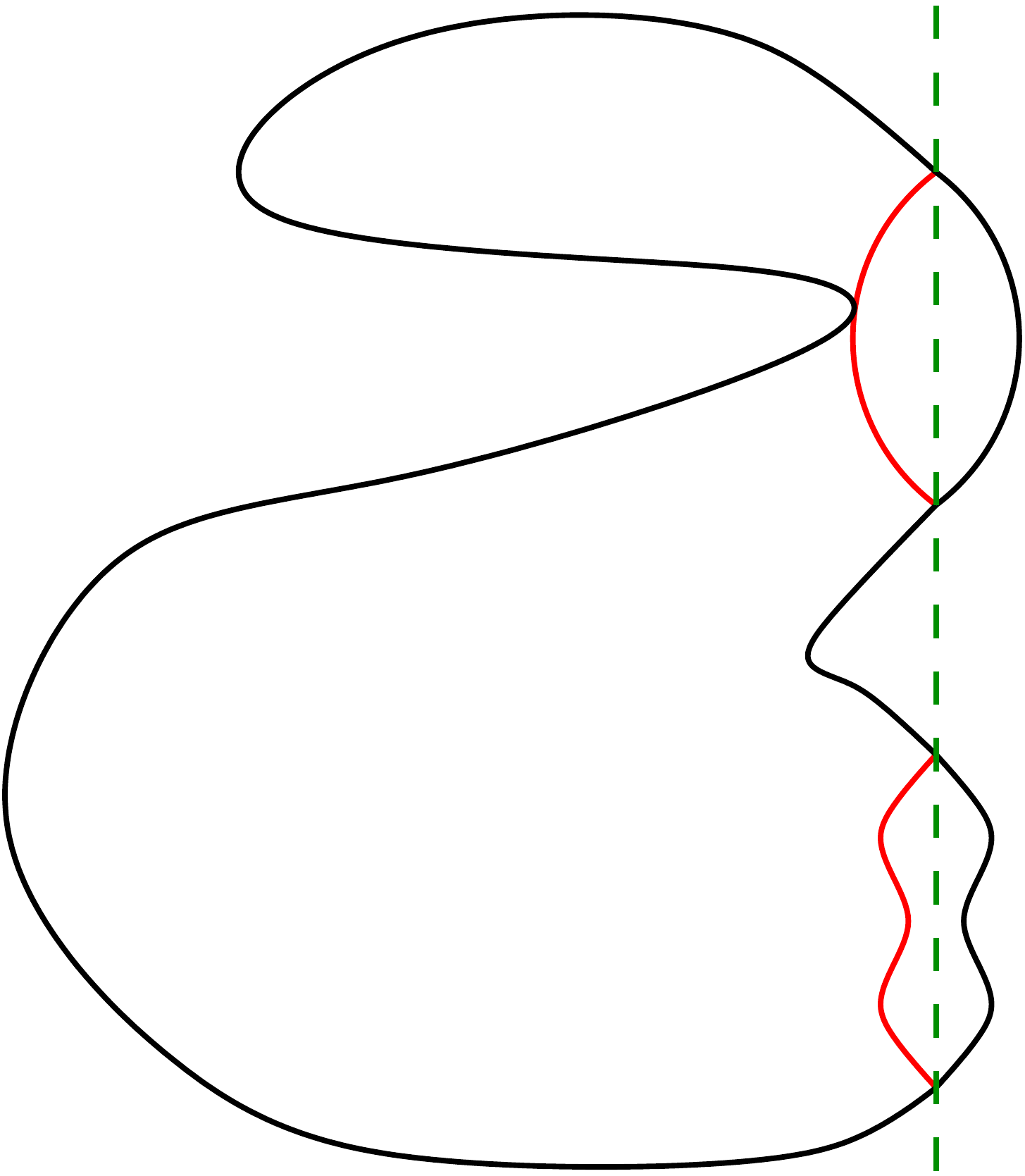}
 \caption{\sl $\Sigma'_0$ is not necessarily convex, nor connected.}\label{JMDSJORGFDFN8SIGMA0h5fdKSLE-21ItangeFIrfd72iqwhrjfngpoewjnikdhgebidkkwdbbXVSdng}
\end{figure}

\begin{lemma}\label{KNS-kwjfluoejdbjfjfj93itkSY} Let~$x=(x_1,\dots,x_n)\in\Sigma'_{\lambda_0}$. Then, $(\tau,x_2,\dots,x_n)\in\Sigma_{\lambda_0}'$
for all~$\tau\in[x_1,\lambda_0)$.
\end{lemma}

\begin{proof} Since~$\Sigma'_{\lambda_0}$ is open and bounded, we can move horizontally from~$x$
towards the left, till we hit the boundary of~$\Sigma'_{\lambda_0}$: that is, we take~$\zeta_0<x_1$
such that~$(\zeta_0,x_2,\dots,x_n)\in\partial \Sigma'_{\lambda_0}$ and \begin{equation}\label{9oj2rlmnfwzeta0}
{\mbox{$
(\zeta,x_2,\dots,x_n)\in \Sigma'_{\lambda_0}$ for all~$\zeta\in(\zeta_0,x_1]$.}}\end{equation}

Now, to prove the desired claim in Lemma~\ref{KNS-kwjfluoejdbjfjfj93itkSY}, we will actually show the following strongest result:
\begin{equation*}
{\mbox{$(\tau,x_2,\dots,x_n)\in\Sigma_{\lambda_0}'$
for all~$\tau\in(\zeta_0,\lambda_0)$.}}
\end{equation*}
For this, it suffices to show that
\begin{equation*}
{\mbox{$(2\lambda_0-\tau,x_2,\dots,x_n)\in\Sigma_{\lambda_0}$
for all~$\tau\in(\zeta_0,\lambda_0)$.}}
\end{equation*}
or equivalently, calling~$\vartheta:=2\lambda_0-\tau$,
\begin{equation}\label{taumbda0foraltauambda0}
{\mbox{$(\vartheta,x_2,\dots,x_n)\in\Omega$
for all~$\vartheta\in(\lambda_0,2\lambda_0-\zeta_0)$.}}
\end{equation}
To check this we start by noticing that, in light of~\eqref{9oj2rlmnfwzeta0},
we have that~$(2\lambda_0-\zeta,x_2,\dots,x_n)\in \Omega$ for all~$\zeta\in(\zeta_0,x_1]$,
whence~$(\vartheta,x_2,\dots,x_n)\in \Omega$ for all~$\zeta\in[2\lambda_0-x_1,2\lambda_0-\zeta_0)$.
As a result, we can define
$$ \lambda^\star:=\inf\Big\{
{\mbox{$t<2\lambda_0-x_1$ s.t.
$(\vartheta,x_2,\dots,x_n)\in\Omega$
for all~$\vartheta\in(t,2\lambda_0-\zeta_0)$}}
\Big\},$$
we observe that~$\lambda^\star\le2\lambda_0-x_1$
and, to prove~\eqref{taumbda0foraltauambda0},
we aim at showing that~$\lambda^\star\le\lambda_0$.

For this, suppose not: then
\begin{equation}\label{robVMSHcsta-d0k4-0}
\lambda^\star>\lambda_0.\end{equation}
We know that
\begin{equation}\label{robVMSHcsta-d0k4-1}
{\mbox{$(\vartheta,x_2,\dots,x_n)\in\Omega$ for all~$\vartheta\in(\lambda^\star,2\lambda_0-\zeta_0)$}}.
\end{equation}
%%but
%%\begin{equation}\label{robVMSHcsta-d0k4-2}
%%{\mbox{$(\lambda^\star,x_2,\dots,x_n)\in\partial\Omega$.}}\end{equation}
{F}rom~\eqref{800ntWAdw3G} and~\eqref{robVMSHcsta-d0k4-1}, we see that
\begin{equation}\label{robVMSHcsta-d0k4-2}
{\mbox{$(2t-\vartheta,x_2,\dots,x_n)\in\Sigma'_t\subseteq\Omega$
for all~$\vartheta\in(\lambda^\star,2\lambda_0-\zeta_0)$ and~$t>\lambda_0$.}}\end{equation}
In particular, we take~$\e>0$ small enough such that~$\vartheta_\e:=\lambda^\star+2\e\in(\lambda^\star,2\lambda_0-\zeta_0)$
and~$t_\e:=\lambda^\star-\e>\lambda_0$ (we stress that we can fulfill the latter condition
owing to~\eqref{robVMSHcsta-d0k4-0}): with this choices, we infer from~\eqref{robVMSHcsta-d0k4-2} that
$$\Omega\ni(2t_\e-\vartheta_\e,x_2,\dots,x_n)
=\big(2(\lambda^\star-\e)-(\lambda^\star+2\e),x_2,\dots,x_n\big)=(\lambda^\star-4\e,x_2,\dots,x_n),
$$
which is in contradiction with the minimality of~$\lambda^\star$.
\end{proof}

\begin{corollary}\label{ssitingtotsendirow}
Let~$\Sigma''_{\lambda_0}$ be the interior of the closure of~$\Sigma_{\lambda_0}\cup\Sigma_{\lambda_0}'$.
Assume that~$(a,x_2,\dots,x_n)$, $(b,x_2,\dots,x_n)\in\Sigma''_{\lambda_0}$ for some~$a<b$.
Then, $(\tau,x_2,\dots,x_n)\in\Sigma''_{\lambda_0}$ for all~$\tau\in[a,b]$.
\end{corollary}

\begin{figure}
  \centering
  \includegraphics[width=.45\linewidth]{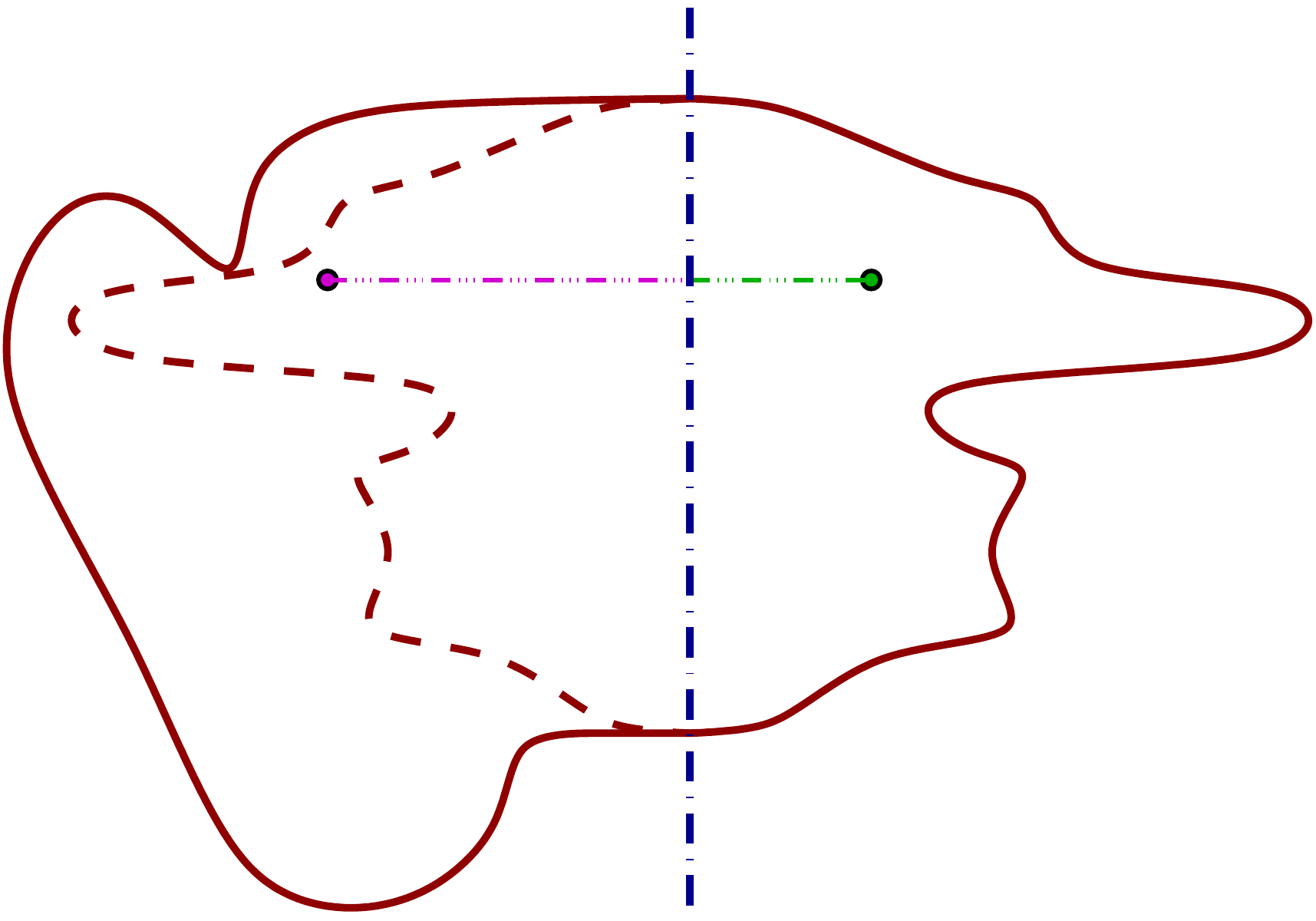}
 \caption{\sl $\Sigma''_{\lambda_0}$ is convex in the~$e_1$ direction.}\label{CONVE1JMDSJORGFDFN8h5fdKSLE-21ItangeFIrfd72iqwhrjfngpoewjnikdhgebidkkwdbbXVSdng}
\end{figure}

\begin{proof} By Lemma~\ref{KNS-kwjfluoejdbjfjfj93itkSY}, we have that~$\Sigma'_{\lambda_0}$
is convex in the~$e_1$ direction, hence so is~$\Sigma_{\lambda_0}$ by symmetry.
This gives the desired result when both~$(a,x_2,\dots,x_n)$ and~$(b,x_2,\dots,x_n)$
lie either in~$\Sigma'_{\lambda_0}$ or in~$\Sigma_{\lambda_0}$.

It remains to consider the case in which~$a\le\lambda_0$ and~$b\ge\lambda_0$.
In this case, one uses Lemma~\ref{KNS-kwjfluoejdbjfjfj93itkSY}
to connect~$(a,x_2,\dots,x_n)$ to~$(\lambda_0,x_2,\dots,x_n)$
(possibly connecting it first to~$(\lambda_0-\e,x_2,\dots,x_n)$ and then sending~$\e\searrow0$)
and similarly to connect~$(b,x_2,\dots,x_n)$ to~$(\lambda_0,x_2,\dots,x_n)$
(possibly connecting it first to~$(\lambda_0+\e,x_2,\dots,x_n)$ and then sending~$\e\searrow0$).
See Figure~\ref{CONVE1JMDSJORGFDFN8h5fdKSLE-21ItangeFIrfd72iqwhrjfngpoewjnikdhgebidkkwdbbXVSdng} for
a sketch of this configuration.
\end{proof}

With this preliminary work, we can now give the proof of the Soap Bubble Theorem
by using the original approach of Aleksandrov~\cite{MR0150710},
via the moving plane method and the reflection principle.

\begin{proof}[Another proof of Theorem~\ref{SOAPTH}]
Suppose first that we are in case~(i) of~\eqref{IDUECASIMPS}.
Then, in the vicinity of~$p_0$, we can write~$\partial\Omega$
and~$\partial\Sigma_{\lambda_0}'$ as the graphs of two functions,
say~$u_1$ and~$u_2$ respectively, that satisfy the mean curvature equation in~\eqref{LAMEANCURVA}.
That is, up to a rigid motion, if~$H\in\R$ denotes the constant mean curvature of~$\partial\Omega$, there exists~$\rho>0$
such that for all~$x\in\R^{n-1}$ with~$|x|<\rho$ we have that
\begin{equation}\label{HSETBDSOB0lk4:1}
\div\left(\frac{\nabla u_1(x)}{\sqrt{1+|\nabla u_1(x)|^2}}\right)=\div\left(\frac{\nabla u_2(x)}{\sqrt{1+|\nabla u_2(x)|^2}}\right)=H,\end{equation}
with
\begin{equation}\label{HSETBDSOB0lk4:2}
u_1(0)=u_2(0),\qquad \nabla u_1(0)=\nabla u_2(0)=0 \qquad{\mbox{and}}\qquad
u_1\ge u_2.\end{equation} Possibly reducing~$\rho$,
we can also suppose that
\begin{equation}\label{HSETBDSOB0lk4:3}
|\nabla u_1|+|\nabla u_2|\le1.\end{equation}
We also remark that, for~$m\in\{1,2\}$,
\begin{eqnarray*}
&&\div\left(\frac{\nabla u_m(x)}{\sqrt{1+|\nabla u_m(x)|^2}}\right)\\
&=&\frac{\Delta u_m(x)}{\sqrt{1+|\nabla u_m(x)|^2}}-
\frac{D^2 u_m(x)\nabla u_m(x)\cdot\nabla u_m(x)}{(1+|\nabla u_m(x)|^2)^{\frac32}}\\
&=&\sum_{i,j=1}^n \alpha_{ij}(\nabla u_m(x))\partial_{ij} u_m(x),
\end{eqnarray*}
where
\begin{equation}\label{ldfgrgite98875231456217485968765}
\alpha_{ij}(\zeta):=
\frac{\delta_{ij}}{\sqrt{1+|\zeta|^2}}-
\frac{\zeta_i\,\zeta_j}{(1+|\zeta|^2)^{\frac32}}.\end{equation}

Thus, if~$w:=u_1-u_2$ we have that~$w\ge0$ and~$w(0)=0$. Also,
\begin{eqnarray*}
0&=&\div\left(\frac{\nabla u_1(x)}{\sqrt{1+|\nabla u_1(x)|^2}}\right)-\div\left(\frac{\nabla u_2(x)}{\sqrt{1+|\nabla u_2(x)|^2}}\right)\\
&=& \sum_{i,j=1}^n \alpha_{ij}(\nabla u_1(x))\partial_{ij} u_1(x)-\sum_{i,j=1}^n \alpha_{ij}(\nabla u_2(x))\partial_{ij} u_2(x)\\
&=& \sum_{i,j=1}^n \alpha_{ij}(\nabla u_1(x))\partial_{ij} w(x)
+\sum_{i,j=1}^n \Big(\alpha_{ij}(\nabla u_1(x))-\alpha_{ij}(\nabla u_2(x))\Big)\partial_{ij} u_2(x)\\
&=& \sum_{i,j=1}^n \alpha_{ij}(\nabla u_1(x))\partial_{ij} w(x)
+\sum_{i,j=1}^n \partial_{ij} u_2(x)\int_0^1
\frac{d}{dt}\Big( \alpha_{ij}\big(t\nabla u_1(x)+(1-t)\nabla u_2(x)\big)\Big)\,dt\\
&=&\sum_{i,j=1}^n \alpha_{ij}(\nabla u_1(x))\partial_{ij} w(x)\\&&\qquad
+\sum_{i,j=1}^n \partial_{ij} u_2(x)\int_0^1
\nabla\alpha_{ij}\big(t\nabla u_1(x)+(1-t)\nabla u_2(x)\big)\,dt
\cdot\big(\nabla u_1(x)-\nabla u_2(x)\big)
\end{eqnarray*}
and therefore
\begin{equation}\label{orollarHOPF12ELLIPTIC-PERFVDJ-LEM-32-CORo}
\sum_{i,j=1}^n a_{ij}(x)\partial_{ij} w(x)+b(x)
\cdot \nabla w(x)=0,
\end{equation}
where
\begin{eqnarray*}
a_{ij}(x)&:=&\alpha_{ij}(\nabla u_1(x))\\
{\mbox{and }}\;b(x)&:=&\sum_{i,j=1}^n \partial_{ij} u_2(x)\int_0^1
\nabla\alpha_{ij}\big(t\nabla u_1(x)+(1-t)\nabla u_2(x)\big)\,dt.
\end{eqnarray*}
We point out that 
\begin{equation}\label{sjxewr45gwiheugherujeri496y98}
{\mbox{the coefficients~$a_{ij}$ satisfies the ellipticity condition
in~\eqref{ELLIPTIC-PERFVDJ}.}}\end{equation}
Indeed, from~\eqref{ldfgrgite98875231456217485968765}
we see that, for all~$\xi\in\partial B_1$,
\begin{eqnarray*}
&& \sum_{i,j=1}^n a_{ij}\xi_i\xi_j=\sum_{i,j=1}^n\alpha_{ij}(\nabla u_1)\xi_i\xi_j=
\sum_{i,j=1}^n
\left(\frac{\delta_{ij}}{\sqrt{1+|\nabla u_1|^2}}-
\frac{\partial_i u_{1}\partial_j u_1}{(1+|\nabla u_1|^2)^{\frac32}}\right)\xi_i\xi_j\\&&\qquad=
\sum_{i,j=1}^n
\left(\frac{\delta_{ij}(1+|\nabla u_1|^2)-
\partial_i u_{1}\partial_j u_1}{(1+|\nabla u_1|^2)^{\frac32}}\right)\xi_i\xi_j
= \frac{|\xi|^2(1+|\nabla u_1|^2)-(\nabla u_1\cdot\xi)^2}{(1+|\nabla u_1|^2)^{\frac32}
}\\&&\qquad\le \frac{|\xi|^2 }{\sqrt{1+|\nabla u_1|^2}}\le |\xi|^2=1
\end{eqnarray*}
and similarly, using~\eqref{HSETBDSOB0lk4:3}
\begin{eqnarray*}
&&  \sum_{i,j=1}^n a_{ij}\xi_i\xi_j= \frac{|\xi|^2(1+|\nabla u_1|^2)-(\nabla u_1\cdot\xi)^2}{(1+|\nabla u_1|^2)^{\frac32}}\ge
\frac{|\xi|^2(1+|\nabla u_1|^2)-|\nabla u_1|^2|\xi|^2}{(1+|\nabla u_1|^2)^{\frac32}}\\&&\qquad
=\frac{|\xi|^2 }{(1+|\nabla u_1|^2)^{\frac32}}\ge \frac{|\xi|^2}{2^{\frac32}}
=\frac{1}{2^{\frac32}}
,\end{eqnarray*}
thus completing the proof of~\eqref{sjxewr45gwiheugherujeri496y98}.

{F}rom these considerations
and the Maximum Principle in Theorem~\ref{MAXPLECON45678CG-PM},
we conclude that~$w$ vanishes identically.
As a result, we have that~$\partial\Omega$ and~$\partial\Sigma'_{\lambda_0}$ coincide in
a neighborhood of~$p_0$.

Suppose now that we are in case~(ii) of~\eqref{IDUECASIMPS}.
Then, in the vicinity of~$p_0$, we can write~$\partial\Omega$
and~$\partial\Sigma_{\lambda_0}'$ as the graphs of two functions,
say~$u_1$ and~$u_2$ respectively, which touch each other at the boundary of their domain of definition:
that is, the setting in~\eqref{HSETBDSOB0lk4:1}, \eqref{HSETBDSOB0lk4:2}, \eqref{HSETBDSOB0lk4:3}
holds true, but the functions~$u_1$ and~$u_2$ are defined in~$D:=\big\{x\in\R^{n-1} {\mbox{ s.t. $|x|<\rho$ and $x_1\ge0$}}\big\}$
and the touching between~$u_1$ and~$u_2$, occurring at the origin, lies on the boundary of the domain~$D$.
We can therefore define~$w:=u_1-u_2$, recall~\eqref{orollarHOPF12ELLIPTIC-PERFVDJ-LEM-32-CORo}
and invoke Corollary~\ref{HOPF12ELLIPTIC-PERFVDJ-LEM-32-CORo}:
the latter result, since~$\nabla w(0)=0$, yields that~$w$ vanishes identically.
Therefore, also in this case we have proved that~$\partial\Omega$ and~$\partial\Sigma'_{\lambda_0}$ coincide in
a neighborhood of~$p_0$.

These considerations show that, in both cases~(i) and~(ii) of~\eqref{IDUECASIMPS}, the
set in which~$\partial\Omega$ and~$\partial\Sigma'_{\lambda_0}$ coincide
is open (as well as closed, by continuity, and nonempty, due to the presence of~$p_0$).
This gives that a connected component of~$(\partial\Sigma'_{\lambda_0})\cap\{x_1<\lambda_0\}$ (namely, the one
containing~$p_0$ in its closure) coincides with a connected component of~$(\partial\Omega)\cap\{x_1<\lambda_0\}$.
By symmetry, and recalling Corollary~\ref{ssitingtotsendirow}, we have thereby identified a connected component
of~$\Omega$ which is invariant under the reflection through~$T_{\lambda_0}$.
Since~$\Omega$ is assumed to be connected, we deduce that~$\Omega$ is invariant under the reflection through~$T_{\lambda_0}$.

We can apply this argument in every direction (not only~$e_1$).
In this way, we find that~$\Omega$ is invariant under hyperplane reflections in every direction
and we can therefore invoke Lemma~\ref{H3O7M8A5C2O21M1EF5K7A90654}
and conclude that~$\Omega$ is spherically symmetric.
This says that~$\Omega$ is the union of disjoint rings.
But applying again Corollary~\ref{ssitingtotsendirow} we deduce that~$\Omega$ is convex, hence it is necessarily a ball.\index{Soap Bubble Theorem|)}
\end{proof}

\section{Back to Serrin's overdetermined problem}\label{BACKSE}

We retake now Theorem~\ref{SERR} and present a different proof of it, based on the moving plane method,\index{Serrin's overdetermined problem|(}
in line with the original article~\cite{MR333220}.
The strategy is to slide a hyperplane in a given direction, say the first coordinate direction from right to left,
and reach a critical configuration, as stated in~\eqref{800ntWAdw3G}.
{F}rom that, one needs to analyze cases~(i) and~(ii) in~\eqref{IDUECASIMPS}.
Case~(i) is technically easier to deal with, relying on the Hopf Lemma
(e.g., Lemma~\ref{JJS:PA} would do the job). Case~(ii) turns out to be more elaborate
and hings on a specific result, known in the literature as Serrin's Corner Lemma\index{Serrin's Corner Lemma},
which goes as follows:

\begin{figure}
  \centering
  \includegraphics[width=.55\linewidth]{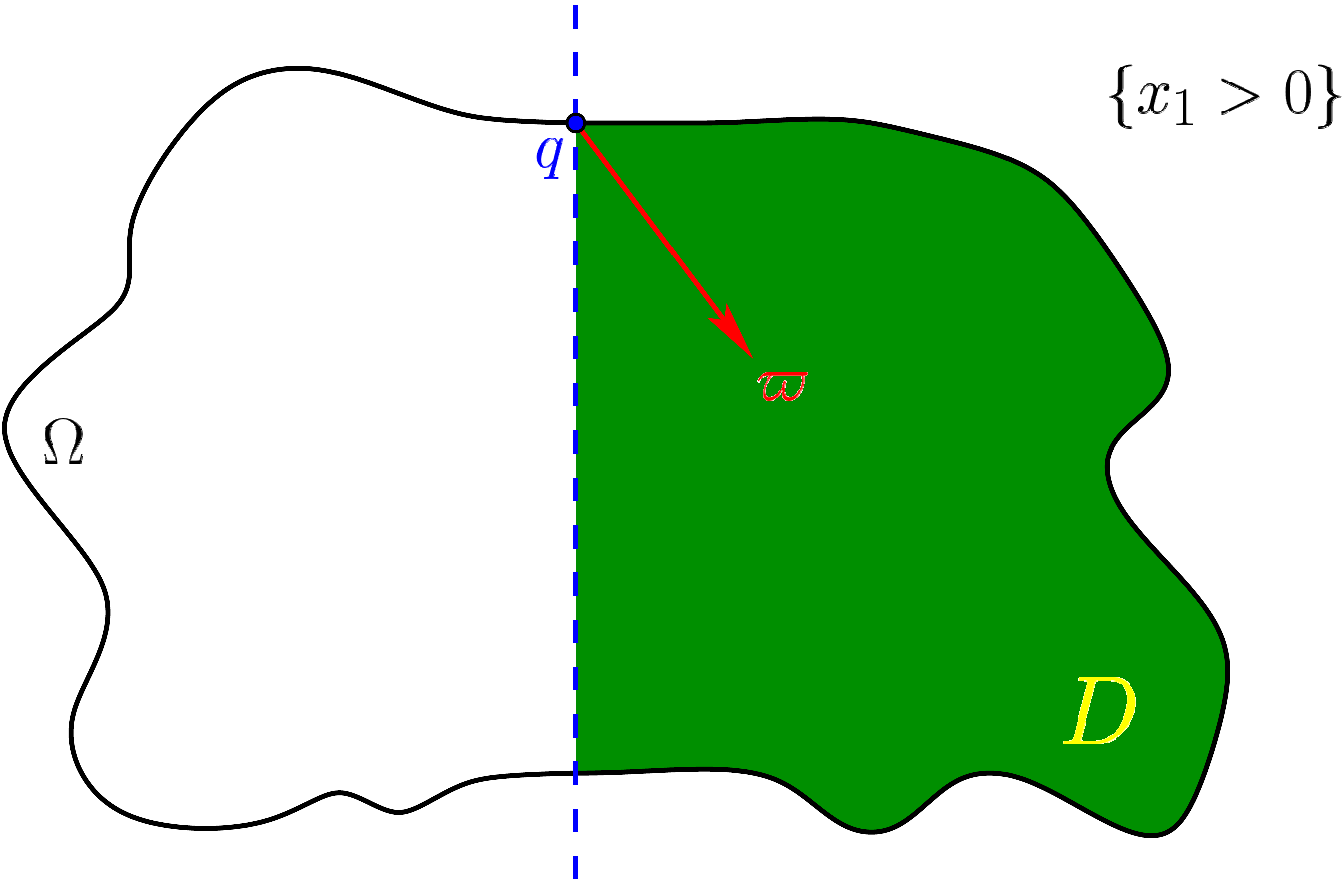}
 \caption{\sl The geometry involved in Serrin's Corner Lemma.}\label{J25376MDSJORGFDFN8h5fdKSLE-21ItangeFIrfd7CORN2iSEpokjgHGSKDqwhrjfngpoewjnikdhgebidkkwdbbXVSdng}
\end{figure}

\begin{lemma}\label{SerrinsCornerLemma}
Let~$\Omega$ be a bounded, open and connected set
of class~$C^2$ with exterior normal~$\nu$. Let~$D:=\Omega\cap\{x_1>0\}$. Let~$q\in (\partial\Omega)\cap\{x_1=0\}$
and suppose that~$\nu(q)=e_n$. Let also~$\varpi=(\varpi_1,\dots,\varpi_n)\in\partial B_1$, with~$\varpi_1>0>\varpi_n$.

Assume that~$w\in C^2(\overline{D})$ satisfies $$
\begin{dcases}
\Delta w\le0 &{\mbox{ in }}D,\\
w\ge0 &{\mbox{ in }}D,\\
w(q)=0 .&\end{dcases}$$
Then, either~$\partial_\varpi w(q)>0$ or~$\partial_{\varpi\varpi} w(q)>0$,
unless~$w$ vanishes identically in~$D$.
\end{lemma}

See Figure~\ref{J25376MDSJORGFDFN8h5fdKSLE-21ItangeFIrfd7CORN2iSEpokjgHGSKDqwhrjfngpoewjnikdhgebidkkwdbbXVSdng} for a sketch of the geometry involved in Lemma~\ref{SerrinsCornerLemma}.
In a sense, Lemma~\ref{SerrinsCornerLemma} addresses one of the cases left out by the theory
developed in Section~\ref{HOPFSECTH}. Namely, the classical setting for the Hopf Lemma
requires an interior ball condition (see Figure~\ref{HOPFLAGIJ7soloDItangeFIHOP}),
which is violated in the assumptions of Lemma~\ref{SerrinsCornerLemma},
since the domain~$D$ presents a right angled corner at the point~$q$.
In this scenario, normal derivatives of harmonic functions reaching their extrema at~$q$ may
vanish. As an example of this phenomenon, consider the case in which~$n=2$, $\Omega:=B_1(-e_2)$, $D:=\Omega\cap\{x_1>0\}$,
and~$w(x):=-x_1x_2$. We have that~$\Delta w=0$, $w(0)=0$ and~$w(x)\ge0$ in~$D$,
but~$\nabla w(0)=0$ thus violating the normal derivative nondegeneracy of the classical Hopf Lemma.

With respect to this occurrence, the power of Lemma~\ref{SerrinsCornerLemma}
is to detect that this example does provide ``the worst possible scenario'',
since the degeneracy of the first derivative (if it occurs) is compensated by a nondegeneracy of second derivatives: indeed, in the example above,~$D^2 w=\left(\begin{matrix}0&-1\\-1&0\end{matrix}\right)$.

We stress that the vector~$\varpi$ in the statement of Lemma~\ref{SerrinsCornerLemma} plays the role of a direction
entering~$D$ nontangentially from~$q$: as a matter of fact, since~$\nu(q)=e_n$, the condition~$\varpi_n\ne0$ is a nontangential
requirement and if~$\Omega$ is parameterized in a neighborhood of~$q$ by the sublevel sets of a function~$\Phi$, say with~$\Phi(q)=0$,
$\nabla \Phi(q)=e_n$
and~$\Phi<0$ in~$\Omega$ in the vicinity of~$q$,
for small~$t>0$ we have that
$$ \Phi(q+t\varpi)=t\nabla \Phi(q)\cdot \varpi+o(t)=t\varpi_n+o(t)<0,$$
hence~$q+t\varpi\in\Omega$, and also
$$ (q+t\varpi)\cdot e_1=0+t\varpi_1>0,$$
giving that~$q+t\varpi\in D$. This shows that
\begin{equation}\label{FMDMIDLAFFJZimnPO83n}
{\mbox{$q+t\varpi$ enters~$D$ for small~$t>0$.}}\end{equation}

In a nutshell, the proof of
Lemma~\ref{SerrinsCornerLemma}
consists in revisiting the argument used to establish the classical Hopf Lemma by carefully
taking into account an additional possible degeneracy
of linear type. Namely, we will still look at configurations such as the one described in Figure~\ref{LwAS1IewN24XGIJ7soloD354hItangeFIrfd72iqwhrjfngpoewjnikdhgebidkkwdbbXVSdng}
(here, up to an intersection with a halfspace)
but in the present framework we will ``correct'' the barrier in~\eqref{LwAS1IewN24XGIJ7soloD354hItangeFIrfd72iqwhrjfngpoewjnikdhgebidkkwdbbXVSdng9028032utDFVBmeKTRA920irkjtd}
by a linear factor.

The details go as follows.

\begin{proof}[Proof of Lemma~\ref{SerrinsCornerLemma}] We use the regularity of~$\Omega$ to pick a ball contained in~$\Omega$ and tangent to~$\partial\Omega$ at~$q$, that is we take~$\zeta=(\zeta_1,\dots,\zeta_n)\in\Omega$ and~$r>0$ such that~$B_r(\zeta)\subseteq\Omega$ and
\begin{equation}\label{PJ-LKS-MPLEMQUIS2}
\{q\}=(\partial B_r(\zeta))\cap(\partial\Omega).\end{equation} Notice that~$\zeta_1=0$.
We let~$U:=B_r(\zeta)\cap B_{r/2}(q)\cap\{x_1>0\}$,
see Figure~\ref{J25372fzxfguicvbhjkopmjhynbgHY6MDSJORGFDFN8h5fdKSLE-21ItangeFIrfd7CORN2iSEpokjgHGSKDqwhrjfngpoewjnikdhgebidkkwdbbXVSdng}.

We point out that~$U\subseteq D$ and we define
$$ z(x)=z(x_1,\dots,x_n):=x_1\big(e^{-\alpha |x-\zeta|^2}-e^{-\alpha r^2}\big),$$
with~$\alpha>1$ to be conveniently chosen in the following calculation.
If~$x\in U$, then~$|x-\zeta|\ge|\zeta-q|-|q-x|\ge r-\frac{r}{2}=\frac{r}2$ and therefore
\begin{equation}\label{PJ-LKS-MPLEMQUISMA}
\Delta z(x)=2\alpha x_1 e^{-\alpha |x-\zeta|^2}\big(2\alpha|x-\zeta|^2-(n+2)\big)
\ge2\alpha x_1 e^{-\alpha |x-\zeta|^2}\left(\frac\alpha2-(n+2)\right)
\ge \alpha x_1 e^{-\alpha |x-\zeta|^2},
\end{equation}
so long as~$\alpha$ is fixed sufficiently large.

Additionally, we observe that~$z$ vanishes along~$\{x_1=0\}$, as well as along~$\partial B_r(\zeta)$.
We take~$\e\in(0,1)$, to be chosen appropriately small here below, and we define~$v:=w-\e z$.
{F}rom the previous observations we arrive at
\begin{equation}\label{PJ-LKS-MPLEMQUIS4}
{\mbox{along $\big(\Omega\cap B_{r/2}(q)\cap\{x_1=0\}\big)\cup\big((\partial B_r(\zeta))\cap B_{r/2}(q)\big)$ it holds that }}
v\ge-\e z=0.
\end{equation}
Moreover, we can suppose that~$w$ does not vanish identically and consequently, by the Strong Maximum Principle in
Theorem~\ref{STRONGMAXPLE1}(ii)
it holds that 
\begin{equation}\label{PJ-LKS-MPLEMQUIS}
w>0\;{\mbox{ in }}\;D.
\end{equation}
We claim that there exists~$c>0$ such that
\begin{equation}\label{PJ-LKS-MPLEMQUIS3}
w(x)\ge cx_1\;{\mbox{ for all }}\;x\in (\partial B_{r/2}(q))\cap B_r(\zeta)\cap\{x_1>0\}.
\end{equation}
Indeed, suppose not, then there exists a sequence of points~$x^{(j)}=
\big(x^{(j)}_1,\dots,x^{(j)}_n\big)\in(\partial B_{r/2}(q))\cap B_r(\zeta)\cap\{x_1>0\}$ such that~$w(x^{(j)})\le\frac{x_1^{(j)}}j$.
Up to a subsequence, we can suppose that~$x^{(j)}$ converges to some point~$\widetilde{x}\in
(\partial B_{r/2}(q))\cap \overline{B_r(\zeta)}\cap\{x_1\ge0\}$ (see Figure~\eqref{J25372fzxfguicvbhjkopmjhynbgHY6MDSJORGFDFN8h5fdKSLE-21ItangeFIrfd7CORN2iSEpokjgHGSKDqwhrjfngpoewjnikdhgebidkkwdbbXVSdng})
and~$w(\widetilde x)\le0$.
{F}rom this and~\eqref{PJ-LKS-MPLEMQUIS}, we find that~$
w(\widetilde x)=0$ and~$\widetilde x\in\{x_1=0\}\setminus(\partial\Omega)$. 

Hence, on the account of~\eqref{PJ-LKS-MPLEMQUIS} and the Hopf Lemma (see Lemma~\ref{JJS:PA}), we infer that~$\frac{\partial w}{\partial x_1}(\widetilde x)\ge c_0$,
for some~$c_0>0$.

As a result, if~$x$ lies in a sufficiently small neighborhood of~$\widetilde x$ we have that~$\frac{\partial w}{\partial x_1}(x)\ge \frac{c_0}2$ and consequently, if~$j$ is large enough,
$$ \frac{x_1^{(j)}}j\ge w(x^{(j)})=
w(x^{(j)})-w(0,(x^{(j)}_2,\dots,x^{(j)}_n)=
\int_0^{x_1^{(j)}}\frac{\partial w}{\partial x_1}(\tau,(x^{(j)}_2,\dots,x^{(j)}_n)\,d\tau\ge\frac{c_0 \,{x_1^{(j)}}}2.$$
We have therefore reached a contradiction which ends the proof of~\eqref{PJ-LKS-MPLEMQUIS3}.

\begin{figure}
  \centering
  \includegraphics[width=.55\linewidth]{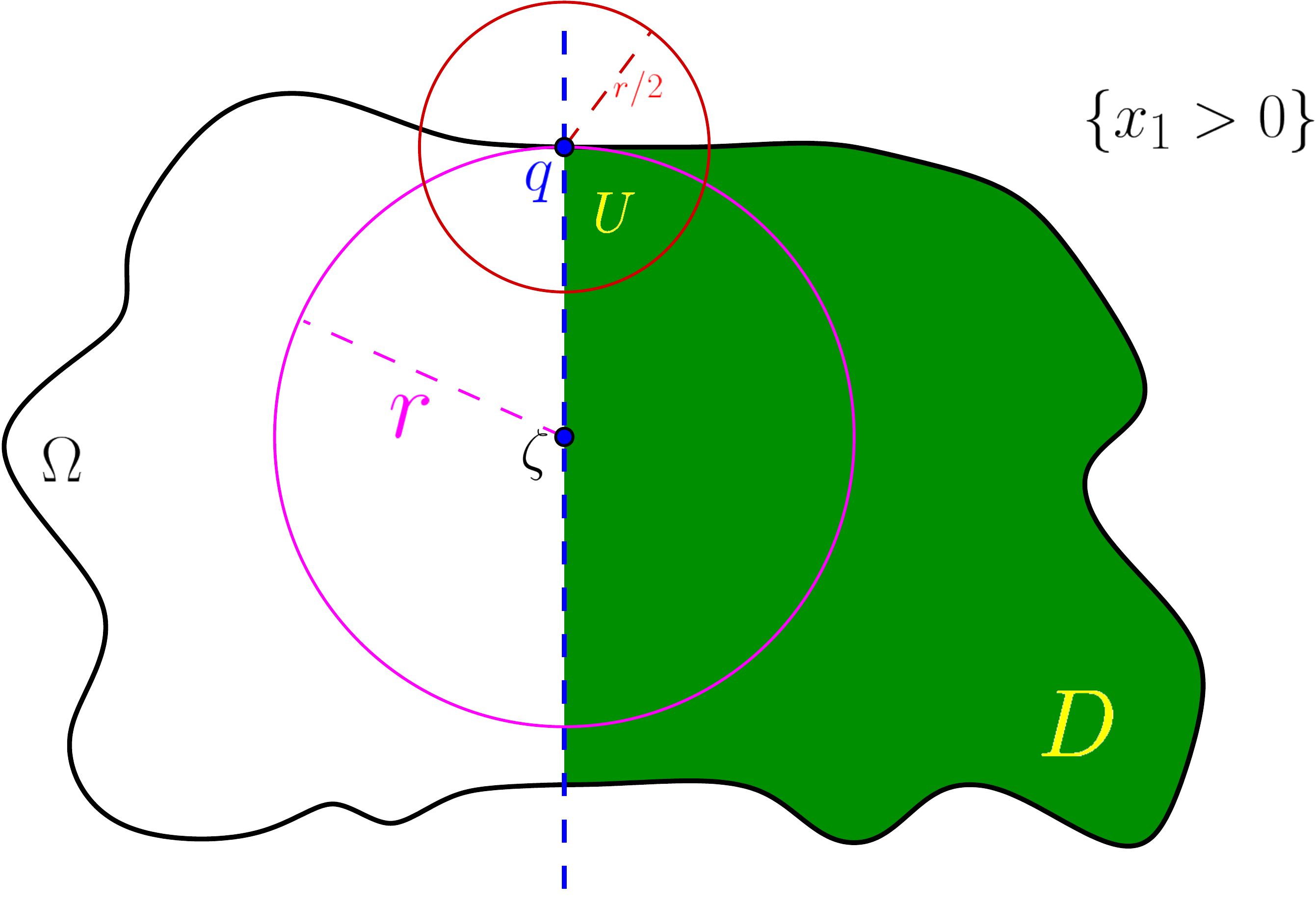}
 \caption{\sl The geometry involved in the proof of Lemma~\ref{SerrinsCornerLemma}.}\label{J25372fzxfguicvbhjkopmjhynbgHY6MDSJORGFDFN8h5fdKSLE-21ItangeFIrfd7CORN2iSEpokjgHGSKDqwhrjfngpoewjnikdhgebidkkwdbbXVSdng}
\end{figure}

{F}rom~\eqref{PJ-LKS-MPLEMQUIS3} we deduce that, for all~$x\in (\partial B_{r/2}(q))\cap B_r(\zeta)\cap\{x_1>0\}$,
we have
$$ v(x)\ge x_1\left(c-\e\sup_{X\in\Omega}
\Big(X_1\big(e^{-\alpha |X-\zeta|^2}-e^{-\alpha r^2}\big)\Big)\right)\ge0$$
if~$\e$ is sufficiently small.
Combining this observation and~\eqref{PJ-LKS-MPLEMQUIS4} we conclude that
\begin{equation}\label{PJ-LKS-MPLEMQUIS9}
v\ge0\;{\mbox{ along }}\;\partial U.
\end{equation}
Besides, by~\eqref{PJ-LKS-MPLEMQUISMA}, in~$U$ we have that~$\Delta v<0$.
Using this, \eqref{PJ-LKS-MPLEMQUIS9} and again the Strong Maximum Principle in
Theorem~\ref{STRONGMAXPLE1}(ii),
we infer that~$v>0$ in~$U$.

As a result, since~$q+t\varpi$ lies in~$D$, and hence in~$U$,
for~$t>0$ sufficiently small (recall~\eqref{FMDMIDLAFFJZimnPO83n}),
\begin{equation}\label{PJ-LKS-MPLEMQUIS9KSg904SA}
0<v(q+t\varpi)=v(q+t\varpi)-v(q)=t\partial_\varpi v(q)+\frac{t^2}2\partial_{\varpi\varpi}v(q)+o(t^2).\end{equation}

Now we distinguish two cases. If~$\partial_\varpi v(q)\ne0$, we write~\eqref{PJ-LKS-MPLEMQUIS9KSg904SA} in the form
$$ 0<t\partial_\varpi v(q)+o(t)$$
as~$t\searrow0$,
whence we deduce that
$$ 0<\partial_\varpi v(q)=\partial_\varpi w(q)-\e\partial_\varpi z(q)
=\partial_\varpi w(q)-\e\varpi_1\big(e^{-\alpha |q-\zeta|^2}-e^{-\alpha r^2}\big)=\partial_\varpi w(q),$$
which gives the desired result.

If instead~$\partial_\varpi v(q)=0$ then~\eqref{PJ-LKS-MPLEMQUIS9KSg904SA} boils down to
$$ 0<\partial_{\varpi\varpi}v(q)+o(1)$$
as~$t\searrow0$, and accordingly
\begin{eqnarray*}&&0\le\partial_{\varpi\varpi}v(q)
=\partial_{\varpi\varpi} w(q)-\e\partial_{\varpi\varpi} z(q)
=\partial_{\varpi\varpi} w(q)+
2\e\alpha \varpi_1 (q-\zeta)\cdot\varpi e^{-\alpha |q-\zeta|^2}\\&&\qquad\qquad
=\partial_{\varpi\varpi} w(q)+
2\e\alpha r\varpi_1\varpi_n e^{-\alpha r^2}
<\partial_{\varpi\varpi} w(q),\end{eqnarray*}
concluding the proof of the desired result.
\end{proof}

With this, we give a proof of Theorem~\ref{SERR} via the moving plane method, as in the original
article~\cite{MR333220}, by proceeding as follows.

\begin{proof}[Another proof of Theorem~\ref{SERR}] 
Recalling the notation in~\eqref{RIFLECONVE-eMPM3675869E242566THTERRAeFI2}
and~\eqref{800ntWAdw3G},
we claim that
\begin{equation}\label{SEHNp0jleSMPolnfWMojdblEDDPWD-R}
{\mbox{$\Omega$ is symmetric with respect to the reflection across~$T_{\lambda_0}$.}}
\end{equation}
Once this is proved, we are done: indeed, one would first repeat the argument in any given direction
and deduce that~$\Omega$ is invariant under hyperplane reflections in every direction;
from this, one can exploit Lemma~\ref{H3O7M8A5C2O21M1EF5K7A90654}
and conclude that~$\Omega$ is spherically symmetric, hence the union of disjoint rings;
thus applying Corollary~\ref{ssitingtotsendirow} one obtains that~$\Omega$ is a ball.

Therefore, we focus on the proof of~\eqref{SEHNp0jleSMPolnfWMojdblEDDPWD-R}.
To this end, we note that from~\eqref{HESE}
and the Weak Maximum Principle in Corollary~\ref{WEAKMAXPLE}(i) it follows that~$u\le0$ in~$\Omega$.

Now, let
$$ \Sigma_{\lambda_0}\ni x\longmapsto w_{\lambda_0}(x):=u(x)-u(x_{\lambda_0}).$$
We note that~$\Delta w_{\lambda_0}=1-1=0$ in~$\Sigma_{\lambda_0}$.
Furthermore, if~$x\in \Omega\cap T_{\lambda_0}$, it follows that~$x_{\lambda_0}=x$ and thus~$w_{\lambda_0}(x)=0$.
If instead~$x\in(\partial \Omega)\cap\{x_1>\lambda_0\}$ then~$w_{\lambda_0}(x)=-u(x_{\lambda_0})\ge0$.
These considerations and the Weak Maximum Principle in Corollary~\ref{WEAKMAXPLE}(ii) yield that~$w_{\lambda_0}\ge0$
in~$\Sigma_{\lambda_0}$.

This, in combination with the Strong Maximum Principle in Theorem~\ref{STRONGMAXPLE1}(ii), gives that
\begin{equation}\label{SEHNp0jleSMPolnfWMojdblEDDPWD}
{\mbox{either~$w_{\lambda_0}>0$
in~$\Sigma_{\lambda_0}$, or~$w_{\lambda_0}$
vanishes identically in~$\Sigma_{\lambda_0}$.}}
\end{equation}
The case in which~$w_{\lambda_0}$
vanishes identically in~$\Sigma_{\lambda_0}$ is easier to deal with,
hence we briefly discuss this situation now.
If~$w_{\lambda_0}$
vanishes identically in~$\Sigma_{\lambda_0}$
we have that~$u(x)=u(x_{\lambda_0})$
for all~$x\in\Sigma_{\lambda_0}$. In this circumstance, the claim in~\eqref{SEHNp0jleSMPolnfWMojdblEDDPWD-R}
must hold true, otherwise there would exist a point~$\overline x\in\partial\Sigma_{\lambda_0}\cap\{x_1>\lambda_0\}$ with~$\overline{x}_{\lambda_0}\in\Omega$, whence~$0=u(\overline x)=u(\overline{x}_{\lambda_0})$. That is, the subharmonic function~$u$ would
have an interior maximum at~$\overline{x}_{\lambda_0}$, and thus it must be constant, thanks to
the Strong Maximum Principle in Theorem~\ref{STRONGMAXPLE1}(i): but this would produce a contradiction with
the equation in~\eqref{HESE}.

{F}rom these observations and~\eqref{SEHNp0jleSMPolnfWMojdblEDDPWD} it follows that we can focus
on the case in which
\begin{equation}\label{SEHNp0jleSMPolnfWMojdblEDDPWD-087}
{\mbox{$w_{\lambda_0}>0$
in~$\Sigma_{\lambda_0}$.}}
\end{equation}
Actually, we will show that this case cannot occur at all, reaching a contradiction in what follows.
To understand this situation, we distinguish between
cases~(i) and~(ii) in~\eqref{IDUECASIMPS}. In all cases,
we notice that both~$p_0$ and its reflection across~$T_{\lambda_0}$ belong to~$\partial\Omega$
and therefore
\begin{equation}\label{092h9832hytjc9h98t0p-098765iut789y0tjjo2outjg7GS}
w_{\lambda_0}(p_0)=0-0=0.
\end{equation}
Now we show that, in the setting of~\eqref{SEHNp0jleSMPolnfWMojdblEDDPWD-087}, case~(i) cannot occur.
To check this we argue for a contradiction.
If case~(i) in~\eqref{IDUECASIMPS} holds true, we let~$p_0$ be the internal point of tangency
given in~\eqref{IDUECASIMPS}(i).
We point out that, in light of~\eqref{SEHNp0jleSMPolnfWMojdblEDDPWD-087}, we have that~$
w_{\lambda_0}<0$ in~$\Sigma_{\lambda_0}'$.
This, equation~\eqref{092h9832hytjc9h98t0p-098765iut789y0tjjo2outjg7GS}
and the version of the Hopf Lemma given in Lemma~\ref{JJS:PA} (applied here
to a ball~$B\subset\Sigma_{\lambda_0}'$ that is tangent to~$\partial \Sigma_{\lambda_0}'$ at the point~$p_0$)
entail that either~$ \partial_\nu w_{\lambda_0}(p_0)\ne0$ or~$w_{\lambda_0}$ is constantly zero in~$\Sigma_{\lambda_0}$.
But the latter situation cannot occur, thanks to~\eqref{SEHNp0jleSMPolnfWMojdblEDDPWD-087}, therefore necessarily
$$ 0\ne\partial_\nu w_{\lambda_0}(p_0)=c-c=0,$$
which is a contradiction (notice that we have used here the normal prescription along~$\partial\Omega$
given in~\eqref{HESE}).

This says that case~(i) of~\eqref{IDUECASIMPS}, in the setting of~\eqref{SEHNp0jleSMPolnfWMojdblEDDPWD-087}, does not occur.
We will show now that also case~(ii) of~\eqref{IDUECASIMPS}, in the setting
of~\eqref{SEHNp0jleSMPolnfWMojdblEDDPWD-087}, does not occur, and this will complete this proof
of Theorem~\ref{SERR}.

To this end, suppose that case~(ii) in~\eqref{IDUECASIMPS} holds true. In this scenario, our goal is to show that~$w_{\lambda_0}$
has a second order zero at~$p_0$ (this will allow us to employ Serrin's Corner Lemma
in Lemma~\ref{SerrinsCornerLemma}). Namely, we claim that
\begin{equation}\label{092h9832hytjc9h98t0p-098765iut789y0tjjo2outjg7GS-2}
D^k w_{\lambda_0}(p_0)=0\qquad{\mbox{for all }}\;k\in\{0,1,2\}.
\end{equation}
As a matter of fact, the case~$k=0$ is given in~\eqref{092h9832hytjc9h98t0p-098765iut789y0tjjo2outjg7GS},
hence we can focus here on the first and second derivatives.
To this end, we observe that, for all~$i$, $j\in\{1,\dots,n\}$,
$$ \partial_i w_{\lambda_0}(x)=\partial_i u(x)-(-1)^{\delta_{i1}}\partial_i u(x_{\lambda_0})$$
and
$$ \partial_{ij} w_{\lambda_0}(x)=\partial_{ij} u(x)-(-1)^{\delta_{i1}+\delta_{j1}}\partial_{ij} u(x_{\lambda_0}).$$
Thus, since~$p_0$ lies on~$T_{\lambda_0}$ (recall that we are dealing here
with case~(ii) in~\eqref{IDUECASIMPS}), and therefore the reflection of~$p_0$ coincides with~$p_0$ itself,
$$ \partial_i w_{\lambda_0}(p_0)=
\big(1-(-1)^{\delta_{i1}}\big)\partial_i u(p_0)=
\begin{dcases}
0 & {\mbox{ if }}i\in\{2,\dots,n\},\\
2\partial_1 u(p_0) & {\mbox{ if }}i=1,
\end{dcases}
$$
and
$$ \partial_{ij} w_{\lambda_0}(p_0)=\big(1-(-1)^{\delta_{i1}+\delta_{j1}}\big)\partial_{ij} u(p_0)=
\begin{dcases}
0 & {\mbox{ if }}i,j\in\{2,\dots,n\},\\
0& {\mbox{ if }}i=j=1,\\
2\partial_{1i} u(p_0) & {\mbox{ if $i\in\{2,\dots,n\}$ and~$j=1$,}}\\
2\partial_{1j} u(p_0) & {\mbox{ if $j\in\{2,\dots,n\}$ and~$i=1$.}}
\end{dcases}$$
Therefore, to prove~\eqref{092h9832hytjc9h98t0p-098765iut789y0tjjo2outjg7GS-2},
it remains to show that
\begin{equation}\label{092h9832hytjc9h98t0p-098765iut789y0tjjo2outjg7GS-2-90}
\partial_1 u(p_0)=0\qquad{\mbox{and}}\qquad\partial_{1i} u(p_0)=0\qquad{\mbox{for all }}\;
i\in\{2,\dots,n\}.
\end{equation}
To check these items of information, up to a rigid motion
we parameterize~$\partial\Omega$ in the vicinity of~$p_0$
by the graph of a function~$\phi:\{x'\in\R^{n-1}$ s.t. $|x'|<\rho\}\to\R$, for some~$\rho>0$. Given the assumption on~$\Omega$ in Theorem~\ref{SERR},
we know that~$\phi$ is of class~$C^2$. Since the normal derivative of~$\Omega$ at~$p_0$ is vertical,
we can suppose that the point~$(0,\phi(0))\in\R^{n-1}\times\R$ is~$p_0$ and~$\nabla\phi(0)=0$.

Since~$\Sigma_{\lambda_0}'\subseteq\Omega$ (up to reverting the vertical direction) we have that
$$\phi(x_1,x_2,\dots,x_{n-1})\le
\phi(-x_1,x_2,\dots,x_{n-1})$$ provided that~$x_1>0$ and~$|(x_1,x_2,\dots,x_{n-1})|<\rho$ and accordingly
\begin{eqnarray*}
0&\ge&\phi(x_1,x_2,\dots,x_{n-1})-\phi(-x_1,x_2,\dots,x_{n-1})\\
&=&\big(\phi(x_1,x_2,\dots,x_{n-1})-\phi(0)\big)-\big(\phi(-x_1,x_2,\dots,x_{n-1})-\phi(0)\big)\\
&=& \frac{1}{2} D^2\phi(0)(x_1,x_2,\dots,x_{n-1})\cdot(x_1,x_2,\dots,x_{n-1})\\&&\qquad
-\frac{1}{2} D^2\phi(0)(-x_1,x_2,\dots,x_{n-1})\cdot(-x_1,x_2,\dots,x_{n-1})+o(|x'|^2).
\end{eqnarray*}
In particular, given~$\ell\in\{2,\dots,n-1\}$, picking~$x_j:=0$ when~$j\not\in\{1,\ell\}$ we find that
\begin{eqnarray*}
0&\ge& \frac12\left[\partial_{11}\phi(0) x_1^2+\partial_{\ell\ell}\phi(0) x_\ell^2+2
\partial_{1\ell}\phi(0) x_1 x_\ell\right]\\&&\qquad-\frac12\left[\partial_{11}\phi(0) x_1^2+\partial_{\ell\ell}\phi(0) x_\ell^2-2
\partial_{1\ell}\phi(0) x_1 x_\ell\right]\\
&=&2\partial_{1\ell}\phi(0) x_1 x_\ell.
\end{eqnarray*}
Thus, choosing~$x_\ell:=\pm x_1$ and dividing by~$2x_1^2$,
$$ 0\ge\pm\partial_{1\ell}\phi(0).$$
Using the arbitrariness of choice in the above sign, we thus obtain that
\begin{equation}\label{LAMDS-KODNpoi1234ty5u67dsr98032hytg9u4jt-2902}
\partial_{1\ell}\phi(0)=0 \qquad{\mbox{for every }}\,\ell\in\{2,\dots,n-1\}.\end{equation}
In this configuration, we make use again of
the boundary prescriptions along~$\partial\Omega$
given in~\eqref{HESE}, thus writing that
$$ 0=u\big(x',\phi(x')\big)\qquad{\mbox{and}}\qquad
c=\partial_\nu u\big(x',\phi(x')\big)=
\nabla u\big(x',\phi(x')\big)\cdot\frac{\big( \nabla\phi(x'),-1\big)}{\sqrt{ 1+|\nabla\phi(x')|^2}}.$$
Differentiating these expressions and then setting~$x':=0$ we thus find that~$\partial_i u(p_0)=0$
for all~$i\in\{1,\dots,n-1\}$ (which, as an aside, establishes
the first claim in~\eqref{092h9832hytjc9h98t0p-098765iut789y0tjjo2outjg7GS-2-90},
and also gives that~$-\partial_n u(p_0)=\partial_\nu u(p_0)=c$) and that
\begin{eqnarray*}
0&=&\left.\frac{\partial}{\partial x_i}\left(
\nabla u\big(x',\phi(x')\big)\cdot\frac{\big( \nabla\phi(x'),-1\big)}{\sqrt{ 1+|\nabla\phi(x')|^2}}\right)\right|_{x'=0}\\
&=&
-\left.\frac{\partial}{\partial x_i}\Big(
\nabla u\big(x',\phi(x')\big)\Big)\right|_{x'=0}\cdot e_n
+
\nabla u(p_0)\cdot
\left.\frac{\partial}{\partial x_i}\left(
\frac{\big( \nabla\phi(x'),-1\big)}{\sqrt{ 1+|\nabla\phi(x')|^2}}\right)\right|_{x'=0}\\&=&
-\left.\frac{\partial}{\partial x_i}\Big(
\partial_n u\big(x',\phi(x')\big)\Big)\right|_{x'=0}
-\partial_n u(p_0)
\left.\frac{\partial}{\partial x_i}\left(
\frac{1}{\sqrt{ 1+|\nabla\phi(x')|^2}}\right)\right|_{x'=0}\\&=&
-\partial_{in} u(p_0).
\end{eqnarray*}
In particular, we have that~$\partial_{1n} u(p_0)=0$. For this reason,
to complete the proof of~\eqref{092h9832hytjc9h98t0p-098765iut789y0tjjo2outjg7GS-2-90},
and thus of~\eqref{092h9832hytjc9h98t0p-098765iut789y0tjjo2outjg7GS-2},
it only remains to show that
\begin{equation}\label{LAMDS-KODNpoir98032hytg9u4jt-290298i2oqejwrlfegbolewdshi}
\partial_{1i} u(p_0)=0\qquad{\mbox{for all }}\;
i\in\{2,\dots,n-1\}.\end{equation}
To this end, we observe that, for all~$i$, $j\in\{1,\dots,n-1\}$,
\begin{eqnarray*}
0&=&\partial_{ij}\Big(u\big(x',\phi(x')\big)\Big)\Big|_{x'=0}\\
&=&\partial_{i}\Big(
\partial_j u\big(x',\phi(x')\big)
+\partial_n u\big(x',\phi(x')\big)\partial_i\phi(x')\Big)\Big|_{x'=0}\\&=&
\partial_{ij} u(p_0)
+\partial_n u(p_0)\partial_{ij}\phi(0)\\
&=&\partial_{ij} u(p_0)
-c\partial_{ij}\phi(0).
\end{eqnarray*}
In particular, recalling~\eqref{LAMDS-KODNpoi1234ty5u67dsr98032hytg9u4jt-2902},
we have that~$\partial_{1i} u(p_0)=c\partial_{1i}\phi(0)=
0$ for all~$i\in\{2,\dots,n-1\}$ and the proof of~\eqref{LAMDS-KODNpoir98032hytg9u4jt-290298i2oqejwrlfegbolewdshi}, and hence of~\eqref{092h9832hytjc9h98t0p-098765iut789y0tjjo2outjg7GS-2}, is thus complete.

For our purposes, the result in~\eqref{092h9832hytjc9h98t0p-098765iut789y0tjjo2outjg7GS-2}
is instrumental in allowing us to use Serrin's Corner Lemma (in the version given here
in Lemma~\ref{SerrinsCornerLemma}). With this, we obtain that~$w_{\lambda_0}$ must vanish identically
in~$\Sigma_{\lambda_0}$, which provides a contradiction with~\eqref{SEHNp0jleSMPolnfWMojdblEDDPWD-087}.
\end{proof}\index{Serrin's overdetermined problem|)}

\section{The problem of Gidas, Ni and Nirenberg}\label{DEDXOJD-GIDANSNI}

The\index{Gidas, Ni and Nirenberg's problem|(} following beautiful result was discovered (or\footnote{Inventions and discoveries
have much in common, but some philosophers make a point in distinguishing
between a mathematical realisms, in which
mathematical entities are understood as existing independently of the human mind (in which case we do not invent mathematics, but rather discover it), and an opposed tendency
of mathematical anti-realism, holding that mathematical statements
do not correspond to a special entity (in which case we really invent mathematics, literally out of the blue).
Personally, we do not find such a distinction amazingly thrilling, but if we were forced to take part in this ontological
debate,
we would opt for ``discovered'', for the result,
and ``invented'', for the proof. Because, as Ennio De Giorgi wisely pointed out~\cite{MR1470169},
``a theorem is something discovered; its proof
something invented''.}
invented?) in~\cite{MR544879} (see
also~\cite{MR1662746}, \cite[Section~8.2]{MR2356201},
\cite[Theorem~2.34]{MR2777537},
\cite[Section~3.3]{MR2866937}
and~\cite[Sections~6.2 and~6.3]{MR2796831} for additional information):

\begin{theorem}\label{GININITH}
Let~$f\in C^1(\R)$ and~$u\in C^2(\overline{B_1})$ be a solution of
\begin{equation}\label{32254236equation-GNN-onwf}\begin{dcases}
\Delta u=f(u) & {\mbox{ in }}B_1,\\
u>0& {\mbox{ in }}B_1,\\
u=0 & {\mbox{ on }}\partial B_1.\end{dcases}\end{equation}
Then, $u$ is radially symmetric and radially strictly decreasing.
\end{theorem}

This is a very elegant statement,
especially comparing the very mild hypotheses assumed with the very strong result obtained.
The result is also of concrete interest, since solutions of the equation under consideration
appear in physics (see e.g.~\cite[equation~(12)]{MR436814},
in the context of the Yang-Mills field theory) and geometry (see e.g.~\cite{MR0358078},
in the search of Riemannian metrics that are covariantly connected under the group of M\"obius Transformations).
\medskip

The proof of Theorem~\ref{GININITH} relies on the reflection technique, but several specific
aspects in the application of the moving plane method arise.
For example, while the core of the idea is still to continuously slide a hyperplane to detect symmetry,
differently from the sketch in Figure~\ref{MPMETHTERRAeFI}, here
the critical situation does not occur when the reflected domain touches the original one but rather
when ``the reflected function  touches the original one'' (in a sense that will be made precise in the course of the proof).
Moreover, differently from the situation pointed out in Figure~\ref{MPME242566THTERRAeFI},
in the course of the forthcoming proof the reflection of the domain never becomes tangential to the original one, since
the domain under consideration is a ball and therefore the sliding hyperplanes do not cut the domain perpendicularly
before reaching the center of the ball.\medskip

Roughly speaking, the use of the moving plane method
in the proof of Theorem~\ref{GININITH} consists in the following conceptual steps. First off, one considers the difference between the solution~$u$ and its reflection: since~$u$ vanishes along
the boundary and it is positive elsewhere, when the hyperplane is very close
to the initial position one is able to detect a sign for such a difference.
Then, one extends the validity of this sign up to a critical position.
It thus remains to analyze this critical position in order to detect
the symmetry of the solution.
Several forms of the Maximum Principle will come in handy in
the above procedure, since the difference between the original
solution and its reflection ends up satisfying a nice partial differential equation.\medskip

Some further remarks about Theorem~\ref{GININITH}:
first of all, when~$f$ is only H\"older continuous, the result in Theorem~\ref{GININITH} does not
hold anymore, since in this case solutions are not necessarily radial, see e.g.~\cite[Exercise~3.18]{MR1751289}.
Also, when the domain is an annulus (instead of the ball), Theorem~\ref{GININITH} does not hold,
see e.g.~\cite{MR709644, MR1167503} for the construction of nonradial positive solutions (but see also~\cite[Theorem 8.2.4]{MR2356201} for a radial symmetry result in the annulus under suitable structural assumptions on the nonlinearity~$f$).
These observations showcase the rather subtle nature of the radial symmetry results
and the several patterns which can lead to symmetry breaking after apparently small modifications
of the problem under consideration.\medskip

To implement the proof of Theorem~\ref{GININITH} we start by presenting some auxiliary results.

\begin{lemma}\label{JMDSJORGFDFN8h5fdKSLE-21}
Let~$\Omega$ be a bounded open set with boundary of class~$C^{2}$, with exterior normal~$\nu$. 

Let~$x_0\in\partial\Omega$, and suppose that
\begin{equation}\label{BUJMGKCH0KDSON01}\nu(x_0)\cdot e_1>0.\end{equation}

For all~$r>0$, let $\Omega_r:=\Omega\cap B_r(x_0)$.

Let~$\e>0$, $f\in C^1(\R)$ and~$u\in C^2(\overline{\Omega_\e})$ be such that~$\Delta u =f(u)$ in~$\Omega_\e$,
$u>0$ in~$\Omega_\e$ and~$u=0$ on~$(\partial\Omega)\cap B_\e(x_0)$.

Then, there exists~$\delta\in(0,\e)$ such that~$\partial_1u(x)<0$ for all~$x\in\Omega_\delta$.
\end{lemma}

\begin{figure}
  \centering
  \includegraphics[width=.45\linewidth]{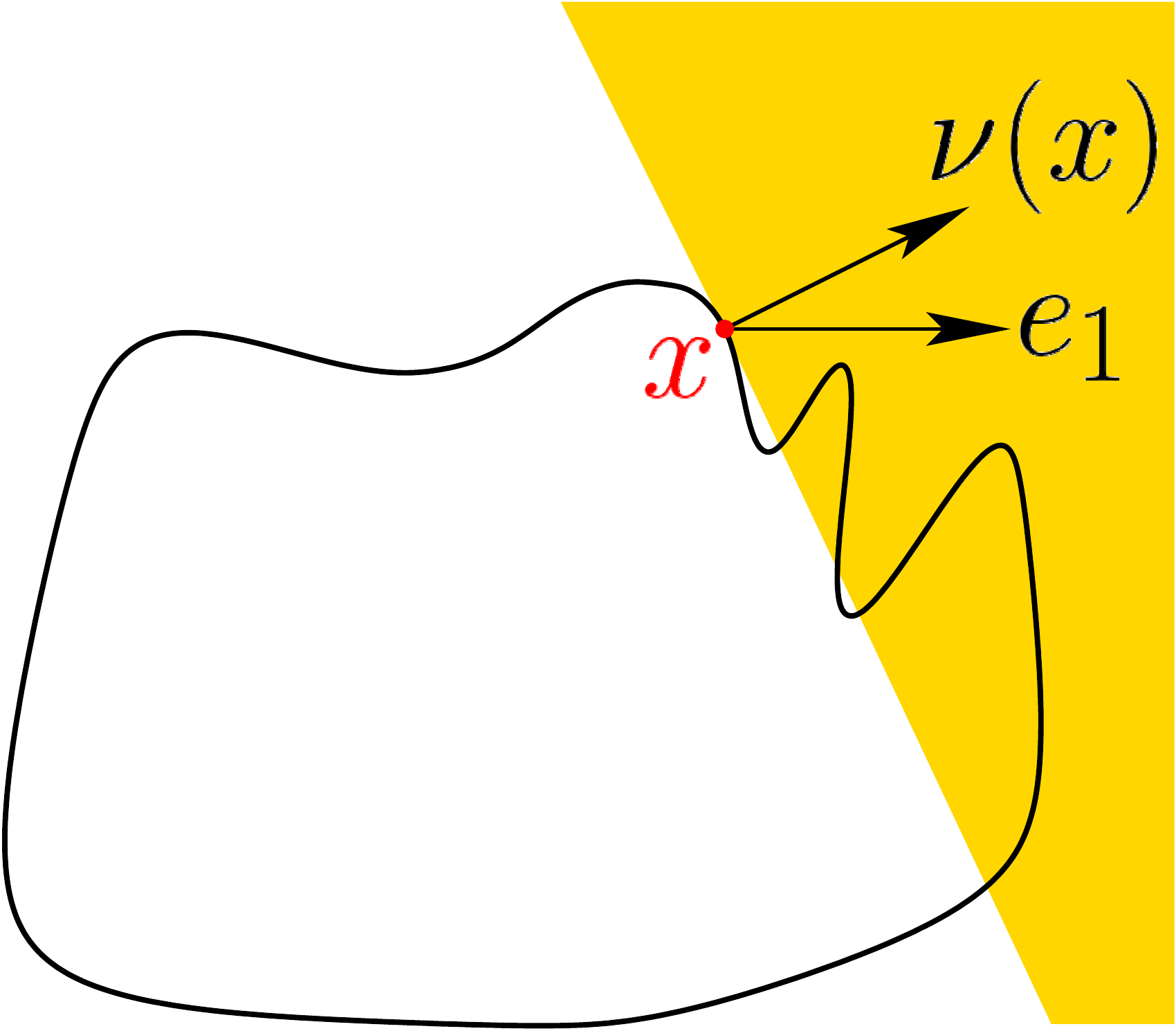}
 \caption{\sl The geometry involved in the proof of Lemma~\ref{JMDSJORGFDFN8h5fdKSLE-21}.}\label{JMDSJORGFDFN8h5fdKSLE-21ItangeFIrfd72iqwhrjfngpoewjnikdhgebidkkwdbbXVSdng}
\end{figure}

\begin{proof} Let~$S_r:=(\partial\Omega)\cap B_r(x_0)$.
In light of~\eqref{BUJMGKCH0KDSON01},
we can find~$\e'\in(0,\e)$ small enough such that~$\nu(x)\cdot e_1>0$
for all~$x\in S_{\e'}$. 
As a consequence (see Figure~\ref{JMDSJORGFDFN8h5fdKSLE-21ItangeFIrfd72iqwhrjfngpoewjnikdhgebidkkwdbbXVSdng})
for all~$x\in S_{\e'}$ there exists~$h_x>0$ such that for all~$h\in(0,h_x)$
we have that~$x-he_1\in\Omega_\e$ and accordingly
\begin{equation} \label{WDFONBAMASSVPLDLMLHISKNB}\partial_1 u(x)=\lim_{h\searrow0}\frac{u(x)-u(x-he_1)}{h}=\lim_{h\searrow0}\frac{0-u(x-he_1)}{h}\le0.
\end{equation}
Now, to prove the desired result, we argue by contradiction and we suppose
that for every~$j\in\N$ there exists a point~$x^{(j)}\in\Omega_{1/j}$
with
\begin{equation}\label{WDFONBAMASSVPLDLMLHISKNB2}\partial_1u(x^{(j)})\ge0.\end{equation}
Notice that
$$ \lim_{j\to+\infty}{\rm dist}(x^{(j)},\partial\Omega)={\rm dist}(x_0,\partial\Omega)=0,$$
therefore we can find~$t^{(j)}>0$ such that~$x^{(j)}+t^{(j)}e_1\in\partial\Omega$, with~$t^{(j)}\to0$ as~$j\to +\infty$.
In particular,
$$ \lim_{j\to+\infty}\big|x^{(j)}+t^{(j)}e_1-x_0\big|=0$$
and accordingly, when~$j$ is sufficiently large, we have that~$x^{(j)}+t^{(j)}e_1\in S_{\e'}$.

{F}rom this and~\eqref{WDFONBAMASSVPLDLMLHISKNB} we arrive at
$$ \partial_1 u\big(x^{(j)}+t^{(j)}e_1\big)\le0.$$
Comparing with~\eqref{WDFONBAMASSVPLDLMLHISKNB2}, we thereby find~$\tau^{(j)}\in\big[0,t^{(j)}\big]$ such that
$$ \partial_1 u\big(x^{(j)}+\tau^{(j)}e_1\big)=0.$$
Taking the limit as~$j\to+\infty$, we conclude that
\begin{equation}\label{CBSAMCTUNMLPLMJH-0ol434} \partial_1 u(x_0)=0.\end{equation}
Moreover,
\begin{equation}\label{CBSAMCTUNMLPLMJH-0ol43435} \partial_{11}u(x_0)=\lim_{j\to+\infty}\frac{\partial_1 u\big(x^{(j)}+\tau^{(j)}e_1\big)
-\partial_1 u(x_0)}{x^{(j)}+\tau^{(j)}e_1-x_0}=\lim_{j\to+\infty}\frac{0-0}{x^{(j)}+\tau^{(j)}e_1-x_0}=0.
\end{equation}
Now let~$\{\zeta_1,\dots,\zeta_{n-1}\}$ be an orthonormal basis for the tangent hyperplane to~$\partial\Omega$ at~$x_0$.
Since~$u=0$ on~$(\partial\Omega)\cap B_\e(x_0)$ we have that
\begin{equation}\label{CBSAMCTUNMLPLMJH-0ol431}
\nabla u(x_0)\cdot\zeta_i=0\qquad{\mbox{for every }}\,i\in\{1,\dots,n-1\}.
\end{equation} 
We also notice that
\begin{equation}\label{CBSAMCTUNMLPLMJH-0ol433}
{\mbox{$e_1$ is linearly independent from~$\zeta_1,\dots,\zeta_{n-1}$,}}\end{equation} because if,
for some~$a_1,\dots,a_n\in\R$,
\begin{equation}\label{CBSAMCTUNMLPLMJH-0ol432}
a_1\zeta_1+\dots+a_{n-1}\zeta_{n-1}+a_n e_1=0\end{equation}
then
$$ 0=\nu(x_0)\cdot\left(a_1\zeta_1+\dots+a_{n-1}\zeta_{n-1}+a_n e_1
\right)=a_n\nu(x_0)\cdot e_1,$$
whence one would deduce, owing to~\eqref{BUJMGKCH0KDSON01}, that~$a_n=0$.
Plugging this information back into~\eqref{CBSAMCTUNMLPLMJH-0ol432}, we obtain that~$a_1\zeta_1+\dots+a_{n-1}\zeta_{n-1}=0$ and hence, by the independence of the basis, necessarily~$a_1=\dots=a_{n-1}=0$.
This proves~\eqref{CBSAMCTUNMLPLMJH-0ol433}.

Hence, by~\eqref{CBSAMCTUNMLPLMJH-0ol434}, \eqref{CBSAMCTUNMLPLMJH-0ol431} and~\eqref{CBSAMCTUNMLPLMJH-0ol433} we conclude that
\begin{equation}\label{CBSAMCTUNMLPLMJH-0ol43X}
\nabla u(x_0)=0.
\end{equation}
Now, for all tangent vectors~$\zeta$ we let~$\gamma_\zeta\in C^2((-1,1),\partial\Omega\cap B_\e(x_0))$ 
with~$\gamma_\zeta(0)=x_0$ and~$\dot\gamma_\zeta(0)=\zeta$.
We have that~$u(\gamma_\zeta(t))=0$ and therefore, by taking derivatives,
$$ 0=\sum_{k=1}^n \partial_k u(\gamma_\zeta(t))\,\dot\gamma_\zeta(t)\cdot e_k$$
and $$
0=\sum_{j,k=1}^n \partial_{jk} u(\gamma_\zeta(t))\,(\dot\gamma_\zeta(t)\cdot e_j)\,(\dot\gamma_\zeta(t)\cdot e_k)
+\sum_{k=1}^n \partial_k u(\gamma_\zeta(t))\,\ddot\gamma_\zeta(t)\cdot e_k.$$
Combining this information with~\eqref{CBSAMCTUNMLPLMJH-0ol43X} we infer that
\begin{equation*}
0=\sum_{j,k=1}^n \partial_{jk} u(x_0)\,(\zeta\cdot e_j) (\zeta\cdot e_k)=
\sum_{j,k=1}^n \partial_{jk} u(x_0)\,(\zeta\cdot e_j) (\zeta\cdot e_k)
=D^2u(x_0)\zeta\cdot\zeta.\end{equation*}
This is valid for all tangent vectors~$\zeta$. In particular, applying this
with~$\zeta:=\zeta_i-\zeta_j$, $\zeta:=\zeta_i$ and~$\zeta:=\zeta_j$,
\begin{equation}\label{90oihwqgoujbfoiu3fbliuy32fr8ogfywal6bc56875tbx-294utujgeHJS-3}
\begin{split}&0=D^2u(x_0)(\zeta_i-\zeta_j)\cdot(\zeta_i-\zeta_j)=
D^2u(x_0)\zeta_i \cdot\zeta_i+D^2u(x_0)\zeta_j\cdot\zeta_j-2
D^2u(x_0)\zeta_i\cdot\zeta_j\\&\qquad\qquad\qquad=-
2D^2u(x_0)\zeta_i\cdot\zeta_j.\end{split}\end{equation}
%%%%  Now, letting~$b_i:=\zeta_i\cdot e_1$ for all~$i\in\{1,\dots,n-1\}$ and~$b_n:=\nu(x_0)\cdot e_1>0$, 
%%%%  since~$\{\zeta_1,\dots,\zeta_{n-1},\nu(x_0)\}$ is an orthonormal frame of~$\R^n$ we have that
%%%%  $$ e_1=\sum_{i=1}^{n-1} b_i\zeta_i+b_n \nu(x_0)$$
%%%%  and thus, by~\eqref{CBSAMCTUNMLPLMJH-0ol43435}
%%%%  and~\eqref{90oihwqgoujbfoiu3fbliuy32fr8ogfywal6bc56875tbx-294utujgeHJS-3},
%%%%  \begin{eqnarray*} &&0=\partial_{11}u(x_0)=D^2u(x_0)e_1\cdot e_1=
%%%%  \sum_{i,j=1}^{n-1} b_i b_j D^2u(x_0)\zeta_i\cdot \zeta_j
%%%%  +2\sum_{i=1}^{n-1} b_i b_n D^2u(x_0)\zeta_i\cdot \nu(x_0)
%%%%  +b_n^2\partial_{\nu\nu}u(x_0)\\&&\qquad\qquad\qquad=
%%%%  2\sum_{i=1}^{n-1} b_i b_n D^2u(x_0)\zeta_i\cdot \nu(x_0)
%%%%  +b_n^2\partial_{\nu\nu}u(x_0).\end{eqnarray*}
Hence, by the rotational invariance of the Laplacian (recall Corollary~\ref{KSMD:ROTAGSZOKA})
and~\eqref{90oihwqgoujbfoiu3fbliuy32fr8ogfywal6bc56875tbx-294utujgeHJS-3},
\begin{equation}\label{CBSAMCTUNMLPLMJH-0ol43XANP} f(0)=f(u(x_0))=\Delta u(x_0)=
\sum_{i=1}^{n-1}D^2u(x_0)\zeta_i\cdot\zeta_i+\partial_{\nu\nu}u(x_0)=\partial_{\nu\nu}u(x_0).\end{equation}
%%%%   As a result,
%%%%   $$0=2\sum_{i=1}^{n-1} b_i b_n D^2u(x_0)\zeta_i\cdot \nu(x_0)
%%%%   +b_n^2 f(0).$$

Now we show that %%, for all~$i\in\{1,\dots,n-1\}$,
\begin{equation}\label{pa24ftrw21e2542s35met3li2r3i88tr2154on12i}
%% D^2u(x_0)\zeta_i\cdot\nu(x_0)=0
f(0)=0.
\end{equation}
To this end, we let~$d$ be the signed
distance function to~$\partial\Omega$, with the convention that~$d>0$ in~$\Omega$ and~$d<0$ in~$\R^n\setminus\Omega$
(see the discussion on page~\pageref{PROIE2}
about the signed distance function and in particular equation~\eqref{PROIE2}).
We let~$t>0$, to be taken suitably small, and~$s\in[-t,t]$. Then, for small~$t$,
\begin{eqnarray*}&& d\big(x_0-t\nu(x_0)+s\zeta_i\big)=
d\big(x_0-t\nu(x_0)\big)+s\nabla d\big(x_0-t\nu(x_0)\big)\cdot \zeta_i+O(s^2)\\&&\quad=
t-s\nu(x_0)\cdot \zeta_i+O(s^2)=t+O(s^2)=t+O(t^2)>0.\end{eqnarray*}
That is, we have that~$x_0-t\nu(x_0)+s\zeta_i\in\Omega$ and accordingly,
by~\eqref{CBSAMCTUNMLPLMJH-0ol43X},
\begin{eqnarray*} 0&<&u\big(x_0-t\nu(x_0)+s\zeta_i\big)\\&=&
u(x_0)+\nabla u(x_0)\cdot\big(-t\nu(x_0)+s\zeta_i\big)
+\frac12 D^2u(x_0)\big(-t\nu(x_0)+s\zeta_i\big)\cdot\big(-t\nu(x_0)+s\zeta_i\big)
+o(t^2)\\
&=&
0+0+\frac{t^2}2 \partial_{\nu\nu}u(x_0)+\frac{s^2}2 \partial_{\zeta_i\zeta_i}u(x_0)
-st D^2u(x_0)\zeta_i\cdot \nu(x_0)+o(t^2).
\end{eqnarray*}
Hence, in view of~\eqref{CBSAMCTUNMLPLMJH-0ol43XANP},
$$ 0=\frac{t^2 f(0)}2-st D^2u(x_0)\zeta_i\cdot \nu(x_0)+o(t^2)$$
and therefore, choosing~$s:=\pm t$ and dividing by~$t^2$,
$$ \frac{f(0)}2\mp D^2u(x_0)\zeta_i\cdot \nu(x_0)=o(1).$$
By taking the limit as~$t\searrow0$, we conclude that
$$ \frac{f(0)}2=\pm D^2u(x_0)\zeta_i\cdot \nu(x_0).$$
The possibility of choosing the sign in the above right hand side gives~\eqref{pa24ftrw21e2542s35met3li2r3i88tr2154on12i},
as desired.

{F}rom~\eqref{pa24ftrw21e2542s35met3li2r3i88tr2154on12i}, letting
$$ c(x):=\int_0^{u(x)}f'(t)\,dt,$$
we find that, in~$\Omega_\e$,
$$ \Delta u(x)=f(u(x))=f(u(x))-f(0)=c(x)u(x).$$
We can therefore invoke the version of the Hopf Lemma given in Corollary~\ref{MS-LEHOPF2-M-NEWCORO}
(applied here to~$v:=-u$)
and deduce that~$\partial_\nu u(x_0)\ne0$.
But this is in contradiction with~\eqref{CBSAMCTUNMLPLMJH-0ol43X}
and the desired result is thus established.
\end{proof}

\begin{lemma}\label{JNSMPERADVBDDOBYHGGJUINMRARYS3123}
Let~$f\in C^1(\R)$ and~$u\in C^2(\overline{B_1})$ be a solution of
$$\begin{dcases}
\Delta u=f(u) & {\mbox{ in }}B_1,\\
u>0& {\mbox{ in }}B_1,\\
u=0 & {\mbox{ on }}\partial B_1.\end{dcases}$$

Then,
\begin{equation}\label{LAPIR-HNG6JEnIERNsTrnE4gt-10431}
{\mbox{$u(-x_1,x_2,\dots,x_n)=u(x_1,x_2,\dots,x_n)$\; for all~$x=(x_1,x_2,\dots,x_n)\in B_1$.}}\end{equation}

Moreover, \begin{equation}\label{LAPIR-HNG6JEnIERNsTrnE4gt-10432}
{\mbox{$\partial_1u(x)<0\,$ for all~$x=(x_1,x_2,\dots,x_n)\in B_1$ with~$x_1>0$.}}\end{equation}
\end{lemma}

\begin{proof}
The moving plane notation in~\eqref{RIFLECONVE-eMPM3675869E242566THTERRAeFI2}
will be adopted here with~$\Omega:=B_1$.
For all~$\lambda\in(0,1)$, we define
\begin{equation}\label{GNNSHLIMSBXNSKt6-PRE-SEGBXDDUEv4JKS} \Sigma_\lambda\ni x\longmapsto w_\lambda(x):=u(x)-u(x_\lambda).\end{equation}
We note that if~$x\in T_\lambda$ then~$x_\lambda=x$ and consequently~$w_\lambda(x)=0$.
If instead~$x\in (\partial\Sigma_\lambda)\cap\{x_1>\lambda\}\subseteq\partial B_1$
we have that~$w_\lambda(x)=0-u(x_\lambda)<0$. {F}rom these observations we arrive at
\begin{equation}\label{GNNSHLIMSBXNSKt6-PRE}
{\mbox{$w_\lambda\le0$ on~$\Sigma_\lambda$, but~$w_\lambda$ does not vanish identically.}}
\end{equation}

Moreover,
we have that~$\partial_1 w_\lambda(x)=\partial_1u(x)+\partial_1u(2\lambda-x_1,x_2,\dots,x_n)$, and thus
\begin{equation}\label{equa0urjoewh3r73855-d}
{\mbox{$\partial_1 w_\lambda=2\partial_1 u$ on~$B_1\cap T_\lambda$.}}\end{equation}
We now let
$$ c_\lambda(x):=\int_0^1 f'\big(tu(x)+(1-t)u_\lambda(x) \big)\,dt $$
and we point out that
\begin{equation}\label{equa0urjoewh3r73855-d-fingivis}\|c_\lambda\|_{L^\infty(B_1)}
\le\|f'\|_{L^\infty([0,3\sup_{B_1}u])}.\end{equation}
Also, for all~$x\in\Sigma_\lambda$,
\begin{equation}\label{GNNSHLIMSBXNSKt6-PRE2}
\Delta w_\lambda(x)=f(u(x))-f(u(x_\lambda))=c_\lambda(x)\big(
u(x)-u(x_\lambda)\big)=c_\lambda(x)\,w_\lambda(x).
\end{equation}

We claim that
\begin{equation}\label{GNNSHLIMSBXNSKt6}
{\mbox{for all~$\lambda\in(0,1)$ and all~$x\in\Sigma_\lambda$ we have that~$w_\lambda(x)<0$.}}
\end{equation}
To this end, we observe that when~$\e_0\in(0,1)$ is conveniently small and~$\lambda\in(1-\e_0,1)$,
we can utilize the Maximum Principle for narrow domains given in
Theorem~\ref{M:PLE:NARROW} (or, alternatively, the
Maximum Principle for small volume domains given in Theorem~\ref{M:PLE:SMALLV})
in combination with~\eqref{GNNSHLIMSBXNSKt6-PRE} and~\eqref{GNNSHLIMSBXNSKt6-PRE2}: this yields that~$w_\lambda<0$ in~$\Sigma_\lambda$ for all~$\lambda\in(1-\e_0,1)$.

We can therefore consider\footnote{In a sense, the definition in~\eqref{800ntWAdw3}
is an analytic counterpart of~\eqref{800ntWAdw3G} in which instead the critical situation arose from
a geometric situation.}
\begin{equation}\label{800ntWAdw3} \lambda_0:=\inf\Big\{\lambda\in(0,1) {\mbox{ s.t. $w_\lambda<0$ 
in $\Sigma_\lambda$}}
\Big\}\end{equation}
and we know that~$\lambda_0\le1-\e_0$.

We observe that~\eqref{GNNSHLIMSBXNSKt6} is proved once we show that~$\lambda_0=0$.
Hence, we argue by contradiction and we suppose that~$\lambda_0>0$.
We recall that, by construction, 
\begin{equation}\label{DATUMWODKDIWPD}
{\mbox{for all~$\lambda\in[\lambda_0,1)$,
$w_{\lambda}\le0$ in~$\Sigma_{\lambda}$
and~$w_{\lambda}=0$ on~$T_{\lambda}$.}}\end{equation}
Moreover,
\begin{equation}\label{CA7CA}
{\mbox{for all~$\lambda\in[\lambda_0,1)$, $w_{\lambda}<0$ in~$(\partial B_1)
\cap\{x_1>\lambda\}$,}}
\end{equation}
since if~$p\in(\partial B_1)
\cap\{x_1>\lambda\}$,
then~$p_\lambda\in B_1$ (since~$\lambda\ge\lambda_0>0$)
and accordingly~$w_{\lambda}(p)=-u(p_{\lambda})<0$, thus proving~\eqref{CA7CA}.

Hence, by~\eqref{GNNSHLIMSBXNSKt6-PRE2},
\eqref{DATUMWODKDIWPD}, \eqref{CA7CA}
and the version of Hopf Lemma given in Corollary~\ref{MS-LEHOPF2-M-NEWCORO},
for all~$\lambda\in[\lambda_0,1)$ either
\begin{equation}\label{OHSn0ujfiHDiegb-0ujhr-00}
{\mbox{$\partial_1w_{\lambda}<0$
along~$B_1\cap T_{\lambda}$}}\end{equation} or~$w_{\lambda}$ vanishes identically.
But the latter possibility cannot hold, owing to~\eqref{CA7CA}.
These considerations yield that necessarily~\eqref{OHSn0ujfiHDiegb-0ujhr-00}
holds true for all~$\lambda\in[\lambda_0,1)$.

In the same way, using the Maximum Principle in Theorem~\ref{MAXPLECON45678CG-PM}, we find that
\begin{equation}\label{OHSn0ujfiHDiegb-0ujhr-00CASILIM}
w_{\lambda_0}<0\, {\mbox{ in~$\Sigma_{\lambda_0}$.}}\end{equation}

{F}rom~\eqref{equa0urjoewh3r73855-d} and~\eqref{OHSn0ujfiHDiegb-0ujhr-00},
it follows that
\begin{equation}\label{OHSn0ujfiHDiegb-0ujhr}
{\mbox{for all~$\lambda\in[\lambda_0,1)$, $\partial_1u<0$
along~$B_1\cap T_{\lambda}$}}.\end{equation}
%%%and accordingly
%%%\begin{equation}\label{OHSn0ujfiHDiegb-0ujhr-tutta}
%%%{\mbox{$\partial_1u<0$
%%%in~$\Sigma_{\lambda_0}$.}}\end{equation}

Now we show that
\begin{equation}\label{OHSn0ujfiHDiegb-0ujhrST}
{\mbox{$\partial_1u<0$
along~$\overline{B_1}\cap T_{\lambda_0}$.}}\end{equation}
Indeed, by~\eqref{OHSn0ujfiHDiegb-0ujhr} (applied with~$\lambda:=\lambda_0$), to check~\eqref{OHSn0ujfiHDiegb-0ujhrST}
we only need to take~$x_0\in(\partial B_1)\cap T_{\lambda_0}$ and prove that~$\partial_1u(x_0)<0$.
But this is a consequence of Lemma~\ref{JMDSJORGFDFN8h5fdKSLE-21} (which can be used here since~$\nu(x_0)\cdot e_1=
x_0\cdot e_1=\lambda_0>0$) and therefore~\eqref{OHSn0ujfiHDiegb-0ujhrST} holds true.

%%%%  By~\eqref{OHSn0ujfiHDiegb-0ujhrST}, there exists~$\delta_0>0$ such that
%%%%  \begin{equation*}
%%%%  {\mbox{$\partial_1 u<0$
%%%%  in~$B_1\cap\big\{x_1\in(\lambda_0-\delta_0,\lambda_0]\big\}$}}.\end{equation*}
%%%%  {F}rom this and~\eqref{OHSn0ujfiHDiegb-0ujhr-tutta} we infer that
%%%%  \begin{equation}\label{OHSn0ujfiHDiegb-0ujhr-2}
%%%%  {\mbox{$\partial_1 u<0$
%%%%  in~$B_1\cap\big\{x_1>\lambda_0-\delta_0\big\}$}}.\end{equation}

Moreover, it follows from the definition of~$\lambda_0$ in~\eqref{800ntWAdw3} that there exist
a sequence of parameters~$\lambda^{(j)}<\lambda_0$ such that~$\lambda^{(j)}\nearrow\lambda_0$ as~$j\to+\infty$ and
a sequence of points~$x^{(j)}\in \Sigma_{\lambda^{(j)}}$ such that~$w_{\lambda^{(j)}}(x^{(j)})\ge0$.
Up to a subsequence, we suppose that~$x^{(j)}$ converges to some point~$x^\star$ as~$j\to+\infty$,
and necessarily, given the convergence of~$\lambda^{(j)}$, we have that~$x^\star\in\overline{\Sigma_{\lambda_0}}$.

Using the notation~$x^{(j)}=\big(x^{(j)}_1,y^{(j)}\big)$ with~$y^{(j)}:=\big(x^{(j)}_2,\dots,x^{(j)}_n\big)$, we find that
\begin{eqnarray*}&& 0\le\lim_{j\to+\infty} w_{\lambda^{(j)}}(x^{(j)})=
\lim_{j\to+\infty} u\big(x^{(j)}_1,y^{(j)}\big)-u\big(2\lambda^{(j)}-x^{(j)}_1,y^{(j)}\big)\\&&\qquad\qquad
=u(x^\star)-u(2\lambda_0-x^\star_1,x^\star_2,\dots,x^\star_n)=w_{\lambda_0}(x^\star).\end{eqnarray*}
Since the opposite inequality holds in~$\overline{\Sigma_{\lambda_0}}$, thanks to the definition of~$\lambda_0$
in~\eqref{800ntWAdw3}, we have that
\begin{equation}\label{OJLDPLOmfSKheann9emaoRK1} w_{\lambda_0}(x^\star)=0.\end{equation}
By comparing~\eqref{CA7CA} and~\eqref{OJLDPLOmfSKheann9emaoRK1} we deduce that
\begin{equation}\label{OJLDPLOmfSKheann9emaoRK}
x^\star\not\in\partial B_1.
\end{equation}
As a consequence of~\eqref{OJLDPLOmfSKheann9emaoRK}, we have that
\begin{equation}\label{OJLDPLOmfSKheann9emaoRK-ZZ}
x^\star\in \big( {B_1}\cap T_{\lambda_0}\big)\cup\Sigma_{\lambda_0}.
\end{equation}
But~$x^\star\not\in\Sigma_{\lambda_0}$,
due to~\eqref{OHSn0ujfiHDiegb-0ujhr-00CASILIM} and~\eqref{OJLDPLOmfSKheann9emaoRK1}.
Therefore we deduce from~\eqref{OJLDPLOmfSKheann9emaoRK-ZZ}
that
\begin{equation}\label{OJLDPLOmfSKheann9emaoRK-ZZ-lAB}
x^\star\in {B_1}\cap T_{\lambda_0}.\end{equation}

We also remark that
$$ u\big(x^{(j)}_1,y^{(j)}\big)=u(x^{(j)})=w_{\lambda^{(j)}}(x^{(j)})+u(x^{(j)}_{\lambda^{(j)}})\ge0+u(x^{(j)}_{\lambda^{(j)}})=u
\big(2\lambda^{(j)} -x^{(j)}_1,y^{(j)}\big).$$
Therefore, by
the Calculus Mean Value Theorem,
there must exist~$t^{(j)}\in\left[2\lambda^{(j)} -x^{(j)}_1,\,x^{(j)}_1\right]$ such that
\begin{equation}\label{KDS-PASTTLQEQCSPE}
\partial_1u
\big(t^{(j)},y^{(j)}\big)\ge0.\end{equation}
We stress that, in light of~\eqref{OJLDPLOmfSKheann9emaoRK-ZZ-lAB},
$$ \lim_{j\to+\infty}x^{(j)}_1=\lambda_0,\qquad{\mbox{hence}}\qquad
\lim_{j\to+\infty}t^{(j)}=\lambda_0.$$
We can therefore take the limit as~$j\to+\infty$ in~\eqref{KDS-PASTTLQEQCSPE} and conclude that~$\partial_1u(x_\star)\ge0$.
This and~\eqref{OJLDPLOmfSKheann9emaoRK-ZZ-lAB} produce a contradiction with~\eqref{OHSn0ujfiHDiegb-0ujhrST}.

Thanks to this contradiction, we have proved that
\begin{equation}\label{800ntWAdw3-95ott}
\lambda_0=0\end{equation}
and therefore~$w_\lambda\le0$ in~$\Sigma_\lambda$ for all~$\lambda\in[0,1)$.
Using this with~$\lambda:=0$ we obtain that, for all~$x\in B_1\cap\{x_1>0\}$,
\begin{equation}\label{so03r8vnhuyeri8tuytyyhriueio}
0\ge w_0(x)= u(x_1,x_2,\dots,x_n)-u(-x_1,x_2,\dots,x_n).\end{equation}
On the other hand, setting~$v(x_1,x_2,\dots,x_n):=u(-x_1,x_2,\dots,x_n)$, we see that, for all~$x\in B_1$,
$$\Delta v(x)=\Delta u(-x_1,x_2,\dots,x_n)=f(u(-x_1,x_2,\dots,x_n))=f(v(x))$$
and also~$v>0$ in~$B_1$ and~$v=0$ on~$\partial B_1$. Accordingly,
from~\eqref{so03r8vnhuyeri8tuytyyhriueio} applied to~$v$ we deduce that, for all~$x\in B_1\cap\{x_1>0\}$,
$$0\ge v(x_1,x_2,\dots,x_n)-v(-x_1,x_2,\dots,x_n) =u(-x_1,x_2,\dots,x_n)-u(x_1,x_2,\dots,x_n) .$$
Putting together this and~\eqref{so03r8vnhuyeri8tuytyyhriueio} we finally obtain~\eqref{LAPIR-HNG6JEnIERNsTrnE4gt-10431}.

Additionally, for all~$x=(\lambda,y)\in\R\times\R^{n-1}$ with~$x\in B_1$ and~$\lambda>0$,
$$ \partial_1u(\lambda,y)=\lim_{h\searrow0}\frac{u(\lambda+h,y)-u(\lambda-h,y)}{2h}=
\lim_{h\searrow0}\frac{w_\lambda(\lambda+h,y)}{2h}\le0.
$$
This gives that~$\partial_1u\le0$ in~$B_1\cap\{x_1>0\}$.
Hence, since~$\Delta(\partial_1u)=f'(u)\partial_1u$ in~$B_1\cap\{x_1>0\}$,
we deduce from the Maximum Principle in Theorem~\ref{MAXPLECON45678CG-PM}
that either~$\partial_1u<0$ or~$\partial_1u$ vanish identically. But the latter cannot be
(since~$u(0)>0=u(e_1)$) and therefore~\eqref{LAPIR-HNG6JEnIERNsTrnE4gt-10432}
is proved.
\end{proof}

Now we utilize Lemmata~\ref{JNSMPERADVBDDOBYHGGJUINMRARYS3123-4} and~\ref{JNSMPERADVBDDOBYHGGJUINMRARYS3123}
to complete the proof of Theorem~\ref{GININITH}.

\begin{proof}[Proof of Theorem~\ref{GININITH}] Since the domain under consideration is a ball,
we can make use of Lemma~\ref{JNSMPERADVBDDOBYHGGJUINMRARYS3123} in every direction,
thus finding that for all~$\omega\in\partial B_1$ and~$x\in\R^n$ we have that~$u({\mathcal{R}}_\omega x)=u(x)$. Thus,
Lemma~\ref{JNSMPERADVBDDOBYHGGJUINMRARYS3123-4} gives that~$u(x)=u_0(|x|)$ for some~$u_0:\R\to\R$.
Hence, retaking the last claim in Lemma~\ref{JNSMPERADVBDDOBYHGGJUINMRARYS3123},
for all~$i\in\{1,\dots,n\}$ we have that~$\partial_iu(x)<0$ for all~$x\in B_1$ with~$x_i>0$.
Also, if~$x_i<0$ then
\begin{eqnarray*}&&\partial_iu(x)=\partial_i\big( u_0(|x|)\big)=\partial_i\Big( u_0\big(|(x_1,\dots,x_{i-1},-x_i,x_{i+1},\dots,x_n)|\big)\Big)\\&&\quad=\partial_i \Big(u
\big(x_1,\dots,x_{i-1},-x_i,x_{i+1},\dots,x_n\big)\Big)=-\partial_i u
(x_1,\dots,x_{i-1},-x_i,x_{i+1},\dots,x_n)>0.\end{eqnarray*}
As a result, we have that~$x_i\partial_iu(x)<0$ for all~$i\in\{1,\dots,n\}$ and~$x\in B_1\setminus\{0\}$, whence~$x\cdot\nabla u(x)<0$ for all~$x\in B_1\setminus\{0\}$, showing that~$u$ is radially strictly decreasing.
\end{proof}

The proof of Theorem~\ref{GININITH} presented here above follows the main ideas of the original article in~\cite{MR544879}. We showcase also a different proof, following~\cite{MR1159383, MR1662746}.

\begin{proof}[Another proof of Theorem~\ref{GININITH}] We provide a different proof of
Lemma~\ref{JNSMPERADVBDDOBYHGGJUINMRARYS3123} (from that, the original proof of Theorem~\ref{GININITH} remains unchanged). We utilize the moving plane method and the notation in~\eqref{GNNSHLIMSBXNSKt6-PRE-SEGBXDDUEv4JKS} and~\eqref{800ntWAdw3}, aiming at proving~\eqref{800ntWAdw3-95ott},
but with a different approach here.
The gist is to exploit the Maximum Principle for small volume domains presented in Theorem~\ref{M:PLE:SMALLV}
(in the first proof of Theorem~\ref{GININITH}, the Maximum Principle for small volume domains
is not essential, since its use can be replaced by the one in narrow domains given in
Theorem~\ref{M:PLE:NARROW}).

The details of this alternative proof go as follows. Suppose, by contradiction, that~\eqref{800ntWAdw3-95ott} does not hold, hence~$\lambda_0>0$.
Let~$\lambda_1\in(0,\lambda_0)$. Let also~$K\Subset\Sigma_{\lambda_0}$ be an open set with $C^\infty$ boundary
and such that
\begin{equation}\label{7iugP3RAFP3i0r3p4jr9POKSjdfh2765r83g5f4d}
\|f'\|_{L^{\infty}(\Omega)}\;|\Sigma_{\lambda_1}\setminus \overline{K}|^{\frac2n}
<c_0,\end{equation}
being~$c_0$ as in~\eqref{22qwrf2KMSDFGHJSIJD9iry93uwgfoiqwt7893uowig9fuwegfiGBSIMfjgmVkfRNKu4gl3a2teft6tejgs},
see Figure~\ref{TER24357iugP3RAFP3i0r3p4jr9POKSjdfh2765r83g5f4dRAeFI}.
In our framework, \eqref{equa0urjoewh3r73855-d-fingivis},
\eqref{GNNSHLIMSBXNSKt6-PRE2} and~\eqref{7iugP3RAFP3i0r3p4jr9POKSjdfh2765r83g5f4d}
allow us to use the Maximum Principle for small volume domains presented in Theorem~\ref{M:PLE:SMALLV},
applied here to the function~$w_{\lambda}$
and the domain~$\Sigma_{\lambda}\setminus \overline{K}$, with~$\lambda\in(\lambda_0,\lambda_1)$
conveniently close to~$\lambda_0$.
To this end, recalling~\eqref{OHSn0ujfiHDiegb-0ujhr-00CASILIM}, we have that~$\displaystyle\max_{\overline{K}}w_{\lambda_0}<0$,
that is~$\displaystyle\max_{\overline{K}}w_{\lambda_0}\le-\delta_0$, for some~$\delta_0>0$.
Hence, for all~$\lambda\in(\lambda_0,\lambda_1)$
sufficiently close to~$\lambda_0$, we have that
$$\max_{\overline{K}}w_{\lambda}\le0.$$
Since~$w_\lambda=0$ on~$T_\lambda\cap\overline{B_1}$
and~$w_\lambda<0$ on~$(\partial B_1)\cap\Sigma_\lambda$, we thus conclude that~$w_\lambda\le0$
on~$\partial(\Sigma_{\lambda}\setminus \overline{K})$. Notice also that~$\|f'\|_{L^{\infty}(\Omega)}\;|\Sigma_{\lambda}\setminus \overline{K}|^{\frac2n}
<c_0$, owing to~\eqref{7iugP3RAFP3i0r3p4jr9POKSjdfh2765r83g5f4d}, hence, by
the Maximum Principle for small volume domains in Theorem~\ref{M:PLE:SMALLV}, we infer that~$w_\lambda<0$
in~$\Sigma_\lambda\setminus \overline{K}$.

\begin{figure}
  \centering
  \includegraphics[width=.4\linewidth]{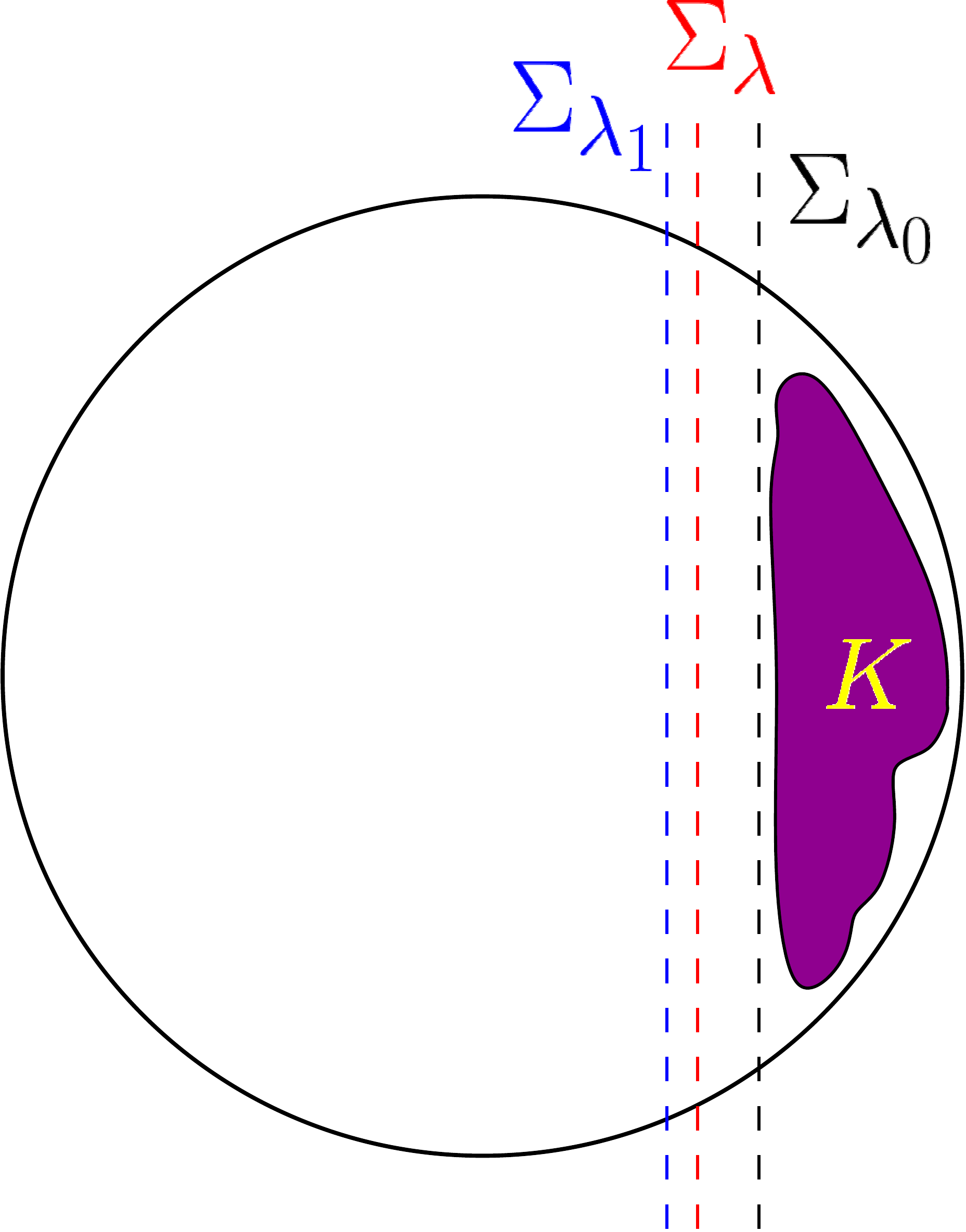}
 \caption{\sl The geometry involved in~\eqref{7iugP3RAFP3i0r3p4jr9POKSjdfh2765r83g5f4d}.}\label{TER24357iugP3RAFP3i0r3p4jr9POKSjdfh2765r83g5f4dRAeFI}
\end{figure}

Given the arbitrariness in the choice of~$K$ (we can freely move it around in~$\Sigma_{\lambda_0}$)
we thereby conclude that~$w_\lambda<0$
in~$\Sigma_\lambda$, which is in contradiction with the infimum property in~\eqref{800ntWAdw3}.\index{Gidas, Ni and Nirenberg's problem|)}\index{moving plane method|)}
\end{proof}

It is also interesting to recall that a different approach to the proof of
Theorem~\ref{GININITH} is possible under the additional assumption that
\begin{equation}\label{ADHYJGNN-A1}
{\mbox{$f$ is nonpositive and nondecreasing.}}
\end{equation}
For this, we follow~\cite{MR3003296} and\footnote{The technique in~\cite{MR3003296}
was also motivated by~\cite{MR653200, MR1382205}, see~\cite{DPV-ipkwd-0243}
for a general approach. Here we actually consider only
the simplest possible case, namely that of Theorem~\ref{GININITH} under the additional assumption~\eqref{ADHYJGNN-A1},
but the method introduced in~\cite{MR3003296} is much more general
and covers also the case of bounded but possibly discontinuous nonlinearities~$f$
and that of possibly singular or degenerate nonlinear operators. The monotonicity
assumption in~\eqref{ADHYJGNN-A1} can also be weakened (and dropped when~$n=2$,
since the parameter~$\beta$ in the forthcoming equation~\eqref{HSknJOPAMRAGN3MEt2r} vanishes
in this case, thus making the analysis of the function~$J$ in~\eqref{VSJIL-G222}
not necessary in dimension~$2$). For the sake of simplicity, here we only
confine ourselves to a model case and we refer to~\cite{MR3003296}
for the method and the results in their full generality.}
we argue as follows:

\begin{proof}[Another proof of Theorem~\ref{GININITH} under the additional assumption~\eqref{ADHYJGNN-A1}]
We let
$$M:=\max_{\overline{B_1}}u.$$
Given~$j\in\{1,\dots,n\}$, we recall that~$\nabla \partial_j u=0$ a.e. in~$\{\partial_ju=0\}$,
see e.g.~\cite[Theorem~6.19]{MR1817225}. In particular, we have that~$\partial_j^2 u=0$ a.e. in~$\{\partial_ju=0\}$
and accordingly 
\begin{equation}\label{VSJIL-G1}
{\mbox{for a.e.~$x\in\{\nabla u=0\}$ we have that~$f(u(x))=\Delta u(x)=0$.}}\end{equation}
Now we recall that~$f$ is assumed to be nonpositive and nondecreasing by~\eqref{ADHYJGNN-A1}.
Hence either~$f<0$ or there exists~$t_0\in\R$ such that
\begin{equation*}
{\mbox{$f<0$ in~$(-\infty,t_0)$ and~$f=0$ in~$[t_0,+\infty)$.}}
\end{equation*}
This observation and~\eqref{VSJIL-G1} give that either~$\{\nabla u=0\}=\varnothing$
or
\begin{equation}\label{VSJIL-G2}
\varnothing\ne\{\nabla u=0\}\subseteq\big\{x\in B_1 {\mbox{ s.t. }}f(u(x))=0\big\}\cup{\mathcal{Z}}\subseteq
\big\{x\in B_1 {\mbox{ s.t. }} u(x)\ge t_0\big\}\cup{\mathcal{Z}},\end{equation}
with~$|{\mathcal{Z}}|=0$.
As a consequence, when~$\big|\{\nabla u=0\}\big|>0$ we have that~$t_0\le M$.

But, more precisely, in this case it holds that~$t_0=M$, because otherwise~$\{ u>t_0\}\ne\varnothing$ and for all~$x\in\{ u>t_0\}$
one would have that~$\Delta u(x)=f(u(x))=0$ and thus the
Weak Maximum Principle in Corollary~\ref{WEAKMAXPLE} would entail that~$u$
is constant in any connected component of~$\{ u>t_0\}$, leading to a contradiction, since at least one of these connected
components must contain a point at which the value of~$u$ equals~$M>t_0$.

We can thus reconsider~\eqref{VSJIL-G2} and write that
\begin{equation*}
\{\nabla u=0\}\subseteq\big\{x\in B_1 {\mbox{ s.t. }} u(x)=M\big\}\cup{\mathcal{Z}}\end{equation*}
and consequently, since~$\nabla u(x)=0$ whenever~$u(x)=M$,
\begin{equation}\label{VSJIL-G22}
\big|\{\nabla u=0\}\big|=\big|\big\{x\in B_1 {\mbox{ s.t. }} u(x)=M\big\}\big|.\end{equation}

It is also convenient to use the Coarea Formula (see e.g.~\cite{MR3409135}) to note that
\begin{eqnarray*}&& 0=\int_{\R^n} |\nabla u(x)|\,\chi_{\{|\nabla u|=0\}}(x)\,dx=\int_\R
\left[\int_{x\in\{u=t\}} \chi_{\{|\nabla u|=0\}}(x)
\,d{\mathcal{H}}^{n-1}_x
\right]\,dt\\&&\qquad\qquad\qquad=\int_\R
\left[
{\mathcal{H}}^{n-1}\Big(
\big\{u=t\big\}\cap\big\{|\nabla u|=0\big\}
\Big)
\right]\,dt
\end{eqnarray*}
and therefore, for a.e.~$t\in\R$,
\begin{equation}\label{VSJIL-G33}
{\mathcal{H}}^{n-1}\Big(
\big\{u=t\big\}\cap\big\{|\nabla u|=0\big\}
\Big)=0.
\end{equation}
Now we set
\begin{equation}\label{VSJIL-G222} J(t):=\big| \{u>t\}\big| \qquad{\mbox{and}}\qquad
J_0(t):=\big| \{u>t\}\cap\{\nabla u\ne0\}\big|.\end{equation}
Note that~$J$ and~$J_0$ are monotone nonincreasing functions of~$t$, hence their derivative is well-defined for a.e.~$t\in(0,M)$.
We claim that
\begin{equation}\label{LjS0paA99i0jt4q38ygfc0hnf-20urjfvio3iy021-PRE}\begin{split}&
{\mbox{$J_0$ is an absolutely continuous function}}\\&{\mbox{and, for a.e.~$t\in(0,M)$, }}\quad
J'_0(t)= -\int_{\{ u=t\}\cap\{\nabla u\ne0\}} \frac{d{\mathcal{H}}^{n-1}_x}{|\nabla u(x)|}.\end{split}
\end{equation}
To check this, let~$\e>0$. We use the Coarea Formula (see e.g.~\cite{MR3409135}) to see that
\begin{equation}\label{LjS0paA99i0jt4q38ygfc0hnf-20urjfvio3iy021-PRE2}
\begin{split}&J_\e(t):=\big|\{u>t\}\cap\{|\nabla u|>\e\}\big|=
\int_{\R^n}\frac{|\nabla u(x)|
\,\chi_{\{u>t\}\cap\{|\nabla u|>\e\}}(x)
}{|\nabla u(x)|} \,dx\\&\qquad\qquad\qquad=
\int_\R\left[\int_{x\in\{u=\lambda\}}
\frac{\chi_{\{u>t\}\cap\{|\nabla u|>\e\}}(x)}{ |\nabla u(x)| }\,d{\mathcal{H}}^{n-1}_x
\right]\,d\lambda\\&\qquad\qquad\qquad=\int_t^M\left[\int_{x\in\{u=\lambda\}\cap\{|\nabla u|>\e\}}
\frac{d{\mathcal{H}}^{n-1}_x}{|\nabla u(x)|}
\right]\,d\lambda.\end{split}
\end{equation}
Therefore,
by the Monotone Convergence Theorem (see e.g.~\cite[Theorem 10.27]{MR3381284}) we can pass~$\e\searrow0$
in~\eqref{LjS0paA99i0jt4q38ygfc0hnf-20urjfvio3iy021-PRE2} and obtain that
\begin{equation*}
J_0(t)=\big|\{u>t\}\cap\{|\nabla u|>0\}\big|=\int_t^M\left[\int_{x\in\{u=\lambda\}\cap\{|\nabla u|>0\}}
\frac{d{\mathcal{H}}^{n-1}_x}{ |\nabla u(x)| }
\right]\,d\lambda
\end{equation*}
and this entails~\eqref{LjS0paA99i0jt4q38ygfc0hnf-20urjfvio3iy021-PRE}, as desired.

As a matter of fact, the result in~\eqref{LjS0paA99i0jt4q38ygfc0hnf-20urjfvio3iy021-PRE} can be sharpened in light of~\eqref{VSJIL-G33}, thus producing that, for a.e.~$t\in(0,M)$,
\begin{equation}\label{LjS0paA99i0jt4q38ygfc0hnf-20urjfvio3iy021-MEG}
J'_0(t)= -\int_{\{ u=t\}} \frac{d{\mathcal{H}}^{n-1}_x}{|\nabla u(x)|}.
\end{equation}

Also, owing to~\eqref{VSJIL-G22} and \eqref{VSJIL-G222}, for all~$t\in(0,M)$ we have that
$$ J(t)=J_0(t)+\big|\{ u=M\}\big|.$$
{F}rom this and~\eqref{LjS0paA99i0jt4q38ygfc0hnf-20urjfvio3iy021-MEG} one deduces that
\begin{equation}\label{LjS0paA99i0jt4q38ygfc0hnf-20urjfvio3iy021}
\begin{split}&
{\mbox{$J$ is an absolutely continuous function}}\\&{\mbox{and, for a.e.~$t\in(0,M)$, }}\quad
J'(t)= -\int_{\{ u=t\}} \frac{d{\mathcal{H}}^{n-1}_x}{|\nabla u(x)|}.\end{split}
\end{equation}

Now we define
$$ I(t):=-\int_{\{u>t\}} f(u(x))\,dx$$
and we observe that~$I$ is a monotone nonincreasing function of~$t$, thanks to~\eqref{ADHYJGNN-A1},
hence its derivative is well-defined for a.e.~$t\in(0,M)$.
Also, as~$\tau\nearrow0$,
\begin{equation}\label{LjS0paA99i0jt4q38ygfc0hnf-20urjfvio3iy022-ABSLCON2}\begin{split}& I(t+\tau)-I(t)=-
\int_{\{u>t+\tau\}} f(u(x))\,dx+\int_{\{u>t\}} f(u(x))\,dx
=\int_{\{t+\tau\ge u>t\}} f(u(x))\,dx\\&\qquad=-\int_{\{t+\tau\ge u>t\}} \big(f(t)+o(1)\big)\,dx=
\big(f(t)+o(1)\big)\big|\{t+\tau\ge u>t\}\big|\\
&\qquad= \big(f(t)+o(1)\big)\big(J(t)-J(t+\tau)\big).
\end{split}\end{equation}
This, the boundedness of~$f$ in the range of the solution~$u$ and~\eqref{LjS0paA99i0jt4q38ygfc0hnf-20urjfvio3iy021} give that
\begin{equation}\label{LjS0paA99i0jt4q38ygfc0hnf-20urjfvio3iy022-ABSLCON}
{\mbox{$I$ is an absolutely continuous function in~$(0,M)$}}.
\end{equation}
Furthermore, dividing by~$\tau$ in~\eqref{LjS0paA99i0jt4q38ygfc0hnf-20urjfvio3iy022-ABSLCON2}
and taking the limit as~$\tau\searrow0$, we thus conclude that, for a.e.~$t\in(0,M)$,
\begin{equation}\label{LjS0paA99i0jt4q38ygfc0hnf-20urjfvio3iy022} I'(t)=-f(t)\,J'(t).\end{equation}

Given~$\alpha$, $\beta\in\R$, to be specified later on, we now take into consideration the function
$$ K(t):= \big(I(t)\big)^\alpha\big(J(t)\big)^\beta.$$
In view of~\eqref{LjS0paA99i0jt4q38ygfc0hnf-20urjfvio3iy021} and~\eqref{LjS0paA99i0jt4q38ygfc0hnf-20urjfvio3iy022-ABSLCON} we know that
\begin{equation}\label{LjS0paA99i0jt4q38ygfc0hnf-20urjfvio3iy022-ABSLCON-LAK}
{\mbox{$K$ is an absolutely continuous function in~$(0,M)$}}.
\end{equation}
Also, we make use of~\eqref{LjS0paA99i0jt4q38ygfc0hnf-20urjfvio3iy022} to find that, a.e. in~$(0,M)$,
\begin{equation}\label{KINCRE2}
K'= \alpha I^{\alpha-1} J^\beta I'+\beta I^\alpha J^{\beta-1} J'=\Big(-
\alpha I^{\alpha-1} J^\beta f+\beta I^\alpha J^{\beta-1} \Big) J'.\end{equation}
Now we choose
\begin{equation}\label{HSknJOPAMRAGN3MEt2r} \alpha:=2\qquad{\mbox{and}}\qquad\beta:=-\frac{n-2}{n}\end{equation}
and we deduce from~\eqref{KINCRE2} that, a.e. in~$(0,M)$,
\begin{equation}\label{KINCRE3}
\begin{split}
K'\,&=\left(-
2I J^{-\frac{n-2}{n}} f-\frac{n-2}{n} I^2 J^{-\frac{n-2}{n}-1} \right) J'
\\&= \left(-
2J f-\frac{n-2}{n} I \right) I J^{-\frac{n-2}{n}-1} J'.\end{split}
\end{equation}
In addition, since~$f$ is supposed to be nondecreasing (recall~\eqref{ADHYJGNN-A1}),
$$ I(t)\le-\int_{\{u>t\}} f(t)\,dx =-f(t)\,\big|\{u>t\}\big|=-f(t)\,J(t).$$
Plugging this information into~\eqref{KINCRE3}, and recalling~\eqref{LjS0paA99i0jt4q38ygfc0hnf-20urjfvio3iy021}, we obtain that, a.e. in~$(0,M)$,
\begin{equation}\label{KINCRE3-bis}
\begin{split}
K'\le \left(-
2J f+\frac{n-2}{n} J
f \right) I J^{-\frac{n-2}{n}-1} J'
=-
\frac{n+2}{n} I J^{-\frac{n-2}{n}} f J'\le0.
\end{split}
\end{equation}
Thanks to~\eqref{LjS0paA99i0jt4q38ygfc0hnf-20urjfvio3iy022-ABSLCON-LAK},
we can thus reconstruct the oscillations of~$K$ from~$K'$ (see e.g.~\cite[Theorem~7.32(ii)]{MR3381284},
used here with~$f:=K$ and~$g:=1$) and we have that, for all~$a$, $b\in(0,M)$ with~$a<b$,
$$ \int_a^b K'(t)\,dt=K(b)-K(a).$$
In particular, by~\eqref{KINCRE3-bis}, we deduce that~$K$ is nonincreasing in~$(0,M)$. Also, recalling~\eqref{KINCRE3}
we see that
\begin{equation}\label{MSNSEMIMNED2454I5T744OMFD}\begin{split}& K(0^+)-K(M^-):=\lim_{{a\searrow0}\atop{b\nearrow M}}K(a)-K(b)=
-\int_0^M K'(t)\,dt\\ &\qquad= \int_0^M
\left(
2J(t) f(t)+\frac{n-2}{n} I(t) \right) I(t) \left(J(t)\right)^{-\frac{n-2}{n}-1} J'(t)
\,dt
.\end{split}\end{equation}

Now we utilize the Isoperimetric Inequality
(see~\cite[Theorem~1 in Section~5.5 and Theorem~2 in Section~5.6.2]{MR3409135}) to write that
\begin{equation}\label{MSNSEMIMNED2454I5T744OMFD7} \frac{ {\mathcal{H}}^{n-1}\big( \{ u=t\}\big)}{\big| \{u>t\}\big|^{\frac{n-1}n}}
\ge \frac{ {\mathcal{H}}^{n-1}(\partial B_1)}{| B_1|^{\frac{n-1}n}}=:c_\star.
\end{equation}
This, in combination with the Cauchy-Schwarz Inequality, gives that
\begin{equation}\label{KDIVBETHNC1O42M21E214VAJMDSONS}\begin{split}&
\left(J(t)\right)^{\frac{n-2}{n}+1}= \big| \{u>t\}\big|^{\frac{2(n-1)}{n}}\le
\left( \frac{ {\mathcal{H}}^{n-1}\big( \{ u=t\}\big)}{c_\star}\right)^2=
\frac1{c_\star^2}\left( \int_{\{u=t\}} \frac{\sqrt{|\nabla u(x)|}}{\sqrt{|\nabla u(x)|}}\,d{\mathcal{H}}^{n-1}_x\right)^2\\&\qquad\le
\frac1{c_\star^2} \int_{\{u=t\}} |\nabla u(x)|\,d{\mathcal{H}}^{n-1}_x\,
\int_{\{u=t\}} \frac{d{\mathcal{H}}^{n-1}_x}{|\nabla u(x)|}.
\end{split}\end{equation}
Furthermore, if~$\nu$ is the external normal to~$\{u>t\}$ we know that~$
\nu=-\frac{\nabla u}{|\nabla u|}$ outside a set of vanishing $(n-1)$-dimensional Hausdorff measure, thanks to~\eqref{VSJIL-G33}, whence by the Divergence Theorem (see e.g.~\cite[Theorem~1 in Section~5.8]{MR3409135}) we find that
\begin{equation}\label{MJM-1ierjfn2urt}
I(t)=-\int_{\{u>t\}} f(u(x))\,dx=-\int_{\{u>t\}} \div(\nabla u(x))\,dx
=\int_{\{u=t\}} |\nabla u(x)|\,d{\mathcal{H}}^{n-1}_x.
\end{equation}
Plugging this into~\eqref{KDIVBETHNC1O42M21E214VAJMDSONS} and recalling~\eqref{LjS0paA99i0jt4q38ygfc0hnf-20urjfvio3iy021} we conclude that
\begin{equation*}
\left(J(t)\right)^{\frac{n-2}{n}+1}\le
\frac{I(t)}{c_\star^2} \,
\int_{\{u=t\}} \frac{d{\mathcal{H}}^{n-1}_x}{|\nabla u(x)|}=
-\frac{I(t)\,J'(t)}{c_\star^2}.
\end{equation*}
We insert this information into~\eqref{MSNSEMIMNED2454I5T744OMFD} and we obtain
\begin{equation}\label{MSNSEMIMNED2454I5T744OMFD-2}
K(0^+)-K(M^-)\ge-
c_\star^2\int_0^M
\left(
2J(t) f(t)+\frac{n-2}{n} I(t) \right) \,dt.\end{equation}
In addition, using the notation
$$ F(r):=\int_0^r f(\tau)\,d\tau,$$
we have that
\begin{eqnarray*}&&
\int_0^M J(t) f(t)\,dt=\int_0^M \big| \{u>t\}\big| f(t)\,dt=
\int_0^M \left[ \int_{B_1}\chi_{\{u>t\}}(x)\,f(t)\,dx\right]\,dt\\&&\qquad=
\int_{B_1}\left[ \int_{0}^{u(x)} f(t)\,dt\right]\,dx=\int_{B_1}F(u(x))\,dx
\end{eqnarray*}
and
\begin{eqnarray*}&&-
\int_0^M I(t) \,dt=\int_0^M\left[
\int_{\{u>t\}} f(u(x))\,dx\right]\,dt=\int_{B_1}\left[
\int_0^{u(x)} f(u(x))\,dt\right]\,dx=\int_{B_1}u(x) f(u(x))\,dx.
\end{eqnarray*}
These observations and~\eqref{MSNSEMIMNED2454I5T744OMFD-2} lead to
\begin{equation*}
\frac{n\,\big(K(0^+)-K(M^-)\big)}{2c_\star^2}\ge- n\int_{B_1}F(u(x))\,dx+\frac{n-2}{2}\int_{B_1} u(x)f(u(x))
\,dx.
\end{equation*}
{F}rom this and the Pohozaev Identity in Theorem~\ref{Pohozaev Identity} we arrive at
\begin{eqnarray*} &&\frac{n\,\big(K(0^+)-K(M^-)\big)}{2c^2_\star}\ge
\frac12\,\int_{\partial B_1} (\partial_\nu u(x))^2\,(x\cdot\nu(x))\,d{\mathcal{H}}^{n-1}_x
=\frac12\,\int_{\partial B_1} (\partial_\nu u(x))^2\,d{\mathcal{H}}^{n-1}_x.
\end{eqnarray*}
This and the Cauchy-Schwarz Inequality lead to
\begin{eqnarray*} &&
\left(\int_{ B_1}f(u(x))\,dx\right)^2=
\left(\int_{ B_1}\div(\nabla u(x))\,dx\right)^2=
\left(\int_{\partial B_1}\nabla u(x)\cdot\nu(x)\,d{\mathcal{H}}^{n-1}_x\right)^2\\&&\qquad=
\left(\int_{\partial B_1} |\nabla u(x)|\,d{\mathcal{H}}^{n-1}_x\right)^2
=
\left(\int_{\partial B_1} |\partial_\nu u(x)|\,d{\mathcal{H}}^{n-1}_x\right)^2\\&&\qquad\le
{\mathcal{H}}^{n-1}(\partial B_1)\int_{\partial B_1}(\partial_\nu u(x))^2\,d{\mathcal{H}}^{n-1}_x
\le\frac{n {\mathcal{H}}^{n-1}(\partial B_1)\,\big(K(0^+)-K(M^-)\big)}{c^2_\star}.
\end{eqnarray*}
Thus, since~$K$ is nonnegative,
\begin{equation}\label{MSNSEMIMNED2454I5T744OMFD7C}
\begin{split}&
\left(\int_{ B_1}f(u(x))\,dx\right)^2
\le 
\frac{n {\mathcal{H}}^{n-1}(\partial B_1)\,K(0^+)}{c^2_\star}\\&\qquad=
\frac{n {\mathcal{H}}^{n-1}(\partial B_1)}{c^2_\star}
\lim_{t\searrow0}
\big(I(t)\big)^\alpha\big(J(t)\big)^\beta\\
&\qquad=\frac{n {\mathcal{H}}^{n-1}(\partial B_1)}{c^2_\star}\left(\int_{\{u>0\}} f(u(x))\,dx\right)^\alpha
\big| \{u>0\}\big|^\beta\\&\qquad=\frac{n {\mathcal{H}}^{n-1}(\partial B_1)\,|B_1|^\beta}{c^2_\star}\left(\int_{B_1} f(u(x))\,dx\right)^\alpha
\\&\qquad=
\frac{n {\mathcal{H}}^{n-1}(\partial B_1)\,|B_1|^{-\frac{n-2}{n}}}{c^2_\star}\left(\int_{B_1} f(u(x))\,dx\right)^2.\end{split}
\end{equation}
But, recalling~\eqref{B1} and~\eqref{MSNSEMIMNED2454I5T744OMFD7},
\begin{eqnarray*}&&
\frac{n {\mathcal{H}}^{n-1}(\partial B_1)\,|B_1|^{-\frac{n-2}{n}}}{c^2_\star}=
\frac{n {\mathcal{H}}^{n-1}(\partial B_1)\,|B_1|^{-\frac{n-2}{n}}\,
| B_1|^{\frac{2(n-1)}n} }{\big( {\mathcal{H}}^{n-1}(\partial B_1)\big)^2}=\frac{n\,
| B_1| }{{\mathcal{H}}^{n-1}(\partial B_1)}=1.
\end{eqnarray*}
With this, we have that all the inequalities used towards~\eqref{MSNSEMIMNED2454I5T744OMFD7C} must be equalities.
In particular, for a.e.~$t\in(0,M)$ the set~$\{u>t\}$ must
optimize the Isoperimetric Inequality in~\eqref{MSNSEMIMNED2454I5T744OMFD7}, hence
it must be a ball, say~$\{u>t\}=B_{r(t)}(x(t))$ for some~$r(t)>0$ and~$x(t)\in B_1$.

Also, the last inequality in~\eqref{KDIVBETHNC1O42M21E214VAJMDSONS} must be an equality and correspondingly~$|\nabla u(x)|$ must be proportional to~$\frac1{|\nabla u(x)|}$ on~$\{u=t\}$. This gives
that~$|\nabla u|$ is constant on~$\{u=t\}$, for a.e.~$t\in(0,M)$, and we thus write
in this setting that~$|\nabla u|=c(t)$ on~$\{u=t\}$.

Hence, for a.e.~$t\in(0,M)$, one can use~\eqref{MJM-1ierjfn2urt} and infer that
\begin{equation}\label{MSNSEMIMNED2454I5T744OMFD791}-
\int_{\{u>t\}} f(u(x))\,dx=\int_{\{u=t\}} |\nabla u(x)|\,d{\mathcal{H}}^{n-1}_x
=c(t)\,{\mathcal{H}}^{n-1}\big(\{u=t\}\big)=
|\nabla u(p)|\,{\mathcal{H}}^{n-1}\big(\{u=t\}\big),
\end{equation}
for all~$p\in\{ u=t\}$.

We now claim that if~$u(x)\in(0,M)$ then
\begin{equation}\label{MSNSEMIMNED2454I5T744OMFD790}
\nabla u(x)\ne0.
\end{equation}
To check this, suppose by contradiction that~$\nabla u(q)=0$ for some~$q\in B_1$ with~$u(q)\in(0,M)$ and let~$p_k$ be a sequence of points such that~\eqref{MSNSEMIMNED2454I5T744OMFD791} holds true and~$p_k\to q$ as~$k\to+\infty$. In this way, passing to the limit~\eqref{MSNSEMIMNED2454I5T744OMFD790} in this case we find that
\begin{equation*}-
\int_{\{u>u(q)\}} f(u(x))\,dx=
|\nabla u(q)|\,{\mathcal{H}}^{n-1}\big(\{u=t\}\big)=0.
\end{equation*}
Thus, since~$f$ is nonpositive (recall~\eqref{ADHYJGNN-A1}), necessarily~$f(u(x))=0$ for all~$x\in\{u>u(q)\}$
and accordingly~$\Delta u=f=0$ in~$\{u>u(q)\}$.
This and the Weak Maximum Principle in Corollary~\ref{WEAKMAXPLE} give that~$u$
is constant in any connected component of~$\{u>u(q)\}$, leading to a contradiction since
one of these connected components must contain a point at which the value of~$u$ is~$M>u(q)$.
The proof of~\eqref{MSNSEMIMNED2454I5T744OMFD790} is thereby complete.

{F}rom~\eqref{MSNSEMIMNED2454I5T744OMFD790} it follows that for all~$t\in(0,M)$
\begin{equation}\label{ODEMD0-pqwjrf-5}
\min_{\{u=t\}} |\nabla u|>0\end{equation}
and accordingly, by~\eqref{LjS0paA99i0jt4q38ygfc0hnf-20urjfvio3iy021}, $J$ is locally Lipschitz
in~$(0,M)$. Hence, since
$$ J(t)=\big|\{u>t\}\big|=\big|B_{r(t)}\big|=|B_1|\,(r(t))^n,$$
we also have that
\begin{equation}\label{MSNSEMIMNED2454I5T744OMFD7908}
{\mbox{$(0,M)\ni t\mapsto r(t)$ is locally Lipschitz.}}
\end{equation}
Now, if~$T$, $t\in(0,M)$ with~$T>t$, we claim that
\begin{equation}\label{MSNSEMIMNED2454I5T744OMFD790889}
|x(T)-x(t)|\le r(t)-r(T).
\end{equation}
For this, we can suppose that~$x(T)\ne x(t)$, otherwise we are done.
Also, we have that~$B_{r(t)}(x(t))=
\{u>t\}\supseteq\{u>T\}=B_{r(T)}(x(T))$
and therefore, for all~$\rho\in(0,r(T))$ and all~$\varpi\in\partial B_1$,
we have that~$x(T)+\rho\varpi\in B_{r(T)}(x(T))\subseteq B_{r(t)}(x(t))$, giving that
$$ |x(T)+\rho\varpi-x(t)|<r(t).
$$
Picking~$\varpi:=\frac{x(T)-x(t)}{|x(T)-x(t)|}$ we get that
\begin{eqnarray*}&& |x(T)-x(t)|+\rho=\left( 1+\frac{\rho}{|x(T)-x(t)|}\right)\,|x(T)-x(t)|\le
\left| 1+\frac{\rho}{|x(T)-x(t)|}\right|\,|x(T)-x(t)|\\&&\qquad\qquad=
|x(T)-x(t)+\rho\varpi|
<r(t).
\end{eqnarray*}
Thus, by taking the limit~$\rho\nearrow r(T)$ we prove~\eqref{MSNSEMIMNED2454I5T744OMFD790889},
as desired.

Then, by~\eqref{MSNSEMIMNED2454I5T744OMFD7908} and~\eqref{MSNSEMIMNED2454I5T744OMFD790889} we obtain that
\begin{equation}\label{MSNSEMIMNED2454I5T744OMFD7908sm4irkf}
{\mbox{$(0,M)\ni t\mapsto x(t)$ is locally Lipschitz.}}\end{equation}

Now we show that
\begin{equation}\label{MSNSEMIMNED2454I5T744OMFD7908sm4irkfX}
{\mbox{$(0,M)\ni t\mapsto x(t)$ is constant.}}\end{equation}
Indeed, suppose not: then there exist~$t_1$ and~$t_2$ in~$(0,M)$ with~$t_1<t_2$
such that~$x(t_1)\ne x(t_2)$. Hence, by~\eqref{MSNSEMIMNED2454I5T744OMFD7908sm4irkf},
$$ 0\ne x(t_2)-x(t_1)=\int_{t_1}^{t_2} \frac{dx}{dt}(t)\,dt$$
and consequently there exists~$t_\star\in(t_1,t_2)$ such that~$\frac{dx}{dt}(t_\star)$ exists and is different from zero.

\begin{figure}
                \centering
                \includegraphics[width=.35\linewidth]{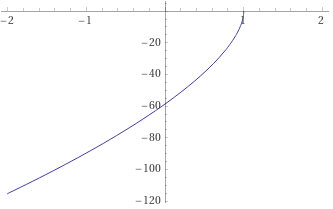}
        \caption{\sl Sketch of the function~$f$ discussed in footnote~\ref{234yu346HAFGo0k-OUMSldRGIRA4ihwjonmfegfAXELHARLRPASDfeDNOJHNFOJEDFO}.}\label{234yu346HAFGo0k-OUMSldRGIRA4ihwjonmfegfAXELHARLRPASDfeDNOJHNFOJED}
\end{figure}

We can thereby define~$y:=\frac{dx}{dt}(t_\star)\ne0$ and~$z:=\frac{y}{|y|}\in\partial B_1$.
Let also~$P_\pm (t):=x(t)\pm r(t)z$. We stress that~$P_\pm(t)\in \partial B_{r(t)}(x(t))=\partial\{u>t\}$.
Therefore, we have that~$u(P_\pm(t))=t$ and~$\big|\nabla u(P_\pm (t))\big|=c(t)$.
Hence, since the exterior normal of~$B_{r(t)}(x(t))$ at~$P_\pm(t)$ is~$\pm z$, we find that
$$ \pm\nabla u(P_\pm (t))\cdot z=-\big|\nabla u(P_\pm(t)\big|=-c(t).
$$
As a result, a.e.~$t\in(0,M)$ we have that
\begin{eqnarray*}&& 1=\frac{d}{dt} t=\frac{d}{dt}
u(P_\pm(t))=\frac{d}{dt}u\big(x(t)\pm r(t)z\big)= \nabla u\big(x(t)\pm r(t)\,z\big)\cdot
\left( 
\frac{dx}{dt}(t)\pm \frac{dr}{dt}(t)z\right)
\\&&\qquad=\nabla u(P_\pm(t))\cdot
\left( 
\frac{dx}{dt}(t)\pm \frac{dr}{dt}(t)z\right)=
\nabla u(P_\pm(t))\cdot\frac{dx}{dt}(t)-c(t)\frac{dr}{dt}(t)
\end{eqnarray*}
and therefore
\begin{eqnarray*}&& 0=1-1
=\left[\nabla u(P_+(t))\cdot\frac{dx}{dt}(t)- c(t)\frac{dr}{dt}(t)\right]
-\left[\nabla u(P_-(t))\cdot\frac{dx}{dt}(t)- c(t)\frac{dr}{dt}(t)\right]\\&&\qquad\qquad\qquad=
\Big(\nabla u(P_+(t))-\nabla u(P_+(t))\Big)\cdot\frac{dx}{dt}(t).
\end{eqnarray*}
Consequently, choosing~$t:=t_\star$,
$$ 0=\Big(\nabla u(P_+(t_\star))-\nabla u(P_+(t_\star))\Big)\cdot y=
|y|\Big(\nabla u(P_+(t_\star))-\nabla u(P_+(t_\star))\Big)\cdot z=-2|y|c(t_\star)
$$
and accordingly~$c(t_\star)=0$. This is in contradiction with~\eqref{MSNSEMIMNED2454I5T744OMFD790}
and thus we have completed the proof of~\eqref{MSNSEMIMNED2454I5T744OMFD7908sm4irkfX}.

Then, by~\eqref{MSNSEMIMNED2454I5T744OMFD7908sm4irkfX}, we deduce that~$\{u>t\}=B_{r(t)}(x_\star)$
for some~$x_\star\in\Omega$ and therefore~$u$ is radially symmetric with respect to the point~$x_\star$
(which, up to a translation, we can suppose to be the origin).

It remains to show that~$u$ is radially strictly decreasing. To this end, it suffices to show that
\begin{equation}\label{ODEMD0-pqwjrf-3}
\{u=M\}=\{0\}.
\end{equation}
Indeed, if this is true, suppose by contradiction that there exist~$a$, $b\in[0,1)$ such that~$a<b$ and~$u(a e_1)\le u(be_1)$.
We have that~$u(be_1)\ne M$ (otherwise, by~\eqref{ODEMD0-pqwjrf-3}, $a<b=0$, which is a contradiction)
and also~$u(ae_1)\ne M$ (otherwise~$u(be_1)=M$, which we have just excluded). Let~$m\in[a,1]$ be such that
$$ u(me_1)=\max_{s\in[a,1]} u(se_1)$$
and note that~$m\ne1$ since, by assumption, $0=u(e_1)<u(be_1)$.
Additionally, we can take~$m>a$, since~$u(a e_1)\le u(be_1)$. As a result, since~$u$ is radial,
\begin{equation}\label{ODEMD0-pqwjrf-4} 0=\left|
\frac{d}{ds} u(se_1)\Big|_{s=m}\right|=\left|\nabla u(se_1)\cdot e_1\right|=|\nabla u(se_1)|.\end{equation}
Moreover, we have that~$u(me_1)\in(0,M)$, due to~\eqref{ODEMD0-pqwjrf-3}. This and~\eqref{ODEMD0-pqwjrf-4}
are in contradiction with~\eqref{ODEMD0-pqwjrf-5}. 

\begin{figure}
                \centering
                \includegraphics[width=.35\linewidth]{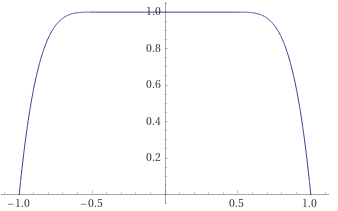}
        \caption{\sl Sketch of the radial profile of the function~$u$ discussed in footnote~\ref{234yu346HAFGo0k-OUMSldRGIRA4ihwjonmfegfAXELHARLRPASDfeDNOJHNFOJEDFO}.}\label{234yu346HAFGo0k-OUMSldRGIRA4ihwjonmfegfAXELHARLRPASDfeDNOJHNFOJEE}
\end{figure}

This argument shows that~$u$ has to be radially strictly decreasing once~\eqref{ODEMD0-pqwjrf-3} is proved.
Hence, we now establish~\eqref{ODEMD0-pqwjrf-3}. To this end, we observe that it is sufficient to
check that
\begin{equation}\label{ODEMD0-pqwjrf-3b}
\{u=M\}\subseteq\{0\},
\end{equation}
since this, together with the fact that~$\{u=M\}\ne\varnothing$, would entail~\eqref{ODEMD0-pqwjrf-3}.
Thus, to prove~\eqref{ODEMD0-pqwjrf-3b} we argue by contradiction, supposing that there exist~$\omega\in\partial B_1$
and~$\ell\in(0,1)$ such that~$u(\ell\omega)=M$. Note that necessarily~$\nabla u(\ell\omega)=0$.
Moreover, we have that~$u(t\omega)=M$ for all~$t\in[0,\ell]$, thanks to~\eqref{MSNSEMIMNED2454I5T744OMFD790}.
As a consequence,
\begin{equation}\label{aiw9c8h374358v7984ifjhkgheria}
0=\Delta u(0)=f(u(0))=f(M).\end{equation}
Now, writing~$u(x)=u_0(|x|)$ for some~$u_0:[0,1]\to\R$,
we use the formula for the Laplace operator in spherical coordinates in
Theorem~\ref{SPH} and we find that~$u_0$ is a solution of the following Cauchy problem:
\begin{equation}\label{vgy0-poj3r1tf90nhoitellomega}
\begin{dcases}
\displaystyle u_0''(r)+\frac{n-1}r\,u'_0(r)=f(u(r)),& {\mbox{ for all }}r\in(0,1),\\
u_0(\ell)=M,\\
u_0'(\ell)=0.
\end{dcases}\end{equation}
We observe that the function constantly equal to~$M$ is also a solution of~\eqref{vgy0-poj3r1tf90nhoitellomega},
in light of~\eqref{aiw9c8h374358v7984ifjhkgheria}.
Therefore, by the uniqueness\footnote{We stress that we are exploiting here the fact
that~$f$ was taken to be of class~$C^1(\R)$ in the statement of Theorem~\ref{GININITH},
in order to use the classical uniqueness result for the Cauchy problem (see e.g.~\cite[page~74]{MR1801796}).
This condition on~$f$ can be relaxed to a Lipschitz assumption, but cannot be completely dropped without losing
the property that~$u$ is radially strictly decreasing. As an example, we observe that, given~$\varrho\in(0,1)$ and~$\alpha>2$, the function
$$ u(x):=\begin{dcases}
1-\displaystyle \left(\frac{|x|^2-\varrho^2}{1-\varrho^2}\right)^\alpha,& {\mbox{ for all }}x\in B_1\setminus B_\varrho,\\
1 & {\mbox{ in }} B_\varrho
\end{dcases}$$
satisfies~\eqref{32254236equation-GNN-onwf} with
$$ f(t):=\begin{dcases}\displaystyle
-\frac{2 \alpha (1-t)^{\frac{\alpha - 2}\alpha}}{(1 - \varrho^2)^{2}}
\Big((n+2\alpha - 2)(1-\varrho^2)(1-t)^{\frac1\alpha}+2 (\alpha - 1)\varrho^2
\Big) & {\mbox{ if }} t<1\\0& {\mbox{ if }}t\ge1
.\end{dcases}$$
Notice that this~$u$ is radial but not radially strictly decreasing, since~$\{u=1\}=\overline{B_\varrho}$.
Note also that in this case~$f\in C^{\frac{\alpha - 2}\alpha}(\R)$ but~$
f\not\in C^1(\R)$. See Figure~\ref{234yu346HAFGo0k-OUMSldRGIRA4ihwjonmfegfAXELHARLRPASDfeDNOJHNFOJED} \label{234yu346HAFGo0k-OUMSldRGIRA4ihwjonmfegfAXELHARLRPASDfeDNOJHNFOJEDFO}
for a sketch of such a function~$f$ when~$\varrho:=\frac12$, $\alpha:=3$ and~$n:=2$. See also Figure~\ref{234yu346HAFGo0k-OUMSldRGIRA4ihwjonmfegfAXELHARLRPASDfeDNOJHNFOJEE} for a sketch of the radial profile of the corresponding solution~$u$.} of solutions of initial data problem for ordinary differential equations, we conclude that
necessarily~$u_0$ is constant, hence so is~$u$, in contradiction with the setting of Theorem~\ref{GININITH}.
The proof of~\eqref{ODEMD0-pqwjrf-3b} is thereby complete.
\end{proof}

\chapter{Local existence theory in the real analytic setting}

While most of these pages are devoted to the study of partial differential equations,
it is interesting now to revisit the classical theory of ordinary differential equations to
recognize its influence on a rather delicate aspect of the theory that we presented.
Namely, we have developed so far a series of rather specific existence and uniqueness
results for some class of elliptic equations (see e.g.
Corollary~\ref{S-coroEXIS-M023} and Theorems~\ref{Theorem6.14GT},
\ref{oind-94984u5f-14UNSpas88saleohi4a2T}, \ref{21980987y4h2nf9325teuationdetildealphauation-2TH}
and~\ref{NOWI}) and it is therefore natural to wonder about a ``general'' approach
guaranteeing existence results for partial differential equations in rather general scenarios.

One can be hopeful about it, since the case of ordinary differential equations
is settled in a general form, ensuring the existence and uniqueness of solutions
via the Picard-Lindel\"of-Cauchy-Lipschitz Theorem (see e.g.~\cite{MR0064934, MR1336820}).
For example, in the case of partial differential equations, one may want to consider an equation,
or a system of equations, of a given order~$m$ (i.e., in which the highest derivative involved is of order~$m$),
assign on a hypersurface the (say, normal) derivatives up to the order~$m-1$, and investigate the existence
and uniqueness of solutions.

\section{The Cauchy-Kowalevsky Theorem}

To formalize the above mentioned setting, we can consider a smooth hypersurface~$\Gamma$ described as a level set
of a smooth function~$\Phi$. That is, we suppose that~$\Phi\in C^\infty(\R^n)$ with~$\nabla\Phi(x)\ne0$
for all~$x\in\R^n$ and set
\begin{equation}\label{GABsaeM-034-001}
\Gamma_0:=\{x\in\R^n{\mbox{ s.t. }}\Phi(x)=0\}.\end{equation}

\begin{figure}
  \centering
  \includegraphics[width=.75\linewidth]{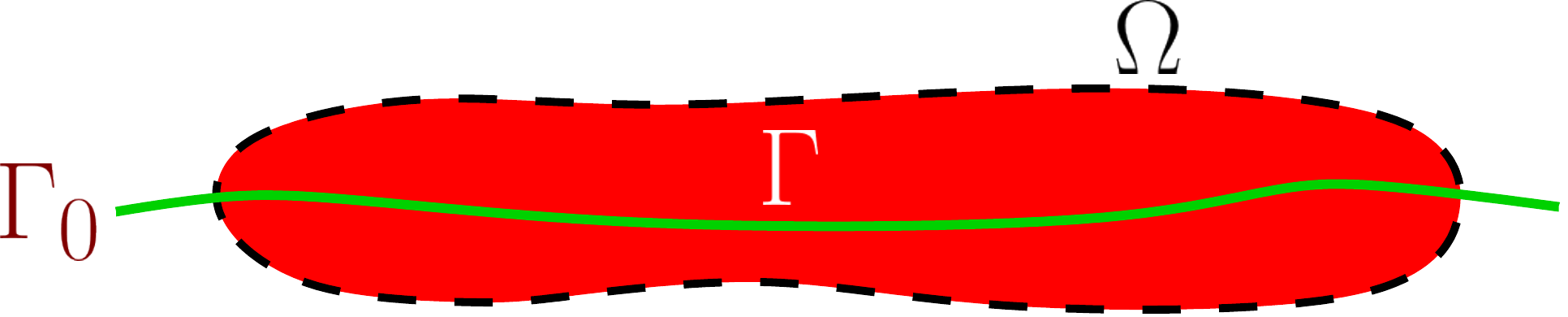}
 \caption{\sl The geometry in~\eqref{GABsaeM-034-001} and~\eqref{GABsaeM-034-001CAK}.}\label{SHOCACURAeFICAK}
\end{figure}

Let~$\Omega$ be an open subset of~$\R^n$, with
\begin{equation}\label{GABsaeM-034-001CAK}
\Gamma:=\Omega\cap\Gamma_0\ne\varnothing,\end{equation} see Figure~\ref{SHOCACURAeFICAK}.

Let also~$N\in\N$ with~$N\ge1$
and~$u:\Omega\to\R^N$. We use the notation~$u=(u_1,\dots,u_N)$, with~$u_\ell:\Omega\to\R$
for each~$\ell\in\{1,\dots,N\}$.

Given~$\mu\in\N^N$, we denote by~$D^\mu_\star u$ the collection of all the partial
derivatives of~$u_\ell$ up to order~$\mu_\ell$, for all~$\ell\in\{1,\dots,N\}$ that is
\begin{equation*}
D^\mu_\star u:=\Big\{ \partial^\alpha u_\ell \Big\}_{{{\alpha\in\N^n}\atop{|\alpha|\le\mu_\ell}}\atop{\ell\in\{1,\dots,N\}}}.\end{equation*}
Thus, a system of partial differential equations of order~$\mu\in\N^N$
is described by a set of equations of the form
\begin{equation}\label{GABsaeM-034-002}
F_\ell(x,D^\mu_\star u(x))=0\end{equation} for~$x\in\Omega$ and~$\ell\in\{1,\dots,N\}$, for some smooth functions~$F_1,\dots,F_N$. We will use the short notation~$F:=(F_1,\dots,F_N)$.

The prescription along~$\Gamma$ reads~$\partial^j_\nu u_\ell=\phi_{j,\ell}$ for all~$j\in\N$ with~$j\le \mu_\ell-1$
and~$\ell\in\{1,\dots,N\}$, for suitable assigned functions~$\phi_{j,\ell}$.

Here, we are denoting by~$\nu$ a smooth unit normal vector field along~$\Gamma$ and by~$
\partial^j_\nu f$ the $j$th derivative of a function~$f$ in the direction of~$\nu$, i.e., for each~$x\in\Gamma$,
\begin{equation}\label{INFEEFEBARRA} \begin{split} \partial^j_\nu f(x)\,&:=\left.
\frac{d^j}{d\tau^j} f(x+\tau\nu(x))\right|_{\tau=0}\\&=
\frac{d^j}{d\tau^j}\left. \left(\sum_{{0\le m\le j}\atop{1\le k_1,\dots,k_m\le n}} \frac{\tau^m}{m!}\,
\frac{\partial^m f}{\partial x_{k_1}\dots\partial x_{k_m}}(x)\,\nu_{k_1}(x)\dots\nu_{k_m}(x)
+o(\tau^j)\right)\right|_{\tau=0}\\&=
\sum_{1\le k_1,\dots,k_j\le n}
\frac{\partial^j f}{\partial x_{k_1}\dots\partial x_{k_j}}(x)\,\nu_{k_1}(x)\dots\nu_{k_j}(x).\end{split}\end{equation}
With this notation, and inspired by the case of ordinary differential equations, one can write \index{Cauchy problem}
a Cauchy problem for partial differential equations in the form
\begin{equation}\label{GEN:CAUC}
\begin{dcases}
F(x,D^\mu_\star u(x))=0 & {\mbox{ for all }}x\in\Omega,\\
\partial^j_\nu u_\ell(x)=
\phi_{j,\ell}(x)& {\mbox{ for all $x\in\Gamma$, $\ell\in\{1,\dots,N\}$
and $j\in\{0,\dots,\mu_\ell-1\}$.}}
\end{dcases}
\end{equation}
The investigation of this problem was pioneered by
Cauchy in~\cite{C:CAUCHY1, C:CAUCHY2, C:CAUCHY3, C:CAUCHY4, C:CAUCHY5, C:CAUCHY6}
for nonlinear differential equations of second order
and then generalized by Kowalevsky in~\cite{MR1579652}
for general analytic nonlinear systems of differential
equations (independently, a rather general case was also treated by Darboux~\cite{GASTO}
and a simplified strategy was put forth by Goursat in~\cite{MR1504315}).
Nowadays, the problems as the one in~\eqref{GEN:CAUC}, in which one seeks
a solution of a partial differential equation satisfying certain conditions that are prescribed on a hypersurface in the domain,
are called ``Cauchy problems''. \index{Cauchy problem}

\begin{figure}
                \centering
                \includegraphics[height=.36\textwidth]{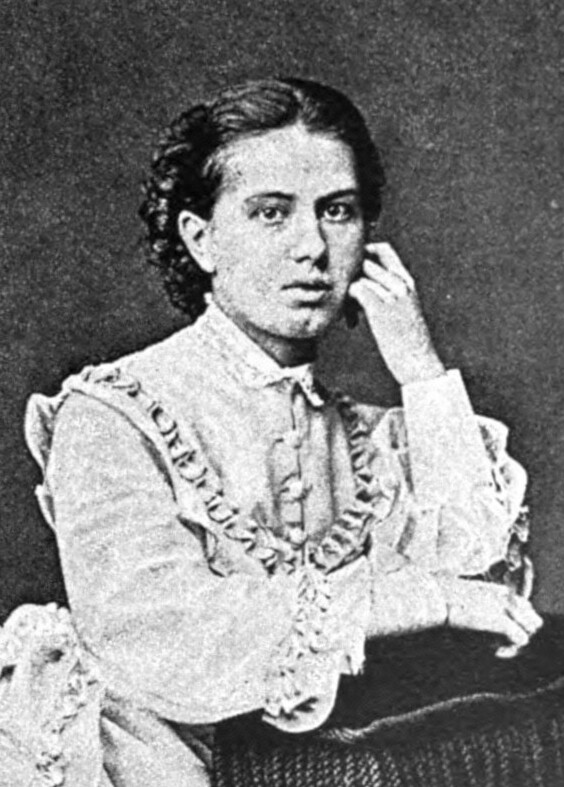} $\quad$
                \includegraphics[height=.36\textwidth]{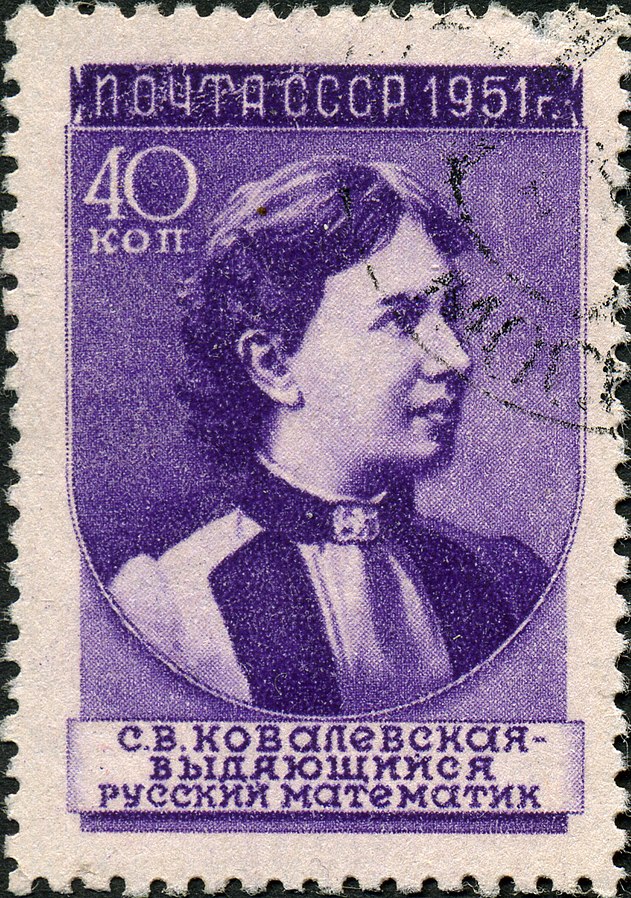} $\quad$
                \includegraphics[height=.36\textwidth]{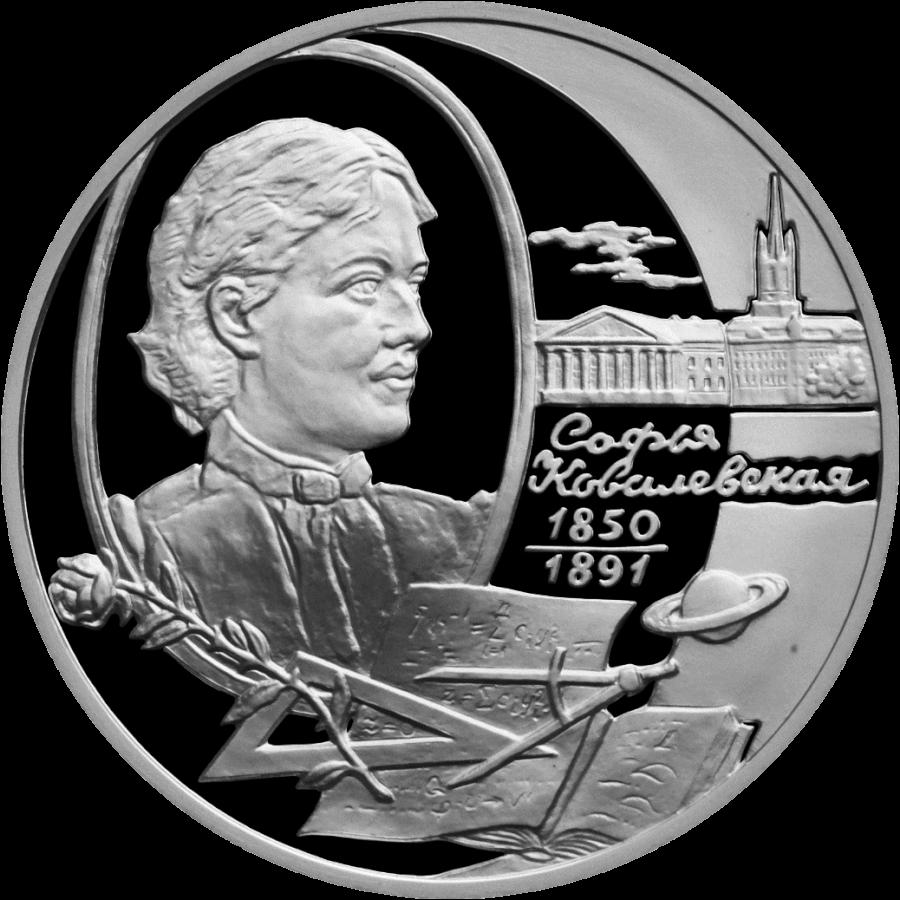}
        \caption{\sl Sophie Kowalevski (Public Domain images from
        Wikipedia).}\label{HAFNEWYDH-3-D-SOFKO}
\end{figure}

The main result related to~\eqref{GEN:CAUC}
goes under the name\footnote{Sof'ja Vasil'evna Kovalevskaja, born Korvin-Krukovskaja,
often used the name of Sophie Kowalevski (or Kowalevsky, or von Kowalevsky) 
in her academic publications
and her first name is ofttimes replaced by the nickname Sonya, hence different spellings and
transliterations are overabundant.
At that time, Russia happened to be a very conservative country and
women were not allowed to attend universities.
To study abroad, Kowalevsky needed written permission from her husband, so, in 1868 she contracted a ``fictitious marriage'' with a paleontology student
(and translator of the works of Charles Darwin).
She obtained her doctorate (first woman ever) in Germany under the supervision of Weierstra{\ss}. The
Cauchy-Kowalevsky Theorem is one of the results of her PhD thesis.

The name of Kowalevsky is also linked to the study of the precession of a spinning top with special symmetry regarding moments of inertia and center of gravity which make this mechanical system integrable (for this study, Kowalevsky was awarded the Prix Bordin from the French Academy of Sciences in 1888).

Besides mathematics, Kowalevsky was very active in the radical movement of Russian nihilism
(defined by Russian anarchist Peter Kropotkin
as ``the symbol of struggle against all forms of tyranny, hypocrisy, and artificiality''),
she wrote plays and autobiographical novels,
and she served as a nurse in the 1871 Paris Commune.

The Motherland, Russia, never recognized Kowalevsky's academic degrees obtained abroad, and
never offered her an academic position (but she became full professor at  Stockholm University,
as well as an editor of the very prestigious Swedish journal Acta Mathematica). The
Russian Academy of Sciences changed its rules to appoint her Corresponding Member.
See Figure~\ref{HAFNEWYDH-3-D-SOFKO} for a portrait of Kowalevsky at 18 years and for a commemorative coin and stamp.}
of Cauchy-Kowalevsky Theorem, which, in a suitable form, will be presented below\footnote{We stress that Theorem~\ref{CKTEO} is certainly not the ultimate or the most general form
of the Cauchy-Kowalevsky Theorem (yet, it is general enough to consider nonlinear systems
of various order).} in Theorem~\ref{CKTEO}.
\medskip

As a special feature, the Cauchy-Kowalevsky Theorem deals with a real analytic setting (the complex analytic case can be treated as well, just with minor modifications in the notation).
Before going into the details of the Cauchy-Kowalevsky Theorem, however, it is instructive to look at some negative examples,
to discover how delicate this matter is. We start with a one-dimensional\footnote{We observe that
the equation in~\eqref{GEN:CAUC:HEAT} is precisely the heat equation in~\eqref{DAGB-ADkrVoiweLL4re2346ytmngrrUj7}
in absence of sources and in one spatial dimension,
with the notation~$x_2:=t$.

It is worth remarking that a solution of the heat equation~\eqref{GEN:CAUC:HEAT} for small~$x_2>0$ and small~$|x_1|$ does exist, just it is not real analytic (thus confirming that analyticity is not always the most convenient setting for a partial differential equation). Indeed, the function
$$ u(x_1,x_2):=\int_{-1/2}^{1/2}\frac{e^{-\frac{(x_1-y)^2}{4x_2}}}{\sqrt{4\pi x_2}\,(1-y)}\,dy$$is a solution as desired, but it is not real analytic in~$x_2$ at~$x_2=0$.} heat equation:

\begin{theorem}
The Cauchy problem\footnote{Notice that any prescription of~$\partial_{x_2}u$ along~$\left(-\frac12,\frac12\right)\times\{0\}$
would cast~\eqref{GEN:CAUC:HEAT} into the general setting of~\eqref{GEN:CAUC}, here with~$n:=2$, $N:=1$,
$\Omega:=(-\epsilon,\epsilon)\times(-\epsilon, \epsilon)$, $\Gamma_0:=(-1,1)\times\{0\}$
and~$\Gamma:=\left(-\epsilon,\epsilon\right)\times\{0\}$.}
\begin{equation}\label{GEN:CAUC:HEAT}
\begin{dcases}
\partial_{x_1x_1} u(x)-\partial_{x_2} u(x)=0 & {\mbox{ for all $x=(x_1,x_2)\in(-\epsilon,\epsilon)\times
(-\epsilon,\epsilon)$,}}\\
u(x_1,0)=\frac{1}{1-x_1}& {\mbox{ for all $x_1\in \left(-\epsilon,\epsilon\right)$.}}
\end{dcases}
\end{equation}
does not
admit any real analytic solution. 

More precisely,
\begin{equation}\label{GEN:CAUC:HEAT:NO}
{\mbox{there is no solution of~\eqref{GEN:CAUC:HEAT} which is real analytic in a neighborhood of the origin.}}
\end{equation}\end{theorem}

\begin{proof} Suppose by contradiction that such a solution exists and let us write it, in a neighborhood of the origin,
as
$$ u(x_1,x_2)=\sum_{i,j=0}^{+\infty} c_{i,j} \,x_1^i \,x_2^j.$$
We claim that
\begin{equation}\label{GEN:CAUC:HEAT:NO2}
c_{i,j}=\frac{\displaystyle\prod_{k=1}^{2j} (i+k)}{j!}.
\end{equation}
The proof of~\eqref{GEN:CAUC:HEAT:NO2} is by induction over~$j$.
Indeed, for small values of~$x_1$ we know that
$$ \sum_{i=0}^{+\infty} c_{i,0} \,x_1^i =u(x_1,0)=\frac{1}{1-x_1}=\sum_{i=0}^{+\infty}x_1^i,$$
showing that~$c_{i,0}=1$, thus confirming~\eqref{GEN:CAUC:HEAT:NO2} when~$j=0$.

Suppose now that~\eqref{GEN:CAUC:HEAT:NO2} is true for the index~$j$ and let us prove it for the index~$j+1$.
To this end, we observe that
\begin{eqnarray*}&&
\sum_{i,j=0}^{+\infty} (i+2)(i+1) c_{i+2,j} \,x_1^i \,x_2^j
=\sum_{i,j=0}^{+\infty} i(i-1) c_{i,j} \,x_1^{i-2} \,x_2^j
=\partial_{x_1x_1} u(x)\\&&\qquad=\partial_{x_2} u(x)
=\sum_{i,j=0}^{+\infty} j c_{i,j} \,x_1^i \,x_2^{j-1}=\sum_{i,j=0}^{+\infty} (j+1) c_{i,j+1} \,x_1^i \,x_2^{j},
\end{eqnarray*}
showing that~$(i+2)(i+1) c_{i+2,j}=(j+1) c_{i,j+1}$.

The inductive hypothesis thereby entails that
$$ c_{i,j+1}=\frac{(i+2)(i+1) c_{i+2,j}}{(j+1) }=
\frac{(i+2)(i+1) \displaystyle\prod_{k=1}^{2j} (i+2+k)}{(j+1) \,j!}
=\frac{ \displaystyle\prod_{k=1}^{2(j+1)} (i+k)}{(j+1)!}
,$$
which ends the proof of~\eqref{GEN:CAUC:HEAT:NO2}.

{F}rom this, we infer that
$$ c_{0,j}=\frac{\displaystyle\prod_{k=1}^{2j} k}{j!}=
\frac{(2j)!}{j!}$$
and accordingly, for small~$x_2$,
$$ u(0,x_2)=\sum_{j=0}^{+\infty} c_{0,j}\,x_2^j=
\sum_{j=0}^{+\infty} \frac{(2j)!}{j!}\,x_2^j,
$$
which is a divergent series and the proof of~\eqref{GEN:CAUC:HEAT:NO} is complete.\end{proof}

\begin{figure}
                \centering
                \includegraphics[width=.43\linewidth]{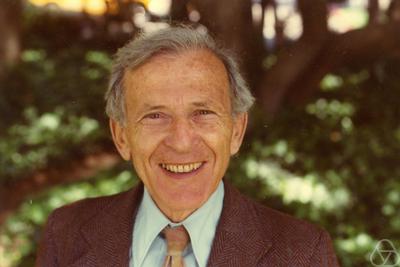}
        \caption{\sl Hans Lewy (photo by  George M. Bergman, image from
        Wikipedia,
        licensed under the GNU Free Documentation License, Version 1.2).}\label{LEWUCVNXIkmdff-4MHDNOJHNFOJED}
\end{figure}

This is another classical example, due to\footnote{Lewy had also an interesting life.
After having studied mathematics and physics at G\"ottingen
and having worked extensively on partial differential equations also in collaboration with Richard Courant,
he traveled to Rome to study algebraic geometry with Tullio Levi-Civita and Federigo Enriques.
He escaped Europe to avoid racial persecutions and, after
a short-term position at Brown University, he moved to the University of
California, Berkeley and devoted himself
to partial differential equations and complex analysis.

After fifteen years, Lewy was fired from Berkeley for refusing to sign a loyalty oath imposed on the faculty
by the University of California's Board of Regents
(he was later reinstated to his job by the California Supreme Court).

See Figure~\ref{LEWUCVNXIkmdff-4MHDNOJHNFOJED} for a picture of Lewy.} Hans Lewy (see~\cite{MR88629}). The example
is in the complex setting: 
in this framework we denote by~$i:=\sqrt{-1}$
the imaginary unit (no confusion should arise when, in these pages, we will use~$i$ as an index). Also, variables in~$\cOMPL$ correspond to~$x+iy$, often identified with~$(x,y)\in\R^2$,
and an additional real variable will be denoted by~$t\in\R$.

\begin{theorem}\label{LEWU:THJSMDD}
Let~$f\in C^\infty(\R,\R)$. Assume that~$f$ is not analytic at the origin.

Then, there cannot be a smooth, complex valued, solution~$u=u(x,y,t)$ of
\begin{equation}\label{TOGDHENVSKD-20345DSmwDRHPIMGMHMS801243K} \partial_x u+i\partial_y u-2i(x+iy)\partial_t u=f(t)\end{equation}
in any neighborhood of the origin.
\end{theorem}

\begin{proof} Suppose by contradiction that such a solution~$u$ exists in a small neighborhood of the origin.
Let~${\mathcal{C}}_r$ be the circle of radius~$r>0$ in the coordinates~$(x,y)\in\R^2$,
traveled counterclockwise. 

Then, for small~$r$ and~$|t|$ we can integrate~\eqref{TOGDHENVSKD-20345DSmwDRHPIMGMHMS801243K} over~$(x,y)\in{\mathcal{C}}_r$ (for a given~$t$) and obtain
\begin{equation}\label{83i4k5r65egmetartiialy} 2\pi rf(t)=\oint_{{\mathcal{C}}_r} \Big(\partial_x u+i\partial_y u-2i(x+iy)\partial_t u\Big).\end{equation}

We introduce polar coordinates~$(\rho,\theta)$ in the~$(x,y)$ variables (in the vicinity of~${\mathcal{C}}_r$) 
by considering the function
$$ u_0(\rho,\theta,t):=u(\rho\cos\theta,\rho\sin\theta,t).$$
We also look at the quantity~$\sigma:=\ln \rho$ and the corresponding function
$$ v(\sigma,\theta,t):=u_0(e^{\sigma},\theta,t)=u(e^{\sigma}\cos\theta,e^{\sigma}\sin\theta,t).$$
Notice that
\begin{equation}\label{FFRP02oth9iuj8yg-098765-pokj9ijh-wefghny-0kj-UH}
\lim_{\rho\searrow0} v(\ln \rho, \theta,t)=
\lim_{\rho\searrow0} u(\rho\cos\theta,\rho\sin\theta,t)=u(0,0,t).
\end{equation}
We remark that
\begin{eqnarray*}
\frac{\partial v}{\partial \sigma}+i\frac{\partial v}{\partial \theta}
&=& 
e^{\sigma}\cos\theta\,\frac{\partial u}{\partial x}+e^{\sigma}\sin\theta\,\frac{\partial u}{\partial y}
-i e^{\sigma}\sin\theta\,\frac{\partial u}{\partial x}+ie^{\sigma}\cos\theta\,\frac{\partial u}{\partial y}\\&=&
e^{\sigma}(\cos\theta-i\sin\theta)\frac{\partial u}{\partial x}
+ie^{\sigma}(\cos\theta-i\sin\theta)\frac{\partial u}{\partial y}\\&=&
e^{\sigma}(\cos\theta-i\sin\theta)\left(\frac{\partial u}{\partial x}
+i\frac{\partial u}{\partial y}\right)\\&=&
e^{\sigma-i\theta}\left(\frac{\partial u}{\partial x}
+i\frac{\partial u}{\partial y}\right),
\end{eqnarray*}
where the derivatives of~$v$ are computed at~$(\sigma,\theta,t)$ and those of~$u$ at~$(e^{\sigma}\cos\theta,e^{\sigma}\sin\theta,t)$.

As a consequence,
\begin{eqnarray*}&& \frac{\partial u}{\partial x}(e^{\sigma}\cos\theta,e^{\sigma}\sin\theta,t)
+i\frac{\partial u}{\partial y}(e^{\sigma}\cos\theta,e^{\sigma}\sin\theta,t)\\&&\qquad=
e^{-\sigma+i\theta}\left(
\frac{\partial v}{\partial \sigma}(\sigma,\theta,t)+i\frac{\partial v}{\partial \theta}(\sigma,\theta,t)\right)\end{eqnarray*}
and accordingly
\begin{eqnarray*}
\oint_{{\mathcal{C}}_r} \Big(\partial_x u+i\partial_y u\Big)&=&
r\int_0^{2\pi} \Big(\partial_x u(r\cos\theta,r\sin\theta,t)+i\partial_y u(r\cos\theta,r\sin\theta,t)\Big)\,d\theta\\&=&
\int_0^{2\pi} e^{i\theta}\Big( \partial_\sigma v(\ln r,\theta,t)+i \partial_\theta v(\ln r,\theta,t)\Big)\,d\theta.
\end{eqnarray*}
Also,
$$\partial_\theta\Big(ie^{i\theta}v(\ln r,\theta,t)\Big)+e^{i\theta}v(\ln r,\theta,t)=ie^{i\theta}\partial_\theta v(\ln r,\theta,t)$$
and therefore, by the periodicity of~$v$ in the variable~$\theta$,
$$ \int_0^{2\pi}e^{i\theta}v(\ln r,\theta,t)\,d\theta=i\int_0^{2\pi} e^{i\theta}\partial_\theta v(\ln r,\theta,t)\,d\theta.$$
These observations lead to
\begin{equation*}
\oint_{{\mathcal{C}}_r} \Big(\partial_x u+i\partial_y u\Big)=
\int_0^{2\pi} e^{i\theta}\Big( \partial_\sigma v(\ln r,\theta,t)+v(\ln r,\theta,t)\Big)\,d\theta.
\end{equation*}
Combining this with~\eqref{83i4k5r65egmetartiialy}, we infer that
\begin{equation}\label{20-3pk5yt-3XFCF4intXftdt} 2\pi rf(t)=
\int_0^{2\pi} e^{i\theta}\Big( \partial_\sigma v(\ln r,\theta,t)+v(\ln r,\theta,t)\Big)\,d\theta
-2ir^2\int_0^{2\pi}e^{i\theta}\partial_t v(\ln r,\theta,t)\,d\theta.\end{equation}

Now we define (for small~$|t|$)
$$ F(t):=\pi\int_0^t f(\tau)\,d\tau$$
and (for~$Y>0$ small and small~$|X|$)
$$ V(X,Y):= i\int_0^{2\pi} e^{i\theta} \sqrt{Y}\,v(\ln \sqrt{Y}, \theta,X)\,d\theta+F(X).$$
We stress that the image of~$f$ (and therefore of~$F$) is supposed to be in the reals, therefore the imaginary part of the complex valued function~$V$ satisfies
$$ \Im V(X,Y)=\Im\left( i\int_0^{2\pi} e^{i\theta} \sqrt{Y}\,v(\ln \sqrt{Y}, \theta,X)\,d\theta\right).$$
As a result, recalling~\eqref{FFRP02oth9iuj8yg-098765-pokj9ijh-wefghny-0kj-UH}, we see that
\begin{equation}\label{20-3pk5yt-3XFCF4intXftdt-7}
\lim_{Y\searrow0} \Im V(X,Y)=0.\end{equation}

Our goal is to deduce from~\eqref{20-3pk5yt-3XFCF4intXftdt} that
\begin{equation}\label{20-3pk5yt-3XFCF4intXftdt2}
{\mbox{$V$ satisfies the Cauchy-Riemann equations in the variables~$(X,Y)$.}}\end{equation}
To this end, we calculate that
\begin{eqnarray*}&&
\partial_X V+i\partial_Y V\\&=&
i\int_0^{2\pi} e^{i\theta} \sqrt{Y}\,\partial_t v(\ln \sqrt{Y}, \theta,X)\,d\theta+\pi f(X)
-\partial_Y\int_0^{2\pi} e^{i\theta} \sqrt{Y}\,v(\ln \sqrt{Y}, \theta,X)\,d\theta\\
&=&i\int_0^{2\pi} e^{i\theta} \sqrt{Y}\,\partial_t v(\ln \sqrt{Y}, \theta,X)\,d\theta+\pi f(X)
-\int_0^{2\pi} e^{i\theta} \frac1{2\sqrt{Y}}\,v(\ln \sqrt{Y}, \theta,X)\,d\theta\\&&\qquad\qquad
-\int_0^{2\pi} e^{i\theta} \frac1{2\sqrt{Y}}\,\partial_\sigma v(\ln \sqrt{Y}, \theta,X)\,d\theta\\&=&0
\end{eqnarray*}
which establishes~\eqref{20-3pk5yt-3XFCF4intXftdt2}.

Thanks to~\eqref{20-3pk5yt-3XFCF4intXftdt-7} and~\eqref{20-3pk5yt-3XFCF4intXftdt2},
we can apply the Schwarz Reflection Principle (see e.g.~\cite[Theorem~11.14]{MR0210528}), and conclude that~$V$
can be extended to a holomorphic function (still denoted by~$V$ for simplicity) in a neighborhood of the origin.
In this way, we can think of~$V=V(X,Y)$ as a holomorphic function in~$(-r_0,r_0)^2$ for some~$r_0>0$ small enough.

As a byproduct, the function~$V(X,0)=F(X)$ is real analytic for~$X\in(-r_0,r_0)^2$.
But then also~$F'=\pi f$ is real analytic for~$X\in(-r_0,r_0)^2$, against our assumptions.
\end{proof}

{F}rom Theorem~\ref{LEWU:THJSMDD} one immediately obtains an interesting example for a system of two real equations:

\begin{corollary}
Let~$f\in C^\infty(\R,\R)$. Assume that~$f$ is not analytic at the origin.

Then, there cannot be a smooth solution~$(u_1(x,y,t),u_2(x,y,t))$, with~$u_1$ and~$u_2$ real valued, of
\begin{equation*} \begin{dcases}
\partial_x u_1-\partial_y u_2 +2y\partial_tu_1+2x\partial_t u_2=f(t),\\
\partial_x u_2+i\partial_y u_1 -2x\partial_t u_1+2y\partial_tu_2=0\end{dcases}\end{equation*}
in any neighborhood of the origin.
\end{corollary}

\begin{proof} If such a solution existed, the complex valued function~$u:=u_1+iu_2$ would be
a solution of~\eqref{TOGDHENVSKD-20345DSmwDRHPIMGMHMS801243K}, in contradiction with Theorem~\ref{LEWU:THJSMDD}.
\end{proof}

The following variation of the example in Theorem~\ref{LEWU:THJSMDD} is inspired by an example that
was put forth by
Fran\c{c}ois Tr\`eves (see~\cite{MR1566203}) to comprise the case of a single real equation:

\begin{theorem}\label{LEWU:THJSMDD-2}
Let~$f\in C^\infty(\R)$. Assume that~$f$ is not analytic at the origin.

Then, there cannot be a smooth solution~$u=u(x,y,t)$ of
\begin{equation}\label{NVINRKJHGFDFGHJ0oihbvI-08gvc3rfIJHNS} \Big(
\partial_x^2+\partial_y^2+4(x^2+y^2)\partial_t^2 +4y\partial_{xt}-4x\partial_{yt}
\Big)^2u+16\,\partial^2_t u=f(t)\end{equation}
in any neighborhood of the origin.
\end{theorem}

\begin{proof} Suppose by contradiction that such a solution exists. We denote by~${\mathcal{L}}_0$ the operator
on the left hand side of~\eqref{NVINRKJHGFDFGHJ0oihbvI-08gvc3rfIJHNS}, thus writing~${\mathcal{L}}_0u=f$.

Let also
$$ {\mathcal{L}}:=\partial_x +i\partial_y -2i(x+iy)\partial_t.$$
We observe that~${\mathcal{L}}$ is the operator dealt with in Theorem~\ref{LEWU:THJSMDD}.
Given a differential operator with complex coefficients, we denote with a superscript~``$\star$''
the operator obtained by replacing the coefficients by their complex conjugates. For example,
$$ {\mathcal{L}}^\star=\partial_x -i\partial_y +2i(x-iy)\partial_t.$$
It is worth pointing out that
$$ {\mathcal{L}}={\mathcal{L}}_1+i{\mathcal{L}}_2\qquad{\mbox{and}}\qquad {\mathcal{L}}^\star={\mathcal{L}}_1-i{\mathcal{L}}_2,$$ where
$$ {\mathcal{L}}_1=\partial_x+2y\partial_t\qquad{\mbox{and}}\qquad
{\mathcal{L}}_2=\partial_y-2x\partial_t.$$
We also observe that
\begin{eqnarray*} &&{\mathcal{M}}:={\mathcal{L}}{\mathcal{L}}^\star
=\Big({\mathcal{L}}_1+i{\mathcal{L}}_2\Big)\Big({\mathcal{L}}_1-i{\mathcal{L}}_2\Big)= {\mathcal{L}}_1^2+{\mathcal{L}}_2^2+i\Big({\mathcal{L}}_2{\mathcal{L}}_1-{\mathcal{L}}_1{\mathcal{L}}_2\Big)
\end{eqnarray*}
and therefore, since the coefficients of~${\mathcal{L}}_1$ and~${\mathcal{L}}_2$ are real,
\begin{eqnarray*} &&{\mathcal{M}}^\star=
{\mathcal{L}}_1^2+{\mathcal{L}}_2^2-i\Big({\mathcal{L}}_2{\mathcal{L}}_1-{\mathcal{L}}_1{\mathcal{L}}_2\Big).\end{eqnarray*}

As a result, setting
$$ {\mathcal{L}}_3:={\mathcal{L}}_1^2+{\mathcal{L}}_2^2\qquad{\mbox{and}}\qquad
{\mathcal{L}}_4:={\mathcal{L}}_2{\mathcal{L}}_1-{\mathcal{L}}_1{\mathcal{L}}_2,$$
we find that
$${\mathcal{M}}={\mathcal{L}}_3+i{\mathcal{L}}_4
\qquad{\mbox{and}}\qquad{\mathcal{M}}^\star={\mathcal{L}}_3-i{\mathcal{L}}_4.$$ 

Furthermore,
\begin{eqnarray*} &&{\mathcal{L}}_4=
\Big(\partial_y-2x\partial_t\Big)\Big(\partial_x+2y\partial_t\Big)-\Big(\partial_x+2y\partial_t\Big)\Big(\partial_y-2x\partial_t\Big)=4\partial_t.\end{eqnarray*}
Consequently, since both~${\mathcal{L}}_1$ and~${\mathcal{L}}_2$ commute with~$\partial_t$,
also~${\mathcal{L}}_3$ commutes with~$\partial_t$ and, for this reason,
\begin{eqnarray*}
{\mathcal{M}}{\mathcal{M}}^\star=\Big({\mathcal{L}}_3+i{\mathcal{L}}_4\Big)\Big({\mathcal{L}}_3-i{\mathcal{L}}_4
\Big)=\Big({\mathcal{L}}_3+4i\partial_t\Big)\Big({\mathcal{L}}_3-4i\partial_t
\Big)={\mathcal{L}}_3^2+16\,\partial_t^2.
\end{eqnarray*}

Additionally,
$$ {\mathcal{L}}_3=\Big(\partial_x+2y\partial_t\Big)^2+\Big(\partial_y-2x\partial_t\Big)^2=
\partial_x^2+\partial_y^2+4(x^2+y^2)\partial_t^2 +4y\partial_{xt}-4x\partial_{yt},
$$
leading to
\begin{eqnarray*}
{\mathcal{M}}{\mathcal{M}}^\star=\Big(
\partial_x^2+\partial_y^2+4(x^2+y^2)\partial_t^2 +4y\partial_{xt}-4x\partial_{yt}
\Big)^2+16\,\partial_t^2={\mathcal{L}}_0.
\end{eqnarray*}

This gives that, defining~$v:={\mathcal{L}}^\star{\mathcal{M}}^\star u$,
$$ f={\mathcal{L}}_0u={\mathcal{M}}{\mathcal{M}}^\star u=
{\mathcal{L}}{\mathcal{L}}^\star{\mathcal{M}}^\star u={\mathcal{L}}v,$$
but this is in contradiction with Theorem~\ref{LEWU:THJSMDD}.
\end{proof}

Now we introduce one of the fundamental notions for the local existence theory,
namely the concept of noncharacteristic hypersurface.
Roughly speaking, we say that a hypersurface is noncharacteristic \index{noncharacteristic hypersurface}
for the Cauchy problem~\eqref{GEN:CAUC} if the knowledge of the normal derivatives of order up to~$m-1$
along~$\Gamma$ combined with the equation involving the derivatives of order up to~$m$
allows one to determine uniquely all the derivatives of order up to~$m$ along~$\Gamma$. A more precise definition
goes as follows: for all~$\ell\in\{1,\dots,N\}$ and~$\alpha\in\N^n$ with~$|\alpha|\le\mu_\ell$
we denote by~$\xi_{\alpha \ell}$ the variable corresponding to~$\partial^\alpha u_\ell$
in the function~$F$ used in~\eqref{GABsaeM-034-002},
i.e. we write
\begin{equation}\label{BjerfFFAIKSFNaSDOFKa}
F_\ell=F_\ell\left(x, \{ \xi_{\alpha j}\}_{{{\alpha\in\N^n}\atop{|\alpha|\le\mu_j}}\atop{j\in\{1,\dots,N\}}}\right).\end{equation}
Recalling the setting in~\eqref{GABsaeM-034-001},
for every~$i$, $j\in\{1,\dots,N\}$ we define
\begin{equation}\label{LAMATRICECK} {\mathcal{M}}_{ij}:=
\sum_{{\alpha\in\N^n}\atop{|\alpha|=\mu_j}}\frac{ \partial F_i}{\partial{\xi_{\alpha j}}}\,
\left(\frac{ \partial \Phi}{\partial x_1}\right)^{\alpha_1}\dots\left(\frac{ \partial \Phi}{\partial x_n}\right)^{\alpha_n}.\end{equation}
In this framework, we can consider the matrix~${\mathcal{M}}:=\{{\mathcal{M}}_{ij}\}_{i,j\in\{1,\dots,N\}}$.
Note also that
$$ {\mathcal{M}}={\mathcal{M}}(x,\xi),\qquad{\mbox{where }}\;
\xi:=\{ \xi_{\alpha \ell}\}_{{{\alpha\in\N^n}\atop{|\alpha|\le\mu_\ell}}\atop{\ell\in\{1,\dots,N\}}}.$$
Then:

\begin{definition}\label{noncharacteristic-DE} We say that~$\Gamma$ 
is noncharacteristic for the Cauchy problem~\eqref{GEN:CAUC} if, when~$x\in\Gamma$,
\begin{equation*}
\det{\mathcal{M}}\ne0.
\end{equation*}
\end{definition}

A special, but very interesting, situation occurs when the system of equations depends linearly on the derivatives of higher order (the dependence on lower order derivatives may be nonlinear).
This is called a quasilinear system of partial differential equations
\index{quasilinear system of partial differential equations}
and in this case
$$ F_\ell=
\sum_{{{{\alpha\in\N^n}\atop{|\alpha|=\mu_j}}\atop{j\in\{1,\dots,N\}}}}
A_{\alpha \ell j}
\left(x, \{ \xi_{\beta k}\}_{{{\beta\in\N^n}\atop{|\beta|\le\mu_k-1}}\atop{k\in\{1,\dots,N\}}}\right)
\xi_{\alpha j}
+G_\ell\left(x, \{ \xi_{\beta k}\}_{{{\beta\in\N^n}\atop{|\beta|\le\mu_k-1}}\atop{k\in\{1,\dots,N\}}}\right),$$
for some~$A_{\alpha \ell j}$ and~$G_\ell$.

Accordingly, in the quasilinear case
$$ {\mathcal{M}}=
\sum_{{\alpha\in\N^n}\atop{|\alpha|=\mu_j}} A_{\alpha \ell j}\,
\left(\frac{ \partial \Phi}{\partial x_1}\right)^{\alpha_1}\dots\left(\frac{ \partial \Phi}{\partial x_n}\right)^{\alpha_n}.$$

The case of systems of second order linear elliptic equations
(see~\eqref{NONDIFORM}) can be written in the form $$
\sum_{k,h=1}^n a_{kh}(x)\partial_{kh}u(x)+\sum_{k=1}^n b_{k}(x)\partial_k u(x)+c(x)u(x)-f(x)=0.
$$
In this case, $N=1$, $\mu_1=2$ and
$$ F=\sum_{{\alpha\in\N^n}\atop{|\alpha|=2}} a_{\alpha}\xi_\alpha+
\sum_{{\alpha\in\N^n}\atop{|\alpha|=1}}b_{\alpha}\xi_\alpha+c\xi_0-f,$$
where
\begin{eqnarray*}
a_\alpha:=\begin{dcases}
a_{kh} & {\mbox{ if $|\alpha|=2$ and~$\alpha=e_k+e_h$,}}\\
a_{kk} & {\mbox{ if $|\alpha|=2$ and~$\alpha=2e_k$}}
\end{dcases}
\end{eqnarray*}
and~$b_\alpha:=b_k$ when~$|\alpha|=1$ and~$\alpha=e_k$.

In this case, ${\mathcal{M}}$ is just a scalar function, namely
\begin{eqnarray*} {\mathcal{M}}&=&\sum_{{\alpha\in\N^n}\atop{|\alpha|=2}}a_\alpha \,
\left(\frac{ \partial \Phi}{\partial x_1}\right)^{\alpha_1}\dots\left(\frac{ \partial \Phi}{\partial x_n}\right)^{\alpha_n}
\\&=&
\sum_{1\le k\ne h\le n}
a_{kh} \frac{ \partial \Phi}{\partial x_k}
\frac{ \partial \Phi}{\partial x_h}
+\sum_{k=1}^n
a_{kk} \left(\frac{ \partial \Phi}{\partial x_k}\right)^{2}
\\&=&
\sum_{1\le k, h\le n}
a_{kh} \frac{ \partial \Phi}{\partial x_k}
\frac{ \partial \Phi}{\partial x_h}.
\end{eqnarray*}
In particular, in this case the ellipticity condition in~\eqref{ELLIPTIC-PERFVDJ} guarantees that~${\mathcal{M}}\ge\lambda|\nabla\Phi|>0$, therefore
\begin{equation}\label{SEMPRECK}
{\mbox{second order linear elliptic equations are always noncharacteristic}}.
\end{equation}\medskip

Having got acquainted with the basic notation,
we can now state the local existence and uniqueness result in the analytic class
for noncharacteristic problems:

\begin{theorem}\label{CKTEO}
Let~$\Gamma$ and~$F$ be real analytic.
Let also~$\partial^j_\nu u_\ell$
be real analytic
for all~$\ell\in\{1,\dots,N\}$
and~$j\in\{0,\dots,m-1\}$.

Assume that~$\Gamma$ 
is noncharacteristic for the Cauchy problem~\eqref{GEN:CAUC}.

Then, the Cauchy problem~\eqref{GEN:CAUC} possesses
one and only one solution which is real analytic in a neighborhood of~$\Gamma$.
\end{theorem}

That is, under real analyticity and noncharacteristic assumptions,
the Cauchy problem~\eqref{GEN:CAUC} admits
one and only one solution which is real analytic in~$\Omega$ as long as~$\Omega$ is ``sufficiently small''.

We stress that Theorem~\ref{CKTEO} always applies to
second order linear elliptic equations with analytic coefficients, source and data, thanks to~\eqref{SEMPRECK}.

The importance of Theorem~\ref{CKTEO} is that properly defined, but very general,
initial value problems for partial differential equations end up admitting
a solution, at least at a local scale, and such a solution is real analytic.

However, the applicability of Theorem~\ref{CKTEO} suffers some structural limitations.
First of all, it necessitates the stringent requirement of working in the class of real analytic functions
(and this setting cannot be, in general, relaxed, given the counterexample in
Theorem~\ref{LEWU:THJSMDD-2},
see also~\cite{MR88629, MR125304, MR130574, MR142873, MR1566203, MR222800, MR257550, MR707206, MR1257229, MR2344206}
for further details).

Furthermore, while Theorem~\ref{CKTEO} guarantees the existence of a unique real analytic function,
it does not preclude the existence of other solutions which are not real analytic.
For linear systems, the analytic solution is however the unique one, see~\cite{zbMATH02662678}
and~\cite[Theorem~5.3.1]{MR1996773}, but for nonlinear equations in real analytic settings there
exist examples of other, non real analytic solutions, see~\cite{MR809715, MR1207483, MR1809287}; see also~\cite{MR217426, MR320486} for generalizations and further information on this uniqueness problem.

\begin{figure}
  \centering
  \includegraphics[width=.45\linewidth,height=.25\textwidth]{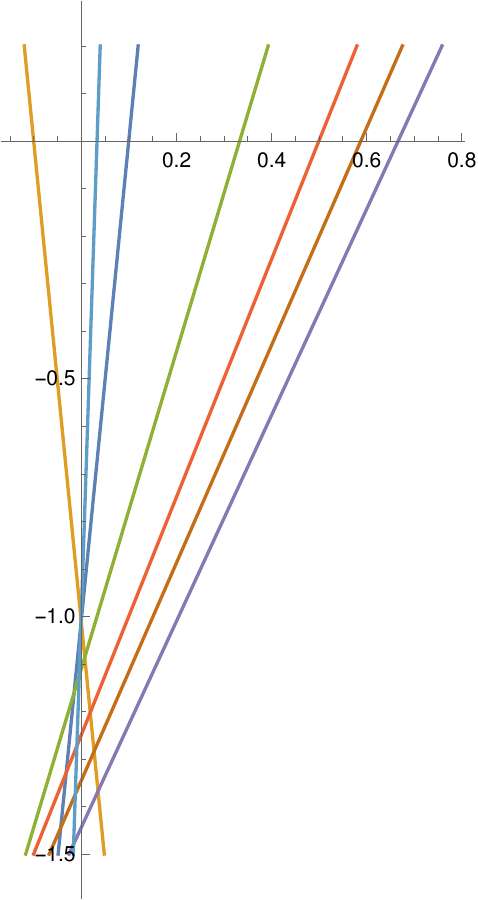}
 \caption{\sl A shock produced by the intersection of characteristic curves.}\label{SHOCACURAeFI}
\end{figure}

Also, the solution provided by Theorem~\ref{CKTEO} is only of local type (i.e., in a small neighborhood
of the hypersurface~$\Gamma$) and cannot
be, in general, taken as a global solution: for instance, even for the case of the Laplacian,
we know that overdetermined boundary value problems do not admit a global solution
(e.g., unless~$\Gamma$
is a sphere, see Theorem~\ref{SERR}). Another interesting example showing that
local solutions do not admit a global extension, not even when~$\Gamma$ is a hyperplane,
and not even when the equation is first order,
is\footnote{The example in~\eqref{GEN:CAUC:HEAT-NOTRds} provides a general paradigm 
for shock formation due to the intersection of ``characteristic curves'',
i.e. curves along which the solution maintains a constant value (see Figure~\ref{SHOCACURAeFI}
for a sketch of this phenomenon),
and coincides, up to renaming variables, with the example discussed for 
the inviscid Bateman-Burgers equation in~\eqref{VERBURG3}.} the following:

\begin{theorem}
The Cauchy problem
\begin{equation}\label{GEN:CAUC:HEAT-NOTRds}
\begin{dcases}
\partial_{x_2} u(x)=u(x)\,\partial_{x_1} u(x) & {\mbox{ for all $x=(x_1,x_2)$ in a neighborhood of }} \R\times\{0\},\\
u(x_1,0)=-\frac{x_1}{1+x_1^2}& {\mbox{ for all $x_1\in\R$.}}
\end{dcases}
\end{equation}
does not admit a real analytic solution in the whole of~$\R^2$.\end{theorem}

\begin{proof}
Notice that the setting of~\eqref{GEN:CAUC:HEAT-NOTRds} is real analytic. Also, in this situation we have that~$N:=1$, $n:=2$, $\mu_1:=1$, $\Gamma=\R\times\{0\}$,
$\Phi(x_1,x_2):=x_2$ and~$F:=\xi_{(0,1)}-\xi_{(0,0)}\xi_{(1,0)}$. Consequently, by~\eqref{LAMATRICECK},
we know that~$ {\mathcal{M}}$ is equal to~$\frac{ \partial F}{\partial{\xi_{(0,1)}}}=1\ne0$, hence the assumptions
of Theorem~\ref{CKTEO} are satisfied.

However, the solution constructed via Theorem~\ref{CKTEO} cannot be extended arbitrarily far from~$\R\times\{0\}$, because, if it were possible to extend it, setting
$$ v(t):= u\left( x_1+\frac{tx_1}{1+x_1^2},t\right)+\frac{x_1}{1+x_1^2}
\qquad{\mbox{and}}\qquad a(t):=\partial_{x_1}u\left( x_1+\frac{tx_1}{1+x_1^2},t\right),$$
then~$v(0)=0$ and, in the domain of validity of the equation,
\begin{eqnarray*}
\dot v(t)&=&\frac{x_1}{1+x_1^2} \,\partial_{x_1}u\left( x_1+\frac{tx_1}{1+x_1^2},t\right)
+\partial_{x_2}u\left( x_1+\frac{tx_1}{1+x_1^2},t\right)\\&=&
\left[\frac{x_1}{1+x_1^2} +u\left( x_1+\frac{tx_1}{1+x_1^2},t\right)\right]
\partial_{x_1}u\left( x_1+\frac{tx_1}{1+x_1^2},t\right)
\\&=& v(t)\partial_{x_1}u\left( x_1+\frac{tx_1}{1+x_1^2},t\right)\\
&=& a(t)v(t).
\end{eqnarray*}
The uniqueness result for ordinary differential equations thus entails that~$v$ vanishes identically, and therefore
\begin{equation}\label{0Suqodihw-FgtasSKYaq23dsc} u\left( x_1+\frac{tx_1}{1+x_1^2},t\right)=-\frac{x_1}{1+x_1^2}.\end{equation}
That is,
if the solution could be extended in~$(-2\delta,2\delta)\times(-1-2\delta^2,0]$
for some~$\delta>0$, then we could pick~$x_1:=\delta$ and~$t:=-1-\delta^2$ in~\eqref{0Suqodihw-FgtasSKYaq23dsc} finding that
\begin{eqnarray*} 
u(0,-1-\delta^2)=
u\left( \delta+\frac{-(1+\delta^2)\delta}{1+\delta^2},-1-\delta^2\right)=-\frac{\delta}{1+\delta^2}<0.
\end{eqnarray*}
But picking~$x_1:=0$ and~$t:=-1-\delta^2$ in~\eqref{0Suqodihw-FgtasSKYaq23dsc} we see that~$u(0,-1-\delta^2)=
0$ which is a contradiction.
\end{proof}

We also stress that the requirement that the hypersurface~$\Gamma$ is noncharacteristic cannot in general be dropped from Theorem~\ref{CKTEO}: indeed, the example in~\eqref{GEN:CAUC:HEAT}, which
does not allow for a local solution, violates the noncharacteristic requirement, since in this case~$N:=1$, $n:=2$,
$\mu_1:=2$, $\Phi(x_1,x_2):=x_2$ and~$F=\xi_{(2,0)}-\xi_{(0,1)}$. In this situation,
we deduce from~\eqref{LAMATRICECK}~$ {\mathcal{M}}$ is equal to~$
\frac{ \partial F }{\partial{\xi_{(2,0)}}}\,
\left(\frac{ \partial \Phi}{\partial x_1}\right)^2=0$, which shows that
the noncharacteristic assumption is not satisfied in this case.
\medskip

\begin{figure}
                \centering
                \includegraphics[width=.55\linewidth]{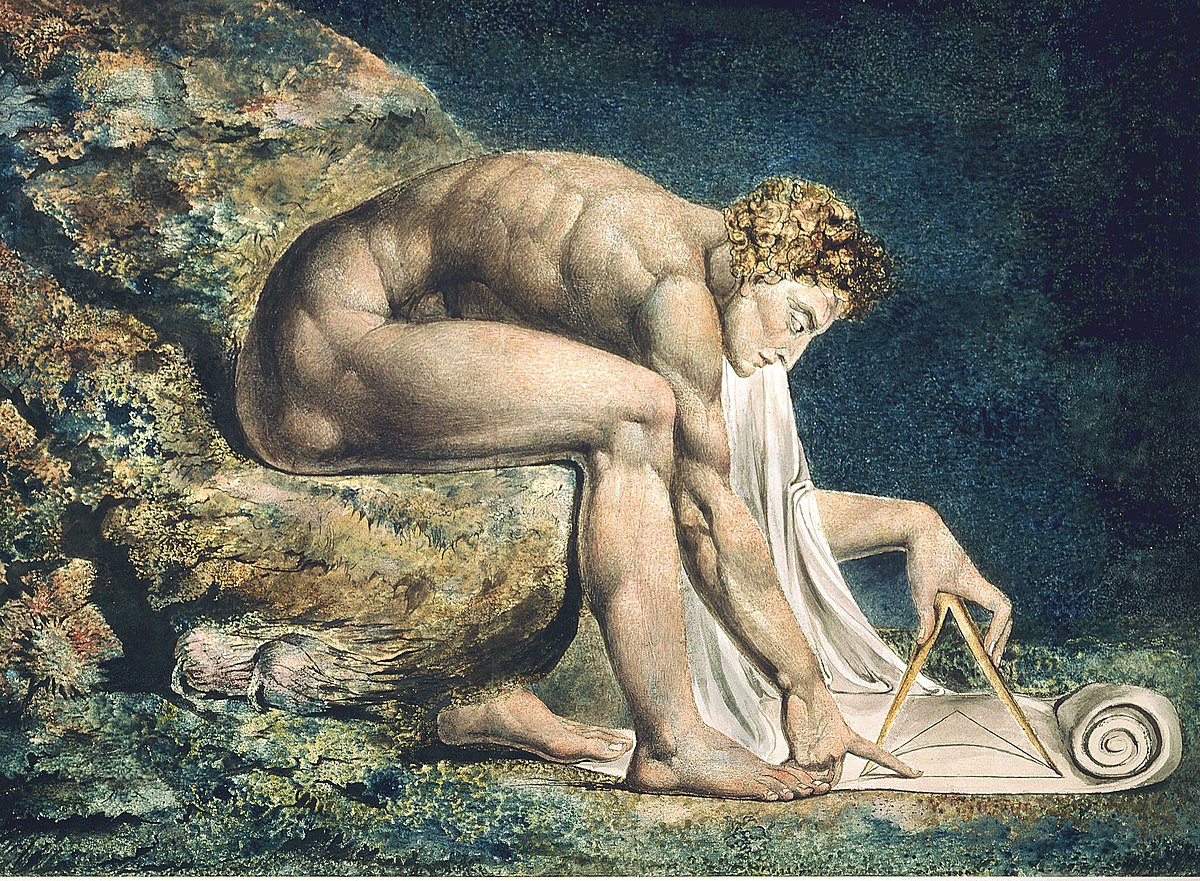}
        \caption{\sl A strapping Isaac Newton as viewed by William Blake (Public Domain image from
        Wikipedia).}\label{HAFNEWYDH-3-D}
\end{figure}

Now we present a proof of Theorem~\ref{CKTEO}. This proof will leverage Cauchy's original technique,
namely the \index{method of majorants} ``method of majorants'', which Cauchy\footnote{The word ``limites''
in this terminology used by Cauchy should be understood as ``bounds'' rather than ``limits''.

The term ``method of majorants'' was instead introduced by Poincar\'e.} called ``calcul de
limites''. The gist of this method is that, on the one hand, one can use formal series expansions
to obtain recursive formulas\footnote{The idea of finding a particular solution by\label{PALISOJA0987654AIKJ} recurrence relations
dates back to Isaac Newton and it is sometimes called the ``method of undetermined coefficients''.
At Newton's times the common belief was that analytic expressions always make sense,
the precise discussion about convergence came much later, and Cauchy was one of the prominent figures
in the introduction of a rigorous and modern approach to calculus and analysis.

See Figure~\ref{HAFNEWYDH-3-D}
for a monotype of Newton by the poet and painter William Blake:
actually, Blake clearly held a grudge against Newton, who is depicted as sitting naked and crouched
at the bottom of the sea, looking down, with his gaze directed only at some blunt geometrical diagrams,
so absorbed in his work
to become oblivious to the beautiful rocks around him and unaware of the similarity
between his drawing and the arc of his own body, as well as the shape of his own eye
(thus, Blake seems to suggest that math is incapable of understanding nature
and lacks introspection).

About Augustin-Louis Cauchy: he was a very prolific mathematician,
maybe second only to Euler, and so many theorems and concepts are named after him.

By the way, the scientific publishers at the time could have been flooded by Cauchy's
overwhelming productivity: for instance, the journal Compte Rendus de l'Acad\'emie des Sciences
instituted a four page limit for papers precisely to try to dam Cauchy's 	exuberance:
this page limit was in effect for a long time, but presently the journal seems to have become more flexible
and states in its webpage that
``there is no page limit for articles submitted to the Comptes Rendus Math\'ematique. However, the maximum recommended length for a manuscript is about 50,000 characters including spaces [...], which is equivalent to about 10 pages''.

Cauchy was also a very staunch royalist and Catholic
(in fact, a ``bigoted Catholic''
according to the prominent mathematician
Niels Henrik Abel, who however highly appreciated Cauchy's mathematical stature).
Cauchy's loyalty to the Bourbons was so firm that he never
swear the required oath of allegiance to the new regimes (preferring to
go into exile at some point).

Cauchy's cast iron opinions were broadly known
and they possibly
made him unpopular in some academic circles. For instance,
a Florentine nobleman called Guglielmo Brutus Icilius Timeleone Libri Carucci dalla Sommaja (known in France simply as
Guillaume Libri) was made chair in mathematics 
at the Coll\`ege de France
before him. However, the name of this guy is nowadays mostly recalled
for the theft of ancient and precious manuscripts:
appointed as the Inspector of Libraries in France, Libri 
fled to England in the company of more than~$3\times 10^4$ precious books. Some of the loot was
sold on the English antiques market: with a couple of auctions, Libri managed to obtain over~$10^6$ francs, in a period in which the average daily wage of a worker was of four francs
(on the top of that, Libri also had the hideous habit
of tearing off from books pages that he found of particular interest).
Well, {\em Nomen est Omen}, since ``Libri'' means ``Books'' in Italian.

One last curious anecdote about Libri: he is the one credited for popularizing the name of medieval mathematician Leonardo Bonacci, or Leonardo Bigollo Pisano, in the form that is most widespread nowadays, namely ``Fibonacci'', short for ``filius Bonacci'', son of Bonacci in Latin (even this is not that original anyway, since the name of Fibonacci to refer to Leonardo Bonacci was already present in documents from the sixteenth century).

See Figure~\ref{HAFOUMSLi:bblFUMHDNOJHNFOJED2423686900-4}
for the portraits of Cauchy and Libri (unsurprisingly, with Libri in the position of nonchalantly
grabbing a book).}
for the coefficients of the solution (if any) and, on the other hand, one can estimate these coefficients 
to rigorously establish convergence: to perform this estimate, rather than looking for specific
direct bounds, one can construct an explicit ``majorant problem'' whose relatively simple solution
automatically provides the necessary bounds for the formal solution of the original problem.

To appreciate the method of majorants in action, let us first review Cauchy's approach to ordinary differential
equations, by recalling the proof of the local existence and uniqueness result for ordinary
differential equations in a real analytic setting:

\begin{theorem}\label{CAUCHYANALODE}
Let~$(t_0,u_0)\in\R\times\R^N$ and let~$f:\R\times\R^N\to\R^N$.
Assume that~$f$ is real analytic in a neighborhood of~$(t_0,u_0)$.

Then, the initial value problem
\begin{equation}\label{OJLMS-o2456o3o}
\begin{dcases}
\dot u(t)=f(t,u(t)) & {\mbox{ for all }}t\in(t_0-\tau,t_0+\tau),\\
u(t_0)=u_0\end{dcases}
\end{equation}
possesses one and only one solution which is real analytic in a neighborhood of~$t_0$.
\end{theorem}

\begin{proof} Up to replacing~$u(t)$ with~$u(t+t_0)$, we can suppose that~$t_0=0$.
It is also convenient to define
\begin{equation}\label{INdipeott}
v(t):=(t,u(t)) \qquad{\mbox{and}}\qquad g(v):=(1,f(v)),\end{equation} thus recasting~\eqref{OJLMS-o2456o3o}
into the form
\begin{equation}\label{OJLMS-o2456o3o-b}
\begin{dcases}
\dot v(t)=g(v(t)) & {\mbox{ for all }}t\in(-\tau,\tau),\\
v(0)=v_0,\end{dcases}
\end{equation}
where~$v_0:=(0,u_0)$.

Actually, up to replacing~$v$ with~$v-v_0$ and~$g(v)$ with~$g(v-v_0)$, we can additionally assume that~$v_0=0$.

Now we establish the local existence and uniqueness theory of~\eqref{OJLMS-o2456o3o-b} by the method
of majorants. To this end, we let~$L:=N+1$ and we formally look for a solution~$v=(v_1,\dots,v_L)$ of~\eqref{OJLMS-o2456o3o-b}, where, for every~$k\in\{1,\dots,L\}$, $v_k$ is in the form
\begin{equation}\label{NLDODEWEQQ} v_k(t):=\sum_{j=0}^{+\infty} \vartheta_{k,j}\, t^j,\end{equation}
for suitable~$\vartheta_{k,j}\in\R$.

Note that, by the initial condition in~\eqref{OJLMS-o2456o3o-b},
\begin{equation}\label{NLDODEWEQQ-IN} 0=v_k(0)=\vartheta_{k,0}.
\end{equation}

We also expand~$g=(g_1,\dots,g_L)$ in its Taylor series near the origin by writing
\begin{equation}\label{NLDODEWEQQ-2}
g_k(v)=\sum_{\beta\in\N^L}\mu_{k,\beta} \,v^\beta,\end{equation}
for suitable~$\mu_{k,\beta}\in\R$.

Now we claim that
\begin{equation}\label{l0qowjfVq92LvezPOLJSdf-1p2orj3tg-304-03-ODE}
\begin{split}&
\vartheta_{k,j}={\mathcal{P}}_{k,j}\left(
\Big\{\mu_{k,\beta}\Big\}_{{{\beta\in\N^L}}\atop{|\beta| \le j-1}}, \Big\{ \vartheta_{m,i}\Big\}_{{1\le m\le L}\atop{{1\le i\le j-1}}}
\right),\\ &{\mbox{where~${\mathcal{P}}_{k,j}$ is a polynomial with {\em
nonnegative coefficients}.}}\end{split}
\end{equation}
We remark that when~$j=0$ this claim follows from~\eqref{NLDODEWEQQ-IN}, so we can suppose that~$j\ge1$.
Thus, to prove~\eqref{l0qowjfVq92LvezPOLJSdf-1p2orj3tg-304-03-ODE},
we plug~\eqref{NLDODEWEQQ} and~\eqref{NLDODEWEQQ-2} into the first line of~\eqref{OJLMS-o2456o3o-b} and we see that
\begin{eqnarray*}
\sum_{j=0}^{+\infty} (j+1) \vartheta_{k,j+1} \,t^{j}
=\sum_{j=0}^{+\infty} j \vartheta_{k,j} \,t^{j-1}=\dot v_k=g_k(v)=
\sum_{\beta\in\N^L}\mu_{k,\beta} \, v^\beta.
\end{eqnarray*}
Hence, since
\begin{eqnarray*}
v_k^{\beta_k}=\left( \sum_{j_k=0}^{+\infty} \vartheta_{k,j_k} \,t^{j_k}\right)^{\beta_k}=
\sum_{j_{k,1},\dots,j_{k,\beta_k}=0}^{+\infty} \vartheta_{k,j_{k,1}} \dots\vartheta_{k,j_{k,\beta_k}} \,t^{j_{k,1}+\dots+j_{k,\beta_k}},
\end{eqnarray*}
we conclude that
\begin{eqnarray*}
&&\sum_{j=0}^{+\infty} (j+1) \vartheta_{k,j+1} \,t^{j}
\\&&\qquad=\sum_{{\beta\in\N^L}\atop{ j_{1,1},\dots,j_{1,\beta_1},\dots,j_{L,1},\dots,j_{L,\beta_L}\in\N}}
\mu_{k,\beta} \, \vartheta_{1,j_{1,1}} \dots\vartheta_{1,j_{1,\beta_1}}\dots
\vartheta_{L,j_{L,1}} \dots\vartheta_{L,j_{L,\beta_L}} \,t^{j_{1,1}+\dots+j_{1,\beta_1}+\dots+
j_{L,1}+\dots+j_{L,\beta_L}}.\end{eqnarray*}
{F}rom this, we arrive at\footnote{Notice that, when~$j=0$, the only term in the right-hand side of~\eqref{l0qowjfVq92LvezPOLJSdf-1p2orj3tg-304-03-ODEq} that is not zero is given by~$\beta=(0,\dots,0)$,
thanks to~\eqref{NLDODEWEQQ-IN}. Hence, when~$j=0$, we have that~$\vartheta_{k,1}=\mu_{k,0}$.}
\begin{equation}\label{l0qowjfVq92LvezPOLJSdf-1p2orj3tg-304-03-ODEq} (j+1) \vartheta_{k,j+1}
=\sum_{{{\beta\in\N^L}\atop{ {j_{1,1},\dots,j_{1,\beta_1},\dots,j_{L,1},\dots,j_{L,\beta_L}}\in\N}}\atop{
{j_{1,1}+\dots+j_{1,\beta_1}+\dots+j_{L,1}+\dots+j_{L,\beta_L}}=j}}
\mu_{k,\beta} \, \vartheta_{1,j_{1,1}} \dots\vartheta_{1,j_{1,\beta_1}}\dots
\vartheta_{L,j_{L,1}} \dots\vartheta_{L,j_{L,\beta_L}},\end{equation}
for all~$j\in\N$.

Hence, to complete the proof of~\eqref{l0qowjfVq92LvezPOLJSdf-1p2orj3tg-304-03-ODE}
(with~$j$ replaced by~$j+1$), we only need to show that, with respect to the indices above,
\begin{equation}\label{l0qowjfVq92LvezPOLJSdf-1p2orj3tg-304-03-ODEp}
|\beta|\le j\qquad{\mbox{and}}\qquad j_{i,\ell}\le j,
\end{equation}
for each~$i\in\{1,\dots,L\}$ and~$\ell\in\{1,\dots,\beta_i\}$.

The second inequality in~\eqref{l0qowjfVq92LvezPOLJSdf-1p2orj3tg-304-03-ODEp} is in fact a consequence
of the fact that, in~\eqref{l0qowjfVq92LvezPOLJSdf-1p2orj3tg-304-03-ODEq},
$$j_{i,\ell}\le j_{1,1}+\dots+j_{1,\beta_1}+\dots+
j_{L,1}+\dots+j_{L,\beta_L}=j.$$
As for the first inequality in~\eqref{l0qowjfVq92LvezPOLJSdf-1p2orj3tg-304-03-ODEp},
we observe that the indices in~\eqref{l0qowjfVq92LvezPOLJSdf-1p2orj3tg-304-03-ODEq}
also satisfy~$j_{i,\ell}\ge1$ (otherwise, the corresponding term in the sum vanishes, due to~\eqref{NLDODEWEQQ-IN}) and therefore
$$|\beta|=\beta_1+\dots+\beta_L\le
j_{1,1}+\dots+j_{1,\beta_1}+\dots+
j_{L,1}+\dots+j_{L,\beta_L}=j.$$
This establishes~\eqref{l0qowjfVq92LvezPOLJSdf-1p2orj3tg-304-03-ODEp},
and thus~\eqref{l0qowjfVq92LvezPOLJSdf-1p2orj3tg-304-03-ODE}, as desired.\medskip

\begin{figure}
                \centering
                \includegraphics[height=.35\linewidth]{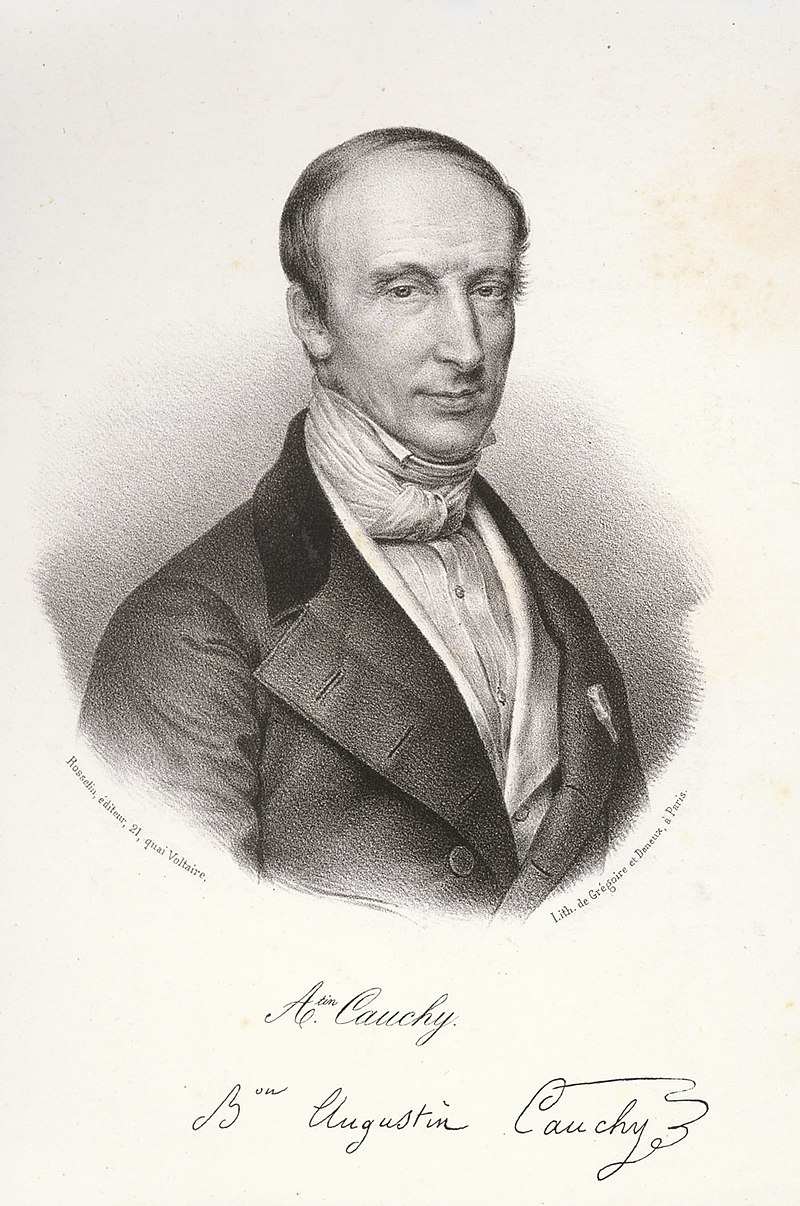}$\quad$
                \includegraphics[height=.35\linewidth]{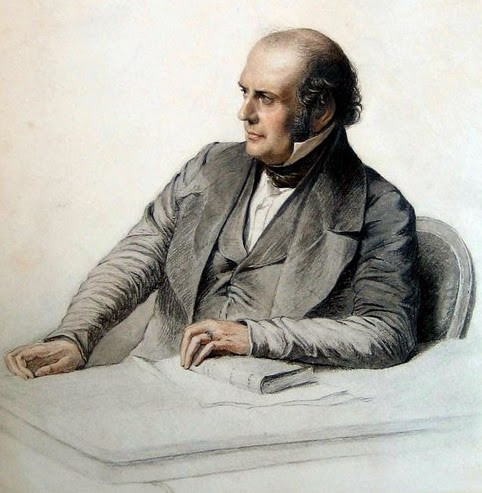}
        \caption{\sl Cauchy and Libri (Public Domain images from
        Wikipedia).}\label{HAFOUMSLi:bblFUMHDNOJHNFOJED2423686900-4}
\end{figure}

We stress that~\eqref{l0qowjfVq92LvezPOLJSdf-1p2orj3tg-304-03-ODE} determines all the coefficients~$\vartheta_{k,j}$ uniquely (i.e., by arguing recursively, starting from~\eqref{NLDODEWEQQ-IN}.
This entails that, if a real analytic solution of~\eqref{OJLMS-o2456o3o-b}
exists, it is uniquely determined by~\eqref{l0qowjfVq92LvezPOLJSdf-1p2orj3tg-304-03-ODE}.\medskip

Having established this uniqueness result for~\eqref{OJLMS-o2456o3o-b}, and thus for~\eqref{OJLMS-o2456o3o}, we
now aim at proving the existence result claimed in Theorem~\ref{CAUCHYANALODE}. To this end, we show that the formal series in~\eqref{NLDODEWEQQ} is indeed convergent. This task will be accomplished
via the method of majorants. For this, we first use~\eqref{NLDODEWEQQ-IN}
and~\eqref{l0qowjfVq92LvezPOLJSdf-1p2orj3tg-304-03-ODE} to see that
\begin{equation}\label{ILPRYHNSOIJSTGISAmeBYHEEL210-ODE} 
\vartheta_{k,j}={\mathcal{Q}}_{k,j}\left(
\Big\{\mu_{k,\beta}\Big\}_{{{\beta\in\N^L}}\atop{|\beta| \le j-1}}
\right),\end{equation}
where~${\mathcal{Q}}_{k,j}$ is a polynomial with {\em nonnegative coefficients}.

As a consequence,
\begin{equation} \label{a32-p0yrt3yt32hg943SA-09uy25zeA}
|\vartheta_{k,j}|\le{\mathcal{Q}}_{k,j}\left(
\Big\{|\mu_{k,\beta}|\Big\}_{{{\beta\in\N^L}}\atop{|\beta| \le j-1}}
\right).\end{equation}
We stress that this inequality is a consequence of the fact that the coefficients of~${\mathcal{Q}}_{k,j}$
are nonnegative. It is interesting to observe how
this monotonicity property will come in handy by allowing us to avoid the explosion of
combinatorial calculations and reduce ourselves to an explicit solution which will bound all the terms~$\vartheta_{k,j}$.

To this end, we employ the analyticity of~$g_k$ in~\eqref{NLDODEWEQQ-2}
and we deduce (see e.g.~\cite[Corollary~1.1.7]{MR1916029}) that
\begin{equation}\label{ILPRYHNSOIJSTGISAmeBYHEEL212q-NQ}
|\mu_{k,\beta}|\le \frac{C}{R^{|\beta|}} ,\end{equation}
for some~$C$, $R>0$.\medskip

To obtain a more manageable problem (see the forthcoming footnote~\ref{0ujor06-2456}), it is convenient to observe that,
for every~$\omega=(\omega_1,\dots,\omega_K)\in\N^K$,
\begin{equation}\label{ILPRYHNSOIJSTGISAmeBYHEEL212q}
\frac{(\omega_1+\dots+\omega_K)!}{\omega_1!\dots\omega_K!}\ge1.
\end{equation}
This inequality can be proved by induction over~$K$. When~$K=1$, it is actually an identity.
When~$K=2$, we can check it by setting~$\varpi:=\omega_1+\omega_2$
and noticing that~$\frac{(\omega_1+\omega_2)!}{\omega_1!\omega_2!}=
\frac{\varpi!}{\omega_1!(\varpi-\omega_1)!}=\left({{\varpi}\atop{\omega_1}}\right)\ge1$.

To perform the inductive step, we assume~\eqref{ILPRYHNSOIJSTGISAmeBYHEEL212q} for the index~$K$ and
we prove it for~$K+1$ with the following calculation:
\begin{eqnarray*}&&
\frac{(\omega_1+\dots+\omega_{K+1})!}{\omega_1!\dots\omega_{K+1}!}=
\frac{(\omega_1+\dots+\omega_{K-1}+(\omega_K+\omega_{K+1}))!}{\omega_1!\dots\omega_{K-1}!\omega_K!\omega_{K+1}!}\\&&\qquad\ge
\frac{\omega_1!\dots!\omega_{K-1}!(\omega_K+\omega_{K+1})!}{\omega_1!\dots\omega_{K-1}!\omega_K!\omega_{K+1}!}=\frac{(\omega_K+\omega_{K+1})!}{\omega_K!\omega_{K+1}!}
\ge1.
\end{eqnarray*}
The proof of~\eqref{ILPRYHNSOIJSTGISAmeBYHEEL212q} is thus complete.

Therefore, owing to~\eqref{ILPRYHNSOIJSTGISAmeBYHEEL212q-NQ} and~\eqref{ILPRYHNSOIJSTGISAmeBYHEEL212q},
\begin{equation}\label{ILPRYHNSOIJSTGISAmeBYHEEL212q-00k4rf}
|\mu_{k,\beta}|\le \frac{C\;|\beta|!}{\beta!\,R^{|\beta|}}=:\mu_{k,\beta}^\star,\end{equation}
and we can thereby majorize once again~\eqref{a32-p0yrt3yt32hg943SA-09uy25zeA},
finding that
\begin{equation} \label{a32-p0yrt3yt32hg943SA-09uy25zeAB}
|\vartheta_{k,j}|\le{\mathcal{Q}}_{k,j}\left(
\Big\{ \mu_{k,\beta}^\star\Big\}_{{{\beta\in\N^L}}\atop{|\beta| \le j-1}}
\right).\end{equation}\medskip

We now cook up a majorization problem for~\eqref{OJLMS-o2456o3o-b} by defining
\begin{equation}\label{OJLMS-o2456o3o-bMMA}
g_k^\star (v):=\sum_{\beta\in\N^L}\mu_{k,\beta}^\star \,v^\beta\end{equation}
and looking for a solution of the system of ordinary differential equations
\begin{equation}\label{OJLMS-o2456o3o-b-MAH}
\begin{dcases}
\dot v^\star(t)=g^\star(v^\star(t)) & {\mbox{ for all }}t\in(-\tau,\tau),\\
v^\star(0)=0.\end{dcases}
\end{equation}
If we seek a real analytic solution in the form
\begin{equation}\label{OJLMS-o2456o3o-b-MAH2} v_k^\star(t):=\sum_{j=0}^{+\infty} \vartheta_{k,j}^\star\, t^j,\end{equation}
we have that all the coefficients~$\vartheta_{k,j}^\star$ are determined by the prescription in~\eqref{ILPRYHNSOIJSTGISAmeBYHEEL210-ODE} , whence
\begin{equation*}
\vartheta_{k,j}^\star={\mathcal{Q}}_{k,j}\left(
\Big\{\mu_{k,\beta}^\star\Big\}_{{{\beta\in\N^L}}\atop{|\beta| \le j-1}}
\right).\end{equation*}
Combining this and~\eqref{a32-p0yrt3yt32hg943SA-09uy25zeAB} it follows that
\begin{equation}\label{901o42-X253-iy44varthetkstar}
|\vartheta_{k,j}|\le\vartheta_{k,j}^\star.\end{equation}
Thus, to complete the proof of Theorem~\ref{CAUCHYANALODE} it suffices to show that
a real analytic solution of the majorization problem in~\eqref{OJLMS-o2456o3o-b-MAH} exists:
indeed, by the uniqueness result that we have already proved, this solution necessarily coincides with the formal power series in~\eqref{OJLMS-o2456o3o-b-MAH2}; furthermore, the analyticity of this solution entails (see e.g.~\cite[Corollary~1.1.7]{MR1916029}) that~$|\vartheta_{k,j}^\star|\le\frac{C_\star}{\rho^{j}}$ for some~$C_\star$, $\rho>0$, and therefore, in view of~\eqref{901o42-X253-iy44varthetkstar}, we have that~$|\vartheta_{k,j}|\le\frac{C_\star}{\rho^{j}}$,
which in turn shows that the formal series in~\eqref{NLDODEWEQQ} is also convergent, thus providing the desired solution of~\eqref{OJLMS-o2456o3o-b} (this is the core of the method of majorants in action!).
\medskip

The advantage of the majorization problem in~\eqref{OJLMS-o2456o3o-b-MAH}
with respect to the original problem in~\eqref{OJLMS-o2456o3o-b} is that the majorization problem is explicitly solvable
(and this would complete the proof of Theorem~\ref{CAUCHYANALODE}, as we have just discussed).
To see this, we remark that
$$ g_k^\star (v)=\sum_{\beta\in\N^L} \frac{C\;|\beta|!}{\beta!\,R^{|\beta|}} \,v^\beta.$$
We also use the notation
$$ \zeta:=\frac{v_1+\dots+v_L}{R}.$$
Since, by the Multinomial Theorem, for every~$\kappa\in\N$ we have that
$$ \zeta^{\kappa}=\sum_{{\beta\in\N^L}\atop{\beta_1+\dots+\beta_L=\kappa}}
\frac{\kappa!}{\beta!\,R^\kappa}\,v_1^{\beta_1}\dots v_L^{\beta_L},$$
we find that
\begin{equation}\label{cc-19425jl0qowjfVq92LvezPOLJSdf-1p2orj3tg-304-01-ODE}
g_k^\star (v)=C\sum_{\kappa\in\N}\sum_{{\beta\in\N^L}\atop{|\beta|=\kappa}} \frac{\kappa!}{\beta!\,R^{\kappa}} \,v^\beta
=C\sum_{\kappa\in\N} \zeta^{\kappa}=\frac{C}{1-\zeta}=\frac{CR}{R-(v_1+\dots+v_L)}.
\end{equation}
As a result, the majorization problem in~\eqref{OJLMS-o2456o3o-b-MAH} becomes
\begin{equation*}
\begin{dcases}
\dot v^\star(t)=\frac{CR}{R-(v_1+\dots+v_L)} & {\mbox{ for all }}t\in(-\tau,\tau),\\
v^\star(0)=0.\end{dcases}
\end{equation*}
Since the role played by each~$v_k^\star$ is the same in this system, it is conceivable to look for solutions in which all the entries are the same, i.e. seek solutions in the form
$$ v_1^\star(t)=\dots=v_L^\star(t)=V(t),$$
for some scalar function~$V$, thus reducing the problem to
\begin{equation}\label{0-1piXDSf35erBJmd2irjfnfrv}
\begin{dcases}
\dot V(t)=\frac{CR}{R-LV} & {\mbox{ for all }}t\in(-\tau,\tau),\\
V(0)=0.\end{dcases}
\end{equation}
This scalar equation can now be solved explicitly\footnote{The \label{0ujor06-2456}
possibility of solving this problem explicitly
was indeed motivating the step of the proof between~\eqref{ILPRYHNSOIJSTGISAmeBYHEEL212q-NQ}
and~\eqref{ILPRYHNSOIJSTGISAmeBYHEEL212q-00k4rf}.} by the method of separation of variables, leading to, for~$t$ close to~$0$,
$$
RV(t)- \frac{LV^2(t)}2=\int_0^t \frac{d}{ds}\left(RV(s)- \frac{LV^2(s)}2\right)\,ds=
\int_0^s \big(R-LV(s)\big) \dot V(s)\,ds=CRt.$$
This is a quadratic equation for~$V(t)$ which can be solved (assuming that~$V(0)=0$), thus obtaining that
\begin{equation}\label{l0qowjfVq9020urVHUICVGHYUoo-01-ODE}
V(t)=\frac{R - \sqrt{R (R - 2 C L t)}}{L}.
\end{equation}

In this way, we have found explicitly a solution of~\eqref{0-1piXDSf35erBJmd2irjfnfrv}, which happens to be real analytic
for~$t$ in a small neighborhood of the origin.
The proof of Theorem~\ref{CAUCHYANALODE} is thereby complete.
\end{proof}

Having reviewed our knowledge of real analytic ordinary differential equations,
we are now ready to dive into the technical details of the proof of Theorem~\ref{CKTEO}.
In a nutshell, we will first perform some reductions of the Cauchy problem in~\eqref{GEN:CAUC}
that exploits the noncharacteristic assumption described in Definition~\ref{noncharacteristic-DE}
(in this, we utilize the invariance of the problem under real analytic transformations that
symplify the geometry under consideration)
and then we take advantage of the real analyticity of the problem.

\begin{proof}[Proof of Theorem~\ref{CKTEO}]
First of all, we utilize the noncharacteristic assumption and the Implicit
Function Theorem to write the system of equations in terms of
the derivatives of highest order (this is actually a special case of a quasilinear system). To this end,
we focus our local analysis in a small neighborhood of a point~$p\in\Gamma$
and, up to changing the order of the spatial coordinates,
we can assume that the hypersurface~$\Gamma$ is locally
described as a graph in the $n$th direction
of a real analytic function~$\Psi:\R^{n-1}\to\R$, i.e.~$\Phi(x',x_n)=x_n-\Psi(x')$,
with~$\Psi(p')=p_n$.

Thus, up to a sign change,
the normal vector can be written in the form
\begin{equation}\label{NUSUGAMMAK} \nu(x'):=-\frac{\nabla \Phi(x',x_n)}{|\nabla\Phi(x',x_n)|}=
\frac{(\nabla\Psi(x'),-1)}{\sqrt{1+|\nabla\Psi(x')|^2}}.\end{equation}
With a slight abuse of notation, we also identify normal and tangent vectors on~$\Gamma$
with their extensions as unit vector fields in a neighborhood of~$\Gamma$
by normal projection (recall~\eqref{DEEST}).\medskip

Moreover, we adopt an alternative notation for the variables of~$F$ corresponding to the derivatives of highest order. Namely,
for every~$j\in\{1,\dots,N\}$ and every~$\alpha\in\N^n$ with~$|\alpha|=\mu_j$,
\begin{equation}\label{IDENEWVACKTH}
\begin{split}&
{\mbox{we identify~$\xi_{\alpha j}$ with~$\zeta_{i_1,\dots,i_{\mu_j},j}$ with~$1\le i_1\le\dots\le i_{\mu_j}\le n$, where~$\alpha_k$
corresponds}}\\ &{\mbox{to the number of elements of~$\{i_1,\dots,i_{\mu_j}\}$ which are equal to~$k$.}}\end{split}\end{equation}
In terms of derivatives, this corresponds to the identification of~$\partial^\alpha=\partial^{\alpha_1}_{x_1}\dots\partial^{\alpha_n}_{x_n}$ with~$\partial_{i_1} \dots\partial_{ i_{\mu_j}}$.

It is also convenient to look at dummy variables~$
\zeta_{i_1,\dots,i_{\mu_j},j}$ with~$1\le i_1,\dots, i_{\mu_j}\le n$:
for these variables, we do not require anymore that the indices are ordered, we simply
identify~$\zeta_{i_1,\dots,i_{\mu_j},j}$ with~$\zeta_{i_1',\dots,i_{\mu_j}',j}$
whenever the indices~$i_1',\dots,i_{\mu_j}'$ constitute a permutation of~$i_1,\dots,i_{\mu_j}$
(for instance, with this notation, if~$j:=1$ and~$\mu_1:=3$, we have that~$\zeta_{2,2,1,1}=\zeta_{1,2,2,1}$).
\medskip

Additionally, we consider a real analytic orthonormal vector frame~$\{\tau_1,\dots,\tau_n\}$
such that~$\tau_n=\nu$ (i.e., an orthonormal frame having~$\nu$ as its last vector).
As usual, for each~$i\in\{1,\dots,n\}$, we write~$\tau_i=(\tau_{i,1},\dots,\tau_{i,n})\in\R^n$. The orthonormality
condition thus gives\footnote{More explicitly, to obtain~\eqref{ROTAORTOCK}, one can consider the matrix~$P$
with entries~$P_{ij}:=\tau_{i,j}$ and we observe that its transpose~$P^T$ coincides with its inverse because
$$ (P P^T)_{ij}=\sum_{k=1}^n \tau_{i,k}\tau_{j,k}=\tau_i\cdot\tau_j=\delta_{i,j}.$$
As a result,
$$ \delta_{i,j}=(P^T P)_{ij}=\sum_{k=1}^n \tau_{k,i}\tau_{k,j},$$
which is~\eqref{ROTAORTOCK}.} that, for all~$m$, $k\in\{1,\dots,n\}$,
\begin{equation}\label{ROTAORTOCK}
\sum_{i=1}^n \tau_{i,m} \tau_{i,k}=\delta_{m,k} .
\end{equation}

For every~$j\in\{1,\dots,N\}$ we consider  the change of variable
\begin{equation}\label{JMSdfvfsdanairfTTi6K90Pdw} \overline\zeta_{i_1,\dots,i_{\mu_j},j}:=\sum_{
1\le k_1,\dots, k_{\mu_j}\le n}\zeta_{k_1,\dots,k_{\mu_j},j} \;\tau_{i_1,k_1}\dots\tau_{i_{\mu_j},k_{\mu_j}}.\end{equation}
The interest in this relation is that
\begin{eqnarray*} \overline\zeta_{n,\dots,n,j}&=&
\sum_{1\le k_1,\dots, k_{\mu_j}\le n}\zeta_{k_1,\dots,k_{\mu_j},j} \;\tau_{n,k_1}\dots\tau_{n,k_{\mu_j}}\\&=&
\sum_{1\le k_1,\dots, k_{\mu_j}\le n}\zeta_{k_1,\dots,k_{\mu_j},j} \;\nu_{k_1}\dots\nu_{k_{\mu_j}}.
\end{eqnarray*}
That is, in the notation of~\eqref{INFEEFEBARRA},
\begin{equation}\label{INFEEFEBARRA-nh}
{\mbox{the variable~$\overline\zeta_{n,\dots,n,j}$ corresponds to the normal derivative
of order~$\mu_j$, namely~$\partial^{\mu_j}_\nu$.}}
\end{equation}
We also stress that~\eqref{JMSdfvfsdanairfTTi6K90Pdw} is indeed an invertible change of variable, since, by~\eqref{ROTAORTOCK},
\begin{equation}\label{90oihf932y03275yhtgfvcb98vnct0987654v867h843u9uyhCVBXnceir-94} \begin{split}
\zeta_{m_1,\dots,m_{\mu_j},j}&=\sum_{{1\le k_1,\dots, k_{\mu_j}\le n}}
\zeta_{k_1,\dots,k_{\mu_j},j} \;\delta_{m_1,k_1}\dots\delta_{m_{\mu_j},k_{\mu_j}}\\
&=\sum_{{1\le k_1,\dots, k_{\mu_j}\le n}\atop{1\le i_1,\dots, i_{\mu_j}\le n}}
\zeta_{k_1,\dots,k_{\mu_j},j} \;\tau_{i_1,k_1}\dots\tau_{i_{\mu_j},k_{\mu_j}}
\;\tau_{i_1,m_1}\dots\tau_{i_{\mu_j},m_{\mu_j}}
\\&=\sum_{1\le i_1,\dots, i_{\mu_j}\le n}\overline\zeta_{i_1,\dots,i_{\mu_j},j}
\;\tau_{i_1,m_1}\dots\tau_{i_{\mu_j},m_{\mu_j}}.\end{split}\end{equation}
An interesting byproduct of~\eqref{90oihf932y03275yhtgfvcb98vnct0987654v867h843u9uyhCVBXnceir-94}
is that
\begin{equation}\label{90oihf932y03275yhtgfvcb98vnct0987654v867h843u9uyhCVBXnceir-94-B} \begin{split}
\frac{\partial \zeta_{m_1,\dots,m_{\mu_j},j}}{\partial\overline\zeta_{n,\dots,n,j}}
=\tau_{n,m_1}\dots\tau_{n,m_{\mu_j}}=\nu_{m_1}\dots\nu_{m_{\mu_j}}
.\end{split}\end{equation}

With these observations, we can define~$\overline F$ as~$F$ in the new set of coordinates given by~\eqref{JMSdfvfsdanairfTTi6K90Pdw}.

Therefore, recalling~\eqref{INFEEFEBARRA-nh},
in order to make it explicit the equation under consideration in terms of
the highest order normal derivatives, we need to use the Implicit Function Theorem on the function~$\overline F$
with respect to the set of variables~$(\overline\zeta_{n,\dots,n,1},\dots,\overline\zeta_{n,\dots,n,N})$.
To this end, in view of~\eqref{90oihf932y03275yhtgfvcb98vnct0987654v867h843u9uyhCVBXnceir-94-B} we observe that, for all~$i$, $j\in\{1,\dots,N\}$,
\begin{eqnarray*} \frac{\partial \overline F_i}{\partial \overline\zeta_{n,\dots,n,j}}=
\sum_{1\le m_1,\dots,m_{\mu_j}\le n} \frac{\partial  F_i}{\partial \zeta_{m_1,\dots,m_{\mu_j},j}}\,
\frac{\partial \zeta_{m_1,\dots,m_{\mu_j},j}}{\partial \overline\zeta_{n,\dots,n,j}}=
\sum_{1\le m_1,\dots,m_{\mu_j}\le n} \frac{\partial  F_i}{\partial \zeta_{m_1,\dots,m_{\mu_j},j}}\,
\nu_{m_1}\dots\nu_{m_{\mu_j}}.\end{eqnarray*}
Hence, recalling~\eqref{IDENEWVACKTH},
\begin{equation*}
\frac{\partial \overline F_i}{\partial \overline\zeta_{n,\dots,n,j}}=
\sum_{{\alpha\in\N^n}\atop{|\alpha|=\mu_j}} \frac{\partial  F_i}{\partial \xi_{\alpha j}}\,
\nu_1^{\alpha_1} \dots\nu_n^{\alpha_n}.
\end{equation*}
As a result, in light of~\eqref{LAMATRICECK}, Definition~\ref{noncharacteristic-DE}
and~\eqref{NUSUGAMMAK},
we infer that the matrix~$\left\{\frac{\partial \overline F_i}{\partial \overline\zeta_{n,\dots,n,j}}\right\}_{i,j\in\{1,\dots ,N\}}$ is invertible and, consequently, by the real analytic version of the Implicit Function Theorem
(see e.g.~\cite{MR1916029}), the system of equations~$F(x,D^\mu_\star u(x))=0$ is equivalent
to
\begin{equation}\label{78o00-4k3-2tjmg}
\partial^{\mu_\ell}_\nu u_\ell=G_\ell,\end{equation} for~$\ell\in\{1,\dots,N\}$,
where~$G_\ell$ is a real analytic function of~$x$ and of the derivatives of~$(u_1,\dots,u_N)$
in a tangential/normal frame of reference, being the order of the derivatives of~$u_m$ less than or equal to~$\mu_m$ and the order of the derivatives of~$u_m$ in the normal direction less than or equal to~$\mu_m-1$.
Recall also that~\eqref{78o00-4k3-2tjmg} is complemented with the boundary conditions in~\eqref{GEN:CAUC}.
\medskip

Now, we want to reduce to the case in which~$\Gamma$ is a hyperplane.
For this, we aim at straightening~$\Gamma$ by using a new coordinate system~$(y_1,\dots,y_{n-1},t)$
in which the coordinates~$(y_1,\dots,y_{n-1})$ parameterize~$\Gamma$ and correspond to
the hyperplane~$\{t=0\}$ and the variable~$t$ identifies displacements in the direction normal to~$\Gamma$.
This reference system is precisely that of \label{FECOORD}
``Fermi Coordinates'', which \index{Fermi Coordinates} were mentioned on page~\pageref{PAGFermiCoordinates}.

Namely, we consider the change of coordinates near the origin described by
\begin{equation}\label{OASK:PAS} x=x(y,t)= (y,\Psi(y))+t\nu(y).\end{equation}
We remark that
$$ \det \left(\left.\frac{\partial x}{\partial (y,t)}\right|_{(y,t)=(0,0)}\right)=-\sqrt{1+|\nabla \Psi(0) |^2}\neq0
,$$
guaranteeing that~\eqref{OASK:PAS} is indeed a real analytic change of variable in a neighborhood of the origin
(see e.g.~\cite{MR1916029} for the Inverse Function Theorem in the real analytic setting).

Therefore, we can consider the function~$\bar{u}(y,t):=u(x(y,t))$ and set~\eqref{78o00-4k3-2tjmg} and the boundary conditions in~\eqref{GEN:CAUC}
into the equivalent problem
\begin{equation}\label{GEN:CA:PIAT1}\begin{dcases}
\partial^{\mu_\ell}_t u_\ell=H_\ell\left( y,t,
\Big\{ \partial^\alpha_{y}\partial^j_t u_m \Big\}_{{m\in\{1,\dots,N\}}\atop{{\alpha\in\N^{n-1},j\in\N}\atop{{|\alpha|+j\le\mu_m}\atop{j\le\mu_m-1}}}}
\right) & {\mbox{ for all~$|y|<\rho$, $|t|<\rho$ and~$\ell\in\{1,\dots,N\}$,}}\\
\partial^j_t u_\ell(y,0)=
\psi_{j,\ell}(y)& {\mbox{ for all $|y|<\rho$, $\ell\in\{1,\dots,N\}$
and $j\in\{0,\dots,\mu_\ell-1\}$,}}
\end{dcases}
\end{equation}
for some~$\rho>0$ (to be chosen conveniently small) and given real analytic functions~$H_\ell$ and~$\psi_{j,\ell}$.
\medskip

Actually, the boundary conditions in~\eqref{GEN:CA:PIAT1} naturally imply the knowledge of
all the derivatives
$$\Big\{ \partial^\alpha_{y}\partial^j_t u_m \Big\}_{{m\in\{1,\dots,N\}}\atop{{\alpha\in\N^{n-1},j\in\N}\atop{{|\alpha|+j\le\mu_m}\atop{j\le\mu_m-1}}}}$$ at~$t=0$, simply by differentiating in the direction of~$y_i$, for each~$i\in\{1,\dots,n-1\}$: more explicitly,
$$ \partial^\alpha_{y}\partial^j_t u_\ell(y,0)=
\partial^\alpha_{y}\psi_{j,\ell}(y)=:\Psi_{\ell,j,\alpha}(y)\qquad {\mbox{ for all $|y|<\rho$, $\ell\in\{1,\dots,N\}$
and $j\in\{0,\dots,\mu_\ell-1\}$.}}$$

As a consequence, the Cauchy problem in~\eqref{GEN:CA:PIAT1} is equivalent to
\begin{equation}\label{GEN:CA:PIAT2p}\begin{dcases}
\partial^{\mu_\ell}_t u_\ell=H_\ell\left( y,t,
\Big\{ \partial^\alpha_{y}\partial^j_t u_m \Big\}_{{m\in\{1,\dots,N\}}\atop{{\alpha\in\N^{n-1},j\in\N}\atop{{|\alpha|+j\le\mu_m}\atop{j\le\mu_m-1}}}}
\right) & {\mbox{ for all~$|y|<\rho$, $|t|<\rho$ and~$\ell\in\{1,\dots,N\}$,}}\\
\partial^\alpha_y\partial^j_t u_\ell(y,0)=
\Psi_{\ell,j,\alpha}(y)& {\mbox{ for all $|y|<\rho$, $\ell\in\{1,\dots,N\}$
and $j\in\{0,\dots,\mu_\ell-1\}$,}}
\end{dcases}
\end{equation}
and we can limit ourselves to consider in the latter equation the indices~$\alpha\in\N^{n-1}$ with~$|\alpha|\le\mu_\ell$.

As a matter of fact, we can also write
\begin{equation}\label{0o2erh0834yt0348ieoh943iegbhhIHBD} \partial^{\mu_\ell}_t u_\ell\big|_{t=0}=\left. H_\ell\left( y,t,
\Big\{ \partial^\alpha_{y}\partial^j_t u_m \Big\}_{{m\in\{1,\dots,N\}}\atop{{\alpha\in\N^{n-1},j\in\N}\atop{{|\alpha|+j\le\mu_m}\atop{j\le\mu_m-1}}}}
\right)\right|_{t=0}=:\Psi_{\ell,\mu_\ell,0},\end{equation}
and we stress that the latter quantity is determined by the boundary data, therefore~\eqref{GEN:CA:PIAT2p}
can also be rewritten in the form
\begin{equation}\label{GEN:CA:PIAT2x}\begin{dcases}
\partial^{\mu_\ell}_t u_\ell=H_\ell\left( y,t,
\Big\{ \partial^\alpha_{y}\partial^j_t u_m \Big\}_{{m\in\{1,\dots,N\}}\atop{{\alpha\in\N^{n-1},j\in\N}\atop{{|\alpha|+j\le\mu_m}\atop{j\le\mu_m-1}}}}
\right) & {\mbox{ for all~$|y|<\rho$, $|t|<\rho$ and~$\ell\in\{1,\dots,N\}$,}}\\
\partial^\alpha_y\partial^j_t u_\ell(y,0)=
\Psi_{\ell,j,\alpha}(y)& {\mbox{ for all $|y|<\rho$, $\ell\in\{1,\dots,N\}$
and $|\alpha|+j\le\mu_\ell$.}}
\end{dcases}
\end{equation}
\medskip

One can actually reduce~\eqref{GEN:CA:PIAT2x} to a quasilinear system of equations.
To this end, we let~$a_{\ell,m,j,\alpha}$ to be the derivative of~$H_\ell$ with respect to the variable corresponding to~$\partial^\alpha_{y}\partial^j_t u_m$ and we consider the system
\begin{equation}\label{GEN:CA:PIAT2l}\begin{dcases}
\partial^{\mu_\ell+1}_t u_\ell=
\partial_t H_\ell+
\sum_{{m\in\{1,\dots,N\}}\atop{{\beta\in\N^{n-1},j\in\N}\atop{{|\beta|+j\le\mu_m}\atop{j\le\mu_m-1}}}} a_{\ell,m,j,\beta}\,\partial^\beta_{y}\partial^{j+1}_t u_m & {\mbox{ for all~$|y|<\rho$, $|t|<\rho$ and~$\ell\in\{1,\dots,N\}$,}}\\
\partial^\alpha_y\partial^j_t u_\ell(y,0)=
\Psi_{\ell,j,\alpha}(y)& {\mbox{ for all $|y|<\rho$, $\ell\in\{1,\dots,N\}$
and $|\alpha|+j\le\mu_\ell$.}}
\end{dcases}
\end{equation}
Indeed, if~\eqref{GEN:CA:PIAT2x} holds true, then we obtain~\eqref{GEN:CA:PIAT2l} simply by taking a derivative in~$t$ of the first line of~\eqref{GEN:CA:PIAT2x}.

Conversely, if~\eqref{GEN:CA:PIAT2l} holds true, then its first line says that~$ \partial_t\left(\partial^{\mu_\ell}_t
u_\ell-H_\ell
\right)=0$ and accordingly, for all~$|t|<\rho$ and~$|y|<\rho$,
we infer that~$\partial^{\mu_\ell}_t u_\ell-H_\ell$ is equal to some function of~$y$ only. However, this function is necessarily zero, due to~\eqref{0o2erh0834yt0348ieoh943iegbhhIHBD} and the boundary condition
contained in the second line of~\eqref{GEN:CA:PIAT2l}.

The equivalence between~\eqref{GEN:CA:PIAT2x} and~\eqref{GEN:CA:PIAT2l} is thereby established.
\medskip

Additionally, we can reduce~\eqref{GEN:CA:PIAT2l} to a Cauchy problem with homogeneous data.
Indeed, up to replacing~$u_\ell(y,t)$ with
$$ u_\ell(y,t)-\sum_{i=0}^{\mu_\ell}\frac{\Psi_{\ell,i,0}(y)}{i!}\,t^i,$$
and eventually modifying~$H_\ell$,
we reduce~\eqref{GEN:CA:PIAT2l} to 
\begin{equation}\label{GEN:CA:PIAT2}\begin{dcases}
\partial^{\mu_\ell+1}_t u_\ell=
\partial_t H_\ell+
\sum_{{m\in\{1,\dots,N\}}\atop{{\beta\in\N^{n-1},j\in\N}\atop{{|\beta|+j\le\mu_m}\atop{j\le\mu_m-1}}}} a_{\ell,m,j,\beta}\,\partial^\beta_{y}\partial^{j+1}_t u_m & {\mbox{ for all~$|y|<\rho$, $|t|<\rho$ and~$\ell\in\{1,\dots,N\}$,}}\\
\partial^\alpha_y\partial^j_t u_\ell(y,0)=0& {\mbox{ for all $|y|<\rho$, $\ell\in\{1,\dots,N\}$
and $|\alpha|+j\le\mu_\ell$.}}
\end{dcases}
\end{equation}
\medskip

It is now convenient to reduce~\eqref{GEN:CA:PIAT2} to a first order system.
For this, if~$u_\ell$ is as in~\eqref{GEN:CA:PIAT2}, whenever~$|\alpha|+j\le\mu_\ell$ we
introduce the functions
\begin{equation}\label{VUELLEALFA}
v_{\ell,j,\alpha}:=\partial^\alpha_y\partial^j_t u_\ell\end{equation}
and consider as the new unknown vectorial solution the function
\begin{equation*} v:=\Big\{ v_{\ell,j,\alpha}\Big\}_{{\ell\in\{1,\dots,N\}}\atop{{\alpha\in\N^{n-1},j\in\{1,\dots,\mu_\ell\}}\atop{{|\alpha|\le\mu_\ell}}}}.\end{equation*}
For each~$\alpha=(\alpha_1,\dots,\alpha_n)\in\N^n\setminus\{0\}$ we also let~$i(\alpha)\in\{1,\dots,N\}$ the smallest index for which~$\alpha_i\ne0$.

In this notation, we claim that~\eqref{GEN:CA:PIAT2} is equivalent to the first order quasilinear system\footnote{Yep, the notation in~\eqref{GEN:CA:PIAT2LANUO} is quite cumbersome. However, let us try to understand what's going on in a specific example. Suppose that we are considering the scalar equation
$$ \partial^3_t u=\partial^2_y \partial_t u+\partial_y u+\partial^2_t u,$$
with~$y\in\R$.

The strategy would be to define
$$ v_1:=u,\qquad v_2:=\partial_yu,\qquad v_3:=\partial^2_yu,\qquad
v_4:=\partial_t u,\qquad v_5:=\partial^2_tu\qquad{\mbox{and}}\qquad v_6:=\partial_y\partial_t u.$$
This reduces the scalar equation above to the first order system
$$\begin{dcases}
\partial_t v_1=\partial_t u=v_4,\\
\partial_t v_2=\partial_y\partial_t u=v_6,\\
\partial_t v_3=\partial_y^2 \partial_t u=\partial_y v_6,\\
\partial_t v_4=\partial_t^2u=v_5,\\
\partial_t v_5=\partial_t^3 u=
\partial^2_y \partial_t u+\partial_y u+\partial^2_t u=
\partial_y v_6+v_2+v_5,\\
\partial_tv_6=\partial_y\partial_t^2 u=\partial_y v_5.
\end{dcases}$$
The framework in~\eqref{GEN:CA:PIAT2LANUO} is nothing more than this strategy,
applied to a general system of equations.}
\begin{equation}\label{GEN:CA:PIAT2LANUO}\begin{dcases}
\partial_t v_{\ell,j,\alpha}=\begin{dcases}
v_{\ell, j+1,\alpha} &\qquad {\mbox{ if }}j\le\mu_\ell-1\\
\\
\begin{matrix}
\displaystyle\partial_t H_\ell+\sum_{{m\in\{1,\dots,N\}}\atop{{\beta\in\N^{n-1},j\in\{1,\dots,\mu_\ell\}}\atop{{|\beta|+j\le\mu_m-1}}}} a_{\ell,m,j,\beta}\, v_{m,j+1,\beta} \cr\displaystyle\qquad
+\sum_{{m\in\{1,\dots,N\}}\atop{{\beta\in\N^{n-1},j\in\{1,\dots,\mu_\ell\}}\atop{{|\beta|+j=\mu_m}\atop{j\le\mu_m-1}}}} a_{\ell,m,j,\beta}\,\partial_{y_{i(\beta)}} v_{m,j+1,\beta-e_{i(\beta)}}\end{matrix}
& \qquad{\mbox{ if }}j=\mu_\ell, \end{dcases} \\ \qquad\qquad\qquad {\mbox{ for all~$|y|<\rho$, $|t|<\rho$, $\ell\in\{1,\dots,N\}$
and $|\alpha|+j\le\mu_\ell$,}}\\
\\
\\
v_{\ell,j,\alpha}(y,0)=0\qquad {\mbox{ for all $|y|<\rho$, $\ell\in\{1,\dots,N\}$
and $|\alpha|+j\le\mu_\ell$,}}
\end{dcases}
\end{equation}
where, in light of~\eqref{GEN:CA:PIAT2x}, the variables of~$\partial_t H_\ell$ and~$a_{\ell,m,i,\beta}$ are now
$$ \left( y,t,\Big\{v_{m,i,\gamma} \Big\}_{{m\in\{1,\dots,N\}}\atop{{\gamma\in\N^{n-1},i\in\N}\atop{{|\gamma|+i\le\mu_m}\atop{i\le\mu_m-1}}}}
\right).$$

Indeed, if we have a solution of~\eqref{GEN:CA:PIAT2} then we obtain a solution of~\eqref{GEN:CA:PIAT2LANUO}
simply by the setting in~\eqref{VUELLEALFA}.

Suppose now that~\eqref{GEN:CA:PIAT2LANUO} holds true and choose
\begin{equation}\label{DGGEN:CA:PIAT2LANUO}
u_\ell(y,t):= \int_0^t\left( \int_0^{\tau_{\mu_\ell}}\left(\dots\left(\int_0^{\tau_2} v_{\ell,{\mu_\ell},0}(y,\tau_1)\,d\tau_1\right)\dots\right)\,d\tau_{{\mu_\ell}-1}\right)\,d\tau_{\mu_\ell}.\end{equation}
We claim that
\begin{equation}\label{DGGEN:CA:PIAT2LANUO2}
{\mbox{whenever~$\ell\in\{1,\dots,N\}$, $j\in\{0,\dots,\mu_\ell\}$ and~$|\alpha|+j\le\mu_\ell$, equation~\eqref{VUELLEALFA} is satisfied.}}\end{equation} 
Notice that once~\eqref{DGGEN:CA:PIAT2LANUO2} is proven, we obtain as a byproduct that~\eqref{GEN:CA:PIAT2} holds true.

We prove~\eqref{DGGEN:CA:PIAT2LANUO2} by induction over~$p:=\mu_\ell-j\in\{0,\dots,\mu_\ell\}$. For this, let us first suppose that~$p=0$. Then, $j=\mu_\ell$ and,
by taking ${\mu_\ell}$ derivatives in~$t$ in~\eqref{DGGEN:CA:PIAT2LANUO},
we have that~$\partial_t^{\mu_\ell} u_\ell(y,t)=v_{\ell,\mu_\ell,0}(y,t)$,
which is~\eqref{VUELLEALFA} for~$j=\mu_\ell$.

Hence, to perform the inductive step, we now suppose that~\eqref{DGGEN:CA:PIAT2LANUO2}
holds true for all~$p\le p_0$ with~$p_0\le\mu_\ell-1$ and we aim at proving it for~$p_0+1$, i.e. for~$j=\mu_\ell-p_0-1$. To this end, we take~$\alpha\in\N^n$ with~$|\alpha|+j\le\mu_\ell$ and we differentiate~\eqref{DGGEN:CA:PIAT2LANUO} to see that
\begin{equation}\label{9.1.32}
\partial^\alpha_y\partial^j_t u_\ell(y,t)= \int_0^t\left( \int_0^{\tau_{\mu_\ell-j}}\left(\dots\left(\int_0^{\tau_2} \partial^\alpha_y v_{\ell,{\mu_\ell},0}(y,\tau_1)\,d\tau_1\right)\dots\right)\,d\tau_{{\mu_\ell}-j-1}\right)\,d\tau_{\mu_\ell-j}.
\end{equation}
We also have that~$j=\mu_\ell-p_0-1\le\mu_\ell-1$. Hence, by~\eqref{GEN:CA:PIAT2LANUO} and the inductive assumption we see that
$$ \partial_t v_{\ell,j,\alpha}=v_{\ell, j+1,\alpha}=\partial^\alpha_y\partial^{j+1}_{t} u_\ell=
\partial_t(\partial^\alpha_y\partial^j_{t} u_\ell).
$$
Accordingly, the function~$v_{\ell,j,\alpha}-\partial_y^{\alpha}\partial^{j}_t u_{\ell}$ is constant in~$t$.
Also, this function  vanishes at~$t=0$, due to the boundary condition in~\eqref{GEN:CA:PIAT2LANUO}
and the expression in~\eqref{9.1.32}. As a consequence, the function~$v_{\ell,j,\alpha}-\partial_y^{\alpha}\partial^{j}_t u_{\ell}$ vanishes identically, giving~\eqref{VUELLEALFA} in this case.
The proof of~\eqref{DGGEN:CA:PIAT2LANUO2} is thereby complete,
hence we have reduced the proof of  Theorem~\ref{CKTEO}
to the analysis of~\eqref{GEN:CA:PIAT2LANUO}.\medskip

So, owing to~\eqref{GEN:CA:PIAT2LANUO}, to prove Theorem~\ref{CKTEO} (using now an easier notation) it suffices to establish, in the real analytic framework, a local existence and uniqueness theory for quasilinear systems of first order with homogeneous boundary data of the form
\begin{equation}\label{9ojknTRvbdfLON0oejgDAnmdcVEujdfrtitghbsdfAJSdfwQUW}\begin{dcases}
\partial_t w_k=\sum_{{1\le i\le n-1}\atop{1\le m\le M}} b_{i,m,k}(y,t,w)\partial_i w_m+c_k(y,t,w)
& {\mbox{ for all~$|y|<\rho$, $|t|<\rho$, $k\in\{1,\dots,M\}$,}}\\
w_k(y,0)=0& {\mbox{ for all $|y|<\rho$ and $k\in\{1,\dots,M\}$,}}
\end{dcases}
\end{equation}
where the unknown vectorial solution is~$w(y,t)=(w_1(y,t),\dots,w_M(y,t))$.\medskip

One can actually suppose, up to adding one more equation and unknown, that the coefficients~$b_{i,m,k}$
and~$c_k$ do not depend on~$t$ (this trick\footnote{Here,
one could also get rid of the dependence on~$y$, again by adding~$n-1$ new equations and unknowns,
but this would introduce a nonhomogeneous boundary term at~$t=0$, therefore
we preferred to maintain the dependence on~$y$ in the coefficients.} was already exploited in~\eqref{INdipeott}
for the case of ordinary differential equations). Namely, by defining~$L:=M+1$
and~$w_L(y,t):=t$, we have that~$\partial_t w_L=1$ and the system in~\eqref{9ojknTRvbdfLON0oejgDAnmdcVEujdfrtitghbsdfAJSdfwQUW} can be written in the form
\begin{equation}\label{9ojknTRvbdfLON0oejgDAnmdcVEujdfrtitghbsdfAJSdfwQUW-2}\begin{dcases}
\partial_t w_k=\sum_{{1\le i\le n-1}\atop{1\le m\le L}} B_{i,m,k}(y,w)\partial_i w_m+C_k(y,w)
& {\mbox{ for all~$|y|<\rho$, $|t|<\rho$, $k\in\{1,\dots,L\}$,}}\\
w_k(y,0)=0& {\mbox{ for all $|y|<\rho$ and $k\in\{1,\dots,L\}$,}}
\end{dcases}
\end{equation}
where the unknown vectorial solution is~$w(y,t)=(w_1(y,t),\dots,w_L(y,t))$.\medskip

We now address the analytic part of the proof of Theorem~\ref{CKTEO}, which relies on the method of majorants
and establishes the local unique solvability of~\eqref{9ojknTRvbdfLON0oejgDAnmdcVEujdfrtitghbsdfAJSdfwQUW-2} in the real analytic setting. To this end,
we use the notation~$Y:=(y,t)$ and seek~$w=(w_1,\dots,w_L)$, with~$w_k$ in the form\footnote{The counterpart of~\eqref{l0qowjfVq92LvezPOLJSdf-1p2orj3tg-304-01} in the setting of ordinary differential equations
was considered in~\eqref{NLDODEWEQQ}.}
\begin{equation}\label{l0qowjfVq92LvezPOLJSdf-1p2orj3tg-304-01} w_k(Y)=\sum_{\alpha\in\N^n} \vartheta_{k,\alpha}\, Y^\alpha,\end{equation}
for suitable~$\vartheta_{k,\alpha}\in\R$. Note also that, by the boundary condition in~\eqref{9ojknTRvbdfLON0oejgDAnmdcVEujdfrtitghbsdfAJSdfwQUW-2}, for every~$\alpha'\in\N^{n-1}$,
\begin{equation}\label{Ojnd-1-Sdfgo0o2ruiytghb-0oritgkfkS}
\vartheta_{k,(\alpha',0)}=0.
\end{equation}
We also expand~$B_{i,m,k}$ and~$C_k$ as\footnote{The counterpart of~\eqref{l0qowjfVq92LvezPOLJSdf-1p2orj3tg-304-02} in the setting of ordinary differential equations
was considered in~\eqref{NLDODEWEQQ-2}.}
\begin{equation}\label{l0qowjfVq92LvezPOLJSdf-1p2orj3tg-304-02} B_{i,m,k}(y,w)=\sum_{{\beta\in\N^L}\atop{\sigma\in\N^{n-1}}} \lambda_{i,m,k,\beta,\sigma} \,y^\sigma w^\beta
\qquad{\mbox{and}}\qquad C_k(y,w)=\sum_{{\beta\in\N^L}\atop{\sigma\in\N^{n-1}}} \mu_{k,\beta,\sigma} \,y^\sigma w^\beta,\end{equation}
for suitable~$\lambda_{i,m,k,\beta,\sigma}$, $\mu_{k,\beta,\sigma}\in\R$.

We claim that\footnote{The counterpart of~\eqref{l0qowjfVq92LvezPOLJSdf-1p2orj3tg-304-03} in the setting of ordinary differential equations
was considered in~\eqref{l0qowjfVq92LvezPOLJSdf-1p2orj3tg-304-03-ODE}.}
\begin{equation}\label{l0qowjfVq92LvezPOLJSdf-1p2orj3tg-304-03}
\begin{split}&
\vartheta_{k,\alpha}={\mathcal{P}}_{k,\alpha}\left(
\Big\{\lambda_{i,m,k,\beta,\sigma}\Big\}_{{{1\le i\le n-1}\atop{1\le m\le L}}\atop{{{\beta\in\N^L}\atop{\sigma\in\N^{n-1}}}\atop{|\beta|+|\sigma|\le|\alpha|-1}}}, \Big\{\mu_{k,\beta,\sigma}\Big\}_{{{\beta\in\N^L}\atop{\sigma\in\N^{n-1}}}\atop{|\beta|+|\sigma|\le|\alpha|-1}}, \Big\{ \vartheta_{m,\delta}\Big\}_{{1\le m\le L}\atop{{\delta\in\N^n}\atop{{|\delta|\le |\alpha|-1}}}}
\right),\\ &{\mbox{where~${\mathcal{P}}_{k,\alpha}$ is a polynomial with {\em
nonnegative coefficients}.}}\end{split}
\end{equation}
Note that this statement is true when~$\alpha_n=0$, thanks to~\eqref{Ojnd-1-Sdfgo0o2ruiytghb-0oritgkfkS},
hence we can focus on the proof of~\eqref{l0qowjfVq92LvezPOLJSdf-1p2orj3tg-304-03} when~$\alpha_n\ne0$.

For this we argue as follows.
Plugging~\eqref{l0qowjfVq92LvezPOLJSdf-1p2orj3tg-304-01} and~\eqref{l0qowjfVq92LvezPOLJSdf-1p2orj3tg-304-02} into~\eqref{9ojknTRvbdfLON0oejgDAnmdcVEujdfrtitghbsdfAJSdfwQUW-2} we find that
\begin{eqnarray*}
&&
\sum_{\alpha\in\N^n} \alpha_n\vartheta_{k,\alpha} Y^{\alpha-e_n}
\\&=&
\partial_t w_k\\&=&\sum_{{1\le i\le n-1}\atop{1\le m\le L}} B_{i,m,k}(y,w)\partial_i w_m+C_k(y,w)
\\&=&
\sum_{{{{1\le i\le n-1}\atop{1\le m\le L}}\atop{{\beta\in\N^L}\atop{\sigma\in\N^{n-1}}}}\atop{{\alpha\in\N^n}}} \lambda_{i,m,k,\beta,\sigma}
\, \alpha_i\vartheta_{m,\alpha} Y^{\alpha+(\sigma,0)-e_i}\,w^\beta
+\sum_{{\beta\in\N^L}\atop{\sigma\in\N^{n-1}}} \mu_{k,\beta,\sigma} \,y^\sigma w^\beta.
\end{eqnarray*}
Hence, since
\begin{eqnarray*}w^\beta&=&w_1^{\beta_1}\dots w_L^{\beta_L}\\&
=&\left(\sum_{\alpha^{(1)}\in\N^n} \vartheta_{1,\alpha^{(1)}}\, Y^{\alpha^{(1)}}\right)^{\beta_1}
\dots\left(\sum_{\alpha^{(L)}\in\N^n} \vartheta_{L,\alpha^{(L)}}\, Y^{\alpha^{(L)}}\right)^{\beta_L}\\
&=&\sum_{{{\alpha^{(1,1)},\dots,\alpha^{(1,\beta_1)}\in\N^n}\atop{\cdots}}\atop{\alpha^{(L,1)},\dots,\alpha^{(L,\beta_L)}\in\N^n}} \vartheta_{1,\alpha^{(1,1)}}\dots\vartheta_{1,\alpha^{(1,\beta_1)}}
\dots\vartheta_{L,\alpha^{(L,1)}}\dots\vartheta_{L,\alpha^{(1,\beta_L)}}\,
Y^{\alpha^{(1,1)}+\dots+\alpha^{(1,\beta_1)}+\dots+\alpha^{(L,1)}+\dots+\alpha^{(L,\beta_L)}}
,\end{eqnarray*}
we obtain that
\begin{eqnarray*}
&&
\sum_{\alpha\in\N^n} \alpha_n\vartheta_{k,\alpha} Y^{\alpha-e_n}
\\&=&
\sum_{{{{{1\le i\le n-1}\atop{1\le m\le L}}\atop{{\beta\in\N^L}\atop{\sigma\in\N^{n-1}}}}\atop{{\alpha\in\N^n}}}\atop{{{\alpha^{(1,1)},\dots,\alpha^{(1,\beta_1)}\in\N^n}\atop{\cdots}}\atop{\alpha^{(L,1)},\dots,\alpha^{(L,\beta_L)}\in\N^n}}}
\alpha_i\lambda_{i,m,k,\beta,\sigma} \,\vartheta_{1,\alpha^{(1,1)}}\dots\vartheta_{1,\alpha^{(1,\beta_1)}}
\dots\vartheta_{L,\alpha^{(L,1)}}\dots\vartheta_{L,\alpha^{(1,\beta_L)}}\vartheta_{m,\alpha}\\&&\qquad\qquad\qquad\qquad\qquad\times
Y^{\alpha^{(1,1)}+\dots+\alpha^{(1,\beta_1)}+\dots+\alpha^{(L,1)}+\dots+\alpha^{(L,\beta_L)}+\alpha+(\sigma,0)-e_i}\\&&\qquad
+\sum_{{{\beta\in\N^L}\atop{\sigma\in\N^{n-1}}}\atop{{{\alpha^{(1,1)},\dots,\alpha^{(1,\beta_1)}\in\N^n}\atop{\cdots}}\atop{\alpha^{(L,1)},\dots,\alpha^{(L,\beta_L)}\in\N^n}}} \mu_{k,\beta,\sigma} \,\vartheta_{1,\alpha^{(1,1)}}\dots\vartheta_{1,\alpha^{(1,\beta_1)}}
\dots\vartheta_{L,\alpha^{(L,1)}}\dots\vartheta_{L,\alpha^{(1,\beta_L)}}\\&&\qquad\qquad\qquad\qquad\qquad\times Y^{\alpha^{(1,1)}+\dots+\alpha^{(1,\beta_1)}+\dots+\alpha^{(L,1)}+\dots+\alpha^{(L,\beta_L)}+(\sigma,0)}.\end{eqnarray*}
As a result,
\begin{equation}\label{Ojnd-1-Sdfgo0o2ruiytghb-0oritgkfkS-XY}\begin{split}
&\alpha_n\vartheta_{k,\alpha} \\ =\,&
\sum_{{{{{{1\le i\le n-1}\atop{1\le m\le L}}\atop{{\beta\in\N^L}\atop{\sigma\in\N^{n-1}}}}\atop{{\gamma\in\N^n}}}\atop{{{\alpha^{(1,1)},\dots,\alpha^{(1,\beta_1)}\in\N^n}\atop{\cdots}}\atop{\alpha^{(L,1)},\dots,\alpha^{(L,\beta_L)}\in\N^n}}}\atop{\alpha^{(1,1)}+\dots+\alpha^{(1,\beta_1)}+\dots+\alpha^{(L,1)}+\dots+\alpha^{(L,\beta_L)}+\gamma-e_i+(\sigma,0)=\alpha-e_n}}
\gamma_i\lambda_{i,m,k,\beta,\sigma} \,\vartheta_{1,\alpha^{(1,1)}}\dots\vartheta_{1,\alpha^{(1,\beta_1)}}
\dots\vartheta_{L,\alpha^{(L,1)}}\dots\vartheta_{L,\alpha^{(1,\beta_L)}}\vartheta_{m,\gamma}\\&\qquad
+\sum_{{{{\beta\in\N^L}\atop{\sigma\in\N^{n-1}}}\atop{{{\alpha^{(1,1)},\dots,\alpha^{(1,\beta_1)}\in\N^n}\atop{\cdots}}\atop{\alpha^{(L,1)},\dots,\alpha^{(L,\beta_L)}\in\N^n}}}\atop{{\alpha^{(1,1)}+\dots+\alpha^{(1,\beta_1)}+\dots+\alpha^{(L,1)}+\dots+\alpha^{(L,\beta_L)}}+(\sigma,0)=\alpha-e_n}} \mu_{k,\beta,\sigma} \,\vartheta_{1,\alpha^{(1,1)}}\dots\vartheta_{1,\alpha^{(1,\beta_1)}}
\dots\vartheta_{L,\alpha^{(L,1)}}\dots\vartheta_{L,\alpha^{(1,\beta_L)}}.\end{split}\end{equation}
This determines~$\vartheta_{k,\alpha}$, but to complete the proof of~\eqref{l0qowjfVq92LvezPOLJSdf-1p2orj3tg-304-03} we need to check the order of the indices in the right hand side of~\eqref{l0qowjfVq92LvezPOLJSdf-1p2orj3tg-304-03}.

To this end, we point out that the conditions~$\gamma_i\ne0$ and
$$ \alpha^{(1,1)}+\dots+\alpha^{(1,\beta_1)}+\dots+\alpha^{(L,1)}+\dots+\alpha^{(L,\beta_L)}+\gamma-e_i+(\sigma,0)=\alpha-e_n$$
give that~$|\gamma|\ge1$ and
$$ |\alpha^{(1,1)}|+\dots+|\alpha^{(1,\beta_1)}|+\dots+|\alpha^{(L,1)}|+\dots+|\alpha^{(L,\beta_L)}|+|\gamma|+|\sigma|=|\alpha|,$$
whence
$$|\alpha^{(1,1)}|+\dots+|\alpha^{(1,\beta_1)}|+\dots+|\alpha^{(L,1)}|+\dots+|\alpha^{(L,\beta_L)}|+|\sigma|
\le|\alpha|-1.$$
Furthermore, in~\eqref{Ojnd-1-Sdfgo0o2ruiytghb-0oritgkfkS-XY} we have that~$|\gamma|\le|\alpha|-1$, due to~\eqref{Ojnd-1-Sdfgo0o2ruiytghb-0oritgkfkS}.

These observations take care of the order of the indices of~$\lambda_{i,m,k,\beta,\sigma}$ in the right hand side of~\eqref{l0qowjfVq92LvezPOLJSdf-1p2orj3tg-304-03}, as well as of the indices of~$\vartheta_{m,\delta}$
in the right hand side of~\eqref{l0qowjfVq92LvezPOLJSdf-1p2orj3tg-304-03} which come from the first sum
in the right hand side of~\eqref{Ojnd-1-Sdfgo0o2ruiytghb-0oritgkfkS-XY}.

Thus, to finish the proof of~\eqref{l0qowjfVq92LvezPOLJSdf-1p2orj3tg-304-03}, we need to take
care, in the right hand side of~\eqref{l0qowjfVq92LvezPOLJSdf-1p2orj3tg-304-03}, of
the order of the indices of~$\mu_{k,\beta,\sigma}$ and that of the indices of~$\vartheta_{m,\delta}$
which come from the second sum
in the right hand side of~\eqref{Ojnd-1-Sdfgo0o2ruiytghb-0oritgkfkS-XY}.

With that in mind, we observe that the condition
\begin{equation}\label{id90438bv0kjhgfhygtfredjhgfdsujhgfd765432534r5t}
\alpha^{(1,1)}+\dots+\alpha^{(1,\beta_1)}+\dots+\alpha^{(L,1)}+\dots+\alpha^{(L,\beta_L)}+(\sigma,0)=\alpha-e_n
\end{equation}
combined with~\eqref{Ojnd-1-Sdfgo0o2ruiytghb-0oritgkfkS} leads to the fact that, in
the second sum
in the right hand side of~\eqref{Ojnd-1-Sdfgo0o2ruiytghb-0oritgkfkS-XY} one has that
$$|\beta|+|\sigma|= \beta_1+\dots+\beta_L+|\sigma|
\le|\alpha^{(1,1)}|+\dots+|\alpha^{(1,\beta_1)}|+\dots+|\alpha^{(L,1)}|+\dots+|\alpha^{(L,\beta_L)}|+|\sigma|=|\alpha|-1,$$
which takes care of the indices of~$\mu_{k,\beta,\sigma}$.

Moreover, from~\eqref{id90438bv0kjhgfhygtfredjhgfdsujhgfd765432534r5t} we also deduce that
$$ |\alpha^{(1,1)}|,\dots,|\alpha^{(1,\beta_1)}|,\dots,|\alpha^{(L,1)}|,\dots,|\alpha^{(L,\beta_L)}|\le|\alpha|-1,$$
which takes care of the indices of~$\vartheta_{m,\delta}$.
This completes the proof of~\eqref{l0qowjfVq92LvezPOLJSdf-1p2orj3tg-304-03}.

We stress that~\eqref{l0qowjfVq92LvezPOLJSdf-1p2orj3tg-304-03} determines all the coefficients~$\vartheta_{k,\alpha}$ uniquely (i.e., by arguing recursively, starting from~\eqref{Ojnd-1-Sdfgo0o2ruiytghb-0oritgkfkS}).
This entails that, if a real analytic solution exists, it is uniquely determined by~\eqref{l0qowjfVq92LvezPOLJSdf-1p2orj3tg-304-01}.\medskip

Having established the uniqueness result in Theorem~\ref{CKTEO}, we
now aim at proving the existence result. To this end, we show that the formal series in~\eqref{l0qowjfVq92LvezPOLJSdf-1p2orj3tg-304-01} is indeed convergent. This task will be accomplished
via the method of majorants. For this, we first use~\eqref{Ojnd-1-Sdfgo0o2ruiytghb-0oritgkfkS}
and~\eqref{l0qowjfVq92LvezPOLJSdf-1p2orj3tg-304-03} to see that\footnote{The counterpart of~\eqref{ILPRYHNSOIJSTGISAmeBYHEEL210} in the setting of ordinary differential equations
was considered in~\eqref{ILPRYHNSOIJSTGISAmeBYHEEL210-ODE}.}
\begin{equation}\label{ILPRYHNSOIJSTGISAmeBYHEEL210} \vartheta_{k,\alpha}={\mathcal{Q}}_{k,\alpha}\left(
\Big\{\lambda_{i,m,k,\beta,\sigma}\Big\}_{{{1\le i\le n-1}\atop{1\le m\le L}}\atop{{{\beta\in\N^L}\atop{\sigma\in\N^{n-1}}}\atop{|\beta|+|\sigma|\le|\alpha|-1}}}, \Big\{\mu_{k,\beta,\sigma}\Big\}_{{{\beta\in\N^L}\atop{\sigma\in\N^{n-1}}}\atop{|\beta|+|\sigma|\le|\alpha|-1}}
\right),\end{equation}
where~${\mathcal{Q}}_{k,\alpha}$ is a polynomial with {\em nonnegative coefficients}.

As a consequence,
\begin{equation}\label{ILPRYHNSOIJSTGISAmeBYHEEL211} |\vartheta_{k,\alpha}|\le{\mathcal{Q}}_{k,\alpha}\left(
\Big\{|\lambda_{i,m,k,\beta,\sigma}|\Big\}_{{{1\le i\le n-1}\atop{1\le m\le L}}\atop{{{\beta\in\N^L}\atop{\sigma\in\N^{n-1}}}\atop{|\beta|+|\sigma|\le|\alpha|-1}}}, \Big\{|\mu_{k,\beta,\sigma}|\Big\}_{{{\beta\in\N^L}\atop{\sigma\in\N^{n-1}}}\atop{|\beta|+|\sigma|\le|\alpha|-1}}
\right).\end{equation}
We stress that in this inequality we have employed the fact that the above coefficients are nonnegative:
this is an interesting case in which one can use
a monotonicity property encoded in an abstract
mathematical structure to deduce useful information hidden behind\footnote{There are however alternative proofs
of the Cauchy-Kowalevsky Theorem that do not rely on the
method of majorants but focus instead on
a direct estimation of all the derivatives of the solution,
see e.g.~\cite[pages~26--35]{MR1306729}.}
incendiary combinatorial calculations.

Furthermore, by the analyticity of~$B_{i,m,k}$ and~$C_k$ as
written in~\eqref{l0qowjfVq92LvezPOLJSdf-1p2orj3tg-304-02},
we know (see e.g.~\cite[Corollary~1.1.7]{MR1916029}) that
\begin{equation}\label{ILPRYHNSOIJSTGISAmeBYHEEL212p}|\lambda_{i,m,k,\beta,\sigma}|+|\mu_{k,\beta,\sigma}|\le \frac{C}{R^{|\beta|+|\sigma|}},\end{equation}
for some~$C$, $R>0$.

As a byproduct of~\eqref{ILPRYHNSOIJSTGISAmeBYHEEL212p}
and~\eqref{ILPRYHNSOIJSTGISAmeBYHEEL212q} we obtain that
\begin{equation}\label{ILPRYHNSOIJSTGISAmeBYHEEL212}|\lambda_{i,m,k,\beta,\sigma}|+|\mu_{k,\beta,\sigma}|\le \frac{C\, (|\beta|+|\sigma|)!}{\beta!\,\sigma!\,R^{|\beta|+|\sigma|}}.\end{equation}
\medskip

We have thus reduced ourselves to the ``worst possible scenario''. Indeed, on the one hand,
if we define
\begin{equation}\label{ILPRYHNSOIJSTGISAmeBYHEEL219}\begin{split}
&\lambda_{i,m,k,\beta,\sigma}^\star:= \frac{C\, (|\beta|+|\sigma|)!}{\beta!\,\sigma!\,R^{|\beta|+|\sigma|}},\\
&\mu_{k,\beta,\sigma}^\star:= \frac{C\, (|\beta|+|\sigma|)!}{\beta!\,\sigma!\,R^{|\beta|+|\sigma|}},\\
& B_{i,m,k}^\star(y,w)=\sum_{{\beta\in\N^L}\atop{\sigma\in\N^{n-1}}} \lambda_{i,m,k,\beta,\sigma}^\star \,y^\sigma w^\beta\\
{\mbox{and}}\qquad& C_k^\star(y,w):=\sum_{{\beta\in\N^L}\atop{\sigma\in\N^{n-1}}} \mu_{k,\beta,\sigma}^\star \,y^\sigma w^\beta
,\end{split}\end{equation}
we infer from~\eqref{ILPRYHNSOIJSTGISAmeBYHEEL211} and~\eqref{ILPRYHNSOIJSTGISAmeBYHEEL212} that
\begin{equation}\label{ILPRYHNSOIJSTGISAmeBYHEEL213} |\vartheta_{k,\alpha}|\le{\mathcal{Q}}_{k,\alpha}\left(
\Big\{\lambda_{i,m,k,\beta,\sigma}^\star\Big\}_{{{1\le i\le n-1}\atop{1\le m\le L}}\atop{{{\beta\in\N^L}\atop{\sigma\in\N^{n-1}}}\atop{|\beta|+|\sigma|\le|\alpha|-1}}}, \Big\{\mu_{k,\beta,\sigma}^\star\Big\}_{{{\beta\in\N^L}\atop{\sigma\in\N^{n-1}}}\atop{|\beta|+|\sigma|\le|\alpha|-1}}
\right).\end{equation}
On the other hand, in light of~\eqref{ILPRYHNSOIJSTGISAmeBYHEEL210}, we also have that if
there exists a real analytic solution 
\begin{equation} \label{ILPRYHNSOIJSTGISAmeBYHEEL21M6}
w_k^\star(Y)=\sum_{\alpha\in\N^n} \vartheta_{k,\alpha}^\star\, Y^\alpha\end{equation}
of the Cauchy problem
\begin{equation}\label{ILPRYHNSOIJSTGISAmeBYHEEL21M}\begin{dcases}
\partial_t w_k^\star=\sum_{{1\le i\le n-1}\atop{1\le m\le L}} B_{i,m,k}^\star(y,w^\star)\partial_i w_m^\star+C_k^\star(y,w^\star)
& {\mbox{ for all~$|y|<\rho$, $|t|<\rho$, $k\in\{1,\dots,L\}$,}}\\
w_k^\star(y,0)=0& {\mbox{ for all $|y|<\rho$ and $k\in\{1,\dots,L\}$,}}
\end{dcases}
\end{equation}
then necessarily
\begin{equation}\label{ILPRYHNSOIJSTGISAmeBYHEEL214} \vartheta_{k,\alpha}^\star={\mathcal{Q}}_{k,\alpha}\left(
\Big\{\lambda_{i,m,k,\beta,\sigma}^\star\Big\}_{{{1\le i\le n-1}\atop{1\le m\le L}}\atop{{{\beta\in\N^L}\atop{\sigma\in\N^{n-1}}}\atop{|\beta|+|\sigma|\le|\alpha|-1}}}, \Big\{\mu_{k,\beta,\sigma}^\star\Big\}_{{{\beta\in\N^L}\atop{\sigma\in\N^{n-1}}}\atop{|\beta|+|\sigma|\le|\alpha|-1}}
\right).\end{equation}
Therefore, we have found a ``majorization problem'' for our original situation, in the sense that,
by~\eqref{ILPRYHNSOIJSTGISAmeBYHEEL213} and~\eqref{ILPRYHNSOIJSTGISAmeBYHEEL214},
\begin{equation}\label{ILPRYHNSOIJSTGISAmeBYHEEL215} |\vartheta_{k,\alpha}|\leq\vartheta_{k,\alpha}^\star.
\end{equation}
Hence, to show that the formal series in~\eqref{l0qowjfVq92LvezPOLJSdf-1p2orj3tg-304-01} defines
in fact an analytic funtion (and thus complete the proof of Theorem~\ref{CKTEO}) it suffices to establish
suitable bounds on the solution (if any) of our majorization problem~\eqref{ILPRYHNSOIJSTGISAmeBYHEEL21M}. This strategy is advantageous, since the
majorization problem can be explicitly
solved.

Indeed, from~\eqref{ILPRYHNSOIJSTGISAmeBYHEEL219} one deduces that
\begin{eqnarray*}&& B_{i,m,k}^\star(y,w)=C_k^\star(y,w)=\sum_{{\beta\in\N^L}\atop{\sigma\in\N^{n-1}}}
\frac{C\,(|\beta|+|\sigma|)!}{\beta!\,\sigma!\,R^{|\beta|+|\sigma|}}\,y^\sigma w^\beta.\end{eqnarray*}
We also set
$$\zeta:=\frac{y_1+\dots+y_{n-1}+w_1+\dots+w_L}{R}.$$
Using the Multinomial Theorem, we observe that, for every~$\kappa\in\N$,
\begin{eqnarray*} \zeta^{\kappa}&
=&\sum _{{\omega\in\N^{L+n-1}}\atop{|\omega|=\kappa}}
\frac{\kappa!}{\omega_1!\dots\omega_{L+n-1}!\,R^{|\omega|}}\,
y_1^{\omega_1}\dots y_{n-1}^{\omega_{n-1}} w_1^{\omega_n} \dots w_L^{\omega_{L+n-1}}\\&=&
\sum _{{\omega=(\beta,\sigma)\in\N^{L}\times\N^{n-1}}\atop{|\beta|+|\sigma|=\kappa}}
\frac{(|\beta|+|\sigma|)!}{\beta!\sigma!\,R^{|\beta|+|\sigma|}}\,
y_1^{\sigma_1}\dots y_{n-1}^{\sigma_{n-1}} w_1^{\beta_1} \dots w_L^{\beta_{L}}
.\end{eqnarray*}
Therefore\footnote{The counterpart of~\eqref{cc-19425jl0qowjfVq92LvezPOLJSdf-1p2orj3tg-304-01} in the setting of ordinary differential equations
was considered in~\eqref{cc-19425jl0qowjfVq92LvezPOLJSdf-1p2orj3tg-304-01-ODE}.} if~$|\zeta|<1$ we find that
\begin{equation} \label{cc-19425jl0qowjfVq92LvezPOLJSdf-1p2orj3tg-304-01}
B_{i,m,k}^\star(y,w)=C_k^\star(y,w)=C\sum_{\kappa\in\N}\zeta^\kappa=\frac{C}{1-\zeta}=\frac{CR}{R-(y_1+\dots+y_{n-1}+w_1+\dots+w_L)}.\end{equation}

Consequently, the majorization problem~\eqref{ILPRYHNSOIJSTGISAmeBYHEEL21M} can be written in the form
\begin{equation}\label{098765r4POJL-234}\begin{dcases}
\frac{R-(y_1+\dots+y_{n-1}+w_1+\dots+w_L)}{CR}\partial_t w_k^\star=\sum_{{1\le i\le n-1}\atop{1\le m\le L}}\partial_i w_m^\star+1
& {\mbox{ for all~$|y|<\rho$, $|t|<\rho$, $k\in\{1,\dots,L\}$,}}\\
w_k^\star(y,0)=0& {\mbox{ for all $|y|<\rho$ and $k\in\{1,\dots,L\}$.}}
\end{dcases}
\end{equation}
We seek a solution for this problem in the form
\begin{equation}\label{098765r4POJL-234bis} w_1^\star(y,t)=\dots=w_L^\star(y,t):=W(\eta,t),\qquad{\mbox{where }}\eta:=y_1+\dots+y_{n-1}.\end{equation}
With this ansatz, up to renaming~$\rho$, we reduce~\eqref{098765r4POJL-234} to
\begin{equation}\label{ICSWNCOMKK893TP}\begin{dcases}
\frac{R-(\eta+LW)}{CR}\partial_t W=(n-1)L\partial_\eta W+1
& {\mbox{ for all~$|\eta|<\rho$, $|t|<\rho$,}}\\
W(\eta,0)=0& {\mbox{ for all $|\eta|<\rho$.}}
\end{dcases}
\end{equation}
To find a solution of this, we use the so-called method of characteristics, i.e. \index{method of characteristics}
we look for a curve~$\tau\mapsto(\eta(\tau),t(\tau))$ along which we control the solution.
Since
$$ \frac{d}{d\tau}\big(W(\eta(\tau),t(\tau))\big)=\partial_\eta W(\eta(\tau),t(\tau))\dot\eta(\tau)+\partial_t W(\eta(\tau),t(\tau))\dot{t}(\tau),$$
in light of~\eqref{ICSWNCOMKK893TP} it is convenient to choose
\begin{equation}\label{IJMINTHBAQG} \dot\eta:=-(n-1)L\qquad{\mbox{and}}\qquad\dot{t}=\frac{R-(\eta+LW)}{CR},\end{equation}
so that
\begin{eqnarray*}&&\frac{d}{d\tau}\big(W(\eta(\tau),t(\tau))\big)=1
\end{eqnarray*}
and accordingly
\begin{equation}\label{IJMINTHBAQG2}W(\eta(\tau),t(\tau))=\tau.\end{equation}
Given~$\eta_0\in\R$, from~\eqref{IJMINTHBAQG} and~\eqref{IJMINTHBAQG2} we arrive at
\begin{equation}\label{IJMINTHBAQG-01}\eta(\tau)=\eta_0-(n-1)L\tau\end{equation}
and
$$ \dot{t}(\tau)=\frac{R-\eta_0+(n-1)L\tau-L\tau}{CR}=\frac{R-\eta_0+(n-2)L\tau}{CR},$$
therefore
\begin{equation}\label{IJMINTHBAQG-02} t(\tau)=\frac{R-\eta_0}{CR}\tau+\frac{(n-2)L}{2CR}\tau^2.\end{equation}
Our goal is now to invert these parametric representations, i.e., rather than having~$\eta$ and~$t$ as functions of~$\tau$ (and~$\eta_0$),
we wish to write~$\tau$ in dependence of~$\eta$ and~$t$. For this, owing to~\eqref{IJMINTHBAQG-01} and~\eqref{IJMINTHBAQG-02}, we note that
$$  t=\frac{R-(\eta+(n-1)L\tau)}{CR}\tau+\frac{(n-2)L}{2CR}\tau^2.$$
Hence, solving for~$\tau$ (and picking the solution for which~$\tau=0$ when~$t=0$), we find
$$ \tau=\frac{ R-\eta-\sqrt{ (R- \eta)^2-2Cn L R t} }{n L}.$$
All in all,
\begin{equation}\label{l0qowjfVq9020urVHUICVGHYUoo-01}
W(\eta,t)=\tau=\frac{ R-\eta-\sqrt{ (R- \eta)^2-2Cn L R t} }{n L},
\end{equation}
which provides\footnote{The counterpart of~\eqref{l0qowjfVq9020urVHUICVGHYUoo-01} in the setting of ordinary differential equations
was considered in~\eqref{l0qowjfVq9020urVHUICVGHYUoo-01-ODE}.} a real analytic solution for the Cauchy problem in~\eqref{ICSWNCOMKK893TP} as long as~$\rho$ is small enough.

Recalling~\eqref{098765r4POJL-234bis}, we have thereby constructed explicitly a real analytic solution of~\eqref{098765r4POJL-234}. {F}rom its analyticity (see e.g.~\cite[Corollary~1.1.7]{MR1916029}), and recalling the notation in~\eqref{ILPRYHNSOIJSTGISAmeBYHEEL21M6}, we infer that
$$ \vartheta_{k,\alpha}^\star\le|\vartheta_{k,\alpha}^\star|\le\frac{C^\star}{\rho^{|\alpha|}},$$
for some~$C^\star$, $\rho>0$, and consequently, by~\eqref{ILPRYHNSOIJSTGISAmeBYHEEL215},
$$ |\vartheta_{k,\alpha}|\le\frac{C^\star}{\rho^{|\alpha|}}.$$
This shows that the formal series in~\eqref{l0qowjfVq92LvezPOLJSdf-1p2orj3tg-304-01} is
convergent near the origin,
thus providing locally a real analytic solution of the Cauchy problem in~\eqref{9ojknTRvbdfLON0oejgDAnmdcVEujdfrtitghbsdfAJSdfwQUW-2} and finishing the proof of Theorem~\ref{CKTEO}.
\end{proof}

For further readings about the Cauchy-Kowalevsky Theorem, see e.g.~\cite[Chapter~6]{MR678605},
\cite[Chapter~1, Section~7]{MR1013360}, \cite[Chapter~2]{MR1141630}, \cite[Chapter~3, Section~3]{MR1185075},
\cite[Sections~2.3, 2.4 and~2.5]{MR1916029},
\cite[Chapter~2, Sections~17, 18, 19, 20 and~21]{MR2301309},
\cite[Sections~II, III and IV]{MR3409063} and~\cite[Chapter~5]{MR4436039}.

For additional information about the method of majorants see~\cite{MR1452105}
and the references therein.

See also~\cite{MR15186, MR76161, MR0091421, MR93634, MR104912, MR115010, MR166462, MR0285941, MR322321,  MR512931, MR696869, MR2301259} and the references therein for different proofs and generalizations
of the Cauchy-Kowalevsky Theorem.

\begin{appendix}

\chapter{An interesting example}\label{MAHH}

\section{The importance of the interior ball condition for the Hopf Lemma}
In this appendix, we showcase
that the interior ball condition in Lemma~\ref{JJS:PA} cannot be dropped
by recalling the following interesting example in which all the hypotheses of Lemma~\ref{JJS:PA}
are met with the exception of the interior ball condition and the thesis of Lemma~\ref{JJS:PA} is violated.

Given~$(x,y)\in\R^2$, we use the notation~$z=x+iy\in \cOMPL$.
We observe that a branch of the complex logarithmic function~$\operatorname{Log}$
in the halfplane~$[0,+\infty)\times\R$ can be uniquely defined: for instance, using polar coordinates~$z=\rho e^{i\vartheta}$ (corresponding to~$x=\rho\cos\vartheta$ and~$y=\rho\sin\vartheta$), when~$\vartheta\in\left[-\frac\pi2,\frac\pi2\right]$ we can write that~$\operatorname{Log} z:=\ln \rho+i\vartheta$. We also set
$$ {\mathcal{U}}:=\big\{
z=x+iy\in\cOMPL {\mbox{ s.t. }}x\in (0,1)
\big\}.$$
In this way, since~$\operatorname{Log} z\ne0$ in~${\mathcal{U}}$, we have that the function~${\mathcal{U}}\ni z\mapsto
\frac{z}{\operatorname{Log} z}$ is holomorphic and accordingly, by~\eqref{KS:RUDVA0o1keLLetyh439}, the function
\begin{equation}\label{USCRPOBIS} u(x,y):=\Re\left(\frac{z}{\operatorname{Log} z}\right)\end{equation}
is harmonic when~$x+iy\in{\mathcal{U}}$.

\begin{figure}
                \centering
                \includegraphics[width=.4\linewidth]{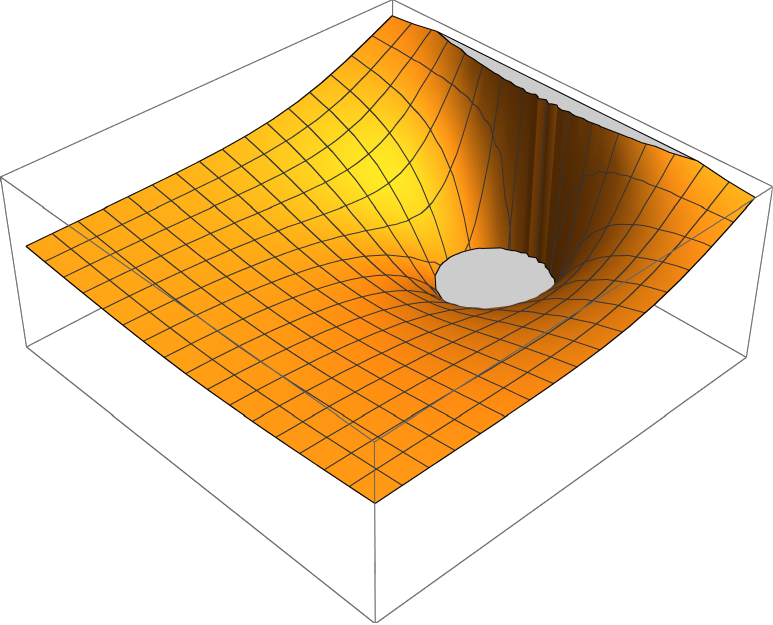}$\quad$\includegraphics[width=.4\linewidth]{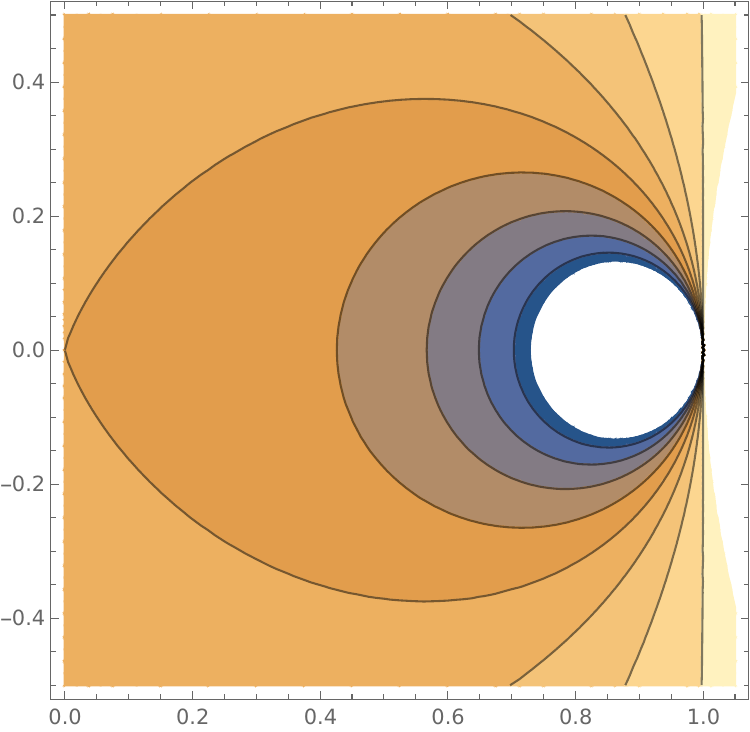}
        \caption{\sl The function~$u$ in~\eqref{USCRPO} and its level sets.}\label{la-funHAFOUMSldRGIRApataNOJHNFOJED}
\end{figure}

We observe that~$u$ can also be written in polar coordinates as
\begin{equation}\label{USCRPO}
u(x,y)=\Re\left(\frac{\rho e^{i\vartheta}}{\ln \rho+i\vartheta}\right)=
\Re\left(\frac{\rho(\cos\vartheta+i\sin\vartheta)(\ln \rho-i\vartheta)}{(\ln \rho)^2+\vartheta^2}\right)=
\frac{\rho\, (\ln \rho\cos\vartheta+\vartheta\sin\vartheta)}{(\ln \rho)^2+\vartheta^2}.
\end{equation}
See Figure~\ref{la-funHAFOUMSldRGIRApataNOJHNFOJED} for a sketch of the graph of~$u$.

Notice also that, for sufficiently small~$\rho$, by~\eqref{USCRPO},
$$ |u(x,y)|\le\frac{2\rho}{ \ln \rho }=o(\rho)$$
and therefore~$u(0,0)=0$ and~$\nabla u(0,0)=(0,0)$.

Furthermore,
we also infer from~\eqref{USCRPOBIS} that, if~$x>0$,
\begin{equation}\label{90:pk8i2hrplieorem} 
\nabla  u(x,y)=\Re\left(\nabla \frac{z}{\operatorname{Log} z}\right)=\Re\left(\frac{(1,i)}{\operatorname{Log} z}
+\frac{z}{\operatorname{Log}^2 z}\nabla \operatorname{Log} z\right)
= \Re\left(\frac{(1,i)}{\operatorname{Log} z}+\frac{(1,i)}{\operatorname{Log}^2 z}
\right).
\end{equation}
Therefore,
if~$x>0$, for small~$\rho$,
\begin{eqnarray*}
|\nabla u(x,y)|=O\left( \frac{1}{|\operatorname{Log} z|}\right)=O\left(\frac{1}{|\ln\rho|}\right)
\end{eqnarray*}
yielding that
\begin{equation}\label{CON0-01-lao}
\lim_{{(x,y)\to(0,0)}\atop{x+iy\in{\mathcal{U}}}}\nabla u(x,y)=0=\nabla u(0,0).
\end{equation}

We now define
\begin{equation}\label{POER-l3e-1pata0}
\Omega_o:=\big\{(x,y)\in\R^2 {\mbox{ s.t. }}x+iy\in{\mathcal{U}}
{\mbox{ and }} u(x,y)<0
\big\}.\end{equation}
Then, we have that
\begin{equation}\label{POER-l3e-1}
{\mbox{$\Omega_o$ is a bounded domain with boundary of class~$C^1$, with~$(0,0)$, $(1,0)\in\partial\Omega$.}}
\end{equation}
Indeed, we can exploit~\eqref{USCRPO} to write~$\Omega_o$ in polar coordinates in the form
\begin{equation}\label{LOSTU}\begin{split}&
\left\{ (\rho,\vartheta)\in(0,+\infty)\times\left(-\frac\pi2,\frac\pi2\right) {\mbox{ s.t. }}
\ln \rho\cos\vartheta+\vartheta\sin\vartheta<0
\right\}\\&\qquad=
\left\{ (\rho,\vartheta)\in(0,+\infty)\times\left(-\frac\pi2,\frac\pi2\right) {\mbox{ s.t. }}
\rho<\frac1{\exp\left(\vartheta\tan\vartheta\right)}
\right\}.\end{split}
\end{equation}
In particular, if~$\rho$ is as above, we have that~$\rho<1$, giving that~$\Omega_o$ is bounded.

Furthermore, since
$$ \lim_{\vartheta\to\pm\pi/2}\vartheta\tan\vartheta=+\infty,$$
we deduce from~\eqref{LOSTU} that
$$ \partial\Omega_o=(0,0)\cup\left\{ (\rho,\vartheta)\in(0,+\infty)\times\left(-\frac\pi2,\frac\pi2\right) {\mbox{ s.t. }}
\rho=\frac1{\exp\left(\vartheta\tan\vartheta\right)}
\right\}.$$
This yields that\footnote{A weaker form of~\eqref{ALS:YHETERNna7TY} can be proved as follows.
{F}rom~\eqref{90:pk8i2hrplieorem}, we have that
that~$\{\nabla u=0\}\cap\{x>0\}=\{
\operatorname{Log}z=-1\}$. {F}rom this and~\eqref{USCRPOBIS} we obtain that~$\{\nabla u=0\}\cap\{x>0\}\subseteq
\{ u=-\Re z\}=\{ u=-x\}$ and thus~$\{\nabla u=0\}\cap\{u=0\}\cap\{x>0\}=\varnothing$.
This and the Implicit Function Theorem yield that
\begin{equation*}
{\mbox{$\partial\Omega_o\setminus (B_\e\cup B_\e(1,0))$ is of class~$C^1$ for each~$\e>0$.}}\end{equation*}
This is weaker than~\eqref{ALS:YHETERNna7TY}, since the point~$(1,0)$ is not
captured by this argument (and it cannot be,
since, as remarked in~\eqref{LKZZpr0p3o4u4z0247ol3}, the point~$(1,0)$ does not belong to~$\{u=0\}$,
though it lies in~$\partial\Omega_o$).

To boot, an additional argument is needed to prove~\eqref{POER-l3e-1},
since~\eqref{CON0-01-lao} prevents us from using the Implicit Function Theorem at the origin.}
\begin{equation}\label{ALS:YHETERNna7TY}
{\mbox{$\partial\Omega_o\setminus B_\e$ is of class~$C^1$ for each~$\e>0$.}}\end{equation} Hence, to complete the proof of~\eqref{POER-l3e-1},
we only need to take into consideration the regularity of~$\partial\Omega_o$ in the vicinity of the origin. To this end,
given a small~$\delta>0$, we consider~$\vartheta\in\left(-\frac\pi2,-\frac\pi2+\delta\right)\cup\left(\frac\pi2-\delta,\frac\pi2\right)$
and we write~$\rho=\frac1{\exp\left(\vartheta\tan\vartheta\right)}$ as
\begin{equation}\label{KMS:CURVESM} x=\frac{\cos\vartheta}{\exp\left(\vartheta\tan\vartheta\right)}
\qquad{\mbox{and}}\qquad y=\frac{\sin\vartheta}{\exp\left(\vartheta\tan\vartheta\right)}.\end{equation}
Let also
\begin{equation}\label{SHNOT7} (-\delta,\delta)\setminus\{0\}\ni \tau\mapsto \varphi(\tau):=\begin{dcases}
\displaystyle\frac1{\exp\left(\left(\tau+\frac\pi2\right)\tan\left(\tau+\frac\pi2\right)\right)}&{\mbox{ if }}\tau\in(-\delta,0),\\
\\
\displaystyle\frac1{\exp\left(\left(\tau-\frac\pi2\right)\tan\left(\tau-\frac\pi2\right)\right)}&{\mbox{ if }}\tau\in(0,\delta).
\end{dcases}\end{equation}
We note that~$\varphi$ can be continuously extended at~$\tau=0$, since
\begin{eqnarray*}&&\lim_{\tau\to0^+}\varphi(\tau)=
\lim_{\tau\to0^+}\frac1{\exp\left(\left(\tau-\frac\pi2\right)\tan\left(\tau-\frac\pi2\right)\right)}=0
\\{\mbox{and}}\qquad&&
\lim_{\tau\to0^-}\varphi(\tau)=
\lim_{\tau\to0^-}\frac1{\exp\left(\left(\tau+\frac\pi2\right)\tan\left(\tau+\frac\pi2\right)\right)}=0,\end{eqnarray*}
so from now on we will consider~$\varphi$ as defined for~$\tau\in(-\delta,\delta)$, with~$\varphi(0)=0$.

\begin{figure}
                \centering
                \includegraphics[width=.45\linewidth]{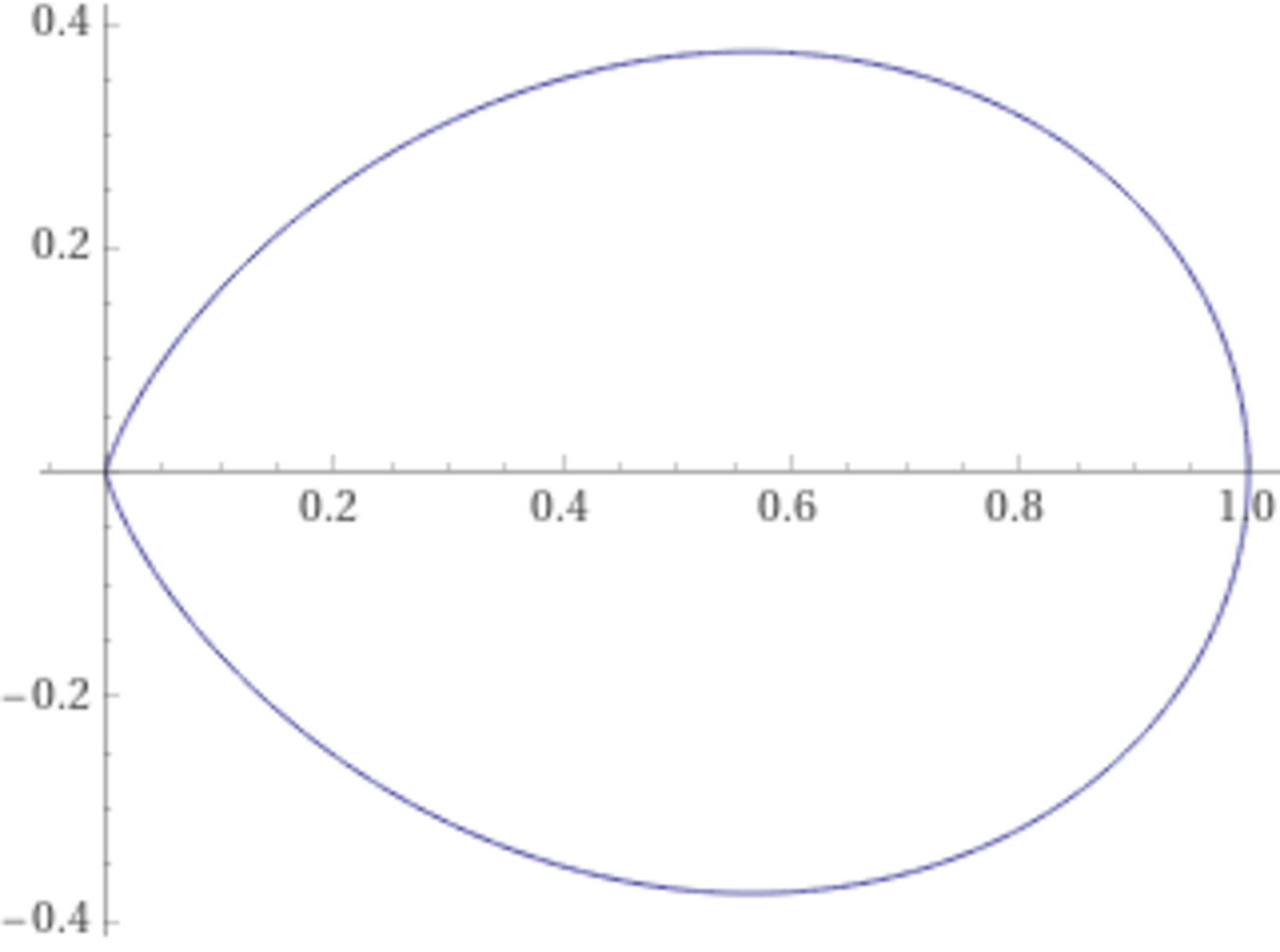}
        \caption{\sl The set~$\Omega_o$ in~\eqref{POER-l3e-1pata0}.}\label{HAFOUMSldRGIRApataNOJHNFOJED}
\end{figure}

We also set
$$ (-\delta,\delta)\ni\tau\mapsto\zeta(\tau):=
\begin{dcases}
\tau+\frac\pi2&{\mbox{ if }}\tau\in(-\delta,0),\\ \\
\tau-\frac\pi2&{\mbox{ if }}\tau\in(0,\delta),
\end{dcases}$$
so that~\eqref{SHNOT7} allows the short notation
$$\varphi(\tau)=\frac1{\exp\left( \zeta(\tau)\tan(\zeta(\tau))\right)},$$ with a continuous extension at~$\tau=0$. Thus, by~\eqref{KMS:CURVESM}, we see that~$\partial\Omega_o$ near the origin can be
identified with the parametric curve
\begin{equation}\label{PAMS-CUYRFVENLKTTA}
\begin{dcases} x(\tau)=\frac{\cos(\zeta(\tau))}{\exp\left((\zeta(\tau))\tan(\zeta(\tau))\right)}=
\varphi(\tau)\,\cos(\zeta(\tau))=\sign\tau \,\sin\tau\, \varphi(\tau)
,\\
y(\tau)=\frac{\sin(\zeta(\tau))}{\exp\left((\zeta(\tau))\tan(\zeta(\tau))\right)}=\varphi(\tau)\,\sin(\zeta(\tau))=
-\sign \tau\,\cos\tau\,\varphi(\tau)
,\end{dcases}\end{equation}
where the sign function has been used.

We also recall that, by L'H\^opital's Rule,
\begin{equation}\label{LHOPCSC1v}
\lim_{\tau\to0}\tau\tan\left(\tau-\frac{\sign\tau\;\pi}2\right)=
\lim_{\tau\to0}\frac{\tan\left(\tau-\frac{\sign\tau\;\pi}2\right)}{1/\tau}=
\lim_{\tau\to0}\frac{1/\sin^2\tau}{-1/\tau^2}=-1.
\end{equation}

We now denote by~$\Xi$ a function that is smooth outside the origin and such that, for every~$k\in\N$,
$$ \limsup_{\tau\to0^+} |D^k\Xi(\tau)|+\limsup_{\tau\to0^-} |D^k\Xi(\tau)|<+\infty.$$
For typographical convenience, such a function~$\Xi$
can vary from line to line of the forthcoming computation.
We claim that, for small~$\tau\ne0$, the function~$\frac{\tau\tan\left(\tau-\frac{\sign\tau\;\pi}2\right)+1}{\tau^2}$
satisfies the above properties, hence we can write that
\begin{equation}\label{MAGBVYHLuevR9-2}
\frac{\tau\tan\left(\tau-\frac{\sign\tau\;\pi}2\right)+1}{\tau^2}=\Xi(\tau).
\end{equation}
To check this, we employ~\eqref{LHOPCSC1v} and standard Taylor expansions for trigonometric functions, thus obtaining that, for~$\tau\ne0$,
\begin{eqnarray*}&&
\tau\tan\left(\tau-\frac{\sign\tau\;\pi}2\right)+1=
\int_0^\tau \frac{d}{d\ell}\left[ \ell\tan\left(\ell-\frac{\sign\ell\;\pi}2\right)\right]\,d\ell\\&&\qquad=
\int_0^\tau \left[ \frac\ell{\sin^2\ell}-\frac{\cos\ell}{\sin\ell}\right]\,d\ell=
\int_0^\tau \left[ \frac{\ell-\sin\ell\cos\ell}{\sin^2\ell}\right]\,d\ell\\&&\qquad=
\int_0^\tau \left[ \frac{\ell-(\ell+\ell^3\;\Xi(\ell))(1+\ell^2\;\Xi(\ell))}{\sin^2\ell}\right]\,d\ell=
\int_0^\tau \frac{\ell^3\;\Xi(\ell)}{\sin^2\ell}\,d\ell\\&&\qquad=\int_0^\tau \ell \;\Xi(\ell)\,d\ell=
\tau^2\int_0^1 \lambda \;\Xi(\tau\lambda)\,d\lambda=\tau^2\;\Xi(\tau),
\end{eqnarray*}
where the possibility of renaming~$\Xi$ at each step of the calculation was repeatedly utilized.
This proves~\eqref{MAGBVYHLuevR9-2}.

Now, by~\eqref{MAGBVYHLuevR9-2},
\[ \tau\tan\left(\tau-\frac{\sign\tau\;\pi}2\right)=-1 + \tau^2 \;\Xi(\tau)\]
and consequently, by~\eqref{SHNOT7},
\begin{eqnarray*}
\varphi(\tau)&=&\frac1{\exp\left(\left(\tau-\frac{\sign\tau\;\pi}2\right)\left(-\frac1\tau +\tau \;\Xi(\tau)\right)\right)}\\
&=&\frac1{\exp\left(
\frac{\sign\tau\;\pi}{2\tau}-1+ \tau \;\Xi(\tau)\right)}
\\&=&\frac1{\exp\left(
\frac{\pi}{2|\tau|}-1+ \tau \;\Xi(\tau)\right)}.
\end{eqnarray*}
It is thereby convenient to reparameterize this curve by setting~$t:=\frac{\sign\tau}{\exp\left(
\frac{\pi}{2|\tau|}\right)}$. In this way, we have that~$\sign t=\sign \tau$ and
$$ \tau=-\frac{\sign t\;\pi}{2 \ln|t|}.$$
Hence, with a slight abuse of notation, for small~$t$ we can write
\begin{eqnarray*}&&
\sign\tau\,\varphi(\tau)=
\frac{\sign\tau}{\exp\left(\frac{\pi}{2|\tau|}\right)
\exp\left(-1+ \tau \;\Xi(\tau)\right)}=
\frac{t}{\exp\left(-1+ \frac1{\ln|t|}\;\Xi\left(\frac1{\ln|t|}\right)\right)}\\&&\qquad=
\frac{et}{1+ \frac1{\ln|t|}\;\Xi\left(\frac1{\ln|t|}\right)}=et \left(1+ \frac1{\ln|t|}\;\Xi\left(\frac1{\ln|t|}\right)\right).
\end{eqnarray*}
Moreover,
\begin{eqnarray*}
\sin\tau=\tau+ \tau^3\;\Xi(\tau)=
-\frac{\sign t\;\pi}{2 \ln|t|}+
\frac{1}{\ln^3|t|}\;\Xi
\left(\frac{1}{\ln|t|}\right)
\end{eqnarray*}
and
\begin{eqnarray*}
\cos\tau=1+\tau^2\;\Xi(\tau)=
1+\frac{1}{\ln^2|t|}\;\Xi\left(\frac{1}{\ln|t|}\right).
\end{eqnarray*}

These observations allow us to recast the parametric curve in~\eqref{PAMS-CUYRFVENLKTTA} into
\begin{equation}\label{9omnSPLvoeDLL}
\begin{dcases} x(t)=et \left(1+
\frac1{\ln|t|}\Xi\left(\frac1{\ln|t|}\right)\right)\left(-\frac{\sign t\;\pi}{2 \ln|t|}+\frac{1}{\ln^3|t|}\;\Xi\left(\frac{1}{\ln|t|}\right)\right)=
et\left(-\frac{\sign t\;\pi}{2 \ln|t|}+\frac{1}{\ln^2|t|}\;\Xi\left(\frac{1}{\ln|t|}\right)\right)
,\\
y(t)=et \left(1+ \frac1{\ln|t|}\;\Xi\left(\frac1{\ln|t|}\right)\right)\left(1+
\frac1{\ln^2|t|}\;\Xi\left(\frac{1}{\ln|t|}\right)\right)
=et \left(1+ \frac1{\ln|t|}\;\Xi\left(\frac1{\ln|t|}\right)\right).
\end{dcases}\end{equation}
{F}rom the latter equation, we deduce that
$$ \frac1e\lim_{t\to0} \dot y(t)=
\lim_{t\to0}\frac{d}{dt}\left[t\left(1+ \frac1{\ln|t|}\;\Xi\left(\frac1{\ln|t|}\right)\right)\right]
=1$$
and, as a result, by the Inverse Function Theorem, for small~$y$,
$$ t=\frac{y}{e}\left(1+\frac{1}{\ln|y|}\;\Xi\left(\frac{1}{\ln|y|}\right)\right).$$
Plugging this information into the first equation in~\eqref{9omnSPLvoeDLL},
we obtain that~$\partial\Omega_o$, in the vicinity of the origin, can be written as
a graph in the horizontal direction in the form
\begin{equation}\label{POER-l3e-10-090vo}
x=y\left(-\frac{\sign y\;\pi}{2 \ln|y|}+\frac{1}{\ln^2|y|}\;\Xi\left(\frac{1}{\ln|y|}\right)\right)
=-\frac{\pi |y|}{2 \ln|y|}+\frac{1}{\ln^2|y|}\;\Xi\left(\frac{|y|}{\ln|y|}\right)
=:F(y).
\end{equation}
Since~$F$ is of class~$C^1$ near~$y=0$, this proves~\eqref{POER-l3e-1}.

Notice also that
\begin{equation}\label{POER-l3e-10}
{\mbox{the interior ball condition for~$\Omega_o$ is not satisfied at the origin,}}\end{equation}
because, by~\eqref{POER-l3e-10-090vo},
$$ \lim_{y\to0}|F''(y)|=+\infty.$$
See Figure~\ref{HAFOUMSldRGIRApataNOJHNFOJED} for a sketch of the domain~$\Omega_o$.

Besides, we point out that, by~\eqref{USCRPO},
\begin{equation}\label{LKZZpr0p3o4u4z0247ol3}
\liminf_{\Omega_o\ni (x,y)\to (1,0)} u(x,y)\le
\lim_{x\nearrow1}u(x,0)
=\lim_{\rho\nearrow1}\frac{\rho\,(\ln \rho+0)}{(\ln \rho)^2+0}
=\lim_{\rho\nearrow1}\frac{1}{\ln \rho}=-\infty.
\end{equation}
Therefore, it is convenient to consider a domain~$\Omega\subset\Omega_o$ such that~$\Omega$ is of class~$C^1$ and
$\Omega\supset\Omega_o\cap\{x\le1/2\}$. In this way, we have that~$u\in C^2(\Omega)\cap C^1(\overline\Omega)$, due to~\eqref{CON0-01-lao}, and, by construction,
$$ \overline\Omega\subseteq\overline{\Omega_o}\subseteq\overline{\{u<0\}}.$$
{F}rom this we infer that if~$x\in\partial\Omega$ then~$u(x)\le0=u(0,0)$, whence
$$ u(0,0)=\max_{\partial\Omega}u.$$
These observations yield that all the assumptions of Lemma~\ref{JJS:PA} are satisfied, with the exception of the inner ball condition.
But the thesis of Lemma~\ref{JJS:PA} is violated, since~$\partial_\nu u(0,0)=0$, owing to~\eqref{CON0-01-lao},
whis showcases the fact that the inner ball condition cannot be removed from
Lemma~\ref{JJS:PA}.

\chapter{An application to physical geodesy}\label{GOEAPL}

\section{Determination of the gravitational potential}

In what follows, relying\index{gravitation}
on combined tools from partial differential equations,
functional analysis and statistics,
we will briefly discuss the following
two classical problems of physical geodesy
(see~\cite{MR0478314, MR414029, MR927099}):
\begin{itemize}
\item given some experimental data on the gravity potential measured at
some specific locations, is it possible to have some {\em
knowledge of the gravity potential at other, possibly unaccessible to experiments, sites?}
\item is it possible to extrapolate from these data
some {\em
knowledge of the mass density of the earth in unreachable locations,
such as the center of the earth?}
\end{itemize}
Concrete answers to these two very fascinating, but quite difficult,
questions will be provided respectively
in~\eqref{JS:3456789JHSmegjSieldilefni-PRE}
and~\eqref{dd-023-929-mea-KSM}-\eqref{dd-023-929-mea-KSM-POSTO}.\medskip

Let us now dive deeper into details,
following~\cite{MR0478314, MR414029, MR927099},
for an application of elliptic partial differential equations,
spherical harmonics,  reproducing kernels
and statistics to the determination of the gravity field.
In spite of the scholastic approach to the problem, which gets
dismissed by stating that the gravity acceleration is everywhere constantly
equal to $9.8$ m/s${}^2$, the problem is technically more complicated,
also in view of the lack of homogeneity of the earth.
Namely,
anomalous density distributions within the earth can produce
significant gravity anomalies that need to be taken into account
in engineering and technological applications.
Roughly speaking,
negative gravity anomalies can be produced
by local deficiencies in mass due to rocks which are less dense
than the average (such as
granite plutons, sedimentary basins, etc.)
and, conversely, positive gravity anomalies
can be the outcome of the presence of dense rocks
(such as basalt).
One of the objectives of physical geodesy is to have a better understanding
of the gravity field starting from the millions of experimental measurements
taken with an uneven distribution
(though nowadays
measurements of satellite orbits make everything much more precise)
to obtain a reliable ``educated guess''
at points where measurements are not available.
For this, mathematics seems indeed to play a pivotal\footnote{{\em``It is a lively
world, which can almost certainly benefit from some more mathematics''},
see~\cite[page 170]{MR927099}.
{\em``It seems to me dangerous to give too much attention
to the non-theoretical approach to physical geodesy [...] An exact
theoretical study of the probabilistic background of the method, as I have
tried to carry through here, will often tell you, how the problems have
been changed and hidden in place, where they are difficult to discover''},
see~\cite[pages~92-93]{MR0478314}.} role.
Here,
we try to provide an exposition accessible to a broad audience,
possibly by skating over some technical details, but trying to
highlight some of the brilliant ideas that emerged in the analysis of this fascinating topic
(see also~\cite{PGsfdgj576E} for the basics of physical geodesy).\medskip

The gravity force at a point~$x\in\R^3$
produced by a point mass~$\mu$ located at a point~$x_0$
is proportional to~$-\frac{\,\mu\,(x-x_0)}{|x-x_0|^2} $, which
corresponds to a potential proportional to~$\frac{\mu}{|x-x_0|} $
(being the gravitational force obtained as the gradient of the potential).
Remarkably, such a gravity potential is harmonic outside the given mass source
(recall~\eqref{GAMMAFU} and the fact that~$n=3$).
The gravity force produced by
an extensive mass can be seen as the superposition of the forces
produced by all ``infinitesimal point masses'' that constitute it,
hence this composite force can also be seen as the gradient of
a gravitation potential, which is harmonic away from the mass itself.
Recalling Section~\ref{SLLD:SPJEDHHA}, we will write this potential as a superposition of spherical harmonics,
according to the so-called ``Laplace expansion of the gravitational potential''.
For this, we point out that, if~$N_k$ is the dimension of
the space of harmonic homogeneous polynomials
of degree~$k$,
as given in~\eqref{ALEHD:DJFJFJ}, the physical case~$n=3$ produces~$N_k=2k+1$.
Furthermore, if~$|y|<|x|$, using the notation~$\tau:=\frac{|y|}{|x|}$
and~$t:=\frac{x\cdot y}{|x|\,|y|}$,
the generating function expression\footnote{Concretely, \label{LISTACON}
using for instance~\eqref{LEGGEB},
\eqref{LEGGE1}
or~\eqref{TRVSBlatsds},
one can list the first ten Legendre polynomials when~$n=3$
as follows:
$$ \label{AEDSFesakns3s24sa098}
{\displaystyle {\begin{array}{r|r}k&P_{k}(t)\\\hline 0&1\\1&t\\2&{\tfrac {1}{2}}\left(3t^{2}-1\right)\\3&{\tfrac {1}{2}}\left(5t^{3}-3t\right)\\4&{\tfrac {1}{8}}\left(35t^{4}-30t^{2}+3\right)\\5&{\tfrac {1}{8}}\left(63t^{5}-70t^{3}+15t\right)\\6&{\tfrac {1}{16}}\left(231t^{6}-315t^{4}+105t^{2}-5\right)\\7&{\tfrac {1}{16}}\left(429t^{7}-693t^{5}+315t^{3}-35t\right)\\8&{\tfrac {1}{128}}\left(6435t^{8}-12012t^{6}+6930t^{4}-1260t^{2}+35\right)\\9&{\tfrac {1}{128}}\left(12155t^{9}-25740t^{7}+18018t^{5}-4620t^{3}+315t\right)\\10&{\tfrac {1}{256}}\left(46189t^{10}-109395t^{8}+90090t^{6}-30030t^{4}+3465t^{2}-63\right)\\\hline \end{array}}}
$$
and higher degrees Legendre polynomials can be computed similarly.} of Legendre polynomials in~\eqref{TRVSBlatsds}
when~$n=3$ gives that
\begin{equation}\label{kjncld cvnkGDBndov}
\begin{split}&
\frac{1}{|x-y|}=\frac1{\sqrt{|x|^2+|y|^2-2x\cdot y}}=
\frac1{|x|\,\sqrt{1-2t\tau+\tau^2}}\\&\quad
=\frac1{|x|}\,\sum _{k=0}^{+\infty }P_k(t)\,\tau^k
=\sum _{k=0}^{+\infty }P_k\left( \frac{x\cdot y}{|x|\,|y|}\right)\,
\frac{|y|^k}{|x|^{k+1}}
.\end{split}\end{equation}
Also, the additional identity for Legendre polynomials
in Theorem~\ref{CARRLegendre} says that
\begin{eqnarray*}&& P_k\left( \frac{x\cdot y}{|x|\,|y|}\right)
=\frac{{\mathcal{H}}^{2}(\partial B_1)}{N_k} \,F_k\left( \frac{x}{|x|},\frac{y}{|y|}\right)=\frac{4\pi}{2k+1}\sum_{j=1}^{2k+1} Y_{k,j}\left( \frac{x}{|x|}\right)\,Y_{k,j}\left( \frac{y}{|y|}\right)
,\end{eqnarray*}
where~$\{Y_{k,1},\dots,Y_{k,2k+1}\}$ is a basis of
spherical harmonics
that are orthonormal with respect to the scalar product in~$L^2(\partial B_1)$.
This and~\eqref{kjncld cvnkGDBndov} lead to
\begin{equation}\label{MNSX cblastSGrAVSUgbvascgO}
\frac{\mu}{|x-y|}=
\sum _{{k\ge0}\atop{1\le j\le 2k+1}}
\frac{4\pi\mu}{2k+1}\,
Y_{k,j}\left( \frac{x}{|x|}\right)\,Y_{k,j}\left( \frac{y}{|y|}\right)\,\frac{|y|^k}{|x|^{k+1}}.
\end{equation}
If, up to scaling, we assume the earth to be a (not necessarily homogeneous) ball
of radius
equal to~$1$,
we can have that the gravity potential at a point~$x\in\R^3\setminus B_1$
is generated by the superposition
of infinitesimal masses at~$y$, distributed according to some
measure~$d\mu_y$.

\begin{figure}
  \centering
  \includegraphics[width=.4\linewidth]{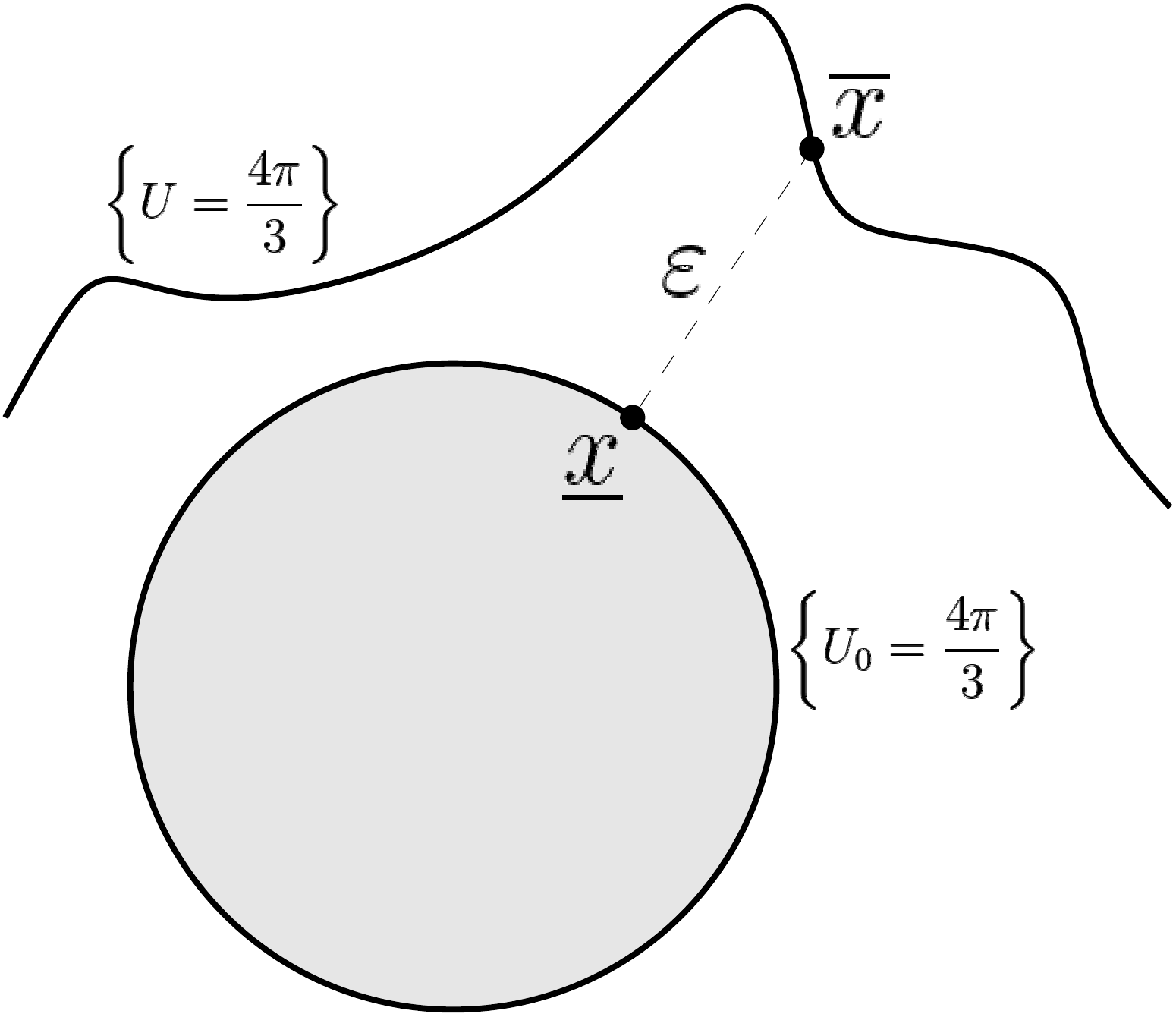}
 \caption{\sl Computing the gravity anomaly.}\label{TERRAeFI}
\end{figure}

We observe that
the homogeneous case, corresponding, up to normalization constants, to~$d\mu_y=dy$, produces
a potential at~$x\in\R^3\setminus B_1$ of the type
\begin{equation}\label{OJHSN-IOU0} U_0(x):=\int_{B_1}\frac{dy}{|x-y|}=
\frac{4\pi}{3\,|x|},
\end{equation}
see Lemma~\ref{QUESTOL}.
In general, without any homogeneity assumption, formula~\eqref{MNSX cblastSGrAVSUgbvascgO}
produces a gravity potential at~$x\in\R^3\setminus B_1$
of the form
\begin{equation}\label{PGSBCrand-1}\begin{split}
U(x)\,&:=\,\int_{B_1}\frac{d\mu_y}{|x-y|}
\\&=\,
\sum _{{k\ge0}\atop{1\le j\le 2k+1}}
\frac{4\pi}{(2k+1)\,|x|^{k+1}}\,
Y_{k,j}\left( \frac{x}{|x|}\right)\,\int_{B_1}
Y_{k,j}\left( \frac{y}{|y|}\right)\,|y|^k\,d\mu_y.\end{split}
\end{equation}
Since the mass distribution inside the earth 
is something we do not know much about,
we can formalize our uncertainty about the true value
of~$d\mu_y$ by treating this measure as a random quantity.
For instance, we could replace~$d\mu_y$ with~$dy+\mu(y,\varpi)\,dy$,
where~$dy$ takes care of the first approximation in which the earth is a uniformly distributed round ball
as in~\eqref{OJHSN-IOU0},
the term~$\mu(y,\varpi)\,dy$
is addressing the correction needed to this uniform model
and the parameter~$\varpi$ belongs to
a sample space of possible outcomes which somehow encodes
our uncertainty on the values of the mass distribution.
We will take suitable
Gau{\ss}ian
ansatz on the random variables,
which are normalized to have expected value zero and a suitable variance.
Namely, we suppose that
$$ {\mathbb{E}}\big(\mu(y,\cdot)\big)=0,$$ that is the expected value
of the mass distribution of the earth coincides with that of
the uniform density model, and
\begin{equation}\label{KMSAbsdstabtz4567mamrouP99s}
{\mathbb{E}}\big(\mu(X,\cdot)\,dX \times \mu(Y,\cdot)\,dY\big)=
\overline{\mu}(|X|)\,\delta_{0}\big(X-Y\big)\,dX\,dY,
\end{equation}
for some function~$\overline{\mu}$,
where~$\delta_0$ is the Dirac Delta Function\footnote{More precisely,
the meaning of~$\delta_{0}\big(X-Y\big)\,dX\,dY$ is that, for every~$\varphi\in C(B_1\times B_1)$,
$$ \iint_{B_1\times B_1}\varphi(X,Y)\,\delta_{0}\big(X-Y\big)\,dX\,dY
=\int_{B_1}\varphi(X,X)\,dX.$$}
at the origin. Roughly speaking, the ansatz in~\eqref{KMSAbsdstabtz4567mamrouP99s}
states that the mass distribution covariance
is expected to be uniform along each sphere~$\partial B_\rho$ with~$\rho\in(0,1)$,
and the distributions at different layers
are expected to be independent.
In this way, the gravity potential in~\eqref{PGSBCrand-1}
also becomes a random quantity
of the form~$U(x,\varpi)=U_0(x)+\xi(x,\varpi)$, with~$U_0$ corresponding to the homogeneous case in~\eqref{OJHSN-IOU0}
and
\begin{equation}\label{234MD:SDKMDJOOS-D}
\xi(x,\varpi)=
\sum _{{k\ge1}\atop{1\le j\le 2k+1}}
\frac{4\pi}{(2k+1)\,|x|^{k+1}}\,
Y_{k,j}\left( \frac{x}{|x|}\right)\,\int_{B_1}
Y_{k,j}\left( \frac{y}{|y|}\right)\,|y|^k\,\mu(y,\varpi)\,dy.
\end{equation}
It is also convenient to introduce the ``covariance function''
\begin{equation}\label{MAGI3} r(x,y):=
{\mathbb{E}}\big(\xi(x,\cdot)\,\xi(y,\cdot)\big).\end{equation}
The idea is now to take suitable ansatz on~$r(x,y)$ that can be compared
with the available data, in order to develop ``educated guesses''
on the gravity potential at places where direct measurements are not available.
For this, combining~\eqref{KMSAbsdstabtz4567mamrouP99s}, \eqref{234MD:SDKMDJOOS-D}
and~\eqref{MAGI3}, we see that
\begin{equation}\label{JS:ASDFGHJ983984SIJN}
\begin{split}
r(x,y)\,&=\,
\sum _{{{m\ge1}\atop{1\le i\le 2m+1}}\atop{{k\ge1}\atop{1\le j\le 2k+1}}}
\frac{16\pi^2}{(2m+1)\,(2k+1)\,|x|^{m+1}\,|y|^{k+1}}\,
Y_{m,i}\left( \frac{x}{|x|}\right)\,Y_{k,j}\left( \frac{y}{|y|}\right)\\&\qquad\quad\times
\iint_{B_1\times B_1}
Y_{m,i}\left( \frac{X}{|X|}\right)\,
Y_{k,j}\left( \frac{Y}{|Y|}\right)\,
|X|^m\,|Y|^k \,{\mathbb{E}} \big(\mu(X,\cdot)\,dX\times
\mu(Y,\cdot)\,dY
\big)\\&=\,
\sum _{{{m\ge1}\atop{1\le i\le 2m+1}}\atop{{k\ge1}\atop{1\le j\le 2k+1}}}
\frac{16\pi^2}{(2m+1)\,(2k+1)\,|x|^{m+1}\,|y|^{k+1}}\,
Y_{m,i}\left( \frac{x}{|x|}\right)\,Y_{k,j}\left( \frac{y}{|y|}\right)\\&\qquad\quad\times
\int_{B_1}
Y_{m,i}\left( \frac{X}{|X|}\right)\,
Y_{k,j}\left( \frac{X}{|X|}\right)\,
|X|^{m+k}\,\overline{\mu}(|X|)\,dX.
\end{split}
\end{equation}
It is also convenient to recall the orthogonality condition in~\eqref{LORZRPO} and use polar coordinates, to see that
\begin{eqnarray*}
&&\int_{B_1}
Y_{m,i}\left( \frac{X}{|X|}\right)\,
Y_{k,j}\left( \frac{X}{|X|}\right)\,
|X|^{m+k}\,\overline{\mu}(|X|)\,dX\\&=&
\int_0^1\left(
\int_{\partial B_1}
Y_{m,i}(\vartheta)\,
Y_{k,j}(\vartheta)\,
r^{m+k+2}\,\overline{\mu}(r)\,d{\mathcal{H}}^{2}_\vartheta
\right)\,dr\\&=&
\delta_{m k}\,\delta_{ij}
\int_0^1
r^{2k+2}\,\overline{\mu}(r)\,dr.
\end{eqnarray*}
Therefore, defining
\begin{equation}\label{SIGMALK} \sigma_k:=\sqrt{\int_0^1
r^{2k+2}\,\overline{\mu}(r)\,dr},\end{equation}
we can rewrite~\eqref{JS:ASDFGHJ983984SIJN} as
\begin{equation*}
r(x,y)=
\sum _{{{k\ge1}\atop{1\le j\le 2k+1}}}
\frac{16\pi^2\,\sigma_k^2}{(2k+1)^2\,|x|^{k+1}\,|y|^{k+1}}\,
Y_{k,j}\left( \frac{x}{|x|}\right)\,Y_{k,j}\left( \frac{y}{|y|}\right).
\end{equation*}
Hence, using again the 
additional identity for Legendre polynomials
in Theorem~\ref{CARRLegendre} and recalling the setting in~\eqref{ASSDM},
\begin{equation}\label{DAVDkjsdstudkg03}
r(x,y)=
\sum _{{{k\ge1}}}
\frac{4\pi\,\sigma_k^2}{(2k+1)\,|x|^{k+1}\,|y|^{k+1}}\,P_k\left( \frac{x}{|x|}\cdot \frac{y}{|y|}\right).
\end{equation}
The values of~$\sigma_k$, which are still to be determined, can be inferred by the data of the measurements
of the gravity anomalies\footnote{In the classical physical geodesy,
the gravity anomaly at a point~$\overline{x}$
is computed as the difference between the intensity of the actual
gravity force field at~$\overline{x}$
and the ideal one at the point~$\underline{x}$ which is the projection of~$\overline{x}$
on the surface~$\{U_0=U(\overline{x})\}$, omitting the random dependence on~$\varpi$ for notational simplicity,
see Figure~\ref{TERRAeFI}.
Concretely, in view of~\eqref{OJHSN-IOU0},
let us consider, for instance~$\overline{x}\in \left\{U=\frac{4\pi}3\right\}$,
$\underline{x}\in\partial B_1=\left\{U_0=\frac{4\pi}3\right\}$ and suppose that~$\e:=|\overline{x}-\underline{x}|$ is small.
Then the gravity anomaly~${\mathcal{A}}$ at~$\overline{x}$ is, expanding in
the perturbation parameters~$\e$ and~$\xi$,
\begin{eqnarray*}&&{\mathcal{A}}(\overline{x}):=
|\nabla U(\overline{x})|-|\nabla U_0(\underline{x})|
=|\nabla U_0(\overline{x})+\nabla\xi(\overline{x})|-|\nabla U_0(\underline{x})|
\\&&\qquad=\left|-\frac{4\pi\,\underline{x}}{3(1+\e)^2}
+\nabla\xi(\overline{x})
\right|-\frac{4\pi}{3}=
\frac{4\pi}{3}\left(
\left|-\frac{\underline{x}}{(1+\e)^2}
+\frac{3\nabla\xi(\overline{x})}{4\pi}
\right|-1
\right)
\\&&\qquad=
\frac{4\pi}{3}\left(
\sqrt{
\frac{1}{(1+\e)^4}+\frac{9|\nabla\xi(\overline{x})|^2}{16\pi^2}-
\frac{3\nabla\xi(\overline{x})\cdot\underline{x}}{2(1+\e)\pi}
}-1
\right)
\simeq
\frac{4\pi}{3}\left(-2\e-\frac{3\nabla\xi(\overline{x})\cdot\underline{x}}{4\pi}
\right)\simeq
\frac{4\pi}{3}\left(-\frac{3\nabla\xi(\overline{x})\cdot\underline{x}}{4\pi}-\frac{2\e}{|\overline{x}|}
\right).
\end{eqnarray*}
Since
$$ \frac{4\pi}{3}=U(\overline{x})=U_0(\overline{x})+\xi(\overline{x})=
\frac{4\pi}{3(1+\e)}+\xi(\overline{x})\simeq
\frac{4\pi(1-\e)}{3}+\xi(\overline{x}), $$
we find that~$\frac{4\pi\e}{3}\simeq\xi(\overline{x})$
and accordingly the gravity anomaly in spherical approximation is outlined by
$$ {\mathcal{A}}(\overline{x})\simeq-\nabla\xi(\overline{x})\cdot\underline{x}-\frac{2\xi(\overline{x})}{|\overline{x}|}
$$
The identity
$$ -\partial_h \xi-\frac{2\xi}{h} ={\mathcal{A}},$$
where~$h$ is the vertical displacement from the center of the earth
is sometimes considered as the ``fundamental equation of physical geodesy
in spherical approximation''
and can be concretely
used to infer the values of the covariance function in terms of the observed data,
see e.g.~\cite[formulas~(3.1) and~(3.2)]{MR414029}.

Interestingly, since~$\xi$ is harmonic outside~$B_1$,
the fundamental equation of physical geodesy in spherical approximation
along~$\partial B_1$ can be considered as a (Robin type) boundary
condition for an elliptic problem.

Let us remark that with the use of GPS for gravity measurements, it is expected
that the use of the gravity anomaly will be replaced by that of ``gravity disturbance''
(with the simplification
of one term in the fundamental equation, thus reducing
to a Neumann type boundary condition,
see e.g.~\cite[pages~93 and~95]{PGsfdgj576E}
for further details on this technical point).

We also stress that, for the sake of simplicity, in our notes we took the liberty of following a purely
spherical approximation rather than an ellipsoidal approximation
(though of course this approach is too crude for accurate geodesic analysis, see~\cite[page~97]{PGsfdgj576E}
for explicit warnings about the careless use of this approximation).} (that are the
scalar discrepancies between the measured gravity accelerations and the ones
predicted by the homogeneous model)
at different points~$x_1,\dots,x_N$.

Comparing with the observations, it was proposed in~\cite{MR0478314}
to consider the values
\begin{equation}\label{DAVDkjsdstudkg04}
\sigma^2_k:=\begin{dcases}
0 & {\mbox{ if }}k\in\{1,2\},\\
\frac{2k+1}{4\pi\,(k-1)(k-2)}& {\mbox{ if }}k\ge3,\end{dcases}\end{equation}
up to normalizing constants that we omit here for the sake of simplicity.
Remarkably, this choice, which appeared to be in good agreement with the measures,
also provides a significant mathematical simplification,
reducing the series in~\eqref{DAVDkjsdstudkg03}
to a simpler closed form. Indeed, replacing~\eqref{DAVDkjsdstudkg03}
with the choice in~\eqref{DAVDkjsdstudkg04} one obtains
\begin{equation}\label{JS:3456789JHSmegjSieldilefni-PRE}
r(x,y)=
\sum _{{{k\ge3}}}
\frac{1}{(k-1)(k-2)\,|x|^{k+1}\,|y|^{k+1}}\,P_k\left( \frac{x}{|x|}\cdot \frac{y}{|y|}\right),
\end{equation}
and this has a closed expression since the quantity
\begin{equation}\label{CKLOSMNDGFOR}
\sum _{{{k\ge3}}}
\frac{\tau^{k}}{(k-1)(k-2)}\,P_k(t)
\end{equation}
can be obtained by exploiting the generating function definition\footnote{Namely,
according to~\eqref{TRVSBlatsds2}, we
notice also that~$\mu_k=1$ when~$n=3$ and therefore
$$ \frac1{\sqrt{1-2t\tau+\tau^{2}}}
=\sum _{k=0}^{+\infty }P_k(t)\,\tau^k.$$
Consequently, recalling that~$P_0(t)=1$, $P_1(t)=t$
and~$P_2(t)=\frac12(3t^2-1)$ (see e.g. the table on page~\pageref{AEDSFesakns3s24sa098}),
one can obtain a closed form expression for~\eqref{CKLOSMNDGFOR}
by algebraic manipulations and
integrating in~$\tau$.}
of
Legendre polynomials in~\eqref{TRVSBlatsds}.
\medskip

It is also interesting to remark that a precise knowledge of the gravitational potential
can be inferred by that of the covariance function according to the following
reproducing kernel procedure.
We use the notation~$r^y(x):=r(x,y)$ and define
$$ \langle r^y, r^z\rangle:=r(y,z).$$
We extend this ``scalar product'' by linearity, namely,
for every~$\alpha_1,\dots,\alpha_N,\beta_1,\dots,\beta_N\in\R$,
$$ \left\langle \sum_{i=1}^N \alpha_i r^{x_i},\,\sum_{j=1}^N \beta_j r^{x_j}\right\rangle:=
\sum_{1\le i,j\le N}\alpha_i \beta_j
\langle r^{x_i}, r^{x_j}\rangle=
\sum_{1\le i,j\le N}\alpha_i \beta_j
r(x_i,x_j).$$
With this, we can introduce a notion of ``norm'' by setting
$$ \left\| \sum_{i=1}^N \alpha_i r^{x_i}\right\|:=
\sqrt{
\left\langle \sum_{i=1}^N \alpha_i r^{x_i},\,\sum_{j=1}^N \alpha_j r^{x_j}\right\rangle
}=\sqrt{
\sum_{1\le i,j\le N}\alpha_i \alpha_j
r(x_i,x_j)
}.$$
The advantage of this setting is that it preserves
the (expected value of the) size of the original random variable, in the sense that
\begin{equation*}\begin{split}&
{\mathbb{E}}\left(
\left|\xi(x,\varpi)-
\sum_{i=1}^N \alpha_i\,\xi(x_i,\varpi) 
\right|^2\right)-
\left\|r^x-\sum_{i=1}^N \alpha_i r^{x_i}\right\|^2
\\&\qquad=
{\mathbb{E}}\left(
\xi(x,\varpi) \,\xi(x,\varpi)+
\sum_{1\le i,j\le N} \alpha_i \alpha_j\,
\xi(x_i,\varpi) \,\xi(x_j,\varpi)-
2\sum_{i=1}^N\alpha_i\, \xi(x,\varpi)\,\xi(x_i,\varpi)
\right)\\&\qquad\quad\qquad-
r(x,x)-\sum_{1\le i,j\le N} \alpha_i \alpha_j r(x_i,x_j)
+2\sum_{i=1}^N\alpha_i r(x,x_i)\\&\qquad
=0.
\end{split}
\end{equation*}
Therefore, given~$x,x_1,\dots,x_N$, 
\begin{equation}\label{MI7}\begin{split}
&{\mbox{the coefficients~$\alpha_1,\dots,\alpha_N$
minimize
the expected value of }}\left|\xi(x,\varpi)
-\sum_{i=1}^N \alpha_i\,\xi(x_i,\varpi) 
\right|^2\\
&{\mbox{if and only if they minimize }}\left\|r^x-\sum_{i=1}^N \alpha_i r^{x_i}\right\|^2.
\end{split}
\end{equation}
In our setting, the points~$x_1,\dots,x_N$ correspond to the sites where the measurements
have been taken (say, for the sake of simplicity, the values~$\xi(x_1,\varpi)$ are known),
while the location~$x$ is inaccessible to experiments and the value of~$\xi(x,\varpi)$ has to be inferred from
the available data.

The minimal~$\alpha_i$ in~\eqref{MI7}
can be determined by
a first variation of the problem, that is, fixing~$j\in\{1,\dots,N\}$, the minimal~$\alpha_i$ satisfy
\begin{eqnarray*}
0&=&\lim_{\e\to0}\frac1{2\e}\,\left[
\left\|r^x-\sum_{i=1}^N (\alpha_i+\e\delta_{ij}) r^{x_i}\right\|^2
-\left\|r^x-\sum_{i=1}^N \alpha_i r^{x_i}\right\|^2
\right]
\\&=& -r(x,x_j)+
\sum_{i=1}^N \alpha_i r(x_j,x_i).
\end{eqnarray*}
Supposing that the symmetric
square matrix~$\{ r(x_j,x_i)\}_{1\le i,j\le N}$ is invertible,
and denoting by~$r^{ij}$ its inverse, we thus
obtain that the minimal coefficients satisfy
\begin{equation}\label{ASIN} \alpha_k=\sum_{j=1}^N r^{jk}\,r(x,x_j).\end{equation}
Accordingly, recalling~\eqref{MI7},
one takes as ``best approximation''~$\xi_\star(x)$, for the value of~$\xi(x,\varpi)$
the quantity~$ \sum_{k=1}^N \alpha_k\,\xi(x_k,\varpi) $ with~$\alpha_k$ as in~\eqref{ASIN}.
As a result, if the observed values of the potential
are denoted by
\begin{equation}\label{EMI}m_i:=\xi(x_i,\varpi),\end{equation} we set
$$ b_j:=\sum_{k=1 }^N
r^{jk}\,m_k$$
and we find that
$$ \xi_\star(x)=\sum_{1\le j,k\le N }
r^{jk}\,r(x,x_j)\,\xi(x_k,\varpi) =\sum_{1\le j,k\le N }
r^{jk}\,r(x,x_j)\,m_k=\sum_{j=1}^N b_j\,r(x,x_j).$$
The problem of an ``educated guess'' on~$\xi$ at the point~$x$ is therefore
reduced to that of determining with a sufficiently good approximation the
covariance function matrix and its inverse, for which explicit expressions
can be of great use.\medskip

Following~\cite{MR414029}, it is also interesting to observe that, as a byproduct of
the choice of the coefficients of the covariance function
in~\eqref{DAVDkjsdstudkg04},
the measurements of the gravity anomaly seem to indicate a mass distribution within the earth which becomes more and more irregular as one approaches the center. To check this intriguing statement,
one can consider the case in which
\begin{equation}\label{dd-023-929-mea-KSM}
\overline{\mu}(r):= \frac1{2\pi}\,\left(\frac5{r^7}-\frac3{r^5}\right)\end{equation}
and observe that, by~\eqref{SIGMALK}, one obtains in this situation
\[
\sigma_k^2= \frac1{2\pi}\, \int_0^1\left( 5 r^{2k-5}- 3 r^{2k-3}\right)\,dr=
\frac1{2\pi}\,\left( \frac{5}{2(k-2)}-\frac{3}{2(k-1)}\right)=\frac{2k+1}{4\pi(k-1)(k-2)}
\]
when~$k\ge3$, which coincides precisely with the choice in~\eqref{DAVDkjsdstudkg04}.

Namely, the covariance choice in~\eqref{DAVDkjsdstudkg04} that nicely
fits the experimental data is closely related to a covariance measure as in~\eqref{dd-023-929-mea-KSM} which is singular as~$r\searrow0$.
This is an indication of the possibility that
\begin{equation}\begin{split}\label{dd-023-929-mea-KSM-POSTO}&{\mbox{the variation of the mass distribution is significantly
different between}}\\&{\mbox{the surface of the earth and
in the vicinity of its center.}}\end{split}\end{equation}
\medskip

A slightly different approach to the problem has been sketched
in~\cite[Section~B]{MR927099}, using a purely ``Bayesian approach''.
Roughly speaking, one takes a development of the gravity potential 
on the surface of the earth of the form
$$ \xi(x,\varpi)=\sum_{k=0}^{+\infty}c_k(\varpi)\, u_k(x),$$
for suitable spherical harmonics~$u_k$ and random functions~$c_k$
that are assumed to be Gau{\ss}ian, with zero mean and variance~$\varsigma_k^2$, namely
\begin{equation}\label{EXAHNDcchir}
{\mathbb{E}}(c_k)=0\qquad{\mbox{and}}\qquad{\mathbb{E}}(c_ic_j)=\delta_{ij}\varsigma_i^2.
\end{equation}
This set of assumptions
can be seen as a ``prior'', somewhat based on a subjective choice,
from which, in light of the~$N$ known measurements, one can obtain a ``posterior''
estimating the potential by Bayes' Law, which in turn
computes the probability of an event, based on some
knowledge of other related events. Since the inference is based on actual data and the predictions
can be confronted with newly collected measurements, one has some control on the reliability of this method.
The core of this is the idea that
the probability that two events, say~$A$
and~$B$,
have both occurred is the ``conditional probability of~$A$
given~$B$''
(that is, the probability that~$A$
has occurred, given that~$B$ has occurred)
multiplied by the probability that~$B$ has occurred.
That is, denoting by~$P(A)$ and~$P(B)$ the probabilities\footnote{Notice that,
exchanging the roles of~$A$ and~$B$ in~\eqref{0oijnhefv3t83848PA},
$$ P(A\cap B)=P(B|A)\,P(A),$$
leading to
\begin{equation*}
P(A|B)=\frac{P(B|A)\,P(A)}{P(B)}.
\end{equation*}
See e.g.~\cite[pages 63-70]{MR2352885} for an accessible
introduction to Bayesian statistics.
In our setting, the left hand side of this identity
represents the posterior probability (that is the estimated range
of values for the gravitational potential at a point~$x$ given the
evidence of the measurements at~$x_1,\dots,x_N$) and the conditional
probability on the right hand side
can be considered as a ``likelihood function''.} of~$A$
and~$B$ respectively, and by~$P(A|B)$ the
conditional probability of~$A$
given~$B$,
\begin{equation}\label{0oijnhefv3t83848PA} P(A\cap B)=P(A|B)\,P(B).\end{equation}
Namely, given~$\delta>0$
and~$v,v_1,\dots,v_N\in\R$, we set~$I:=(v-\delta,v+\delta)$,
$I_1:=(v_1-\delta,v_1+\delta),\dots,I_N:=(v_N-\delta,v_N+\delta)$,
we take~$A$ as the event ``the value of the potential at~$x$
lies in~$I$'' and~$B$ as the event
``for all~$i\in\{1,\dots,N\}$ the value of the potential at~$x_i$
lies in~$I_i$''. 
We also recall that the linear combination of independent
Gau{\ss}ian random variables is
also Gau{\ss}ian (see e.g.~\cite[Example~2.46 and page~73]{ROSS5})
and in particular, by~\eqref{EXAHNDcchir},
the random variable~$\xi(x,\omega)$ is Gau{\ss}ian 
with zero mean and variance
$$S(x):=\sum_{k=1}^N( u_k(x))^2\varsigma_i^2.$$
We stress that this quantity is known in dependence on
the prior variances~$\varsigma_i$ and on the point~$x$ (since the spherical harmonics
are given and we suppose we can compute their values at
a given point with sufficient accuracy).
That is, if, for a given~$i$ we denote by~$B_i$ the event ``the value of the potential at~$x_i$
lies in~$I_i$'', we have that~$B=B_1\cap\dots \cap B_N$.
We consider the joint probability density function~$\gamma_{B_1,\dots,B_N}$
(which, under nondegeneracy assumptions, is explicitly known given the covariances of~$\xi(x_1,\omega),\dots,\xi(x_N,\omega)$,
see e.g.~\cite[Theorem~16.1]{MR1956867},
\cite[formula~(2.27) in Chapter~2]{MR2118904}, or~\cite[formula~(5.7)]{MR2484222}) thus obtaining
\begin{equation}\label{89ke:LAOKS1} P(B)=\int_{I_1\times I_N} \gamma_{B_1,\dots,B_N}(\varpi_1,\dots,\varpi_N)\,d\varpi_1\dots d\varpi_N.\end{equation}
Similarly, considering the joint probability density function~$\gamma_{A,B_1,\dots,B_N}$,
\begin{equation}\label{89ke:LAOKS2} P(A\cap B)=\int_{I\times I_1\times I_N} \gamma_{A,B_1,\dots,B_N}(\varpi_0,\varpi_1,\dots,\varpi_N)\,d\varpi_0\dots d\varpi_N.\end{equation}
In this way, one can exploit~\eqref{0oijnhefv3t83848PA}, \eqref{89ke:LAOKS1}
and~\eqref{89ke:LAOKS2}
to determine the probability~$
P(A|B)$ for the gravity potential at the point~$x$ to lie in the given range~$I$ of values
given the known measures 
that determine the potential at points~$x_1,\dots,x_N$
in the intervals~$I_1,\dots,I_N$.

\end{appendix}

\section*{Acknowledgments}

The authors are members of AustMS and INdAM.
This work was supported by the Australian Research Council DECRA DE180100957 and the Australian
Laureate Fellowship FL190100081.
It is a pleasure to thank Dario Bambusi, Giancarlo Benettin, Hardy Chan, H\'{e}ctor Chang-Lara,
Lorenzo D'Ambrosio, Alberto Farina, Giovanni Gallavotti, Massimo Grossi,
Damon Haddon,
Sergei Kuzenko,
Steffen Lauritzen, Jun Lim, Andrew Munyard,
Gopalan Nair, David Perrella, David Pfefferl\'e, Giorgio Poggesi, Fabio Pusateri, Alfonso Sorrentino,
Thomas Stemler, Jack Thompson, Dami\`a Torres-Latorre,
Josh Troy, Callum Vukovich,
and Jingshi Xu
for interesting comments.

\printindex

\end{document}